\documentclass[UTF8,openright,12pt,a4paper,fontset=fandol,scheme=plain]{ctexbook}
\usepackage{amsmath,amssymb,amsthm,mathtools,mathrsfs}
\usepackage{geometry,multicol,booktabs,array,longtable,enumitem,needspace,emptypage}
\usepackage{tikz}
\usepackage{placeins}
\usetikzlibrary{calc,arrows.meta}
\usepackage{fancyhdr,microtype,xcolor}
\usepackage[unicode,bookmarksnumbered,hidelinks]{hyperref}
\usepackage{bookmark}
\definecolor{CoverPaper}{HTML}{FFFAF2}
\definecolor{CoverMint}{HTML}{ECF7F3}
\definecolor{InkSlate}{HTML}{253F4D}
\definecolor{Jade}{HTML}{19786D}
\definecolor{Clay}{HTML}{AE634A}
\definecolor{MutedSlate}{HTML}{59716E}
\colorlet{inkblue}{Jade}
\xeCJKsetup{
 CJKecglue={\hskip 0.39em plus 0.08em minus 0.02em\relax},
 xCJKecglue=false}
\ctexset{
 chapter={format=\raggedright\fontsize{26}{34}\selectfont\bfseries\color{Jade},
   nameformat=\fontsize{16}{22}\selectfont\bfseries\color{MutedSlate},
   aftername={\par\nobreak\vspace{8pt}},
   beforeskip=18pt,afterskip=28pt},
 section={format=\fontsize{17}{23}\selectfont\bfseries\color{Jade},beforeskip=22pt,afterskip=12pt},
 subsection={format=\fontsize{14}{20}\selectfont\bfseries\color{InkSlate},beforeskip=16pt,afterskip=9pt}}
\hypersetup{bookmarksdepth=2}
\hypersetup{pdftitle={Finite Fields: Root Sets, Frobenius, and Algebraic Structure},pdfauthor={Hongfeng Wu and Li Zhu},pdfsubject={Root sets, orbits, duality, and operator structure in finite fields; Version 153}}
\author{Hongfeng Wu and Li Zhu}

\fancypagestyle{plain}{\fancyhf{}\fancyfoot[C]{\thepage}}
\setlist[enumerate]{label=\textup{(\arabic*)},leftmargin=2.1em,itemsep=0.35em,topsep=0.35em}
\setlist[itemize]{leftmargin=2em,itemsep=0.3em,topsep=0.35em}
\allowdisplaybreaks[2]
\newtheoremstyle{cnthm}{0.65em}{0.55em}{\normalfont}{}{\bfseries}{}{0.7em}{}
\theoremstyle{cnthm}
\newtheorem{theorem}{Theorem}[chapter]
\newtheorem{lemma}[theorem]{Lemma}
\newtheorem{proposition}[theorem]{Proposition}
\newtheorem{corollary}[theorem]{Corollary}
\newtheorem{definition}[theorem]{Definition}
\newtheorem{example}[theorem]{Example}
\newtheorem{remark}[theorem]{Remark}
\newtheoremstyle{cnexercise}{0.3em}{0.3em}{\normalfont}{}{\bfseries}{}{0.7em}
{\thmname{#1}\thmnumber{ \hyperlink{solution.\theexercise}{#2}}\thmnote{ {\normalfont(#3)}}}
\theoremstyle{cnexercise}
\newtheorem{exercise}{Exercise}[chapter]
\theoremstyle{cnthm}

\newenvironment{solution}{\par\noindent\textbf{Solution}\quad}{\par\medskip}
\newcommand{\solutionheading}[1]{\par\Needspace{4\baselineskip}\noindent\hypertarget{solution.\getrefnumber{#1}}{\textbf{Exercise~\ref{#1}}}\par\nobreak}
\newcommand{\F}{\mathbb F}
\newcommand{\Z}{\mathbb Z}
\newcommand{\N}{\mathbb N}
\newcommand{\Q}{\mathbb Q}
\newcommand{\C}{\mathbb C}
\newcommand{\Aut}{\operatorname{Aut}}
\newcommand{\Fix}{\operatorname{Fix}}
\newcommand{\ord}{\operatorname{ord}}
\newcommand{\lcm}{\operatorname{lcm}}
\newcommand{\Char}{\operatorname{char}}
\newcommand{\id}{\operatorname{id}}
\newcommand{\Irr}{\operatorname{Irr}}
\newcommand{\Tr}{\operatorname{Tr}}
\newcommand{\Nm}{\operatorname{N}}
\newcommand{\End}{\operatorname{End}}
\newcommand{\im}{\operatorname{im}}
\newcommand{\rank}{\operatorname{rank}}
\newcommand{\rad}{\operatorname{rad}}
\newcommand{\GR}{\operatorname{GR}}
\numberwithin{equation}{chapter}
\AddToHook{env/theorem/before}{\par\Needspace{6\baselineskip}}
\AddToHook{env/lemma/before}{\par\Needspace{6\baselineskip}}
\AddToHook{env/proposition/before}{\par\Needspace{6\baselineskip}}
\AddToHook{env/corollary/before}{\par\Needspace{6\baselineskip}}
\AddToHook{env/definition/before}{\par\Needspace{6\baselineskip}}
\AddToHook{env/example/before}{\par\Needspace{6\baselineskip}}
\AddToHook{env/remark/before}{\par\Needspace{6\baselineskip}}
\AddToHook{env/exercise/before}{\par\Needspace{6\baselineskip}}

\makeatletter
\newcommand*\l@chapterunit[2]{%
  \addpenalty{-200}\addvspace{0.65em}%
  \begingroup\small\bfseries\color{MutedSlate}%
  \noindent\hspace*{1.5em}#1\par\endgroup\nobreak}
\providecommand*\toclevel@chapterunit{1}
\makeatother
\newcommand{\VersionDate}{6 November 2026}
\newcommand{\CompilationDate}{%
\number\day\space\ifcase\month\or January\or February\or March\or April\or May\or June\or July\or August\or September\or October\or November\or December\fi\space\number\year}

\makeatletter
\renewcommand{\@pnumwidth}{2.8em}
\renewcommand{\@tocrmarg}{3.8em}
\makeatother

\newcommand{\GuideRef}[1]{\ref{#1} (p.~\pageref{#1})}
\newcommand{\GuideSection}[2]{\par\Needspace{7\baselineskip}\section*{#1}\phantomsection\label{#2}\addcontentsline{toc}{section}{#1}}

\begin{document}
\hypersetup{pageanchor=false}
\begin{titlepage}
\begin{tikzpicture}[remember picture,overlay]
\fill[CoverPaper] (current page.south west) rectangle (current page.north east);
\fill[CoverMint] ($(current page.north east)+(-2.25cm,0)$) rectangle (current page.south east);
\draw[Jade,line width=1.1pt] ($(current page.north west)+(3cm,-3.7cm)$) -- ++(10.7cm,0);
\end{tikzpicture}
\vspace*{23mm}
\noindent{\color{MutedSlate}\large FINITE FIELDS\par}
\vspace{20mm}
\noindent{\fontsize{48}{60}\selectfont\bfseries\color{InkSlate}Finite Fields\par}
\vspace{13mm}
\noindent{\LARGE\color{Jade}Root Sets, Frobenius, and Algebraic Structure\par}
\vspace{15mm}
\vspace{18mm}
\noindent{\Large\color{InkSlate}Hongfeng Wu and Li Zhu\par}
\vspace{5mm}
\noindent{\color{MutedSlate}Compiled on \VersionDate\par}
\vfill
\vspace{22mm}
\end{titlepage}
\thispagestyle{empty}
\null\vfill
\noindent\begin{minipage}{0.88\textwidth}
\setlength{\parindent}{0pt}
{\large\bfseries\color{Jade}Edition and Contact Information\par}
\medskip
{\normalsize\itshape Finite Fields: Root Sets, Frobenius, and Algebraic Structure}\\[0.6em]
Hongfeng Wu and Li Zhu\\
\medskip

This English version of the book was translated from the Chinese version. The Chinese version may contain additional 
content and be updated more promptly. For further details, please refer to the Chinese version of the book.
\\[0.4em]
Compiled on \CompilationDate

\textbf{Hongfeng Wu}\\
College of Science\\
North China University of Technology\\
Beijing, China\\[0.3em]
Email:
\href{mailto:whfmath@gmail.com}
     {\nolinkurl{whfmath@gmail.com}}

\bigskip

\textbf{Li Zhu}\\
School of Mathematical Sciences\\
Guizhou Normal University\\
Guiyang, China\\[0.3em]
Email:
\href{mailto:lizhumath@pku.edu.cn}
     {\nolinkurl{lizhumath@pku.edu.cn}}

\end{minipage}
\vspace*{18mm}
\clearpage
\frontmatter
\hypersetup{pageanchor=true}
\chapter*{Preface}
\addcontentsline{toc}{chapter}{Preface}
\markboth{Preface}{Preface}

Polynomial roots can form a field, while an orbit's length can determine an irreducible polynomial's degree.
Both phenomena arise from Frobenius \(x\mapsto x^p\) in characteristic \(p\), which preserves addition and multiplication and organizes algebraic elements by their periods.
This book develops finite-field theory and its constructions through root sets, conjugate orbits, and linear actions.

Existence comes first. After finitely many root adjunctions, \(X^{p^n}-X\) splits.
The Frobenius identities give the root set a field structure; degree and the absence of repeated roots determine its size.
Neither an algebraic closure nor irreducible polynomials of prescribed degree are assumed.
Conversely, every finite subfield \(S\) satisfies
\[
\prod_{a\in S}(X-a)=X^{|S|}-X.
\]
This identity expresses closure under the field operations and a common period.

Over \(k=\F_q\), relative Frobenius \(\sigma:x\mapsto x^q\) gives conjugacy through its orbits, degrees through their lengths, and subfields through its iterated fixed points. Extending the base field changes the orbit step size.
As a \(k\)-linear operator, \(\sigma\) becomes a cyclic shift in a normal basis.
The root bound excludes low-degree relations, providing a common tool for multiplicative cyclicity, nonvanishing of trace, and existence of cyclic vectors.

These descriptions retain different information. Orbit length determines field degree but not the dimension spanned by conjugates; distinct roots recover monic squarefree polynomials but not multiplicities.
Vandermonde matrices require distinct nodes, whereas Moore matrices require linear independence over the base field.
These distinctions identify the additional hypotheses needed when passing from orbits to linear relations or from root sets to quotient algebras.

Descent asks when data over an extension come from the base field.
Frobenius stability gives coefficient descent for root polynomials; equivariant functions recover the corresponding reduced quotient algebra.
Cyclic codes' Fourier coordinates satisfy the same equivariance condition, so their spectral components, organized by cyclotomic cosets, use the structure of root-set descent.
Compatibility of local restrictions and preservation of relations under coordinate changes require examination in the relevant spaces.

Trace connects multiplication with linear duality; linearized polynomials and cyclic modules connect normal bases, subspaces, and recurrences.
A subfield's minimum representation degree can depend on the power basis. Relative trace turns low-degree representation conditions into orthogonality within the subfield; short trace recurrences recover root-block polynomials, connecting representations with functional decomposition.
For repeated factors, coefficient pairings and derivatives supply information lost by trace pairings.
Base-field translation of a generator uses addition in the fixed subspace; reciprocation combines traces and defining equations to produce new normality conditions.

Construction and counting ask how local data determine a whole structure.
Return maps express total orbit lengths as products of base-orbit lengths and internal return periods, recovering factorizations and supporting compression and irreducible constructions. Recursion also requires preservation of the next step's hypotheses.
Character orthogonality turns trace and multiplicative conditions into sums. Gaussian-period moments record short additive relations in multiplicative subgroups, while joint additive and multiplicative character expansions connect primitivity with normality.
Lifting root orbits to $p$-adic precision, Witt carries, coefficient valuations and Newton slopes recover information lost upon reduction; equal-characteristic formal series retain coefficient-field operations while equations and roots vary with a parameter. In the reverse direction, residue fields of unramified cyclotomic extensions recover finite fields; integral relations among the same roots connect the root-set construction with higher arithmetic precision.
Extension-field character sums and curve point counts return to finite-order recurrences; syndrome recurrence recovery turns code root constraints into error location and codeword recovery.

Exceptional parameters in reciprocal families, degree obstructions in subfield representations, preservation of return conditions, and empty trace fibers turn one criterion's remaining conditions into further questions.
Examples and exercises compare representations, test hypotheses, and trace subsequent uses of the same structure.

The book assumes abstract algebra and linear algebra, and addresses advanced undergraduates, graduate students, and researchers.
Basic structure theorems and main constructions are proved; elementary tools appear in the appendices. General references are \cite{LN,Feng,Conrad}.

The book may still contain omissions and errors. The authors would be grateful for readers' comments and corrections.

\chapter*{Acknowledgments}
\addcontentsline{toc}{chapter}{Acknowledgments}
\markboth{Acknowledgments}{Acknowledgments}
We thank GPT-6 Astra for assistance in preparing this book, especially in translating the Chinese version into English.

\chapter*{Reading Guide and Notation}
\addcontentsline{toc}{chapter}{Reading Guide and Notation}
\markboth{Reading Guide and Notation}{Reading Guide and Notation}

\noindent\textbf{Prerequisites.}
The reader should know the basic theory of groups, rings, fields, and finite-dimensional vector spaces. Field extensions and generated subfields are defined when first used. Appendix~\ref{chap:polynomial-tools} supplies field embeddings, restriction of scalars, and polynomial tools; the other appendices may be consulted as needed.

\noindent\textbf{Reading routes.}
Begin with Chapters~\ref{chap:foundations}--\ref{chap:orders}. Then use the roots of unity and cyclotomic cosets in Chapter~\ref{chap:cyclotomic} alongside Chapters~\ref{chap:trace}--\ref{chap:normal}, which develop the field structure, trace duality, and linear tools needed later.

Chapters~\ref{chap:algorithms} and~\ref{chap:root-fibers} treat algorithms and root fibers. Chapters~\ref{chap:root-symmetry}--\ref{chap:root-constructions} continue with root-set symmetry, subfield representations, and recursive constructions. Begin Chapter~\ref{chap:subfield-representations} with its opening power-basis example.

Chapters~\ref{chap:functions}--\ref{chap:codes} treat polynomial functions, linear recurrences, and codes. Finite-observation recovery in Chapter~\ref{chap:sequences} can be read together with the decoding example in Section~\ref{v118:sec-bch-decoding}.

Chapter~\ref{chap:characters} uses trace and the multiplicative group for counting. Exact-period layers in Chapter~\ref{chap:period-layers} may also be read after the first four chapters; its character-sum and recurrence material can be revisited after the relevant chapters. The quick reference and indexes locate individual results and their hypotheses.
The thematic reading guide on page~\pageref{v122:guide} connects chapters through their tools and problems.

\par\Needspace{8\baselineskip}
\noindent\phantomsection\label{v119:par-argument-organization}\textbf{Organization of the arguments.}
The table lists the principal structures established in each part and their subsequent uses.
\begingroup\small\linespread{1.08}\selectfont
\renewcommand{\arraystretch}{1.15}
\begin{longtable}{@{}>{\raggedright\arraybackslash}p{0.10\textwidth}>{\raggedright\arraybackslash}p{0.39\textwidth}>{\raggedright\arraybackslash}p{0.42\textwidth}@{}}
\toprule
Chapters&Structure established&Uses in subsequent arguments\\
\midrule
\endhead
\bottomrule
\endfoot
\ref{chap:foundations}--\ref{chap:frobenius}&
Adjoining roots, fixed points, conjugate orbits, and equivariant functions&
Existence and uniqueness of finite fields, subfield classification, root-set descent, and reduced-algebra decompositions\\
\ref{chap:orders}, \ref{chap:cyclotomic}&
Multiplicative orders, cyclotomic cosets, orbit step sizes, and compatible lifting&
Binomial criteria, factorization over extensions, lifting factors, and root constraints for codes\\
\ref{chap:trace}--\ref{chap:normal}&
Trace and norm, trace duality, linearized root spaces, and cyclic modules&
Coordinate recovery, normal bases, subspace criteria, and trace representations of recurrences\\
\ref{chap:algorithms}, \ref{chap:root-fibers}&
Fixed algebras, reciprocal families, and fiber return&
Polynomial factorization, normality tests, and recovery of a full factorization from one fiber\\
\ref{chap:root-symmetry}&
Root multiplicities, affine stabilizers, and normal forms&
Classification of symmetric root sets, root-set partitions, and functional decomposition\\
\ref{chap:subfield-representations}&
Subfield degree profiles, root-block coefficients, and relative traces&
Low-degree representations, tests for inner factors, and recovery of root blocks from short trace prefixes\\
\ref{chap:root-constructions}&
Orbit composition, resultants, and conditions preserved recursively&
Factorization of composed products, irreducible iterates, and recursive constructions\\
\ref{chap:functions}&
Difference conditions and additive and multiplicative fibers&
Permutation criteria, fiber lifting, and structure-preserving function constructions\\
\ref{chap:sequences}&
Shift operators, quotient algebras, and finite-prefix recovery&
Periods and transients of recurrences, syndrome recurrences, and error location\\
\ref{chap:codes}&
Root constraints, evaluation spaces, and spectral coordinates&
Code distances and duality, support structure, coefficient-field descent, and error correction\\
\ref{chap:characters}&
Character orthogonality, Gauss sums, and weighted orbit products&
Trace distributions, primitive normal bases, and recurrences for extension-field character sums and curve point counts\\
\ref{chap:period-layers}&
Exact-period layers, compression, and return maps&
Separation of exact periods, lifting compressed factors, and orbit counts\\
\ref{revadic:chapter}&
Unramified cyclotomic extensions, integer rings, and reduction&
Finite-field construction, normal integral bases, coefficient valuations, and point-count slopes\\
\end{longtable}
\endgroup

\noindent\phantomsection\label{def:basic-notation}\ignorespaces \textbf{Notation.}
Throughout, \(p\) is prime and \(q=p^r\); we usually write \(k=\F_q\subseteq E=\F_{q^n}\),
\(\sigma(x)=x^q\), and \(\phi(x)=x^p\).
The notation \([E:k]\) denotes the extension degree, and \(\Aut_k(E)\) denotes the group of automorphisms fixing \(k\) pointwise.
The variable \(X\) is usually used for polynomials in field elements, and \(T\) for polynomials in operators.
All field homomorphisms preserve the identity element; uniqueness of a finite field of a given order means uniqueness of its isomorphism type.
Intersections and composita of subfields, and comparisons of root sets, are taken in a specified common extension field.
Statements allowing polynomials of degree less than a given positive integer include the zero polynomial.

Sections~\ref{v135:sec-foundations}--\ref{revadic:sec-frob} of Chapter~\ref{revadic:chapter} can be read after the first four chapters. Trace and normal integral bases connect with Chapters 6--9, and character sums with Chapter 18. Route eight gives a staged reading order.

\tableofcontents

\chapter*{Thematic Reading Guide}
\phantomsection\label{v122:guide}
\addcontentsline{toc}{chapter}{Thematic Reading Guide}
\markboth{Thematic Reading Guide}{Thematic Reading Guide}

The same structure serves different purposes in different chapters. Cyclic modules turn normality into generation and later describe shift relations in sequences; trace pairings recover coordinates and turn linear equations into restrictions on frequencies. The eight routes below start from particular problems and connect these uses. References and page numbers point to the main text. Each route can be read independently after the indicated prerequisites.

\begingroup\small\linespread{1.08}\selectfont
\renewcommand{\arraystretch}{1.18}
\begin{longtable}{@{}>{\raggedright\arraybackslash}p{0.30\textwidth}>{\raggedright\arraybackslash}p{0.35\textwidth}>{\raggedright\arraybackslash}p{0.27\textwidth}@{}}
\toprule
Topic&Starting question&Main locations\\
\midrule
\hyperref[v122:route-normal]{Normal bases and cyclic modules}&When do conjugates span the extension?&Chapters 8, 9, 16, 17\\
\hyperref[v122:route-fibers]{Root fibers and recursive irreducible constructions}&When can an irreducible lifting be repeated?&Chapters 6, 11, 14, 19\\
\hyperref[v122:route-fourier]{Trace duality and Fourier analysis}&How does linear solvability become a counting formula?&Chapters 7, 8, 18\\
\hyperref[v122:route-decoding]{Linear recurrences and codes}&How do finitely many syndromes determine errors?&Chapters 13, 16, 17\\
\hyperref[v122:route-subfields]{Subfield representations and relative trace recovery}&How early does a subfield appear in a power basis?&Chapters 7, 11, 13, 16\\
\hyperref[v122:route-reciprocals]{Reciprocal families and normality obstructions}&Which linear relations can a change of generator remove?&Chapters 8, 9, 10\\
\hyperref[v125:route-orbits]{Root exponents, phases, and coefficient lifting}&Which root-set structures persist under successive lifts?&Chapters 5, 10, 14, 17, 19\\
\hyperref[v132:route-adic]{Finite fields and unramified constructions}&How do roots, integral bases, and precision correspond?&Chapters 5--9, 18, 20\\
\bottomrule
\end{longtable}
\endgroup

We retain \(k=\F_q\subseteq E=\F_{q^n}\) and \(\sigma(x)=x^q\). Fields, conjugacy, and subfields from the first four chapters are common prerequisites; routes involving traces also use Chapter 6. The special trace sequences of Chapter 13 can first be read through their direct proofs, before the general recurrence theory of Chapter 16.

\GuideSection{1. Normal bases and cyclic modules}{v122:route-normal}
\noindent\textbf{Route:} Section 8.3, Sections 9.1--9.5 and 9.7; then Section 16.2 and Sections 17.3--17.4.

Example~\GuideRef{ex:normal-f8-example} provides a concrete starting point: an element has three distinct conjugates but is not normal because their sum is zero. Field degree records orbit length; the degree of the Frobenius annihilator \(\mu_a\) records the dimension spanned by the conjugates. Section 9.1 separates these quantities, and the normal basis theorem, Theorem~\GuideRef{thm:normal-basis}, turns the problem into finding a cyclic vector for \(\sigma\).

There are two successive steps. Section 8.3 uses the root bound to prove that the minimal polynomial of \(\sigma\) is \(T^n-1\); Section 9.2 then applies the cyclic vector lemma. Once a normal element \(\eta\) is chosen, the map
\[
k[T]/(T^n-1)\longrightarrow E,\qquad
\bar a\longmapsto a(\sigma)\eta
\]
identifies normal elements with units of the quotient ring. Proposition~\GuideRef{prop:normal-test} and the subsequent count thus use the same criterion: every primary component must be fully generated.

The same proposition gives \(U_P=\ker((F/P)(\sigma))\) and its trace dual \(V_P\). A nonzero detecting-operator value, avoidance of an obstruction kernel, and a nonzero weighted trace express the same local generation condition. Section 10.6 intersects these kernels with reciprocal families; Section 18.4 turns the same conditions into the exact Fourier support in Proposition~\GuideRef{v128:prop-local-spectrum}. The fibre-sum criterion~\GuideRef{v131:thm-fibre-normal} and component product~\GuideRef{v131:thm-component-product} then group intermediate-field normality conditions into independent components; Proposition~\GuideRef{v131:prop-component-exact} recovers their exact probabilities through joint ranks in trace duality. Corollary~\GuideRef{v131:cor-base-complete-count} classifies all base-field translates; Proposition~\GuideRef{v131:prop-intermediate-translate-count} reduces intermediate-field translation counts to affine joint ranks and image compatibility.

When \(p\mid n\), the polynomial \(T^n-1\) has repeated factors. The cyclic module remains valid, whereas individual eigenvalues do not record the nilpotent chains. Read the \(\F_4/\F_2\) example in Section 9.2 together with the invariant-subspace classification in Section 9.7. The former shows why failure of diagonalizability does not prevent a normal basis; the latter translates kernels, intersections, sums, and trace orthogonal complements into divisors, greatest common divisors, least common multiples, and reciprocal polynomials. Proposition~\GuideRef{v150:prop-row-intersection} then uses primary components to decide whether any two nonzero row spaces in the operator algebra intersect.

After coefficient extension, the relevant action on a linearized root space is the Frobenius return fixing the coefficients. Theorem~\GuideRef{v145:thm-skew-roots} identifies monic right factors with return-invariant subspaces; Proposition~\GuideRef{v145:prop-return-matrix} computes the return characteristic polynomial through a Moore matrix. The return need not be cyclic, but irreducibility of its characteristic polynomial characterizes skew irreducibility. Proposition~\GuideRef{v146:prop-spectrum-descent} extends two-dimensional returns to any Frobenius generator; Lemma~\GuideRef{v146:lem-semilinear-normal} then converts a coefficient fiber with prescribed return spectrum into a count of normalized cyclic vectors. Theorem~\GuideRef{v148:thm-all-fibers} counts all return-norm fibers through local matrix modules; Theorem~\GuideRef{v148:thm-relative-fibers} then explains fixed-constant uniformity by common subfields. To recover ordinary factorization patterns, retain the local partitions as well as the characteristic polynomial; the orbit formulas are in Theorem~\GuideRef{v148:thm-ordinary-return}.

Section 16.2 uses a quotient algebra to model a recurrence space, with \(T\) acting as the shift. In Sections 17.3--17.4, \(X\) acts as cyclic or constacyclic shift on codeword coordinates, and stable subspaces become ideals. These models explain the operators, but the normal-coordinate isomorphism is generally not a ring isomorphism. It does not turn field multiplication into coordinatewise multiplication.

For self-duality, continue with Sections 7.6 and 9.15: a self-dual basis for a nondegenerate pairing and a self-dual normal basis formed by one Frobenius orbit satisfy different constraints. For multiplicative primitivity, turn to Sections 9.16 and 18.9. The unit condition in the cyclic module has already settled normality; what remains is generation of the multiplicative group.

\GuideSection{2. Root fibers and recursive irreducible constructions}{v122:route-fibers}
\noindent\textbf{Route:} Section 11.1, Sections 14.1--14.2 and 14.4--14.6; then Sections 19.6--19.7 and 19.10. Traces and norms use Chapter 6.

For a monic irreducible polynomial \(P\), the study of \(P(h(X))\) can begin with one root \(\alpha\) of \(P\) and the smaller equation \(h(X)=\alpha\). The two return examples in Section 11.1 have the same base orbit and fiber size but different total factorizations: the difference is the permutation induced inside a fiber after one circuit of the base orbit.

Lemma~\GuideRef{lem:fiber-return} gives
\[
\text{total orbit length}=\text{base-orbit length}\times\text{least return period}.
\]
Theorem~\GuideRef{thm:fiber-factorization} recovers all factors and their multiplicities; Capelli's criterion is the case of an irreducible fiber polynomial. Fixed-point counts of the return iterates suffice for factor degrees and numbers, whereas explicit factors require the corresponding root orbits. Corollary~\GuideRef{v129:cor-return-base-change} separates base extension into splitting of the base orbit and splitting of the internal return; these are distinct operations.

For images, the same period formula records information loss. Section 14.1 first synchronizes Frobenius on root pairs and then takes their sums or products. Proposition~\GuideRef{v129:prop-phase-images} expresses the composed polynomial through the minimal polynomials of the phase images: shortening an orbit or obtaining conjugate images from different phases can produce repeated factors. Example~\GuideRef{v129:ex-composed-cubic} distinguishes both sources using the sum and product compositions of one cubic. Once phases and periods determine the factorization type, the matrices and resultants of Section 14.2 supply coefficients.

Recursion also requires the next fiber criterion to remain valid. Section 14.4 preserves a nonzero trace after reciprocation; Section 14.5 preserves a trace or square class for quadratic fibers; Section 14.6 uses norms to descend extension-field criteria to the base field. For \(h=X^2+c\) and target \(\tau\), Theorem~\GuideRef{v55:thm-quadratic-stability} requires
\[
\tau-c\text{ nonsquare},\qquad
h^{\circ j}(0)-\tau\text{ nonsquare}\quad(j\ge2).
\]
Corollary~\GuideRef{v55:cor-quadratic-finite-test} shows that repetition of the critical orbit ends the appearance of new conditions. Example~\GuideRef{v55:ex-quadratic-stability} also shows why affine conjugacy must transform the target, and why the leading coefficient of a general quadratic enters the square-class test. Trace descent must separately exclude shorter periods by a low-degree relation; membership in the target field alone does not establish generation.

Section 19.7 first compresses affine orbits and then recovers the internal period from a norm and an orbit sum of the compression value; see Proposition~\GuideRef{prop:explicit-affine-return}. Corollary~\GuideRef{v129:cor-affine-factor-orbits} identifies each compressed factor with an affine orbit of irreducible polynomials, with return order equal to the order of the affine stabilizer over the base field. Formula~\GuideRef{v129:eq-return-euler} in Section 19.10 records all return periods in an Euler product; it retains the degrees and numbers of distinct factors, but does not recover root coordinates or multiplicities.

Roots of unity provide another continuation: Corollary~\GuideRef{v124:cor-binomial-extension} determines the least binomial-splitting extension from exponent blocks, and Proposition~\GuideRef{v47:prop-galois-med} uses Hensel uniqueness to preserve composition structure at every coefficient precision.

\GuideSection{3. Trace duality and Fourier analysis}{v122:route-fourier}
\noindent\textbf{Route:} Sections 7.1--7.2, Section 8.6, Sections 18.1--18.4; then Sections 18.8 and 18.15--18.20; Section 4.7 supplies quadratic target counts.

Nondegeneracy of the trace pairing in Section 7.1 expresses every \(k\)-linear functional as \(x\mapsto\Tr_{E/k}(ax)\). Section 7.2 uses this to recover power-basis coordinates; Section 8.6 expresses dual maps as adjoint operators. For a linear map \(A:E\to E\), the relation
\[
\im A=(\ker A^*)^\perp
\]
turns solvability of \(A(x)=b\) into trace conditions on \(b\). Knowing the adjoint kernel suffices here; finding every preimage is unnecessary.

In Chapter 18, take \(K=\F_{p^r}\) and compose the absolute trace with a nontrivial complex character of \(\F_p\). Linear functionals become additive characters. In the four-element example of Section 18.2, two trace fibers have identical Fourier support, while translation changes the phase. In general,
\[
\widehat{1_{b+U}}(a)
=|U|\psi(-ab)1_{U^\perp}(a),
\]
where \(U\) is an \(\F_p\)-subspace and orthogonality is defined by the absolute trace pairing. The trace orthogonality condition has become a frequency restriction.

For an \(\F_p\)-linear map \(A:K\to K\), let \(A^*\) denote its adjoint for the absolute trace pairing. Proposition~\GuideRef{prop:fourier-adjoint} further gives
\[
\widehat{A_*F}(a)=\widehat F(A^*a),
\qquad (A_*F)(y)=\sum_{A(x)=y}F(x).
\]
Taking \(F=1\) recovers the number of solutions of a linear equation. When linear algebra already determines the fibers, this formula chiefly explains their frequency structure. Adding multiplicative-subgroup or primitivity constraints makes joint character expansions provide further counting information. Section 18.8 first determines scaling symmetries on the nonzero trace values, then expands on the corresponding quotient group, avoiding coefficients already forced to vanish by symmetry.

Normality supplies local conditions whose expansion is explicit. In Proposition~\GuideRef{v128:prop-local-spectrum}, the direct sum of dual spaces makes the frequency representation unique, preventing coefficient cancellation. Corollary~\GuideRef{v128:cor-affine-spectrum} then restricts an affine-fiber count to \(V_M\cap W^\perp\). The eight-element example~\GuideRef{v128:ex-affine-spectrum-eight} recovers the four counts from Chapter 9 and shows why independent linear coordinates need not yield independent normality conditions. Proposition~\GuideRef{v130:prop-annihilator-spectrum} extends the exact spectrum to any prescribed annihilator stratum; its support dimension is precisely the minimum amount of linear information in Chapter 9. Section 18.9 uses this expansion directly and adds multiplicative characters.

Translate intersections of the norm-one group connect additive decomposition with character-sum moments. Proposition~\GuideRef{v151:prop-normone-intersections} determines all intersection counts from the norm of the translation parameter, and Corollary~\GuideRef{v151:cor-normone-moment} gives the exact fourth moment by orthogonality. Proposition~\GuideRef{v153:prop-normone-short-sums} turns the same norm condition into short zero-sum relations.

Computing character sums from field operations also has a reverse direction: Proposition~\GuideRef{v150:prop-jacobi-recovery} recovers addition from a faithful Jacobi row, while Example~\GuideRef{v150:ex-partial-jacobi} shows that a partial table retains only joint coset counts. Propositions~\GuideRef{v150:prop-subset-products} and~\GuideRef{v150:prop-joint-products} count fibers with one or two prescribed subset products; Theorem~\GuideRef{v150:thm-ers-design} turns a nonzero Fourier coefficient into a three-design obstruction for an extended RS code.

Support alone usually does not determine function values. Equality in the uncertainty principle, Theorem~\GuideRef{v90:thm-fourier-uncertainty}, is an exception: when the product of the support sizes of a nonzero function and its Fourier transform equals the group order, the function is a nonzero scalar multiple of a coset indicator times a character. For the quadratic phases of Section 18.15, the polarization radical determines the direction and amplitude of support, while values on the radical are also needed to determine its translation.

A second counting route starts from return norms. Theorem~\GuideRef{v148:thm-relative-fibers} in Chapter~9 first determines relative-norm symmetry of constant coefficients: coprime degrees already give every fiber, without character sums. Hasse--Davenport lifting handles noncoprime quadratic returns in Theorems~\GuideRef{v145:thm-return-uniform} and~\GuideRef{v148:thm-quadratic-even}; the matching expansion~\GuideRef{v146:prop-return-matching} converts these counts into trace--norm distributions in low-degree extensions.

Arbitrary even degrees over a quadratic coefficient field require another cancellation. Lemma~\GuideRef{v148:lem-hankel-shift} makes singular Hankel strata depend only on the norm, while full-rank Hermitian sums leave a correction determined by commutative factor constants. Theorem~\GuideRef{v148:thm-all-even-fibers} thus treats irreducible, reducible, and repeated-factor targets together. Proposition~\GuideRef{v149:prop-uniformity} further reduces uniformity to divisibility conditions on the constant-coefficient phases of the factors. Proposition~\GuideRef{v148:prop-fiber-fourier} relates the same phase to an elliptic trace.

For prescribed return orders, Theorem~\GuideRef{v148:thm-even-order} gives the full spectrum. The Lucas atoms of Chapter~6 pull back in Theorem~\GuideRef{v148:thm-optimal-pullback} to optimal coefficient criteria for skew quadratics; exact cancellation at the zero parameter determines the degree. Fixing another coefficient leads instead to cubic value sets and norm circles, in Theorems~\GuideRef{v148:thm-cubic-value-set} and~\GuideRef{v148:thm-cubic-circles}. Finally, norm transfer~\GuideRef{v148:prop-sparse-transfer} reduces equally spaced powers to target selection and weighted fibers; the coprime count~\GuideRef{v148:thm-sparse-coprime} and primitivity obstruction~\GuideRef{v148:cor-sparse-order} distinguish the two goals.

The finite-field-valued Fourier coordinates of Chapter 17 provide a useful comparison, with the value fields kept distinct. Spectral decomposition of cyclic codes takes place in a suitable finite extension; counting, absolute values, and Parseval in this route use complex characters. Orthogonality organizes both theories, but their characteristic assumptions differ.

\GuideSection{4. Linear recurrences and codes}{v122:route-decoding}
\noindent\textbf{Route:} Sections 16.1--16.3 and 16.6, then Section 17.5; Section 13.4 supplies an earlier instance of Hankel recovery.

Chapter 16 first expresses recurrence of an infinite sequence as an annihilating relation for the shift, then turns to finite prefixes. Passing between them requires a bound on recurrence order. Without it, one prefix admits different continuations. The Berlekamp--Massey theorem, Theorem~\GuideRef{thm:berlekamp-massey}, finds a shortest connection relation for the prefix; Corollary~\GuideRef{cor:bm-recovery} states when that relation recovers the minimal polynomial of the whole sequence.

Section 17.5 supplies a natural order bound. Let the code length be \(N\), with \(p\nmid N\) and \(1\le2t<N\), and suppose every codeword vanishes at \(\alpha,\ldots,\alpha^{2t}\), where \(\alpha\) is a primitive \(N\)-th root of unity. If at most \(t\) errors occur in \(r=c+e\), write \(e(X)=\sum_{i=1}^w e_iX^{j_i}\) with \(e_i\ne0\) and distinct \(0\le j_i<N\). The computable syndromes are
\[
S_a=r(\alpha^a)=\sum_{i=1}^w e_i x_i^a,
\qquad x_i=\alpha^{j_i},\quad 1\le a\le2t.
\]
The root constraints eliminate the unknown codeword; the number of errors becomes the recurrence order of this exponential-sum sequence.

Proposition~\GuideRef{v118:prop-syndrome-recurrence} applies the algorithm of Chapter 16 to these data. Distinct nodes and nonzero weights make the minimal order exactly \(w\), and the recovered connection polynomial is
\[
\Lambda(Z)=\prod_{i=1}^w(1-x_iZ).
\]
Its roots identify the inverse locations \(\alpha^{-j_i}\). A Vandermonde system then gives the error values and hence \(c=r-e\). Example~\GuideRef{v118:ex-bch-fifteen-decoding} carries out two-error decoding for a binary BCH code of length fifteen, verifying the result by re-encoding and the distance bound.

The distance bound and recovery algorithm use the same evaluation equations for different ends. In the BCH bound, a few nonzero coefficients contradict a consecutive block of zeros; in decoding, a finite list of syndromes determines the positions and values of the coefficients. Vandermonde invertibility suffices for the distance bound. Recovering unknown error locations is what calls for a recurrence and a locator polynomial.

Lemma~\GuideRef{v127:lem-pairing-recovery} reduces this route to a nondegenerate product pairing; Proposition~\GuideRef{v127:prop-observation-kernel} explains a drop in order through the observation kernel. The dual-number example in Chapter 16 shows that ordinary trace can lose the entire recurrence, whereas the coefficient pairing still recovers repeated factors. In Chapter 17, the syndrome pairing is nondegenerate and diagonal in evaluation coordinates at the error locations, so its Hankel rank equals the number of errors. Each recovery problem also has limits. A general finite prefix requires the connection length as well as the polynomial, since the polynomial degree can be smaller. Decoding requires the error-count bound, and evaluation of the received polynomial at other powers need not continue the error sequence. The recovery conclusion cannot be carried beyond these data without further justification.

To ask when the distance bound is attained, Proposition~\GuideRef{v123:prop-bch-equality} identifies the kernel of a Vandermonde matrix with one row removed through interpolation weights. Corollary~\GuideRef{v123:cor-paired-bch-word} then turns equality of power-trace prefixes of two root orbits into a word attaining the bound; in odd characteristic the signs must be retained. Example~\GuideRef{v123:ex-bch-equality-small} gives complete binary and ternary examples.

Proposition~\GuideRef{v142:prop-full-fibers} also turns complete orbit fibers into irreducible composition. Corollary~\GuideRef{v142:cor-irreducible-interleaving} then gives independent short-code blocks. With only a subgroup in a proper subfield, Proposition~\GuideRef{v142:prop-subfield-distance} still bounds both the code distance and its dual distance. Example~\GuideRef{v142:ex-twenty-five-code} uses the same twenty-fifth root to obtain exact parameters $[25,20,2]$. Corollary~\GuideRef{v153:cor-normone-binary-code} obtains a binary distance bound from short zero-sum relations in a norm-one group. Example~\GuideRef{v153:ex-normone-binary-codes} then uses a recursive factor from Chapter 14 to obtain a $[17,9,5]$ code.

Theorem~\GuideRef{v140:thm-code-group} connects defining sets to the full permutation group of a canonical lift. The second digit already reduces the symmetries of all nonzero check columns from the general linear group to multiplication and Frobenius. Corollary~\GuideRef{v144:cor-subgroup-code} treats proper-subgroup check columns: separation of the mixed exponents in the first carry determines every higher-precision permutation group, whereas complete root blocks may retain additional permutations. Proposition~\GuideRef{v140:prop-common-symmetry} shows that common permutation groups and multiplier stabilizers persist at every level.

\GuideSection{5. Subfield representations and relative trace recovery}{v122:route-subfields}
\noindent\textbf{Route:} The opening example of Chapter 13 and Sections 13.1--13.4, together with Sections 7.2, 11.1, and 16.6.

Example~\GuideRef{sfr:ex-quartic-comparison} compares the same four-element subfield of \(\F_{16}\) in two power bases. Using a root of \(X^4+X+1\), a subfield generator has a quadratic representative; using its inverse, the least polynomial degree becomes three. Subfield classification has identified the object, but not its position in the low-degree power spaces.

The coefficient field of a root block can first be determined from its period. Proposition~\GuideRef{v50:prop-multiplier-blocks} does not require the multiplier group to contain Frobenius: the first return of the root block gives its exact coefficient field. In the subfield representation problem, one must further decide whether that block is a full fiber of a polynomial of low degree.

Section 13.1 records that position through \(K_d\cap\mathcal P_t(\alpha)\), where \(1<d<n\), \(d\mid n\), \(K_d=\F_{q^d}\), and \(\mathcal P_t(\alpha)=\operatorname{span}_k\{1,\alpha,\ldots,\alpha^t\}\). The minimum-degree theorem, Theorem~\GuideRef{sfr:thm-optimal-composition}, identifies attainment of the bound \([E:K_d]\) for a subfield generator with the corresponding polynomial decomposition. This leads back to Chapter 11: a subfield supplies a Frobenius-stable root partition, but realizing its blocks as full fibers of one low-degree polynomial imposes further relations on their coefficients.

Section 13.2 brings in the leading-coefficient formula of Chapter 7. Let \(f\) be the monic minimal polynomial of \(\alpha\), let \(K=K_d\), put \(m=[E:K]\), and set
\[
c_j=\Tr_{E/K}\left(\frac{\alpha^j}{f'(\alpha)}\right).
\]
The weight \(1/f'(\alpha)\) makes the orthogonal complement of a low-degree space another low-degree space; relative trace then descends the conditions to \(K\). Theorem~\GuideRef{str:thm-profile-trace-duality} recovers the degree profile from the ranks of these \(c_j\) over \(k\). For existence of a representative at a specified degree, the intersection kernel or the corresponding trace rank already provides the required information.

To recover a root-block equation, view the same \(c_j\) as a \(K\)-valued sequence and allow recurrence coefficients in \(K\). Corollary~\GuideRef{str:thm-minimal-recurrence} identifies its minimal recurrence polynomial, and Corollary~\GuideRef{str:thm-hankel-recovery} recovers all coefficients from the first \(2m\) terms. Both are direct applications of Lemma~\GuideRef{v127:lem-pairing-recovery}. The Hankel matrix is the Gram matrix of a nondegenerate weighted trace pairing; it remains invertible even if \(p\mid m\) or \(c_0=0\). Rank over the base field and recurrence over the subfield extract different information from the same data. The coefficient field belongs to the specification of the problem.

Unweighted power traces pose a further recovery question: how much information is lost when the last trace is omitted? Theorem~\GuideRef{v123:thm-power-trace-pencil} shows that, among monic square-free polynomials of degree \(m\), the fibers of the first \(2m-2\) power sums are singletons or affine lines, according to \(f''\). Example~\GuideRef{v123:ex-f3-trace-pencil} exhibits a nonsingleton fiber in odd characteristic. These data are the actual unweighted power sums, distinct from the ranks used in weighted relative-trace recovery.

If polynomial representations fail to attain the bound, Section 13.3 considers rational representations. The quartic example admits a quadratic rational representative, whereas the root-block coefficient rank prevents the same degree reduction in the sextic example. Whether a denominator helps is determined by the coefficient space; it is not an automatic remedy for failure of polynomial representation. The recovered root block therefore answers questions about representation, composition, and fibers at once.

Theorem~\GuideRef{v137:thm-jump-realization} turns degree profiles into realization conditions: once the jump positions are fixed, multiplicative closure and Frobenius fixed-point equations test the same lower coefficients. Example~\GuideRef{v137:ex-quartic-classification} lists every binary quartic defining polynomial. Degree constraints retain leading terms; deciding realizability also requires compatibility of the lower coefficients.

Theorem~\GuideRef{v140:thm-profile-smith} extends the degree profile across precisions: positive Smith exponents of a Frobenius difference determine when additional representation directions disappear. Example~\GuideRef{v140:ex-quartic-depth} returns to the opening quartic example, distinguishing canonical lifting of a subfield from persistence of its low-degree representations.

\GuideSection{6. Reciprocal families and normality obstructions}{v122:route-reciprocals}
\noindent\textbf{Route:} The translation criterion in Section 9.3 and reciprocal families in Section 10.6; then compare reciprocal rank in Section 8.2.

A change of generator preserves extension degree but may alter linear relations among conjugates. Proposition~\GuideRef{v82:prop-normal-translates} first determines the limits of translation. If \(p\nmid n\), base-field translation changes only the \(T-1\) normality condition; if \(p\mid n\), whether any of the detecting expressions vanishes is unchanged. When the obstruction lies in another primary component, further adjustment of the translation parameter cannot produce a normal element.

For a generator \(\alpha\) of degree \(n>1\), Section 10.6 considers \((\alpha+t)^{-1}\), with \(t\in k\). Proposition~\GuideRef{v115:prop-primary-reciprocal} clears denominators and applies \(((Y^n-1)/P)(\sigma)\) for every \(P\mid Y^n-1\). The greatest common divisor of one family of coordinate numerators records one local obstruction; a further greatest common divisor with \(T^q-T\) locates all exceptional base-field parameters. This retains nonsplit factors and the multiplicities required when the characteristic divides the degree.

Corollary~\GuideRef{v117:cor-primary-rank-bound} uses operator ranks to bound the common-factor degrees; the bound \((n-1)^2\) in Corollary~\GuideRef{v68:prop-reciprocal-normal-family} follows. The subsequent determinant proof remains available without factoring the operator polynomial. For existence, a sufficient bound ends the argument; if it is inconclusive, the exact common zeros or relations between numerators determine the next step.

In Example~\GuideRef{v68:ex-reciprocal-normal-f343}, where \(\theta^3=2\) in \(\F_{7^3}\), every simple translate lacks the component containing \(\theta^2\), but reciprocation supplies it. When \(n\mid q-1\), the general formula becomes nonvanishing of individual eigen-coordinates, and quartic discriminant constraints improve the separate degree bounds. When \(n\) is a power of the characteristic, formula~\GuideRef{v128:eq-power-degree-exceptions} requires only the base-field zeros of \(f'\). The example \(X^3-X-1\) over \(\F_3\) shows that all three simple translates can fail while all three reciprocals succeed, even though Frobenius is not diagonalizable.

Section 8.2 asks a related but different question: what is the dimension spanned by the reciprocals of all nonzero elements of a binary subspace? For a nonzero finite-dimensional space \(\mathcal U\subseteq\F_2[X]\), let \(M\) be the least common multiple of its nonzero elements, and let \(W\) be the \(\F_2\)-span of the numerators \(M/u\). If the minimal polynomial \(m_\alpha\) over \(\F_2\) is coprime to \(M\), Proposition~\GuideRef{v117:prop-specialization-loss} gives
\[
r\bigl(\mathcal U(\alpha)\bigr)
=\dim_{\F_2}W-\dim_{\F_2}\bigl(W\cap(m_\alpha)\bigr).
\]
At sufficiently large degree, the intersection vanishes and all rational-function relations are retained. When the degree condition fails, the evaluation kernel still records the extra relations exactly. Proposition~\GuideRef{v150:prop-power-intersection} uses the same degree comparison for scalar intersections of power spaces; Corollary~\GuideRef{v150:cor-power-distribution} gives their full distribution. Both routes begin by clearing denominators; their goals then determine which information to keep. Normality requires generation in every primary component, whereas reciprocal rank requires all linear relations and their loss under evaluation.

\Needspace{30\baselineskip}
\GuideSection{7. Root exponents, phases, and coefficient lifting}{v125:route-orbits}
\noindent\textbf{Main route:} Sections 5.6 and 5.8--5.10; compare defining equations in Chapter 10 and sparse multiples in Chapter 17.

Theorem~\GuideRef{v125:thm-orbit-structure} reads two quantities from the additive stabilizer of an orbit: \(\kappa_C\) is the maximum composition exponent over the original field, while \(\omega_C\) is the least extension degree making all irreducible factors binomials. Coefficients may still cancel after maximal compression. Example~\GuideRef{v125:ex-extra-cancellation} separates these two sources of sparsity using polynomials with equal root orders and stabilizers but different supports. Example~\GuideRef{v142:ex-twenty-five-roots} compresses the orbit of a twenty-fifth root to fifth roots: four fibers of size five give $\Phi_{25}(X)=\Phi_5(X^5)$. Chapter 17 uses this composition to identify independent code-coordinate blocks, with relative trace providing the explicit map between the representations.

After decomposing the modulus, relative phase is still needed to assemble local orbits. Along a tower of moduli, Proposition~\GuideRef{v125:prop-phase-lifting} shows that the product of local growth factors equals the global period growth factor times the phase-fiber size. Corollary~\GuideRef{v126:cor-phase-splitting} determines the phase kernel and the existence of additive sections from the local periods. In Example~\GuideRef{v125:ex-160-tower}, unit orbits modulo \(2^a\) grow stably from modulus eight, but adjoining the period modulo five delays unique global lifting until modulus 80. Thereafter, power substitutions in a single quartic give the entire irreducible factor tower.

Example~\GuideRef{v125:ex-160-hensel} lifts this tower's coefficients: power substitution commutes with precision reduction and preserves the maximum composition exponent, but additional zeros may disappear. Proposition~\GuideRef{v125:prop-chain-composition} extends this preservation to finite commutative local rings: the maximal ideal is nilpotent, and coprime factors lift uniquely at each precision. When the residue characteristic divides the exponent, canonical comaximal blocks replace irreducible factors.

Theorem~\GuideRef{v134:thm-subset-periods} at the end of Chapter 5 carries phases into canonical lifted coefficients: subset-sum counts give the weights of orbit periods, while high-order cancellation governs their positive valuations. Proposition~\GuideRef{v129:prop-phase-images} in Section 14.1 uses the same synchronized phase quotient, but images of root pairs require additional checks for shorter periods and collisions between phases. Section 19.7 compresses affine orbits into root coordinates, and the return order recovers the stabilizer size. Phases record synchronization; they label distinct irreducible factors directly only when the image map is injective.

Proposition~\GuideRef{v137:prop-norm-depth} constrains local cancellation depths by a resultant norm. When the multiplicative stabilizer leaves two components, the minimum and sum determine both depths. Proposition~\GuideRef{v137:prop-quadratic-period-depth} realizes this situation through quadratic periods, giving exact dyadic and odd-characteristic valuations.

Theorem~\GuideRef{v141:thm-index-budget} combines an integer discriminant with finite-field multiplication kernels. At minimal positive excess, Theorem~\GuideRef{v141:thm-minimal-excess} confines all exceptional pairs to a matching or triangles at the last layer. Proposition~\GuideRef{v141:prop-period-duality} recovers positions by trace duality. Example~\GuideRef{v141:ex-period-indices} determines a three-component depth multiset from a norm and differences, while retaining the distinction between depths and labels. Proposition~\GuideRef{v144:prop-profile-compression} and Theorem~\GuideRef{v144:thm-invariant-smith} compress complete fibers and coefficient fixed algebras in succession. Example~\GuideRef{v144:ex-profile-tower} uses the same three periods to recover phase valuations at modulus $43^2$.

\GuideSection{8. Finite fields and unramified constructions: roots, integral bases, and precision}{v132:route-adic}
\noindent\textbf{Route:} Chapters 1--4; Sections~\GuideRef{v135:sec-foundations}--\GuideRef{revadic:sec-frob} of Chapter 20; Galois rings and Witt carries in Chapter 5; trace duality and normal bases in Chapters 6--9; character sums in Chapter 18 and coefficient valuations and slopes in Chapter 20.

The first four sections of Chapter 20 can be read after Chapters 1--4. Compatible congruences define $\mathbb Z_p$, and Hensel lifting carries factors modulo $p$ into characteristic zero. Theorem~\GuideRef{revadic:thm-main} constructs a field of $p^n$ elements as the residue field of $\mathbb Q_p(\zeta_{p^n-1})$. Multiplication and Frobenius agree on the roots, while Proposition~\GuideRef{revadic:prop-teich-projection} recovers addition through a limit. After Chapter 5, Corollary~\GuideRef{revadic:cor-witt} identifies this integer ring with the inverse limit of Galois rings.

Example~\GuideRef{v137:ex-construction-comparison} compares four constructions for the fields of sizes eight and nine, distinguishing retained information, known inputs, and precision costs.

Equation~\GuideRef{v135:eq-finite-projection} gives a finite power formula for root projection at each precision; Theorem~\GuideRef{v135:thm-hom-lifting} uniquely lifts every residue-field embedding to compatible ring maps at all precisions. Base extension changes the Frobenius step, whereas truncation fixes the base field and changes the retained digits.

In the integer ring, a unit determinant controls both normality and normal integral bases. Theorem~\GuideRef{revadic:thm-normal} and Proposition~\GuideRef{revadic:prop-dual} describe arbitrary lifts of normal elements and reduction of trace duality. Proposition~\GuideRef{revadic:prop-fixed-precision} extends the fixed-field criterion to every finite precision. Corollary~\GuideRef{v136:cor-integral-hilbert90} identifies the image of Frobenius difference with the relative-trace kernel at every precision; the same cyclic-coordinate proof does not require the extension degree to be invertible. Seventh roots compare primitivity, normality, self-duality, and coefficient valuation in the same example.

Theorem~\GuideRef{v140:thm-teich-rigidity} recovers multiplication from the first carry, and Theorem~\GuideRef{v140:thm-digit-subspaces} identifies scalar multiples of subfields. Theorem~\GuideRef{v147:thm-carry-directions} further classifies nonzero translation directions: in odd characteristic one direction forces a scaled subfield, whereas in characteristic two the direction space is determined by the maximal field of scalars. Corollary~\GuideRef{v144:cor-carry-closure} recovers a subfield by closing a partial subspace under carries. Proposition~\GuideRef{v144:prop-partial-lift} retains mixed-term cancellations on partial root sets and distinguishes liftable permutations from those excluded by the second digit. Precision profiles of representations and preservation depths of code permutations are then determined by Smith exponents of Frobenius difference maps and valuations of spectral-projector matrices.

Theorem~\GuideRef{v138:thm-periodic-cosets} also records a cyclotomic coset as the purely periodic digits of $-a/m$. Proposition~\GuideRef{v138:prop-period-steps} distinguishes the base-$p$ period from the base-$q$ period obtained by grouping $f$ digits, and connects digit rotation with Frobenius on roots of unity. Corollary~\GuideRef{v138:cor-word-gauss} reads Gauss-sum valuations through the Stickelberger formula of Chapter~18; equation~\GuideRef{v138:eq-seven-gauss-period} gives the exact connection with the seventh-root periods. Digit encoding, canonical lifting, and coefficient cancellation thus carry the same orbit into different problems.

Proposition~\GuideRef{v139:prop-digit-leaders} studies leaders through minimal rotations subject to divisibility. Theorem~\GuideRef{v139:thm-period-compression} then identifies additive stabilizers, exact period compression, and monomial composition of canonical factors. Proposition~\GuideRef{v139:prop-digit-crt} records local phases in the same digit model. For codes, Proposition~\GuideRef{v139:prop-digit-bch-gap} detects consecutive roots by addition with carries, while Corollary~\GuideRef{v139:cor-compressed-code} recovers interleaved shorter codes from complete fibers.

Theorem~\GuideRef{v134:thm-subset-periods} in Chapter 5 expresses canonical coefficients through integral subset counts. Theorem~\GuideRef{revadic:thm-integral-resolvent} uses cyclotomic remainders to detect exact zeros and local factors to test initial cancellation. Theorem~\GuideRef{revadic:thm-local-min} determines the minimum valuation across all local components. For information modulo $p^2$, use Proposition~\GuideRef{v133:prop-first-carry} and Proposition~\GuideRef{revadic:prop-first-profile}; to distinguish exact zero from higher divisibility, combine integral remainders with local factors at successive precisions. Theorem~\GuideRef{v143:thm-valuation-smith} recovers the full multiset of local valuations from the Smith exponents of an integer multiplication matrix; local factors and trace duality supply the phase labels. Once the multiset is known, Proposition~\GuideRef{v144:prop-label-precision} attaches all phase labels at precision one above the second-largest distinct valuation, without lifting roots through the largest valuation.

The seventh-root period identities are established in Chapter~5, Example~\GuideRef{v132:ex-seven-period}; Examples~\GuideRef{revadic:ex-eta}, \GuideRef{revadic:ex-sharp-vs-cancel}, and \GuideRef{revadic:ex-period-count} in Chapter~20 continue with trace pairings, local components, and precision loss in Newton division.

Atom factorization holds for every index, while exact root-ratio detection requires \(p\nmid n\) (Corollary~\GuideRef{v151:cor-lucas-frobenius}, Example~\GuideRef{v151:ex-lucas-characteristic}). For a primitive \(\bar c\in\F_q^\times\), Proposition~\GuideRef{v146:prop-adic-lucas} places all roots of \(\Lambda_{q+1}(B,[\bar c])\) in an unramified integer ring. In characteristic two, the two branches over one residue root differ by an element of valuation exactly one, determining the discriminant valuation.

Character orthogonality and Hasse--Davenport lifting from Chapter 18 enter the Fourier and Gauss sums of Chapter 20. Unramified Frobenius invariance identifies the local coefficient field of a Gauss sum; Stickelberger determines its valuation separately. Proposition~\GuideRef{v150:prop-jacobi-carries} uses two-state carry matrices to recover the full Jacobi valuation distribution for a fixed product character; Example~\GuideRef{v150:ex-carry-distribution} distinguishes parameters with equal digit sums but different digit positions. The curve section uses the classical zeta-function theorem. Theorem~\GuideRef{revadic:thm-genus-counts} reconstructs the polynomial from the first $g$ point counts, while Theorem~\GuideRef{revadic:thm-fermat-slopes} obtains Fermat slopes through Jacobi sums. The same theorem gives factors $1+J_DT^{|D|}$ over the original base field by $q$-orbits. Corollary~\GuideRef{v141:cor-orbit-carries} expresses return valuations as carry frequencies, and Example~\GuideRef{v141:ex-fermat-orbits} exhibits equal orbit lengths with different slopes.

\mainmatter

\chapter{Existence of Finite Fields: Characteristic, Adjoining Roots, and Fixed Points}
\label{chap:foundations}
The order of a finite field is determined by its characteristic and its dimension over its prime field. Existence follows by showing that a Frobenius fixed-point set is closed under field operations. Successively adjoining finitely many roots provides a suitable splitting extension, and the root bound determines the size of the fixed-point set.

\section{Characteristic and the Prime Field}

Let \(K\) be a finite field. The sequence
\[
0,\quad 1,\quad 1+1,\quad 1+1+1,\quad\ldots
\]
must contain a repetition, so \(m\cdot1=0\) for some positive integer \(m\).
\leavevmode\phantomsection\label{def:characteristic}\ignorespaces The least positive integer with this property is the \textbf{characteristic} of \(K\), denoted by \(\Char K\).

\begin{proposition}\label{prop:characteristic}
The characteristic of a finite field is prime. If its characteristic is \(p\), it contains a subfield naturally isomorphic to \(\F_p=\Z/p\Z\).
\end{proposition}
\begin{proof}
Since \(1\ne0\), the characteristic \(m\ge2\). If \(m=ab\) with \(1<a,b<m\), then
\[
(a\cdot1)(b\cdot1)=m\cdot1=0.
\]
A field has no zero divisors, so \(a\cdot1=0\) or \(b\cdot1=0\), contrary to the minimality of \(m\).
Thus \(m=p\) is prime.

The map \(\Z\to K\), \(t\mapsto t\cdot1\), preserves addition and multiplication and has kernel \(p\Z\).
Hence \(0,1,2\cdot1,\ldots,(p-1)\cdot1\) form a subfield isomorphic to \(\Z/p\Z\).
\leavevmode\phantomsection\label{def:prime-subfield}\ignorespaces Every subfield contains \(1\), and therefore all these elements. This is the \textbf{prime field} of \(K\).
\end{proof}

We henceforth denote this prime field by \(\F_p\). Its embedding into \(K\), that is, an injective field homomorphism, is uniquely determined: it sends the residue class \(\overline t\) modulo \(p\) to \(t\cdot1\).

\section{Dimension over the Prime Field and the Order of a Finite Field}

\phantomsection\label{def:field-extension}
If \(k\) is a subfield of a field \(E\), we call \(E\) an \textbf{extension field} of \(k\) and \(E/k\) a \textbf{field extension}; \(k\) is the base field.
The addition and multiplication of \(E\), with elements of \(k\) acting as scalars, make \(E\) a \(k\)-vector space:
vector addition is field addition, scalar multiplication is multiplication of an element of \(k\) by an element of \(E\), and the vector-space axioms follow from the field axioms.

\phantomsection\label{def:extension-degree}
The dimension of this vector space is the \textbf{degree} of the extension \(E/k\), denoted by
\[
[E:k]:=\dim_k E.
\]
If \([E:k]=n<\infty\), then \(E/k\) is a \textbf{finite extension}, or an extension of degree \(n\).
A basis of \(E/k\) means a basis of \(E\) as a \(k\)-vector space.
Thus \([E:k]=n\) if and only if there exist \(e_1,\ldots,e_n\in E\) such that every \(x\in E\) has a unique expression
\[
x=a_1e_1+\cdots+a_ne_n,\qquad a_1,\ldots,a_n\in k.
\]
The degree depends on the specified base field. Finiteness of an extension refers to its dimension over that field; when the base field is also finite, these coordinates determine the number of elements of the extension field.
For general remarks on specified embeddings, bases, and restriction of scalars, see the first two sections of Appendix~\ref{chap:polynomial-tools}.

\begin{theorem}\label{thm:cardinality}
If a finite field \(K\) has characteristic \(p\), then for some positive integer \(n\),
\[
|K|=p^n,\qquad n=[K:\F_p].
\]
More generally, if \(k\) is a subfield of a finite field \(K\) and \(|k|=q\), then \(|K|=q^{[K:k]}\).
\end{theorem}
\begin{proof}
Regard \(K\) as a vector space over \(k\). Since \(K\) is a finite set, the process of successively choosing linearly independent elements terminates,
yielding a finite basis \(e_1,\ldots,e_n\). Every element has a unique expression
\[
a_1e_1+\cdots+a_ne_n,\qquad a_i\in k.
\]
Each coefficient has \(q\) possible values, giving \(q^n\) elements. Taking \(k=\F_p\) proves the first assertion.
\end{proof}

\begin{lemma}[The tower law]\label{lem:tower}
Let \(k\subseteq L\subseteq E\) be fields, with both \([L:k]\) and \([E:L]\) finite.
Then \(E/k\) is finite and
\[
[E:k]=[E:L][L:k].
\]
\end{lemma}
\begin{proof}
Let \(u_1,\ldots,u_a\) be a basis of \(L/k\), and \(v_1,\ldots,v_b\) a basis of \(E/L\).
Expand any \(x\in E\) in the \(v_j\), and then expand each coefficient in \(L\) in the \(u_i\), to obtain
\[
x=\sum_{j=1}^{b}b_jv_j
 =\sum_{j=1}^{b}\sum_{i=1}^{a}c_{ij}u_iv_j,
\qquad b_j\in L,\quad c_{ij}\in k.
\]
Thus these \(ab\) products span \(E\) over \(k\).
If they satisfy a \(k\)-linear relation
\[
\sum_{j=1}^{b}\left(\sum_{i=1}^{a}c_{ij}u_i\right)v_j=0,
\]
the linear independence of \(v_1,\ldots,v_b\) over \(L\) gives \(\sum_i c_{ij}u_i=0\) for each \(j\).
The linear independence of \(u_1,\ldots,u_a\) over \(k\) then gives \(c_{ij}=0\) for all \(i,j\).
Hence \(\{u_iv_j:1\le i\le a,\ 1\le j\le b\}\) is a basis of \(E/k\), proving the formula.
\end{proof}

\begin{corollary}\label{cor:subdegree-necessary}
If \(K\) has \(p^n\) elements and a subfield \(L\) has \(p^d\) elements, then \(d\mid n\).
\end{corollary}

This divisibility is necessary for the existence of a subfield. Chapter~\ref{chap:frobenius} proves sufficiency and uniqueness of the corresponding subfield.

\section{Two Identities Arising from Finiteness}

\begin{lemma}[Power maps in characteristic \(p\)]\label{lem:freshman}
In a field of characteristic \(p\), for any \(a,b\) and any integer \(s\geq0\),
\[
(a+b)^{p^s}=a^{p^s}+b^{p^s},\qquad
(ab)^{p^s}=a^{p^s}b^{p^s}.
\]
\end{lemma}
\begin{proof}
For \(s=1\), all intermediate binomial coefficients \(\binom pi\) are divisible by \(p\), and hence vanish in the field.
Iteration gives the additive identity for every \(s\). The multiplicative identity follows from commutativity of multiplication.
\end{proof}

\phantomsection\label{def:frobenius-first}
These identities show that \(\phi:x\mapsto x^p\) preserves the field operations. It is called the \textbf{Frobenius map}.
Its \(s\)-fold iterate is \(\phi^s(x)=x^{p^s}\). An element satisfying \(\phi^s(x)=x\) is a \textbf{fixed point} of this iterate;
the set of all such elements is denoted by \(\Fix(\phi^s)\). We shall use this set directly in the root-set construction.

\begin{proposition}[The power identity in a finite field]\label{prop:fermat}
If \(|K|=Q\), then every \(a\in K\) satisfies \(a^Q=a\).
\end{proposition}
\begin{proof}
If \(a\ne0\), it belongs to the group \(K^\times\) of order \(Q-1\). Lagrange's theorem gives
\(a^{Q-1}=1\), so \(a^Q=a\). The assertion also holds for \(a=0\).
\end{proof}

\begin{example}
\(\Z/4\Z\) is not a field, since \(\overline2\ne0\) but \(\overline2^2=0\).
The field of order four has characteristic \(2\); its additive group is isomorphic to \((\Z/2\Z)^2\), rather than \(\Z/4\Z\).
A field with \(p^n\) elements and the ring of integers modulo \(p^n\) are distinct concepts.
\end{example}

We have shown that the order of a finite field must be a prime power. To construct a finite field of a prescribed prime-power order, it suffices to find an extension over which a specified polynomial splits and then take the fixed points of a Frobenius iterate.
We use polynomial division, the existence of irreducible factors, and B\'ezout's identity. These tools over general fields, together with their proofs, are collected in Appendix~\ref{chap:polynomial-tools}.

\section{Adjoining One Root and the Field Structure of a Quotient Ring}

Fix a monic irreducible polynomial \(f\in k[X]\), with \(\deg f=d\).
Define \(A\equiv B\pmod f\) if and only if \(f\mid A-B\).
Each congruence class has a unique representative of degree less than \(d\). Addition and multiplication of congruence classes define
\[
k[X]/(f)=\{\overline A:A\in k[X]\}.
\]
To perform these operations, carry out ordinary polynomial arithmetic and then take the remainder modulo \(f\).
Replacing representatives by congruent ones leaves sums and products congruent, so the operations are well defined.

\begin{proposition}[Adjoining one root]\label{prop:adjoin}
The quotient \(k[X]/(f)\) is a field. The residue classes of constant polynomials give an embedding of \(k\); if \(\alpha=\overline X\), then
\[
f(\alpha)=0,\qquad
k[X]/(f)=\{a_0+a_1\alpha+\cdots+a_{d-1}\alpha^{d-1}:a_i\in k\}.
\]
This expression is unique, so the extension has degree \(d\). If \(k\) has \(q\) elements, the extension field has \(q^d\) elements.
\end{proposition}
\begin{proof}
The ring axioms follow from those of \(k[X]\), and \(\overline1\ne\overline0\), since a polynomial of positive degree cannot divide \(1\).
If \(\overline A\ne0\), then \(f\nmid A\). Since \(f\) is irreducible, \(\gcd(A,f)=1\).
By B\'ezout's identity, there exist \(U,V\) with \(UA+Vf=1\), whence \(\overline U\,\overline A=\overline1\). Every nonzero element is invertible, so the quotient is a field.

No nonzero constant is divisible by \(f\), so the map from constants to their residue classes is injective.
We have \(f(\overline X)=\overline f=0\), and uniqueness of remainders shows that \(1,\alpha,\ldots,\alpha^{d-1}\) is a \(k\)-basis. The counting assertion follows.
\end{proof}

The constant polynomials identify \(k\) with a subfield of the quotient. Under this identification, we write \(k\subseteq k[X]/(f)\). We shall continue to identify a base field with its image under a specified embedding.

\par\Needspace{12\baselineskip}
\section{Adjoining Finitely Many Roots and Complete Splitting}

\begin{lemma}[Constructing a splitting extension]\label{lem:finite-splitting}
Let \(k\) be a field and \(P(X)\in k[X]\) a nonconstant polynomial. By adjoining finitely many roots, one obtains an extension field \(L\)
in which \(P\) factors completely into linear factors in \(L[X]\). The extension may be chosen to have finite degree.
If \(k\) is finite, so is \(L\).
\end{lemma}
\begin{proof}
Induct on \(D=\deg P\), proving at the same time that the extension can be chosen to have degree at most \(D!\).
For \(D=1\), the root already lies in \(k\), so take \(L=k\).

Suppose \(D>1\). By Lemma \ref{lem:factor}, \(P\) has a monic irreducible factor \(f\) of degree \(d\), where \(1\le d\le D\).
Adjoining one root gives the field \(k_1=k[X]/(f)\) and the element \(\alpha=\bar X\).
Using the embedding defined by constant polynomials to identify \(k\subseteq k_1\), we have
\[
[k_1:k]=d,\qquad P(\alpha)=0.
\]
Division by \(X-\alpha\) in \(k_1[X]\) gives remainder \(P(\alpha)=0\), so
\[
P(X)=(X-\alpha)P_1(X),\qquad \deg P_1=D-1.
\]
The induction hypothesis, applied to the base field \(k_1\) and the polynomial \(P_1\), gives an extension \(L/k_1\)
over which \(P_1\) splits completely, with \([L:k_1]\le(D-1)!\).
Thus \(P\) splits completely over \(L\), and the tower law gives
\[
[L:k]=[L:k_1][k_1:k]\le (D-1)!\,D=D!.
\]
If \(k\) is finite, then \(|L|=|k|^{[L:k]}<\infty\). Each induction step removes one linear factor, so the argument also applies to polynomials with repeated roots.
\end{proof}

\begin{remark}
The construction adjoins roots of the specified polynomial, using an irreducible factor of the remaining polynomial at each step.
The existence of such a factor follows by induction on degree, and its degree is determined by the factorization.
\end{remark}

\par\Needspace{12\baselineskip}
\section{The Field Structure of Frobenius Fixed Points}

\begin{lemma}[The root bound]\label{lem:root-bound}
A nonzero polynomial of degree \(D\) over a field has at most \(D\) distinct roots in any extension field.
\end{lemma}
\begin{proof}
Induct on \(D\). A nonzero constant has no roots. If \(D>0\) and \(a\) is a root, polynomial division gives \(P(X)=(X-a)Q(X)\), with \(\deg Q=D-1\). Every root \(b\ne a\) satisfies \((b-a)Q(b)=0\), and is therefore a root of \(Q\). The induction hypothesis bounds the total number of roots by \(D\).
\end{proof}

\phantomsection\label{def:separable-squarefree}
A nonzero polynomial with no repeated roots in a splitting extension is called \textbf{separable};
if its irreducible factorization has no repeated factor, it is called \textbf{square-free}.
The first condition always implies the second. Over a finite field they are equivalent, since irreducible polynomials have no repeated roots, as we shall prove in Corollary \ref{cor:separable}.

\begin{corollary}[Separability of the power polynomial]\label{lem:simple-root}
In characteristic \(p\), for every integer \(n\ge1\), the polynomial \(X^{p^n}-X\) has no repeated roots.
\end{corollary}
\begin{proof}
By the derivative criterion for repeated roots (Proposition \ref{prop:derivative-roots}), it suffices to note that \((X^{p^n}-X)'=-1\).
\end{proof}

\Needspace{6\baselineskip}
\begin{theorem}[The root-set construction]\label{thm:root-construction}
For every prime \(p\) and positive integer \(n\), there exists a field with exactly \(p^n\) elements.
\end{theorem}
\begin{proof}
Consider \(P=X^{p^n}-X\) over \(\F_p\). Lemma~\ref{lem:finite-splitting} gives a finite extension \(L/\F_p\)
over which \(P\) splits completely. Set
\[
E=\{a\in L:a^{p^n}=a\}.
\]
Since the roots are simple, \(E\) has exactly \(p^n\) elements.

To verify that \(E\) is a subfield, first note that \(0,1\in E\), since \(0^{p^n}=0\) and \(1^{p^n}=1\). If \(a,b\in E\), Lemma \ref{lem:freshman} gives
\[
(a-b)^{p^n}=a^{p^n}-b^{p^n}=a-b,\qquad
(ab)^{p^n}=a^{p^n}b^{p^n}=ab.
\]
If \(a\ne0\), then
\[
(a^{-1})^{p^n}=(a^{p^n})^{-1}=a^{-1}.
\]
Thus \(E\) is closed under subtraction, multiplication, and inversion of nonzero elements, and is a subfield of \(L\).
\end{proof}

The splitting extension \(L\) supplies all the roots; closure under field operations makes their set a subfield. Its order is determined by the degree and separability of \(P\), independently of the chosen splitting extension.

\begin{corollary}[Construction over a specified base field]\label{cor:relative-existence}
Let \(k\) be a given field with \(q\) elements. For every \(m\geq1\), there exists an extension field \(E\supseteq k\)
with \(|E|=q^m\) and \([E:k]=m\).
\end{corollary}
\begin{proof}
Choose an extension of \(k\) over which \(X^{q^m}-X\) splits, and let \(E\) be the set of all its roots.
Since \(q^m\) is a power of \(p\), the preceding closure argument applies, and \(|E|=q^m\).
For \(a\in k\), Proposition \ref{prop:fermat} gives \(a^q=a\); iteration yields \(a^{q^m}=a\), so \(k\subseteq E\).
Finally, \(|E|=q^{[E:k]}\) gives \([E:k]=m\).
\end{proof}

\begin{example}[The field of order four]
\(X^2+X+1\) has no root in \(\F_2\), and is therefore irreducible.
Adjoin a root \(\alpha\), and write \(\F_2[\alpha]\) for the set of all \(A(\alpha)\), with \(A\in\F_2[X]\). Polynomial division gives
\[
\F_2[\alpha]=\{0,1,\alpha,1+\alpha\},\qquad \alpha^2=\alpha+1.
\]
Thus \(\alpha^3=1\) and \(\alpha^4=\alpha\), so all four elements are roots of \(X^4-X\).
In this example, the explicit model obtained by adjoining a root agrees directly with the root-set construction.
\end{example}

\section{The Root Polynomial of a Finite Subfield}

For a finite subset \(S\) of a field \(K\), the monic polynomial with simple roots
\[
R_S(X)=\prod_{a\in S}(X-a)
\]
records precisely all elements of \(S\). When \(S\) is itself a subfield, this polynomial is completely determined by \(|S|\).

\begin{proposition}[The root polynomial of a finite subfield]\label{prop:subfield-root-polynomial}
Let \(K\) be a field, and let \(S\subseteq K\) be a finite subfield with \(|S|=p^n\). Then
\[
\boxed{\ \prod_{a\in S}(X-a)=X^{p^n}-X.\ }
\]
In particular, finite subfields of the same order in a common extension field coincide.
\end{proposition}
\begin{proof}
By Proposition \ref{prop:fermat}, every \(a\in S\) satisfies \(a^{p^n}=a\).
Thus \(R_S\) divides \(X^{p^n}-X\). Both are monic of degree \(p^n\), so they are equal.
If two finite subfields have the same order, the formula identifies both as the complete root set of the same polynomial in \(K\), so they coincide.
\end{proof}

\begin{example}[Root sets, multiplicities, and the standard root polynomial]\label{eg:rootset-multiplicity}
In characteristic two,
\[
X^7-X=X(X^3-1)^2.
\]
In a splitting extension, its distinct roots are zero and the three cube roots of unity, that is, the elements satisfying \(x^3=1\). Together they form \(\F_4\).
The root zero has multiplicity one, whereas each of the three nonzero roots has multiplicity two.
Replacing every multiplicity by one gives the monic polynomial \(X(X^3-1)=X^4-X\).
Thus the intrinsic description of a finite field concerns the set of all distinct roots and its monic polynomial with simple roots;
specifying the root set alone does not recover the multiplicities of the original polynomial.
\end{example}

\begin{remark}[Closure of fixed points and determination of their number]
If \(\tau:K\to K\) is a field homomorphism preserving the identity, then \(\Fix(\tau)=\{x:\tau(x)=x\}\) is a subfield.
Indeed, for \(x,y\in\Fix(\tau)\), we have \(\tau(x-y)=x-y\) and \(\tau(xy)=xy\),
and, when \(x\ne0\), \(\tau(x^{-1})=\tau(x)^{-1}=x^{-1}\).
Neither finiteness of \(K\) nor surjectivity of \(\tau\) is required.
In characteristic \(p\), take \(\tau=\phi^n\) for an integer \(n\ge1\). Its fixed points satisfy a nonzero polynomial equation,
so there are only finitely many. If this polynomial splits completely, separability further shows that their number is \(p^n\).
Closure under the operations, the root bound, and complete splitting play three distinct roles in the argument.
\end{remark}

\par\Needspace{11\baselineskip}
\section*{Exercises}
\phantomsection
\addcontentsline{toc}{section}{Exercises}
\begin{exercise}\label{ex:orders}
Determine whether fields of orders \(12,16,27,36\) can exist. In each possible case, give the characteristic and the dimension over the prime field.
Distinguish the necessary condition imposed by the order from the sufficiency supplied by the existence theorem for finite fields.
\end{exercise}
\begin{exercise}\label{ex:additive}
Prove that the additive group of a finite field of order \(p^n\) is isomorphic to \((\Z/p\Z)^n\).
Explain why this additive group is not cyclic when \(n>1\).
\end{exercise}
\begin{exercise}\label{ex:tower}
What orders can subfields of a field with \(2^{12}\) elements have?
Explain whether the tower law gives a necessary or a sufficient condition here.
\end{exercise}

\begin{exercise}[Power identities and field structure]\label{ex:wu-ring-power}
Prove that every element of \(\Z/6\Z\) satisfies \(x^3=x\). Explain why this identity alone does not imply that the ring is a field.
Next, determine whether a finite nonzero commutative ring \(R\) with identity must be a field if \(x^m=1\) for every nonzero \(x\), where \(m\ge1\) is fixed.
\end{exercise}

\begin{exercise}\label{ex:bezout}
Verify \((X+1)(X+2)-(X^2+1)=1\) in \(\F_3[X]\). Use this identity to find the inverse of \(\overline X+1\) in \(\F_3[X]/(X^2+1)\), and explain why this quotient ring is a field.
\end{exercise}
\begin{exercise}\label{ex:reducible}
If \(f\) is a reducible polynomial of positive degree, prove that \(k[X]/(f)\) is not a field.
\end{exercise}
\begin{exercise}[A degree bound for a splitting extension]\label{ex:split-bound}
Let \(P\in k[X]\) have degree \(D\), and suppose that \(P\) has a root in \(k\).
Prove that there is a finite extension \(L/k\) over which \(P\) splits, with \([L:k]\le(D-1)!\).
If \(k\) is finite, does this imply that every extension over which \(P\) splits satisfies this bound?
\end{exercise}
\begin{exercise}\label{ex:roots-needs-splitting}
Consider the root set of \(X^4-X\) in \(\F_2\).
Explain why closure of the root set under the field operations does not by itself ensure that the set has four elements.
\end{exercise}

\begin{exercise}\label{ex:supp-crt}
In \(\F_3[X]\), find a polynomial \(A\) of degree less than two such that \(A(1)=2\) and \(A(-1)=1\).
Then find two idempotents in \(\F_3[X]/(X^2-1)\), that is, elements satisfying \(e^2=e\),
whose values at \(1,-1\) are \((1,0)\) and \((0,1)\), respectively.
\end{exercise}

\begin{exercise}[When a subspace is a subfield]\leavevmode\phantomsection\label{ex:wu-rootfield}
Let \(V\) be a finite \(\F_p\)-subspace of a field of characteristic \(p\), and suppose that \(1\in V\).
Prove that \(V\) is a subfield if and only if it is closed under multiplication.
If \(\dim_{\F_p}V=n\), prove that this condition is also equivalent to \(\prod_{v\in V}(X-v)=X^{p^n}-X\).
\end{exercise}

\chapter{Uniqueness and Generators of Finite Fields}
\label{chap:generators}
Matching roots determines the isomorphism type of a finite field, while the root bound forces its multiplicative group to be cyclic. These facts yield field generators, irreducible polynomials of every degree, and quotient-ring coordinates, connecting intrinsic structure with computational models.

\section{Minimal Polynomials and Fields Generated by One Element}

\phantomsection\label{def:generated-subfield}
Let \(k\subseteq E\) be fields and \(S\subseteq E\). The intersection of all subfields of \(E\) containing \(k\cup S\) is denoted by \(k(S)\) and is called the \textbf{subfield generated by \(k\) and \(S\)}.
This family contains \(E\), and its intersection is closed under the field operations. Thus \(k(S)\) exists and is the smallest subfield containing \(k\cup S\).
For \(S=\{\alpha_1,\ldots,\alpha_r\}\), we write \(k(\alpha_1,\ldots,\alpha_r)\); when \(S\) is empty, \(k(S)=k\).
If \(E=k(\alpha)\), then \(E/k\) is a \textbf{simple extension}, and \(\alpha\) is a field generator.

\phantomsection\label{def:polynomial-subring}
For \(\alpha\in E\), also write
\[
k[\alpha]=\{A(\alpha):A\in k[X]\}.
\]
This is the subring generated by \(k\) and \(\alpha\), namely the smallest subring of \(E\) containing both. Hence \(k[\alpha]\subseteq k(\alpha)\).
The two notations record polynomial evaluation and field generation, respectively. The next proposition shows that they coincide when \(\alpha\) is algebraic over \(k\).

\begin{definition}\label{def:algebraic-minpoly}\label{def:annihilating-polynomial}
An element \(\alpha\) of an extension field of \(k\) is \textbf{algebraic} over \(k\) if \(A(\alpha)=0\) for some nonzero \(A\in k[X]\).
A nonzero polynomial \(A\) satisfying \(A(\alpha)=0\) is an \textbf{annihilating polynomial} of \(\alpha\) over \(k\).
Choose an annihilating polynomial of least degree and divide it by its leading coefficient. The resulting monic polynomial is the \textbf{minimal polynomial} of \(\alpha\) over \(k\), denoted by \(m_{\alpha,k}\).
\end{definition}

\begin{proposition}\label{prop:minpoly}
Suppose that \(\alpha\) is algebraic over \(k\) and that \(m=m_{\alpha,k}\) has degree \(d\). Then:
\begin{enumerate}
  \item \(m\) is irreducible, and for \(A\in k[X]\), \(A(\alpha)=0\) if and only if \(m\mid A\);
  \item the minimal polynomial is unique;
  \item \(k(\alpha)=k[\alpha]\);
  \item evaluation induces
  \[
  k[X]/(m)\ \cong\ k(\alpha),\qquad \overline A\longmapsto A(\alpha),
  \]
  and \(1,\alpha,\ldots,\alpha^{d-1}\) is a basis of \(k(\alpha)/k\). In particular,
  \[
  [k(\alpha):k]=\deg m_{\alpha,k}=d.
  \]
\end{enumerate}
\end{proposition}
\begin{proof}
If \(m=AB\), with both \(A,B\) of positive degree, then \(A(\alpha)B(\alpha)=0\). A field has no zero divisors, so one of these values is zero.
The corresponding factor has degree less than \(\deg m\), contradicting minimality. Thus \(m\) is irreducible.

For any \(A\in k[X]\), division with remainder gives \(A=Qm+R\), where \(R=0\) or \(\deg R<d\).
If \(A(\alpha)=0\), then \(R(\alpha)=0\), and minimality forces \(R=0\), so \(m\mid A\).
The converse follows immediately by evaluation. If another monic polynomial \(m_1\) of least degree vanishes at \(\alpha\),
then \(m\mid m_1\). Since both are monic of the same degree, \(m_1=m\).

Consider the evaluation homomorphism \(\operatorname{ev}_\alpha:k[X]\to k[\alpha]\), \(A\mapsto A(\alpha)\).
The divisibility criterion shows that its kernel is precisely the ideal \((m)\). Consequently,
\[
\overline A\longmapsto A(\alpha)
\]
is well defined, surjective, and injective on \(k[X]/(m)\), and preserves addition, multiplication, and the identity.
Since \(m\) is irreducible, the quotient is a field, so \(k[\alpha]\) is also a field.
It contains \(k\) and \(\alpha\), and hence contains the smallest such subfield \(k(\alpha)\).
Conversely, \(A(\alpha)\in k(\alpha)\), so \(k[\alpha]=k(\alpha)\).

Division with remainder shows that \(1,\alpha,\ldots,\alpha^{d-1}\) span this field.
A nontrivial \(k\)-linear relation among them would give a nonzero annihilating polynomial of degree less than \(d\), contradicting minimality.
Thus they form a basis.
\end{proof}

\phantomsection\label{def:operator-annihilator}
The same terminology applies to evaluating polynomials at linear operators. Let \(A:V\to V\) be a linear operator on a finite-dimensional \(k\)-space.
For \(g(T)=\sum_i c_iT^i\), define \(g(A)=\sum_i c_iA^i\), with \(A^0=\id_V\).
A nonzero \(g\in k[T]\) satisfying \(g(A)=0\) is an annihilating polynomial of \(A\).
If only \(g(A)v=0\) is required, it is an annihilating polynomial of the vector \(v\) relative to \(A\).
The monic polynomials of least degree satisfying these two conditions are called the minimal polynomial of \(A\) and the minimal annihilating polynomial of \(v\) relative to \(A\), respectively.
Such polynomials exist: if \(\dim_kV=N\), the linear operators form a space of dimension \(N^2\), so \(\id_V,A,\ldots,A^{N^2}\) are linearly dependent. The resulting operator relation also annihilates every specified vector.
Division with remainder, as in Proposition~\ref{prop:minpoly}, proves uniqueness of the monic polynomial of least degree.
The first condition requires vanishing on the entire space, and the second only on a specified vector. Both minimal polynomials must be distinguished from general annihilating polynomials.

\section{Extending a Base-Field Isomorphism After Adjoining a Root}

Let \(\tau:k\to k'\) be a field isomorphism. For \(A=\sum a_iX^i\in k[X]\), write
\[
A^\tau=\sum\tau(a_i)X^i\in k'[X].
\]
Applying \(\tau\) coefficientwise preserves addition, multiplication, and degree, and hence also irreducibility and divisibility.

\begin{lemma}[Extending an isomorphism in one step]\leavevmode\phantomsection\label{lem:extend}
Let \(\alpha\) be algebraic over \(k\), with minimal polynomial \(m\).
If \(\beta\) is a root of \(m^\tau\) in an extension field of \(k'\), there is a unique field isomorphism
\[
\widetilde\tau:k(\alpha)\longrightarrow k'(\beta)
\]
extending \(\tau\) and sending \(\alpha\) to \(\beta\). Explicitly,
\[
\widetilde\tau\!\left(\sum_{i=0}^{d-1}a_i\alpha^i\right)
=\sum_{i=0}^{d-1}\tau(a_i)\beta^i,
\qquad d=\deg m.
\]
\end{lemma}
\begin{proof}
Applying the isomorphism coefficientwise preserves irreducibility, so \(m^\tau\) is monic and irreducible. Since it vanishes at \(\beta\), it is the minimal polynomial of \(\beta\) over \(k'\).

Define a map of quotient rings by \(\overline A\mapsto\overline{A^\tau}\).
If \(A-B\) is divisible by \(m\), then \(A^\tau-B^\tau\) is divisible by \(m^\tau\), so this map is independent of representatives.
It preserves the ring operations, and applying \(\tau^{-1}\) coefficientwise gives its inverse. Thus it is an isomorphism
\[
k[X]/(m)\longrightarrow k'[X]/(m^\tau).
\]
Composing with the evaluation isomorphisms at both ends gives an isomorphism \(k(\alpha)\to k'(\beta)\).
It restricts to \(\tau\) on the constants and sends \(\alpha\) to \(\beta\), as required.
Any homomorphism with these properties must send \(A(\alpha)\) to \(A^\tau(\beta)\), proving uniqueness.
\end{proof}

\section{Matching Roots and Uniqueness of Splitting Fields}

\begin{definition}\label{def:splitting-field}
Let \(P\in k[X]\) be nonconstant. An extension field \(E\) is a \textbf{splitting field} of \(P\) over \(k\) if \(P\) splits completely over \(E\)
and \(E\) is generated by \(k\) and all the roots of \(P\).
\end{definition}

Lemma~\ref{lem:finite-splitting} constructs an extension over which the polynomial splits. The subfield generated by all its roots in that extension is a splitting field.
Successively matching roots proves uniqueness up to isomorphism.

\begin{proposition}[Matching roots between splitting fields]\label{prop:splitting-unique}
Let \(\tau:k\to k'\) be an isomorphism, let \(E/k\) be a splitting field of \(P\), and let \(E'/k'\) be a splitting field of \(P^\tau\).
Then \(\tau\) extends to an isomorphism \(E\to E'\).
\end{proposition}
\begin{proof}
Write the distinct roots of \(P\) as \(\alpha_1,\ldots,\alpha_t\), so that \(E=k(\alpha_1,\ldots,\alpha_t)\).
Starting with \(\tau\), suppose that an embedding
\[
\tau_i:k(\alpha_1,\ldots,\alpha_i)\longrightarrow E'
\]
has been constructed. Regard it as an isomorphism from its domain onto its image. Let \(m\) be the minimal polynomial of \(\alpha_{i+1}\) over the current domain.
Since \(P(\alpha_{i+1})=0\), Proposition \ref{prop:minpoly} gives \(m\mid P\). Applying \(\tau_i\) coefficientwise yields
\[
m^{\tau_i}\mid P^\tau.
\]
The right-hand side splits completely in \(E'\), so \(m^{\tau_i}\), which has positive degree, has a root \(\beta\) in \(E'\).
The one-step extension lemma gives the next embedding \(\tau_{i+1}\).
After finitely many steps, we obtain \(\sigma:E\to E'\).

It remains to prove surjectivity. Over \(E\), write \(P\) as a product of linear factors, retaining multiplicities.
Applying \(\sigma\) to the coefficients and roots gives a complete factorization of \(P^\tau\) over \(\sigma(E)\).
Thus every root of \(P^\tau\) belongs to \(\sigma(E)\). Since \(E'\) is generated by \(k'\) and these roots,
we have \(E'=\sigma(E)\).
\end{proof}

\begin{theorem}[Uniqueness up to isomorphism of finite fields of the same order]\label{thm:unique}
Any two fields with \(p^n\) elements are isomorphic. Such an isomorphism restricts to the canonical isomorphism between their prime fields; after identifying both prime fields with \(\F_p\), it is an \(\F_p\)-isomorphism.
Thus each prime power \(p^n\) determines exactly one isomorphism type of finite field, denoted by \(\F_{p^n}\).
\end{theorem}
\begin{proof}
If \(|E|=p^n\), Proposition \ref{prop:subfield-root-polynomial} identifies all elements of \(E\) with the complete root set of \(X^{p^n}-X\). Hence \(E\) is a splitting field of this polynomial over \(\F_p\).
The same holds for any other field with \(p^n\) elements, and Proposition \ref{prop:splitting-unique} gives an isomorphism.
Existence was established in Theorem \ref{thm:root-construction}.
\end{proof}

\begin{corollary}[Uniqueness over a specified base field]\label{cor:relative-unique}
Let \(k\) be a specified field with \(q\) elements, and let \(E,E'\) be fields with \(q^m\) elements containing \(k\).
Then there is an isomorphism \(E\to E'\) fixing \(k\) pointwise.
\end{corollary}
\begin{proof}
Both are splitting fields of \(X^{q^m}-X\) over \(k\). Apply the root-matching proposition to \(\id_k\).
\end{proof}

\begin{example}[An isomorphism need not be unique]
In the field with four elements, let \(\alpha\) satisfy \(\alpha^2+\alpha+1=0\).
The other root of this polynomial is \(\alpha+1\).
The one-step extension lemma gives an automorphism fixing \(\F_2\) and sending \(\alpha\) to \(\alpha+1\).
This differs from the identity. Thus uniqueness of the isomorphism type does not imply uniqueness of the isomorphism.
\end{example}

Having determined the isomorphism type, we turn to generators within the field. The root bound constrains the least common multiple of the element orders in a finite multiplicative group, and this constraint suffices to prove cyclicity.

\section{Combining Elements of Coprime Orders}

\begin{lemma}\label{lem:coprime-product}
If commuting elements \(a,b\) of a group satisfy \(\ord(a)=r\), \(\ord(b)=s\), and \(\gcd(r,s)=1\),
then \(\ord(ab)=rs\).
\end{lemma}
\begin{proof}
First, \((ab)^{rs}=1\). If \((ab)^t=1\), then \(a^t=b^{-t}\) lies in both \(\langle a\rangle\) and \(\langle b\rangle\). Its order divides both \(r\) and \(s\), so it is the identity.
Thus \(r\mid t\) and \(s\mid t\), and coprimality gives \(rs\mid t\).
\end{proof}

\begin{lemma}[Realizing the least common multiple as an element order]\label{lem:exponent}
Let \(G\) be a finite abelian group, and set
\[
M=\lcm\{\ord(g):g\in G\},
\]
called the \textbf{exponent} of \(G\). Then \(G\) contains an element of order exactly \(M\).
\end{lemma}
\begin{proof}
If \(G=\{1\}\), then \(M=1\), the order of the identity. Otherwise write \(M=\ell_1^{a_1}\cdots\ell_t^{a_t}\), where the \(\ell_i\) are distinct primes.
For each \(i\), the definition of least common multiple allows us to choose \(g_i\in G\) such that
\[
\ord(g_i)=\ell_i^{a_i}u_i,\qquad \gcd(\ell_i,u_i)=1.
\]
Then \(h_i=g_i^{u_i}\) has order \(\ell_i^{a_i}\),
since \(h_i^v=1\) is equivalent to \(\ell_i^{a_i}u_i\mid u_iv\), or \(\ell_i^{a_i}\mid v\).
These orders are pairwise coprime, and all the \(h_i\) commute. Repeated application of Lemma \ref{lem:coprime-product}
shows that \(h_1\cdots h_t\) has order \(M\).
\end{proof}

\section{The Root Bound and Multiplicative Cyclicity}

\begin{theorem}[Cyclicity of finite multiplicative subgroups of a field]\label{thm:cyclic}
Every finite subgroup of the multiplicative group \(K^\times\) of a field \(K\) is cyclic.
In particular, if \(|E|=Q\), then \(E^\times\) is cyclic of order \(Q-1\).
\end{theorem}
\begin{proof}
Let \(G\subseteq K^\times\) be a finite subgroup, and let \(M\) be the least common multiple of all its element orders.
Every \(g\in G\) satisfies \(g^M=1\), so all elements of \(G\) are roots of \(X^M-1\).
The root bound gives \(|G|\leq M\).

On the other hand, \(G\) is abelian, so Lemma \ref{lem:exponent} gives an element \(h\) of order \(M\).
Its cyclic subgroup has \(M\) elements, so \(M\leq|G|\). Therefore
\[
|G|=M=|\langle h\rangle|,
\]
and hence \(G=\langle h\rangle\).
\end{proof}

This proof uses a special property of fields: a polynomial of degree \(M\) cannot have more than \(M\) roots.
A general finite abelian group has no analogous restriction. For example, \((\Z/2\Z)^2\) has four elements, each of order at most two.

\begin{definition}\label{def:primitive-element}
A generator of the multiplicative group of a finite field \(E\) is a \textbf{primitive element} of \(E\).
If \(|E|=Q\), a primitive element \(\zeta\) satisfies
\[
E^\times=\{1,\zeta,\zeta^2,\ldots,\zeta^{Q-2}\},\qquad \ord(\zeta)=Q-1.
\]
\end{definition}

\begin{corollary}\label{cor:power-roots}
Let \(|E|=Q\) and \(t\geq1\). The equation \(x^t=1\) has exactly \(\gcd(t,Q-1)\) solutions in \(E\).
\end{corollary}
\begin{proof}
Write \(x=\zeta^j\), with \(j\) taken modulo \(Q-1\). Then \(x^t=1\) is equivalent to \((Q-1)\mid tj\), which holds for exactly \(\gcd(t,Q-1)\) residue classes.
\end{proof}

\begin{example}[Primitive elements in the field with nine elements]
Let \(E=\F_3(\alpha)\), where \(\alpha^2=-1=2\).
Since \(X^2+1\) has no root in \(\F_3\), this is indeed a field with nine elements.
We have \(\alpha^2=2\) and \(\alpha^4=1\), so \(\alpha\) has order four and is not primitive.
However,
\[
(1+\alpha)^2=2\alpha,\qquad
(1+\alpha)^4=2,\qquad (1+\alpha)^8=1,
\]
so \(1+\alpha\) has order eight and is primitive.
Thus an element generating the whole extension field need not generate its multiplicative group.
\end{example}

\section{Successively Increasing the Multiplicative Order: Gauss's Method}
The cyclicity proof also gives a finite construction. Suppose that \(a\in E^\times\) has order \(m\).
If \(m<|E^\times|\), choose \(b\notin\langle a\rangle\), and denote its order by \(s\).
Then \(s\nmid m\): otherwise \(b\) and all \(m\) elements of \(\langle a\rangle\) would be roots of \(X^m-1\), contradicting the root bound.

For each prime dividing \(m\) or \(s\), choose the larger prime power appearing in their factorizations; if they are equal, choose the one from \(m\).
Let \(d\) be the product of those chosen from \(m\), and \(e\) the product of the remaining ones, with an empty product equal to \(1\).
Then \(d\mid m\), \(e\mid s\), \(\gcd(d,e)=1\), and \(de=\lcm(m,s)\).
The element
\[
c=a^{m/d}b^{s/e}
\]
has order \(de>m\). Replace \(a\) by \(c\) and repeat. The order strictly increases, so the process eventually yields a primitive element.
This is a finite construction; its computational complexity also depends on determining element orders and factoring integers.

\begin{example}
In \(\F_{13}\), the elements \(3\) and \(5\) have orders \(3\) and \(4\), respectively. These orders are coprime,
so \(3\cdot5=2\) has order \(12\), giving a primitive element of \(\F_{13}\).
\end{example}

A multiplicative primitive element also generates the whole extension field. Its minimal polynomial therefore converts an already constructed finite field into a quotient-ring model and yields irreducible polynomials of every degree.

\section{Primitive Elements Make Finite Extensions Simple}

\begin{proposition}\label{prop:primitive-generates}
If \(k\subseteq E\) are finite fields and \(\zeta\) is a primitive element of \(E\), then \(E=k(\zeta)\).
\end{proposition}
\begin{proof}
The subfield \(k(\zeta)\) contains \(0\) and all nonnegative integral powers of \(\zeta\).
These powers comprise all nonzero elements of \(E\), so the subfield is \(E\).
\end{proof}

\begin{theorem}[Existence of irreducible polynomials of every degree]\label{thm:irreducible-existence}
Let \(k\) be any finite field, with \(|k|=q\). For every \(n\geq1\), there is a monic irreducible polynomial of degree \(n\) in \(k[X]\).
\end{theorem}
\begin{proof}
By Corollary \ref{cor:relative-existence}, first construct a field \(E\) with \(q^n\) elements containing \(k\).
By Theorem \ref{thm:cyclic}, choose a primitive element \(\zeta\) of \(E\). Then \(E=k(\zeta)\). Set \(f=m_{\zeta,k}\). This polynomial is monic and irreducible, and
\[
\deg f=[k(\zeta):k]=[E:k]=n.
\]
\end{proof}

The root-set construction establishes the existence of finite fields, while the minimal polynomial of a primitive element supplies quotient-ring coordinates. These describe the intrinsic field structure and a concrete computational model, respectively.

\begin{theorem}[Quotient-ring representations of finite fields]\label{thm:quotient-model}
Let \(k\) be a field with \(q\) elements.
\begin{enumerate}
  \item If \(f\in k[X]\) is monic and irreducible of degree \(n\), then \(k[X]/(f)\) is a field with \(q^n\) elements containing \(k\).
  \item Every field with \(q^n\) elements containing \(k\) is isomorphic to such a quotient. In fact, for every monic irreducible \(f\) of degree \(n\),
  it is isomorphic to \(k[X]/(f)\) over \(k\).
\end{enumerate}
\end{theorem}
\begin{proof}
The first assertion is Proposition \ref{prop:adjoin}.
For the second, choose a primitive element of \(E\) and its minimal polynomial. Proposition \ref{prop:minpoly} gives a quotient-ring representation.
For any other monic irreducible \(f\) of degree \(n\), the quotient ring and \(E\) are both extensions of \(k\) with \(q^n\) elements,
so Corollary \ref{cor:relative-unique} gives a \(k\)-isomorphism between them.
\end{proof}

\section{Quotient-Ring Coordinates and Finite-Field Arithmetic}\label{v48:sec-quotient-coordinates}

Let \(f\in k[X]\) be monic and irreducible, with \(\deg f=n\ge1\). In \(k[X]/(f)\), write \(\alpha=\overline X\) for the residue class of \(X\).
The elements of the quotient ring are residue classes of polynomials. Choosing the unique representative of degree less than \(n\) for computation gives coordinates
\[
a_0+a_1\alpha+\cdots+a_{n-1}\alpha^{n-1}.
\]
Existence of the representative follows from division with remainder. Uniqueness follows because the difference of two distinct representatives cannot have degree less than \(n\) while also being divisible by \(f\).
Thus the vector space \(k^n\) supplies coordinates for all elements, while multiplication also depends on the chosen polynomial \(f\).

Writing \(\operatorname{rem}_f A\) for the remainder of \(A\) upon division by \(f\), the coordinate operations are
\[
A\mathbin{+}B=\operatorname{rem}_f(A+B),\qquad
A\mathbin{\cdot}B=\operatorname{rem}_f(AB).
\]
If \(A'=A+Uf\) and \(B'=B+Vf\), then the differences between \(A'+B'\) and \(A+B\), and between \(A'B'\) and \(AB\), are divisible by \(f\).
Thus these operations are independent of the representatives. Associativity and distributivity descend from the corresponding polynomial identities.

Let \(0\ne A\) with \(\deg A<n\). Irreducibility gives \(\gcd(A,f)=1\), and the extended Euclidean algorithm yields
\[
U(X)A(X)+V(X)f(X)=1.
\]
Evaluation at \(\alpha\) gives \(U(\alpha)A(\alpha)=1\), so the inverse is represented by \(\operatorname{rem}_f U\).
The field operations are therefore reduced to arithmetic in the base field, polynomial division, and B\'ezout's identity.

\begin{example}[A field with sixteen elements]\leavevmode\phantomsection\label{exm:quartic-model}
Take \(f=X^4+X+1\in\F_2[X]\). It is nonzero at \(0,1\), so it has no linear factor.
The only monic irreducible quadratic in \(\F_2[X]\) is \(X^2+X+1\):
a monic quadratic without a root must have constant term \(1\), and evaluation at \(1\) then forces the linear coefficient to be \(1\).
Modulo \(X^2+X+1\), we have \(X^3\equiv1\) and \(X^4\equiv X\), so \(f\equiv1\).
A reducible quartic without a linear factor must have an irreducible quadratic factor. Hence \(f\) is irreducible.

Consequently,
\[
E=\F_2[\alpha],\qquad \alpha^4=\alpha+1
\]
has sixteen elements, each uniquely expressible as \(a+b\alpha+c\alpha^2+d\alpha^3\), with coefficients in \(\F_2\).
For example,
\[
(1+\alpha^2)(\alpha+\alpha^3)
=\alpha+\alpha^5=\alpha^2,
\]
since \(\alpha^5=\alpha^2+\alpha\).
Moreover, the order of \(\alpha\) divides fifteen and is not \(1,3,5\):
the degree-four minimal polynomial implies \(\alpha^3\ne1\),
and the same degree argument gives \(\alpha^5=\alpha^2+\alpha\ne1\).
Thus \(\alpha\) has order fifteen and is also a primitive element.
\end{example}

\begin{example}[Inversion by successive polynomial divisions]\label{v48:ex-bezout-coordinates}
With \(f=X^4+X+1\) as in Example~\ref{exm:quartic-model}, take \(h=X^2+1\). In \(\F_2[X]\),
\[
f=(X^2+1)h+X,\qquad h=X\cdot X+1.
\]
Substituting the first identity into the second gives
\[
1=h+X\bigl(f+(X^2+1)h\bigr)
  =(X^3+X+1)h+Xf.
\]
Hence
\[
(1+\alpha^2)^{-1}=\alpha^3+\alpha+1.
\]
Direct multiplication gives \((1+\alpha^2)(\alpha^3+\alpha+1)=\alpha^5+\alpha^2+\alpha+1=1\),
in agreement with the B\'ezout calculation.
\end{example}

\begin{remark}\label{def:primitive-polynomial}
An irreducible polynomial ensures that a root generates an extension of the prescribed degree, but does not ensure that the root generates the multiplicative group.
A monic irreducible polynomial whose root is a primitive element of the resulting finite field is called a \textbf{primitive polynomial}.
The proof of Theorem \ref{thm:irreducible-existence} in fact produces primitive polynomials of every degree.
\end{remark}

By uniqueness of finite fields of the same order, irreducible polynomials of the same degree give quotient-ring models with the same \(k\)-isomorphism type,
although their reduction relations generally give different multiplication rules in coordinates.
Let \(f,g\in k[X]\) be monic irreducible polynomials of degree \(n\). Choose a root \(b\) of \(f\) in \(k[Y]/(g)\). Then
\[
k[X]/(f)\longrightarrow k[Y]/(g),\qquad
\overline{A(X)}\longmapsto A(b)
\]
is a \(k\)-isomorphism by the one-step extension lemma. Conversion between these coordinate systems requires such a choice of root; simply identifying the respective polynomial bases need not define a field isomorphism.

\begin{example}[Two reduction relations for the same field with nine elements]\label{v48:ex-coordinate-model-change}
In \(E=\F_3[\alpha]\), with \(\alpha^2+1=0\), set \(\beta=1+\alpha\). Then
\[
\beta^2+\beta+2=0.
\]
The three values of \(g(Y)=Y^2+Y+2\) on \(\F_3\) are \(2,1,2\), so \(g\) is irreducible.
Thus
\[
\F_3[Y]/(g)\longrightarrow E,\qquad
u+v\overline Y\longmapsto (u+v)+v\alpha
\]
is a field isomorphism, with inverse \(a+b\alpha\mapsto(a-b)+b\overline Y\).
The coordinate-change matrices are
\[
\begin{pmatrix}1&1\\0&1\end{pmatrix},\qquad
\begin{pmatrix}1&-1\\0&1\end{pmatrix}.
\]
Here the generator's minimal polynomial has changed: \(\overline Y\) satisfies \(g\), while \(\alpha\) satisfies \(X^2+1\).
Sending \(\alpha\) to \(1+\alpha\) within the same \(X^2+1\) model, however, does not preserve the defining quadratic relation; see Exercise~\ref{ex:wu-embedding-compatibility}.
\end{example}

\par\Needspace{11\baselineskip}
\section*{Exercises}
\phantomsection
\addcontentsline{toc}{section}{Exercises}
\begin{exercise}\label{ex:f8iso}
Prove that \(f=X^3+X+1\) and \(g=X^3+X^2+1\) are both irreducible over \(\F_2\).
Let \(\alpha\) be a root of \(f\) and \(\beta\) a root of \(g\).
Verify that \(g(\alpha+1)=0\), and give an explicit isomorphism between the two fields \(\F_2(\beta)\) and \(\F_2(\alpha)\) of order eight.
\end{exercise}
\begin{exercise}\label{ex:iso-condition}
Explain why the one-step extension lemma does not allow \(\alpha\) to be sent to an arbitrary element of the target extension.
What condition must its image satisfy?
\end{exercise}
\begin{exercise}\label{ex:finite-domain}
Prove that every finite commutative integral domain is a field. If one starts with closure of a root set under subtraction and multiplication,
can this result replace the separate verification of closure under inverses in the root-set construction?
\end{exercise}

\begin{exercise}[Images of roots and algebraic relations]\label{ex:wu-embedding-compatibility}
Let \(E=\F_3(\alpha)\), with \(\alpha^2=-1\).
Prove that every field homomorphism \(E\to E\) fixing \(\F_3\) sends \(\alpha\) to \(\alpha\) or \(-\alpha\),
and that both choices define homomorphisms. Why is the linear map sending \(\alpha\) to \(1+\alpha\) not a field homomorphism?
\end{exercise}

\begin{exercise}\label{ex:primitive-test}
Let \(G=\langle\zeta\rangle\) have order \(N\). Prove \(\ord(\zeta^j)=N/\gcd(N,j)\), and deduce that \(G\) has \(\varphi(N)\) generators.
Here \(\varphi(N)\) counts the invertible residue classes modulo \(N\), with \(\varphi(1)=1\).
\end{exercise}
\begin{exercise}\label{ex:power-eq}
Let \(E\) have \(Q\) elements, let \(a\in E^\times\) and \(t\geq1\), and set \(g=\gcd(t,Q-1)\).
Prove that \(x^t=a\) is solvable if and only if \(a^{(Q-1)/g}=1\), and that it has exactly \(g\) solutions when solvable.
\end{exercise}
\begin{exercise}\label{ex:cosets}
Let \(k\subseteq E\) be finite fields, with \(|k|=q\) and \(|E|=q^m\).
Using only the coset decomposition of a group, prove
\[
E^\times=\bigsqcup_{i=1}^{s}z_i k^\times,
\qquad q^m-1=s(q-1).
\]
Determine \(s\), and explain how the subgroup \(k^\times\) arises from the subfield \(k\subseteq E\). Then give an irreducible polynomial whose root set is not closed under multiplication.
\end{exercise}

\begin{exercise}[Generating the extension and generating the multiplicative group]\label{ex:wu-powers-groups}
In \(\F_9=\F_3(\alpha)\), with \(\alpha^2=-1\), prove that \(\alpha\) generates the extension but is not primitive.
Find all primitive elements, and decompose \(\F_9^\times\) into cosets of \(\F_3^\times\).
\end{exercise}

\begin{exercise}\label{ex:phi5}
Prove that \(h=X^4+X^3+X^2+X+1\) is irreducible over \(\F_2\).
If \(\beta\) is a root, prove that \(\F_2(\beta)\) has sixteen elements but \(\beta\) has multiplicative order five.
\end{exercise}
\begin{exercise}\label{ex:inverse16}
In \(\F_2[\alpha]\), with \(\alpha^4=\alpha+1\), find \(\alpha^{-1}\) and \((\alpha+1)^{-1}\).
\end{exercise}
\begin{exercise}\label{ex:primitive-polynomial}
Explain why Theorem \ref{thm:irreducible-existence} also proves the existence of primitive polynomials of every degree.
Identify the step ensuring that the minimal polynomial has degree equal to the extension degree.
\end{exercise}

\chapter{Frobenius: Subfields, Conjugates, and Descent of Root Sets}
\label{chap:frobenius}
The fixed points of Frobenius determine subfields, and its orbits determine minimal polynomials. Relating functions on root sets to the Frobenius action makes factorization the algebraic expression of orbit decomposition.

Write \(q=p^r\), and fix a field \(k\) with \(q\) elements. A specified finite extension is denoted by \(E/k\), with \(|E|=q^m\) and \([E:k]=m\) recording its order and degree.
When comparing subfields, intersections, and composita, we first specify a common finite extension containing the fields in question; any change of this ambient extension will be stated.

\section{The Order of Frobenius and the Extension Degree}

\begin{definition}\label{def:frobenius-maps}
On a field of characteristic \(p\), the map \(\phi:a\mapsto a^p\) is called the \textbf{absolute Frobenius map}.
For \(E/k\), the map
\[
\sigma:E\longrightarrow E,\qquad a\longmapsto a^q=\phi^r(a)
\]
is the relative Frobenius map over \(k\).
\end{definition}

\begin{proposition}\label{prop:frobenius-order}
The map \(\sigma\) is an automorphism of \(E\), fixes \(k\) pointwise, and has order exactly \(m\).
Its fixed elements are precisely those of \(k\):
\[
\Fix(\sigma)=\{a\in E:a^q=a\}=k.
\]
\end{proposition}
\begin{proof}
Lemma \ref{lem:freshman} shows that \(\sigma\) preserves addition, multiplication, and \(1\). If \(\sigma(a)=\sigma(b)\),
then \((a-b)^q=0\), so \(a=b\). Thus \(\sigma\) is injective; an injection on a finite set is surjective, so \(\sigma\) is an automorphism.
For \(a\in k\), we have \(a^q=a\), so \(\sigma\) fixes \(k\).

Every element of \(E\) satisfies \(a^{q^m}=a\), and hence \(\sigma^m=\id\).
If \(\sigma^t=\id\) for some \(1\leq t<m\), all \(q^m\) elements of \(E\) would be roots of \(X^{q^t}-X\). Since \(q^m>q^t\), this contradicts the root bound.
Therefore \(\ord(\sigma)=m\).

Finally, \(X^q-X\) has at most \(q\) roots, and all \(q\) elements of \(k\) are roots. Thus the fixed elements of \(\sigma\) are precisely the elements of \(k\).
\end{proof}

\section{Cyclicity of the Automorphism Group}

\begin{theorem}\label{thm:automorphisms}
\[
\Aut_k(E)=\langle\sigma\rangle
=\{\id,\sigma,\ldots,\sigma^{m-1}\},
\qquad \Aut_k(E)\cong\Z/m\Z.
\]
\end{theorem}
\begin{proof}
These \(m\) powers of Frobenius are distinct \(k\)-automorphisms.
Conversely, choose a primitive element \(\zeta\) of \(E\), so that \(E=k(\zeta)\).
The minimal polynomial of \(\zeta\) over \(k\) has degree \(m\).
Any \(k\)-automorphism is uniquely determined by the image of \(\zeta\), which must be a root of this minimal polynomial.
The polynomial has at most \(m\) roots, so there are at most \(m\) such automorphisms. Comparing the two bounds proves the assertion.
\end{proof}

\begin{corollary}
If \(|E|=p^N\), every field automorphism fixes the prime field, and
\[
\Aut(E)=\langle a\mapsto a^p\rangle\cong\Z/N\Z.
\]
For the same field with \(q^m\) elements, \(\Aut(E)\) has order \(rm\), whereas \(\Aut_k(E)\) has order \(m\).
\end{corollary}

\par\Needspace{12\baselineskip}
\section{Subfields and Divisors of the Extension Degree}

\begin{lemma}[Divisibility of power polynomials]\label{lem:power-divisibility}
If \(d\mid m\), then
\[
X^{q^d}-X\mid X^{q^m}-X.
\]
\end{lemma}
\begin{proof}
Write \(m=ds\) and set \(A=X^{q^d}-X\). In characteristic \(p\),
\[
\sum_{j=0}^{s-1}A^{q^{dj}}
=\sum_{j=0}^{s-1}\bigl(X^{q^{d(j+1)}}-X^{q^{dj}}\bigr)
=X^{q^{ds}}-X.
\]
Every term on the left is divisible by \(A\), proving the claim.
\end{proof}

\begin{theorem}[Classification of subfields]\leavevmode\phantomsection\label{thm:subfields}
For each \(d\mid m\), the set
\[
E_d=\{a\in E:a^{q^d}=a\}=\Fix(\sigma^d)
\]
is the unique subfield of \(E\) containing \(k\) and having \(q^d\) elements.
Conversely, every subfield of \(E\) containing \(k\) arises this way.
\end{theorem}
\begin{proof}
If \(d\mid m\), Lemma \ref{lem:power-divisibility} shows that \(X^{q^d}-X\) divides
\(X^{q^m}-X\), which splits completely in \(E\). Hence the former also splits completely in \(E\).
It has no repeated roots, so it has exactly \(q^d\) roots. The closure argument from Chapter~\ref{chap:foundations} shows that its root set \(E_d\) is a subfield.
Every element of \(k\) is fixed by the \(q\)-th power map, and hence by the \(q^d\)-th power map. Thus \(k\subseteq E_d\).

If \(B\) is any subfield containing \(k\), set \([B:k]=d\). Counting gives \(|B|=q^d\),
and the tower law gives \(d\mid m\). Each \(b\in B\) satisfies \(b^{q^d}=b\), so \(B\subseteq E_d\).
The two sets have the same cardinality, so \(B=E_d\). This proves both uniqueness and completeness of the classification.
\end{proof}

Within this specified extension \(E\), we often abbreviate \(E_d\) to \(\F_{q^d}\).
For \(d,e\mid m\), the tower law and the root-set description immediately give
\[
E_d\subseteq E_e\quad\Longleftrightarrow\quad d\mid e.
\]
Thus field inclusion corresponds precisely to divisibility of integers.

\par\Needspace{12\baselineskip}
\section{Fixed Fields of Iterates, Intersections, and Composita}

\begin{proposition}\label{prop:fixed-gcd}
For any integer \(t\geq0\), set \(g=\gcd(m,t)\), with \(\gcd(m,0)=m\). Then
\[
\Fix(\sigma^t)=E_g.
\]
\end{proposition}
\begin{proof}
For \(t=0\), both sides are \(E\). Suppose \(t>0\).
B\'ezout's identity for integers gives integers \(u,v\) such that \(g=um+vt\).
Since \(\sigma\) is an automorphism, its negative powers are defined. If \(\sigma^t(a)=a\), then, together with \(\sigma^m(a)=a\), this gives
\[
\sigma^g(a)=(\sigma^m)^u(\sigma^t)^v(a)=a.
\]
Conversely, since \(g\mid t\), the equality \(\sigma^g(a)=a\) implies \(\sigma^t(a)=a\).
Therefore \(\Fix(\sigma^t)=\Fix(\sigma^g)=E_g\).
\end{proof}

B\'ezout's identity for integers used here follows from integer division by the same argument as Lemma \ref{lem:bezout}:
take the least positive value among the nonzero integer linear combinations and apply division with remainder.

\begin{corollary}\leavevmode\phantomsection\label{cor:intersection-compositum}
Let \(d,e\mid m\). Within \(E\),
\[
E_d\cap E_e=E_{\gcd(d,e)},\qquad
E_dE_e=E_{\lcm(d,e)}.
\]
Here \(E_dE_e\) denotes the smallest subfield containing both subfields, called their \textbf{compositum}.
\end{corollary}
\begin{proof}
An element is fixed by both \(\sigma^d\) and \(\sigma^e\) if and only if it is fixed by \(\sigma^{\gcd(d,e)}\).
The same B\'ezout argument as above gives the intersection formula.

Set \(\ell=\lcm(d,e)\). Since \(E_\ell\) contains \(E_d\) and \(E_e\), it contains their compositum.
By the subfield classification, the compositum is some \(E_t\). Its containment of \(E_d,E_e\) gives \(d\mid t\) and \(e\mid t\),
so \(\ell\mid t\). Conversely, \(E_t\subseteq E_\ell\) implies \(t\mid\ell\), and hence \(t=\ell\).
\end{proof}

\begin{remark}[Subfields and automorphism subgroups]
A subfield \(E_d\) containing \(k\) corresponds to the group of automorphisms fixing it pointwise,
\[
\Aut_{E_d}(E)=\langle\sigma^d\rangle,
\qquad |\Aut_{E_d}(E)|=m/d.
\]
Indeed, \(\sigma^d\) is the relative Frobenius over \(E_d\), so Theorem \ref{thm:automorphisms} applies.
Every subgroup of \(\langle\sigma\rangle\) has the form \(\langle\sigma^d\rangle\) for some \(d\mid m\).
This gives an inclusion-reversing bijection between subfields and automorphism subgroups: the Galois correspondence for finite fields.
\end{remark}

\begin{example}[Subfields of a degree-twelve extension]
In \(E=\F_{2^{12}}\), over \(\F_2\), the possible subfield degrees are \(1,2,3,4,6,12\), with exactly one subfield of each degree. For example,
\[
\F_{2^4}\cap\F_{2^6}=\F_{2^2},\qquad
\F_{2^4}\F_{2^6}=\F_{2^{12}}.
\]
The absolute Frobenius \(\phi(a)=a^2\) has order twelve, and \(\Fix(\phi^8)=\F_{2^4}\), since \(\gcd(12,8)=4\).
\end{example}

\section{Orbits, Conjugates, and Minimal Polynomials}

\phantomsection\label{def:orbit-conjugacy}
Let \(\tau\) be a permutation of a finite set. The \textbf{orbit} of an element \(a\) is \(\{a,\tau(a),\tau^2(a),\ldots\}\), obtained by repeatedly applying \(\tau\).
The least positive integer \(d\) for which \(\tau^d(a)=a\) is the period of \(a\); its orbit has exactly \(d\) elements.
A subset preserved setwise by the permutation is called stable. If it consists of one orbit, the action on that subset is transitive.
On the other hand, two algebraic elements are \textbf{conjugate} over a base field if they have the same minimal polynomial over that field.
The following theorem identifies this algebraic condition with membership in the same Frobenius orbit.

\begin{theorem}[Frobenius orbits and conjugate roots]\label{thm:orbit}
Let \(\alpha\in E\), and set \(d=[k(\alpha):k]\). Then \(d\mid m\), and the least positive period of \(\alpha\) under \(\sigma\) is exactly \(d\). Its minimal polynomial is
\[
m_{\alpha,k}(X)=\prod_{i=0}^{d-1}(X-\alpha^{q^i}).
\]
These \(d\) roots are distinct.
\end{theorem}
\begin{proof}
The tower law gives \(d\mid m\). Since \(\sigma^m=\id\), there is a least positive \(t\) such that \(\sigma^t(\alpha)=\alpha\).
If \(0\le i<j<t\) and \(\sigma^i(\alpha)=\sigma^j(\alpha)\),
applying \(\sigma^{-i}\) gives \(\alpha=\sigma^{j-i}(\alpha)\), contrary to the minimality of \(t\).
Thus the first \(t\) orbit elements are distinct.

Set \(G(X)=\prod_{i=0}^{t-1}(X-\sigma^i(\alpha))\).
Applying \(\sigma\) coefficientwise shifts each root to the next, with the last returning to \(\alpha\), so the resulting polynomial is again \(G\).
Every coefficient of \(G\) therefore lies in \(\Fix(\sigma)=k\).
Since \(G(\alpha)=0\), we have \(m_{\alpha,k}\mid G\), and hence \(d\le t\).

Conversely, because the coefficients of \(m_{\alpha,k}\) lie in \(k\),
\[
m_{\alpha,k}(\sigma^i(\alpha))
=\sigma^i\bigl(m_{\alpha,k}(\alpha)\bigr)=0
\qquad(0\le i<t).
\]
This exhibits \(t\) distinct roots of \(m_{\alpha,k}\). The root bound gives \(t\le d\), so \(t=d\).
The monic polynomial \(m_{\alpha,k}\) divides the monic polynomial \(G\) of the same degree; hence they are equal.
\end{proof}

Together with Proposition~\ref{prop:minpoly}, this theorem identifies three independently defined quantities:
\[
\boxed{\dim_k k(\alpha)
 =\deg m_{\alpha,k}
 =\min\{j\ge1:\sigma^j(\alpha)=\alpha\}.}
\]
The first is the number of coordinates in a basis, the second is the least degree of a nonzero polynomial relation, and the third is the length of the Frobenius orbit.
Thus the extension degree, defined as a vector-space dimension, is also a Frobenius period over a finite field.

\phantomsection\label{def:normal-element}
An element whose distinct Frobenius conjugates form a basis of the whole extension over the base field is called a \textbf{normal element}.
This notion is developed in Chapter~\ref{chap:normal}. We first compare it with field generators and primitive elements.

\begin{example}[Three kinds of generation for an element]\label{v100:ex-three-generation-nine}
In \(E=\F_9=\F_3(\alpha)\), with \(\alpha^2=-1\), we have \(\alpha\notin\F_3\),
so \(E=\F_3(\alpha)\). However, \(\alpha^4=1\), and its multiplicative order is four, less than eight.
Also, \(\alpha^3=-\alpha\), so its two conjugates are linearly dependent. Thus \(\alpha\) is neither primitive nor normal.

Set \(\beta=1+\alpha\). Since \(\beta^4=-1\), it has multiplicative order eight and is primitive.
Its two conjugates are \(1+\alpha,1-\alpha\), whose coordinate matrix in the basis \((1,\alpha)\) is
\[
\begin{pmatrix}1&1\\1&-1\end{pmatrix},
\qquad \det=-2=1\ne0\quad\text{in }\F_3.
\]
Hence \(\beta\) is also normal.
\end{example}

\leavevmode\phantomsection\label{def:perfect-field}\ignorespaces A field over which every irreducible polynomial is separable is called \textbf{perfect}.

\begin{corollary}[Irreducible polynomials over finite fields have no repeated roots]\label{cor:separable}
Let \(f\in k[X]\) be monic and irreducible of degree \(d\). It splits completely in an extension with \(q^d\) elements,
and its roots \(\alpha,\alpha^q,\ldots,\alpha^{q^{d-1}}\) are distinct.
Consequently, finite fields are perfect.
\end{corollary}
\begin{proof}
Take \(\alpha=\overline X\) in \(k[X]/(f)\), so that \(f=m_{\alpha,k}\).
Apply Theorem \ref{thm:orbit} to this extension with \(q^d\) elements.
\end{proof}

\begin{example}[Subfields and periods in the field with sixty-four elements]\label{exm:f64-structure}
Let \(E=\F_{64}\) be obtained by the root-set construction, and set \(\sigma(x)=x^2\). Its subfields are precisely \(\F_2,\F_4,\F_8,E\), with inclusions shown below. Edges indicate inclusion, and neither intermediate subfield contains the other.
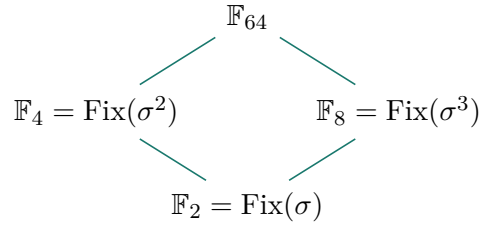
\begin{center}
\begin{tikzpicture}[every node/.style={font=\small},line width=0.6pt]
\node (top) at (0,2.5) {\(\F_{64}\)};
\node (left) at (-2,1.25) {\(\F_4=\Fix(\sigma^2)\)};
\node (right) at (2,1.25) {\(\F_8=\Fix(\sigma^3)\)};
\node (bottom) at (0,0) {\(\F_2=\Fix(\sigma)\)};
\draw[Jade] (bottom)--(left)--(top)--(right)--(bottom);
\end{tikzpicture}
\end{center}
Their intersection is \(\F_2\), and their compositum is \(E\). The possible least Frobenius periods of elements are \(1,2,3,6\).
The elements of period six are exactly \(E\setminus(\F_4\cup\F_8)\), numbering \(64-4-8+2=54\).
For any such \(\alpha\), the orbit length over \(\F_4\) is three, and that over \(\F_8\) is two:
the relative Frobenius over \(\F_4\) advances two steps along the absolute Frobenius orbit, while the relative Frobenius over \(\F_8\) advances three.
These conclusions follow first from fixed points and periods. In Section~\ref{sec:f64-computation}, we choose a specific irreducible polynomial and compute coordinates, traces, and dual bases.
\end{example}

\par\Needspace{12\baselineskip}
\section{Orbit Periods and Irreducible Factors of Power Polynomials}

\begin{theorem}\label{thm:divisibility-main}
Let \(f\in k[X]\) be monic and irreducible, with \(\deg f=d\). For every \(n\geq1\),
\[
\boxed{\ f(X)\mid X^{q^n}-X\quad\Longleftrightarrow\quad d\mid n.\ }
\]
\end{theorem}
\begin{proof}
Take \(\alpha=\overline X\) in \(L=k[X]/(f)\), and set \(\rho(a)=a^q\).
Since \(f\) is the minimal polynomial of \(\alpha\),
\[
f\mid X^{q^n}-X\quad\Longleftrightarrow\quad \rho^n(\alpha)=\alpha.
\]
By Theorem \ref{thm:orbit}, the least positive period of \(\alpha\) is \(d\).
Divide \(n\) by \(d\): \(n=ad+b\), with \(0\leq b<d\). Then \(\rho^n(\alpha)=\rho^b(\alpha)\), which equals \(\alpha\) if and only if \(b=0\).
Thus the required condition is exactly \(d\mid n\).
\end{proof}

\begin{remark}[The full collection of orbits and a single orbit]
The polynomial \(X^{q^n}-X\) collects all orbits whose lengths divide \(n\), whereas an irreducible polynomial records just one orbit.
An individual orbit is generally not closed under addition or multiplication. For example, the two roots of \(X^2+X+1\) in \(\F_4\) have product one,
which does not belong to that orbit.
Taking the union of all orbits whose lengths divide \(n\) gives the full fixed-point set of \(\sigma^n\).
Chapter~\ref{chap:period-layers} further decomposes this set into exact-period layers and separates them by polynomial elimination.
\end{remark}

In any common finite splitting extension, \(\F_{q^n}\) consists exactly of all elements whose periods divide \(n\).
The set of elements of period exactly \(n>1\) does not contain one, and hence is not itself a field.
Both the field structure of a root set and divisibility of degrees are governed by Frobenius: preservation of the operations makes the fixed-point set closed,
while the orbit theorem identifies the extension degree of an element with its least period.

\begin{corollary}[Irreducible factorization of the power polynomial]\label{cor:factor-product}
Let \(\Irr_q(d)\) denote the set of all monic irreducible polynomials of degree \(d\) in \(k[X]\). Then
\[
X^{q^n}-X=\prod_{d\mid n}\ \prod_{f\in\Irr_q(d)}f(X).
\]
Each factor occurs exactly once.
\end{corollary}
\begin{proof}
The preceding theorem determines precisely which irreducible polynomials occur as factors.
Since the derivative is \(-1\), there are no repeated roots, and hence no repeated irreducible factors.
Both sides are monic, so the equality follows.
\end{proof}

\begin{corollary}\label{cor:gcd-poly}
For positive integers \(a,b\), the monic greatest common divisor in \(k[X]\) is
\[
\gcd\bigl(X^{q^a}-X,\,X^{q^b}-X\bigr)=X^{q^{\gcd(a,b)}}-X.
\]
\end{corollary}
\begin{proof}
By the preceding corollary, the common irreducible factors on the left are exactly those whose degrees divide both \(a\) and \(b\),
equivalently, those whose degrees divide \(\gcd(a,b)\). Their product is the right-hand side.
\end{proof}

\section{Descent of Root Sets and Decompositions of Function Algebras}
\label{sec:root-descent}

Preservation of the field operations makes the fixed points of Frobenius a subfield. For a monic polynomial with simple roots, whether its coefficients lie in the base field is determined by Frobenius stability of its root set.
The latter principle both identifies irreducible polynomials and turns factorization into a direct-product decomposition of a function algebra.
Here ``descent'' means that a polynomial initially written over an extension has all its coefficients in the base field.

\begin{theorem}[Descent of root sets]\label{thm:root-descent}
Let \(k=\F_q\subseteq L=\F_{q^N}\), let \(\sigma(a)=a^q\), and let \(S\subseteq L\) be a nonempty finite set.
Set \(R_S(X)=\prod_{a\in S}(X-a)\). Then
\[
R_S\in k[X]\quad\Longleftrightarrow\quad \sigma(S)=S.
\]
When these conditions hold, the monic irreducible factors of \(R_S\) correspond bijectively to the Frobenius orbits in \(S\).
\end{theorem}
\begin{proof}
Apply \(\sigma\) to every coefficient, and denote the resulting polynomial by \(\sigma(R_S)\). Then
\[
\sigma(R_S)(X)=\prod_{a\in S}(X-\sigma(a))=R_{\sigma(S)}(X).
\]
A monic polynomial with simple roots is uniquely determined by its root set. Hence \(\sigma(R_S)=R_S\) if and only if \(\sigma(S)=S\). The elements of \(L\) fixed by \(\sigma\) are precisely those of \(k\),
so the former equality holds if and only if every coefficient belongs to \(k\). Finally, apply Theorem \ref{thm:orbit} to each orbit:
the root polynomial \(\prod_{\alpha\in\mathcal O}(X-\alpha)\) of an orbit \(\mathcal O\) is the minimal polynomial over \(k\) of any of its elements.
Distinct orbits are disjoint, so the resulting factors are distinct.
\end{proof}

Stability here means setwise preservation; it does not require each root to be fixed. For a monic nonconstant separable polynomial,
having coefficients in the base field is equivalent to stability of the root set; under this condition, irreducibility is equivalent to transitivity of Frobenius on the roots.
The assumption that the polynomial is monic fixes its leading coefficient, which the root set alone cannot determine. For example, if \(c\in L\setminus k\), the root set \(\{0\}\) of \(cX\) is stable,
but \(cX\notin k[X]\).

For a monic polynomial with repeated roots, coefficient descent is still characterized by the root multiset: Frobenius must preserve both roots and their multiplicities.
A single root orbit, however, no longer implies irreducibility. For example, the square \(P^2\) of an irreducible polynomial \(P\) still has just one orbit of distinct roots,
but is reducible. Irreducible polynomials over finite fields are automatically separable, so the simple-root condition in the orbit criterion cannot be omitted.

\Needspace{8\baselineskip}
\phantomsection\label{def:equivariant-map}
Let sets \(S,T\) carry maps \(\tau:S\to S\) and \(\rho:T\to T\), respectively.
A map \(h:S\to T\) satisfying
\[
h\circ\tau=\rho\circ h,
\]
is called \textbf{equivariant} with respect to \(\tau,\rho\).
This equality says that \(h\) respects the specified iteration; induction gives \(h\circ\tau^j=\rho^j\circ h\) for every \(j\ge0\).
If both \(\tau,\rho\) are permutations, it holds for every integer \(j\).
More generally, if a group \(G\) acts on both \(S,T\), a map \(h\) is \(G\)-equivariant if
\[
h(g\cdot s)=g\cdot h(s)\qquad(g\in G,\ s\in S).
\]
For the action of a cyclic group generated by a specified permutation, it suffices to check the equality for the generator.
When both maps are denoted by \(\sigma\), the condition is abbreviated to \(h\sigma=\sigma h\), with the two occurrences of \(\sigma\) acting on their respective sets.
Equivariance expresses compatibility with the specified actions; it does not itself require \(h\) to preserve addition or multiplication.

\phantomsection\label{def:algebra-language}
We shall use both maps of sets and algebras of functions. A \(k\)-vector space equipped with an associative multiplication and an identity,
with multiplication \(k\)-linear in each variable, is called a \(k\)-\textbf{algebra}.
An algebra homomorphism is a \(k\)-linear map preserving multiplication and the identity. The function algebras below are commutative,
with pointwise operations; operations in a direct product of algebras are likewise performed coordinatewise.

\par
\begin{figure}[!htbp]
\centering
\begin{tikzpicture}[
  vpoint/.style={circle,draw=Jade,fill=CoverMint,minimum size=8mm,inner sep=2pt},
  varrow/.style={-{Stealth[length=2mm]},Jade,semithick},
  vbox/.style={draw=MutedSlate,rounded corners=2pt,inner sep=5pt,align=center},
  every node/.style={font=\small}]

\node[vpoint] (z) at (0,0) {\(0\)};
\node[vpoint] (o) at (3,0) {\(1\)};
\node[vpoint] (w) at (6,0) {\(\omega\)};
\node[vpoint] (ww) at (8.5,0) {\(\omega^2\)};
\draw[varrow] (z) edge[loop above] node[above] {\(\sigma\)} (z);
\draw[varrow] (o) edge[loop above] node[above] {\(\sigma\)} (o);
\draw[varrow] (w) to[bend left=25] node[above] {\(\sigma\)} (ww);
\draw[varrow] (ww) to[bend left=25] (w);
\node at (0,-.85) {\(a\in\F_2\)};
\node at (3,-.85) {\(b\in\F_2\)};
\node at (6,-.85) {\(c\in\F_4\)};
\node at (8.5,-.85) {\(c^2\)};
\node[vbox] (f0) at (0,-1.8) {\(\F_2\)};
\node[vbox] (f1) at (3,-1.8) {\(\F_2\)};
\node[vbox] (f2) at (7.25,-1.8) {\(\F_4\)};
\draw[varrow,dashed,Clay] (0,-1.1)--(f0);
\draw[varrow,dashed,Clay] (3,-1.1)--(f1);
\draw[varrow,dashed,Clay] (7.25,-1.1)--(f2);

\end{tikzpicture}
\caption{Solid arrows indicate squaring on roots; dashed arrows record the independent values allowed on each orbit. The two-cycle contributes one copy of the four-element field.}\label{v102:fig-four-root-components}
\end{figure}
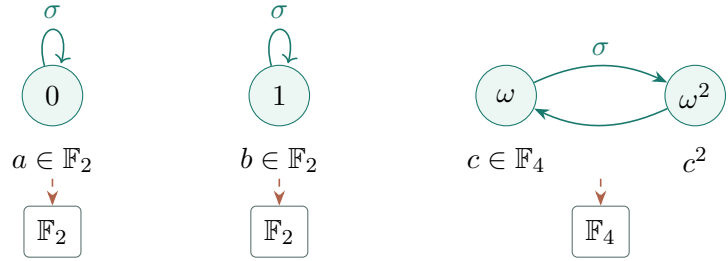

\Needspace{7\baselineskip}
\noindent\phantomsection\label{v102:par-four-root-components}\textbf{Four roots and three components.}\quad
Take \(k=\F_2\) and \(L=\F_4=\F_2(\omega)\), where \(\omega^2+\omega+1=0\).
Squaring fixes \(0,1\) and interchanges \(\omega,\omega^2\). Thus
\[
R_S(X)=X^4-X=X(X+1)(X^2+X+1),\qquad S=L.
\]
Before introducing the general function model, consider the values of a polynomial over the base field.
For \(P\in\F_2[X]\), we have \(P(\omega^2)=P(\omega)^2\), so its four values have the form
\[
\bigl(P(0),P(1),P(\omega),P(\omega^2)\bigr)=(a,b,c,c^2),
\quad a,b\in\F_2,\quad c\in\F_4.
\]
Figure~\ref{v102:fig-four-root-components} displays this restriction orbit by orbit.

\FloatBarrier
The Chinese remainder theorem allows \(a,b,c\) to be prescribed independently. Hence
\[
\F_2[X]/(X^4-X)\longrightarrow\F_2\times\F_2\times\F_4,
\qquad \overline P\longmapsto(P(0),P(1),P(\omega))
\]
is an algebra isomorphism. For example, \(\overline X\) corresponds to \((0,1,\omega)\), and \(\overline{X^2+X}\) to \((0,0,1)\).
There are four roots and three field components; the quotient algebra has dimension \(1+1+2=4\) and cardinality \(2\cdot2\cdot4=16\).
One initial value determines all values on an orbit, while distinct orbits allow independent choices. The next theorem extends this calculation to any stable root set.

\begin{theorem}[The equivariant function model]\label{thm:equivariant-functions}
Under the hypotheses of Theorem \ref{thm:root-descent}, suppose that \(\sigma(S)=S\), and define
\[
\mathcal A(S)=\{u:S\to L:\ u(\sigma(a))=\sigma(u(a))\text{ for every }a\in S\}.
\]
Give this set its \(k\)-algebra structure by pointwise addition and multiplication.\phantomsection\label{v85term:def-evaluation}\label{v85term:def-evaluation-map}%
Sending a polynomial \(P\) to the function \(a\mapsto P(a)\) is called \textbf{evaluation}; the resulting \textbf{evaluation map} induces an isomorphism
\[
k[X]/(R_S)\ \xrightarrow{\ \sim\ }\ \mathcal A(S).
\]
If the orbit lengths are \(d_1,\ldots,d_s\), then, after choosing orbit representatives,
\[
\mathcal A(S)\simeq\prod_{j=1}^s\F_{q^{d_j}}.
\]
\phantomsection\label{v85term:def-evaluation-set}\label{v85term:def-evaluation-point}\ignorespaces
Here \(S\) is called the \textbf{evaluation set}, or set of evaluation points, and its elements are the \textbf{evaluation points}.
The elements satisfying \(u^q=u\) form a subalgebra isomorphic to \(k^s\).
\end{theorem}
\begin{proof}
For \(f\in k[X]\), the identity \(f(a^q)=f(a)^q\) gives equivariance of its evaluation function. The kernel consists of polynomials vanishing on all of \(S\),
which form precisely the ideal \((R_S)\). Thus the induced map is injective.
Given \(u\in\mathcal A(S)\), Lagrange interpolation (Corollary~\ref{cor:interpolation-basic}) gives a unique \(P\in L[X]\) with \(\deg P<|S|\) and \(P(a)=u(a)\) for every \(a\in S\). For each \(a\in S\),
\[
\sigma(P)(\sigma(a))=\sigma(P(a))=\sigma(u(a))=u(\sigma(a)).
\]
Since \(\sigma\) permutes \(S\), the polynomials \(P\) and \(\sigma(P)\) agree on \(S\), and both have degree less than \(|S|\). Uniqueness of interpolation gives \(\sigma(P)=P\), so \(P\in k[X]\). This proves surjectivity.

On an orbit of length \(d_j\), an equivariant function is determined by its value at one representative \(a_j\).
Compatibility after one complete traversal of the orbit requires \(u(a_j)^{q^{d_j}}=u(a_j)\). Since \(d_j\mid N\),
the allowed values are exactly the elements of \(\F_{q^{d_j}}\) in \(L\). Different orbits admit independent choices, giving the direct-product isomorphism.
Finally, \(u^q=u\) requires every value to lie in \(k\); equivariance then makes the value constant along each orbit.
Each orbit contributes one copy of \(k\), so the fixed subalgebra is \(k^s\).
\end{proof}

For the set of all \(N\)th roots of unity with \(p\nmid N\), this equivariant-function model is the finite-field-valued Fourier coordinate model on a cyclic group. Proposition~\ref{v96:prop-fourier-descent} will write the same condition explicitly for cyclic codes.

\begin{example}\label{exm:descent-product}
In \(\F_2[X]\), let \(f=(X^2+X+1)(X^3+X+1)\). Both factors are irreducible, since their degrees are two and three, respectively, and neither has a root in the base field.
Their roots can all be placed in \(\F_{64}\), where they form two orbits of lengths two and three. Thus
\[
\F_2[X]/(f)\simeq\F_4\times\F_8,
\qquad \dim_{\F_2}\{a:a^2=a\}=2.
\]
The root set itself has only five elements and is not a field, but its function algebra records both orbits and their lengths.
\end{example}

This distinction reappears in factorization algorithms: first find the Frobenius-fixed elements of the quotient algebra, then use them to separate the root orbits.
The dimension of the fixed subalgebra therefore directly measures the number of irreducible factors.

\section{Finite Reduced Algebras and Abstract Frobenius Orbits}
\label{sec:finite-reduced-algebras}

The equivariant function model of a root set can be developed without choosing a polynomial. This unifies different quotient-ring representations
and also treats products of finite fields that cannot be generated by one element. Throughout this section, algebras are nonzero commutative finite-dimensional \(k=\F_q\)-algebras with identity,
and homomorphisms preserve the identity and fix \(k\) pointwise.

\begin{definition}\leavevmode\phantomsection\label{def:reduced-algebra}
An element \(a\) is \textbf{nilpotent} if some positive power of \(a\) is zero.
An algebra with no nonzero nilpotent elements is called \textbf{reduced}.
An element satisfying \(e^2=e\) is an \textbf{idempotent}; two idempotents whose product is zero are \textbf{orthogonal}.
\end{definition}

\par\Needspace{7\baselineskip}
\noindent\phantomsection\label{v102:par-component-idempotents}\textbf{Selecting the components.}\quad
Continue in \(A=\F_2[X]/(X^4-X)\), interpreting the following polynomials as residue classes. Put
\[
e_0=1+X^3,\qquad e_1=X+X^2+X^3,\qquad e_2=X+X^2.
\]
Their values at the three orbit representatives are

\[
\begin{array}{c|ccc}
 &0&1&\omega\\\hline
e_0&1&0&0\\
e_1&0&1&0\\
e_2&0&0&1
\end{array}
\qquad
e_i^2=e_i,\quad e_i e_j=0\ (i\ne j),\quad e_0+e_1+e_2=1.
\]
The identities can be checked by evaluation, since the evaluation map is injective.
Multiplication by \(e_i\) retains the \(i\)-th component and sets the others to zero.
Thus \(e_0A\simeq\F_2\), \(e_1A\simeq\F_2\), and \(e_2A\simeq\F_4\); the identity of each component is its own \(e_i\).
This is the concrete meaning of splitting an algebra by orthogonal idempotents.

For an eight-element root set, take \(\beta^3+\beta+1=0\). Then \(\F_8^\times=\langle\beta\rangle\), and the squaring orbits are
\[
\{0\},\quad\{1\},\quad\{\beta,\beta^2,\beta^4\},\quad
\{\beta^3,\beta^6,\beta^5\}.
\]
Consequently, \(\F_2[X]/(X^8-X)\simeq\F_2\times\F_2\times\F_8\times\F_8\).
The component degrees sum to eight, but the algebra has \(2^8\) elements.
The next theorem shows that the same component description applies to any finite-dimensional commutative reduced algebra, even before a polynomial presentation is chosen.

\begin{theorem}[Structure of finite reduced algebras]\label{thm:finite-reduced-algebra}
Let \(A\ne0\) be a finite-dimensional commutative \(k=\F_q\)-algebra with identity. The following conditions are equivalent:
\begin{enumerate}
\item \(A\) is reduced;
\item \(\sigma_A:a\mapsto a^q\) is an algebra automorphism;
\item there is a positive integer \(N\) such that \(a^{q^N}=a\) for every \(a\in A\);
\item there are positive integers \(s,d_1,\ldots,d_s\) such that
\[
A\simeq\prod_{i=1}^s\F_{q^{d_i}}.
\]
\end{enumerate}
These degrees, including their multiplicities, are uniquely determined by \(A\), up to reordering.
\end{theorem}
\begin{proof}
Commutativity makes \(\sigma_A\) a \(k\)-algebra homomorphism. If \(A\) is reduced, then \(a^q=0\) implies \(a=0\),
so \(\sigma_A\) is injective. Since \(A\) is a finite set, it is bijective and has finite order as a permutation.
This proves the first implication and then the second.
If the third condition holds and \(a\) is nilpotent, repeatedly using \(a=a^{q^N}\) gives \(a=a^{q^{jN}}=0\) for sufficiently large \(j\).
Thus the first three conditions are equivalent.

To deduce the fourth from the third, induct on \(\dim_k A\). If every nonzero element is invertible, then \(A\) is a finite field and the conclusion holds.
Otherwise choose a nonzero nonunit \(a\). The kernel of evaluation \(k[X]\to A\), \(X\mapsto a\), is generated by a monic polynomial \(h\),
and \(k[a]\simeq k[X]/(h)\). Since \(a^{q^N}=a\), we have \(h\mid X^{q^N}-X\), so \(h\) is square-free.
Moreover, \(h(0)=0\); otherwise \(h(a)=0\) would express an inverse of \(a\) in \(k[a]\).
Thus \(h=Xh_1\), where \(h_1(0)\ne0\), and \(a\ne0\) implies \(\deg h_1\ge1\).
The Chinese remainder theorem (Theorem~\ref{thm:crt-basic}) gives
\[
k[a]\simeq k\times k[X]/(h_1).
\]
Take the idempotent \(e\) corresponding to \((0,1)\), so that \(e\ne0,1\). The map
\[
A\longrightarrow eA\times(1-e)A,\qquad x\longmapsto(ex,(1-e)x)
\]
is an algebra isomorphism with inverse \((u,v)\mapsto u+v\); the identities of the two factors are \(e\) and \(1-e\), respectively.
Each factor is reduced and has dimension strictly less than \(\dim_k A\). Induction gives the required product of finite fields.
A product of finite fields has no nonzero nilpotent elements, so the fourth condition implies the first.

In a product of fields, every coordinate of an idempotent is zero or one.
The nonzero idempotents that cannot be expressed as a sum of two nonzero orthogonal idempotents are exactly the coordinate identities of the factors.
Every algebra isomorphism permutes these idempotents and induces isomorphisms between the corresponding field factors.
The orders of these finite fields determine the \(d_i\) and their multiplicities.
\end{proof}

\leavevmode\phantomsection\label{def:finite-etale-algebra}\ignorespaces A finite-dimensional commutative reduced algebra over a finite field is also called a \textbf{finite \'etale algebra}.
Over a finite field, this terminology is equivalent to the preceding description as a product of finite fields.
The correspondence below is the concrete finite-field form of the general relation between finite \'etale algebras and Galois actions; see \cite[Lemma 58.2.2]{StacksEtale}.

Let \(S\) be a nonempty finite set and \(\tau:S\to S\) a permutation.
Choose \(N\) divisible by all orbit lengths, and define, in \(L=\F_{q^N}\),
\[
\mathcal A(S,\tau)=\{u:S\to L:\ u(\tau s)=u(s)^q\}.
\]
By the proof of Theorem \ref{thm:equivariant-functions}, an orbit of length \(d\) contributes one copy of \(\F_{q^d}\).
Enlarging \(L\) therefore introduces no new function values: compatibility around the orbit already requires every value to lie in the corresponding \(\F_{q^d}\).
The set \(S\) need not be embedded in a field; it records the orbits and their lengths.

We say that a finite reduced algebra \(A\) splits over a finite extension \(L/k\) if every field factor of \(A\) embeds in \(L\) as a \(k\)-extension.
\phantomsection\label{v85term:def-valued-point}\ignorespaces
A \(k\)-algebra homomorphism \(A\to L\) is called an \textbf{\(L\)-valued point} of \(A\). We also call it an evaluation point of the algebra, since each algebra element can be evaluated under this homomorphism; this is the algebraic counterpart of choosing a field element as a polynomial evaluation point.
For a polynomial quotient, a homomorphism \(k[X]/(f)\to L\) is determined by the image of \(\overline X\), which is a root of \(f\).

\begin{theorem}[The contravariant correspondence between algebras and orbits]\label{thm:finite-orbit-duality}
Let \(A\) be a finite reduced \(k\)-algebra, and choose \(L=\F_{q^N}\) so that the degrees of all its field factors divide \(N\). Write
\[
S_L(A)=\operatorname{Hom}_{k\text{-alg}}(A,L),\qquad
\tau(\iota)=\sigma_L\circ\iota.
\]
Then \(|S_L(A)|=\dim_k A\), and evaluation induces an isomorphism
\[
A\xrightarrow{\ \sim\ }\mathcal A(S_L(A),\tau),\qquad
a\longmapsto\bigl(\iota\mapsto\iota(a)\bigr).
\]
If \(A,B\) split over the same \(L\), composition gives a bijection
\[
\operatorname{Hom}_{k\text{-alg}}(A,B)
\ \simeq\ 
\operatorname{Hom}_{\tau}\bigl(S_L(B),S_L(A)\bigr).
\]
The right-hand side consists of set maps commuting with the Frobenius permutations. This correspondence reverses the order of composition.
\end{theorem}
\begin{proof}
Write \(A=\prod_i E_i\), with \([E_i:k]=d_i\), and let \(e_i\) be the identity in the \(i\)-th coordinate.
These idempotents satisfy
\[
e_i^2=e_i,\qquad e_i e_j=0\quad(i\ne j),\qquad \sum_i e_i=1.
\]
For any homomorphism \(\iota:A\to L\), the only idempotents in the field \(L\) are zero and one, so \(\iota(e_i)\in\{0,1\}\).
If \(i\ne j\), then \(\iota(e_i)\iota(e_j)=\iota(e_i e_j)=0\), so at most one image is one.
Since \(\sum_i\iota(e_i)=\iota(1)=1\), exactly one image is one; denote its index by \(i\).
Then \(\iota\) vanishes on all other components, while its restriction to \(E_i\) preserves the identity and is therefore a field embedding.
Thus \(\iota\) is the composition of the projection \(A\to E_i\) and an embedding \(E_i\to L\).
There are \(d_i\) such embeddings, cyclically permuted by Frobenius. Hence \(S_L(A)\) decomposes into orbits of lengths \(d_i\).

Evaluation functions satisfy equivariance. If \(a\ne0\), some coordinate is nonzero; any embedding of that coordinate into \(L\) gives a nonzero evaluation.
Thus evaluation is injective. Its target has dimension equal to the sum of the orbit lengths, \(\sum_i d_i=\dim_k A\), so it is also surjective.

A homomorphism \(v:A\to B\) induces \(v^*(\iota)=\iota\circ v\). The identity
\((\iota\circ v)(a)^q=\iota(v(a))^q\) shows that this map commutes with \(\tau\).
Conversely, given an equivariant map \(h:S_L(B)\to S_L(A)\), send a function \(u\) to \(u\circ h\).
This operation, called pullback along \(h\), is a homomorphism of equivariant function algebras.
The evaluation isomorphisms already proved then give a homomorphism \(A\to B\). Evaluation at each point shows that the two constructions are mutually inverse and that the correspondence reverses composition.
\end{proof}

\begin{example}[Evaluating a homomorphism point by point]\label{v103:ex-contravariant-map}
Take \(k=\F_2\), \(L=\F_4=\F_2(\omega)\), where \(\omega^2+\omega+1=0\), and write \(\sigma(x)=x^2\).
Let \(A=L\), \(B=L\times L\), and define the \(k\)-algebra homomorphism
\[
v:A\longrightarrow B,\qquad a\longmapsto(a,a^2).
\]
We have \(S_L(A)=\{\id,\sigma\}\). If \(\pi_1,\pi_2\) are the coordinate projections, the four evaluation points of \(B\) are
\(\iota_{j,e}=\sigma^e\circ\pi_j\), where \(j=1,2\) and \(e=0,1\). Composition gives
\[
\begin{array}{c|cccc}
\iota&\iota_{1,0}&\iota_{1,1}&\iota_{2,0}&\iota_{2,1}\\\hline
\iota\circ v&\id&\sigma&\sigma&\id\\
(\iota\circ v)(\omega)&\omega&\omega^2&\omega^2&\omega
\end{array}
\]
For instance, \(\iota_{2,0}(v(\omega))=\pi_2(\omega,\omega^2)=\omega^2\).
To evaluate on \(A\), first apply \(v\) to reach \(B\), then use an evaluation point of \(B\).
Thus the point map goes from \(S_L(B)\) to \(S_L(A)\).
Frobenius exchanges each pair \(\iota_{j,0},\iota_{j,1}\), and also exchanges their images \(\id,\sigma\); hence this reverse map is equivariant.
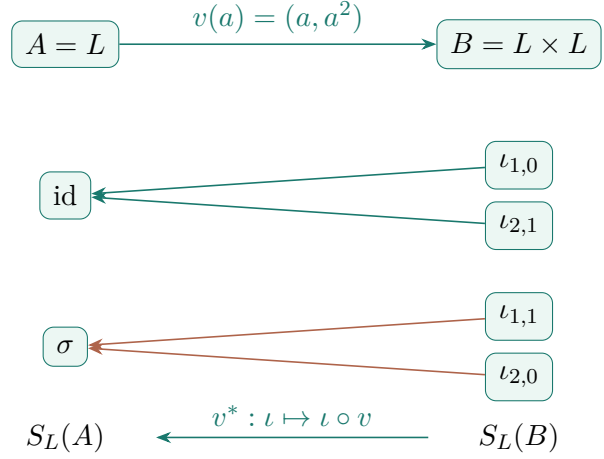
\begin{figure}[!htbp]
\centering
\begin{tikzpicture}[every node/.style={font=\small},
 box/.style={draw=Jade,fill=CoverMint,rounded corners,inner sep=5pt},
 ar/.style={-{Stealth[length=2mm]},semithick,Jade}]
\node[box] (A) at (-3,2) {\(A=L\)};
\node[box] (B) at (3,2) {\(B=L\times L\)};
\draw[ar] (A) -- node[above]{\(v(a)=(a,a^2)\)} (B);
\node[box] (id) at (-3,0) {\(\id\)};
\node[box] (sg) at (-3,-2) {\(\sigma\)};
\node[box] (p10) at (3,.4) {\(\iota_{1,0}\)};
\node[box] (p21) at (3,-.4) {\(\iota_{2,1}\)};
\node[box] (p11) at (3,-1.6) {\(\iota_{1,1}\)};
\node[box] (p20) at (3,-2.4) {\(\iota_{2,0}\)};
\draw[ar] (p10) -- (id); \draw[ar] (p21) -- (id);
\draw[ar,Clay] (p11) -- (sg); \draw[ar,Clay] (p20) -- (sg);
\node at (-3,-3.2) {\(S_L(A)\)};
\node at (3,-3.2) {\(S_L(B)\)};
\draw[ar] (1.8,-3.2) -- node[above]{\(v^*: \iota\mapsto\iota\circ v\)} (-1.8,-3.2);
\end{tikzpicture}
\caption{The algebra homomorphism goes right and the evaluation-point map goes left. The image of each of the four points is computed by composition.}\label{v103:fig-contravariant-map}
\end{figure}
\FloatBarrier
\end{example}

\begin{corollary}\label{cor:orbit-algebra-dictionary}
In the notation of Theorem \ref{thm:finite-orbit-duality}, a finite reduced \(k\)-algebra \(A\) is a field
if and only if \(\tau:\iota\mapsto\sigma_L\circ\iota\) acts transitively on the set \(S_L(A)\) of evaluation points.
Write \(A^{\sigma_A}=\{a\in A:\sigma_A(a)=a\}\) for the fixed subalgebra of \(\sigma_A\).
If \(S_L(A)\) has \(s\) orbits under \(\tau\), then \(A^{\sigma_A}\simeq k^s\), and \(A\) has exactly \(2^s\) idempotents.
\end{corollary}
\begin{proof}
The evaluation isomorphism identifies the field factors of \(A\) with the orbits in \(S_L(A)\). Thus \(A\) is a field if and only if there is one orbit.
In the equivariant function model, \(\sigma_A\) acts pointwise on values. A fixed function takes values in \(k\) and is constant on each orbit,
so \(A^{\sigma_A}\simeq k^s\). An idempotent function takes only the values zero and one. Choosing any collection of orbits on which the value is one
uniquely specifies an idempotent, giving \(2^s\) choices.
\end{proof}

\begin{example}[Actions on elements and on evaluation points]\label{exm:two-frobenius-actions}
Take \(A=L=\F_{q^3}\). On the set of elements of \(A\), the map \(\sigma_A\) fixes every element of \(k\),
so it cannot be transitive. On the set of evaluation points, however,
\[
S_L(A)=\{\id,\sigma_A,\sigma_A^2\},
\qquad \tau:\id\longmapsto\sigma_A\longmapsto\sigma_A^2\longmapsto\id.
\]
This three-cycle has no fixed point but has just one orbit. Thus \(\dim_k A^{\sigma_A}=1\) counts
the orbits on the evaluation points. The two actions are on algebra elements and algebra homomorphisms, respectively, and must be distinguished.
\end{example}

\begin{example}[Abstract orbits beyond univariate root sets]\label{exm:three-fixed-points}
The algebra \(\F_2^3\) corresponds to three fixed points, but cannot be written as \(\F_2[X]/(f)\).
Indeed, every \(a\in\F_2^3\) satisfies \(a^2=a\), so \(\F_2[a]\) has at most four elements and cannot equal the eight-element algebra \(\F_2^3\).
Moreover, in any extension field of characteristic two, squaring has only the two fixed points zero and one, so these three abstract fixed points cannot simultaneously be realized as distinct roots.
A general criterion for generation by one element is given in Proposition~\ref{v104:prop-monogenic-reduced} below.

\begin{figure}[!htbp]
\centering
\begin{tikzpicture}[
  vpoint/.style={circle,draw=Jade,fill=CoverMint,minimum size=8mm,inner sep=2pt},
  varrow/.style={-{Stealth[length=2mm]},Jade,semithick},
  vbox/.style={draw=MutedSlate,rounded corners=2pt,inner sep=5pt,align=center},
  every node/.style={font=\small}]

\node[vbox] (A) at (0,0) {\(\F_2\times\F_2\times\F_2\)};
\foreach \j/\x in {1/4,2/6,3/8} {
  \node[vpoint] (p\j) at (\x,0) {\(\pi_{\j}\)};
  \draw[varrow] (p\j) edge[loop above] (p\j);
}
\draw[varrow,Clay,dashed] (A.east)--(3.15,0);
\node at (6,-.85) {\(\tau(\pi_j)=\sigma_{\F_2}\circ\pi_j=\pi_j\)};

\end{tikzpicture}
\caption{The three coordinate projections are distinct evaluation points, each fixed by Frobenius. They are algebra homomorphisms, not three distinct fixed elements in a field.}\label{v102:fig-abstract-fixed-points}
\end{figure}
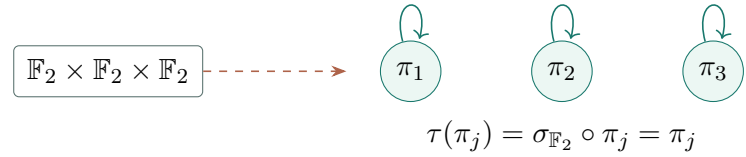
\FloatBarrier
\end{example}

\leavevmode\phantomsection\label{def:irreducible-count}\ignorespaces
Write \(I_q(d)=|\Irr_q(d)|\) for the number of monic irreducible polynomials of degree \(d\) over \(\F_q\).
Generation of a finite reduced algebra by one element also depends on how many root orbits of each length are available.
\begin{proposition}[A criterion for a finite reduced algebra to be monogenic]\label{v104:prop-monogenic-reduced}
Let \(A\) be a nonzero finite-dimensional commutative reduced \(\F_q\)-algebra, written as
\[
A\simeq\prod_{d\ge1}(\F_{q^d})^{m_d},\qquad m_d\ge0,
\]
with only finitely many nonzero \(m_d\). Then \(A\) is a \textbf{monogenic algebra}, meaning that \(A=\F_q[a]\) for some \(a\in A\), if and only if
\[
\boxed{\ m_d\le I_q(d)\quad\text{for every }d\ge1.\ }
\]
Equivalently, precisely under these inequalities, \(A\) can be written as \(\F_q[X]/(f)\) for a monic squarefree polynomial \(f\).
\end{proposition}
\begin{proof}
First suppose that \(A=\F_q[a]\), and let \(a_{d,j}\) be the coordinate of \(a\) in the \(j\)th factor \(\F_{q^d}\).
Projection gives \(\F_q[a_{d,j}]=\F_{q^d}\), so its monic minimal polynomial \(P_{d,j}\) has degree \(d\).
The kernel of the evaluation map \(\F_q[X]\to A\) is the intersection of the ideals \((P_{d,j})\), hence is generated by their least common multiple. Therefore
\[
\dim_{\F_q}\F_q[a]=\deg\lcm_{d,j}P_{d,j}.
\]
If two polynomials \(P_{d,j}\) were equal, the right-hand side would be strictly smaller than \(\sum_d dm_d=\dim_{\F_q}A\), a contradiction.
Thus each degree requires \(m_d\) distinct irreducible polynomials, proving the inequalities.
Conversely, for every degree choose \(m_d\) distinct monic irreducible polynomials, and choose a root of each in its corresponding field \(\F_{q^d}\) as \(a_{d,j}\).
All the \(P_{d,j}\) are pairwise coprime, so the Chinese remainder theorem gives
\[
\F_q[X]\Big/\Big(\prod_{d,j}P_{d,j}\Big)
\simeq\prod_{d,j}\F_q[X]/(P_{d,j})\simeq A,
\]
where \(\overline X\) corresponds to \(a=(a_{d,j})\). Hence \(A=\F_q[a]\). The product is squarefree, which also proves the final equivalent formulation.
\end{proof}
Proposition~\ref{v107:prop-generating-tuples} will use tuple orbits to give the existence criterion and exact count for several generators.

The assumption that the algebra is reduced is also essential. Every homomorphism from \(k[\varepsilon]/(\varepsilon^2)\) to a field sends \(\varepsilon\) to zero,
so this algebra has the same evaluation points as \(k\), but additional nilpotent structure.
The root set and its permutation record separable points; multiplicities require further algebraic structure.
Corollary \ref{cor:reduced-trace-pairing} will show that reducedness is also equivalent to nondegeneracy of the algebra trace pairing.
For a monogenic quotient algebra with repeated factors, a derivative factor records the degeneracy of this pairing.

\section{Recovering Algebra Decompositions from Frobenius Fixed Spaces}
\label{fra:sec-recovery}

Theorem~\ref{thm:finite-reduced-algebra} expresses a reduced algebra as a product of finite fields. Conversely, the fixed spaces of Frobenius iterates determine those field components and their multiplicities.
The relevant data are the integer dimensions of the fixed spaces, not merely the characteristic polynomial of Frobenius.

Let \(A\) be a nonzero finite-dimensional commutative reduced \(k=\F_q\)-algebra, set \(N=\dim_k A\), and write \(\sigma_A(x)=x^q\).
Grouping components of the same degree, write
\[
 A\simeq\prod_{j=1}^{N}\F_{q^j}^{\,a_j},\qquad
 a_j\ge0,\qquad \sum_{j=1}^{N}ja_j=N.
\]
The equality \(a_j=0\) means that the corresponding component does not occur.
\leavevmode\phantomsection\label{fra:def-fixed-dimensions}\ignorespaces
For \(r\ge1\), define the \textbf{Frobenius fixed-space dimension} by
\[
 b_r(A)=\dim_k\ker(\sigma_A^r-\id).
\]

\begin{proposition}\leavevmode\phantomsection\label{fra:prop-gcd-dimensions}
We have
\begin{equation}\label{fra:eq-gcd-dimensions}
 b_r(A)=\sum_{j=1}^{N}a_j\gcd(j,r)\qquad(r\ge1).
\end{equation}
In particular, \(b_1(A)=\sum_j a_j\) is the number of field components.
\end{proposition}
\begin{proof}
On the component \(\F_{q^j}\), the fixed field of \(\sigma_A^r\) is \(\F_{q^{\gcd(j,r)}}\), of \(k\)-dimension \(\gcd(j,r)\).
Fixed points can be chosen independently in each component. Adding the dimensions gives the formula.
\end{proof}

The inversions below are performed over the integers or rational numbers; dimensions are not reduced modulo the characteristic \(p\).
Let \(\varphi(e)\) count the integers among \(1,\ldots,e\) that are coprime to \(e\).
\leavevmode\phantomsection\label{def:mobius-function}\ignorespaces
The number-theoretic M\"obius function \(\mu\) is defined by \(\mu(1)=1\), by \(\mu(e)=(-1)^h\) when \(e\) is a product of \(h\) distinct primes,
and by \(\mu(e)=0\) if \(e\) has a squared prime factor. It satisfies
\[
 \sum_{t\mid e}\mu(t)=\begin{cases}1,&e=1,\\0,&e>1.\end{cases}
\]
For \(e>1\), this follows by summing over subsets of its distinct prime factors and expanding \((1-1)^h\).

\begin{theorem}[Recovering field components by two M\"obius inversions]\label{fra:thm-recovery}
For \(1\le e\le N\), set
\begin{equation}\label{fra:eq-first-inversion}
 B_e(A)=\frac1{\varphi(e)}\sum_{t\mid e}\mu(e/t)b_t(A).
\end{equation}
Then
\begin{equation}\label{fra:eq-second-inversion}
 B_e(A)=\sum_{\substack{j\le N\\e\mid j}}a_j,
 \qquad
 \boxed{\quad a_j=\sum_{h=1}^{\lfloor N/j\rfloor}\mu(h)B_{jh}(A).\quad}
\end{equation}
Thus \(N\) and \(b_1(A),\ldots,b_N(A)\) uniquely determine all field components of \(A\).
\end{theorem}
\begin{proof}
Grouping the integers \(1,\ldots,v\) by their greatest common divisor with \(v\) gives \(\sum_{e\mid v}\varphi(e)=v\).
Applying this to \(v=\gcd(j,r)\), equation \eqref{fra:eq-gcd-dimensions} becomes
\[
 b_r(A)=\sum_{e\mid r}\varphi(e)\sum_{\substack{j\le N\\e\mid j}}a_j.
\]
Replace \(r\) by \(t\), multiply by \(\mu(r/t)\), and sum over \(t\mid r\). Interchanging the finite sums leaves the inner sum \(\sum_{h\mid r/e}\mu(h)\),
which vanishes unless \(e=r\). This gives the first expression for \(B_e\).
Substitute it into the right-hand side of the second inversion:
\[
 \sum_{h\le N/j}\mu(h)B_{jh}(A)
 =\sum_{\ell\le N/j}a_{j\ell}\sum_{h\mid\ell}\mu(h)=a_j.
\]
Thus all multiplicities can be recovered, and Theorem~\ref{thm:finite-reduced-algebra} then determines the algebra structure.
\end{proof}

On the set of evaluation points \(S_L(A)\) from Theorem~\ref{thm:finite-orbit-duality}, let \(z_r\) be the number of fixed points of \(\tau^r\). An orbit of length \(j\) contributes \(j\) fixed points if \(j\mid r\), and none otherwise. Thus
\[
 z_r=\sum_{\substack{j\mid r\\j\le N}}ja_j,\qquad
 a_j=\frac1j\sum_{d\mid j}\mu(j/d)z_d\quad(1\le j\le N).
\]
By contrast, \(b_r(A)=\sum_j\gcd(j,r)a_j\) in \eqref{fra:eq-gcd-dimensions} counts the orbits of \(\tau^r\) on the evaluation points, or equivalently the dimension of the fixed subspace of the algebra.
Both sets of integers recover the field components, but they record different fixed objects; Exercise~\ref{fra:ex-two-counts} compares them for an algebra with quadratic and cubic components.

\begin{corollary}[Frobenius similarity type and algebra isomorphism type]\label{fra:cor-similarity}
Two nonzero finite-dimensional commutative reduced \(k\)-algebras are isomorphic if and only if their Frobenius maps are similar as \(k\)-linear operators.
\end{corollary}
\begin{proof}
An algebra isomorphism commutes with taking \(q\)-th powers, so it conjugates one Frobenius operator to the other.
Conversely, similar operators have equal underlying dimensions and equal fixed-space dimensions for every iterate. Theorem~\ref{fra:thm-recovery} then recovers the same field components and multiplicities.
\end{proof}

What is recovered here is the algebra's isomorphism type. A given Frobenius-equivariant linear isomorphism need not preserve multiplication;
Example~\ref{exm:linear-not-field} gives a concrete counterexample.

\begin{example}[Equal characteristic polynomials and different fixed spaces]\leavevmode\phantomsection\label{fra:ex-characteristic}
Over \(k=\F_2\), compare \(A=\F_{16}\) and \(A'=\F_4\times\F_4\).
Both algebras have dimension four, and their Frobenius maps satisfy \(\sigma_A^4=\id\) and \(\sigma_{A'}^2=\id\), respectively.
In characteristic two, \(T^4-1=(T-1)^4\) and \(T^2-1=(T-1)^2\). Hence all eigenvalues of both operators are one, and both characteristic polynomials are \((T-1)^4\).
Their first four fixed-space dimensions, however, are
\[
 (b_1,b_2,b_3,b_4)=(1,2,1,4),\qquad(2,4,2,4).
\]
By \eqref{fra:eq-first-inversion}, the first data set gives \((B_1,B_2,B_3,B_4)=(1,1,0,1)\), and hence \(a_4=1\).
The second gives \((2,2,0,0)\), and hence \(a_2=2\). All other multiplicities are zero.
The fixed-space dimensions distinguish the two algebras even though the characteristic polynomial does not.
\end{example}

\begin{example}[Nilpotent multiplication cannot be recovered this way]\label{fra:ex-nonreduced}
Let \(q\ge3\), and take
\[
 A=k[\varepsilon]/(\varepsilon^3)\times k,\qquad
 A'=k[\varepsilon]/(\varepsilon^2)\times k[\varepsilon]/(\varepsilon^2).
\]
Both Frobenius maps fix the constants in each component and annihilate all nilpotent terms. As operators on four-dimensional \(k\)-spaces, both are therefore similar to \(\operatorname{diag}(1,1,0,0)\).
The algebra \(A\) contains a nilpotent element \((\varepsilon,0)\) with nonzero square, whereas every nilpotent element of \(A'\) has square zero.
Thus the algebras are not isomorphic, and the assumption of reducedness cannot be omitted.
\end{example}

If the matrix of \(\sigma_A\) in a \(k\)-basis is \(M\), then \(b_r(A)=N-\rank_k(M^r-I)\).
For \(A=k[X]/(g)\), with \(g\) monic and square-free, use the basis \(1,\overline X,\ldots,\overline X^{\deg g-1}\). With columns indexed from zero, column \(i\) is the coefficient column of \(X^{iq}\bmod g\), for \(0\le i<\deg g\).
These integer data can thus be recovered directly from \(g\), giving the number of irreducible factors of each degree. The individual factors can then be separated by the methods in Chapter~\ref{chap:algorithms}.

\section{Frobenius Descent of Multivariate Linear Factors}
\label{sec:multivariate-linear-descent}

The essence of root-set descent is that a Frobenius-stable collection of factors gives a polynomial over the base field, while one orbit gives an irreducible factor.
For multivariate polynomials, the same idea applies when the polynomial splits completely into linear factors over an extension.
The phrase \textbf{polynomial splitting into linear factors} describes the form of the factorization; the linear factors are not assumed to be distinct.
The following orbit formulation gives one form of Carlitz's criterion; see \cite[Chapter 13]{Agou}.

Let \(E=\F_{q^N}\), and let \(\sigma\) act coefficientwise on \(E[X_1,\ldots,X_m]\), fixing the variables.
For a nonconstant affine linear polynomial \(M=c_0+\sum_{j=1}^m c_jX_j\), take the smallest index \(j\in\{1,\ldots,m\}\) for which \(c_j\ne0\),
and call \(L=c_j^{-1}M\) its \textbf{normalized linear factor}.
Linear polynomials differing by a nonzero scalar give the same normalized factor, and normalization commutes with \(\sigma\).

\begin{theorem}[Frobenius descent of linear factors]\label{thm:multivariate-linear-descent}
Let \(k=\F_q\), let \(L=a_0+\sum_{j=1}^m a_jX_j\) be a normalized linear factor as above, and suppose that
\[
 E_0=k(a_0,\ldots,a_m)=\F_{q^s}.
\]
Then the coefficientwise Frobenius orbit of \(L\) has length exactly \(s\), and
\begin{equation}\label{eq:multivariate-orbit-product}
 F_L=\prod_{i=0}^{s-1}L^{(q^i)}
 =\prod_{i=0}^{s-1}\left(a_0^{q^i}+\sum_{j=1}^m a_j^{q^i}X_j\right)
\end{equation}
belongs to \(k[X_1,\ldots,X_m]\), is irreducible in this ring, and has total degree \(s\). Parenthesized superscripts indicate powers applied only to coefficients.
Moreover,
\[
 s=\lcm\bigl([k(a_0):k],\ldots,[k(a_m):k]\bigr).
\]
Conversely, suppose that a nonconstant \(f\in k[X_1,\ldots,X_m]\) splits completely into linear factors over some finite extension.
Then \(f\) is irreducible in \(k[X_1,\ldots,X_m]\) if and only if \(f=cF_L\) for some normalized linear factor \(L\) and \(c\in k^\times\).
\end{theorem}
\begin{proof}
The equality \(\sigma^r(L)=L\) holds if and only if every \(a_j\) is fixed by \(\sigma^r\).
By the subfield classification, the least positive such \(r\) is \(s\). The compositum formula also shows that \(s\) is the least common multiple of the coefficient degrees.
Thus \(L,\sigma(L),\ldots,\sigma^{s-1}(L)\) are distinct.
The coefficient of the same variable is one in each factor, so they are pairwise nonassociate: no two differ by a nonzero scalar.
Since \(\sigma\) cyclically permutes these factors, \(\sigma(F_L)=F_L\), showing that its coefficients lie in \(k\).

Here a nonzero nonunit \(P\) is called a \textbf{prime element} if \(P\mid GH\) implies \(P\mid G\) or \(P\mid H\).
Every nonconstant linear polynomial is a prime element of the multivariate polynomial ring: after dividing by a nonzero variable coefficient,
passing to the quotient by that polynomial eliminates the corresponding variable. The quotient is a polynomial ring in the remaining variables, hence an integral domain.
If \(F_L=GH\), with \(G,H\in k[X_1,\ldots,X_m]\), then \(L\) divides \(G\) or \(H\).
Suppose without loss of generality that \(L\mid G\). Repeatedly applying \(\sigma\), and using \(\sigma(G)=G\), shows that every \(\sigma^i(L)\) divides \(G\).
These prime elements are pairwise nonassociate, so successive use of primality gives \(F_L\mid G\).
Therefore \(H\) is constant, proving irreducibility of \(F_L\).

Conversely, let \(f\) be irreducible. Choose one of its linear factors, normalize it to \(L\), and work in a finite extension \(E\) containing all coefficients.
Since \(\sigma(f)=f\), all distinct conjugates of \(L\) divide \(f\). Thus \(f=F_LQ\), with \(Q\in E[X_1,\ldots,X_m]\).
Applying \(\sigma\) and canceling the nonzero factor \(F_L\) gives \(\sigma(Q)=Q\), so \(Q\in k[X_1,\ldots,X_m]\).
Irreducibility forces \(Q=c\in k^\times\). The reverse implication follows from the irreducibility of \(F_L\) already proved.
\end{proof}

In one variable, a normalized linear factor is \(X-\alpha\), and \eqref{eq:multivariate-orbit-product} is the orbit product for the minimal polynomial.
In several variables, Frobenius instead organizes the coefficients of linear factors; the product of conjugates gives the same descent mechanism.
In particular, the linear factors of an irreducible polynomial splitting into linear factors are necessarily distinct in a splitting extension.

\begin{example}\leavevmode\phantomsection\label{exm:multivariate-linear-descent}
In \(\F_4\), take \(\omega^2+\omega+1=0\) and set \(L=X+\omega Y+\omega^2\). Then
\[
 F_L=(X+\omega Y+\omega^2)(X+\omega^2Y+\omega)
 =X^2+XY+Y^2+X+Y+1.
\]
The normalized coefficients generate \(\F_4\), so \(F_L\) is irreducible in \(\F_2[X,Y]\).
Nevertheless, \(F_L(1,1)=0\): a multivariate polynomial with a zero defined over the base field need not have a linear factor over that field.
\end{example}

\section{Counting Irreducible Polynomials}

Continue to write \(I_q(d)\) for the number of monic irreducible polynomials of degree \(d\). Comparing degrees in Corollary \ref{cor:factor-product} gives
\begin{equation}\label{eq:counting}
q^n=\sum_{d\mid n}dI_q(d).
\end{equation}
This identity also has a direct interpretation: group the elements of \(\F_{q^n}\) by their minimal polynomials.
An irreducible polynomial of degree \(d\) contributes \(d\) distinct roots, and the possible degrees are precisely the positive divisors of \(n\).

Inversion using the M\"obius function \(\mu\) yields
\begin{equation}\label{eq:mobius}
I_q(n)=\frac1n\sum_{d\mid n}\mu(d)q^{n/d}.
\end{equation}
Here \(\mu\) is the number-theoretic M\"obius function defined on page~\pageref{def:mobius-function}. We next express the same inversion in the general form of a periodic-point count.

\begin{example}
\[
I_q(1)=q,\qquad I_q(2)=\frac{q^2-q}{2},\qquad
I_q(3)=\frac{q^3-q}{3},\qquad I_2(4)=\frac{16-4}{4}=3.
\]
Thus \(\F_2\) has exactly one monic irreducible quadratic and two monic irreducible cubics,
in agreement with the explicit polynomials listed earlier.
\end{example}

Formula \eqref{eq:mobius} counts irreducible polynomials of a prescribed degree; the denominator \(n\) is the number of elements in each orbit of exact period \(n\).

\subsection{M\"obius Inversion for Periodic-Point Counts}
The divisor-sum identity for the M\"obius function implies, for functions with values in any additive abelian group,
\[
A(n)=\sum_{d\mid n}B(d)\quad\Longleftrightarrow\quad
B(n)=\sum_{d\mid n}\mu(n/d)A(d).
\]
Substitute the first relation into the right-hand side of the second and interchange the finite sums:
\[
\sum_{d\mid n}\mu(n/d)A(d)
=\sum_{e\mid n}B(e)\sum_{\substack{d\mid n\\ e\mid d}}\mu(n/d)
=\sum_{e\mid n}B(e)\sum_{t\mid n/e}\mu(t)=B(n).
\]
Substitution in the other direction likewise recovers \(A(n)\), so the relations are equivalent.
For values in a multiplicative abelian group, the same computation gives the multiplicative form
\[
A(n)=\prod_{d\mid n}B(d)\quad\Longleftrightarrow\quad
B(n)=\prod_{d\mid n}A(d)^{\mu(n/d)}.
\]
Objects whose periods divide \(n\) decompose into disjoint layers of exact periods \(d\mid n\). M\"obius inversion recovers an individual layer from these cumulative data.

\subsection{Several coordinates and the number of algebra generators}
The monogenic criterion~\ref{v104:prop-monogenic-reduced} assigns each field component to a different orbit of elements.
When one coordinate supplies too few orbits, we may use several coordinates. The requirement is still to separate all evaluation points.

\begin{proposition}[Generators of a finite reduced algebra: existence and exact counts]\label{v107:prop-generating-tuples}
Let \(A\simeq\prod_{d\ge1}(\F_{q^d})^{m_d}\) be a nonzero finite-dimensional commutative reduced \(k=\F_q\)-algebra, and let \(r\ge1\). Put
\begin{equation}\label{v107:eq-tuple-orbit-count}
J_{q,r}(d)=\frac1d\sum_{e\mid d}\mu(d/e)q^{re}=I_{q^r}(d).
\end{equation}
Then \(A\) is generated by \(r\) elements as a unital \(k\)-algebra if and only if
\[
\boxed{m_d\le J_{q,r}(d)\quad\text{for every }d.}
\]
The number of ordered \(r\)-tuples satisfying \(k[a_1,\ldots,a_r]=A\) is
\begin{equation}\label{v107:eq-generating-tuple-count}
\boxed{\prod_{d\ge1}d^{m_d}\bigl(J_{q,r}(d)\bigr)_{m_d}.}
\end{equation}
Here \((J)_m=J(J-1)\cdots(J-m+1)\), \((J)_0=1\), and the value is zero when \(m>J\).
Thus the least positive number of generators is the least \(r\) satisfying all these inequalities; no convention about an empty generating tuple is needed.
\end{proposition}
\begin{proof}
Frobenius acts coordinatewise on \((\F_{q^d})^r\). A tuple has exact period \(e\mid d\) if and only if its coordinates jointly generate \(\F_{q^e}\).
If \(B(e)\) counts tuples of exact period \(e\), then \(q^{rd}=\sum_{e\mid d}B(e)\). M\"obius inversion and division by the \(d\) points in each orbit give~\eqref{v107:eq-tuple-orbit-count}. Its equality with \(I_{q^r}(d)\) is an equality of counting formulas, not an additional identification of field models.

Collect the coordinates of \(a_1,\ldots,a_r\) in the \(j\)-th degree-\(d\) field component into a tuple \(z_{d,j}\).
Generation of \(A\) requires each \(z_{d,j}\) to generate that field jointly.
For such a tuple, evaluation \(k[X_1,\ldots,X_r]\to\F_{q^d}\) is surjective and has a maximal ideal as kernel.
Two kernels coincide if and only if the tuples are related by an isomorphism over the base field. Such finite-field isomorphisms are powers of Frobenius, so this means precisely that the tuples lie in one orbit.
Components of different degrees cannot have the same kernel.
If two kernels coincide, their values are constrained by that isomorphism and cannot be prescribed independently. If all kernels differ, they are pairwise comaximal, and the Chinese remainder theorem makes joint evaluation surjective onto \(A\).
Indeed, two different maximal ideals have sum equal to the whole ring, and the Chinese remainder construction is the same as in one variable.

For each \(d\), we must therefore choose distinct length-\(d\) orbits for the \(m_d\) fixed field components, then choose a point on each orbit.
There are \((J_{q,r}(d))_{m_d}\) ordered orbit choices and \(d^{m_d}\) point choices. Choices for different degrees are independent, proving the count and the existence criterion.
The count is for tuples in a fixed algebra; it is not divided by algebra automorphisms.
\end{proof}

\begin{example}[Three components and two coordinates]\label{v107:ex-three-components}
The algebra \(A=\F_2^3\) is not monogenic, since \(J_{2,1}(1)=2<3\). However,
\[
a=(0,0,1),\qquad b=(0,1,0)
\]
generate \(A\): the three coordinate idempotents are \(1+a+b,b,a\). The coordinate pairs \((0,0),(0,1),(1,0)\) are distinct.
The formula gives exactly \((4)_3=24\) ordered generating pairs. Also \(J_{2,2}(2)=(16-4)/2=6\), so \(\F_4\times\F_4\) has \(2^2(6)_2=120\) ordered generating pairs, although it has no single generator.
The number of abstract field components constrains a single coordinate. With several coordinates, the constraint is the supply of tuple orbits of each length.
\end{example}

\subsection{From Orbit Counts to Euler Products}
\label{sec:counting-euler}
For a commutative ring \(R\), the \textbf{formal power series ring} \(R[[u]]\) consists of expressions
\(\sum_{j\ge0}a_j u^j\), with \(a_j\in R\).
Equality means equality of coefficients in every degree. Addition is termwise, and the coefficient of \(u^m\) in a product is the finite sum
\(\sum_{i=0}^m a_i b_{m-i}\). No convergence is involved in these operations.
In \(\Q[[u]]\), for \(v\) with constant term zero, one also defines the \textbf{formal logarithm} and \textbf{formal exponential}:
\[
\log(1+v)=\sum_{j\ge1}\frac{(-1)^{j-1}v^j}{j},
\qquad \exp(v)=\sum_{j\ge0}\frac{v^j}{j!}.
\]
Since \(v\) has constant term zero, every coefficient is still determined by finitely many operations.
Denominators in the counting identities below are interpreted in \(\Q\). Although the polynomials are defined over \(\F_q\), the counting series have coefficients in the characteristic-zero field \(\Q\).
The same definitions apply when formal logarithms and exponentials are later used in \(\C[[u]]\).

\begin{proposition}\label{prop:finite-field-euler}
In the formal power series ring \(\Q[[u]]\),
\begin{equation}\label{eq:finite-field-euler}
\boxed{\quad
\prod_{d\ge1}(1-u^d)^{-I_q(d)}
=\exp\left(\sum_{r\ge1}\frac{q^r u^r}{r}\right)
=\frac1{1-qu}.\quad}
\end{equation}
\end{proposition}
\begin{proof}
A monic irreducible polynomial \(P\) of degree \(d\) can occur with any nonnegative integer exponent in a factorization,
and therefore contributes \((1-u^d)^{-1}=\sum_{j\ge0}u^{dj}\). Multiplying over all \(P\),
unique factorization ensures that every monic polynomial is counted exactly once. There are \(q^r\) monic polynomials of degree \(r\), so
\[
\prod_{P\text{ monic irreducible}}(1-u^{\deg P})^{-1}
=\sum_{r\ge0}q^r u^r=(1-qu)^{-1}.
\]
For any fixed coefficient of \(u^r\), only the finitely many irreducible polynomials of degrees at most \(r\) contribute.
Thus the infinite product is well defined as a formal power series, without a convergence condition.
Expanding \(-\log(1-z)=\sum_{j\ge1}z^j/j\) in \(\Q[[u]]\) gives
\[
\log\prod_{d\ge1}(1-u^d)^{-I_q(d)}
=\sum_{r\ge1}\frac{u^r}{r}\sum_{d\mid r}dI_q(d)
=\sum_{r\ge1}\frac{q^r u^r}{r}.
\]
The final step uses \eqref{eq:counting}; taking the formal exponential gives the middle expression.
\end{proof}

An orbit of length \(d\) contributes \(d\) fixed points of the \(r\)-th iterate exactly when \(d\mid r\),
explaining the factor \(d\) in the formal logarithm. The identity connects unique factorization, degrees of irreducible factors, and counts of Frobenius fixed points.
Section~\ref{sec:hasse-davenport} assigns trace- and norm-dependent weights to the orbits and obtains the corresponding lifting formula for Gauss sums.
Section~\ref{sec:orbit-euler} then expresses this product as the determinant of a permutation operator.

\subsection{Synchronized Orbits and the Proportion of Irreducible Polynomials}
\label{v105:subsec-necklace-density}

Formula~\eqref{eq:mobius} also allows us to compare the proportions of irreducible polynomials over different base fields.
Replace the field size temporarily by a formal parameter and define the \textbf{necklace polynomial}
\phantomsection\label{v105:def-necklace-polynomial}
\begin{equation}\label{v105:eq-necklace-polynomial}
M_n(Z)=\frac1n\sum_{d\mid n}\mu(d)Z^{n/d}\in\Q[Z]
\qquad(n\ge1).
\end{equation}
The name comes from counting orbits of cyclic words; see Exercise~\ref{v40:ex-necklace-models}.
For a prime power \(q\), we have \(M_n(q)=I_q(n)\).
For now, we regard \(M_n\) as a polynomial with rational coefficients, without requiring its argument to be a field size.

\begin{lemma}[Synchronized products and a prime recursion]\label{v105:lem-necklace-identities}
For every \(n\ge1\), the following polynomial identity holds:
\begin{equation}\label{v105:eq-necklace-product}
M_n(XY)=
\sum_{\substack{i,j\ge1\\ \lcm(i,j)=n}}
\gcd(i,j)M_i(X)M_j(Y).
\end{equation}
If \(\ell\) is a prime and \(d\ge1\), then
\begin{equation}\label{v105:eq-necklace-prime-recursion}
M_{\ell d}(X)=
\begin{cases}
\bigl(M_d(X^\ell)-M_d(X)\bigr)/\ell,&\ell\nmid d,\\[2pt]
M_d(X^\ell)/\ell,&\ell\mid d.
\end{cases}
\end{equation}
\end{lemma}
\begin{proof}
M\"obius inversion gives
\[
X^n=\sum_{i\mid n}iM_i(X),\qquad
Y^n=\sum_{j\mid n}jM_j(Y).
\]
Multiply these identities, group the terms by \(e=\lcm(i,j)\), and use \(ij=e\gcd(i,j)\):
\[
(XY)^n
=\sum_{e\mid n}e
\sum_{\lcm(i,j)=e}\gcd(i,j)M_i(X)M_j(Y).
\]
Comparing this with \((XY)^n=\sum_{e\mid n}eM_e(XY)\) and applying M\"obius inversion once more
proves~\eqref{v105:eq-necklace-product}.

If \(\ell\mid d\), then \(\ell d\) and \(d\) have the same squarefree divisors.
Keeping only terms with \(\mu(e)\ne0\) in~\eqref{v105:eq-necklace-polynomial} gives
\[
M_{\ell d}(X)
=\frac1\ell\left(\frac1d\sum_{e\mid d}\mu(e)(X^\ell)^{d/e}\right)
=\frac1\ell M_d(X^\ell).
\]
If \(\ell\nmid d\), the divisors of \(\ell d\) are of the two forms \(e\) and \(\ell e\), with \(e\mid d\),
and \(\mu(\ell e)=-\mu(e)\). Separating the two sums gives the other recursion.
\end{proof}

The least common multiple and greatest common divisor in~\eqref{v105:eq-necklace-product} have a direct orbit interpretation:
moving simultaneously around cycles of lengths \(i,j\) produces \(\gcd(i,j)\) orbits of length \(\lcm(i,j)\).
Proposition~\ref{prop:orbit-algebra-operations} gives the same mechanism for tensor products of finite reduced algebras.
We now apply it to proportions of irreducible polynomials; see~\cite{CheboluEtAl2026Necklace} for this result and further properties.

\begin{proposition}[The proportion of irreducibles increases with the field size]\label{v105:prop-irreducible-density}
For every \(n\ge2\), \(M_n(x)>0\) for \(x>1\), and the function
\[
x\longmapsto\frac{M_n(x)}{x^n}
\]
is strictly increasing on \([1,\infty)\).
Consequently, for prime powers \(q_1<q_2\),
\begin{equation}\label{v105:eq-irreducible-density}
\frac{I_{q_1}(n)}{q_1^n}
<
\frac{I_{q_2}(n)}{q_2^n}.
\end{equation}
For \(n=1\), the proportion is identically one.
\end{proposition}
\begin{proof}
We prove positivity and monotonicity simultaneously by induction on \(n\).
Since \(M_1(x)=x\), positivity holds and the ratio is identically one.
Moreover, \eqref{v105:eq-necklace-polynomial} gives \(M_n(1)=0\) for every \(n>1\).
Let \(n\ge2\), and assume the assertions for all smaller indices.
Choose a prime \(\ell\mid n\) and write \(n=\ell d\).
If \(\ell\mid d\), the recursion gives \(M_n(x)>0\) immediately for \(x>1\).
If \(\ell\nmid d\), monotonicity of the ratio for the smaller index \(d\) gives
\[
M_d(x^\ell)\ge x^{(\ell-1)d}M_d(x)>M_d(x)\qquad(x>1).
\]
The final inequality uses \(M_d(x)>0\) and \(x^{(\ell-1)d}>1\);
it therefore also covers \(d=1\), when the ratio is constant.
The recursion again yields \(M_n(x)>0\).

We now know that all \(M_i\) with \(i\le n\) are nonnegative on \([1,\infty)\)
and positive on \((1,\infty)\).
For \(x\ge1\) and \(y>1\), retain in~\eqref{v105:eq-necklace-product}
all terms with \(i=n\), together with the additional term \((i,j)=(1,n)\):
\[
\begin{aligned}
M_n(xy)
&\ge M_n(x)\sum_{j\mid n}jM_j(y)+M_1(x)M_n(y)\\
&>y^nM_n(x).
\end{aligned}
\]
The additional term is positive and, because \(n>1\), does not belong to the terms with \(i=n\).
Divide by \((xy)^n\) and put \(x=a,y=b/a\).
For any \(1\le a<b\), this gives
\[
\frac{M_n(b)}{b^n}>\frac{M_n(a)}{a^n}.
\]
This completes the induction. Substituting \(M_n(q)=I_q(n)\) gives the assertion for finite fields.
\end{proof}

\begin{example}[Splitting of one orbit and the overall proportion]\label{v105:ex-density-base-extension}
The polynomial \(X^2+X+1\) is irreducible over \(\F_2\).
If \(\F_4=\F_2(\omega)\), with \(\omega^2+\omega+1=0\), then over \(\F_4\)
\[
X^2+X+1=(X-\omega)(X-\omega^2).
\]
Nevertheless, the proportion of irreducible monic quadratics increases from
\[
\frac{I_2(2)}{2^2}=\frac14
\qquad\text{to}\qquad
\frac{I_4(2)}{4^2}=\frac6{16}=\frac38.
\]
The factorization follows one fixed root orbit as the base field is enlarged;
the proportions compare all polynomials in two different coefficient spaces.
An increase in the proportion does not require each original irreducible polynomial to remain irreducible after extension.
\end{example}

\par\Needspace{11\baselineskip}
\section*{Exercises}
\phantomsection
\addcontentsline{toc}{section}{Exercises}
\begin{exercise}\label{ex:subfields18}
List all subfields of \(\F_{3^{18}}\), and compute \(\F_{3^6}\cap\F_{3^9}\) and their compositum.
Set \(\phi(a)=a^3\). Find \(\Fix(\phi^{12})\) and the automorphism group fixing \(\F_{3^6}\).
\end{exercise}
\begin{exercise}\label{ex:embedding}
Prove that an \(\F_q\)-embedding from \(\F_{q^a}\) to \(\F_{q^b}\) exists if and only if \(a\mid b\).
When it exists, prove that there are exactly \(a\) such embeddings.
Both fields here are equipped with specified \(\F_q\)-extension structures.
\end{exercise}
\begin{exercise}[The counting formula and exclusion of subfields]\label{ex:counting}
Compute \(I_2(6)\) and the number of elements of \(\F_{64}\) having degree exactly six over \(\F_2\).
Compute the latter number independently using unions and intersections of subfields, and explain why the two calculations agree.
\end{exercise}
\begin{exercise}[Compatibility of successive base-field extensions]\label{ex:factor-extension}
Let \(f\in\F_q[X]\) be monic and irreducible of degree \(d\), and suppose that \(r\mid s\).
Compare \(\F_q\subseteq\F_{q^r}\subseteq\F_{q^s}\) inside one finite extension, and set \(g=\gcd(d,r)\) and \(h=\gcd(d/g,s/r)\).
Prove that each irreducible factor obtained over \(\F_{q^r}\) splits into exactly \(h\) irreducible factors over \(\F_{q^s}\),
each of degree \(d/(gh)\), and that \(gh=\gcd(d,s)\).
Use root indices to show that successive extension and direct extension give the same partition into factors. For \(d=12,r=2,s=8\), write out the orbits at both stages.
\end{exercise}

\begin{exercise}[Affine transformations preserve degree]\label{ex:wu-affine-degree}
For \(a\in k^\times\) and \(b\in k\), the invertible substitution \(X\mapsto aX+b\) is a one-dimensional \textbf{affine transformation}.
Let \(\alpha\) have degree \(d\) over \(k=\F_q\), and let \(a\in k^\times,b\in k\).
Prove \(k(a\alpha+b)=k(\alpha)\), and determine the monic minimal polynomial of \(a\alpha+b\).
Explain why this transformation need not preserve multiplicative order.
\end{exercise}

\begin{exercise}[Normalization and repeated factors in products of conjugates]\label{ex:multivariate-normalization}
Let \(E=\F_{q^N}\), suppose that the coefficients of a normalized linear factor \(L\) generate \(\F_{q^s}\subseteq E\), and define \(F_L\) by \eqref{eq:multivariate-orbit-product}.
For \(u\in E^\times\), prove
\[
 \prod_{i=0}^{N-1}(uL)^{(q^i)}
 =u^{(q^N-1)/(q-1)}F_L^{N/s}.
\]
In \(\F_4\), take \(\omega^2+\omega+1=0\), and compute the product of the two coefficientwise conjugates of \(M=\omega(X+Y)\).
Explain why the degree of the coefficient field before normalization does not directly give an irreducibility criterion.
\end{exercise}

\begin{exercise}[When a product of finite fields is generated by one element]\label{ex:wu2-monogenic}
Let \(A=\prod_{d\ge1}(\F_{q^d})^{m_d}\), with only finitely many \(m_d\) nonzero and at least one positive.
Prove that \(A=\F_q[a]\) for some \(a\in A\) if and only if \(m_d\le I_q(d)\) for every \(d\).
Deduce that \(\F_2^3\) cannot be represented as a quotient algebra of \(\F_2[X]\).
\end{exercise}

\begin{exercise}[Maps between orbits and algebra homomorphisms]\leavevmode\label{ex:wu3-algebra-hom}
Let \(A=\prod_{i=1}^s\F_{q^{d_i}}\) and \(B=\prod_{j=1}^t\F_{q^{e_j}}\). Prove
\[
|\operatorname{Hom}_{k\text{-alg}}(A,B)|
=\prod_{j=1}^t\left(\sum_{i:d_i\mid e_j}d_i\right).
\]
If \(A\) has \(m_d\) field factors of degree \(d\), prove that \( |\Aut_k(A)|=\prod_d d^{m_d}m_d!\).
\end{exercise}
\begin{exercise}[Fixed subalgebras of the iterates]\label{ex:wu3-algebra-fixed}
For \(A=\prod_i\F_{q^{d_i}}\), prove
\[
A^{\sigma_A^r}\simeq\prod_i\F_{q^{\gcd(d_i,r)}},\qquad
\dim_k A^{\sigma_A^r}=\sum_i\gcd(d_i,r).
\]
Explain why this dimension is the number of orbits of \(\tau^r\) on the evaluation points, rather than its number of fixed points.
\end{exercise}

\begin{exercise}[Recovering an algebra from five dimensions]\label{fra:ex-five-dimensions}
Let \(A\) be a five-dimensional commutative reduced \(\F_q\)-algebra satisfying
\[
 (b_1(A),b_2(A),b_3(A),b_4(A),b_5(A))=(2,3,4,3,2).
\]
Compute \(B_1(A),\ldots,B_5(A)\), and recover the direct-product decomposition of \(A\) into fields.
\end{exercise}

\begin{exercise}[Root-period counts and fixed spaces of the algebra]\label{fra:ex-two-counts}
Let \(P_2,P_3\in k[X]\) be monic irreducible polynomials of degrees two and three, respectively, and set \(g=P_2P_3\), \(A=k[X]/(g)\).
For \(r\ge1\), compare
\[
 z_r=\deg\gcd(g,X^{q^r}-X),\qquad b_r(A)=\dim_k\ker(\sigma_A^r-\id).
\]
Give formulas for both quantities, explain which Frobenius action each counts, and compute both when \(r=1\).
\end{exercise}

\chapter{Element Representations, Multiplicative Orders, and Polynomial Orders}
\label{chap:orders}
The multiplicative order and Frobenius period of an element describe the cyclic subgroup and the field extension it generates, respectively. Quotient-ring coordinates and multiplication matrices make these invariants accessible to polynomial computations, while repeated factors record the contribution of powers of the characteristic.

\section{Coordinates for Computing with Field Elements}\label{sec:field-coordinates}
Fix \(k=\F_q\subseteq E=\F_{q^n}\). If \(E=k(\alpha)\), every element has a unique expression
\[
x=x_0+x_1\alpha+\cdots+x_{n-1}\alpha^{n-1},\qquad x_i\in k.
\]
This is called a \textbf{polynomial basis representation}.
The relation \(m_{\alpha,k}(\alpha)=0\) reduces higher powers to lower ones.
Addition is coordinatewise, multiplication consists of polynomial multiplication followed by reduction, and inversion uses the extended Euclidean algorithm. This algorithm records the coefficients of linear combinations during successive divisions, thereby obtaining the coefficients in B\'ezout's identity.
The chosen \(\alpha\) need only generate the extension, not the multiplicative group.

\leavevmode\phantomsection\label{def:discrete-logarithm}\ignorespaces If a primitive element \(\zeta\) is chosen, every nonzero \(x\) is uniquely expressible as \(\zeta^j\), where \(j\in\Z/(q^n-1)\Z\). This is the \textbf{exponent representation}; the exponent class \(j\) is the \textbf{discrete logarithm} of \(x\) to base \(\zeta\), denoted by \(\log_\zeta x\).
In this representation, multiplication adds exponents, inversion negates them, and taking an integer power multiplies the exponent by that integer. Addition of field elements requires further information.
\leavevmode\phantomsection\label{def:zech-logarithm}\ignorespaces For example, when \(1+\zeta^t\ne0\), define the \textbf{Zech logarithm}
\[
Z(t)=\log_\zeta(1+\zeta^t).
\]
Then \(\zeta^i+\zeta^j=\zeta^{i+Z(j-i)}\). If \(1+\zeta^{j-i}=0\), the sum is zero.
In practice, addition can use a Zech logarithm table or a conversion to polynomial basis coordinates.

If \(k\subseteq L\subseteq E\), choose a basis \(u_i\) of \(L/k\) and a basis \(v_j\) of \(E/L\),
and write elements as \(\sum_{i,j}a_{ij}u_iv_j\). Lemma \ref{lem:tower} guarantees uniqueness; this is a \textbf{tower representation}.
These representations describe the same field. Choosing coordinates does not change minimal polynomials, conjugacy orbits, or element orders.

\section{The Regular Representation by Multiplication}
\label{sec:regular-representation}
Choose a basis of \(E/k\). For \(a\in E\), let \(M_a\) denote the matrix of the linear map \(x\mapsto ax\).
Then
\[
M_{a+b}=M_a+M_b,\qquad M_{ab}=M_aM_b,\qquad M_1=I.
\]
This representation is injective, since \(M_a=0\) implies \(a\cdot1=0\). For nonzero \(M_a\), its inverse is \(M_{a^{-1}}\).
Thus \(a\mapsto M_a\) embeds \(E\) into a matrix algebra, with image a subfield isomorphic to \(E\). This embedding is the \textbf{regular representation} of the field.

If \(E=k(\alpha)\) and \(m_\alpha(X)=X^n+c_{n-1}X^{n-1}+\cdots+c_0\),
then, in the basis \((1,\alpha,\ldots,\alpha^{n-1})\), the matrix \(M_\alpha=A\) is called the \textbf{companion matrix} of \(m_\alpha\):
its first \(n-1\) columns record \(\alpha\cdot\alpha^j=\alpha^{j+1}\) for \(0\le j<n-1\), and its last column is \((-c_0,\ldots,-c_{n-1})^{\mathsf t}\).
Consequently, \(E\cong k[A]\). If \(B\) is the matrix of Frobenius, then
\[
BM_aB^{-1}=M_{a^q}.
\]
Thus conjugate field elements give similar multiplication matrices. Chapter~\ref{chap:trace} uses this relation to determine traces, norms, and characteristic polynomials.

\begin{proposition}[Conjugacy of matrix models in the same dimension]\label{v48:prop-matrix-models}
Let \([E:k]=n\), and let \(V\) be an \(n\)-dimensional \(k\)-space.
If \(\rho:E\to\End_k(V)\) is a unital \(k\)-algebra homomorphism, there is a \(k\)-linear isomorphism \(S:E\to V\) such that
\[
\rho(a)=S M_a S^{-1}\qquad(a\in E).
\]
Thus every matrix model of a field in dimension equal to its extension degree is obtained from the regular representation by one simultaneous change of basis.
\end{proposition}
\begin{proof}
Choose \(0\ne v\in V\), and define \(S(x)=\rho(x)v\). If \(x\ne0\), then \(\rho(x)\) has inverse \(\rho(x^{-1})\),
so \(S(x)=0\) would force \(v=0\), a contradiction. Hence \(S\) is injective, and equality of dimensions makes it an isomorphism.
For \(a,x\in E\), we have \(\rho(a)S(x)=\rho(ax)v=S(M_a x)\), proving the conjugacy relation.
\end{proof}

The regular representation therefore preserves all multiplication relations simultaneously. For every \(h\in k[T]\), we have \(h(M_a)=M_{h(a)}\),
so the annihilating polynomials of the multiplication matrix are exactly those of the field element. For the resulting equality of minimal polynomials, see Exercise~\ref{ex:wu-matrix-polynomial}.
A change of basis alters matrix entries but preserves the same system of matrix relations.

\begin{example}
Reuse \(\F_9=\F_3(\alpha)\), \(\alpha^2=-1\), from Example~\ref{v100:ex-three-generation-nine}. Here multiplication is represented by matrices: in the basis \((1,\alpha)\),
\[
a+b\alpha\ \longleftrightarrow\ \begin{pmatrix}a&-b\\b&a\end{pmatrix},\qquad a,b\in\F_3.
\]
These nine matrices are closed under the field operations. For a nonzero matrix, the determinant \(a^2+b^2\) is nonzero in \(\F_3\).
\end{example}

\begin{example}[Companion matrices and the coordinate action of conjugacy]\label{v48:ex-companion-frobenius}
Take \(E=\F_2[\alpha]\), with \(\alpha^3=\alpha+1\). Since \(X^3+X+1\) is nonzero at \(0,1\),
it is irreducible, and \(E\) has eight elements. In the basis \((1,\alpha,\alpha^2)\), multiplication by \(\alpha\) and squaring have respective matrices
\[
A=\begin{pmatrix}0&0&1\\1&0&1\\0&1&0\end{pmatrix},\qquad
B=\begin{pmatrix}1&0&0\\0&0&1\\0&1&1\end{pmatrix}.
\]
The last column of \(A\) comes from \(\alpha^3=1+\alpha\); that of \(B\) comes from \(\alpha^4=\alpha+\alpha^2\).
Thus the multiplication matrix of \(a+b\alpha+c\alpha^2\) is \(aI+bA+cA^2\). The specific relations are
\[
A^3=A+I,\qquad A^7=I,\qquad B^3=I,\qquad BA=A^2B.
\]
The first three record the reduction relation, multiplicative order, and Frobenius period, respectively; the last records compatibility between squaring and multiplication.
Since \(A\ne I\) and \(7\) is prime, \(A\) has order seven. Since \(B\ne I\), it has order three.
Multiplication by the chosen element and Frobenius are therefore represented by distinct operators in the same coordinates.
\end{example}

\section{Multiplicative Order and Frobenius Period}\label{v48:sec-multiplicative-period}
\begin{proposition}\label{prop:degree-from-order}
Let \(\alpha\ne0\) lie in a finite extension of \(k\), with multiplicative order \(e\). Then \(\gcd(e,q)=1\), and
\[
[k(\alpha):k]=\ord_e(q).
\]
For \(e>1\), \(\ord_e(q)\) is the least positive integer \(d\) satisfying \(q^d\equiv1\pmod e\); set \(\ord_1(q)=1\).
\end{proposition}
\begin{proof}
The integer \(e\) divides the order of the ambient finite field's multiplicative group, and hence is coprime to the characteristic.
For every \(t\geq1\),
\[
\alpha^{q^t}=\alpha
\quad\Longleftrightarrow\quad \alpha^{q^t-1}=1
\quad\Longleftrightarrow\quad e\mid q^t-1.
\]
By Theorem \ref{thm:orbit}, the least positive \(t\) on the left is the degree of the element; the least positive \(t\) on the right is \(\ord_e(q)\).
If \(e=1\), then \(\alpha=1\), which has degree one.
\end{proof}

\begin{remark}
The multiplicative order of \(\alpha\) is the least positive integer \(e\) satisfying \(\alpha^e=1\).
Its Frobenius period is the least positive integer \(d\) satisfying \(\alpha^{q^d}=\alpha\).
These numbers generally differ and are related by \(d=\ord_e(q)\).
Zero has Frobenius period one but has no multiplicative order.
\end{remark}

\begin{example}
In the field with sixteen elements, a primitive element \(\alpha\) has multiplicative order fifteen and degree four over \(\F_2\).
The element \(\alpha^3\) has order five and still has degree four; \(\alpha^5\) has order three and degree two.
Thus generating the whole extension, generating the multiplicative group, and determining subfield membership are related but distinct questions.
\end{example}

\begin{proposition}[Two periods of a multiplication matrix]\label{v48:prop-two-matrix-periods}
Let \(E/k\) be finite, let \(0\ne a\in E\) have order \(e\) and degree \(d\), and let \(B\) be the Frobenius matrix in the same basis.
Then \(M_a\) has order \(e\), while the least positive integer \(t\) satisfying \(B^tM_aB^{-t}=M_a\) is \(d\).
Moreover, every orbit of \(x\mapsto ax\) on \(E^\times\) has length \(e\).
\end{proposition}
\begin{proof}
Since \(M_a^j=M_{a^j}\) and the regular representation is injective, \(M_a^j=I\) if and only if \(a^j=1\).
Also, \(B^tM_aB^{-t}=M_{a^{q^t}}\), so the second period is the Frobenius period of \(a\).
Finally, for \(x\ne0\), the equality \(a^jx=x\) is equivalent to \(a^j=1\). Hence every nonzero orbit has length \(e\).
\end{proof}

\begin{corollary}[Degrees of powers]\leavevmode\phantomsection\label{v48:cor-degree-power}
Let \(a\ne0\) have order \(e\). For any integer \(t\), set \(e_t=e/\gcd(e,t)\). Then
\[
\ord(a^t)=e_t,\qquad [k(a^t):k]=\ord_{e_t}(q).
\]
In particular, within the same algebraic closure, \(a^t\in\F_{q^r}\) if and only if \(e_t\mid q^r-1\).
\end{corollary}
\begin{proof}
The condition \((a^t)^j=1\) is equivalent to \(e\mid tj\), whose least positive solution \(j\) is \(e/\gcd(e,t)\).
Applying Proposition~\ref{prop:degree-from-order} gives the degree formula.
The subfield condition follows from \((a^t)^{q^r}=a^t\) and the fixed-point characterization of finite fields.
This includes \(t=0\), when \(e_t=1\) and \(a^t=1\).
\end{proof}

\begin{example}[Subfield descent by taking powers]\label{v48:ex-power-subfields}
Choose a primitive element \(\zeta\) in \(\F_{81}\), so that \(\ord(\zeta)=80\). Over the base field \(\F_3\),
\[
\begin{array}{c|rrrrrrr}
t&1&2&5&10&20&40&80\\\hline
\ord(\zeta^t)&80&40&16&8&4&2&1\\
{[\F_3(\zeta^t):\F_3]}&4&4&4&2&2&1&1
\end{array}
\]
For example, \(3^2\equiv9\not\equiv1\pmod{16}\) and \(3^4\equiv1\pmod{16}\),
so \(\zeta^5\) still generates all of \(\F_{81}\), although its multiplicative order is only sixteen.
By contrast, \(\zeta^{10}\) has order eight, and \(\ord_8(3)=2\), so it generates the subfield \(\F_9\) and is primitive there.
The element \(\zeta^{20}\) also generates \(\F_9\), but has order only four in that field.
Finally, \(\zeta^{40}=-1\in\F_3\). The multiplicative order of a power is first determined by an integer greatest common divisor; the subfield generated by that power is then determined by the Frobenius period modulo its order.
\end{example}

\section{Polynomial Orders and the Contribution of Repeated Factors}
\begin{definition}\label{def:poly-order}
Let \(f\in k[X]\) be monic with \(f(0)\ne0\). The least positive integer \(N\) such that \(f\mid X^N-1\)
is called the \textbf{order} of \(f\), often also its period, and is denoted by \(\ord(f)\). Set \(\ord(1)=1\).
\end{definition}
For \(f\ne1\), this least integer exists: \(\overline X\) is invertible in the finite ring \(k[X]/(f)\),
since \(\gcd(X,f)=1\). It has finite order in the finite unit group, and \(\overline X^N=1\) is equivalent to \(f\mid X^N-1\).

\begin{proposition}\leavevmode\phantomsection\label{prop:irreducible-order}
If \(f\) is monic and irreducible, \(f\ne X\), and \(\alpha\) is a root, then
\[
\ord(f)=\ord(\alpha),\qquad
\deg f=\ord_{\ord(f)}(q).
\]
All roots of \(f\) have the same multiplicative order.
\end{proposition}
\begin{proof}
The divisibility criterion for the minimal polynomial gives \(f\mid X^N-1\Longleftrightarrow\alpha^N=1\).
The roots are related by Frobenius, and automorphisms preserve element orders.
Finally apply Proposition \ref{prop:degree-from-order}.
\end{proof}

\begin{theorem}[The order of a general polynomial]\label{thm:order-multiplicity}
Let \(f=\prod_{i=1}^t f_i^{a_i}\), where \(t\ge1\), \(a_i\ge1\), and the \(f_i\ne X\) are distinct monic irreducible polynomials.
Set
\[
e=\lcm_i\ord(f_i),\qquad A=\max_i a_i,
\qquad b=\min\{j\geq0:p^j\geq A\}.
\]
Then \(\ord(f)=ep^b\).
\end{theorem}
\begin{proof}
Write \(e_i=\ord(f_i)\). Since the nonzero roots of \(f_i\) lie in finite fields, \(p\nmid e_i\), and hence \(p\nmid e\).
Write any candidate exponent uniquely as \(N=p^cN_0\), with \(c\ge0\) and \(p\nmid N_0\).
The Frobenius identity gives
\[
X^N-1=(X^{N_0}-1)^{p^c}.
\]
The derivative of \(X^{N_0}-1\) is \(N_0X^{N_0-1}\). Its roots are nonzero, and \(N_0\ne0\) in \(k\),
so it has no repeated roots. Thus each distinct irreducible factor of \(X^N-1\) has multiplicity exactly \(p^c\).

By definition, \(f_i\mid X^{N_0}-1\) if and only if \(e_i\mid N_0\);
this is also equivalent to any root of \(f_i\) having \(N_0\)-th power one.
Together with unique factorization, this gives
\[
f\mid X^N-1
\quad\Longleftrightarrow\quad
e_i\mid N_0,\ a_i\le p^c\ \text{for all }i
\quad\Longleftrightarrow\quad
e\mid N_0,\ b\le c.
\]
Every admissible \(N\) is therefore divisible by \(ep^b\), and \(N=ep^b\) itself satisfies all the conditions.
Hence it is the least positive exponent.
\end{proof}

The part coprime to the characteristic is determined by the distinct roots, while the power of \(p\) is determined by root multiplicities.
Looking only at the set of distinct roots loses this multiplicity information in the polynomial order.

\begin{corollary}\leavevmode\phantomsection\label{cor:primitive-count}
The number of monic primitive polynomials of degree \(n\) in \(k[X]\) is
\[
\frac{\varphi(q^n-1)}{n}.
\]
More generally, for a positive integer \(e\) with \(\gcd(e,q)=1\), there are \(\varphi(e)/\ord_e(q)\) monic irreducible polynomials of order \(e\), each of degree \(\ord_e(q)\).
\end{corollary}
\begin{proof}
In \(\F_{q^{\ord_e(q)}}\), exactly \(\varphi(e)\) elements have order \(e\).
All their Frobenius orbits have length \(\ord_e(q)\), and each orbit corresponds to one monic irreducible polynomial.
For primitive elements, take \(e=q^n-1\); the orbit length is \(n\), giving the first formula.
\end{proof}

\section{Euler Totient Functions and Multiplicative Orders in Quotient Rings}
\label{sec:euler-polynomial}
Euler's totient function \(\varphi(m)=|(\Z/m\Z)^\times|\) satisfies
\[
\varphi(m)=m\prod_{\ell\mid m}(1-\ell^{-1}),\qquad
\sum_{d\mid m}\varphi(d)=m.
\]
The first formula excludes residue classes divisible by the prime factors; the second classifies the elements of a cyclic group of order \(m\) by their exact orders.
The Chinese remainder theorem for integers gives \(\varphi(ab)=\varphi(a)\varphi(b)\) when \(a,b\) are coprime.

\leavevmode\phantomsection\label{def:polynomial-euler}\ignorespaces For a nonzero monic polynomial \(F\in k[T]\), define the \textbf{polynomial Euler totient function}
\[
\varphi_q(F)=|(k[T]/(F))^\times|,
\qquad \varphi_q(1)=1.
\]
B\'ezout's identity shows that \(\bar h\) is invertible if and only if \(\gcd(h,F)=1\). If \(F=\prod_i P_i^{e_i}\), the Chinese remainder theorem gives
\begin{equation}\label{eq:polynomial-euler}
\varphi_q(F)=q^{\deg F}\prod_{P\mid F}(1-q^{-\deg P}),
\qquad \sum_{D\mid F}\varphi_q(D)=q^{\deg F}.
\end{equation}
The product runs over the distinct monic irreducible factors, and the sum over all monic divisors. The first formula also follows directly by inclusion--exclusion.
For the second, classify by \(D=F/\gcd(h,F)\): write \(h=(F/D)u\), with \(u\bmod D\) coprime to \(D\), giving exactly \(\varphi_q(D)\) possibilities.
This has the same structure as classification of a cyclic group by element orders.

The quantity \(\varphi_q(F)\) also counts vectors whose successive images under a linear operator form a basis.
Specifically, let \(m=\deg F\ge1\), and consider multiplication by \(\bar T\) as a linear operator on \(R=k[T]/(F)\).
For \(\bar h\in R\), the sequence
\[
\bar h,\quad \bar T\bar h,\quad\ldots,\quad\bar T^{m-1}\bar h
\]
is a \(k\)-basis if and only if \(\gcd(h,F)=1\): it is the image of the standard basis \(1,\bar T,\ldots,\bar T^{m-1}\) under multiplication by \(\bar h\),
and this map is invertible if and only if \(\bar h\) is a unit.
Chapter~\ref{chap:normal} represents the Frobenius action on an extension of degree \(n\) by the quotient ring \(k[T]/(T^n-1)\),
thereby reducing the count of normal elements to this count of units.

\begin{proposition}\leavevmode\phantomsection\label{prop:maximal-order-primitive}
Let \(f\in k[X]\) be monic, with \(\deg f=d\ge1\) and \(f(0)\ne0\). Then
\[
\ord(f)\le\varphi_q(f)\le q^d-1.
\]
Moreover, \(\ord(f)=q^d-1\) if and only if \(f\) is primitive; irreducibility need not be assumed in advance.
\end{proposition}
\begin{proof}
The class \(\bar X\) is a unit of \(k[X]/(f)\), of order \(\ord(f)\). Lagrange's theorem bounds its order by the number of units.
If its order is \(q^d-1\), every nonzero element of the quotient ring is a unit. Thus the quotient is a field, \(f\) is irreducible, and \(\bar X\) is primitive.
The converse follows from the definition of a primitive polynomial.
\end{proof}

\section{Projective Orders Modulo Base-Field Scalars}
\begin{definition}\label{def:projective-order}
Let \(f\in k[X]\) be monic, with \(\deg f=d\ge1\) and \(f(0)\ne0\). Its \textbf{projective order} is the least positive integer \(r\) such that \(X^r\equiv a\pmod f\) for some \(a\in k^\times\).
\end{definition}
With \(R=k[X]/(f)\), this is the order of the coset \(\bar Xk^\times\) in \(R^\times/k^\times\).
Units differing by a nonzero base-field scalar are identified as one point. The projective order records the first return of the power sequence \(1,\bar X,\bar X^2,\ldots\) under this identification.
If \(\bar X^r=a\), then
\[
\ord(f)=r\ord(a),\qquad r\le\frac{q^d-1}{q-1}.
\]
For the first formula, \(\bar X^j=1\) first requires \(r\mid j\), and then \(a^{j/r}=1\).
The second follows by grouping units modulo nonzero base-field scalars, with \(q-1\) elements in each group.

\begin{proposition}\label{prop:projective-primitive}
Set \(S=(q^d-1)/(q-1)\). The polynomial \(f\) is primitive if and only if its projective order is \(S\)
and \((-1)^df(0)\) generates \(k^\times\).
\end{proposition}
\begin{proof}
If the projective order is \(S\), the quotient ring already has \(S(q-1)=q^d-1\) units, and hence is a field.
Set \(\alpha=\bar X\). Its distinct conjugates are \(\alpha,\alpha^q,\ldots,\alpha^{q^{d-1}}\), so \(\alpha^S=\prod_{i=0}^{d-1}\alpha^{q^i}=(-1)^df(0)\). This product of conjugates is the norm introduced in Chapter~\ref{chap:trace}.
The preceding formula then gives the necessary and sufficient condition for \(\ord(\alpha)=S(q-1)\).
The converse follows from the same formula, or from the formula for orders in quotients of cyclic groups.
\end{proof}

\subsection{Projective Order and Minimal Binomial Multiples}
\label{v40:sec-minimal-binomial}

Projective order can also be computed directly from roots and their multiplicities. Ordinary order requires powers of all roots to equal one; projective order allows them to take a common nonzero base-field value. For the structure of minimal binomial multiples, see \cite{ZWMinimalBinomial}.
\leavevmode\phantomsection\label{v40:def-minimal-binomial}
For a nonconstant monic polynomial \(f\in k[X]\) with \(f(0)\ne0\), a monic binomial \(X^L-c\), \(c\in k^\times\), of least degree among those divisible by \(f\) is called its \textbf{minimal binomial multiple}. Its degree is the projective order of \(f\).

\begin{proposition}[Root ratios and minimal binomial multiples]\label{v40:prop-minimal-binomial}
In a finite splitting extension, let the distinct roots of \(f\) be \(\alpha_1,\ldots,\alpha_s\), and let their maximum multiplicity be \(M\). Set
\[
L_0=\lcm_{1\le i\le s}\ord(\alpha_i/\alpha_1),
\qquad b=\min\{j\ge0:p^j\ge M\},
\qquad \beta=\alpha_1^{L_0}.
\]
Then \(\beta\in k^\times\), and the minimal binomial multiple is unique, namely
\[
B_f(X)=X^{L_0p^b}-\beta^{p^b}.
\]
All monic binomial multiples are precisely
\[
X^{uL_0p^b}-\beta^{up^b},\qquad u\ge1.
\]
\end{proposition}
\begin{proof}
The roots are nonzero, and their ratios have orders coprime to \(p\), so \(p\nmid L_0\). By the definition of \(L_0\), all \(\alpha_i^{L_0}\) equal \(\beta\). Frobenius permutes the roots of \(f\), so \(\beta^q=(\alpha_1^q)^{L_0}=\beta\), and hence \(\beta\in k^\times\). The derivative of \(X^{L_0}-\beta\) is nonzero at each \(\alpha_i\), so all these roots have multiplicity one in it. Taking the \(p^b\)-th power gives
\[
(X^{L_0}-\beta)^{p^b}=X^{L_0p^b}-\beta^{p^b},
\]
where each root has multiplicity \(p^b\ge M\). Thus \(f\mid B_f\).

Conversely, if \(f\mid X^L-c\), then \(\alpha_i^L=\alpha_1^L\), so every \(\ord(\alpha_i/\alpha_1)\) divides \(L\), and hence \(L_0\mid L\). Write \(L=p^vL'\), with \(p\nmid L'\). In the finite field, \(c\) has a unique \(p^v\)-th root \(c'\), and
\[
X^L-c=(X^{L'}-c')^{p^v}.
\]
All nonzero roots of the inner polynomial are simple, so every root of the original binomial has multiplicity exactly \(p^v\). This forces \(p^v\ge M\), or \(v\ge b\). Together with \(p\nmid L_0\), it gives \(L_0p^b\mid L\). The constant must be \(c=\alpha_1^L\), so every candidate has the stated form. Conversely, \(X^{uL_0p^b}-\beta^{up^b}\) is divisible by \(B_f\), and hence by \(f\). This proves minimality, uniqueness, and the description of all binomial multiples.
\end{proof}

\begin{corollary}[Projective order of an irreducible polynomial]\label{v40:cor-irreducible-projective}
If \(f\ne X\) is monic and irreducible, and a root has multiplicative order \(N\), then
\[
\text{projective order of \(f\)}=\frac{N}{\gcd(N,q-1)}.
\]
\end{corollary}
\begin{proof}
The roots are \(\alpha^{q^i}\). Every root ratio \(\alpha^{q^i-1}\) is a power of \(\alpha^{q-1}\), and \(\alpha^{q-1}\) is itself a root ratio. These ratios therefore generate \(\langle\alpha^{q-1}\rangle\), of order \(N/\gcd(N,q-1)\). Since an irreducible polynomial has no repeated roots, \(b=0\) in the preceding proposition.
\end{proof}

The group \(H=\langle\alpha_i/\alpha_1:1\le i\le s\rangle\) generated by root ratios gives the smallest multiplicative coset \(\alpha_1H\) containing the root set: the subgroup underlying any such coset must contain all these ratios. Since the multiplicative group of a finite field is cyclic, \(|H|=L_0\). This parallels the affine envelope spanned by root differences in Proposition~\ref{prop:affine-q-multiple}.
The minimal binomial multiple is determined by root ratios and maximum multiplicity. Writing the original polynomial as \(G(X^r)\), by contrast, requires the corresponding multiplicative symmetry of its whole root set. These are different conditions.

\begin{example}\label{v40:ex-minimal-binomial}
Take \(f=X^3+X+1\) in \(\F_2[X]\). Its roots have multiplicative order seven, so its projective order is seven and its minimal binomial multiple is \(X^7+1\). However, the positive exponents in the support of \(f\) are \(1,3\), so the largest \(r\) allowing an expression of the form \(G(X^r)\) is one. For \(f^3\), the maximum root multiplicity is three, and the least required power of two is four. Hence its minimal binomial multiple is \(X^{28}+1\), whereas \(X^{21}+1\) is square-free and cannot be divisible by \(f^3\).
\end{example}

For an irreducible quadratic, projective order can also be computed from the first zero term of positive index in a coefficient recurrence; see Section~\ref{v45:sec-quadratic-criterion}.

A binomial multiple allows only two nonzero terms. With more terms allowed but the divisor fixed, Subsection~\ref{v97:sec-sparse-relations} records the first possible degree by \(L_s(f)\); the multiplicative period alone no longer determines this threshold.

\section{Reciprocal Polynomials and Inverse-Root Orbits}\label{sec:reciprocal-polynomials}
\leavevmode\phantomsection\label{def:reciprocal-polynomial}\ignorespaces For a monic polynomial \(f\) of degree \(d\) with \(f(0)\ne0\), define its monic \textbf{reciprocal polynomial} by
\[
f^\vee(X)=\frac{X^df(X^{-1})}{f(0)}.
\]
Its roots are the inverses of the roots of \(f\), with the same respective multiplicities, and \((f^\vee)^\vee=f\).
\phantomsection\label{def:self-reciprocal-polynomial}
A polynomial satisfying \(f=f^\vee\) is called \textbf{self-reciprocal}.
This operation preserves irreducibility, since taking reciprocals of a nontrivial factorization again gives a nontrivial factorization.

\begin{proposition}\leavevmode\phantomsection\label{prop:self-reciprocal}
Let \(f\ne X\) be monic and irreducible. If \(f=f^\vee\), then either \(f=X-1\) or \(X+1\),
or \(\deg f=2s\) and every root \(\alpha\) satisfies \(\alpha^{q^s+1}=1\).
Conversely, if a root of a monic irreducible polynomial of degree \(2s\) satisfies this equality, then \(f=f^\vee\).
\end{proposition}
\begin{proof}
Self-reciprocity is equivalent to \(\alpha^{-1}\) lying in the same Frobenius orbit as \(\alpha\), for some root \(\alpha\).
If \(\alpha^{-1}=\alpha\), then \(\alpha=\pm1\), giving a linear polynomial.
Otherwise write \(\alpha^{-1}=\alpha^{q^j}\), where \(0<j<d=\deg f\).
Taking the \(q^j\)-th power again gives \(\alpha^{q^{2j}}=\alpha\), so \(d\mid2j\).
Since \(0<2j<2d\), we have \(2j=d\). Thus \(d=2s\), \(j=s\), and \(\alpha^{q^s+1}=1\).
Conversely, the stated equality places the inverse in the original orbit, so the two monic minimal polynomials coincide.
\end{proof}

In exponent coordinates, inversion is \(i\mapsto-i\), and Frobenius is \(i\mapsto qi\).
They commute, so inversion sends an entire Frobenius orbit to another orbit.
If \(\alpha\) has order \(e\), its conjugates and the conjugates of its inverse are, respectively,
\[
\{\alpha^{q^j}:j\ge0\},\qquad \{\alpha^{-q^j}:j\ge0\}.
\]
These roots determine \(f\) and \(f^\vee\), respectively. Thus reciprocity is an involution on irreducible polynomials,
whose fixed points are precisely the self-reciprocal irreducibles; all remaining irreducibles occur in pairs.

\begin{corollary}[An order criterion for self-reciprocity]\label{v48:cor-reciprocal-order}
Let \(f\ne X\) be a monic irreducible polynomial whose roots have multiplicative order \(e>2\), and set \(d=\ord_e(q)\). Then
\[
f=f^\vee\quad\Longleftrightarrow\quad
-1\in\langle q\rangle\subseteq(\Z/e\Z)^\times
\quad\Longleftrightarrow\quad
d\text{ is even and }q^{d/2}\equiv-1\pmod e.
\]
\end{corollary}
\begin{proof}
The inverse is a conjugate if and only if \(q^j\equiv-1\pmod e\) for some \(j\), giving the first equivalence.
Since \(e>2\), the element \(-1\) has order two. The cyclic group \(\langle q\rangle\) of order \(d\) has an element of order two if and only if \(d\) is even,
in which case its unique element of order two is \(q^{d/2}\). This proves the second equivalence.
\end{proof}

The criterion requires inversion exactly halfway around the orbit; even degree alone is insufficient.
For example, \(\ord_8(3)=2\), but \(3\not\equiv-1\pmod8\).
Thus irreducible quadratics over \(\F_3\) whose roots have multiplicative order eight are not self-reciprocal.
The polynomial \(X^2+X+2\) from Example~\ref{v48:ex-coordinate-model-change} and its reciprocal \(X^2+2X+2\) form such a pair.

\leavevmode\phantomsection\label{v48:def-norm-one-group}\ignorespaces
In \(\F_{q^{2s}}^\times\), set
\[
U_s=\{x:x^{q^s+1}=1\}.
\]
This is a cyclic subgroup of order \(q^s+1\), called the \textbf{norm-one subgroup} relative to \(\F_{q^s}\):
the norm of Chapter~\ref{chap:trace} is here \(x\mapsto x\,x^{q^s}=x^{q^s+1}\).
Its order follows from Corollary~\ref{cor:power-roots}, since \(q^s+1\mid q^{2s}-1\).
Proposition~\ref{prop:self-reciprocal} shows that self-reciprocal irreducible polynomials of degree \(2s\) correspond precisely to the conjugacy orbits of degree-\(2s\) elements in \(U_s\).
The degree condition is essential: \(U_s\) also contains \(1\) and may contain elements of proper subfields.

\begin{example}[Inversion at half-period for roots of order ten]\leavevmode\phantomsection\label{v48:ex-reciprocal-ten}
Choose \(\gamma\) of order ten in \(\F_{81}\). Then
\[
\ord_{10}(3)=4,\qquad 3^2\equiv-1\pmod{10}.
\]
Thus the conjugates of \(\gamma\) occur in the order \(\gamma,\gamma^3,\gamma^9,\gamma^7\),
and \(\gamma^{-1}=\gamma^9\) is the conjugate halfway around the Frobenius orbit.
Since \(\gamma^5=-1\) and \(\gamma\ne-1\), the identity
\[
X^5+1=(X+1)(X^4-X^3+X^2-X+1)
\]
shows that \(\gamma\) is a root of \(h=X^4-X^3+X^2-X+1\).
Its degree is four, so \(h\) is its minimal polynomial, both irreducible and self-reciprocal.
On the other hand, \(\gamma^{10}=1\) places it in \(U_2\), and its conjugate over \(\F_9\) is \(\gamma^9=\gamma^{-1}\).
Self-reciprocity thus appears both as a relation between opposite points of the base-field Frobenius orbit and as the statement that the product of an element with its conjugate over the intermediate field is one.
\end{example}

\subsection{Inverse pairs and counts of full-degree elements}
\label{v145:sec-inverse-count}
Pairing elements of the norm-one group with their inverses gives quadratic polynomials; grouping them into Frobenius orbits gives self-reciprocal irreducibles. The two counts differ according to whether the degree over the original base field has been specified.

Let $S_q(2n)$ be the number of monic self-reciprocal irreducible polynomials of degree $2n$ over $\F_q$, and put
\begin{equation}\label{v145:eq-target-count}
 C_n(q)=\#\left\{b\in\F_{q^n}:
 \begin{array}{l}Z^2+bZ+1\text{ is irreducible over }\F_{q^n},\\
 {}[\F_q(b):\F_q]=n
 \end{array}\right\}.
\end{equation}

\begin{proposition}[Three counts of full inverse orbits]\label{v145:prop-inverse-count}
For $n\ge1$, set $\epsilon_q=0$ if $q$ is even and $\epsilon_q=1$ if $q$ is odd. Then
\begin{equation}\label{v145:eq-inverse-count}
\begin{aligned}
 C_n(q)
 &=\frac12\sum_{\substack{d\mid q^n+1\\\ord_d(q)=2n}}\varphi(d)\\
 &=\frac12\sum_{\substack{e\mid n\\e\text{ odd}}}
       \mu(e)\bigl(q^{n/e}-\epsilon_q\bigr)
 =nS_q(2n),
\end{aligned}
\end{equation}
with the convention $\ord_1(q)=1$.
\end{proposition}
\begin{proof}
On $U_n$, the map $\alpha\mapsto b=-(\alpha+\alpha^{-1})$ pairs elements outside $\F_{q^n}$. Since
\[
 U_n\cap\F_{q^n}^{\times}=\{\alpha:\alpha^2=1\},
\]
the resulting quadratics are precisely the irreducible ones, with roots $\alpha,\alpha^{-1}$. If $b$ has degree $n$ over $\F_q$, then $\F_q(b)=\F_{q^n}$ and the quadratic extension gives $[\F_q(\alpha):\F_q]=2n$. Conversely, if $\alpha$ has degree $2n$, then $b\in\F_{q^n}$ and $[\F_q(\alpha):\F_q(b)]\le2$, forcing $b$ to have degree $n$. There are $\varphi(d)$ elements of order $d$, each of degree $\ord_d(q)$. Dividing by two proves the first formula.

There is one element of $U_j$ satisfying $\alpha^2=1$ when $q$ is even, and two when $q$ is odd. Thus the number of irreducible quadratics $Z^2+bZ+1$ over $\F_{q^j}$ is
\[
 A_j=\frac{q^j-\epsilon_q}{2}.
\]
Classify $b$ by its exact degree $d\mid j$. An irreducible quadratic over $\F_{q^d}$ stays irreducible over $\F_{q^j}$ exactly when $j/d$ is odd: its root generates $\F_{q^{2d}}$ and belongs to $\F_{q^j}$ exactly when $2d\mid j$. Therefore
\[
 A_j=\sum_{\substack{d\mid j\\j/d\text{ odd}}}C_d(q).
\]
M\"obius inversion on the odd divisors gives the second formula.

By Proposition~\ref{prop:self-reciprocal}, the degree-$2n$ elements of $U_n$ are exactly the roots of the self-reciprocal irreducibles being counted. Each Frobenius orbit of size $2n$ contains $n$ inverse pairs, proving $C_n(q)=nS_q(2n)$.
\end{proof}

For example, $C_2(3)=4$, so $S_3(4)=2$. Example~\ref{v48:ex-reciprocal-ten} gives $X^4-X^3+X^2-X+1$ from roots of order ten; roots of order five give $X^4+X^3+X^2+X+1$. These are the complete list.
For $q=2,n=3$, one has $C_3(2)=(8-2)/2=3$ and $S_2(6)=1$. The unique polynomial is $X^6+X^3+1$, whose roots have order nine.

These target coefficients also enter a return-map count: in Theorem~\ref{v145:thm-sparse-count}, $C_n(q)$ counts the allowable return characteristic polynomials, and a second factor counts the preimages of each target.

\par\Needspace{11\baselineskip}
\section*{Exercises}
\phantomsection
\addcontentsline{toc}{section}{Exercises}
\begin{exercise}\label{ex:order12}
Find the order of \((X^2+X+1)^3(X+1)^2\) in \(\F_2[X]\), and explain its odd part and power-of-two part separately.
\end{exercise}
\begin{exercise}\label{ex:order-degree}
Let \(\alpha\) be an element of order thirteen in \(\F_{2^{12}}\). Determine its degrees over \(\F_2\), \(\F_4\), and \(\F_8\).
\end{exercise}
\begin{exercise}\label{ex:reciprocal3}
Find the reciprocal of \(X^3+X+1\) over \(\F_2\), and explain why it is not self-reciprocal.
Prove that \(X^4+X^3+X^2+X+1\) is self-reciprocal.
\end{exercise}
\begin{exercise}\label{ex:primitive4count}
Count the monic primitive quartics over \(\F_2\), and compare this with the number of all monic irreducible quartics.
\end{exercise}

\begin{exercise}\leavevmode\phantomsection\label{ex:supp-projective}
Let \(f=X^2+1\in\F_3[X]\). Determine its ordinary and projective orders, and decide whether it is primitive. Then explain why knowing only that the quotient-ring unit \(\bar X\) has order \(q^d-1\) suffices to prove that \(f\) is irreducible.
\end{exercise}

\begin{exercise}[The minimal polynomial of a multiplication matrix]\label{ex:wu-matrix-polynomial}
Let \(E/k\) be an extension of finite fields, let \(a\in E\), and set \(M_a(x)=ax\).
Prove that the minimal polynomial of the linear operator \(M_a\) is exactly the minimal polynomial \(m_{a,k}\) of the field element.
Compare this with the monic annihilating polynomial of least degree for the vector \(a\) under the Frobenius operator.
\end{exercise}
\par

\chapter{Roots of Unity, Base Extension, and Tensor Product Decompositions}
\label{chap:cyclotomic}
On roots of unity, Frobenius acts by multiplication on exponents; extending the base field changes the step size. Cyclotomic factorization, intersections, and composita are therefore governed by the same period calculation, while tensor products record all components of the factorization.
For Galois rings, residue-field orbits and Hensel lifting further determine the factors associated with roots of unity and the direct product decomposition of quotient rings.

Once root periods have determined factor degrees, four questions remain. The first concerns a single orbit: what additive symmetries does its exponent set have, and can they force the root polynomial to take the form \(G(X^r)\)? Section~\ref{sec:cyclotomic-med} uses the additive stabilizer to determine maximal blocks and turns their period into the extension degree needed for binomial factorization.

The second question involves different moduli: do local orbits at each prime power determine the global orbit? Section~\ref{v44:sec-cyclotomic-local-global} expresses the missing information as relative phase. The third follows a divisibility chain of moduli: how many orbits occur in one lift, and when does splitting stop? Section~\ref{v44:sec-cyclotomic-systems} records branches by return maps and determines stable factor towers from local multiplicative orders.

Finally, when coefficients are lifted from a finite field to a Galois ring, to what precision does the composition structure determined by the root set persist? Section~\ref{v47:sec-galois-rings} answers this through Hensel uniqueness, distinguishing complete root blocks, additional zero coefficients, and repeated-root blocks. The same family of roots with \(q=3\) and moduli \(5\cdot2^a\) will run through these questions.

\section{The Splitting Field of the Roots of Unity}
\leavevmode\phantomsection\label{def:roots-of-unity}\ignorespaces Let \(N\) be a positive integer with \(\gcd(N,q)=1\). In an extension in which \(X^N-1\) splits, put
\[
\mu_N=\{\zeta:\zeta^N=1\}.
\]
The elements of \(\mu_N\) are called \(N\)-th \textbf{roots of unity}. This set is a multiplicative group. The derivative \(NX^{N-1}\) shows that the roots are simple, so \(|\mu_N|=N\).
By Theorem~\ref{thm:cyclic}, \(\mu_N\) is cyclic; an element of order \(N\) is called a \textbf{primitive \(N\)-th root of unity}.

\begin{proposition}\label{prop:unity-splitting}
The splitting field of \(X^N-1\) over \(\F_q\) is \(\F_{q^s}\), where \(s=\ord_N(q)\).
For \(N=1\), we set \(s=1\).
\end{proposition}
\begin{proof}
Choose a generator \(\zeta\) of \(\mu_N\). Every root is a power of \(\zeta\), so the splitting field is \(\F_q(\zeta)\).
Proposition~\ref{prop:degree-from-order} gives its degree \(s\).
\end{proof}

If \(N=p^aM\), where \(p\nmid M\), then
\[
X^N-1=(X^M-1)^{p^a}.
\]
The distinct roots are still \(\mu_M\), but each has multiplicity \(p^a\).
Thus there are no nontrivial \(p\)-th roots of unity in characteristic \(p\).

\section{Cyclotomic Polynomials over the Integers: A Root-Set Definition}
To clarify the relation between cyclotomic polynomials and reduction modulo \(p\), we first define them in characteristic zero.
The construction in Lemma~\ref{lem:finite-splitting} applies over \(\Q\). In an extension in which \(X^N-1\) splits, the group \(\mu_N\) is again cyclic of order \(N\).
\leavevmode\phantomsection\label{def:cyclotomic-polynomial}\ignorespaces Define
\[
\Phi_N(X)=\prod_{\substack{\zeta\in\mu_N\\\ord(\zeta)=N}}(X-\zeta).
\]
This is the \(N\)-th \textbf{cyclotomic polynomial}.

\begin{proposition}\label{prop:cyclotomic-integral}
We have \(\Phi_N\in\Z[X]\), \(\deg\Phi_N=\varphi(N)\), and
\[
X^N-1=\prod_{d\mid N}\Phi_d(X).
\]
These identities uniquely determine all the \(\Phi_N\).
\end{proposition}
\begin{proof}
The order of an \(N\)-th root of unity is a positive divisor of \(N\). Grouping the roots by order gives the product formula.
A cyclic group has exactly \(\varphi(N)\) elements of order \(N\), giving the degree.

We prove that \(\Phi_N\) has integer coefficients by induction on \(N\). We have \(\Phi_1=X-1\). Suppose \(\Phi_d\in\Z[X]\) for every proper divisor \(d\mid N\).
Their product \(D\) is a monic polynomial with integer coefficients. On dividing \(X^N-1\) by \(D\), the quotient and remainder have integer coefficients because the divisor is monic.
In the extension containing the roots, the quotient is \(\Phi_N\) and the remainder is zero. Uniqueness of division gives \(\Phi_N\in\Z[X]\).
The same recursion proves uniqueness.
\end{proof}

Over a field of characteristic \(p\), \(\Phi_N\) denotes the polynomial obtained by reducing its integer coefficients modulo \(p\).
If \(p\nmid N\), its roots are precisely the primitive \(N\)-th roots of unity: reduce the product formula above,
use the fact that \(X^N-1\) has no repeated roots, and induct over the divisors of \(N\).

\begin{proposition}[When the characteristic divides the index]\leavevmode\phantomsection\label{prop:cyclotomic-p}
If \(p\nmid M\) and \(a\geq1\), then in characteristic \(p\),
\[
\Phi_{Mp^a}(X)=\Phi_M(X)^{p^{a-1}(p-1)}.
\]
\end{proposition}
\begin{proof}
In characteristic zero, the roots of \(\Phi_M(X^{p^a})\) have multiplicative orders exactly \(M,Mp,\ldots,Mp^a\). Hence
\[
\Phi_M(X^{p^a})=\prod_{j=0}^{a}\Phi_{Mp^j}(X).
\]
Comparing this with the formula for \(a-1\) gives \(\Phi_M(X^{p^a})=\Phi_M(X^{p^{a-1}})\Phi_{Mp^a}(X)\).
After reduction to characteristic \(p\), the integer coefficients lie in the prime field, so the polynomials \(\Phi_M(X^{p^a})\) and \(\Phi_M(X^{p^{a-1}})\) become \(\Phi_M(X)^{p^a}\) and \(\Phi_M(X)^{p^{a-1}}\), respectively. Cancellation in the polynomial ring, which is an integral domain, proves the claim.
\end{proof}

\section{Cyclotomic Cosets and Exponent Orbits}
Fix \(\gcd(N,q)=1\) and a primitive \(N\)-th root of unity \(\zeta\).
The map \(\sigma(a)=a^q\) satisfies
\[
\sigma(\zeta^i)=\zeta^{qi}.
\]
Thus conjugacy becomes a question about exponents modulo \(N\).

\begin{definition}\label{def:cyclotomic-coset}
A \(q\)-\textbf{cyclotomic coset} modulo \(N\) is an orbit of the permutation \(i\mapsto qi\):
\[
C_i=\{i,qi,q^2i,\ldots\}\subseteq\Z/N\Z.
\]
The word ``coset'' is conventional; in general, this is an orbit of a group action rather than a coset of the additive group modulo \(N\).
\end{definition}

\begin{theorem}\label{thm:coset-factorization}
For each cyclotomic coset \(C\), put
\[
M_C(X)=\prod_{i\in C}(X-\zeta^i).
\]
Then \(M_C\in\F_q[X]\) is monic and irreducible, and
\[
X^N-1=\prod_C M_C(X),\qquad
\deg M_{C_i}=|C_i|=\ord_{N/\gcd(N,i)}(q).
\]
For \(i=0\), we have \(N/\gcd(N,0)=1\), with the convention \(\ord_1(q)=1\).
\end{theorem}
\begin{proof}
The set \(\{\zeta^i:i\in C\}\) consists of all Frobenius conjugates of one element, so Theorem~\ref{thm:orbit} gives irreducibility and degree.
The orbits partition \(\mu_N\), giving the product formula.
The multiplicative order of \(\zeta^i\) is \(N/\gcd(N,i)\); Proposition~\ref{prop:degree-from-order} now gives the last identity.
\end{proof}

\begin{corollary}\label{cor:cyclotomic-equal-degree}
If \(\gcd(N,q)=1\), every irreducible factor of \(\Phi_N\) over \(\F_q\) has degree \(\ord_N(q)\),
and there are \(\varphi(N)/\ord_N(q)\) such factors.
In particular, \(\Phi_N\) is irreducible if and only if \(\ord_N(q)=\varphi(N)\).
\end{corollary}
\begin{proof}
The roots of \(\Phi_N\) correspond to exponents coprime to \(N\). All have multiplicative order \(N\), so their orbits have equal length.
\end{proof}

\begin{proposition}[A cyclotomic criterion for a common Frobenius period]\label{v124:prop-uniform-cyclotomic}
Let \(N>1\), \(\gcd(N,q)=1\), and \(m\ge1\). The following conditions are equivalent:
\begin{enumerate}
\item every monic irreducible factor of \(S_N(X)=(X^N-1)/(X-1)\) over \(\F_q\) has degree \(m\);
\item \(N\mid q^m-1\), and \(\ord_\ell(q)=m\) for every prime \(\ell\mid N\);
\item \(N\mid\Phi_m(q)\) and \(\gcd(N,m)=1\).
\end{enumerate}
Here \(\Phi_m(q)\) is the value of the integer cyclotomic polynomial at the integer \(q\). Under these conditions, every nonzero cyclotomic coset has length \(m\), there are \((N-1)/m\) factors, and \(m\mid\ell-1\) for every prime divisor \(\ell\mid N\).
\end{proposition}
\begin{proof}
The roots of \(S_N\) are \(\mu_N\setminus\{1\}\). Every divisor \(d>1\) of \(N\) occurs as a root order, with Frobenius period \(\ord_d(q)\). The first condition therefore says that all these orders equal \(m\), and implies the second. Conversely, under the second condition, choose a prime \(\ell\mid d\) for any divisor \(d>1\) of \(N\). Then
\[
m=\ord_\ell(q)\mid\ord_d(q)\mid m.
\]
Thus every nonidentity root has degree \(m\). The orbit decomposition also gives the number of factors.

Under the second condition, \(m=\ord_\ell(q)\mid\ell-1\), so \(\gcd(N,m)=1\). For every proper divisor \(a\) of \(m\), one has \(\ell\nmid q^a-1\). Hence in
\[
q^m-1=\Phi_m(q)\prod_{\substack{a\mid m\\a<m}}\Phi_a(q)
\]
the entire power of \(\ell\) comes from \(\Phi_m(q)\). Since \(N\mid q^m-1\), it follows that \(N\mid\Phi_m(q)\).

Under the third condition, \(N\mid q^m-1\). Fix a prime \(\ell\mid N\). Since \(\ell\nmid m\), the reduction of \(X^m-1\) modulo \(\ell\) is square-free, so distinct factors in its cyclotomic product have no common root. The residue of \(q\) is a root of \(\Phi_m\). If its order were a proper divisor \(a\) of \(m\), it would also be a root of \(X^a-1\), contradicting this coprimality. Thus \(\ord_\ell(q)=m\), proving the second condition.
\end{proof}

For fixed \(q,m\), the admissible moduli are exactly the divisors greater than one of \(\Phi_m(q)\) that are coprime to \(m\). For \(m=1\), this is simply \(N\mid q-1\). The criterion replaces separate orbit checks by one integer cyclotomic evaluation.

\begin{example}[Repeated prime divisors and a common degree]\label{v124:ex-uniform-121}
Since \(\Phi_5(3)=1+3+3^2+3^3+3^4=121=11^2\) and \(\gcd(121,5)=1\), the polynomial \((X^{121}-1)/(X-1)\) has exactly twenty-four irreducible factors of degree five over \(\F_3\). Roots of orders eleven and 121 contribute \(\varphi(11)/5=2\) and \(\varphi(121)/5=22\) factors, respectively.

The coprimality condition cannot be omitted. Although \(3\mid\Phi_6(2)=3\), the only nonzero binary cyclotomic coset modulo three is \(\{1,2\}\), of length two rather than six.
\end{example}

\begin{example}[Seventh roots of unity]\label{v50:ex-seven-orbits}
The binary cyclotomic cosets modulo seven are
\[
\{0\},\quad\{1,2,4\},\quad\{3,6,5\}.
\]
Thus \(X^7-1\) factors over \(\F_2\) into one linear and two cubic factors. Explicitly,
\[
X^7-1=(X+1)(X^3+X+1)(X^3+X^2+1).
\]
Choosing \(\zeta\) as a root of the first cubic identifies the two nonzero orbits with the two cubic factors in the displayed order.
Changing the primitive root of unity may interchange the labels of the factors, but leaves the factorization unchanged.
\end{example}

\begin{example}[Fifteenth roots of unity]
The binary cyclotomic cosets modulo fifteen are
\[
\{0\},\ \{1,2,4,8\},\ \{3,6,12,9\},\ \{5,10\},\ \{7,14,13,11\}.
\]
Take \(\zeta\) to be a root of \(X^4+X+1\). The corresponding irreducible factors, in order, are
\[
X+1,\quad X^4+X+1,\quad X^4+X^3+X^2+X+1,
\]
\[
X^2+X+1,\quad X^4+X^3+1.
\]
The two orbits of primitive fifteenth roots give two reciprocal quartic factors.
The roots of order five give the self-reciprocal factor \(\Phi_5\).
\end{example}

\begin{proposition}[Frobenius on a multiplicative quotient and descent of root sets]\label{v85:prop-coset-frobenius}
Let \(k=\F_q\), \(K=\F_{q^m}\), and let \(H\le K^\times\) have order \(n\) and index \(d\). Choose a generator \(\gamma\) of \(K^\times\).
Then \(H=\langle\gamma^d\rangle\), and Frobenius acts on the multiplicative quotient group by
\[
 K^\times/H\simeq\Z/d\Z,\qquad
 \gamma^jH\longmapsto\gamma^{qj}H\quad\longleftrightarrow\quad j\longmapsto qj.
\]
The orbit length of \(\gamma^jH\) is \(\ord_{d/\gcd(d,j)}(q)\), with \(\ord_1(q)=1\). Its root polynomial is
\[
 p_{\gamma^jH}(T)=T^n-\gamma^{jn},
\]
and belongs to \(k[T]\) if and only if \((q-1)j\equiv0\pmod d\).
\end{proposition}
\begin{proof}
The power map \(x\mapsto x^q\) is invertible on both \(K^\times\) and \(H\), so it induces the stated quotient action.
Its \(s\)th iterate fixes \(\gamma^jH\) exactly when \(d\mid(q^s-1)j\), giving the orbit length.
The elements of \(H\) are all the roots of \(T^n-1\); scaling this root set gives \(p_{\gamma^jH}\).
Its coefficients belong to \(k\) exactly when \((\gamma^{jn})^{q-1}=1\), equivalently \(d\mid(q-1)j\).
This is also Theorem~\ref{thm:root-descent} applied to a multiplicative coset.
\end{proof}
A multiplicative coset here is a set of roots, not necessarily a single Frobenius orbit. The quotient orbit records how the entire set changes under Frobenius.
When \(K=k\), the relative \(q\)-Frobenius is the identity; to study cosets relative to the prime field, use \(x\mapsto x^p\) instead.

\subsection{Coset Leaders and Order-Preserving Filtering}
\label{v50:sec-coset-leaders}

Any element of an exponent orbit can serve as its representative. Choosing the least nonnegative integer gives a uniform index for the irreducible factors;
the indices for different moduli can also be related by order-preserving scaling~\cite{WZGeneralizedCosets2026}.

\begin{definition}[Cyclotomic coset leader]\label{v50:def-coset-leader}
For an integer \(a\) and a positive integer \(n\), let \([a]_n\) denote the least nonnegative residue of \(a\) modulo \(n\).
When \(\gcd(n,q)=1\), write
\[
C_a^{(n)}=\{[aq^j]_n:j\ge0\},\qquad
\lambda_n(a)=\min C_a^{(n)}.
\]
We call \(\lambda_n(a)\) the \textbf{coset leader} of \(C_a^{(n)}\), and put
\[
\Gamma_n(q)=\{\lambda_n(a):0\le a<n\}.
\]
Thus \(\Gamma_n(q)\) contains one representative of each cyclotomic coset, including the leader \(0\) of the zero coset.
\end{definition}

\begin{theorem}[Order-preserving scaling and coset-leader filtering]\label{v50:thm-leader-scaling}
Let \(N=cn\), where \(c,n\ge1\) and \(\gcd(N,q)=1\). The map
\[
\iota_c:\Z/n\Z\longrightarrow\Z/N\Z,\qquad a\longmapsto ca
\]
is an isomorphism onto the additive subgroup \(c\Z/N\Z\), and commutes with multiplication by \(q\). For every integer \(a\), it gives
\begin{equation}\label{v50:eq-leader-scaling}
C_{ca}^{(N)}=cC_a^{(n)},\qquad
|C_{ca}^{(N)}|=|C_a^{(n)}|,\qquad
\lambda_N(ca)=c\lambda_n(a).
\end{equation}
Corresponding elements also have the same additive order. The two sets of coset leaders satisfy
\begin{equation}\label{v50:eq-leader-filter}
\boxed{\quad
\Gamma_n(q)=\{A/c:A\in\Gamma_N(q),\ c\mid A\}.
\quad}
\end{equation}
This correspondence preserves strict inequalities between coset leaders. Thus the \(j\)-th largest leader for the smaller modulus is obtained from the leaders for the larger modulus
by retaining those divisible by \(c\), selecting the \(j\)-th entry in descending order, and dividing by \(c\), provided that entry exists.
\end{theorem}
\begin{proof}
The congruence \(ca\equiv cb\pmod{cn}\) is equivalent to \(a\equiv b\pmod n\), so \(\iota_c\) is well defined and injective,
with image \(c\Z/N\Z\). It preserves addition and satisfies \(\iota_c(qa)=q\iota_c(a)\).
Writing \(aq^j=un+r\), with \(0\le r<n\), gives
\[
[caq^j]_N=cr=c[aq^j]_n.
\]
Thus \(\iota_c\) maps each orbit for the smaller modulus bijectively onto its corresponding orbit for the larger modulus, preserving orbit length.
Since \(0\le cr<N\) and \(c>0\), least nonnegative representatives scale by the same factor, giving~\eqref{v50:eq-leader-scaling}.
For every positive integer \(t\), \(tca\equiv0\pmod N\) is equivalent to \(ta\equiv0\pmod n\), so additive orders agree.

Each \(L\in\Gamma_n(q)\) corresponds to \(cL\in\Gamma_N(q)\). Conversely, if \(A\in\Gamma_N(q)\) and \(c\mid A\),
then \(a=A/c\) satisfies \(0\le a<n\). Equation~\eqref{v50:eq-leader-scaling} gives
\(A=\lambda_N(A)=c\lambda_n(a)\), so \(a=\lambda_n(a)\in\Gamma_n(q)\).
This proves the filtering formula. Finally, \(L_1<L_2\) is equivalent to \(cL_1<cL_2\), proving strict order preservation and the corresponding ranking.
The integer cancellation above does not require \(c\) and \(n\) to be coprime.
\end{proof}

\begin{example}[Filter before ranking]\label{v50:ex-leader-filter-twentyone}
Take \(q=2\) and \(N=21=3\cdot7\). The coset leaders modulo twenty-one and their orbit lengths are
\[
(0,1),\quad(1,6),\quad(3,3),\quad(5,6),\quad(7,2),\quad(9,3).
\]
The nonzero orbits are
\[
\begin{aligned}
C_1^{(21)}&=\{1,2,4,8,16,11\},&
C_3^{(21)}&=\{3,6,12\},\\
C_5^{(21)}&=\{5,10,20,19,17,13\},&
C_7^{(21)}&=\{7,14\},\\
C_9^{(21)}&=\{9,18,15\}.&&
\end{aligned}
\]
Retaining the leaders \(0,3,9\) divisible by three and dividing by three gives the leaders \(0,1,3\) modulo seven,
with the same lengths \(1,3,3\). The second largest leader modulo twenty-one is seven, which is removed by the filter;
the second largest retained leader is three, corresponding to the second largest leader \(1\) modulo seven.

Here scaling is order preserving, whereas natural reduction need not preserve coset leaders.
For example, the leader \(9\) modulo twenty-one reduces to \(2\) modulo seven, whose coset has leader \(1\).

For a primitive twenty-first root of unity \(\zeta\), the element \(\eta=\zeta^3\) is a primitive seventh root, and
\[
\prod_{b\in C_a^{(7)}}(X-\eta^b)
=\prod_{A\in C_{3a}^{(21)}}(X-\zeta^A).
\]
Thus filtering preserves factor degrees and selects precisely the irreducible factors of \(X^7-1\) from the factorization of \(X^{21}-1\).
\end{example}

For an integer \(0\le b<n\), the monic polynomial of least degree in \(\F_q[X]\) vanishing at \(\zeta^i\) for \(1\le i\le b\) is
\begin{equation}\label{v50:eq-leader-interval-polynomial}
\begin{aligned}
G_b(X)&=\mathop{\operatorname{lcm}}_{1\le i\le b}M_{C_i^{(n)}}(X)
=\prod_{\substack{L\in\Gamma_n(q)\\1\le L\le b}}M_{C_L^{(n)}}(X),
\\[0.3em]
\deg G_b&=\sum_{\substack{L\in\Gamma_n(q)\\1\le L\le b}}|C_L^{(n)}|.
\end{aligned}
\end{equation}
Here \(\zeta\) is the primitive \(n\)-th root used to define the \(M_C\), and an empty product or least common multiple is one.
Indeed, a nonzero orbit meets \(\{1,\ldots,b\}\) if and only if its leader is at most \(b\).
Theorem~\ref{thm:coset-factorization} shows that distinct orbits give distinct irreducible factors, yielding the product and degree formulas.
This converts vanishing at consecutive roots of unity into a range condition on coset leaders.
In the cyclic-code model of Section~\ref{v50:sec-cyclic-trace-model},
if codeword polynomials must vanish at these roots of unity, \(G_b\) is the generator polynomial. This gives the root-set construction of BCH codes.
Consecutive roots also bound the minimum support size of a nonzero codeword; see Theorem~\ref{v56:thm-bch-bound} and Example~\ref{v56:ex-bch-seven}.

\noindent\textbf{Uniform cyclotomic-coset lengths versus short relations.} Let \(H\) be a multiplicative subgroup of order \(n\) and put \(d=\ord_n(q)\ge2\). If every nonidentity element \(h\in H\setminus\{1\}\) has minimal-polynomial degree \(d\), then all nonzero \(q\)-cyclotomic cosets modulo \(n\) have size \(d\), and \(n\equiv1\pmod d\). Equal orbit lengths, however, do not exclude short additive relations involving different orbits. Section~\ref{v95:sec-subgroup-periods} separates these issues through the linear relation length, cyclotomic intersection numbers, and Gaussian periods.

\section{Irreducibility of Cyclotomic Polynomials}
For a prime \(\ell\) and a positive integer \(a\), let \(v_\ell(a)\) be the largest nonnegative integer \(j\) with \(\ell^j\mid a\), called the \(\ell\)-adic valuation.
The notation \(\ell^j\Vert a\) means \(v_\ell(a)=j\), that is, \(\ell^j\mid a\) but \(\ell^{j+1}\nmid a\).
We use this notation in the group-theoretic argument over the integers below and the binomial criterion in Chapter~\ref{chap:root-fibers}.

\begin{lemma}[Lifting periods at a prime divisor]\label{lem:valuation-lift}
Let \(Q>1\) be an integer and \(\ell\mid Q-1\) a prime. If \(\ell\) is odd, or if \(\ell=2\) and \(Q\equiv1\pmod4\), then
\[
v_\ell(Q^s-1)=v_\ell(Q-1)+v_\ell(s)\qquad(s\ge1).
\]
\end{lemma}
\begin{proof}
Put \(b=v_\ell(Q-1)\), so \(\ell^b\mid Q-1\) but \(\ell^{b+1}\nmid Q-1\).
First suppose \(\ell\nmid u\). Since \(Q\equiv1\pmod\ell\),
\[
\frac{Q^u-1}{Q-1}=1+Q+\cdots+Q^{u-1}\equiv u\not\equiv0\pmod\ell,
\]
we obtain \(v_\ell(Q^u-1)=b\).

Next, multiplying the exponent by \(\ell\) increases the valuation by exactly one.
Write \(U=1+c\ell^b\), with \(\ell\nmid c\); for \(\ell=2\), we also have \(b\ge2\).
The binomial expansion gives
\[
U^\ell-1=c\ell^{b+1}
+\sum_{i=2}^{\ell-1}\binom{\ell}{i}c^i\ell^{bi}
+c^\ell\ell^{b\ell}.
\]
For odd \(\ell\), each intermediate term has valuation at least \(1+2b\ge b+2\),
and the last term has valuation \(b\ell\ge b+2\). For \(\ell=2\), the intermediate sum is empty and the last term has valuation \(2b\ge b+2\).
Only the first term has valuation \(b+1\), so no cancellation at that valuation is possible. Hence \(v_\ell(U^\ell-1)=b+1\).

Write \(s=\ell^c u\), where \(\ell\nmid u\).
Apply the first calculation to the exponent \(u\), then multiply the exponent by \(\ell\) successively \(c\) times. This gives \(v_\ell(Q^s-1)=b+c\), as required.
\end{proof}

The condition \(Q\equiv1\pmod4\) in the case \(\ell=2\) is essential. For example, when \(Q=3\) and \(s=2\),
\(v_2(3^2-1)=3\), whereas \(v_2(3-1)+v_2(2)=2\).

\Needspace{5\baselineskip}
\begin{proposition}\label{prop:cyclotomic-irreducible}
Let \(\gcd(M,q)=1\). Then \(\Phi_M\) is irreducible over \(\F_q\) if and only if \(\ord_M(q)=\varphi(M)\), equivalently, \(q\) generates \((\Z/M\Z)^\times\).
This can occur only for
\[
M=1,\ 2,\ 4,\ \ell^a,\ 2\ell^a\qquad(\ell\text{ an odd prime},\ a\ge1).
\]
For these moduli, one must still check whether the particular \(q\) is a generator.
\end{proposition}
\begin{proof}
By Corollary~\ref{cor:cyclotomic-equal-degree}, \(\Phi_M\) has a single irreducible factor if and only if \(\ord_M(q)=\varphi(M)\).
The unit group modulo \(M\) has order \(\varphi(M)\), so this is equivalent to \(q\) generating that group.

We prove the required classification of cyclic unit groups modulo positive integers. For an odd prime \(\ell\),
\((\Z/\ell\Z)^\times\) is the multiplicative group of a finite field, so some integer \(g\) has order \(\ell-1\) modulo \(\ell\).
For \(0\le c<\ell\), we have
\[
(g+c\ell)^{\ell-1}
\equiv g^{\ell-1}+(\ell-1)g^{\ell-2}c\ell\pmod{\ell^2}.
\]
Subtract one from the right side and divide by \(\ell\). Its residue modulo \(\ell\) is a nonconstant linear function of \(c\).
At most one value of \(c\) makes it zero. Choosing another value and setting \(u=g+c\ell\) gives \(v_\ell(u^{\ell-1}-1)=1\).
Applying Lemma~\ref{lem:valuation-lift} to \(Q=u^{\ell-1}\) yields
\[
v_\ell\bigl(u^{(\ell-1)s}-1\bigr)=1+v_\ell(s)\qquad(s\ge1).
\]
If \(u^t\equiv1\pmod{\ell^a}\), reduction modulo \(\ell\) first gives \(\ell-1\mid t\);
putting \(s=t/(\ell-1)\), the valuation formula then gives \(\ell^{a-1}\mid s\).
Conversely, the exponent \((\ell-1)\ell^{a-1}\) satisfies this congruence, so \(u\) has order \((\ell-1)\ell^{a-1}=\varphi(\ell^a)\).
Thus the unit group modulo an odd prime power is cyclic.
The Chinese remainder theorem (CRT) over the integers identifies the unit group modulo \(2\ell^a\) with that modulo \(\ell^a\).
The unit groups modulo one, two, and four have one, one, and two elements, respectively, and are cyclic.

Finally, we exclude the other moduli. If a modulus has two distinct odd prime divisors, its unit-group decomposition under the CRT has at least two components of even order;
each contains an element of order two, so their direct product has at least three elements of order two and cannot be cyclic.
The same argument applies when the modulus is divisible by both four and an odd prime.
If \(2^a\mid M\) with \(a\ge3\), the four distinct units modulo \(2^a\)
\[
1,\quad -1,\quad 1+2^{a-1},\quad -1+2^{a-1}
\]
all square to one, so this unit group is not cyclic.
The CRT makes it a direct factor of the unit group modulo \(M\). Since a quotient of a cyclic group is cyclic, the latter group is not cyclic either.
The remaining moduli are exactly those listed.
\end{proof}

\section{Base Extension and Orbit Step Size}
\begin{theorem}\label{thm:factor-base-extension}
Let \(f\in\F_q[X]\) be monic irreducible of degree \(d\). For a positive integer \(r\), put \(g=\gcd(d,r)\). Then in \(\F_{q^r}[X]\),
\[
\boxed{\quad f=f_0f_1\cdots f_{g-1},\qquad \deg f_j=\frac d g,\quad}
\]
where the \(f_j\) are distinct monic irreducible polynomials. If \(\alpha\) is a root of \(f\), we may take
\begin{equation}\label{eq:base-extension-orbits}
f_j(X)=\prod_{t=0}^{d/g-1}\bigl(X-\alpha^{q^{j+tr}}\bigr),\qquad 0\leq j<g.
\end{equation}
For a nonmonic irreducible polynomial \(f\) with leading coefficient \(c\), the factorization is \(f=c\,f_0f_1\cdots f_{g-1}\), with the monic factors \(f_j\) still given by~\eqref{eq:base-extension-orbits}.
\end{theorem}
\begin{proof}
In the common extension \(M=\F_{q^{\lcm(d,r)}}\), the roots of \(f\) are the distinct elements \(\alpha_i=\alpha^{q^i}\), with indices in \(\Z/d\Z\).
Frobenius relative to the new base field \(\F_{q^r}\) is \(\tau=\sigma^r\), so \(\tau(\alpha_i)=\alpha_{i+r}\). Starting at any index, the least positive return time \(t\) satisfies \(d\mid tr\). Write \(d=gd'\), \(r=gr'\), where \(\gcd(d',r')=1\). This condition is equivalent to \(d'\mid t\).
Hence each \(\tau\)-orbit has length exactly \(d/g\).

In \(\Z/d\Z\), the subgroup generated by \(r\) is \(g\Z/d\Z\), whose cosets have representatives \(0,1,\ldots,g-1\). Thus there are exactly \(g\) disjoint orbits covering all roots.
The product in \eqref{eq:base-extension-orbits} is taken over one such orbit. Its coefficients are fixed by \(\tau\), hence belong to \(\F_{q^r}\).
Applying the orbit and minimal-polynomial theorem over the new base field, each product is the minimal polynomial of any of its roots, and is therefore irreducible.
Disjoint orbits give distinct factors; together they cover all roots, and monicity gives the product identity.
\end{proof}

\begin{figure}[!htbp]
\centering
\begin{tikzpicture}[
  point/.style={circle,draw=Jade,fill=CoverMint,minimum size=7mm,inner sep=1pt,font=\small},
  orbit/.style={-{Stealth[length=2mm]},Jade,semithick},
  every node/.style={font=\small}]
\node[point] (a0) at (0,1.45) {\(0\)};
\node[point] (a1) at (1.26,.725) {\(1\)};
\node[point] (a2) at (1.26,-.725) {\(2\)};
\node[point] (a3) at (0,-1.45) {\(3\)};
\node[point] (a4) at (-1.26,-.725) {\(4\)};
\node[point] (a5) at (-1.26,.725) {\(5\)};
\foreach \a/\b in {0/1,1/2,2/3,3/4,4/5,5/0}
  \draw[orbit] (a\a)--(a\b);
\node at (0,0) {\(\sigma\)};
\draw[-{Stealth[length=2.4mm]},Clay,thick] (2.05,0)--(3.25,0)
  node[midway,above=3pt] {\(\sigma\mapsto\sigma^4\)};
\node[point] (b0) at (4.8,1.0) {\(0\)};
\node[point] (b4) at (3.9,-.7) {\(4\)};
\node[point] (b2) at (5.7,-.7) {\(2\)};
\draw[orbit] (b0)--(b4);
\draw[orbit] (b4)--(b2);
\draw[orbit] (b2)--(b0);
\node at (4.8,-.1) {\(\sigma^4\)};
\node[point] (c1) at (8.1,1.0) {\(1\)};
\node[point] (c5) at (7.2,-.7) {\(5\)};
\node[point] (c3) at (9.0,-.7) {\(3\)};
\draw[orbit] (c1)--(c5);
\draw[orbit] (c5)--(c3);
\draw[orbit] (c3)--(c1);
\node at (8.1,-.1) {\(\sigma^4\)};
\node[align=center] at (0,-2.15) {One orbit of length six};
\node[align=center] at (6.45,-2.15) {\(\gcd(6,4)=2\) orbits, each of length \(6/2=3\)};
\end{tikzpicture}
\caption{Changing the Frobenius step size splits an orbit. Arrows represent one specified iteration, not field inclusions.}
\label{fig:frobenius-step}
\end{figure}
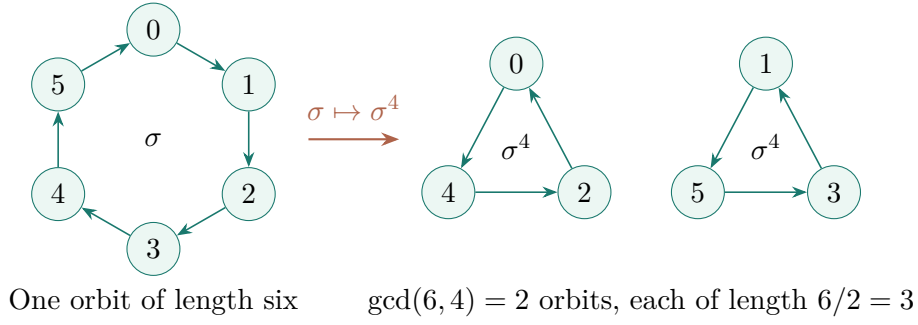

Figure~\ref{fig:frobenius-step} illustrates the case \(d=6\) and \(r=4\). The index \(i\) at each vertex denotes the root \(\alpha_i=\alpha^{q^i}\), read modulo six.
The original action increases \(i\) by one; after base extension, it increases \(i\) by four, splitting the roots into two three-element orbits with even and odd indices.
Accordingly, a degree-six irreducible polynomial splits over \(\F_{q^4}\) into two irreducible cubic factors.

\begin{corollary}\label{cor:base-extension-field-degree}
For the root \(\alpha\) above, we have
\[
\begin{aligned}
\F_{q^d}\cap\F_{q^r}&=\F_{q^g},&
\F_{q^r}(\alpha)&=\F_{q^{\lcm(d,r)}},\\
[\F_{q^r}(\alpha):\F_{q^r}]&=\frac d g.&
\end{aligned}
\]
In particular, \(f\) remains irreducible over \(\F_{q^r}\) if and only if \(\gcd(d,r)=1\), and splits completely if and only if \(d\mid r\).
\end{corollary}
\begin{proof}
Since \(\F_q(\alpha)=\F_{q^d}\), its intersection and compositum with the new base field follow from the subfield structure theorem.
The relative degree of the compositum is \(\lcm(d,r)/r=d/g\), in agreement with the orbit calculation above.
\end{proof}

\begin{proposition}[The exact coefficient field of a factor]\label{prop:factor-coefficient-field}
Under the hypotheses of Theorem~\ref{thm:factor-base-extension}, the coefficients of each factor \(f_j\)
generate precisely the intersection \(\F_{q^g}=\F_{q^d}\cap\F_{q^r}\) over \(\F_q\). In particular, these factors are already defined and irreducible over \(\F_{q^g}\).
\end{proposition}
\begin{proof}
On the common extension, let \(\sigma\) act coefficientwise on polynomials. The indices of the roots of \(f_j\) are the classes congruent to \(j\) modulo \(g\),
so \(\sigma(f_j)=f_{j+1}\), with factor indices read modulo \(g\).
The factors are distinct, so \(\sigma^t(f_j)=f_j\) if and only if \(g\mid t\).
In particular, \(\sigma^g\) fixes every coefficient, so all coefficients lie in \(\F_{q^g}\).
If they generate a subfield \(\F_{q^e}\subseteq\F_{q^g}\), then \(e\mid g\) and \(\sigma^e(f_j)=f_j\).
Thus \(g\mid e\) as well, and hence \(e=g\).
Finally, \(f_j\) is irreducible over the larger field \(\F_{q^r}\), and therefore also over its coefficient field \(\F_{q^g}\).
\end{proof}

The intersection of the enlarged base field with the original splitting field determines the number of factors; the coefficients of each factor belong to this intersection and generate it over \(\F_q\).
Thus the intersection and factorization formulas not only give the same degrees, but also identify the least coefficient field over which the factorization occurs.

\begin{example}\leavevmode\phantomsection\label{exm:base-extension-four}
Example~\ref{exm:quartic-model} proves that \(f=X^4+X+1\) is irreducible over \(\F_2\). Put \(\F_4=\F_2(\omega)\), where \(\omega^2+\omega+1=0\). Then
\[
X^4+X+1=(X^2+X+\omega)(X^2+X+\omega^2).
\]
The theorem ensures that both quadratic factors are irreducible over \(\F_4\). Over \(\F_8\), the original polynomial remains irreducible;
over \(\F_{16}\), it splits into four linear factors. The extent of factorization does not vary monotonically with the size of the base field:
it is governed by divisibility relations between extension degrees.
\end{example}

\section{Multiple Equal-Difference Structure of Cyclotomic Cosets and Binomial Descent}
\label{sec:cyclotomic-med}

Extending the base field splits a Frobenius orbit into suborbits with a larger step size.
For roots of unity, some suborbits also form complete equal-difference blocks, whose root polynomials are binomials.
This relation explains the additive symmetry within exponent orbits, the forced zero coefficients of the factors, and the appearance of lower-degree polynomials suited to computation.
Throughout this section, fix \(\gcd(N,q)=1\), and retain the notation \(C_i=\{iq^j\bmod N:j\ge0\}\) and \(M_{C_i}\).

For multiple equal-difference structure and its correspondence with binomial factorization, see \cite{ZZLW2025}. We describe these decompositions uniformly using complete additive cosets.

\par\Needspace{12\baselineskip}
\subsection{Complete Equal-Difference Root Blocks and the Coarsest Partition}
\begin{definition}\label{def:med-cosets}
For \(r\mid N\), a set
\[
B=b+\frac Nr\Z/N\Z
=\left\{b,b+\frac Nr,\ldots,b+(r-1)\frac Nr\right\}
\]
is called a \textbf{complete equal-difference block} of size \(r\).
A partition of a nonempty set into disjoint complete equal-difference blocks of the same size is called a \textbf{multiple equal-difference decomposition}.
Each block must be a complete coset of the corresponding additive subgroup; an arbitrary segment of a finite arithmetic progression need not satisfy this condition.
The decomposition records neither the order of the blocks nor the starting point of a block. The case \(r=1\) gives the singleton partition.
A partition is finer than another if each of its blocks is contained in a block of the latter. A partition coarser than every admissible partition is called the \textbf{coarsest partition}.
\leavevmode\phantomsection\label{def:additive-stabilizer}\ignorespaces For a nonempty set \(C\subseteq\Z/N\Z\), its \textbf{additive stabilizer} is
\[
H_+(C)=\{u\in\Z/N\Z:C+u=C\}.
\]
\end{definition}

\begin{proposition}[Binomials and complete equal-difference root sets]\label{prop:ed-binomial-correspondence}
Let \(\zeta\) be a primitive \(N\)-th root of unity, let \(\varnothing\ne B\subseteq\Z/N\Z\), and put
\(F_B(X)=\prod_{a\in B}(X-\zeta^a)\). Then \(F_B\) is a binomial if and only if \(B\) is a complete equal-difference block.
More precisely, if \(B=b+(N/r)\Z/N\Z\), then
\[
\boxed{\quad F_B(X)=X^r-\zeta^{br}.\quad}
\]
In particular, \(M_C\) is an irreducible binomial if and only if the cyclotomic coset \(C\) itself is a complete equal-difference block.
\end{proposition}
\begin{proof}
If \(B=b+(N/r)\Z/N\Z\), its root set is \(\zeta^b\mu_r\). Since \(p\nmid N\), the group \(\mu_r\) has exactly \(r\) elements.
These are therefore all the roots of \(X^r-\zeta^{br}\), and comparison of leading coefficients gives the identity.

Conversely, suppose \(F_B=X^r-c\). Its constant term is nonzero and its roots are distinct, so \(c\ne0\) and \(p\nmid r\).
Choose a root \(\alpha=\zeta^b\). The full root set is then \(\alpha\mu_r\).
Every root, as well as \(\alpha\), belongs to \(\mu_N\); taking their ratios gives \(\mu_r\subseteq\mu_N\), and hence \(r\mid N\).
The cyclic group \(\mu_N\) has a unique subgroup of order \(r\), namely \(\mu_r=\langle\zeta^{N/r}\rangle\).
Thus \(B=b+(N/r)\Z/N\Z\). The final assertion follows immediately from the irreducibility of \(M_C\).
\end{proof}

A complete equal-difference root block determines a binomial.
For all irreducible factors after base extension to be binomials, the corresponding Frobenius suborbits must be precisely such blocks.
We first use the additive stabilizer to determine all possible decompositions, then calculate the least iteration step that produces them.

\begin{theorem}[All multiple equal-difference decompositions of a cyclotomic coset]\label{thm:med-coset-partitions}
Let \(C=C_i\), \(\tau=|C|\), and \(\kappa_C=|H_+(C)|\). Then \(\kappa_C\mid\gcd(N,\tau)\), and
\[
H_+(C)=\frac N{\kappa_C}\Z/N\Z.
\]
Every complete equal-difference block in \(C\) is contained in an additive coset of \(H_+(C)\).
These cosets, each of size \(\kappa_C\), form the unique coarsest multiple equal-difference decomposition.
The admissible block sizes are precisely the positive divisors of \(\kappa_C\). For each such \(r\), the decomposition is unique and consists of the nonempty residue classes of \(C\) modulo \(N/r\).

Put \(\omega_C=\tau/\kappa_C\), and for \(t\ge1\) define \(C_{a;t}=\{aq^{tj}\bmod N:j\ge0\}\).
All multiple equal-difference decompositions can also be written as
\begin{equation}\label{eq:med-step-partitions}
\boxed{\quad
C_i=\bigsqcup_{j=0}^{t-1}C_{iq^j;t},
\qquad \omega_C\mid t\mid\tau.
\quad}
\end{equation}
The corresponding block size is \(\tau/t\). For any positive integer \(m\), all \(q^m\)-suborbits are complete equal-difference blocks if and only if \(\omega_C\mid m\).
\end{theorem}
\begin{proof}
Translations preserving \(C\) are closed under composition and inversion, so \(H_+(C)\) is an additive subgroup.
The set \(C\) is a disjoint union of cosets of this subgroup, so the subgroup order divides both \(N\) and \(\tau\).
Uniqueness of a subgroup of a given order in a cyclic group gives the displayed expression for \(H_+(C)\).

Let \(b+K\subseteq C\) be a complete equal-difference block, where \(K\) is the corresponding additive subgroup.
For any \(a\in C\), there is a \(j\) such that \(a=q^jb\). Multiplication by the unit \(q^j\) preserves \(K\), so
\[
a+K=q^j(b+K)\subseteq C.
\]
Taking the union over all \(a\in C\) gives \(C+K=C\), and hence \(K\subseteq H_+(C)\).
Thus every complete equal-difference block is contained in a maximal block.
Conversely, for any subgroup \(K\subseteq H_+(C)\), the \(K\)-cosets contained in \(C\) form a partition of \(C\).
The classification of subgroups of a cyclic group now gives all admissible block sizes, uniqueness for each size, and uniqueness of the coarsest decomposition.

Fix a block size \(r\mid\kappa_C\). Multiplication by \(q\) permutes the blocks transitively,
since its action on \(C\) is a single \(\tau\)-cycle. Thus the blocks form a cycle of length \(t=\tau/r\).
Multiplication by \(q^t\) preserves the block containing \(i\), and the orbit of \(i\) under this map has length \(\tau/t=r\), so this suborbit is exactly that block.
The other blocks are likewise the suborbits in \eqref{eq:med-step-partitions}.
The condition \(r\mid\kappa_C\) is equivalent to \(\omega_C\mid t\mid\tau\), giving all decompositions.

For arbitrary \(m\), put \(g=\gcd(m,\tau)\). The \(q^m\)-suborbits and \(q^g\)-suborbits of the original cycle coincide.
They are complete equal-difference blocks if and only if \(\omega_C\mid g\). Since \(\omega_C\mid\tau\), this is equivalent to \(\omega_C\mid m\).
\end{proof}

Transitivity of multiplication by \(q\) carries any complete equal-difference block to blocks covering the orbit, thereby partitioning the entire coset into blocks of equal size.

\begin{proposition}[An additive-stabilizer criterion for a multiplicative orbit]\label{v74:prop-orbit-additive-stabilizer}
Let \(R=\Z/N\Z\), \(H\le R^\times\), and \(C=uH\) with \(u\in R^\times\). For any additive subgroup \(I\) of \(R\),
\begin{equation}\label{v74:eq-additive-principal-unit}
C+I=C\quad\Longleftrightarrow\quad1+I\subseteq H.
\end{equation}
Thus \(H_+(C)\) is the largest additive subgroup satisfying the right-hand condition, independently of the unit representative \(u\).
\end{proposition}
\begin{proof}
Every additive subgroup of \(R\) is an ideal and is mapped onto itself by multiplication by a unit.
If \(C+I=C\), then \(u\in C\) gives \(u+I\subseteq uH\); multiplying by \(u^{-1}\) gives \(1+I\subseteq H\).
Conversely, suppose \(1+I\subseteq H\). For \(h\in H\) and \(i\in I\),
\[
h+i=h(1+h^{-1}i)\in H.
\]
Hence \(H+I=H\), and \(uH+I=u(H+I)=uH\). The last assertion follows from the definition of the additive stabilizer.
\end{proof}

The same condition has two meanings: additively, the orbit is a union of complete cosets; multiplicatively, \(H\) contains \(1+I\).
For \(H=\langle q\rangle\), the multiplier group therefore identifies the multiple equal-difference structure directly.
Theorem~\ref{v74:thm-med-nonaffine} will also show that, at odd prime-power moduli, a nontrivial structure of this kind corresponds exactly to extra non-affine permutations preserving difference orbits.

\begin{example}[Successive reduction and complete equal-difference orbits]\leavevmode\phantomsection\label{exm:med-nine}
Take \(q=4,N=9\). Starting at one and repeatedly multiplying by four and reducing modulo nine gives
\[
1\longmapsto4\longmapsto7\longmapsto1.
\]
Starting at two gives \(2\mapsto8\mapsto5\mapsto2\). The full list of cyclotomic cosets modulo nine is
\[
\{0\},\quad C_1=\{1,4,7\},\quad C_2=\{2,5,8\},\quad
C_3=\{3\},\quad C_6=\{6\}.
\]
Here the elements are listed in increasing order; the order in which successive multiplication by four visits them need not be retained. Both nonsingleton cosets are complete equal-difference blocks:
\[
C_1=1+3\Z/9\Z,\qquad C_2=2+3\Z/9\Z.
\]
For \(C_1\), for example, \(H_+(C_1)=\{0,3,6\}\), so \(\tau=3\), \(\kappa_C=3\), and \(\omega_C=1\).
The only admissible block sizes are three and one, giving the whole coset and the singleton partition, respectively.

Choose a primitive ninth root of unity \(\zeta\), and put \(\omega=\zeta^3\). Then \(\omega\in\F_4\) and \(\omega^2+\omega+1=0\). The root-product formula for complete equal-difference blocks gives
\[
M_{C_1}(X)=X^3-\omega,\qquad M_{C_2}(X)=X^3-\omega^2.
\]
Both are irreducible over the base field \(\F_4\). The five cosets together give
\[
X^9-1=(X-1)(X-\omega)(X-\omega^2)(X^3-\omega)(X^3-\omega^2).
\]
Thus the complete equal-difference orbits obtained by successive reduction correspond to irreducible binomials through their root products.
\end{example}

\begin{example}[A coset with no nonsingleton equal-difference block]\label{exm:med-fifteen-trivial}
Take \(q=2,N=15\). Then
\[
1\longmapsto2\longmapsto4\longmapsto8\longmapsto1,
\qquad C_1=\{1,2,4,8\}.
\]
Since \(\kappa_C\mid\gcd(15,4)=1\), we have \(H_+(C_1)=\{0\}\), \(\kappa_C=1\), and \(\omega_C=4\).
This coset contains no complete equal-difference block of size greater than one; its only multiple equal-difference decomposition is the singleton partition.
Although \(1,2\) form a finite arithmetic progression, they are not a complete equal-difference block modulo fifteen: their common difference one generates all of \(\Z/15\Z\), not this two-element set.
The \(q^2\)-suborbits \(\{1,4\}\) and \(\{2,8\}\) are not complete equal-difference blocks either. All suborbits become singletons only when the step size is a multiple of four.
\end{example}

\subsection{An Explicit Formula for the Least Step and the Dyadic Exception}
\leavevmode\phantomsection\label{def:squarefree-kernel}\ignorespaces For a positive integer \(d\), its \textbf{square-free kernel} is \(\rad(d)=\prod_{\ell\mid d}\ell\), with \(\rad(1)=1\).
\begin{proposition}[An arithmetic criterion for a complete equal-difference orbit]\label{prop:ed-coset-criterion}
Let \(C=C_i\), and put \(d=N/\gcd(N,i)\). Then \(C\) itself is a complete equal-difference block if and only if
\begin{equation}\label{eq:ed-coset-criterion}
\rad(d)\mid q-1,
\qquad 8\mid d\ \Longrightarrow\ q\equiv1\pmod4.
\end{equation}
When these conditions hold, \(\tau=\ord_d(q)=d/\gcd(d,q-1)\).
\end{proposition}
\begin{proof}
For \(d=1\), the coset \(C\) is a singleton and the assertion follows from the conventions. Assume henceforth that \(d>1\).
Put \(g=N/d\), and write \(i=gu\), where \(\gcd(u,d)=1\).
The scaling map \(x\mapsto gx\) sends the unit orbit \(u\langle q\rangle\) modulo \(d\) to \(C\).
If a complete equal-difference block lies in \(g\Z/N\Z\), differences of its elements show that its additive subgroup also lies in this subgroup, so the scaling can be reversed.
Multiplication by the unit \(u^{-1}\) also preserves complete equal-difference blocks. It therefore suffices to consider \(\langle q\rangle\) modulo \(d\).

If this orbit is a complete equal-difference block, it contains one and hence equals \(1+\Delta\Z/d\Z\), where \(\Delta=d/\tau\).
Every element of this block is a unit, so \(\rad(d)\mid\Delta\):
otherwise, for some \(\ell\mid d\) with \(\ell\nmid\Delta\), one could choose an integer \(j\) such that \(1+j\Delta\equiv0\pmod\ell\), contradicting the unit condition.
Also, \(q\) belongs to the block, so \(\Delta\mid q-1\), and hence \(\rad(d)\mid q-1\).

To verify the dyadic condition, suppose \(2^a\Vert d\), \(a\ge3\), and \(q\equiv3\pmod4\). Put \(b=v_2(q+1)\ge2\).
An odd exponent cannot return \(q\) to one modulo \(2^a\). For an even exponent \(2s\), applying Lemma~\ref{lem:valuation-lift} to \(q^2\) gives
\[
v_2(q^{2s}-1)=b+1+v_2(s),
\qquad \ord_{2^a}(q)=2^{\max(1,a-b)}.
\]
On the other hand, \(\rad(d)\mid q-1\) makes every local order at an odd prime power odd.
The Chinese remainder theorem therefore gives \(v_2(\tau)=\max(1,a-b)\), so
\[
v_2(\Delta)=a-\max(1,a-b)=\min(a-1,b)\ge2.
\]
This contradicts \(\Delta\mid q-1\) and \(v_2(q-1)=1\), proving necessity.

Conversely, assume \eqref{eq:ed-coset-criterion}. For every odd prime power \(\ell^a\Vert d\),
Lemma~\ref{lem:valuation-lift} gives
\[
\ord_{\ell^a}(q)=\ell^{\max(0,a-v_\ell(q-1))}.
\]
At the dyadic component, the same formula holds if \(q\equiv1\pmod4\).
If \(q\equiv3\pmod4\), the condition requires \(a\le2\); direct calculation modulo two and four again gives the same formula.
These local orders are powers of distinct primes, so their least common multiple equals their product, giving
\[
\tau=\frac d{\gcd(d,q-1)}.
\]
Put \(\Delta=\gcd(d,q-1)\). Since \(q\equiv1\pmod\Delta\), all powers of \(q\) lie in \(1+\Delta\Z/d\Z\).
Both sets have \(d/\Delta=\tau\) elements, so they coincide, proving sufficiency.
\end{proof}

\leavevmode\phantomsection\label{v124:def-small-modulus}\ignorespaces
To retain the information at the odd primes together with the dyadic exception, put
\[
b(d)=\begin{cases}
\rad(d),&8\nmid d,\\
2\rad(d),&8\mid d.
\end{cases}
\]
Thus \(b(d)\mid d\), with \(b(1)=1\). When \(8\mid d\), the two-primary part of \(b(d)\) retains reduction modulo four.

\begin{theorem}[The least Frobenius step for multiple equal-difference structure]\label{thm:med-minimal-step}
Let \(C=C_i\), and put
\[
d=\frac N{\gcd(N,i)},\qquad \tau=\ord_d(q),\qquad
\rho=\ord_{\rad(d)}(q).
\]
The least step in Theorem~\ref{thm:med-coset-partitions} is
\begin{equation}\label{v124:eq-small-step}
\boxed{\quad
\omega_C=\ord_{b(d)}(q),\qquad
\kappa_C=\frac{\ord_d(q)}{\ord_{b(d)}(q)}.
\quad}
\end{equation}
Equivalently,
\begin{equation}\label{eq:med-minimal-step}
\omega_C=
\begin{cases}
2\rho,&8\mid d\text{ and }q^\rho\equiv3\pmod4,\\
\rho,&\text{otherwise},
\end{cases}
\qquad \kappa_C=\frac\tau{\omega_C}.
\end{equation}
In particular, the common difference of a maximal complete equal-difference block is \(N/\kappa_C\), and the coarsest decomposition has exactly \(\omega_C\) blocks.
\end{theorem}
\begin{proof}
The elements of every \(q^m\)-suborbit have the same reduced modulus \(d\). Applying Proposition~\ref{prop:ed-coset-criterion} with multiplier \(q^m\) gives precisely
\[
q^m\equiv1\pmod{\rad(d)},\qquad
8\mid d\ \Longrightarrow\ q^m\equiv1\pmod4.
\]
These conditions are equivalent to \(q^m\equiv1\pmod{b(d)}\), so the admissible steps are exactly the multiples of \(\ord_{b(d)}(q)\). Theorem~\ref{thm:med-coset-partitions} gives \(\omega_C=\ord_{b(d)}(q)\) and \(\kappa_C=\tau/\omega_C\).

If \(8\nmid d\), then \(b(d)=\rad(d)\). If \(8\mid d\), first write \(m=\rho u\). When \(q^\rho\equiv1\pmod4\), every \(u\) is allowed; when \(q^\rho\equiv3\pmod4\), precisely the even values of \(u\) are allowed. This proves the expanded formula.
\end{proof}

If \(d\mid D\), then \(b(d)\mid b(D)\), and hence \(\ord_{b(d)}(q)\mid\ord_{b(D)}(q)\). The least step at the largest reduced modulus therefore works at all of its divisor moduli.

The decompositions and parameter calculations above use only \(\gcd(N,q)=1\), so the arithmetic applies to any integer \(q>1\).
When \(q\) is a prime power, Proposition~\ref{prop:ed-binomial-correspondence} converts each complete equal-difference block into a binomial,
and Frobenius acting on the blocks determines the coefficient fields of these binomials and their product.

\begin{example}[Two nine-element blocks and three levels of decomposition]\label{exm:med-twentyseven}
Take \(q=2,N=27,i=1\). Successive multiplication by two visits the elements in the following order; the second line continues from the last entry of the first:
\[
\begin{aligned}
1&\longmapsto2\longmapsto4\longmapsto8\longmapsto16
\longmapsto5\longmapsto10\longmapsto20\longmapsto13,\\
13&\longmapsto26\longmapsto25\longmapsto23\longmapsto19
\longmapsto11\longmapsto22\longmapsto17\longmapsto7\longmapsto14\longmapsto1.
\end{aligned}
\]
These eighteen entries are precisely all units modulo twenty-seven, so \(\tau=18\). Also,
\[
\rho=\ord_3(2)=2,\qquad \omega_C=2,\qquad \kappa_C=9.
\]
Reordering the same set by residue modulo three displays the coarsest decomposition:
\[
\begin{aligned}
B_1&=\{1,4,7,10,13,16,19,22,25\}=1+3\Z/27\Z,\\
B_2&=\{2,5,8,11,14,17,20,23,26\}=2+3\Z/27\Z,\\
C_1&=B_1\sqcup B_2.
\end{aligned}
\]
Multiplication by two interchanges \(B_1,B_2\), while multiplication by four traverses each block; for example,
\[
1\longmapsto4\longmapsto16\longmapsto10\longmapsto13
\longmapsto25\longmapsto19\longmapsto22\longmapsto7\longmapsto1.
\]
Thus arithmetic-progression order describes the additive structure of the set, while arrows describe the multiplicative orbit; the two orders usually differ.

Since \(\kappa_C=9\), the admissible block sizes are \(9,3,1\), corresponding to Frobenius steps \(t=2,6,18\):
\begin{center}
\begin{tabular}{@{}cccc@{}}
\toprule
Block size \(r\)&Blocks \(t=18/r\)&Common difference \(27/r\)&Decomposition\\
\midrule
9&2&3&\(B_1,B_2\)\\
3&6&9&The six triples below\\
1&18&27&Eighteen singletons\\
\bottomrule
\end{tabular}
\end{center}
The intermediate level is
\[
\begin{aligned}
C_1={}&\{1,10,19\}\sqcup\{2,11,20\}\sqcup\{4,13,22\}\\
&\sqcup\{8,17,26\}\sqcup\{7,16,25\}\sqcup\{5,14,23\}.
\end{aligned}
\]
These blocks have representatives \(1,2,4,8,16,5\), respectively, and are precisely the \(q^6\)-suborbits; here \(2^6\equiv10\pmod{27}\).
For example, \(1\mapsto10\mapsto19\mapsto1\) and \(16\mapsto25\mapsto7\mapsto16\).
These decompositions correspond exactly to the three subgroups of \(H_+(C_1)=3\Z/27\Z\).
\end{example}

\begin{example}[Nonunit indices for the same modulus]\label{exm:med-nonunit}
Again take \(q=2,N=27\), but start at three. Direct calculation gives
\[
3\longmapsto6\longmapsto12\longmapsto24\longmapsto21\longmapsto15\longmapsto3.
\]
Here the reduced modulus is \(d=27/\gcd(27,3)=9\), so
\[
\tau=\ord_9(2)=6,\qquad\omega_{C_3}=2,\qquad\kappa_{C_3}=3.
\]
The corresponding coarsest decomposition is
\[
C_3=\{3,12,21\}\sqcup\{6,15,24\},\qquad
H_+(C_3)=9\Z/27\Z.
\]
For the index nine, \(C_9=\{9,18\}\), \(d=3\), and \(\tau=2\), so \(\omega_{C_9}=2\) and \(\kappa_{C_9}=1\).
This gives all cosets modulo twenty-seven:
\begin{center}
\begin{tabular}{@{}ccccc@{}}
\toprule
Coset&Reduced modulus \(d\)&Length \(\tau\)&Least step \(\omega_C\)&Maximal block size \(\kappa_C\)\\
\midrule
\(C_0=\{0\}\)&1&1&1&1\\
\(C_1\)&27&18&2&9\\
\(C_3\)&9&6&2&3\\
\(C_9\)&3&2&2&1\\
\bottomrule
\end{tabular}
\end{center}
The lengths sum to \(1+18+6+2=27\). For the same modulus, different indices may have different maximal equal-difference block sizes; the parameters depend on \(N/\gcd(N,i)\), not only on \(N\).
\end{example}

\begin{proposition}[Compressing coset-leader calculations by maximal equal-difference blocks]\label{v50:prop-med-leader-compression}
Let \(C=C_i\subseteq\Z/N\Z\), with \(\tau=|C|\), \(\omega_C\), and \(\kappa_C=\tau/\omega_C\)
as in Theorem~\ref{thm:med-minimal-step}, and put \(D=N/\kappa_C=N\omega_C/\tau\). Then
\begin{equation}\label{v50:eq-med-leader-compression}
\boxed{\quad
\lambda_N(i)=\min_{0\le j<\omega_C}[iq^j]_D.
\quad}
\end{equation}
In particular, if \(C\) itself is a complete equal-difference block, then \(\lambda_N(i)=[i]_{N/\tau}\).
\end{proposition}
\begin{proof}
Theorem~\ref{thm:med-coset-partitions} partitions \(C\) into \(\omega_C\) maximal complete equal-difference blocks,
the \(j\)-th being \(iq^j+D\Z/N\Z\). Put \(r_j=[iq^j]_D\). The least nonnegative representatives of the elements in this block are
\[
r_j,\quad r_j+D,\quad\ldots,\quad r_j+(\kappa_C-1)D.
\]
These integers are distinct and less than \(N\); the smallest is \(r_j\). Taking the minimum over all blocks gives the formula.
If the whole coset is a complete equal-difference block, then \(\omega_C=1\) and \(\kappa_C=\tau\), giving the last assertion.
\end{proof}
Theorem~\ref{v139:thm-period-compression} expresses the same compression through integer multiplication of periodic digits: a complete fiber has size equal to the exact period-loss factor, and the canonical root polynomial is a monomial composition. Proposition~\ref{v144:prop-profile-compression} applies the compression to all coefficient valuations. Example~\ref{v144:ex-profile-tower} reduces a 602-element orbit modulo $43^2$ to a fourteen-element orbit modulo 43, then recovers its phase profile in a three-dimensional fixed algebra.

The input has two parts: \(q,N,i\) specify the problem, while \(\tau=\ord_d(q)\) and \(\omega_C=\ord_{b(d)}(q)\) are the required arithmetic data. Once both orders are known, successive multiplication by \(q\), starting at \([i]_D\), uses \(\omega_C-1\) modular multiplications and at most as many comparisons; direct orbit traversal visits \(\tau\) residues. Computing the orders is a separate cost. Given a factored multiple \(M\) of an order, one can successively remove prime factors, using modular exponentiation to test whether the result remains a multiple of that order. A factorization of \(N\) alone does not automatically supply these orders or factorizations of their bounds.

\begin{example}[Two residue calculations for a 324-element orbit]\label{v50:ex-med-leader-compression}
Take \(q=5\), \(N=3888=2^4\cdot3^5\), and \(i=1001\). Since \(\gcd(1001,3888)=1\), we have
\[
\ord_{16}(5)=4,\qquad\ord_{243}(5)=162,\qquad
\tau=\operatorname{lcm}(4,162)=324.
\]
Also, \(\rad(3888)=6\), \(\rho=\ord_6(5)=2\), and \(5^2\equiv1\pmod4\),
so \(\omega_C=2\), \(\kappa_C=162\), and \(D=24\). Hence
\[
\lambda_{3888}(1001)=\min\{[1001]_{24},[5005]_{24}\}
=\min\{17,13\}=13.
\]
The coset leader is therefore found by comparing the least representatives of just two maximal blocks. Here \(\tau\) and \(\omega_C\) are calculated from multiplicative orders and valuations;
once these parameters are known, the leader search reduces from 324 orbit elements to two modular reductions.
\end{example}

Order-preserving scaling preserves the entire Frobenius orbit and factor degree across different moduli; maximal-block compression uses the additive symmetry within one orbit.
The latter gives both the coset-leader formula and the monomial composition decomposition in the next subsection.

\par\Needspace{12\baselineskip}
\subsection{Multiple Equal-Difference Decompositions, Monomial Composition, and Binomial Descent}
\leavevmode\phantomsection\label{def:coefficient-extraction}\ignorespaces
For a positive integer \(r\), write \(k[X^r]=\{G(X^r):G\in k[X]\}\), and let \([X^j]P\) denote the coefficient of \(X^j\) in \(P\).
For a nonconstant polynomial \(P\), the largest positive integer \(r\) such that \(P\in k[X^r]\) is its \textbf{maximum monomial composition exponent}.

\begin{theorem}[Stabilizers, monomial composition, and binomial descent]
\label{v125:thm-orbit-structure}\label{prop:med-binomial-descent}\label{cor:med-monomial-composition}\label{cor:med-binomial-extension}
Let \(C=C_i\subseteq\Z/N\Z\), \(\alpha=\zeta^i\), and put
\[
d=N/\gcd(N,i),\qquad \tau=|C|=\ord_d(q),\qquad
\kappa_C=|H_+(C)|,\qquad \omega_C=\tau/\kappa_C.
\]
Then
\[
H_+(C)=(N/\kappa_C)\Z/N\Z,\qquad
\omega_C=\ord_{b(d)}(q),\qquad
\kappa_C=\frac{\ord_d(q)}{\ord_{b(d)}(q)}.
\]
These numbers determine the following structures simultaneously.
\begin{enumerate}
\item For every positive integer \(r\), one has \(M_C\in k[X^r]\) if and only if \(r\mid\kappa_C\). In particular,
\[
\kappa_C=\gcd\{j\ge1:[X^j]M_C\ne0\}.
\]
\item If \(r\mid\kappa_C\) and \(t=\tau/r\), the minimal polynomial \(G_r\in k[Y]\) of \(\alpha^r\) has degree \(t\), and
\begin{equation}\label{eq:med-binomial-descent}
M_C(X)=G_r(X^r)=\prod_{j=0}^{t-1}(X^r-\alpha^{rq^j}).
\end{equation}
Each binomial on the right is irreducible in \(\F_{q^t}[X]\); its root exponents form one of the complete equal-difference blocks of size \(r\) in \(C\).
\item For \(m\ge1\), every irreducible factor of \(M_C\) over \(\F_{q^m}\) is a binomial if and only if \(\omega_C\mid m\). The least extension degree is \(\omega_C\); over that extension there are exactly \(\omega_C\) irreducible binomials of degree \(\kappa_C\).
\end{enumerate}
\end{theorem}
\begin{proof}
The form of the stabilizer and the arithmetic formula for \(\omega_C\) follow from Theorems~\ref{thm:med-coset-partitions} and~\ref{thm:med-minimal-step}, respectively.
If \(r\mid\kappa_C\), the blocks of size \(r\) are represented by \(iq^j\), \(0\le j<t\).
Proposition~\ref{prop:ed-binomial-correspondence} gives the product in~\eqref{eq:med-binomial-descent}.
Two roots have the same \(r\)-th power precisely when their exponents agree modulo \(N/r\), so the power map collapses exactly these blocks.
Frobenius acts on the blocks as a cycle of length \(t\); hence \(G_r\) is the degree-\(t\) minimal polynomial of \(\alpha^r\).
The tower formula gives \([k(\alpha):k(\alpha^r)]=r\), so \(X^r-\alpha^r\) is the relative minimal polynomial. Its conjugates give the remaining irreducible binomials.

Conversely, suppose \(M_C=G(X^r)\). Separability gives \(p\nmid r\). For every \(\varepsilon\in\mu_r\), the element \(\alpha\varepsilon\) is again a root.
Both roots lie in \(\mu_N\), so \(\mu_r\subseteq\mu_N\) and \(r\mid N\). Applying this to every root gives
\((N/r)\Z/N\Z\subseteq H_+(C)\), or \(r\mid\kappa_C\).
Moreover, \(M_C\in k[X^r]\) means exactly that \(r\) divides every positive exponent with nonzero coefficient; their greatest common divisor is therefore \(\kappa_C\).

After base extension, irreducible factors correspond to \(q^m\)-suborbits.
Proposition~\ref{prop:ed-binomial-correspondence} identifies the binomial condition with these suborbits being complete equal-difference blocks; Theorem~\ref{thm:med-coset-partitions} identifies the latter condition with \(\omega_C\mid m\). At \(m=\omega_C\), the identity \(\tau=\kappa_C\omega_C\) gives the number and degrees of the factors.
\end{proof}

Thus the stabilizer order \(\kappa_C\) measures maximal composition compression over the original field, while the quotient orbit length \(\omega_C\) is the least extension degree for binomial descent.
For \(\beta=\alpha^{\kappa_C}\),
\[
[k(\beta):k]=\omega_C,\qquad [k(\alpha):k(\beta)]=\kappa_C.
\]
Return maps for power and affine compression are compared in Section~\ref{sec:power-affine-comparison}.

\begin{proposition}[Complete orbit fibers and irreducible composition]\label{v142:prop-full-fibers}
Let $n>1$, $\gcd(n,q)=1$, $d\mid n$, $d>1$, and put $s=n/d$.
Choose a primitive $n$th root $\xi$, set $\eta=\xi^s$, and write
\[
 m=\ord_n(q),\qquad r=\ord_d(q),\qquad
 h=m_\xi,\qquad h_0=m_\eta.
\]
Reduction modulo $d$ maps the coset $C_n(1)$ modulo $n$ onto $C_d(1)$, with $m/r\le s$ elements in each fiber. The following conditions are equivalent:
\begin{enumerate}
\item Every nonempty fiber is a complete fiber of reduction modulo $d$.
\item $m=sr$.
\item $h(X)=h_0(X^s)$.
\item $h_0(X^s)$ is irreducible in $k[X]$.
\end{enumerate}
They are also equivalent to $s\mid |H_+(C_n(1))|$. Under these conditions,
\[
 [k(\xi):k(\eta)]=s,
 \qquad X^s-\eta\text{ is the minimal polynomial of $\xi$ over $k(\eta)$.}
\]
\end{proposition}
\begin{proof}
The two orbits have lengths $m,r$, and $r\mid m$. The exponents $q^i,q^j$ agree modulo $d$ exactly when $r\mid i-j$, so each fiber has $m/r$ elements. Each ambient fiber has $s$ residue classes, proving the equivalence of the first two conditions.
Since $h_0(\xi^s)=0$, we have $h\mid h_0(X^s)$. Both polynomials are monic, of degrees $m,sr$, respectively, proving the equivalence of the last three conditions.
Complete fibers are cosets of $d\mathbb Z/n\mathbb Z$, so the first condition says precisely that this subgroup of order $s$ lies in the additive stabilizer.
Finally, the tower law gives relative degree $m/r=s$; the degree-$s$ polynomial $X^s-\eta$ vanishes at $\xi$ and is therefore its relative minimal polynomial.
\end{proof}

This criterion tests composition along a specified reduction of the modulus, whereas Theorem~\ref{v125:thm-orbit-structure} lists all admissible compression exponents at once. Complete fibers will also yield independence of the trace-code coordinate blocks in Corollary~\ref{v142:cor-irreducible-interleaving}.

\begin{example}[Five-element fibers of twenty-fifth roots]\label{v142:ex-twenty-five-roots}
Take $q=23$, $n=25$, and $d=s=5$. Since $23\equiv3\pmod5$, we have $\ord_5(23)=4$. Moreover,
\[
 23^4\equiv16=1+3\cdot5\pmod{25},
 \qquad (1+3\cdot5)^j\equiv1+3j\cdot5\pmod{25},
\]
so $\ord_{25}(23)=20=5\cdot4$. The coset $C_{25}(1)$ consists of all twenty unit residue classes. Its fibers over $C_5(1)=\{1,2,3,4\}$ are
\[
 a+5\mathbb Z/25\mathbb Z,\qquad a=1,2,3,4.
\]
Thus $H_+(C_{25}(1))=5\mathbb Z/25\mathbb Z$, and the maximum composition exponent is five. For a primitive root $\xi$ and $\eta=\xi^5$,
\[
 \begin{aligned}
 h_0(Y)&=\Phi_5(Y)=Y^4+Y^3+Y^2+Y+1,\\
 h(X)&=\Phi_{25}(X)=X^{20}+X^{15}+X^{10}+X^5+1=h_0(X^5).
 \end{aligned}
\]
Both polynomials are irreducible over $\F_{23}$ and give the tower $\F_{23}\subset\F_{23^4}\subset\F_{23^{20}}$. Example~\ref{v142:ex-twenty-five-code} will use this same composition to determine a cyclic code's distance.
\end{example}

\begin{remark}[Coefficients after compression]\label{v125:rem-compressed-coefficients}
Retain $\alpha=\zeta^i$ and $\beta=\alpha^{\kappa_C}$. Write \(M_C(X)=G(X^{\kappa_C})\), where
\(G(Y)=Y^{\omega_C}+c_1Y^{\omega_C-1}+\cdots+c_{\omega_C}\). Then
\[
c_j=(-1)^j\!\sum_{0\le i_1<\cdots<i_j<\omega_C}
\beta^{q^{i_1}+\cdots+q^{i_j}},\qquad
c_1=-\sum_{i=0}^{\omega_C-1}\beta^{q^i}.
\]
Substitution of \(Y=X^{\kappa_C}\) preserves the number of nonzero terms. The stabilizer forces coefficients at exponents not divisible by \(\kappa_C\) to vanish; further zeros are determined by these symmetric functions. The first is the negative trace of the compressed root.
If \(\omega_C>1\), then \(G\) has at least three nonzero terms: otherwise \(M_C\) itself would be a binomial, contradicting its maximum composition exponent \(\kappa_C<\tau\).
Section~\ref{v79:sec-defining-equations} studies representation complexity after changing the generator; sparse multiples in Chapter~\ref{chap:codes} also allow the degree of the defining relation to vary.
\end{remark}

\begin{example}[Equal stabilizers and different coefficient supports]\label{v125:ex-extra-cancellation}
In \(\F_3[Y]\), the polynomials
\[
G_1=Y^3+2Y+1,\qquad G_2=Y^3+Y^2+2Y+1
\]
have no roots in the base field and are therefore irreducible. If \(\beta\) is a root of either polynomial, then
\(\beta^{13}=\beta\beta^3\beta^9=-1\), while \(\beta\notin\F_3\); hence \(\beta\) has multiplicative order 26.
Choose \(\alpha^2=\beta\). Its order is 52, and \(\ord_{52}(3)=6\) shows that both
\[
G_1(X^2)=X^6+2X^2+1,\qquad
G_2(X^2)=X^6+X^4+2X^2+1
\]
are irreducible. Their root-exponent orbits have the same \(d=52\), \(\kappa_C=2\), and \(\omega_C=3\), but the polynomials have three and four nonzero terms, respectively.
The missing \(X^4\) term in the first polynomial comes from the zero trace of \(\beta\), not from a larger additive stabilizer.
\end{example}

\begin{remark}[Equal-difference orbits and sparsest defining equations]\label{v79:rem-med-defining-equation}
Let \(C\) be a \(q\)-cyclotomic coset with \(|C|=n\ge2\). Proposition~\ref{prop:ed-binomial-correspondence} gives
\[
C\text{ is a complete equal-difference block}
\quad\Longleftrightarrow\quad
M_C(X)=X^n-a\quad(a\in k^\times).
\]
Every irreducible polynomial of degree \(n\ge2\) has both a leading term and a nonzero constant term, so such a binomial has the minimum possible support size.
This links an additive property of the exponent orbit to an optimal defining equation; the two complexity measures in Section~\ref{v79:sec-defining-equations} make this precise.
Starting from an arbitrary degree-\(n\) generator \(\alpha\), one may therefore ask whether a low-degree rational function \(R\in k(X)\), defined at \(\alpha\), sends its Frobenius orbit bijectively to an orbit whose exponent set is a complete equal-difference block.
The target modulus is the multiplicative order of \(R(\alpha)\), which may differ from that of \(\alpha\).
Proposition~\ref{v79:prop-generator-orbit-image} records the image orbit and any loss of degree; Exercise~\ref{v79:ex-binomial-defining-equation} gives an explicit example.
\end{remark}

\Needspace{5\baselineskip}
\begin{example}[Two levels of binomial descent for twenty-seventh roots of unity]\label{exm:med-binomial-twentyseven}
In Example~\ref{exm:med-twentyseven}, choose a primitive twenty-seventh root of unity \(\zeta\). All primitive roots form one eighteen-element orbit, so
\[
M_{C_1}(X)=\Phi_{27}(X)=X^{18}+X^9+1
\]
is irreducible over \(\F_2\). Put \(\omega=\zeta^9\). Then \(\omega\in\F_4\setminus\F_2\) and \(\omega^2+\omega+1=0\).
The two nine-element equal-difference blocks give
\[
M_{C_1}(X)=(X^9-\omega)(X^9-\omega^2)
\quad\text{in }\F_4[X].
\]
Both degree-nine polynomials are irreducible: on each root set, Frobenius relative to \(\F_4\), namely \(x\mapsto x^4\), has orbit length nine.
The composition over the base field is \(M_{C_1}(X)=G_9(X^9)\), where \(G_9(Y)=Y^2+Y+1\).

For the six three-element blocks, put \(\delta=\zeta^3\). Then \(\delta\) has order nine and \([\F_2(\delta):\F_2]=6\). The corresponding factorization is
\[
M_{C_1}(X)=\prod_{j=0}^{5}(X^3-\delta^{2^j})
\quad\text{in }\F_{64}[X].
\]
Every cubic is irreducible, and the composition over the base field is
\[
M_{C_1}(X)=G_3(X^3),\qquad G_3(Y)=Y^6+Y^3+1.
\]
Thus block sizes nine and three for the same coset correspond to outer degrees two and six and binomial degrees nine and three, respectively.
The eighteen singleton blocks give all linear factors over \(\F_{2^{18}}\), the splitting field of \(M_{C_1}\) over \(\F_2\).
The other two nonzero cosets in Example~\ref{exm:med-nonunit} give, respectively,
\[
M_{C_3}(X)=\Phi_9(X)=X^6+X^3+1,
\qquad M_{C_9}(X)=\Phi_3(X)=X^2+X+1.
\]
These expressions agree with the maximal block sizes three and one in the table.
\end{example}

\begin{corollary}[When every factor of a roots-of-unity polynomial is binomial]\label{cor:med-all-binomial}
Let \(N\ge2\), \(\gcd(N,q)=1\), and \(t=\ord_N(q)\). The following conditions are equivalent:
\begin{enumerate}
\item Every irreducible factor of \(X^N-1\) in \(k[X]\) is a binomial.
\item Every \(q\)-cyclotomic coset modulo \(N\) is a complete equal-difference block.
\item \(\rad(N)\mid q-1\), and if \(8\mid N\), then \(q\equiv1\pmod4\).
\item \(N\mid t(q-1)\).
\end{enumerate}
More generally, put \(\rho_N=\ord_{\rad(N)}(q)\), and define
\[
\omega_N=\ord_{b(N)}(q)=
\begin{cases}
2\rho_N,&8\mid N\text{ and }q^{\rho_N}\equiv3\pmod4,\\
\rho_N,&\text{otherwise}.
\end{cases}
\]
For any positive integer \(m\), all irreducible factors of \(X^N-1\) over \(\F_{q^m}\) are binomials if and only if \(\omega_N\mid\nobreak m\).
\end{corollary}
\begin{proof}
Applying Proposition~\ref{prop:ed-binomial-correspondence} to each orbit gives equivalence of the first two conditions.
If the second holds, applying Proposition~\ref{prop:ed-coset-criterion} to \(C_1\) gives the third.
Conversely, the reduced modulus \(d=N/\gcd(N,i)\) of any \(C_i\) divides \(N\), so the third condition implies the corresponding criterion for \(d\).
Hence every \(C_i\) is a complete equal-difference block.

Choose a primitive \(N\)-th root of unity \(\zeta\), which has degree \(t\) over \(k\).
If the first condition holds, its minimal polynomial has the form \(X^t-a\), where \(a=\zeta^t\in k^\times\).
Thus \(\zeta^{t(q-1)}=1\), or equivalently \(N\mid t(q-1)\).
Conversely, the fourth condition makes \(a=\zeta^t\) satisfy \(a^{q-1}=1\), so \(a\in k\).
The degree-\(t\) polynomial \(X^t-a\) has \(\zeta\) as a root and is therefore its minimal polynomial.
Proposition~\ref{prop:ed-binomial-correspondence} shows that \(C_1\) is a complete equal-difference block,
and Proposition~\ref{prop:ed-coset-criterion} then gives the third condition. By the equivalence already proved, the first condition follows as well.

Finally, apply the third condition with \(q^m\) in place of \(q\). The resulting congruence is \(q^m\equiv1\pmod{b(N)}\), or equivalently \(\omega_N\mid m\).
\end{proof}

When the equivalent conditions of Corollary~\ref{cor:med-all-binomial} hold, \(C_1=1+(N/t)\Z/N\Z\). Taking Frobenius conjugates of \(\zeta^j\) amounts to multiplying this exponent set by \(j\),
giving exactly \(t/\gcd(t,j)\) distinct exponents. Thus, for \(d\mid t\), the number of irreducible factors of degree \(d\) is
\[
\frac1d\,\#\{0\le j<N:\gcd(t,j)=t/d\}.
\]
Here the roots are first grouped by orbit length, then their number is divided by the length \(d\) of each orbit.

\begin{example}
For \(q=5,N=8\), we have \(t=2\) and \(8\mid2(5-1)\); indeed,
\[
X^8-1=(X-1)(X+1)(X-2)(X+2)(X^2-2)(X^2+2).
\]
In contrast, when \(q=2\) and \(N>1\) is odd, \(t\le\varphi(N)<N\), so \(N\mid t\) is impossible and at least one factor is not a binomial.
\end{example}

\begin{example}[The dyadic exception and two four-element equal-difference blocks]\label{exm:med-thirtytwo}
Take \(q=3,N=32,i=1\). Its cyclotomic coset is
\[
C_1=\{1,3,9,27,17,19,25,11\},\qquad \tau=8.
\]
Although \(\rho=\ord_2(3)=1\), we have \(8\mid32\) and \(3\equiv3\pmod4\), so \(\omega_C=2\) and \(\kappa_C=4\).
The coarsest decomposition is
\[
C_1=\{1,9,17,25\}\sqcup\{3,11,19,27\}
=(1+8\Z/32\Z)\sqcup(3+8\Z/32\Z).
\]
Multiplication by three interchanges the two blocks, and multiplication by nine traverses each:
\[
1\mapsto9\mapsto17\mapsto25\mapsto1,
\qquad3\mapsto27\mapsto19\mapsto11\mapsto3.
\]
The admissible block sizes are \(4,2,1\), with steps \(t=2,4,8\). The intermediate decomposition into blocks of size two is
\[
C_1=\{1,17\}\sqcup\{3,19\}\sqcup\{9,25\}\sqcup\{11,27\}.
\]
Each block has common difference sixteen, and \(3^4\equiv17\pmod{32}\) interchanges its two elements.
The cyclic gaps between consecutive elements in increasing order are not constant, so \(C_1\) itself is not a complete equal-difference block. This shows the necessity of the additional dyadic condition in \eqref{eq:med-minimal-step}.

In \(\F_3[X]\),
\[
\Phi_{32}(X)=X^{16}+1
=(X^8+X^4+2)(X^8+2X^4+2).
\]
Since \(\ord_{32}(3)=8\), both degree-eight factors are irreducible.
Choose \(\zeta\) to be a root of the first factor. Then \(M_{C_1}=X^8+X^4+2=G(X^4)\), where \(G(Y)=Y^2+Y+2\).
The element \(\gamma=\zeta^4\) belongs to \(\F_9\setminus\F_3\), and
\[
M_{C_1}(X)=(X^4-\gamma)(X^4-\gamma^3)
\quad\text{in }\F_9[X].
\]
Both quartics are irreducible. Equal-difference blocks, Frobenius step two, and binomial factorization over the quadratic field thus correspond in this example.
\end{example}

\begin{example}[From roots of order 80 to roots of order 160]\label{v125:ex-160-composition}
Take \(q=3\). The polynomial \(G(Y)=Y^4+Y+2\) is irreducible over \(\F_3\): it has no linear factor, and its remainders upon division by the three monic irreducible quadratics
\(Y^2+1,Y^2+Y+2,Y^2+2Y+2\) are \(Y,Y+1,Y+1\), respectively.
If \(\beta\) is a root, repeated squaring using \(\beta^4=2\beta+1\) gives
\[
\beta^{40}=-1,\qquad \beta^{16}=2\beta^3+\beta+2\ne1.
\]
Since \(\beta\in\F_{81}^\times\), these two identities exclude every proper divisor of 80 as its order. Thus \(\ord(\beta)=80\).
Choose \(\alpha^2=\beta\). Then \(\ord(\alpha)=160\), and \(\ord_{160}(3)=8\) shows that
\[
M(X)=G(X^2)=X^8+X^2+2
\]
is irreducible. Taking \(\alpha\) as the primitive 160th root, its root-exponent orbit is
\[
C_1=\{1,3,9,27,81,83,89,107\}
=\{1,3,9,27\}+\{0,80\}.
\]
Here \(b(160)=20\) and \(\ord_{20}(3)=4\), so \(\kappa_C=2\) and \(\omega_C=4\). The four two-element blocks give
\[
M(X)=\prod_{j=0}^{3}(X^2-\beta^{3^j})
\quad\text{in }\F_{81}[X],
\]
with every quadratic irreducible. Section~\ref{v44:sec-cyclotomic-local-global} will recover this \(C_1\) from its local orbits modulo 32 and five.
\end{example}

\subsection{Three-Sparse Factorizations and Images under Root Squaring}
\label{v49:subsec-three-sparse}

\leavevmode\phantomsection\label{v49:def-three-sparse}\ignorespaces
For a positive integer \(n\), call \(X^n-1\) \textbf{3-sparse} over \(k\) if each of its monic irreducible factors over \(k\) has at most three nonzero terms.
This concerns the number of terms of each irreducible factor, not of \(X^n-1\) itself.
The constant terms of these factors are nonzero, so 3-sparsity means that every irreducible factor is a binomial or trinomial.
Let \(q=p^e\), and write \(n=p^am\), where \(\gcd(m,p)=1\). Since
\[
 X^n-1=(X^m-1)^{p^a},
\]
3-sparsity depends only on \(m\), not on the multiplicities of the factors.

Multiple equal-difference structure gives a sufficient condition for a 3-sparse factorization: if the coarsest equal-difference decomposition has at most two blocks,
the outer polynomial in \eqref{eq:med-binomial-descent} has degree at most two and retains at most three nonzero terms after power substitution.
In odd characteristic, this arithmetic condition is also necessary~\cite{ChengOddSparse2025}.

\begin{theorem}[A criterion for 3-sparsity in odd characteristic]\label{v49:thm-odd-three-sparse}
Let \(q\) be an odd prime power, \(p=\operatorname{char}k\), and write \(n=p^am\), where \(\gcd(m,p)=1\). Then
\[
 \boxed{\quad X^n-1\text{ is 3-sparse over }k
 \quad\Longleftrightarrow\quad \rad(m)\mid q^2-1.\quad}
\]
\end{theorem}
\begin{proof}
We first prove sufficiency. Multiplicities do not change the form of the irreducible factors, so it suffices to consider \(X^m-1\).
The case \(m=1\) is immediate. Suppose \(\rad(m)\mid q^2-1\). For any \(q\)-cyclotomic coset \(C_i\) modulo \(m\),
put \(d=m/\gcd(m,i)\). If \(d=1\), its factor is linear.
If \(d>1\), then \(\rho=\ord_{\rad(d)}(q)\) divides two.
The value \(\omega_{C_i}\) given by Theorem~\ref{thm:med-minimal-step} also divides two:
the additional dyadic factor two can occur only when \(\rho=1\), since \(q^2\equiv1\pmod4\).
Taking \(r=\kappa_{C_i}\) in Theorem~\ref{prop:med-binomial-descent} gives
\[
 M_{C_i}(X)=G_{\kappa_{C_i}}(X^{\kappa_{C_i}}),
 \qquad \deg G_{\kappa_{C_i}}=\omega_{C_i}\le2.
\]
Thus every irreducible factor has at most three nonzero terms.

For necessity, suppose \(X^m-1\) is 3-sparse, but some prime \(\ell\mid m\) does not divide \(q^2-1\).
Since \(q\) is odd, so is \(\ell\). Put \(t=\ord_\ell(q)\). Then \(t\ge3\),
and all irreducible factors of \(\Phi_\ell\) have degree \(t\).
None of these factors is a binomial: otherwise Propositions~\ref{prop:ed-binomial-correspondence} and~\ref{prop:ed-coset-criterion}
would give \(\ell\mid q-1\), a contradiction.
Thus we may write
\[
 \Phi_\ell(X)=\prod_j\bigl(X^t+a_jX^{s_j}+b_j\bigr),
 \qquad a_jb_j\ne0,\quad 1\le s_j<t.
\]
Since the coefficient of \(X\) in \(\Phi_\ell=1+X+\cdots+X^{\ell-1}\) is nonzero,
we cannot have \(s_j>1\) for every \(j\). Hence there is an irreducible factor
\[
 g(X)=X^t+aX+b,\qquad ab\ne0.
\]
Let \(\xi\) be a root. It has order \(\ell\), as does \(\xi^2\),
so the minimal polynomial \(h\) of \(\xi^2\) also has degree \(t\).
Squaring each conjugate root of \(g\) gives
\[
 h(T^2)=\prod_{j=0}^{t-1}\bigl(T^2-\xi^{2q^j}\bigr)
       =(-1)^t g(T)g(-T).
\]
Expanding and replacing \(T^2\) by \(X\), we obtain
\[
 h(X)=
 \begin{cases}
 X^t+2bX^{t/2}-a^2X+b^2,&t\text{ even},\\[2pt]
 X^t+2aX^{(t+1)/2}+a^2X-b^2,&t\text{ odd}.
 \end{cases}
\]
Odd characteristic and \(ab\ne0\) ensure that all four coefficients are nonzero, while \(t\ge3\) ensures that the four exponents are distinct.
Thus \(h\) is an irreducible four-term factor of \(X^\ell-1\), contradicting 3-sparsity.
Consequently every \(\ell\mid m\) divides \(q^2-1\), or equivalently \(\rad(m)\mid q^2-1\).
\end{proof}

This criterion restricts the prime divisors of \(m\), not their exponents.
The degrees of irreducible factors can be arbitrarily large; their sparse form arises from substituting a power of \(X\) into a linear or quadratic outer polynomial, as determined by the equal-difference root blocks.
For example, over \(\F_3\), the 3-sparse polynomials \(X^n-1\) are precisely those with \(n=3^a2^b\), \(a,b\ge0\).
The two degree-eight trinomials in Example~\ref{exm:med-thirtytwo} illustrate this structure.
Over \(\F_9\), the corresponding lengths are precisely \(n=3^a2^b5^c\), \(a,b,c\ge0\).

\begin{example}[Squaring roots destroys the trinomial form]\label{v49:ex-square-four-terms}
The polynomial \(g=X^3-X-1\) has no root in \(\F_3\), so it is irreducible.
If \(\xi\) is a root, then
\[
 \xi^3=\xi+1,\qquad \xi^9=\xi-1,\qquad
 \xi^{13}=\xi(\xi+1)(\xi-1)=1.
\]
Since \(\xi\ne1\), it is a primitive thirteenth root of unity.
By the root-squaring formula above, the minimal polynomial of \(\xi^2\) is
\[
 h(X)=X^3+X^2+X-1.
\]
It is still an irreducible cubic, but has four nonzero terms.
Thus \(X^{13}-1\) is not 3-sparse over \(\F_3\), in agreement with \(13\nmid3^2-1\).
The complete conjugate orbits of these two factors, together with the other two factors, are given in Example~\ref{fas:ex-thirteen}.
\end{example}

\subsection{The Exceptional Family of Powers of Seven in Characteristic Two}

\begin{theorem}[Classification of 3-sparsity in characteristic two]\label{v49:thm-binary-three-sparse}
Let \(q=2^e\) and \(n=2^\lambda m\), where \(m\) is odd and \(\lambda\ge0\).
Every irreducible factor of \(X^n-1\) has at most three nonzero terms if and only if one of the following conditions holds:
\begin{enumerate}
\item \(\rad(m)\mid q^2-1\);
\item \(3\nmid e\), and we can write
\[
m=7^A s_0,\qquad A\ge1,\qquad \gcd(s_0,7)=1,
\]
with
\[
\rad(s_0)\mid q-1,\qquad
3\nmid\frac{s_0}{\gcd(s_0,q-1)},\qquad
\ord_{7^a}(q)=3\cdot7^{a-1}\quad(1\le a\le A).
\]
The final set of order conditions is equivalent to \(A=1\) or \(7\nmid e\).
\end{enumerate}
\end{theorem}
This classification is taken from \cite{ChengEvenSparse2026,ChengSparseCorrigendum2026}, with the corrected conditions for roots of unity of order a power of seven.
The power of the characteristic prime in the exponent changes only factor multiplicities; the exceptional conditions restrict both the remaining prime components and Frobenius periods at each power-of-seven level.

\Needspace{5\baselineskip}
\begin{corollary}[The binary and quaternary fields]\label{v49:cor-small-binary-sparse}
Write \(n=2^\lambda m\), where \(m\) is odd.
\begin{enumerate}
\item Over \(\F_2\), \(X^n-1\) is 3-sparse if and only if \(m=3^u\) or \(m=7^v\), where \(u,v\ge0\).
\item Over \(\F_4\), \(X^n-1\) is 3-sparse if and only if one of the three cases
\[
m=3^u5^v\quad(u,v\ge0),\qquad
m=7^A\quad(A\ge1),\qquad
m=3\cdot7^A\quad(A\ge1)
\]
holds.
\end{enumerate}
\end{corollary}
\begin{proof}
For \(q=2\), the main condition permits only the prime three, and \(s_0\) in the exceptional condition must equal one.
For \(q=4\), the main condition permits only the primes three and five; the exceptional condition requires \(s_0=3^u\),
and \(3\nmid s_0/\gcd(s_0,3)\) requires precisely \(u\le1\). Both base fields satisfy \(7\nmid e\), so \(A\) is unrestricted.
\end{proof}
These two classifications are due to Oliveira and Reis~\cite{OliveiraReis2021}.
For roots of unity of order a power of three over \(\F_2\), see Section~\ref{v49:sec-sparse-families}; the power-of-seven family has the following orbit interpretation.

\begin{example}[Orders 21 and 63 over the four-element field]\label{v104:ex-four-sparse-boundary}
Take \(\F_4=\F_2(\omega)\), where \(\omega^2+\omega+1=0\).
In the exceptional family \(m=7s_0\), the orders \(21\) and \(63\) correspond to \(s_0=3\) and \(s_0=9\), respectively.
Both satisfy all the other conditions, but \(s_0/\gcd(s_0,3)\) is \(1\) and \(3\), respectively. The shapes of the factors exhibit the difference.
Scaling the two cubic factors for seventh roots of unity gives
\[
X^{21}-1=(X^3-1)\prod_{j=0}^{2}
(X^3+\omega^{2j}X+1)(X^3+\omega^jX^2+1).
\]
Here \(X^3-1\) splits into three linear factors over \(\F_4\), and all six cubic factors are irreducible.
Indeed, \(X^3+X+1\) and \(X^3+X^2+1\) are irreducible over \(\F_2\); their degree three is coprime to the extension degree two, so they remain irreducible over \(\F_4\).
The other factors result from \(X\mapsto\omega^{-j}X\). Thus every irreducible factor here has at most three terms.

In contrast, put \(P=X^3+X^2+X+\omega\). Since
\[
P(0)=\omega,\quad P(1)=\omega^2,\quad
P(\omega)=P(\omega^2)=\omega,
\]
the polynomial \(P\) is irreducible. A root \(\theta\in\F_{64}\) has norm \(\theta^{21}=\omega\), so \(\ord(\theta)\) is nine or 63.
If it were nine, then \(\theta^3=\theta^{21}=\omega\), and \(P(\theta)=0\) would give \(\theta^2+\theta=0\), a contradiction.
Thus \(\theta\) has order 63, and \(P\) is a four-term irreducible factor of \(X^{63}-1\).
The example shows why the restriction \(3\nmid s_0/\gcd(s_0,q-1)\) cannot be omitted even when all the other conditions hold.
\end{example}

\begin{proposition}[The two maximal orbits of roots of unity of power-of-seven order]\label{v49:prop-seven-power-trinomials}
Let \(q=2^e\), \(a\ge1\), and \(r=7^{a-1}\). In \(\F_q[X]\), we always have
\[
\Phi_{7^a}(X)
=(X^{3r}+X^r+1)(X^{3r}+X^{2r}+1).
\]
The two displayed trinomials are both irreducible if and only if
\[
\gcd(3\cdot7^{a-1},e)=1.
\]
In particular, when \(3\nmid e\), this condition is equivalent to \(a=1\) or \(7\nmid e\).
\end{proposition}
\begin{proof}
Direct multiplication gives \(\Phi_7(X)=(X^3+X+1)(X^3+X^2+1)\).
For roots of unity of power-of-seven order, \(\Phi_{7^a}(X)=\Phi_7(X^{7^{a-1}})\), yielding the displayed identity.
We have \(\ord_7(2)=3\) and \(v_7(2^3-1)=1\). By Lemma~\ref{lem:valuation-lift},
\[
\ord_{7^a}(2)=3\cdot7^{a-1},\qquad
\ord_{7^a}(2^e)
=\frac{3\cdot7^{a-1}}{\gcd(3\cdot7^{a-1},e)}.
\]
Corollary~\ref{cor:cyclotomic-equal-degree} shows that every irreducible factor of \(\Phi_{7^a}\) has degree equal to the latter order.
Both displayed trinomials have degree \(3r\), so they are irreducible precisely when that order equals \(3r\).
\end{proof}

The sparse product identity holds over every finite field of characteristic two, while irreducibility depends on whether the orbits retain maximal length after base extension.
When \(7\mid e\) and \(a\ge2\), the original trinomials split further, and the number of terms in their irreducible factors may increase.

\begin{example}[A four-term factor after base extension]\label{v49:ex-binary-nonmaximal-seven}
Take \(q=128=2^7\) and \(m=49\). We have
\[
\ord_{49}(128)=\frac{21}{\gcd(21,7)}=3.
\]
Thus the irreducible factors of \(\Phi_{49}\) over \(\F_{128}\) are cubics, rather than the degree-twenty-one trinomials in the preceding proposition.
First put \(F(T)=T^7+T+1\). Its values at 0 and 1 are both one, and
\[
\begin{aligned}
 F&\equiv1\pmod{T^2+T+1},\\
 F&\equiv T\pmod{T^3+T+1},\qquad
 F\equiv T\pmod{T^3+T^2+1}.
\end{aligned}
\]
The displayed divisors are all the monic irreducible quadratics and cubics over \(\F_2\), and a reducible polynomial of degree seven has an irreducible factor of degree at most three.
Hence \(F\) is irreducible, so \(\F_{128}=\F_2(\tau)\), with \(\tau^7=\tau+1\), does give a field model.
Put
\[
 A=\tau^4+\tau^3+\tau^2+\tau,\qquad B=\tau^4+\tau,
 \qquad H(X)=X^3+AX^2+BX+1.
\]
Let \(R_j\) be the remainder of \(X^j\) modulo \(H\). The relation \(X^3\equiv AX^2+BX+1\) gives
\[
 R_0=1,\quad R_1=X,\quad R_2=X^2,\qquad
 R_{j+3}=AR_{j+2}+BR_{j+1}+R_j.
\]
Set
\[
\begin{aligned}
 u&=\tau^3+\tau,&v&=\tau^5+\tau^4+\tau^3+1,\\
 w&=\tau^6+\tau^5+\tau^4+\tau^3+1,&z&=\tau^6+\tau^4+\tau^2+\tau.
\end{aligned}
\]
Applying the recurrence up to \(R_7\), then reducing \(R_{14}\equiv R_7^2\) and \(R_{21}\equiv R_7R_{14}\) using \(\tau^7=\tau+1\), gives
\[
\begin{aligned}
 R_7&=uX^2+vX+w,\\
 R_{14}&=vX^2+BX+z,\\
 R_{21}&=uX^2+vX+w+1.
\end{aligned}
\]
Thus \(H\mid X^{21}+X^7+1\). In characteristic two,
\[
\begin{aligned}
 &(X^{21}+X^7+1)(X^{21}+X^{14}+1)\\
 &\qquad=X^{42}+X^{35}+X^{28}+X^{21}+X^{14}+X^7+1
 =\Phi_{49}(X).
\end{aligned}
\]
Consequently \(H\mid\Phi_{49}\); all its irreducible factors have degree three, so \(H\) itself is irreducible.
Both \(A\) and \(B\) are represented by nonzero polynomials in \(\tau\) of degree less than seven, and hence are nonzero.
Therefore \(H\) has four nonzero terms and \(X^{49}-1\) is not 3-sparse over this base field.
This example appears in \cite{ChengSparseCorrigendum2026}; it shows that \(3\nmid e\) alone cannot replace the maximal-orbit condition at power-of-seven levels.
\end{example}

\subsection{Cyclotomic Polynomials of Low Order and Coefficient Recurrences}
\label{v49:sec-low-order-cyclotomic}

When the square-free kernel of the index divides \(q+1\), roots of unity of low order can first be treated in the quadratic extension, where Frobenius acts by inversion.
When the dyadic part remains unchanged, coefficients of factors of higher-index cyclotomic polynomials can therefore be obtained recursively from coefficients of quadratic factors at lower indices.
For the stabilization that produces factors at higher indices from those at lower indices by power substitution, see \cite{WuZhuFengYang2017} and Proposition~\ref{v50:prop-cyclotomic-uniform-stabilization}.
The following quadratic coefficient recurrence is given in \cite{SehrawatSingh2025}; its proof uses multiplicative orders of roots and orbit lengths.

\begin{proposition}[Monomial lifting of quadratic factors of low order]
\label{v49:prop-low-order-cyclotomic}
Let \(n>1\), \(\rad(n)\mid q+1\), and put
\[
 m=\gcd(n,q+1)>2,\qquad t=n/m.
\]
Assume \(v_2(n)\le v_2(q+1)\). Then \(t\) is odd.
Choose a primitive \(m\)-th root of unity \(u\in\F_{q^2}\), put \(c=u+u^{-1}\in k\), and define recursively in \(k\)
\begin{equation}\label{v49:eq-coefficient-recurrence}
 c_0=2,\qquad c_1=c,\qquad c_{j+1}=c\,c_j-c_{j-1}\quad(j\ge1).
\end{equation}
Let \(J_m=\{j:1\le j<m/2,\ \gcd(j,m)=1\}\). Then
\begin{equation}\label{v49:eq-low-order-cyclotomic}
 \boxed{\quad
 \Phi_n(X)=\Phi_m(X^t)
 =\prod_{j\in J_m}\bigl(X^{2t}-c_jX^t+1\bigr).
 \quad}
\end{equation}
The right side has precisely \(\varphi(m)/2\) distinct monic irreducible factors, each of degree \(2t\).
\end{proposition}
\begin{proof}
Since \(m\mid q+1\), we have \(u^q=u^{-1}\), so \(c\in k\). Direct multiplication gives
\[
 (u+u^{-1})(u^j+u^{-j})-(u^{j-1}+u^{-(j-1)})
 =u^{j+1}+u^{-(j+1)}.
\]
Comparison of initial values and recurrences gives \(c_j=u^j+u^{-j}\). If \(j\) is coprime to \(m\), then \(u^j\) has order \(m>2\),
so \(u^j\ne u^{-j}\) and \(u^j\notin k\). Its Frobenius orbit is exactly \(\{u^j,u^{-j}\}\),
with minimal polynomial \(X^2-c_jX+1\). The set \(J_m\) chooses one representative from each pair of opposite unit indices, so
\[
 \Phi_m(X)=\prod_{j\in J_m}(X^2-c_jX+1).
\]

Since \(\rad(n)\mid q+1\), we have \(\rad(m)=\rad(n)\), and hence \(\rad(t)\mid m\).
If \(\beta^t\) has order \(m\), put \(d=\ord(\beta)\). Then \(d\mid mt\) and
\(d/\gcd(d,t)=m\). For every \(\ell\mid m\), the latter equality gives
\[
 v_\ell(d)-\min\{v_\ell(d),v_\ell(t)\}=v_\ell(m)>0,
\]
so \(v_\ell(d)=v_\ell(m)+v_\ell(t)=v_\ell(n)\). Every prime divisor of \(n\) divides \(m\), so \(d=n\).
Conversely, the \(t\)-th power of a primitive \(n\)-th root of unity has order \(m\). Thus the roots of \(\Phi_m(X^t)\) are precisely all primitive \(n\)-th roots of unity.
Since \(p\nmid n\), these roots are distinct; comparison of leading coefficients gives \(\Phi_n(X)=\Phi_m(X^t)\).

Finally, we calculate their orbit lengths. Since \(q\equiv-1\pmod m\) and \(m>2\), we have \(\ord_m(q)=2\),
so \(n\mid q^a-1\) requires \(a\) to be even. For every odd prime \(\ell\mid n\), put \(b=v_\ell(q+1)\).
Since \(\ell\nmid q-1\), Lemma~\ref{lem:valuation-lift} gives
\[
 v_\ell(q^{2s}-1)=b+v_\ell(s),\qquad
 v_\ell(t)=\max\{0,v_\ell(n)-b\}.
\]
Thus the divisibility conditions at odd primes are precisely \(t\mid s\). If \(n\) is even, its dyadic part divides \(q+1\) by hypothesis,
and therefore divides every \(q^{2s}-1\), imposing no further condition. Hence
\(\ord_n(q)=2t\). Every degree-\(2t\) factor in \eqref{v49:eq-low-order-cyclotomic}
has roots of order \(n\) and Frobenius period \(2t\), and is therefore irreducible. Distinct factors have root sets whose \(t\)-th power images lie in different two-element orbits, so they are distinct.
\end{proof}

Equation~\eqref{v49:eq-coefficient-recurrence} is also the Dickson recurrence of Section~\ref{sec:dickson-polynomials}:
with \(D_j\) as in Definition~\ref{def:dickson-polynomial}, we have \(c_j=D_j(c)\).
The recurrence calculates trace parameters of reciprocal root pairs, while the multiple equal-difference decomposition explains why these quadratics lift by monomial composition.
Once the initial value \(c\) is chosen as specified, all coefficient calculations take place in the base field.

\leavevmode\phantomsection\label{v49:criterion-trace-order}\ignorespaces The initial value must also record the exact multiplicative order of the roots. In odd characteristic, the condition that \(c^2-4\) is a nonsquare says only that the roots of \(X^2-cX+1\) lie in the quadratic extension and are Frobenius conjugates;
the two roots have the same order, which divides \(q+1\) and may vary with \(c\). For any nonzero root \(u\),
\[
 c_r=2\quad\Longleftrightarrow\quad
 (u^r-1)^2=0\quad\Longleftrightarrow\quad u^r=1.
\]
Consequently, assuming \(X^2-cX+1\) is irreducible, its roots have order exactly \(m\) if and only if
\[
 c_m=2,\qquad c_{m/\ell}\ne2\quad\text{for every prime }\ell\mid m.
\]
This criterion distinguishes traces of primitive \(m\)-th roots of unity from traces of roots of smaller order. The coefficients \(c_j\) may vanish, so \eqref{v49:eq-low-order-cyclotomic} ensures at most three terms, not necessarily exactly three.

\begin{example}[From trace parameters of order seven to factors of order forty-nine]
\label{v49:ex-low-order-49}
Take \(q=13\) and \(n=49\). Then \(m=7\) and \(t=7\). Choose \(c=7\) in \(\F_{13}\). We have
\(c^2-4=6\), a nonsquare, so \(X^2-7X+1\) is irreducible.
Equation~\eqref{v49:eq-coefficient-recurrence} gives
\[
 (c_0,c_1,\ldots,c_7)=(2,7,8,10,10,8,7,2).
\]
Since \(c_7=2\) and \(c_1\ne2\), the roots of this quadratic have order seven. Thus
\[
 \Phi_7(X)=(X^2-7X+1)(X^2-8X+1)(X^2-10X+1),
\]
and
\[
 \Phi_{49}(X)
 =(X^{14}-7X^7+1)(X^{14}-8X^7+1)(X^{14}-10X^7+1).
\]
All three factors are irreducible because \(\ord_{49}(13)=14\). Expanding their product gives
\[
 \Phi_{49}(X)=X^{42}+X^{35}+X^{28}+X^{21}+X^{14}+X^7+1,
\]
in agreement with reduction of the integer cyclotomic polynomial.

The dyadic condition cannot be omitted. Again take \(q=13\), but now let \(n=28\). Then \(m=14\), \(t=2\),
and \(v_2(28)>v_2(14)\). The identity \(\Phi_{28}(X)=\Phi_{14}(X^2)\) still holds, but its composed factors may split; for example,
\[
 X^4-3X^2+1=(X^2-X-1)(X^2+X-1)
 \quad\text{in }\F_{13}[X].
\]
Both quadratic factors have discriminant five, a nonsquare, and are therefore irreducible.
Thus the monomial lifting identity for a cyclotomic polynomial of low order and irreducibility of its lifted factors are separate matters to check.
\end{example}

\subsection{Unions of Orbits, Root Power Sums, and the Maximum Composition Exponent}
\label{v40:sec:root-power-composition}

Transitivity of the multiplicative action on a single cyclotomic coset ensures that each complete equal-difference block extends to a partition of the whole coset.
For a union of cosets, existing blocks may merge and produce translation periods absent from each constituent coset.
The maximum monomial composition exponent of the entire root set is still determined by its additive stabilizer; root power sums give another way to compute this stabilizer.

\begin{proposition}[Composition exponents and root power sums for unions of orbits]
\label{v40:prop:root-power-composition}
Let \(p\nmid N\), let \(\zeta\) be a primitive \(N\)-th root of unity, and let \(\varnothing\ne T\subseteq\Z/N\Z\) satisfy \(qT=T\). Put
\[
 P_T(X)=\prod_{a\in T}(X-\zeta^a)\in k[X],
 \qquad \kappa_T=|H_+(T)|.
\]
Then, for every positive integer \(r\),
\[
 P_T\in k[X^r]\quad\Longleftrightarrow\quad r\mid\kappa_T.
\]
Thus \(\kappa_T\) is the maximum monomial composition exponent of \(P_T\). Put
\[
 s_j=\sum_{a\in T}\zeta^{aj}.
\]
Then
\begin{equation}\label{v40:eq:root-power-composition}
\boxed{\quad
 \kappa_T
 =\gcd\{j\ge1:[X^j]P_T\ne0\}
 =\gcd\bigl(N,\{j:0\le j<N,\ s_j\ne0\}\bigr).
\quad}
\end{equation}
We call \(s_j\) a \textbf{power sum of roots} of this root set. Each \(s_j\) belongs to \(k\), and the set of frequencies with \(s_j\ne0\) is nonempty.
If \(T\) is a single Frobenius orbit of \(\alpha=\zeta^a\), of length \(\tau\), then
\[
 s_j=\sum_{i=0}^{\tau-1}\alpha^{jq^i}
 =\Tr_{\F_{q^\tau}/k}(\alpha^j).
\]
The notation on the right denotes the sum of the \(\tau\) Frobenius conjugate values; the general definition and properties of trace are given in Chapter~\ref{chap:trace}.
\end{proposition}
\begin{proof}
Write \(H_+(T)=(N/\kappa_T)\Z/N\Z\). The set \(T\) is a disjoint union of cosets of this subgroup.
Applying Proposition~\ref{prop:ed-binomial-correspondence} to each coset gives
\(P_T(X)=Q(X^{\kappa_T})\). Since \(P_T\in k[X]\), extracting coefficients by exponent gives \(Q\in k[X]\), so \(P_T\in k[X^r]\) whenever \(r\mid\kappa_T\).

Conversely, if \(P_T=G(X^r)\), the absence of repeated roots in \(P_T\) implies \(p\nmid r\).
For every root \(\alpha\) and \(\varepsilon\in\mu_r\), the element \(\varepsilon\alpha\) is again a root; both lie in \(\mu_N\), so \(\mu_r\subseteq\mu_N\) and \(r\mid N\).
Thus \(T+(N/r)\Z/N\Z=T\), giving \((N/r)\Z/N\Z\subseteq H_+(T)\) and \(r\mid\kappa_T\).
On the other hand, \(P_T\in k[X^r]\) means precisely that \(r\) divides every positive exponent with nonzero coefficient, giving the first equality in \eqref{v40:eq:root-power-composition}.

To prove the second equality, in a finite extension \(E/k\) containing \(\zeta\), define for a function \(v:\Z/N\Z\to E\)
\[
 \widehat v(j)=\sum_{a\bmod N}v(a)\zeta^{aj}.
\]
The orthogonality relation for roots of unity is
\[
 \sum_{j=0}^{N-1}\zeta^{j(a-b)}
 =\begin{cases}N,&a=b\pmod N,\\0,&a\ne b\pmod N.\end{cases}
\]
The second case follows by summing a geometric progression. Since \(p\nmid N\), the element \(N\) is invertible in \(E\), and hence
\[
 v(a)=N^{-1}\sum_{j=0}^{N-1}\widehat v(j)\zeta^{-aj}.
\]
Taking \(v=1_T\) gives \(\widehat v(j)=s_j\). The inverse transform shows that these \(s_j\) are not all zero, and that
\[
 T+u=T
 \quad\Longleftrightarrow\quad
 \zeta^{ju}s_j=s_j\quad\text{for every }j.
\]
Put \(d=\gcd(N,\{j:s_j\ne0\})\). The right-hand condition is equivalent to \(N\mid ju\) for every frequency with nonzero coefficient; by B\'ezout's identity, this is equivalent to \(N\mid du\), or \(N/d\mid u\).
Thus \(H_+(T)=(N/d)\Z/N\Z\) and \(\kappa_T=d\).
Finally, \(qT=T\) gives \(s_j^q=s_j\), so \(s_j\in k\). For a single orbit, the formula is the trace expression for the power sum of conjugate roots.
\end{proof}

\leavevmode\phantomsection\label{v40:def:cyclic-fourier}\ignorespaces
This transform is called the \textbf{finite-field-valued Fourier transform on a cyclic group}.
Its values lie in a finite extension containing the roots of unity, and its inverse requires \(p\nmid N\); the additive Fourier transform in Chapter~\ref{chap:characters} uses complex-valued additive characters.
Both recover functions through character orthogonality, but their value fields and underlying groups differ.
In particular, \(s_0=|T|\cdot1_k\) may vanish because of the characteristic. The zero frequency imposes no restriction on translations, so its vanishing does not affect the stabilizer criterion.

\Needspace{5\baselineskip}
\begin{example}[A union produces a new composition exponent]
\label{v40:ex:union-emergent-stabilizer}
Take \(q=2\) and \(N=15\). Both cyclotomic cosets \(C_0=\{0\}\) and \(C_5=\{5,10\}\) have trivial additive stabilizers: their orders divide \(\gcd(15,1)\) and \(\gcd(15,2)\), respectively.
Their union \(T=\{0,5,10\}=5\Z/15\Z\), however, has an additive stabilizer of order three and gives
\[
 P_T(X)=(X-1)(X^2+X+1)=X^3-1.
\]
If \(\eta=\zeta^5\), then \(s_j=1+\eta^j+\eta^{2j}\) is one when \(3\mid j\) and zero otherwise; the frequencies with nonzero coefficients are \(0,3,6,9,12\), whose greatest common divisor is three.
Thus the stabilizer of the entire root set is not obtained merely by intersecting the stabilizers of the irreducible factors: translations may interchange roots of different factors without preserving each part separately.
\end{example}

\section{Base Extension and Tensor Product Decompositions}
\label{sec:tensor-fields}

Extending the base field replaces \(\sigma\) by an iterate and may therefore split one orbit into several.
The tensor product gives an algebraic expression of this phenomenon. For finite extensions \(E/k\) and \(F/k\),
\(E\otimes_kF\) is the vector space generated by symbols \(x\otimes y\) subject to bilinear relations, with multiplication defined by \((x\otimes y)(x'\otimes y')=xx'\otimes yy'\). If \((e_i)\) and \((f_j)\) are bases of the two fields,
then \((e_i\otimes f_j)\) is a basis of the tensor product, whose dimension is therefore \([E:k][F:k]\).

\begin{example}[Compute the tensor product of two four-element fields first]\label{v103:ex-tensor-four}
Let \(E=\F_4=\F_2(\alpha)\), where \(\alpha^2+\alpha+1=0\). Then
\[
E\otimes_{\F_2}E\simeq E[T]/\bigl((T-\alpha)(T-\alpha^2)\bigr)\simeq E\times E.
\]
In the middle quotient ring, \(e_0=T+\alpha^2\) and \(e_1=T+\alpha\) satisfy \(e_0+e_1=1\), \(e_0e_1=0\), and \(e_i^2=e_i\), and correspond to \((1,0)\) and \((0,1)\), respectively.
Both arise from the splitting of the same root orbit after base extension.
The first isomorphism sends \(1\otimes\alpha\) to \(T\) and can be checked directly on the \(E\)-bases \((1\otimes1,1\otimes\alpha)\) and \((1,T)\).
The second sends \(P(T)\) to \((P(\alpha),P(\alpha^2))\), so their composite is
\[
x\otimes y\longmapsto(xy,xy^2).
\]
In particular, \(e_0(\alpha)=1,e_0(\alpha^2)=0\), while \(e_1\) has the opposite values. The two components and their projections have thus been computed explicitly.
\end{example}

\begin{theorem}\leavevmode\phantomsection\label{thm:tensor-fields}
In a common extension, put \(E_a=\F_{q^a}\), \(E_b=\F_{q^b}\), \(g=\gcd(a,b)\), and \(\ell=\lcm(a,b)\). There is an \(E_a\)-algebra isomorphism
\[
E_a\otimes_kE_b\simeq\prod_{j=0}^{g-1}E_\ell,
\qquad x\otimes y\longmapsto (x\,y^{q^j})_{0\leq j<g}.
\]
In particular, the multiplication map \(E_a\otimes_kE_b\to E_\ell\) is an isomorphism if and only if \(a,b\) are coprime.
\end{theorem}
\begin{proof}
Write \(E_b=k(\beta)\), and let \(m_\beta\in k[T]\) be its minimal polynomial. The map \(x\otimes P(\beta)\mapsto xP(T)\) induces
\[
E_a\otimes_kE_b\simeq E_a[T]/(m_\beta).
\]
This can be checked directly using the respective \(E_a\)-bases \(1\otimes1,1\otimes\beta,\ldots,1\otimes\beta^{b-1}\) and \(1,T,\ldots,T^{b-1}\).
The roots of \(m_\beta\) are \(\beta^{q^i}\), \(i\in\Z/b\Z\). Frobenius relative to \(E_a\) sends the index \(i\) to \(i+a\); its orbits are precisely the residue classes modulo \(g\), giving \(g\) orbits, each of length \(b/g\).
Thus \(m_\beta\) splits in \(E_a[T]\) into \(g\) distinct irreducible factors.
Evaluation identifies the quotient by the \(j\)-th factor with the field generated by \(E_a\) and \(\beta^{q^j}\), namely the compositum \(E_\ell\).
The Chinese remainder theorem gives the displayed direct-product isomorphism, whose \(j\)-th component is evaluation at \(T=\beta^{q^j}\).

The multiplication map is the zeroth component. It is always surjective; its domain has dimension \(ab\), while its target has dimension \(\ell=ab/g\), so it is injective if and only if \(g=1\).
\end{proof}

\leavevmode\phantomsection\label{def:linearly-disjoint}\ignorespaces If the multiplication map above is injective, the two extensions are called \textbf{linearly disjoint} over \(k\).
Equivalently, a \(k\)-basis of one field remains linearly independent over the other. For finite fields, this is determined entirely by coprimality of their degrees.
The compositum corresponds to just one component of the tensor product; the tensor product also records the other ways in which the two fields can agree on their common subfield.

\subsection*{Direct Products, Tensor Products, and Changes of Orbit Step}

Theorem~\ref{thm:finite-orbit-duality} expresses tensor product decomposition as a general operation on orbits.
We compare evaluation points in a sufficiently large common finite extension \(L\).

\begin{proposition}\leavevmode\phantomsection\label{prop:orbit-algebra-operations}
For finite reduced \(k\)-algebras \(A,B\), there are natural Frobenius-equivariant bijections
\[
S_L(A\times B)\simeq S_L(A)\sqcup S_L(B),\qquad
S_L(A\otimes_k B)\simeq S_L(A)\times S_L(B).
\]
The tensor product is again reduced. If the base field is extended to \(k_r=\F_{q^r}\subseteq L\), the \(k_r\)-algebra homomorphisms \(A\otimes_k k_r\to L\) can be identified with \(S_L(A)\), while relative Frobenius becomes \(\tau^r\).
\end{proposition}
\begin{proof}
A homomorphism from a direct product to a field selects exactly one component, giving the first bijection.
The second sends \(\iota:A\otimes B\to L\) to the maps \(a\mapsto\iota(a\otimes1)\) and \(b\mapsto\iota(1\otimes b)\).
Its inverse sends \((u,v)\) to \(a\otimes b\mapsto u(a)v(b)\); the bilinear relations ensure that this map is well defined.
Tensor products distribute over finite direct products, and Theorem~\ref{thm:tensor-fields} expresses the tensor product of each pair of field factors as a direct product of finite fields. Hence the tensor product is reduced.

For base extension, fix \(k_r\subseteq L\). Every \(k\)-homomorphism \(u:A\to L\) extends uniquely to a \(k_r\)-homomorphism \(A\otimes k_r\to L\), given by \(a\otimes c\mapsto u(a)c\).
Relative Frobenius for the new base field is \(x\mapsto x^{q^r}\), corresponding to \(\tau^r\) on evaluation points.
\end{proof}

For example, the direct product of a cycle of length \(a\) and a cycle of length \(b\)
has \(\gcd(a,b)\) orbits under the simultaneous cyclic action, each of length \(\lcm(a,b)\).
This explains both the number of components of a finite-field tensor product and the degree of each component.
Changing the step of a cycle of length \(d\) to \(r\) gives the factorization law in Theorem~\ref{thm:factor-base-extension}.

\begin{figure}[!htbp]
\centering
\begin{tikzpicture}[
  point/.style={draw=Jade,fill=CoverMint,rounded corners=2pt,minimum width=12mm,minimum height=7mm,inner sep=2pt},
  orbit/.style={-{Stealth[length=2mm]},Jade,semithick},
  every node/.style={font=\small}]
\node[point] (p00) at (0,0) {\((0,0)\)};
\node[point] (p11) at (2.5,0) {\((1,1)\)};
\node[point] (p02) at (2.5,-1.5) {\((0,2)\)};
\node[point] (p13) at (0,-1.5) {\((1,3)\)};
\draw[orbit] (p00)--(p11);
\draw[orbit] (p11)--(p02);
\draw[orbit] (p02)--(p13);
\draw[orbit] (p13)--(p00);
\node at (1.25,-.75) {\(\mathcal O_0\)};
\node[point] (p10) at (6,0) {\((1,0)\)};
\node[point] (p01) at (8.5,0) {\((0,1)\)};
\node[point] (p12) at (8.5,-1.5) {\((1,2)\)};
\node[point] (p03) at (6,-1.5) {\((0,3)\)};
\draw[orbit] (p10)--(p01);
\draw[orbit] (p01)--(p12);
\draw[orbit] (p12)--(p03);
\draw[orbit] (p03)--(p10);
\node at (7.25,-.75) {\(\mathcal O_1\)};
\draw[-{Stealth[length=2mm]},Clay] (1.25,-1.95)--(1.25,-2.8)
  node[midway,right=3pt] {Equivariant functions};
\draw[-{Stealth[length=2mm]},Clay] (7.25,-1.95)--(7.25,-2.8)
  node[midway,right=3pt] {Equivariant functions};
\node at (1.25,-3.15) {\(\mathcal A(\mathcal O_0,\tau)\simeq\F_{q^4}\)};
\node at (7.25,-3.15) {\(\mathcal A(\mathcal O_1,\tau)\simeq\F_{q^4}\)};
\node at (4.25,-3.15) {\(\times\)};
\end{tikzpicture}
\caption{Simultaneous Frobenius action on the evaluation points has two four-element orbits, each corresponding to a field component of the tensor product.}
\label{fig:tensor-evaluation}
\end{figure}
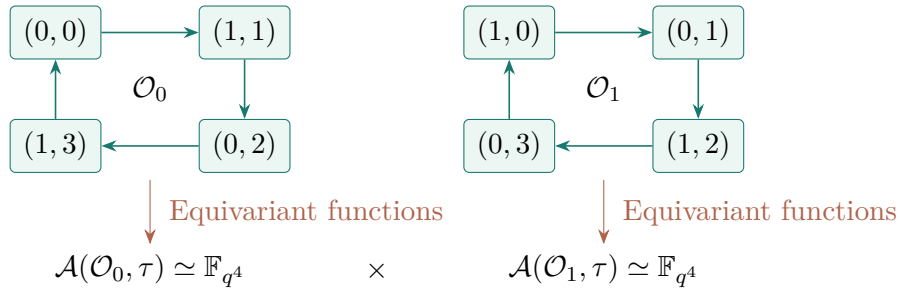

\begin{example}[Evaluation-point orbits and tensor product components]\leavevmode\phantomsection\label{exm:tensor-evaluation-orbits}
Let \(K=\F_{q^2}\subseteq M=\F_{q^4}\), and choose \(L=M\) as a common splitting field.
Write the two \(k\)-embeddings of \(K\) and the four \(k\)-embeddings of \(M\) as \(a\mapsto a^{q^i}\) and \(b\mapsto b^{q^j}\), respectively. Thus
\[
S_L(K\otimes_kM)=\{(i,j):i\in\Z/2\Z,\ j\in\Z/4\Z\},
\]
where \((i,j)\) corresponds to the evaluation homomorphism \(a\otimes b\mapsto a^{q^i}b^{q^j}\).
Frobenius advances both coordinates by one: \(\tau(i,j)=(i+1,j+1)\).
The two cycles in Figure~\ref{fig:tensor-evaluation} are all its orbits.
On either four-element orbit, an equivariant function is determined by its value at a representative; one full circuit requires that value to satisfy \(z^{q^4}=z\),
so each orbit contributes one \(M\)-component.

Choosing representatives \((0,0)\) and \((1,0)\), evaluation gives the explicit \(M\)-algebra isomorphism
\[
K\otimes_k M\longrightarrow M\times M,
\qquad a\otimes b\longmapsto (ab,a^q b).
\]
Here \(M\) acts through the second tensor factor, and each direct-product component records the equivariant functions on one orbit of evaluation points.
\end{example}

\section{Synchronization of Local Orbits and Explicit Factorization of Roots-of-Unity Polynomials}
\label{v44:sec-cyclotomic-local-global}

Conjugate orbits of roots of unity can first be calculated at pairwise coprime moduli and then combined by the Chinese remainder theorem.
In the combination, the local orbits advance simultaneously; even when all local orbits are specified, they may yield more than one global orbit.
The distinctions between these orbits are determined by the relative phases of their components.
We give nonredundant phase representatives and apply them to cyclotomic cosets for odd moduli.
This local-to-global description corresponds to the explicit construction of representatives for odd moduli in \cite{ZLW2026}.

\subsection{Simultaneous Action and the Phase Quotient}

Let \(s\ge1\), let \(m_1,\ldots,m_s\ge2\) be pairwise coprime, put \(N=m_1\cdots m_s\), and assume \(\gcd(N,q)=1\).
The Chinese remainder theorem for integers gives
\[
 \Z/N\Z\simeq\prod_{i=1}^s\Z/m_i\Z,
\]
and multiplication by \(q\) modulo \(N\) corresponds to simultaneous multiplication by \(q\) in every component.
Fix an orbit in the \(i\)-th component,
\[
 C_i=\{a_iq^t\bmod m_i:t\ge0\},\qquad |C_i|=\tau_i.
\]
The map \(j_i\mapsto a_iq^{j_i}\) is a bijection from \(\Z/\tau_i\Z\) to \(C_i\).
\leavevmode\phantomsection\label{v44:def-local-phase}\ignorespaces
We call \(j_i\) the \textbf{phase} relative to the representative \(a_i\).

\begin{theorem}[Synchronized decomposition of local orbits]\label{v44:thm-local-orbit-synchronization}
Put \(L=\lcm(\tau_1,\ldots,\tau_s)\). Under the CRT identification above,
every global \(q\)-orbit in \(C_1\times\cdots\times C_s\) has length \(L\),
and the global orbits are indexed bijectively by the \textbf{phase quotient}
\begin{equation}\label{v44:eq-phase-quotient}
 \mathcal P=\left(\prod_{i=1}^s\Z/\tau_i\Z\right)
 \big/\langle(1,\ldots,1)\rangle.
\end{equation}
Thus the number of global orbits is
\[
 \frac{\tau_1\cdots\tau_s}{\lcm(\tau_1,\ldots,\tau_s)}.
\]
\end{theorem}
\begin{proof}
In phase coordinates, multiplication by \(q\) sends \((j_1,\ldots,j_s)\) to
\((j_1+1,\ldots,j_s+1)\). After \(t\) iterations, it returns to its starting point if and only if
\(\tau_i\mid t\) for every \(i\). Hence the least positive return time is \(L\).
Two phase tuples belong to the same orbit if and only if their difference belongs to
\(\langle(1,\ldots,1)\rangle\), giving \eqref{v44:eq-phase-quotient}.
This subgroup has \(L\) elements, while the phase space has \(\tau_1\cdots\tau_s\), giving the orbit count.
\end{proof}

For \(s=2\), there are \(\gcd(\tau_1,\tau_2)\) orbits, each of length \(\lcm(\tau_1,\tau_2)\).
This is the synchronized decomposition of two Frobenius cycles from Section~\ref{sec:tensor-fields}.
If the lengths are pairwise coprime, the phase quotient has one element, and the chosen local orbits combine into one global orbit.
In general, phase representatives must also be chosen.

\begin{proposition}[Recursive phase representatives and CRT reconstruction]\label{v44:prop-recursive-phase-representatives}
Put
\[
 L_0=1,\qquad L_i=\lcm(\tau_1,\ldots,\tau_i),\qquad
 \delta_i=\gcd(L_{i-1},\tau_i)\quad(1\le i\le s).
\]
Each class of the phase quotient~\eqref{v44:eq-phase-quotient} has a unique representative satisfying
\begin{equation}\label{v44:eq-phase-box}
 j_1=0,\qquad 0\le j_i<\delta_i\quad(2\le i\le s).
\end{equation}
More explicitly, put \(N_i=N/m_i\), choose \(c_i\) with \(c_iN_i\equiv1\pmod{m_i}\), and set \(E_i=c_iN_i\).
Then the elements
\begin{equation}\label{v44:eq-crt-orbit-representatives}
 \gamma(\mathbf j)
 =\sum_{i=1}^s E_i a_iq^{j_i}\pmod N,
 \qquad \mathbf j\text{ satisfying \eqref{v44:eq-phase-box}},
\end{equation}
represent all global orbits in \(C_1\times\cdots\times C_s\), exactly once each.
\end{proposition}
\begin{proof}
We first prove existence. A common phase shift makes the first coordinate zero.
Suppose the first \(i-1\) coordinates already lie in their prescribed ranges.
The common shifts preserving these coordinates are precisely the multiples of \(L_{i-1}\);
the subgroup they generate in \(\Z/\tau_i\Z\) is \(\delta_i\Z/\tau_i\Z\).
Thus the \(i\)-th coordinate can be reduced uniquely to its representative modulo \(\delta_i\)
in \(0,\ldots,\delta_i-1\), without changing the previously reduced coordinates.
Continuing gives a representative in \eqref{v44:eq-phase-box}.

For uniqueness, suppose two such representatives \(\mathbf j,\mathbf j'\) differ by a common shift \(t\).
Their first coordinates are both zero, so \(\tau_1\mid t\).
If the first \(i-1\) coordinates agree, then \(L_{i-1}\mid t\), and
\[
 j_i'-j_i\equiv t\pmod{\tau_i}.
\]
Thus \(\delta_i\mid j_i'-j_i\). Both coordinates lie in \([0,\delta_i-1]\), so they agree, and \(\tau_i\mid t\).
Inductively, all coordinates agree and \(L\mid t\).
Finally, \eqref{v44:eq-crt-orbit-representatives} is exactly the CRT preimage of the specified local phases,
so existence and uniqueness of phase representatives give all the assertions.
\end{proof}

Given local representatives \(a_i\) and lengths \(\tau_i\), the phase box directly outputs \(P=\prod_i\tau_i/L\) global representatives. Precompute the Chinese remainder coefficients and \(a_iq^j\) for \(0\le j<\delta_i\). Each representative then uses \(s\) modular operations; including successive computation of these powers, the count is \(O(sP+\sum_i\delta_i)\). This excludes obtaining local representatives, orders, and prime factorizations, as well as expanding root-polynomial coefficients. If every orbit element is requested, the output already contains \(PL=\prod_i\tau_i\) entries. With ordinary integer multiplication and division, each operation modulo \(N\) uses at most \(O((\log N)^2)\) bit operations; the number of explicit representatives must still be charged to the output size.

The number of phase representatives can also be checked directly from
\[
 L_i=\frac{L_{i-1}\tau_i}{\delta_i}
 \quad\Longrightarrow\quad
 \prod_{i=1}^s\delta_i=\frac{\tau_1\cdots\tau_s}{L}.
\]
A CRT orbit representative need not be the coset leader of Definition~\ref{v50:def-coset-leader}.
If indexing by coset leaders is desired, Proposition~\ref{v50:prop-med-leader-compression} gives the least representative of each resulting orbit.

\subsection{Local Representatives at Odd Prime Powers}

Let \(\ell\) be an odd prime with \(\gcd(\ell,q)=1\).
Cyclicity of the unit group modulo \(\ell^e\) was established in the proof of Proposition~\ref{prop:cyclotomic-irreducible}.
Nonunit indices are separated by \(\ell\)-adic valuation; each layer reduces to the unit group for a smaller modulus.

\begin{proposition}[Local lengths and valuation layers]\leavevmode\phantomsection\label{v44:prop-prime-power-local-orbits}\ignorespaces
Put
\[
 h=\ord_\ell(q),\qquad b=v_\ell(q^h-1).
\]
For every \(a\ge1\),
\begin{equation}\label{v44:eq-prime-power-order}
 \ord_{\ell^a}(q)=h\ell^{\max(0,a-b)}.
\end{equation}
Fix \(e\ge1\), and choose a primitive root \(g\) modulo \(\ell^e\).
For \(0\le y<e\), put
\[
 f_y=\ord_{\ell^{e-y}}(q),\qquad
 G_y=\frac{\varphi(\ell^{e-y})}{f_y}.
\]
Every orbit of indices with valuation exactly \(y\) has length \(f_y\), and a nonredundant set of representatives is
\begin{equation}\label{v44:eq-prime-power-local-representatives}
 \ell^yg^j\pmod{\ell^e},\qquad 0\le j<G_y.
\end{equation}
Adjoining the singleton orbit of the zero index gives all \(q\)-cyclotomic cosets modulo \(\ell^e\).
\end{proposition}
\begin{proof}
A return modulo \(\ell\) requires the exponent to have the form \(hu\). Applying Lemma~\ref{lem:valuation-lift} to \(Q=q^h\) gives
\[
 v_\ell(q^{hu}-1)=b+v_\ell(u).
\]
Thus the least such exponent returning modulo \(\ell^a\) is \(h\ell^{\max(0,a-b)}\), proving \eqref{v44:eq-prime-power-order}.

Write an index of valuation \(y\) as \(\ell^yu\), where \(u\) is a unit modulo \(\ell^{e-y}\).
The map \(u\mapsto\ell^yu\) is an equivariant bijection from that unit group to this valuation layer,
and \(\ell^yuq^t\equiv\ell^yu\pmod{\ell^e}\) if and only if \(q^t\equiv1\pmod{\ell^{e-y}}\).
Hence every orbit in the layer has length \(f_y\).

The image of the primitive root \(g\) modulo \(\ell^{e-y}\) still generates the unit group.
A cyclic group has a unique subgroup of order \(f_y\), so
\[
 \langle q\rangle=\langle g^{G_y}\rangle
 \quad\text{in }(\Z/\ell^{e-y}\Z)^\times.
\]
Its multiplicative cosets have representatives \(g^j\), \(0\le j<G_y\).
Multiplication by \(\ell^y\) gives \eqref{v44:eq-prime-power-local-representatives}.
Distinct valuation layers are disjoint, and zero forms its own orbit, so every orbit has been listed.
\end{proof}

For example, when \(q=2\) and \(\ell=3\), we have \(h=2\) and \(b=v_3(2^2-1)=1\), so
\[
 \ord_{3^a}(2)=2\cdot3^{a-1}\qquad(a\ge1).
\]
Each nonzero valuation layer has one orbit.
The three nonzero orbits modulo twenty-seven, of lengths eighteen, six, and two, can therefore also be obtained by this local calculation,
in agreement with Example~\ref{exm:med-nonunit}.

\par\Needspace{12\baselineskip}
\subsection{From Local Representatives to Factors of Roots-of-Unity Polynomials}

\begin{corollary}[Explicit representatives and irreducible factors]\label{v44:cor-crt-cyclotomic-factorization}
Let \(N>1\), \(N=m_1\cdots m_s\), where the \(m_i\ge2\) are pairwise coprime, and suppose \(\gcd(N,q)=1\).
Choose one representative from each orbit modulo \(m_i\), forming a set \(\mathcal R_i\).
For each tuple \(a_i\in\mathcal R_i\), use its local lengths to define the phase representatives in \eqref{v44:eq-phase-box},
then obtain \(\gamma\) from \eqref{v44:eq-crt-orbit-representatives}.
These \(\gamma\) form a nonredundant set \(\mathcal R\) of representatives of all orbits modulo \(N\).
If \(q\) is a prime power and \(\zeta\) is a primitive \(N\)-th root of unity, then
\begin{equation}\label{v44:eq-crt-cyclotomic-factorization}
 X^N-1=\prod_{\gamma\in\mathcal R}M_{C_\gamma}(X),
 \qquad
 \Phi_N(X)=\prod_{\substack{\gamma\in\mathcal R\\\gcd(\gamma,N)=1}}M_{C_\gamma}(X),
\end{equation}
and the factor corresponding to \(\gamma\) has degree equal to the least common multiple of its local orbit lengths.
For odd \(N\), the \(m_i\) may be chosen as prime powers, and Proposition~\ref{v44:prop-prime-power-local-orbits}
then gives each \(\mathcal R_i\).
\end{corollary}
\begin{proof}
Products of different local orbits are disjoint and exhaust the CRT product space.
Proposition~\ref{v44:prop-recursive-phase-representatives} gives all global orbit representatives in each product.
The factorization formulas and factor degrees follow from Theorem~\ref{thm:coset-factorization}.
Primitive roots of unity correspond exactly to indices whose local components are all units, giving the formula for \(\Phi_N\).
\end{proof}

\Needspace{7\baselineskip}
\begin{example}[Local phases for fifteenth roots of unity]\label{v44:ex-local-phases-fifteen}
Take \(q=2\) and \(N=15=3\cdot5\). Each local unit group has just one orbit:
\[
 C_1\pmod3=\{1,2\},\qquad C_1\pmod5=\{1,2,4,3\},
\]
of lengths two and four. Their product has \(2\cdot4/4=2\) global orbits, each of length four.
The recursive phase representatives are \((0,0)\) and \((0,1)\).
We may take CRT coefficients \(E_1=10\) and \(E_2=6\), giving global representatives
\[
 10+6\equiv1\pmod{15},\qquad
 10+6\cdot2\equiv7\pmod{15}.
\]
Thus all unit indices split into
\[
 C_1=\{1,2,4,8\},\qquad C_7=\{7,14,13,11\}.
\]
The four-element orbit in Example~\ref{exm:med-fifteen-trivial} is the first of these two synchronized phases.
Specifying the local orbits does not distinguish it from the second; the difference is precisely one extra step in the modulo-five component relative to the modulo-three component.

Choose \(\zeta\) to be a primitive fifteenth root of unity that is a root of \(X^4+X+1\). Then
\[
 M_{C_1}(X)=X^4+X+1,\qquad M_{C_7}(X)=X^4+X^3+1.
\]
The second identity also follows from \(C_7=C_{-1}\) and inversion of roots for reciprocal polynomials.
The other three combinations of local orbits are zero--zero, nonzero--zero, and zero--nonzero,
with global orbit lengths one, two, and four.
Together with the two unit orbits above, they recover all five irreducible factors in the earlier example of fifteenth roots of unity.
\end{example}

\begin{corollary}[Phases and composition exponents]\label{v125:cor-phase-composition}
Fix the local orbits above, and write
\[
C_i=C_{m_i}(a_i),\qquad d_i=m_i/\gcd(m_i,a_i),\qquad d=\prod_i d_i.
\]
Every global orbit \(C\) assembled from them satisfies
\[
|C|=L,\qquad \omega_C=\ord_{b(d)}(q),\qquad
\kappa_C=L/\ord_{b(d)}(q).
\]
Thus relative phase changes the global orbit and its root polynomial, but not its maximum composition exponent or its least extension degree for binomial factorization.
\end{corollary}
\begin{proof}
Every index in the \(i\)-th local orbit has additive order \(d_i\). The Chinese remainder theorem gives the global index additive order \(\lcm(d_1,\ldots,d_s)=d\), independently of phase. Apply Theorem~\ref{v125:thm-orbit-structure}.
\end{proof}

\begin{example}[Phases modulo 160]\label{v125:ex-160-phases}
Continue Example~\ref{v125:ex-160-composition}. The units modulo 32 form two orbits \(A=\langle3\rangle\) and \(5A\), each of length eight; modulo five there is just one unit orbit \(B=\{1,3,4,2\}\), of length four.
For fixed \(A,B\), the phase quotient has \(8\cdot4/8=4\) elements. Using the Chinese remainder idempotents 65 and 96, the four global orbits have representatives
\[
65+96\cdot3^j\pmod{160},\qquad 0\le j<4,
\]
namely \(1,33,129,97\); \(j=0\) gives the preceding \(C_1\). The pair \(5A,B\) yields four more orbits, so \(\Phi_{160}\) has eight irreducible factors of degree eight.
All eight orbits have \(\kappa_C=2\) and \(\omega_C=4\). The maximum local block sizes are four for \(A\) and one for \(B\), whereas the global block size is two: local block sizes cannot simply be multiplied.
\end{example}

For each orbit determined in \eqref{v44:eq-crt-cyclotomic-factorization},
the multiple equal-difference structure of Section~\ref{sec:cyclotomic-med} can further compress the root product
and determine the least extension degree over which all its factors are binomials.
Local representatives enumerate and label the orbits; multiple equal-difference structure describes the additive symmetry within each orbit.
Together they give an explicit factorization of roots-of-unity polynomials.
For \(N=1\), the exponent set consists of one singleton orbit, with factor \(X-1\).

\subsection{Multiplier Orbits and Blocks of Factors}
\label{v50:sec-multiplier-blocks}
Frobenius orbits determine irreducible factors. Requiring root sets to be invariant under additional power maps groups these factors into larger blocks.
The exact coefficient field of a root block is determined by its period under Frobenius acting on the set of blocks~\cite{WZGeneralizedCosets2026,WZMultiplierOrbits2026}.

\begin{proposition}[Exact coefficient fields and factor blocks of multiplier orbits]\label{v50:prop-multiplier-blocks}
Suppose \(\gcd(N,q)=1\) and \(H\le(\Z/N\Z)^\times\). Write
\[
C_H(i;N)=\{hi\bmod N:h\in H\}.
\]
For \(d\mid N\), let \(H_d\) be the image of \(H\) under reduction in \((\Z/d\Z)^\times\). The indices of additive order \(d\) form \(\varphi(d)/|H_d|\) orbits under \(H\), each of length \(|H_d|\). Fix a primitive \(N\)-th root of unity \(\zeta\), and put
\[
B_{H,i}(X)=\prod_{j\in C_H(i;N)}(X-\zeta^j),\qquad
d=\frac N{\gcd(N,i)}.
\]
Define
\begin{equation}\label{v124:eq-coefficient-degree}
e_{q,H}(d)=\min\{r\ge1:q^r\bmod d\in H_d\},
\end{equation}
with \(e_{q,H}(1)=1\). The coefficients of \(B_{H,i}\) generate exactly \(\F_{q^{e_{q,H}(d)}}\) over \(\F_q\), and for every \(t\ge1\),
\[
B_{H,i}\in\F_{q^t}[X]\quad\Longleftrightarrow\quad e_{q,H}(d)\mid t.
\]
Put \(e=e_{q,H}(d)\) and \(J=\langle H,q\bmod N\rangle\). Coefficientwise Frobenius gives \(e\) distinct conjugate blocks, with
\begin{equation}\label{v124:eq-block-closure}
B_{H,i}^{[q^r]}=B_{H,q^ri},\qquad
\prod_{r=0}^{e-1}B_{H,q^ri}=B_{J,i}\in\F_q[X].
\end{equation}
Here \(P^{[q^r]}\) denotes the polynomial obtained by raising every coefficient of \(P\) to the \(q^r\)-th power.

In particular, \(B_{H,i}\in\F_q[X]\) if and only if \(q\bmod d\in H_d\). In that case it has exactly \(|H_d|/\ord_d(q)\) distinct monic irreducible factors, each of degree \(\ord_d(q)\), and is irreducible if and only if \(H_d=\langle q\bmod d\rangle\). If \(q\bmod N\in H\), all blocks are defined over the base field and
\[
X^N-1=\prod_{C_H(i;N)}B_{H,i}(X).
\]
\end{proposition}
\begin{proof}
An index of additive order \(d\) has a unique expression \((N/d)u\), with \(u\in(\Z/d\Z)^\times\). This scaling converts the \(H\)-action into multiplication by \(H_d\) on the unit group. Its orbits are cosets \(uH_d\), so each has length \(|H_d|\), and there are \(\varphi(d)/|H_d|\) orbits.

Frobenius sends the root exponents \(Hi\) to \(qHi\), giving the first identity in \eqref{v124:eq-block-closure}. A monic square-free polynomial is determined by its root set, so
\[
B_{H,i}^{[q^t]}=B_{H,i}
\iff q^tuH_d=uH_d
\iff q^t\bmod d\in H_d.
\]
The last condition is equivalent to \(e\mid t\): \(e\) is the order of \(qH_d\) in \((\Z/d\Z)^\times/H_d\). If the coefficients generate \(\F_{q^a}\), their least positive coefficientwise Frobenius return time is \(a\). Thus \(a=e\), proving the exact coefficient field.

The cosets \(uH_d,uqH_d,\ldots,uq^{e-1}H_d\) are disjoint and their union is \(u\langle H_d,q\rangle\). Multiplication over these root sets proves the product formula for the conjugate blocks.

When \(e=1\), one has \(\langle q\bmod d\rangle\le H_d\). Its cosets in \(uH_d\) are the Frobenius orbits, each of length \(\ord_d(q)\). Theorem~\ref{thm:coset-factorization} gives the stated irreducible factors and their number. If \(q\bmod N\in H\), every stratum has \(e=1\); the \(H\)-orbits partition the exponent set, giving the total product formula.
\end{proof}

Taking \(H=\langle q^t\bmod N\rangle\) and \(\tau=\ord_d(q)\) gives \(e_{q,H}(d)=\gcd(\tau,t)\). This recovers the coefficient-field formula of Proposition~\ref{prop:factor-coefficient-field} for roots of unity: base extension and multiplier fusion are both controlled by the Frobenius period of a root block.

Taking \(H=\langle q,-1\rangle\) includes inversion of roots.
If a \(q\)-cyclotomic coset is already stable under negation, its block is a self-reciprocal irreducible polynomial;
otherwise, two opposite cosets merge, giving the block \(M_C M_C^\vee\).
This corresponds to the inverse-orbit description in Section~\ref{sec:reciprocal-polynomials}.
For example, take \(N=7\) over \(\F_2\). Since \(\langle2,-1\rangle=(\Z/7\Z)^\times\),
the two nonzero Frobenius orbits of Example~\ref{v50:ex-seven-orbits} merge into the block \(\Phi_7\).
The block records invariance of the root set under the whole multiplier group, while the irreducible factors are still determined by its Frobenius orbits.

\section{Cyclotomic Systems and Compatible Factor Lifting}
\label{v44:sec-cyclotomic-systems}

When moduli increase along a divisibility chain, power maps between roots of unity relate cyclotomic cosets at different levels.
A cyclotomic system records this compatibility and corresponds to a compatible tower of irreducible factors.
Fix \(q=p^e\); all moduli below are coprime to \(p\).
For cyclotomic systems and their classification in prime-power directions, see \cite{ZWCosetsII}. We describe their lifting structure by return permutations on finite fibers~\cite{WuCyclotomicSystems}.

\subsection{Reduction Maps and Compatible Systems}

To compare different moduli, write an orbit modulo \(m\) as
\[
C_m(\gamma)=\{q^j\gamma\bmod m:j\ge0\},\qquad
\tau_m(\gamma)=|C_m(\gamma)|,
\]
and let \(\mathcal C_m\) be the finite set of all such orbits.
Theorem~\ref{thm:coset-factorization} gives
\[
\tau_m(\gamma)=\ord_{m/\gcd(m,\gamma)}(q),\qquad \ord_1(q)=1.
\]
If \(m\mid n\), reduction commutes with multiplication by \(q\), inducing a map
\[
\pi_{n,m}:\mathcal C_n\longrightarrow\mathcal C_m,
\qquad C_n(\gamma)\longmapsto C_m(\gamma).
\]
It is surjective: choose a representative of a parent orbit, lift it to an integer modulo \(n\), and take its orbit.
For \(l\mid m\mid n\), we also have \(\pi_{n,l}=\pi_{m,l}\pi_{n,m}\).

\Needspace{5\baselineskip}
\par
\begin{figure}[!htbp]
\centering
\begin{tikzpicture}[every node/.style={font=\small},
 nd/.style={draw=MutedSlate,rounded corners,fill=white,inner sep=3pt},
 ar/.style={-{Stealth[length=1.7mm]},draw=MutedSlate},
 path/.style={-{Stealth[length=2mm]},draw=Jade,line width=1.1pt}]
\node[nd,draw=Jade,fill=CoverMint] (a) at (0,0) {\(C_{1}(0)\)};
\node[nd] (b0) at (-2,-1.3) {\(C_{2}(0)\)};
\node[nd,draw=Jade,fill=CoverMint] (b1) at (4,-1.3) {\(C_{2}(1)\)};
\node[nd] (c0) at (-4,-2.6) {\(C_{4}(0)\)};
\node[nd] (c2) at (0,-2.6) {\(C_{4}(2)\)};
\node[nd,draw=Jade,fill=CoverMint] (c1) at (4,-2.6) {\(C_{4}(1)\)};
\node[nd] (d0) at (-5,-3.9) {\(C_{8}(0)\)};
\node[nd] (d4) at (-2,-3.9) {\(C_{8}(4)\)};
\node[nd] (d2) at (0.5,-3.9) {\(C_{8}(2)\)};
\node[nd,draw=Jade,fill=CoverMint] (d1) at (3,-3.9) {\(C_{8}(1)\)};
\node[nd] (d5) at (5.5,-3.9) {\(C_{8}(5)\)};
\node[nd] (e0) at (-6,-5.2) {\(C_{16}(0)\)};
\node[nd] (e8) at (-4,-5.2) {\(C_{16}(8)\)};
\node[nd] (e4) at (-2,-5.2) {\(C_{16}(4)\)};
\node[nd] (e2) at (0,-5.2) {\(C_{16}(2)\)};
\node[nd] (e10) at (2,-5.2) {\(C_{16}(10)\)};
\node[nd,draw=Jade,fill=CoverMint] (e1) at (4,-5.2) {\(C_{16}(1)\)};
\node[nd] (e5) at (6,-5.2) {\(C_{16}(5)\)};
\draw[ar] (a) -- (b0);
\draw[path] (a) -- node[right,font=\scriptsize]{\(1\)} (b1);
\draw[ar] (b0) -- (c0);
\draw[ar] (b0) -- (c2);
\draw[path] (b1) -- node[right,font=\scriptsize]{\(0\)} (c1);
\draw[ar] (c0) -- (d0);
\draw[ar] (c0) -- (d4);
\draw[ar] (c2) -- (d2);
\draw[path] (c1) -- node[right,font=\scriptsize]{\(0\)} (d1);
\draw[ar] (c1) -- (d5);
\draw[ar] (d0) -- (e0);
\draw[ar] (d0) -- (e8);
\draw[ar] (d4) -- (e4);
\draw[ar] (d2) -- (e2);
\draw[ar] (d2) -- (e10);
\draw[path] (d1) -- node[right,font=\scriptsize]{\(0\)} (e1);
\draw[ar] (d5) -- (e5);
\end{tikzpicture}
\caption{The full finite lifting tree from \(C_1(0)\). Digits on the green edges are admissible digits for the chosen compatible representatives. Extending a path requires choosing a further edge.}\label{v103:fig-cyclotomic-tree}
\end{figure}
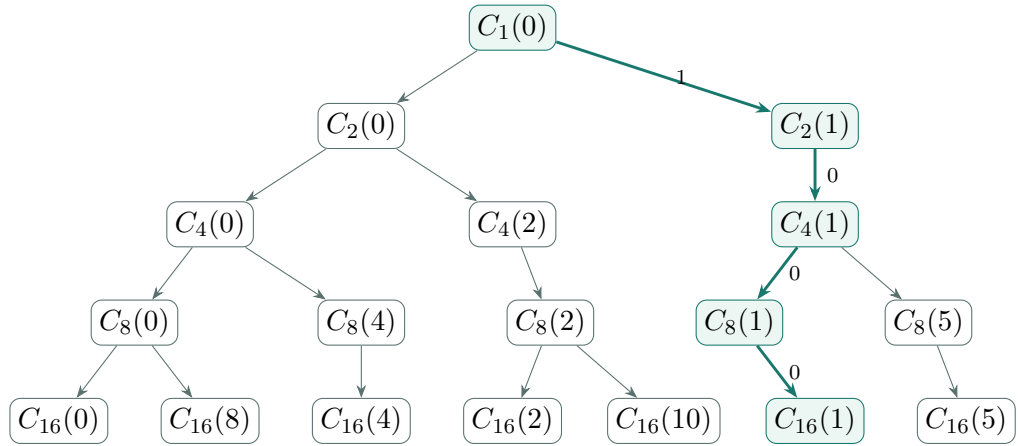

\Needspace{7\baselineskip}
\noindent\phantomsection\label{v103:par-cyclotomic-tree}\textbf{Start with a finite tree of five levels.}\quad
Take \(q=3\) and \(N_i=2^i\). Figure~\ref{v103:fig-cyclotomic-tree} lists all orbits modulo one, two, four, eight, and 16.
Each edge goes from a parent orbit to a lift; reduction follows the edge in reverse.
For example, \(C_4(1)=\{1,3\}\) has two children, \(C_8(1)=\{1,3\}\) and \(C_8(5)=\{5,7\}\).

\FloatBarrier
Extending the tree to arbitrarily high levels, a choice of one vertex at each level connected by successive edges is the compatible system defined below.
The green path uses representatives \(0,1,1,1,1\), with admissible digits \(1,0,0,0\).
Choosing \(C_8(5)\) instead changes these digits to \(1,0,1,0\).
The all-zero path is the principal system; the green path first leaves it at level 1, giving deviation depth \(u=0\).
Its orbit lengths are \(1,1,2,2,4\): growth can still be followed by splitting.
Theorem~\ref{v46:thm-cyclotomic-stable-growth} will show that the stabilization index is \(s=3\): from modulus eight onward, each lift is unique and doubles the length.
The finite diagram itself stops at modulus 16; the assertion for all later levels needs that theorem.

\begin{definition}[Compatible cyclotomic systems]\label{v44:def-cyclotomic-system}
Let \(N_0\mid N_1\mid\cdots\) be a divisibility chain.
A sequence \((C_i)_{i\ge0}\) satisfying
\[
C_i\in\mathcal C_{N_i},\qquad
\pi_{N_{i+1},N_i}(C_{i+1})=C_i\quad(i\ge0)
\]
is called a \textbf{cyclotomic system} on this chain.
The set of all compatible sequences is denoted by \(\varprojlim_i\mathcal C_{N_i}\); this \textbf{inverse limit} notation abbreviates precisely the compatibility condition above.

\leavevmode\phantomsection\label{v44:def-local-cyclotomic-system}\ignorespaces
If \(\ell\ne p\) is prime and \(\gcd(m,p\ell)=1\), take \(N_i=\ell^i m\). A cyclotomic system on this chain is called a \textbf{local cyclotomic system} in the \(\ell\)-direction with base level \(m\).
The set of all such systems is denoted by
\[
\mathsf{PC}_{q;\ell,m}=\varprojlim_{i\ge0}\mathcal C_{\ell^i m}.
\]
For \(C\in\mathcal C_m\), local systems whose initial term is \(C\) are called systems above \(C\).
\leavevmode\phantomsection\label{v44:def-total-cyclotomic-system}\ignorespaces
A \textbf{total cyclotomic system} is a family of orbits, one for each modulus coprime to \(p\), compatible with reduction whenever one modulus divides another.
The set of all such systems is denoted by
\[
\mathsf{PC}_q=
\left\{(C_n)_{\gcd(n,p)=1}:
C_n\in\mathcal C_n,\ \pi_{n,m}(C_n)=C_m\text{ for all }m\mid n\right\}.
\]
\end{definition}

Adjacent equal moduli in a chain give identity reduction maps, so deleting repetitions does not change the system.
Compatibility is imposed on orbits, not on arbitrarily chosen representatives; a recursive choice of compatible representatives is given by~\eqref{v44:eq-cyclotomic-admissible-digits}.

\par\Needspace{12\baselineskip}
\subsection{Affine Return on a Lifting Fiber}

A fiber of a map consists of all elements mapped to a specified point; see Section~\ref{sec:fibers}.

\begin{theorem}[Lifting cyclotomic cosets by return maps]\label{v44:thm-cyclotomic-affine-return}
Let \(m\mid n\), \(r=n/m\), choose \(0\le\gamma<m\), and put \(\tau=\tau_m(\gamma)\).
Define an affine permutation on \(\Z/r\Z\) by
\begin{equation}\label{v44:eq-cyclotomic-return}
T(k)=q^\tau k+\frac{\gamma(q^\tau-1)}m\pmod r.
\end{equation}
Let \(D(T)\) consist of the least nonnegative representative of each orbit of \(T\). The orbits modulo \(n\) reducing to \(C_m(\gamma)\) are precisely
\begin{equation}\label{v44:eq-cyclotomic-lifts}
\{C_n(\gamma+mk):k\in D(T)\},
\end{equation}
and are pairwise distinct. If the return orbit of \(k\) has length \(s\), then
\[
\tau_n(\gamma+mk)=\tau s.
\]
Consequently the lengths of these lifted orbits sum to \(r\tau\).
\end{theorem}
\begin{proof}
Since \(q^\tau\gamma\equiv\gamma\pmod m\), the numerator in \eqref{v44:eq-cyclotomic-return} is divisible by \(m\).
Also, \(\gcd(q,r)=1\), so \(T\) is indeed a permutation.

Among the residue classes modulo \(n\), consider the slice
\[
S=\{\gamma+mk:k\in\Z/r\Z\}.
\]
Every lifted orbit has image \(C_m(\gamma)\) modulo \(m\), and hence meets \(S\).
For \(x\in S\), if \(q^j x\in S\), then \(q^j\gamma\equiv\gamma\pmod m\), so \(\tau\mid j\).
The first-return action on the slice satisfies
\[
q^\tau(\gamma+mk)
=\gamma+m\left(q^\tau k+\frac{\gamma(q^\tau-1)}m\right).
\]
Thus the intersection of a lifted orbit with \(S\) is exactly one \(T\)-orbit, and each \(T\)-orbit produces a unique lifted orbit.
This proves completeness and nonredundancy in \eqref{v44:eq-cyclotomic-lifts}.

If the least positive period of \(k\) under \(T\) is \(s\), any return exponent of \(\gamma+mk\) must be a multiple of \(\tau\);
the first return occurs at \(\tau s\). Finally, the lengths of the \(T\)-orbits sum to \(r\); multiplication by \(\tau\) gives the total length.
\end{proof}

\begin{corollary}[Three branching types for a prime lift]\label{v44:cor-cyclotomic-prime-lifts}
If \(r=\ell\ne p\) is prime, write the return permutation as \(T(k)=ak+b\) on \(\F_\ell\), where
\[
a=q^\tau\bmod\ell,\qquad
b=\frac{\gamma(q^\tau-1)}m\bmod\ell.
\]
Then the lifts have the following three types:
\begin{enumerate}
\item \(a=1,b=0\): there are \(\ell\) lifted orbits, each of length \(\tau\).
\item \(a=1,b\ne0\): there is a unique lifted orbit, of length \(\ell\tau\).
\item \(a\ne1\): put \(k_*=-b/(a-1)\) and \(d=\ord_\ell(a)\).
There is one lift of length \(\tau\), and \((\ell-1)/d\) lifts of length \(d\tau\).
\end{enumerate}
\end{corollary}
\begin{proof}
In the first case, the return map is the identity. In the second, it is a nonzero translation, first returning after \(\ell\) iterations, so all points form a single \(\ell\)-cycle.
In the third, \(k_*\) is the unique fixed point; the change of variable \(u=k-k_*\) turns the return into \(u\mapsto au\).
Every nonzero point has period \(d\), giving \((\ell-1)/d\) nonzero orbits. Apply the return-lifting theorem.
\end{proof}

For explicit representatives in the third case, choose a generator \(g\) of \(\F_\ell^\times\), and put \(h=(\ell-1)/d\). Then
\[
\{k_*\}\cup\{k_*+g^j:0\le j<h\}
\]
gives representatives of all return orbits, since \(\langle a\rangle=\langle g^h\rangle\).
These may differ from the least integer representatives, but either choice makes \eqref{v44:eq-cyclotomic-lifts} nonredundant.
When \(\ell=2\), necessarily \(a=1\), so the two translation types exhaust all cases.

\par\Needspace{12\baselineskip}
\subsection{Compatible Roots of Unity and Towers of Irreducible Factors}
\begin{lemma}[Lifting primitive roots of unity]
\label{v44:lem-primitive-root-lift}
Let \(d,N\) be positive integers with \(d\mid N\) and \(p\nmid N\). For any primitive \(d\)-th root of unity \(\eta\),
one can choose, in a common finite extension, a primitive \(N\)-th root of unity \(\zeta\) such that
\[
 \zeta^{N/d}=\eta.
\]
\end{lemma}
\begin{proof}
Choose a finite extension containing \(\eta\) and all \(N\)-th roots of unity, and a primitive \(N\)-th root of unity \(\xi\) in it.
Then \(\xi^{N/d}\) has order \(d\), so \(\eta=(\xi^{N/d})^b\), where \(\gcd(b,d)=1\).
Choose an integer \(h\) with \(h\equiv b\pmod d\),
and impose \(h\equiv1\pmod\ell\) for every prime \(\ell\mid N\) with \(\ell\nmid d\).
These moduli are pairwise coprime, so the Chinese remainder theorem gives such an \(h\).
No prime divisor of \(d\) divides \(h\), and no prime dividing \(N\) but not \(d\) divides \(h\) either.
Thus \(\gcd(h,N)=1\). Taking \(\zeta=\xi^h\) gives \(\ord(\zeta)=N\) and \(\zeta^{N/d}=\eta\).
For \(d=1\), the condition modulo \(d\) is vacuous, and the same construction applies.
\end{proof}

Let \(m\mid n\) and \(r=n/m\). In a common finite extension, primitive roots of unity can be chosen so that
\begin{equation}\label{v44:eq-compatible-primitive-roots}
\zeta_n^r=\zeta_m.
\end{equation}
Lemma~\ref{v44:lem-primitive-root-lift} ensures that every primitive \(m\)-th root has such a primitive \(n\)-th preimage.
For an infinite divisibility chain, apply the lemma successively in one algebraic closure to obtain roots of unity satisfying this power relation at every pair of adjacent levels;
for the definition and existence of an algebraic closure, see Appendix~\ref{app:algebraic-closure}.

For \(C\in\mathcal C_m\), use the chosen primitive root to define
\[
M_{C,m}(X)=\prod_{a\in C}(X-\zeta_m^a).
\]
This is a monic irreducible factor of \(X^m-1\), of degree \(|C|\).

\begin{proposition}[Compatible factor lifting]\label{v44:prop-cyclotomic-factor-tower}
With the choice in \eqref{v44:eq-compatible-primitive-roots}, we have
\begin{equation}\label{v44:eq-cyclotomic-factor-lift}
\boxed{\quad
M_{C,m}(X^r)=
\prod_{\substack{D\in\mathcal C_n\\\pi_{n,m}(D)=C}}M_{D,n}(X).
\quad}
\end{equation}
The irreducible factors on the right are distinct. If \(C=C_m(\gamma)\), their degrees are determined by the return periods in Theorem~\ref{v44:thm-cyclotomic-affine-return}.
\end{proposition}
\begin{proof}
If \(x\) is a root of the left side, then \(x^r=\zeta_m^a\) for some \(a\in C\), and hence \(x^n=(x^r)^m=1\).
Thus \(x=\zeta_n^b\) for some \(b\bmod n\).
Compatibility of the roots of unity gives \(x^r=\zeta_m^b\), so \(x\) is a root precisely when \(b\bmod m\in C\).
These indices form the disjoint union of all orbits modulo \(n\) reducing to \(C\).

Since \(p\nmid r\), and \(M_{C,m}\) is separable with nonzero roots, substitution of \(X^r\) introduces no repeated roots: at any root \(x\), the derivative is
\[
r x^{r-1}M_{C,m}'(x^r)\ne0.
\]
Taking root products over the individual orbits gives the distinct factors on the right. Both sides are monic with the same root set, proving \eqref{v44:eq-cyclotomic-factor-lift}.
Finally, factor degrees equal exponent-orbit lengths; apply the return-lifting theorem.
\end{proof}

After choosing compatible primitive roots along the divisibility chain, a system \((C_i)\) corresponds to factors \(F_i=M_{C_i,N_i}\) satisfying
\[
F_i\mid X^{N_i}-1,\qquad
F_{i+1}\mid F_i\!\left(X^{N_{i+1}/N_i}\right),
\]
with each \(F_i\) monic and irreducible. Conversely, by \eqref{v44:eq-cyclotomic-factor-lift}, a sequence of factors satisfying these conditions determines a unique compatible system.
Thus lifting the modulus describes both the branching of exponent orbits and the factorization of power substitutions in the corresponding irreducible polynomials. Base extension instead changes the Frobenius step at a fixed modulus; the two operations have different effects.

\Needspace{5\baselineskip}
\begin{example}[Growth followed by further splitting in the dyadic direction]\label{v44:ex-cyclotomic-two-adic-three}
Take \(q=3\). The orbit \(C_2(1)=\{1\}\) modulo two has length one. From modulus two to four, the return is \(T(k)=k+1\) on \(\F_2\), so its unique lift
\[
C_4(1)=\{1,3\}
\]
has length two. Lifting next from modulus four to eight, we have \(\tau=2\) and
\[
a=3^2\bmod2=1,\qquad b=(3^2-1)/4\bmod2=0.
\]
It therefore splits into \(C_8(1)=\{1,3\}\) and \(C_8(5)=\{5,7\}\), each of length two.
Thus a length increase may still be followed by splitting at the next level.

From modulus eight to sixteen, the representatives one and five give, respectively,
\[
b=(3^2-1)/8=1,\qquad b=5(3^2-1)/8\equiv1\pmod2.
\]
Both returns are nonzero translations, so each has a unique lift whose length doubles to four:
\[
C_{16}(1)=\{1,3,9,11\},\qquad
C_{16}(5)=\{5,15,13,7\}.
\]
The corresponding polynomial relations are
\[
X^2+1\text{ is irreducible over }\F_3,\qquad
X^4+1=(X^2+X+2)(X^2+2X+2).
\]
Choose \(\zeta_8\) as a root of the first quadratic factor, then choose \(\zeta_{16}^2=\zeta_8\). The next level is
\[
X^8+1=(X^4+X^2+2)(X^4+2X^2+2),
\]
and both quartic factors are irreducible. They are the unique lifts obtained by substituting \(X^2\) into the respective quadratic factors above.
\end{example}

\subsection{Admissible Digits and Compatible Representatives}

If each step of the divisibility chain multiplies by one prime, so \(N_{i+1}=\ell_iN_i\) with \(\ell_i\ne p\), the return permutations give coordinates at successive levels.
Fix an initial representative \(0\le\gamma_0<N_0\), and put \(\tau_0=\tau_{N_0}(\gamma_0)\).
At step \(i\), denote \eqref{v44:eq-cyclotomic-return} by \(T_i\), let \(D(T_i)\) be the set of least nonnegative representatives of its orbits, and define recursively
\begin{equation}\label{v44:eq-cyclotomic-admissible-digits}
k_i\in D(T_i),\qquad
\gamma_{i+1}=\gamma_i+N_i k_i,\qquad
\tau_{i+1}=\tau_i\,|\operatorname{Orb}_{T_i}(k_i)|.
\end{equation}
\leavevmode\phantomsection\label{v44:def-admissible-digit-path}\ignorespaces
Call \((k_i)_{i\ge0}\) an \textbf{admissible digit path}.
Every path gives a system with initial term \(C_{N_0}(\gamma_0)\), and every such system has exactly one path:
given the parent orbit and its representative, the return-lifting theorem puts all child orbits in bijection with \(D(T_i)\), so both the successive choices and their inverses are unique.
The recursion also ensures
\[
0\le\gamma_i<N_i,\qquad \gamma_{i+1}\equiv\gamma_i\pmod{N_i}.
\]
Thus it produces compatible representatives, not just compatible orbits.

Any divisibility chain can be refined by factoring each successive quotient into primes and inserting finitely many intermediate levels.
The orbit at an intermediate level is uniquely determined by reduction from a higher level, so this refinement preserves the original compatibility information.
For local systems, take \(N_0=m\) and \(\ell_i=\ell\).
For total systems, take \(N_0=1\) and a sequence of lifting primes in which every \(\ell\ne p\) occurs infinitely often.
Then every \(n\) coprime to \(p\) divides some \(N_i\), and its orbit is obtained by reduction from that level; all sufficiently high levels give the same result.

A compatible representative \(\gamma_i\) need not be the coset leader \(\lambda_{N_i}(\gamma_i)\) of the entire orbit \(C_{N_i}(\gamma_i)\).
Choosing the coset leader independently at each level still represents the same system, but may violate \(\gamma_{i+1}\equiv\gamma_i\pmod{N_i}\), so the original digit recursion cannot simply be continued with those representatives.
For subsequent returns, retain the representatives chosen in \eqref{v44:eq-cyclotomic-admissible-digits}.

Once the local lifts are fixed, the remaining branches are described by the kernel of phase reduction.

\begin{proposition}[Phase reduction, its kernel, and global branches]\label{v125:prop-phase-lifting}
Let \(m_i\mid n_i\) for \(1\le i\le s\), where the \(n_i\) are pairwise coprime and coprime to \(q\).
Choose local orbits \(D_i\subseteq\Z/n_i\Z\) reducing to \(C_i\subseteq\Z/m_i\Z\), of lengths \(T_i,\tau_i\), respectively. Put
\[
L'=\lcm_i T_i,\quad L=\lcm_i\tau_i,\qquad
P_{\boldsymbol T}=\frac{\prod_i\Z/T_i\Z}{\langle(1,\ldots,1)\rangle},\quad
P_{\boldsymbol\tau}=\frac{\prod_i\Z/\tau_i\Z}{\langle(1,\ldots,1)\rangle}.
\]
After compatible local representatives are chosen, reduction of global orbits is described by the commutative diagram
\[
\begin{tikzpicture}[baseline=(current bounding box.center),>=stealth]
\node (a) at (0,1.35) {$\prod_i\Z/T_i\Z$};
\node (b) at (4.3,1.35) {$P_{\boldsymbol T}$};
\node (c) at (0,0) {$\prod_i\Z/\tau_i\Z$};
\node (d) at (4.3,0) {$P_{\boldsymbol\tau}$};
\draw[->] (a)--(b);\draw[->] (c)--(d);
\draw[->] (a)--(c);\draw[->] (b)--node[right] {$\rho$} (d);
\end{tikzpicture}
\]
where the horizontal maps quotient by common translation and the vertical maps reduce each coordinate.
Put \(u_i=T_i/\tau_i\), and let \(c=(L/\tau_i)_i\) in \(B=\prod_i\Z/u_i\Z\). There is a \textbf{short exact sequence}
\begin{equation}\label{v126:eq-phase-exact}
0\longrightarrow B/\langle c\rangle
\xrightarrow{\ \iota\ }P_{\boldsymbol T}
\xrightarrow{\ \rho\ }P_{\boldsymbol\tau}\longrightarrow0,
\qquad
\iota([z])=[(\tau_i z_i)_i].
\end{equation}
In particular, every phase fiber is a coset of \(\ker\rho\cong B/\langle c\rangle\), with exactly
\begin{equation}\label{v125:eq-phase-fiber}
K=|\ker\rho|=\frac{\prod_i(T_i/\tau_i)}{L'/L}
\end{equation}
elements. Thus, within the chosen local lifts \((D_i)\), every global orbit at the lower level has exactly \(K\) lifts, each of length \(L'\); the corresponding factor degrees are multiplied by \(L'/L\).
\end{proposition}
\begin{proof}
Choose \(b_i\in D_i\), and use its reduction \(a_i\in C_i\) as representative. Since \(\tau_i\mid T_i\), coordinate reduction is well defined and surjective. It maps the diagonal generator to the diagonal generator, and hence induces a surjection \(\rho\).
Theorem~\ref{v44:thm-local-orbit-synchronization} identifies the two phase quotients with the corresponding sets of global orbits. Compatibility of representatives makes these identifications commute with reduction of moduli.

The kernel of coordinate reduction is identified with \(B\) by \(z\mapsto(\tau_i z_i)_i\).
If a phase class lies in \(\ker\rho\), subtracting a suitable common translation gives a representative in this coordinate kernel.
The intersection of the coordinate kernel with the upper diagonal subgroup is generated by \((L,\ldots,L)\). It corresponds to \(c\) in \(B\), and has order \(L'/L\).
Thus \(\iota\) is injective with image \(\ker\rho\), proving exactness and the kernel order.
Every fiber of a group homomorphism is a kernel coset. Orbit lengths give factor degrees, and Proposition~\ref{v44:prop-cyclotomic-factor-tower} gives the divisibility relations between factors.
\end{proof}

A surjection always permits a separate choice of lift for each phase. An additive choice requires a homomorphism \(\jmath:P_{\boldsymbol\tau}\to P_{\boldsymbol T}\) with \(\rho\jmath=\id\); the short exact sequence is then called \textbf{split}. The following criterion uses only prime valuations of the local periods.

\begin{corollary}[Phase kernels and a splitting criterion]\label{v126:cor-phase-splitting}
Retain the notation of Proposition~\ref{v125:prop-phase-lifting}. For each prime \(\ell\mid L'\), put
\[
a_i=v_\ell(\tau_i),\quad b_i=v_\ell(T_i),\quad
a=\max_i a_i,\quad b=\max_i b_i,\quad c_i=b_i-a_i.
\]
For \(r\ge0\), define
\begin{equation}\label{v126:eq-phase-kernel-torsion}
\begin{aligned}
\eta_\ell(r)&=\max\bigl(0,\max_i\{\min(r,c_i)-a+a_i\}\bigr),\\
u_\ell(r)&=\sum_i\min(r,c_i)-\eta_\ell(r).
\end{aligned}
\end{equation}
The \(\ell\)-primary component of \(\ker\rho\) has exactly \(\ell^{u_\ell(r)}\) elements annihilated by \(\ell^r\), and
\[
(\ker\rho)_\ell\cong
\prod_{r=1}^{b}(\Z/\ell^r\Z)^{\,2u_\ell(r)-u_\ell(r-1)-u_\ell(r+1)}.
\]
The short exact sequence~\eqref{v126:eq-phase-exact} splits if and only if, for every \(\ell\mid L'\) and \(1\le r\le b\),
\begin{equation}\label{v126:eq-phase-split-test}
u_\ell(r)=\sum_i\bigl(\min(r,b_i)-\min(r,a_i)\bigr)
-\min(r,b)+\min(r,a).
\end{equation}
\end{corollary}
\begin{proof}
In the \(\ell\)-primary component of \(B\), multiplication by a unit in each coordinate transforms the relation element \(c\) into
\((\ell^{a-a_i})_i\in B_\ell=\prod_i\Z/\ell^{c_i}\Z\), of order \(\ell^{b-a}\).
The least nonnegative integer \(t\) such that \(\ell^t c\in\ell^rB_\ell\) is \(\eta_\ell(r)\), since the \(i\)-th coordinate requires
\(a-a_i+t\ge\min(r,c_i)\).
Hence \(\langle c\rangle\cap\ell^rB_\ell\) has order \(\ell^{b-a-\eta_\ell(r)}\). Comparing the orders of
\(B_\ell/\langle c\rangle\) and its image under multiplication by \(\ell^r\) gives \(u_\ell(r)\).
For a finite abelian \(\ell\)-group, \(u_\ell(r)-u_\ell(r-1)\) counts the cyclic direct factors of order at least \(\ell^r\). A second difference gives the stated decomposition.

For a finite abelian group \(G\), write \(G[\ell^r]=\{x:\ell^r x=0\}\).
The sequence splits if and only if all maps
\(P_{\boldsymbol T}[\ell^r]\to P_{\boldsymbol\tau}[\ell^r]\) are surjective. A section proves necessity. Conversely, lift a generator of each prime-power cyclic direct factor of the lower group to an element annihilated by the same prime power, and extend over the direct sum.
The kernels of these restricted maps are \((\ker\rho)[\ell^r]\), so surjectivity is equivalent to multiplication of the corresponding element counts.
In the \(\ell\)-primary component of \(P_{\boldsymbol T}\), choose a coordinate with \(b_i=b\). Subtracting its value eliminates the diagonal direction and identifies the quotient with the product of the remaining cyclic groups. Thus
\[
\log_\ell|P_{\boldsymbol T}[\ell^r]|=\sum_i\min(r,b_i)-\min(r,b),\qquad
\log_\ell|P_{\boldsymbol\tau}[\ell^r]|=\sum_i\min(r,a_i)-\min(r,a).
\]
Substitution of the kernel count gives~\eqref{v126:eq-phase-split-test}. For \(r>b\), all the groups are already annihilated, so no further test is needed.
\end{proof}

When the sequence splits, each phase has a unique expression \(\iota(k)+\jmath(\theta)\), with \(k\in B/\langle c\rangle\) and \(\theta\in P_{\boldsymbol\tau}\). This separates the branch coordinate from the lower phase as a direct sum of groups. Choosing sections at successive levels and composing them gives compatible additive lifts along the tower.

For \(s=2\), put \(g=\gcd(\tau_1,\tau_2)\) and \(g'=\gcd(T_1,T_2)\). Phase differences identify \(P_{\boldsymbol\tau}\cong\Z/g\Z\) and \(P_{\boldsymbol T}\cong\Z/g'\Z\), with \(\rho\) the reduction map. The sequence therefore splits if and only if \(\gcd(g,g'/g)=1\). Coprimeness gives a section by the Chinese remainder theorem; if a prime divides both numbers, a generator of the lower primary component has no lift of the same order.

Successive choices of local orbits are recorded by the admissible digits in~\eqref{v44:eq-cyclotomic-admissible-digits}. Once compatible representatives are chosen, phase classes satisfying \(\rho(\theta_{j+1})=\theta_j\) determine the global system and its factor tower uniquely.
Changing representatives only changes the origins of the phase coordinates. To count all global branches, sum~\eqref{v125:eq-phase-fiber} over the possible tuples of local child orbits.

In particular, suppose that a single local length changes from \(\tau_i\) to \(\ell\tau_i\), where \(\ell\) is prime, and all other lengths remain fixed. Then
\[
\begin{array}{c|cc}
&L'/L&K\\\hline
v_\ell(\tau_i)<v_\ell(L)&1&\ell\\
v_\ell(\tau_i)=v_\ell(L)&\ell&1
\end{array}
\]
Local growth can initially produce global phase splitting. It forces global length growth only after catching up with the \(\ell\)-parts of the other periods.

\begin{example}[Dyadic growth and phase branching]\label{v125:ex-160-tower}
Continue with \(q=3\), and let \(N_a=5\cdot2^a\), \(a\ge3\), restricting to unit indices.
Each unit orbit modulo \(2^a\) has length \(t_a=2^{a-2}\) and lifts uniquely at every step from \(a=3\) onward; the period modulo five remains four.
For fixed local orbits, therefore,
\[
\begin{array}{c|ccc}
a&3&4&a\ge5\\\hline
N_a&40&80&5\cdot2^a\\
t_a&2&4&2^{a-2}\\
L_a=\lcm(t_a,4)&4&4&2^{a-2}\\
|P_a|=\gcd(t_a,4)&2&4&4
\end{array}
\]
From modulus 40 to modulus 80, the local length doubles while the global length stays fixed, so each global orbit splits into two. From modulus 80 onward, each step has a unique global lift and doubles its length.
Thus the local stabilization index \(a=3\) differs from the global index \(a=4\).

Phase differences identify reduction at this step with
\[
0\longrightarrow\Z/2\Z\xrightarrow{\,1\mapsto2\,}
\Z/4\Z\longrightarrow\Z/2\Z\longrightarrow0.
\]
This sequence does not split: both lifts of the nonzero lower phase are odd classes of order four, so neither can be the homomorphic image of a generator of order two. Compatible orbit lifts still exist, but no additive section can be chosen at this step.

Choose compatible primitive roots so that \(C_1\) modulo 80 corresponds to \(G=X^4+X+2\) from Example~\ref{v125:ex-160-composition}. The factor-tower identity gives the corresponding factor at modulus \(N_a\):
\begin{equation}\label{v125:eq-160-tower}
F_a(X)=G(X^{2^{a-4}})
=X^{2^{a-2}}+X^{2^{a-4}}+2\qquad(a\ge4).
\end{equation}
Every such polynomial is irreducible, with \(\omega_{C_a}=4\) and \(\kappa_{C_a}=2^{a-4}\).
The degree-eight polynomial at modulus 160 is the case \(a=5\). Valuation formulas in the next subsection determine the finite stabilization index in general.
\end{example}

\subsection{Stable Growth of Local Systems}
\label{v46:sec-cyclotomic-stable-growth}

Once a prime direction is fixed, the lengths in a compatible system are determined by the first deviation and the local multiplicative orders.
We give the exact level from which all lifts are unique and translate it into a substitution formula for irreducible factors.
This local classification is given in \cite{ZWCosetsII,WuCyclotomicSystems}.

Let \(\ell\ne p\) be prime, let \(\gcd(m,p\ell)=1\), and fix a base-level orbit
\(C_m(\gamma)\), with \(0\le\gamma<m\). Put
\[
 N_i=m\ell^i,\qquad \tau=\tau_m(\gamma),\qquad Q=q^\tau.
\]
The Chinese remainder theorem uniquely determines an integer \(0\le\eta_i<N_i\) such that
\[
 \eta_i\equiv\gamma\pmod m,\qquad \eta_i\equiv0\pmod{\ell^i}.
\]
The orbits \(P_i=C_{N_i}(\eta_i)\) are compatible and have length \(\tau\) at every level.
\leavevmode\phantomsection\label{v46:def-principal-cyclotomic-system}\ignorespaces
We call \((P_i)_{i\ge0}\) the \textbf{principal cyclotomic system} above \(C_m(\gamma)\).
It is the unique system whose coordinate modulo \(\ell^i\) is zero at every level.

For any system \((C_i)_{i\ge0}\) different from the principal system, there is a unique \(u\ge0\) such that
\[
 C_i=P_i\quad(0\le i\le u),\qquad C_{u+1}\ne P_{u+1}.
\]
\leavevmode\phantomsection\label{v46:def-cyclotomic-deviation-degree}\ignorespaces
We call \(u\) its \textbf{deviation depth}.
At level \(u+1\), every index in the orbit is divisible by \(\ell^u\) but not by \(\ell^{u+1}\);
compatibility preserves this valuation at all higher levels.

\Needspace{5\baselineskip}
\begin{theorem}[Local lengths and eventual uniqueness of lifting]
\label{v46:thm-cyclotomic-stable-growth}
With this notation, a system of deviation depth \(u\) satisfies
\begin{equation}\label{v46:eq-cyclotomic-local-length}
 |C_i|=\tau\,\ord_{\ell^{\max(i-u,0)}}(Q)\qquad(i\ge0).
\end{equation}
Put \(b=i-u\). The relevant multiplicative orders and the starting levels of eventual stabilization are as follows.
\begin{enumerate}
\item If \(\ell\) is odd, put
\[
 d=\ord_\ell(Q),\qquad v=v_\ell(Q^d-1).
\]
Then, for \(b\ge1\),
\[
 \ord_{\ell^b}(Q)=d\ell^{\max(0,b-v)},\qquad s=u+v.
\]
\item If \(\ell=2\) and \(Q\equiv1\pmod4\), put \(v=v_2(Q-1)\). Then
\[
 \ord_{2^b}(Q)=2^{\max(0,b-v)}\quad(b\ge1),\qquad s=u+v.
\]
\item If \(\ell=2\) and \(Q\equiv3\pmod4\), put \(v=v_2(Q+1)\). Then
\[
 \ord_2(Q)=1,\qquad
 \ord_{2^b}(Q)=2^{\max(1,b-v)}\quad(b\ge2),\qquad s=u+v+1.
\]
\end{enumerate}
In each case, for every \(i\ge s\), the orbit \(C_i\) has a unique lift, and
\[
 |C_{i+1}|=\ell|C_i|.
\]
Thus the orbit at level \(s\) uniquely determines all higher levels of the system.
\leavevmode\phantomsection\label{v46:def-cyclotomic-stable-degree}\ignorespaces
We call \(s\) the \textbf{stabilization index} of the system.
It is the least level from which the stated uniqueness of lifting holds at every subsequent level.
\end{theorem}
\begin{proof}
Every \(C_i\) reduces to \(C_m(\gamma)\), so it contains a representative \(a_i\) with
\(a_i\equiv\gamma\pmod m\). A return exponent for \(a_i\) must have the form \(\tau t\).
When \(i\le u\), we have \(a_i\equiv0\pmod{\ell^i}\), so the least return exponent is \(\tau\).
When \(i>u\), we have \(v_\ell(a_i)=u\), and
\[
 q^{\tau t}a_i\equiv a_i\pmod{N_i}
 \quad\Longleftrightarrow\quad
 Q^t\equiv1\pmod{\ell^{i-u}}.
\]
The condition modulo \(m\) follows from \(q^\tau\gamma\equiv\gamma\pmod m\); the Chinese remainder theorem now gives
\eqref{v46:eq-cyclotomic-local-length}.

If \(\ell\) is odd, every exponent for which \(Q^t\equiv1\pmod{\ell^b}\) is divisible by \(d\).
Write \(t=dh\). Applying Lemma~\ref{lem:valuation-lift} to \(Q^d\) gives
\[
 v_\ell(Q^{dh}-1)=v+v_\ell(h).
\]
The least admissible \(h\) is \(\ell^{\max(0,b-v)}\), yielding the first formula.
The second formula follows directly from the same lemma.
In the third case, odd \(t\) satisfies \(v_2(Q^t-1)=1\);
for even \(t=2h\), the congruence \(Q^2\equiv1\pmod4\) and the lemma give
\[
 v_2(Q^{2h}-1)=v_2(Q^2-1)+v_2(h)=v+v_2(t).
\]
For \(b\ge2\), the least admissible \(t\) is \(2^{\max(1,b-v)}\), proving the third formula.

When \(i>u\), every index in every child orbit of \(C_i\) still has valuation \(u\).
Thus the length formulas above apply to all child orbits.
In the stated range \(i\ge s\), each child orbit has length \(\ell|C_i|\),
while Theorem~\ref{v44:thm-cyclotomic-affine-return} says that their lengths also sum to \(\ell|C_i|\).
There is therefore just one child orbit, and successive lifting uniquely determines the entire tail.

Finally, we check minimality. In the first two cases, if \(s-1>u\),
the length does not change from level \(s-1\) to level \(s\), so there are \(\ell\) child orbits.
If \(s-1=u\), the orbit at that level is principal, and its lifts include the principal child and at least one nonprincipal child.
In the third case, \(s-1=u+v>u\), and the length again remains unchanged from that level to level \(s\), giving two child orbits.
These splittings show that \(s\) is the least starting level for eventual uniqueness of lifting.
\end{proof}

The stabilization index depends on the deviation depth \(u\). For any finite level \(j\), the principal system and a system first deviating at level \(j+1\)
have the same first \(j\) levels. Hence no fixed finite level simultaneously determines all systems in the entire fiber above the base orbit.

\Needspace{5\baselineskip}
\begin{corollary}[Irreducible substitutions in stable factor towers]
\label{v46:cor-cyclotomic-stable-factor-substitution}
Choose compatible primitive roots of unity, and put \(F_i=M_{C_i,N_i}\).
If the system has stabilization index \(s\), then for every \(i\ge s\) and \(t\ge0\),
\begin{equation}\label{v46:eq-cyclotomic-stable-factor-substitution}
 F_{i+t}(X)=F_i(X^{\ell^t}),\qquad
 \deg F_{i+t}=\ell^t\deg F_i.
\end{equation}
In particular, \(F_i(X^{\ell^t})\) is irreducible in \(\F_q[X]\).
\end{corollary}
\begin{proof}
From level \(i\) onward, each level has a unique child orbit.
The factor product in Proposition~\ref{v44:prop-cyclotomic-factor-tower} therefore has just one factor, giving
\(F_{i+1}(X)=F_i(X^\ell)\). Iteration gives the displayed formula.
Each \(F_{i+t}\) is the root polynomial of one Frobenius orbit and is therefore irreducible.
The degree formula also follows from the growth in lengths.
\end{proof}

The finite input consists of the base orbit, deviation depth \(u\), and the relevant local orders and valuations. These determine the stabilization index \(s\) and all subsequent lengths; after choosing the orbit at level \(s\), no further factorization is needed for the tail. Before stabilization, the prime-step return map acts on only \(\ell\) digits, so child representatives do not require listing all indices modulo \(m\ell^{i+1}\). A single system can be stored as a finite initial segment and a power-substitution formula; expanded polynomials at each level require a separate allowance for coefficient output.

\begin{example}[Deeper splitting and sparse irreducible families after stabilization]
\label{v46:ex-cyclotomic-depth-eleven}
Take \(q=3\), \(m=1\), and \(\ell=11\). Here \(\tau=1\), and
\[
 d=\ord_{11}(3)=5,\qquad 3^5-1=242=2\cdot11^2,
 \qquad v=2.
\]
A system of deviation depth zero already has a nonzero coordinate modulo eleven,
but eventual uniqueness of lifting begins only at modulus 121.
Separating all orbits by the valuation of their indices gives the first levels:
\[
\begin{array}{c|ccc|c}
 \text{Modulus}&\text{Length 1}&\text{Length 5}&\text{Length 55}&\text{Total orbits}\\ \hline
 1&1&0&0&1\\
 11&1&2&0&3\\
 121&1&24&0&25\\
 1331&1&24&22&47
\end{array}
\]
For example, modulo 121 there are \(110/5=22\) orbits of valuation zero and \(10/5=2\) of valuation one.
Modulo 1331, there are \(1210/55=22\) unit orbits of length \(55\).
There are also \(22\) orbits of valuation \(1\) and \(2\) orbits of valuation \(2\), all of length \(5\).
At every level, the sum of orbit lengths equals the modulus.

On the polynomial side, \(\Phi_{11}\) has two irreducible quintic factors and \(\Phi_{121}\) has twenty-two.
Substitution of \(X^{11}\) into each factor of \(\Phi_{11}\) gives eleven quintic factors,
whereas every factor \(G\) of \(\Phi_{121}\) satisfies
\[
 G(X^{11^t})\text{ is irreducible in }\F_3[X],\qquad
 \deg G(X^{11^t})=5\cdot11^t\quad(t\ge0).
\]
For example, polynomial division gives \(G(X)=X^5+X^4+2\mid\Phi_{121}(X)\).
All factors of \(\Phi_{121}\) have degree five, so \(G\) is irreducible.
This gives the irreducible family
\[
 X^{5\cdot11^t}+X^{4\cdot11^t}+2\qquad(t\ge0).
\]
Here the stabilization threshold is determined by the valuation \(v=2\); the order \(d=5\) modulo eleven gives only the initial orbit length.
\end{example}

Primitive roots of unity of a fixed order all have the same Frobenius orbit length.
The valuation criteria then give a common finite stabilization threshold.
For the method of determining higher-index cyclotomic factorizations from finitely many lower-index factorizations, see \cite{WuZhuFengYang2017}.
We give an odd-prime formulation using root orders and Frobenius periods.

\begin{proposition}[A uniform stabilization threshold for primitive roots of unity]
\label{v50:prop-cyclotomic-uniform-stabilization}
Let \(\ell\) be an odd prime, \(\ell\nmid q\), let \(r\ge1\) with \(\gcd(r,q\ell)=1\), and put
\[
 h=\ord_{\ell r}(q)
 =\lcm\bigl(\ord_r(q),\ord_\ell(q)\bigr),\qquad
 b=v_\ell(q^h-1)\ge1.
\]
For every \(a\ge1\), we have
\begin{equation}\label{v50:eq-cyclotomic-uniform-order}
 \ord_{\ell^a r}(q)=h\ell^{\max(0,a-b)}.
\end{equation}
If \(\Phi_{\ell^b r}=\prod_{j=1}^s f_j\) is the monic irreducible factorization in \(\F_q[X]\),
then for all \(a\ge b\),
\begin{equation}\label{v50:eq-cyclotomic-uniform-substitution}
 \boxed{\quad
 \Phi_{\ell^a r}(X)
 =\prod_{j=1}^s f_j\bigl(X^{\ell^{a-b}}\bigr).
 \quad}
\end{equation}
Every factor on the right remains irreducible, of degree \(h\ell^{a-b}\), and
\(s=\varphi(\ell^b r)/h\) is independent of \(a\).
\end{proposition}
\begin{proof}
For \(a\ge1\), if \(\ell^a r\mid q^t-1\), then \(\ell r\mid q^t-1\), so \(h\mid t\).
Write \(t=hu\). Lemma~\ref{lem:valuation-lift} gives
\[
 v_\ell(q^{hu}-1)=b+v_\ell(u).
\]
Since \(r\mid q^h-1\), the condition \(\ell^a r\mid q^{hu}-1\) is equivalent to
\(v_\ell(u)\ge\max\{0,a-b\}\). The least admissible \(u\) is \(\ell^{\max(0,a-b)}\),
giving \eqref{v50:eq-cyclotomic-uniform-order}.

Put \(L=\ell^{a-b}\). If \(\beta^L\) has multiplicative order \(\ell^b r\), write \(d=\ord(\beta)\). Then
\[
 d\mid\ell^a r,\qquad \frac{d}{\gcd(d,L)}=\ell^b r.
\]
Since \(b\ge1\), comparison of \(\ell\)-adic valuations gives \(v_\ell(d)=a\);
comparison at the other primes gives \(d=\ell^a r\). Conversely, the \(L\)-th power of an element of order \(\ell^a r\) has order \(\ell^b r\).
Thus the roots of \(\Phi_{\ell^b r}(X^L)\) are precisely all primitive \(\ell^a r\)-th roots of unity.
Since \(\ell^a r\) is coprime to the characteristic, these roots are distinct; comparison of leading coefficients proves \eqref{v50:eq-cyclotomic-uniform-substitution}.

By Corollary~\ref{cor:cyclotomic-equal-degree}, \(\deg f_j=h\).
Every root of \(f_j(X^L)\) has multiplicative order \(\ell^a r\)
and Frobenius period \(hL\), exactly the degree of \(f_j(X^L)\), so this polynomial is irreducible.
The roots of different factors map to distinct base-level orbits under the \(L\)-th power map, so the displayed factors are distinct.
Their number is the base-level degree \(\varphi(\ell^b r)\) divided by \(h\).
\end{proof}

When \(q\equiv1\pmod\ell\), one may take \(h=\ord_r(q)\). In general, the period modulo \(\ell\) must also be included;
the stabilization threshold is determined by the valuation of \(q^h-1\), not merely that of \(q-1\).
Thus, for fixed \(r\), only the finitely many initial factorizations of \(\Phi_r,\Phi_{\ell r},\ldots,\Phi_{\ell^b r}\) are needed;
all higher levels follow by power substitution. Here we restrict to primitive \(\ell^a r\)-th roots, whose local indices all have valuation zero.
For the entire fiber above the base orbit, containing arbitrary deviation depths, the earlier conclusion that no common finite determining level exists still applies.

\begin{example}[From factors of order fifteen to order \(3^a5\)]
\label{v50:ex-cyclotomic-three-power-five}
Take \(q=7\), \(\ell=3\), and \(r=5\). Since \(\ord_{15}(7)=4\) and
\[
 7^4-1=2400=3\cdot800,
\]
we have \(h=4\) and \(b=1\). Direct multiplication in \(\F_7[X]\) gives
\[
 \Phi_{15}(X)=f_1(X)f_2(X),
\]
where
\[
 f_1=X^4+2X^3+4X^2+X+2,\qquad
 f_2=X^4+4X^3+2X^2+X+4.
\]
Both are irreducible, since all irreducible factors of \(\Phi_{15}\) have degree four.
Thus, for \(a\ge1\),
\[
 \Phi_{3^a5}(X)
 =f_1\bigl(X^{3^{a-1}}\bigr)f_2\bigl(X^{3^{a-1}}\bigr).
\]
There are therefore exactly two irreducible factors at every level, each of degree \(4\cdot3^{a-1}\).
In particular, both irreducible factors of \(\Phi_{45}\) have degree twelve.

The initial level must reach the stabilization threshold. If instead \(q=19\), \(\ell=3\), and \(r=1\), then \(h=1\) and \(b=2\).
Although \(X-7\mid\Phi_3\), we have
\[
 X^3-7=(X-4)(X-6)(X-9)\quad\text{in }\F_{19}[X],
\]
so cubic power substitution at the first level still splits; irreducible lifting must start with factors of \(\Phi_9\).
\end{example}

\subsection{Binomial Root Sets and Factorization with Multiplicities}
\label{v44:sec-binomial-orbits}

Cyclotomic orbits also describe irreducible factorization of general binomials.
Let \(q=p^e\), \(\lambda\in\F_q^\times\), and \(g=\ord(\lambda)\), so \(g\mid q-1\).
First assume \(n\ge1\) and \(p\nmid n\). The roots of \(X^n-\lambda\) embed in the group of \(ng\)-th roots of unity;
their exponent set is a Frobenius-stable additive coset modulo \(ng\). For explicit representatives and lengths in the corresponding cyclotomic systems, see \cite{ZWCosetsII}.

\begin{proposition}[Cyclotomic-orbit factorization of a binomial]
\label{v44:prop-binomial-orbits}
One can choose a primitive \(ng\)-th root of unity \(\zeta\) with \(\zeta^n=\lambda\). Put
\[
 S_{n,g}=1+g\Z/ng\Z
 =\{1,1+g,\ldots,1+(n-1)g\}\subseteq\Z/ng\Z.
\]
Then \(qS_{n,g}=S_{n,g}\). For each \(q\)-cyclotomic orbit \(C\) in this set, put
\[
 M_C(X)=\prod_{a\in C}(X-\zeta^a).
\]
These polynomials are monic, irreducible, and pairwise distinct, and
\begin{equation}\label{v44:eq-binomial-orbits}
 \begin{aligned}
 X^n-\lambda&=\prod_{C\subseteq S_{n,g}}M_C(X),\\
 \deg M_C&=|C|=\ord_{ng/\gcd(ng,a)}(q)\quad(a\in C).
 \end{aligned}
\end{equation}
\end{proposition}
\begin{proof}
Since \(p\nmid ng\), apply Lemma~\ref{v44:lem-primitive-root-lift} with \(d=g\), \(N=ng\), and \(\eta=\lambda\) to choose the required \(\zeta\).
The full root set is \(\zeta\mu_n\), with \(\mu_n=\langle\zeta^g\rangle\), so its exponent set is exactly \(S_{n,g}\).
Since \(q\equiv1\pmod g\), multiplication by \(q\) preserves this set; as \(q\) is invertible modulo \(ng\), we have \(qS_{n,g}=S_{n,g}\).
The map \(\zeta^a\mapsto\zeta^{qa}\) is precisely Frobenius on the root set.
By Theorem~\ref{thm:orbit}, the root product for each orbit is monic and irreducible, of degree equal to the orbit length.
Distinct orbits are disjoint and together exhaust all \(n\) roots, giving the product formula.
Finally, \(\zeta^a\) has multiplicative order \(ng/\gcd(ng,a)\),
so Proposition~\ref{prop:degree-from-order} gives the stated length.
\end{proof}

\begin{corollary}[The least extension for binomial factorization of a general binomial]\label{v124:cor-binomial-extension}
Let \(n\ge1\), \(p\nmid n\), \(\lambda\in\F_q^\times\), and \(g=\ord(\lambda)\). Put
\begin{equation}\label{v124:eq-binomial-index}
\Omega_{n,\lambda}=\ord_{b(ng)}(q),
\end{equation}
where \(b\) is as defined in Section~\ref{sec:cyclotomic-med}. For every \(t\ge1\), all irreducible factors of \(X^n-\lambda\) over \(\F_{q^t}\) are binomials if and only if \(\Omega_{n,\lambda}\mid t\). Thus the smallest such extension is \(\F_{q^{\Omega_{n,\lambda}}}\), and its degree depends on \(\lambda\) only through \(g\).
\end{corollary}
\begin{proof}
Proposition~\ref{v44:prop-binomial-orbits} identifies the root exponents with \(S_{n,g}=1+g\Z/ng\Z\). This set contains one, and hence contains an orbit with reduced modulus \(ng\). By Theorem~\ref{cor:med-binomial-extension} and \eqref{v124:eq-small-step}, all irreducible descendants of that orbit are binomials if and only if \(\ord_{b(ng)}(q)\mid t\). This proves necessity.

Every other orbit has reduced modulus \(d\mid ng\). Since \(b(d)\mid b(ng)\), one has \(\ord_{b(d)}(q)\mid\ord_{b(ng)}(q)\). The same divisibility condition therefore makes all irreducible descendants of every orbit binomials, proving sufficiency.
\end{proof}

\begin{corollary}[Characteristic powers and factor multiplicities]
\label{v44:cor-binomial-multiplicities}
Let \(N=p^v n\), where \(v\ge0\) and \(p\nmid n\).
There is a unique \(\lambda_0\in\F_q^\times\) with \(\lambda_0^{p^v}=\lambda\), and \(\ord(\lambda_0)=g\).
Applying Proposition~\ref{v44:prop-binomial-orbits} to \(X^n-\lambda_0\) gives
\begin{equation}\label{v44:eq-binomial-multiplicities}
 X^N-\lambda=(X^n-\lambda_0)^{p^v}
 =\prod_{C\subseteq S_{n,g}}M_C(X)^{p^v}.
\end{equation}
Thus orbit lengths still determine the degrees of the distinct monic irreducible factors, each with multiplicity exactly \(p^v\).
\end{corollary}
\begin{proof}
The map \(x\mapsto x^{p^v}\) is an automorphism of \(\F_q\), so \(\lambda_0\) exists uniquely and has the same multiplicative order.
Taking powers in characteristic \(p\) gives the first equality in \eqref{v44:eq-binomial-multiplicities}.
Since \((X^n-\lambda_0)'=nX^{n-1}\) and \(\lambda_0\ne0\), this polynomial has no repeated roots.
Its irreducible factors are distinct, and raising it to the \(p^v\)-th power gives each multiplicity exactly \(p^v\).
\end{proof}

Since \(\ord(\lambda_0)=g\), Corollary~\ref{v124:cor-binomial-extension} also gives the least extension degree over which all distinct irreducible factors of \(X^N-\lambda\) are binomials; the multiplicity \(p^v\) does not change this condition.

\begin{example}[A quartic binomial and its repeated-root lift]
\label{v44:ex-binomial-four}
Over \(\F_3\), the binomial \(X^4+1=X^4-(-1)\) has \(g=2\) and \(ng=8\).
The exponent set \(S_{4,2}=\{1,3,5,7\}\) splits into two \(3\)-cyclotomic orbits,
\[
 \{1,3\},\qquad\{5,7\}.
\]
Thus there are exactly two irreducible quadratic factors, namely
\[
 X^4+1=(X^2+X+2)(X^2-X+2).
\]
Both quadratics have discriminant \(2\), a nonsquare in \(\F_3\), and are therefore irreducible.
Here \(b(8)=4\) and \(\ord_4(3)=2\), so Corollary~\ref{v124:cor-binomial-extension} gives \(\F_9\) as the smallest extension for binomial factorization; both quadratics split over this field. By contrast, \(X^4-1\) has \(g=1\) and \(b(4)=2\), so the least extension degree is one, and
\[
X^4-1=(X-1)(X+1)(X^2+1)
\]
is already a factorization into irreducible binomials over \(\F_3\).

Also, \((-1)^3=-1\), so
\[
 X^{12}+1=(X^4+1)^3
 =(X^2+X+2)^3(X^2-X+2)^3.
\]
There are still just two distinct root orbits, while each irreducible factor now has multiplicity three.
\end{example}

Besides the root-exponent orbits, a norm class gives a common divisibility constraint on all factor degrees; see Proposition~\ref{v97:prop-norm-factor-degrees}. The support period of a single factor constrains the scalar action remaining in the splitting of the whole binomial; see Proposition~\ref{v97:prop-support-splitting}.

\section{Galois Rings and Lifting of Roots-of-Unity Factors}
\label{v47:sec-galois-rings}

Does the composition structure determined by root-set symmetry over a finite field persist when coefficient precision increases?

Over a finite field, roots-of-unity factors are determined by Frobenius orbits. Over a Galois ring of characteristic \(p^s\),
roots of unity of order coprime to \(p\) retain the same orbit structure, while the coefficients of their factors are uniquely determined by coprime Hensel lifting.
When \(p\mid n\), the residual root orbits correspond to canonical comaximal blocks.
For the standard theory of Galois rings, see \cite{WanGalois2003}; for lifting cyclotomic cosets and multiple equal-difference structure, see \cite{WuGaloisOrbitReport}.
Below, \(s\) is the exponent in the characteristic \(p^s\), and \(q=p^r\) is fixed.

\subsection{Quotient Ring Models and Coprime Hensel Lifting}

\begin{definition}[Galois rings and basic irreducible polynomials]
\label{v47:def-galois-ring}
Let \(p\) be prime, \(s,r\ge1\), and let \(P\in(\Z/p^s\Z)[T]\) be monic of degree \(r\),
with irreducible reduction \(\overline P\) modulo \(p\). The ring
\[
 A=(\Z/p^s\Z)[T]/(P)
\]
is called a \textbf{Galois ring}, denoted by \(\operatorname{GR}(p^s,r)\).
Its isomorphism type depends only on \(p,s,r\). \leavevmode\phantomsection\label{v47:def-basic-irreducible}\ignorespaces
A polynomial \(G\in A[X]\) is called \textbf{basic irreducible} if it is monic and
\(\overline G\in\F_q[X]\) is irreducible.
Here and below, reduction is coefficientwise modulo \(p\).
The second parameter \(r\) is the degree of the residue field over \(\F_p\).
In particular, \(\GR(p,r)\cong\F_{p^r}\) and \(\GR(p^s,1)\cong\Z/p^s\Z\).
\end{definition}

\leavevmode\phantomsection\label{v47:def-local-ring}\ignorespaces
A commutative ring with a unique maximal ideal is called a \textbf{local ring};
its quotient by that maximal ideal is called its \textbf{residue field}.
\leavevmode\phantomsection\label{v47:def-chain-ring}\ignorespaces
A commutative ring whose ideals are totally ordered by inclusion is called a \textbf{chain ring}.

\begin{proposition}[Units, the residue field, and the ideal chain]
\label{v47:prop-galois-ring-structure}
The ring \(A\) has characteristic \(p^s\), cardinality \(p^{sr}\), and unique maximal ideal \(pA\), with \(A/pA\cong\F_q\).
An element \(u\in A\) is invertible if and only if \(\overline u\ne0\). Every nonzero element determines a unique
\(0\le j<s\) for which it can be written as \(p^ju\) with \(u\) a unit. All ideals form the chain
\[
 A\supset pA\supset\cdots\supset p^{s-1}A\supset p^sA=0,
 \qquad\operatorname{Ann}_A(p^j)=p^{s-j}A\quad(0\le j\le s).
\]
\end{proposition}
\begin{proof}
Monic division shows that the images of \(1,T,\ldots,T^{r-1}\) form a basis of \(A\) over \(\Z/p^s\Z\).
Thus its characteristic is \(p^s\), its cardinality is \(p^{sr}\), and its quotient modulo \(p\) is \(\F_p[T]/(\overline P)\cong\F_q\).
If \(\overline u\ne0\), lifting its inverse gives \(uv=1+pw\). Since \((pw)^s=0\),
a finite geometric sum gives the inverse of \(1+pw\), so \(u\) is a unit. Conversely, a unit reduces to a unit.
The least \(p\)-adic valuation among the basis coordinates of a nonzero element is \(j<s\); after extracting \(p^j\), the remaining element has nonzero reduction.
Thus \(j\) is unique, and the least such \(j\) in a nonzero ideal determines that ideal as \(p^jA\).
The same basis coordinates directly give the annihilator identity. The nonunits are exactly \(pA\), so it is the unique maximal ideal.
\end{proof}

\begin{lemma}[Hensel lifting of simple roots]
\label{v47:lem-hensel-simple-root}
Let \(f\in A[X]\), and suppose \(\overline a\in\F_q\) satisfies
\(\overline f(\overline a)=0\), \(\overline f'(\overline a)\ne0\).
There is a unique \(a\in A\) with \(f(a)=0\) and \(a\bmod p=\overline a\).
This process is called \textbf{Hensel lifting}.
\end{lemma}
\begin{proof}
If \(a_j\) satisfies \(f(a_j)\equiv0\pmod{p^j}\), \(1\le j<s\),
then \(f(a_j)/p^j\bmod p\) is well defined: any two quotients differ by an element of
\(\operatorname{Ann}_A(p^j)=p^{s-j}A\subseteq pA\). Put
\[
 a_{j+1}=a_j+p^j\delta_j,\qquad
 \overline{\delta_j}=-\frac{f(a_j)}{p^j}\,
       \overline f'(\overline a)^{-1}\quad\text{in }\F_q.
\]
On expansion, all correction terms of degree at least two contain \(p^{2j}\) and vanish modulo \(p^{j+1}\), so
\(f(a_{j+1})\equiv0\pmod{p^{j+1}}\). Every lift must satisfy the same linear correction equation,
so each level is unique. Iteration to \(j=s\) gives the desired root; for \(s=1\), the result is immediate.

\end{proof}

\begin{theorem}[Hensel lifting of coprime factors]
\label{v47:thm-hensel-coprime}
Let \(f\in A[X]\) be monic, and suppose
\(\overline f=g_1\cdots g_v\), where the \(g_i\in\F_q[X]\) are monic and pairwise coprime.
There is a unique family of monic polynomials \(G_i\in A[X]\) satisfying
\[
 \overline{G_i}=g_i,\qquad f=\prod_{i=1}^vG_i.
\]
These factors are pairwise \textbf{comaximal}: \((G_i,G_j)=A[X]\) whenever \(i\ne j\).
\end{theorem}
\begin{proof}
First consider two coprime factors. Let \(\overline f=gh\), and suppose monic lifts
\(G_j,H_j\) of degrees \(\deg g,\deg h\) have been chosen with \(f\equiv G_jH_j\pmod{p^j}\). Put
\[
 D_j=\frac{f-G_jH_j}{p^j}\bmod p.
\]
The leading terms cancel, so \(\deg D_j<\deg g+\deg h\). Since \(h\) is invertible modulo \(g\), there is a unique polynomial
\(U_j\) with \(\deg U_j<\deg g\) and \(hU_j\equiv D_j\pmod g\). Put
\(V_j=(D_j-hU_j)/g\). Then \(\deg V_j<\deg h\), and
\[
 D_j=hU_j+gV_j.
\]
Choose coefficient lifts of \(U_j,V_j\) to \(A[X]\) without increasing degrees, and retain the same notation.
Then \(G_{j+1}=G_j+p^jU_j\) and \(H_{j+1}=H_j+p^jV_j\) are monic factors at the next level.
Subtracting two correction equations, reducing modulo \(g\), and using the degree bounds gives uniqueness of \(U_j\), and then of \(V_j\).
Iteration proves the two-factor case; grouping and induction give the general case, with constant factors \(1\) omitted.
Finally, lifting a B\'ezout identity for \(g_i,g_j\) gives
\(G_iU+G_jV=1+pW\). The polynomial \(pW\) is nilpotent, so the right side is invertible, proving comaximality.
\end{proof}

\leavevmode\phantomsection\label{v47:prop-comaximal-crt}\ignorespaces
Comaximality allows the Chinese remainder theorem to be applied to factors over a ring.
In general, if ideals \(I,J\) in a commutative ring \(R\) satisfy \(I+J=R\), choose \(u\in I\) and \(v\in J\) with \(u+v=1\).
For \(x\in I\cap J\), we have \(x=xu+xv\in IJ\), so \(I\cap J=IJ\).
Given \(a+I\) and \(b+J\), the element \(av+bu\) maps to these two residue classes. Hence
\[
 R/(IJ)\cong R/I\times R/J.
\]
Induction gives the corresponding direct-product decomposition for finitely many pairwise comaximal ideals.

\par\Needspace{12\baselineskip}
\subsection{Teichm\"uller Representatives and the Frobenius Automorphism}

\begin{definition}[The Teichm\"uller set]
\label{v47:def-teichmuller}
The set
\[
 \mathcal T(A)=\{t\in A:t^q=t\}
\]
is called the \textbf{Teichm\"uller set}, and its elements are called
\textbf{Teichm\"uller representatives}.
\end{definition}

\begin{theorem}[Digit expansions and the Frobenius automorphism]
\label{v47:thm-galois-frobenius}
Reduction gives a bijection \(\mathcal T(A)\to\F_q\) and an isomorphism of cyclic groups
\(\mathcal T(A)\setminus\{0\}\cong\F_q^\times\).
Every \(z\in A\) has a unique
\textbf{\(p\)-adic expansion}
\begin{equation}\label{v47:eq-teichmuller-expansion}
 z=t_0+pt_1+\cdots+p^{s-1}t_{s-1},\qquad t_j\in\mathcal T(A).
\end{equation}
There is a unique ring automorphism \(\phi\) lifting the \(p\)-th power automorphism of \(\F_q\),
called the \textbf{absolute Frobenius automorphism}, and
\begin{equation}\label{v47:eq-galois-frobenius}
 \phi(z)=\sum_{j=0}^{s-1}p^jt_j^p.
\end{equation}
Its order is \(r\), and its fixed subring is the embedded copy of \(\Z/p^s\Z\).
\end{theorem}
\begin{proof}
Modulo \(p\), the polynomial \(X^q-X\) has precisely \(\F_q\) as its root set, and its derivative is the constant \(-1\).
Lemma~\ref{v47:lem-hensel-simple-root} gives a unique Teichm\"uller representative for each residue class.
Nonzero representatives are units satisfying \(t^{q-1}=1\). They are closed under multiplication and inversion,
and reduction preserves multiplication, so they form a cyclic group isomorphic to \(\F_q^\times\).
Match residue classes modulo \(p\) successively, subtract the representative, and extract a factor \(p\); this gives \eqref{v47:eq-teichmuller-expansion}.
If two expansions first differ at digit \(j\), their difference is \(p^j\) times a unit, with \(j<s\), and cannot vanish. Thus the expansion is unique.

To show that \eqref{v47:eq-galois-frobenius} defines a ring automorphism, use the quotient model and put
\(\alpha=T+(P)\). The element \(\overline\alpha^{p}\) is a simple root of \(\overline P\),
so there is a unique \(\beta\in A\) with \(P(\beta)=0\) and \(\overline\beta=\overline\alpha^{p}\).
Substitution \(\alpha\mapsto\beta\) defines a ring homomorphism \(\phi\) fixing \(\Z/p^s\Z\),
whose reduction is the \(p\)-th power automorphism. Similarly, lifting the \(p^{r-1}\)-th power of \(\overline\alpha\) to a root of \(P\)
gives an inverse homomorphism: both composites send \(\alpha\) to a root of \(P\) reducing to \(\overline\alpha\).
By uniqueness, that root is \(\alpha\). Any ring automorphism with the specified reduction must act this way on \(\alpha\), so the lift is unique.

For \(t\in\mathcal T(A)\), both \(\phi(t)\) and \(t^p\) are Teichm\"uller representatives of the same residue class,
so \(\phi(t)=t^p\). A ring homomorphism preserves \(p\), giving the digit formula.
The map \(\phi^r\) reduces to the identity, and the same root uniqueness gives \(\phi^r=\id\).
The reduction has order \(r\), so \(\phi\) has order exactly \(r\). The digit formula shows that every digit of a fixed element reduces to an element of \(\F_p\).
Applying the same simple-root lifting and digit expansion within \(\Z/p^s\Z\) shows that these elements form precisely that subring.
\end{proof}

This proof also shows that the isomorphism type of a Galois ring is independent of the choice of \(P\).
In another quotient ring with the same parameters, send \(\overline\alpha\) by a residue-field isomorphism to a root of \(\overline P\),
then lift it uniquely to a root of \(P\), obtaining a ring homomorphism. Its reduction modulo \(p\) is an isomorphism,
and the image of every nonzero element \(p^ju\), \(j<s\), remains nonzero, so it is injective; equality of cardinalities makes it an isomorphism.

For a positive integer \(d\), the same method embeds \(A\) in \(S=\operatorname{GR}(p^s,rd)\):
first embed the residue fields \(\F_q\subseteq\F_{q^d}\), then lift the corresponding simple root of \(P\).
\leavevmode\phantomsection\label{v47:def-unramified-extension}\ignorespaces
Such an extension of Galois rings with the same characteristic is called an \textbf{unramified extension}.
\leavevmode\phantomsection\label{v47:def-galois-relative-frobenius}\ignorespaces
Let \(\phi_S\) be the absolute Frobenius automorphism of \(S\). Then
\[
 \sigma=\phi_S^r,\qquad
 \sigma\left(\sum_jp^jt_j\right)=\sum_jp^jt_j^q
\]
is the \textbf{relative Frobenius automorphism},
of order \(d\). In the embedded model, \(\alpha\) and \(\sigma(\alpha)\) are roots of \(P\)
with the same residue class, so uniqueness of simple-root lifting gives \(\sigma(\alpha)=\alpha\).
Thus \(\sigma\) fixes \(A\). The digit formula shows that it has exactly \(q^s=|A|\) fixed elements,
so its fixed subring is precisely \(A\).
On Teichm\"uller elements, \(\sigma(t)=t^q\), but in general \(\sigma(z)=z^q\) does not hold on the whole ring.

\Needspace{7\baselineskip}
\begin{example}[Frobenius in characteristic \(4\)]
\label{v47:ex-gr-four}
In \(A=\operatorname{GR}(4,2)=(\Z/4\Z)[\omega]\), take
\(\omega^2+\omega+1=0\). Then \(\omega^3=1\),
\(\mathcal T(A)=\{0,1,\omega,\omega^2\}\), and
\[
 X^3-1=(X-1)(X-\omega)(X-\omega^2).
\]
Absolute Frobenius interchanges \(\omega\) and \(\omega^2\) and fixes \(\Z/4\Z\).
In particular, \(\phi(2)=2\), whereas \(2^2=0\). The Teichm\"uller set is not closed under addition,
since \(1+1=2\notin\mathcal T(A)\).
\end{example}

\subsection{Roots-of-Unity Orbits, Factor Lifting, and Quotient Ring Decompositions}

Assume \(p\nmid n\), put \(L=\ord_n(q)\), and choose in
\(S=\operatorname{GR}(p^s,rL)\) a Teichm\"uller element \(\xi\) of order exactly \(n\).
Such an element exists because the Teichm\"uller unit group of \(S\) is cyclic of order \(q^L-1\).
For a \(q\)-cyclotomic coset \(C\) modulo \(n\), put
\[
 F_C(X)=\prod_{b\in C}(X-\xi^b),\qquad
 M_C(X)=\prod_{b\in C}(X-\overline\xi^{b}).
\]

\begin{theorem}[Roots-of-unity factors over Galois rings]
\label{v47:thm-galois-cyclotomic-factorization}
Each \(F_C\) lies in \(A[X]\), is basic irreducible, and satisfies \(\deg F_C=|C|\) and \(\overline{F_C}=M_C\).
These factors are pairwise comaximal and give
\begin{equation}\label{v47:eq-galois-factorization}
 X^n-1=\prod_C F_C(X),\qquad
 A[X]/(X^n-1)\cong\prod_C A[X]/(F_C).
\end{equation}
Each quotient ring \(A[X]/(F_C)\) is isomorphic to \(\operatorname{GR}(p^s,r|C|)\).
\end{theorem}
\begin{proof}
Relative Frobenius satisfies \(\sigma(\xi)=\xi^q\), so it permutes all linear factors of \(F_C\).
The coefficients are fixed by \(\sigma\) and hence belong to \(A\). The reduced roots form a complete finite-field Frobenius orbit.
By Theorem~\ref{thm:orbit}, \(M_C\) is irreducible of degree \(|C|\), so \(F_C\) is basic irreducible.

Distinct \(\xi^b\) have distinct reductions, so their differences are units; hence the factors \(X-\xi^b\) are pairwise comaximal.
Each linear factor divides \(X^n-1\), and the intersection of comaximal ideals equals their product.
Thus the product of all \(n\) linear factors divides \(X^n-1\). Both are monic of the same degree, so they are equal.
Grouping by orbits gives the first identity. Distinct \(M_C\) are coprime;
the lifted B\'ezout argument of Theorem~\ref{v47:thm-hensel-coprime} shows that distinct \(F_C\) are comaximal,
and the Chinese remainder theorem gives the direct-product isomorphism.

Let \(d=|C|\) and \(a\in C\). Since \(d\mid L\), the extension-ring construction above embeds
\(\operatorname{GR}(p^s,rd)\) in \(S\). Its image is fixed by \(\sigma^d\),
and the digit formula shows that the fixed subring of \(\sigma^d\) has exactly \(q^{ds}\) elements, so the two coincide.
Denote this subring by \(S_d\). Then \(\xi^a\in S_d\).
Evaluation gives \(A[X]/(F_C)\to S_d\). Both are free \(A\)-modules of rank \(d\):
for \(S_d\), lift an \(\F_q\)-basis of the residue field and successively eliminate \(p\)-adic digits to obtain a generating set.
If a nontrivial linear relation exists, extract the least \(p\)-adic valuation \(j<s\) among its coefficients;
the remaining combination then lies in \(\operatorname{Ann}_{S_d}(p^j)=p^{s-j}S_d\subseteq pS_d\).
Reduction contradicts linear independence of the residue-field basis.
Modulo \(p\), evaluation is the isomorphism \(\F_q[X]/(M_C)\to\F_{q^d}\),
so its matrix determinant is a unit of \(A\), and evaluation itself is an isomorphism.
\end{proof}

By Theorem~\ref{v47:thm-hensel-coprime}, \(F_C\) is also the unique monic factor of \(X^n-1\) reducing to \(M_C\).
Thus Frobenius orbits determine factor labels and degrees, while Hensel uniqueness determines coefficients at each precision.
For \(h\mid n\), the polynomial \(\Phi_h\) still has
\(\varphi(h)/\ord_h(q)\) basic irreducible factors, each of degree \(\ord_h(q)\):
they are the orbit factors of roots of order \(h\), and equality of their root product with \(\Phi_h\) also follows from uniqueness of lifting.
On extending the base ring to \(\operatorname{GR}(p^s,rt)\), the Frobenius step changes from \(q\) to \(q^t\),
so a factor of degree \(d\) splits into \(\gcd(d,t)\) basic irreducible factors of degree \(d/\gcd(d,t)\).

\begin{proposition}[Monomial composition of root blocks and Hensel lifting]\label{v47:prop-galois-med}
Let \(p\nmid n\), let \(\varnothing\ne T\subseteq\Z/n\Z\) satisfy \(qT=T\), and put
\[
P_T(X)=\prod_{a\in T}(X-\overline\xi^{a}),\qquad
\kappa_T=|H_+(T)|,\qquad D=n/\kappa_T.
\]
Let \(T_0\) be the image of \(T\) modulo \(D\), and put
\[
Q_{T_0}(Y)=\prod_{b\in T_0}(Y-\overline\xi^{\kappa_Tb}).
\]
Let \(F_T\in A[X]\) and \(G_{T_0}\in A[Y]\) be the unique monic divisor lifts of \(P_T\) and \(Q_{T_0}\) in \(X^n-1\) and \(Y^D-1\), respectively. For every positive integer \(h\),
\begin{equation}\label{v124:eq-root-block-lift}
F_T(X)=G_{T_0}(X^{\kappa_T}),\qquad
F_T\in A[X^h]\quad\Longleftrightarrow\quad h\mid\kappa_T.
\end{equation}
For every \(1\le j\le s\), the greatest common divisor of all positive exponents of nonzero terms in \(F_T\bmod p^j\) is still \(\kappa_T\).

If \(T=C\) is a single cyclotomic coset, put \(\tau=|C|\) and \(\omega_C=\tau/\kappa_C\). Then \(H_C=G_{T_0}\) is basic irreducible of degree \(\omega_C\), and
\begin{equation}\label{v47:eq-galois-med}
F_C(X)=H_C(X^{\kappa_C}).
\end{equation}
The factor \(F_C\) is a binomial if and only if \(C\) is a complete equal-difference block; in that case, for \(a\in C\), one has \(F_C(X)=X^\tau-\xi^{a\tau}\).
\end{proposition}
\begin{proof}
The set \(T\) is a union of complete cosets of \(D\Z/n\Z\), so Proposition~\ref{prop:ed-binomial-correspondence} gives
\[
P_T(X)=Q_{T_0}(X^{\kappa_T}).
\]
Reduction commutes with multiplication by \(q\), so \(Q_{T_0}\in\F_q[Y]\); it is square-free and divides \(Y^D-1\). Since \(p\nmid n\), each selected factor is coprime to its complementary factor. Theorem~\ref{v47:thm-hensel-coprime} therefore gives unique lifts.

Write \(Y^D-1=G_{T_0}(Y)K(Y)\). Substitution of \(Y=X^{\kappa_T}\) yields
\[
X^n-1=G_{T_0}(X^{\kappa_T})K(X^{\kappa_T}).
\]
The first factor is monic and reduces to \(P_T\). Uniqueness of the lift proves the first identity in \eqref{v124:eq-root-block-lift}, and hence \(F_T\in A[X^h]\) whenever \(h\mid\kappa_T\). Conversely, \(F_T\in A[X^h]\) implies \(P_T\in\F_q[X^h]\) after reduction, so Proposition~\ref{v40:prop:root-power-composition} gives \(h\mid\kappa_T\). The same argument applies at every precision \(p^j\); the gcd of the nonzero positive coefficient exponents characterizes the maximal composition exponent.

If \(T=C\), its image modulo \(D\) is again a single Frobenius orbit, of length \(\tau/\kappa_C=\omega_C\). Thus \(Q_{T_0}\) is irreducible and \(H_C\) is basic irreducible. The root product of a complete equal-difference block is \(X^\tau-\xi^{a\tau}\); unit differences between distinct roots ensure that this product identity remains valid in the splitting Galois ring. Conversely, \(F_C\) has unit constant term. If it is a binomial, its reduction remains a binomial, so \(C\) is a complete equal-difference block.
\end{proof}

Here \(F_T\) is the canonical divisor lift in \(X^n-1\). An arbitrary coefficient lift of \(P_T\) may acquire new nonzero coefficients divisible by \(p\), destroying this composition structure.

\begin{proposition}[Composition structure over finite local rings]\label{v125:prop-chain-composition}
Let \(R\) be a finite commutative local ring with maximal ideal \(\mathfrak m\) and residue field \(k=\F_q\), and suppose \(p\nmid n\).
Retain \(T,P_T,Q_{T_0},\kappa_T,D\) from Proposition~\ref{v47:prop-galois-med}.
There are unique monic factors \(F_T,G_{T_0}\) of \(X^n-1,Y^D-1\) over \(R\) reducing to \(P_T,Q_{T_0}\), respectively. For every positive integer \(h\),
\[
F_T(X)=G_{T_0}(X^{\kappa_T}),\qquad
F_T\in R[X^h]\quad\Longleftrightarrow\quad h\mid\kappa_T.
\]
These canonical factors commute with reduction \(R\to R/\mathfrak m^j\); for every \(j\ge1\), the maximum monomial composition exponent of \(F_T\bmod\mathfrak m^j\) is \(\kappa_T\).
\end{proposition}
\begin{proof}
First, \(\mathfrak m\) is nilpotent. For \(x\in\mathfrak m\), finiteness gives \(x^a=x^{a+b}\) for some \(a,b\ge1\). Since \(1-x^b\) is a unit, \(x^a=0\). Choose generators \(\mathfrak m=(x_1,\ldots,x_t)\). If \(x_i^{e_i}=0\), then
\[
\mathfrak m^{\,1+\sum_i(e_i-1)}=0.
\]
The case \(\mathfrak m=0\) is immediate.

For \(j\ge1\), \(V_j=\mathfrak m^j/\mathfrak m^{j+1}\) is a finite-dimensional \(k\)-space, not necessarily of dimension one.
Let a monic polynomial \(f\in R[X]\) satisfy \(\overline f=gh\), where \(g,h\in k[X]\) are monic and coprime, of degrees \(d,e\).
Starting with a monic factorization \(f\equiv G_jH_j\pmod{\mathfrak m^j}\), its error gives \(E_j\in V_j[X]_{<d+e}\).
The correction is governed by
\[
V_j[X]_{<d}\oplus V_j[X]_{<e}\longrightarrow V_j[X]_{<d+e},
\qquad (U,V)\longmapsto hU+gV,
\]
where \(V_j[X]_{<d}\) denotes polynomials of degree less than \(d\) with coefficients in \(V_j\).
In any \(k\)-basis of \(V_j\), this map is a direct sum of copies of the same scalar map.
The scalar map has zero kernel: \(hU+gV=0\) and \(\gcd(g,h)=1\) imply \(g\mid U\), and the degree bounds give \(U=V=0\). Equality of dimensions makes it an isomorphism.
Thus \(E_j=hU+gV\) has a unique solution. Lift \(U,V\) arbitrarily to \(\mathfrak m^j[X]\) and add them to \(G_j,H_j\). Since \(\mathfrak m^{2j}\subseteq\mathfrak m^{j+1}\), this gives the unique monic factorization at the next precision.
Iteration until \(\mathfrak m^j=0\) proves existence and uniqueness of coprime monic factor lifts. This is also the general Hensel property of nilpotent ideals~\cite{StacksHensel}.

Since \(p\nmid n\), each of \(P_T,Q_{T_0}\) is coprime to its complementary factor. Substitute \(Y=X^{\kappa_T}\) into the canonical lifted factorization of \(Y^D-1\). Uniqueness gives \(F_T=G_{T_0}(X^{\kappa_T})\).
Hence \(h\mid\kappa_T\) implies \(F_T\in R[X^h]\); conversely, reduction to \(k\) and Proposition~\ref{v40:prop:root-power-composition} give \(h\mid\kappa_T\).
The reduction of a canonical factor is still a monic lift of the prescribed residue factor, so uniqueness gives compatibility at every precision. The same two implications apply over each \(R/\mathfrak m^j\).
\end{proof}

For example,
\[
R=\F_3[\varepsilon_1,\varepsilon_2]/(\varepsilon_1,\varepsilon_2)^2
\]
has maximal ideal \(\mathfrak m=(\varepsilon_1,\varepsilon_2)\) with \(\mathfrak m^2=0\). The ideals \((\varepsilon_1)\) and \((\varepsilon_2)\) are incomparable, so \(R\) is not a chain ring.
Here \(\dim_{\F_3}\mathfrak m=2\), and the correction acts independently in the two coefficient directions.
Since \(\F_3\subset R\), a coprime factorization of \(X^n-1\) over \(\F_3\) is already its canonical lift. In particular, \(X^8+X^2+2\) from Example~\ref{v125:ex-160-composition} remains the canonical factor of \(X^{160}-1\), with maximum composition exponent two.
The map \(z\mapsto z^3\) annihilates \(\mathfrak m\) and is not an automorphism; also \(\mathfrak m\ne3R=0\). Composition preservation uses only coprimeness of the residue factors and uniqueness of nilpotent lifting. The Galois-ring digit expansion and Frobenius description retain their respective hypotheses.

\subsection{Local Corrections Determined by the Derivative}

\begin{proposition}[Local lifting of roots-of-unity factors]
\label{v47:prop-galois-local-lift}
Let \(p\nmid n\), \(1\le j<s\), and let \(G_j\in A[X]\) be monic of positive degree and divide \(X^n-1\) modulo \(p^j\).
If \(X^n-1=Q_jG_j\pmod{p^j}\), then
\begin{equation}\label{v47:eq-galois-local-inverse}
 Q_j^{-1}=n^{-1}XG_j'(X)\quad\text{in }(A/p^jA)[X]/(G_j).
\end{equation}
Choose any monic coefficient lift of \(G_j\), and compute modulo \(p^{j+1}\)
\[
 R_j=(X^n-1)\bmod G_j,\qquad
 E_j=R_j/p^j\bmod p,\qquad
 \Delta_j=E_jn^{-1}X\overline{G_j}'\bmod\overline{G_j}.
\]
Take \(\Delta_j\) as the remainder of degree less than \(\deg G_j\), and lift its coefficients arbitrarily. Then
\[
 G_{j+1}=G_j+p^j\Delta_j\pmod{p^{j+1}}
\]
is the unique monic factor lift of the same degree.
\end{proposition}
\begin{proof}
Differentiate \(X^n-1=Q_jG_j\pmod{p^j}\) and reduce modulo \(G_j\) to obtain
\(Q_jG_j'=nX^{n-1}\). In the quotient ring, \(X^n=1\) and \(n\) is a unit;
multiplication by \(n^{-1}X\) gives \eqref{v47:eq-galois-local-inverse}.
Monic division gives \(X^n-1=QG_j+p^jE_j\pmod{p^{j+1}}\).
The correction satisfies \(\overline Q\Delta_j=E_j\pmod{\overline{G_j}}\),
so there is a \(\Theta\in\F_q[X]\) with \(E_j=\overline Q\Delta_j+\overline{G_j}\Theta\).
Choose any coefficient lift of \(\Theta\), retaining its notation. Then
\[
 (Q+p^j\Theta)(G_j+p^j\Delta_j)
 =X^n-1\pmod{p^{j+1}},
\]
since \(2j\ge j+1\). Uniqueness follows from Theorem~\ref{v47:thm-hensel-coprime};
the condition \(p\nmid n\) ensures that the reduced factor is coprime to the complementary factor.
\end{proof}

\begin{example}[\(X^8-1\) over \(\Z/9\Z\)]
\label{v47:ex-z9-eight}
Over \(\F_3\), we have
\[
 X^8-1=(X-1)(X+1)(X^2+1)(X^2+X+2)(X^2+2X+2).
\]
All three quadratic factors have nonsquare discriminant \(2\). For \(X^4+1\), choose integer-coefficient initial lifts
\(G=X^2+X+2\) and \(H=X^2+2X+2\). Then
\[
 \frac{X^4+1-GH}{3}=2X^3+X^2+X+2\pmod3.
\]
Put \(g=\overline G\) and \(h=\overline H\). In the two-factor correction equation, take \(U=V=X+2\), giving \(hU+gV=2X^3+X^2+X+2\).
Their lifts modulo \(9\) are therefore \(G+3U=X^2+4X+8\) and \(H+3V=X^2+5X+8\), so
\[
 \begin{aligned}
 X^8-1={}&(X-1)(X+1)(X^2+1)\\
 &\cdot(X^2+4X+8)(X^2+5X+8)
 \quad\text{in }(\Z/9\Z)[X].
 \end{aligned}
\]
Each factor is basic irreducible; its degree is unchanged from \(\F_3\), while lifting changes its coefficients.
\end{example}

\begin{example}[Coefficient precision along the same factor tower]\label{v125:ex-160-hensel}
Let \(A_s=\Z/3^s\Z\), and let \(\widetilde G_s\) be the canonical monic factor of \(Y^{80}-1\) reducing to \(G=Y^4+Y+2\).
Every \(F_a\) in~\eqref{v125:eq-160-tower} has a canonical lift, and
\[
\widetilde F_{a,s}(X)=\widetilde G_s(X^{2^{a-4}})
\qquad(a\ge4,\ s\ge1).
\]
Indeed, \(\widetilde G_s\mid Y^{80}-1\), so the right side divides \(X^{5\cdot2^a}-1\) and reduces to \(F_a\); uniqueness gives the equality. Consequently,
\[
\widetilde F_{a+1,s}(X)=\widetilde F_{a,s}(X^2),\qquad
\widetilde F_{a,s+1}\bmod3^s=\widetilde F_{a,s}.
\]
Lifting the modulus commutes with increasing coefficient precision. The maximum composition exponent remains \(2^{a-4}\) at every precision.

At precision nine, use the same integer polynomial \(G\) as initial lift. Monic division gives
\[
(Y^{80}-1)\bmod G=3+6Y^2+6Y^3\pmod9.
\]
The correction in Proposition~\ref{v47:prop-galois-local-lift} is
\[
\Delta=(1+2Y^2+2Y^3)\,80^{-1}Y\overline G'\bmod\overline G
=Y^3+Y^2+2.
\]
Hence
\[
\widetilde G_2=Y^4+3Y^3+3Y^2+Y+8,\qquad
\widetilde F_{5,2}=X^8+3X^6+3X^4+X^2+8.
\]
The two-element root blocks modulo 160 still force every odd-degree coefficient to vanish, whereas the \(X^6,X^4\) coefficients change from zero modulo three to nonzero values.
Thus the maximum composition exponent persists at every precision, but additional coefficient zeros need not persist.
\end{example}

\subsection{Comaximal Blocks When the Residue Characteristic Divides the Exponent}

Now let \(n=p^am\), where \(a\ge1\) and \(p\nmid m\).
Choose a primitive \(m\)-th Teichm\"uller root of unity \(\xi\) in a suitable Galois extension ring, and define \(F_C\) by the orbit product above.
For a \(q\)-cyclotomic coset \(C\) modulo \(m\), put \(p^aC=\{p^ab\bmod m:b\in C\}\).
Multiplication by \(p^a\) permutes these cosets and preserves their sizes.

\begin{proposition}[Canonical comaximal blocks]
\label{v47:prop-galois-ramified-blocks}
Let \(Y^m-1=\prod_C F_C(Y)\) be the basic irreducible factorization, and define
\[
 B_C(X)=F_{p^aC}(X^{p^a}).
\]
Then
\begin{equation}\label{v47:eq-galois-ramified}
 X^n-1=\prod_C B_C(X),\qquad
 \overline{B_C}=M_C^{p^a}.
\end{equation}
The \(B_C\) are pairwise comaximal, and each quotient \(A[X]/(B_C)\) is a local ring.
They are the unique monic factors of \(X^n-1\) reducing to \(M_C^{p^a}\).
They are called the \textbf{canonical comaximal blocks} and give the corresponding Chinese remainder decomposition.
\end{proposition}
\begin{proof}
The permutation of cosets and substitution \(Y=X^{p^a}\) give the product identity.
Making the same substitution in a B\'ezout identity for \(F_C,F_D\) proves comaximality of the blocks.
After reduction to the splitting field, taking powers in characteristic \(p\) gives
\[
 \overline{F_{p^aC}}(X^{p^a})
 =\prod_{b\in C}(X^{p^a}-\overline\xi^{p^ab})
 =\prod_{b\in C}(X-\overline\xi^{b})^{p^a}=M_C^{p^a}.
\]
The polynomials \(M_C^{p^a}\) are pairwise coprime, so uniqueness of the lifts follows from Theorem~\ref{v47:thm-hensel-coprime}.
Modulo \(p\), the quotient ring is \(\F_q[X]/(M_C^{p^a})\), where every element coprime to \(M_C\) is invertible,
so it has a unique maximal ideal. Since \(p\) is nilpotent, every maximal ideal of the original quotient contains \(p\),
and hence the original quotient is also local.
\end{proof}

A block cannot split further into two comaximal monic factors of positive degree: its quotient ring would then be a direct product of two nonzero rings,
contradicting locality. A block may still admit a factorization into noncomaximal factors.
When \(s>1\), one generally cannot write \(B_C=F_C(X)^{p^a}\);
the power identity in \eqref{v47:eq-galois-ramified} holds only after reduction modulo \(p\).
When \(r=1\), multiplication by \(p^a\) preserves every \(p\)-cyclotomic coset, so \(B_C=F_C(X^{p^a})\).

\begin{example}[The coset twist in residue degree two]\label{v104:ex-gr-four-twist}
Take \(R=\operatorname{GR}(4,2)=(\Z/4\Z)[\omega]/(\omega^2+\omega+1)\).
Here \(\omega^3=1\), and the residue field is \(\F_4\). For \(X^6-1\), we have \(p^a=2\) and \(m=3\).
The four-cyclotomic cosets modulo three are \(C_0=\{0\},C_1=\{1\},C_2=\{2\}\), and multiplication by two exchanges \(C_1,C_2\).
Since \(F_{C_j}(Y)=Y-\omega^j\), the canonical blocks are
\[
\begin{aligned}
B_{C_0}&=X^2-1,&
B_{C_1}&=X^2-\omega^2,&
B_{C_2}&=X^2-\omega,\\
X^6-1&=(X^2-1)(X^2-\omega^2)(X^2-\omega).
\end{aligned}
\]
Every difference of two blocks is a unit, since their constant terms have distinct residues. Thus the three blocks are pairwise comaximal.
The twist ensures the correct labeling:
\[
\overline{B_{C_1}}=(X-\overline\omega)^2,
\qquad
\overline{F_{C_1}(X^2)}=(X-\overline\omega^{\,2})^2.
\]
Omitting \(2C_1=C_2\) assigns the block for \(C_2\) to \(C_1\), although the product of all blocks remains unchanged.
Also, \(B_{C_1}=(X-\omega)(X+\omega)\); these two linear factors have the same reduction modulo two and are not comaximal. Moreover,
\[
B_{C_1}-(X-\omega)^2=2\omega(X-\omega)\ne0.
\]
Thus this smallest nontrivial residue-field extension displays the coset twist, Chinese remainder decomposition into blocks, and repeated-root factorization within a block at once.
\end{example}

\begin{example}[\(X^6-1\) over \(\Z/9\Z\)]
\label{v47:ex-z9-six}
From \(Y^2-1=(Y-1)(Y+1)\), we obtain the canonical block decomposition
\[
 X^6-1=(X^3-1)(X^3+1).
\]
The blocks differ by the unit \(2\), so they are comaximal; modulo \(3\), they become \((X-1)^3\) and \((X+1)^3\), respectively.
However, in \((\Z/9\Z)[X]\),
\[
 X^3-1-(X-1)^3=3X(X-1)\ne0.
\]
Inside a block, we still have \(X^3-1=(X-1)(X^2+X+1)\), but these two factors are not comaximal:
their reductions modulo \(3\) are \(X-1\) and \((X-1)^2\).
Thus Chinese remainder decomposition into repeated-root blocks and factorization within a block have different structures.
\end{example}

\medskip
\noindent\textbf{From root orbits to coefficient precision.}
When \(p\nmid n\), each finite-field root block has a unique canonical factor lift.
The additive stabilizer determines its monomial composition structure, and Hensel uniqueness preserves that structure at every coefficient precision.
When \(n=p^am\), \(p\nmid m\), the same uniqueness applies to the coprime primary blocks over the residue field, rather than to individual repeated factors inside a block.

\begingroup\small
\renewcommand{\arraystretch}{1.18}
\begin{longtable}{@{}>{\raggedright\arraybackslash}p{0.20\textwidth}>{\raggedright\arraybackslash}p{0.34\textwidth}>{\raggedright\arraybackslash}p{0.38\textwidth}@{}}
\toprule
Structure&Finite field \(\F_q\)&Galois ring \(A=\GR(p^s,r)\)\\
\midrule
\endhead
\bottomrule
\endfoot
Roots and Frobenius&For \(p\nmid n\), a \(q\)-orbit \(C\)&Teichmüller roots have the same exponent orbit; relative Frobenius is the \(q\)-th power on these roots\\
Factors and degrees&\(M_C\) irreducible of degree \(|C|\)&Its canonical lift \(F_C\) is basic irreducible of the same degree\\
Composition&\(P_T=Q_{T_0}(X^{\kappa_T})\)&\(F_T=G_{T_0}(X^{\kappa_T})\); the maximum composition exponent is \(\kappa_T\) at every precision\\
Additional zero coefficients&Determined by symmetric functions of compressed roots&May become nonzero multiples of \(p\); see Example~\ref{v125:ex-160-hensel}\\
Stable factor towers&Unique orbit lifts give power substitution&Power substitution commutes with reduction of coefficient precision\\
\(n=p^am\)&Distinct factors \(M_C\), each of multiplicity \(p^a\)&Canonical comaximal blocks \(F_{p^aC}(X^{p^a})\), generally different from \(F_C^{p^a}\)\\
\end{longtable}
\endgroup
Composition preservation also holds over every finite commutative local ring, by Proposition~\ref{v125:prop-chain-composition}. Thus coprime exponents preserve root orbits and all their composition exponents. With a characteristic-power factor in the exponent, the preserved structure is the comaximal decomposition between distinct residual root orbits; further factorization inside a block requires separate analysis.

\par\Needspace{11\baselineskip}

\section{From finite precision to \texorpdfstring{$p$}{p}-adic lifts: root orbits and coefficient valuations}
\label{v132:sec-adic-lift}
The preceding discussion lifted roots and factors at every fixed precision $\operatorname{GR}(p^s,r)$. We now assemble the compatible levels. The point is not to replace the orbit structure already determined over the finite field, but to distinguish a coefficient that vanishes identically from a nonzero coefficient that vanishes upon reduction modulo $p$.

\begin{theorem}[Unramified integers and $p$-adic digits]
\label{v132:thm-unramified-limit}
Let $q=p^r$. Choose a monic irreducible polynomial $\bar F\in\F_p[X]$ of degree $r$ and fix a monic coefficientwise lift $F\in\Z[X]$. Put $A_s=(\Z/p^s\Z)[X]/(F\bmod p^s)\simeq\operatorname{GR}(p^s,r)$. Reduction yields a complete ring
\[
 \Z_q:=\varprojlim_s A_s,\qquad \Z_q/p^s\Z_q\simeq A_s,
 \qquad \Z_q/p\Z_q\simeq\F_q.
\]
It is a complete discrete valuation ring with uniformizer $p$; write $\mathbb Q_q$ for its fraction field. Every $a\in\F_q$ has a unique \textbf{Teichmüller lift} $[a]\in\Z_q$ with $[a]^q=[a]$ and residue $a$. Each $x\in\Z_q$ has a unique convergent expansion
\begin{equation}\label{v132:eq-adic-expansion}
 x=\sum_{j\ge0}p^j[a_j],\qquad a_j\in\F_q.
\end{equation}
The compatible absolute Frobenius automorphisms induce
\[
 \phi(x)=\sum_{j\ge0}p^j[a_j^p],\qquad \phi^r=\operatorname{id}.
\]
For $x\ne0$, define $v_p(x)$ as the largest $s$ for which $x\in p^s\Z_q$, and let $v_p(0)=+\infty$. Then $v_p(p)=1$ and the extension of this valuation to $\mathbb Q_q$ satisfies the multiplicative and ultrametric laws.
\end{theorem}
\begin{proof}
Coefficient reduction of the fixed defining equation supplies compatible quotient rings and surjective reduction maps. At each level the digit expansion is unique, so a compatible sequence can be constructed, and uniquely recovered, one digit at a time. If its first nonzero digit has index $j$, then $x=p^ju$ with $u$ having nonzero residue and hence invertible. Consequently every nonzero ideal is $(p^j)$. The inverse limit is complete and separated because $\bigcap_{j\ge0}p^j\Z_q=0$. Products and valuations are controlled by the first nonzero digit. The existence and uniqueness of Teichmüller representatives and the Frobenius action follow level by level from Theorem~\ref{v47:thm-galois-frobenius}, and commute with reduction.
\end{proof}

Teichmüller lifting is multiplicative, $[ab]=[a][b]$, but not generally additive. For instance, in $\Z_2$ one has $[1]+[1]=2$ whereas $[1+1]=0$. Thus an additive identity over the residue field need not lift to an identity in characteristic zero.

\subsection{Additive carries and Witt vectors}
\label{v133:subsec-witt-carry}
Put $K=\F_q$, $A=\Z_q$, and let $\pi:A\to K$ denote reduction.
There is no unital ring homomorphism $s:K\to A$ with $\pi\circ s=\id$.
Indeed, additivity would give $ps(1)=s(p\cdot1)=0$, whereas $A$ has no $p$-torsion,
so $s(1)=0$. The obstacle therefore comes from mixed characteristic, rather than
from a defect in the choice of Teichmüller representatives.

Although $a\mapsto[a]$ is not additive, it determines both ring operations on $A$.
Write $x=\sum_{j\ge0}p^j[a_j]$ as its sequence of digits $(a_0,a_1,\ldots)$.
Sums and products of such sequences are uniquely induced by the operations in $A$;
they involve corrections depending on preceding digits, not coordinatewise addition or multiplication.
This realizes the ring of \textbf{$p$-typical Witt vectors} over the perfect field $K$,
written $W(K)\simeq\Z_q$; truncation to $s$ digits gives
$W_s(K)\simeq\GR(p^s,r)$.
The digits $a_j$ and the standard \textbf{Witt coordinates} $x_j$ are related by
\begin{equation}\label{v136:eq-witt-coordinates}
 x_j=a_j^{p^j},\qquad
 x=\sum_{j\ge0}p^j[x_j^{p^{-j}}].
\end{equation}
Here $p^{-j}$ denotes the $j$th inverse iterate of Frobenius on the perfect field $K$; the standard coordinates use the Witt polynomials
$w_j(X_0,\ldots,X_j)=\sum_{i=0}^j p^iX_i^{p^{j-i}}$.
Thus the \emph{Witt ring itself} has both ring operations, but the Teichmüller map
from $K$ into it does not preserve the original addition and multiplication simultaneously.

\begin{proposition}[Polynomial formula for the first carry]\label{v133:prop-first-carry}
Let $q=p^r$ and put
\[
 C_p(X,Y)=\frac{(X+Y)^p-X^p-Y^p}{p}
 =\sum_{i=1}^{p-1}\frac1p\binom pi X^iY^{p-i}\in\Z[X,Y].
\]
Write $\overline C_p$ for coefficientwise reduction modulo $p$, and define
\[
 \kappa(a,b)=\bigl(-\overline C_p(a,b)\bigr)^{p^{r-1}}
 \quad(a,b\in K).
\]
Then
\begin{equation}\label{v133:eq-first-carry}
 [a]+[b]\equiv[a+b]+p[\kappa(a,b)]\pmod{p^2A}.
\end{equation}
In particular, for $p=2$ one has $\kappa(a,b)=(ab)^{2^{r-1}}$, and hence
$[1]+[1]\equiv2[1]\pmod4$.
\end{proposition}
\begin{proof}
In the complete unramified integer ring, write $[a]+[b]=[a+b]+p\delta$. Taking $p$th powers modulo $p^2$ gives
\[
 [a^p]+[b^p]+pC_p([a],[b])\equiv[(a+b)^p].
\]
Applying arithmetic Frobenius to the original equality instead gives
$[a^p]+[b^p]=[(a+b)^p]+p\phi(\delta)$.
Comparing, dividing by $p$, and reducing yields $\bar\delta^p=-\overline C_p(a,b)$.
The inverse of $p$th-power Frobenius on $K$ is $z\mapsto z^{p^{r-1}}$, so $\bar\delta=\kappa(a,b)$, proving the congruence. Reduction gives the same identity in every finite-precision Galois ring.
\end{proof}

Theorem~\ref{v140:thm-teich-rigidity} recovers multiplicative rigidity from all first carries, and Theorem~\ref{v140:thm-digit-subspaces} classifies linear subspaces whose repeated digit sets remain additively closed.

\begin{proposition}[Obstructions to lifting zero-sum relations]\label{v133:prop-zero-sum-obstruction}
Suppose $a_1,\ldots,a_m\in K$ satisfy $\sum_i a_i=0$. Then
$S=\sum_i[a_i]\in pA$, and
\[
 \mathcal O_1(a_1,\ldots,a_m)=\overline{S/p}\in K
\]
is independent of the placement of parentheses. Moreover,
\[
 \mathcal O_1=0\quad\Longleftrightarrow\quad S\in p^2A.
\]
Whenever $S\in p^jA$, the $j$th-digit obstruction $\overline{S/p^j}$ is defined;
this requires all preceding obstructions to vanish.
\end{proposition}
\begin{proof}
Reduction gives $\pi(S)=\sum_i a_i=0$, so $S\in pA$.
As $A$ has no $p$-torsion, $S/p$ is unique; its reduction is zero exactly when
$S/p\in pA$. The higher-digit assertion is identical, while independence
of parentheses follows from associativity in $A$.
\end{proof}

For an extension $E=\F_{q^d}$, norm commutes exactly with Teichmüller lifting:
\begin{equation}\label{v133:eq-norm-exact}
 \Nm_{\Z_{q^d}/\Z_q}([a])=[\Nm_{E/K}(a)].
\end{equation}
Trace, however, is only guaranteed to commute with reduction.
When $\Tr_{E/K}(a)=0$, we may define
\[
 \eta_1(a)=\overline{\Tr_{\Z_{q^d}/\Z_q}([a])/p}\in K.
\]
For example, with $E=\F_4$ and $K=\F_2$, one has $\Tr_{E/K}(1)=0$,
but its lifted trace equals $2$, so $\eta_1(1)=1$;
by contrast, the trace of the lifted zero element vanishes identically.
For general zero-sum identities, carry digits must satisfy compatibility
constraints across different parenthesizations; an exact identity for the first digit
is derived in the exercises.

\subsection{Ring sections and formal deformations in equal characteristic}
\label{v133:subsec-equichar}
An algebraic section impossible in mixed characteristic can exist in equal characteristic.
Let $R=K[[t]]$. The constant embedding $\iota:K\hookrightarrow R$ preserves
both addition and multiplication and is a section of $R\to R/(t)=K$.
It is the unique unital ring section: if $s:K\to R$ is such a section,
then $s(a)^q=s(a)$; however, the only roots of $z^q=z$ in $R$ are constants.
Indeed, if $z=a_0+u$ with $u\in tR\setminus\{0\}$, then
$v_t(u^q)=qv_t(u)>v_t(u)$, so $u^q=u$ is impossible.
This preserves the operations, but constant embedding alone creates
no new elements.

A more substantive construction lets \emph{equations and their roots} vary with $t$,
while retaining the field operations on the coefficients.
\begin{proposition}[Formal Hensel lifting in equal characteristic]\label{v133:prop-equal-hensel}
Let $f\in K[X]$ and let $\alpha\in E=\F_{q^d}$ satisfy
$f(\alpha)=0$ and $f'(\alpha)\ne0$. If
$F(X,t)\in K[[t]][X]$ satisfies $F(X,0)=f(X)$, there is a unique
$\alpha(t)\in E[[t]]$ such that
\[
 \alpha(t)\equiv\alpha\pmod t,\qquad
 F(\alpha(t),t)=0.
\]
\end{proposition}
\begin{proof}
Suppose $u_j\in E[t]$ satisfies $u_j\equiv\alpha\pmod t$ and
$F(u_j,t)\equiv0\pmod{t^j}$. Write
$F(u_j,t)\equiv b_jt^j\pmod{t^{j+1}}$ and set
\[
 u_{j+1}=u_j-b_jf'(\alpha)^{-1}t^j.
\]
Formal Taylor expansion gives $F(u_{j+1},t)\equiv0\pmod{t^{j+1}}$.
The same calculation proves uniqueness of the next coefficient.
Completeness assembles the compatible truncations into a unique $\alpha(t)$.
\end{proof}

\begin{example}[A nontrivial deformation of an Artin--Schreier equation]
\label{v133:ex-artin-schreier-deformation}
Choose $c\in K$ with $\Tr_{K/\F_p}(c)\ne0$.
Then $f(X)=X^p-X-c$ is irreducible over $K$.
Let $\alpha\in E=\F_{q^p}$ be a root and put
\[
 B(t)=-\sum_{j\ge0}t^{p^j}\in K[[t]],\qquad
 \alpha(t)=\alpha+B(t).
\]
Since $B(t)^p-B(t)=t$, the identity
\begin{equation}\label{v133:eq-AS-deform}
 \alpha(t)^p-\alpha(t)=c+t
\end{equation}
holds exactly. The conjugates of $\alpha(t)$ in the constant extension
$E[[t]]/K[[t]]$ are $\alpha(t)+j$ for $j\in\F_p$, whence
\[
 \Tr_{E[[t]]/K[[t]]}(\alpha(t))=\Tr_{E/K}(\alpha),
 \qquad
 \Nm_{E[[t]]/K[[t]]}(\alpha(t))=c+t.
\]
The trace remains constant while the norm varies with $t$.
It is essential to distinguish the ring homomorphism $K\hookrightarrow K[[t]]$
from the assignment $\alpha\mapsto\alpha(t)$, which deforms a root of a chosen
equation and is \emph{not} a previously defined ring homomorphism on $E$.
\end{example}

The two constructions preserve different information. Witt carries measure
divisibility digits concealed by reduction modulo $p$; equal-characteristic deformations
preserve coefficient-field operations and use $t$-adic precision to follow the roots
of varying equations. The former serves lifted coefficient valuations, traces and Gauss
periods; the latter serves deformations of separable roots, Frobenius conjugates, and
parameter-dependent trace and norm.

\begin{proposition}[Compatible lifts of trace, norm, and cyclotomic factors]
\label{v132:prop-coset-adic-factor}
Let $E=\F_{q^d}$ and $\Z_{q^d}=\varprojlim_s\operatorname{GR}(p^s,rd)$. Its relative Frobenius is $\phi^r$, and the trace and norm of the ring extension $\Z_{q^d}/\Z_q$ satisfy
\[
 \operatorname{Tr}_{\Z_{q^d}/\Z_q}(x)=\sum_{i=0}^{d-1}\phi^{ri}(x),\quad
 \operatorname{Nm}_{\Z_{q^d}/\Z_q}(x)=\prod_{i=0}^{d-1}\phi^{ri}(x).
\]
Their reductions are the usual trace and norm for $E/\F_q$. If $p\nmid N$ and $N\mid q^d-1$, choose a primitive $N$th Teichmüller root of unity $\widetilde\zeta$. For a $q$-cyclotomic coset $C$ modulo $N$,
\[
 \widetilde M_C(X)=\prod_{i\in C}(X-\widetilde\zeta^{\,i})\in\Z_q[X],
 \qquad \widetilde M_C\bmod p=M_C.
\]
\end{proposition}
\begin{proof}
The relative Frobenius has order $d$ on every finite unramified level. Its orbit sums and products commute with reduction, and its fixed ring is $\operatorname{GR}(p^s,r)$; hence the limit has fixed ring $\Z_q$. Since $\phi^r(\widetilde\zeta^i)=\widetilde\zeta^{qi}$, the coset root product is fixed by $\phi^r$ and has coefficients in $\Z_q$. Its reduction is exactly $M_C$.
\end{proof}

\begin{example}[Seventh roots of unity: the first nonzero lift of a zero trace]
\label{v132:ex-seven-period}
Let $\zeta$ be a primitive seventh Teichmüller root of unity in $\Z_{8}=W(\F_8)$. The nonzero binary cyclotomic cosets modulo seven are $C_1=\{1,2,4\}$ and $C_3=\{3,5,6\}$. Put
\[
 \eta=\zeta+\zeta^2+\zeta^4,\qquad
 \eta'=\zeta^3+\zeta^5+\zeta^6.
\]
Both are fixed by $\phi$ and hence belong to $\mathbb Z_2$. The sum of the seven roots is zero; in the nine terms of the product, three equal $1$ and the other six run through the nontrivial seventh roots. Hence
\begin{equation}\label{v137:eq-seven-periods}
 \eta+\eta'=-1,\qquad \eta\eta'=2,
 \qquad \eta^2+\eta+2=0.
\end{equation}
Choose $\zeta$ whose reduction is a root of $X^3+X+1$. Then $\overline\eta=0$ and $\overline{\eta'}=1$. Consequently $\eta'$ is a $2$-adic unit, and $\eta\eta'=2$ gives
\[
 \boxed{v_2(\eta)=1,\qquad \eta\equiv2\pmod4.}
\]
The vanishing trace in the finite field is not an exact zero in characteristic zero: it retains one digit of divisibility, visible modulo four.
\end{example}

\subsection{From exponent orbits to valuation profiles of canonical coefficients}
\label{v134:subsec-coset-profile}
Retain the notation $q=p^r$, $p\nmid N$, with a $q$-cyclotomic coset
$C\subseteq\Z/N\Z$ of size $d$. Let $\widetilde\zeta$ be a primitive
$N$th Teichm\"uller root of unity and put $A=\Z_q$. Write
\[
 \widetilde M_C(X)=\prod_{i\in C}(X-\widetilde\zeta^i)
 =\sum_{j=0}^d a_jX^{d-j},\qquad a_0=1,
\]
with $v_p(0)=+\infty$. Define the \textbf{valuation profile}
\begin{equation}\label{v136:eq-valuation-profile}
 \mathcal V_p(C)=\bigl(v_p(a_1),\ldots,v_p(a_d)\bigr).
\end{equation}
The coefficient support over the residue field detects only
whether $v_p(a_j)>0$; the valuation profile distinguishes positive finite valuation
from exact vanishing. We first separate coefficients determined by an exact root symmetry.

\begin{proposition}[Additive stabilizers and exact zeros in canonical lifts]
\label{v134:prop-stabilizer-zero}
Put
\[
 H_C=\{h\in\Z/N\Z:C+h=C\},\qquad\kappa_C=|H_C|.
\]
Then $\kappa_C\mid d$, and
\begin{equation}\label{v134:eq-composition-zero}
 \widetilde M_C(X)\in A[X^{\kappa_C}],\qquad
 \kappa_C\nmid j\ \Longrightarrow\ a_j=0.
\end{equation}
When $\kappa_C\mid j$, this symmetry alone asserts neither that $a_j$ is a unit
nor that it vanishes exactly.
\end{proposition}
\begin{proof}
The additive stabilizer $H_C$ is a subgroup of $\Z/N\Z$ acting freely on $C$ by
translation, so $\kappa_C\mid d$. For $h\in H_C$ let $u=\widetilde\zeta^h$.
Multiplication by $u$ permutes the roots, hence
$\widetilde M_C(uX)=u^d\widetilde M_C(X)$.
Comparison of coefficients of $X^{d-j}$ gives $(u^j-1)a_j=0$.
Since the root-containing ring $A[\widetilde\zeta]$ is an integral domain, $a_j\ne0$ requires $u^j=1$ for every $h\in H_C$;
this is equivalent to $\kappa_C\mid j$.
The conclusion is exact vanishing in characteristic zero, not merely vanishing modulo $p$.
\end{proof}

Proposition~\ref{v141:prop-affine-profile} also uses affine exponent symmetries to identify equal valuations at different phases and complementary coefficients, distinguishing symmetry of exact zeros from symmetry of positive valuations.

\begin{theorem}[Orbit-period formula from subset-sum counts]
\label{v134:thm-subset-periods}
For $0\le j\le d$ and $u\in\Z/N\Z$ define the integer
\[
 m_j(u)=\#\left\{J\subseteq C:|J|=j,\ \sum_{i\in J}i\equiv u\pmod N\right\}.
\]
Then $m_j(qu)=m_j(u)$. Let $\mathscr D$ denote all $q$-cyclotomic cosets in
$\Z/N\Z$, including $\{0\}$, and set
\[
 \eta_D=\sum_{u\in D}\widetilde\zeta^u\in A,\qquad
 m_j(D)=m_j(u)\quad(u\in D).
\]
Then
\begin{equation}\label{v134:eq-subset-periods}
 \boxed{a_j=(-1)^j\sum_{D\in\mathscr D}m_j(D)\eta_D.}
\end{equation}
In particular, $a_d=(-1)^d\widetilde\zeta^{\sum_{i\in C}i}$ is a unit, so
$v_p(a_d)=0$. Identity~\eqref{v134:eq-subset-periods} holds at every finite
Galois-ring precision $p^s$. Thus orbit counts, relative phases, and orbit-period
values modulo $p^s$ determine cancellations at that precision; the full profile requires compatible data at all precisions.
\end{theorem}
\begin{proof}
Expanding the elementary symmetric function of degree $j$ gives
$a_j=(-1)^j\sum_{|J|=j}\widetilde\zeta^{\sum_{i\in J}i}$.
Multiplication by $q$ permutes $C$ and maps subsets $J$ bijectively to subsets $qJ$.
Hence $m_j(qu)=m_j(u)$, and grouping exponents by $q$-orbits proves the formula.
The constant coefficient is a product of units.
\end{proof}

\begin{proposition}[First Witt-carry criterion for elementary symmetric sums]
\label{v134:prop-first-coefficient-carry}
Let $\zeta$ be the reduction of $\widetilde\zeta$. For
$J\subseteq C$, $|J|=j$, put $b_J=\zeta^{\sum_{i\in J}i}\in\F_{q^d}$.
If $a_j\equiv0\pmod p$, then $\sum_{|J|=j}b_J=0$, and
\[
 \mathcal B_{1,j}(C)=
 \overline{\frac1p\sum_{|J|=j}[b_J]}\in\F_q
\]
is well defined. Moreover,
\[
 \boxed{v_p(a_j)=1\ \Longleftrightarrow\ \mathcal B_{1,j}(C)\ne0,\qquad
 v_p(a_j)\ge2\ \Longleftrightarrow\ \mathcal B_{1,j}(C)=0.}
\]
Exact zeros are included in the second equivalence. The first carry is
computable using Proposition~\ref{v133:prop-zero-sum-obstruction}.
\end{proposition}
\begin{proof}
Multiplicativity of Teichm\"uller lifts gives
$a_j=(-1)^j\sum_{|J|=j}[b_J]$. The sum is fixed by relative Frobenius,
so it belongs to $\Z_q$. If $\overline a_j=0$, it lies in $p\Z_q$.
Dividing by $p$ and reducing gives $(-1)^j\mathcal B_{1,j}(C)$;
this scalar is nonzero in the residue field precisely when the carry is nonzero.
\end{proof}

The $\eta_D$ in~\eqref{v134:eq-subset-periods} are
\textbf{Frobenius-orbit periods}. If $N=\ell$ is prime, $\ell\ne p$,
and $C=aH$ with $H=\langle q\rangle\le(\Z/\ell\Z)^\times$,
then the nonzero orbit periods are Gauss periods in the usual sense.
Write $e=[(\Z/\ell\Z)^\times:H]$. Character orthogonality on the quotient gives
\begin{equation}\label{v134:eq-phase-gauss-fourier}
 \eta_{aH}=\frac1e\sum_{\chi\vert_H=1}\chi(a)^{-1}
 \sum_{x\in(\Z/\ell\Z)^\times}\chi(x)\widetilde\zeta^x,
\end{equation}
where the characters take values in a suitable algebraic closure.
This is an exact Fourier relation between phase $aH$ and Gauss sums. It must be
interpreted in a $p$-adic field containing the character values, and $1/e$ carries
negative valuation if $p\mid e$. In particular, the digit-sum Stickelberger formula
in Chapter 18 concerns additive characters of characteristic-$p$ finite fields
using $p$th roots of unity; it cannot be applied unchanged to the present Gauss
sums using $\ell$th roots of unity.

\begin{proposition}[Newton power sums and a valuation bound]
\label{v134:prop-newton-profile}
Put $P_k(C)=\sum_{i\in C}\widetilde\zeta^{ki}$ for $k\ge1$.
These are orbit-period sums for the exponent image $kC$, counted with multiplicity.
They satisfy
\begin{equation}\label{v134:eq-newton-profile}
 j a_j=-\sum_{k=1}^j a_{j-k}P_k(C)\qquad(1\le j\le d).
\end{equation}
For a partition $\lambda=(1^{m_1}2^{m_2}\cdots)\vdash j$ set
\[
 B_\lambda=\sum_{k\ge1}m_kv_p(P_k(C))
 -\sum_{k\ge1}\bigl(m_kv_p(k)+v_p(m_k!)\bigr).
\]
A summand with $P_k=0$ and $m_k>0$ is interpreted as zero with valuation $+\infty$.
Then
\begin{equation}\label{v134:eq-partition-bound}
 \boxed{v_p(a_j)\ge\min_{\lambda\vdash j}B_\lambda.}
\end{equation}
If the minimum is finite and attained by exactly one partition, equality holds.
For $p\nmid j$, equation~\eqref{v134:eq-newton-profile} permits division by $j$
in $A$; when $p\mid j$, that division loses $v_p(j)$ digits of precision.
\end{proposition}
\begin{proof}
Let $E(T)=\prod_{i\in C}(1-\widetilde\zeta^iT)=\sum_{j=0}^d a_jT^j$.
Taking a logarithmic derivative in $\mathbb Q_q[[T]]$ gives
$TE'(T)/E(T)=-\sum_{k\ge1}P_k(C)T^k$, proving Newton's recurrence.
Expanding
\[
 E(T)=\exp\!\left(-\sum_{k\ge1}\frac{P_k(C)}{k}T^k\right)
\]
gives the partition formula
\[
 a_j=\sum_{\lambda\vdash j}
 \prod_{k\ge1}\frac{(-P_k(C))^{m_k}}{m_k!\,k^{m_k}}.
\]
The ultrametric inequality yields the bound. A unique minimum-valuation term
cannot cancel with other terms, giving equality. The precision loss at $p\mid j$
follows immediately from $v_p(j)$.
\end{proof}

Actual period values, or sufficiently many $p$-adic digits, are needed to decide
sharpness: a list of the valuations $v_p(P_k)$ alone cannot rule out cancellation
among terms of equal minimum valuation. Put $f=\operatorname{ord}_N(q)$, with $f=1$ when $N=1$. An effective calculation first lifts a primitive $N$th root to $\operatorname{GR}(p^s,rf)$ and then evaluates
\eqref{v134:eq-subset-periods}. This uses only ring addition and multiplication,
with no division by a nonunit. The result is the true reduction $a_j\bmod p^s$.
If the coefficient first becomes nonzero modulo $p^{s+1}$, then $v_p(a_j)=s$;
continued vanishing at finite precision implies only a lower bound, never exact zero.

\begin{example}[The two phases modulo seven and their full profiles]
\label{v134:ex-seven-full-profile}
With $\eta,\eta'$ as in Example~\ref{v132:ex-seven-period},
\[
 \widetilde M_{C_1}(X)=X^3-\eta X^2+\eta'X-1,\qquad
 \widetilde M_{C_3}(X)=X^3-\eta'X^2+\eta X-1.
\]
Using $\eta\eta'=2$, $\eta+\eta'=-1$, and $\overline\eta=0$, we obtain
\[
 \boxed{\mathcal V_2(C_1)=(1,0,0),\qquad
 \mathcal V_2(C_3)=(0,1,0).}
\]
The two orbits have the same length, root orders, and additive stabilizer,
yet different phases place the first positive valuation in different positions.
This follows from the algebraic period relations and the chosen residue root,
not from absolute values of the periods alone.
\end{example}

\begin{proposition}[Compatibility of profiles under base-field extension]
\label{v134:prop-base-change-profile}
Extend the base field from $\F_q$ to $\F_{q^s}$ and put
$g=\gcd(s,d)$. The coset $C$ splits into $g$ orbits
$C_0,\ldots,C_{g-1}$ for $q^s$, each of size $d/g$, with
\[
 \widetilde M_C=\prod_{\nu=0}^{g-1}\widetilde M_{C_\nu}
 \quad\text{in }\Z_{q^s}[X].
\]
Corresponding coefficients of the $g$ factors are unramified Frobenius
conjugates and therefore have \emph{identical valuation profiles}. Their
product may nevertheless exhibit cancellations, so the profile of the product
is not determined solely by the numerical profiles of its factors.
\end{proposition}
\begin{proof}
The $q^s$ action on a cyclic $q$-orbit of length $d$ has $\gcd(s,d)$
cycles of equal length. Grouping the roots gives the factorization.
The $q$-power Frobenius moves transitively among these cycles and preserves
$p$-adic valuations. The ultrametric inequality need not be an equality when
the factor polynomials are multiplied.
\end{proof}

\subsection*{Predicting valuation profiles from orbit structure}
\phantomsection\label{v132:problem-coset-valuations}
Theorems~\ref{revadic:thm-integral-resolvent} and~\ref{revadic:thm-local-min} in Chapter~\ref{revadic:chapter} give the integral-remainder zero criterion and the minimum valuation across all local components.
For each $q$-cyclotomic coset $C$, first remove the exact zeros forced by
$\kappa_C\nmid j$. Can the remaining $v_p(a_j)$ be predicted precisely
from subset-sum orbit counts $m_j(D)$, synchronized local phases,
and known Gauss-period congruences or valuations?
Equation~\eqref{v134:eq-subset-periods} gives an exact expression and
\eqref{v134:eq-partition-bound} a verifiable lower bound. The essential
unresolved issue is cancellation among terms of the same least valuation:
when can orbit structure replace successive Hensel computations?
Proposition~\ref{v137:prop-norm-depth} gives a norm constraint on local depths; the quadratic periods in Proposition~\ref{v137:prop-quadratic-period-depth} yield exact positive valuations without successive lifting. Theorem~\ref{v141:thm-index-budget} relates excess phase-difference depths to a generating-order index; Theorem~\ref{v141:thm-minimal-excess} determines the last congruence layer under minimal positive excess. What further compatibility relations govern larger excess? One may also seek profile patterns uniform along towers of moduli and study
how equal-profile factors recombine after extension of the base field.
Theorem~\ref{v138:thm-periodic-cosets} converts the same exponent orbit into a periodic digit block. Digit sums control Gauss-sum valuations, while coefficient cancellation also depends on joint block distributions.

\section*{Exercises}
\phantomsection
\addcontentsline{toc}{section}{Exercises}
\begin{exercise}\label{ex:phi9}
List the binary cyclotomic cosets modulo nine, and determine the factor degrees of \(\Phi_9\) over \(\F_2\).
Then determine its factor degrees over \(\F_4\) and \(\F_8\).
\end{exercise}
\begin{exercise}\label{ex:cyclop}
Find the reductions of \(\Phi_{12}\) modulo \(2\) and \(3\), and explain the root multiplicities.
\end{exercise}
\begin{exercise}\label{ex:coset-negative}
Show that taking reciprocal polynomials corresponds to the operation \(C\mapsto-C\) on exponent orbits.
Use this to reinterpret the self-reciprocal factor in the fifteenth-root example.
\end{exercise}
\begin{exercise}\label{ex:degree5extension}
Does an irreducible quintic over \(\F_q\) remain irreducible over \(\F_{q^6}\)? How does it factor over \(\F_{q^{10}}\)?
\end{exercise}

\begin{exercise}[The degree of the splitting field of an arbitrary polynomial]\label{ex:wu-splitting-degree}
Let \(F=\prod_i f_i^{a_i}\in\F_q[X]\) be nonconstant, with \(a_i\ge1\) and distinct monic irreducible polynomials \(f_i\),
where \(\deg f_i=d_i\). Prove that its splitting field has degree \(\lcm_i d_i\).
Do the multiplicities \(a_i\) affect this degree?
\end{exercise}

\begin{exercise}[Coefficient fields and further base extension]\label{ex:wu2-factor-coefficients}
Let \(f\in\F_q[X]\) be monic and irreducible of degree six.
Determine the number and degrees of its factors over \(\F_{q^4}\) and \(\F_{q^9}\), and the exact coefficient field of each factor.
Then determine its factorization over \(\F_{q^{12}}\), and explain whether the fact that the coefficients collectively generate a field implies that each nonzero coefficient generates that field.
\end{exercise}

\begin{exercise}[Multiple equal-difference decompositions and the least extension for binomial factorization]\label{ex:med-seventytwo}
Take \(q=5,N=72,i=1\). Compute \(\tau,\rho,\omega_C,\kappa_C\), and list the coarsest multiple equal-difference decomposition and all admissible block sizes.
If \(\alpha\) is a primitive seventy-second root of unity, determine the maximum monomial composition exponent of its minimal polynomial
and the least extension degree over which all its irreducible factors become binomials.
\end{exercise}

\begin{exercise}[Trivial block decompositions and changes of Frobenius step]\label{ex:med-boundaries}
Show that if \(d=N/\gcd(N,i)\) is square-free, then \(C_i\) admits only the singleton multiple equal-difference decomposition.
For the degree-eight irreducible polynomial in Example~\ref{exm:med-thirtytwo}, determine the number and degrees of its factors over \(\F_{3^2}\), \(\F_{3^3}\), and \(\F_{3^4}\),
and decide whether all these irreducible factors are binomials.
\end{exercise}

\begin{exercise}\label{ex:supp-binomial-factors}
Determine whether all irreducible factors of \(X^{16}-1\) over \(\F_5\) and over \(\F_3\) are binomials.
\end{exercise}

\Needspace{5\baselineskip}
\begin{exercise}[Recognizing equal-difference blocks from traversal order]\label{ex:med-eighteen-five}
Take \(q=5,N=18\). List all \(q\)-cyclotomic cosets modulo eighteen and compute \(\tau,\omega_C,\kappa_C\) for each.
List the elements of each nonsingleton maximal equal-difference block in increasing order. Explain why the complete equal-difference structures of \(C_1\) and \(C_3\) differ.
\end{exercise}

\begin{exercise}[Maximal blocks and two levels of binomial factorization]\label{ex:med-twentyfive}
Take \(q=4,N=25\). Compute \(C_1\), write its coarsest multiple equal-difference decomposition, and find the least extension degree over which every irreducible factor of \(M_{C_1}\) is a binomial.
For a primitive twenty-fifth root of unity \(\zeta\), write the corresponding binomial product. Explain how the maximal block size changes when the index is replaced by five.
\end{exercise}

\begin{exercise}[Recovering factorization types from fixed roots of unity]\label{v40:ex-unit-fixed-counts}
Let \(N\ge1\) and \(\gcd(N,q)=1\). Let \(\mu_N\) be the set of all \(N\)-th roots of unity in a common extension, and let
\(\sigma(x)=x^q\). Denote by \(a_r\) the number of monic irreducible factors of degree \(r\) of \(X^N-1\) in \(\F_q[X]\).
\begin{enumerate}
\item Prove that, for every \(t\ge1\),
\[
 |\mu_N\cap\F_{q^t}|=\gcd(N,q^t-1).
\]
Deduce that
\[
 a_r=\frac1r\sum_{t\mid r}\mu(r/t)\gcd(N,q^t-1)
 =\frac1r\sum_{\substack{d\mid N\\\operatorname{ord}_d(q)=r}}\varphi(d),
\]
where \(\operatorname{ord}_1(q)=1\).
\item Without first listing all cyclotomic cosets, compute the number of factors of each degree for \(q=2,N=15\).
\item Put \(k=\F_q\), \(A=k[X]/(X^N-1)\), and \(\sigma_A(a)=a^q\). Prove that
\[
 \dim_k\ker(\sigma_A^t-\id_A)=\sum_{r\ge1}a_r\gcd(r,t).
\]
For \(q=2,N=15,t=1\), compute the number of fixed points on the root set and the dimension of the fixed space in the algebra,
and explain why they are different kinds of data.
\end{enumerate}
\end{exercise}

\begin{exercise}[A common Frobenius period for nontrivial roots of unity]\label{v40:ex-uniform-root-degrees}
Let \(N>1\), \(\gcd(N,q)=1\), and put
\[
S_N(X)=\frac{X^N-1}{X-1}.
\]
\begin{enumerate}
\item If every monic irreducible factor of \(S_N\) has degree \(m\), determine \(\mu_N\cap\F_{q^t}\) for every \(t\ge1\).
\item Factor \(S_N\) for \(q=2,N=9\). Show that agreement of the multiplicative orders at all prime divisors does not allow the condition \(N\mid q^m-1\) in Proposition~\ref{v124:prop-uniform-cyclotomic} to be omitted.
\item Show that for \(q=9,N=91=7\cdot13\), the polynomial \(S_N\) has exactly thirty irreducible cubic factors. Count the factors with root orders seven, thirteen, and ninety-one separately.
\item Take \(N=q^s-1\), with \(s\ge2\). Determine all pairs \((q,s)\) for which the irreducible factors of \(S_N\) have a common degree.
\end{enumerate}
\end{exercise}

\begin{exercise}[Synchronized phases of three local orbits]\label{v44:ex-local-phases-315}
Take \(q=2\) and \(N=315=9\cdot5\cdot7\).
\begin{enumerate}
\item Find representatives and lengths of the unit orbits modulo nine, five, and seven.
\item For the orbits represented by one modulo nine and modulo five, and the orbit represented by \(b\in\{1,3\}\) modulo seven,
use Proposition~\ref{v44:prop-recursive-phase-representatives} to list the recursive phase ranges.
\item Show that \(\Phi_{315}\) has exactly twelve irreducible factors over \(\F_2\), each of degree twelve,
and verify that orbit representatives may be chosen as
\[
 280+126\cdot2^{j_2}+225\cdot b\,2^{j_3}\pmod{315},
 \quad b\in\{1,3\},\quad 0\le j_2<2,\quad 0\le j_3<3.
\]
\end{enumerate}
\end{exercise}

\begin{exercise}[Local and global coefficient fields]\label{v124:ex-local-coefficient-fields}
Choose a primitive twenty-fourth root of unity \(\zeta\) in \(\overline{\F}_{13}\), and put \(H=\{1,5,7,11\}\subseteq(\Z/24\Z)^\times\). Determine the exact coefficient fields of
\[
P_8(X)=\prod_{a\in(\Z/8\Z)^\times}(X-\zeta^{3a}),\qquad
P_3(X)=\prod_{a\in(\Z/3\Z)^\times}(X-\zeta^{8a}),
\]
\[
B(X)=\prod_{a\in H}(X-\zeta^a).
\]
Then find the coefficientwise thirteenth-power conjugate of \(B\) and its product with \(B\).
\end{exercise}

\begin{exercise}[Valuation layers and eighty-first roots of unity]\label{v44:ex-local-orbits-81}
Take \(q=4\) and \(N=81\), and let \(\omega\in\F_4\) satisfy \(\omega^2+\omega+1=0\).
\begin{enumerate}
\item Use Proposition~\ref{v44:prop-prime-power-local-orbits} to calculate the length distribution of all cyclotomic cosets modulo eighty-one,
and hence the number and degree distribution of the irreducible factors of \(X^{81}-1\).
\item Prove that
\[
 X^{81}-1=(X-1)\prod_{a=0}^{3}
 (X^{3^a}-\omega)(X^{3^a}-\omega^2),
\]
and that every factor on the right is irreducible in \(\F_4[X]\).
\item Find the irreducible factorization of \(\Phi_{81}\), and explain which valuation layer contains the roots of each degree-twenty-seven factor.
\end{enumerate}
\end{exercise}

\begin{exercise}[Return permutations and factor lifting]\label{v44:ex-cyclotomic-return-factors}
Over \(\F_2\), take the parent orbit \(C_3(1)=\{1,2\}\) modulo three.
Lift it to moduli nine and fifteen, calculating the return permutations, all lifted orbits, and their lengths.
Deduce that \(X^6+X^3+1\) is irreducible, and give the complete irreducible factorization of \(X^{10}+X^5+1\).
\end{exercise}

\begin{exercise}[Admissible digits and least coset representatives]\label{v44:ex-cyclotomic-compatible-leaders}
Take \(q=2\), and lift \(C_9(1)\) from modulus nine to modulus forty-five.
Using the least nonnegative representative of each return orbit as an admissible digit, find compatible representatives and lengths of all child orbits.
Then find the least integer in each child orbit, test whether it reduces to the parent representative one, and explain whether this contradicts orbit compatibility.
\end{exercise}

\begin{exercise}[Exponent orbits and multiplicities for a binomial]
\label{v44:ex-binomial-six}
Over \(\F_5\), let \(\lambda=2\).
\begin{enumerate}
\item Find \(g=\ord(\lambda)\), list all \(5\)-cyclotomic orbits in \(S_{6,g}=1+g\Z/6g\Z\), and determine the number and degrees of the irreducible factors of \(X^6-2\).
\item Write the monic irreducible factorization of \(X^6-2\), and prove that every listed factor is irreducible.
\item Find the monic irreducible factorization of \(X^{30}-2\), explaining what determines the factor degrees and multiplicities.
\end{enumerate}
\end{exercise}

\begin{exercise}[Dyadic stabilization indices and irreducible substitutions]
\label{v46:ex-cyclotomic-stable-binary}
Take \(q=7\), \(m=1\), and \(\ell=2\), and consider the system \((C_{2^i}(1))_{i\ge0}\).
\begin{enumerate}
\item List its orbits and lengths modulo two, four, eight, sixteen, and thirty-two, and find its deviation depth and stabilization index.
\item Prove that \(f(X)=X^2+X+6\in\F_7[X]\) is irreducible and that its roots have multiplicative order sixteen.
\item Prove that, for every \(t\ge0\), the polynomial
\[
 X^{2^{t+1}}+X^{2^t}+6
\]
is irreducible in \(\F_7[X]\), and find its degree.
\end{enumerate}
\end{exercise}

\begin{exercise}[The obstruction to ring sections in mixed characteristic]\label{v133:ex-no-section}
Let $R$ be a nonzero commutative unital ring without $p$-torsion.
Show that every additive group homomorphism $\F_q\to(R,+)$ is zero,
so that no unital ring homomorphism from $\F_q$ to $R$ exists.
Replacing the target by $\F_q[[t]]$, show that the constant embedding
is the unique unital ring section of reduction modulo $t$.
\end{exercise}

\begin{exercise}[The first Witt carry and its compatibilities]\label{v133:ex-cocycle}
Let $\kappa$ be as in Proposition~\ref{v133:prop-first-carry}.
\begin{enumerate}
\item Use associativity in $W_2(K)$ to show
$\kappa(a,b)=\kappa(b,a)$ and
$\kappa(a,b)+\kappa(a+b,c)=\kappa(b,c)+\kappa(a,b+c)$.
\item Show that $\kappa(ua,ub)=u\kappa(a,b)$ for all $u\in K$.
\item Give a simplified expression for $\kappa(a,b)$ over $K=\F_4$;
explain why the first carry of $[a]+[b]$ is not simply $ab$.
\end{enumerate}
\end{exercise}

\begin{exercise}[Different lifting depths of zero-sum identities]\label{v133:ex-zero-sum-depth}
Over $K=\F_2$, compare the valuations of the Teichmüller lifts
of $1+1=0$ and $0+0=0$.
Now let $\omega\in\F_4$ be a primitive cube root of unity.
Show that their Teichmüller lifts in $W(\F_4)$ satisfy $[1]+[\omega]+[\omega^2]=0$ exactly.
Explain why a zero-sum identity in a finite field does not by itself
determine the valuation of its lift.
\end{exercise}

\begin{exercise}[First trace obstruction and norms]\label{v133:ex-trace-norm}
Let $E=\F_4$, $K=\F_2$. For $a=0,1$, compute the finite-field trace,
the relative trace of the Teichm\"uller representative $[a]$, and the first obstruction (when the finite-field trace is zero).
For general $E=\F_{q^d}$, prove~\eqref{v133:eq-norm-exact}
and explain why the corresponding assertion for traces does not follow.
\end{exercise}

\begin{exercise}[Formal roots and different deformations of trace and norm]\label{v133:ex-formal-deform}
Let $q=p^r$ and choose $c\in\F_q$ with $\Tr_{\F_q/\F_p}(c)\ne0$.
\begin{enumerate}
\item Find the unique $B\in\F_q[[t]]$ with $B(0)=0$ and $B^p-B=t$.
\item Prove that the norm in Example~\ref{v133:ex-artin-schreier-deformation}
is $c+t$ and compute its trace separately for $p=2$ and $p>2$.
\item Prove that $X^p-X-(c+t)$ remains irreducible over $\F_q((t))$.
Hint: if it has a root $u\in\F_q((t))$, use the $t$-adic valuation
to show that $u\in\F_q[[t]]$, then reduce modulo $t$.
\end{enumerate}
\end{exercise}

\begin{exercise}[Stabilizers, exact zeros, and additional cancellation]\label{v134:ex-stabilizer}
Suppose $p\nmid N$, and let $C$ be a $q$-cyclotomic coset of size $d$.
\begin{enumerate}
\item Show that $|H_C|=\kappa_C$ implies $\widetilde M_C\in\Z_q[X^{\kappa_C}]$.
\item Explain why a further coefficient of $M_C$ vanishing modulo $p$ does not imply
that its canonical lift vanishes exactly.
\item In $\widetilde F_{5,2}$ from Example~\ref{v125:ex-160-hensel}, distinguish
coefficients forced to vanish by the composition exponent from those vanishing
modulo three but not modulo nine.
\end{enumerate}
\end{exercise}
\begin{exercise}[Subset-sum counts and orbit periods]\label{v134:ex-subset-periods}
Take $q=2$, $N=7$, and $C=\{1,2,4\}$.
\begin{enumerate}
\item Find all nonzero $m_j(u)$ for $j=1,2,3$ and derive $a_1,a_2,a_3$ from
 the orbit-period formula.
\item Use Example~\ref{v132:ex-seven-period} to determine $\mathcal V_2(C)$.
\item Repeat the calculation for $C'=\{3,5,6\}$ and compare the profiles.
\end{enumerate}
\end{exercise}
\begin{exercise}[Equality and cancellation in the Newton partition bound]\label{v134:ex-newton-bound}
\begin{enumerate}
\item Derive $a_1=-P_1$ and $2a_2=P_1^2-P_2$ from
$E(T)=\prod_{i\in C}(1-\widetilde\zeta^iT)$.
\item If $u=v_p(P_1)$ and $w=v_p(P_2)$ are finite, prove
$v_p(a_2)\ge\min(2u,w)-v_p(2)$, with equality when $2u\ne w$.
\item Explain why $u,w$ alone do not decide equality if $2u=w$;
when $p=2$, explain the precision issue in dividing by two to compute $a_2\bmod 2^s$.
\end{enumerate}
\end{exercise}
\begin{exercise}[Equal-profile factors after extending the base field]\label{v134:ex-base-change}
Let $q=2$, $N=7$, and $C=\{1,2,4\}$, and extend the base field to $\F_8$.
\begin{enumerate}
\item List the three $8$-cyclotomic subcosets and factor the canonical polynomial in $\Z_8[X]$.
\item Determine the profiles of the three linear factors and explain why these
profiles alone do not determine the parent cubic profile $(1,0,0)$.
\item Generalize to a $q$-coset of length $d$ and base extension degree $s$.
\end{enumerate}
\end{exercise}

\begin{exercise}[The first nonzero carry of a zero coefficient]\label{v134:ex-first-coefficient-carry}
Show that $\mathcal B_{1,j}(C)$ in
Proposition~\ref{v134:prop-first-coefficient-carry} belongs to $\F_q$,
not merely to $\F_{q^d}$. For the two binary cosets $C_1,C_3$ modulo seven,
compute $\mathcal B_{1,1}(C_1)$ and $\mathcal B_{1,2}(C_3)$.
\end{exercise}

\begin{exercise}[Teichm\"uller representatives and Frobenius]\label{v47:ex-galois-extensions}
Let \(A=(\Z/4\Z)[\omega]\), where \(\omega^2+\omega+1=0\).
\begin{enumerate}
\item Find the cardinality, number of units, and Teichm\"uller set of \(A\).
\item For \(a,b\in\Z/4\Z\), write the absolute Frobenius \(\phi(a+b\omega)\), and determine its fixed subring.
\item Compare \(\phi(1+\omega)\) with \((1+\omega)^2\), and find the multiplicative order of \(1+\omega\).
\end{enumerate}
\end{exercise}

\begin{exercise}[Lifting roots-of-unity factors by local corrections]\label{v47:ex-hensel-seven}
Over \(\F_2\), we have \(X^7-1=(X-1)(X^3+X+1)(X^3+X^2+1)\).
Starting with \(G_1=X^3+X+1\), use Proposition~\ref{v47:prop-galois-local-lift} to find successively
the unique monic factor lifts modulo \(4\), \(8\), and \(16\).
Then use monic reciprocal polynomials to write the complete basic irreducible factorization of \(X^7-1\) over \(\Z/16\Z\),
and verify the product directly.
\end{exercise}

\begin{exercise}[Coset orbits and coefficient fields]\label{v85:ex-coset-orbits}
Take \(K=\F_{16}\), \(k=\F_2\), a generator \(\gamma\) of \(K^\times\), and \(H=\langle\gamma^3\rangle\).
Determine all orbits of the squaring Frobenius on \(K^\times/H\), and the smallest coefficient field of each of the three coset root polynomials.
Prove also that
\[
 p_{\gamma H}(T)p_{\gamma^2H}(T)=T^{10}+T^5+1.
\]
Explain the distinction between coset orbit length, coset size, and the Frobenius period of an individual element of a coset.
\end{exercise}

\chapter{Trace, Norm, and Conjugate Coordinates}
\label{chap:trace}
Sums and products along conjugate orbits give trace and norm. They are also the trace and determinant of multiplication operators and govern Frobenius difference equations. Equivariant splitting of the trace further reveals the interaction between the characteristic and the extension degree.

\section{Orbit Sums and Orbit Products}
Fix \(k=\F_q\subseteq E=\F_{q^n}\) and \(\sigma(a)=a^q\).
\begin{definition}\label{def:trace-norm}
The \textbf{trace} and \textbf{norm} of \(a\in E\) over \(k\) are defined by
\[
\Tr_{E/k}(a)=\sum_{i=0}^{n-1}\sigma^i(a),\qquad
\Nm_{E/k}(a)=\prod_{i=0}^{n-1}\sigma^i(a).
\]
Thus
\[
\Tr_{E/k}(a)=a+a^q+\cdots+a^{q^{n-1}},\qquad
\Nm_{E/k}(a)=a^{(q^n-1)/(q-1)}.
\]
They are also called the \textbf{relative trace} and \textbf{relative norm} over \(k\). When \(k=\F_p\) is the prime field, \(\Tr_{E/\F_p}\) is called the \textbf{absolute trace}.
\end{definition}
Frobenius cyclically permutes the terms in the sum and product, so both values are fixed by \(\sigma\) and hence belong to \(k\).
The trace is \(k\)-linear, and the norm on \(E^\times\) is a homomorphism of multiplicative groups.
For \(c\in k\), we have \(\Tr_{E/k}(c)=nc\) and \(\Nm_{E/k}(c)=c^n\).
Here \(n\) denotes its image in the field, which may be zero.

If \(d=[k(a):k]\), the conjugate orbit has only \(d\) distinct elements, but the definitions above still have \(n\) terms.
Each conjugate occurs \(n/d\) times, so
\begin{equation}\label{eq:trace-subdegree}
\Tr_{E/k}(a)=\frac nd\Tr_{k(a)/k}(a),\qquad
\Nm_{E/k}(a)=\Nm_{k(a)/k}(a)^{n/d}.
\end{equation}
In the first identity, \(n/d\) counts repeated additions; it does not denote division in the field.
Summing or multiplying only the distinct conjugate roots gives the trace or norm for the smaller extension \(k(a)/k\).

Descent of multivariate linear factors is explained by the same conjugate product. For the
normalized linear factor \(L\) in Theorem~\ref{thm:multivariate-linear-descent}, suppose its coefficients generate \(E_0=\F_{q^s}\). For every \(x\in k^m\),
\begin{equation}\label{eq:linear-form-norm-value}
F_L(x)=\prod_{i=0}^{s-1}L(x)^{q^i}=\Nm_{E_0/k}(L(x)).
\end{equation}
Thus the orbit product of linear factors, evaluated at a base-field point, is the norm of the corresponding linear expression.
Even if \(L(x)\) lies in a proper subfield of \(E_0\), the product must retain all \(s\) terms, in agreement with the repetitions in \eqref{eq:trace-subdegree}.

\begin{proposition}[Transitivity]\label{prop:trace-tower}
If \(k\subseteq L\subseteq E\), then
\[
\Tr_{E/k}=\Tr_{L/k}\circ\Tr_{E/L},\qquad
\Nm_{E/k}=\Nm_{L/k}\circ\Nm_{E/L}.
\]
\end{proposition}
\begin{proof}
Put \([L:k]=d\) and \([E:L]=s\), so \(n=ds\). Frobenius relative to \(L\) is \(\sigma^d\). Hence
\[
\Tr_{L/k}\bigl(\Tr_{E/L}(a)\bigr)
=\sum_{i=0}^{d-1}\sum_{j=0}^{s-1}\sigma^{i+dj}(a)
=\sum_{t=0}^{n-1}\sigma^t(a).
\]
The same grouping of indices in the product gives the norm formula.
\end{proof}

\section{Surjectivity and Fiber Sizes}
\label{sec:fibers}

For a map \(h:S\to T\) and a fixed element \(t\) of the target, the set of all elements mapped to \(t\),
\[
h^{-1}(t)=\{x\in S:h(x)=t\}
\]
is called the \textbf{fiber} of \(h\) over \(t\). Here \(h^{-1}(t)\) is the preimage of the singleton \(\{t\}\);
it does not require \(h\) to have an inverse. Fibers may be empty; the nonempty fibers are disjoint and their union is the whole domain \(S\).
For example, the trace fiber over \(c\in k\) is \(\{x\in E:\Tr_{E/k}(x)=c\}\),
and the fiber of the norm on the multiplicative group over \(c\in k^\times\) is \(\{x\in E^\times:\Nm_{E/k}(x)=c\}\).
A fiber is always understood relative to a specified map and domain.

\begin{theorem}\label{thm:trace-norm-surjective}
The trace \(\Tr_{E/k}:E\to k\) is surjective, with \(q^{n-1}\) elements in each fiber.
The norm \(\Nm_{E/k}:E^\times\to k^\times\) is also surjective, with \((q^n-1)/(q-1)\) preimages of each nonzero value.
\end{theorem}
\begin{proof}
The trace is given by the nonzero polynomial \(X+X^q+\cdots+X^{q^{n-1}}\), whose degree \(q^{n-1}<q^n\) prevents it from vanishing at every element of \(E\).
Thus it is a nonzero linear map to the one-dimensional \(k\)-space \(k\), and is surjective.
Its kernel has dimension \(n-1\). If \(x_0\) has trace \(c\), then \(\Tr(x)=c\) is equivalent to \(x-x_0\in\ker\Tr\), so the fiber is \(x_0+\ker\Tr\), of size \(q^{n-1}\).

Choose a generator \(\zeta\) of \(E^\times\). The element \(\Nm_{E/k}(\zeta)=\zeta^{(q^n-1)/(q-1)}\) has order \(q-1\), and hence generates \(k^\times\).
If \(\Nm(x_0)=c\ne0\), then \(\Nm(x)=c\) is equivalent to \(x_0^{-1}x\in\ker\Nm\).
Thus the fiber is the multiplicative coset \(x_0\ker\Nm\), whose size follows by counting kernel and image of the group homomorphism.
\end{proof}

\begin{remark}
When \(p\mid n\), we have \(\Tr_{E/k}(1)=0\), but the trace map is still nonzero.
A linear operator \(P\) with \(P^2=P\) is called a projection; it maps the whole space onto \(\im P\) and acts identically on that image.
The trace operator itself need not be a projection: viewing \(k\) as a subspace of \(E\), it satisfies \(\Tr^2=n\Tr\). Only when \(p\nmid n\) does \(n^{-1}\Tr\) define a projection onto \(k\).
\end{remark}

\begin{example}
In \(E=\F_3(\alpha)\), with \(\alpha^2=-1\), we have \(\alpha^3=-\alpha\). For \(x=a+b\alpha\),
\[
\Tr_{E/\F_3}(x)=2a,\qquad
\Nm_{E/\F_3}(x)=a^2+b^2.
\]
In \(\F_4/\F_2\), however, \(\Tr(1)=0\); if \(\beta^2+\beta+1=0\), then \(\Tr(\beta)=1\).
\end{example}

\par\noindent\textbf{Norm preservation and linear structure.}\phantomsection\label{v74:note-norm-preservers}
The nonzero norm fibers are multiplicative cosets of \(\ker\Nm_{E/k}\).
If a \(k\)-linear map \(T:E\to E\) preserves the norm of every element, then \(T(x)-T(y)=T(x-y)\), so it also preserves the norm of every nonzero difference.
A condition on norm fibers thus becomes a compatibility condition between addition and multiplication.
When \(|k|>2\), such a map must be multiplication by a norm-one element followed by a \(k\)-automorphism; see Exercise~\ref{v74:ex-norm-preservers}.
This is the finite-field case of Lenstra's norm-form theorem \cite[Theorem~4]{Lenstra1990Automorphisms}; the solution explains the mechanism using the difference-quotient rigidity established later.

\section{Compatibility of Multiple Norms and Joint Fibers}
\label{v58:sec-compatible-norms}

When norms to two subfields are prescribed separately, surjectivity of each norm does not ensure that one element realizes both values.
Their norms to the intersection must agree; multiplicative cyclicity turns this into compatibility of integer congruences.
The following result also appears in \cite[Lemma~3.1]{FPRPrescribedNorms}.

Let \(d_1,\ldots,d_s\mid n\), and put \(L_i=\F_{q^{d_i}}\subseteq E\).
Given nonzero values \(a_i\in L_i^\times\), suppose that, for every pair \(i,j\), putting \(e=\gcd(d_i,d_j)\) gives
\begin{equation}\label{v58:eq-norm-compatibility}
 \Nm_{L_i/\F_{q^e}}(a_i)=\Nm_{L_j/\F_{q^e}}(a_j).
\end{equation}
\leavevmode\phantomsection\label{v58:def-norm-compatibility}\ignorespaces
We then say that the values satisfy the \textbf{norm compatibility condition}.
Here \(\F_{q^e}=L_i\cap L_j\).
\leavevmode\phantomsection\label{v58:def-joint-norm-fiber}\ignorespaces
The elements \(x\in E^\times\) simultaneously satisfying \(\Nm_{E/L_i}(x)=a_i\) form the
\textbf{joint fiber} of the norm map \(x\mapsto(\Nm_{E/L_i}(x))_i\).

\begin{proposition}[Solvability and counting for multiple norms]\label{v58:prop-compatible-norms}
Let \(E=\F_{q^n}\), let \(s\ge1\), and let \(d_i\mid n\), \(L_i=\F_{q^{d_i}}\), and \(a_i\in L_i^\times\) for \(1\le i\le s\). The fiber \(\{x\in E^\times:\Nm_{E/L_i}(x)=a_i\text{ for all }i\}\) is nonempty if and only if the norm compatibility condition~\eqref{v58:eq-norm-compatibility} holds. When the prescribed values are compatible, the fiber is a multiplicative coset of the kernel of the joint norm map, of size
\begin{equation}\label{v58:eq-joint-norm-count}
 \frac{q^n-1}{\lcm(q^{d_1}-1,\ldots,q^{d_s}-1)}
 =\gcd\left(\frac{q^n-1}{q^{d_1}-1},\ldots,
             \frac{q^n-1}{q^{d_s}-1}\right).
\end{equation}
\end{proposition}
\begin{proof}
Necessity follows immediately from transitivity of the norm. For sufficiency and the count, choose a generator \(\zeta\) of \(E^\times\),
and put \(N=q^n-1\), \(m_i=q^{d_i}-1\). The element \(\zeta^{N/m_i}\) generates \(L_i^\times\), so there is a unique residue class \(r_i\) modulo \(m_i\) with
\[
 a_i=\zeta^{(N/m_i)r_i}.
\]
Writing \(x=\zeta^r\), the equation \(\Nm_{E/L_i}(x)=a_i\) is equivalent to \(r\equiv r_i\pmod{m_i}\).
The Euclidean algorithm gives
\(\gcd(q^{d_i}-1,q^{d_j}-1)=q^{\gcd(d_i,d_j)}-1\).
Thus \eqref{v58:eq-norm-compatibility} states precisely that
\[
 r_i\equiv r_j\pmod{\gcd(m_i,m_j)}.
\]

These pairwise compatibility conditions make the whole congruence system solvable. Indeed, for every prime \(\ell\) dividing some \(m_i\),
choose an index \(i\) maximizing \(v_\ell(m_i)\), and prescribe \(r\equiv r_i\pmod{\ell^{v_\ell(m_i)}}\).
Pairwise compatibility ensures that this prescription satisfies the conditions in the \(\ell\)-direction for all other indices.
The powers of distinct primes are coprime, so the Chinese remainder theorem gives a unique solution modulo \(M=\lcm(m_1,\ldots,m_s)\).
Since \(M\mid N\), there are \(N/M\) solutions modulo \(N\), and hence that many field elements.
Comparison of prime-adic valuations also gives the second expression in \eqref{v58:eq-joint-norm-count}.
Finally, the joint norm map is a group homomorphism, so each nonempty fiber is a coset of its kernel.
\end{proof}
The nonzero-value hypothesis distinguishes the group fiber from its zero boundary. If zero values are allowed, all prescribed values zero give the singleton fiber \(\{0\}\); a mixture of zero and nonzero prescribed values gives an empty fiber, since \(\Nm_{E/L_i}(x)=0\) if and only if \(x=0\).

\begin{example}[Two norms in a degree-six extension]\label{v58:ex-two-norms-degree-six}
Take \(E=\F_{q^6}\), \(L_1=\F_{q^2}\), and \(L_2=\F_{q^3}\).
Given \(a\in L_1^\times\) and \(b\in L_2^\times\), the two norms can be prescribed simultaneously if and only if
\[
 a^{q+1}=b^{q^2+q+1}.
\]
When compatible, the joint fiber has size
\[
 \frac{(q^6-1)(q-1)}{(q^2-1)(q^3-1)}=q^2-q+1.
\]
For example, when \(q=5\), each compatible pair of nonzero norm values has exactly \(21\) common preimages.
\end{example}

One norm value determines a multiplicative coset; the fiber determined by several norm values depends on compatibility over subfield intersections and the least common multiple of integer moduli.
Section~\ref{v58:sec-prescribed-normal-traces} gives the corresponding result for relative traces using normal coordinates.

\section{Two Hilbert 90 Identities}
\begin{theorem}[Additive and multiplicative Hilbert 90]\label{thm:hilbert90}
We have
\[
\ker\Tr_{E/k}=\{b^q-b:b\in E\},
\]
\[
\ker\bigl(\Nm_{E/k}:E^\times\to k^\times\bigr)
=\{b^{q-1}:b\in E^\times\}.
\]
\end{theorem}
\begin{proof}
For every \(b\), \(\Tr(\sigma(b)-b)=0\), by cancellation in the sum.
The kernel of \(\sigma-\id\) is \(k\), so its image has \(k\)-dimension \(n-1\).
The trace kernel has the same dimension; hence the two coincide.

Similarly, \(\Nm(\sigma(b)/b)=1\). The homomorphism \(b\mapsto b^{q-1}\) has kernel \(k^\times\),
so its image has \((q^n-1)/(q-1)\) elements, exactly the size of the norm kernel.
\end{proof}

\begin{example}[A Frobenius difference equation with trace-zero right-hand side]
Reuse \(E=\F_9=\F_3(\alpha)\), \(\alpha^2=-1\), from Example~\ref{v100:ex-three-generation-nine}; the new question concerns the image and fibers of the Frobenius difference.
For \(x=u+v\alpha\), we have
\[
x^3-x=(u-v\alpha)-(u+v\alpha)=v\alpha.
\]
Thus the image is the one-dimensional space \(\F_3\alpha\), while \(\ker(x\mapsto x^3-x)=\F_3\).
Also, \(\Tr(u+v\alpha)=2u\), so the trace-zero space is likewise \(\F_3\alpha\).
The solutions of \(x^3-x=\alpha\) are \(\alpha,\alpha+1,\alpha+2\),
while \(x^3-x=1\) has no solution.
This exhibits both the image identity in Hilbert 90 and the coset structure of each nonempty fiber.
\end{example}

Thus \(x^q-x=a\) is solvable in \(E\) if and only if \(\Tr_{E/k}(a)=0\), and then has \(q\) solutions.
For \(a\ne0\), the equation \(x^{q-1}=a\) is solvable if and only if \(\Nm_{E/k}(a)=1\), and then has \(q-1\) solutions.
These equations turn additive orbit differences and multiplicative orbit ratios, respectively, into one condition over the base field.

\begin{proposition}[A fractional-linear parametrization of the norm-one group]\label{v85:prop-norm-one-parametrization}
Let \(E=\F_{q^2}\), \(k=\F_q\), and choose \(\theta\in E\setminus k\). The norm-one group is
\[
 \mu_{q+1}=\{u\in E^\times:u^{q+1}=1\}=\ker\Nm_{E/k}.
\]
The M\"obius transformation
\[
 \phi_\theta:k\longrightarrow\mu_{q+1}\setminus\{1\},\qquad
 a\longmapsto\frac{\theta^q-a}{\theta-a}
\]
is a bijection, with inverse \(u\mapsto(u\theta-\theta^q)/(u-1)\). Setting \(\phi_\theta(\infty)=1\) extends it to a bijection from \(k\cup\{\infty\}\) onto the entire norm-one group.
For any finite set \(S\subseteq k\), put \(p_S(T)=\prod_{a\in S}(T-a)\). Then
\begin{equation}\label{v85:eq-norm-root-product}
 p_S(\theta)^{q-1}=\prod_{a\in S}\phi_\theta(a).
\end{equation}
\end{proposition}
\begin{proof}
The denominator is nonzero, and \(\phi_\theta(a)^q=\phi_\theta(a)^{-1}\), so its value lies in \(\mu_{q+1}\) and is not one.
If two parameters have the same image, cross-multiplication gives \((\theta^q-\theta)(a-b)=0\), hence \(a=b\).
Both sets have \(q\) elements, so this injection is a bijection. Solving for \(a\) gives the inverse formula.
Finally, \(p_S\in k[T]\) implies \(p_S(\theta)^q=p_S(\theta^q)\). Expanding both values as products of root factors gives \eqref{v85:eq-norm-root-product}.
\end{proof}
In \(\F_9=\F_3(\alpha)\), where \(\alpha^2=-1\), we have \(\alpha^3=-\alpha\), and therefore
\[
 \phi_\alpha(0)=-1,\qquad \phi_\alpha(1)=\alpha,\qquad
 \phi_\alpha(-1)=-\alpha.
\]
Together with \(\phi_\alpha(\infty)=1\), these values give \(\mu_4\). Ratios of conjugates produce norm-one elements, and evaluating a root polynomial produces their product.

\begin{proposition}[Small-exponent Frobenius models and additive decompositions]\label{v150:prop-small-frobenius}
Let $H\le\F_q^\times$ have order $t>1$, and let $p=\Char\F_q$. Suppose $e\ge1$ and an integer $r$ satisfy
\[
 p^e\equiv r\pmod t,\qquad r\ne1,\qquad p\nmid r.
\]
Put $M(r)=r-1$ for $r\ge2$, and $M(r)=2|r|$ for $r\le-1$. For distinct $a,b\in\F_q$,
\begin{equation}\label{v150:eq-translate-intersection}
 |(H-a)\cap(H-b)|\le M(r).
\end{equation}
If $t>M(r)^2$, no subsets $A,B\subseteq\F_q$ with $|A|,|B|\ge2$ satisfy $H=A+B$.
\end{proposition}
\begin{proof}
An intersection point $x$ has $x+a,x+b\in H$, and therefore
\[
 (x+a)^r-(x+b)^r=(a-b)^{p^e}.
\]
For $r\ge2$, subtracting the right side gives a polynomial of degree $r-1$ with leading coefficient $r(a-b)\ne0$.
For $r=-h<0$, clearing denominators gives
\[
 (X+b)^h-(X+a)^h-(a-b)^{p^e}(X+a)^h(X+b)^h,
\]
a nonzero polynomial of degree $2h$. The root bound proves~\eqref{v150:eq-translate-intersection}.
If $H=A+B$, two distinct elements of $A$ give $|B|\le M(r)$; reversing the roles gives $|A|\le M(r)$.
Thus $t=|A+B|\le |A||B|\le M(r)^2$, a contradiction.
\end{proof}

\begin{example}[Additive indecomposability of a norm-one group]\label{v150:ex-norm-one-indecomposable}
Let $Q$ be a prime power and $H=\mu_{Q+1}\subseteq\F_{Q^2}^\times$. Then $H$ is the sum of two sets of size at least two if and only if $Q=3$.
For $Q\ge4$, the congruence $Q\equiv-1\pmod{Q+1}$ allows $r=-1$, bounding the intersection of two distinct translates by two. Since $|H|=Q+1>4$, no such decomposition exists.
For $Q=2$, a decomposition $H=A+B$ with $|A|,|B|\ge2$ would force $|A|=|B|=2$ by the same intersection bound. In characteristic two the sum of two two-element sets has either two or four elements, whereas $|H|=3$.
For $Q=3$, in $\F_9=\F_3(i)$ with $i^2=-1$,
\[
 H=\{\pm1,\pm i\}=\{1,i\}+\{0,-1-i\}.
\]
In particular, the subgroup of order five in $\F_{16}$ is additively indecomposable. Frobenius acts on $H$ by inversion; a quadratic equation first bounds translate intersections, which then exclude a global additive decomposition.
\end{example}

The norm of the translation parameter determines the intersection count exactly.
\begin{proposition}[Exact intersections of translates of the norm-one group]\label{v151:prop-normone-intersections}
Let $Q$ be a prime power, $K=\F_{Q^2}$, and $H=\ker\Nm_{K/\F_Q}$. For $\delta\in K^\times$, put
\[
 D=\Nm_{K/\F_Q}(\delta),\qquad
 r_H(\delta)=\#\{(x,y)\in H^2:x-y=\delta\}=|H\cap(H+\delta)|.
\]
When $Q$ is odd, let $\chi$ be the quadratic character of $\F_Q$, with $\chi(0)=0$. Then
\begin{equation}\label{v151:eq-normone-intersections}
 r_H(\delta)=
 \begin{cases}
 1-\chi\bigl(D(D-4)\bigr),&Q\text{ odd},\\[2pt]
 1-(-1)^{\Tr_{\F_Q/\F_2}(D^{-1})},&Q\text{ even}.
 \end{cases}
\end{equation}
In particular, in odd characteristic the count is one precisely when $D=4$; in every other case it is zero or two.
\end{proposition}
\begin{proof}
Write $x=\delta z$ and $y=\delta(z-1)$. Since
\[
 \Nm(z-1)=\Nm(z)-\Tr_{K/\F_Q}(z)+1,
\]
the conditions $x,y\in H$ are equivalent to
\[
 \Tr_{K/\F_Q}(z)=1,\qquad \Nm_{K/\F_Q}(z)=D^{-1}.
\]
Such an element is a root of $Z^2-Z+D^{-1}$. If this polynomial is irreducible, both conjugate roots have the required trace and norm. If it has two distinct roots in $\F_Q$, neither qualifies: a base-field element $z$ has relative trace $2z$, and $2z=1$ can only give a repeated root. In odd characteristic the discriminant is $1-4/D$; a repeated root occurs exactly when $D=4$ and contributes one solution. Since $D^2$ is a nonzero square, $\chi(1-4/D)=\chi(D(D-4))$, proving the first formula.

In even characteristic the quadratic is separable. Additive Hilbert 90 shows that it is irreducible precisely when $\Tr_{\F_Q/\F_2}(D^{-1})=1$, proving the second formula.
\end{proof}

Thus all translation parameters in the same norm fiber give the same intersection count. Section~\ref{v58:sec-character-moments} uses these counts to determine a fourth moment of character sums.

\begin{proposition}[Short zero-sum relations in the norm-one group]\label{v153:prop-normone-short-sums}
Let $Q$ be a prime power and $H=\ker\Nm_{\F_{Q^2}/\F_Q}$.
\begin{enumerate}
\item For $s\ne0$, there is at most one unordered pair $\{x,y\}$ with $x,y\in H$ and $x+y=s$, allowing repetition, and
\[
 xy=s/s^Q.
\]
\item In characteristic three, the zero-sum triples in $H$ are precisely $(x,x,x)$. In other characteristics their entries are distinct and form a coset $x\mu_3\subseteq H$; such triples exist if and only if $3\mid Q+1$.
\item Four distinct elements cannot sum to zero in characteristic two. In odd characteristic a zero-sum four-element subset is precisely the union of two pairs of opposite elements; there are $\binom{(Q+1)/2}{2}$ such subsets.
\end{enumerate}
\end{proposition}
\begin{proof}
Since $x^Q=x^{-1}$,
\[
 s^Q=x^{-1}+y^{-1}=s/(xy).
\]
Thus $x,y$ are determined by $T^2-sT+s/s^Q$, proving the first assertion.
If $x+y+z=0$, put $a=x/z$, $b=y/z$. Raising $a+b=-1$ to the $Q$th power gives $a^{-1}+b^{-1}=-1$, whence $ab=1$. Thus $a,b$ are the two roots of $T^2+T+1$. In characteristic three both equal one; otherwise they are the primitive cube roots of unity, which lie in $H$ if and only if $3\mid Q+1$.

Suppose that $x_1+x_2+x_3+x_4=0$ with distinct entries. In characteristic two, the first assertion applied to $x_1+x_2=x_3+x_4\ne0$ gives a contradiction. In odd characteristic, if $x_1+x_2=0$, the required pairing is already present; otherwise apply the first assertion to $x_1+x_2=(-x_3)+(-x_4)$. The converse is immediate. There are $(Q+1)/2$ opposite pairs in $H$, and choosing two gives the count.
\end{proof}

\section{Trace, Differences, and Equivariant Decomposition}
\label{sec:trace-exact}

Put \(D=\sigma-\id\). Additive Hilbert 90 and surjectivity of the trace together give
\begin{equation}\label{eq:trace-exact}
0\longrightarrow k\longrightarrow E
\xrightarrow{\ D\ }E\xrightarrow{\ \Tr_{E/k}\ }k\longrightarrow0.
\end{equation}
Exactness means that the image of each map is the kernel of the next. This expresses two facts:
Frobenius differences annihilate base-field constants, and an element is a Frobenius difference if and only if its trace is zero.
Let \(\sigma\) act on the two copies of \(E\), and let the identity act on the two copies of \(k\). All maps above are then equivariant.
In particular, equivariance of \(P:E\to k\) means \(P\sigma=P\); viewing its values in \(E\), this is also equivalent to \(P\sigma=\sigma P\).

\begin{proposition}\label{prop:trace-splitting}
Regard \(\Tr=\Tr_{E/k}\) as a linear operator taking values in \(E\). Then \(\Tr^2=n\Tr\).
If \(p\nmid n\), there is a Frobenius-equivariant direct-sum decomposition
\[
E=k\oplus\ker\Tr,
\qquad P=\frac1n\Tr:E\to k,
\]
and \(D\) is invertible on \(\ker\Tr\). If \(p\mid n\), then \(k\subseteq\ker\Tr=\im D\), and no \(k\)-linear map \(E\to k\) commuting with \(\sigma\) restricts to the identity on \(k\).
\end{proposition}
\begin{proof}
For \(c\in k\), we have \(\Tr(c)=nc\). Since \(\Tr(x)\in k\), this gives \(\Tr(\Tr(x))=n\Tr(x)\). If \(p\nmid n\), then \(P^2=P\) and \(P|_k=\id\),
while \(P\sigma=P=\sigma P\), giving the direct sum. The space \(\ker\Tr\) is \(D\)-stable, and
\[
\ker(D|_{\ker\Tr})=k\cap\ker\Tr=\{0\}.
\]
Finite dimensionality gives invertibility.

If \(p\mid n\), every base-field element has trace zero, so \(k\subseteq\im D\).
Suppose such a map \(P\) exists. Since its values lie in \(k\), commutation means \(P\sigma=P\), and hence \(PD=0\).
But \(1\in\im D\), so \(P(1)=0\), contradicting \(P|_k=\id\).
\end{proof}

When \(p\mid n\), the \(k\)-vector space \(E\) still admits a complement to \(k\), but no such complement is Frobenius-stable.
In this case, the nonzero operator \(\Tr\) satisfies \(\Tr^2=0\). This phenomenon will be further explained by the repeated factors of the Frobenius minimal polynomial \(T^n-1\).

\section{Artin--Schreier Equations and Binomials of Prime Degree}
\begin{theorem}\label{thm:artin-schreier}
Let \(k=\F_q\), \(q=p^r\), and \(a\in k\). The polynomial \(X^p-X-a\) satisfies exactly one of the following:
\begin{enumerate}
\item If \(\Tr_{k/\F_p}(a)=0\), it splits completely over \(k\).
\item If \(\Tr_{k/\F_p}(a)\ne0\), it is irreducible in \(k[X]\).
\end{enumerate}
\end{theorem}
\begin{proof}
Choose a root \(\alpha\) in a finite extension. Iterating \(\alpha^p-\alpha=a\) and summing gives
\[
\alpha^q-\alpha=\sum_{i=0}^{r-1}a^{p^i}=c,\qquad c\in\F_p.
\]
If \(c=0\), then \(\alpha\in k\), and all roots \(\alpha+t\), \(t\in\F_p\), lie in \(k\).
If \(c\ne0\), then \(\alpha^{q^j}=\alpha+jc\), whose least positive period is \(p\).
Thus \(\alpha\) has degree \(p\) over \(k\), and the original monic degree-\(p\) polynomial is its minimal polynomial.
\end{proof}

\begin{proposition}\label{prop:prime-binomial}
Let \(\ell\) be prime and \(a\in\F_q^\times\). Then \(X^\ell-a\) is irreducible over \(\F_q\)
if and only if \(\ell\mid q-1\) and \(a\) is not an \(\ell\)-th power in \(\F_q\).
\end{proposition}
\begin{proof}
If \(\ell\nmid q-1\), the map \(x\mapsto x^\ell\) permutes \(\F_q^\times\), so the polynomial has a base-field root and is reducible.
Suppose \(\ell\mid q-1\). All \(\ell\)-th roots of unity lie in the base field.
For a root \(\alpha\), we have \((\alpha^{q-1})^\ell=a^{q-1}=1\), so \(\alpha^q=\zeta\alpha\), where \(\zeta\in\F_q\) and \(\zeta^\ell=1\).
If \(a\) is not an \(\ell\)-th power, then \(\alpha\notin\F_q\), so \(\zeta\ne1\).
Since \(\ell\) is prime, \(\zeta\) has order \(\ell\), and the Frobenius period of \(\alpha\) is \(\ell\), proving irreducibility.
If \(a\) is an \(\ell\)-th power, the polynomial has a base-field root and is not irreducible.
\end{proof}

\par\Needspace{12\baselineskip}
The norm also constrains all irreducible factors of a binomial, without requiring the binomial itself to be irreducible.
\begin{proposition}[A norm class and the degrees of binomial factors]\label{v97:prop-norm-factor-degrees}
Let \(n\ge1\) and \(\lambda\in\F_q^\times\), and put
\[
 d=\ord\bigl([\lambda]\in\F_q^\times/(\F_q^\times)^n\bigr).
\]
The degree of every irreducible factor of \(X^n-\lambda\) is divisible by \(d\). This also holds when \(p\mid n\).
\end{proposition}
\begin{proof}
If \(\alpha\) is a root of an irreducible factor of degree \(r\), then \(\alpha\in\F_{q^r}^\times\) and \(\alpha^n=\lambda\). Taking norms gives
\[
 \Nm_{\F_{q^r}/\F_q}(\alpha)^n=\lambda^r.
\]
Thus \([\lambda]^r=1\), and hence \(d\mid r\).
\end{proof}

This is a common divisibility constraint, not a complete classification of factor degrees.
For example, over \(\F_3\), the class of \(-1\) in \(\F_3^\times/(\F_3^\times)^8\) has order \(2\).
Nevertheless, every root of \(X^8+1\) has multiplicative order \(16\), so \(\ord_{16}(3)=4\) shows that every irreducible factor has degree \(4\).
Example~\ref{v97:ex-binomial-support} will explain the additional splitting extension through coefficient support.

\section{Explicit Solutions of Hilbert 90 Equations}
\begin{proposition}[An explicit formula for additive Hilbert 90]\label{prop:hilbert90-explicit}
Let \(E=\F_{q^n}\), \(\sigma(x)=x^q\), and choose \(\theta\in E\) with \(\Tr_{E/k}(\theta)=1\).
For \(\Tr_{E/k}(a)=0\), put
\[
B=\sum_{i=1}^{n-1}\left(\sum_{j=0}^{i-1}\sigma^j(a)\right)\sigma^i(\theta).
\]
Then \(x=-B\) satisfies \(x^q-x=a\), and all solutions form \(-B+k\). For \(n=1\), the empty sum is zero.
\end{proposition}
\begin{proof}
When \(n=1\), the trace condition gives \(a=0\), and the empty-sum formula is immediate. Assume \(n\ge2\).
Put \(S_0=0\) and \(S_i=\sum_{j=0}^{i-1}\sigma^j(a)\). Then \(\sigma(S_i)=S_{i+1}-a\) and \(S_n=\Tr_{E/k}(a)\).
For an arbitrary value of \(\Tr(a)\), we have
\[
\begin{aligned}
\sigma(B)
&=\sum_{i=1}^{n-1}(S_{i+1}-a)\sigma^{i+1}(\theta)\\
&=\sum_{j=2}^{n-1}S_j\sigma^j(\theta)
  +S_n\theta-a\sum_{j=2}^{n}\sigma^j(\theta).
\end{aligned}
\]
On the other hand, \(B=a\sigma(\theta)+\sum_{j=2}^{n-1}S_j\sigma^j(\theta)\).
Subtracting cancels the intermediate terms; using \(\sigma^n(\theta)=\theta\), we obtain
\[
B-\sigma(B)
=a\sum_{j=1}^{n}\sigma^j(\theta)-\theta S_n
=a\Tr_{E/k}(\theta)-\theta\Tr_{E/k}(a).
\]
Under the two trace hypotheses, the right side equals \(a\). Thus \(\sigma(-B)-(-B)=a\). Any two solutions differ by an element of \(\ker(\sigma-\id)=k\), and adding any element of \(k\) to a solution gives another solution. Hence the full solution set is \(-B+k\).
\end{proof}
The equation \(\sigma(x)-x=a\) is a difference equation on a conjugate orbit. Trace zero means that the differences sum to zero around one circuit; an element of trace one supplies the normalization for the explicit solution formula.

\Needspace{5\baselineskip}
\begin{proposition}[An explicit formula for multiplicative Hilbert 90]\label{prop:multiplicative-hilbert-explicit}
Let \(a\in E^\times\) satisfy \(\Nm_{E/k}(a)=1\), and put \(c_0=1\), \(c_i=\prod_{j=0}^{i-1}\sigma^j(a)\).
One can choose \(\theta\in E\) such that
\[
B=\sum_{i=0}^{n-1}c_i\sigma^i(\theta)\ne0.
\]
Then \(a=B/\sigma(B)\), so \(x=B^{-1}\) satisfies \(x^{q-1}=a\).
\end{proposition}
\begin{proof}
The relations \(a\sigma(c_i)=c_{i+1}\) and \(c_n=1\) give \(a\sigma(B)=B\).
As a polynomial in \(\theta\), the expression \(B\) has linear coefficient one and degree at most \(q^{n-1}<q^n\),
so it cannot vanish on all of \(E\). Thus some \(\theta\) gives \(B\ne0\).
\end{proof}
The additive and multiplicative formulas correspond to orbit differences and orbit ratios; trace zero and norm one are the respective compatibility conditions over a full orbit.

\begin{example}[Half-trace]
Let \(E=\F_{2^n}\), with \(n\) odd. Define
\[
H(a)=\sum_{j=0}^{(n-1)/2}a^{2^{2j}}.
\]
Adding the even and odd Frobenius iterates gives \(H(a)^2+H(a)=a+\Tr_{E/\F_2}(a)\).
Thus, when the absolute trace is zero, \(H(a)\) is a solution of \(X^2+X=a\).
\end{example}

\begin{proposition}[A unified criterion for quadratic equations]\label{prop:quadratic-test}
Let \(A,B,C\in k\), with \(A\ne0\). In odd characteristic, \(AX^2+BX+C\) has a root in \(k\) if and only if \(B^2-4AC\) is a square, including zero.
In characteristic two, if \(B=0\), it has a unique root in \(k\); if \(B\ne0\), it has a root if and only if
\[
\Tr_{k/\F_2}(AC/B^2)=0.
\]
In the latter case, there are two distinct roots.
\end{proposition}
\begin{proof}
In odd characteristic, complete the square; in characteristic two with \(B=0\), use bijectivity of the squaring map.
For \(B\ne0\), put \(Y=AX/B\). The equation becomes \(Y^2+Y=AC/B^2\); apply additive Hilbert 90.
\end{proof}

\section{Trace Fibers and Polynomial Factorization}
\begin{proposition}\label{prop:trace-fiber-factorization}
Let \(r\) be a prime power, \(q=r^m\), and \(b\in\F_r\). Then
\[
X^q-X-b=
\prod_{\substack{\beta\in\F_q\\ \Tr_{\F_q/\F_r}(\beta)=b}}
(X^r-X-\beta).
\]
When \(r=p\) and \(b\ne0\), every factor on the right is irreducible of degree \(p\) over \(\F_q\).
\end{proposition}
\begin{proof}
Compare the root sets in a common finite extension in which the polynomials split. If \(x^q-x=b\), put \(\beta=x^r-x\). Then
\[
\beta^q-\beta
=(x^q-x)^r-(x^q-x)=b^r-b=0.
\]
Thus \(\beta\in\F_q\). We have the telescoping identity
\[
\sum_{i=0}^{m-1}\beta^{r^i}
=\sum_{i=0}^{m-1}(x^{r^{i+1}}-x^{r^i})
=x^q-x=b.
\]
This gives \(\Tr_{\F_q/\F_r}(\beta)=b\), and \(x\) is a root of the corresponding factor on the right.
Conversely, if \(x^r-x=\beta\) and this trace equals \(b\), the same identity gives \(x^q-x=b\).

The derivative of the left side is \(-1\), as is that of every factor on the right, so each has no repeated roots.
Factors with different \(\beta\) cannot share a root, since the same \(x^r-x\) cannot equal two different values.
Surjectivity of the trace gives \(q/r\) indices, so the right side has degree \(r(q/r)=q\).
Both sides are monic of the same degree and have the same distinct roots, so they are equal.
When \(r=p\) and \(b\ne0\), each \(\beta\) has nonzero absolute trace,
so Theorem~\ref{thm:artin-schreier} makes every \(X^p-X-\beta\) irreducible over \(\F_q\).
\end{proof}

\begin{theorem}[Affine equations of prime degree]\label{thm:affine-prime}
Let \(k=\F_q\), \(\Char k=p\), \(a\in k^\times\), and \(b\in k\).
The polynomial \(X^p-aX-b\) is irreducible if and only if there exists \(u\in k^\times\) satisfying
\[
a=u^{p-1},\qquad \Tr_{k/\F_p}(b/u^p)\ne0.
\]
\end{theorem}
\begin{proof}
First determine whether \(u^{p-1}=a\) is solvable in \(k^\times\).
If not, the \(\F_p\)-linear map \(L:k\to k\), \(L(x)=x^p-ax\), has trivial kernel:
a nonzero kernel element would satisfy \(x^{p-1}=a\). An injection of a finite set into itself is surjective, so some \(x\in k\) satisfies \(L(x)=b\).
The original polynomial has a root in \(k\) and degree \(p\ge2\), so it is reducible.

If \(u\in k^\times\) satisfies \(u^{p-1}=a\), fix such a \(u\), and make the invertible substitution \(X=uY\), obtaining
\[
X^p-aX-b=u^p\bigl(Y^p-Y-b/u^p\bigr).
\]
Substitution and multiplication by a nonzero constant preserve irreducibility. By Theorem~\ref{thm:artin-schreier},
the right side is irreducible if and only if
\[
\Tr_{k/\F_p}(b/u^p)\ne0.
\]
This gives the stated necessary and sufficient condition.
\end{proof}
If \(a=0\), the polynomial is the \(p\)-th power of a linear polynomial and is also not irreducible.

\par\Needspace{12\baselineskip}
A trace condition alone gives an additive fiber; combining trace and norm also involves multiplicative structure. Corollary~\ref{v145:cor-cubic-trace-norm} identifies $\Nm(a)-\Tr(a)$ in a cubic extension with the negative trace of a two-dimensional Frobenius return, and character sums determine the exact fibers over irreducible quadratic targets.

\section{Conjugate Coordinates, Tensor Product Decomposition, and Multiplication Operators}
For \(a\in E\), let \(M_a:E\to E\) be multiplication \(x\mapsto ax\), a \(k\)-linear operator.
\begin{theorem}[Conjugate coordinates]\label{thm:conjugate-coordinates}
Let \(b_0,\ldots,b_{n-1}\) be a basis of \(E/k\), and put
\[
V=(b_j^{q^i})_{0\leq i,j<n}.
\]
Then \(V\) is invertible over \(E\). For \(a\in E\), let \(M_a\) denote the matrix of multiplication \(x\mapsto ax\) in this basis. Then
\[
V M_a=\operatorname{diag}(a,a^q,\ldots,a^{q^{n-1}})V.
\]
Thus all multiplication matrices are simultaneously diagonalizable after extending scalars to \(E\).
\end{theorem}
\begin{proof}
First take the power basis \(1,\theta,\ldots,\theta^{n-1}\) of \(E=k(\theta)\). Then \(V\) is a Vandermonde matrix
with distinct nodes \(\theta,\theta^q,\ldots,\theta^{q^{n-1}}\), so its determinant is nonzero.
Any basis is obtained from the power basis by a matrix in \(\mathrm{GL}_n(k)\); the corresponding \(V\) is multiplied on the right by that matrix, and remains invertible.
Apply \(\sigma^i\) to \(ab_j=\sum_r(M_a)_{rj}b_r\). Since the matrix entries belong to \(k\), this gives the asserted matrix identity.
\end{proof}

The tensor product decomposition in Theorem~\ref{thm:tensor-fields} now takes the concrete form
\begin{equation}\label{eq:conjugate-tensor-map}
\mathcal J:E\otimes_kE\longrightarrow E^n,\qquad
c\otimes x\longmapsto\bigl(c\sigma^i(x)\bigr)_{i=0}^{n-1}.
\end{equation}
The domain is an \(E\)-vector space through its first tensor factor. Relative to the basis \((1\otimes b_j)_j\),
the matrix of \(\mathcal J\) is exactly \(V\), so the theorem also directly proves that it is an \(E\)-algebra isomorphism.
Evaluation on pure tensors gives
\begin{equation}\label{eq:conjugate-operator-models}
\begin{aligned}
\mathcal J(\id\otimes M_a)\mathcal J^{-1}
 &=\operatorname{diag}\bigl(a,\sigma(a),\ldots,\sigma^{n-1}(a)\bigr),\\
\mathcal J(\id\otimes\sigma)\mathcal J^{-1}
 &=P,\qquad P(z_0,\ldots,z_{n-1})=(z_1,\ldots,z_{n-1},z_0).
\end{aligned}
\end{equation}
The first identity follows from \(\sigma^i(ax)=\sigma^i(a)\sigma^i(x)\), and the second from \(\sigma^n=\id\).
Here \(\id\otimes\sigma\) is the scalar extension of the \(k\)-linear operator \(\sigma\), fixing coefficients in the first tensor factor.
Thus the same evaluation isomorphism represents field multiplication as coordinatewise multiplication and Frobenius as a cyclic shift of conjugate coordinates.
These models retain the multiplicative relations and the order of conjugation, respectively.
Simultaneous diagonalization of multiplication operators over \(E\) follows from separability of the field extension.
Equation~\eqref{eq:conjugate-operator-models} also directly determines the minimal polynomial of Frobenius as a \(k\)-linear operator.
For a standard basis vector \(e\) of \(E^n\), the vectors \(e,Pe,\ldots,P^{n-1}e\) are a permutation of the full standard basis,
so no nonzero polynomial of degree less than \(n\) annihilates \(P\).
If such a \(k\)-polynomial annihilated \(\sigma\), it would annihilate \(\id\otimes\sigma\), and hence \(P\), a contradiction.
Together with \(\sigma^n=\id\), this shows that the minimal polynomial of \(\sigma\) is exactly \(T^n-1\).
Over the splitting field of \(T^n-1\), it is diagonalizable if and only if \(p\nmid n\).
For diagonalizability over the base field \(k\), all \(n\)-th roots of unity must also belong to \(k\); the combined condition is \(n\mid q-1\).
Thus separability of the field extension, absence of repeated roots in the operator's minimal polynomial, and splitting of that polynomial over a specified scalar field must be considered separately.
For example, Frobenius on \(\F_8/\F_2\) satisfies the square-free polynomial \(T^3-1=(T+1)(T^2+T+1)\),
but this polynomial splits over neither \(\F_2\) nor \(\F_8\), so the operator is diagonalizable over neither scalar field.

\begin{proposition}\leavevmode\phantomsection\label{prop:characteristic-multiplication}
If \(d=[k(a):k]\), the characteristic polynomial of \(M_a\) is
\[
\det(TI-M_a)=m_{a,k}(T)^{n/d}
=\prod_{i=0}^{n-1}(T-a^{q^i}).
\]
Thus its matrix trace and determinant are \(\Tr_{E/k}(a)\) and \(\Nm_{E/k}(a)\), respectively.
\end{proposition}
\begin{proof}
By Theorem~\ref{thm:conjugate-coordinates}, the multiplication matrix is similar over \(E\) to the diagonal matrix of all Frobenius conjugate values,
so its characteristic polynomial is the displayed product. The orbit of \(a\) has length \(d\), and each distinct conjugate occurs exactly \(n/d\) times among the \(n\) terms.
The product is therefore \(m_{a,k}^{n/d}\). Comparing the coefficient of \(T^{n-1}\) and the constant term gives the matrix trace and determinant as the field trace and norm.
\end{proof}

In particular, if \(a\) generates the whole extension, the minimal polynomial \(T^n+c_{n-1}T^{n-1}+\cdots+c_0\) has coefficients \(c_{n-1}=-\Tr(a)\) and \(c_0=(-1)^n\Nm(a)\).
For elements not generating the whole extension, first use \eqref{eq:trace-subdegree}.

\section{Root Ratios and Coefficient Criteria for Primitive Quadratics}
\label{v45:sec-quadratic-criterion}

For a monic irreducible quadratic, the trace and norm of a root determine the linear and constant coefficients, while the multiplicative order of the ratio of its conjugate roots is its projective order.
Translating the multiplicative order of this ratio into coefficient conditions gives primitivity criteria, counts with fixed constant term, and the minimal binomial multiple.
For related descriptions of primitive quadratics, see \cite{VegaQuadratic,VegaEasy};
for background on Lucas atoms, see \cite{SaganTirrellLucasAtoms}.
The following construction of optimal determining polynomials over finite fields is given in \cite{WZOptimalQuadratic}.
Throughout this section, fix \(k=\F_q\) of characteristic \(p\), and put \(N=q+1\).

\subsection{Root Ratios, Recurrence Remainders, and Projective Order}

\leavevmode\phantomsection\label{v45:def-binomial-order}\ignorespaces
In the context of binomial multiples, projective order is also called \textbf{binomial order};
the \emph{quasi-order} in \cite{VegaQuadratic} has the same definition.

\begin{proposition}[Quadratic primitivity and root ratios]\label{v45:prop-quadratic-root-ratio}
Let \(f=X^2+bX+c\in k[X]\), with \(c\ne0\), and let its roots be \(\theta_1,\theta_2\).
Then \(f\) is primitive if and only if \(c\) generates \(k^\times\) and
\[
 \ord(\theta_1/\theta_2)=q+1.
\]
\end{proposition}
\begin{proof}
If \(f\) is irreducible, its roots are \(\theta,\theta^q\), and \(\theta^{q+1}=c\).
Proposition~\ref{v40:prop-minimal-binomial} gives its projective order as \(\ord(\theta/\theta^q)\);
Proposition~\ref{prop:projective-primitive} then gives the equivalent conditions for primitivity.
Conversely, if the root ratio has order \(q+1\), then \(f\) is irreducible: otherwise both roots lie in \(k\),
and their ratio has order dividing \(q-1\), which cannot equal \(q+1\). The equivalence above then applies.
\end{proof}

\begin{definition}\label{v45:def-lucas-polynomials}
Define the \textbf{Lucas polynomials of the first kind} in \(k[B,C]\) by
\[
 U_0=0,\qquad U_1=1,\qquad
 U_n=-BU_{n-1}-CU_{n-2}\quad(n\ge2).
\]
For fixed \((b,c)\), put \(u_n=U_n(b,c)\); the sequence \((u_n)\) is the corresponding Lucas sequence.
\end{definition}

\begin{proposition}[Recurrence remainders and the first zero term of positive index]\label{v45:prop-lucas-first-zero}
Suppose \(f=X^2+bX+c\) is irreducible. For \(n\ge1\),
\begin{equation}\label{v45:eq-lucas-remainder}
 X^n\equiv u_nX-cu_{n-1}\pmod f.
\end{equation}
If \(\theta\) is a root, then
\begin{equation}\label{v45:eq-lucas-binet}
 u_n=\frac{\theta^n-\theta^{nq}}{\theta-\theta^q}\quad(n\ge0).
\end{equation}
Put \(r=\min\{n\ge1:u_n=0\}\). Then \(r\ge2\), \(r\mid q+1\), and \(r\) is the projective order of \(f\).
Writing \(\lambda=-cu_{r-1}\), its minimal binomial multiple and ordinary order are, respectively,
\begin{equation}\label{v45:eq-lucas-binomial-order}
 B_f(X)=X^r-\lambda=X^r+cu_{r-1},
 \qquad \ord(f)=r\ord(\lambda).
\end{equation}
\end{proposition}
\begin{proof}
Equation~\eqref{v45:eq-lucas-remainder} holds for \(n=1\).
If it holds for \(n\), multiply by \(X\) and use \(X^2\equiv-bX-c\pmod f\) to obtain
\[
 X^{n+1}\equiv(-bu_n-cu_{n-1})X-cu_n
 =u_{n+1}X-cu_n\pmod f,
\]
completing the induction. The right side of \eqref{v45:eq-lucas-binet} has initial values zero and one,
and, since both roots satisfy \(\theta^2=-b\theta-c\), it satisfies the same recurrence. Hence it equals \(u_n\).

Put \(\rho=\theta/\theta^q\). Since \(\rho^q=\rho^{-1}\), its order divides \(q+1\).
Equation~\eqref{v45:eq-lucas-binet} gives \(u_n=0\iff\rho^n=1\), so the least positive index at which \(u_n\) vanishes is
\(r=\ord(\rho)\ge2\).
The remainder formula shows that \(X^n\) becomes a base-field constant modulo \(f\) if and only if \(u_n=0\).
At \(n=r\), that constant is \(\lambda=\theta^r=-cu_{r-1}\ne0\).
This gives the stated minimal binomial multiple; the ordinary-order formula follows from the properties of projective order.
\end{proof}

For this irreducible \(f\), the recurrence terms are also coefficients of a long-division quotient. Since \(\theta^{q+1}=c\), we have \(f\mid X^{q+1}-c\), and
\begin{equation}\label{v45:eq-first-zero-quotient}
 \frac{X^{q+2}-cX}{f(X)}
 =h_1X^q+h_2X^{q-1}+\cdots+h_qX+h_{q+1}.
\end{equation}
Comparing coefficients in the product gives \(h_1=1\), \(h_2=-b\), and
\(h_i=-bh_{i-1}-ch_{i-2}\) for \(3\le i\le q+1\); hence \(h_i=u_i\).
Also, \eqref{v45:eq-lucas-binet} gives \(h_q=u_q=-1\) and \(h_{q+1}=u_{q+1}=0\).
Thus the first zero coefficient occurs at \(r\), and the preceding coefficient recovers \(\lambda\).
In particular, \(h_1,\ldots,h_q\) are all nonzero if and only if the projective order is \(q+1\);
then the minimal binomial multiple is \(X^{q+1}-c\), and the ordinary order is \((q+1)\ord(c)\).

\par\Needspace{12\baselineskip}
\subsection{Counting Primitive Quadratics with Fixed Norm}

\begin{proposition}\label{v45:prop-fixed-constant-primitive-count}
Fix a generator \(c\) of \(k^\times\). The number of \(b\in k\) for which \(X^2+bX+c\) is primitive is
\begin{equation}\label{v45:eq-fixed-constant-primitive-count}
 t_q=
 \begin{cases}
 \varphi(q+1),&q\text{ odd},\\
 \varphi(q+1)/2,&q\text{ even}.
 \end{cases}
\end{equation}
\end{proposition}
\begin{proof}
Choose a primitive element \(g\) of \(E=\F_{q^2}\). Then \(h=g^{q+1}\) generates \(k^\times\).
Write \(c=h^a\), with \(\gcd(a,q-1)=1\). All nonzero elements of norm \(c\) are precisely
\[
 \theta_j=g^{a+j(q-1)},\qquad 0\le j<q+1.
\]
Since the exponent is already coprime to \(q-1\), the element \(\theta_j\) is primitive if and only if
\(\gcd(a+j(q-1),q+1)=1\). Modulo \(N=q+1\), this exponent is \(a-2j\).

If \(q\) is even, multiplication by two is invertible on \(\Z/N\Z\), giving exactly \(\varphi(N)\) primitive elements.
If \(q\) is odd, then \(a\) is odd, and \(j\mapsto a-2j\) is exactly two-to-one onto the odd residue classes.
All units modulo \(N\) are odd, giving exactly \(2\varphi(N)\) primitive elements.
None lies in \(k\), and each pairs with \(\theta_j^q\) in a two-element conjugate orbit.
Each orbit gives one quadratic, with constant term the norm \(c\) and linear coefficient \(-\Tr_{E/k}(\theta_j)\).
Dividing the element count by two gives \eqref{v45:eq-fixed-constant-primitive-count}.
\end{proof}

\subsection{Symmetric Cyclotomic Polynomials and Lucas Atoms}

\begin{lemma}[The required two-variable symmetric-polynomial theorem]\label{v103:lem-symmetric-two}
Over any field \(k\), every \(F\in k[U,V]\) satisfying \(F(U,V)=F(V,U)\) can be written uniquely as
\(Q(U+V,UV)\), with \(Q\in k[S,P]\). This includes characteristic two.
\end{lemma}
\begin{proof}
Group monomials under interchange of \(U,V\): diagonal terms are \((UV)^j\), and off-diagonal pairs are
\((UV)^j(U^r+V^r)\), with \(r>0\). Set \(T_0=2,T_1=S\). Then
\[
T_r=ST_{r-1}-PT_{r-2}\quad(r\ge2),\qquad
T_r(U+V,UV)=U^r+V^r.
\]
This proves existence without dividing by two. For uniqueness, use lexicographic order with \(U>V\).
The leading monomial of \((U+V)^a(UV)^b\) is \(U^{a+b}V^b\).
Distinct pairs \((a,b)\) give distinct leading monomials, so the largest such monomial in a nonzero \(Q\) cannot cancel.
The substitution map is therefore injective.
\end{proof}

For every integer \(n\ge3\), map the coefficients of the integer cyclotomic polynomial \(\Phi_n\) into \(k\). Let \(\alpha,\beta\) be independent indeterminates,
and homogenize the cyclotomic polynomial as
\[
 \Phi_n^{\mathrm{hom}}(\alpha,\beta)
 =\beta^{\varphi(n)}\Phi_n(\alpha/\beta).
\]
Since \(\Phi_n\) is self-reciprocal, this homogeneous polynomial is invariant under interchanging \(\alpha,\beta\).
Lemma~\ref{v103:lem-symmetric-two} therefore gives a unique \(\Lambda_n(B,C)\in k[B,C]\) such that
\begin{equation}\label{v45:eq-lucas-atom-definition}
 \Lambda_n\bigl(-\alpha-\beta,\alpha\beta\bigr)
 =\Phi_n^{\mathrm{hom}}(\alpha,\beta).
\end{equation}
\leavevmode\phantomsection\label{v45:def-lucas-atom}\ignorespaces
We call \(\Lambda_n\) the \(n\)-th \textbf{Lucas atom}.
Also put \(\Lambda_2(B,C)=-B\). This definition imposes no restriction on \(p\) and \(n\), and does not require \(\Lambda_n\) to be irreducible in \(k[B,C]\).

\begin{lemma}[Inverse pairing and detection of root ratios]\label{v45:lem-lucas-atom-ratio}
Assume \(n\ge3\) and \(p\nmid n\). The Lucas atom is monic in \(B\), of degree \(\varphi(n)\), and belongs to \(k[B^2,C]\).
In an extension containing all \(n\)-th roots of unity,
\begin{equation}\label{v45:eq-lucas-atom-pairs}
 \Lambda_n(B,C)
 =\prod_{\{\zeta,\zeta^{-1}\}}
 \bigl(B^2-C(2+\zeta+\zeta^{-1})\bigr),
\end{equation}
where the product runs over inverse pairs of primitive \(n\)-th roots of unity.
If \(c\ne0\), and \(X^2+bX+c\) has roots \(\theta_1,\theta_2\), then
\begin{equation}\label{v45:eq-lucas-atom-root-test}
 \Lambda_n(b,c)=0
 \quad\Longleftrightarrow\quad
 \ord(\theta_1/\theta_2)=n.
\end{equation}
Moreover, for every \(c\in k^\times\), \(\Lambda_n(B,c)\) is square-free in odd characteristic, while in characteristic two every root has multiplicity exactly two.
\end{lemma}
\begin{proof}
A primitive \(n\)-th root of unity differs from its inverse because \(n\ge3\). Multiplying the corresponding pair of factors in the homogeneous factorization gives
\[
 (\alpha-\zeta\beta)(\alpha-\zeta^{-1}\beta)
 =(\alpha+\beta)^2-(2+\zeta+\zeta^{-1})\alpha\beta.
\]
This proves \eqref{v45:eq-lucas-atom-pairs} and gives the leading term, degree, and even-exponent property in \(B\).
Substituting the two roots into \eqref{v45:eq-lucas-atom-definition} gives
\[
 \Lambda_n(b,c)=\theta_2^{\varphi(n)}\Phi_n(\theta_1/\theta_2).
\]
Both roots are nonzero, and the roots of \(\Phi_n\) are precisely the primitive \(n\)-th roots of unity, proving \eqref{v45:eq-lucas-atom-root-test}.

Distinct inverse pairs give distinct \(\zeta+\zeta^{-1}\): if two values agree, the corresponding roots of unity are roots of
\(T^2-(\zeta+\zeta^{-1})T+1\), so they belong to the same inverse pair.
Also,
\[
 2+\zeta+\zeta^{-1}=\frac{(\zeta+1)^2}{\zeta}\ne0,
\]
since a primitive \(n\)-th root of unity is not \(-1\).
Thus the nonzero constants in \eqref{v45:eq-lucas-atom-pairs} remain distinct after setting \(C=c\).
In odd characteristic, each quadratic factor has two simple roots, and distinct factors share no root.
In characteristic two, each quadratic factor is the square of a linear polynomial, and distinct factors have distinct roots,
so every root has multiplicity exactly two.
\end{proof}

\begin{proposition}[Atom factorization of the Lucas recurrence]\label{v146:prop-lucas-factorization}
For every integer \(n\ge2\),
\begin{equation}\label{v146:eq-lucas-factorization}
 U_n(B,C)=\prod_{\substack{d\mid n\\d>1}}\Lambda_d(B,C).
\end{equation}
Define the \textbf{Lucas polynomials of the second kind} by
\[
 V_0=2,\qquad V_1=-B,\qquad V_n=-BV_{n-1}-CV_{n-2}\quad(n\ge2).
\]
If \(X^2+bX+c\) has roots \(\theta_1,\theta_2\), then
\[
 V_n(b,c)=\theta_1^n+\theta_2^n,\qquad
 V_n=U_{n+1}-CU_{n-1}\quad(n\ge1).
\]
\end{proposition}
\begin{proof}
For independent variables \(\alpha,\beta\), the recurrence gives
\[
 U_n(-\alpha-\beta,\alpha\beta)
 =\frac{\alpha^n-\beta^n}{\alpha-\beta}
 =\prod_{\substack{d\mid n\\d>1}}\Phi_d^{\mathrm{hom}}(\alpha,\beta).
\]
The factor for \(d=2\) is \(\alpha+\beta\). Injectivity of substitution in Lemma~\ref{v103:lem-symmetric-two} gives~\eqref{v146:eq-lucas-factorization}. The remaining identities follow from their initial values and recurrences.
\end{proof}

For example, \(U_4=-B(B^2-2C)\), whereas \(U_5=B^4-3B^2C+C^2=\Lambda_5\). The second-kind polynomials record the power sums of the same roots. For distinct nonzero roots, the recurrence remainder detects whether the root-ratio order divides \(n\); when \(n\ge3\) and \(p\nmid n\), the atom \(\Lambda_n\) detects whether that order is exactly \(n\).

\begin{corollary}[Characteristic-power indices in the Lucas recurrence]\label{v151:cor-lucas-frobenius}
For $a\ge0$ and $n\ge1$, the following identities hold in $k[B,C]$:
\begin{equation}\label{v151:eq-lucas-frobenius}
 U_{p^a n}(B,C)=
 \begin{cases}
 (B^2-4C)^{(p^a-1)/2}U_n(B,C)^{p^a},&p\ne2,\\[2pt]
 (-B)^{2^a-1}U_n(B,C)^{2^a},&p=2.
 \end{cases}
\end{equation}
Let $b\in k$ and $c\in k^\times$. If $X^2+bX+c$ is separable, then $U_{p^a n}(b,c)=0$ if and only if $U_n(b,c)=0$. If the quadratic has a repeated root $\theta$, then $U_n(b,c)=n\theta^{n-1}$, and its first zero term of positive index occurs at index $p$.
\end{corollary}
\begin{proof}
For independent indeterminates $\alpha,\beta$, the Frobenius identity gives
\[
 \frac{\alpha^{p^a n}-\beta^{p^a n}}{\alpha-\beta}
 =(\alpha-\beta)^{p^a-1}
 \left(\frac{\alpha^n-\beta^n}{\alpha-\beta}\right)^{p^a}.
\]
Under $B=-\alpha-\beta$ and $C=\alpha\beta$, in odd characteristic $(\alpha-\beta)^2=B^2-4C$; in characteristic two $\alpha-\beta=-B$. Injectivity of substitution proves~\eqref{v151:eq-lucas-frobenius}. Separability makes the discriminant, or $b$ in characteristic two, nonzero, so the prefactor does not vanish. The repeated-root formula follows from the recurrence. Since $c\ne0$ implies $\theta\ne0$, its first positive zero index is $p$.
\end{proof}

\begin{example}[Reduction of an atom versus the order of a root ratio]\label{v151:ex-lucas-characteristic}
We have $\Lambda_3(B,C)=B^2-C$. Over $\F_3$, $\Lambda_3(1,1)=0$, whereas
\[
 X^2+X+1=(X-1)^2.
\]
The root ratio is one, whose order is not three. Thus atom factorization holds in every characteristic, while detecting the exact order of a root ratio by $\Lambda_n$ still requires $p\nmid n$. For separable quadratics, characteristic-power indices repeat existing zero conditions; a repeated root instead produces a zero at index $p$.
\end{example}

\subsection{Optimal Coefficient Determining Polynomials}

\leavevmode\phantomsection\label{v45:def-primitive-coefficient-set}\ignorespaces
Fix a primitive element \(c\in k^\times\) as the constant term, and put
\[
 \mathcal P_c=\{b\in k:X^2+bX+c\text{ is primitive}\}.
\]
\leavevmode\phantomsection\label{v45:def-optimal-determining-polynomial}\ignorespaces
A monic polynomial \(P_{q,c}(B)\) is called the \textbf{optimal determining polynomial} of this fiber if
\[
 P_{q,c}(B)=\prod_{b\in\mathcal P_c}(B-b).
\]
It is square-free and has least degree among all nonzero polynomials vanishing on \(\mathcal P_c\).

\begin{theorem}[An optimal coefficient criterion for primitive quadratics]\label{v45:thm-optimal-quadratic-criterion}
Let \(c\) generate \(k^\times\), and put \(N=q+1\).
If \(q\) is odd, then
\begin{equation}\label{v45:eq-optimal-quadratic-odd}
 P_{q,c}(B)=\Lambda_N(B,c),\qquad \deg P_{q,c}=\varphi(N).
\end{equation}
If \(q\) is even, \(\Lambda_N(B,c)\) has a unique square root in \(k[B]\), and
\begin{equation}\label{v45:eq-optimal-quadratic-even}
 P_{q,c}(B)^2=\Lambda_N(B,c),\qquad
 \deg P_{q,c}=\varphi(N)/2.
\end{equation}
In both cases, \(P_{q,c}\) is monic and square-free, all its roots lie in \(k\), and
\[
 X^2+bX+c\text{ is primitive}\quad\Longleftrightarrow\quad P_{q,c}(b)=0.
\]
\end{theorem}
\begin{proof}
Since \(N\ge3\) and \(p\nmid N\), the preceding lemma and Proposition~\ref{v45:prop-quadratic-root-ratio} apply.
They show that, for \(b\in k\), the equality \(\Lambda_N(b,c)=0\) is equivalent to primitivity of the corresponding quadratic.
If \(q\) is odd, \(\Lambda_N(B,c)\) is monic and square-free of degree \(\varphi(N)\).
Proposition~\ref{v45:prop-fixed-constant-primitive-count} gives exactly that many distinct roots in \(k\),
so all roots lie in \(k\), proving \eqref{v45:eq-optimal-quadratic-odd}.

If \(q\) is even, write \(\Lambda_N(B,c)=\sum_i a_iB^{2i}\). The squaring map is bijective on a finite field, so its unique square root is
\[
 P_{q,c}(B)=\sum_i a_i^{q/2}B^i.
\]
The preceding lemma shows that every root of \(\Lambda_N(B,c)\) is exactly double; hence this square root is monic and square-free.
Its degree is \(\varphi(N)/2\), and the fixed-constant-term count gives exactly that many primitive quadratics.
Thus all its roots lie in \(k\), giving \eqref{v45:eq-optimal-quadratic-even} and the criterion.
Finally, every nonzero polynomial vanishing on all of \(\mathcal P_c\) is divisible by the product of these distinct linear factors,
so the constructed monic polynomial is the unique one of least degree.
\end{proof}

In even characteristic, each inverse pair explains the halving of the degree. After setting \(C=c\), every factor in \eqref{v45:eq-lucas-atom-pairs} has the form
\[
B^2-c(\zeta+\zeta^{-1})
=\bigl(B-\sqrt{c(\zeta+\zeta^{-1})}\bigr)^2.
\]
Square roots are unique even in the algebraic closure, and distinct inverse pairs give distinct square roots. Thus the \(\varphi(N)\) roots counted with multiplicity represent only \(\varphi(N)/2\) distinct coefficients.
Taking the polynomial square root removes one repetition of each root and retains every admissible coefficient.

\leavevmode\phantomsection\label{v45:prop-global-optimal-uniqueness}\ignorespaces
All primitive constant terms can be treated together. Put
\[
 A=k[C]/(\Phi_{q-1}(C)).
\]
The roots of \(\Phi_{q-1}\) are precisely the primitive elements of \(k\), all simple, so CRT gives
\(A\simeq\prod_{c\text{ primitive}}k\).
Hence the polynomials \(P_{q,c}\) uniquely determine a polynomial in \(A[B]\).
For odd \(q\), it is the image of \(\Lambda_N\). For even \(q\), the inverse of squaring on \(A\) is \(a\mapsto a^{q/2}\),
so if \(\Lambda_N=\sum_i a_iB^{2i}\) in \(A[B]\), the corresponding image is \(\sum_i a_i^{q/2}B^i\).
In \(k[B,C]\), representatives may differ by a multiple of \(\Phi_{q-1}(C)\); in \(A[B]\), the element is unique.
This uniqueness comes from polynomial equality in each constant-term fiber, not merely equality of values on a finite point set.

Direct calculation may also use only the recurrence. Equation~\eqref{v45:eq-lucas-atom-root-test} and the remainder formula
give the following criterion without precomputing Lucas atoms.
\begin{corollary}[A recursive coefficient criterion]\leavevmode\phantomsection\label{v45:cor-recursive-quadratic-criterion}
Let \(c\) generate \(k^\times\), and let \(f=X^2+bX+c\). Then \(f\) is primitive if and only if
\begin{equation}\label{v45:eq-recursive-quadratic-criterion}
 U_N(b,c)=0,
 \qquad U_{N/\ell}(b,c)\ne0\quad\text{for every prime }\ell\mid N.
\end{equation}
\end{corollary}
\begin{proof}
If \(f\) has a repeated root \(\theta\), the recurrence gives \(U_n(b,c)=n\theta^{n-1}\),
so \(U_N(b,c)=\theta^q\ne0\). Thus the first condition excludes repeated roots.
For two distinct roots, the same Binet formula gives \(U_n(b,c)=0\iff(\theta_1/\theta_2)^n=1\).
Hence \eqref{v45:eq-recursive-quadratic-criterion} is equivalent to the root ratio having order exactly \(N\);
apply Proposition~\ref{v45:prop-quadratic-root-ratio}.
\end{proof}

\begin{example}[Two small fields of odd order]\label{v45:ex-small-odd-quadratics}
We have \(\Lambda_4(B,C)=B^2-2C\) and \(\Lambda_6(B,C)=B^2-3C\).
In \(\F_3\), the only primitive constant term is two, so
\[
 P_{3,2}(B)=B^2+2=(B-1)(B-2).
\]
The corresponding primitive quadratics are \(X^2+X+2\) and \(X^2+2X+2\).
In \(\F_5\), take the primitive constant term two. Then
\[
 P_{5,2}(B)=B^2-1=(B-1)(B+1).
\]
The corresponding primitive quadratics are \(X^2+X+2\) and \(X^2-X+2\).
A primitive constant term does not ensure that a quadratic is primitive; for example, \(X^2+2\) is irreducible over \(\F_5\), but has projective order two and ordinary order eight.
\end{example}

\subsection{Trace Polynomials of the Norm-One Group and Lucas Atoms}
\label{v146:sec-lucas-pullback}
A root ratio \(\rho\) and its inverse give the same coefficient condition. Pairing inverses in the norm-one group and substituting the resulting trace condition into a quadratic explains the characteristic-dependent multiplicities in the optimal criterion.

\begin{proposition}[Exact-order trace polynomials and quadratic pullback]\label{v146:prop-lucas-pullback}
Let \(s>2\), \(s\mid q+1\), and \(d_s=\varphi(s)/2\). Define the \textbf{exact-order trace polynomial}
\begin{equation}\label{v146:eq-order-trace}
 H_{q,s}(Y)=\prod_{\substack{\{\rho,\rho^{-1}\}\\\ord(\rho)=s}}
 (Y+\rho+\rho^{-1}).
\end{equation}
Then \(H_{q,s}\in\F_q[Y]\) is monic and square-free of degree \(d_s\), with all roots in \(\F_q\). It satisfies
\begin{equation}\label{v146:eq-lucas-pullback}
 \Lambda_s(B,C)=(-C)^{d_s}H_{q,s}\!\left(2-\frac{B^2}{C}\right).
\end{equation}
After cancellation of denominators, this is a polynomial identity in \(\F_q[B,C]\).
\end{proposition}
\begin{proof}
Since \(\rho^q=\rho^{-1}\), every \(-\rho-\rho^{-1}\) belongs to \(\F_q\). Distinct inverse pairs have distinct sums, so the roots are simple. Substitution into each factor gives
\[
 -C\left(2-\frac{B^2}{C}+\rho+\rho^{-1}\right)
 =B^2-C(2+\rho+\rho^{-1}).
\]
Multiply these identities and apply~\eqref{v45:eq-lucas-atom-pairs}.
\end{proof}

In particular, put \(H_q=H_{q,q+1}\). For a primitive constant term \(c\), equation~\eqref{v146:eq-lucas-pullback} pulls the roots of \(H_q\) back to the optimal coefficient criterion. In odd characteristic each root gives two opposite nonzero coefficients; Theorem~\ref{v45:thm-optimal-quadratic-criterion} ensures that both lie in \(\F_q\). In characteristic two, squaring is bijective, so each root gives one coefficient with polynomial multiplicity two. Thus \(\Lambda_{q+1}(B,c)\) equals \(P_{q,c}\) in odd characteristic and \(P_{q,c}^2\) in characteristic two.

\begin{example}[The trace polynomial and atom over four elements]\label{v146:ex-lucas-four}
Let \(\F_4=\F_2(\omega)\), where \(\omega^2+\omega+1=0\). The two inverse pairs of primitive fifth roots of unity give
\[
 H_4(Y)=Y^2+Y+1,\qquad
 \Lambda_5(B,C)=B^4+B^2C+C^2.
\]
For the primitive constant term \(c=\omega\), this gives
\[
 P_{4,\omega}(B)=B^2+\omega^2B+\omega
 =(B-1)(B-\omega).
\]
Hence \(X^2+X+\omega\) and \(X^2+\omega X+\omega\) are precisely the primitive quadratics with constant term \(\omega\). Their root ratios have order five; their roots have order fifteen.
\end{example}

In Chapter 18, \(H_{q,s}\) will prescribe the order of the eigenvalues of a return operator. Here it detects the root-ratio order of a quadratic; equation~\eqref{v146:eq-lucas-pullback} connects the two conditions through the same trace polynomial.

\section{Quadratic Root Sets in Additive Cosets and Frobenius Eigenspaces}
\label{v45:sec-quadratic-cosets}

Let \(q\) be an odd prime power, put \(K=\F_{q^2}\), \(E=\F_{q^4}\), and fix
\(\mu\in\F_q^\times\) and \(\alpha\in K\setminus\F_q\). Consider the family of quadratics
\[
 F_a(X)=X^2+\mu X-a,\qquad a\in\alpha+\F_q.
\]
The coefficients of each member generate \(K\) over \(\F_q\). Reducibility and primitivity in this section are relative to \(K\):
a primitive member has roots of order \(q^4-1\) in \(E^\times\).
The roots of this family can be organized using the two eigenspaces of relative Frobenius, each of dimension two over \(\F_q\). The norm counts the members of each type; see \cite{ZWAdditiveQuadratics} for the related structure.

Put
\[
 L=\{d\in K:d^q=-d\},\qquad
 D=\alpha-\alpha^q,\qquad \Delta=\frac D{\mu^2}.
\]
The space \(L\) is the one-dimensional trace-zero space of \(K/\F_q\), and \(D,\Delta\in L\setminus\{0\}\).
Since the Frobenius difference vanishes on base-field constants, \(D\) depends only on the coset \(\alpha+\F_q\).
Put
\begin{equation}\label{v45:eq-quadratic-normalization}
 w=\frac x\mu+\frac12,\qquad
 \Phi(w)=w^2-w^{2q}.
\end{equation}

\begin{proposition}[Two eigenspaces for the normalized roots]
\label{v45:prop-quadratic-planes}
Under \eqref{v45:eq-quadratic-normalization}, all roots of the family \(\{F_a\}\) correspond bijectively to all solutions of
\[
 \Phi(w)=\Delta.
\]
Every solution belongs to \(E\) and lies in one of the spaces
\[
 V_+=\{w\in E:w^{q^2}=w\}=K,\qquad
 V_-=\{w\in E:w^{q^2}=-w\}.
\]
Both are two-dimensional \(\F_q\)-spaces, and \(E=V_+\oplus V_-\).
Solutions in \(V_+\) correspond to roots of reducible members, and those in \(V_-\) to roots of irreducible members.
\end{proposition}
\begin{proof}
Completing the square gives
\[
 F_a(x)=0\quad\Longleftrightarrow\quad
 a=\mu^2\left(w^2-\frac14\right).
\]
If \(a\in\alpha+\F_q\), then \(a-a^q=D\), so \(\Phi(w)=\Delta\).
Conversely, suppose \(w\) satisfies this equation in a common splitting extension. Taking \(q\)-th powers gives
\(w^{2q}-w^{2q^2}=-\Delta\); adding the original equation gives \(w^{2q^2}=w^2\).
Since the characteristic is odd, \(w^{q^2}=w\) or \(w^{q^2}=-w\). Hence \(w^{q^4}=w\), so \(w\in E\).
In both cases, \(w^2\in K\). Thus \(a=\mu^2(w^2-1/4)\in K\) and \(a-a^q=D\),
giving \(a-\alpha\in\F_q\). This also shows that each root determines its member uniquely.

Relative Frobenius \(\tau:E\to E\), \(\tau(w)=w^{q^2}\), satisfies \(\tau^2=\id\); the projections
\[
 P_+=\frac{\id+\tau}{2},\qquad P_-=\frac{\id-\tau}{2}
\]
give the direct-sum decomposition above. This is Proposition~\ref{prop:trace-splitting} applied to \(E/K\).
Since \(V_+=K\) has dimension two over \(\F_q\), so does \(V_-\).

If \(w\in V_+\), then \(x=\mu(w-1/2)\in K\), and the corresponding quadratic is reducible.
If \(w\in V_-\), then
\[
 x^{q^2}=\mu\left(-w-\frac12\right)=-\mu-x.
\]
Since \(\Delta\ne0\) forces \(w\ne0\), we have \(x^{q^2}\ne x\), so \(x\notin K\) and the corresponding quadratic is irreducible.
\end{proof}

On either plane, \(w^2\in K\), so \(\Phi(w)\in L\).
Thus \(Q_\pm(w)=\Phi(w)/\Delta\) are \(\F_q\)-valued functions.
\leavevmode\phantomsection\label{v45:def-binary-quadratic-forms}\ignorespaces
After choosing an \(\F_q\)-basis of each plane, both are homogeneous quadratic polynomials in the coordinates, called binary quadratic forms.
We calculate these two forms and their nonzero fibers directly.

\begin{proposition}[Two types of quadratic forms and counts of members]
\label{v45:prop-quadratic-coset-count}
In suitable coordinates, \(Q_+\) is a nonzero constant times \(uv\), and its zero set is the union of two lines.
The form \(Q_-\) vanishes only at the zero vector, and in suitable coordinates has the form
\[
 A(u^2-\delta v^2),\qquad
 A\in\F_q^\times,\quad\delta\in\F_q^\times\text{ a nonsquare}.
\]
\leavevmode\phantomsection\label{v45:def-binary-quadratic-types}\ignorespaces
The former is called split and the latter anisotropic. The equation \(\Phi(w)=\Delta\) has \(q-1\) solutions in \(V_+\) and \(q+1\) in \(V_-\).
Consequently, the family \(\{F_a\}\) has exactly \((q-1)/2\) reducible and \((q+1)/2\) irreducible members.
\end{proposition}
\begin{proof}
Choose \(\kappa\in L\setminus\{0\}\). Then \(\kappa^q=-\kappa\) and \(K=\F_q\oplus\F_q\kappa\).
Write \(w=u+v\kappa\), \(\Delta=d\kappa\), with \(u,v\in\F_q\) and \(d\in\F_q^\times\). Then
\[
 \Phi(w)=(u+v\kappa)^2-(u-v\kappa)^2=4uv\kappa.
\]
Thus \(Q_+(w)=4uv/d\), and \(Q_+(w)=1\) is equivalent to \(uv=d/4\), with exactly \(q-1\) solutions.

For \(w\in V_-\), the identity
\[
 \Phi(w)=(w-w^q)(w+w^q)
\]
shows that \(\Phi(w)=0\) implies \(w^q=w\) or \(w^q=-w\), either of which gives \(w^{q^2}=w\).
Combining this with \(w^{q^2}=-w\) gives \(w=0\).
Choose a basis of \(V_-\), and write \(Q_-=As^2+Bst+Ct^2\).
Since no nonzero vector has value zero, \(A\ne0\). Completing the square gives
\[
 Q_-=A\left[\left(s+\frac B{2A}t\right)^2-\delta t^2\right],
 \qquad \delta=\frac{B^2-4AC}{4A^2}.
\]
If \(\delta\) were a square or zero, the form would vanish at a nonzero vector, a contradiction. Thus \(\delta\) is a nonsquare.
Choose \(\vartheta\in K\) with \(\vartheta^2=\delta\). Then \(\vartheta^q=-\vartheta\), and
\[
 u^2-\delta v^2=\Nm_{K/\F_q}(u+v\vartheta).
\]
Since \(A^{-1}\ne0\), Theorem~\ref{thm:trace-norm-surjective} gives exactly \(q+1\) elements in this norm fiber.
Its size also equals that of the norm kernel in Theorem~\ref{thm:hilbert90}.
Hence \(Q_-=1\) has exactly \(q+1\) solutions.

Each \(F_a\) has nonzero discriminant \(\mu^2+4a\), since \(a\notin\F_q\).
Its two distinct roots have centered coordinates \(w,-w\) and lie in the same plane.
Each member contributes two roots; dividing the respective root counts by two gives the numbers of the two types of members.
\end{proof}

\leavevmode\phantomsection\label{v45:def-root-pair-norm-descent}\ignorespaces
The norm takes the same value on a conjugate root pair. Expressing it as a function of a parameter for the pair is called \textbf{norm descent} here.
\begin{proposition}[Conjugate root pairs and norm parameters]
\label{v45:prop-quadratic-norm-descent}
The two roots of an irreducible member are \(x\) and \(x^{q^2}=-\mu-x\), so
\begin{equation}\label{v45:eq-quadratic-norm-descent}
 \Tr_{E/K}(x)=-\mu,\qquad
 \Nm_{E/K}(x)=-a=\mu^2\left(\frac14-w^2\right).
\end{equation}
The map \(w\mapsto z=w^2\) sends each pair \(\{w,-w\}\) corresponding to conjugate roots to a single element.
The irreducible members correspond bijectively to the set
\[
 Z_\Delta=\{z\in K:z-z^q=\Delta,\ z\text{ is a nonsquare in }K\}
\]
via \(a=\mu^2(z-1/4)\). If a member is primitive, then \(-a\) is a primitive element of \(K\).
\end{proposition}
\begin{proof}
Relative Frobenius interchanges the two roots; their sum and product give \eqref{v45:eq-quadratic-norm-descent}.
If \(w\in V_-\setminus\{0\}\), then \(w^2\) cannot be a square in \(K\):
otherwise \(w^2=y^2\), \(y\in K\), would force \(w=\pm y\in K\), contradicting \(w^{q^2}=-w\ne w\).
Conversely, if \(z\in Z_\Delta\), then \(T^2-z\) is irreducible over \(K\); its two roots \(w,-w\) belong to \(E\) and satisfy \(w^{q^2}=-w\).
They satisfy \(\Phi(w)=z-z^q=\Delta\), and therefore correspond uniquely to the stated irreducible member.
Finally, if \(x\) has order \(q^4-1\), then \(x^{q^2+1}\) has order \(q^2-1\), so \(-a\) is primitive in \(K\).
\end{proof}

\par\Needspace{11\baselineskip}
\section*{Exercises}
\phantomsection
\addcontentsline{toc}{section}{Exercises}
\begin{exercise}\label{ex:trace-f9}
In \(\F_9=\F_3(\alpha)\), with \(\alpha^2=-1\), determine whether \(X^3-X-1\) and \(X^3-X-\alpha\) are irreducible, and list all roots when they split.
\end{exercise}
\begin{exercise}\label{ex:trace-iterate}
Let \(E=\F_{q^n}\), \(s\geq1\), and \(g=\gcd(n,s)\). Prove that \(x^{q^s}-x=a\) is solvable in \(E\)
if and only if \(\Tr_{E/\F_{q^g}}(a)=0\), and then has \(q^g\) solutions.
\end{exercise}
\begin{exercise}\label{ex:trace-norm-fibers}
Let \(q\) be odd and \(t,u\in\F_q\). The quadratic character \(\chi\) takes value \(1\) on nonzero squares,
value \(-1\) on nonsquares, and \(\chi(0)=0\).
Prove that the number of elements of \(\F_{q^2}\) satisfying both \(\Tr(x)=t\) and \(\Nm(x)=u\) is \(1-\chi(t^2-4u)\).
\end{exercise}
\begin{exercise}\label{ex:q-artin-schreier}
Let \(q=p^r\) and \(a\in\F_q^\times\). Prove that \(X^q-X-a\) factors over \(\F_q\) into \(q/p\) distinct irreducible factors of degree \(p\).
\end{exercise}

\begin{exercise}[Right inverses for a Frobenius difference equation]\label{ex:supp-hilbert-formula}
Let \(D=\sigma-\id:E\to\ker\Tr_{E/k}\), and fix a \(k\)-linear right inverse \(R\)
given by Proposition~\ref{prop:hilbert90-explicit}, so \(D R=\id\).
Show that all \(k\)-linear right inverses are precisely \(R+\ell\), with \(\ell:\ker\Tr_{E/k}\to k\) an arbitrary linear functional,
and count these right inverses. Use this to explain the normalization choice in the explicit solution formula.
\end{exercise}
\begin{exercise}\label{ex:supp-trace-factor}
Explain how \(X^4+X+1\) factors over \(\F_4\) according to absolute-trace fibers, and prove that the resulting quadratic factors are irreducible.
\end{exercise}

\begin{exercise}[Restricting the trace pairing to a subfield]\label{ex:wu-subfield-pairing}
Let \(k\subseteq L\subseteq E\), with \([E:L]=s\).
Show that, for \(x,y\in L\), we have \(\Tr_{E/k}(xy)=s\Tr_{L/k}(xy)\).
The trace pairing is the bilinear form \((x,y)\mapsto\Tr_{E/k}(xy)\); it is nondegenerate if the zero vector is the only vector pairing to zero with every vector.
Determine when its restriction to \(L\) is nondegenerate and when it vanishes identically. Explain compatibility with nondegeneracy on \(E\), proved in Theorem~\ref{thm:trace-nondegenerate}.
\end{exercise}
\par

\begin{exercise}[Equivariant linear functionals are scalar multiples of trace]\leavevmode\phantomsection\label{ex:wu2-equivariant-functionals}
Let \(\ell:E\to k\) be \(k\)-linear and satisfy \(\ell\sigma=\ell\).
Prove that there is a unique \(c\in k\) with \(\ell=c\Tr_{E/k}\). Deduce again when a Frobenius-equivariant projection \(E\to k\) exists.
\end{exercise}

\begin{exercise}[Separability and diagonalization over specified scalar fields]\label{ex:wu10-diagonalization}
Let \(E=\F_8\) and \(\sigma(x)=x^2\). Choose an \(\F_2\)-basis and represent \(\sigma\) by a three-by-three matrix.
Regard this same matrix over \(\F_2\), \(\F_8\), and \(\F_{64}\), and determine its diagonalizability in each case.
Find the least finite scalar extension over which it is diagonalizable, and explain how this differs from separability of \(E/\F_2\).
\end{exercise}

\begin{exercise}[The first zero term of positive index and the order of a quadratic]\label{v45:ex-quadratic-first-zero-f5}
In \(\F_5[X]\), consider \(f_0=X^2+2\) and \(f_1=X^2+X+2\).
\begin{enumerate}
\item Prove that both are irreducible, and calculate \(u_1,\ldots,u_6\) for each.
\item Find their projective orders, minimal binomial multiples, and ordinary orders, and determine primitivity.
\item Use the recurrence terms for \(f_1\) to write \((X^7-2X)/f_1\), checking the position of the first zero coefficient.
\item Knowing only that the constant term two is primitive, is \(b\ne0\) enough to ensure that \(X^2+bX+2\) is primitive? Justify your answer.
\end{enumerate}
\end{exercise}

\begin{exercise}[Square roots of optimal determining polynomials in characteristic two]\label{v45:ex-optimal-quadratics-f4}
Let \(k=\F_4=\F_2(\omega)\), where \(\omega^2+\omega+1=0\).
\begin{enumerate}
\item Use the homogeneous \(\Phi_5\) or the inverse-pair formula to prove
\(\Lambda_5(B,C)=B^4+B^2C+C^2\).
\item In \(A[B]\), where \(A=k[C]/(C^2+C+1)\), prove that the unique square root of \(\Lambda_5(B,C)\) is
\(P(B,C)=B^2+C^2B+C\).
\item List all primitive quadratics with constant terms \(\omega\) and \(\omega^2\), respectively.
\item Explain why \(\Lambda_5(B,c)\) itself detects primitivity but is not the optimal determining polynomial.
\end{enumerate}
\end{exercise}

\begin{exercise}[An additive-coset family over the nine-element field]
\leavevmode\phantomsection\label{v45:ex-quadratic-cosets-nine}
Let \(K=\F_3(\beta)\), with \(\beta^2=-1\), and consider the three quadratics in \(K[X]\)
\[
 F_t(X)=X^2+X-(\beta+t),\qquad t=0,1,2.
\]
\begin{enumerate}
\item Determine reducibility of each member, factor the reducible member, and compare with the count in Proposition~\ref{v45:prop-quadratic-coset-count}.
\item Calculate the norms of the roots of the two irreducible members in \(\F_{81}/K\), and determine which members are primitive polynomials over \(K\).
\end{enumerate}
\end{exercise}

\begin{exercise}[Norm-preserving linear maps (advanced)]\label{v74:ex-norm-preservers}
Let \(q>2\), \(m\ge1\), \(E=\F_{q^m}\), and let \(T:E\to E\) be \(\F_q\)-linear.
\begin{enumerate}
\item Prove that \(\Nm_{E/\F_q}(T(x))=\Nm_{E/\F_q}(x)\) for every \(x\in E\) if and only if
\[
T(x)=ax^{q^j},\qquad \Nm_{E/\F_q}(a)=1,\qquad 0\le j<m.
\]
Invertibility need not be assumed. Explain why linearity alone does not give this classification.
\item Over \(\F_2\), explain why norm preservation imposes no additional condition on a linear bijection.
Take \(\F_8=\F_2(\theta)\), with \(\theta^3+\theta+1=0\). Show that the linear map determined by
\[
T(1)=1,\qquad T(\theta)=\theta,\qquad T(\theta^2)=\theta^2+1
\]
preserves the norm but is not of the form \(ax^{2^j}\).
\end{enumerate}
\emph{Hint:} For the first part, Lenstra's Theorem~4 \cite{Lenstra1990Automorphisms} may be quoted.
After studying Theorem~\ref{v72:thm-carlitz-mcconnel}, one can instead take \(H=\ker\Nm_{E/\F_q}\), turn norm preservation into a difference-quotient condition, and use \(\F_q\)-linearity to determine the Frobenius step.
\end{exercise}

\begin{exercise}[Two adjacent norm fibers]\label{v82:ex-adjacent-norm-fibers}
Let \(E=\F_{q^2}\) and \(k=\F_q\). Fix \(a,b\in k^\times\), and put
\[
 M(a,b)=\#\{x\in E:\Nm_{E/k}(x)=a,\ \Nm_{E/k}(x+1)=b\},
 \qquad t=b-a-1.
\]
\begin{enumerate}
\item Convert the two norm conditions into prescribed trace and norm.
For odd \(q\), use Exercise~\ref{ex:trace-norm-fibers} to prove \(M(a,b)=1-\chi(t^2-4a)\).
\item Let \(q=2^s\). Prove that \(M(a,b)=1\) when \(t=0\), and determine \(M(a,b)\) from \(\Tr_{k/\F_2}(a/t^2)\) when \(t\ne0\).
\item Explain why merely counting the roots of \(X^2-tX+a\) in \(E\) is insufficient: one must also know whether they are exchanged by the \(q\)-Frobenius.
\end{enumerate}
Nonzero norm fibers are multiplicative cosets of \(\ker\Nm_{E/k}\), so \(M(a,b)\) also counts the intersection of one coset with an additive translate of another. For this cyclotomic-number viewpoint, see~\cite{KudoNogataCyclotomic2026}.
\end{exercise}

\begin{exercise}[The sign obstruction in the constant term of a primitive polynomial]\label{v82:ex-primitive-constant-term}
Let \(Q\) be a prime power, and let \(f\in\F_Q[X]\) be monic and primitive of degree \(m\ge1\).
\begin{enumerate}
\item Prove that \((-1)^m f(0)\) generates \(\F_Q^\times\).
\item If \(Q\equiv3\pmod4\) and \(m\) is odd, prove that no monic degree-\(m\) polynomial with primitive constant term can be primitive.
\item Using Example~\ref{v48:ex-primitive-normal-f16} or Exercise~\ref{v45:ex-quadratic-first-zero-f5}, explain why irreducibility, primitivity, and multiplicative-order conditions on coefficients must be distinguished.
\end{enumerate}
\end{exercise}

\begin{exercise}\label{v146:ex-lucas-characteristic}
Let the field have characteristic \(p\), and write \(n=p^v m\), where \(p\nmid m\). Express \(U_n(B,C)\) as a product of a power of \(U_m(B,C)\) and a power of \(B^2-4C\); in characteristic two use \(B\) instead. Describe the change in zeros and multiplicities when the characteristic divides the index, and decide whether it produces any new nonzero root-ratio orders.
\end{exercise}

\chapter{Trace Duality, Bases, and Coordinate Recovery}
\label{chap:dual}
The trace pairing identifies field elements with linear functionals, and conjugate coordinates express this pairing as a matrix. Dual bases recover coordinates, while relative trace is adjoint to subfield inclusion, linking subfield structure to linear orthogonality.

\section{Nondegeneracy of the Trace Pairing}
\leavevmode\phantomsection\label{def:trace-pairing}\ignorespaces Fix \(k=\F_q\subseteq E=\F_{q^n}\). Define the \textbf{trace pairing} by
\[
B(x,y)=\Tr_{E/k}(xy).
\]
It is a symmetric \(k\)-bilinear form.
\phantomsection\label{def:pairing-radical}
A bilinear form is also called a bilinear pairing. The \textbf{radical} of a symmetric pairing \(B\) is
\[
\rad B=\{x\in E:B(x,y)=0\text{ for every }y\in E\}.
\]
A pairing is nondegenerate if its radical is zero. Vectors in the radical pair to zero with every vector.
A linear functional is a linear map to the base field. The space \(E^\vee=\operatorname{Hom}_k(E,k)\) of all such maps is the dual space.

\begin{theorem}\label{thm:trace-nondegenerate}
The trace pairing is nondegenerate: if \(B(x,y)=0\) for every \(y\in E\), then \(x=0\).
Thus every \(k\)-linear functional \(\lambda:E\to k\) has a unique expression \(\lambda(x)=\Tr_{E/k}(ax)\), with \(a\in E\).
\end{theorem}
\begin{proof}
If \(x\ne0\), choose \(z\) with \(\Tr(z)=1\) by surjectivity of the trace, and set \(y=z/x\). Then \(B(x,y)=1\).
Hence the pairing is nondegenerate.
The map \(E\to\operatorname{Hom}_k(E,k)\), \(a\mapsto(x\mapsto\Tr(ax))\), is therefore injective.
Its domain and codomain both have dimension \(n\), so it is an isomorphism.
\end{proof}

\begin{theorem}[Dual bases]\label{thm:dual-basis}
For every ordered \(k\)-basis \(b_0,\ldots,b_{n-1}\), there is a unique ordered basis \(b_0^*,\ldots,b_{n-1}^*\) such that
\[
\Tr(b_i b_j^*)=\delta_{ij}.
\]
Here \(\delta_{ij}\) is the Kronecker delta: it equals \(1\) when \(i=j\) and \(0\) otherwise.
This basis is the \textbf{dual basis} of the original basis. For every \(x\in E\),
\[
\boxed{\ x=\sum_{i=0}^{n-1}\Tr(xb_i^*)b_i
=\sum_{i=0}^{n-1}\Tr(xb_i)b_i^*.\ }
\]
\end{theorem}
\begin{proof}
Apply Theorem~\ref{thm:trace-nondegenerate} to the \(i\)-th coordinate functional \(\lambda_i\) of the original basis
to obtain a unique \(b_i^*\). These elements are linearly independent: pairing \(\sum c_i b_i^*=0\) with \(b_j\) gives \(c_j=0\).
The coordinate recovery formula follows from the definition. Symmetry implies that the dual of the dual basis is the original basis, giving the second formula.
\end{proof}

\leavevmode\phantomsection\label{def:gram-matrix}\ignorespaces The matrix \(G=(\Tr(b_i b_j))\) of pairings of basis vectors is the \textbf{Gram matrix} of the pairing in this basis. Then
\[
b_j^*=\sum_i (G^{-1})_{ij}b_i.
\]
This gives a direct algorithm for computing the dual basis.
\leavevmode\phantomsection\label{def:trace-orthogonal-complement}\ignorespaces For a subspace \(U\subseteq E\), define \(U^\perp=\{y\in E:\Tr(xy)=0\text{ for all }x\in U\}\), the \textbf{orthogonal complement} of \(U\) with respect to the trace pairing. Nondegeneracy gives
\[
\dim_k U^\perp=n-\dim_k U,\qquad (U^\perp)^\perp=U.
\]

\begin{example}[Three traces recover an extension-field linear combination]\label{v150:ex-three-traces}
Let $\F_4=\F_2(\omega)$, where $\omega^2+\omega+1=0$, and write $\Tr=\Tr_{\F_4/\F_2}$. For $a,b\in\F_4$, retain only
\[
 A=\Tr(\omega^2a),\qquad B=\Tr(b),\qquad
 C=\Tr\bigl(\omega(a+b)\bigr).
\]
Writing $a=a_0+a_1\omega$ and $b=b_0+b_1\omega$ gives
$A=a_0$, $B=b_1$, and $C=a_0+a_1+b_0+b_1$. Hence
\begin{equation}\label{v150:eq-three-traces}
 a+\omega b=(A+B)+(A+C)\omega.
\end{equation}
Three base-field values therefore recover this extension-field linear combination. The observation map has kernel
$\{(0,0),(\omega,1)\}$, contained in the kernel of the target map $(a,b)\mapsto a+\omega b$; the observations do not recover $a,b$ themselves.
The observation map is $\F_2$-linear, whereas the target map is also $\F_4$-linear. In general, for linear maps over a common base field, an observation $L:V\to W$ recovers a target $T:V\to U$ if and only if $\ker L\subseteq\ker T$: precisely then $T$ factors uniquely through $L(V)$.
\end{example}

\section{Interpolation and the Dual of a Power Basis}
\begin{theorem}\label{thm:power-dual}
Let \(E=k(\alpha)\), and let \(m=m_{\alpha,k}\) have degree \(n\). Write in \(E[T]\)
\[
\frac{m(T)}{T-\alpha}=\sum_{i=0}^{n-1}b_i(\alpha)T^i.
\]
Then the dual basis of the power basis \(1,\alpha,\ldots,\alpha^{n-1}\) is
\[
\frac{b_0(\alpha)}{m'(\alpha)},\ldots,
\frac{b_{n-1}(\alpha)}{m'(\alpha)}.
\]
\end{theorem}
\begin{proof}
Write \(m(T)=T^n+c_{n-1}T^{n-1}+\cdots+c_0\).
Synthetic division gives the coefficients of the quotient as
\[
b_{n-1}(U)=1,\qquad b_{i-1}(U)=c_i+U b_i(U)\quad(1\le i<n),
\]
where \(b_i(U)\in k[U]\). Thus, at each conjugate root \(\alpha_j=\alpha^{q^j}\),
\[
\frac{m(T)}{T-\alpha_j}=\sum_{i=0}^{n-1}b_i(\alpha_j)T^i.
\]
These roots are distinct, so \(m'(\alpha_j)\ne0\).

For \(0\le t<n\), interpolation at the nodes \(\alpha_0,\ldots,\alpha_{n-1}\) gives
\[
T^t=\sum_{j=0}^{n-1}
\alpha_j^t\,\frac{m(T)}{(T-\alpha_j)m'(\alpha_j)}.
\]
Each interpolation factor is one at its own node and zero at the others. Both sides have degree less than \(n\), so the identity holds.
Comparing coefficients of \(T^i\) gives
\[
\delta_{it}
=\sum_{j=0}^{n-1}\frac{\alpha_j^t b_i(\alpha_j)}{m'(\alpha_j)}.
\]
Since \(b_i,m'\) have coefficients in \(k\), the \(j\)-th term is the \(q^j\)-th power of \(\alpha^t b_i(\alpha)/m'(\alpha)\). Hence the sum equals
\[
\Tr_{E/k}\left(\alpha^t\frac{b_i(\alpha)}{m'(\alpha)}\right).
\]
This holds for all \(0\le i,t<n\) and is precisely the defining property of a dual basis. Uniqueness follows from Theorem~\ref{thm:dual-basis}.
\end{proof}

\begin{proposition}[Top-degree coefficients and trace functionals]\label{prop:coefficient-trace}
Let \(f\in k[T]\) be monic and irreducible, with \(\deg f=d\), \(E=k(\alpha)\), and \(f(\alpha)=0\).
For \(u=P(\alpha)\) with \(\deg P<d\), define \(\lambda(u)=[T^{d-1}]P(T)\). Then
\begin{equation}\label{eq:coefficient-trace}
\boxed{\quad \lambda(u)=\Tr_{E/k}\left(\frac{u}{f'(\alpha)}\right).\quad}
\end{equation}
\end{proposition}
\begin{proof}
Irreducible polynomials over finite fields are separable, so \(f'(\alpha)\ne0\).
In Theorem~\ref{thm:power-dual}, the leading coefficient of \(f(T)/(T-\alpha)\) is one.
Thus \(1/f'(\alpha)\) is dual to the last power-basis vector \(\alpha^{d-1}\).
Writing \(u=\sum_{i=0}^{d-1}u_i\alpha^i\) gives
\[
\Tr_{E/k}\left(\frac{u}{f'(\alpha)}\right)
=\sum_{i=0}^{d-1}u_i\Tr_{E/k}\left(\frac{\alpha^i}{f'(\alpha)}\right)
=\sum_{i=0}^{d-1}u_i\delta_{i,d-1}=u_{d-1}.
\]
This is exactly the definition of \(\lambda(u)\), including the case \(d=1\).
\end{proof}

The denominator \(f'(\alpha)\) is the normalization factor in Lagrange interpolation. Summing over all conjugate roots combines local evaluations into a coefficient in the base field.
The same functional can therefore extract coordinates in the quotient ring or sum conjugate values on the root set.
Section~\ref{str:sec-trace-duality} uses this weight to turn intersections with subfields in low degree into relative-trace orthogonality conditions;
Section~\ref{sec:recurrence-algebra} then identifies the quotient-algebra model of recurrence sequences directly with their trace representation.

\begin{example}
In the same \(\F_9=\F_3(\alpha)\), \(\alpha^2=-1\), as in Example~\ref{v100:ex-three-generation-nine}, interpolation now recovers dual coordinates. We have
\[
\frac{T^2+1}{T-\alpha}=T+\alpha,\qquad m'(\alpha)=2\alpha.
\]
Thus the dual basis of \((1,\alpha)\) is \((2,\alpha)\).
If \(x=a+b\alpha\), then \(a=\Tr(2x)\) and \(b=\Tr(\alpha x)\), recovering \(x\) from two base-field values.
\end{example}

Interpolation also compares evaluations at two different root sets. The roots in either set need not form a single conjugacy orbit.
\begin{proposition}[Changing the nodes of root evaluation]\label{v97:prop-cauchy-node-change}
Let \(f,g\in k[X]\) be monic, squarefree, and coprime, both of degree \(d\ge1\).
In a common splitting field \(E\), list their roots as \(\alpha_1,\ldots,\alpha_d\) and \(\beta_1,\ldots,\beta_d\), and put
\[
 V_\alpha=(\alpha_i^r)_{0\le r<d,\,1\le i\le d},\qquad
 V_\beta=(\beta_j^r)_{0\le r<d,\,1\le j\le d}.
\]
Then \(M=V_\alpha^{-1}V_\beta\) satisfies
\begin{equation}\label{v97:eq-cauchy-node-change}
 \boxed{M_{ij}=\frac{f(\beta_j)}{f'(\alpha_i)(\beta_j-\alpha_i)}.}
\end{equation}
Every nonempty square minor of \(M\) is nonzero.
For \(\deg P<d\), the evaluation column vectors \(v_\alpha=(P(\alpha_i))_i\) and \(v_\beta=(P(\beta_j))_j\) satisfy
\begin{equation}\label{v97:eq-evaluation-transpose}
 v_\beta=M^{\mathsf T}v_\alpha.
\end{equation}
\end{proposition}
\begin{proof}
The Lagrange basis polynomials \(\ell_i(X)=f(X)/((X-\alpha_i)f'(\alpha_i))\) satisfy
\(X^r=\sum_i\alpha_i^r\ell_i(X)\) for \(0\le r<d\).
Evaluating at \(\beta_j\) yields \(V_\beta=V_\alpha(\ell_i(\beta_j))_{i,j}\), proving \eqref{v97:eq-cauchy-node-change}.
Applying the same interpolation formula to \(P\) gives the transpose in \eqref{v97:eq-evaluation-transpose}.

After extracting the nonzero row factors \(1/f'(\alpha_i)\) and column factors \(f(\beta_j)\), it remains to show that every square submatrix of the \textbf{Cauchy matrix} \((1/(\beta_j-\alpha_i))\) is invertible.
Select \(r\) nodes of each kind. A left-kernel vector \((u_i)\) would make the rational function
\(\sum_i u_i/(T-\alpha_i)\) vanish at the selected \(r\) nodes \(\beta_j\).
Multiplication by \(\prod_i(T-\alpha_i)\) gives a polynomial of degree at most \(r-1\) with \(r\) distinct roots, so it is zero.
Evaluation at each \(\alpha_i\) now gives \(u_i=0\).
\end{proof}

For example, if \(f=X^d-a\) and \(g=X^d-b\), where \(a,b\in k^\times\), \(a\ne b\), and \(p\nmid d\), the formula becomes
\[
 M_{ij}=\frac{b-a}{d\alpha_i^{d-1}(\beta_j-\alpha_i)}.
\]
Here \(p\nmid d\) ensures that the binomials have distinct roots; the general node-change statement only requires distinct roots in each set.
The matrices are interpreted over the common splitting field and describe two linear coordinate systems on the same space of low-degree polynomials.
They do not automatically define a ring isomorphism between \(k[X]/(f)\) and \(k[X]/(g)\).

\section{Coefficient Pairings, Derivatives, and Algebra Trace}
\label{sec:algebra-trace-coefficient}

When the conjugate roots are distinct, interpolation expresses the coefficient functional as a trace. Quotient algebras with repeated factors still retain the coefficient functional,
but its relation to the trace involves an additional derivative factor that may fail to be invertible.

\subsection{The Derivative Formula and Nondegeneracy}
\label{sec:trace-derivative-criterion}
\leavevmode\phantomsection\label{def:algebra-trace}\ignorespaces For a finite-dimensional commutative \(k\)-algebra \(A\), define its \textbf{algebra trace} by
\[
\Tr_{A/k}(u)=\operatorname{tr}(M_u),\qquad M_u(v)=uv.
\]
The right-hand side is the matrix trace of a linear operator and is independent of the basis. If \(k\) is finite and \(A/k\) is a finite field extension,
Proposition~\ref{prop:characteristic-multiplication} identifies it with the field trace defined by the conjugate orbit.

\begin{theorem}[Coefficient pairings, derivatives, and algebra trace]\leavevmode\phantomsection\label{thm:algebra-trace-coefficient}
Let \(k\) be any field, let \(f\in k[T]\) be monic with \(\deg f=d\ge1\), and put \(A=k[T]/(f)\) and \(\alpha=\overline T\).
For \(u\in A\), let \(\lambda(u)\) denote the coefficient of \(T^{d-1}\) in its standard remainder. Then
\begin{enumerate}
\item the pairing \((u,v)\mapsto\lambda(uv)\) is nondegenerate;
\item for every \(u\in A\),
\begin{equation}\label{eq:algebra-trace-coefficient}
\boxed{\quad\Tr_{A/k}(u)=\lambda\bigl(f'(\alpha)u\bigr).\quad}
\end{equation}
\end{enumerate}
\end{theorem}
\begin{proof}
If the standard remainder of \(u\ne0\) has degree \(r\) and leading coefficient \(b\ne0\), take \(v=\alpha^{d-1-r}\). The product has degree exactly \(d-1\) and requires no reduction, so \(\lambda(uv)=b\ne0\).
This proves the first assertion.

To compute the algebra trace, write
\[
f(T)=T^d+c_{d-1}T^{d-1}+\cdots+c_0,\qquad
f(Z)=(Z-\alpha)\sum_{i=0}^{d-1}b_iZ^i\quad\text{in }A[Z].
\]
Since \(f(\alpha)=0\), this division has zero remainder. Comparing coefficients gives
\[
b_{d-1}=1,\qquad
\alpha b_i=b_{i-1}-c_i\ (1\le i<d),\qquad
\alpha b_0=-c_0.
\]
Starting from \(b_{d-1}=1\), the recurrence shows that the standard remainder of \(b_i\) is monic of degree \(d-1-i\). Hence \(\lambda(b_i)=\delta_{i0}\). These recurrences imply
\[
\lambda(\alpha^j b_i)=\delta_{ij}\qquad(0\le i,j<d).
\]
Indeed, the case \(j=0\) is established. If \(0\le j\le d-2\), then \(\lambda(\alpha^j)=0\);
for \(i\ge1\),
\[
\lambda(\alpha^{j+1}b_i)
=\lambda(\alpha^j b_{i-1})-c_i\lambda(\alpha^j)
=\delta_{i-1,j}=\delta_{i,j+1}.
\]
For \(i=0\), both sides vanish. This completes the induction and also covers \(d=1\), where no induction is needed.
Thus \((b_i)\) is dual to the power basis under the coefficient pairing.

The \(i\)-th diagonal entry of the multiplication matrix \(M_u\) is \(\lambda(u\alpha^i b_i)\), so
\[
\Tr_{A/k}(u)=\lambda\left(u\sum_{i=0}^{d-1}\alpha^i b_i\right).
\]
Formally differentiating \(f(Z)=(Z-\alpha)\sum_i b_iZ^i\) and substituting \(Z=\alpha\) gives \(\sum_i\alpha^i b_i=f'(\alpha)\), proving \eqref{eq:algebra-trace-coefficient}.
The trace formula therefore also holds when \(f\) has repeated factors.
\end{proof}

\begin{proposition}[The algebra trace pairing and separability]\label{prop:trace-pairing-separability}
Let \(k\) be any field, let \(f\in k[T]\) be monic with \(\deg f\ge1\), and put \(A=k[T]/(f)\) and \(\alpha=\overline T\). For the algebra trace pairing \(B_{\Tr}(u,v)=\Tr_{A/k}(uv)\),
\begin{equation}\label{eq:trace-pairing-separability}
\boxed{\quad
\text{algebra trace pairing nondegenerate}
\iff f'(\alpha)\in A^\times
\iff \gcd(f,f')=1.
\quad}
\end{equation}
These conditions are also equivalent to separability of \(f\), that is, to having no repeated roots in a splitting field.
If moreover \(k=\F_q\), they are equivalent to \(A\) being reduced, and to \(u\mapsto u^q\) being an automorphism of \(A\).
When these conditions hold, the coefficient functional in Theorem~\ref{thm:algebra-trace-coefficient} satisfies
\[
\lambda(u)=\Tr_{A/k}\left(\frac{u}{f'(\alpha)}\right).
\]
\end{proposition}
\begin{proof}
By the trace formula in Theorem~\ref{thm:algebra-trace-coefficient} and nondegeneracy of the coefficient pairing,
an element \(u\) pairs to zero with every \(v\) under the trace pairing if and only if \(f'(\alpha)u=0\).
Thus nondegeneracy of the trace pairing is equivalent to injectivity of multiplication by \(f'(\alpha)\). Finite dimensionality makes this equivalent to surjectivity,
and hence to invertibility of \(f'(\alpha)\). B\'ezout's identity makes this equivalent to \(\gcd(f,f')=1\).
The derivative criterion identifies \(\gcd(f,f')=1\) with separability of \(f\); the base field need not be finite here.
Now assume \(k=\F_q\). Finite fields are perfect, so \(f\) is separable if and only if \(f\) has no repeated irreducible factor.
If \(f\) has no repeated factor, the Chinese remainder theorem expresses \(A\) as a product of finite fields, so \(A\) is reduced.
If \(f\) has repeated factors, let \(r\) be the product of its distinct monic irreducible factors.
Then \(r(\alpha)\ne0\), but a sufficiently high power of this element is zero, so \(A\) is not reduced.
The equivalence between reducedness and invertibility of Frobenius follows from Theorem~\ref{thm:finite-reduced-algebra}.
Finally, replacing \(u\) by \(u/f'(\alpha)\) in \eqref{eq:algebra-trace-coefficient} gives the last formula.
\end{proof}

Write \(B_\lambda(u,v)=\lambda(uv)\) and \(B_{\Tr}(u,v)=\Tr_{A/k}(uv)\). The trace formula gives
\[
B_{\Tr}(u,v)=B_\lambda\bigl(f'(\alpha)u,v\bigr).
\]
Thus passing from the coefficient pairing to the trace pairing amounts to applying the multiplication operator \(M_{f'(\alpha)}\) to the first argument.
Nondegeneracy of the coefficient pairing uses only standard remainders and their leading terms; it does not require distinct roots.
For the trace pairing to remain nondegenerate, this additional multiplication operator must be invertible.
With repeated factors, the operator has a nonzero kernel, so the trace pairing loses information in some directions.
The coefficient pairing therefore applies to every monic nonconstant polynomial, whereas nondegeneracy of the algebra trace pairing is characterized exactly by separability.

Over a finite base field, Proposition~\ref{prop:coefficient-trace} is the case in which \(A\) is a field. In the square-free but reducible case,
the algebra trace is the sum of the traces of the field components, and the derivative-normalized coefficient formula still holds.
Section~\ref{sec:recurrence-algebra} uses the coefficient pairing to parameterize the solution space of a recurrence, retaining the nilpotent structure associated with repeated factors.

\begin{corollary}[A trace criterion for reduced algebras]\label{cor:reduced-trace-pairing}
Let \(A\) be a finite-dimensional commutative \(k=\F_q\)-algebra. Then \(A\) is reduced if and only if its algebra trace pairing is nondegenerate.
\end{corollary}
\begin{proof}
If \(A\) is reduced, Theorem~\ref{thm:finite-reduced-algebra} expresses it as a product \(\prod_i E_i\) of finite fields.
The multiplication matrices are block diagonal, so the algebra trace is the sum of the component field traces. Each field trace pairing is nondegenerate, and so is their direct sum.
Conversely, if \(A\) is not reduced, choose a nonzero nilpotent element \(u\). For every \(v\in A\), commutativity implies that \(uv\) is nilpotent,
so the multiplication operator \(M_{uv}\) is nilpotent and has matrix trace zero. Thus \(u\) pairs to zero with every \(v\), and the pairing is degenerate.
\end{proof}

\subsection{The Periodic Part of Frobenius and the Reduced Quotient}
\label{sec:algebra-periodic-part}

In a finite-dimensional commutative algebra, the ideal of nilpotent elements is the \textbf{nilradical}.
The following proposition also proves that these elements form an ideal and realizes the reduced quotient inside the original algebra by means of periodic points of Frobenius.

We use the definitions of local ring and residue field on page~\pageref{v47:def-local-ring}.
\phantomsection\label{def:local-residue}
An ideal \(\mathfrak m\) is called a \textbf{nilpotent ideal} if \(\mathfrak m^e=0\) for some positive integer \(e\), that is, if every product of \(e\) of its elements is zero.
\relax

\begin{proposition}[The periodic part of Frobenius and the reduced quotient]\label{prop:algebra-periodic-reduction}
Let \(A\ne0\) be a finite-dimensional commutative \(k=\F_q\)-algebra, and put
\[
F:A\longrightarrow A,\qquad F(a)=a^q,
\qquad \mathfrak n=\{a\in A:a^m=0\text{ for some }m\ge1\}.
\]
Write
\[
A_{\mathrm{per}}=\{a\in A:F^m(a)=a\text{ for some }m\ge1\}.
\]
Then \(\mathfrak n\) is an ideal, \(A_{\mathrm{per}}\) is a reduced \(k\)-subalgebra, and
\begin{equation}\label{eq:algebra-periodic-reduction}
\boxed{\quad A=A_{\mathrm{per}}\oplus\mathfrak n,
\qquad A_{\mathrm{per}}\xrightarrow{\ \sim\ }A/\mathfrak n.\quad}
\end{equation}
The first identity is a direct sum of \(k\)-vector spaces; the second is induced by restricting the natural quotient map.
For sufficiently large \(r\),
\[
\ker F^r=\mathfrak n,\qquad \im F^r=A_{\mathrm{per}}.
\]
There is a positive integer \(N\) such that \(F^N\) is an algebra projection onto \(A_{\mathrm{per}}\) with kernel \(\mathfrak n\).
Moreover, \(A_{\mathrm{per}}\) is the unique \(k\)-subalgebra mapped isomorphically onto \(A/\mathfrak n\) by the natural quotient map.
\end{proposition}
\begin{proof}
Since \(A\) is commutative and \(k=\F_q\), \(F\) is both an algebra homomorphism and a \(k\)-linear map.
Finite dimensionality ensures that the ascending chain of kernels and the descending chain of images both stabilize. Choose \(r\ge1\) such that both are stable from step \(r\) onward,
and put \(B=\im F^r\). Each \(\ker F^j\) is an ideal. An element is nilpotent if and only if some iterate of the \(q\)-th power map sends it to zero, so
\[
\mathfrak n=\bigcup_{j\ge1}\ker F^j=\ker F^r.
\]
This also proves that \(\mathfrak n\) is an ideal.

Stability of the image gives \(F(B)=B\). Thus \(F|_B\) is surjective on a finite-dimensional space, hence bijective.
The set \(B\) is finite, so all its elements are periodic under \(F\).
Conversely, if \(F^m(a)=a\), choose \(jm\ge r\); then \(a=F^{jm}(a)\in B\).
Hence \(B=A_{\mathrm{per}}\). In particular, the periodic points form a unital \(k\)-subalgebra.

If \(a\in B\cap\mathfrak n\), then \(F^r(a)=0\). Since \(F^r|_B\) is injective, \(a=0\).
The rank--nullity theorem gives \(\dim B+\dim\mathfrak n=\dim A\), whence \(A=B\oplus\mathfrak n\).
The quotient map restricted to \(B\) is therefore a bijective algebra homomorphism. Since \(B\) has no nonzero nilpotents, it is reduced.
Reducedness of the quotient can also be checked directly: if \(a^m\in\mathfrak n\), then \(a\) is nilpotent, so \(a\in\mathfrak n\).

Let \(M\) be the order of the permutation \(F|_B\), and choose \(N\ge r\) with \(M\mid N\).
For the unique decomposition \(a=b+u\) with \(b\in B,u\in\mathfrak n\), we have \(F^N(a)=b\).
Thus \(F^N\) is the required projection and satisfies \((F^N)^2=F^N\).

Finally, suppose a subalgebra \(C\subseteq A\) maps isomorphically onto \(A/\mathfrak n\) under the quotient map. Then \(C\) is reduced.
Its Frobenius is injective, hence bijective by finiteness. Every element of \(C\) is therefore periodic, giving \(C\subseteq B\).
Both have dimension \(\dim(A/\mathfrak n)\), so \(C=B\).
\end{proof}

\begin{corollary}[Fixed points of Frobenius iterates in algebras with repeated factors]\label{cor:algebra-fixed-points-reduced}
If \(A/\mathfrak n\simeq\prod_{i=1}^s\F_{q^{d_i}}\), then for every \(m\ge1\),
\[
\Fix(F^m)\simeq\prod_{i=1}^s\F_{q^{\gcd(d_i,m)}},
\qquad |\Fix(F^m)|=q^{\sum_i\gcd(d_i,m)}.
\]
\end{corollary}
\begin{proof}
All fixed points lie in \(A_{\mathrm{per}}\). The natural isomorphism in the proposition commutes with Frobenius,
so the fixed-point algebra is the product of the fixed fields in the field components. Apply Proposition~\ref{prop:fixed-gcd} to each component.
\end{proof}

Thus the reduced quotient can be obtained either by factoring out nilpotents or by selecting the periodic points of Frobenius.
Multiplication in~\eqref{eq:algebra-periodic-reduction} still includes the action of the periodic subalgebra on the nilpotent ideal,
so this vector-space direct sum is generally not a direct product of algebras.
The Berlekamp algebra in Chapter~\ref{chap:algorithms} is the case \(\Fix(F)\) here.

The periodic subalgebra reduces multiplication characteristic polynomials to calculations in the field components.
\begin{proposition}[Local components, characteristic polynomials, and the trace pairing]\label{v117:prop-local-dimension-trace}
Let \(A\ne0\) be a finite-dimensional commutative algebra over \(k=\F_q\), let \(p=\Char k\), and write
\[
A_{\mathrm{per}}=\prod_{i=1}^r K_i.
\]
Let \(e_i\) be the corresponding coordinate idempotents, put \(A_i=e_iA\) and \(\mathfrak n_i=e_i\mathfrak n\), and regard \(K_i=e_iA_{\mathrm{per}}\) as a subfield of \(A_i\). Set \(\ell_i=\dim_{K_i}A_i\). If \(a=b+u\), where \(b=(b_i)\in A_{\mathrm{per}}\) and \(u\in\mathfrak n\), let \(M_a\) denote multiplication by \(a\) as a \(k\)-linear operator. Then
\begin{equation}\label{v121:eq-local-characteristic}
\det_k(TI-M_a)=\prod_{i=1}^r
m_{b_i,k}(T)^{\ell_i[K_i:k(b_i)]}.
\end{equation}
Here \(m_{b_i,k}\) is the monic minimal polynomial of \(b_i\) over \(k\). In particular,
\begin{equation}\label{v121:eq-local-determinant}
\det_k M_a=\prod_{i=1}^r\Nm_{K_i/k}(b_i)^{\ell_i},
\end{equation}
and
\begin{equation}\label{v117:eq-local-dimension-trace}
\Tr_{A/k}(a)=\sum_{i=1}^r\ell_i\Tr_{K_i/k}(b_i),\qquad
\rank B_{\Tr}=\sum_{p\nmid\ell_i}[K_i:k].
\end{equation}
Under \(A=\prod_iA_i\), the radical of the trace pairing has component \(A_i\) if \(p\mid\ell_i\), and \(\mathfrak n_i\) if \(p\nmid\ell_i\).
\end{proposition}
\begin{proof}
The orthogonal idempotents give \(A=\prod_iA_i\), with \(A_i/\mathfrak n_i\simeq K_i\). If \(b_i\ne0\), then \(b_i+u_i=b_i(1+b_i^{-1}u_i)\) is invertible by a finite geometric series; if \(b_i=0\), it is nilpotent. Thus \(A_i\) is local, with maximal ideal \(\mathfrak n_i\).

Multiplication by \(u_i\) is nilpotent over \(K_i\). In a \(K_i\)-basis adapted to its kernel filtration, its matrix is strictly upper triangular. Multiplication by \(a_i=b_i+u_i\) therefore has an upper triangular matrix with every diagonal entry equal to \(b_i\). Expanding each coordinate in a \(k\)-basis of \(K_i\) gives a block upper triangular \(k\)-matrix with \(\ell_i\) diagonal blocks, each representing multiplication by \(b_i\) on \(K_i\).
By Proposition~\ref{prop:characteristic-multiplication}, each block has characteristic polynomial \(m_{b_i,k}^{[K_i:k(b_i)]}\), trace \(\Tr_{K_i/k}(b_i)\), and determinant \(\Nm_{K_i/k}(b_i)\). Multiplying or summing over the components proves the three formulas.

Since the quotient map is an algebra homomorphism, the trace pairing on the \(i\)-th component is
\[
(a_i,c_i)\longmapsto\ell_i\Tr_{K_i/k}(\bar a_i\bar c_i).
\]
Nondegeneracy of the field trace pairing gives the asserted radical and rank.
\end{proof}

In general, the weight is the dimension over the residue field, not the least exponent that annihilates the maximal ideal.
\begin{example}[Dimension versus nilpotency index]\label{v117:ex-dimension-index}
Let \(A=k[X,Y]/(X^2,XY,Y^2)\). Its maximal ideal \(\mathfrak n=(\bar X,\bar Y)\) satisfies \(\mathfrak n^2=0\), whereas \(A_{\mathrm{per}}=k\) and \(\dim_kA=3\). Hence
\[
\Tr_{A/k}(a+b\bar X+c\bar Y)=3a,\qquad
\det_k(TI-M_{a+b\bar X+c\bar Y})=(T-a)^3.
\]
The trace pairing has rank one in characteristic two and rank zero in characteristic three; the nilpotency index is two in both cases.
If \((b,c)\ne(0,0)\), the minimal polynomial of this multiplication operator is \((T-a)^2\); if \(b=c=0\), it is \(T-a\). The characteristic polynomial is unchanged, whereas the minimal polynomial distinguishes whether the nilpotent part vanishes.
\end{example}

\subsection{Local Multiplicities and the Rank of the Trace Pairing}
\label{sec:local-trace-rank}

\begin{example}[Two kinds of trace degeneracy in the dual numbers]\label{exm:dual-numbers-trace}
For any finite field \(k\), let \(A=k[\varepsilon]/(\varepsilon^2)\), the \textbf{ring of dual numbers}.
In the basis \((1,\varepsilon)\), multiplication by \(u=a+b\varepsilon\) has matrix
\[
M_u=\begin{pmatrix}a&0\\b&a\end{pmatrix},\qquad
\Tr_{A/k}(u)=2a.
\]
If \(v=c+d\varepsilon\), then \(uv=ac+(ad+bc)\varepsilon\). The trace and coefficient pairings are therefore
\[
B_{\Tr}(u,v)=2ac,\qquad B_\lambda(u,v)=ad+bc,
\quad \lambda(a+b\varepsilon)=b,
\]
with Gram matrices
\[
G_{\Tr}=\begin{pmatrix}2&0\\0&0\end{pmatrix},\qquad
G_\lambda=\begin{pmatrix}0&1\\1&0\end{pmatrix}.
\]
In odd characteristic, \(G_{\Tr}\) has rank one and radical \(k\varepsilon\); in characteristic two it is zero and its radical is all of \(A\).
In both cases \(G_\lambda\) is invertible and \(\lambda(\varepsilon\cdot1)=1\), so the coefficient pairing detects the nilpotent direction.
For \(f(T)=T^2\), the derivative \(f'=2T\) in the trace formula explains both degeneracies.
The following corollary expresses the dimensions in the general trace formula as multiplicities of irreducible factors.
Example~\ref{exm:mixed-algebra-trace} then combines different local components in an eight-dimensional model.
Example~\ref{exm:repeated-recurrence-algebra} computes the period-two sequences produced by multiplication by \(1+\varepsilon\) in the same dual-number algebra in characteristic two.
\end{example}

For a monic irreducible polynomial \(P\in k[T]\) and a positive integer \(e\), put
\[
R_P=k[T]/(P^e),\qquad
\mathfrak m=(\overline P),\qquad K_P=k[T]/(P).
\]
B\'ezout's identity shows that an element of \(R_P\) is invertible if and only if its image modulo \(P\) is nonzero.
Thus the nonunits form exactly \(\mathfrak m\). The quotient \(R_P/\mathfrak m\simeq K_P\) is a field, and no proper ideal contains a unit,
so \(\mathfrak m\) is the unique maximal ideal. Moreover, \(\mathfrak m^e=0\).
This model retains both the field associated with an irreducible factor and the nilpotent structure recorded by its multiplicity.

\[
R_P\supset\mathfrak m\supset\mathfrak m^2\supset\cdots\supset\mathfrak m^e=0.
\]
\phantomsection\label{def:ideal-filtration}
Such a descending chain of ideals is called a \textbf{filtration}. For \(0\le j<e\), the map \(V\bmod P\mapsto P^jV\bmod P^{j+1}\) identifies \(K_P\) with \(\mathfrak m^j/\mathfrak m^{j+1}\). Thus each layer retains the same residue field.

\begin{corollary}[Local multiplicities and algebra trace]\label{prop:local-multiplicity-trace}
Let \(k=\F_q\) have characteristic \(p\), and let \(f=\prod_{i=1}^rP_i^{e_i}\) be the irreducible factorization of a monic nonconstant polynomial,
with the \(P_i\) pairwise distinct. Put
\[
A=k[T]/(f),\qquad A_i=k[T]/(P_i^{e_i}),\qquad K_i=k[T]/(P_i).
\]
Let \(\bar u_i\in K_i\) be the image of \(u\in A\) in the corresponding field \(K_i\). Then
\begin{equation}\label{eq:local-multiplicity-trace}
\boxed{\quad\Tr_{A/k}(u)=\sum_{i=1}^r e_i\Tr_{K_i/k}(\bar u_i).\quad}
\end{equation}
The coefficient \(e_i\) is interpreted as its image in \(k\). The rank of the algebra trace pairing satisfies
\begin{equation}\label{eq:trace-rank-multiplicities}
\boxed{\quad\rank B_{\Tr}=\sum_{p\nmid e_i}\deg P_i
=\deg f-\deg\gcd(f,f').\quad}
\end{equation}
Under the Chinese remainder decomposition \(A\simeq\prod_i A_i\), the radical in component \(i\) is
\(A_i\) if \(p\mid e_i\), and the ideal \((\overline{P_i})\) if \(p\nmid e_i\).
\end{corollary}
\begin{proof}
On the component \(A_i\), we have \(\dim_kA_i=e_i\deg P_i\) and \([K_i:k]=\deg P_i\), hence \(\ell_i=e_i\). Proposition~\ref{v117:prop-local-dimension-trace} gives the trace formula, the radical, and the first rank formula.

Write \(f=P_i^{e_i}R_i\) with \(P_i\nmid R_i\). Since finite fields are perfect, \(P_i\nmid P_i'\), and
\[
f'=e_iP_i^{e_i-1}P_i'R_i+P_i^{e_i}R_i'.
\]
If \(p\nmid e_i\), the exponent of \(P_i\) in \(\gcd(f,f')\) is exactly \(e_i-1\);
if \(p\mid e_i\), it is \(e_i\). Comparing degrees gives the second formula, including the case \(f'=0\).
\end{proof}

Formula~\eqref{eq:local-multiplicity-trace} shows that the algebra trace takes the residue-field trace in each component, weights it by the multiplicity, and sums the results.
The trace pairing vanishes whenever one argument is nilpotent. If the multiplicity vanishes in the base field, the entire residue-field contribution vanishes as well.
Section~\ref{sec:recurrence-algebra} retains the coefficient functional \(\lambda\) precisely to distinguish these directions and their shift action.
The unit norm equation in Chapter~\ref{chap:normal} is likewise first solved over the residue field; lifting errors must then be corrected in the nilpotent ideal
(Lemma~\ref{lem:odd-circulant-norm}).
Thus residue-field traces and multiplicities determine the algebra trace, whereas solving unit norm equations also requires the original local ring.

\begin{remark}[The role of perfectness in the reducedness criterion]\label{rem:perfectness-trace}
The coefficient pairing and derivative formula of Theorem~\ref{thm:algebra-trace-coefficient} hold over any base field,
as does the separability criterion in Proposition~\ref{prop:trace-pairing-separability}.
Identifying separability further with reducedness of the quotient algebra uses perfectness of the base field.
For example, take an indeterminate \(u\) and put \(k=\F_p(u^p)\) and \(E=\F_p(u)\). Example~\ref{ex:inseparable-rational-extension} in the appendix proves that
\[
E\simeq k[T]/(T^p-u^p),\qquad [E:k]=p.
\]
This algebra is a field and hence reduced, but its defining polynomial has zero derivative, so~\eqref{eq:algebra-trace-coefficient} gives \(\Tr_{E/k}=0\).
Surjectivity of Frobenius over finite fields ensures that irreducible polynomials are separable, excluding this reduced but inseparable case.
\end{remark}

\begin{example}[The periodic part and trace degeneracy in an eight-dimensional quotient algebra]\label{exm:mixed-algebra-trace}
Over \(k=\F_2\), put
\[
P=T^2+T+1,\qquad f=(T+1)^2P^3,\qquad A=k[T]/(f).
\]
The Chinese remainder theorem gives
\[
A\simeq A_1\times A_2,
\qquad A_1=k[T]/((T+1)^2),\quad A_2=k[T]/(P^3).
\]
The two components have dimensions two and six. Put
\[
\mathfrak n=(\overline{T+1})A_1\times(\overline P)A_2.
\]
This is exactly the nilradical, of dimension \(1+4=5\), satisfying \(\mathfrak n^3=0\), and
\[
A/\mathfrak n\simeq\F_2\times\F_4.
\]
Thus \(A_{\mathrm{per}}\) has dimension three. If \(F(u)=u^2\), then \(F^2(u)=u^4\) is the projection onto this periodic subalgebra.
Indeed, each element of the reduced quotient satisfies \(\bar u^4=\bar u\), so \(u^4-u\in\mathfrak n\).
Since \(\mathfrak n^3=0\), we have \((u^4-u)^4=0\), or \(u^{16}=u^4\).
Thus \(F^2\) fixes its image pointwise. Its kernel is exactly \(\mathfrak n\), so its image has dimension three and is precisely \(A_{\mathrm{per}}\).

For \(u=(u_1,u_2)\), Corollary~\ref{prop:local-multiplicity-trace} gives
\[
\Tr_{A/k}(u)=\Tr_{\F_4/\F_2}(\bar u_2),
\qquad
\rad B_{\Tr}=A_1\times(\overline P)A_2.
\]
The first component has multiplicity two, which is zero in the base field; the second has multiplicity three, which equals one.
Hence the trace pairing has rank two and radical dimension six, strictly exceeding the nilradical dimension five.
The coordinate idempotent \((1,0)\) lies in the radical of the trace pairing but is not nilpotent.

In contrast, let \(\lambda\) extract the coefficient of \(T^7\) in the standard remainder. By Theorem~\ref{thm:algebra-trace-coefficient}, the pairing \(\lambda(uv)\) is nondegenerate and has rank eight.
It still distinguishes the directions erased by the algebra trace. Example~\ref{exm:mixed-algebra-recurrence} studies multiplication by \(\overline T\) in this same algebra.
\end{example}

\subsection{Quadratic Quotient Algebras and Product Coordinates of Root Sets}\label{v85:sec-quadratic-root-coordinates}
Reduction of a root polynomial modulo a low-degree polynomial retains only part of the root-set data. A quadratic quotient turns these data into a norm-one product, a ratio of two values, or a sum of reciprocals. The factorization type of the modulus determines which form occurs.
Let \(k=\F_q\), let \(Q\in k[T]\) be monic of degree two, and put
\[
 R_Q=k[T]/(Q),\qquad G_Q=R_Q^\times/k^\times.
\]
\leavevmode\phantomsection\label{v85:def-projective-unit-group}\ignorespaces
The group \(G_Q\) is the projective unit group of \(R_Q\): units differing by a nonzero base-field scalar represent the same class, denoted by \([u]\). It is the quadratic instance of the unit group modulo scalars used in Definition~\ref{def:projective-order}.

\begin{proposition}[Projective unit groups of quadratic quotient algebras]\label{v85:prop-quadratic-unit-groups}
The following are group isomorphisms, where \(\overline f\) denotes the residue class of \(f\).
\begin{enumerate}
\item If \(Q\) is irreducible and \(\theta\) is a root, then \(R_Q\simeq\F_{q^2}\), and
\[
 G_Q\xrightarrow{\ \sim\ }\mu_{q+1},\qquad
 [\overline f]\longmapsto f(\theta)^{q-1}.
\]
\item If \(Q=(T-r)(T-s)\), where \(r,s\in k\) and \(r\ne s\), then \(R_Q\simeq k\times k\), and
\[
 G_Q\xrightarrow{\ \sim\ }k^\times,\qquad
 [\overline f]\longmapsto f(r)/f(s).
\]
\item If \(Q=(T-r)^2\), put \(\varepsilon=\overline T-r\). Then \(R_Q\simeq k[\varepsilon]/(\varepsilon^2)\), the ring of dual numbers, and
\[
 G_Q\xrightarrow{\ \sim\ }(k,+),\qquad
 [a+b\varepsilon]\longmapsto b/a\quad(a\ne0).
\]
\end{enumerate}
In these coordinates, the maps induced by \(q\)-Frobenius on \(G_Q\) are inversion, the identity, and the zero map, respectively.
\end{proposition}
\begin{proof}
In the first case, \(u\mapsto u^{q-1}\) has kernel \(k^\times\) and image the norm-one group, by Theorem~\ref{thm:hilbert90}.
In the second case, the Chinese remainder theorem identifies units with pairs \((a,b)\in k^\times\times k^\times\); the kernel of \((a,b)\mapsto a/b\) is the diagonal scalar group.
In the third case, \(a+b\varepsilon\) is a unit exactly when \(a\ne0\), and each projective class has a unique representative \(1+c\varepsilon\). The identity
\[
 (1+c\varepsilon)(1+c'\varepsilon)=1+(c+c')\varepsilon
\]
gives the additive group. Since finite fields are perfect, these exhaust the factorization types of a quadratic.
Finally, \(u^q=u^{-1}\) in the first coordinate group; every element of the second algebra is fixed by the \(q\)th-power map; and \((a+b\varepsilon)^q=a\) in the third algebra. Thus every projective class in the last case maps to the identity class, which has additive coordinate zero.
\end{proof}

\begin{corollary}[Coordinates of a root set in a quadratic quotient]\label{v85:cor-rootset-quotient-coordinates}
Let \(S\subseteq k\), put \(p_S(T)=\prod_{a\in S}(T-a)\), and assume \(\gcd(p_S,Q)=1\).
The coordinates of \([\overline{p_S}]\in G_Q\) in the three groups above are
\begin{align*}
 &\prod_{a\in S}\frac{\theta^q-a}{\theta-a},
 &&Q\text{ irreducible};\\
 &\frac{p_S(r)}{p_S(s)}=\prod_{a\in S}\frac{r-a}{s-a},
 &&Q=(T-r)(T-s),\ r\ne s;\\
 &\frac{p_S'(r)}{p_S(r)}=\sum_{a\in S}\frac1{r-a},
 &&Q=(T-r)^2.
\end{align*}
The coprimality condition is automatic in the irreducible case. In the other two cases it requires \(r,s\notin S\) and \(r\notin S\), respectively.
\end{corollary}
\begin{proof}
The first formula is \eqref{v85:eq-norm-root-product}; the second follows by expanding both values.
In the third case, \(p_S(r+\varepsilon)=p_S(r)+p_S'(r)\varepsilon\), so the coordinate is \(p_S'(r)/p_S(r)\).
Differentiating the product of root factors and dividing by \(p_S(T)\) gives the rational-function identity
\[
 \frac{p_S'(T)}{p_S(T)}=\sum_{a\in S}\frac1{T-a}.
\]
\leavevmode\phantomsection\label{v85:def-logarithmic-derivative}\ignorespaces
The left-hand side is the logarithmic derivative of \(p_S\). Evaluation at \(r\notin S\) gives the result.
\end{proof}
For fixed \(Q\), the coordinates depend on the chosen roots and their order, whereas the projective residue class they describe is unchanged. Disjoint unions of root sets become products in the first two coordinate groups and sums in the third.

\begin{proposition}[A linear condition as one projective unit class]\label{v85:prop-functional-unit-fiber}
Let \(\ell:R_Q\to k\) be linear, and suppose that \((u,v)\mapsto\ell(uv)\) is nondegenerate.
There is a unit \(\kappa\in R_Q^\times\) such that \(\ker\ell=k\kappa\). For every root set with \(\gcd(p_S,Q)=1\),
\[
 \ell(\overline{p_S})=0\quad\Longleftrightarrow\quad
 [\overline{p_S}]=[\kappa]\text{ in }G_Q.
\]
\end{proposition}
\begin{proof}
Nondegeneracy implies \(\ell\ne0\), so its kernel is one-dimensional. If a nonzero nonunit \(v\) were in the kernel, multiplication by \(v\) would not be invertible. Its image \(vR_Q\) contains \(v\), so it would be exactly the one-dimensional space \(kv=\ker\ell\). This would give \(\ell(vu)=0\) for every \(u\), contrary to nondegeneracy.
Thus every nonzero element of the kernel is a unit; choose one as \(\kappa\).
\end{proof}
A linear vanishing condition on the quotient algebra has therefore become a prescribed norm-one product, value ratio, or sum of reciprocals. The nondegeneracy hypothesis is essential:
in the ring of dual numbers, \(\Tr(a+b\varepsilon)=2a\) does not detect \(b\), whereas the coefficient functional \(a+b\varepsilon\mapsto b\) still gives a nondegenerate pairing.
This explains concretely why repeated roots require the coefficient pairing rather than the algebra trace alone.

\subsection{Weighted trace forms and symmetric multiplication matrices}\label{v90:subsec-weighted-trace}
Let $E=\F_{q^n}$ and $k=\F_q$. For $\lambda\in E^\times$, define the nondegenerate symmetric bilinear form
\[
B_\lambda(x,y)=\Tr_{E/k}(\lambda xy).
\]
Relative to a nondegenerate symmetric pairing \(B\), a linear operator \(A\) is called \textbf{self-adjoint} if
\(B(Ax,y)=B(x,Ay)\) for all \(x,y\).
Equivalently, \(A^*=A\); the general definition and existence and uniqueness of the adjoint map are given at \ref{def:adjoint-map} below.
\begin{proposition}[Self-adjoint multiplication and simultaneous symmetric matrices]\label{v90:prop-symmetric-mult}
For all $a,x,y\in E$, one has $B_\lambda(ax,y)=B_\lambda(x,ay)$. Moreover, there exist $\lambda\in E^\times$ and a $k$-basis of $E$ in which the matrices of all multiplication operators $m_a:x\mapsto ax$ are symmetric.
\end{proposition}
\begin{proof}
Self-adjointness follows from \(\lambda axy=\lambda xay\). It suffices to give some \(B_\lambda\) an orthonormal basis.
We use the classification criterion in Lemma~\ref{lem:finite-symmetric-form} below; its proof is independent of this proposition.
Fix a \(k\)-basis, let \(G_\lambda\) be the Gram matrix of \(B_\lambda\), and set \(D_1=\det G_1\in k^\times\).
If \(M_\lambda\) is the matrix of multiplication by \(\lambda\), then \(G_\lambda=G_1M_\lambda\), so
\[
\det G_\lambda=\Nm_{E/k}(\lambda)D_1.
\]
When \(q\) is odd, surjectivity of the norm allows a choice of \(\lambda\) making this determinant a nonzero square; Lemma~\ref{lem:finite-symmetric-form} then supplies an orthonormal basis.
When \(q\) is even, \(B_\lambda(x,x)=\Tr(\lambda x^2)\) is not identically zero, since the trace pairing is nondegenerate and squaring is bijective.
The form is therefore nonalternating, and the same lemma again supplies an orthonormal basis.
In such a basis the Gram matrix is the identity, so self-adjointness makes every \(M_a\) equal to its transpose.
\end{proof}
\begin{remark}
The weighted form need not equal the ordinary trace form $B_1$. Simultaneous symmetry of multiplication matrices does not assert that the ordinary trace pairing admits a self-dual basis. For the standard classification of bilinear forms see \cite{HandbookFiniteFields}.
\end{remark}

\section{Tower Bases and Transitivity of Duality}
\begin{proposition}\label{prop:dual-tower}
Let \(k\subseteq L\subseteq E\). If \((u_i),(u_i^*)\) are dual bases for \(L/k\),
and \((v_j),(v_j^*)\) are dual bases for \(E/L\), then \((u_iv_j)_{i,j}\) and \((u_i^*v_j^*)_{i,j}\) are dual bases for \(E/k\).
\end{proposition}
\begin{proof}
They are bases by the tower law. Transitivity of the trace gives
\[
\Tr_{E/k}(u_i u_r^*v_jv_s^*)
=\Tr_{L/k}\bigl(u_i u_r^*\Tr_{E/L}(v_jv_s^*)\bigr)
=\delta_{js}\delta_{ir}.\qedhere
\]
\end{proof}

\phantomsection\label{def:adjoint-map}
Suppose finite-dimensional spaces \(V,W\) carry nondegenerate symmetric pairings \(B_V,B_W\). The \textbf{adjoint map} of a linear map \(A:V\to W\) is the linear map satisfying
\[
B_W(Av,w)=B_V(v,A^*w)\qquad(v\in V,\ w\in W).
\]
Here \(A^*:W\to V\). It exists uniquely: for fixed \(w\), the left-hand side is a linear functional on \(V\),
and nondegeneracy expresses it uniquely as pairing with some \(A^*w\). This correspondence is linear in \(w\).
When \(V=W\), \(A^*\) is also called the adjoint operator. It depends on the chosen pairing and is distinct from the companion matrix of a polynomial.

\begin{proposition}[Relative trace is adjoint to subfield inclusion]\label{prop:relative-trace-adjoint}
Let \(k\subseteq L\subseteq E\) be finite fields, and equip \(E\) and \(L\) with their trace pairings over \(k\).
For \(x\in E\) and \(y\in L\),
\[
\Tr_{E/k}(xy)=\Tr_{L/k}\bigl(\Tr_{E/L}(x)y\bigr).
\]
Thus, viewing \(L\) as a \(k\)-subspace of \(E\),
\[
L^\perp=\ker\Tr_{E/L}.
\]
\end{proposition}
\begin{proof}
The first identity follows from transitivity of the trace and the \(L\)-linearity of \(\Tr_{E/L}\).
Thus \(x\in L^\perp\) is equivalent to the identity
\[
\Tr_{L/k}\bigl(\Tr_{E/L}(x)y\bigr)=0.
\]
This identity must hold for every \(y\in L\). Nondegeneracy of the trace pairing for \(L/k\) makes this equivalent to \(\Tr_{E/L}(x)=0\).
\end{proof}

In particular, even when \(p\mid[E:L]\) makes the trace pairing for \(E/k\) identically zero on \(L\times L\),
the pairing on all of \(E\) remains nondegenerate. A subfield, viewed as a multiplicatively closed set and as a linear subspace equipped with a pairing,
carries compatible but distinct structures.

\section{Discriminants and Matrices of Conjugates}
\begin{definition}\label{def:discriminant}
For \(b_1,\ldots,b_n\in E\), define the \textbf{trace discriminant} by
\[
\Delta(b_1,\ldots,b_n)=\det\bigl(\Tr_{E/k}(b_ib_j)\bigr)_{i,j=1}^n.
\]
\end{definition}
Let \(V=(b_j^{q^i})_{0\le i<n,\,1\le j\le n}\). Computing the entries of the matrix product gives the trace Gram matrix \(V^{\mathsf t}V\), so
\[
\Delta(b_1,\ldots,b_n)=(\det V)^2.
\]
Since \(\det V\) belongs to \(E\), this identity alone does not imply that the discriminant is a square in \(k\).
Nondegeneracy of the trace pairing shows that \(\Delta\ne0\) if and only if these elements form a basis.
\leavevmode\phantomsection\label{def:square-class}\ignorespaces Two nonzero elements belong to the same \textbf{square class} if their quotient is a square in \(k\).
A change of basis with matrix \(C\) multiplies the discriminant by \((\det C)^2\), so its square class in \(k^\times\) is independent of the basis.

For a basis \((b_j)_{j=1}^n\) with trace-dual basis \((b_j^*)_{j=1}^n\), put
\(W=(\sigma^i(b_j^*))_{0\le i<n,\,1\le j\le n}\). The duality relations give
\[
(V^{\mathsf t}W)_{rs}=\sum_{i=0}^{n-1}\sigma^i(b_r b_s^*)
=\Tr_{E/k}(b_r b_s^*)=\delta_{rs}.
\]
Hence \(VW^{\mathsf t}=I\), and the inverse of~\eqref{eq:conjugate-tensor-map} is
\begin{equation}\label{eq:inverse-conjugate-tensor-map}
\mathcal J^{-1}\bigl((z_i)_{i=0}^{n-1}\bigr)
=\sum_{j=1}^n\left(\sum_{i=0}^{n-1}z_i\sigma^i(b_j^*)\right)\otimes b_j.
\end{equation}
When \(z_i=\sigma^i(x)\), the expression in parentheses is \(\Tr_{E/k}(xb_j^*)\),
so recovering \(1\otimes x\) gives precisely the trace-dual coordinate formula. Conjugate evaluation and trace duality thus describe the forward and inverse directions of the same isomorphism.

\begin{proposition}\leavevmode\phantomsection\label{prop:power-discriminant}
Let \(E=k(\alpha)\), \(f=m_{\alpha,k}\), and \(\deg f=n\). Then
\[
\Delta(1,\alpha,\ldots,\alpha^{n-1})
=(-1)^{n(n-1)/2}\Nm_{E/k}(f'(\alpha)).
\]
\end{proposition}
\begin{proof}
Let the conjugate roots be \(\alpha_i\). The matrix of conjugates is a Vandermonde matrix, so its determinant squared is \(\prod_{i<j}(\alpha_j-\alpha_i)^2\).
On the other hand, \(\prod_i f'(\alpha_i)=\prod_i\prod_{j\ne i}(\alpha_i-\alpha_j)\),
and each pair of indices contributes one minus sign, giving the formula.
\end{proof}
Thus the discriminant, the factor \(f'(\alpha)\) in the dual-basis formula, and the distinctness of conjugate roots are different expressions of the same matrix of conjugates.

The product of root differences also applies to reducible polynomials.
Let \(f\in k[X]\) be monic of degree \(n\ge1\), with roots
\(\alpha_1,\ldots,\alpha_n\) in a common splitting field, counted with multiplicity.
Its \textbf{polynomial discriminant} is
\phantomsection\label{v105:def-polynomial-discriminant}
\begin{equation}\label{v105:eq-polynomial-discriminant}
\operatorname{disc}(f)=\prod_{1\le i<j\le n}(\alpha_j-\alpha_i)^2.
\end{equation}
This expression is independent of the numbering of the roots and is fixed by Frobenius, so it belongs to \(k\).
It is nonzero exactly when \(f\) has no repeated root.
When \(f=m_{\alpha,k}\) is irreducible, the proof of the preceding proposition gives
\(\operatorname{disc}(f)=\Delta(1,\alpha,\ldots,\alpha^{n-1})\).

\begin{proposition}[The discriminant and the parity of a factorization]\label{v105:prop-discriminant-parity}
Let \(q\) be odd, and let \(f\in\F_q[X]\) be monic and squarefree of degree \(n\ge1\),
with exactly \(r\) monic irreducible factors in \(\F_q[X]\).
Let \(\chi:\F_q^\times\to\{1,-1\}\) take the value one on squares and minus one on nonsquares.
Then
\begin{equation}\label{v105:eq-discriminant-parity}
\boxed{\quad
\chi\bigl(\operatorname{disc}(f)\bigr)=(-1)^{n-r}.
\quad}
\end{equation}
\end{proposition}
\begin{proof}
Number the distinct roots of \(f\) in a common finite splitting field and put
\[
V=(\alpha_i^{j-1})_{1\le i,j\le n},\qquad
D=\det V=\prod_{i<j}(\alpha_j-\alpha_i).
\]
The roots are distinct, so \(D\ne0\), and \(\operatorname{disc}(f)=D^2\).
Write the action of Frobenius on the root indices as a permutation \(\pi\), so that
\(\alpha_i^q=\alpha_{\pi(i)}\).
Raising every entry of \(V\) to the \(q\)-th power permutes its rows by \(\pi\), giving
\[
D^q=\operatorname{sgn}(\pi)D.
\]
Each Frobenius orbit of roots corresponds to one irreducible factor.
If the factor degrees are \(n_1,\ldots,n_r\), these are exactly the cycle lengths of \(\pi\).
Consequently,
\[
\operatorname{sgn}(\pi)
=\prod_{j=1}^r(-1)^{n_j-1}
=(-1)^{n-r}.
\]
Cyclicity of the multiplicative group gives \(\chi(a)=a^{(q-1)/2}\), identifying the field elements \(\pm1\) with the corresponding signs.
Hence
\[
\chi\bigl(\operatorname{disc}(f)\bigr)
=(D^2)^{(q-1)/2}=D^{q-1}
=\operatorname{sgn}(\pi)=(-1)^{n-r}.
\]
\end{proof}

Equation~\eqref{v105:eq-discriminant-parity} is \textbf{Stickelberger's discriminant formula}.
It expresses the sign of the Frobenius permutation of the roots as a square class in the base field.
For irreducible \(f\), the sign is \((-1)^{n-1}\);
for a general squarefree polynomial, it determines only the parity of the number of irreducible factors.
This restriction also enters counts of irreducible polynomials with prescribed discriminant;
see~\cite[Appendix~A]{OstafeShparlinski2026Thin}.

\begin{example}[The same discriminant with different cubic factorizations]\label{v105:ex-cubic-same-discriminant}
In \(\F_5[X]\), take
\[
f=X^3+X+1,\qquad g=X^3-X=X(X-1)(X+1).
\]
The values of \(f\) at \(0,1,2,3,4\) are respectively \(1,3,1,1,4\), so it has no root in the base field.
A reducible cubic has a linear factor; therefore \(f\) is irreducible.
If \(\alpha\) is a root of \(f\), then \(\alpha^3=-\alpha-1\).
In the basis \(1,\alpha,\alpha^2\), multiplication by \(f'(\alpha)=1+3\alpha^2\) has matrix
\[
\begin{pmatrix}
1&2&0\\
0&3&2\\
3&0&3
\end{pmatrix},
\qquad
\det\begin{pmatrix}
1&2&0\\
0&3&2\\
3&0&3
\end{pmatrix}=1\quad\text{in }\F_5.
\]
Proposition~\ref{prop:power-discriminant} gives
\(\operatorname{disc}(f)=-\Nm_{\F_{125}/\F_5}(f'(\alpha))=-1=4\).
On the other hand, the roots of \(g\) are \(0,1,-1\), and their differences give
\(\operatorname{disc}(g)=4\) directly.
Thus the two polynomials have the same nonzero discriminant, although one is irreducible and the other splits completely.
They have respectively one and three Frobenius orbits, in agreement with~\eqref{v105:eq-discriminant-parity}.
\end{example}

\section{Existence of Self-Dual Bases}
\leavevmode\phantomsection\label{def:self-dual-basis}\ignorespaces A basis equal to its dual basis, that is, satisfying \(\Tr(b_i b_j)=\delta_{ij}\), is called a \textbf{self-dual basis}.
No normality condition is imposed here.

\begin{lemma}\leavevmode\phantomsection\label{lem:finite-symmetric-form}
On a finite-dimensional space of positive dimension over a finite field \(k\), a nondegenerate symmetric bilinear form admits a basis with identity Gram matrix if and only if:
\begin{enumerate}
\item in characteristic different from two, its Gram determinant is a square in \(k\);
\item in characteristic two, it is not an alternating form, that is, some \(v\) satisfies \mbox{\(B(v,v)\ne0\)}.
\end{enumerate}
\end{lemma}
\begin{proof}
First suppose the characteristic is not two. If \(B(v,v)=0\) for every vector, then
\[
2B(x,y)=B(x+y,x+y)-B(x,x)-B(y,y)=0
\]
for all \(x,y\), contradicting nondegeneracy. Thus one can choose \(v\) with \(B(v,v)\ne0\).
Every vector \(x\) has a unique decomposition
\[
x=\frac{B(x,v)}{B(v,v)}v+
\left(x-\frac{B(x,v)}{B(v,v)}v\right)\in kv\oplus(kv)^\perp.
\]
The restriction to the orthogonal complement remains nondegenerate: a vector there orthogonal to both the complement and \(v\) is orthogonal to the entire space and hence is zero.
Induction gives a diagonal Gram matrix with nonzero entries.

For any \(a,b\in k^\times\), the sets \(\{ax^2:x\in k\}\) and \(\{1-by^2:y\in k\}\) each contain \((q+1)/2\) elements,
so they intersect. Thus the binary diagonal form \(ax^2+by^2\) represents \(1\), and a line with Gram coefficient \(1\) splits off.
Its orthogonal complement is a nondegenerate line whose coefficient has square class \(ab\), since a change of basis multiplies the determinant by a nonzero square.
Applying this operation successively to two as yet unnormalized diagonal entries gives \(\operatorname{diag}(1,\ldots,1,c)\), where \(c\) has the square class of the original determinant.
If the determinant is a square, scaling also gives coefficient \(1\) on the last line.
Necessity follows from the same change-of-basis formula for determinants; the conclusion includes the one-dimensional case.

Now suppose the characteristic is two. The form with identity matrix is not alternating, so this condition is necessary.
If some \(v\) satisfies \(B(v,v)\ne0\), the squaring map is bijective on a finite field, so \(v\) can be scaled to have self-pairing \(1\).
The same orthogonal decomposition successively splits off unit lines until the orthogonal complement is zero or is a nondegenerate alternating space.
A nondegenerate alternating space decomposes into two-dimensional blocks: choose \(u\ne0\) and, by nondegeneracy, choose \(v\) with \(B(u,v)=1\).
The matrix on \((u,v)\) is then \(\left(\begin{smallmatrix}0&1\\1&0\end{smallmatrix}\right)\);
split off this block and apply induction to its orthogonal complement.

At least one unit line \(kw\) has already been obtained. If \(u,v\) is an alternating block orthogonal to \(w\), then
\[
w+u,\qquad w+v,\qquad w+u+v
\]
each have self-pairing \(1\), and their pairwise pairings are \(1+1=0\).
Their Gram matrix is the identity, so they are linearly independent and form a basis of this three-dimensional space.
Thus one unit line and one alternating block can be replaced by three unit lines.
Repeating this for each block gives an orthonormal basis of the entire space.
\end{proof}

\begin{theorem}\label{thm:self-dual-basis}
The extension \(\F_{q^n}/\F_q\) admits a self-dual basis if and only if \(q\) is even or both \(q,n\) are odd.
\end{theorem}
\begin{proof}
If \(q\) is even, \(B(x,x)=\Tr(x^2)=\Tr(x)^2\). The trace is surjective, so the form is not alternating; the preceding lemma gives a self-dual basis.

Suppose \(q\) is odd. Choose any \(k\)-basis \(b_0,\ldots,b_{n-1}\) and put
\[
V=(b_j^{q^i})_{0\leq i,j<n},\qquad D=\det V.
\]
Theorem~\ref{thm:conjugate-coordinates} shows that \(V\) is invertible, so \(D\ne0\).
The Gram matrix is \(V^{\mathsf T}V\) and has determinant \(D^2\).
Raising \(D\) to the \(q\)-th power cyclically permutes the rows of \(V\), so \(D^q=(-1)^{n-1}D\).
If \(n\) is odd, then \(D\in k\) and the determinant is a square.
If \(n\) is even, then \(D\notin k\); but any square root of \(D^2\) in \(k\) would have to equal \(D\) or \(-D\), a contradiction.
The determinant is therefore a nonsquare. The preceding lemma proves the result.
\end{proof}

This criterion comes from the Frobenius action on the determinant of the matrix of conjugates.
Self-duality is a coordinate condition, whereas normality requires the basis vectors to form a single Frobenius orbit.
The conditions under which both hold will be discussed separately in the chapter on normal bases.

\par\Needspace{11\baselineskip}
\section*{Exercises}
\phantomsection
\addcontentsline{toc}{section}{Exercises}
\begin{exercise}\label{ex:dual-f4}
Let \(\F_4=\F_2(\alpha)\) with \(\alpha^2+\alpha+1=0\).
Find the dual basis of \((1,\alpha)\) and prove that \((\alpha,\alpha^2)\) is self-dual.
\end{exercise}
\begin{exercise}\label{ex:trace-orthogonal}
Prove that \(k^\perp=\ker\Tr_{E/k}\), viewing \(k\) as a one-dimensional \(k\)-subspace of \(E\).
Also prove that \((U+V)^\perp=U^\perp\cap V^\perp\).
\end{exercise}
\begin{exercise}\label{ex:no-selfdual-power}
Prove that for \(n>1\), the power basis \((1,\alpha,\ldots,\alpha^{n-1})\) cannot be self-dual.
Hint: for \(n\geq3\), compare \(\Tr(1\cdot\alpha^2)\) and \(\Tr(\alpha\cdot\alpha)\); treat \(n=2\) separately.
\end{exercise}

\begin{exercise}\label{ex:supp-discriminant}
Let \(\F_9=\F_3(\alpha)\) with \(\alpha^2=-1\). Compute the trace Gram matrix and discriminant of \((1,\alpha)\) directly, then check the result using the derivative-norm formula.
\end{exercise}

\begin{exercise}[The sum of a self-dual basis in characteristic two]\label{ex:wu-even-selfdual-sum}
Let \(q\) be even and let \((b_1,\ldots,b_n)\) be a self-dual basis of \(\F_{q^n}/\F_q\).
Prove that \(\Tr(b_i)=1\) for every index and that \(\sum_i b_i=1\).
Explain why, when \(n\) is even, this does not contradict \(\Tr(1)=0\).
\end{exercise}

\begin{exercise}[The orthogonal complement of a sum of subfields]\label{ex:wu2-sum-subfields}
Let \(a,b\mid n\) and put \(L=\F_{q^a}\) and \(M=\F_{q^b}\) in \(E=\F_{q^n}\).
Find the dimension of \(L+M\) and express its trace-orthogonal complement as the intersection of two relative-trace kernels.
Determine when \(L+M\) is itself a subfield.
\end{exercise}

\begin{exercise}[Derivatives, coefficient pairings, and trace degeneracy]\label{ex:wu7-trace-radical}
Let \(f\in k[T]\) be monic with \(\deg f=d\ge1\), and put \(A=k[T]/(f)\) and \(\alpha=\overline T\).
Let \(\lambda\) extract the coefficient of \(T^{d-1}\) in the standard remainder. Write
\[
B_\lambda(u,v)=\lambda(uv),\qquad B_{\Tr}(u,v)=\Tr_{A/k}(uv).
\]
\begin{enumerate}
\item Put \(h=\gcd(f,f')\). Prove that
\[
\rad B_{\Tr}=\ker M_{f'(\alpha)}
=\bigl((f/h)(\alpha)\bigr),\qquad
\rank B_{\Tr}=d-\deg h,
\]
where \(\bigl((f/h)(\alpha)\bigr)\) denotes the ideal of \(A\) generated by \((f/h)(\alpha)\).
Explain how this transforms the nondegenerate coefficient pairing into the potentially degenerate trace pairing.
\item Suppose \(k\) has characteristic \(p>0\). Take \(f=(T-1)^e\) with \(e\ge1\) and put \(\varepsilon=\alpha-1\).
Find the Gram matrices, ranks, and radicals of \(B_\lambda\) and \(B_{\Tr}\) in the basis \((1,\varepsilon,\ldots,\varepsilon^{e-1})\).
Compare the cases \(p\mid e\) and \(p\nmid e\), and compare the radicals with the nilradical; check in particular \(e=p\) and, in odd characteristic, \(e=2\).
\item Return to Example~\ref{exm:mixed-algebra-trace}, where
\(k=\F_2\), \(P=T^2+T+1\), and \(f=(T+1)^2P^3\).
Compute \(f'\), \(h\), and a basis of \(\rad B_{\Tr}\) directly, and check the ranks of both pairings using the derivative formula and the local multiplicity formula.
Prove that \(e_1=\alpha P(\alpha)^3\) is a nonzero idempotent satisfying
\[
B_{\Tr}(e_1,v)=0\quad(v\in A),\qquad B_\lambda(e_1,1)=1.
\]
\item Let \(u\) be an indeterminate and put \(k=\F_p(u^p)\) and \(E=\F_p(u)\).
Example~\ref{ex:inseparable-rational-extension} in the appendix gives \(E\simeq k[T]/(T^p-u^p)\).
Compute the Gram matrices of the two pairings in the basis \((1,u,\ldots,u^{p-1})\), and compare their radicals with the nilradical of \(E\).
Explain which data agree with the case \(e=p\) in part two, which algebraic properties differ, and the role of the perfect-field hypothesis in the criterion for reducedness.
\end{enumerate}
\end{exercise}

\begin{exercise}[Changing multiplicities while preserving the Frobenius periodic part]\label{ex:wu10-local-trace}
Replace the algebra in Example~\ref{exm:mixed-algebra-trace} by
\[
A'=\F_2[T]/((T+1)^3(T^2+T+1)^3).
\]
Find its nilradical, the dimension of its Frobenius periodic part, and the ranks of its algebra-trace pairing and top-coefficient pairing.
Prove that the radical of the trace pairing is now exactly the nilradical.
For each \(r\ge1\), count the solutions of \(u^{2^r}=u\) in \(A'\), and explain which data agree with the original example and which change with the multiplicities.
\end{exercise}

\begin{exercise}[Three quadratic residues of the same root set]\label{v85:ex-quadratic-coordinates}
In \(k=\F_5\), take \(S=\{0,1\}\), and let
\[
 Q_1=T^2+2,\qquad Q_2=(T-2)(T-3),\qquad Q_3=(T-2)^2.
\]
Prove that \(Q_1\) is irreducible. Compute the image of \([p_S]\) in each projective unit group using Proposition~\ref{v85:prop-quadratic-unit-groups};
use \(\theta^2=3\) in the first case and the ordered roots \((2,3)\) in the second.
Determine the effect of the fifth-power Frobenius on each image, and explain why the third action is not a group automorphism.
\end{exercise}

\begin{exercise}[A node-change matrix and multiplication]\label{v97:ex-node-change}
Over \(\F_7\), take \(f=X(X-1)(X-2)\) and \(g=(X-3)(X-4)(X-5)\), ordering their roots as \((0,1,2)\) and \((3,4,5)\).
Compute the matrix \(M\) of Proposition~\ref{v97:prop-cauchy-node-change}, and use \(P=1+2X+X^2\) to check \(v_\beta=M^{\mathsf T}v_\alpha\).
Explain why every square minor of \(M\) is nonzero, whereas the linear map \(k[X]/(f)\to k[X]/(g)\) defined by using the same representative of degree less than three does not preserve multiplication.
\end{exercise}

\chapter{Linearized Polynomials and Operator Algebras}
\label{chap:linearized}
Finite-dimensional subspaces of finite extensions have linearized root polynomials. Composition corresponds to operator multiplication, while the condition that the coefficients lie in the base field determines whether the operator commutes with Frobenius. Trace adjoints relate kernels to images.

\section{Polynomial Representations of Linear Maps}\label{sec:linearized-definition}
Let \(k=\F_q\subseteq E=\F_{q^n}\), as before.
Whenever all roots of a polynomial are considered in this chapter, we work in a finite extension over which it splits completely;
when considering a kernel on \(E\), we take only the roots belonging to \(E\). A polynomial of the form
\[
L(X)=a_0X+a_1X^q+\cdots+a_sX^{q^s},\qquad a_i\in E,
\]
is called a \(q\)-\textbf{linearized polynomial}.
If \(a_s\ne0\), then \(s\) is its \(q\)-degree; its ordinary degree is \(q^s\).
The map \(L\) is \(k\)-linear on every extension field containing its coefficients, since the \(q\)th-power map preserves addition and fixes \(k\).
Its zero set is therefore a linear subspace, called the \textbf{root space} of \(L\).

\begin{proposition}\label{prop:linearized-rootspace}
The roots of a nonzero linearized polynomial of \(q\)-degree \(s\) form a \(k\)-vector space of dimension at most \(s\).
The polynomial has no repeated roots if and only if \(a_0\ne0\).
\end{proposition}
\begin{proof}
The root set is the kernel of a linear map and hence is a subspace.
The usual bound on the number of roots gives at most \(q^s\) roots, so the dimension is at most \(s\).
The derivative is the constant \(a_0\). If it is nonzero, there are no repeated roots. If it is zero, all exponents are divisible by \(p\);
taking \(p\)th roots of the coefficients in the finite coefficient field expresses \(L\) as the \(p\)th power of a nonconstant polynomial, so it has repeated roots.
\end{proof}

If \(a_0\ne0\), then \(L\) has \(q^s\) distinct roots in its splitting field, so its full root space has dimension exactly \(s\); its kernel in a fixed field \(E\) may have smaller dimension. Section~\ref{sec:linearized-division} treats all roots and their multiplicities, while Corollary~\ref{v70:cor-all-extension-kernels} describes the kernels in finite extensions when the coefficients lie in \(k\).

\begin{theorem}[Polynomial representations of all linear maps]\label{thm:linear-map-representation}
Every \(k\)-linear map \(A:E\to E\) has a unique representation
\[
A(x)=a_0x+a_1x^q+\cdots+a_{n-1}x^{q^{n-1}},\qquad a_i\in E.
\]
\end{theorem}
\begin{proof}
The coefficient tuples form an \(n^2\)-dimensional \(k\)-space, and \(\End_k(E)\) also has dimension \(n^2\).
If two representations define the same function, their difference vanishes on all \(q^n\) elements of \(E\),
whereas a nonzero difference polynomial has degree at most \(q^{n-1}\) and cannot have so many roots.
The map from coefficient tuples to linear maps is therefore injective, and comparing dimensions proves surjectivity.
\end{proof}

\begin{theorem}[Moore determinant]\leavevmode\phantomsection\label{thm:moore}
Let \(x_1,\ldots,x_s\in E\) with \(1\leq s\leq n\). The Moore determinant satisfies
\[
\det\bigl(x_j^{q^i}\bigr)_{0\leq i<s,\ 1\leq j\leq s}\ne0
\]
if and only if \(x_1,\ldots,x_s\) are linearly independent over \(k\).
\end{theorem}
\begin{proof}
A \(k\)-linear relation, raised to each \(q^i\)th power, gives the same relation among the matrix columns, so the determinant vanishes.
Conversely, if the determinant vanishes, there exist \(a_0,\ldots,a_{s-1}\in E\), not all zero,
such that the linearized polynomial \(L(X)=\sum_i a_iX^{q^i}\) vanishes at every \(x_j\).
If these elements are independent over \(k\), their span contains \(q^s\) elements, all roots of \(L\),
contradicting \(\deg L\leq q^{s-1}\).
\end{proof}

Given a basis and its images, the linearized representation is obtained by solving the linear system defined by the Moore matrix.
Inverting this matrix thus gives all coefficients of the representation.

\Needspace{10\baselineskip}
\section{Subspace Polynomials}
\begin{theorem}\label{thm:subspace-polynomial}
Let \(V\subseteq E\) be a \(k\)-subspace of dimension \(s\). Define its \textbf{subspace polynomial} by
\[
L_V(X)=\prod_{v\in V}(X-v).
\]
It is a monic linearized polynomial of \(q\)-degree \(s\) with nonzero linear coefficient.
\end{theorem}
\begin{proof}
We induct on \(s=\dim_kV\). If \(s=0\), then \(V=\{0\}\), so \(L_V=X\) and the assertion holds.
Suppose \(s>0\). Choose a subspace \(W\subset V\) of codimension one and an element \(u\in V\setminus W\) such that \(V=W\oplus ku\).
The additive cosets \(cu+W\) are pairwise disjoint, so
\[
L_V(X)=\prod_{c\in k}\prod_{w\in W}(X-cu-w)
=\prod_{c\in k}L_W(X-cu).
\]
By induction, \(L_W\) is \(k\)-linearized, and hence \(L_W(X-cu)=L_W(X)-cL_W(u)\).
Put \(b=L_W(u)\). Since the roots of \(L_W\) are exactly \(W\), we have \(b\ne0\).
In \(\prod_{c\in k}(Y-c)=Y^q-Y\), rescale the variable to \(Y/b\) to obtain
\[
\prod_{c\in k}(Y-cb)=Y^q-b^{q-1}Y.
\]
Substituting \(Y=L_W(X)\) gives
\[
L_V(X)=L_W(X)^q-b^{q-1}L_W(X).
\]
If \(L_W=\sum_{i=0}^{s-1}a_iX^{q^i}\) with \(a_{s-1}=1\) and \(a_0\ne0\),
the resulting exponents are powers of \(q\), the leading term is \(X^{q^s}\), and the linear coefficient is \(-b^{q-1}a_0\ne0\).
This completes the induction.
\end{proof}

\begin{corollary}\label{cor:subspace-coefficients}
All coefficients of \(L_V\) belong to \(k\) if and only if \(\sigma(V)=V\).
If \(U\subseteq V\) are subspaces, then
\[
L_V=L_{L_U(V)}\circ L_U.
\]
\end{corollary}
\begin{proof}
Applying \(\sigma\) to the coefficients of \(L_V\) gives the monic polynomial with root set \(\sigma(V)\),
so uniqueness of the root-set polynomial proves the first assertion.
For the second, the kernel of \(L_U:V\to L_U(V)\) is \(U\), and hence the image has dimension \(\dim V-\dim U\).
The composed polynomial is monic of ordinary degree \(q^{\dim V}\) and vanishes on all of \(V\), so it equals \(L_V\).
\end{proof}

The root polynomial of a finite subfield has the special form \(X^{q^s}-X\),
while that of a finite \(k\)-subspace has the more general linearized form.
Both reflect closure of the root set under the corresponding operations.

\begin{definition}\label{v70:def-scalar-field}
For a \(k\)-subspace \(W\subseteq E\), define its \textbf{maximal field of scalars} by
\[
 H(W)=\{a\in E:aW\subseteq W\}.
\]
In particular, \(H(0)=E\). When the base field needs to be specified, we write \(H_k(W)\).
\end{definition}

\begin{proposition}[Recovering the scalar field from the support of the root polynomial]\label{v70:prop-scalar-support}
The set \(H(W)\) is a subfield containing \(k\), and is the largest intermediate field over which \(W\) is a vector space.
Suppose \(W\ne0\), put \(r=\dim_kW\), and write
\[
 L_W(X)=\sum_{i=0}^r c_iX^{q^i},\qquad c_r=1,\quad c_0\ne0.
\]
Then
\begin{equation}\label{v70:eq-scalar-support}
 H(W)=\F_{q^d},\qquad
 d=\gcd\bigl(n,\{i:0\le i\le r,\ c_i\ne0\}\bigr).
\end{equation}
In particular, \(d\mid r\). Neither Frobenius stability of \(W\) nor membership of the coefficients \(c_i\) in \(k\) is required.
\end{proposition}
\begin{proof}
The assertion is immediate for \(W=0\). Suppose \(W\ne0\). The set \(H(W)\) contains \(k\) and is closed under addition, subtraction, and multiplication.
If \(0\ne a\in H(W)\), multiplication by \(a\) is injective, so finiteness gives \(aW=W\), and hence \(a^{-1}W=W\).
Thus \(H(W)\) is a subfield. Its definition also shows that \(W\) is an \(H(W)\)-space and that every intermediate field over which \(W\) is a vector space is contained in \(H(W)\).

For \(a\in E^\times\), comparison of root sets and leading terms gives
\[
 aW=W\quad\Longleftrightarrow\quad
 L_W(aX)=a^{q^r}L_W(X).
\]
Comparing the linear coefficients and using \(c_0\ne0\) gives \(a^{q^r}=a\). Comparison of the remaining nonzero coefficients gives \(a^{q^i}=a\) whenever \(c_i\ne0\); conversely, these equalities imply the polynomial identity.
Together with \(a^{q^n}=a\), these fixed-point conditions define the common fixed field \(\F_{q^d}\) in \eqref{v70:eq-scalar-support}, by the intersection formula for fixed fields. Zero also satisfies all the conditions.
Since \(c_r=1\), the greatest common divisor includes \(r\), so \(d\mid r\).
\end{proof}

This subfield is determined by the scalar symmetries of the root set and should be distinguished from the field \(k(W)\) generated by \(W\).
Writing \(W/W=\{x/y:x,y\in W,\ y\ne0\}\), the ratio-generated field \(k(W/W)\) is a third object.
For example, if \(W=bk\), with \(b\ne0\), then \(H(W)=k(W/W)=k\), whereas \(k(W)=k(b)\) may be larger.

\subsection{Reciprocal Spans and Rational Function Relations}\label{v116:subsec-inverse-rank}
A subspace polynomial records additive closure; taking reciprocals generally changes the linear relations.
For a binary subspace \(U\) of \(E=\F_{2^n}\), define its \textbf{reciprocal rank} by
\begin{equation}\label{v116:def-inverse-rank}
I(U)=\operatorname{span}_{\F_2}\{u^{-1}:u\in U\setminus\{0\}\},\qquad
r(U)=\dim_{\F_2}I(U).
\end{equation}
For \(c\in E^\times\), we have \(I(cU)=c^{-1}I(U)\), so scaling preserves the rank. If \(\dim U=s\), the \(2^s-1\) distinct reciprocals lie in \(I(U)\), giving \(s\le r(U)\le\min\{n,2^s-1\}\).

\begin{proposition}[Loss of reciprocal rank under evaluation]\label{v117:prop-specialization-loss}
Let \(\mathcal U\subseteq\F_2[X]\) be a nonzero finite-dimensional subspace, and put
\[
M=\operatorname{lcm}\{u:u\in\mathcal U\setminus\{0\}\},\qquad D=\deg M.
\]
Let \(\alpha\) be algebraic over \(\F_2\), with \(f=m_\alpha\) and \(\gcd(f,M)=1\). Put
\[
W=\operatorname{span}_{\F_2}\{M/u:u\in\mathcal U\setminus\{0\}\},
\qquad U=\{u(\alpha):u\in\mathcal U\}.
\]
Then \(\dim U=\dim\mathcal U\), and
\begin{equation}\label{v117:eq-specialization-loss}
r(U)=\dim W-\dim\bigl(W\cap(f)\bigr),
\end{equation}
where \((f)\) is the ideal generated by \(f\) in \(\F_2[X]\). Equivalently, \(r(U)\) is the rank of the matrix whose columns are the coefficient vectors of the remainders of \(M/u\) modulo \(f\).
\end{proposition}
\begin{proof}
Every nonzero \(u\in\mathcal U\) divides \(M\), so \(u(\alpha)\ne0\); evaluation is therefore injective on \(\mathcal U\).
Multiplication by \(M(\alpha)\ne0\) identifies \(I(U)\) with the evaluation space \(W(\alpha)\). The kernel of evaluation \(W\to W(\alpha)\) is exactly \(W\cap(f)\), and rank--nullity gives the formula.
\end{proof}

In general, reduce the numerators modulo \(f\) before computing their rank. The following degree condition makes the evaluation kernel zero.
\begin{corollary}[Preserving all reciprocal relations within a degree bound]\label{v116:prop-rational-specialization}
Let \(\mathcal U\subseteq\F_2[X]\) be a nonzero finite-dimensional subspace, and put
\[
M=\operatorname{lcm}\{u:u\in\mathcal U\setminus\{0\}\},\qquad D=\deg M.
\]
If \(\alpha\) has degree greater than \(D\) over \(\F_2\), then \(U=\{u(\alpha):u\in\mathcal U\}\) has the same dimension as \(\mathcal U\), and
\begin{equation}\label{v116:eq-rational-specialization}
\begin{aligned}
r(U)&=\dim_{\F_2}\operatorname{span}\{1/u:u\in\mathcal U\setminus\{0\}\}\\
&=\rank_{\F_2}\bigl([X^i](M/u)\bigr)_{0\le i\le D,\ u\in\mathcal U\setminus\{0\}}.
\end{aligned}
\end{equation}
\end{corollary}
\begin{proof}
Here \(\gcd(f,M)=1\), and every polynomial in \(W\) has degree at most \(D<\deg f\), so \(W\cap(f)=0\). Apply Proposition~\ref{v117:prop-specialization-loss} and express \(W\) as a column space of coefficients. Multiplication by the common denominator \(M\) also identifies \(\operatorname{span}\{1/u\}\) with \(W\).
\end{proof}

\begin{example}[Equal dimensions and different reciprocal ranks]\label{v116:ex-inverse-models}
For an exponent set \(J\), put \(\mathcal U_J=\operatorname{span}_{\F_2}\{X^j:j\in J\}\), and let \(M_J\) be the monic least common multiple of its nonzero elements. The ranks below are first computed in \(\F_2(X)\). If \(\deg\alpha>\deg M_J\), they also give the reciprocal ranks of \(U_J=\langle\alpha^j:j\in J\rangle_{\F_2}\).
\begin{center}
\begin{tabular}{@{}ccc@{}}
\toprule
\(J\)&Dimension of the reciprocal span&\(\deg M_J\)\\
\midrule
\(\{0,1,2\}\)&6&6\\
\(\{0,1,3\}\)&7&10\\
\(\{0,1,2,3\}\)&11&14\\
\(\{0,1,3,7\}\)&15&42\\
\bottomrule
\end{tabular}
\end{center}

Write \(Q=X^2+X+1\), \(A=X^3+X+1\), and \(B=X^3+X^2+1\).
The first row is spanned by \(1,X^{-1},(X+1)^{-1},X^{-2},(X+1)^{-2},Q^{-1}\), which are independent by their pole orders. The remaining reciprocal satisfies
\[
\frac1X+\frac1{X+1}+\frac1{X(X+1)}=0.
\]
For the second row, poles at \(A,Q,(X+1)^2,X^3\) successively eliminate the coefficients associated with denominators \(A,X^3+1,X^3+X,X^3\). The remaining functions \(1,X^{-1},(X+1)^{-1}\) are independent.
The two least common multiples are \(X^2(X+1)^2Q\) and \(X^3(X+1)^2QA\).

All reciprocals in the third row belong to the following partial-fraction space: constants, principal parts of order at most three at \(0,1\), simple principal parts at \(Q\), and \(A^{-1},B^{-1}\).
Its dimension is \(1+6+2+1+1=11\), and these reciprocals span every component. The two directions at \(Q\) can be taken as \(Q^{-1},X/Q\), using
\[
\frac XQ=\frac1{XQ}+\frac1X+\frac1Q.
\]
The least common multiple is \(X^3(X+1)^3QAB\), of degree fourteen.

The fourth row adds eight degree-seven denominators to the seven in the second row. For each new denominator, the table marks an irreducible factor whose indicated multiplicity occurs in that denominator alone:
\begin{center}\small
\begin{tabular}{@{}ll@{}}
\toprule
Denominator&Factor power used for comparison\\
\midrule
\(X^7\)&\(X^7\)\\
\(X^7+1\)&\(B\)\\
\(X^7+X\)&\(Q^2\)\\
\(X^7+X+1\)&\(C_1=X^7+X+1\)\\
\(X^7+X^3\)&\((X+1)^4\)\\
\(X^7+X^3+1\)&\(C_2=X^7+X^3+1\)\\
\(X^7+X^3+X\)&\(A^2\)\\
\(X^7+X^3+X+1\)&\(D_0=X^4+X+1\)\\
\bottomrule
\end{tabular}
\end{center}
These divisibilities follow by multiplication. Both \(C_1,C_2\) take value one at \(0,1\); their remainders modulo \(Q,A,B\) are \((1,X,X)\) and \((X,X+1,X^2+1)\), respectively. Thus both are irreducible; \(D_0\) has no irreducible factor of degree one or two.
In any reciprocal relation, comparing the highest principal parts at the marked poles eliminates the eight new coefficients, and the second row then eliminates the remaining seven.
The least common multiple is \(X^7(X+1)^4Q^2A^2BC_1C_2D_0\), of degree forty-two.
\end{example}

Thus, for the same odd prime \(n>42\) and the same \(\alpha\in\F_{2^n}\setminus\F_2\), the four-dimensional spaces \(\langle1,\alpha,\alpha^2,\alpha^3\rangle\) and \(\langle1,\alpha,\alpha^3,\alpha^7\rangle\) have reciprocal ranks eleven and fifteen.
Chapter~\ref{chap:algorithms} will use a similar common denominator to test normality in reciprocal families of translates; here it retains every linear relation.

For an arbitrary four-dimensional \(U\), these relations form the space
\begin{equation}\label{v116:def-inverse-relations}
\begin{gathered}
\mathcal C(U)=\left\{(\varepsilon_u)_{u\in U\setminus\{0\}}\in\F_2^{15}:
\sum_{u\ne0}\varepsilon_u/u=0\right\},\\
\dim\mathcal C(U)=15-r(U).
\end{gathered}
\end{equation}
Does there exist an integer \(N\) such that \(r(U)\ge11\) for every odd prime \(n>N\) and every four-dimensional subspace \(U\subseteq\F_{2^n}\)? Equivalently, must \(\dim\mathcal C(U)\le4\) hold? The preceding example attains eleven when \(n>42\).
A four-dimensional space has fifteen three-dimensional subspaces. Do relations supported in those subspaces span \(\mathcal C(U)\)? If not, which further relations crossing these sections are needed?
Exercise~\ref{v116:ex-inverse-sum} connects one such relation with a coefficient of the subspace polynomial.

\section{Composition and Operator Polynomials}
Multiplication of linear maps means composition, rather than pointwise multiplication of their values.
For example, \((aX^q)\circ(bX^q)=ab^qX^{q^2}\), which generally differs from \((bX^q)\circ(aX^q)\).
Thus linearized polynomials with coefficients in \(E\) generally do not commute under composition.

\leavevmode\phantomsection\label{def:linearized-symbol}\ignorespaces When all coefficients lie in \(k\), define, for \(a(T)=\sum a_iT^i\in k[T]\),
\[
\mathcal L_a(X)=\sum a_iX^{q^i}.
\]
Then \(\mathcal L_{a+b}=\mathcal L_a+\mathcal L_b\) and \(\mathcal L_{ab}=\mathcal L_a\circ\mathcal L_b\).
The constant polynomial \(1\) corresponds to the identity map \(X\). On \(E\), we have \(\mathcal L_a=a(\sigma)\), where \(\sigma(x)=x^q\).

\begin{theorem}[Kernel dimension and the greatest common divisor]\label{thm:linearized-gcd}
For \(a\in k[T]\), let \(h=\gcd(a,T^n-1)\) be the monic greatest common divisor. Then
\[
\ker\bigl(a(\sigma):E\to E\bigr)=\ker h(\sigma),\qquad
\dim_k\ker a(\sigma)=\deg h.
\]
In particular, \(a(\sigma)\) is invertible if and only if \(\gcd(a,T^n-1)=1\).
\end{theorem}
\begin{proof}
Since \(T^n-1\ne0\), the monic greatest common divisor \(h\) is defined even when \(a=0\).
B\'ezout's identity gives \(h=ua+v(T^n-1)\). Evaluating at the operator \(\sigma\) and using \(\sigma^n=\id\) yields
\[
h(\sigma)=u(\sigma)a(\sigma).
\]
Thus \(\ker a(\sigma)\subseteq\ker h(\sigma)\).
Conversely, \(a=hc\) gives \(a(\sigma)=c(\sigma)h(\sigma)\), so the kernels agree.

Write \(T^n-1=hb\). The correspondence between ordinary multiplication and linearized composition gives
\[
X^{q^n}-X=\mathcal L_b\bigl(\mathcal L_h(X)\bigr).
\]
A linearized polynomial has no constant term, so \(\mathcal L_b(Y)=Y\,C(Y)\) for some \(C\in k[Y]\);
hence \(\mathcal L_h\) divides \(X^{q^n}-X\) as an ordinary polynomial.
The latter splits over \(E\) with no repeated roots, so all roots of \(\mathcal L_h\) also lie in \(E\) and are simple.
If \(d=\deg h\), then \(\deg_X\mathcal L_h=q^d\), so its root space contains exactly \(q^d\) elements and has dimension \(d\).
This includes \(h=1\): in that case \(\mathcal L_h=X\) and the kernel is \(\{0\}\).

Finally, \(a(\sigma)\) is invertible if and only if its kernel is zero, equivalently \(\deg h=0\), or \(h=1\).
If \(ua+v(T^n-1)=1\), evaluating the same identity at the operator also gives the explicit inverse \(u(\sigma)\).
\end{proof}

For example, when both \(q\) and \(n\) are odd, \(x\mapsto x+x^q\) is invertible, with inverse
\[
x\longmapsto\frac12\bigl(x-x^q+x^{q^2}-\cdots+x^{q^{n-1}}\bigr).
\]
This is the identity \((1+T)(1-T+\cdots+T^{n-1})=1+T^n\) evaluated at \(\sigma\).

\section{The Operator Algebra Generated by Multiplication and Frobenius}
\label{sec:operator-algebra}

\begin{theorem}\label{thm:operator-algebra}
For \(a\in E\), write \(M_a(x)=ax\). Every \(k\)-linear operator \(A:E\to E\) has a unique expression
\[
A=\sum_{i=0}^{n-1}M_{a_i}\sigma^i,\qquad a_i\in E.
\]
Multiplication of these operators is determined by the relations
\[
\sigma M_a=M_{a^q}\sigma,\qquad \sigma^n=\id.
\]
In \(\End_k(E)\) we have
\[
\begin{aligned}
\{A:A\sigma=\sigma A\}=k[\sigma]\simeq k[T]/(T^n-1),
\\
\{A:AM_a=M_aA\text{ for all }a\in E\}=\{M_c:c\in E\}.
\end{aligned}
\]
\end{theorem}
\begin{proof}
The first expression follows from the linearized polynomial representation in Theorem~\ref{thm:linear-map-representation}; the relations follow from \((ax)^q=a^qx^q\). Comparing
\[
A\sigma=\sum_i M_{a_i}\sigma^{i+1},\qquad
\sigma A=\sum_i M_{a_i^q}\sigma^{i+1},
\]
and reducing exponents modulo \(n\), uniqueness of the representation shows that the two agree if and only if \(a_i\in k\) for every index.
The kernel of \(k[T]\to k[\sigma]\) is \((T^n-1)\): first divide by \(T^n-1\) to remove the higher-degree terms.
If a remainder of degree less than \(n\) gives the zero operator, the same uniqueness shows that all its coefficients vanish.

If \(A\) commutes with every \(M_a\), then for any \(x\in E\),
\(A(x)=A(M_x(1))=M_x(A(1))=xA(1)\), so \(A=M_{A(1)}\). The reverse inclusion follows immediately from commutativity of field multiplication.
\end{proof}

Allowing coefficients in \(E\) gives all linear operators; restricting them to \(k\)
gives precisely the operators commuting with Frobenius. The coefficient field thus determines a commutation property.
The operators commuting with both Frobenius and all multiplication operators are precisely multiplication by scalars in \(k\).

\begin{example}[Commuting with Frobenius need not preserve multiplication]\label{exm:linear-not-field}
Let \(E=\F_{16}\) and \(\sigma(x)=x^2\), and take \(A=1+\sigma+\sigma^2\).
Since \(T^4-1=(T+1)^4\) and \((1+T+T^2)|_{T=1}=1\),
the polynomials \(1+T+T^2\) and \(T^4-1\) are coprime. Theorem~\ref{thm:linearized-gcd} therefore shows that \(A\) is invertible.
It commutes with \(\sigma\) and satisfies \(A(1)=1\).
If \(A\) also preserved multiplication, it would be an automorphism of \(E/\F_2\), hence some \(\sigma^j\) with \(0\le j<4\).
The nonzero polynomial \(1+T+T^2-T^j\) would then annihilate \(\sigma\), although its degree is at most three,
contradicting the fact that the kernel of \(k[T]\to k[\sigma]\) in Theorem~\ref{thm:operator-algebra} is \((T^4-1)\).
Thus an invertible linear map fixing the identity and commuting with Frobenius need not be a field isomorphism.
\end{example}

\section{Root Multiplicities and Composition Division}
\label{sec:linearized-division}
If \(L(X)=\sum_{i=h}^s a_iX^{q^i}\) and \(a_ha_s\ne0\), take the unique \(q^h\)th root of each coefficient in the finite coefficient field to obtain
\[
L(X)=M(X)^{q^h},\qquad M'(X)\ne0.
\]
In a splitting extension, every root has multiplicity exactly \(q^h\), and the distinct roots form an \(\F_q\)-space of dimension \(s-h\).
Conversely, if a nonzero polynomial splits completely in a specified extension, its distinct roots form a finite \(\F_q\)-space \(V\), and every root has the same multiplicity \(q^h\), then the polynomial must be \(cL_V(X)^{q^h}\) and hence is linearized.
The condition that all roots have the same multiplicity, a power of \(q\), is essential; the condition that the root set be a subspace is not sufficient by itself.

\begin{theorem}[Right composition division]\label{thm:composition-division}
Let \(L,M\) be \(q\)-linearized polynomials over the same finite extension, with \(M\ne0\). There are unique \(q\)-linearized polynomials \(Q,R\) over that field such that
\[
L=Q\circ M+R,\qquad R=0\ \text{or}\ \deg_q R<\deg_q M.
\]
Moreover, divisibility \(M\mid L\) as ordinary polynomials holds if and only if \(L=Q\circ M\).
\end{theorem}
\begin{proof}
Let the leading terms be \(aX^{q^s}\) and \(bX^{q^t}\), with \(s\ge t\).
Subtracting \((a/b^{q^{s-t}})M(X)^{q^{s-t}}\) cancels the leading term. Repeating finitely many times proves existence.
The \(q\)-degree of a nonzero composition is the sum of the \(q\)-degrees, so subtracting two division identities proves uniqueness.
Every term of \(Q(M(X))\) is divisible by the ordinary polynomial \(M(X)\). Thus if \(M\mid L\), then \(M\mid R\);
the remainder has smaller ordinary degree and must be zero. Conversely, if \(L=Q\circ M\), the same termwise divisibility gives \(M\mid L\).
\end{proof}
Older literature calls composition ``symbolic multiplication.'' For coefficients in \(\F_q\),
\(\mathcal L_a\circ\mathcal L_b=\mathcal L_{ab}\), so symbolic factorization is ordinary factorization of operator polynomials;
over a larger coefficient field, composition is generally noncommutative, and its order must be specified.

To prescribe several right remainders, first identify the data recorded by each remainder.
When the divisor is a subspace polynomial \(L_W\), the identity \(L=Q\circ L_W+R\) shows that \(R\) and \(L\) have the same restriction to \(W\).
Compatibility of several remainders therefore depends on the linear relations among the corresponding subspaces.
This gives a direct form of the linearized Chinese remainder problem for subspace polynomials; more general reconstruction from right remainders is treated in~\cite{GaboritGarnierRuatta2026CRT}.

\begin{proposition}[Compatibility of linearized remainders]\label{v105:prop-linearized-remainder-compatibility}
Let \(E=\F_{q^m}\), and let \(W_1,\ldots,W_s\subseteq E\) be nonzero \(\F_q\)-subspaces. Write
\[
 r_i=\dim_{\F_q}W_i,\qquad W=W_1+\cdots+W_s,
 \qquad r=\dim_{\F_q}W.
\]
For a positive integer \(d\), put
\[
 \mathscr L_{<d}(E)=
 \left\{\sum_{j=0}^{d-1}a_jX^{q^j}:a_j\in E\right\}.
\]
Denote by \(\operatorname{rem}_iL\) the \textbf{right composition remainder} of \(L\) upon division by \(L_{W_i}\), so that
\[
 L=Q_i\circ L_{W_i}+\operatorname{rem}_iL,
 \qquad \operatorname{rem}_iL\in\mathscr L_{<r_i}(E).
\]
Then the map
\begin{equation}\label{v105:eq-linearized-remainder-map}
 \begin{aligned}
 \Phi:\mathscr L_{<r}(E)&\longrightarrow
 \bigoplus_{i=1}^s\mathscr L_{<r_i}(E),\\
 L&\longmapsto(\operatorname{rem}_iL)_{i=1}^s.
 \end{aligned}
\end{equation}
is \(E\)-linear and injective, with
\begin{equation}\label{v105:eq-linearized-remainder-codimension}
 \operatorname{codim}_E\operatorname{im}\Phi
 =\sum_{i=1}^s r_i-r.
\end{equation}
Given \(R_i\in\mathscr L_{<r_i}(E)\), a polynomial \(L\in\mathscr L_{<r}(E)\) satisfying \(\operatorname{rem}_iL=R_i\) for every \(i\) exists if and only if
\begin{equation}\label{v105:eq-linearized-remainder-compatibility}
 w_i\in W_i,\quad \sum_{i=1}^s w_i=0
 \quad\Longrightarrow\quad
 \sum_{i=1}^s R_i(w_i)=0.
\end{equation}
When it exists, \(L\) is unique. In particular, all the right remainders can be prescribed independently if and only if
\[
 W=W_1\oplus\cdots\oplus W_s.
\]
\end{proposition}
\begin{proof}
For \(a\in E\), the identity
\[
 aL=(aQ_i)\circ L_{W_i}+a\operatorname{rem}_iL
\]
and uniqueness of the right remainder give \(\operatorname{rem}_i(aL)=a\operatorname{rem}_iL\). Addition is treated in the same way.
Thus \(\Phi\) is linear for the \(E\)-scalar operations on polynomial coefficients.

If \(\Phi(L)=0\), then \(L\) vanishes on every \(W_i\).
Since \(L\) is \(\F_q\)-linear, it vanishes on their sum \(W\).
A nonzero polynomial in \(\mathscr L_{<r}(E)\) has ordinary degree at most \(q^{r-1}\), so it cannot have \(|W|=q^r\) roots. Hence \(L=0\), proving injectivity.
Now \(\dim_E\mathscr L_{<d}(E)=d\), and dimension counting gives~\eqref{v105:eq-linearized-remainder-codimension}.

If \(L\) has the prescribed remainders, then \(L|_{W_i}=R_i|_{W_i}\), and additivity gives~\eqref{v105:eq-linearized-remainder-compatibility}.
Conversely, suppose that this condition holds. Define
\[
 f:W\longrightarrow E,\qquad
 f(w_1+\cdots+w_s)=R_1(w_1)+\cdots+R_s(w_s).
\]
Subtracting two representations of the same element gives a relation with sum zero, so compatibility makes \(f\) well defined.
The \(\F_q\)-linearity of the \(R_i\) makes \(f\) an \(\F_q\)-linear map.
Choose a basis \(v_1,\ldots,v_r\) of \(W\). By the Moore determinant criterion, the system
\[
 \sum_{j=0}^{r-1}a_jv_\ell^{q^j}=f(v_\ell),
 \qquad 1\le\ell\le r,
\]
has a unique solution \((a_0,\ldots,a_{r-1})\in E^r\).
The resulting polynomial \(L(X)=\sum_{j=0}^{r-1}a_jX^{q^j}\) agrees with \(f\) on all of \(W\).
Consequently, \(\operatorname{rem}_iL-R_i\) vanishes on \(W_i\). It belongs to \(\mathscr L_{<r_i}(E)\), so the root bound forces it to be zero.
This proves existence, and injectivity of \(\Phi\) gives uniqueness.
Finally, \(\Phi\) is surjective if and only if \(r=\sum_i r_i\), which is precisely the direct-sum condition.
\end{proof}

The \(E\)-linearity here concerns the coefficient spaces.
Under composition, linearized polynomials naturally form an \(\F_q\)-algebra, generally not an \(E\)-algebra; the statement is therefore not a direct-product decomposition of quotient \(E\)-algebras.
The moduli in the proposition are also split, separable subspace polynomials, whose \(q\)-degrees equal their root-space dimensions.
For a linearized modulus with repeated roots, its distinct root space alone does not give this criterion.

For a one-dimensional space \(W=\F_qg\), where \(g\ne0\), we have
\[
 L_W(X)=X^q-g^{q-1}X.
\]
The right remainder has the form \(cX\). Evaluating at \(g\) gives \(cg=L(g)\), hence
\begin{equation}\label{v105:eq-linearized-point-remainder}
 \operatorname{rem}_{L_W}L=\frac{L(g)}{g}X.
\end{equation}
Thus a one-dimensional right remainder records a normalized evaluation on a linear direction.

\begin{example}[Three lines and simultaneous compatibility]\label{v105:ex-linearized-crt-three-lines}
Let \(E=\F_4=\F_2(\omega)\), where \(\omega^2+\omega+1=0\), and take
\[
 W_1=\F_2,\qquad W_2=\omega\F_2,
 \qquad W_3=\omega^2\F_2.
\]
The lines meet pairwise only at zero. Their subspace polynomials are
\[
 P_1=X^2+X,\qquad P_2=X^2+\omega X,
 \qquad P_3=X^2+\omega^2X.
\]
As ordinary polynomials, their pairwise monic greatest common divisors are all \(X\).
By Theorem~\ref{thm:composition-division}, their pairwise monic greatest common right factors under composition are also \(X\).
Since \(X\) is the identity for composition, these moduli are \textbf{pairwise right coprime}.

Nevertheless, their three right remainders cannot be prescribed independently.
If \(\operatorname{rem}_1L=\operatorname{rem}_2L=0\), then \(L(1)=L(\omega)=0\), and therefore
\[
 L(\omega^2)=L(1+\omega)=0.
\]
We cannot also require \(\operatorname{rem}_3L=X\), since this would give \(L(\omega^2)=\omega^2\ne0\).
More precisely, writing \(R_i=c_iX\), the complete compatibility condition is
\[
 c_1+\omega c_2+\omega^2c_3=0.
\]
Indeed, \(1,\omega\) form a basis of \(E/\F_2\), and the third value must satisfy exactly the relation given by their sum.
Here \(r_1+r_2+r_3=3\) and \(r=2\), so this single \(E\)-linear relation also agrees with~\eqref{v105:eq-linearized-remainder-codimension}.
\end{example}

\begin{proposition}[Minimal linearized polynomial]\leavevmode\phantomsection\label{prop:minimal-linearized}
Given a finite field \(K\supseteq\F_q\) and an element \(\beta\) algebraic over \(K\), there is a unique monic nonzero \(q\)-linearized polynomial \(M\in K[X]\) of least \(q\)-degree such that \(M(\beta)=0\).
It is called the \textbf{minimal linearized polynomial} \(M\) of \(\beta\) over \(K\).
The \(q\)-linearized polynomials with coefficients in \(K\) vanishing at \(\beta\) are exactly those of the form \(Q\circ M\), where \(Q\in K[X]\) is also \(q\)-linearized; in this case \(M\) is called a right factor under composition.
\end{proposition}
\begin{proof}
The element \(\beta\) lies in a finite field, so some \(X^{q^N}-X\) vanishes at it. Choosing the least degree proves existence.
Right composition division of any other annihilating polynomial leaves a remainder vanishing at \(\beta\), which must be zero by minimality.
Subtracting two monic polynomials of least degree proves uniqueness.
\end{proof}

\begin{remark}[Roots of the minimal linearized polynomial]
Let \(|K|=Q\), and take all \(K\)-conjugates \(\beta,\beta^Q,\ldots\) of \(\beta\) in a finite extension.
Write \(V\) for their \(\F_q\)-linear span.
The \(Q\)th-power map preserves \(V\), so the coefficients of \(L_V\) are fixed by the \(Q\)th-power map and belong to \(K\).
Every nonzero linearized polynomial with coefficients in \(K\) vanishing at \(\beta\)
also vanishes at all its \(K\)-conjugates and hence, by \(\F_q\)-linearity, on \(V\).
Its ordinary degree is therefore at least \(|V|\).
The polynomial \(L_V\) attains this degree, so the minimal linearized polynomial is exactly \(L_V\).
The ordinary minimal polynomial has the conjugate orbit as its root set; the minimal linearized polynomial has the linear span of that orbit as its root set.
\end{remark}

When the coefficient field is exactly \(k=\F_q\), this \(M\) corresponds to the Frobenius annihilator polynomial \(\mu_\beta\) defined in the next chapter, namely \(M=\mathcal L_{\mu_\beta}\).
Its root set is the linear span of the conjugate orbit, whereas the root set of the ordinary minimal polynomial is the conjugate orbit itself.

\begin{proposition}[Invariant subspaces and operator factors]\label{prop:invariant-root-subspaces}
Let \(h\in k[T]\) be monic, with \(\deg h\ge1\) and \(h(0)\ne0\), and let \(W\) be the space of all roots of \(\mathcal L_h\).
The Frobenius-invariant \(k\)-subspaces of \(W\) correspond bijectively to the monic divisors of \(h\):
\[
d\mid h\quad\longleftrightarrow\quad\{x:\mathcal L_d(x)=0\}.
\]
\end{proposition}
\begin{proof}
The monic root polynomial \(L_V\) of an invariant subspace \(V\) has coefficients in \(k\), so it equals some \(\mathcal L_d\).
Since \(V\subseteq W\) and the roots are simple, we have ordinary divisibility \(\mathcal L_d\mid\mathcal L_h\).
Right division under composition over \(k\) gives a linearized quotient with coefficients in \(k\), say \(\mathcal L_c\).
Thus \(\mathcal L_h=\mathcal L_c\circ\mathcal L_d=\mathcal L_{cd}\); comparing the coefficients of each \(X^{q^i}\) gives \(h=cd\).
The same composition identity proves the converse, and uniqueness of the root-set polynomial establishes the bijection.
\end{proof}
In particular, \(W\) has no nonzero proper invariant subspace if and only if \(h\) is irreducible in \(k[T]\).

\begin{proposition}[A root-product formula for the Moore determinant]\label{prop:moore-product}
For arbitrary \(x_1,\ldots,x_s\), we have
\[
\det(x_j^{q^i})_{0\le i<s,\,1\le j\le s}
=x_1\prod_{j=1}^{s-1}\prod_{c_1,\ldots,c_j\in k}
\left(x_{j+1}-\sum_{i=1}^jc_ix_i\right).
\]
\end{proposition}
\begin{proof}
Induct on \(s\). If the first \(s-1\) elements are linearly dependent, the left side has dependent columns and vanishes; the product on the right formed from the first \(s-1\) elements also vanishes by induction.
Otherwise, regard the left side as a polynomial in the last variable \(X=x_s\). Its leading coefficient is the preceding Moore determinant and its degree is \(q^{s-1}\).
Every element in the span of the first \(s-1\) elements makes this determinant vanish, giving exactly \(q^{s-1}\) roots.
Compare degrees and leading coefficients, then apply the induction hypothesis.
\end{proof}

This product formula also computes the subspace polynomial from a basis. Let \(V=\langle v_1,\ldots,v_s\rangle_k\), where \(v_1,\ldots,v_s\) are linearly independent and \(s\ge1\). Then
\begin{equation}\label{eq:subspace-moore-determinant}
L_V(X)=
\frac{\displaystyle\det\begin{pmatrix}
v_1&\cdots&v_s&X\\
v_1^q&\cdots&v_s^q&X^q\\
\vdots&&\vdots&\vdots\\
v_1^{q^s}&\cdots&v_s^{q^s}&X^{q^s}
\end{pmatrix}}
{\displaystyle\det(v_j^{q^i})_{0\le i<s,\,1\le j\le s}}.
\end{equation}
The Moore criterion makes the denominator nonzero; it is precisely the leading coefficient of the numerator, and the root-product formula identifies the quotient as \(\prod_{v\in V}(X-v)\).
Expansion along the last column also gives the individual linearized coefficients of \(L_V\).
Thus a subspace may be specified by a basis or by its monic linearized root polynomial; formula~\eqref{eq:subspace-moore-determinant} converts the first description directly into the second.

Here the Moore determinant both detects linear independence over the base field
and constructs the subspace polynomial from a basis. Chapter~\ref{chap:normal} arranges the
\emph{Frobenius conjugates of a single element} in a Moore matrix; a column
permutation converts it into a multiplication matrix in a cyclic algebra
(Proposition~\ref{v88:prop-moore-circulant}). Thus the trace-power determinant
formula in characteristic-power degree is an exact expression of the same
independence criterion in local cyclic-algebra coordinates, not a new
unrelated test.

\Needspace{7\baselineskip}
Furthermore, if \(a,b\in k[T]\) are not both zero, then
\begin{equation}\label{eq:linearized-euclidean-gcd}
\gcd(\mathcal L_a,\mathcal L_b)=\mathcal L_{\gcd(a,b)},
\end{equation}
where all greatest common divisors are monic. This identity follows directly from the Euclidean algorithm.
If \(b\ne0\), division with remainder in \(k[T]\) gives \(a=cb+r\), and then
\[
\mathcal L_a=\mathcal L_c\circ\mathcal L_b+\mathcal L_r.
\]
Since \(\mathcal L_c(Y)\) has no constant term, the ordinary polynomial \(\mathcal L_b(X)\) divides \(\mathcal L_c(\mathcal L_b(X))\), so
\[
\gcd(\mathcal L_a,\mathcal L_b)=\gcd(\mathcal L_b,\mathcal L_r).
\]
Following the Euclidean algorithm for \(a,b\) reduces the left side to the linearization of the last nonzero remainder.
Multiplication by a nonzero base-field constant commutes with linearization, so monic normalization gives~\eqref{eq:linearized-euclidean-gcd}.
If \(b=0\), the identity follows directly by monic normalization. The argument takes place in the ordinary polynomial ring and thus preserves both roots and their multiplicities;
it also shows that, for base-field coefficients, the ordinary greatest common divisor agrees with the greatest common right divisor under composition.

\section{Adjoints and Solvability of Linear Equations}
\begin{proposition}\label{prop:linearized-adjoint}
Let \(L(x)=\sum_{i=0}^{n-1}a_ix^{q^i}\). Its adjoint with respect to the trace pairing is
\[
L^*(x)=a_0x+\sum_{i=1}^{n-1}a_i^{q^{n-i}}x^{q^{n-i}},
\]
that is, \(\Tr(L(x)y)=\Tr(xL^*(y))\). Consequently,
\[
\im L=(\ker L^*)^\perp.
\]
The equation \(L(x)=b\) is solvable if and only if \(\Tr(bz)=0\) for every \(z\in\ker L^*\);
when solvable, its solutions form a coset of \(\ker L\).
\end{proposition}
\begin{proof}
The trace is Frobenius-invariant, so \(\Tr(a_i x^{q^i}y)=\Tr(x(a_i y)^{q^{n-i}})\). Summing termwise gives the adjoint formula.
An element \(z\) is orthogonal to \(\im L\) if and only if \(\Tr(xL^*(z))=0\) for every \(x\);
by nondegeneracy, this is equivalent to \(L^*(z)=0\). Taking orthogonal complements gives the image formula. The description of the solution set follows from linearity.
\end{proof}

In particular, \(\sigma^*=\sigma^{-1}\).
If \(c\in E\), the linearized polynomial \(X^q-cX\) permutes \(E\) if and only if \(\Nm_{E/k}(c)\ne1\).
For \(c=0\) this is invertibility of Frobenius; for \(c\ne0\), a nonzero kernel element exists if and only if \(x^{q-1}=c\) is solvable, and multiplicative Hilbert 90 applies.

When \(\Nm_{E/k}(c)=1\), choose \(u\in E^\times\) with \(c=u^{q-1}\). Then for \(L(x)=x^q-cx\),
\[
\ker L=uk,\qquad \ker L^*=u^{-q}k.
\]
The first identity follows from \(L(uy)=u^q(y^q-y)\). For the second, substitute \(u^{-q}\) into the adjoint formula, then use equality of the ranks of an operator and its adjoint to determine the kernel dimension.
The general adjoint criterion therefore becomes
\[
 x^q-cx=b\text{ is solvable in }E
 \quad\Longleftrightarrow\quad
 \Tr_{E/k}(b/u^q)=0.
\]
The substitution \(x=uy\) also reduces this directly to additive Hilbert 90; when solvable, all solutions are \(x_0+uk\).
Thus the trace condition in Hilbert 90 is the orthogonality condition determined by the adjoint kernel, and \(u\) gives both the direction of the solution fiber and the weight in the trace condition.

The adjoint also determines the restriction of the trace pairing to an operator image.

\begin{corollary}[The trace pairing on an operator image]\label{v68:cor-image-trace-radical}
Let \(L:E\to E\) be \(k\)-linear, put \(W=\im L\), and let \(B(x,y)=\Tr_{E/k}(xy)\).
The radical of the restricted pairing consists of the elements \(w\in W\) satisfying \(B(w,z)=0\) for every \(z\in W\). It satisfies
\[
\rad(B|_W)=\im L\cap\ker L^*,
\]
and
\[
\dim_k\rad(B|_W)=\rank_k L-\rank_k(L^*\circ L).
\]
\end{corollary}
\begin{proof}
Proposition~\ref{prop:linearized-adjoint} gives \(W^\perp=\ker L^*\), proving the first equality.
The map \(L\) sends \(\ker(L^*\circ L)\) onto \(\im L\cap\ker L^*\), and the kernel of this restriction is \(\ker L\).
Thus the radical has dimension \(\dim_k\ker(L^*\circ L)-\dim_k\ker L\), which is the stated rank difference.
\end{proof}

\begin{example}[The parity of the extension degree]\label{v68:ex-trace-image-parity}
Let \(E=\F_{2^n}\), \(k=\F_2\), and \(L(x)=x^2+x\). Additive Hilbert 90 gives
\[
W=\im L=\ker\Tr_{E/\F_2},\qquad \dim_{\F_2}W=n-1.
\]
Since \(L^*=\sigma^{-1}+I\), its kernel is \(\F_2\). Hence
\[
\rad(B|_W)=W\cap\F_2
=\begin{cases}
\{0\},&n\text{ odd},\\
\F_2,&n\text{ even},
\end{cases}
\]
because \(\Tr_{E/\F_2}(1)=n\cdot1\).
For \(E=\F_8\), the pairing on the two-dimensional image is nondegenerate; for \(E=\F_{16}\), the three-dimensional image has radical \(\F_2\).
The trace pairing on the whole field is nondegenerate in both cases.
\end{example}

\subsection*{The Trace Pairing Identifies the Dual Map with the Adjoint}
Let \(E^\vee=\operatorname{Hom}_k(E,k)\). The trace pairing gives an isomorphism
\[
\kappa:E\longrightarrow E^\vee,\qquad
\kappa(y)(x)=\Tr_{E/k}(xy).
\]
For a \(k\)-linear map \(A:E\to E\), define its dual map by \(A^\vee(\ell)=\ell\circ A\). Figure~\ref{fig:trace-adjoint} expresses
\begin{equation}\label{eq:adjoint-dual-conjugation}
\kappa A^*=A^\vee\kappa,
\qquad A^*=\kappa^{-1}A^\vee\kappa.
\end{equation}
Indeed, for every \(x,y\in E\), the two paths evaluated at \(x\) give \(\Tr(xA^*y)\) and \(\Tr((Ax)y)\), which agree by the definition of the adjoint.

\begin{figure}[!htbp]
\centering
\begin{tikzpicture}[every node/.style={font=\normalsize},
  map/.style={-{Stealth[length=2.2mm]},Jade,semithick}]
\node (E1) at (0,2) {\(E\)};
\node (E2) at (5,2) {\(E\)};
\node (D1) at (0,0) {\(E^\vee\)};
\node (D2) at (5,0) {\(E^\vee\)};
\draw[map] (E1)--(E2) node[midway,above] {\(A^*\)};
\draw[map] (D1)--(D2) node[midway,below] {\(A^\vee:\ell\mapsto\ell\circ A\)};
\draw[map] (E1)--(D1) node[midway,left] {\(\kappa\)} node[midway,right] {\(\simeq\)};
\draw[map] (E2)--(D2) node[midway,right] {\(\kappa\)} node[midway,left] {\(\simeq\)};
\end{tikzpicture}
\caption{The commutative diagram for the trace pairing. Equality of the functionals obtained along the two paths is equivalent to \(\Tr((Ax)y)=\Tr(xA^*y)\).}
\label{fig:trace-adjoint}
\end{figure}

Thus \(\sigma^*=\sigma^{-1}\) corresponds to preservation of the trace pairing by Frobenius:
\(\Tr(\sigma x\,\sigma y)=\Tr(xy)\). For multiplication operators, \(M_a^*=M_a\).
These two adjoint identities determine the adjoint of every linearized operator and yield \(\im A=(\ker A^*)^\perp\). Section~\ref{sec:fourier-adjoint} uses the same commutation relation to describe the transformation of Fourier frequencies.

\subsection{Quotient Multiplication and Trace Duality}

\begin{definition}\label{v70:def-subspace-product}
For \(k\)-subspaces \(U,V\subseteq E\), the \textbf{product of subspaces} is defined by
\[
 UV=\operatorname{span}_k\{uv:u\in U,\ v\in V\}.
\]
The linear span is part of the definition; \(UV\) need not equal the set of individual products.
\end{definition}

\begin{definition}\label{v70:def-quotient-multiplication}
Given a \(k\)-subspace \(W\subseteq E\), define the \textbf{quotient multiplication map} for each \(a\in E\) by
\[
 T_{a,W}:W\longrightarrow E/W,\qquad x\longmapsto ax+W.
\]
These maps depend \(k\)-linearly on \(a\); the resulting parameter map is
\[
 \Phi_W:E\longrightarrow\operatorname{Hom}_k(W,E/W),\qquad
 a\longmapsto T_{a,W}.
\]
Here \(E/W\) is a quotient vector space, generally not a quotient ring.
For \(a\in E^\times\), the number
\[
 I_W(a)=\dim_k(W\cap aW)
\]
is called the \textbf{multiplicative intersection dimension} of \(W\) with its scalar multiple.
\end{definition}

\begin{proposition}[Kernels and duals of quotient multiplication]\label{v70:prop-quotient-multiplication}
Put \(r=\dim_kW\), and take orthogonal complements with respect to \(\Tr_{E/k}(xy)\). Then
\[
 \ker\Phi_W=H(W).
\]
For every \(a\in E^\times\),
\[
 \ker T_{a,W}=W\cap a^{-1}W,\qquad
 \rank T_{a,W}=r-I_W(a).
\]
The trace pairing gives natural isomorphisms
\[
 (E/W)^\vee\simeq W^\perp,\qquad W^\vee\simeq E/W^\perp.
\]
Under these identifications, the dual of \(T_{a,W}\) is \(T_{a,W^\perp}\). Consequently,
\begin{equation}\label{v70:eq-intersection-duality}
 H(W)=H(W^\perp),\qquad
 r-I_W(a)=(n-r)-I_{W^\perp}(a)\quad(a\ne0).
\end{equation}
Moreover, for all \(k\)-subspaces \(U,V\subseteq E\),
\begin{equation}\label{v70:eq-subspace-product-orthogonal}
 (UV)^\perp=\{z\in E:zU\subseteq V^\perp\},
\end{equation}
and in particular,
\begin{equation}\label{v70:eq-product-relative-trace}
 WW^\perp=H(W)^\perp=\ker\Tr_{E/H(W)}.
\end{equation}
\end{proposition}
\begin{proof}
The equality \(T_{a,W}=0\) is equivalent to \(aW\subseteq W\), so the kernel of the parameter map \(\Phi_W\) is \(H(W)\).
For \(a\ne0\), the condition \(T_{a,W}(x)=0\) is equivalent to \(x\in W\cap a^{-1}W\). Multiplication by \(a\) maps this intersection isomorphically onto \(aW\cap W\), proving the rank formula.
The kernel of the parameter map consists of scalars for which the entire map vanishes, whereas the kernel of one \(T_{a,W}\) consists of vectors whose products lie in \(W\).

The two quotient pairings are
\[
 (z+W,y)\longmapsto\Tr_{E/k}(zy),\qquad
 (x,z+W^\perp)\longmapsto\Tr_{E/k}(xz).
\]
In the first, \(y\in W^\perp\); in the second, \(x\in W\). Thus both are independent of the coset representatives.
Nondegeneracy of the trace pairing on \(E\) gives \((W^\perp)^\perp=W\) and \(\dim W^\perp=n-r\). It follows that both quotient pairings are nondegenerate on each side and give the stated duality isomorphisms.
No nondegeneracy assumption on the restriction to \(W\) is needed.
For \(x\in W\) and \(y\in W^\perp\), the identity
\[
 \Tr_{E/k}((ax)y)=\Tr_{E/k}(x(ay))
\]
identifies the dual map with \(y\mapsto ay+W^\perp\). A map and its dual vanish simultaneously and have the same rank, giving \eqref{v70:eq-intersection-duality}.

For every \(z\in E\),
\[
\begin{aligned}
 z\in(UV)^\perp
 &\quad\Longleftrightarrow\quad
 \Tr_{E/k}(zuv)=0\quad\text{for all }u\in U,\ v\in V\\
 &\quad\Longleftrightarrow\quad zU\subseteq V^\perp.
\end{aligned}
\]
Taking \(U=W\) and \(V=W^\perp\) gives \((WW^\perp)^\perp=H(W)\); taking orthogonal complements proves the first equality.
Put \(H=H(W)\). Transitivity of trace gives
\[
 \Tr_{E/k}(zh)=\Tr_{H/k}\bigl(h\Tr_{E/H}(z)\bigr)\qquad(h\in H).
\]
By nondegeneracy of the trace pairing for \(H/k\), this vanishes for every \(h\) if and only if \(\Tr_{E/H}(z)=0\), proving the second equality.
The argument also covers \(W=0\) and \(W=E\), and does not require the characteristic to be coprime to the extension degree.
\end{proof}

Equation~\eqref{v70:eq-subspace-product-orthogonal} also gives a simple product criterion: if \(\dim_kU+\dim_kV>n\), then \(UV=E\).
Otherwise choose \(0\ne z\in(UV)^\perp\). Then \(zU\subseteq V^\perp\), contradicting \(\dim_kU>n-\dim_kV\).

\Needspace{10\baselineskip}
\subsection{Counting Scalar Intersections}

\begin{proposition}\label{v70:prop-intersection-count}
Let \(0\ne W\subseteq E\), put \(r=\dim_kW\), and write \(H(W)=\F_{q^d}\). For \(0\le j\le r\), set
\[
 N_j=\#\{a\in E^\times:I_W(a)=j\}.
\]
Then
\begin{equation}\label{v70:eq-intersection-moments}
\begin{aligned}
 \sum_{j=0}^rN_j&=q^n-1,&
 \sum_{j=0}^r(q^j-1)N_j&=(q^r-1)^2,\\
 N_r&=q^d-1,&
 N_j&=0\quad\text{if }d\nmid j.
\end{aligned}
\end{equation}
For every \(a\ne0\), we have \(I_W(a^{-1})=I_W(a)\). If additionally \(\sigma(W)=W\), then \(I_W(a^q)=I_W(a)\).
Each \(N_j\) is divisible by \(q^d-1\); the quotient \(N_j/(q^d-1)\) counts the distinct scalar multiples \(aW\) satisfying \(\dim_k(W\cap aW)=j\).
\end{proposition}
\begin{proof}
The first identity partitions \(E^\times\) by intersection dimension.
For a fixed \(a\ne0\), the pairs \((x,y)\in(W\setminus\{0\})^2\) satisfying \(x=ay\) are in bijection with \((W\cap aW)\setminus\{0\}\), and hence number \(q^{I_W(a)}-1\).
Conversely, every nonzero ordered pair has a unique ratio \(a=x/y\); summing over \(a\) proves the second identity.
The intersection has dimension \(r\) precisely when \(aW=W\), or equivalently \(a\in H(W)^\times\), giving \(N_r=q^d-1\).
The intersection \(W\cap aW\) is still an \(H(W)\)-space, so its \(k\)-dimension is divisible by \(d\).

Multiplication by \(a\) gives \(W\cap a^{-1}W\simeq aW\cap W\), proving invariance under inversion.
When \(\sigma(W)=W\), Frobenius maps \(W\cap aW\) isomorphically onto \(W\cap a^qW\), proving Frobenius invariance.
Finally, \(aW=bW\) if and only if \(b^{-1}a\in H(W)^\times\). Each scalar multiple is therefore represented by exactly \(q^d-1\) parameters, all having the same intersection dimension.
\end{proof}

\begin{proposition}[Scalar intersections of power spaces and fractional linear transformations]\label{v150:prop-power-intersection}
Let $E=k(\alpha)=\F_{q^t}$, $k=\F_q$, and $r\ge2$, with $t>2r-2$. Put
\[
 W=P_r(\alpha)=\operatorname{span}_k\{1,\alpha,\ldots,\alpha^{r-1}\}.
\]
For nonzero coprime $f,g\in k[X]$ of degrees at most $r-1$,
\begin{equation}\label{v150:eq-power-intersection}
 I_W\left(\frac{f(\alpha)}{g(\alpha)}\right)
 =r-\max\{\deg f,\deg g\}.
\end{equation}
Every $\lambda\in E^\times$ with $I_W(\lambda)>0$ arises uniquely from such a reduced fraction, where a common factor in $k^\times$ in the pair $(f,g)$ is ignored.
For $\beta\in E\setminus k$, moreover,
\begin{equation}\label{v150:eq-power-pgl}
 P_r(\beta)\in\{uP_r(\alpha):u\in E^\times\}
 \iff
 \beta=\frac{a\alpha+b}{c\alpha+d},\quad
 a,b,c,d\in k,\quad ad-bc\ne0.
\end{equation}
\end{proposition}
\begin{proof}
An element of the intersection satisfies $A(\alpha)=f(\alpha)B(\alpha)/g(\alpha)$, with $\deg A,\deg B<r$.
The polynomial $gA-fB$ has degree at most $2r-2<t$, so vanishing at $\alpha$ forces it to be zero.
Coprimality gives $A=fh$, $B=gh$, with $\deg h<r-\max\{\deg f,\deg g\}$.
Evaluation and multiplication by $f(\alpha)$ are injective, proving~\eqref{v150:eq-power-intersection}.
A nonzero intersection point gives $\lambda=A(\alpha)/B(\alpha)$; cancellation preserves the degree bounds. Equality at $\alpha$ of two such fractions similarly forces their cross products to agree identically, proving uniqueness of the reduced fraction.

If $P_r(\beta)=uW$, this space has dimension $r$, so $1,\beta,\ldots,\beta^{r-1}$ are linearly independent. Hence $P_r(\beta)\cap\beta P_r(\beta)$ contains the $r-1$ linearly independent elements $\beta,\ldots,\beta^{r-1}$.
Since $1,\beta\in uW$, $\beta$ is a ratio of evaluations at $\alpha$ of two polynomials of degree less than $r$.
Formula~\eqref{v150:eq-power-intersection} forces its reduced numerator and denominator to have degrees at most one. As $\beta\notin k$, their determinant is nonzero.
Conversely, multiplication by $(c\alpha+d)^{r-1}$ identifies $1,\beta,\ldots,\beta^{r-1}$ with the evaluations of
\[
 (aX+b)^i(cX+d)^{r-1-i},\qquad 0\le i<r.
\]
An invertible linear substitution induces an automorphism on binary homogeneous forms of degree $r-1$, so these polynomials form a basis of $k[X]_{<r}$. This proves the converse.
\end{proof}

\begin{corollary}[The full intersection distribution of a power space]\label{v150:cor-power-distribution}
Under these hypotheses, the frequencies in Proposition~\ref{v70:prop-intersection-count} are
\begin{equation}\label{v150:eq-power-distribution}
 N_r=q-1,\qquad
 N_{r-d}=q^{2d-1}(q^2-1)\ (1\le d<r),\qquad
 N_0=q^t-q^{2r-1}.
\end{equation}
In particular, $H(W)=k$; the multipliers giving the largest nontrivial intersection dimension $r-1$ are exactly the fractional linear images of $\alpha$.
\end{corollary}
\begin{proof}
Let $C_d$ count nonzero reduced fractions $f/g$ with $\max(\deg f,\deg g)\le d$, ignoring a common nonzero scalar.
Grouping all pairs of nonzero polynomials by their monic greatest common divisor gives
\[
 \frac{(q^{d+1}-1)^2}{q-1}=\sum_{j=0}^dq^jC_{d-j}.
\]
Subtracting the preceding identity multiplied by $q$ gives $C_d=q^{2d+1}-1$. Injectivity of evaluation in the proposition now yields $N_r=C_0$, $N_{r-d}=C_d-C_{d-1}$, and $N_0=q^t-1-C_{r-1}$.
\end{proof}

\begin{example}[Power relations at the degree boundary]\label{v150:ex-power-boundary}
In $\F_{16}=\F_2(\alpha)$, take $\alpha^4+\alpha+1=0$. Then
\[
 P_3(\alpha^2)=\langle1,\alpha^2,\alpha^4\rangle_{\F_2}
 =\langle1,\alpha^2,\alpha+1\rangle_{\F_2}=P_3(\alpha).
\]
But $\alpha^2=(a\alpha+b)/(c\alpha+d)$ with $ad-bc\ne0$ would give a nonzero polynomial relation of degree at most three, a contradiction.
Here $t=4=2r-2$: the evaluation kernel can already change scalar equivalence, so the strict degree bound is essential.
\end{example}

\begin{example}[Scalar multiples of subfields and relative-trace hyperplanes]\label{v70:ex-trace-hyperplane-intersection}
Let \(k\subseteq H=\F_{q^d}\subseteq E\) and \(b\in E^\times\). For \(W=bH\), we have \(H(W)=H\) and
\[
 I_W(a)=\begin{cases}d,&a\in H^\times,\\0,&a\notin H.\end{cases}
\]
Indeed, a nonzero element \(x=bu=abv\) in the intersection gives \(a=u/v\in H\); conversely, \(a\in H^\times\) gives \(aW=W\).

Now assume \(H\subsetneq E\) and put \(W=\ker\Tr_{E/H}\). Surjectivity of the relative trace gives \(\dim_kW=n-d\).
Trace transitivity gives \(H\subseteq W^\perp\), and equality follows by dimension.
Thus \(H(W)=H(H)=H\), and \eqref{v70:eq-intersection-duality}, together with the first case, gives
\[
 I_W(a)=\begin{cases}n-d,&a\in H^\times,\\n-2d,&a\notin H.\end{cases}
\]
Consequently, \(N_{n-d}=q^d-1\) and \(N_{n-2d}=q^n-q^d\), with all other \(N_j\) zero.
Here \([E:H]\ge2\), so \(n-2d\ge0\); when \([E:H]=2\), the trace kernel is itself a one-dimensional \(H\)-space.
\end{example}

\begin{example}[Quotient multiplication in the field of order eight]\label{v70:ex-f8-quotient-multiplication}
Take \(E=\F_2(\alpha)\), where \(\alpha^3+\alpha+1=0\). We have
\[
 \Tr_{E/\F_2}(1)=1,\qquad
 \Tr_{E/\F_2}(\alpha)=\Tr_{E/\F_2}(\alpha^2)=0.
\]
Hence \(W=\ker\Tr_{E/\F_2}=\operatorname{span}_{\F_2}\{\alpha,\alpha^2\}\).
Identify \(E/W\) with \(\F_2\) through the trace. For \(a=a_0+a_1\alpha+a_2\alpha^2\), the relations \(\alpha^3=\alpha+1\) and \(\alpha^4=\alpha^2+\alpha\) give
\[
 \Tr(a\alpha)=a_2,\qquad \Tr(a\alpha^2)=a_1,
 \qquad [T_{a,W}]=(a_2,a_1).
\]
The kernel of the parameter map is therefore \(\F_2\), so \(H(W)=\F_2\). For nonzero parameters \(a\), the individual map has rank zero at \(a=1\) and rank one in the remaining six cases.
Thus \(N_2=1\), \(N_1=6\), and \(N_0=0\).
The same trace pairing gives coordinates on the quotient and expresses the multiplicative intersection dimension through the rank of a one-row matrix.
\end{example}

\section{Dickson Matrices: Linear Operators in Conjugate Coordinates}\label{sec:dickson}

A linearized polynomial represents a linear operator by its coefficients. When all coefficients lie in \(k\), Theorem~\ref{thm:linearized-gcd} determines its kernel from the greatest common divisor of operator polynomials;
when the coefficients range over all of \(E\), the polynomial and its Frobenius conjugates determine a matrix from which the root space, composition, and trace adjoint can be computed directly.

Let \(k=\F_q\subseteq E=\F_{q^n}\), as before, and write
\[
L(X)=\sum_{j=0}^{n-1}a_jX^{q^j},\qquad a_j\in E.
\]
Define the conjugate-value vector
\[
\mathbf c(x)=(x,x^q,\ldots,x^{q^{n-1}})^{\mathsf T}\in E^n.
\]
Raising \(L(x)\) to the \(q^i\)th power and using \(x^{q^n}=x\) gives
\[
L(x)^{q^i}=\sum_{j=0}^{n-1}a_{j-i}^{q^i}x^{q^j},
\]
where coefficient indices are taken modulo \(n\). Hence
\begin{equation}\label{eq:dickson-conjugate-action}
\boxed{\qquad \mathbf c(L(x))=D_L\mathbf c(x),\qquad
D_L=(a_{j-i}^{q^i})_{0\le i,j<n}.\qquad}
\end{equation}
\leavevmode\phantomsection\label{def:dickson-matrix}\ignorespaces
The matrix \(D_L\) is called the \textbf{Dickson matrix} of \(L\).
It is determined directly by the coefficients of \(L\) and their Frobenius conjugates, without first choosing a basis of \(E/k\).
The vector \(\mathbf c(x)\) records conjugate values; the Moore matrix below connects it with ordinary basis coordinates.

\begin{proposition}[Dickson matrices and root spaces]\label{prop:dickson-rank}
Let \(e_0,\ldots,e_{n-1}\) be a basis of \(E/k\), let \(B\in k^{n\times n}\) represent \(L\) in this basis,
and let \(M=(e_j^{q^i})_{0\le i,j<n}\) be the corresponding Moore matrix. Then
\[
D_LM=MB,\qquad D_L=MBM^{-1}.
\]
In particular,
\begin{equation}\label{eq:dickson-root-count}
\begin{aligned}
\rank_E D_L&=\rank_k L,&\qquad \det D_L&=\det B\in k,\\
\dim_k\ker L&=n-\rank_E D_L,&
|\{x\in E:L(x)=0\}|&=q^{\,n-\rank_E D_L}.
\end{aligned}
\end{equation}
Thus \(L:E\to E\) is invertible if and only if \(\det D_L\ne0\).
\end{proposition}
\begin{proof}
Column \(j\) of \(M\) is \(\mathbf c(e_j)\). If \(L(e_j)=\sum_rB_{rj}e_r\), then, since \(B_{rj}\in k\),
\[
D_L\mathbf c(e_j)=\mathbf c(L(e_j))
=\sum_rB_{rj}\mathbf c(e_r).
\]
Comparing columns gives \(D_LM=MB\). Theorem~\ref{thm:moore} ensures that \(M\) is invertible, proving the similarity relation.
The rank of \(B\) is unchanged on extending scalars from \(k\) to \(E\), since it is the largest size of a nonzero minor.
This proves the rank and determinant identities. The rank--nullity theorem in dimension \(n\), together with the cardinality of a finite-dimensional \(k\)-space, gives the root count and invertibility criterion.
\end{proof}

This similarity relation extends the conjugate-coordinate description of Theorem~\ref{thm:conjugate-coordinates} to arbitrary linear operators.
Consider the evaluation isomorphism from~\eqref{eq:conjugate-tensor-map}:
\[
\mathcal J:E\otimes_kE\longrightarrow E^n,\qquad
u\otimes x\longmapsto u\mathbf c(x).
\]
Identity~\eqref{eq:dickson-conjugate-action} states that
\[
\mathcal J\circ(\id_E\otimes L)=D_L\circ\mathcal J.
\]
Thus \(D_L\) is precisely the matrix of \(L\) in conjugate coordinates after extending scalars to \(E\).
This also explains why a rank computed over \(E\) still records the \(k\)-rank of the original map.

\begin{proposition}[Composition, inverses, and trace adjoints]\leavevmode\phantomsection\label{prop:dickson-operations}
Let \(L,H:E\to E\) be \(k\)-linear maps represented by linearized polynomials of degree less than \(q^n\). Then
\[
D_{L+H}=D_L+D_H,\qquad
D_{L\circ H}=D_LD_H,\qquad D_{\id_E}=I_n.
\]
Here the composed polynomial is replaced by its standard linearized representative modulo \(X^{q^n}-X\).
If \(L\) is invertible, then \(D_{L^{-1}}=D_L^{-1}\). The adjoint with respect to the trace pairing satisfies
\begin{equation}\label{eq:dickson-adjoint-transpose}
\boxed{\qquad D_{L^*}=D_L^{\mathsf T}.\qquad}
\end{equation}
\end{proposition}
\begin{proof}
The assertions for addition and the identity map follow directly from the definition. For composition,
\[
\mathbf c(L(H(x)))=D_L\mathbf c(H(x))=D_LD_H\mathbf c(x).
\]
Take \(x=e_0,\ldots,e_{n-1}\). Invertibility of the Moore matrix shows that these vectors \(\mathbf c(e_j)\) form a basis of \(E^n\) over \(E\),
so \(D_{L\circ H}=D_LD_H\). Setting \(H\) equal to \(L^{-1}\) gives the inverse formula.

The trace pairing has a particularly simple expression in conjugate coordinates:
\begin{equation}\label{eq:trace-conjugate-dot-product}
\Tr_{E/k}(xy)=\sum_{i=0}^{n-1}x^{q^i}y^{q^i}
=\mathbf c(x)^{\mathsf T}\mathbf c(y).
\end{equation}
Combining this with \(\Tr(L(x)y)=\Tr(xL^*(y))\) gives
\[
\mathbf c(x)^{\mathsf T}D_L^{\mathsf T}\mathbf c(y)
=\mathbf c(x)^{\mathsf T}D_{L^*}\mathbf c(y).
\]
Letting \(x,y\) range over the chosen basis yields \(M^{\mathsf T}(D_L^{\mathsf T}-D_{L^*})M=0\).
Since \(M\) is invertible, this proves~\eqref{eq:dickson-adjoint-transpose}.
\end{proof}

Matrix rank measures the root space, operator composition becomes matrix multiplication, and the trace adjoint becomes transposition.
These correspondences place the operator algebra of Section~\ref{sec:operator-algebra} and the trace duality of the preceding section in the same conjugate coordinates.
When \(L\) is invertible, the first row of \(D_L^{-1}\) gives the linearized coefficients of the inverse map.

\begin{example}[The norm criterion in a quadratic extension]\leavevmode\phantomsection\label{exm:dickson-quadratic}
Take \(E=\F_{q^2}\) and let \(L(x)=ax+bx^q\), where \(a,b\in E\). Then
\[
D_L=\begin{pmatrix}a&b\\b^q&a^q\end{pmatrix},\qquad
\Delta:=\det D_L=a^{q+1}-b^{q+1}
=\Nm_{E/k}(a)-\Nm_{E/k}(b).
\]
Thus \(L\) is invertible if and only if \(\Nm(a)\ne\Nm(b)\). Under this condition, \(\Delta\in k^\times\), and the first row of the inverse matrix gives
\[
L^{-1}(y)=\frac{a^qy-by^q}{\Delta}.
\]
If \(\Delta=0\) but \((a,b)\ne(0,0)\), then \(D_L\) has rank one, so \(L\) has exactly \(q\) roots in \(E\);
if \(a=b=0\), all of \(E\) is the root space.

When both \(a,b\) are nonzero, a nonzero root satisfies \(x^{q-1}=-a/b\). Multiplicative Hilbert 90 reduces solvability to
\(\Nm(-a/b)=1\), which is also equivalent to \(\Nm(a)=\Nm(b)\).
The norm criterion and the Dickson determinant criterion therefore give the same condition in this case.
\end{example}

The fixed spaces of iterates can also be computed from \(D_L^t-I\); Section~\ref{v53:subsec-linear-cycles} uses this to recover all cycles,
while Section~\ref{v53:subsec-linear-parity} uses the Dickson determinant to determine the sign of a permutation in odd characteristic.

\par\Needspace{11\baselineskip}
\section*{Exercises}
\phantomsection
\addcontentsline{toc}{section}{Exercises}
\begin{exercise}\label{ex:subspace-line}
For \(u\in E^\times\), prove that the root polynomial of the one-dimensional subspace \(ku\) is \(X^q-u^{q-1}X\).
Use the recurrence in this chapter to find the root polynomial of the two-dimensional space \(ku+kv\), where \(u,v\) are independent over \(k\).
\end{exercise}
\begin{exercise}\label{ex:kernel-gcd}
Find the kernel and image dimensions of \(x\mapsto x+x^4\) on \(\F_{2^6}\).
\end{exercise}
\begin{exercise}\label{ex:trace-adjoint}
Prove that the trace operator is self-adjoint, and use \(\im L=(\ker L^*)^\perp\) to recover \((\ker\Tr)^\perp=k\).
\end{exercise}
\begin{exercise}\label{ex:linearized-not-ordinary}
Explain why \(L_U\mid L_V\) as ordinary polynomials is equivalent to \(U\subseteq V\),
but the ordinary least common multiple generally differs from \(L_{U+V}\). Give an example in a two-dimensional space.
\end{exercise}

\begin{exercise}\label{ex:supp-symbolic-division}
In \(\F_2[X]\), express \(X^8+X\) as \(Q\circ(X^2+X)\) and explain why divisibility holds both in the ordinary polynomial ring and on the right under composition.
\end{exercise}

\begin{exercise}[Noncommutative composition and reduction of functions]\label{ex:wu-noncommuting}
In \(\F_4=\F_2(\alpha)\) with \(\alpha^2+\alpha+1=0\), let \(L(X)=\alpha X^2\) and \(M(X)=X^2\).
Compare \(L\circ M\) with \(M\circ L\) and determine whether they represent the same function on \(\F_4\).
Give another example showing that distinct linearized polynomials can represent the same function.
\end{exercise}
\begin{exercise}[Operators commuting with an iterate of Frobenius]\label{ex:wu-centralizer}
Let \(s\ge1\) and \(g=\gcd(n,s)\). Prove that a \(k\)-linear map \(A(x)=\sum_{i=0}^{n-1}a_ix^{q^i}\) commutes with \(\sigma^s\) if and only if every coefficient satisfies \(a_i\in\F_{q^g}\).
Find the dimension of the \(k\)-space of these operators.
\end{exercise}

\begin{exercise}[The Ore correspondence and Dickson matrices]\label{ex:ore-dickson-comparison}
Let \(E=\F_{q^n}\) and \(k=\F_q\), and put \(\sigma(x)=x^q\).
For \(a(T)\in k[T]\), the correspondence \(a(T)=\sum_i a_iT^i\mapsto\mathcal L_a(X)=\sum_i a_iX^{q^i}\) is also called the \textbf{Ore correspondence}.
\begin{enumerate}
\item Let \(a_i\in k\). Use the kernel-dimension formula to show that \(\mathcal L_a\) permutes \(E\) if and only if
\(\gcd(a,T^n-1)=1\), and construct the inverse map from B\'ezout's identity.
\item If \(a_0,\ldots,a_{n-1}\in E\), explain why the permutation property of \(L(X)=\sum_{i=0}^{n-1}a_iX^{q^i}\) should be tested by
\(\det D_L\ne0\). Use \(\sigma M_c=M_{c^q}\sigma\) to explain why the commutative polynomial operations in the previous part do not transfer directly to \(E[T]\).
\item Take \(E=\F_4=\F_2(\omega)\) and \(\omega^2+\omega+1=0\), and let \(L(X)=X+\omega X^2\).
Compare \(\gcd_{E[T]}(1+\omega T,T^2-1)\), \(\det D_L\), and \(\ker L\) to show that the base-field coefficient condition cannot be omitted.
\end{enumerate}
\end{exercise}

\begin{exercise}[The complete multiplicative intersection distribution of a two-dimensional space]\label{v70:ex-two-dimensional-intersections}
Let \(n\ge2\), and let \(W\subseteq E=\F_{q^n}\) be a two-dimensional \(k=\F_q\)-subspace. Use the notation \(N_j\) from Proposition~\ref{v70:prop-intersection-count}. Prove that exactly one of the following cases occurs:
\begin{enumerate}
\item \(2\mid n\) and \(W=b\F_{q^2}\) for some \(b\in E^\times\). Then
\[
 N_2=q^2-1,\qquad N_1=0,\qquad N_0=q^n-q^2.
\]
\item The space \(W\) is not such a scalar multiple of a subfield. Then \(H(W)=k\), and
\[
 N_2=q-1,\qquad N_1=q(q^2-1),\qquad N_0=q^n-q^3.
\]
\end{enumerate}
Explain why the second case can occur only for \(n\ge3\), and determine the total number of distinct scalar multiples of \(W\) in each case.
\end{exercise}

\begin{exercise}[Reciprocal sums and additional relations after evaluation]\label{v116:ex-inverse-sum}
Let \(U\subseteq\F_{2^n}\) be a binary subspace of dimension \(s\ge1\), with \(L_U(X)=\sum_{i=0}^s c_iX^{2^i}\) and \(c_0\ne0\). Prove
\[
\sum_{u\in U\setminus\{0\}}u^{-1}=c_1/c_0.
\]
For \(s=4\), interpret \(c_1=0\) as a relation in \(\mathcal C(U)\).
Then take \(\alpha^3+\alpha+1=0\) and compare the reciprocal rank of \(\langle1,\alpha,\alpha^2\rangle_{\F_2}\) with the first row of Example~\ref{v116:ex-inverse-models}. At which step can evaluation create additional relations?
\end{exercise}

\chapter{Normal Bases, Cyclic Modules, and Invariant Subspaces}
\label{chap:normal}
A normal basis represents Frobenius as a cyclic shift, making the extension field a concrete cyclic module. This unifies normality criteria, root-space counts, and the classification of invariant subspaces; reciprocation and the trace pairing further describe duality.
Families of translates apply the unit criterion to changes of generator: the same cyclic-module model determines which primary components can be repaired and which obstructions remain unchanged.

\section{Orbit Length and the Dimension of Its Span}
\phantomsection\label{def:cyclic-module}
A \(k\)-linear operator \(A:V\to V\) makes \(V\) a \(k[T]\)-module: the action of a polynomial \(a(T)\) on a vector \(v\) is defined as \(a(A)v\).
If the iterates \(v,Av,A^2v,\ldots\) of a vector \(v\) span all of \(V\), then \(v\) is called a \textbf{cyclic vector},
and the corresponding \(k[T]\)-module is called a \textbf{cyclic module}.
A subspace \(U\) satisfying \(A(U)\subseteq U\) is called \(A\)-invariant; it is precisely a submodule for this module structure.
In this chapter \(A=\sigma\), so a polynomial acts as a linear combination of Frobenius iterates.

\begin{definition}\label{def:normal-basis}
Let \(E=\F_{q^n}\) and \(k=\F_q\), and let \(\eta\) be a \hyperref[def:normal-element]{normal element}. Its conjugates
\[
\eta,\eta^q,\ldots,\eta^{q^{n-1}}
\]
form a basis of \(E/k\), called a \textbf{normal basis}.
\end{definition}
A normal element has \(n\) distinct conjugates and therefore generates the extension field; the converse fails.
Orbit length counts distinct elements, whereas normality requires these elements to be linearly independent.

\begin{example}\leavevmode\phantomsection\label{ex:normal-f8-example}
Let \(\F_8=\F_2(\alpha)\) with \(\alpha^3+\alpha+1=0\).
The element \(\alpha\) has degree three, but \(\alpha+\alpha^2+\alpha^4=0\), so it is not normal.
Put \(\eta=\alpha+1\). Its conjugates are
\[
\alpha+1,\quad\alpha^2+1,\quad\alpha^2+\alpha+1.
\]
Their coordinate matrix in the basis \((1,\alpha,\alpha^2)\) is
\[
\begin{pmatrix}1&1&1\\1&0&1\\0&1&1\end{pmatrix}.
\]
Its determinant is one in \(\F_2\), so the three conjugates are linearly independent and \(\eta\) is normal.
\end{example}

\phantomsection\label{sec:q-order}
Regard \(\sigma:x\mapsto x^q\) as a \(k\)-linear operator.
For \(a\in E\), define its \textbf{Frobenius annihilator polynomial} \(\mu_a(T)\)
to be the monic polynomial of least degree satisfying \(\mu_a(\sigma)(a)=0\).
It is also called the \(q\)-order polynomial of \(a\). Set \(\mu_0=1\).
Division with remainder shows that every operator polynomial annihilating \(a\) is divisible by \(\mu_a\); in particular, \(\mu_a\mid T^n-1\).

If \(d=\deg\mu_a\), then \(a,\sigma(a),\ldots,\sigma^{d-1}(a)\) are linearly independent and span the entire orbit.
Hence
\[
\deg\mu_a=\dim_k\langle a,\sigma(a),\sigma^2(a),\ldots\rangle_k.
\]
This operator polynomial generally differs from the minimal polynomial \(m_{a,k}(X)\) of the field element.
The former describes linear relations among conjugates; the latter has the distinct conjugates as its roots.

\section{The Normal Basis Theorem}

\par
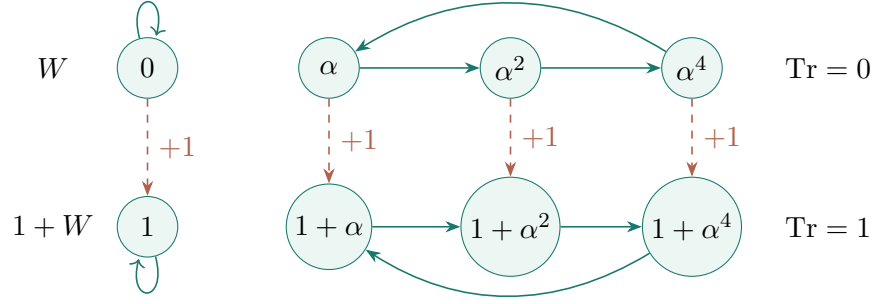
\begin{figure}[!htbp]
\centering
\begin{tikzpicture}[
  vpoint/.style={circle,draw=Jade,fill=CoverMint,minimum size=8mm,inner sep=2pt},
  varrow/.style={-{Stealth[length=2mm]},Jade,semithick},
  vbox/.style={draw=MutedSlate,rounded corners=2pt,inner sep=5pt,align=center},
  every node/.style={font=\small}]

\node at (-1.25,0) {\(W\)};
\node at (-1.25,-2.1) {\(1+W\)};
\node[vpoint] (z) at (0,0) {\(0\)};
\node[vpoint] (o) at (0,-2.1) {\(1\)};
\draw[varrow] (z) edge[loop above] (z);
\draw[varrow] (o) edge[loop below] (o);
\node[vpoint] (a) at (2.4,0) {\(\alpha\)};
\node[vpoint] (b) at (4.8,0) {\(\alpha^2\)};
\node[vpoint] (c) at (7.2,0) {\(\alpha^4\)};
\node[vpoint] (aa) at (2.4,-2.1) {\(1+\alpha\)};
\node[vpoint] (bb) at (4.8,-2.1) {\(1+\alpha^2\)};
\node[vpoint] (cc) at (7.2,-2.1) {\(1+\alpha^4\)};
\draw[varrow] (a)--(b);\draw[varrow] (b)--(c);
\draw[varrow] (c) to[bend right=32] (a);
\draw[varrow] (aa)--(bb);\draw[varrow] (bb)--(cc);
\draw[varrow] (cc) to[bend left=32] (aa);
\foreach \u/\v in {z/o,a/aa,b/bb,c/cc}
  \draw[varrow,dashed,Clay] (\u)--node[right] {\(+1\)} (\v);
\node at (9,0) {\(\Tr=0\)};
\node at (9,-2.1) {\(\Tr=1\)};

\end{tikzpicture}
\caption{Solid arrows are Frobenius; dashed arrows are translation. Both nontrivial cycles have length three, but the upper one spans a two-dimensional space and the lower one is a normal basis.}\label{v102:fig-cyclic-eight}
\end{figure}

\Needspace{7\baselineskip}
\noindent\phantomsection\label{v102:par-cyclic-eight}\textbf{Two invariant subspaces inside the eight-element field.}\quad
Use \(E=\F_2(\alpha)\), \(\alpha^3+\alpha+1=0\), from Example~\ref{ex:normal-f8-example}, and put \(\sigma(x)=x^2\).
Since \(\alpha^4=\alpha^2+\alpha\), we obtain
\[
W=\ker\Tr_{E/\F_2}=\{0,\alpha,\alpha^2,\alpha^4\},\qquad
\Tr_{E/\F_2}(1)=1.
\]
Every \(x\in E\) therefore has a unique expression \(x=c+w\), with \(c=\Tr(x)\in\F_2\) and \(w=x+\Tr(x)\in W\).
This gives the decomposition \(E=\F_2\oplus W\) as an \(\F_2\)-space and an \(\F_2[T]\)-module, with \(T\) acting as \(\sigma\).
Figure~\ref{v102:fig-cyclic-eight} displays all eight elements; the lower row is the coset obtained by adding one to the upper row.

\FloatBarrier
The projections can also be written as operators. Put
\[
P(T)=T^2+T+1,\quad e_0(T)=1+T+T^2,\quad e_1(T)=T+T^2.
\]
In \(\F_2[T]/(T^3-1)\), the residues of \(e_0,e_1\) modulo \(T+1,P\) are \((1,0),(0,1)\), respectively. Thus
\[
e_0(\sigma)x=\Tr(x),\qquad e_1(\sigma)x=x+\Tr(x).
\]
These are the projections onto \(\F_2\) and \(W\). The vectors \(\alpha,\alpha^2\) form a basis of \(W\), so
\[
(E,\sigma)\simeq\F_2[T]/(T+1)\ \oplus\ \F_2[T]/(P).
\]
The following table lists the annihilators of all eight vectors, with \(\mu_0=1\):

\[
\begin{array}{c|cccc}
w&0&\alpha&\alpha^2&\alpha^4\\\hline
\mu_w&1&P&P&P\\
\mu_{1+w}&T+1&(T+1)P&(T+1)P&(T+1)P
\end{array}
\]
A nonzero \(w\in W\) has annihilator \(P\), since \(P(\sigma)w=0\) and \(w\notin\F_2\).
Applying the two projections to \(1+w\) recovers \(1\) and \(w\); for \(w\ne0\), its annihilator must therefore contain both coprime factors \(T+1\) and \(P\).
A cyclic vector must generate both components. This anticipates the proof below: project to the primary components, then add their cyclic vectors.
The displayed isomorphism is a \textbf{module isomorphism}: \(E\) remains a field, and \(W\) is not a nonzero proper ideal of it. It is not a product decomposition of the field \(E\).

\begin{lemma}[Cyclic vector lemma]\label{lem:cyclic-vector}
Let \(A\) be a linear operator on a finite-dimensional \(k\)-space \(V\), and let \(m_A\) be the monic polynomial of least degree satisfying \(m_A(A)=0\).
There exists \(v\in V\) whose monic annihilator polynomial of least degree with respect to \(A\) is exactly \(m_A\).
In particular, if \(\deg m_A=\dim V\), then \(v,A(v),\ldots,A^{\dim V-1}(v)\) form a basis.
\end{lemma}
\begin{proof}
If \(V=\{0\}\), take \(v=0\); the corresponding minimal polynomial is \(1\) and the assertion holds. Henceforth assume \(V\ne0\).
The definition of operator annihilating polynomials on page~\pageref{def:operator-annihilator} and division with remainder show that \(m_A\) and each \(\mu_w\) exist and divide every corresponding annihilating polynomial.

Choose a basis \(e_1,\ldots,e_d\). An operator \(g(A)\) is zero if and only if it vanishes on every \(e_j\);
therefore
\[
m_A=\lcm(\mu_{e_1},\ldots,\mu_{e_d}).
\]
Write the irreducible factorization \(m_A=\prod_{i=1}^r P_i^{a_i}\).
For each \(i\), the definition of the least common multiple provides a basis vector \(u_i\) such that \(\mu_{u_i}=P_i^{a_i}R_i\), where \(\gcd(P_i,R_i)=1\).
Put \(v_i=R_i(A)u_i\). Then
\[
g(A)v_i=0
\ \Longleftrightarrow\ \mu_{u_i}\mid gR_i
\ \Longleftrightarrow\ P_i^{a_i}\mid g.
\]
Thus \(\mu_{v_i}=P_i^{a_i}\).

Put \(v=\sum_i v_i\). To show that distinct components cannot cancel, set \(H_i=\prod_{j\ne i}P_j^{a_j}\) and use B\'ezout's identity to choose \(U_i\) with \(U_iH_i\equiv1\pmod{P_i^{a_i}}\).
The operator \((U_iH_i)(A)\) is the identity on \(v_i\) and zero on \(v_j\) for \(j\ne i\),
so it maps \(v\) to \(v_i\).
If \(g(A)v=0\), the two operators commute because both are polynomials in \(A\), and hence
\[
g(A)v_i=(U_iH_i)(A)g(A)v=0.
\]
Thus \(P_i^{a_i}\mid g\) for every index; pairwise coprimality gives \(m_A\mid g\).
On the other hand, \(m_A(A)v=0\), so \(\mu_v=m_A\).

If \(\deg m_A=d=\dim V\), a nontrivial linear relation among \(v,Av,\ldots,A^{d-1}v\)
would give an annihilating polynomial of degree less than \(d\), which is impossible. These are \(d\) linearly independent vectors and hence form a basis.
\end{proof}

\begin{theorem}[Normal basis theorem]\label{thm:normal-basis}
Every finite-field extension \(\F_{q^n}/\F_q\) admits a normal basis.
\end{theorem}
\begin{proof}
Theorem~\ref{thm:operator-algebra} shows that the map \(k[T]\to k[\sigma]\),
\(a\mapsto a(\sigma)\) has kernel \((T^n-1)\). Thus the minimal polynomial of \(\sigma\) is exactly \(T^n-1\),
whose degree is \(\dim_k E=n\).
The cyclic vector lemma~\ref{lem:cyclic-vector} therefore gives \(\eta\) such that \(\eta,\sigma(\eta),\ldots,\sigma^{n-1}(\eta)\) form a basis.
\end{proof}

The minimal-polynomial fact used above also has a direct verification. If
\(\sum_{i=0}^{n-1}a_i\sigma^i=0\), the linearized polynomial
\(\sum_{i=0}^{n-1}a_iX^{q^i}\) vanishes at all \(q^n\) elements of \(E\).
Unless all \(a_i\) vanish, this is a nonzero polynomial of degree at most
\(q^{n-1}<q^n\), contradicting the root bound. Thus no polynomial of degree
less than \(n\) annihilates \(\sigma\); together with \(\sigma^n=I\), this gives
\(m_\sigma=T^n-1\). The root bound identifies the minimal polynomial;
the cyclic-vector lemma is still needed to obtain a normal basis.

\begin{example}[Normal bases in the presence of repeated factors]
In \(\F_4/\F_2\), take \(\eta^2+\eta+1=0\).
The pair \((\eta,\eta^2)\) is a normal basis, and the matrix of Frobenius is
\[
S=\begin{pmatrix}0&1\\1&0\end{pmatrix}.
\]
In characteristic two, \(S\ne I\), while \((S-I)^2=0\); the minimal polynomial of the operator is therefore \((T-1)^2=T^2-1\), which is not square-free.
The only eigenvalue is one and the operator is not diagonalizable, but \(\eta,S\eta\) still form a cyclic basis.
Thus the normal basis theorem depends on the cyclic vector property, not on diagonalizability of Frobenius.

\par\Needspace{7\baselineskip}
\noindent\phantomsection\label{v102:par-repeated-root-chain}\textbf{Distinguishing a two-cycle from a nilpotent chain.}\quad
In this example put \(N=\sigma-I\). Direct calculation gives \(N(\eta)=\eta^2+\eta=1\) and \(N(1)=0\).
The two sets of arrows in Figure~\ref{v102:fig-repeated-root-chain} therefore represent different operators.

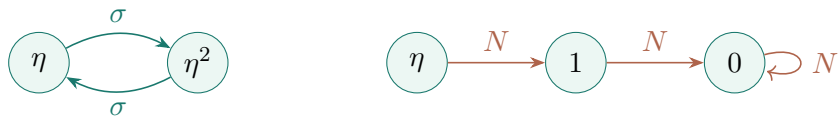
\begin{figure}[!htbp]
\centering
\begin{tikzpicture}[
  vpoint/.style={circle,draw=Jade,fill=CoverMint,minimum size=8mm,inner sep=2pt},
  varrow/.style={-{Stealth[length=2mm]},Jade,semithick},
  vbox/.style={draw=MutedSlate,rounded corners=2pt,inner sep=5pt,align=center},
  every node/.style={font=\small}]

\node[vpoint] (a) at (0,0) {\(\eta\)};
\node[vpoint] (b) at (2.1,0) {\(\eta^2\)};
\draw[varrow] (a) to[bend left=28] node[above] {\(\sigma\)} (b);
\draw[varrow] (b) to[bend left=28] node[below] {\(\sigma\)} (a);
\node[vpoint] (c) at (5,0) {\(\eta\)};
\node[vpoint] (d) at (7.1,0) {\(1\)};
\node[vpoint] (e) at (9.2,0) {\(0\)};
\draw[varrow,Clay] (c)--node[above] {\(N\)} (d);
\draw[varrow,Clay] (d)--node[above] {\(N\)} (e);
\draw[varrow,Clay] (e) edge[loop right] node[right] {\(N\)} (e);

\end{tikzpicture}
\caption{Frobenius on the left interchanges the normal basis vectors. The difference operator on the right satisfies \(N\ne0\), \(N^2=0\); the repeated factor appears as a nilpotent chain.}\label{v102:fig-repeated-root-chain}
\end{figure}
\FloatBarrier
The corresponding module is \(\F_2[T]/((T-1)^2)\). The two occurrences of \(T-1\) are not coprime, so the Chinese remainder decomposition used for the eight-element field does not split this module into two lines.
Indeed, on any invariant one-dimensional \(\F_2\)-space, the invertible operator \(\sigma\) must be the identity. A direct sum of two such lines would force \(\sigma=I\), a contradiction.

\end{example}

In a normal basis, the matrix of \(\sigma\) is a cyclic shift matrix, with both characteristic and minimal polynomial \(T^n-1\).

\subsection{The Root Bound and Linear Nondegeneracy}
\label{sec:root-bound-rigidity}

The key restriction in the normal-basis proof is that a nonzero polynomial of low degree cannot vanish on too large a set.
For linearized polynomials, comparing ordinary degree with the cardinality of a subspace turns this restriction into linear independence.
Specifically, if \(a_0,\ldots,a_s\in E\) are not all zero, then
\[
\left|\left\{x\in E:\sum_{i=0}^s a_i x^{q^i}=0\right\}\right|\le q^s.
\]
Thus this operator cannot annihilate a \(k\)-subspace of dimension greater than \(s\).
Four results already obtained use this same restriction:
\begin{enumerate}
\item \textbf{Cyclicity of the multiplicative group.}
The exponent \(e\) of a finite multiplicative subgroup is attained as the order of an element, and every group element is a root of \(X^e-1\).
The group order is therefore at most \(e\), while that element already generates a cyclic subgroup of \(e\) elements. The two groups coincide (Theorem~\ref{thm:cyclic}).
\item \textbf{Nondegeneracy of the trace pairing.}
The trace polynomial \(X+X^q+\cdots+X^{q^{n-1}}\) is nonzero of degree less than \(q^n\) and cannot vanish on all of \(E\).
Scaling a nonzero trace value by a base-field scalar gives an element of trace one; multiplication by the inverse of any nonzero element then proves nondegeneracy of the trace pairing
(Theorems~\ref{thm:trace-norm-surjective} and~\ref{thm:trace-nondegenerate}).
\item \textbf{The Moore criterion.}
If the Moore determinant of \(s\) elements linearly independent over the base field vanished, a nonzero linearized polynomial of \(q\)-degree at most \(s-1\)
would vanish on their span of \(q^s\) elements, a contradiction (Theorem~\ref{thm:moore}).
\item \textbf{Existence of normal bases.}
The same degree comparison excludes every operator relation for \(\sigma\) of degree less than \(n\), forcing its minimal polynomial to be \(T^n-1\);
the cyclic vector lemma then gives a normal basis. The root bound determines the minimal-polynomial degree, and the cyclic vector lemma converts that degree into a basis.
\end{enumerate}
Preservation of the field operations by Frobenius supplies fixed points, orbits, and linearized structure; the root bound excludes unintended degeneracies within these structures.
Together they show how counting restrictions in finite fields become generation, surjectivity, and linear independence.

\section{Testing and Counting Normal Elements}
Fix a normal element \(\eta\), and put \(F=T^n-1\). The map
\[
k[T]/(F)\longrightarrow E,\qquad \overline a\longmapsto a(\sigma)\eta
\]
is a \(k[T]\)-module isomorphism that identifies multiplication by \(T\) with Frobenius. Since
\(g(\sigma)a(\sigma)\eta=0\) is equivalent to \(F\mid ga\),
\[
\mu_{a(\sigma)\eta}=\frac{F}{\gcd(F,a)}.
\]
Thus normal elements correspond exactly to units of the quotient ring. The unit condition can be expressed both by operator kernels and by trace orthogonality.

\begin{proposition}[Local obstructions to normality and trace duality]\label{prop:normal-test}
Let \(M=\operatorname{rad}(F)\). For each distinct monic irreducible factor \(P\) of \(F\), put
\[
A_P=(F/P)(\sigma),\qquad U_P=P(\sigma)E,\qquad V_P=U_P^\perp,
\]
where orthogonal complements and adjoints refer to \(\Tr_{E/k}(xy)\). Then
\begin{equation}\label{v128:eq-local-obstructions}
\ker A_P=U_P,\qquad \operatorname{rank}_k A_P=\deg P,
\qquad V_P=\im A_P^*.
\end{equation}
For \(x\in E\), the following are equivalent: \(x\) is normal; \(\mu_x=F\); \(A_Px\ne0\) for every \(P\mid F\); for every \(P\mid F\), there exists \(b\in V_P\) with \(\Tr_{E/k}(bx)\ne0\). Moreover,
\begin{equation}\label{v128:eq-obstruction-dual-sum}
\bigcap_{P\mid F}U_P=M(\sigma)E,\qquad
\bigoplus_{P\mid F}V_P=(M(\sigma)E)^\perp.
\end{equation}
These statements allow repeated factors in \(F\); division by \(P\) lowers its multiplicity by only one.
\end{proposition}
\begin{proof}
In the cyclic module above, \((F/P)a\equiv0\pmod F\) is equivalent to \(P\mid a\). Hence \(\ker A_P=U_P\), with codimension \(\deg P\). The image--kernel relation for trace adjoints gives \(V_P=\im A_P^*\). The unit condition requires avoidance of every \(U_P\), and \(x\notin U_P=V_P^\perp\) is equivalent to the stated nonzero trace condition.

The distinct factors \(P\) are coprime, so the ideals \((P)/(F)\) intersect in \((M)/(F)\). Taking orthogonal complements gives the sum of the \(V_P\). Its dimension is \(\deg M=\sum_{P\mid F}\deg P\), so the sum is direct.
\end{proof}

The space \(U_P\) consists of the elements for which one local generation condition fails; \(V_P\) gives all weighted traces detecting that obstruction. For an individual element one can use \(A_P\); for a candidate family one can determine its intersection with \(U_P\). In counting, the same \(V_P\) becomes a space of additive frequencies; see Proposition~\ref{v128:prop-local-spectrum}.

\par\Needspace{5\baselineskip}
\noindent\phantomsection\label{v98:par-prescribed-annihilator}\textbf{Prescribing the Frobenius annihilator.}
The same formula can be used in reverse. Given a monic nonconstant divisor \(D\mid F=T^n-1\), take
\[
\alpha=\left(\frac{F}{D}\right)(\sigma)\eta,
\qquad \mu_\alpha=D.
\]
The second equality follows directly from \(\gcd(F,F/D)=F/D\), even when \(F\) has repeated factors.
Thus the conjugates of \(\alpha\) span a space of dimension exactly \(\deg D\), whereas its degree over \(k\) is the least positive integer \(r\) for which \(D\mid T^r-1\), since this divisibility is equivalent to \(\sigma^r(\alpha)=\alpha\).
Choosing \(D\) determines these two distinct data, and the normal element \(\eta\) then yields the required element; taking \(D=F\) recovers \(\eta\) itself.

\begin{theorem}\label{thm:normal-count}
The number of normal elements is
\[
\varphi_q(T^n-1)=q^n\prod_{P\mid T^n-1}\left(1-q^{-\deg P}\right),
\]
where the product runs over distinct monic irreducible factors and \(\varphi_q(F)\) denotes the number of units of \(k[T]/(F)\).
When normal bases are identified up to cyclic shift, their number is this quantity divided by \(n\).
\end{theorem}
\begin{proof}
The preceding cyclic-module isomorphism identifies normal elements exactly with units of \(k[T]/(T^n-1)\).
Formula~\eqref{eq:polynomial-euler} gives the order of this unit group, yielding the product formula.
The \(n\) conjugates of a normal element are distinct and give different cyclic starting points of the same normal basis, so divide by \(n\) to count bases.
\end{proof}

A normal element has nonzero trace, since otherwise the sum of its conjugates would give a nontrivial linear relation.
In particular, if \(n=p^a\), then \(T^n-1=(T-1)^n\) and there is just one normality condition:
\[
\left(\frac{T^n-1}{T-1}\right)(\sigma)(x)=\Tr_{E/k}(x)\ne0.
\]
Thus, when the extension degree is a power of the characteristic \(p\), the normal elements are exactly those with nonzero trace.

Normality can also be forced by a stronger family of linear-independence conditions. Section~\ref{v84:subsec-trace-mds} proves that, when the irreducible trace code of a multiplicative subgroup is nontrivial and MDS, the Frobenius orbit of every nonidentity element is a normal basis. The subgroup condition first excludes proper subfields, and column independence then recovers generation of the Frobenius cyclic module; see Theorem~\ref{v84:thm-subgroup-rigidity}.

Upgrading the normality of every nonidentity element to joint independence across a whole subgroup requires control of additive relations between different Frobenius orbits; Example~\ref{v96:ex-nineteen} gives an explicit subgroup whose nontrivial orbits are all normal despite a three-term relation.

\subsection{The Effect of Translation on Normality}\label{v82:subsec-normal-translates}
For a fixed \(\theta\in E\), adding a base-field element preserves the generated field: \(k(\theta+t)=k(\theta)\).
Normality, however, requires the conjugates to generate the entire cyclic module, so one must determine which generation conditions such a translation can change.

\begin{proposition}[Normal elements in a base-field translate]\label{v82:prop-normal-translates}
Let \(k=\F_q\), \(E=\F_{q^n}\), \(n\ge1\), and \(p=\Char k\), and put \(F=T^n-1\) and \(\sigma(x)=x^q\). Fix \(\theta\in E\).
\begin{enumerate}
\item Suppose \(p\nmid n\). If \((F/P)(\sigma)\theta=0\) for some monic irreducible factor \(P\mid F\) other than \(T-1\), then \(\theta+k\) contains no normal element.
If \((F/P)(\sigma)\theta\ne0\) for every such \(P\), then
\begin{equation}\label{v82:eq-normal-translate-trace}
 \theta+t\text{ is normal}
 \quad\Longleftrightarrow\quad
 \Tr_{E/k}(\theta)+nt\ne0\qquad(t\in k).
\end{equation}
Thus the family contains exactly \(q-1\) normal elements, and its unique excluded parameter is \(t=-\Tr_{E/k}(\theta)/n\).
\item If \(p\mid n\), then, for every \(t\in k\), the element \(\theta+t\) is normal if and only if \(\theta\) is normal.
\end{enumerate}
In particular, the number of normal elements in \(\theta+k\) is zero or \(q-1\) when \(p\nmid n\), and zero or \(q\) when \(p\mid n\).
\end{proposition}
\begin{proof}
For \(H\in k[T]\) and \(t\in k\), the equality \(\sigma(t)=t\) gives \(H(\sigma)t=H(1)t\). Hence
\begin{equation}\label{v82:eq-translate-annihilator}
 \left(\frac FP\right)(\sigma)(\theta+t)
 =\left(\frac FP\right)(\sigma)\theta
  +\left(\frac FP\right)(1)t.
\end{equation}
If \(P\ne T-1\), then \(T-1\mid F/P\), so the last term vanishes.
For \(P=T-1\), the operator \((F/P)(\sigma)\) is \(\Tr_{E/k}\), and \((F/P)(1)=n\) in \(k\).
When \(p\nmid n\), only this condition changes with \(t\), and it fails at exactly one parameter.
When \(p\mid n\), none of the conditions changes with \(t\).
Proposition~\ref{prop:normal-test} now proves the assertions; the original multiplicities of all factors of \(F\) have been retained throughout.
\end{proof}

When the characteristic does not divide the extension degree, a base-field translation can change only the \(T-1\) component.
When the characteristic divides the extension degree, such a translation cannot change whether any primary component is generated in full.
The latter statement does not mean that every component of the element is unchanged; it means that the nonvanishing conditions deciding normality are unchanged.
Consequently, changing the trace generally cannot repair a failure of generation in another primary component.
Exercise~\ref{v82:ex-universal-normal-translation} classifies the extensions in which every full-degree element can be made normal by a base-field translation;
Exercise~\ref{v82:ex-subfield-translates} considers translation by an intermediate subfield.
The reciprocal change of variable in Section~\ref{sec:construct-normal-dual} is not subject to the same linear restriction; Example~\ref{v68:ex-reciprocal-normal-f343} will compare the two changes of generator.
Subsection~\ref{v114:sec-reciprocal-exceptions} further uses coordinate numerators to locate the failures and to distinguish a fixed family, an average, and a maximum.

\paragraph{Subspace translations with prescribed annihilator data.}
\phantomsection\label{v113:par-annihilator-question}
Suppose that an element generates a given Frobenius cyclic module. After translation along an invariant subspace, how many elements still generate that same module? The dimension of the direction space counts the candidates; retaining the original annihilator also requires each primary component to retain its original generation property.
Base-field translation moves along only one invariant direction. We now extend the question to \(\alpha+\ker H(\sigma)\), asking when a perturbation can cause a primary component to lose that property. The decisive distinction is whether the corresponding irreducible factor occurs once or repeatedly in the original annihilator. The small-field comparisons that follow will show why direction spaces of the same dimension can give different preservation counts.

\begin{proposition}[Affine fibers with a fixed annihilator]\label{v106:prop-annihilator-fiber}
Let \(E=\F_{q^n}\), \(\sigma(x)=x^q\), and \(\mu_\alpha=D\). Let \(H\) be a monic squarefree divisor of \(D\), and put \(W=\ker H(\sigma)\). Then \(\alpha+W\subseteq\ker D(\sigma)\), and the number of elements in this fiber whose annihilator is exactly \(D\) is
\begin{equation}\label{v106:eq-annihilator-fiber}
\boxed{q^{\deg H}\prod_{\substack{P\mid H\\P^2\nmid D}}(1-q^{-\deg P}).}
\end{equation}
The product runs over monic irreducible factors. In particular, if \(H=P\) is irreducible, every element of the fiber retains \(D\) when \(P^2\mid D\); exactly one element fails to retain \(D\) when \(P^2\nmid D\).
\end{proposition}
\begin{proof}
The cyclic-module isomorphism
\[
k[T]/(D)\longrightarrow\ker D(\sigma),\qquad \bar a\longmapsto a(\sigma)\alpha
\]
follows from \(\mu_\alpha=D\) and \(\dim_k\ker D(\sigma)=\deg D\). Under this isomorphism, \(W\) corresponds to the ideal \((D/H)/(D)\). Thus the elements of the fiber have the unique form
\[
\left(1+\frac DH b\right)(\sigma)\alpha,\qquad b\bmod H.
\]
Their annihilators are \(D/\gcd(D,1+(D/H)b)\). We therefore count the parameters for which \(1+(D/H)b\) is a unit modulo \(D\).
For \(P\mid D/H\), this expression is always one modulo \(P\). The remaining factors are precisely those satisfying \(P\mid H\) and \(P^2\nmid D\). For each of them, \(D/H\) is invertible modulo \(P\), so exactly one residue of \(b\bmod P\) is excluded. The Chinese remainder theorem for squarefree \(H\) multiplies these counts. When \(H=1\), the empty product is one, in agreement with \(W=\{0\}\).
\end{proof}
The full annihilator distribution along an arbitrary invariant direction appears in Theorem~\ref{v107:thm-affine-annihilator-distribution}. The squarefree case here first exhibits how multiplicity protects generation.

\begin{example}[The same direction in degrees three and six]\label{v106:ex-affine-f8}
In \(\F_8\), choose \(v^3+v^2+1=0\), and put \(P=T^2+T+1\). Since \(v+v^2+v^4=1\) and \(v\) has degree three, \(v\) is normal. Direct calculation gives
\[
\begin{aligned}
W&=\ker P(\sigma)=\{0,v+1,v^2+1,v^2+v\},\\
v+W&=\{v,1,v^4,v^2\}.
\end{aligned}
\]
The factor \(P\) occurs only once in \(T^3-1\), so exactly three elements are normal; the unique exception is \(1\). After embedding \(\F_8\) in \(\F_{64}\), the same \(W\) is still \(\ker P(\sigma)\), whereas \(T^6-1=(T+1)^2P^2\). For every normal element \(b\) of \(\F_{64}/\F_2\), all four elements of \(b+W\) are therefore normal. Example~\ref{v106:ex-complete-f64} will show that they need not be normal over every intermediate field.
\end{example}
Multiplicity influences both diagonalizability and the preservation of full annihilator data under affine perturbations.

\section{Linearized Root Spaces with Base-Field Coefficients}
\begin{theorem}\label{thm:root-module-count}
Let \(h\in k[T]\) be monic with \(h(0)\ne0\), let \(W\) be the root space of \(\mathcal L_h\), and let \(s=\deg h\).
There exists \(\eta\in W\) such that
\[
k[T]/(h)\longrightarrow W,\qquad \bar a\longmapsto a(\sigma)\eta
\]
is a linear isomorphism carrying multiplication by \(T\) to Frobenius.
For each monic \(d\mid h\), exactly \(\varphi_q(d)\) elements of \(W\) satisfy \(\mu_x=d\).
If \(x\ne0\), its degree over the base field is \(\ord(\mu_x)\).
\end{theorem}
\begin{proof}
Since \(\mathcal L_h\) has no repeated roots, \(\dim_kW=s\). If a nonzero \(g\) of degree less than \(s\) annihilated the entire space,
then \(\mathcal L_g\) would have \(q^s\) roots but degree less than \(q^s\), which is impossible. Thus the minimal polynomial of \(\sigma|_W\) is exactly \(h\).
The cyclic vector lemma gives \(\eta\). As for normal bases, \(\mu_{a(\sigma)\eta}=h/\gcd(h,a)\).
Requiring this polynomial to be \(d\) is equivalent to \(a=(h/d)u\) and \(\gcd(u,d)=1\), with \(\varphi_q(d)\) choices modulo \(d\).
Finally, \(\sigma^j(x)=x\) is equivalent to \(\mu_x\mid T^j-1\), so the least period is \(\ord(\mu_x)\).
\end{proof}

\begin{corollary}[The splitting field of a linearized root space]\label{cor:linearized-splitting-field}
Let \(f\in k[T]\) be monic with \(\deg f\ge1\) and \(f(0)\ne0\), and put \(N=\ord(f)\).
The splitting field of \(\mathcal L_f\) over \(k=\F_q\) is exactly \(\F_{q^N}\).
\end{corollary}
\begin{proof}
By \(N=\ord(f)\), write \(T^N-1=fh\). Linearization with base-field coefficients takes multiplication to composition, so
\[
X^{q^N}-X=\mathcal L_h\bigl(\mathcal L_f(X)\bigr).
\]
Consequently, all roots of \(\mathcal L_f\) belong to \(\F_{q^N}\).
On the other hand, Theorem~\ref{thm:root-module-count} gives a root \(\eta\) whose Frobenius annihilator polynomial is exactly \(f\).
Thus \(\sigma^t\eta=\eta\) is equivalent to \(f\mid T^t-1\), with least positive period \(N\).
The orbit theorem gives \([k(\eta):k]=N\), so \(k(\eta)=\F_{q^N}\).
Every field containing all the roots must contain \(\eta\); the splitting field can therefore be no smaller.
\end{proof}

Here \(\deg f\) determines the dimension of the root space, whereas \(\ord(f)\) determines the least extension degree needed to contain all roots;
these numbers describe linear relations and iteration periods, respectively.

\begin{corollary}[Kernels of a fixed linearized polynomial over all extensions]\label{v70:cor-all-extension-kernels}
Let \(k=\F_q\) have characteristic \(p\), and fix a nonzero polynomial \(h\in k[T]\). Write
\[
 h=cT^a h_0,\qquad c\in k^\times,\quad a\ge0,\quad
 h_0\text{ monic},\quad h_0(0)\ne0.
\]
If \(h_0=1\), then \(\mathcal L_h\) permutes every finite extension \(\F_{q^m}\).
If \(h_0\ne1\), put \(N=\ord(h_0)\), and let \(W\) be the full root space of \(\mathcal L_{h_0}\) in an algebraic closure. Then, for every integer \(m\ge1\),
\begin{equation}\label{v70:eq-all-extension-root-space}
 \ker\bigl(\mathcal L_h|_{\F_{q^m}}\bigr)
 =W\cap\F_{q^m}
 =W\cap\F_{q^{\gcd(m,N)}}.
\end{equation}
Factor
\[
 h_0=\prod_{i=1}^s P_i^{e_i},\qquad o_i=\ord(P_i),
\]
where the \(P_i\) are pairwise distinct monic irreducible polynomials and \(e_i\ge1\). Then
\begin{equation}\label{v70:eq-all-extension-kernel-dimension}
 d_m:=\dim_k\ker\bigl(\mathcal L_h|_{\F_{q^m}}\bigr)
 =\sum_{\substack{1\le i\le s\\o_i\mid m}}
     \deg P_i\min\{e_i,p^{v_p(m)}\}.
\end{equation}
\end{corollary}
\begin{proof}
Base-field coefficients satisfy \(b^{q^a}=b\), so
\[
 \mathcal L_h(X)=c\bigl(\mathcal L_{h_0}(X)\bigr)^{q^a}.
\]
Multiplication by \(c\) and the \(q^a\)th-power map are bijections on every finite extension, so they leave the kernel unchanged. When \(h_0=1\), this proves the permutation assertion.
Otherwise, Corollary~\ref{cor:linearized-splitting-field} gives \(W\subseteq\F_{q^N}\). The intersection formula for finite fields now gives~\eqref{v70:eq-all-extension-root-space}.

By Theorem~\ref{thm:linearized-gcd}, we have \(d_m=\deg\gcd(h_0,T^m-1)\). Write \(m=p^v u\), where \(p\nmid u\). Then
\[
 T^m-1=(T^u-1)^{p^v},
\]
and \(T^u-1\) has no repeated roots. Since \(o_i\) is coprime to \(p\), the condition \(P_i\mid T^u-1\) is equivalent to \(o_i\mid u\), and hence to \(o_i\mid m\).
The multiplicity of \(P_i\) in the greatest common divisor is therefore \(\min\{e_i,p^v\}\) or zero. Taking degrees proves~\eqref{v70:eq-all-extension-kernel-dimension}.
\end{proof}

Here \(h\) is fixed and the extension degree \(m\) varies; all kernels are controlled by a single finite root space and its intersections with subfields.
We next fix the base field and vary the number of compositions, studying the growth of the full root space and its splitting field.

\begin{corollary}[Splitting fields of linearized iterates]\label{v68:cor-linearized-iterate-splitting}
Let \(k=\F_q\) have characteristic \(p\), and let \(f\in k[T]\) be monic and nonconstant, with \(f(0)\ne0\).
Write
\[
f=\prod_{i=1}^s P_i^{e_i},\qquad
f_0=\prod_{i=1}^s P_i,\qquad
e=\max_i e_i,\qquad N=\ord(f_0),
\]
where the \(P_i\) are distinct monic irreducible polynomials. Put \(A(X)=\mathcal L_f(X)\).
For every integer \(r\ge1\), the splitting field of the \(r\)-fold composition \(A^{\circ r}\) over \(k\) is
\begin{equation}\label{v68:eq-linearized-iterate-splitting}
\F_{q^{N p^{b_r}}},\qquad
b_r=\min\{b\in\mathbb Z_{\ge0}:p^b\ge er\}=\lceil\log_p(er)\rceil.
\end{equation}
Its root space has dimension \(r\deg f\) over \(k\).
\end{corollary}
\begin{proof}
The coefficients lie in \(k\), so multiplication of operator polynomials corresponds to composition of linearized polynomials:
\[
A^{\circ r}=\mathcal L_{f^r},\qquad
f^r=\prod_{i=1}^s P_i^{re_i}.
\]
Corollary~\ref{cor:linearized-splitting-field} identifies the splitting-field degree with \(\ord(f^r)\).
The square-free part of \(f^r\) is \(f_0\), and its largest factor multiplicity is \(er\).
The order formula in Theorem~\ref{thm:order-multiplicity} therefore gives \(\ord(f^r)=Np^{b_r}\).
Finally, the coefficient of \(X\) in \(\mathcal L_{f^r}\) is \(f(0)^r\ne0\), so all its roots are simple.
Its ordinary degree is \(q^{r\deg f}\); hence its root space has \(q^{r\deg f}\) elements and dimension \(r\deg f\).
\end{proof}

Thus the root-space dimension grows linearly with the number of compositions, whereas the splitting-field degree changes only when the largest factor multiplicity crosses a power of the characteristic.
This gives the base-field coefficient case of the splitting-field growth studied in~\cite[Theorem~4.1]{Reis2025Iterating}.
The multiplicities of the irreducible factors of \(f^r\) describe the structure of the Frobenius operator; they are not the root multiplicities of \(\mathcal L_{f^r}(X)\) as an ordinary polynomial.
Its coefficient of \(X\) is nonzero, so all its roots remain simple.

\begin{example}[Different root spaces in the same splitting field]\label{v68:ex-linearized-iterate-f2}
Take \(k=\F_2\), \(f(T)=T^2+T+1\), and
\[
A(X)=X^4+X^2+X.
\]
The polynomial \(f\) is irreducible and \(\ord(f)=3\), so \(e=1\) and \(N=3\).
Formula~\eqref{v68:eq-linearized-iterate-splitting} gives
\[
\begin{array}{c|ccccc}
r&1&2&3&4&5\\\hline
\dim_{\F_2}\ker A^{\circ r}&2&4&6&8&10\\
\text{Splitting-field degree}&3&6&12&12&24
\end{array}
\]
Here each kernel is taken in an algebraic closure.
In particular, \(A^{\circ3}\) and \(A^{\circ4}\) have respectively \(2^6\) and \(2^8\) distinct roots, but both have splitting field \(\F_{2^{12}}\).
The root spaces are nested, since \(A(0)=0\); their dimensions may therefore increase strictly while their splitting field remains unchanged.

Alternatively, fix \(A\) and \(A^{\circ3}=\mathcal L_{f^3}\), and vary the ambient field \(\F_{2^m}\). Corollary~\ref{v70:cor-all-extension-kernels} gives
\[
\begin{array}{c|cc}
 &3\nmid m&3\mid m\\\hline
\dim_{\F_2}\ker\bigl(A|_{\F_{2^m}}\bigr)&0&2\\
\dim_{\F_2}\ker\bigl(A^{\circ3}|_{\F_{2^m}}\bigr)&0&2\min\{3,2^{v_2(m)}\}
\end{array}
\]
Thus they permute exactly the same extensions, although their kernels on other extensions may differ.
Exercise~\ref{v70:ex-linearized-values-all-extensions} translates these kernel dimensions into value-set sizes and counts of ordered pairs with the same image.
\end{example}

\leavevmode\phantomsection\label{def:linear-primitive-root}\ignorespaces A root satisfying \(\mu_\eta=h\) is called a \textbf{linearly primitive root} of \(\mathcal L_h\), or a \(q\)-primitive root; when \(h(0)\ne0\), there are \(\varphi_q(h)\) such roots.
If \(T\mid h\), then \(\mathcal L_h\) has repeated roots, and the Frobenius annihilator polynomial of every root lacks the factor \(T\); there are therefore no such roots.
Linear and multiplicative primitivity are distinct generation properties; normal elements are precisely the linearly primitive roots of \(\mathcal L_h\) when \(h=T^n-1\).

\begin{corollary}\leavevmode\phantomsection\label{cor:linearized-primitive-irreducible}
Let \(h\in k[T]\) be monic irreducible, with \(h(0)\ne0\) and \(\deg h=s\).
Every irreducible factor of \(\mathcal L_h(X)/X\) has degree \(\ord(h)\).
The polynomial \(\mathcal L_h(X)/X\) is irreducible if and only if \(h\) is a primitive polynomial.
\end{corollary}
\begin{proof}
The Frobenius annihilator polynomial of every nonzero root must be \(h\), so every nonzero orbit has length \(\ord(h)\).
There are \(q^s-1\) nonzero roots, which form a single orbit precisely when \(\ord(h)=q^s-1\).
\end{proof}

\begin{theorem}[Realizing an invertible linear operator as Frobenius]\label{thm:frobenius-realization}
Let \(V\ne0\) be a finite-dimensional \(k=\F_q\)-space, let \(A\in\operatorname{GL}_k(V)\), and let \(m_A\) be its minimal polynomial.
The following conditions are equivalent:
\begin{enumerate}
\item There are a finite extension \(L/k\) and an injective \(k\)-linear map \(\iota:V\to L\) such that \(\iota(Av)=\iota(v)^q\) for every \(v\in V\);
\item \(\deg m_A=\dim_k V\);
\item \(A\) has a cyclic vector, that is, there exists \(v\) such that \(v,Av,\ldots,A^{\dim V-1}v\) form a basis of \(V\).
\end{enumerate}
When these conditions hold, one may take \(\iota(V)\) to be the full root space of \(\mathcal L_{m_A}\).
\end{theorem}
\begin{proof}
Put \(s=\dim_kV\). First prove that the first condition implies the second. Write \(W=\iota(V)\). Then \(W\) is a \(k\)-subspace of \(L\) of dimension \(s\) satisfying \(\sigma(W)=W\). The theorem on subspace polynomials and its corollary on coefficient descent give
\[
L_W(X)=\prod_{w\in W}(X-w)=\mathcal L_h(X)
\]
for some monic \(h\in k[T]\) with \(\deg h=s\) and \(h(0)\ne0\).
Equivariance gives \(h(A)=0\), hence \(m_A\mid h\) and therefore \(\deg m_A\le s\).
If \(\deg m_A<s\), the nonzero polynomial \(\mathcal L_{m_A}\) would have ordinary degree less than \(q^s\)
but vanish at all \(q^s\) elements of \(W\), contradicting the root bound. Thus \(\deg m_A=s\), and \(m_A=h\).

The cyclic vector lemma~\ref{lem:cyclic-vector} proves that the second condition implies the third.
Conversely, if \(v,Av,\ldots,A^{s-1}v\) is a basis, no nonzero polynomial of degree less than \(s\) annihilates \(v\).
Expanding \(A^sv\) in this basis gives a monic polynomial \(h\) of degree \(s\) with \(h(A)v=0\).
Since \(h(A)\) commutes with \(A\), it vanishes on every basis vector \(A^iv\), so \(h(A)=0\).
Thus \(m_A=h\), proving the equivalence of the second and third conditions.

Finally, the second and third conditions imply the first. Invertibility of \(A\) implies \(m_A(0)\ne0\): otherwise \(m_A=Tg\) would give \(A g(A)=0\) and hence \(g(A)=0\), contradicting minimality.
Take the full root space \(W\) in a finite splitting extension \(L/k\) of \(\mathcal L_{m_A}\).
Theorem~\ref{thm:root-module-count} gives a Frobenius-compatible isomorphism \(k[T]/(m_A)\simeq W\). On the other hand, a cyclic vector \(v\) gives \(k[T]/(m_A)\simeq V\), \(\bar a\mapsto a(A)v\).
Compose the inverse of the second isomorphism with the first and then with \(W\subseteq L\) to obtain the required \(\iota\).
\end{proof}

\begin{example}[An identity action that cannot be realized]\label{exm:noncyclic-realization}
The identity operator on a two-dimensional space \(k^2\) cannot be realized this way. An equivariant injection would have every element of its image satisfying \(x^q=x\), so the image would lie in \(k\), contradicting its dimension two.
In general, the root polynomial of a Frobenius-stable subspace attains the root bound, forcing the degree of the minimal polynomial of the restricted Frobenius operator to equal the dimension of the space.
This is why the realization criterion requires cyclicity, not merely invertibility.
\end{example}

\subsection{Coefficient extension and Frobenius return on root spaces}
\label{v145:sec-skew-roots}
With base-field coefficients, composition factors are described by factors of an ordinary operator polynomial. After extending the coefficient field, root spaces still record composition order, but the Frobenius step fixing the coefficients changes.

Let $Q$ be a prime power, $K=\F_{Q^s}$, and $\sigma(c)=c^Q$. The \textbf{skew polynomial ring} $K[t;\sigma]$ is defined by $tc=\sigma(c)t$, with coefficients written on the left. The correspondence
\[
 f(t)=\sum_i c_it^i\quad\longleftrightarrow\quad
 L_f(X)=\sum_i c_iX^{Q^i}
\]
gives $L_{gh}=L_g\circ L_h$, so a right factor is the linearized polynomial applied first. Skew irreducibility means that no factorization into two factors of positive degree exists in $K[t;\sigma]$.

\begin{theorem}[Right factors and return-invariant subspaces]\label{v145:thm-skew-roots}
Let $f\in K[t;\sigma]$ be monic of degree $r\ge1$, with $f(0)\ne0$, and put
\[
 V_f=\ker L_f\subset\overline{\F}_Q,
 \qquad \Theta(x)=x^{Q^s}.
\]
Then $\dim_{\F_Q}V_f=r$, and there is an order-preserving bijection
\begin{equation}\label{v145:eq-right-factor-lattice}
 \begin{gathered}
 \{\text{monic right factors of $f$}\}
 \longleftrightarrow
 \{\text{$\Theta$-invariant $\F_Q$-subspaces of $V_f$}\},\\
 h\longmapsto\ker L_h.
 \end{gathered}
\end{equation}
Factor degree equals subspace dimension, and right divisibility corresponds to inclusion. In particular,
\begin{equation}\label{v145:eq-skew-criterion}
\begin{aligned}
 f\text{ is skew irreducible}
 &\iff V_f\text{ has no nonzero proper $\Theta$-invariant subspace}\\
 &\iff\chi_{\Theta|V_f}(Z)\text{ is irreducible in }\F_Q[Z].
\end{aligned}
\end{equation}
\end{theorem}
\begin{proof}
Since $L_f'=f(0)\ne0$, its $Q^r$ distinct roots form an $r$-dimensional $\F_Q$-space. Its coefficients lie in $K$, so $\Theta$ preserves this space. If $f=gh$, both factors have nonzero constant term, and composition gives the $\Theta$-invariant subspace $\ker L_h\subseteq V_f$.

Conversely, for such a subspace $W$, Theorem~\ref{thm:subspace-polynomial} and coefficient descent give the monic linearized polynomial
\[
 L_W(X)=\prod_{w\in W}(X-w)\in K[X].
\]
Let $h$ be the corresponding skew polynomial. Right composition division, Theorem~\ref{thm:composition-division}, gives $f=gh+R$, where $R=0$ or $\deg R<\dim W$. Now $L_R$ vanishes on $W$. If $R\ne0$, it has at most $Q^{\deg R}<|W|$ roots, a contradiction. Thus $h$ right-divides $f$. A monic separable polynomial is determined by its root set, so the constructions are inverse. The same division argument identifies right divisibility with inclusion.

If no nonzero proper invariant subspace exists, every nonzero vector generates the whole space as a cyclic module. Thus $V_f\simeq\F_Q[Z]/(g)$ with $g$ irreducible, and $g=\chi_\Theta$. Conversely, if $\chi_\Theta$ is irreducible, the monic annihilator of every nonzero vector must equal $\chi_\Theta$, so its cyclic submodule already has dimension $r$. This proves the last equivalence.
\end{proof}

\begin{proposition}[Moore matrices and the return characteristic polynomial]\label{v145:prop-return-matrix}
Write $f=t^r+c_{r-1}t^{r-1}+\cdots+c_0$. Let $A_f$ be the column companion matrix: its first $r-1$ columns are $e_2,\ldots,e_r$, and its last column is $(-c_0,\ldots,-c_{r-1})^{\mathsf T}$. Then
\begin{equation}\label{v145:eq-return-polynomial}
 \chi_{\Theta|V_f}(Z)
 =\det\!\left(ZI_r-A_f\sigma(A_f)\cdots\sigma^{s-1}(A_f)\right),
\end{equation}
where $\sigma$ acts entrywise. The right side belongs to $\F_Q[Z]$ and is called the \textbf{return characteristic polynomial} of $f$.
\end{proposition}
\begin{proof}
Choose an $\F_Q$-basis $u_1,\ldots,u_r$ of $V_f$, put $u=(u_1,\ldots,u_r)^{\mathsf T}$, and write $F(u)=u^Q$. The Moore matrix $B=[u,Fu,\ldots,F^{r-1}u]$ is invertible, and the root relation gives $F(B)=BA_f$. Choose $C\in\operatorname{GL}_r(\F_Q)$ with $F^s(u)=Cu$; the matrix of $\Theta$ in the chosen basis is $C^{\mathsf T}$. On the other hand,
\[
 CB=F^s(B)=B A_f\sigma(A_f)\cdots\sigma^{s-1}(A_f).
\]
The product matrix is therefore similar to $C$, which has the same characteristic polynomial as $C^{\mathsf T}$.
\end{proof}

When the return characteristic polynomial is square-free, its quotient algebra is a finite étale algebra as in Chapter~3; when it is irreducible, that algebra is a field. The factor criterion has become irreducibility of a commutative polynomial, while the Frobenius order in the matrix product remains essential.

\begin{example}[Two factor structures for the same coefficient expression]\label{v145:ex-skew-four}
Take $K=\F_4$ and $\sigma(c)=c^2$. For $f=t^2+t+1$,
\[
 A_f=\begin{pmatrix}0&1\\1&1\end{pmatrix},\qquad
 A_f\sigma(A_f)=\begin{pmatrix}1&1\\1&0\end{pmatrix}.
\]
The return characteristic polynomial is $Z^2+Z+1$, irreducible over $\F_2$. Thus $f$ is irreducible in $K[t;\sigma]$, although the same coefficient expression splits in the ordinary polynomial ring $K[t]$.

For $t^2+1$, the root space is $\F_4$ and the return map is the identity, with characteristic polynomial $(Z-1)^2$. The three $\F_2$-lines give three monic right factors $t-c$, for $c\in\F_4^\times$; directly,
\[
 (t+c^2)(t+c)=t^2+1.
\]
The return characteristic polynomial detects irreducibility, whereas the full lattice of right factors is recorded by the invariant subspaces of the return operator.
Here the identity is a second Frobenius return; the realization problem in Example~\ref{exm:noncyclic-realization} uses a single Frobenius step.
\end{example}

\section{Return norms, common subfields, and fibers}
\label{v148:sec-norm-fibers}
The return characteristic polynomial converts right factorization into factorization in a commutative polynomial ring. Conversely, how many skew polynomials realize a prescribed return polynomial? If the constant coefficient is also prescribed, does the answer still depend only on degrees? Local modules and common subfields govern these two questions.

Let $K=\F_{q^k}$, $\tau(x)=x^{q^s}$, $\gcd(s,k)=1$, and $R=K[t;\tau]$, with central variable $Z=t^k$. For a monic $f\in R$ of degree $n$, use the column companion matrix $A_f$ and define its \textbf{return norm} by
\[
 \mathcal N(f)=\det\bigl(ZI_n-A_f\tau(A_f)\cdots\tau^{k-1}(A_f)\bigr).
\]
This is the monic normalized reduced norm. When $\tau(x)=x^q$, it is the return characteristic polynomial of Proposition~\ref{v145:prop-return-matrix}. Cyclically permuting a matrix product preserves its characteristic polynomial, so the coefficients lie in $\F_q$, also when $f(0)=0$.

\begin{proposition}[Products and constant coefficients]\label{v148:prop-norm-basic}
For monic skew polynomials $f,g$,
\begin{equation}\label{v148:eq-norm-basic}
 \mathcal N(fg)=\mathcal N(f)\mathcal N(g),\qquad
 \deg_Z\mathcal N(f)=\deg_t f,\qquad \mathcal N(t)=Z.
\end{equation}
If $\deg f=n$ and $f(0)=c$, then
\[
 \mathcal N(f)(0)=(-1)^{n(k+1)}\Nm_{K/\F_q}(c).
\]
\end{proposition}
\begin{proof}
Left multiplication by $t$ on $R/Rf$ is the semilinear operator $A_f\tau$; its $k$th power has the matrix in the definition. For $fg$, use the exact sequence
\[
 0\longrightarrow R/Rf\xrightarrow{\ u\mapsto ug\ }R/R(fg)
 \longrightarrow R/Rg\longrightarrow0.
\]
In a basis adapted to the subspace, the matrix of $t^k$ is block upper triangular, so its characteristic polynomial is the product of the two diagonal characteristic polynomials. The constant coefficient follows from $\det A_f=(-1)^nc$.
\end{proof}

For monic $h\in\F_q[Z]$, put
\[
 \begin{aligned}
 U_k(h)&=\#\{f\in R\text{ monic}:\mathcal N(f)=h\},\\
 F_{h,k}(c)&=\#\{f\in R\text{ monic}:\mathcal N(f)=h,\ f(0)=c\}.
 \end{aligned}
\]
We suppress $\tau$ from the notation, writing $F_{h,k}^{\tau}$ when the generator varies.

\subsection{Local modules and unrestricted fibers}
\begin{theorem}[All return-norm fibers]\label{v148:thm-all-fibers}
Write
\[
 h=Z^{e_0}\prod_i\pi_i^{e_i},\qquad
 \pi_i\ne Z\text{ monic irreducible},\quad d_i=\deg\pi_i.
\]
Then
\begin{equation}\label{v148:eq-all-fibers}
 U_k(h)=\prod_i{e_i+k-1\brack k-1}_{q^{d_i}},\qquad
 {e+k-1\brack k-1}_Q=\prod_{j=1}^e\frac{Q^{k+j-1}-1}{Q^j-1}.
\end{equation}
In particular, for irreducible $h\ne Z$ of degree $n$,
\[
 U_k(h)=\frac{q^{kn}-1}{q^n-1}.
\]
Moreover, $f$ is skew irreducible if and only if $\mathcal N(f)$ is irreducible.
\end{theorem}
\begin{proof}
First suppose $h(0)\ne0$. Each monic $f$ determines a unique left ideal $Rf$. Cayley--Hamilton gives $Rh(t^k)\subseteq Rf$, allowing the count to be made in a finite quotient algebra.

We describe the matrix model at a nonzero central prime. Put $A=\F_q[Z]/(\pi^a)$, where $\pi\ne Z$, and $E=K\otimes_{\F_q}A$. The cyclic norm $E^\times\to A^\times$ is surjective. On residue fields this is finite-field norm surjectivity; at successive precisions,
\[
 \prod_{i=0}^{k-1}\tau^i(1+\pi^jx)
 \equiv1+\pi^j\sum_{i=0}^{k-1}\tau^i(x)\pmod{\pi^{j+1}},
\]
where the corresponding trace is surjective. Choose $b$ of norm $Z$ and let $t$ act as $b\tau$. This gives
\[
 R/(\pi(t^k)^a)\simeq\End_A(E)\simeq M_k(A).
\]
To verify surjectivity, reduce modulo $\pi$ and extend to an algebraic closure: $E$ splits into $k$ components, multiplication gives all diagonal matrices, and $b\tau$ cyclically permutes the components. These generate every matrix unit. Both sides are free of rank $k^2$ over $A$, so the isomorphism lifts to $A$.

In the finite matrix model above, a left ideal is determined by the common submodule formed by its rows. Pulling this submodule back over the local principal ideal ring $\mathcal O=\F_q[Z]_{(\pi)}$ gives a lattice $\pi^a\mathcal O^k\subseteq L\subseteq\mathcal O^k$. Put $M=\mathcal O^k/L$, of colength $e$. The quotient of the matrix left module, viewed as a central module, is $M^{\oplus k}$. Every composition factor of $M$ is $\F_q[Z]/(\pi)$, so the $\F_q$-characteristic polynomial of $Z$ on this quotient is $\pi^{ke}$.
On the other hand, restriction of scalars from $K$ to $\F_q$ on $R/Rf$ gives the characteristic polynomial $\mathcal N(f)^k$ for $t^k$, since $\mathcal N(f)$ has coefficients in $\F_q$. Comparing the exponents of $\pi$ identifies its multiplicity in the return norm with $e$. A lattice of colength $e$ contains $\pi^e\mathcal O^k$, so the finite quotient model suffices.

Let $L_k(e;Q)$ count these lattices, where $Q=q^{\deg\pi}$. Projection onto the last coordinate has image $\pi^a\mathcal O$; the intersection with the first $k-1$ coordinates has colength $e-a$. There are $Q^{e-a}$ lifts of the last-coordinate generator. Hence
\[
 L_k(e;Q)=\sum_{a=0}^e Q^{e-a}L_{k-1}(e-a;Q),\qquad
 \sum_{e\ge0}L_k(e;Q)u^e=\prod_{j=0}^{k-1}(1-Q^ju)^{-1}.
\]
Its coefficients give~\eqref{v148:eq-all-fibers}. The Chinese remainder theorem combines distinct primes independently. If $h=Z^{e_0}h_1$ with $h_1(0)\ne0$, there is a unique expression $f=g t^{e_0}$ with $\mathcal N(g)=h_1$; the zero prime introduces no further choice.

Finally, $f$ is irreducible precisely when $R/Rf$ is simple. At a nonzero prime, the matrix model makes this equivalent to one irreducible central factor of colength one. At the zero prime, the only simple case is $f=t$. This proves the irreducibility criterion as well.
\end{proof}

\begin{proposition}[Local types within a fiber]\label{v148:prop-local-types}
Let $\pi\ne Z$ be a monic irreducible central factor. At $\pi^e$, the quotient type is a partition $\lambda\vdash e$ with $\ell(\lambda)\le k$. Let $\lambda'$ be the conjugate partition, let $m_a(\lambda)$ be the multiplicity of the part $a$, and put $Q=q^{\deg\pi}$. The number of this type is
\begin{equation}\label{v148:eq-local-type}
 T_k(\lambda;Q)=
 \frac{Q^{ke-\sum_j(\lambda'_j)^2}
       \prod_{j=0}^{\ell(\lambda)-1}(1-Q^{j-k})}
      {\prod_{a\ge1}\prod_{j=1}^{m_a(\lambda)}(1-Q^{-j})}.
\end{equation}
Multiplication over the distinct primes counts a prescribed collection of local types.
\end{proposition}
\begin{proof}
Put $M_\lambda=\bigoplus_i\mathcal O/(\pi^{\lambda_i})$. The desired lattices are kernels of surjections $\mathcal O^k\to M_\lambda$. A tuple generates $M_\lambda$ exactly when its reduction generates the residual vector space: successive subtraction of lifted linear combinations proves sufficiency. Thus the number of surjections is
$Q^{ke}\prod_{j=0}^{\ell(\lambda)-1}(1-Q^{j-k})$.
There are $Q^{\min(a,b)}$ homomorphisms between cyclic summands of lengths $a,b$, so the endomorphism ring has size $Q^{\sum_j(\lambda'_j)^2}$. Its invertible fraction is
$\prod_a\prod_{j=1}^{m_a(\lambda)}(1-Q^{-j})$, from the full-rank condition on equal-length summands. Dividing by the automorphism count proves the formula.
\end{proof}

\subsection{Common subfields and a fixed constant coefficient}
\begin{theorem}[Relative-norm symmetry of fibers]\label{v148:thm-relative-fibers}
Suppose $h(0)\ne0$, $\deg h=n\ge1$, and $g=\gcd(n,k)$, and take $c,c'\in K^\times$. If
\[
 \Nm_{K/\F_{q^g}}(c)=\Nm_{K/\F_{q^g}}(c'),
\]
then $F_{h,k}(c)=F_{h,k}(c')$. This equality persists after fixing coefficient support or local module type. If $\gcd(n,k)=1$, then
\begin{equation}\label{v148:eq-coprime-fiber}
 F_{h,k}(c)=\frac{q-1}{q^k-1}U_k(h)
 \quad\text{when}\quad
 \Nm_{K/\F_q}(c)=(-1)^{n(k+1)}h(0),
\end{equation}
and the fiber is empty otherwise. In particular, the nonempty fibers of irreducible targets have size
\begin{equation}\label{v148:eq-Cnk}
 C_{n,k}(q)=\frac{(q^{kn}-1)(q-1)}{(q^n-1)(q^k-1)}.
\end{equation}
\end{theorem}
\begin{proof}
For $v\in K^\times$, use the monic transformation
\[
 f^{[v]}=\tau^n(v)^{-1}fv,\qquad
 a_i\longmapsto\tau^n(v)^{-1}a_i\tau^i(v).
\]
It preserves the return norm, support, and skew-module isomorphism type, while replacing the constant coefficient by $cv/\tau^n(v)$. Finite-field Hilbert 90 gives
\[
 \{v/\tau^n(v):v\in K^\times\}
 =\ker\Nm_{K/\F_{q^g}}.
\]
Thus constants with the same relative norm are connected by such bijections. When $g=1$, all $(q^k-1)/(q-1)$ compatible constants have equal fibers; summing proves the formula.
\end{proof}

In general, a fixed norm to the base field leaves at most $(q^g-1)/(q-1)$ relative norm values to distinguish. A trivial common subfield gives uniformity; otherwise the residual problem lies in this smaller multiplicative group. Fixing support preserves equality between constants for the same target, but need not give equality between targets. For instance, in $\F_4[t;\bar{\phantom a}]$,
\[
 \mathcal N(t^3+at+1)=Z^3+\Tr(a)Z^2+\Nm(a)Z+1.
\]
The irreducible targets $Z^3+Z+1$ and $Z^3+Z^2+1$ have fibers of sizes one and zero on this support.

Choose $c_0$ in a compatible norm coset and put
\[
 H_0=\ker\Nm_{K/\F_q},\qquad H_g=\ker\Nm_{K/\F_{q^g}}.
\]
The function $u\mapsto F_{h,k}(c_0u)$ descends to $H_0/H_g$. Its multiplicative Fourier coefficients are therefore supported on characters trivial on $H_g$: translation along that subgroup kills every other character sum. Thus the residual nonuniformity has only $(q^g-1)/(q-1)$ parameters, rather than all of $K^\times$.

\begin{corollary}[A prescribed constant and root order]\label{v148:cor-fixed-order}
Let $n\ge2$, $\gcd(n,k)=1$, $c\ne0$, and $\nu=(-1)^{nk}\Nm(c)$. The number of skew irreducible polynomials of degree $n$ and constant coefficient $c$ is
\begin{equation}\label{v148:eq-fixed-irred}
 \frac{C_{n,k}(q)}{n(q-1)}
 \sum_{d\mid n}\mu(d)(q^{n/d}-1)R_d(\nu),
\end{equation}
where $R_d(\nu)=\#\{u\in\F_q^\times:u^d=\nu\}$.
For $\ell\mid q^n-1$ with $\ord_\ell(q)=n$, put
$d_\ell=\ell/\gcd(\ell,(q^n-1)/(q-1))$. The number with return-root order exactly $\ell$ is
\[
 C_{n,k}(q)\frac{\varphi(\ell)}{n\varphi(d_\ell)}
 \quad\text{if }\ord(\nu)=d_\ell,
\]
and zero otherwise. In particular, the number whose return norm is primitive is
\[
 C_{n,k}(q)\frac{\varphi(q^n-1)}{n\varphi(q-1)}
 \quad\text{if }\nu\text{ is primitive},
\]
and zero otherwise.
\end{corollary}
\begin{proof}
An irreducible target root $\alpha$ satisfies $\Nm_{\F_{q^n}/\F_q}(\alpha)=\nu$. Apply Möbius inversion over subfields. If $\alpha\in\F_{q^{n/d}}$, its norm from that smaller field has $d$th power $\nu$, giving $(q^{n/d}-1)R_d(\nu)/(q-1)$ choices. Divide by the orbit length $n$ and multiply by the common fiber size.

The norm image of an element of order $\ell$ has order $d_\ell$. Generators of a cyclic group have equally many preimages over each generator of the image, namely $\varphi(\ell)/\varphi(d_\ell)$. Here $R_d(\nu)=\gcd(d,q-1)$ if $\nu^{(q-1)/\gcd(d,q-1)}=1$, and is zero otherwise.
\end{proof}

\subsection{Return orders and ordinary factorization}
A monic $f$ is \textbf{norm-primitive} if $\mathcal N(f)$ is primitive. For skew irreducible $f$, the endomorphism field of the simple module $R/Rf$ is $\F_{q^n}$, and central $Z$ corresponds to a root $\alpha$. Norm-primitivity is therefore equivalent to the action of $Z$ having order $q^n-1$, or to the least positive $\ell$ in the central return relation
\[
 f\text{ right divides }t^{k\ell}-1
\]
being $q^n-1$. These relations have exponents divisible by $k$. Without fixing the constant coefficient, the norm-primitive count is
\[
 \frac{q^{kn}-1}{q^n-1}\frac{\varphi(q^n-1)}n.
\]
If the return-root order is $\ell$, the semilinear operator $A_f\tau$ has order $k\ell$: an identity power must first eliminate the semilinear part of order $k$, and then the return operator of order $\ell$.

Now take $\tau(x)=x^q$, $\deg f=n\ge1$, and $c=f(0)\ne0$. Write
\[
 S_j=1+q+\cdots+q^{j-1},\qquad
 P_f(X)=\frac{L_f(X)}X,\qquad
 Q_f(Y)=Y^{S_n}+\sum_{i=1}^{n-1}a_iY^{S_i}+c.
\]
Then $P_f(X)=Q_f(X^{q-1})$. Their roots respectively represent nonzero vectors and $\F_q$-lines in the root space.

\begin{theorem}[Return orbits determine ordinary factorization]\label{v148:thm-ordinary-return}
All irreducible factorizations below are over the coefficient field $K=\F_{q^k}$.
If $\mathcal N(f)$ is irreducible with root order $\ell$, then $P_f$ has $(q^n-1)/\ell$ distinct irreducible factors of degree $\ell$. Every irreducible factor of $Q_f$ has degree
\[
 e=\frac{\ell}{\gcd(\ell,q-1)}.
\]
In particular, $P_f$ is irreducible if and only if $f$ is norm-primitive; $Q_f$ is irreducible if and only if the return action is transitive on lines.

More generally, let the elementary divisors of the return operator be $\pi^{\lambda_{\pi,j}}$, and put
\[
 d_s(a)=\sum_\pi\deg\pi\sum_j
 \min\{\lambda_{\pi,j},v_\pi(Z^s-a)\},\qquad a\in\F_q^\times.
\]
The numbers of irreducible factors of degree $e$ in $P_f$ and $Q_f$ are respectively
\begin{align}
 b_e&=\frac1e\sum_{s\mid e}\mu(e/s)(q^{d_s(1)}-1),\label{v148:eq-ordinary-types}\\
 b_e^{\rm proj}&=\frac1{e(q-1)}\sum_{s\mid e}\mu(e/s)
       \sum_{a\in\F_q^\times}(q^{d_s(a)}-1).\label{v148:eq-projective-types}
\end{align}
All factor degrees divide $rp^j$, where $r$ is the least common multiple of the root orders of the $\pi$, and $p^j$ is the least power of $p$ at least as large as the largest part $\lambda_{\pi,j}$.
\end{theorem}
\begin{proof}
The $K$-Frobenius on the root space is the return operator. In the irreducible case, identify this space with $\F_{q^n}$, with return action multiplication by $\alpha$. Every nonzero vector has orbit length $\ell$, and every line has orbit length the least positive $e$ such that $\alpha^e\in\F_q^\times$. Transitivity on nonzero vectors or on lines precludes a nonzero proper invariant subspace, so the converse also forces the return polynomial to be irreducible.

On general elementary modules, $d_s(a)=\dim\ker(T^s-aI)$. Thus the fixed nonzero-vector count is $q^{d_s(1)}-1$, while the fixed-line count is
$(q-1)^{-1}\sum_a(q^{d_s(a)}-1)$. Möbius inversion gives the two formulas. The semisimple part has order $r$, and the largest unipotent block has order $p^j$, giving the degree bound. Separability follows from $L_f'=c\ne0$.
\end{proof}

When $\gcd(n,k)=1$, the frequency of a prescribed elementary-divisor type and a compatible constant coefficient is
\begin{equation}\label{v148:eq-type-frequency}
 \frac{q-1}{q^k-1}\prod_\pi T_k(\lambda_\pi;q^{\deg\pi}).
\end{equation}
Local partitions therefore distinguish different actions with the same return polynomial and give exact frequencies of ordinary factorization patterns.

\begin{proposition}[Counting projectively irreducible members]\label{v148:prop-projective-count}
Let $n\ge2$, $\gcd(n,k)=1$, $S=(q^n-1)/(q-1)$, and $d=\ord((-1)^{nk}\Nm(c))$. If
$\gcd(S,(q-1)/d)=1$, the number of irreducible $Q_f$ with fixed $c$ is
\[
 \frac{C_{n,k}(q)}n\,
 S\prod_{\substack{\rho\mid S\text{ prime}\\\rho\nmid q-1}}
 (1-\rho^{-1});
\]
otherwise it is zero.
\end{proposition}
\begin{proof}
Choose a generator $\gamma$ of $\F_{q^n}^\times$ and write $\alpha=\gamma^u$. Maximal projective order is equivalent to $\gcd(u,S)=1$. A prescribed norm fixes $u$ modulo $q-1$. This residue class is invertible at the common prime divisors of $S$ and $q-1$ precisely under the stated gcd condition. Excluding the remaining prime divisors by the Chinese remainder theorem gives the displayed product. Maximal projective order forces degree $n$; divide by $n$ and multiply by the fiber size.
\end{proof}

\begin{example}[One return target and several ordinary patterns]\label{v148:ex-return-patterns}
In $\F_4[t;\bar{\phantom a}]$, take $f=t^3+a_2t^2+a_1t+1$. For target $(Z+1)^3$, the allowed partitions are $(3)$ and $(2,1)$. Formula~\eqref{v148:eq-local-type} gives totals $12$ and $3$, hence frequencies $4$ and $1$ after fixing the constant coefficient. The factor degrees of $P_f=X^7+a_2X^3+a_1X+1$ are respectively $(1,2,4)$ and $(1,1,1,2,2)$. The two irreducible cubic targets give six irreducible degree-seven members; target $(Z+1)(Z^2+Z+1)$ gives five members of pattern $(1,3,3)$. The four frequencies sum to $16$.

Over $\F_8$, take degree four and constant coefficient one. Since $C_{4,3}(2)=39$ and two of the three irreducible quartic targets are primitive, $78$ polynomials $P_f$ are irreducible of degree fifteen; the remaining $39$ split into three degree-five factors. Over $\F_9$, take degree three. Here $C_{3,2}(3)=7$, and each nonzero constant coefficient gives $28$ irreducible polynomials
$Q_f=Y^{13}+a_2Y^4+a_1Y+c$. When $\Nm(c)=1$, the corresponding $P_f$ split into two degree-thirteen factors; when $\Nm(c)=-1$, these $28$ polynomials $P_f$ are irreducible of degree twenty-six.
\end{example}

\section{Norms and Cyclic Decomposition of Weighted Frobenius}
\label{wfr:sec-weighted}

Combinations of multiplication operators and Frobenius give all linear operators. This section considers a single-term combination and reduces its cyclic structure to the relative norm of the weight.

Fix \(k=\F_q\) and \(E=\F_{q^n}\), choose \(1\le s<n\), and put
\[
 d=\gcd(n,s),\qquad m=n/d,\qquad K=\F_{q^d},\qquad \tau(x)=x^{q^s}.
\]
The fixed field of \(\tau\) is \(K\), and this automorphism generates the automorphism group of \(E/K\), which has order \(m\).
\begin{definition}\label{wfr:def-weighted}
For \(a\in E^\times\), call
\[
 S_a(x)=a\tau(x)
\]
the \textbf{weighted Frobenius operator}. Write
\[
 c=\Nm_{E/K}(a),\qquad r=[k(c):k],
\]
and let \(P_c\in k[T]\) denote the monic minimal polynomial of \(c\) over \(k\).
\end{definition}

The map \(S_a\) is invertible and linear over \(K\), and is also linear over \(k\). Regarding \(E\) as a one-dimensional space over itself, it satisfies
\(S_a(ux)=\tau(u)S_a(x)\) for \(u,x\in E\). \leavevmode\phantomsection\label{def:semilinear-map}\ignorespaces An additive map obeying this scalar transformation rule is called a \textbf{semilinear map} relative to \(\tau\).
For background on cyclic decompositions of semilinear transformations, see \cite{DFH}. We examine the linear structure over the fixed field \(K\), then restrict scalars to \(k\).

\begin{lemma}[A cyclic model over the fixed subfield]\label{wfr:lem-relative-model}
With \(T\) acting as \(S_a\), there is an isomorphism of \(K[T]\)-modules
\begin{equation}\label{wfr:eq-relative-model}
 (E,S_a)\simeq K[T]/(T^m-c).
\end{equation}
In particular, the minimal and characteristic polynomials of \(S_a\) over \(K\) are both \(T^m-c\).
\end{lemma}
\begin{proof}
Put \(a_0=1\) and \(a_j=\prod_{i=0}^{j-1}\tau^i(a)\). Iterated composition gives
\[
 S_a^j=a_j\tau^j,\qquad S_a^m=c\id.
\]
If \(\sum_{j=0}^{m-1}u_j S_a^j=0\) with \(u_j\in K\), then
\(\sum_j u_j a_j x^{q^{sj}}=0\) for every \(x\in E\).
Reduce \(sj\) modulo \(n\) to \(e_j\in\{0,\ldots,n-1\}\); since \(m=n/\gcd(n,s)\), these \(e_j\) are distinct.
If some \(u_j\ne0\), the resulting \(\sum_j u_j a_jX^{q^{e_j}}\) would be nonzero of degree at most \(q^{n-1}\) but vanish on all \(q^n\) elements of \(E\), a contradiction.
Thus \(\id,S_a,\ldots,S_a^{m-1}\) are linearly independent over \(K\), and the minimal polynomial is exactly \(T^m-c\).
Its degree is \(\dim_K E=m\), so the cyclic vector lemma~\ref{lem:cyclic-vector} gives~\eqref{wfr:eq-relative-model}. The characteristic polynomial has the same degree and hence is equal to it.
\end{proof}

\begin{theorem}[Cyclic decomposition determined by the norm]\label{wfr:thm-norm-decomposition}
As \(k[T]\)-modules,
\begin{equation}\label{wfr:eq-norm-decomposition}
 \boxed{\quad (E,S_a)\simeq
 \bigl(k[T]/(P_c(T^m))\bigr)^{\oplus d/r}.\quad}
\end{equation}
Consequently, the minimal and characteristic polynomials of \(S_a\) as a \(k\)-linear operator are, respectively,
\begin{equation}\label{wfr:eq-minimal-characteristic}
 \mu_{S_a}(T)=P_c(T^m),\qquad
 \chi_{S_a}(T)=P_c(T^m)^{d/r}.
\end{equation}
This also applies when \(p\mid m\). Formula~\eqref{wfr:eq-norm-decomposition} is a module isomorphism, not a direct-product decomposition of the field \(E\).
\end{theorem}
\begin{proof}
Put \(L=k(c)\subseteq K\). Choose a basis of \(K/L\). Since \(T^m-c\in L[T]\), restriction of scalars gives
\[
 K[T]/(T^m-c)\simeq
 \bigl(L[T]/(T^m-c)\bigr)^{\oplus d/r}
\]
as \(L[T]\)-modules.
In \(L[T]/(T^m-c)\) we have \(\overline T^{\,m}=c\), while \(k[c]=L\), so \(\overline T\) already generates the entire quotient algebra as a \(k\)-algebra.
The map \(k[T]\to L[T]/(T^m-c)\) is therefore surjective, and \(P_c(T^m)\) belongs to its kernel.
The quotient \(k[T]/P_c(T^m)\) of the source and the target both have \(k\)-dimension \(mr\), so the induced map is an isomorphism.
Together with Lemma~\ref{wfr:lem-relative-model}, this proves the cyclic decomposition.
Each component has both minimal and characteristic polynomial \(P_c(T^m)\). The minimal polynomial of a direct sum is the least common multiple, and the characteristic polynomial is the product, giving the last formulas.
\end{proof}

\begin{corollary}[Cyclicity and realization inside a field]\label{wfr:cor-cyclic-realization}
The following conditions are equivalent:
\begin{enumerate}
\item \(S_a\) has a cyclic vector over \(k\);
\item \(k(\Nm_{E/K}(a))=K\);
\item There are a finite extension \(L'/k\) and an injective \(k\)-linear map \(\iota:E\to L'\) such that
\(\iota(S_a x)=\iota(x)^q\) for every \(x\in E\).
\end{enumerate}
\end{corollary}
\begin{proof}
By~\eqref{wfr:eq-minimal-characteristic}, \(\deg\mu_{S_a}=mr\), which equals \(n=md\) if and only if \(r=d\).
The cyclic vector lemma proves equivalence of the first two conditions, and Theorem~\ref{thm:frobenius-realization} proves equivalence of the first and third.
\end{proof}

Here the norm must generate all of \(K\) over \(k\); it need not generate the multiplicative group \(K^\times\). These generation properties are measured by field degree and multiplicative order, respectively.

\begin{corollary}[Similarity classification for a fixed Frobenius step]\label{wfr:cor-similarity}
Fix \(E,k,s\). For \(a,a'\in E^\times\), the operators \(S_a\) and \(S_{a'}\) are similar over \(k\) if and only if
\(\Nm_{E/K}(a)\) and \(\Nm_{E/K}(a')\) are conjugate over \(k\).
\end{corollary}
\begin{proof}
Conjugate norms have the same minimal polynomial, so~\eqref{wfr:eq-norm-decomposition} gives the same cyclic modules.
Conversely, similarity implies \(P_c(T^m)=P_{c'}(T^m)\); injectivity of substitution gives \(P_c=P_{c'}\), so the norms are conjugate.
\end{proof}

\begin{example}[Weights can change cyclicity]\leavevmode\phantomsection\label{wfr:ex-sixteen}
Take \(k=\F_2\), \(E=\F_{16}\), and \(s=2\), so \(K=\F_4\) and \(m=d=2\).
For \(a=1\), we have \(c=1\) and \(r=1\), hence
\[
 (E,S_1)\simeq\bigl(\F_2[T]/((T+1)^2)\bigr)^{\oplus2},
 \qquad \mu_{S_1}=(T+1)^2,\quad\chi_{S_1}=(T+1)^4.
\]
This operator has no cyclic vector over \(\F_2\).
Choose \(\omega\in\F_4\) with \(\omega^2+\omega+1=0\) and put \(a=\omega^2\).
Since \(\Nm_{E/K}(a)=a^5=\omega\), we now have \(r=d=2\), giving
\[
 \mu_{S_a}(T)=\chi_{S_a}(T)=T^4+T^2+1=(T^2+T+1)^2.
\]
This operator has a cyclic vector, but its minimal polynomial still has repeated factors, so it remains nondiagonalizable after extension of scalars.
The weight combines two cyclic components into one; having a cyclic vector and becoming diagonalizable after scalar extension are still distinct conditions.
\end{example}

Formula~\eqref{v53:eq-cyclic-gcd-count} also uses this cyclic-module decomposition to compute fixed spaces of the iterates of \(S_a\) and hence recover all cycles.

Section~\ref{afr:sec-affine-periods} adds a translation term to \(S_a\). The norm still controls the linear part, while a relative trace determines whether a change of origin can remove the translation.

\subsection{Return Spectra of Two-Dimensional Semilinear Operators}
\label{v146:sec-return-spectrum}
A one-dimensional weighted Frobenius becomes multiplication by a norm after iteration. In dimension two, the trace and determinant of the resulting linear operator together record its return characteristic polynomial.

Let \(P\) be a prime power, \(K=\F_{P^k}\), \(r\ge1\), \(\gcd(r,k)=1\), and \(\tau(x)=x^{P^r}\). For \(a,c\in K\), define
\begin{equation}\label{v146:eq-return-spectrum}
\begin{aligned}
 A_{a,c}&=\begin{pmatrix}0&-c\\1&-a\end{pmatrix},
 \qquad \Phi_{a,c}(v)=A_{a,c}\tau(v),\\
 M_{k,r}(a,c)&=A_{a,c}\tau(A_{a,c})\cdots\tau^{k-1}(A_{a,c}),\\
 B_{k,r}(a,c)&=-\operatorname{tr}M_{k,r}(a,c).
\end{aligned}
\end{equation}
Here \(\tau\) acts on coordinates or matrix entries, and \(\Phi_{a,c}^k=M_{k,r}(a,c)\) is \(K\)-linear.

\begin{proposition}[Descent of the return characteristic polynomial]\label{v146:prop-spectrum-descent}
For all \(a,c\in K\), including \(c=0\),
\begin{equation}\label{v146:eq-spectrum-descent}
 \det(ZI_2-M_{k,r}(a,c))
 =Z^2+B_{k,r}(a,c)Z+\Nm_{K/\F_P}(c)\in\F_P[Z].
\end{equation}
\end{proposition}
\begin{proof}
Applying \(\tau\) to the trace cyclically permutes the matrix factors. The identity \(\operatorname{tr}(UV)=\operatorname{tr}(VU)\) therefore makes the trace fixed by \(\tau\), without requiring invertibility. The determinant is \(\prod_{i=0}^{k-1}\tau^i(c)=\Nm_{K/\F_P}(c)\).
\end{proof}

\leavevmode\phantomsection\label{v146:def-return-spectrum}\ignorespaces
The multiset of roots of~\eqref{v146:eq-spectrum-descent} is the \textbf{return spectrum} of \(\Phi_{a,c}\). For \(c\ne0\), the return polynomial over \(\F_P\) is irreducible, split with distinct roots, or repeated-root. In the irreducible case its roots are \(\alpha,\alpha^P\), the ratio \(\rho=\alpha/\alpha^P\) belongs to \(\mu_{P+1}\), and
\[
 \frac{B_{k,r}(a,c)^2}{\Nm_{K/\F_P}(c)}=2+\rho+\rho^{-1}.
\]
This is the coefficient substitution in Proposition~\ref{v146:prop-lucas-pullback}. A repeated-root spectrum does not record the Jordan structure: the same \((Z-\lambda)^2\) can arise from a scalar operator or a nontrivial Jordan block.

Put
\begin{equation}\label{v146:eq-theta-definition}
 \Theta_{k,r}(a)=B_{k,r}(a,1),\qquad
 \mathcal I_P=\{b\in\F_P:Z^2+bZ+1\text{ is irreducible over }\F_P\}.
\end{equation}
Irreducibility here always refers to the fixed field \(\F_P\). For even \(k\), these quadratics already split over \(K\). When \(r=1\), Proposition~\ref{v145:prop-return-matrix} identifies them with the return characteristic polynomials on linearized root spaces. Chapter 18 will count the coefficient fibers of the different spectral types.

\section{The Lattice of Invariant Subspaces and Trace Duality}
\label{sec:invariant-duality}

A family of subspaces ordered by inclusion is a \textbf{lattice} if any two members have, within the family, both a largest common subspace and a smallest subspace containing them.
For invariant subspaces these operations are intersection and sum; this section identifies them with polynomial greatest common divisors and least common multiples.

Put \(F(T)=T^n-1\). A normal element \(\theta\) gives an isomorphism
\begin{equation}\label{eq:normal-module-model}
k[T]/(F)\longrightarrow E,\qquad \overline a\longmapsto a(\sigma)\theta.
\end{equation}
This is an isomorphism of \(k[T]\)-modules, with \(T\) acting on the right as \(\sigma\). In general it is not a field or ring isomorphism:
the element \(\overline1\) on the left maps to \(\theta\), and when \(n>1\), the polynomial \(T-1\) is a proper factor of \(F\), so the ring on the left has zero divisors.
This model describes the linear action of Frobenius; it does not replace field multiplication.

\begin{theorem}\label{thm:invariant-duality}
For a monic divisor \(h\) of \(F\), write \(U_h=\ker h(\sigma)\). Then \(h\mapsto U_h\) gives a bijection between the monic divisors of \(F\) and the Frobenius-invariant \(k\)-subspaces of \(E\), with
\[
\dim_kU_h=\deg h,\qquad
U_h\cap U_g=U_{\gcd(h,g)},\qquad
U_h+U_g=U_{\lcm(h,g)}.
\]
Use the monic reciprocal notation \(h^\vee\) from Section~\ref{sec:reciprocal-polynomials}. For the trace pairing of \(E/k\),
\begin{equation}\label{eq:invariant-orthogonal}
U_h^\perp=U_{(F/h)^\vee}.
\end{equation}
\end{theorem}
\begin{proof}
The bijection follows from Proposition~\ref{prop:invariant-root-subspaces}, and the dimension from Theorem~\ref{thm:linearized-gcd}.
Moreover, \(h\mid g\) implies \(U_h\subseteq U_g\). Conversely, in model~\eqref{eq:normal-module-model},
\(U_h\) corresponds to the ideal generated by \(F/h\): the equation \(ha\equiv0\pmod F\) is exactly equivalent to \(F/h\mid a\).
Thus \(U_h\subseteq U_g\) means that the ideal generated by \(F/h\) is contained in that generated by \(F/g\).
Taking preimages in \(k[T]\) gives \((F/h)\subseteq(F/g)\), equivalently \(F/g\mid F/h\), and hence \(h\mid g\).
The intersection formula follows from B\'ezout's identity.
The sum is contained in \(U_{\lcm(h,g)}\) and has dimension \(\deg h+\deg g-\deg\gcd(h,g)=\deg\lcm(h,g)\), so the two spaces agree.

Frobenius invariance of the trace shows that \(\sigma^{-1}\) is the adjoint of \(\sigma\).
If \(d=\deg h\), the adjoint of \(h(\sigma)\) is
\[
h(\sigma^{-1})=h(0)\sigma^{-d}h^\vee(\sigma).
\]
For any operator \(A\) on a space with nondegenerate pairing, \((\ker A)^\perp=\im A^*\):
the right side is contained in the left, and both have dimension \(n-\dim\ker A\). Thus \(U_h^\perp=\im h^\vee(\sigma)\). The latter space lies in \(\ker(F/h^\vee)(\sigma)\), and both have dimension \(n-d\), so they agree.
Finally, monic reciprocation preserves multiplication and \(F^\vee=F\), so \(F/h^\vee=(F/h)^\vee\).
\end{proof}

A linear operator over a finite field is called \textbf{semisimple} if it becomes diagonalizable after extending scalars to a splitting field of its minimal polynomial.
\begin{corollary}\label{cor:frobenius-semisimple}
Write \(n=p^e m\) with \(p\nmid m\). Then \(T^n-1=(T^m-1)^{p^e}\). Frobenius becomes diagonalizable after extending scalars if and only if \(p\nmid n\).
\end{corollary}
\begin{proof}
The derivative of \(T^m-1\) is \(mT^{m-1}\), so it has no repeated roots.
The minimal polynomial of Frobenius is \(T^n-1\) and has no repeated roots if and only if \(e=0\).
A finite-dimensional linear operator is diagonalizable over a splitting field if and only if its minimal polynomial has no repeated roots:
the B\'ezout decomposition for distinct linear factors gives a direct sum of eigenspaces, and the converse follows directly from a diagonal matrix.
\end{proof}

\begin{corollary}\label{cor:selfdual-invariant-space}
There is a Frobenius-invariant subspace \(U\subseteq E\) satisfying \(U=U^\perp\) if and only if both \(q\) and \(n\) are even.
In this case, the subfield \(\F_{q^{n/2}}\) is such a subspace.
\end{corollary}
\begin{proof}
If \(U=U_h=U_h^\perp\), the theorem gives \(h=(F/h)^\vee\), equivalently \(hh^\vee=F\).
The factor \(T-1\) has the same multiplicity in \(h\) and \(h^\vee\), so its multiplicity in \(F\) must be even.
That multiplicity is \(p^{v_p(n)}\); hence \(p=2\) and \(2\mid n\) are necessary.
Conversely, put \(L=\F_{q^{n/2}}\). When \(p=2\), every \(x,y\in L\) satisfies \(\Tr_{E/L}(xy)=2xy=0\), so \(L\subseteq L^\perp\); both spaces have dimension \(n/2\) and therefore agree.
\end{proof}

\begin{example}\label{exm:unipotent-sixteen}
Take \(E=\F_{16}\) and \(k=\F_2\), so \(F=(T+1)^4\). Its invariant subspaces are precisely
\[
U_j=\ker(\sigma-\id)^j\quad(0\leq j\leq4),\qquad \dim_{\F_2}U_j=j.
\]
Here \(U_1=\F_2\) and \(U_2=\F_4\), with \(U_3=\ker\Tr_{E/k}\) because \(1+T+T^2+T^3=(T+1)^3\).
The three-dimensional space \(U_3\) is not a subfield, whereas \(U_2=U_2^\perp\). The normal elements are precisely \(E\setminus U_3\).
Here every eigenvalue of Frobenius is one, but its action is not the identity; its minimal polynomial records four nontrivial linear layers.
\end{example}

The equation \(U=U^\perp\) for a self-dual subspace differs from the equation \(\Tr(b_ib_j)=\delta_{ij}\) for a self-dual basis.
The former requires the pairing to vanish identically on \(U\times U\); the latter requires the Gram matrix of the basis to be the identity.
The self-dual normal basis criterion below cannot be replaced by the subspace criterion of this section.

Intersections and sums are controlled by greatest common divisors and least common multiples of factors. For the product of subspaces in Definition~\ref{v70:def-subspace-product}, the semisimple case is instead controlled by multiplication of eigenvalue sets.
\begin{theorem}[Products of invariant subspaces in the semisimple case]\label{v70:thm-semisimple-product}
Let \(E=\F_{q^n}\), \(k=\F_q\), and \(\Char k=p\), with \(p\nmid n\). Take monic divisors \(h,g\) of \(T^n-1\), and put
\[
 U_h=\ker h(\sigma),\qquad U_g=\ker g(\sigma).
\]
In a finite extension \(K\) containing both \(E\) and all roots of \(T^n-1\), let \(A_h,A_g\) be the root sets of \(h,g\), respectively. Put
\[
 C=A_hA_g=\{\lambda\mu:\lambda\in A_h,\ \mu\in A_g\},
 \qquad c(T)=\prod_{\nu\in C}(T-\nu).
\]
Then \(c\in k[T]\) is a monic divisor of \(T^n-1\), and
\begin{equation}\label{v70:eq-semisimple-product}
 U_hU_g=U_c.
\end{equation}
The product root set retains only distinct elements; an empty root set gives \(c=1\) and \(U_c=0\).
\end{theorem}
\begin{proof}
The root sets \(A_h,A_g\) are stable under the \(q\)th-power map, hence so is \(C\). The coefficients of \(c\) therefore lie in \(k\).
Moreover, \(\lambda^n=\mu^n=1\) implies that every element of \(C\) is an \(n\)th root of unity, so \(c\mid T^n-1\).

Extend the conjugate coordinates in~\eqref{eq:conjugate-tensor-map} further to \(K\). This gives a \(K\)-algebra isomorphism
\[
 \mathcal J_K:K\otimes_k E\longrightarrow K^n,
 \qquad u\otimes x\longmapsto\bigl(u\sigma^i(x)\bigr)_{i=0}^{n-1}.
\]
Field multiplication becomes coordinatewise multiplication, and \(\id_K\otimes\sigma\) becomes the cyclic shift
\(P(z_0,\ldots,z_{n-1})=(z_1,\ldots,z_{n-1},z_0)\).
For each \(n\)th root of unity \(\lambda\), take
\[
 e_\lambda=(1,\lambda,\ldots,\lambda^{n-1}).
\]
Then \(Pe_\lambda=\lambda e_\lambda\), and coordinatewise multiplication gives
\[
 e_\lambda e_\mu=e_{\lambda\mu}.
\]
Since \(p\nmid n\), the polynomial \(T^n-1\) has \(n\) distinct roots. The corresponding Vandermonde matrix is invertible, so these vectors form a basis of \(K^n\).
Consequently,
\[
 \mathcal J_K(K\otimes_k U_h)=\bigoplus_{\lambda\in A_h}Ke_\lambda,
 \qquad
 \mathcal J_K(K\otimes_k U_g)=\bigoplus_{\mu\in A_g}Ke_\mu.
\]
The product of these spaces is exactly \(\bigoplus_{\nu\in C}Ke_\nu\), which is \(\mathcal J_K(K\otimes_k U_c)\).
Products of subspaces commute with extension of scalars: expanding products of pure tensors gives
\((K\otimes_k U_h)(K\otimes_k U_g)=K\otimes_k(U_hU_g)\).
Finally, extension of scalars preserves inclusions and dimensions of subspaces, so equality after extension implies~\eqref{v70:eq-semisimple-product}.
\end{proof}

Conjugate coordinates first preserve field multiplication, after which the cyclic shift is decomposed into eigenspaces.
The field \(E\) itself need not contain all roots of \(T^n-1\), so the common splitting field \(K\) is needed. This coordinatewise product is not componentwise multiplication in normal basis coordinates.

\begin{example}[Squares and mixed products of two three-dimensional spaces]\label{v70:ex-f128-products}
Take \(E=\F_{128}\), \(k=\F_2\), and
\[
 h_1=T^3+T+1,\qquad h_2=T^3+T^2+1.
\]
Neither polynomial has a root in \(\F_2\), so both are irreducible cubics, and
\[
 (T-1)h_1h_2=T^7-1.
\]
Choose a root \(\zeta\) of \(h_1\) in a common splitting field; its multiplicative order is seven.
The roots of \(h_1\) are \(\zeta^i\) with \(i\in A=\{1,2,4\}\). The remaining nontrivial seventh roots of unity are the roots of \(h_2\), corresponding to \(B=\{3,5,6\}\).
Addition modulo seven gives
\[
 A+A=B+B=\{1,2,3,4,5,6\},\qquad
 A+B=\mathbb Z/7\mathbb Z.
\]
The theorem therefore yields
\[
 U_{h_1}^2=U_{h_2}^2=U_{h_1h_2}=\ker\Tr_{E/k},
 \qquad U_{h_1}U_{h_2}=E.
\]
Here \(U_h^2\) denotes the subspace product \(U_hU_h\), not merely the set of squares of its elements.
Each three-dimensional space has square equal to the six-dimensional trace kernel, whereas their mixed product is the entire seven-dimensional field.
\end{example}

\begin{example}[Eigenvalues alone do not suffice in the repeated-root case]\label{v70:ex-f16-product-boundary}
As in Example~\ref{exm:unipotent-sixteen}, take \(E=\F_{16}\), \(k=\F_2\), and
\(W=U_3=\ker\Tr_{E/k}\). Then
\[
 W^2=E.
\]
Indeed, if \(W^2\ne E\), nondegeneracy of the trace pairing gives \(0\ne z\in(W^2)^\perp\).
For all \(x,y\in W\), we have \(\Tr_{E/k}(zxy)=0\), so \(zW\subseteq W^\perp\).
Multiplication by \(z\) is invertible, however, and hence \(\dim_k zW=3\), whereas \(\dim_k W^\perp=1\), a contradiction.
Here \(T^4-1=(T-1)^4\), so every nonzero invariant subspace has only the eigenvalue one. These root sets distinguish neither the dimensions nor the products of the spaces.
Thus the hypothesis \(p\nmid n\) in Theorem~\ref{v70:thm-semisimple-product} cannot be omitted.
\end{example}

\Needspace{10\baselineskip}
\subsection{Intersecting row spaces in an operator algebra}
\begin{proposition}[A row-space criterion for a singly generated algebra]\label{v150:prop-row-intersection}
Let $B\in M_n(k)$, $n\ge1$, with minimal polynomial $m_B$. Every two nonzero matrices in $k[B]$ have row spaces with nonzero intersection if and only if $m_B$ has just one distinct monic irreducible factor.
\end{proposition}
\begin{proof}
If $m_B$ has at least two coprime primary components, the Chinese remainder theorem gives nonzero idempotents $e_1,e_2\in k[B]$ with $e_1e_2=e_2e_1=0$.
A row vector $v$ in both row spaces satisfies $ve_1=v$, whereas $v\in\operatorname{Row}(e_2)$ gives $ve_1=0$. Thus their intersection is zero.
Conversely, suppose $m_B=f^e$, with $f$ irreducible. Every nonzero $a(B)$ has the form $f(B)^ju(B)$, where $0\le j<e$ and $u(B)$ is a unit. The exponent $j$ is unique, although the unit need not be.
The factors commute, and multiplication on the left by an invertible matrix preserves row space. Hence
\[
 \operatorname{Row}(a(B))=\operatorname{Row}(f(B)^j).
\]
These nonzero row spaces form a descending chain, so any two have nonzero intersection.
\end{proof}

If $B$ represents Frobenius on $\F_{q^n}/\F_q$, its minimal polynomial is $T^n-1$. Writing $n=p^a m$ with $p\nmid m$ gives
$T^n-1=(T^m-1)^{p^a}$. Since $T^m-1$ is square-free and divisible by $T-1$, the row-space intersection property holds exactly when $m=1$. It has the same source as the single primary normality condition for extensions of characteristic-power degree.

\begin{example}[Equal rank distributions with different row-space intersections]\label{v150:ex-transpose-intersection}
In the matrix space
\[
 \mathcal A=\left\{
 \begin{pmatrix}a&0&0\\b&a&0\\c&0&a\end{pmatrix}:a,b,c\in\F_q
 \right\},
\]
every matrix with $a\ne0$ is invertible. When $a=0$, every nonzero matrix has row space $\langle(1,0,0)\rangle$. Thus any two nonzero matrices have row spaces with nonzero intersection.
Transposition preserves each matrix rank, but $\mathcal A^{\mathsf T}$ contains $E_{12},E_{13}$, whose row spaces are $\langle(0,1,0)\rangle$ and $\langle(0,0,1)\rangle$, with zero intersection.
The rank distribution alone therefore does not recover the relative positions of the images of several operators.
\end{example}

\section{Cyclic Correlation Polynomials and Normality}
\begin{proposition}\label{prop:normal-circulant-test}
For any \(\alpha\in E\), put
\[
g_\alpha(T)=\sum_{i=0}^{n-1}\Tr_{E/k}(\alpha\alpha^{q^i})T^i\in k[T].
\]
Call \(g_\alpha\) the \textbf{cyclic correlation polynomial} of \(\alpha\); its coefficients record the trace pairings of \(\alpha\) with its Frobenius conjugates.
Then \(\alpha\) is normal if and only if \(\gcd(g_\alpha,T^n-1)=1\).
If it is normal, take \(u g_\alpha\equiv1\pmod{T^n-1}\); then \(u(\sigma)\alpha\) generates its dual normal basis.
Moreover, putting \(h_\alpha(T)=\sum_{i=0}^{n-1}\alpha^{q^i}T^i\in E[T]\) gives the equivalent criterion \(\gcd_{E[T]}(h_\alpha,T^n-1)=1\).
\end{proposition}
\begin{proof}
Put \(c_i=\Tr(\alpha\sigma^i(\alpha))\), with indices modulo \(n\).
The trace is invariant under \(\sigma\), so
\[
c_{-i}=\Tr(\sigma^i(\alpha)\alpha)=c_i,\qquad
\Tr(\sigma^i(\alpha)\sigma^j(\alpha))=c_{j-i}.
\]
The Gram matrix of the conjugates is therefore \(G=(c_{j-i})\).
If they form a basis, nondegeneracy of the trace pairing makes \(G\) invertible;
if they are linearly dependent, the same coefficient relation holds among the columns of \(G\), so \(G\) is singular.
Normality is therefore equivalent to invertibility of \(G\).

In the basis \(1,T,\ldots,T^{n-1}\) of \(R=k[T]/(T^n-1)\),
multiplication by \(g_\alpha=\sum_i c_iT^i\) has matrix \((c_{i-j})\), which equals \(G\) by \(c_{-i}=c_i\).
This multiplication map is invertible if and only if \(g_\alpha\) is a unit of \(R\), equivalently \(\gcd(g_\alpha,T^n-1)=1\).

If \(u(T)=\sum_i u_iT^i\) is the inverse modulo \(T^n-1\), put \(\beta=u(\sigma)\alpha\).
For every \(j\),
\[
\Tr(\beta\sigma^j(\alpha))
=\sum_i u_i c_{j-i}
=[T^j](u g_\alpha\bmod(T^n-1))
=\delta_{0j}.
\]
Thus \(\beta\) is the first dual vector; applying \(\sigma^r\) to both factors gives the remaining dual vectors \(\sigma^r(\beta)\).

Finally, multiplication by \(h_\alpha\) on \(E[T]/(T^n-1)\) has matrix \((\sigma^{i-j}(\alpha))_{i,j}\).
Relabeling column \(j\) as \(-j\bmod n\) gives \((\sigma^{i+j}(\alpha))_{i,j}\), the Moore matrix of the conjugates.
A column permutation preserves invertibility, so the Moore criterion and the unit criterion give the final equivalent condition.
\end{proof}

Cyclic correlation polynomials also give a direct construction of linearized operators. Write \(c_i=\Tr_{E/k}(\alpha\alpha^{q^i})\) and put
\[
\mathcal C_\alpha(X)=\mathcal L_{g_\alpha}(X)
=\sum_{i=0}^{n-1}c_iX^{q^i}.
\]
Since \(c_i\in k\), its Dickson matrix is exactly the trace Gram matrix of the conjugates:
\begin{equation}\label{eq:normal-correlation-dickson}
D_{\mathcal C_\alpha}=(c_{j-i})_{0\le i,j<n}
=\bigl(\Tr_{E/k}(\alpha^{q^i}\alpha^{q^j})\bigr)_{0\le i,j<n}.
\end{equation}
Proposition~\ref{prop:dickson-rank} therefore expresses normality as
\[
\alpha\text{ is normal}\quad\Longleftrightarrow\quad
\mathcal C_\alpha:E\longrightarrow E\text{ is invertible}.
\]
In the normal case, \(u g_\alpha\equiv1\pmod{T^n-1}\) gives \(\mathcal C_\alpha^{-1}=u(\sigma)\),
so the first dual normal element is exactly \(\mathcal C_\alpha^{-1}(\alpha)\).
Proposition~\ref{prop:dickson-operations} also represents this inverse operator by the inverse Gram matrix.
The same trace-correlation coefficients thus test whether the conjugates form a basis and give the inverse operator needed to recover dual coordinates.

\begin{proposition}[Product identity for conjugate coefficients and cyclic correlation]\label{v88:prop-conjugate-correlation}
Keep the notation \(h_\alpha,g_\alpha\) from this section, and put
\(h_\alpha^*(T)=h_\alpha(T^{-1})\) in \(E[T]/(T^n-1)\). Then
\begin{equation}\label{v88:eq-conjugate-correlation}
\boxed{\quad g_\alpha(T)\equiv h_\alpha(T)h_\alpha^*(T)\pmod{T^n-1}.\quad}
\end{equation}
In particular, \(g_\alpha(1)=\Tr_{E/k}(\alpha)^2\).
\end{proposition}
\begin{proof}
Collecting exponents modulo \(n\), the coefficient of \(T^r\) in the product is
\[
\sum_{j=0}^{n-1}\alpha^{q^{j+r}}\alpha^{q^j}
=\Tr_{E/k}(\alpha^{q^r}\alpha),
\]
which is the corresponding coefficient of \(g_\alpha\). Substituting \(T=1\)
then gives \(h_\alpha(1)=\Tr(\alpha)\) and the asserted square identity.
\end{proof}

\begin{proposition}[Moore determinants and cyclic multiplication matrices]\label{v88:prop-moore-circulant}
Let \(M_\alpha=(\alpha^{q^{i+j}})_{0\le i,j<n}\), and let \(C_\alpha\) be the matrix
of multiplication by \(h_\alpha\) on \(E[T]/(T^n-1)\) in its standard basis.
Then
\begin{equation}\label{v88:eq-moore-circulant}
\boxed{\quad \det M_\alpha=(-1)^{\lfloor(n-1)/2\rfloor}\det C_\alpha.\quad}
\end{equation}
If \(n=p^a\), then more precisely
\begin{equation}\label{v88:eq-moore-trace-power}
\boxed{\quad\det(\alpha^{q^{i+j}})_{0\le i,j<n}
=(-1)^{\lfloor(n-1)/2\rfloor}\Tr_{E/k}(\alpha)^n.\quad}
\end{equation}
\end{proposition}
\begin{proof}
The \((i,j)\)-entry of \(C_\alpha\) is \(\alpha^{q^{i-j}}\), with indices taken modulo
\(n\). Relabel its columns by \(j\mapsto-j\), obtaining \(M_\alpha\). This permutation
fixes column zero and exchanges \(j\) and \(n-j\) in
\(\lfloor(n-1)/2\rfloor\) pairs, giving the first formula.
If \(n=p^a\), then \(T^n-1=(T-1)^n\). Put \(\varepsilon=T-1\), so that
\(E[T]/(T^n-1)=E[\varepsilon]/(\varepsilon^n)\). In the basis
\(1,\varepsilon,\ldots,\varepsilon^{n-1}\), multiplication by
\(h_\alpha(1+\varepsilon)\) has a triangular matrix with every diagonal entry equal to
\(h_\alpha(1)\). Indeed, for \(0\le j<n\),
\[
h_\alpha(1+\varepsilon)\varepsilon^j
\equiv h_\alpha(1)\varepsilon^j\pmod{(\varepsilon^{j+1})}.
\]
Since the determinant is independent of the chosen basis,
\(\det C_\alpha=h_\alpha(1)^n=\Tr(\alpha)^n\), proving the second formula.
\end{proof}

The coefficients of \(h_\alpha\) are the Frobenius conjugates, whereas the
coefficients of \(g_\alpha\) are their trace pairings. The former test generation
through a Moore matrix; the latter test and recover duality through a cyclic Gram
matrix. Formula~\eqref{v88:eq-conjugate-correlation} connects them.
Here \(h_\alpha\) is a \emph{conjugate-coefficient polynomial}, not the polynomial
whose roots form the Frobenius orbit.

\section{When Nonzero Trace Implies Normality}
\begin{theorem}\label{thm:trace-normal-classification}
Let \(n\ge1\). The property that \emph{every element of degree \(n\) over \(k\) with nonzero trace is normal}
holds if and only if one of the following conditions is satisfied:
\begin{enumerate}
\item \(n=p^a\), where \(a\ge0\);
\item \(n=\ell\ne p\) is prime and \(\ord_\ell(q)=\ell-1\).
\end{enumerate}
\end{theorem}
\begin{proof}
Write \(n=p^a m\) with \(p\nmid m\) and put \(F=T^n-1=(T^m-1)^{p^a}\).
For \(\mu=\mu_\alpha\mid F\), the trace operator is \(F/(T-1)\) evaluated at \(\sigma\), so
\[
\Tr(\alpha)\ne0\quad\Longleftrightarrow\quad (T-1)^{p^a}\mid\mu.
\]
Every nonconstant monic divisor \(\mu\mid F\) is the Frobenius annihilator polynomial of some nonzero element, whose degree over \(k\) is \(\ord(\mu)\).
The question is therefore: when do \((T-1)^{p^a}\mid\mu\) and \(\ord(\mu)=n\) for a monic polynomial \(\mu\mid F\) force \(\mu=F\)?

If \(m=1\), the trace condition already forces \(\mu=F\).
If \(a=0\), \(m=\ell\) is prime, and \(\Phi_\ell\) is irreducible, then \(F=(T-1)\Phi_\ell\).
The trace condition requires the first factor and degree \(\ell\) requires the second, forcing \(\mu=F\).
The cyclotomic orbit theorem shows that \(\Phi_\ell\) is irreducible if and only if \(\ord_\ell(q)=\ell-1\).

Conversely, if \(a>0,m>1\), take the proper divisor
\[
\mu=(T-1)^{p^a}\frac{T^m-1}{T-1}.
\]
The second factor contains a primitive \(m\)th root of unity, of order \(m\); the multiplicity contributes \(p^a\).
Thus \(\ord(\mu)=n\) but \(\mu\ne F\), giving a nonnormal counterexample.
If \(a=0,m>1\) and \(m\) is either composite or prime with \(\Phi_m\) reducible, choose an irreducible factor \(P\mid\Phi_m\) and put \(\mu=(T-1)P\).
Its order is \(m\), but its degree satisfies \(1+\ord_m(q)<m\): for composite degree use \(1+\varphi(m)<m\), and for prime degree use reducibility of \(\Phi_m\).
This again gives the required counterexample and exhausts the remaining cases.
\end{proof}
Nonzero trace detects the multiplicity of \(T-1\) in the Frobenius annihilator polynomial. The classification identifies the parameters for which this condition, together with full degree, forces normality.

Proposition~\ref{v82:prop-normal-translates} allows base-field translation to alter only the trace condition. Exercise~\ref{v82:ex-universal-normal-translation} uses this to compare two different quantifiers: whether nonzero trace forces normality, and whether every full-degree element can be made normal by translation.

\begin{corollary}[Nonzero trace without a degree hypothesis]\label{v88:cor-trace-universal-normal}
For \(E=\F_{q^n}\) and \(k=\F_q\),
\[
\boxed{\quad
\bigl(\forall\alpha\in E,\ \Tr_{E/k}(\alpha)\ne0\Rightarrow\alpha\text{ is normal}\bigr)
\iff n=p^a\ (a\ge0).\quad}
\]
\end{corollary}
\begin{proof}
If \(n=p^a\), then \(T^n-1=(T-1)^n\). In the Frobenius cyclic module,
a residue class is a unit exactly when its reduction modulo \(T-1\) is nonzero;
this is detected by nonzero trace. If \(n\) is not a power of \(p\), let \(P\)
be an irreducible factor of \(T^n-1\) other than \(T-1\). Fix a normal element
\(\eta\) and set \(\alpha=P(\sigma)\eta\). Since \(P(1)\ne0\),
\(\Tr(\alpha)=P(1)\Tr(\eta)\ne0\), whereas \(P\) is not a unit in
\(k[T]/(T^n-1)\), so \(\alpha\) is not normal.
\end{proof}

The two statements have different quantifiers: Theorem~\ref{thm:trace-normal-classification} only
considers elements of \emph{degree exactly \(n\)}, whereas this corollary concerns
\emph{all} elements. For example, \(\Tr_{8/2}(1)=1\) but \(1\) is not normal;
every element of degree three and nonzero trace in that extension is normal.

\begin{example}[Full degree, nonzero trace, and normality]\label{exm:full-degree-nonnormal-f64}
Take \(E=\F_{64}\) and \(k=\F_2\), write \(\sigma(x)=x^2\), and choose a normal element \(\eta\) of \(E/k\). Put
\[
P(T)=T^2+T+1,\qquad
\alpha=P(\sigma)\eta=\eta+\eta^2+\eta^4.
\]
Since \(\eta\) is normal, its trace is nonzero, so \(\Tr_{E/k}(\eta)=1\). Frobenius invariance of the trace gives
\[
\Tr_{E/k}(\alpha)=\Tr_{E/k}(\eta)+\Tr_{E/k}(\eta^2)+\Tr_{E/k}(\eta^4)=1.
\]
On the other hand, \(T^6-1=(T+1)^2P(T)^2\), and the annihilator formula in the cyclic module gives
\[
\mu_\alpha(T)=\frac{T^6-1}{\gcd(T^6-1,P)}
=(T+1)^2P=T^4+T^3+T+1.
\]
Here \((T+1)^2\) has order two and \(P\) has order three; they are coprime, so
\[
[k(\alpha):k]=\ord(\mu_\alpha)=6,
\qquad
\dim_k\langle\alpha,\sigma\alpha,\ldots,\sigma^5\alpha\rangle
=\deg\mu_\alpha=4.
\]
Thus \(\alpha\) has full Frobenius period and nonzero trace but is not normal. It satisfies the nontrivial relation
\(\alpha^{16}+\alpha^8+\alpha^2+\alpha=0\).
Orbit length records the first return, while the dimension of the orbit span records linear independence; these are controlled by \(\ord(\mu_\alpha)\) and \(\deg\mu_\alpha\), respectively.

All such counterexamples can also be determined. Trace one requires the annihilator polynomial to contain \((T+1)^2\);
degree six requires it to contain \(P\), while nonnormality requires it to lack \(P^2\). The annihilator polynomial must therefore be the above \(\mu_\alpha\).
By Theorem~\ref{thm:root-module-count}, the number of such elements is exactly
\[
\varphi_2\bigl((T+1)^2P\bigr)=(2^2-2)(2^2-1)=6
\]
and they form a single Frobenius orbit.

The equation \(\mu_\alpha=(T+1)(T^3+1)\) also gives \(\alpha^8+\alpha\in\F_2\).
This value cannot be zero, since then the orbit length of \(\alpha\) would divide three; hence \(\alpha^8+\alpha=1\).
These six elements are therefore exactly the roots of the degree-six factor in Example~\ref{exm:additive-spectrum-two-six},
\(X^6+X^5+X^3+X^2+1\).
That example will reinterpret the same six roots through the return action on an additive fiber.
\end{example}

\subsection{Prescribed trace and conjugate-span dimension}
Nonzero trace detects only one primary component. Retaining the degree of the annihilator gives an exact count for every conjugate-span dimension.
\begin{theorem}[$k$-normal elements with prescribed trace]\label{v109:thm-knormal-trace}
Let $E=\F_{q^n}$, $p=\Char\F_q$, and $0\le k\le n$. An element $\alpha\in E$ is \textbf{$k$-normal} if its Frobenius conjugates span a space of dimension $n-k$.
Write $N_k(a)$ for the number of such elements with trace $a$, and put
\[
e=p^{v_p(n)},\qquad F=T^n-1,\qquad Q=F/(T-1),\qquad P_0=F/(T-1)^e.
\]
Summing over monic divisors gives
\[
\boxed{N_k(0)=\sum_{\substack{h\mid Q\\\deg h=n-k}}\varphi_q(h),\qquad
N_k(a)=q^{e-1}\sum_{\substack{h\mid P_0\\\deg h=n-k-e}}\varphi_q(h)\quad(a\ne0).}
\]
A sum with negative degree is empty, and $\varphi_q(1)=1$. With this convention, zero is the unique $n$-normal element.
\end{theorem}
\begin{proof}
In the Frobenius cyclic module $E\simeq\F_q[T]/(F)$, we have $\deg\mu_\alpha=n-k$, and exactly $\varphi_q(h)$ elements have annihilator $h$.
Moreover, $\Tr(\alpha)=Q(\sigma)\alpha$, so trace zero is equivalent to $\mu_\alpha\mid Q$, proving the first formula.
Since $T^{n/e}-1$ is squarefree, $P_0$ is coprime to $T-1$. Nonzero trace is equivalent to
$\mu_\alpha=(T-1)^e h$ with $h\mid P_0$. For fixed $h$, there are
$(q-1)q^{e-1}\varphi_q(h)$ such elements. Multiplication by a nonzero base-field scalar preserves the annihilator and gives bijections between the nonzero trace fibers. Each fiber therefore contains a fraction $1/(q-1)$ of these elements. Summing by degree proves the second formula.
\end{proof}
\begin{example}[All trace layers in the field of 27 elements]\label{v109:ex-knormal27}
For $\F_{27}/\F_3$, we have $F=(T-1)^3$, and hence
\[
\begin{array}{c|c|c|c}
k&\deg\mu_\alpha&\text{number}&\text{trace distribution}\\\hline
0&3&18&\text{nine each for }1,2\\
1&2&6&\text{all zero}\\
2&1&2&\text{all zero}\\
3&0&1&\text{the zero element}
\end{array}
\]
For instance, take $\theta^3-\theta-1=0$. Then $\sigma(\theta)-\theta=1$, so $\mu_\theta=(T-1)^2$.
Its field degree is three, its conjugate-span dimension is two, and its trace is zero. In contrast, $\Tr(\theta^2)=2$, so $\theta^2$ is normal.
Trace zero permits a substantial nilpotent chain, but excludes the full $(T-1)^3$ primary component.
\end{example}

\section{Normal Coordinates, Duality, and Multiplication}
Let \(\eta_i=\eta^{q^i}\). If \(x=\sum_i x_i\eta_i\) and \(x_i\in k\), then
\[
x^q=\sum_i x_i\eta_{i+1},
\qquad \Tr(x)=\left(\sum_i x_i\right)\Tr(\eta),
\]
with indices modulo \(n\). In normal coordinates, Frobenius acts as a cyclic shift.

\begin{proposition}\label{prop:dual-normal}
Let \((\eta,\eta^q,\ldots,\eta^{q^{n-1}})\) be a normal basis, and let \(\theta\) be the element dual to \(\eta\). Its ordered dual basis is \((\theta,\theta^q,\ldots,\theta^{q^{n-1}})\), and \(\Tr(\eta)\Tr(\theta)=1\).
\end{proposition}
\begin{proof}
By trace invariance, with indices taken modulo \(n\),
\[
\Tr(\eta^{q^i}\theta^{q^j})=\Tr(\eta^{q^{i-j}}\theta)=\delta_{ij}.
\]
These are exactly the defining relations for the ordered dual basis; uniqueness proves the assertion.
Using the duality relations again,
\[
\Tr(\eta)\Tr(\theta)
=\Tr\left(\left(\sum_i\eta_i\right)\theta\right)
=\sum_i\Tr(\eta_i\theta)=1.
\]
\end{proof}

In the expansion
\[
\eta_i\eta_j=\sum_k c_{ijk}\eta_k
\]
\leavevmode\phantomsection\label{def:multiplication-constants}\ignorespaces the coefficients \(c_{ijk}\in k\) are the \textbf{structure constants} of field multiplication in this basis.
If \((\theta_k)_k\) is the dual basis, the coordinate-recovery formula gives
\[
c_{ijk}=\Tr(\eta_i\eta_j\theta_k).
\]
They satisfy \(c_{i+1,j+1,k+1}=c_{ijk}\). Frobenius symmetry therefore reduces the number of constants that must be stored independently.
Normal coordinates make the \(q\)th-power operation particularly simple; the cost of general multiplication still depends on the sparsity of these constants.
The coordinatewise multiplication in Theorem~\ref{v70:thm-semisimple-product} takes place in conjugate evaluation coordinates after extending scalars; here we use the coefficients of field elements in a normal basis. These are different coordinate systems.

The cyclic shift above describes the \(q\)-Frobenius. When \(q=p^e\), the absolute Frobenius \(x\mapsto x^p\) also acts on coordinate coefficients, so its semilinearity must be retained.
\begin{proposition}[Normal coordinates of the absolute Frobenius]\label{v97:prop-absolute-frobenius-coordinates}
Let \(E=\F_{q^n}\), where \(q=p^e\), and choose a normal basis \(\mathcal B=(\alpha,\alpha^q,\ldots,\alpha^{q^{n-1}})\). Write
\[
 \alpha^p=\sum_{i=0}^{n-1}c_i\alpha^{q^i},\qquad
 C=(c_{i-j})_{i,j\bmod n},\qquad c_i\in\F_q.
\]
Then \(x\mapsto x^p\) acts on coordinate columns by \(v\mapsto Cv^{[p]}\), where \(^{[p]}\) denotes the \textbf{entrywise Frobenius}.
If \(S e_j=e_{j+1}\) is the cyclic shift matrix, with indices modulo \(n\), then
\begin{equation}\label{v97:eq-semilinear-norm}
 \boxed{CC^{[p]}\cdots C^{[p^{e-1}]}=S.}
\end{equation}
Put \(c(T)=\sum_i c_iT^i\), and let \(^{[p^j]}\) act only on polynomial coefficients. Equivalently,
\begin{equation}\label{v97:eq-semilinear-circulant}
 c(T)c^{[p]}(T)\cdots c^{[p^{e-1}]}(T)\equiv T\pmod{T^n-1}.
\end{equation}
These identities do not require \(p\nmid n\).
\end{proposition}
\begin{proof}
Write \(b_j=\alpha^{q^j}\). The absolute Frobenius commutes with the \(q\)-Frobenius, and \(c_i\in\F_q\), so
\[
 b_j^p=\sum_i c_i b_{i+j}.
\]
For \(x=\sum_j v_jb_j\), one has \(x^p=\sum_j v_j^p b_j^p\), giving the coordinate action \(Cv^{[p]}\).
Iterating \(e\) times successively applies the \(p\)th power to the coefficient matrices.
Since \(v^{[p^e]}=v\) and \(x^{p^e}=x^q\) acts as \(Sv\), this proves \eqref{v97:eq-semilinear-norm}.
Finally, \(C=c(S)\), and \(I,S,\ldots,S^{n-1}\) are linearly independent.
The algebra of circulant matrices is therefore isomorphic to \(k[T]/(T^n-1)\), yielding \eqref{v97:eq-semilinear-circulant}.
No diagonalization by roots of unity is used.
\end{proof}

Formula~\eqref{v97:eq-semilinear-norm} uses the same norm mechanism as Section~\ref{wfr:sec-weighted}: repeated composition of a semilinear map also transforms the coefficients at each step.
Here the final result is the cyclic shift, not the ordinary matrix power \(C^e\).

\begin{example}[Semilinear coordinates when the characteristic divides the basis length]\label{v97:ex-f16-semilinear}
Take \(K=\F_4=\F_2(\omega)\), where \(\omega^2+\omega+1=0\), and let \(\theta^2+\omega^2\theta+1=0\).
The quadratic takes the values \(1,\omega^2,\omega^2,1\) at \(0,1,\omega,\omega^2\), so \(K(\theta)=\F_{16}\).
Its defining equation gives
\[
 \theta^4=\theta+\omega^2,\qquad
 \theta^2=\theta+\omega\theta^4.
\]
The first identity shows that \((\theta,\theta^4)\) is linearly independent over \(K\), hence is a normal basis.
The absolute Frobenius has matrix
\[
 C=\begin{pmatrix}1&\omega\\\omega&1\end{pmatrix},\qquad
 CC^{[2]}=\begin{pmatrix}0&1\\1&0\end{pmatrix},\qquad C^2=\omega I.
\]
The two products differ. Both the characteristic and the basis length are two, but the semilinear identity remains valid.
\end{example}

\section{Descent to Subfields and Composition in Coprime Degrees}
\begin{proposition}\label{prop:normal-descent}
If \(\eta\) is normal in \(E/k\) and \(L=\F_{q^d}\subseteq E\), then \(b=\Tr_{E/L}(\eta)\) is normal in \(L/k\).
If the original normal basis is self-dual, so is the normal basis generated by \(b\).
\end{proposition}
\begin{proof}
We have \(b=\sum_{j=0}^{n/d-1}\sigma^{dj}(\eta)\), so \(b,\sigma(b),\ldots,\sigma^{d-1}(b)\) are sums over \(d\) disjoint groups partitioning the original normal basis.
These sums are linearly independent and hence form a basis of \(L/k\).
In the self-dual case, for \(0\leq i<d\), transitivity of the trace gives
\[
\Tr_{L/k}(b\sigma^i(b))
=\Tr_{E/k}(\eta\sigma^i(b))
=\sum_j\Tr_{E/k}(\eta\sigma^{i+dj}(\eta))=\delta_{i0}.
\]
\end{proof}

\begin{proposition}\label{prop:normal-coprime}
Let \(\gcd(a,b)=1\). If \(u\in\F_{q^a}\) and \(v\in\F_{q^b}\) are both normal over \(k\),
then \(uv\) is normal in the extension \(\F_{q^{ab}}/k\). If both original normal bases are self-dual, the resulting normal basis is also self-dual.
\end{proposition}
\begin{proof}
By Theorem~\ref{thm:tensor-fields}, multiplication gives \(\F_{q^a}\otimes_k\F_{q^b}\simeq\F_{q^{ab}}\). A tensor-product basis therefore maps to \(u^{q^i}v^{q^j}\), where \(0\le i<a\) and \(0\le j<b\).
The Chinese remainder theorem shows that \(t\mapsto(t\bmod a,t\bmod b)\) runs through all index pairs,
so these products are precisely the \(ab\) conjugates of \(uv\), proving normality.
Summing over the same index pairs, for \(x\in\F_{q^a}\) and \(y\in\F_{q^b}\) we have
\[
\Tr_{\F_{q^{ab}}/k}(xy)=\Tr_{\F_{q^a}/k}(x)\Tr_{\F_{q^b}/k}(y).
\]
Apply this to the product of any two basis vectors. Each factor trace is the corresponding Kronecker delta, so the product basis is self-dual.
\end{proof}

\paragraph{Normality over every intermediate base field.}
\leavevmode\phantomsection\label{v82:def-completely-normal}\ignorespaces
An element \(\alpha\in E\) is \textbf{completely normal} for \(E/k\) if it is normal for \(E/K\) for every intermediate field \(k\subseteq K\subseteq E\).
This differs from trace descent in Proposition~\ref{prop:normal-descent}: there the relative trace is normal in the subfield, whereas here the original element must be normal over each intermediate base field.

\begin{proposition}[Complete normality in characteristic-power degree]\label{v82:prop-completely-normal-p-power}
If \([E:k]=p^a\), where \(p=\Char k\) and \(a\ge0\), then the following conditions on \(\alpha\in E\) are equivalent:
\(\alpha\) is normal; \(\alpha\) is completely normal; \(\Tr_{E/k}(\alpha)\ne0\).
\end{proposition}
\begin{proof}
The preceding criterion for characteristic-power degree identifies normality with nonzero total trace.
If that trace is nonzero, transitivity of trace gives \(\Tr_{E/K}(\alpha)\ne0\) for every intermediate field \(K\).
Since \([E:K]\) is again a power of \(p\), the same criterion applied to \(E/K\) makes \(\alpha\) normal over \(K\).
Conversely, complete normality includes normality for \(K=k\).
\end{proof}

\paragraph{Base-field translations of completely normal elements.}
\begin{lemma}\label{v106:lem-complete-translate}
Suppose that \(\alpha\) is completely normal in \(E=\F_{q^n}\) over \(k=\F_q\), and let \(c\in k\). Then
\[
\boxed{\alpha+c\text{ is completely normal}\quad\Longleftrightarrow\quad
\alpha+c\text{ is normal over }k.}
\]
If \(p\mid n\), every element of \(\alpha+k\) is completely normal. If \(p\nmid n\), the unique exceptional parameter is \(c=-\Tr_{E/k}(\alpha)/n\).
\end{lemma}
\begin{proof}
Only the reverse implication needs proof. Take an intermediate field \(K=\F_{q^d}\), and put \(m=n/d\). The original element \(\alpha\) is normal over \(K\), and \(c\in K\). If \(p\mid m\), Proposition~\ref{v82:prop-normal-translates} preserves normality over \(K\). If \(p\nmid m\), the same proposition says that the only possible obstruction is \(\Tr_{E/K}(\alpha+c)=0\). Applying \(\Tr_{K/k}\) would then give \(\Tr_{E/k}(\alpha+c)=0\), contrary to normality over the base field. Thus every intermediate-field condition holds. The final assertions follow by applying the translation proposition to \(E/k\).
\end{proof}

\begin{corollary}[Exact count of completely normal base-field translates]\label{v131:cor-base-complete-count}
Let $E=\F_{q^n}$, $k=\F_q$, $n\ge2$, and $p=\Char k$. For any $\theta\in E$, put
\[
 C_{q,n}(\theta)=\#\{t\in k:\theta+t\text{ is completely normal in }E/k\}.
\]
If $p\mid n$, then $C_{q,n}(\theta)\in\{0,q\}$; if $p\nmid n$, then
$C_{q,n}(\theta)\in\{0,q-1\}$. In the latter case, whenever the count is positive,
the unique excluded parameter is $-\Tr_{E/k}(\theta)/n$.
\end{corollary}
\begin{proof}
If no translate is completely normal, the count is zero. Otherwise choose $t_0\in k$
such that $\alpha=\theta+t_0$ is completely normal. By the preceding lemma,
$\theta+t=\alpha+(t-t_0)$ is completely normal precisely when it is normal over $k$.
Proposition~\ref{v82:prop-normal-translates} makes normality constant under base-field
translation when $p\mid n$, and excludes exactly one parameter when $p\nmid n$.
The latter parameter is determined by
$\Tr_{E/k}(\theta+t)=\Tr_{E/k}(\theta)+nt$.
\end{proof}

\begin{example}[Three ranges of preservation in the field of order 64]\label{v106:ex-complete-f64}
Use the model verified in full in Section~\ref{sec:f64-computation}:
\[
\alpha^6+\alpha+1=0,\qquad u=\alpha^{21},\qquad v=\alpha^9.
\]
Then \(u^2+u+1=0\) and \(v^3+v^2+1=0\); these elements generate \(\F_4\) and \(\F_8\), respectively. Both are normal over \(\F_2\), so Proposition~\ref{prop:normal-coprime} makes \(a=uv\) normal over \(\F_2\). Over \(\F_4\), its conjugates are \(uv,uv^4,uv^2\), which are independent by extending scalars in the normal basis of \(\F_8\). Over \(\F_8\), its conjugates are \(uv,u^2v\), independent by extending scalars in the normal basis of \(\F_4\). These exhaust the proper intermediate fields, and \(a\ne0\), so \(a\) is completely normal. However, \(\ord(a)=\lcm(3,7)=21\).

Put \(b=1+uv\). Since \(2\mid6\), the lemma makes \(b\) completely normal. Using the two low-degree relations, first compute \(b^3=u(v^2+v)\) and \(b^7=uv+u+v^2+1\), and then obtain
\[
b^{21}=u\ne1,\qquad b^9=1+v+v^2\ne1.
\]
The second inequality follows from \(v\notin\{0,1\}\). As the only distinct prime divisors of \(63\) are \(3,7\), the element \(b\) is primitive.

Now put \(c=uv+v^2=b+(v^2+1)\). Here \(v^2+1\in W=\ker(\sigma^2+\sigma+1)\), so Proposition~\ref{v106:prop-annihilator-fiber} makes \(c\) normal over \(\F_2\). Yet
\[
c+u c^4+u^2c^{16}=0.
\]
Indeed, \(c^4=uv^4+v\) and \(c^{16}=uv^2+v^4\); substitution cancels the coefficients of \(v,v^2,v^4\) separately. This is a nontrivial conjugate relation over \(\F_4\). On the other hand, \(c^{21}=u+1=u^2\ne1\) and \(c^9=v^2+v\ne1\), so \(c\) is also primitive. Within this same field, primitive normality therefore still fails to imply complete normality.
\begin{center}\small
\begin{tabular}{@{}lll@{}}\toprule
Element&Preserved property&Boundary\\\midrule
\(uv\)&Completely normal&Multiplicative order 21\\
\(1+uv\)&Completely normal and primitive&Multiplicative order 63\\
\(uv+v^2\)&Primitive and normal over \(\F_2\)&Not normal over \(\F_4\)\\\bottomrule
\end{tabular}\end{center}
Starting from a completely normal element, base-field translation controls all levels. Translation by a larger invariant subspace guarantees only the annihilator data over the base field specified in the proposition.
\end{example}

\paragraph{A different behaviour for intermediate-field translations.}
The corollary concerns only base-field parameters. When $t$ varies in a proper
intermediate field, complete normality can change even if $p\mid n$;
the base-field dichotomy does not extend to this larger parameter set.

\begin{example}[Intermediate-field translates in a binary extension of degree six]\label{v131:ex-intermediate-complete-f64}
Retain $u\in\F_4$ and $v\in\F_8$ from Example~\ref{v106:ex-complete-f64},
and put $\theta=uv$. It is completely normal over $\F_2$, and
Corollary~\ref{v131:cor-base-complete-count} shows that both elements of
$\theta+\F_2$ are completely normal. Now take $t\in\F_4$ and write
$t=Au+Bu^2$ in the normal basis $(u,u^2)$ over $\F_2$, with $A,B\in\F_2$.
Since $1=v+v^2+v^4$ and $\gcd(2,3)=1$, the six Frobenius conjugates of $uv$
form a normal basis of $E/\F_2$. In that basis the coordinate polynomial of
$\theta+t$ is
\[
 h(T)=1+A(1+T^2+T^4)+B(T+T^3+T^5).
\]
As $T^6-1=(T+1)^2(T^2+T+1)^2$, normality is equivalent to nonvanishing modulo
each of the two distinct irreducible factors. Modulo $T^2+T+1$ one has
$h\equiv1$, while $h(1)=1+A+B$. Consequently,
\[
 \theta+t\text{ is normal over }\F_2\iff A=B\iff t\in\{0,1\}.
\]
Both of these translates are completely normal by the base-field lemma, hence
\[
 \boxed{\#\{t\in\F_4:\theta+t\text{ is completely normal}\}=2.}
\]
In addition, $t=u$ gives $\Tr_{E/\F_4}(\theta+t)=u+t=0$, so it also destroys
normality at another level. This shows why the overlaps between the
intermediate-field failure conditions must be counted rather than presumed independent.
\end{example}

\paragraph{Partial resolvents and recovery.}
\phantomsection\label{v113:par-resolvent-question}
Weighted sums of conjugates can provide a normality test, but can the same data recover the element itself? If two elements have identical full resolvents, in which subspace must their difference lie? This is a further question beyond testing generation.
The proposition below shows that the full resolvents retain exactly a relative trace. In the presence of repeated factors, different elements can therefore have identical full resolvent data. To distinguish those elements, we retain positional information in the sums: separate sums are recorded according to positions modulo the characteristic-power part of the extension degree. Partial resolvents thereby become coordinates for recovering the conjugates, rather than merely another normality test.

\begin{proposition}\label{v106:prop-partial-resolvents}
Write \(n=p^s n_0\), where \(p\nmid n_0\), and choose a common finite extension \(\Omega\) containing \(E\) and \(\mu_{n_0}\). For \(\eta\in\mu_{n_0}\) and \(j\in\Z/p^s\Z\), define the \textbf{partial resolvent}
\[
F_\eta(w;j)=\sum_{\substack{0\le t<n\\t\equiv j\ ({\rm mod}\ p^s)}}\eta^{-t}\sigma^t(w).
\]
For \(0\le t<n\), the recovery formula is
\begin{equation}\label{v106:eq-partial-inversion}
\boxed{\sigma^t(w)=\frac1{n_0}\sum_{\eta\in\mu_{n_0}}\eta^tF_\eta(w;t\bmod p^s).}
\end{equation}
Put \(R_\eta(w)=\sum_jF_\eta(w;j)\). The element \(w\) is normal if and only if every \(R_\eta(w)\ne0\).
More precisely, put \(K_0=\F_{q^{n_0}}\) and \(\tau=\Tr_{E/K_0}\). The full resolvents determine exactly the relative trace \(\tau(w)\):
\begin{equation}\label{v112:eq-resolvent-trace}
\begin{aligned}
R_\eta(w)&=\sum_{j=0}^{n_0-1}\eta^{-j}\sigma^j(\tau(w)),\\
\tau(w)&=\frac1{n_0}\sum_{\eta\in\mu_{n_0}}R_\eta(w).
\end{aligned}
\end{equation}
Consequently,
\begin{equation}\label{v112:eq-resolvent-kernel}
\bigcap_{\eta\in\mu_{n_0}}\ker R_\eta
=\ker\tau=(\sigma^{n_0}-1)E.
\end{equation}
As a \(k\)-linear map, \(w\mapsto(R_\eta(w))_\eta\) has rank \(n_0\) and kernel dimension \(n-n_0\). Thus it recovers \(w\) precisely when \(s=0\); the partial resolvents recover \(w\) for every \(s\).
\end{proposition}
\begin{proof}
Expand the right-hand side of the recovery formula. Root-of-unity orthogonality retains only indices congruent to \(t\) modulo \(n_0\); the partial sum already fixes their residue modulo \(p^s\). The Chinese remainder theorem leaves just \(\sigma^t(w)\).

Put \(A_w(Z)=\sum_{t=0}^{n-1}\sigma^t(w)Z^t\). Reordering the rows of the Moore matrix \((\sigma^{i+j}w)_{i,j}\) gives a circulant matrix. Its invertibility is equivalent to \(A_w\) being a unit in \(\Omega[Z]/(Z^n-1)\), hence to \(A_w(\eta^{-1})=R_\eta(w)\ne0\) for every \(\eta\). Here \(Z^n-1=(Z^{n_0}-1)^{p^s}\), and units are detected at the distinct roots. The Moore criterion proves the normality assertion.
For the remaining assertions, write each index uniquely as \(t=j+an_0\), with \(0\le j<n_0\) and \(0\le a<p^s\). Since \(\eta^{n_0}=1\), grouping the full sum gives the first identity in~\eqref{v112:eq-resolvent-trace}. Summing over \(\eta\) and using root-of-unity orthogonality gives the second. Thus all resolvents vanish exactly when \(\tau(w)=0\).
Additive Hilbert 90 for \(E/K_0\), whose generating Frobenius is \(\sigma^{n_0}\), identifies this kernel with \((\sigma^{n_0}-1)E\). Surjectivity of the relative trace gives its dimension \(n-n_0\), and rank--nullity proves the rank assertion.
\end{proof}
In the preceding example, \(n=6\), \(p^s=2\), and \(n_0=3\); we may take \(\eta\in\{1,u,u^2\}\) inside \(E\). For each \(\eta\), retain separately the even positions \(0,2,4\) and the odd positions \(1,3,5\). Combining them retains exactly \(\Tr_{E/\F_8}(w)=w+w^8\), with kernel \((\sigma^3-1)E=\F_8\); keeping them separate recovers all conjugates by the formula. Testing generation can require less information than recovering the element.

\begin{remark}[Shared coordinates and dependent conditions]\label{v106:rem-resolvent-research}
Normality conditions over different intermediate fields come from the same conjugates; their success proportions cannot simply be multiplied. Example~\ref{v106:ex-complete-f64} already shows the insufficiency of testing only base-field normality. The dependency pattern is determined by independent pair-orbit coordinates, fibre sums at each level, and the supports of the conditions. The following section establishes the corresponding exact product formula.
\end{remark}

\section{Independent resolvent components and the count of completely normal elements}
\label{v131:sec-resolvent-components}

Normality over different intermediate fields is determined by the same conjugates, but the conditions need not involve the same coordinates. The partial resolvents of the preceding section retain the positions lost when the Frobenius polynomial has repeated factors. We now use them to separate the conditions into groups that can be counted independently. This requires distinguishing independence of linear coordinates from independence of nonvanishing conditions that may share those coordinates. The fibre sums and condition reduction follow the structure in~\cite[Sections~2--3]{Hellwig2026Resolvents}; the proofs below are given in the language of the preceding CRT, relative-trace, and cyclic-module theory.

\subsection{Pair orbits and independent linear coordinates}
\label{v131:subsec-pair-orbits}
Retain $E=\F_{q^n}$, $k=\F_q$, $n=p^a n_0$, $p\nmid n_0$, where $p=\Char k$, and choose a common finite extension $\Omega$ containing $E$ and $\mu_{n_0}$. We use the partial resolvents $F_\eta(w;u)$ of~\eqref{v106:eq-partial-inversion}. Consider the permutation
\[
 (\eta,u)\longmapsto (\eta^q,u+1)
\]
of $\mu_{n_0}\times\Z/p^a\Z$. Its orbits will be called \emph{pair orbits}. Denote their set by $\Pi$. For each $\pi\in\Pi$ choose a representative $(\eta_\pi,u_\pi)$ and let $s_\pi=|\pi|$.

\begin{proposition}[Independent coordinates from partial resolvents]\label{v131:prop-pair-coordinates}
For all $j\ge0$,
\begin{equation}\label{v131:eq-pair-frobenius}
 F_{\eta^{q^j}}(w;u+j)=\eta^{-jq^j}F_\eta(w;u)^{q^j}.
\end{equation}
Put
\[
 X_\pi=\{x\in\Omega:x=\eta_\pi^{-s_\pi}x^{q^{s_\pi}}\}.
\]
Then $\dim_kX_\pi=s_\pi$, and
\begin{equation}\label{v131:eq-pair-product}
 E\longrightarrow\prod_{\pi\in\Pi}X_\pi,\qquad
 w\longmapsto(F_{\eta_\pi}(w;u_\pi))_{\pi\in\Pi}
\end{equation}
is a $k$-linear isomorphism. Hence, for uniformly distributed $w\in E$, the representative coordinates belonging to distinct pair orbits are independent and uniformly distributed in the corresponding spaces $X_\pi$.
\end{proposition}
\begin{proof}
Raise the defining sum for $F_\eta(w;u)$ to the $q$-th power and reindex $t$ as $t+1$. This gives $F_\eta(w;u)^q=\eta^qF_{\eta^q}(w;u+1)$, and iteration proves~\eqref{v131:eq-pair-frobenius}. All coordinates on one orbit are determined by one representative. Applying the relation with $j=s_\pi$ shows that the representative belongs to $X_\pi$. The root bound for the linearized polynomial $X^{q^{s_\pi}}-\eta_\pi^{s_\pi}X$ gives $\dim_kX_\pi\le s_\pi$. The recovery formula~\eqref{v106:eq-partial-inversion} makes~\eqref{v131:eq-pair-product} injective. Thus
\[
 n\le\sum_{\pi}\dim_kX_\pi\le\sum_{\pi}s_\pi
 =|\mu_{n_0}|p^a=n.
\]
All inequalities are equalities, proving surjectivity and the dimension assertions. A uniform distribution on a finite direct product has independent uniform coordinates.
\end{proof}

This is the coordinate independence supplied by CRT and roots-of-unity orthogonality. Full resolvents alone do not give such an isomorphism when $p\mid n$.

\subsection{Fibre sums and normality across intermediate fields}
\label{v131:subsec-fibre-sums}
Fix $d\mid n$, write $d=p^cd_0$ with $p\nmid d_0$ and $0\le c\le a$, and put $N=n_0/d_0$. Let $\zeta$ generate $\mu_{n_0}$. For $k_0\in\Z/N\Z$ define
\begin{equation}\label{v131:eq-fibre-sum}
 C_{d,k_0}(w)=\sum_{t=0}^{d_0-1}
 \ \sum_{u\in p^c\Z/p^a\Z}F_{\zeta^{k_0+tN}}(w;u).
\end{equation}
The summation index $t$ labels a fibre of roots of unity; it is unrelated to the base field $k$. Set
\[
 e(d,k_0)=\frac{N}{\gcd(N,k_0)},\qquad
 o(d,k_0)=\operatorname{ord}_{e(d,k_0)}(q),
\]
with $\operatorname{ord}_1(q)=1$. As shown below, the kernel depends only on the orbit of $k_0$ under multiplication by $q^d$; we identify indices in this orbit.

\begin{theorem}[Fibre-sum criterion for intermediate-field normality]\label{v131:thm-fibre-normal}
An element $w\in E$ is normal over $\F_{q^d}$ if and only if
\begin{equation}\label{v131:eq-fibre-normal}
 \boxed{C_{d,k_0}(w)\ne0\quad\text{for every }k_0\in\Z/N\Z.}
\end{equation}
If $V_{d,k_0}=\ker C_{d,k_0}$, then
\begin{equation}\label{v131:eq-fibre-codim}
 \boxed{\operatorname{codim}_kV_{d,k_0}
 =\operatorname{lcm}\bigl(d,o(d,k_0)\bigr).}
\end{equation}
When $d\mid o(d,k_0)$, this codimension is $o(d,k_0)$.
\end{theorem}
\begin{proof}
Let $m=n/d=p^{a-c}N$ and $\widehat\xi=\zeta^{dk_0}$. Its order is $e(d,k_0)$, and $\widehat\xi\in\mu_N$. Consider the full resolvent of $E/\F_{q^d}$,
\[
 R^{(d)}_{\widehat\xi}(w)
 =\sum_{s=0}^{m-1}\widehat\xi^{-s}\sigma^{ds}(w).
\]
Substituting the recovery formula~\eqref{v106:eq-partial-inversion} gives
\[
 R^{(d)}_{\widehat\xi}(w)=\frac1{n_0}
 \sum_{\eta\in\mu_{n_0}}\sum_{s=0}^{m-1}
 (\widehat\xi^{-1}\eta^d)^s F_\eta(w;ds\bmod p^a).
\]
Apply CRT to the coprime moduli $N$ and $p^{a-c}$ in the index $s\bmod m$. The roots-of-unity factor depends only on $s\bmod N$, and the position $ds\bmod p^a$ only on $s\bmod p^{a-c}$. Orthogonality cancels the sum over $s\bmod N$ unless $\eta^d=\widehat\xi$, in which case it contributes $N$. The remaining positions run through $p^c\Z/p^a\Z$. Since raising to the $p^c$-th power permutes $\mu_{n_0}$, the equation $\eta^d=\zeta^{dk_0}$ is equivalent to $\eta=\zeta^{k_0+tN}$ for $0\le t<d_0$. Consequently
\begin{equation}\label{v131:eq-fibre-resolvent}
 R^{(d)}_{\widehat\xi}(w)=d_0^{-1}C_{d,k_0}(w).
\end{equation}
The cyclic-module unit criterion from Proposition~\ref{prop:normal-test}, applied to $E/\F_{q^d}$, characterizes normality by nonvanishing of these resolvents. Since $p\nmid d_0$, this proves~\eqref{v131:eq-fibre-normal}.

Let $g_{\widehat\xi}$ be the minimal polynomial of $\widehat\xi$ over $\F_{q^d}$. The resolvent kernel is $g_{\widehat\xi}(\sigma^d)E$; its codimension over $\F_{q^d}$ equals
\[
 \deg g_{\widehat\xi}=\operatorname{ord}_{e(d,k_0)}(q^d)
 =\frac{o(d,k_0)}{\gcd(d,o(d,k_0))}.
\]
Multiplying by $d$ converts this to the codimension over $k$ and proves~\eqref{v131:eq-fibre-codim}.
\end{proof}

The elements $\zeta^{d_0k_0}$ and $\widehat\xi=\zeta^{dk_0}$ must be distinguished. The former is the target of the fibre, while the latter is a root of the irreducible factor of the cyclic-module polynomial over $\F_{q^d}$. They need not coincide when $p\mid d$; this distinction is essential in repeated-root cases.

\subsection{Essential conditions and components}
\label{v131:subsec-essential-components}
Call $(d,k_0)$ an \emph{essential condition} if
\[
 d\mid o(d,k_0).
\]
Its failure kernel then has codimension $o(d,k_0)$ over $k$.

\begin{proposition}[Reduction to essential conditions]\label{v131:prop-essential-reduction}
An element $w\in E$ is completely normal if and only if it satisfies the nonvanishing conditions~\eqref{v131:eq-fibre-normal} for all essential conditions.
\end{proposition}
\begin{proof}
It suffices to show that failure of a nonessential condition forces failure of some essential condition. Let $C=(d,k_0)$ be nonessential, put $e=e(d,k_0)$, $o=\operatorname{ord}_e(q)$, and $h=\gcd(d,o)<d$. Set $\widehat\xi=\zeta^{dk_0}$. Since $\widehat\xi^{q^o}=\widehat\xi$, the element $\widehat\xi^{q^h}$ belongs to the same $q^d$-orbit as $\widehat\xi$.

For every $j\ge0$, the resolvents satisfy
\[
\begin{aligned}
 R^{(d)}_{\widehat\xi}(\sigma^d v)&=\widehat\xi R^{(d)}_{\widehat\xi}(v),\\
 R^{(d)}_{\widehat\xi^{q^j}}(\sigma^jv)&=R^{(d)}_{\widehat\xi}(v)^{q^j}.
\end{aligned}
\]
The first identity makes the kernel stable under $\sigma^d$. Together with the second, it shows that $q^d$-conjugate parameters have the same kernel. Since $\widehat\xi^{q^h}$ and $\widehat\xi$ are such conjugates, the second identity makes $V_C=\ker R^{(d)}_{\widehat\xi}$ stable under $\sigma^h$. Write $h=p^\gamma h_0$, $p\nmid h_0$, and take the level-$h$ target $\widehat\xi_1=\zeta^{hk_0}$. Since $\widehat\xi_1^{d/h}=\widehat\xi$, decomposing the summation index as $t=(d/h)s+r$, $0\le r<d/h$, gives
\[
 R^{(h)}_{\widehat\xi_1}(w)
 =\sum_{r=0}^{d/h-1}\widehat\xi_1^{-r}
 R^{(d)}_{\widehat\xi}(\sigma^{hr}w).
\]
Thus $V_C\subseteq V_{C'}$, where $C'$ denotes this level-$h$ condition. Moreover $e\mid\operatorname{ord}(\widehat\xi_1)$, whence $h\mid o\mid\operatorname{ord}_{\operatorname{ord}(\widehat\xi_1)}(q)$; therefore $C'$ is essential. Every nonessential failure is thus absorbed by an essential one.
\end{proof}

For an essential condition $C=(d,k_0)$, let $\operatorname{Supp}(C)\subseteq\Pi$ consist of the pair orbits of the coordinates used by~\eqref{v131:eq-fibre-sum}. Two essential conditions are \emph{adjacent} if their supports intersect; the equivalence classes generated by adjacency are the \emph{resolvent components}. Write them as $\mathcal C_1,\ldots,\mathcal C_s$, and define
\[
 G_{\mathcal C}=\{w\in E:C(w)\ne0\text{ for every }C\in\mathcal C\},
 \qquad w_{\mathcal C}=\sum_{C\in\mathcal C}q^{-o(C)}.
\]

\begin{theorem}[Product of independent components]\label{v131:thm-component-product}
Let $CN_q(n)$ be the number of elements of $E$ completely normal over $k$. For a uniformly distributed $w\in E$, the events $G_{\mathcal C_i}$ are independent, and
\begin{equation}\label{v131:eq-component-product}
 \boxed{\frac{CN_q(n)}{q^n}
 =\prod_{i=1}^s\Pr(G_{\mathcal C_i}).}
\end{equation}
If $w_{\mathcal C_i}<1$ for all $i$, then
\begin{equation}\label{v131:eq-component-union}
 \frac{CN_q(n)}{q^n}\ge\prod_{i=1}^s(1-w_{\mathcal C_i}).
\end{equation}
\end{theorem}
\begin{proof}
By~\eqref{v131:eq-pair-frobenius}, each fibre sum is determined by the representative coordinates indexed by $\operatorname{Supp}(C)$. The supports of distinct components are disjoint, for an overlap would make two conditions adjacent. Proposition~\ref{v131:prop-pair-coordinates} makes the coordinates on disjoint sets of pair orbits independent, and essential reduction shows that satisfying every component is equivalent to complete normality. This yields the exact product.

For each essential condition $C$, formula~\eqref{v131:eq-fibre-codim} gives failure probability $q^{-o(C)}$. Apply the union bound within each component and then multiply across components to obtain~\eqref{v131:eq-component-union}.
\end{proof}

It is probabilities of success for \emph{distinct components} that can be multiplied; conditions within one component still share coordinates. Hellwig~\cite[Section~7]{Hellwig2026Resolvents} uses closure of unit-root fibres under Frobenius and multiplicative-coset operations to show that every component without a codimension-one condition has weight at most $3/4$, with sharper bounds when $q$ is larger. This uniform inequality also requires combinatorial counting of admissible target degrees; it does not follow from CRT alone.

\begin{remark}[Primitivity together with complete normality]
The product of independent components concerns the additive generating conditions of complete normality. Requiring an element to generate the multiplicative group also calls for the multiplicative-character sieve used in the primitive normal basis theorem. Hellwig~\cite{Hellwig2026Resolvents} proves that a primitive completely normal element exists in \(\F_{q^n}\) for every prime power \(q\) and \(n\ge1\). That general proof also uses arithmetic closure of components, uniform combinatorial weight bounds, and explicit verification of finitely many parameter pairs. The direct-product decomposition by itself does not imply the primitive completely normal existence theorem.
\end{remark}

\begin{example}[Two independent components in the binary extension of degree six]\label{v131:ex-f64-components}
Take $(q,n)=(2,6)$, so $n_0=3$ and $p^a=2$, and let $\zeta\in\F_4$ be a primitive cube root of unity. The pair orbits are
\[
\begin{aligned}
 \pi_0&=\{(1,0),(1,1)\},\\
 \pi_1&=\{(\zeta,0),(\zeta^2,1)\},\\
 \pi_2&=\{(\zeta,1),(\zeta^2,0)\}.
\end{aligned}
\]
The essential conditions are $(1,0),(1,1),(2,1),(2,2)$. The first reads only $\pi_0$ and succeeds with probability $1/2$; the other three involve $\pi_1\cup\pi_2$ and form a second component. Put $x=F_\zeta(w;0)$ and $y=F_\zeta(w;1)$. They are independent and uniformly distributed in the same four-element $\F_2$-linear space. The three nonvanishing requirements become
\[
 x\ne0,\qquad y\ne0,\qquad x+y\ne0.
\]
There are $3\cdot2=6$ successful pairs among all $4^2=16$. Thus
\[
 \boxed{CN_2(6)=2^6\cdot\frac12\cdot\frac38=12.}
\]
The second component has weight $3\cdot2^{-2}=3/4$. Its union-bound success probability is only $1/4$, whereas the exact probability is $3/8$. Hence the product across components is exact but a union bound within a component may lose information.
\end{example}

\begin{example}[Four independent components in the binary extension of degree twenty-one]\label{v131:ex-f2-21-components}
Take $(q,n)=(2,21)$ and put $\sigma(x)=x^2$. Since
\[
 21=3\cdot7,\qquad \ord_7(2)=3,\qquad \ord_3(2)=2,
\]
the preceding reduction to essential conditions shows that it suffices to test normality over $\F_2$ and $\F_8$; the conditions over $\F_{2^7}$ follow from these. The successful coordinates in each independent component can be counted without enumerating the elements of $\F_{2^{21}}$.

Let $K=\F_8$, choose $c\in K^\times$ of order seven, and choose $z\in E=\F_{2^{21}}$ with $z^7=c$. Since $K$ contains the seventh roots of unity and $c$ is not a seventh power, $X^7-c$ is irreducible over $K$. Thus $E=K(z)$ and
\[
 \sigma^3(z)=z^8=cz.
\]
Grouping the power basis by doubling orbits on its exponents gives an $\F_2$-direct sum of $\sigma$-invariant spaces:
\begin{equation}\label{v131:eq-f21-decomposition}
 \begin{aligned}
 E&=K\oplus U_1\oplus U_2,\\
 U_1&=Kz\oplus Kz^2\oplus Kz^4,\\
 U_2&=Kz^3\oplus Kz^6\oplus Kz^5.
 \end{aligned}
\end{equation}
Their $\F_2$-dimensions are $3,9,9$, and their $\sigma$-annihilating polynomials are pairwise coprime. The number of normal elements in $K/\F_2$ is
$\varphi_2(T^3-1)=(2-1)(2^2-1)=3$; this part further splits into independent primary components of dimensions one and two, with respectively $1$ and $3$ successful vectors.

Consider $U_1$. Write $x=az+bz^2+ez^4$ with $a,b,e\in K$. The eigenvalues of $\sigma^3$ on these three basis vectors are $c,c^2,c^4$, so normality over $K$ requires
\[
 \boxed{abe\ne0}
\]
in this component. On the other hand put
\[
 \Pi=1+\sigma^7+\sigma^{14}\qquad\text{on }U_1.
\]
As $(\sigma^7)^3=1$ in characteristic two, $\Pi$ is a projector; its image and kernel are the cubic and sextic irreducible primary components for $\sigma$. Consequently, normality over $\F_2$ requires both $\Pi x\ne0$ and $x-\Pi x\ne0$.

Using $z^7=c$ and $c^7=1$ one obtains
\[
 \sigma^7(x)=c^3e^2z+c^4a^2z^2+cb^2z^4,
\]
and hence
\begin{equation}\label{v131:eq-f21-projection}
 \Pi x=Lz+c^4L^2z^2+c^2L^4z^4,
 \qquad L=a+c^3e^2+c^5b^4.
\end{equation}
Among the $7^3$ triples with $a,b,e\in K^\times$, the condition $\Pi x=0$ is equivalent to $L=0$. Given nonzero $b,e$, this determines $a$ uniquely; precisely seven pairs yield $a=0$, namely those with $e=cb^2$. Thus there are
\[
 7^2-7=42
\]
such failures. The condition $x=\Pi x$ has precisely seven nonzero solutions, corresponding to
\[
 (a,b,e)=(t,c^4t^2,c^2t^4),\qquad t\in K^\times.
\]
The two kinds of failure are disjoint because $\Pi$ is a projector. Therefore the number of elements in $U_1$ that satisfy the three conditions from $E/K$ and the two conditions from $E/\F_2$ is
\begin{equation}\label{v131:eq-f21-local-count}
 \boxed{7^3-42-7=294.}
\end{equation}
In $U_2$, use $z^3$ as a Kummer generator with seventh power $c^3$; the same calculation again gives $294$ successful elements.

Thus the dimensions, successful-vector counts and failure weights of the four independent components are
\[
\begin{array}{c|cccc}
\text{component dimension}&1&2&9&9\\ \hline
\text{successful vectors}&1&3&294&294\\
\text{failure weight}&\frac12&\frac14&\frac{33}{64}&\frac{33}{64}
\end{array}
\]
Each nine-dimensional component contains three conditions of codimension three over $\F_2$ arising from normality over $K$, and two conditions of codimensions three and six arising from normality over $\F_2$. Its weight is therefore $4\cdot2^{-3}+2^{-6}=33/64$. The exact product over independent components yields
\[
 \boxed{CN_2(21)=1\cdot3\cdot294^2=259308.}
\]
However, the sum of the failure probabilities of all essential conditions is
\[
 \frac12+\frac14+2\cdot\frac{33}{64}=\frac{57}{32}>1,
\]
so a global union bound gives no positive lower bound. Applying the union bound only within independent components gives instead
\[
 \frac{CN_2(21)}{2^{21}}
 \ge\frac12\cdot\frac34\left(\frac{31}{64}\right)^2>0.
\]
Example~\ref{v131:ex-f64-components} shows why conditions within one component cannot be treated as independent events. The present example shows why an undifferentiated global union bound loses the multiplicative structure across separate components.
\end{example}

\subsection{Exact component probabilities from local Fourier spectra}
\label{v131:subsec-component-exact}
The component decomposition determines which probabilities may be multiplied. To compute $\Pr(G_{\mathcal C})$ exactly requires the intersections of the linear failure kernels within the component. For $\mathcal J\subseteq\mathcal C$, put
\[
 V_{\mathcal J}=\bigcap_{C\in\mathcal J}V_C,\qquad
 V_\emptyset=E,\qquad
 r_{\mathcal J}=\operatorname{codim}_kV_{\mathcal J}.
\]

\begin{proposition}[Exact inclusion--exclusion and joint ranks]\label{v131:prop-component-exact}
For every resolvent component $\mathcal C$,
\begin{equation}\label{v131:eq-component-exact}
 \boxed{\Pr(G_{\mathcal C})
 =\sum_{\mathcal J\subseteq\mathcal C}(-1)^{|\mathcal J|}q^{-r_{\mathcal J}}.}
\end{equation}
If the $k$-linear maps $C:E\to\Omega$ are assembled into a joint map $(C)_{C\in\mathcal J}$, then $r_{\mathcal J}$ is its $k$-rank. If orthogonal complements are taken in $E$ relative to the nondegenerate trace pairing $\Tr_{E/k}(xy)$, then
\begin{equation}\label{v131:eq-component-dual-rank}
 r_{\mathcal J}=\dim_k\left(\sum_{C\in\mathcal J}V_C^\perp\right).
\end{equation}
\end{proposition}
\begin{proof}
Average $\prod_{C\in\mathcal C}(1-\boldsymbol1_{V_C})$ over $E$ and expand. Since $|V_{\mathcal J}|=q^{n-r_{\mathcal J}}$, the terms yield~\eqref{v131:eq-component-exact}. The intersection is the kernel of the joint map; rank--nullity gives the rank statement. Nondegeneracy of the trace pairing gives $(\cap_CV_C)^\perp=\sum_CV_C^\perp$ and the last formula.
\end{proof}

The local Fourier expansions in Section~\ref{v128:sec-normal-spectrum} give a frequency interpretation. For every $k$-subspace $V\subseteq E$, the indicator has Fourier transform
\[
 \widehat{\boldsymbol1_V}(b)=|V|\boldsymbol1_{V^\perp}(b).
\]
Thus intersections of failure kernels can be studied through sums of their dual frequency subspaces. Formula~\eqref{v131:eq-component-exact} is exact for all parameters, but a direct enumeration of $\mathcal J$ is exponential in the component size and does not itself give a closed formula.

\paragraph{Cyclotomic parameters and joint ranks.}
For a given resolvent component, can one compute every $r_{\mathcal J}$ directly from the levels of the conditions, the orders of their target roots of unity, their $q$-cyclotomic cosets, and their fibre phases? Can the resulting inclusion--exclusion formula be reorganised into a closed formula or a recurrence without enumerating all subsets? The codimension of each individual kernel is given by~\eqref{v131:eq-fibre-codim}, but the joint rank also depends on the phases and the linear relations among conditions. In the example $(2,6)$, each of the three kernels has codimension two, whereas their joint rank is four; this dependence explains why the exact success probability $3/8$ exceeds the union-bound value $1/4$.

This exact-counting problem complements the independent-component theory of~\cite[Sections~3 and~12]{Hellwig2026Resolvents}, and joins the earlier treatments of joint trace fibres, trace duality, local Fourier spectra, and cyclotomic orbits. The product formula identifies what is independent; the joint-rank formula identifies what remains to be computed. A uniform weight bound instead provides an existence guarantee without determining all those intersections.

\paragraph{Exact counts for intermediate-field translation families.}
Let $e\mid n$, $1\le e\le n$, and put $L=\F_{q^e}$. Define
\[
 C_{q,n;e}(\theta)=\#\{t\in L:\theta+t\text{ is completely normal over }k\}.
\]
For $e=1$ this is the count in Corollary~\ref{v131:cor-base-complete-count}.
For $e>1$ the failure conditions can overlap upon restriction to the affine
space $\theta+L$, and the independent-component product for uniform elements
of $E$ does not automatically remain valid on that affine space.

\begin{proposition}[Joint-rank formula for restricted translation families]\label{v131:prop-intermediate-translate-count}
Let $\mathcal N$ denote the essential conditions of
Proposition~\ref{v131:prop-essential-reduction}. For $C\in\mathcal N$, write
$\Phi_C:E\to\Omega$ for the corresponding $k$-linear fibre-sum map, with
$\ker\Phi_C=V_C$. Given $\mathcal J\subseteq\mathcal N$, put
$\Phi_{\mathcal J}=(\Phi_C)_{C\in\mathcal J}$ and
\[
 r_{\mathcal J,L}=\operatorname{rank}_k(\Phi_{\mathcal J}|_L),\qquad
 \varepsilon_{\mathcal J}(\theta)=
 \begin{cases}
 1,&-\Phi_{\mathcal J}(\theta)\in\Phi_{\mathcal J}(L),\\
 0,&\text{otherwise}.
 \end{cases}
\]
Use rank zero and $\varepsilon_\emptyset=1$ for the empty set. Then
\begin{equation}\label{v131:eq-intermediate-translate-count}
 \boxed{C_{q,n;e}(\theta)=
 \sum_{\mathcal J\subseteq\mathcal N}(-1)^{|\mathcal J|}
 \varepsilon_{\mathcal J}(\theta)q^{e-r_{\mathcal J,L}}.}
\end{equation}
\end{proposition}
\begin{proof}
On the affine space $\theta+L$, expand the complete-normality indicator
$\prod_{C\in\mathcal N}(1-\boldsymbol1_{\ker\Phi_C})$.
For a given $\mathcal J$, the parameters that fail all its conditions solve
$\Phi_{\mathcal J}(t)=-\Phi_{\mathcal J}(\theta)$.
The system either has no solution or has an affine fibre of size
$q^{e-r_{\mathcal J,L}}$. Summing these terms proves the formula.
\end{proof}

Thus the exact problem reduces to two computable quantities: the rank of each
restricted joint map, and membership of the prescribed value in its image.
Can the phases, cyclotomic cosets, and local cyclic-module data in the partial
resolvent coordinates provide formulas for these quantities without enumerating
all subsets? Can one classify the possible values of $C_{q,n;e}(\theta)$?
Example~\ref{v131:ex-intermediate-complete-f64} already shows that the simple
base-field dichotomy need not persist for $e>1$, even if $p\mid n$. If multiplicative primitivity is imposed as well, a multiplicative-character sieve may be applied on this affine family; compatibility between additive joint ranks and multiplicative character sums then requires separate analysis.

\section{Normal Elements with Prescribed Intermediate Traces}
\label{v58:sec-prescribed-normal-traces}

The relative trace of a normal element is normal. The converse question asks whether prescribed normal elements in subfields can be realized as traces of a single normal element in the larger field.
Normal coordinates turn this into polynomial congruences: traces on intersections determine compatibility, and units in quotient rings determine normality.
Write \(p=\Char k\). Under the assumption \(p\nmid n\), \cite[Theorem 3.6]{MKBPrescribedTraces} counts normal elements with prescribed intermediate traces;
the cyclic-module proof below retains the multiplicities of irreducible factors and therefore also applies when \(p\mid n\).

Let \(d_1,\ldots,d_s\mid n\), put \(L_i=\F_{q^{d_i}}\subseteq E\), and prescribe \(a_i\in L_i\).
For each pair \(i,j\), put \(e=\gcd(d_i,d_j)\). If
\begin{equation}\label{v58:eq-trace-compatibility}
 \Tr_{L_i/\F_{q^e}}(a_i)=\Tr_{L_j/\F_{q^e}}(a_j)
\end{equation}
\leavevmode\phantomsection\label{v58:def-trace-compatibility}\ignorespaces
holds for every pair, the values satisfy the \textbf{trace compatibility condition}.
Write
\[
 F=T^n-1,\qquad m_i=T^{d_i}-1,\qquad
 G=\lcm(m_1,\ldots,m_s),
\]
with the least common multiple monic. Use the notation of Theorem~\ref{thm:normal-count},
\(\varphi_q(H)=|(k[T]/(H))^\times|\).

\begin{theorem}[Joint trace fibers and normal lifts]\label{v58:thm-prescribed-normal-traces}
The conditions
\begin{equation}\label{v58:eq-prescribed-trace-system}
 \Tr_{E/L_i}(x)=a_i\qquad(1\le i\le s)
\end{equation}
admit a solution \(x\in E\) if and only if the prescribed values satisfy trace compatibility. When compatible, all solutions form an affine \(k\)-subspace of size
\[
 q^{n-\deg G}.
\]
This fiber contains a normal element of \(E/k\) if and only if every \(a_i\) is normal in \(L_i/k\).
In that case, the number of normal solutions is exactly
\begin{equation}\label{v58:eq-prescribed-normal-count}
 \frac{\varphi_q(F)}{\varphi_q(G)}.
\end{equation}
These conclusions do not require \(p\nmid n\).
\end{theorem}
\begin{proof}
Necessity of compatibility follows from transitivity of the trace; every intermediate trace of a normal solution is normal by Proposition~\ref{prop:normal-descent}.
We prove sufficiency and the counts.

Fix a normal element \(\eta\in E\). For each \(d\mid n\), put
\[
 \eta_d=\Tr_{E/\F_{q^d}}(\eta),\qquad
 H_d=\frac{F}{T^d-1}.
\]
The trace definition gives \(\eta_d=H_d(\sigma)\eta\), and normal descent ensures that \(\eta_d\) is normal in \(\F_{q^d}/k\).
Thus every \(a_i\) has a unique expression \(a_i=f_i(\sigma)\eta_{d_i}\), with \(f_i\) taken modulo \(m_i\).
Similarly, every \(x\in E\) has a unique expression \(x=f(\sigma)\eta\), with \(f\) taken modulo \(F\).
Since the coefficients of \(f\) lie in \(k\), the trace commutes with \(f(\sigma)\), so
\[
 \Tr_{E/L_i}(x)=f(\sigma)\eta_{d_i}.
\]
The Frobenius annihilator polynomial of \(\eta_{d_i}\) is \(m_i\), so~\eqref{v58:eq-prescribed-trace-system} is equivalent to
\begin{equation}\label{v58:eq-trace-polynomial-congruences}
 f\equiv f_i\pmod{m_i}\qquad(1\le i\le s).
\end{equation}

If \(e\mid d_i\), trace transitivity gives \(\Tr_{L_i/\F_{q^e}}(\eta_{d_i})=\eta_e\),
and hence \(\Tr_{L_i/\F_{q^e}}(a_i)=f_i(\sigma)\eta_e\). Thus~\eqref{v58:eq-trace-compatibility} is equivalent to
\[
 f_i\equiv f_j\pmod{T^{\gcd(d_i,d_j)}-1}.
\]
Applying the Euclidean algorithm with \(T^a-1-T^{a-b}(T^b-1)=T^{a-b}-1\) gives
\[
 \gcd(T^a-1,T^b-1)=T^{\gcd(a,b)}-1,
\]
so these are pairwise compatibility conditions modulo \(\gcd(m_i,m_j)\). This identity retains multiplicities when the characteristic divides an exponent.

To combine all congruences, for each monic irreducible polynomial \(P\) dividing \(G\),
choose an index \(i\) maximizing the multiplicity of \(P\) in \(m_i\), and use \(f_i\) to prescribe the residue class of \(f\) modulo that highest power.
Pairwise compatibility ensures that this prescription meets every other condition in the \(P\)-component.
Powers of distinct irreducibles are coprime, so the Chinese remainder theorem gives a unique residue class \(f_0\) modulo \(G\),
and this is precisely a solution of~\eqref{v58:eq-trace-polynomial-congruences}.
Since \(G\mid F\), all solutions modulo \(F\) are \(f=f_0+Gh\), with \(h\) ranging modulo \(F/G\).
There are \(q^{n-\deg G}\) such solutions, forming an affine \(k\)-subspace under the cyclic-module isomorphism.

Finally, count the normal solutions. Normality of each \(a_i\) is equivalent to \(\gcd(f_i,m_i)=1\).
Together these conditions are equivalent to \(\gcd(f_0,G)=1\): each \(P\mid G\) divides at least one \(m_i\),
and the corresponding \(f_i\) and \(f_0\) agree modulo \(P\).
Write \(F=\prod_P P^{u_P}\) and \(G=\prod_P P^{v_P}\), where \(0\le v_P\le u_P\).
Normal elements \(f(\sigma)\eta\) correspond to units modulo \(F\).
If \(v_P>0\), the fixed \(f_0\bmod P^{v_P}\) is already a unit;
it has \(q^{(u_P-v_P)\deg P}\) lifts modulo \(P^{u_P}\), all still units.
If \(v_P=0\), this component is unrestricted, giving \(\varphi_q(P^{u_P})\) choices of units.
Multiply these counts by the Chinese remainder theorem and use the formula for the polynomial Euler totient function~\eqref{eq:polynomial-euler}.
The total is exactly \(\varphi_q(F)/\varphi_q(G)\) and is strictly positive.
The argument retains all irreducible-factor multiplicities throughout, so it also applies when \(p\mid n\).
\end{proof}

\begin{corollary}[Normal elements of prescribed nonzero trace]\label{v58:cor-normal-fixed-trace}
For each \(c\in k^\times\), the number of normal elements of trace \(c\) is exactly
\[
 \frac{\varphi_q(T^n-1)}{q-1}.
\]
There are no normal elements of trace zero.
\end{corollary}
\begin{proof}
Take \(s=1\) and \(d_1=1\) in the theorem. The normal elements of the one-dimensional extension \(k/k\) are precisely the nonzero elements,
and \(G=T-1\), \(\varphi_q(T-1)=q-1\).
\end{proof}

\begin{example}[Two intermediate traces and one free field component]\label{v58:ex-normal-traces-five-six}
Take \(k=\F_5\), \(E=\F_{5^6}\), \(L_1=\F_{25}\), and \(L_2=\F_{125}\).
For prescribed \(a\in L_1\) and \(b\in L_2\), compatibility means
\[
 \Tr_{L_1/k}(a)=\Tr_{L_2/k}(b).
\]
Here
\[
 \begin{aligned}
 F&=(T-1)(T+1)(T^2+T+1)(T^2-T+1),\\
 G&=(T-1)(T+1)(T^2+T+1).
 \end{aligned}
\]
Both quadratic factors have discriminant \(2\), a nonsquare in \(\F_5\), and are therefore irreducible.
Thus \(\deg G=4\), and each compatible pair of target values has \(5^2=25\) common preimages.
If \(a,b\) are both normal over \(k\), the three components modulo \(G\) are already units,
and the remaining free component is \(k[T]/(T^2-T+1)\cong\F_{25}\).
Its \(24\) choices of nonzero values give exactly all normal lifts.

Such compatible normal target values do exist: for each prescribed nonzero trace, the two subfields contain, respectively,
\(\varphi_5(T^2-1)/4=4\) and \(\varphi_5(T^3-1)/4=24\) normal elements.
Choose the same nonzero trace and then a normal element in each subfield to obtain such \((a,b)\).
\end{example}

\begin{example}[Restrictions from normality and a prescribed norm]\label{v58:ex-normal-norm-nine}
In \(E=\F_9=\F_3(\alpha)\), \(\alpha^2=-1\), from Example~\ref{v100:ex-three-generation-nine}, we now compare all normal elements with the norm fibers.
For \(x=a+b\alpha\), the conjugate is \(a-b\alpha\).
In the basis \((1,\alpha)\), the coefficient matrix of these two vectors has determinant \(-2ab=ab\), so
\[
 x\text{ is normal}\quad\Longleftrightarrow\quad ab\ne0.
\]
The four normal elements are exactly \(\pm1\pm\alpha\), and all have norm
\[
 \Nm_{E/\F_3}(x)=a^2+b^2=2=-1.
\]
Thus the norm-\(1\) fiber is nonempty but contains no normal elements; the norm-\(-1\) fiber consists precisely of the four normal elements.
Prescribing intermediate traces amounts to polynomial reduction in normal coordinates, whereas prescribing a norm imposes a condition on field multiplication. They impose different restrictions on normality.
\end{example}

Prescribed intermediate traces specify quotient coordinates, whereas adding elements of an intermediate subfield gives a coset of an ideal in the cyclic module.
Exercise~\ref{v82:ex-subfield-translates} compares these two types of fibers through the same unit condition, and shows that a larger subfield need not offer more ways to repair normality.

\subsection{All annihilator data along an invariant direction}
Counting normal elements distinguishes units from nonunits. To determine which primary components an element loses within an affine fiber, we must also retain multiplicities.
Only the direction subspace is required to be Frobenius-invariant below. The coset \(\alpha+U\) itself is invariant only if \(\sigma(\alpha)-\alpha\in U\) as well.

\begin{theorem}[Annihilator distribution in an affine fiber]\label{v107:thm-affine-annihilator-distribution}
Let \(E=\F_{q^n}\), \(k=\F_q\), \(\sigma(x)=x^q\), and write
\[
F=T^n-1=\prod_P P^{e_P},\qquad
D=\mu_\alpha=\prod_P P^{d_P},\qquad H=\prod_P P^{h_P}\mid F.
\]
Products run over the distinct monic irreducible factors of \(F\), allowing zero exponents. Put \(U=\ker H(\sigma)\) and \(Q_P=q^{\deg P}\).
For \(G=\prod_P P^{g_P}\mid F\), the number of \(x\in\alpha+U\) with \(\mu_x=G\) is
\begin{equation}\label{v107:eq-affine-annihilator-distribution}
\boxed{\prod_P c_P(g_P),\qquad
c_P(g)=\begin{cases}
Q_P^{h_P},&d_P>h_P,\ g=d_P,\\
\varphi_q(P^g),&d_P\le h_P,\ 0\le g\le h_P,\\
0,&\text{otherwise}.
\end{cases}}
\end{equation}
We use \(\mu_0=1\) and \(\varphi_q(1)=1\). In particular, every \(\mu_x\) divides \(\lcm(D,H)\), and the number attaining this maximum is
\begin{equation}\label{v107:eq-maximal-annihilator-count}
q^{\deg H}\prod_{\substack{P\mid F\\h_P\ge d_P,\ h_P>0}}(1-q^{-\deg P}).
\end{equation}
\end{theorem}
\begin{proof}
Choose a normal element to identify \(E\), as a \(k[T]\)-module, with \(k[T]/(F)\), and apply the Chinese remainder decomposition.
In the local component \(R_P=k[T]/(P^{e_P})\), the direction \(U\) corresponds to the ideal \((P^{e_P-h_P})\), which contains \(Q_P^{h_P}\) elements.
If \(\alpha_P\ne0\), its \(P\)-valuation is \(e_P-d_P\); assign valuation \(e_P\) to zero by convention.

When \(d_P>h_P\), a perturbation has strictly higher valuation and cannot cancel the first nonzero layer of \(\alpha_P\). The local annihilator therefore remains \(P^{d_P}\).
When \(d_P\le h_P\), the center already lies in the ideal, so its translate is the whole ideal.
Within it, \(Q_P^g\) elements are killed by \(P^g\). For \(g\ge1\), subtracting the \(Q_P^{g-1}\) elements killed by \(P^{g-1}\) gives \(\varphi_q(P^g)\). For \(g=0\), only zero occurs.
The primary components are independent, giving the product distribution. Taking the largest exponent \(\max(d_P,h_P)\) in each component gives the final formula.
The local counts sum to \(Q_P^{h_P}\), so the total remains \(q^{\deg H}=|U|\).
\end{proof}

\begin{corollary}[Preserving a prescribed annihilator]\label{v107:cor-annihilator-preservation}
If also \(H\mid D\), the number of elements of \(\alpha+\ker H(\sigma)\) whose annihilator remains \(D\) is
\[
\boxed{q^{\deg H}\prod_{\substack{P\mid H\\v_P(H)=v_P(D)}}(1-q^{-\deg P}).}
\]
A loss can occur only when the perturbation reaches the full annihilator length of a primary component. Perturbations in lower layers all preserve that component.
Squarefree \(H\) recovers Proposition~\ref{v106:prop-annihilator-fiber}.
\end{corollary}
\begin{proof}
Take \(G=D\) in the theorem. Since \(h_P\le d_P\), a nonzero proportion \(1-q^{-\deg P}\) is needed only when \(h_P=d_P>0\).
\end{proof}

\begin{example}[Four-point and sixteen-point fibers in the field of 64 elements]\label{v107:ex-affine-sixteen}
Use \(b=1+uv\) from Example~\ref{v106:ex-complete-f64}, and put \(P=T^2+T+1\). Since \(b\) is normal, \(D=(T+1)^2P^2\).
Along \(\ker P(\sigma)\), all four elements remain normal. Enlarging the direction to \(U=\ker(P^2)(\sigma)\) gives
\[
\begin{array}{c|ccc}
\mu_x&(T+1)^2&(T+1)^2P&(T+1)^2P^2\\\hline
\text{number}&1&3&12\\
{[\F_2(x):\F_2]}&2&6&6
\end{array}
\]
Here \(P^2(\sigma)=P(\sigma)\circ P(\sigma)\).
The unique quadratic element is \(u^2\). Indeed, \(v+v^2+v^4=1\) gives \((1+\sigma^2+\sigma^4)b=u^2\), while \((1+\sigma^2+\sigma^4)u^2=u^2\), so \(b+u^2\in U\).
The three elements in the middle column still generate the entire field, but their conjugates span only four dimensions. Only the last column consists of normal elements.
A larger direction permits more annihilator data, while the proportion preserving the original data drops from one to \(3/4\).
\end{example}

\begin{corollary}[Normal elements in affine spaces with invariant direction]\label{v88:thm-invariant-affine-normal-count}
Write \(F=T^n-1=\prod_{i=1}^sP_i^{e_i}\), where the \(P_i\) are distinct
monic irreducible polynomials and \(d_i=\deg P_i\). Let \(U\subseteq E\)
be a Frobenius-invariant \(k\)-subspace, take \(\alpha_0\in E\), and put
\[
\pi_i:E\longrightarrow E/P_i(\sigma)E,
\qquad S(U)=\{i:\pi_i(U)\ne0\}.
\]
Then \(\dim_k E/P_i(\sigma)E=d_i\), and the number of normal elements in
\(\alpha_0+U\) is
\begin{equation}\label{v88:eq-invariant-affine-normal-count}
\boxed{\quad
\begin{cases}
0,&\text{if some }i\notin S(U)\text{ has }\pi_i(\alpha_0)=0,\\[2mm]
|U|\displaystyle\prod_{i\in S(U)}(1-q^{-d_i}),&\text{otherwise}.
\end{cases}\quad}
\end{equation}
An empty product is one. No assumption that \(p\nmid n\) is needed.
\end{corollary}
\begin{proof}
By Theorem~\ref{thm:invariant-duality}, write \(U=\ker H(\sigma)\) uniquely with \(H\mid F\).
In the component \(k[T]/(P_i^{e_i})\), put \(h_i=v_{P_i}(H)\); the direction corresponds to \((P_i^{e_i-h_i})\).
Its projection modulo \(P_i\) is nonzero exactly when \(h_i=e_i\), in which case it is the whole residue field. The quotient dimension is therefore \(d_i\).
Apply Theorem~\ref{v107:thm-affine-annihilator-distribution} with \(G=F\). A free component \(h_i=e_i\) contributes proportion \(1-q^{-d_i}\); every other component must already have \(v_{P_i}(\mu_{\alpha_0})=e_i\), equivalently \(\pi_i(\alpha_0)\ne0\).
Multiplying the local sizes gives the asserted formula.
\end{proof}

\begin{example}[Normal elements in an affine subset of the nine-element field]\label{v88:ex-affine-normal-nine}
Take \(E=\F_9=\F_3(\xi)\), \(\xi^2=-1\), \(U=k=\F_3\), and \(\alpha_0=\xi\).
Example~\ref{v58:ex-normal-norm-nine} already shows that \(a+b\xi\) is normal exactly when \(ab\ne0\).
Here we restrict this criterion to the affine fiber: the point \(\xi\) is not
normal, while \(\xi+1\) and \(\xi-1\) are, giving two elements.
Alternatively, apply Corollary~\ref{v88:thm-invariant-affine-normal-count}.
Since \(T^2-1=(T-1)(T+1)\), the affine set \(\xi+U\) takes every value in the
one-dimensional component corresponding to \(T-1\), while its component
corresponding to \(T+1\) is fixed and nonzero. The number of normal elements is therefore
\(|U|(1-q^{-1})=3(1-1/3)=2\).

Here \(\sigma(\xi)=-\xi\), so \(\sigma(\xi+k)=-\xi+k\ne\xi+k\). The direction \(k\) is invariant, but the coset is not. The count does not require the latter condition.
\end{example}

With \(U=\ker\Tr_{E/k}=(\sigma-1)E\) and \(\Tr(\alpha_0)=c\ne0\),
Corollary~\ref{v88:thm-invariant-affine-normal-count} gives \(q^{n-1}\prod_{P\ne T-1}(1-q^{-\deg P})\), exactly
\(\varphi_q(F)/(q-1)\) in Corollary~\ref{v58:cor-normal-fixed-trace}.
Prescribed trace is therefore a special case of the count for affine spaces with invariant direction. Proposition~\ref{v110:prop-affine-count} will treat arbitrary directions by intersections of invariant kernels.

\section{Relative Trace Reduction and Affine Spaces of Normal Elements}\label{v93:sec-normal-affine}

An intermediate trace fixes certain local residue classes, whereas normality requires every local residue class to be nonzero.
We first determine when ordinary intermediate traces detect all local components, then use trace duality to construct affine spaces consisting entirely of normal elements.
Throughout this section, write
\[
k=\F_q,\qquad E=\F_{q^n},\qquad Q=q^n,\qquad N=Q-1,\qquad
\sigma(x)=x^q,
\]
and let \(n=p^a m\), where \(p=\Char k\) and \(p\nmid m\). Put
\[
F=T^n-1=M^{p^a},\qquad
M=T^m-1=\prod_{i=1}^sP_i,\qquad d_i=\deg P_i.
\]
The polynomials \(P_i\) are distinct, monic, and irreducible, so \(\sum_i d_i=m\).
Fix a normal element \(\eta\). The cyclic-module isomorphism
\(A=k[T]/(F)\simeq E\), given by \(u\mapsto u(\sigma)\eta\), yields
\begin{equation}\label{v93:eq-normal-local}
u(\sigma)\eta\text{ is normal}
\quad\Longleftrightarrow\quad u\in A^\times
\quad\Longleftrightarrow\quad u\not\equiv0\pmod{P_i}\quad(1\le i\le s).
\end{equation}

\subsection{What Ordinary Intermediate Traces Detect}

\begin{theorem}[A single-trace criterion for the prime-to-characteristic part]\label{v93:thm-trace}
Let \(K=\F_{q^m}\) and \(\tau=\Tr_{E/K}\). For every \(x\in E\),
\begin{equation}\label{v93:eq-trace-normal}
\boxed{\quad x\text{ is normal in }E/k
\quad\Longleftrightarrow\quad
\tau(x)\text{ is normal in }K/k.\quad}
\end{equation}
Consequently, for each normal \(b\in K\), the trace fiber \(\tau^{-1}(b)\) contains \(q^{n-m}\) elements, all of which are normal in \(E/k\).
\end{theorem}
\begin{proof}
The relative trace is the operator \(S(\sigma)\), where
\[
S(T)=\frac{T^n-1}{T^m-1}=M(T)^{p^a-1}.
\]
Put \(\eta_m=\tau(\eta)=S(\sigma)\eta\).
Its Frobenius annihilating polynomial is \(F/\gcd(F,S)=M\), so \(\eta_m\) is normal in \(K/k\).
If \(x=u(\sigma)\eta\), then \(\tau(x)=u(\sigma)\eta_m\).
By the unit criterion, normality of \(x\) is equivalent to \(\gcd(u,F)=1\), while normality of \(\tau(x)\) is equivalent to \(\gcd(u,M)=1\).
These conditions are equivalent because \(F=M^{p^a}\).
Surjectivity of the trace and the dimension formula give the fiber size.
\end{proof}

When \(p\nmid n\), we have \(m=n\), and this trace is the identity map.
When \(n\) is a power of \(p\), we have \(m=1\), recovering the criterion that normality is equivalent to a nonzero trace to \(k\).

\begin{theorem}[The exact boundary for ordinary joint trace fibers]\label{v93:thm-ordinary-boundary}
Let \(r_1,\ldots,r_t\mid n\), where \(t\ge1\), and prescribe
\(c_j\in K_j=\F_{q^{r_j}}\), each normal in \(K_j/k\),
satisfying the trace compatibility condition~\eqref{v58:eq-trace-compatibility}. Put
\[
G=\lcm_{1\le j\le t}(T^{r_j}-1).
\]
Every element of the joint trace fiber
\(\{x:\Tr_{E/K_j}(x)=c_j\ (1\le j\le t)\}\) is normal if and only if
\begin{equation}\label{v93:eq-ordinary-boundary}
\boxed{\quad M\mid G.\quad}
\end{equation}
This is also equivalent to \(m\mid r_j\) for at least one \(j\).
If the condition fails, the same fiber contains both normal and nonnormal elements.
\end{theorem}
\begin{proof}
By Theorem~\ref{v58:thm-prescribed-normal-traces}, the joint trace conditions become
\(u\equiv u_0\pmod G\) in normal coordinates, and normality of the target values ensures that \(\gcd(u_0,G)=1\).
If \(M\mid G\), the residue class modulo every \(P_i\) is fixed and nonzero, so every lift is a unit.
If some \(P_i\nmid G\), the Chinese remainder theorem allows the residue modulo \(P_i\) to be chosen either zero or nonzero while preserving the condition modulo \(G\).
Taking all other unfixed components nonzero gives a nonnormal lift and a normal lift.

Since \(p\nmid m\), the polynomial \(T^m-1\) has \(m\) distinct roots in a splitting field.
These roots form a finite multiplicative subgroup of order \(m\), hence a cyclic group; a generator \(\zeta\) is a primitive \(m\)th root of unity.
If \(M\mid G\), then \(\zeta\) is a root of some \(T^{r_j}-1\), so \(m\mid r_j\).
Conversely, \(m\mid r_j\) immediately gives \(M\mid T^{r_j}-1\mid G\).
This argument includes \(m=1\).
\end{proof}

\begin{example}[Local components undetected by the ordinary trace]\label{v93:ex-trace-f128}
In \(E=\F_{2^7}\),
\[
T^7-1=(T+1)(T^3+T+1)(T^3+T^2+1).
\]
Neither cubic factor has a root in \(\F_2\), so both are irreducible.
The only proper subfield is \(\F_2\).
The trace-one fiber contains \(64\) elements, of which \((2^3-1)^2=49\) are normal.
Since \(127\) is prime, every element of \(E\) other than zero and one is primitive.
The element one also has trace one, so this fiber contains \(63\) primitive elements, \(14\) of which are not normal.
The ordinary trace fixes the component corresponding to \(T+1\), but does not detect the other two cubic components.
\end{example}

\begin{proposition}[Linear information determining an exact annihilator]\label{v93:prop-normal-linear-info}\label{v110:prop-annihilator-info}
For a monic divisor $D\mid F$, put
\[
U_D=\ker D(\sigma),\qquad \mathcal S_D=\{x\in E:\mu_x=D\},
\qquad R_D=\prod_{P\mid D}P,
\]
where the product uses distinct monic irreducibles and $R_1=1$. The \textbf{additive translation stabilizer} of this set is
\begin{equation}\label{v110:eq-annihilator-stabilizer}
\{v\in E:\mathcal S_D+v=\mathcal S_D\}=R_D(\sigma)U_D.
\end{equation}
Membership in $\mathcal S_D$ is constant on every fiber of a $k$-linear map $L:E\to V$ if and only if
$\ker L\subseteq R_D(\sigma)U_D$. Consequently, the smallest possible rank of such a map is
\begin{equation}\label{v110:eq-annihilator-information}
\boxed{n-\deg D+\deg R_D.}
\end{equation}
In particular, the normal set $\mathcal N(E/k)=\mathcal S_F$ satisfies
\begin{equation}\label{v93:eq-normal-translation-stabilizer}
\{v\in E:\mathcal N(E/k)+v=\mathcal N(E/k)\}
=M(\sigma)E=\ker\tau.
\end{equation}
Thus determining normality on every fiber requires rank at least $m$, attained by $L=\tau$.
\end{proposition}
\begin{proof}
If $D=1$, then $\mathcal S_D=U_D=\{0\}$, its stabilizer is zero, and membership in this singleton requires an injective linear map. The assertions follow.
For $D\ne1$, choose $\eta_D$ with annihilator $D$ as in Section~\ref{sec:q-order}. The cyclic-module isomorphism
\[
k[T]/(D)\longrightarrow U_D,\qquad \bar u\longmapsto u(\sigma)\eta_D
\]
identifies $\mathcal S_D$ with the unit set. Any translation preserving $\mathcal S_D$ lies in $U_D$, because it is the difference of two elements of this nonempty set.
In $k[T]/(D)$, translation by an element of $(R_D)$ preserves every nonzero residue modulo $P\mid D$, hence preserves all units. Conversely, if $v\notin(R_D)$, some residue $v\bmod P$ is nonzero. The Chinese remainder theorem gives a unit $u$ with $u\equiv-v\pmod P$. Then $u+v$ is not a unit, so this translation does not preserve the unit set.
Its stabilizer is therefore $(R_D)$, which corresponds to $R_D(\sigma)U_D$ and has dimension $\deg D-\deg R_D$.

Fibers of $L$ are cosets of $\ker L$. Membership in $\mathcal S_D$ is constant on them precisely when every translation in $\ker L$ preserves $\mathcal S_D$. Rank--nullity gives the lower bound; the quotient map $E\to E/R_D(\sigma)U_D$ attains it.
For $D=F$, we have $R_D=M$ and $U_D=E$. Since $F=M^{p^a}$ and $\tau$ is represented by multiplication by $M^{p^a-1}$, its kernel corresponds to $(M)$. This proves the final assertions.
\end{proof}

The rank in~\eqref{v110:eq-annihilator-information} counts both exclusion of points outside $U_D$ and detection of its generating elements. If the input is already known to lie in $U_D$, the corresponding minimum is only $\deg R_D$.
For $E=\F_{27}$ over $\F_3$, the strata $D=(T-1)^j$, $j=0,1,2,3$, require ranks $3,3,2,1$, respectively, when the input ranges over all of $E$.
Thus testing the largest conjugate span can require less linear information than testing a smaller exact span. This concerns membership in one prescribed stratum, not recovery of the element or of all its annihilator data at once. Proposition~\ref{v130:prop-annihilator-spectrum} recovers this amount of information from the exact Fourier support of the stratum indicator.

\subsection{Weighted Traces and Affine Spaces of Normal Elements}

The rank bound above requires every fiber to determine normality.
To construct a single fiber consisting entirely of normal elements, it suffices to select one nonzero linear condition from each local component.
Put
\[
U_i=P_i(\sigma)E,\qquad
V_i=U_i^\perp=\{b\in E:\Tr_{E/k}(bu)=0\text{ for all }u\in U_i\}.
\]
Proposition~\ref{prop:normal-test} gives \(\dim_kV_i=d_i\) and
\begin{equation}\label{v93:eq-local-trace-dual}
V_1\oplus\cdots\oplus V_s=(M(\sigma)E)^\perp.
\end{equation}
Consequently, any choice of one nonzero vector from each \(V_i\) gives \(s\) linearly independent vectors over \(k\).
\phantomsection\label{v93:def-weighted-trace}
The linear functional \(x\mapsto\Tr_{E/k}(bx)\) is a \textbf{weighted trace}, with weight \(b\in E\).
\phantomsection\label{v93:def-normal-affine-space}
A nonempty affine \(k\)-subspace is called an \textbf{affine space of normal elements}
if every one of its elements is normal in \(E/k\).
This condition requires every element of the affine space to be normal over the fixed base field \(k\).
It does not require normality over every intermediate base field, which is the condition defining a completely normal element (see p.~\pageref{v82:def-completely-normal}).

\begin{theorem}[Weighted trace fibers of normal elements and their optimal dimension]\label{v93:thm-patch}
Choose any \(b_i\in V_i\setminus\{0\}\), and define
\begin{equation}\label{v93:eq-patch}
\mathcal A_{\boldsymbol b}
=\{x\in E:\Tr_{E/k}(b_ix)=1\quad(1\le i\le s)\}.
\end{equation}
Then \(\mathcal A_{\boldsymbol b}\) is an affine space of normal elements of dimension \(n-s\).
Every affine space of normal elements is contained in some \(\mathcal A_{\boldsymbol b}\).
In particular, every normal element belongs to such a fiber, every affine space of normal elements has dimension at most \(n-s\),
and those of dimension \(n-s\) are precisely the fibers above.
\end{theorem}
\begin{proof}
By~\eqref{v93:eq-local-trace-dual} and nondegeneracy of the trace pairing,
the \(k\)-linear map \(x\mapsto(\Tr(b_1x),\ldots,\Tr(b_sx))\) has rank \(s\) and is therefore onto \(k^s\).
Thus \(\mathcal A_{\boldsymbol b}\) is nonempty and has dimension \(n-s\).
An element \(x\) of this fiber lies outside every \(U_i\), since membership in \(U_i\) would make the \(i\)th trace zero.
Equation~\eqref{v93:eq-normal-local} then shows that \(x\) is normal.

Let \(B=x_0+W\) be any affine space of normal elements, and let \(\pi_i:E\to E/U_i\) be the quotient map.
Its image \(\pi_i(x_0)+\pi_i(W)\) does not contain zero, so \(\pi_i(x_0)\notin\pi_i(W)\).
Choose a linear functional on the quotient that vanishes on \(\pi_i(W)\) and takes the value one at \(\pi_i(x_0)\).
Pulling it back to \(E\) and using nondegeneracy of the trace pairing, write it uniquely as \(x\mapsto\Tr(b_ix)\),
where \(b_i\in V_i\setminus\{0\}\). It follows that \(B\subseteq\mathcal A_{\boldsymbol b}\).
Taking \(B\) to be a singleton proves the covering assertion.
Comparing dimensions gives the bound, and inclusion is equality when the dimensions agree.
\end{proof}

Determining normality on every fiber requires at least \(m=\sum_i d_i\) linear coordinates,
whereas constructing one affine space of normal elements requires only \(s\) weighted trace conditions.
The former retains every local residue class; the latter fixes just one nonzero linear value on each local component.
These direction spaces need not be Frobenius-invariant, so they are not subject to the direction restriction in the preceding result for affine spaces with invariant direction.

\paragraph{Exact annihilator counts under arbitrary linear constraints.}
Invariant directions allow independent counts in primary components. A general direction may couple those components; the intersections with invariant kernels still give an exact answer.
\begin{proposition}[Annihilator counts in an arbitrary affine space]\label{v110:prop-affine-count}
Let $D\mid F=T^n-1$ be monic, let $\mathcal P(D)$ be its set of distinct monic irreducible factors, and let $B=x_0+W$, where $W\subseteq E$ is any $k$-subspace. For $h\mid F$, put $U_h=\ker h(\sigma)$ and
\[
c_B(h)=
\begin{cases}
q^{\dim_k(W\cap U_h)},&x_0\in W+U_h,\\
0,&x_0\notin W+U_h.
\end{cases}
\]
Then
\begin{equation}\label{v110:eq-affine-count}
\boxed{\#\{x\in B:\mu_x=D\}
=\sum_{J\subseteq\mathcal P(D)}(-1)^{|J|}
c_B\left(D/\prod_{P\in J}P\right).}
\end{equation}
Empty products equal one. In particular, for a $k$-linear map $L:E\to V$ and a prescribed value $b\in V$, replace $c_B(h)$ by
\begin{equation}\label{v110:eq-linear-constraint-count}
c_{L,b}(h)=
\begin{cases}
q^{\deg h-\operatorname{rank}_k(L|_{U_h})},&b\in L(U_h),\\
0,&b\notin L(U_h).
\end{cases}
\end{equation}
This gives the same count in $L^{-1}(b)$, including when that fiber is empty.
\end{proposition}
\begin{proof}
The elements with annihilator exactly $D$ are those in $U_D$ outside every $U_{D/P}$, for $P\in\mathcal P(D)$. Indeed, every proper divisor of $D$ divides $D/P$ for at least one such $P$.
For a subset $J$, the intersection of these excluded kernels is
\[
\bigcap_{P\in J}U_{D/P}=U_{D/\prod_{P\in J}P},
\]
where for $J=\varnothing$ the intersection is taken inside $U_D$ and equals $U_D$.
This follows from the gcd identity for invariant kernels; dividing by $P$ removes only one occurrence, so repeated factors are retained.
Inclusion--exclusion within $B\cap U_D$ now proves~\eqref{v110:eq-affine-count}: the intersection $B\cap U_h$ is nonempty precisely when $x_0\in W+U_h$, and then it is a coset of $W\cap U_h$.
For a prescribed linear value, apply rank--nullity to $L|_{U_h}$ and use $\dim_k U_h=\deg h$. This proves~\eqref{v110:eq-linear-constraint-count}.
\end{proof}

Taking $B=E$ recovers $\varphi_q(D)$. Summing over the divisors $D\mid F$ of degree $n-j$ counts $j$-normal elements under the same constraints.
The earlier product formula is more efficient for invariant directions; the formula here also applies when a direction is not invariant, without assuming independence of the local conditions.

\begin{example}[Independent trace coordinates and dependent generation conditions]\label{v110:ex-trace-coordinates-eight}
Let $E=\F_8=\F_2(v)$, with $v^3+v^2+1=0$, as in Example~\ref{v106:ex-affine-f8}. Put $P=T^2+T+1$ and $F=(T+1)P$.
The identities $v+v^2+v^4=1$, $\Tr(v)=1$ and $\Tr(v^3)=0$ show that $(v,v^2,v^4)$ is self-dual: its diagonal trace pairings are one, and its off-diagonal pairings are zero by Frobenius invariance.
Consequently, for $x=rv+sv^2+tv^4$, the two weighted traces $\Tr(vx)$ and $\Tr(v^2x)$ are precisely $r$ and $s$.
The four fibers of this rank-two map have the following annihilator distributions:
\begin{center}\small
\begin{tabular}{@{}ccrrrr@{}}\toprule
$(r,s)$&Elements&$1$&$T+1$&$P$&$F$\\\midrule
$(0,0)$&$0,\ v^4$&1&0&0&1\\
$(1,0)$&$v,\ v+v^4$&0&0&1&1\\
$(0,1)$&$v^2,\ v^2+v^4$&0&0&1&1\\
$(1,1)$&$v+v^2,\ 1$&0&1&1&0\\\bottomrule
\end{tabular}\end{center}
Indeed, zero has annihilator one, the element one has annihilator $T+1$, the three nonzero trace-zero elements have annihilator $P$, and $v,v^2,v^4$ have annihilator $F$.
Every fiber has direction $\langle v^4\rangle_{\F_2}$, which is not Frobenius-invariant because $\sigma(v^4)=v$.
The trace coordinates are independent on $E$, but on its three normal elements each event $r=1$ or $s=1$ has proportion $1/3$, while their intersection is empty.
The linear rank controls the size of the ambient fibers; the intersections with annihilator kernels control the generating elements within them.
\end{example}

\subsection{Imposing Primitivity as Well}

\phantomsection\label{def:primitive-normal-element}
An element that generates \(E^\times\) and is normal in \(E/k\) is called a \textbf{primitive normal element}.
Let \(\operatorname{PN}_q(n)\) be the number of such elements.
For \(b\in K\), let \(P_E(b;\tau)\) be the number of elements of multiplicative order \(N\) in \(\tau^{-1}(b)\).
Theorem~\ref{v93:thm-trace} gives
\begin{equation}\label{v93:eq-one-trace-count}
\operatorname{PN}_q(n)=
\sum_{b\in\mathcal N(K/k)}P_E(b;\tau).
\end{equation}
Weighted trace fibers can overlap, but every normal element occurs with the same multiplicity, giving another exact counting formula.

\begin{corollary}[Counting primitive normal elements through weighted trace fibers]\label{v93:cor-exact-pn}
Let \(P(\mathcal A_{\boldsymbol b})\) be the number of primitive elements in \(\mathcal A_{\boldsymbol b}\). Then
\begin{equation}\label{v93:eq-exact-pn}
\boxed{\quad
\operatorname{PN}_q(n)=q^{s-m}
\sum_{\substack{b_i\in V_i\setminus\{0\}\\1\le i\le s}}
P(\mathcal A_{\boldsymbol b}).\quad}
\end{equation}
\end{corollary}
\begin{proof}
Fix a normal element \(x\).
Since \(x\notin U_i=V_i^\perp\), the map \(b\mapsto\Tr(bx)\) is a nonzero linear functional on \(V_i\).
The equation \(\Tr(b_ix)=1\) has exactly \(q^{d_i-1}\) solutions, all automatically nonzero.
Thus \(x\) lies in exactly \(\prod_iq^{d_i-1}=q^{m-s}\) fibers indexed by weight tuples.
A nonnormal element belongs to none of these fibers.
Restricting to primitive elements and counting in the two orders proves the formula.
\end{proof}

Equations~\eqref{v93:eq-one-trace-count} and~\eqref{v93:eq-exact-pn} express the joint count through primitive-element counts in ordinary and weighted trace fibers, respectively.
Section~\ref{v93:sec-primitive-normal-proof} will use Gauss sums to estimate the multiplicative conditions and count over the local components together, proving existence for all parameters.

\subsection{Scalar Symmetry of Exact-Order Elements and Trace Distributions}
\phantomsection\label{v113:par-trace-question}
The preceding counts impose multiplicative and linear generation on the same elements. Even once such elements are known to exist, two further questions remain: can their trace be prescribed, and do different traces occur equally often? The primitive normal basis theorem guarantees ordinary primitive normal elements, but does not by itself determine the size of each prescribed-trace fiber or cover arbitrary combinations of multiplicative order and conjugate-span dimension. An ordinary normal element has nonzero trace, so its zero-trace fiber must be distinguished from the nonzero ones.

A total count determines only an average. Example~\ref{v111:ex-order52-traces} will give a concrete comparison in \(\F_{5^4}\): elements of multiplicative order 52 and conjugate-span dimension three have counts
\[
(4,0,0,4)
\]
at the nonzero traces \(1,2,3,4\), respectively. The average is two elements per nonzero trace, yet two fibers are empty. Thus a positive total does not ensure that every fiber is nonempty, and the average does not establish uniformity. The example will verify these counts from six irreducible factors.

We first seek a way to compare fibers without counting them separately. If a nonzero base-field scalar preserves both the prescribed multiplicative order and the conjugate-span dimension, multiplication by it gives a bijection between the corresponding trace fibers. We therefore determine the multiplicative stabilizer of the exact-order set and then intersect it with the base-field scalars. This identifies fibers that must have equal sizes before asking whether different cosets still have different counts.

In this subsection, reset the notation to $K=\F_q$, $E=\F_{q^n}$, $N=q^n-1$ and $M=N/(q-1)$, with $n\ge1$.
\phantomsection\label{v111:def-r-primitive}
For $r\mid N$, an element of multiplicative order exactly $N/r$ is called an \textbf{$r$-primitive element}; $r=1$ recovers an ordinary primitive element.
This specifies only the multiplicative order, without imposing field degree or normality. Using the definition of $k$-normality in Theorem~\ref{v109:thm-knormal-trace}, put
\[
\mathcal A_{r,k}=\{\alpha\in E^\times:\ord(\alpha)=N/r,\ \deg\mu_\alpha=n-k\},
\qquad 0\le k<n,
\]
\begin{equation}\label{v111:eq-trace-count}
N_{r,k}(a)=\#\{\alpha\in\mathcal A_{r,k}:\Tr_{E/K}(\alpha)=a\}\qquad(a\in K).
\end{equation}
We call $a\mapsto N_{r,k}(a)$ the corresponding \textbf{trace distribution}. For the following symmetry, see~\cite{ChoudharyPisolkarTrace2026}. We prove it directly from exact-order elements in a cyclic group and the Chinese remainder theorem.

\begin{theorem}[The stabilizer of an exact-order set and trace cosets]\label{v111:thm-order-trace-symmetry}
Put $R=N/r$, with $\rad(1)=1$, and set
\[
D_r=r\rad(R),\qquad t=\frac{D_r}{\gcd(D_r,M)}.
\]
If $E^\times=\langle g\rangle$, the full multiplicative stabilizer of the set $\mathcal P_r$ of elements of order $R$ is
\begin{equation}\label{v111:eq-order-stabilizer}
\{\lambda\in E^\times:\lambda\mathcal P_r=\mathcal P_r\}
=\langle g^{D_r}\rangle.
\end{equation}
Its intersection with $K^\times$ is
\[
H_r=\langle g^{tM}\rangle=(K^\times)^t,\qquad t\mid q-1.
\]
For every $0\le k<n$, $a\in K$ and $\lambda\in H_r$,
\begin{equation}\label{v111:eq-trace-cosets}
N_{r,k}(\lambda a)=N_{r,k}(a).
\end{equation}
Thus the distribution of nonzero traces is constant on each multiplicative coset of $H_r$. If $D_r\mid M$, it is constant on all of $K^\times$.
\end{theorem}
\begin{proof}
We have $\mathcal P_r=\{g^{ru}:\gcd(u,R)=1\}$. If $D_r\mid c$, write $c=rv$. Every prime $\ell\mid R$ then divides $v$, so $\gcd(u,R)=1$ if and only if $\gcd(u+v,R)=1$. Hence $g^c\mathcal P_r=\mathcal P_r$.
Conversely, if this equality holds, $g^r\in\mathcal P_r$ implies $r\mid c$; write $c=rv$.
If a prime $\ell\mid R$ fails to divide $v$, the Chinese remainder theorem supplies an integer $u$ such that
\[
u\equiv-v\pmod\ell,\qquad
u\equiv1\pmod{\ell'}\quad(\ell'\mid R,\ \ell'\ne\ell\text{ prime}).
\]
Then $u$ is coprime to $R$ but $u+v$ is not, a contradiction. Thus $\rad(R)\mid v$, proving~\eqref{v111:eq-order-stabilizer}.
When $R=1$, both $\mathcal P_r$ and its stabilizer are $\{1\}$, and the same formula applies.

Since $K^\times=\langle g^M\rangle$, a scalar $g^{jM}$ preserves $\mathcal P_r$ precisely when $D_r\mid jM$, equivalently $t\mid j$.
Moreover, $D_r\mid N=(q-1)M$ gives $t\mid q-1$, proving the formula for $H_r$.
For $\lambda\in K^\times$, the identities $(\lambda\alpha)^{q^i}=\lambda\alpha^{q^i}$ show that scalar multiplication preserves the entire Frobenius annihilator, and hence $k$-normality.
If $\lambda\in H_r$, it also preserves the prescribed multiplicative order, while $\Tr(\lambda\alpha)=\lambda\Tr(\alpha)$.
Multiplication by $\lambda$ and by $\lambda^{-1}$ therefore gives inverse bijections between the two trace fibers, proving~\eqref{v111:eq-trace-cosets}.
\end{proof}

\begin{corollary}[Uniform nonzero traces of primitive $k$-normal elements]\label{v111:cor-primitive-trace-uniform}
For $r=1$,
\[
t=\prod_{\substack{\ell\mid q-1\\\ell\nmid n}}\ell
=\frac{\rad(q-1)}{\gcd(\rad(q-1),n)}.
\]
In particular, if $\rad(q-1)\mid n$, then $N_{1,k}(a)$ is constant on $K^\times$; under the same hypothesis,
\[
N_{1,0}(a)=\frac{\operatorname{PN}_q(n)}{q-1}\qquad(a\ne0).
\]
\end{corollary}
\begin{proof}
For $\ell\mid q-1$, we have $M=1+q+\cdots+q^{n-1}\equiv n\pmod\ell$.
If $\ell\mid N$ but $\ell\nmid q-1$, then $\ell\mid M$.
Together with $D_1=\rad(N)$, these observations give the formula for $t$. The stated divisibility condition makes $H_1=K^\times$.
A normal element has nonzero trace: its annihilator $T^n-1$ cannot divide the trace-operator polynomial $(T^n-1)/(T-1)$.
The primitive normal elements therefore split into $q-1$ equally large nonzero trace fibers, giving the last formula.
\end{proof}

The divisibility condition guarantees uniformity but does not characterize it. The stabilizer action does not compare the zero-trace fiber with the nonzero fibers, and the common count may be zero.

\begin{example}[Equal multiplicative order and trace, different conjugate-span dimensions]\label{v111:ex-order52-traces}
Take $q=5$, $n=4$ and $r=12$. Then $N=624$, $R=52$, $M=156$, $D_{12}=312$ and $t=2$, so $H_{12}=\{1,4\}$.
Over $\F_5$, the polynomial $\Phi_{52}=\sum_{j=0}^{12}(-1)^jX^{2j}$ is the product of the six factors below.
Multiplication verifies the factorization; since $\ord_{52}(5)=4$, every factor is irreducible.
Its roots have order 52 and trace equal to the negative of its cubic coefficient. The last column is the rank of $X,X^5,X^{25},X^{125}$ in the corresponding quotient ring.
\begin{center}\small
\begin{tabular}{@{}lcc@{}}\toprule
Factor&Trace&Conjugate-span dimension\\\midrule
$X^4+X^3-X+1$&4&3\\
$X^4+X^3-X^2-X+1$&4&4\\
$X^4+2X^3+X^2-2X+1$&3&4\\
$X^4-2X^3+X^2+2X+1$&2&4\\
$X^4-X^3+X+1$&1&3\\
$X^4-X^3-X^2+X+1$&1&4\\\bottomrule
\end{tabular}\end{center}
For a concrete calculation, take roots $\alpha,\beta$ of the last two factors. In their respective power bases $(1,z,z^2,z^3)$, the conjugate coordinate columns form
\[
C_\alpha=\begin{pmatrix}0&4&1&1\\1&3&0&1\\0&4&4&2\\0&1&1&3\end{pmatrix},\qquad
C_\beta=\begin{pmatrix}0&4&1&1\\1&3&4&2\\0&0&4&1\\0&2&1&2\end{pmatrix}.
\]
The upper-left $3\times3$ minor of $C_\alpha$ is 3, and its fourth row is four times its third; meanwhile $\det C_\beta=3$.
Their ranks are therefore three and four, respectively, and
\[
\alpha^{125}-2\alpha^{25}-\alpha^5+2\alpha=0,\qquad
\mu_\alpha=\frac{T^4-1}{T-3},\quad \mu_\beta=T^4-1.
\]
Each row supplies four conjugate roots. Thus all 24 elements of order 52 have the trace distribution
\[
\begin{array}{c|rrrrr}
a&0&1&2&3&4\\\hline
N_{12,0}(a)&0&4&4&4&4\\
N_{12,1}(a)&0&4&0&0&4
\end{array}
\]
All other values of $k$ give zero counts. Both functions are constant on $\{1,4\}$ and on $\{2,3\}$, but only the first is constant on all nonzero traces.
Multiplicative order, field degree, conjugate-span dimension and trace record different information. Even the full stabilizer of an exact-order set need not be the full symmetry group of its trace-counting function.
\end{example}

Section~\ref{v111:sec-trace-fourier} will express these coset symmetries as Fourier support constraints: the group action first restricts the permitted frequencies, and the actual counts determine which of their coefficients vanish.

\section{The Complete Criterion for Self-Dual Normal Bases}
This section treats normality and the trace pairing within one algebra. The key reduction expresses changes of normal basis as units in the cyclic-module algebra
and self-duality as a unit norm equation for an involution. Solvability is determined by residue-field norms and lifting through nilpotent ideals.
\phantomsection\label{def:involution}
A map whose square is the identity is an \textbf{involution}; the identity map itself is allowed.
\leavevmode\phantomsection\label{def:quotient-involution}\ignorespaces Put \(R=k[T]/(T^n-1)\) and define \(u^*(T)=u(T^{-1})\) on \(R\).
This map preserves the algebra operations and satisfies \((u^*)^*=u\), so it is an algebra involution.
Since \(T^n=1\), we have \(T^{-1}=T^{n-1}\) here, making the definition meaningful.
This involution does not normalize to monic form; if \(h\in k[T]\) is monic and \(h(0)\ne0\), then in \(R\),
\[
\overline{h^\vee}=h(0)^{-1}T^{\deg h}\bar h^*.
\]
Also, \(\sigma^*=\sigma^{-1}\) shows that the trace adjoint of \(u(\sigma)\) is \(u^*(\sigma)\).
Monic reciprocation, the inversion involution in the quotient ring, and operator adjoints thus arise from the same inversion action, each with its own normalization.

Fix a normal element \(\eta\) and put
\[
g(T)=\sum_{i=0}^{n-1}\Tr(\eta\sigma^i(\eta))T^i\in R.
\]
Then \(g^*=g\) and \(g(1)=\Tr(\eta)^2\ne0\).
Multiplication by \(g\) has the symmetric circulant Gram matrix of the normal basis, so \(g\) is a unit of \(R\).
If \(u=\sum_a u_aT^a\in R^\times\) and \(\xi=u(\sigma)(\eta)\), write \(c_i=\Tr(\eta\sigma^i(\eta))\). Then
\[
\Tr(\xi\sigma^i(\xi))=\sum_{a,b}u_au_b c_{i+b-a}.
\]
The right side is exactly the coefficient of \(T^i\) in \(uu^*g_\eta\) reduced modulo \(T^n-1\), so
\begin{equation}\label{eq:normal-gram-change}
g_\xi=u u^* g_\eta.
\end{equation}
Existence of a self-dual normal element is therefore equivalent to solvability of the unit norm equation \(uu^*=g^{-1}\).
We first solve modulo irreducible factors and then handle the errors arising in lifts to their powers.
The inversion involution may exchange two components or preserve one; direct pairing and finite-field norms handle these cases, respectively.

\begin{lemma}\leavevmode\phantomsection\label{lem:odd-circulant-norm}
Suppose \(n\) is odd, \(g\in R^\times\), \(g^*=g\), and \(g(1)\) is a square in \(k^\times\).
There exists \(u\in R^\times\) such that \(uu^*=g^{-1}\).
\end{lemma}
\begin{proof}
First suppose \(q\) is even. Since \(n\) is odd, \(T^n-1\) has no repeated factors,
and the Chinese remainder theorem decomposes \(R\) as a finite product of finite fields. Squaring is bijective in every component,
so there is a unique \(u\in R\) satisfying \(u^2=g^{-1}\), necessarily a unit.
Applying the involution gives \((u^*)^2=g^{-1}\). Uniqueness implies \(u^*=u\), hence \(uu^*=g^{-1}\).

Now suppose \(q\) is odd. Write \(n=p^a m\) with \(p\nmid m\) and put \(e=p^a\). Then both \(e,m\) are odd, and
\[
T^n-1=(T^m-1)^e,\qquad
R\cong\prod_{P\mid T^m-1}R_P,\quad R_P=k[T]/(P^e).
\]
The product runs over all distinct monic irreducible factors of \(T^m-1\).
By the local-ring model of Section~\ref{sec:local-trace-rank}, each \(R_P\) has unique maximal ideal
\(\mathfrak m=(\overline P)\), residue field \(K_P=k[T]/(P)\), and \(\mathfrak m^e=0\).
An element is a unit if and only if its image in \(K_P\) is nonzero.
The involution \(T\mapsto T^{-1}\) takes the component for \(P\) to that for its monic reciprocal \(P^\vee\).
We solve along the component orbits of this involution.

\noindent\textbf{Step 1: Pair the exchanged components directly.}
If \(P\ne P^\vee\), the involution induces mutually inverse isomorphisms between the two components.
Since \(g^*=g\), we have \(g_{P^\vee}=g_P^*\). In this pair of components, choose
\[
u_P=g_P^{-1},\qquad u_{P^\vee}=1.
\]
The two components of \(uu^*\) are then \(g_P^{-1}\) and \((g_P^{-1})^*=g_{P^\vee}^{-1}\). The equation already holds on the pair of local rings; no further lifting is needed.

\noindent\textbf{Step 2: Solve the residue-field equation on fixed components.}
Suppose \(P=P^\vee\) and first solve in the residue field \(K_P=k[T]/(P)\).
If \(P=T-1\), the residue-field involution is the identity and the equation becomes \(\bar u^2=g(1)^{-1}\), solvable by the square hypothesis.
Since \(m\) and the characteristic are odd, \(T+1\) is not a factor of \(T^m-1\).
For the other self-reciprocal factors, Proposition~\ref{prop:self-reciprocal} gives \(\deg P=2d\),
and the inverse of the residue class of \(T\) is its \(q^d\)th power. The involution on all of \(K_P\) is therefore the \(q^d\)th-power map,
with fixed field \(\F_{q^d}\). The condition \(g^*=g\) ensures that \(\bar g\) belongs to this fixed field.
Thus
\[
\bar u\bar u^*=\bar g^{-1}
\]
is the norm equation for a quadratic finite-field extension and is solvable by norm surjectivity.

\noindent\textbf{Step 3: Lift through the nilpotent ideal.}
In a fixed component \(R_P\), lift the nonzero residue-field solution from Step 2 arbitrarily to \(u_0\in R_P\).
Its nonzero image modulo \(P\) makes \(u_0\) a unit. Put
\[
h=\frac{g^{-1}}{u_0u_0^*}.
\]
Commutativity and \(g^*=g\) give \(h^*=h\). The residue-field equation says \(h\equiv1\pmod P\),
so \(h=1+r\) for some \(r^e=0\). Since \(e=p^a\) is odd,
the Frobenius identity gives \(h^e=1+r^e=1\), and hence
\[
v=h^{(e+1)/2},\qquad v^2=h^{e+1}=h,\qquad v^*=v.
\]
Thus \(u=u_0v\) is a unit satisfying \(uu^*=u_0u_0^*h=g^{-1}\).
This corrects the error in the arbitrary lift without changing its residue-field solution.

Combine the solutions from Step 1 on exchanged components and Step 3 on fixed components by the Chinese remainder theorem
to obtain an element of \(R^\times\) satisfying \(uu^*=g^{-1}\).
\end{proof}

\begin{theorem}[Realizing prescribed cyclic correlations in odd degree]\label{v88:thm-prescribed-correlation}
Assume that \(n\) is odd and put \(R=k[T]/(T^n-1)\). For \(b\in R\),
there exists a normal element \(\alpha\in E/k\) with \(g_\alpha=b\) if and only if
\begin{equation}\label{v88:eq-prescribed-correlation}
\boxed{\quad b\in R^\times,\qquad b^*=b,\qquad
b(1)\in(k^\times)^2.\quad}
\end{equation}
If \(q\) is even, every nonzero base-field element is a square; the last
condition then follows from the unit condition.
\end{theorem}
\begin{proof}
If \(b=g_\alpha\), Proposition~\ref{prop:normal-circulant-test} shows that \(b\)
is a unit. Symmetry of the trace pairing gives \(b^*=b\), while
Proposition~\ref{v88:prop-conjugate-correlation} yields
\(b(1)=\Tr(\alpha)^2\in(k^\times)^2\).
Conversely, fix a normal element \(\eta\) and let \(g=g_\eta\).
Then \(g\) is a self-adjoint unit with \(g(1)=\Tr(\eta)^2\).
The element \(c=gb^{-1}\) is again a self-adjoint unit, and
\(c(1)=g(1)/b(1)\) is a nonzero square. By
Lemma~\ref{lem:odd-circulant-norm}, there is \(u\in R^\times\) satisfying
\(uu^*=c^{-1}=bg^{-1}\). The element \(\alpha=u(\sigma)\eta\) remains normal,
and~\eqref{eq:normal-gram-change} gives \(g_\alpha=uu^*g=b\).
\end{proof}

This theorem prescribes all \(n\) trace-pairing values, not merely the relative trace $\Tr_{E/k}(\alpha)$. Taking \(b=1\) recovers the existence of self-dual normal bases in odd
degree; other admissible \(b\) describe all invertible symmetric trace
autocorrelations of normal Frobenius orbits.
\begin{example}[Every normal basis of the eight-element field is self-dual]\label{v88:ex-f8-normal-selfdual}
Example~\ref{ex:normal-f8-example} constructs a normal basis of the eight-element field; here we determine its trace pairing. Let \(E=\F_8\), \(k=\F_2\). Every normal element \(\alpha\) has
\(\Tr(\alpha)=1\). Set \(c=\Tr(\alpha\alpha^2)\).
By cyclic symmetry and \(\Tr(\alpha^2)=1\),
\(g_\alpha=1+c(T+T^2)\).
If \(c=1\), then \(g_\alpha=1+T+T^2\) is not coprime to \(T^3-1\),
contradicting the normality test. Hence \(c=0\) and \(g_\alpha=1\),
so every normal basis of \(\F_8/\F_2\) is self-dual.
\end{example}

\begin{theorem}[The self-dual normal basis criterion]\label{thm:self-dual-normal}
The extension \(\F_{q^n}/\F_q\) admits a self-dual normal basis if and only if
\[
\boxed{\quad n\text{ is odd, or }q\text{ is even and }n\equiv2\pmod4.\quad}
\]
\end{theorem}
\begin{proof}
Suppose \(n\) is odd. Choose a normal element \(\eta\) by the normal basis theorem.
Its cyclic correlation polynomial \(g_\eta\) is a unit satisfying \(g_\eta^*=g_\eta\) and \(g_\eta(1)=\Tr(\eta)^2\ne0\).
Lemma~\ref{lem:odd-circulant-norm} gives a unit \(u\) with \(uu^*=g_\eta^{-1}\).
Put \(\xi=u(\sigma)\eta\).
The quotient-ring model of a normal basis shows that the action of a unit preserves normality; formula~\eqref{eq:normal-gram-change} gives \(g_\xi=1\), meaning \(\Tr(\xi\sigma^i(\xi))=\delta_{0i}\).
Shifting both factors gives \(\Tr(\sigma^j(\xi)\sigma^l(\xi))=\delta_{jl}\), so the resulting normal basis is self-dual.

Suppose \(q\) is even, \(n=2m\), and \(m\) is odd.
In \(\F_{q^2}\), trace surjectivity permits a choice of \(b\) with \(b+b^q=1\).
Base-field elements have relative trace zero, so \(b\notin k\).
If \(cb+db^q=0\) with \(c,d\in k\), taking traces gives \(c+d=0\);
the original equation then becomes \(c(b+b^q)=c=0\), so \(c=d=0\).
Thus \((b,b^q)\) is a normal basis. Moreover,
\[
\Tr(b^2)=\Tr(b)^2=1,\qquad
\Tr(bb^q)=2\,bb^q=0.
\]
Frobenius invariance shows that the other diagonal pairing is also one, so this basis is self-dual.
Take a self-dual normal basis in \(\F_{q^m}/k\) from the odd-degree case
and combine it using \(\gcd(2,m)=1\) and Proposition~\ref{prop:normal-coprime} to obtain the required basis.

For necessity, there are two cases. If \(q\) is odd and \(n\) even, Theorem~\ref{thm:self-dual-basis} already excludes a self-dual basis.
If \(q\) is even and \(4\mid n\), suppose a self-dual normal basis generator exists.
Take its relative trace to the subfield \(\F_{q^4}\). Proposition~\ref{prop:normal-descent} gives a self-dual normal element \(b\) of \(\F_{q^4}/k\).
Self-duality gives \(\Tr(b^2)=1\), while \(\Tr(b^2)=\Tr(b)^2\). Squaring is injective in characteristic two, so \(\Tr(b)=1\).
Put \(c=b+b^{q^2}\). Then \(c\in\F_{q^2}\) and \(c+c^q=1\).
On the other hand, expanding and regrouping the trace gives
\[
\begin{aligned}
0&=\Tr_{\F_{q^4}/k}(bb^q)\\
 &=bb^q+b^qb^{q^2}+b^{q^2}b^{q^3}+b^{q^3}b\\
 &=(b+b^{q^2})(b^q+b^{q^3})=cc^q.
\end{aligned}
\]
A field has no zero divisors, so \(c=0\), contradicting \(c+c^q=1\).
These cases exhaust the parameters not satisfying the criterion.
\end{proof}

For corresponding generalizations to abelian Galois extensions, see \cite{BL}.

\section{Primitive Normal Elements: Two Generation Properties}\label{v48:sec-primitive-normal}
A primitive normal element generates both the multiplicative group and the Frobenius cyclic module.
Let \(k=\F_q\) and \(E=\F_{q^n}\), and write \(\sigma(x)=x^q\), \(N=q^n-1\), and \(F=T^n-1\).
A primitive normal element \(\eta\) simultaneously satisfies
\[
E^\times=\langle\eta\rangle,\qquad
E=\bigoplus_{i=0}^{n-1}k\,\sigma^i(\eta).
\]
The first equality expresses every nonzero element as a multiplicative power; the second expresses every element of the extension as a linear combination of conjugates.
Both imply \(E=k(\eta)\). Field generation only requires Frobenius orbit length \(n\); it does not require multiplicative order \(N\) or linear independence of the orbit elements.

\begin{theorem}[Primitive normal basis theorem]\leavevmode\phantomsection\label{thm:primitive-normal}
For every prime power \(q\) and positive integer \(n\), the extension \(\F_{q^n}/\F_q\) has a primitive normal element.
\end{theorem}
For complete proofs, see Lenstra and Schoof~\cite{LS} and Cohen and Huczynska~\cite{CohenHuczynska2003}.
Section~\ref{v93:sec-primitive-normal-proof} gives a proof based on the local decomposition and trace duality developed here, together with character orthogonality and Gauss sums.
The arithmetic argument covering all parameters is supplied in full in Appendix~\ref{v93:chap-primitive-normal-arithmetic}.

\begin{proposition}[Primitive elements in quadratic extensions are normal]\label{v88:prop-primitive-quadratic-normal}
For every prime power \(q\), every primitive element of \(\F_{q^2}\) is
normal over \(\F_q\).
\end{proposition}
\begin{proof}
If a primitive \(\alpha\) were not normal, then \(\alpha^q=c\alpha\) for
some \(c\in k^\times\), so \(\alpha^{q-1}=c\).
Since \(c^q=c\) and \(\alpha^{q^2-1}=1\), we have
\(c^{q+1}=c^2=1\), hence \(\alpha^{2(q-1)}=1\).
This contradicts \(\ord(\alpha)=(q-1)(q+1)>2(q-1)\).
\end{proof}

\begin{corollary}[Primitive normal and completely normal elements in two degree families]\label{v88:cor-primitive-normal-special-degrees}
Let \(q=p^m\). If \(n=p^a\), then \(\F_{q^n}/\F_q\) has a primitive
completely normal element. If \(n=\ell\ne p\) is prime with
\(\ord_\ell(q)=\ell-1\), then it has a primitive normal element.
\end{corollary}
\begin{proof}
Theorem~\ref{v88:thm-primitive-nonzero-trace} supplies a primitive element
\(\alpha\) of nonzero trace. Being primitive, \(\alpha\) has full extension
degree, so Theorem~\ref{thm:trace-normal-classification} makes it normal.
For \(n=p^a\), Proposition~\ref{v82:prop-completely-normal-p-power}
makes it completely normal. The case \(n=1\) follows as well.
\end{proof}

In these special degrees, the factor structure of the Frobenius module reduces normality to a trace condition,
so the character-sum argument of Chapter~\ref{chap:characters} directly produces a primitive normal element.
In general, all local components must be treated together; the joint count in Section~\ref{v93:sec-primitive-normal-proof} completes this argument.

\subsection{Two Unit Coordinate Systems and a Joint Criterion}
Fix a primitive element \(\zeta\) and a normal element \(\theta\). The two generation problems have coordinate models
\begin{equation}\label{v48:eq-two-generator-models}
\begin{aligned}
\Z/N\Z&\longrightarrow E^\times,\qquad j\longmapsto\zeta^j,\\
R=k[T]/(F)&\longrightarrow E,\qquad\bar u\longmapsto u(\sigma)\theta.
\end{aligned}
\end{equation}
The first is an isomorphism from an additive group to a multiplicative group; the second is a \(k[T]\)-module isomorphism, with \(T\) acting on \(E\) as \(\sigma\).
The order formula in a cyclic group and the cyclic-module model in this chapter give
\[
\begin{aligned}
\zeta^j\text{ is primitive}&\ \Longleftrightarrow\ \gcd(j,N)=1,\\
u(\sigma)\theta\text{ is normal}&\ \Longleftrightarrow\ \gcd(u,F)=1.
\end{aligned}
\]
The two types of generators therefore correspond to \((\Z/N\Z)^\times\) and \(R^\times\), respectively.
The models preserve different operations: exponent coordinates turn field multiplication into addition of exponents, while normal coordinates turn Frobenius into a cyclic shift.
The second isomorphism generally does not preserve field multiplication, and the first does not describe field addition.
Coprimality in the two coordinate systems is therefore not interchangeable.

\begin{proposition}[A joint criterion for primitive normal elements]\label{v48:prop-primitive-normal-test}
For \(a\in E^\times\), the following two conditions hold simultaneously if and only if \(a\) is primitive and normal:
\begin{enumerate}
\item For every distinct prime divisor \(\ell\mid N\), we have \(a^{N/\ell}\ne1\);
\item For every distinct monic irreducible factor \(P\mid F\), we have \((F/P)(\sigma)(a)\ne0\).
\end{enumerate}
When \(N=1\), the first set of conditions is empty.
\end{proposition}
\begin{proof}
Let \(e=\ord(a)\), so \(e\mid N\). If \(e<N\), choose a prime divisor \(\ell\) of \(N/e\). Then \(e\mid N/\ell\), hence \(a^{N/\ell}=1\).
Conversely, if this equality holds, \(e\mid N/\ell<N\), so \(a\) is not primitive.
The second set is exactly the normality criterion of Proposition~\ref{prop:normal-test}. Combining them proves the assertion.
\end{proof}

The first set excludes elements of multiplicative order less than \(N\); the second excludes elements whose conjugates span a proper subspace.
The multiplicative condition excludes the subgroup \(\{a\in E^\times:a^{N/\ell}=1\}\), of order \(N/\ell\);
the normality condition excludes the subspace \(\ker(F/P)(\sigma)\). In the cyclic-module model, this corresponds to the ideal \((\bar P)\) of \(R\), of dimension \(n-\deg P\).
When \(p\mid n\), the polynomial \(F\) may have repeated factors, but \(F/P\) reduces the multiplicity of \(P\) by only one, so the second set still tests the full primary structure.
Taking \(P=T-1\) gives \(\Tr_{E/k}(a)\ne0\); the trace condition is generally only one normality condition. Theorem~\ref{thm:trace-normal-classification} describes when nonzero trace, together with degree \(n\) over \(k\), is sufficient for normality.

\subsection{Changing Generators and Intersecting the Two Conditions}
If \(\gcd(j,N)=1\), the power map \(a\mapsto a^j\) is an automorphism of \(E^\times\) and therefore takes primitive elements to primitive elements.
If \(\bar u\in R^\times\), the operator \(u(\sigma)\) is an invertible \(k\)-linear map taking normal elements to normal elements.
Indeed, if \(\mu_a=F\), the cyclic-module model gives
\(\mu_{u(\sigma)a}=F/\gcd(F,u)=F\).
The first type of change need not preserve normality, and the second need not preserve multiplicative order.
Frobenius preserves both: in exponent coordinates it is multiplication by \(q\), and in normal coordinates multiplication by \(T\); each multiplier is a unit in its coordinate ring.
The identity \(\sigma M_a\sigma^{-1}=M_{a^q}\) in the regular representation of Section~\ref{sec:regular-representation} also shows how these actions are compatible within one extension field.

\begin{example}[Two generation properties in the field of order 16]\leavevmode\phantomsection\label{v48:ex-primitive-normal-f16}
In \(E=\F_{16}=\F_2(\alpha)\) with \(\alpha^4+\alpha+1=0\), Exercise~\ref{ex:wu-three-generators} and its solution give three four-element Frobenius orbits:
\begin{center}
\begin{tabular}{@{}ccccc@{}}
Representative&Multiplicative order&Degree&Primitive&Normal\\\hline
\(\alpha\)&15&4&Yes&No\\
\(\alpha^3\)&5&4&No&Yes\\
\(\alpha^{-1}\)&15&4&Yes&Yes
\end{tabular}
\end{center}
The other three nonzero elements lie in \(\F_4^\times\) and are neither primitive nor normal.
In particular, inversion takes the primitive but nonnormal \(\alpha\) to the primitive normal \(\beta=\alpha^{-1}\), so a power map preserving primitivity can change normality.
On the other hand, \(F=(T+1)^4\), and \(u=T+T^2+T^3\) satisfies \(u(1)=1\), hence \(\bar u\in R^\times\). By \(\Tr(\beta)=1\) and \(\alpha^{-1}=\alpha^3+1\),
\[
u(\sigma)\beta=\beta^2+\beta^4+\beta^8
=1+\beta=\alpha^3.
\]
This change of normal basis preserves normality but lowers multiplicative order from 15 to five.
\end{example}

Counting primitive and normal elements separately gives
\[
\#\mathcal P=\varphi(N),\qquad
\#\mathcal N=\varphi_q(F),
\]
where \(\mathcal P,\mathcal N\subseteq E^\times\) denote the respective sets.
Both counts come from coprimality in unit coordinates, but these two individual counts do not determine the size of the intersection.
For example, in the preceding example \(\#\mathcal P=\#\mathcal N=8\), whereas \(\#(\mathcal P\cap\mathcal N)=4\).
For a uniformly chosen element of \(E^\times\), independence of the two conditions would give intersection size \(8\cdot8/15=64/15\), different from the actual count.
Computing the intersection requires further control of intersections between multiplicative subgroups and kernels of Frobenius operators.
Separate existence of the two types of generators only says that both sets are nonempty; the primitive normal basis theorem guarantees that a single element can satisfy both sets of exclusion conditions.
Such an element provides exponent coordinates for the nonzero elements and normal coordinates for the whole field. The relation between the two coordinate systems still depends on the specific field operations.

\par\Needspace{11\baselineskip}

\paragraph{Testing normality before sieving for primitivity.}
\phantomsection\label{v82:rem-primitive-normal-translates}
When candidates are restricted to \(\theta+k\), Proposition~\ref{v82:prop-normal-translates} first determines whether an obstruction to normality persists throughout the translate.
If the family contains no normal element, it contains no primitive normal element. If it does contain a normal element, put
\[
 P_\theta=\#\{t\in k:\theta+t\text{ is primitive}\},\qquad
 PN_\theta=\#\{t\in k:\theta+t\text{ is primitive and normal}\}.
\]
When \(p\nmid n\), also put \(\beta_\theta=\theta-\Tr_{E/k}(\theta)/n\). Then
\begin{equation}\label{v82:eq-primitive-normal-translates}
 PN_\theta=
 \begin{cases}
 P_\theta,&p\mid n,\\
 P_\theta-\mathbf1_{\{\beta_\theta\text{ is primitive}\}},&p\nmid n.
 \end{cases}
\end{equation}
The indicator is one when \(\beta_\theta\) is primitive and zero otherwise; zero is not counted as a primitive element.
The formula assumes that the translate \emph{contains a normal element}; it does not dispense with checking the normality obstruction.
It turns the joint condition into a concrete procedure: first test generation in the cyclic module, and then apply the multiplicative-order test to the remaining candidates.
It gives no estimate for \(P_\theta\), nor does it assert that every translate contains a primitive element.

\paragraph{The common structure of cyclic modules, trace duality, and local counts.}
The commutative algebra \(A=k[T]/(T^n-1)\) explains several results of this
chapter. In normal coordinates, \(T\) acts by Frobenius cyclic shift, and
normal elements correspond to units of \(A\). After factoring
\(T^n-1=\prod_P P^{e_P}\), normality means that every local component is
nonzero modulo \(P\). This explains normal-element counts, prescribed
intermediate traces, and the exact count for affine spaces with invariant direction.
The trace pairing gives the involution \(T\mapsto T^{-1}\) and the correlation
identity \(g_\alpha=h_\alpha h_\alpha^*\); self-dual normal bases and
prescribed correlations are governed by unit norm equations \(uu^*=b\).
The additive and multiplicative structures must not be conflated:
primitivity is still a multiplicative-order condition in \(E^\times\) and
does not follow from being a unit in the cyclic module. Their simultaneous
satisfaction is precisely the difficulty addressed by the primitive normal
basis theorem.

\section*{Exercises}
\phantomsection
\addcontentsline{toc}{section}{Exercises}

\begin{exercise}[A cyclic formula for the Moore determinant]\label{v88:ex-moore-sign}
In \(E=\F_4/\F_2\), compute
\(\det\begin{pmatrix}x&x^2\\x^2&x\end{pmatrix}\) for arbitrary \(x\in E\),
express it using \(\Tr_{4/2}(x)\), and explain why this is a special case
of~\eqref{v88:eq-moore-trace-power}.
\end{exercise}
\begin{exercise}[Quantifiers in the nonzero-trace criterion]\label{v88:ex-trace-quantifier}
Prove that \(1\) has nonzero trace but is not normal in \(\F_8/\F_2\).
Then prove that every element of degree exactly three and nonzero trace
in that extension is normal. Identify the hypothesis that distinguishes
the two statements.
\end{exercise}
\begin{exercise}[Normal elements in translates over the nine-element field]\label{v88:ex-affine-nine-cosets}
In \(E=\F_3(\xi)\), with \(\xi^2=-1\), determine all normal elements and their number in \(b\xi+\F_3\) for each \(b\in\F_3\).
\end{exercise}
\begin{exercise}[All correlation polynomials in a binary cubic extension]\label{v88:ex-correlation-f8}
Determine all invertible \(b\in\F_2[T]/(T^3-1)\) satisfying
\(b^*=b\) and \(b(1)=1\). Deduce without enumerating \(\F_8\)
that every normal basis of \(\F_8/\F_2\) is self-dual.
\end{exercise}
\begin{exercise}[Primitive normal elements in the nine-element field]\label{v88:ex-primitive-normal-nine}
Let \(E=\F_3(\xi)\), \(\xi^2=-1\). Determine every primitive element
and every normal element of \(E\), show that they form the same set of four,
and compare with Proposition~\ref{v88:prop-primitive-quadratic-normal}.
\end{exercise}

\begin{exercise}[Two normal traces with repeated factors]\label{v58:ex-normal-traces-binary-six}
Take \(E=\F_{64}\) and \(k=\F_2\), and let \(a\in\F_4\) and \(b\in\F_8\) both be normal over \(k\).
Prove that they automatically satisfy trace compatibility, and count all solutions of
\[
 \Tr_{E/\F_4}(x)=a,\qquad \Tr_{E/\F_8}(x)=b
\]
and the normal elements among them.
Prove also that \(x\in E\) is normal if and only if both intermediate traces above are normal.
Explain why this sufficiency does not contradict the nondiagonalizability of Frobenius on \(E/k\).
\end{exercise}

\begin{exercise}\label{ex:normal-count8}
Count the normal elements of \(\F_8/\F_2\) and verify that the three conjugates in Example~\ref{ex:normal-f8-example} give all of them.
\end{exercise}
\begin{exercise}[Frobenius-invariant hyperplanes]\label{ex:normal-p-power}
Let \(E=\F_{q^n}\) and \(k=\F_q\). Prove that \(E\) has exactly \(\gcd(n,q-1)\) Frobenius-invariant \(k\)-hyperplanes, and describe all of them using the trace pairing.
In particular, if \(n=p^a\), prove that the only invariant hyperplane is \(\ker\Tr_{E/k}\).
\end{exercise}
\begin{exercise}\label{ex:normal-selfdual8}
Verify that the normal basis in Example~\ref{ex:normal-f8-example} is self-dual.
Then determine whether \(\F_{2^4}/\F_2\) and \(\F_{3^4}/\F_3\) have self-dual bases and self-dual normal bases.
\end{exercise}
\begin{exercise}\label{ex:minpoly-vs-qorder}
In \(\F_8/\F_2\), find the ordinary minimal polynomials and Frobenius annihilator polynomials of \(\alpha\) and \(\eta=\alpha+1\) from Example~\ref{ex:normal-f8-example},
and explain the distinction between the two types of polynomial.
\end{exercise}

\begin{exercise}\label{ex:supp-linear-primitive}
Let \(h(T)=T^2+T+1\in\F_2[T]\). Find \(\mathcal L_h\) and use Frobenius periods on the root space to prove that \(\mathcal L_h/X\) is irreducible.
\end{exercise}
\begin{exercise}[Do a full primary factor and maximal period force normality?]\label{ex:supp-nonnormal-six}
Let \(k=\F_q\), \(E=\F_{q^n}\), and \(n=p^a m\), where \(p\nmid m\), and write \(F=T^n-1\).
Fix a normal element \(\eta\in E\).
\begin{enumerate}
\item For each nonconstant monic divisor \(\mu\mid F\), put \(\alpha_\mu=(F/\mu)(\sigma)\eta\). Prove that
\[
\mu_{\alpha_\mu}=\mu,\qquad
[k(\alpha_\mu):k]=\ord(\mu),\qquad
\dim_k k[\sigma]\alpha_\mu=\deg\mu.
\]
Then prove that \(\Tr_{E/k}(\alpha_\mu)\ne0\) if and only if \((T-1)^{p^a}\mid\mu\).
\item Determine whether the following implication holds for every \(q,n\):
\[
\mu\mid F,\qquad (T-1)^{p^a}\mid\mu,\qquad \ord(\mu)=n
\quad\Longrightarrow\quad \mu=F.
\]
Use Theorem~\ref{thm:trace-normal-classification} to determine all parameters for which it holds for every \(\mu\),
and explain how the proper divisors constructed in its proof yield explicit nonnormal elements.
Check this conversion against Example~\ref{exm:full-degree-nonnormal-f64}.
\end{enumerate}
\end{exercise}

\begin{exercise}[Three types of generators in the field of order 16]\leavevmode\phantomsection\label{ex:wu-three-generators}
Let \(\alpha^4+\alpha+1=0\) in \(\F_{16}\).
Prove that \(\alpha\) is primitive but not normal, \(\alpha^3\) is normal but not primitive, and \(\alpha^{-1}\) is both primitive and normal.
Count all normal, primitive, and primitive normal elements.
\end{exercise}

\begin{exercise}[Existence of invariant complements]\label{ex:wu2-invariant-complements}
Factor \(T^n-1\) as \(\prod_iP_i^{e_i}\). Prove that the total number of Frobenius-invariant subspaces is \(\prod_i(e_i+1)\).
Then prove that every invariant subspace has an invariant direct complement if and only if \(p\nmid n\).
\end{exercise}

\begin{exercise}[Two extremes for linearized root sets]\label{ex:wu7-rootfield-extremes}
Let \(f\in\F_q[T]\) be monic of degree \(d\ge1\) with \(f(0)\ne0\), and let \(W_f\) be the space of all roots of \(\mathcal L_f\).
Prove that \(W_f\) is a subfield if and only if \(f=T^d-1\).
Then prove that its splitting-field degree is at most \(q^d-1\), with equality if and only if \(f\) is a primitive polynomial.
In the latter case, explain why \(W_f\setminus\{0\}\) forms a single Frobenius orbit.
Compare the information conveyed in these two cases by \(d\), the root-set cardinality, and the splitting-field degree.
\end{exercise}

\begin{exercise}[Norm fibers and similarity classes]\label{wfr:ex-norm-classes}
Take \(E=\F_{64}\), \(k=\F_2\), and \(\tau(x)=x^4\). For \(a\in E^\times\), put \(S_a=a\tau\).
Determine the \(k\)-linear similarity classes of these operators, the number of weights \(a\) in each class, and the total number of weights giving a cyclic vector.
\end{exercise}

\begin{exercise}[The norm need only generate the field]\label{wfr:ex-nonprimitive-norm}
Take \(k=\F_3\), \(E=\F_{81}\), \(\tau(x)=x^9\), and \(K=\F_9\), and choose \(i\in K\) with \(i^2=-1\).
Prove that there exists \(a\in E^\times\) with \(\Nm_{E/K}(a)=i\). Find the minimal and characteristic polynomials of \(S_a=a\tau\) over \(k\),
and explain why it has a cyclic vector although the norm \(i\) is not a multiplicative generator of \(K\).
\end{exercise}

\begin{exercise}[Digit rotation and normal coordinates]\label{v40:ex-necklace-models}
Let \(m\ge2\). Equivalence classes of words of length \(m\) under cyclic shift are called \textbf{necklaces}\leavevmode\phantomsection\label{v40:def-necklace};
a word that is not a repeated concatenation of a shorter word is called a \textbf{primitive word}\leavevmode\phantomsection\label{v40:def-primitive-word}.
\begin{enumerate}
\item Write \(0\le a<q^m-1\) as a fixed-length base-\(q\) digit word. Prove that multiplication by \(q\) modulo \(q^m-1\) corresponds to cyclic rotation.
Explain why the all-\(q-1\) word must be excluded, and prove that orbit length equals the least cyclic period of the digit word.
\item Choose a normal basis of \(\F_{q^m}/\F_q\) and a bijection of sets between \(\F_q\) and \(\{0,\ldots,q-1\}\).
Prove that Frobenius likewise corresponds to cyclic rotation of coordinate words.
For every \(r\mid m\) and \(r>1\), prove that the number of length-\(r\) orbits in both models is
\[
I_q(r)=\frac1r\sum_{e\mid r}\mu(r/e)q^e,
\]
whereas the numbers of length-one orbits are \(q-1\) and \(q\), respectively.
\item For \(q=2,m=6\), list the orbit counts at every length in both models and check the total number of elements.
Explain why the comparison depends on the chosen normal basis and alphabet labeling and therefore does not establish an algebraic isomorphism between the digit and field-coordinate models.
\end{enumerate}
\end{exercise}

\begin{exercise}[Extensions in which every full-degree element has a normal translate]\label{v82:ex-universal-normal-translation}
Let \(k=\F_q\), \(E=\F_{q^n}\), \(n\ge2\), and \(p=\Char k\). Prove that the property
\[
 \text{for every }\theta\text{ with }k(\theta)=E,
 \text{ some }t\in k\text{ makes }\theta+t\text{ normal}
\]
holds if and only if one of the following occurs:
\[
 n=\ell\ne p\text{ is prime and }\ord_\ell(q)=\ell-1;
 \qquad n=p=2.
\]
\begin{enumerate}
\item For \(p\nmid n\), compare this property with Theorem~\ref{thm:trace-normal-classification}: translation preserves full degree, and the trace is the only normality condition that can change.
\item For \(p\mid n\), use the cyclic module to construct full-degree elements that are not normal, excluding every case other than \(n=p=2\).
\item Explain directly why \(n=p=2\) works, and state the additional boundary case when \(n=1\) is allowed.
\end{enumerate}
\end{exercise}

\begin{exercise}[Primary-component obstructions to subfield translation]\label{v82:ex-subfield-translates}
Let \(k=\F_q\), \(E=\F_{q^n}\), \(d\mid n\), and \(L=\F_{q^d}\), and put
\[
 F=T^n-1,\qquad D=T^d-1,\qquad M=F/D.
\]
Fix a normal element \(\eta\), and write \(\theta=A(\sigma)\eta\) in the normal coordinates \(k[T]/(F)\simeq E\).
\begin{enumerate}
\item Prove that \(L=\ker(\sigma^d-\id)\) corresponds to the ideal \((\overline M)\) in this cyclic module.
\item For a monic irreducible factor \(P\mid F\), determine whether \(\theta+L\) can change the unit condition modulo \(P\).
Prove that the translate contains a normal element if and only if \(\gcd(A,M)=1\), and count its normal elements exactly, retaining all multiplicities in \(F\) and \(D\).
\item With \(p=\Char k\), prove that if \(v_p(d)<v_p(n)\), then
\[
 \theta+u\text{ is normal}\quad\Longleftrightarrow\quad\theta\text{ is normal}
 \qquad(u\in L).
\]
\item For \(E=\F_{64}\) and \(k=\F_2\), compare \(d=2\) with \(d=3\).
In particular, let \(\theta=\eta+\eta^2\). Prove that \(\theta\) has full degree but is not normal, and count the normal elements in \(\theta+\F_4\) and in \(\theta+\F_8\).
\end{enumerate}
\end{exercise}

\begin{exercise}\label{v145:ex-right-factors}
Let $r\mid s$, $K=\F_{Q^s}$, and $\sigma(c)=c^Q$. For each $0\le j\le r$, determine the number of monic degree-$j$ right factors of $t^r-1$ in $K[t;\sigma]$.
\end{exercise}

\begin{exercise}[Right-factor lattices over one return target]\label{v148:ex-factor-lattices}
Let $\pi\ne Z$ be monic irreducible over $\F_q$, and let $e,k\ge1$. Among monic skew polynomials with $\mathcal N(f)=\pi^e$, determine when every $f$ has its monic right factors totally ordered by right divisibility. Determine also when at least one such $f$ has this property.
\end{exercise}

\chapter{Frobenius Tests, Factorization, and Changes of Representation}
\label{chap:algorithms}
Irreducibility criteria test orbit lengths, fixed subalgebras separate distinct orbits, and the minimal polynomial of a fixed element determines the actual splitting values.
Root-set descent thus becomes a finite sequence of linear-algebra and polynomial operations. Value fibers of fixed elements also yield idempotent projections; constructing normal elements compares the different effects of translation and reciprocation on linear relations among conjugates.

\section{Frobenius Criteria for Irreducibility}
For the following criterion and its algorithmic background, see Rabin~\cite{Rabin}.

\begin{theorem}[Rabin's criterion]\label{thm:rabin}
Let \(f\in\F_q[X]\) be monic with \(\deg f=n\geq1\). Then \(f\) is irreducible if and only if
\[
X^{q^n}\equiv X\pmod f
\]
and, for every distinct prime divisor \(\ell\) of \(n\),
\[
\gcd\bigl(f,X^{q^{n/\ell}}-X\bigr)=1.
\]
When \(n=1\), the second set of conditions is empty.
\end{theorem}
\begin{proof}
If \(f\) is irreducible, Theorem~\ref{thm:divisibility-main} gives all the conditions directly.
Conversely, the first condition says \(f\mid X^{q^n}-X\), so \(f\) is square-free and the degree \(d\) of every irreducible factor divides \(n\).
If \(d<n\), choose a prime divisor \(\ell\) of \(n/d>1\); then \(d\mid n/\ell\).
This factor would divide \(X^{q^{n/\ell}}-X\), contradicting the greatest-common-divisor condition.
Thus every factor has degree \(n\); since \(\deg f=n\), there can be only one.
\end{proof}

\begin{proposition}[Successive degree tests]\label{prop:ben-or}
A monic polynomial \(f\) of degree \(n\) with \(n\ge1\) is irreducible if and only if
\[
\gcd(f,X^{q^i}-X)=1\qquad(1\leq i\leq\lfloor n/2\rfloor).
\]
\end{proposition}
\begin{proof}
If \(f\) is irreducible, its degree \(n\) divides none of \(1\le i\le\lfloor n/2\rfloor\), so \(f\nmid X^{q^i}-X\) and every indicated greatest common divisor is one.
If \(f\) is reducible, its irreducible factorization, counted with multiplicity, contains a factor of degree \(i\le n/2\).
That factor divides \(X^{q^i}-X\), making the corresponding greatest common divisor nontrivial.
\end{proof}

These criteria do not require expanding polynomials of degree \(q^n\).
Starting with \(u_0=X\bmod f\), compute \(u_{i+1}=u_i^q\bmod f\) successively to obtain the remainder of \(X^{q^i}\).
\leavevmode\phantomsection\label{v101:def-binary-exponentiation}\ignorespaces Expanding an exponent in binary and using repeated squaring and multiplication gives \textbf{binary exponentiation}.
Computing a \(q\)th power this way requires \(O(\log q)\) modular polynomial multiplications.
Here \(O(g)\) means that, for sufficiently large parameters, the operation count is bounded by a constant multiple of \(g\); Section~\ref{sec:construct-normal-dual} specifies the cost model.
Rabin's criterion requires greatest common divisors only at the specified \(i=n/\ell\) and a congruence check at \(i=n\).

\begin{example}
For \(f=X^4+X+1\in\F_2[X]\), reduction modulo \(f\) gives
\[
X^4\equiv X+1,\qquad X^8\equiv X^2+1,\qquad X^{16}\equiv X.
\]
Since \(X^4-X\equiv1\pmod f\), the required greatest common divisor is one. Rabin's criterion proves that \(f\) is irreducible.
\end{example}

\begin{proposition}[Testing primitive polynomials]\label{prop:primitive-test-poly}
Let \(f\) be monic irreducible of degree \(n\), with \(f(0)\ne0\). Put \(N=q^n-1\).
Then \(f\) is primitive if and only if, for every distinct prime divisor \(\ell\) of \(N\),
\[
X^{N/\ell}\not\equiv1\pmod f.
\]
\end{proposition}
\begin{proof}
In \(\F_q[X]/(f)\), the order of \(\overline X\) divides \(N\).
Every proper divisor \(e\mid N\) divides some \(N/\ell\), so the stated conditions say precisely that this order is \(N\).
\end{proof}

\begin{corollary}[Testing irreducibility and primitivity simultaneously]\label{cor:primitive-test-unconditional}
Let \(f\) be monic of degree \(d\ge1\) with \(f(0)\ne0\), and put \(N=q^d-1\). Then \(f\) is primitive if and only if
\[
X^N\equiv1\pmod f,\qquad
X^{N/\ell}\not\equiv1\pmod f\quad\text{for every prime divisor }\ell\mid N.
\]
\end{corollary}
\begin{proof}
These conditions say exactly that the quotient-ring unit \(\bar X\) has order \(N\); now apply Proposition~\ref{prop:maximal-order-primitive}.
\end{proof}
If irreducibility has not been checked beforehand, the first congruence cannot be omitted: excluding periods \(N/\ell\) alone does not ensure that the period divides \(N\).

\section{Square-Free and Distinct-Degree Factorization}\label{sec:preliminary-factorization}
Square-free factorization first reduces the factorization of a general polynomial to that of square-free polynomials. One may then apply Berlekamp's method directly, or first group the irreducible factors by degree and then split each equal-degree part by a randomized method.

\subsection{Separating Multiplicities}\label{sec:squarefree-factorization}
Finite fields are perfect, so their irreducible polynomials are separable.
If \(f'=0\), every nonzero term has exponent divisible by \(p\), and there is a unique expression \(f=g^p\).
Over the coefficient field \(\F_{p^r}\), the \(p\)th root of a coefficient \(a\) is \(a^{p^{r-1}}\).

The following gives the \textbf{square-free factorization} of a monic nonconstant polynomial, namely \(f=\prod_i z_i^i\), where each \(z_i\) is monic and the \(z_i\) different from \(1\) are square-free and pairwise coprime.
\begin{enumerate}
\item Put \(c=\gcd(f,f')\), \(w=f/c\), and \(i=1\).
\item While \(w\ne1\), compute \(y=\gcd(w,c)\) and \(z=w/y\); if \(z\ne1\), record the factor \(z\) with multiplicity \(i\). Then put \(w=y\), \(c=c/y\), and \(i=i+1\).
\item After the loop, if the remaining \(c\ne1\), it is a \(p\)th power. Apply the same procedure recursively to \(c^{1/p}\) and multiply every resulting multiplicity by \(p\).
\end{enumerate}
Correctness can be checked separately for each irreducible factor \(P\) and its multiplicity \(e\):
if \(p\nmid e\), it occurs with multiplicity \(e-1\) in \(c\) and once in \(w\); the loop extracts it exactly when \(i=e\).
If \(p\mid e\), it initially remains entirely in \(c\) and is handled by the final \(p\)th-root step.
Every recursive call lowers the degree, so the algorithm terminates.

\subsection{Grouping by Orbit Length}\label{sec:distinct-degree-factorization}
Let \(f\) be monic, nonconstant, and square-free. \textbf{Distinct-degree factorization} uses the fact that the roots of \(X^{q^i}-X\) are exactly the elements whose periods divide \(i\).

Put \(r=f\), \(u=X\bmod r\), and \(i=1\). While \(2i\leq\deg r\), compute
\[
u\leftarrow u^q\bmod r,\qquad g=\gcd(r,u-X).
\]
If \(g\ne1\), record \(g\) as the product of all irreducible factors of degree \(i\), then put \(r\leftarrow r/g\) and reduce \(u\) modulo the new \(r\).
If \(r=1\), stop; otherwise put \(i\leftarrow i+1\) and continue.
On exit, if \(r\ne1\), it is itself an irreducible factor.

It suffices to prove the invariant: on entering step \(i\), the polynomial \(r\) has no irreducible factor of degree less than \(i\), and \(u\) records the previous Frobenius-power remainder.
Every common irreducible factor at this step has degree both dividing \(i\) and at least \(i\), hence exactly \(i\).
If two irreducible factors remained on exit, each would have degree at least the current \(i\), forcing \(\deg r\geq2i\), a contradiction.

\section{Berlekamp's Method: Separating Orbits with the Fixed Subalgebra}
Berlekamp's method~\cite{Berlekamp} computes fixed elements in the quotient algebra by linear equations and uses them to separate factors.

\begin{theorem}\leavevmode\phantomsection\label{thm:berlekamp}
Let \(f=f_1\cdots f_s\) be the irreducible factorization of a square-free polynomial, and put \(A=k[X]/(f)\). Then
\[
B=\{h\in A:h^q=h\}
\]
is a \(k\)-space of dimension \(s\) and is isomorphic as a \(k\)-algebra to \(k^s\).
\end{theorem}
\begin{proof}
In a finite splitting extension, let \(S\) be the root set of \(f\). It has \(s\) Frobenius orbits.
Theorem~\ref{thm:equivariant-functions} identifies \(A\) with the equivariant function algebra
and \(B\) with functions taking a constant base-field value on each orbit. Thus \(B\simeq k^s\), of dimension \(s\).
The \textbf{indicator function} of a subset is one on the subset and zero on its complement.
The orbit indicators give pairwise orthogonal idempotents in \(B\) summing to one; these are precisely the algebra elements used to separate factors.
\end{proof}

\begin{example}[Separating two root orbits with one fixed element]\label{eg:berlekamp-worked}
In \(\F_2[X]\), put
\[
f=(X+1)(X^2+X+1)=X^3+1,\qquad h=X+X^2.
\]
Modulo \(f\), we have \(X^3=1\) and \(X^4=X\), hence \(h^2=X^2+X^4=h\).
Thus \(h\) belongs to the Berlekamp subalgebra.
On the component corresponding to \(X+1\), we have \(h(1)=0\);
on the other component, \(X^2+X+1=0\) gives \(h=1\).
Therefore
\[
\gcd(f,h)=X+1,\qquad \gcd(f,h-1)=X^2+X+1.
\]
A fixed element takes a constant base-field value on each root orbit; distinct values separate the orbits into different factors.
\end{example}

In the basis \(1,\overline X,\ldots,\overline X^{\deg f-1}\), compute the matrix of \(h\mapsto h^q\)
and solve \(h^q-h=0\) to obtain a basis of \(B\).
If \(\dim B=1\), then \(f\) is irreducible. If the dimension exceeds one, choose \(h\in B\) outside the constant subspace.
For every \(c\in k\), compute \(\gcd(f,h-c)\) to group the irreducible factors of \(f\) by the values on their components.
Since these values are not all equal, at least two nontrivial groups result; continuing on the resulting factors gives the full factorization.

The fixed subalgebra is isomorphic to a finite product of copies of the base field and is generally not a field.
It records exactly how many Frobenius orbits make up the root set.
For square-free \(f\), computing the fixed-space dimensions of all Frobenius iterates and applying Theorem~\ref{fra:thm-recovery} also recovers the number of irreducible factors of each degree.

\subsection{Value Polynomials: Avoiding Enumeration of the Base Field}\label{sec:berlekamp-value-polynomial}
When \(q\) is much larger than \(\deg f\), computing every \(\gcd(f,h-c)\) gives many results equal to one.
Linear algebra can first determine the values actually attained by a fixed element; the resultant also records the number of roots corresponding to each value. See also \cite{Feng} for this value-polynomial method.
\phantomsection\label{def:resultant-root-product}
When the first polynomial is monic, the \textbf{resultant} is the product of the values of the second polynomial at all roots of the first: if
\(F(X)=\prod_{i=1}^m(X-\alpha_i)\), then
\[
\operatorname{Res}_X(F,H)=\prod_{i=1}^m H(\alpha_i).
\]
Roots are counted with multiplicity, and the coefficients of \(H\) may contain independent parameters. The section on \hyperref[sec:composed-resultants]{resultants and composed constructions} further explains its relation to composed root products.

The root-product definition of the resultant also gives an exact expression
for the Moore determinant in Chapter~\ref{chap:normal}. In the notation of
Proposition~\ref{v88:prop-moore-circulant},
\begin{equation}\label{v88:eq-moore-resultant}
\boxed{\quad
\det(\alpha^{q^{i+j}})_{0\le i,j<n}
=(-1)^{\lfloor(n-1)/2\rfloor}
\operatorname{Res}_T(T^n-1,h_\alpha(T)).\quad}
\end{equation}
Indeed, \(C_\alpha\) represents multiplication by \(h_\alpha\) in the quotient,
so its determinant is this root product; the column permutation from Chapter~\ref{chap:normal}
supplies the sign. Roots are counted with multiplicity even when \(p\mid n\).

\begin{proposition}[Values of fixed elements and multiplication operators]\label{prop:berlekamp-value-polynomial}
Let \(f\in k[X]\) be monic, nonconstant, and square-free, let \(A=k[X]/(f)\), and let \(h\in A\) satisfy \(h^q=h\).
In the product decomposition given by \(f=\prod_{i=1}^s f_i\), write the coordinates of \(h\) as \(c_i\in k\), and put \(C=\{c_1,\ldots,c_s\}\).
The minimal polynomial of \(h\) over \(k\) is
\[
m_h(Y)=\prod_{c\in C}(Y-c),\qquad m_h\mid Y^q-Y.
\]
In this section, \(m_h\) is called the \textbf{value polynomial} of \(h\).
Choose any polynomial representative of \(h\) and put \(f_c=\gcd(f,h-c)\), taking the greatest common divisor monic. Then
\begin{equation}\label{eq:berlekamp-value-factors}
f=\prod_{c\in C}f_c,\qquad
\chi_h(Y):=\det(YI-M_h)
=\operatorname{Res}_X(f(X),Y-h(X))
=\prod_{c\in C}(Y-c)^{\deg f_c},
\end{equation}
where \(M_h:A\to A\) is multiplication by \(h\). Thus \(m_h\) records distinct values, while \(\chi_h\) also records the total degree of each value fiber.
\end{proposition}
\begin{proof}
On field component \(i\), the equation \(h^q=h\) forces \(h=c_i\in k\).
Hence, for \(a(Y)\in k[Y]\), we have \(a(h)=0\) if and only if \(a(c_i)=0\) for every \(i\), equivalently \(\prod_{c\in C}(Y-c)\mid a\).
This proves the minimal-polynomial formula.
Also, \(f_i\mid h-c\) if and only if \(c_i=c\), so \(f_c\) is exactly the product of these \(f_i\), proving the first identity.
On \(k[X]/(f_i)\), the operator \(M_h\) is the scalar operator \(c_iI\), with characteristic polynomial \((Y-c_i)^{\deg f_i}\).
Multiplying over the components gives the formula for \(\chi_h\).
Finally, the root-product formula for the resultant gives
\(\operatorname{Res}_X(f,Y-h)=\prod_{f(\alpha)=0}(Y-h(\alpha))\);
the value on each root orbit is its corresponding \(c_i\), so this product is also \(\chi_h(Y)\).
\end{proof}

\begin{proposition}[Idempotent projections onto value fibers]\label{v82:prop-value-idempotents}
Let \(A\) be a finite-dimensional commutative unital \(k=\F_q\)-algebra, and let \(h\in A\) satisfy \(h^q=h\). For \(c\in k\), define
\begin{equation}\label{v82:eq-value-idempotents}
 e_c=1-(h-c)^{q-1}.
\end{equation}
Each \(e_c\) is idempotent, and
\[
 e_ce_d=0\ (c\ne d),\qquad
 \sum_{c\in k}e_c=1,\qquad h=\sum_{c\in k}c e_c.
\]
In particular, for the preceding square-free quotient \(A=k[X]/(f)\), the element \(e_c\) is one on field components where \(h\) has value \(c\), and zero on all other components.
If \(f_c=\gcd(f,h-c)\ne1\), then \(Ae_c\simeq k[X]/(f_c)\), with identity element \(e_c\) in \(Ae_c\).
\end{proposition}
\begin{proof}
In \(k[Y]/(Y^q-Y)\), the polynomial \(1-(Y-c)^{q-1}\) is one at \(Y=c\) and zero at every other base-field element.
The Chinese remainder theorem therefore gives the stated idempotence, orthogonality, sum, and weighted-sum identities.
Since \(h^q=h\), substitution \(Y=h\) induces an algebra homomorphism to \(A\), so the identities hold in \(A\) without any assumption that \(A\) is reduced.
For square-free \(f\), evaluation in the field component associated with each irreducible factor gives the description of \(e_c\) and the product decomposition of \(Ae_c\).
\end{proof}

Thus separating factors by \(\gcd(f,h-c)\) and projecting onto the corresponding components by multiplication by \(e_c\) are the polynomial and algebra descriptions of the same splitting.
A nonzero \(e_c\) need not be primitive: if one value combines several root orbits, further splitting is required.
For example, in Example~\ref{eg:berlekamp-worked}, we have \(q=2\), so \(e_0=1-h\) and \(e_1=h\) separate the two orbits.

Let \(t=\deg m_h\). Compute successively the coordinates of \(1,h,h^2,\ldots\) in the power basis of \(A\).
At the first linear dependence there is a unique expression
\[
h^t+b_{t-1}h^{t-1}+\cdots+b_0=0,
\qquad m_h(Y)=Y^t+b_{t-1}Y^{t-1}+\cdots+b_0.
\]
Indeed, the preceding powers are linearly independent. If this relation had degree greater than \(\deg m_h\), the minimal polynomial would give an earlier relation, a contradiction.
Thus \(m_h\) can be computed without knowing \(f_i\) or \(C\), and \(t\le s\le\deg f\).
It remains to find the base-field roots of this completely split square-free polynomial and compute the greatest common divisors in~\eqref{eq:berlekamp-value-factors} for just these \(t\) values.
Root finding itself need not enumerate the base field.

\begin{proposition}[Randomized root finding for completely split polynomials]\label{prop:split-values-random}
Let \(M\in k[Y]\) be monic, square-free, and completely split, with \(\deg M=t\ge2\).
Choose a polynomial \(b\) of degree less than \(t\) uniformly at random. Each attempt of the following procedure finds a nontrivial factor of \(M\) with probability at least \(1/2\).
\begin{enumerate}
\item If \(q\) is odd, first compute \(g_0=\gcd(M,b)\). If \(0<\deg g_0<t\), return \(g_0\); if \(g_0=M\), choose \(b\) again; if \(g_0=1\), compute
\[
g_1=\gcd\bigl(M,b^{(q-1)/2}-1\bigr).
\]
If \(0<\deg g_1<t\), return \(g_1\); otherwise choose \(b\) again.
\item If \(q=2^r\), compute
\[
B=b+b^2+\cdots+b^{2^{r-1}}\pmod M,
\qquad g_0=\gcd(M,B).
\]
If \(0<\deg g_0<t\), return \(g_0\); otherwise choose \(b\) again.
\end{enumerate}
Apply the procedure recursively to the factors until all have degree one, obtaining every root of \(M\). Randomness affects only the number of attempts, not correctness of the output.
\end{proposition}
\begin{proof}
The Chinese remainder theorem gives \(k[Y]/(M)\simeq k^t\), so the values of \(b\) at the roots are independent and uniform in \(k\).
In odd characteristic an attempt fails exactly when the values are all zero, all nonzero squares, or all nonsquares.
Thus the failure probability is
\[
q^{-t}+2\left(\frac{q-1}{2q}\right)^t
\le\frac1{q^2}+\frac{(q-1)^2}{2q^2}\le\frac12.
\]
In even characteristic the absolute trace is a surjection \(k\to\F_2\), so the components of \(B\) are independent, with zero and one equally likely.
Splitting fails only when all components agree, with probability \(2^{1-t}\le1/2\).
A nontrivial greatest common divisor and its quotient remain square-free and completely split, so recursion applies.
The binary splitting tree has \(t-1\) successful splits; the expected total number of attempts is therefore at most \(2(t-1)\).
\end{proof}

This root-finding procedure uses only polynomial operations of degree at most \(t\), random base-field elements, and modular exponentiation;
the exponents have bit length \(O(\log q)\), and no enumeration of the \(q\) base-field elements is needed.
If \(h\notin k\), then \(t\ge2\), yielding at least two nontrivial factor groups.
If a group is still reducible, recompute its fixed space and continue splitting.
The original polynomial must first be decomposed into square-free factors; neither a fixed-space dimension nor a single derivative computation replaces this step.
In particular, some root multiplicities of \(\chi_h\) may be divisible by the characteristic, so taking \(\chi_h/\gcd(\chi_h,\chi_h')\) directly generally does not give \(m_h\); see Exercise~\ref{ex:berlekamp-value-multiplicity}.

\begin{example}[Equal values merge distinct root orbits]\label{exm:berlekamp-values-101}
Over \(k=\F_{101}\), take
\[
f=(X-1)(X-2)(X-3),\qquad h=(X-1)(X-2)\pmod f.
\]
Since \(f\) splits completely over \(k\), we have \(h^{101}=h\) in \(A=k[X]/(f)\). Also,
\[
h^2-2h=Xf\equiv0\pmod f,
\]
and \(h\) is nonconstant, so the first linear relation gives \(m_h(Y)=Y(Y-2)\).
Only two greatest common divisors are needed:
\[
f_0=\gcd(f,h)=(X-1)(X-2),\qquad
f_2=\gcd(f,h-2)=X-3.
\]
Thus \(\chi_h(Y)=Y^2(Y-2)\). The value zero combines two singleton root orbits;
the degree of \(m_h\), which is two, counts distinct values rather than irreducible factors.
For \(f_0\), choose \(h_0=X\) instead. Its value polynomial is \((Y-1)(Y-2)\), and another split completes the factorization.
In the same coordinates the idempotents are \(e_2=h/2\) and \(e_0=1-h/2\). The first projects onto the root three; the second still combines the roots one and two, so it is not primitive.
\end{example}

\paragraph{Constructing a multiple from the affine span of the roots.}
\phantomsection\label{def:affine-space}
One may also approach root finding by first embedding the root set into an \textbf{affine subspace}, a translate of a linear subspace \(V\),
\(\beta+V=\{\beta+v:v\in V\}\). We also call this an affine space; its dimension and affine maps are defined in Section~\ref{sec:affine-root-space}.
A constant plus a \(q\)-linearized polynomial is called an \textbf{affine \(q\)-polynomial}.
The space spanned by root differences gives the following degree bound.

\begin{proposition}[An affine \(q\)-polynomial multiple]\label{prop:affine-q-multiple}
Let \(f\in k[X]\) be monic and square-free with \(\deg f=n\ge1\), and write its root set in a finite splitting extension as \(S\).
Choose \(\beta\in S\) and put
\[
V=\operatorname{span}_k\{\alpha-\beta:\alpha\in S\},\qquad v=\dim_kV.
\]
Then \(v\le n-1\), both \(V\) and \(\beta+V\) are Frobenius-stable, and
\[
A_f(X)=L_V(X)-L_V(\beta)\in k[X],\qquad
f\mid A_f,\qquad \deg A_f=q^v\le q^{n-1}.
\]
Here \(\beta+V\) is the smallest affine \(k\)-subspace containing \(S\). Thus \(A_f\) is independent of the choice of \(\beta\).
\end{proposition}
\begin{proof}
Subtracting the chosen root from every element of \(S\) gives \(n\) distinct differences, including \(\beta-\beta=0\), so \(v\le n-1\).
Write \(\sigma(x)=x^q\). Since \(\sigma\) permutes \(S\),
\[
\sigma(\alpha-\beta)=(\sigma(\alpha)-\beta)-(\sigma(\beta)-\beta)\in V.
\]
Thus \(\sigma(V)\subseteq V\); finiteness and injectivity of \(\sigma\) give equality.
Also \(\sigma(\beta)-\beta\in V\), so \(\sigma(\beta+V)=\beta+V\).
By Theorem~\ref{thm:subspace-polynomial}, the polynomial \(L_V\) is \(q\)-linearized, and Frobenius stability of its root set gives \(L_V\in k[X]\).
Thus
\[
L_V(\beta)^q=L_V(\sigma(\beta))
=L_V(\beta)+L_V(\sigma(\beta)-\beta)=L_V(\beta),
\]
so \(L_V(\beta)\in k\). The root set of \(A_f\) is exactly \(\beta+V\), containing \(S\); square-freeness of \(f\) gives \(f\mid A_f\).
Every affine space containing \(S\) contains \(\beta\), and its direction space contains all root differences and their \(k\)-linear combinations. It therefore contains \(\beta+V\).
This proves minimality and independence of the base point; the degree formula follows from \(|V|=q^v\).
\end{proof}

The polynomial \(A_f\) can also be computed without knowing any roots.
In \(k[X]/(f)\), compute successively the coordinates of
\[
1,\ X,\ X^q,\ \ldots,\ X^{q^{n-1}}\pmod f.
\]
These \(n+1\) vectors are necessarily dependent. Take the first relation and normalize its leading coefficient to one to obtain a nonzero affine \(q\)-polynomial multiple.
It is exactly \(A_f\): an affine \(q\)-polynomial vanishing on \(S\) has, by additivity and \(k\)-homogeneity of its linearized part, a root set containing \(\beta+V\), so its degree is at least \(q^v\);
but \(A_f\) already provides a relation of that degree. The difference of two monic relations of least degree has lower degree and hence must be zero.
Using root differences lowers the upper bound for the span dimension from \(n\) to \(n-1\).
After constructing the multiple, one must still find its roots and test the original polynomial; this degree bound alone does not guarantee efficient root finding for every input.

For example, the two roots of \(f=X^2+X+1\in\F_2[X]\) are \(\alpha,\alpha+1\), and the root-difference space is \(V=\F_2\).
Thus \(L_V=X^2+X\) and \(L_V(\alpha)=1\), and the affine multiple is \(f\) itself.
The roots themselves span the linear space \(\F_4\), with subspace polynomial \(X^4-X\).
Here retaining the translation constant lowers the multiple's degree from four to two.

\subsection*{Repeated Factors Do Not Increase the Number of Fixed Components}
\begin{proposition}\label{prop:berlekamp-repeated}
Let \(f=\prod_{i=1}^sP_i^{e_i}\in k[X]\), where \(s\ge1\), the \(P_i\) are pairwise distinct monic irreducible polynomials, and the \(e_i\) are positive integers. Put \(B_f=\{\bar h\in k[X]/(f):\bar h^q=\bar h\}\). Then \(B_f\cong k^s\), so \(\dim_kB_f=s\).
\end{proposition}
\begin{proof}
The reduced quotient of this quotient algebra is \(\prod_{i=1}^s k[X]/(P_i)\).
Apply Corollary~\ref{cor:algebra-fixed-points-reduced} with \(m=1\). Each field component has fixed part \(k\), so \(B_f\simeq k^s\).
The natural isomorphism follows from Proposition~\ref{prop:algebra-periodic-reduction}: realize the reduced quotient by Frobenius-periodic points, then take the Frobenius-fixed points.
\end{proof}
Thus \(\dim B_f=1\) only says that there is one distinct irreducible factor.
It implies irreducibility of \(f\) only when \(f\) is already square-free.

\subsection*{A Berlekamp Basis on Roots of Unity}
Let \(p\nmid N\) and \(f=X^N-1\). In the quotient ring, Frobenius takes \(X^i\) to \(X^{qi\bmod N}\).
A basis of \(B_f\) is therefore
\[
\sum_{i\in C}X^i,
\]
where \(C\) runs through all \(q\)-cyclotomic cosets modulo \(N\).
Compare coefficients: a polynomial is fixed exactly when its coefficients are constant on each exponent orbit.
When \(N=7,q=2\), this gives \(1,\ X+X^2+X^4,\ X^3+X^5+X^6\).
Conjugate root orbits control factorization, while monomial exponent orbits control the functions separating the factors; the two orbit counts coincide.

\section{Equal-Degree Factorization and Random Splitting}\label{sec:equal-degree-factorization}
Let \(f\in k[X]\) be monic and square-free, a product of \(s\ge2\) distinct irreducibles of degree \(d\), and write \(Q=q^d\). The Chinese remainder theorem gives
\[
A=k[X]/(f)\cong\prod_{i=1}^s k[X]/(f_i)\cong(\F_Q)^s.
\]
This isomorphism explains correctness; during computation the factors \(f_i\) are unknown, and the field components on the right need not be constructed.
Choosing \(sd\) independent uniform base-field coefficients to form \(h\) of degree less than \(sd\) is equivalent to sampling uniformly from \(A\).
Since the isomorphism is a bijection of finite sets, the component values are independent and uniform in \(\F_Q\).

Equal degrees allow every component to use the same exponent or trace map.
These operations map values in unknown field components into a small set, after which greatest common divisors separate factors with different values.
The resulting elements satisfy \(b^q=b\) and thus lie in the Berlekamp subalgebra of the preceding section; sampling and modular operations produce the separating elements here.

\subsection*{Splitting Steps and Success Probabilities}
An attempt follows these rules, with all powers and sums computed in the quotient ring modulo \(f\).
\begin{enumerate}
\item If \(q\) is odd, first put \(g_0=\gcd(f,h)\). If \(0<\deg g_0<sd\), return \(g_0\) and \(f/g_0\); if \(g_0=f\), sample again.
If \(g_0=1\), compute
\[
b=h^{(Q-1)/2}\pmod f,\qquad g_+=\gcd(f,b-1).
\]
If \(0<\deg g_+<sd\), return \(g_+\) and \(f/g_+\); otherwise sample again.
\item If \(q=2^r\), compute the absolute trace to \(\F_2\),
\begin{equation}\label{v48:eq-edf-absolute-trace}
b=h+h^2+h^{2^2}+\cdots+h^{2^{rd-1}}\pmod f,
\qquad g_0=\gcd(f,b).
\end{equation}
If \(0<\deg g_0<sd\), return \(g_0\) and \(f/g_0\); otherwise sample again.
\end{enumerate}
All greatest common divisors are monic. In odd characteristic, when \(h\) is a unit, each component of \(b\) equals \(1\) or \(-1\), corresponding to nonzero squares and nonsquares; in even characteristic, each component's absolute trace is \(0\) or \(1\).
Every split therefore partitions whole irreducible factors into two groups, and each output degree is a positive multiple of \(d\).

\begin{proposition}[Success probability of equal-degree splitting]\label{v48:prop-edf-success}
The success probability of one attempt above is
\[
P_{\mathrm{odd}}=1-Q^{-s}-2\left(\frac{Q-1}{2Q}\right)^s
\quad(q\text{ odd}),\qquad
P_{\mathrm{even}}=1-2^{1-s}
\quad(q\text{ even}).
\]
Both are at least \(1/2\). In odd characteristic, conditional on \(\gcd(f,h)=1\), the subsequent sign split succeeds with probability exactly \(1-2^{1-s}\).
\end{proposition}
\begin{proof}
The componentwise proof of Proposition~\ref{prop:split-values-random} applies to \(s\) components isomorphic to \(\F_Q\).
In odd characteristic, failure occurs only when all values are zero, all nonzero squares, or all nonsquares, giving the formula.
Since \(s\ge2\) and \(Q\ge3\), the failure probability is at most
\[
\frac1{Q^2}+\frac{(Q-1)^2}{2Q^2}\le\frac12.
\]
The condition \(\gcd(f,h)=1\) means every component is nonzero. Conditional on this, the components remain independent and uniform in \(\F_Q^*\), with the two signs equally likely; failure occurs only if all signs agree.
In even characteristic, Theorem~\ref{thm:trace-norm-surjective} makes the absolute traces independent with zero and one equally likely; failure occurs only if they all agree.
\end{proof}

In odd characteristic, the unconditional success probability also includes the opportunity to separate zero components using \(\gcd(f,h)\); the conditional sign-splitting probability counts only the subsequent stage.
In even characteristic the two trace values are equally likely in each component. Using only the relative trace \(\F_Q\to\F_q\) and then the greatest common divisor for value zero would give success probability
\[
1-q^{-s}-\left(\frac{q-1}{q}\right)^s,
\]
which tends to zero as \(q\) grows with \(s\) fixed. Transitivity also computes the absolute trace: first compute
\[
u=h+h^q+\cdots+h^{q^{d-1}}\pmod f,
\]
and then \(b=u+u^2+\cdots+u^{2^{r-1}}\pmod f\), giving exactly~\eqref{v48:eq-edf-absolute-trace}.

\subsection*{Recursive Factorization and Correctness of the Output}
Repeat the procedure on each output factor \(F\), sampling independently and uniformly from polynomials of degree less than \(\deg F\).
If \(\deg F=d\), return \(F\): it must be irreducible.
If \(\deg F>d\), the polynomial \(F\) is still a product of distinct irreducibles of degree \(d\) with at least two components, so the same success-probability bound applies.
Failure only triggers resampling; success strictly lowers the degrees of both subproblems.

\leavevmode\phantomsection\label{v48:def-las-vegas}\ignorespaces
A randomized algorithm whose running time depends on random choices, whose output is always correct, and whose expected running time is finite is called a \textbf{Las Vegas algorithm}.
This procedure is of that type: each trial uses finitely many polynomial operations, and each recursive node needs at most two trials in expectation. A full factorization has \(s-1\) successful splitting nodes, so the expected total number of trials is at most \(2(s-1)\).
Independent retries terminate with probability one. Every split uses a greatest common divisor and exact division, and the final factors have degree exactly \(d\); randomness therefore does not affect output correctness.
This gives Cantor--Zassenhaus equal-degree factorization; see \cite{CZ}.

\begin{example}[Splitting by quadratic character values in odd characteristic]\label{v48:ex-edf-odd}
Over \(\F_3[X]\), take
\[
f=X^4+X^3+X+2,\qquad d=2,\qquad h=X.
\]
Modular computation gives \(X^9\equiv X\pmod f\) and \(\gcd(f,X^3-X)=1\), so \(f\) is square-free with every irreducible factor of degree two.
Since \(f(0)=2\), we have \(\gcd(f,h)=1\).
Here \(Q=9\), so the quadratic character value is given by \(h^4\), and
\[
b=X^4\bmod f=2X^3+2X+1,
\qquad b-1=2X(X^2+1).
\]
Hence
\[
g_+=\gcd(f,b-1)=X^2+1,
\qquad \frac f{g_+}=X^2+X+2.
\]
Neither quadratic has a root in \(\F_3\), so the factorization is complete.
In the two residue fields, \(X^4=1\) modulo \(X^2+1\) and \(X^4=2=-1\) modulo \(X^2+X+2\), so the two components have different quadratic character values.
Here \(s=2\), and the unconditional success probability of one trial is
\[
1-\frac1{81}-2\left(\frac49\right)^2=\frac{16}{27},
\]
while, conditional on \(h\) being a unit, sign splitting succeeds with probability \(1/2\).
\end{example}

\begin{example}[Splitting by the absolute trace in even characteristic]\label{v48:ex-edf-even}
Over \(\F_2[X]\), take
\[
f=X^6+X^5+X^4+X^3+X^2+X+1,\qquad d=3,\qquad h=X.
\]
Since \((X+1)f=X^7+1\), we have \(X^8\equiv X\pmod f\); also \(\gcd(f,X^2-X)=1\), so \(f\) is square-free with every irreducible factor of degree three.
Each component lies in \(\F_8\), with absolute trace
\[
b=X+X^2+X^4=X(X^3+X+1).
\]
Computing the greatest common divisor and the corresponding quotient gives
\[
\gcd(f,b)=X^3+X+1,
\qquad \frac f{\gcd(f,b)}=X^3+X^2+1.
\]
Neither cubic has a root in \(\F_2\), so both are irreducible.
On the first component, \(X^3=X+1\) gives \(X^4=X^2+X\), hence \(b=0\);
on the second, \(X^3=X^2+1\) gives \(X^4=X^2+X+1\), hence \(b=1\).
Distinct trace values separate the two root orbits, by the same mechanism as the fixed-element separation in Example~\ref{eg:berlekamp-worked}.
If \(h=1\) is chosen, both component traces are one and no split is obtained; a new sample is then selected.
This example has two field components, so each uniformly sampled attempt succeeds with probability \(1/2\).
\end{example}

\section{From Quotient-Ring Representations to Explicit Isomorphisms}
\label{sec:explicit-field-isomorphism}

The structure theorem determines the isomorphism type of a finite field; computation also requires the image of a specified generator.
Root finding directly connects these two levels of information.

\begin{proposition}\leavevmode\phantomsection\label{prop:construct-field-isomorphism}
Let \(f,g\in k[X]\) be monic irreducible polynomials of degree \(n\), and put
\(E=k[X]/(f)=k(\alpha)\) and \(E'=k[Y]/(g)=k(\beta)\).
Choose a root \(b\) of \(f\) in \(E'\). Then
\[
\Phi:E\longrightarrow E',\qquad P(\alpha)\longmapsto P(b)
\]
is a \(k\)-isomorphism. All \(k\)-isomorphisms are exactly the \(n\) maps taking \(\alpha\) to \(b,b^q,\ldots,b^{q^{n-1}}\).
\end{proposition}
\begin{proof}
Since \(f\mid X^{q^n}-X\), it splits into \(n\) distinct linear factors in \(E'\), so a root \(b\) exists.
Evaluation \(P\mapsto P(b)\) has kernel \((f)\) and therefore induces an injection from \(E\) to \(E'\).
Both sides have dimension \(n\), so it is an isomorphism.
Every \(k\)-isomorphism is determined by the image of \(\alpha\), which must be a root of \(f\); its roots are exactly the listed Frobenius orbit.
\end{proof}

Algorithmically, apply this chapter's factorization methods to \(f\) over \(E'[X]\) and take a linear factor.
Enumerating the elements of \(E'\) also terminates, but its \(q^n\) candidates grow exponentially with the representation length \(n\log q\).
Once a root is chosen, computing the isomorphism requires only evaluation.
Precompute the coordinates of \(1,b,\ldots,b^{n-1}\) in the power basis of \(\beta\) and use them as the columns of a matrix \(P\).
Then
\[
[\Phi(x)]_{(1,\beta,\ldots,\beta^{n-1})}
=P[x]_{(1,\alpha,\ldots,\alpha^{n-1})}.
\]
The matrix \(P\) is invertible, and computing \(P^{-1}\) gives the coordinate formula for the inverse isomorphism.

\subsection*{Defining Equations and Changes of Generator}
\phantomsection\label{v79:sec-defining-equations}
Fix the extension \(E/k\), where \(E=\F_{q^n}\), \(k=\F_q\), and \(n\ge2\). A \textbf{defining polynomial} for this extension is the monic minimal polynomial of a generator of \(E/k\).
Every member of \(\Irr_q(n)\) gives such a representation, but the coefficients can have very different support patterns.
Thus the existence of a generator leads to a further question: can it be chosen so that its defining equation has a simpler form?

\leavevmode\phantomsection\label{v79:def-defining-complexity}\ignorespaces
For a polynomial \(f(X)=\sum_i a_iX^i\), its \textbf{support size} is \(\operatorname{wt}(f)=\#\{i:a_i\ne0\}\).
Two measures of the simplest defining equations of degree \(n\) are
\begin{equation}\label{v79:eq-defining-complexity}
\begin{aligned}
w_q(n)&=\min_{f\in\Irr_q(n)}\operatorname{wt}(f),\\
\theta_q(n)&=\min\{\deg g:0\ne g\in k[X],\ \deg g<n,\ X^n+g\in\Irr_q(n)\}.
\end{aligned}
\end{equation}
The second quantity minimizes the highest degree of a nonleading term. These minima exist because \(\Irr_q(n)\) is finite and nonempty.
An irreducible polynomial of degree \(n\ge2\) has a nonzero constant term, so
\[
2\le w_q(n)\le\theta_q(n)+2\le n+1.
\]
For the middle inequality, choose \(g\) of minimum degree and count its at most \(\deg g+1\) nonzero coefficients.
In particular,
\begin{equation}\label{v79:eq-binomial-optimum}
\begin{aligned}
w_q(n)=2\quad&\Longleftrightarrow\quad\theta_q(n)=0\\
&\Longleftrightarrow\quad X^n+a\in\Irr_q(n)\text{ for some }a\in k^\times.
\end{aligned}
\end{equation}
Indeed, a monic irreducible polynomial with two nonzero terms must consist of its leading term and its nonzero constant term.
Theorem~\ref{thm:binomial-full} and Exercise~\ref{ex:supp-binomial} determine exactly when this happens and which constants are possible.

The two objectives need not have the same minimizers. Over \(\F_2\), both \(X^3+X+1\) and \(X^3+X^2+1\) are irreducible, since neither has a root in \(\F_2\).
No irreducible cubic binomial exists: its constant term would have to be \(1\), and \(X^3+1\) vanishes at \(1\).
Thus both cubics attain \(w_2(3)=3\), but only the first attains \(\theta_2(3)=1\).

Here the defining equation is optimized over all field generators. For a fixed generator, sparse relations among its powers also involve the allowed exponent range; Subsection~\ref{v97:sec-sparse-relations} treats that problem separately.

\noindent\textbf{Recovering minimal polynomials from orbits.}
Given an element \(a\in E\) in the chosen representation, compute \(a,a^q,a^{q^2},\ldots\) successively until the first return to \(a\).
If the orbit length is \(d\), then
\[
m_a(X)=\prod_{i=0}^{d-1}(X-a^{q^i})\in k[X].
\]
Frobenius stability places the coefficients in \(k\), and the orbit theorem gives minimality of the degree.
After enlarging the base field, change the step size in the same root list and take the product over each new orbit.
Minimal polynomials, factorization over extensions, and relative traces and norms can therefore share the same conjugate data.

Taking the image of a generator also requires retaining the multiplicities of its conjugate values.

\begin{proposition}[Rational images of a generating orbit]\label{v79:prop-generator-orbit-image}
Let \(k=\F_q\), \(E=k(\alpha)\), \([E:k]=n\ge2\), and \(f=m_\alpha\), and put
\(\beta=R(\alpha)\), where \(R=A/B\in k(X)\) has relatively prime numerator and denominator \(A,B\in k[X]\), and \(B(\alpha)\ne0\).
Write \(d=[k(\beta):k]\). Then \(R\) maps the Frobenius orbit of \(\alpha\) onto that of \(\beta\), each image value occurs \(n/d\) times, and
\begin{equation}\label{v70:eq-image-conjugate-product}
\prod_{i=0}^{n-1}(Y-\beta^{q^i})=m_\beta(Y)^{n/d}.
\end{equation}
Equivalently,
\begin{equation}\label{v70:eq-image-resultant}
\frac{\operatorname{Res}_X\bigl(f(X),B(X)Y-A(X)\bigr)}
     {\operatorname{Res}_X\bigl(f(X),B(X)\bigr)}
=m_\beta(Y)^{n/d},
\end{equation}
where the denominator is \(\Nm_{E/k}(B(\alpha))\in k^\times\).
The conjugate product is a new defining polynomial for \(E/k\) if and only if \(d=n\), equivalently, if and only if \(R\) is injective on the original root orbit.
\end{proposition}
\begin{proof}
Because the coefficients of \(A,B\) lie in \(k\), for every \(i\) we have
\[
B(\alpha^{q^i})=B(\alpha)^{q^i}\ne0,
\qquad R(\alpha^{q^i})=R(\alpha)^{q^i}=\beta^{q^i}.
\]
Thus the orbit map is defined on every conjugate of \(\alpha\) and commutes with Frobenius.
Its image is the full orbit of \(\beta\), whose length is \(d\mid n\). Each of its values therefore occurs exactly \(n/d\) times, giving \eqref{v70:eq-image-conjugate-product}.
The root-product formula for the resultant gives
\[
\operatorname{Res}_X\bigl(f(X),B(X)Y-A(X)\bigr)
=\prod_{i=0}^{n-1}B(\alpha^{q^i})\,
 \prod_{i=0}^{n-1}(Y-\beta^{q^i}).
\]
The first product is both \(\operatorname{Res}_X(f,B)\) and \(\Nm_{E/k}(B(\alpha))\), and is nonzero.
Division proves \eqref{v70:eq-image-resultant}. Finally, the image has \(n\) distinct elements precisely when \(d=n\); in that case its root polynomial is the minimal polynomial of a generator of \(E/k\).
\end{proof}

For a nonconstant rational function \(R=A/B\) in lowest terms, its degree is \(\deg R=\max\{\deg A,\deg B\}\).
Lemma~\ref{v70:lem-low-degree-generation} gives a low-degree condition that preserves generation;
Proposition~\ref{v70:prop-rational-reconstruction} addresses the converse task of representing a specified element by a rational function with bounded numerator and denominator degrees.
The minimum support size concerns all defining polynomials of the extension. Starting from a prescribed generator adds a different constraint: the degree of the rational function needed to reach a chosen defining equation.
If an irreducible binomial of degree \(n\) exists, one of its roots lies in \(E=k[\alpha]\) and is a polynomial in \(\alpha\) of degree less than \(n\); this proves reachability without imposing a smaller degree bound.
Exercises~\ref{v79:ex-rational-generation} and~\ref{v79:ex-binomial-defining-equation} compare these degree restrictions and exhibit a change to a binomial defining equation.

\subsection{Sparse Relations for a Fixed Defining Polynomial}\label{v97:sec-sparse-relations}
Changing a generator changes its minimal polynomial. For a fixed generator, one may instead allow a relation to have higher degree in order to reduce its number of nonzero terms.
These are different problems. We now fix a defining polynomial and express the second problem as linear relations among powers in a quotient algebra.

\begin{definition}[Minimum degree of a sparse multiple]\label{v97:def-sparse-degree}
Let \(k=\F_q\), and let \(f\in k[X]\) be monic and nonconstant, with \(f(0)\ne0\). For an integer \(s\ge2\), define
\[
 L_s(f)=\min\{\deg P:0\ne P\in k[X],\ f\mid P,\ \operatorname{wt}(P)\le s\}.
\]
We call \(L_s(f)\) the \textbf{minimum degree of a sparse multiple} of \(f\) with at most \(s\) terms.
The degree is that of the multiple \(P\), not of the multiplier \(P/f\).
\end{definition}

In the finite quotient algebra \(k[X]/(f)\), the element \(\overline X\) is a unit. Hence \(f\) divides some \(X^B-1\), and all the minima above are finite.
Every nonzero multiple has degree at least \(\deg f\), and increasing the permitted number of terms enlarges the admissible set. Therefore,
\[
 \deg f\le L_{s+1}(f)\le L_s(f).
\]
A multiple attaining the minimum degree has nonzero constant term: otherwise one could divide it by \(X\), retain divisibility by \(f\), and lower its degree.
Thus a degree-minimal binomial, after monic normalization, has the form \(X^B-c\), with \(c\in k^\times\).
By Proposition~\ref{v40:prop-minimal-binomial}, \(L_2(f)\) is precisely the order of \(\overline X\) in the unit group modulo base-field scalars.
When \(f\) is irreducible and \(\alpha\) is a root, divisibility \(f\mid\sum_i a_iX^i\) is equivalent to \(\sum_i a_i\alpha^i=0\).
Thus \(L_s(f)\) controls both the number of terms and the exponent range of a relation, whereas \(w_q(n)\) allows the generator to vary. The two quantities are not interchangeable.
The quantity \(\lambda_q(H)\) in Section~\ref{v95:sec-subgroup-periods} instead fixes the entire multiplicative subgroup and minimizes the number of terms in a relation, without an additional restriction on their exponent positions.

\begin{proposition}[Coefficient cancellation by a linear multiplier]\label{v97:prop-linear-cancellation}
Let \(d=\deg f\ge2\), and write \(f=\sum_{i=0}^d a_iX^i\). If some coefficient is zero, then \(L_d(f)=d\).
If all coefficients are nonzero, define the \textbf{adjacent coefficient ratios} by \(r_i=a_i/a_{i-1}\) for \(1\le i\le d\). For \(r\in k^\times\),
\begin{equation}\label{v97:eq-linear-cancellation}
 \operatorname{wt}\bigl((1-rX)f\bigr)
 =d+2-\#\{i:1\le i\le d,\ r_i=r\}.
\end{equation}
Consequently, when all coefficients are nonzero,
\[
 \boxed{L_d(f)=d+1\quad\Longleftrightarrow\quad r_1,\ldots,r_d\text{ contain a repeated value}.}
\]
\end{proposition}
\begin{proof}
A nonzero multiple of degree \(d\) is a scalar multiple of \(f\). This proves the first case and shows that \(L_d(f)>d\) when every coefficient is nonzero.
The constant and leading coefficients of \((1-rX)f\) are \(a_0\) and \(-ra_d\), both nonzero.
Its coefficient of \(X^i\) is \(a_i-ra_{i-1}\), which vanishes precisely when \(r_i=r\). This proves \eqref{v97:eq-linear-cancellation}.
A repeated ratio therefore cancels at least two coefficients, producing a multiple of degree \(d+1\) with at most \(d\) terms.
Conversely, a multiple of degree \(d+1\) has a linear multiplier.
If the multiplier has zero constant term, the product is a scalar multiple of \(Xf\) and still has \(d+1\) terms.
Otherwise the multiplier can be normalized to \(1-rX\), and \eqref{v97:eq-linear-cancellation} forces a repeated ratio.
\end{proof}

The ratio \(r_i=r\) corresponds to the multiplier \(1-rX\), whose monic form is \(X-r^{-1}\), not \(X-r\).
The direction of the reciprocal is fixed by the coefficient \(a_i-ra_{i-1}\).

\begin{corollary}[Cancellation over a smaller coefficient field]\label{v97:cor-small-coefficient-field}
If the coefficients of \(f\) lie in a subfield \(\F_Q\subseteq k\) and \(d\ge Q\), then \(L_d(f)\le d+1\).
\end{corollary}
\begin{proof}
If a coefficient is zero, then \(L_d(f)=d\). Otherwise the \(d\) adjacent ratios lie in the \(Q-1\) element set \(\F_Q^\times\), so two coincide. Apply the preceding proposition.
\end{proof}

\begin{example}[The same multiplicative period, but different three-term relations]\label{v97:ex-f5-sparse}
In \(\F_5[X]\), take
\[
 f_1=X^3+3X^2+X+1,\qquad f_2=X^3+4X^2+X+1.
\]
Their values at \(0,1,2,3,4\) are \((1,1,3,3,2)\) and \((1,2,2,2,3)\), respectively. Hence both are irreducible.
For a root \(\alpha\) of either polynomial, the norm gives \(\alpha^{31}=-1\), so its multiplicative order is \(2\) or \(62\).
Degree three excludes order \(2\). Both roots therefore have multiplicative order \(62\) and projective order \(31\), giving
\[
 L_2(f_1)=L_2(f_2)=31.
\]
On the other hand, direct multiplication gives
\[
 (X^2+2X+3)f_1=X^5+2X^2+3,\qquad
 (X+1)f_2=X^4+2X+1.
\]
The adjacent ratios for \(f_1\) are \(1,3,2\), with no repetition, so it has no three-term multiple of degree four.
Those for \(f_2\) are \(1,4,4\). Since each polynomial itself has four terms,
\[
 \boxed{L_3(f_1)=5,\qquad L_3(f_2)=4.}
\]
The multiplicative and projective periods do not determine where the first three-term relation occurs.
\end{example}

For general \(f\), writing each \(X^i\bmod f\) as a column vector turns sparse relations into linear relations among specified remainder columns.
This also explains why the monomial coordinates must remain fixed: a translation of the variable preserves irreducibility but usually changes the number of terms.
Section~\ref{ffs:sec-coefficient-recovery} will describe the same relations by evaluations that retain root multiplicities.

\subsection{Gauss periods, roots of unity, and optimal normal bases}\label{v90:subsec-gauss-optimal}
The following construction turns cyclotomic irreducibility into sparsity of normal-basis multiplication. Let $E=\F_{q^n}$ and write $\alpha_i=\alpha^{q^i}$ for a normal basis. If
\[
\alpha\alpha_i=\sum_{j=0}^{n-1}t_{ij}\alpha_j\qquad(0\le i<n),
\]
the number of nonzero entries of \((t_{ij})\) is its \emph{multiplication complexity}. We first prove the general lower bound, then give an explicit construction attaining it.
\begin{lemma}[The lower bound for normal-basis multiplication complexity]\label{v104:lem-normal-complexity}
For every \(n\ge1\), the multiplication complexity of a normal basis of \(\F_{q^n}/\F_q\) is at least \(2n-1\).
\end{lemma}
\begin{proof}
In the notation above, multiplication by \(\alpha\ne0\) is an invertible linear map. Thus \((t_{ij})\) is nonsingular and every column is nonzero.
Since \(\alpha_0=\alpha\), we also have
\[
\sum_{i=0}^{n-1}\alpha\alpha_i=\Tr_{E/\F_q}(\alpha)\alpha_0,
\qquad
\sum_{i=0}^{n-1}t_{ij}=
\begin{cases}\Tr_{E/\F_q}(\alpha),&j=0,\\0,&j\ne0.\end{cases}
\]
Linear independence of the normal basis also ensures that \(\Tr(\alpha)=\sum_i\alpha_i\ne0\).
The zeroth column therefore has at least one nonzero entry. Every other column is nonzero but has sum zero, so each has at least two nonzero entries.
Their total number is at least \(1+2(n-1)=2n-1\).
\end{proof}
\begin{theorem}[Optimal normal bases of Gauss-period type I]\label{v90:thm-gauss-optimal}
Suppose $r=n+1$ is prime and $\ord_r(q)=n$. If $\zeta$ is a primitive $r$th root of unity, then
\[
\zeta,\zeta^q,\ldots,\zeta^{q^{n-1}}
\]
is a normal basis of $\F_{q^n}/\F_q$ with multiplication complexity $2n-1$; hence it is optimal.
\end{theorem}
\begin{proof}
The hypothesis implies that $\Phi_r(X)=1+X+\cdots+X^n$ is irreducible over $\F_q$ and that the conjugates listed are a permutation of $\zeta,\zeta^2,\ldots,\zeta^n$. A nontrivial linear relation would yield a nonzero polynomial of degree at most $n$, with zero constant term, vanishing at $\zeta$. This would have to be a scalar multiple of $\Phi_r$, contradicting $\Phi_r(0)=1$. Thus the conjugates form a normal basis.
Order the basis as $\zeta,\zeta^2,\ldots,\zeta^n$. For $1\le j<n$, the product $\zeta\zeta^j=\zeta^{j+1}$ contributes one nonzero coefficient, whereas
\[
\zeta\zeta^n=1=-\sum_{j=1}^n\zeta^j
\]
contributes $n$. Consequently the complexity is $2n-1$, which attains the general lower bound.
\end{proof}
\begin{remark}[Cyclotomic cosets and trace correlation]
The basis comes from the transitive action of $q$ on nonzero exponents modulo $r$; its sparse multiplication comes from $\Phi_r(\zeta)=0$. The cyclic correlation polynomial from Chapter~\ref{chap:normal} can also recover its trace-dual basis. A primitive root of unity need not be primitive in $E^\times$: optimality here refers only to multiplication complexity.
\end{remark}

Period sums of roots of unity can also serve as observations from which irreducible-factor coefficients are recovered; Corollary~\ref{v97:cor-period-recovery} gives the corresponding Hankel system.

\section{Finding Normal Elements, Dual Bases, and Computational Costs}
\label{sec:construct-normal-dual}

Let \(F(T)=T^n-1\), with its factorization into powers of distinct monic irreducibles known.
Proposition~\ref{prop:normal-test} tests normality without first computing the full Frobenius annihilator polynomial:
\[
a\text{ is normal}\quad\Longleftrightarrow\quad
(F/P)(\sigma)a\ne0\quad\text{for every distinct irreducible factor }P\mid F.
\]
Here \(F/P\) means the quotient of the original \(F=T^n-1\) by \(P\): if \(P\) has multiplicity \(a\) in \(F\), its multiplicity in \(F/P\) is \(a-1\), and all other multiplicities are unchanged.

Enumerating and testing elements eventually finds a normal element. Independent uniform random candidates succeed with probability exactly
\[
\rho=\prod_{P\mid F}(1-q^{-\deg P}),
\]
so the expected number of draws is \(\rho^{-1}\). This exact probability follows from the normal-element count and does not require \(p\nmid n\).

\subsection{Normality and Exceptional Parameters among Reciprocals of Translates}\label{v114:sec-reciprocal-exceptions}
Starting from a fixed field generator, can a normal element be found simply by translating over the base field and taking a reciprocal? Which parameters succeed, and which conjugate components account for failure?
Let \(E=\F_{q^n}\), \(k=\F_q\), \(n\ge2\), and \([k(\alpha):k]=n\). For every \(t\in k\),
\[
 k\bigl((\alpha+t)^{-1}\bigr)=k(\alpha),
 \qquad \alpha=\bigl((\alpha+t)^{-1}\bigr)^{-1}-t.
\]
Every element in this family therefore generates the field; what changes is the linear span of its conjugates. Put
\begin{equation}\label{v114:def-bad-count}
 b_{q,n}(\alpha)=\#\{t\in k:(\alpha+t)^{-1}\text{ is not normal over }k\}.
\end{equation}
These values of \(t\) are called \textbf{exceptional parameters}. We first determine the exact parameter set from the local obstructions and use component ranks to bound its size. When \(n\mid q-1\), these conditions become nonvanishing of individual eigen-coordinates.

Proposition~\ref{v82:prop-normal-translates} shows that base-field translation may never remove an obstruction to normality.
Taking reciprocals need not preserve the same linear relations among conjugates, and gives a different one-parameter family of candidates from a power-basis generator.
For a fixed generator \(\theta\notin k\), the \textbf{family of reciprocals of translates} is obtained by translating \(\theta\) to \(\theta-t\), with \(t\in k\), and then taking reciprocals.
This family is used to construct normal elements in~\cite[Lemma~3.1]{GiesbrechtJamshidpeySchost2026}.

First apply the detecting operators of Proposition~\ref{prop:normal-test} to a common numerator. The common zeros of the resulting coordinate polynomials determine which parameters meet each obstruction kernel. Use \(Y\) for the operator variable, \(T\) for the translation parameter, and write \(\sigma(x)=x^q\).

\begin{proposition}[Exceptional parameters from primary components]\label{v115:prop-primary-reciprocal}
Let \(E=\F_{q^n}=k(\alpha)\), with \(n>1\), and put
\[
F(Y)=Y^n-1,\qquad D(T)=\prod_{i=0}^{n-1}(T+\alpha^{q^i}),\qquad
U(T)=\frac{D(T)}{T+\alpha}\in E[T].
\]
Choose a basis \(e_1,\ldots,e_n\) of \(E/k\). For each distinct monic irreducible factor \(P\) of \(F\), apply \((F/P)(\sigma)\) coefficientwise and write
\[
(F/P)(\sigma)U(T)=\sum_{i=1}^n u_{P,i}(T)e_i,
\qquad G_P(T)=\gcd_i u_{P,i}(T),
\]
with monic greatest common divisor. Each family \(u_{P,i}\) is not identically zero, and
\begin{equation}\label{v115:eq-primary-exceptions}
 b_{q,n}(\alpha)=\deg\gcd\!\left(T^q-T,\prod_{P\mid F}G_P(T)\right).
\end{equation}
The formula allows \(p\mid n\).
\end{proposition}
\begin{proof}
Write \(D=T^n+d_{n-1}T^{n-1}+\cdots+d_0\). Successive coefficient comparisons in \((T+\alpha)U=D\) give the coefficients of \(U\), from highest degree to constant term, as
\[
1,\quad d_{n-1}-\alpha,\quad
d_{n-2}-d_{n-1}\alpha+\alpha^2,\quad\ldots.
\]
They form an invertible triangular change from \(1,\alpha,\ldots,\alpha^{n-1}\), and hence span \(E\).
The minimal polynomial of Frobenius is \(F\). Since \(\deg(F/P)<n\), the operator \((F/P)(\sigma)\) is nonzero. It cannot annihilate every coefficient of \(U\), so \(G_P\) is defined.

For \(t\in k\), we have \(D(t)\ne0\), and base-field scaling commutes with \(\sigma\). Thus
\[
(F/P)(\sigma)(\alpha+t)^{-1}=0
\iff u_{P,i}(t)=0\text{ for every }i
\iff G_P(t)=0.
\]
Apply Proposition~\ref{prop:normal-test} and remove multiplicities by taking the greatest common divisor with \(T^q-T\). Throughout, \(F/P\) retains all other multiplicities.
\end{proof}

A change of basis makes invertible base-field linear combinations of the \(u_{P,i}\), so \(G_P\) is unchanged.
In general, first compute the coordinate greatest common divisors, then take a greatest common divisor with \(T^q-T\). Modular exponentiation computes the required residue of \(T^q\) without enumerating the base field.

The larger the image of a testing operator, the less room its coordinate numerators have for a common factor. This gives a uniform bound without requiring \(Y^n-1\) to split.
\begin{corollary}[Primary component ranks and reciprocal-family bounds]\label{v117:cor-primary-rank-bound}
Retain the notation of Proposition~\ref{v115:prop-primary-reciprocal}. Write \(n=p^a m\), with \(p\nmid m\), and let \(s\) be the number of distinct monic irreducible factors of \(Y^m-1\). Then
\begin{equation}\label{v117:eq-primary-rank-bound}
b_{q,n}(\alpha)\le
\begin{cases}
\min\{q,(s-1)(n-1)\},&p\nmid n,\\
\min\{q,s(n-1)-m\},&p\mid n.
\end{cases}
\end{equation}
In particular, if \(n>1\) is a power of the characteristic, then \(b_{q,n}(\alpha)\le\min\{q,n-2\}\).
\end{corollary}
\begin{proof}
For \(P\mid F\), put \(A_P=(F/P)(\sigma)\) and \(d_P=\deg P\). Proposition~\ref{prop:normal-test} gives \(\rank A_P=d_P\). The preceding proof shows that the coefficients of \(U(T)\) span \(E\), so those of \(A_PU(T)\) span \(\im A_P\).

Write \(A_PU=G_PR_P\). Multiplication by a nonzero base-field polynomial preserves the span of the coefficients: one inclusion follows by multiplication, and the reverse inclusion by successive elimination from the leading term. Thus the coefficients of \(R_P\) span a \(d_P\)-dimensional space, and
\[
d_P\le\deg R_P+1=\deg(A_PU)-\deg G_P+1.
\]
The leading coefficient of \(U\) is one, and \(A_P(1)=(F/P)(1)\). This is nonzero only when \(p\nmid n\) and \(P=Y-1\); in that case \(d_P=1\) and \(\deg G_P\le n-1\). In every other case \(\deg(A_PU)\le n-2\), whence
\[
\deg G_P\le n-1-d_P.
\]
Finally, \(F=(Y^m-1)^{p^a}\) and \(\sum_Pd_P=m\). Summing over \(P\) and applying \eqref{v115:eq-primary-exceptions} gives both bounds.
\end{proof}
When \(n\mid q-1\), we have \(s=n\) and recover the bound \((n-1)^2\); the common quartic discriminant constraint can still improve it. If \(q=n=5\), every full-degree element has at least two normal reciprocals among its five translates, since \(b_{5,5}(\alpha)\le3\).

\begin{corollary}[Normal elements among reciprocals of translates]\leavevmode\phantomsection\label{v68:prop-reciprocal-normal-family}
Let \(f\in k[X]\), with \(k=\F_q\), be monic irreducible of degree \(n>1\), and let \(E=k(\theta)\), where \(f(\theta)=0\).
For each \(t\in k\), put
\[
\beta_t=(\theta-t)^{-1}.
\]
At most \((n-1)^2\) parameters \(t\in k\) give an element \(\beta_t\) that is not normal over \(k\).
In particular, if \(q>(n-1)^2\), some member of this family is normal.
\end{corollary}
\begin{proof}
Replacing \(t\) by \(-t\) preserves the number of exceptions. If \(p\nmid n\), then \(s\le n\) in Corollary~\ref{v117:cor-primary-rank-bound}, giving \((s-1)(n-1)\le(n-1)^2\). If \(p\mid n\), then \(s\le m\le n-1\), so \(s(n-1)-m\le m(n-2)\le(n-1)^2\).
\end{proof}

\paragraph{A proof without factoring the operator polynomial.}
Since \(n>1\), we have \(\theta\notin k\), so every \(\beta_t\) is defined and \(f(t)\ne0\).
Write \(\theta_i=\theta^{q^i}\) for \(0\le i<n\), with subscripts read modulo \(n\).
For an indeterminate \(U\), define
\[
b_i(U)=\frac{f(U)}{\theta_i-U}
=-\prod_{j\ne i}(U-\theta_j)\in E[U],
\qquad B(U)=\bigl(b_{j-i}(U)\bigr)_{0\le i,j<n}.
\]
For \(t\in k\), the \((i,j)\)-entry of the Moore matrix of
\(\beta_t,\beta_t^q,\ldots,\beta_t^{q^{n-1}}\) is \((\theta_{i+j}-t)^{-1}\).
Permuting its rows by \(i\mapsto-i\) and multiplying every entry by \(f(t)\) gives \(B(t)\).
The Moore criterion (Theorem~\ref{thm:moore}) therefore makes normality of \(\beta_t\) equivalent to \(\det B(t)\ne0\).

All \(b_i(U)\) have the same leading term \(-U^{n-1}\). Subtract row zero from each of the other \(n-1\) rows. This preserves the determinant and lowers every entry in those rows to degree at most \(n-2\). Hence
\[
\deg\det B(U)\le(n-1)+(n-1)(n-2)=(n-1)^2.
\]
At \(U=\theta\), we have \(b_i(\theta)=0\) for \(i\ne0\) and \(b_0(\theta)=-f'(\theta)\ne0\), since \(f\) is separable.
Thus \(B(\theta)=-f'(\theta)I_n\), and its determinant is nonzero.
The polynomial \(\det B(U)\) is therefore nonzero and has at most \((n-1)^2\) roots in \(k\), proving both assertions.

The bound uses no restriction on whether the characteristic divides \(n\).
When \(q>(n-1)^2\), testing any \((n-1)^2+1\) distinct base-field parameters must produce a normal element.
For smaller base fields, this bound alone does not guarantee a successful parameter in \(k\).

For \(n=2\), the bound is one, so every quadratic extension has a normal element in this reciprocal family.
For example, in \(\F_9=\F_3(\xi)\) with \(\xi^2=-1\), parameters \(0,1,-1\) give \(-\xi,\xi+1,\xi-1\), respectively. Only the first is not normal, attaining the bound.

\begin{example}[A reciprocal change of generator over \(\F_7\)]\label{v68:ex-reciprocal-normal-f343}
Let \(k=\F_7\) and \(E=k(\theta)\), where \(\theta^3=2\).
The cubes in \(\F_7\) are \(0,1,-1\), so \(X^3-2\) has no root in \(\F_7\) and is irreducible.
Since \(\theta^7=4\theta\) and \((\theta^2)^7=2\theta^2\), the three power-basis vectors lie in distinct Frobenius eigenspaces.
First consider translation alone. For every \(t\in\F_7\), all conjugates of \(\theta+t\) lie in the two-dimensional space \(\langle1,\theta\rangle_{\F_7}\), so none of these translates is normal.
The missing \(\theta^2\)-component cannot be supplied by base-field translation.

For \(t\in\F_7\),
\[
\beta_t=\frac{\theta^2+t\theta+t^2}{2-t^3}.
\]
The coordinate matrix of its three conjugates in the basis \((1,\theta,\theta^2)\) is
\[
\frac1{2-t^3}
\begin{pmatrix}
t^2&t^2&t^2\\
t&4t&2t\\
1&2&4
\end{pmatrix},
\qquad
\det=\frac{t^3}{(2-t^3)^3}.
\]
Hence \(\beta_t\) is normal exactly when \(t\ne0\).
For \(t=1\), this gives the normal basis
\[
1+\theta+\theta^2,\qquad
1+4\theta+2\theta^2,\qquad
1+2\theta+4\theta^2.
\]
For \(t=0\), the element \(\beta_0=4\theta^2\) still has degree three, but its conjugates span only the line \(k\theta^2\).
The reciprocal change of variable thus alters the linear relations among the conjugates while preserving the degree of the generator.
When \(t\ne0\), the coefficients of \(1,\theta,\theta^2\) in the reciprocal expression are all nonzero, so all three eigencomponents occur.
Simple translation remains in a fixed proper linear subspace, whereas taking reciprocals changes the Frobenius components. This explains why the reciprocal family constructs normal elements in this example.
\end{example}

Returning to the notation of Corollary~\ref{v68:prop-reciprocal-normal-family}, taking reciprocals after translation can also control the trace. Equation~\eqref{v70:eq-reciprocal-trace} gives
\(\Tr_{E/k}((\theta-t)^{-1})=-f'(t)/f(t)\), so a parameter with \(f'(t)=0\) cannot yield a normal element:
zero trace means that the sum of all conjugates is zero, making them linearly dependent when \(n>1\).
Proposition~\ref{v70:prop-trace-zero-generator} preserves field generation while making the trace zero, a different requirement from the normality sought here.

If \(n>1\) is a power of the characteristic, normality is already equivalent to nonzero trace. Thus, for \(f=m_\theta\), the exact exceptional set with the minus-sign parameter is
\begin{equation}\label{v128:eq-power-degree-exceptions}
\{t\in k:(\theta-t)^{-1}\text{ is not normal}\}
=\{t\in k:f'(t)=0\},\qquad
b_{q,n}(\theta)=\deg\gcd(T^q-T,f').
\end{equation}
For example, \(f=X^3-X-1\in\F_3[X]\) takes the value \(-1\) throughout the base field and is therefore irreducible. If \(f(\theta)=0\), then \(\Tr(\theta+t)=0\) for all \(t\in\F_3\), so no simple translate is normal. Yet
\[
\Tr((\theta-t)^{-1})=-\frac{f'(t)}{f(t)}=-1,
\]
so all three reciprocals are normal. Here \(Y^3-1=(Y-1)^3\); diagonalizing Frobenius is unnecessary.

\paragraph{Eigen-coordinates in the split case.}
If \(n\mid q-1\), every local quotient is one-dimensional, and the exact exceptional-set formula becomes a zero condition on individual eigen-coordinates.
Choose a primitive \(n\)-th root of unity \(\zeta\in k\), and write \(\sigma(x)=x^q\). By the cyclic-module isomorphism \(E\simeq k[T]/(T^n-1)\) of Chapter~\ref{chap:normal}, each \(\zeta^j\)-eigenspace is one-dimensional. Choose a nonzero \(\theta\) with \(\theta^q=\zeta\theta\).
The elements \(\theta^j\) have distinct eigenvalues \(\zeta^j\), so \(1,\theta,\ldots,\theta^{n-1}\) form a \(k\)-basis, and \(d=\theta^n\in k^\times\). In particular, \(X^n-d\) is the minimal polynomial of \(\theta\).
For \(x=\sum_{j=0}^{n-1}x_j\theta^j\), the conjugate-coordinate matrix is \(\operatorname{diag}(x_0,\ldots,x_{n-1})(\zeta^{ij})_{j,i}\). The second matrix is an invertible Vandermonde matrix. Hence
\begin{equation}\label{v114:eq-eigen-normal}
 x\text{ is normal}\quad\Longleftrightarrow\quad x_j\ne0\quad(0\le j<n).
\end{equation}
An orbit of length \(n\) need not involve every eigenspace; normality requires every component to be nonzero.

\begin{corollary}[Coordinate numerators and exact exception counts]\label{v114:prop-reciprocal-count}
Suppose \(n\ge2\), \(n\mid q-1\), and \([k(\alpha):k]=n\), and choose the eigenbasis above. Put
\[
 D_\alpha(T)=\prod_{i=0}^{n-1}(T+\alpha^{q^i})\in k[T],
 \qquad \frac{D_\alpha(T)}{T+\alpha}=\sum_{j=0}^{n-1}p_j(T)\theta^j,
 \quad p_j\in k[T].
\]
Every \(p_j\) is nonzero, and
\[
 p_0=\frac{D_\alpha'}n,\qquad \deg p_0=n-1,
 \qquad \deg p_j\le n-2\quad(j>0).
\]
The exceptional parameters are exactly the union of the base-field root sets of the \(p_j\). Thus
\begin{equation}\label{v114:eq-bad-gcd}
 \boxed{b_{q,n}(\alpha)=\deg\gcd\!\left(T^q-T,\prod_{j=0}^{n-1}p_j(T)\right)
 \le\min\{q,(n-1)^2\}.}
\end{equation}
\end{corollary}
\begin{proof}
In Proposition~\ref{v115:prop-primary-reciprocal}, every detecting operator now has a one-dimensional eigenspace as its image, so \(G_P\) is the monic normalization of the corresponding coordinate numerator \(p_j\). It remains to determine their degrees.

Write \(\alpha_i=\alpha^{q^i}\). Frobenius acts on coefficients and fixes \(T\), so
\[
 \prod_{r\ne i}(T+\alpha_r)=\sum_{j=0}^{n-1}p_j(T)\zeta^{ij}\theta^j.
\]
Root-of-unity orthogonality gives
\begin{equation}\label{v114:eq-coordinate-projection}
 n\theta^j p_j(T)=\sum_{i=0}^{n-1}\zeta^{-ij}\prod_{r\ne i}(T+\alpha_r).
\end{equation}
At \(T=-\alpha\), only the term \(i=0\) remains, yielding
\[
 n\theta^j p_j(-\alpha)=\prod_{r=1}^{n-1}(\alpha_r-\alpha)\ne0.
\]
Here \(n\ne0\) in \(k\), and full degree makes the conjugates distinct. Thus every \(p_j\) is nonzero.
For \(j=0\), the right-hand side of~\eqref{v114:eq-coordinate-projection} is \(D_\alpha'\), with leading coefficient \(n\); hence \(p_0\) is monic of degree \(n-1\). For \(j>0\), the leading coefficient \(\sum_i\zeta^{-ij}\) vanishes, giving degree at most \(n-2\).

For \(t\in k\), we have \(D_\alpha(t)\ne0\): otherwise a conjugate would belong to \(k\), contrary to \(n>1\). Consequently,
\[
 (\alpha+t)^{-1}=D_\alpha(t)^{-1}\sum_{j=0}^{n-1}p_j(t)\theta^j.
\]
Equation~\eqref{v114:eq-eigen-normal} gives the asserted union of root sets. The roots of \(T^q-T\) are exactly the elements of \(k\), all simple, so the greatest common divisor both restricts the roots to the base field and removes repeated counting.
Finally, the sum of the numerator degrees is at most \((n-1)+(n-1)(n-2)=(n-1)^2\).
\end{proof}

In these coordinates, the trace reads only the constant component:
\[
 \Tr_{E/k}\bigl((\alpha+t)^{-1}\bigr)
 =\frac{n p_0(t)}{D_\alpha(t)}=\frac{D_\alpha'(t)}{D_\alpha(t)}.
\]
If \(f=m_\alpha\), then \(D_\alpha(T)=(-1)^n f(-T)\), agreeing with the earlier trace formula after replacing \(t\) by \(-t\). Nonzero trace guarantees only \(p_0(t)\ne0\); the other components still have to be checked.
Failure to attain the bound can result from lower numerator degrees, too few distinct base-field roots, or roots shared by different numerators. Moreover, the numerators arise from one reciprocal expression, so their root counts cannot always attain their individual degree bounds simultaneously. Proposition~\ref{v115:prop-quartic-bound} gives a concrete quartic constraint.

\begin{example}[One and four exceptions in the same extension]\label{v114:ex-cubic-comparison}
Use \(\F_7(\theta)\), \(\theta^3=2\), from Example~\ref{v68:ex-reciprocal-normal-f343}. That example gives \(b_{7,3}(\theta)=1\); with the plus-sign parameter the exceptional set is still \(\{0\}\).
Take \(\alpha=\theta+\theta^2\). It is outside \(\F_7\), and the extension degree three is prime, so it again has full degree. Multiplication and reduction by \(\theta^3=2\) give
\[
 (t+\theta+\theta^2)\bigl((t^2-2)+(2-t)\theta+(1-t)\theta^2\bigr)
 =t^3+t+6.
\]
Thus
\[
 \frac1{\alpha+t}
 =\frac{(t^2-2)+(2-t)\theta+(1-t)\theta^2}{t^3+t+6}.
\]
The three numerator root sets in the base field are \(\{3,4\}\), \(\{2\}\), and \(\{1\}\), respectively. Therefore
\[
 b_{7,3}(\theta+\theta^2)=4=(3-1)^2.
\]
The two generators give the comparison
\[
\begin{array}{c|c|c}
 \text{Generator}&\text{Exceptional parameters}&b_{7,3}\\\hline
 \theta&\{0\}&1\\
 \theta+\theta^2&\{1,2,3,4\}&4
\end{array}
\]
At \(t=1,2\), the constant numerator is nonzero, as is the trace of the reciprocal, but the \(\theta^2\)- or \(\theta\)-component, respectively, is absent. Thus not every failure of normality in this example is caused by zero trace.
\end{example}

\begin{example}[A reciprocal family containing no normal element]\label{v114:ex-quartic-failure}
Over \(\F_5\), take \(\theta^4=2\) in an algebraic closure. Since \(\theta^5=2\theta\) and \(2\) has order four, the Frobenius orbit of \(\theta\) has length four; hence \(E=\F_5(\theta)=\F_{5^4}\).
Put \(\alpha=\theta+\theta^3\). Independence of the power basis gives \(\alpha\ne0\), while \(\alpha^{25}=-\alpha\ne\alpha\). Thus \(\alpha\) lies outside the unique quadratic subfield and also has full degree.
Multiplication and reduction by \(\theta^4=2\) verify
\[
 \frac1{\alpha+t}
 =\frac{(t^3+t)+(3-t^2)\theta+3t\theta^2+(1-t^2)\theta^3}
 {t^4+2t^2+3}.
\]
The denominator has no root in \(\F_5\). The following numerator root sets cover the base-field parameters:
\[
\begin{array}{c|c}
 \text{Parameters}&\text{One vanishing coordinate numerator}\\\hline
 t=0&3t\\
 t=\pm1&1-t^2\\
 t=\pm2&t^3+t
\end{array}
\]
Hence \(b_{5,4}(\alpha)=5\). The normal basis theorem guarantees normal elements in \(E\), but does not guarantee one in this candidate family from an arbitrarily prescribed generator.
The example shows that complete failure can occur over a small base field; it does not assert that \((n-1)^2\) is the optimal threshold for guaranteed success.
\end{example}

\paragraph{A fixed generator, an average, and a maximum.}\label{v114:par-extremal-question}
Base-field affine changes preserve the failure count: if \(c\in k^\times\) and \(a\in k\), then
\[
 (c\alpha+a+t)^{-1}=c^{-1}\bigl(\alpha+(a+t)/c\bigr)^{-1},
 \qquad b_{q,n}(c\alpha+a)=b_{q,n}(\alpha).
\]
The parameter changes bijectively, and multiplication by a nonzero base-field scalar preserves normality. Thus the constant coordinate of a generator may first be removed in low-degree calculations.
Put \(\mathcal G_{q,n}=\{\alpha\in E:[k(\alpha):k]=n\}\), and define
\begin{equation}\label{v114:def-max-bad}
 B_q(n)=\max_{\alpha\in\mathcal G_{q,n}}b_{q,n}(\alpha).
\end{equation}
Exercise~\ref{v114:ex-reciprocal-average} computes the average failure count over all generators by double counting. This average does not determine the successful parameters of a fixed family, nor does it replace the maximum.
Exercises~\ref{v114:ex-quadratic-count}--\ref{v114:ex-quartic-bound} resolve successive low-degree cases: an exact quadratic count, a cubic classification by characteristic and square class, and a quartic improvement from nine to seven using a relation between two numerator discriminants.
These computations turn the existence argument into a testable extremal question; the end of the chapter asks further about simultaneous splitting of the numerators.

The bound estimates each numerator separately, without using relations among them. In degree four, the square classes of two discriminants already force a strict improvement.
\begin{proposition}[A common constraint in quartic reciprocal families]\label{v115:prop-quartic-bound}
If \(4\mid q-1\) and \(\alpha\) has full degree in \(\F_{q^4}/\F_q\), then
\[
b_{q,4}(\alpha)\le\min\{q,7\}.
\]
\end{proposition}
\begin{proof}

Choose a primitive fourth root of unity \(\zeta\in k\) and an eigenbasis \(1,\theta,\theta^2,\theta^3\) with \(\theta^q=\zeta\theta\), and put \(d=\theta^4\in k\). Base-field affine changes preserve \(b_{q,4}\); after removing the constant coordinate, write \(\alpha=a\theta+b\theta^2+c\theta^3\), where \(a,b,c\in k\). Then
\[
 d^{(q-1)/2}=\theta^{2(q-1)}=\zeta^2=-1.
\]
Thus \(d\) is a nonsquare. If \(a=c=0\), then \(\alpha=b\theta^2\in\F_{q^2}\), contrary to full degree.
Write \(D\) for the common norm denominator. Multiplication gives
\[
 \frac{D(T)}{T+\alpha}=p_0(T)+p_1(T)\theta+p_2(T)\theta^2+p_3(T)\theta^3,
\]
where
\begin{equation}\label{v114:eq-quartic-numerators}
\begin{aligned}
 p_0&=T^3-d(2ac+b^2)T+bd(a^2+c^2d),\\
 p_1&=-aT^2+2bcdT+d(a^2c-ab^2-c^3d),\\
 p_2&=-bT^2+(a^2+c^2d)T+bd(b^2-2ac),\\
 p_3&=-cT^2+2abT-a^3+ac^2d-b^2cd.
\end{aligned}
\end{equation}
These formulas can also be checked directly: after multiplication by \(T+\alpha\), all nonconstant \(\theta\)-coordinates vanish, leaving a monic quartic \(S(T)\). Substitution \(T=-\alpha\) gives \(S(-\alpha)=0\); since \(-\alpha\) has degree four, \(S=D\).

First suppose \(ac\ne0\). Both \(p_1,p_3\) are quadratic, with discriminants
\begin{equation}\label{v114:eq-quartic-discriminants}
\begin{aligned}
 \Delta_1&=4d(a^2-c^2d)(ac-b^2),\\
 \Delta_3&=-4(a^2-c^2d)(ac-b^2),\qquad \Delta_1=-d\Delta_3.
\end{aligned}
\end{equation}
Since \(4\mid q-1\), \(-1\) is a square, so \(-d\) is a nonsquare. If the discriminants are nonzero, exactly one is a square, and the two numerators have a total of two distinct base-field roots, counted separately. If the discriminants vanish, each numerator has one distinct root, again giving total two.

If \(a=0,c\ne0\), then
\[
 p_1=cd(2bT-c^2d),\qquad p_3=-c(T^2+b^2d).
\]
For \(b\ne0\), the first has one root and the second has none, because \(-d\) is a nonsquare. For \(b=0\), the first is a nonzero constant and the second has only the distinct root zero.
If \(c=0,a\ne0\), then
\[
 p_1=-a(T^2+b^2d),\qquad p_3=a(2bT-a^2),
\]
and their total root count is again one. These cases exhaust \((a,c)\ne(0,0)\).
By Corollary~\ref{v114:prop-reciprocal-count}, all numerators are nonzero, and \(p_0,p_2\) have at most three and two roots, respectively. Their combined root set therefore has size at most \(3+2+2=7\), and also at most \(q\).

\end{proof}

\begin{example}[Seven exceptional parameters]\label{v115:ex-quartic-seven}
If \(\theta^4=3\), then \(\theta^{17}=13\theta\), and \(13\) has order four in \(\F_{17}^\times\), so \(\theta\) has degree four.
For \(\alpha=\theta+3\theta^2-4\theta^3\), \(\sigma^2\) fixes the even powers and negates the odd powers. Independence of the power basis gives \(\alpha^{17^2}\ne\alpha\), so \(\alpha\) also has full degree.
Equation~\eqref{v114:eq-quartic-numerators} becomes
\[
\begin{array}{c|c|c}
 j&p_j(T)&\text{Roots in }\F_{17}\\\hline
 0&T^3-3T-1&3,4,10\\
 1&-T^2-4T-7&\varnothing\\
 2&-3T^2-2T&0,5\\
 3&4T^2+6T+2&8,16
\end{array}
\]
The first row factors as \((T-3)(T-4)(T-10)\); the second has nonsquare discriminant \(5\); the last two factor directly. The seven roots are distinct, proving \(b_{17,4}(\alpha)=7\) and hence \(B_{17}(4)=7\).
\end{example}

Degrees bound the roots of each numerator; repeated roots, common roots, and splitting restrict their union further. Writing \(Z_j=\{t\in k:p_j(t)=0\}\), we have
\[
\left|\bigcup_j Z_j\right|\le\sum_j|Z_j|\le\sum_j\deg p_j\le(n-1)^2.
\]
The quartic proposition improves the second inequality by a relation between discriminants. In higher degree, the relevant numerators still come from one reciprocal expression; they are not arbitrary polynomials of bounded degree.

\subsection*{Computing Dual Bases}

Choose any basis \(\mathcal B=(b_0,\ldots,b_{n-1})\) and compute the Gram matrix \(G=(\Tr(b_i b_j))\). If \(H=G^{-1}\), its dual basis is
\[
b_j^*=\sum_i H_{ij}b_i,
\qquad x=\sum_j\Tr(xb_j^*)b_j.
\]
If \(\mathcal B\) is normal, only \(b_0^*\) needs to be computed; Frobenius iteration gives the remaining dual vectors.
Frobenius is a cyclic shift in normal coordinates, but multiplication still requires a multiplication table or conversion to a power basis; it is not coordinatewise multiplication.

\subsection*{Measuring Arithmetic Costs}
Count each addition, subtraction, multiplication, or division in \(k\) as one field operation,
and let \(M(n)\) bound the number of base-field operations needed for one multiplication in \(k[X]/(f)\).
Naive multiplication and reduction give \(M(n)=O(n^2)\). Binary exponentiation computes a \(q\)th power in \(O(M(n)\log q)\) base-field operations.
The table lists upper bounds for these direct implementations.
\begin{center}\small
\begin{tabular}{@{}p{0.43\textwidth}p{0.50\textwidth}@{}}
\toprule
Task and known input&Base-field operation cost\\
\midrule
Generate a conjugate orbit of length \(d\)&\(O(dM(n)\log q)\)\\
Multiply root factors successively along an orbit&\(O(d^2M(n))\)\\
Given a target root, precompute the isomorphism matrix&\(O(nM(n))\)\\
Change coordinates of one element using a known matrix&\(O(n^2)\)\\
Invert a change-of-basis matrix or a known Gram matrix&\(O(n^3)\)\\
Test a normal-element candidate given the factorization of \(F\)&\(O(nM(n)\log q+\nu n^2)\)\\
\bottomrule
\end{tabular}
\end{center}
Here \(\nu\) is the number of distinct irreducible factors of \(F\); the last row first computes all conjugates and then forms a linear combination for each \(F/P\).
Root finding, factoring \(F\), and constructing the Gram matrix are separate costs, not included in the linear-algebra costs once these data are known.
Testing multiplicative primitivity requires the distinct prime factors of \(q^n-1\); computing those integer factors is also outside finite-field arithmetic costs.
A bit-operation model must additionally account for the base-field representation and the bit complexity of each base-field operation.

\section{Random Computation and Canonical Output}\label{v58:sec-canonical-randomness}
Random factorization may discover the same factors in different orders; sorting the complete factor set still gives a fixed output.
First fix a field representation: use the integer representatives \(0,1,\ldots,p-1\) in \(\F_p\),
and compare elements of the quotient-ring model \(\F_p[Z]/(g)\) by lexicographic order on coefficient vectors \((a_0,\ldots,a_{r-1})\).
This order depends on the representation and is not compatible with field operations.

\begin{definition}[Pseudo-deterministic algorithms]\label{v58:def-pseudo-deterministic}
A randomized algorithm is \textbf{pseudo-deterministic} if, for every input, there is a valid output determined by that input that the algorithm returns with probability at least \(2/3\).

This fixed output is called the \textbf{canonical output} for the input.
\end{definition}
A fixed output is compatible with the use of randomness during computation.
A Las Vegas algorithm requires correct results but need not return the same result on different runs;
selecting from the complete candidate set by a fixed rule can leave randomness affecting only running time.
Rai~\cite{RaiPseudo2024} gives a pseudo-deterministic construction of irreducible polynomials of any prescribed degree.
The quadratic construction below illustrates canonical selection.

\begin{proposition}[A fixed nonsquare from successive root extraction]\label{v58:prop-canonical-nonsquare}
Let \(p\) be an odd prime and write \(p-1=2^s m\) with \(m\) odd.
Start from \(a=-1\). If \(a\) is a square, find both roots of \(X^2-a\)
and replace \(a\) by the root with the smaller integer representative; if \(a\) is a nonsquare, return \(X^2-a\).
The procedure stops after exactly \(s-1\) root extractions, and the quadratic irreducible polynomial returned depends only on \(p\).
\end{proposition}
\begin{proof}
The initial element \(-1\) has order two.
If \(a\) has order \(2^j\) and \(b^2=a\), then
\(\ord(b)/\gcd(\ord(b),2)=2^j\), so \(\ord(b)=2^{j+1}\).
The square subgroup of the cyclic group \(\F_p^\times\) has order \(2^{s-1}m\),
so an element of order \(2^j\) is a square if and only if \(j<s\).
Each root extraction therefore doubles the order, and the procedure stops exactly when the order reaches \(2^s\).
Each step makes a fixed choice from the complete pair of roots, so the whole sequence and final output are independent of random choices.
The final \(a\) is a nonsquare, making \(X^2-a\) irreducible.
\end{proof}

Test whether an element is a square using \(a^{(p-1)/2}\) and find both roots using the random splitting in Proposition~\ref{prop:split-values-random}.
If splitting fails, sample again; choose the smaller root only after obtaining two roots verified by substitution.
The root-finding subroutine is Las Vegas, so it terminates with probability one and finite expected running time.
Combining at most \(\log_2p\) root extractions with binary exponentiation gives expected bit complexity polynomial in \(\log p\) for the entire construction.
Suppose a known bound \(T(p)\) on the expected bit-operation count is available, where \(T(p)\) is a polynomial in \(\log p\).
Report failure if execution has not finished after \(3T(p)\) steps. Markov's inequality bounds the timeout probability by \(1/3\),
so the probability of returning the fixed valid output within this limit is at least \(2/3\).

\begin{example}[A fixed quadratic construction in \(\F_{73}\)]\label{v58:ex-canonical-f73}
Since \(72=2^3\cdot9\), the procedure is
\[
72=-1\longmapsto27\longmapsto10,
\qquad\ord(72)=2,\quad\ord(27)=4,\quad\ord(10)=8.
\]
Indeed, the two roots of \(-1\) are \(27,46\), and those of \(27\) are \(10,63\);
choosing the smaller root gives the displayed sequence.
Also \(10^{36}=(10^4)^9=(-1)^9=-1\), so the final output is \(X^2-10\).
This example appears in \cite[Example 2.1]{RaiPseudo2024}; its increasing-order mechanism also parallels multiplicative period extension in Chapter~\ref{chap:generators}.
\end{example}

\begin{remark}[Root partitions and elliptic-curve root-finding]\label{v105:rem-elliptic-root-finding}
The root partitions obtained by greatest common divisors can also begin a change of the space in which roots are studied.
Suppose that the characteristic of \(\F_q\) exceeds three and that \(f\in\F_q[X]\) is monic, squarefree, and split over \(\F_q\).
Choose \(a,b\in\F_q\), put \(h(X)=X^3+aX+b\) with \(4a^3+27b^2\ne0\), and fix a nonsquare \(\beta\in\F_q^\times\). Set
\[
 f_0=\gcd(f,h),\qquad g=f/f_0,
 \qquad f_\pm=\gcd\bigl(g,h^{(q-1)/2}\mp1\bigr).
\]
Then
\begin{equation}\label{v105:eq-elliptic-root-partition}
 f=f_0f_+f_-,\qquad
 \gcd(f_0,f_+)=\gcd(f_0,f_-)=\gcd(f_+,f_-)=1.
\end{equation}
Indeed, at each root \(x\), the value \(h(x)\) is either zero, a nonzero square, or a nonsquare.
Euler's criterion describes the last two cases by \(h(x)^{(q-1)/2}=1\) and \(-1\), respectively, and squarefreeness gives the factorization.
Each root of \(f_+\) is the first coordinate of an \(\F_q\)-point on \(Y^2=h(X)\).
Each root of \(f_-\) is the first coordinate of an \(\F_q\)-point on the quadratic twist \(\beta Y^2=h(X)\), since the quotient of two nonsquares is a square.
The roots of \(f_0\) correspond to \(Y=0\) and can be treated separately.
Only modular exponentiation modulo \(g\) is required; this partition does not require the individual roots of \(f\) to be known.

Spaenlehauer~\cite{Spaenlehauer2026Roots} develops this approach by encoding point sets with root polynomials and organizing deterministic root-finding through the group operations on an elliptic curve and its quadratic twist.
The method separates field-dependent preprocessing from the subsequent computation for a given polynomial.
The preprocessing data include the curves, a nonsquare, the two rational-point group orders with their prime factorizations, and group bases.
Given these data, the online algorithm is deterministic; its cost also depends on the largest prime dividing the two group orders.
Constructing suitable preprocessing data efficiently therefore remains a central problem in this approach.
This connects with the preceding discussion of canonical outputs: determinism of the output, determinism of the root-finding procedure, and the cost of preprocessing are separate questions.
The complete curve-based algorithm and its complexity are given in the cited paper.
\end{remark}

\subsection*{Sparse Representations and Complexity}
\phantomsection\label{v58:remark-sparse-complexity}
The preceding costs are measured by degree or by the number of base-field operations.
In a \textbf{sparse representation}, only nonzero coefficients and their binary-encoded exponents are stored,
so input length can be much smaller than degree.
For example, \(X^{2^u}-1\) has just two terms, and its exponent needs only \(u+1\) bits;
a dense representation storing all coefficients needs \(2^u+1\) coefficient positions.
A polynomial bound in the degree need not therefore be polynomial in sparse input length.

For sparse greatest common divisors, Qiu et al.~\cite{QiuSparseGCD2026} prove a conditional complexity obstruction in a model with the finite-field size also part of the input and exponents encoded in binary:
if general sparse greatest common divisors could be computed in randomized polynomial bit time in the sparse lengths of the input and output,
then \(\mathrm{NP}\subseteq\mathrm{BPP}\).
\leavevmode\phantomsection\label{v58:def-np-bpp}\ignorespaces
Here \(\mathrm{NP}\) (nondeterministic polynomial time) consists of decision problems admitting a deterministic polynomial-time verifier and a polynomial \(p\)
such that every yes-instance \(x\) has an accepted certificate of length at most \(p(|x|)\), whereas no no-instance has an accepted certificate.
The class \(\mathrm{BPP}\) (bounded-error probabilistic polynomial time) consists of decision problems solvable by randomized polynomial-time algorithms
with error probability at most \(1/3\) for every input, taken over the algorithm's internal random choices.
The conclusion uses \(\mathrm{NP}\nsubseteq\mathrm{BPP}\) as the complexity assumption and concerns finite fields varying with the input;
it does not establish the same obstruction over every fixed finite field.
The preceding algorithms measure cost by polynomial degree, whereas this statement measures bit complexity by sparse input and output length.

\section{A Running Example: Representations and Orbits in the Field of Order 64}
\label{sec:f64-computation}

Example~\ref{exm:f64-structure} determined the subfields and periods of the field of order 64. We now choose a concrete quotient-ring model to make those structures computationally explicit.

\begin{example}\leavevmode\phantomsection\label{exm:f64-computation}
Take
\[
f(X)=X^6+X+1\in\F_2[X],\qquad
E=\F_2[X]/(f),\qquad\alpha=\overline X.
\]
First compute in the quotient ring, then verify that it is a field. The relation \(\alpha^6=\alpha+1\) gives
\[
\begin{array}{c|cccc}
j&8&16&32&64\\\hline
\alpha^j&\alpha^3+\alpha^2&\alpha^4+\alpha+1&\alpha^3+1&\alpha
\end{array}
\]
We also have \(\gcd(f,X^4+X)=\gcd(f,X^8+X)=1\).
The first identity follows using \(X^4+X=X(X^3+1)\) and \(f\bmod(X^3+1)=X\);
for the second, first reduce \(X^8+X\) modulo \(f\) to \(X(X^2+X+1)\), then use \(f\bmod(X^2+X+1)=X\).
Rabin's criterion proves that \(f\) is irreducible.

Further computation gives
\[
\alpha^9=\alpha^4+\alpha^3\ne1,\qquad
\alpha^{21}=\alpha^5+\alpha^4+\alpha^3+\alpha+1\ne1.
\]
Since \(63=3^2\cdot7\), the primitivity criterion shows that \(\alpha\) has multiplicative order \(63\).
But \(\Tr_{E/\F_2}(\alpha)=0\), so \(\alpha\) is not normal.
\end{example}

\subsection*{Subfields and Factors After Enlarging the Base Field}
Put \(u=\alpha^{21}\). Then \(u\) has order three, so \(\F_4=\F_2(u)\) and \(u^2+u+1=0\). Similarly, \(v=\alpha^9\) has order seven,
and \(2\) has order three modulo seven, so \(\F_8=\F_2(v)\).
These subfields have intersection \(\F_2\) and compositum \(E\).

Over \(\F_4\), the orbit of \(\alpha\) is \(\alpha,\alpha^4,\alpha^{16}\), so
\[
\Tr_{E/\F_4}(\alpha)=1,\qquad \Nm_{E/\F_4}(\alpha)=u.
\]
The sum of the pairwise products is \(\alpha^5+\alpha^{17}+\alpha^{20}=u^2\); hence
\begin{equation}\label{eq:f64-factor-f4}
f(X)=(X^3+X^2+u^2X+u)(X^3+X^2+uX+u^2).
\end{equation}
Squaring the coefficients interchanges the two cubics. Both are irreducible, since they correspond to two relative Frobenius orbits of length three.
Nonzero relative trace and zero absolute trace are compatible, since \(\Tr_{\F_4/\F_2}(1)=0\).

\subsection*{Two compatible block matrices for one operator}
\phantomsection\label{v106:par-tower-matrices}
The relation \(v^3=v^2+1\) used in Example~\ref{v106:ex-complete-f64} follows directly from \(v=\alpha^9\). Moreover,
\[
\alpha=u+(u+1)v^2=v^2+uv^3.
\]
In the basis \((1,v,v^2)\) of \(E/\F_4\) and the basis \((1,u)\) of \(E/\F_8\), multiplication by \(\alpha\) has the respective matrices
\begin{equation}\label{v106:eq-tower-matrices}
M_\alpha=\begin{pmatrix}u&u+1&u+1\\0&u&u+1\\u+1&u+1&1\end{pmatrix},
\qquad N_\alpha=\begin{pmatrix}v^2&v^3\\v^3&1\end{pmatrix}.
\end{equation}
Column \(j\) records the coordinates of \(\alpha\) times basis vector \(j\). For instance, \(\alpha v=(u+1)+uv+(u+1)v^2\) gives the second column of \(M_\alpha\).

The binary matrices for multiplication by \(u\) in \((1,u)\) and by \(v\) in \((1,v,v^2)\) are
\[
U=\begin{pmatrix}0&1\\1&1\end{pmatrix},\qquad
V=\begin{pmatrix}0&0&1\\1&0&0\\0&1&1\end{pmatrix}.
\]
Replace \(u\) and \(1\) in \(M_\alpha\) by \(U\) and \(I_2\), and replace \(v^j\) in \(N_\alpha\) by \(V^j\). The resulting six-by-six base-field matrices correspond to
\[
\begin{aligned}
B&=(1,u,v,uv,v^2,uv^2),\\
C&=(1,v,v^2,u,uv,uv^2)=(B_1,B_3,B_5,B_2,B_4,B_6).
\end{aligned}
\]
Thus simultaneous row and column reordering by \((1,3,5,2,4,6)\) carries one matrix to the other. This permutation compatibility comes from the chosen product bases; a general basis change gives a similarity transformation.

The traces, determinants, and characteristic polynomials give
\[
\begin{aligned}
\operatorname{tr}M_\alpha&=1,&\det M_\alpha&=u,&
\det(XI_3-M_\alpha)&=X^3+X^2+u^2X+u,\\
\operatorname{tr}N_\alpha&=v^3,&\det N_\alpha&=v,&
\det(XI_2-N_\alpha)&=X^2+v^3X+v.
\end{aligned}
\]
The first polynomial agrees with the relative factor in~\eqref{eq:f64-factor-f4}. Continuing down either tower gives absolute trace \(0\) and absolute norm \(1\). The regular representation of Section~\ref{sec:regular-representation} thus displays two subfield structures within the same six-dimensional operator.

\subsection*{Power Bases, Normal Bases, and Dual Coordinates}
By the power-basis dual formula, or by direct trace computation, the dual basis of \((1,\alpha,\ldots,\alpha^5)\) is
\[
(\alpha^5+1,\alpha^4,\alpha^3,\alpha^2,\alpha,1).
\]
For example, the constant coordinate of \(x=\sum_{i=0}^5x_i\alpha^i\) is \(x_0=\Tr(x(\alpha^5+1))\) and its highest-degree coordinate is \(x_5=\Tr(x)\).

Here \(f'(X)=1\), so the latter identity directly illustrates Proposition~\ref{prop:coefficient-trace}.

Take \(\gamma=\alpha^5\). Its successive conjugates are
\[
\begin{aligned}
\gamma&=\alpha^5,&\gamma^2&=\alpha^5+\alpha^4,\\
\gamma^4&=\alpha^5+\alpha^4+\alpha^3+\alpha^2,&
\gamma^8&=\alpha^5+\alpha^3+\alpha^2+\alpha+1,\\
\gamma^{16}&=\alpha^5+\alpha^2+\alpha,&
\gamma^{32}&=\alpha^5+\alpha^2.
\end{aligned}
\]
Here \(T^6-1=(T+1)^2(T^2+T+1)^2\). The two normality tests give, respectively,
\[
\Tr(\gamma)=1,\qquad
(1+\sigma+\sigma^3+\sigma^4)(\gamma)=\alpha^4+\alpha^3+1\ne0.
\]
These six conjugates therefore form a normal basis. To compute its dual normal element, first find the cyclic correlation coefficients:
\[
\bigl(\Tr(\gamma\gamma^{2^i})\bigr)_{i=0}^5=(1,1,1,0,1,1).
\]
The cyclic correlation polynomial of Proposition~\ref{prop:normal-circulant-test} is thus
\[
g_\gamma(T)=1+T+T^2+T^4+T^5.
\]
Squaring in characteristic two and using \(T^6\equiv1\pmod{T^6-1}\) gives
\[
g_\gamma(T)^2=1+T^2+T^4+T^8+T^{10}\equiv1\pmod{T^6-1}.
\]
Thus \(g_\gamma\) is its own inverse modulo \(T^6-1\). By the same proposition, the dual normal element is
\[
\begin{aligned}
\gamma^*&=g_\gamma(\sigma)(\gamma)
=\Tr(\gamma)+\sigma^3(\gamma)=1+\gamma^8\\
&=\alpha^5+\alpha^3+\alpha^2+\alpha=\alpha^{55}.
\end{aligned}
\]
The six corresponding pairing values are
\(\bigl(\Tr(\gamma^*\gamma^{2^j})\bigr)_{j=0}^5=(1,0,0,0,0,0)\).
Frobenius invariance then gives the full dual basis \((\gamma^*,(\gamma^*)^2,\ldots,(\gamma^*)^{32})\) and the normal coordinates of every element,
\[
x=\sum_{i=0}^5\Tr\bigl(x(\gamma^*)^{2^i}\bigr)\gamma^{2^i}.
\]
Although Frobenius is a cyclic shift in this basis, its minimal polynomial \(T^6-1\) has repeated factors, so it is not semisimple.
Separability of the field extension and semisimplicity of this linear operator are different properties.

\subsection*{Isomorphisms Between Two Quotient-Ring Models}
Put
\[
g(Y)=f(Y+1)=Y^6+Y^4+Y^2+Y+1,\qquad E'=\F_2[Y]/(g),\quad\beta=\overline Y.
\]
Translation of the variable is invertible, so \(g\) is irreducible. An explicit isomorphism is
\[
P(\alpha)\longmapsto P(\beta+1).
\]
Its inverse takes \(\beta\) to \(\alpha+1\). The other five isomorphisms take \(\alpha\) to \((\beta+1)^{2^j}\) for \(1\le j<6\).
Example~\ref{exm:bm-f64} will use these computations to recover minimal polynomials from trace sequences,
and Example~\ref{exm:f64-quadratic} will use them to connect linearized root spaces with counts of solutions of quadratic equations.

\par\Needspace{11\baselineskip}
\section*{Exercises}
\phantomsection
\addcontentsline{toc}{section}{Exercises}
\begin{exercise}\label{ex:rabin-prime-degree}
When \(n\) is prime, simplify Rabin's criterion to two conditions.
Explain why having no base-field root does not ensure irreducibility in degree four.
\end{exercise}
\begin{exercise}\label{ex:squarefree-example}
Perform square-free factorization on \(f=(X+1)^3(X^2+X+1)^2\) over \(\F_2[X]\), giving \(c,w,y,z\) at each step.
\end{exercise}

\begin{exercise}[Complete factorization from multiplicities to root orbits]\label{v56:ex-factorization-pipeline}
Over \(\F_3[X]\), put
\[
F=X^8+2X^7+X^6+X^4+X^3+2X^2+1.
\]
Determine its complete irreducible factorization, including every multiplicity, and prove completeness.
\end{exercise}
\begin{exercise}\label{ex:berlekamp7}
Let \(f=X^7-1\in\F_2[X]\). Find the dimension of the Berlekamp fixed space and explain its relation to the number of cyclotomic cosets modulo seven.
\end{exercise}

\begin{exercise}[Values and multiplicities in positive characteristic]\label{ex:berlekamp-value-multiplicity}
Over \(\F_2[X]\), put
\[
P=X^2+X+1,\qquad Q=X^4+X+1,\qquad f=PQ,\qquad h=X^4+X\pmod f.
\]
Prove that \(P,Q\) are irreducible and \(h^2=h\). Find the value polynomial \(m_h\) of \(h\) and the characteristic polynomial \(\chi_h\) of multiplication by it,
and compute \(\chi_h/\gcd(\chi_h,\chi_h')\). Explain why this last quotient cannot recover \(m_h\) even when \(f\) is square-free.
Also take \(F=X^2(X+1)\in\F_2[X]\), compute \(F/\gcd(F,F')\), and determine whether it retains every distinct root of \(F\).
How does this relate to the recursive \(p\)th-root step in Subsection~\ref{sec:squarefree-factorization}?
\end{exercise}

\begin{exercise}\label{ex:supp-berlekamp-repeated}
Let \(f=(X^2+X+1)^2\in\F_2[X]\). Find the dimension of the Berlekamp fixed algebra and explain why it alone cannot determine whether this \(f\) is irreducible.
\end{exercise}

\begin{exercise}[Necessary conditions in irreducibility tests]\label{ex:wu-rabin-missing}
Let \(f=(X^2+X+1)(X^3+X+1)\in\F_2[X]\).
Prove that \(f\) has prime degree and no base-field root but is reducible.
Without expanding \(X^{32}\), identify the condition in Rabin's criterion that excludes this polynomial.
\end{exercise}

\begin{exercise}[Idempotents and subsets of root orbits]\label{ex:wu2-idempotents}
Let a monic nonconstant square-free polynomial \(f\in k[X]\) have \(s\) distinct irreducible factors.
Prove that \(k[X]/(f)\) has exactly \(2^s\) idempotents.
Explain their correspondence with arbitrary subsets of root orbits, and give a formula factoring \(f\) from an idempotent.
\end{exercise}

\begin{exercise}[Two coordinate systems for one element]\label{ex:wu3-f64-coordinates}
In Example~\ref{exm:f64-computation}, use the dual basis to find the normal coordinates of \(\alpha\)
and compute \(\Tr_{\F_{64}/\F_4}(\gamma)\). Explain why this relative trace is still normal in \(\F_4/\F_2\).
\end{exercise}
\begin{exercise}[Counting linearly compatible maps and field automorphisms]\label{ex:wu3-isomorphism-certificate}
Let \(E=\F_{16}\) and \(\sigma(x)=x^2\).
Count the invertible \(\F_2\)-linear maps \(A:E\to E\) satisfying both \(A\sigma=\sigma A\) and \(A(1)=1\).
How many preserve multiplication? Compare with Example~\ref{exm:linear-not-field}.
\end{exercise}

\begin{exercise}[The uniqueness boundary for degree-minimal sparse multiples]\label{v97:ex-sparse-uniqueness}
Let \(f\in\F_q[X]\) be monic and nonconstant, with \(f(0)\ne0\).
Let \(\mathcal M_s(f)\) consist of its monic multiples of degree \(L_s(f)\) having at most \(s\) terms.
Prove that the support of any \(P\in\mathcal M_s(f)\) is a minimally dependent set of columns \(X^i\bmod f\): it is dependent, while every proper subset is independent. For a fixed such support, prove uniqueness of the monic relation.
If distinct \(P,Q\) belong to \(\mathcal M_s(f)\), prove
\[
 s+4\le\operatorname{wt}(P)+\operatorname{wt}(Q).
\]
Deduce that \(\mathcal M_2(f)\) and \(\mathcal M_3(f)\) are singletons.
Finally, over \(\F_3\), take \(f=X^4-X^3+X^2+X+1\) and use \((X+1)f\) and \((X+2)f\) to show that uniqueness need not hold for four terms.
\end{exercise}

\begin{exercise}[Scalar extension and the scale of the support]\label{v97:ex-sparse-base-change}
Let \(k=\F_q\), let \(f\in k[X]\) be monic and nonconstant with \(f(0)\ne0\), and let \(s\ge2\).
Temporarily write \(L_{s,k}(f)\) to indicate the coefficient field.
\begin{enumerate}
\item For every extension \(K/k\), prove \(L_{s,K}(f)=L_{s,k}(f)\). The polynomial \(f\) remains fixed; it is not replaced by a root's new minimal polynomial over \(K\).
\item For a positive integer \(a\), prove \(L_s(f(X^a))=aL_s(f)\). Hint: group the multiplier's exponents by residue modulo \(a\).
\item If \(p=\Char k\) and \(u\ge0\), prove \(L_s(f^{p^u})=p^uL_s(f)\). Distinguish taking a polynomial power from applying Frobenius only to its coefficients.
\end{enumerate}
\end{exercise}

\begin{exercise}[Conjugate reciprocity and norm-triggered cancellation]\label{v97:ex-norm-cancellation}
Let \(r\) be a prime power, and let \(F=\sum_{i=0}^d a_iX^i\in\F_{r^2}[X]\) be monic of degree \(d\ge2\), with every \(a_i\ne0\).
Put \(F^{[r]}=\sum_i a_i^rX^i\), and suppose
\[
 F^{[r]}=F^\vee,\qquad F^\vee=a_0^{-1}X^dF(X^{-1}).
\]
For \(t_i=a_i/a_{i-1}\), prove \(t_{d+1-i}=t_i^{-r}\).
If some \(1\le i\le d\) satisfies \(2i\ne d+1\) and
\(\Nm_{\F_{r^2}/\F_r}(a_i)=\Nm_{\F_{r^2}/\F_r}(a_{i-1})\), construct an \(F\)-multiple of degree \(d+1\) with at most \(d\) terms.
Explain why a norm equality at the central position \(2i=d+1\) alone does not ensure two cancellations.
\end{exercise}

\begin{exercise}[Three-dimensional projection and the existence of sparse relations]\label{v97:ex-projection-sparse}
Let \(q\) be odd, and let \(f\in\F_q[X]\) be monic with \(f(0)\ne0\) and \(d=\deg f\ge2\). Prove
\[
 L_2(f)\le q+1\quad(d=2),\qquad
 L_d(f)\le d+q-2\quad(d\ge3).
\]
Hint: quotient \(k[X]/(f)\) by the vector subspace spanned by its first \(d-3\) powers to obtain a three-dimensional space.
If \(q+2\) nonzero, pairwise nonproportional vectors in that space had every three independent, consider all two-dimensional subspaces through a given one-dimensional subspace and derive a parity contradiction.
Finally, take \(\alpha\in\F_{q^2}\) of order \(q+1\), and prove
\[
 f=(X-1)(X-\alpha)(X-\alpha^{-1})\in\F_q[X]
 \quad\Longrightarrow\quad L_3(f)=q+1.
\]
Thus the cubic bound is attained.
\end{exercise}

\begin{exercise}[Exact counts for quadratic reciprocal families]\label{v114:ex-quadratic-count}
Let \(\alpha\in\F_{q^2}\setminus\F_q\). Prove
\[
 b_{q,2}(\alpha)=
 \begin{cases}1,&q\text{ odd},\\0,&q\text{ even}.\end{cases}
\]
Find the unique exceptional parameter in odd characteristic. The even-characteristic assertion does not assume \(2\mid q-1\).
\end{exercise}

\begin{exercise}[Complete classification of cubic reciprocal families]\label{v114:ex-cubic-classification}
Suppose \(3\mid q-1\), \(X^3-d\in\F_q[X]\) is irreducible, and \(\theta^3=d\). Write \(\alpha=a+u\theta+v\theta^2\), where \(a,u,v\in\F_q\) and \((u,v)\ne(0,0)\), and put \(w=t+a\).
\begin{enumerate}
\item Verify that the three coordinate numerators of the reciprocal are
\[
 w^2-uvd,\qquad v^2d-uw,\qquad u^2-vw,
\]
with common denominator \(w^3+u^3d+v^3d^2-3uvdw\).
\item If \(uv=0\), prove that the unique exception satisfies \(w=0\). If \(uv\ne0\), prove that the two linear numerators have distinct roots, neither a root of the quadratic numerator. Deduce
\[
 b_{q,3}(\alpha)=2+\#\{w\in\F_q:w^2=uvd\}\quad(uv\ne0).
\]
List all possible values in odd and even characteristic, and determine \(B_q(3)\).
\end{enumerate}
\end{exercise}

\begin{exercise}[A common constraint on quartic numerators]\label{v114:ex-quartic-bound}
Suppose \(4\mid q-1\), and let \(\alpha\) have full degree in \(\F_{q^4}/\F_q\). Prove
\[
 b_{q,4}(\alpha)\le\min\{q,7\}.
\]
Hint: choose \(\theta^q=\zeta\theta\), with \(\zeta\) a primitive fourth root of unity. Remove the constant coordinate and write \(\alpha=a\theta+b\theta^2+c\theta^3\), with \(d=\theta^4\). Prove that \(d\) is a nonsquare and \((a,c)\ne(0,0)\).
Compute the four coordinate numerators of the reciprocal. When \(ac\ne0\), compare the discriminants of the \(\theta\)- and \(\theta^3\)-numerators; also handle the degree drops when \(a=0\) or \(c=0\).
Finally, in \(\F_{17}(\theta)\), with \(\theta^4=3\), take \(\alpha=\theta+3\theta^2-4\theta^3\) and verify that \(b_{17,4}(\alpha)=7\).
\end{exercise}

\begin{exercise}[Average failure count over all generators]\label{v114:ex-reciprocal-average}
Suppose \(n\ge2\), \(n\mid q-1\), and put
\[
 G_{q,n}=|\mathcal G_{q,n}|=\sum_{d\mid n}\mu(d)q^{n/d}.
\]
Use double counting to prove
\[
 \frac1{G_{q,n}}\sum_{\alpha\in\mathcal G_{q,n}}b_{q,n}(\alpha)
 =q\left(1-\frac{(q-1)^n}{G_{q,n}}\right).
\]
Explain why every normal element belongs to \(\mathcal G_{q,n}\), and why, for fixed \(t\), the map \(\alpha\mapsto(\alpha+t)^{-1}\) is a bijection on this set.
Compute the averages for \((q,n)=(7,3)\) and \((5,4)\), and compare them with Examples~\ref{v114:ex-cubic-comparison} and~\ref{v114:ex-quartic-failure}.
\end{exercise}

\begin{exercise}[An operator polynomial that does not split]\label{v115:ex-reciprocal-eight}
Let \(E=\F_2(\alpha)\), where \(\alpha^3+\alpha+1=0\). Use
\(Y^3-1=(Y+1)(Y^2+Y+1)\) to determine all normal elements among \((\alpha+t)^{-1}\), \(t\in\F_2\).
Describe what each local test detects on this family, and check the result using explicit reciprocals.
\end{exercise}

\subsection*{Extremal Exception Counts}
\phantomsection\label{v114:problem-extremal}
For \(n\ge2\) and \(n\mid q-1\), determine or estimate \(B_q(n)\). The preceding exercises give
\[
 B_q(2)=1,\qquad
 B_q(3)=\begin{cases}4,&q\text{ odd},\\3,&q\text{ even},\end{cases}
 \qquad B_q(4)\le\min\{q,7\},
\]
Example~\ref{v114:ex-quartic-failure} gives \(B_5(4)=5\), and Example~\ref{v115:ex-quartic-seven} gives \(B_{17}(4)=7\). The general bound is \(\min\{q,(n-1)^2\}\).
When \(q\ge(n-1)^2\), attaining \((n-1)^2\) requires every numerator to have its maximum allowed degree, to split over the base field with simple roots, and to have a root set disjoint from those of all the other numerators. These conditions are also sufficient.
Which generators can satisfy them simultaneously? Which extension degrees force identities among the numerators that obstruct simultaneous splitting, as in degree four?
The question is how compatibility among the numerators sharpens the bound supplied by a common denominator.

\chapter{Root Fibers, Return Periods, and Irreducible Constructions}
\label{chap:root-fibers}
\phantomsection\label{v113:par-fiber-question}
Suppose that the outer polynomial \(P\) is irreducible. Why can polynomial substitution, or rational substitution followed by clearing denominators, nevertheless produce a reducible polynomial? Irreducibility of the outer polynomial says that the image roots form one Frobenius orbit; it does not yet determine how Frobenius permutes the preimages above each root.
Starting at a root in one fiber and following the image orbit for one full turn guarantees a return to that fiber, but not necessarily to the starting root. We must therefore study the return map induced inside the fiber by that turn. The base orbit length and the internal return period together determine the full orbit length; different internal cycles can yield different irreducible factors.

This chapter formulates that observation through \hyperref[def:equivariant-map]{equivariant maps}, which commute with Frobenius: \(\pi\sigma=\sigma\pi\). Together with the multiplicity data of the fiber polynomials, return maps characterize the irreducibility of composite polynomials and the degrees and multiplicities of their irreducible factors.
Divisibility relations among return periods also compare the fields generated by the roots; for linearized compositions, periods in affine fibers control the least root field.

\section{Root Fibers and Equivariant Return Maps}
\label{sec:root-fiber-language}
Let \(P(Y)\in k[Y]\) be monic and irreducible, with \(\deg P=n\), and let \(\alpha\) be a root.
Take relatively prime \(A,B\in k[X]\), with \(B\ne0\) and \(A/B\) nonconstant, and assume that \(A-\alpha B\) is nonconstant.
In a finite extension \(L/k\) containing all roots of \(P\) and over which \(Q(X)=B(X)^nP(A(X)/B(X))\) splits completely, consider the map
\[
h:\{x\in L:B(x)\ne0\}\longrightarrow L,\qquad h(x)=A(x)/B(x).
\]
\leavevmode\phantomsection\label{def:root-fiber}\ignorespaces By the definition in Section~\ref{sec:fibers}, the fiber over \(\alpha^{q^i}\) is
\[
F_{\alpha^{q^i}}=h^{-1}(\alpha^{q^i})
=\{x\in L:A(x)-\alpha^{q^i}B(x)=0\}.
\]
The second equality holds because coprimality of \(A,B\) prevents a root of the equation from making the denominator zero.
These fibers lie over the roots of \(P\) and are called \textbf{root fibers}.
The corresponding polynomial \(A-\alpha^{q^i}B\) is called the fiber polynomial.
A root fiber is a set and records only distinct roots; its fiber polynomial also records their multiplicities.
Below, calling a fiber separable or irreducible means that its fiber polynomial has that property; factor multiplicities within a fiber are likewise recorded by that polynomial.

\leavevmode\phantomsection\label{def:fiber-return-map}\ignorespaces Put \(\mathcal O=\{\alpha,\alpha^q,\ldots,\alpha^{q^{n-1}}\}\), and let \(S\) be the set of distinct roots of \(Q\).
The restriction \(\pi=h|_S:S\to\mathcal O\) partitions \(S\) into these root fibers and satisfies \(\pi\sigma=\sigma\pi\).
Thus Frobenius maps each fiber bijectively onto the next, and its \(n\)th iterate is the first to preserve the initial fiber.
The permutation obtained by restricting \(\sigma^n\) to this fiber is called the \textbf{return map}:
it records how the original action moves elements within the same fiber after one circuit of the base orbit.
The entire root set is therefore a single orbit precisely when the return map acts transitively on a fiber.
The change of step shown in Figure~\ref{fig:frobenius-step} is here restricted to a fiber:
\(\sigma\) moves successively through the fibers, whereas \(\sigma^n\) preserves the chosen fiber.

Theorem~\ref{v44:thm-cyclotomic-affine-return} gave a concrete form of this mechanism in fibers arising from modulus lifting; we now formulate it as an orbit correspondence for general finite sets.

\Needspace{5\baselineskip}

\par\Needspace{7\baselineskip}
\noindent\phantomsection\label{v102:par-two-fiber-returns}\textbf{One base orbit, two return maps.}\quad
Take \(k=\F_3\), \(P(Y)=Y^2+1\), and \(i\in\F_9\) with \(i^2=-1\).
The base orbit is \(i\mapsto-i\mapsto i\), of length two. Compare \(h_0(X)=X^2\) with \(h_1(X)=X^2+1\).
In both cases each fiber has two distinct points, and the return map is \(\sigma^2:x\mapsto x^9\).

For \(h_0\), put \(a=1-i\). Then \(a^2=i\), so \(F_i=\{a,-a\}\subseteq\F_9\), and the return map fixes both points.
With \(b=a^3\), the other fiber is \(F_{-i}=\{b,-b\}\); the full orbits are \(\{a,b\}\) and \(\{-a,-b\}\).
Direct multiplication verifies
\[
P(h_0(X))=X^4+1=(X^2+X+2)(X^2-X+2).
\]
Neither quadratic has a root in \(\F_3\), so both are irreducible.

For \(h_1\), take \(a^2=i-1\) and again put \(b=a^3\). Since
\[
(i-1)^2=i,\qquad (i-1)^4=-1,\qquad
a^9=a(a^2)^4=-a,
\]
the return map interchanges \(a,-a\). All four roots, in order \(a,b,-a,-b\), form a single Frobenius orbit.
Thus \(P(h_1(X))=X^4+2X^2+2\) is irreducible. The symbols \(a,b\) are chosen separately for each \(h_j\); they do not denote the same pair in the two panels.

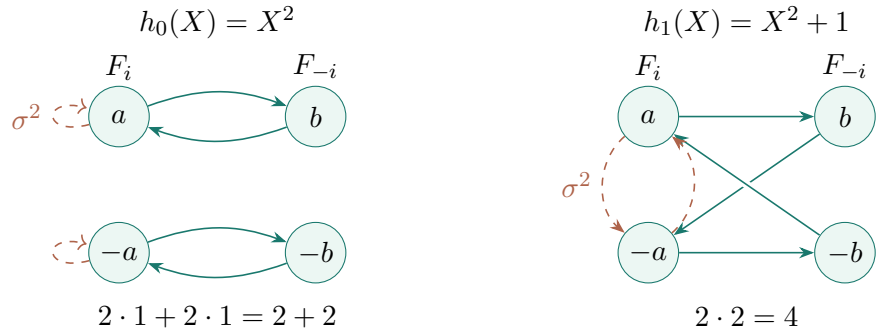
\begin{figure}[!htbp]
\centering
\begin{tikzpicture}[
  vpoint/.style={circle,draw=Jade,fill=CoverMint,minimum size=8mm,inner sep=2pt},
  varrow/.style={-{Stealth[length=2mm]},Jade,semithick},
  vbox/.style={draw=MutedSlate,rounded corners=2pt,inner sep=5pt,align=center},
  every node/.style={font=\small}]

\begin{scope}
\node at (1.3,1.25) {\(h_0(X)=X^2\)};
\node at (0,.65) {\(F_i\)};\node at (2.6,.65) {\(F_{-i}\)};
\node[vpoint] (a) at (0,0) {\(a\)};\node[vpoint] (b) at (2.6,0) {\(b\)};
\node[vpoint] (aa) at (0,-1.8) {\(-a\)};\node[vpoint] (bb) at (2.6,-1.8) {\(-b\)};
\draw[varrow] (a) to[bend left=20] (b);\draw[varrow] (b) to[bend left=20] (a);
\draw[varrow] (aa) to[bend left=20] (bb);\draw[varrow] (bb) to[bend left=20] (aa);
\draw[varrow,Clay,dashed] (a) edge[loop left] node[left] {\(\sigma^2\)} (a);
\draw[varrow,Clay,dashed] (aa) edge[loop left] (aa);
\node at (1.3,-2.65) {\(2\cdot1+2\cdot1=2+2\)};
\end{scope}
\begin{scope}[xshift=7cm]
\node at (1.3,1.25) {\(h_1(X)=X^2+1\)};
\node at (0,.65) {\(F_i\)};\node at (2.6,.65) {\(F_{-i}\)};
\node[vpoint] (c) at (0,0) {\(a\)};\node[vpoint] (d) at (2.6,0) {\(b\)};
\node[vpoint] (cc) at (0,-1.8) {\(-a\)};\node[vpoint] (dd) at (2.6,-1.8) {\(-b\)};
\draw[varrow] (c)--(d);\draw[varrow] (cc)--(dd);
\draw[varrow] (d)--(cc);
\draw[varrow,preaction={draw=white,line width=3pt}] (dd)--(c);
\draw[varrow,Clay,dashed] (c) to[bend right=50] node[left] {\(\sigma^2\)} (cc);
\draw[varrow,Clay,dashed] (cc) to[bend right=50] (c);
\node at (1.3,-2.65) {\(2\cdot2=4\)};
\end{scope}

\end{tikzpicture}
\caption{Every solid arrow is one Frobenius step; dashed arrows show return on the initial fiber only. There are two two-cycles on the left and one four-cycle on the right. Crossings are not vertices. The numbers below are irreducible factor degrees.}\label{v102:fig-two-fiber-returns}
\end{figure}
\FloatBarrier
The base orbit length and fiber size are identical in the two examples; the return permutation determines their different factorization types.
In general, a base orbit of length \(n\) and return cycles of lengths \(d_1,\ldots,d_s\) give full orbits of distinct roots of lengths \(nd_1,\ldots,nd_s\).
Factor multiplicities still require the root multiplicities in the fiber polynomial. The following lemma proves the orbit-length correspondence.

\begin{lemma}[Return map of an orbit fiber]\label{lem:fiber-return}
Let each of the finite sets \(S\) and \(O\) carry a permutation denoted by \(\sigma\), and suppose that \(O\) is a single orbit of length \(n\).
If a surjection \(\pi:S\to O\) satisfies \(\pi\sigma=\sigma\pi\), then, for every \(a\in O\),
the \(\sigma\)-orbits in \(S\) correspond bijectively to the \(\sigma^n\)-orbits in the fiber \(\pi^{-1}(a)\).
A return orbit of length \(d\) corresponds to a full orbit of length \(nd\).
\end{lemma}
\begin{proof}
The image of each full orbit is a nonempty \(\sigma\)-invariant subset of \(O\), hence is all of \(O\); thus the orbit meets the chosen fiber.
If \(x\) lies in this fiber, then \(\sigma^t(x)\) lies in it if and only if \(\sigma^t(a)=a\),
that is, \(n\mid t\). The intersection of the full orbit with the fiber is therefore exactly one \(\sigma^n\)-orbit.
If the latter has least period \(d\), the return times of \(x\) to itself are precisely the multiples of \(nd\), and its least positive period is \(nd\).
Distinct full orbits are disjoint, so this correspondence is bijective.
\end{proof}

\begin{theorem}[Recovering the full factorization from one fiber]\leavevmode\phantomsection\label{thm:fiber-factorization}
Let \(P\in k[Y]\) be monic irreducible of degree \(n>1\), and let \(A,B\in k[X]\) be coprime with \(A/B\) nonconstant. Choose a root \(\alpha\) of \(P\), and put \(E=k(\alpha)\) and
\[
Q(X)=B(X)^nP\bigl(A(X)/B(X)\bigr).
\]
Then \(\deg Q=n\max(\deg A,\deg B)\).
Normalize \(A-\alpha B\) to be monic and factor it in \(E[X]\) as
\[
\frac{A-\alpha B}{\operatorname{lc}(A-\alpha B)}
=\prod_{j=1}^s h_j^{e_j},
\]
where the \(h_j\) are distinct monic irreducible polynomials. Apply \(\sigma(a)=a^q\) to coefficients and define
\[
H_j(X)=\prod_{i=0}^{n-1}\sigma^i(h_j)(X).
\]
Then the \(H_j\in k[X]\) are pairwise distinct, monic, and irreducible, and
\[
\deg H_j=n\deg h_j,\qquad
\frac{Q}{\operatorname{lc}(Q)}=\prod_{j=1}^sH_j^{e_j}.
\]
\end{theorem}
\begin{proof}
Since \(\alpha\notin k\), each conjugate fiber \(A-\alpha^{q^i}B\) has degree \(h=\max(\deg A,\deg B)\): cancellation of leading terms of equal degree would imply \(\alpha^{q^i}\in k\). The identity
\[
Q=\prod_{i=0}^{n-1}(A-\alpha^{q^i}B)
\]
gives the degree formula.
In a finite splitting extension, take the set of distinct roots of \(Q\), and put \(\pi(\beta)=A(\beta)/B(\beta)\).
Coprimality ensures that the denominator is nonzero, and the fibers over distinct \(\alpha^{q^i}\) are disjoint.
The map \(\pi\) commutes with Frobenius, and its image is the full orbit of \(\alpha\).
Within the initial fiber, the roots of \(h_j\) form a single orbit relative to \(E\), with return map \(\sigma^n\).
By Lemma~\ref{lem:fiber-return}, their orbit in the full root set has length \(n\deg h_j\).
The root polynomial obtained by applying Frobenius successively to this root set through the first \(n\) fibers is exactly \(H_j\), so it is the minimal polynomial over the base field.
Distinct return orbits correspond to distinct full orbits, so the \(H_j\) are distinct.

Applying Frobenius to coefficients preserves factor multiplicities. Multiplying the monic factorizations of the conjugate fibers gives the final identity.
Even when \(Q\) has repeated roots, this argument determines both the orbit structure and all factor multiplicities.
\end{proof}

Irreducibility is the case of a single root orbit with multiplicity one.
\begin{corollary}[Root fiber criterion, or Capelli's criterion]\label{thm:capelli}
Let \(P\in k[Y]\) be monic and irreducible with \(\deg P=n\), and let \(A,B\in k[X]\) be relatively prime with \(A/B\) nonconstant.
Take a root \(\alpha\) of \(P\), and assume that \(A-\alpha B\) is nonconstant. Then
\[
Q(X)=B(X)^nP\bigl(A(X)/B(X)\bigr)
\]
is irreducible over \(k\) if and only if \(A(X)-\alpha B(X)\) is irreducible over \(k(\alpha)\).
If \(n>1\), then \(\deg Q=n\max(\deg A,\deg B)\).
\end{corollary}
\begin{proof}
For \(n>1\), Theorem~\ref{thm:fiber-factorization} makes \(Q\) irreducible precisely when the fiber has a single irreducible factor of multiplicity one, which is exactly the irreducibility of \(A-\alpha B\). The same theorem gives the degree formula.
For \(n=1\), we have \(P(Y)=Y-\alpha\), \(\alpha\in k\), and \(Q=A-\alpha B\), so the assertion is immediate.
\end{proof}

\leavevmode\phantomsection\label{v98:par-fiber-construction}\ignorespaces
For a construction, first choose an irreducible \(P\), and then choose \(A,B\) satisfying the hypotheses of the theorem. Their base-field coefficients make the fibers over conjugate roots compatible with Frobenius; these fibers cannot be prescribed independently.
A desired factor degree is the product of the base-orbit length and the \emph{least} return period; knowing that one iterate returns gives only a multiple of the period.
For the entire composite polynomial to be irreducible, the fiber polynomial must also have no repeated roots, and the return map must traverse the whole fiber.
In a recursive construction, continuing to the next step requires the data enforcing these conditions to be preserved; Chapter~\ref{chap:root-constructions} controls them through traces, square classes, and multiplicative orders.
The same return lemma also applies to linear state spaces and their quotients; Section~\ref{sec:recurrence-algebra} gives the corresponding root-space realization.

If only factor degrees and counts are required, count fixed points of iterates of the return map before separating individual factors in the fiber.
\begin{corollary}[Return fixed points and factor degrees]\label{v115:cor-return-counts}
In the notation of Theorem~\ref{thm:fiber-factorization}, let \(H\in E[X]\) be the monic fiber polynomial, of degree \(h\), where \(|E|=q^n\). For \(1\le e\le h\), put
\[
r_e=\deg\gcd\bigl(H,X^{q^{ne}}-X\bigr).
\]
The number of distinct monic irreducible factors of the full polynomial \(Q\) having degree \(nj\) is
\begin{equation}\label{v115:eq-return-counts}
a_j=\frac1j\sum_{e\mid j}\mu(j/e)r_e\qquad(1\le j\le h).
\end{equation}
If \(H\) is squarefree, these numbers give the complete factorization type of \(Q\). Otherwise they count distinct factors without their multiplicities.
\end{corollary}
\begin{proof}
On the distinct roots of \(H\), the return map is \(\tau=\sigma^n\). Since \(X^{q^{ne}}-X\) has simple roots, its greatest common divisor with \(H\) retains exactly the fixed points of \(\tau^e\).
If \(c_j\) counts return cycles of length \(j\), then
\[
r_e=\sum_{j\mid e}jc_j.
\]
M\"obius inversion gives~\eqref{v115:eq-return-counts}. The return lemma sends each cycle of length \(j\) to a base-field orbit of length \(nj\), so \(a_j=c_j\).
Multiplicities do not enter the permutation of distinct roots. Nothing is lost in the squarefree case; otherwise the fiber multiplicities are still required.
\end{proof}

The two-point fibers in Figure~\ref{v102:fig-two-fiber-returns} have return fixed-point counts \((r_1,r_2)=(2,2)\) and \((0,2)\), respectively.
These give \(a_1=2\) and \(a_2=1\), hence full factor degrees \(2,2\) and \(4\).
Compute \(r_e\) by successive modular powers \(X^{q^{ne}}\bmod H\). To obtain the factors themselves, return to the fiber factorization theorem: fixed-point counts discard the coordinates distinguishing factors of the same degree.

\begin{corollary}[Fiber returns after a change of base field]\label{v129:cor-return-base-change}
Retain the notation of Lemma~\ref{lem:fiber-return}, and put \(F=\pi^{-1}(a)\) and \(\tau=\sigma^n|_F\).
For an integer \(r\ge1\), put \(g=\gcd(n,r)\) and \(s=r/g\). Under \(\sigma^r\), the base orbit splits into \(g\) orbits of length \(n/g\), and the first return on \(F\) is \(\tau^s\).
An original return cycle of length \(j\) gives, above each new base orbit, \(\gcd(j,s)\) total cycles, each of length
\[
\frac{n}{g}\,\frac{j}{\gcd(j,s)}.
\]
Thus an original total cycle of length \(nj\) splits into \(\gcd(nj,r)\) cycles altogether.
\end{corollary}
\begin{proof}
The least period of \(\sigma^r\) on the base orbit is \(n/g\), so its first return on \(F\) is
\[
(\sigma^r)^{n/g}|_F=(\sigma^n)^{r/g}|_F=\tau^s.
\]
The \(s\)-th power of a \(j\)-cycle has \(\gcd(j,s)\) cycles of length \(j/\gcd(j,s)\); the return lemma gives the asserted total lengths.
The fibers \(\sigma^i(F)\), \(0\le i<g\), may be chosen over the other new base orbits, and \(\sigma^i\) conjugates their return permutations.
Finally, \(\gcd(n/g,s)=1\), whence
\[
\gcd(nj,r)=g\gcd((n/g)j,s)=g\gcd(j,s).
\]
\end{proof}

For root fibers, this describes factorization after extending the base field from \(\F_q\) to \(\F_{q^r}\). Splitting of the base orbit and splitting of the internal return must be computed separately; irreducible-factor multiplicities are unchanged.
For example, a four-cycle with \(n=2,j=2\) splits into two two-cycles when \(r=2\): the base orbit becomes two fixed points, but \(s=1\), so the return in each fiber remains a two-cycle.

For inverse images, return cycles produce longer root orbits; for images, the same relation records repeated conjugate values.
Let \(S\) be a finite Frobenius-stable root set, and let \(R\in k(X)\) be defined on \(S\). For each root orbit \(\mathcal O\subseteq S\), choose \(x_{\mathcal O}\in\mathcal O\), and put
\(n_{\mathcal O}=|\mathcal O|\) and \(d_{\mathcal O}=[k(R(x_{\mathcal O})):k]\). Then
\begin{equation}\label{v129:eq-orbit-image-product}
\prod_{x\in S}(Y-R(x))
=\prod_{\mathcal O\subseteq S}
m_{R(x_{\mathcal O})}(Y)^{n_{\mathcal O}/d_{\mathcal O}}.
\end{equation}
Each original orbit maps onto an image orbit with equal fiber sizes, proving the formula; Proposition~\ref{v79:prop-generator-orbit-image} gives the resultant expression for a single orbit.
A single orbit may shorten, and different orbits may have the same image orbit; in the latter case, exponents of equal factors must be added.
Thus distinct roots in \(S\) can give repeated factors in the image product, which is squarefree exactly when \(R|_S\) is injective.

\begin{corollary}[Fractional linear transformations]\label{cor:mobius-transform}
Let \(P\in k[X]\) be irreducible with \(\deg P=n\ge2\), and let \(a,b,c,d\in k\) satisfy \(ad-bc\ne0\). Then
\[
(cX+d)^nP\left(\frac{aX+b}{cX+d}\right)
\]
is again an irreducible polynomial of degree \(n\).
\end{corollary}
\begin{proof}
Denote the displayed polynomial by \(Q\), and take a root \(\alpha\) of \(P\). Since \(n\ge2\), we have \(\alpha\notin k\).
If \(c\ne0\), then \(a-c\alpha\ne0\); if \(c=0\), the determinant condition gives \(a\ne0\), so this denominator is again nonzero.
We may therefore define
\[
\beta=\frac{d\alpha-b}{a-c\alpha}.
\]
Direct calculation gives
\[
c\beta+d=\frac{ad-bc}{a-c\alpha}\ne0,\qquad
\alpha=\frac{a\beta+b}{c\beta+d}.
\]
Thus \(Q(\beta)=0\), and the two expressions give \(\beta\in k(\alpha)\) and \(\alpha\in k(\beta)\), respectively,
so \(k(\alpha)=k(\beta)\) and \([k(\beta):k]=n\).

The condition \(ad-bc\ne0\) also ensures that \(aX+b\) and \(cX+d\) are relatively prime and their ratio is nonconstant.
After making \(P\) monic, apply the degree formula of Corollary~\ref{thm:capelli} to obtain \(\deg Q=n\).
Hence the minimal polynomial of \(\beta\) divides \(Q\) and has the same degree, so \(Q\) is irreducible.
All coefficients of the root transformation lie in \(k\), so it commutes with Frobenius and maps the single root orbit of \(P\) bijectively onto that of \(Q\).
\end{proof}

\subsection*{Eliminating the Trace while Preserving Generation}
The preceding composite polynomial takes preimages under a rational function: its new root \(\beta\) satisfies \(R(\beta)=\alpha\).
We now take the image \(\beta=R(\alpha)\) of a given generator and ask when a trace condition can be imposed without shortening its Frobenius orbit.
Here a generator means a field generator, with no additional requirement of primitivity or normality.

\begin{lemma}[Low-degree changes and field generation]\leavevmode\phantomsection\label{v70:lem-low-degree-generation}
Let \(k=\F_q\), \(E=k(\alpha)\), and \([E:k]=n\ge2\), and let \(\ell\) be the least prime divisor of \(n\).
Let \(R=A/B\in k(X)\) be nonconstant, where \(A,B\) are relatively prime and \(B(\alpha)\ne0\), and put
\(\deg R=\max\{\deg A,\deg B\}\). If \(\beta=R(\alpha)\), then
\[
[E:k(\beta)]\le\deg R,\qquad [E:k(\beta)]\mid n.
\]
In particular, if \(\deg R<\ell\), then \(k(\beta)=E\).
\end{lemma}
\begin{proof}
The polynomial \(A(X)-\beta B(X)\in k(\beta)[X]\) is nonzero.
Otherwise, comparison at any nonzero coefficient of \(B\) would give \(\beta\in k\) and \(A=\beta B\), contrary to the nonconstancy of \(R\).
Since it vanishes at \(\alpha\), the minimal polynomial of \(\alpha\) over \(k(\beta)\) divides it, giving
\[
[E:k(\beta)]\le\deg(A-\beta B)\le\deg R.
\]
The tower formula gives \([E:k(\beta)]\mid n\). Every divisor of \(n\) greater than one is at least \(\ell\), so the strict inequality \(\deg R<\ell\) forces this relative degree to be one.
\end{proof}

The strict inequality cannot be replaced by \(\deg R\le\ell\). For example, if \(\alpha^4+\alpha+1=0\) with \(\alpha\in\F_{16}\), then
\(\beta=\alpha^2+\alpha\) satisfies \(\beta^2+\beta+1=0\) and generates only \(\F_4\);
this is also the quadratic-subfield representation in Example~\ref{sfr:ex-quartic-comparison}.
Thus a quadratic change in an extension of even degree requires a separate generation argument.

\begin{proposition}[A trace criterion for fractional linear changes]\leavevmode\phantomsection\label{v70:prop-mobius-trace-zero}
Let \(k=\F_q\) have characteristic \(p\), let \(E=k(\alpha)\) with \([E:k]=n\ge2\), and put
\(f=m_\alpha\) and \(T=\Tr_{E/k}(\alpha)\). If \(p\nmid n\), then \(\alpha-T/n\) is a generator of trace zero.
If \(p\mid n\), there is a fractional linear transformation \(R\) over \(k\) such that \(R(\alpha)\) still generates \(E\) and has trace zero if and only if
\[
T=0\quad\text{or}\quad f'(t)=0\text{ for some }t\in k.
\]
\end{proposition}
\begin{proof}
A fractional linear transformation and its inverse are defined over \(k\), so it preserves field generation; its denominator does not vanish at \(\alpha\), since \(\alpha\notin k\).
If \(p\nmid n\), then \(\Tr_{E/k}(T/n)=T\), and translation gives the required element.
Now suppose \(p\mid n\). An affine change \(u\alpha+v\), with \(u\in k^\times\) and \(v\in k\), has trace \(uT\), so such a change works exactly when \(T=0\).
Every nonaffine fractional linear function has the form \(u+v/(X-t)\), where \(u,t\in k\) and \(v\in k^\times\).
The logarithmic derivative of \(f(X)=\prod_{i=0}^{n-1}(X-\alpha^{q^i})\) gives
\begin{equation}\label{v70:eq-reciprocal-trace}
\Tr_{E/k}\left(\frac1{\alpha-t}\right)
=\sum_{i=0}^{n-1}\frac1{\alpha^{q^i}-t}
=-\frac{f'(t)}{f(t)}.
\end{equation}
Here \(f(t)\ne0\). Since \(\Tr_{E/k}(u)=nu=0\), the trace of this change is zero exactly when \(f'(t)=0\), proving both necessity and sufficiency.
\end{proof}

\begin{proposition}[A trace-zero generator of rational degree at most two]\leavevmode\phantomsection\label{v70:prop-trace-zero-generator}
Let \(k=\F_q\), \(E=k(\alpha)\), and \([E:k]=n\ge3\). There is a nonconstant rational function \(R\in k(X)\),
whose reduced numerator and denominator both have degree at most two and whose denominator does not vanish at \(\alpha\), such that \(\beta=R(\alpha)\) satisfies
\[
k(\beta)=E,\qquad \Tr_{E/k}(\beta)=0.
\]
Such a \(\beta\) is called a \textbf{trace-zero generator}.
\end{proposition}
\begin{proof}
Put \(p=\Char k\), \(T=\Tr_{E/k}(\alpha)\), and \(S=\Tr_{E/k}(\alpha^2)\).
If \(p\nmid n\), take \(\beta=\alpha-T/n\); if \(T=0\), take \(\beta=\alpha\).
It remains to consider \(p\mid n\) and \(T\ne0\). Put \(c=S/T\) and set
\[
\beta_0=\alpha^2-c\alpha.
\]
Linearity of the trace gives \(\Tr_{E/k}(\beta_0)=S-cT=0\). Since \(\alpha\) satisfies \(X^2-cX-\beta_0=0\),
we have \([E:k(\beta_0)]\le2\). If \(\beta_0\) does not generate \(E\), this relative degree is two, \(n\) is even,
and \(k(\beta_0)=\F_{q^{n/2}}\). The nonidentity automorphism of this quadratic extension sends \(\alpha\) to the other root, so
\[
\alpha^{q^{n/2}}=c-\alpha,\qquad
T=\Tr_{\F_{q^{n/2}}/k}\bigl(\alpha+\alpha^{q^{n/2}}\bigr)=(n/2)c.
\]
If \(p\) is odd, then \(p\mid n\) implies \(p\mid n/2\); if \(p=2\) and \(4\mid n\), again \(n/2=0\) in \(k\).
Both cases contradict \(T\ne0\), so \(\beta_0\) already gives the required generator.

Only the case \(p=2\), \(n\equiv2\pmod4\), with \(\beta_0\) failing to generate \(E\), remains.
Here \(S=T^2\), so \(c=T\). Put \(\gamma=\alpha^{-1}\) and \(U=\Tr_{E/k}(\gamma)\);
we have \(\alpha\ne0\), and \(\gamma\) still generates \(E\). If \(U=0\), take \(\beta=\gamma\).
If \(U\ne0\), set
\[
\beta=\gamma^2+U\gamma=\frac{1+U\alpha}{\alpha^2}.
\]
Since \(\Tr_{E/k}(\gamma^2)=U^2\), this element again has trace zero.
If \(\beta\) also fails to generate \(E\), the relation \(\gamma^2+U\gamma+\beta=0\) again gives \([E:k(\beta)]=2\).
The subfield of degree \(n/2\) in a finite field is unique, so \(k(\beta)=k(\beta_0)=\F_{q^{n/2}}\),
and the same relative Frobenius satisfies \(\gamma^{q^{n/2}}=\gamma+U\). Together with \(\alpha^{q^{n/2}}=\alpha+T\), this gives
\[
U=\frac1\alpha+\frac1{\alpha+T}
=\frac{T}{\alpha(\alpha+T)}=\frac{T}{\beta_0}.
\]
Thus \(\beta_0=T/U\in k\), contradicting \([k(\beta_0):k]=n/2>1\).
The second construction therefore preserves generation. The functions used are a linear polynomial, a quadratic polynomial, \(1/X\), or \((1+UX)/X^2\);
each is nonconstant and satisfies the stated degree and denominator conditions.
\end{proof}

In a quadratic extension of characteristic two, \(k\) lies in the trace kernel, which has dimension one over \(k\); hence the kernel is exactly \(k\), and no trace-zero generator exists.
Thus the hypothesis \(n\ge3\) cannot be replaced unconditionally by \(n\ge2\); a quadratic extension of odd characteristic is already handled by translation.
The zero-trace condition removes the next-to-leading term of the minimal polynomial, while generation ensures that the new equation still defines the original extension.
When \(n\ge3\) and the criterion for a fractional linear change fails, the quadratic construction attains the least possible rational degree.

\section{Binomials and Lifting Multiplicative Periods}
Lemma~\ref{lem:valuation-lift} describes how valuations change when exponents increase.
We apply it to lifting multiplicative orders and thereby determine when the roots of a binomial form a single Frobenius orbit.
For the classical binomial irreducibility criterion, see \cite[Chapter 3, Section 5]{LN}.

\begin{theorem}[Binomial irreducibility criterion]\label{thm:binomial-full}
Let \(a\in\F_q^\times\), \(e=\ord(a)\), and \(t\ge2\). Then \(X^t-a\) is irreducible if and only if
\begin{enumerate}
\item every prime \(\ell\mid t\) divides \(e\) but does not divide \((q-1)/e\);
\item if \(4\mid t\), then \(q\equiv1\pmod4\).
\end{enumerate}
\end{theorem}
\begin{proof}
We first prove necessity. For each prime \(\ell\mid t\), if \(Y^\ell-a\) were reducible,
substituting \(Y=X^{t/\ell}\) into a nontrivial factorization would make \(X^t-a\) reducible.
Thus every \(X^\ell-a\) is irreducible.
By Proposition~\ref{prop:prime-binomial}, this requires \(\ell\mid q-1\) and that \(a\) is not an \(\ell\)th power.
Take a generator \(z\) of \(k^\times\), and write \(a=z^j\).
Since \(\ell\mid q-1\), the element \(a\) is an \(\ell\)th power if and only if \(\ell\mid j\).
Also, \(e=(q-1)/\gcd(q-1,j)\), so the condition of not being an \(\ell\)th power is equivalent to
\[
v_\ell(e)=v_\ell(q-1)>0.
\]
This is precisely condition (1).

If \(4\mid t\), then \(q\) is odd. Suppose \(q\equiv3\pmod4\); condition (1) gives \(v_2(e)=1\).
Take \(\beta^4=a\), and write \(M=\ord(\beta)\).
The identity \(\ord(\beta^4)=M/\gcd(M,4)=e\) shows that \(v_2(M)=3\) and that the odd part of \(M\) equals that of \(e\), so \(M=4e\).
Since \(e\mid q-1\) and \(8\mid q^2-1\), we have \(4e\mid q^2-1\).
Thus \(\beta^{q^2}=\beta\), so \(\beta\) has degree at most two over \(k\), and \(X^4-a\) is reducible.
Substitution of \(X^{t/4}\) now contradicts the original irreducibility. This proves condition (2).

Conversely, suppose that both conditions hold. Condition (1) ensures \(p\nmid t\), so \(X^t-a\) has no repeated roots.
Take a root \(\beta\), and put \(M=\ord(\beta)\). Again,
\[
\frac{M}{\gcd(M,t)}=e.
\]
If \(\ell\mid t\), then \(v_\ell(e)>0\); comparing \(\ell\)-adic valuations in this identity gives
\[
v_\ell(M)-\min\{v_\ell(M),v_\ell(t)\}=v_\ell(e)>0,
\]
hence \(v_\ell(M)=v_\ell(e)+v_\ell(t)\).
If \(\ell\nmid t\), then \(v_\ell(M)=v_\ell(e)\).
Combining all primes gives \(M=et\).

We now determine the Frobenius period of \(\beta\), namely \(\ord_{et}(q)\).
For every odd prime \(\ell\mid t\), and also for \(\ell=2\) when \(q\equiv1\pmod4\),
we have \(v_\ell(e)=v_\ell(q-1)\). The valuation-lifting lemma gives
\[
v_\ell(q^s-1)\ge v_\ell(et)
\quad\Longleftrightarrow\quad v_\ell(s)\ge v_\ell(t).
\]
The only remaining case is \(q\equiv3\pmod4\) and \(2\mid t\).
Condition (2) then requires \(v_2(t)=1\), while \(v_2(e)=1\),
so the 2-primary part requires only \(4\mid q^s-1\), which is equivalent to \(s\) being even.
For primes not dividing \(t\), the corresponding prime powers in \(et\) already divide \(e\mid q-1\) and impose no further condition.
Therefore \(et\mid q^s-1\) if and only if \(t\mid s\), and the least positive \(s\) is exactly \(t\).
Thus \(\beta\) has degree \(t\), its minimal polynomial has the same degree as \(X^t-a\), and the latter is irreducible.
\end{proof}

The same criterion controls return roots of equally spaced skew polynomials. Proposition~\ref{v148:prop-sparse-transfer} converts the exponent spacing into $X^m-\alpha$; Theorem~\ref{v148:thm-sparse-coprime} combines the root-order conditions with fixed-constant fibers to give existence and exact counts.

For a fixed degree, the binomial criterion also determines every admissible constant term.
\leavevmode\phantomsection\label{v77:def-binomial-constants}\ignorespaces For \(n\ge1\), write
\[
\mathcal A_q(n)=\{a\in\F_q:X^n+a\text{ is irreducible over }\F_q\},
\]
the \textbf{constant-term set of irreducible binomials} of degree \(n\).
Since \(X^n+a=X^n-(-a)\), the multiplicative-order and power-residue conditions in the criterion apply to \(-a\), not directly to \(a\).
Exercise~\ref{ex:supp-binomial} gives a parametrization by a primitive element, equivalent power-residue and multiplicative-order tests, and an exact count.

\begin{example}[All constant terms, arbitrarily large degrees, and the exception modulo four]\label{v77:ex-binomial-constant-families}
In \(\F_5\), the nonsquares are exactly \(2,3\), and negation preserves this set.
Since \(5\equiv1\pmod4\), the binomial criterion gives, for every \(r\ge1\),
\[
\mathcal A_5(2^r)=\{2,3\}.
\]
Thus \(X^{2^r}+2\) and \(X^{2^r}+3\) are irreducible for all \(r\ge1\); the latter is \(X^{2^r}-2\).
The degrees can be arbitrarily large while the admissible coefficient set remains unchanged.
Irreducibility concerns ordinary polynomials, so the exponents cannot be reduced using identities of functions on the base field.

In \(\F_7\), the nonzero cubes are \(1,6\), and negation preserves the set of remaining nonzero elements. Hence
\[
\mathcal A_7(3^r)=\{2,3,4,5\}\qquad(r\ge1).
\]
In particular, \(X^{3^r}+5=X^{3^r}-2\) is irreducible, although \(2\) has multiplicative order three rather than six.
Thus \(-a\) need not be primitive: its order need only contain the full prime-power part of \(q-1\) at each prime dividing the degree.

Over \(\F_3\), in contrast,
\[
\mathcal A_3(2)=\{1\},\qquad \mathcal A_3(4)=\varnothing.
\]
These two degrees have the same prime divisor, but degree four activates the restriction modulo four, and \(3\not\equiv1\pmod4\). Indeed,
\[
X^4+1=(X^2+X+2)(X^2+2X+2).
\]
When a prime exponent is increased, one must therefore check whether the exceptional condition \(4\mid n\) becomes relevant for the first time.
\end{example}

\begin{corollary}\leavevmode\phantomsection\label{cor:power-composition}
Let \(P\in\F_q[X]\) be monic and irreducible, with \(\deg P=n\) and \(P(0)\ne0\), let its roots have multiplicative order \(e\), and let \(t\ge2\) be an integer.
Then \(P(X^t)\) is irreducible if and only if both conditions of Theorem~\ref{thm:binomial-full} hold after replacing \(q\) by \(q^n\) and leaving \(e,t\) unchanged.
When these conditions hold, the resulting polynomial has degree \(nt\) and order \(et\).
\end{corollary}
\begin{proof}
Choose a root \(\alpha\) of \(P\), and put \(E=\F_q(\alpha)=\F_{q^n}\). The fiber polynomial over \(\alpha\) is \(X^t-\alpha\). Apply the binomial criterion over \(E\), followed by the root fiber theorem.
The proof of the binomial criterion shows that the multiplicative order of each root is \(et\).
\end{proof}
For fixed admissible \(n,e,t\), this construction gives a bijection between the monic irreducible polynomials of order \(e\) and degree \(n\) and those of order \(et\) and degree \(nt\).
It is injective, and the respective numbers are \(\varphi(e)/n\) and \(\varphi(et)/(nt)\); condition (1) gives \(\varphi(et)=t\varphi(e)\), so it is also surjective.

The criterion above decides irreducibility of the binomial itself. The norm argument of Proposition~\ref{v97:prop-norm-factor-degrees} still constrains every factor degree when the binomial is reducible. Further restrictions on the splitting field from coefficient support appear in Proposition~\ref{v97:prop-support-splitting}.

\section{Quadratic Fibers and Self-Reciprocal Constructions}
The map \(x\mapsto x+x^{-1}\) sends reciprocal roots \(x\) and \(x^{-1}\) to the same point.
For an irreducible polynomial \(P\) of degree \(n\) with root \(\alpha\), the fiber polynomial over \(\alpha\) is \(X^2-\alpha X+1\).
When this quadratic has no repeated root, \(\sigma^n\) either fixes both roots or interchanges them.
By Lemma~\ref{lem:fiber-return}, the first case yields two irreducible factors of degree \(n\), and the second yields one irreducible factor of degree \(2n\).
Thus the trace and square-class conditions below determine exactly when the return map interchanges a reciprocal pair; self-reciprocity follows from invariance of the full root set under inversion.

\begin{theorem}\label{thm:reciprocal-transform}
Let \(P(X)=X^n+c_{n-1}X^{n-1}+\cdots+c_0\) be irreducible with \(c_0\ne0\). Define its \textbf{\(Q\)-transform} by
\[
\mathcal Q(P)(X)=X^nP(X+X^{-1}).
\]
It is a monic self-reciprocal polynomial of degree \(2n\), and:
\begin{enumerate}
\item in characteristic two, it is irreducible if and only if \(\Tr_{k/\F_2}(c_1/c_0)=1\), where \(c_1=1\) when \(n=1\);
\item in odd characteristic, it is irreducible if and only if \(P(2)P(-2)\) is a nonsquare in \(k\).
\end{enumerate}
\end{theorem}
\begin{proof}
The fiber polynomial over the root \(\alpha\) is \(X^2-\alpha X+1\).
In characteristic two, it is irreducible if and only if \(\Tr_{E/\F_2}(\alpha^{-2})=1\), where \(E=k(\alpha)\).
The absolute trace is invariant under squaring, and \(\Tr_{E/k}(\alpha^{-1})=-c_1/c_0\); transitivity of the trace gives (1).
In odd characteristic, the condition is that \(\alpha^2-4\) be a nonsquare in \(E\).
An element \(z\in E^\times\) is a square if and only if \(\Nm_{E/k}(z)\) is a square in \(k\), since \(\Nm(z)^{(q-1)/2}=z^{(q^n-1)/2}\). Also,
\[
\Nm_{E/k}(\alpha^2-4)=P(2)P(-2),
\]
which gives (2). If this product is zero, the fiber polynomial has a repeated root and again cannot be irreducible.
Self-reciprocity follows by substituting \(X^{-1}\).
\end{proof}
\leavevmode\phantomsection\label{def:r-transform}\ignorespaces In odd characteristic, one also often uses the normalized \textbf{\(R\)-transform}
\[
\mathcal R(P)(X)=(2X)^nP\left(\frac{X+X^{-1}}2\right).
\]
Replacing the parameter \(\alpha\) in the root fiber by \(2\alpha\), the same criterion shows that \(\mathcal R(P)\) is irreducible if and only if \(P(1)P(-1)\) is a nonsquare.

\section{Additive Fibers and Trace Conditions}
\begin{proposition}[Irreducibility criterion for a general additive composition]\leavevmode\phantomsection\label{prop:general-additive-composition}
Let \(k=\F_q\) have characteristic \(p\), let \(P\in k[X]\) be monic and irreducible with \(\deg P=n\),
let \(\alpha\) be a root of \(P\), put \(E=k(\alpha)\), and take \(a\in k^\times\).
Then \(P(X^p-aX)\) is irreducible if and only if there exists \(u\in E^\times\) such that
\[
u^{p-1}=a,\qquad \Tr_{E/\F_p}(\alpha/u^p)\ne0.
\]
\end{proposition}
\begin{proof}
By Corollary~\ref{thm:capelli}, \(P(X^p-aX)\) is irreducible if and only if \(X^p-aX-\alpha\) is irreducible over \(E\).
Apply Theorem~\ref{thm:affine-prime} over the base field \(E\), with parameters \(a,\alpha\), to obtain the stated conditions.
\end{proof}

The power equation \(u^{p-1}=a\) has a solution in \(E^\times\) if and only if \(a^{(q^n-1)/(p-1)}=1\). This is the power-equation criterion in the cyclic group \(E^\times\); it does not yet include the trace condition.
Since \(a\in k^\times\), we have \(a^{q-1}=1\). Write \(M=(q-1)/(p-1)\). Then
\[
\gcd\!\left(q-1,\frac{q^n-1}{p-1}\right)
=M\gcd\bigl(p-1,1+q+\cdots+q^{n-1}\bigr)
=M\gcd(p-1,n).
\]
Thus solvability of the power equation is also equivalent to \(a^{M\gcd(n,p-1)}=1\).

When the power equation is solvable, any two solutions satisfy \(u'=cu\), where \(c\in\F_p^\times\). Hence
\[
\Tr_{E/\F_p}\!\left(\frac{\alpha}{(u')^p}\right)
=c^{-1}\Tr_{E/\F_p}\!\left(\frac{\alpha}{u^p}\right).
\]
Whether the trace is nonzero is therefore independent of the chosen solution. Write \(P(X)=X^n+c_{n-1}X^{n-1}+\cdots+c_0\). If \(u\) can be chosen in \(k\), then \(\Tr_{E/k}(\alpha)=-c_{n-1}\), and the second condition becomes \(\Tr_{k/\F_p}(c_{n-1}/u^p)\ne0\).

\begin{corollary}\label{cor:artin-schreier-composition}
Let \(P=X^n+c_{n-1}X^{n-1}+\cdots+c_0\) be monic and irreducible, and let \(b\in k\). Then
\[
P(X^p-X-b)\text{ is irreducible}\quad\Longleftrightarrow\quad
\Tr_{k/\F_p}(nb-c_{n-1})\ne0.
\]
\end{corollary}
\begin{proof}
Put \(Q(Y)=P(Y-b)\). This is monic irreducible with root \(\alpha+b\). Apply Proposition~\ref{prop:general-additive-composition} to \(Q(X^p-X)\), taking \(u=1\). The required condition is \(\Tr_{E/\F_p}(\alpha+b)\ne0\), where \(E=k(\alpha)\). Now
\[
\Tr_{E/k}(\alpha+b)=-c_{n-1}+nb,
\]
and transitivity of trace gives the result.
\end{proof}

\begin{example}[A solvable power equation with a failed trace condition]\label{exm:affine-trace-obstruction}
Reuse the irreducible polynomial \(P=X^3+X+1\) and its trace-zero root \(\alpha\in\F_8\) from Example~\ref{ex:normal-f8-example}, and take \(a=1\).
The power equation allows \(u=1\), but the required nonzero trace condition fails, so this composition is reducible. Explicitly,
\[
P(X^2+X)=(X^3+X+1)(X^3+X^2+1).
\]
This factorization shows that solvability of the power equation cannot replace the nonzero trace condition required for irreducibility of the fiber polynomial.
\end{example}

\Needspace{5\baselineskip}
\begin{proposition}[Primitive Artin--Schreier polynomials over the prime field]\label{prop:primitive-as}
Let \(b\in\F_p^\times\), and put \(S=(p^p-1)/(p-1)\). Then \(X^p-X-b\) is primitive if and only if \(b\) generates \(\F_p^\times\) and \(\ord(X^p-X-1)=S\).
\end{proposition}
\begin{proof}
If \(\beta^p-\beta=1\), then \(\alpha=b\beta\) is a root of \(X^p-X-b\).
We have \(\Nm(\beta)=1\), so \(\ord(\beta)\mid S\). Also, \(\gcd(S,p-1)=1\),
hence \(\ord(b\beta)=\ord(b)\ord(\beta)\), which reaches \(p^p-1\) if and only if both factors reach their upper bounds.
\end{proof}
The condition \(\ord(X^p-X-1)=S\) is part of the criterion and must be checked separately.

\section{Restrictions on Irreducibility of Affine Root Spaces}\label{sec:affine-root-space}
For an \hyperref[def:affine-space]{affine subspace} \(\beta+V\), its dimension is defined as \(\dim V\).
After choosing a base point \(\beta\), the correspondence \(v\leftrightarrow\beta+v\) gives linear coordinates, but the affine subspace itself has no distinguished zero vector.
In these coordinates, a map of the form \(v\mapsto Av+c\), with \(A\) linear and \(c\) fixed, is called an affine map;
when \(A\) is invertible, it is an affine transformation. We shall write nonempty fibers of linearized equations as such translated spaces.

\begin{lemma}[Least dividing period of an affine permutation]\label{lem:affine-least-period}
Let \(V\) be a finite-dimensional \(\F_p\)-space, and let \(\phi:V\to V\) be an affine permutation,
that is, \(\phi(v)=A(v)+b\), where \(A\) is an invertible \(\F_p\)-linear map and \(b\in V\).
Then \(\phi\) has an orbit of length \(p^e\) whose length divides every orbit length.
More precisely, if \(v_p(\ell)\) denotes the exponent of \(p\) in the integer \(\ell\), then
\[
e=\min\{v_p(\ell):\ell\text{ is an orbit length of }\phi\}.
\]
\end{lemma}
\begin{proof}
Choose an orbit whose length \(\ell=p^e m\) has the least possible \(p\)-adic exponent, where \(p\nmid m\), and take a point \(x\) on this orbit.
Since \(m\) is invertible in \(\F_p\), we may define
\[
y=\frac1m\sum_{j=0}^{m-1}\phi^{jp^e}(x).
\]
Affine maps preserve linear combinations whose coefficients sum to one, so
\[
\phi^{p^e}(y)
=\frac1m\sum_{j=0}^{m-1}\phi^{(j+1)p^e}(x)=y;
\]
the last equality follows by shifting the sum cyclically and using \(\phi^{mp^e}(x)=x\).
Hence the period of \(y\) is some \(p^a\), with \(a\le e\).
The minimality of \(e\) also gives \(a\ge e\), so the period of \(y\) is exactly \(p^e\).
Every other orbit length has \(p\)-adic exponent at least \(e\), and is therefore divisible by \(p^e\).
\end{proof}

\begin{proposition}[Computing the least period of an affine permutation]\leavevmode\phantomsection\label{prop:affine-least-period-computation}
With the notation of Lemma~\ref{lem:affine-least-period}, put \(r=\dim_{\F_p}V\), \(N=A-I\), and
\(h=\lceil\log_p(r+1)\rceil\). The least period dividing all orbit lengths is \(p^e\), where
\begin{equation}\label{eq:affine-least-period-image}
e=\min\{0\le j\le h:N^{p^j-1}b\in\operatorname{im}N^{p^j}\}.
\end{equation}
In particular, \(e\le h\). After choosing a basis of \(V\), the condition in the formula is equivalent to
\[
\operatorname{rank}\bigl[N^{p^j}\mid N^{p^j-1}b\bigr]
=\operatorname{rank}N^{p^j},
\]
so finitely many matrix operations determine \(e\). Here \(N^0=I\); when \(j=0\), the condition is \(b\in\operatorname{im}(A-I)\), that is, \(\phi\) has a fixed point.
\end{proposition}
\begin{proof}
For \(m=p^j\), polynomial identities in characteristic \(p\) give
\[
A^m=(I+N)^m=I+N^m,\qquad
I+A+\cdots+A^{m-1}=N^{m-1}.
\]
The second identity follows from \(\sum_{i=0}^{m-1}(1+T)^i=T^{m-1}\), so it does not require \(N\) to be invertible. Iterating the affine map gives
\[
\phi^{p^j}(v)=v+N^{p^j}v+N^{p^j-1}b.
\]
Thus \(\phi^{p^j}\) has a fixed point if and only if \(N^{p^j-1}b\in\operatorname{im}N^{p^j}\).
By Lemma~\ref{lem:affine-least-period}, there is an orbit of length \(p^e\), and its length divides all orbit lengths;
therefore \(\phi^{p^j}\) has a fixed point exactly when \(j\ge e\). This proves the minimum formula and the rank criterion.

It remains to show that a suitable index always exists in \(0\le j\le h\).
The descending chain \(V\supseteq\operatorname{im}N\supseteq\operatorname{im}N^2\supseteq\cdots\) stabilizes by step \(r\) at the latest:
once two consecutive terms are equal, further applications of \(N\) preserve equality; if every preceding inclusion is strict, the dimension has already fallen to zero.
Since \(p^h-1\ge r\), we have
\[
N^{p^h-1}b\in\operatorname{im}N^{p^h-1}
=\operatorname{im}N^{p^h}.
\]
Hence \(e\le h\). Changing the affine base point to \(v_0\) changes the translation term to \(b+Nv_0\),
so the vector on the left of the membership condition in~\eqref{eq:affine-least-period-image} changes only by \(N^{p^j}v_0\). The criterion is therefore independent of the base point.
\end{proof}
Theorem~\ref{v107:thm-affine-cycles} expresses this image condition through the translation obstruction height and combines it with elementary divisors to count fixed points of every iterate and all cycles.

\begin{theorem}\label{thm:affine-linearized-obstruction}
Let \(L\in k[X]\) be a \(p\)-linearized polynomial of ordinary degree \(p^v\), where \(v\ge1\), and let \(b\in k\).
If \(L(X)-b\) is irreducible in \(k[X]\), then \(v=1\) or \(p=v=2\).
\end{theorem}
\begin{proof}
If \(L'=0\), then \(L-b\) is the \(p\)th power of a polynomial over the finite coefficient field,
and cannot be irreducible of positive degree. Hence \(L'\ne0\).
In a common finite splitting extension, the \(p^v\) distinct roots of \(L\) form a \(v\)-dimensional \(\F_p\)-space \(V\).
Take a root \(\beta\) of \(L-b\). Additivity gives \(L(\beta+x)-b=L(x)\), so the full root set of \(L-b\) is \(\beta+V\).

Put \(\sigma:x\mapsto x^q\). Since the coefficients lie in \(k\), we have \(\sigma(V)=V\) and \(\sigma(\beta)-\beta\in V\).
In the coordinates \(\beta+x\leftrightarrow x\), the map \(\sigma\) is therefore the affine permutation of \(V\) given by
\[
x\longmapsto Ax+c,\qquad A=\sigma|_V,\quad c=\sigma(\beta)-\beta.
\]
Irreducibility means that the full root set is a single orbit of length \(p^v\).
Applying Proposition~\ref{prop:affine-least-period-computation} to this \(v\)-dimensional affine space gives
\[
v\le\lceil\log_p(v+1)\rceil.
\]
When \(v\ge3\), we have \(2^{v-1}\ge v+1\), so the right-hand side is at most \(v-1\), a contradiction.
When \(v=2\) and \(p\ge3\), the right-hand side is one, again a contradiction.
Thus only \(v=1\) or \(p=v=2\) is possible.
\end{proof}
This restriction compares the unique orbit length required by irreducibility with the upper bound for the shortest orbit length in an affine space of the given dimension.

\begin{proposition}[Complete form of the quartic exception]\label{prop:affine-quartic}
Let \(\Char k=2\), let \(L\in k[X]\) be monic and \(2\)-linearized with \(\deg L=4\), and let \(b\in k\). Then \(L(X)+b\) is irreducible if and only if there exist \(u,v\in k^\times\) such that
\[
L(X)=(X^2+uX)^2+v(X^2+uX)
=X(X+u)(X^2+uX+v),
\]
and both \(X^2+uX+v\) and \(X^2+vX+b\) are irreducible.
\end{proposition}
\begin{proof}
First suppose that \(L+b\) is irreducible. Its derivative is nonzero, so the roots of \(L\) form a two-dimensional \(\F_2\)-space \(V\).
As in the preceding proof, Frobenius acts as a four-cycle on the affine root space; denote its linear part by \(A\).
We have \(A^4=I\), hence \((A-I)^4=0\). The nilpotency index on a two-dimensional space is at most two, so \(A^2=I\).
If \(A=I\), the affine transformation is a translation of order at most two, a contradiction.
Thus \(A-I\ne0\) and \((A-I)^2=0\), and its kernel has dimension one.
The Frobenius-fixed elements of \(V\) are exactly \(\{0,u\}\), with \(u\in k^\times\).

Write \(L=X^4+aX^2+cX\), where \(c\ne0\).
Since \(L(u)=0\), we have \(c=u^3+au=u(a+u^2)\).
Putting \(v=a+u^2\in k^\times\) gives
\[
L=(X^2+uX)^2+v(X^2+uX)
=X(X+u)(X^2+uX+v).
\]
The roots of the final quadratic factor are the other two elements of \(V\), neither of which lies in \(k\),
so \(X^2+uX+v\) is irreducible over \(k\).
If the outer quadratic \(Y^2+vY+b\) were reducible, substituting \(Y=X^2+uX\) would give a nontrivial factorization of \(L+b\).
It too is therefore irreducible, proving necessity.

Conversely, suppose that the stated decomposition and both quadratic irreducibility conditions hold.
Take a root \(\gamma\) of \(Y^2+vY+b\). Then \(E=k(\gamma)=\F_{q^2}\) and \(\Tr_{E/k}(\gamma)=v\).
The fiber polynomial is \(X^2+uX+\gamma\), which is irreducible over \(E\) if and only if
\[
\Tr_{E/\F_2}(\gamma/u^2)
=\Tr_{k/\F_2}(v/u^2)=1.
\]
The right-hand condition is exactly the quadratic trace criterion for irreducibility of \(X^2+uX+v\) over \(k\).
Thus the fiber polynomial is irreducible, and the root fiber theorem gives the irreducibility of \(L+b\).
\end{proof}

\begin{corollary}\label{cor:affine-composition-classification}
Let \(P\) be monic and irreducible of degree \(n\), with coefficient \(c_{n-1}\) of the next-to-leading term, and let \(L\) be monic and \(p\)-linearized of degree \(p^v\), with \(v\ge2\).
Then \(P(L)\) is irreducible if and only if \(p=v=2\), \(n\) is odd, and
\[
L=X(X+u)(X^2+uX+v_0),\qquad u,v_0\in k^\times,
\]
where both \(X^2+uX+v_0\) and \(X^2+v_0X+c_{n-1}\) are irreducible.
\end{corollary}
\begin{proof}
Take a root \(\alpha\) of \(P\), and put \(E=k(\alpha)=\F_{q^n}\).
By the root fiber criterion~\ref{thm:capelli}, \(P(L)\) is irreducible if and only if \(L(X)-\alpha\) is irreducible over \(E\).

First prove necessity. Apply Theorem~\ref{thm:affine-linearized-obstruction} to the fiber over \(E\).
Since \(v\ge2\), we obtain \(p=v=2\). Thus the root space \(V\) of \(L\) is a two-dimensional \(\F_2\)-space,
and \(V\setminus\{0\}\) has exactly three elements.
Put \(\sigma(x)=x^q\). Since \(L\in k[X]\), the map \(\sigma\) permutes these three nonzero roots;
Frobenius relative to \(E\) is \(\sigma^n\).
The necessity proof in Proposition~\ref{prop:affine-quartic} shows that \(\sigma^n\) acts on the three nonzero roots as a transposition.
A permutation of three elements is the identity, a transposition, or a three-cycle. If its \(n\)th power is a transposition,
the original permutation must itself be a transposition and \(n\) must be odd.

Denote the unique nonzero root fixed by this transposition by \(u\). Then \(u\in k^\times\).
Write \(L=X^4+aX^2+cX\), where \(a,c\in k\) and \(c\ne0\). Since \(L(u)=0\),
putting \(v_0=a+u^2=c/u\in k^\times\) gives
\[
L=X(X+u)(X^2+uX+v_0).
\]
The other two nonzero roots are interchanged by \(\sigma\), so \(X^2+uX+v_0\) is irreducible over \(k\).
Applying Proposition~\ref{prop:affine-quartic} to \(L+\alpha\) also shows that \(Y^2+v_0Y+\alpha\) is irreducible over \(E\), that is,
\[
\Tr_{E/\F_2}(\alpha/v_0^2)
=\Tr_{k/\F_2}\!\left(\frac{\Tr_{E/k}(\alpha)}{v_0^2}\right)
=\Tr_{k/\F_2}(c_{n-1}/v_0^2)=1.
\]
Here we used \(\Tr_{E/k}(\alpha)=-c_{n-1}=c_{n-1}\) in characteristic two.
The last equality is equivalent to irreducibility of \(X^2+v_0X+c_{n-1}\) over \(k\).

Now prove sufficiency. Suppose all the stated conditions hold. Since \(n\) is odd, \(\gcd(2,n)=1\),
and the base-extension theorem~\ref{thm:factor-base-extension} ensures that \(X^2+uX+v_0\) remains irreducible over \(E\).
Irreducibility of the second quadratic over \(k\) gives \(\Tr_{k/\F_2}(c_{n-1}/v_0^2)=1\). The trace equality above shows that
\(Y^2+v_0Y+\alpha\) is also irreducible over \(E\).
Proposition~\ref{prop:affine-quartic} therefore implies that \(L+\alpha\) is irreducible over \(E\),
and the root fiber criterion gives the irreducibility of \(P(L)\) over \(k\).
\end{proof}

\section{Root-Field Spectra and Least Root Fields}
\label{sec:root-field-spectrum}

Irreducibility requires all roots to form a single Frobenius orbit. For a reducible polynomial, divisibility among orbit lengths also records inclusions among the fields generated by individual roots.
The fiber method can determine whether the roots of a composite polynomial form a single orbit and can also compare the extensions required by different orbits.

\begin{definition}[Root-field spectrum and least root field]\label{def:root-field-spectrum}
Let \(k=\F_q\), and let \(f\in k[X]\) be nonconstant. Define the \textbf{root-field spectrum} of \(f\) over \(k\) by
\[
\mathcal D_k(f)=\{\deg P:P\in k[X]\text{ is monic irreducible and }P\mid f\}.
\]
In a common finite splitting extension, if a root \(\beta_0\) of \(f\) satisfies
\[
k(\beta_0)\subseteq k(\beta)\qquad\text{for every }\beta\text{ with }f(\beta)=0,
\]
then \(k(\beta_0)\) is called the \textbf{least root field} of \(f\). A polynomial with this property is also called a \textbf{subnormal polynomial}; see Agou~\cite{Agou}.
\end{definition}

The root-field spectrum records only the degrees of distinct irreducible factors, not their numbers or multiplicities.
What is required here is a \emph{least element} under inclusion: a root field contained in every other root field,
not merely a \emph{minimal element} containing no smaller root field. The latter exists for every nonempty finite collection of root fields and does not characterize the property under discussion.

\begin{proposition}\label{prop:least-root-field}
Let \(\mathcal D_k(f)=\{d_1,\ldots,d_s\}\), and put \(d_* =\gcd(d_1,\ldots,d_s)\). The following conditions are equivalent:
\begin{enumerate}
\item \(f\) has a least root field;
\item there exists \(d_0\in\mathcal D_k(f)\) such that \(d_0\mid d\) for every \(d\in\mathcal D_k(f)\);
\item \(d_*\in\mathcal D_k(f)\).
\end{enumerate}
When these conditions hold, the least root field is \(\F_{q^{d_*}}\). In general, the intersection of all root fields is \(\F_{q^{d_*}}\), but this intersection need not be generated by a root of \(f\).
\end{proposition}
\begin{proof}
Every root of an irreducible factor of degree \(d\) generates \(\F_{q^d}\). Within a common extension,
\(\F_{q^a}\subseteq\F_{q^b}\) if and only if \(a\mid b\), so the first two conditions are equivalent.
If the second holds, then \(d_0\) divides all \(d_i\), and some \(d_i=d_0\), so \(d_*=d_0\).
Conversely, if \(d_*\) belongs to the root-field spectrum, it itself satisfies the second condition.
The subfield intersection formula proves the final assertion.
\end{proof}

\begin{proposition}[Root-field spectrum of a composite polynomial]\label{prop:composition-root-spectrum}
Let \(P\in k[Y]\) be monic and irreducible with \(\deg P=n\), let \(\alpha\) be a root, and put \(E=k(\alpha)\).
For every nonconstant \(G\in k[X]\), we have
\begin{equation}\label{eq:composition-root-spectrum}
\mathcal D_k(P(G))
=n\mathcal D_E(G-\alpha)
=\{nd:d\in\mathcal D_E(G-\alpha)\}.
\end{equation}
Thus \(P(G)\) has a least root field over \(k\) if and only if \(G-\alpha\) has a least root field over \(E\).
If the latter has degree \(d_0\) over \(E\), then the former has degree \(nd_0\) over \(k\).
\end{proposition}
\begin{proof}
When \(n>1\), apply Theorem~\ref{thm:fiber-factorization} with \(A=G\) and \(B=1\):
each irreducible factor \(h_j\) in the fiber corresponds to an irreducible factor of degree \(n\deg h_j\) over the base field, and all factors arise this way.
This gives~\eqref{eq:composition-root-spectrum}; multiplicities in the fiber do not affect the identity.
When \(n=1\), we have \(E=k\) and \(P(G)=G-\alpha\), so the equality is immediate.
Finally, \(nd_0\mid nd\) if and only if \(d_0\mid d\). Multiplication by \(n\) therefore preserves and reflects the divisibility order, proving the equivalence of the least-root-field conditions.
\end{proof}

Thus, to find the least root field of a composite polynomial, one may first seek in one fiber the shortest return period dividing all return periods.
For linearized polynomials, Lemma~\ref{lem:affine-least-period} guarantees its existence, and Proposition~\ref{prop:affine-least-period-computation} gives a method of computation and an upper bound.

\begin{corollary}[Least root field of a linearized composition]\label{cor:additive-least-root-field}
Let \(P\in k[Y]\) be monic and irreducible with \(\deg P=n\), let \(k=\F_q\) have characteristic \(p\), and suppose that
\[
L(X)=\sum_{i=0}^r a_iX^{p^i}\in k[X],\qquad a_0a_r\ne0.
\]
Put \(h=\lceil\log_p(r+1)\rceil\). Then there exists \(0\le e\le h\) such that \(np^e\) belongs to \(\mathcal D_k(P(L))\) and divides each of its elements.
In particular, \(P(L)\) has least root field \(\F_{q^{np^e}}\).
If \(\alpha\) is a root of \(P\) and \(E=k(\alpha)\), this exponent can also be determined directly by polynomial operations:
\begin{equation}\label{eq:additive-least-root-field-gcd}
e=\min\left\{0\le j\le h:
\gcd_{E[X]}\bigl(L(X)-\alpha,\,X^{q^{np^j}}-X\bigr)\ne1\right\}.
\end{equation}
\end{corollary}
\begin{proof}
In a finite extension over which \(L(X)-\alpha\) splits completely, take a root \(\beta\) and put \(V=\ker L\).
Since \(L'=a_0\ne0\), the kernel \(V\) is an \(\F_p\)-space of \(p^r\) elements, and the root set of \(L(X)-\alpha\) is exactly \(\beta+V\).
Put \(\sigma(x)=x^q\). The condition that the coefficients lie in \(k\) gives \(\sigma^n(V)=V\), and since \(\sigma^n(\alpha)=\alpha\), we have
\[
c=\sigma^n(\beta)-\beta\in V.
\]
Thus, in the coordinates \(\beta+v\leftrightarrow v\), the fiber return map is
\[
v\longmapsto\sigma^n(v)+c.
\]
Its orbit lengths are exactly the degrees in \(\mathcal D_E(L-\alpha)\).
Lemma~\ref{lem:affine-least-period} gives the least dividing period \(p^e\),
and Proposition~\ref{prop:affine-least-period-computation} gives \(e\le h\).
Applying Proposition~\ref{prop:composition-root-spectrum} now yields the conclusion on the least root field.

The greatest common divisor in~\eqref{eq:additive-least-root-field-gcd} is nonconstant if and only if the fiber contains a root fixed by \(\sigma^{np^j}\).
Such a root exists if and only if the return map has an orbit whose length divides \(p^j\).
Since \(p^e\) divides all return periods and is itself a return period, this is equivalent to \(j\ge e\).
Hence the minimum is precisely the exponent \(e\) above.
For computation, one can first find the remainder of \(X^{q^{np^j}}\) modulo \(L-\alpha\), then take the greatest common divisor, without expanding the high-degree polynomial.
\end{proof}

For \(L=X^{p^r}-aX\), with \(a\ne0\), the kernel above is an \(r\)-dimensional \(\F_p\)-space;
the same argument handles any separable \(p\)-linearized polynomial.
The least root field exists because of the affine return action; this does not mean that all return orbits have the same length.
The least root field does not determine the full set of factor degrees. These depend on all return-cycle lengths, which can also be recovered from fixed-point counts of iterates by M\"obius inversion.

This least-root-field conclusion also applies to inseparable nonconstant \(p\)-linearized polynomials.
Indeed, perfection of finite fields allows us to write
\[
L=L_0^{p^u},\qquad L_0'\ne0,
\]
where \(L_0\in k[X]\) is again \(p\)-linearized.
Take the unique \(\alpha_0\) in \(E\) with \(\alpha_0^{p^u}=\alpha\). Then
\[
L(X)-\alpha=(L_0(X)-\alpha_0)^{p^u}.
\]
Since \(\alpha_0\in E\), the distinct roots on the right still form an affine space stable under the relative Frobenius \(\sigma^n:x\mapsto x^{q^n}\),
so the same averaging argument gives a least root field;
the common multiplicity \(p^u\) changes only factor multiplicities, not the root-field spectrum.

\begin{example}[A least root field with unequal factor degrees]\label{exm:additive-spectrum-two-six}
Put \(F=X^8+X+1\) in \(\F_2[X]\). Its roots satisfy
\[
\sigma^3(\beta)=\beta+1,\qquad \sigma^6(\beta)=\beta,
\qquad \sigma(\beta)=\beta^2.
\]
Thus every root orbit has length dividing six but not three, so its length is either two or six.
In \(\F_4\), we have \(x^8=x^2\), so the roots of \(F\) in \(\F_4\) are exactly the two roots of \(X^2+X+1\).
Also, \(F'=1\), so the remaining six roots must form a single orbit of length six. Polynomial division therefore gives the full irreducible factorization
\[
X^8+X+1=(X^2+X+1)(X^6+X^5+X^3+X^2+1).
\]
Hence \(\mathcal D_{\F_2}(F)=\{2,6\}\); the least root field is \(\F_4\), whereas the splitting field is \(\F_{64}\).
Here the root fiber is \(\beta+\F_8\), and the third iterate of Frobenius acts on it as the nonidentity translation \(x\mapsto x+1\);
the two degrees, two and six, reflect both the additive fiber and the Frobenius periods.
\end{example}

\par\Needspace{11\baselineskip}
\section*{Exercises}
\phantomsection
\addcontentsline{toc}{section}{Exercises}
\begin{exercise}[A trace-zero change in an extension of degree eight]\label{v70:ex-octic-trace-zero}
Over \(\F_2\), put
\[
f(X)=X^8+X^7+X^3+X+1,\qquad
 g(Y)=Y^8+Y^5+Y^3+Y^2+1.
\]
\begin{enumerate}
\item Prove that \(f\) is irreducible by excluding irreducible factors of degree at most four.
\item Let \(\alpha\) be a root of \(f\). Prove that no fractional linear transformation over \(\F_2\) sends \(\alpha\) to a trace-zero generator.
\item Prove that \(\beta=\alpha^2+\alpha\) still generates \(\F_{2^8}\) and has trace zero. Verify
\[
g(X^2+X)=f(X)f(X+1),
\]
and deduce the minimal polynomial of \(\beta\).
\end{enumerate}
\end{exercise}

\begin{exercise}[Classification and counting of constant terms of irreducible binomials]\label{ex:supp-binomial}
Let \(q\) be a prime power and \(n\ge1\). Use the notation \(\mathcal A_q(n)\), and let \(\rad(n)\) denote the squarefree kernel of \(n\), with \(\rad(1)=1\).
\begin{enumerate}
\item Treat \(n=1\) and the zero coefficient separately. For \(n\ge2\), prove that \(\mathcal A_q(n)\ne\varnothing\) if and only if
\begin{equation}\label{v77:eq-binomial-existence}
\rad(n)\mid q-1,\qquad
4\mid n\ \Longrightarrow\ q\equiv1\pmod4.
\end{equation}
Explain why \(n\mid q-1\) is not required.
\item Suppose \(n\ge2\) and \eqref{v77:eq-binomial-existence} holds. Fix a primitive element \(g\) of \(\F_q^\times\), and prove the complete parametrization
\begin{equation}\label{v77:eq-binomial-constant-exponents}
\mathcal A_q(n)=\{-g^k:0\le k\le q-2,\ \gcd(k,n)=1\}.
\end{equation}
For \(a\in\F_q^\times\) and \(e=\ord(-a)\), prove the equivalent prime-by-prime tests:
\begin{equation}\label{v77:eq-binomial-coefficient-tests}
\begin{aligned}
a\in\mathcal A_q(n)
&\iff (-a)^{(q-1)/\ell}\ne1
&&\text{for every prime }\ell\mid n,\\
&\iff v_\ell(e)=v_\ell(q-1)>0
&&\text{for every prime }\ell\mid n.
\end{aligned}
\end{equation}
\item Derive the exact count
\begin{equation}\label{v77:eq-binomial-count}
|\mathcal A_q(n)|=
\begin{cases}
q,&n=1,\\[2pt]
(q-1)\dfrac{\varphi(n)}n,
&n\ge2\text{ and \eqref{v77:eq-binomial-existence} holds},\\[5pt]
0,&\text{otherwise}.
\end{cases}
\end{equation}
Explain why, among degrees \(n\ge2\) satisfying the existence conditions, the admissible coefficient set depends only on the distinct prime divisors of \(n\). The condition modulo four must still be retained when comparing different degrees.
\item Determine \(\mathcal A_5(8)\), \(\mathcal A_7(9)\), and \(\mathcal A_3(4)\), and use them to decide the irreducibility of
\(X^8-2\in\F_5[X]\), \(X^9-2\in\F_7[X]\), and \(X^4+1\in\F_3[X]\).
Use the classification to verify the two infinite families in Example~\ref{v77:ex-binomial-constant-families}.
Prove that \(\mathcal A_2(n)=\varnothing\) for every \(n\ge2\), and explain why this does not extend to all finite fields of characteristic two.
Finally, compute \(\mathcal A_7(2)\) to explain why the tests apply to \(-a\), and show that requiring \(-a\) merely not to be an \(n\)th power does not in general ensure irreducibility.
\end{enumerate}
\end{exercise}
\begin{exercise}\label{ex:supp-fiber-binary}
Prove that \(X^4+X^3+X^2+X+1\) is irreducible over \(\F_2\), and determine the degrees and number of its irreducible factors over \(\F_4\).
\end{exercise}
\begin{exercise}\label{ex:supp-affine-quartic}
Explain why \(X^4+X+1\) is irreducible over \(\F_2\), whereas no \(2\)-linearized polynomial of degree eight plus a constant can be irreducible.
\end{exercise}
\begin{exercise}\label{ex:supp-composite-as}
Let \(P\in\F_3[X]\) be monic and irreducible of degree three, with nonzero next-to-leading coefficient. Prove that \(P(X^3-X-b)\) is irreducible for every \(b\in\F_3\).
\end{exercise}

\begin{exercise}[Necessity of the linearized hypothesis]\label{ex:wu-quartic-hypothesis}
Prove that \(f=X^4+X^3+1\in\F_2[X]\) is irreducible.
Prove that its four roots cannot form an affine coset of a two-dimensional \(\F_2\)-linear space.
Which hypothesis in the affine quartic criterion does this show to be indispensable?
\end{exercise}

\begin{exercise}[Recovering factor degrees and multiplicities from one root fiber]\label{ex:wu2-fiber-multiplicity}\label{ex:return-factor-capstone}
Let \(k=\F_2\), \(P(Y)=Y^2+Y+1\), take \(\alpha\in\F_4\) with \(P(\alpha)=0\), and write \(\sigma(x)=x^2\).
\begin{enumerate}
\item Compute the root fibers, return maps, and root multiplicities of \(X^2=\alpha\) and \(X^2+X=\alpha\), respectively, and hence find the irreducible factorizations of \(P(X^2)\) and \(P(X^2+X)\).
\item Put \(L(X)=X^8+X\), and, in a splitting extension, put
\[
S_\alpha=\{x:L(x)=\alpha\},\qquad R=\sigma^2|_{S_\alpha}.
\]
Prove that \(|S_\alpha|=8\), \(R^3(x)=x+1\), and \(R^6=\id\). Prove that \(R^2\) has exactly two fixed points, and determine all cycle lengths of \(R\) and their numbers.
\item Without first factoring \(P(L)\), use these return cycles to determine the degrees, numbers, and multiplicities of its irreducible factors over \(k\), and then verify
\[
\begin{aligned}
P(X^8+X)
&=X^{16}+X^8+X^2+X+1\\
&=(X^4+X+1)
  (X^{12}+X^9+X^8+X^6+X^3+X^2+1).
\end{aligned}
\]
Explain why the return periods in the second part are not directly the degrees of irreducible factors over \(k\), and why return cycles alone do not contain multiplicity information.
\end{enumerate}
\end{exercise}

\begin{exercise}[Least root fields and base extension]\label{ex:root-spectrum-base-extension}
Let \(f\in\F_q[X]\) be nonconstant, and let \(r\ge1\).
\begin{enumerate}
\item Prove that
\[
\mathcal D_{\F_{q^r}}(f)
=\left\{\frac d{\gcd(d,r)}:d\in\mathcal D_{\F_q}(f)\right\}.
\]
\item Prove that if \(f\) has a least root field over \(\F_q\), then it also has one over \(\F_{q^r}\).
\item Take irreducible polynomials \(P_2,P_3\) over \(\F_q\) of degrees two and three, respectively.
Compare the existence of a least root field for \(P_2P_3\) over \(\F_q\) and over \(\F_{q^2}\), showing that the converse of the preceding assertion fails.
Also explain why existence of the greatest common divisor of the root-field spectrum does not suffice to ensure existence of a least root field.
\end{enumerate}
\end{exercise}

\begin{exercise}[Preserving generation and reaching a prescribed defining equation]\label{v79:ex-rational-generation}
Let \(k=\F_q\), \(E=k(\alpha)\), and \([E:k]=n\ge2\). Let \(\ell\) be the least prime divisor of \(n\), and put \(r=\lfloor n/2\rfloor\).
\begin{enumerate}
\item Let \(R=A/B\in k(X)\) be nonconstant and in lowest terms, with \(\deg R<\ell\). Show that \(B(\alpha)\ne0\) and \(k(R(\alpha))=E\).
\item Prove that every \(\beta\in E\setminus k\) can be written as \(\beta=R(\alpha)\) for a nonconstant rational function \(R=A/B\) in lowest terms satisfying
\[
\deg A,\deg B\le r,\qquad B(\alpha)\ne0.
\]
\emph{Hint.} Express \(A(\alpha)-\beta B(\alpha)=0\) in the power basis of \(\alpha\), and compare the numbers of equations and unknown coefficients.
\item Deduce that, for any prescribed monic irreducible binomial of degree \(n\) over \(k\), every given generator of \(E/k\) can be transformed by a rational function of degree at most \(r\) into a generator with that binomial as its minimal polynomial.
Explain why this conclusion does not assert that degree at most two always suffices.
\end{enumerate}
\end{exercise}

\begin{exercise}[Resultants and binomial defining equations]\label{v79:ex-binomial-defining-equation}
Let \(k=\F_q\), let \(f\in\Irr_q(n)\), with \(n\ge2\), and let \(\alpha\) be a root of \(f\). Put \(E=k(\alpha)\).
Let \(R=A/B\in k(X)\) be nonconstant and in lowest terms, with \(B(\alpha)\ne0\), and put \(\beta=R(\alpha)\) and \(d=[k(\beta):k]\).
\begin{enumerate}
\item Prove
\[
\operatorname{Res}_X\bigl(f(X),B(X)Y-A(X)\bigr)
=\Nm_{E/k}\bigl(B(\alpha)\bigr)\,m_\beta(Y)^{n/d}.
\]
Explain how the formula gives a new monic defining polynomial for \(E/k\) when \(d=n\), and what the repeated factors record when \(d<n\).
\item Take \(k=\F_5\), \(f(X)=X^2+X+1\), and \(R(X)=X+3\). Prove that \(f\) is irreducible, that \(\beta^2=3\), and that
\[
\operatorname{Res}_X\bigl(f(X),Y-X-3\bigr)=Y^2-3.
\]
Compare the numbers of nonzero terms of the two defining polynomials. Compute \(\ord(\alpha)\) and \(\ord(\beta)\), and identify the complete equal-difference exponent orbit of the roots of \(Y^2-3\).
\end{enumerate}
\end{exercise}

\begin{exercise}[Irreducible-polynomial counts classified by trace]\label{v109:ex-binary-trace}
Let $A_i(n)$ count monic irreducibles of degree $n\ge1$ over $\F_2$ whose next-to-leading coefficient is $i$; include $X$ in degree one.
Let $f$ be monic irreducible of degree $n$. Prove that $f(X^2+X)$ is irreducible exactly when the next-to-leading coefficient of $f$ is one. If that coefficient is zero, prove that the composition has two distinct irreducible factors of degree $n$.
Use M\"obius inversion to show
\[
A_1(n)=\frac1{2n}\sum_{\substack{d\mid n\\d\text{ odd}}}\mu(d)2^{n/d},\qquad
A_1(n)-A_0(n)=\begin{cases}0,&n\text{ odd},\\A_1(n/2),&n\text{ even}.\end{cases}
\]
Compute both counts for $n=1,2,3,4,6$. This is an elementary instance of recursive counting induced by a translation quotient; see \cite{SchulzRecursive2026}.
\end{exercise}

\begin{exercise}[Equal fixed-point counts and different return periods]\label{v115:ex-return-observations}
Two return permutations on a six-point fiber have cycle lengths \((1,1,4)\) and \((1,1,2,2)\).
Compute their \(r_1,r_2,r_3,r_4\), and recover the cycle counts using~\eqref{v115:eq-return-counts}.
If the base orbit has length \(n\), list the full factor degrees. Can the fiber size and \(r_1\) alone distinguish these two factorization types?
\end{exercise}

\chapter{Affine Symmetries of Root Sets}
\label{chap:root-symmetry}
Frobenius orbits determine irreducible factors; translations, scalings, and affine transformations preserving root multiplicities determine another partition.
Translation orbits yield additive decompositions, while full affine orbits give a normal form that also records the special orbit and its multiplicity.
An affine symmetry of the quotient lifts to an affine transformation of the original root set only when its slope satisfies a compatibility condition.
\leavevmode\phantomsection\label{ffs:def-affine-slope}\ignorespaces Here the \textbf{slope} of an affine transformation \(X\mapsto aX+b\) is its linear coefficient \(a\); compatibility requires scalar multiplication \(v\mapsto av\) to preserve the translation subspace.
Conversely, polynomial differences and coefficient comparison recover these symmetries without first finding the individual roots.
This chapter recovers affine stabilizers from root multiplicities and coefficient identities, then describes factorization through orbits and Frobenius return.
General polynomial decomposition need not arise from affine symmetry; Chapter~\ref{chap:subfield-representations} turns instead to low-degree representations of subfields.

Fix \(k=\F_q\), with \(q=p^r\), and write \(\sigma(x)=x^q\) for relative Frobenius.
Root sets and root multiplicities are compared in a finite extension \(E/k\) over which the polynomials concerned split completely; enlarge \(E\) when necessary.
We use \(X\) for polynomials in the original coordinate and \(Y,B\) as auxiliary variables.

\leavevmode\phantomsection\label{v58:def-multiset}\ignorespaces
A \textbf{multiset} is an unordered collection in which elements may occur with multiplicity.
A root multiset thus records both the roots and their multiplicities.

\par\smallskip
\Needspace{3\baselineskip}
\begingroup\small
\noindent\phantomsection\label{box:descent-factor-composition}\textbf{Descent, factorization, and compositional blocks of root sets}
\par\nobreak\smallskip
\noindent\textbf{Descent concerns stability.} Frobenius preserves root multiplicities exactly when the monic root polynomial is defined over the base field. This is Theorem~\ref{thm:root-descent} with multiplicities retained.
\par\smallskip
\noindent\textbf{Factors correspond to orbits.} Each Frobenius orbit gives a monic irreducible factor; its length and common root multiplicity determine the degree and exponent, respectively.
\par\smallskip
\noindent\textbf{Composition requires complete fibers.} A Frobenius-compatible partition is not enough: in the irreducible case, a candidate root block must be a complete fiber of a suitable base-field polynomial; see Theorem~\ref{sfr:thm-optimal-composition}. Stabilizers describe the symmetry of certain fibers, whereas low-degree representations do not require such symmetry.
\par\endgroup
\medskip

\phantomsection\label{unit:root-symmetry-return}
Sections~\ref{ffs:sec-stabilizers}--\ref{ffs:sec-irreducible} recover translation and affine symmetries from root multiplicities, then determine factor degrees and multiplicities by fiber return; irreducibility also yields conjugate-difference criteria, recovery on quotients, and low-degree classifications.

\section{Root Multiplicities and Three Stabilizers}
\label{ffs:sec-stabilizers}

\begin{definition}\label{ffs:def-root-multiplicity}
Let \(f\in k[X]\) be monic and nonconstant, and put \(n=\deg f\). Define its \textbf{root multiplicity function} by
\[
\nu_f(\alpha)=\max\{j\ge0:(X-\alpha)^j\mid f\}\qquad(\alpha\in E).
\]
The set \(S_f=\{\alpha\in E:\nu_f(\alpha)>0\}\) consists of the distinct roots of \(f\).
At a point that is not a root, the multiplicity is zero; the element zero itself may or may not be a root.
\end{definition}

\begin{definition}\label{ffs:def-stabilizers}
Define the \textbf{additive stabilizer}, \textbf{multiplicative stabilizer}, and \textbf{affine stabilizer} of \(f\), respectively, by
\begin{align*}
V_f&=\{b\in E:\nu_f(\alpha+b)=\nu_f(\alpha)\text{ for all }\alpha\in E\},\\
M_f&=\{a\in E^\times:\nu_f(a\alpha)=\nu_f(\alpha)\text{ for all }\alpha\in E\},\\
\mathcal A_f&=\{g_{a,b}:X\mapsto aX+b:\ a\in E^\times,\ b\in E,\
\nu_f(g_{a,b}(\alpha))=\nu_f(\alpha)\text{ for all }\alpha\in E\}.
\end{align*}
The group operation is composition of transformations, namely
\[
g_{a,b}\circ g_{a',b'}=g_{aa',\,b+ab'}.
\]
For any field \(K\), all transformations \(X\mapsto aX+b\), with \(a\in K^\times,b\in K\), form the \textbf{one-dimensional affine group}, denoted by \(\operatorname{AGL}_1(K)\).
We call \(\mathcal A_f(k)=\mathcal A_f\cap\operatorname{AGL}_1(k)\) the \textbf{affine stabilizer over the base field}.
\end{definition}

These groups preserve the entire root multiplicity function and may permute distinct roots; \(\mathcal A_f(k)\) additionally requires the transformation coefficients to lie in the base field.
\leavevmode\phantomsection\label{ffs:def-point-stabilizer}\ignorespaces For an action of a group \(G\) on a set,
the \textbf{point stabilizer} of \(\alpha\) is \(G_\alpha=\{g\in G:g(\alpha)=\alpha\}\).
An action is called \textbf{free} if every point stabilizer is trivial.
The base-field affine stabilizer restricts transformation coefficients, whereas a point stabilizer requires a specified point to be fixed; these are generally different.
In particular, the map \(a\mapsto g_{a,0}\) gives a group isomorphism
\[
M_f\simeq(\mathcal A_f)_0=\{g_{a,0}:a\in M_f\}.
\]
\leavevmode\phantomsection\label{ffs:def-slope-image}\ignorespaces
The map \(\mathcal A_f\to E^\times\), \(g_{a,b}\mapsto a\), is called the \textbf{slope map}.
Its image is the \textbf{image of the slope map}, namely the set
\[
\{a\in E^\times:\text{there exists }b\in E\text{ with }g_{a,b}\in\mathcal A_f\}.
\]
Here \(b\ne0\) is allowed, so this image generally differs from \(M_f\), which includes only coefficients of scalings about the origin.

\begin{proposition}[Substitution identity]\label{ffs:prop-affine-identity}
For \(a\ne0\), we have
\[
g_{a,b}\in\mathcal A_f
\quad\Longleftrightarrow\quad
f(aX+b)=a^n f(X).
\]
In particular, \(b\in V_f\iff f(X+b)=f(X)\), and \(a\in M_f\iff f(aX)=a^nf(X)\).
All these equalities are identities of ordinary polynomials.
\end{proposition}
\begin{proof}
The substitution \(X\mapsto aX+b\) is invertible and preserves the multiplicity of every linear factor.
Thus the root multiplicities agree if and only if the two polynomials differ by a nonzero constant factor; comparing leading terms gives the factor \(a^n\).
\end{proof}

The factor \(a^n\) generally cannot be omitted. For example, scaling \(X\mapsto aX\), with \(a\in\F_{q^d}^\times\), changes \(X^{q^d}-X\) into \(a(X^{q^d}-X)\).
Equality of values at all base-field points also cannot replace a polynomial identity, since the nonzero polynomial \(X^q-X\) vanishes throughout the base field.

\begin{proposition}[Finiteness and the field of the parameters]\label{ffs:prop-finiteness}
The set \(V_f\) is a finite \(\F_p\)-linear subspace, and \(|V_f|\mid |S_f|\).
If \(|S_f|\ge2\), then \(\mathcal A_f\) is finite, all its parameters already lie in the splitting field of \(f\), and enlarging the splitting extension does not enlarge \(V_f,M_f,\mathcal A_f\).
\end{proposition}
\begin{proof}
Composition and inversion of translations make \(V_f\) an additive subgroup, hence an \(\F_p\)-subspace.
Fix \(\alpha\in S_f\). The map \(b\mapsto\alpha+b\) injects \(V_f\) into \(S_f\), proving finiteness.
A nonzero translation fixes no point, so each translation orbit has size \(|V_f|\), giving the divisibility assertion.

Take two distinct roots \(\alpha,\beta\). If an affine transformation sends them to \(\gamma,\eta\), respectively, then
\[
a=\frac{\eta-\gamma}{\beta-\alpha},\qquad b=\gamma-a\alpha.
\]
The parameters therefore lie in the splitting field of the polynomial and are uniquely determined by these two images.
There are only finitely many such ordered image pairs, and the same reasoning applies in every larger extension.
\end{proof}

\begin{remark}[The case of a single distinct root]\label{ffs:rem-single-root}
If \(|S_f|=1\), then \(f=(X-c)^n\); Frobenius fixes the unique root, so \(c\in k\).
In this case \(V_f=\{0\}\), and all affine stabilizing transformations have the form \(X\mapsto c+a(X-c)\).
Allowing a larger field of parameters also allows more slopes \(a\).
Thus, whenever we discuss the full finite affine structure independently of the splitting extension below, we assume \(|S_f|\ge2\).
The results on translation differences and the common \(p\)-power do not require this assumption.
\end{remark}

\begin{proposition}[Compatibility with Frobenius]\label{ffs:prop-frobenius-normalizes}
We have \(\nu_f(\alpha^q)=\nu_f(\alpha)\), \(\sigma(V_f)=V_f\), and
\[
\sigma\circ g_{a,b}\circ\sigma^{-1}=g_{a^q,b^q}.
\]
Thus Frobenius maps \(\mathcal A_f\) to itself and commutes with every transformation in \(\mathcal A_f(k)\).
The latter conclusion does not in general apply to the whole of \(\mathcal A_f\).
\end{proposition}
\begin{proof}
Applying coefficientwise Frobenius to the linear factorization of \(f\) leaves \(f\) unchanged, so root multiplicities are preserved.
The conjugation formula follows by direct substitution; taking \(q\)th powers of the coefficients in the affine identity also shows that the right-hand transformation remains in the stabilizer.
The transformation commutes with \(\sigma\) if and only if \(a^q=a,b^q=b\).
\end{proof}

\begin{remark}[Addition of root-of-unity indices]\label{ffs:rem-index-addition}
When \(f(0)\ne0\), we may choose \(N\) with \(p\nmid N\) and a primitive \(N\)th root of unity \(\zeta\) such that \(S_f\subseteq\mu_N\).
Write \(\Omega(i)=\nu_f(\zeta^i)\) on \(\Z/N\Z\). Then
\[
\{j:\Omega(i+j)=\Omega(i)\text{ for all }i\}
\longrightarrow M_f,\qquad j\longmapsto\zeta^j
\]
is a group isomorphism. Indeed, \(\zeta^{i+j}=\zeta^j\zeta^i\), and a multiplier preserving the root set is a ratio of two nonzero roots, hence lies in \(\mu_N\).
Thus the additive stabilizer of the root-of-unity indices is another representation of the multiplicative stabilizer of the root values; it differs from \(V_f\), which is defined by translations of root values.
\end{remark}

\section{Recovering the Common \texorpdfstring{\(p\)}{p}-Power and Stabilizers from Coefficients}
\label{ffs:sec-coefficient-recovery}

\subsection{The Common \texorpdfstring{\(p\)}{p}-Power and Hasse Derivatives}

\begin{proposition}[Common \(p\)-power]\label{ffs:prop-common-power}
Write \(f=\sum_{i=0}^n c_iX^i\), with \(c_n=1\). Then
\[
u=\min_{\substack{i>0\\c_i\ne0}}v_p(i)
 =\min_{\alpha\in S_f}v_p(\nu_f(\alpha)).
\]
Put \(e=p^u\), called the \textbf{common \(p\)-power} of \(f\): the largest common \(p\)-power divisor of its root multiplicities. There is a unique monic \(h\in k[X]\) such that
\[
f=h^e,\qquad h'\ne0.
\]
Moreover, \(\nu_f=e\nu_h\), so \(f,h\) have the same three stabilizer groups defined above, with transformation parameters in a field.
\end{proposition}
\begin{proof}
The \(p\)th-power map is bijective on a finite field. Thus \(f\in k[X^{p^v}]\) if and only if \(f\) is the \(p^v\)th power of a monic polynomial in \(k[X]\), and that monic polynomial is unique.
From the coefficients, the largest \(v\) is the minimum valuation in the first expression.
From the root factorization, being a \(p^v\)th power is equivalent to every root multiplicity being divisible by \(p^v\):
the latter condition first gives a unique monic polynomial over the splitting field with the required power; uniqueness makes it fixed under coefficientwise Frobenius, so it belongs to \(k[X]\).
The two maxima are therefore equal. After extracting the largest common power, \(h\) is no longer a \(p\)th power, that is, \(h'\ne0\).
Dividing all multiplicities by \(e\) does not change the transformations preserving them.
\end{proof}

The condition \(h'\ne0\) means that \(h\) is not a \(p\)th power; it does not imply that \(h\) is square-free.
For example, in characteristic three, \((X-1)(X+1)^2\) has common \(p\)-power one but still has a double root.

\begin{definition}\leavevmode\phantomsection\label{ffs:def-hasse}
The \(j\)th \textbf{Hasse derivative} of \(f=\sum_i c_iX^i\) is defined by
\[
\partial^{[j]}f(X)=\sum_{i\ge j}\binom ij c_iX^{i-j}.
\]
Equivalently, \(f(X+B)=\sum_j\partial^{[j]}f(X)B^j\).
\end{definition}

\begin{proposition}\label{ffs:prop-hasse-multiplicity}
For \(m\ge1\), we have
\[
\nu_f(\alpha)\ge m
\quad\Longleftrightarrow\quad
\partial^{[j]}f(\alpha)=0\quad(0\le j<m).
\]
\end{proposition}
\begin{proof}
Substitute \(X=\alpha+Z\) and expand in powers of \(Z\). The coefficient of \(Z^j\) is exactly \(\partial^{[j]}f(\alpha)\);
divisibility of this expansion by \(Z^m\) is equivalent to vanishing of its first \(m\) coefficients.
\end{proof}

Hasse derivatives require no division by \(j!\), so the criterion holds for arbitrary multiplicities.
For example, in characteristic \(p\), all positive-order ordinary derivatives of \(f=X^p\) vanish, whereas \(\partial^{[p]}f=1\); the criterion still gives the correct multiplicity \(p\) of the zero root.
Thus the vanishing of ordinary derivatives caused by a common \(p\)-power does not prevent recovery of root multiplicities from the coefficients of the translation expansion.
\begin{proposition}[Remainder coordinates and root evaluations with multiplicities]\label{v97:prop-hasse-coordinates}
Let \(f\in k[X]\) be monic of degree \(d\ge1\), and write \(f=\prod_\alpha(X-\alpha)^{m_\alpha}\) over a splitting field \(E\).
The map
\[
 E[X]_{<d}\longrightarrow\prod_\alpha E^{m_\alpha},\qquad
 P\longmapsto\bigl(\partial^{[j]}P(\alpha)\bigr)_{\alpha,\,0\le j<m_\alpha}
\]
is a linear isomorphism.
For a finite exponent set \(I\subseteq\Z_{\ge0}\), let the \(i\)th column of \(R_I\) be the power-basis coordinate vector of \(X^i\bmod f\), and put
\begin{equation}\label{v97:eq-hasse-coordinates}
 V_I=\bigl(\partial^{[j]}X^i\big|_{X=\alpha}\bigr)_{(\alpha,j),\,i\in I}.
\end{equation}
There is an invertible matrix \(J\in\mathrm{GL}_d(E)\), independent of \(I\), such that \(V_I=JR_I\).
The two matrices therefore describe the same divisibility relations, and \(\rank_E V_I=\rank_k R_I\).
\end{proposition}
\begin{proof}
If all the indicated Hasse derivatives of a polynomial of degree less than \(d\) vanish, the multiplicity criterion makes it divisible by \(f\), hence zero.
Both spaces have dimension \(d\), so the map is an isomorphism; denote its matrix by \(J\).
Subtracting the remainder modulo \(f\) from any polynomial changes none of the indicated values or derivatives.
Column by column, this gives \(V_I=JR_I\). Invertible left multiplication and scalar extension preserve rank.
\end{proof}

In \eqref{v97:eq-hasse-coordinates}, the entry is \(\binom ij\alpha^{i-j}\) when \(\alpha\ne0\) and \(i\ge j\), and is zero when \(i<j\).
At \(\alpha=0\), the entry is \(\delta_{ij}\), so the multiplicity of a zero root is retained as well.
Ordinary root evaluation records only \(j=0\) and cannot replace these coordinates in the presence of repeated roots.

We next collect the same bivariate difference in powers of \(X\) to obtain all conditions for translation invariance.

\subsection{The Translation Stabilizer Polynomial}

\begin{definition}\label{ffs:def-translation-polynomial}
Write
\[
f(X+B)-f(X)=\sum_{j=0}^{n-1}D_{f,j}(B)X^j,
\qquad
D_{f,j}(B)=\sum_{i=j+1}^{n}\binom ij c_iB^{i-j}.
\]
These coefficients are not all zero, since \(D_{f,0}(B)=f(B)-f(0)\) is nonconstant.
Their monic greatest common divisor
\[
\Psi_f(B)=\gcd\{D_{f,j}(B):D_{f,j}\ne0\}
\]
is called the \textbf{translation stabilizer polynomial} of \(f\).
\end{definition}

The set \(V_f\) is an \(\F_p\)-vector space, but need not be a \(k\)-vector space.
Applying Theorem~\ref{thm:subspace-polynomial} over the prime field \(\F_p\) gives
\[
L_{V_f}(X)=\prod_{v\in V_f}(X-v)
 =X^{p^s}+a_{s-1}X^{p^{s-1}}+\cdots+a_0X,
\qquad a_0\ne0,
\]
where \(s=\dim_{\F_p}V_f\).
Since \(\sigma(V_f)=V_f\), uniqueness of the monic root polynomial gives \(L_{V_f}\in k[X]\).
It satisfies \(L_{V_f}(X+Y)=L_{V_f}(X)+L_{V_f}(Y)\).
When \(V_f=\{0\}\), we use \(L_{V_f}=X\) and \(|V_f|=1\).

\begin{theorem}[Coefficient recovery formula for translations]\leavevmode\phantomsection\label{ffs:thm-psi}
For any monic nonconstant \(f\in k[X]\), let \(e\) be its common \(p\)-power from Proposition~\ref{ffs:prop-common-power}. Then
\begin{equation}\label{ffs:eq-psi}
\boxed{\ \Psi_f(B)=L_{V_f}(B)^e.\ }
\end{equation}
Thus the distinct roots of \(\Psi_f\) are exactly \(V_f\), each with multiplicity \(e\) in \(\Psi_f\).
\end{theorem}
\begin{proof}
First suppose \(f'\ne0\). Any common root \(b\), possibly in a larger finite extension, satisfies \(f(X+b)=f(X)\).
Fix \(\alpha\in S_f\). Then \(\alpha+b\in S_f\), so \(b=(\alpha+b)-\alpha\) already lies in the splitting field of \(f\).
Hence the common roots of all \(D_{f,j}\) are exactly \(V_f\).
Choose \(i\ge1\) with \(ic_i\ne0\). Then
\[
D'_{f,i-1}(0)=ic_i\ne0.
\]
All difference coefficients vanish at zero, and at least one has only a simple root there, so their greatest common divisor has multiplicity one at zero.
For any \(b\in V_f\), we have
\[
f(X+B+b)-f(X)=f(X+B)-f(X).
\]
Comparing coefficients of each power of \(X\) gives \(D_{f,j}(B+b)=D_{f,j}(B)\).
Every common root therefore has the same multiplicity as zero, and the monic greatest common divisor is exactly \(L_{V_f}\).

In general, write \(f=h^e\), with \(h'\ne0\). The Frobenius identity gives
\[
f(X+B)-f(X)=\bigl(h(X+B)-h(X)\bigr)^e.
\]
The nonzero coefficients of powers of \(X\) on the left are exactly the \(e\)th powers of the difference coefficients on the right, so the greatest common divisor is also raised to the \(e\)th power.
Using \(V_h=V_f\) now yields the formula.
\end{proof}

\begin{corollary}\label{ffs:cor-psi-recovery}
The coefficients of \(f\) recover
\[
e=\operatorname{ord}_{B=0}\Psi_f,
\qquad |V_f|=\frac{\deg\Psi_f}{e},
\qquad L_{V_f}=\Psi_f^{\,1/e}.
\]
Here \(\operatorname{ord}_{B=0}\Psi_f\) is the multiplicity of zero as a root of \(\Psi_f\); the final expression denotes the unique monic polynomial \(G\) satisfying \(G^e=\Psi_f\).
\end{corollary}

This recovers the polynomial describing the actual subspace, not just the size of the stabilizer.
Even when \(V_f=\{0\}\), we have \(\Psi_f=B^e\), from which the common \(p\)-power \(e\) can still be recovered.

\subsection{Recovering additive decompositions from constant differences}
Translation invariance requires a zero difference. Allowing a constant difference reveals a larger additive structure. Here translations are restricted to the base field, unlike the root stabilizers in the splitting field considered above.
\begin{proposition}[The maximal additive translation space]\label{v109:prop-additive-directions}
Let $P\in\F_q[X]$, $p=\Char\F_q$, and $P_0=P-P(0)$. Define
\[
\mathcal A_P=\{v\in\F_q:P_0(X+v)-P_0(X)=P_0(v)\text{ in }\F_q[X]\}.
\]
Then $\mathcal A_P$ is an $\F_p$-subspace, and $P_0|_{\mathcal A_P}$ is $\F_p$-linear.
Put $L=L_{\mathcal A_P}$. There are unique $H\in\F_q[X]$ and a $p$-linearized polynomial $M$ such that
\[
P(X)=P(0)+H(L(X))+M(X),\qquad H(0)=0,\quad\deg M<\deg L.
\]
The zero polynomial $M=0$ is allowed. If a decomposition of this kind uses the vanishing polynomial of a base-field subspace $U$, then $U\subseteq\mathcal A_P$.
\end{proposition}
\begin{proof}
Apply two translation identities successively and evaluate at $X=0$. This proves additive closure and $P_0(u+v)=P_0(u)+P_0(v)$. A finite additive subgroup is an $\F_p$-space.
Let $d=\dim_{\F_p}\mathcal A_P$. Moore interpolation gives a unique
$M=\sum_{i=0}^{d-1}m_iX^{p^i}$ with $M(v)=P_0(v)$ for $v\in\mathcal A_P$; take $M=0$ if $d=0$.
Thus $R=P_0-M$ is invariant under these translations.
Repeated monic division uniquely writes $R=\sum_i R_i(X)L(X)^i$, where $\deg R_i<\deg L$.
Since $L(X+v)=L(X)$, uniqueness of the expansion gives $R_i(X+v)=R_i(X)$. Therefore $R_i-R_i(0)$ has $\deg L$ roots in $\mathcal A_P$ and must vanish. Hence $R=H(L)$, and $R(0)=0$ gives $H(0)=0$.
Any such decomposition restricts to $M$ on $\mathcal A_P$. Moore interpolation proves uniqueness of $M$, and then of $H$.
If the same type of decomposition uses $U$, taking differences for $v\in U$ gives $P_0(X+v)-P_0(X)=M(v)=P_0(v)$, so $U\subseteq\mathcal A_P$.
\end{proof}
If $P_0(X+B)-P_0(X)-P_0(B)=\sum_{j\ge1}c_j(B)X^j$, then $L_{\mathcal A_P}(B)=\gcd(B^q-B,c_1(B),c_2(B),\ldots)$, with monic gcd. Both sides have exactly the same base-field roots, and $B^q-B$ is squarefree.

This is the additive-decomposition viewpoint of \cite[Theorem 3.2]{ReisWangAdditive2021}. Base-field translations preserving $P$ itself form
$\ker(P_0|_{\mathcal A_P})$. A nonzero kernel prevents injectivity along those directions. For the connection with the AGW criterion, see Section~\ref{sec:agw}.
\begin{example}[Constant differences and zero differences]\label{v109:ex-additive-directions}
Over $\F_{p^n}$ of odd characteristic, take $P=(X^p-X)^2+X$. Direct calculation gives
\[
P(X+v)-P(X)-P(v)=2(v^p-v)(X^p-X).
\]
Thus $\mathcal A_P=\F_p$, with canonical decomposition $L=X^p-X$, $H=X^2$, $M=X$.
However, $P(X+v)-P(X)=v$, so only the zero base-field translation preserves $P$ itself.
Polynomial identities cannot be replaced by function identities: $P=X(X^q-X)$ is the zero function on $\F_q$, but its difference is $v(X^q-X)$, giving $\mathcal A_P=\{0\}$.
\end{example}
\subsection{Coefficient Equations for Scalings and the Full Affine Stabilizer}

\begin{proposition}[Exponent criterion for scalings]\label{ffs:prop-multiplier-support}
Suppose \(f\) is not a monomial. Put \(m=\nu_f(0)=\min\{i:c_i\ne0\}\), and define
\[
s_0=\gcd\{i-m:i>m,\ c_i\ne0\},\qquad
s_*=\frac{s_0}{p^{v_p(s_0)}}.
\]
Then \(M_f=\mu_{s_*}\), so \(|M_f|=s_*\).
There is also a unique \(G\in k[Y]\) such that
\[
f(X)=X^mG(X^{s_0}),\qquad G(0)\ne0.
\]
\end{proposition}
\begin{proof}
Write \(f=X^m g\), where \(g(0)\ne0\). The affine identity gives
\(g(aX)=a^{n-m}g(X)\). Comparing constant terms first gives \(a^{n-m}=1\); comparing the remaining coefficients then gives exactly the conditions
\(a^{i-m}=1\) for all \(c_i\ne0\).
These are equivalent to \(a^{s_0}=1\). In characteristic \(p\), the distinct roots of \(X^{s_0}-1\) are exactly \(\mu_{s_*}\).
The compositional expression follows directly because all exponent differences are multiples of \(s_0\); uniqueness follows from injectivity of nonconstant substitution.
\end{proof}

\begin{definition}[Support period and support-primitive core]\label{v44:def-support-period}
Let \(f\in k[X]\) be monic and nonconstant with \(f(0)\ne0\). Its \textbf{support period} is
\[
\operatorname{sp}(f)=\gcd\{j\ge1:[X^j]f\ne0\}.
\]
\leavevmode\phantomsection\label{v44:def-support-primitive}\ignorespaces
If \(\operatorname{sp}(f)=1\), then \(f\) is called a \textbf{support-primitive polynomial}.
\leavevmode\phantomsection\label{v44:def-support-core}\ignorespaces
Put \(s=\operatorname{sp}(f)\). There is a unique monic support-primitive polynomial \(f_{\mathrm{sp}}\in k[Y]\) such that
\[
f(X)=f_{\mathrm{sp}}(X^s).
\]
The polynomial \(f_{\mathrm{sp}}\) is called the \textbf{support-primitive core} of \(f\).
\end{definition}

Every positive exponent with nonzero coefficient is divisible by \(s\), so extracting the coefficients and dividing their exponents by \(s\) gives \(f_{\mathrm{sp}}\);
the positive exponents in the resulting support have greatest common divisor one, and injectivity of the substitution \(Y\mapsto X^s\) ensures uniqueness.
Thus the support period is the maximum monomial composition exponent defined on page~\pageref{def:coefficient-extraction}, and
\(\deg f_{\mathrm{sp}}=\deg f/s\).
Proposition~\ref{ffs:prop-common-power} and the preceding proposition together describe it through root multiplicities; see Remark~\ref{ffs:rem-index-addition} for its expression in root-of-unity indices.
For this terminology and its applications to decompositions of cyclic codes, see~\cite{ZWSupportPrimitive}.
Support-primitivity requires only that the exponents in the support have greatest common divisor one; it requires neither irreducibility nor that the roots generate the multiplicative group.

When \(f(0)\ne0\), we have \(m=0\), and Proposition~\ref{ffs:prop-common-power} gives \(p^{v_p(s_0)}=e\), so \(s_0=e|M_f|\).
When zero is a root, the equality \(s_0=e|M_f|\) need not hold.
For example, over \(\F_3\), the polynomial \(f=X(X^3-1)\) has \(m=1,s_0=3,e=1\), but \(M_f=\{1\}\).
The multiplicity of zero and the common \(p\)-power of the other root multiplicities must be recorded separately.

The support period also records a scalar action that may remain on the roots of a binomial.
It can turn coefficient information about a single factor into a splitting condition for the entire binomial.
\begin{proposition}[Coefficient support and binomial splitting]\label{v97:prop-support-splitting}
Let \(n\ge2\), \(p\nmid n\), and \(\lambda\in\F_q^\times\), and put
\[
 m=\ord_n(q),\qquad c=\lambda^{(q^m-1)/n}\in\F_q^\times.
\]
Every monic nonconstant polynomial \(f\in\F_q[X]\) dividing \(X^n-\lambda\) satisfies \(c^{\operatorname{sp}(f)}=1\).
In particular, if one factor is support-primitive, the entire binomial splits over \(\F_{q^m}\).
More generally, if \(e=\ord(c)\), its splitting field is exactly
\begin{equation}\label{v97:eq-binomial-splitting}
 \boxed{\F_{q^{me}}.}
\end{equation}
\end{proposition}
\begin{proof}
Every root \(\alpha\) of the binomial satisfies
\[
 \alpha^{q^m}=\alpha(\alpha^n)^{(q^m-1)/n}=c\alpha.
\]
The root set of \(f\) is stable under the \(q^m\)-Frobenius, and hence under multiplication by \(c\).
Since the binomial has distinct roots, \(f(cX)=c^{\deg f}f(X)\).
Comparing nonzero constant terms gives \(c^{\deg f}=1\); comparison at each nonzero coefficient then gives \(c^i=1\) for every positive support exponent \(i\).
Taking their greatest common divisor yields \(c^{\operatorname{sp}(f)}=1\).

A field containing all binomial roots contains all their ratios, hence all \(n\)th roots of unity.
Its degree over \(\F_q\) is therefore a multiple of \(m\).
Since \(c\in\F_q\), iteration gives \(\alpha^{q^{mt}}=c^t\alpha\).
This fixes every root precisely when \(e\mid t\), proving the exact splitting degree \(me\).
For a support-primitive factor, \(c=1\), so \(e=1\).
\end{proof}

Having every coefficient nonzero is sufficient for support-primitivity, but is not necessary.
What matters is whether the support exponents permit a common nontrivial scaling, not how many coefficients are nonzero.
Formula~\eqref{v97:eq-binomial-splitting} also separates splitting into two steps: adjoining the roots of unity, then removing the remaining scalar action.

\begin{example}[A support period retains an additional splitting degree]\label{v97:ex-binomial-support}
In \(\F_3[X]\),
\[
 X^8+1=(X^4+X^2-1)(X^4-X^2-1).
\]
Every root has order \(16\), and \(\ord_{16}(3)=4\), so both quartic factors are irreducible.
Here \(m=\ord_8(3)=2\), \(c=(-1)^{(9-1)/8}=-1\), and \(e=2\).
The splitting field is therefore \(\F_{81}\), not \(\F_9\).
Both factors have support period \(2\), allowing the scaling \(c=-1\).
Proposition~\ref{v97:prop-norm-factor-degrees} only forces even factor degrees; the multiplicative root order and the remaining scalar action determine degree four in this case.
\end{example}

\begin{proposition}[Coefficient criterion for the full affine stabilizer]\leavevmode\phantomsection\label{ffs:prop-coefficient-equations}
We have \(g_{a,b}\in\mathcal A_f\) if and only if \(a\ne0\) and
\begin{equation}\label{ffs:eq-affine-coefficients}
a^j\sum_{i=j}^{n}\binom ijc_i b^{i-j}=a^nc_j
\qquad(0\le j<n).
\end{equation}
\end{proposition}
\begin{proof}
Expand the identity in Proposition~\ref{ffs:prop-affine-identity} and compare coefficients.
The leading-term equality holds automatically.
\end{proof}

The same coefficient identities can be tested over a ring containing a nonzero element whose square is zero.
This connects affine symmetries with the dependence of the algebra trace on root multiplicities in Corollary~\ref{prop:local-multiplicity-trace}.

\begin{definition}[First-order affine symmetries]\leavevmode\phantomsection\label{v67:def-first-order-affine}
Take the ring of dual numbers \(D=k[\varepsilon]/(\varepsilon^2)\) (see Proposition~\ref{v85:prop-quadratic-unit-groups}), and let \(f\in k[X]\) be monic of degree \(n\ge1\).
Define the \textbf{first-order affine symmetry space} of \(f\) by
\[
I_f=\left\{(u,v)\in k^2:
f\bigl((1+\varepsilon u)X+\varepsilon v\bigr)
=(1+\varepsilon u)^n f(X)\text{ in }D[X]\right\}.
\]
Here \(1+\varepsilon u\) is a unit, with inverse \(1-\varepsilon u\), so each substitution is invertible and reduces to the identity modulo \((\varepsilon)\).
\end{definition}

\begin{proposition}[First-order affine symmetries and the rank of the trace pairing]\leavevmode\phantomsection\label{v67:prop-first-order-affine}
With the notation above,
\begin{equation}\label{v67:eq-first-order-affine}
I_f=\{(u,v)\in k^2:(uX+v)f'(X)=nu\,f(X)\}.
\end{equation}
In particular, \(I_f\) is a \(k\)-linear space.
Let \(A=k[X]/(f)\), and let \(r\) be the rank of the algebra trace pairing
\(B_{\Tr}(x,y)=\Tr_{A/k}(xy)\). Then
\begin{equation}\label{v67:eq-first-order-trace-rank}
\boxed{\quad\dim_k I_f=\max\{0,2-r\}.\quad}
\end{equation}
\end{proposition}
\begin{proof}
Since \(\varepsilon^2=0\), expansion gives
\[
f\bigl((1+\varepsilon u)X+\varepsilon v\bigr)
=f(X)+\varepsilon(uX+v)f'(X),\qquad
(1+\varepsilon u)^n=1+\varepsilon nu.
\]
Comparison of the coefficients of \(\varepsilon\) proves~\eqref{v67:eq-first-order-affine}; its coefficient equations are linear in \(u,v\).

Extend coefficients to the splitting field \(E\), and write \(m_\alpha=\nu_f(\alpha)\) for \(\alpha\in S_f\).
The logarithmic derivative and the equality \(n=\sum_{\alpha\in S_f}m_\alpha\) give
\[
(uX+v)\frac{f'(X)}{f(X)}-nu
=\sum_{\alpha\in S_f}\frac{m_\alpha(u\alpha+v)}{X-\alpha},
\]
where the integer multiplicities are interpreted in \(E\).
The rational function on the right vanishes if and only if every numerator vanishes:
multiply by \(\prod_{\alpha\in S_f}(X-\alpha)\) and evaluate successively at its distinct roots.
Thus the nonzero constraints are exactly
\[
u\alpha+v=0\qquad(\alpha\in S_f,\ p\nmid m_\alpha).
\]
Corollary~\ref{prop:local-multiplicity-trace} identifies the number of these roots with \(r\).
If there are none, there is no constraint; one such root gives one independent equation; two distinct such roots give two independent equations, forcing \(u=v=0\).
The dimension over \(E\) is therefore \(\max\{0,2-r\}\).
Extension of scalars preserves the rank of a linear system over \(k\), so this is also \(\dim_k I_f\).
\end{proof}

\begin{example}[A trivial affine stabilizer with a nonzero first-order symmetry]\leavevmode\phantomsection\label{v67:ex-first-order-affine}
Over \(\F_5\), take
\[
h=(X-1)(X-2)^5=X^6-X^5-2X+2.
\]
The two roots have different multiplicities, so every affine transformation preserving these multiplicities fixes both roots and is the identity.
Nevertheless,
\[
h'=X^5-2,\qquad (X-1)h'=h.
\]
As \(6=1\) in \(\F_5\), equation~\eqref{v67:eq-first-order-affine} becomes
\((u+v)h'=0\), and hence
\[
I_h=\{(u,-u):u\in\F_5\},\qquad \dim_{\F_5}I_h=1.
\]
For example, over the dual numbers,
\[
h\bigl((1+\varepsilon)X-\varepsilon\bigr)
=(1+\varepsilon)^6h(X).
\]
Only the root \(1\) has multiplicity not divisible by five, so the trace pairing of \(\F_5[X]/(h)\) has rank one, in agreement with~\eqref{v67:eq-first-order-trace-rank}.
In the centered coordinate \(Y=X-1\), the polynomial is \(h(Y+1)=Y^6-Y\).
\end{example}

First-order information does not determine the common \(p\)-power or the coefficient identities over rings with higher-order nilpotents:
if \(g\) is monic and nonconstant with \(g'\ne0\), then \(g^p\) and \(g^{p^2}\) both have zero derivative and degree divisible by \(p\), so
\(I_{g^p}=I_{g^{p^2}}=k^2\), although their common \(p\)-powers are \(p\) and \(p^2\), respectively.

\section{Translation Orbits and Additive Decompositions}
\label{ffs:sec-additive-composition}

\begin{theorem}[Translation invariance and subspace composition]\label{ffs:thm-additive-composition}
Let \(W\) be a finite Frobenius-stable \(\F_p\)-subspace, and let \(f\in k[X]\). Then
\[
f(X+w)=f(X)\text{ for all }w\in W
\quad\Longleftrightarrow\quad
f(X)=F(L_W(X))\text{ for some }F\in k[Y].
\]
The outer polynomial \(F\) is unique; if \(f\) is monic and nonconstant, so is \(F\), and \(\deg f=|W|\deg F\).
\end{theorem}
\begin{proof}
Sufficiency follows immediately from \(L_W(X+w)=L_W(X)\).
For necessity, take a finite extension containing \(W\).
Since \(f(w)=f(0)\) for every \(w\in W\), the monic square-free polynomial \(L_W\) divides \(f-f(0)\), so
\[
f(X)=f(0)+L_W(X)Q(X).
\]
For every \(w\in W\), compare this identity with its translate and cancel the nonzero polynomial \(L_W\). It follows that \(Q\) is also translation-invariant.
Induction on the degree gives \(f=F(L_W)\). Nonconstant substitution \(Y\mapsto L_W(X)\) maps nonzero polynomials to nonzero polynomials, so is injective and ensures uniqueness.

Since \(L_W\in k[X]\), applying coefficientwise Frobenius to this identity gives another decomposition with the same inner polynomial;
uniqueness forces every coefficient of \(F\) to be Frobenius-fixed, so \(F\in k[Y]\).
\end{proof}

On root sets, the same conclusion appears as a product over complete additive cosets. For every \(\alpha\),
\[
\prod_{w\in W}(X-\alpha-w)=L_W(X)-L_W(\alpha).
\]
If \(W\subseteq V_f\), the multiplicities on each root coset are equal, and the multiplicity of \(L_W(\alpha)\) as a root of the outer polynomial \(F\) equals this common multiplicity.
This compositional characterization of translation-invariant polynomials also appears in~\cite{FFSReis2016}.

\begin{corollary}[The maximal translation inner component and its outer component]\label{ffs:cor-maximal-additive}
Put \(V=V_f\), and write \(f=h^e\), with \(h'\ne0\). There is a unique monic \(F\in k[Y]\) such that
\[
h=F(L_V),\qquad F'\ne0.
\]
Moreover, \(F\) is invariant under no nonzero translation: \(V_F=\{0\}\) in every finite splitting extension.
For every finite Frobenius-stable \(\F_p\)-subspace \(W\), we have
\[
f\in k[L_W]\quad\Longleftrightarrow\quad W\subseteq V_f.
\]
Thus \(L_V\) has maximal degree among such additive inner components.
\end{corollary}
\begin{proof}
The decomposition and the final equivalence follow from Theorem~\ref{ffs:thm-additive-composition}.
Since \(L_V'\) is a nonzero constant, the chain rule gives \(h'=F'(L_V)L_V'\), hence \(F'\ne0\).
If \(F(Y+d)=F(Y)\), take \(b\) in a finite extension with \(L_V(b)=d\). Then
\[
h(X+b)=F(L_V(X)+d)=h(X).
\]
Thus \(b\in V_h=V\), and hence \(d=L_V(b)=0\).
\end{proof}

\leavevmode\phantomsection\label{v98:par-exact-translation}\ignorespaces
To start from a prescribed translation space \(W\), first require \(W\) to be Frobenius-stable so that \(L_W\in k[X]\), and then choose a nonconstant outer polynomial \(G\in k[Y]\).
The composite \(f=G(L_W)\) guarantees \(W\subseteq V_f\). For the full translation stabilizer to be exactly \(W\), nonzero translation parameters of \(G\) must also be excluded, as in the preceding proof's examination of additional translations.
The absence of nonzero translation symmetries of the outer polynomial does not exclude its other functional decompositions; the full affine structure considered below also requires studying slopes beyond translations.

\section{The Full Affine Stabilizer, Center Coset, and Normal Form}
\label{ffs:sec-full-affine}

\leavevmode\phantomsection\label{ffs:def-affine-data}\ignorespaces Henceforth assume \(|S_f|\ge2\), and write \(V=V_f\) and \(Q=|V|\).
The slope map \(g_{a,b}\mapsto a\) is a group homomorphism \(\mathcal A_f\to E^\times\) with kernel exactly the translation group \(V\).
Its image is cyclic of order coprime to \(p\); denote it by \(\mu_t\).
Here \(t\) need not equal \(|M_f|\): \(M_f\) consists only of scalings about zero, whereas the slope image also allows a nonzero translation part.

\begin{theorem}[Structure of a finite affine group]\leavevmode\phantomsection\label{ffs:thm-affine-group}
If \(t=1\), then \(\mathcal A_f\) is the translation group \(V\).
If \(t>1\), there exists \(c\in E\) such that
\begin{equation}\label{ffs:eq-centered-group}
\mathcal A_f
 =\{X\mapsto c+a(X-c)+v:\ a\in\mu_t,\ v\in V\}.
\end{equation}
Here \(aV=V\) and \(L_V(aX)=aL_V(X)\) for every \(a\in\mu_t\).
In particular, \(|\mathcal A_f|=Qt\); if \(Q>1\), then \(t\mid Q-1\).
\end{theorem}
\begin{proof}
When \(t=1\), every element has slope one, proving the assertion.
Suppose \(t>1\). Take a generator \(a\) of \(\mu_t\), and an element \(g(X)=aX+b\) of slope \(a\).
\leavevmode\phantomsection\label{ffs:def-affine-center}\ignorespaces Its unique fixed point \(c=b/(1-a)\) is called the \textbf{center of the affine transformation}; in the coordinate \(Y=X-c\), the transformation \(g\) becomes the scaling \(Y\mapsto aY\), and
\[
g^j(X)=c+a^j(X-c),\qquad g^t=\id.
\]
Every group element has the same slope as some \(g^j\) and differs from it by a translation, giving~\eqref{ffs:eq-centered-group}.

Conjugating the translation \(X\mapsto X+v\) gives \(X\mapsto X+av\), so \(aV=V\).
The root product gives \(L_V(aX)=a^QL_V(X)\).
The coefficient of \(X\) in \(L_V\) is nonzero; comparing linear terms gives \(a^Q=a\), hence \(L_V(aX)=aL_V(X)\).
When \(Q>1\), we have \(a^{Q-1}=1\), so the order \(t\) of the generator divides \(Q-1\).
\end{proof}

In the coordinate \(Y=X-c\), each group element is uniquely written as \(Y\mapsto aY+v\), and the group operation is
\[
(v,a)(v',a')=(v+av',aa').
\]
\leavevmode\phantomsection\label{def:semidirect-affine}\ignorespaces This group structure, determined by the scalar action \(a:v\mapsto av\), is called the \textbf{semidirect product} of the translation group \(V\) and the cyclic slope group \(\mu_t\), denoted by \(V\rtimes\mu_t\).
The translation group \(V\) is normal, its intersection with \(\{(0,a):a\in\mu_t\}\) is the identity, and every group element is uniquely a product of elements from these two subgroups.

\begin{proposition}[Special and free orbits]\leavevmode\phantomsection\label{ffs:prop-affine-orbits}
Suppose \(t>1\). The coset \(C=c+V\) is called the \textbf{center coset}; as an orbit under \(\mathcal A_f\), it is called the \textbf{special orbit}.
It is independent of the choice of center \(c\) and consists exactly of the points of \(E\) with nontrivial point stabilizers.
The set \(C\) is one orbit of size \(Q\); every other orbit has size \(Qt\), with free action.
Put \(U(X)=L_V(X-c)\). For \(\alpha\notin C\), its orbit polynomial is
\[
\prod_{\beta\in\mathcal A_f\alpha}(X-\beta)
 =U(X)^t-U(\alpha)^t.
\]
\end{proposition}
\begin{proof}
In centered coordinates, \(C\) is \(V\), on which the translation group acts transitively.
If \(ay+v=y\) with \(a\ne1\), then \((1-a)y=v\).
We have \((1-a)V\subseteq V\), and multiplication by \(1-a\ne0\) is injective on the finite space \(V\), hence is also surjective on \(V\).
Thus \(y\in V\). Conversely, every \(y\in V\) is fixed by \(z\mapsto az+(1-a)y\).
Nontrivial point stabilizers therefore occur exactly on \(C\), which also proves that this coset is independent of the chosen center.
The action outside \(C\) is free, so each orbit has \(Qt\) points.

Taking products first along translation orbits and then over \(\mu_t\) gives
\[
\prod_{a\in\mu_t}\bigl(U(X)-aU(\alpha)\bigr)
 =U(X)^t-U(\alpha)^t.
\]
Here \(U(\alpha)\ne0\), so the roots are distinct and indeed form the entire free orbit.
\end{proof}

\noindent In characteristic two, the same distinction can be seen through derivatives.
When the translation part is trivial, a root set consisting only of nonzero orbits of a multiplicative group of order \(d\ge5\)
has root polynomial \(P(X^d)\), whose derivative has a square root divisible by \(X^2\).
Including the fixed point zero may remove this obstruction.
Proposition~\ref{v83:prop-two-mult-cosets} gives a completely split construction, illustrating that fixed points can affect not only orbit sizes but also repeated-root behavior of derivatives.

\begin{theorem}[Descent of the center coset to the base field]\label{ffs:thm-center-descent}
When \(t>1\), Frobenius preserves \(C=c+V\), and
\[
c^q-c\in V,\qquad \delta:=L_V(c)\in k.
\]
Thus \(\delta\) and
\begin{equation}\label{ffs:eq-centered-polynomial}
U(X)=L_V(X)-\delta=\prod_{\gamma\in C}(X-\gamma)\in k[X]
\end{equation}
are uniquely determined by \(\mathcal A_f\), although a center \(c\in k\) need not exist.
\end{theorem}
\begin{proof}
Proposition~\ref{ffs:prop-frobenius-normalizes} shows that Frobenius conjugates nontrivial point stabilizers to nontrivial point stabilizers.
Proposition~\ref{ffs:prop-affine-orbits} therefore gives \(\sigma(C)=C\), that is, \(c^q-c\in V\).
Since \(L_V\in k[X]\) and is additive,
\[
L_V(c)^q=L_V(c^q)=L_V(c)+L_V(c^q-c)=L_V(c).
\]
Hence \(\delta\in k\). Changing the center adds an element of \(V\) to \(c\), leaving \(L_V(c)\) unchanged.
The root-product expression follows from the translation formula for the subspace polynomial.
\end{proof}

\begin{corollary}[Scalings about the origin and the full slope group]\label{v54:cor-origin-slopes}
Suppose \(t>1\). The stabilizer \(M_f\) of scalings about the origin satisfies
\[
M_f=\begin{cases}
\mu_t,&0\in C\quad(\delta=0),\\
\{1\},&0\notin C\quad(\delta\ne0).
\end{cases}
\]
\end{corollary}
\begin{proof}
By~\eqref{ffs:eq-centered-group}, the scaling \(aX\) preserves root multiplicities if and only if
\(a\in\mu_t\) and \((1-a)c\in V\). If \(c\in V\), this holds for every slope.
If \(a\ne1\), the map \(v\mapsto(1-a)v\) maps the finite space \(V\) into itself injectively, hence bijectively;
thus \((1-a)c\in V\) conversely forces \(c\in V\). This proves both cases.
Finally, \(0\in c+V\iff c\in V\iff L_V(c)=0\), that is, \(\delta=0\).
\end{proof}

\begin{corollary}[Recovering the base-field stabilizer]\label{ffs:cor-rational-group}
Suppose \(t>1\), and put \(V(k)=V\cap k\). The center coset contains a base-field point if and only if \(\delta\in L_V(k)\).
If there is a base-field center \(c_0\), then
\[
\mathcal A_f(k)
 =\{X\mapsto c_0+a(X-c_0)+v:
 a\in\mu_t\cap k^\times,\ v\in V(k)\}.
\]
If there is no base-field center, \(\mathcal A_f(k)\) consists only of translations by elements of \(V(k)\). Consequently,
\[
|\mathcal A_f(k)|=
\begin{cases}
|V(k)|\gcd(t,q-1),&\delta\in L_V(k),\\
|V(k)|,&\delta\notin L_V(k).
\end{cases}
\]
\end{corollary}
\begin{proof}
The center coset is defined by \(L_V(X)=\delta\), so the first assertion directly expresses solvability over \(k\).
Choosing \(c_0\in k\), both coefficients of the transformation in~\eqref{ffs:eq-centered-group} lie in \(k\) exactly when \(a,v\in k\), yielding the first group expression.
If a base-field affine transformation has slope \(a\ne1\), its fixed point \(b/(1-a)\) lies both in \(k\) and in \(C\).
Such a slope therefore cannot occur when no base-field center exists.
Finally use \(|\mu_t\cap k^\times|=\gcd(t,q-1)\) in the multiplicative group of a finite field.
\end{proof}

\begin{theorem}[Affine normal form preserving root multiplicities]\label{ffs:thm-normal-form}
Suppose \(|S_f|\ge2\), write \(f=h^e\), with \(h'\ne0\), and take the full stabilizer data \(V,Q,t\).
\begin{enumerate}
\item If \(t=1\), then \(h=F(L_V)\), where \(F\in k[Y]\) is monic and \(F'\ne0\).
\item If \(t>1\), let \(U=L_V-\delta\) be the polynomial in~\eqref{ffs:eq-centered-polynomial}.
There are a unique integer \(w\ge0\) and a unique monic \(H\in k[Y]\), with \(H(0)\ne0\), such that
\begin{equation}\label{ffs:eq-normal-form}
\boxed{\quad f(X)=\bigl(U(X)^wH(U(X)^t)\bigr)^e.\quad}
\end{equation}
Here \(w=\nu_h(c)\). Every point of the special orbit \(C\) has root multiplicity \(ew\), with \(w=0\) meaning that \(C\) is disjoint from the root set.
For a nonzero root \(\eta\) of \(H\), the fiber \(U(X)^t=\eta\) is a free orbit, and each of its roots has multiplicity \(e\nu_H(\eta)\).
Moreover,
\[
\deg f=eQ\bigl(w+t\deg H\bigr),
\qquad p\nmid w\ \text{or}\ H'\ne0.
\]
\end{enumerate}
Conversely, given \(\mu_t\) satisfying \(aV=V\) and a center \(c\), every polynomial of the form~\eqref{ffs:eq-normal-form} has its root multiplicities preserved by the group in~\eqref{ffs:eq-centered-group}.
This group may be only a subgroup of its full stabilizer.
\end{theorem}
\begin{proof}
The first assertion was proved in Corollary~\ref{ffs:cor-maximal-additive}.
Suppose \(t>1\). The polynomials \(h\) and \(f\) have the same stabilizer, and \(\nu_h\) is constant on each affine orbit.
The root polynomial of the special orbit is \(U\), and those of the other orbits are given by Proposition~\ref{ffs:prop-affine-orbits}.
Multiplying over the orbits gives \(h=U^wH(U^t)\), with \(H(0)\ne0\).
The factor \(U^w\) records the zero fiber exactly, so \(w\) is unique; injectivity of nonconstant substitution then gives uniqueness of \(H\).

Since \(h,U\in k[X]\), coefficientwise Frobenius gives another expression of the same form; uniqueness implies \(H\in k[Y]\).
The derivative \(U'\) is a nonzero constant and \(p\nmid t\), so \(U\) is square-free, as is \(U^t-\eta\) when \(\eta\ne0\).
This gives all the multiplicity and degree formulas.
If \(p\mid w\), the chain rule shows \(h'\ne0\iff H'\ne0\);
if \(p\nmid w\), then \(H(0)\ne0\) ensures that the derivative of the lowest-degree term in \(U\) does not vanish.
Thus the final condition expresses exactly \(h'\ne0\).

Conversely, for a group element, \(U(c+a(X-c)+v)=aU(X)\), and \(a^t=1\).
The normal form is therefore multiplied by a nonzero constant, preserving root multiplicities.
However, some choices of the outer polynomial may introduce further symmetries, so the converse guarantees only containment.
\end{proof}

\begin{corollary}[The square-free case]\label{ffs:cor-squarefree-normal}
When \(t>1\), the polynomial \(f\) is square-free if and only if \(e=1\), \(w\in\{0,1\}\), and \(H\) is square-free.
\end{corollary}
\begin{proof}
Distinct fibers in the normal form are disjoint and every fiber polynomial is square-free; the assertion follows pointwise from the multiplicity formulas in the theorem.
\end{proof}

\leavevmode\phantomsection\label{def:orbit-quotient-free}\ignorespaces
For an action of a group \(G\) on a set \(S\), regard each orbit \(Gx=\{g(x):g\in G\}\) as one element. The resulting set
\(S/G=\{Gx:x\in S\}\) is called the \textbf{orbit space}.
This is a quotient of sets; it does not automatically inherit any field operations the original set may possess.

\begin{proposition}[Orbit space and invariant polynomials]\leavevmode\phantomsection\label{ffs:prop-invariant-ring}
Suppose \(t>1\), and put \(\pi(X)=U(X)^t\). Then
\[
\pi(\alpha)=\pi(\beta)
\quad\Longleftrightarrow\quad
\beta\in\mathcal A_f\alpha.
\]
Furthermore, if
\[
k[X]^{\mathcal A_f}
 =\{P\in k[X]:P(g(X))=P(X)\text{ for all }g\in\mathcal A_f\},
\]
then \(k[X]^{\mathcal A_f}=k[\pi(X)]\).
For \(f\) in Theorem~\ref{ffs:thm-normal-form}, we have
\[
f\in k[\pi]\quad\Longleftrightarrow\quad t\mid w.
\]
\end{proposition}
\begin{proof}
The zero fiber is exactly the special orbit \(C\). In a nonzero fiber, \(U(\beta)^t=U(\alpha)^t\) gives \(U(\beta)=aU(\alpha)\) for some \(a\in\mu_t\).
Thus \(L_V(\beta-c-a(\alpha-c))=0\), so the two points lie in the same affine orbit.
The converse follows directly from \(U(g(X))=aU(X)\).

If \(P\) is invariant under all translations, Theorem~\ref{ffs:thm-additive-composition} and \(\delta\in k\) give \(P=R(U)\), with \(R\in k[Y]\).
By injectivity of substitution and coefficient comparison, the additional requirement \(R(aU)=R(U)\) is equivalent to every nonzero exponent of \(R\) being divisible by \(t\).
This is exactly \(P\in k[U^t]\).
The scalar multiplier of the normal form under slope \(a\) is \(a^{ew}\). Since \(\gcd(e,t)=1\), it equals one for every \(a\in\mu_t\) if and only if \(t\mid w\).
\end{proof}

The factor \(U^{ew}\) in the normal form records root multiplicities on the special orbit, and \(H(U^t)^e\) records those in the nonzero fibers.
The polynomial \(f\) itself is invariant under the affine group if and only if \(t\mid w\).
For the general relation between polynomial invariants and group orbit decompositions, see~\cite{FFSSchulz2023}.

\subsection{Polynomial Spaces with a Fixed Scalar Multiplier}
\label{fas:subsec-semi-invariants}

In the substitution identity, preserving root multiplicities requires only multiplication of the polynomial by a nonzero constant; an invariant polynomial requires that constant to be one.
Fixing a scalar multiplier allows us to describe all polynomials satisfying the condition at once, without finding their normal forms individually.

\begin{definition}\leavevmode\phantomsection\label{fas:def-semi-invariant}
Suppose \(t>1\), and retain \(\Gamma=\mathcal A_f\), \(Q=|V_f|\), and \(U=L_V-\delta\).
For \(0\le r<t\), a polynomial \(P\in k[X]\) satisfying
\[
P(aX+b)=a^rP(X)\qquad(g_{a,b}\in\Gamma)
\]
is called a \textbf{semi-invariant} with the corresponding scalar multiplier.
For an integer \(N\ge0\), denote the \(k\)-linear space of such polynomials of degree at most \(N\), including the zero polynomial, by \(\mathcal W_{\Gamma,r}(N)\).
\end{definition}

\begin{proposition}[A basis of the semi-invariant space]\leavevmode\phantomsection\label{fas:prop-semi-invariant-basis}
All the semi-invariants above form \(U^r k[U^t]\), and
\begin{equation}\label{fas:eq-semi-invariant-basis}
\mathcal W_{\Gamma,r}(N)
=\bigoplus_{\substack{j\ge0\\Q(r+tj)\le N}}k\,U^{r+tj}.
\end{equation}
Consequently,
\[
\dim_k\mathcal W_{\Gamma,r}(N)
=\max\left\{0,\,1+\left\lfloor\frac{N-Qr}{Qt}\right\rfloor\right\}.
\]
\end{proposition}
\begin{proof}
For translations, whose slope is \(a=1\), the semi-invariant condition is translation invariance. Theorem~\ref{ffs:thm-additive-composition} and \(\delta\in k\) give the unique expression \(P=R(U)\), with \(R\in k[Y]\).
Using \(U(aX+b)=aU(X)\), the condition becomes \(R(aY)=a^rR(Y)\) for every \(a\in\mu_t\).
Coefficient comparison shows that every exponent of a nonzero term of \(R\) must be congruent to \(r\) modulo \(t\), so \(P\in U^rk[U^t]\).
The converse follows by substitution. The powers \(U^{r+tj}\) have distinct degrees \(Q(r+tj)\), so are linearly independent; imposing the degree bound gives the basis and dimension formula.
\end{proof}

\leavevmode\phantomsection\label{fas:def-common-support}\ignorespaces
Write \(\operatorname{Supp}_X(P)=\{i:[X^i]P\ne0\}\), and define the support of the zero polynomial to be empty.
The \textbf{common support set} of \(\mathcal W_{\Gamma,r}(N)\) is exactly
\[
\bigcup_{\substack{j\ge0\\Q(r+tj)\le N}}
\operatorname{Supp}_X(U^{r+tj}).
\]
This is the smallest set containing the support of every polynomial in the space: the support of any linear combination lies in the union of the basis supports, and every basis vector belongs to the space.
This union need not be the support of any single polynomial in the space.
For example, over \(\F_2\), when \(U=X+1\), the common support of \(\operatorname{span}_{\F_2}\{1,U^3,U^6\}\) is \(\{0,1,2,3,4,6\}\).
If a member has nonzero cubic and sextic terms, the coefficients of both \(U^3,U^6\) must be one; its quadratic terms then cancel. Thus no polynomial in the space has the entire common support as its support.
The power-basis expression for semi-invariants also converts coefficient restrictions into finitely many linear equations.
For example, given a set of allowed exponents \(I\subseteq\{0,\ldots,N\}\), write
\[
P=\sum_{\substack{j\ge0\\Q(r+tj)\le N}}\lambda_jU^{r+tj},\qquad \lambda_j\in k.
\]
Then \(\operatorname{Supp}_X(P)\subseteq I\) if and only if
\[
\sum_{\substack{j\ge0\\Q(r+tj)\le N}}
\lambda_j[X^i]U^{r+tj}=0
\qquad(0\le i\le N,\ i\notin I).
\]
One can therefore first determine the power basis from the group symmetry, then impose support restrictions by solving linear equations over the base field.
If certain coefficients must also be nonzero, this must be checked separately on the solutions.
Section~\ref{v54:subsec-affine-test-codes} combines these spaces with divisibility conditions from a generator polynomial to obtain subspaces of codes.

\begin{corollary}[Common support bounds from affine symmetries]\label{v67:cor-affine-support-bounds}
Let \(f\in k[X]\) be monic with \(f(0)\ne0\) and \(|S_f|\ge2\). Let \(e\) be the common \(p\)-power of \(f\), write \(f=h^e\) with \(h'\ne0\), and put
\[
d=\deg h,\qquad \Gamma=\mathcal A_f,\qquad V=V_f,\qquad
Q=|V|,\qquad m=|M_f|,\qquad s=\operatorname{sp}(f)=em.
\]
Let \(t\) be the order of the slope image. Use \(U=L_V-\delta\) when \(t>1\), and set \(U=L_V\) when \(t=1\).
Let \(r\) be the least nonnegative residue of \(d\) modulo \(t\), and define
\begin{align*}
E_\Gamma&=\bigcup_{\substack{j\ge0\\Q(r+tj)\le d}}
\operatorname{Supp}_X(U^{r+tj}),\\
E_0&=\bigcup_{\substack{j\ge0\\Qmj\le d}}
\operatorname{Supp}_X(L_V^{mj}).
\end{align*}
For a set of exponents \(I\), write \(eI=\{ei:i\in I\}\); this multiplies the exponents by \(e\) and leaves the cardinality unchanged.
Then
\begin{equation}\label{v67:eq-affine-support-chain}
\operatorname{Supp}_X(f)\subseteq eE_\Gamma\subseteq eE_0
\subseteq\{0,s,2s,\ldots,\deg f\},
\end{equation}
and consequently
\begin{equation}\label{v67:eq-affine-support-cardinalities}
|\operatorname{Supp}_X(f)|\le |E_\Gamma|\le |E_0|
\le \frac{\deg f}{s}+1.
\end{equation}
\end{corollary}
\begin{proof}
The identity \(s=em\) follows from the common \(p\)-power and the exponent criterion for scalings, since \(f(0)\ne0\).
The polynomial \(h\) has the same affine stabilizer as \(f\), and satisfies \(h(aX+b)=a^dh(X)\).
For \(t>1\), Proposition~\ref{fas:prop-semi-invariant-basis} therefore gives \(\operatorname{Supp}_X(h)\subseteq E_\Gamma\).
For \(t=1\), translation invariance gives \(h\in k[L_V]\), with the same conclusion because \(r=0\).
The Frobenius identity \(f=h^e\) gives \(\operatorname{Supp}_X(f)=e\operatorname{Supp}_X(h)\), proving the first inclusion.

Each power defining \(E_\Gamma\) is translation-invariant and has scalar multiplier \(a^d\) under \(\Gamma\).
For \(a\in M_f\), evaluation of \(h(aX)=a^dh(X)\) at zero gives \(a^d=1\), since \(h(0)\ne0\).
Thus each such power is also invariant under scalings about the origin by \(M_f=\mu_m\).
Translation invariance first expresses it as \(R(L_V)\); using \(L_V(aX)=aL_V(X)\), invariance under \(\mu_m\) then forces every exponent in \(R\) to be divisible by \(m\).
The degree bound gives a linear combination of the powers defining \(E_0\), and hence \(E_\Gamma\subseteq E_0\).

Each \(L_V^{mj}\) is invariant under \(X\mapsto aX\) for every \(a\in\mu_m\), so its nonzero coefficient exponents are multiples of \(m\) and are at most \(d\).
Multiplication by \(e\) therefore gives multiples of \(s=em\) between zero and \(ed=\deg f\).
Taking cardinalities proves the final inequalities.
\end{proof}

\section{Recovering the Normal-Form Data: Three Examples}
\label{ffs:sec-computation}

\begin{proposition}[A finite coefficient recovery procedure]\label{ffs:prop-recovery-procedure}
Let \(f\in k[X]\) be monic with \(|S_f|\ge2\). The following procedure determines the full affine normal form from the coefficients of \(f\).
\begin{enumerate}
\item Find \(e\) from the \(p\)-adic valuations of the positive exponents occurring in the polynomial, and extract the unique monic root \(h=f^{1/e}\).
\item Take the monic greatest common divisor of the coefficients of powers of \(X\) in \(h(X+B)-h(X)\), obtaining \(L_V=\Psi_h\).
Use Proposition~\ref{ffs:prop-multiplier-support} separately to find the multiplicative stabilizer \(M_f\) about the origin.
\item Solve the coefficient equations~\eqref{ffs:eq-affine-coefficients} in a finite splitting extension to find \(\mathcal A_f\) and its slope image \(\mu_t\).
These equations may instead be solved for the lower-degree polynomial \(h\). Once \(h=F(L_V)\) has been recovered, Proposition~\ref{v70:prop-outer-center-recovery} also obtains the slopes from outer centers and support.
\item If \(t>1\), take a generator \(a\ne1\) of the slope image and an element \(g_{a,b}\in\mathcal A_f\) of slope \(a\), and compute
\[
c=\frac b{1-a},\qquad \delta=L_V(c),\qquad U=L_V-\delta.
\]
The resulting \(\delta\) lies in \(k\) and is independent of the chosen \(g_{a,b}\).
\item If \(t=1\), recover the unique expression \(h=F(L_V)\). If \(t>1\), first find the largest \(w\) with \(U^w\mid h\), then recover the unique expression \(h/U^w=H(U^t)\).
\end{enumerate}
\end{proposition}
\begin{proof}
The first four steps follow respectively from the common \(p\)-power, the translation coefficient formula, the affine coefficient criterion, and the center descent theorem.
Divisibility in the fifth step is tested by ordinary polynomial division.
Given a monic nonconstant \(P\) and a polynomial \(R\) known to belong to \(k[P]\), the outer polynomial can likewise be recovered by removing leading terms:
if \(\deg R=m\deg P\) and its leading coefficient is \(a\), replace \(R\) by \(R-aP^m\), and repeat until zero.
At each step the degree strictly decreases and the remaining polynomial still belongs to \(k[P]\), so all outer coefficients are recovered.
\end{proof}

This procedure always terminates, although the enumeration in the third step may be large.
Without first factoring \(f\), one may take
\(N=\lcm(1,2,\ldots,n)\) and \(E=\F_{q^N}\):
each irreducible factor has degree at most \(n\), so all its roots lie in \(E\).
One can enumerate \(E^\times\times E\) and test the coefficient equations, but this general bound on the extension degree is usually much larger than the degree actually needed.
By contrast, the first two steps take place entirely over \(k\), requiring neither factorization of \(f\) nor enumeration of the elements of \(V_f\).
When \(f\) is irreducible, Subsections~\ref{fas:subsec-frobenius-differences}--\ref{fas:subsec-centered-recovery} give another route to the full group using specified conjugate differences and a base-field center.

Once the maximal translation inner component is known, the remaining slopes can also be recovered from centers and support of the outer polynomial.

\begin{proposition}[Recovering the full affine group from outer centers]\label{v70:prop-outer-center-recovery}
Let \(f\in k[X]\) be monic with \(|S_f|\ge2\), and let \(E\) be a finite splitting extension. Write
\[
f=h^e,\qquad h'\ne0,\qquad h=F(L_V),\qquad V=V_f,
\]
where \(F\) is the monic outer polynomial of the maximal translation quotient.
\begin{enumerate}
\item Suppose \(F\) has at least two distinct roots. For each \(z\in k\), write
\[
F_z(Y)=F(Y+z)=\sum_i d_{z,i}Y^i,
\qquad w_z=\min\{i:d_{z,i}\ne0\},
\]
and put
\[
g_z=\gcd\{i-w_z:i>w_z,\ d_{z,i}\ne0\},
\qquad r_z=\frac{g_z}{p^{v_p(g_z)}}.
\]
The set of exponent differences here is nonempty. Define
\[
t_z=
\begin{cases}
\gcd\bigl(r_z,|H_{\F_p}(V)|-1\bigr),&V\ne\{0\},\\
r_z,&V=\{0\}.
\end{cases}
\]
At most one \(z\in k\) satisfies \(t_z>1\). If none does, \(\mathcal A_f\) consists only of translations from \(V\).
Otherwise put \(\delta=z\) and \(t=t_z\). Then
\begin{equation}\label{v70:eq-outer-center-group}
\mathcal A_f=
\{X\mapsto aX+b:a\in E^\times, a^t=1,
\ b\in E, L_V(b)=(1-a)\delta\}.
\end{equation}
Each slope has exactly \(|V|\) lifts.
\item If \(F\) has only one distinct root, then \(F=(Y-\delta)^w\), where \(\delta\in k\), \(p\nmid w\), and \(V\ne\{0\}\).
In this case put \(t=|H_{\F_p}(V)|-1\); formula~\eqref{v70:eq-outer-center-group} again gives the full affine group.
\end{enumerate}
In either case, if \(t>1\), then \(\delta\in L_V(E)\); a center can then be chosen in \(k\) precisely when \(\delta\in L_V(k)\).
\end{proposition}
\begin{proof}
Corollary~\ref{ffs:cor-maximal-additive} gives \(F'\ne0\) and \(V_F=\{0\}\).
First suppose \(F\) has at least two distinct roots. Its affine stabilizer is finite, and all its parameters belong to \(E\), since an affine transformation is determined by the images of two distinct roots.
As \(V_F=\{0\}\), the slope map embeds this group into the cyclic group \(E^\times\), so it too is cyclic.
If it is nontrivial, take a generator \(Y\mapsto aY+d\), with \(a\ne1\). Its unique fixed point is \(\delta=d/(1-a)\), which is fixed by the whole group.
Since \(F\in k[Y]\), Frobenius preserves the group and hence fixes its unique common center. Thus \(\delta\in k\).

For every \(z\in k\), the polynomial \(F_z\) is not a monomial. Proposition~\ref{ffs:prop-multiplier-support} gives its multiplicative stabilizer about the origin as \(\mu_{r_z}\).
Returning to the original coordinate gives exactly the outer affine transformations centered at \(z\).
If \(V\ne\{0\}\), Proposition~\ref{v70:prop-scalar-support} identifies the compatible slopes with \(H_{\F_p}(V)^\times\), whose intersection with \(\mu_{r_z}\) has order \(t_z\).
If \(V=\{0\}\), then \(L_V=X\), and there is no further compatibility restriction.
In that case \(\mu_{r_z}\) is already contained in \(E^\times\): these scalings preserve the roots of \(F_z\), and \(F_z\) has a nonzero root, so each scalar is the ratio of its image to itself.
There is consequently no need for an additional greatest common divisor with \(|E|-1\).
Every \(t_z>1\) gives the center of a nonidentity outer transformation. The nontrivial outer group has only one common center, so at most one such \(z\) exists.

We next verify lifting inside the fixed field \(E\).
The slope \(a\) of an original affine stabilizer satisfies \(aV=V\) and
\(L_V(aX)=aL_V(X)\). For \(Q=|V|>1\), the latter follows from the support criterion for the maximal field of scalars; for \(Q=1\), it is immediate.
Thus the original transformation induces \(Y\mapsto aY+L_V(b)\) on the quotient and preserves the root multiplicities of \(F\).
Conversely, take a compatible outer transformation \(Y\mapsto aY+d\) and a root \(\alpha\in E\) of \(h\).
It sends \(L_V(\alpha)\) to another root of \(F\). Since \(h=F(L_V)\) splits completely in \(E\), there is a root \(\alpha'\in E\) with
\[
L_V(\alpha')=aL_V(\alpha)+d.
\]
Put \(b=\alpha'-a\alpha\); then \(L_V(b)=d\). Every root fiber of \(L_V\) is separable, and its roots have in \(h\) the multiplicity of the corresponding root of \(F\). Hence \(X\mapsto aX+b\) preserves the root multiplicities of \(h\), and therefore of \(f\).
All lifting translation parameters form \(b+V\), giving exactly \(|V|\) lifts.
This proves~\eqref{v70:eq-outer-center-group} and shows that only the translation group remains when there is no nontrivial compatible slope.
The argument uses splitting of all root fibers in \(E\), not surjectivity of \(L_V:E\to E\).

If \(t>1\), choose \(a\ne1\) and a lift \(b\) satisfying \(L_V(b)=(1-a)\delta\).
For \(V\ne\{0\}\), we have \((1-a)^{-1}\in H_{\F_p}(V)\), and therefore
\[
L_V\left(\frac{b}{1-a}\right)=\delta.
\]
For \(V=\{0\}\), the same equality follows directly from \(L_V=X\). This produces a center in \(E\), without ensuring that it lies in \(k\).

Finally, suppose \(F=(Y-\delta)^w\). Frobenius fixes its unique root, so \(\delta\in k\), and \(F'\ne0\) gives \(p\nmid w\).
Choose a root \(c\in E\) of \(h\). Then \(L_V(c)=\delta\), and the original root support is exactly \(c+V\). Since \(|S_f|\ge2\), we have \(V\ne\{0\}\).
All these roots have the same multiplicity, and every \(a\in H_{\F_p}(V)^\times\) gives all transformations of that slope preserving the coset:
\[
X\mapsto aX+(1-a)c+v,\qquad v\in V.
\]
Conversely, preserving the coset forces \(aV=V\), so no other slopes occur. This proves the second case. The group formula also holds when \(t=1\), in which case the group consists of translations.

In either case, if \(t>1\), formula~\eqref{v70:eq-outer-center-group} also shows that \(c\in E\) can serve as a center representative in~\eqref{ffs:eq-centered-group} if and only if \(L_V(c)=\delta\).
Sufficiency follows from \(L_V((1-a)c+v)=(1-a)\delta\). Conversely, choose a slope \(a\ne1\), substitute \(b=(1-a)c\) in the group formula, and divide by \(1-a\) to obtain \(L_V(c)=\delta\).
Thus a center in the original coordinate can be chosen in \(k\) if and only if \(\delta\in L_V(k)\).
\end{proof}

Thus, once \(V\) and \(F\) are known, the case of at least two distinct outer roots requires testing only the \(q\) possible centers in the base field.
The support criterion uses exponent differences, and the \(p\)-part of their greatest common divisor must be removed. For example, over \(\F_2\), the scaling group of \(Y^5+Y^2\) about zero is \(\mu_3\), whereas \(Y^3+Y\) has exponent-difference greatest common divisor two but trivial scaling group.
When the outer polynomial has only one distinct root, the second case applies; no greatest common divisor of an empty set of exponent differences is needed.

\begin{example}[Scalings about the origin do not recover all slopes]\leavevmode\phantomsection\label{ffs:ex-off-center}
Over \(\F_5\), take
\[
f(X)=(X-1)^2-2=X^2+3X+4.
\]
Since \(2\) is a nonsquare, \(f\) is irreducible and has two distinct roots.
Its root set has size two, whereas \(|V_f|\) is a power of five, so \(V_f=\{0\}\) and \(\Psi_f=B\).
The nonzero exponents have greatest common divisor one, so \(M_f=\{1\}\).
However, \mbox{\(f(2-X)=f(X)\)}, and the nonidentity transformation \(X\mapsto2-X\) interchanges the roots.
An affine transformation is determined by the images of two roots, so the full stabilizer is exactly
\[
\mathcal A_f=\{X\mapsto X,\ X\mapsto2-X\}.
\]
Thus \(t=2,c=1,U=X-1,w=0,H(Y)=Y-2\).
Here \(t=2>|M_f|=1\); the difference comes precisely from the nonzero affine center.
\end{example}

\begin{example}[A center outside the base field and a center-coset polynomial over the base field]\label{ffs:ex-quartic-center}
Take the irreducible polynomial \(f=X^4+X+1\in\F_2[X]\) from Example~\ref{exm:quartic-model}.
Direct calculation gives
\[
f(X+B)-f(X)=B^4+B,
\]
hence \(e=1,V_f=\F_4,L_V=X^4+X\).
Take a root \(c\in\F_{16}\); then the root set is \(c+\F_4\).
Its full affine stabilizer is
\[
\mathcal A_f=\{X\mapsto c+a(X-c)+v:
 a\in\F_4^\times,\ v\in\F_4\},
\qquad |\mathcal A_f|=12.
\]
Indeed, all these 12 transformations preserve the root set; conversely, every slope must preserve \(\F_4\), so lies in \(\F_4^\times\).
Thus \(t=3\), whereas the multiplicative stabilizer about the origin is still \(M_f=\{1\}\).

Here \(\delta=L_V(c)=1\), so \(U=L_V-1=X^4+X+1\in\F_2[X]\).
The normal-form parameters are \(w=1,H=1\), and all roots form the special orbit.
Since \(U\) has no base-field root, the center coset is disjoint from the base field.
The base-field affine stabilizer consists of \(X\mapsto X\) and \(X\mapsto X+1\), and has order two.
For any root \(\alpha\), its point stabilizer is
\[
(\mathcal A_f)_\alpha
=\{X\mapsto\alpha+a(X-\alpha):a\in\F_4^\times\},
\]
and hence has order three. The four stabilizers are compared in the following table:
\begin{center}\small
\renewcommand{\arraystretch}{1.15}
\begin{tabular}{@{}llc@{}}
\toprule
Group&Restriction&Order\\
\midrule
\(\mathcal A_f\)&Preserves the full root multiplicity function&\(12\)\\
\(\mathcal A_f(\F_2)\)&Affine coefficients belong to \(\F_2\)&\(2\)\\
\((\mathcal A_f)_\alpha\)&Fixes the root \(\alpha\)&\(3\)\\
\((\mathcal A_f)_0\simeq M_f\)&Fixes the origin&\(1\)\\
\bottomrule
\end{tabular}
\end{center}
Example~\ref{exm:compression-return-four} further compares how this root-orbit period is distributed between quotient and fiber when taking quotients by translations from \(\F_2\) and \(\F_4\).
\end{example}

\begin{example}[Residual unequal multiplicities]\label{ffs:ex-unequal-multiplicity}
Over \(\F_3\), put \(L=X^3-X\), and take
\[
f=(L-1)(L+1)^2.
\]
Neither \(X^3-X-1\) nor \(X^3-X+1\) has a root in \(\F_3\), so they are distinct irreducible cubics.
The six distinct roots form two additive \(\F_3\)-cosets, with root multiplicities one and two, respectively.
Thus \(V_f=\F_3\), \(e=1\), and \(\Psi_f=B^3-B\).
The translation group cannot have order nine, since its order must divide the number of distinct roots, six.

Every affine slope must preserve \(V_f\), so is either \(1\) or \(-1\).
In the quotient coordinate \(L\), the roots \(1,-1\) have different multiplicities and must therefore be fixed individually;
the only affine transformation fixing both points is the identity.
Hence \(\mathcal A_f\) is exactly the translation group of \(\F_3\), and \(t=1\).
The outer polynomial \(F(Y)=(Y-1)(Y+1)^2\) satisfies \(F'\ne0\) but is not square-free.

Replacing \(f\) by \(f^3\) leaves all three stabilizers unchanged and changes the root multiplicities to three and six, while
\[
e_{f^3}=3,\qquad \Psi_{f^3}(B)=(B^3-B)^3.
\]
The common \(p\)-power and the residual outer multiplicities are thus recorded separately.
\end{example}

\section{Frobenius Return and Recovery of Irreducible Factors}
\label{ffs:sec-frobenius-factorization}

Symmetry orbits and Frobenius orbits generally give different partitions.
The former determine quotient polynomials and the corresponding normal forms, while the latter determine irreducible factors.
Lemma~\ref{lem:fiber-return} connects them: after one circuit of the base orbit, the return period in a fiber, multiplied by the base orbit length, is the degree of the original root.
We now express this return map as a computable finite permutation, retaining multiplicities as well.

\subsection{Affine Return on Additive Fibers}

\begin{theorem}[Complete factorization type of an additive composition]\label{ffs:thm-additive-return}
Let \(V\) be a finite Frobenius-stable \(\F_p\)-subspace, put \(Q=|V|\), and let \(L=L_V\in k[X]\).
Let \(P\in k[Y]\) be monic and irreducible of degree \(d\), take a root \(\beta\), and then take \(\alpha\) with \(L(\alpha)=\beta\).
Define on \(V\)
\[
A(v)=v^{q^d},\qquad b_\beta=\alpha^{q^d}-\alpha\in V,
\qquad R_\beta(v)=A(v)+b_\beta.
\]
If the cycle lengths of \(R_\beta\) are \(\ell_1,\ldots,\ell_s\), then \(P(L(X))\) has exactly \(s\) distinct monic irreducible factors, of degrees
\[
d\ell_1,\ldots,d\ell_s.
\]
Each factor has multiplicity one. Choosing a different \(\alpha\) in the same fiber gives a conjugate finite permutation, so the cycle lengths and their numbers are unchanged; \leavevmode\phantomsection\label{def:permutation-cycle-type}\ignorespaces these data are called the \textbf{cycle type} of the permutation.
\end{theorem}
\begin{proof}
Since \(L(\alpha^{q^d})=\beta^{q^d}=\beta\), additivity gives \(b_\beta\in V\).
The map \(A\) is an invertible \(\F_p\)-linear map on \(V\).
The fiber \(L(X)=\beta\) is \(\alpha+V\); in the coordinates \(\alpha+v\leftrightarrow v\), the \(d\)th iterate of Frobenius is exactly
\[
\alpha+v\longmapsto\alpha^{q^d}+v^{q^d}
 =\alpha+R_\beta(v).
\]
By the fiber return lemma, each return cycle of length \(\ell_j\) corresponds to a Frobenius orbit of length \(d\ell_j\), whose root polynomial is the corresponding irreducible factor.
Since \(L'\) is a nonzero constant and \(P\) is square-free, \(P(L)\) is square-free too.
Replacing \(\alpha\) by \(\alpha+v_0\) translates the coordinates, so the return permutations are conjugate.
\end{proof}

\begin{corollary}[Recovering multiplicities]\label{ffs:cor-additive-multiplicities}
Suppose
\[
f=F(L)^e,\qquad F=\prod_i P_i^{h_i},
\]
where the \(P_i\) are distinct monic irreducible polynomials. Factor each \(P_i(L)\) by Theorem~\ref{ffs:thm-additive-return};
each irreducible factor occurring there has multiplicity exactly \(eh_i\) in \(f\).
Factors arising from different \(P_i\) are distinct.
\end{corollary}
\begin{proof}
The root sets of the \(P_i\) are disjoint, so their inverse images under \(L\) are disjoint as well.
Each inverse-image polynomial is square-free; raising it to the \(eh_i\)th power gives exactly the stated multiplicities.
\end{proof}

\begin{proposition}[Counting factors by linear equations]\label{ffs:prop-return-counts}
Use the notation of Theorem~\ref{ffs:thm-additive-return}. For \(j\ge1\), put
\[
b_j=\sum_{i=0}^{j-1}A^i(b_\beta),\qquad
N_j=\#\Fix(R_\beta^j).
\]
Then
\[
N_j=
\begin{cases}
p^{\dim_{\F_p}\ker(A^j-I)},&-b_j\in\im(A^j-I),\\
0,&-b_j\notin\im(A^j-I).
\end{cases}
\]
The number of irreducible factors of \(P(L)\) of degree \(d\ell\) is
\begin{equation}\label{ffs:eq-return-mobius}
C_\ell=\frac1\ell\sum_{j\mid\ell}\mu(\ell/j)N_j.
\end{equation}
\end{proposition}
\begin{proof}
Induction gives \(R_\beta^j(v)=A^jv+b_j\). The fixed-point condition is the linear equation
\((A^j-I)v=-b_j\). If solvable, its solution set is a coset of \(\ker(A^j-I)\); otherwise there are no fixed points.
If \(C_\ell\) is the number of cycles of length \(\ell\), then \(N_j=\sum_{\ell\mid j}\ell C_\ell\).
M\"obius inversion gives~\eqref{ffs:eq-return-mobius}; the return factorization theorem then identifies cycle counts with factor counts.
\end{proof}

\begin{corollary}[Two equal-degree factorization types for base-field translations]\label{ffs:cor-base-translations}
If \(0\ne V\subseteq k\) and \(|V|=p^s\), then \(P(L_V)\) has only two possible factorization types:
\[
\begin{array}{c|c|c}
\text{Return translation}&\text{Number of distinct irreducible factors}&\text{Degree of each factor}\\ \hline
b_\beta=0&p^s&d\\
b_\beta\ne0&p^{s-1}&pd
\end{array}
\]
In particular, \(P(L_V)\) is reducible whenever \(s\ge2\).
More generally, the same conclusion requires only \(V\subseteq\F_{q^d}\).
\end{corollary}
\begin{proof}
This hypothesis gives \(A=I\). If \(b_\beta=0\), the return permutation is the identity; otherwise it is a nonzero translation, with all cycles of length \(p\).
The respective cycle counts are \(p^s\) and \(p^{s-1}\); apply Theorem~\ref{ffs:thm-additive-return}.
\end{proof}

\begin{remark}[Distinguishing the two factorizations by a remainder]\label{ffs:rem-return-remainder}
Under the assumption \(V\subseteq\F_{q^d}\), division of \(X^{q^d}-X\) by \(L_V(X)-\beta\) leaves the constant remainder \(b_\beta\).
Indeed, at every root \(\alpha+v\), the former polynomial has value \(\alpha^{q^d}-\alpha\),
while the remainder has degree less than \(|V|\); the root bound forces it to be this constant.
Thus the two factorizations can be distinguished without first finding a root in the fiber: modular exponentiation in \(k(\beta)[X]\) suffices.
\end{remark}

For \(L=X^p-X\), this constant is \(\Tr_{\F_{q^d}/\F_p}(\beta)\).
Indeed, successively taking \(p\)th powers of \(\alpha^p-\alpha=\beta\) and adding gives
\[
\alpha^{q^d}-\alpha
 =\sum_{i=0}^{rd-1}\beta^{p^i}
 =\Tr_{\F_{q^d}/\F_p}(\beta).
\]
Thus the earlier Artin--Schreier composition criterion is the case in which the return translation is nonzero and the fiber is a single \(p\)-cycle.

\subsection{Return on Full Affine Fibers}

Now suppose \(t>1\), and take \(G=\mathcal A_f\), \(U\), and \(\pi=U^t\) from Theorem~\ref{ffs:thm-normal-form}.
Let \(P\mid H\) be monic and irreducible, put \(d=\deg P\), and let \(\beta\ne0\) be a root.
Take \(\alpha\) with \(\pi(\alpha)=\beta\); its fiber is a free \(G\)-orbit.

\begin{theorem}[Factorization over a full affine fiber]\leavevmode\phantomsection\label{ffs:thm-full-return}
There is a unique \(g_0\in G\) such that \(\sigma^d(\alpha)=g_0(\alpha)\).
Put
\[
\tau(g)=\sigma^d\circ g\circ\sigma^{-d},\qquad
\mathcal R(g)=\tau(g)g_0\qquad(g\in G).
\]
Here \(\tau\) is an automorphism of \(G\), and multiplication on the right means composition of transformations.
If the finite permutation \(\mathcal R\) has cycle lengths \(\ell_1,\ldots,\ell_s\), then the distinct irreducible factors of \(P(\pi(X))\) have degrees exactly
\[
d\ell_1,\ldots,d\ell_s.
\]
If \(P\) has multiplicity \(h_P\) in \(H\), all these factors have multiplicity \(eh_P\) in \(f\).
\end{theorem}
\begin{proof}
Since \(\pi\in k[X]\), the map \(\sigma^d\) preserves the chosen fiber; freeness gives the existence and uniqueness of \(g_0\).
Conjugation by Frobenius maps \(G\) to itself, so \(\tau\) is an automorphism of \(G\).
Under the bijection \(G\to\pi^{-1}(\beta)\), \(g\mapsto g(\alpha)\), we have
\[
\sigma^d(g(\alpha))
 =\tau(g)(\sigma^d(\alpha))
 =(\tau(g)g_0)(\alpha).
\]
Thus \(\mathcal R\) is the return map of this fiber.
The fiber return lemma gives the degrees, and square-freeness of \(U^t-\beta\) gives the multiplicity formula.
\end{proof}

The conjugation action can also be calculated directly from parameters.
Writing \(g=g_{a,b}\), we have \(\tau(g)=g_{a^{q^d},b^{q^d}}\).
In centered coordinates \(g(X)=c+a(X-c)+v\), put \(\lambda_d=c^{q^d}-c\in V\). The parameters then change by
\[
(a,v)\longmapsto
\bigl(a^{q^d},\ v^{q^d}+(1-a^{q^d})\lambda_d\bigr).
\]
This formula computes the return permutation directly from the affine parameters.

\begin{corollary}[Coefficients fixed by the return Frobenius]\leavevmode\phantomsection\label{ffs:cor-constant-group-return}
If every transformation in \(G\) is defined over \(\F_{q^d}\), then \(\tau=\id\).
In this case \(\mathcal R\) is right multiplication by \(g_0\), and all cycles have length \(\ord(g_0)\).
Consequently, the number of distinct irreducible factors of \(P(\pi)\) is
\[
\frac{|G|}{\ord(g_0)}
\]
and each has degree \(d\ord(g_0)\).
\end{corollary}
\begin{proof}
The \(j\)th iterate of right multiplication by \(g_0\) fixes some \(g\) if and only if \(g_0^j=1\), so every cycle has the same length \(\ord(g_0)\).
\end{proof}

The return for a base-field affine group in Section~\ref{sec:compression-return-lifting} has Frobenius-fixed parameters;
thus all return cycles in the same free fiber there have equal length.

\begin{example}[Frobenius on the parameters cannot be omitted]\label{ffs:ex-nonuniform}
Over \(\F_2\), take \(f=X^3-1\). Its root set is \(\mu_3\subseteq\F_4\), \(V_f=\{0\}\), and its full affine stabilizer is the scaling group \(G=\mu_3\).
Here \(U=X,\pi=X^3,H(Y)=Y-1\). Choose \(\beta=1,\alpha=1\); then \(g_0=1\), but \(\tau(a)=a^2\).
On the three group elements, \(\tau\) fixes one and interchanges the other two, so the return cycle lengths are one and two.
Correspondingly,
\[
X^3-1=(X+1)(X^2+X+1)\quad\text{in }\F_2[X].
\]
Here the return map is the squaring permutation induced by Frobenius on \(\mu_3\); its orbits of lengths one and two give the linear and quadratic factors, respectively.
\end{example}

Factors from the special orbit can also be recovered directly.
On \(C=c+V\), Frobenius has coordinate form
\[
v\longmapsto v^q+(c^q-c).
\]
Its cycle lengths are the degrees of the irreducible factors of \(U\), all of which have multiplicity \(ew\) in \(f\).
Together with the free fibers for each \(P\mid H\), Theorem~\ref{ffs:thm-full-return} therefore gives all irreducible factorization data of the normal form.

\section{Stabilizers of Irreducible Polynomials: Restrictions, Recovery, and Classification}
\label{ffs:sec-irreducible}

Irreducibility requires Frobenius to act as a single cycle on all roots. This restricts the sizes of the base-field affine group and translation blocks, and also allows specified conjugate differences to identify all translations.
After passing to the translation quotient, compatible slopes and centered multiplicative orders recover the full affine group; a centered multiplicative order is the multiplicative order obtained by subtracting a candidate center from a quotient root.
This yields explicit classifications in prime degree and in the quartic case in characteristic two.

\subsection{Restrictions on the Base-Field Group and Translation Group}
\label{fas:subsec-group-restrictions}

\begin{theorem}[Classification of the base-field affine stabilizer]\label{ffs:thm-arithmetic-cyclic}
If \(f\in k[X]\) is monic and irreducible with \(n=\deg f\ge2\), then \(\mathcal A_f(k)\) is cyclic and its order divides \(n\).
More precisely, it is of one of the following types:
\begin{enumerate}
\item the trivial group;
\item a group of order \(p\) generated by a nonzero translation;
\item a group that becomes the scaling group \(\mu_d\subseteq k^\times\) after a base-field translation of coordinates, where \(d>1\) and \(d\mid\gcd(n,q-1)\).
\end{enumerate}
In particular, \(|V_f\cap k|\le p\).
\end{theorem}
\begin{proof}
Take a root \(\alpha\) of \(f\). For each \(g\in\mathcal A_f(k)\), there is a unique \(j\bmod n\) such that \(g(\alpha)=\alpha^{q^j}\).
Since \(g\) commutes with Frobenius,
\[
g(\alpha^{q^i})=\alpha^{q^{i+j}}\qquad(0\le i<n).
\]
Thus \(g\) acts on all roots as the \(j\)th power of the same cycle, and composition corresponds to addition of exponents.
An affine transformation fixing two distinct points is the identity, giving a group embedding
\(\mathcal A_f(k)\hookrightarrow\Z/n\Z\).
This proves that the group is cyclic of order dividing \(n\).

Take a generator \(g(X)=aX+b\). If \(a=1\), then \(b=0\) gives the identity and \(b\ne0\) gives a translation of order \(p\).
If \(a\ne1\), then \(c=b/(1-a)\in k\), and \(g(X)-c=a(X-c)\).
Its order is therefore \(\ord(a)\mid q-1\), and conjugation by a base-field translation makes it a scaling group.
Finally, \(V_f\cap k\) is both a subgroup of a cyclic group and an additive group in which every nonzero element has order \(p\), so its order is at most \(p\).
\end{proof}

The corresponding root structures are also explicit.
The translation type gives \(f=F(X^p-b^{p-1}X)\), and the scaling type gives \(f=F((X-c)^d)\).
In the latter case there is no factor \(X-c\), since an irreducible polynomial of degree at least two cannot have a root \(c\in k\).
The outer polynomial \(F\) must be irreducible; in a converse construction, irreducibility of the inner fiber must still be checked by the root fiber criterion, since symmetry alone does not ensure irreducibility.

\begin{theorem}[Size restrictions on all root translations]\label{ffs:thm-geometric-bound}
If \(f\in k[X]\) is irreducible with \(\deg f\ge2\), then
\[
|V_f|\in
\begin{cases}
\{1,p\},&p\text{ is an odd prime},\\
\{1,2,4\},&p=2.
\end{cases}
\]
Here \(V_f\) allows all translation parameters in the splitting extension, not just base-field parameters.
\end{theorem}
\begin{proof}
The polynomial \(f\) is square-free, so \(e=1\). If \(|V_f|=p^s>1\), the maximal translation decomposition gives
\(f=F(L_{V_f})\). The outer polynomial \(F\) is irreducible, since otherwise its factorization would yield a nontrivial factorization of \(f\) by substitution.
Take a root \(\beta\) of \(F\). By Corollary~\ref{thm:capelli}, \(L_{V_f}(X)-\beta\) is irreducible over the finite field \(k(\beta)\).
Applying Theorem~\ref{thm:affine-linearized-obstruction} to this fiber gives \(s=1\) or \(p=s=2\).
These are exactly the listed possibilities.
\end{proof}

This argument also explains why the restriction applies to general irreducible polynomials, without requiring the original polynomial itself to be a linearized polynomial plus a constant:
take the maximal translation block first, then apply the affine root-space theorem over the corresponding outer root field.

\begin{corollary}[Exact conditions for a four-element translation block]\leavevmode\phantomsection\label{ffs:cor-four-block}
Suppose \(p=2\), and let irreducible \(f\) satisfy \(|V_f|=4\). Write
\[
f=F(L_V),\qquad \deg F=m,\qquad
F(Y)=Y^m+c_{m-1}Y^{m-1}+\cdots+c_0.
\]
Then \(m\) is odd, and there exist \(u,v\in k^\times\) such that
\[
L_V=X(X+u)(X^2+uX+v),
\]
where both \(X^2+uX+v\) and \(X^2+vX+c_{m-1}\) are irreducible over \(k\).
Conversely, if \(F\) is monic and irreducible, \(m\) is odd, and \(u,v\) satisfy these conditions, then \(F(L_V)\) is irreducible and its full additive stabilizer has exactly four elements.
Moreover, \(V_f\subseteq\F_{q^2}\), \(|V_f\cap k|=2\), Frobenius has order two on \(V_f\), and \(v_2(\deg f)=2\).
\end{corollary}
\begin{proof}
The necessity and sufficiency statements follow from Corollary~\ref{cor:affine-composition-classification} applied to this maximal translation inner component.
The irreducible polynomial obtained in the converse construction has four translation parameters, and Theorem~\ref{ffs:thm-geometric-bound} excludes more, so this stabilizer is full.
All roots of \(L_V=X(X+u)(X^2+uX+v)\) lie in \(\F_{q^2}\), and its base-field roots are exactly \(0,u\); hence \(|V_f\cap k|=2\).
Frobenius fixes these two points and interchanges the other two, so has order exactly two on \(V_f\).
The last assertion follows from \(\deg f=4m\) with \(m\) odd.
\end{proof}

\begin{corollary}[An irreducible root set and the special orbit]\label{ffs:cor-exceptional-irreducible}
Let \(f\) be irreducible, and suppose that the slope image of its full affine stabilizer has \(t>1\).
If the special orbit meets the root set, then all roots form this special orbit, and \(f=U\).
Otherwise \(f=H(U^t)\), where \(H\) is monic and irreducible with \(H(0)\ne0\).
\end{corollary}
\begin{proof}
The special orbit is Frobenius-stable. Once it meets the irreducible root set, transitivity on the latter forces it to contain all roots.
Conversely, since the special orbit is itself an \(\mathcal A_f\)-orbit, a nonempty intersection forces the whole orbit into the root set; the two sets are therefore equal.
Since \(f\) is monic and square-free, \(f=U\). When they are disjoint, the normal form has \(w=0,e=1\), giving \(f=H(U^t)\), whose outer polynomial must be irreducible.
\end{proof}

The bounds in Theorem~\ref{ffs:thm-geometric-bound} are attained.
For an odd prime \(p\), the polynomial \(X^p-X-1\in\F_p[X]\) is irreducible: every root \(\alpha\) satisfies \(\alpha^{p^j}=\alpha+j\), with least positive return time \(p\).
Its additive stabilizer is exactly \(\F_p\).
In characteristic two, \(X^4+X+1\) in Example~\ref{ffs:ex-quartic-center} attains the four-element bound.

\subsection{Frobenius Differences and Recovery of the Translation Stabilizer}
\label{fas:subsec-frobenius-differences}

Fix a monic irreducible \(f\in k[X]\) with \(n=\deg f\ge2\), take a root \(\alpha\), and compute in \(E=k(\alpha)=\F_{q^n}\).
This is also the splitting field of \(f\). The size restriction on the translation group allows all translation parameters to be recovered from specified conjugate differences.

\begin{theorem}[Detecting nonzero translations by conjugate differences]\label{fas:thm-translation-difference}
If \(p\nmid n\), then \(V_f=\{0\}\). If \(p\mid n\), put
\begin{equation}\label{fas:eq-prime-difference}
b_p=\sigma^{n/p}(\alpha)-\alpha\ne0.
\end{equation}
Then
\begin{equation}\label{fas:eq-translation-test}
\boxed{\quad V_f\ne\{0\}
\iff b_p^{q-1}\in\F_p^\times
\iff b_p^{(q-1)(p-1)}=1.\quad}
\end{equation}
If \(p\) is odd and the condition holds, then \(V_f=\F_p b_p\).

When \(p=2\), the case \(|V_f|=4\) is possible only when \(n=4r\) with \(r\) odd. In this case put
\[
b=\sigma^{2r}(\alpha)+\alpha=b_2,
\qquad u=\sigma^r(\alpha)+\alpha.
\]
Then
\begin{equation}\label{fas:eq-four-translation-test}
\boxed{\quad |V_f|=4\iff b^q=b\ \text{and}\ u^q=u+b.\quad}
\end{equation}
When these conditions hold, \(V_f=\{0,b,u,b+u\}\). For any even \(n\), if the degree condition and both equalities of this four-element criterion do not all hold, then \(V_f=\{0,b_2\}\) when \(b_2^q=b_2\), and \(V_f=\{0\}\) otherwise.
\end{theorem}
\begin{proof}
The number \(|V_f|\) is a power of \(p\) dividing \(n\), so equals one when \(p\nmid n\).
Now suppose \(p\mid n\), and put \(r=n/p\). The least positive Frobenius period of \(\alpha\) is \(n\), so \(b_p\ne0\).

First suppose \(|V_f|=p\). The root set partitions into \(r\) cosets of \(V_f\), on which Frobenius acts transitively.
Thus \(\sigma^r\) acts on \(\alpha+V_f\) as a \(p\)-cycle.
Identify this coset with a one-dimensional \(\F_p\)-space, taking \(\alpha\) as origin. The linear part of the return transformation lies in \(\F_p^\times\), and its order divides both \(p\) and \(p-1\), so it is the identity.
The return transformation is therefore a nonzero translation by \(b_p\), so \(V_f=\F_pb_p\).
Frobenius stability of the direction space \(V_f=\F_pb_p\) is equivalent to \(b_p^q=\lambda b_p\), with \(\lambda\in\F_p^\times\), giving necessity.

Conversely, if \(b_p^q=\lambda b_p\) with \(\lambda\in\F_p^\times\), then \(b_p\in E^\times\) gives \(\lambda^n=1\).
The order of \(\lambda\) is coprime to \(p\), so \(\lambda^r=1\), that is, \(\sigma^r(b_p)=b_p\). Inductively,
\[
\sigma^{jr}(\alpha)=\alpha+jb_p\qquad(0\le j<p).
\]
These \(p\) roots form \(\alpha+\F_pb_p\). Its direction space \(\F_pb_p\) is Frobenius-stable, so
\[
S_f=\bigcup_{i=0}^{r-1}\bigl(\sigma^i(\alpha)+\F_pb_p\bigr).
\]
Each coset is preserved by translations from \(\F_pb_p\), hence \(\F_pb_p\subseteq V_f\). In odd characteristic, Theorem~\ref{ffs:thm-geometric-bound} immediately gives equality.
Moreover, a nonzero element has \((p-1)\)st power one if and only if it lies in \(\F_p^\times\), so the last two conditions in~\eqref{fas:eq-translation-test} are equivalent.

It remains to treat the four-element case in characteristic two. Corollary~\ref{ffs:cor-four-block} gives \(n=4r\), with \(r\) odd, and \(A=\sigma|_{V_f}\) of order two.
On \(\alpha+V_f\), the map \(\sigma^r\) is a four-cycle with linear part \(A\). Thus \(u\in V_f\), and
\[
b=\sigma^{2r}(\alpha)+\alpha=u+A(u),\qquad A(b)=b.
\]
This gives \(b^q=b\) and \(u^q=u+b\), and also shows that \(b_2\) satisfies~\eqref{fas:eq-translation-test}.
Conversely, if these two equalities hold, then \(b\ne0\) and \(u\notin\{0,b\}\), so \(W=\{0,b,u,b+u\}\) is a four-element \(\F_2\)-space.
Frobenius has order two on \(W\); since \(r\) is odd, \(\sigma^r|_W=\sigma|_W\). Hence
\[
\alpha,\ \sigma^r(\alpha),\ \sigma^{2r}(\alpha),\ \sigma^{3r}(\alpha)
=\alpha,\ \alpha+u,\ \alpha+b,\ \alpha+u+b.
\]
Since the direction space \(W\) is Frobenius-stable, the full root set is the union of the Frobenius conjugate cosets of \(\alpha+W\), each preserved by translations from \(W\). Thus \(W\subseteq V_f\), and the size bound forces equality.
Together with the one-dimensional case, this exhausts all possibilities.
\end{proof}

Condition~\eqref{fas:eq-translation-test} requires only Frobenius stability of \(\F_pb_p\), not \(b_p\in k\).
For example, \(b_p^q=-b_p\) in odd characteristic still permits nontrivial translations.
Replacing \(\alpha\) by a conjugate applies Frobenius to all the differences above and does not change the recovered \(V_f\).
Given \(f\), one may take \(\alpha=X\bmod f\) directly and compute these differences in \(k[X]/(f)\), without listing all conjugate roots.

\subsection{The Translation Quotient and Lifting Compatible Slopes}
\label{fas:subsec-lifting}

\leavevmode\phantomsection\label{fas:def-translation-quotient}\ignorespaces
For a decomposition \(f=G(L_V)\), the map \(x\mapsto L_V(x)\) identifies each translation coset \(x+V\) in the root set with one point. The resulting orbit space is the \textbf{translation quotient}, with quotient coordinate \(Y=L_V(X)\).
Taking \(V=V_f\), the outer polynomial of the maximal translation quotient has no nonzero stabilizing translations, but may still have scaling symmetries or affine symmetries with nonzero centers.
Whether these symmetries lift depends on whether the slope preserves the original translation space, that is, whether \(aV=V\); such slopes are called compatible with \(V\).
The following proposition also applies to polynomials with multiplicities.

\leavevmode\phantomsection\label{fas:def-scalar-normalizer}\ignorespaces
In a fixed finite splitting extension \(E\), write
\[
\mathcal N_E(V)=\{a\in E^\times:aV=V\},\qquad h_V=|\mathcal N_E(V)|,
\]
and call \(\mathcal N_E(V)\) the \textbf{scalar stabilizer} of \(V\).
If \(V=\{0\}\), it is \(E^\times\); in this case the choice of \(E\) must be retained.

If \(V\ne\{0\}\), regard it as an \(\F_p\)-subspace. The maximal field of scalars from Definition~\ref{v70:def-scalar-field} gives
\[
\mathcal N_E(V)=H_{\F_p}(V)^\times.
\]
The base field here must be \(\F_p\), since \(V_f\) need not be a \(k\)-subspace.
Put \(Q=|V|=p^s\). Proposition~\ref{v70:prop-scalar-support} shows that \(aV=V\) implies
\[
a^Q=a,\qquad L_V(aX)=aL_V(X).
\]
Thus the slope in the quotient coordinate is still \(a\). The maximal field of scalars supplies the compatible slopes; whether a compatible slope actually occurs also depends on whether some accompanying translation preserves the outer roots and their multiplicities.
If \(V=\{0\}\), then \(L_V=X\) and \(Q=1\), so the scalar identities hold directly.

\begin{proposition}[Affine lifting from the translation quotient]\label{fas:prop-compatible-lifting}
Let \(f\in k[X]\) be monic with at least two distinct roots, put \(V=V_f\) and \(L=L_V\), and write \(f=G(L)\).
Take \(E\) over which \(f\) splits completely, and define
\begin{equation}\label{fas:eq-compatible-quotient}
\mathcal A_{G,E}^{[V]}
=\{Y\mapsto aY+d:\ a\in\mathcal N_E(V),\ d\in E,
\ G(aY+d)=a^{\deg G}G(Y)\}.
\end{equation}
Then the map
\[
\Phi:\mathcal A_f\longrightarrow\mathcal A_{G,E}^{[V]},
\qquad g_{a,b}\longmapsto (Y\mapsto aY+L(b))
\]
is a surjective homomorphism whose kernel consists exactly of translations from \(V\). Any inverse image under \(\Phi\) of a specified quotient transformation is called a \textbf{lift}; each quotient transformation has exactly \(|V|\) lifts. Thus there is an exact sequence
\begin{equation}\label{fas:eq-lifting-exact}
1\longrightarrow V\longrightarrow\mathcal A_f
\xrightarrow{\ \Phi\ }\mathcal A_{G,E}^{[V]}\longrightarrow1.
\end{equation}
\end{proposition}
\begin{proof}
Conjugating translations shows that \(g_{a,b}\in\mathcal A_f\) implies \(aV=V\). Proposition~\ref{v70:prop-scalar-support}, together with the direct calculation when \(V=\{0\}\), gives \(L(aX)=aL(X)\). Additivity then gives
\[
L(aX+b)=aL(X)+L(b).
\]
The derivative \(L'\) is a nonzero constant, so each root fiber is separable, and \(\nu_f(x)=\nu_G(L(x))\) for \(x\in E\).
Thus \(g_{a,b}\) induces a quotient transformation preserving the root multiplicities of \(G\), with scalar multiplier determined by \(\deg G\).
Additivity and this scalar compatibility also show that \(\Phi\) preserves composition. The quotient transformation is the identity if and only if \(a=1,L(b)=0\), that is, \(b\in V\).

Conversely, given \(a,d\) in~\eqref{fas:eq-compatible-quotient}, take a root \(\alpha\) of \(f\).
Then \(aL(\alpha)+d\) is again a root of \(G\); since \(f=G(L)\) splits completely in \(E\), there is a root \(\alpha'\) of \(f\) with
\(L(\alpha')=aL(\alpha)+d\). Put \(b=\alpha'-a\alpha\); then \(L(b)=d\). For every \(x\in E\),
\[
\nu_f(ax+b)=\nu_G(aL(x)+d)=\nu_G(L(x))=\nu_f(x).
\]
Hence \(g_{a,b}\in\mathcal A_f\). All solutions of \(L(b')=d\) form \(b+V\), giving exactly \(|V|\) lifts.
The inclusion in~\eqref{fas:eq-lifting-exact} sends \(v\in V\) to \(X\mapsto X+v\); its image is exactly \(\ker\Phi\), and \(\Phi\) is surjective.
\end{proof}

\begin{example}[Quotient scalings need not lift]\label{fas:ex-nonlifting}
Over \(\F_2\), take \(L=X^2+X\) and \(G=Y^3+1\), and put
\[
f=G(L)=X^6+X^5+X^4+X^3+1.
\]
Since \(f'=L^2\) and \(f\) equals one where \(L=0\), the polynomial \(f\) has six distinct roots.
The group \(V_f\) contains \(\F_2\), and its size is a power of two dividing six, so \(V_f=\F_2\) and \(\mathcal N_E(V_f)=\{1\}\).
The root set of \(G\) is \(\mu_3\), with root sum zero. If \(Y\mapsto aY+d\) preserves these three roots, comparing root sums gives \(3d=d=0\); thus its full affine stabilizer consists exactly of three scalings.
The two nonidentity scalings do not preserve \(V_f\), so cannot lift. The full affine stabilizer of the original polynomial consists only of the two translations from \(\F_2\).
\end{example}

Return now to irreducible \(f\), and take \(E=k(\alpha)\).
Put \(Q=|V_f|\), \(D=n/Q\), and \(\beta=L_V(\alpha)\). Then
\begin{equation}\label{fas:eq-irreducible-quotient}
\beta\ne0,\qquad [k(\beta):k]=D,\qquad f=G(L_V),
\end{equation}
where \(G\) is the monic minimal polynomial of \(\beta\), and \(V_G=\{0\}\).
Indeed, \(\beta=0\) would give \(\alpha\in V_f\), so a translation would send the root \(\alpha\) to zero, contradicting irreducibility of \(f\) of degree at least two.
The fibers of \(L_V\) on the root set are exactly the translation cosets. Frobenius is a single cycle on the \(D\) cosets, so \(\beta\) has degree \(D\).
Multiplying their root polynomials gives \(f=G(L_V)\), and \(V_G=\{0\}\) was proved in Corollary~\ref{ffs:cor-maximal-additive}.

\begin{lemma}[The scalar stabilizer in the irreducible case]\label{fas:lem-scalar-normalizer}
With the notation above, \(\mathcal N_E(V_f)=E^\times\) when \(Q=1\), and \(\mathcal N_E(V_f)=\F_p^\times\) when \(Q=p\).
In characteristic two, if \(Q=4\), write
\[
L_V=X^4+A_2X^2+A_1X,\qquad A_1\ne0.
\]
Then
\[
h_V=\begin{cases}3,&A_2=0,\\1,&A_2\ne0.\end{cases}
\]
When \(A_2=0\), we have \(V_f=b\F_4\) for some \(b\ne0\).
\end{lemma}
\begin{proof}
The case \(Q=1\) follows from the convention for the fixed field \(E\). When \(Q=p\), the space \(V_f\) has dimension one over \(\F_p\), so its maximal field of scalars is \(\F_p\).
If \(Q=4\), the support criterion in Proposition~\ref{v70:prop-scalar-support} gives maximal field of scalars \(\F_2\) when \(A_2\ne0\). When \(A_2=0\), it remains only to check that \(\F_4\subseteq E\).
Take \(b\in V_f\setminus\{0\}\). Then \(b^3=A_1\), and
\(L_V(bY)=b^4(Y^4+Y)\), so \(V_f=b\F_4\).
Since \(V_f\subseteq E\), we have \(\F_4\subseteq E\). The maximal field of scalars is therefore \(\F_4\), whose nonzero elements give all three scalars.
\end{proof}

\subsection{Centered Multiplicative Orders and Recovery of the Full Stabilizer}
\label{fas:subsec-centered-recovery}

The residual group on the translation quotient has no nonzero translations, so if nontrivial it has a unique fixed center.
To test candidate centers using root values, first express the multiple equal-difference structure of Section~\ref{sec:cyclotomic-med} as a finite divisor formula.

\begin{proposition}[Counting stabilizing multipliers from multiplicative orders]\leavevmode\phantomsection\label{fas:prop-order-multiplier}
For \(\gcd(N,q)=1\), put \(d_N=\ord_N(q)\), with the convention \(\ord_1(q)=1\), and define
\begin{equation}\label{fas:eq-order-multiplier}
s_q(N)=\max\left\{r:
 r\mid\gcd(N,d_N),\quad \ord_{N/r}(q)=d_N/r\right\}.
\end{equation}
If \(\gamma\ne0\) has multiplicative order \(N\), then the multiplicative stabilizer of its minimal polynomial over \(k\) has exactly \(s_q(N)\) elements.
This number is \(\kappa_C\) for the corresponding cyclotomic coset and can be computed from multiplicative orders at the small modulus in \eqref{v124:eq-small-step}.
\end{proposition}
\begin{proof}
Take \(\gamma\) as a primitive \(N\)th root of unity. The root index set is \(C=\langle q\rangle\subseteq\Z/N\Z\), of size \(d_N\).
For \(r\mid N\), reduction modulo \(N/r\) maps \(C\) to an orbit of size \(\ord_{N/r}(q)\), with all nonempty fibers of equal size.
Every fiber of the reduction map \(\Z/N\Z\to\Z/(N/r)\Z\) has \(r\) elements, so \(C\) is a union of complete fibers if and only if
\(d_N=r\ord_{N/r}(q)\). This is equivalent to \(C\) being preserved by the additive subgroup of order \(r\) in \(\Z/N\Z\).
By Theorem~\ref{thm:med-coset-partitions}, such \(r\) are exactly the divisors of \(\kappa_C\), with maximum \(\kappa_C\).
The correspondence between index translations and multiplication of root values proves the assertion. When \(N=1\), both sides are one.
\end{proof}

\begin{theorem}[Recovering the full affine stabilizer from root values]\label{fas:thm-full-recovery}
Let \(f\in k[X]\) be monic and irreducible with \(n=\deg f\ge2\), and determine \(V=V_f\) in \(E=k(\alpha)\) by Theorem~\ref{fas:thm-translation-difference}.
Put \(L=L_V\), \(Q=|V|\), \(\beta=L(\alpha)\), and \(D=n/Q\), and take \(\mathcal N_E(V),h_V\) from Lemma~\ref{fas:lem-scalar-normalizer}.
\begin{enumerate}
\item If \(D=1\), then \(\beta\in k^\times\), \(f=L-\beta\), and
\begin{equation}\label{fas:eq-one-fiber-group}
\mathcal A_f
=\{X\mapsto aX+(1-a)\alpha+v:
 a\in\mathcal N_E(V),\ v\in V\},
\quad |\mathcal A_f|=Qh_V.
\end{equation}
\item If \(D>1\), define, for every \(z\in k\),
\begin{equation}\label{fas:eq-center-test}
N_z=\ord_{E^\times}(\beta-z),\qquad
 t_z=\gcd\bigl(h_V,s_q(N_z)\bigr).
\end{equation}
\leavevmode\phantomsection\label{fas:def-centered-order}\ignorespaces
Here \(N_z\) is obtained by first translating the quotient root \(\beta\) to \(\beta-z\), then taking its multiplicative order; it is called the \textbf{multiplicative order after centering} at the candidate center \(z\).
At most one \(z\in k\) has \(t_z>1\). If no such candidate exists, \(\mathcal A_f\) consists only of translations from \(V\).
Otherwise, denote it by \(\delta\) and put \(t=t_\delta\). Then
\begin{equation}\label{fas:eq-recovered-group}
\boxed{\quad
\mathcal A_f=\{g_{a,b}:a,b\in E,\ a^t=1,\ L(b)=(1-a)\delta\}.
\quad}
\end{equation}
Each slope has exactly \(Q\) translation parameters, and \(|\mathcal A_f|=Qt\), with \(t\mid D\).
There is also a monic irreducible \(H\in k[Y]\) such that
\begin{equation}\label{fas:eq-recovered-composition}
f=H\bigl((L-\delta)^t\bigr),\qquad H(0)\ne0,\quad \deg H=D/t.
\end{equation}
If \(p\nmid D\), the only possible center for nontrivial slopes is
\begin{equation}\label{fas:eq-trace-center}
\delta=D^{-1}\Tr_{k(\beta)/k}(\beta).
\end{equation}
In this case only this value needs to be tested for \(t_\delta>1\).
\end{enumerate}
\end{theorem}
\begin{proof}
Equation~\eqref{fas:eq-irreducible-quotient} gives irreducible \(G\), \(f=G(L)\), \(\deg G=D\), and \(V_G=\{0\}\).
If \(D=1\), the root set is \(\alpha+V\). Any slope preserving it must belong to \(\mathcal N_E(V)\), and the image of \(\alpha\) must be \(\alpha+v\), giving~\eqref{fas:eq-one-fiber-group}; conversely, every listed transformation preserves this coset.

Suppose \(D>1\), and denote the root set of \(G\) by \(B\). Proposition~\ref{fas:prop-compatible-lifting} identifies \(\mathcal A_f/V\) with the compatible affine stabilizer of \(B\).
It has no nonzero translations, so the slope map embeds it into \(E^\times\); hence it is cyclic.
If nontrivial, its generator has slope different from one and thus a unique fixed point; this is also the unique common fixed point of the whole group, its center.
Frobenius preserves both \(B\) and \(V\), so normalizes this group. Uniqueness forces the center to be Frobenius-fixed, hence to lie in \(k\).

For \(z\in k\), the element \(\beta-z\) still has degree \(D>1\), so it is nonzero and \(\ord_{N_z}(q)=D\).
The group of scalings centered at \(z\) that preserve \(B\) is the multiplicative stabilizer of \(G(Y+z)\), of size \(s_q(N_z)\).
It and \(\mathcal N_E(V)\) are subgroups of the cyclic group \(E^\times\), so their intersection has size \(t_z\).
Proposition~\ref{fas:prop-compatible-lifting} shows that each quotient transformation in this intersection lifts.
Thus every \(t_z>1\) gives the center of the same nontrivial quotient group, so there is at most one such \(z\).
If there is no such \(z\), the quotient group is trivial. If one exists, its full slope group is \(\mu_t\), and the lifting condition is \(L(b)=(1-a)\delta\).
This proves the full group formula and the number of lifts for each slope.

Since \(G\) is irreducible with \(D>1\), the base-field point \(\delta\) is not a root. In the coordinate \(Y-\delta\), the group \(\mu_t\) acts freely on all roots,
so \(t\mid D\), and \(G(Y)=H((Y-\delta)^t)\), with \(H(0)\ne0\).
The polynomial \(H\) must be irreducible, since otherwise substitution gives a nontrivial factorization of \(G\). Substituting \(L\) then yields~\eqref{fas:eq-recovered-composition}.
Finally, take a stabilizing slope \(a\ne1\). Comparing root sums gives
\[
\Tr_{k(\beta)/k}(\beta)
=a\Tr_{k(\beta)/k}(\beta)+D(1-a)\delta.
\]
If \(p\nmid D\), cancel \(1-a\) and divide by the nonzero image of the integer \(D\) in \(k\), obtaining~\eqref{fas:eq-trace-center}.
\end{proof}

Here \(\delta\) is the center in the quotient coordinate \(L\); a center \(c\) in the original coordinate satisfies \(L(c)=\delta\) and need not lie in \(k\).
The trace formula gives only the unique candidate to test; it does not ensure that the corresponding slope group is nontrivial.
When \(p\mid D\), test all \(q\) values in \(k\) using~\eqref{fas:eq-center-test}.
Compared with general parameter enumeration in Section~\ref{ffs:sec-computation}, this route first uses irreducibility to reduce the translation candidates, then computes the centered multiplicative orders of quotient roots.
Its input can be represented entirely in \(k[X]/(f)\); Frobenius powers, multiplicative orders, and finite-field operations still have to be computed.
The resulting group is independent of the chosen conjugate root, since changing the root applies Frobenius simultaneously to all intermediate data.

\leavevmode\phantomsection\label{fas:coefficient-recovery}\ignorespaces If the coefficients of the outer polynomial \(G\) have already been recovered from \(f=G(L)\), the same \(t_z\) can also be computed directly from coefficient support.
Suppose \(D>1\), and, for \(z\in k\), put
\[
r_z=\gcd\{i:1\le i\le D,\ [Y^i]G(Y+z)\ne0\}.
\]
Irreducibility gives \(G(z)\ne0\), so \(G(Y+z)\) has nonzero constant term.
Also, \(G\) is separable, so \(G'\ne0\). The translated polynomial still has at least one nonzero term whose exponent is not divisible by \(p\), hence \(p\nmid r_z\).
Proposition~\ref{ffs:prop-multiplier-support} now gives \(M_{G(Y+z)}=\mu_{r_z}\).
It and \(\mathcal N_E(V)\) are subgroups of the cyclic group \(E^\times\); their intersection is precisely the compatible slope group centered at \(z\). Thus
\[
\boxed{\quad t_z=\gcd(h_V,r_z).\quad}
\]
On the other hand, the minimal polynomial of \(\beta-z\) is \(G(Y+z)\), and Proposition~\ref{fas:prop-order-multiplier} gives \(r_z=s_q(N_z)\), so this agrees with~\eqref{fas:eq-center-test}.
When the outer coefficients are known, each candidate center can therefore be tested simply by translating the polynomial and taking the greatest common divisor of its nonzero-term exponents, without separately computing the multiplicative order of \(\beta-z\).

\subsection{Possible Group Types and Classification in Prime Degree}
\label{fas:subsec-prime-classification}

Below, \(C_j\) denotes a cyclic group of order \(j\), \(S_j\) the symmetric group of all permutations of \(j\) elements, and \(A_j\) its alternating subgroup of even permutations; these symbols are distinguished from cyclotomic-coset notation by context.

\begin{corollary}[Full group types and the four-element exception]\leavevmode\phantomsection\label{fas:cor-group-types}
Retain \(Q,D\) from Theorem~\ref{fas:thm-full-recovery}, and write \(t=|\mathcal A_f|/Q\) for the order of the slope image.
\begin{enumerate}
\item If \(Q=1\), then \(\mathcal A_f\) is cyclic of order \(t\mid n\), with \(p\nmid t\).
\item If \(Q=p\), then \(\mathcal A_f\simeq C_p\rtimes C_t\), where \(t\mid p-1\); slope \(a\) acts on the translation group by \(v\mapsto av\), and only \(a=1\) fixes that group pointwise.
If \(D>1\), we also have \(t\mid D\); if \(D=1\), then \(t=p-1\).
\item If \(p=2,Q=4\), then \(\mathcal A_f\) is isomorphic to \(C_2\times C_2\) or \(A_4\).
The latter occurs only when \(n=4\) and \(a\) in \(q=2^a\) is odd.
\end{enumerate}
\end{corollary}
\begin{proof}
The first two assertions follow directly from the finite affine group structure, the scalar stabilizer, and Theorem~\ref{fas:thm-full-recovery}.
When \(Q=4\), the slope order is either one or three. A slope of order three cyclically permutes the three nonzero translations; the resulting group \((C_2)^2\rtimes C_3\) acts on four points by even permutations and has order 12, so is \(A_4\).

If a slope of order three occurs, Lemma~\ref{fas:lem-scalar-normalizer} gives \(V=b\F_4\).
Corollary~\ref{ffs:cor-four-block} also shows that \(\sigma|_V\) has order two and \(D=n/4\) is odd.
If \(a\) were even, the \(q\)th-power map would be the identity on \(\F_4\), and \(\sigma\) on \(b\F_4\) would correspond to multiplication by \(b^{q-1}\in\F_4^\times\), of order one or three, a contradiction.
Thus \(a\) is odd.
If \(D>1\), then \(3\mid D\); write \(D=3j\), with \(j\) odd.
The quotient polynomial is \(G(Y)=H((Y-\delta)^3)\), where \(H\) is irreducible of degree \(j\).
Take a root \(\theta\in\F_{q^j}^\times\) of \(H\). Since \(q^j\equiv2\pmod3\), cubing is bijective on \(\F_{q^j}^\times\), so some \(\xi\in\F_{q^j}\) satisfies \(\xi^3=\theta\).
Then \(\delta+\xi\) is a root of \(G\) of degree at most \(j<D\), contradicting irreducibility of \(G\).
Hence \(D=1,n=4\).
\end{proof}

\begin{theorem}[Symmetric irreducible polynomials of degree equal to the characteristic]\label{fas:thm-prime-characteristic}
Let \(f\in k[X]\) be monic and irreducible with \(\deg f=p\). Its full affine stabilizer is nontrivial if and only if there exist \(b\in k^\times\) and \(a_0\in k\) such that
\begin{equation}\label{fas:eq-artin-schreier-classification}
f(X)=X^p-b^{p-1}X-a_0,
\qquad \Tr_{k/\F_p}(a_0/b^p)\ne0.
\end{equation}
In this case \(V_f=b\F_p\), \(\mathcal A_f\simeq\operatorname{AGL}_1(\F_p)\), and \(|\mathcal A_f|=p(p-1)\), while \(\mathcal A_f(k)\) consists exactly of the \(p\) translations from \(b\F_p\).
The number of such monic irreducible polynomials is
\begin{equation}\label{fas:eq-as-count}
\frac{q(q-1)}p.
\end{equation}
\end{theorem}
\begin{proof}
If \(V_f=\{0\}\), the preceding corollary gives \(t\mid p\) and \(p\nmid t\), so the full group is trivial.
In the nontrivial case, necessarily \(|V_f|=p\), and all roots form a single translation coset. Frobenius acts on this coset as a \(p\)-cycle, so the order of its linear part divides \(p\); this linear part also belongs to \(\F_p^\times\), so its order divides \(p-1\) and it is the identity. Thus Frobenius fixes \(V_f\) pointwise, that is, \(V_f\subseteq k\).
Writing \(V_f=b\F_p\), with \(b\in k^\times\), the root polynomial must be \(X^p-b^{p-1}X-a_0\), with \(a_0\in k\).

For a root \(\alpha\), put \(y=\alpha/b\). The equation becomes \(y^p-y=a_0/b^p\). Iterating \(p\)th powers gives
\[
y^q-y=\Tr_{k/\F_p}(a_0/b^p).
\]
If the right-hand side is nonzero, \(q\)-Frobenius is a \(p\)-cycle on \(y+\F_p\), so the polynomial is irreducible; if it is zero, all roots lie in \(k\).
The full group follows from~\eqref{fas:eq-one-fiber-group}. The base-field group contains \(p\) translations, and its order divides \(p\) by Theorem~\ref{ffs:thm-arithmetic-cyclic}, so it consists exactly of these translations.

There are \((q-1)/(p-1)\) distinct one-dimensional subspaces \(b\F_p\subseteq k\) over the prime field, giving distinct linear coefficients \(-b^{p-1}\).
Choose one \(b\) for each line. The nonzero trace condition allows exactly \(q(p-1)/p\) choices of \(a_0\).
Multiplying \(b\) by \(\lambda\in\F_p^\times\) merely multiplies the trace by \(\lambda^{-1}\), preserving nonvanishing, so this count is independent of the representative.
Multiplying the two counts gives~\eqref{fas:eq-as-count}.
\end{proof}

\begin{theorem}[Symmetric irreducible polynomials of prime degree different from the characteristic]\label{fas:thm-prime-tame}
Let \(\ell\ne p\) be prime, and let \(f\in k[X]\) be monic and irreducible with \(\deg f=\ell\). Its full affine stabilizer is nontrivial if and only if
\begin{equation}\label{fas:eq-prime-binomial-classification}
\ell\mid q-1,\qquad
f(X)=(X-c)^\ell-a_0,
\quad c\in k,\quad a_0\in k^\times\setminus(k^\times)^\ell.
\end{equation}
In this case \(V_f=\{0\}\), the full group is \(\mu_\ell\) centered at \(c\), and all its parameters lie in \(k\).
The number of such polynomials is
\begin{equation}\label{fas:eq-prime-binomial-count}
\begin{cases}
q(q-1)(\ell-1)/\ell,&\ell\mid q-1,\\
0,&\ell\nmid q-1.
\end{cases}
\end{equation}
\end{theorem}
\begin{proof}
Since \(p\nmid\ell\), the group \(V_f\) is trivial. If the full group is nontrivial, its slope order divides \(\ell\), hence equals \(\ell\).
Its unique center \(c\) is Frobenius-fixed, so lies in \(k\); it is not a root of the irreducible polynomial.
All roots therefore form a single multiplicative orbit centered at \(c\), and \(f=(X-c)^\ell-a_0\), with \(a_0\in k^\times\).
If \(\ell\nmid q-1\), the \(\ell\)th-power map on \(k^\times\) is bijective, giving a base-field root, a contradiction.

Suppose \(\ell\mid q-1\). If \(a_0\) is an \(\ell\)th power, there is plainly a base-field root.
If not, take \(z^\ell=a_0\). Then \(z^q/z\in\mu_\ell\subseteq k\), and this ratio differs from one.
Since \(\ell\) is prime, the ratio has order \(\ell\), so the Frobenius period of \(z\) is \(\ell\), proving irreducibility.
The \(\ell\) scalings already exhibited form the full group, since its slope order cannot exceed \(\ell\).
Finally, \(c\) is uniquely determined by the coefficient \(-\ell c\) of \(X^{\ell-1}\), after which \(a_0\) is unique as well.
There are \((q-1)(\ell-1)/\ell\) nonzero elements that are not \(\ell\)th powers, giving the stated count.
\end{proof}

These two theorems classify irreducible polynomials of prime degree with nontrivial full affine stabilizer.
Additive cosets give the Artin--Schreier type, and centered multiplicative orbits give translated binomials, providing a symmetry characterization of the two constructions in Chapter~\ref{chap:trace}.

\begin{corollary}[Classification of irreducible quartics in characteristic two]\leavevmode\phantomsection\label{fas:cor-quartic-coefficients}
Let \(k\) have characteristic two, and let \(f=X^4+a_3X^3+a_2X^2+a_1X+a_0\in k[X]\) be irreducible. Then
\[
\mathcal A_f\simeq
\begin{cases}
\{1\},&a_3\ne0,\\
C_2\times C_2,&a_3=0,\ a_2\ne0,\\
A_4,&a_3=a_2=0.
\end{cases}
\]
In the last two cases, \(|V_f|=4\) and \(a_1a_0\ne0\); the final case requires \([k:\F_2]\) to be odd.
\end{corollary}
\begin{proof}
A nonzero stabilizing translation partitions the four roots into two pairs, each with sum equal to the translation parameter. Thus the sum of all four roots is zero, that is, \(a_3=0\).
Conversely, if four distinct roots sum to zero, translating one fixed root to zero gives a set \(\{0,b,u,b+u\}\), where \(b,u\) are linearly independent over \(\F_2\).
Thus the original root set is a coset of a four-element translation space, and \(|V_f|=4\).
When \(a_3\ne0\), the group \(V_f\) is trivial and the slope order is an odd divisor of four, so the full group is trivial.
When \(a_3=0\), we have \(f=L_V-\beta\), with \(L_V=X^4+a_2X^2+a_1X\).
Separability and irreducibility give \(a_1a_0\ne0\); Lemma~\ref{fas:lem-scalar-normalizer} and~\eqref{fas:eq-one-fiber-group} give the two group types. The restriction on the base-field degree follows from Corollary~\ref{fas:cor-group-types}.
\end{proof}

\subsection{Four Frobenius Orbits of Roots of Order 13}
\label{fas:subsec-thirteen}

\begin{example}[Equal multiplicative orders with different affine groups]\label{fas:ex-thirteen}
Over \(\F_3\), take \(\alpha^3-\alpha=1\). The polynomial \(X^3-X-1\) has no base-field root and is therefore irreducible; hence
\[
\alpha^9=\alpha-1,\qquad
\alpha^{13}=\alpha(\alpha+1)(\alpha-1)=1.
\]
Since \(\alpha\ne1\), it is a primitive 13th root of unity.
Relative to \(\zeta=\alpha\), the nonzero \(3\)-cyclotomic cosets modulo 13 and their root polynomials are as follows.
The groups in the table are full affine stabilizers allowing parameters in the splitting field.
\begin{center}\small
\renewcommand{\arraystretch}{1.18}
\begin{tabular}{@{}llll@{}}
\toprule
Cyclotomic coset&Monic irreducible polynomial&\(V_f\)&Group type of \(\mathcal A_f\)\\
\midrule
\(\{1,3,9\}\)&\(X^3-X-1\)&\(\F_3\)&\(S_3\)\\
\(\{2,6,5\}\)&\(X^3+X^2+X-1\)&\(\{0\}\)&Trivial\\
\(\{4,12,10\}\)&\(X^3+X^2-1\)&\(\{0\}\)&Trivial\\
\(\{7,8,11\}\)&\(X^3-X^2-X-1\)&\(\{0\}\)&Trivial\\
\bottomrule
\end{tabular}
\end{center}

The polynomial column follows by successive reductions using \(\alpha^{j+3}=\alpha^{j+1}+\alpha^j\).
For example, the elementary symmetric sums of the three roots in the first row are \(0,-1,1\); after squaring each root, the sums become \(-1,1,1\), giving the second row.
The root indices in the third row are the negatives of those in the first, so its polynomial is the monic reciprocal of the first polynomial.
The fourth row likewise has the negatives of the second row's indices and is obtained by taking the monic reciprocal.
All four cosets have length three, so their polynomials are irreducible and their product is \(\Phi_{13}\).

For the first row, \(b_3=\alpha^3-\alpha=1\), so \(V_f=\F_3\), \(L_V=X^3-X\), and \(\beta=1,D=1\).
Theorem~\ref{fas:thm-full-recovery} gives
\[
\mathcal A_f=\{X\mapsto aX+(1-a)\alpha+v:
 a\in\F_3^\times,\ v\in\F_3\}\simeq S_3.
\]
The base-field group, however, contains only the three translations. For the second row, take the root \(\alpha^2\); its test difference is
\[
(\alpha^2)^3-\alpha^2=1-\alpha,
\qquad(1-\alpha)^3=-\alpha\ne\pm(1-\alpha).
\]
The translation criterion therefore fails, and the classification in degree equal to the characteristic shows that the full group is trivial.
The third and fourth rows have nonzero quadratic coefficients, so are not of the form~\eqref{fas:eq-artin-schreier-classification}; their full groups are also trivial.

All four rows have \(q=3\), root order 13, and irreducible degree three; since \(\gcd(13,3)=1\), their multiplicative stabilizers are all trivial as well.
Their full affine groups nevertheless differ, because additive relations among root values are not determined by these multiplicative-order data.
The first and third rows are reciprocal but have different group types; the third has a trivial group despite having a zero coefficient between its leading and constant terms.
Thus inversion need not preserve affine stability, and the full affine group does not explain every coefficient cancellation.
\end{example}

An inner component obtained from a translation or affine stabilizer is maximal only within the corresponding symmetry class and does not exhaust general polynomial decompositions.
For irreducible polynomials, the next chapter gives a criterion independent of affine symmetry, using low-degree representations of subfields.

\par\Needspace{6\baselineskip}
\section*{Exercises}
\phantomsection
\addcontentsline{toc}{section}{Exercises}

\begin{exercise}[Root sets and root multiplicities]\label{ffs:ex-weighted-two}
In \(\F_3[X]\), take \(f=X(X-1)^2\). Find \(V_f,M_f,\mathcal A_f\), and \(\Psi_f\).
Explain why the affine stabilizer of the root set \(\{0,1\}\), with multiplicities ignored, is larger.
\end{exercise}

\begin{exercise}[Recovering the common \(p\)-power from differences]\label{ffs:ex-power-four}
In \(\F_2[X]\), take \(f=(X^2+X+1)^4\). Without first finding roots, compute \(e,\Psi_f,L_{V_f}\), and \(M_f\), and find the full affine stabilizer.
\end{exercise}

\begin{exercise}[A slope group of order six with a nonzero center]\label{ffs:ex-six-slopes}
In \(\F_5[X]\), take \(f=(X-2)^6-1\). Find \(e,V_f,M_f,t,\delta\), and the normal form.
Prove that the full set of slopes is \(\mu_6\), and find the order of \(\mathcal A_f(\F_5)\).
Use Corollary~\ref{v67:cor-affine-support-bounds} to compute \(E_\Gamma\) and \(E_0\), and compare their cardinalities with \(|\operatorname{Supp}_X(f)|\) and \(\deg f/\operatorname{sp}(f)+1\).
\end{exercise}

\begin{exercise}[Two factorization types with the same inner polynomial]\label{ffs:ex-two-additive-types}
Over \(\F_3\), put \(L=X^3-X\), \(P_1=Y^2+1\), and \(P_2=Y^2+Y+2\).
Use return translations to determine the irreducible factor degrees of \(P_1(L)\) and \(P_2(L)\), and write the full factorization of \(P_1(L)\).
\end{exercise}

\begin{exercise}[Linear coordinates for a four-cycle]\label{ffs:ex-four-cycle}
On \(V=\F_2^2\), put
\[
A(x,y)=(x+y,y),\qquad b=(0,1),\qquad R(v)=Av+b.
\]
Find the cycles of \(R\) and \(N_1,N_2,N_4\), then recover the cycle counts using~\eqref{ffs:eq-return-mobius}.
Relate this to the root fiber of \(X^4+X+1\).
\end{exercise}

\begin{exercise}[Splitting-field and base-field translations]\label{ffs:ex-degree-restrictions}
Let \(f\in\F_2[X]\) be irreducible.
For \(\deg f=6,12,24\), respectively, list the values of \(|V_f|\) allowed by Theorem~\ref{ffs:thm-geometric-bound} and Corollary~\ref{ffs:cor-four-block}.
Only list the stabilizer sizes allowed by these necessary conditions; there is no need to prove that every size is attained.
\end{exercise}

\begin{exercise}[The center equation and the base-field group]\label{ffs:ex-rational-centers}
Suppose \(t>1\). Prove that if the special orbit contains a base-field point, it contains exactly \(|V_f\cap k|\) such points.
Use this to explain why \(U\) in Example~\ref{ffs:ex-quartic-center} is defined over the base field but has no base-field affine center.
\end{exercise}

\begin{exercise}[A product may acquire new symmetries]\label{ffs:ex-product-new-symmetry}
In \(\F_3[X]\), put \(L=X^3-X\), \(f_1=L-1\), and \(f_2=L+1\).
Find \(\mathcal A_{f_1}\), \(\mathcal A_{f_2}\), and \(\mathcal A_{f_1f_2}\), and explain why the stabilizer of a product is generally not the intersection of the stabilizers of its factors.
\end{exercise}

\begin{exercise}[Trivial affine stabilizer with nontrivial composition]\label{ffs:ex-not-all-compositions}
In \(\F_7[X]\), put \(R=X^3+X\) and \(f=R^2+R=X^6+2X^4+X^3+X^2+X\).
Prove that \(f\) has six distinct roots and trivial \(\mathcal A_f\), yet has a nontrivial composition \(f=F(R)\) with outer degree two and inner degree three.
\end{exercise}

\begin{exercise}[A semi-invariant basis and common support]\leavevmode\phantomsection\label{fas:ex-semi-support}
Over \(\F_3\), put \(L=X^3-X\), and take the group
\(\Gamma=\{X\mapsto aX+b:a\in\{1,-1\},\ b\in\F_3\}\).
Find a basis, dimension, and common support set for each of the two spaces of polynomials of degree at most 12 satisfying \(P(aX+b)=P(X)\) and \(P(aX+b)=aP(X)\), respectively.
\end{exercise}

\begin{exercise}[A translation direction outside the base field]\label{fas:ex-nonrational-line}
In \(\F_3[X]\), take
\[
f=(X^3+X)^2+1=X^6-X^4+X^2+1.
\]
Prove that \(f\) is irreducible. Take a root \(\alpha\), put \(i=\alpha^3+\alpha\), and compute \(b_3=\alpha^9-\alpha\), \(V_f\), \(\mathcal A_f\), and \(\mathcal A_f(\F_3)\).
Then recover the same full group from the center and multiplicative orders of the translation quotient.
\end{exercise}

\begin{exercise}[Counts in prime degree and the quartic exception]\label{fas:ex-prime-counts}
Count the monic irreducible cubics over \(\F_3\) with nontrivial full affine stabilizer, and list them all.
Then count such monic irreducible cubics over \(\F_4\).
Finally, classify the three monic irreducible quartics over \(\F_2\) by their full affine groups.
\end{exercise}

\begin{exercise}[Four-element translations with a nontrivial outer component]\label{fas:ex-twelve-four-translations}
In \(\F_2[X]\), put
\[
L=X^4+X,\qquad G=Y^3+Y^2+1,\qquad
f=G(L)=X^{12}+X^9+X^8+X^6+X^3+X^2+1.
\]
Prove that \(f\) is irreducible, find its full affine stabilizer, and explain why the order-three slopes preserving \(\F_4\) do not preserve this root set.
\end{exercise}

\chapter{Subfield Representations and Polynomial Decomposition}
\label{chap:subfield-representations}
\phantomsection\label{unit:subfield-trace-blocks}
The central question of this chapter is concrete: why can the same subfield require polynomial representatives of different degrees in different power bases?
We first compute a quartic example, then use degree profiles, root blocks, and relative traces to explain it.

\begin{example}[One four-element subfield, quadratic and cubic representatives]\label{sfr:ex-quartic-comparison}
Take
\[
f_1=X^4+X+1,\qquad f_2=X^4+X^3+1
\]
in \(\F_2[X]\). Neither has a root in \(\F_2\), and their remainders modulo the unique monic irreducible quadratic \(X^2+X+1\) are \(1,X\), respectively.
Thus both are irreducible. If \(f_1(\alpha)=0\), then \(\gamma=\alpha^{-1}\) satisfies \(f_2(\gamma)=0\).
The two quotient rings therefore represent the same field \(E=\F_{16}\), with the same subfield \(K=\F_4\).

In the power basis \((1,\alpha,\alpha^2,\alpha^3)\), put \(\beta=\alpha^2+\alpha\). Since \(\alpha^4=\alpha+1\),
\[
\beta^2+\beta+1=0,\qquad
K=\{0,1,\beta,\beta+1\}.
\]
Here \(\beta\notin\F_2\), since otherwise \(\alpha\) would satisfy a quadratic equation over the base field.
Likewise, no nonconstant linear expression \(a_0+a_1\alpha\) can belong to \(K\), since that would imply \(\alpha\in K\).
The smallest degree of a representative of a generator of \(K/\F_2\) is therefore two. We can already see that
\[
f_1(X)=(X^2+X)^2+(X^2+X)+1.
\]

Now use the power basis \((1,\gamma,\gamma^2,\gamma^3)\). The relations
\[
\gamma^4=\gamma^3+1,\quad
\gamma^6=\gamma^3+\gamma^2+\gamma+1,\quad
\gamma^8=\gamma^3+\gamma^2+\gamma
\]
show that \(\eta=\gamma^3+\gamma\) satisfies \(\eta^2+\eta+1=0\).
Linear independence of the power basis gives \(\eta\notin\F_2\), so \(K=\{0,1,\eta,\eta+1\}\).
In fact, \(\alpha=\gamma^3+\gamma^2\) and \(\alpha^2=\gamma^2+\gamma\), so \(\beta=\eta\): we are even comparing two representatives of the same subfield element.
To exclude quadratic representatives, write \(z=a_0+a_1\gamma+a_2\gamma^2\), with \(a_i\in\F_2\). Direct calculation gives
\[
z^4+z=(a_1+a_2)\gamma^3+(a_1+a_2)\gamma+a_1.
\]
Hence \(z\in K\), equivalently \(z^4=z\), if and only if \(a_1=a_2=0\).
In this basis a subfield generator requires a representative of degree at least three, and \(\eta\) attains that degree.

\[
\begin{array}{c|c|c}
\text{Power-basis generator}&\text{A generator of }\F_4/\F_2&\text{Smallest degree}\\\hline
\alpha&\alpha^2+\alpha&2\\
\gamma=\alpha^{-1}&\gamma^3+\gamma&3
\end{array}
\]
These minimum degrees were obtained directly, without a composition criterion.
The abstract subfield is unchanged; its position in the successive low-degree power spaces has changed.
\end{example}

Example~\ref{v140:ex-quartic-depth} continues with the canonical lift of this same four-element subfield: its quadratic representative disappears modulo four, whereas a different Hensel lift preserving the composition equation retains a quadratic representative.

The example raises four questions that organize this chapter:
\begin{itemize}
\item Why does the minimum representative degree depend on the power basis when the subfield is the same?
\item When is the lower bound \([E:K]\) attained in a general extension?
\item Why is equality equivalent to a composition of the defining polynomial?
\item How can relative traces detect equality and recover root blocks from finitely many data?
\end{itemize}
The degree profile records the subfield's intersections with low-degree spaces; the defect measures the departure from the optimal distribution.
Hankel recovery turns a finite relative-trace prefix back into a root-block equation.
We first establish the degree bound and its equality criterion, and then give the full table of tools.

\section{Low-Degree Representations of Subfields and Polynomial Decomposition}
\label{sfr:sec-subfield-representation}

Subfield classification determines which subfields a finite field contains, but does not specify the degrees needed to represent their elements in given coordinates.
Once a generator is fixed, the power basis assigns a degree to each field element; intersections of a subfield with the low-degree spaces record its distribution in these coordinates.
This section studies such representation degrees and connects them with partitions of irreducible root sets and polynomial decomposition.
Here the measure of complexity is the degree of a polynomial representative, rather than the number of field operations or bit operations required for computation.

\subsection{Power-Basis Degrees and the Degree Profile of a Subfield}
\label{sfr:subsec-filtration}

Fix \(k=\F_q\), let \(f\in k[X]\) be monic and irreducible with \(\deg f=n\ge2\), take a root \(\alpha\), and write
\[
E=k(\alpha)=\F_{q^n},\qquad \sigma(x)=x^q.
\]
Every \(u\in E\) has a unique polynomial representative \(h_u\in k[X]\) of degree less than \(n\), satisfying \(u=h_u(\alpha)\).
\leavevmode\phantomsection\label{sfr:def-representative-degree}\ignorespaces For \(u\ne0\), define its \textbf{degree in the power basis} by \(\deg_\alpha(u)=\deg h_u\), and set \(\deg_\alpha(0)=-\infty\).
For \(0\le t<n\), put
\[
\mathcal P_t(\alpha)=\operatorname{span}_k\{1,\alpha,\ldots,\alpha^t\},
\qquad \mathcal P_{-1}(\alpha)=\{0\}.
\]
\leavevmode\phantomsection\label{sfr:def-degree-filtration}\ignorespaces 
These spaces form the \textbf{power-basis degree filtration} of \(E\), the increasing sequence of subspaces
\[
\{0\}\subset\mathcal P_0(\alpha)\subset\cdots\subset\mathcal P_{n-1}(\alpha)=E.
\]
The function \(\deg_\alpha\) records the degree of the polynomial representative in the given power basis; its multiplicative behavior depends on whether reduction modulo \(f\) is needed.
If \(u,v\ne0\) and \(\deg_\alpha(u)+\deg_\alpha(v)<n\), no degree reduction modulo \(f\) is needed after multiplication, so
\begin{equation}\label{sfr:eq-degree-product}
\deg_\alpha(uv)=\deg_\alpha(u)+\deg_\alpha(v).
\end{equation}
Once the sum of degrees reaches \(n\), this equality may fail.

\leavevmode\phantomsection\label{sfr:def-subfield-coordinate}\ignorespaces By the subfield classification theorem~\ref{thm:subfields}, for every \(d\mid n\), the field \(E\) contains a unique subfield
\[
K_d=\F_{q^d}=\{u\in E:\sigma^d(u)=u\}.
\]
\begin{definition}\leavevmode\phantomsection\label{sfr:def-profile}
For \(d\mid n\) and \(-1\le t<n\), define
\[
U_{f,d}(t)=K_d\cap\mathcal P_t(\alpha),\qquad
\rho_{f,d}(t)=\dim_k U_{f,d}(t).
\]
The function \(\rho_{f,d}\) is called the \textbf{degree profile} of \(K_d\) in this power basis.
\leavevmode\phantomsection\label{sfr:def-generator-degree}\ignorespaces For \(d>1\), define its \textbf{minimum generator-representation degree} by
\begin{equation}\label{sfr:eq-kappa-definition}
\kappa_f(d)=\min\{\deg h:h\in k[X],\ \deg h<n,\ k(h(\alpha))=K_d\}.
\end{equation}
\end{definition}

The extension \(K_d/k\) has a generator, and every element has a unique representative of degree less than \(n\), so the minimum exists.
Allowing arbitrary degrees for \(h\) gives the same minimum: first take the remainder modulo \(f\).
The value \(\rho_{f,d}(t)\) records the dimension of the whole intersection space, whereas \(\kappa_f(d)\) requires a value generating the specified subfield; the two must not be confused.
When \(d\) is composite, the presence of an element of the intersection space outside the base field does not by itself imply that the space contains a generator of \(K_d/k\).

\begin{proposition}[Coordinate dependence and invariance]\leavevmode\phantomsection\label{sfr:prop-coordinate-invariance}
The degree profile and minimum generator-representation degree are independent of the choice of a conjugate root of \(f\).
They are also unchanged upon replacing the generator by \(\alpha'=a\alpha+b\), with \(a\in k^\times,b\in k\).
In general, they depend on the generator, not merely on the abstract field \(E\).
\end{proposition}
\begin{proof}
If \(\alpha'=\sigma^j(\alpha)\), then \(\sigma^j\) fixes \(k\) pointwise, preserves \(K_d\), and maps
\(\mathcal P_t(\alpha)\) bijectively onto \(\mathcal P_t(\alpha')\).
It also preserves degrees of elements over \(k\), so both kinds of data are unchanged.
If \(\alpha'=a\alpha+b\), binomial expansion gives \(\mathcal P_t(\alpha')\subseteq\mathcal P_t(\alpha)\),
and \(\alpha=a^{-1}(\alpha'-b)\) gives the reverse inclusion. The degree filtrations are therefore identical.
Example~\ref{sfr:ex-quartic-comparison} at the start of this chapter already shows that a general change of generator may change these data.
\end{proof}

When \(d=1\), we have \(\rho_{f,1}(t)=1\) for \(0\le t<n\); when \(d=n\),
\(\rho_{f,n}(t)=t+1\) and \(\kappa_f(n)=1\). Whenever a nontrivial decomposition is considered below, assume \(1<d<n\), and write
\[
m=n/d=[E:K_d].
\]

\subsection{The Least Representation Degree, Conjugate Root Blocks, and Decomposition}
\label{sfr:subsec-optimal-composition}

\begin{theorem}[Least degree and functional decomposition]\label{sfr:thm-optimal-composition}
Suppose \(1<d<n\), \(d\mid n\), and put \(m=n/d\). Then
\begin{equation}\label{sfr:eq-kappa-lower}
\kappa_f(d)\ge m.
\end{equation}
Equality holds if and only if there exist monic \(P,h\in k[X]\) such that
\begin{equation}\label{sfr:eq-optimal-composition}
f=P\circ h,\qquad \deg P=d,\quad\deg h=m.
\end{equation}
In this case \(P\) is irreducible, and \(h(\alpha)\) generates \(K_d/k\).
\end{theorem}
\begin{proof}
Suppose \(\beta=h(\alpha)\) generates \(K_d/k\), and let \(P_\beta\) be its monic minimal polynomial over \(k\).
Since \(P_\beta(h(\alpha))=0\),
\[
f\mid P_\beta(h),\qquad n\le d\deg h.
\]
Thus \(\deg h\ge m\). If \(\deg h=m\), divide \(h\) by its leading coefficient to make it monic;
this does not change the subfield generated by its value. Denote the monic minimal polynomial of the new value by \(P\).
Then \(f\) and \(P(h)\) are monic of the same degree, and the former divides the latter, so they are equal.

Conversely, if~\eqref{sfr:eq-optimal-composition} holds, then \(P\) must be irreducible,
since otherwise a nontrivial factorization, after substitution of \(h\), would make \(f\) reducible.
Thus \(h(\alpha)\) has degree \(d\), and uniqueness of the subfield implies that it generates \(K_d\).
This gives \(\kappa_f(d)\le m\), and the lower bound yields equality.
\end{proof}

The equality theorem gives the first central criterion. Now put \(K=K_d=\F_{q^d}\) and \(m=n/d\).
The remaining tools are organized below by the questions they answer; notation not yet introduced is defined in the cited sections.
\begingroup\small
\linespread{1.08}\selectfont
\setlength{\tabcolsep}{4pt}
\renewcommand{\arraystretch}{1.04}
\begin{longtable}{@{}>{\raggedright\arraybackslash}p{0.23\textwidth}>{\raggedright\arraybackslash}p{0.53\textwidth}>{\raggedright\arraybackslash}p{0.16\textwidth}@{}}
\caption{From defining equations and root blocks to subfield representations}\label{v100:tab-subfield-route}\\
\toprule
Data and question&Structure and information recovered&Main reference\\
\midrule
\endfirsthead
\toprule
Data and question&Structure and information recovered&Main reference\\
\midrule
\endhead
\bottomrule
\endfoot
Given \(f\), which low-degree values lie in \(K\)?&The Frobenius fixed-point matrix \(M_{f,d,t}\) gives \(\rho_{f,d}(t)=t+1-\rank M_{f,d,t}\); kernel vectors give a basis of the intersection.&\ref{sfr:prop-rank-recovery}\\
\addlinespace
Does a low-degree space contain a generator of \(K/k\)?&M\"obius inversion over divisor subfields gives exact generator counts and \(\kappa_f(d)\). A value outside \(k\) alone does not suffice.&\ref{sfr:prop-exact-count}\\
\addlinespace
Is there a polynomial inner factor of degree \(m\)?&The minimum representation degree equals \(m\) exactly when \(f=P\circ h\), with \(\deg P=d\), \(\deg h=m\); fixed-point equations then recover the inner factor.&\ref{sfr:thm-optimal-composition}, \ref{sfr:prop-rank-recovery}\\
\addlinespace
Given a relative root block, can rational degree \(m\) be attained?&Coefficient rank \(\eta_f(d)=2\) characterizes equality. Whether its coordinate-polynomial space contains \(1\) determines whether the denominator can be removed at the same degree.&\ref{str:thm-rational-rank}, \ref{str:cor-polynomial-direction}\\
\addlinespace
A finite prefix of weighted relative traces&Hankel equations recover the relative root block from \(2m\) traces. Retaining their \(K/k\) coordinates then recovers the degree profiles by trace duality.&\ref{str:thm-hankel-recovery}, \ref{str:thm-profile-trace-duality}\\
\end{longtable}
\endgroup

Irreducibility here already implies separability, so equality entails \(P'\ne0\) and \(h'\ne0\).
The theorem requires neither \(p\nmid n\) nor \(p\nmid m\).
A general low-degree generator representation \(h\) guarantees only \(f\mid P_\beta(h)\);
precisely when \(d\deg h=n\), equal degrees turn this divisibility into a compositional identity.

To interpret this equality on root sets, write the roots of \(f\) as \(\alpha,\alpha^q,\ldots,\alpha^{q^{n-1}}\), and put
\begin{equation}\label{sfr:eq-root-blocks}
S_j=\{\alpha^{q^{j+\ell d}}:0\le\ell<m\},
\qquad 0\le j<d.
\end{equation}
By Theorem~\ref{thm:factor-base-extension}, these are the root blocks of the irreducible factors of \(f\) over \(K_d\), equivalently its \(\sigma^d\)-orbits; \(\sigma\) maps \(S_j\) successively to \(S_{j+1}\), with indices modulo \(d\).
If \(k(h(\alpha))=K_d\), then
\[
h(\alpha^{q^i})=h(\alpha)^{q^i},
\]
so two roots have equal images if and only if their indices are congruent modulo \(d\).
Thus~\eqref{sfr:eq-root-blocks} is also the fiber partition of \(h\) on the original root set.

Each \(S_j\) contains \(m\) distinct roots. If a nonconstant polynomial \(h\) has constant value \(b\) on this block, then \(h(X)-b\) has at least \(m\) distinct roots, so \(\deg h\ge m\).
When \(\deg h=m\), the block \(S_j\) is exactly the full root set of the corresponding fiber polynomial, explaining why a least-degree representation yields a compositional identity.
If \(\deg h>m\), the block \(S_j\) is still contained in the fiber but need not be its full root set.

\begin{proposition}[A coefficient criterion for root-block polynomials]\label{sfr:prop-block-coefficients}
Put
\[
B_j(X)=\prod_{\gamma\in S_j}(X-\gamma).
\]
Then \(B_0\in K_d[X]\) and \(B_j=\sigma^j(B_0)\), where \(\sigma\) acts on coefficients.
The following conditions are equivalent:
\begin{enumerate}
\item \(\kappa_f(d)=m\);
\item all nonconstant coefficients of \(B_0\) belong to \(k\);
\item the nonconstant parts of all \(B_j\) agree, that is, every \(B_j-B_0\) is constant.
\end{enumerate}
When these conditions hold, one may take
\[
h(X)=B_0(X)-B_0(0),\qquad \beta=-B_0(0),
\]
and then \(h\in k[X]\), \(h(\alpha)=\beta\), and \([k(\beta):k]=d\).
\end{proposition}
\begin{proof}
The map \(\sigma^d\) preserves \(S_0\), hence fixes the coefficients of \(B_0\); the transformation of root sets directly gives \(B_j=\sigma^j(B_0)\).
If (1) holds, Theorem~\ref{sfr:thm-optimal-composition} gives \(f=P(h)\).
The polynomial \(h(X)-h(\alpha)^{q^j}\) is monic of degree \(m\) and vanishes on \(S_j\), so equals \(B_j\).
Thus the \(B_j\) differ only by constants, and their nonconstant coefficients lie in \(k\).

Condition (2) plainly implies (3). Conversely, \(B_1=\sigma(B_0)\) and (3) imply that
each nonconstant coefficient of \(B_0\) is fixed by \(\sigma\), hence lies in \(k\).
Under (2), the stated \(h\) belongs to \(k[X]\) and satisfies \(B_j=h-\beta^{q^j}\).
The \(B_j\) have distinct, pairwise disjoint root sets, so \(\beta,\beta^q,\ldots,\beta^{q^{d-1}}\) are distinct.
Since \(\beta\in K_d\), its degree is exactly \(d\), giving \(\kappa_f(d)=m\).
\end{proof}

\par
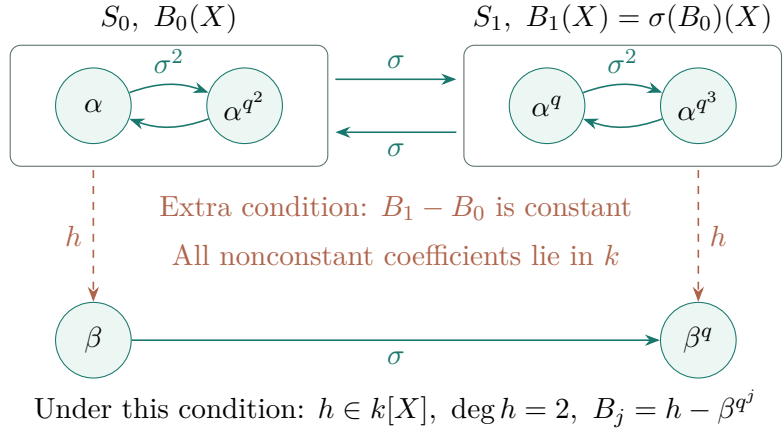
\begin{figure}[!htbp]
\centering
\begin{tikzpicture}[every node/.style={font=\small},
 pt/.style={circle,draw=Jade,fill=CoverMint,minimum size=10mm,inner sep=2pt},
 ar/.style={-{Stealth[length=2mm]},semithick,Jade}]
\draw[draw=MutedSlate,rounded corners] (-5.1,-.8) rectangle (-.9,.8);
\draw[draw=MutedSlate,rounded corners] (.9,-.8) rectangle (5.1,.8);
\node[pt] (a0) at (-4,0) {\(\alpha\)};
\node[pt] (a2) at (-2,0) {\(\alpha^{q^2}\)};
\node[pt] (a1) at (2,0) {\(\alpha^q\)};
\node[pt] (a3) at (4,0) {\(\alpha^{q^3}\)};
\draw[ar] (a0) to[bend left=20] node[above]{\(\sigma^2\)} (a2);
\draw[ar] (a2) to[bend left=20] (a0);
\draw[ar] (a1) to[bend left=20] node[above]{\(\sigma^2\)} (a3);
\draw[ar] (a3) to[bend left=20] (a1);
\node at (-3,1.15) {\(S_0,\ B_0(X)\)};
\node at (3,1.15) {\(S_1,\ B_1(X)=\sigma(B_0)(X)\)};
\draw[ar] (-.8,.35) -- node[above]{\(\sigma\)} (.8,.35);
\draw[ar] (.8,-.35) -- node[below]{\(\sigma\)} (-.8,-.35);
\node[align=center,text=Clay] at (0,-1.65) {Extra condition: \(B_1-B_0\) is constant\\All nonconstant coefficients lie in \(k\)};
\node[pt] (b0) at (-4,-3.1) {\(\beta\)};
\node[pt] (b1) at (4,-3.1) {\(\beta^q\)};
\draw[ar,Clay,dashed] (-4,-.85) -- node[left]{\(h\)} (b0);
\draw[ar,Clay,dashed] (4,-.85) -- node[right]{\(h\)} (b1);
\draw[ar] (b0) -- node[below]{\(\sigma\)} (b1);
\node at (0,-4) {Under this condition: \(h\in k[X],\ \deg h=2,\ B_j=h-\beta^{q^j}\)};
\end{tikzpicture}
\caption{Stable blocks always exist. The quadratic complete fibers shown by dashed arrows additionally require the two monic block polynomials to differ only by a constant.}\label{v103:fig-root-blocks}
\end{figure}

\Needspace{7\baselineskip}
\noindent\phantomsection\label{v103:par-quartic-blocks}\textbf{Root blocks and traces in the quartic example.}\quad
Take \(n=4,d=m=2\). Figure~\ref{v103:fig-root-blocks} first shows the two \(\sigma^2\)-orbits that always exist,
then identifies the additional condition needed to realize them as complete fibers of a single quadratic polynomial.

\FloatBarrier

For \(f_1\) in Example~\ref{sfr:ex-quartic-comparison}, the block polynomials are
\[
B_0=X^2+X+\beta,\qquad B_1=X^2+X+\beta^2.
\]
They differ only in their constant terms, giving the common inner polynomial \(h=X^2+X\).
The relative trace reads \(\Tr_{E/K}(\alpha)=\alpha+\alpha^4=1\).
For \(f_2\), taking \(\gamma,\gamma^4\) as the first block gives
\[
B_0=X^2+(\eta+1)X+(\eta+1),\qquad
B_1=X^2+\eta X+\eta.
\]
Their linear coefficients differ, so there is no common quadratic inner polynomial.
Here \(\Tr_{E/K}(\gamma)=\gamma+\gamma^4=\eta+1\notin\F_2\) detects exactly this obstruction.
The cubic representative \(h_2=X^3+X\) from the opening example still maps each root block to a single value, but each complete fiber has an additional root:
\[
\begin{aligned}
h_2(X)+\eta&=B_0(X)(X+\eta+1),\\
h_2(X)+\eta^2&=B_1(X)(X+\eta).
\end{aligned}
\]
The extra roots \(\eta+1,\eta\) lie in \(K\), outside the four roots of \(f_2\). Consequently,
\[
h_2(X)^2+h_2(X)+1=f_2(X)(X^2+X+1).
\]
This is a proper divisibility relation. A fiber partition restricted to the root set must not be mistaken for the identity \(f_2=P\circ h_2\).
For larger \(m\), the first relative trace detects only the coefficient of \(X^{m-1}\).
The other nonconstant coefficients must also be controlled: a trace in the base field alone does not imply a polynomial inner factor.

Every \(d\mid n\) gives these conjugate root blocks.
Proposition~\ref{prop:factor-coefficient-field} further shows that all coefficients of \(B_0\) together generate \(K_d\) over \(k\);
when its nonconstant coefficients lie in \(k\), this generating property is concentrated in the constant term.
Expressing these blocks as complete fibers of a single base-field polynomial additionally requires their nonconstant coefficients to agree.
Thus existence of the subfield and existence of the corresponding polynomial decomposition are different conditions.

\subsection{Dimension Bounds, Degree Jumps, and Representation Defect}
\label{sfr:subsec-profile-bound}

\begin{theorem}[A dimension bound for the degree profile]\label{sfr:thm-profile-bound}
For any \(d\mid n\), put \(m=n/d\). Then
\begin{equation}\label{sfr:eq-profile-bound}
\boxed{\quad \rho_{f,d}(t)\le 1+\left\lfloor\frac tm\right\rfloor
\qquad(0\le t<n).\quad}
\end{equation}
\end{theorem}
\begin{proof}
Put \(U=U_{f,d}(t)\) and \(r=\dim_k U\).
Since \(E=K_d(\alpha)\) and \([E:K_d]=m\), the elements \(1,\alpha,\ldots,\alpha^{m-1}\) are linearly independent over \(K_d\).
As \(U\subseteq K_d\), the following sum is a direct sum of \(k\)-spaces:
\[
U\oplus\alpha U\oplus\cdots\oplus\alpha^{m-1}U.
\]
Its dimension is \(mr\). If \(t\le n-m\), each summand has representatives of degree at most \(t+m-1<n\), so this direct sum lies in
\(\mathcal P_{t+m-1}(\alpha)\), of dimension \(t+m\). Hence
\[
mr\le t+m,\qquad r\le 1+\left\lfloor\frac tm\right\rfloor.
\]
If \(n-m\le t<n\), the right-hand side is \(d\), and the assertion follows directly from \(r\le\dim_k K_d=d\).
\end{proof}

\begin{definition}\label{sfr:def-jumps-defect}
Since \(\rho_{f,d}(-1)=0\), \(\rho_{f,d}(n-1)=d\), and each step increases by at most one, there is a unique sequence
\[
0=\lambda_0<\lambda_1<\cdots<\lambda_{d-1}<n,
\qquad
\lambda_i=\min\{t:\rho_{f,d}(t)\ge i+1\}.
\]
Call \((\lambda_0,\ldots,\lambda_{d-1})\) the \textbf{degree-jump sequence},
and write \(\Lambda_{f,d}=\{\lambda_0,\ldots,\lambda_{d-1}\}\).
\leavevmode\phantomsection\label{sfr:def-representation-defect}\ignorespaces Define the \textbf{representation defect} by
\begin{equation}\label{sfr:eq-defect}
\Delta_f(d)=\sum_{i=0}^{d-1}(\lambda_i-im),\qquad m=n/d.
\end{equation}
\end{definition}

By~\eqref{sfr:eq-profile-bound}, we have \(\lambda_i\ge im\), so \(\Delta_f(d)\) is a nonnegative integer.
The set \(\Lambda_{f,d}\) consists exactly of all possible representation degrees of nonzero elements of \(K_d\):
if \(u\in K_d^\times\) and \(\deg_\alpha(u)=a\), then \(u\in U_{f,d}(a)\setminus U_{f,d}(a-1)\), so \(a\) is a jump position; conversely, at each jump position \(a\), an element of \(U_{f,d}(a)\setminus U_{f,d}(a-1)\) has representation degree \(a\).
Choosing one new element at each jump position gives a basis of \(K_d\) whose representative degrees are the successive \(\lambda_i\).

\Needspace{10\baselineskip}
\begin{theorem}[Equivalent characterizations of zero defect]\label{sfr:thm-zero-defect}
Suppose \(1<d<n\) and \(m=n/d\). The following conditions are equivalent:
\begin{enumerate}
\item \(\kappa_f(d)=m\);
\item \(\rho_{f,d}(m)=2\);
\item \(\rho_{f,d}(t)=1+\lfloor t/m\rfloor\) for every \(0\le t<n\);
\item \(\lambda_i=im\) for \(0\le i<d\), or equivalently \(\Delta_f(d)=0\).
\end{enumerate}
They are also equivalent to the decomposition in Theorem~\ref{sfr:thm-optimal-composition}.
\end{theorem}
\begin{proof}
If (1) holds, take \(h\) of degree \(m\) such that \(\beta=h(\alpha)\) generates \(K_d/k\).
Then \(1,\beta,\ldots,\beta^{d-1}\) is a basis of \(K_d\), and its polynomial representatives
\(1,h,\ldots,h^{d-1}\) have successive degrees \(0,m,\ldots,(d-1)m<n\).
Thus \(\rho_{f,d}(t)\ge1+\lfloor t/m\rfloor\); the upper bound gives (3).
Conditions (3) and (4) are equivalent by the definition of the jump sequence; a sum of nonnegative terms is zero if and only if each term is zero.
Both imply (2).

Finally, suppose (2) holds. Take \(u=h(\alpha)\in U_{f,d}(m)\setminus k\), with \(\deg h\le m\).
Put \(c=[k(u):k]\), so \(1<c\mid d\). Denoting the minimal polynomial of \(u\) by \(P_u\), we have
\[
f\mid P_u(h),\qquad n\le c\deg h\le cm\le dm=n.
\]
All inequalities are therefore equalities, giving \(c=d\) and \(\deg h=m\), which is (1).
\end{proof}

This proof shows why belonging to a subfield and generating it must generally be distinguished:
in this argument, attainment of the lower bound \(m\) is what allows degree comparison to force every such element outside the base field to generate \(K_d/k\).
The defect \(\Delta_f(d)\) measures the total displacement of the entire jump sequence, not the single difference \(\kappa_f(d)-m\).
It can also be written as the sum of dimension deficits at all degree thresholds:
\[
\Delta_f(d)=\sum_{t=0}^{n-1}
\left(1+\left\lfloor\frac{t}{m}\right\rfloor-\rho_{f,d}(t)\right).
\]
Indeed, the upper-bound profile and the actual profile count the indices \(i\) satisfying \(im\le t\) and \(\lambda_i\le t\), respectively;
for each fixed \(i\), they differ at exactly the \(\lambda_i-im\) thresholds with \(im\le t<\lambda_i\).
Thus the representation defect records the subfield directions still absent from all low-degree spaces, not only the first occurrence of a generator.

\begin{corollary}[Uniqueness of the normalized inner component]\leavevmode\phantomsection\label{sfr:cor-normalized-inner}
If \(\Delta_f(d)=0\), there is a unique monic \(h\in k[X]\) satisfying
\[
\deg h=m,\qquad h(0)=0,\qquad f=P(h)\text{ for some monic }P\in k[X].
\]
The corresponding \(P\) is also unique. All polynomials of degree \(m\) whose values generate \(K_d\) are exactly \(ah+b\), with \(a\in k^\times,b\in k\).
\end{corollary}
\begin{proof}
Since \(\rho_{f,d}(m)=2\), choosing one \(h\) with the stated monicity and zero constant term gives
\[
U_{f,d}(m)=k\cdot1\oplus k\cdot h(\alpha).
\]
Representatives of degree less than \(n\) are unique, so every generator representation of degree \(m\) is \(ah+b\).
Monicity forces \(a=1\), and the zero constant term then forces \(b=0\).
The outer polynomial \(P\) is the monic minimal polynomial of \(h(\alpha)\), so is unique as well.
\end{proof}

Thus, for a fixed monic irreducible polynomial, each specified inner degree gives at most one such normalized inner component.
This uniqueness assertion applies to each fixed inner degree separately.

\subsection{Recovering Profiles from Coefficients and Counting Exact Degrees}
\label{sfr:subsec-coefficient-recovery}

\leavevmode\phantomsection\label{sfr:def-fixed-point-matrix}\ignorespaces Write \(k[X]_{\le t}\) for the space of polynomials of degree at most \(t\), including zero.
For \(d\mid n\), define the \(k\)-linear map
\[
T_{f,d,t}:k[X]_{\le t}\longrightarrow k[X]/(f),
\qquad h\longmapsto h^{q^d}-h\pmod f.
\]
Linearity follows because the \(q^d\)th-power map is additive and fixes \(k\) pointwise.
Represent the target in the basis \(1,X,\ldots,X^{n-1}\), and denote its matrix by \(M_{f,d,t}\). Column \(i\),
for \(0\le i\le t\), is the coefficient column of
\[
X^{iq^d}-X^i\pmod f.
\]
In particular, column \(i=0\) is zero.

\begin{proposition}[Recovery by linear equations]\leavevmode\phantomsection\label{sfr:prop-rank-recovery}
For \(0\le t<n\), we have
\begin{equation}\label{sfr:eq-rank-recovery}
\boxed{\quad\rho_{f,d}(t)=\dim_k\ker T_{f,d,t}
=t+1-\rank_k M_{f,d,t}.\quad}
\end{equation}
Evaluating a basis of the kernel at \(\alpha\) gives a basis of \(U_{f,d}(t)\).
Thus all degree profiles, jump sequences, and representation defects can be recovered from the coefficients of \(f\).
\end{proposition}
\begin{proof}
The evaluation map \(k[X]_{\le t}\to E\) is injective when \(t<n\), and
\[
h\in\ker T_{f,d,t}
\quad\Longleftrightarrow\quad
h(\alpha)^{q^d}=h(\alpha)
\quad\Longleftrightarrow\quad
h(\alpha)\in K_d.
\]
The kernel is therefore linearly isomorphic to \(U_{f,d}(t)\); rank--nullity gives the displayed equalities.
\end{proof}

To compute these columns, first find \(Z_d=X^{q^d}\bmod f\), then compute the successive remainders of \(1,Z_d,Z_d^2,\ldots\) modulo \(f\).
There is no need to expand polynomials of degree \(iq^d\) or first find individual roots in an extension.
To test only a specified decomposition degree, take \(t=m\) and check whether the kernel has dimension two.
If it does, making a nonconstant kernel vector monic and subtracting its constant term recovers \(h\) from Corollary~\ref{sfr:cor-normalized-inner}.
Then find the minimal polynomial of \(h(X)\bmod f\), or compare coefficients successively in powers of \(h\), to recover \(P\).

\begin{proposition}[Counting generators of an exact subfield]\label{sfr:prop-exact-count}
For \(d\mid n\) and \(0\le t<n\), put
\[
N_f(d,t)=\#\{h\in k[X]_{\le t}:[k(h(\alpha)):k]=d\}.
\]
Then
\begin{equation}\label{sfr:eq-mobius-count}
N_f(d,t)=\sum_{c\mid d}\mu(d/c)\,q^{\rho_{f,c}(t)},
\end{equation}
where \(\mu\) is the integer M\"obius function. For \(d>1\),
\[
\kappa_f(d)=\min\{t:N_f(d,t)>0\}.
\]
\end{proposition}
\begin{proof}
The space \(U_{f,d}(t)\) has exactly \(q^{\rho_{f,d}(t)}\) elements.
The degree over \(k\) of each element divides \(d\); conversely, every element of degree \(c\mid d\) belongs to \(K_d\).
Since representatives are unique when \(t<n\),
\[
q^{\rho_{f,d}(t)}=\sum_{c\mid d}N_f(c,t).
\]
M\"obius inversion gives~\eqref{sfr:eq-mobius-count}, and the final equality follows from the definition.
\end{proof}

Here \(\rho_{f,d}(t)>1\) cannot simply replace \(N_f(d,t)>0\).
When \(d\) is prime, they are equivalent because \(N_f(d,t)=q^{\rho_{f,d}(t)}-q\).

\subsection{Multiplicative Constraints and the Least Positive Defect}
\label{sfr:subsec-small-defect}

\begin{proposition}[Truncated additive closure of the degree set]\label{sfr:prop-truncated-closure}
If \(a,b\in\Lambda_{f,d}\) and \(a+b<n\), then \(a+b\in\Lambda_{f,d}\).
\end{proposition}
\begin{proof}
Take \(u,v\in K_d^\times\) with \(\deg_\alpha(u)=a\) and \(\deg_\alpha(v)=b\).
Since \(uv\in K_d^\times\), equation~\eqref{sfr:eq-degree-product} shows that \(uv\) has representation degree \(a+b\), which is therefore a jump position.
\end{proof}

This condition comes from multiplicative closure of the subfield, not merely from its structure as a \(k\)-linear subspace.
The truncation condition \(a+b<n\) is essential: once \(n\) is reached, reduction modulo \(f\) is needed.

\begin{proposition}[Parameter restrictions for defect one]\label{sfr:prop-defect-one}
Suppose \(1<d<n\) and \(d\mid n\). If \(\Delta_f(d)=1\), then
\[
d=2\qquad\text{or}\qquad(n,d)=(6,3).
\]
This is a necessary condition for a given polynomial to have defect one, not a necessary and sufficient criterion determined by the parameters alone.
\end{proposition}
\begin{proof}
Put \(m=n/d\ge2\). If \(\lambda_1=m\), Theorem~\ref{sfr:thm-zero-defect} gives \(\Delta_f(d)=0\), a contradiction.
Since each \(\lambda_i-im\) is a nonnegative integer, a total of one forces
\[
\Lambda_{f,d}=\{0,m+1,2m,3m,\ldots,(d-1)m\}.
\]
When \(d=2\), we have the first listed possibility. Now suppose \(d\ge3\).
If \(m\ge3\), then \(2(m+1)<dm=n\), but \(2m+2\) is neither \(m+1\) nor a multiple of \(m\),
so it does not belong to \(\Lambda_{f,d}\), contrary to truncated additive closure.
Hence \(m=2\). If also \(d\ge4\), then \(3,4\in\Lambda_{f,d}\) and \(3+4=7<n\),
but \(7\notin\Lambda_{f,d}\), again a contradiction. Thus only \(m=2,d=3\) remains, that is, \((n,d)=(6,3)\).
\end{proof}

\subsection{Different Representation Degrees in the Same Finite Field}
\label{sfr:subsec-examples}

\par\Needspace{7\baselineskip}
\noindent\phantomsection\label{v103:par-quartic-profile}\textbf{Return to the opening quartic example.}\quad
The direct calculations in Example~\ref{sfr:ex-quartic-comparison} now give the degree profiles and defects:

\[
\begin{array}{c|ccc}
&\kappa_f(2)&\Lambda_{f,2}&\Delta_f(2)\\\hline
f_1=X^4+X+1&2&\{0,2\}&0\\
f_2=X^4+X^3+1&3&\{0,3\}&1
\end{array}
\qquad
\begin{array}{c|cccc}
t&0&1&2&3\\\hline
\rho_{f_1,2}(t)&1&1&2&2\\
\rho_{f_2,2}(t)&1&1&1&2
\end{array}
\]
For the second defining polynomial the subfield and the two-point root blocks still exist; the quadratic complete-fiber representation fails.
One can also check directly that, in characteristic two, \(P=Y^2+cY+e\) and \(h=X^2+aX+b\) give a composition \(P(h)\) with no cubic term, so it cannot equal \(f_2\).
Example~\ref{exm:compression-return-four} compares translation-quotient returns for the same \(f_1\); here the subfield is fixed and its representations in different power bases are compared.

\begin{example}[Defect one in degree six]\label{sfr:ex-six-defect-one}
Use \(f=X^6+X+1\in\F_2[X]\), proved irreducible in Example~\ref{exm:f64-computation}.
Let \(f(\alpha)=0\), and put
\[
u=\alpha^3+\alpha^2+\alpha,\qquad v=\alpha^4+\alpha^3.
\]
Direct computation modulo \(f\) gives
\[
u^2=u+v+1,\qquad v^2=u+1.
\]
Thus \(u^8=u\) and \(v^8=v\), while \(1,u,v\) have distinct representative degrees \(0,3,4\) and are therefore linearly independent.
They form a basis of \(K_3=\F_8\). In any linear combination, the leading term can occur only at the largest chosen degree, so
\[
\Lambda_{f,3}=\{0,3,4\},\qquad
\rho_{f,3}=(1,1,1,2,3,3),\qquad
\Delta_f(3)=1.
\]
Since \(3\) is prime, \(u\notin\F_2\) generates \(K_3\), giving \(\kappa_f(3)=3>6/3\).
This realizes the degree-six exception in Proposition~\ref{sfr:prop-defect-one}.
\end{example}

\par\Needspace{12\baselineskip}
The restriction on defect one leaves a realization question: after one has been excluded, does the first remaining positive integer occur? The following degree-eight example has defect two, certified entirely by quotient-ring identities.
\begin{example}[The least positive defect in degree eight]\label{v115:ex-defect-two}
Over \(\F_2\), take
\[
f=X^8+X^7+X^3+X^2+1,\qquad
r=X^5+X^4+X^3+X^2+X+1.
\]
Repeated squaring modulo \(f\) gives \(X^{16}+X\equiv r\) and \(X^{256}\equiv X\), while
\[
(X^3+X^2+1)f+(X^6+X^5+X^4)r=1.
\]
Thus \(\gcd(f,X^{16}-X)=1\). The Frobenius irreducibility criterion proves that \(f\) is irreducible: two is the only prime divisor of eight, so only the degree-four proper-subfield test is required.

Let \(f(\alpha)=0\), and put
\[
u=\alpha^3+\alpha,\qquad
v=\alpha^5+\alpha^4+\alpha^2,\qquad
w=\alpha^6+\alpha^2.
\]
Reduction gives
\[
u^2=w,\qquad w^2=1+u,\qquad v^2=u+v.
\]
Hence \(u^{16}=u\) and \(w^{16}=w\); also \(v^4=w+u+v\) and \(v^8=1+w+v\), so \(v^{16}=v\).
The representatives of \(1,u,v,w\) have degrees \(0,3,5,6\). They are therefore independent and form a basis of \(K_4=\F_{16}\).
Their distinct leading terms cannot cancel, giving
\[
\Lambda_{f,4}=\{0,3,5,6\},\qquad
\Delta_f(4)=(3-2)+(5-4)+(6-6)=2.
\]
Since \(u^4=1+u\ne u\), the element \(u\) generates \(K_4/\F_2\), and \(\kappa_f(4)=3\).
Proposition~\ref{sfr:prop-defect-one} excludes defect one for \((n,d)=(8,4)\). Thus two is the least positive defect for this base field and these parameters.
\end{example}

\subsection{Small Positive Defect and Completing Fibers}\label{v116:subsec-small-defect}
The restriction on defect one can be pushed further by the same argument. Write \(n=dm\). When the total defect is less than \(m\), each interval of length \(m\) contains just one jump, whose position is then constrained by multiplication.
\begin{theorem}[Arithmetic progressions of jumps at small defect]\label{v116:thm-small-defect}
Suppose \(d\ge4\). If \(0<\Delta_f(d)<m\), there is an integer \(a\ge1\) such that
\begin{equation}\label{v116:eq-small-defect}
\lambda_i=i(m+a)\quad(0\le i<d),\qquad
\Delta_f(d)=a\binom d2.
\end{equation}
Consequently, positive defect satisfies \(\Delta_f(d)\ge\min\{m,\binom d2\}\).
\end{theorem}
\begin{proof}
Write \(\lambda_i=im+r_i\), where \(r_i\ge0\) and \(\sum_i r_i<m\). Each interval \([im,(i+1)m)\) contains exactly one jump.
Put \(a=r_1\). If \(a=0\), Theorem~\ref{sfr:thm-zero-defect} gives zero defect; hence \(a>0\).
If \(2a\ge m\), the distinct numbers \(2\lambda_1\) and \(\lambda_1+\lambda_2\) both lie in \([3m,4m)\): use \(a<m\) for the first and \(a+r_2<m\) for the second.
Both are less than \(dm\), so truncated additive closure makes both jumps, a contradiction. Thus \(2a<m\), giving \(\lambda_2=2\lambda_1\).
For \(2\le i<d-1\), we have \(a+r_i<m\), so \(\lambda_1+\lambda_i\in[(i+1)m,(i+2)m)\).
Closure and uniqueness now give \(\lambda_{i+1}=\lambda_1+\lambda_i\) successively. Summing proves the formula. If \(\Delta_f(d)\ge m\), the last bound is immediate.
\end{proof}

The range cannot simply be extended to \(\Delta=m\). In Example~\ref{v115:ex-defect-two}, \(m=\Delta=2\), and the jumps \(0,3,5,6\) do not form an arithmetic progression.

\begin{corollary}[Small defect and the residual factor of a composition]\label{v116:cor-residual-factor}
Under the hypotheses of Theorem~\ref{v116:thm-small-defect}, there is a unique monic polynomial \(h\) with zero constant term such that \(\deg h=m+a\) and \(h(\alpha)\in K_d\setminus k\).
Its value generates \(K_d/k\). If \(P\) is the monic minimal polynomial of that value, then
\begin{equation}\label{v116:eq-residual-factor}
P(h(X))=f(X)R(X),\qquad \deg P=d,\quad\deg R=da.
\end{equation}
Every irreducible factor of \(R\) has degree divisible by \(d\). In particular, when \(a=1\), \(R\) is irreducible of degree \(d\).
\end{corollary}
\begin{proof}
The space \(U_{f,d}(m+a)\) is two-dimensional and contains the constants, giving the unique monic \(h\) with zero constant term.
If \(\beta=h(\alpha)\) generated only a subfield of degree \(c<d\), then \(c\mid d\) and \(c\le d/2\).
Since \(f\mid P_\beta(h)\) and \(m+a<2m\), we would have \(dm\le c(m+a)<2cm\le dm\), a contradiction.
Thus \(\deg P=d\), and divisibility gives~\eqref{v116:eq-residual-factor}.
For any root \(\gamma\) of \(R\), the identity \(P(h(\gamma))=0\) shows that \([k(h(\gamma)):k]=d\) divides \([k(\gamma):k]\).
When \(a=1\), the factor degrees, counted with multiplicity, sum to \(d\). There can therefore be just one irreducible factor of degree \(d\), with multiplicity one.
\end{proof}

When \(a=1\), each root block of \(f\), of size \(m\), needs exactly one more root to become a full fiber of \(h\); the added \(d\) roots form another Frobenius orbit.
We now construct a fiber whose return cycles have lengths \(1,m\).

\begin{proposition}[Realizing the first positive defect by linearized fibers]\label{v116:prop-defect-construction}
Let \(C(T)\) be a monic primitive irreducible polynomial of degree \(s\) over \(k=\F_q\), and put
\[
m=q^s-1,\qquad h=\mathcal L_C,\qquad d\ge4,\qquad
\gcd(d,m)=1,\quad m>\binom d2.
\]
For every monic irreducible \(P\in k[X]\) of degree \(d\), we have
\[
P(h)=Rf,
\]
where \(R,f\) are distinct monic irreducible polynomials of degrees \(d,dm\), respectively, and
\[
\Lambda_{f,d}=\{0,m+1,\ldots,(d-1)(m+1)\},\qquad
\Delta_f(d)=\binom d2.
\]
\end{proposition}
\begin{proof}
In a finite splitting extension put \(V=\ker \mathcal L_C\) and \(\sigma(x)=x^q\). Since \(C(0)\ne0\), the space \(V\) has \(q^s\) elements.
Irreducibility of \(C\) makes \(V\) isomorphic to \(k[T]/(C)\) as a \(k[T]\)-module. Primitivity means that \(\sigma\) is a single cycle of length \(m\) on \(V\setminus\{0\}\).
On \(K=\F_{q^d}\), the operator \(C(\sigma)\) is invertible: the roots of \(C\) have order \(m\), and \(m\nmid d\), so \(\gcd(C,T^d-1)=1\).
For a root \(\beta\in K\) of \(P\), there is a unique \(\gamma\in K\) with \(h(\gamma)=\beta\).
The fiber is \(\gamma+V\); its return \(\sigma^d\) fixes \(\gamma\) and is a single cycle on the other \(m\) points.
Theorem~\ref{thm:fiber-factorization} gives the two factors of degrees \(d,dm\). Both \(P\) and \(h\) are separable, so there are no repeated factors.

For a root \(\alpha\) of \(f\), the element \(h(\alpha)\) generates \(K\). Since \(m>\binom d2\ge d-1\), the polynomials
\(1,h,\ldots,h^{d-1}\) have distinct degrees below \(dm\).
Their values form a basis of \(K/k\); distinct leading terms cannot cancel in a linear combination. These degrees are therefore exactly the jumps, and the defect follows by substitution.
\end{proof}

\begin{example}[Defect six in a degree-twenty-eight representation]\label{v116:ex-defect-six}
Over \(\F_2\), take \(C=T^3+T+1\), \(h=X^8+X^2+X\), and \(P=X^4+X+1\).
The cubic \(C\) is irreducible, and its roots are nonidentity elements of \(\F_8^\times\), hence have order seven. Irreducibility of \(P\) was checked in Example~\ref{sfr:ex-quartic-comparison}.
Here \(m=7,d=4\). If \(P(\gamma)=0\), then \(\gamma^8=\gamma^2+1\), so \(h(\gamma)=\gamma+1\) is again a root of \(P\).
Thus \(R=P\) in the proposition, and
\[
f(X)=\frac{X^{32}+X^4+X^2+X+1}{X^4+X+1}
\]
is irreducible, with \(\Lambda_{f,4}=\{0,8,16,24\}\) and \(\Delta_f(4)=6\). The small-defect bound shows that six is the least positive defect for these parameters.
\end{example}

Arithmetic progressions give the necessary degree pattern, and the fiber construction realizes one layer. For fixed \(q,d,m\), which values \(a\binom d2<m\) occur?
The factor degrees of \(R\) in~\eqref{v116:eq-residual-factor}, divided by \(d\) and counted with multiplicity, form a partition of \(a\).
Can every such partition occur, and how many polynomials \(f\) realize each one? Here both the degree profile and the fiber return cycles must be retained.

\Needspace{9\baselineskip}
\subsection{Realization conditions for jump sequences}
\label{v137:subsec-jump-realization}
Degree inequalities and truncated additive closure retain only leading terms. To decide whether a candidate sequence actually comes from a subfield, the lower coefficients must satisfy multiplication relations and Frobenius equations simultaneously.

\begin{theorem}[Multiplication and fixed-point tests for jump sequences]
\label{v137:thm-jump-realization}
Let $k=\mathbb F_q$, let $f\in k[X]$ be monic irreducible of degree $n$, let $d\mid n$, and put $E=k[X]/(f)$ and $\alpha=X\bmod f$. Let $K_d$ be the subfield of $E$ with $q^d$ elements. Given
\[
 \Lambda=\{0=\lambda_0<\lambda_1<\cdots<\lambda_{d-1}<n\},
\]
choose coefficients $b_{ia}\in k$ and set
\[
 h_i(X)=X^{\lambda_i}
 +\sum_{\substack{0\le a<\lambda_i\\a\notin\Lambda}}b_{ia}X^a
 \quad(0\le i<d).
\]
Write $\operatorname{rem}_f P$ for the remainder of $P$ on division by $f$, and $[X^a]P$ for its coefficient of $X^a$. Define
\begin{equation}\label{v137:eq-jump-projection}
 \Pi_\Lambda(P)=\operatorname{rem}_f P
 -\sum_{i=0}^{d-1}\bigl([X^{\lambda_i}]\operatorname{rem}_f P\bigr)h_i.
\end{equation}
The following conditions are equivalent:
\begin{enumerate}[label=\textup{(\roman*)}]
\item $\{h_i(\alpha):0\le i<d\}$ is a $k$-basis of $K_d$; in this case $\Lambda_{f,d}=\Lambda$.
\item $\Pi_\Lambda(h_ih_j)=0$ for all $0\le i,j<d$.
\item $\operatorname{rem}_f(h_i^{q^d}-h_i)=0$ for all $0\le i<d$.
\end{enumerate}
Thus $\Lambda$ is realizable if and only if some choice of the $b_{ia}$ satisfies \textup{(ii)} or \textup{(iii)}; when such a choice exists, it is unique.
\end{theorem}
\begin{proof}
The coefficient matrix of the $h_i$ at the jump positions is the identity, so they are linearly independent, and the kernel of~\eqref{v137:eq-jump-projection} in $k[X]/(f)$ is their span $W$. The possible degrees of its nonzero elements are exactly $\Lambda$.
Since $h_0=1$, condition \textup{(ii)} says that $W$ is a subalgebra containing the identity. A finite subalgebra of a field is a field; since $\dim_k W=d$, it must be $K_d$.
The fixed points of $x\mapsto x^{q^d}$ in $E$ are exactly $K_d$, so \textup{(iii)} says that $W\subseteq K_d$; equality follows from the dimensions.
Finally, elimination by leading terms on the polynomial representatives of $K_d$, with zero coefficients at the other jump positions in each basis vector, gives the stated unique form.
\end{proof}

For fixed $f$, condition \textup{(iii)} is linear in the $b_{ia}$, because the $q^d$th power fixes $k$; it connects with the kernel formula in Proposition~\ref{sfr:prop-rank-recovery}. Condition \textup{(ii)} exhibits the multiplication compatibility beyond degree constraints.
When these conditions hold, the matrix of the $q$th-power Frobenius in this basis is
\[
 S_{ij}=[X^{\lambda_i}]\operatorname{rem}_f(h_j^q).
\]
By the normal basis theorem, its minimal polynomial is $T^d-1$, and $\dim_k\ker(S^a-I)=\gcd(d,a)$ for every $a\ge1$. These restrictions also hold when $p\mid d$.

\begin{example}[All jump sequences for binary quartic representations]
\label{v137:ex-quartic-classification}
There are exactly $(2^4-2^2)/4=3$ monic irreducible quartics over $\mathbb F_2$. Their data for the quadratic subfield are
\[
\begin{array}{c|c|c|c}
 f&h&\Lambda_{f,2}&\Delta_f(2)\\\hline
 X^4+X+1&X^2+X&\{0,2\}&0\\
 X^4+X^3+1&X^3+X&\{0,3\}&1\\
 X^4+X^3+X^2+X+1&X^3+X^2&\{0,3\}&1
\end{array}
\]
The first two rows are Example~\ref{sfr:ex-quartic-comparison}. The third polynomial is $\Phi_5$, irreducible because $\operatorname{ord}_5(2)=4$. If $\alpha$ is a root, then $\alpha^5=1$, and
\[
 (\alpha^3+\alpha^2)^2+(\alpha^3+\alpha^2)+1
 =1+\alpha+\alpha^2+\alpha^3+\alpha^4=0.
\]
Thus each row gives a basis $1,h(\alpha)$ of $\mathbb F_4$. This lists every possible sequence and counts its occurrences among irreducible defining polynomials.
\end{example}

For fixed $q,n,d$, the degree bounds and truncated multiplication first give candidate sets $\Lambda$; which satisfy the coefficient equations above? As $f$ varies, which jump sequences occur, and how many irreducible defining polynomials realize each one? Composition in the zero-defect case and fiber constructions for small positive defect provide two classes of realizations.

\subsection{Optimal Inner Components and Degree Transfer in Subfield Towers}
\label{sfr:subsec-towers}

\begin{proposition}[Tower transfer of optimal representations]\label{sfr:prop-tower-scaling}
Suppose \(d\mid n\) and \(1<d\le n\), and put \(m=n/d\). If \(f=P(h)\), where \(P,h\in k[X]\) are monic with
\(\deg P=d\) and \(\deg h=m\), then, putting \(\beta=h(\alpha)\), we have \(K_d=k(\beta)\).
For every \(c\mid d\) and \(0\le t<n\),
\begin{equation}\label{sfr:eq-tower-profile}
\rho_{f,c}(t)=\rho_{P,c}\!\left(\left\lfloor\frac tm\right\rfloor\right).
\end{equation}
Correspondingly,
\[
\lambda_{f,c,i}=m\lambda_{P,c,i}\quad(0\le i<c),
\qquad\Delta_f(c)=m\Delta_P(c).
\]
For \(c>1\), we also have
\begin{equation}\label{sfr:eq-tower-kappa}
\kappa_f(c)=m\kappa_P(c).
\end{equation}
The data on the right are defined within \(K_d\), using the root \(\beta\) of \(P\) as generator.
\end{proposition}
\begin{proof}
Every \(u\in K_d\) is uniquely \(A(\beta)\), where \(A\in k[Y]\) and \(\deg A<d\).
If \(u\ne0\), then \(A(h(X))\) has degree \(m\deg A<n\), so is its unique representative in the \(\alpha\)-coordinates.
Therefore
\[
\deg_\alpha(u)=m\deg_\beta(u)\qquad(u\in K_d^\times).
\]
Restricting to \(K_c\) gives a bijection of the filtration spaces, yielding~\eqref{sfr:eq-tower-profile} and the jump-sequence formula.
The reference spacing in the definition of \(\Delta_f(c)\) is \(n/c=m(d/c)\); subtracting term by term and summing gives the defect formula.
Finally, minimizing over generators of \(K_c/k\) gives~\eqref{sfr:eq-tower-kappa}.
\end{proof}

In particular, once \(f\) has an optimal inner component relative to \(K_d\), further representation in \(K_c\) becomes exactly the same problem for the outer polynomial \(P\).
We have \(\Delta_f(c)=0\) if and only if \(\Delta_P(c)=0\), so when both stages attain the least degree, the normalized inner components can be composed successively.

The earlier decomposition \(f=F(L_{V_f})\) of an irreducible polynomial, obtained from its translation stabilizer, is a special case:
if \(1<\deg F<n\), take \(d=\deg F\); then \(\kappa_f(d)=|V_f|\) and \(\Delta_f(d)=0\).
More explicitly, after Theorem~\ref{fas:thm-full-recovery} recovers \(f=H((L_V-\delta)^t)\) from root values,
if \(d=\deg H\) satisfies \(1<d<n\), then \(\kappa_f(d)=|V_f|t\) and \(\Delta_f(d)=0\).
If the outer degree is one, the quotient value already lies in the base field, outside the proper intermediate-subfield case considered here; in particular, the group order in the special-orbit case cannot simply be taken as the inner degree.
The polynomial \(h\) in this section, however, need not arise from any affine stabilizer: it is characterized by the least representation degree of the subfield.
Thus conjugate root blocks provide all candidate partitions, and representation degree determines when they become complete polynomial fibers;
stabilizers explain why some of these fibers have additive or affine symmetry.

\section{Duality Between Subfield Degree Profiles and Relative Traces}
\label{str:sec-trace-duality}

Low-degree representations of a subfield are determined by its intersections with the power-basis degree spaces; trace pairing converts these intersections into orthogonality conditions.
Relative trace expresses these conditions over the subfield, recovering the entire degree profile from the ranks of spans of successive trace values.
Continue to assume that \(f\in k[X]\) is monic and irreducible with \(\deg f=n\), \(E=k(\alpha)\), and \(f(\alpha)=0\), and fix
\[
1<d<n,\qquad d\mid n,\qquad K=K_d=\F_{q^d},\qquad m=n/d.
\]
The degree spaces \(\mathcal P_t(\alpha)\), profile \(\rho_{f,d}\), and defect \(\Delta_f(d)\) retain their definitions from the preceding section.

\subsection{Weighted-Trace Orthogonal Complements of Degree Spaces}

\leavevmode\phantomsection\label{str:def-weighted-pairing}\ignorespaces
Define the \textbf{weighted trace pairing} on \(E\) by
\[
\langle u,v\rangle_\alpha
=\Tr_{E/k}\left(\frac{uv}{f'(\alpha)}\right).
\]
An irreducible polynomial over a finite field is separable, so \(f'(\alpha)\ne0\); invertibility of nonzero scalar multiplication and nondegeneracy of the trace pairing show that this pairing is nondegenerate.
Write \(U^{\perp_\alpha}\) for the orthogonal complement with respect to it, to distinguish it from that for the ordinary trace pairing.

\begin{proposition}[Reverse orthogonality of the degree filtration]\label{str:prop-filter-orthogonal}
For \(-1\le t<n\), we have
\begin{equation}\label{str:eq-filter-orthogonal}
\mathcal P_t(\alpha)^{\perp_\alpha}
=\mathcal P_{n-2-t}(\alpha).
\end{equation}
Here \(\mathcal P_{-1}(\alpha)=\{0\}\).
\end{proposition}
\begin{proof}
If \(0\le t\le n-2\), \(u\in\mathcal P_t(\alpha)\), and \(v\in\mathcal P_{n-2-t}(\alpha)\),
then \(uv\) has a polynomial representative of degree at most \(n-2\). The leading-coefficient formula~\eqref{eq:coefficient-trace} gives \(\langle u,v\rangle_\alpha=0\).
Thus the right side is contained in the left. Nondegeneracy gives both sides dimension \(n-t-1\), so they are equal.
For \(t=-1\) and \(t=n-1\), the equalities are \(\{0\}^{\perp_\alpha}=E\) and \(E^{\perp_\alpha}=\{0\}\), respectively.
\end{proof}

The weight \(1/f'(\alpha)\) turns the trace functional into the leading-coefficient functional, making the orthogonal complement of a low-degree space another low-degree space.
For the ordinary trace pairing, a scalar factor \(1/f'(\alpha)\) must be inserted on the right; formula~\eqref{str:eq-filter-orthogonal} cannot be retained unchanged.

\subsection{Relative Trace Sequences and Recovery of the Profile}

\begin{definition}\label{str:def-relative-trace-sequence}
Define the \textbf{weighted relative trace sequence} associated with \((f,d)\) by
\begin{equation}\label{str:eq-relative-trace-sequence}
 c_j=\Tr_{E/K}\left(\frac{\alpha^j}{f'(\alpha)}\right)\in K
 \qquad(j\ge0).
\end{equation}
\leavevmode\phantomsection\label{str:def-trace-rank}\ignorespaces
Define the spans of its initial segments and its \textbf{trace-rank function} by
\[
\mathcal C_{f,d}(s)=\operatorname{span}_k\{c_0,\ldots,c_s\},\qquad
r_{f,d}(s)=\dim_k\mathcal C_{f,d}(s)\quad(s\ge0),
\]
and set \(\mathcal C_{f,d}(-1)=\{0\}\) and \(r_{f,d}(-1)=0\).
\end{definition}

Each \(c_j\) lies in the subfield \(K\), but the span and dimension here are taken over the base field \(k\).

\begin{theorem}[Duality of degree profiles and relative traces]\leavevmode\phantomsection\label{str:thm-profile-trace-duality}
For \(-1\le t<n\), we have
\begin{equation}\label{str:eq-profile-trace-duality}
\boxed{\quad \rho_{f,d}(t)+r_{f,d}(n-2-t)=d.\quad}
\end{equation}
More precisely, if \(\perp_K\) denotes orthogonal complement for the ordinary trace pairing of \(K/k\), then
\begin{equation}\label{str:eq-subfield-annihilator}
U_{f,d}(t)=\mathcal C_{f,d}(n-2-t)^{\perp_K}.
\end{equation}
\end{theorem}
\begin{proof}
For \(u\in K\), Proposition~\ref{str:prop-filter-orthogonal} and the double-orthogonal-complement identity for a nondegenerate pairing give
\[
 u\in\mathcal P_t(\alpha)
 \quad\Longleftrightarrow\quad
 \Tr_{E/k}\left(\frac{u\alpha^j}{f'(\alpha)}\right)=0
 \quad(0\le j\le n-2-t).
\]
By transitivity of trace and \(K\)-linearity of \(\Tr_{E/K}\), the right-hand condition is equivalent to
\[
 \Tr_{K/k}(u c_j)=0\qquad(0\le j\le n-2-t).
\]
These are precisely the conditions in~\eqref{str:eq-subfield-annihilator}. Taking dimensions in the nondegenerate trace pairing of \(K/k\) gives~\eqref{str:eq-profile-trace-duality}.
When the upper bound is \(-1\), there are no orthogonality conditions, and the right side is all of \(K\), so the argument still applies.
\end{proof}

Taking \(t=-1\) gives \(r_{f,d}(n-1)=d\), so the first \(n\) trace values already span all of \(K\).
Taking \(t=0\) gives \(r_{f,d}(n-2)=d-1\). These two boundary cases also follow from the leading-coefficient formula:
\[
\Tr_{K/k}(c_j)=0\quad(0\le j<n-1),\qquad
\Tr_{K/k}(c_{n-1})=1.
\]
Thus the first \(n-1\) trace values span the trace kernel of \(K/k\), and the last trace value completes it to all of \(K\).

\subsection{Two Kinds of Jumps and a Trace-Rank Criterion for Decomposition}

\leavevmode\phantomsection\label{str:def-trace-jumps}\ignorespaces
Since the trace rank increases by at most one per step and reaches \(d\) by \(n-1\), define the \textbf{trace-rank jump positions} by
\[
\mu_i=\min\{s\ge0:r_{f,d}(s)\ge i+1\},\qquad 0\le i<d.
\]
Here the \(\mu_i\) are integers listed in increasing order, not the integer M\"obius function \(\mu(\cdot)\).

\begin{corollary}[Reversed jumps and representation defect]\label{str:cor-reversed-jumps}
Let \(\lambda_0,\ldots,\lambda_{d-1}\) be the degree-jump sequence of \(K\). Then
\begin{equation}\label{str:eq-reversed-jumps}
\mu_i=n-1-\lambda_{d-1-i}\qquad(0\le i<d),
\end{equation}
and hence
\begin{equation}\label{str:eq-trace-defect}
\Delta_f(d)=\sum_{i=0}^{d-1}\bigl((i+1)m-1-\mu_i\bigr).
\end{equation}
Every summand is nonnegative.
\end{corollary}
\begin{proof}
For \(0\le t<n\), subtract~\eqref{str:eq-profile-trace-duality} at \(t-1\) from the identity at \(t\) to obtain
\[
\rho_{f,d}(t)-\rho_{f,d}(t-1)
=r_{f,d}(n-1-t)-r_{f,d}(n-2-t).
\]
Thus a degree jump at \(t\) corresponds exactly to a trace-rank jump at \(n-1-t\); ordering them gives~\eqref{str:eq-reversed-jumps}.
Since \(\lambda_{d-1-i}\ge(d-1-i)m\), we have \(\mu_i\le(i+1)m-1\). Summing and using \(n=dm\) gives the defect formula.
\end{proof}

\begin{corollary}[A trace-rank criterion for polynomial composition]\label{str:cor-trace-composition}
Let \(1<d<n\), \(d\mid n\), and \(m=n/d\), with \(f\) as above. The following conditions are equivalent:
\begin{enumerate}
\item there exist monic \(P,h\in k[X]\) with \(f=P(h)\), \(\deg P=d\), and \(\deg h=m\);
\item \(r_{f,d}(n-m-2)=d-2\);
\item \(r_{f,d}(s)=\lfloor(s+1)/m\rfloor\) for every \(-1\le s<n\);
\item \(\mu_i=(i+1)m-1\) for \(0\le i<d\).
\end{enumerate}
\end{corollary}
\begin{proof}
Taking \(t=m\) in the duality formula makes (2) equivalent to \(\rho_{f,d}(m)=2\).
Theorem~\ref{sfr:thm-zero-defect} equates this with (1) and with the full least-degree profile.
Substitution of that profile into~\eqref{str:eq-profile-trace-duality} gives (3); conditions (3) and (4) are equivalent by the definition of jump positions.
\end{proof}

\begin{example}[Reversed jumps in the quartic examples]\label{str:ex-trace-jumps-four}
In Example~\ref{sfr:ex-quartic-comparison}, the degree jumps for the four-element subfield are \((0,2)\) for \(f_1=X^4+X+1\) and \((0,3)\) for \(f_2=X^4+X^3+1\).
The corresponding trace-rank jumps are therefore \((1,3)\) and \((0,3)\), and the trace-rank functions at \(s=0,1,2,3\) are, respectively,
\[
(0,1,1,2),\qquad(1,1,1,2).
\]
The first polynomial has a quadratic inner component, whereas the second does not; on the trace side, the distinction is whether the first nonzero trace value already occurs at \(s=0\).
\end{example}

\section{Root-Block Coefficient Rank and Rational Representations of Subfields}
\label{str:sec-rational-representation}

A polynomial representation attaining the degree lower bound requires the conjugate root blocks to differ only in their constant terms.
If a rational generator representation attains the same bound, the corresponding root-block polynomial is a linear combination of two base-field polynomials.
This section characterizes that condition by root-block coefficient rank and compares polynomial and rational representations.
Retain \(k,E,f,\alpha,K,d,m\), always assuming \(1<d<n\).

\subsection{Two Coefficient Spaces of a Relative Root Block}

\leavevmode\phantomsection\label{str:def-block-polynomial}\ignorespaces
Denote the monic minimal polynomial of \(\alpha\) over \(K\) by
\begin{equation}\label{str:eq-block-polynomial}
 B_{f,d}(X)=\prod_{\ell=0}^{m-1}(X-\alpha^{q^{d\ell}})
 =X^m+b_{m-1}X^{m-1}+\cdots+b_0.
\end{equation}
This is \(B_0\) from Section~\ref{sfr:sec-subfield-representation}. Proposition~\ref{prop:factor-coefficient-field} gives
\begin{equation}\label{str:eq-block-coefficient-field}
k(b_0,\ldots,b_{m-1})=K.
\end{equation}
The field generated by these coefficients must be distinguished from their linear span, defined below.

\begin{definition}\label{str:def-coefficient-rank}
Put
\[
\mathcal V_{f,d}=\operatorname{span}_k\{1,b_0,\ldots,b_{m-1}\}\subseteq K,
\qquad \eta_f(d)=\dim_k\mathcal V_{f,d}.
\]
Call \(\eta_f(d)\) the \textbf{root-block coefficient rank}.
\leavevmode\phantomsection\label{str:def-coordinate-space}\ignorespaces
Choose a basis \(e_1,\ldots,e_d\) of \(K/k\), and write uniquely
\[
B_{f,d}(X)=\sum_{\ell=1}^d e_\ell B_\ell(X),\qquad B_\ell\in k[X],
\]
and define the \textbf{root-block coordinate-polynomial space}
\[
\mathcal W_{f,d}=\operatorname{span}_k\{B_1,\ldots,B_d\}\subseteq k[X]_{\le m}.
\]
\end{definition}

\begin{proposition}\label{str:prop-coefficient-spaces}
The space \(\mathcal W_{f,d}\) is independent of the chosen basis of \(K/k\), and
\[
\dim_k\mathcal W_{f,d}=\eta_f(d),\qquad
2\le\eta_f(d)\le\min\{d,m+1\}.
\]
Replacing the root of \(f\) by any conjugate leaves both \(\eta_f(d)\) and \(\mathcal W_{f,d}\) unchanged.
\end{proposition}
\begin{proof}
Expand the \(m+1\) coefficients of \(B_{f,d}\) in the basis \((e_\ell)\), obtaining a \(d\times(m+1)\) matrix.
Its column rank is \(\eta_f(d)\), and its row space corresponds to \(\mathcal W_{f,d}\). Equality of row and column ranks proves the dimension identity.
Changing basis multiplies the matrix on the left by an invertible matrix, preserving its row space.
The upper bound follows from the number of coefficients and \(\dim_k K=d\). If the rank were one, the leading coefficient one would force all coefficients to lie in \(k\),
contradicting~\eqref{str:eq-block-coefficient-field} and \(d>1\). Thus the lower bound is two.
Changing the conjugate root applies some \(\sigma^j\) simultaneously to all coefficients. This is an invertible \(k\)-linear map on \(K\), so again preserves rank and row space.
\end{proof}

The trace pairing also gives a direct description of the coordinate-polynomial space:
\[
\mathcal W_{f,d}
=\{\Tr_{K/k}(zB_{f,d}(X)):z\in K\},
\]
where trace acts coefficientwise. If \(e_1^\vee,\ldots,e_d^\vee\) is the trace-dual basis to \(e_1,\ldots,e_d\),
then \(\Tr_{K/k}(e_\ell^\vee B_{f,d})=B_\ell\); the trace for any \(z\) is a \(k\)-linear combination of these \(B_\ell\), proving the equality.
Thus \(\mathcal W_{f,d}\) collects all base-field polynomials extracted from the coefficients of the same relative root block,
and \(\eta_f(d)\) is the least number of independent base-field polynomials needed to reconstruct that block.
Rank two yields a numerator and denominator. While retaining representation degree \(m\), the denominator can be chosen constant if and only if this space contains the constant polynomial \(1\).

\subsection{Rational Representations Attaining the Degree Lower Bound}

\begin{definition}\label{str:def-rational-degree}
For relatively prime \(A,C\in k[X]\) with \(C\ne0\), define the \textbf{degree} of the nonconstant rational function \(R=A/C\) by
\[
\deg R=\max\{\deg A,\deg C\}.
\]
Common factors of numerator and denominator must first be canceled; multiplying both by a nonzero constant does not change the degree.
\end{definition}

First consider representations of an arbitrary prescribed element, without restricting the subfield to which it belongs. For nonnegative integers \(a,b\), let \(\mathcal R_{a,b}\) be the set of rational functions \(A/B\) satisfying
\[
A,B\in k[X],\qquad B\text{ monic},\qquad
\gcd(A,B)=1,\qquad \deg A\le a,\quad\deg B\le b,
\]
where the zero function has the unique normalized representation \(0/1\).

\begin{proposition}[Rational reconstruction with separate degree bounds]\label{v70:prop-rational-reconstruction}
Let \(E=k(\alpha)\), with \([E:k]=n\), and let \(0\le a,b<n\).
If \(a+b\ge n-1\), every \(\beta\in E\) has a representation
\[
\beta=\frac{A(\alpha)}{B(\alpha)},\qquad A/B\in\mathcal R_{a,b}.
\]
If \(a+b<n\), evaluation defines an injection \(\mathcal R_{a,b}\to E\).
In particular, evaluation is a bijection when \(a+b=n-1\).
\end{proposition}
\begin{proof}
Since \(\deg B<n\), every nonzero denominator satisfies \(B(\alpha)\ne0\), so evaluation is defined throughout.
Suppose \(a+b\ge n-1\). Expressing \(A(\alpha)-\beta B(\alpha)=0\) in a basis of \(E/k\) gives at most \(n\) homogeneous linear conditions on \(a+b+2\ge n+1\) unknown coefficients. There is therefore a nonzero solution \((A,B)\).
If \(B=0\), then \(A(\alpha)=0\) and \(\deg A<n\) force \(A=0\), a contradiction. Thus \(B\ne0\) and \(B(\alpha)\ne0\).
Canceling common factors and making the denominator monic do not increase the degrees, giving the required representation.

If \(A/B,A'/B'\in\mathcal R_{a,b}\) have the same value at \(\alpha\), then
\[
(AB'-A'B)(\alpha)=0,\qquad\deg(AB'-A'B)\le a+b<n.
\]
Hence \(AB'-A'B=0\), so the rational functions are equal. Coprimality and the monic denominator then make their normalized representatives equal as well.
When \(a+b=n-1\), existence and injectivity hold simultaneously, proving bijectivity.
\end{proof}

The uniqueness here concerns reduced rational functions, not unreduced solution vectors of the original homogeneous system.
Writing \(P_j=\operatorname{span}_k\{1,\alpha,\ldots,\alpha^j\}\), a representation with these bounds exists for \(\beta\ne0\) if and only if \(P_a\cap\beta P_b\ne\{0\}\); for \(\beta=0\), use \(0/1\).
The cardinality of \(\mathcal R_{a,b}\) itself is computed in Exercise~\ref{v70:ex-rational-rectangle-count}.

Representing a prescribed element, representing a generator of a prescribed subfield, and preserving a generator of the whole extension while making its trace zero impose different requirements.
The preceding proposition treats the first question through intersections of low-degree power spaces. The root-block coefficient rank below determines when the second attains the bound \(n/d\), while Proposition~\ref{v70:prop-trace-zero-generator} makes the trace zero without shortening the orbit.

\begin{definition}

\leavevmode\phantomsection\label{str:def-rational-generator-degree}\ignorespaces
Define the \textbf{minimum rational generator-representation degree} of \(K\) by
\begin{equation}\label{str:eq-rational-kappa}
\kappa_f^{\mathrm{rat}}(d)
=\min\left\{\deg(A/C):
\begin{array}{l}
 A,C\in k[X],\ \gcd(A,C)=1,\ C(\alpha)\ne0,\\
 k(A(\alpha)/C(\alpha))=K
\end{array}\right\}.
\end{equation}
\end{definition}

Taking \(C=1\) includes all polynomial representations, so the minimum exists and \(\kappa_f^{\mathrm{rat}}(d)\le\kappa_f(d)\).
For rational-function degree and substitution after clearing denominators, see also~\cite{FFSSchulz2023}; here we determine from a given relative root block whether a rational generator representation attains the degree lower bound.

\begin{theorem}[Root-block coefficient rank and the rational representation lower bound]\leavevmode\phantomsection\label{str:thm-rational-rank}
We have
\begin{equation}\label{str:eq-rational-lower}
\kappa_f^{\mathrm{rat}}(d)\ge m,
\end{equation}
with equality if and only if \(\eta_f(d)=2\). When this holds, one may choose \(\beta\in K\), with \(k(\beta)=K\), and relatively prime \(A,C\in k[X]\) such that
\begin{equation}\label{str:eq-normalized-rational-block}
B_{f,d}=A-\beta C,\qquad
A\text{ is monic},\quad\deg A=m,\quad0\le\deg C<m.
\end{equation}
If \(P_\beta\) is the monic minimal polynomial of \(\beta\) over \(k\), then
\begin{equation}\label{str:eq-rational-transform}
\boxed{\quad f(X)=C(X)^d P_\beta\left(\frac{A(X)}{C(X)}\right).\quad}
\end{equation}
The right side is the polynomial obtained by clearing denominators, not an ordinary polynomial composition \(P_\beta\circ h\).
\end{theorem}
\begin{proof}
Suppose \(\beta=A(\alpha)/C(\alpha)\) generates \(K/k\). Since \(\beta\notin k\), the polynomial \(A-\beta C\) is nonzero;
otherwise comparison at any nonzero coefficient of \(C\) would give \(\beta\in k\).
It vanishes at \(\alpha\), so \(B_{f,d}\mid A-\beta C\), and hence
\[
m\le\deg(A-\beta C)\le\max\{\deg A,\deg C\}.
\]
This proves the lower bound. If equality holds, let \(\ell\in K^\times\) be the leading coefficient of \(A-\beta C\). Then
\[
B_{f,d}=\ell^{-1}(A-\beta C).
\]
All coefficients lie in the two-dimensional space \(\ell^{-1}(k+k\beta)\), so \(\eta_f(d)\le2\); Proposition~\ref{str:prop-coefficient-spaces} gives \(\eta_f(d)=2\).

Conversely, if \(\eta_f(d)=2\), choose a basis \(1,\beta\) of \(\mathcal V_{f,d}\).
Writing coefficientwise \(B_{f,d}=A-\beta C\) gives \(A,C\in k[X]\). The leading coefficient one makes \(A\) monic with \(\deg A=m\) and \(\deg C<m\);
rank two also ensures \(C\ne0\). From \(\mathcal V_{f,d}\subseteq k(\beta)\) and~\eqref{str:eq-block-coefficient-field}, we get \(k(\beta)=K\).
If \(A,C\) had a nonconstant common factor, it would divide the irreducible polynomial \(B_{f,d}\) in \(K[X]\),
but would have degree at most \(\deg C<m\), a contradiction. Thus they are relatively prime.
Also, \(C\ne0\) and \(\deg C<n\), so the degree \(n\) of \(\alpha\) over \(k\) gives \(C(\alpha)\ne0\).
Evaluation at \(\alpha\) yields \(\beta=A(\alpha)/C(\alpha)\), attaining the lower bound \(m\).

Finally, apply \(\sigma^j\) to the coefficients for \(0\le j<d\). The root blocks are distinct and exhaust all roots of \(f\), so
\[
f(X)=\prod_{j=0}^{d-1}\bigl(A(X)-\beta^{q^j}C(X)\bigr).
\]
Factoring \(C(X)\) out of each factor gives~\eqref{str:eq-rational-transform}.
Under the chosen monic normalization, both sides have leading term \(X^{dm}\), so no additional scalar factor is needed.
\end{proof}

When the lower bound is attained, each fiber of \(A/C\) on the original root set is exactly one \(\sigma^d\)-root block.
Coprimality ensures that roots of the fiber equation \(A-\beta^{q^j}C=0\) are not zeros of the denominator; its degree \(m\) ensures that the root block exhausts all roots of this equation.
The root block need not be an orbit of an affine group.

\subsection{Removal of the Denominator}

\begin{corollary}[The constant direction and a polynomial inner component]\label{str:cor-polynomial-direction}
If \(\eta_f(d)=2\), then
\begin{equation}\label{str:eq-polynomial-direction}
\kappa_f(d)=m\quad\Longleftrightarrow\quad1\in\mathcal W_{f,d}.
\end{equation}
In general, attainment of the degree lower bound by the two kinds of representation falls into three cases:
\begin{enumerate}
\item \(\eta_f(d)=2\) and \(1\in\mathcal W_{f,d}\): both polynomial and rational representations attain \(m\);
\item \(\eta_f(d)=2\) and \(1\notin\mathcal W_{f,d}\): only rational representations attain \(m\);
\item \(\eta_f(d)>2\): neither kind attains \(m\).
\end{enumerate}
\end{corollary}
\begin{proof}
In~\eqref{str:eq-normalized-rational-block}, we have \(\mathcal W_{f,d}=\operatorname{span}_k\{A,C\}\).
It contains a polynomial of degree \(m\). If it also contains \(1\), choose a basis consisting of a monic \(h\) and \(1\), with \(\deg h=m\).
Then \(B_{f,d}=\xi h+\zeta\), with \(\xi,\zeta\in K\); comparison of leading coefficients gives \(\xi=1\).
Thus all its nonconstant coefficients lie in \(k\), and Proposition~\ref{sfr:prop-block-coefficients} gives \(\kappa_f(d)=m\).
Conversely, a polynomial representation attaining the lower bound gives \(B_{f,d}=h-\beta\), with \(k(\beta)=K\), hence \(\mathcal W_{f,d}=\operatorname{span}_k\{h,1\}\).
The three cases now follow from Theorem~\ref{str:thm-rational-rank} and \(\kappa_f^{\mathrm{rat}}(d)\le\kappa_f(d)\).
\end{proof}

\begin{corollary}[Fractional linear classification of rational representations attaining the lower bound]\label{str:cor-rational-uniqueness}
Suppose \(\eta_f(d)=2\), and take \(A,C\) from~\eqref{str:eq-normalized-rational-block}.
All rational functions of degree \(m\) whose values at \(\alpha\) generate \(K/k\) are exactly
\begin{equation}\label{str:eq-rational-mobius-class}
\frac{aA+bC}{cA+eC},\qquad a,b,c,e\in k,\quad ae-bc\ne0.
\end{equation}
Thus it is the two-dimensional space \(\mathcal W_{f,d}\) that is uniquely determined; changing its ordered basis corresponds to an outer fractional linear transformation.
\end{corollary}
\begin{proof}
In~\eqref{str:eq-rational-mobius-class}, the numerator and denominator are relatively prime, since inverse linear combinations recover \(A,C\).
Their degrees cannot both be less than \(m\), so the rational function has degree \(m\).
Its value at \(\alpha\) is \((a\beta+b)/(c\beta+e)\); the denominator is nonzero and the inverse transformation recovers \(\beta\), so this value still generates \(K\).

Conversely, if \(\widetilde A/\widetilde C\) attains the lower bound, the proof of necessity in the preceding theorem gives
\(B_{f,d}=\ell^{-1}(\widetilde A-\gamma\widetilde C)\), where \(k(\gamma)=K\).
Multiplying all coefficients in \(K\) simultaneously by \(\ell\ne0\) is an invertible \(k\)-linear transformation, so does not change the coordinate-polynomial space.
Since \(1,\gamma\) are linearly independent, \(\mathcal W_{f,d}=\operatorname{span}_k\{\widetilde A,\widetilde C\}\).
The invertible change between the two bases gives~\eqref{str:eq-rational-mobius-class}.
\end{proof}

\begin{remark}[Rational functions and generator values]\label{v70:rem-rational-functions-values}
Under the assumptions \(1<d<n\), \(m=n/d\), and \(\eta_f(d)=2\), the classification above counts distinct functions in \(k(X)\), not automatically distinct values at \(\alpha\).
If \(d>2\), then \(2m<n\), and the denominator-clearing argument in Proposition~\ref{v70:prop-rational-reconstruction} makes evaluation injective.
If \(d=2\), then \(2m=n\), and distinct functions can have the same value: the \(q(q^2-1)\) rational functions of minimum degree attain all \(q^2-q\) generator values in \(K_2\setminus k\), each exactly \(q+1\) times. See Exercise~\ref{str:ex-count-rational} for the proof.
\end{remark}

\begin{corollary}[Quadratic subfields]\label{str:cor-quadratic-rational}
If \(n>2\) is even, then every monic irreducible polynomial \(f\) of degree \(n\) satisfies
\[
\kappa_f^{\mathrm{rat}}(2)=n/2.
\]
\end{corollary}
\begin{proof}
Take \(d=2\). Proposition~\ref{str:prop-coefficient-spaces} forces \(\eta_f(2)=2\); apply Theorem~\ref{str:thm-rational-rank}.
\end{proof}

\subsection{Obstructions That Denominators Can and Cannot Remove}

\begin{example}[The two-dimensional coefficient space of a quartic root block]\label{str:ex-rational-four}
Use \(f_2=X^4+X^3+1\) from Example~\ref{sfr:ex-quartic-comparison}, take a root \(\gamma\), and put
\[
b=\frac{\gamma^2}{\gamma+1}=\gamma^3+\gamma+1=\eta+1.
\]
Here \(\eta=\beta\) retains its meaning from the opening example; \(b\) is the other non-base-field element of the four-element subfield.
The relation \(\gamma^4=\gamma^3+1\) verifies this equality and \(b^2+b+1=0\), so \(b\) generates the four-element subfield.
The relative root-block polynomial is
\[
B_{f_2,2}(X)=X^2+b X+b
=X^2-b(X+1).
\]
Thus \(\mathcal W_{f_2,2}=\operatorname{span}_{\F_2}\{X^2,X+1\}\), which does not contain \(1\).
Taking \(A=X^2\), \(C=X+1\), and \(P=Y^2+Y+1\) gives
\[
f_2=A^2+AC+C^2=C^2P(A/C),
\qquad
\kappa_{f_2}^{\mathrm{rat}}(2)=2<3=\kappa_{f_2}(2).
\]
This root block is a complete fiber of a quadratic rational function, but cannot be a complete fiber of a quadratic polynomial.
The functions \(R=A/C\) and \(1+1/R=(A+C)/A\) are distinct rational functions of degree two, but
\[
R-\left(1+\frac1R\right)=\frac{A^2+AC+C^2}{AC}=\frac{f_2}{AC}
\]
shows that they have the same value at \(\gamma\). Here \(A(\gamma)C(\gamma)\ne0\), and the cross product has degree exactly four, at the boundary where the strict degree condition fails.
\end{example}

\begin{example}[The three-dimensional coefficient space of a sextic root block]\label{str:ex-rational-six}
Use \(f=X^6+X+1\), \(K=\F_8\), and \(\alpha,u,v\) from Example~\ref{sfr:ex-six-defect-one}.
Since \(\alpha^8=\alpha^3+\alpha^2\) and \(\alpha^9=\alpha^4+\alpha^3\), we have
\[
B_{f,3}(X)=(X-\alpha)(X-\alpha^8)=X^2+uX+v.
\]
The elements \(1,u,v\) are linearly independent over \(\F_2\), so \(\eta_f(3)=3\), and there is no quadratic rational generator representation.
The existing cubic polynomial representation gives an upper bound, hence
\[
\kappa_f^{\mathrm{rat}}(3)=\kappa_f(3)=3>2=n/3.
\]
The obstruction is the three-dimensional base-field span of the root-block coefficients; allowing a denominator therefore cannot lower the generator-representation degree to two.
\end{example}

\section{Recovering Root Blocks and Representation Structure from Relative Traces}
\label{str:sec-trace-recovery}

The trace-rank function over \(k\) records the subfield degree profile; over the subfield \(K\), the same trace values satisfy a linear recurrence.
Its minimal polynomial is exactly the relative root-block polynomial, so the relative trace sequence also recovers the coefficient rank and rational representation of the preceding section.
Continue to use \(c_j\) from~\eqref{str:eq-relative-trace-sequence} and \(B_{f,d}\) from~\eqref{str:eq-block-polynomial}.

\subsection{The Shortest Recurrence from a Relative Root Block}

\leavevmode\phantomsection\label{str:def-trace-recurrence}\ignorespaces
If a nonzero \(R(T)=\sum_{i=0}^{\ell}r_iT^i\in K[T]\) satisfies
\[
\sum_{i=0}^{\ell}r_i c_{j+i}=0\qquad(j\ge0),
\]
it is called a \textbf{recurrence polynomial} of the sequence \((c_j)\).
The monic one of least degree is called the minimal polynomial of the sequence over \(K\).
Here recurrence is expressed by this identity; Chapter~\ref{chap:sequences} develops the general theory of linear recurring sequences using the shift operator.

The role of relative trace can be expressed through a general pairing. The following result requires neither a field quotient nor a squarefree defining polynomial.

\leavevmode\phantomsection\label{str:def-hankel}\ignorespaces
A matrix whose entries depend only on the sum of the row and column indices is a \textbf{Hankel matrix}; its entries are constant along each antidiagonal.

\begin{lemma}[Nondegenerate pairings and recovery from finite observations]\label{v127:lem-pairing-recovery}
Let \(K\) be a field, let \(g\in K[T]\) be monic of degree \(d\ge1\), and put \(A=K[T]/(g)\) and \(t=\overline T\).
Suppose a \(K\)-linear functional \(\ell:A\to K\) makes the pairing \((u,v)\mapsto\ell(uv)\) nondegenerate. Then the sequence \(s_j=\ell(t^j)\) has minimal polynomial \(g\).
Put
\[
H=(s_{i+j})_{0\le i,j<d},\qquad H^+=(s_{i+j+1})_{0\le i,j<d}.
\]
The matrix \(H\) is invertible. If \(g=T^d+\sum_{i=0}^{d-1}g_iT^i\), then
\begin{equation}\label{v127:eq-pairing-recovery}
H(g_0,\ldots,g_{d-1})^{\mathsf t}=-(s_d,\ldots,s_{2d-1})^{\mathsf t},
\qquad M_t=H^{-1}H^+,
\end{equation}
where \(M_t\) represents multiplication by \(t\) in the power basis \(1,t,\ldots,t^{d-1}\). Thus \(d\) and the first \(2d\) terms determine both \(g\) and \(M_t\).
\end{lemma}
\begin{proof}
A polynomial \(R\) gives a recurrence if and only if \(\ell(t^jR(t))=0\) for all \(j\ge0\). Since the powers span \(A\) and the pairing is nondegenerate, this is equivalent to \(R(t)=0\), or \(g\mid R\).
The matrix \(H\) is the Gram matrix of the pairing in the power basis, so it is invertible. Multiplying \(g(t)=0\) successively by \(1,t,\ldots,t^{d-1}\) and applying \(\ell\) gives the first identity.
Column \(j\) of \(M_t\) contains the coordinates of \(t^{j+1}\), so \(HM_t=H^+\), proving the second identity.
\end{proof}

\begin{corollary}[Minimal polynomial of the relative trace sequence]\label{str:thm-minimal-recurrence}
The minimal polynomial of \((c_j)\) over \(K\) is \(B_{f,d}(T)\). Explicitly,
\begin{equation}\label{str:eq-relative-recurrence}
c_{j+m}+b_{m-1}c_{j+m-1}+\cdots+b_0c_j=0
\qquad(j\ge0),
\end{equation}
and every recurrence polynomial \(R\in K[T]\) is divisible by \(B_{f,d}\).
\end{corollary}
\begin{proof}
Apply Lemma~\ref{v127:lem-pairing-recovery} with \(A=E\simeq K[T]/(B_{f,d})\), \(t=\alpha\), and
\[
\ell(u)=\Tr_{E/K}\left(\frac{u}{f'(\alpha)}\right).
\]
Since \(f'(\alpha)\ne0\), nondegeneracy of the trace pairing makes \(\ell(uv)\) nondegenerate, giving the result.
\end{proof}

The same \(c_j\) may be viewed as a vector over \(k\) or as a term of a sequence over \(K\).
The first viewpoint studies linear independence of trace values over \(k\); the second allows recurrence coefficients in \(K\).
These two linear structures yield the degree profile and the relative root block, respectively; their coefficient fields cannot be interchanged.

\subsection{Hankel Matrices and Recovery from a Finite Prefix}

\phantomsection\label{v113:par-hankel-question}
Given \(m=[E:K]\) and \(c_0,\ldots,c_{2m-1}\), take the recurrence at \(m\) consecutive positions to recover the root-block coefficients.
Its coefficient matrix is the Gram matrix of the nondegenerate pairing above, so the system has a unique solution. Put
\begin{equation}\label{str:eq-hankel}
\mathcal H_{f,d}=(c_{i+j})_{0\le i,j<m}\in M_m(K).
\end{equation}

\begin{corollary}[Recovering a root block from the first \(2m\) trace values]\leavevmode\phantomsection\label{str:thm-hankel-recovery}
The matrix \(\mathcal H_{f,d}\) is invertible, and the coefficients of \(B_{f,d}\) are uniquely determined by
\begin{equation}\label{str:eq-hankel-recovery}
\boxed{\quad
\mathcal H_{f,d}
\begin{pmatrix}b_0\\b_1\\\vdots\\b_{m-1}\end{pmatrix}
=-\begin{pmatrix}c_m\\c_{m+1}\\\vdots\\c_{2m-1}\end{pmatrix}.
\quad}
\end{equation}
Thus \(c_0,\ldots,c_{2m-1}\) recover the relative root block without first factoring \(f\) in \(K[X]\).
\end{corollary}
\begin{proof}
Use the pairing in Corollary~\ref{str:thm-minimal-recurrence}. Equation~\eqref{v127:eq-pairing-recovery} gives the stated system and its unique solution.
\end{proof}

To test only whether a polynomial inner component of degree \(m\) exists, first use the kernel of \(M_{f,d,m}\) to test \(\rho_{f,d}(m)=2\); the entire root block is unnecessary.
To compare polynomial and rational representations, its coefficient space is also needed, and the Hankel equation supplies it.
The length \(2m\) is sufficient for this recovery procedure; no minimality is asserted. A shorter-observation question must specify the available information. If the full polynomial \(f\) and quotient-field operations are given, the traces can be recomputed. If only a finite trace prefix is given, a prior bound on the recurrence order is needed; see Corollary~\ref{cor:bm-recovery}.

This invertibility follows from nondegeneracy of the weighted trace pairing and remains valid when \(p\mid m\) or \(c_0=0\).

The same lemma also gives
\[
C_B=\mathcal H_{f,d}^{-1}\mathcal H_{f,d}^{+},\qquad
\mathcal H_{f,d}^{+}=(c_{i+j+1})_{0\le i,j<m},
\]
where \(C_B\) represents multiplication by \(\alpha\) in the power basis. The first \(2m\) traces thus also recover the multiplication representation of Section~\ref{sec:regular-representation}.

\begin{theorem}[The Hankel--Vandermonde determinant for exponential moments]\label{v88:thm-hankel-vandermonde}
Let \(z_1,\ldots,z_w\in K^\times\) be pairwise distinct and
\(a_1,\ldots,a_w\in K^\times\), and put
\(s_r=\sum_{j=1}^w a_jz_j^r\) for \(r\ge0\). Then for every \(r\ge0\),
\begin{equation}\label{v88:eq-hankel-vandermonde}
\boxed{\quad\det(s_{r+i+j})_{0\le i,j<w}
=\left(\prod_{j=1}^w a_jz_j^r\right)
 \prod_{1\le i<j\le w}(z_j-z_i)^2\ne0.\quad}
\end{equation}
Consequently the sequence has minimum recurrence order \(w\) over \(K\), with
monic minimum recurrence polynomial \(\prod_{j=1}^w(T-z_j)\).
\end{theorem}
\begin{proof}
Put \(V=(z_j^i)_{0\le i<w,\,1\le j\le w}\). Then
\[
(s_{r+i+j})_{i,j}=V\,\operatorname{diag}(a_1z_1^r,\ldots,a_wz_w^r)\,V^{\mathsf T}.
\]
Taking determinants gives the formula. Put \(g(T)=\prod_j(T-z_j)\). Under the evaluation isomorphism \(K[T]/(g)\simeq K^w\), the product pairing associated with
\[
\ell(u)=\sum_{j=1}^w a_j u(z_j)
\]
has diagonal matrix \(\operatorname{diag}(a_1,\ldots,a_w)\). It is nondegenerate, so Lemma~\ref{v127:lem-pairing-recovery} gives the minimal recurrence polynomial \(g\).
\end{proof}

Weighted relative trace and weighted evaluation at distinct nodes both give nondegenerate product pairings: the former uses a power basis over the subfield, the latter evaluation coordinates over a splitting field.
When some weights vanish, or the ordinary algebra trace cannot detect nilpotent directions, the ideal removed by observation must first be determined; see Proposition~\ref{v127:prop-observation-kernel}.
The minimal recurrence order of an entire sequence still cannot be inferred from the rank of an arbitrary short prefix.

\begin{example}[Recovering traces and a root block in the quartic example]\leavevmode\phantomsection\label{str:ex-hankel-four}
First recover the root block for \(f_1\) from the opening Example~\ref{sfr:ex-quartic-comparison}.
Since \(f_1'=1\), direct calculation gives

\[
(c_0,c_1,c_2,c_3)=(0,1,1,\beta+1),\qquad
\mathcal H_{f_1,2}=\begin{pmatrix}0&1\\1&1\end{pmatrix}.
\]
The Hankel equation gives \((b_0,b_1)=(\beta,1)\), recovering \(B_{f_1,2}=X^2+X+\beta\).
Although \(c_0=0\), the determinant is one; the base-field ranks of the trace prefixes are \((0,1,1,2)\).

For \(f_2,\gamma,b\) in Example~\ref{str:ex-rational-four}, we have \(f_2'(\gamma)=\gamma^2\).
Using \(\Tr_{\F_{16}/\F_4}(z)=z+z^4\) and \(\gamma^4=\gamma^3+1\), we obtain
\[
(c_0,c_1,c_2,c_3)=(1,1,0,b),\qquad
\mathcal H_{f_2,2}=\begin{pmatrix}1&1\\1&0\end{pmatrix}.
\]
For example, \(\gamma^{-1}=\gamma^3+\gamma^2\) and \(\gamma^{-2}=\gamma^2+\gamma\);
substitution into the relative trace gives the first two terms, while the last two are \(\Tr(1)=0\) and \(\Tr(\gamma)=b\).
Equation~\eqref{str:eq-hankel-recovery} then gives \(b_0=b_1=b\), recovering
\[
B_{f_2,2}(X)=X^2+b X+b.
\]
Meanwhile, the \(\F_2\)-ranks of the four trace-value prefixes grow as \((1,1,1,2)\); duality recovers the degree profile \((1,1,1,2)\).
The root-block coefficient rank is two, but the coordinate space contains no nonzero constant polynomial, giving minimum rational generator-representation degree two and minimum polynomial generator-representation degree three.
\end{example}

\begin{example}[Four trace values recover the representation structure of a sextic root]\label{str:ex-hankel-six}
Use \(f=X^6+X+1\), \(K=\F_8\), and \(\alpha,u,v\) from Example~\ref{str:ex-rational-six}.
Here \(n=6,d=3,m=2\), so \(2m=4<n\). Using \(f'=1\) and \(\alpha^6=\alpha+1\), compute directly
\(c_j=\alpha^j+\alpha^{8j}\), obtaining the first four terms
\[
 (c_0,c_1,c_2,c_3)=(0,u,1+u+v,1+v).
\]
We already know \(u^2=1+u+v\) and \(v^2=1+u\); the relation \(\alpha^6=\alpha+1\) also gives \(uv=u+v\). The Hankel equation is
\[
 \begin{pmatrix}0&u\\u&1+u+v\end{pmatrix}
 \begin{pmatrix}b_0\\b_1\end{pmatrix}
 =\begin{pmatrix}1+u+v\\1+v\end{pmatrix}.
\]
Its determinant is \(u^2\ne0\), and its unique solution is \(b_0=v,b_1=u\), recovering
\[
 B_{f,3}(X)=X^2+uX+v,\qquad c_{j+2}=u c_{j+1}+v c_j.
\]
The recurrence recovered from the first four terms then gives \(c_4,c_5\). The successive \(\F_2\)-ranks are:
\begin{center}
\renewcommand{\arraystretch}{1.15}
\begin{tabular}{cclc}
\toprule
\(j\)&\(c_j\)&How the term is obtained&\(r_{f,3}(j)\)\\
\midrule
0&\(0\)&Relative trace&0\\
1&\(u\)&Relative trace&1\\
2&\(1+u+v\)&Relative trace&2\\
3&\(1+v\)&Relative trace&2\\
4&\(1+v\)&\(u c_3+v c_2\)&2\\
5&\(1+u\)&\(u c_4+v c_3\)&3\\
\bottomrule
\end{tabular}
\end{center}
The three trace-rank jump positions are \(\mu_0=1,\mu_1=2,\mu_2=5\).
Equation~\eqref{str:eq-profile-trace-duality} and the reversed-jump formula give
\[
 (\rho_{f,3}(0),\ldots,\rho_{f,3}(5))=(1,1,1,2,3,3),
 \qquad (\lambda_0,\lambda_1,\lambda_2)=(0,3,4),
\]
hence \(\Delta_f(3)=(3-2)+(4-4)=1\).
On the other hand, linear independence of \(1,u,v\) gives \(\eta_f(3)=3\), so a rational generator representation cannot attain the degree lower bound two.
The same short prefix recovers the root block, the degree profile, and the degree obstruction to rational representation.
\end{example}

Power sums along a cyclotomic orbit provide another set of observations from which a polynomial can be recovered linearly.
\begin{corollary}[Recovering an irreducible factor from period sums]\label{v97:cor-period-recovery}
Suppose \(p\nmid n\), let \(\zeta\) be a primitive \(n\)th root of unity, and let \(C\) be a \(q\)-cyclotomic coset modulo \(n\), of size \(d\). Put
\[
 h_C(X)=\prod_{u\in C}(X-\zeta^u)=X^d+\sum_{i=0}^{d-1}h_iX^i,
 \qquad s_j=\sum_{u\in C}\zeta^{ju}\quad(j\ge0).
\]
These \textbf{period sums of roots of unity} belong to \(\F_q\), the matrix \(H=(s_{i+j})_{0\le i,j<d}\) is invertible, and
\begin{equation}\label{v97:eq-period-recovery}
 H(h_0,\ldots,h_{d-1})^{\mathsf T}
 =-(s_d,\ldots,s_{2d-1})^{\mathsf T}.
\end{equation}
Thus the first \(2d\) period sums determine \(h_C\), even when \(p\mid d\).
\end{corollary}
\begin{proof}
Multiplication by \(q\) permutes \(C\), so \(s_j^q=s_j\).
Over an extension containing \(\zeta\), apply Theorem~\ref{v88:thm-hankel-vandermonde} to the distinct nodes \(\zeta^u\), with every weight equal to one, to obtain \(\det H\ne0\).
This determinant already lies in \(\F_q\), so \(H\) is invertible over the base field as well.
Multiply \(h_C(\zeta^u)=0\) by \(\zeta^{ju}\) and sum over \(u\) to obtain the recurrence.
Taking \(j=0,\ldots,d-1\) gives \eqref{v97:eq-period-recovery}.
\end{proof}

The root product, Newton identities, and Hankel recurrence recovery use the same conjugate data.
When Newton identities are used to determine coefficients successively, division by \(j\) is available only if \(p\nmid j\); the linear recovery above has no such restriction.

\begin{example}[Period-sum recovery when the degree equals the characteristic]\label{v97:ex-f9-period-recovery}
In an algebraic closure of characteristic three, take a primitive seventh root \(\zeta\) and put \(\eta=\zeta+\zeta^2+\zeta^4\).
The sum over the other nonzero exponent orbit is \(-1-\eta\). Expanding the square gives
\[
 \eta^2=\eta+2(-1-\eta)=1-\eta.
\]
Since \(Y^2+Y+2\) has no root in \(\F_3\), one has \(K=\F_3(\eta)=\F_9\).
As \(9\equiv2\pmod7\), the set \(C=\{1,2,4\}\) is a \(9\)-cyclotomic orbit of length three.
Its first six period sums are
\[
 s_0=0,\qquad s_1=s_2=s_4=\eta,\qquad s_3=s_5=-1-\eta.
\]
Hence
\[
 H=\begin{pmatrix}
 0&\eta&\eta\\
 \eta&\eta&-1-\eta\\
 \eta&-1-\eta&\eta
 \end{pmatrix},\qquad \det H=-1\ne0.
\]
Using \(\eta^2=1-\eta\) in \eqref{v97:eq-period-recovery} gives the unique solution
\((h_0,h_1,h_2)=(-1,-1-\eta,-\eta)\), and therefore recovers
\[
 h_C(X)=X^3-\eta X^2+(-1-\eta)X-1.
\]
Here \(s_0=3=0\), and the third Newton identity does not allow division by three, yet the Hankel matrix remains invertible.
A zero trace value does not mean that the trace pairing is degenerate; the failure of one division does not mean that the full power-sum data no longer determine the polynomial.
\end{example}

The period sums here take values in finite extension fields.
The Gaussian periods defined through complex-valued additive characters in Section~\ref{v95:sec-subgroup-periods} are a different numerical realization.
Their value fields and operations differ, so the latter's absolute-value estimates cannot be applied directly to these finite-field-valued sequences.

\subsection{Power-Trace Prefixes and Polynomial Identification}\label{v123:sec-power-trace-fibers}

For a monic square-free polynomial \(f\in k[X]\), write
\[
s_j(f)=\sum_{f(\alpha)=0}\alpha^j\qquad(j\ge0).
\]
If \(k=\F_q\) and \(f\) is irreducible of degree \(m\), then \(s_j(f)=\Tr_{\F_{q^m}/k}(\alpha^j)\).
In Corollary~\ref{v97:cor-period-recovery}, the degree already determines \(s_0=m\), so the first \(2m-1\) power traces \(s_1,\ldots,s_{2m-1}\) identify the corresponding root orbit.
Which polynomials remain indistinguishable if the last entry is omitted?

\begin{lemma}[Vanishing power sums and a constant derivative]\label{v123:lem-moment-derivative}
Let \(x_1,\ldots,x_w\) be distinct elements of a field, where \(w\ge2\), and put \(D(X)=\prod_{i=1}^w(X-x_i)\). Then
\[
\sum_{i=1}^w x_i^j=0\quad(0\le j\le w-2)
\quad\Longleftrightarrow\quad D'(X)=c\ne0
\]
for some constant \(c\). In this case, \(\sum_i x_i^{w-1}=c\).
\end{lemma}
\begin{proof}
Apply Lagrange interpolation to \(X^j\) and compare coefficients of \(X^{w-1}\). This gives
\begin{equation}\label{v123:eq-lagrange-moments}
\sum_{i=1}^w\frac{x_i^j}{D'(x_i)}=
\begin{cases}0,&0\le j\le w-2,\\1,&j=w-1.\end{cases}
\end{equation}
The matrix \((x_i^j)_{0\le j\le w-2,\,1\le i\le w}\) has rank \(w-1\), so its kernel is spanned by \((D'(x_i)^{-1})_i\).
If the first \(w-1\) unweighted power sums vanish, then \((1,\ldots,1)\) also lies in this kernel. Thus all \(D'(x_i)\) equal the same nonzero constant \(c\).
Since \(D'-c\) has degree less than \(w\) and vanishes at every \(x_i\), it is zero. The converse and the final identity follow from \eqref{v123:eq-lagrange-moments}.
\end{proof}

\begin{theorem}[Fibers of a power-sum prefix]\label{v123:thm-power-trace-pencil}
Let \(k=\F_q\), and let \(f\in k[X]\) be monic and square-free of degree \(m\ge1\). Among all monic square-free polynomials of degree \(m\), consider those sharing the first \(2m-2\) power sums of \(f\).
\begin{enumerate}
\item If \(f''\ne0\), the only such polynomial is \(f\).
\item If \(f''=0\), let \(h_f\) be the unique polynomial satisfying
\[
f'h_f\equiv1\pmod f,\qquad \deg h_f<m.
\]
The required polynomials are exactly \(f_c=f+c h_f\), for \(c\in k\), and
\begin{equation}\label{v123:eq-last-power-trace}
s_{2m-1}(f_c)=s_{2m-1}(f)-c.
\end{equation}
They are monic and square-free, with \(h_{f_c}=h_f\); for \(a\ne b\),
\[
\gcd(f_a,f_b)=1,\qquad f_a'f_b-f_af_b'=b-a.
\]
\end{enumerate}
For \(m=1\), the condition on the first \(2m-2\) power sums is empty.
\end{theorem}
\begin{proof}
First suppose that \(g\ne f\) is another monic square-free polynomial of degree \(m\) with the same prefix.
If they share a root, remove all common roots, leaving \(r\) distinct roots with \(1\le r\le2m-2\). Give the remaining roots of \(f\) weight one and those of \(g\) weight minus one. Their weighted power sums vanish for \(j=0,\ldots,2m-2\).
The first \(r\) equations and Vandermonde invertibility force every weight to vanish, a contradiction. Thus \(f\) and \(g\) are coprime.

Apply \eqref{v123:eq-lagrange-moments} to the \(2m\) roots of \(D=fg\). The corresponding Vandermonde matrix has a one-dimensional kernel, so the signed weights above equal \(c/D'(x_i)\) for some \(c\ne0\).
Thus \(f'g-fg'\) takes the value \(c\) at every root of either polynomial. Its degree is less than \(2m\), so it is the constant \(c\), which belongs to \(k\). The same interpolation identity gives the converse and the difference of the next terms. Hence, for distinct monic square-free polynomials of degree \(m\),
\begin{equation}\label{v123:eq-wronskian-prefix}
\begin{gathered}
s_j(f)=s_j(g)\quad(1\le j\le2m-2)\\
\Longleftrightarrow\quad f'g-fg'=c\in k^\times,\\
s_{2m-1}(f)-s_{2m-1}(g)=c\quad\text{under these conditions.}
\end{gathered}
\end{equation}
Differentiating \(f'g-fg'=c\) gives \(f''g-fg''=0\). Since \(\gcd(f,g)=1\) and \(\deg f''<m\), we obtain \(f''=0\). This proves the first assertion.

Now assume \(f''=0\), and write \(f'h_f=1+fu\). Differentiation gives
\[
f'(h_f'-u)=fu'.
\]
Since \(\gcd(f,f')=1\) and \(\deg(h_f'-u)<m\), we have \(h_f'=u\), and therefore
\begin{equation}\label{v123:eq-derivative-inverse}
f'h_f-fh_f'=1.
\end{equation}
Consequently, \(f_a'f_b-f_af_b'=b-a\). For \(a\ne b\), this nonzero constant excludes common roots and repeated roots of either factor. Thus the \(f_c\) are square-free and pairwise coprime.
Reducing \(f_c'h_f-f_ch_f'=1\) modulo \(f_c\) gives \(h_{f_c}=h_f\).
Equation~\eqref{v123:eq-wronskian-prefix} shows that the \(f_c\) share the same prefix and gives \eqref{v123:eq-last-power-trace}.
Conversely, if \(g\ne f\) shares this prefix, then \(f'g-fg'=c\ne0\). Reduction modulo \(f\) gives \(g\equiv c h_f\pmod f\); comparison of degrees and leading coefficients yields \(g=f+c h_f\).
\end{proof}

Every polynomial in a nonsingleton fiber has zero second derivative, and \((f_a)_b=f_{a+b}\). Thus each nonempty fiber is either a singleton or an affine line of \(q\) polynomials; the final power sum \(s_{2m-1}\) distinguishes the points of the line.
Restricting to irreducible or primitive polynomials leaves at most \(q\) elements in each fiber, but need not preserve every element.

\paragraph{Characteristic two.}
Here \(f''=0\) always holds, and \(f'g-fg'=(fg)'\). Every fiber has \(q\) elements; for \(q=2\), this gives a pairing.

\begin{corollary}[Binary derivative pairing]\label{v123:cor-binary-derivative-pair}
On the monic square-free polynomials of degree \(m\) in \(\F_2[X]\),
\[
\mathcal P(f)=f+h_f
\]
defines a fixed-point-free involution, where \(h_f\) is the inverse representative in the theorem.
Distinct monic irreducible polynomials \(f,g\) of degree \(m\) have the same first \(2m-2\) power traces if and only if \(g=\mathcal P(f)\), equivalently, \((fg)'=1\).
In this case there is a unique monic polynomial \(H\in\F_2[X]\) of degree \(m\) such that \(fg=H^2+X\).
\end{corollary}
\begin{proof}
The first two assertions follow from the theorem with \(q=2\). In characteristic two, \((fg)'=1\) means precisely that \(fg-X\) has only even powers. Taking square roots term by term over \(\F_2\) gives the unique \(H\).
\end{proof}

The pairing need not preserve irreducibility. For example, \(f=X^3+X+1\) has no root in \(\F_2\), so it is irreducible. Since \(Xf'=f+1\),
\[
\mathcal P(f)=f+X=X^3+1=(X+1)(X^2+X+1).
\]
For which degrees does this involution pair two irreducible polynomials? How does the answer change if their roots must also generate the full multiplicative group?

\begin{example}[Two primitive quadratics over the four-element field]\label{v123:ex-f4-trace-pair}
Take \(k=\F_4=\F_2(\omega)\), where \(\omega^2+\omega+1=0\), and set
\[
f=X^2+X+\omega,\qquad g=X^2+X+\omega^2=f+1.
\]
Both polynomials are irreducible because \(\Tr_{k/\F_2}(\omega)=\Tr_{k/\F_2}(\omega^2)=1\).
If \(\alpha\) is a root of \(f\), then \(\alpha^5=\Nm_{\F_{16}/k}(\alpha)=\omega\), so its order is three or fifteen. Since \(\alpha\notin k\), its order cannot be three. The roots of \(g\) likewise have order fifteen.
Here \(f'=1\), \(h_f=1\), and
\[
fg=X^4+X+1.
\]
This is the first quartic in Example~\ref{sfr:ex-quartic-comparison}. Its roots form two orbits over \(k\), whose power traces satisfy
\[
s_1(f)=s_1(g)=s_2(f)=s_2(g)=1,\qquad
s_3(f)=\omega^2,\quad s_3(g)=\omega.
\]
The first two power traces do not distinguish the orbits, whereas the third does. The other two polynomials in the same fiber are \(X^2+X\) and \(X^2+X+1\), both of which split over \(k\).
\end{example}

\paragraph{Odd characteristic.}\label{v123:rem-odd-characteristic}
Let \(p=\Char k\) be odd. Then \(f''=0\) is equivalent to
\[
f(X)=A(X)^p+XB(X)^p\qquad(A,B\in k[X]),
\]
 since only exponents congruent to zero or one modulo \(p\) can occur. In particular, if \(m\not\equiv0,1\pmod p\), the first \(2m-2\) power sums already determine \(f\) uniquely.
Even when the degree satisfies this congruence condition, one must still check \(f''=0\) to obtain a full affine fiber. On each nonsingleton fiber, every orbit of \(f\mapsto f+h_f\) has length \(p\).
In odd characteristic, the required constant is \(f'g-fg'\); the derivative \((fg)'\) of the product cannot replace it.

\begin{example}[A cubic fiber over the three-element field]\label{v123:ex-f3-trace-pencil}
In \(\F_3[X]\), put \(f_a=X^3-X-a\), for \(a\in\F_3\). All have derivative \(-1\) and zero second derivative.
The polynomial \(f_0=X(X-1)(X+1)\) splits completely, whereas \(f_1,f_2\) have no root in \(\F_3\) and are therefore irreducible.
The relation \(\alpha^3=\alpha+a\) and the root--coefficient identities give
\[
(s_1,s_2,s_3,s_4)(f_a)=(0,2,0,2),\qquad s_5(f_a)=2a.
\]
Thus the first four power sums do not distinguish these three polynomials, while the fifth distinguishes all three.
If \(\alpha\) is a root of \(f_1\), then \(\alpha^{13}=\Nm_{\F_{27}/\F_3}(\alpha)=1\), and \(\alpha\ne1\), so its order is thirteen.
A root \(\beta\) of \(f_2\) satisfies \(\beta^{13}=2\), so its order is two or twenty-six. Since \(\beta\notin\F_3\), its order is twenty-six.
The two irreducible polynomials share a short power-trace prefix, yet only the latter is primitive.
\end{example}

\subsection{Two Linear Interpretations of the Same Data}

\begin{corollary}[Recovering subfield representations from trace prefixes]\label{str:cor-unified-recovery}\label{v70:cor-trace-prefix-generators}
Fix \(f,d\), with \(1<d<n\), \(d\mid n\), and \(m=n/d\), and regard the trace values as elements of the given quotient field \(E=k[X]/(f)\).
The first \(2m\) terms \(c_0,\ldots,c_{2m-1}\) determine
\(B_{f,d}\), \(\eta_f(d)\), and \(\mathcal W_{f,d}\), as well as the entire degree profile \(\rho_{f,d}\) and the defect \(\Delta_f(d)\).
In particular, they decide whether a polynomial or rational generator representation of degree \(m\) exists, and recover the corresponding representation when it does.
The same trace values also determine, for every \(c\mid d\), all degree profiles \(\rho_{f,c}(t)\), for \(-1\le t<n\), and the counts \(N_f(c,t)\) of representatives of degree at most \(t\) that generate \(K_c/k\), for \(0\le t<n\), where \(K_c=\F_{q^c}\). They therefore determine the minimum generator-representation degrees
\[
\kappa_f(c)=\min\{0\le t<n:N_f(c,t)>0\}
\qquad(1<c\mid d).
\]
The data here are trace values with coordinates in the quotient field, not merely their linear ranks.
\end{corollary}
\begin{proof}
First recover \(B_{f,d}\) by Corollary~\ref{str:thm-hankel-recovery}, then expand its coefficients over \(k\) to obtain \(\eta_f(d)\) and \(\mathcal W_{f,d}\).
If \(c_{2m},\ldots,c_{n-1}\) are also needed, compute them successively from~\eqref{str:eq-relative-recurrence}.
Since \(d\ge2\), we have \(2m\le n\); when \(d=2\), the original prefix already contains the required first \(n\) terms.
Then compute prefix ranks over \(k\); Theorem~\ref{str:thm-profile-trace-duality} and Corollary~\ref{str:cor-reversed-jumps} recover the profile and defect.

When \(\eta_f(d)=2\), choose a basis \(1,\beta\) of \(\mathcal V_{f,d}\); separation of coefficients gives \(B_{f,d}=A-\beta C\).
If \(1\in\mathcal W_{f,d}\), choose instead a basis \(h,1\), with \(h\) monic and \(h(0)=0\), recovering the normalized polynomial inner component.
The respective outer polynomials are the minimal polynomials of the corresponding subfield generators, and can also be computed directly by products of coefficient conjugates.

To recover the other divisor subfields, write \(K_c=\F_{q^c}\) and denote the existing sequence by
\(c_j^{(d)}=\Tr_{E/K_d}(\alpha^j/f'(\alpha))\).
For every \(c\mid d\), transitivity of the trace gives
\[
\begin{aligned}
a_j^{(c)}
&=\Tr_{K_d/K_c}(c_j^{(d)})
=\sum_{\ell=0}^{d/c-1}(c_j^{(d)})^{q^{c\ell}}\\
&=\Tr_{E/K_c}\left(\frac{\alpha^j}{f'(\alpha)}\right),
\qquad 0\le j<n.
\end{aligned}
\]
The same weight \(1/f'(\alpha)\) is retained at every level. Put
\[
r_c(s)=\dim_k\operatorname{span}_k\{a_0^{(c)},\ldots,a_s^{(c)}\},
\qquad r_c(-1)=0.
\]
For \(1<c\mid d\), trace duality gives
\[
\rho_{f,c}(t)=c-r_c(n-2-t),\qquad -1\le t<n.
\]
At \(t=-1\), the largest required index is \(n-1\); at \(t=n-1\), the formula uses \(r_c(-1)=0\). Thus both endpoints are covered by the recovered data.
For \(c=1\), directly \(\rho_{f,1}(-1)=0\) and \(\rho_{f,1}(t)=1\) for \(0\le t<n\).
Proposition~\ref{sfr:prop-exact-count} now gives
\[
N_f(c,t)=\sum_{e\mid c}\mu(c/e)q^{\rho_{f,e}(t)},
\qquad c\mid d,\quad0\le t<n.
\]
Since \(t<n\), evaluation of polynomial representatives at \(\alpha\) is injective. Hence \(N_f(c,t)>0\) means precisely that a representative of degree at most \(t\) generates \(K_c/k\).
When \(t=n-1\), every element of \(K_c\) has a representative, and a finite-field extension has a generator. For \(c>1\), the displayed minimum therefore exists and equals \(\kappa_f(c)\).
\end{proof}

For composite \(c\), knowing only the single profile \(\rho_{f,c}\) does not in general decide generation. Here the actual trace values recover the profiles of all divisor subfields through intermediate traces, allowing exact counting.
For the base field, \(N_f(1,t)=q\) for \(0\le t<n\), with no need to define \(\kappa_f(1)\).
At the endpoint of the whole extension, \(\rho_{f,n}(t)=t+1\) and \(\kappa_f(n)=1\) are already known; the corollary retains the proper intermediate-subfield hypothesis.

The procedure takes place entirely in the given quotient-ring coordinates. First use a B\'ezout identity to invert \(f'(X)\) modulo \(f\), then compute
\begin{equation}\label{str:eq-quotient-trace-computation}
 c_j=\sum_{\ell=0}^{m-1}
 \left(\frac{\alpha^j}{f'(\alpha)}\right)^{q^{d\ell}}.
\end{equation}
These values are fixed by \(\sigma^d\), hence lie in \(K\).
The Hankel equation can be solved in \(E\); its coefficients and right side lie in \(K\), so its unique solution also lies in \(K\).
To find \(k\)-linear ranks, write each element as a coefficient column in the power basis \(1,\alpha,\ldots,\alpha^{n-1}\).
There is therefore no need first to choose a separate quotient-ring model for \(K\) or to find the conjugate roots of \(f\) individually.

\medskip
The opening example leads to the following connections within this chapter:
\begin{itemize}
\item Defining equations and Frobenius fixed-point matrices recover intersections of subfields with low-degree power spaces.
\item Degree jumps, exact generator counts, and defects quantify the subfield's position in the chosen power basis.
\item A minimum representative degree of \(n/d\) is equivalent to a polynomial composition of the corresponding degrees.
\item Nonconstant root-block coefficients and coefficient rank distinguish polynomial representations from rational representations allowing denominators.
\item Base-field ranks of relative traces and shortest recurrences over the subfield recover degree profiles and root-block equations, respectively.
\end{itemize}
The same relative trace sequence gives different information over two coefficient fields: prefix ranks over the base field \(k\) recover the degree profile of the subfield;
the shortest recurrence over the subfield \(K\) recovers the relative root-block polynomial.
The former records the position of the subfield in the power basis; the latter determines the root block and its coefficient space, and hence tests low-degree generator representations and functional decompositions.

\par\Needspace{6\baselineskip}
\section*{Exercises}
\phantomsection
\addcontentsline{toc}{section}{Exercises}

\begin{exercise}[Recovering two exponential nodes]\label{v88:ex-hankel-two-nodes}
Over \(\F_5\), put \(s_r=1+2^r\) for \(r\ge0\). Compute
\(s_0,s_1,s_2,s_3\), check the order-two Hankel determinant
in~\eqref{v88:eq-hankel-vandermonde}, and recover the minimum
recurrence polynomial from these four terms.
\end{exercise}

\begin{exercise}[A quadratic subfield and the first degree jump]\label{sfr:ex-quadratic-count}
Let \(n>2\) be even, and let \(f\in\F_q[X]\) be monic and irreducible of degree \(n\). Prove that
\[
\kappa_f(2)=\lambda_1,\qquad
\Delta_f(2)=\kappa_f(2)-n/2.
\]
For \(0\le t<n\), determine \(N_f(2,t)\) when \(t<\kappa_f(2)\) and when \(t\ge\kappa_f(2)\).
\end{exercise}

\begin{exercise}[Coefficient matrices for the quartic examples]\label{sfr:ex-quartic-matrices}
For \(f_1,f_2\) in Example~\ref{sfr:ex-quartic-comparison}, write \(M_{f_i,2,3}\), and recover a basis of the four-element subfield from its kernel.
Use the first three columns to test whether a quadratic inner component exists, without first finding roots.
\end{exercise}

\begin{exercise}[Coordinates from fifth roots of unity]\label{sfr:ex-cyclotomic-four}
Prove that \(f=X^4+X^3+X^2+X+1\in\F_2[X]\) is irreducible.
For a root \(\alpha\), prove that \(\alpha^3+\alpha^2\) generates the four-element subfield, and find \(\rho_{f,2}\), \(\kappa_f(2)\), and \(\Delta_f(2)\).
\end{exercise}

\begin{exercise}[Representation defect after an optimal inner component]\label{sfr:ex-tower-eight}
In \(\F_2[X]\), put
\[
P=Y^4+Y^3+1,\qquad h=X^2+X,
\qquad f=P(h)=X^8+X^6+X^5+X^3+1.
\]
Prove that \(f\) is irreducible, and find the degree-jump sequences, minimum generator-representation degrees, and representation defects of \(K_4\) and \(K_2\) in the root coordinates of \(f\).
Explain why \(f\) has an inner component of degree two but none of degree four.
\end{exercise}

\begin{exercise}[Truncation and a nonoptimal tower bound]\label{sfr:ex-truncation-tower}
\begingroup\emergencystretch=6em
\begin{enumerate}
\item In the \(f_1\)-\allowbreak coordinates of Example~\ref{sfr:ex-quartic-comparison}, take \(\beta=\alpha^2+\alpha\).
Compare \(\deg_\alpha(\beta^2)\) with \(2\deg_\alpha(\beta)\), and explain why Proposition~\ref{sfr:prop-truncated-closure} requires truncation.
\item Suppose \(1<c\mid d\mid n\), take \(h\) attaining \(\kappa_f(d)\), put \(\beta=h(\alpha)\), and let \(P\) be the monic minimal polynomial of \(\beta\).
Prove that \(\kappa_f(c)\le\kappa_f(d)\kappa_P(c)\). Explain why this generally gives only an inequality, whereas Proposition~\ref{sfr:prop-tower-scaling} gives equality.
\end{enumerate}
\par\endgroup
\end{exercise}

\begin{exercise}[Orthogonal complements for ordinary and weighted traces]\label{str:ex-trace-weight}
With the notation of Section~\ref{str:sec-trace-duality}, write \(\perp_{\Tr}\) for the ordinary trace orthogonal complement of \(E/k\).
Prove that, for \(-1\le t<n\),
\[
\mathcal P_t(\alpha)^{\perp_{\Tr}}
=\frac{1}{f'(\alpha)}\mathcal P_{n-2-t}(\alpha).
\]
For \(f=X^4+X^3+1\in\F_2[X]\) and \(d=2\), show that replacing \(c_j\) by the unweighted trace \(\Tr_{E/K}(\alpha^j)\)
makes~\eqref{str:eq-profile-trace-duality} fail at \(t=2\).
\end{exercise}

\begin{exercise}[An invertible Hankel matrix with first trace term zero]\leavevmode\phantomsection\label{str:ex-hankel-zero-first}
Let \(f=X^4+X+1\in\F_2[X]\), with \(f(\alpha)=0\), and put \(\beta=\alpha^2+\alpha\) and \(K=\F_4\).
Compute \(c_0,c_1,c_2,c_3\) and \(\mathcal H_{f,2}\), and recover \(B_{f,2}\).
From the same data, find the trace-rank jumps, degree jumps, and polynomial inner component of least degree.
\end{exercise}

\begin{exercise}[Counting reduced rational functions with separate degree bounds]\label{v70:ex-rational-rectangle-count}
For arbitrary nonnegative integers \(a,b\), prove that the set defined in Proposition~\ref{v70:prop-rational-reconstruction} satisfies
\[
|\mathcal R_{a,b}|=q^{a+b+1}.
\]
Obtain a direct count, independent of any field of evaluation, by grouping unreduced numerator--denominator pairs according to the degree of their monic greatest common divisor.
Explain how this count agrees with the reconstruction bijection when \(a+b=n-1\).
\end{exercise}

\begin{exercise}[Counting rational representations attaining the lower bound]\label{str:ex-count-rational}
Suppose \(\eta_f(d)=2\). Prove that the number of distinct rational functions of degree \(m=n/d\) whose values at \(\alpha\) generate \(K_d/k\) is exactly
\[
q(q^2-1).
\]
If \(1\in\mathcal W_{f,d}\), exactly \(q(q-1)\) of these are polynomials; if \(1\notin\mathcal W_{f,d}\), none is a polynomial.
In counting, representations obtained by multiplying numerator and denominator by the same nonzero scalar are regarded as the same rational function.
Also prove that these functions have distinct values at \(\alpha\) when \(d>2\), whereas each generator of \(K_2/k\) has exactly \(q+1\) such rational-function representations when \(d=2\).
\end{exercise}

\begin{exercise}[Two subfields of the same sextic root field]\label{str:ex-six-two-subfields}
Let \(f=X^6+X+1\in\F_2[X]\), with \(f(\alpha)=0\), and put
\[
\beta=\alpha^5+\alpha^4+\alpha^3+\alpha.
\]
Verify \(\beta^2+\beta+1=0\), and prove that, relative to \(K_2=\F_4\),
\[
B_{f,2}(X)=X^3+X^2+\beta X+\beta+1.
\]
Find \(\kappa_f(2)\), \(\kappa_f^{\mathrm{rat}}(2)\), and \(\Delta_f(2)\), and write the corresponding rational-substitution identity.
Compare this with the representation of the same root relative to \(K_3=\F_8\) in Example~\ref{str:ex-rational-six}.
\end{exercise}

\begin{exercise}[Consecutive trace windows and arbitrary nonzero weights]\label{str:ex-shifted-hankel}
\leavevmode\par
\begin{enumerate}
\item For every integer \(s\ge0\), prove that the matrix
\[
(c_{s+i+j})_{0\le i,j<m}
\]
is invertible. Deduce that any \(2m\) consecutive trace values recover \(B_{f,d}\).
\item Replace the weight \(1/f'(\alpha)\) by any \(w\in E^\times\), and put
\[
a_j=\Tr_{E/K}(w\alpha^j).
\]
Prove that the minimal polynomial of this sequence over \(K\) is still \(B_{f,d}\), and its Hankel matrix of size \(m\) remains invertible.
Explain why the degree-profile duality formula cannot in general be retained unchanged.
\end{enumerate}
\end{exercise}

\begin{exercise}[Affine changes of generator and degree profiles]\label{v115:ex-profile-affine}
Let \(E=k(\alpha)\), with \([E:k]=n\), and put \(\beta=c\alpha+a\), where \(c\in k^\times\) and \(a\in k\). Prove, for \(0\le t<n\), that
\[
\operatorname{span}_k\{1,\alpha,\ldots,\alpha^t\}
=\operatorname{span}_k\{1,\beta,\ldots,\beta^t\}.
\]
Compare the profiles, jumps, defects, and minimum generator-representation degrees of every subfield in these two power bases.
Do the same conclusions hold for \(\beta=\alpha^{-1}\)?
\end{exercise}

\begin{exercise}[Residual degrees, returns, and multiplicities]\label{v116:ex-residual-types}
Take \(a=2\) in Corollary~\ref{v116:cor-residual-factor}. Which factorization types of \(R\) are permitted by degree divisibility alone?
If \(R\) is squarefree, determine the return cycles and fixed-point counts \((r_1,r_2)\) on the added part over one root of \(P\). If \(R\) is not squarefree, what information do these counts still retain?
\end{exercise}

\chapter{Orbit Composition, Recursive Lifting, and Trace Descent}
\label{chap:root-constructions}
Direct products and tensor products of root orbits provide composition constructions, while return maps provide recursive lifting. Trace can also be used to descend constructions from an auxiliary extension to the base field; combining these local period conditions yields irreducible polynomials of prescribed degrees.
The resultant expresses the product of conjugate values as a norm, and its symmetry under interchange yields polynomial reciprocity laws.
Polynomial iteration requires the root fibers to remain irreducible at every level; critical orbits and norm descent for power residues can reduce certain infinite conditions to finite tests in the base field.

The following table places the one-step fiber criteria of Chapter~\ref{chap:root-fibers} alongside the recursive and descent constructions of this chapter.
Let \(P\) be monic irreducible of degree \(n\), let \(\alpha\) be a root, and put \(E=k(\alpha)\).
The key to a one-step construction is that the return map has a single orbit on the entire fiber; a recursive construction must also preserve the data needed for the next step.
Trace descent additionally requires ruling out membership in a smaller subfield. In the table, power substitution assumes \(\alpha\ne0\), additive lifting assumes \(b\in k\),
and trace descent assumes that \(\ell\ne p\) is an odd prime and \(e\ge1\) is an integer.

\begingroup\small
\setlength{\tabcolsep}{4pt}
\renewcommand{\arraystretch}{1.2}
\begin{longtable}{@{}>{\raggedright\arraybackslash}p{0.13\textwidth}>{\raggedright\arraybackslash}p{0.25\textwidth}>{\raggedright\arraybackslash}p{0.28\textwidth}>{\raggedright\arraybackslash}p{0.25\textwidth}@{}}
\caption{Fibers, return conditions, and recursive data in irreducible constructions}\label{tab:construction-map}\\
\toprule
Method&Fiber polynomial or auxiliary equation&Return period or descent condition&Data preserved or lifted recursively\\
\midrule
\endfirsthead
\toprule
Method&Fiber polynomial or auxiliary equation&Return period or descent condition&Data preserved or lifted recursively\\
\midrule
\endhead
\bottomrule
\endfoot
Power substitution&\(X^t-\alpha\)&Irreducible over \(E\); tested using the prime valuations of \(\ord(\alpha)\) and the binomial criterion&The fiber period is \(t\), and the root order increases from \(\ord(\alpha)\) to \(t\ord(\alpha)\); see Corollary~\ref{cor:power-composition}\\
\addlinespace
Additive lifting&\(X^p-X-(\alpha+b)\)&\(\Tr_{E/\F_p}(\alpha+b)\ne0\)&The fiber is a single translation orbit of length \(p\); nonzero trace at one step does not automatically ensure the next step\\
\addlinespace
Reciprocal additive recursion&\(F\mapsto\widehat F(X^p-X+\delta)\)&The next-to-leading coefficient after reciprocation has nonzero absolute trace&The values on \(\F_p\) are a common \(A\in\F_p^\times\), and the derivative values are a common \(D\), with \(\Tr_{k/\F_p}(D)\ne0\); see Theorem~\ref{thm:recursive-as}\\
\addlinespace
Q-transform (characteristic two)&\(X^2-\alpha X+1\)&\(\Tr_{E/\F_2}(\alpha^{-1})=1\)&From the first step onward, the constant term is one, and the coefficients of the linear and next-to-leading terms are equal and remain unchanged; see Theorem~\ref{thm:q-transform-iteration} for the initial two-trace condition\\
\addlinespace
R-transform (odd characteristic)&\(X^2-2\alpha X+1\)&\(P(1)P(-1)\) is a nonsquare&This square class is preserved under the stated parity condition on the degree; see Theorem~\ref{thm:r-transform-iteration}\\
\addlinespace
Iteration of a fixed polynomial&\(h(X)-\beta_n\), \(h(\beta_{n+1})=\beta_n\)&Every fiber is irreducible, with return period \(\deg h\)&For quadratic iteration, square classes are controlled by norms along the critical orbit; a cubic family in characteristic two preserves noncubicity; see Section~\ref{v55:sec-dynamical-irreducibility}\\
\addlinespace
Auxiliary field and trace descent&\(\beta^{\ell^e}=a\), with the base field first extended to \(\F_{q^s}\)&\(s=\ord_\ell(q)\), \(a\) is not an \(\ell\)-th power; the exponents in the trace sum are distinct modulo \(\ell\)&A drop in degree after descent is ruled out, retaining the required degree \(\ell^e\); see Theorem~\ref{thm:trace-descent-construction}\\
\end{longtable}
\endgroup

In a recursive tower, suppose that a map defined over the base field sends the new root \(\beta_{j+1}\) to the old root \(\beta_j\),
and put \(n_j=[k(\beta_j):k]\). The return map on the fiber containing the new root is then \(\sigma^{n_j}\).
If its period at the new root is \(r_j\), Lemma~\ref{lem:fiber-return} gives
\[
n_{j+1}=n_jr_j.
\]
Thus degree growth in the recursion results from multiplying the internal return periods at successive levels.
Inductive conditions on traces, function values, and derivative values ensure that the return map at the next level still traverses every element of its fiber,
rather than merely ensuring the irreducibility of the polynomial already constructed.

\section{Composed Sums and Products of Roots}
\leavevmode\phantomsection\label{def:composed-products}\ignorespaces For monic polynomials \(f,g\in k[X]\), define their \textbf{composed sum} and \textbf{composed product} in a common finite splitting field by
\[
(f\oplus g)(X)=\prod_{f(\alpha)=0}\prod_{g(\beta)=0}(X-\alpha-\beta),
\quad
(f\odot g)(X)=\prod_{f(\alpha)=0}\prod_{g(\beta)=0}(X-\alpha\beta).
\]
All roots are counted with multiplicity. Both polynomials have degree \((\deg f)(\deg g)\).
More generally, if the set \(G\) containing the roots has an operation \(\diamond\) satisfying \(\sigma(\alpha\diamond\beta)=\sigma(\alpha)\diamond\sigma(\beta)\), define
\[
(f\diamond g)(X)=\prod_{\alpha,\beta}(X-\alpha\diamond\beta).
\]
Frobenius permutes these roots together with their multiplicities, so the resulting coefficients belong to \(k\).
The evaluation points in Figure~\ref{fig:tensor-evaluation} are precisely such pairs of roots, with Frobenius acting on both coordinates simultaneously.
The operation on roots then sends each pair to its sum or product; distinct evaluation points may have the same image, so the resulting multiplicities must be retained.
Moreover, \(f\diamond(g_1g_2)=(f\diamond g_1)(f\diamond g_2)\), and the same holds in the other variable.

\begin{theorem}[Composition of coprime orbits]\leavevmode\phantomsection\label{thm:composed-products}
Suppose that \((G,\diamond)\) is a Frobenius-stable group and that Frobenius preserves its group operation.
If the roots of monic polynomials \(f,g\) all lie in \(G\), and their degrees are \(m,n\ge1\), respectively, then \(f\diamond g\) is irreducible
if and only if both \(f,g\) are irreducible and \(\gcd(m,n)=1\).
\end{theorem}
\begin{proof}
Distributivity of the composed-product operation over ordinary products shows that irreducibility of the result forces both \(f,g\) to be irreducible.
If \(d=\gcd(m,n)>1\), every root pair \((\alpha,\beta)\) has period at most \(\lcm(m,n)=mn/d\),
and its composed value has period at most this number, so it cannot be a root of an irreducible polynomial of degree \(mn\).

Conversely, suppose that \(f\) and \(g\) are irreducible and that \(m,n\) are coprime. Let \(\gamma=\alpha\diamond\beta\) have period \(t\).
Since \(\sigma^{mn}\) fixes both \(\alpha\) and \(\beta\), it fixes \(\gamma\); hence its orbit length \(t\) divides \(mn\). Since \(\sigma^{mt}(\gamma)=\gamma\) and \(\sigma^{mt}(\alpha)=\alpha\),
cancellation gives \(\sigma^{mt}(\beta)=\beta\). Hence \(n\mid mt\), and therefore \(n\mid t\).
Similarly, \(m\mid t\), so \(mn\mid t\). Thus \(t=mn\), and the degree of the composed value over the base field equals the degree of the composed polynomial.
\end{proof}
The additive group allows zero roots; the multiplicative group requires \(f(0)g(0)\ne0\).
For example, \(X\oplus g=g\), whereas \(X\odot g=X^{\deg g}\), so the nonzero condition cannot be omitted in the multiplicative case.
Any binary polynomial operation compatible with Frobenius defines a composition, but need not satisfy the irreducibility statement above: the proof uses cancellation in a group.

\begin{example}\label{ex:composed-six}
Take \(f=X^2+X+1\) and \(g=X^3+X+1\in\F_2[X]\), and let \(\alpha^2+\alpha+1=0\).
For a fixed root \(\alpha\) of \(f\), the polynomial of the product roots is \(\alpha^3g(X/\alpha)=X^3+\alpha^2X+1\), and that of the sum roots is \(g(X+\alpha)\).
Multiplying by the respective conjugate polynomials gives
\[
\begin{aligned}
f\odot g&=(X^3+\alpha^2X+1)(X^3+\alpha X+1)
         =X^6+X^4+X^2+X+1,\\
f\oplus g&=g(X+\alpha)g(X+\alpha^2)
         =X^6+X^5+X^3+X^2+1.
\end{aligned}
\]
Both are irreducible, since the original degrees two and three are coprime. The roots of the first have multiplicative order \(3\cdot7=21\),
whereas those of the second have order \(63\) (as can be checked by the primitivity criterion); thus only the second is primitive.
If \(\gcd(s,t)=1\) and \(p\nmid st\), we also have \(\Phi_s\odot\Phi_t=\Phi_{st}\): products of roots of unity of coprime orders run through precisely all primitive \(st\)-th roots of unity.
\end{example}

\subsection{Synchronized Phases and Factorization of Composed Products}\label{v40:sec-composed-phases}
After choosing initial points \(\alpha,\beta\) in two root orbits, a root pair
\((\alpha^{q^i},\beta^{q^j})\) is recorded by its two iteration indices. We call \(i,j\) its
\textbf{phases}: the numbers of Frobenius iterations from the chosen initial points, taken modulo the respective orbit lengths.
\leavevmode\phantomsection\label{v40:def-composed-phase}
Frobenius increases both phases by one. The phase quotient first separates root-pair orbits; the composed values then determine whether these orbits shorten or merge.

\begin{proposition}[Synchronized phases, shorter periods, and collisions]\label{v129:prop-phase-images}
Let \(f,g\in k[X]\) be monic irreducible polynomials of degrees \(d,e\), with roots \(\alpha,\beta\), and put
\(h=\gcd(d,e)\) and \(L=\lcm(d,e)\). For \(\Psi\in k[U,V]\), define
\[
C_\Psi(X)=\prod_{f(u)=0}\prod_{g(v)=0}(X-\Psi(u,v)).
\]
For \(0\le r<h\), put \(\gamma_r=\Psi(\alpha,\beta^{q^r})\) and \(s_r=[k(\gamma_r):k]\). Then \(s_r\mid L\), and
\begin{equation}\label{v129:eq-phase-images}
C_\Psi(X)=\prod_{r=0}^{h-1}m_{\gamma_r}(X)^{L/s_r}.
\end{equation}
In particular, \(C_\Psi\) is squarefree if and only if every \(s_r=L\) and the polynomials \(m_{\gamma_r}\) are pairwise distinct; equivalently, \(\Psi\) is injective on all root pairs.
\end{proposition}
\begin{proof}
The phases of root pairs lie in \(\Z/d\Z\times\Z/e\Z\), and Frobenius adds one to both coordinates.
Every orbit has length \(L\), and the phase difference \(j-i\bmod h\) separates the \(h\) orbits, with representatives \((0,r)\), \(0\le r<h\).
Indeed, two phase pairs lie in the same orbit exactly when the two congruences for a common connecting step are compatible, that is, when their phase differences agree modulo \(h\).
Since \(\Psi\) commutes with Frobenius, the image of the \(r\)-th root-pair orbit is the orbit of \(\gamma_r\), of length \(s_r\mid L\), and each image value occurs \(L/s_r\) times.
Multiplying over the orbits gives the formula. Squarefreeness requires the absence of collisions both within and between orbits.
\end{proof}

Factor exponents have two sources: \(L/s_r\) records the loss of period within one root-pair orbit, whereas equal polynomials \(m_{\gamma_r}\) record the merging of distinct phase orbits.
These are the two kinds of collision in~\eqref{v129:eq-orbit-image-product}, now with a product of root sets as the original set.
For the factorization type, these periods and conjugacy relations suffice; explicit coefficients can then be obtained from the determinants or resultants in the next section.

\begin{example}[Two compositions of the same cubic]\label{v129:ex-composed-cubic}
Take \(f=X^3+X+1\in\F_2[X]\), with root \(\alpha\). It has no root in the base field and is therefore irreducible; since \(\alpha\ne0,1\), its multiplicative order is seven.
Here \(h=L=3\). For the composed sum, the three phase representatives give
\[
0,\quad \alpha+\alpha^2=\alpha^4,\quad
\alpha+\alpha^4=\alpha^2.
\]
Their periods are one, three, and three, and the last two values are conjugate. Thus
\[
f\oplus f=X^3f(X)^2.
\]
For the composed product, the representatives give \(\alpha^2,\alpha^3,\alpha^5\). Doubling the exponents gives the orbits
\(\{1,2,4\}\) and \(\{3,6,5\}\); the latter is the negative of the former, so its root polynomial is \(f^\vee=X^3+X^2+1\). Hence
\[
f\odot f=f(X)\bigl(f^\vee(X)\bigr)^2.
\]
No individual phase orbit shortens in the product case, but two phase orbits have the same image, still producing a squared factor.
\end{example}

\begin{corollary}[Composed products with coprime root orders]\leavevmode\phantomsection\label{v40:prop-coprime-root-orders}
Let \(f,g\in\F_q[X]\) be monic irreducible polynomials of degrees \(d,e\), respectively, and choose nonzero roots \(\alpha,\beta\) in a common extension.
Suppose that the multiplicative orders of \(\alpha,\beta\) are coprime, and put
\[
 h=\gcd(d,e),\qquad L=\lcm(d,e).
\]
Then
\[
 f\odot g=\prod_{r=0}^{h-1}F_r,
 \qquad
 F_r(X)=\prod_{j=0}^{L-1}
 \bigl(X-\alpha^{q^j}\beta^{q^{r+j}}\bigr).
\]
The polynomials \(F_r\) all belong to \(\F_q[X]\), are pairwise distinct and irreducible, and each has degree \(L\).
In particular, \(f\odot g\) has no repeated roots.
\end{corollary}
\begin{proof}
Put \(u=\ord(\alpha)\) and \(v=\ord(\beta)\). Equality of products of two root pairs gives
\[
\alpha^{q^i}\alpha^{-q^{i'}}=\beta^{q^{j'}}\beta^{-q^j}
\in\mu_u\cap\mu_v=\{1\},
\]
so the two coordinates are respectively equal. Multiplication is injective on the root-pair set, and Proposition~\ref{v129:prop-phase-images} gives distinct irreducible phase factors, all of degree \(L\).
\end{proof}

\begin{example}[Noncoprime degrees with distinct product roots]\leavevmode\phantomsection\label{v40:ex-phase-fifteen}
In \(\F_2[X]\), take
\[
 f=\Phi_3=X^2+X+1,\qquad
 g=\Phi_5=X^4+X^3+X^2+X+1.
\]
Both are irreducible. Their roots have multiplicative orders three and five, respectively, but their degrees are two and four.
Coprimality of the multiplicative orders ensures that the eight product roots are distinct; the phase difference modulo two splits them into two Frobenius orbits of size four.
Consequently,
\[
 f\odot g=\Phi_{15}
 =(X^4+X+1)(X^4+X^3+1),
\]
and both factors are irreducible. The two phase classes here correspond to the two primitive-root orbits in the example of fifteenth roots of unity in Chapter~\ref{chap:cyclotomic}.
Coprime degrees ensure that the product of the root sets is a single Frobenius orbit; coprime multiplicative orders ensure that different root pairs have different products.
These two conditions address different questions.
\end{example}

\par\Needspace{12\baselineskip}
\section{Matrix Representations of Composed Constructions}
\begin{proposition}\label{prop:composed-matrices}
Let \(A,B\) be the companion matrices of monic polynomials \(f,g\), of sizes \(m,n\), respectively. Then
\begin{align*}
f\odot g&=\det(XI_{mn}-A\otimes B),\\
f\oplus g&=\det\bigl(XI_{mn}-(A\otimes I_n+I_m\otimes B)\bigr).
\end{align*}
Here \(A\otimes B=(a_{ij}B)\) is the \textbf{Kronecker product}, and \(A\otimes I_n+I_m\otimes B\) is the \textbf{Kronecker sum}.
More generally, the polynomial whose roots are obtained by applying \(\Psi(U,V)=\sum c_{ij}U^iV^j\) is
\[
\det\left(XI_{mn}-\sum_{i,j}c_{ij}A^i\otimes B^j\right).
\]
\end{proposition}
\begin{proof}
Over a finite splitting field, put \(A,B\) separately into upper triangular form: choose an eigenvector and then argue inductively on the quotient space.
Their diagonal entries are the roots of \(f,g\), respectively, listed with algebraic multiplicity.
The corresponding Kronecker matrices are also upper triangular, with diagonal entries \(\alpha_i\beta_j\), \(\alpha_i+\beta_j\), or \(\Psi(\alpha_i,\beta_j)\), respectively.
Their determinants are the required root products. This proof also applies to polynomials with repeated roots and does not assume diagonalizability.
\end{proof}

These composed polynomials also admit determinant expressions of smaller size. If \(g(X)=\sum_{j=0}^n b_jX^j\), then
\begin{equation}\label{eq:small-composed-determinants}
(f\oplus g)(X)=\det\bigl(g(XI_m-A)\bigr),\qquad
(f\odot g)(X)=\det\left(\sum_{j=0}^n b_jX^jA^{n-j}\right).
\end{equation}
The proof again uses triangularization: the diagonal entries are \(g(X-\alpha_i)\) and \(\alpha_i^ng(X/\alpha_i)\), respectively,
where the latter expression is interpreted as the expanded polynomial, so zero roots cause no difficulty.
By interchanging \(f,g\), the determinant size can be chosen as \(\min(m,n)\).
\phantomsection\label{sec:composed-resultants}
In the \hyperref[def:resultant-root-product]{root-product definition of the resultant}, roots are counted with multiplicity, and the coefficients of the second polynomial may contain an independent parameter \(X\).
The product is a symmetric polynomial in the roots and therefore belongs to the corresponding coefficient ring over the base field.
In~\eqref{eq:small-composed-determinants}, this gives
\[
\begin{aligned}
f\oplus g&=\operatorname{Res}_Y\bigl(f(Y),g(X-Y)\bigr),\\
f\odot g&=\operatorname{Res}_Y\left(f(Y),\sum_{j=0}^n b_jX^jY^{n-j}\right).
\end{aligned}
\]
The second expression is written as an ordinary polynomial and includes the case where \(f\) has a zero root.

For irreducible \(f,g\), these matrices also have a common algebraic interpretation.
In \(k[X]/(f)\otimes_k k[Y]/(g)\), denote the two coordinate elements by \(\alpha,\beta\).
The matrix of multiplication by \(\alpha\otimes1+1\otimes\beta\) is the Kronecker sum of the corresponding multiplication matrices,
and that of multiplication by \(\alpha\otimes\beta\) is their Kronecker product.
Thus the composed-sum and composed-product polynomials are the characteristic polynomials of these two multiplication operators.
Theorem~\ref{thm:tensor-fields} shows that, when the degrees are coprime, this tensor product is itself the compositum field.
Theorem~\ref{thm:composed-products} further ensures that \(\alpha\otimes1+1\otimes\beta\) generates this field;
if \(f(0)g(0)\ne0\), then \(\alpha\otimes\beta\) also generates it.
When the degrees are not coprime, several field components occur, and different root pairs may give the same sum or product, so multiplicity information must be retained.

\section{Power-Residue Symbols, Norms, and Polynomial Reciprocity}
\label{sec:polynomial-reciprocity}
Throughout this section, \(k=\F_q\). The resultant multiplies the values of one polynomial at all roots of another.
When these roots form a single Frobenius orbit, this product is a norm.
Taking a suitable power of the norm then determines whether an element of the residue field is a power of the prescribed exponent.
Interchanging the two root sets changes the resultant only by a sign, yielding a polynomial power-residue reciprocity law~\cite[Chapter 9]{Agou}.

\begin{definition}[Power-residue symbol]\leavevmode\phantomsection\label{def:power-residue-symbol}
Let \(d\mid q-1\), \(d\ge1\), and let \(P\in k[X]\) be monic irreducible, with \(\deg P=n\).
For \(g\in k[X]\) with \(P\nmid g\), choose a root \(\theta\) of \(P\) and define
\[
\left(\frac{g}{P}\right)_d=g(\theta)^{(q^n-1)/d}.
\]
This element is called the \textbf{\(d\)-th power-residue symbol} of \(g\) modulo \(P\).
\end{definition}
The symbol takes values in \(\mu_d=\{u\in k^\times:u^d=1\}\) and does not depend on the chosen root.
Indeed, its \(d\)-th power is one, and \(d\mid q-1\) ensures that all \(d\)-th roots of unity lie in \(k\).
Replacing \(\theta\) by another root \(\theta^{q^j}\) simply raises the value of the symbol to the \(q^j\)-th power, leaving it unchanged.
It also depends only on \(g\bmod P\) and is multiplicative in \(g\).

In the residue field \(E=k[X]/(P)\simeq\F_{q^n}\), write \(\overline g=\zeta^a\),
where \(\zeta\) generates \(E^\times\). Then
\[
\left(\frac{g}{P}\right)_d=1
\quad\Longleftrightarrow\quad d\mid a
\quad\Longleftrightarrow\quad
\overline g\in(E^\times)^d.
\]
Thus the power-residue symbol reduces the test for a \(d\)-th power in the residue field to a value in the group of roots of unity in the base field.

\begin{proposition}[Norm expression for the power-residue symbol]\label{prop:residue-norm-resultant}
Under the hypotheses of Definition~\ref{def:power-residue-symbol}, put \(E=k(\theta)\). Then
\begin{equation}\label{eq:residue-norm-resultant}
\left(\frac{g}{P}\right)_d
=\Nm_{E/k}\bigl(g(\theta)\bigr)^{(q-1)/d}
=\operatorname{Res}(P,g)^{(q-1)/d}.
\end{equation}
\end{proposition}
\begin{proof}
Since the coefficients of \(g\) lie in \(k\), we have \(g(\theta)^{q^j}=g(\theta^{q^j})\). The conjugate-product formula for the norm gives
\[
\Nm_{E/k}\bigl(g(\theta)\bigr)
=\prod_{j=0}^{n-1}g(\theta^{q^j})
=g(\theta)^{1+q+\cdots+q^{n-1}}
=\operatorname{Res}(P,g).
\]
Raising this equality to the \((q-1)/d\)-th power and using
\((1+q+\cdots+q^{n-1})(q-1)=q^n-1\) proves the claim.
\end{proof}

The connection between norms and resultants also retains repeated factors. Applying Proposition~\ref{v117:prop-local-dimension-trace} to a monogenic quotient algebra gives the following formula.

\begin{corollary}[Local multiplicities and determinants of multiplication operators]\label{prop:local-multiplicity-determinant}
Let \(F=\prod_{i=1}^rP_i^{e_i}\in k[X]\) be monic and nonconstant, where the \(P_i\) are pairwise distinct monic irreducible polynomials.
Put \(A=k[X]/(F)\) and \(K_i=k[X]/(P_i)\). For any \(g\in k[X]\),
write \(u=g\bmod F\) and \(\bar u_i=g\bmod P_i\), and let \(M_u:A\to A\) be the \(k\)-linear operator of multiplication by \(u\). Then
\begin{equation}\label{eq:local-multiplicity-determinant}
\boxed{\quad
\det_k M_u=\operatorname{Res}(F,g)
=\prod_{i=1}^r\Nm_{K_i/k}(\bar u_i)^{e_i}.
\quad}
\end{equation}
Here \(u\) need not be invertible; if some \(P_i\mid g\), both sides are zero.
\end{corollary}
\begin{proof}
The component \(A_i=k[X]/(P_i^{e_i})\) has residue field \(K_i\) and \(\ell_i=e_i\). Equation~\eqref{v121:eq-local-determinant} gives
\[
\det_k M_u=\prod_i\Nm_{K_i/k}(\bar u_i)^{e_i}.
\]
On the other hand, the roots of \(P_i\) form a single Frobenius orbit, so the root product with multiplicities gives
\[
\operatorname{Res}(F,g)
=\prod_i\operatorname{Res}(P_i,g)^{e_i}
=\prod_i\Nm_{K_i/k}(\bar u_i)^{e_i}.
\]
This proves the required equality.
\end{proof}

The same multiplication operator occurs on every layer of the filtration: taking traces gives the weighted sum
\(\sum_i e_i\Tr_{K_i/k}(\bar u_i)\), while taking determinants gives the product of norms with multiplicities.
In the trace formula, the integers \(e_i\) are interpreted as scalars in the base field; in the determinant formula, they remain integer exponents.

For any monic nonconstant polynomial \(F=\prod_iP_i^{e_i}\), with the \(P_i\) monic irreducible and pairwise distinct,
and any \(g\in k[X]\) satisfying \(\gcd(F,g)=1\), set
\[
\left(\frac{g}{F}\right)_d
=\prod_i\left(\frac{g}{P_i}\right)_d^{e_i}.
\]
This extends the symbol multiplicatively in its modulus. Listing the roots of \(F\) with multiplicity, formula~\eqref{eq:residue-norm-resultant} gives
\begin{equation}\label{eq:residue-composite-resultant}
\left(\frac{g}{F}\right)_d=\operatorname{Res}(F,g)^{(q-1)/d}.
\end{equation}
Repeated factors of \(F\) are allowed here; the formula retains the multiplicity of each irreducible factor.
When \(F\) is reducible, however, \((g/F)_d=1\) does not in general imply that \(g\bmod F\) is a \(d\)-th power in the quotient ring:
the symbol values for different factors may cancel in the product, even if \(F\) is square-free.
Thus this multiplicative extension retains the product of the component criteria, rather than each criterion separately; see Exercise~\ref{ex:power-residue-boundaries}.

\begin{theorem}[Polynomial reciprocity]\label{thm:polynomial-reciprocity}
Let \(d\mid q-1\) be a positive integer, and let \(F,G\in k[X]\) be coprime monic nonconstant polynomials,
with \(\deg F=m\) and \(\deg G=n\). Then
\begin{equation}\label{eq:polynomial-reciprocity}
\boxed{
\left(\frac{G}{F}\right)_d
=(-1)^{mn(q-1)/d}\left(\frac{F}{G}\right)_d.
}
\end{equation}
In particular, when \(q\) is odd, the quadratic residue symbol satisfies
\[
\left(\frac{G}{F}\right)_2
=(-1)^{mn(q-1)/2}\left(\frac{F}{G}\right)_2.
\]
\end{theorem}
\begin{proof}
In a common splitting field, write, with multiplicities,
\(F=\prod_{i=1}^m(X-\alpha_i)\) and \(G=\prod_{j=1}^n(X-\beta_j)\).
Since both polynomials are monic,
\[
\operatorname{Res}(F,G)=\prod_{i,j}(\alpha_i-\beta_j)
=(-1)^{mn}\prod_{j,i}(\beta_j-\alpha_i)
=(-1)^{mn}\operatorname{Res}(G,F).
\]
Raising both sides to the \((q-1)/d\)-th power and applying~\eqref{eq:residue-composite-resultant} proves the result.
\end{proof}

Reciprocity interchanges the numerator and modulus of the power-residue symbol; in the root-product expression, it reverses the order of the two terms in each root difference.
This is a different operation from the change of roots \(\alpha\mapsto\alpha^{-1}\) in a reciprocal polynomial.
The condition \(d\mid q-1\) places the power-residue symbols in all extension fields in the same group of roots of unity in the base field,
so that norm descent and interchange in the resultant can be connected directly.

\begin{example}[Detecting squares in an extension from base-field values]\label{exm:polynomial-reciprocity-cubic}
In \(\F_3[X]\), take \(F=X^3-X-1\) and \(G=X\).
For every \(a\in\F_3\), we have \(F(a)=-1\), so the cubic \(F\) is irreducible.
Direct computation gives \((F/X)_2=F(0)=-1\). By reciprocity,
\[
\left(\frac{X}{F}\right)_2=(-1)^{3}\left(\frac{F}{X}\right)_2=1.
\]
Thus a root \(\theta\in\F_{27}\) of \(F\) is a square.
The norm formula gives the same conclusion: \(\Nm_{\F_{27}/\F_3}(\theta)=-F(0)=1\),
that is, \(\theta^{13}=1\), whence \(\theta=(\theta^7)^2\).
\end{example}

\section{Recursive Additive Lifting Preserving Trace Conditions}
\label{sec:recursive-as}
For a polynomial \(F\) of degree \(N\) with nonzero constant term, write \(\widehat F(X)=X^NF(X^{-1})\) in this section for its reciprocal polynomial, \emph{without monic normalization}.
When \(F\) is monic, this is related to the monic reciprocal polynomial of Section~\ref{sec:reciprocal-polynomials} by \(\widehat F=F(0)F^\vee\).

\begin{example}[One irreducible step and the second return]\label{v115:ex-as-failure}
Take \(k=\F_3\) and \(h=X^3-X-1\). It has no root in \(k\), so it is irreducible. If \(h(\alpha)=0\), then \(E=k(\alpha)=\F_{27}\) and \(\Tr_{E/k}(\alpha)=0\).
The fiber for the next composition is
\[
h(X)-\alpha=X^3-X-(\alpha+1).
\]
Since \(\Tr_{E/k}(\alpha+1)=0+3=0\), Hilbert 90 gives a root in \(E\). Adding the three elements of \(\F_3\) gives all roots, and the derivative is \(-1\).
The return \(\sigma^3\) is therefore the identity on this fiber, and \(h\circ h\) has three distinct irreducible cubic factors.
The first irreducible step uses \(\Tr_{k/\F_3}(1)=1\), but the trace needed for the second step has vanished.
\end{example}

Iteration thus requires a quantity whose nonvanishing survives the change of variable. The next construction first takes a reciprocal, turning the linear coefficient into the next-to-leading coefficient, then controls its trace through identities for values and derivatives.

\begin{theorem}\label{thm:recursive-as}
Let \(P=X^n+c_{n-1}X^{n-1}+\cdots+c_0\in k[X]\) be monic irreducible.
Suppose that there exists \(\delta_0\in\F_p\) such that
\[
P(\delta_0)\in\F_p^\times,\quad
\Tr_{k/\F_p}(n\delta_0+c_{n-1})\ne0,\quad
\Tr_{k/\F_p}(P'(\delta_0))\ne0.
\]
Define
\[
F_0=P(X^p-X+\delta_0),\qquad
F_{j+1}=\widehat F_j(X^p-X+\delta_{j+1}),\quad \delta_{j+1}\in\F_p^\times.
\]
Then every \(F_j\) is irreducible, of degree \(np^{j+1}\).
\end{theorem}
\begin{proof}
First let \(\alpha\) be a root of \(P\).
The fiber polynomial for \(F_0\) is \(X^p-X-(\alpha-\delta_0)\), whose irreducibility condition is
\[
\Tr_{k(\alpha)/\F_p}(\alpha-\delta_0)
=-\Tr_{k/\F_p}(c_{n-1}+n\delta_0)\ne0.
\]
Thus \(F_0\) is irreducible of degree \(np\).
For any \(c\in\F_p\), we have \(c^p-c=0\), so
\[
F_0(c)=A_0=P(\delta_0)\in\F_p^\times,\qquad
F_0'(c)=D_0=-P'(\delta_0).
\]
By hypothesis, \(\Tr_{k/\F_p}(D_0)\ne0\).

Assume inductively that \(F_j\) is irreducible of degree \(N=np^{j+1}\), and that at every element of \(\F_p\),
it takes the same value \(A\in\F_p^\times\), and its derivative takes the same value \(D\in k\), where \(\Tr_{k/\F_p}(D)\ne0\).
In particular, \(F_j(0)=A\) and \(F_j'(0)=D\), so its constant and linear coefficients are \(A,D\), respectively.
Reciprocation preserves irreducibility, and the leading and next-to-leading coefficients of \(\widehat F_j\) are \(A,D\), respectively.
After monic normalization, its next-to-leading coefficient is \(D/A\).
Since \(A\in\F_p^\times\),
\(\Tr_{k/\F_p}(D/A)=A^{-1}\Tr_{k/\F_p}(D)\ne0\).
Also, \(p\mid N\), so the additive fiber condition for the composition \(\widehat F_j(X^p-X+\delta)\) is independent of \(\delta\) and is always ensured by this nonzero trace.
Taking \(\delta=\delta_{j+1}\in\F_p^\times\), we obtain an irreducible \(F_{j+1}\) of degree \(Np\).

It remains to verify that the constant-value conditions and the nonzero-trace condition persist at the next step.
Put \(Y=X^p-X+\delta\). At \(c\in\F_p\), we have \(Y(c)=\delta\) and \(Y'(c)=-1\), while
\[
\widehat F_j(Y)=Y^NF_j(Y^{-1}).
\]
Since \(N=0\) in \(k\), differentiation with respect to \(Y\) gives
\[
\widehat F_j'(Y)=-Y^{N-2}F_j'(Y^{-1}).
\]
Since \(\delta^{-1}\in\F_p\), we have
\[
F_{j+1}(c)=\delta^N A,\qquad
F_{j+1}'(c)=\delta^{N-2}D.
\]
The former lies in \(\F_p^\times\), and the absolute trace of the latter is \(\delta^{N-2}\Tr_{k/\F_p}(D)\ne0\).
Thus all inductive conditions persist, proving the result.
\end{proof}

\leavevmode\phantomsection\label{v98:par-recursive-conditions}\ignorespaces
In this construction, each additive fiber is an affine line over \(\F_p\); a nonzero absolute trace makes the return map a nonidentity translation, with least period exactly \(p\).
The initial parameters must ensure both irreducibility of the first layer and the possibility of lifting again after reciprocation.
The preserved conditions \(A\in\F_p^\times\) and \(\Tr_{k/\F_p}(D)\ne0\) ensure that the next-to-leading coefficient \(D/A\) continues to have nonzero trace; the values and derivatives at the next layer are multiplied only by nonzero prime-field scalars.
The recursion therefore depends on preservation of these conditions, not merely on irreducibility at the first step.

\begin{corollary}\label{cor:recursive-prime-field}
Over the prime field, if \(a\in\F_p^\times\) and \((na+c_{n-1})P'(a)\ne0\),
then starting with \(F_0=P(X^p-X+a)\) and repeatedly taking the reciprocal followed by substitution of \(X^p-X+a\) gives irreducible polynomials of degrees \(np^{j+1}\).
In particular, \(P=X,a=-1\) gives irreducible polynomials of degrees \(p^{j+1}\) (\(j\ge0\)).
\end{corollary}
\begin{proof}
When \(n\ge2\), irreducibility ensures \(P(a)\ne0\); when \(n=1\), this value is \(na+c_{n-1}\), which is also nonzero.
The trace from the prime field to itself is the identity map, so all conditions of the preceding theorem hold.
\end{proof}

\section{Recursive Quadratic Lifting: Self-Reciprocal Orbits}\label{v48:sec-recursive-quadratic}
\begin{theorem}[The Q-transform in characteristic two]\label{thm:q-transform-iteration}
Let \(k\) have characteristic two, and let \(P=X^n+c_{n-1}X^{n-1}+\cdots+c_0\in k[X]\) be monic irreducible with \(c_0\ne0\). Retain the convention \(c_1=1\) when \(n=1\), and suppose that
\[
\Tr_{k/\F_2}(c_1/c_0)=1,\qquad \Tr_{k/\F_2}(c_{n-1})=1.
\]
Put \(P_0=P\) and \(P_{j+1}=\mathcal Q(P_j)\). Then each \(P_j\) is irreducible of degree \(2^jn\), and is self-reciprocal for \(j\ge1\).
\end{theorem}
\begin{proof}
By Theorem~\ref{thm:reciprocal-transform}, the first trace condition ensures that \(P_1\) is irreducible.
For monic \(F\) of degree \(N\), the constant term of \(\mathcal Q(F)\) is one, and its linear coefficient is the next-to-leading coefficient of \(F\).
Since \(\mathcal Q(F)\) is self-reciprocal, its linear and next-to-leading coefficients are equal. Thus, from \(P_1\) onward, both coefficients remain equal to \(c_{n-1}\).
The second trace condition then ensures irreducibility of the quadratic fiber polynomial at every step.
\end{proof}

\begin{example}\label{v153:ex-qtransform-seventeen}
Starting from \(P_0=X+1\in\F_2[X]\), we obtain successively
\[
\begin{aligned}
P_1&=X^2+X+1,\\
P_2&=X^4+X^3+X^2+X+1,\\
P_3&=X^8+X^7+X^6+X^4+X^2+X+1.
\end{aligned}
\]
At each step, Frobenius acts transitively on the union of the two-element fibers over the preceding orbit, producing an orbit of twice its length.
\end{example}

\begin{theorem}[The R-transform in odd characteristic]\label{thm:r-transform-iteration}
Let \(q\) be odd, and let \(P\in k[X]\) be monic irreducible of degree \(n\), with \(P(1)P(-1)\) a nonsquare.
If \(q\equiv3\pmod4\), assume in addition that \(n\) is even.
Put \(P_0=P\) and \(P_{j+1}=\mathcal R(P_j)\). Then every \(P_j\) is irreducible of degree \(2^jn\), and is monic and self-reciprocal for \(j\ge1\).
\end{theorem}
\begin{proof}
If \(F\) has degree \(N\), then
\[
\mathcal R(F)(1)\mathcal R(F)(-1)=(-1)^N2^{2N}F(1)F(-1).
\]
When \(-1\) is a square or \(N\) is even, the square class is unchanged. The hypotheses ensure that the first transform preserves the square class; all subsequent degrees are even, so every later transform also preserves it.
Successive application of the \hyperref[def:r-transform]{R-transform criterion} in Chapter~\ref{chap:root-fibers} proves irreducibility at every iteration.
For a monic polynomial \(F\) of degree \(N\), the leading and constant coefficients of \(\mathcal R(F)\) are both one, and
\[
X^{2N}\mathcal R(F)(X^{-1})=\mathcal R(F)(X).
\]
Hence all polynomials from the first step onward are monic and self-reciprocal.

\end{proof}
\leavevmode\phantomsection\label{v98:par-quadratic-conditions}\ignorespaces
When choosing the initial polynomial, the nonsquare product of its values controls the first layer; when \(q\equiv3\pmod4\), the additional even-degree condition keeps this square class nonsquare at the next layer.
Over \(\F_5\), one may start from \(X-2\), since the product of its values at \(1\) and \(-1\) is \(3\), a nonsquare.
We next give a recursive example over \(\F_3\), together with a comparison where the first step succeeds but a later step splits.

\begin{example}[A square-class invariant and a recursive tower of quadratic fibers]\label{v48:ex-r-transform-tower}
In \(\F_3[X]\), take
\[
P_0=X^2+X-1,\qquad P_{j+1}=\mathcal R(P_j).
\]
Since \(P_0\) has no root in the base field, it is irreducible, and
\[
P_0(1)=1,\qquad P_0(-1)=2.
\]
Put \(u=X+X^{-1}\). Since \(1/2=-1\) in \(\F_3\), the first two transforms can be written explicitly as
\begin{equation}\label{v48:eq-r-transform-polynomials}
\begin{aligned}
P_1&=X^2(u^2-u-1)
=X^4+2X^3+X^2+2X+1,\\
P_2&=X^4(u^4+u^3+u^2+u+1)\\
&=X^8+X^7+2X^6+X^5+X^3+2X^2+X+1.
\end{aligned}
\end{equation}
These expansions display the coefficients; we now explain their irreducibility through fiber return.

Choose a root \(\alpha\in\F_9\) of \(P_0\). The map
\[
h(x)=\frac{x+x^{-1}}2
\]
has fiber over \(\alpha\) given by \(X^2-2\alpha X+1\), with discriminant \(4(\alpha^2-1)\). By the norm product formula,
\[
\Nm_{\F_9/\F_3}(\alpha^2-1)=P_0(1)P_0(-1)=2.
\]
A nonzero element of a finite field is a square if and only if its norm is a square; see the proof of Theorem~\ref{thm:reciprocal-transform}.
Thus \(\alpha^2-1\) is a nonsquare in \(\F_9\), and the quadratic fiber polynomial is irreducible.
If \(\beta\) is one root, the other is \(\beta^{-1}\), and
\[
\beta^9=\beta^{-1}\ne\beta.
\]
The base orbit has length two, and the return map \(x\mapsto x^9\) interchanges the two roots in the fiber.
By Lemma~\ref{lem:fiber-return}, the total orbit has length \(2\cdot2=4\), so \(P_1\) is irreducible.

The criterion preserved for the second lift can be read directly from base-field values:
\[
P_1(1)=1,\qquad P_1(-1)=2,\qquad
\Nm_{\F_{81}/\F_3}(\beta^2-1)=2.
\]
Thus \(X^2-2\beta X+1\) remains irreducible over \(\F_{81}\).
Choose a root \(\gamma\). Then
\[
\gamma^{81}=\gamma^{-1}\ne\gamma.
\]
Now the base orbit has length four and the return map interchanges the two-element fiber, so the total orbit of \(\gamma\) has length eight, and \(P_2\) is irreducible.
In each calculation, inversion agrees with the Frobenius iterate at half the full orbit length, which also explains the self-reciprocity of \(P_1\) and \(P_2\).

The recursion continues because the square class controlling the quadratic fiber is preserved throughout.
For \(N_j=\deg P_j=2^{j+1}\), we have
\[
P_{j+1}(1)=2^{N_j}P_j(1)=P_j(1),\qquad
P_{j+1}(-1)=(-2)^{N_j}P_j(-1)=P_j(-1).
\]
Hence \(P_j(1)=1\) and \(P_j(-1)=2\) for all \(j\ge0\); the nonsquare-norm criterion holds at every step.
Each doubling of the degree corresponds to a nontrivial return in a new two-element fiber, and the two base-field values ensure that this return can be repeated indefinitely.

The parity condition on the degree is essential.
If instead we start from \(F_0=X\in\F_3[X]\), then \(F_0(1)F_0(-1)=2\), and the first step does give the irreducible polynomial
\(F_1=\mathcal R(F_0)=X^2+1\). However, \(F_1(1)F_1(-1)=1\), and the second step gives
\[
\mathcal R(F_1)=X^4+1
=(X^2+X+2)(X^2+2X+2).
\]
The odd initial degree causes the first transform to change the square class.
Thus one-step irreducibility and a continuing recursion impose different conditions: the former requires the current quadratic fiber to be irreducible, while the latter also requires this criterion to be preserved by the transform.
\end{example}

\section{Polynomial Iteration and Dynamical Irreducibility}\label{v55:sec-dynamical-irreducibility}
The preceding recursive constructions preserve a trace, a square class, or a multiplicative order at every step, so that each new root fiber continues to have a single return orbit.
Repeated substitution of one fixed polynomial gives another class of towers of irreducible polynomials.
The root-fiber criterion still applies, but we must explain why the conditions controlling the fibers hold at every level.

\subsection{Prescribed Targets and Successive Root Fibers}
\begin{definition}[Dynamical irreducibility]\label{v55:def-dynamical-irreducibility}
Let \(h\in k[X]\), with \(\deg h=d\ge2\), and set
\[
h^{\circ0}(X)=X,\qquad h^{\circ(n+1)}=h\circ h^{\circ n}.
\]
If \(h^{\circ n}\) is irreducible in \(k[X]\) for every \(n\ge1\),
then \(h\) is called \textbf{dynamically irreducible}, or \textbf{stable}.
More generally, if \(\tau\in k\) and every \(h^{\circ n}(X)-\tau\) is irreducible,
then the pair \((h,\tau)\) is called dynamically irreducible.
\end{definition}
Ordinary stability corresponds to the target \(\tau=0\). Iteration here takes place in the polynomial ring;
the degree is \(d^n\), without reduction modulo \(X^q-X\).

\begin{corollary}[Fiber criterion for iteration]\label{v55:cor-iterated-capelli}
Suppose that \(h^{\circ n}(X)-\tau\) is irreducible over \(k\), where \(n\ge0\), and choose a root \(\beta_n\).
Then \(h^{\circ(n+1)}(X)-\tau\) is irreducible if and only if
\(h(X)-\beta_n\) is irreducible over \(k(\beta_n)\).
For \(n=0\), take \(\beta_0=\tau\); for \(n\ge1\), we have \([k(\beta_n):k]=d^n\).
\end{corollary}
\begin{proof}
Apply Corollary~\ref{thm:capelli} to the monic normalization of the outer polynomial \(h^{\circ n}(Y)-\tau\) and the inner polynomial \(h(X)\).
\end{proof}
When every level is irreducible, we may choose \(h(\beta_{n+1})=\beta_n\), obtaining
\[
k\subseteq k(\beta_1)\subseteq k(\beta_2)\subseteq\cdots,
\qquad [k(\beta_{n+1}):k(\beta_n)]=d.
\]
Each multiplication of the degree by \(d\) corresponds to a full return orbit in a new root fiber.

The iteration level \(n\) and a separate base-extension degree \(r\) have different roles: the former changes the polynomial to \(h^{\circ n}-\tau\) while keeping \(k\) fixed;
the latter keeps that polynomial fixed and replaces Frobenius \(\sigma\) by \(\sigma^r\).
If this level is irreducible over \(k\), Theorem~\ref{thm:factor-base-extension} gives
\(\gcd(d^n,r)\) distinct irreducible factors over \(\F_{q^r}\), each of degree \(d^n/\gcd(d^n,r)\).
Thus a dynamically irreducible tower remains dynamically irreducible after base extension if and only if \(\gcd(d,r)=1\); see Exercise~\ref{v55:ex-stability-base-extension}.

\Needspace{7\baselineskip}

\begin{proposition}[Transforming the polynomial and target together]\label{v55:prop-pointed-conjugacy}
Let \(\phi(X)=uX+v\), where \(u\in k^\times\) and \(v\in k\), and put
\(g=\phi\circ h\circ\phi^{-1}\). Then
\[
g^{\circ n}(X)-\phi(\tau)
=u\bigl(h^{\circ n}(\phi^{-1}(X))-\tau\bigr).
\]
Thus \((h,\tau)\) is dynamically irreducible if and only if \((g,\phi(\tau))\) is.
\end{proposition}
\begin{proof}
Conjugation commutes with iteration, so \(g^{\circ n}=\phi\circ h^{\circ n}\circ\phi^{-1}\).
Substitution of \(\phi(\tau)=u\tau+v\) gives the formula; invertible linear substitution and multiplication by a nonzero constant preserve irreducibility.
\end{proof}
We refer to this simultaneous change of polynomial and target as \textbf{pointed affine conjugacy}.
Prescribing the target keeps every root fiber after an affine change of variables in correspondence with an original fiber.
For example, stability corresponds to the new target \(v\) after conjugation by a translation; see~\cite{DayDynamical2025}.

\par\Needspace{12\baselineskip}
\subsection{Iteration of Power Maps and the Binomial Criterion}
\begin{theorem}[Stability of translated power maps]\label{v55:thm-shifted-power-stability}
Let \(d\ge2\), \(\zeta\in k^\times\), and put \(f(X)=(X+\zeta)^d-\zeta\). Then
\[
f^{\circ n}(X)=(X+\zeta)^{d^n}-\zeta.
\]
The polynomial \(f\) is stable if and only if \(X^d-\zeta\) is irreducible and, when \(d\) is even, \(q\equiv1\pmod4\).
Equivalently, putting \(e=\ord(\zeta)\), the conditions are:
\begin{enumerate}
\item Every prime \(\ell\mid d\) divides \(e\) but does not divide \((q-1)/e\).
\item If \(2\mid d\), then \(q\equiv1\pmod4\).
\end{enumerate}
\end{theorem}
\begin{proof}
The map \(\phi(X)=X+\zeta\) gives \(\phi\circ f\circ\phi^{-1}=X^d\), proving the iteration formula.
A linear substitution reduces irreducibility at level \(n\) to that of \(X^{d^n}-\zeta\).
Since \(d^n\) and \(d\) have the same prime divisors, the first condition of Theorem~\ref{thm:binomial-full} is the same at every level.
When \(d\) is odd, the condition modulo four in that criterion never arises;
when \(d\) is even, \(4\mid d^2\), so \(q\equiv1\pmod4\) is required from the second level onward.
This gives the stated necessary and sufficient conditions.
\end{proof}
Consequently, when \(d\) is odd, the first level suffices to test stability in this family; when \(d\) is even, the first two levels must be tested.
This result also appears in~\cite{AnWu2026}; the closed iteration formula reduces it to the multiplicative-period criterion already established.

\begin{example}\label{v55:ex-shifted-power-stability}
The binomial examples in Chapter~\ref{chap:root-fibers} give two stable polynomials:
\[
(X+2)^2-2\in\F_5[X],\qquad (X+2)^3-2\in\F_7[X].
\]
Over \(\F_3\), the polynomial \(f=(X+2)^2-2=X^2+X+2\) has no root and is therefore irreducible, but its second iterate satisfies
\[
f^{\circ2}=(X+2)^4+1=(X^2+1)(X^2+2X+2).
\]
The congruence condition modulo four is absent at the first level but fails at the second.
\end{example}

\subsection{Quadratic Iteration and Finite Critical Orbits}
\leavevmode\phantomsection\label{def:critical-point-value}\ignorespaces
For a nonconstant polynomial \(H\) whose formal derivative \(H'\) is nonzero, a root of \(H'\) in a splitting extension is called a \textbf{critical point},
and the corresponding value of \(H\) is called a \textbf{critical value}.
The successive images \(H^{\circ j}(\gamma)\), \(j\ge0\), of a critical point \(\gamma\) form its \textbf{critical orbit}.
We now take \(h=X^2+c\) in odd characteristic; its unique critical point is \(0\).

\Needspace{7\baselineskip}

\begin{theorem}[Square-class criterion for quadratic iteration]\label{v55:thm-quadratic-stability}
Let \(q\) be odd, \(h=X^2+c\in k[X]\), and \(\tau\in k\).
Then \((h,\tau)\) is dynamically irreducible if and only if
\[
\tau-c\text{ is a nonsquare},\qquad
h^{\circ j}(0)-\tau\text{ is a nonsquare}\quad(j\ge2).
\]
Here a nonsquare means an element of \(k^\times\) that is not a square.
\end{theorem}
\begin{proof}
The first level is \(X^2-(\tau-c)\), whose irreducibility condition is precisely that \(\tau-c\) is a nonsquare.
Suppose that level \(n\) is irreducible, with \(n\ge1\), choose a root \(\beta_n\), and put \(E=k(\beta_n)\).
By Corollary~\ref{v55:cor-iterated-capelli}, the next level is irreducible if and only if \(\beta_n-c\) is a nonsquare in \(E\).
The power-residue expression for the norm, formula~\eqref{eq:residue-norm-resultant}, gives
\[
z^{(q^{2^n}-1)/2}
=\Nm_{E/k}(z)^{(q-1)/2}\qquad(z\in E^\times).
\]
Thus nonsquareness can be tested in \(k\) by taking the norm; the norm of zero is also zero.
Since \(2^n\) is even,
\[
\Nm_{E/k}(\beta_n-c)
=h^{\circ n}(c)-\tau=h^{\circ(n+1)}(0)-\tau.
\]
Thus level \(n+1\) is controlled precisely by the listed critical-orbit value of index \(n+1\).
Induction on the levels proves the result.
\end{proof}
The critical-orbit criterion for quadratics is a basic result in the study of dynamical irreducibility; see~\cite{JonesBoston2012,DayDynamical2025}.
In this proof, the relevant critical-orbit value minus the target is the norm of the extension-field fiber parameter, so one square-class condition controls the entire two-element fiber.

\begin{corollary}[A finite test for quadratic stability]\leavevmode\phantomsection\label{v55:cor-quadratic-finite-test}
Retain the hypotheses of the preceding theorem. First test \(\tau-c\), then iterate point by point starting from \(h^{\circ2}(0)\).
If every value \(h^{\circ j}(0)-\tau\) encountered is a nonsquare until an orbit point repeats, then \((h,\tau)\) is dynamically irreducible.
In particular, it suffices that \(h^{\circ n}(X)-\tau\) be irreducible for \(1\le n\le(q+3)/2\).
\end{corollary}
\begin{proof}
Once an orbit point repeats, all subsequent points repeat, so the preceding theorem proves the first assertion.
For the stated bound, put \(B=(q+3)/2\) and suppose that \(h^{\circ n}(X)-\tau\) is irreducible for \(1\le n\le B\).
The successive norm tests in the preceding proof show that \(\tau-c\) and all values \(h^{\circ j}(0)-\tau\), \(2\le j\le B\), are nonsquares.
Hence the \(B-1=(q+1)/2\) orbit points \(h^{\circ j}(0)\), \(2\le j\le B\), all lie in
\[
\tau+\{u\in k^\times:u\text{ is a nonsquare}\},
\]
a set of \((q-1)/2\) elements. Two of these points therefore coincide.
Every subsequent point repeats within this segment, so all later square-class conditions hold as well.
Together with the first-level condition, the preceding theorem gives dynamical irreducibility.
\end{proof}
This finite test examines base-field elements in the critical orbit, without expanding iterate polynomials of degree \(2^n\).

\Needspace{7\baselineskip}
\begin{corollary}[Normalization and stability of a general quadratic]\label{v56:cor-general-quadratic-stability}
Let \(q\) be odd, and let \(f(X)=aX^2+bX+c\in k[X]\), where \(a\ne0\). Put
\[
\begin{aligned}
\phi(X)&=aX+\frac b2,& h(X)&=X^2+C,\\
C&=ac+\frac b2-\frac{b^2}{4},& \tau&=\frac b2.
\end{aligned}
\]
Then \(\phi\circ f\circ\phi^{-1}=h\), and \(f\) is stable if and only if \((h,\tau)\) is dynamically irreducible.
Writing \(\Delta=b^2-4ac\) and \(\gamma=-b/(2a)\), this condition is equivalent to
\[
\Delta\text{ is a nonsquare},\qquad
 a f^{\circ j}(\gamma)\text{ is a nonsquare}\quad(j\ge2).
\]
Here \(\gamma\) is the critical point of \(f\), and all square classes are tested in \(k\).
\end{corollary}
\begin{proof}
Substituting \(\phi^{-1}(X)=(X-b/2)/a\) directly gives the stated conjugacy identity.
Since \(\phi(0)=\tau\), Proposition~\ref{v55:prop-pointed-conjugacy} gives
\[
h^{\circ n}(X)-\tau
=a f^{\circ n}(\phi^{-1}(X)),
\]
so irreducibility agrees at every level on the two sides. Since \(\phi(\gamma)=0\), we also have
\[
\tau-C=\frac{\Delta}{4},\qquad
h^{\circ j}(0)-\tau=\phi(f^{\circ j}(\gamma))-\tau
=a f^{\circ j}(\gamma).
\]
As \(4\) is a nonzero square, Theorem~\ref{v55:thm-quadratic-stability} gives the necessary and sufficient conditions.
\end{proof}
Affine conjugacy puts the root fibers of every quadratic in odd characteristic into monic form while transforming the target \(0\) to \(\tau\).
Thus the finite test of Corollary~\ref{v55:cor-quadratic-finite-test} also applies to general quadratic polynomials.

\begin{example}[Preservation and failure along critical orbits]\leavevmode\phantomsection\label{v55:ex-quadratic-stability}
Over \(\F_3\), the critical orbit of \(h=X^2+1\) is
\[
0\longmapsto1\longmapsto2\longmapsto2.
\]
The first-level condition is that \(-1=2\) is a nonsquare, and every later value to be tested is also \(2\), so \(h\) is stable.
Taking \(\phi=X+1\), the conjugate polynomial is \(g=X^2+X\); the pair \((g,1)\) is dynamically irreducible.

Over the same field, the nonmonic polynomial \(f=2X^2+2\) has critical point \(\gamma=0\), discriminant \(\Delta=2\), and critical orbit
\[
0\longmapsto2\longmapsto1\longmapsto1.
\]
For \(j\ge2\), the value \(f^{\circ j}(\gamma)=1\) is a square, whereas multiplication by the leading coefficient gives the nonsquare \(a f^{\circ j}(\gamma)=2\).
By Corollary~\ref{v56:cor-general-quadratic-stability}, \(f\) is stable.
With \(\eta(X)=2X\), we have \(\eta\circ f\circ\eta^{-1}=h\) and \(\eta(0)=0\), so conjugacy with the same target gives the same conclusion.
This example shows that the leading-coefficient factor \(a\) cannot be omitted from the square-class criterion for \(j\ge2\).

On the other hand, over \(\F_5\), the polynomial \(X^2+2\) is irreducible, but the second iterate of its critical point is \(1\), and
\[
(X^2+2)^{\circ2}=X^4-X^2+1
=(X^2+2X-1)(X^2-2X-1).
\]
The change of square class in the second fiber causes the second level to split.
\end{example}

\subsection{Cubic Iteration in Characteristic Two and Reciprocal Coordinates}
In characteristic two, the identity
\[
(u+u^{-1})^3+(u+u^{-1})=u^3+u^{-3}
\]
connects \(h=X^3+X\) with the cubing map.
Reciprocal coordinates turn successive fibers into lifts by cube roots, and the norm can preserve noncubicity throughout.

\Needspace{7\baselineskip}

\begin{theorem}[Cubic iteration with a prescribed target]\label{v55:thm-char2-target-cubic}
Let \(q=2^m\), \(h=X^3+X\), and \(\tau\in\F_q\).
Then \((h,\tau)\) is dynamically irreducible if and only if \(h(X)-\tau\) is irreducible over the base field.
\end{theorem}
\begin{proof}
Only sufficiency needs proof. Put \(E=k(\omega)\), where \(\omega\) is a primitive cube root of unity.
If \(q\equiv1\pmod3\), then \(E=k\); if \(q\equiv-1\pmod3\), then \(E=\F_{q^2}\).
Irreducibility ensures \(\tau\ne0\), so \(h-\tau\) and \(h'=X^2+1\) are coprime. Choose
\[
t^2+\tau t+1=0,\qquad u^3=t.
\]
The three distinct roots of \(h(X)-\tau\) are
\[
x_i=\omega^iu+\omega^{-i}u^{-1},\qquad i\in\Z/3\Z.
\]
If \(q\equiv1\pmod3\) but \(t\notin k\), then \(t^q=t^{-1}\) and \(u^q=\omega^a/u\).
Frobenius acts on root indices as \(i\mapsto-i-a\), which has a fixed point, contradicting irreducibility.
If \(q\equiv-1\pmod3\) but \(t\in k\), then cubing is a bijection on \(k^\times\),
so we may take \(u\in k\), again producing a base-field root.
Thus \(t\in E\); in the latter congruence class, we also have \(t^q=t^{-1}\).
The element \(t\) must be a noncube in \(E\); otherwise we could take \(u\in E\), producing a root of degree at most two over \(k\).

Put \(K_j=\F_{q^{3^j}}\) and \(L_j=EK_j\), and induct starting from \(\beta_0=\tau=t_0+t_0^{-1}\), \(t_0=t\).
Suppose that \(\beta_j\) has degree \(3^j\) over \(k\), that \(t_j\in L_j\) is a noncube, and that \(\beta_j=t_j+t_j^{-1}\).
Choose \(U^3=t_j\). The prime-degree binomial criterion gives \([L_j(U):L_j]=3\), hence \(L_j(U)=L_{j+1}\).
Put \(\beta_{j+1}=U+U^{-1}\); then \(h(\beta_{j+1})=\beta_j\).
When \(q\equiv1\pmod3\), we have \(\beta_{j+1}\in K_{j+1}\).
When \(q\equiv-1\pmod3\), we also preserve \(t_j^{q^{3^j}}=t_j^{-1}\) in the induction.
Since three is odd,
\[
\bigl(U^{q^{3^{j+1}}}U\bigr)^3=1.
\]
The expression in parentheses is the norm of \(U\) from \(L_{j+1}\) to \(K_{j+1}\), and hence lies in \(K_{j+1}\).
The multiplicative group of this field has order not divisible by three, so the expression equals one.
Thus \(U^{q^{3^{j+1}}}=U^{-1}\), and again \(\beta_{j+1}\in K_{j+1}\).

If \(\beta_{j+1}\in L_j\), then \(U\) satisfies \(X^2+\beta_{j+1}X+1\),
contradicting its degree three over \(L_j\).
Since \(\beta_{j+1}\in K_{j+1}\), \([K_{j+1}:K_j]=3\), and \(\beta_{j+1}\notin L_j\supseteq K_j\),
we have \([K_j(\beta_{j+1}):K_j]=3\); then \(\beta_j=h(\beta_{j+1})\) gives
\([k(\beta_{j+1}):k]=3^{j+1}\).
Finally, the three conjugates of \(U\) are \(U,\omega U,\omega^2U\), so
\[
\Nm_{L_{j+1}/L_j}(U)=U^3=t_j.
\]
If \(U\) were a cube, its norm would also be a cube, a contradiction. Put \(t_{j+1}=U\); all inductive conditions continue to hold.
Thus \(\beta_n\) is a root of \(h^{\circ n}(X)-\tau\) of degree \(3^n\), so this polynomial is irreducible.
\end{proof}

\begin{corollary}[A family of stable cubics]\label{v55:cor-linwang-stable-cubic}
Let \(q\) be a power of two and \(a,b\in\F_q^\times\).
Then \(f(X)=b^2(X+a)^3+X\) is stable whenever it is irreducible.
\end{corollary}
\begin{proof}
Taking \(\phi(X)=b(X+a)\), direct calculation gives \(\phi\circ f=h\circ\phi\) and \(\phi(0)=ab\).
Apply Proposition~\ref{v55:prop-pointed-conjugacy} and the preceding theorem.
\end{proof}
This cubic family appears in~\cite{LinWangStable2024}; translated power families of degrees two and four are controlled by the preceding binomial criterion.

\Needspace{7\baselineskip}
\begin{example}[A stable cubic over the binary field]\label{v55:ex-binary-stable-cubic}
The polynomial \(h(X)-1=X^3+X+1\) has no root over \(\F_2\), so it is irreducible.
Taking \(a=b=1\) gives
\[
f(X)=(X+1)^3+X=X^3+X^2+1,
\qquad f^{\circ n}(X)=h^{\circ n}(X+1)+1.
\]
Thus every positive iterate of \(f\) is irreducible, of degree \(3^n\).

The first two lifts of the reciprocal parameter can be made explicit. Choose
\[
t^2+t+1=0,\qquad U^3=t,\qquad V^3=U,
\]
and put \(\beta_0=1\), \(\beta_1=U+U^{-1}\), and \(\beta_2=V+V^{-1}\).
Since \(t\ne1\) and \(t^3=1\), we have
\[
\ord(t)=3,\qquad \ord(U)=9,\qquad \ord(V)=27.
\]
From \(\ord_3(2)=2\), \(\ord_9(2)=6\), and \(\ord_{27}(2)=18\), the least extension fields containing the reciprocal parameters are, successively,
\[
\F_2(t)=\F_4,\qquad
\F_2(U)=\F_{64},\qquad
\F_2(V)=\F_{2^{18}}.
\]
On the other hand, \(2^3\equiv-1\pmod9\) and \(2^9\equiv-1\pmod{27}\), so
\[
U^{2^3}=U^{-1},\qquad V^{2^9}=V^{-1},
\qquad
\beta_1\in\F_8,\qquad \beta_2\in\F_{512}.
\]
These are precisely the fields generated by \(\beta_1,\beta_2\), respectively, over the prime field. Indeed, for \(W=U\) or \(V\),
\(W+W^{-1}\) is fixed by the \(2^r\)-Frobenius if and only if
\(W^{2^r}\in\{W,W^{-1}\}\): these two elements are
the two roots of \(X^2+(W+W^{-1})X+1\).
Hence its Frobenius period is the least positive integer satisfying \(2^r\equiv\pm1\pmod{\ord(W)}\), namely \(3\) and \(9\), respectively.
Thus
\[
\F_2(\beta_1)=\F_8,\qquad
\F_2(\beta_2)=\F_{512}.
\]
The reciprocal-coordinate identity also gives
\[
h(\beta_1)=t+t^{-1}=1,\qquad h(\beta_2)=\beta_1.
\]
Therefore \(\beta_1\) and \(\beta_2\) are roots of \(h(X)-1\) and \(h^{\circ2}(X)-1\), respectively,
with degrees exactly \(3\) and \(9\); the corresponding roots of the iterates of \(f\) are \(\beta_1+1\) and \(\beta_2+1\), and generate the same respective fields.

In the notation of Theorem~\ref{v55:thm-char2-target-cubic}, fix the auxiliary field \(E=\F_4\). The three types of fields are
\[
\begin{array}{c|c|c|c}
j&\text{Auxiliary field }E&\text{Iterate root field }K_j&\text{Reciprocal parameter field }L_j=EK_j\\
\hline
0&\F_4&\F_2&\F_4\\
1&\F_4&\F_8&\F_{64}\\
2&\F_4&\F_{512}&\F_{2^{18}}
\end{array}
\]
At every level, the parameter field is a quadratic extension of the corresponding root field, and the reciprocal sum descends to the root field.
Also, \(\Nm_{\F_{64}/\F_4}(U)=t\) and \(\Nm_{\F_{2^{18}}/\F_{64}}(V)=U\),
showing that norms preserve noncubicity at successive levels, while Frobenius return periods determine the degrees of the iterate roots.
Section~\ref{sec:dickson-polynomials} places this reciprocal-coordinate identity in the general root-parameter model of Dickson polynomials.
\end{example}

\subsection{Mixed Compositions and Forward-Invariant Sets in the Base Field}
\begin{definition}[Dynamically irreducible set]\label{v55:def-dynamical-set}
Let \(\mathcal A=\{f_1,\ldots,f_s\}\subseteq k[X]\), with \(\deg f_i\ge2\) for every \(i\).
If, for every length \(r\ge1\) and every choice of indices \(i_1,\ldots,i_r\),
the composition \(f_{i_1}\circ\cdots\circ f_{i_r}\) is irreducible, then \(\mathcal A\) is called a \textbf{dynamically irreducible set}.
\end{definition}
For translated power families, norm descent for power residues allows mixed compositions of arbitrary length to be controlled by a forward-invariant set in the base field; see~\cite{SharmaDynamical2026}.

\begin{theorem}[A finite criterion for translated power families of odd degree]\label{v55:thm-mixed-power-stability}
Let \(d\ge3\) be odd, suppose that every prime \(\ell\mid d\) divides \(q-1\), and put
\[
f_i(X)=(X-a_i)^d+b_i\quad(1\le i\le s),\qquad a_i,b_i\in k.
\]
Let \(S\subseteq k\) be the smallest set containing all \(b_i\) and satisfying \(f_i(S)\subseteq S\).
Then \(\{f_1,\ldots,f_s\}\) is dynamically irreducible if and only if every element of \(S\) is nonzero
and is not an \(\ell\)-th power in \(k\) for any prime \(\ell\mid d\).
\end{theorem}
\begin{proof}
Let \(G\) be an irreducible mixed composition, let \(D=\deg G\) be odd, choose a root \(\alpha\), and put \(E=k(\alpha)\).
The fiber polynomial for the next substitution of \(f_i\) on the right is
\[
f_i(X)-\alpha=(X-a_i)^d-(\alpha-b_i).
\]
By the binomial criterion, it is irreducible if and only if \(\alpha-b_i\ne0\) and it is not an \(\ell\)-th power for any \(\ell\mid d\).
The norm expression gives
\[
(\alpha-b_i)^{(q^D-1)/\ell}
=\Nm_{E/k}(\alpha-b_i)^{(q-1)/\ell},\qquad
\Nm_{E/k}(\alpha-b_i)=(-1)^D G(b_i)=-G(b_i).
\]
All \(\ell\) are odd, and \(-1=(-1)^\ell\), so changing the sign does not change \(\ell\)-th-power status.
Thus \(G\circ f_i\) is irreducible precisely when \(G(b_i)\) satisfies the stated power-residue conditions.
The same argument applies to the initial outer polynomial \(G=X\), with \(D=1\) and \(\alpha=0\), giving the irreducibility condition for each \(f_i\).

The set \(S\) consists exactly of all \(G(b_i)\), where \(G\) runs through all mixed compositions and the identity polynomial \(X\).
If \(S\) satisfies the conditions, the one-step equivalence above and induction on composition length ensure that all mixed compositions are irreducible.
Conversely, if all mixed compositions are irreducible, applying the equivalence to each \(G\) and \(f_i\)
shows that every \(G(b_i)\) satisfies the conditions, proving necessity as well.
\end{proof}
Starting from \(S_0=\{b_1,\ldots,b_s\}\), repeatedly adjoining every \(f_i(x)\) produces \(S\).
At most \(q\) distinct elements occur; by recording the elements already processed, each element need only be evaluated once under each \(f_i\).
Thus the condition on infinitely many mixed compositions in this family admits a finite test in the base field.

\Needspace{7\baselineskip}

\begin{example}\label{v55:ex-mixed-cubics-seven}
In \(\F_7\), the nonzero cubes are \(1,6\). Put
\[
f_1=X^3+3,\qquad f_2=X^3+4,
\qquad S=\{2,3,4,5\}.
\]
For \(x=2,4\), we have \((f_1(x),f_2(x))=(4,5)\);
for \(x=3,5\), we have \((f_1(x),f_2(x))=(2,3)\).
Thus \(S\) contains both constant terms, is closed under both polynomials, and consists entirely of noncubes.
The preceding theorem ensures that every mixed composition is irreducible; a composition of length \(r\) has degree \(3^r\).
\end{example}

\par\Needspace{12\baselineskip}
\subsection{Denominators of Rational Iterates and Inverse Stability}
\begin{definition}[Inverse stability]\label{v55:def-inverse-stability}
Let \(\psi\in k[X]\), with \(\deg\psi\ge2\), and put \(\Psi=1/\psi\). For each integer \(n\ge1\), write its rational iterate as
\[
\Psi^{\circ n}(X)=\frac{A_n(X)}{B_n(X)},\qquad
\gcd(A_n,B_n)=1,\quad B_n\text{ monic}.
\]
If all \(B_n\) are irreducible and pairwise distinct, then \(\psi\) is called \textbf{inversely stable}.
\end{definition}
The roots of these denominators are the finite poles of the corresponding rational iterates: the points where the denominator vanishes and the numerator does not.
The reciprocal additive recursion of Section~\ref{sec:recursive-as} gives an example for a fixed rational map.

\begin{proposition}[An inversely stable Artin--Schreier family]\leavevmode\phantomsection\label{v55:prop-inverse-as}
Let \(s\in\F_p^\times\) and \(q=p^e\). Then \(\psi=X^p-X+s\) is inversely stable over \(\F_q\)
if and only if \(p\nmid e\). When this holds, the \(n\)-th denominator has degree \(p^n\) for every integer \(n\ge1\).
\end{proposition}
\begin{proof}
First work over \(\F_p\). Put \(D_1=\psi\), and, following Corollary~\ref{cor:recursive-prime-field}, define
\[
D_{n+1}(X)=\widehat D_n(\psi(X)),\qquad
\widehat D_n(Y)=Y^{p^n}D_n(Y^{-1}).
\]
Taking \(P=X\) and \(a=s\) in that corollary shows that each \(D_n\) is irreducible and \(\deg D_n=p^n\).
We now show that these are also denominators of the iterates of \(\Psi=1/\psi\).
This is clear at the first level. Suppose inductively that \(\Psi^{\circ n}=C_n/D_n\) is in lowest terms and \(\deg C_n\le p^n\). Then
\[
\Psi^{\circ(n+1)}=\Psi^{\circ n}\circ\Psi
=\frac{\psi(X)^{p^n}C_n(\psi(X)^{-1})}
{\psi(X)^{p^n}D_n(\psi(X)^{-1})}.
\]
Both expressions are polynomials, and the denominator is exactly \(D_{n+1}\).
At a zero of \(\psi\), the denominator equals the nonzero leading coefficient of \(D_n\);
at any other point, simultaneous vanishing would give a common root of \(C_n,D_n\). Hence the two expressions have no common factor.
The new numerator has degree at most \(p^{n+1}\), so the degree condition in the induction persists.
Thus the monic normalization of the reduced denominator equals the monic normalization of \(D_n\).

After extending the base field to \(\F_{p^e}\), Theorem~\ref{thm:factor-base-extension} gives:
\(D_n\) remains irreducible if and only if \(\gcd(p^n,e)=1\).
If \(p\nmid e\), every level is irreducible, and the differing degrees also ensure that the denominators are pairwise distinct;
if \(p\mid e\), the first level is already reducible.
\end{proof}
This family also belongs to the inverse-stability constructions studied in~\cite{AnWu2026}; its proof follows directly from the reciprocal additive recursion already established.
The ordinary polynomial iterates of \(\psi\), however, are affine linearized polynomials:
when \(p\) is odd, the second iterate has degree \(p^2\); when \(p=2\), the third iterate has degree eight.
Theorem~\ref{thm:affine-linearized-obstruction} rules out irreducibility at these respective levels.

The constructions above preserve multiplicative orders, square classes, noncube classes, or nonzero traces after reciprocation.
The quadratic critical-orbit criterion and the forward-invariant set for translated power families of odd degree give two explicit finite tests.
General cubic fibers contain further information: in characteristic different from \(2,3\), a nonzero square discriminant allows both irreducibility and complete splitting,
so the square class of the discriminant gives only partial factorization information.

Proposition~\ref{v105:prop-discriminant-parity} expresses this limitation through the parity of the root permutation; Example~\ref{v105:ex-cubic-same-discriminant} gives two cubic factorization types with the same discriminant.
For cubic factorization models on short critical orbits and their connection with iterated Galois groups over the rationals, see~\cite{SanMartinCubic2026}.

\section{Roots of Unity and Sparse Trinomial Families}\label{v49:sec-sparse-families}
\begin{proposition}\label{prop:sparse-trinomial-family}
Suppose that \(q\equiv3\pmod4\), and put \(r=v_2(q+1)\). Choose an element \(\zeta\) of order \(2^{r+1}\) in \(\F_{q^2}\), and put \(c=\zeta+\zeta^q\in k\).
Then \(c\ne0\), and, for every \(j\ge1\),
\[
X^{2^j}-cX^{2^{j-1}}-1
\]
is an irreducible trinomial over \(k\).
\end{proposition}
\begin{proof}
We have \(\zeta^{q+1}=-1\), so the minimal polynomial of \(\zeta\) is \(P=X^2-cX-1\).
The element \(\zeta\) does not lie in \(k\), since the \(2\)-part of its order exceeds the \(2\)-part of \(q-1\).
If \(c=0\), then \(\zeta^q=-\zeta\); together with \(\zeta^{q+1}=-1\), this gives \(\zeta^2=1\), a contradiction.
In \(\F_{q^2}\), the order \(\ord(\zeta)=2^{r+1}\) contains the full \(2\)-part of \(q^2-1\), and \(q^2\equiv1\pmod4\).
For \(j=1\), the root fiber is linear; for \(j\ge2\), the binomial criterion ensures that \(X^{2^{j-1}}-\zeta\) is irreducible. The root-fiber criterion then makes \(P(X^{2^{j-1}})\) irreducible.
\end{proof}

Over \(\F_2\), there is another simple family of trinomials:
\[
X^{2\cdot3^j}+X^{3^j}+1=\Phi_{3^{j+1}}(X),\qquad j\ge0.
\]
These polynomials are irreducible because \(\ord_{3^{j+1}}(2)=2\cdot3^j\).
The period equality follows immediately from \(2^2=4\) and Lemma~\ref{lem:valuation-lift}.
For the classification of three-sparse irreducible factors of root-of-unity polynomials, and the family involving powers of \(7\) in characteristic two, see Subsection~\ref{v49:subsec-three-sparse} and Proposition~\ref{v49:prop-seven-power-trinomials}.

\subsection{Counting Sparse Families and Conjugate Roots}
Every irreducible \(f\) of degree at least two satisfies \(f(0)\ne0\), \(f(1)\ne0\), and \(\gcd(f,f')=1\).
Hence an irreducible polynomial over \(\F_2\) of degree at least two has an odd number of nonzero terms and at least one term of odd degree. These are only necessary conditions.

\begin{proposition}[A square-root estimate for a sparse family]\label{prop:sparse-square-root}
Let \(n\ge2\), \(q=p^r\), and \(p\nmid2n(n-1)\), and put
\[
T_n(q)=\#\{b\in\F_q:X^n+X+b\text{ is irreducible}\}.
\]
Then there is a constant \(B_n\), depending only on \(n\), such that
\[
\left|T_n(q)-\frac qn\right|\le B_n\sqrt q.
\]
\end{proposition}
This is a special case of Cohen's polynomial distribution theorem; for a proof, see~\cite[Theorem 1 and p.~257]{Cohen}.
Critical points and critical values are defined in Section~\ref{v55:sec-dynamical-irreducibility}.
The restriction \(p\nmid n(n-1)\) on the characteristic has the following interpretation. Put \(H(X)=X^n+X\); then \(H'(X)=nX^{n-1}+1\). Its roots are all nonzero,
and \(H''(X)=n(n-1)X^{n-2}\) is nonzero at these roots, so there are exactly \(n-1\) distinct critical points.
If \(H'(\beta)=0\), then
\[
\beta^{n-1}=-\frac1n,\qquad H(\beta)=\frac{n-1}{n}\beta.
\]
The coefficient \((n-1)/n\) is nonzero, so the critical values are also pairwise distinct.
The characteristic condition also excludes \(p=2\), as required by the cited distribution theorem.

For \(X^n+X+b\) with distinct roots, irreducibility is equivalent to Frobenius acting on its root set as an \(n\)-cycle.
The factor \(1/n\) in the main term agrees with the proportion \((n-1)!/n!\) of \(n\)-cycles in the symmetric group \(S_n\);
turning this group-theoretic proportion into a distribution of parameters \(b\), and controlling the square-root error, still requires the distribution theorem in the cited paper.
For fixed \(n\), the estimate ensures the existence of such irreducible trinomials over sufficiently large finite fields satisfying the characteristic condition; it does not ensure existence over every small finite field.

\section{Combined Constructions for Arbitrary Degrees}
The binomial criterion for power substitution requires the multiplicative group of an auxiliary field to contain roots of unity of suitable orders.
For an odd prime \(\ell\ne p\), first take \(s=\ord_\ell(q)\), so that \(\ell\mid q^s-1\); a binomial extension of degree \(\ell^e\) can then be constructed over \(K=\F_{q^s}\).
Next take the trace in \(\F_{q^{s\ell^e}}/\F_{q^{\ell^e}}\) to eliminate the auxiliary degree \(s\).
That the trace lies in the target field does not automatically mean that it generates that field; the key to the following proof is to use the distinct powers in the trace sum to rule out membership in a proper subfield.

\begin{theorem}[Construction by an auxiliary extension and trace descent]\leavevmode\phantomsection\label{thm:trace-descent-construction}
Let \(\ell\ne p\) be an odd prime, put \(s=\ord_\ell(q)\), and choose \(a\in K^\times\), where \(K=\F_{q^s}\), such that \(a\) is not an \(\ell\)-th power in \(K\).
For \(e\ge1\), choose \(\beta^{\ell^e}=a\), and put
\[
\gamma=\Tr_{\F_{q^{s\ell^e}}/\F_{q^{\ell^e}}}(\beta)
=\sum_{i=0}^{s-1}\beta^{q^{i\ell^e}}.
\]
Then \(\gamma\) has degree exactly \(\ell^e\) over \(\F_q\).
\end{theorem}
\begin{proof}
Since \(s=\ord_\ell(q)\) divides \(\ell-1\), we have \(\gcd(s,\ell)=1\) and \(\ell\mid q^s-1\).
The condition that \(a\) is not an \(\ell\)-th power in \(K\) is equivalent to its multiplicative order containing the full \(\ell\)-part of \(q^s-1\).
Since \(\ell\) is an odd prime, the binomial criterion shows that \(X^{\ell^e}-a\) is irreducible over \(K\).
Hence
\[
[K(\beta):K]=\ell^e,\qquad K(\beta)=\F_{q^{s\ell^e}}.
\]
The element \(\beta^\ell\) is a root of \(X^{\ell^{e-1}}-a\).
When \(e>1\), the same criterion shows that it has degree \(\ell^{e-1}\) over \(K\);
when \(e=1\), we have \(\beta^\ell=a\in K\), and the conclusion still holds.
Thus
\[
D=K(\beta^\ell)=\F_{q^{s\ell^{e-1}}},\qquad [K(\beta):D]=\ell.
\]

The definition of relative trace ensures that \(\gamma\in\F_{q^{\ell^e}}\).
If its degree over \(k\) were less than \(\ell^e\), the classification of subfields would make that degree a proper divisor of \(\ell^e\),
so \(\gamma\in\F_{q^{\ell^{e-1}}}\subseteq D\).
For \(0\le i<s\), use integer division to write
\[
q^{i\ell^e}=\ell u_i+t_i,\qquad 1\le t_i<\ell.
\]
The remainder is nonzero because \(\ell\nmid q\).
Also, \(\ord_\ell(q^{\ell^e})=s/\gcd(s,\ell^e)=s\),
so \(t_0,\ldots,t_{s-1}\) are pairwise distinct.
If \(\gamma\in D\), the trace identity would become
\[
\sum_{i=0}^{s-1}(\beta^\ell)^{u_i}\beta^{t_i}-\gamma=0.
\]
The left-hand side is the value at \(\beta\) of a polynomial over \(D\);
its degree is less than \(\ell\), and its positive-degree terms have pairwise distinct exponents and nonzero coefficients, so this polynomial is nonzero.
This contradicts \([D(\beta):D]=\ell\). Thus the degree of \(\gamma\) must be \(\ell^e\).
\end{proof}

The trace descent places \(\gamma\) in \(\F_{q^{\ell^e}}\), giving only an upper bound on its orbit length. Distinct residues turn a hypothetical shorter orbit into a nonzero relation of too low a degree, yielding the reverse bound. In terms of~\eqref{v129:eq-orbit-image-product}, this excludes loss of period after taking the trace.

Let \(N=\prod_i\ell_i^{e_i}\). The prime-power degrees are handled by the following constructions.
\begin{itemize}
\item For the characteristic prime \(p\), successive Artin--Schreier equations with nonzero trace construct extensions of degree \(p\); over the prime field, the explicit reciprocal recursion above is also available.
\item For an odd prime \(\ell\ne p\), the trace-descent construction above gives an element of degree \(\ell^{e_i}\).
\item For \(\ell=2\ne p\), use the binomial family from a nonsquare if \(q\equiv1\pmod4\), and the trinomial family above if \(q\equiv3\pmod4\).
\end{itemize}
To construct an extension whose degree is a power of the characteristic, let \(e\ge1\) and \(K_0=k\), choose successively
\(a_j\in K_{j-1}\) with \(\Tr_{K_{j-1}/\F_p}(a_j)\ne0\), and put
\[
\beta_j^p-\beta_j=a_j,\qquad K_j=K_{j-1}(\beta_j)\qquad(1\le j\le e).
\]
The Artin--Schreier criterion gives \([K_j:K_{j-1}]=p\), so \(K_j=\F_{q^{p^j}}\).
The final root satisfies \(\beta_e\notin K_{e-1}\), and the subfield classification (Theorem~\ref{thm:subfields}) shows that every proper intermediate field of \(K_e/k\) is contained in \(K_{e-1}\).
Hence \(k(\beta_e)=K_e\). Taking its Frobenius orbit polynomial directly,
\[
\prod_{i=0}^{p^e-1}(X-\beta_e^{q^i})
\]
gives the required irreducible polynomial of degree \(p^e\).
Finally, composed sums or composed products of irreducible polynomials of distinct prime-power degrees give an irreducible polynomial of degree \(N\).
The computational complexity of the different methods depends on the sizes of the auxiliary fields, the computation of element orders, and the evaluation of trace sums.

\par\Needspace{11\baselineskip}

\begin{example}[Carrying the arbitrary-degree construction through to degree six]\label{v103:ex-composite-six}
We construct degree \(6=2\cdot3\) over \(\F_2\). First take the quadratic component
\[
K=\F_4=\F_2(u),\qquad u^2+u+1=0,
\]
given by an Artin--Schreier equation. For the cubic component, \(s=\ord_3(2)=2\), so the auxiliary field can also be \(K\).
The group \(K^\times\) has order three, and every nonzero cube in it is one; thus \(u\ne1\) is not a cube.
The polynomial \(X^3-u\) is consequently irreducible over \(K\).
Taking \(\beta^3=u\) gives \(K(\beta)=\F_{64}\). Since \(\beta^9=1\) and \(\beta^3\ne1\), the multiplicative order of \(\beta\) is nine.

Now perform trace descent:
\[
\gamma=\Tr_{\F_{64}/\F_8}(\beta)=\beta+\beta^8=\beta+\beta^{-1}.
\]
A direct expansion gives more than an existence statement:
\[
\gamma^3=\beta^3+\beta^{-3}+\beta+\beta^{-1}
=u+u^2+\gamma=1+\gamma.
\]
Thus \(\gamma\) satisfies \(g(Y)=Y^3+Y+1\). This cubic has no root in \(\F_2\), so it is irreducible and gives the required cubic component \(\F_2(\gamma)=\F_8\).

Finally form the composed sum of the components of coprime degrees, taking \(\theta=u+\gamma\).
Write \(A(X)=X^2+X+1\). Then
\[
g(X+u)=X^3+uA(X),\qquad g(X+u^2)=X^3+u^2A(X).
\]
Multiplying eliminates the auxiliary-field coefficients and produces the final output:
\[
\begin{aligned}
F(X)&=g(X+u)g(X+u^2)\\
&=X^6+X^3A(X)+A(X)^2\\
&=\boxed{X^6+X^5+X^3+X^2+1}.
\end{aligned}
\]
Since \(g(\theta+u)=g(\gamma)=0\), we have \(F(\theta)=0\). The degree still needs checking:
\[
\theta^8+\theta=u^2+u=1,\qquad
\theta^4+\theta=\gamma^4+\gamma\ne0.
\]
For the last inequality, \(\gamma^4=\gamma\) would put the degree-three element \(\gamma\) in \(\F_4\cap\F_8=\F_2\), a contradiction.
Thus \(\theta\) lies in neither \(\F_8\) nor \(\F_4\), excluding every proper subfield of \(\F_{64}\).
Its degree is exactly six, and \(F\) is its monic irreducible defining polynomial.
The quadratic component, the cubic binomial in the auxiliary field, trace descent, and the final coefficient elimination have all been given as explicit elements and equations.
\end{example}

\section*{Exercises}
\phantomsection
\addcontentsline{toc}{section}{Exercises}
\begin{exercise}[Tensor-product components and trace duality]\label{ex:supp-composed-six}
Let \(E_m=\F_{q^m}\), \(E_n=\F_{q^n}\), and \(C=E_m\otimes_k E_n\).
\begin{enumerate}
\item Prove that \(\Tr_{C/k}(u\otimes v)=\Tr_{E_m/k}(u)\Tr_{E_n/k}(v)\),
and explain why, under \(C\simeq\prod_{j=1}^{g}\F_{q^\ell}\), it becomes the sum of the field traces of the components, where \(g=\gcd(m,n)\) and \(\ell=\lcm(m,n)\).
\item If \((b_i),(b_i^*)\) and \((c_j),(c_j^*)\) are dual bases of the two extensions, respectively, prove that \((b_i\otimes c_j)\) and \((b_i^*\otimes c_j^*)\) are trace-dual bases of \(C\).
\item Using \(\F_4\otimes_{\F_2}\F_4\), explain why neither the failure of the algebra to be a field nor \(\Tr_{C/k}(1)=0\) prevents its trace pairing from being nondegenerate.
\end{enumerate}
\end{exercise}
\begin{exercise}\label{ex:supp-recursive}
Over \(\F_2\), start from \(F_0=X^2+X+1\) and define \(F_{j+1}=\widehat F_j(X^2+X+1)\). Compute \(F_1\) and explain why every \(F_j\) is irreducible.
\end{exercise}
\begin{exercise}\label{ex:supp-trace-descent}
In the trace-descent construction, explain how the choice \(s=\ord_\ell(q)\) gives both \(\ell\mid q^s-1\) and \(\gcd(s,\ell)=1\). How is each fact used in the proof?
\end{exercise}

\begin{exercise}[Collisions for noncoprime orbits]\label{ex:wu-composed-collisions}
Let \(f=X^2+X+1\in\F_2[X]\). Determine all irreducible factors of \(f\oplus f\) and \(f\odot f\), with their multiplicities. Why do both composed polynomials have repeated factors even though \(f\) is squarefree?
\end{exercise}

\begin{exercise}[Limits of the hypotheses for power-residue symbols]\label{ex:power-residue-boundaries}
\begin{enumerate}
\item In \(\F_2[X]\), take \(P=X^2+X+1\), \(g=X\), and \(d=3\).
Although \(d\mid2^2-1\), show that \(g(\theta)^{(2^2-1)/d}\) depends on the chosen root \(\theta\),
and identify the role of \(d\mid q-1\) in Definition~\ref{def:power-residue-symbol}.
\item In \(\F_5[X]\), take \(F=X(X-1)\), \(g=2\), and \(d=2\).
Prove that \((g/F)_2=1\), but \(g\) is not a square in the quotient ring \(\F_5[X]/(F)\).
Explain why, for a reducible modulus polynomial, a single symbol value cannot replace the power-residue test in each irreducible component.
\end{enumerate}
\end{exercise}

\begin{exercise}[Dynamical irreducibility after extending the base field]\label{v55:ex-stability-base-extension}
Suppose that \((h,\tau)\) is dynamically irreducible over \(\F_q\), with \(\deg h=d\ge2\), and let \(r\ge1\).
Prove that \((h,\tau)\) is dynamically irreducible over \(\F_{q^r}\) if and only if \(\gcd(d,r)=1\).
Explain how this condition relates to the degree \(d^n\) at each level.
\end{exercise}

\begin{exercise}[The same function with different iterates]\label{v55:ex-function-stability}
Let \(f=X^3+X^2+1\) and \(g=X^2+X+1\in\F_2[X]\).
\begin{enumerate}
\item Prove that \(f,g\) define the same constant function on \(\F_2\), but \(f\) is stable.
\item Compute \(g^{\circ2}\) and \(g^{\circ3}\), and prove that \(g\) is not stable.
\item Explain why polynomials cannot be reduced modulo \(X^2-X\) in questions of irreducibility under iteration.
\end{enumerate}
\end{exercise}

\begin{exercise}[The second-level boundary for inverse stability]\label{v55:ex-inverse-quartic-boundary}
In \(\F_2[X]\), put \(\psi=X^4+X+1\) and \(\Psi=1/\psi\).
\begin{enumerate}
\item Prove that \(\psi\) is irreducible and write \(\Psi^{\circ2}\) as a fraction in lowest terms.
\item Verify that its denominator satisfies
\[
\psi^4+\psi^3+1
=(X^8+X^6+X^5+X+1)(X^8+X^6+X^5+X^2+1),
\]
and hence explain why irreducibility at the first level does not ensure inverse stability.
\end{enumerate}
\end{exercise}

\chapter{Polynomial Functions and Fiber Permutations}
\label{chap:functions}
\phantomsection\label{unit:function-fibers}
Functions on a finite field are determined by their values on a finite root set. Starting from low-degree representatives, this chapter studies permutations, value sets, difference preservation, and compatibility between maps on fibers and on their base sets.
Traces and power maps provide two natural projections: a map on the smaller base, together with control within each fiber, constructs a permutation and determines its cycles.
Difference rigidity identifies functions forced to have linear or semilinear structure; Section~\ref{v73:sec-common-cyclic-symmetry} will apply it to common permutation groups of cyclic codes.
The prerequisites here are root-set evaluation, traces and norms, and linearized polynomials, not coding theory.

\section{Function Algebras and Low-Degree Polynomial Representatives}\label{sec:polynomial-functions}
Polynomial representatives of functions in this section are reduced modulo \(X^q-X\). Dynamical irreducibility in Section~\ref{v55:sec-dynamical-irreducibility} concerns the irreducibility of \(h^{\circ n}(X)-\tau\) for a fixed target \(\tau\), and the successive root fibers arising in extension fields.
Exercise~\ref{v55:ex-function-stability} compares two families of iterates with the same function values on the base field.
\begin{theorem}\label{thm:polynomial-functions}
Every function \(h:\F_q\to\F_q\) is represented by a unique polynomial of degree less than \(q\).
Explicitly, we may take
\[
H(X)=\sum_{a\in\F_q}h(a)\bigl(1-(X-a)^{q-1}\bigr).
\]
Two polynomials define the same function if and only if they are congruent modulo \(X^q-X\).
\end{theorem}
\begin{proof}
\leavevmode\phantomsection\label{def:point-indicator}\ignorespaces The polynomial \(\delta_a(X)=1-(X-a)^{q-1}\) is one at \(a\) and zero at every other element of the base field.
Thus the given \(H\) has exactly the prescribed function values and has degree at most \(q-1\).
If the difference of two such polynomials vanishes everywhere, it has \(q\) roots and degree less than \(q\), and must therefore be zero.
An arbitrary polynomial vanishes everywhere if and only if it is divisible by \(\prod_{a\in\F_q}(X-a)=X^q-X\).
\end{proof}

The function algebra under pointwise addition and multiplication is therefore
\[
\F_q[X]/(X^q-X)\cong\prod_{a\in\F_q}\F_q.
\]
The root set \(\F_q\) is a field, and all \(\F_q\)-valued functions on it form a direct product of \(q\) fields under pointwise operations.
The \(\delta_a\) are pairwise orthogonal idempotents whose sum is one.

\section{Permutation Polynomials and Inverse Functions}
\leavevmode\phantomsection\label{def:permutation-polynomial}\ignorespaces A polynomial \(F\in\F_q[X]\) is called a \textbf{permutation polynomial} if it induces a bijection on \(\F_q\).
Its inverse function is represented by a unique polynomial \(G\) of degree less than \(q\), satisfying
\[
G(F(X))\equiv F(G(X))\equiv X\pmod{X^q-X}.
\]
These inverse relations hold modulo \(X^q-X\); they need not hold as identities in the ordinary polynomial ring.

In terms of root sets, \(F\) is a permutation polynomial if and only if
\[
\prod_{a\in\F_q}(T-F(a))=T^q-T,
\]
or, equivalently, if \(\gcd(F(X)-b,X^q-X)\) has degree exactly one for every \(b\in\F_q\).
These conditions say, respectively, that distinct inputs have distinct images and that each fiber has exactly one element.

\begin{example}
For a positive integer \(d\), \(X^d\) is a permutation polynomial over \(\F_q\) if and only if \(\gcd(d,q-1)=1\).
This follows by writing each nonzero element as a power of a primitive element.
Over \(E=\F_{q^n}\), the polynomial \(\sum a_iX^{q^i}\), with \(a_i\in\F_q\), is a permutation polynomial
if and only if \(\gcd(\sum a_iT^i,T^n-1)=1\). These are parallel criteria for the multiplicative-group and Frobenius-operator actions.
\end{example}

For a linearized polynomial over \(E=\F_{q^n}\) with arbitrary coefficients in \(E\), the Dickson matrix of Section~\ref{sec:dickson} applies.
Take its canonical function representative \(L(X)=\sum_{i=0}^{n-1}a_iX^{q^i}\), and write \(r=\operatorname{rank}_E D_L\).
Proposition~\ref{prop:dickson-rank} shows that every nonempty fiber \(L^{-1}(b)\) has \(q^{n-r}\) elements,
since any solution \(x_0\) identifies this fiber as \(x_0+\ker L\). Thus \(L\) is a permutation if and only if \(\det D_L\ne0\);
the linearized polynomial of its inverse can be read from the first row of \(D_L^{-1}\), as in Proposition~\ref{prop:dickson-operations}.
A matrix rank computation therefore determines both fiber sizes and the permutation property.

\subsection{Redei root polynomials, secant directions, and a lacunarity obstruction}\label{v90:subsec-redei}
For a function $f:\F_q\to\F_q$, define the set of \emph{secant directions}
\[
D(f)=\left\{\frac{f(x)-f(y)}{x-y}:x\ne y\right\}.
\]
Directly from the definition,
\[
\boxed{\ f(X)-mX\text{ permutes }\F_q\iff m\notin D(f).\ }
\]
In particular, if the expression is $f(X)+mX$, the condition is $-m\notin D(f)$.
Encode the graph by the \emph{Redei root polynomial}
\[
R_f(U,W)=\prod_{x\in\F_q}(W+Ux+f(x)).
\]
Every $R_f(m,W)$ splits completely over $\F_q$. If $-m\notin D(f)$, its $q$ roots are distinct and $R_f(m,W)=W^q-W$; otherwise it has a repeated root. Thus a question about directions becomes one about multiplicities and coefficients of a root polynomial.
\begin{theorem}[Lacunarity obstruction over a prime field]\label{v90:thm-redei-lacunary}
Let $p$ be prime and $g\in\F_p[X]$, $\deg g<p$. If $F(X)=X^p+g(X)$ splits completely over $\F_p$, with multiplicity, then either $g$ is constant, $g=-X$, or
\[
\deg g\ge\frac{p+1}{2}.
\]
\end{theorem}
\begin{proof}
Put $s=\gcd(F,X^p-X)$, the squarefree polynomial with the same distinct roots as $F$. Since $F-(X^p-X)=X+g$, we have $s\mid X+g$. Write $F=sr$. At each root of $F$, its multiplicity in $r$ is one less than in $F$, hence $r\mid F'=g'$. Therefore
\[
F=sr\mid (X+g)g'.
\]
If the right-hand side is zero, either $g'=0$, forcing $g$ to be constant because $\deg g<p$, or $g=-X$. Otherwise, with $d=\deg g$, one obtains
\[
p=\deg F\le\deg(X+g)+\deg g'\le\max(1,d)+d-1.
\]
Here $d\ge2$, so $p\le2d-1$, as claimed.
\end{proof}
\begin{remark}
The same root set is constrained by the squarefree polynomial $X^p-X$ and by the multiplicities detected by $F'=g'$. For the direction-set classification and its projective-geometric connections, see \cite{HandbookFiniteFields}.
\end{remark}
\subsection{Difference fibers of Gold power maps}\label{v90:subsec-gold}
\begin{proposition}[Difference fibers of a Gold power map]\label{v90:prop-gold}
Let $E=\F_{2^n}$, $f(X)=X^{2^k+1}$, and $d=\gcd(k,n)$. For every $a\in E^\times$ and $b\in E$, the equation
\[
f(x+a)+f(x)=b
\]
has either no solutions or exactly $2^d$ solutions. In particular, if $d=1$ then every difference fiber has at most two points.
\end{proposition}
\begin{proof}
Write $x=ay$. In characteristic two,
\[
f(x+a)+f(x)=a^{2^k+1}(y^{2^k}+y+1).
\]
The equation becomes $y^{2^k}+y=c$. The kernel of this $\F_2$-linear map is $\F_{2^{\gcd(k,n)}}$, of size $2^d$, and each nonempty fiber is its coset.
\end{proof}

\subsection{Formal Derivatives and Lifting Permutations}\label{v105:subsec-polynomial-lifts}
Difference fibers depend only on function values, whereas a formal derivative depends on the chosen polynomial representative.
Indeed, if \(f,H\in\F_q[X]\) and \(g=f+(X^q-X)H\), then \(f\) and \(g\) define the same function on \(\F_q\), but
\[
g'(a)=f'(a)-H(a)\qquad(a\in\F_q).
\]
Reduction modulo \(X^q-X\) therefore preserves function values without necessarily preserving derivative values.
This distinction affects permutation behavior when coefficients are lifted to a finite local ring.

Let \(R\) be a finite commutative local ring with nonzero maximal ideal \(\mathfrak m\) and residue field \(R/\mathfrak m=\F_q\).
\leavevmode\phantomsection\label{v105:def-coefficientwise-lift}\ignorespaces
If the coefficientwise reduction modulo \(\mathfrak m\) of \(F\in R[X]\) is \(f\in\F_q[X]\), then \(F\) is called a
\textbf{coefficientwise lift} of \(f\).
Here the residue polynomial \(f\) is fixed, rather than just the function that it defines on \(\F_q\).

\begin{proposition}[A permutation lifting criterion over finite local rings]\label{v105:prop-permutation-lifting}
Let \(f\in\F_q[X]\) be the reduction modulo \(\mathfrak m\) of \(F\in R[X]\). Then
\begin{equation}\label{v105:eq-permutation-lifting}
F\text{ permutes }R
\quad\Longleftrightarrow\quad
\begin{cases}
f\text{ permutes }\F_q,\\
f'(a)\ne0\quad\text{for every }a\in\F_q.
\end{cases}
\end{equation}
\end{proposition}
\begin{proof}
If \(F\) is surjective, the induced map \(f\) is surjective and hence permutes \(\F_q\).
Let \(s\ge2\) be minimal with \(\mathfrak m^s=0\), and choose \(0\ne h\in\mathfrak m^{s-1}\).
If \(f'(\bar a)=0\), then \(F'(a)\in\mathfrak m\), so \(hF'(a)=h^2=0\). Consequently
\[
 F(a+h)=F(a)+hF'(a)=F(a),
\]
contrary to injectivity.

Conversely, assume both conditions on the right and take \(x\ne y\).
If \(\bar x\ne\bar y\), injectivity of \(f\) gives \(F(x)\ne F(y)\).
If \(\bar x=\bar y\), put \(\delta=y-x\in\mathfrak m\). Polynomial expansion gives
\[
 F(y)-F(x)=\delta\bigl(F'(x)+\delta c\bigr)\qquad(c\in R).
\]
The factor in parentheses has nonzero residue \(f'(\bar x)\), so it is a unit. Since \(\delta\ne0\), the product is nonzero.
Thus \(F\) is injective and hence permutes the finite set \(R\).
\end{proof}

The first condition separates the images of distinct residue classes; the second excludes collisions within each residue class.
For fixed \(f\), the derivative condition is independent of the choices of lifted coefficients within their residue classes.
In particular, it gives the same criterion for every \(\GR(p^k,m)\) with \(k\ge2\).
The nonzero maximal-ideal hypothesis is necessary: \(X^p\) permutes a finite field but has zero derivative.
For the relation with two-element difference fibers in characteristic two, see \cite{BartoliStanica2026Lifts}.

For example, in the non-chain ring
\[
 R=\F_3[\varepsilon_1,\varepsilon_2]/(\varepsilon_1,\varepsilon_2)^2,
\]
\(X+X^3\) is a permutation: writing \(x=a+u\) with \(a\in\F_3\) and \(u\in\mathfrak m\) gives \(x+x^3=2a+u\).
In contrast, \(x^3=a\) collapses each residue class to one point, although \(X^3\) permutes the residue field.

\begin{example}[Function representatives and lifting behavior over the field of order 32]\label{v105:ex-F32-permutation-lifts}
Take \(E=\F_{32}\), put
\[
A(X)=X^2+X^4+X^{16},\qquad G(X)=X^5,
\]
and define
\begin{equation}\label{v105:eq-F32-liftable-polynomial}
P(X)=X+X^4+X^6+X^{10}+X^{12}+X^{20}+X^{24}.
\end{equation}
This seven-term polynomial is the specialization to \(E\) of the construction in \cite[Theorem 2.3]{BartoliStanica2026Lifts}.

First, \(A\) is \(\F_2\)-linear. For \(a\in E\), the identity \(a^{32}=a\) gives
\[
A(a)^4+A(a)=a^8+a^4.
\]
If \(A(a)=0\), then \(a^4\) is fixed by squaring, so \(a^4\in\F_2\).
The fourth-power map is an automorphism of \(E\) restricting to the identity on \(\F_2\), and hence \(a\in\F_2\).
Then \(A(a)=a\), so \(a=0\). Thus \(A\) is invertible.

Next, in the function algebra modulo \(X^{32}-X\),
\[
A(X)^4\equiv X^8+X^{16}+X^2.
\]
Multiplication by \(A(X)\) cancels the two \(X^{18}\) terms, and reduction replaces \(X^{32}\) by \(X\), giving
\begin{equation}\label{v105:eq-F32-linear-input}
A(X)^5\equiv P(X)\pmod{X^{32}-X}.
\end{equation}
Since \(\gcd(5,31)=1\), \(G\) permutes \(E\).
As functions on \(E\), we have \(P=G\circ A\), so \(P\) also permutes \(E\).
Applying Proposition~\ref{v90:prop-gold} to \(G(X)=X^{2^2+1}\) and using \(\gcd(2,5)=1\),
we see that every nonempty difference fiber of \(G\) in a nonzero direction has exactly two elements.
For \(a\ne0\),
\[
P(x+a)+P(x)=G(A(x)+A(a))+G(A(x)).
\]
As \(A\) is invertible, the substitution \(y=A(x)\) ranges over \(E\), and \(A(a)\ne0\).
Hence every nonempty difference fiber of \(P\) in a nonzero direction also has exactly two elements.

The formal derivatives, however, are
\[
P'(X)=1,\qquad G'(X)=X^4,\qquad\bigl(A(X)^5\bigr)'=0.
\]
It follows that, for every \(k\ge2\), every coefficientwise lift of \(P\) permutes \(\GR(2^k,5)\),
whereas no coefficientwise lift of either \(G\) or \(A^5\) permutes that ring.
In particular, interpreting each coefficient \(1\) in \eqref{v105:eq-F32-liftable-polynomial} as the identity element of the ring gives such a permutation lift.

There are two distinct comparisons here.
The polynomials \(A^5\) and \(P\) define the same function on \(E\), but the former has ordinary degree 80,
whereas the latter, of degree 24, is its unique representative of degree less than 32.
The polynomials \(G\) and \(P\), on the other hand, are both already reduced to degree less than 32 and are related by the invertible \(\F_2\)-linear input map \(A\).
This input change preserves difference fiber sizes but does not preserve the coefficientwise lifting behavior of the reduced representative.
\end{example}

\section{Difference Preservation and Algebraic Rigidity}
\label{v73:sec-difference-rigidity}
How much of a function is determined by its local differences? If parallel difference quotients are independent of the base point, the differences recover linearity over the prime field.
Preservation of the multiplicative cosets of nonzero differences instead forces rigidity through low-degree function representatives and Frobenius automorphisms.
A multiplicative subgroup constraint on second divided differences recovers ordinary quadratics through compatibility between base points; constant values for a fixed shape recover linearity over its generated subfield.
Over an integer residue ring, the same question is expressed by cyclotomic difference orbits and reduced by the Chinese remainder theorem to prime-power moduli.
These conditions link constraints on pairs or triples to global algebraic form. Over residue rings, however, certain additive blocks also allow extra non-affine permutations, so rigidity and its exceptions must both be identified.

\subsection{Parallel Difference Quotients and Linearized Maps}
\phantomsection\label{v72:subsec-parallel-difference}
In a finite field, the difference quotient for a fixed nonzero displacement usually depends on the base point.
When it depends only on the displacement, the global additive structure can be recovered from local differences.

\begin{theorem}[Parallel difference quotients characterize linearized polynomials]\label{v72:thm-parallel-difference}
Let \(q=p^n\), where \(p\) is prime. For a function \(f:\F_q\to\F_q\), the following are equivalent:
\begin{enumerate}
\item For every \(t\in\F_q^\times\) and all \(x,y\in\F_q\),
\begin{equation}\label{v72:eq-parallel-difference}
\frac{f(x+t)-f(x)}{t}=\frac{f(y+t)-f(y)}{t}.
\end{equation}
\item There are unique \(b,a_0,\ldots,a_{n-1}\in\F_q\) such that, as functions on \(\F_q\),
\[
f(X)=b+\sum_{i=0}^{n-1}a_iX^{p^i}.
\]
\end{enumerate}
Neither bijectivity nor nonvanishing of the difference quotient is assumed.
\end{theorem}
\begin{proof}
Assume~\eqref{v72:eq-parallel-difference}. Setting \(y=0\) and canceling \(t\ne0\) gives
\[
f(x+t)-f(x)=f(t)-f(0)\qquad(t\ne0).
\]
The same identity holds when \(t=0\). Put \(L(x)=f(x)-f(0)\). Then
\[
L(x+t)=L(x)+L(t)\qquad(x,t\in\F_q).
\]
Thus \(L\) is additive and hence \(\F_p\)-linear, since \(L(cx)=cL(x)\) for all \(c\in\F_p\).
Apply Theorem~\ref{thm:linear-map-representation} of Chapter~\ref{chap:linearized} with \(k=\F_p\) and \(E=\F_q\).
It yields the unique representation \(L(X)=\sum_{i=0}^{n-1}a_iX^{p^i}\); take \(b=f(0)\).
Conversely, \((x+t)^{p^i}=x^{p^i}+t^{p^i}\) implies
\[
\frac{f(x+t)-f(x)}{t}=\sum_{i=0}^{n-1}a_it^{p^i-1},
\]
which is independent of \(x\), proving the first condition.
\end{proof}

\begin{example}[Parallel quotients without a permutation]
Over \(\F_9=\F_3(u)\), where \(u^2=-1\), take \(f(X)=X+X^3\). For every \(t\ne0\),
\[
\frac{f(x+t)-f(x)}{t}=1+t^2,
\]
independent of \(x\). However, \(f(u)=0=f(0)\), so \(f\) is not a permutation.
Thus parallel difference-quotient invariance characterizes all \(\F_p\)-affine maps, not merely the bijective ones.
\end{example}

\subsection{Multiplicative Coset Preservation and Rigidity}
\phantomsection\label{v72:subsec-coset-rigidity}
The preceding theorem recovers all \(\F_p\)-affine maps from parallel difference quotients. We now impose a different condition---preservation of multiplicative cosets of nonzero differences---which yields a stronger rigidity conclusion without assuming parallel difference-quotient invariance.
Let \(H\) be a proper subgroup of \(\F_q^\times\). A function preserves the multiplicative \(H\)-coset of every nonzero difference if
\[
(f(x)-f(y))H=(x-y)H\qquad(x\ne y).
\]
Equivalently, every difference quotient belongs to \(H\).
The condition preserves the image of each difference in the quotient group \(\F_q^\times/H\).
These multiplicative cosets are distinct from the \(q\)-cyclotomic cosets of exponents in Chapter~\ref{chap:cyclotomic}.
The classical Carlitz--McConnel theorem shows that this condition determines the function up to a scalar, a translation, and a power of Frobenius
\cite{Carlitz1960Permutations,McConnel1963PseudoOrdered}. The proof uses the function algebra of Section~\ref{sec:polynomial-functions} and follows Lenstra~\cite[Section~2]{Lenstra1990Automorphisms}.

\begin{theorem}[Carlitz--McConnel]\label{v72:thm-carlitz-mcconnel}
Let \(q=p^n\), where \(p\) is prime, let \(H\) be a proper subgroup of \(\F_q^\times\), and put
\(h=|H|\) and \(d=(q-1)/h\). A function \(f:\F_q\to\F_q\) satisfies
\begin{equation}\label{v72:eq-coset-difference}
\frac{f(x)-f(y)}{x-y}\in H\qquad(x\ne y)
\end{equation}
if and only if
\[
f(x)=ax^{p^j}+b,\qquad
a\in H,\quad b\in\F_q,\quad 0\le j<n,\quad d\mid p^j-1.
\]
Every such function is a permutation of \(\F_q\), and the parameters \(a,b,j\) are unique.
\end{theorem}
\begin{proof}
The condition separates the images of distinct elements, so \(f\) is a permutation.
Subtracting \(f(0)\) does not change any difference quotient; assume first that \(f(0)=0\).
Since \(h\mid q-1\) and \(H\) is proper, we have \(p\nmid h\) and \(2h\le q-1\).
The subgroup \(H\) is exactly the root set of \(X^h-1\) in \(\F_q\), so
\[
(f(x)-f(y))^h=(x-y)^h,\qquad f(x)^h=x^h.
\]
The second identity alone does not readily control \(f\). The first compares all base points, suggesting that we consider all translates of \(X^h\) together.

In the function algebra \(A=\F_q[X]/(X^q-X)\), write \(X\) for the identity function.
Translations act by \(T_cg(X)=g(X-c)\), with \(T_cT_b=T_{c+b}\). Put
\begin{equation}\label{v74:eq-translation-power-space}
V=\operatorname{span}_{\F_q}\{(X-c)^h:c\in\F_q\}.
\end{equation}
This is the span of the translation orbit of \(X^h\); each of its elements has a polynomial representative of degree at most \(h\).
We recover \(X\) and \(X^{h-1}\) from these translates because their product is precisely \(X^h\).
Cyclicity of the multiplicative group and a geometric sum give, for every positive integer \(m\),
\(\sum_{c\in\F_q}c^m=-1\) if \(q-1\mid m\), and zero otherwise; the same formula will be used in Lemma~\ref{lem:power-sums}.
Binomial expansion now gives
\[
\begin{aligned}
\sum_{c\in\F_q}c^{q-h}(X-c)^h&=(-1)^h hX,\\
\sum_{c\in\F_q}c^{q-2}(X-c)^h&=hX^{h-1}.
\end{aligned}
\]
Indeed, the term \(\binom hkX^{h-k}(-c)^k\) gives exponents \(q-h+k\) and \(q-2+k\), respectively.
All these exponents lie strictly between \(0\) and \(2(q-1)\), so only exponent \(q-1\) contributes; it occurs at \(k=h-1\) and \(k=1\), respectively.
As \(h\ne0\) in \(\F_q\), we obtain \(X,X^{h-1}\in V\). This also covers \(h=1\), with \(X^{h-1}\) interpreted as the constant function one.

Composition with the permutation \(f\) is an algebra automorphism \(U_f:g\mapsto g\circ f\) of \(A\).
For \(b=f(c)\),
\[
U_f\bigl((X-b)^h\bigr)=(f(X)-f(c))^h=(X-c)^h.
\]
Because \(f\) is bijective, \(U_f\) permutes these generators of \(V\), so \(U_f(V)=V\).
In particular, \(U_f(X)\) and \(U_f(X^{h-1})\) have representatives \(P,Q\) of degree at most \(h\), with
\(P(x)=f(x)\) and \(Q(x)=f(x)^{h-1}\); the latter is a pointwise power, not an iterate.
Thus \(P(x)Q(x)=x^h\) for all \(x\in\F_q\).
Since \(\deg(PQ-X^h)\le2h<q\), the root bound turns this equality of functions into the polynomial identity
\[
P(X)Q(X)=X^h.
\]
Unique factorization and \(P(0)=0\) force \(P(X)=aX^k\), where \(a\ne0\) and \(1\le k\le h\).

Translation compatibility now determines the exponent. For every \(y\in\F_q\), the map
\(f_y(x)=f(x+y)-f(y)\) satisfies the original condition and fixes zero.
The same argument supplies a monomial representative \(a_yX^{k_y}\), with \(1\le k_y\le h\).
Both sides have degree less than \(q\), so equality of their values again yields a polynomial identity:
\[
a\bigl((X+y)^k-y^k\bigr)=a_yX^{k_y}.
\]
Comparing degrees and leading coefficients gives \(k_y=k\) and \(a_y=a\), hence
\((X+y)^k=X^k+y^k\).
The map \(x\mapsto x^k\) therefore preserves addition, multiplication, and one, and is bijective; it is a field automorphism.
Apply Theorem~\ref{thm:automorphisms} over the prime field \(\F_p\) to obtain \(x\mapsto x^{p^j}\), where \(0\le j<n\).
Uniqueness of representatives of degree less than \(q\) also gives \(k=p^j\).

Restoring the constant term yields \(f(x)=ax^{p^j}+b\).
The pair \((x,y)=(1,0)\) gives \(a\in H\).
For a primitive element \(\zeta\), the condition \(a\zeta^{p^j-1}\in H=\langle\zeta^d\rangle\) gives \(d\mid p^j-1\).
Conversely, these parameters ensure
\[
\frac{f(x)-f(y)}{x-y}=a(x-y)^{p^j-1}\in H.
\]
Finally, \(b=f(0)\), \(a=f(1)-f(0)\), and the \(n\) Frobenius powers are distinct, so the parameters are unique.
\end{proof}

The proof uses two distinct conversions: coset preservation leaves the translation-power space invariant, and low-degree representatives with their product relation force a monomial; changing the base point then forces that monomial to preserve addition.
The inequality \(2h<q\) turns pointwise equality into a polynomial identity. Additive translations, multiplicative constraints, and the root bound together identify Frobenius; bijectivity alone does not do so.

\begin{corollary}[Counting the permitted affine Frobenius permutations]\label{v72:cor-carlitz-mcconnel-count}
With the notation of Theorem~\ref{v72:thm-carlitz-mcconnel}, put \(r=\ord_d(p)\).
There are exactly
\[
qh\frac nr
\]
functions satisfying the difference-quotient condition. Those fixing \(0\) and \(1\) are exactly \(x\mapsto x^{p^{tr}}\), where \(0\le t<n/r\).
\end{corollary}
\begin{proof}
Since \(d\mid p^n-1\), we have \(r\mid n\), and the allowed exponents are \(j=0,r,\ldots,n-r\).
By uniqueness of the parameters, there are \(q\) choices for \(b\), \(h\) for \(a\), and \(n/r\) for \(j\).
The normalization fixes \(b=0\) and \(a=1\).
\end{proof}

\begin{example}[Square classes of differences over \(\F_9\)]\label{v72:ex-carlitz-mcconnel-f9}
Let \(H=(\F_9^\times)^2\). Then \(h=4\), \(d=2\), and \(\ord_2(3)=1\).
The functions preserving the square class of every nonzero difference are exactly
\[
f(x)=ax+b\quad\text{or}\quad f(x)=ax^3+b,
\qquad a\in H,\quad b\in\F_9.
\]
There are \(9\cdot4\cdot2=72\) such functions.
Denote the quadratic generator of Example~\ref{v100:ex-three-generation-nine} by \(u\). Its Frobenius action is already determined by \(u\mapsto-u\); the difference-quotient condition restricts all functions to these two families.
Thus the normalized functions are the identity and this nontrivial conjugation.
\end{example}

\begin{remark}[Permuting cosets and the hidden subfield]\label{v74:rem-hidden-subfield}
Fixing each multiplicative coset differs from allowing a global permutation of the cosets.
Continue to assume that \(H<\F_q^\times\) is a proper subgroup. Write \(D=H\), \(G=\F_q^\times/D\), let \(\pi:\F_q^\times\to G\) be the quotient map, and put \(K=\F_p(D)\).
Lenstra's Theorem~2 \cite{Lenstra1990Automorphisms} states that there exists a single permutation \(\kappa\) of the quotient group for which a permutation \(f\) satisfies
\[
\pi(f(x)-f(y))=\kappa\bigl(\pi(x-y)\bigr)\qquad(x\ne y)
\]
if and only if \(f(x)=\beta(x)+b\), where \(\beta\) is a \(K\)-semilinear automorphism.
As in the definition on page~\pageref{def:semilinear-map}, this means that \(\beta\) is an additive bijection and
\(\beta(cx)=\gamma(c)\beta(x)\) for some \(\gamma\in\Aut(K)\), all \(c\in K\), and all \(x\in\F_q\).
The same \(\kappa\) must work for every pair; choosing a different relabeling for each pair is not the same condition.

The subfield arises from scalar multiplication: if \(M_c(x)=cx\), the operator algebra
\(\F_p[M_c:c\in D]\) consists exactly of the \(M_c\) with \(c\in K\), and is therefore isomorphic to \(K\). Thus
\[
D\longmapsto K=\F_p(D)\longmapsto K\text{-semilinear symmetries}.
\]
When \(K=\F_q\), this vector space has dimension one, forcing \(\beta(x)=ax^{p^j}\).
For a proper subfield \(K\), semilinearity is instead understood on a higher-dimensional \(K\)-space; a single Frobenius monomial cannot be assumed in advance.
The theorem thus identifies a scalar field from the multiplicative subgroup: freedom under coset relabeling is governed by this subfield and the relative dimension.
\end{remark}

\subsection{A Few Failed Difference Quotients and Local Stability}\label{v116:subsec-stability}
The exact difference-quotient condition determines all permitted permutations. If only a few pairs fail, one can instead compare the permutation with these models point by point.
Fix \(q=p^e\), a proper subgroup \(H<\F_q^\times\), and put \(h=|H|\), \(d=(q-1)/h\). Let
\[
\Gamma_H=\{x\mapsto ax^{p^j}+b:a\in H,\ b\in\F_q,\ 0\le j<e,\ d\mid p^j-1\}.
\]
For a permutation \(f\), define
\begin{equation}\label{v116:def-stability}
\begin{aligned}
\mathcal E_H(f)&=\#\{(x,y)\in\F_q^2:x\ne y,\ (f(x)-f(y))/(x-y)\notin H\},\\
D_H(f)&=\min_{g\in\Gamma_H}\#\{x\in\F_q:f(x)\ne g(x)\}.
\end{aligned}
\end{equation}
Pairs are ordered. The Carlitz--McConnel theorem says that \(\mathcal E_H(f)=0\) if and only if \(D_H(f)=0\).
If \(f\) differs from a model at \(r\) points, only pairs involving those points can fail. Hence
\[
\mathcal E_H(f)\le2(q-1)D_H(f).
\]

\begin{example}[Exchanging two points]\label{v116:ex-transposition}
Let \(f\) exchange \(0,1\) and fix every other point. Then
\[
\mathcal E_H(f)=4(q-h-1)+2\mathbf1_{\{-1\notin H\}}.
\]
Indeed, for \(y\notin\{0,1\}\), the two quotients are \((y-1)/y\) and \(y/(y-1)\).
The former ranges over \(\F_q\setminus\{0,1\}\), and the latter is its inverse. Each group therefore has \(q-h-1\) failures; including reversed pairs gives the first term.
The remaining pairs \((0,1),(1,0)\) both have quotient \(-1\).
For example, \(q=7\) and \(H=\{1,2,4\}\) give fourteen bad pairs.
Exchanging \(1,2\) instead and fixing the other points gives \(f(x)/x\in H\) for every \(x\ne0\), yet some pairs still fail. Conditions at one origin do not replace the full difference-quotient condition.
\end{example}

\begin{proposition}[A lower bound when an approximate model is known]\label{v116:prop-local-stability}
If \(f\) is a permutation and \(g\in\Gamma_H\) agrees with \(f\) at exactly \(q-r\) points, then
\begin{equation}\label{v116:eq-local-stability}
\mathcal E_H(f)\ge2r\max\{0,q-r-h+1\}.
\end{equation}
In particular, if \(r\le(q-h+1)/2\), then \(\mathcal E_H(f)\ge r(q-h+1)\).
\end{proposition}
\begin{proof}
Permutations satisfying the exact condition are closed under composition and inversion, since the corresponding difference quotients multiply within \(H\).
Replacing \(f\) by \(g^{-1}\circ f\) preserves \(\mathcal E_H\) and the number of disagreements with the model. We may therefore take \(g\) to be the identity.
Let \(G=\{y:f(y)=y\}\), so \(|G|=q-r\). Fix \(x\notin G\) and write \(a=f(x)\ne x\). The map
\[
y\longmapsto\frac{a-y}{x-y}
\]
is a bijection from \(\F_q\setminus\{x\}\) to \(\F_q\setminus\{1\}\), of whose values exactly \(h-1\) belong to \(H\).
Thus \(x\) forms a bad pair with at least \(\max\{0,q-r-h+1\}\) points of \(G\).
Sum over the \(r\) disagreement points and include both directions.
\end{proof}

A transposition affects only \(O(q)\) pairs, so a small failure proportion does not force membership in \(\Gamma_H\). Once a sufficiently close model has been found, however, the lower bound converts errors on pairs into errors at points.
For fixed index \(d\), do there exist constants \(\varepsilon_d>0\) and \(0<\eta_d<(1-1/d)/2\), independent of \(q\), such that over sufficiently large prime fields
\[
\mathcal E_H(f)\le\varepsilon_dq^2\quad\Longrightarrow\quad D_H(f)\le\eta_dq?
\]
Such a global gap would imply \(\mathcal E_H(f)\ge c_dqD_H(f)\), as Exercise~\ref{v116:ex-stability-gap} shows.
The invariant space of translated powers in the exact proof provides an entry point. With errors allowed, one must still recover a single global model from approximate invariance before applying the local bound.

\subsection{Three-Point Difference Determinants and Second-Order Rigidity}
\label{v76:subsec-second-difference-rigidity}
A two-point difference quotient compares the multiplicative cosets of two nonzero differences. Comparing two secant slopes leads to a three-point condition.
Throughout this subsection, let \(q=p^n\ge3\). For pairwise distinct \(x,y,z\in\F_q\), define the \textbf{three-point difference determinant}
\begin{equation}\label{v76:def-three-point-determinant}
\mathcal A_f(x,y,z)=
\det\begin{pmatrix}
x-y&f(x)-f(y)\\
z-y&f(z)-f(y)
\end{pmatrix}.
\end{equation}
It vanishes exactly when the three graph points are collinear. Put \(V(x,y,z)=(x-y)(z-y)(z-x)\), and define the \textbf{second divided difference}
\begin{equation}\label{v76:def-second-divided-difference}
\begin{aligned}
\mathcal K_f(x,y,z)
&=\frac{\mathcal A_f(x,y,z)}{V(x,y,z)}\\
&=\frac{1}{z-x}
\left(\frac{f(z)-f(y)}{z-y}-\frac{f(x)-f(y)}{x-y}\right).
\end{aligned}
\end{equation}
By three-point Lagrange interpolation, \(\mathcal K_f(x,y,z)\) is the coefficient of \(X^2\) in the interpolating polynomial, so it is symmetric in the nodes.
In particular, adding an affine polynomial to \(f\) leaves it unchanged, and
\[
\mathcal K_{aX^2+bX+c}(x,y,z)=a,\qquad
\mathcal K_{X^3}(x,y,z)=x+y+z.
\]
Here and below, polynomial expressions describe functions on \(\F_q\); whenever ordinary degrees are compared, we use the unique representative of degree less than \(q\).
A multiplicative subgroup constraint on the raw determinant generally differs from one on the divided difference, since the normalizing factor \(V(x,y,z)\) varies with the nodes.

\Needspace{9\baselineskip}
\begin{theorem}[Multiplicative subgroup rigidity for second divided differences]\label{v76:thm-second-difference-rigidity}
Let \(H\) be a proper subgroup of \(\F_q^\times\). A function \(f:\F_q\to\F_q\) satisfies
\begin{equation}\label{v76:eq-all-triples-subgroup}
\mathcal K_f(x,y,z)\in H
\qquad(x,y,z\text{ pairwise distinct})
\end{equation}
if and only if
\[
f(X)=aX^2+bX+c,\qquad a\in H,\quad b,c\in\F_q.
\]
The parameters \(a,b,c\) are unique, so there are \(q^2|H|\) such functions. These functions need not be permutations.
\end{theorem}
\begin{proof}
Let \(F\in\F_q[X]\) be the representative of \(f\) of degree less than \(q\). Fix \(u\in\F_q\), and put
\[
g_u(T)=\frac{F(T)-F(u)}{T-u}\in\F_q[T],\qquad \deg g_u\le q-2.
\]
This is polynomial division, so \(g_u\) is also defined at \(u\), with \(g_u(u)=F'(u)\). For distinct \(v,w\ne u\), equation~\eqref{v76:def-second-divided-difference} gives
\begin{equation}\label{v76:eq-punctured-difference}
\frac{g_u(v)-g_u(w)}{v-w}=\mathcal K_f(u,v,w)\in H.
\end{equation}
Thus \(g_u\) is injective on \(\F_q\setminus\{u\}\), and its image omits a unique element \(c_u\).
The finite-field power-sum formula shows that the sum of the values of every polynomial of degree at most \(q-2\) is zero; a constant term also contributes \(q=0\).
Moreover, \(\sum_{a\in\F_q}a=0\), since \(q>2\). Hence
\[
0=\sum_{t\in\F_q}g_u(t)
=g_u(u)+\sum_{a\in\F_q\setminus\{c_u\}}a
=g_u(u)-c_u.
\]
It follows that \(g_u(u)=c_u\), so \(g_u\) permutes all of \(\F_q\).

The difference-quotient condition must still be extended to pairs involving \(u\). Fix \(v\ne u\), and multiply equation~\eqref{v76:eq-punctured-difference} over \(w\notin\{u,v\}\).
Since \(g_u\) is a permutation, the numerator and denominator each run through all nonzero field elements with one element omitted. Therefore
\[
\prod_{w\notin\{u,v\}}
\frac{g_u(w)-g_u(v)}{w-v}
=\frac{u-v}{g_u(u)-g_u(v)}\in H.
\]
Taking inverses gives \((g_u(u)-g_u(v))/(u-v)\in H\). Together with equation~\eqref{v76:eq-punctured-difference}, this proves the Carlitz--McConnel condition for every pair of distinct elements.
Theorem~\ref{v72:thm-carlitz-mcconnel} now yields \(a_u\in H\), \(b_u\in\F_q\), and \(0\le j_u<n\) such that
\[
g_u(T)=a_uT^{p^{j_u}}+b_u.
\]
Both sides have degree less than \(q\), so this is an identity of ordinary polynomials, not merely an equality of functions. Consequently,
\begin{equation}\label{v76:eq-basepoint-polynomial}
F(T)=F(u)+(T-u)(a_uT^{p^{j_u}}+b_u).
\end{equation}
The right-hand side has degree \(p^{j_u}+1\) and leading coefficient \(a_u\). Thus \(j_u=j\) and \(a_u=a\) are independent of \(u\).
If \(j\ge1\), then \(p^j\ge2\), and the coefficient of \(T^{p^j}\) in equation~\eqref{v76:eq-basepoint-polynomial} is \(-au\).
Taking \(u=0\) and \(u=1\) gives different coefficients for the same polynomial \(F\), a contradiction.
Hence \(j=0\), and \(F=aT^2+bT+c\), with \(a\in H\).
Conversely, the second divided difference of such a quadratic is constantly \(a\), as required.
Finally, \(2<q\) makes the representative of degree at most two unique, which also proves the count.
\end{proof}

The proof reduces a three-point condition to a one-point-deleted version of a two-point condition. Power sums complete the permutation, a product identity completes the difference-quotient condition, and compatibility between base points excludes nontrivial Frobenius exponents.
The required first-order rigidity was proved above by Lenstra's function-algebra method; no classification theorem for finite projective planes is needed.

If only one relative position of three points is prescribed, the subfield generated by the position parameter governs the answer.
Fix \(\lambda\in\F_q\setminus\{0,1\}\), called the \textbf{shape parameter} of \(x,x+t,x+\lambda t\), and write
\begin{equation}\label{v76:def-fixed-shape-difference}
\begin{aligned}
K_{\lambda,f}(x,t)
&=\mathcal K_f(x,x+t,x+\lambda t)\\
&=\frac{f(x+\lambda t)-\lambda f(x+t)+(\lambda-1)f(x)}
{\lambda(\lambda-1)t^2},\qquad t\ne0.
\end{aligned}
\end{equation}

\Needspace{9\baselineskip}
\begin{theorem}[Constant divided differences of one shape and subfield linearity]\label{v76:thm-fixed-shape-constant}
Let \(K=\F_p(\lambda)=\F_{p^r}\), and fix \(a\in\F_q\). Then
\[
K_{\lambda,f}(x,t)=a\qquad(x\in\F_q,\ t\ne0)
\]
if and only if
\begin{equation}\label{v76:eq-fixed-shape-classification}
f(X)=aX^2+L(X)+b,\qquad
L(X)=\sum_{i=0}^{n/r-1}\ell_iX^{p^{ri}},\qquad \ell_i,b\in\F_q.
\end{equation}
Here \(L\) ranges over all \(K\)-linear maps, not necessarily invertible.
In particular, when \(K=\F_q\), the solutions are precisely \(aX^2+b_1X+b_0\).
\end{theorem}
\begin{proof}
The divided difference \(K_{\lambda,X^2}\) is constantly one. Put \(h(x)=f(x)-ax^2-f(0)\). Then \(h(0)=0\), and the condition is equivalent to
\[
h(x+\lambda t)=(1-\lambda)h(x)+\lambda h(x+t).
\]
This identity also holds for \(t=0\). Setting \(y=x+t\) gives
\begin{equation}\label{v76:eq-weighted-affine-identity}
h((1-\lambda)x+\lambda y)
=(1-\lambda)h(x)+\lambda h(y)\qquad(x,y\in\F_q).
\end{equation}
Taking \(x=0\) and \(y=0\), respectively, yields
\(h(\lambda y)=\lambda h(y)\) and \(h((1-\lambda)x)=(1-\lambda)h(x)\).
For arbitrary \(u,v\), substitute \(x=u/(1-\lambda)\) and \(y=v/\lambda\) to obtain
\[
h(u+v)=h(u)+h(v).
\]
Thus \(h\) is \(\F_p\)-linear and commutes with multiplication by \(\lambda\).
It therefore commutes with all scalar multiplications from \(\F_p[\lambda]=K\), so it is \(K\)-linear.
Theorem~\ref{thm:linear-map-representation} gives equation~\eqref{v76:eq-fixed-shape-classification}, with \(b=f(0)\).
Conversely, a \(K\)-linear map satisfies equation~\eqref{v76:eq-weighted-affine-identity}, so its divided difference of this shape is zero. Adding \(aX^2+b\) gives the required constant.
\end{proof}

For a prescribed constant \(a\in H\), the constant-value condition is stronger than requiring membership in \(H\). In particular, a subgroup condition for one fixed shape does not meet the hypothesis of the preceding all-triples theorem.

\begin{example}[A fixed shape admitting a nonquadratic function]\label{v76:ex-fixed-shape-f16}
In \(\F_{16}\), take \(H=(\F_{16}^\times)^3\) and \(\omega\in\F_4\setminus\F_2\), with \(\omega^2+\omega+1=0\).
Let \(f(X)=X^8\). The Frobenius identity in characteristic two gives
\[
K_{\omega,f}(x,t)
=\frac{\omega^8-\omega}{\omega(\omega-1)}t^6=t^6\in H.
\]
The coefficient is one because \(\omega^3=1\), \(\omega^8=\omega^2\), and \(\omega^2+\omega=1\).
As \(t\) runs through \(\F_{16}^\times\), the values \(t^6\) run through \(H\), so these divided differences are not constant.
Since \(X^8\) has degree less than \(16\), it cannot define the same function as a polynomial of degree at most two.
Exercise~\ref{v76:ex-frobenius-shape-factor} gives the general exponent criterion behind this example; Exercise~\ref{v76:ex-multishape-recovery} distinguishes constant-value recovery from subgroup-value recovery.
\end{example}

\begin{remark}[Sidon sets and non-Sidon four-term relations]\label{v76:rem-sidon-shape-relations}
In characteristic two, a fixed shape also connects three-point constraints to four-term additive relations.
A subset \(S\subseteq\F_{2^n}\) is a \textbf{Sidon set} if the sum of an unordered pair of distinct elements determines that pair uniquely.
Equivalently, no four pairwise distinct elements of \(S\) sum to zero. For this additive property of multiplicative subgroups, see \cite{CarletMesnager2022Sidon}.
Define the \textbf{shape-specific four-term relation count}
\begin{equation}\label{v76:def-shape-four-term-count}
N_\lambda(H)=\#\{(r,z)\in\F_q\times\F_q^\times:
 r,r+z,r+\lambda z,r+(1+\lambda)z\in H\}.
\end{equation}
Consider \(x,x+t,x+\lambda t,x+(1+\lambda)t\), and let \(\Delta\) be the sum of the four function values.
Each of their four triples is a translate of the original triple. Expanding equation~\eqref{v76:def-second-divided-difference} shows that their divided differences can be ordered as
\[
k,\quad k+z,\quad k+\lambda z,\quad k+(1+\lambda)z,
\qquad z=\frac{\Delta}{\lambda(1+\lambda)t^2}.
\]
If \(K_{\lambda,f}(x,t)\in H\) everywhere and \(N_\lambda(H)=0\), then \(z=0\), giving
\begin{equation}\label{v76:eq-shape-parallelogram}
f(x)+f(x+t)+f(x+\lambda t)+f(x+(1+\lambda)t)=0.
\end{equation}
The Sidon property makes \(N_\lambda(H)=0\) for every shape, but a zero count for a single shape does not require \(H\) to be Sidon.
For example, take \(H=\F_8^\times\) inside \(\F_{64}\). From an \(\F_2\)-basis \(1,u,v\) of \(\F_8\), the four elements \(1,u,v,1+u+v\) show that \(H\) is not Sidon.
Nevertheless, \(N_\lambda(H)=0\) whenever \(\lambda\notin\F_8\): the first two entries of a proposed relation would force \(z\in\F_8^\times\), and the first and third would then give \(\lambda z\in\F_8\), a contradiction.
Equation~\eqref{v76:eq-shape-parallelogram} asserts only a four-point identity in the specified directions. It does not by itself imply global additivity or quadratic rigidity; compatibility between shapes must still be examined.
\end{remark}

\subsection{Difference Orbits and Local Rigidity at Prime Powers}
\label{v73:subsec-local-rigidity}
Every nonzero difference in a field can be used as a divisor. Over a residue ring, a nonzero difference need not be a unit, so we retain the orbit condition rather than form a quotient.
This connects multiplicative coset preservation over fields with the exponent orbits of Chapter~\ref{chap:cyclotomic}.

\begin{definition}[Preservation of cyclotomic difference orbits]\label{v73:def-difference-orbit}
Let \(n\ge2\), \(R=\Z/n\Z\), and \(K\le R^\times\). For \(z\in R\), put \(Kz=\{uz:u\in K\}\).
A function \(f:R\to R\) \textbf{preserves cyclotomic difference orbits} if
\begin{equation}\label{v73:eq-difference-orbit}
f(x)-f(y)\in K(x-y)\qquad(x,y\in R).
\end{equation}
Write \(\mathcal R_K(n)\) for the set of these functions, and put
\[
\operatorname{Aff}_K(n)=\{x\mapsto ux+b:u\in K,\ b\in R\}.
\]
\end{definition}
Here \(K\) need not be generated by a field size; when \(K=\langle q\rangle\), the sets \(Kz\) are the \(q\)-cyclotomic cosets modulo \(n\).
Multiplication by a unit cannot send a nonzero element to zero, so condition~\eqref{v73:eq-difference-orbit} makes \(f\) a permutation automatically.
It gives the equality \(K(f(x)-f(y))=K(x-y)\), so these permutations are closed under composition and inversion.
Thus \(\mathcal R_K(n)\) is a group containing \(\operatorname{Aff}_K(n)\).

If \(n=p\) is prime and \(K\ne\F_p^\times\), the Carlitz--McConnel theorem gives
\(\mathcal R_K(p)=\operatorname{Aff}_K(p)\), since the prime field has no nontrivial automorphism.
If \(K=\F_p^\times\), every permutation satisfies the condition; all of them are affine precisely when \(p=2,3\).
Indeed, the affine group has \(p(p-1)\) elements, equal to \(2!\) and \(3!\) in these two cases.
For \(p\ge5\), a transposition fixing every other point is not affine.

At a prime-power modulus we must also compare congruence layers. The following lemmas show that excluding certain deepest-layer units forces a common local multiplier.

\begin{lemma}[Local rigidity at odd prime powers]\label{v73:lem-odd-local-rigidity}
Let \(p\) be an odd prime, \(a\ge2\), \(R=\Z/p^a\Z\), and \(L\le R^\times\). If
\[
1+p^{a-1}R\not\subseteq L,
\]
then every \(f\in\mathcal R_L(p^a)\) has the form \(f(x)=ux+b\), with \(u\in L\) and \(b\in R\).
\end{lemma}
\begin{proof}
The proof of Proposition~\ref{prop:cyclotomic-irreducible} establishes that \(R^\times\) is cyclic.
Its unique subgroup of order \(p\) is \(1+p^{a-1}R\). The hypothesis therefore implies \(p\nmid|L|\).
Reduction \(L\to\F_p^\times\) is injective, since its kernel lies in the group \(1+pR\), of order \(p^{a-1}\).

The difference condition preserves congruence modulo \(p\), so \(f\) induces a permutation \(\bar f\) of \(\F_p\).
For distinct \(r,s\in\F_p\) and any lifts \(x\equiv r\), \(y\equiv s\pmod p\), the quotient
\((f(x)-f(y))/(x-y)\) belongs to \(L\), and its reduction is
\((\bar f(r)-\bar f(s))/(r-s)\). Injectivity of reduction makes it an element \(u_{r,s}\in L\) depending only on \(r,s\). Thus
\[
f(x)-f(y)=u_{r,s}(x-y).
\]
Replace \(x\) by \(x+p\) and subtract to obtain
\[
f(x+p)-f(x)=p\,u_{r,s}.
\]
The left side is independent of \(s\), so \(u_{r,s}\equiv u_{r,t}\pmod{p^{a-1}}\) for \(s,t\ne r\).
Injectivity of reduction modulo \(p\) gives \(u_{r,s}=u_{r,t}\).
Since also \(u_{r,s}=u_{s,r}\), all multipliers between distinct residue classes equal one \(u\in L\).
For two points in the same class, compare each with a point in another class and subtract.
This gives \(f(x)-f(y)=u(x-y)\) for all pairs; set \(y=0\).
\end{proof}

\begin{lemma}[Local rigidity at powers of two]\label{v73:lem-two-local-rigidity}
Let \(a\ge1\), \(R=\Z/2^a\Z\), and let \(L\le R^\times\) be cyclic.
No further assumption is needed if \(a\le2\); if \(a\ge3\), assume
\[
1+2^{a-2}R\not\subseteq L.
\]
Then every \(f\in\mathcal R_L(2^a)\) has the form \(f(x)=ux+b\), with \(u\in L\) and \(b\in R\).
\end{lemma}
\begin{proof}
Since \(f(1)-f(0)\in L\), subtracting \(f(0)\) and dividing by this unit reduces to
\(f(0)=0\), \(f(1)=1\), without changing the difference-orbit condition.
We prove that the normalized map is the identity.
This is immediate modulo \(2\); modulo \(4\), the condition gives \(f(2)\in2L=\{2\}\), determining the last point as well.

Modulo \(8\), the cyclic subgroups of the unit group are
\(\{1\},\{1,3\},\{1,5\},\{1,7\}\). For every \(x\in\Z/8\Z\), the difference condition gives
\[
f(x)\in S_L(x):=xL\cap\bigl(1+(x-1)L\bigr).
\]
If \(L=\{1\}\), this already determines \(f\). For the other three groups, the required intersections are
\[
\begin{array}{c|ccc}
L&S_L(2)&S_L(3)&S_L(4)\\ \hline
\{1,3\}&\{2\}&\{3\}&\{4\}\\
\{1,5\}&\{2\}&\{3\}&\{4\}\\
\{1,7\}&\{2\}&\{3\}&\{4\}
\end{array}
\]
The identity \(S_L(1-x)=1-S_L(x)\) also determines \(x=7,6,5\), completing this case.

Now let \(a\ge4\) and argue by induction on \(a\). Write \(\bar L\) for the image of \(L\) modulo \(2^{a-1}\).
We first check the induction hypothesis for \(\bar L\). If \(a=4\), its image is cyclic and cannot contain the noncyclic group \((\Z/8\Z)^\times\).
If \(a\ge5\) and \(\bar L\) contains \(1+2^{a-3}(\Z/2^{a-1}\Z)\), lift an element to \(v\in L\) with
\[
v\equiv1+2^{a-3}\pmod{2^{a-1}}.
\]
For an odd integer representative, \(v_2(v-1)=a-3\) and \(v_2(v+1)=1\), hence \(v_2(v^2-1)=a-2\).
Thus \(v^2\) modulo \(2^a\) has the form \(1+c2^{a-2}\) with \(c\) odd. Its square is \(1+2^{a-1}\), so it has order four and generates \(1+2^{a-2}R\), a contradiction.

The difference condition makes \(f\) descend modulo \(2^{a-1}\). This induced permutation fixes \(0,1\), so induction makes it the identity.
Consequently,
\[
f(x)=x+2^{a-1}\varepsilon(x),\qquad
\varepsilon:R\to\{0,1\},\quad\varepsilon(0)=\varepsilon(1)=0.
\]
If \(\varepsilon(x+2)\ne\varepsilon(x)\), the condition for the pair \(x+2,x\) gives \(v\in L\) such that
\[
2v\equiv2+2^{a-1}\pmod{2^a},\qquad
v\equiv1+2^{a-2}\pmod{2^{a-1}}.
\]
Such a \(v\) again generates \(1+2^{a-2}R\), a contradiction.
Therefore \(\varepsilon(x+2)=\varepsilon(x)\). It is constant on each parity class and vanishes on both because \(\varepsilon(0)=\varepsilon(1)=0\).
Undoing the normalization proves the claim.
\end{proof}

\subsection{Chinese Remainder Decomposition and Global Rigidity}
\label{v73:subsec-global-rigidity}
The local rigidity argument uses not the full projection of \(K\) at a prime power, but the multipliers that can act there while every other coordinate is fixed.
This distinction parallels orbit synchronization in Section~\ref{v44:sec-cyclotomic-local-global}: an iterate allowed by one local projection need not satisfy the phase conditions on the other components.

\leavevmode\phantomsection\label{v73:def-reduction-kernel}\ignorespaces
For \(d\mid n\), define the \textbf{unit reduction kernel}
\[
\Delta_n(d)=\ker\bigl((\Z/n\Z)^\times\longrightarrow(\Z/d\Z)^\times\bigr)
=\{u\in(\Z/n\Z)^\times:u\equiv1\pmod d\}.
\]
The unit group modulo one is taken to be trivial, so \(\Delta_n(1)=(\Z/n\Z)^\times\).

\par\Needspace{12\baselineskip}
\begin{theorem}[Cyclotomic difference-orbit rigidity]\label{v73:thm-ring-rigidity}
Let \(n\ge2\), and let \(K\le(\Z/n\Z)^\times\) be cyclic.
The group \(\mathcal R_K(n)\) strictly contains \(\operatorname{Aff}_K(n)\) if and only if at least one of the following occurs:
\begin{enumerate}[label=\textup{(R\arabic*)}]
\item There is a prime \(p\ge5\) with \(p\parallel n\) and \(\Delta_n(n/p)\subseteq K\).
\item There is an odd prime \(p\) with \(p^2\mid n\) and \(\Delta_n(n/p)\subseteq K\).
\item We have \(8\mid n\) and \(\Delta_n(n/4)\subseteq K\).
\end{enumerate}
In particular, if none occurs, condition~\eqref{v73:eq-difference-orbit} holds exactly for
\[
f(x)=ux+b,\qquad u\in K,\quad b\in\Z/n\Z.
\]
Here \(p\parallel n\) means \(p\mid n\) and \(p^2\nmid n\).
\end{theorem}
\begin{proof}
\textbf{Rigidity in the absence of the three exceptions.}
Write \(n=\prod_i p_i^{a_i}\) and identify, by the Chinese remainder theorem,
\[
R=\Z/n\Z\simeq\prod_iR_i,\qquad R_i=\Z/p_i^{a_i}\Z.
\]
Multiplication by a unit preserves every ideal. Thus the difference condition gives, for \(d\mid n\),
\[
x\equiv y\pmod d\quad\Longleftrightarrow\quad f(x)\equiv f(y)\pmod d.
\]
In particular, the \(i\)th output coordinate depends only on the \(i\)th input coordinate and induces a permutation \(f_i:R_i\to R_i\). Hence
\[
f(x_1,\ldots,x_t)=(f_1(x_1),\ldots,f_t(x_t)).
\]
Define the \textbf{independent local multiplier group}
\begin{equation}\label{v73:eq-independent-local-group}
L_i=\{u\in R_i^\times:(1,\ldots,1,u,1,\ldots,1)\in K\}.
\end{equation}
It is isomorphic to a subgroup of \(K\), and is therefore cyclic.

Fix \(i\). In every other component, fix the pair \(1,0\), and put
\(d_j=f_j(1)-f_j(0)\in R_j^\times\) for \(j\ne i\).
First use the pair \(1,0\) in all components. The difference condition gives \(k^0\in K\) whose \(j\)th component is \(f_j(1)-f_j(0)\).
Now let the \(i\)th pair be arbitrary \(x,y\), leaving all other pairs unchanged.
The corresponding \(k\in K\) must satisfy \(k_j=d_j=k_j^0\) for \(j\ne i\).
Thus \(k(k^0)^{-1}\) is one in every other component. Writing \(u_i=k_i^0\), we obtain
\[
f_i(x)-f_i(y)\in u_iL_i(x-y).
\]
Consequently \(u_i^{-1}f_i\in\mathcal R_{L_i}(p_i^{a_i})\).
This argument does not divide by \(x-y\), and also applies when \(n\) has only one prime divisor.

If \(a_i=1\) and \(p_i\ge5\), exclusion of the first case gives \(L_i\ne\F_{p_i}^\times\), so the prime-field case of Carlitz--McConnel makes the local map affine.
For \(p_i=2,3\), use the elementary observation above.
If \(p_i\) is odd and \(a_i\ge2\), exclusion of the second case gives
\(1+p_i^{a_i-1}R_i\not\subseteq L_i\), so Lemma~\ref{v73:lem-odd-local-rigidity} applies.
For \(p_i=2\), exclusion of the third case is precisely the hypothesis of Lemma~\ref{v73:lem-two-local-rigidity}; no extra hypothesis is needed when \(a_i\le2\).
In each instance the relevant unit reduction kernel acts only in the specified CRT component.

All \(f_i\) are therefore affine. The Chinese remainder theorem combines them into \(f(x)=ux+b\), with \(u\in R^\times\).
Taking \((x,y)=(1,0)\) gives \(u=f(1)-f(0)\in K\). Conversely, these affine maps preserve all difference orbits.

\textbf{A nonaffine permutation in each exceptional case.}
In each construction, change only the specified prime-power coordinate and leave the others fixed.
In the first case, \(\Delta_n(n/p)\) consists of all nonzero multipliers in that \(\F_p\) coordinate.
There we may exchange \(0,1\) and fix every other point: every quotient of corresponding nonzero differences is allowed.
This permutation fixes at least two distinct points but is not the identity, so it is not affine.

In the second case, write \(p^a\parallel n\), \(a\ge2\), and \(R_p=\Z/p^a\Z\).
The contained local kernel is \(1+p^{a-1}R_p\).
Let \(\bar x\in\{0,\ldots,p-1\}\) represent \(x\) modulo \(p\), and set
\begin{equation}\label{v73:eq-odd-local-shear}
g(x)=x+p^{a-1}\bar x^{\,2}.
\end{equation}
This translates each residue class modulo \(p\), so it is a permutation.
Differences within a class are unchanged; between classes they are units, and the ratio of new to old difference belongs to \(1+p^{a-1}R_p\).
If \(g(x)=vx+b\), the values at \(0,1\) force \(b=0\), \(v=1+p^{a-1}\).
But \(g(2)-2v=2p^{a-1}\ne0\), a contradiction.

In the third case, write \(2^a\parallel n\).
For \(a=3\), the kernel would be isomorphic to the noncyclic group \((\Z/8\Z)^\times\), and could not lie in the cyclic group \(K\). Thus \(a\ge4\).
On \(R_2=\Z/2^a\Z\), set
\begin{equation}\label{v73:eq-two-local-shear}
g(x)=x+2^{a-2}(x\bmod2).
\end{equation}
It fixes the even class and translates the odd class, so it is a permutation.
Differences of equal parity are unchanged; differences of opposite parity are units, and their ratios belong to \(1+2^{a-2}R_2\).
It fixes \(0,2\) but sends \(1\) to \(1+2^{a-2}\), so it cannot be affine.

In all three cases the local multipliers extend independently to elements of \(K\), giving permutations of \(R\) that preserve every difference orbit.
An affine map of \(R\) induces an affine map on each CRT component, contradicting these constructions. This proves the classification.
\end{proof}

\begin{theorem}[Multiple equal-difference structure and non-affine symmetries]\label{v74:thm-med-nonaffine}
Let \(p\) be an odd prime, \(a\ge2\), \(R=\Z/p^a\Z\), and \(K=\langle q\rangle\le R^\times\), where \(\gcd(p,q)=1\).
Put \(\tau=\ord_{p^a}(q)\) and \(t=v_p(\tau)\). For every unit \(u\), the cyclotomic coset \(C=uK\) satisfies
\[
H_+(C)=p^{a-t}R,\qquad \kappa_C=p^t,\qquad \omega_C=\frac{\tau}{p^t}.
\]
In particular, the following are equivalent:
\begin{enumerate}
\item \(p\mid\tau\);
\item \(C\) admits a multiple equal-difference decomposition with blocks larger than one;
\item There is a non-affine permutation \(f:R\to R\) such that \(f(x)-f(y)\in K(x-y)\) for all \(x,y\in R\).
\end{enumerate}
The third condition is equivalent to \(\mathcal R_K(p^a)\supsetneq\operatorname{Aff}_K(p^a)\).
It asserts the existence of extra permutations, not that every preserving map is non-affine.
\end{theorem}
\begin{proof}
The unit group at an odd prime power is cyclic of order \(p^{a-1}(p-1)\), so \(0\le t\le a-1\).
For \(1\le s\le a\), the group \(1+p^sR\) has order \(p^{a-s}\).
Subgroups of a cyclic group are determined by their orders; hence
\[
1+p^sR\subseteq K\quad\Longleftrightarrow\quad p^{a-s}\mid\tau
\quad\Longleftrightarrow\quad a-s\le t.
\]
The additive stabilizer \(H_+(C)\) is proper, so it is some \(p^sR\).
Proposition~\ref{v74:prop-orbit-additive-stabilizer} shows that the largest admissible ideal is \(p^{a-t}R\).
Now use \(\kappa_C=|H_+(C)|\) and \(\omega_C=\tau/\kappa_C\) from Theorem~\ref{thm:med-coset-partitions} to obtain the formulas and the equivalence of the first two conditions.

For the modulus \(p^a\), only the second exception of Theorem~\ref{v73:thm-ring-rigidity} can occur. Its condition is
\[
\Delta_{p^a}(p^{a-1})=1+p^{a-1}R\subseteq K.
\]
This subgroup has order \(p\), so the containment is equivalent to \(p\mid\tau\), and hence to the third condition.
When it holds, formula~\eqref{v73:eq-odd-local-shear} already supplies the non-affine permutation; no further existence argument is needed.
\end{proof}

Multiple equal-difference structure does not always strengthen rigidity.
Here it retains a complete additive direction, allowing perturbations within the corresponding congruence layers to preserve each difference orbit.
The result also follows from the odd-characteristic Galois-ring classification of Evdokimov--Ponomarenko \cite[Theorem~1.4]{EP2006NormalCyclotomic}; the unit-group structure and local rigidity already established here give a direct proof.
The equivalence is restricted to the stated odd prime-power case: modulo \(8\), the group \(K=\langle5\rangle=\{1,5\}\) has \(H_+(K)=\{0,4\}\),
but Lemma~\ref{v73:lem-two-local-rigidity} gives \(\mathcal R_K(8)=\operatorname{Aff}_K(8)\).
Thus nontrivial equal-difference blocks and extra non-affine symmetries cannot be identified without further conditions at powers of two.

\begin{example}[The field of nine elements and the residue ring modulo nine]\label{v73:ex-field-ring-rigidity}
The field \(\F_9\) in Example~\ref{v72:ex-carlitz-mcconnel-f9} admits a nontrivial Frobenius automorphism.
The ring \(\Z/9\Z\) is different: its additive group is cyclic, so an additive automorphism is determined by the image of \(1\).
For \(K=\{1,-1\}\le(\Z/9\Z)^\times\), we have
\(\Delta_9(3)=\{1,4,7\}\not\subseteq K\).
Thus its difference-orbit-preserving maps are exactly \(\pm x+b\), numbering \(18\).
For \(K=\langle4\rangle=\{1,4,7\}\), the second exception holds and
\(x\mapsto x+3(x\bmod3)^2\) is nonaffine.
The residue-ring exceptions therefore arise from local perturbations of congruence layers, not from field Frobenius automorphisms.
\end{example}

These results recover prime-field linearity from parallel differences, Frobenius from multiplicative cosets of field differences, quadratic polynomials from second divided differences on all triples, and, after excluding the stated exceptions, a common multiplier from residue-ring difference orbits.
Their underlying sets and permitted local variations differ; in particular, \(\F_{p^a}\) and \(\Z/p^a\Z\) are distinct algebraic structures when \(a>1\).
Section~\ref{v73:sec-common-cyclic-symmetry} will apply the last rigidity theorem to the coordinate symmetries common to all cyclic codes, using root-of-unity Fourier projections.

\section{Value Sets, Fibers, and Minimal Value Set Polynomials}
\label{v58:sec-value-sets}
Every fiber of a permutation polynomial has one element. For a general polynomial, computing the nonempty fibers and their sizes
determines the number of values actually attained. The root bound gives a general restriction, and subspace polynomials provide constructions attaining it.

\begin{definition}\leavevmode\phantomsection\label{v58:def-value-set}
Let \(K=\F_Q\) and \(F\in K[X]\). The set
\[
\mathcal V_F=\{F(x):x\in K\}
\]
is called the \textbf{value set} of \(F\) over \(K\). For \(b\in K\), write
\(N_b(F)=\#\{x\in K:F(x)=b\}\) for the \textbf{fiber size} over \(b\).
The fiber sizes are integers, and \(\sum_bN_b(F)=Q\).
\end{definition}

\begin{proposition}[The root bound and a lower bound for the value set]\label{v58:prop-value-set-bound}
If \(F\) has ordinary degree \(d\ge1\), then
\[
|\mathcal V_F|\ge\left\lceil\frac Qd\right\rceil
=\left\lfloor\frac{Q-1}{d}\right\rfloor+1.
\]
\end{proposition}
\begin{proof}
For every \(b\), the nonzero polynomial \(F-b\) has at most \(d\) roots, so \(N_b(F)\le d\).
Summing the sizes of the nonempty fibers gives \(Q\le d|\mathcal V_F|\).
\end{proof}

\leavevmode\phantomsection\label{v58:def-minimal-value-set}\ignorespaces
A nonconstant polynomial attaining this lower bound is called a \textbf{minimal value set polynomial}.
The degree here is the ordinary polynomial degree: reduction modulo \(X^Q-X\) preserves the function and its value set, but may change the degree and whether this bound is attained.
For example, \(X\) and \(X^Q\) define the same function on \(\F_Q\), but the former is a minimal value set polynomial and the latter is not.

\subsection{From Additive Cosets to Power Images of Subspaces}
A subspace polynomial determines all its nonempty fibers, not just its zero fiber.
Postcomposing such a polynomial with a power map yields further minimal value set polynomials.
\leavevmode\phantomsection\label{v58:def-power-image}\ignorespaces
For a set \(W\) and a positive integer \(v\), the set \(W^v=\{w^v:w\in W\}\) is called the \(v\)th \textbf{power image} of \(W\).

\begin{proposition}[Powers of subspace polynomials and their value sets]\label{v58:prop-subspace-power-values}
Let \(k=\F_q\subseteq E=\F_{q^n}\), let \(U\subsetneq E\) be an \(s\)-dimensional \(k\)-subspace,
and put
\[
L_U(X)=\prod_{u\in U}(X-u),\qquad W=L_U(E).
\]
Take \(v\mid q-1\), \(a\in E^\times\), and \(b\in E\). Then
\[
F(X)=aL_U(X)^v+b
\]
is a minimal value set polynomial over \(E\), and
\begin{equation}\label{v58:eq-subspace-power-values}
\mathcal V_F=aW^v+b,\qquad
|\mathcal V_F|=1+\frac{q^{n-s}-1}{v}.
\end{equation}
The fiber over \(b\) has \(q^s\) elements; all other nonempty fibers have \(vq^s\) elements.
\end{proposition}
\begin{proof}
Theorem~\ref{thm:subspace-polynomial} shows that \(L_U\) is a \(k\)-linearized polynomial of ordinary degree \(q^s\),
with kernel in \(E\) exactly \(U\). Hence \(W\) is an \(n-s\)-dimensional \(k\)-subspace,
and each nonempty fiber is an additive coset of \(U\) of size \(q^s\).

Since \(v\mid q-1\), the group \(\mu_v\) of \(v\)th roots of unity is contained in \(k^\times\).
For \(w\in W\setminus\{0\}\), the elements satisfying \(z^v=w^v\) are exactly \(w\mu_v\);
they all belong to \(W\), and there are \(v\) of them. Thus each nonempty fiber of the power map on the nonzero part of \(W\) has \(v\) elements,
while zero contributes a separate value. Multiplying by \(a\) and adding \(b\) gives the asserted value-set and fiber formulas.

Finally, \(\deg F=vq^s\) and \(v\mid q^{n-s}-1\). Therefore
\[
\left\lceil\frac{q^n}{\deg F}\right\rceil
=\left\lceil\frac{q^{n-s}}v\right\rceil
=1+\frac{q^{n-s}-1}{v}=|\mathcal V_F|.
\]
This proves minimality.
\end{proof}

The proof used \(|\ker L_U|=\deg L_U\). An arbitrary linearized polynomial need not have this property:
its kernel may have cardinality smaller than its ordinary degree, so it need not be a minimal value set polynomial; see Exercise~\ref{v58:ex-linearized-minimal-values}.
When \(U=0\), the construction reduces to \(aX^v+b\); when \(v=1\), the value set is an affine subspace.

\Needspace{12\baselineskip}
\begin{example}[A five-element value set in the field of order 27]\label{v58:ex-F27-five-values}
Take \(E=\F_{27}=\F_3(\theta)\), where \(\theta^3-\theta-1=0\).
This cubic has no root in \(\F_3\) and is therefore irreducible. Put
\[
L(X)=X^3-X,\qquad F(X)=L(X)^2.
\]
Additive Hilbert 90 gives \(L(E)=\ker\Tr_{E/\F_3}\).
From \(\theta^3=\theta+1\) and \(\theta^9=\theta-1\), we obtain
\(\Tr(1)=\Tr(\theta)=0\). The trace kernel is two-dimensional, so
\[
W=L(E)=\operatorname{span}_{\F_3}\{1,\theta\}.
\]
Its eight nonzero elements form four pairs \(\{w,-w\}\), with representatives \(1,\theta,\theta+1,\theta-1\).
Squaring gives
\[
\boxed{\quad\mathcal V_F=
\{0,1,\theta^2,\theta^2-\theta+1,\theta^2+\theta+1\}.\quad}
\]
These five values are distinct, since two nonzero elements have the same square if and only if they differ by a sign.
Since \(\ker L=\F_3\), the zero fiber has three elements, and each of the other four nonempty fibers has six elements.
Thus
\[
|\mathcal V_F|=5=\left\lceil\frac{27}{6}\right\rceil.
\]
This value set is neither a subfield nor an \(\F_3\)-affine subspace: either type of set must have size a power of three.
It is determined jointly by additive cosets and the nonzero fibers of the squaring map.
\end{example}

\begin{remark}
The construction above gives a class of minimal value set polynomials with a direct proof.
Further questions ask which sets can occur as value sets of minimal value set polynomials and which minimal value set polynomials realize a given such set.
Borges and Reis proved that, in \(\F_{p^N}\), the realizable sets of size greater than two are exactly
\(aW^v+b\), where \(a\ne0\), \(W\) is a subspace containing \(1\) over some intermediate field \(\F_{p^h}\),
and \(v\mid p^h-1\); see \cite[Theorem 2.2]{BorgesReisMinimal2025}.
This classifies the value sets. The subsequent preprint \cite{BorgesReisComplete2026} resolves the classification of realizing polynomials proposed in Conjecture 1.1 of that paper.
\end{remark}
\subsection{The structural classification of minimal value set polynomials}
Can the additive and multiplicative structure of a value set be lifted to a compositional structure of every polynomial realizing it? The subsequent classification of Borges--Reis answers this question affirmatively.
\begin{theorem}[Structural classification of MVSPs, {\cite[Theorem~1.1]{BorgesReisComplete2026}}]\label{v109:thm-mvsp-classification}
Let $q=p^n$ and $F\in\F_q[X]$ with $|\mathcal V_F|>2$. Then $F$ is a minimal value set polynomial if and only if, after invertible affine changes of variable and value, it has the form
\[
F(X)=A(P(X))^v,
\]
where $P$ is an MVSP whose value set is exactly a subfield $K=\F_{q'}\subseteq\F_q$,
$U\subseteq K$ is a vector subspace over a subfield $k_0=\F_{p^h}\subseteq K$, and
\[
A(X)=\prod_{u\in U}(X-u),\qquad v\mid p^h-1,\qquad 2v\deg A<q'.
\]
Then $|\mathcal V_F|=1+(q'/\deg A-1)/v$. The composition is an identity of ordinary polynomials.
\end{theorem}
\begin{proof}
Necessity follows from \cite[Theorem~1.1]{BorgesReisComplete2026}. We prove sufficiency below.
Put $a=\deg A=|U|$ and $d=\deg P$. The subspace polynomial $A$ is $k_0$-linear and has kernel exactly $U$ on $K$, so its image $B=A(K)$ has $q'/a$ elements.
Since $v\mid|k_0^\times|$, every nonzero $v$-th-power fiber in $B$ has exactly $v$ elements. The image therefore has size $r+1$, where $r=(q'/a-1)/v$.
Minimality of $P$ gives $(q'-1)d<q\le q'd$. Consequently,
\[
r\deg F=(q'-a)d<q\le(q'+a(v-1))d=(r+1)\deg F.
\]
Thus $\lceil q/\deg F\rceil=r+1$. Invertible affine changes preserve degree and value-set size.
\end{proof}
The restriction $2va<q'$ gives $q'/a>2v$. Since $q'/a\equiv1\pmod v$, there are at least three values. One-value and two-value cases are excluded here.

The power-value condition also has a direct rigidity consequence in low degree: if $F(\F_q)\subseteq\{0\}\cup(\F_q^\times)^v$, with $v\ge2$, $v\mid q-1$ and $\deg F<\sqrt q$, then $F$ itself is a $v$th power. Exercise~\ref{v111:ex-low-degree-power} proves this from the multiplicative Weil bound and uses trace polynomials at the square-root boundary to explain why the exponent cannot be prescribed arbitrarily. It does not assume a minimal value set and does not replace the general classification above. Its final part strengthens the degree condition to a bound on the number of distinct roots in a splitting field, allowing arbitrary multiplicities.

\begin{example}[Five values and all fibers over the field of 729 elements]\label{v109:ex-mvsp729}
Take $E=\F_{729}$ and $K=\F_{27}$, and put
\[
F(X)=(X^{84}-X^{28})^2,\qquad N(x)=x^{28},\quad L(y)=y^3-y.
\]
The successive maps are $E\xrightarrow{N}K\xrightarrow{L}W\xrightarrow{z\mapsto z^2}W^2$, where $N=\Nm_{E/K}$ and $W=\ker\Tr_{K/\F_3}$.
The norm is surjective, with nonzero fibers of size 28. The kernel of $L$ is $\F_3$, and its image $W$ has nine elements. Pairing its eight nonzero elements by sign gives
\[
|\mathcal V_F|=5=\left\lceil\frac{729}{168}\right\rceil,
\qquad N_0(F)=1+2\cdot28=57,
\]
\[
N_b(F)=2\cdot3\cdot28=168\quad(b\in W^2\setminus\{0\}).
\]
Hence $57+4\cdot168=729$. Using $\theta$ from Example~\ref{v58:ex-F27-five-values}, the five values are again
$0,1,\theta^2,\theta^2-\theta+1,\theta^2+\theta+1$.
In the classification theorem take $P=X^{28}$, $A=X^3-X$, and $v=2$. The polynomial $P$ itself satisfies $27=\lceil729/28\rceil$, and $2\cdot2\cdot3<27$.
Compared with Subsection~\ref{v105:subsec-norm-trace-values}, the final square merges nonzero trace-kernel fibers further; the zero fiber still requires separate treatment.
\end{example}

\subsection{From Norm Fibers to a Trace Kernel}\label{v105:subsec-norm-trace-values}
The preceding construction first combines additive cosets and then takes a power image.
One may instead first combine multiplicative cosets by a norm and then take a Frobenius difference in its target field.
The image of the latter linear map is a trace kernel, so both the value set and all fiber sizes can be determined directly.
The following binomials belong to the minimal value set family in \cite[Theorem A(v)]{AprigioGuardieiro2026Binomials}.

\begin{proposition}[A Frobenius difference after a norm]\label{v105:prop-norm-trace-values}
Let \(Q\) be a prime power, let positive integers \(n,s\) satisfy \(s\mid n\) and \(s\ge2\), and put
\[
k=\F_Q,\qquad K=\F_{Q^s}\subseteq E=\F_{Q^n},\qquad
t=\frac{Q^n-1}{Q^s-1}.
\]
Then the binomial
\[
F(X)=X^{Qt}-X^t
\]
is a minimal value set polynomial over \(E\), with
\begin{equation}\label{v105:eq-norm-trace-values}
\mathcal V_F=\ker\Tr_{K/k},\qquad
|\mathcal V_F|=Q^{s-1}.
\end{equation}
Its fiber sizes are
\begin{equation}\label{v105:eq-norm-trace-fibers}
N_b(F)=
\begin{cases}
1+(Q-1)t,&b=0,\\
Qt,&b\in\ker\Tr_{K/k}\setminus\{0\},\\
0,&b\notin\ker\Tr_{K/k},
\end{cases}
\qquad(b\in E).
\end{equation}
\end{proposition}
\begin{proof}
Extend the norm to zero, writing \(N:E\to K\), \(N(x)=x^t\).
It is surjective on the nonzero elements, with \(t\) preimages of each nonzero value, whereas \(N^{-1}(0)=\{0\}\).
On \(K\), put \(L(y)=y^Q-y\). Additive Hilbert~90 gives
\[
L(K)=\ker\Tr_{K/k},\qquad\ker L=k.
\]
Every nonempty fiber of \(L\) is therefore an additive coset of \(k\), of size \(Q\),
and its image has \(Q^{s-1}\) elements.
Since \(F=L\circ N\) and \(N\) is surjective, this proves \eqref{v105:eq-norm-trace-values}.

In the zero fiber \(L^{-1}(0)=k\), zero has one \(N\)-preimage and each of the remaining \(Q-1\) elements has \(t\) preimages.
Hence \(N_0(F)=1+(Q-1)t\).
If \(b\ne0\) belongs to \(L(K)\), then \(L^{-1}(b)\) does not contain zero;
each of its \(Q\) elements has \(t\) preimages under \(N\), giving \(N_b(F)=Qt\).
All other fibers are empty.

Finally, \(\deg F=Qt\), and
\[
\left\lfloor\frac{Q^n-1}{Qt}\right\rfloor+1
=\left\lfloor\frac{Q^s-1}{Q}\right\rfloor+1
=Q^{s-1}=|\mathcal V_F|.
\]
Thus \(F\) attains the lower bound in Proposition~\ref{v58:prop-value-set-bound}.
\end{proof}

The zero fiber has \(t-1\) fewer elements than any other nonempty fiber, a difference also visible in the root multiplicities.
If \(p=\Char(k)\), then
\[
t=\sum_{j=0}^{n/s-1}Q^{sj}\equiv1\pmod p,
\qquad
F(X)=X^t\bigl(X^{(Q-1)t}-1\bigr).
\]
Since \(p\nmid(Q-1)t\) and \((Q-1)t\mid Q^n-1\), the second factor has \((Q-1)t\) distinct nonzero roots in \(E\).
The root zero has multiplicity \(t\), so the number of distinct roots in the zero fiber is \(t-1\) less than the ordinary degree \(Qt\).
When \(t=1\), the construction is the linearized polynomial \(X^Q-X\), and all nonempty fibers have the same size.

\begin{example}[A three-element value set in the field of order 81]\label{v105:ex-F81-norm-values}
Take \(E=\F_{81}\) and its subfield \(K=\F_9\).
Since \(X^2+1\) is irreducible over \(\F_3\), write \(K=\F_3(u)\) with \(u^2=-1\).
For \(a,b\in\F_3\),
\[
(a+bu)^3=a-bu,\qquad
\Tr_{K/\F_3}(a+bu)=2a,\qquad
(a+bu)^3-(a+bu)=bu.
\]
Consequently, \(\ker\Tr_{K/\F_3}=\F_3u\).
The norm \(N_{E/K}(x)=x^{10}\) has ten preimages of each nonzero value, so
\[
F(X)=X^{30}-X^{10},\qquad
\mathcal V_F=\{0,u,-u\}.
\]
The zero fiber has \(1+2\cdot10=21\) elements, and the other two nonempty fibers have \(3\cdot10=30\) elements each.
Their sizes sum to \(21+30+30=81\), while
\[
|\mathcal V_F|=3=\left\lfloor\frac{80}{30}\right\rfloor+1.
\]
The five-element value set in Example~\ref{v58:ex-F27-five-values} is the square image of a trace kernel.
Here the three-element value set is itself a trace kernel, but the input first passes through the multiplicative fibers of a norm.
\end{example}

\subsection{Pairs with the Same Image and Collision Counts}
\leavevmode\phantomsection\label{v58:def-collision-number}\ignorespaces
For a function \(F:K\to K\), \(K=\F_Q\), define its \textbf{collision count} by
\[
C_F=\#\{(x,y)\in K^2:F(x)=F(y)\}
=\sum_{b\in K}N_b(F)^2.
\]
The pairs are ordered, and the count includes the \(Q\) diagonal pairs with \(x=y\).

\begin{proposition}[Collision counts and value-set sizes]\label{v58:prop-value-collision}
For every function \mbox{\(F:K\to K\)},
\[
|\mathcal V_F|\ge\frac{Q^2}{C_F}.
\]
Equality holds if and only if all nonempty fibers have equal size. In particular, \(C_F=Q\) if and only if \(F\) is a permutation.
\end{proposition}
\begin{proof}
Applying the Cauchy--Schwarz inequality to the sizes of the nonempty fibers gives
\[
Q^2=\left(\sum_{b\in\mathcal V_F}N_b(F)\right)^2
\le |\mathcal V_F|\sum_{b\in\mathcal V_F}N_b(F)^2
=|\mathcal V_F|C_F.
\]
Equality holds precisely when all these positive integers are equal.
Moreover, \(C_F-Q=\sum_bN_b(F)(N_b(F)-1)\ge0\),
which vanishes precisely when every fiber has at most one element; for finite sets of equal size, this is equivalent to a permutation.
\end{proof}

The collision count in Example~\ref{v58:ex-F27-five-values} is
\[
C_F=3^2+4\cdot6^2=153.
\]
Hence \(|\mathcal V_F|\ge\lceil729/153\rceil=5\), again giving the exact value-set size.
The inequality before rounding is strict here because the nonempty fibers have unequal sizes.
Proposition~\ref{v58:prop-collision-fourier} expresses the collision count as the mean squared absolute value of additive character sums,
thereby connecting value-set counting with Fourier duality.

For a fixed linearized polynomial with coefficients in the base field, every nonempty fiber on each finite extension is a coset of its kernel.
Thus the kernel dimension in Corollary~\ref{v70:cor-all-extension-kernels} determines both the value-set size and the collision count, with equality in the preceding bound.
Exercise~\ref{v70:ex-linearized-values-all-extensions} illustrates this using one polynomial and its third iterate.

\paragraph{Descent of collision factors and extension of the base field.}
A nondiagonal factor of the collision equation over an extension need not have a nondiagonal point over the base field. The following quartic rational function expresses this distinction as a root-orbit degree condition. This family occurs among the low-degree rational permutations in~\cite[Theorem 1.4]{DingZieve2022Rational}.
\begin{proposition}\label{v106:prop-quartic-rational}
Let \(q\) be an odd prime power and let \(h(T)=T^3+aT+b\in\F_q[T]\) be irreducible. Put
\[
N(T)=T^4-2aT^2-8bT+a^2,\qquad R=N/h.
\]
Then \(R\) permutes \(\mathbb P^1(\F_{q^m})\) if and only if \(3\nmid m\). Here \(\mathbb P^1(K)=K\cup\{\infty\}\), and a rational function is evaluated as its reduced projective map.
\end{proposition}
\begin{proof}
Irreducible polynomials over finite fields are separable. For each root \(\gamma\) of \(h\), the relation \(\gamma^3=-a\gamma-b\) gives \(N(\gamma)=h'(\gamma)^2\ne0\). Thus the numerator and denominator are coprime, and \(R(\infty)=\infty\).
In the splitting field we have the exact identity
\begin{equation}\label{v106:eq-quartic-collision}
\begin{split}
N(X)h(Y)-N(Y)h(X)={}&(X-Y)\\[-0.4em]
&\cdot\prod_{h(\gamma)=0}\bigl[XY-\gamma(X+Y)-2\gamma^2-a\bigr].
\end{split}
\end{equation}
Substitute \(Y=(\gamma X+2\gamma^2+a)/(X-\gamma)\) in the left-hand side, clear denominators, and reduce using \(\gamma^3+a\gamma+b=0\). The remainder is zero, so the corresponding factor of bidegree \((1,1)\) divides the left-hand side. Its coefficient matrix has determinant \(-3\gamma^2-a=-h'(\gamma)\ne0\), hence the factor is irreducible. The three factors are distinct and are also distinct from \(X-Y\). The product of these three factors and \(X-Y\) has bidegree \((4,4)\), the bidegree of the left-hand side. Comparing the coefficient of \(X^4Y^3\) gives scalar factor one.

If \(3\nmid m\), the change-of-step formula for Frobenius shows that \(h\) stays irreducible over \(K=\F_{q^m}\). Distinct finite points \(x,y\in K\) with the same image would, by the collision identity, force a root \(\gamma\) to satisfy
\[
2\gamma^2+(x+y)\gamma+(a-xy)=0,
\]
contradicting \([K(\gamma):K]=3\). Since \(h\) has no root in \(K\), no finite point maps to \(\infty\). The projective map is injective and hence a permutation. If \(3\mid m\), the three distinct denominator roots belong to \(K\), and all three, together with \(\infty\), map to \(\infty\). Injectivity fails.
\end{proof}
\begin{example}[A five-cycle and a cubic root orbit]\label{v106:ex-quartic-f5}
Reuse \(h=X^3+X+1\in\F_5[X]\) from Example~\ref{v105:ex-cubic-same-discriminant}. Then
\[
R(X)=\frac{X^4+3X^2+2X+1}{X^3+X+1}.
\]
\begin{center}\begin{tabular}{c|ccccc}
\(x\)&0&1&2&3&4\\\hline
\(R(x)\)&1&4&3&0&2
\end{tabular}\end{center}
The finite points form the five-cycle \((0,1,4,2,3)\), and \(\infty\) is fixed. The function still permutes \(\F_{25}\), whereas over \(\F_{125}\) four distinct points map to \(\infty\). Its extension behavior is controlled by the cubic root orbit, rather than by the length of the original five-cycle.
\end{example}

\subsection{Exact Value Counts from an Affine Linearized Equation}
\label{v68:sec-affine-values}
A change of coordinate can turn the nonzero fibers of an ordinary polynomial into fibers of a linearized map.
For the following family, this determines the entire value distribution from the intersection of two finite fields.
The formula is the full-affine specialization of~\cite[Corollary 5.4]{Kumallagov2026Value}; we give a direct proof.

\begin{proposition}[A value distribution controlled by an intersection field]\label{v68:prop-affine-values}
Let \(p\) be prime, let \(e,m\ge1\), and put
\[
N=p^e,\qquad Q=p^m,\qquad s=p^{\gcd(e,m)},\qquad
P_N(X)=X(X-1)^{N-1}\in\F_p[X].
\]
Thus \(\F_N\cap\F_Q=\F_s\) in a fixed algebraic closure.
The zero fiber of \(P_N\) on \(\F_Q\) is \(\{0,1\}\).
Among nonzero target values, the only possible fiber sizes are \(0,1,s\), and their numbers are
\begin{equation}\label{v68:eq-affine-fiber-counts}
\begin{aligned}
\#\{t\in\F_Q^\times:N_t(P_N)=0\}&=\frac Qs,\\
\#\{t\in\F_Q^\times:N_t(P_N)=1\}&=\frac{(Q-1)(s-2)}{s-1},\\
\#\{t\in\F_Q^\times:N_t(P_N)=s\}&=\frac{Q-s}{s(s-1)}.
\end{aligned}
\end{equation}
Consequently,
\begin{equation}\label{v68:eq-affine-value-count}
|\mathcal V_{P_N}|=Q-\frac Qs.
\end{equation}
No containment assumption between \(\F_N\) and \(\F_Q\) is required.
\end{proposition}
\begin{proof}
Since \(N\ge2\), the zeros are exactly \(0,1\).
Fix \(t\in\F_Q^\times\). Every solution of \(P_N(u)=t\) has \(u\ne0,1\), so the substitution
\[
v=\frac1{u-1},\qquad u=1+\frac1v
\]
gives
\begin{equation}\label{v68:eq-affine-fiber-linearization}
P_N(u)=t\quad\Longleftrightarrow\quad tv^N-v=1.
\end{equation}
Conversely, a solution of the right-hand equation is neither \(0\) nor \(-1\):
zero does not satisfy it, and \((-1)^N=-1\) would give \(t=0\) if \(v=-1\).
Thus the substitution is a bijection between the two solution sets.

The map \(A_t(v)=tv^N-v\) is \(\F_s\)-linear on \(\F_Q\).
A nonzero element lies in its kernel precisely when
\[
v^{N-1}=t^{-1}.
\]
The cyclic group \(\F_Q^\times\) and the identity
\(\gcd(N-1,Q-1)=s-1\) show that this equation has either no solutions or exactly \(s-1\) solutions.
Consequently, \(|\ker A_t|\) is either \(1\) or \(s\), and the second case occurs for exactly \((Q-1)/(s-1)\) nonzero values of \(t\).
If the kernel is trivial, \(A_t\) is bijective, so \(A_t(v)=1\) has one solution.
If the kernel has size \(s\), its fiber over \(1\) is either empty or a coset of the kernel and therefore has \(s\) elements.

Write \(a_j=\#\{t\in\F_Q^\times:N_t(P_N)=j\}\) for \(j=0,1,s\).
The trivial-kernel cases give
\[
a_1=(Q-1)-\frac{Q-1}{s-1}=\frac{(Q-1)(s-2)}{s-1}.
\]
Counting the nonzero targets and then all source elements yields
\[
a_0+a_1+a_s=Q-1,\qquad 2+a_1+sa_s=Q.
\]
These two identities give the remaining formulas in~\eqref{v68:eq-affine-fiber-counts}.
Since the zero value is attained, precisely \(a_0=Q/s\) values are omitted, proving~\eqref{v68:eq-affine-value-count}.
\end{proof}

\begin{example}[The same quartic over two extension fields]\label{v68:ex-affine-values-small}
In characteristic two, take
\[
P_4(X)=X(X-1)^3=X^4+X^3+X^2+X.
\]
Over \(\F_8\), the intersection \(\F_4\cap\F_8=\F_2\) gives \(s=2\).
There are four omitted values and four nonempty fibers, each of size two.
Over \(\F_{16}\), the intersection is \(\F_4\), so \(s=4\).
There are four omitted values, ten singleton fibers, the zero fiber of size two, and one nonzero fiber of size four.
Thus
\[
\begin{array}{c|rrrr|r}
&N_t=0&N_t=1&N_t=2&N_t=4&|\mathcal V_{P_4}|\\\hline
\F_8&4&0&4&0&4\\
\F_{16}&4&10&1&1&12
\end{array}
\]
where each entry in the middle four columns counts target values.
The ordinary degree is four in both cases, but the intersection field changes the kernel sizes in~\eqref{v68:eq-affine-fiber-linearization} and hence the value distribution.
\end{example}

\subsection{Two-Point Coefficient Fibers and Fractional-Linear Involutions}\label{v85:sec-pair-coefficient-fibers}
A value-set problem can also use a root set as its variable. Let \(H\le k^\times\), where \(k=\F_q\), and write \(\binom H2\) for the family of two-element subsets of \(H\).
For \(\{x,y\}\in\binom H2\), the root polynomial is
\[
 p_{\{x,y\}}(T)=T^2-e_1T+e_2,
 \qquad e_1=x+y,\quad e_2=xy.
\]
Fix \(A,B\in k\). A coefficient fiber is a fiber of the map
\[
 \Phi_{A,B}:\binom H2\longrightarrow k,\qquad
 \{x,y\}\longmapsto B(x+y)-Axy.
\]
Put
\begin{equation}\label{v85:eq-pair-intersection}
 N_H(A,B,C)=\#\{(x,y)\in H^2:x\ne y,\ Axy-B(x+y)+C=0\}.
\end{equation}
Then \(\#\Phi_{A,B}^{-1}(C)=N_H(A,B,C)/2\). This division is in the integers, so it remains valid in characteristic two.

\leavevmode\phantomsection\label{v85:def-projective-line}\ignorespaces
To include poles of fractional-linear maps, adjoin \(\infty\) to \(k\), obtaining the projective line \(\mathbb P^1(k)\).
Equivalently, its points are the one-dimensional subspaces of \(k^2\): \(x\) corresponds to \([x:1]\), and \(\infty\) to \([1:0]\).
An invertible two-by-two matrix acts on these subspaces, and nonzero scalar multiples induce the same transformation.

\begin{proposition}[Two-root coefficient conditions and projective involutions]\label{v85:prop-pair-involution}
Suppose \(\Delta=B^2-AC\ne0\). The matrix
\[
 J=\begin{pmatrix}B&-C\\A&-B\end{pmatrix},\qquad J^2=\Delta I,
\]
induces a transformation of order at most two on \(\mathbb P^1(k)\), given by
\[
 J(x)=\frac{Bx-C}{Ax-B}.
\]
For \(x,y\in k\), the two-variable equation is equivalent to \(y=J(x)\), restricting to finite coordinates. Consequently,
\begin{equation}\label{v85:eq-involution-intersection}
 N_H(A,B,C)=|H\cap J(H)|-
 \#\{x\in H:Ax^2-2Bx+C=0\}.
\end{equation}
Except when the characteristic is two and \(A=C=0\), the map \(J\) is nonidentity and is a projective involution; the subtracted fixed-point count is then at most two.
In the exceptional case, \(B\ne0\), \(J\) is the identity, and \(N_H=0\).
\end{proposition}
\begin{proof}
Matrix multiplication gives \(J^2=\Delta I\), and \(\det J=-\Delta\ne0\).
Rewrite the equation as \((Ax-B)y=Bx-C\). If \(Ax-B=0\), then \(A\ne0\) and \(Bx-C=\Delta/A\ne0\), so no finite \(y\) exists. Otherwise division gives the formula.
Since \(J^2\) is the identity on the projective line, the conditions \(x\in H\) and \(J(x)\in H\) are equivalent to \(x\in H\cap J(H)\). Each such \(x\) determines exactly one \(y\); subtracting the solutions with \(y=x\) gives the count.
The matrix \(J\) is scalar exactly when \(A=C=0\) and \(2B=0\). Under \(\Delta\ne0\), this is precisely the stated exception. In the nonscalar case the fixed-point equation is not the zero polynomial, so it has at most two roots.
\end{proof}

\begin{example}[Coset obstructions and degenerate fibers]\label{v85:ex-pair-degenerate}
Let \(n=|H|\ge2\). If \(A\ne0\) and \(B=0\), put \(\beta=-C/A\). The equation is \(xy=\beta\), hence
\[
 N_H(A,0,C)=n\mathbf1_H(\beta)-\nu_H(\beta),\qquad
 \nu_H(\beta)=\#\{x\in H:x^2=\beta\}.
\]
The square map on a cyclic group gives \(\nu_H(\beta)=\gcd(2,n)\) if \(\beta\in H^2=\{h^2:h\in H\}\), and zero otherwise.
For example, in \(H=\{1,2,4\}\le\F_7^\times\), the equation \(xy=3\) has no solution, whereas \(xy=1\) has the two distinct-coordinate solutions \((2,4),(4,2)\).

If \(A\ne0\) and \(\Delta=0\), then
\[
 Axy-B(x+y)+C=A(x-B/A)(y-B/A),
\]
so \(N_H=2(n-1)\mathbf1_H(B/A)\). If \(A=0\) and \(B\ne0\), the equation is \(x+y=C/B\).
In particular, \(x+y=0\) has \(n\) distinct-coordinate solutions exactly when the characteristic is not two and \(-1\in H\); otherwise it has none.
These cases exhibit the constraints imposed by multiplicative cosets, linear factors, and the diagonal.
\end{example}
For \(q=p^r\), taking \(p\)th powers of both coordinates also gives
\[
 N_H(A^p,B^p,C^p)=N_H(A,B,C),
\]
since Frobenius preserves \(H\) and distinctness of coordinates. Proposition~\ref{v85:prop-coset-frobenius} thus acts on these fiber parameters as well.
An exact Jacobi-sum expression for nonzero sum fibers is given in Exercise~\ref{v85:ex-restricted-sums}.

\section{Power Sums and Hermite's Criterion}
The permutation property requires each value to have exactly one preimage. Let \(k=\F_q\) and put \(N_b=\#\{a\in k:F(a)=b\}\). Then
\[
\sum_{a\in k}F(a)^t=\sum_{b\in k}N_b b^t.
\]
Thus the power sums of \(F\) record information about fiber sizes. We first express a power sum as a coefficient of the canonical remainder,
then use the Vandermonde matrix of the distinct base-field elements to determine \(N_b\) modulo the characteristic.
Under the criterion's hypotheses, we obtain \(N_b\equiv1\pmod p\) for every fiber; together with \(\sum_bN_b=q\), this implies \(N_b=1\) for every fiber.

\begin{lemma}\label{lem:power-sums}
For a positive integer \(t\),
\[
\sum_{a\in\F_q}a^t=
\begin{cases}
-1,&q-1\mid t,\\
0,&q-1\nmid t.
\end{cases}
\]
For \(t=0\), interpreting the summand as the constant polynomial \(X^0=1\) gives the sum \(q=0\).
Consequently, for \(\deg H<q\), we have \(\sum_aH(a)=-[X^{q-1}]H\).
\end{lemma}
\begin{proof}
If \(q-1\mid t\), the \(t\)th power of every nonzero element is one, so the sum is \(q-1=-1\).
Otherwise, choose a primitive element \(\zeta\). Then \(\zeta^t\ne1\), and multiplication by \(\zeta\) permutes the nonzero elements,
so the power sum equals itself multiplied by \(\zeta^t\) and must be zero. The last assertion follows by applying this formula to each term of \(H\).
\end{proof}

\begin{theorem}[Hermite's criterion]\label{thm:hermite}
A polynomial \(F\in\F_q[X]\) is a permutation polynomial if and only if:
\begin{enumerate}
\item \(F\) has exactly one zero in \(\F_q\);
\item for \(1\leq t\leq q-2\), the remainder of \(F^t\) modulo \(X^q-X\) has no \(X^{q-1}\) term.
\end{enumerate}
For \(q=2\), the second family of conditions is empty.
\end{theorem}
\begin{proof}
If \(F\) induces a permutation, it has exactly one zero, and for every positive integer \(t\),
\(\sum_aF(a)^t=\sum_b b^t\).
For \(1\le t\le q-2\), the right-hand side is zero;
by Lemma~\ref{lem:power-sums}, this is equivalent to the coefficient of \(X^{q-1}\) in the canonical remainder of \(F^t\) being zero.

Conversely, for each \(b\in k\), let the integer \(N_b\ge0\) be the number of solutions of \(F(a)=b\).
We have \(\sum_bN_b=q\), and the hypothesis gives \(N_0=1\).
For \(1\le t\le q-2\), the power-sum condition gives \(\sum_bN_b b^t=\sum_b b^t=0\).
For \(t=q-1\), since \(F\) has exactly one zero, \(\sum_a F(a)^{q-1}=q-1=-1=\sum_b b^{q-1}\).
For \(t=0\), interpret every \(b^0\) as one, including \(b=0\); both sides then equal \(q=0\).
Thus, in \(k\),
\[
\sum_{b\in k}(N_b-1)b^t=0\qquad(0\le t\le q-1).
\]
The coefficient matrix is the Vandermonde matrix formed by all \(q\) distinct base-field elements;
its determinant is the product of their pairwise differences and is nonzero.
Hence every \(N_b-1\) is zero in \(k\), giving the integer congruence \(N_b\equiv1\pmod p\).
All \(N_b\) are therefore at least one.
There are \(q\) of them, and their sum is \(q\), so every \(N_b=1\).
Thus each base-field value has exactly one preimage, and \(F\) is a permutation polynomial.
\end{proof}

\section{The AGW Criterion and Additive and Multiplicative Fibers}
\label{sec:trace-permutations}
\label{sec:agw}
Descending a map to a smaller set separates the permutation property into two conditions: images of distinct fibers must remain separate, and no collisions may occur within a fiber.
This set-theoretic principle is called the \textbf{AGW criterion}; see \cite{AGW,Feng}.
Trace fibers are cosets of an additive subgroup, whereas nonzero power-map fibers are cosets of a multiplicative subgroup; the two constructions therefore have the same form of criterion.

\subsection{A Fiber Criterion for Finite Sets}
\begin{lemma}[AGW criterion]\label{lem:fiber-permutation}
Let \(A,S,\overline S\) be finite sets with \(|S|=|\overline S|\), and suppose that
\[
\lambda:A\longrightarrow S,\qquad
\overline\lambda:A\longrightarrow\overline S
\]
are surjective. If maps \(F:A\to A\) and \(\varphi:S\to\overline S\) satisfy
\(\overline\lambda\circ F=\varphi\circ\lambda\), then \(F\) is a permutation if and only if \(\varphi\) is bijective and \(F\) is injective on every fiber \(\lambda^{-1}(s)\).
When these conditions hold, each restriction
\[
F:\lambda^{-1}(s)\longrightarrow
\overline\lambda^{-1}(\varphi(s))
\]
is bijective.
\end{lemma}
\begin{proof}
If \(F\) is bijective, then \(\overline\lambda F\) is surjective, and commutativity implies that \(\varphi\) is surjective.
Since the two base sets are finite and have equal size, \(\varphi\) is bijective; the restriction of \(F\) to each fiber is certainly injective.
Conversely, if \(F(x)=F(y)\), applying \(\overline\lambda\) to both sides gives
\(\varphi(\lambda(x))=\varphi(\lambda(y))\). By injectivity of \(\varphi\), the elements \(x,y\) lie in the same fiber;
injectivity within that fiber then gives \(x=y\). Hence \(F\) is injective and therefore bijective on the finite set \(A\).
Finally, each point in the target fiber has an \(F\)-preimage, which commutativity and injectivity of \(\varphi\) force to lie in \(\lambda^{-1}(s)\). Thus the restriction is also surjective.
\end{proof}

The commutativity relation is expressed by the following diagram:
\begin{center}
\begin{tikzpicture}[>=Stealth,every node/.style={font=\normalsize},line width=0.6pt]
\node (a) at (0,1.5) {$A$};
\node (b) at (3.8,1.5) {$A$};
\node (s) at (0,0) {$S$};
\node (t) at (3.8,0) {$\overline S$};
\draw[->] (a)--node[above] {$F$}(b);
\draw[->] (s)--node[above] {$\varphi$}(t);
\draw[->] (a)--node[left] {$\lambda$}(s);
\draw[->] (b)--node[right] {$\overline\lambda$}(t);
\end{tikzpicture}
\end{center}
When \(S=\overline S\) and \(\lambda=\overline\lambda=\pi\), this is the fiber permutation criterion compatible with \(\pi F=\varphi\pi\).
To compute the cycles of \(F\) as well, choose one fiber over each cycle of \(\varphi\) and apply Lemma~\ref{lem:fiber-return}:
the total cycle length is the product of the base-cycle length and the fiber return period.

\subsection{Constructing Permutations Along Trace Fibers}
\begin{theorem}\label{thm:trace-permutations}
Let \(E=\F_{q^n}\), \(k=\F_q\), and \(a\in k^\times\), and let \(g:k\to E\) be any function. Put
\[
F(x)=ax+g(\Tr_{E/k}(x)),\qquad
\varphi(t)=at+\Tr_{E/k}(g(t)).
\]
Then \(F\) permutes \(E\) if and only if \(\varphi\) permutes \(k\). In this case,
\[
F^{-1}(z)=a^{-1}\bigl(z-g(\varphi^{-1}(\Tr_{E/k}(z)))\bigr).
\]
\end{theorem}
\begin{proof}
Linearity of the trace gives \(\Tr\circ F=\varphi\circ\Tr\). The trace is surjective, and every fiber has \(q^{n-1}\) elements.
On the fiber \(\Tr(x)=t\), \(F\) is the affine map \(x\mapsto ax+g(t)\), which is injective because \(a\ne0\),
and hence maps bijectively onto the corresponding target fiber. Lemma~\ref{lem:fiber-permutation} gives the permutation criterion.
If \(F(x)=z\), then \(\Tr(x)=\varphi^{-1}(\Tr(z))\); substituting into the definition and solving for \(x\) gives the inverse formula.
\end{proof}

By interpolation, \(g\) can always be represented by a polynomial with coefficients in \(E\), so this construction does give permutation polynomials.
If \(g\) takes values in \(k\), then \(\varphi(t)=at+ng(t)\). In particular, when \(p\mid n\), every such \(g\) gives a permutation. This is a nonlinear application of \(k\subseteq\ker\Tr\) in Proposition~\ref{prop:trace-splitting}.

\par\Needspace{5\baselineskip}
\noindent\phantomsection\label{v98:par-prescribed-trace-permutation}\textbf{Choosing parameters from a base-field permutation.}
One may instead prescribe a permutation \(\varphi\) of \(k\) and a scalar \(a\in k^\times\), and then construct \(g\). By surjectivity of the trace, choose \(\theta\in E\) with \(\Tr_{E/k}(\theta)=1\). For any function \(u:k\to\ker\Tr_{E/k}\), put
\[
g(t)=\theta\bigl(\varphi(t)-at\bigr)+u(t).
\]
Then \(\Tr_{E/k}(g(t))=\varphi(t)-at\), so Theorem~\ref{thm:trace-permutations} gives the permutation \(F(x)=ax+g(\Tr_{E/k}(x))\).
For fixed \(\theta,a,\varphi\), every \(g\) satisfying this trace condition has the displayed form; \(u\) changes the actual maps between fibers without changing \(\varphi\) on the base field.
This choice does not require division by \(n\), and therefore also applies when \(p\mid n\).

Along a cycle of \(\varphi\) of length \(\ell\), choose an initial trace fiber and an origin in it. The first return then has the form
\(v\mapsto a^\ell v+\delta\), where \(v,\delta\in\ker\Tr_{E/k}\).
The translation term records the accumulated values of \(g\) along the base cycle, so the slope \(a^\ell\) alone does not determine all cycles.
Exercise~\ref{v70:ex-additive-return} gives this term explicitly and determines the resulting cycle decomposition.

\begin{example}
Over \(\F_{16}/\F_2\), \(F(x)=x+\Tr(x)=x^2+x^4+x^8\) is an involutory permutation.
Indeed, \(\Tr(F(x))=\Tr(x)\), so \(F(F(x))=x\). If the extension degree four is replaced by any odd integer,
the same expression \(x+\Tr(x)\) always has trace zero and cannot permute the extension field.
\end{example}

\subsection{Multiplicative Fibers and Explicit Return Maps}
\begin{proposition}[Multiplicative fiber criterion]\label{prop:multiplicative-agw}
Let \(q-1=ds\), where \(d,s\) are positive integers, and take \(r\ge1\) and \(g\in\F_q[X]\). Put
\[
F(X)=X^r g(X^s),\qquad H(Z)=Z^r g(Z)^s.
\]
Then \(F\) permutes \(\F_q\) if and only if
\[
\gcd(r,s)=1,\qquad H\text{ permutes }\mu_d.
\]
The second condition includes the requirement that \(g(z)\ne0\) for every \(z\in\mu_d\).
\end{proposition}
\begin{proof}
Since \(r\ge1\), we have \(F(0)=0\). If some \(z\in\mu_d\) satisfies \(g(z)=0\), choose \(x\in\F_q^\times\) with \(x^s=z\). Then \(F(x)=F(0)\), so \(F\) is not a permutation.
Now assume that \(g\) is nonzero everywhere on \(\mu_d\). Then \(F\) preserves \(\F_q^\times\), and
\[
\pi:\F_q^\times\longrightarrow\mu_d,\qquad \pi(x)=x^s
\]
is surjective and satisfies \(\pi F=H\pi\). For \(z\in\mu_d\), choose \(x_0\) with \(x_0^s=z\). Then
\[
\pi^{-1}(z)=x_0\mu_s,\qquad
F(x_0u)=F(x_0)u^r\quad(u\in\mu_s).
\]
Thus \(F\) is injective on this fiber if and only if \(u\mapsto u^r\) is injective on the cyclic group \(\mu_s\) of order \(s\), equivalently, \(\gcd(r,s)=1\).
Lemma~\ref{lem:fiber-permutation} gives the permutation criterion on the nonzero elements; adjoining the fixed point \(0\) proves the result.
\end{proof}

Choosing a representative \(x_0\) of the source fiber and a representative \(F(x_0)\) of the target fiber, we may compare the two constructions as follows:
\begin{center}
\small
\begin{tabular}{@{}lll@{}}
\toprule
Projection and base set&Source fiber&Map in relative coordinates\\
\midrule
\(\Tr_{E/k}:E\to k\)&\(x_0+\ker\Tr_{E/k}\)&\(v\mapsto av\)\\
\(x\mapsto x^s:\F_q^\times\to\mu_d\)&\(x_0\mu_s\)&\(u\mapsto u^r\)\\
\bottomrule
\end{tabular}
\end{center}
The AGW criterion determines whether either type of map is a permutation; return maps further determine its cycles.
In the multiplicative case, \(H\) determines how the fibers are permuted, and \(\gcd(r,s)=1\) ensures that the power action within each fiber is invertible.
After iteration along a base cycle, the map within the fiber is again a composition of multiplication by a fixed element and a power map.

\begin{proposition}[Return formula on multiplicative fibers]\label{prop:multiplicative-permutation-return}
Use the notation of Proposition~\ref{prop:multiplicative-agw} and assume that \(F\) is a permutation.
Let \(z_0,z_1,\ldots,z_{\ell-1}\) be a cycle of \(H\) of length \(\ell\), ordered by \(z_{j+1}=H(z_j)\), with indices taken modulo \(\ell\).
Take \(x_0^s=z_0\) and put
\[
C=\prod_{j=0}^{\ell-1}g(z_j)^{r^{\ell-1-j}},
\qquad
\delta=Cx_0^{r^\ell-1}.
\]
Then \(\delta\in\mu_s\), and in the coordinate \(x=x_0u\), the first-return map \(F^\ell\) is
\begin{equation}\label{eq:multiplicative-permutation-return}
R(u)=\delta u^{r^\ell}\qquad(u\in\mu_s).
\end{equation}
If \(\zeta\) generates \(\mu_s\) and we write \(\delta=\zeta^b\), then \(u=\zeta^v\) transforms \(R\) into
\[
v\longmapsto r^\ell v+b\pmod s.
\]
Each cycle of \(R\) of length \(t\) corresponds to a cycle of \(F\) of length \(\ell t\); the total number of \(F\)-cycles above this base cycle equals the number of cycles of \(R\).
\end{proposition}
\begin{proof}
On the fiber defined by \(x^s=z_0\), repeated substitution into \(F(x)=x^r g(x^s)\) gives
\[
F^j(x)=x^{r^j}\prod_{i=0}^{j-1}g(z_i)^{r^{j-1-i}}.
\]
Taking \(j=\ell\) gives \(F^\ell(x)=Cx^{r^\ell}\). Since the base cycle returns at step \(\ell\), \(F^\ell(x_0)^s=x_0^s\),
so \(\delta=F^\ell(x_0)/x_0\in\mu_s\). Substituting \(x=x_0u\) now gives~\eqref{eq:multiplicative-permutation-return} and the formula in exponent coordinates.
The cycle correspondence and length formula follow directly from Lemma~\ref{lem:fiber-return}.
\end{proof}

The affine action in exponent coordinates is defined over the residue class ring \(\Z/s\Z\); \(s\) need not be prime. When \(s\) is composite, this ring is not a field.
Changing the fiber representative \(x_0\) changes only the coordinates of the return map, not its cycle lengths.

\begin{example}[Two six-cycles above the same base cycle]\leavevmode\phantomsection\label{exm:agw-f13-return}
In \(\F_{13}\), take \(d=3\), \(s=4\), \(r=3\), and
\[
g(Z)=-Z^2+5Z+11,\qquad
F(X)=-X^{11}+5X^7+11X^3.
\]
Since \(\mu_3=\{1,3,9\}\), and \(g(1)=2\), \(g(3)=4\), \(g(9)=1\), we have
\[
H(1)=3,\qquad H(3)=9,\qquad H(9)=1.
\]
As \(\gcd(3,4)=1\), the map \(F\) is a permutation. Take \(z_0=x_0=1\); then \(\ell=3\) and
\[
\delta=2^9\cdot4^3=8\quad\text{in }\F_{13}.
\]
The group \(\mu_4\) is generated by \(5\), and \(8=5^3\), so the return map in exponent coordinates is
\[
v\longmapsto3v+3\pmod4.
\]
Its cycles are \((0,3)\) and \((1,2)\), so the nonzero elements form two \(F\)-cycles of length six.
Explicitly, they are
\[
(1,2,6,8,10,9),\qquad (3,4,12,11,7,5),
\]
and there is also the fixed point \(0\). These lengths are the products of the three-step base cycle and the two-step return cycles.
Figure~\ref{v103:fig-f13-return} arranges the first six-cycle in three columns according to \(x^4\). Three steps return to the starting fiber; the points \(1\) and \(8\) each need two such internal returns to return to themselves.
\begin{figure}[!htbp]
\centering
\begin{tikzpicture}[every node/.style={font=\small},
 pt/.style={circle,draw=Jade,fill=white,minimum size=8mm},
 ar/.style={-{Stealth[length=2mm]},semithick,Jade}]
\foreach \x in {0,3.5,7} {\fill[CoverMint,rounded corners] (\x-.65,-.65) rectangle (\x+.65,2.65);}
\node[pt] (a) at (0,2) {\(1\)};
\node[pt] (b) at (3.5,2) {\(2\)};
\node[pt] (c) at (7,2) {\(6\)};
\node[pt] (d) at (0,0) {\(8\)};
\node[pt] (e) at (3.5,0) {\(10\)};
\node[pt] (f) at (7,0) {\(9\)};
\draw[ar] (a)--node[above]{\(F\)}(b);
\draw[ar] (b)--node[above]{\(F\)}(c);
\draw[ar] (c.south)--(7,1)--(0,1)--(d.north);
\draw[ar] (d)--node[below]{\(F\)}(e);
\draw[ar] (e)--node[below]{\(F\)}(f);
\draw[ar] (f.east)--(8.1,0)--(8.1,3.2)--(0,3.2)--(a.north);
\draw[<->,Clay,dashed,semithick] (a.west) to[bend right=40] node[left]{\(F^3\)} (d.west);
\node[pt] (h1) at (0,-2) {\(1\)};
\node[pt] (h3) at (3.5,-2) {\(3\)};
\node[pt] (h9) at (7,-2) {\(9\)};
\foreach \x/\n in {0/h1,3.5/h3,7/h9} {\draw[ar,MutedSlate,dashed] (\x,-.7)--node[right]{\(x\mapsto x^4\)} (\n);}
\draw[ar] (h1)--node[above]{\(H\)}(h3);
\draw[ar] (h3)--node[above]{\(H\)}(h9);
\draw[ar] (h9.south)--(7,-2.9)--node[below]{\(H\)}(0,-2.9)--(h1.south);
\end{tikzpicture}
\caption{A base cycle of length three and an internal return of length two produce the six-cycle \((1,2,6,8,10,9)\). Each green column contains its two points in one multiplicative fiber.}\label{v103:fig-f13-return}
\end{figure}
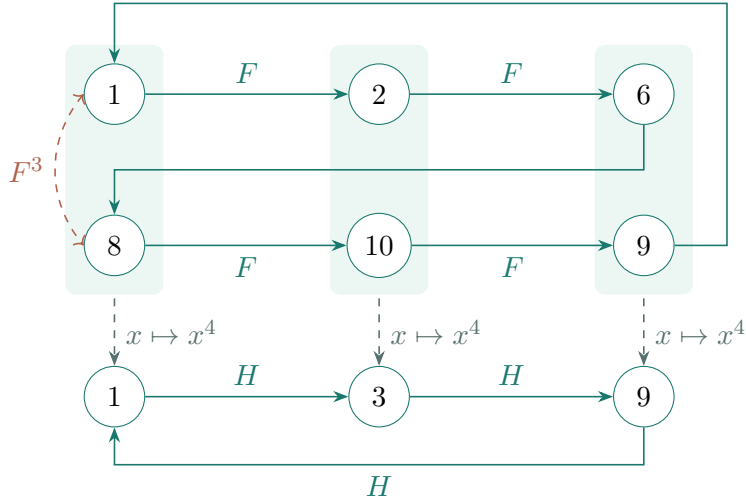
\FloatBarrier

\end{example}

An invertible linear map on the quotient can also produce a nonlinear permutation with a sparse inverse.
The following example specializes the lifting construction of~\cite[Theorem 3.1]{Li2026SparseLifts}.

\begin{example}[A nonlinear trinomial involution from a trace map]\label{v68:ex-sparse-involution}
Take \(\F_8\subset\F_{64}\), and denote the trace-zero root of Example~\ref{ex:normal-f8-example} by \(\eta\), so that \(\eta^3+\eta+1=0\).
Its zero trace is now used to construct an involution and then lift it along multiplicative fibers.
On \(\F_8\), put
\[
L(z)=z+\Tr_{\F_8/\F_2}(\eta z)
=(1+\eta)z+\eta^2z^2+\eta^4z^4.
\]
If \(c=\Tr_{\F_8/\F_2}(\eta z)\in\F_2\), then
\(\Tr_{\F_8/\F_2}(\eta(z+c))=c+c\Tr_{\F_8/\F_2}(\eta)=c\).
Thus \(L(L(z))=z+c+c=z\), so \(L\) is an involution.

Now define
\begin{equation}\label{v68:eq-sparse-involution}
F(X)=(1+\eta)X+\eta^2X^{37}+\eta^4X^{46}\in\F_{64}[X].
\end{equation}
For \(x\in\F_{64}^\times\), set \(t=x^{36}\).
The power map \(\pi:x\mapsto x^{36}\) has image \(\F_8^\times\), since its image has order \(63/\gcd(63,36)=7\).
Moreover, \(36\equiv1\pmod7\), so \(\pi\) restricts to the identity on \(\F_8^\times\).
For \(t\in\F_8^\times\), write
\[
R(t)=\frac{L(t)}t=(1+\eta)+\eta^2t+\eta^4t^3.
\]
This value is nonzero because \(L\) is bijective and \(L(0)=0\).
As \(x^{109}=x^{46}\) on \(\F_{64}^\times\), formula~\eqref{v68:eq-sparse-involution} gives
\[
F(x)=xR(t),\qquad
\pi(F(x))=tR(t)^{36}=tR(t)=L(t).
\]
Hence
\[
F(F(x))
=x\frac{L(t)}t\frac{L(L(t))}{L(t)}=x.
\]
The same identity holds at \(x=0\), so \(F^{-1}=F\) as functions on \(\F_{64}\).
Its canonical representative has degree less than \(64\), with nonzero terms of degrees \(37\) and \(46\), neither a power of two.
By the representation of linear maps in Chapter~\ref{chap:linearized}, \(F\) is therefore not \(\F_2\)-linear.
The inverse is obtained by transporting the linear involution \(L\) through the multiplicative quotient \(\pi\), while each fiber is mapped by multiplication by the nonzero scalar \(R(t)\).
\end{example}

\paragraph{Hidden linearization in permutation binomials.}
\begin{proposition}\label{v106:prop-hidden-linearization}
Let \(E=\F_{q^e}\), \(e\ge2\), \(a\in E^\times\), and \(r\ge1\), and put \(\ell_j=(q^j-1)/(q-1)\). Suppose that
\[
\gcd(r,q-1)=1,\qquad(-a)^{\ell_e}\ne1,
\]
and that some \(1\le h<e\) satisfies \(\gcd(h,e)=1\) and \(r\ell_h\equiv1\pmod{\ell_e}\). Then
\[
F(X)=X^r(X^{q-1}+a)
\equiv (X^{q^h}+aX)\circ X^r\pmod{X^{q^e}-X}
\]
permutes \(E\). If \(\rho r\equiv1\pmod{q^e-1}\), \(1\le\rho<q^e-1\), and \(L(X)=X^{q^h}+aX\), its inverse function is \(F^{-1}(y)=(L^{-1}(y))^\rho\).
\end{proposition}
\begin{proof}
The exponent congruence gives \(r(q^h-1)\equiv q-1\pmod{q^e-1}\), so the two polynomials agree at every nonzero element, and both vanish at zero. Also \(r\) is coprime to \(\ell_e\); together with \(\gcd(r,q-1)=1\), this gives \(\gcd(r,q^e-1)=1\), making the power map invertible.
The nonzero kernel equation \(L(z)=0\) is \(z^{q^h-1}=-a\). Since \(\gcd(q^h-1,q^e-1)=q-1\), the power-image criterion in a cyclic group says that this equation has a solution exactly when \((-a)^{\ell_e}=1\). Thus \(L\) has no nonzero kernel and is an invertible linearized map. Reversing the composition gives the inverse formula.
\end{proof}
\begin{example}[A two-sided inverse over the field of order 27]\label{v106:ex-binomial-f27}
Take \(q=3,e=3,h=2,r=23,a=1\). We have \(\ell_2=4\), \(\ell_3=13\), \(23\cdot4\equiv1\pmod{13}\), and \((-1)^{13}\ne1\). Hence
\[
F(X)=X^{25}+X^{23},\qquad F(x)=L(x^{23}),\qquad L(z)=z^9+z.
\]
On \(\F_{27}\), put \(\sigma(z)=z^3\). Since \(\sigma^3=1\), the operator identity
\[
(1+\sigma^2)(2+2\sigma+\sigma^2)=1
\]
gives \(L^{-1}(y)=2y+2y^3+y^9\). Now \(23\cdot17\equiv1\pmod{26}\) yields
\[
\boxed{F^{-1}(y)=(2y+y^9+2y^3)^{17}.}
\]
Both operators are polynomials in \(\sigma\), so the identity also holds in the opposite order. The power maps are likewise mutual inverses, proving both inverse compositions.
\end{example}
These are congruences in the finite-field function algebra. A factorization into two permutations does not make \(F\) and \(L\) conjugate as dynamical systems, and does not identify their cycle lengths.
The displayed conditions are sufficient for permutation. Fan's preprint~\cite{Fan2026Binomials} further reports that, within the same family with \(e\ge2\) and \(a\ne0\), they are also necessary.

\section{The Trace Obstruction and Cycle Decomposition of Affine Frobenius Maps}
\label{afr:sec-affine-periods}

The preceding section constructed permutations on additive and multiplicative fibers and determined their cycles through return maps. We now compute all cycle lengths of affine Frobenius permutations, use cyclic modules to treat general linear permutations, and express their permutation signs in terms of determinants.
Fix \(k=\F_q\), \(E=\F_{q^n}\), and \(1\le s<n\), and retain from Section~\ref{wfr:sec-weighted} the notation
\[
 d=\gcd(n,s),\quad m=n/d,\quad K=\F_{q^d},\quad \tau(x)=x^{q^s}.
\]
\leavevmode\phantomsection\label{afr:def-affine-frobenius}\ignorespaces
For \(a\in E^\times\) and \(b\in E\), define the \textbf{affine Frobenius permutation}
\[
 F_{a,b}(x)=a\tau(x)+b=S_a(x)+b,\qquad c=\Nm_{E/K}(a).
\]
Here affine refers to the \(K\)-vector space structure of \(E\); we study the iterates of this permutation on \(E\).

\par\Needspace{12\baselineskip}
\subsection{Norms and Changes of Origin}

\begin{lemma}[Hilbert 90 for any generating automorphism]\label{afr:lem-hilbert90}
Let \(H\subseteq E\) be finite fields, and let \(\gamma\) generate \(\Aut_H(E)\). Then
\[
 \im(\gamma-\id)=\ker\Tr_{E/H},\qquad
 \{\gamma(u)/u:u\in E^\times\}=\ker\Nm_{E/H}.
\]
\end{lemma}
\begin{proof}
The powers of \(\gamma\) run through all conjugates, so the difference has trace zero and the quotient has norm one.
The kernel of the first map is \(H\), and that of the second homomorphism is \(H^\times\). Comparing dimensions and group orders, respectively,
and using surjectivity of the trace and norm proves both equalities. This is the proof of Theorem~\ref{thm:hilbert90} and does not depend on the chosen generating automorphism.
\end{proof}

\begin{proposition}[Linearization and the trace obstruction]\label{afr:prop-normalization}
The affine permutation \(F_{a,b}\) is described by the following three cases.
\begin{enumerate}
\item If \(c\ne1\), then \(F_{a,b}\) has a unique fixed point; taking this point as the origin conjugates \(F_{a,b}\) to \(S_a\).
\item If \(c=1\), choose \(u\in E^\times\) such that \(a\tau(u)=u\), and put
\begin{equation}\label{afr:eq-trace-obstruction}
 \beta=b/u,\qquad \theta=\Tr_{E/K}(\beta).
\end{equation}
The substitution \(x=uy\) transforms \(F_{a,b}\) into \(G(y)=\tau(y)+\beta\). If \(\theta=0\), then \(G\) is conjugate to \(\tau\) by a translation.
\item In the notation of the second case, if \(\theta\ne0\), then \(G\) has no fixed point, and
\begin{equation}\label{afr:eq-iterate-order}
 G^m(y)=y+\theta,\qquad \ord(G)=pm.
\end{equation}
\end{enumerate}
Whether \(\theta\) vanishes is independent of the choice of \(u\) satisfying the condition.
\end{proposition}
\begin{proof}
If \(S_a(x)=x\) and \(x\ne0\), then \(a=x/\tau(x)\), and taking norms gives \(c=1\).
Thus, when \(c\ne1\), the map \(\id-S_a\) is invertible, and the fixed-point equation \((\id-S_a)x=b\) has a unique solution \(x_0\).
The first assertion follows from \(F_{a,b}(x_0+y)=x_0+S_a(y)\).

If \(c=1\), Lemma~\ref{afr:lem-hilbert90} supplies the required \(u\), and
\(u^{-1}F_{a,b}(uy)=\tau(y)+b/u\).
The latter has a fixed point if and only if \(y-\tau(y)=\beta\) is solvable, equivalently, \(\theta=0\); taking such a fixed point as the origin yields \(\tau\).
Two admissible choices of \(u,u'\) satisfy \(\tau(u'/u)=u'/u\), so \(u'/u\in K^\times\); the corresponding traces differ only by a nonzero \(K\)-scalar.

Finally, \(G^m(y)=\tau^m(y)+\sum_{i=0}^{m-1}\tau^i(\beta)=y+\theta\).
When \(\theta\ne0\), this translation has order exactly \(p\). If \(G^t=\id\), comparison of the linear parts gives \(m\mid t\);
writing \(t=mh\) and then using \(G^{mh}(y)=y+h\theta\) gives \(p\mid h\). Hence \(\ord(G)=pm\).
\end{proof}

\subsection{Fixed Points of Iterates and All Cycles}

Below let \(t\ge1\) and \(g=\gcd(m,t)\), and write \(K_g=\F_{q^{dg}}\).
The automorphism \(\tau^t\) has fixed field \(K_g\) and generates the automorphism group of \(E/K_g\), which has order \(m/g\).

\begin{lemma}[Relative sums and products in iterates]\label{afr:lem-iteration-count}
For \(a\in E^\times\) and \(\beta\in E\), write \(c=\Nm_{E/K}(a)\) and \(\theta=\Tr_{E/K}(\beta)\), and put
\[
 a_t=\prod_{i=0}^{t-1}\tau^i(a),\qquad
 \beta_t=\sum_{i=0}^{t-1}\tau^i(\beta).
\]
Then
\begin{equation}\label{afr:eq-iteration-count}
 \Nm_{E/K_g}(a_t)=c^{t/g},\qquad
 \Tr_{E/K_g}(\beta_t)=(t/g)\theta.
\end{equation}
The integer \(t/g\) in the second formula is interpreted in the characteristic of the field.
\end{lemma}
\begin{proof}
Expanding the relative norm or trace using the generator \(\tau^t\) gives the exponents
\(jt+i\pmod m\), where \(0\le j<m/g\) and \(0\le i<t\).
For fixed \(i\), \(jt\pmod m\) runs through all multiples of \(g\) exactly once.
In \(0,\ldots,t-1\), each residue class modulo \(g\) occurs exactly \(t/g\) times,
so each exponent modulo \(m\) occurs \(t/g\) times in total. Multiplication and addition give the respective formulas.
\end{proof}

\begin{theorem}[Fixed-point counts for affine Frobenius maps]\label{afr:thm-fixed-counts}
Put \(N_t(F)=\#\{x\in E:F^t(x)=x\}\). Then
\begin{equation}\label{afr:eq-linear-fixed}
 N_t(S_a)=
 \begin{cases}
 q^{dg},&c^{t/g}=1,\\
 1,&c^{t/g}\ne1.
 \end{cases}
\end{equation}
For \(F=F_{a,b}\), all fixed-point counts are given by
\begin{equation}\label{afr:eq-affine-fixed}
 N_t(F)=
 \begin{cases}
 N_t(S_a),&c\ne1,\\
 q^{dg},&c=1,\ \theta=0,\\
 q^{dg},&c=1,\ \theta\ne0,\ p\mid t/g,\\
 0,&c=1,\ \theta\ne0,\ p\nmid t/g.
 \end{cases}
\end{equation}
In the last three cases, \(\theta\) is defined by~\eqref{afr:eq-trace-obstruction}.
\end{theorem}
\begin{proof}
Since \(S_a^t=a_t\tau^t\), a nonzero fixed point exists if and only if \(\Nm_{E/K_g}(a_t)=1\),
by the multiplicative form of Lemma~\ref{afr:lem-hilbert90}. Lemma~\ref{afr:lem-iteration-count} reduces this condition to \(c^{t/g}=1\).
If there is a nonzero solution \(x_0\), the nonzero solutions are exactly \(x_0K_g^\times\); together with zero, there are \(q^{dg}\) solutions. Otherwise zero is the only solution. This proves~\eqref{afr:eq-linear-fixed}.

Proposition~\ref{afr:prop-normalization} already handles the conjugacies for \(c\ne1\) and \(c=1,\theta=0\).
In the remaining case, \(G^t(y)=\tau^t(y)+\beta_t\). By additive Hilbert 90,
the fixed-point equation has a solution if and only if \(\Tr_{E/K_g}(\beta_t)=0\), that is, \((t/g)\theta=0\).
Since \(\theta\ne0\), this is equivalent to \(p\mid t/g\). When solvable, its solution set is an additive coset of \(K_g\), so there are \(q^{dg}\) solutions.
\end{proof}

\begin{corollary}[The characteristic-power part of periods with nonzero trace]\label{afr:cor-period-jump}
If \(c=1\) and \(\theta\ne0\), then every cycle length \(\ell\) of \(F_{a,b}\) satisfies
\begin{equation}\label{afr:eq-period-jump}
 \boxed{\quad v_p(\ell)=v_p(m)+1.\quad}
\end{equation}
\end{corollary}
\begin{proof}
From \(\ord(F_{a,b})=pm\), we have \(\ell\mid pm\). For \(t=\ell\), the fixed-point set is nonempty,
so~\eqref{afr:eq-affine-fixed} gives \(p\mid\ell/\gcd(m,\ell)\), that is, \(v_p(\ell)>v_p(m)\).
Together with the divisibility bound, this gives equality.
\end{proof}

\leavevmode\phantomsection\label{afr:def-cycle-count}\ignorespaces
Let \(C_t(F)\) denote the number of cycles of length exactly \(t\). A cycle of length \(e\) contributes \(e\) fixed points to the \(t\)th iterate if and only if \(e\mid t\). Thus
\begin{equation}\label{afr:eq-cycle-count}
 N_t(F)=\sum_{e\mid t}eC_e(F),\qquad
 C_t(F)=\frac1t\sum_{e\mid t}\mu(t/e)N_e(F).
\end{equation}
Theorem~\ref{afr:thm-fixed-counts} therefore determines all cycle lengths and their multiplicities.

\begin{corollary}[Shortest cycles and their number]\label{afr:cor-shortest-cycles}
Retain the notation of this section. If \(c=1\) and \(\theta\ne0\), put \(\nu=v_p(m)\). Then the shortest cycle length of \(F_{a,b}\) and the number of such cycles are, respectively,
\begin{equation}\label{afr:eq-shortest-cycles}
\boxed{\qquad
\ell_{\min}=p^{\nu+1},\qquad
C_{\ell_{\min}}(F_{a,b})=\frac{q^{d p^\nu}}{p^{\nu+1}}.
\qquad}
\end{equation}
If \(c\ne1\), or if \(c=1\) and \(\theta=0\), the shortest cycle length is one.
\end{corollary}
\begin{proof}
In the case \(c=1\) and \(\theta\ne0\), take \(t=p^{\nu+1}\). Then
\(g=\gcd(m,t)=p^\nu\), hence \(p\mid t/g\).
Theorem~\ref{afr:thm-fixed-counts} gives
\[
N_t(F_{a,b})=q^{d p^\nu}>0.
\]
The cycle lengths \(\ell\) of these fixed points divide \(t\), whereas Corollary~\ref{afr:cor-period-jump} gives
\(v_p(\ell)=\nu+1\). Thus all these cycle lengths equal \(t\).
Cycles of length \(t\) therefore exist, their number is \(N_t/t\), and all other cycle lengths are also multiples of \(t\).
Hence \(t\) is the shortest cycle length.
In the other two cases, Proposition~\ref{afr:prop-normalization} gives a fixed point, so the shortest cycle length is one.
\end{proof}

Lemma~\ref{lem:affine-least-period} gives a least period dividing every cycle length of a general affine permutation;
the trace obstruction here further determines that period and the corresponding number of cycles.
Here we compare all cycles in the entire \(\F_p\)-affine space;
the affine compression of Chapter~\ref{chap:period-layers} instead examines a return transformation within a single fiber.
When the group acts freely on that fiber, the return period equals the order of the corresponding group element; see Proposition~\ref{prop:affine-return-order}.
These two computations use, respectively, the cycle decomposition of the whole space and a group action on a single fiber.

\begin{example}[Two eight-cycles on 16 elements]\label{afr:ex-eight-cycles}
Take \(E=\F_{16}\), \(k=\F_2\), and \(\tau(x)=x^2\), and choose \(\beta\) satisfying \(\Tr_{E/k}(\beta)=1\).
Then \(G(x)=x^2+\beta\) satisfies \(G^4(x)=x+1\) and \(G^8(x)=x\).
Here \(m=4\) and \(\nu=2\). By Corollary~\ref{afr:cor-shortest-cycles}, the shortest cycle length is \(2^3=8\),
and the number of such cycles is \(2^{2^2}/2^3=2\). These two cycles already contain all 16 elements.
The same squaring map without the translation term has order four; the nonzero trace raises the order of the affine permutation to eight, and every point in this example has period eight.
\end{example}

The norm determines the fixed-point condition for the weighted linear part, and the trace determines whether the translation can be removed when the norm is one.
These conditions come from products and sums of conjugates, respectively, and~\eqref{afr:eq-cycle-count} converts them into the complete cycle decomposition.

\subsection{Cycle Data for General Linear Permutations}
\label{v53:subsec-linear-cycles}

Let \(V\ne0\) be a finite-dimensional \(k=\F_q\)-space, let \(A:V\to V\) be invertible, and suppose that \(p=\operatorname{char}k\).
Letting \(T\) act as \(A\) makes \(V\) a \(k[T]\)-module.
A cyclic vector identifies an invariant subspace with a cyclic module given by a polynomial quotient ring. After decomposition into an invariant direct sum of such subspaces,
the fixed-point counts of all iterates can be computed directly from the multiplicities of irreducible factors.

\begin{lemma}[Invariant cyclic decomposition]\label{v53:lem-cyclic-decomposition}
There exist monic irreducible polynomials \(P_1,\ldots,P_r\in k[T]\) and positive integers \(e_1,\ldots,e_r\) such that
\begin{equation}\label{v53:eq-elementary-divisors}
 (V,A)\simeq\bigoplus_{i=1}^r
 \bigl(k[T]/(P_i^{e_i}),M_{\overline T}\bigr).
\end{equation}
The \(P_i\) may repeat, and \(P_i\ne T\). Such a decomposition is called a \textbf{cyclic decomposition};
the \(P_i^{e_i}\) are the \textbf{elementary divisors} of \(A\), and their multiset is uniquely determined by \(A\).
\end{lemma}
\begin{proof}
Put \(m=m_A\). By Lemma~\ref{lem:cyclic-vector}, choose \(v\in V\) such that \(\mu_v=m\).
Then \(U=k[A]v\) and \(R=k[T]/(m)\) are isomorphic via \(\bar a\mapsto a(A)v\).
Write \(\alpha=\overline T\), and let \(\lambda:R\to k\) extract the coefficient of \(T^{\deg m-1}\) from the canonical remainder.
Extend the functional \(a(A)v\mapsto\lambda(\bar a)\) on \(U\) to a \(k\)-linear functional \(\ell:V\to k\).

The nondegenerate coefficient pairing of Theorem~\ref{thm:algebra-trace-coefficient} ensures that
for every \(x\in V\) there is a unique \(\pi(x)\in R\) such that
\[
 \lambda(\bar a\,\pi(x))=\ell(a(A)x)\qquad(\bar a\in R).
\]
The right-hand side is independent of the representative of \(a\), since \(m(A)=0\). Uniqueness implies that \(\pi\) is linear and
\[
 \pi(Ax)=\alpha\pi(x),\qquad \pi(a(A)v)=\bar a.
\]
Thus the restriction of \(\pi\) to \(U\) is the inverse of the preceding isomorphism, \(\ker\pi\) is an \(A\)-invariant subspace,
and \(V=U\oplus\ker\pi\). Applying induction to \(\ker\pi\) gives a direct sum of cyclic modules.
Factoring each defining polynomial into irreducibles and applying the Chinese remainder theorem yields~\eqref{v53:eq-elementary-divisors}.
Finally, invertibility of \(A\) forces each defining polynomial to have nonzero constant term, so \(P_i=T\) does not occur.
The multiset of elementary divisors is independent of the decomposition: for each irreducible polynomial \(P\),
\[
 \frac{\dim_k\ker P(A)^j-\dim_k\ker P(A)^{j-1}}{\deg P}
 =\#\{i:P_i=P,\ e_i\ge j\}\qquad(j\ge1).
\]
Here \(P(A)^0=I\). The equality follows by direct computation on each component, and its right-hand side determines all exponents and their multiplicities.
\end{proof}

\begin{theorem}[Elementary divisors and cycle counts]\label{v53:thm-linear-cycle-count}
With the notation of~\eqref{v53:eq-elementary-divisors}, put \(\theta_i=\ord(P_i)\).
For a positive integer \(t\), put
\begin{equation}\label{v53:eq-linear-fixed-dimension}
 \beta_A(t)=\sum_{\theta_i\mid t}
 (\deg P_i)\min\{e_i,p^{v_p(t)}\}.
\end{equation}
Then the number of fixed points of \(A^t\) and the number of cycles of \(A\) of length exactly \(t\) are, respectively,
\begin{equation}\label{v53:eq-linear-cycle-count}
 N_t(A)=q^{\beta_A(t)},\qquad
 C_t(A)=\frac1t\sum_{d\mid t}\mu(t/d)q^{\beta_A(d)}.
\end{equation}
Moreover, putting \(b=\min\{j\ge0:p^j\ge\max_i e_i\}\) gives
\begin{equation}\label{v53:eq-linear-order}
 \ord(A)=p^b\lcm_i\theta_i.
\end{equation}
\end{theorem}
\begin{proof}
First compute the kernel of multiplication by \(T^t-1\) on a component \(k[T]/(P^e)\).
If \(h=\min\{e,v_P(T^t-1)\}\), where \(v_P\) denotes the multiplicity of the irreducible factor \(P\),
then this kernel is the ideal \((\overline{P^{e-h}})\), of dimension \(h\deg P\).
Indeed, write \(T^t-1=P^s u\) and \(P\nmid u\). Then \(\bar u\) is a unit,
and the annihilation condition \(P^e\mid P^s a\) is equivalent to \(P^{e-h}\mid a\).

Write \(t=p^a u\) and \(p\nmid u\). Then
\[
 T^t-1=(T^u-1)^{p^a}.
\]
The polynomial \(T^u-1\) has no repeated factor; by \(\theta=\ord(P)\) and \(p\nmid\theta\),
\(P\mid T^u-1\) if and only if \(\theta\mid t\).
Thus \(h=\min\{e,p^{v_p(t)}\}\) or zero, according as \(\theta\mid t\) holds or fails.
Summing over all components gives \(\dim_k\ker(A^t-I)=\beta_A(t)\) and hence the fixed-point count.
Substitution into~\eqref{afr:eq-cycle-count} gives the cycle count.

Finally, \(A^t=I\) if and only if every elementary divisor divides \(T^t-1\).
Theorem~\ref{thm:order-multiplicity} gives~\eqref{v53:eq-linear-order}.
\end{proof}

If a cyclic module decomposition \((V,A)\simeq\bigoplus_j k[T]/(f_j)\) with \(f_j(0)\ne0\) is known, the same kernel computation gives
\begin{equation}\label{v53:eq-cyclic-gcd-count}
 \beta_A(t)=\sum_j\deg\gcd(f_j,T^t-1).
\end{equation}
Indeed, putting \(g=\gcd(f_j,T^t-1)\) makes \(f_j\mid(T^t-1)u\) equivalent to \(f_j/g\mid u\),
and these residue classes form a \(\deg g\)-dimensional space. This formula retains repeated-factor information and does not require \(f_j\) to be square-free.

If \(x\in V\) has minimal annihilating polynomial \(\mu_x\), then \(A^t x=x\) if and only if \(\mu_x\mid T^t-1\),
so its least period is \(\ord(\mu_x)\); the zero vector has period one.
Consequently, removing all powers of \(p\) from any cycle length \(t\) leaves the least common multiple of some of the \(\theta_i\);
the least common multiple of the empty set is taken to be one. The \(p\)-power part is determined by the multiplicities in the annihilating polynomial of that vector.
Lemma~\ref{lem:cyclic-vector} also gives a cycle of length \(\ord(A)\).

\begin{example}[The same characteristic polynomial with different cycles]\label{v53:ex-characteristic-not-cycles}
Take \(k=\F_3\) and \(E=\F_{27}\), and compare
\[
 L(x)=x+\Tr_{E/k}(x),\qquad \sigma(x)=x^3.
\]
The trace map \(N=\Tr_{E/k}\) has rank one, and \(N^2=0\):
for \(c\in k\), we have \(\Tr_{E/k}(c)=3c=0\).
Thus the elementary divisors of \(L=I+N\) are \((T-1)^2,T-1\),
whereas \(\sigma\) is represented by the elementary divisor \(T^3-1=(T-1)^3\) of the normal-basis cyclic module.
Both have characteristic polynomial \((T-1)^3\) and order three, but
\[
\begin{array}{c|cc|cc}
 &N_1&N_3&C_1&C_3\\\hline
 L&9&27&9&6\\
 \sigma&3&27&3&8
\end{array}
\]
Here the fixed space of \(L\) is \(\ker\Tr\), of dimension two, whereas that of \(\sigma\) is \(k\), of dimension one.
The order of the entire map gives only the least common multiple of all periods, and the characteristic polynomial gives only the total multiplicity of each irreducible factor;
the distribution of those multiplicities among the cyclic modules determines the full cycle structure.
\end{example}

For an invertible linearized map \(L\) on \(E=\F_{q^n}\),
Proposition~\ref{prop:dickson-rank} relates its base-field matrix to \(D_L\). Therefore
\begin{equation}\label{v53:eq-dickson-iterate-kernel}
 \beta_L(t)=n-\rank_E(D_L^t-I).
\end{equation}
Extending scalars does not change matrix rank. The rank on the right is computed over \(E\), while the left side remains the fixed-space dimension over \(k\).
Substituting these integers into~\eqref{v53:eq-linear-cycle-count} gives all cycle counts.
If the linearized coefficients of \(L=a(\sigma)\) belong to \(k\), no matrix is needed.
Since \(L^t-I=(a(T)^t-1)(\sigma)\), Theorem~\ref{thm:linearized-gcd} directly gives
\begin{equation}\label{v53:eq-base-coefficient-cycles}
 \beta_L(t)=\deg\gcd(T^n-1,a(T)^t-1).
\end{equation}
Substituting this into~\eqref{v53:eq-linear-cycle-count} again gives all cycle counts.

\subsection{Translation obstructions and general affine cycles}
\phantomsection\label{v113:par-affine-question}
The preceding argument determined all cycle counts from the elementary divisors of the linear part. Are these data still sufficient after adding a translation? On \(\F_p\), the maps \(x\mapsto x\) and \(x\mapsto x+1\) have the same linear part, but the former has \(p\) fixed points whereas the latter is one \(p\)-cycle: its \(t\)th iterate is \(x+t\), and the first return occurs at \(t=p\). Thus even the full classification of the linear part must be supplemented by information about the translation.

For \(\Phi(x)=Ax+b\), changing the origin by \(x=y+c\) replaces the translation term by \(b+(A-I)c\). With \(A\) fixed, the information that cannot be removed by changing the origin is therefore the class of \(b\) in \(V/(A-I)V\). This class is zero exactly when \(\Phi\) has a fixed point, in which case a change of origin turns \(\Phi\) into a linear map. For a nonzero class, determining cycle lengths further requires knowing which iterates have fixed points.

The trace condition for weighted Frobenius and the image conditions of Chapter~\ref{chap:root-fibers} are instances of this elimination problem. For a general linear part, we examine the images of \(\ker(A-I)^j\) in the quotient above and record the first level containing the obstruction class. This integer will control compatibility for all iterates at once.

\begin{definition}\label{v107:def-translation-height}
Let \(A\in\operatorname{GL}(V)\), where \(V\) is a finite-dimensional \(\F_q\)-space, and put \(N=A-I\), \(b\in V\). Define the \textbf{obstruction height} of the translation relative to \(A\) by
\[
h_A(b)=\min\{j\ge0:b\in\im N+\ker N^j\},\qquad \ker N^0=\{0\}.
\]
By elementary divisors, \(N\) is invertible on the primary components other than \(T-1\), and nilpotent on the \((T-1)\)-primary component. Thus the minimum exists and is at most \(\dim V\).
It is unchanged by a change of basis or by replacing \(b\) with \(b+Nc\). The condition \(h_A(b)=0\) says exactly that the affine map has a fixed point.
\end{definition}

\begin{theorem}[Exact cycle counts for a general affine permutation]\label{v107:thm-affine-cycles}
Let \(\Phi(x)=Ax+b\), \(p=\Char\F_q\), and use \(\beta_A(t)\) from~\eqref{v53:eq-linear-fixed-dimension}. For \(t\ge1\),
\begin{equation}\label{v107:eq-affine-cycles}
\boxed{N_t(\Phi)=
\begin{cases}
q^{\beta_A(t)},&p^{v_p(t)}>h_A(b),\\
0,&p^{v_p(t)}\le h_A(b),
\end{cases}\qquad
C_t(\Phi)=\frac1t\sum_{d\mid t}\mu(t/d)N_d(\Phi).}
\end{equation}
Consequently, the least cycle length, which divides every cycle length, is
\[
p^{\lceil\log_p(h_A(b)+1)\rceil}.
\]
\end{theorem}
\begin{proof}
Decompose \(V=V_1\oplus V_*\), with \(V_1\) the \((T-1)\)-primary component. Since \(N\) is invertible on \(V_*\), a change of origin removes the translation in that component.
Compatibility need therefore be tested only on \(V_1\), and \(h_A(b)\) depends only on this component as well.

Write \(t=ru\), where \(r=p^{v_p(t)}\) and \(p\nmid u\). On \(V_1\),
\[
A^t-I=N^rU(N),\qquad
I+A+\cdots+A^{t-1}=N^{r-1}U(N),
\]
where \(U(Z)=(((1+Z)^u-1)/Z)^r\) and \(U(0)=u^r\ne0\).
Nilpotence of \(N\) makes \(U(N)\) invertible. Hence \(\Phi^t(x)=x\) is solvable exactly when
\[
N^{r-1}b\in\im N^r
\quad\Longleftrightarrow\quad
b\in\im N+\ker N^{r-1}
\quad\Longleftrightarrow\quad h_A(b)<r.
\]
For the middle equivalence, write \(N^{r-1}b=N^ry\) and rearrange to obtain \(N^{r-1}(b-Ny)=0\).
When solvable, the solution set is an affine coset of \(\ker(A^t-I)\), with size given by the preceding theorem; otherwise it is empty. M\"obius inversion gives the cycle count.
Put \(e=\lceil\log_p(h_A(b)+1)\rceil\). Compatibility shows that every cycle length is divisible by \(p^e\). Since \(\Phi^{p^e}\) has a fixed point, some cycle length also divides \(p^e\), and must therefore equal it.
\end{proof}

\begin{corollary}[When a translation increases the order]\label{v110:cor-affine-order}
Put $r=\ord(A)$ and $h=h_A(b)$. Then
\[
\ord(\Phi)=\operatorname{lcm}\left(r,p^{\lceil\log_p(h+1)\rceil}\right)
=\begin{cases}
pr,&h=p^{v_p(r)},\\
r,&h<p^{v_p(r)}.
\end{cases}
\]
Thus translation can increase the order only by a factor of the characteristic. A nonremovable translation need not increase it.
\end{corollary}
\begin{proof}
The identity $\Phi^t=\id$ holds precisely when every point is fixed. By Theorem~\ref{v107:thm-affine-cycles}, this is equivalent to $A^t=I$ and $p^{v_p(t)}>h$. Their least common solution gives the lcm formula.
The height $h$ is at most the largest exponent of a $(T-1)$ elementary divisor, or zero when none occurs. Every elementary-divisor exponent is at most $p^{v_p(r)}$ by~\eqref{v53:eq-linear-order}. Hence $0\le h\le p^{v_p(r)}$, and the two cases follow.
The condition that translation cannot be removed is only $h>0$; increasing the order requires the stronger equality displayed above.
\end{proof}

\begin{example}[Three cycle types with one linear part]\label{v107:ex-affine-three-types}
On \(V=\F_2^3\), take
\[
A=\begin{pmatrix}1&1&0\\0&1&0\\0&0&1\end{pmatrix},\qquad
\im(A-I)=\langle e_1\rangle,\quad
\ker(A-I)=\langle e_1,e_3\rangle.
\]
Since \((A-I)^2=0\), three translations give
\[
\begin{array}{c|c|rrr|l}
b&h_A(b)&N_1&N_2&N_4&\text{cycle decomposition}\\\hline
0&0&4&8&8&\text{four fixed points, two 2-cycles}\\
e_3&1&0&8&8&\text{four 2-cycles}\\
e_2&2&0&0&8&\text{two 4-cycles}
\end{array}
\]
The linear part is identical in all three cases. The obstruction height decides which linear fixed spaces of iterates become nonempty affine fibers.
\end{example}

For \(r=p^j\), the compatibility condition in the proof is exactly the image criterion of Proposition~\ref{prop:affine-least-period-computation}.
The present theorem also uses \(\beta_A(t)\) to count every compatible fiber, passing from the least period to the full cycle decomposition.
If \(A\) is a relative Frobenius of order \(m\) and the translation has nonzero relative trace, all \((T-1)\)-primary blocks have length \(p^{v_p(m)}\), and at least one translation residue is nonzero. Thus \(h_A(b)=p^{v_p(m)}\). This recovers the characteristic-power jump in Corollary~\ref{afr:cor-period-jump}; the elementary divisors still determine the individual fixed-point counts.

\subsection{The Sign and Determinant of a Linear Permutation}
\label{v53:subsec-linear-parity}

\leavevmode\phantomsection\label{v53:def-permutation-sign}\ignorespaces
A permutation of a finite set is a product of disjoint cycles. The \textbf{permutation sign} of a cycle of length \(r\) is \((-1)^{r-1}\);
the sign of the whole permutation is the product of its cycle signs, with value \(1\) or \(-1\), corresponding to an even or odd permutation, respectively.
If the finite set has \(Q\) elements and \(c_r\) denotes the number of cycles of length \(r\), then
\[
\operatorname{sgn}(f)=(-1)^{\sum_r(r-1)c_r}
=(-1)^{Q-\sum_rc_r}.
\]
The sign also equals the determinant over \(\mathbb Q\) of the associated permutation matrix, and is therefore multiplicative under composition.
For linear permutations, this cycle invariant can also be determined directly by linear algebra.

\begin{theorem}[The sign of a linear permutation]\label{v53:thm-linear-sign}
Let \(V\) be an \(n\)-dimensional vector space over \(\F_q\), with \(n\ge1\) and \(A\in\operatorname{GL}(V)\).
Regard \(A\) as a permutation of the \(q^n\) elements of \(V\).
If \(q\) is odd, define the quadratic character \(\chi_q:\F_q^\times\to\{1,-1\}\)
to have value \(1\) on nonzero squares and \(-1\) on nonsquares. Then
\begin{equation}\label{v53:eq-linear-sign}
\boxed{\qquad\operatorname{sgn}(A)=\chi_q(\det A).\qquad}
\end{equation}
If \(q\) is even, all invertible linear maps are even permutations, except when \((q,n)=(2,2)\);
for \((q,n)=(2,2)\), the odd permutations are exactly the nonidentity linear maps of order \(2\).
\end{theorem}
\begin{proof}
Choose a basis of \(V\). Elementary row operations show that invertible matrices are generated by two types of maps:
shears \(x_i\mapsto x_i+c x_j\) (\(i\ne j\)) and diagonal scalings \(x_i\mapsto a x_i\) (\(a\ne0\)).
Coordinate transpositions can also be expressed using these two types; for example, on the corresponding two coordinates,
\[
\begin{pmatrix}0&1\\1&0\end{pmatrix}
=\begin{pmatrix}1&1\\0&1\end{pmatrix}
 \begin{pmatrix}1&0\\-1&1\end{pmatrix}
 \begin{pmatrix}1&1\\0&1\end{pmatrix}
 \begin{pmatrix}-1&0\\0&1\end{pmatrix}.
\]
Both permutation sign and determinant are multiplicative under composition, so it suffices to compute them on these generators.

Suppose that \(q\) is odd. A nonidentity shear has order equal to the characteristic \(p\), so its cycle lengths are \(1\) or \(p\), its sign is \(1\), and its determinant is also \(1\).
Choose a generator \(z\) of \(\F_q^\times\). The multiplication map \(x\mapsto zx\) fixes zero and forms one cycle of length \(q-1\) on the nonzero elements, with sign \(-1\).
Thus the sign of multiplication by \(a=z^j\) is \((-1)^j=\chi_q(a)\).
A diagonal scaling induces this multiplication permutation on each of \(q^{n-1}\) parallel coordinate lines.
Since \(q^{n-1}\) is odd, its sign is again \(\chi_q(a)\), the quadratic character of its determinant.
This proves~\eqref{v53:eq-linear-sign}.

Suppose that \(q\) is even. A diagonal scaling has order dividing the odd integer \(q-1\), so it is an even permutation.
When \(n\ge2\) and \(c\ne0\), a shear fixes all points satisfying \(x_j=0\), and the remaining points form cycles of length \(2\).
The number of these cycles is
\[
\frac{q^n-q^{n-1}}2=\frac{q^{n-1}(q-1)}2.
\]
This number is odd if and only if \((q,n)=(2,2)\).
For \(n=1\), only diagonal scalings occur, and all are even permutations.
Finally, the action of \(\operatorname{GL}_2(\F_2)\) on the three nonzero vectors gives an isomorphism with the symmetric group on three points:
the action is faithful, and both groups have order \(6\). Its odd permutations are therefore exactly its nonidentity elements of order \(2\).
\end{proof}

\begin{corollary}[The Dickson determinant and permutation sign]\label{v53:cor-dickson-sign}
Let \(L:E\to E\) be an invertible \(q\)-linearized map, where \(E=\F_{q^n}\). If \(q\) is odd, then
\[
\operatorname{sgn}(L)=\chi_q(\det D_L).
\]
In particular, for the multiplication map \(M_a:x\mapsto ax\),
\[
\operatorname{sgn}(M_a)
=\chi_q\bigl(\Nm_{E/\F_q}(a)\bigr)=\chi_{q^n}(a),\qquad a\in E^\times.
\]
\end{corollary}
\begin{proof}
Proposition~\ref{prop:dickson-rank} shows that \(\det D_L\) is the determinant of \(L\) over the base field.
Applying Proposition~\ref{prop:characteristic-multiplication} to the multiplication map gives the norm formula.
Finally, \(\Nm(a)=a^{(q^n-1)/(q-1)}\), so
\(\Nm(a)^{(q-1)/2}=a^{(q^n-1)/2}\), proving that the two quadratic character values are equal.
\end{proof}

Thus, in odd characteristic, the parity information in the complete cycle decomposition depends only on the square class of the determinant.
For the permutation of \(\F_{q^2}\) induced by \(L(x)=ax+bx^q\), the norm difference in Example~\ref{exm:dickson-quadratic} also gives
\[
\operatorname{sgn}(L)=\chi_q\bigl(\Nm(a)-\Nm(b)\bigr).
\]
Computing the sign therefore does not require finding all cycle lengths first.

If \(L=\sum_{i=0}^{n-1}a_i\sigma^i=g(\sigma)\) and all \(a_i\in\F_q\),
then in a normal basis the matrix of \(\sigma\) is the cyclic shift matrix \(S\) of length \(n\), so the matrix of \(L\) is \(g(S)\).
In the cyclic module \(k[T]/(T^n-1)\), \(g(S)\) corresponds to multiplication by \(\overline{g(T)}\).
Corollary~\ref{prop:local-multiplicity-determinant} and the \hyperref[def:resultant-root-product]{root-product formula for the resultant} give
\[
\det D_L=\det g(S)=\operatorname{Res}(T^n-1,g(T))=\prod_{\omega^n=1}g(\omega),
\]
where roots of unity are counted with their multiplicities in \(T^n-1\), so the formula also applies when the characteristic divides \(n\).
In odd characteristic, taking the quadratic character of this product gives the sign of \(L\).
In particular, \(\sigma^s\) forms \(\gcd(n,s)\) cycles on the normal-basis coordinates, so
\[
\det\sigma^s=(-1)^{n-\gcd(n,s)},\qquad
\operatorname{sgn}(\sigma^s)=\chi_q\bigl((-1)^{n-\gcd(n,s)}\bigr).
\]
Here the determinant of \(\sigma\) as a linear operator and its permutation sign on all field elements are related by the quadratic character.

\begin{corollary}[The sign of an affine permutation]\label{v53:cor-affine-sign}
If \(|V|>2\), then for every \(A\in\operatorname{GL}(V)\) and \(b\in V\),
\[
\operatorname{sgn}(x\mapsto Ax+b)=\operatorname{sgn}(A).
\]
In particular, if \(q\) is odd, the affine Frobenius map of this section satisfies
\[
\operatorname{sgn}(F_{a,b})
=\chi_q\bigl((-1)^{n-d}\Nm_{E/k}(a)\bigr),\qquad d=\gcd(n,s).
\]
\end{corollary}
\begin{proof}
If \(b\ne0\), all cycles of the translation \(x\mapsto x+b\) have length equal to the characteristic \(p\), and there are \(|V|/p\) of them.
For odd \(p\) its sign is \(1\); for \(p=2\) its sign is \((-1)^{|V|/2}\), which equals \(-1\) only when \(|V|=2\).
The identity translation also has sign \(1\), so the first formula follows from multiplicativity under composition.
Now express \(F_{a,b}\) as the composition of a translation, a multiplication map, and \(\tau=\sigma^s\).
In a normal basis, the matrix of \(\sigma^s\) forms \(d\) cycles on the \(n\) basis vectors, so its determinant is \((-1)^{n-d}\), proving the second formula.
\end{proof}

\begin{example}[All linear permutations of \(\F_4\)]\label{v53:ex-F4-sign}
Write \(\F_4=\F_2(\alpha)\), where \(\alpha^2+\alpha+1=0\).
Every \(\F_2\)-linear map has a unique expression \(L(x)=ax+bx^2\), with Dickson determinant \(a^3+b^3\).
Every nonzero element of \(\F_4\) satisfies \(u^3=1\), so \(L\) is invertible if and only if exactly one of \(a,b\) is nonzero.
The three maps \(x\mapsto ax\) are even permutations, and the three maps \(x\mapsto bx^2\) are odd permutations.
The latter satisfy \(L^{\circ2}(x)=b^3x^4=x\), and each induces a transposition on the four field elements.
Here the base field of the linear structure is \(\F_2\); with base field \(\F_4\) instead, the only one-dimensional linear permutations are \(x\mapsto ax\), all even.
\end{example}

\section{Dickson Polynomials: Power Maps Descending Along Conjugate Fibers}
\label{sec:dickson-polynomials}
The behavior of a power map on a cyclic group is determined by its exponent. Identifying two mutually inverse roots 
as a single value allows this group action to descend to a nonlinear permutation of a finite field.

Given a quadratic polynomial, raise both of its roots, counted with
multiplicity, to the same power. The sum of the resulting roots,
expressed as a polynomial in the sum and product of the original roots,
is a Dickson polynomial.
We first consider the case in which the product of the two roots is one.

\begin{definition}\leavevmode\phantomsection\label{def:dickson-polynomial}
The \textbf{Dickson polynomials} of the first kind with parameter one are defined by the recurrence
\[
D_0(X)=2,\qquad D_1(X)=X,\qquad
D_{m+1}(X)=XD_m(X)-D_{m-1}(X)\quad(m\ge1).
\]
After reducing the integer coefficients, we use the same notation over a base field of any characteristic.
For \(m\ge1\), \(D_m\) is monic of degree \(m\).
\end{definition}
The first few terms are \(D_2=X^2-2\), \(D_3=X^3-3X\), and \(D_4=X^4-4X^2+2\).
The recurrence directly gives, for every nonzero element \(u\) of any extension field,
\begin{equation}\label{eq:dickson-root-parameter}
\boxed{\quad D_m(u+u^{-1})=u^m+u^{-m}.\quad}
\end{equation}
Indeed, the identity holds for \(m=0,1\), and substitution of two consecutive terms into the recurrence cancels the two middle terms.

\begin{lemma}[Root parameters for inversion orbits]\label{lem:dickson-quotient}
In \(\F_{q^2}^\times\), put \(U=\mu_{q-1}\cup\mu_{q+1}\) and define \(\pi(u)=u+u^{-1}\).
Then \(\pi:U\to k\) is surjective, and
\[
\pi(u)=\pi(v)\quad\Longleftrightarrow\quad v=u\ \text{or}\ v=u^{-1}.
\]
Thus each fiber is exactly one orbit of the involution \(u\mapsto u^{-1}\); the fiber is a singleton when \(u=u^{-1}\).
\end{lemma}
\begin{proof}
If \(u\in\mu_{q-1}\), then \(u^q=u\); if \(u\in\mu_{q+1}\), then \(u^q=u^{-1}\).
In both cases, \(\pi(u)^q=\pi(u)\), so \(\pi(u)\in k\).
Conversely, for \(x\in k\), the quadratic polynomial \(Y^2-xY+1\) splits in \(\F_{q^2}\),
with nonzero roots \(u,u^{-1}\). Frobenius permutes these two roots, so \(u^q=u\) or \(u^q=u^{-1}\). Hence \(u\in U\) and \(\pi(u)=x\).
Finally, multiplying \(\pi(u)=\pi(v)\) by \(uv\) gives the equivalent equation \((u-v)(uv-1)=0\), proving the fiber description.
This calculation also covers the case of a repeated root of the quadratic polynomial.
\end{proof}

\begin{theorem}[Dickson permutation criterion]\label{thm:dickson-permutation}
For a positive integer \(m\), the polynomial \(D_m\) permutes \(\F_q\) if and only if
\[
\boxed{\gcd(m,q^2-1)=1.}
\]
If \(a\ge1\) satisfies \(ma\equiv1\pmod{q^2-1}\), then \(D_a\) represents its inverse function.
\end{theorem}
\begin{proof}
Equation~\eqref{eq:dickson-root-parameter} states that \(D_m\pi=\pi\rho_m\), where \(\rho_m(u)=u^m\).
If \(\gcd(m,q^2-1)=1\), choose \(a\) as in the statement. Every element of \(U\) satisfies \(u^{q^2-1}=1\),
so \(\rho_m\) and \(\rho_a\) are mutually inverse on \(U\), and both commute with inversion.
Consequently, for every \(x=\pi(u)\in k\),
\[
D_a(D_m(x))=\pi(u^{ma})=\pi(u)=x,
\]
and the composition in the reverse order gives the same result, proving the permutation property and the inverse formula.

If \(\gcd(m,q^2-1)>1\), choose a common prime divisor \(\ell\). It divides \(q-1\) or \(q+1\),
so the corresponding cyclic group contains an element \(\zeta\ne1\) of order \(\ell\).
By \(\zeta^m=1\) and the fiber description, \(\pi(\zeta)\ne\pi(1)=2\), but
\[
D_m(\pi(\zeta))=2=D_m(2).
\]
Thus \(D_m\) is not injective. The argument remains valid in characteristic two, where \(2=0\).
\end{proof}

For example, over \(\F_8\), the polynomial \(D_5=X^5+X^3+X\) is a permutation polynomial because \(\gcd(5,63)=1\).
Here \(q-1\) comes from the multiplicative group of the base field, and \(q+1\) from the norm-one group of the quadratic extension.
On the root parameters, Frobenius either fixes \(u\) or exchanges \(u,u^{-1}\), and the power map descends simultaneously along both types of fibers.
The self-reciprocal constructions of Chapter~\ref{chap:root-fibers} use the same root parameter \(u+u^{-1}\): there the question is transitivity of the return map,
whereas here it is whether taking powers preserves all inversion orbits.
The polynomial \(X^3+X\) in Theorem~\ref{v55:thm-char2-target-cubic} is exactly \(D_3\) in characteristic two; the noncube condition in the same reciprocal coordinate ensures dynamical irreducibility with the prescribed target, as stated in that theorem.
For related classical results, see \cite[Chapter 7, Section 2]{LN}.
The Dickson polynomials of this section, the Dickson matrices of Section~\ref{sec:dickson},
and the Dickson invariants of Chapter~\ref{chap:period-layers} are distinct objects.

\subsection{Reciprocal-Orbit Root Polynomials and a Frobenius Equation}\label{v97:sec-reciprocal-dickson}
The root parameter above identifies a root of unity with its inverse.
Forming a new root polynomial from these parameter values allows us to compare the Frobenius equations arising from roots of unity of different orders.

For odd \(m=2r+1\), define
\begin{equation}\label{v97:def-reciprocal-orbit-polynomial}
 P_m(Y)=1+\sum_{j=1}^{r}D_j(Y),\qquad P_1(Y)=1.
\end{equation}
Formula~\eqref{eq:dickson-root-parameter} gives the Laurent polynomial identity
\begin{equation}\label{v97:eq-reciprocal-orbit-polynomial}
 \frac{X^m-1}{X-1}=X^rP_m(X+X^{-1}).
\end{equation}
When \(p\nmid m\), the \(r\) distinct roots of \(P_m\) are the sums over the two-element reciprocal orbits of the nonidentity \(m\)th roots of unity, with each pair counted once.
We therefore call \(P_m\) the corresponding \textbf{reciprocal-orbit root polynomial}.
For \(m>1\), \(D_r\) is monic, so \(P_m\) is monic of degree \(r\); moreover, \(z+z^{-1}=w+w^{-1}\) holds precisely when \(w=z\) or \(z^{-1}\).
Since \(m\) is odd, inversion fixes no nonidentity root, giving exactly \(r\) roots.

\begin{theorem}[Root-set decomposition of a Dickson--Frobenius equation]\label{v97:thm-dickson-frobenius}
For every odd prime power \(q\), the following identity holds in \(\F_q[Y]\):
\[
 \boxed{Y^q-Y^2+2=(Y-2)P_{q+2}(Y)P_{q-2}(Y).}
\]
The nonconstant factors on the right are squarefree and pairwise coprime; they need not be irreducible.
\end{theorem}
\begin{proof}
If \(z^{q+2}=1\), then \(z^q=z^{-2}\); if \(z^{q-2}=1\), then \(z^q=z^2\).
In either case, \(y=z+z^{-1}\) satisfies
\[
 y^q=z^2+z^{-2}=y^2-2=D_2(y).
\]
After excluding \(z=1\), the two sets of reciprocal-orbit sums have \((q+1)/2\) and \((q-3)/2\) elements, none equal to \(2\).
They are disjoint: equality between two parameter values would make the corresponding roots equal or inverse, with order dividing both \(q+2\) and \(q-2\).
These odd integers are coprime, leaving only the excluded root \(1\).
Including \(2\), we obtain \(q\) distinct roots. Both sides are monic of degree \(q\), proving the identity.
Squarefreeness also follows directly: the derivative of the left side is \(-2Y\), while its value at \(0\) is \(2\ne0\).
\end{proof}

For \(q=3\), one has \(P_{q-2}=P_1=1\), and the identity and root count still apply.
All roots are considered in a common algebraic closure; the equation \(y^q=D_2(y)\) does not require \(y\) to lie in the base field.
Likewise, \(z+z^{-1}\) is first a reciprocal-orbit sum.
It is the trace from the quadratic extension only when \(z\in\F_{q^2}\) and \(z^q=z^{-1}\).

\begin{corollary}[Explicit interpolation weights]\label{v97:cor-dickson-weights}
Put \(H=P_{q+2}\) and \(J=P_{q-2}\). For every root \(t\) of \(H\),
\[
 \boxed{H'(t)=-\frac{2t}{(t-2)J(t)}\ne0.}
\]
\end{corollary}
\begin{proof}
Differentiate the preceding identity and substitute \(H(t)=0\) to obtain
\(-2t=(t-2)H'(t)J(t)\).
Coprimality of the factors and \(t\ne0\) show that the denominators are nonzero.
\end{proof}

The factor \(1/H'(t)\) is exactly the normalization weight in the Lagrange interpolation of Chapter~\ref{chap:dual}.
Thus the root-set decomposition supplies not only the factors but also the weights by which root evaluations recover the highest coefficient; see Exercise~\ref{v97:ex-dickson-functional}.

\begin{example}[The roots need not lie in a quadratic extension]\label{v97:ex-dickson-five}
For \(q=5\), one has \(P_3=Y+1\) and \(P_7=Y^3+Y^2-2Y-1\), so
\[
 Y^5-Y^2+2=(Y-2)(Y+1)(Y^3+Y^2-2Y-1).
\]
The cubic takes the values \(4,4,2,4,1\) at \(0,1,2,3,4\), and is therefore irreducible.
Its roots lie in \(\F_{125}\), not in \(\F_{25}\).
On the root-of-unity parameters, \(5^3\equiv-1\pmod7\), whereas neither \(5\) nor \(5^2\) is congruent to \(\pm1\).
The reciprocal orbits of seventh roots thus have Frobenius period three, in agreement with the cubic factor.
\end{example}

\par\Needspace{6\baselineskip}
\section*{Exercises}
\phantomsection
\addcontentsline{toc}{section}{Exercises}

\begin{exercise}[Recovering subfield linearity from slopes]\label{v72:exer-parallel-subfield}
Let \(q=p^n\), let \(r\mid n\), and put \(K=\F_{p^r}\subseteq\F_q\).
Suppose \(f:\F_q\to\F_q\) satisfies~\eqref{v72:eq-parallel-difference}, and write
\(s(t)=(f(t)-f(0))/t\) for \(t\ne0\). Prove that
\(s(ct)=s(t)\) for every \(c\in K^\times\) and \(t\in\F_q^\times\) if and only if
\[
f(X)=b+\sum_{j=0}^{n/r-1}a_jX^{p^{rj}},\qquad a_j,b\in\F_q.
\]
\end{exercise}

\begin{exercise}[The subgroup determines the permitted Frobenius powers]\label{v72:exer-carlitz-mcconnel}
Let \(\zeta\) generate \(\F_{16}^\times\), and put
\(H_5=\langle\zeta^3\rangle\) and \(H_3=\langle\zeta^5\rangle\).
For each of these two subgroups \(H\), determine all functions \(f:\F_{16}\to\F_{16}\) satisfying
\[
\frac{f(x)-f(y)}{x-y}\in H\qquad(x\ne y).
\]
Count the functions and identify those satisfying \(f(0)=0\) and \(f(1)=1\).
Explain why a nontrivial power of Frobenius occurs for one subgroup but not the other.
\end{exercise}

\begin{exercise}[Difference quotients, difference products, and coset rigidity]\label{v74:ex-exponent-rigidity}
Let \(q=p^n\), let \(H\) be a proper subgroup of \(\F_q^\times\) of index \(d>1\), and take integers \(u\ge1\) and \(v\in\Z\) with \(\gcd(u,d)=1\).
Consider the condition
\begin{equation}\label{v74:eq-mixed-difference}
(f(x)-f(y))^u(x-y)^v\in H\qquad(x\ne y).
\end{equation}
Negative \(v\) is allowed because \(x-y\ne0\).
\begin{enumerate}
\item Prove that \(f\) is a permutation and that existence of such a function implies \(\gcd(v,d)=1\).
\item When \(v\equiv-u\pmod d\), reduce the condition directly to Carlitz--McConnel and determine all functions.
\item In general, prove that~\eqref{v74:eq-mixed-difference} holds if and only if
\[
\begin{gathered}
f(x)=ax^{p^j}+b,\\
a\in H,\qquad b\in\F_q,\qquad0\le j<n,\qquad d\mid up^j+v.
\end{gathered}
\]
Take \((u,v)=(1,-1)\) and \((1,1)\) to compare the difference-quotient and difference-product conditions.
For the subgroup of index three in \(\F_{16}^\times\), compare the allowed Frobenius powers.
\item Drop \(\gcd(u,d)=1\) and take \(u=d,v=0\). Show that every permutation then satisfies~\eqref{v74:eq-mixed-difference}, and explain the loss of rigidity.
\end{enumerate}
\emph{Hint:} Put \(G=\F_q^\times/H\), write \(\pi\) for the quotient map, and choose \(ur\equiv1\pmod d\).
First show that \(\pi(f(x)-f(y))=\pi(x-y)^{-rv}\).
For part (3), use the following classification of Lenstra \cite[Theorem~3]{Lenstra1990Automorphisms}:
if a permutation \(f\) transforms the cosets of nonzero differences by one fixed affine automorphism of \(G\),
that is, \(\pi(f(x)-f(y))=e_0\chi(\pi(x-y))\) with \(e_0\in G\) and \(\chi\in\Aut(G)\), then \(f(x)=ax^{p^j}+b\).
Here an affine automorphism of the quotient group means \(z\mapsto e_0\chi(z)\), not an affine map of the field.
Part (2) follows directly from Carlitz--McConnel.
\end{exercise}

\begin{exercise}[Frobenius shape factors]\label{v76:ex-frobenius-shape-factor}
Let \(q=p^n\), \(1\le j<n\), and let \(H<\F_q^\times\) have index \(d\). Put \(f(X)=aX^{p^j}+bX+c\), where \(a\ne0\).
For \(\lambda\notin\{0,1\}\), define the \textbf{Frobenius shape factor}
\begin{equation}\label{v76:def-frobenius-shape-factor}
C_j(\lambda)=\frac{\lambda^{p^j}-\lambda}{\lambda(\lambda-1)}
=\frac{\lambda^{p^j-1}-1}{\lambda-1}.
\end{equation}
\begin{enumerate}
\item Prove that \(K_{\lambda,f}(x,t)=aC_j(\lambda)t^{p^j-2}\), and identify the subfield corresponding to \(C_j(\lambda)=0\).
\item For a nonempty \(\Lambda\subseteq\F_q\setminus\{0,1\}\), prove that the same \(f\) satisfies the subgroup condition for every shape in \(\Lambda\) if and only if
\[
d\mid p^j-2,\qquad aC_j(\lambda)\in H\quad(\lambda\in\Lambda).
\]
Deduce that such an \(a\ne0\) exists if and only if the divisibility condition holds, all \(C_j(\lambda)\) are nonzero, and the shape factors lie in one \(H\)-coset. Explain why this is only a criterion within the specified Frobenius monomial family.
\end{enumerate}
\end{exercise}

\begin{exercise}[Recovery from several shapes: constant values and subgroup values]\label{v76:ex-multishape-recovery}
Let \(q=p^n\ge3\), let \(\varnothing\ne\Lambda\subseteq\F_q\setminus\{0,1\}\), and put \(K=\F_p(\Lambda)\).
\begin{enumerate}
\item Suppose that, for each \(\lambda\in\Lambda\), the value \(K_{\lambda,f}(x,t)\) is independent of \(x,t\), while its constant value \(a_\lambda\) may initially depend on \(\lambda\).
Prove that all \(a_\lambda\) are equal and that the solutions are exactly \(f(X)=aX^2+L(X)+b\), with \(L\) a \(K\)-linear map.
Prove that this condition forces every such \(f\) to have a polynomial representative of degree at most two if and only if \(K=\F_q\).
\item Let \(H<\F_q^\times\), and assume only \(K_{\lambda,f}(x,t)\in H\). Prove that the condition for one shape automatically covers
\[
\lambda,\quad1-\lambda,\quad\lambda^{-1},\quad
(1-\lambda)^{-1},\quad\frac{\lambda-1}{\lambda},\quad\frac{\lambda}{\lambda-1}.
\]
If \(\Lambda\) meets every orbit of these transformations on \(\F_q\setminus\{0,1\}\), prove that \(f(X)=aX^2+bX+c\), with \(a\in H\).
\item In \(\F_{16}\), take \(H=(\F_{16}^\times)^3\) and a \(\lambda\) of degree four over \(\F_2\). Use the preceding exercise to construct a nonquadratic function satisfying the subgroup condition for this shape. Explain why the subfield-generation criterion in part~(1) does not apply directly to part~(2).
\end{enumerate}
\emph{Hint:} In part~(1), compare the coefficients of \(X^2\) in the unique low-degree representative of \(f\); note that \(\F_p(\lambda)\) has at least three elements. In part~(2), use symmetry of second divided differences in the nodes.
\end{exercise}

\begin{exercise}[Linearization does not guarantee a minimal value set]\label{v58:ex-linearized-minimal-values}
Let \(E=\F_{q^n}\), and let \(L\in E[X]\) be a nonzero \(q\)-linearized polynomial of ordinary degree \(q^t\), with \(0\le t<n\).
Prove that the following conditions are equivalent: \(L\) is a minimal value set polynomial over \(E\);
\(|\ker(L|_E)|=q^t\); \(L\mid X^{q^n}-X\).
Then show that \(X^3-X\) is a minimal value set polynomial over \(\F_{27}\), whereas \(X^3+X\), of the same degree, is not.
\end{exercise}

\begin{exercise}[Value sets of one linearized polynomial over all extensions]\label{v70:ex-linearized-values-all-extensions}
Fix \(k=\F_q\) and \(0\ne h=\sum_i a_iT^i\in k[T]\), and put \(\mathcal L_h(X)=\sum_i a_iX^{q^i}\).
For each \(m\ge1\), write \(K_m=\F_{q^m}\) and \(d_m=\deg\gcd(h,T^m-1)\).
\begin{enumerate}
\item Prove that
\[
|\ker(\mathcal L_h|_{K_m})|=q^{d_m},\qquad
|\mathcal L_h(K_m)|=q^{m-d_m},\qquad
C_{\mathcal L_h|_{K_m}}=q^{m+d_m},
\]
and show that equality holds in the bound of Proposition~\ref{v58:prop-value-collision}.
\item For \(A=X^4+X^2+X\) from Example~\ref{v68:ex-linearized-iterate-f2} and the fixed polynomial \(A^{\circ3}\),
determine all \(m\) for which they permute \(\F_{2^m}\), and compute their kernel dimensions and value-set sizes for \(m=6,12\).
\item Explain why multiplying \(h\) by \(T^a\), where \(a\ge0\), changes none of these kernels, value-set sizes, or collision counts,
whereas a factor \(T^a\) in the minimal polynomial of a one-sided recurrence sequence can produce a transient.
\end{enumerate}
\end{exercise}

\begin{exercise}[Recovering a value set from multiple collisions]\label{v58:ex-factorial-value-count}
Let \(F\in\F_Q[X]\) have ordinary degree \(d\ge1\), and put \(m=\min\{d,Q\}\).
For each positive integer \(j\), let \(M_j(F)\) count the ordered \(j\)-tuples satisfying \(F(x_1)=\cdots=F(x_j)\),
with the \(x_i\) pairwise distinct; for \(j>m\) the count is zero. Prove the integer counting formula
\[
|\mathcal V_F|=\sum_{j=1}^m\frac{(-1)^{j-1}}{j!}M_j(F),
\qquad M_1(F)=Q,\qquad M_2(F)=C_F-Q.
\]
Compute these numbers in Example~\ref{v58:ex-F27-five-values} and verify the formula.
\end{exercise}

\begin{exercise}\label{ex:function-inverse}
Find the inverse polynomial of the permutation polynomial \(X^5\) over \(\F_7\), and explain why it represents the inverse function.
\end{exercise}

\begin{exercise}\label{ex:delta-idempotent}
In \(\F_q[X]/(X^q-X)\), prove \(\delta_a^2=\delta_a\), \(\delta_a\delta_b=0\) (\(a\ne b\)), and \(\sum_a\delta_a=1\).
\end{exercise}

\begin{exercise}[Frobenius in a function algebra]\label{ex:wu-function-frobenius}
Prove that every element \(h\) of \(R=\F_q[X]/(X^q-X)\) satisfies \(h^q=h\).
When \(q>1\), why is \(R\) not a field? Compare it with the field structure on the root set \(\F_q\) itself.
\end{exercise}

\begin{exercise}[A permutation group on the same trace fibers]\label{ex:wu2-trace-shears}
Let \(p\mid n\). For each function \(g:k\to k\), put \(F_g(x)=x+g(\Tr(x))\).
Prove \(F_g\circ F_h=F_{g+h}\), determine when these permutations are the identity, and find their possible orders.
\end{exercise}

\begin{exercise}[Affine return on a trace fiber]\label{v70:ex-additive-return}
Use the notation of Theorem~\ref{thm:trace-permutations}, assume that \(\varphi\) permutes \(k\), and put \(p=\Char(k)\).
Choose a cycle \(t_0,\ldots,t_{\ell-1}\) of \(\varphi\) of length \(\ell\), with \(t_{j+1}=\varphi(t_j)\), where indices are read modulo \(\ell\).
Choose \(x_0\in E\) with \(\Tr_{E/k}(x_0)=t_0\), and put
\[
W=\ker\Tr_{E/k},\qquad
B=\sum_{j=0}^{\ell-1}a^{\ell-1-j}g(t_j),\qquad
\delta=(a^\ell-1)x_0+B.
\]
Prove that \(\delta\in W\) and that, in the fiber coordinate \(x=x_0+v\), the first return map is
\[
R(v)=a^\ell v+\delta\qquad(v\in W).
\]
Determine all cycle lengths and their multiplicities for \(F\) on these \(\ell\) trace fibers in the three cases
\(a^\ell\ne1\), \(a^\ell=1\) with \(\delta=0\), and \(a^\ell=1\) with \(\delta\ne0\).
Explain why changing \(x_0\) preserves the cycle data, and treat the boundary case \(n=1\).
\end{exercise}

\begin{exercise}[The same base cycle with different return cycles]\label{ex:agw-return-comparison}
In Example~\ref{exm:agw-f13-return}, replace \(g\) by
\[
\widetilde g(Z)=-2Z^2+9Z+8,\qquad
\widetilde F(X)=X^3\widetilde g(X^4).
\]
Prove that \(\widetilde F\) is still a permutation of \(\F_{13}\) and induces the same three-cycle on \(\mu_3\).
Compute the return map on the fiber containing \(1\) and use it to determine all cycle lengths of \(\widetilde F\) and their multiplicities.
Explain why the permutation of the base set does not determine the cycle structure of the whole permutation.
\end{exercise}

\begin{exercise}[Dickson composition and inverse functions]\label{ex:wu10-dickson-inverse}
Prove \(D_m(D_r(X))=D_{mr}(X)\). Over \(\F_8\), find the polynomial representative of degree less than eight of the inverse function of \(D_5=X^5+X^3+X\).
Compare the identity \(D_m\circ D_r=D_{mr}\) with the inverse-function equality modulo \(X^8-X\).
\end{exercise}

\begin{exercise}[A translation of nonzero trace on 64 elements]\label{afr:ex-sixtyfour-cycles}
On \(E=\F_{64}\), take \(G(x)=x^2+\beta\), where \(\Tr_{E/\F_2}(\beta)=1\).
Without enumerating points, find all cycle lengths of \(G\) and the number of cycles of each length.
\end{exercise}

\begin{exercise}[All affine periods with the same linear part]\label{afr:ex-weighted-sixteen}
Retain Example~\ref{wfr:ex-sixteen} and take \(E=\F_{16}\), \(\tau(x)=x^4\), and \(a=\omega^2\), where \(\omega^2+\omega+1=0\).
For each \(b\in E\), find the fixed-point count and all cycle types of \(F_{a,b}(x)=ax^4+b\), and explain why the answer is independent of \(b\).
\end{exercise}

\begin{exercise}[Cyclic modules, root spaces, and permutation signs]\label{v53:ex-cycles-root-sign}
Let \(k=\F_3\) and \(f(T)=(T-1)^2(T+1)\), and let \(A\) be multiplication by \(\overline T\) on \(k[T]/(f)\).
Without enumerating points, find the order of \(A\), all its cycle lengths and their multiplicities, and compute its permutation sign using both cycle counts and the determinant.
Let \(W_f\) be the root space of \(\mathcal L_f(X)=X^{27}-X^9-X^3+X\) in a finite splitting extension.
Explain how these period data become Frobenius orbit data on \(W_f\),
and find the splitting field of \(\mathcal L_f\) over \(k\).
\end{exercise}

\begin{exercise}[Reciprocal-orbit roots and a coefficient functional]\label{v97:ex-dickson-functional}
Let \(q\) be odd, and put \(H=P_{q+2}\), \(J=P_{q-2}\), and \(N=(q+1)/2\), with \(P_m\) as in \eqref{v97:def-reciprocal-orbit-polynomial}.
For \(A=\F_q[Y]/(H)\), define \(\Lambda(R)=[Y^{N-1}]R\) using the unique representative \(R\) of degree less than \(N\).
As \(t\) ranges over all roots of \(H\) in a splitting field, prove
\[
 \Lambda(R)=-\frac12\sum_{H(t)=0}\frac{(t-2)J(t)}t\,R(t).
\]
Deduce the value of \(\sum_{H(t)=0}(t-2)J(t)t^{j-1}\) for \(0\le j<N\), and explain why the sum lies in \(\F_q\), although it is computed in a splitting field.
\end{exercise}

\begin{exercise}[R\'edei rational functions and the parity of extension degrees]\label{v106:ex-redei-rational}
Let \(q\) be odd, let \(\delta\in\F_q^\times\) be a nonsquare, take \(\theta^2=\delta\), and let \(n\ge1\). Put
\[
\phi(X)=\frac{X+\theta}{X-\theta},\qquad
R_n=\phi^{-1}\circ(Z\mapsto Z^n)\circ\phi.
\]
Prove that \(R_n\in\F_q(X)\) and that, for every \(m\ge1\),
\[
R_n\text{ permutes }\mathbb P^1(\F_{q^m})
\quad\Longleftrightarrow\quad\gcd\bigl(n,q^m-(-1)^m\bigr)=1.
\]
For \(q=7,\delta=3,n=3\), compute \(R_3\) and determine the extensions on which it remains a permutation. These R\'edei rational functions are distinct from the R\'edei root polynomials discussed earlier.
\end{exercise}

\begin{exercise}[From a global gap to linear stability]\label{v116:ex-stability-gap}
Fix \(d\ge2\). Suppose that uniform constants \(\varepsilon>0\), \(0<\eta<(1-1/d)/2\) satisfy, for all permutations over the prime fields under consideration,
\[
\mathcal E_H(f)\le\varepsilon q^2\quad\Longrightarrow\quad D_H(f)\le\eta q,
\qquad |H|=(q-1)/d.
\]
Use Proposition~\ref{v116:prop-local-stability} to prove
\[
\mathcal E_H(f)\ge\min\{\varepsilon,1-1/d\}\,qD_H(f).
\]
Conversely, how does a uniform linear lower bound give such a gap condition?
\end{exercise}

\chapter{Linear Recurrences and Recovery from Finite Observations}
\label{chap:sequences}
\phantomsection\label{unit:recurrence-models}
A linear recurrence is governed by a shift on a finite-dimensional state space. Identifying the shift with multiplication by a coordinate element in a quotient algebra gives a common model for minimal polynomials and periods: the zero-root component produces transients, while multiplicities of nonzero roots affect the period.
In the separable case the coefficient pairing becomes a trace representation; finite observations, in the reverse direction, recover an unknown recurrence. Root-block recovery in Section~\ref{str:sec-trace-recovery} thus leads to two questions: which quotient algebra does a general observation retain, and when does a finite prefix determine the entire recurrence?
We proceed through shifts, quotient algebras and root spaces, trace models, maximum periods, transients, and generating functions, ending with Berlekamp--Massey recovery.
The underlying tools are polynomial orders from Chapter~\ref{chap:orders}, the coefficient pairing from Chapter~\ref{chap:dual}, and cyclic modules from Chapter~\ref{chap:normal}; see also \cite[Chapter 8]{LN}.

\section{The Shift Operator and Linear Recurrence Relations}
\label{sec:linear-recurrences}
\leavevmode\phantomsection\label{def:sequence-shift}\ignorespaces Let \(s=(s_0,s_1,\ldots)\) be a sequence with values in \(k=\F_q\), and define the shift operator \(S(s)_t=s_{t+1}\). For a monic polynomial
\[
f(T)=T^d+c_{d-1}T^{d-1}+\cdots+c_0,\qquad d\ge1,
\]
the equation \(f(S)s=0\) is exactly the recurrence relation
\[
s_{t+d}+c_{d-1}s_{t+d-1}+\cdots+c_0s_t=0\qquad(t\geq0).
\]
A sequence satisfying some recurrence of this form is called a \textbf{linear recurring sequence}. For a given \(f\), a solution is uniquely determined by its first \(d\) terms, so the solution space is a \(d\)-dimensional \(k\)-space.
The tuple \((s_t,\ldots,s_{t+d-1})\) is called the \textbf{state} at time \(t\).
If \(c_0\ne0\), the recurrence also recovers each state uniquely from the next one.
The state transition is therefore a bijection of a finite set, so every state orbit is cyclic from the start.

\begin{definition}\label{def:sequence-minpoly}
The \textbf{minimal polynomial} \(h_s(T)\) of a linear recurring sequence \(s\) is the monic polynomial of least degree satisfying \(h_s(S)s=0\).
The zero sequence is assigned minimal polynomial \(1\). If \(s_{t+N}=s_t\) holds for all \(t\geq0\), then \(N\) is called a period of the sequence;
the smallest positive period is its \textbf{least period}.
\end{definition}
Division with remainder again shows that the polynomials annihilating the sequence \(s\) are exactly the multiples of \(h_s\).
This is the same operator principle as for the Frobenius annihilator polynomial of an element, with Frobenius replaced by the shift.

\Needspace{6\baselineskip}
\begin{theorem}\label{thm:sequence-period}
Suppose that \(f(S)s=0\) and \(f(0)\ne0\). Then \(h_s\mid f\), and the least period of the sequence is
\[
\boxed{\quad \operatorname{per}(s)=\ord(h_s).\quad}
\]
If \(f\) is irreducible, every nonzero solution sequence has minimal polynomial \(f\) and least period \(\ord(f)\).
\end{theorem}
\begin{proof}
The first assertion follows from division with remainder. Moreover, \(N\) is a period if and only if \((S^N-1)s=0\), equivalently, \(h_s\mid T^N-1\).
Thus the least period is the polynomial order defined in Chapter~\ref{chap:orders}.
If \(f\) is irreducible, \(h_s\mid f\) allows only \(h_s=1\) or \(f\); the former corresponds precisely to the zero sequence.
\end{proof}

The orders of irreducible factors determine the part of the period coprime to the characteristic, and their multiplicities determine the additional \(p\)-power part.
For a nonzero recurring sequence satisfying \(h_s(0)\ne0\), if \(h_s=\prod_i f_i^{a_i}\) and \(b\) is the least nonnegative integer satisfying \(p^b\geq\max_i a_i\), then
\[
\operatorname{per}(s)=p^b\lcm_i\ord(f_i).
\]
Thus the multiplicities of irreducible factors affect the sequence period through an additional power of \(p\).

\section{A Quotient-Algebra Model of the Recurrence Space}
\label{sec:recurrence-algebra}

A normal basis represents Frobenius as an action on a cyclic module. For recurring sequences, the corresponding operator is multiplication by the coordinate element in a quotient algebra.
The coefficient pairing provides a common linear model for these two representations while retaining the nilpotent structure of repeated factors.

\par\Needspace{12\baselineskip}
\subsection{The Coefficient Pairing and the Recurrence-Space Isomorphism}
\label{sec:recurrence-coefficient-model}
\begin{theorem}[The coefficient-pairing model of the recurrence space]\label{thm:recurrence-algebra}
Let \(f\in k[T]\) be monic, with \(\deg f=d\ge1\), and put \(A=k[T]/(f)\) and \(\alpha=\overline T\).
Define \(\lambda:A\to k\) to be the coefficient of \(T^{d-1}\) in the canonical remainder. Then the pairing
\[
(u,v)\longmapsto\lambda(uv)
\]
is nondegenerate, and the map
\[
c\longmapsto s(c),\qquad s(c)_t=\lambda(c\alpha^t),\quad t\ge0,
\]
is a linear isomorphism from \(A\) to the solution space of the recurrence \(f(S)s=0\). Under this isomorphism, the shift \(S\) corresponds to multiplication by \(\alpha\). If \(c\) is represented by a remainder of degree less than \(d\), then
\[
h_{s(c)}=\frac{f}{\gcd(f,c)}.
\]
\end{theorem}
\begin{proof}
Nondegeneracy of the coefficient pairing was proved in Theorem~\ref{thm:algebra-trace-coefficient}.
This pairing gives a concrete model of the recurrence solution space.

Since \(f(\alpha)=0\), the sequence \(s(c)\) satisfies the recurrence.
If its first \(d\) terms vanish, then \(\lambda(c\alpha^i)=0\) for \(0\le i<d\).
These \(\alpha^i\) form a basis of \(A\), so nondegeneracy gives \(c=0\). Both the domain and the solution space have dimension \(d\),
hence the map is an isomorphism. The identity \(s(c)_{t+1}=\lambda((\alpha c)\alpha^t)\) shows that the shift corresponds to multiplication.

For \(g\in k[T]\), the condition \(g(S)s(c)=0\) is equivalent to \(\lambda(g(\alpha)c\alpha^t)=0\) for all \(t\ge0\); by nondegeneracy of the pairing, this is equivalent to \(g(\alpha)c=0\), that is, \(f\mid gc\). Write \(b=\gcd(f,c)\), \(f=bf_1\), and \(c=bc_1\),
where \(\gcd(f_1,c_1)=1\). Then \(f\mid gc\) if and only if \(f_1\mid g\). The monic annihilating polynomial of least degree is therefore \(f_1\).
The zero element corresponds to \(b=f\) and \(f_1=1\), and is included in the formula.
\end{proof}

The coefficient pairing also measures how much recurrence information a linear observation retains. For \(\ell:A\to k\), the equality \(\ell(u)=0\) alone does not make \(u\) permanently invisible; all its shifts must be considered.

\begin{proposition}[The observation kernel and recurrence order]\label{v127:prop-observation-kernel}
Retain \(A,\alpha,\lambda\) from Theorem~\ref{thm:recurrence-algebra}. Every \(k\)-linear functional \(\ell\) has a unique expression \(\ell(u)=\lambda(cu)\), with \(c\in A\); use its canonical remainder in the formula below. Put
\begin{equation}\label{v127:eq-observation-kernel}
I_\ell=\{u\in A:\ell(u\alpha^j)=0\text{ for all }j\ge0\},
\qquad h=\frac{f}{\gcd(f,c)},\qquad r=\deg h.
\end{equation}
Then \(I_\ell=\operatorname{ann}_A(c)=(h)/(f)\), and
\[
A/I_\ell\simeq k[T]/(h).
\]
The sequence \(s_j=\ell(\alpha^j)\) has minimal polynomial \(h\). For every \(N\ge r\) with \(N\ge1\),
\begin{equation}\label{v127:eq-observation-rank}
\operatorname{rank}(s_{i+j})_{0\le i,j<N}=r.
\end{equation}
If \(r>0\), the leading \(r\)-by-\(r\) Hankel matrix is invertible, and the first \(2r\) terms recover \(h\). In particular, the \(d\)-by-\(d\) Hankel matrix is invertible if and only if \(c\) is a unit of \(A\).
\end{proposition}
\begin{proof}
Nondegeneracy of the coefficient pairing makes \(c\mapsto\lambda(c\,\cdot)\) an isomorphism \(A\to A^*\).
Since the powers span \(A\), membership in \(I_\ell\) is equivalent to \(\lambda(cuv)=0\) for all \(v\in A\), or \(cu=0\). The equivalence \(f\mid cu\iff h\mid u\) identifies the ideal and quotient algebra.
The product pairing induced by \(\ell\) on this quotient is nondegenerate: its radical is exactly what \(I_\ell\) removes. For \(r>0\), Lemma~\ref{v127:lem-pairing-recovery} gives the minimal polynomial and invertibility of the \(r\)-by-\(r\) Hankel matrix.
For \(N\ge r\), the recurrence expresses every remaining column in terms of the first \(r\), so the rank is exactly \(r\). If \(r=0\), then \(c=0\) and the sequence vanishes, giving the same conclusion.
Finally, \(r=d\) is equivalent to \(\gcd(f,c)=1\), hence to \(c\) being a unit.
\end{proof}

A known annihilating polynomial \(f\) supplies the order bound \(d\). Even when \(r\) is unknown, the first \(2d\) terms suffice to determine \(r\) from the Hankel rank and then recover \(h\); if \(r=0\), then \(h=1\). Section~\ref{sec:berlekamp-massey} performs the same recovery by updating connection relations term by term, using only an order bound rather than a known \(f\).

\begin{corollary}\label{cor:recurrence-count}
For each monic divisor \(h\mid f\), exactly \(\varphi_q(h)\) solution sequences have minimal polynomial \(h\).
If \(f(0)\ne0\), all these sequences have least period \(\ord(h)\).
\end{corollary}
\begin{proof}
The condition \(f/\gcd(f,c)=h\) is equivalent to \(c=(f/h)u\), where \(u\) is taken modulo \(h\) and is coprime to \(h\).
There are exactly \(\varphi_q(h)\) choices. The period assertion follows from Theorem~\ref{thm:sequence-period}.
\end{proof}

This uses the same Euler totient function for polynomial quotient rings as the count of normal elements: both count vectors in a cyclic module with a prescribed minimal annihilating polynomial.
Multiplicative generation of the field and the cyclic-vector property for Frobenius must still be checked separately.

If \(f(0)\ne0\), formula~\eqref{v53:eq-cyclic-gcd-count} also gives
\[
 \#\{s:f(S)s=0,\ S^t s=s\}=q^{\deg\gcd(f,T^t-1)}\qquad(t\ge1).
\]
Substitution into~\eqref{v53:eq-linear-cycle-count} yields all cycle counts for recurrence states.

\par\Needspace{12\baselineskip}
\subsection{Linearized Root Spaces and Compatible Quotient Maps}
\label{sec:recurrence-root-quotient}
\begin{proposition}[The recurrence space and the linearized root space]\label{prop:recurrence-root-correspondence}
Let \(f\in k[T]\) be monic, with \(\deg f=d\ge1\) and \(f(0)\ne0\), and put
\[
\mathcal R_f=\{s:f(S)s=0\},\qquad
W_f=\{x:\mathcal L_f(x)=0\},
\]
where \(W_f\) is taken in a finite splitting extension of \(\mathcal L_f\). After choosing a cyclic vector \(\eta\) in \(W_f\), there are compatible \(k[T]\)-module isomorphisms
\begin{equation}\label{eq:recurrence-root-models}
(\mathcal R_f,S)\simeq\bigl(k[T]/(f),M_{\overline T}\bigr)\simeq(W_f,\sigma).
\end{equation}
Under these isomorphisms, the minimal polynomial of a sequence equals the Frobenius annihilator polynomial of the corresponding element, and the least period of the sequence equals the Frobenius orbit length of that element.
Moreover, if \(N=\ord(f)\), the splitting field of \(\mathcal L_f\) over \(k\) is exactly \(\F_{q^N}\).
\end{proposition}
\begin{proof}
Write \(A=k[T]/(f)\) and \(\alpha=\overline T\), and retain the functional \(\lambda\) of Theorem~\ref{thm:recurrence-algebra}.
That theorem gives the isomorphism \(\Psi:A\to\mathcal R_f\), namely \(\Psi(\bar a)_t=\lambda(a(\alpha)\alpha^t)\), satisfying \(\Psi M_\alpha=S\Psi\).
On the other hand, Theorem~\ref{thm:root-module-count} gives a root \(\eta\) satisfying \(\mu_\eta=f\) and an isomorphism
\[
\Phi_\eta:A\longrightarrow W_f,\qquad
\bar a\longmapsto a(\sigma)\eta.
\]
This isomorphism satisfies \(\Phi_\eta M_\alpha=\sigma\Phi_\eta\). Thus \(\Phi_\eta\Psi^{-1}\) is the required operator-compatible isomorphism.
For \(s=\Psi(\bar a)\) and \(x=\Phi_\eta(\bar a)\),
\[
h_s=\mu_x=\frac{f}{\gcd(f,a)}.
\]
Also, \(S^t s=s\) if and only if \(\sigma^t x=x\), so the two least positive periods coincide; the zero sequence and zero element both have period one by convention.

The splitting-field assertion follows from Corollary~\ref{cor:linearized-splitting-field}.
\end{proof}

Formula~\eqref{eq:recurrence-root-models} describes linear spaces equipped with operator actions;
it depends on the choice of cyclic vector \(\eta\) and generally does not preserve multiplication. The linear subspace \(W_f\) need not be a subfield either.
The condition \(f(0)\ne0\) ensures that the shift is invertible and \(\mathcal L_f\) has no repeated roots. If \(T\mid f\),
sequences in the recurrence space may have a transient phase, an initial segment before entering a periodic cycle, so the entire recurrence space cannot be realized in this way as a Frobenius-stable root space.
The correspondence remains valid even if \(f\) has other repeated factors; these record the nonsemisimple part of the linear operator.

\begin{proposition}[Compatibility of the coefficient pairing with three quotient maps]\label{prop:quotient-model-compatibility}
Let \(f=hv\), where \(h,v\in k[T]\) are monic and \(\deg h\ge1\). Write
\[
A_f=k[T]/(f),\qquad A_h=k[T]/(h),
\qquad \rho:A_f\longrightarrow A_h,\quad [a]_f\longmapsto[a]_h.
\]
Let \(\lambda_f,\lambda_h\) be the functionals extracting the coefficients of \(T^{\deg f-1}\) and \(T^{\deg h-1}\), respectively, from the standard remainders,
and let \(\mathcal C_f:A_f\to\mathcal R_f\) and \(\mathcal C_h:A_h\to\mathcal R_h\) denote the isomorphisms in Theorem~\ref{thm:recurrence-algebra}.
Then, for every \(a\in k[T]\),
\begin{equation}\label{eq:coefficient-quotient-compatibility}
\lambda_f([va]_f)=\lambda_h([a]_h),
\qquad v(S)\mathcal C_f=\mathcal C_h\rho.
\end{equation}
If, in addition, \(f(0)\ne0\), choose a cyclic vector \(\eta_f\) in \(W_f\) and put
\[
\eta_h=\mathcal L_v(\eta_f),\qquad
\Phi_f([a]_f)=a(\sigma)\eta_f,\quad
\Phi_h([a]_h)=a(\sigma)\eta_h.
\]
Then the Frobenius annihilator polynomial of \(\eta_h\) is \(h\), all the maps \(\Phi_f,\Phi_h\) above are module isomorphisms, and
\begin{equation}\label{eq:recurrence-quotient-compatibility}
\mathcal L_v\Phi_f=\Phi_h\rho.
\end{equation}
Thus the natural quotient map \(\rho\), the map \(v(S)\) on recurrence spaces, and the map \(\mathcal L_v\) on root spaces
represent the same quotient map under these isomorphisms.
\end{proposition}
\begin{proof}
Reduce \(a\) modulo \(h\) to \(a_0\), where \(\deg a_0<\deg h\).
Since \(f=hv\), we have \([va]_f=[va_0]_f\) and \(\deg(va_0)<\deg f\).
As \(v\) is monic, the coefficient of \(T^{\deg f-1}\) in \(va_0\) is exactly the coefficient of \(T^{\deg h-1}\) in \(a_0\),
proving the first identity; the zero polynomial also satisfies it.

For every \(j\ge0\), the correspondence between shifting and multiplication by the coordinate element gives
\[
\bigl(v(S)\mathcal C_f([a]_f)\bigr)_j
=\lambda_f([vaT^j]_f)
=\lambda_h([aT^j]_h)
=\bigl(\mathcal C_h([a]_h)\bigr)_j.
\]
Equality term by term proves \(v(S)\mathcal C_f=\mathcal C_h\rho\). This part does not require \(f(0)\ne0\) and therefore also applies to recurrences with a transient phase.

Now assume that \(f(0)\ne0\). By \(\mu_{\eta_f}=f\) and the formula for minimal annihilating polynomials in a cyclic module,
\[
\mu_{\eta_h}=\frac{f}{\gcd(f,v)}=h.
\]
Thus \(\eta_h\) is cyclic in \(W_h\), and \(\Phi_h\) is a module isomorphism.
Operator polynomials with base-field coefficients commute, so
\[
\mathcal L_v\bigl(\Phi_f([a]_f)\bigr)
=v(\sigma)a(\sigma)\eta_f
=a(\sigma)\eta_h
=\Phi_h\bigl(\rho([a]_f)\bigr).
\]
This proves the last compatibility relation.
\end{proof}

If \(h=1\), take \(A_h\), \(\mathcal R_h\), and \(W_h\) to be zero spaces;
all the compatible quotient maps above are zero, and \(\lambda_h\) need not be defined.

\Needspace{8\baselineskip}
\begin{proposition}[Linearized fibers and affine returns for recurrence quotients]\label{prop:recurrence-fiber-return}
Let \(h,v\in k[T]\) be monic, and let \(f=hv\) satisfy \(\deg f\ge1\) and \(f(0)\ne0\). Place \(W_f,W_h,W_v\) in one common finite splitting extension.
Then the map
\[
\pi=\mathcal L_v:W_f\longrightarrow W_h
\]
is surjective, commutes with Frobenius, and satisfies \(\ker\pi=W_v\).
Take \(x\in W_f\) and \(y=\pi(x)\), and let \(y\) have Frobenius period \(t\).
The fiber is \(x+W_v\), and its first-return map in the coordinate \(z\in W_v\) is
\[
z\longmapsto \sigma^t(z)+\delta,\qquad \delta=\sigma^t(x)-x\in W_v.
\]
If zero has period \(e\) under this affine map, then \(x\) has Frobenius period \(te\).
Under the isomorphisms in~\eqref{eq:recurrence-root-models}, this is also the shift period of the corresponding recurrence state.
\end{proposition}
\begin{proof}
Composition of linearized polynomials with base-field coefficients satisfies \(\mathcal L_h\circ\mathcal L_v=\mathcal L_f\),
so \(\pi(W_f)\subseteq W_h\). Since \(\mathcal L_f=\mathcal L_h\circ\mathcal L_v\),
we have \(W_v\subseteq W_f\), hence \(\ker(\pi|_{W_f})=W_v\). The dimension of each root space equals the degree of the corresponding operator polynomial, so
\[
\dim_k\pi(W_f)=\deg f-\deg v=\deg h=\dim_kW_h.
\]
This proves surjectivity; coefficients in \(k\) also give \(\pi\sigma=\sigma\pi\).
Proposition~\ref{prop:quotient-model-compatibility} further identifies \(\pi\) with the natural quotient map \(\rho\) in cyclic-module coordinates
and with \(v(S)\) in recurrence coordinates. For \(h=1\), the corresponding target is the zero space \(\{0\}\), and the conclusion still holds.

Since \(\pi\) is linear, \(\pi^{-1}(y)=x+W_v\). Equivariance and \(\sigma^t(y)=y\) give
\[
\pi(\sigma^t(x)-x)=\sigma^t(y)-y=0,
\]
so \(\delta\in W_v\). Applying \(\sigma^t\) to \(x+z\) and subtracting the base point \(x\) gives
\[
\sigma^t(x+z)-x=\sigma^t(z)+\delta.
\]
The base point \(x\) corresponds to zero in the fiber coordinate, so its internal return period is exactly the period \(e\) of zero under the stated affine map.
Restricting \(\pi\) over the entire orbit of \(y\) and applying Lemma~\ref{lem:fiber-return} gives the total period \(te\).
Finally, compatibility of the shift and Frobenius in Proposition~\ref{prop:recurrence-root-correspondence} gives the interpretation as a recurrence-state period.
\end{proof}

This proposition is also a concrete form of Theorem~\ref{ffs:thm-additive-return} in the linearized-root-space model of a recurrence quotient.
Retaining its notation, take \(V=W_v\). Since \(v(0)\ne0\), the polynomial \(\mathcal L_v\) is monic and has no repeated roots, with root set exactly \(W_v\). Thus
\[
L_V=\mathcal L_v.
\]
Let \(P\in k[Y]\) be the minimal polynomial of the element \(y\) over \(k\). Then \(\deg P=t\) and \(P\mid\mathcal L_h\).
Here \(P\) is determined by the conjugates of \(y\), whereas \(h\) specifies the operator relation \(h(\sigma)y=0\).
Since \(\mathcal L_h\circ\mathcal L_v=\mathcal L_f\), in the chosen splitting extension we have
\[
\{u:P(\mathcal L_v(u))=0\}
=\pi^{-1}\bigl(\{y,y^q,\ldots,y^{q^{t-1}}\}\bigr)\subseteq W_f.
\]
Thus the root set of \(P(\mathcal L_v)\) is exactly the preimage under \(\pi\) of the entire Frobenius orbit of \(y\).
In Theorem~\ref{ffs:thm-additive-return}, take \(\beta=y\), \(\alpha=x\), and \(d=t\). The resulting fiber return is
\[
z\longmapsto z^{q^t}+(x^{q^t}-x),\qquad z\in W_v,
\]
which is exactly the internal return of the recurrence quotient above.
Proposition~\ref{prop:quotient-model-compatibility} also identifies \(\pi\) with \(v(S)\) on the recurrence space,
so the return period of the root fiber determines the shift period of the recurrence state.
The quotient state has period \(t\), and the original state has return period \(e\) within the corresponding fiber, so the original state has period \(te\).
The additive recursions of Chapter~\ref{chap:root-constructions} use a nonzero trace to make this internal period maximal;
for a general recurrence, the possible internal periods are determined by the operator polynomial and its factors.

\begin{example}[A nilpotent direction and a sequence of period two]\label{exm:repeated-recurrence-algebra}
Over \(k=\F_2\), take \(f=(T+1)^2\). The quotient algebra can be written as \(A=\F_2[\varepsilon]/(\varepsilon^2)\), where \(\alpha=1+\varepsilon\).
This is the algebra of Example~\ref{exm:dual-numbers-trace}, but now the shift corresponds to multiplication by the unit \(1+\varepsilon\).
Since \(a+b\varepsilon=(a+b)+b\alpha\), the linear coefficient of the canonical remainder is still \(\lambda(a+b\varepsilon)=b\). From \(\alpha^t=1+t\varepsilon\), the parameter \(c=1\) gives the sequence \(s_t=t\pmod2\), namely \(0,1,0,1,\ldots\). Its minimal polynomial is \((T+1)^2\), and its least period is two.
For the first four terms,
\[
H=\begin{pmatrix}0&1\\1&0\end{pmatrix},\qquad
H^+=\begin{pmatrix}1&0\\0&1\end{pmatrix},\qquad
H^{-1}H^+=\begin{pmatrix}0&1\\1&0\end{pmatrix}.
\]
Thus finite observations recover the multiplication operator in the basis \(1,\alpha\); its minimal polynomial remains \((T+1)^2\).

Example~\ref{exm:dual-numbers-trace} showed that \(\Tr_{A/k}\) is identically zero.
If \(\lambda\) in \(s(c)_t=\lambda(c\alpha^t)\) is replaced by the algebra trace,
every parameter \(c\) produces only the zero sequence, so this two-dimensional solution space cannot be parametrized.
If one first passes to the reduced quotient \(A/(\varepsilon)\simeq\F_2\), the image of \(\alpha\) is one, and the multiplicative transition becomes the identity,
representing only constant sequences and again losing period two. Both losses of information show why the nilpotent structure must be retained.

The corresponding linearized polynomial, however, is \(\mathcal L_f(X)=X^4+X\), with root space \(\F_4\).
Choose \(\eta\in\F_4\setminus\F_2\). Then \(\eta^2+\eta=1\) and \(\mu_\eta=(T+1)^2\).
Under the isomorphism of Proposition~\ref{prop:recurrence-root-correspondence}, \(0,1,0,1,\ldots\) corresponds to \(\eta\),
and its shift \(1,0,1,0,\ldots\) corresponds to \(\eta^2\). The least period two is also realized by Frobenius interchanging these two roots.
\end{example}

\begin{example}[Two root sets of size 64 and two types of cycles]\label{exm:two-sixtyfour-rootsets}
Retain Example~\ref{exm:f64-computation}, put \(k=\F_2\) and \(f(T)=T^6+T+1\), and take \(E=\F_{64}=k(\alpha)\) and \(f(\alpha)=0\). That example proved that \(f\) is primitive, so \(\ord(f)=\ord(\alpha)=63\). In a common finite extension, compare
\[
E=\{x:x^{64}-x=0\},\qquad
W=\{x:x^{64}+x^2+x=0\}.
\]
The second polynomial is \(\mathcal L_f\), with linear coefficient one, so all its roots form a six-dimensional \(k\)-space.
Theorem~\ref{thm:root-module-count} shows that every nonzero root has Frobenius annihilator polynomial \(f\),
and therefore degree \(\ord(f)=63\) over \(k\).
These 63 nonzero roots consequently form a single Frobenius orbit; Corollary~\ref{cor:linearized-splitting-field}
gives the splitting field \(\F_{2^{63}}\).

\begin{center}\small
\setlength{\tabcolsep}{5pt}
\begin{tabular}{@{}lcc@{}}
\toprule
Structure compared&\(E\)&\(W\)\\
\midrule
Cardinality and \(k\)-dimension&\(64,\ 6\)&\(64,\ 6\)\\
A subfield?&Yes&No\\
Order of the squaring operator \(\sigma\)&\(6\)&\(63\)\\
Minimal polynomial of squaring&\(T^6-1\)&\(T^6+T+1\)\\
Field generated by the root set&\(\F_{2^6}\)&\(\F_{2^{63}}\)\\
\bottomrule
\end{tabular}
\end{center}
The set \(W\) in the table is not a subfield, since \(1^{64}+1^2+1=1\ne0\) and hence \(1\notin W\).
Equal cardinalities and linear dimensions do not determine either multiplicative closure of a root set or its Frobenius periods.

Take \(0\ne\eta\in W\). The multiplication operator \(M_\alpha:E\to E\) has minimal polynomial \(f\),
and \(\sigma|_W\) has the same minimal polynomial \(f\). An explicit isomorphism is
\[
\Phi_\eta:E\longrightarrow W,\qquad
a(\alpha)\longmapsto a(\sigma)\eta,\qquad
\Phi_\eta M_\alpha=\sigma\Phi_\eta.
\]
It comes from two realizations of the same cyclic module \(k[T]/(f)\) and preserves the linear structure and the specified operators;
since \(\Phi_\eta(1)=\eta\) and \(W\) is not a field, it cannot be interpreted as a field isomorphism.
Under this correspondence, \(\alpha^j\) maps to \(\eta^{2^j}\): the former runs through the multiplicative orbit in \(E^\times\),
and the latter through the Frobenius orbit in \(W\setminus\{0\}\).

Moreover, since \(f'(\alpha)=1\), the coefficient functional on \(E\) is the trace, so Proposition
\ref{prop:recurrence-root-correspondence} associates the sequence
\[
s(c)_j=\Tr_{E/k}(c\alpha^j)
\]
with \(\Phi_\eta(c)\). The shift period 63 of a nonzero sequence becomes Frobenius period 63 in \(W\),
whereas \(\alpha\) itself still has Frobenius period six in \(E\).
\end{example}

When \(f\) is irreducible, \(A=k(\alpha)\) is a field. By Proposition~\ref{prop:coefficient-trace}, the preceding isomorphism becomes
\[
s(c)_t=\lambda(c\alpha^t)
=\Tr_{A/k}\left(\frac{c}{f'(\alpha)}\alpha^t\right).
\]
Since multiplication by \(1/f'(\alpha)\) is an invertible linear map, the quotient-algebra parameter and the trace parameter differ only by this fixed factor.
If \(f\) is square-free but reducible, the same formula holds, with the trace interpreted as the sum of the traces on the field components.
If \(f=\prod_iP_i^{e_i}\) has repeated factors, Corollary~\ref{prop:local-multiplicity-trace} gives
\[
\Tr_{A/k}(c\alpha^j)=
\sum_i e_i\Tr_{K_i/k}(\bar c_i\bar\alpha_i^{\,j}).
\]
This depends only on the images of the parameter in the fields \(K_i\); when \(p\mid e_i\), the corresponding component disappears altogether.
Thus replacing \(\lambda\) by the algebra trace loses the nilpotent directions and their shift data.
Example~\ref{exm:repeated-recurrence-algebra} shows that the same direction can still produce a nonconstant periodic sequence under multiplicative transition.
Nondegeneracy of the coefficient pairing retains all states in the general recurrence model, while the derivative formula describes exactly which states are lost when the trace is used instead.

\section{Trace Representations and Conjugate Decomposition of Recurrences}\label{v100:sec-recurrence-trace}
In the irreducible or squarefree case, the coefficient functional of the quotient algebra has a trace representation. We express recurrence parameters as elements of the root fields and retain each separable component through conjugate decomposition.
\begin{theorem}\label{thm:sequence-trace}
Let \(f\in k[T]\) be an irreducible polynomial of degree \(d\), with \(f(0)\ne0\), and let \(\alpha\) be a root in \(E=\F_{q^d}\).
Then every solution of the recurrence \(f(S)s=0\) has a unique representation of the form
\[
\boxed{\quad s_t=\Tr_{E/k}(c\alpha^t),\qquad c\in E.\quad}
\]
\end{theorem}
\begin{proof}
Making \(f\) monic changes neither its roots nor its recurrence solution space, so we may assume that \(f\) is monic.
Theorem~\ref{thm:recurrence-algebra} gives a unique quotient-algebra parameter \(b\in E\) expressing each solution as \(s_t=\lambda(b\alpha^t)\). By~\eqref{eq:coefficient-trace}, this equals \(\Tr_{E/k}((b/f'(\alpha))\alpha^t)\). Put \(c=b/f'(\alpha)\);
since \(f'(\alpha)\ne0\), this parameter change is bijective, proving existence and uniqueness of the trace representation.
Conversely, for any \(c\in E\), taking \(b=cf'(\alpha)\) gives a solution by the same theorem.
\end{proof}

This connects two kinds of orbits: repeated multiplication by \(\alpha\) in the extension field, followed by the trace to the base field, produces a recurring sequence.
The weighted relative trace sequence of Section~\ref{str:sec-trace-recovery} is an instance with the base field replaced by the subfield \(K_d\);
its Hankel matrix uses the nondegeneracy of the trace pairing to give a nonsingular linear system that recovers the relative root block from a finite prefix.
Replacing the root by a conjugate \(\alpha^{q^i}\) does not change the set of represented sequences, since the trace is Frobenius-invariant.
To recover \(c\) from the first \(d\) terms, use the dual basis \((\beta_0,\ldots,\beta_{d-1})\) of the power basis:
\[
c=\sum_{i=0}^{d-1}s_i\beta_i.
\]

\begin{corollary}[Conjugate decomposition of separable recurrences]\label{cor:squarefree-recurrence-trace}
Suppose that \(f=f_1\cdots f_s\) is a product of distinct monic irreducible polynomials and \(f(0)\ne0\).
Choose a root \(\alpha_i\in E_i=\F_{q^{\deg f_i}}\) of \(f_i\). Then each solution sequence has a unique expression
\[
s_t=\sum_{i=1}^s\Tr_{E_i/k}(c_i\alpha_i^t),\qquad c_i\in E_i.
\]
Its minimal polynomial is the product of all \(f_i\) for which \(c_i\ne0\).
\end{corollary}
\begin{proof}
By the Chinese remainder theorem, \(A=k[T]/(f)\simeq\prod_iE_i\), with \(\overline T\) corresponding to \((\alpha_i)_i\).
The nondegenerate pairing of Theorem~\ref{thm:recurrence-algebra} shows that every solution sequence also has a unique expression \(s_t=\ell(\overline T^t)\), where \(\ell\) is a linear functional on \(A\).
Every functional on the direct product decomposes uniquely as a sum of functionals on its components. The trace pairing represents the functional on the \(i\)th component uniquely as \(x\mapsto\Tr_{E_i/k}(c_ix)\). This proves the representation and its uniqueness.
A polynomial \(g\) annihilates the sequence if and only if \(c_i g(\alpha_i)=0\) for every \(i\),
so its monic minimal annihilating polynomial is exactly the stated product.
\end{proof}

\section{Primitive Polynomials and Maximum-Length Sequences}
\label{v100:sec-maximum-period}
Suppose that a nonzero purely periodic sequence has a minimal polynomial of degree \(d\). No \(d\) consecutive terms can all vanish,
since the recurrence would make all subsequent terms zero, and pure periodicity would then make the whole sequence zero.
States cannot repeat within one least period either: equal states determine the same subsequent terms, and pure periodicity would give a shorter period.
There are only \(q^d-1\) nonzero states, so the least period is at most \(q^d-1\). The following definition describes sequences attaining this bound.

\begin{definition}\leavevmode\phantomsection\label{def:m-sequence}
A linear recurring sequence whose minimal polynomial has degree \(d\ge1\) and whose least period is \(q^d-1\) is called a \textbf{maximum-length sequence}, or an \(m\)-sequence.
\end{definition}

\begin{theorem}\label{thm:m-sequence}
If \(f\) is a primitive polynomial of degree \(d\), every nonzero solution sequence is an \(m\)-sequence;
the states formed by \(d\) consecutive terms run through \(k^d\setminus\{0\}\) exactly once in one period.
In one period, the symbol \(0\) occurs \(q^{d-1}-1\) times, and each nonzero symbol occurs \(q^{d-1}\) times.
Conversely, if a sequence with minimal polynomial of degree \(d\) has least period \(q^d-1\), its minimal polynomial is primitive.
\end{theorem}
\begin{proof}
Let \(\alpha\) be a root of \(f\). For nonzero \(c\), the elements \(c\alpha^t\) run through \(E^\times\).
The linear isomorphism from the trace representation associates them with all nonzero states, whose first coordinate is \(\Tr(c\alpha^t)\).
Each trace fiber has \(q^{d-1}\) elements; removing the zero element of the extension field gives the symbol distribution.

Conversely, let the sequence have minimal polynomial \(h\), with \(\deg h=d\). Periodicity gives
\(h\mid T^{q^d-1}-1\), hence \(h(0)\ne0\). By Theorem~\ref{thm:sequence-period},
\(\ord(h)=q^d-1\). Proposition~\ref{prop:maximal-order-primitive} therefore shows that \(h\) is primitive.
This reduces the converse for maximum-length sequences directly to the criterion for multiplicative order in a quotient algebra to attain its upper bound.
\end{proof}

\begin{example}\leavevmode\phantomsection\label{exm:m-sequence-seven}
Take \(f(T)=T^3+T+1\) from Example~\ref{ex:normal-f8-example}. Its roots are nonidentity elements of \(\F_8^\times\), hence have order seven, so \(f\) is primitive.
Here the same root determines recurrence states and a symbol distribution.
The recurrence is \(s_{t+3}=s_{t+1}+s_t\). The initial state \((1,0,0)\) gives
\[
1,0,0,1,0,1,1,\quad1,0,0,1,0,1,1,\quad\ldots
\]
One period contains four ones and three zeros, and its seven consecutive triples are exactly the nonzero vectors of \(\F_2^3\).
\end{example}

\section{Transients, Eventual Periods, and Rational Generating Functions}
\label{sec:transient-period}
\subsection{The Exact Decomposition into Transient and Periodic Parts}
\label{sec:transient-decomposition}
\leavevmode\phantomsection\label{def:transient-period}\ignorespaces For a sequence \(s=(s_j)_{j\ge0}\), suppose that there exist \(u\ge0\) and \(N\ge1\) such that
\(s_{j+N}=s_j\) for all \(j\ge u\). Then \(s\) is called eventually periodic.
The smallest such \(u\) is its \textbf{minimal preperiod};
the least positive integer \(N\) for which the equality holds on some tail is its \textbf{least eventual period}.
A sequence with preperiod zero is purely periodic. Below we adopt the convention \(\ord(1)=1\).

\begin{theorem}[Preperiods and eventual periods of general recurrences]\label{thm:transient-period}
Let \(s\) be a linear recurring sequence over \(k=\F_q\) with minimal polynomial
\[
h_s(T)=T^a g(T),\qquad a\ge0,\quad g(0)\ne0,
\]
where \(g\) is monic. Then \(s\) is eventually periodic, with
\[
\boxed{\quad\text{minimal preperiod}=a,\qquad
\text{least eventual period}=\ord(g).\quad}
\]
More precisely, for every \(u\ge0\) and \(N\ge1\),
\begin{equation}\label{eq:transient-divisibility}
s_{j+N}=s_j\ (j\ge u)
\quad\Longleftrightarrow\quad u\ge a\ \text{and}\ g\mid T^N-1.
\end{equation}
\end{theorem}
\begin{proof}
The tail-period relation is equivalent to \(S^u(S^N-1)s=0\) for the whole sequence. By the definition of the minimal polynomial and division with remainder,
this is equivalent to \(h_s\mid T^u(T^N-1)\).
Since \(T\nmid T^N-1\) and \(g(0)\ne0\), the last divisibility splits exactly into \(u\ge a\) and \(g\mid T^N-1\).
If \(g\ne1\), the quotient ring \(k[T]/(g)\) is finite and \(\overline T\) is a unit, so some positive integer \(N\) satisfies \(g\mid T^N-1\);
for \(g=1\), take \(N=1\). Thus the least possible \(u\) is \(a\) and the least possible \(N\) is \(\ord(g)\), and both are attained simultaneously.
The zero sequence corresponds to \(h_s=g=1\) and \(a=0\) and is included in the conclusion.
\end{proof}

A B\'ezout identity gives explicit projections for this decomposition. Choose \(U,V\in k[T]\) such that
\[
U(T)T^a+V(T)g(T)=1,
\]
and define
\begin{equation}\label{eq:transient-projectors}
\begin{aligned}
s_{\mathrm{nil}}&=V(S)g(S)s,\\
s_{\mathrm{per}}&=U(S)S^a s.
\end{aligned}
\qquad s=s_{\mathrm{nil}}+s_{\mathrm{per}}.
\end{equation}
In the cyclic module generated by \(s\), the isomorphism
\[
k[T]/(T^ag)\longrightarrow k[T]s,\qquad [A]\longmapsto A(S)s
\]
sends \(1\) to \(s\). The Chinese remainder decomposition
\(k[T]/(T^ag)\simeq k[T]/(T^a)\times k[T]/(g)\)
then sends \(Vg\) and \(UT^a\) to \((1,0)\) and \((0,1)\), respectively. Hence
\[
h_{s_{\mathrm{nil}}}=T^a,\qquad h_{s_{\mathrm{per}}}=g.
\]
When a modulus polynomial is one, its component is zero, interpreted as having minimal polynomial one.
Since \(S^a s_{\mathrm{nil}}=0\), the first component vanishes from index \(a\) onward;
if \(a>0\), then \(S^{a-1}s_{\mathrm{nil}}\ne0\), so its \(a-1\)th term is nonzero.
The second component has least period \(\ord(g)\) from index zero.
This gives the unique decomposition of \(s\) as the sum of a finitely supported sequence, meaning one with finitely many nonzero terms, and a purely periodic sequence:
the difference of two such decompositions would be both finitely supported and purely periodic and must therefore vanish.

Thus only the zero-root factor causes a transient phase. Repeated factors at nonzero roots remain in the invertible part
and may enlarge the \(p\)-power part of the period as in Theorem~\ref{thm:order-multiplicity}.
Proposition~\ref{prop:recurrence-root-correspondence} realizes the invertible component as Frobenius on a linearized root space,
while~\eqref{eq:transient-projectors} explicitly gives the map from the original sequence to that component.
A nonzero \(T\)-primary component on which the shift is nilpotent cannot be realized as Frobenius on a Frobenius-stable subspace of an extension field,
because the latter is always injective.

This restriction depends on the specified linear action. Again take, as in Examples~\ref{exm:dual-numbers-trace} and~\ref{exm:repeated-recurrence-algebra},
\(A=\F_2[\varepsilon]/(\varepsilon^2)\), and let \(M_a\) denote the linear operator of multiplication by \(a\). Then
\[
M_\varepsilon^2=0,
\qquad M_{1+\varepsilon}^2=I,\quad M_{1+\varepsilon}\ne I.
\]
The nonzero orbits of the former eventually reach zero, whereas all orbits of the latter are periodic from the start and the operator has order two.
Take \(\eta\in\F_4\) satisfying \(\eta^2+\eta+1=0\). The \(\F_2\)-linear isomorphism
\[
\Psi:A\longrightarrow\F_4,\qquad a+b\varepsilon\longmapsto a\eta+b
\]
satisfies \(\Psi\circ M_{1+\varepsilon}=\sigma\circ\Psi\), where \(\sigma(x)=x^2\). Indeed,
\[
\Psi\bigl((1+\varepsilon)(a+b\varepsilon)\bigr)
=a\eta+(a+b)=(a\eta+b)^2.
\]
This is an isomorphism of linear spaces with their operators and does not preserve multiplication.
Thus nonzero nilpotents in a quotient algebra do not prevent an invertible multiplication operator on it from having a Frobenius model;
the transient phase comes from the zero-root primary component of the shift, while repeated factors at nonzero roots may remain in a periodic action.

\begin{example}[Squaring and multiplication on the same quotient algebra]\leavevmode\phantomsection\label{exm:mixed-algebra-recurrence}
Retain from Example~\ref{exm:mixed-algebra-trace} the algebra
\[
f=(T+1)^2(T^2+T+1)^3,\qquad A=\F_2[T]/(f),\qquad \alpha=\overline T.
\]
First consider its own Frobenius map \(F(u)=u^2\). That example showed that \(F^2\) is a projection onto a three-dimensional periodic subalgebra.
On this subalgebra, \(F\) corresponds to componentwise squaring in \(\F_2\times\F_4\), so all eventual periods divide two,
and all preperiods are at most two. The bound is attained: in Chinese remainder coordinates,
\(u=(0,\overline P)\) and \(P=T^2+T+1\) satisfy \(F(u)\ne0\) and \(F^2(u)=0\).
Since the periodic subalgebra contains an \(\F_4\) component, elements of least period two also occur.

Now consider the multiplication operator \(M_\alpha\). The condition \(f(0)=1\) ensures that \(\alpha\) is invertible.
The orders of \(T+1\) and \(P\) are one and three, respectively, and the largest multiplicity is three; the smallest power of two at least this multiplicity is four.
Theorem~\ref{thm:order-multiplicity} therefore gives
\[
\ord(M_\alpha)=\ord(f)=4\lcm(1,3)=12.
\]
The first equality can also be checked directly: \(M_\alpha^j=I\) if and only if \(\alpha^j=1\);
necessity follows by evaluation at \(1\in A\).
Thus all orbits of this multiplication action are purely periodic, with periods dividing 12.
Let \(\lambda\) extract the coefficient of \(T^7\). By Theorem~\ref{thm:recurrence-algebra}, the sequence \(s_j=\lambda(\alpha^j)\)
has minimal polynomial \(f\) and hence least period exactly 12.

\begin{center}\small
\renewcommand{\arraystretch}{1.2}
\begin{tabular}{@{}p{0.35\textwidth}p{0.58\textwidth}@{}}
\toprule
Object in the same algebra&Corresponding data\\
\midrule
Frobenius \(F(u)=u^2\)&\(F^2\) projects onto a three-dimensional periodic part;\newline preperiod at most two, eventual period dividing two\\
Algebra trace pairing&Rank two, radical dimension six; nilradical dimension five\\
Coefficient pairing&Nondegenerate, rank eight\\
Multiplication operator \(M_\alpha\)&Order 12; all orbits purely periodic, with an orbit of period 12\\
\bottomrule
\end{tabular}
\end{center}

Proposition~\ref{prop:recurrence-root-correspondence} also realizes \(M_\alpha\) as Frobenius on the eight-dimensional linearized root space \(W_f\),
whose splitting field is \(\F_{2^{12}}\). The operator realized is the specified \(M_\alpha\),
whereas the squaring map \(F\) of \(A\) itself has nonzero kernel.
Thus trace degeneracy, nilpotents, and recurrence transients must be determined from the pairing, algebra multiplication, and specified transition operator, respectively.
\end{example}

\subsection{Rational Generating Functions and Minimal Polynomials}
\phantomsection\label{def:sequence-generating-function}
The \textbf{generating function} of a sequence is the formal power series \(G_s(Z)\in k[[Z]]\) having the \(t\)th term as its coefficient of \(Z^t\).
Operations on formal power series are described in Section~\ref{sec:counting-euler}; here the coefficients lie in the base field \(k\) of the sequence.
If \(s\) satisfies \(f(S)s=0\), write
\[
G_s(Z)=\sum_{t\geq0}s_tZ^t,
\qquad C_f(Z)=Z^d f(Z^{-1})=1+c_{d-1}Z+\cdots+c_0Z^d.
\]
Direct coefficient comparison makes the recurrence equivalent to
\[
C_f(Z)G_s(Z)=P(Z),\qquad \deg P<d.
\]
If \(h_s(0)\ne0\), reduce this rational expression and normalize the constant term of its denominator to one. The resulting denominator is exactly \(C_{h_s}(Z)=Z^{\deg h_s}h_s(Z^{-1})\).
Indeed, the numerator obtained from the minimal recurrence has degree below the denominator; any cancellable common factor would give a recurrence of smaller order, a contradiction.
Conversely, the recurrence from the reduced denominator annihilates the entire sequence and must therefore correspond to the minimal polynomial.

For the general case \(h_s=T^ag\), formula~\eqref{eq:transient-projectors} further gives
\[
G_s(Z)=N(Z)+\frac{P(Z)}{C_g(Z)},\qquad
C_g(Z)=Z^{\deg g}g(Z^{-1}),
\]
where \(N\) is the generating polynomial of \(s_{\mathrm{nil}}\). If \(a>0\), then \(\deg N=a-1\); if \(a=0\), then \(N=0\).
When \(g\ne1\), the second term is a reduced proper fraction with \(\deg P<\deg g\); when \(g=1\), the second term is zero.
Thus the polynomial part records the transient exactly, and the reduced denominator records the minimal polynomial of the invertible component exactly.
Combining the two parts gives the entire minimal polynomial directly from the reduced numerator and denominator.

\begin{proposition}[Recovering the minimal polynomial from a reduced generating function]\label{v70:prop-rational-generating-annihilator}
Suppose that a nonzero sequence \(s\) has the reduced rational generating function
\[
G_s(Z)=\frac{A(Z)}{B(Z)},\qquad A,B\in k[Z],\qquad \gcd(A,B)=1,\quad B(0)=1.
\]
Put
\[
b=\deg B,\qquad g(T)=T^bB(T^{-1}),\qquad
 a=\max\{0,\deg A-b+1\}.
\]
Then the minimal polynomial of \(s\) is
\[
h_s(T)=T^ag(T).
\]
Its least transient length is \(a\), and its least eventual period is \(\ord(g)\), with \(\ord(1)=1\).
The formula includes the case \(B=1\); the zero sequence has \(h_s=1\).
\end{proposition}
\begin{proof}
Division with remainder gives \(A=NB+P\), where \(\deg P<b\), with \(P=0\) when \(B=1\).
If \(\deg A<b\), then \(N=0\); otherwise, \(\deg N=\deg A-b=a-1\).
Let \(u,v\) be the coefficient sequences of \(N\) and \(P/B\), respectively, so that \(s=u+v\).
If \(u\ne0\), then \(S^au=0\) and \(S^{a-1}u\ne0\), giving minimal polynomial \(T^a\);
if \(u=0\), this polynomial is \(1=T^0\).

If \(b>0\), reducedness gives \(\gcd(P,B)=1\), and \(P/B\) is a proper fraction.
The preceding coefficient comparison and cancellation argument show that the minimal polynomial of \(v\) is exactly \(g\), with \(g(0)\ne0\).
If \(B=1\), then \(v=0\) and \(g=1\).
Since \(T^a\) and \(g\) are coprime, choose \(c,d\in k[T]\) with \(cT^a+dg=1\). Then
\[
d(S)g(S)s=u,\qquad c(S)S^as=v.
\]
Every polynomial annihilating \(s\) therefore annihilates both \(u\) and \(v\), and is divisible by both \(T^a\) and \(g\).
Conversely, \(T^ag\) annihilates their sum, so it is exactly \(h_s\).
The transient and eventual-period assertions follow from Theorem~\ref{thm:transient-period}.
\end{proof}

For example, the binary sequence in Exercise~\ref{ex:wu-transient} has reduced generating function
\[
G_s(Z)=\frac{1+Z+Z^2}{1+Z}.
\]
Here \(b=1\), \(a=2\), and \(g=T+1\), giving \(h_s=T^2(T+1)\).
The denominator records the periodic component, while the degree of the numerator relative to the denominator supplies the zero-root factor.
This recovery assumes that the complete rational generating function is known; a finite initial segment alone does not determine the minimal annihilating polynomial of the entire sequence.
The next section allows a factor \(T\) in the minimal polynomial when studying finite prefixes, and records both a connection polynomial and a connection length.

\section{Recovering Recurrences from Finite Observations: The Berlekamp--Massey Algorithm}
\label{sec:berlekamp-massey}

A minimal polynomial describes the shift relations of the entire sequence; computations usually begin with finitely many terms.
A finite prefix determines connection relations over a finite range; recovering a recurrence for the entire sequence also requires an order bound.
Proposition~\ref{v127:prop-observation-kernel} already recovers the recurrence from a Hankel rank and a linear system.
When observations arrive one at a time, one can retain a shortest connection relation and correct it only when a new term violates it.
The discrepancy bound below determines when this correction must increase the length; it also applies to sequences with transients.

\begin{definition}\leavevmode\phantomsection\label{def:finite-linear-complexity}
Let \(N\ge0\). For a finite prefix \(s_0,\ldots,s_{N-1}\in k\), if \(0\le L\le N\) and \(C(Z)=1+c_1Z+\cdots+c_LZ^L\) satisfy
\[
s_j+c_1s_{j-1}+\cdots+c_Ls_{j-L}=0\qquad(L\le j<N),
\]
then \((L,C)\) is called a \textbf{connection relation} for this prefix, where \(L\) is the \textbf{connection length} and \(C\) the \textbf{connection polynomial}.
The smallest possible \(L\) is the \textbf{linear complexity} of the prefix.
Since \(c_L=0\) is allowed, the length \(L\) must be retained; recording only \(\deg C\) does not suffice.
For \(L=N\) there are no equations to check, so the minimum always exists.
\end{definition}

\leavevmode\phantomsection\label{def:recurrence-discrepancy}\ignorespaces For a connection relation \((L,C)\), the value computed at a new index \(j\ge L\),
\(\delta=s_j+\sum_{i=1}^L c_i s_{j-i}\), is called the \textbf{discrepancy} there. It vanishes precisely when the previous connection relation remains valid at that position.

\begin{lemma}[A lower bound from one discrepancy]\leavevmode\phantomsection\label{lem:bm-discrepancy-bound}
Suppose that \(0\le L\le j\) and that \((L,C)\) satisfies the connection relations for the first \(j\) terms but has a nonzero discrepancy at index \(j\),
\[
\delta=s_j+\sum_{i=1}^L c_i s_{j-i}\ne0.
\]
If \((L',D)\) satisfies the first \(j+1\) terms, then \(L+L'\ge j+1\).
\end{lemma}
\begin{proof}
Suppose to the contrary that \(L+L'\le j\), and write \(c_0=d_0=1\). Consider the same double sum in two orders:
\[
\sum_{u=0}^{L'}d_u\sum_{i=0}^{L}c_i s_{j-u-i}
=\sum_{i=0}^{L}c_i\sum_{u=0}^{L'}d_u s_{j-i-u}.
\]
On the left, the term for \(u=0\) is \(\delta\); in every other inner sum the largest sequence index \(j-u\) lies in the range \(L,\ldots,j-1\), so that sum vanishes.
Each inner sum on the right is a valid connection relation for \(D\), since \(j-i\ge L'\). The right-hand side is therefore zero, a contradiction.
\end{proof}

The algorithm below retains the current relation and the relation just before the most recent increase in complexity.
The quantities \(B,b,t\) record the latter relation, its nonzero discrepancy, and the discrepancy position, respectively.
The Berlekamp--Massey algorithm solves the shortest linear feedback shift register synthesis problem: it finds a shortest linear recurrence generating a given prefix; see \cite{Massey}.

\Needspace{18\baselineskip}
\begin{quote}\small
\leavevmode\phantomsection\label{alg:berlekamp-massey}\ignorespaces Input: \(s_0,\ldots,s_{N-1}\in k\).

Initialize \(C=B=1\), \(L=0\), \(b=1\), and \(t=-1\).

For \(j=0,1,\ldots,N-1\) in order:
\begin{enumerate}
\item Compute \(\delta=s_j+\sum_{i=1}^{L}c_i s_{j-i}\); when \(i>\deg C\), take \(c_i=0\).
\item If \(\delta=0\), proceed to the next iteration.
\item If \(\delta\ne0\), first store \(D=C\) and \(L_0=L\), then set
\[
C\leftarrow C-\frac{\delta}{b}Z^{j-t}B.
\]
\item If \(2L_0\le j\), also set
\[
L\leftarrow j+1-L_0,\qquad B\leftarrow D,\qquad b\leftarrow\delta,\qquad t\leftarrow j.
\]
\end{enumerate}
Output \((L,C)\) and the monic polynomial \(H(T)=T^L C(T^{-1})\).
\end{quote}

\noindent\phantomsection\label{v101:par-bm-boundaries}\ignorespaces
At the start of iteration \(j\), one has \(L\le j\), so the discrepancy uses no negative sequence index. The stored \(b\) is nonzero after initialization and after every update, so division by \(b\) is defined. For \(N=0\) or an all-zero input, the output is \((0,1)\). A zero final coefficient does not change the connection length \(L\); it must not be replaced by \(\deg C\).

\begin{theorem}\label{thm:berlekamp-massey}
The algorithm outputs the minimum linear complexity of the given prefix and a corresponding connection relation.
If the final complexity is \(L_N\), a straightforward implementation uses \(O(N(L_N+1))\) arithmetic operations in \(k\);
apart from storing the input, the stored polynomials contain \(O(L_N+1)\) coefficients.
\end{theorem}
\begin{proof}
Induct on the number of processed terms. The zero prefix has minimum length zero; if no discrepancy ever occurs, the algorithm keeps \(L=0,C=1\) and the conclusion follows.
If the first nonzero term occurs at position \(j\), the update gives \(C=1-s_jZ^{j+1}\) and \(L=j+1\).
There are then no recurrence relations to check, so this is feasible; Lemma~\ref{lem:bm-discrepancy-bound} requires length at least \(j+1\), proving minimality.

Thereafter, suppose that the stored \(B\) first failed at position \(t\) with discrepancy \(b\ne0\), when the minimum length was \(\ell\).
The update at the most recent complexity increase ensures that
\[
L=t+1-\ell,\qquad \deg B\le\ell,
\]
and \(B\) satisfies all connection relations at indices from \(\ell\) through \(t-1\).
If the current discrepancy is zero, the old relation remains valid. Any shorter relation would also apply to the old prefix, contradicting the induction hypothesis.

If the current discrepancy \(\delta\ne0\), every new relation has length at least the old minimum \(L\) and, by the discrepancy lemma, at least \(j+1-L\).
Its length is therefore at least
\[
L_1=\max\{L,j+1-L\}.
\]
The polynomial \(C'=C-(\delta/b)Z^{j-t}B\) constructed by the algorithm has degree at most \(L_1\), since \(j-t+\ell=j+1-L\). For every \(L_1\le v<j\), the old relation \(C\) is valid at position \(v\),
while the shifted \(B\) corresponds to position \(v-j+t\), satisfying
\[
\ell=t+1-L\le v-j+t<t.
\]
Thus the correction term vanishes at these positions. At position \(j\), its discrepancy is \((\delta/b)b=\delta\), exactly cancelling the old discrepancy.
Hence \((L_1,C')\) is feasible and attains the lower bound.
When \(L_1>L\), equivalently \(2L\le j\), storing the old relation and the current discrepancy restores the invariant; otherwise the previous record is retained.
This proves correctness of all updates and of the final output.

Complexity is nondecreasing. Each discrepancy computation and polynomial update costs \(O(L_N+1)\) field operations,
and polynomial lengths are also bounded in terms of \(L_N\), giving the stated costs.
\end{proof}

\begin{corollary}[Information needed to recover an infinite sequence]\label{cor:bm-recovery}
Let \(D\ge0\), and suppose that the minimal polynomial of an infinite sequence is known to have degree at most \(D\). Run the algorithm on its first \(2D\) terms.
The output \(H(T)\) is the minimal polynomial of the entire sequence.
\end{corollary}
\begin{proof}
Let the true minimal polynomial be \(h\), of degree \(d\le D\). If \(d=0\), then \(h=1\) and the sequence is zero, so the algorithm also returns \(H=1\). Assume henceforth that \(d\ge1\). It gives a prefix relation of length \(d\), so the output length satisfies \(L\le d\).
Put \(w=H(S)s\). Since \(H(S)\) and \(h(S)\) commute, \(w\) also satisfies \(h(S)w=0\).
The connection relation ensures that the first \(2D-L\ge d\) terms of \(w\) vanish.
A monic recurrence of degree \(d\) is determined by its first \(d\) terms, so \(w=0\); thus \(H\) annihilates the entire sequence.
Minimality gives \(h\mid H\), which together with \(\deg H=L\le d\) yields \(H=h\).
\end{proof}

Recovering the minimal polynomial of the entire sequence requires a known degree bound. A finite prefix can have arbitrary continuations,
so its shortest connection relation need not annihilate the entire sequence.
The recovery statement also applies when \(H(0)=0\). For example, the sequence \(1,0,0,\ldots\) has minimal polynomial \(T\), with connection length one and constant connection polynomial one.

Section~\ref{v118:sec-bch-decoding} applies this recovery result to syndrome sequences to determine error locations in codewords.

\begin{example}[Recovering the recurrence of a degree-six orbit]\label{exm:bm-f64}
Use the field and element \(\alpha\) of Example~\ref{exm:f64-computation}, and put \(s_j=\Tr_{\F_{64}/\F_2}(\alpha^j)\).
The first 12 terms are
\[
0,0,0,0,0,1,0,0,0,0,1,1.
\]
Grouping consecutive steps with the same connection length \(L\) and connection polynomial \(C\) gives the following computation:
\begin{center}
\begin{tabular}{ccl}
\toprule
Processed positions \(j\)&Complexity \(L\)&Connection polynomial \(C(Z)\)\\
\midrule
\(0\)--\(4\)&\(0\)&\(1\)\\
\(5\)--\(9\)&\(6\)&\(1+Z^6\)\\
\(10\)--\(11\)&\(6\)&\(1+Z^5+Z^6\)\\
\bottomrule
\end{tabular}
\end{center}
Thus \(H(T)=T^6+T+1\). This agrees with the trace representation of the multiplicative orbit:
\(s\) is a nonzero solution of an irreducible recurrence, with minimal polynomial equal to that of \(\alpha\) and period \(63\).
\end{example}

\par\Needspace{6\baselineskip}
\section*{Exercises}
\phantomsection
\addcontentsline{toc}{section}{Exercises}

\begin{exercise}\label{ex:sequence-trace8}
Suppose that \(\alpha^3+\alpha+1=0\) in \(\F_8\). Prove that the sequence of Example~\ref{exm:m-sequence-seven} is \(s_t=\Tr_{\F_8/\F_2}(\alpha^t)\), and show that all nonzero solution sequences are its cyclic shifts.
\end{exercise}

\begin{exercise}\label{ex:sequence-repeated}
If a sequence over \(\F_2\) has minimal polynomial \((T^2+T+1)^3\), find its least period.
Then prove that, for every monic polynomial \(h\) with \(h(0)\ne0\), there exists a sequence with minimal polynomial \(h\).
\end{exercise}

\begin{exercise}\label{ex:sequence-generating}
Find the generating function of the binary period-seven sequence in this chapter, and explain its minimal polynomial using the denominator.
\end{exercise}

\begin{exercise}[Eventual periodicity and the zero-root factor]\leavevmode\phantomsection\label{ex:wu-transient}
Consider the sequence \(1,0,1,1,1,\ldots\) over \(\F_2\).
Prove that its minimal polynomial is \(T^2(T+1)\), whereas the reduced generating-function denominator is \(1+Z\).
Explain why recovering a minimal polynomial from a generating-function denominator requires accounting for a finite initial segment.
\end{exercise}

\begin{exercise}[Separating transients and computing periods in a cyclic module]\leavevmode\phantomsection\label{ex:wu10-transient-repeated}
Over \(k=\F_2\), put
\[
f(T)=T^9+T^3=T^3(T+1)^2(T^2+T+1)^2,
\qquad A=k[T]/(f),\qquad \alpha=\overline T.
\]
Let \(\lambda:A\to k\) extract the coefficient of \(T^8\) from the canonical remainder.
For \(c\in k[T]\), define the sequence \(s^{(c)}_j=\lambda(c(\alpha)\alpha^j)\).
\begin{enumerate}[label=\textup{(\arabic*)}]
\item Use nondegeneracy of the coefficient pairing to prove
\[
h_{s^{(c)}}=\frac{f}{\gcd(f,c)}.
\]
For \(c=T(T^2+T+1)\), compute this minimal polynomial and the first six terms,
find the minimal preperiod and least eventual period, and explain why neither can be shortened.
\item Write the resulting minimal polynomial as \(T^2g(T)\) and find an explicit B\'ezout identity
\(U(T)T^2+V(T)g(T)=1\).
Use it to express the finitely supported and purely periodic components of \(s^{(c)}\) in terms of the shift operator,
and give the values of both components explicitly.
\item Now take \(c=1\) and \(c=T^3(T+1)^2(T^2+T+1)\),
and compute the minimal polynomial, minimal preperiod, and least eventual period of each sequence.
All these sequences satisfy \(f(S)s=0\). Explain why the multiplicity of the zero root and the order of the nonzero-root part of an arbitrary annihilating polynomial
cannot be identified directly with the minimal preperiod and least eventual period of the sequence.
\end{enumerate}
\end{exercise}

\begin{exercise}[Prescribing recurrence complexity and period simultaneously]\label{ex:wu2-recurrence-profile}
Over \(\F_2\), take \(f=(T+1)^2(T^2+T+1)\). Find the number of sequences with minimal polynomial exactly \(f\) and their common least period;
then find the number of all solution sequences with least period three.
\end{exercise}

\begin{exercise}[Finite observations and connection length]\label{ex:wu3-bm-short}
For the binary prefix \(1,0,0,1,0,1\), find a shortest connection relation for every nonempty initial segment and prove that its length is minimal.
Explain why the connection polynomial for the first two terms may be the constant one, while its connection length cannot be zero.
\end{exercise}

\begin{exercise}[The worst case for the recovery bound]\label{ex:wu3-bm-sharp}
For every \(D\ge1\), construct over any finite field two sequences whose minimal polynomials both have degree \(D\),
whose first \(2D-1\) terms agree, but whose \(2D\)th terms differ.
\end{exercise}

\begin{exercise}[Three types of period data for the same element]\label{ex:wu7-f64-two-polynomials}
Retain \(\alpha\) from Example~\ref{exm:two-sixtyfour-rootsets}. Prove
\[
m_{\alpha,\F_2}(T)=T^6+T+1,\qquad
\mu_\alpha(T)=\frac{T^6-1}{T+1}=(T+1)(T^2+T+1)^2.
\]
Use this to find the Frobenius orbit length of \(\alpha\), the dimension of the span of that orbit, and the multiplicative order of \(\alpha\).
Then explain why the sequence \(s_j=\Tr_{\F_{64}/\F_2}(\alpha^j)\) has minimal polynomial \(m_{\alpha,\F_2}\)
rather than \(\mu_\alpha\).
\end{exercise}

\chapter{Root Sets and Coding}
\label{chap:codes}
Root sets and finite-field linear algebra turn distance requirements into vanishing conditions, column relations, or bounds on root-space dimensions.
We begin with linear codes, two weights, and distance bounds, then study linearized evaluation, cyclic and constacyclic structures, consecutive root constraints, support root polynomials, and projective evaluation.
Trace representations of recurrences were established in Chapter~\ref{chap:sequences}. Here a period becomes a coordinate vector: root orbits determine code spaces, while root sets recover coefficient fields and permutation symmetries.
Consecutive-root constraints first exclude low-weight codewords and then eliminate the unknown codeword. The remaining syndromes are weighted evaluations at error locations; the pairing recovery of Section~\ref{str:sec-trace-recovery} and the recurrence algorithm of the preceding chapter determine the locations and values.

\section{Linear Codes, Weights, and Distance Bounds}
\label{v100:sec-code-foundations}
\subsection{Hamming Distance, Duality, and the Column Criterion}
\leavevmode\phantomsection\label{v40:def-linear-code}\ignorespaces
A \textbf{linear code} over a field \(K\) is a linear subspace of \(K^N\). The integer \(N\) is the \textbf{code length}, and the vector-space dimension is the \textbf{dimension of the code}. Its elements are called \textbf{codewords}. A code of length \(N\) and dimension \(r\) is briefly called an \([N,r]\) code.

\par\Needspace{3\baselineskip}\noindent\phantomsection\label{v44:def-hamming-weight}\ignorespaces
The \textbf{Hamming weight} of a codeword \(c=(c_0,\ldots,c_{N-1})\) is
\(\operatorname{wt}_H(c)=\#\{j:c_j\ne0\}\), and the \textbf{Hamming distance} between two codewords is \(\operatorname{wt}_H(c-c')\).
\leavevmode\phantomsection\label{v44:def-minimum-distance}\ignorespaces
The \textbf{minimum Hamming distance} of a nonzero linear code \(C\) is
\[
d_H(C)=\min_{0\ne c\in C}\operatorname{wt}_H(c),
\]
equivalently, the minimum distance between two distinct codewords. A linear code of length \(N\), dimension \(K\), and minimum distance \(d\) is denoted by \([N,K,d]\).

In the following discussion take \(k=\F_q\), and write \(N,K\) for length and dimension.
\par\Needspace{3\baselineskip}\noindent\phantomsection\label{v100:def-generator-matrix}\ignorespaces
A matrix whose rows form a basis of a code \(C\) is its \textbf{generator matrix}.
\par\Needspace{3\baselineskip}\noindent\phantomsection\label{v100:def-code-dual}\ignorespaces
With respect to the standard inner product, its \textbf{dual code} is
\[
 C^\perp=\{y\in k^N:\textstyle\sum_{j=0}^{N-1}c_jy_j=0\text{ for all }c\in C\}.
\]
\par\Needspace{3\baselineskip}\noindent\phantomsection\label{v100:def-information-set}\ignorespaces
A set of \(K\) coordinates is an \textbf{information set} if restriction to those coordinates gives an isomorphism \(C\to k^K\).
The next proposition gives the Hamming Singleton bound and the column criterion for equality.

\leavevmode\phantomsection\label{v84:def-mds}\ignorespaces
An \([N,K,d]\) linear code is \textbf{maximum distance separable}, or \textbf{MDS}, if \(d=N-K+1\).
It is called nontrivial when \(2\le K\le N-2\).
The following linear-algebraic criterion translates this equality into independence of a root set.

\begin{proposition}[Distance, column ranks, and scalar extension]\label{v84:prop-mds-linear}
Let \(G\in k^{K\times N}\) have full row rank, with \(1\le K\le N\), and let \(C\) be its row space. Then
\begin{equation}\label{v84:eq-distance-rank}
 d_H(C)=N-\max\{|I|:\rank G_I<K\}\le N-K+1.
\end{equation}
Equality holds if and only if every \(K\) columns are independent, equivalently every set of at most \(K\) columns is independent.
For every extension \(L/k\), the code \(C_L\) generated by the same \(G\) over \(L\) has the same dimension and distance, and
\(C_L\cap k^N=C\) and \((C^\perp)_L=(C_L)^\perp\).
If \(C\) is MDS and \(K<N\), then \(C^\perp\) is an \([N,N-K,K+1]\) MDS code.
If \(C\) is MDS, imposing zeros on any \(0\le s<K\) coordinates and then deleting these coordinates gives a \textbf{shortened code} with parameters \([N-s,K-s,N-K+1]\).
\end{proposition}
\begin{proof}
The condition \(\rank G_I<K\) says precisely that some \(0\ne u\in k^K\) gives a word \(uG\) vanishing on \(I\), which proves the distance formula.
Any \(K-1\) columns have rank less than \(K\), giving the Singleton bound. A smaller dependent set can be enlarged to a dependent set of \(K\) columns, proving the equivalence of the independence conditions.
The rank of every submatrix is determined by its nonzero minors, which are unchanged by scalar extension. Hence dimension and distance are preserved.
Choose an invertible submatrix \(G_I\) of \(K\) columns. If \(uG\in k^N\), then \(u=(uG)_IG_I^{-1}\in k^K\), proving the intersection formula.
The extended dual contains \((C^\perp)_L\), and both spaces have dimension \(N-K\), so they agree.

For an MDS code, any \(K\) coordinates allow a generator matrix of the form \((I_K\mid A)\); a generator of the dual is \((-A^{\mathsf T}\mid I_{N-K})\).
Since the original \(K\) coordinates were arbitrary, their complement is always an information set for the dual, which therefore satisfies the column criterion.
For shortening, the selected \(s\) columns impose zero-coordinate conditions of rank \(s\), leaving dimension \(K-s\).
Deleting zero coordinates preserves weight. The resulting distance is at least \(N-K+1\), and the Singleton bound gives the reverse inequality.
\end{proof}

\subsection{Rank Weight and the Rank-Metric Singleton Bound}
Let \(E=\F_{q^m}\), with \(1\le K\le n\le m\).
\leavevmode\phantomsection\label{v99:def-rank-weight}\ignorespaces
For \(c=(c_1,\ldots,c_n)\in E^n\), its \textbf{rank weight} is
\[
\operatorname{wt}_R(c)=\dim_{\F_q}\langle c_1,\ldots,c_n\rangle.
\]
\leavevmode\phantomsection\label{v99:def-rank-distance}\ignorespaces
The \textbf{rank distance} between two vectors is \(\operatorname{wt}_R(c-c')\). The \textbf{minimum rank distance} of a nonzero \(E\)-linear code \(\mathcal C\) is
\[
d_R(\mathcal C)=\min_{0\ne c\in\mathcal C}\operatorname{wt}_R(c).
\]
\leavevmode\phantomsection\label{v99:par-rank-singleton}\ignorespaces
If \(\dim_E\mathcal C=K\), imposing zero values in any \(K-1\) coordinates gives a linear map \(\mathcal C\to E^{K-1}\) with nonzero kernel.
A nonzero vector in this kernel has at most \(n-K+1\) nonzero coordinates, so the span of its coordinates has dimension at most \(n-K+1\). Hence
\[
d_R(\mathcal C)\le n-K+1.
\]
This is the \textbf{rank-metric Singleton bound} for \(E\)-linear codes when \(n\le m\).
\leavevmode\phantomsection\label{v99:def-mrd}\ignorespaces
A code attaining it is called a \textbf{maximum rank distance code}.

\section{Linearized Evaluation and Rank-Metric Codes}
\label{v90:subsec-rankmetric}
Continue with \(E=\F_{q^m}\) and \(1\le K\le n\le m\).
Choose \(\F_q\)-linearly independent \(g_1,\ldots,g_n\in E\).
\begin{theorem}[Maximum rank distance from linearized evaluations]\label{v90:thm-rank-metric}
The set
\[
\mathcal G_K=\left\{(L(g_1),\ldots,L(g_n)):\ L(X)=\sum_{i=0}^{K-1}a_iX^{q^i},\ a_i\in E\right\}
\]
is an $E$-linear code of dimension $K$, with minimum nonzero rank weight
\[
\boxed{d_R(\mathcal G_K)=n-K+1.}
\]
It therefore meets the rank-metric Singleton bound.
\end{theorem}
\begin{proof}
Let $V=\langle g_1,\ldots,g_n\rangle_{\F_q}$. For any nonzero $L$ above, its root space has dimension at most $K-1$, so
\[
\dim_{\F_q}L(V)=n-\dim_{\F_q}(V\cap\ker L)\ge n-K+1.
\]
The coordinates of the evaluation vector span exactly \(L(V)\), giving the lower bound on rank weight. In particular, evaluation is injective on these polynomials, so \(\dim_E\mathcal G_K=K\).
Conversely, choose an \((K-1)\)-dimensional \(\F_q\)-subspace \(U\subseteq V\). By Theorem~\ref{thm:subspace-polynomial},
\[
L_U(X)=\prod_{u\in U}(X-u)
\]
is a linearized polynomial of \(q\)-degree \(K-1\), with coefficients in \(E\) and root set exactly \(U\). It therefore belongs to the family being evaluated and satisfies
\[
\dim_{\F_q}L_U(V)=n-\dim_{\F_q}U=n-K+1.
\]
Its nonzero evaluation vector attains the lower bound. When \(K=1\), take \(U=\{0\}\) and \(L_U=X\); the same argument applies.
\end{proof}
\begin{remark}
These are commonly called Gabidulin codes. Their construction and distance follow from Moore independence, the linearized root-space bound, and subspace polynomials.
Rank weight measures the dimension of the coordinate span over \(\F_q\); Hamming weight, defined above, counts the nonzero coordinates. These are different measurements.
\end{remark}

\begin{example}[Length three and dimension two over the eight-element field]\label{v103:ex-rank-eight}
Let \(E=\F_8=\F_2(\alpha)\), where \(\alpha^3+\alpha+1=0\).
Take \((g_1,g_2,g_3)=(1,\alpha,\alpha^2)\), which are linearly independent over \(\F_2\).
Evaluating \(L(X)=aX+bX^2\), with \(a,b\in E\), gives the generator matrix
\[
G=\begin{pmatrix}
1&\alpha&\alpha^2\\
1&\alpha^2&\alpha^4
\end{pmatrix}
=\begin{pmatrix}
1&\alpha&\alpha^2\\
1&\alpha^2&\alpha^2+\alpha
\end{pmatrix}.
\]
This is a Moore matrix. The minor on the first two columns is \(\alpha^2+\alpha\ne0\), so the code \(C\) has dimension two over \(E\).
For \(L(X)=X^2+(\alpha+1)X\), the resulting word is
\[
c=(L(1),L(\alpha),L(\alpha^2))=(\alpha,\alpha,1).
\]
All three entries are nonzero, so \(\operatorname{wt}_H(c)=3\).
Their span over the base field is only \(\langle1,\alpha\rangle_{\F_2}\), giving \(\operatorname{wt}_R(c)=2\).
Writing the entries as coordinate columns in the basis \((1,\alpha,\alpha^2)\) gives
\[
\begin{pmatrix}0&0&1\\1&1&0\\0&0&0\end{pmatrix},
\]
of binary rank two. The general theorem gives \(d_R(C)\ge3-2+1=2\), and this word attains equality.
Thus the minimum rank distance is witnessed explicitly: Hamming weight counts nonzero positions, while rank weight measures the base-field dimension spanned by the entries.
\end{example}

\begin{example}[Column normalization and the two distances]\label{v105:ex-rank-column-normalization}
Continue with the notation of Example~\ref{v103:ex-rank-eight}, and set
\[
 D=\operatorname{diag}(1,\alpha^{-1},\alpha^{-2}),
 \qquad C'=CD=\{cD:c\in C\}.
\]
Dividing each column by its first entry is called \textbf{column normalization}. Here it gives
\begin{equation}\label{v105:eq-rank-column-normalization}
 GD=
 \begin{pmatrix}
 1&1&1\\
 1&\alpha&\alpha^2
 \end{pmatrix}.
\end{equation}
Since \(D\) is invertible, \(C'\) still has length three and dimension two over \(E\).
The preceding example gives \(d_R(C)=2\).
The dimension of the coordinate span of a vector is at most its number of nonzero entries, so
\(d_H(C)\ge d_R(C)=2\). The Hamming Singleton bound gives \(d_H(C)\le3-2+1=2\), hence \(d_H(C)=2\).
Multiplying the coordinates separately by nonzero elements does not change the zero positions, so \(d_H(C')=d_H(C)=2\).
However, \((1,1,1)\in C'\) has binary rank weight one. Every nonzero vector has rank weight at least one, and therefore
\[
 \boxed{d_H(C)=d_H(C')=2,\qquad
 d_R(C)=2,\quad d_R(C')=1.}
\]

The matrix \(G\) is the Moore matrix formed from \(1,\alpha,\alpha^2\), whereas \(GD\) is a two-row Vandermonde matrix on the same nodes.
Explicitly, the words of \(C'\) are
\[
 \bigl(a+b,\ a+b\alpha,\ a+b\alpha^2\bigr),
 \qquad a,b\in E.
\]
Thus \(C'\) is the Reed--Solomon code obtained by evaluating ordinary polynomials of degree at most one at \(1,\alpha,\alpha^2\).
The invertible column scaling preserves Hamming distance, but it does not preserve the joint linear relations among the coordinates over the base field.

More generally, for \(E=\F_{q^m}\), \(c\in E^n\), \(a\in E^\times\), and \(B\in\operatorname{GL}_n(\F_q)\), we have
\[
 \operatorname{wt}_R(acB)=\operatorname{wt}_R(c).
\]
Indeed, applying \(B\) and \(B^{-1}\) shows that the coordinates of \(cB\) and \(c\) span the same \(\F_q\)-space, and multiplication by the common scalar \(a\) gives a linear isomorphism from that span to its scalar multiple.
Separate multiplication by arbitrary nonzero extension-field elements need not preserve rank weight, as~\eqref{v105:eq-rank-column-normalization} shows.
The one-dimensional right-remainder coordinates in~\eqref{v105:eq-linearized-point-remainder} are exactly \(L(g_j)/g_j\).
Their relation to the evaluation coordinates \(L(g_j)\) is another instance of this column normalization, for which the two distances must be distinguished.
\end{example}

\section{Root Orbits, Cyclic Codes, and Fourier Coordinates}
\label{v40:sec-cyclic-trace}
\label{v50:sec-cyclic-trace-model}
Restricting a purely periodic sequence to one period gives a vector of finite length, and the shift becomes a cyclic shift on these vectors.
In this model, the decomposition into root orbits also gives the decomposition of a shift-invariant space.

Here the coordinate field is \(k=\F_q\).
\leavevmode\phantomsection\label{v40:def-cyclic-code}
A linear code invariant under the cyclic right shift
\[
R_N(c_0,\ldots,c_{N-1})=(c_{N-1},c_0,\ldots,c_{N-2})
\]
is called a \textbf{cyclic code}. Representing a codeword by \(\sum_{j=0}^{N-1}c_jX^j\) identifies the cyclic right shift, in
\(A_N=k[X]/(X^N-1)\), with multiplication by \(X\); thus cyclic codes are exactly the ideals of \(A_N\).
The transformations induced by \(R_N\) and the sequence shift \((Ss)_j=s_{j+1}\) on one period are mutually inverse.
Assume henceforth that \(p\nmid N\), and fix a primitive \(N\)th root of unity \(\zeta\).

\begin{proposition}[Root constraints and cyclic codes]\label{v40:prop-cyclic-root-constraints}
Every cyclic code has a unique expression \(C=(g)\), where \(g\in k[X]\) is monic and \(g\mid X^N-1\). Its dimension is \(N-\deg g\).
\leavevmode\phantomsection\label{v40:def-generator-check}
The polynomial \(g\) is called the \textbf{generator polynomial}, and \(h=(X^N-1)/g\) the \textbf{check polynomial}.
\leavevmode\phantomsection\label{v40:def-defining-set}
The \textbf{defining set} of the code is
\[
T(C)=\{i\in\Z/N\Z:g(\zeta^i)=0\}.
\]
Cyclic codes correspond bijectively to \(q\)-stable exponent sets, explicitly by
\[
g_T(X)=\prod_{i\in T}(X-\zeta^i),\qquad \dim_k C=N-|T|.
\]
\end{proposition}
\begin{proof}
The preimage of an ideal under the quotient map \(k[X]\to A_N\) is generated by a unique monic polynomial \(g\) and contains \(X^N-1\), so \(g\mid X^N-1\).
Put \(h=(X^N-1)/g\). The residue of any \(ag\) modulo \(X^N-1\) depends only on \(a\bmod h\).
Moreover, \(ag=0\) if and only if \(h\mid a\). If \(h=1\), then \(C=0\); otherwise \(g,Xg,\ldots,X^{\deg h-1}g\) is a basis of \(C\).
Since \(p\nmid N\), the polynomial \(X^N-1\) is square-free; Theorem~\ref{thm:coset-factorization} identifies its irreducible factors with complete Frobenius orbits.
Hence monic factors are determined precisely by \(q\)-stable root sets, and their degrees equal their numbers of roots, proving the last two equalities.
\end{proof}

This is also a concrete application of the structure of finite reduced algebras. The Chinese remainder theorem gives
\[
A_N\simeq\prod_C k[X]/(M_C)\simeq\prod_C\F_{q^{|C|}},
\]
where \(C\) runs through the cyclotomic cosets modulo \(N\). An ideal is either zero or the whole field on each field component: the generator polynomial contains \(M_C\) precisely when that component is set to zero.
The orthogonal idempotents of Section~\ref{sec:finite-reduced-algebras} give the projections onto these components.
The root-set description records which components vanish, while a trace representation gives actual coordinates in a retained component.

\begin{proposition}[The trace model of a single root orbit]\label{v40:prop-cyclic-trace-model}
Let \(d=[k(\zeta):k]=\ord_N(q)\) and \(E=k(\zeta)=\F_{q^d}\), and let \(m_\zeta\) be the minimal polynomial of \(\zeta\) over \(k\).
Then
\[
\Phi(a)=\bigl(\Tr_{E/k}(a\zeta^{-j})\bigr)_{j=0}^{N-1},\qquad a\in E,
\]
is a linear isomorphism from \(E\) to a \(d\)-dimensional cyclic code \(C_\zeta\), and
\[
R_N\Phi(a)=\Phi(a\zeta),\qquad
C_\zeta=\left(\frac{X^N-1}{m_\zeta(X)}\right).
\]
Thus the check polynomial of \(C_\zeta\) is \(m_\zeta\).
\end{proposition}
\begin{proof}
If \(\Phi(a)=0\), then the trace pairings of \(a\) with all \(1,\zeta^{-1},\ldots,\zeta^{-(d-1)}\) vanish.
These elements form a power basis of \(E/k\), so nondegeneracy of the trace pairing gives \(a=0\). Thus \(\Phi\) is injective and its image has dimension \(d\).
The identity \(\zeta^N=1\) gives \(R_N\Phi(a)=\Phi(a\zeta)\) coordinatewise, so the image is a cyclic code.

For every \(P\in k[X]\), this equality gives \(P(R_N)\Phi(a)=\Phi(P(\zeta)a)\).
By injectivity, \(P(R_N)\) vanishes on the whole image if and only if \(P(\zeta)=0\), equivalently, \(m_\zeta\mid P\).
On the other hand, if a cyclic code has generator polynomial \(g\) and check polynomial \(h\), then
\(P(R_N)C=0\) if and only if \(X^N-1=gh\mid Pg\), that is, \(h\mid P\).
Comparing the two annihilator ideals yields \(h=m_\zeta\) and the stated generator polynomial.
\end{proof}

The Frobenius orbit length \(d\) of the root gives the code dimension, while its multiplicative order \(N\) gives the cyclic shift period of each nonzero codeword:
by injectivity of \(\Phi\), the equality \(R_N^r\Phi(a)=\Phi(a)\) holds for \(a\ne0\) if and only if \(\zeta^r=1\).

\leavevmode\phantomsection\label{v40:def-irreducible-cyclic-code}
A nonzero cyclic code containing no proper nonzero cyclic subcode is called an \textbf{irreducible cyclic code}.
In the preceding product of fields, it retains exactly one field component, so this is equivalent to irreducibility of the check polynomial.
Proposition~\ref{v40:prop-cyclic-trace-model} gives trace coordinates for this component; its generator polynomial is generally reducible.
More generally, any \(N\)th root of unity \(\eta\) gives, by the same trace construction over \(k(\eta)/k\), a cyclic code of length \(N\) with check polynomial \(m_\eta\), the minimal polynomial of \(\eta\) over \(k\).
If \(e=\ord(\eta)<N\), the coordinate block of length \(e\) repeats \(N/e\) times, while the check polynomial remains \(m_\eta\).
Every nonzero codeword has least period \(\ord(\eta)\) under the cyclic shift.

\leavevmode\phantomsection\label{v40:def-finite-field-dft}
The evaluation map \(c\mapsto(c(\zeta^i))_{i\in\Z/N\Z}\) gives an extension-field-valued \textbf{discrete Fourier transform}: it applies the Fourier transform on a cyclic group defined on page~\pageref{v40:def:cyclic-fourier} to codeword coefficients.
\begin{proposition}[Frobenius descent of Fourier coordinates]\label{v96:prop-fourier-descent}
Suppose \(p\nmid N\), and let \(E/k\) be a finite extension containing a primitive \(N\)th root of unity \(\zeta\).
For \(c=(c_j)_{j=0}^{N-1}\in E^N\), put
\[
\widehat c(i)=\sum_{j=0}^{N-1}c_j\zeta^{ij}\qquad(i\in\Z/N\Z).
\]
Then
\begin{equation}\label{v96:eq-fourier-descent}
\boxed{c\in k^N\quad\Longleftrightarrow\quad
\widehat c(qi)=\widehat c(i)^q\quad\text{for every }i\in\Z/N\Z.}
\end{equation}
On a \(q\)-cyclotomic coset \(O\) of length \(d\), spectral data satisfying this condition are uniquely determined by one value in \(\F_{q^d}\).
Thus an orbit of length \(d\) contributes a \(d\)-dimensional \(k\)-space, and the full set of frequencies allowed to be nonzero must be a union of cyclotomic cosets.
\end{proposition}
\begin{proof}
If every \(c_j\in k\), then
\[
\widehat c(i)^q=\sum_jc_j\zeta^{qij}=\widehat c(qi).
\]
Conversely, Fourier inversion gives
\[
c_j=N^{-1}\sum_i\widehat c(i)\zeta^{-ij}.
\]
Since \(N\in k^\times\), raising this expression to the \(q\)th power and substituting \(r=qi\), which permutes \(\Z/N\Z\), yields
\[
c_j^q=N^{-1}\sum_i\widehat c(qi)\zeta^{-qij}
     =N^{-1}\sum_r\widehat c(r)\zeta^{-rj}=c_j.
\]
Hence \(c_j\in k\). After one circuit of an orbit of length \(d\), its representative value \(u\) must satisfy \(u^{q^d}=u\).
Conversely, any such value uniquely determines equivariant data on that orbit.
Since \(d\mid[E:k]\), these values form exactly \(\F_{q^d}\).
One nonzero value makes every value on the orbit nonzero, and different orbits may be assigned values independently.
\end{proof}
This is precisely the equivariance condition in Theorem~\ref{thm:equivariant-functions}, applied to the root set
\(\{\zeta^i:i\in\Z/N\Z\}\): the root \(\zeta^i\) is sent to \(\zeta^{qi}\), and its value is raised to the \(q\)th power.
Root-set descent and spectral descent for cyclic codes are therefore two descriptions of the same evaluation isomorphism.
Cyclic right shift becomes multiplication by \(\zeta^i\) in these coordinates, so selecting the root orbits to retain selects the corresponding cyclic-code components.

For the preceding single-orbit trace code, write its coefficient polynomial as
\(c_a(X)=\sum_{j=0}^{N-1}\Tr_{E/k}(a\zeta^{-j})X^j\). Expanding the trace and using the root-of-unity sum gives
\[
c_a(\zeta^i)=\sum_{t=0}^{d-1}a^{q^t}
\sum_{j=0}^{N-1}\zeta^{j(i-q^t)}
=\begin{cases}
N a^{q^t},&i\equiv q^t\pmod N\quad(0\le t<d),\\
0,&i\notin C_1.
\end{cases}
\]
The inner sum is \(N\) when \(i\equiv q^t\), and otherwise is zero by the geometric-series formula; here \(N\in k^\times\).
Thus a single root orbit is exactly the set of evaluation coordinates that may be nonzero in this code, and the trace parameter \(a\) expands into its Frobenius conjugates in these coordinates.
In particular, \(c_a(\zeta)=Na\), so the explicit map from this trace model to a single CRT field component is \(a\mapsto Na\).
This identifies the multiplication model of cyclic right shift, the trace model of recurrences, and the evaluation model of root sets.

\begin{example}\label{v40:ex-cyclic-seven}
Take \(k=\F_2\), \(N=7\), and denote the root of order seven from Example~\ref{ex:normal-f8-example} by \(\zeta\), so that \(\zeta^3+\zeta+1=0\). Here a trace period is written as a codeword and its retained root orbit is identified.
Then \(C_1=\{1,2,4\}\), and the associated three-dimensional irreducible cyclic code has
\[
h=X^3+X+1,\qquad g=\frac{X^7-1}{h}=X^4+X^2+X+1.
\]
Indeed, \(\Phi(1)=(1,1,1,0,1,0,0)\), and its coefficient polynomial is \(g\).
The generator polynomial contains the other two root orbits, whereas the check polynomial records the orbit retained by this trace model.
\end{example}

\section{Constacyclic Codes, Duality, and Sparse Multiples}
\label{v100:sec-constacyclic}
\begin{definition}[Constacyclic codes]\label{v54:def-constacyclic-code}
Let \(\lambda\in k^\times\). A linear code \(C\subseteq k^N\) invariant under
\[
(c_0,\ldots,c_{N-1})\longmapsto(\lambda c_{N-1},c_0,\ldots,c_{N-2})
\]
is called a \textbf{\(\lambda\)-constacyclic code}.
Under the coefficient correspondence, it is an ideal of \(R_{N,\lambda}=k[X]/(X^N-\lambda)\).
Thus there is a unique monic \(f\mid X^N-\lambda\) such that \(C=(f)\).
\end{definition}
The proof of this ideal description is the same as in the cyclic case of Proposition~\ref{v40:prop-cyclic-root-constraints}:
take the monic generator in \(k[X]\) of the preimage ideal, and use its containment of \(X^N-\lambda\) to obtain the divisibility.
Below take \(0<\tau=\deg f<N\). Every codeword has a unique representative of degree less than \(N\), and
\begin{equation}\label{v54:eq-code-polynomial-realization}
C=\{fq:q\in k[X],\ \deg q<N-\tau\},\qquad \dim_k C=N-\tau.
\end{equation}
The zero polynomial is understood to satisfy the degree bound. This expression also applies to generator polynomials with repeated roots.

\begin{theorem}[Full-support generator and check polynomials of constacyclic MDS codes]\label{v87:thm-constacyclic-mds-reciprocal}
Let \(C=(g)\subseteq k[X]/(X^N-\lambda)\) be a \(\lambda\)-constacyclic code,
where \(\lambda\in k^\times\) and \(1\le K=\dim_k C<N\). Write
\[
 X^N-\lambda=g(X)h(X),\qquad
 \deg g=N-K,\quad\deg h=K,
\]
with \(g,h\) monic. Since \(h(0)\ne0\), its \textbf{monic reciprocal polynomial}
\[
 h^\vee(X)=\frac{X^Kh(X^{-1})}{h(0)}
\]
is well defined (also commonly denoted by \(h^*(X)\)). Then \(C^\perp\) is the
\(\lambda^{-1}\)-constacyclic code generated by \(h^\vee\), and
\begin{equation}\label{v87:eq-constacyclic-dual-multiple}
 \boxed{\quad C\text{ is MDS}\ \Longleftrightarrow\
 \min_{\substack{0\ne F\in k[X]\\h^\vee\mid F,\ \deg F<N}}\operatorname{wt}(F)=K+1.\quad}
\end{equation}
In particular, if \(C\) is MDS, then
\[
 \boxed{\operatorname{wt}(g)=N-K+1,\qquad \operatorname{wt}(h)=K+1.}
\]
Thus every coefficient from the constant term to the leading term of each polynomial is nonzero.
\end{theorem}
\begin{proof}
The polynomial \(h^\vee\) divides \(X^N-\lambda^{-1}\) and generates a
\(\lambda^{-1}\)-constacyclic code \(D\) of dimension \(N-K\).
Take \(c=ga\in C\), \(d=h^\vee b\in D\), with \(\deg a<K\) and \(\deg b<N-K\).
The usual inner product of the coefficient vectors is the constant coefficient of the Laurent polynomial
\(c(X)d(X^{-1})\). Since \(gh=X^N-\lambda\),
\[
 c(X)d(X^{-1})
 =\frac{X^N-\lambda}{h(0)}\,X^{-K}a(X)b(X^{-1}).
\]
Every exponent in \(X^{-K}a(X)b(X^{-1})\) lies between \(-(N-1)\) and \(-1\).
Multiplication by \(X^N-\lambda\) cannot produce a constant term; thus \(D\subseteq C^\perp\).
Both codes have dimension \(N-K\), so \(D=C^\perp\).
By Proposition~\ref{v84:prop-mds-linear}, \(C\) is MDS if and only if \(C^\perp\) is
an \([N,N-K,K+1]\) MDS code. Its standard polynomial representatives of degree less than \(N\)
are exactly the corresponding multiples of \(h^\vee\), proving~\eqref{v87:eq-constacyclic-dual-multiple}.

If \(C\) is MDS, then its nonzero codeword \(g\), of degree \(N-K\), satisfies
\(N-K+1\le\operatorname{wt}(g)\le N-K+1\).
Likewise, \(h^\vee\) is a nonzero word of the MDS dual and has degree \(K\), so
\(\operatorname{wt}(h^\vee)=K+1\). Taking the reciprocal preserves coefficient weight,
and hence \(\operatorname{wt}(h)=K+1\).
\end{proof}

For any \(K\times N\) generator matrix \(G\), if \(F=\sum_{j=0}^{N-1}f_jX^j\), then
\[
 h^\vee\mid F\quad\Longleftrightarrow\quad
 G(f_0,\ldots,f_{N-1})^{\mathsf T}=0.
\]
Thus the divisibility condition describes \emph{all column relations} of the original generator matrix,
with no irreducibility assumption on \(h\).

The \emph{dual code itself} should be distinguished from a \emph{weight-preserving reversal of coordinates}.
For \(\deg F<N\), put \(\operatorname{rev}_N(F)=X^{N-1}F(X^{-1})\).
If \(F=hQ\) with \(\deg Q<N-K\), then
\[
 \operatorname{rev}_N(hQ)
 =h(0)h^\vee(X)X^{N-1-K}Q(X^{-1}).
\]
Reversal is an involution preserving Hamming weight. Thus the minimum in
\eqref{v87:eq-constacyclic-dual-multiple} also equals the minimum weight among the
nonzero multiples of \(h\) of degree less than \(N\).
\textbf{However}, it is the multiples of \(h^\vee\), not generally the multiples of \(h\),
that directly represent words in \(C^\perp\). This distinction remains valid for repeated-root constacyclic codes.

The root convention for power vectors is equally important. For the right constacyclic shift \(T_\lambda\),
\[
 T_\lambda(1,\theta,\ldots,\theta^{N-1})
 =\theta^{-1}(1,\theta,\ldots,\theta^{N-1})
 \quad\Longleftrightarrow\quad\theta^N=\lambda^{-1}.
\]
Thus, if \(\alpha\) is a root of \(h\), then \(\alpha^N=\lambda\), and the corresponding
power vector uses \(\theta=\alpha^{-1}\), not \(\theta=\alpha\).
When \(p\nmid N\), choose a primitive \(N\)-th root of unity \(\zeta\) and write
\(\alpha^{q-1}=\zeta^s\). On the reciprocal root sets
\(\alpha\zeta^j\) and \(\alpha^{-1}\zeta^j\), Frobenius acts on the exponents by
\(j\mapsto qj+s\) and \(j\mapsto qj-s\), respectively.
These actions are conjugate via \(j\mapsto-j\); in particular, the orbit lengths are unchanged.

\begin{example}[Why the inverse root matters for constacyclic codes]\label{v87:ex-constacyclic-inverse-root}
Over \(k=\F_5\), take \(N=3\) and \(\lambda=2\). We have
\[
 X^3-2=(X-3)(X^2+3X+4),\qquad
 g=X-3,\quad h=X^2+3X+4,\quad h^\vee=X^2+2X+4.
\]
The dual code generated by \(h^\vee\) has dimension one and distance three, so
\(C=(g)\) is a \([3,2,2]\) MDS code.
The coefficient weights of \(g\) and \(h\) are two and three, respectively, and
\(h^\vee\ne h\). If \(\alpha\) is a root of \(h\), then \(\alpha^3=2\), whereas
\(\theta=\alpha^{-1}\), the parameter of the right-shift eigenvector, satisfies
\(\theta^3=3=\lambda^{-1}\).
\end{example}

\section{Consecutive Root Constraints and the BCH Bound}
\phantomsection\label{v56:subsec-bch-bound}
\leavevmode\phantomsection\label{v56:def-bch-code}\ignorespaces
Taking the union of the Frobenius orbits meeting an interval of consecutive exponents and using the corresponding roots of unity as the zeros of the generator polynomial gives a \textbf{BCH code}.
Explicitly, given integers \(b\) and \(2\le\delta\le N\), its generator polynomial is
\[
g_{b,\delta}(X)=\mathop{\operatorname{lcm}}_{0\le r\le\delta-2}m_{\zeta^{b+r}}(X),
\]
where \(m_{\zeta^{b+r}}\) is the monic minimal polynomial of \(\zeta^{b+r}\) over \(k\), and \(\delta\) is called the \textbf{designed distance}.
The \(G_b\) in~\eqref{v50:eq-leader-interval-polynomial} corresponds to \(b\) consecutive zeros starting at exponent one.
Consecutive zeros give the following distance bound.

\Needspace{7\baselineskip}
\begin{theorem}[BCH distance bound]\label{v56:thm-bch-bound}
Let \(p\nmid N\), let \(\zeta\) be a primitive \(N\)th root of unity, and let \(C\subseteq k^N\) be a nonzero cyclic code with generator polynomial \(g\).
If, for some integer \(b\) and \(2\le\delta\le N\),
\[
g(\zeta^{b+r})=0\qquad(0\le r\le\delta-2),
\]
then \(d_H(C)\ge\delta\).
\end{theorem}
\begin{proof}
Suppose there is a nonzero codeword of weight \(w<\delta\), and write its canonical remainder as
\[
c(X)=\sum_{i=1}^w a_iX^{j_i},\qquad
0\le j_1<\cdots<j_w<N,\qquad a_i\in k^\times.
\]
Since \(g\mid c\), the same zero conditions hold for \(c\).
Put \(x_i=\zeta^{j_i}\) and \(y_i=a_ix_i^b\). Taking the first \(w\) of these conditions gives
\[
\sum_{i=1}^w y_ix_i^r=0\qquad(0\le r<w).
\]
The \(x_i\) are distinct, so the Vandermonde determinant of the coefficient matrix is
\[
\det(x_i^{r})_{\substack{0\le r<w\\1\le i\le w}}
=\prod_{1\le i<j\le w}(x_j-x_i)\ne0.
\]
Thus all \(y_i\) vanish, contradicting \(a_i\ne0\) and \(x_i\ne0\).
\end{proof}
Proposition~\ref{v139:prop-digit-bch-gap} reads the same bound from the largest cyclic gap between omitted exponents in periodic-digit coordinates. Corollary~\ref{v139:cor-compressed-code} identifies complete root fibers with interleaved copies of a shorter code.

\Needspace{6\baselineskip}

\paragraph{A Hankel-rank interpretation of the BCH bound.}
In the proof of Theorem~\ref{v56:thm-bch-bound}, write
\(s_r=c(\zeta^{b+r})=\sum_{j\in\mathrm{supp}(c)}a_j\zeta^{bj}(\zeta^j)^r\).
The exponential nodes are nonzero and distinct. By
Theorem~\ref{v88:thm-hankel-vandermonde}, the minimum recurrence order and
the corresponding Hankel rank of this sequence equal the weight \(w\).
It cannot have \(w\) consecutive zero terms: the recurrence and the
invertible Hankel matrix exclude this. This is another interpretation of
the BCH bound; the direct Vandermonde system remains its shortest proof.

\begin{example}[A binary cyclic code from consecutive roots]\label{v56:ex-bch-seven}
Retain \(k=\F_2\), \(N=7\), and \(\zeta\) from Example~\ref{v40:ex-cyclic-seven}, but take
\(M(X)=X^3+X+1\) as the generator polynomial.
Its roots are \(\zeta,\zeta^2,\zeta^4\), so its defining set is \(\{1,2,4\}\).
The first two roots correspond to consecutive exponents one and two, and Theorem~\ref{v56:thm-bch-bound} gives \(d_H(C)\ge3\).
On the other hand, \(M\) itself is a nonzero codeword of weight three, so \(d_H(C)=3\).
Since \(\dim_k C=7-\deg M=4\), this gives a \([7,4,3]\) code.

Example~\ref{v40:ex-cyclic-seven} uses the same polynomial as a check polynomial, retaining the field component corresponding to \(C_1\) and giving dimension three;
here it is the generator polynomial, removing that component and giving dimension four.
\end{example}

\leavevmode\phantomsection\label{v68:rem-singular-fourier-minor}
The same example shows why the consecutive exponents in the BCH argument matter.
The polynomial \(c(X)=1+X+X^3\) has three nonzero coefficients and vanishes at \(\zeta,\zeta^2,\zeta^4\). Thus the Fourier submatrix
\[
\bigl(\zeta^{ij}\bigr)_{\substack{i\in\{1,2,4\}\\j\in\{0,1,3\}}}
\]
annihilates \((1,1,1)^{\mathsf T}\) and is singular. The seven-point Fourier transform of the coefficient vector of \(c\) has exactly four nonzero coordinates, since these three points are all the roots of \(c\).
Even at prime transform length, arbitrary square submatrices of a finite-field Fourier matrix need not be nonsingular.

Frobenius closure determines which root constraints descend to the base field, while consecutive exponents make the evaluation equations form an invertible Vandermonde submatrix.
Thus root sets determine not only cyclic-code dimensions but also the minimum number of nonzero coefficients.
The distinction from the Moore matrices of trace codes is made precise at~\ref{v96:rem-vandermonde-moore}.

\begin{corollary}[The BCH bound for square-free constacyclic codes]\label{v88:cor-constacyclic-bch-bound}
Assume \(p\nmid N\) and \(\lambda\in k^\times\). In a splitting field,
choose \(\omega\) with \(\omega^N=\lambda\) and a primitive \(N\)-th root
of unity \(\zeta\). Suppose a nonzero \(\lambda\)-constacyclic code
\(C=(f)\subseteq k[X]/(X^N-\lambda)\) has
\[
f(\omega\zeta^{b+r})=0\qquad(0\le r\le\delta-2),\qquad2\le\delta\le N.
\]
Then \(d_H(C)\ge\delta\). There is no requirement that \(\lambda\)
have an \(N\)-th root in \(k\).
\end{corollary}
\begin{proof}
If a nonzero codeword had weight \(w<\delta\), take its representative of
degree below \(N\), \(c(X)=\sum_{j=1}^w a_jX^{e_j}\), with nonzero
coefficients and distinct exponents. By \eqref{v54:eq-code-polynomial-realization}, this standard representative is divisible
by \(f\) in \(k[X]\), so it vanishes at the listed roots. For \(r=0,\ldots,w-1\), we get
\[
\sum_{j=1}^w\bigl(a_j\omega^{e_j}\zeta^{be_j}\bigr)(\zeta^{e_j})^r=0.
\]
The nodes \(\zeta^{e_j}\) are pairwise distinct, so the Vandermonde matrix
is invertible, forcing all \(a_j=0\), a contradiction.
\end{proof}

\subsection{Equality in the BCH Bound}\label{v123:sec-bch-equality}
The proof of the BCH bound gives an invertible square matrix when the weight is smaller than the designed distance. At equality, one row is missing; interpolation identifies the resulting one-dimensional kernel.

\begin{proposition}[Support weights under consecutive-root constraints]\label{v123:prop-bch-equality}
Let \(k=\F_q\), with \(p=\Char k\) and \(p\nmid N\), and let \(\zeta\) be a primitive \(N\)-th root of unity in a finite extension \(E/k\). For \(2\le w\le N\), write
\[
u(Z)=\sum_{i=1}^w a_iZ^{j_i},\qquad
a_i\in k^\times,\quad0\le j_1<\cdots<j_w<N,
\]
and put \(x_i=\zeta^{j_i}\), \(D(X)=\prod_i(X-x_i)\). For every integer \(b\),
\begin{equation}\label{v123:eq-bch-equality-weights}
u(\zeta^{b+r})=0\ (0\le r\le w-2)
\quad\Longleftrightarrow\quad
a_i=\frac{\lambda x_i^{-b}}{D'(x_i)}\ (1\le i\le w)
\end{equation}
for some \(\lambda\in E^\times\).
If such a word \(u\) belongs to a cyclic code \(C\) having these consecutive zeros, then \(d_H(C)=w\).
\end{proposition}
\begin{proof}
The zero conditions say that \((a_ix_i^b)_i\) lies in the kernel of \((x_i^r)_{0\le r\le w-2,\,1\le i\le w}\).
By \eqref{v123:eq-lagrange-moments} and the Vandermonde rank, this kernel is spanned by \((D'(x_i)^{-1})_i\), proving the equivalence.
Finally, the BCH bound gives \(d_H(C)\ge w\), while \(u\) has weight \(w\).
\end{proof}

When constructing a word from its support, one must also require every right-hand side of \eqref{v123:eq-bch-equality-weights} to lie in the base field \(k\).
For initial exponent \(b=0\) and all nonzero coefficients equal to one, the condition becomes
\[
D'(X)=\lambda\ne0.
\]
This criterion holds in every characteristic and forces \(p\mid w\): the first zero condition is \(u(1)=w=0\).
In odd characteristic it differs from the signed weights expressing equality of the power sums of two root sets, which lead to \eqref{v123:eq-wronskian-prefix}.

\begin{corollary}[Trace-prefix fibers and words attaining the bound]\label{v123:cor-paired-bch-word}
Let \(\gcd(q,N)=1\), let \(\zeta\) be a primitive \(N\)-th root of unity, and let \(C,D\) be distinct \(q\)-cyclotomic cosets modulo \(N\), with \(|C|=|D|=m\) and \(2m\le N\). Set
\[
f=\prod_{i\in C}(X-\zeta^i),\qquad
 g=\prod_{i\in D}(X-\zeta^i),\qquad
 U(Z)=\sum_{i\in C}Z^i-\sum_{i\in D}Z^i.
\]
The following are equivalent:
\begin{enumerate}
\item \(s_j(f)=s_j(g)\) for \(1\le j\le2m-2\);
\item \(f'g-fg'\) is a nonzero constant;
\item \(U\) belongs to the \(q\)-ary BCH code with consecutive zeros \(1,\zeta,\ldots,\zeta^{2m-2}\).
\end{enumerate}
When these conditions hold, the code has minimum distance \(2m\). If \(q=2\), the second condition is also equivalent to \(g=\mathcal P(f)\), and \(U^2\equiv U\pmod{Z^N-1}\).
\end{corollary}
\begin{proof}
The first two conditions are equivalent by \eqref{v123:eq-wronskian-prefix}.
Moreover, \(U(\zeta^j)=s_j(f)-s_j(g)\) and \(U(1)=m-m=0\), so the first condition is equivalent to the third.
The weight \(2m\) attains the BCH bound. In the binary case, Corollary~\ref{v123:cor-binary-derivative-pair} gives the pairing, and multiplication by two permutes \(C\cup D\), yielding \(U^2\equiv U\pmod{Z^N-1}\).
\end{proof}

Neither orbit is required to arise from primitive exponents. If \(N=q^m-1\), then both \(f\) and \(g\) are primitive if and only if the exponents in both orbits are coprime to \(N\).
A word attaining the BCH bound need not have support consisting of two equally long primitive orbits.

\begin{example}[Binary and ternary words attaining the bound]\label{v123:ex-bch-equality-small}
First take \(X^4+X+1\) from Example~\ref{v123:ex-f4-trace-pair}, with a root \(\alpha\) of order fifteen. Its roots over \(\F_2\) are \(\alpha,\alpha^2,\alpha^4,\alpha^8\), and the root polynomial has derivative one. Hence
\[
U(Z)=Z+Z^2+Z^4+Z^8
\]
vanishes at \(1,\alpha,\alpha^2\) and belongs to the binary BCH code of length fifteen generated by
\[
(Z+1)(Z^4+Z+1).
\]
The code has dimension ten; the BCH bound and \(\operatorname{wt}_H(U)=4\) give parameters \([15,10,4]\).
Here \(U^2\equiv U\pmod{Z^{15}-1}\). The support is one orbit of length four over \(\F_2\); over \(\F_4\), it splits into the two orbits of length two in the earlier example.

Next take a root \(\zeta\) of order twenty-six of \(f_2=X^3-X-2\) from Example~\ref{v123:ex-f3-trace-pencil}. Since \(\zeta^{13}=-1\), the element \(-\zeta=\zeta^{14}\) is a root of \(f_1=X^3-X-1\).
The corresponding ternary cyclotomic cosets modulo twenty-six are
\[
C_1=\{1,3,9\},\qquad C_{14}=\{14,16,22\}.
\]
The two polynomials share their first four power traces, so
\[
V(Z)=Z^{14}+Z^{16}+Z^{22}-Z-Z^3-Z^9
\]
vanishes at \(1,\zeta,\zeta^2,\zeta^3,\zeta^4\).
The ternary BCH code of length twenty-six defined by these roots has defining set
\[
\{0\}\cup\{1,3,9\}\cup\{2,6,18\}\cup\{4,12,10\}.
\]
Its dimension is therefore \(26-10=16\); the BCH bound and \(\operatorname{wt}_H(V)=6\) give parameters \([26,16,6]\).
Here \(f_1'f_2-f_1f_2'=1\), whereas \((f_1f_2)'=X^3-X\) is not constant. The signs of the weights are essential.
\end{example}

For the study of minimum-weight BCH codewords through power sums of locators and their root polynomials, see~\cite{AugotCharpinSendrier1992}.

\subsection{Recovering Error Locations from Syndromes}\label{v118:sec-bch-decoding}
The evaluation equations supplied by consecutive roots also recover errors. Write the received polynomial as \(r=c+e\), where \(c\) is a codeword and \(e\) is the error polynomial. If \(c(\alpha^j)=0\), then
\[
S_j=r(\alpha^j)=e(\alpha^j).
\]
\phantomsection\label{v118:def-syndrome}
These values, computed from the received word, are called \textbf{syndromes}. They eliminate the unknown codeword from the equations.

\begin{proposition}[Syndrome recurrences and error recovery]\label{v118:prop-syndrome-recurrence}
Let \(k=\F_q\), with \(p=\Char k\) and \(p\nmid N\), and let \(\alpha\) be a primitive \(N\)th root of unity in a finite extension \(E\). Suppose \(1\le2t<N\) and every codeword of \(C\subseteq k[X]/(X^N-1)\) vanishes at \(\alpha,\ldots,\alpha^{2t}\).
If \(r=c+e\), with \(c\in C\) and \(\operatorname{wt}_H(e)=w\le t\), write
\[
e(X)=\sum_{i=1}^w e_iX^{j_i},\qquad x_i=\alpha^{j_i},
\]
where \(e_i\ne0\) and the positions \(0\le j_i<N\) are distinct. Applied over \(E\) to \(S_1,\ldots,S_{2t}\), the Berlekamp--Massey algorithm returns connection length \(w\) and connection polynomial
\begin{equation}\label{v118:def-error-locator}
\Lambda(Z)=\prod_{i=1}^w(1-x_iZ).
\end{equation}
The polynomial \(\Lambda\) is called the \textbf{error-locator polynomial}. Its roots are \(\alpha^{-j_i}\). Once these positions are known, the first \(w\) syndromes uniquely determine \(e_1,\ldots,e_w\), and hence \(c=r-e\).
Moreover, the number of errors is determined directly by
\begin{equation}\label{v127:eq-syndrome-rank}
w=\operatorname{rank}(S_{i+j+1})_{0\le i,j<t}.
\end{equation}
For \(w>0\), the leading \(w\)-by-\(w\) matrix is invertible; solving the recurrence coefficient system with the first \(2w\) syndromes also recovers \(\Lambda\).
\end{proposition}
\begin{proof}
If \(w=0\), all syndromes vanish and both recoveries give \(\Lambda=1\). Suppose \(w>0\), and put
\[
h(T)=\prod_{i=1}^w(T-x_i),\qquad A=E[T]/(h),\qquad
\ell(u)=\sum_{i=1}^w e_ix_i u(x_i).
\]
The evaluation isomorphism \(A\simeq E^w\) makes \(\ell(uv)\) diagonal, with every diagonal entry \(e_ix_i\) nonzero.
Lemma~\ref{v127:lem-pairing-recovery} shows that \(s_a=\ell(\overline T^a)\) has minimal polynomial \(h\).
Its first \(2t\) terms are \(S_1,\ldots,S_{2t}\), and \(\deg h=w\le t\). Corollary~\ref{cor:bm-recovery} therefore recovers \(h\), with connection polynomial \(Z^wh(Z^{-1})=\Lambda(Z)\).
The \(w\)-by-\(w\) Hankel matrix is invertible, and the recurrence expresses the remaining columns of every larger matrix in terms of its first \(w\) columns. This proves~\eqref{v127:eq-syndrome-rank} and the coefficient recovery.
Finally, the coefficient matrix of \(S_a=\sum_i e_ix_i^a\), for \(1\le a\le w\), is a Vandermonde matrix times the invertible diagonal matrix \(\operatorname{diag}(x_i)\). Hence the error values are unique.
\end{proof}
The infinite sequence here is defined by the error polynomial; evaluations of the received polynomial at other powers need not agree with it. Recovery uses the first \(2t\) syndromes together with the bound on the number of errors.

\begin{example}[Correcting two errors in a binary BCH code of length fifteen]\label{v118:ex-bch-fifteen-decoding}
Take \(E=\F_2[a]\), where \(a^4+a+1=0\). This quartic has no linear factor and is not divisible by \(X^2+X+1\), so it is irreducible. Moreover, \(a^3\ne1\) and \(a^5=a^2+a\ne1\), so \(a\) has order fifteen.
The binary cyclotomic cosets \(\{1,2,4,8\}\) and \(\{3,6,12,9\}\) modulo fifteen give the generator polynomial
\[
\begin{aligned}
g(X)&=(X^4+X+1)(X^4+X^3+X^2+X+1)\\
&=X^8+X^7+X^6+X^4+1.
\end{aligned}
\]
The second factor is the fifth cyclotomic polynomial and has the order-five element \(a^3\) as a root. Thus \(C=(g)\subseteq\F_2[X]/(X^{15}-1)\) has consecutive zeros \(a,a^2,a^3,a^4\) and dimension seven. The BCH bound gives \(d_H(C)\ge5\), while \(g\) itself has weight five, so the code has parameters \([15,7,5]\).

Suppose at most two errors occur, and the received polynomial is
\[
r(X)=X^9+X^7+X^4+X^3+X^2+1.
\]
Successive reduction using \(a^4=a+1\) gives
\[
(S_1,S_2,S_3,S_4)=(a^{12},a^9,a^7,a^3).
\]
For example, \(S_1=a^9+a^7+a^4+a^3+a^2+1=a^{12}\); binary coefficients also give \(S_2=S_1^2\) and \(S_4=S_2^2\).
The same syndromes give
\[
H_2=\begin{pmatrix}a^{12}&a^9\\a^9&a^7\end{pmatrix},\qquad
\det H_2=a^7\ne0,
\]
so there are exactly two errors. Solving
\[
H_2\begin{pmatrix}b_0\\b_1\end{pmatrix}
=\begin{pmatrix}a^7\\a^3\end{pmatrix}
\]
gives \(b_0=a^{13}\) and \(b_1=a^{12}\), hence \(h(T)=T^2+a^{12}T+a^{13}\). The following successive recovery yields the same result.

Feed \(s_0=S_1,\ldots,s_3=S_4\) into Algorithm~\ref{alg:berlekamp-massey}. Initialize \(C=B=1,L=0,b=1,t=-1\). The state after each iteration is
\begin{center}
\begin{tabular}{@{}cccc@{}}
\toprule
\(j\)&Discrepancy \(\delta\)&\(L\)&\(C(Z)\)\\
\midrule
0&\(a^{12}\)&1&\(1+a^{12}Z\)\\
1&0&1&\(1+a^{12}Z\)\\
2&\(a^{10}\)&2&\(1+a^{12}Z+a^{13}Z^2\)\\
3&0&2&\(1+a^{12}Z+a^{13}Z^2\)\\
\bottomrule
\end{tabular}
\end{center}
At \(j=2\), we have \(\delta=S_3+a^{12}S_2=a^{10}\) and stored discrepancy \(b=a^{12}\), so the correction term is \((\delta/b)Z^2=a^{13}Z^2\).

The output is
\[
\Lambda(Z)=1+a^{12}Z+a^{13}Z^2=(1+a^3Z)(1+a^{10}Z).
\]
Evaluation at \(a^{-j}\), for \(0\le j<15\), gives zeros exactly at \(a^{12}=a^{-3}\) and \(a^5=a^{-10}\). Thus the error positions are three and ten. A nonzero binary error value must be one, so
\[
e=X^3+X^{10},\qquad
c=r+e=X^{10}+X^9+X^7+X^4+X^2+1.
\]
Finally, \(c=(X^2+1)g\) verifies that the recovered word belongs to the code. If another codeword were within distance two of \(r\), the two codewords would be within distance four of each other, contradicting \(d_H(C)=5\). The recovered word is therefore unique.
\end{example}

\section{Support Root Polynomials and the Distance of Binary Goppa Codes}\label{v82:subsec-goppa-support}
The BCH bound and its equality conditions express codeword supports through node polynomials and interpolation weights.
For binary Goppa codes, the logarithmic derivative of the support root polynomial further turns the parity-check condition into divisibility.
The following codes need not be cyclic; the argument still uses only root sets, quotient rings, and square-freeness.
Further connections between split auxiliary polynomials and exact distances are developed in~\cite{ChenChenLiuLaoGoppa2026}.

\leavevmode\phantomsection\label{v82:def-binary-goppa-code}\ignorespaces
Let \(E=\F_{2^m}\), let \(G\in E[X]\) be monic and square-free of degree \(r\ge1\), and let \(L\subseteq E\) be a nonempty set containing no root of \(G\).
The corresponding \textbf{binary Goppa code} is
\begin{equation}\label{v82:eq-goppa-definition}
 \Gamma(L,G)=\left\{(c_a)_{a\in L}\in\F_2^L:
       \sum_{a\in L}c_a(X-a)^{-1}=0\text{ in }E[X]/(G)\right\}.
\end{equation}
Since \(G(a)\ne0\), each \(X-a\) is a unit in the quotient ring, so the definition is meaningful.
This is a binary linear code; the target of its check map has dimension \(mr\) over \(\F_2\), so its dimension is at least \(|L|-mr\).
\leavevmode\phantomsection\label{v82:def-support-root-polynomial}\ignorespaces
For nonempty \(S\subseteq L\), define the \textbf{support root polynomial}
\[
 F_S(X)=\prod_{a\in S}(X-a).
\]
A binary vector is uniquely determined by its support. Write \(\mathbf1_S\) for the indicator vector, whose weight is \(|S|=\deg F_S\).

\begin{proposition}[A derivative test for supports and a distance bound]\label{v82:prop-goppa-support}
For nonempty \(S\subseteq L\),
\begin{equation}\label{v82:eq-goppa-support-test}
 \mathbf1_S\in\Gamma(L,G)
 \quad\Longleftrightarrow\quad G\mid F_S'
 \quad\Longleftrightarrow\quad G^2\mid F_S'.
\end{equation}
Consequently, every nonzero codeword of \(\Gamma(L,G)\) has weight at least \(2r+1\).
\end{proposition}
\begin{proof}
The product rule gives
\[
 \frac{F_S'}{F_S}=\sum_{a\in S}\frac1{X-a}.
\]
Since \(S\subseteq L\), we have \(\gcd(F_S,G)=1\). Multiplication by the unit \(F_S\) in the quotient modulo \(G\) proves the first equivalence.
In characteristic two, the derivative \(F_S'\) contains only even powers. Each coefficient has a unique square root in \(E\), so \(F_S'=H^2\) for some \(H\in E[X]\).
Square-freeness of \(G\) makes \(G\mid H^2\) equivalent to \(G\mid H\), hence to \(G^2\mid F_S'\).
Finally, \(F_S\) is nonconstant and square-free, so \(F_S'\ne0\). With \(w=|S|\),
\[
 2r\le\deg F_S'\le w-1,
\]
which proves the weight bound. The statement also covers the zero code without assigning it a separate minimum-distance convention.
\end{proof}

\begin{corollary}[Attaining the distance bound]\label{v82:cor-goppa-equality}
The code \(\Gamma(L,G)\) contains a word of weight \(2r+1\) if and only if some \(U\in E[X]\) with \(\deg U\le r\) satisfies
\begin{equation}\label{v82:eq-goppa-split-polynomial}
 XG(X)^2+U(X)^2\ \mid\ \prod_{a\in L}(X-a).
\end{equation}
Here \(U=0\) is allowed. When this condition holds, the code is nonzero and its minimum distance is exactly \(2r+1\).
\end{corollary}
\begin{proof}
If \(|S|=2r+1\), both \(F_S'\) and \(G^2\) are monic of degree \(2r\); divisibility therefore gives \(F_S'=G^2\).
Thus \(F_S-XG^2\) has zero derivative and degree at most \(2r\), so it equals \(U^2\) for some \(U\) of degree at most \(r\).
Conversely, condition~\eqref{v82:eq-goppa-split-polynomial} makes its monic left-hand side, of degree \(2r+1\), split into distinct roots in \(L\).
If \(S\) is this root set, its root polynomial has derivative \(G^2\), and Proposition~\ref{v82:prop-goppa-support} supplies a codeword of weight \(2r+1\).
\end{proof}

\begin{example}[A distance-nine code over the field of 64 elements]\label{v82:ex-goppa-sixtyfour}
Take \(E=\F_{64}\), \(G=X^4+X^3+1\), and
\[
 F=(X^7+1)(X^2+X+1).
\]
Over \(\F_2\), the polynomial \(G\) has no linear factor, and its remainder modulo the unique monic irreducible quadratic \(X^2+X+1\) is \(X\); hence it is irreducible.
Its roots have degree four over \(\F_2\) and therefore do not belong to \(E\). Also, \(G'=X^2\) and \(G(0)=1\), so \(G\) is square-free and we may take \(L=E\).
The roots of \(X^7+1\) are \(\F_8^\times\), and those of \(X^2+X+1\) are \(\F_4\setminus\F_2\).
Both sets lie in \(E\) and are disjoint, so \(F\) has nine distinct roots. Direct differentiation gives
\[
 F'=X^8+X^6+1=G^2.
\]
These nine roots therefore give a word of weight nine, while the general bound is also \(2\cdot4+1=9\).
Thus \(\Gamma(E,G)\) has length 64, dimension at least 40, and minimum distance exactly nine.
The distance is determined by a split root set and its derivative, without enumerating codewords.
\end{example}

\medskip
\noindent\textbf{Multiplicative cosets and square roots of derivatives.}
The preceding example starts with a split root polynomial and uses its derivative to verify attainment of the distance bound.
Conversely, multiplicative orbits can be used to choose a root set whose derivative has a prescribed form.
The construction uses only roots of unity from Chapter~\ref{chap:cyclotomic} and multiplicative orbits from Chapter~\ref{chap:root-symmetry}.

\begin{proposition}[A derivative construction from two multiplicative cosets]\label{v83:prop-two-mult-cosets}
Let \(E=\F_{2^m}\), let \(d\mid 2^m-1\), and choose distinct \(a,b\in(E^\times)^d\). Then
\[
 F(X)=X(X^d-a)(X^d-b)
\]
has \(2d+1\) distinct roots in \(E\) and satisfies
\begin{equation}\label{v83:eq-two-cosets-derivative}
 F'(X)=\bigl(X^d+\sqrt{ab}\bigr)^2.
\end{equation}
The monic polynomial \(G(X)=X^d+\sqrt{ab}\) is square-free.
For \(L=\{x\in E:G(x)\ne0\}\), the code \(\Gamma(L,G)\) has minimum distance \(2d+1\).
\end{proposition}
\begin{proof}
The integer \(d\) is odd. Each of \(X^d-a\) and \(X^d-b\) has \(d\) distinct nonzero roots, and their root sets are disjoint.
Expanding and differentiating gives
\[
 F=X^{2d+1}+(a+b)X^{d+1}+abX,
 \qquad F'=X^{2d}+ab,
\]
which proves~\eqref{v83:eq-two-cosets-derivative}.
Also \(G'=X^{d-1}\) and \(G(0)=\sqrt{ab}\ne0\), so \(\gcd(G,G')=1\).
No root of \(F\) is a root of \(G\), since \(F\) is square-free and \(F'=G^2\).
Corollary~\ref{v82:cor-goppa-equality} now gives the claimed distance.
\end{proof}

\begin{remark}[Nonzero cosets and the fixed point]\label{v83:rem-coset-fixed-point}
If \(d\ge5\) and \(S\subseteq E^\times\) is a union of full \(\mu_d\)-cosets, its root polynomial has the form \(F_S=P(X^d)\). Thus
\[
 \sqrt{F_S'(X)}=X^{(d-1)/2}\sqrt{P'(X^d)}.
\]
The right-hand side is divisible by \(X^2\) and cannot be square-free.
This obstruction does not extend to \(d=3\): for example, the square root of the derivative of \(F=X^3+1\) is \(X\).
Proposition~\ref{v83:prop-two-mult-cosets} also includes zero among its roots and thereby uses the fixed point of the multiplicative action, in accord with the distinction between special and free orbits in Chapter~\ref{chap:root-symmetry}.
\end{remark}

\medskip
\noindent\textbf{Parameters and splitting for a fixed derivative.}
The two-coset construction belongs to a more general parameter family.
Here \(Y\) records a multiplicative coset, whereas \(X\) records the roots within that coset.

\begin{proposition}[Fixed-derivative parameters and a splitting test]\label{v83:prop-derivative-fibers}
Let \(K=\F_Q\), where \(Q=2^m\), let \(d\mid Q-1\), and put \(h=(Q-1)/d\).
Let \(U\in K[Y]\) be monic and square-free, with \(\deg U=s\ge1\) and \(U(0)\ne0\).
For every \(V\) with \(\deg V<s\), define
\begin{equation}\label{v83:eq-derivative-fiber}
 P_V(Y)=U(Y)^2+YV(Y)^2,
 \qquad F_V(X)=XP_V(X^d),\qquad G(X)=U(X^d).
\end{equation}
Then \(G\) is monic and square-free, \(\deg G=ds\), and
\[
 F_V'(X)=G(X)^2,\qquad
 F_V\text{ is square-free}\ \Longleftrightarrow\ \gcd(U,V)=1.
\]
Conversely, if a monic polynomial \(F=XP(X^d)\) has degree \(2ds+1\) and \(F'=G^2\), there is a unique \(V\) of degree less than \(s\) such that \(F=F_V\).
The number of \(V\) satisfying this square-freeness condition is
\begin{equation}\label{v83:eq-euler-derivative-fibers}
 \varphi_Q(U)=Q^s\prod_{\pi\mid U}\bigl(1-Q^{-\deg\pi}\bigr),
\end{equation}
where the product is over the distinct monic irreducible factors of \(U\).
If \(H=(K^\times)^d\), then also
\begin{equation}\label{v83:eq-fiber-splitting-test}
 F_V\mid X^Q-X
 \quad\Longleftrightarrow\quad
 P_V\mid Y^h-1.
\end{equation}
\end{proposition}
\begin{proof}
Since \(d\mid Q-1\), \(d\) is odd. In characteristic two,
\[
 F_V'=\bigl(XU(X^d)^2\bigr)'+\bigl(X^{d+1}V(X^d)^2\bigr)'=U(X^d)^2.
\]
Conversely, if \(F=XP(X^d)\) has the stated degree and derivative, then \(P(Y)+YP'(Y)=U(Y)^2\). Comparison of coefficients shows that the even-degree terms of \(P\) are fixed by \(U^2\), while the odd-degree terms are free. Thus \(P=U^2+YV^2\) for a unique \(V\) of degree less than \(s\).
At zero, \(F_V'(0)=U(0)^2\ne0\). For nonzero \(x\), if \(G(x)=0\), then
\(F_V(x)=x^{d+1}V(x^d)^2\).
Hence \(F_V\) and its derivative have a common root precisely when \(U\) and \(V\) have a common root.
The assumptions \(U(0)\ne0\), odd \(d\), and square-freeness of \(U\) also imply that \(G\) is square-free.
The polynomials \(V\) of degree less than \(s\) represent all residue classes of \(K[Y]/(U)\), and the coprime classes are precisely its units.
Their number is therefore the polynomial Euler function of Section~\ref{sec:euler-polynomial}, proving~\eqref{v83:eq-euler-derivative-fibers}.
Finally, the image of the \(d\)-th power map on \(K^\times\) is the subgroup \(H\) of order \(h\).
For \(y\ne0\), the equation \(x^d=y\) has \(d\) distinct roots in \(K\) precisely when \(y\in H\).
Since \(P_V(0)=U(0)^2\ne0\), \(F_V\) splits into distinct linear factors over \(K\) if and only if the \(2s\) roots of \(P_V\) are distinct and lie in \(H\), which is equivalent to~\eqref{v83:eq-fiber-splitting-test}.
\end{proof}

\begin{example}[Algebraic candidates and split root sets with the same derivative]\label{v83:ex-sixtyfour-fibers}
Take \(K=\F_{64}\), \(d=3\), and \(U(Y)=Y+1\).
For \(v\in K\), equation~\eqref{v83:eq-derivative-fiber} becomes
\[
 F_v(X)=X^7+v^2X^4+X,\qquad G(X)=X^3+1.
\]
We have \(F_v'=G^2\), and \(F_v\) is square-free exactly when \(v\ne0\).
Thus there are \(\varphi_{64}(Y+1)=63\) algebraic candidates.
However, \(H=(K^\times)^3\) has only \(21\) elements. If \(F_v\) splits completely, the roots of
\(P_v(Y)=Y^2+v^2Y+1\) must be two distinct elements \(y,y^{-1}\in H\).
Because \(|H|=21\) is odd, \(y=y^{-1}\) in \(H\) occurs only at \(y=1\).
There are therefore exactly \((21-1)/2=10\) such unordered pairs, each determining a unique \(v\).
Thus exactly \(10\) of the \(63\) algebraic candidates split completely over \(K\).
Moreover, \(G\) has three roots in \(K\), so its full support \(L=\{x\in K:G(x)\ne0\}\) has \(61\) elements.
The ten distinct root sets yield at least ten words of weight seven, and \(\Gamma(L,G)\) has minimum distance seven.
These are counts of distinct locator root sets, not of different Goppa codes.
\end{example}

\medskip
\noindent\textbf{Adjoining roots while preserving conjugate orbits.}
Multiplicative cosets supply one kind of symmetry. To keep coefficients in a prescribed subfield, added roots should instead form complete relative Frobenius orbits.

\begin{proposition}[Adjoining a pair of relative Frobenius conjugates]\label{v83:prop-conjugate-pair-extension}
Let \(K=\F_Q\), where \(Q=2^s\), and let \(E=\F_{Q^{2t}}\), with \(s,t,r\ge1\).
Suppose that \(F\in K[X]\) is monic and square-free of degree \(2r+1\), and splits completely over \(E\).
If \(Q>4r\), there exists \(\beta\in\F_{Q^2}\setminus K\) such that
\[
 \widetilde F(X)=(X-\beta)(X-\beta^Q)F(X)
\]
splits completely and is square-free over \(E\), while the square root of \(\widetilde F'\) is also square-free.
This square root may be taken monic; after fixing the nonzero parameter \(b\) below, there are at least \(Q/2-2r\) distinct such square roots.
\end{proposition}
\begin{proof}
Since \(F\) is monic of degree \(2r+1\), there are \(G,u\in K[X]\) such that
\[
 F=XG^2+u^2,\qquad \deg G=r,\quad\deg u\le r.
\]
Square-freeness of \(F\) implies \(\gcd(G,u)=1\).
Fix \(b\in K^\times\). For \(a\in K\), put
\[
 M_a(X)=X^2+b^2X+a^2,\qquad
 \widetilde G_a(X)=(X+a)G(X)+bu(X).
\]
The quadratic trace criterion of Chapter~\ref{chap:trace} shows that \(M_a\) is irreducible over \(K\) precisely when
\(\Tr_{K/\F_2}(a/b^2)=1\); there are exactly \(Q/2\) such values of \(a\).
The roots of each such \(M_a\) form a relative Frobenius conjugate pair in \(\F_{Q^2}\).
Differentiation gives
\[
 (M_aF)'=\bigl((X+a)G+bu\bigr)^2=\widetilde G_a^2.
\]
To exclude repeated roots of \(\widetilde G_a\), put \(A=XG+bu\) and \(B=G\).
Since \(\gcd(A,B)=1\), a repeated root \(\xi\) of \(A+aB\) must satisfy
\[
 D(\xi)=0,\qquad D=A'B+AB'
 =G^2+b(u'G+uG').
\]
The polynomial \(D\) is monic of degree \(2r\) and satisfies \(D'=0\); it is therefore a square and has at most \(r\) distinct roots.
Each such root determines at most one parameter \(a\) for which \(A+aB\) has a repeated root.
We must also exclude parameters for which \(M_a\) and \(F\) have a common root.
The \(M_a\) are distinct irreducible quadratics, while \(\deg F=2r+1\), so there are at most \(r\) such parameters.
At least \(Q/2-2r>0\) parameters therefore make both \(M_aF\) and \(\widetilde G_a\) square-free.
Different values of \(a\) give different \(\widetilde G_a\), since their difference is a nonzero constant multiple of \(G\).
\end{proof}

This construction preserves symmetry under the relative Frobenius \(x\mapsto x^Q\).
When \(K\ne\F_2\), it does not by itself imply that the new root set is a union of full absolute Frobenius orbits.
This distinction is consistent with the descent principles of Chapter~\ref{chap:root-symmetry}.
Thus the distance criterion for binary Goppa codes can be used not only to test a given root set but also to construct one from admissible orbits.

The support root polynomial has \emph{coordinate positions} as roots and degree equal to the codeword weight. It differs from the polynomial whose coefficients form the codeword.
In the binary case, the coefficients are zero or one, so the root set already records the entire word. This is why the logarithmic derivative directly expresses the check condition.
We next examine minimal root completion within a fixed base field, before returning to coefficient representations of cyclic codes.

\par\Needspace{12\baselineskip}
\section{Constant-Derivative Completion of a Prescribed Root Set}\label{v116:sec-completion}
Parameterizing a fixed derivative also addresses another question: how few base-field elements must be added while retaining prescribed roots? Proposition~\ref{v123:prop-bch-equality} shows that root polynomials with nonzero constant derivative also give vanishing consecutive power sums and words attaining the BCH bound. Here all roots must additionally lie in the prescribed base field.
For \(A\subseteq k=\F_q\) with \(|A|\ge2\), write \(R_A=\prod_{a\in A}(X-a)\) and define the \textbf{minimum constant-derivative completion size}
\begin{equation}\label{v116:def-completion}
c_q(A)=\min\{|B|:A\subseteq B\subseteq k,\ R_B'\in k^\times\}.
\end{equation}
The full field has root polynomial \(X^q-X\), so the minimum exists. If the affine hull of \(A\) over \(\F_p\) has dimension \(r\), its subspace polynomial gives \(c_q(A)\le p^r\).
An invertible affine change preserves this minimum, since \(R_{uB+v}(X)=u^{|B|}R_B((X-v)/u)\).

\begin{proposition}[Interpolation candidates and a degree gap]\label{v116:prop-completion-gap}
Let \(q=p^e\), \(s=|A|\ge2\), and let \(H_A\) be the unique polynomial of degree less than \(s\) with \(H_A(a)=a^{q/p}\).
All polynomials satisfying \(F'=1\) and \(F|_A=0\) are precisely
\begin{equation}\label{v116:eq-completion-family}
F_Q=X-(H_A+R_AQ)^p,\qquad Q\in k[X].
\end{equation}
Writing \(\rho=\deg H_A\), we have \(\deg F_0=p\rho\), whereas \(\deg F_Q=p(s+\deg Q)\ge ps\) for \(Q\ne0\). Consequently,
\[
F_0\mid X^q-X\ \Longrightarrow\ c_q(A)=p\rho;
\qquad F_0\nmid X^q-X\ \Longrightarrow\ c_q(A)\ge ps.
\]
\end{proposition}
\begin{proof}
Perfectness and formal differentiation show that \(F'=1\) is equivalent to \(F=X-G^p\), with unique \(G\in k[X]\).
The root condition is equivalent to \(G(a)=a^{q/p}\), or \(G\equiv H_A\pmod{R_A}\), giving the parameterization.
The polynomial \(H_A\) is nonconstant because inverse Frobenius has different values at distinct points; thus \(\rho\ge1\).
Since \(\deg H_A<s\), there is no leading-term cancellation when \(Q\ne0\), proving both degree formulas.
The root polynomial of any completion, divided by its nonzero constant derivative, belongs to this family. Conversely, if \(F_Q\mid X^q-X\), its distinct roots form a completion of size \(\deg F_Q\). Compare the degrees.
\end{proof}

The derivative and prescribed-root conditions have now been solved. What remains is splitting over the fixed base field, which can be tested by computing \(X^q\bmod F_Q\) without finding each root.

\begin{corollary}[A size restriction for proper subsets]\label{v116:cor-completion-size}
If \(B\subsetneq\F_q\), \(|B|\ge2\), and \(R_B'\in\F_q^\times\), then \(p\mid|B|\) and \(|B|\le q/p\). In particular, \(|A|>q/p\) implies \(c_q(A)=q\).
\end{corollary}
\begin{proof}
Normalize to write \(R_B/R_B'=X-G^p\). Since the degree exceeds one, \(|B|=p\deg G\).
As \(|B|<q\), we have \(\deg G<q/p\), and the nonzero polynomial \(G-X^{q/p}\) vanishes on \(B\). The root bound gives \(|B|\le q/p\).
\end{proof}

\begin{proposition}[Completing four points in characteristic two]\label{v116:prop-four-completion}
Let \(q\) be a power of two. Affinely normalize four distinct points to \(A=\{0,1,a,b\}\), and put
\[
e_1=1+a+b,\qquad e_2=a+b+ab,\qquad e_3=ab.
\]
If \(e_1\ne0\), put \(u=e_2+e_3/e_1\). Then
\begin{equation}\label{v116:eq-four-completion}
c_q(A)=
\begin{cases}
4,&e_1=0,\\
6,&e_1\ne0,\ \Tr_{q/2}(u/e_1^2)=0,\\
8,&e_1\ne0,\ \Tr_{q/2}(u/e_1^2)=1.
\end{cases}
\end{equation}
\end{proposition}
\begin{proof}
We have \(R_A=X^4+e_1X^3+e_2X^2+e_3X\). If \(e_1=0\), its derivative is \(e_3\ne0\); also \(b=1+a\), so the four points already form an affine plane over \(\F_2\).
If \(e_1\ne0\), the four-point set itself does not have constant derivative, and every completion has even size.
The root polynomial of a six-point completion must be \(R_A(X)(X^2+vX+w)\). Comparing the coefficients of \(X^5,X^3\), the constant-derivative condition forces
\[
v=e_1,\qquad w=e_2+e_3/e_1=u.
\]
This unique candidate \(S=R_A(X)(X^2+e_1X+u)\) satisfies \(S'=e_3u\), and
\[
e_1u=e_1e_2+e_3=(a+b)(a+1)(b+1)\ne0.
\]
Thus \(S\) has no repeated root, and its quadratic factor has no root in common with \(R_A\). That factor splits if and only if \(\Tr_{q/2}(u/e_1^2)=0\).
If the trace is one, six points are impossible; the affine hull of the four points has at most eight points and supplies a completion. The minimum is therefore eight.
\end{proof}

\begin{example}[Six- and eight-point completions with the same affine hull]\label{v116:ex-completion-six-eight}
Let \(k=\F_2(\theta)\), where \(\theta^4+\theta+1=0\), and take
\[
A_6=\{0,1,\theta,\theta^2\},\qquad
A_8=\{0,1,\theta,\theta^2+\theta\}.
\]
Both affine hulls are \(\langle1,\theta,\theta^2\rangle_{\F_2}\).
For \(A_6\), the two roots of the quadratic factor above are \(\theta^3+\theta\) and \(\theta^3+\theta^2+1\). The resulting six-point set \(B_6\) satisfies
\[
R_{B_6}=X^6+(\theta+1)X^4+(\theta^2+1)X^2+(\theta^2+\theta+1)X.
\]
For \(A_8\), we have \(e_1=1+\theta^2\), \(u=1+\theta^2\), and \(u/e_1^2=\theta^3+\theta+1\), of absolute trace one.
Hence \(c_{16}(A_6)=6\) and \(c_{16}(A_8)=8\). The common affine hull gives the same upper bound; splitting of the quadratic factor decides whether it can be improved strictly.
\end{example}

Continue to assume \(\Char k=2\). For \(|A|=5\), we have \(3\le\deg H_A\le4\). If the degree is three, the degree-six candidate already has five base-field roots, and the remaining linear factor splits automatically. If the degree is four, first test the residual cubic in the degree-eight candidate.
If it does not split, nonzero constant \(Q\) gives all degree-ten candidates. When can they all fail to split? On increasing \(\deg Q\), can the residual factor types determine the minimum completion size in advance?
The affine hull has at most sixteen points, so no degree beyond sixteen is needed; every layer must still obey Corollary~\ref{v116:cor-completion-size}.
Corollary~\ref{v116:cor-unweighted-self-dual} will also interpret constant-derivative completion as constructing a self-dual evaluation code while retaining prescribed evaluation points.

\section{Support Blocks and Affine Coefficient Constraints}
\label{v100:sec-code-support}
\subsection{Support Periods and Interleaved Direct Sums}
\phantomsection\label{v44:subsec-interleaving}
Here \(N\) may be any positive integer, including the case \(p\mid N\).
\begin{proposition}[Support periods and interleaved decompositions of cyclic codes]\leavevmode\phantomsection\label{v44:prop-interleaving}
Let \(f\in k[X]\) be monic and nonconstant, with \(f\mid X^N-1\), and put
\(s=\operatorname{sp}(f)\), \(f(X)=f_{\mathrm{sp}}(X^s)\).
Then \(s\mid N\); putting \(M=N/s\) gives \(f_{\mathrm{sp}}\mid Y^M-1\). Write
\[
\begin{aligned}
A_N&=k[X]/(X^N-1),& A_M&=k[Y]/(Y^M-1),\\
C&=(f)\subseteq A_N,& C_{\mathrm{sp}}&=(f_{\mathrm{sp}})\subseteq A_M.
\end{aligned}
\]
The map
\[
\Phi:A_M^{\oplus s}\longrightarrow A_N,\qquad
(a_0,\ldots,a_{s-1})\longmapsto\sum_{j=0}^{s-1}X^j a_j(X^s)
\]
is a \(k\)-linear isomorphism and restricts to
\[
\boxed{\quad C_{\mathrm{sp}}^{\oplus s}\simeq C.\quad}
\]
Computing weights from the coefficients of canonical remainders gives
\[
\operatorname{wt}_H\bigl(\Phi(a_0,\ldots,a_{s-1})\bigr)
=\sum_{j=0}^{s-1}\operatorname{wt}_H(a_j),
\qquad \dim_k C=s\dim_k C_{\mathrm{sp}}.
\]
If \(C\ne\{0\}\), then \(d_H(C)=d_H(C_{\mathrm{sp}})\).
\end{proposition}
\begin{proof}
The polynomial ring decomposes by exponent residue classes as
\[
k[X]=\bigoplus_{j=0}^{s-1}X^j k[X^s].
\]
First write \(N=sM_0+r\) and \(0\le r<s\), and choose \(h\) such that \(X^N-1=f_{\mathrm{sp}}(X^s)h(X)\).
If \(r>0\), decompose \(h\) as above and compare the component of residue class zero modulo \(s\). This gives
\(-1=f_{\mathrm{sp}}(X^s)h_0(X^s)\), contradicting the assumption that \(f_{\mathrm{sp}}\) is nonconstant.
Thus \(r=0\), that is, \(s\mid N\) and \(M_0=M\). Comparing the other components gives \(h_j=0\) for \(j>0\),
and hence \(Y^M-1=f_{\mathrm{sp}}(Y)h_0(Y)\).

Substitution of \(Y=X^s\) sends \(Y^M-1\) to \(X^N-1\), so \(\Phi\) is well-defined on the quotient spaces.
Each canonical representative of degree less than \(N\) has a unique expression
\(\sum_{j=0}^{s-1}X^j a_j(X^s)\), where \(\deg a_j<M\). Thus \(\Phi\) is a linear isomorphism.
Different components occupy disjoint coordinate positions, so their weights add.

Since \(f\mid X^N-1\), a canonical representative \(c\) belongs to \(C\) if and only if \(f\mid c\).
If \(c=f h\), decompose \(h\) by exponent residue classes and compare the components to obtain
\(a_j(Y)=f_{\mathrm{sp}}(Y)h_j(Y)\); conversely, these divisibilities directly imply \(f\mid c\).
Thus \(\Phi(C_{\mathrm{sp}}^{\oplus s})=C\), and dimensions add.
When \(C\ne\{0\}\), every nonzero element has a nonzero component, so its weight is at least \(d_H(C_{\mathrm{sp}})\);
taking one nonzero component of minimum weight and setting all others to zero attains the bound.
\end{proof}

\leavevmode\phantomsection\label{v44:def-support-decomposition}\ignorespaces
The code \(C_{\mathrm{sp}}\) is called the \textbf{support-primitive core} of \(C\), and the decomposition above is called the \textbf{support-primitive decomposition}; see \cite{ZWSupportPrimitive}.
It is a linear isomorphism given by coordinate interleaving. When \(s>1\), the map \(\Phi\) is not a ring isomorphism:
in the product ring the product of \((1,0,\ldots)\) and \((0,1,0,\ldots)\) is zero, while the product of their images \(1\) and \(X\) is nonzero.
This decomposition therefore differs from the ring-product decomposition of the Chinese remainder theorem.
One cyclic right shift successively permutes the exponent-class components and induces a cyclic right shift of the shorter code when the last component returns to the first;
each component is invariant under \(s\) right shifts and is generally not itself a cyclic subcode of length \(N\).

\begin{corollary}[Complete orbit fibers and interleaved irreducible codes]\label{v142:cor-irreducible-interleaving}
Use $n,d,s,m,r,\xi,\eta,h,h_0$ from Proposition~\ref{v142:prop-full-fibers}, and assume $m=sr$. Put
\[
 g_0(Y)=\frac{Y^d-1}{h_0(Y)},\qquad C_0=(g_0)\subseteq k[Y]/(Y^d-1).
\]
The length-$n$ irreducible cyclic code $C_\xi$ with check polynomial $h$ satisfies
\begin{equation}\label{v142:eq-orbit-code-interleaving}
 g(X)=g_0(X^s),\qquad
 C_\xi\cong C_0^{\oplus s},\qquad
 \dim_k C_\xi=sr,\qquad d_H(C_\xi)=d_H(C_0).
\end{equation}
The isomorphism groups coordinates by their indices modulo $s$. If $s>1$, then
\[
 d_H(C_\xi)\le d-r+1<s(d-r)+1=n-m+1,
\]
so $C_\xi$ is not MDS.
\end{corollary}
\begin{proof}
The identity $h(X)=h_0(X^s)$ gives
\[
 g(X)=\frac{X^{ds}-1}{h_0(X^s)}=g_0(X^s).
\]
The decomposition by exponent classes in the proof of Proposition~\ref{v44:prop-interleaving} shows that $C_\xi$ consists of independent choices from $X^iC_0(X^s)$, $0\le i<s$. This does not require $s$ to be the maximum support period of $g$. The blocks have disjoint supports, giving the dimension and distance identities. Finally, $1\le r\le\varphi(d)<d$, so the Singleton bound for $C_0$ gives the strict inequality.
\end{proof}

Relative trace realizes the same decomposition in trace parameters. Set $E=k(\xi)$ and $F=k(\eta)$. For $0\le i<s$ and $0\le j<d$,
\begin{equation}\label{v142:eq-relative-trace-blocks}
 \Tr_{E/k}(a\xi^{-(i+sj)})
 =\Tr_{F/k}\bigl(\Tr_{E/F}(a\xi^{-i})\eta^{-j}\bigr).
\end{equation}
Since $1,\xi^{-1},\ldots,\xi^{-(s-1)}$ is an $F$-basis of $E$, nondegeneracy of the trace pairing makes
\[
 E\longrightarrow F^s,\qquad a\longmapsto
 \bigl(\Tr_{E/F}(a\xi^{-i})\bigr)_{i=0}^{s-1}
\]
an $F$-linear isomorphism. The short-code trace parameters can therefore be prescribed independently. Without $m=sr$, identity~\eqref{v142:eq-relative-trace-blocks} still holds. The first $m/r$ inverse powers already form an $F$-basis of $E$, so the parameter map remains injective, but its image has $F$-dimension only $m/r<s$. The blocks therefore satisfy linear constraints.
A coordinate direct sum does not contradict cyclic irreducibility: one cyclic right shift permutes the blocks, and an individual block is not a length-$n$ cyclic subcode.

\begin{example}[Five single parity-check blocks in twenty-five coordinates]\label{v142:ex-twenty-five-code}
Continuing Example~\ref{v142:ex-twenty-five-roots}, the check polynomial $h=\Phi_{25}$ gives
\[
 g_0(Y)=Y-1,\qquad g(X)=X^5-1.
\]
Thus $C_0$ is the $[5,4,2]$ single parity-check code over $\F_{23}$, and the words of $C_\xi$ are precisely those satisfying
\[
 \sum_{j=0}^4 c_{i+5j}=0\qquad(0\le i<5).
\]
Hence $C_\xi$ has parameters $[25,20,2]$. The generator $X^5-1$ itself is a word of weight two, while no word of weight one satisfies these checks. Its dual consists of five disjoint repetition blocks of length five and has parameters $[25,5,5]$.
Here $h$ is irreducible, yet the coordinate space splits into five blocks. Complete fibers of the root orbit become polynomial composition in one representation and independent short-code parameters in the other.
\end{example}

\begin{example}[Two shorter codes in a repeated-root cyclic code]\label{v44:ex-interleaving-six}
Over \(\F_2\), take \(N=6\) and \(f=X^2+1=(X+1)^2\). Since \(f\mid X^6+1\), it generates a cyclic code \(C\).
Here \(s=2\), \(f_{\mathrm{sp}}=Y+1\), and \(M=3\). The shorter code \(C_{\mathrm{sp}}=(Y+1)\subseteq\F_2[Y]/(Y^3+1)\) is
\[
C_{\mathrm{sp}}=\{(0,0,0),(1,1,0),(1,0,1),(0,1,1)\}.
\]
It consists of length-three vectors of even weight and has parameters \([3,2,2]\). Interleaving two copies in the even and odd positions gives
\[
C=\{(c_0,\ldots,c_5):c_0+c_2+c_4=0,\ c_1+c_3+c_5=0\}.
\]
Thus \(C\) has \(16\) codewords and parameters \([6,4,2]\). For example, \((1,0,1,0,0,0)\) has weight two.
The root \(1\) has multiplicity two, the common \(p\)-power is two, and the multiplicative stabilizer is trivial; the support period two comes here from root multiplicity.
\end{example}

\subsection{Affine Symmetries and Code Subspaces}
\label{v54:subsec-affine-test-codes}

A support period gives a direct sum by coordinate interleaving. Affine symmetry also involves a center and translation directions,
which give rise to code subspaces through the semi-invariant spaces of Section~\ref{fas:subsec-semi-invariants}.
We first test divisibility in the polynomial space, then return to the original coefficient coordinates to compute Hamming weights.

Continue with \(k=\F_q\), and fix \(C=(f)\subseteq k[X]/(X^N-\lambda)\), where \(\lambda\ne0\) and \(0<\tau=\deg f<N\).
Equation~\eqref{v54:eq-code-polynomial-realization} specifies the polynomial representatives and their degree range.

\begin{proposition}[A criterion for affine substitutions to preserve a code space]\label{v54:prop-affine-code-stability}
Let \(K/k\) be an extension field, with \(a\in K^\times\) and \(b\in K\). In the space of polynomials of degree less than \(N\),
the substitution \(g(X)\mapsto g(aX+b)\) preserves \(C\otimes_k K\) if and only if
\[
f(aX+b)=a^\tau f(X).
\]
Its affine stabilizer is therefore exactly the root-multiplicity-preserving affine stabilizer of its generator polynomial.
\end{proposition}
\begin{proof}
If the identity holds, \(fq\) is sent to \(a^\tau f(X)q(aX+b)\), and the multiplier still satisfies the original degree bound;
the inverse substitution gives equality of the code spaces. Conversely, \(f(aX+b)\) belongs to the code space and has degree \(\tau\),
so by~\eqref{v54:eq-code-polynomial-realization} it is a constant multiple of \(f\); comparing leading terms gives the constant \(a^\tau\).
\end{proof}
The action here is in the specified polynomial coordinates. If \(b\ne0\), then \(X\) is sent to \(aX+b\),
and its Hamming weight changes from one to two. Thus a substitution generally does not preserve weight.
Put \(s=\operatorname{sp}(f)\) and write \(f(X)=f_{\mathrm{sp}}(X^s)\). Replacing \(-1\) by \(-\lambda\) in the proof of Proposition~\ref{v44:prop-interleaving} gives
\(s\mid N\), \(f_{\mathrm{sp}}\mid Y^{N/s}-\lambda\), and the same coordinate interleaving with additive Hamming weights.
The resulting blocks also determine every decomposition along disjoint sets of coordinates.

\begin{definition}[Coordinate direct sums]\label{v67:def-coordinate-direct-sum}
For \(B\subseteq\{0,\ldots,N-1\}\), let \(k^B\subseteq k^N\) denote the subspace of vectors supported on \(B\).
A partition \(B_1,\ldots,B_m\) of the coordinate set gives a \textbf{coordinate direct-sum decomposition} of \(C\) if
\[
C=\bigoplus_{i=1}^m(C\cap k^{B_i}).
\]
A nonzero summand is \textbf{coordinate-indecomposable} if it cannot be split further into two nonzero subcodes supported on disjoint coordinate sets.
\end{definition}

\begin{proposition}[The indecomposable coordinate blocks]\label{v67:prop-coordinate-blocks}
Let \(C=(f)\subseteq k[X]/(X^N-\lambda)\), where \(f\) is monic, \(f\mid X^N-\lambda\), \(\lambda\ne0\), and \(0<\deg f<N\).
Put \(s=\operatorname{sp}(f)\). The sets
\[
R_j=\{i:0\le i<N,\ i\equiv j\pmod s\},\qquad 0\le j<s,
\]
are precisely the indecomposable coordinate blocks of \(C\).
Every coordinate direct-sum decomposition of \(C\) is obtained by grouping these blocks.
\end{proposition}
\begin{proof}
A nonzero codeword is called support-minimal if no nonzero codeword has strictly smaller support contained in its support.
Since \(f\mid X^N-\lambda\), we have \(f(0)\ne0\).
If a nonzero codeword \(g\) is supported inside \(\operatorname{Supp}_X(f)\), its standard representative has degree at most \(\deg f\).
Divisibility by \(f\) then makes \(g\) a nonzero scalar multiple of \(f\).
Thus \(f\) is support-minimal. The constacyclic shift is invertible and permutes the coordinates while multiplying one coordinate by a nonzero scalar,
so every constacyclic shift of \(f\) is also support-minimal.

Identify the coordinate set with \(\mathbb Z/N\mathbb Z\), and join two coordinates whenever they occur in the support of the same shifted generator.
These supports are all translates of \(\operatorname{Supp}_X(f)\).
Because \(0\in\operatorname{Supp}_X(f)\), the connected component of zero is the subgroup generated by the positive support indices:
each such index gives edges from \(a\) to its translate by that index, and every edge difference belongs to the subgroup they generate.
Since their greatest common divisor is \(s\) and \(s\mid N\), the connected components are exactly the \(R_j\).

The shifted generators span \(C\), and each lies within one \(R_j\), so these blocks give a coordinate direct sum.
In any coordinate direct sum, a support-minimal word must lie in one summand;
otherwise projection onto a nonzero summand would give a nonzero codeword with strictly smaller support.
Consequently no support of a shifted generator, and hence no connected component \(R_j\), can be split between two summands.
This proves both assertions.
\end{proof}

\begin{corollary}[A distance bound from the compressed code]\label{v67:cor-core-distance}
Under the hypotheses of Proposition~\ref{v67:prop-coordinate-blocks}, put \(\tau=\deg f\) and
\[
D=(f_{\mathrm{sp}})\subseteq k[Y]/(Y^{N/s}-\lambda).
\]
If \(C\) has parameters \([N,k_C,d_H(C)]\), then \(D\) has parameters
\[
[N/s,k_C/s,d_H(C)],\qquad k_C=N-\tau,
\]
and
\[
d_H(C)\le \frac{\tau}{s}+1.
\]
\end{corollary}
\begin{proof}
Interleaving gives \(C\simeq D^{\oplus s}\), so dimensions add and the two minimum distances agree.
For a nonzero linear code of length \(M\) and dimension \(K\), imposing zero values in any \(K-1\) coordinate positions leaves a nonzero word,
whose weight is at most \(M-K+1\). Applying this Singleton bound to \(D\) gives the stated inequality.
\end{proof}
For the \([6,4,2]\) code in Example~\ref{v44:ex-interleaving-six}, the even and odd positions are therefore its two indecomposable coordinate blocks.
The bound obtained from the \([3,2,2]\) compressed code is attained.

\begin{definition}[The effective support of a subcode]\label{v54:def-effective-support}
For a linear subspace \(D\subseteq C\), write
\[
\operatorname{Supp}_X(D)=\bigcup_{g\in D}\operatorname{Supp}_X(g),
\]
and call it the \textbf{effective support} of \(D\).
It is the union of the coefficient supports of any basis; no single codeword is required to be nonzero at every position in it.
\end{definition}

\begin{proposition}[The divisibility matrix and effective support]\label{v54:prop-test-subcode}
Let \(W\subseteq k[X]_{<N}\) have basis \(w_1,\ldots,w_l\), with \(D=C\cap W\).
Form the matrix \(M\) whose columns are the coefficient vectors of the remainders of \(w_i\) modulo \(f\). Put \(\rho=\rank_k M\) and write
\(r_j=([X^j]w_1,\ldots,[X^j]w_l)\). Then
\begin{equation}\label{v54:eq-effective-support}
\dim_kD=l-\rho,\qquad
\operatorname{Supp}_X(D)=\{0\le j<N:r_j\notin\operatorname{row}_k(M)\}.
\end{equation}
If \(l>\rho\), then
\begin{equation}\label{v54:eq-test-distance}
d_H(C)\le d_H(D)\le |\operatorname{Supp}_X(D)|-(l-\rho)+1.
\end{equation}
\end{proposition}
\begin{proof}
Basis coordinates identify \(D\) with \(\ker M\), so the dimension formula follows from rank--nullity.
The \(j\)th coefficient vanishes identically on \(D\) if and only if the linear functional \(r_j\) vanishes on \(\ker M\),
equivalently, \(r_j\) belongs to the row space of \(M\). This proves the support formula.
Put \(d=l-\rho>0\). Requiring any \(d-1\) coordinates in the effective support to vanish on \(D\)
imposes at most \(d-1\) linear conditions, leaving a nonzero solution of weight at most \(|\operatorname{Supp}_X(D)|-d+1\).
Finally, \(D\subseteq C\) gives the first inequality.
\end{proof}

\subsection{Reducing Divisibility by the Affine Normal Form}
\label{v54:subsec-outer-divisibility}

Let \(f\) have at least two distinct roots, and retain the data of Theorem~\ref{ffs:thm-normal-form}.
When \(t>1\), write \(f=(U^wH(U^t))^e\), \(U=L_V-\delta\), and \(H(0)\ne0\).
When \(t=1\), take \(U=L_V\) and write \(F(Y)\) uniquely as \(Y^wH(Y)\), \(H(0)\ne0\), giving the same expression.
Write
\[
A(X)=U(X)^e,\qquad P(X)=A(X)^t,\qquad
B(Y)=H^{[e]}(Y),
\]
where \(H^{[e]}(Y)=\sum_i h_i^eY^i\) means raising the coefficients to the \(e\)th power.
By the power identity in characteristic \(p\),
\begin{equation}\label{v54:eq-powered-normal-form}
f=A^wB(P),\qquad B(0)\ne0,\qquad \deg A=eQ.
\end{equation}
Here \(H(Y)^e=H^{[e]}(Y^e)\), and both coefficient powers and exponent powers must be retained.

\begin{theorem}[Outer divisibility retaining the special orbit]\label{v54:thm-outer-divisibility}
For \(0\le r<t\), put
\[
u_r=\max\left\{0,\left\lceil\frac{w-r}{t}\right\rceil\right\},
\qquad B_r(Y)=Y^{u_r}B(Y).
\]
Then every \(R\in k[Y]\) satisfies
\begin{equation}\label{v54:eq-outer-divisibility}
f\mid A^rR(P)\quad\Longleftrightarrow\quad B_r\mid R.
\end{equation}
Suppose that \(J\subseteq\Z_{\ge0}\) is finite and nonempty, and \(eQ(r+t\max J)<N\). Take
\[
W_{r,J}=\operatorname{span}_k\{A^{r+tj}:j\in J\}.
\]
Let \(\rho_{r,J}\) be the rank of the remainder map
\[
\operatorname{span}_k\{Y^j:j\in J\}\longrightarrow k[Y]/(B_r).
\]
Then \(\dim_k(C\cap W_{r,J})=|J|-\rho_{r,J}\).
For \(B_r=1\), the target is the zero space and the rank is zero.
In particular, if \(J=\{0,\ldots,L\}\), this dimension is
\(\max\{0,L+1-\deg B_r\}\), and when nonzero the space has basis
\[
A^rB_r(P)P^j\qquad(0\le j\le L-\deg B_r).
\]
\end{theorem}
\begin{proof}
The zero polynomial case is immediate. For nonzero \(R\), write \(R=Y^uR_0\) and \(R_0(0)\ne0\).
At every root of \(A\), both \(B(P)\) and \(R_0(P)\) are nonzero,
so the multiplicity required by \(A^w\) is equivalent to \(r+tu\ge w\), that is, \(Y^{u_r}\mid R\).
Since \(A\) and \(B(P)\) are coprime, the remaining condition is \(B(P)\mid R(P)\).
Divide \(R\) by the monic polynomial \(B\), writing \(R=BS+T\) and \(\deg T<\deg B\).
If \(T\ne0\), then \(\deg T(P)<\deg B(P)\), so \(B(P)\) cannot divide \(T(P)\).
Thus the remaining condition is equivalent to \(B\mid R\); this also holds when \(B=1\).
Since \(B(0)\ne0\), the two conditions combine to give \(B_r\mid R\).

The powers \(A^{r+tj}\) have distinct degrees and are therefore linearly independent, so the parameter space has dimension \(|J|\).
Formula~\eqref{v54:eq-outer-divisibility} identifies the intersection with the code as the kernel of the outer remainder map, proving the dimension formula.
For a consecutive exponent set, \(R\) is exactly \(B_rS\) with \(\deg S\le L-\deg B_r\);
taking \(S=1,Y,\ldots,Y^{L-\deg B_r}\) gives the stated basis.
\end{proof}

When \(t>1\), the identity \(U(aX+b)=aU(X)\) makes these powers semi-invariants with scaling factor \(a^{er}\);
when \(t=1\), they are translation invariants.
The factor \(Y^{u_r}\) retains the root multiplicity of the special orbit on \(U=0\) and cannot be omitted when \(w>r\).
The polynomial defining the outer quotient ring may have repeated roots; the proof does not require \(B_r\) to be separable.
The effective support of the subcode must still be computed in \(X\)-coordinates using Proposition~\ref{v54:prop-test-subcode};
the number of terms of the outer polynomial cannot be read directly as the codeword weight.

\begin{corollary}[A recurrence matrix for outer divisibility]\label{v67:cor-hankel-divisibility}
Retain \(r,J\), and \(B_r\) from Theorem~\ref{v54:thm-outer-divisibility}.
Put
\[
H(Y)=B_r(Y)=Y^d+\sum_{a=0}^{d-1}h_aY^a,
\]
and first assume \(d\ge1\). In \(A_H=k[Y]/(H)\), write \(y=\overline Y\), and let \(\ell_H\) extract the coefficient of \(Y^{d-1}\) from the standard remainder.
The sequence \(s_j=\ell_H(y^j)\) is given by
\[
s_0=\cdots=s_{d-2}=0,\qquad s_{d-1}=1,\qquad
s_{j+d}+\sum_{a=0}^{d-1}h_as_{j+a}=0\quad(j\ge0),
\]
where the initial zero string is empty when \(d=1\).
The matrix
\begin{equation}\label{v67:eq-hankel-divisibility}
\mathcal H_J=(s_{i+j})_{\substack{0\le i<d\\j\in J}}
\end{equation}
has the same kernel as the outer remainder matrix. Consequently,
\[
\rho_{r,J}=\rank_k\mathcal H_J,\qquad
\dim_k(C\cap W_{r,J})=|J|-\rank_k\mathcal H_J.
\]
No consecutiveness assumption on \(J\) or separability assumption on \(H\) is required.
If \(H=1\), the remainder map has zero target; its kernel is all of \(k^J\), and its rank is zero.
\end{corollary}
\begin{proof}
The initial values follow from the definition of \(\ell_H\), and the recurrence follows from \(H(y)=0\).
By Theorem~\ref{thm:algebra-trace-coefficient}, the coefficient pairing on \(A_H\) is nondegenerate.
Since \(1,y,\ldots,y^{d-1}\) is a basis, the map
\[
A_H\longrightarrow k^d,\qquad
u\longmapsto\bigl(\ell_H(u),\ell_H(yu),\ldots,\ell_H(y^{d-1}u)\bigr)
\]
is an isomorphism. Composing the outer remainder map with this isomorphism gives the matrix~\eqref{v67:eq-hankel-divisibility},
so its kernel and rank are unchanged. The dimension formula follows from Theorem~\ref{v54:thm-outer-divisibility}.
This is the sequence corresponding to \(1\in A_H\) in Theorem~\ref{thm:recurrence-algebra}; in particular, its minimal polynomial is \(H\).
\end{proof}

The matrix~\eqref{v67:eq-hankel-divisibility} is obtained by selecting the columns indexed by \(J\) from the Hankel matrix \((s_{i+j})\).
Thus one recurrence computes the divisibility tests for different exponent sets.
Its kernel still parametrizes the subcode in the original powers \(A^{r+tj}\), so the coefficient-row test of Proposition~\ref{v54:prop-test-subcode} computes effective support from this matrix as well.

The same rank has a root description that retains multiplicities.
Let \(K\) split \(H\), and write \(H(Y)=\prod_\beta(Y-\beta)^{m_\beta}\) in \(K[Y]\).
Using the Hasse derivatives of Definition~\ref{ffs:def-hasse}, form
\begin{equation}\label{v67:eq-hasse-divisibility}
\mathcal D_J(H)=
\bigl(\partial^{[a]}Y^j\big|_{Y=\beta}\bigr)_{\substack{(\beta,a),\ 0\le a<m_\beta\\j\in J}}.
\end{equation}
Proposition~\ref{v97:prop-hasse-coordinates} relates this matrix to the outer remainder matrix by an invertible left factor, retaining all root multiplicities, including that of a zero root.
Its kernel is therefore the scalar extension to \(K\) of the remainder-map kernel, and
\[
\rank_K\mathcal D_J(H)=\rank_k\mathcal H_J=\rho_{r,J}.
\]
For \(H=1\), both tests have no rows and rank zero.
The recurrence matrix uses coefficient coordinates over \(k\); the Hasse matrix uses root values and their multiplicities over a splitting field.
Both describe the same quotient-algebra constraint.

For example, the algebra of Example~\ref{exm:repeated-recurrence-algebra} has \(H=(Y+1)^2\) over \(\F_2\).
Its sequence is \(0,1,0,1,\ldots\), and for \(J=\{0,1\}\),
\[
\mathcal H_J=\begin{pmatrix}0&1\\1&0\end{pmatrix}.
\]
This matrix retains both quotient coordinates, even though the ordinary algebra trace is identically zero.

\Needspace{24\baselineskip}
\begin{example}[A three-term codeword from a nonzero center]\label{v54:ex-centered-code}
Over \(\F_5\), take \(U=X-1\) and \(f=U(U^2-2)=X^3+2X^2+X+1\).
Put \(\beta^2=2\). The roots of \(f\) are \(1,1+\beta,1-\beta\), all distinct.
Since \((1\pm\beta)^3=2\) and \(2\) has order \(4\), the order of \(1\pm\beta\) can only be \(4\) or \(12\).
The former is impossible: elements of order \(4\) in \(\F_{25}^\times\) belong to \(\F_5^\times\),
but \(1\pm\beta\notin\F_5\). Thus both roots have order \(12\), and hence \(f\mid X^{12}-1\).
It generates a cyclic code \(C\) of length \(12\) and dimension \(9\).

The mean of the roots is \(1\), so every element of the affine stabilizer fixes \(1\). In the coordinate \(U\),
it must preserve \(\{\beta,-\beta\}\), so its slope is \(\pm1\).
Thus \(V=0\), \(Q=e=1\), \(t=2\), \(w=1\), and \(H(Y)=Y-2\).
The corresponding outer factors are
\[
B_0(Y)=Y(Y-2),\qquad B_1(Y)=Y-2.
\]
The extra factor \(Y\) in residue class zero comes precisely from the root \(X=1\); testing only \(Y-2\) does not retain this root.

Take \(r=1\) and \(J=\{0,2\}\), that is, \(W=\operatorname{span}\{U,U^5\}\).
Modulo \(Y-2\), the remainder of \(1,Y^2\) is \(1,4\), so the outer remainder matrix is \((1\ \ 4)\),
whose kernel is generated by \((1,1)\). Hence
\[
C\cap W=k(U^5+U),\qquad
U^5+U=X^5+X-2=f(U^2+2).
\]
The effective support of this subcode is \(\{0,1,5\}\), so \(d_H(C)\le3\).
If a codeword of weight two existed, removing a nonzero scalar and a monomial factor would give
\(f\mid X^j-a\), \(0<j<12\), and \(a\in\F_5^\times\).
The root \(1\) gives \(a=1\), while the other root, of order \(12\), gives \(12\mid j\), a contradiction.
A codeword of weight one cannot vanish at a nonzero root either. Thus \(C\) has parameters \([12,9,3]\).
Here \(\operatorname{sp}(f)=1\), and the original generator polynomial has four terms; the centered power basis and the outer relation produce a codeword of weight three.
\end{example}

\section{Projective Evaluation and Interpolation Duality}
\label{v100:sec-projective-codes}
\subsection{Projective Root Forms and Evaluation Spaces}\label{v84:subsec-projective-evaluation}
A univariate root polynomial records affine points; after adding the point at infinity, a binary homogeneous form records the root set.
Homogenization preserves the root bound while relating changes of coordinates, Frobenius descent, and evaluation spaces.
\phantomsection\label{v85term:def-evaluation-space}\label{v85term:def-evaluation-code}\ignorespaces
For a fixed set of evaluation points, the image of a polynomial vector space under pointwise evaluation is called an \textbf{evaluation space}; viewed as a linear code, it is an \textbf{evaluation code}.
We continue to write \(k=\F_q\) for the base field and use \(N,K\) for code length and dimension, respectively.

\leavevmode\phantomsection\label{v84:def-projective-line}\ignorespaces
Recall the projective line \(\mathbb P^1(k)\) from Subsection~\ref{v85:sec-pair-coefficient-fibers}: its points are the one-dimensional subspaces of \(k^2\).
The point represented by a nonzero vector \((u,v)\) is denoted by \([u:v]\), with \([u:v]=[\lambda u:\lambda v]\) for \(\lambda\in k^\times\).
Every point has a unique representative \([a:1]\), \(a\in k\), or \(\infty=[1:0]\), so there are \(q+1\) points.
The \emph{zeros} of a homogeneous form \(F(U,V)\) do not depend on the representatives; its \emph{values} generally do, since
\(F(\lambda u,\lambda v)=\lambda^{\deg F}F(u,v)\). \phantomsection\label{v85term:def-projective-evaluation}\ignorespaces
Evaluation at these fixed representatives is called \textbf{projective evaluation}; on forms of a fixed degree, the corresponding map is a \textbf{homogeneous evaluation map}. Throughout, projective evaluation uses the representatives fixed above.

\begin{proposition}[The root form of the projective line]\label{v84:prop-projective-root}
Define the \textbf{projective root form}
\begin{equation}\label{v84:eq-projective-root}
 B_q(U,V)=U^qV-UV^q
 =V\prod_{a\in k}(U-aV).
\end{equation}
Its projective zeros are exactly \(\mathbb P^1(k)\), and every linear factor occurs once.
A nonzero binary homogeneous form of degree \(d\) has at most \(d\) distinct projective zeros.
If \(F\) vanishes at every point of \(\mathbb P^1(k)\), then \(B_q\mid F\); in particular, such an \(F\) of degree less than \(q+1\) is zero.
Moreover, for every \(A=\left(\begin{smallmatrix}a&b\\c&d\end{smallmatrix}\right)\in\mathrm{GL}_2(k)\),
\begin{equation}\label{v84:eq-projective-covariance}
 B_q(aU+bV,cU+dV)=\det(A)B_q(U,V).
\end{equation}
\end{proposition}
\begin{proof}
Homogenizing \(X^q-X=\prod_{a\in k}(X-a)\) and multiplying by \(V\) gives \eqref{v84:eq-projective-root}.
If \(F(a,1)=0\), division in \(U\) by the monic linear form \(U-aV\) leaves the remainder \(F(aV,V)=V^dF(a,1)=0\), so \(U-aV\mid F\).
Likewise, \(F(1,0)=0\) is equivalent to \(V\mid F\).
Distinct projective points give nonproportional linear factors. Their product divides \(F\), and comparison of degrees proves both the root bound and the divisibility assertion.
Finally, \(a^q=a,b^q=b,c^q=c,d^q=d\) makes the \(U^{q+1}\) and \(V^{q+1}\) terms cancel in the expansion; the remaining terms are \((ad-bc)U^qV\) and \(-(ad-bc)UV^q\).
\end{proof}

Equation~\eqref{v84:eq-projective-covariance} shows that a change of basis multiplies the form by a nonzero scalar: its roots describe \(k\)-rational directions rather than a special choice of coordinates.
For a nonzero vector \((u,v)\) over an extension field,
\[
 B_q(u,v)=0
 \quad\Longleftrightarrow\quad [u^q:v^q]=[u:v]
 \quad\Longleftrightarrow\quad [u:v]\in\mathbb P^1(k).
\]
On the chart \(v\ne0\), this is \((u/v)^q=u/v\).
Thus \(B_q\) is also the determinant of the two vectors \((u,v)\) and \((u^q,v^q)\).
This binary form is sometimes called a Dickson form; it is not the univariate Dickson polynomial of Section~\ref{sec:dickson-polynomials}.

\begin{proposition}[Descent and period layers of projective root sets]\label{v84:prop-projective-descent}
Let \(L/k\) be a finite extension and \(S\subseteq\mathbb P^1(L)\) a finite set.
Put \(A=\{a\in L:[a:1]\in S\}\), and let \(\varepsilon\in\{0,1\}\) indicate whether \(\infty\in S\). Define
\[
 R_S(U,V)=V^{\varepsilon}\prod_{a\in A}(U-aV).
\]
Then \(R_S\in k[U,V]\) if and only if \(S\) is stable under \([u:v]\mapsto[u^q:v^q]\).
An orbit of finite points corresponds to the homogenization of its univariate minimal polynomial; the point at infinity corresponds to the factor \(V\).
For \(r,s\ge1\), in a common finite extension one has
\[
 \gcd(B_{q^r},B_{q^s})=B_{q^{\gcd(r,s)}}
 \quad\text{up to a nonzero scalar}.
\]
In particular, if \(r\mid s\), then \(B_{q^r}\mid B_{q^s}\). The quotient is square-free, has degree \(q^s-q^r\), and has root set
\(\mathbb P^1(\F_{q^s})\setminus\mathbb P^1(\F_{q^r})\).
\end{proposition}
\begin{proof}
Frobenius fixes \(\infty\) and sends \([a:1]\) to \([a^q:1]\).
After setting \(V=1\), descent is precisely Theorem~\ref{thm:root-descent}; conversely, homogenize the resulting polynomial.
This also identifies the factor of a single orbit.
The intersection of the two sets of rational projective points consists of the affine points in the intersection of the two fields together with their common point at infinity. It is therefore
\(\mathbb P^1(\F_{q^{\gcd(r,s)}})\).
All roots are simple, so comparison of linear factors proves the assertions about the greatest common divisor and the quotient.
\end{proof}

The form \(B_{q^r}\) is not the least-degree homogenization of \(X^{q^r}-X\): the latter is
\(U^{q^r}-UV^{q^r-1}\), which does not vanish at infinity.
Multiplication by \(V\) adds this last point and gives all \(q^r+1\) projective points.
For example, \(B_2=UV(U+V)\), whereas \(B_4/B_2=U^2+UV+V^2\) records the Frobenius orbit of length two in \(\F_4\setminus\F_2\).

For a monic polynomial \(f\) of degree \(d\), write \(f^{\mathrm h}(U,V)=V^df(U/V)\). The orbit description also gives
\begin{equation}\label{v84:eq-projective-orbit-factorization}
 B_{q^r}(U,V)=V\prod_{d\mid r}\ \prod_{f\in\Irr_q(d)}f^{\mathrm h}(U,V).
\end{equation}
Every factor occurs once; this follows by homogenizing the irreducible factorization of \(X^{q^r}-X\).
Thus the rational projective points give \(q+1\) orbits of length one, while projective Frobenius orbits of length \(d>1\) still correspond bijectively to monic irreducible polynomials of degree \(d\).
Adding the point at infinity adds just one fixed orbit and leaves the higher-degree conjugate orbits unchanged. The factorization of projective root forms therefore returns to the earlier irreducible-factor counts.

\begin{proposition}[Kernel and dimension of homogeneous evaluation]\label{v84:prop-projective-evaluation}
Let \(R_d=k[U,V]_d\) be the space of homogeneous forms of degree \(d\), with \(R_d=0\) for \(d<0\).
For \(d\ge0\), the map
\[
 \operatorname{ev}_d:R_d\longrightarrow k^{\mathbb P^1(k)},\qquad
 F\longmapsto\bigl((F(a,1))_{a\in k};F(1,0)\bigr)
\]
satisfies
\[
 \ker\operatorname{ev}_d=B_qR_{d-q-1},\qquad
 \dim\im\operatorname{ev}_d=\min\{d+1,q+1\}.
\]
It is injective for \(d\le q\), and is an isomorphism for \(d=q\).
Nonzero homogeneous vanishing relations first occur in degree \(q+1\).
\end{proposition}
\begin{proof}
Proposition~\ref{v84:prop-projective-root} describes the kernel.
Since \(\dim R_d=d+1\) and multiplication by the nonzero form \(B_q\) is injective, rank--nullity gives the dimension formula.
\end{proof}

There is an important distinction from the affine function algebra \(k[X]/(X^q-X)\).
On the affine line, \(X^q-X\) already represents the zero function; under projective evaluation, its degree-\(q\) homogenization takes the value one at \(\infty\), so it does not belong to the kernel.
The extra coordinate retains exactly this highest-degree information.
These assertions concern each fixed homogeneous degree; forms of different degrees are not identified as the same kind of function without choosing representatives.

\begin{corollary}[Elementary parameter restrictions]\label{v84:cor-elementary-mds}
A nontrivial \([N,K]\) MDS code satisfies
\[
 N\le q+K-1,\qquad K\le q-1,\qquad N\le2q-2.
\]
\end{corollary}
\begin{proof}
Shorten on \(K-2\) coordinates to obtain a two-dimensional MDS code.
Its generator columns represent distinct points of \(\mathbb P^1(k)\), so \(N-K+2\le q+1\).
Applying this inequality to the dual gives \(K\le q-1\), and substitution gives the last bound.
\end{proof}

\begin{proposition}[The root bound and Reed--Solomon evaluation codes]\label{v84:prop-rs}
Let \(S\subseteq k\), with \(|S|=n\ge1\).
For \(1\le K\le n\), the \textbf{Reed--Solomon code}
\[
 \operatorname{RS}_K(S)=\{(f(a))_{a\in S}:f\in k[X],\ \deg f<K\}
\]
is an \([n,K,n-K+1]\) MDS code.
For \(1\le K\le n+1\), the \textbf{extended Reed--Solomon code}
\[
 \operatorname{ERS}_K(S)=
 \{((f(a))_{a\in S};f_{K-1}):\deg f\le K-1\}
\]
is an \([n+1,K,n-K+2]\) MDS code, where \(f_{K-1}\) denotes the coefficient of \(X^{K-1}\).
\end{proposition}
\begin{proof}
In ordinary evaluation, a nonzero \(f\) vanishes at at most \(K-1\) points.
Thus evaluation is injective and the distance is at least \(n-K+1\).
Extended evaluation is the evaluation of \(F(U,V)=V^{K-1}f(U/V)\) on
\(\{[a:1]:a\in S\}\cup\{\infty\}\).
The projective root bound gives injectivity and distance at least \(n-K+2\).
The Singleton bound gives equality in both cases.
\end{proof}

\phantomsection\label{v85term:def-full-evaluation-set}\label{v85term:def-full-length-rs}\label{v85term:def-subgroup-evaluation-set}\ignorespaces
When \(S=k=\F_q\), the set \(S\) is the \textbf{full evaluation set}, and \(\operatorname{RS}_K(k)\) is a \textbf{full-length Reed--Solomon code}, of length \(q\). When \(S=H\le k^\times\), we call \(H\) a \textbf{multiplicative subgroup evaluation set}; the code length is \(|H|\). The term ``full-length'' modifies the code, not the evaluation set. The extended code obtained by adjoining the coordinate at infinity is distinct from the full-length code of length \(q\).

\leavevmode\phantomsection\label{v84:def-projective-arc}\ignorespaces
More generally, \(\mathbb P^{K-1}(k)\) is the set of one-dimensional subspaces of \(k^K\).
A set of at least \(K\) points for which every \(K\) representative vectors are independent is a \textbf{projective arc}.
For \(K\ge2\), the columns of an MDS generator matrix give precisely such an arc.
Changing the representatives only rescales the corresponding code coordinates by nonzero scalars and does not change distance.
For example, the columns of \(\operatorname{ERS}_K(k)\) are \((1,a,\ldots,a^{K-1})^{\mathsf T}\) and \((0,\ldots,0,1)^{\mathsf T}\).
The BCH bound starts with a small coefficient support and obtains evaluation equations; here a small root set of an evaluation polynomial controls the zero coordinates of a word.
Both arguments use root bounds or Vandermonde independence, but the polynomials record different objects.

\begin{example}[An extra evaluation direction in characteristic two]\label{v84:ex-hyperoval}
Here an evaluation direction means the projective class of the linear functional defining a code coordinate; an additional coefficient functional need not be evaluation at another ordinary point.
Suppose that \(q\) is a power of two. Take the columns
\[
 (1,t,t^2)^{\mathsf T}\ (t\in\F_q),\qquad
 (0,0,1)^{\mathsf T},\qquad(0,1,0)^{\mathsf T}.
\]
Any three are independent.
For three finite parameters, the determinant is a Vandermonde product.
For two finite parameters \(a\ne b\) and one of the last two columns, it equals, up to sign, \(b-a\) or \(b^2-a^2=(b-a)^2\), respectively.
For one finite parameter and both special columns, it equals one.
The columns therefore generate an \([q+2,3,q]\) MDS code, whose dual has parameters \([q+2,q-1,4]\).
Such an arc of \(q+2\) points in the projective plane is called a \textbf{hyperoval}.
The independence is verified directly, without a classification of general projective arcs.
For \(q=2\), the dimensions are three and one, respectively, so both codes are trivial boundary cases.
\end{example}

Quotienting by the span of \(s\) points of an arc, with \(0\le s\le K-2\), leaves an arc in the quotient: a dependence there lifts to a dependence involving at most \(K\) original columns.
For \(s=K-2\), the quotient is two-dimensional and contains \(N-K+2\) distinct occupied directions.
The remaining \(q+K-1-N\) directions can be recorded by the product of their linear forms.
The root forms of the occupied and missing directions multiply to a nonzero scalar multiple of \(B_q\).
This is a factorization within each two-dimensional quotient; it does not assert that the missing-direction forms from different quotients automatically come from one global polynomial.

\subsection{Minimum-weight layers and prescribed coordinate values}
\leavevmode\phantomsection\label{v150:def-qary-design}\ignorespaces
A family of equal-weight $q$-ary vectors is a \textbf{$q$-ary $3$-design} if prescribing any three nonzero values at any three distinct coordinates gives the same number of vectors with those values. This definition retains the nonzero symbols; a set design recording only supports omits this information.

\begin{theorem}[A three-design obstruction for the minimum-weight layer of an extended RS code]\label{v150:thm-ers-design}
Let $q\ge3$ and $2\le K\le q-1$. The minimum-weight codewords of $\operatorname{ERS}_K(\F_q)$ do not form a $q$-ary $3$-design.
\end{theorem}
\begin{proof}
Put $r=K-1$ and $N=q-1$, so $1\le r<N$. Fix the coordinate at infinity to be $a\ne0$.
Minimum weight is $q+1-r$, so the corresponding polynomial has precisely the form
\[
 f(X)=a\prod_{s\in S}(X-s),\qquad S\subseteq\F_q,\quad |S|=r.
\]
Prescribing $f(0)=b\ne0$ and $f(1)=c\ne0$ forces $S\subseteq\F_q\setminus\{0,1\}$ and
$\prod S=(-1)^rb/a$, $\prod(1-S)=c/a$. The design property would therefore make
\begin{equation}\label{v150:eq-joint-products}
 B_r(u,v)=\#\{S\subseteq\F_q\setminus\{0,1\}:|S|=r,
 \ \textstyle\prod S=u,\ \prod(1-S)=v\}
\end{equation}
constant on all $(u,v)\in(\F_q^\times)^2$.

Choose a generator $g$ of $\F_q^\times$ and set $\chi(g^j)=\exp(2\pi i j/N)$. If $B_r$ were constant, a geometric sum would give $\sum_{u,v}B_r(u,v)\chi(u)=0$.
Expanding by subsets instead gives
\[
 \sum_{u,v}B_r(u,v)\chi(u)
 =[T^r]\prod_{s\ne0,1}(1+T\chi(s))
 =[T^r]\frac{1-(-T)^N}{1+T}=(-1)^r\ne0,
\]
a contradiction. The last step uses $r<N$; the deleted element one contributes exactly the factor $1+T$.
\end{proof}

Proposition~\ref{v150:prop-subset-products} compares subset products on the full multiplicative group; Proposition~\ref{v150:prop-joint-products} further recovers these joint fibers through Jacobi sums.

\subsection{Interpolation Weights, Duality, and Squaring in Characteristic Two}\label{v84:subsec-interpolation-dual}
Duality of evaluation spaces is not a separate construction.
The weights that extract the highest coefficient in Lagrange interpolation give the orthogonality relations for low-degree evaluation spaces.
They have the same source as the trace-dual weights \(1/f'(\alpha)\) for power bases in Chapter~\ref{chap:dual}.

\begin{proposition}[Dual weights of evaluation codes]\label{v84:prop-rs-dual}
Let \(S\subseteq k\), with \(|S|=n\ge1\), and put \(P_S(X)=\prod_{a\in S}(X-a)\).
For \(1\le K\le n\),
\[
 \operatorname{RS}_K(S)^\perp
 =\left\{\left(\frac{g(a)}{P_S'(a)}\right)_{a\in S}:\deg g<n-K\right\},
\]
where the right side is the zero space for \(K=n\).
The dual of extended evaluation is
\begin{equation}\label{v84:eq-ers-dual}
 \operatorname{ERS}_K(S)^\perp
 =\left\{\left(\left(\frac{g(a)}{P_S'(a)}\right)_{a\in S};-g_{n-K}\right):
               \deg g\le n-K\right\}.
\end{equation}
For $K=n+1$, one has $\operatorname{ERS}_{n+1}(S)=k^{n+1}$, whose dual is the zero space.
\end{proposition}
\begin{proof}
The coefficient identity~\eqref{v84:eq-euler-lagrange} in the appendix gives
\[
 \sum_{a\in S}\frac{A(a)}{P_S'(a)}=[X^{n-1}]A
 \qquad(\deg A\le n-1).
\]
For ordinary evaluation take \(A=fg\); then \(\deg A\le n-2\), so the inner product is zero.
For extended evaluation, \(\deg A\le n-1\), and the coefficient of \(X^{n-1}\) is \(f_{K-1}g_{n-K}\), which is canceled by the product of the coordinates at infinity.
Both parameter maps on the right are injective. For ordinary evaluation this is the root bound.
For extended evaluation, if all coordinates vanish, the highest coefficient is zero, and the remaining polynomial has degree less than \(n\) and vanishes throughout \(S\), so it is zero.
The dimensions are \(n-K\) and \(n+1-K\), respectively, proving equality in each orthogonal inclusion.
Finally, for $K=n+1$ the extended evaluation map has zero kernel: vanishing of the highest coefficient leaves a polynomial of degree less than $n$ that vanishes at all $n$ finite nodes. Equal dimensions then give surjectivity.
\end{proof}

\begin{corollary}[Constant derivative and self-duality without column multipliers]\label{v116:cor-unweighted-self-dual}
Let \(q\) be even and \(B\subseteq\F_q\), with \(|B|=2k\ge2\). Then
\[
\operatorname{RS}_k(B)=\operatorname{RS}_k(B)^\perp
\quad\Longleftrightarrow\quad R_B'\in\F_q^\times.
\]
\end{corollary}
\begin{proof}
If \(R_B'\) is a nonzero constant, all dual weights are equal, and the duality formula gives equality immediately.
Conversely, put \(w_b=1/R_B'(b)\). Self-duality makes \((w_b)\) the evaluation of some \(h\) with \(\deg h<k\), and makes \((w_bb^{k-1})\) the evaluation of some \(g\) with \(\deg g<k\).
Thus \(hX^{k-1}-g\) vanishes at all \(2k\) points and has degree at most \(2k-2\), so it is zero.
Since \(h\ne0\), we obtain \(\deg h=0\). Hence \(R_B'\) has the same nonzero value on \(B\); as \(\deg R_B'<|B|\), it is itself constant.
\end{proof}
Thus \(c_q(A)\) is the shortest length of such an ordinary self-dual evaluation code retaining the evaluation set \(A\).
There are no independent column multipliers here; allowing them changes the self-duality condition.

In characteristic two, squaring also turns a differential interpolation condition into ordinary low-degree evaluation.
The following result gives the exact correspondence between this condition and the evaluation space\cite[Lemma 3.1 and Proposition 3.2]{ChenHuang2026Tools}.

\begin{proposition}[Differential interpolation and evaluation at squares]\label{v84:prop-differential-interpolation}
Let \(q=2^e\ge4\), \(m=q/2\), and let \(H\subseteq\F_q\) satisfy \(1\le h=|H|\le m-1\). Put \(S=\F_q\setminus H\) and \(P_H=\prod_{\tau\in H}(X-\tau)\).
For data \(z=((z_s)_{s\in S};a)\), define
\[
 F_z(X)=\sum_{s\in S}z_s\frac{X^q-X}{X-s}.
\]
This is the unique polynomial of degree less than \(q\) that vanishes on \(H\) and takes the values \(z_s\) on \(S\).
Let \(\mathcal E_H\) be the space of data for which some \(Q\) of degree less than \(h\) satisfies
\begin{equation}\label{v84:eq-differential-data}
 P_H(s)\bigl(F_z'(s)+a\bigr)=Q(s)\qquad(s\in S).
\end{equation}
Then
\begin{equation}\label{v84:eq-square-evaluation}
 \mathcal E_H=
 \{((E(s)+as)_{s\in S};a):
   E\in\operatorname{span}_k\{1,X^2,\ldots,X^{q-2}\},\ a\in k\}.
\end{equation}
After relabeling finite coordinates by \(s\leftrightarrow s^2\), one has \(\mathcal E_H=\operatorname{ERS}_{m+1}(S^2)\), of dimension \(m+1\).
Proposition~\ref{v84:prop-rs} and duality for MDS codes give the parameters
\[
 \begin{gathered}
 \boxed{\mathcal E_H:\ [q-h+1,\ q/2+1,\ q/2-h+1]},\\
 \boxed{\mathcal E_H^\perp:\ [q-h+1,\ q/2-h,\ q/2+2]}.\end{gathered}
\]
The hypothesis \(h\le q/2-1\) keeps the dual dimension positive, so both minimum distances have their usual meaning in terms of nonzero codewords.
Writing \(t=m-h\), its dual is
\begin{equation}\label{v84:eq-square-dual}
 \mathcal E_H^\perp=
 \{((P_H(s)^2g(s^2))_{s\in S};g_{t-1}):\deg g\le t-1\}.
\end{equation}
\end{proposition}
\begin{proof}
Since \((X^q-X)'=1\), interpolation gives the values and uniqueness of \(F_z\).
Suppose \(z\in\mathcal E_H\), and put \(A=F_z'+a\), \(P_S=(X^q-X)/P_H\).
Condition~\eqref{v84:eq-differential-data} implies
\[
 P_HA-Q=P_SG,\qquad \deg G\le2h-2.
\]
Multiply by \(P_H\) and differentiate. In characteristic two \(A'=0\), so
\[
 P_H^2A-P_HQ=(X^q-X)G,
 \qquad (P_HQ)'=G+(X^q-X)G'.
\]
The left side and \(G\) have degree at most \(2h-2<q\).
If \(G'\ne0\), the highest term on the right in the last identity has degree at least \(q\) and cannot cancel.
Hence \(G'=0\) and \(G=(P_HQ)'\), giving
\[
 A=\frac Q{P_H}+(X^q-X)\left(\frac Q{P_H}\right)'.
\]
Since \(P_H\) has distinct roots and \(\deg Q<h\), partial fractions give
\(Q/P_H=\sum_{\tau\in H}c_\tau/(X-\tau)\).
Substitution and \(\tau^q=\tau\) yield
\[
 F_z'=a+\sum_{\tau\in H}c_\tau(X-\tau)^{q-2}.
\]
Thus \(F_z-aX-\sum_\tau c_\tau(X-\tau)^{q-1}\) has zero derivative and contains only even powers.
For \(s\in S\), each \((s-\tau)^{q-1}=1\); their constant sum can be absorbed into the even polynomial, giving the form in \eqref{v84:eq-square-evaluation}.

Conversely, given \(E,a\) on the right, set
\[
 c_\tau=E(\tau)+a\tau,\quad c_0=\sum_{\tau\in H}c_\tau,\quad
 F=E+c_0+aX+\sum_{\tau\in H}c_\tau(X-\tau)^{q-1}.
\]
Direct evaluation shows that \(F\) vanishes on \(H\) and equals \(E(s)+as\) on \(S\), hence \(F=F_z\).
Taking \(Q=P_H\sum_\tau c_\tau/(X-\tau)\) gives \(\deg Q<h\); differentiation and evaluation on \(S\) give \eqref{v84:eq-differential-data}.

Write \(E(X)=\sum_{j=0}^{m-1}e_jX^{2j}\). With \(x=s^2\), we have \(x^m=s^q=s\), and therefore
\[
 E(s)+as=\sum_{j=0}^{m-1}e_jx^j+ax^m.
\]
The last coordinate is the highest coefficient \(a\), so this is the extended evaluation space.
Its parameter map is injective by Proposition~\ref{v84:prop-rs}.
Finally, \(P_HP_S=X^q-X\) gives \(P_S'(s)=P_H(s)^{-1}\), hence
\[
 P_{S^2}'(s^2)=P_S'(s)^2=P_H(s)^{-2}.
\]
Apply \eqref{v84:eq-ers-dual} to \(S^2\), noting that minus and plus agree in characteristic two, to obtain \eqref{v84:eq-square-dual}.
\end{proof}

\begin{example}[A seven-coordinate MDS code from differential conditions]\label{v104:ex-differential-eight}
Take \(q=8\) and \(H=\{0,1\}\). Then \(P_H=X^2+X\) and \(|S|=6\).
Proposition~\ref{v84:prop-differential-interpolation} turns the differential conditions into an \([7,5,3]\) MDS code, with an \([7,2,6]\) MDS dual.
For example, the dual has the two-parameter expression
\[
 \left(\bigl((s^2+s)^2(b_0+b_1s^2)\bigr)_{s\in S};b_1\right),
 \qquad b_0,b_1\in\F_8.
\]
All weights \((s^2+s)^2\) on \(S\) are nonzero. A nonzero linear polynomial has at most one zero among the squared evaluation points;
when \(b_1=0\), only the last coordinate is zero. This also gives the dual distance six directly.
\end{example}

Squaring not only lowers the effective degree; it also explains the dual scaling of projective representatives.
On a two-dimensional space, let \(H_0,f\) be homogeneous of degrees \(h,t-1\), with \(h+t=q/2\).
Then \((H_0f)^2\) has homogeneous degree \(q-2\), and for \(\lambda\in\F_q^\times\),
\[
 (H_0f)^2(\lambda v)=\lambda^{q-2}(H_0f)^2(v)
 =\lambda^{-1}(H_0f)^2(v).
\]
This is exactly the inverse scaling needed for dual coordinates when the original columns are rescaled.
Also, \(f\mapsto(H_0f)^2\) is semilinear for the squaring Frobenius.
Since squaring is bijective on \(\F_q\), its image is still an \(\F_q\)-linear space.
Semilinearity differs from linearity, although the usual dimension formula still holds for kernels and images.

\subsection{Low-Degree and Critical-Degree Rigidity of the Frobenius Derivation}\label{v84:subsec-frobenius-derivation}
Formal differentiation detects characteristic powers, whereas Frobenius identifies directions rational over the base field.
Combining the two turns a differential identity into restrictions on polynomial support and root sets, rather than merely an equality of values on a finite set.

\leavevmode\phantomsection\label{v84:def-frobenius-derivation}\ignorespaces
Let \(k=\F_q\), with \(q=p^e\). On \(R=k[X_1,\ldots,X_r]\), define the \textbf{Frobenius derivation}
\[
 \mathcal D_q=\sum_{i=1}^rX_i^q\frac{\partial}{\partial X_i}.
\]
A derivation is a \(k\)-linear operator satisfying the Leibniz product rule; \(\mathcal D_q\) is not the Frobenius power map \(F\mapsto F^q\).
Write \(R_d\) for the \(k\)-vector space of homogeneous forms of degree \(d\) in \(R\), for \(d\ge0\).

\begin{theorem}[Low-degree and critical-degree rigidity of the Frobenius derivation]\label{v84:prop-frobenius-derivation}\label{v86:thm-frobenius-kernel}
Every \(k\)-linear form \(L\) satisfies \(\mathcal D_qL=L^q\), so \(\mathcal D_q((L))\subseteq(L)\).
The derivation has the same expression after any \(k\)-linear change of coordinates, and \(\mathcal D_q^p=0\).
For every \(0\le d\le q\), one has
\begin{equation}\label{v86:eq-low-degree-homogeneous-kernel}
 \boxed{\ker(\mathcal D_q|_{R_d})=
 \begin{cases}
 \{Q^p:Q\in R_{d/p}\},&p\mid d,\\
 \{0\},&p\nmid d.
 \end{cases}}
\end{equation}
More generally, for any polynomial \(P\) of total degree at most \(q\),
\begin{equation}\label{v84:eq-low-degree-p-power}
 \mathcal D_qP=0
 \quad\Longleftrightarrow\quad
 P\in k[X_1^p,\ldots,X_r^p]
 \quad\Longleftrightarrow\quad P=Q^p\text{ for some }Q\in R.
\end{equation}
At the first degree beyond this range, the kernel is exactly
\begin{equation}\label{v86:eq-critical-kernel}
 \boxed{\ker(\mathcal D_q|_{R_{q+1}})
 =\operatorname{span}_{k}
 \{X_i^qX_j-X_iX_j^q:1\le i<j\le r\}.}
\end{equation}
These generators are linearly independent. Consequently,
\begin{equation}\label{v86:eq-critical-dimension}
 \boxed{\dim_k\ker(\mathcal D_q|_{R_{q+1}})=\binom r2.}
\end{equation}
All assertions hold in arbitrary characteristic. In particular, the critical kernel is zero when \(r=1\); for \(r\ge2\), degree \(q+1\) is the first degree in which a kernel element that is not a \(p\)th power can occur.
\end{theorem}
\begin{proof}
For \(L=\sum_i a_iX_i\), with \(a_i\in k\), one has \(\mathcal D_qL=\sum_i a_iX_i^q=L^q\).
The product rule gives preservation of \((L)\). New linear coordinates \(Y_i\) still satisfy \(\mathcal D_qY_i=Y_i^q\).
Since a derivation is determined by its values on the variables, its coordinate expression is unchanged.
The operators \(X_i^q\partial_{X_i}\) commute, and \(\partial_{X_i}(X_j^q)=0\).
Expanding the \(p\)th power of their sum in characteristic \(p\) gives
\(\mathcal D_q^p=\sum_iX_i^{pq}\partial_{X_i}^p=0\), since any product of \(p\) consecutive integers vanishes in \(k\).

Suppose \(\deg P\le q\). Every nonzero monomial in \(X_i^q\partial_{X_i}P\) has \(X_i\)-exponent at least \(q\), while its total degree is at most \(2q-1\); it cannot have a second variable whose exponent is at least \(q\).
Thus the nonzero supports for distinct \(i\) are disjoint, and \(\mathcal D_qP=0\) if and only if every partial derivative vanishes.
Termwise differentiation shows that this happens precisely when all exponents in all monomials are divisible by \(p\). Since the \(p\)th-power map on \(k\) is bijective, this is equivalent to \(P=Q^p\).
If \(P\) is homogeneous of degree \(d\), then \(Q\) must be homogeneous of degree \(d/p\); when \(p\nmid d\), no such nonzero \(Q\) exists. This proves the low-degree assertions.

Now let \(P\in R_{q+1}\) satisfy \(\mathcal D_qP=0\). Euler's homogeneous identity gives
\[
 \sum_iX_i\frac{\partial P}{\partial X_i}=(q+1)P=P,
\]
for \(q=0\) in \(k\). For every \(x\in k^r\), since \(x_i^q=x_i\), we obtain
\[
 P(x)=\sum_i x_i\frac{\partial P}{\partial X_i}(x)
      =(\mathcal D_qP)(x)=0.
\]
We determine the first possible homogeneous vanishing degree and its entire space, without invoking a general theorem on projective vanishing ideals.

First we need a low-degree affine observation: if \(f\in k[Y_1,\ldots,Y_s]\) has total degree at most \(q\) and vanishes on \(k^s\), then
\begin{equation}\label{v86:eq-affine-vanishing}
 f=\sum_{i=1}^s c_i(Y_i^q-Y_i),\qquad c_i\in k.
\end{equation}
Indeed, subtract from \(f\) the coefficients of its \(Y_i^q\)-terms multiplied by \(Y_i^q-Y_i\).
The remaining polynomial vanishes on \(k^s\) and has degree less than \(q\) in every variable; successive applications of the univariate root bound make it identically zero.

First, there are no nonzero homogeneous vanishing forms in degrees at most \(q\). The degree-zero case is immediate. For \(1\le d\le q\), proceed by induction on the number of variables: the one-variable case is clear. If the assertion holds for \(r-1\) variables, a homogeneous vanishing form of degree \(d\) restricts to zero at \(X_r=0\), so it has the form \(X_rG\), where \(\deg G=d-1<q\). As \(G(x_1,\ldots,x_{r-1},1)=0\) throughout the affine chart, the successive univariate root bound makes \(G=0\).

For degree \(q+1\), proceed by induction on \(r\). When \(r=1\), a homogeneous form \(aX_1^{q+1}\) vanishing on \(k\) has \(a=0\).
Assume the assertion for \(r-1\) variables. The restriction \(P(X_1,\ldots,X_{r-1},0)\) vanishes on \(k^{r-1}\); by induction it is a linear combination of
\(X_i^qX_j-X_iX_j^q\), with \(i<j<r\).
Subtract this linear combination from \(P\). The remainder is still a homogeneous vanishing form and is divisible by \(X_r\), say \(X_rF\), with \(\deg F=q\).
Since \(F(x_1,\ldots,x_{r-1},1)=0\) for all \((x_1,\ldots,x_{r-1})\in k^{r-1}\), the affine observation followed by homogenization yields
\[
 F=\sum_{i=1}^{r-1}c_i(X_i^q-X_iX_r^{q-1}),
 \qquad
 X_rF=\sum_{i=1}^{r-1}c_i(X_i^qX_r-X_iX_r^q).
\]
This completes the induction and shows that the right side of~\eqref{v86:eq-critical-kernel} spans every homogeneous vanishing form of degree \(q+1\).
Conversely, direct differentiation gives
\[
 \mathcal D_q(X_i^qX_j-X_iX_j^q)
 =X_i^qX_j^q-X_i^qX_j^q=0.
\]
Thus these forms belong to the derivation kernel and give the stated equality.
Finally, for \(i<j\), the monomial \(X_i^qX_j\) occurs in no generator other than the one with indices \((i,j)\).
As \(q\ge2\), the generators are linearly independent, giving dimension \(\binom r2\).
\end{proof}

Let \(I(\mathbb P^{r-1}(k))\) denote the ideal generated by all homogeneous forms vanishing on every vector in \(k^r\), called the \textbf{homogeneous vanishing ideal of projective rational points}.
This ideal is homogeneous, and its degree-\(d\) homogeneous component consists precisely of the degree-\(d\) forms vanishing at every rational point.
The argument above proves the \textbf{critical-degree rigidity} identity
\begin{equation}\label{v86:eq-vanishing-kernel}
 \boxed{\ker(\mathcal D_q|_{R_{q+1}})
       =I(\mathbb P^{r-1}(k))_{q+1}.}
\end{equation}
For \(r=2\), this space is spanned by the projective root form \(B_q(U,V)=U^qV-UV^q\) already introduced in this chapter.
In general, it is spanned by the pairwise \textbf{projective Frobenius relations} \(X_i^qX_j-X_iX_j^q\).
Before the critical degree, the derivation kernel comes from characteristic powers; at the critical degree, new relations arise from all rational directions of the finite projective space.
This is a uniform phenomenon in odd characteristic and characteristic two.

\begin{example}[The first projective relations in three variables, in odd characteristic]\label{v86:ex-odd-critical-kernel}
Take \(k=\F_3\) and \(r=3\). Put
\[
 B_{ij}=X_i^3X_j-X_iX_j^3\quad(1\le i<j\le3).
\]
The linearly independent forms \(B_{12},B_{13},B_{23}\) form a basis of \(\ker(\mathcal D_3|_{R_4})\).
For instance,
\[
 \begin{aligned}
 P&=(X_1+X_2)^3(X_2+X_3)-(X_1+X_2)(X_2+X_3)^3\\
  &=B_{12}+B_{13}+B_{23}.
 \end{aligned}
\]
Since \(x^3=x\) for all \(x\in\F_3\), the form \(P\) vanishes on \(\F_3^3\), and \(\mathcal D_3P=0\).
It is nevertheless nonzero and has degree four, not divisible by three, so it is not a cube.
By contrast, below degree four there are no kernel elements that are not cubes.
\end{example}

All identities involving \(\mathcal D_q\) here are polynomial identities; equality of the corresponding functions on rational points is insufficient.
For example, \(X^q-X\) represents the zero function on \(k\), but
\(\mathcal D_q(X^q-X)=-X^q\) is not the zero polynomial.
Therefore \(\mathcal D_q\) does not descend directly to \(k[X]/(X^q-X)\).

\subsection{Compatible Restrictions and Recovery of Homogeneous Polynomials}\label{v84:subsec-homogeneous-gluing}
Univariate interpolation recovers a polynomial from point values.
In three-dimensional space, a homogeneous form can also be recovered from compatible restrictions to several two-dimensional subspaces.
The restrictions below are polynomial restrictions in coordinate rings, not merely equality as functions on the rational points of a finite field.

\begin{proposition}[Compatible recovery on linear sections]\label{v84:prop-gluing}
Let \(k\) be any field and let \(L_1,\ldots,L_N\in k[X_1,X_2,X_3]_1\setminus\{0\}\) be pairwise nonproportional, with any three distinct forms linearly independent.
Put \(H_i=\ker L_i\).
For \(0\le d<N\), suppose that each \(H_i\) carries a homogeneous binary form \(g_i\) of degree \(d\), and that
\[
 g_i|_{H_i\cap H_j}=g_j|_{H_i\cap H_j}\qquad(i\ne j).
\]
There is a unique \(G\in k[X_1,X_2,X_3]_d\) satisfying \(G|_{H_i}=g_i\) for every \(i\).
\end{proposition}
\begin{proof}
First lift \(g_1\) to a ternary homogeneous form.
Suppose that \(G\) has already been constructed correctly on \(H_1,\ldots,H_{j-1}\).
The residual form \(g_j-G|_{H_j}\) vanishes on each \(H_i\cap H_j\), \(i<j\), so it is divisible by \(L_i|_{H_j}\).
Independence of any three original forms makes these restricted factors pairwise nonproportional. Their product therefore divides the residual form.
If \(j-1>d\), the residual form is zero by degrees.
Otherwise lift the quotient to a ternary homogeneous form \(Q\) of degree \(d-j+1\), and replace \(G\) by
\(G+(L_1\cdots L_{j-1})Q\).
This preserves all previous restrictions and corrects the one on \(H_j\), proving existence inductively.
If \(G,\widetilde G\) have the same restrictions, each \(L_i\) divides their difference.
Their product of degree \(N\) therefore divides a form of degree \(d<N\), forcing the difference to be zero.
\end{proof}

The proof applies over an arbitrary field and uses only divisibility by linear forms and a degree comparison.
Over a finite field, homogeneity makes one value at a nonzero representative determine the entire restriction to a one-dimensional intersection.
Without homogeneity, equality at one point does not justify this conclusion.
The condition \(d<N\) controls uniqueness: at \(d=N\), adding \(cL_1\cdots L_N\) to any solution preserves every restriction.

\begin{example}[A binary code glued from three local quadratic forms]\label{v104:ex-projective-gluing-code}
Over \(k=\F_2\), order the seven points of \(\mathbb P^2(k)\) as
\[
 e_1,\ e_2,\ e_3,\ e_1+e_2,\ e_1+e_3,\ e_2+e_3,\ e_1+e_2+e_3.
\]
Let \(C_2\) be the evaluation code of all ternary quadratic homogeneous forms.
The three coordinate lines \(X_i=0\) cover the first six points.
On each line, a binary quadratic form \(aU^2+bUV+cV^2\) has values
\((a,c,a+b+c)\) at \([1:0],[0:1],[1:1]\); these values can be prescribed arbitrarily and determine the local form uniquely.
Arbitrary values on the first six coordinates therefore give three local forms agreeing on their pairwise intersections.
Proposition~\ref{v84:prop-gluing}, with \(d=2<3=N\), glues them to a unique ternary quadratic form.
Thus the first six coordinates form an information set of \(C_2\).
Each square monomial is one at four of the seven points, and each mixed monomial is one at two, so every evaluation vector has coordinate sum zero.
Consequently
\[
 C_2=\{c\in\F_2^7:\ \textstyle\sum_{i=1}^7c_i=0\},\qquad
 C_2:[7,6,2],\quad C_2^\perp:[7,1,7].
\]
This constructs codewords from local data and gives the check matrix \((1,1,1,1,1,1,1)\).

The critical derivation kernel also determines the change at the next degree.
The space of ternary cubic forms has dimension ten. By~\eqref{v86:eq-vanishing-kernel}, its evaluation kernel is precisely the span of the three forms
\(X_i^2X_j-X_iX_j^2\), \(i<j\).
Hence the cubic evaluation code \(C_3\) has dimension \(10-3=7\): it is all of \(\F_2^7\), with parameters \([7,7,1]\).
Homogeneous gluing constructs the quadratic code, while the critical Frobenius kernel counts the redundant cubic representations exactly.
\end{example}

\section{Trace Codes, Sparse Multiples, and Multiplicative Subgroup Rigidity}\label{v84:subsec-trace-mds}
The trace model of an irreducible cyclic code already identifies cyclic shift with multiplication by a field element.
For an MDS code of dimension \(K\), every set of \(K\) generator columns must be linearly independent.
In this model, independence means that a multiplicative subgroup contains no short linear relation.

\begin{example}[A five-coordinate irreducible cyclic MDS code]\leavevmode\phantomsection\label{v104:ex-trace-five}
Over \(k=\F_4=\F_2(\omega)\), where \(\omega^2+\omega+1=0\), take
\[
 h=X^2+\omega X+1,\qquad g=X^3+\omega X^2+\omega X+1.
\]
Multiplication gives \(gh=X^5-1\). The polynomial \(h\) has no root in \(\F_4\), so it is irreducible.
The cyclic code \(C=(g)\subseteq k[X]/(X^5-1)\) has generating rows \(g,Xg\), hence generator matrix
\[
 G=\begin{pmatrix}1&\omega&\omega&1&0\\0&1&\omega&\omega&1\end{pmatrix}.
\]
Its five columns represent all five points of \(\mathbb P^1(\F_4)\); every two columns are independent, making \(C\) an \([5,2,4]\) MDS code.
If \(\alpha\) is a root of \(h\), then \(\alpha^4=\alpha^{-1}\) and \(\ord(\alpha)=5\).
Proposition~\ref{v40:prop-cyclic-trace-model} expresses the same code as
\[
 C=\{(\Tr_{\F_{16}/\F_4}(a\alpha^{-j}))_{j=0}^{4}:a\in\F_{16}\}.
\]
Also \(\langle\alpha\rangle\cap\F_4^\times=\{1\}\), so any two distinct subgroup elements are linearly independent over the base field.
The short-relation and rigidity conditions below can first be compared with this equality case; Exercise~\ref{v84:ex-dimension-two} treats the general two-dimensional situation.
\end{example}

\begin{lemma}[Full coordinate period of an irreducible MDS code]\label{v84:lem-faithful}
Let a cyclic code \(C=(g)\) of length \(N\) have irreducible check polynomial \(h=(X^N-1)/g\), with \(\deg h=K\ge2\).
If \(C\) is MDS, every root \(\alpha\) of \(h\) has multiplicative order exactly \(N\); in particular, \(p\nmid N\).
\end{lemma}
\begin{proof}
Put \(s=\ord(\alpha)\). Then \(s\mid N\), \(p\nmid s\), and \(h\mid X^s-1\).
Writing \(N=rs\), we obtain
\[
 g=(1+X^s+\cdots+X^{(r-1)s})\frac{X^s-1}{h}.
\]
Every word is uniquely represented by \(ag\), with \(\deg a<K\).
The polynomial \(b=a(X)(X^s-1)/h\) has degree less than \(s\), so the coefficient vector of \(ag\) consists of \(r\) identical blocks of length \(s\).
If \(r>1\), a generator matrix has two equal columns, contradicting the MDS column criterion for \(K\ge2\).
\end{proof}

Fix henceforth \(\ord(\alpha)=N\), \(E=k(\alpha)=\F_{q^K}\), \(K=\ord_N(q)\), and \(H=\langle\alpha\rangle\).
With the right-shift convention of Proposition~\ref{v40:prop-cyclic-trace-model}, the code with check polynomial \(m_\alpha\) is
\begin{equation}\label{v84:eq-trace-code-convention}
 C_\alpha=\left\{\bigl(\Tr_{E/k}(a\alpha^{-j})\bigr)_{j=0}^{N-1}:a\in E\right\}.
\end{equation}
The inverse root comes from the right-shift convention: it cannot silently be replaced by \(\alpha^j\) while keeping the check polynomial unchanged.
Nevertheless, \(\{\alpha^{-j}:0\le j<N\}=H\), so independence of the entire column set is unaffected.

\begin{proposition}[The relation kernel and the dual distance of a trace code]\label{v96:prop-relation-kernel}
Retain the notation \(E,k,N,K,\alpha,H\) above, and define the \textbf{linear-relation map}
\begin{equation}\label{v96:eq-relation-map}
\mathcal R_\alpha:k^N\longrightarrow E,\qquad
(c_0,\ldots,c_{N-1})\longmapsto\sum_{j=0}^{N-1}c_j\alpha^{-j}.
\end{equation}
Then \(\mathcal R_\alpha\) is surjective and
\begin{equation}\label{v96:eq-dual-distance}
\boxed{\ker\mathcal R_\alpha=C_\alpha^\perp,\qquad
\lambda_q(H)=d_H(C_\alpha^\perp).}
\end{equation}
Here \(\lambda_q(H)\) is the smallest support size of a nontrivial \(k\)-linear relation among elements of \(H\); Section~\ref{v95:sec-subgroup-periods} studies this quantity directly as a finite-field invariant.
\end{proposition}
\begin{proof}
The subgroup \(H\) spans \(E\), so the map is surjective. For \(c\in k^N\) and \(a\in E\),
\[
\sum_j c_j\Tr_{E/k}(a\alpha^{-j})
 =\Tr_{E/k}\bigl(a\mathcal R_\alpha(c)\bigr).
\]
By nondegeneracy of the trace pairing, the left side vanishes for every \(a\) exactly when \(\mathcal R_\alpha(c)=0\), proving the kernel identity.
The elements \(1,\alpha^{-1},\ldots,\alpha^{-(N-1)}\) are distinct and exhaust \(H\), so the support size of a relation is the Hamming weight of its kernel vector.
Finally, \(K=\ord_N(q)\le\varphi(N)<N\), so the kernel is nonzero and its minimum nonzero weight is defined.
\end{proof}
Thus the shortest relation length is the minimum distance of the \emph{dual code}, not of the trace code itself.
The weights of the trace code are instead determined by intersections of \(H\) with trace hyperplanes; Subsection~\ref{v95:subsec-trace-hyperplanes} connects the two descriptions.

\begin{corollary}[The binary relation code of a norm-one group]\label{v153:cor-normone-binary-code}
Let $f\ge1$ be an integer, $Q=2^f$, $N=Q+1$, and index the coordinates by $H=\mu_N\subseteq\F_{Q^2}^{\times}$. Set
\[
 \mathcal D_H=\left\{(c_x)_{x\in H}\in\F_2^H:\sum_{x\in H}c_xx=0\right\}.
\]
Ordering the coordinates by consecutive powers of a generator of $H$ makes $\mathcal D_H$ a cyclic code of dimension $N-2f$. It has no word of weight four. If $f$ is odd, its minimum distance is three, and its weight-three supports are exactly the $N/3$ cosets of $\mu_3$ in $H$. If $f$ is even, its minimum distance is at least five.
\end{corollary}
\begin{proof}
Since $2^f\equiv-1\pmod N$, one has $\ord_N(2)\mid2f$. If this order were smaller than $2f$, it would be a proper divisor of $2f$ and hence at most $f$, contradicting $N>2^f$. Thus $\ord_N(2)=2f$, so $H$ spans $\F_{Q^2}$ over $\F_2$. The relation map is surjective, and its kernel has dimension $N-2f$. A cyclic shift multiplies the relation by a generator of $H$, preserving the kernel.

All nonzero coefficients of a binary word equal one, so its weight is the size of its zero-sum support. Weights one and two are impossible. Proposition~\ref{v153:prop-normone-short-sums} excludes weight four and identifies the weight-three supports as cosets of $\mu_3$. Finally, $3\mid2^f+1$ if and only if $f$ is odd.
\end{proof}

\begin{example}[Relation codes of lengths five and seventeen]\label{v153:ex-normone-binary-codes}
For $f=2$, the code $\mathcal D_H$ has dimension one. Since $\sum_{x\in\mu_5}x=0$, it contains the all-one word and has parameters $[5,1,5]$.
For $f=4$, the factorization over $\F_2$ is
\[
 \Phi_{17}(X)=(X^8+X^5+X^4+X^3+1)
 (X^8+X^7+X^6+X^4+X^2+X+1).
\]
Both factors are irreducible because $\ord_{17}(2)=8$. If $\alpha$ is a root of the first factor, then $\ord(\alpha)=17$ and
\[
 1+\alpha^3+\alpha^4+\alpha^5+\alpha^8=0.
\]
This gives a word of weight five, so the preceding lower bound yields a $[17,9,5]$ code. The second degree-eight factor is precisely $P_3$ obtained by the $\mathcal Q$ transform in Example~\ref{v153:ex-qtransform-seventeen}; the same irreducible orbit now supplies the input for a distance calculation.
\end{example}

\begin{proposition}[Restriction to subgroup coordinates and distance loss]\label{v142:prop-subfield-distance}
Use the notation $E,k,N,K,\alpha,H$. Take $d\mid N$, $d>1$, and put
\[
 s=N/d,\qquad \eta=\alpha^s,\qquad r=\ord_d(q),\qquad
 F=k(\eta)=\F_{q^r}.
\]
Restriction of $C_\alpha$ to the coordinate set $J_d=\{0,s,\ldots,(d-1)s\}$ has image exactly the length-$d$ irreducible trace code $C_\eta$, of dimension $r$. Its kernel is
\begin{equation}\label{v142:eq-subgroup-restriction}
 D_d=\left\{c_a:a\in\ker\Tr_{E/F}\right\},
 \qquad \dim_k D_d=K-r,
\end{equation}
where $c_a$ denotes the trace word in \eqref{v84:eq-trace-code-convention}.
The subcode of $C_\alpha^\perp$ supported on $J_d$ becomes $C_\eta^\perp$ after the other zero coordinates are deleted, and has dimension $d-r$. In particular,
\begin{equation}\label{v142:eq-subfield-distance}
 \begin{aligned}
 d_H(C_\alpha^\perp)&\le r+1,\\
 r<K\quad&\Longrightarrow\quad
 d_H(C_\alpha)\le N-d-K+r+1
                =N-K+1-(d-r).
 \end{aligned}
\end{equation}
Thus, when $r<K$, the minimum distance of $C_\alpha$ falls below the Singleton bound by at least $d-r>0$.
\end{proposition}
\begin{proof}
Transitivity of trace gives
\[
 c_a(st)=\Tr_{F/k}\bigl(\Tr_{E/F}(a)\eta^{-t}\bigr)
 \qquad(0\le t<d).
\]
The relative trace $E\to F$ is surjective, and the trace parametrization $F\to C_\eta$ is injective. Hence the restriction image is $C_\eta$, and its kernel is parametrized exactly by $\ker\Tr_{E/F}$, giving the dimension formula.
A vector supported on $J_d$ is orthogonal to $C_\alpha$ exactly when its restriction is orthogonal to $C_\eta$. The corresponding dual subcode therefore has dimension $d-r$.
Since $r\le\varphi(d)<d$, the Singleton bound for $C_\eta^\perp$ yields $d_H(C_\alpha^\perp)\le r+1$.

If $r<K$, then $D_d\ne0$. Delete the $d$ coordinates in $J_d$, which vanish identically on $D_d$. The resulting code has length $N-d$ and dimension $K-r$, with Hamming weights unchanged. Applying the Singleton bound and using $D_d\subseteq C_\alpha$ gives the second distance bound.
\end{proof}

The first bound records short linear relations within $\langle\eta\rangle$; the second uses trace functionals vanishing on the whole subfield $F$. The same restriction map thus controls the dual distance and the original distance. These bounds require a nontrivial subgroup in a proper subfield, without requiring complete orbit fibers.
For Example~\ref{v142:ex-twenty-five-code}, take $d=5$ and $r=4$. The proposition gives $d_H(C_\xi^\perp)\le5$ and $d_H(C_\xi)\le5$; independent blocks supplied by complete fibers sharpen the latter bound to the exact value two. Exercise~\ref{v142:ex-partial-fibers} considers a distance obstruction with incomplete fibers.

\begin{theorem}[An MDS criterion through sparse multiples]\label{v84:thm-sparse-multiple}
Suppose \(2\le K\le N-2\). The following conditions are equivalent:
\begin{enumerate}
\item \(C_\alpha\) is MDS;
\item every set of at most \(K\) distinct elements of \(H\) is linearly independent over \(k\);
\item one has
\begin{equation}\label{v84:eq-sparse-mds}
 \boxed{\quad
 \min_{\substack{0\ne F\in k[X]\\m_\alpha\mid F,\ \deg F<N}}
       \operatorname{wt}(F)=K+1.\quad}
\end{equation}
\end{enumerate}
In this case all \(K+1\) coefficients of \(m_\alpha\) are nonzero.
\end{theorem}
\begin{proof}
Choose a basis \(a_1,\ldots,a_K\) of \(E/k\). Nondegeneracy of the trace pairing shows that
\[
 E\longrightarrow k^K,\qquad
 x\longmapsto\bigl(\Tr(a_1x),\ldots,\Tr(a_Kx)\bigr)^{\mathsf T}
\]
is a linear isomorphism.
The generator columns of \(C_\alpha\) are the images of the elements of \(H\), so the first two conditions are equivalent by the column criterion.

For \(F=\sum_{i=1}^r b_iX^{e_i}\), where \(b_i\ne0\) and \(0\le e_1<\cdots<e_r<N\),
\[
 m_\alpha\mid F
 \quad\Longleftrightarrow\quad
 \sum_{i=1}^r b_i\alpha^{e_i}=0.
\]
The exponent range makes these subgroup elements distinct.
Thus the second condition says that the minimum weight on the left is at least \(K+1\).
On the other hand, \(m_\alpha\) itself has degree \(K<N\) and weight at most \(K+1\), so equality follows.
Taking \(F=m_\alpha\) proves that all its coefficients are nonzero.
\end{proof}

\begin{remark}[Consecutive powers and Frobenius powers]\label{v96:rem-vandermonde-moore}
Both the BCH bound and the MDS criterion above reduce to nonsingularity, but they require different kinds of independence.
In Theorem~\ref{v56:thm-bch-bound}, a word of weight \(w\) gives distinct nodes \(x_j=\zeta^{e_j}\), and consecutive zero conditions give
\[
V=(x_j^i)_{\substack{0\le i<w\\1\le j\le w}},\qquad
\det V=\prod_{u<v}(x_v-x_u)\ne0.
\]
The initial exponent merely multiplies each column by a nonzero factor, so distinctness of the nodes suffices.

For the trace code \(C_\alpha\), choose \(K\) coordinates \(e_1,\ldots,e_K\), put \(x_j=\alpha^{-e_j}\), and choose a basis \(a_1,\ldots,a_K\) of \(E/k\).
The corresponding generator submatrix factors as
\begin{equation}\label{v96:eq-trace-moore-factorization}
G=(\Tr_{E/k}(a_i x_j))_{i,j}=AM,\qquad
A=(a_i^{q^t})_{i,t},\quad M=(x_j^{q^t})_{t,j},\quad 0\le t<K.
\end{equation}
By the Moore criterion (Theorem~\ref{thm:moore}), \(A\) is invertible, whereas \(M\) is invertible exactly when \(x_1,\ldots,x_K\) are linearly independent over \(k\).
Thus \(G\) is nonsingular over the base field exactly when its Moore model over the extension is nonsingular.
Distinct nodes alone do not suffice here: consecutive powers ensure Vandermonde nonsingularity, whereas Frobenius powers test base-field linear independence.
The former rules out codewords with few nonzero coefficients; the latter requires ruling out every short relation in a specified multiplicative subgroup.
\end{remark}

An irreducible polynomial can be sparse, but the check polynomial of an irreducible cyclic MDS code cannot be sparse.
Moreover, having all coefficients nonzero is not enough: every low-weight multiple of degree less than the code length must be excluded.
The invariant \(w_q(K)\) at \ref{v79:def-defining-complexity} allows a choice among all generators to minimize the support of a minimal polynomial.
In contrast, \eqref{v84:eq-sparse-mds} fixes a multiplicative orbit and tests \emph{all} defining relations within the length range.
One is an optimization over generators; the other excludes short relations for a specified root. They are not interchangeable.
For a general constacyclic code, the dual generator is the monic reciprocal of the check polynomial,
as proved in Theorem~\ref{v87:thm-constacyclic-mds-reciprocal}.
Here we are in the cyclic case \(\lambda=1\): although the dual itself is generated by
\(m_\alpha^\vee\), Formula~\eqref{v84:eq-sparse-mds} may use multiples of \(m_\alpha\),
since reversing the \(N\) coefficient positions preserves Hamming weight and
bijectively relates these two families of multiples.
The trace-model columns use \(\alpha^{-j}\), not \(\alpha^j\); the two unordered sets of
subgroup elements coincide, so the column-independence criterion is unchanged.

\begin{example}[A full-support minimal polynomial with a sparse multiple]\label{v84:ex-full-support-not-mds}
The polynomial \(m=X^2+X+2\) has no root in \(\F_3\), so it is irreducible.
Its root \(\alpha\) satisfies \(\alpha^2=2\alpha+1\) and \(\alpha^4=-1\), giving \(\ord(\alpha)=8\).
Although all three coefficients of \(m\) are nonzero,
\[
 (X^2+X+2)(X^2+2X+2)=X^4+1.
\]
The right side has degree less than eight and weight only two.
Hence the length-eight, dimension-two cyclic code with check polynomial \(m\) is not MDS.
\end{example}

\begin{theorem}[Full-degree rigidity of multiplicative subgroups]\label{v84:thm-subgroup-rigidity}
If \(C_\alpha\) is a nontrivial MDS code, then
\begin{equation}\label{v84:eq-all-divisor-orders}
 \ord_d(q)=K\qquad(d\mid N,\ d>1).
\end{equation}
Consequently, \(K\mid\ell-1\) for every prime \(\ell\mid N\), and
\[
 N\text{ is odd},\qquad
 \gcd(N,q-1)=\gcd(N,K)=1,\qquad N\equiv1\pmod K,\qquad N\ge2K+1.
\]
Every \(x\in H\setminus\{1\}\) is normal in \(E/k\), and all coefficients of its degree-\(K\) minimal polynomial are nonzero.
\end{theorem}
\begin{proof}
For $d\mid N$, $d>1$, we have $r=\ord_d(q)\mid K$. If $r<K$, Proposition~\ref{v142:prop-subfield-distance} gives $d_H(C_\alpha)<N-K+1$, contrary to the MDS condition. Hence $r=K$. Taking \(d=\ell\) gives \(K\mid\ell-1\); all prime divisors are larger than \(K\ge2\).
This proves oddness and \(\gcd(N,K)=1\).
If \(\ell\mid q-1\), then \(\ord_\ell(q)=1\), again contradicting \(K\ge2\).
Multiplying the prime factors modulo \(K\) gives \(N\equiv1\pmod K\), and \(N\ge K+2\) then implies \(N\ge2K+1\).

For \(x\ne1\), put \(d=\ord(x)>1\).
Equation~\eqref{v84:eq-all-divisor-orders} gives a Frobenius orbit of length \(K\).
Its \(K\) conjugates belong to \(H\), so they are independent and form a normal basis.
Also \(d\ge K+1\), making \(1,x,\ldots,x^K\) distinct.
A zero coefficient in the degree-\(K\) minimal polynomial of \(x\) would give a relation involving at most \(K\) of these elements, again contradicting column independence.
\end{proof}

Condition~\eqref{v84:eq-all-divisor-orders} is also equivalent to
\[
 H\cap\F_{q^e}^{\times}=\{1\},\qquad
 \gcd(N,q^e-1)=1\qquad(e\mid K,\ e<K).
\]
Indeed, the intersection has order equal to this gcd, and any nonidentity element in it has order $d>1$ with $\ord_d(q)\mid e$. Conversely, if $\ord_d(q)<K$, the subgroup of order $d$ lies in that proper subfield.

Subgroup orders and subfield degrees provide the first step of this rigidity; normality is the second.
An orbit of length \(K\) need not span a space of dimension \(K\); the MDS column condition supplies this additional independence.
The subgroup restrictions and normality conclusion also appear in\cite[Theorem 3.1 and Proposition 3.7]{Zhang2026IrreducibleTools}.
In particular, \(\Tr_{E/k}(x)\ne0\) for \(x\in H\setminus\{1\}\), since a zero conjugate sum would be a nontrivial relation.
The word with trace parameter \(a=1\) therefore has weight \(N\) or \(N-1\), according as \(p\nmid K\) or \(p\mid K\).
The only potentially zero coordinate is \(\Tr(1)=K\cdot1\).
Normality tests particular Frobenius-orbit subsets and does not replace independence of every \(K\)-subset.
For \(K\ge3\), one also has \(\gcd(N,q^2-1)=1\): a common prime divisor would make the corresponding multiplicative order at most two.

\begin{corollary}[Restrictions from prime divisors and cyclotomic values]\label{v84:cor-cyclotomic-restrictions}
Under these hypotheses, write \(\Omega(N)=\sum_{\ell\mid N}v_\ell(N)\). Then
\[
 N\ge(K+1)^{\Omega(N)},\qquad q\ge(K+1)^{\Omega(N)}-K+1.
\]
If \(K\) is odd, every prime divisor even satisfies \(\ell\ge2K+1\).
Moreover,
\[
 \begin{aligned}
 K\mid q-1&\ \Longrightarrow\ N\le q-K,\\
 K\mid q,\ K\ge3&\ \Longrightarrow\ N\le q-K+1.
 \end{aligned}
\]
If \(N=q+t>q+1\), then
\begin{equation}\label{v84:eq-cyclotomic-value}
 2\le t\le K-1,\qquad N\mid\Phi_K(-t).
\end{equation}
For fixed \(K\), only finitely many pairs \((q,N)\) satisfy this last set of conditions.
\end{corollary}
\begin{proof}
Every prime divisor is at least \(K+1\); combine this with \(q\ge N-K+1\) to obtain the first two inequalities.
If \(K\) is odd and \(\ell=aK+1\) is an odd prime, then \(a\) is even, hence \(a\ge2\).
For \(q=uK+1\), the congruence \(N\equiv1\pmod K\) and \(N\le q+K-1\) give \(N\le q\).
Since \(p\nmid N\), equality is impossible, and the next admissible value is \(q-K\).
For \(q=uK\), the same argument gives \(N\le q+1\); equality would imply \(\ord_N(q)=2\), contrary to \(K\ge3\).

Finally, the elementary length bound gives \(2\le t\le K-1\), while \(q\equiv-t\pmod N\) gives \(N\mid(-t)^K-1\).
For each \(\ell\mid N\), one has \(\ord_\ell(-t)=K\).
In \((-t)^K-1=\prod_{d\mid K}\Phi_d(-t)\), the prime \(\ell\) divides no factor with \(d<K\), since otherwise its order would divide \(d\).
Thus the entire prime-power contribution of \(\ell\) in \(N\) lies in \(\Phi_K(-t)\), proving integer divisibility, not merely divisibility by the distinct prime factors.
There are finitely many choices of \(t\), each nonzero integer \(\Phi_K(-t)\) has finitely many positive divisors, and then \(q=N-t\) is determined.
\end{proof}

These conclusions constrain the linear relations allowed within a specified multiplicative orbit.
In the same trace model, \(x\mapsto\alpha x\) and \(x\mapsto x^q\) give a standard permutation subgroup \(C_N\rtimes_q C_K\) of order \(NK\).
Under the trace-pairing identification of columns with \(H\), a linear map preserving \(H\) permutes the generator columns and hence induces a coordinate automorphism of the code. The two maps have orders \(N,K\), respectively, and satisfy \(\sigma M_\alpha\sigma^{-1}=M_\alpha^q\).
If \(M_\alpha^a\sigma^j=M_\alpha^b\sigma^i\), evaluation at one first gives \(a=b\), and the order of Frobenius then gives \(i=j\).
Only the existence of this subgroup is asserted; it is not thereby identified with the full permutation automorphism group of an individual code.

Beyond its coding-theoretic interpretation, the finite-field question is: how many distinct subgroup elements are needed in the shortest nontrivial base-field linear relation? Section~\ref{v95:sec-subgroup-periods} studies the invariant \(\lambda_q(H)\) through multiplicative orders, cyclotomic intersections, and trace hyperplanes.

\section{Residue Layers and Balance of Fourier Supports}\label{v84:subsec-fourier-balance}
An irreducible trace code retains a single root orbit; a general cyclic code may retain several.
Expanding these orbits into a Fourier basis over a splitting field preserves their Frobenius descent information and makes it possible to compare support sizes in different residue layers.
Throughout this subsection, assume \(\gcd(N,q)=1\), fix a primitive \(N\)th root \(\zeta\), and put \(E=k(\zeta)\).

\leavevmode\phantomsection\label{v84:def-simple-root-code}\ignorespaces A cyclic code satisfying this coprimality condition is called a \textbf{simple-root cyclic code}: the derivative \(NX^{N-1}\) is coprime to \(X^N-1\).

\leavevmode\phantomsection\label{v84:def-fourier-support}\ignorespaces
For a cyclic code \(C\) with defining set \(T(C)\), retain the earlier evaluation convention
\(\widehat c(i)=c(\zeta^i)\), and define its \textbf{Fourier support} by
\[
 S(C)=(\Z/N\Z)\setminus T(C).
\]
This is also called the \textbf{spectral support}; equation~\eqref{v73:eq-code-spectrum} uses the same notation \(S(C)\). It consists of the frequencies at which at least one codeword has a nonzero Fourier coefficient; it is not the Fourier support of a particular word.
Write \(S=S(C)\), \(K=\dim_k C\).
Then \(|S|=K\), \(qS=S\), and
\begin{equation}\label{v84:eq-fourier-basis}
 C_E=\operatorname{span}_E\{v_d:d\in S\},\qquad
 v_d=(\zeta^{-jd})_{j=0}^{N-1}.
\end{equation}
Indeed, \(\widehat v_d(i)=N\) for \(i=d\), and is zero otherwise; Fourier inversion gives \eqref{v84:eq-fourier-basis}.
Since \(\langle v_d,v_e\rangle=N\) exactly when \(d+e=0\), and is zero otherwise,
\begin{equation}\label{v84:eq-dual-support}
 S(C^\perp)=(\Z/N\Z)\setminus(-S).
\end{equation}
The negative exponent in the basis makes \(S\) equal to the Fourier support of the code in the existing evaluation model.
With a positive-exponent basis the support must simultaneously be replaced by \(-S\); changing just one sign would mix conventions.

\begin{theorem}[Balance of Fourier supports]\label{v84:thm-balanced-support}
If \(C\) is a nontrivial cyclic MDS code, then for every \(t\mid N\) and every residue class \(j\) modulo \(t\),
\begin{equation}\label{v84:eq-balanced-support}
 \left\lfloor\frac Kt\right\rfloor
 \le\#\{d\in S:d\equiv j\pmod t\}
 \le\left\lceil\frac Kt\right\rceil.
\end{equation}
Thus the Fourier support of a cyclic MDS code must be distributed as evenly as possible in every residue layering by a divisor of the length.
\end{theorem}
\begin{proof}
Write \(N=t\ell\), \(S_j=\{d\in S:d\equiv j\pmod t\}\), and \(a_j=|S_j|\).
If \(a_j>0\), let \(W_j=\operatorname{span}_E\{v_d:d\in S_j\}\subseteq C_E\), of dimension \(a_j\).
Every \(w\in W_j\) satisfies
\[
 w_{x+y\ell}=\zeta^{-yj\ell}w_x
 \qquad(0\le x<\ell,\ 0\le y<t).
\]
Projection onto the first block of length \(\ell\) is therefore injective, and all blocks have the same weight.
Apply the Singleton bound to this \(a_j\)-dimensional projection and use that \(C_E\) is still MDS:
\[
 N-K+1\le d_H(W_j)
 \le t(\ell-a_j+1).
\]
Hence \(a_j\le1+(K-1)/t\), which for integers is \(a_j\le\lceil K/t\rceil\).
This upper bound also holds when \(a_j=0\).
Apply the same upper bound to the dual code.
By \eqref{v84:eq-dual-support}, its residue class \(-j\) modulo \(t\) has size \(\ell-a_j\), so
\[
 \ell-a_j\le\left\lceil\frac{N-K}{t}\right\rceil
 =\ell-\left\lfloor\frac Kt\right\rfloor,
\]
which gives the lower bound.
\end{proof}

This balance condition concerns only the exponent set and can be tested before a generator matrix is constructed\cite[Lemma 2.8]{LuoCaoLiLi2026Tools}.
It is necessary, not sufficient.
For example, the length-seven support \(S=\{1,2,4\}\subseteq\Z/7\Z\) is balanced modulo every divisor.
Over \(\F_2\), however, it gives the three-dimensional trace code of Example~\ref{v40:ex-cyclic-seven}, which already has a nonzero word of weight four, whereas MDS equality would require distance five.
For prime length, residue balance imposes no restriction on a set of prescribed size.
It must not be mistaken for nonsingularity of all square submatrices of a finite-field Fourier matrix.

\begin{proposition}[Inheritance under zero-class compression]\label{v84:prop-fourier-compression}
Let \(C\) be a nontrivial \([N,K]\) cyclic MDS code, let \(h\mid N\), put \(M=N/h\), and define
\[
 A=\{a\in\Z/M\Z:ha\in S\},\qquad c=|A|.
\]
Then \(qA=A\), so it defines a cyclic code \(B\) over \(k\) of length \(M\) and dimension \(c\), with
\(c\in\{\lfloor K/h\rfloor,\lceil K/h\rceil\}\).
If \(c>0\), then \(B\) is MDS; if \(c=0\), it is the zero code.
\end{proposition}
\begin{proof}
The identity \(qS=S\) gives \(qA=A\), and balance gives the two possible dimensions.
With \(\theta=\zeta^h\), the vector \(v_{ha}\) consists of \(h\) repetitions of \((\theta^{-ja})_{j=0}^{M-1}\).
The subcode they span in \(C_E\) is exactly the \(h\)-fold repetition of \(B_E\).
Thus, when \(c>0\),
\[
 h\,d_H(B_E)\ge N-K+1,
 \qquad d_H(B_E)\ge M-\lceil K/h\rceil+1.
\]
If \(c=\lceil K/h\rceil\), the Singleton bound gives equality immediately.
Otherwise \(h\nmid K\) and \(c=\lfloor K/h\rfloor>0\).
The zero residue class of the dual support is \((\Z/M\Z)\setminus(-A)\), precisely the support of \(B^\perp\), of size
\(M-c=\lceil(N-K)/h\rceil\).
Apply the larger-dimension case already proved to \(C^\perp\); it makes \(B^\perp\) MDS, and duality makes \(B\) MDS.
Finally, scalar extension preserves distance, so the conclusion over \(E\) descends to \(k\).
\end{proof}

The zero residue class produces genuine repeated blocks, so compression remains a cyclic code over the base field.
Other residue classes produce blocks repeated with nonzero proportionality factors.
These factors can change under Frobenius, so the other classes cannot automatically be treated as the same kind of base-field cyclic code.
This distinction parallels the interleaved direct sums obtained from support periods: exponent layers and coefficient blocks are related but are not the same decomposition.

\section{The Smallest Field of Definition of a Cyclic Code}\label{v84:subsec-code-descent}
Root-set descent determines not only the coefficient field of a polynomial but also whether a cyclic code can already be defined over a smaller alphabet.
Coordinates are fixed here: no coordinate scalings or permutations are allowed before fields of definition are compared.

\leavevmode\phantomsection\label{v84:def-code-field}\ignorespaces
Let \(q=p^f\), \(N>1\), and \(\gcd(N,p)=1\), and let \(C\subseteq\F_q^N\) be cyclic.
For \(b\mid f\), say that \(C\) is \textbf{defined over \(\F_{p^b}\) in the fixed coordinates} if
\(C=\operatorname{span}_{\F_q}C_0\) for some \(C_0\subseteq\F_{p^b}^N\).
The least such subfield is the \textbf{smallest field of definition} of this cyclic realization.

\begin{theorem}[The field of definition and the return period of the root set]\label{v84:thm-code-field}
Put \(K=\dim_{\F_q}C\), let \(g\) be the monic generator polynomial of \(C\), put \(T=T(C)\), \(e=\ord_N(p)\), and define
\[
 s=\min\{r\ge1:p^rT=T\}.
\]
Then \(s\mid f\), \(s\mid e\), and for every \(b\mid f\),
\[
 \begin{aligned}
 C\text{ is defined over }\F_{p^b}
 &\quad\Longleftrightarrow\quad g\in\F_{p^b}[X]\\
 &\quad\Longleftrightarrow\quad p^bT=T
 \quad\Longleftrightarrow\quad s\mid b.
 \end{aligned}
\]
The smallest field of definition is therefore \(\F_{p^s}\).
Computing with the complement \(S(C)\) instead of \(T\) gives the same \(s\).
If the check polynomial is irreducible and its root has order \(N\), then
\begin{equation}\label{v84:eq-faithful-coefficient-field}
 s=\gcd(e,f),\qquad K=\ord_N(q)=\frac e{s}.
\end{equation}
In this case, descent to \(\F_{p^s}\) preserves irreducibility and dimension \(K\).
\end{theorem}
\begin{proof}
For the zero code and the whole space, the defining sets are \(\Z/N\Z\) and the empty set, respectively; both have \(s=1\), and the assertions are immediate. We prove the remaining cases below.
Multiplication by \(p\) permutes \(\Z/N\Z\); the return times of \(T\) under this action are precisely the multiples of \(s\).
The identities \(p^fT=T\) and \(p^e\equiv1\pmod N\) give \(s\mid f,e\).
Coefficient Frobenius sends the root set \(\{\zeta^i:i\in T\}\) of \(g\) to \(\{\zeta^{p^bi}:i\in T\}\).
Monicity and root-set descent therefore give \(g\in\F_{p^b}[X]\) if and only if \(p^bT=T\).

If \(g\) has coefficients in this subfield, the generator matrix with rows corresponding to \(g,Xg,\ldots,X^{N-\deg g-1}g\) is also defined there.
Conversely, suppose \(C=\operatorname{span}_{\F_q}C_0\), where \(C_0\subseteq\F_{p^b}^N\) is merely a generating set. First put
\[
C_b=\operatorname{span}_{\F_{p^b}}C_0.
\]
Choose a reduced row-echelon basis \(v_1,\ldots,v_r\) of \(C_b\), with pivot positions \(j_1,\ldots,j_r\).
The pivot submatrix is the identity, so these vectors remain linearly independent over \(\F_q\). They also span \(C\), and hence form an \(\F_q\)-basis of \(C\).
If \(x=\sum_i a_iv_i\in C\cap\F_{p^b}^N\), the pivot coordinates give \(a_i=x_{j_i}\in\F_{p^b}\), so \(x\in C_b\). The reverse inclusion is immediate, and therefore
\[
C_b=C\cap\F_{p^b}^N.
\]
Cyclic shift preserves both \(C\) and \(\F_{p^b}^N\), so \(C_b\) is a cyclic code over \(\F_{p^b}\).
Its monic generator still generates \(C\) after scalar extension; uniqueness makes it \(g\). Thus \(g\) has the required coefficients.
A subset and its complement have the same return period, so \(S(C)\) may replace \(T\).

In the irreducible case, choose a check root \(\alpha\) of order \(N\).
Its root orbit corresponds to a coset \(u\langle p^f\rangle\) in the multiplicative group modulo \(N\), with \(u\) a unit.
In the cyclic group \(\langle p\rangle\) of order \(e\),
\(\langle p^f\rangle=\langle p^{\gcd(e,f)}\rangle\).
Hence this orbit is preserved by \(p^b\) exactly when \(\gcd(e,f)\mid b\).
The generator and check polynomials have the same coefficient field: they are monic with product \(X^N-1\), and monic division recovers either from the other over that field.
This proves \(s=\gcd(e,f)\).
Over \(\F_{p^s}\), the orbit of \(\alpha\) has length \(e/s=K\), so the check polynomial remains irreducible and the dimension stays \(K\).
\end{proof}

This unifies descent of a general root set with descent of a single Frobenius orbit; the precise irreducible case also appears in\cite[Theorem 5.1]{Zhang2026IrreducibleTools}.
The smallest field belongs to the \emph{fixed cyclic realization}; it has not thereby been shown to be smallest among all equivalent codes.

\begin{corollary}[Distance restrictions from proper-subfield descent]\label{v84:cor-code-descent-obstruction}
If \(C\) is a nontrivial MDS code with smallest field of definition \(\F_{p^s}\), then
\[
 N\le2p^s-2.
\]
In particular, \(s<f\) implies \(N\le2\sqrt q-2\).
If \(N>q+1\), one necessarily has
\[
 s=f,\qquad f\mid\ord_N(p),\qquad
 \ord_N(p)=f\,\ord_N(q).
\]
In the irreducible case, the last equality is \(\ord_N(p)=fK\), and one further requires \((f-1)K\le q-2\).
\end{corollary}
\begin{proof}
Scalar extension preserves dimension and distance, so the descended code is still nontrivial and MDS.
Apply Corollary~\ref{v84:cor-elementary-mds} over its actual coefficient field to obtain the first bound.
If \(s<f\), divisibility \(s\mid f\) gives \(2s\le f\), proving the second.
Thus \(N>q+1\) forces \(s=f\).
Now \(s\mid e\) and \(\ord_N(p^f)=e/\gcd(e,f)\) give the order identity.
In the irreducible case, \(fK=e\le\varphi(N)\le N-1\le q+K-2\); rearranging gives the dimension restriction.
\end{proof}

\begin{example}[Coefficient-field descent without changing coordinates]\label{v84:ex-descent-eleven}
Take \(q=9\) and \(N=11\).
Since \(\ord_{11}(3)=5\) and \(\ord_{11}(9)=5\), an irreducible cyclic code whose check root is a primitive eleventh root has dimension five and smallest field of definition \(\F_3\).
Were it MDS, the descended code would still have length eleven, whereas the elementary restriction would give \(N\le2\cdot3-2=4\), a contradiction.
No search through codewords is required.
Nor is the dimension of a general cyclic code being identified with \(\ord_N(q)\): that equality applies here because the check polynomial consists of one irreducible orbit.
\end{example}

\section{Common Permutation Symmetries of Cyclic Codes}
\label{v73:sec-common-cyclic-symmetry}
Root constraints determine the vectors in a cyclic code. We now ask the reverse question: which coordinate permutations preserve every cyclic code simultaneously?
This is a problem of recovering common symmetries from a family of linear subspaces.
Cyclotomic orbits, idempotent decompositions of quotient algebras, and the finite-field-valued Fourier transform provide the three descriptions needed to answer it.
Throughout this section, assume \(n\ge2\), \(\gcd(n,q)=1\), and put \(R=\Z/n\Z\), \(K=\langle q\rangle\le R^\times\). 
Theorems~\ref{v73:thm-common-permutation-orbits} and~\ref{v73:thm-common-permutation-arithmetic} 
correspond to Theorems~A and~B of Wu and Xiang \cite{WuXiang2026CyclicCodes}, respectively.
The proofs use finite-field-valued Fourier spectral projections and the cyclotomic difference-orbit rigidity established earlier.

\subsection{Recovering Difference Orbits from Spectral Projections}
\label{v73:subsec-code-projections}
\leavevmode\phantomsection\label{v73:def-common-permutation-group}\ignorespaces
A coordinate permutation \(\sigma\in\operatorname{Sym}(R)\) acts on \(\F_q^R\) by
\((\sigma c)_i=c_{\sigma^{-1}(i)}\). Define the \textbf{permutation automorphism group} of a code \(C\), and the \textbf{common permutation group} at length \(n\), by
\[
\operatorname{PAut}(C)=\{\sigma:\sigma C=C\},\qquad
\mathcal G_q(n)=\bigcap_C\operatorname{PAut}(C),
\]
where the intersection runs over all length-\(n\) cyclic codes over \(\F_q\), including the zero code and the full space.
Each code is preserved individually; merely permuting the family of cyclic codes is not enough.
Coordinate permutations preserve Hamming weight and differ from substitutions in the coefficient-polynomial variable, as in Proposition~\ref{v54:prop-affine-code-stability}.

Fix a primitive \(n\)th root of unity \(\zeta\) in a finite extension \(E/\F_q\).
With the defining-set notation, retain the Fourier-support notation from Section~\ref{v84:subsec-fourier-balance}:
\begin{equation}\label{v73:eq-code-spectrum}
S(C)=R\setminus T(C).
\end{equation}
This Fourier support, also called spectral support, records the Fourier evaluation coordinates that are allowed to be nonzero, and \(\dim C=|S(C)|\).
By Proposition~\ref{v40:prop-cyclic-root-constraints}, it is a union of \(q\)-cyclotomic cosets.
The defining set records vanishing coordinates, whereas the spectral support records the retained coordinates.

\begin{theorem}[Common permutation symmetries characterized by difference orbits]\label{v73:thm-common-permutation-orbits}
Under these assumptions,
\[
\boxed{\quad\mathcal G_q(n)=\mathcal R_{\langle q\rangle}(n).\quad}
\]
Thus a coordinate permutation \(\sigma\) preserves every cyclic code if and only if
\[
\sigma(i)-\sigma(j)\in\langle q\rangle(i-j)\qquad(i,j\in R).
\]
\end{theorem}
\begin{proof}
The finite-field-valued Fourier transform defined on page~\pageref{v40:def:cyclic-fourier} and its inverse are
\[
\widehat c(a)=\sum_{j\in R}c_j\zeta^{aj},\qquad
c_j=\frac1n\sum_{a\in R}\widehat c(a)\zeta^{-aj}.
\]
The inverse follows from geometric sums of roots of unity, using \(n\in\F_q^\times\).
For each \(q\)-cyclotomic orbit \(O\subseteq R\), define an orbit sum and its \textbf{Fourier spectral projection} by
\begin{equation}\label{v73:eq-fourier-orbit-projector}
h_O(z)=\sum_{a\in O}\zeta^{-az},\qquad
(P_O)_{ij}=\frac1n h_O(i-j).
\end{equation}
Since \(qO=O\), we have \(h_O(z)^q=h_O(z)\). Hence \(h_O(z)\in\F_q\) and \(P_O\in M_n(\F_q)\).
In Fourier coordinates, the inverse formula shows that \(P_O\) retains the coordinates in \(O\) and sets all other coordinates to zero.
Therefore
\[
P_O^2=P_O,\qquad P_OP_{O'}=0\ (O\ne O'),\qquad\sum_OP_O=I.
\]
Its image \(C_O=P_O\F_q^R\) is the irreducible cyclic code with spectral support \(O\), and
\[
\F_q^R=\bigoplus_OC_O,\qquad C=\bigoplus_{O\subseteq S(C)}C_O.
\]
These are the earlier CRT field-component projections written explicitly in coefficient coordinates.

Use the same symbol \(\sigma\) for the permutation matrix.
If \(\sigma\) preserves every cyclic code, it preserves each \(C_O\) and the sum of all other components, hence commutes with the corresponding projection.
Conversely, commuting with every \(P_O\) preserves each image and therefore every cyclic code. Thus
\[
\sigma\in\mathcal G_q(n)
\quad\Longleftrightarrow\quad
\sigma P_O\sigma^{-1}=P_O\quad\text{for every }O.
\]
Comparing matrix entries gives the equivalent condition
\begin{equation}\label{v73:eq-orbit-sums-preserved}
h_O(\sigma(i)-\sigma(j))=h_O(i-j)\qquad\text{for all }i,j,O.
\end{equation}

It remains to show that the orbit sums separate the orbits.
Each \(h_O\) is constant on \(K\)-orbits, because \(h_O(qz)=h_O(z)\).
The indicators \(1_O\) are linearly independent, and \(n^{-1}h_O\) are their inverse Fourier transforms.
Thus the \(h_O\) are linearly independent over \(E\), and in particular over \(\F_q\).
The space of \(\F_q\)-valued functions on \(R\) that are constant on \(K\)-orbits has dimension equal to the number of orbits, so the \(h_O\) form a basis of that space.
It contains the indicator of each orbit, and consequently
\[
\bigl(h_O(z)=h_O(w)\text{ for every }O\bigr)
\quad\Longleftrightarrow\quad w\in Kz.
\]
Apply this equivalence to~\eqref{v73:eq-orbit-sums-preserved} to finish the proof.
\end{proof}

Proposition~\ref{v140:prop-common-symmetry} shows that this common group persists at every unramified precision. Additional permutations of a single code are controlled by the matrix valuations in Proposition~\ref{v140:prop-permutation-depth}.

No orbit length was divided out or assumed prime to the characteristic; the only required invertibility is \(\gcd(n,q)=1\).
The proof expresses the same idempotents as field-component projections, Fourier coordinate restrictions, and functions of difference orbits.
In this way, common invariance of linear subspaces becomes an orbit condition on pairs of coordinates.

\subsection{An Arithmetic Criterion for the Common Affine Group}
\label{v73:subsec-code-arithmetic}
For \(K=\langle q\rangle\), the group \(\operatorname{Aff}_K(n)\) is generated by the cyclic shift \(i\mapsto i+1\) and the Frobenius multiplier \(i\mapsto qi\).
The preceding theorem places it inside \(\mathcal G_q(n)\).
To determine whether further common symmetries occur, it remains to translate the kernel containments of Theorem~\ref{v73:thm-ring-rigidity} into multiplicative-order conditions.

\leavevmode\phantomsection\label{v73:def-multiplier-index}\ignorespaces
For \(d\mid n\), define the \textbf{multiplier index}
\[
\iota_q(d)=[(\Z/d\Z)^\times:\langle q\bmod d\rangle]
=\frac{\varphi(d)}{\ord_d(q)}.
\]
We take \(\ord_1(q)=\varphi(1)=1\), so \(\iota_q(1)=1\).
The symbol \(\iota_q\) distinguishes this index from \(I_q\), which counts irreducible polynomials.

\begin{theorem}[An arithmetic criterion for the common permutation group]\label{v73:thm-common-permutation-arithmetic}
Let \(n\ge2\) and \(\gcd(n,q)=1\). Then
\[
\mathcal G_q(n)\ne\operatorname{Aff}_{\langle q\rangle}(n)
\]
if and only if at least one of the following holds:
\begin{enumerate}[label=\textup{(E\arabic*)}]
\item There is a prime \(p\ge5\) with \(p\parallel n\) and \(\iota_q(n/p)=\iota_q(n)\).
\item There is an odd prime \(p\) with \(p^2\mid n\) and \(\iota_q(n/p)=\iota_q(n)\).
\item We have \(8\mid n\) and \(\iota_q(n/4)=\iota_q(n)\).
\end{enumerate}
If none holds, the common permutations are exactly
\[
i\longmapsto q^ji+b,
\qquad 0\le j<\ord_n(q),\quad b\in R,
\]
numbering \(n\ord_n(q)\).
\end{theorem}
\begin{proof}
For \(d\mid n\), reduction \(\pi_d:R^\times\to(\Z/d\Z)^\times\) is surjective.
Indeed, in each prime-power CRT component a unit lifts to a unit; components at primes absent from \(d\) can be chosen arbitrarily among units.
The image \(\pi_d(K)\) is generated by \(q\bmod d\). Counting orders gives
\[
|\Delta_n(d)|=\frac{\varphi(n)}{\varphi(d)},\qquad
|K\cap\Delta_n(d)|=\frac{\ord_n(q)}{\ord_d(q)}.
\]
Consequently,
\begin{equation}\label{v73:eq-kernel-index-identity}
\frac{|\Delta_n(d)|}{|K\cap\Delta_n(d)|}
=\frac{\iota_q(n)}{\iota_q(d)},\qquad
\Delta_n(d)\subseteq K\ \Longleftrightarrow\ \iota_q(d)=\iota_q(n).
\end{equation}
Take \(d=n/p\) or \(d=n/4\) as appropriate and apply Theorems~\ref{v73:thm-common-permutation-orbits} and~\ref{v73:thm-ring-rigidity}.
In the absence of exceptions, the image of zero determines \(b\), and the difference of the images of one and zero determines \(q^j\).
Thus the listed parameters are distinct, yielding the count.
\end{proof}

\begin{example}[Length seven and an extra symmetry modulo sixteen]\label{v73:ex-common-symmetry}
Take \(q=2\), \(n=7\). Since \(\ord_7(2)=3\), we have \(\iota_2(7)=2\ne\iota_2(1)=1\).
The common permutations of all binary length-seven cyclic codes are therefore \(i\mapsto2^ji+b\), numbering \(21\).
However, the repetition code spanned by the all-one vector is preserved by all \(7!\) coordinate permutations.
Thus \(\mathcal G_q(n)\) records symmetries shared by every code, not necessarily the full permutation automorphism group of any specified code.

Next take \(q=5\), \(n=16\). As \(\ord_{16}(5)=4\) and \(\ord_4(5)=1\),
\(\iota_5(16)=2=\iota_5(4)\), so the third condition holds.
Here \(\Delta_{16}(4)=\{1,5,9,13\}=\langle5\rangle\).
The map \(i\mapsto i+4(i\bmod2)\) from~\eqref{v73:eq-two-local-shear} is not affine, yet it preserves every length-sixteen cyclic code over \(\F_5\).
\end{example}

Root-set theory thus determines not only the construction and dimension of cyclic codes, but also symmetries common to the entire family.
Root orbits give field components, these give idempotent projections, and the matrix entries of the projections separate difference orbits.
Local rigidity then turns the orbit constraints into a common affine multiplier.
For repeated-root cyclic codes, \(n\) is not invertible in the base field, so these Fourier projections are unavailable.
For an individual code, the common-group theorem likewise does not replace an analysis of its full automorphism group.

\subsection{Multiplier Preservation and Fixed-Point Rigidity for One Code}\label{v84:subsec-multiplier-rigidity}
The preceding two subsections studied permutations preserving every cyclic code.
For a fixed code, the MDS condition constrains the multipliers that preserve its Fourier support.
The quantifiers are different: the full automorphism group of one code cannot generally be read off from the common group.
Continue to assume \(\gcd(N,q)=1\) and use the support \(S=S(C)\) and negative-exponent basis of Section~\ref{v84:subsec-fourier-balance}.

\leavevmode\phantomsection\label{v84:def-multiplier-stabilizer}\ignorespaces
Define the \textbf{multiplier stabilizer}
\[
 H(S)=\{u\in(\Z/N\Z)^\times:uS=S\}.
\]
It is a subgroup of the unit group and contains \(\langle q\rangle\).
If \((P_ux)_j=x_{uj}\), then \(P_uv_d=v_{ud}\).
The Fourier model and the intersection formula for scalar extension show that \(u\in H(S)\) if and only if \(P_u\) preserves \(C\).

\begin{theorem}[A nonidentity multiplier has at most two fixed points]\label{v84:thm-multiplier-fixed}
If \(C\) is a nontrivial cyclic MDS code, then every \(u\in H(S)\), \(u\ne1\), satisfies
\begin{equation}\label{v84:eq-multiplier-fixed}
 \#\{j\in\Z/N\Z:uj=j\}=\gcd(u-1,N)\le2.
\end{equation}
\end{theorem}
\begin{proof}
First suppose that \(u\) has prime order \(r\).
Its coordinate orbits have length one or \(r\).
Let \(\nu\) be the number of fixed points and \(b\) the number of other orbits, so \(N=\nu+rb\).
Since \(uS=S\), the support contains \(a\) of the fixed points and \(t\) of the orbits of length \(r\), giving \(K=a+rt\).
One has \(0<t<b\).
Indeed, if \(t=0\), every \(d\in S\) is fixed by \(u\) and hence divisible by \(N/\nu\).
Columns zero and \(\nu\) of the Fourier generator matrix are equal, with \(1\le\nu<N\), contradicting the column criterion for \(K\ge2\).
If \(t=b\), the dual support contains no orbit of length \(r\), and the same argument applied to \(C^\perp\) gives a contradiction.

Over the splitting field \(E\), put \(P=P_u\) and \(W=(P-I)C_E\).
On the space spanned by a Fourier orbit of length \(r\), the kernel of \(P-I\) consists of vectors with all coefficients equal and has dimension one; on a fixed basis vector the map is zero.
Thus \(\dim_EW=(r-1)t\).
On the coordinate side, every \(w=(P-I)c\) is zero at the fixed coordinates and satisfies \(\sum_{j\in O}w_j=0\) on each orbit of length \(r\), since \(P\) permutes that orbit.
These facts remain valid when the characteristic divides \(r\); no division by \(r\) has occurred.

Impose zeros on \(t-1\) whole coordinate orbits of length \(r\), and on \(r-2\) coordinates of one more orbit.
Each whole orbit contributes at most \(r-1\) independent conditions because its coordinate sum is zero.
The total rank is at most
\[
 (t-1)(r-1)+(r-2)=t(r-1)-1<\dim_EW.
\]
A nonzero \(w\in W\) therefore still satisfies these conditions, and
\[
 \operatorname{wt}_H(w)
 \le rb-\bigl(r(t-1)+r-2\bigr)=r(b-t)+2.
\]
On the other hand, \(W\subseteq C_E\) and the MDS property give
\[
 \operatorname{wt}_H(w)\ge N-K+1
 =\nu-a+r(b-t)+1.
\]
Hence \(\nu-a\le1\).
The dual support \((\Z/N\Z)\setminus(-S)\) contains \(\nu-a\) fixed points and \(b-t\) long orbits.
Applying the same argument to the dual gives \(a\le1\), so \(\nu\le2\).

For general order \(h>1\), choose a prime \(r\mid h\) and put \(v=u^{h/r}\).
Then \(v\in H(S)\) has order \(r\).
Every point fixed by \(u\) is fixed by \(v\), so there are at most two.
Finally, the congruence \((u-1)j\equiv0\pmod N\) has \(\gcd(u-1,N)\) solutions, as required.
\end{proof}

This fixed-point restriction appears in\cite[Lemma 3.1]{LuoCaoLiLi2026Tools}.
The subcode distance estimate here follows directly from the zero-sum relations on the orbits and rank--nullity.
The restriction therefore comes not from the multiplier action alone, but from its preservation of a linear space with strong column independence.

\begin{corollary}[Synchronization of orders on divisors]\label{v84:cor-multiplier-orders}
Under these hypotheses, for every \(d\mid N\), \(d\ge3\), reduction induces an injection
\[
 H(S)\longrightarrow(\Z/d\Z)^\times.
\]
In particular, writing \(m=\ord_N(q)\), one has
\begin{equation}\label{v84:eq-general-divisor-orders}
 \ord_d(q)=m\qquad(d\mid N,\ d\ge3).
\end{equation}
For odd \(N\), all nonzero \(q\)-cyclotomic cosets have size \(m\).
For even \(N\), all cosets other than \(\{0\}\) and \(\{N/2\}\) have size \(m\).
The dimension \(K\) of a general cyclic code may be the sum of the sizes of several such orbits and need not equal \(m\).
\end{corollary}
\begin{proof}
If \(u\) belongs to the kernel of reduction, then \(d\mid u-1\), hence \(\gcd(u-1,N)\ge d\ge3\).
Theorem~\ref{v84:thm-multiplier-fixed} forces \(u=1\), proving injectivity.
Its restriction to \(\langle q\rangle\) preserves the order of \(q\), giving \eqref{v84:eq-general-divisor-orders}.
The coset of an exponent \(i\) has size \(\ord_{N/\gcd(N,i)}(q)\).
This modulus is less than three only for \(i=0\), or for even \(N\) and \(i=N/2\).
Apply synchronization of orders in every other case.
\end{proof}

The subgroup argument for an irreducible trace code gives order \(K\) for every \(d>1\), and also excludes even length.
The multiplier argument for a general cyclic code only requires \(d\ge3\), and does not identify \(K\) with \(m\).
The two rigidity statements therefore share a Frobenius-period mechanism but have different hypotheses and boundary cases.
In particular, the restriction \(d\ge3\) in \eqref{v84:eq-general-divisor-orders} is essential.

\begin{corollary}[Two fixed directions and parity of dimension]\label{v84:cor-even-length-parity}
Let \(C\) be a nontrivial \([N,K]\) cyclic MDS code, with \(\gcd(N,q)=1\). 
If \(\gcd(N,q^2-1)=2\) and \(m=\ord_N(q)\), then
\[
 K\equiv N-K\equiv1\pmod{2m},\qquad
 \min\{K,N-K\}\ge2m+1.
\]
\end{corollary}
\begin{proof}
Here \(N\) is even and \(q\) is odd, and \(N/2\) is odd: otherwise \(4\mid N\) and \(4\mid q^2-1\) would contradict the hypothesis.
Also \(m>2\), since otherwise \(N\mid q^2-1\) would force \(N=2\), contrary to nontriviality.
Choose a nonidentity element \(u\) of prime order in \(\langle q\rangle\).
It fixes \(0,N/2\), so it has exactly two fixed points.
The prime-order proof of Theorem~\ref{v84:thm-multiplier-fixed} also shows that \(S\) contains exactly one of them.

Every other \(q\)-orbit has size \(m\), and since \(q\) is odd it preserves exponent parity.
Suppose that \(S\) contains \(x\) even-exponent orbits and \(y\) odd-exponent orbits.
Put \(\varepsilon=1\) if \(0\in S\), and \(\varepsilon=0\) otherwise.
Balance modulo two gives
\[
 |mx+\varepsilon-(my+1-\varepsilon)|
 =|m(x-y)+2\varepsilon-1|\le1.
\]
Since \(m>2\), one has \(x=y\), and hence \(K=1+2mx\).
Nontriviality gives \(x\ge1\).
The dual code satisfies the same hypotheses, so the same conclusion applies to its dimension.
\end{proof}

Root-set descent determines which frequencies combine to form a code over the base field.
Residue layers constrain their distribution, and multiplier preservation then constrains the common periods of those layers.
For a single irreducible orbit, these conditions also return to full-degree multiplicative subgroups, normality, and the obstruction from sparse multiples.
They arise from the same Frobenius cyclic structure, tested through coefficients, orbits, and linear relations, respectively.

\par\Needspace{6\baselineskip}
\section*{Exercises}
\phantomsection
\addcontentsline{toc}{section}{Exercises}

\begin{exercise}[A constacyclic code without a base-field twisting root]\label{v88:ex-constacyclic-bch-five}
Over \(k=\F_5\), take \(N=4\) and \(\lambda=-1\). Let \(C\) be the
constacyclic code generated by \(f=X^2+2\) in \(k[X]/(X^4+1)\).
Show that \(X^4+1=(X^2+2)(X^2+3)\), while \(\lambda\) has no fourth
root in \(k\). Apply Corollary~\ref{v88:cor-constacyclic-bch-bound}
to give \(d_H(C)\ge2\), and use the weight of \(f\) to prove equality.
\end{exercise}

\begin{exercise}[Recovering a derivative fiber by interpolation on prescribed roots]\label{v83:ex-interpolated-split-fiber}
Let \(K=\F_Q\), where \(Q\) is a power of two, let \(d\mid Q-1\), and put \(h=(Q-1)/d\) and \(H=(K^\times)^d\).
Fix a monic polynomial \(U\in K[Y]\) of degree \(s\), with \(U(0)\ne0\), and write
\(P_V=U^2+YV^2\), where \(\deg V<s\).
\begin{enumerate}
\item For \(a\in H\) with \(U(a)\ne0\), prove that
\[
 P_V(a)=0\quad\Longleftrightarrow\quad V(a)=U(a)/\sqrt a.
\]
\item Specify \(s\) distinct elements \(a\in H\) at which \(U(a)\ne0\).
Prove that exactly one polynomial \(V\) of degree less than \(s\) makes them all roots of \(P_V\).
For their set \(A\), put \(L_A(Y)=\prod_{a\in A}(Y-a)\) and write \(P_V=L_AR_A\).
Give a necessary and sufficient condition for \(P_V\mid Y^h-1\).
\end{enumerate}
\end{exercise}

\begin{exercise}[Counting split multiplicative-coset fibers over a fixed field]\label{v83:ex-sixteen-coset-fibers}
Take \(K=\F_{16}\), \(d=3\), and \(H=(K^\times)^3\). For \(u\in K^\times\) and \(v\in K\), define
\[
 U(Y)=Y+u,\qquad
 P_v(Y)=(Y+u)^2+Yv^2,\qquad F_v(X)=XP_v(X^3).
\]
Prove that the number of parameters \(v\) satisfying \(F_v\mid X^{16}-X\) is two when \(u\in H\) and zero when \(u\notin H\).
When \(u=1\), how many polynomials \(F_v\) are square-free? Find the length and minimum distance of the full-support Goppa code defined by \(G=X^3+1\).
\end{exercise}

\begin{exercise}[The largest interleaving parameter]\label{v44:ex-interleaving-intrinsic}
Let \(C=(f)\subseteq k[X]/(X^N-1)\), where \(f\) is monic, \(f\mid X^N-1\), and \(1\le\deg f<N\).
For \(d\mid N\), let \(V_j^{(d)}\) be the coordinate subspace supported only at exponents congruent to \(j\pmod d\),
and let \(\pi_j^{(d)}\) be the coordinate projection onto this subspace. Prove that the following conditions are equivalent:
\begin{enumerate}
\item \(f\in k[X^d]\);
\item \(\pi_j^{(d)}(C)\subseteq C\) holds for \(0\le j<d\);
\item \(C=\bigoplus_{j=0}^{d-1}(C\cap V_j^{(d)})\).
\end{enumerate}
Deduce that the largest \(d\) satisfying these conditions is \(\operatorname{sp}(f)\).
\end{exercise}

\begin{exercise}[Coordinates vanishing throughout a subcode]\label{v54:ex-effective-support}
Over \(\F_2\), take \(U=X+1\) and \(f=U^6+U^3+1\).
Prove \(f\mid X^{63}-1\), and put \(C=(f)\subseteq\F_2[X]/(X^{63}-1)\) and
\(W=\operatorname{span}\{1,U^3,U^6\}\).
Use the outer remainder matrix to compute the dimension and effective support of \(C\cap W\).
Explain why the quadratic term allowed in \(W\) vanishes throughout this subcode, and specify the scope of the resulting upper bound on distance.
\end{exercise}

\begin{exercise}[Independent multipliers versus local projections]\label{v73:exer-independent-multipliers}
Take \(K=\langle2\rangle\le(\Z/15\Z)^\times\).
\begin{enumerate}
\item List the CRT coordinates of \(K\). Find its projections modulo three and five, and the two independent local multiplier groups in~\eqref{v73:eq-independent-local-group}.
\item Determine all functions satisfying \(f(x)-f(y)\in K(x-y)\) and count them. Explain why the full projection onto the unit group modulo five does not allow arbitrary local permutations.
\item Determine the common permutation group of all binary cyclic codes of length fifteen.
\end{enumerate}
\end{exercise}

\begin{exercise}[Preserving each code versus permuting the family]\label{v73:exer-family-versus-common}
Let \(\gcd(n,q)=1\), \(u\in(\Z/n\Z)^\times\), and \(\mu_u(i)=ui\).
\begin{enumerate}
\item Show that \(\mu_u\) sends cyclic codes to cyclic codes and that their spectral supports satisfy \(S(\mu_u C)=u^{-1}S(C)\).
\item Take \(q=2\), \(n=7\), and \(u=3\). Show that \(\mu_3\) permutes the family of all cyclic codes but does not belong to their common permutation group.
Identify the spectral supports of the two irreducible cyclic codes that it exchanges.
\end{enumerate}
\end{exercise}

\begin{exercise}[Cyclic covering and full-weight codewords]\label{v82:ex-cyclic-covering}
Let \(V=\F_q^N\), with \(N\ge1\), and let \(R_N\) be the cyclic right shift.
\leavevmode\phantomsection\label{v82:def-cyclically-covering}\ignorespaces
A subspace \(U\subseteq V\) is \textbf{cyclically covering} if \(\bigcup_{i=0}^{N-1}R_N^iU=V\).
\leavevmode\phantomsection\label{v82:def-full-weight-codeword}\ignorespaces
A word of weight \(N\), with all coordinates nonzero, is called a \textbf{full-weight codeword}.
Prove that a proper cyclically covering subspace exists if and only if a nonzero cyclic code without a full-weight codeword exists.
Do not assume \(\gcd(N,q)=1\).
\emph{Hint:} In one direction, enlarge a proper subspace to \(w^\perp\) and consider
\(x\mapsto(\langle x,R_N^iw\rangle)_{i=0}^{N-1}\); in the other, use the generator polynomial of a cyclic code.
This equivalence appears in~\cite{LiYuanCyclicCovering2026}.
\end{exercise}

\begin{exercise}[An explicit interpolation basis on the projective line]\label{v84:ex-projective-basis}
For \(a\in\F_q\), put
\[
 F_a(U,V)=V^q-V(U-aV)^{q-1},\qquad
 F_\infty(U,V)=U^q-UV^{q-1}.
\]
Prove that the evaluation vectors of these \(q+1\) forms of degree \(q\), at the representatives \((b,1)\) and \((1,0)\), form the standard basis.
Write the unique homogeneous interpolant of degree \(q\) for arbitrary prescribed projective point values, and explain why \(F_\infty\) is not the zero form.
\end{exercise}

\begin{exercise}[Linear-form identities and coordinate invariance of the critical kernel]\label{v84:ex-ternary-kernel}
Let \(q=p^e\), \(R=\F_q[X_1,\ldots,X_r]\), and \(r\ge2\). For \(\F_q\)-linear forms
\(L=\sum_i a_iX_i\) and \(M=\sum_i b_iX_i\), prove that
\[
 L^qM-LM^q
 =\sum_{1\le i<j\le r}(a_ib_j-a_jb_i)(X_i^qX_j-X_iX_j^q).
\]
Deduce that every \(\F_q\)-linear change of coordinates preserves the critical kernel of Theorem~\ref{v86:thm-frobenius-kernel}, and explain why linearly independent \(L,M\) yield a nonzero element of that kernel.
Finally, take \(q=5\), \(r=3\), \(L=X_1+2X_2\), and \(M=X_2+X_3\). Write out the resulting form of degree six and explain why it is neither a fifth power nor the zero polynomial.
\end{exercise}

\begin{exercise}[Two-dimensional irreducible trace codes and the norm-one subgroup]\label{v84:ex-dimension-two}
Let \(\alpha\in\F_{q^2}^\times\), with \(\ord(\alpha)=N>1\) and \(\ord_N(q)=2\).
Prove that \(C_\alpha\) is MDS if and only if \(\gcd(N,q-1)=1\), and that the possible lengths are exactly the odd divisors greater than one of \(q+1\).
State the additional restriction on \(N\) imposed by nontriviality.
Then take \(q=4,N=5\), prove that \(\omega=\alpha+\alpha^{-1}\) satisfies \(\omega^2+\omega+1=0\), find the check polynomial and code parameters, and verify normality of the nonidentity subgroup elements.
\end{exercise}

\begin{exercise}[Detecting low-weight multiples]\label{v84:ex-low-degree-multiplier}
Let \(m_\alpha(X)=\sum_{i=0}^K a_iX^i\), under the hypotheses of Theorem~\ref{v84:thm-sparse-multiple}. For \(0\le s\le N-K-1\), give a finite criterion involving only these coefficients for the existence of \(0\ne B\), with \(\deg B\le s\), such that \(m_\alpha B\) has weight at most \(K\). Deduce a characterization of when \(C_\alpha\) is MDS.
In Example~\ref{v84:ex-full-support-not-mds}, find the least possible degree of such a \(B\) and give a corresponding multiple.
\end{exercise}

\begin{exercise}[A distance obstruction with incomplete fibers]\label{v142:ex-partial-fibers}
Choose a primitive twenty-first root $\xi$ in $\F_{64}$, let $h=m_\xi\in\F_2[X]$, and let $C_\xi$ be the length-twenty-one cyclic code with check polynomial $h$.
Determine whether $h$ can be written as $G(X^3)$, and prove that $C_\xi$ is not MDS without listing all its words. Find the dimensions of the image and kernel of restriction to $\{0,3,6,\ldots,18\}$, and compare the resulting distance upper bounds for the code and its dual.
\end{exercise}

\begin{exercise}[A trace code from an element of order fifteen]\label{v84:ex-fifteen-obstruction}
Choose \(\alpha\in\F_{16}\) of order fifteen and consider the trace code \(C_\alpha\) over \(\F_2\). Find its dimension and Fourier support, decide whether it is MDS, and prove your conclusion.
\end{exercise}

\begin{exercise}[A code after scalar extension]\label{v84:ex-extension-splitting}
Take the code \(C\) with \(q=4,N=5\) from Exercise~\ref{v84:ex-dimension-two}, and extend scalars to \(\F_{16}\), obtaining \(C_{16}\). Find its parameters and decide whether it is an irreducible cyclic code. Using the same primitive fifth root for evaluation, find its Fourier support and its smallest field of definition in the fixed coordinates.
\end{exercise}

\begin{exercise}[Orders on divisors and short arithmetic restrictions]\label{v84:ex-arithmetic-orders}
Let \(C\) be a nontrivial simple-root cyclic MDS code over \(\F_q\), of length \(N\) and dimension \(K\).
\begin{enumerate}
\item Prove that \(N\mid q^2-1\) implies \(N\le q+1\).
\item Prove the same conclusion when \(3\mid N\) or \(4\mid N\), and when \(N\) has distinct odd prime divisors \(\ell_1,\ell_2\) with \(\gcd(\ell_1-1,\ell_2-1)=2\).
\item Prove that \(N\le q+1\) also holds if \(\gcd(N,q^2-1)>2\), or if \(K\) is even and \(\gcd(N,q^2-1)>1\).
\end{enumerate}
Use only the multiplier fixed-point restriction, synchronization of orders, and elementary parameter bounds proved in this chapter.
\end{exercise}

\begin{exercise}[Two boundary cases for homogeneous recovery]\label{v84:ex-gluing-boundaries}
Explain why uniqueness in Proposition~\ref{v84:prop-gluing} fails at \(d=N\).
Over any field, then take \(H_1=\{X=0\}\), \(H_2=\{Y=0\}\), and \(H_3=\{X-Y=0\}\).
Set \(g_1=g_2=0\) and \(g_3=X|_{H_3}\).
Prove that these homogeneous linear forms agree on every pairwise intersection but that no ternary homogeneous linear form has all three restrictions.
Identify the hypothesis of the proposition that fails.
\end{exercise}

\begin{exercise}[Five-point completion and a residual cubic]\label{v116:ex-five-completion}
Let \(q\) be a power of two and \(|A|=5\), using the notation of Proposition~\ref{v116:prop-completion-gap}. Prove \(3\le\deg H_A\le4\), and show that \(c_q(A)=6\) if and only if \(\deg H_A=3\), whereas \(c_q(A)=8\) if and only if \(\deg H_A=4\) and \(F_0/R_A\) splits completely.
Express the latter condition as a greatest-common-divisor degree test. For \(q=16\), prove that \(c_{16}(A)\in\{6,8,16\}\) and give a criterion for each value.
\end{exercise}

\chapter{Fourier Duality, Quadratic Forms, and Character Sums}
\label{chap:characters}
The trace pairing makes field elements into parameters for additive characters. The Fourier transform describes changes of frequency parameters through adjoint operators and converts solvability into orthogonality.
Polynomial indicator functions first give root-count congruences for low-degree equations; polarization and linearized root spaces determine the structure of quadratic forms, and character sums then determine their exact value distributions.
Joint expansions in additive and multiplicative characters count primitive elements of prescribed trace and, combined with the local Frobenius decomposition, prove the primitive normal basis theorem; weighted orbit products organize extension-field character sums into finite-order recurrences, and trace fibers recover extension-field point counts for a family of curves in characteristic two.
The same orthogonality relation also converts moments of character sums into counts of solutions of power-sum equations, and correlations of orbit sequences into character sums.

\section{Additive and Multiplicative Characters}
\leavevmode\phantomsection\label{def:characters}\ignorespaces Throughout this chapter write \(K=\F_Q\) and \(Q=p^m\), and put \(\zeta_p=\exp(2\pi i/p)\).
An \textbf{additive character} is a map \(K\to\C^\times\) satisfying \(\psi(x+y)=\psi(x)\psi(y)\).
A \textbf{multiplicative character} is a group homomorphism \(\chi:K^\times\to\C^\times\).
A character identically equal to one is called trivial.

\par\Needspace{7\baselineskip}
\noindent\phantomsection\label{v102:par-four-character-table}\textbf{Writing trace duality as four rows of signs.}\quad
Temporarily take \(K=\F_4=\F_2(\omega)\), with \(\omega^2+\omega+1=0\).
The absolute trace is \(\Tr(x)=x+x^2\), so
\[
\Tr(0)=\Tr(1)=0,\qquad\Tr(\omega)=\Tr(\omega^2)=1.
\]
For \(a\in K\), first compute the linear functional \(x\mapsto\Tr(ax)\), then replace \(0,1\) by the complex numbers \(+1,-1\), respectively. This gives \(\psi_a(x)=(-1)^{\Tr(ax)}\):

\[
\begin{array}{c|rrrr}
\psi_a(x)&x=0&x=1&x=\omega&x=\omega^2\\\hline
a=0&1&1&1&1\\
a=1&1&1&-1&-1\\
a=\omega&1&-1&-1&1\\
a=\omega^2&1&-1&1&-1
\end{array}
\]
For instance, \(\psi_\omega(\omega)=(-1)^{\Tr(\omega^2)}=-1\), whereas \(\psi_\omega(\omega^2)=(-1)^{\Tr(1)}=1\).
The four rows are distinct. Each nontrivial row has two plus signs and two minus signs, and the entrywise product of any two distinct rows sums to zero.
This displays both the trace parametrization of characters and their orthogonality. Arithmetic with the table entries takes place in \(\C\); its \(-1\) must not be identified with \(1\) in characteristic two.

\begin{theorem}\label{thm:characters}
All additive characters of \(K\) are exactly
\[
\psi_a(x)=\zeta_p^{\Tr_{K/\F_p}(ax)},\qquad a\in K.
\]
Fixing a generator \(g\) of \(K^\times\), all multiplicative characters are uniquely determined by \(\chi_j(g)=\exp(2\pi i j/(Q-1))\), \(0\leq j<Q-1\).
\end{theorem}
\begin{proof}
The homomorphism property of \(\psi_a\) follows from additivity of the trace; nondegeneracy of the trace pairing ensures that distinct \(a\) give distinct characters.
Every additive character takes values among the \(p\)th roots of unity and therefore has a unique expression \(\zeta_p^{\ell(x)}\),
where \(\ell:K\to\F_p\) is an \(\F_p\)-linear functional. Chapter~\ref{chap:dual} shows that every such functional has the form \(x\mapsto\Tr(ax)\).
The multiplicative statement follows directly because the image of a cyclic-group generator can be any \((Q-1)\)th root of unity.
\end{proof}

\begin{proposition}[Orthogonality relations]\leavevmode\phantomsection\label{prop:character-orthogonality}
A nontrivial additive character satisfies \(\sum_{x\in K}\psi(x)=0\), and a nontrivial multiplicative character satisfies \(\sum_{x\in K^\times}\chi(x)=0\). Moreover,
\[
\sum_{a\in K}\psi_a(x)=
\begin{cases}Q,&x=0,\\0,&x\ne0,\end{cases}
\qquad
\sum_{j=0}^{Q-2}\chi_j(x)=
\begin{cases}Q-1,&x=1,\\0,&x\in K^\times\setminus\{1\}.\end{cases}
\]
\end{proposition}
\begin{proof}
If \(\psi(b)\ne1\), translating the summation variable by \(b\) leaves the sum unchanged while multiplying it by \(\psi(b)\), so the sum must be zero.
For the multiplicative case, multiply the variable by a suitable \(b\). The formulas summing over characters follow from the trace pairing and finite geometric series, respectively.
\end{proof}

Passing from a polynomial function to its trace can identify further parameters. The relevant coefficient space is determined by the field generated by the power values; this also gives nonredundant parameters for the corresponding additive-character phases \cite[Lemma~4.1 and Proposition~4.2]{ManganielloRandrianarivo2026}.

\begin{proposition}[Parameters of trace evaluation]\label{v68:prop-trace-evaluation-parameters}
Let \(K=\F_{p^m}\), let \(g\) generate \(K^\times\), and let \(1\le j\le p^m-1\). Put
\[
K_j=\F_p(g^j),\qquad d_j=[K_j:\F_p],
\]
and define the \(\F_p\)-linear map
\[
T_j:K\longrightarrow\{\text{functions }K\to\F_p\},\qquad
T_j(a)(x)=\Tr_{K/\F_p}(ax^j).
\]
Then
\[
\ker T_j=\ker\Tr_{K/K_j},\qquad
\dim_{\F_p}\im T_j=d_j.
\]
If \(u\in K\) satisfies \(\Tr_{K/K_j}(u)=1\), every function in \(\im T_j\) has a unique expression \(T_j(ub)\) with \(b\in K_j\).
\end{proposition}
\begin{proof}
For every \(x\in K\), we have \(x^j\in K_j\). Trace transitivity therefore gives
\[
T_j(a)(x)=\Tr_{K_j/\F_p}\bigl(x^j\Tr_{K/K_j}(a)\bigr).
\]
Thus \(\ker\Tr_{K/K_j}\subseteq\ker T_j\).
Conversely, if \(T_j(a)=0\), the element \(\Tr_{K/K_j}(a)\) is trace-orthogonal to every \(x^j\).
The values \(x^j\) include all powers of \(g^j\), which span \(K_j\) over \(\F_p\). Nondegeneracy of the trace pairing on \(K_j\) gives \(\Tr_{K/K_j}(a)=0\).
The relative trace is surjective, so rank--nullity gives \(\dim_{\F_p}\im T_j=d_j\).
It also supplies such an element \(u\). Since \(\Tr_{K/K_j}(ub)=b\), these coefficients give each function exactly once.
\end{proof}

The integer \(d_j\) is the length of the \(p\)-Frobenius orbit of \(g^j\), equivalently the size of the \(p\)-cyclotomic coset of \(j\) modulo \(p^m-1\).
Reading the values at \(x=g^r\) gives the trace sequences of Section~\ref{v40:sec-cyclic-trace}; passing to \(K_j\) removes the redundancy caused by using a larger coefficient field.
The hypothesis \(j\ge1\) keeps the value at zero fixed at zero: in particular, \(X^{p^m-1}\) and the constant polynomial one must remain distinct as functions on \(K\).

\begin{example}[A nonzero polynomial function with zero trace]\label{v68:ex-trace-x5}
Take \(K=\F_{16}\) and \(j=5\). Then \(K_j=\F_4\), and \(x^5\in\F_4\) for every \(x\in\F_{16}\). Consequently,
\[
\Tr_{\F_{16}/\F_2}(ax^5)
=\Tr_{\F_4/\F_2}\bigl((a+a^4)x^5\bigr).
\]
For \(a=1\), this is identically zero, although \(X^5\) is a nonzero polynomial function and its exponent is odd. Its additive-character phase \(x\mapsto\zeta_2^{\Tr_{\F_{16}/\F_2}(x^5)}\) is therefore identically one.
As \(a\) ranges over \(\F_{16}\), there are exactly four distinct trace functions, each represented by four coefficients; the kernel is \(\F_4\).
Thus excluding exponents divisible by the characteristic does not by itself make trace-evaluation parameters unique.
Theorem~\ref{v70:thm-phase-reduction} will give a phase reduction valid over every extension; this example shows the additional degeneracy that can arise from trace evaluation on one fixed field.
\end{example}

\section{Fourier Transforms and Root-Count Formulas}
\leavevmode\phantomsection\label{def:fourier-transform}\ignorespaces Fix \(\psi(x)=\zeta_p^{\Tr_{K/\F_p}(x)}\). For any function \(F:K\to\C\), define its \textbf{additive Fourier transform} by
\[
\widehat F(a)=\sum_{x\in K}F(x)\psi(-ax).
\]
\leavevmode\phantomsection\label{def:fourier-support}\ignorespaces The parameter \(a\) is called a frequency, and the set \(\operatorname{supp}\widehat F=\{a\in K:\widehat F(a)\ne0\}\) is the \textbf{Fourier support} of \(F\).

\par
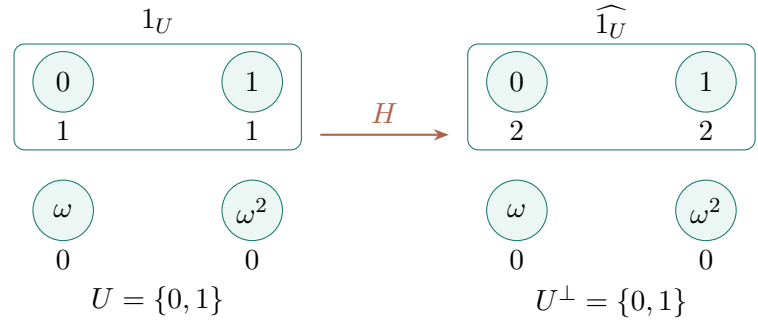
\begin{figure}[!htbp]
\centering
\begin{tikzpicture}[
  vpoint/.style={circle,draw=Jade,fill=CoverMint,minimum size=8mm,inner sep=2pt},
  varrow/.style={-{Stealth[length=2mm]},Jade,semithick},
  vbox/.style={draw=MutedSlate,rounded corners=2pt,inner sep=5pt,align=center},
  every node/.style={font=\small}]

\node at (1.25,1.5) {\(1_U\)};\node at (7.25,1.5) {\(\widehat{1_U}\)};
\foreach \name/\x/\y/\lab/\val in {a/0/.7/0/1,b/2.5/.7/1/1,c/0/-1/\omega/0,d/2.5/-1/\omega^2/0} {
  \node[vpoint] (\name) at (\x,\y) {\(\lab\)};
  \node at (\x,\y-.65) {\(\val\)};
}
\foreach \name/\x/\y/\lab/\val in {e/6/.7/0/2,f/8.5/.7/1/2,g/6/-1/\omega/0,h/8.5/-1/\omega^2/0} {
  \node[vpoint] (\name) at (\x,\y) {\(\lab\)};
  \node at (\x,\y-.65) {\(\val\)};
}
\draw[Jade,rounded corners] (-.65,-.2) rectangle (3.15,1.2);
\draw[Jade,rounded corners] (5.35,-.2) rectangle (9.15,1.2);
\draw[varrow,Clay,thick] (3.4,0)--node[above] {\(H\)} (5.1,0);
\node at (1.25,-2.2) {\(U=\{0,1\}\)};
\node at (7.25,-2.2) {\(U^\perp=\{0,1\}\)};

\end{tikzpicture}
\caption{The left points are field elements; the right points are frequencies labelled by the trace pairing. Boxes mark the two supports. The value two comes from adding two equal signs.}\label{v102:fig-fourier-four}
\end{figure}

\Needspace{7\baselineskip}
\noindent\phantomsection\label{v102:par-fourier-four}\textbf{Fourier transforms of the two trace fibers.}\quad
Return to the four-element example, ordering both function values and frequencies as \((0,1,\omega,\omega^2)\).
Let \(H=(\psi_a(x))_{a,x}\) be the preceding table. Since \(-ax=ax\) in characteristic two, \(\widehat F=HF\).
Put \(U=\ker\Tr=\{0,1\}\), and let \(1_U\) be its indicator. Summing the appropriate entries in each row gives

\[
H\begin{pmatrix}1\\1\\0\\0\end{pmatrix}
=\begin{pmatrix}2\\2\\0\\0\end{pmatrix},\qquad
H\begin{pmatrix}0\\0\\1\\1\end{pmatrix}
=\begin{pmatrix}2\\-2\\0\\0\end{pmatrix}.
\]
The frequencies with nonzero Fourier coefficients in the first identity form precisely
\[
U^\perp=\{a:\Tr(au)=0\text{ for all }u\in U\}=\{0,1\}.
\]
The second identity concerns the other trace fiber, \(\omega+U=\{\omega,\omega^2\}\).
Translation preserves the support and changes the sign at frequency one.
Figure~\ref{v102:fig-fourier-four} labels each point with its function or transform value, including zeros.

\FloatBarrier
Row orthogonality gives \(H^2=4I\), so inversion applies \(H\) once more and divides by four.
The equality \(U=U^\perp\) is special to this example. In general, a subspace indicator satisfies
\(\widehat{1_U}=|U|1_{U^\perp}\): summing a trivial character over \(U\) gives \(|U|\), while a nontrivial character cancels.
Returning to general \(K=\F_Q\), the same cancellation gives the following inversion formula.

Orthogonality gives Fourier inversion and Parseval's identity:
\begin{equation}\label{eq:fourier}
F(x)=\frac1Q\sum_{a\in K}\widehat F(a)\psi(ax),
\qquad
\sum_x|F(x)|^2=\frac1Q\sum_a|\widehat F(a)|^2.
\end{equation}
Substituting the definition of the transform into the right-hand side of the inversion formula gives
\[
\frac1Q\sum_a\sum_yF(y)\psi(a(x-y))
=\sum_yF(y)\,\frac1Q\sum_a\psi(a(x-y))=F(x).
\]
The inner sum is \(Q\) when \(x=y\), and zero otherwise.
Similarly, since \(\overline{\psi(t)}=\psi(-t)\),
\[
\sum_a|\widehat F(a)|^2
=\sum_{x,y}F(x)\overline{F(y)}
  \sum_a\psi(a(y-x))
=Q\sum_x|F(x)|^2,
\]
which proves Parseval's identity.

If \(P:K^r\to K\) is a polynomial function, the number of solutions of \(P(x)=b\) is
\begin{equation}\label{eq:character-root-count}
N_P(b)=\frac1Q\sum_{a\in K}\psi(-ab)
\sum_{x\in K^r}\psi(aP(x)).
\end{equation}
For each \(x\), the corresponding indicator equals one when \(P(x)=b\) and zero otherwise, so this is an exact equality.
The point indicators of Chapter~\ref{chap:functions} take values in a field of characteristic \(p\), whereas the character orthogonality formulas here take values in \(\C\).
Section~\ref{sec:chevalley-warning} uses the former to determine root-count congruences modulo \(p\); the latter retain integer solution counts and give exact enumerations.

\subsection{The Fourier support uncertainty principle}\label{v90:subsec-fourier-uncertainty}
\begin{theorem}[Support uncertainty on a finite abelian group]\label{v90:thm-fourier-uncertainty}
Let \(G\) be a finite abelian group, let \(\widehat G\) be its group of complex-valued characters, and let \(f:G\to\C\) be nonzero. Define the unnormalized Fourier transform by
\[
\widehat f(\gamma)=\sum_{x\in G}f(x)\overline{\gamma(x)}
\qquad(\gamma\in\widehat G).
\]
For the additive group of a finite field, the trace pairing parametrizes these characters and recovers the earlier definition. Then
\[
\boxed{|\operatorname{supp}f|\,|\operatorname{supp}\widehat f|\ge |G|.}
\]
\leavevmode\phantomsection\label{v121:par-uncertainty-equality}
Write \(G\) additively. Equality holds if and only if there are a subgroup \(H\le G\), an element \(x_0\in G\), a character \(\gamma_0\in\widehat G\), and \(c\in\C^\times\) such that
\[
f(x)=c\gamma_0(x)1_{x_0+H}(x).
\]
In this case, \(\operatorname{supp}\widehat f=\gamma_0H^\perp\), where \(H^\perp=\{\gamma\in\widehat G:\gamma|_H=1\}\).
\end{theorem}
\begin{proof}
Put $A=\operatorname{supp}f$, $B=\operatorname{supp}\widehat f$. By Cauchy--Schwarz,
\[
\|\widehat f\|_\infty\le\|f\|_1\le\sqrt{|A|}\,\|f\|_2.
\]
By character orthogonality and Parseval,
\[
|G|\,\|f\|_2^2=\|\widehat f\|_2^2\le|B|\,\|\widehat f\|_\infty^2
\le|A||B|\,\|f\|_2^2.
\]
Cancel the positive factor $\|f\|_2^2$.

Suppose equality holds. Every inequality above is then an equality. Fix \(\gamma_0\in B\). Equality in Cauchy--Schwarz makes \(|f|\) constant on \(A\); equality in the triangle inequality gives \(f(x)=c\gamma_0(x)\) for \(x\in A\), with \(c\ne0\).
For each \(\gamma\in B\), the same equality condition makes \(\gamma\gamma_0^{-1}\) constant on \(A\). Choose \(x_0\in A\) and let \(H\) be the subgroup generated by \(A-x_0\). Then
\[
A\subseteq x_0+H,\qquad B\subseteq\gamma_0H^\perp.
\]
Character orthogonality gives \(|H^\perp|=|G|/|H|\). Thus \(|A||B|=|G|\) forces equality in both inclusions.
Conversely, direct summation for the stated function gives
\[
\widehat f(\gamma)
=c\gamma_0(x_0)\overline{\gamma(x_0)}\,|H|\,
1_{\gamma_0H^\perp}(\gamma),
\]
so the two support sizes have product \(|G|\).
\end{proof}
\begin{remark}
For this equality characterization, see also \cite[Section~1]{BonamiGhobber2013}. For \(G=(K,+)\), it takes the form
\[
f(x)=c\psi(ax)1_{b+U}(x),\qquad
\operatorname{supp}\widehat f=a+U^\perp,
\]
where \(U\) is an \(\F_p\)-subspace of \(K\) and \(U^\perp\) is defined by the absolute trace pairing. Thus the preceding subspace and trace-fiber indicators attain equality.
Here the Fourier transform is complex valued and uses additive-character orthogonality. It must not be conflated with a transform in modular characteristic dividing the group order. The inequality complements the spectral projections and cyclic-code representations established earlier.
\end{remark}

\section{Moments of Character Sums and Power-Sum Equations}
\label{v58:sec-character-moments}

Character orthogonality counts roots of a single equation and also converts higher moments of a family of character sums into counts of common solutions of several equations.
Fiber collisions in value-set problems give the simplest case.

\begin{proposition}[A Fourier expression for the collision count]\leavevmode\phantomsection\label{v58:prop-collision-fourier}
Let \(F:K\to K\) and \(K=\F_Q\), and fix a nontrivial additive character \(\psi\).
Write \(C_F=\#\{(x,y)\in K^2:F(x)=F(y)\}\). Then
\begin{equation}\label{v58:eq-collision-fourier}
C_F=\frac1Q\sum_{a\in K}
\left|\sum_{x\in K}\psi(aF(x))\right|^2.
\end{equation}
\end{proposition}
\begin{proof}
Expanding the squared absolute value and interchanging sums makes the right-hand side equal to
\[
\sum_{x,y\in K}\frac1Q\sum_{a\in K}
\psi\bigl(a(F(x)-F(y))\bigr).
\]
By Proposition~\ref{prop:character-orthogonality}, the normalized inner sum is one when \(F(x)=F(y)\) and zero otherwise.
\end{proof}

Proposition~\ref{v58:prop-value-collision} thus relates the value-set lower bound to this sum of squares.
For example, the collision count of \(F=(X^3-X)^2\) in Example~\ref{v58:ex-F27-five-values} is 153,
so the unnormalized sum of squares on the right of~\eqref{v58:eq-collision-fourier} is \(27\cdot153=4131\).
The fiber sizes already determine this sum of squares, without evaluating the individual character sums.

\leavevmode\phantomsection\label{v58:def-character-moment}\ignorespaces
Let \(A\) be a finite set, with \(S:A\to\C\) and \(r\ge1\). We call \(\sum_{a\in A}|S(a)|^{2r}\) its \textbf{moment of order \(2r\)},
with the convention that the sum is not normalized by the number of parameters.
For a subset $S\subseteq K$, its \textbf{additive energy} is
\begin{equation}\label{v151:def-additive-energy}
 \mathcal E_+(S)=\#\{(x_1,x_2,x_3,x_4)\in S^4:x_1+x_2=x_3+x_4\}.
\end{equation}
Grouping these quadruples by $x_1-x_3=x_4-x_2$, or expanding the Fourier fourth moment, gives
\begin{equation}\label{v151:eq-energy-fourier}
 \mathcal E_+(S)=\sum_{\delta\in K}|S\cap(S+\delta)|^2
 =\frac1{|K|}\sum_{a\in K}\bigl|\widehat{1_S}(a)\bigr|^4.
\end{equation}
The second equality again follows from additive-character orthogonality: only quadruples satisfying $x_1+x_2=x_3+x_4$ survive the inner sum.

\begin{corollary}[Additive energy and a fourth moment for the norm-one group]\label{v151:cor-normone-moment}
Let $Q$ be a prime power, $K=\F_{Q^2}$, and $H=\ker\Nm_{K/\F_Q}$. Then
\begin{equation}\label{v151:eq-normone-energy}
 \mathcal E_+(H)=
 \begin{cases}
 3Q(Q+1),&Q\text{ odd},\\
 (Q+1)(3Q+1),&Q\text{ even}.
 \end{cases}
\end{equation}
Consequently, for the Fourier transform defined by any nontrivial additive character,
\[
 \sum_{a\in K^\times}\bigl|\widehat{1_H}(a)\bigr|^4
 =Q^2\mathcal E_+(H)-(Q+1)^4.
\]
\end{corollary}
\begin{proof}
The zero translation has $Q+1$ intersection points, and each nonzero norm fiber also has $Q+1$ elements. Proposition~\ref{v151:prop-normone-intersections} therefore reduces the energy to a sum over $D\in\F_Q^\times$.
For odd $Q$,
\[
 \sum_{D\ne0}\chi(D(D-4))=-1,\qquad
 \sum_{D\ne0}\chi(D(D-4))^2=Q-2.
\]
To prove the first identity, count all pairs $(D,z)\in\F_Q^2$ satisfying $z^2=D(D-4)$. The substitution
$(D-2-z)(D-2+z)=4$ gives exactly $Q-1$ solutions, so the total quadratic-character sum is $-1$. The term at $D=0$ is zero, proving the first identity. The second identity only removes the zero at $D=4$. Hence
\[
 \sum_{D\ne0}\bigl(1-\chi(D(D-4))\bigr)^2=2Q-1,
\]
and the energy is $(Q+1)^2+(Q+1)(2Q-1)=3Q(Q+1)$.
For even $Q$, exactly $Q/2$ nonzero values of $D$ satisfy $\Tr_{\F_Q/\F_2}(D^{-1})=1$; each gives intersection count two. The energy is therefore $(Q+1)^2+2Q(Q+1)$.
Finally, apply~\eqref{v151:eq-energy-fourier} and remove the zero-frequency contribution $(Q+1)^4$ to obtain the fourth-moment formula.
\end{proof}

For example, when $Q=4$ the energy is $65$, and the fourth moment over nonzero frequencies is $16\cdot65-5^4=415$. The trace--norm equation determines the intersections, and orthogonality then determines the moment without evaluating individual Fourier coefficients.

Return to a general field \(K=\F_Q\). Put \(d,r\ge1\). For \(a=(a_1,\ldots,a_d)\in K^d\), define
\[
S_d(a)=\sum_{t\in K}\psi(a_1t+\cdots+a_dt^d).
\]
Expanding \(S_d(a)^r\overline{S_d(a)}^{\,r}\) and applying orthogonality separately to each \(a_j\) gives
\begin{equation}\label{v58:eq-moment-power-systems}
\sum_{a\in K^d}|S_d(a)|^{2r}
=Q^d\#\left\{(x,y)\in K^r\times K^r:
\sum_{i=1}^r x_i^j=\sum_{i=1}^r y_i^j\ (1\le j\le d)\right\}.
\end{equation}
Thus exact moment computation reduces to fiber counts for power-sum equations; for related work in finite-field Fourier analysis, see \cite[Section 2.1]{BiswasMoment2025}.

We compare roots with multiplicity using the \hyperref[v58:def-multiset]{multisets} defined in Chapter~\ref{chap:root-symmetry}.

\begin{lemma}[Newton identities]\label{v58:lem-newton}
In a field, take \(x_1,\ldots,x_m\) and write
\[
\prod_{i=1}^m(X-x_i)=X^m+c_1X^{m-1}+\cdots+c_m,
\qquad s_j=\sum_{i=1}^m x_i^j.
\]
Then the \textbf{Newton identities} are
\begin{equation}\label{eq:root-power-newton}
s_j+\sum_{i=1}^{j-1}c_i s_{j-i}+jc_j=0
\qquad(1\le j\le m).
\end{equation}
\end{lemma}
\begin{proof}
Put \(H(u)=\prod_i(1-x_i u)=1+c_1u+\cdots+c_mu^m\). In the formal power series ring,
\[
-\frac{uH'(u)}{H(u)}
=\sum_i\frac{x_i u}{1-x_i u}
=\sum_{j\ge1}s_j u^j.
\]
Multiplication by \(H(u)\) and comparison of the coefficient of \(u^j\) gives the formula.
This derivation is valid in every characteristic and requires no division by \(j\).
\end{proof}

\begin{proposition}[Power sums determine the root multiset]\label{v58:prop-power-sum-multiset}
Suppose that \(p=\Char K>d\). If \((x_1,\ldots,x_d),(y_1,\ldots,y_d)\in K^d\) satisfy
\[
\sum_{i=1}^d x_i^j=\sum_{i=1}^d y_i^j\qquad(1\le j\le d),
\]
then the two lists of elements form the same multiset, counting multiplicities.
\end{proposition}
\begin{proof}
Write \(\prod_i(X-x_i)=X^d+c_1X^{d-1}+\cdots+c_d\) and put \(s_j=\sum_i x_i^j\).
By Lemma~\ref{v58:lem-newton},
\[
c_j=-\frac1j\left(s_j+\sum_{i=1}^{j-1}c_i s_{j-i}\right)
\qquad(1\le j\le d).
\]
Under the assumption \(p>d\), all \(1,\ldots,d\) are invertible,
so the first \(d\) power sums determine all coefficients successively and uniquely.
The two lists therefore give the same monic root polynomial, and unique factorization gives the same roots and multiplicities.
\end{proof}

The characteristic restriction is used to recover coefficients from power sums.
In \(\F_3\), the first three power sums of both \((0,0,0)\) and \((1,1,1)\) are zero,
but their root polynomials are \(X^3\) and \((X-1)^3=X^3-1\), respectively, and are different.

For distinct roots, Theorem~\ref{v123:thm-power-trace-pencil} requires no assumption that the characteristic exceeds the degree: the first \(2m-1\) power sums always recover a monic square-free polynomial of degree \(m\), while the second derivative and the affine fibers describe all ambiguity in the first \(2m-2\) terms. This conclusion does not apply directly to repeated roots.

\begin{example}[The fourth moment of quadratic-phase character sums]\label{v58:ex-quadratic-fourth-moment}
Let \(Q\) be odd, and put
\[
W(a,b)=\sum_{t\in K}\psi(at+bt^2)\qquad(a,b\in K).
\]
By~\eqref{v58:eq-moment-power-systems},
\[
\sum_{a,b\in K}|W(a,b)|^4
=Q^2\#\left\{(x,y,z,w)\in K^4:
\begin{aligned}
x+y&=z+w,\\[-2pt]
x^2+y^2&=z^2+w^2
\end{aligned}\right\}.
\]
Since two is invertible, the two equations also give \(xy=zw\), so the two pairs are roots of the same quadratic polynomial:
\(\{x,y\}=\{z,w\}\), with multiplicities counted.
Fix \((x,y)\). When \(x\ne y\), there are two permutations for \((z,w)\); when \(x=y\), there is only one.
There are \(2Q(Q-1)+Q=2Q^2-Q\) such quadruples. Therefore
\begin{equation}\label{v58:eq-quadratic-fourth-moment}
\boxed{\quad\sum_{a,b\in K}|W(a,b)|^4=2Q^4-Q^3.\quad}
\end{equation}
\end{example}

This identity follows by counting fibers of root multisets, without first computing the individual \(W(a,b)\).
Higher moments are still governed by~\eqref{v58:eq-moment-power-systems}. When \(r=d<p\),
each nonempty fiber consists of all permutations of one root multiset, its size determined by the root multiplicities.

\section{Fourier Transforms and Adjoints of Linear Operators}
\label{sec:fourier-adjoint}

The trace pairing recovers coordinates and also determines how linear maps act on frequencies. Let \(A:K\to K\) be an \(\F_p\)-linear map, and let \(A^*\) be its adjoint for the absolute trace pairing, so that \(\Tr(aA(x))=\Tr(A^*(a)x)\). If \(A\) is a linearized operator over an intermediate base field,
transitivity of the trace shows that its adjoint for the relative trace pairing is also its adjoint here.
The commutative relation \(\kappa A^*=A^\vee\kappa\) in Figure~\ref{fig:trace-adjoint} now applies over \(\F_p\):
composing the functional \(x\mapsto\Tr(ax)\) with \(A\) changes its parameter to \(A^*(a)\).
Composing these functionals with the character \(t\mapsto\zeta_p^t\) of \(\F_p\) gives the following formula for the action on Fourier frequencies.

\par\Needspace{7\baselineskip}
\noindent\phantomsection\label{v102:par-adjoint-four}\textbf{Computing the adjoint in the same four-element field.}\quad
Take \(A(x)=x^2+x=\Tr(x)\), regarding its values \(0,1\) as elements of \(K=\F_4\).
Here \(A^*=A\), since \(\Tr(ax^2)=\Tr(a^2x)\). Its nonempty fibers are
\[
A^{-1}(0)=\{0,1\},\qquad A^{-1}(1)=\{\omega,\omega^2\}.
\]
Write \(A_*1\) for the sum of the constant function one over each fiber. This gives \(A_*1=(2,2,0,0)^{\mathsf T}\), in the same ordering.
The table gives \(H(A_*1)=(4,4,0,0)^{\mathsf T}=4\,1_{\ker A^*}\).
Thus \(x^2+x=b\) has two solutions for \(b=0,1\), and none for \(b=\omega,\omega^2\).
Solvability is exactly \(b\in(\ker A^*)^\perp\). The next proposition extends this fiber count to an arbitrary linear operator.

\begin{proposition}\label{prop:fourier-adjoint}
For a function \(F:K\to\C\), define its fiber-sum function by
\[
(A_*F)(y)=\sum_{A(x)=y}F(x).
\]
Then
\[
\widehat{A_*F}(a)=\widehat F(A^*(a)).
\]
In particular, the number of solutions of \(A(x)=b\) is
\begin{equation}\label{eq:linear-fourier-count}
\#\{x:A(x)=b\}=\sum_{a\in\ker A^*}\psi(-ab)
=\begin{cases}
|\ker A|,&b\in(\ker A^*)^\perp,\\
0,&b\notin(\ker A^*)^\perp.
\end{cases}
\end{equation}
\end{proposition}
\begin{proof}
Summing over fibers and combining indices gives
\[
\widehat{A_*F}(a)
=\sum_{x\in K}F(x)\psi(-aA(x))
=\sum_{x\in K}F(x)\psi(-A^*(a)x).
\]
Take \(F\equiv1\); its Fourier transform is \(Q\) at zero and vanishes elsewhere.
Thus \(\widehat{A_*1}=Q\,1_{\ker A^*}\). Inversion gives \((A_*1)(b)=\sum_{a\in\ker A^*}\psi(ab)\); replacing \(a\) by \(-a\) gives the sign in the formula.
If \(b\) is orthogonal to \(\ker A^*\), each term is one; otherwise the sum is that of a nontrivial character on this additive subgroup and is zero.
Finally, \(A\) and \(A^*\) have the same rank, so their kernels have the same size.
\end{proof}

This again gives \(\im A=(\ker A^*)^\perp\). When \(A=\sigma-\id\),
it is the trace condition of additive Hilbert 90. The kernel of a linearized polynomial, a trace-orthogonal complement, and the support of a Fourier transform
are thus three descriptions of the same linear equation.

\subsection{The local spectrum of the normality indicator}
\label{v128:sec-normal-spectrum}

Counting solutions of a linear equation uses orthogonality for one subspace. Normality requires avoidance of several subspaces; the trace-dual direct sum from Chapter~\ref{chap:normal} gives their product an explicit spectrum.

\begin{proposition}[Fourier expansion of local generation conditions]\label{v128:prop-local-spectrum}
Let \(k=\F_q\), \(E=\F_{q^n}\), \(F=T^n-1\), and \(M=\operatorname{rad}(F)\). For \(P\mid F\), use \(U_P,V_P\) from Proposition~\ref{prop:normal-test}. Choose a nontrivial additive character \(\psi\) of \(k\), put \(\Psi(x)=\psi(\Tr_{E/k}(x))\), and define the Fourier transform using \(\Psi(-bx)\).
For a monic divisor \(g\mid M\), put
\[
\boldsymbol1_g(x)=\prod_{P\mid g}\boldsymbol1_{x\notin U_P},\qquad
\Theta(g)=\prod_{P\mid g}(1-q^{-\deg P}),\qquad
V_g=\bigoplus_{P\mid g}V_P.
\]
Empty products are one, and \(V_1=\{0\}\). For \(b=\sum_{P\mid g}b_P\in V_g\), define
\[
c_g(b)=\prod_{b_P\ne0}\frac{-1}{q^{\deg P}-1}.
\]
Then
\begin{equation}\label{v93:eq-local-character-expansion}
\boldsymbol1_g(x)=\Theta(g)\sum_{b\in V_g}c_g(b)\Psi(bx),
\qquad c_g(0)=1,
\end{equation}
and its exact Fourier support is given by
\begin{equation}\label{v128:eq-normal-spectrum}
\widehat{\boldsymbol1_g}(b)=
\begin{cases}
q^n\Theta(g)c_g(b),&b\in V_g,\\
0,&b\notin V_g.
\end{cases}
\qquad \operatorname{supp}\widehat{\boldsymbol1_g}=V_g.
\end{equation}
Writing \(W_q(g)=2^{\omega_q(g)}\), where \(\omega_q(g)\) counts distinct irreducible factors, we also have
\begin{equation}\label{v93:eq-local-coefficient-sum}
\sum_{0\ne b\in V_g}|c_g(b)|=W_q(g)-1.
\end{equation}
\end{proposition}
\begin{proof}
Subspace orthogonality gives
\[
\boldsymbol1_{x\notin U_P}
=1-q^{-\deg P}\sum_{b\in V_P}\Psi(bx)
=(1-q^{-\deg P})\left(1-\frac{\sum_{0\ne b\in V_P}\Psi(bx)}{q^{\deg P}-1}\right).
\]
Here \(U_P\) is a \(k\)-subspace. If \(\Tr_{E/k}(bu)\ne0\) for some \(u\in U_P\), then \(\Tr_{E/k}(btu)\) ranges over \(k\) as \(t\) ranges over \(k\), so the nontrivial character \(\psi\) cannot be identically one on these values. Thus the relative trace orthogonal complement parametrizes exactly the characters trivial on \(U_P\). Multiply these identities. The direct sum gives each frequency a unique component representation, proving the expansion. Fourier orthogonality gives the coefficient formula; every \(c_g(b)\ne0\), so the support is exactly \(V_g\).

Fix a nonempty set \(J\subseteq\{P:P\mid g\}\), and require \(b_P\ne0\) precisely for \(P\in J\). There are \(\prod_{P\in J}(q^{\deg P}-1)\) such frequencies, and their coefficients have total absolute value one. Summing over all nonempty \(J\) gives the final identity.
\end{proof}

Taking \(g=M\) gives the indicator of the normal elements \(\mathcal N\). Write \(n=p^a m\), with \(p\nmid m\). Then
\[
\operatorname{supp}\widehat{\boldsymbol1_{\mathcal N}}
=(M(\sigma)E)^\perp,\qquad
\left|\operatorname{supp}\widehat{\boldsymbol1_{\mathcal N}}\right|=q^m.
\]
Every allowed frequency has a nonzero coefficient. Thus a translation preserves normality exactly when it belongs to \(M(\sigma)E\). The spectrum recovers the full translation stabilizer and the minimum linear information in Proposition~\ref{v93:prop-normal-linear-info}. If \(n\) is a power of the characteristic, the frequency space is the base field \(k\), reflecting the sufficiency of the total trace; if \(p\nmid n\), the support is all of \(E\).

\begin{proposition}[The exact spectrum of a prescribed annihilator stratum]\label{v130:prop-annihilator-spectrum}
Retain \(E,k,F,\Psi\) above. For a monic divisor \(D\mid F\), put
\[
U_D=\ker D(\sigma),\qquad
\mathcal S_D=\{x\in E:\mu_x=D\},\qquad
J_D=\operatorname{rad}(D)(\sigma)U_D.
\]
Write \(D=\prod_{i=1}^sP_i^{a_i}\), with distinct irreducible \(P_i\), \(d_i=\deg P_i\), and \(a_i\ge1\), and set
\[
W_i=\ker P_i^{a_i}(\sigma),\qquad
W_i^-=\ker P_i^{a_i-1}(\sigma).
\]
All orthogonal complements are taken in \(E\) with respect to \(\Tr_{E/k}(xy)\). Then
\begin{equation}\label{v130:eq-annihilator-spectrum}
\widehat{\boldsymbol1_{\mathcal S_D}}(b)
=\prod_{i=1}^s\left(
q^{a_i d_i}\boldsymbol1_{W_i^\perp}(b)
-q^{(a_i-1)d_i}\boldsymbol1_{(W_i^-)^\perp}(b)\right).
\end{equation}
In particular,
\begin{equation}\label{v130:eq-annihilator-support}
\operatorname{supp}\widehat{\boldsymbol1_{\mathcal S_D}}=J_D^\perp,
\qquad
\left|\operatorname{supp}\widehat{\boldsymbol1_{\mathcal S_D}}\right|
=q^{\,n-\deg D+\deg\operatorname{rad}(D)}.
\end{equation}
For \(D=1\), the empty product is one and \(\mathcal S_1=J_1=\{0\}\), so the formulas remain valid.
\end{proposition}
\begin{proof}
The primary decomposition gives \(U_D=\bigoplus_iW_i\). An element \(x=\sum_i x_i\) has annihilator exactly \(D\) if and only if every \(x_i\) lies in \(W_i\setminus W_i^-\). Thus the character sum over \(\mathcal S_D\) factors as
\[
\prod_i\left(\sum_{x\in W_i}\Psi(-bx)
-\sum_{x\in W_i^-}\Psi(-bx)\right).
\]
Subspace orthogonality and \(\dim_k\ker P_i^j(\sigma)=jd_i\) give~\eqref{v130:eq-annihilator-spectrum}.

Since \(W_i^-\subseteq W_i\), each factor vanishes outside \((W_i^-)^\perp\). Within that space it takes one of the nonzero values \(q^{(a_i-1)d_i}(q^{d_i}-1)\) and \(-q^{(a_i-1)d_i}\). The Fourier support is therefore exactly \(\bigcap_i(W_i^-)^\perp\).
On \(W_i\), every \(P_j(\sigma)\) with \(j\ne i\) is invertible, whereas \(P_i(\sigma)W_i=W_i^-\). Hence
\[
J_D=\bigoplus_iW_i^-,\qquad
\bigcap_i(W_i^-)^\perp=J_D^\perp.
\]
The support size follows from \(\dim_kJ_D=\deg D-\deg\operatorname{rad}(D)\).
\end{proof}

Formula~\eqref{v130:eq-annihilator-support} expresses the minimum amount of linear information in Proposition~\ref{v110:prop-annihilator-info} as the dimension of the exact spectrum. Taking \(D=F\) recovers the support of the normality indicator; taking \(D=1\) gives the full spectrum of the indicator of zero.
For example, in \(\F_{27}/\F_3\), the four strata for \(D_j=(T-1)^j\) satisfy
\[
\begin{array}{c|rrrr}
j&0&1&2&3\\\hline
|\mathcal S_{D_j}|&1&2&6&18\\
|\operatorname{supp}\widehat{\boldsymbol1_{\mathcal S_{D_j}}}|&27&27&9&3
\end{array}
\]
A larger stratum can therefore require fewer frequencies to determine membership.

\begin{corollary}[Normal elements in affine fibers and frequency intersections]\label{v128:cor-affine-spectrum}
Let \(B=x_0+W\), where \(W\subseteq E\) is any \(k\)-subspace. Then
\begin{equation}\label{v128:eq-affine-spectrum}
\#(B\cap\mathcal N)
=|W|\Theta(M)\sum_{b\in V_M\cap W^\perp}c_M(b)\Psi(bx_0).
\end{equation}
\end{corollary}
\begin{proof}
Sum~\eqref{v93:eq-local-character-expansion} over \(x_0+W\). The sum \(\sum_{w\in W}\Psi(bw)\) is \(|W|\) for \(b\in W^\perp\) and zero otherwise.
\end{proof}

For any monic \(D\mid F\), the same summation gives
\begin{equation}\label{v130:eq-stratum-affine-count}
\#\bigl((x_0+W)\cap\mathcal S_D\bigr)
=\frac{|W|}{q^n}\sum_{b\in J_D^\perp\cap W^\perp}
\widehat{\boldsymbol1_{\mathcal S_D}}(b)\Psi(bx_0).
\end{equation}
The coefficients are given by~\eqref{v130:eq-annihilator-spectrum}; no Frobenius invariance of \(W\) is required.

Proposition~\ref{v110:prop-affine-count} intersects invariant kernels in the original space; this formula intersects frequencies in the dual space. Both compute the same integer. For a Frobenius-invariant direction, Corollary~\ref{v88:thm-invariant-affine-normal-count} already gives a product formula. A general direction can link different local conditions, so their individual proportions cannot simply be multiplied.

\begin{example}[Two dependent trace conditions in the eight-element field]\label{v128:ex-affine-spectrum-eight}
Use the self-dual normal basis \((v,v^2,v^4)\) from Example~\ref{v110:ex-trace-coordinates-eight}, where \(v^3+v^2+1=0\). For fixed \(r,s\in\F_2\), put
\[
B_{r,s}=rv+sv^2+\langle v^4\rangle_{\F_2}.
\]
Here \(M=(T+1)(T^2+T+1)\), \(\Theta(M)=3/8\), and \(V_M=E\), while \(W^\perp=\langle v,v^2\rangle\). The two local dual spaces are \(\F_2\) and \(\ker\Tr\), respectively. Hence
\[
c_M(v)=c_M(v^2)=\frac13,\qquad c_M(v+v^2)=-\frac13.
\]
Self-duality gives \(\Psi(vx_0)=(-1)^r\) and \(\Psi(v^2x_0)=(-1)^s\), so
\[
\#(B_{r,s}\cap\mathcal N)
=\frac{3+(-1)^r+(-1)^s-(-1)^{r+s}}4.
\]
The four fiber counts are \(1,1,1,0\), recovering the earlier distribution. Two trace conditions, each separately nonzero on some normal element, can together exclude every normal element.
\end{example}

\section{Gauss and Jacobi Sums}
\label{sec:gauss-jacobi}
From now on extend every multiplicative character by \(\chi(0)=0\), including the trivial multiplicative character \(\varepsilon\).
This convention must be fixed when zero is involved.
\leavevmode\phantomsection\label{def:gauss-sum}\ignorespaces For a nontrivial additive character \(\psi\), define the \textbf{Gauss sum}
\[
G(\chi,\psi)=\sum_{x\in K^\times}\chi(x)\psi(x).
\]
If \(\chi=\varepsilon\), then \(G(\varepsilon,\psi)=-1\).

Gauss sums describe the additive Fourier transform of multiplicative characters. For the fixed additive character \(\psi\) of this chapter,
for \(a\ne0\), the substitution \(y=-ax\) gives
\[
\widehat\chi(a)=\sum_{x\in K}\chi(x)\psi(-ax)
=\chi(-a)^{-1}G(\chi,\psi).
\]
Thus one Gauss sum determines the Fourier coefficients of \(\chi\) at all nonzero frequencies.
If \(\chi\) is nontrivial, its zero-frequency coefficient vanishes; the next theorem determines the absolute values of all remaining coefficients.
This allows multiplicative power conditions to enter equation counts through additive character orthogonality.

\begin{theorem}\label{thm:gauss}
If both \(\chi\) and \(\psi\) are nontrivial, then
\[
|G(\chi,\psi)|=\sqrt Q,
\qquad G(\chi,\psi)G(\chi^{-1},\psi)=\chi(-1)Q.
\]
For every \(a\in K^\times\), we also have \(G(\chi,x\mapsto\psi(ax))=\chi(a)^{-1}G(\chi,\psi)\).
\end{theorem}
\begin{proof}
Both additive and multiplicative character values are complex roots of unity, so complex conjugation equals inversion.
Expanding the squared absolute value and substituting \(x=ty\) gives
\[
\begin{aligned}
|G(\chi,\psi)|^2
&=\sum_{x,y\in K^\times}\chi(x)\chi(y)^{-1}\psi(x-y)\\
&=\sum_{t\in K^\times}\chi(t)\sum_{y\in K^\times}\psi((t-1)y).
\end{aligned}
\]
When \(t=1\), the inner sum is \(Q-1\).
When \(t\ne1\), multiplication by \(t-1\) permutes \(K^\times\), so the inner sum is \(\sum_{z\in K^\times}\psi(z)=-1\).
A nontrivial multiplicative character also satisfies \(\sum_{t\in K^\times}\chi(t)=0\), hence
\[
|G(\chi,\psi)|^2
=(Q-1)-\sum_{t\ne0,1}\chi(t)
=(Q-1)-(0-1)=Q.
\]
Thus \(|G(\chi,\psi)|=\sqrt Q\).

Now compute the complex conjugate, putting \(y=-x\):
\[
\overline{G(\chi,\psi)}
=\sum_{x\ne0}\chi(x)^{-1}\psi(-x)
=\chi(-1)^{-1}G(\chi^{-1},\psi).
\]
Since \((-1)^2=1\) and \(\chi(-1)^2=1\), we have \(\chi(-1)^{-1}=\chi(-1)\).
Multiplying by \(G(\chi,\psi)\) and using the squared-absolute-value formula gives \(G(\chi,\psi)G(\chi^{-1},\psi)=\chi(-1)Q\).
Finally, for \(a\ne0\), the substitution \(y=ax\) gives
\[
\sum_{x\ne0}\chi(x)\psi(ax)
=\chi(a)^{-1}\sum_{y\ne0}\chi(y)\psi(y),
\]
which is the character-scaling formula.
\end{proof}

\subsection{Complex absolute values and \texorpdfstring{$p$}{p}-adic valuations of Gauss sums}
\label{v132:subsec-adic-gauss}
The preceding complex absolute-value formula controls the size of a character sum, not the divisibility of the corresponding algebraic integer. Let $q=p^r$, choose compatible $p$th and $(q-1)$st roots of unity in $\overline{\mathbb Q}_p$, and fix the Teichmüller character $\omega:\F_q^\times\to\mu_{q-1}$. Gauss sums can now also be evaluated in a $p$-adic cyclotomic field. Normalize $v_p(p)=1$; the ramification from $p$th roots of unity permits fractional valuations.

\begin{theorem}[Digit-sum valuation of Gauss sums (Stickelberger)]
\label{v132:thm-gauss-digit-valuation}
Use the standard additive character $\psi(x)=\zeta_p^{\operatorname{Tr}_{\F_q/\F_p}(x)}$ and evaluate in the chosen $p$-adic cyclotomic field. For $1\le a\le q-2$, let $s_p(a)$ denote the sum of the base-$p$ digits of $a$. Then
\begin{equation}\label{v132:eq-gauss-digit-valuation}
 \boxed{v_p\bigl(G(\omega^{-a},\psi)\bigr)=\frac{s_p(a)}{p-1}.}
\end{equation}
Here $\omega^{-a}$ is nontrivial. The valuation depends on the chosen $p$-adic embedding, whereas the complex absolute value is always $\sqrt q$. This valuation formula also follows from Gross--Koblitz\cite{GrossKoblitz1979,Robert2001GrossKoblitz}.
\end{theorem}

Formula~\eqref{v134:eq-phase-gauss-fourier} in Chapter 5 uses a different type of Gauss sum: its additive exponent has order $\ell\ne p$, rather than order $p$ as here. These must be distinguished; the present digit-sum formula does not directly apply to the former.

For example, let $q=8$, $p=2$. The binary digit sums of $1\le a\le6$ are $1,1,2,1,2,2$. The corresponding nontrivial Gauss sums all have complex absolute value $\sqrt8$, but their $2$-adic valuations are respectively $1,1,2,1,2,2$. This uses the stated arithmetic valuation theorem; it cannot be inferred from the complex absolute value alone.

\subsection*{From cyclotomic orbits to slopes of exponential sums}
\phantomsection\label{v132:problem-gauss-slopes}
Section~\ref{revadic:sec-newton} reconstructs curve polynomials from point counts and determines Fermat slopes through Jacobi sums.
Fix a family of extension-field exponential sums $S_m$ arising from a multiplicative subgroup, root set, or polynomial, and form the formal series
\[
 L(T)=\exp\left(\sum_{m\ge1}\frac{S_mT^m}{m}\right).
\]
When polynomiality or rationality has been established and the numerator and denominator have suitable $p$-adic integral coefficients, one may study their Newton polygons. Which cyclotomic orbits, phases, or root-set symmetries force given slopes of Frobenius eigenvalues? When are the resulting lower bounds sharp? The precise family of sums and its coefficient field must be specified before making a slope assertion; the formal expression by itself does not imply rationality.

\leavevmode\phantomsection\label{def:jacobi-sum}\ignorespaces Define the \textbf{Jacobi sum}
\[
J(\chi,\lambda)=\sum_{x\in K}\chi(x)\lambda(1-x).
\]
\begin{proposition}\label{prop:jacobi}
If all \(\chi,\lambda,\chi\lambda\) are nontrivial, then
\[
G(\chi,\psi)G(\lambda,\psi)
=J(\chi,\lambda)G(\chi\lambda,\psi),
\qquad |J(\chi,\lambda)|=\sqrt Q.
\]
If \(\chi\ne\varepsilon\), then \(J(\chi,\chi^{-1})=-\chi(-1)\),
\(J(\varepsilon,\chi)=J(\chi,\varepsilon)=-1\), and \(J(\varepsilon,\varepsilon)=Q-2\).
\end{proposition}
\begin{proof}
Expand the product of the two Gauss sums and group terms by \(z=x+y\).
For \(z=0\), the total contribution is \(\lambda(-1)\sum_{x\ne0}(\chi\lambda)(x)=0\).
For \(z\ne0\), put \(x=zt\) and \(y=z(1-t)\) to separate the expression into the product of \(J\) and \(G(\chi\lambda,\psi)\).
The absolute-value formula follows from the preceding theorem.

For inverse characters, the variable \(t=x/(1-x)\) maps \(x\ne0,1\) bijectively onto \(t\ne0,-1\), so
\[
J(\chi,\chi^{-1})=\sum_{t\ne0,-1}\chi(t)=-\chi(-1).
\]
The formulas involving trivial characters follow by removing zero and one from the complete character sums.
\end{proof}

\subsection{Recovering addition from one row of Jacobi data}
Jacobi sums can also determine field operations in the reverse direction\cite{Snowden2024Jacobi}. With the labeled multiplicative group fixed, one faithful character row already recovers addition.

\begin{proposition}[A faithful character row recovers addition]\label{v150:prop-jacobi-recovery}
Let $\chi$ be a faithful complex character of $\F_q^\times$. Given this labeled multiplicative group, the row
$\{J(\chi,\lambda):\lambda\in\widehat{\F_q^\times}\}$ uniquely determines addition in $\F_q$.
\end{proposition}
\begin{proof}
Extend every multiplicative character by zero at zero. For $y\ne0$, orthogonality gives
\begin{equation}\label{v150:eq-jacobi-recovery}
 \frac1{q-1}\sum_{\lambda}J(\chi,\lambda)\lambda(y)^{-1}
 =\chi(1-y).
\end{equation}
Faithfulness and $\chi(0)=0$ uniquely determine $R(y)=1-y$; also $R(0)=1$.
The multiplicative group determines $-1$ as well: it is one in even characteristic and the unique element of order two in odd characteristic. Thus
\[
 x+y=xR(-y/x)\qquad(x\ne0),
\]
together with $0+y=y$, recovers all of addition.
\end{proof}

\begin{corollary}[Power-permutation stabilizers of a Jacobi row]\label{v150:cor-jacobi-stabilizer}
Let $q=p^f>2$ and $N=q-1$. For a fixed faithful $\chi$,
\begin{equation}\label{v150:eq-jacobi-stabilizer}
 \left\{a\in(\mathbb Z/N\mathbb Z)^\times:
 J(\chi^a,\lambda^a)=J(\chi,\lambda)\text{ for all }\lambda\right\}
 =\langle p\rangle.
\end{equation}
\end{corollary}
\begin{proof}
Frobenius substitution gives inclusion from right to left. Conversely, let $a$ belong to the set on the left and choose $b$ with $ab\equiv1\pmod N$.
Extend $F(x)=x^b$ by $F(0)=0$, and transport addition by
$x\oplus y=F^{-1}(F(x)+F(y))$. Substitution $z=F(x)$ shows that the Jacobi row for this new addition is $J(\chi^a,\lambda^a)$, equal to the original row.
Proposition~\ref{v150:prop-jacobi-recovery} forces $\oplus=+$, so $F$ is a field automorphism. Hence $b\in\langle p\rangle$ and therefore $a\in\langle p\rangle$.
\end{proof}

\begin{example}[Relations retained by a partial character table]\label{v150:ex-partial-jacobi}
On the same multiplicative group of $\F_{13}$, put
\[
 x\oplus y=(x^5+y^5)^5.
\]
Since $5^2\equiv1\pmod{12}$, this defines another field addition, but $1\oplus1=6\ne2$.
Every character of order dividing four satisfies $\lambda^5=\lambda$. Substitution shows that the two additions have identical Jacobi data on these characters.
For $C_i=2^i\langle2^4\rangle$, double Fourier inversion of these data recovers
\[
 \left(\#\{x\in C_i:1-x\in C_j\}\right)_{0\le i,j<4}
 =\begin{pmatrix}0&1&0&1\\1&0&1&1\\0&1&2&0\\1&1&0&1\end{pmatrix}.
\]
The entries sum to $11=q-2$. Equal joint counts between cosets do not require equality of addition on individual elements.
\end{example}

\Needspace{18\baselineskip}

\Needspace{15\baselineskip}
\section{Short Linear Relations in Multiplicative Subgroups and Gaussian Periods}\label{v95:sec-subgroup-periods}
Multiplicative order first determines the least field containing one element, and Frobenius orbits then organize its conjugates; neither determines all additive relations among several elements.
Even when every nontrivial orbit in a multiplicative subgroup is a normal basis, elements chosen from different orbits may be dependent.
Passing from orbit structure to the linear position of the entire subgroup therefore requires recording the support of a shortest relation.
We first obtain necessary restrictions from multiplicative orders, then detect joint independence using trace hyperplanes, and finally count intersections and relations through moments of Gaussian periods.
The trace-code interpretation has appeared in Section~\ref{v84:subsec-trace-mds}; here the main objects are linear relations in finite fields.

\subsection{Linear relation length and multiplicative orders}\label{v95:subsec-short-relations}
Let $n>1$, $\gcd(n,q)=1$, and $k=\ord_n(q)\ge2$. Put $L=\F_{q^k}$ and let $H=\langle\beta\rangle$ be the unique subgroup of $L^\times$ of order $n$. Since $\F_q(\beta)=L$, the $\F_q$-span of $H$ is $L$.
\leavevmode\phantomsection\label{v95:def-linear-relation-length}\ignorespaces
Define the \textbf{linear relation length} of $H$ by
\begin{equation}\label{v95:eq-lambda-definition}
\lambda_q(H)=\min\left\{\#\{h\in H:c_h\ne0\}:\ 
 \sum_{h\in H}c_hh=0,\quad(c_h)\in\F_q^H\setminus\{0\}\right\}.
\end{equation}
Its support counts distinct subgroup elements, independently of the coefficient values. Any $k+1$ elements of $H$ are dependent, while $0\notin H$; hence
\[2\le\lambda_q(H)\le k+1.\]
The extreme value $\lambda_q(H)=k+1$ says that every set of at most $k$ distinct elements of $H$ is linearly independent. We call this \textbf{maximal linear general position}. It is stronger than requiring every nonidentity element to have degree $k$, since the latter merely excludes elements contained in proper subfields.

\begin{proposition}[The first restrictions from multiplicative orders]\label{v95:prop-order-lambda}
Set $g=\gcd(n,q-1)$. Then
\begin{equation}\label{v95:eq-order-bound}
\boxed{\lambda_q(H)\le1+\min_{\substack{d\mid n\\d>1}}\ord_d(q)},\qquad
\boxed{\lambda_q(H)=2\iff g>1.}
\end{equation}
If $\lambda_q(H)\ge s+1$ with $2\le s\le k$, then $\ord_d(q)\ge s$ for every $d\mid n$, $d>1$, and every prime divisor $\ell\mid n$ satisfies $\ell\ge s+1$.
\end{proposition}
\begin{proof}
Choose an element $x\in H$ of order $d$, and write $r=\ord_d(q)$. Its minimal polynomial of degree $r$ supplies a nontrivial relation among $1,x,\ldots,x^r$. Since $r\le\varphi(d)\le d-1$, these powers are distinct, giving $\lambda_q(H)\le r+1$. Two distinct elements $x,y\in H$ are linearly dependent precisely when $x/y\in H\cap\F_q^\times\setminus\{1\}$. The intersection of these cyclic groups has order $\gcd(n,q-1)=g$, proving the two-term criterion. The remaining assertions follow from the order bound and $\ord_\ell(q)\mid\ell-1$.
\end{proof}

Proposition~\ref{v142:prop-subfield-distance} expresses this order obstruction through restriction of a trace code: if the subgroup of order $d$ generates a subfield of degree $r<k$, the original code loses at least $d-r$ from the Singleton distance bound, while short relations bound the dual distance.

\subsection{When every Frobenius orbit is a normal basis}\label{v95:subsec-frobenius-normal}
The equality $\lambda_q(H)=k+1$ first forces every nontrivial subgroup to generate the full extension and then restricts the linear relations within each Frobenius orbit.
\begin{proposition}[Cyclotomic and normality consequences of general position]\label{v95:prop-normal-cyclotomic}
If $\lambda_q(H)=k+1$, then:
\begin{enumerate}
\item $\ord_d(q)=k$ for every $d\mid n$, $d>1$, and $n\mid\Phi_k(q)$;
\item every nonzero $q$-cyclotomic coset modulo $n$ has length $k$, and $n\equiv1\pmod k$;
\item $H\setminus\{1\}$ is partitioned into $(n-1)/k$ Frobenius orbits, each of which is a normal basis of $L/\F_q$. In particular, $\Tr_{L/\F_q}(h)\ne0$ for $h\ne1$;
\item the minimal polynomial $m_h(X)$ of every $h\in H\setminus\{1\}$ has all its coefficients nonzero.
\end{enumerate}
\end{proposition}
\begin{proof}
The order equality follows from Proposition~\ref{v95:prop-order-lambda}. Proposition~\ref{v124:prop-uniform-cyclotomic} then gives $n\mid\Phi_k(q)$. Every nonidentity $n$th root of unity has degree $k$, giving the uniform Frobenius orbit lengths. Each orbit consists of $k$ distinct elements of $H$ and is linearly independent by general position; hence it is a normal basis, and its sum, the trace, cannot vanish. Finally, the order of $h$ is at least $k+1$. If its minimal polynomial had a zero coefficient, $m_h(h)=0$ would give a nontrivial relation involving at most $k$ of the distinct elements $1,h,\ldots,h^k$, a contradiction.
\end{proof}
One must distinguish \textbf{uniform Frobenius orbit lengths} from \textbf{linear independence of the elements of an orbit}. Even when every orbit is a normal basis, one cannot infer that arbitrary $k$ elements drawn from different orbits are independent. The primitive normal basis theorem guarantees an element that generates \(L^\times\) and is normal; here the assertion is that \emph{every nonidentity element} of a specified subgroup is normal. The objects and quantifiers are different.

\par\Needspace{6\baselineskip}

\begin{example}[Normal orbits and a short relation at order nineteen]\label{v96:ex-nineteen}
Take \(q=11\) and \(n=19\). Since \(11^2\equiv7\pmod{19}\) and \(11^3\equiv1\pmod{19}\), we have \(\ord_{19}(11)=3\).
In \(\F_{11}[X]\),
\begin{equation}\label{v96:eq-nineteen-factorization}
\begin{gathered}
\Phi_{19}(X)=\prod_{(a,b)\in\mathcal B}(X^3+aX^2+bX-1),\\
\mathcal B=\{(1,2),(2,-3),(3,-5),(3,-2),(5,-3),(-2,-1)\}.
\end{gathered}
\end{equation}
This identity follows by direct multiplication. Every irreducible factor of \(\Phi_{19}\) has degree \(\ord_{19}(11)=3\), so all six displayed cubic factors are irreducible.
Choose a root \(\beta\) of \(f(X)=X^3+3X^2-2X-1\), and put \(H=\langle\beta\rangle\subseteq\F_{11^3}^\times\).
The order of \(\beta\) is nineteen, and the eighteen elements of \(H\setminus\{1\}\) are exactly the roots of these cubics.

Their quadratic coefficients are nonzero, so all eighteen elements have nonzero trace; each also has degree three.
Since \(\ord_3(11)=2\), Theorem~\ref{thm:trace-normal-classification} makes every one of them normal.
Moreover, all four coefficients of each cubic are nonzero, and \(19\mid\Phi_3(11)=133\).
Thus all the order conditions, orbitwise normality, and full-support minimal-polynomial properties in Proposition~\ref{v95:prop-normal-cyclotomic} hold.

Nevertheless, direct multiplication gives
\[
(X-3)f(X)=X^4+5X+3,
\]
and hence
\begin{equation}\label{v96:eq-nineteen-relation}
3+5\beta+\beta^4=0,\qquad
3\beta+5\beta^2+\beta^5=0.
\end{equation}
Since \(\gcd(19,10)=1\), Proposition~\ref{v95:prop-order-lambda} excludes two-term relations, and therefore
\[\lambda_{11}(H)=3<4.\]
The exponents 1, 2, and 5 in the second relation lie in the three distinct nonzero Frobenius orbits
\(\{1,7,11\}\), \(\{2,3,14\}\), and \(\{5,16,17\}\).
The failure is thus not normality within an orbit, but joint independence across different orbits: even together, the listed necessary properties do not imply maximal linear general position.
\end{example}

\begin{proposition}[Unique representations by short weighted sums]\label{v95:prop-additive-unique}
If $\lambda_q(H)\ge2r+1$, all weighted sums $\sum_{h\in H}c_hh$ with coefficient support of size at most $r$ are distinct, and
\begin{equation}\label{v95:eq-additive-packing}
\boxed{\sum_{j=0}^{r}\binom nj(q-1)^j\le q^k.}
\end{equation}
In particular, $\lambda_q(H)\ge4$ excludes every weighted three-term relation among distinct elements, while $\lambda_q(H)\ge5$ gives uniqueness of weighted two-term sums.
\end{proposition}
\begin{proof}
Equality of two such sums would give a nontrivial relation supported on at most $2r$ elements. The sum on the left of the inequality counts all the eligible coefficient vectors.
\end{proof}

\subsection{Trace hyperplanes and Singer difference sets}\label{v95:subsec-trace-hyperplanes}
Using the trace pairing, for $a\in L^\times$ define
\[
V_a=\{x\in L:\Tr_{L/\F_q}(ax)=0\},\qquad N_H(a)=\#(H\cap V_a).
\]
By nondegeneracy of the trace pairing, $V_a$ runs through all $\F_q$-hyperplanes, with $V_a=V_b$ precisely when $a/b\in\F_q^\times$.
\begin{proposition}[The hyperplane criterion for general position]\label{v95:prop-hyperplane-criterion}
\[
\boxed{\lambda_q(H)=k+1\iff\max_{a\in L^\times}N_H(a)\le k-1.}
\]
\end{proposition}
\begin{proof}
Any $k$ distinct elements of $H$ are dependent precisely when their span is proper, in which case a nonzero trace-linear functional vanishes simultaneously on them.
\end{proof}

\leavevmode\phantomsection\label{v96:def-projective-arc}\ignorespaces
Now suppose \(\gcd(n,q-1)=1\). The \textbf{projective image} of the subgroup is
\[
\mathcal P(H)=\{\F_qh:h\in H\}\subseteq\operatorname{PG}(k-1,q).
\]
Here \(\operatorname{PG}(k-1,q)\) is the projective space of one-dimensional \(\F_q\)-subspaces of \(L\); the coprimality condition makes the \(n\) displayed points distinct.
An \(n\)-point set is a \textbf{projective arc}, or an \(n\)-\textbf{arc}, if representatives of every \(k\) of its points are linearly independent.
This condition is unchanged by nonzero rescaling of representatives and is equivalent to every projective hyperplane containing at most \(k-1\) of the points.

\begin{corollary}[Short relations, trace codes, and projective arcs]\label{v96:cor-arc-trace-equivalence}
If \(\gcd(n,q-1)=1\), the following conditions are equivalent:
\begin{enumerate}
\item \(\lambda_q(H)=k+1\);
\item the \(k\)-dimensional trace code \(C_\beta\) constructed as in \eqref{v84:eq-trace-code-convention} is MDS;
\item \(\mathcal P(H)\) is an \(n\)-arc in \(\operatorname{PG}(k-1,q)\);
\item \(\max_{a\in L^\times}N_H(a)\le k-1\).
\end{enumerate}
\end{corollary}
\begin{proof}
The trace pairing maps the elements of \(H\) to the generator columns of \(C_\beta\) by a linear isomorphism.
Thus independence of every \(k\) subgroup elements is equivalent both to independence of every \(k\) columns and to independence of representatives of every \(k\) projective points.
The last condition follows from Proposition~\ref{v95:prop-hyperplane-criterion}.
\end{proof}
For \(k=3\), the arc condition simply excludes three collinear points.
Each nontrivial Frobenius orbit in Example~\ref{v96:ex-nineteen} gives a noncollinear triangle, but
\(3\beta+5\beta^2+\beta^5=0\) gives collinear points chosen from three different triangles.
This is the geometric distinction between normality within each orbit and general position of the entire subgroup.

\begin{definition}[Difference set]\label{v103:def-difference-set}
Let \(G\) be a finite group of order \(v\). A subset \(D\subseteq G\) of size \(b\) is a
\((v,b,\lambda)\)-\textbf{difference set} if every \(g\ne1\) has exactly \(\lambda\) representations
\(g=d_1d_2^{-1}\) by ordered pairs \((d_1,d_2)\in D^2\).
Equivalently, \(|D\cap gD|=\lambda\) for every \(g\ne1\).
Counting ordered pairs of distinct elements also gives \(b(b-1)=\lambda(v-1)\).
\end{definition}

\begin{proposition}[Trace-zero points form a Singer difference set]\label{v103:prop-singer-difference-set}
Let \(L=\F_{q^k}\), with \(k\ge3\), and \(G=L^\times/\F_q^\times\). Then
\[
D=\{[x]\in G:\Tr_{L/\F_q}(x)=0\}
\]
is a \textbf{Singer difference set} in the cyclic group \(G\), with parameters
\[
\left(\frac{q^k-1}{q-1},\ \frac{q^{k-1}-1}{q-1},\ \frac{q^{k-2}-1}{q-1}\right).
\]
\end{proposition}
\begin{proof}
The trace-zero condition is unchanged by a nonzero base-field scalar, so \(D\) is well defined.
The trace is a nonzero linear functional with a kernel of dimension \(k-1\).
Each projective point has \(q-1\) nonzero representatives, giving the stated values of \(|G|\) and \(|D|\).
For \([a]\ne1\), membership \([x]\in D\cap[a]D\) is equivalent to
\[
\Tr(x)=0,\qquad\Tr(a^{-1}x)=0.
\]
These two trace functionals are linearly independent. Indeed, if \(c,d\in\F_q\) and \(c\Tr(x)+d\Tr(a^{-1}x)=0\) for every \(x\),
nondegeneracy of the trace pairing gives \(c+da^{-1}=0\).
Since \(a\notin\F_q\), this forces \(c=d=0\).
Their common kernel therefore has dimension \(k-2\), and
\[
|D\cap[a]D|=\frac{q^{k-2}-1}{q-1}.
\]
This is the required constant intersection number for every nonidentity translate.
\end{proof}

If \(\lambda_q(H)\ge3\), then \(H\cap\F_q^\times=\{1\}\), so its image in \(G\) is isomorphic to \(H\).
With this identification, \(N_H(a)=|H\cap[a]^{-1}D|\).
Thus general position for the subgroup is an intersection constraint with all translates of the Singer difference set.
For cyclic arcs of Singer type in the projective plane, see also \cite{Giulietti2002Subgroups,MomiharaXiang2022Subgroups}.

\subsection{Three-term relations from one cyclotomic intersection number}\label{v95:subsec-cyclotomic-triples}
To detect three-term relations with \emph{arbitrary base-field coefficients}, counting $H\cap(1+H)$ is insufficient: we must first adjoin all nonzero base-field scalars. Put
\[
g=\gcd(n,q-1),\qquad U=\F_q^\times H,\qquad M=|U|=\frac{(q-1)n}{g},
\]
and define the \textbf{scalar-saturated cyclotomic intersection number}
\[\mathcal C(U)=\#\{u\in U:1+u\in U\}.\]
\begin{theorem}[Exact enumeration of dependent triples]\label{v95:thm-triple-cyclotomic}
Suppose $g=1$ and let $T_3(H)$ denote the number of linearly dependent three-element subsets of $H$. Then
\begin{equation}\label{v95:eq-triple-exact}
\boxed{6T_3(H)=n\bigl(\mathcal C(U)-(q-2)\bigr).}
\end{equation}
In particular, for $k\ge3$,
\begin{equation}\label{v95:eq-lambda-4-cyclo}
\boxed{\lambda_q(H)\ge4\iff \mathcal C(U)=q-2.}
\end{equation}
\end{theorem}
\begin{proof}
Since $g=1$, each $u\in U$ has a unique representation $th$ with $t\in\F_q^\times$ and $h\in H$. The elements $u\in\F_q^\times$ satisfying $1+u\in U$ are precisely the $q-2$ elements of $\F_q^\times\setminus\{-1\}$. Every remaining solution has unique expressions $u=tx$ and $1+u=sy$, with $x,y\in H\setminus\{1\}$ and $t,s\in\F_q^\times$. This yields a relation $1+tx-sy=0$ among three distinct elements. Conversely, divide any ordered dependent triple by its first element, then normalize the coefficient of the constant element to one. This gives a unique solution of the stated type. There are $n$ choices for the first element, and every unordered triple has six orderings, proving \eqref{v95:eq-triple-exact}. Finally, $g=1$ already excludes two-term dependencies, so the absence of dependent triples is equivalent to $\lambda_q(H)\ge4$.
\end{proof}
The term $q-2$ is not an error term: it counts the unavoidable \emph{scalar solutions} in $U$. For $q=2$, we have $U=H$, and the criterion reduces to $(1+H)\cap H=\varnothing$. For $q>2$, ordinary sum-freeness does not exclude all weighted relations over the base field.

\subsection{Gaussian periods and their first three moments}\label{v95:subsec-period-moments}
Let $Q=q^k$ and $e=(Q-1)/M$. Fix a nontrivial additive character $\psi$ of $\F_q$, and put $\Psi=\psi\circ\Tr_{L/\F_q}$. Choose representatives $a_1,\ldots,a_e$ for $L^\times/U$. Define the \textbf{Gaussian periods}
\begin{equation}\label{v95:eq-gaussian-period}
\eta(a)=\sum_{u\in U}\Psi(au),\qquad \eta_i=\eta(a_i)\quad(1\le i\le e).
\end{equation}
Because $\F_q^\times\subseteq U$, each sum over a base-field scalar orbit equals either $q-1$ or $-1$; hence every $\eta_i$ is an integer. More precisely, additive-character orthogonality gives, for $a\ne0$,
\begin{equation}\label{v95:eq-period-hyperplane}
\boxed{qN_H(a)=n+g\eta(a).}
\end{equation}
The factor $g$ is indispensable: under the multiplication map $\F_q^\times\times H\to U$, every element of $U$ has exactly $g$ preimages.

With this character choice, \(\sum_{t\in\F_q^\times}\psi(t\Tr(ax))\) equals \(q-1\) when \(\Tr(ax)=0\) and \(-1\) otherwise, proving both integrality and \eqref{v95:eq-period-hyperplane}.
Replacing \(\Psi(x)\) by \(\Psi(bx)\), with \(b\in L^\times\), also replaces the intersection number in that identity by \(N_H(ab)\); one cannot change only the character while leaving its parameter fixed.
For the trace-code word \(c_a=(\Tr_{L/\F_q}(a\beta^{-j}))_{j=0}^{n-1}\),
\begin{equation}\label{v96:eq-period-weight}
\operatorname{wt}(c_a)=n-N_H(a)
 =\frac{(q-1)n-g\eta(a)}q\qquad(a\in L^\times).
\end{equation}
Thus individual period values give the weights of the trace code itself, whereas moments of the intersection numbers detect short subgroup relations, hence the minimum distance of the dual code; the latter identification is Proposition~\ref{v96:prop-relation-kernel}.

\begin{theorem}[The third period moment of a cyclotomic intersection]\label{v95:thm-period-third}
With the above notation,
\begin{equation}\label{v95:eq-period-third}
\boxed{\sum_{i=1}^e\eta_i=-1,\qquad
\sum_{i=1}^e\eta_i^2=Q-M,\qquad
\sum_{i=1}^e\eta_i^3=Q\mathcal C(U)-M^2.}
\end{equation}
Consequently, if $g=1$ and $k\ge3$, the identity
\[\boxed{\sum_i\eta_i^3=Q(q-2)-M^2}\]
is equivalent to $\lambda_q(H)\ge4$.
\end{theorem}
\begin{proof}
Write $\eta(a)=\sum_{u\in U}\Psi(au)$ and set $\eta(0)=M$. For each integer $r\ge1$, expanding the power and applying additive-character orthogonality gives
\[
\sum_{a\in L}\eta(a)^r=Q\,\#\{(u_1,\ldots,u_r)\in U^r:u_1+\cdots+u_r=0\}.
\]
Each $\eta_i$ occurs $M$ times as $a$ ranges over $L^\times$. For $r=1$ the count on the right is zero; for $r=2$, since $-1\in U$, there are $M$ zero-sum pairs; for $r=3$, after choosing the first element $u$, write the other two as $uv$ and $-u(1+v)$, obtaining $M\mathcal C(U)$ triples. Subtract the $a=0$ term $M^r$ and divide by $M$ to obtain the three identities. The final assertion follows from Theorem~\ref{v95:thm-triple-cyclotomic}.
\end{proof}
Consider the \textbf{Cayley graph} on the additive group $(L,+)$ with connection set $U$. Its nontrivial adjacency eigenvalues are precisely these $\eta_i$, each occurring $M$ times when indexed by cosets. Thus the third moment counts both closed walks of length three and three-term relations within the multiplicative subgroup. No strongly regular graph assumption is needed. Further connections among cyclic arcs of Singer type, three-valued Gaussian periods, and strongly regular graphs are developed in \cite{MomiharaXiang2022Subgroups}.

\paragraph{What real and integral spectra detect.}
The preceding integrality argument used \(\F_q^\times\subseteq U\). For a general multiplicative subgroup, the exact condition involves only prime-field scalars.
\begin{proposition}[Real and integral spectra of Gaussian periods]\label{v106:prop-real-integral-spectrum}
Let \(K=\F_Q\) have characteristic \(p\), let \(H\le K^\times\), and fix a nontrivial additive character \(\psi\). The translation operator
\[
(T_Hf)(x)=\sum_{h\in H}f(x+h)
\]
has eigenvalues \(\lambda_a=\sum_{h\in H}\psi(ah)\), \(a\in K\), in the additive-character basis. Moreover,
\[
\begin{aligned}
\text{every }\lambda_a\in\mathbb R&\quad\Longleftrightarrow\quad-1\in H,\\
\text{every }\lambda_a\in\Z&\quad\Longleftrightarrow\quad\F_p^\times\subseteq H.
\end{aligned}
\]
\end{proposition}
\begin{proof}
Substitution of \(x\mapsto\psi(ax)\) gives the eigenvalue formula. Complex conjugation sends \(\lambda_a\) to \(\lambda_{-a}\). If every value is real, the indicator functions of \(H\) and \(-H\) have the same Fourier transform. Inversion gives \(H=-H\), equivalently \(-1\in H\). The converse follows by changing variables.

All values belong to \(\Q(\zeta_p)\) and are algebraic integers. The automorphism \(\zeta_p\mapsto\zeta_p^c\), for \(c\in\F_p^\times\), sends \(\lambda_a\) to \(\lambda_{ca}\). If every value is an integer, each is fixed; Fourier inversion gives \(cH=H\), hence \(\F_p^\times\subseteq H\). Conversely, this inclusion fixes each value under all these automorphisms, so it is a rational algebraic integer, and therefore an integer. Alternatively, partition \(H\) into \(\F_p^\times\)-orbits; each orbit sum is \(p-1\) or \(-1\).
\end{proof}
\begin{example}[Real spectra, integral spectra, and additive closure]\label{v106:ex-small-spectra}
For \(H=\{1,4\}\subseteq\F_5\), the spectrum is
\[
2^{(1)},\qquad\left(\frac{-1+\sqrt5}{2}\right)^{(2)},\qquad
\left(\frac{-1-\sqrt5}{2}\right)^{(2)}.
\]
Parenthesized superscripts indicate multiplicities. Since \(-1\in H\) but \(\F_5^\times\nsubseteq H\), the spectrum is real but not integral.
In \(\F_9=\F_3(i)\), with \(i^2=-1\), the square subgroup is \(H=\{\pm1,\pm i\}\). Write \(a=s+ti\). Using the canonical additive character and \(\Tr(s+ti)=2s\), we get
\[
\lambda_a=A(2s)+A(t),\qquad
A(z)=\begin{cases}2,&z=0,\\-1,&z\ne0.\end{cases}
\]
Thus the spectrum is \(4^{(1)},1^{(4)},(-2)^{(4)}\). This time \(\F_3^\times\subseteq H\), so the spectrum is integral. Yet \(H\cup\{0\}\) has five elements and cannot be a subfield. The spectrum detects stability under prime-field scaling, but does not establish additive closure.
\end{example}

\subsection{Higher factorial moments and recovery of relation length}\label{v95:subsec-factorial-moments}
The third moment detects all three-term relations. For larger subsets, it is more useful to count \emph{factorial moments} of hyperplane intersections than to list relations one by one. For $0\le r\le k$ define
\[
B_r=\frac{q^{k-r}-1}{q-1},\qquad
\mathfrak M_s(H)=\sum_{[a]\in L^\times/\F_q^\times}\binom{N_H(a)}s,
\qquad 0\le s\le k.
\]
The number $B_r$ counts the hyperplanes containing a fixed $r$-dimensional subspace; in particular, $B_k=0$.
\begin{theorem}[Factorial hyperplane moments and subset ranks]\label{v95:thm-factorial-rank}
For every $0\le s\le k$,
\begin{equation}\label{v95:eq-factorial-rank}
\boxed{\mathfrak M_s(H)=\sum_{\substack{S\subseteq H\\|S|=s}}B_{\dim_{\F_q}\langle S\rangle}.}
\end{equation}
Suppose $2\le s\le k$ and $\lambda_q(H)\ge s$. If $R_s(H)$ is the number of $s$-element subsets of rank $s-1$, then
\begin{equation}\label{v95:eq-rank-step}
\boxed{R_s(H)=q^{s-k}\left[\mathfrak M_s(H)-\binom ns B_s\right].}
\end{equation}
Thus one may test $R_s$ successively for $s=2,3,\ldots,k$. The first index $s$ with $R_s>0$ equals $\lambda_q(H)$; if all these numbers vanish, $\lambda_q(H)=k+1$.
\end{theorem}
\begin{proof}
Double-count pairs $(S,V)$ where $S$ is an $s$-element subset of $H$ and $V$ is a trace hyperplane containing $S$. For fixed $S$ of rank $r$, there are exactly $B_r$ such hyperplanes, proving the first identity. If every $(s-1)$-element subset is independent, each $s$-element subset has rank either $s$ or $s-1$. Since $B_{s-1}-B_s=q^{k-s}$, the first identity becomes
\[\mathfrak M_s(H)=\binom ns B_s+q^{k-s}R_s(H).\]
As shorter dependencies have been excluded, $R_s>0$ is equivalent to the shortest relation having length $s$.
\end{proof}
\begin{corollary}[The first $k$ period moments determine $\lambda_q(H)$]\label{v95:cor-period-determines}
Put $h=n/g$ and $N_i=(n+g\eta_i)/q$. Then
\begin{equation}\label{v95:eq-moments-from-periods}
\boxed{\mathfrak M_s(H)=h\sum_{i=1}^{e}\binom{N_i}s\qquad(0\le s\le k).}
\end{equation}
Consequently the full distribution of Gaussian period values, and even just the power sums $\sum_i\eta_i^j$ for $0\le j\le k$, determines $\lambda_q(H)$ exactly.
\end{corollary}
\begin{proof}
Each $U$-coset contains $M/(q-1)=h$ projective parameters in $L^\times/\F_q^\times$, all with hyperplane intersection number $N_i$ by \eqref{v95:eq-period-hyperplane}. This gives the formula. The binomial coefficient $\binom{(n+gT)/q}s$ is a polynomial in $T$ of degree at most $s$, so its total sum is determined by power sums of degrees $0$ through $s$. Now apply Theorem~\ref{v95:thm-factorial-rank}.
\end{proof}
This exact recovery principle does \emph{not} imply that high period moments are easy to evaluate: they encode substantial additive information. Unlike the easily computed multiplicative orders, high moments may be as difficult to determine as the short relations themselves. It is therefore natural to begin with the order restrictions, the third moment, and special structure before resorting to higher-order data.

\subsection{Base-field descent, a seven-element example, and Fourier minors}\label{v95:subsec-examples-fourier}
The length of a shortest relation is unchanged by certain extensions of scalars. If $q=p^f$, $e_0=\ord_n(p)$, and $a=\gcd(e_0,f)$, then
\[
\F_{p^{e_0}}\cap\F_q=\F_{p^a},\qquad
\F_{p^{e_0}}\F_q=\F_{q^k},\qquad k=e_0/a.
\]
The two finite fields are linearly disjoint over their intersection. Consequently, any finite subset of the subgroup of order $n$ is linearly independent over $\F_{p^a}$ if and only if it is independent over $\F_q$. Hence
\begin{equation}\label{v95:eq-descent-lambda}
\boxed{\lambda_{p^a}(H)=\lambda_q(H).}
\end{equation}
This corresponds to the smallest coefficient field in Section~\ref{v84:subsec-code-descent}.

\leavevmode\phantomsection\label{v96:rem-scalar-extension}\ignorespaces
This identity compares the same subgroup inside an algebraic closure of one fixed characteristic \(p\); linear disjointness ensures that extending scalars introduces no new dependent subsets.
The following example of order seven instead changes the characteristic from two to eleven, comparing subgroups in different characteristics rather than a scalar extension of this kind.
It shows that the same multiplicative orders and cyclotomic-coset partition can accompany different additive relations, without contradicting \eqref{v95:eq-descent-lambda}.

\begin{example}[The same cyclotomic order data with different three-term relations]\label{v95:ex-seven}
Take $n=7$. For both $q=2$ and $q=11$ one has $k=3$, and the two nonzero Frobenius cosets modulo seven are $\{1,2,4\}$ and $\{3,5,6\}$.
Over $\F_8$, the subgroup is $H=\F_8^\times$. The sum of two distinct nonzero elements is again in $H$, so $\lambda_2(H)=3$.
Over $\F_{11}$ one has
\[
X^7-1=(X-1)(X^3+5X^2+4X-1)(X^3-4X^2-5X-1).
\]
Let $\beta$ be a root of $f(X)=X^3+5X^2+4X-1$. It has order seven. From $\beta^3=1-4\beta-5\beta^2$, the coordinate vectors of $1,\beta,\ldots,\beta^6$ in the basis $(1,\beta,\beta^2)$ are
\[
(1,0,0),\ (0,1,0),\ (0,0,1),\ (1,7,6),\ (6,10,10),\ (10,10,4),\ (4,5,1).
\]
Exponent translation and the $11$-Frobenius action $j\mapsto4j\pmod7$ partition the $\binom73=35$ exponent triples into three orbits. The determinants corresponding to the representative triples $(0,1,2)$, $(0,1,3)$, and $(0,1,5)$ are $1,6,4$ modulo eleven, respectively. All are nonzero, and hence $\lambda_{11}(H)=4$.
Moreover, $Q=1331$, $M=70$, $e=19$, and the Gaussian period values and their multiplicities are
\[
\begin{array}{c|rrr}
\eta_i&-7&4&15\\\hline
\text{multiplicity}&10&6&3
\end{array}
\]
They satisfy $\sum\eta_i^3=7079=1331\cdot9-70^2$, giving $\mathcal C(U)=9=q-2$. Thus \emph{the same multiplicative orders and Frobenius coset partition do not determine $\lambda$}: what changes is the joint additive position of the subgroup elements.
\end{example}

\leavevmode\phantomsection\label{v96:rem-two-fourier}\ignorespaces
Two Fourier descriptions must also be distinguished here.
The Gaussian periods above use complex-valued characters of the additive group \((L,+)\): orthogonality counts intersections, and the periods are eigenvalues of an additive Cayley graph.
The matrices below instead use root-of-unity characters of the cyclic group \(\Z/n\Z\), with minors first considered in a characteristic-zero cyclotomic field and then reduced to characteristic \(p\).
The former frequencies are parameterized by the trace pairing; the latter indices are organized into cyclotomic cosets, and finite-field-valued inversion requires \(p\nmid n\).
When the latter row indices form a Frobenius orbit, the reduced submatrix is a Moore matrix testing base-field linear independence.
This connects with the comparison at~\ref{v96:rem-vandermonde-moore}; it does not identify complex character sums with nonsingularity of finite-field minors.

For fixed $n$, characteristic-dependent behavior can also be explained through Fourier minors over a cyclotomic number field. Put $\Gamma=\langle q\bmod n\rangle$, with $|\Gamma|=k$, and let $S\subset\Z/n\Z$ have size $k$. Consider
\[
D_{\Gamma,S}(\zeta_n)=\det(\zeta_n^{\gamma j})_{\gamma\in\Gamma,\ j\in S}.
\]
Choose a prime ideal above $p$ in the cyclotomic integer ring and identify the image of $\zeta_n$ with a primitive $n$th root of unity in characteristic $p$. Ordering the rows as $1,q,\ldots,q^{k-1}$, the reduced minor is the Moore determinant of the corresponding $k$ subgroup elements and hence tests linear independence. If the cyclotomic integer is nonzero in characteristic zero, its integer norm is nonzero. Each conjugate determinant has absolute value at most $k!$, giving an absolute norm bound $(k!)^{\varphi(n)}$. Accordingly, \emph{if all these minors are nonzero in characteristic zero}, the inequality $p>(k!)^{\varphi(n)}$ ensures $\lambda_q(H)=k+1$. When $n=\ell$ is prime, Chebotarëv's nonvanishing theorem for all minors of a prime-order Fourier matrix guarantees the hypothesis \cite{Tao2005Subgroups,EmmrichKunis2025Subgroups}, yielding the sufficient, though very coarse, condition
\begin{equation}\label{v95:eq-chebotarev-bound}
\boxed{n=\ell\text{ prime},\quad p>(k!)^{\ell-1}\quad\Longrightarrow\quad\lambda_q(H)=k+1.}
\end{equation}
This deduction uses the classical Fourier-minor theorem, which is not proved here. At composite $n$, the possibility of minors vanishing already in characteristic zero must not be ignored.

From Theorem~\ref{v95:thm-factorial-rank}, short linear relations can be analyzed in successive layers: multiplicative orders provide quick necessary bounds; one cyclotomic intersection number detects all three-term dependencies; and factorial moments of Gaussian periods recover dependencies of greater length. Group theory, additive counting, and character sums are linked by the linear position of the same subgroup.

\paragraph{Subfield translation covers and quotient fibers.}
Trace hyperplanes record intersections of multiplicative subgroups with additive subspaces. When a subfield is the translation set, covering is likewise determined by projection fibers.
\begin{proposition}\label{v106:prop-subfield-cover}
Let \(E=\F_{q^2}\), \(K=\F_q\), and \(d\mid q^2-1\), and put \(H=(E^\times)^d\), \(H_0=H\cup\{0\}\). Choose a nonzero \(K\)-linear map \(\pi:E\to K\) with kernel \(K\). Then
\[
\boxed{K+H_0=E\quad\Longleftrightarrow\quad\pi(H_0)=K.}
\]
A necessary condition for covering is
\begin{equation}\label{v106:eq-cover-bound}
\boxed{d+\gcd(d,q+1)\le q+1.}
\end{equation}
If covering holds and equality holds in this bound, every \(x\notin K\) has a unique expression \(x=a+h\), where \(a\in K\) and \(h\in H\).
\end{proposition}
\begin{proof}
The condition \(x\in K+H_0\) means that \(\pi(x)=\pi(h)\) for some \(h\in H_0\), proving the first assertion. The intersection formula for cyclic subgroups gives
\[
|H\cap K^\times|=\gcd\left(\frac{q^2-1}{d},q-1\right)
=\frac{q-1}{d}\gcd(d,q+1).
\]
Every nonzero projection value must be attained on \(H\setminus K\), so
\[
q-1\le |H\setminus K|=\frac{q-1}{d}\bigl(q+1-\gcd(d,q+1)\bigr).
\]
This is the stated bound. Under the covering hypothesis, equality means that each nonzero projection fiber contains exactly one \(h\in H\), which also uniquely determines \(a=x-h\).
\end{proof}
This projection bound and the following construction appear in~\cite[Proposition 5.3 and Example 5.4]{ZhouWu2026Covering}; the proof uses only quotient spaces and cyclic groups.
\begin{example}[The same counting margin, different fiber distributions]\label{v106:ex-cover-5-11}
Suppose that \(q\) is odd and \(q\equiv2\pmod3\). Choose a cube root of unity \(\zeta\notin K\), and put \(T=(K^\times)^2\), \(H=T\langle\zeta\rangle\). Since \(3\nmid q-1\), we have \(|H|=3(q-1)/2\), with index \(d=2(q+1)/3\). In the basis \((1,\zeta)\), take \(\pi(a+b\zeta)=b\). The relation \(\zeta^2=-1-\zeta\) gives
\[
\pi(H_0)=\{0\}\cup T\cup(-T).
\]
Thus covering holds exactly when \(-1\) is a nonsquare, equivalently \(q\equiv3\pmod4\).
\begin{center}\small
\begin{tabular}{@{}cccccc@{}}\toprule
\(q\)&\(d\)&\(|H|\)&\(d+\gcd(d,q+1)\)&\(|\pi(H_0)|\)&Points missed\\\midrule
11&8&15&12&11&0\\
5&4&6&6&3&10\\\bottomrule
\end{tabular}\end{center}
Both examples attain equality in the necessary bound. In the second, \(T=-T=\{1,4\}\), and the projection misses \(2,3\), each accounting for an additive coset of size five. Having enough preimages in total does not distribute them among distinct fibers. Example~\ref{v106:ex-coset-character} will compare the positions of subfield cosets again through character sums.
\end{example}

\section{Multiple Jacobi Sums and Weil-Type Bounds for Diagonal Equations}
\label{sec:diagonal-jacobi}

In this section take \(k=\F_q\), retaining the convention that every multiplicative character is zero at zero.

\subsection{Exact Counts for Two-Variable Equations}

\begin{proposition}\label{prop:diagonal-jacobi-count}
Let \(a,b\in k^\times\) and \(u,v\ge1\), and put \(d=\gcd(u,q-1)\) and \(e=\gcd(v,q-1)\).
Let \(\mathcal X_d\) be the set of multiplicative characters satisfying \(\chi^d=1\), and let \(\varepsilon\) be the trivial character.
For \(c\ne0\), the number of solutions of \(ax^u+by^v=c\) is
\[
q+\sum_{\substack{\chi\in\mathcal X_d\setminus\{\varepsilon\}\\
                 \lambda\in\mathcal X_e\setminus\{\varepsilon\}}}
\chi(c/a)\lambda(c/b)J(\chi,\lambda).
\]
For \(c=0\), the number of solutions is
\[
q+(q-1)\sum_{\chi\in\mathcal X_{\gcd(d,e)}\setminus\{\varepsilon\}}\chi(-b/a).
\]
\end{proposition}
\begin{proof}
Character orthogonality on a cyclic group gives, for every \(z\in k\), including \(z=0\),
\[
|\{x:x^u=z\}|=1+\sum_{\chi\in\mathcal X_d\setminus\{\varepsilon\}}\chi(z).
\]
Thus the required solution count is
\[
\sum_{t\in k}
\left(1+\sum_{\chi\ne\varepsilon}\chi(t/a)\right)
\left(1+\sum_{\lambda\ne\varepsilon}\lambda((c-t)/b)\right),
\]
where \(\chi\) and \(\lambda\) run through the nontrivial characters in the two respective sets above.
The constant term contributes \(q\), and the sums of terms with only one nontrivial character vanish.
If \(c\ne0\), substituting \(t=cz\) in the cross terms gives the first formula.
If \(c=0\), a cross term becomes \(\chi(a^{-1})\lambda(-b^{-1})\sum_{t\ne0}(\chi\lambda)(t)\),
which is nonzero only when \(\lambda=\chi^{-1}\) and then equals \((q-1)\chi(-b/a)\).
These common characters form exactly \(\mathcal X_{\gcd(d,e)}\), proving the second formula.
\end{proof}

Although the solution count is an integer, it can be computed first in the character basis and then recovered through cancellation between conjugate characters.
Section~\ref{sec:hasse-davenport} proves the Hasse--Davenport formula describing lifts of existing characters under base-field extension;
if \(\gcd(u,q^n-1)\) or \(\gcd(v,q^n-1)\) increases, newly appearing characters must also be included.
Thus a solution count over the base field alone does not determine the counts over all extensions.

Restricting variables to a multiplicative subgroup gives a Jacobi-sum count by expanding its indicator function. Removing the diagonal then imposes distinctness of roots; see Exercise~\ref{v85:ex-restricted-sums}. Products of distinct elements can instead be treated directly by sums of cyclic-group exponents, as in Exercise~\ref{v85:ex-distinct-products}.

\subsection{Uniform subset products and joint fibers}
\begin{proposition}[Subset-product distributions in a cyclic group]\label{v150:prop-subset-products}
Let $G$ be cyclic of order $N>1$, let $1\le r<N$, and put
\[
 A_r(c)=\#\{S\subseteq G:|S|=r,\ \textstyle\prod_{s\in S}s=c\}.
\]
Then $A_r(c)$ is independent of $c\in G$ if and only if $\gcd(r,N)=1$; its uniform value is $N^{-1}\binom Nr$. In general,
\begin{equation}\label{v150:eq-subset-products}
 A_r(c)=\frac1N\sum_{d\mid\gcd(N,r)}
 (-1)^{r+r/d}\binom{N/d}{r/d}
 \sum_{\substack{\chi\in\widehat G\\\ord\chi=d}}\chi(c)^{-1}.
\end{equation}
\end{proposition}
\begin{proof}
If $\gcd(r,N)=1$, choose $u^r=c'/c$ for any $c,c'$. Then $S\mapsto uS$ is a bijection between the two fibers.
For a complex character $\chi$ of order $d$, each $d$th root of unity occurs $N/d$ times. Hence
\[
 \sum_cA_r(c)\chi(c)
 =[T^r]\prod_{x\in G}(1+T\chi(x))
 =[T^r](1-(-T)^d)^{N/d}.
\]
This is zero when $d\nmid r$, and otherwise equals
$(-1)^{r+r/d}\binom{N/d}{r/d}$. Fourier inversion gives~\eqref{v150:eq-subset-products}.
If $\gcd(r,N)>1$, choose a divisor $d>1$ of that greatest common divisor. The corresponding nontrivial Fourier coefficient is nonzero, excluding uniformity.
\end{proof}

\begin{example}[Nonempty fibers and uniform fibers]\label{v150:ex-subset-products-five}
In $\F_5^\times$, two-element subsets have product $1$ or $4$ once each, and product $2$ or $3$ twice each. Three-element subsets attain every nonzero product once.
Exercise~\ref{v85:ex-distinct-products} ensures that every specified product occurs; the Fourier coefficients above further decide whether the fibers have equal size.
\end{example}

Deleting one point deletes one factor from the generating function. Theorem~\ref{v150:thm-ers-design} in Chapter~17 uses precisely
$[T^r](1-(-T)^{q-1})/(1+T)=(-1)^r$, for $1\le r<q-1$.
Prescribing the root polynomial at both zero and one instead gives the joint fibers $B_r(u,v)$ of~\eqref{v150:eq-joint-products}.

\begin{proposition}[A Jacobi recurrence for joint-product fibers]\label{v150:prop-joint-products}
For $0\le r\le q-2$ and multiplicative characters $\chi,\lambda$, define
\[
 \mathcal F_r(\chi,\lambda)
 =\sum_{u,v\in\F_q^\times}B_r(u,v)\chi(u)\lambda(v).
\]
Then
\begin{align}
 \mathcal F_r(\chi,\lambda)
 &=[T^r]\prod_{x\ne0,1}(1+T\chi(x)\lambda(1-x)),
 \label{v150:eq-joint-generating}\\
 \mathcal F_0&=1,\qquad
 r\mathcal F_r=\sum_{j=1}^r(-1)^{j-1}
 J(\chi^j,\lambda^j)\mathcal F_{r-j}.
 \label{v150:eq-joint-newton}
\end{align}
Double Fourier inversion therefore determines every $B_r(u,v)$.
\end{proposition}
\begin{proof}
In the product expansion, a subset contributes $\chi(\prod x)\lambda(\prod(1-x))$, proving the first identity.
The $j$th power sum of these weights is $J(\chi^j,\lambda^j)$, including trivial power characters, which are still extended by zero at zero. The second identity is Newton's identity over the complex numbers.
\end{proof}

For $r=2$, $B_2(u,v)$ is one or zero according as
\[
 X^2-(1+u-v)X+u
\]
has two distinct roots in $\F_q$; $u,v\ne0$ already exclude zero and one as roots.
In odd characteristic, this means that the discriminant is a nonzero square. In characteristic two, writing $s=1+u-v$, it means $s\ne0$ and $\Tr_{q/2}(u/s^2)=0$.
For $r=3$, $B_3(u,v)$ counts the $w\in\F_q$ for which
\[
 X^3-(1+w-u-v)X^2+wX-u
\]
splits over $\F_q$ with distinct roots.
Root sets turn the two product constraints into low-degree splitting conditions; the Jacobi recurrence organizes all fibers at once.

\subsection{Multiple Jacobi Sums and the Product Character}

A two-variable Jacobi sum combines two multiplicative conditions with an additive constraint. The same structure holds for several variables;
the square-root bound is determined by whether the product of the multiplicative characters is trivial.

\leavevmode\phantomsection\label{def:multiple-jacobi-sum}\ignorespaces
For multiplicative characters \(\chi_1,\ldots,\chi_n\), define the \textbf{multiple Jacobi sum}
\[
J(\chi_1,\ldots,\chi_n)
=\sum_{\substack{y_1,\ldots,y_n\in k\\y_1+\cdots+y_n=1}}
  \prod_{i=1}^n\chi_i(y_i).
\]
For \(n=1\) this sum equals one; for \(n=2\) it is the two-variable Jacobi sum of Section~\ref{sec:gauss-jacobi}.

\begin{proposition}[Gauss-sum expressions for multiple Jacobi sums]\label{prop:multiple-jacobi-gauss}
Let \(n\ge1\), suppose that all \(\chi_1,\ldots,\chi_n\) are nontrivial, and write \(\rho=\chi_1\cdots\chi_n\).
Fix a nontrivial additive character \(\psi\) of \(k\). Then
\[
J(\chi_1,\ldots,\chi_n)=
\begin{cases}
\displaystyle\frac{\prod_iG(\chi_i,\psi)}{G(\rho,\psi)},&\rho\ne\varepsilon,\\[6pt]
\displaystyle-\frac1q\prod_iG(\chi_i,\psi),&\rho=\varepsilon.
\end{cases}
\]
Consequently,
\begin{equation}\label{eq:multiple-jacobi-size}
\left|J(\chi_1,\ldots,\chi_n)\right|=
\begin{cases}
q^{(n-1)/2},&\rho\ne\varepsilon,\\
q^{(n-2)/2},&\rho=\varepsilon.
\end{cases}
\end{equation}
The second case can occur only when \(n\ge2\).
\end{proposition}
\begin{proof}
Put
\[
S(t)=\sum_{y_1+\cdots+y_n=t}\prod_i\chi_i(y_i),\qquad J=S(1).
\]
For \(t\ne0\), substituting \(y_i=tz_i\) in all coordinates gives \(S(t)=\rho(t)J\).
If \(\rho\ne\varepsilon\), choose \(u\in k^\times\) such that \(\rho(u)\ne1\).
The simultaneous substitution \(y_i=uz_i\) in \(S(0)\) gives \(S(0)=\rho(u)S(0)\), hence \(S(0)=0\).
Grouping the product of Gauss sums by \(t=\sum_i y_i\) gives
\[
\prod_iG(\chi_i,\psi)=\sum_{t\in k}S(t)\psi(t)=JG(\rho,\psi).
\]
If \(\rho=\varepsilon\), then \(S(t)=J\) for every \(t\ne0\). On the other hand,
\[
\sum_{t\in k}S(t)=\prod_i\sum_{y\in k}\chi_i(y)=0,
\]
so \(S(0)=-(q-1)J\). Since \(\sum_{t\ne0}\psi(t)=-1\), in this case
\[
\prod_iG(\chi_i,\psi)=S(0)+J\sum_{t\ne0}\psi(t)=-qJ.
\]
The absolute values of both expressions follow from Theorem~\ref{thm:gauss}.
\end{proof}

\Needspace{16\baselineskip}
\subsection{Exact Counts and Weil-Type Bounds for Diagonal Equations}

\begin{theorem}[Character expansion and square-root bounds for diagonal equations]\label{thm:diagonal-weil-bound}
Let \(n\ge1\) and \(a_1,\ldots,a_n\in k^\times\), and let the positive integers \(d_i\) all divide \(q-1\).
Let the number of solutions of
\[
a_1x_1^{d_1}+\cdots+a_nx_n^{d_n}=b
\]
in \(k^n\) be \(N_b\). If \(b\ne0\), then
\begin{equation}\label{eq:multiple-diagonal-exact}
N_b=q^{n-1}+
\sum_{\substack{\chi_i\in\mathcal X_{d_i}\setminus\{\varepsilon\}\\1\le i\le n}}
\left(\prod_{i=1}^n\chi_i(b/a_i)\right)
J(\chi_1,\ldots,\chi_n).
\end{equation}
In particular,
\begin{equation}\label{eq:diagonal-weil-bound}
\left|N_b-q^{n-1}\right|
\le q^{(n-1)/2}\prod_{i=1}^n(d_i-1).
\end{equation}
More precisely, put \(M=\prod_i(d_i-1)\), and let \(M_0\) be the number of character tuples in this sum satisfying
\(\chi_1\cdots\chi_n=\varepsilon\). Then
\[
\left|N_b-q^{n-1}\right|
\le (M-M_0)q^{(n-1)/2}+M_0q^{(n-2)/2}.
\]
\end{theorem}
\begin{proof}
By the power-map fiber formula in the proof of Proposition~\ref{prop:diagonal-jacobi-count}, for each \(y\in k\),
\[
\#\{x:a_ix^{d_i}=y\}
=1+\sum_{\chi_i\in\mathcal X_{d_i}\setminus\{\varepsilon\}}\chi_i(y/a_i).
\]
The constant function one also equals one at zero and cannot be replaced by the trivial character \(\varepsilon\), which is zero there.
Hence
\[
N_b=\sum_{y_1+\cdots+y_n=b}
\prod_{i=1}^n\left(1+
\sum_{\chi_i\in\mathcal X_{d_i}\setminus\{\varepsilon\}}\chi_i(y_i/a_i)\right).
\]
Expand the product. The term choosing the constant in every coordinate contributes \(q^{n-1}\).
If nontrivial characters are chosen in a nonempty proper subset \(I\subsetneq\{1,\ldots,n\}\) of coordinates and constants in the others,
then, for fixed \((y_i)_{i\in I}\), the remaining coordinates have \(q^{n-|I|-1}\) choices. This part of the sum is
\[
q^{n-|I|-1}\prod_{i\in I}\sum_{y_i\in k}\chi_i(y_i/a_i)=0.
\]
Thus only the terms choosing nontrivial characters in every coordinate remain. In these terms, simultaneously substituting \(y_i=bz_i\)
gives~\eqref{eq:multiple-diagonal-exact}.
Each coefficient \(\prod_i\chi_i(b/a_i)\) has absolute value one,
and there are \(M\) character tuples. Applying~\eqref{eq:multiple-jacobi-size} term by term
first gives the bound involving \(M_0\) and then~\eqref{eq:diagonal-weil-bound}.
\end{proof}

If some \(d_i=1\), the sum is empty, so \(N_b=q^{n-1}\); this also follows directly because that variable is uniquely determined by the others.
For \(n=1\), the multiple Jacobi sum equals one, and~\eqref{eq:multiple-diagonal-exact} recovers the power-map fiber size.
For arbitrary positive exponents \(u_i\), take \(d_i\) to be \(\gcd(u_i,q-1)\):
the two power maps have the same image and fiber sizes on \(k\), so the diagonal equations have the same solution count.

These Weil-type bounds for diagonal equations \cite[Chapter 11]{Agou} come entirely from the square-root magnitude of Gauss sums.
Multiplicative characters expand the fibers of power maps, and the additive constraint combines these characters into Jacobi sums;
whether their product character is trivial determines the two scales of the error term.
Exact counts and square-root estimates are therefore two consequences of the same character expansion.

\section{Primitive Elements of Prescribed Trace}
\label{sec:primitive-prescribed-trace}

Primitivity is characterized by order in the multiplicative group, whereas prescribing the trace specifies an additive affine subspace.
Expanding both conditions in characters reduces their common solution count to multiplicative character sums over trace fibers.
These fiber sums are determined by a Gauss sum over the extension field and the restriction of the character to the base field.
This application connects the trace fibers of Chapter~\ref{chap:trace} with the square-root estimates of this chapter;
for a classical discussion, see \cite[Chapter 8]{Feng}.

\subsection{Multiplicative Symmetry and Fourier Reduction of Trace Distributions}\label{v111:sec-trace-fourier}
\phantomsection\label{v113:par-fourier-question}
Theorem~\ref{v111:thm-order-trace-symmetry} has already supplied bijections between certain trace fibers, making the counts constant on the corresponding multiplicative cosets. What remains is whether different cosets have the same count. The example in Chapter~\ref{chap:normal} shows that the total or average over nonzero traces alone neither excludes empty fibers nor determines uniformity.

With the base field again denoted by \(K=\F_q\), view the nonzero trace distribution as a function on \(K^\times\), using the Fourier convention of this chapter. Its coefficient at the trivial character is the sum of these counts; division by \(q-1\) gives the average. The nontrivial coefficients record differences that the average does not determine. Uniformity therefore requires all those coefficients to vanish, rather than merely a positive total. The known coset symmetry also rules out some frequencies in advance, expressing the remaining question on a smaller quotient group.

We first prove this connection in general and then apply it to exact-order and normality counts. The algebraic argument of Chapter~\ref{chap:normal} and character orthogonality in this chapter thus give the following chain:
\[
\boxed{\begin{gathered}
\text{Stabilizer of an exact-order set in the multiplicative group}\\
\Downarrow\\
\text{Scalar action preserving exact order and $k$-normality}\\
\Downarrow\\
\text{Multiplicative coset symmetry of the trace distribution}\\
\Downarrow\\
\text{Fourier support constraints and a uniformity criterion}
\end{gathered}}
\]
The first three steps use the cyclic group and the linearity of trace; the last uses only character orthogonality from this chapter.
Here Fourier transformation is on the multiplicative group of nonzero trace values, which must be distinguished from transformation on the additive group of the field.

\begin{proposition}[Coset invariance and Fourier support]\label{v111:prop-fourier-cosets}
Let $G$ be a finite abelian group, $H\le G$, and $f:G\to\C$. Use the unnormalized transform of Theorem~\ref{v90:thm-fourier-uncertainty}, and put
\[
\widehat f(\chi)=\sum_{x\in G}f(x)\overline{\chi(x)},\qquad
H^\perp=\{\chi\in\widehat G:\chi|_H=1\}.
\]
The following are equivalent: $f(hx)=f(x)$ for all $h\in H,x\in G$; $f$ is constant on each $H$-coset; and $\operatorname{supp}\widehat f\subseteq H^\perp$.
Moreover, $|H^\perp|=[G:H]$. If $c_1,\ldots,c_t$ represent the cosets, $b_j=f(c_j)$, and these equivalent conditions hold, then
\begin{equation}\label{v111:eq-quotient-fourier}
\widehat f(\chi)=
\begin{cases}
|H|\displaystyle\sum_{j=1}^t b_j\overline{\chi(c_j)},&\chi\in H^\perp,\\
0,&\chi\notin H^\perp.
\end{cases}
\end{equation}
\end{proposition}
\begin{proof}
The first two conditions are equivalent. If $f(hx)=f(x)$, substituting $x=hy$ gives $\widehat f(\chi)=\overline{\chi(h)}\widehat f(\chi)$.
When $\chi|_H\ne1$, choose $h$ with $\chi(h)\ne1$; the coefficient must vanish.
Conversely, character orthogonality gives the inversion formula
\[
f(x)=\frac1{|G|}\sum_{\chi\in\widehat G}\widehat f(\chi)\chi(x).
\]
If the support lies in $H^\perp$, replacing $x$ by $hx$ leaves the right-hand side unchanged.
Characters in $H^\perp$ are precisely the characters that descend through the quotient map to $G/H$, so $|H^\perp|=|G/H|$.
Finally, grouping the sum by $x=c_jh$ gives
\[
\widehat f(\chi)=\sum_{j=1}^t b_j\overline{\chi(c_j)}
\sum_{h\in H}\overline{\chi(h)}.
\]
The inner sum is $|H|$ or zero, proving~\eqref{v111:eq-quotient-fourier}.
\end{proof}
\begin{corollary}[Recovering the full symmetry group and measuring nonuniformity]\label{v112:cor-spectral-stabilizer}
Let \(f:G\to\C\) be as in Proposition~\ref{v111:prop-fourier-cosets}, and let \(\varepsilon\) be the trivial character. Define its \textbf{translation stabilizer} by
\[
\operatorname{Stab}(f)=\{h\in G:f(hx)=f(x)\text{ for every }x\in G\}.
\]
Then
\begin{equation}\label{v112:eq-spectral-stabilizer}
\boxed{\operatorname{Stab}(f)=
\bigcap_{\chi\in\operatorname{supp}\widehat f}\ker\chi.}
\end{equation}
The intersection over an empty support is \(G\). The quotient \(G/\operatorname{Stab}(f)\) has the smallest order among the group quotients through which \(f\) factors. If \(G\) is cyclic, then
\[
[G:\operatorname{Stab}(f)]
=\operatorname{lcm}_{\chi\in\operatorname{supp}\widehat f}\operatorname{ord}(\chi),
\]
with the empty least common multiple equal to one.
Writing \(\mu_f=|G|^{-1}\sum_{x\in G}f(x)\), one also has
\begin{equation}\label{v112:eq-fourier-nonuniformity}
\boxed{\sum_{x\in G}|f(x)-\mu_f|^2
=\frac1{|G|}\sum_{\chi\ne\varepsilon}|\widehat f(\chi)|^2.}
\end{equation}
If \(f\) is \(H\)-invariant, the last sum may be restricted to \(H^\perp\setminus\{\varepsilon\}\).
\end{corollary}
\begin{proof}
For each \(h\in G\), the substitution and inversion arguments in Proposition~\ref{v111:prop-fourier-cosets} show that \(f(hx)=f(x)\) for all \(x\) if and only if \(\chi(h)=1\) at every nonzero Fourier coefficient. This proves~\eqref{v112:eq-spectral-stabilizer}.
The function factors through \(G/H\) exactly when \(H\subseteq\operatorname{Stab}(f)\), proving minimality of the quotient. In a cyclic group, \(\ker\chi\) has index \(\operatorname{ord}(\chi)\); the index of an intersection of subgroups is the least common multiple of their indices.
Finally, Fourier inversion expresses \(f-\mu_f\) as \(|G|^{-1}\sum_{\chi\ne\varepsilon}\widehat f(\chi)\chi\). Expanding its squared absolute value and summing over \(G\), character orthogonality cancels all terms with distinct characters and contributes \(|G|\) to each diagonal term. This gives~\eqref{v112:eq-fourier-nonuniformity}.
\end{proof}
The actual support therefore recovers the full symmetry of the distribution, while the squared magnitudes of its nontrivial coefficients measure its deviation from uniformity. The number of nonzero coefficients alone need not equal the order of the smallest quotient.

\begin{corollary}[Quotient expansion of exact-order and normality counts]\label{v111:cor-trace-fourier}
Use the notation of Theorem~\ref{v111:thm-order-trace-symmetry}.
Set $f_{r,k}=N_{r,k}|_{K^\times}$. Then
\[
\operatorname{supp}\widehat f_{r,k}\subseteq H_r^\perp,
\qquad |\operatorname{supp}\widehat f_{r,k}|\le t.
\]
The nonzero trace distribution is uniform if and only if
\[
\widehat f_{r,k}(\chi)=0\qquad
(\chi\in H_r^\perp\setminus\{\varepsilon\}),
\]
where $\varepsilon$ is the trivial character. Thus at most $t-1$ nontrivial coefficients need to be tested.
For representatives $c_j$ of $K^\times/H_r$, put $b_j=N_{r,k}(c_j)$; formula~\eqref{v111:eq-quotient-fourier} gives every coefficient.
Furthermore, $b_j$ is exactly the number of $H_r$-orbits in $\mathcal A_{r,k}$ whose traces lie in $c_jH_r$.
\end{corollary}
\begin{proof}
Theorem~\ref{v111:thm-order-trace-symmetry} supplies $H_r$-invariance, so the preceding proposition gives the support constraint and the quotient formula.
Inversion shows that having no nontrivial frequency is equivalent to being constant, including the zero function.
The scalar action on nonzero elements is free, so every $H_r$-orbit has $|H_r|$ elements.
If $\Tr(\alpha)\ne0$, the map
\[
h\alpha\longmapsto h\Tr(\alpha)\qquad(h\in H_r)
\]
is a bijection from that orbit to its trace coset. Each such orbit contributes exactly one element to the fiber with trace $c_j$, giving the orbit count $b_j$.
\end{proof}

This reduces the calculation on $q-1$ nonzero trace values exactly to a quotient group of order $t$.
The support relation is an inclusion, not necessarily equality: even when $H_r$ comes from the full stabilizer of an exact-order set, the trace counts may have a larger symmetry group.
Zero trace is fixed under scalar multiplication and is absent from this quotient group; its fiber must be counted separately.

\begin{example}[The same permitted frequencies, different actual supports]\label{v111:ex-two-fourier-supports}
In Example~\ref{v111:ex-order52-traces}, $K=\F_5$, $n=4$ and $r=12$, with $H_{12}=\{1,4\}$ and $t=2$.
Let $\eta$ be the quadratic character of $K^\times$. Then $H_{12}^\perp=\{\varepsilon,\eta\}$, and the previously computed trace distribution gives
\[
f_{12,0}=4,\qquad f_{12,1}=2+2\eta.
\]
Consequently,
\[
\begin{array}{c|rr}
&\varepsilon&\eta\\\hline
\widehat f_{12,0}&16&0\\
\widehat f_{12,1}&8&8
\end{array}
\]
Both coefficients at characters of order four vanish. With coset representatives $1,2$, the quotient counts are $(4,4)$ and $(4,0)$; their two-point Fourier transforms, multiplied by $|H_{12}|=2$, recover the table.
The stabilizer specifies which frequencies are permitted, while the actual orbit counts determine whether their coefficients vanish.
The full stabilizers of these two functions are, respectively, \(K^\times\) and \(\{1,4\}\), by Corollary~\ref{v112:cor-spectral-stabilizer}. Their means are four and two, and their squared deviations from those means are zero and \(16=8^2/4\). Thus the same permitted support yields different smallest quotients and different deviations from uniformity.
\end{example}

Coset symmetry reduces the parameter set. If the question is only whether every nonzero trace occurs, a second moment can first bound the number of empty cosets.
\begin{proposition}[Second moments and empty trace cosets]\label{v115:prop-empty-cosets}
Let \(G=K^\times\), with \(H\le G\), and let a nonnegative counting function be constant on \(H\)-cosets, with values \(b_1,\ldots,b_t\), where \(t=[G:H]\). Put
\[
\mu=\frac1t\sum_{j=1}^t b_j>0,\qquad
V=\sum_{j=1}^t(b_j-\mu)^2,\qquad
z=\#\{j:b_j=0\}.
\]
Then
\begin{equation}\label{v115:eq-empty-bound}
z\le\left\lfloor\frac{tV}{t\mu^2+V}\right\rfloor.
\end{equation}
In particular, if \(t>1\) and \(V<t\mu^2/(t-1)\), every coset value is positive. If \(t=1\), positivity follows directly from \(\mu>0\).
For the unnormalized Fourier transform on \(G/H\),
\begin{equation}\label{v115:eq-quotient-variance}
V=\frac1t\sum_{\chi\ne\varepsilon}|\widehat b(\chi)|^2.
\end{equation}
\end{proposition}
\begin{proof}
There are \(t-z\) positive values. Cauchy--Schwarz gives
\[
(t\mu)^2=\left(\sum_{b_j>0}b_j\right)^2
\le(t-z)\sum_jb_j^2=(t-z)(t\mu^2+V).
\]
Rearrange and use integrality of \(z\). Under the stated strict inequality, the bound is less than one.
Equation~\eqref{v115:eq-quotient-variance} is Parseval's identity on the quotient group, also obtained directly from Corollary~\ref{v112:cor-spectral-stabilizer}.
\end{proof}

For \(N_{r,k}|_{K^\times}\), take \(H=H_r\). Each empty coset accounts for \(|H_r|\) nonzero trace values with no preimage; zero trace must still be treated separately.
The two quotient distributions in Example~\ref{v111:ex-two-fourier-supports} are \((4,4)\) and \((4,0)\). The former has \(V=0\); the latter has \(\mu=2,V=8\), and an empty coset occurs exactly at \(t\mu^2/(t-1)=8\).
The strictness of the threshold is therefore necessary.

When \(t=2\), the total and the single nontrivial Fourier coefficient recover both values. With more cosets, the second moment retains only the sum of squared coefficient magnitudes.
If the bound does not exclude an empty coset, one can compute the actual quotient counts or nontrivial coefficients; a positive average alone does not imply nonempty fibers.
The scalar stabilizer of the exact-order set specifies permitted frequencies. Further vanishing, and hence additional symmetry of the trace distribution, still depends on these counts.

We now return to ordinary primitive elements, using their multiplicative character expansion and Gauss sums to control the sizes of the fibers.

\subsection{The Primitivity Indicator}

Let \(k=\F_q\) and \(E=\F_{q^n}\), where \(n\ge1\), and write
\[
Q=q^n,\qquad M=Q-1,\qquad \vartheta(M)=\frac{\varphi(M)}{M}.
\]
All multiplicative characters, including the trivial character \(\varepsilon\), are assigned value zero at zero.
Let \(1_{\mathrm{prim}}\) denote the indicator of the primitive elements of \(E\).

\begin{proposition}[Multiplicative character expansion of primitivity]\leavevmode\phantomsection\label{prop:primitive-character-indicator}
For every \(x\in E\),
\begin{equation}\label{eq:primitive-character-indicator}
1_{\mathrm{prim}}(x)
=\vartheta(M)\sum_{d\mid M}\frac{\mu(d)}{\varphi(d)}
   \sum_{\substack{\chi\in\widehat{E^\times}\\\ord(\chi)=d}}\chi(x).
\end{equation}
Here \(\widehat{E^\times}\) denotes the multiplicative character group of \(E^\times\). When \(M=1\),
the right-hand side contains only the trivial-character term.
\end{proposition}
\begin{proof}
Both sides vanish at \(x=0\). For \(x\ne0\), choose \(E^\times=\langle z\rangle\) and write \(x=z^j\).
Since the character group is also cyclic of order \(M\), every character of square-free order decomposes uniquely as a product of characters of distinct prime orders. Hence
\[
\sum_{d\mid M}\frac{\mu(d)}{\varphi(d)}\sum_{\ord(\chi)=d}\chi(x)
=\prod_{\ell\mid M}
 \left(1-\frac1{\ell-1}\sum_{\ord(\chi)=\ell}\chi(x)\right),
\]
where \(\ell\) runs through the prime divisors of \(M\). Fix a character of order \(\ell\); its nontrivial powers give all characters of order \(\ell\).
The finite geometric-series formula gives
\[
\sum_{\ord(\chi)=\ell}\chi(z^j)
=\begin{cases}\ell-1,&\ell\mid j,\\-1,&\ell\nmid j.\end{cases}
\]
Thus, if \(\gcd(j,M)>1\), the product vanishes; if \(\gcd(j,M)=1\),
the product equals \(\prod_{\ell\mid M}\ell/(\ell-1)=\vartheta(M)^{-1}\).
These are precisely the cases where \(x\) is nonprimitive or primitive, respectively. When \(M=1\), the empty product is one, and the conclusion still holds.
\end{proof}

\subsection{Character Sums over Trace Fibers}

Fix a nontrivial additive character \(\psi\) of \(k\), and put
\(\Psi=\psi\circ\Tr_{E/k}\). Surjectivity of the trace makes \(\Psi\) nontrivial.
\leavevmode\phantomsection\label{def:trace-fiber-character-sum}\ignorespaces For \(a\in k\) and a multiplicative character \(\chi\) of \(E\), write
\[
U_a=\{x\in E:\Tr_{E/k}(x)=a\},\qquad
S_a(\chi)=\sum_{x\in U_a}\chi(x),\qquad \rho=\chi|_{k^\times}.
\]
Here \(U_a\) is a trace fiber, and \(S_a(\chi)\) is called a \textbf{character sum over a trace fiber}.
For the trivial character,
\begin{equation}\label{eq:trace-fiber-principal}
S_a(\varepsilon)=q^{n-1}-\delta_{a,0},
\end{equation}
where \(\delta_{a,0}\) is one when \(a=0\) and zero otherwise.

\begin{proposition}[Gauss-sum expressions for trace-fiber sums]\label{prop:trace-fiber-gauss}
If \(\chi\ne\varepsilon\), then
\begin{equation}\label{eq:trace-fiber-gauss}
S_a(\chi)=\frac{G_E(\chi,\Psi)}q
 \sum_{t\in k^\times}\rho(t)^{-1}\psi(-at).
\end{equation}
More explicitly, let \(\varepsilon_k\) be the trivial character of \(k^\times\). Then
\[
S_a(\chi)=\frac{G_E(\chi,\Psi)}q
\begin{cases}
q-1,&a=0,\ \rho=\varepsilon_k,\\
0,&a=0,\ \rho\ne\varepsilon_k,\\
-1,&a\ne0,\ \rho=\varepsilon_k,\\
\rho(-a)G_k(\rho^{-1},\psi),&a\ne0,\ \rho\ne\varepsilon_k.
\end{cases}
\]
The corresponding absolute values are therefore
\begin{equation}\label{eq:trace-fiber-absolute}
|S_a(\chi)|=
\begin{cases}
(q-1)q^{(n-2)/2},&a=0,\ \rho=\varepsilon_k,\\
0,&a=0,\ \rho\ne\varepsilon_k,\\
q^{(n-2)/2},&a\ne0,\ \rho=\varepsilon_k,\\
q^{(n-1)/2},&a\ne0,\ \rho\ne\varepsilon_k.
\end{cases}
\end{equation}
\end{proposition}
\begin{proof}
Expanding the indicator of the trace condition gives
\[
S_a(\chi)=\frac1q\sum_{t\in k}\psi(-at)
 \sum_{x\in E^\times}\chi(x)\Psi(tx).
\]
Since \(\chi\) is nontrivial, the inner sum for \(t=0\) vanishes; for \(t\ne0\), putting \(y=tx\)
makes the inner sum equal to \(\rho(t)^{-1}G_E(\chi,\Psi)\). This proves~\eqref{eq:trace-fiber-gauss}.

If \(a=0\), the base-field sum is determined by multiplicative character orthogonality.
If \(a\ne0\) and \(\rho\) is trivial, it equals \(-1\);
if \(a\ne0\) and \(\rho\) is nontrivial, substitution of \(u=-at\) turns it into
\(\rho(-a)G_k(\rho^{-1},\psi)\).
Finally, applying the Gauss-sum absolute-value formulas over \(E\) and \(k\) gives~\eqref{eq:trace-fiber-absolute}.
\end{proof}

The two error scales have a specific source: when the restriction of the character to \(k^\times\) is nontrivial,
the base-field sum contributes an additional \(\sqrt q\); when the restriction is trivial, this sum reduces to \(-1\) or \(q-1\).
Only characters of the latter type contribute to the zero-trace fiber, so it must be counted separately from nonzero-trace fibers.

\subsection{Counting Formulas and Existence Conditions}

\leavevmode\phantomsection\label{def:primitive-trace-count}\ignorespaces Define
\[
N_a=\#\{x\in E^\times:\ord(x)=Q-1,\ \Tr_{E/k}(x)=a\}.
\]
For a positive integer \(r\), let \(\omega(r)\) be its number of distinct prime factors,
and \(W(r)=2^{\omega(r)}=\sum_{d\mid r}|\mu(d)|\) its number of square-free divisors, with the convention \(W(1)=1\).
Also write
\[
R=\frac{Q-1}{q-1},\qquad B=W(R)-1,\qquad C=W(M)-W(R).
\]

\begin{theorem}[Counting primitive elements of prescribed trace]\label{thm:primitive-prescribed-trace}
For every \(a\in k\), the exact formula is
\begin{equation}\label{eq:primitive-trace-exact}
N_a=\vartheta(M)\left(q^{n-1}-\delta_{a,0}
 +\sum_{\substack{d\mid M\\d>1}}\frac{\mu(d)}{\varphi(d)}
       \sum_{\ord(\chi)=d}S_a(\chi)\right).
\end{equation}
The contribution of the trivial character is
\[
\vartheta(M)(q^{n-1}-\delta_{a,0})
=\frac{\varphi(M)}q+\vartheta(M)\left(\frac1q-\delta_{a,0}\right).
\]
The total contribution of nontrivial characters satisfies
\begin{align}
\bigl|N_0-\vartheta(M)(q^{n-1}-1)\bigr|
 &\le \vartheta(M)B(q-1)q^{(n-2)/2},\label{eq:primitive-trace-zero-bound}\\
\bigl|N_a-\vartheta(M)q^{n-1}\bigr|
 &\le \vartheta(M)\bigl(Bq^{(n-2)/2}+Cq^{(n-1)/2}\bigr)
 \quad(a\ne0).\label{eq:primitive-trace-nonzero-bound}
\end{align}
\end{theorem}
\begin{proof}
Sum~\eqref{eq:primitive-character-indicator} over \(U_a\)
and separate the trivial character using~\eqref{eq:trace-fiber-principal} to obtain the exact formula.

To estimate the remaining terms, first determine which characters are trivial on \(k^\times\).
Take \(E^\times=\langle z\rangle\). Then \(k^\times=\langle z^R\rangle\).
A character of order \(d\) satisfies \(\chi(z^R)=1\) if and only if \(d\mid R\),
because \(\chi(z)\) is a primitive \(d\)th root of unity.
Thus, for each \(d\mid R\), the restrictions of all \(\varphi(d)\) characters of order \(d\) are trivial;
for the remaining \(d\mid M\), all these restrictions are nontrivial.

Apply the triangle inequality to the exact formula, followed by~\eqref{eq:trace-fiber-absolute}.
For each order \(d\), the number \(\varphi(d)\) of characters cancels the denominator in the coefficient.
The total weight of the orders of nontrivial characters with trivial restriction is
\(\sum_{d\mid R,\,d>1}|\mu(d)|=B\), and that of the remaining orders is \(C\).
This gives the bounds for zero and nonzero trace, respectively.
\end{proof}

\begin{corollary}[Directly testable sufficient conditions]\label{cor:primitive-trace-existence}
If
\[
q^{n-1}>Bq^{(n-2)/2}+Cq^{(n-1)/2},
\]
then every \(a\in k^\times\) is the trace of a primitive element of \(E\).
In particular, for \(a\ne0\) there is the simpler estimate
\[
\bigl|N_a-\vartheta(M)q^{n-1}\bigr|
\le\vartheta(M)(W(M)-1)q^{(n-1)/2},
\]
so the sufficient condition can be simplified to
\begin{equation}\label{eq:primitive-trace-simple-criterion}
\boxed{\ q^{(n-1)/2}>W(q^n-1)-1.\ }
\end{equation}
For zero trace, if
\begin{equation}\label{eq:primitive-trace-zero-criterion}
q^{n-1}-1>(W(R)-1)(q-1)q^{(n-2)/2},
\end{equation}
then \(E\) contains a primitive element of trace zero.
\end{corollary}
\begin{proof}
Each strict inequality makes the lower bound for the count in the theorem positive.
Moreover, \(q^{(n-2)/2}\le q^{(n-1)/2}\) and \(B+C=W(M)-1\),
so the refined bound for nonzero trace implies the displayed simpler estimate; cancelling the positive factor gives~\eqref{eq:primitive-trace-simple-criterion}.
\end{proof}

These are sufficient conditions for existence. When they fail, the exact formula or other methods are still needed.
In particular, the zero-trace condition cannot be replaced by the nonzero-trace condition; see Exercise~\ref{ex:primitive-zero-trace-boundary} for the zero-trace obstruction in quadratic extensions.
When \(n=1\), the trace is the identity and \(N_a\) is the primitivity indicator of \(a\);
the formulas above remain valid.

\begin{lemma}[An elementary bound for square-free divisor counts]\label{v88:lem-squarefree-divisor-inequality}
For every positive integer \(r\), \((W(r)-1)^2<2r\); if \(r\) is odd,
then the stronger inequality \((W(r)-1)^2<r\) holds.
\end{lemma}
\begin{proof}
Put \(s=\omega(r)\), and let \(P_s\) be the product of the \(s\) smallest primes, with \(P_0=1\).
Since \(r\ge P_s\), it suffices to prove \((2^s-1)^2<2P_s\).
For \(s=0,1,2,3\), substitute \(P_s=1,2,6,30\).
For \(s\ge3\), one has \(P_{s+1}/P_s\ge7\), whereas
\((2^{s+1}-1)^2<5(2^s-1)^2\), so induction proves the inequality.
If \(r\) is odd, instead use the product of the \(s\) smallest odd primes.
The first four products are \(1,3,15,105\), all satisfying \((2^s-1)^2<P_s\).
Thereafter \(P_{s+1}/P_s\ge11\), and the same ratio estimate completes the induction.
\end{proof}

\begin{theorem}[Universal existence of primitive elements with nonzero trace]\label{v88:thm-primitive-nonzero-trace}
For every prime power \(q\) and every \(n\ge1\), there exists
\(\alpha\in\F_{q^n}\) satisfying
\[
\boxed{\quad\ord(\alpha)=q^n-1,\qquad\Tr_{q^n/q}(\alpha)\ne0.\quad}
\]
More precisely, with \(R=(Q-1)/(q-1)\) and
\(P_{\ne0}=\sum_{a\in k^\times}N_a\),
\begin{equation}\label{v88:eq-primitive-nonzero-trace-lower}
\boxed{\quad P_{\ne0}\ge\vartheta(M)\frac{q-1}{q}
\bigl(Q-(W(R)-1)\sqrt Q\bigr)>0.\quad}
\end{equation}
\end{theorem}
\begin{proof}
Sum the primitive indicator~\eqref{eq:primitive-character-indicator} over
all elements of nonzero trace. The trivial multiplicative character contributes
\(\vartheta(M)(q-1)Q/q\).
For a nontrivial multiplicative character \(\chi\), summing the trace-fiber
Gauss expression over \(a\ne0\) gives
\[
\sum_{\Tr(x)\ne0}\chi(x)
=-\frac{G_E(\chi,\Psi)}q\sum_{t\in k^\times}\chi(t)^{-1}.
\]
If \(\chi|_{k^\times}\ne1\), this vanishes. Otherwise its absolute value
is \((q-1)\sqrt Q/q\). Choose a generator \(g\) of \(E^\times\); then
\(k^\times=\langle g^R\rangle\). Hence \(\chi|_{k^\times}=1\) is equivalent to
\(\chi(g)^R=1\), or equivalently to \(\ord(\chi)\mid R\).
Grouping the primitive-indicator expansion by square-free orders cancels the
number \(\varphi(d)\) of characters of each order against its denominator.
Thus the total error has absolute value at most
\(\vartheta(M)(q-1)(W(R)-1)\sqrt Q/q\), proving the first inequality.
If \(q=2\), \(R=Q-1\) is odd, and the preceding lemma gives
\((W(R)-1)^2<R<Q\). If \(q\ge3\), then \(Q=(q-1)R+1>2R\), and the lemma
gives \((W(R)-1)^2<2R<Q\). The displayed lower bound is therefore positive.
\end{proof}

The theorem asks for \emph{some} nonzero trace, not a prescribed value.
Summing all nonzero trace fibers is easier than controlling each fiber
separately, and supplies exactly the existence statement used for the
special-degree primitive normal constructions in Chapter~\ref{chap:normal}.
For stronger classical results on prescribing a trace, see
\cite{CohenPrimitiveTrace1990}; the proof here uses only the character-sum
formula already established in this section.

\begin{example}[An improvement from restricted characters and an exact count]\leavevmode\phantomsection\label{exm:primitive-trace-small-fields}
Take \(q=3,n=3\). Then \(M=26\) and \(R=13\), so \(B=1\) and \(C=2\).
The sufficient condition for nonzero trace is \(9>\sqrt3+6\), and that for zero trace is \(8>2\sqrt3\).
Both hold, so every trace fiber of \(\F_{27}\) contains a primitive element.
Here the simplified sufficient condition would require \(3>W(26)-1=3\), so it is inconclusive;
retaining the character restrictions to the base field completes the argument.

Now take \(E=\F_8\) and \(k=\F_2\). Since \(|E^\times|=7\) is prime,
one is the only nonprimitive nonzero element. The two trace fibers each contain four elements, and
\(\Tr(0)=0\) and \(\Tr(1)=1\). Removing zero and one from their respective fibers gives
\[
N_0=N_1=3.
\]
Here \(\vartheta(7)=6/7\), and the two main terms in the theorem are \(18/7\) and \(24/7\),
with exact errors \(3/7\) and \(-3/7\), both bounded by \(6\sqrt2/7\).
\end{example}

Finally, a primitive element generates the entire extension: if \(x\) belonged to a proper subfield \(\F_{q^d}\),
then \(\ord(x)\le q^d-1<q^n-1\).
Its minimal polynomial therefore has degree \(n\), with coefficient \(-\Tr_{E/k}(x)\) at \(X^{n-1}\).
Counting primitive elements of prescribed trace also counts primitive polynomials with a prescribed coefficient of the next-to-leading term:
each Frobenius orbit has \(n\) roots, so the number of such monic primitive polynomials is exactly \(N_a/n\).

\section{The primitive normal basis theorem: a joint count of two generation conditions}
\label{v93:sec-primitive-normal-proof}

Section~\ref{v93:sec-normal-affine} expresses normality as nonvanishing in the local Frobenius components and uses weighted traces to describe affine fibers consisting entirely of normal elements.
Primitivity, in turn, requires an element to generate the multiplicative group of the field. We expand the two conditions in characters and use Gauss sums to control the number of elements satisfying both.
A normality-preserving scaling group first reduces the multiplicative prime factors that must be tested; a local density sieve retains the density of each remaining condition and handles the boundary cases of the basic estimate.
This approach is related to the classical character-sum proofs of the primitive normal basis theorem\cite{LS,CohenHuczynska2003}.

Throughout this section, write
\[
k=\F_q,\qquad E=\F_{q^n},\qquad Q=q^n,\qquad N=Q-1,\qquad \sigma(x)=x^q.
\]
Write \(n=p^a m\), where \(p=\Char k\) and \(p\nmid m\), and put
\[
F(T)=T^n-1=M(T)^{p^a},\qquad
M(T)=T^m-1=\prod_{i=1}^sP_i(T),\qquad d_i=\deg P_i.
\]
The \(P_i\) are distinct monic irreducible polynomials, so \(M=\operatorname{rad}(F)\) and \(\sum_i d_i=m\).
Let \(\operatorname{PN}_q(n)\) denote the number of primitive normal elements of \(E/k\), and retain the notation
\(\vartheta(u)=\varphi(u)/u\) and \(W(u)=2^{\omega(u)}\). Here \(s\) counts the distinct irreducible local components, whereas \(\omega(u)\) counts the distinct prime factors of an integer.

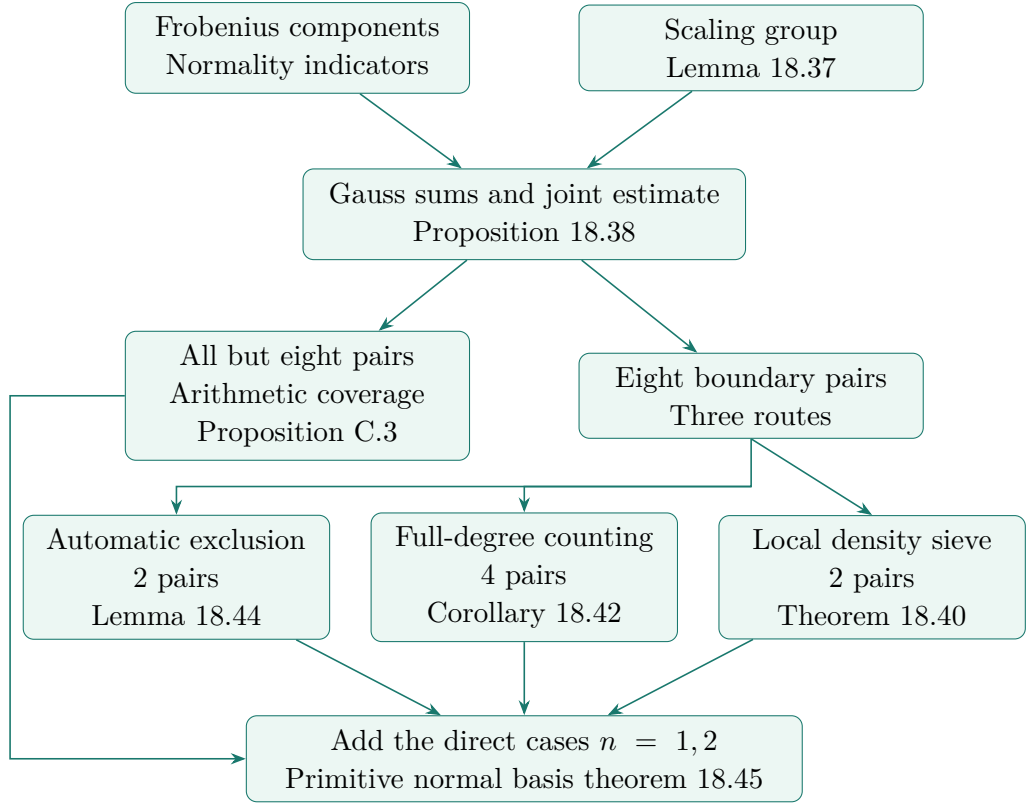
\begin{figure}[!htbp]
\centering
\begin{tikzpicture}[every node/.style={font=\small\linespread{1.05}\selectfont},
 bx/.style={draw=Jade,fill=CoverMint,rounded corners,align=center,text width=4.2cm,inner sep=5pt},
 ar/.style={-{Stealth[length=2mm]},semithick,Jade}]
\node[bx] (local) at (-3,0) {Frobenius components\\Normality indicators};
\node[bx] (scale) at (3,0) {Scaling group\\Lemma~\ref{v93:lem-scaling}};
\node[bx,text width=5.5cm] (basic) at (0,-2.2) {Gauss sums and joint estimate\\Proposition~\ref{v93:prop-basic}};
\node[bx] (generic) at (-3,-4.6) {All but eight pairs\\Arithmetic coverage\\Proposition~\ref{v93:prop-arithmetic-complete}};
\node[bx] (eight) at (3,-4.6) {Eight boundary pairs\\Three routes};
\node[bx,text width=3.7cm] (auto) at (-4.6,-7) {Automatic exclusion\\2 pairs\\Lemma~\ref{v93:lem-automatic}};
\node[bx,text width=3.7cm] (union) at (0,-7) {Full-degree counting\\4 pairs\\Corollary~\ref{v93:cor-quotient-full-degree-union}};
\node[bx,text width=3.7cm] (sieve) at (4.6,-7) {Local density sieve\\2 pairs\\Theorem~\ref{v93:thm-sieve}};
\node[bx,text width=7cm] (end) at (0,-9.4) {Add the direct cases \(n=1,2\)\\Primitive normal basis theorem~\ref{v93:thm-PNBT}};

\draw[ar] (local)--(basic); \draw[ar] (scale)--(basic);
\draw[ar] (basic)--(generic); \draw[ar] (basic)--(eight);
\draw[ar] (eight.south)--(3,-5.8)--(-4.6,-5.8)--(auto.north);
\draw[ar] (3,-5.8)--(0,-5.8)--(union.north);
\draw[ar] (eight.south)--(sieve.north);
\draw[ar] (auto)--(end); \draw[ar] (union)--(end); \draw[ar] (sieve)--(end);
\draw[ar] (generic.west)--(-6.8,-4.6)--(-6.8,-9.4)--(end.west);
\end{tikzpicture}
\caption{Dependencies of the primitive normal basis theorem. The estimate in a fixed weighted-trace fiber is a separate application, not a required step in this proof.}\label{v103:fig-pnbt-dependencies}
\end{figure}

Figure~\ref{v103:fig-pnbt-dependencies} lists the dependencies of the proof in this section. The basic estimate is only sufficient; the theorem for all parameters also needs the arithmetic coverage in the appendix and the treatment of all eight boundary pairs.
\FloatBarrier

\subsection{A scaling group that preserves normality}

\begin{lemma}[A normality-preserving scaling group]\label{v93:lem-scaling}
Put
\[
\mathcal C=\{c\in E^\times:c^{q-1}\in k^\times\},\qquad
\kappa=\gcd(n,q-1),\qquad R=\frac{Q-1}{(q-1)\kappa}.
\]
Then \(|\mathcal C|=(q-1)\kappa\), and multiplication by any \(c\in\mathcal C\) preserves normality.
A \(\mathcal C\)-coset contains a primitive element if and only if it generates the cyclic quotient group \(E^\times/\mathcal C\).
Each such generating coset contains \(\varphi(N)/\varphi(R)\) primitive elements.
\end{lemma}
\begin{proof}
If \(c^q=\lambda c\), with \(\lambda\in k^\times\), then
\[
(cx)^{q^i}=c\lambda^i x^{q^i}.
\]
The two families of conjugates are therefore linearly independent simultaneously.
Moreover, \(\mathcal C\) is the subgroup of solutions of \(c^{(q-1)^2}=1\) in \(E^\times\), so its order is
\[
\gcd(N,(q-1)^2)=(q-1)\gcd\!\left(\frac{q^n-1}{q-1},q-1\right)=(q-1)\kappa.
\]
Fix a generator \(z\) of \(E^\times\). The coset \(z^j\mathcal C\) generates the quotient if and only if \(\gcd(j,R)=1\).
The reduction map on units \((\mathbb Z/N\mathbb Z)^\times\to(\mathbb Z/R\mathbb Z)^\times\) is surjective: choose lifts separately at the prime factors of \(N\), then apply the Chinese remainder theorem.
Its fibers all have the same size, namely \(\varphi(N)/\varphi(R)\).
\end{proof}

\leavevmode\phantomsection\label{v93:def-e-free}\ignorespaces
An element \(x=z^j\ne0\) is called \textbf{\(e\)-free} if \(\gcd(j,e)=1\), where \(e\mid R\).
This condition is independent of the choice of the generator \(z\).
\leavevmode\phantomsection\label{v93:def-local-normality}\ignorespaces
For a nonconstant squarefree polynomial \(g\mid M\), say that \(x\) satisfies the \(g\)-condition if
\[
x\notin P(\sigma)E\qquad(P\mid g\text{ irreducible}).
\]
Let \(\mathscr N(e,g)\) count the elements that are \(e\)-free and satisfy the \(g\)-condition. The lemma gives
\begin{equation}\label{v93:eq-quotient-pn}
\boxed{\quad\operatorname{PN}_q(n)=\frac{\vartheta(N)}{\vartheta(R)}\mathscr N(R,M).\quad}
\end{equation}
Thus it suffices to find a \emph{normal element} whose coset generates \(E^\times/\mathcal C\), then scale within that coset.
Scaling uses field multiplication, while the Frobenius module structure continues to detect normality.

\subsection{Local trace indicators and the joint estimate}

Choose a nontrivial additive character \(\psi\) of \(k\), and put \(\Psi=\psi\circ\Tr_{E/k}\). Proposition~\ref{v128:prop-local-spectrum} gives the indicator expansion for the \(g\)-condition and the total absolute value \(W_q(g)-1\) of its nonzero-frequency coefficients. These additive conditions are now intersected with multiplicative generation conditions.

The cyclic-group argument in Proposition~\ref{prop:primitive-character-indicator} also applies to the \(e\)-free condition:
in \(E^\times=\langle z\rangle\), we only need to test whether the exponent \(j\) is coprime to \(e\). Thus
\[
\boldsymbol1_{e\text{-free}}(x)
=\vartheta(e)\sum_{d\mid e}\frac{\mu(d)}{\varphi(d)}
\sum_{\ord\chi=d}\chi(x)\qquad(x\ne0).
\]
Extend each multiplicative character by \(\chi(0)=0\).
Since \(g\ne1\), the \(g\)-condition automatically excludes zero, so no additional correction at zero is needed.
\begin{proposition}[The basic joint estimate]\label{v93:prop-basic}
For \(e\mid R\) and nonconstant \(g\mid M\),
\begin{equation}\label{v93:eq-basic}
\left|\mathscr N(e,g)-\vartheta(e)\Theta(g)Q\right|
\le\vartheta(e)\Theta(g)(W(e)-1)(W_q(g)-1)\sqrt Q.
\end{equation}
In particular, if
\begin{equation}\label{v93:eq-easy}
\boxed{\quad q^{n/2}>(2^{\omega(R)}-1)(2^s-1),\quad}
\end{equation}
then \(\operatorname{PN}_q(n)>0\).
\end{proposition}
\begin{proof}
Multiply the two indicators and sum over \(E^\times\).
The inner sum associated with the trivial multiplicative character is exactly the number of elements satisfying the \(g\)-condition.
In the local Frobenius decomposition, each \(P_i\mid g\) requires the coordinate in a quotient space of dimension \(d_i\) to be nonzero.
These conditions are independent, so their solution count is \(Q\Theta(g)\).
Since \(g\ne1\), none of these solutions is zero.
For a nontrivial multiplicative character, the sum at additive frequency zero vanishes.
Writing \(\Psi_b(x)=\Psi(bx)\), we obtain the exact identity
\begin{equation}\label{v93:eq-mixed-character-count}
\frac{\mathscr N(e,g)}{\vartheta(e)\Theta(g)}
=Q+\sum_{\substack{d\mid e\\d>1}}\frac{\mu(d)}{\varphi(d)}
\sum_{\ord\chi=d}\sum_{0\ne b\in V_g}
c_g(b)G_E(\chi,\Psi_b).
\end{equation}
For \(b\ne0\), nondegeneracy of the trace pairing ensures that \(\Psi_b\) is nontrivial,
so Theorem~\ref{thm:gauss} gives \(|G_E(\chi,\Psi_b)|=\sqrt Q\).
For each \(d\mid e\), there are \(\varphi(d)\) multiplicative characters of order \(d\). Thus
\[
\sum_{\substack{d\mid e\\d>1}}\frac{|\mu(d)|}{\varphi(d)}
\#\{\chi:\ord\chi=d\}=W(e)-1.
\]
Combine this with~\eqref{v93:eq-local-coefficient-sum} and apply the triangle inequality to obtain~\eqref{v93:eq-basic}.
Taking \(e=R,g=M\) and using~\eqref{v93:eq-quotient-pn} proves the existence assertion.
\end{proof}

\begin{proposition}[Primitive elements in a fixed weighted-trace fiber]\label{v93:prop-fixed-fiber}
For any weighted-trace fiber \(\mathcal A_{\boldsymbol b}\) from Theorem~\ref{v93:thm-patch}, let \(P(\mathcal A_{\boldsymbol b})\) denote its number of primitive elements. Then
\[
P(\mathcal A_{\boldsymbol b})\ge\vartheta(N)
\left[q^{n-s}-(W(N)-1)(1-q^{-s})\sqrt Q\right].
\]
In particular, if \(q^{n/2}>(q^s-1)(W(N)-1)\), every such fiber contains a primitive element.
\end{proposition}
\begin{proof}
Expand the fiber indicator using \(s\) additive character orthogonality relations.
For a nontrivial multiplicative character, the contribution at frequency zero vanishes.
The other \(q^s-1\) frequencies are all nonzero, and the corresponding Gauss sums have absolute value \(\sqrt Q\) by Theorem~\ref{thm:gauss}.
Each frequency coefficient has absolute value \(q^{-s}\).
Expanding the primitivity indicator, the normalized coefficients of the nontrivial multiplicative characters have total absolute value \(W(N)-1\), giving the stated lower bound.
\end{proof}

This criterion directly estimates the number of primitive elements in one fiber consisting entirely of normal elements, but is not a necessary step in proving general existence.
We now return to the joint count over all normal local components and sieve only the conditions that remain to be imposed.

\subsection{Sieving only the remaining local conditions}
The basic sufficient condition does not hold for every pair of parameters.
We use a \textbf{local density sieve}: within the set already satisfying some of the conditions, apply the union bound while retaining the density of each remaining condition.
Let \(e\mid R\) be a squarefree integer and let \(g\mid M\) be nonconstant.
Write the remaining prime factors as \(\ell_1,\ldots,\ell_u\) and the remaining irreducible polynomials as \(H_1,\ldots,H_v\), with \(f_j=\deg H_j\). Put
\[
\delta=1-\sum_{i=1}^u\frac1{\ell_i}-\sum_{j=1}^v q^{-f_j},\qquad
A_0=\sum_{i=1}^u(1-\ell_i^{-1}),\qquad B_0=\sum_{j=1}^v(1-q^{-f_j}).
\]
When \(\delta>0\), define
\begin{equation}\label{v93:eq-C}
\begin{split}
C(e,g)={}&(W(e)-1)(W_q(g)-1)\\
&+\frac{W(e)(W_q(g)-1)A_0+W_q(g)(W(e)-1)B_0}{\delta}.
\end{split}
\end{equation}

\begin{theorem}[The local density sieve]\label{v93:thm-sieve}
If \(\delta>0\), then
\begin{equation}\label{v93:eq-sieve}
\mathscr N(R,M)\ge\vartheta(e)\Theta(g)\delta
\bigl(Q-C(e,g)\sqrt Q\bigr).
\end{equation}
Consequently, \(q^{n/2}>C(e,g)\) suffices for the existence of a primitive normal element.
\end{theorem}
\begin{proof}
Within the set satisfying both the \(e\)-condition and the \(g\)-condition, exclude the remaining obstructions one at a time. This gives
\[
\mathscr N(R,M)\ge\sum_i\mathscr N(e\ell_i,g)+\sum_j\mathscr N(e,gH_j)
-(u+v-1)\mathscr N(e,g).
\]
If there are no remaining conditions, this is interpreted as an identity.
By the preceding expansions, adding one multiplicative prime factor introduces only new multiplicative characters, so
\[
\begin{split}
&\left|\mathscr N(e\ell,g)-(1-\ell^{-1})\mathscr N(e,g)\right|\\
&\hspace{1em}\le\vartheta(e)\Theta(g)(1-\ell^{-1})W(e)(W_q(g)-1)\sqrt Q.
\end{split}
\]
When a local polynomial factor \(H\) is added, all new additive frequencies are nonzero, and the difference contributed by the trivial multiplicative character is zero. Hence
\[
\begin{split}
&\left|\mathscr N(e,gH)-(1-q^{-\deg H})\mathscr N(e,g)\right|\\
&\hspace{1em}\le\vartheta(e)\Theta(g)(1-q^{-\deg H})W_q(g)(W(e)-1)\sqrt Q.
\end{split}
\]
Substitute these difference estimates into the union bound.
The coefficient of \(\mathscr N(e,g)\) is exactly \(\delta\).
Applying the lower bound from~\eqref{v93:eq-basic} to this term gives~\eqref{v93:eq-sieve}.
\end{proof}

\subsection{Two generation conditions among elements of full degree}

Generating the multiplicative quotient also forces an element to have full extension degree.
For some parameters, this allows a direct comparison of the numbers of the two types of generators within the set of elements of full degree.

\begin{lemma}[Generating the quotient forces full extension degree]\label{v93:lem-quotient-full-degree}
Let \(n\ge3\), and write
\[
\mathcal C=\{c\in E^\times:c^{q-1}\in k^\times\},\qquad
R=\frac{q^n-1}{(q-1)\gcd(n,q-1)}.
\]
If the image of \(x\in E^\times\) generates \(E^\times/\mathcal C\), then \(x\) has degree exactly \(n\) over \(k\).
\end{lemma}
\begin{proof}
Take any proper divisor \(d\) of \(n\), so \(d\le n/2\). Write \(\kappa=\gcd(n,q-1)\le n\). Then
\[
\frac{R}{q^d-1}
=\frac{1+q^d+q^{2d}+\cdots+q^{n-d}}{(q-1)\kappa}
>\frac{q^{n-d}}{n(q-1)}
\ge\frac{q^{\lceil n/2\rceil}}{n(q-1)}\ge1.
\]
The last inequality follows directly. Put \(h=\lceil n/2\rceil\ge2\). Since \((q-2)^2\ge0\) and \(q\ge2\),
\[
\frac{q^h}{q-1}=\frac{q^2}{q-1}q^{h-2}
\ge4\cdot2^{h-2}=2^h\ge n.
\]
For the final step, start with \(n=3,4\) and induct by increasing \(n\) by two. Thus
\[
 R>q^d-1\qquad(d\mid n,\ d<n).
\]
If \(x\in\F_{q^d}^\times\), its multiplicative order is at most \(q^d-1\).
Its image in \(E^\times/\mathcal C\) therefore also has order at most \(q^d-1<R\), and cannot generate that quotient.
\end{proof}

Let \(D_q(n)\) be the number of elements of degree exactly \(n\) over \(k\). Counting Frobenius orbits gives
\[
D_q(n)=\sum_{d\mid n}\mu(d)q^{n/d}.
\]
Let \(\mathcal N\) be the set of normal elements and \(\mathcal P_R\) the set of elements whose images generate \(E^\times/\mathcal C\).
For \(n\ge3\), the preceding lemma and normality give, respectively,
\[
\mathcal P_R\subseteq\{x:\deg_k(x)=n\},\qquad
\mathcal N\subseteq\{x:\deg_k(x)=n\}.
\]
This yields an existence criterion based on counting elements of full degree.

\begin{corollary}[Inclusion--exclusion for two generation conditions in full degree]\label{v93:cor-quotient-full-degree-union}
If \(n\ge3\) and
\begin{equation}\label{v93:eq-quotient-full-degree-union}
\boxed{\quad
(q^n-1)\frac{\varphi(R)}R
+q^n\Theta(M)>D_q(n),\quad}
\end{equation}
then \(E/k\) has a primitive normal element. Moreover,
\[
\boxed{\quad
\operatorname{PN}_q(n)\ge
\frac{\vartheta(q^n-1)}{\vartheta(R)}
\left[
(q^n-1)\vartheta(R)+q^n\Theta(M)-D_q(n)
\right]>0.
\quad}
\]
\end{corollary}
\begin{proof}
Apply inclusion--exclusion within the set of elements of full degree.
Both \(\mathcal P_R\) and \(\mathcal N\) lie in this set, so
\[
|\mathcal P_R\cap\mathcal N|
\ge |\mathcal P_R|+|\mathcal N|-D_q(n).
\]
Here \(|\mathcal P_R|=(q^n-1)\vartheta(R)\) and \(|\mathcal N|=q^n\Theta(M)\).
In each \(\mathcal C\)-coset that generates the quotient, the proportion of primitive elements is \(\vartheta(q^n-1)/\vartheta(R)\), and normality is constant throughout the coset.
Multiplying the intersection count by this proportion therefore gives the number of primitive normal elements.
\end{proof}

\begin{example}[Direct intersection in a quartic extension]
Take \((q,n)=(5,4)\). Here \(R=39\), and \(T^4-1\) has four linear factors over \(\F_5\), so
\[
 |\mathcal P_R|=624\cdot\frac{24}{39}=384,\quad
 |\mathcal N|=5^4(1-5^{-1})^4=256,\quad
 D_5(4)=625-25=600.
\]
Hence
\[
 |\mathcal P_R\cap\mathcal N|\ge40,
 \qquad \operatorname{PN}_5(4)\ge 40\cdot\frac{\vartheta(624)}{\vartheta(39)}=20.
\]
This lower bound comes directly from counting the two types of generators.
\end{example}

The inclusion--exclusion condition is only sufficient.
For example, when \((q,n)=(3,4)\), we have \(|\mathcal P_R|=|\mathcal N|=32\), whereas \(D_3(4)=72\).
Cardinalities alone therefore do not establish that the intersection is nonempty.
We next use the relation between normality and particular multiplicative obstructions.

\subsection{Multiplicative obstructions automatically excluded by normality}
\begin{lemma}[Short Frobenius relations exclude multiplicative subgroups]\label{v93:lem-automatic}
Let \(x\in E\) be normal, and let \(1\le d<n\). Then
\[
 x^{q^d-1}\notin k^\times.
\]
In particular, if a prime \(\ell\mid R\) satisfies
\begin{equation}\label{v93:eq-auto-divisibility}
 \frac{Q-1}{\ell}\ \bigm|\ (q^d-1)(q-1)
\end{equation}
for some \(1\le d<n\), every normal element is automatically \(\ell\)-free.
\end{lemma}
\begin{proof}
If \(x^{q^d-1}=c\in k^\times\), then \((\sigma^d-cI)x=0\).
The nonzero operator polynomial \(T^d-c\) has degree less than \(n\), so it cannot annihilate a normal element.
If \(x\) is not \(\ell\)-free, then \(x^{(Q-1)/\ell}=1\).
Equation~\eqref{v93:eq-auto-divisibility} gives \((x^{q^d-1})^{q-1}=1\), hence \(x^{q^d-1}\in k^\times\), a contradiction.
\end{proof}
Prime factors whose obstructions are already excluded by normality need not be included in the multiplicative sieve.

\paragraph{Direct applications.}
For \((q,n)=(3,4)\), we have \(R=20\). Take \(\ell=5,d=2\); then
\[
 (Q-1)/5=16\mid(3^2-1)(3-1)=16.
\]
Only the prime factor \(2\) remains to be sieved, and its normalized error coefficient is \(7<3^{4/2}=9\).
For \((q,n)=(5,8)\), we have \(R=24414\). Take \(\ell=313,d=4\); then
\[
 (Q-1)/313=1248\mid(5^4-1)(5-1)=2496.
\]
The remaining multiplicative prime factors are \(2,3,13\), giving error coefficient \(441<5^{8/2}=625\).
The pair \((5,4)\) can also be handled this way, with \(\ell=13,d=2\), but the preceding inclusion--exclusion argument in full degree already gives a more direct proof.
In the order \((3,4),(5,4),(5,8)\), the values of \((Q-1)/\ell\) are \(16,48,1248\), which divide \(16,96,2496\), respectively.
Applying the basic joint estimate to \((e,M)\) first gives an \(e\)-free normal element.
Its automatic \(\ell\)-freeness then makes it \(R\)-free, and the scaling lemma gives a primitive normal element.
Primality of \(313\) requires trial division only by primes at most 17.

\subsection{Eight boundary cases and general existence}

Proposition~\ref{v93:prop-arithmetic-complete} in the appendix identifies the exceptions to the basic sufficient condition as exactly
\[
(2,3),(2,6),(2,15),(3,4),(3,8),(5,4),(5,8),(7,6).
\]
The pairs \((3,4)\) and \((5,8)\) have already been handled by the preceding automatic exclusion of obstructions.

Corollary~\ref{v93:cor-quotient-full-degree-union} also handles \((q,n)=(2,3),(2,6),(2,15),(5,4)\).
The full calculation for \((5,4)\) was given above.

For the other three binary cases, \(n=3,6,15\), respectively, we have
\[
\begin{aligned}
\bigl(|\mathcal P_R|,|\mathcal N|,D_2(n)\bigr)
&=(6,3,6),\quad(36,24,54),\quad(27000,10125,32730).
\end{aligned}
\]
The corresponding intersections therefore have at least \(3,6,4395\) elements, respectively, and are all nonempty.
Here \(2^{15}-1=7\cdot31\cdot151\), and the irreducible factors of \(M\) have degrees \(1,2,4,4,4\).
The other numbers follow directly from the standard counts for cyclic groups and Frobenius orbits.

\paragraph{Two local density sieves.}
Only \((3,8)\) and \((7,6)\) remain. Apply Theorem~\ref{v93:thm-sieve} with
\[
\begin{aligned}
(q,n)=(3,8):&\quad e=2,\quad g=(T-1)(T+1),\\
(q,n)=(7,6):&\quad e=2,\quad g=T-1.
\end{aligned}
\]
We compute the densities and error coefficients explicitly; both satisfy the required strict inequalities for a positive lower bound.
In the first case, the remaining multiplicative prime factors are \(5,41\), and the remaining three polynomial factors all have degree two. Hence
\[
 \delta=1-\frac15-\frac1{41}-\frac39=\frac{272}{615},\quad
 A_0=\frac45+\frac{40}{41},\quad B_0=\frac83.
\]
Substitution into~\eqref{v93:eq-C} gives \(C=3+(6A_0+4B_0)/\delta=1741/34\).
In the second case, the remaining multiplicative prime factors are \(19,43\), and the remaining five polynomial factors are all linear. Hence
\[
 \delta=1-\frac1{19}-\frac1{43}-\frac57=\frac{1200}{5719},\quad
 A_0=\frac{18}{19}+\frac{42}{43},\quad B_0=\frac{30}{7},
\]
so \(C=1+2(A_0+B_0)/\delta=6019/100\).
The final comparisons are \(1741<34\cdot81\) and \(6019<100\cdot343\), respectively, and both are strict.

\begin{theorem}[The primitive normal basis theorem]\label{v93:thm-PNBT}
For every prime power \(q\) and positive integer \(n\), the extension \(\F_{q^n}/\F_q\) has a primitive normal element.
\end{theorem}
\begin{proof}
For \(n=1\), a generator of \(k^\times\) is nonzero and hence is a basis of the one-dimensional \(k\)-space \(k\). For \(n=2\), Proposition~\ref{v88:prop-primitive-quadratic-normal} shows that every primitive element is normal.

For \(n\ge3\), Proposition~\ref{v93:prop-arithmetic-complete} in the appendix and the basic joint estimate cover every pair except the eight boundary cases.
Automatic exclusion handles \((3,4),(5,8)\); inclusion--exclusion in full degree handles \((2,3),(2,6),(2,15),(5,4)\); and the local density sieve handles \((3,8),(7,6)\).
The scaling lemma turns the resulting generators of the quotient into primitive elements while preserving normality.
This exhausts all cases.
\end{proof}

The proof places the two generation conditions in a single counting formula:
the local Frobenius decomposition provides additive frequencies, the multiplicative quotient provides multiplicative characters, and Gauss sums control their interaction.
Appendix~\ref{v93:chap-primitive-normal-arithmetic} gives the complete arithmetic coverage from the uniform estimates to the eight boundary cases and is part of the proof for all parameters.
For classical proofs of the primitive normal basis theorem and sieve methods, see~\cite{LS,CohenHuczynska2003}.

\section{The Quadratic Character and Exact Counts}
Let \(Q\) be odd. Denote the unique nontrivial multiplicative character of order two by \(\chi\); it is one on nonzero squares and minus one on nonsquares, with \(\chi(0)=0\).
Thus \(y^2=b\) has exactly \(1+\chi(b)\) solutions.

\begin{proposition}\leavevmode\phantomsection\label{prop:quadratic-sum}
If \(a\ne0\), then
\[
\sum_{x\in K}\chi(x^2-a)=-1.
\]
Consequently, the number of solutions of \(x^2+y^2=a\) is
\[
\begin{cases}
Q-\chi(-1),&a\ne0,\\
Q+(Q-1)\chi(-1),&a=0.
\end{cases}
\]
\end{proposition}
\begin{proof}
The equation \(x^2-y^2=a\) can be written as \((x-y)(x+y)=a\).
Since two is invertible, any nonzero \(u=x-y\) determines \(v=x+y=a/u\) uniquely, and then uniquely recovers \(x,y\).
Hence there are \(Q-1\) solutions. On the other hand, their number is \(\sum_x(1+\chi(x^2-a))\), proving the first formula.
When \(a\ne0\), the number of solutions of \(x^2+y^2=a\) is \(Q+\sum_x\chi(a-x^2)=Q-\chi(-1)\).
When \(a=0\), the zero solution contributes one, and each nonzero \(x\) contributes \(1+\chi(-1)\), giving the second assertion.
\end{proof}

\begin{example}[Different character sums on affine cosets of equal dimension]\label{v106:ex-coset-character}
Let \(q\) be odd and \(E=\F_{q^2}\), and let \(\chi_E,\chi_q\) be quadratic characters, both extended by zero at zero. Then
\[
\boxed{\sum_{x\in\F_q}\chi_E(x+a)=
\begin{cases}q-1,&a\in\F_q,\\-1,&a\notin\F_q.\end{cases}}
\]
Indeed, \(\chi_E(z)=\chi_q(\Nm_{E/\F_q}z)\). If \(a\in\F_q\), the coset contains exactly one zero, and all other base-field elements are squares in \(E\). If \(a\notin\F_q\), then
\[
\Nm(x+a)=x^2+\Tr(a)x+\Nm(a)
\]
has nonzero discriminant \((a-a^q)^2\). Completing the square and applying Proposition~\ref{prop:quadratic-sum} gives sum \(-1\). Equal additive dimension does not force equal character sums: position relative to the subfield affects multiplicative cancellation. The projection fibers in Example~\ref{v106:ex-cover-5-11} show how the same position issue affects covering.
\end{example}

\subsection{Restricted Coefficients and Translation Orbits}
\label{v105:subsec-restricted-coefficients}

To construct irreducible polynomials with restricted coefficients,
one can first choose a family that preserves irreducibility and then impose the coefficient conditions within it.
Here we consider coefficients in prescribed multiplicative cosets.
Let \(k=\F_q\), and let \(f\in k[X]\) be monic and irreducible of degree \(n\ge2\).
For every \(t\in k\), the polynomial \(f(X+t)\) remains irreducible:
substitution of \(X+t\) is an automorphism of \(k[X]\), with inverse given by substitution of \(X-t\).
If the roots of \(f\) are \(\alpha_1,\ldots,\alpha_n\), those of \(f(X+t)\) are
\(\alpha_1-t,\ldots,\alpha_n-t\).
All root differences are unchanged, and hence
\[
\operatorname{disc}(f(X+t))=\operatorname{disc}(f).
\]
Both properties hold in every characteristic.

When \(p=\Char k>n\), each \(i!\) is invertible in \(k\), and Taylor expansion gives
\begin{equation}\label{v105:eq-translate-coefficients}
f(X+t)=\sum_{i=0}^n c_i(t)X^i,\qquad
c_i(t)=\frac{f^{(i)}(t)}{i!}.
\end{equation}
If \(f=X^n+a_{n-1}X^{n-1}+\cdots+a_0\), then \(c_{n-1}(t)=a_{n-1}+nt\),
so these \(q\) translated polynomials are distinct.
A translation orbit therefore fixes irreducibility and discriminant,
while conditions on the coefficients become multiplicative conditions on a single parameter \(t\).
The hypothesis \(p>n\) is used for the expansion with ordinary derivatives;
it is not needed to preserve the two invariants.

More explicitly, choose \(H_i\le k^\times\) and \(h_i\in k^\times\), and put \(e_i=[k^\times:H_i]\).
Extend every multiplicative character, including the trivial character, by zero at zero.
Character orthogonality then gives
\begin{equation}\label{v105:eq-translate-coset-count}
\begin{aligned}
&\#\{t\in k:c_i(t)\in h_iH_i\ \ (0\le i<n)\}\\
&\qquad=
\frac1{\prod_{i=0}^{n-1}e_i}
\sum_{\substack{\rho_i^{\,e_i}=1\\0\le i<n}}
\ \sum_{t\in k}\prod_{i=0}^{n-1}\rho_i\bigl(h_i^{-1}c_i(t)\bigr).
\end{aligned}
\end{equation}
Indeed, for nonzero \(z\), the average
\(e_i^{-1}\sum_{\rho_i^{e_i}=1}\rho_i(h_i^{-1}z)\)
is the indicator of \(h_iH_i\); at \(z=0\), both sides are zero as well.
Multiplying these indicators and summing over \(t\) proves the formula.
Thus counting restricted coefficients reduces to joint character sums of derivative values.
For subsequent estimates, common or repeated roots of the derivatives can make the multiplicative conditions dependent,
and must be separated from the general case.
For this translation method and its counting applications to thin multiplicative subgroups,
see~\cite{OstafeShparlinski2026Thin}.
The following quadratic case can be evaluated exactly using the preceding proposition.

\begin{proposition}[Irreducible quadratics with two nonleading square coefficients]\label{v105:prop-square-coefficient-quadratics}
Let \(q\) be odd, and let \(\chi\) be the quadratic character of \(\F_q\), with \(\chi(0)=0\).
Write \(S=(\F_q^\times)^2\) for the subgroup of nonzero squares, and set
\[
N_q^\square=\#\{(b,c)\in S^2:X^2+bX+c\text{ irreducible}\}.
\]
Then
\begin{equation}\label{v105:eq-square-coefficient-count}
\boxed{\quad
N_q^\square=\frac{(q-1)(q+\chi(-1)-2)}8.
\quad}
\end{equation}
\end{proposition}
\begin{proof}
Fix \(b\in S\) and put \(t=c/b^2\).
Multiplication by the nonzero square \(b^{-2}\) gives a bijection between \(c\in S\) and \(t\in S\).
The quadratic discriminant is
\[
b^2-4c=b^2(1-4t).
\]
In odd characteristic, completing the square shows that the quadratic is irreducible exactly when its discriminant is a nonsquare.
Thus the required parameters satisfy \(\chi(t)=1\) and \(\chi(1-4t)=-1\).
Let \(A\) denote their number.
Except at \(t=0,1/4\), the expression
\((1+\chi(t))(1-\chi(1-4t))\)
is four when these conditions hold and zero otherwise.
It is zero at \(t=0\) and two at \(t=1/4\), so
\[
4A=\sum_{t\in\F_q}
(1+\chi(t))(1-\chi(1-4t))-2.
\]
Both character sums with linear arguments vanish.
Moreover, \(u=8t-1\) runs through \(\F_q\), and
\[
t(1-4t)=\frac{1-u^2}{16}.
\]
Proposition~\ref{prop:quadratic-sum} therefore gives
\[
\sum_{t\in\F_q}\chi\bigl(t(1-4t)\bigr)
=\sum_{u\in\F_q}\chi(1-u^2)
=\chi(-1)\sum_{u\in\F_q}\chi(u^2-1)
=-\chi(-1).
\]
Hence \(4A=q+\chi(-1)-2\).
This number does not depend on the chosen \(b\).
Multiplying by \(|S|=(q-1)/2\) proves the formula.
\end{proof}

\begin{example}[Small fields and the correction for zero discriminant]\label{v105:ex-square-coefficient-boundary}
For \(q=3\), formula~\eqref{v105:eq-square-coefficient-count} gives zero.
Here \(S=\{1\}\), and the only candidate is
\[
X^2+X+1=(X-1)^2.
\]
For \(q=5\), we have \(S=\{1,4\}\), and the formula gives exactly two polynomials:
\[
X^2+X+1,\qquad X^2+4X+1.
\]
Both have discriminant two, a nonsquare.
The other two candidates have constant term four and discriminant zero.
The subtraction of two in the proof removes exactly the contribution from \(t=1/4\);
this zero-discriminant case must not be counted as irreducible.
\end{example}

Among all \(|S|^2\) candidates, the proportion that are irreducible is
\[
\frac12+\frac{\chi(-1)-1}{2(q-1)}
=
\begin{cases}
\frac12,&q\equiv1\pmod4,\\[2pt]
\frac12-\frac1{q-1},&q\equiv3\pmod4.
\end{cases}
\]
This exact formula gives both the main proportion under the coefficient restrictions and the correction for small fields.

\section{The Chevalley--Warning Congruence for Root Counts}
\label{sec:chevalley-warning}
Classifying quadratic forms requires criteria for the existence of a nontrivial zero. We first use the polynomial indicators and power sums of Chapter~\ref{chap:functions}
to prove a root-count congruence for low-degree equations in several variables.
\begin{theorem}\leavevmode\phantomsection\label{thm:chevalley-warning}
Let \(F_1,\ldots,F_s\in\F_q[X_1,\ldots,X_m]\) be nonzero polynomials with \(\sum_i\deg F_i<m\). Their number of common zeros in \(\F_q^m\) is divisible by \(p=\Char\F_q\).
\end{theorem}
\begin{proof}
The function
\[
I(x)=\prod_{i=1}^s\bigl(1-F_i(x)^{q-1}\bigr)
\]
is one at a common zero and zero elsewhere, so \(\sum_{x\in\F_q^m}I(x)\) is the image in the field of the number of common zeros.
Expand \(I\). Each monomial has total degree less than \(m(q-1)\), so some variable has exponent less than \(q-1\).
Summing over that variable gives \(q=0\) when its exponent is zero, and also zero by the power-sum lemma when its exponent lies between one and \(q-2\).
Thus every monomial sums to zero over \(\F_q^m\). The number of common zeros is zero in the field, equivalently, divisible by \(p\).
\end{proof}

\begin{corollary}
If all \(F_i\) are homogeneous of positive degree and \(\sum_i\deg F_i<m\), they have a nonzero common zero.
In particular, a homogeneous quadratic form in three or more variables over a finite field always has a nontrivial zero.
\end{corollary}
\begin{proof}
The zero vector is already a solution. Since the solution count is divisible by \(p\), it cannot be the only solution.
\end{proof}

This proof uses \(1-a^{q-1}\) to express vanishing as a polynomial function and then uses power sums to eliminate all terms of low total degree.
Unlike the univariate root bound, it controls a counting congruence for common roots in several variables.

\begin{proposition}[Summation and the top coefficient]\label{v91:prop-sum-top-coefficient}
Let \(m\ge1\). For \(P\in\F_q[X_1,\ldots,X_m]\), let \(\overline P\) be its unique representative with degree less than \(q\) in each variable in the function algebra
\[
A_m=\F_q[X_1,\ldots,X_m]/(X_1^q-X_1,\ldots,X_m^q-X_m).
\]
Then the following identity holds in \(\F_q\):
\begin{equation}\label{v91:eq-sum-top-coefficient}
\sum_{a\in\F_q^m}P(a)
=(-1)^m[X_1^{q-1}\cdots X_m^{q-1}]\overline P,
\end{equation}
where brackets denote extraction of the indicated coefficient.
\end{proposition}
\begin{proof}
Taking products of the point indicators from Chapter~\ref{chap:functions} gives
\[
\overline P(X)=\sum_{a\in\F_q^m}P(a)
\prod_{j=1}^m\bigl(1-(X_j-a_j)^{q-1}\bigr).
\]
Both sides have the same value at every point and degree less than \(q\) in each variable, so repeated application of the univariate root bound makes them equal.
The coefficient of \(X_1^{q-1}\cdots X_m^{q-1}\) in each point indicator is \((-1)^m\). Comparing coefficients proves the identity.
\end{proof}

Thus summation over the full evaluation set reads a single top coefficient of the standard representative.
In the Chevalley--Warning theorem, \(I=\prod_i(1-F_i^{q-1})\) has total degree less than \(m(q-1)\).
Taking remainders modulo \(X_j^q-X_j\), one variable at a time, does not increase total degree, so \(\overline I\) has no such top term.
Equation~\eqref{v91:eq-sum-top-coefficient} again shows that the common-zero count has zero image in \(\F_q\).
This explains exactly how the degree condition forces the root-count congruence modulo \(p\). The sum here is taken in \(\F_q\), whereas complex character orthogonality retains the integer solution count.

\section{Polarization, Radicals, and Classification of Quadratic Forms}
\label{sec:quadratic-forms}

In this section take \(k=\F_q\), and let \(V\) be a \(k\)-space of dimension \(m\ge1\).
Solution counts for quadratic equations are controlled by the polar form and the remaining quadratic information; the trace pairing further turns the former into the kernel of a linearized operator.
For the classification of quadratic forms and solution-count formulas, see \cite[Chapter 6, Section 2]{LN}.

\begin{definition}\label{def:quadratic-form}
A map \(Q:V\to k\) is called a \textbf{quadratic form} if \(Q(cx)=c^2Q(x)\) and
\[
B_Q(x,y)=Q(x+y)-Q(x)-Q(y)
\]
is bilinear. The form \(B_Q\) is called the \textbf{polar form} of \(Q\). Write
\[
R=\rad B_Q=\{x\in V:B_Q(x,y)=0\text{ for all }y\in V\}.
\]
A nonzero vector satisfying \(Q(x)=0\) is called \textbf{isotropic};
a quadratic form with no such vector is called \textbf{anisotropic}.
\leavevmode\phantomsection\label{def:quadratic-isometry}\ignorespaces Quadratic forms \(Q:V\to k\) and \(Q':V'\to k\) are called \textbf{isometric} if there is a \(k\)-linear isomorphism \(T:V\to V'\) such that \(Q'(Tx)=Q(x)\).
\end{definition}

For any basis \((e_i)\), the quadratic-form axioms give
\[
Q\Bigl(\sum_i x_i e_i\Bigr)
=\sum_i Q(e_i)x_i^2+\sum_{i<j}B_Q(e_i,e_j)x_ix_j.
\]
Thus the definition is equivalent to a homogeneous quadratic polynomial in coordinates. In odd characteristic, \(B_Q(x,x)=2Q(x)\),
so the polar form determines the quadratic form completely; in characteristic two, \(B_Q\) is always alternating, and the square terms cannot be recovered from it.

\begin{example}[A trace quadratic form with known polarization and undetermined type]\leavevmode\phantomsection\label{exm:f64-polar}
Over \(E=\F_{64}\), take \(Q(x)=\Tr_{E/\F_2}(x^3)\), retaining the notation \(\sigma(x)=x^2\).
For \(x,y\in E\),
\[
B_Q(x,y)=\Tr(x^2y+xy^2)
=\Tr\bigl(x(\sigma^{-1}+\sigma)y\bigr).
\]
The second equality follows from invariance of the trace under \(\sigma\): \(\Tr(x^2y)=\Tr(xy^{32})\).
The trace pairing is nondegenerate and \(\sigma\) is invertible, so
\[
\rad B_Q=\ker(\sigma+\sigma^{-1})
=\ker(\sigma^2-1)=\F_4.
\]
Thus \(B_Q\) has rank \(6-2=4\). For \(x\in\F_4\), we have \(x^3\in\F_2\), hence \(Q(x)=6x^3=0\). For every \(r\in\F_4\), \(Q(x+r)=Q(x)+Q(r)+B_Q(x,r)=Q(x)\),
so \(Q\) descends to the four-dimensional quotient \(E/\F_4\) with nondegenerate polar form.
Polar rank still does not determine the solution count: the classification below leaves two quadratic types,
and Example~\ref{exm:f64-quadratic} uses Gauss sums to determine the type and value distribution in this example.
\end{example}

\phantomsection\label{def:orthogonal-quadratic-language}
If \(V=U\oplus W\) and \(B_Q(U,W)=0\), then \(Q(u+w)=Q(u)+Q(w)\). This is an \textbf{orthogonal sum} of quadratic forms, denoted by \(Q|_U\perp Q|_W\).
The binary quadratic form \(H(x,y)=xy\) is called a \textbf{hyperbolic plane}; \(H^{\perp s}\) denotes the orthogonal sum of \(s\) hyperbolic planes, with \(H^{\perp0}\) the zero form on the zero-dimensional space.
In characteristic two also write \(N_\delta(x,y)=x^2+xy+\delta y^2\).
For a nondegenerate alternating bilinear form, a basis \((e_1,f_1,\ldots,e_s,f_s)\) satisfying \(B(e_i,e_j)=B(f_i,f_j)=0\) and \(B(e_i,f_j)=\delta_{ij}\) is called a \textbf{symplectic basis}.

\begin{theorem}[Canonical forms of quadratic forms over finite fields]\label{thm:quadratic-classification}
Let \(R=\rad B_Q\).
\begin{enumerate}
\item If \(q\) is odd, then \(Q|_R=0\). In suitable coordinates,
\[
Q=a_1x_1^2+\cdots+a_rx_r^2,\qquad a_i\ne0,\quad r=\rank B_Q.
\]
When \(r>0\), the isometry type is determined by \(m,r\) and the square class of \(\Delta=a_1\cdots a_r\). Here \(\sum_{i=1}^r a_ix_i^2\) is called the nondegenerate part of \(Q\), and the square class of \(\Delta\) is its discriminant square class.
\item If \(q\) is even, write \(\rank B_Q=2s\). If \(Q|_R\ne0\), the canonical form is
\[
x_1y_1+\cdots+x_sy_s+z^2,
\]
with the remaining \(m-2s-1\) coordinates absent.
If \(Q|_R=0\), the canonical form is one of the following two types:
\[
\begin{aligned}
H^{\perp s}&=x_1y_1+\cdots+x_sy_s,\\
H^{\perp(s-1)}\perp N_\delta
&=x_1y_1+\cdots+x_{s-1}y_{s-1}+x_s^2+x_sy_s+\delta y_s^2,
\end{aligned}
\]
where \(\Tr_{k/\F_2}(\delta)=1\); the second type requires \(s\ge1\), and the remaining \(m-2s\) coordinates do not occur.
\end{enumerate}
\end{theorem}
\begin{proof}
In odd characteristic, \(2Q(x)=B_Q(x,x)\) gives \(Q|_R=0\). Choose any linear complement of \(R\);
the restriction of \(B_Q\) to it is nondegenerate, and orthogonal decomposition gives a diagonal form.
The binary change of variables in the proof of Lemma~\ref{lem:finite-symmetric-form} applies equally to these quadratic coefficients:
every \(ax^2+by^2\) represents one, so successive changes reduce the form to \(x_1^2+\cdots+x_{r-1}^2+\Delta' x_r^2\), where \(\Delta'\) has the same square class as \(\Delta\).
An invertible change of basis multiplies the determinant by a square, so this square class must also be an isometry invariant.

In even characteristic, a nondegenerate alternating space decomposes by successively choosing two-dimensional blocks with \(B_Q(e_i,f_i)=1\).
Thus \(V=U\oplus R\), where \(U\) has a symplectic basis \(e_1,f_1,\ldots,e_s,f_s\) and distinct two-dimensional blocks are orthogonal.
On \(R\), the function \(Q\) is additive, and \(\ell(x)=Q(x)^{q/2}\) is a \(k\)-linear functional because \(\ell(cx)=c\ell(x)\).
If \(Q|_R\ne0\), choose \(z\in R\) with \(Q(z)=1\), and form a basis of \(R\) from \(z\) and a basis of \(\ker\ell\).
Replace each \(e_i\) by \(e_i+Q(e_i)^{q/2}z\), and make the same adjustment to each \(f_i\);
their quadratic values then vanish without changing the polar pairings. This gives the first form in the second assertion.

Now assume that \(Q|_R=0\). It suffices to classify the quadratic form on \(U\).
If \(u\ne0\) and \(Q(u)=0\), nondegeneracy of the polar form gives \(v\) such that \(B_Q(u,v)=1\);
replacing \(v\) by \(v+Q(v)u\) makes \(Q(v)=0\). This splits off a hyperbolic block \(H(x,y)=xy\).
Take its orthogonal complement and repeat. If the remaining dimension is at least three, the Chevalley--Warning theorem~\ref{thm:chevalley-warning}
ensures another nontrivial zero; the remaining dimension is even, so the process ends with either a zero-dimensional space or an anisotropic two-dimensional space.

Write a two-dimensional block as \(ax^2+xy+by^2\). If \(a=0\), it has an isotropic vector and is therefore \(H\).
If \(a\ne0\), the substitution \(u=ax\) transforms it into \(a^{-1}(u^2+uy+ab\,y^2)\). This form is anisotropic if and only if \(T^2+T+ab\) is irreducible over \(k\), that is, \(\Tr_{k/\F_2}(ab)=1\).
If \(\theta^2+\theta+ab=0\), then
\[
u^2+uy+ab\,y^2=\Nm_{k(\theta)/k}(u+y\theta).
\]
Surjectivity of the norm allows multiplication by an element of norm \(a^{-1}\) to absorb the coefficient \(a^{-1}\).
The difference of any two parameters of absolute trace one has the form \(c^2+c\); substitution of \(u\mapsto u+cy\) transforms their norm forms into one another.
This gives the listed canonical forms. The solution counts in Theorem~\ref{thm:quadratic-even-count} will also distinguish the two types directly.
\end{proof}

\section{Value Distributions of Quadratic Forms: Odd Characteristic}
\label{sec:quadratic-odd-count}

\leavevmode\phantomsection\label{def:quadratic-value-count}\ignorespaces Write \(N_Q(b)=|\{x\in V:Q(x)=b\}|\). If \(Q=0\), then \(N_Q(0)=q^m\), and every other value has no preimage.
Below assume polar rank \(r>0\), and let \(\eta\) be the quadratic multiplicative character of \(k\), with \(\eta(0)=0\).

\begin{theorem}\leavevmode\phantomsection\label{thm:quadratic-odd-count}
Suppose that \(q\) is odd and that the nondegenerate part of \(Q\) is \(\sum_{i=1}^r a_i x_i^2\), with \(\Delta=\prod_i a_i\).
If \(r\) is even, put \(\epsilon=\eta((-1)^{r/2}\Delta)\). Then
\[
N_Q(b)=q^{m-1}+\epsilon\,v(b)q^{m-r/2-1},
\qquad
v(b)=\begin{cases}q-1,&b=0,\\-1,&b\ne0.\end{cases}
\]
If \(r\) is odd, then
\[
N_Q(b)=q^{m-1}+
\eta\bigl((-1)^{(r-1)/2}\Delta\bigr)\eta(b)q^{m-(r+1)/2}.
\]
\end{theorem}
\begin{proof}
Fix a nontrivial additive character \(\psi\) and write \(G=G(\eta,\psi)\). Counting square roots and using character orthogonality gives, for \(c\ne0\),
\[
\sum_{x\in k}\psi(cx^2)
=\sum_{y\in k}(1+\eta(y))\psi(cy)=\eta(c)G.
\]
By Theorem~\ref{thm:gauss}, we also have \(G^2=\eta(-1)q\).
Separating the coordinates and allowing the \(m-r\) radical coordinates to vary freely gives, for \(t\ne0\),
\[
S_Q(t):=\sum_{x\in V}\psi(tQ(x))
=q^{m-r}\eta(t)^r\eta(\Delta)G^r.
\]
Additive character orthogonality yields
\[
N_Q(b)=q^{m-1}+
q^{m-r-1}\eta(\Delta)G^r\sum_{t\ne0}\eta(t)^r\psi(-tb).
\]
For even \(r\), the last sum is \(v(b)\); substituting \(G^r=(\eta(-1)q)^{r/2}\) gives the first formula.
For odd \(r\), the last sum vanishes if \(b=0\); if \(b\ne0\), a substitution gives its value as \(\eta(-b)G\).
Using \(G^{r+1}=(\eta(-1)q)^{(r+1)/2}\) again and simplifying the signs gives the second formula.
\end{proof}

For example, over \(\F_3\) the ternary form \(x^2+y^2+z^2\) has \(\Delta=1\) and \(r=3\),
so the numbers of vectors with values zero, one, and two are \(9,6,12\), respectively. Quadratic forms of the same rank may still have different value distributions, determined by the discriminant square class.

\section{Value Distributions of Quadratic Forms: Characteristic Two}
\label{sec:quadratic-even-count}

Now let \(q=2^e\), \(\rank B_Q=2s\), and \(R=\rad B_Q\).
Choose a complement \(U\) to \(R\) and a symplectic basis in \(U\), and write the restriction as
\[
Q|_U=\sum_{i=1}^s(a_i x_i^2+x_i y_i+b_i y_i^2).
\]
\leavevmode\phantomsection\label{def:quadratic-type-sign}\ignorespaces When \(Q|_R=0\), define
\[
\epsilon(Q)=(-1)^{\Tr_{k/\F_2}(\sum_i a_i b_i)}.
\]
The proof below also shows that this is independent of the complement \(U\) and symplectic basis.

\begin{theorem}\label{thm:quadratic-even-count}
If \(Q|_R\ne0\), then every \(b\in k\) satisfies
\[
N_Q(b)=q^{m-1}.
\]
If \(Q|_R=0\), then
\[
N_Q(b)=q^{m-1}+\epsilon(Q)v(b)q^{m-s-1},
\]
where \(v(0)=q-1\) and \(v(b)=-1\) for \(b\ne0\).
The sign \(\epsilon=1\) corresponds to the canonical form \(H^{\perp s}\), and \(\epsilon=-1\) to \(H^{\perp(s-1)}\perp N_\delta\).
\end{theorem}
\begin{proof}
If \(Q|_R\ne0\), then \(Q(r)=\ell(r)^2\), where \(\ell:R\to k\) is a nonzero linear functional.
Thus \(Q|_R\) is surjective, and each fiber has \(q^{\dim R-1}\) elements.
On each coset \(x+R\), the identity \(Q(x+r)=Q(x)+Q(r)\) shows that all values occur equally often.
Multiplication by the number of cosets \(q^{m-\dim R}\) gives \(q^{m-1}\).

Suppose that \(Q|_R=0\), and take the canonical additive character \(\psi(z)=(-1)^{\Tr_{k/\F_2}(z)}\).
The trace is invariant under squaring, so \(\psi(cy^2)=\psi(c^{q/2}y)\). For \(t\ne0\), the sum over one two-dimensional block is
\[
\sum_x\psi(ta x^2)\sum_y\psi\bigl(((tb)^{q/2}+tx)y\bigr).
\]
The inner sum is nonzero only when \(x=t^{-1}(tb)^{q/2}\), when it equals \(q\); substituting this \(x\) gives \(ta x^2=ab\).
Thus this two-dimensional sum equals \(q\psi(ab)\), independently of the choice of \(t\ne0\).
Multiplying over the blocks and including the free coordinates of \(R\) gives
\[
S_Q(t)=\epsilon(Q)q^{m-s}\qquad(t\ne0).
\]
The left-hand side depends only on the quadratic form, and \(s\) is determined by the polar rank. Hence \(\epsilon(Q)\) is independent of both complement and symplectic basis.
The solution-count formula follows from \(N_Q(b)=q^{-1}\sum_t\psi(-tb)S_Q(t)\).
For \(H\), we have \(a=b=0\) and sign one; for \(N_\delta\), we have \(a=1,b=\delta\) and sign minus one.
Thus the sign distinguishes exactly the two canonical types.
\end{proof}

\begin{definition}\label{def:arf-invariant}
When \(Q|_R=0\), the \textbf{Arf invariant} of the quadratic form on the quotient \(V/R\) is
\[
\operatorname{Arf}(Q)=\sum_i a_i b_i\pmod{\{c^2+c:c\in k\}}.
\]
Since \(\{c^2+c\}=\ker\Tr_{k/\F_2}\), this quotient group has two elements;
the invariance of the sign just proved ensures that the definition is independent of the symplectic basis.
\end{definition}

For example, \(Q_1(x,y)=xy\) and \(Q_2(x,y)=x^2+xy+y^2\) have the same polar form over \(\F_2^2\),
but have three and one zeros, respectively. Also, \(Q(z)=z^2\) has zero polar form but exactly one preimage of each value.
These show, respectively, that in characteristic two neither the Arf invariant nor the restriction of the quadratic form to the polar radical can be omitted.

\section{Trace Quadratic Forms and Adjoints of Linearized Operators}
\label{sec:trace-quadratic}

Identifying an \(m\)-dimensional space with the additive space of \(E=\F_{q^m}\) allows the trace pairing to represent all quadratic forms uniformly.

\begin{proposition}\leavevmode\phantomsection\label{prop:trace-quadratic-operator}
For any \(k\)-linear operator \(L:E\to E\), the map
\[
Q_L(x)=\Tr_{E/k}(xL(x))
\]
is a quadratic form, and
\[
B_{Q_L}(x,y)=\Tr_{E/k}\bigl(x(L+L^*)(y)\bigr),\qquad
\rad B_{Q_L}=\ker(L+L^*).
\]
Conversely, every quadratic form on the additive space of \(E\) is some \(Q_L\).
\end{proposition}
\begin{proof}
Linearity gives the scaling condition. Expanding \(Q_L(x+y)\) and rewriting the cross terms by the adjoint definition gives the polarization formula.
Nondegeneracy of the trace pairing makes the radical of the polar form \(B_{Q_L}\) exactly the stated kernel.

Conversely, choose a basis and write any quadratic form as \(Q(x)=\sum_i a_i x_i^2+\sum_{i<j}b_{ij}x_ix_j\).
Define the bilinear form \(D(x,y)=\sum_i a_i x_i y_i+\sum_{i<j}b_{ij}x_i y_j\). Then \(D(x,x)=Q(x)\).
For each \(y\), the trace pairing represents the functional \(x\mapsto D(x,y)\) uniquely as \(x\mapsto\Tr(xL(y))\).
Uniqueness makes \(y\mapsto L(y)\) linear, giving \(Q=Q_L\).
\end{proof}

Theorem~\ref{thm:linear-map-representation} also expresses \(L\) as a linear combination of Frobenius powers.
If \(D_L\) is the Dickson matrix of this reduced linearized representation, Propositions~\ref{prop:dickson-operations} and~\ref{prop:dickson-rank} give
\[
D_{L+L^*}=D_L+D_L^{\mathsf T},\qquad
\rank_k B_{Q_L}=\rank_E(D_L+D_L^{\mathsf T}),
\]
and
\[
\dim_k\rad B_{Q_L}=m-\rank_E(D_L+D_L^{\mathsf T}).
\]
The second equality also uses nondegeneracy of the trace pairing, which ensures that \(B_{Q_L}\) and \(L+L^*\) have equal rank.
Thus the polar rank can be computed directly from the coefficients of the linearized polynomial. In characteristic two, \(D_L+D_L^{\mathsf T}\) is alternating, so its rank is even,
also explaining the evenness of polar rank through its matrix representation.
The linear operator \(L\) representing a given quadratic form \(Q\) is generally not unique; for a fixed \(L\), its reduced linearized polynomial representation is uniquely determined by Theorem~\ref{thm:linear-map-representation}.
Polarization retains only \(L+L^*\); in characteristic two, the values of \(Q_L\) on this kernel must also be retained.

For the common one-term trace quadratic form
\[
Q_a(x)=\Tr_{E/k}(a x^{q^h+1}),\qquad 0\le h<m,
\]
take \(L=M_a\sigma^h\). The adjoint formula of Proposition~\ref{prop:linearized-adjoint} gives
\[
\rad B_{Q_a}
=\{x\in E:a x^{q^h}+a^{q^{m-h}}x^{q^{m-h}}=0\}.
\]
After finding the radical, compute the discriminant of the nondegenerate part in odd characteristic.
In characteristic two, first test the restriction of \(Q_a\) to the radical: if nonzero, all values have equal numbers of preimages;
if zero, compute the Arf invariant of the nondegenerate quadratic form on the quotient.
The corresponding counting theorem then determines the value distribution of the quadratic equation.

\begin{example}[The norm form in a quadratic extension]\leavevmode\phantomsection\label{exm:quadratic-norm}
Let \(E=\F_{q^2}\) and \(c=\Tr_{E/k}(a)\). Since \(x^{q+1}\in k\),
\[
\Tr_{E/k}(a x^{q+1})=c\Nm_{E/k}(x).
\]
If \(c=0\), this is the zero quadratic form; if \(c\ne0\), its only zero is zero, and each nonzero value has exactly \(q+1\) preimages.
This follows directly from norm-fiber sizes and agrees with the anisotropic binary counts above in both odd and even characteristic.
\end{example}

\paragraph{Diagonal recovery of a sesquilinear form.}
\phantomsection\label{v113:par-diagonal-question}
The preceding trace quadratic forms show that an ordinary bilinear representation need not be determined uniquely by its diagonal values: polarization gives only \(B(u,v)+B(v,u)\), which need not separate the two cross terms. Does the same loss persist if the second variable instead scales through a nontrivial field conjugation?
After subtracting the known diagonal terms, comparing the diagonal values at \(u+v\) and \(u+cv\) then assigns the coefficients \(c^q\) and \(c\) to the two cross terms. Choosing \(c\notin\F_q\) makes the two observations independent, suggesting recovery of the whole form; the next proposition gives a single formula valid in every characteristic.
The data retained here are the full field-valued diagonal function. Taking a further trace can destroy the information needed for recovery. The \(\F_4\) example following the proposition will compare these two kinds of observation explicitly.

\begin{proposition}\label{v106:prop-sesquilinear-recovery}
Let \(K=\F_{q^2}\), let \(V\) be a \(K\)-vector space, and let \(h:V\times V\to K\) be \(K\)-linear in the first variable and additive in the second, with \(h(u,cv)=c^qh(u,v)\). This is a \textbf{\(q\)-sesquilinear form}. Put \(d(x)=h(x,x)\), choose \(c\in K\setminus\F_q\), and set
\[
\begin{aligned}
S&=d(u+v)-d(u)-d(v),\\
T&=d(u+cv)-d(u)-c^{q+1}d(v).
\end{aligned}
\]
In every characteristic,
\begin{equation}\label{v106:eq-sesquilinear-recovery}
\boxed{h(u,v)=\frac{T-cS}{c^q-c}.}
\end{equation}
If \(V=E=\F_{q^{2n}}\), there is also a unique representation
\[
h(u,v)=\Tr_{E/K}\bigl(uL(v^q)\bigr),\qquad
L(Z)=\sum_{i=0}^{n-1}a_iZ^{q^{2i}},\quad a_i\in E.
\]
\end{proposition}
\begin{proof}
Expansion in the two variables gives \(S=h(u,v)+h(v,u)\) and \(T=c^qh(u,v)+ch(v,u)\). Subtracting and dividing by the nonzero \(c^q-c\) proves recovery, without dividing by two or assuming Hermitian symmetry.
For \(V=E\), nondegeneracy of the trace pairing gives a unique \(A(v)\in E\) for each \(v\), such that \(h(u,v)=\Tr_{E/K}(uA(v))\). Uniqueness makes \(A\) \(q\)-semilinear. Put \(\tau(v)=v^q\); then \(L=A\circ\tau^{-1}\) is \(K\)-linear and \(A(v)=L(v^q)\). Theorem~\ref{thm:linear-map-representation} gives its unique reduced \(q^2\)-linearized representative. Note that \(\tau^{-1}(z)=z^{q^{2n-1}}\); for \(n>1\), \(\tau\) must not be mistaken for an involution of \(E\).
\end{proof}
\begin{example}[The full diagonal and a further trace]\label{v106:ex-hermitian-f4}
On \(K=\F_4\), \(V=K^2\), take
\[
h(x,y)=x_1y_1^2+x_2y_2^2,\qquad d(x)=x_1^3+x_2^3.
\]
This is a Hermitian form, since \(h(y,x)=h(x,y)^2\). Every nonzero coordinate has cube one, giving
\[
\#d^{-1}(0)=1+3^2=10,\qquad \#d^{-1}(1)=2\cdot3=6.
\]
The diagonal values of any Hermitian form belong to \(\F_q\), so \(\Tr_{K/\F_q}d(x)=2d(x)\). In characteristic two this diagonal trace function is identically zero, although the two-variable trace form can be nonzero. If \(\omega^2+\omega+1=0\), \(x=(1,0)\), and \(y=(\omega,0)\), then
\[
\Tr_{\F_4/\F_2}h(x,y)=\Tr_{\F_4/\F_2}(\omega^2)=1.
\]
The full field-valued diagonal recovers the form by~\eqref{v106:eq-sesquilinear-recovery}; a further trace can erase all diagonal information. These are different kinds of data.
\end{example}

\subsection{Fourier Support of Quadratic Phases}
\label{sec:quadratic-fourier-support}

Over the base field \(\F_2\), encoding a quadratic form \(Q\) as the function \(x\mapsto(-1)^{Q(x)}\) gives a \(\{1,-1\}\)-valued \textbf{quadratic phase function}.
The polar radical determines the direction of the Fourier support, its dimension determines the magnitudes of nonzero Fourier coefficients, and the restriction of the quadratic form to the radical determines the affine coset containing the support.

\begin{proposition}[Fourier support and magnitudes of quadratic phases]\leavevmode\phantomsection\label{prop:quadratic-fourier-support}
Let \(E=\F_{2^n}\) and let \(Q:E\to\F_2\) be a quadratic form. Put \(R=\rad B_Q\) and \(d=\dim_{\F_2}R\), and write \(\Tr=\Tr_{E/\F_2}\).
For \(F_Q(x)=(-1)^{Q(x)}\), the Fourier convention of this chapter gives
\[
\widehat F_Q(b)=\sum_{x\in E}(-1)^{Q(x)+\Tr(bx)}.
\]
Define
\[
\mathcal S_Q=\{b\in E:\Tr(br)=Q(r)\text{ for all }r\in R\}.
\]
Then \(\mathcal S_Q\) is nonempty and is an affine coset of \(R^\perp\), and
\begin{equation}\label{eq:quadratic-fourier-support}
\boxed{\quad \operatorname{supp}\widehat F_Q=\mathcal S_Q,
\qquad
|\widehat F_Q(b)|=
\begin{cases}
2^{(n+d)/2},&b\in\mathcal S_Q,\\
0,&b\notin\mathcal S_Q.
\end{cases}\quad}
\end{equation}
If \(Q(x)=\Tr(xL(x))\), where \(L\) is an \(\F_2\)-linear map, put \(A=L+L^*\).
Then \(R=\ker A\) and \(R^\perp=\im A\), so \(\mathcal S_Q=b_0+\im A\), where \(b_0\) is any element satisfying the restriction condition above.
In particular, if \(Q|_R=0\), then \(\mathcal S_Q=\im A\), and for \(b=Az\),
\begin{equation}\label{eq:quadratic-fourier-completion}
\widehat F_Q(b)=(-1)^{Q(z)}\widehat F_Q(0).
\end{equation}
The right-hand side is independent of the choice of solution \(z\) of \(Az=b\).
\end{proposition}
\begin{proof}
Since \(B_Q\) vanishes on \(R\times R\), the restriction \(Q|_R\) is additive and hence an \(\F_2\)-linear functional.
Extend it to a linear functional on \(E\) and use nondegeneracy of the trace pairing to represent it as \(x\mapsto\Tr(b_0x)\), giving \(b_0\in\mathcal S_Q\).
Two parameters have the same restriction if and only if their difference lies in \(R^\perp\). Hence \(\mathcal S_Q=b_0+R^\perp\).

Write \(S_b=\widehat F_Q(b)\). All terms are real; substituting \(y=x+z\) in \(S_b^2\) gives
\[
\begin{aligned}
S_b^2
&=\sum_{z\in E}(-1)^{Q(z)+\Tr(bz)}
       \sum_{x\in E}(-1)^{B_Q(x,z)}\\
&=2^n\sum_{z\in R}(-1)^{Q(z)+\Tr(bz)}.
\end{aligned}
\]
In the second step, if \(z\in R\), the inner sum is \(2^n\); if \(z\notin R\),
then \(x\mapsto B_Q(x,z)\) is a nonzero linear functional and its character sum is zero by orthogonality.
On \(R\), the map \(Q(z)+\Tr(bz)\) is also linear, so the last sum is \(2^d\) when \(b\in\mathcal S_Q\) and zero otherwise.
This proves the support and magnitude formulas. The polar form is alternating, so its rank \(n-d\) is even, consistent with the integer exponent in the magnitude.

If \(Q=Q_L\), Proposition~\ref{prop:trace-quadratic-operator} gives \(R=\ker A\).
Also, \(A^*=A\), so kernel--image orthogonality for adjoints gives \(R^\perp=\im A\).
When \(Q|_R=0\), we may take \(b_0=0\). For \(b=Az\), polarization gives
\[
Q(x)+\Tr(bx)=Q(x+z)+Q(z).
\]
Translating the summation variable in the Fourier sum proves~\eqref{eq:quadratic-fourier-completion}.
If another solution is \(z+r\) with \(r\in R\), then \(Q(z+r)=Q(z)+Q(r)+B_Q(z,r)=Q(z)\), so the sign is independent of the chosen solution.
\end{proof}

Thus the kernel of \(L+L^*\) determines the direction of the Fourier support, and its dimension determines the magnitude of every nonzero Fourier coefficient.
In characteristic two, \(Q|_R\) must also be retained: it determines whether the support is \(\im(L+L^*)\) itself or a nonzero coset of that subspace.
Writing \(N_{Q,b}(c)=|\{x:Q(x)+\Tr(bx)=c\}|\), for \(c\in\F_2\) we have
\[
N_{Q,b}(c)=\frac{2^n+(-1)^c\widehat F_Q(b)}2.
\]
Thus the same proposition controls the solution counts of all linear perturbations of the quadratic equation.

\subsection{Differential Kernels and Trace Components of a Field-Valued Quadratic Map}

The preceding arguments fix a quadratic form and examine its polar operator. When a quadratic map takes values in an extension field, the direction and trace parameters give two families of operators whose kernels are related by the trace pairing.

\begin{definition}\label{v70:def-field-valued-quadratic}
Let \(E=\F_{2^n}\), with \(n\ge1\). A map \(F:E\to E\) is a \textbf{field-valued quadratic map} if \(F(0)=0\) and
\[
B_F(a,x)=F(a+x)+F(a)+F(x)
\]
is \(\F_2\)-bilinear. Linear terms are allowed in this definition.
For \(a,b\in E\), put
\[
L_a(x)=B_F(a,x),\qquad Q_b(x)=\Tr_{E/\F_2}(bF(x)),\qquad R_b=\rad B_{Q_b}.
\]
The quadratic form \(Q_b\) is a \textbf{trace component} of \(F\). All adjoints below are taken with respect to the absolute trace pairing \(\Tr_{E/\F_2}(xy)\), and all dimensions are over \(\F_2\).
\end{definition}

\begin{proposition}[Incidence between differential kernels and trace-component radicals]\label{v70:prop-dual-kernel-incidence}
With this notation, for every \(a\in E\), each nonempty solution set of \(F(x+a)+F(x)=c\) is a coset of \(\ker L_a\). For every \(b\in E\),
\begin{equation}\label{v70:eq-dual-kernel-incidence}
R_b=\{a\in E:L_a^*(b)=0\}.
\end{equation}
If \(r_a=\dim_{\F_2}\ker L_a\) and \(s_b=\dim_{\F_2}R_b\), then
\begin{equation}\label{v70:eq-dual-kernel-count}
\sum_{a\in E^\times}(2^{r_a}-1)
=\sum_{b\in E^\times}(2^{s_b}-1).
\end{equation}
\end{proposition}
\begin{proof}
By definition, \(F(x+a)+F(x)=L_a(x)+F(a)\). Thus, if the equation has one solution \(x_0\), its full solution set is \(x_0+\ker L_a\).
The polar form of \(Q_b\) is
\[
B_{Q_b}(a,x)=\Tr_{E/\F_2}(bL_a(x))
=\Tr_{E/\F_2}(L_a^*(b)x).
\]
Nondegeneracy of the trace pairing shows that this vanishes for all \(x\) if and only if \(L_a^*(b)=0\), proving~\eqref{v70:eq-dual-kernel-incidence}.

Now count the set
\[
\{(a,b)\in E^\times\times E^\times:L_a^*(b)=0\}
\]
in two ways. For fixed \(a\ne0\), an operator and its adjoint have the same rank, so \(\dim\ker L_a^*=r_a\); there are \(2^{r_a}-1\) choices of nonzero \(b\).
For fixed \(b\ne0\), equation~\eqref{v70:eq-dual-kernel-incidence} gives \(2^{s_b}-1\) choices of nonzero \(a\). Equating the counts proves~\eqref{v70:eq-dual-kernel-count}.
\end{proof}

The kernel \(\ker L_a\) fixes the direction \(a\) and varies the input, whereas \(R_b\) fixes the trace parameter \(b\) and varies the direction; the two objects cannot be identified merely by using the same subscript. The zero parameters give \(L_0=0\) and \(R_0=E\); the counting identity removes the zero parameter on both sides. A map with a constant term should first be normalized as \(F(x)-F(0)\).

\begin{example}[Differential kernels and all trace-component radicals of the cubic map]\label{v70:ex-cubic-dual-kernels}
Take \(F(x)=x^3\) on \(E=\F_{2^n}\), and abbreviate \(\Tr=\Tr_{E/\F_2}\). For \(a\ne0\),
\[
L_a(x)=a^2x+ax^2,\qquad L_a(at)=a^3(t+t^2),\qquad
\ker L_a=\{0,a\}.
\]
Proposition~\ref{prop:linearized-adjoint} gives
\[
L_a^*(b)=a^2b+(ab)^{2^{n-1}},\qquad
(L_a^*(b))^2=ab(a^3b+1).
\]
Since squaring is injective, for \(b\ne0\),
\begin{equation}\label{v70:eq-cubic-radical}
R_b=\{0\}\cup\{a\in E^\times:a^3b=1\}.
\end{equation}

If \(n\) is odd, then \(\gcd(3,2^n-1)=1\), so every \(b\ne0\) has exactly one nonzero radical element \(u\), and \(\dim R_b=1\). However,
\[
Q_b(u)=\Tr(bu^3)=\Tr(1)=1.
\]
Thus \(Q_b(x+u)=Q_b(x)+1\), and translation exchanges the two value fibers, each of size \(2^{n-1}\). In this case the radical is not contained in the zero set.

If \(n\) is even, then \(E\) contains \(\F_4\). Exactly \((2^n-1)/3\) nonzero parameters \(b\) have \(b^{-1}\) a cube. For any such parameter, choose \(u^3b=1\); then \(R_b=u\F_4\) and \(\dim R_b=2\). The remaining \(2(2^n-1)/3\) nonzero parameters have \(R_b=\{0\}\). In the former case, for \(\lambda\in\F_4^\times\),
\[
Q_b(u\lambda)=\Tr(\lambda^3)=\Tr(1)=0.
\]
Hence every \(Q_b\) vanishes on its radical and descends to \(E/R_b\); the exact value distribution also depends on the quadratic type on the quotient.
These conclusions include \(n=1\) and \(n=2\): in the former case \(Q_1(x)=x\), while in the latter \(Q_1=0\) and the components corresponding to the other two nonzero parameters have zero radical.

Every nonzero direction has a one-dimensional differential kernel. For odd \(n\), both sides of~\eqref{v70:eq-dual-kernel-count} equal \(2^n-1\); for even \(n\), its right-hand side is \(((2^n-1)/3)(2^2-1)=2^n-1\). Examples~\ref{exm:f64-polar} and~\ref{exm:f64-quadratic} determine the radical and quadratic type, respectively, when \(n=6,b=1\).
\end{example}

\subsection{Mixed traces, eigenlines, and semilinear returns}
Trace and determinant determine the eigenvalues of a two-by-two matrix. Adding $\operatorname{tr}(SS^\tau)$ detects the position of its eigenlines relative to Frobenius.

\begin{lemma}[Mixed traces and eigenlines]\label{v150:lem-mixed-trace}
Let $K$ be a finite field of characteristic two, let $\tau$ be an automorphism, and write
$S=Q\operatorname{diag}(\lambda,\mu)Q^{-1}$, where $Q\in\operatorname{GL}_2(K)$ and $\lambda\ne\mu$. Put
\[
 \begin{gathered}
 Q^{-1}Q^\tau=\begin{pmatrix}\alpha&\beta\\\gamma&\delta\end{pmatrix},
 \qquad\Delta=\alpha\delta+\beta\gamma,\\
 A=\lambda\lambda^\tau+\mu\mu^\tau,
 \qquad B=\lambda\mu^\tau+\mu\lambda^\tau.
 \end{gathered}
\]
Then
\begin{equation}\label{v150:eq-mixed-trace}
 \operatorname{tr}(SS^\tau)=A+\frac{\beta\gamma}{\Delta}(A+B),
 \qquad A+B=(\lambda+\mu)(\lambda^\tau+\mu^\tau)\ne0.
\end{equation}
In particular, the mixed trace equals $A$ if and only if $\beta\gamma=0$; then at least one eigenline is preserved by $\tau$.
Moreover, nonzero two-by-two matrices $R_1,R_2$ satisfying $R_i^2=0$ and $\operatorname{tr}(R_1R_2)=0$ have equal kernels.
\end{lemma}
\begin{proof}
Conjugate $SS^\tau$ to $DCD^\tau C^{-1}$, where $D=\operatorname{diag}(\lambda,\mu)$ and $C=Q^{-1}Q^\tau$.
Its trace is $(\alpha\delta A+\beta\gamma B)/\Delta$, proving~\eqref{v150:eq-mixed-trace}.
If $\gamma=0$, the first column of $Q$ spans a $\tau$-stable line; if $\beta=0$, the second column does.

For the last assertion, put $R_1$ in the form $\left(\begin{smallmatrix}0&1\\0&0\end{smallmatrix}\right)$ and write
$R_2=\left(\begin{smallmatrix}u&v\\w&u\end{smallmatrix}\right)$. The mixed trace gives $w=0$, and $\det R_2=0$ gives $u^2=0$. Thus $u=0$, $v\ne0$, and the kernels agree.
\end{proof}

\begin{proposition}[Restrictions on mixed traces from an irreducible return]\label{v150:prop-mixed-return}
Let $K$ be a finite field of characteristic two, let $\tau\in\operatorname{Aut}(K)$ have order $h$, and suppose $S\in\operatorname{GL}_2(K)$ and
$\Pi=S S^\tau\cdots S^{\tau^{h-1}}$ has irreducible characteristic polynomial over $K$.
If $S$ has distinct eigenvalues $\lambda,\mu$ in $K$, then
\[
 \operatorname{tr}(SS^\tau)\ne\lambda\lambda^\tau+\mu\mu^\tau.
\]
If the characteristic polynomial of $S$ has a repeated root, then $\operatorname{tr}(SS^\tau)\ne0$.
\end{proposition}
\begin{proof}
Failure of the first assertion would give, by the lemma, a $K$-line preserved by both $S$ and $\tau$. The semilinear operator $S\tau$ and its return $\Pi$ would preserve that line, contradicting irreducibility.
In the repeated-root case, write $S=\lambda I+R$, with $R^2=0$. If $R=0$, the return is scalar and cannot have irreducible characteristic polynomial.
If the mixed trace were zero, then $\operatorname{tr}(RR^\tau)=0$, so the lemma gives $\ker R=\ker R^\tau=\tau(\ker R)$.
This again supplies a line preserved by both $S$ and $\tau$, with the same contradiction.
\end{proof}

For example, let $K=\F_{2^m}$, with $m$ odd and $3\nmid m$, and let $\tau$ be any Frobenius generator of $K/\F_2$. Put
\[
 M=\begin{pmatrix}0&1\\1&1\end{pmatrix},\qquad S=P^{-1}MP^\tau,
 \quad P\in\operatorname{GL}_2(K).
\]
Cancellation of adjacent factors gives
\begin{equation}\label{v150:eq-mixed-telescope}
 S S^\tau\cdots S^{\tau^{m-1}}=P^{-1}M^mP.
\end{equation}
Since $M^3=I$ and $Z^2+Z+1$ is irreducible over an odd-degree binary extension, the characteristic polynomial on the right is irreducible. The restrictions above can therefore test mixed traces forced by local difference equations.

\begin{example}[Trace coefficients and additional difference-kernel elements]\label{v150:ex-apn-boundary}
For a finite-dimensional $\F_2$-vector space $V$, a map $G:V\to V$ is \textbf{almost perfect nonlinear} (APN) if $G(X+A)+G(X)=B$ has at most two solutions for every $A\ne0$ and every $B$.
For a quadratic map over $\F_2$, this is equivalent to requiring the linear map
\[
 D_AG(X)=G(X+A)+G(X)+G(A)+G(0)
\]
to have kernel $\{0,A\}$ for every $A\ne0$.

In $K=\F_{32}=\F_2(\alpha)$, where $\alpha^5+\alpha^2+1=0$, write $\Tr=\Tr_{K/\F_2}$, $t_2(z)=z+z^2$, and put
\[
 q_j(x,y)=x^{2^j+1}+xy^{2^j}+y^{2^j+1},
\]
\[
 \begin{aligned}
 G_c(x,y)=\bigl(&q_2(x,y)+t_2(xy)+c\Tr(xy),\\
              &q_1(x,y)+xy+\Tr(xy)\bigr).
 \end{aligned}
\]
Take $c=\alpha$, $A=(0,\alpha)$, and $X=(\alpha^4+\alpha,\alpha^3+\alpha^2+1)$. Using $\alpha^5=\alpha^2+1$ gives
\[
 \det(A,X)=1,\qquad D_Aq_2(X)=\alpha,\qquad D_Aq_1(X)=0.
\]
Since $t_2(1)=0$ and $\Tr(1)=1$, we obtain $D_AG_\alpha(X)=0$.
Its kernel thus contains the four distinct elements $0,A,X,X+A$, so $G_\alpha$ is not APN.
Here the difference-kernel size depends on the extension-field coefficient of the trace term, and cannot be inferred from $q_2,q_1$ alone.
\end{example}

\section{Character Lifts to Extensions and Weighted Orbit Products}
\label{sec:hasse-davenport}

The trace and norm lift additive and multiplicative characters, respectively, from the base field to an extension.
In an extension-field sum, the number of repetitions of a base-field Frobenius orbit controls both lifts simultaneously.

Thus the geometric-series contribution of each orbit in~\eqref{eq:finite-field-euler} can be weighted by its trace and norm;
the lifting formula below follows from the same unique-factorization expansion and formal logarithm calculation.

\begin{theorem}[Hasse--Davenport lifting formula]\label{thm:hasse-davenport}
Let \(k=\F_q\) and \(E_n=\F_{q^n}\), let \(\chi\) be a multiplicative character of \(k^\times\), and let \(\psi\) be a nontrivial additive character of \(k\).
Put
\[
\chi_n=\chi\circ\Nm_{E_n/k},\qquad
\psi_n=\psi\circ\Tr_{E_n/k}.
\]
All multiplicative characters, including the trivial one, take value zero at zero. Then
\[
\boxed{\quad G_{E_n}(\chi_n,\psi_n)=(-1)^{n-1}G_k(\chi,\psi)^n.\quad}
\]
\end{theorem}
\begin{proof}
For a monic polynomial \(F(X)=X^d+c_{d-1}X^{d-1}+\cdots+c_0\) of positive degree, define the weight
\[
w(F)=\chi((-1)^d c_0)\psi(-c_{d-1}),\qquad w(1)=1.
\]
Multiplying two polynomials multiplies their constant terms and adds their root sums, so \(w(FH)=w(F)w(H)\).
This remains valid when a constant term vanishes, since both sides are then zero.

In the formal power series ring \(\C[[u]]\), put
\[
\mathcal L(u)=\sum_{F\text{ monic}}w(F)u^{\deg F}.
\]
The linear polynomial \(X-a\) has weight \(\chi(a)\psi(a)\), so the linear coefficient is \(G_k(\chi,\psi)\).
For \(d\ge2\), the coefficients \(c_0\) and \(c_{d-1}\) can be chosen independently, and the other \(d-2\) coefficients are free. Hence
\[
\sum_{\deg F=d}w(F)
=q^{d-2}\left(\sum_{a\in k}\chi((-1)^d a)\right)
\left(\sum_{b\in k}\psi(-b)\right)=0.
\]
Thus
\begin{equation}\label{eq:gauss-linear-L}
\mathcal L(u)=1+G_k(\chi,\psi)u.
\end{equation}

On the other hand, unique factorization and multiplicativity of the weight give the Euler product
\[
\mathcal L(u)=\prod_{P\text{ monic irreducible}}
\bigl(1-w(P)u^{\deg P}\bigr)^{-1}.
\]
Only finitely many factors contribute to any fixed degree, so this identity is entirely formal.
Taking formal logarithms and collecting by degree gives
\begin{equation}\label{eq:weighted-orbit-log}
\log\mathcal L(u)
=\sum_{n\ge1}\frac{u^n}{n}
\sum_{d\mid n}d\sum_{\deg P=d}w(P)^{n/d}.
\end{equation}
Let \(\alpha\in E_n\) have minimal polynomial \(P\) of degree \(d\mid n\). The root sum and product give
\[
w(P)=\chi(\Nm_{\F_{q^d}/k}(\alpha))
\psi(\Tr_{\F_{q^d}/k}(\alpha)).
\]
When \(\alpha\ne0\), the list of conjugates in the extension repeats this orbit \(n/d\) times, so
\[
\Tr_{E_n/k}(\alpha)=\frac nd\Tr_{\F_{q^d}/k}(\alpha),\qquad
\Nm_{E_n/k}(\alpha)=\Nm_{\F_{q^d}/k}(\alpha)^{n/d}.
\]
Using \(\psi(u)^j=\psi(ju)\) gives
\[
w(P)^{n/d}=\chi_n(\alpha)\psi_n(\alpha).
\]
When \(\alpha=0\), both sides vanish. Each irreducible polynomial of degree \(d\) contributes exactly \(d\) roots,
so the inner sum in~\eqref{eq:weighted-orbit-log} is precisely \(G_{E_n}(\chi_n,\psi_n)\).
Hence
\[
\mathcal L(u)=\exp\left(\sum_{n\ge1}G_{E_n}(\chi_n,\psi_n)\frac{u^n}{n}\right).
\]
Compare with~\eqref{eq:gauss-linear-L}, expand \(\log(1+Gu)=\sum_{n\ge1}(-1)^{n-1}G^n u^n/n\), and compare coefficients to obtain the formula.
\end{proof}

This is weighted orbit counting: an irreducible polynomial corresponds to an orbit, whose product and sum determine its weight.
The unweighted form was established in Proposition~\ref{prop:finite-field-euler}; here the number of repetitions of an orbit also becomes the exponent of its weight.
For the \(L\)-function interpretation of this weighted Euler product, see \cite[Section 3]{ConradGaussL}.

The lifting formula also gives a product relation among Gauss sums over the same field.
The key is to combine the fibers of a power map with power sums of roots, then compare the leading coefficient of one weighted orbit product in two ways.

\begin{theorem}[Davenport--Hasse multiplication formula]\label{v91:thm-dh-product}
Let \(k=\F_q\), let \(m\mid q-1\), and let \(\eta\) be a multiplicative character of \(k^\times\) of order \(m\).
Let \(\chi\) be any multiplicative character and \(\psi\) a nontrivial additive character.
Write \(G(\rho)=G_k(\rho,\psi)\), retaining the convention that every multiplicative character takes the value zero at zero. Then
\begin{equation}\label{v91:eq-dh-product}
\boxed{\;
\prod_{j=0}^{m-1}G(\chi\eta^j)
=-\chi(m^{-m})G(\chi^m)
\prod_{j=0}^{m-1}G(\eta^j).
\;}
\end{equation}
Here \(m\) is regarded as a nonzero element of \(k\); the assumption \(m\mid q-1\) ensures that \(p\nmid m\).
\end{theorem}
\begin{proof}
Put \(E_n=\F_{q^n}\), using the subscript \(n\) for characters lifted by norm or trace. Define
\[
S_n(\chi)=\sum_{x\in E_n^\times}\chi_n(x)^m\psi_n(x^m).
\]
Surjectivity of the norm ensures that \(\eta_n\) still has order \(m\).
For \(y\in E_n^\times\), the equation \(x^m=y\) has
\(\sum_{j=0}^{m-1}\eta_n^j(y)\) solutions, and each solution satisfies \(\chi_n(x)^m=\chi_n(y)\). Hence
\[
S_n(\chi)=\sum_{j=0}^{m-1}G_{E_n}(\chi_n\eta_n^j,\psi_n)
=(-1)^{n-1}\sum_{j=0}^{m-1}G(\chi\eta^j)^n,
\]
where the last equality follows from Theorem~\ref{thm:hasse-davenport}.
In \(\C[[u]]\), put
\[
\mathcal L_\chi(u)=
\exp\left(\sum_{n\ge1}S_n(\chi)\frac{u^n}{n}\right).
\]
Expanding the formal logarithm gives
\begin{equation}\label{v91:eq-dh-product-L}
\mathcal L_\chi(u)=\prod_{j=0}^{m-1}\bigl(1+G(\chi\eta^j)u\bigr).
\end{equation}

We now compute the coefficient of \(u^m\) from polynomial coefficients.
For a monic polynomial \(F\in k[X]\) with roots \(\alpha_1,\ldots,\alpha_d\), listed with multiplicity, define
\[
w_\chi(F)=
\chi^m\!\left(\prod_{i=1}^d\alpha_i\right)
\psi\!\left(\sum_{i=1}^d\alpha_i^m\right),\qquad w_\chi(1)=1.
\]
Both the product and the sum lie in \(k\), and \(w_\chi(FH)=w_\chi(F)w_\chi(H)\).
Unique factorization gives
\[
\sum_{F\text{ monic}}w_\chi(F)u^{\deg F}
=\prod_{P\text{ monic irreducible}}
\bigl(1-w_\chi(P)u^{\deg P}\bigr)^{-1}.
\]
If \(P\) has degree \(d\mid n\) and \(\alpha\) is one of its roots, its orbit is repeated \(n/d\) times in \(E_n\), so
\[
w_\chi(P)^{n/d}=\chi_n(\alpha)^m\psi_n(\alpha^m).
\]
For \(P=X\), both sides are zero. Every other \(P\) contributes exactly \(d\) roots.
Taking the formal logarithm as in~\eqref{eq:weighted-orbit-log}, we obtain
\begin{equation}\label{v91:eq-dh-product-weight}
\mathcal L_\chi(u)=\sum_{F\text{ monic}}w_\chi(F)u^{\deg F}.
\end{equation}

Write \(F=X^m+a_1X^{m-1}+\cdots+a_m\).
Newton's identities (Lemma~\ref{v58:lem-newton}) give
\[
\sum_{i=1}^m\alpha_i^m
=R_m(a_1,\ldots,a_{m-1})-ma_m,
\]
where \(R_m\) is independent of \(a_m\). Put
\[
A_m=\sum_{a_1,\ldots,a_{m-1}\in k}
\psi\bigl(R_m(a_1,\ldots,a_{m-1})\bigr).
\]
In the coefficient of \(u^m\) in~\eqref{v91:eq-dh-product-weight}, substitute \(t=-ma_m\).
Since \(m(m+1)\) is even, this gives
\begin{align*}
[u^m]\mathcal L_\chi(u)
&=A_m\sum_{a_m\in k}\chi^m((-1)^ma_m)\psi(-ma_m)\\
&=A_m\chi(m^{-m})G(\chi^m).
\end{align*}
The auxiliary quantity \(A_m\) does not depend on \(\chi\).
Taking \(\chi=\varepsilon\) and using \(G(\varepsilon)=-1\), equation~\eqref{v91:eq-dh-product-L} yields
\[
\prod_{j=0}^{m-1}G(\eta^j)=-A_m.
\]
Substituting this into the identity for general \(\chi\) proves~\eqref{v91:eq-dh-product}.
\end{proof}

The proof allows \(\chi^m=\varepsilon\); this case need not be excluded.
When \(m=1\), we have \(A_1=1\), and the formula remains valid.
For the lifting formula, the freedom of the coefficient next to the leading term and the constant coefficient makes the weighted generating function linear.
For the multiplication formula, Newton's identities isolate the constant coefficient in the degree-\(m\) term, where it contributes a single Gauss sum.
Both arguments pass from root orbits to coefficients, but compare different coefficients.

\begin{example}[A product formula involving a quadratic character]\label{v91:ex-dh-quadratic}
Suppose that \(q\) is odd and \(\eta\) is the quadratic character.
Taking \(m=2\) and using \(G(\varepsilon)=-1\), we obtain
\[
G(\chi)G(\chi\eta)=\chi(4^{-1})G(\chi^2)G(\eta).
\]
The factor \(4^{-1}\) is \(2^{-2}\).
It can also be read directly from the root power sum
\(\alpha_1^2+\alpha_2^2=a_1^2-2a_2\) for \(F=X^2+a_1X+a_2\), together with the substitution above.
\end{example}

\begin{corollary}[Lifting Jacobi sums]\leavevmode\phantomsection\label{cor:jacobi-lifting}
If all \(\chi,\lambda,\chi\lambda\) are nontrivial, then
\[
J_{E_n}(\chi_n,\lambda_n)=(-1)^{n-1}J_k(\chi,\lambda)^n.
\]
\end{corollary}
\begin{proof}
Surjectivity of the norm keeps all three lifted multiplicative characters nontrivial, and surjectivity of the trace makes \(\psi_n\) nontrivial.
Use \(J(\chi,\lambda)=G(\chi,\psi)G(\lambda,\psi)/G(\chi\lambda,\psi)\) over each of the two fields,
then substitute the lifting formula. The numerator has two sign factors and the denominator one, leaving \((-1)^{n-1}\).
\end{proof}

\begin{example}[A trace quadratic equation over the field of order 64]\label{exm:f64-quadratic}
Continue Example~\ref{exm:f64-polar}, taking \(Q(x)=\Tr_{E/\F_2}(x^3)\) over \(E=\F_{64}\).
Its polar radical is known to be \(\F_4\), its polar rank is four, and \(Q|_{\F_4}=0\).
This also agrees with the general formula of Proposition~\ref{prop:trace-quadratic-operator}.
The corresponding character sum determines the remaining quadratic type.

Take \(k=\F_4=\F_2(\omega)\) and \(\omega^2+\omega+1=0\),
put \(\chi(\omega)=\zeta\), where \(\zeta\) is a primitive complex cube root of unity,
and take \(\psi(z)=(-1)^{\Tr_{k/\F_2}(z)}\). Since \(\psi(1)=1\) and \(\psi(\omega)=\psi(\omega^2)=-1\),
\[
G_k(\chi,\psi)=G_k(\chi^2,\psi)=1-\zeta-\zeta^2=2.
\]
The extension \(E/k\) has degree three, so the lifting formula gives eight for both corresponding Gauss sums.
The lifted characters \(\chi_3,\chi_3^2\) are exactly the two nontrivial cubic characters of \(E\),
and \(\psi_3(z)=(-1)^{\Tr_{E/\F_2}(z)}\). Using the number \(1+\chi_3(y)+\chi_3^2(y)\) of preimages in the equation \(x^3=y\) gives
\[
\sum_{x\in E}(-1)^{Q(x)}
=G_E(\chi_3,\psi_3)+G_E(\chi_3^2,\psi_3)=16.
\]
Thus \(\epsilon(Q)=1\), and
\[
|\{Q=0\}|=40,\qquad |\{Q=1\}|=24.
\]
The radical has four elements, whereas the zero set has 40: the radical is the kernel of the polar operator, not the zero set of the quadratic equation.
Proposition~\ref{prop:quadratic-fourier-support} further gives the full support of \(F_Q(x)=(-1)^{Q(x)}\).
Here \(L(x)=x^2\) and \(A=L+L^*\) is \(A(x)=x^2+x^{32}\), and
\[
\operatorname{supp}\widehat F_Q=\im A=(\F_4)^\perp
=\ker\Tr_{\F_{64}/\F_4}.
\]
This support has 16 elements, and all Fourier magnitudes on it are 16.
For \(b=Az\), we have \(\widehat F_Q(b)=16(-1)^{Q(z)}\); coefficients outside the support vanish.
Thus the linearized operator determines the polar rank and frequency support, Gauss sums over the extension determine the sign of the zero-frequency Fourier coefficient,
and the translation formula for the quadratic form determines the signs of the remaining Fourier coefficients on the support.
\end{example}

\subsection{Quadratic Frobenius returns and exact fibers}
\label{v145:sec-return-fibers}
A linearized root space turns irreducibility into a condition on a return operator. We now fix its characteristic polynomial and count the coefficients producing it. Normalized root bases turn this question into a power-sum equation, whose fiber size is determined by Hasse--Davenport lifting.

Let $P$ be a prime power, $K=\F_{P^k}$, and $F(x)=x^P$. For $a\in K$, put
\begin{equation}\label{v145:eq-quadratic-return}
\begin{split}
 L_a(X)&=X^{P^2}+aX^P+X,
 \qquad A_a=\begin{pmatrix}0&-1\\1&-a\end{pmatrix},\\
 M_k(a)&=A_aA_{a^P}\cdots A_{a^{P^{k-1}}},
 \qquad b_k(a)=-\operatorname{tr}M_k(a).
\end{split}
\end{equation}
Thus $M_k(a)=M_{k,1}(a,1)$ and $b_k(a)=\Theta_{k,1}(a)$, in the return-spectrum notation of Chapter 9.

Proposition~\ref{v145:prop-return-matrix} gives $b_k(a)\in\F_P$. The characteristic polynomial of $F^k$ on the two-dimensional $\F_P$-root space is $Z^2+b_k(a)Z+1$.

\begin{lemma}[Normalized root bases and target fibers]\label{v145:lem-normalized-bases}
Suppose that $G_b(Z)=Z^2+bZ+1$ is irreducible over $\F_P$, and choose $C\in\operatorname{SL}_2(\F_P)$ with characteristic polynomial $G_b$. For a column vector $u=(u_1,u_2)^{\mathsf T}$, put
\[
 \Delta(u)=\det(u,Fu),\qquad
 W_C=\{u\in(\overline{\F}_P)^2:F^ku=Cu\}.
\]
Then $W_C$ is a two-dimensional $\F_{P^k}$-space, $\Delta(W_C)\subseteq\F_{P^k}$, and
\begin{equation}\label{v145:eq-fiber-bases}
 (P+1)\#\{a\in K:b_k(a)=b\}
 =\#\{u\in W_C:\Delta(u)=1\}.
\end{equation}
Moreover, $\Delta(u)\ne0$ for $u\in W_C\setminus\{0\}$, and
\begin{equation}\label{v145:eq-delta-homogeneous}
 \Delta(\lambda u)=\lambda^{P+1}\Delta(u)
 \qquad(\lambda\in K).
\end{equation}
\end{lemma}
\begin{proof}
If $\Delta(u)=1$, then $\det(Fu,F^2u)=1$, so
\[
 F^2u=-u-a(u)Fu,\qquad a(u)=-\det(u,F^2u).
\]
Thus $u_1,u_2$ form a root basis of $L_{a(u)}$. Conversely, every root basis of $L_a$ has $\Delta^P=\Delta\ne0$. A base-field change of basis of determinant $\Delta^{-1}$ makes $\Delta=1$. These \textbf{normalized root bases} form a free transitive set under $\operatorname{SL}_2(\F_P)$.

When $a\in K$, each normalized basis determines $F^ku=C_u u$, where $C_u\in\operatorname{SL}_2(\F_P)$ has the same characteristic polynomial as $M_k(a)$. Matrices with a specified irreducible quadratic characteristic polynomial form one conjugacy class in $\operatorname{SL}_2(\F_P)$: first conjugate in the general linear group, then adjust the determinant of the conjugating matrix using the centralizer $\F_{P^2}^{\times}$. Its determinant is the surjective quadratic norm. Hence the centralizer of $C$ in the special linear group is the norm-one group of order $P+1$. Each $a$ with $b_k(a)=b$ gives exactly $P+1$ normalized bases with $C_u=C$. Conversely, $F^ku=Cu$ and $\det C=1$ give $a(u)^{P^k}=a(u)$, proving~\eqref{v145:eq-fiber-bases}.

Write $Q=P^k$. Choose eigenvalues $\alpha,\alpha^{-1}\in\F_{P^2}$ of $C$, with eigenvectors $e,e^P$. If $k$ is even, the condition on $u=ex+e^Py$ is $x^Q=\alpha x$, $y^Q=\alpha^{-1}y$, each with $Q$ solutions. If $k$ is odd, it is $x^{Q^2}=\alpha^2x$, $y=\alpha x^Q$, giving $Q^2$ solutions. Thus $W_C$ has $Q^2$ elements and is a two-dimensional $K$-space.

The value field follows from $\Delta(u)^Q=\det(Cu,CFu)=\Delta(u)$. If a nonzero $u$ had $\Delta(u)=0$, the projective point $[u]$ would be fixed by $F$ and hence lie in $\mathbb P^1(\F_P)$. It would also be fixed by $C$, contradicting irreducibility of $G_b$. Homogeneity follows directly from the determinant.
\end{proof}

\begin{lemma}[Character sums on the quadratic norm quotient]\label{v145:lem-norm-one-gauss}
Choose an isomorphism $\iota:\mu_{P+1}(\F_{P^2})\to\mu_{P+1}(\C)$, and put
\[
 \lambda(z)=\iota(z^{P-1}),\qquad
 \varepsilon_P=\iota(-1)=
 \begin{cases}1,&P\text{ even},\\-1,&P\text{ odd}.
 \end{cases}
\]
Let $Q=P^{2s}$, $\chi=\lambda\circ\Nm_{\F_Q/\F_{P^2}}$, and lift a nontrivial additive character from $\F_P$. For $1\le i\le P$,
\begin{equation}\label{v145:eq-semi-primitive-gauss}
 G_{\F_{P^2}}(\lambda^i)=P\varepsilon_P^i,
 \qquad
 G_{\F_Q}(\chi^i)=(-1)^{s-1}P^s\varepsilon_P^{si}.
\end{equation}
For $1\le i,j\le P$, the corresponding Jacobi sums are
\begin{equation}\label{v145:eq-return-jacobi}
 J(\chi^i,\chi^j)=
 \begin{cases}
 (-1)^{s-1}P^s,&i+j\not\equiv0\pmod{P+1},\\
 -1,&i+j\equiv0\pmod{P+1}.
 \end{cases}
\end{equation}
\end{lemma}
\begin{proof}
Choose $0\ne\beta\in\F_{P^2}$ with $\Tr_{\F_{P^2}/\F_P}(\beta)=0$. Then $\beta^{P-1}=-1$, so $\lambda^i(\beta)=\varepsilon_P^i$. Put
\[
 S_i=\sum_{\Tr(x)=1}\lambda^i(x).
\]
The sum on the nonzero trace-zero elements is $(P-1)\varepsilon_P^i$. Since $\lambda$ is trivial on $\F_P^\times$, every nonzero trace fiber has sum $S_i$. The sum over all nonzero elements is zero, so $S_i=-\varepsilon_P^i$. A nontrivial additive character sums to $-1$ on $\F_P^\times$, giving the Gauss sum $(P-1)\varepsilon_P^i-S_i=P\varepsilon_P^i$. Theorem~\ref{thm:hasse-davenport} gives the extension-field formula.

For nontrivial product character, Proposition~\ref{prop:jacobi} and $\varepsilon_P^{P+1}=1$ give the first case of~\eqref{v145:eq-return-jacobi}. Otherwise that proposition gives $-\chi^i(-1)=-1$: in odd characteristic $(Q-1)/(P+1)$ is even, and in even characteristic $-1=1$.
\end{proof}

\begin{lemma}[Semilinear normal forms and normalized cyclic vectors]\label{v146:lem-semilinear-normal}
Let \(\tau(x)=x^{P^r}\), \(\gcd(r,k)=1\), and fix \(b\in\mathcal I_P\). The \(\tau\)-semilinear transformations of \(K^2\) satisfying
\[
 \det(\Phi v,\Phi w)=\tau(\det(v,w)),\qquad
 \chi_{\Phi^k}(Z)=Z^2+bZ+1
\]
form a nonempty single \(\operatorname{SL}_2(K)\)-conjugacy class, with centralizer of order \(P+1\). A vector satisfying \(\det(v,\Phi v)=1\) is a \textbf{normalized cyclic vector}; denote their set by \(\mathcal C(\Phi)\). Then
\begin{equation}\label{v146:eq-cyclic-fiber}
 \#\{a\in K:\Theta_{k,r}(a)=b\}
 =\frac{\#\mathcal C(\Phi)}{P+1}.
\end{equation}
\end{lemma}
\begin{proof}
Choose a root \(\alpha\) of \(Z^2+bZ+1\), so \(\alpha^P=\alpha^{-1}\ne\alpha\).

If \(k=2s\), then \(r\) is odd and \(\Phi\) exchanges the two eigenlines of \(\Phi^k\). In an ordered eigenbasis of determinant one, it has the form
\begin{equation}\label{v146:eq-even-normal}
\begin{aligned}
 \Phi_w(x,y)&=(w\tau(y),-w^{-1}\tau(x)),\\
 \Phi_w^k&=\operatorname{diag}\bigl((-1)^sN^{1-P},(-1)^sN^{P-1}\bigr),
\end{aligned}
\end{equation}
where \(w\in K^\times\) and \(N=\Nm_{K/\F_{P^2}}(w)\). The second identity uses the fact that \(\tau^2\) generates the Galois group of \(K/\F_{P^2}\). The norm is surjective, and \(N\mapsto N^{1-P}\) has image \(\mu_{P+1}\), so the first eigenvalue can be prescribed as \(\alpha\).

If \(w,w'\) give the same first eigenvalue, put \(\delta=w'/w\). Then \(\Nm_{K/\F_{P^2}}(\delta)\in\F_P^\times\), and Hilbert 90 gives \(u\ne0\) such that
\[
 \frac{\tau^2(u)}u=\frac{\tau(\delta)}\delta,
 \qquad
 \operatorname{diag}(u,\delta^{-1}\tau(u))\Phi_w
 =\Phi_{w'}\operatorname{diag}(u,\delta^{-1}\tau(u)).
\]
An invertible \(K\)-linear map commuting with \(\Phi_w\) is diagonal because it commutes with \(\Phi_w^k\). The semilinear relation then makes it \(\operatorname{diag}(t,\tau(t))\), with \(t\in\F_{P^2}^\times\). Its determinant \(t^{P+1}\) ranges over \(\F_P^\times\).

If \(k\) is odd, put \(Q=P^k\) and \(E=\F_{P^{2k}}\). The quadratic remains irreducible over \(K\), so identify \(K^2\) with \(E\) such that \(\Phi^k\) is multiplication by \(\alpha\). Let \(j\in[0,2k)\) be the unique even integer with \(j\equiv r\pmod k\), and put \(\rho(z)=z^{P^j}\). Its restriction to \(K\) is \(\tau\), and it generates the Galois group of \(E/\F_{P^2}\). Commutation with multiplication by \(\alpha\) gives
\begin{equation}\label{v146:eq-odd-normal}
 \Phi_w(z)=w\rho(z),\qquad
 \Phi_w^k(z)=\Nm_{E/\F_{P^2}}(w)z.
\end{equation}
Choose \(\vartheta\in\F_{P^2}^\times\) with \(\vartheta^P=-\vartheta\), and put
\[
 \omega_\vartheta(z,z')=\vartheta(zz'^Q-z^Qz').
\]
This is a nondegenerate alternating form with values in \(K\). Any two such forms differ by a scalar in \(K^\times\). Multiplication maps have determinant equal to the surjective norm \(\Nm_{E/K}\), so the identification can be adjusted to carry \(\omega_\vartheta\) exactly to the standard determinant. Now
\[
 \omega_\vartheta(\Phi_wz,\Phi_wz')
 =w^{Q+1}\rho(\omega_\vartheta(z,z')),
\]
so the determinant-one condition is \(w^{Q+1}=1\). The image of \(\mu_{Q+1}\) under \(\Nm_{E/\F_{P^2}}\) is precisely \(\mu_{P+1}\): the identities
\[
 \frac{P^{2k}-1}{P^2-1}
 =\frac{Q+1}{P+1}\frac{Q-1}{P-1},\qquad
 \gcd\!\left(P+1,\frac{Q-1}{P-1}\right)=1
\]
show that its kernel has order \((Q+1)/(P+1)\). Thus \(w^{Q+1}=1\) and \(\Nm_{E/\F_{P^2}}(w)=\alpha\) can be imposed simultaneously.

Hilbert 90 supplies a multiplication map intertwining two parameters with the same norm. The invertible \(K\)-linear centralizer is \(\F_{P^2}^\times\): commutation with multiplication by \(\alpha\) first puts it in \(E^\times\), and commutation with \(w\rho\) then makes it fixed by \(\rho\). Its determinant is again \(t^{P+1}\), ranging over \(\F_P^\times\).

In both cases, an intertwining matrix \(g\) between determinant-one semilinear transformations satisfies \(\tau(\det g)=\det g\). Surjectivity of the determinant on the centralizer adjusts \(g\) to determinant one. There is therefore a single special linear conjugacy class, and its centralizer is a norm-one group of order \(P+1\).

Finally, if \(v\in\mathcal C(\Phi)\), then \((v,\Phi v)\) is a determinant-one basis and \(\det(\Phi v,\Phi^2v)=1\). There is a unique \(a\in K\) such that
\[
 \Phi^2v=-v-a\Phi v.
\]
In this basis \(\Phi=A_a\tau\). Conversely, every \(A_a\tau\) with the prescribed return polynomial is special-linearly conjugate to \(\Phi\). The bases giving a fixed \(a\) form a coset of the centralizer, of size \(P+1\). Each basis is determined by its first vector, proving~\eqref{v146:eq-cyclic-fiber}.
\end{proof}

\begin{theorem}[Uniform fibers of irreducible quadratic return polynomials]\label{v145:thm-return-uniform}
Let $\gcd(r,k)=1$. If $Z^2+bZ+1$ is irreducible over $\F_P$, then
\begin{equation}\label{v145:eq-return-fiber-size}
 \#\{a\in\F_{P^k}:\Theta_{k,r}(a)=b\}=D_k(P),
 \qquad
 D_k(P)=
 \begin{cases}
 \dfrac{P^k+1}{P+1},&k\text{ odd},\\[2mm]
 \dfrac{P^k+1-2(-P)^{k/2}}{P+1},&k\text{ even}.
 \end{cases}
\end{equation}
\end{theorem}
\begin{proof}
First take $r=1$. Use the notation of Lemma~\ref{v145:lem-normalized-bases}, with $Q=P^k$. It remains to count $u\in W_C$ with $\Delta(u)=1$.

If $k$ is odd, then $\gcd(P+1,Q-1)=\gcd(P+1,2)$. For even $P$, equation~\eqref{v145:eq-delta-homogeneous} gives a unique solution on each $K$-line, hence $Q+1$ solutions. For odd $P$, each line gives two solutions or none according as its nonzero $\Delta$-values are squares or nonsquares. Choose $T$ in the general linear centralizer of $C$ with determinant a nonsquare in $\F_P$; surjectivity of the quadratic norm permits this choice. An odd-degree extension preserves nonsquares, while $T$ preserves $W_C$ and satisfies
\[
 \Delta(Tu)=(\det T)\Delta(u).
\]
It therefore bijects the square and nonsquare classes of lines, again giving $Q+1$ solutions.

If $k=2s$, choose eigenvectors $e,e^P$ as above and put $\delta=\det(e,e^P)$, so $\delta^P=-\delta$. Choose $x_0\ne0$ with $x_0^{Q-1}=\alpha$. Every element of $W_C$ has the form
\[
 u=e x_0\xi+e^P x_0^P\eta,\qquad \xi,\eta\in\F_Q.
\]
Expansion of the determinant gives
\begin{equation}\label{v145:eq-return-diagonal}
 \Delta(u)=A\xi^{P+1}+A^P\eta^{P+1},
 \qquad A=\delta x_0^{P+1}\in\F_Q^\times.
\end{equation}
For the character $\chi$ of order $P+1$ from Lemma~\ref{v145:lem-norm-one-gauss},
\[
 \chi(A)=\varepsilon_P^s\iota(\alpha)=:\upsilon,
 \qquad \chi(A^P)=\upsilon^{-1},\qquad \upsilon^2\ne1.
\]
The first identity follows from $\delta^{P-1}=-1$ and $x_0^{Q-1}=\alpha$; the others follow from $\alpha^{P+1}=1$ and $\alpha\notin\F_P$.

Extend nontrivial multiplicative characters by zero at zero. The power-map fiber function is
\[
 R(z)=\#\{x\in\F_Q:x^{P+1}=z\}
     =1+\sum_{i=1}^{P}\chi^i(z),
\]
where the constant term is one also at zero. The required count is
\[
 \sum_{x+y=1}R(x/A)R(y/A^P)
 =Q+\sum_{i,j=1}^{P}\upsilon^{j-i}J(\chi^i,\chi^j).
\]
Terms containing only one nontrivial character vanish by orthogonality. Since $\upsilon\ne1$ and $\upsilon^2\ne1$,
\[
 \sum_{i,j=1}^{P}\upsilon^{j-i}=1,\qquad
 \sum_{\substack{1\le i,j\le P\\i+j\equiv0\ (P+1)}}\upsilon^{j-i}
 =\sum_{i=1}^{P}\upsilon^{-2i}=-1.
\]
Equation~\eqref{v145:eq-return-jacobi} therefore gives the total
\[
 Q+2(-1)^{s-1}P^s+1=Q+1-2(-P)^s.
\]
Dividing the two respective counts by $P+1$ proves the formulas for $r=1$.

Now let \(\gcd(r,k)=1\). By Lemma~\ref{v146:lem-semilinear-normal}, it suffices to compare normalized cyclic-vector counts for the same target.

If \(k=2s\), fix \(w\) in~\eqref{v146:eq-even-normal} and use \(\tau(x)=x^{P^r}\) and \(F(x)=x^P\), respectively. Their return eigenvalues agree, and the cyclic-vector conditions are
\[
 -w^{-1}x^{P^r+1}-wy^{P^r+1}=1,
 \qquad -w^{-1}x^{P+1}-wy^{P+1}=1.
\]
Since \(\gcd(P^r+1,P^k-1)=P+1\), the two power maps have the same image and fiber sizes on \(K^\times\); each has one preimage of zero. The two equations therefore have the same number of solutions.

If \(k\) is odd, retain \(E,Q,j,w,\omega_\vartheta\) from~\eqref{v146:eq-odd-normal} and put \(\ell=j+k\). Then \(\ell\) is odd, \(\gcd(\ell,2k)=1\), and
\[
 \omega_\vartheta(z,\Phi_wz)=\vartheta(w^QS-wS^Q),
 \qquad S=z^{P^\ell+1}.
\]
Replacing \(\rho\) by \(\rho_0(z)=z^{P^{k+1}}\) gives a transformation restricting to standard Frobenius on \(K\), with the same return multiplication by \(\alpha\); its corresponding \(S\) is \(z^{P+1}\). Since \(\gcd(P^\ell+1,P^{2k}-1)=P+1\), these two maps to \(S\) have identical images and fiber sizes on \(E\). Hence the equation \(\omega_\vartheta(z,\Phi_wz)=1\) has the same number of solutions. Equation~\eqref{v146:eq-cyclic-fiber} now gives \(D_k(P)\) for every generator.
\end{proof}

\subsection{Split Spectra, Repeated-Root Fibers, and Constant-Term Normalization}
\label{v146:sec-other-spectra}
An irreducible target has a quadratic norm-one group as centralizer. When both eigenvalues lie in the fixed field, normalized cyclic vectors still count the coefficient fiber, but the centralizer and the vector equation change.

\begin{theorem}[Uniform fibers of split return spectra with distinct roots]\label{v146:thm-split-fibers}
Let \(\gcd(r,k)=1\), \(\lambda\in\F_P^\times\), and \(\lambda^2\ne1\). Then
\begin{equation}\label{v146:eq-split-fibers}
 \#\{a\in\F_{P^k}:\Theta_{k,r}(a)=-(\lambda+\lambda^{-1})\}
 =\frac{P^k-1}{P-1}.
\end{equation}
\end{theorem}
\begin{proof}
Put \(\tau(x)=x^{P^r}\). A determinant-one semilinear transformation with return eigenvalues \(\lambda,\lambda^{-1}\) preserves both eigenlines. In a determinant-one eigenbasis it is
\[
 \Phi_w(x,y)=(w\tau(x),w^{-1}\tau(y)),\qquad
 \Nm_{K/\F_P}(w)=\lambda.
\]
Norm surjectivity gives existence, and Hilbert 90 makes all parameters of the same norm special-linearly conjugate. The special linear centralizer is \(\{\operatorname{diag}(t,t^{-1}):t\in\F_P^\times\}\), of order \(P-1\). The basis-counting argument at the end of Lemma~\ref{v146:lem-semilinear-normal} therefore gives fiber size \(\#\mathcal C(\Phi_w)/(P-1)\).

The normalization condition is \(w^{-1}x\tau(y)-wy\tau(x)=1\), so \(x\ne0\). Write \(y=xz\); then
\[
 w^{-1}\tau(z)-wz=(x\tau(x))^{-1}.
\]
The left side is an \(\F_P\)-linear map with zero kernel: a nonzero kernel element would give \(\tau(z)/z=w^2\), whose norm forces \(1=\lambda^2\). Thus each \(x\ne0\) gives exactly one \(z\), for a total of \(P^k-1\) normalized cyclic vectors.
\end{proof}

\begin{corollary}[Repeated-root fibers from the numbers of spectral types]\label{v146:cor-repeated-fibers}
Let \(\gcd(r,k)=1\). If \(k\) is odd, then for every \(\lambda\in\F_P\) with \(\lambda^2=1\),
\begin{equation}\label{v146:eq-repeated-odd}
 \#\{a:\Theta_{k,r}(a)=-2\lambda\}
 =\frac{P^{k+1}-1}{P^2-1}.
\end{equation}
If \(k\) is even, put \(h=(-P)^{k/2}\). In characteristic two the unique repeated-root target is zero, and
\begin{equation}\label{v146:eq-repeated-even-char}
 \#\{a:\Theta_{k,r}(a)=0\}
 =\frac{P^{k+1}-1+P(P-1)h}{P^2-1}.
\end{equation}
In odd characteristic the sum of the two repeated-root fiber sizes is
\begin{equation}\label{v146:eq-repeated-total}
 \#\Theta_{k,r}^{-1}(2)+\#\Theta_{k,r}^{-1}(-2)
 =\frac{2(P^{k+1}-1)+(P-1)^2h}{P^2-1}.
\end{equation}
\end{corollary}
\begin{proof}
For odd \(P\), the numbers of irreducible, split distinct-root, and repeated-root targets are \((P-1)/2,(P-3)/2,2\); for even \(P\), they are \(P/2,(P-2)/2,1\). These counts follow by inverse pairing in \(\mu_{P+1}\) and \(\F_P^\times\). Subtracting the total sizes of the first two types from \(P^k\) gives the third.

If both \(k\) and \(P\) are odd, put \(J=\operatorname{diag}(1,-1)\). The identity \(A_{-a}=-JA_aJ^{-1}\) gives \(\Theta_{k,r}(-a)=-\Theta_{k,r}(a)\), so the two repeated-root fibers have equal size. Substitution from Theorems~\ref{v145:thm-return-uniform} and~\ref{v146:thm-split-fibers}, followed by simplification, gives the formulas.
\end{proof}

For odd \(k\), all three spectral types now have determined fiber sizes. For even \(k\) in odd characteristic, equation~\eqref{v146:eq-repeated-total} gives only their sum; separating the two targets requires finer information than spectral type.

\begin{proposition}[Semilinear scaling of the constant term]\label{v146:prop-return-scaling}
Suppose \(c=u\tau(u)\), with \(u\in K^\times\). Put \(a'=a/\tau(u)\) and \(d=\Nm_{K/\F_P}(u)\). Then
\begin{equation}\label{v146:eq-return-scaling}
\begin{aligned}
 B_{k,r}(a,c)&=d\Theta_{k,r}(a'),\\
 \det(ZI_2-M_{k,r}(a,c))&=Z^2+d\Theta_{k,r}(a')Z+d^2.
\end{aligned}
\end{equation}
Such a \(u\) exists if and only if
\begin{equation}\label{v146:eq-scaling-image}
 c\in(K^\times)^g,\qquad
 g=\begin{cases}P+1,&k\text{ even},\\\gcd(P+1,2),&k\text{ odd}.\end{cases}
\end{equation}
For a normalizable \(c\), the fiber over \(b\) therefore has the same size as the fiber of \(\Theta_{k,r}\) over \(b/d\).
\end{proposition}
\begin{proof}
Put \(T=\operatorname{diag}(1,u^{-1})\). Direct multiplication gives
\[
 A_{a,c}=uT A_{a',1}\tau(T)^{-1},\qquad
 M_{k,r}(a,c)=dT M_{k,r}(a',1)T^{-1}.
\]
Adjacent factors cancel during iteration, giving the second identity and hence~\eqref{v146:eq-return-scaling}. On the cyclic group \(K^\times\), the image of \(u\mapsto u^{P^r+1}\) is \((K^\times)^g\), where \(g=\gcd(P^r+1,P^k-1)\). Since \(\gcd(r,k)=1\), this greatest common divisor has the value in~\eqref{v146:eq-scaling-image}. Finally, \(a\mapsto a'\) is bijective.
\end{proof}

\begin{proposition}[All constant-term fibers in a quadratic extension]\label{v146:prop-constant-quadratic}
Let \(K=\F_{P^2}\), \(\gcd(r,2)=1\), and fix \(c\in K\). Then
\[
 B_{2,r}(a,c)=\Tr_{K/\F_P}(c)-\Nm_{K/\F_P}(a),
\]
and for \(b\in\F_P\),
\begin{equation}\label{v146:eq-constant-quadratic}
 \#\{a:B_{2,r}(a,c)=b\}
 =\begin{cases}1,&b=\Tr(c),\\P+1,&b\ne\Tr(c).\end{cases}
\end{equation}
\end{proposition}
\begin{proof}
Here \(\tau(x)=x^P\), and
\[
 A_{a,c}A_{a^P,c^P}
 =\begin{pmatrix}-c&ca^P\\-a&-c^P+aa^P\end{pmatrix}.
\]
The norm has only zero over zero and exactly \(P+1\) preimages of each nonzero value.
\end{proof}

\begin{example}[A constant term that cannot be normalized]\label{v146:ex-constant-obstruction}
First take \(K=\F_4\), \(P=2\), and \(c=\omega\), where \(\omega^2+\omega+1=0\). The irreducible target \(Z^2+Z+1\) has only one preimage, \(a=0\), whereas the fiber size for constant term one is \(D_2(2)=3\).

Next use \(\F_4\) as the fixed field, take \(K=\F_{16}\), and let \(c\) be a root of \(Z^2+\omega Z+1\). Then \(\Tr_{K/\F_4}(c)=\omega\) and \(\Nm_{K/\F_4}(c)=1\). For this fixed \(c\), the two irreducible targets
\[
 Z^2+\omega Z+1,\qquad Z^2+\omega^2Z+1
\]
have fiber sizes one and five, respectively. Their irreducibility follows from the product in Example~\ref{v145:ex-sparse-four}. Here \(c\) has order five, whereas \((K^\times)^5=\F_4^\times\), so condition~\eqref{v146:eq-scaling-image} fails.
\end{example}

For \(c=1,k=2\), the target \(b=2\) has one preimage, and every other target has \(P+1\). Thus the two repeated-root targets \(2,-2\) in odd characteristic need not have equal fibers.

If \(c=0\), direct multiplication gives
\begin{equation}\label{v146:eq-singular-return}
 \det(ZI_2-M_{k,r}(a,0))
 =Z^2+(-1)^{k+1}\Nm_{K/\F_P}(a)Z.
\end{equation}
The corresponding coefficient map has one preimage of zero and \((P^k-1)/(P-1)\) preimages of every nonzero value. Invertible normalization does not apply, but the norm already determines every fiber.

\subsection{Cyclic Matchings, Trace--Norm Distributions, and Exact Orders}
\label{v146:sec-return-applications}
The return matrix product can also be expanded by pairing adjacent factors. This makes one fiber theorem apply to a family of different trace--norm expressions.

\begin{proposition}[Cyclic matching expansion of the return coefficient]\label{v146:prop-return-matching}
Let \(k\ge3\), \(\gcd(r,k)=1\), and \(a_i=a^{P^{ri}}\). Let \(C_k\) be the cycle on vertices \(0,\ldots,k-1\). A \textbf{matching} is a set of edges with pairwise disjoint endpoints; write \(\mathscr M(C_k)\) for all matchings and \(V(M)\) for the covered vertices. Then
\begin{equation}\label{v146:eq-return-matching}
 \Theta_{k,r}(a)=(-1)^{k+1}
 \sum_{M\in\mathscr M(C_k)}(-1)^{|M|}
 \prod_{i\notin V(M)}a_i.
\end{equation}
The empty matching and empty products are included with their usual conventions.
\end{proposition}
\begin{proof}
Expand \(\operatorname{tr}(A_{a_0}\cdots A_{a_{k-1}})\) as a sum over matrix indices with equal initial and final indices. Since the upper-left entry is zero, every visit to the first index must be followed immediately by departure. These two off-diagonal transitions occupy adjacent factors and contribute \(-1\). The pairs are disjoint, so they form a matching of the cycle. Every remaining factor contributes its lower-right entry \(-a_i\). For \(|M|\) edges the sign in the negative trace is \((-1)^{1+|M|+k-2|M|}=(-1)^{k+1+|M|}\), proving the formula.
\end{proof}

For standard Frobenius, let each trace and norm below be from \(\F_{P^k}\) to \(\F_P\). Then
\begin{equation}\label{v146:eq-small-return-polynomials}
\begin{aligned}
 \Theta_{1,1}(a)&=a,&\Theta_{2,1}(a)&=2-\Nm(a),\\
 \Theta_{3,1}(a)&=\Nm(a)-\Tr(a),\\
 \Theta_{4,1}(a)&=\Tr(a^{P+1})-\Nm(a)-2,\\
 \Theta_{5,1}(a)&=\Nm(a)-\Tr(a^{1+P+P^2})+\Tr(a).
\end{aligned}
\end{equation}
The first two identities follow by multiplication; the others group matchings by their number of edges. The constant \(2\) is interpreted in the field.

\begin{corollary}[Trace--norm fibers in a cubic extension]\label{v145:cor-cubic-trace-norm}
If $Z^2+bZ+1$ is irreducible over $\F_P$, then
\begin{equation}\label{v145:eq-cubic-trace-norm}
 \#\{a\in\F_{P^3}:\Nm_{\F_{P^3}/\F_P}(a)
                      -\Tr_{\F_{P^3}/\F_P}(a)=b\}
 =P^2-P+1.
\end{equation}
\end{corollary}
\begin{proof}
Equation~\eqref{v146:eq-small-return-polynomials} gives $b_3(a)=\Nm(a)-\Tr(a)$. Substitute $D_3(P)=(P^3+1)/(P+1)$.
\end{proof}

\begin{example}[Cubic fibers in even and odd characteristic]\label{v145:ex-cubic-fibers}
Over $\F_8/\F_2$, the polynomial $Z^2+Z+1$ is irreducible, so $\Nm(a)-\Tr(a)=1$ has three solutions. Every nonzero element has norm one, so these are exactly the nonzero trace-zero elements. The remaining five elements lie over zero.

Over $\F_{27}/\F_3$, the polynomial $Z^2+1$ is irreducible, so zero has seven preimages. The map $a\mapsto-a$ negates $\Nm(a)-\Tr(a)$, hence one and two each have ten preimages. These different fiber sizes belong to different quadratic spectral types.
\end{example}

\begin{corollary}[A mixed distribution in a quartic extension]\label{v146:cor-quartic-distribution}
If \(Z^2+bZ+1\) is irreducible over \(\F_P\), then
\[
 \#\{a\in\F_{P^4}:\Tr(a^{P+1})-\Nm(a)-2=b\}
 =(P-1)^2(P+1).
\]
\end{corollary}
\begin{proof}
Apply~\eqref{v146:eq-small-return-polynomials} and \(D_4(P)=(P^2-1)^2/(P+1)\).
\end{proof}

The degree-five expression likewise has \(D_5(P)=P^4-P^3+P^2-P+1\) preimages of each irreducible target. For other targets, Theorem~\ref{v146:thm-split-fibers} and Corollary~\ref{v146:cor-repeated-fibers} give the corresponding spectral distributions.

\begin{corollary}[Coefficient counts for return orders and Lucas atoms]\label{v146:cor-return-order}
Let \(\gcd(r,k)=1\), \(s>2\), and \(s\mid P+1\). Then
\begin{equation}\label{v146:eq-return-order}
 \#\{a\in\F_{P^k}:H_{P,s}(\Theta_{k,r}(a))=0\}
 =\frac{\varphi(s)}2D_k(P).
\end{equation}
For these \(a\), every nonzero vector under the return operator \(M_{k,r}(a,1)\) has period \(s\). In particular, if \(c\) is primitive in \(\F_P\) and \(t_{P,c}\) counts primitive quadratics with constant term \(c\), then
\begin{equation}\label{v146:eq-primitive-return-comparison}
 \#\{a:H_P(\Theta_{k,r}(a))=0\}
 =\begin{cases}
 \dfrac{D_k(P)}2t_{P,c},&P\text{ odd},\\[1mm]
 D_k(P)t_{P,c},&P\text{ even}.
 \end{cases}
\end{equation}
\end{corollary}
\begin{proof}
The \(\varphi(s)\) elements of order \(s\) form \(\varphi(s)/2\) inverse pairs. Each gives a root of \(H_{P,s}\), hence an irreducible return target. Summing Theorem~\ref{v145:thm-return-uniform} proves~\eqref{v146:eq-return-order}. The return operator diagonalizes over a splitting field, with both eigenvalues of order \(s\). Every nonzero vector has a nonzero eigen-coordinate, so its period is exactly \(s\). Finally, set \(s=P+1\) and apply~\eqref{v45:eq-fixed-constant-primitive-count}.
\end{proof}

This period belongs to the linear return \(\Phi_{a,1}^k\), not to the original semilinear transformation \(\Phi_{a,1}\). Equation~\eqref{v146:eq-lucas-pullback} explains how the same \(H_{P,s}\) converts a return-eigenvalue order into the root-ratio order of an ordinary quadratic. The polynomial \(H_4\) in Example~\ref{v146:ex-lucas-four} has two roots. For \(k=2\), each has five return coefficients, giving the ten sparse skew-irreducible coefficients of Example~\ref{v145:ex-sparse-four}.

\begin{theorem}[Counting sparse skew polynomials in divisible extension degrees]\label{v145:thm-sparse-count}
Let $q$ be a prime power, $n,k\ge1$, $K=\F_{q^{nk}}$, and $\sigma(c)=c^q$. Write
\[
 \mathcal N_n(q,nk)=\#\{a\in K:
 t^{2n}+at^n+1\text{ is irreducible in }K[t;\sigma]\}.
\]
Then
\begin{equation}\label{v145:eq-sparse-count}
 \mathcal N_n(q,nk)=C_n(q)D_k(q^n)
                  =nS_q(2n)D_k(q^n).
\end{equation}
\end{theorem}
\begin{proof}
Put $P=q^n$. The corresponding $q$-linearized polynomial is
\[
 L_a(X)=X^{q^{2n}}+aX^{q^n}+X=X^{P^2}+aX^P+X.
\]
Its root space has dimension two over $\F_P$ and dimension $2n$ over $\F_q$. The map $\Theta(x)=x^{q^{nk}}=x^{P^k}$ is $\F_P$-linear, with characteristic polynomial $G_a(Z)=Z^2+b_k(a)Z+1$ over $\F_P$. Restricting scalars to $\F_q$, conjugate coordinates give
\begin{equation}\label{v145:eq-restricted-return}
 \chi_{\Theta,\F_q}(Z)
 =\prod_{i=0}^{n-1}\bigl(Z^2+b_k(a)^{q^i}Z+1\bigr).
\end{equation}
Indeed, after extension to the algebraic closure, $\F_{q^n}\otimes_{\F_q}\overline{\F}_q$ splits according to its $n$ embeddings. The corresponding two-dimensional matrices are the entrywise conjugates of the original matrix, so their determinants multiply.

By Theorem~\ref{v145:thm-skew-roots}, the skew polynomial is irreducible exactly when~\eqref{v145:eq-restricted-return} is irreducible. If $G_a$ is irreducible and $b_k(a)$ has degree $n$ over $\F_q$, the root-degree argument in Proposition~\ref{v145:prop-inverse-count} gives degree $2n$ for a root. The product is therefore its minimal polynomial. Conversely, if $G_a$ is reducible, its roots lie in $\F_{q^n}$ and have degree at most $n$; if $b_k(a)$ has degree $d<n$, the product is the $(n/d)$th power of a degree-$2d$ polynomial. Neither case gives an irreducible polynomial of degree $2n$.

There are $C_n(q)$ allowable targets $b$, each with $D_k(q^n)$ preimages, proving~\eqref{v145:eq-sparse-count}.
\end{proof}

If $n$ is a power of two, its only odd divisor is one, so
\[
 \mathcal N_n(q,nk)=\frac{q^n-\epsilon_q}{2}D_k(q^n).
\]
For $k=1$, equation~\eqref{v145:eq-sparse-count} becomes $\mathcal N_n(q,n)=nS_q(2n)$.

\begin{example}[Quartic skew polynomials and fifth roots of unity]\label{v145:ex-sparse-four}
Take $q=2,n=2$, so $C_2(2)=2$. For $k=1,2,3$, the values of $D_k(4)$ are $1,5,13$, giving $2,10,26$ irreducible coefficients, respectively.
For $k=1$, write $\F_4=\F_2(\omega)$ with $\omega^2+\omega+1=0$. The allowable coefficients are exactly $a=\omega,\omega^2$, since
\[
 (Z^2+\omega Z+1)(Z^2+\omega^2Z+1)
 =Z^4+Z^3+Z^2+Z+1
\]
is irreducible over $\F_2$, with primitive fifth roots of unity as roots.
These two polynomials $t^4+at^2+1$ are irreducible in the skew polynomial ring; as ordinary polynomials, they equal $(t^2+a^2t+1)^2$. The multiplication counted in~\eqref{v145:eq-sparse-count} is determined by $tc=c^2t$.
\end{example}

If the coefficient field is $\F_{q^h}$ with $n\nmid h$, the $q^h$-power map on the two-dimensional $\F_{q^n}$-root space is generally only semilinear. The root-space criterion still applies, but the restriction-of-scalars step in~\eqref{v145:eq-restricted-return} must change. Proposition~\ref{v148:prop-sparse-transfer} uses the common fixed field and norm transfer; Corollary~\ref{v148:cor-trinomial-spacing} gives the count for arbitrary spacing.

\subsection{Unequal Fibers of Four-Dimensional Return Spectra}
\label{v146:sec-higher-return}
Irreducible quadratic return targets have equal fibers. In larger root spaces, a return characteristic polynomial still detects skew irreducibility, but irreducible targets of the same degree need not have equal numbers of coefficient preimages.

\begin{example}[Nonuniform irreducible return spectra in dimension four]\label{v146:ex-quartic-counterexample}
Let \(K=\F_4=\F_2(\omega)\), \(\omega^2+\omega+1=0\), and \(\tau(z)=z^2\). For
\[
 f=t^4+a_3t^3+a_2t^2+a_1t+1\in K[t;\tau],
 \qquad
 A_f=\begin{pmatrix}
 0&0&0&1\\1&0&0&a_1\\0&1&0&a_2\\0&0&1&a_3
 \end{pmatrix},
\]
the return matrix is \(M_f=A_f\tau(A_f)\). The three monic irreducible quartics over \(\F_2\) have the following fibers:
\begin{center}
\begin{tabular}{lc}
\toprule
Return characteristic polynomial & Number of \((a_1,a_2,a_3)\)\\\midrule
\(Z^4+Z+1\)&3\\
\(Z^4+Z^3+1\)&3\\
\(Z^4+Z^3+Z^2+Z+1\)&11\\\bottomrule
\end{tabular}
\end{center}
Put \(T(x)=x+x^2\) and \(N(x)=x^3\). Multiplying the matrices and expanding the determinant gives
\begin{equation}\label{v146:eq-quartic-return-coefficients}
\begin{split}
 \det(ZI_4-M_f)&=Z^4+c_3Z^3+c_2Z^2+c_1Z+1,\\
 c_1&=N(a_1)+T(a_2),\\
 c_2&=N(a_2)+T(a_1^2a_3),\\
 c_3&=N(a_3)+T(a_2).
\end{split}
\end{equation}
The joint values \((T,N)\) are \((0,0)\) at zero, \((0,1)\) at one, and \((1,1)\) at both \(\omega,\omega^2\).

For \((c_3,c_2,c_1)=(0,0,1)\), if \(T(a_2)=0\), then \(a_1\ne0,a_3=0\), and \(c_2=0\) forces \(a_2=0\), giving three triples. If \(T(a_2)=1\), then \(a_1=0,a_3\ne0\), but \(c_2=N(a_2)=1\), which is impossible. Interchanging \(a_1,a_3\) gives three triples for \((1,0,0)\) as well; here \(T(a_1^2a_3)=T(a_3^2a_1)\).

For \((1,1,1)\), if \(T(a_2)=1\), then \(a_1=a_3=0\), with two choices for \(a_2\). If \(a_2=0\), then \(a_1,a_3\ne0\) and \(T(a_1^2a_3)=1\), giving \(3\cdot2=6\) triples. If \(a_2=1\), the condition becomes \(T(a_1^2a_3)=0\), giving three triples. The total is therefore \(2+6+3=11\).
\end{example}

The three targets have the same degree, determinant, and isomorphism type of general linear centralizer, but their eigenvalue orders are fifteen, fifteen, and five. Degree or centralizer type alone therefore cannot recover the fiber size. Whether prescribing eigenvalue order restores uniformity also depends on the conditions that the coefficient family imposes on normalized cyclic bases.

\section{Higher-degree return fibers over a common subfield}
\label{v148:sec-higher-fibers}
The three fiber sizes in Example~\ref{v146:ex-quartic-counterexample} are unequal, whereas Theorem~\ref{v148:thm-relative-fibers} guarantees uniformity when the degree and coefficient-extension degree are coprime. We now fix a quadratic coefficient extension. Odd degrees are already settled by coprimality; in even degrees, a commutative norm equation describes the entire deviation.

Fix $K=\F_{q^2}$, $\bar a=a^q$, $R=K[t;\bar{\phantom a}]$, and $Z=t^2$. For an ordinary polynomial $P\in K[Z]$, write $N_2(P)=P\bar P$, with conjugation acting only on coefficients. Every skew polynomial has a unique expression $P(Z)+Q(Z)t$. The representation
\[
 a\longmapsto\begin{pmatrix}a&0\\0&\bar a\end{pmatrix},\qquad
 t\longmapsto\begin{pmatrix}0&Z\\1&0\end{pmatrix}
\]
gives determinant $N_2(P)-ZN_2(Q)$. This is monic in even skew degree; in odd degree it must be multiplied by $-1$.

\begin{proposition}[Factor allocation in a commutative norm]\label{v148:prop-commutative-norm}
Let $V(H)=\#\{Q\in K[Z]:Q\bar Q=H\}$. Then $V(0)=1$. If $H=a\prod_\pi\pi^{e_\pi}\ne0$, with $a\in\F_q^\times$, then
\[
 V(H)=
 \begin{cases}
 (q+1)\displaystyle\prod_{2\mid\deg\pi}(e_\pi+1),
   &e_\pi\text{ is even for every odd-degree }\pi,\\
 0,&\text{otherwise}.
 \end{cases}
\]
For a monic target $h$ of degree $2m$, if $h(0)=\Nm(c)$, then
\begin{equation}\label{v148:eq-direct-even-fiber}
 F_{h,2}(c)=
 \sum_{\substack{P\text{ monic},\ \deg P=m\\P(0)=c}}
 V\!\left(\frac{P\bar P-h}{Z}\right).
\end{equation}
For a monic target of degree $2m+1$, if $h(0)=-\Nm(c)$, then
\[
 F_{h,2}(c)=\frac1{q+1}
 \sum_{\substack{\deg P\le m\\P(0)=c}}
 V\!\left(\frac{h+P\bar P}{Z}\right).
\]
Incompatible constant coefficients give zero.
\end{proposition}
\begin{proof}
An odd-degree irreducible factor remains irreducible over $K$ and is fixed by conjugation, so its exponent must be even. An even-degree factor splits as $g\bar g$; its two exponents in $Q$ sum to $e_\pi$, giving $e_\pi+1$ choices. The norm equation for the leading coefficient has $q+1$ solutions. The last two formulas follow from $P\bar P-ZQ\bar Q$; in odd degree, $Q$ is monic, removing the $q+1$ leading-coefficient choices. The quotient degree automatically imposes the required bound on $\deg Q$.
\end{proof}

The direct formula retains every factor choice. To compare different constant coefficients, most of this sum can be eliminated. For monic $h$ of degree $2m$, define
\[
 E_h(c)=\#\{P\in K[Z]\text{ monic}:\deg P=m,\ P(0)=c,\ P\bar P=h\},
 \qquad E(h)=\sum_cE_h(c).
\]

\subsection{Hankel kernel shifts and singular strata}
Fix a nontrivial additive character $\psi$. For $\lambda=(\lambda_0,\ldots,\lambda_{2m-1})\in\F_q^{2m}$, set $\Lambda(Z^j)=\lambda_j$ and $\Lambda(Z^{2m})=0$, and extend $\Lambda$ $K$-linearly. Put
\[
 J_\lambda=(\lambda_{i+j+1})_{0\le i,j<m},\qquad
 A_\lambda(c)=\sum_{\substack{P\text{ monic},\ \deg P=m\\P(0)=c}}
                 \psi\bigl(\Lambda(P\bar P)\bigr).
\]
The entries of the Hankel matrix $J_\lambda$ are constant along antidiagonals.

\begin{lemma}[Norm invariance on a singular Hankel stratum]\label{v148:lem-hankel-shift}
If $J_\lambda$ is singular, then $A_\lambda(c)$ depends only on $\Nm(c)$. More precisely, either $A_\lambda(c)=0$ for every $c$, or there is a monic $S\in\F_q[Z]$ with $\deg S=m$, $S(0)=0$, and
\[
 \Lambda(SV)=0\qquad(\deg V<m).
\]
\end{lemma}
\begin{proof}
Write a nonzero kernel vector as $r\in\F_q[Z]$, $\deg r<m$. Its kernel condition is
\[
 \Lambda(ZrV)=0\qquad(\deg V<m).
\]
If $\deg r=m-1$, monic normalization of $Zr$ gives $S$. Otherwise put $R_0=Zr$ and $\delta=\Lambda(R_0Z^m)$. For every summation polynomial $P$,
$\Lambda(R_0\bar P)=\delta$ and $\Lambda(R_0^2)=0$. Translation $P\mapsto P+aR_0$ preserves degree, monicity, and constant coefficient, and adds $\Tr_{K/\F_q}(a\delta)$ to the phase. If $\delta\ne0$, trace surjectivity gives a translation multiplying the sum by a character value other than one; the sum is therefore zero.

If $\delta=0$, then $Zr$ is again a kernel polynomial. For $\deg V\le m-2$, apply the original condition to $ZV$; the endpoint $V=Z^{m-1}$ is supplied by $\delta=0$. Thus the kernel degree increases by one. Repetition either gives the zero-sum case or produces $S$.

In the latter case, write $P=S+V$. The cross terms vanish, so
\[
 \Lambda(P\bar P)=\Lambda(S^2)+\Lambda(V\bar V).
\]
For $\Nm(\zeta)=1$, the map $V\mapsto\zeta V$ preserves this expression and changes $c$ to $\zeta c$. This proves norm invariance, including in characteristic two.
\end{proof}

\begin{theorem}[Even-degree fibers over a quadratic coefficient field]\label{v148:thm-all-even-fibers}
Let $h\in\F_q[Z]$ be monic, with $\deg h=2m\ge2$ and $h(0)\ne0$. For $\Nm(c)=h(0)$,
\begin{equation}\label{v148:eq-all-even-fibers}
 F_{h,2}(c)=\frac{U_2(h)-(-q)^mE(h)}{q+1}+(-q)^mE_h(c).
\end{equation}
The irreducible factors of $h$ may have arbitrary multiplicities. In particular, $F_{h,2}(c)$ is uniform on compatible constants if and only if $E_h(c)$ is uniform there.
\end{theorem}
\begin{proof}
Identify monic degree-$2m$ targets with their lower coefficient space $\F_q^{2m}$ and take the unnormalized Fourier transform in the target. Since $\mathcal N(f)=P\bar P-ZQ\bar Q$,
\[
 \widehat F_c(\lambda)=A_\lambda(c)G_\lambda,
 \qquad \widehat E_c(\lambda)=A_\lambda(c),\qquad
 G_\lambda=\sum_{\deg Q<m}\psi\bigl(-\Lambda(ZQ\bar Q)\bigr).
\]
For $r=\rank J_\lambda$, diagonalization of the Hermitian form gives
\[
 G_\lambda=(-1)^rq^{2m-r}.
\]
Indeed, each nonzero diagonal direction contributes
$1+(q+1)\sum_{u\ne0}\psi(au)=-q$, and each zero direction contributes $q^2$.

Put $a_m=(-q)^m$. On the full-rank stratum, $G_\lambda=a_m$, so
\[
 \widehat F_c(\lambda)-a_m\widehat E_c(\lambda)
 =A_\lambda(c)(G_\lambda-a_m)
\]
is zero. On singular strata, Lemma~\ref{v148:lem-hankel-shift} shows that it depends only on $\Nm(c)$. Fourier inversion therefore makes $F_{h,2}(c)-a_mE_h(c)$ constant on each norm coset. Summing over its $q+1$ elements gives $U_2(h)-a_mE(h)$ and proves the formula.
\end{proof}

The singular strata need not be evaluated separately: a kernel shift first removes their dependence on the constant-coefficient phase, while the full-rank stratum gives the surviving correction. Hermitian rank and factor constants thus play different roles in the count.

\begin{corollary}[Three-value formula for irreducible targets]\label{v148:cor-three-values}
If $h$ is irreducible, write $h=g\bar g$, where $g\in K[Z]$ is monic of degree $m$, and let $b=g(0)$. Put
\[
 B_m(q)=\frac{q^{2m}+1-2(-q)^m}{q+1}.
\]
For $\Nm(c)=h(0)$,
\begin{equation}\label{v148:eq-three-values}
 F_{h,2}(c)=B_m(q)+(-q)^m
       \bigl(\mathbf1_{c=b}+\mathbf1_{c=b^q}\bigr).
\end{equation}
\end{corollary}
\begin{proof}
Here $U_2(h)=q^{2m}+1$, and the monic solutions of $P\bar P=h$ are exactly $g,\bar g$. Hence $E(h)=2$, and the sum of indicators is $E_h(c)$.
\end{proof}

The three possible numbers are $B_m$, $B_m+(-q)^m$, and $B_m+2(-q)^m$; each fixed target realizes at most two of them. If $b\in\F_q$, the corrections coincide at one constant; otherwise they occur at two conjugate constants.

\begin{corollary}[All targets over a quadratic coefficient extension]\label{v148:cor-quadratic-classification}
Let $\deg h=n\ge1$. If $h(0)=0$, only $c=0$ is possible, and $F_{h,2}(0)=U_2(h)$. If $h(0)\ne0$, incompatible constants give zero; compatible fibers have size $U_2(h)/(q+1)$ for odd $n$, and are given by~\eqref{v148:eq-all-even-fibers} for even $n$.

For monic $h=\prod\pi^{e_\pi}$ with $h(0)\ne0$,
\[
 E(h)=
 \begin{cases}
 \displaystyle\prod_{2\mid\deg\pi}(e_\pi+1),
   &\text{every odd-degree factor has even exponent},\\
 0,&\text{otherwise}.
 \end{cases}
\]
The number $E_h(c)$ is obtained by retaining the factor allocations whose constant coefficient is $c$.
\end{corollary}
\begin{proof}
A zero constant coefficient is equivalent to a zero constant coefficient of the return norm. Use Theorem~\ref{v148:thm-relative-fibers} in odd degree and the preceding theorem in even degree. The product is the monic case of Proposition~\ref{v148:prop-commutative-norm}.
\end{proof}

\begin{proposition}[Factor phases and a divisibility criterion for uniform fibers]
\label{v149:prop-uniformity}
Let $h\in\F_q[Z]$ be monic, with irreducible factorization $h=\prod_\pi\pi^{e_\pi}$, $\deg h=2m\ge2$, and $h(0)\ne0$. If an odd-degree factor has odd exponent, the fibers over compatible constant coefficients are uniform.
Otherwise, for each even-degree factor choose $\pi=g_\pi\bar g_\pi$ with $g_\pi\in K[Z]$ monic, and put
\[
 P_* =\prod_{\deg\pi\text{ odd}}\pi^{e_\pi/2}
       \prod_{\deg\pi\text{ even}}g_\pi^{e_\pi},
 \qquad c_*=P_*(0).
\]
Choose a generator $\theta$ of the norm-one group $U=\ker\Nm_{K/\F_q}$, write $N=q+1$, and define
\[
 \rho_\pi=\frac{g_\pi(0)^q}{g_\pi(0)}=\theta^{r_\pi},
 \qquad r_\pi\in\mathbb Z/N\mathbb Z.
\]
Then, in $\mathbb Z[T]/(T^N-1)$,
\begin{equation}\label{v149:eq-phase-distribution}
 \sum_{a=0}^{N-1}E_h(c_*\theta^a)T^a
 =\prod_{2\mid\deg\pi}
       \left(1+T^{r_\pi}+\cdots+T^{e_\pi r_\pi}\right).
\end{equation}
The fibers over compatible constant coefficients are uniform if and only if, for every $d\mid N$ with $d>1$, some even-degree factor $\pi$ satisfies
\begin{equation}\label{v149:eq-uniform-divisibility}
 d\nmid r_\pi,\qquad d\mid(e_\pi+1)r_\pi.
\end{equation}
These conditions are independent of the integer representatives of $r_\pi$, the generator $\theta$, and the choices of $g_\pi$.
\end{proposition}
\begin{proof}
In the first case $E_h(c)=0$, so Theorem~\ref{v148:thm-all-even-fibers} applies directly.
Otherwise, every monic solution of $P\bar P=h$ comes from independent choices
$g_\pi^{e_\pi-j}\bar g_\pi^j$, $0\le j\le e_\pi$, at the even-degree factors. Relative to $P_*$, its constant coefficient is multiplied by $\prod_\pi\rho_\pi^j$. This proves~\eqref{v149:eq-phase-distribution}.

Cyclic Fourier inversion over $\mathbb C$ shows that all coefficients on the left are equal exactly when the right side vanishes at every nontrivial $N$th root of unity. At a root $\xi$ of exact order $d$, a single factor takes the value
\[
 \sum_{j=0}^{e_\pi}\xi^{jr_\pi}.
\]
If $d\mid r_\pi$, this is the positive integer $e_\pi+1$. Otherwise the geometric-sum formula makes it zero exactly when $d\mid(e_\pi+1)r_\pi$. The product vanishes exactly when at least one factor does, yielding~\eqref{v149:eq-uniform-divisibility}. The even-degree fiber formula then transfers uniformity of $E_h$ to uniformity of $F_{h,2}$.

Conjugating $g_\pi$ negates $r_\pi$, and changing the generator multiplies all exponents by a unit modulo $N$. Both operations preserve the divisibility conditions.
\end{proof}

Thus factor degrees and multiplicities determine $E(h)$, while their constant-coefficient phases determine its distribution. Uniformity requires $q+1\mid E(h)$, but this total-count condition is generally insufficient. If $h=\pi^e$ with $\deg\pi$ even, uniformity holds exactly when $\rho_\pi$ generates $U$ and $q+1\mid e+1$. In particular, an irreducible target of even degree is never uniform, since then $E(h)=2<q+1$.

Suppose $\deg h=4$ and $h(0)\ne0$, and let $\rho(h)$ count the distinct roots of $h$ in $\F_q^\times$. For compatible constants, the formula also becomes
\begin{equation}\label{v148:eq-quartic-all-targets}
 F_{h,2}(c)=(q-1)^2(q+1)+q(q+1)\rho(h)+q^2E_h(c).
\end{equation}
Indeed, substituting the possible degree-four factor degrees and multiplicities into the two product formulas gives
$U_2(h)-q^2E(h)=(q+1)^2((q-1)^2+q\rho(h))$; apply the even-degree fiber formula. Thus the baseline for a general quartic target also records its number of distinct base-field roots.

\begin{example}[Repeated factors and restored uniformity]\label{v148:ex-repeated-fibers}
Take $q=2$ and $K=\F_4=\F_2(\omega)$. For $h=(Z^2+Z+1)^2$, we have $U_2(h)=21$. The three monic commutative norm preimages are
$(Z+\omega)^2$, $(Z+\omega)(Z+\omega^2)$, and $(Z+\omega^2)^2$.
Their constants are $\omega^2,1,\omega$, so all three fibers have size $7$. Here $E(h)\ne0$, yet the correction remains uniform.
Here $\rho_\pi$ has order three and $e_\pi+1=3$, exactly as required by Proposition~\ref{v149:prop-uniformity}.

In contrast, take $h=(Z^4+Z^3+Z^2+Z+1)^2$. Its irreducible quartic factor splits over $K$ as
\[
 (Z^2+\omega Z+1)(Z^2+\omega^2Z+1).
\]
Both constant coefficients are one. Thus $E(h)=3$, but all three commutative norm preimages have constant coefficient one. Since $U_2(h)=1+16+256=273$ and $(-2)^4=16$, the fibers for constants $1,\omega,\omega^2$ have sizes $123,75,75$. Hence $q+1\mid E(h)$ alone does not imply uniformity.

For $h=(Z+1)^6$, we have $U_2(h)=127$, and the sole commutative norm preimage is $(Z+1)^3$. Since $(-2)^3=-8$, the fibers for constants $1,\omega,\omega^2$ have sizes
\[
 37,\qquad45,\qquad45.
\]
The quartic irreducible distribution $3,3,11$ in Example~\ref{v146:ex-quartic-counterexample} is likewise determined by the constants of the two conjugate quadratic factors.
\end{example}

\subsection{Return-order spectra and character sums}
\begin{theorem}[Complete return-order spectrum in even degree]\label{v148:thm-even-order}
Put $Q=q^{2m}$. Let $\ell\mid Q-1$, $\ord_\ell(q)=2m$, and
\[
 d=\frac{\ell}{\gcd(\ell,(Q-1)/(q-1))},\qquad
 e=\frac{\ell}{\gcd(\ell,(Q-1)/(q^2-1))}.
\]
If $\ord(\Nm(c))=d$, the number of skew irreducible polynomials of degree $2m$, constant coefficient $c$, and return-root order $\ell$ is
\begin{equation}\label{v148:eq-even-order}
 J_{2m,\ell}(c)=B_m(q)\frac{\varphi(\ell)}{2m\varphi(d)}
 +(-q)^m\frac{\varphi(\ell)}{m\varphi(e)}
                 \mathbf1_{\ord((-1)^mc)=e}.
\end{equation}
An incompatible norm order gives zero. In particular, if $\Nm(c)$ is primitive, the norm-primitive count is
\begin{equation}\label{v148:eq-even-primitive}
 B_m(q)\frac{\varphi(Q-1)}{2m\varphi(q-1)}
 +(-q)^m\frac{\varphi(Q-1)}{m\varphi(q^2-1)}
                    \mathbf1_{c\text{ primitive in }K};
\end{equation}
otherwise it is zero.
\end{theorem}
\begin{proof}
Roots of order $\ell$ have norm-image orders $d,e$ over $\F_q$ and $K$. Prescribing the base-field norm leaves $\varphi(\ell)/(2m\varphi(d))$ targets. The constant of the minimal polynomial over $K$ of a target root $\alpha$ is
$(-1)^m\Nm_{\F_Q/K}(\alpha)$. Each such factor corresponds to $m$ roots, so the sum of $E_h(c)$ over all targets is
$\varphi(\ell)/(m\varphi(e))$ when $\ord((-1)^mc)=e$, and zero otherwise. This also accounts for coincident constants $b=b^q$. Substitution into the three-value formula gives the first identity.

Set $\ell=Q-1$ for the second. Multiplication by $-1$ preserves primitivity in $K=\F_{q^2}$: in odd characteristic $4\mid q^2-1$, and in even characteristic $-1=1$.
\end{proof}

\begin{example}[Sextic return orders and ordinary factorization]\label{v148:ex-sixth-orders}
Take $q=2$ and $m=3$. Irreducible sextic targets have root orders $9,21,63$, and $B_3(2)=27$. For constants $1,\omega,\omega^2$, formula~\eqref{v148:eq-even-order} gives
\[
 \begin{array}{c|rrr}
 \ell\backslash c&1&\omega&\omega^2\\\hline
 9&27&19&19\\
 21&22&54&54\\
 63&162&114&114
 \end{array}
\]
By Theorem~\ref{v148:thm-ordinary-return}, the corresponding $P_f$ have seven degree-nine factors, three degree-twenty-one factors, or one irreducible degree-sixty-three factor, respectively.
\end{example}

\begin{proposition}[Multiplicative Fourier correction and elliptic traces]\label{v148:prop-fiber-fourier}
Let $\chi$ be a multiplicative character of $K^\times$ nontrivial on the norm-one group. For any monic even-degree $h$ with $h(0)\ne0$,
\begin{equation}\label{v148:eq-fiber-fourier}
 \sum_{\mathcal N(f)=h}\chi(f(0))
 =(-q)^m\sum_{\substack{P\text{ monic}\\P\bar P=h}}\chi(P(0)).
\end{equation}
If $q=2$, $h$ is irreducible, and $\alpha\in\F_{2^{2m}}$ is a root, put
\[
 E_\alpha:\quad Y^2+Y=\alpha X^3,\qquad
 a_{2^{2m}}(\alpha)=2^{2m}+1-\#E_\alpha(\F_{2^{2m}}).
\]
Then
\begin{equation}\label{v148:eq-elliptic-fiber}
 3F_{h,2}(1)-(2^{2m}+1)=2a_{2^{2m}}(\alpha).
\end{equation}
\end{proposition}
\begin{proof}
The sum of $\chi$ on each nonzero norm coset is zero. Multiply~\eqref{v148:eq-all-even-fibers} by $\chi(c)$ and sum; the constant term vanishes.

Over $\F_4$, both nontrivial cubic characters have Gauss sum $2$, as a direct sum over the two trace fibers shows. Hasse--Davenport lifting to $\F_{4^m}$ gives $G=(-1)^{m-1}2^m$. Hence
\[
 a_{2^{2m}}(\alpha)=(-2)^m\epsilon(\alpha),\qquad
 \epsilon(\alpha)=
 \begin{cases}2,&\alpha\text{ is a cube},\\-1,&\alpha\text{ is not a cube}.
 \end{cases}
\]
Indeed, write the point-count trace as $-\sum_x\psi(\alpha x^3)$ and expand with cubic characters. On the other hand, $\Nm_{\F_{4^m}/\F_4}(\alpha)=1$ exactly when $\alpha$ is a cube. Thus $E_h(1)=2$ or $0$ in the corresponding cases; the three-value formula proves the identity.
\end{proof}

The identity relates two counts governed by the same relative-norm phase: fixed-constant skew-polynomial fibers and an elliptic trace obtained from cubic character sums.

\section{Quadratic returns, Lucas atoms, and coefficient criteria}
\label{v148:sec-quadratic-lucas}
For high degree over a quadratic coefficient field, constants of commutative factors govern the fiber formula. Interchanging the two degrees leads to quadratic returns over arbitrary even coefficient extensions. The same common subfield appears, but the phases must be computed anew.

\begin{theorem}[Quadratic fibers over an even coefficient extension]\label{v148:thm-quadratic-even}
Let $K=\F_{q^{2s}}$, $\tau(x)=x^q$, and let $h=Z^2+BZ+D\in\F_q[Z]$ be irreducible, with roots $\alpha,\alpha^q\in\F_{q^2}$. For $c\in K^\times$, put $\rho=\Nm_{K/\F_{q^2}}(c)$. If $\Nm_{K/\F_q}(c)=D$, then
\begin{equation}\label{v148:eq-quadratic-even}
 F_{h,2s}(c)=B_s(q)+(-q)^s
 \bigl(\mathbf1_{\rho=(-1)^s\alpha}
       +\mathbf1_{\rho=(-1)^s\alpha^q}\bigr).
\end{equation}
Otherwise the fiber is zero. The two indicators cannot both be one, and every compatible fiber is positive.
\end{theorem}
\begin{proof}
Put $Q=q^{2s}$. Choose $C\in\operatorname{GL}_2(\F_q)$ with characteristic polynomial $h$, and an $\alpha$-eigenvector $v\in\F_{q^2}^2$. Set
\[
 d_0=\det(v,v^q),\qquad d_0^q=-d_0,\qquad \zeta^{Q-1}=\alpha.
\]
Every solution of $u^Q=Cu$ is uniquely $u=\zeta xv+\zeta^qyv^q$, with $x,y\in K$. Calculation gives
\[
 \Delta(u)=\det(u,u^q)=A x^{q+1}+A^q y^{q+1},\qquad
 A=d_0\zeta^{q+1}.
\]
Choose $\delta^{q-1}=c$. Compatibility implies $A/\delta,A^q/\delta\in K$. The normalized-root-basis count is as in Lemma~\ref{v145:lem-normalized-bases}: $\Delta(u)=\delta$ forces companion determinant $c$, and the $q+1$ determinant-one elements of the centralizer act transitively on the normalized bases of one polynomial. Thus
\[
 F_{h,2s}(c)=\frac1{q+1}
 \#\{(x,y)\in K^2:(A/\delta)x^{q+1}+(A^q/\delta)y^{q+1}=1\}.
\]
For $E=(Q-1)/(q+1)$, the two coefficient phases modulo $(q+1)$st powers are
\[
 U=(A/\delta)^E=\frac{(-1)^s\alpha}{\rho},\qquad
 V=(A^q/\delta)^E=\frac{(-1)^s\alpha^q}{\rho}.
\]
Here $d_0^E=(-1)^s$ and
$\delta^E=c^{1+q^2+\cdots+q^{2s-2}}=\rho$.

For a character $\chi$ of $K^\times$ of order $q+1$, and $1\le i,j\le q$,
\[
 J(\chi^i,\chi^j)=
 \begin{cases}(-1)^{s-1}q^s,&i+j\not\equiv0\pmod{q+1},\\-1,&i+j\equiv0\pmod{q+1}.
 \end{cases}
\]
Indeed, choose a nonzero trace-zero $\eta\in\F_{q^2}$. Summation by trace fibers gives $G(\chi_0^i)=q\chi_0^i(\eta)$. Hasse--Davenport lifting and the Gauss--Jacobi identity cancel the phases; for inverse characters use $\chi(-1)=1$.

Fix an isomorphism from $\mu_{q+1}$ to the complex $(q+1)$st roots of unity, and use $U,V$ also for these complex phases in the character sums. Irreducibility gives $U\ne V$. With $e=\mathbf1_{U=1}+\mathbf1_{V=1}$,
\[
 \sum_{i,j=1}^qU^{-i}V^{-j}=1-(q+1)e,\qquad
 \sum_{i=1}^qU^{-i}V^i=-1.
\]
Character orthogonality now gives
$Q+1-2(-q)^s+(q+1)(-q)^se$ solutions, proving the formula. For even $s$, the baseline is $(q^s-1)^2/(q+1)>0$; for odd $s$, the smaller value is
$(q^{2s}-(q-1)q^s+1)/(q+1)>0$.
\end{proof}

There is no exceptional target when $\rho\in\F_q$. Otherwise exactly the target with root $(-1)^s\rho$ receives a correction. In particular, $c=1$ recovers Theorem~\ref{v145:thm-return-uniform}. The correction in both counts is determined by the relative norm to the quadratic common subfield.

\begin{corollary}[Quadratic primitive fibers and the full order spectrum]\label{v148:cor-quadratic-orders}
Put $d_q=\varphi(q^2-1)/(2\varphi(q-1))$. If $D=\Nm(c)$ is primitive, the number of norm-primitive quadratics is
\[
 B_s(q)d_q+(-q)^s\mathbf1_{\rho\text{ primitive in }\F_{q^2}};
\]
otherwise it is zero. More generally, let $\ell\mid q^2-1$, $\ord_\ell(q)=2$, and $d=\ell/\gcd(\ell,q+1)$. If $\ord(D)=d$, the number of return order $\ell$ is
\[
 B_s(q)\frac{\varphi(\ell)}{2\varphi(d)}
 +(-q)^s\mathbf1_{\ord((-1)^s\rho)=\ell};
\]
it is zero otherwise.
\end{corollary}
\begin{proof}
Sum over quadratic targets with prescribed norm and root order. The sole possible correction target has root $(-1)^s\rho$; the target count follows from counting generators under a cyclic-group norm.
\end{proof}

\subsection{Pullback and exact cancellation in an optimal criterion}
Return to arbitrary $K=\F_{q^k}$. For $f=t^2+at+c$, write
\[
 \mathcal N(f)=Z^2+B_{k,c}(a)Z+D,\qquad D=\Nm(c).
\]
The Lucas atoms and optimal criterion of Chapter~6 apply to the return coefficient $B_{k,c}(a)$. For primitive $D$, norm-primitivity is equivalent to
$P_{q,D}(B_{k,c}(a))=0$; for nonprimitive $D$, there are no such members.

\begin{theorem}[Optimal pullback over a quadratic coefficient field]\label{v148:thm-optimal-pullback}
Let $k=2$, let $D=\Nm(c)$ be primitive, and put $T=\Tr(c)$ and $\epsilon(c)=\mathbf1_{c\text{ primitive}}$. The unique monic polynomial of least degree vanishing on all $a$ for which $t^2+at+c$ is norm-primitive is
\begin{equation}\label{v148:eq-optimal-pullback}
 H_c(A)=(-1)^{d_q}\frac{P_{q,D}(T-A^{q+1})}{A^{q\epsilon(c)}}.
\end{equation}
The division is exact in the polynomial ring. The result is squarefree, splits over $K$, and has degree $(q+1)d_q-q\epsilon(c)$.
\end{theorem}
\begin{proof}
Direct calculation gives $B_{2,c}(a)=T-\Nm(a)$. For each primitive target coefficient $b$, the equation $A^{q+1}=T-b$ has $q+1$ distinct roots if $T-b\ne0$, and only the zero root otherwise. The latter case occurs exactly when
$(Z+c)(Z+\bar c)=Z^2+TZ+D$ is primitive, equivalently when $c$ is primitive. The pullback gives this zero root multiplicity $q+1$; division by $A^q$ reduces it to one. All other roots are simple, and different targets give disjoint root sets. The root bound proves minimality and uniqueness.
\end{proof}

For example, over $\F_4$ we have $H_1=A^3+1$, while $H_c=A$ for $c\in\F_4\setminus\F_2$. A zero middle coefficient gives every norm-primitive member in the latter case.

\begin{proposition}[Parameter strata and optimal degree]\label{v148:prop-optimal-strata}
For arbitrary $k$ and primitive $D$, the optimal polynomial is
\[
 H_{k,c}(A)=\operatorname{monic}\gcd\bigl(A^{q^k}-A,
                                   P_{q,D}(B_{k,c}(A))\bigr).
\]
Its degree is $C_{2,k}(q)d_q$ for odd $k$; for even $k=2s$, it is given by Corollary~\ref{v148:cor-quadratic-orders}. For nonprimitive $D$, put $H_{k,c}=1$.

For $k=2$, let
\[
 \mathcal A=K[C]/(\Phi_{q-1}(C^{q+1})),\qquad
 e(C)=1-\Phi_{q^2-1}(C)^{q^2-1}.
\]
This is the product of copies of $K$ indexed by constants with primitive norm, and $e(c)=\epsilon(c)$. In this product define componentwise
$G(A,C)=P_{q,C^{q+1}}(C+C^q-A^{q+1})$. Then
\[
 (-1)^{d_q}\left((1-e)G+\frac{eG}{A^q}\right)\in\mathcal A[A]
\]
specializes to $H_c$ on every component. If both degree strata are nonempty, no globally monic polynomial can have all these optimal specializations.
\end{proposition}
\begin{proof}
The monic gcd with $A^{q^k}-A$ selects exactly the eligible parameters as simple roots, so its degree is the fiber count. The defining polynomial of $\mathcal A$ is squarefree and splits over $K$, allowing the finite étale decomposition of Chapter~3. The element $e$ is the indicator idempotent of the primitive stratum. Since $A^{q+1}$ divides $eG$, the displayed division is valid. A globally monic polynomial has the same degree on every field component and cannot realize two distinct optimal degrees.
\end{proof}

\begin{proposition}[The first zero Lucas coefficient and short criteria]\label{v148:prop-lucas-return-order}
Let $f=t^2+at+c$ be skew irreducible with return norm $Z^2+BZ+D$. Put
$U_0=0$, $U_1=1$, and $U_{j+1}=-BU_j-DU_{j-1}$. If $r$ is the least positive index with $U_r=0$, then
\[
 r=\min\{j\ge1:f\text{ right divides }t^{kj}-\lambda
                         \text{ for some }\lambda\in\F_q^\times\},
 \qquad\lambda=-DU_{r-1}.
\]
The central return order is $r\ord(\lambda)$. For primitive $D$, norm-primitivity is equivalent to $r=q+1$.

Fix any $k\ge1$. For every (equivalently, for some) $c$ of primitive norm, all skew irreducible quadratics with $B\ne0$ are norm-primitive if and only if $q+1$ is a power of two, an odd prime, or twice an odd prime. If $q+1=4\pi$ with $\pi$ an odd prime, then for irreducible targets of primitive $D$, norm-primitivity is equivalent to $B^2\ne2D$.
\end{proposition}
\begin{proof}
Cayley--Hamilton gives $M_f^j=U_jM_f-DU_{j-1}I$. Since $M_f$ is nonscalar, its first scalar power is detected by the first zero of $U_j$. For a return root $\alpha$, we have $r=\ord(\alpha/\alpha^q)$ and $\alpha^r=\lambda$. This gives the order formula and the root-ratio criterion of Chapter~6.

Every compatible irreducible target has a positive fiber, by coprime uniformity for odd $k$ and Theorem~\ref{v148:thm-quadratic-even} for even $k$. It therefore suffices to compare ordinary quadratic targets. For even $q$, their total is $q/2$ and their primitive count is $\varphi(q+1)/2$; equality holds exactly when $q+1$ is prime. If $q\equiv1\pmod4$, removing the zero middle coefficient leaves $(q-1)/2$ targets, equal to $\varphi(q+1)$ exactly when $q+1$ is twice an odd prime. If $q\equiv3\pmod4$, the zero-middle-coefficient target is already reducible, and the condition becomes $\varphi(q+1)=(q+1)/2$, equivalent to $q+1$ being a power of two.

If $q+1=4\pi$, primitive norm forces the root-ratio order to contain the full power of two. The nonprimitive irreducible case has root-ratio order four, detected by $U_4=-B(B^2-2D)$. Here $B=0$ is already reducible, leaving precisely the exclusion $B^2=2D$.
\end{proof}

Right division at arbitrary exponents also gives
\[
 \begin{gathered}
 t^j\equiv u_jt+v_j\pmod{Rf},\qquad (u_0,v_0)=(0,1),\\
 u_{j+1}=\tau(v_j)-a\tau(u_j),\qquad v_{j+1}=-c\tau(u_j).
 \end{gathered}
\]
Here $u_j=0$ detects a right binomial relation with scalar in $K$, whereas the Lucas recurrence detects central return relations.

\section{Cubic fibers with a fixed middle coefficient}
\label{v148:sec-cubic-slices}
A fixed-constant total does not show whether every middle coefficient admits an irreducible member. In degree three, reducibility first becomes a value-set problem; trace--norm fibers then describe the overlaps in that value set.

\begin{theorem}[Cubic right roots and nonempty slices]\label{v148:thm-cubic-value-set}
Let $K=\F_{q^k}$, let $\tau$ generate $\operatorname{Gal}(K/\F_q)$, and let $b\in K$, $c\ne0$. Then
\begin{equation}\label{v148:eq-cubic-value-set}
 \begin{aligned}
 &\#\{a\in K:t^3+at^2+bt+c\text{ skew irreducible}\}\\
 &\quad=q^k-\left|\left\{-\tau^2(u)-\frac b{\tau(u)}
                         -\frac c{u\tau(u)}:u\in K^\times\right\}\right|\ge1.
 \end{aligned}
\end{equation}
\end{theorem}
\begin{proof}
A skew cubic with nonzero constant coefficient is reducible exactly when it has a linear right factor. Sufficiency is immediate. For necessity, its return norm is a reducible ordinary cubic and has a nonzero root $\lambda\in\F_q$. Choose $d\in K^\times$ with $\Nm(d)=\lambda$ and represent $t$ by $d\tau$. The matrix model of Theorem~\ref{v148:thm-all-fibers} at this central point is $\End_{\F_q}(K)$. Vanishing of the norm at $\lambda$ makes $f(d\tau)$ singular. For a nonzero kernel vector $x$, put $u=d\tau(x)/x$. Then
\[
 (d\tau)^j(x)/x=u\tau(u)\cdots\tau^{j-1}(u),
\]
so $u$ is a right root. Right evaluation of the cubic is
$u\tau(u)\tau^2(u)+a u\tau(u)+bu+c$; solving for $a$ gives the stated value set. Its domain has only $q^k-1$ elements, so at least one $a$ is omitted.
\end{proof}

Now take $K=\F_{q^2}$, $\tau(x)=x^q$, and put
\[
 w_q(T,N)=\#\{x\in K:\Tr(x)=T,\ \Nm(x)=N\}.
\]
The quadratic trace--norm criterion gives
\begin{equation}\label{v148:eq-trace-norm-weight}
 w_q(T,N)=
 \begin{cases}
 1-\chi(T^2-4N),&q\text{ odd},\\
 1,&q\text{ even and }T=0,\\
 1-(-1)^{\Tr_{\F_q/\F_2}(N/T^2)},&q\text{ even and }T\ne0,
 \end{cases}
\end{equation}
where $\chi(0)=0$. Two distinct base-field roots do not have the required extension trace and norm; a repeated root contributes one, and an irreducible quadratic contributes two.

\begin{theorem}[Norm-circle intersections and exact slices]\label{v148:thm-cubic-circles}
Let $M_b=\#\{a\in K:t^3+at^2+bt-1\text{ skew irreducible}\}$.
For $x,y\in\F_q^\times$, put $z=(xy)^{-1}$ and
\[
 W_b(x,y)=w_q\bigl(xy+xz+yz-\Nm(b),\ x+y+z+\Tr(b)\bigr),
\]
\[
 A_b=\sum_{\substack{x,y\in\F_q^\times\\x\ne y}}W_b(x,y),\qquad
 D_b=\sum_{\substack{x\in\F_q^\times\\x^3\ne1}}W_b(x,x^{-2}),\qquad
 \epsilon_b=\mathbf1_{b\in\F_q^\times}.
\]
Then
\begin{equation}\label{v148:eq-cubic-circles}
 M_b=1+q\epsilon_b+\frac{A_b+D_b}{3}.
\end{equation}
Division by three takes place in the integers, including in characteristic three. For any $c\ne0$, the corresponding irreducible slice satisfies
\[
 I_b(c)\ge1+q\epsilon_b.
\]
\end{theorem}
\begin{proof}
Calculation gives
\begin{equation}\label{v148:eq-cubic-return}
 \mathcal N(t^3+at^2+bt+c)
 =Z^3+(\Tr(b)-\Nm(a))Z^2
       +(\Nm(b)-\Tr(a\bar c))Z-\Nm(c).
\end{equation}
For $c=-1$, the event that $x\in\F_q^\times$ is a return root is
\[
 \Nm(a-x^{-1})=x^{-1}\Nm(b+x).
\]
There are $q+1$ solutions when the right side is nonzero, and one otherwise. The sum of all root-event sizes is therefore
$S_1=q^2-1-q\epsilon_b$. Two distinct roots $x,y$ force the third root to be $z=(xy)^{-1}$; coefficient comparison gives intersection size $W_b(x,y)$.

Let $S_3$ be the sum of intersections for three distinct roots. Targets with exactly two distinct roots are $(Z-x)^2(Z-x^{-2})$, $x^3\ne1$, and are counted by $D_b$. Hence $A_b=6S_3+2D_b$. Inclusion--exclusion gives
\[
 M_b=q^2-S_1+A_b/2-S_3
     =1+q\epsilon_b+(A_b+D_b)/3.
\]
For general $c$, the root event becomes
$\Nm(a+c/x)=\Nm(b+x)/x$. Its total size is still $S_1$, and the union bound proves the lower bound.
\end{proof}

Thus, for $b\in\F_q^\times$, at least $q$ nonzero values of $a$ give irreducible cubics with all four displayed coefficients nonzero. Formula~\eqref{v148:eq-cubic-circles} depends only on $t=\Tr(b)$ and $n=\Nm(b)$; denote its value by $M_q(t,n)$. The complete frequency distribution is
\[
 \#\{b:M_b=j\}=\sum_{t,n\in\F_q}w_q(t,n)\mathbf1_{M_q(t,n)=j}.
\]
For nonzero weight, $\epsilon_b=1$ exactly when $w_q(t,n)=1$ and $n\ne0$. The fixed-norm irreducible cubic count also gives
\[
 \sum_bM_b=(q^2-q+1)\frac{q^2+q+1-\gcd(3,q-1)}3.
\]
For $K=\F_4$ and $b=0,1,\omega,\omega^2$, the slice sizes are $1,3,1,1$.

\subsection{Trace--norm strata of primitive slices}
Fix $c\ne0$, put $D=\Nm(c)$, and let
$P_b(c)=\#\{a:t^3+at^2+bt+c\text{ norm-primitive}\}$.
For nonprimitive $D$, every slice is zero. For primitive $D$, define
\[
 \mathcal T_D=\{(A,B)\in\F_q^2:Z^3+AZ^2+BZ-D\text{ primitive}\},\qquad
 T_q=\frac{\varphi(q^3-1)}{3\varphi(q-1)}.
\]
This set has size $T_q$. It can also be parameterized by primitive roots:
$A=-\Tr(\alpha)$ and $B=D\Tr(\alpha^{-1})$, where
$\alpha\in\F_{q^3}$ is primitive with $\Nm(\alpha)=D$; each target occurs three times.

\begin{theorem}[Primitive cubic slices and moments]\label{v148:thm-cubic-primitive-slices}
For primitive $D$,
\begin{equation}\label{v148:eq-cubic-primitive-slice}
 P_b(c)=P_D(\Tr(b),\Nm(b)),\qquad
 P_D(t,n)=\sum_{(A,B)\in\mathcal T_D}w_q(n-B,D(t-A)).
\end{equation}
The frequency distribution, zero-slice count, and second moment are
\[
 \mathcal H_D(j)=\sum_{t,n}w_q(t,n)\mathbf1_{P_D(t,n)=j},\quad
 Z_D=\mathcal H_D(0),\quad
 M_2(D)=\sum_{t,n}w_q(t,n)P_D(t,n)^2,
\]
with all sums over $\F_q^2$. The first moment is
\[
 \sum_bP_b(c)=(q^2-q+1)T_q.
\]
\end{theorem}
\begin{proof}
For each target $(A,B)$, set $u=a\bar c$. Formula~\eqref{v148:eq-cubic-return} becomes
$\Tr(u)=\Nm(b)-B$ and $\Nm(u)=D(\Tr(b)-A)$. The map $a\mapsto u$ is bijective, proving the slice formula. Stratifying $b$ by trace and norm gives the frequencies and moments. Each primitive cubic target has fiber size $C_{3,2}(q)=q^2-q+1$, giving the first moment.
\end{proof}

In odd characteristic the slice formula is also
\[
 P_D(t,n)=T_q-\sum_{(A,B)\in\mathcal T_D}
             \chi\bigl((n-B)^2-4D(t-A)\bigr),
\]
a quadratic character sum over the base field.

\begin{proposition}[Cube classes and slice distributions]\label{v148:prop-cubic-classes}
Suppose $D,D'$ are primitive and $D'=D/s^3$ with $s\in\F_q^\times$. Then
\[
 P_{D'}(\Tr(b/s),\Nm(b/s))=P_D(\Tr(b),\Nm(b)).
\]
Consequently $\mathcal H_D,Z_D,M_2(D)$ depend only on the cube class of $D$ and are invariant under $D\mapsto D^p$, where $q=p^e$. If $3\nmid q-1$ or $p\equiv2\pmod3$, all primitive $D$ have the same distribution; otherwise there are at most two distributions.

If $q^2+q+1$ is prime, then for every extension $K=\F_{q^k}$ and every $c$ with primitive $(-1)^{3k}\Nm_{K/\F_q}(c)$, each fixed-$b$ cubic family contains a norm-primitive member. Over a quadratic coefficient extension, the lower bound is $1+q\epsilon_b$.
\end{proposition}
\begin{proof}
Choose $d\in K$ of norm $s$, substitute $t\mapsto dt$, and normalize to obtain
\[
 (a,b,c)\longmapsto(a/d,b/s,c/(ds)).
\]
The return target becomes $s^{-3}h(sZ)$, replacing $\alpha$ by $\alpha/s$. Primitivity is equivalent to primitive base-field norm and maximal order in $\F_{q^3}^\times/\F_q^\times$. Base-field scaling preserves the latter condition, and both prescribed norms are primitive. Thus the transformation preserves primitivity. Frobenius also gives a bijection. There is one stratum when cubing is surjective; otherwise primitive elements lie in the two nontrivial cube classes, exchanged by Frobenius when $p\equiv2\pmod3$.

Finally, if the quotient group has prime order $q^2+q+1$, the nontrivial image of an irreducible cubic root already has maximal order. With primitive norm, irreducibility then implies primitivity. Apply Theorem~\ref{v148:thm-cubic-value-set} and the norm-circle bound.
\end{proof}

Nonempty irreducible slices do not in general imply nonempty primitive slices. Formula~\eqref{v148:eq-cubic-primitive-slice} keeps the distinction explicit: total fibers may be uniform while fixing one more coefficient leaves a distribution governed by trace--norm correlations.

\section{Norm transfer for equally spaced powers}
\label{v148:sec-sparse-transfer}
Fiber formulas also apply to coefficient families with large exponent gaps. First change the Frobenius step to lower the skew degree; then determine whether the resulting ordinary polynomial remains irreducible after descent.

Let $K=\F_{q^k}$, $\tau(x)=x^q$, and $R=K[t;\tau]$. Consider
\[
 F(t)=t^{rn}+a_{r-1}t^{(r-1)n}+\cdots+a_1t^n+c,
 \qquad r\ge2,\quad n\ge1,\quad c\ne0.
\]
Put
\[
 d=\gcd(n,k),\quad m=n/d,\quad H=k/d,\quad P=q^d,
 \quad u=t^n,\quad S=K[u;\tau^n].
\]
The fixed field of $S$ is $\F_P$, and its central variable is $W=u^H$. In the original ring, $W=Z^m$, where $Z=t^k$. Write
$g=u^r+a_{r-1}u^{r-1}+\cdots+a_1u+c$ and $h=\mathcal N_S(g)$.

\begin{proposition}[Power embedding and norm descent]\label{v148:prop-sparse-transfer}
We have
\begin{equation}\label{v148:eq-sparse-transfer}
 \mathcal N_R(F)(Z)=\prod_{j=0}^{d-1}h^{(q^j)}(Z^m),\qquad
 h(0)=(-1)^{r(H+1)}\Nm_{K/\F_P}(c).
\end{equation}
The superscript acts only on coefficients. The right side is irreducible exactly when $h$ is irreducible and its root $\alpha$ satisfies
\begin{equation}\label{v148:eq-admissible-root}
 [\F_q(\alpha):\F_q]=rd,\qquad
 \Nm_{\F_{P^r}/\F_P}(\alpha)=\nu,\qquad
 X^m-\alpha\text{ is irreducible over }\F_{P^r},
\end{equation}
where $\nu=(-1)^{rH}\Nm_{K/\F_P}(c)$.
\end{proposition}
\begin{proof}
Insert $S'=K[v;\tau^d]$ with $v=t^d$, so $u=v^m$ and $v^H=Z$. Over the rational function field, where central variables are invertible, extend scalars by $W\mapsto Z^m$. Since $\gcd(m,H)=1$, choose $Am+BH=1$; then $v=u^AZ^B$, giving an isomorphism of the two cyclic representations. Determinants of split matrices give the power substitution for norms. The monic-normalization sign difference is $(-1)^{r(m-1)(H-1)}=1$.

Next group the original representation by exponents modulo $d$, obtaining $d$ blocks of size $H$. Their determinants are the $d$ coefficient conjugates of $h(Z^m)$. The other normalization sign has exponent $rmHd(d-1)$, again even. This proves the transfer formula. Since both sides are polynomials, the rational-function identity is a polynomial identity.

If the right side is irreducible, $h$ is irreducible and its $d$ coefficient conjugates are distinct. Equivalently, the $q$-orbit of $\alpha$ has length $rd$, with minimal polynomial $J(Z)=\prod_jh^{(q^j)}(Z)$. For $\beta^m=\alpha$, the degree tower gives
\[
 [\F_q(\beta):\F_q]
 =rd\,[\F_{P^r}(\beta):\F_{P^r}],
\]
so irreducibility of $J(Z^m)$ is equivalent to that of the indicated binomial. The constant coefficient gives the prescribed norm $\nu$.
\end{proof}

Let $\mathcal A_c$ be the set of $\F_P$-minimal polynomials satisfying~\eqref{v148:eq-admissible-root}. Without any coprimality assumption, the irreducible coefficient-vector count is
\begin{equation}\label{v148:eq-weighted-sparse}
 \#\{(a_1,\ldots,a_{r-1}):F\text{ skew irreducible}\}
 =\sum_{h\in\mathcal A_c}F_{h,H}^{\tau^n}(c).
\end{equation}
This separates target selection from fiber weights. The two counts may be multiplied directly only when the weights are equal.

\begin{theorem}[Exact counts after coprime degree reduction]\label{v148:thm-sparse-coprime}
Suppose $\gcd(r,H)=1$. Put $Q=P^r$, $S_r=(Q-1)/(P-1)$, $e=\ord(\nu)$, and $C=C_{r,H}(P)$.
If $4\mid m$ and $Q\not\equiv1\pmod4$, the count is zero. Otherwise define
\[
 \mathcal T_e=\left\{\ell\mid Q-1:
 \begin{array}{l}
 \ell>1,\ \ord_\ell(q)=rd,\quad \ell/\gcd(\ell,S_r)=e,\\
 \rad(m)\mid\ell,\quad\gcd(m,(Q-1)/\ell)=1
 \end{array}\right\}.
\]
Then
\begin{equation}\label{v148:eq-sparse-coprime}
 \#\{\mathbf a:F\text{ skew irreducible}\}
 =\frac C{r\varphi(e)}\sum_{\ell\in\mathcal T_e}\varphi(\ell).
\end{equation}
An irreducible member exists if and only if
\begin{equation}\label{v148:eq-sparse-existence}
 \rad(m)\mid S_re,\qquad
 \gcd\bigl(m,(P-1)/e\bigr)=1,\qquad
 4\nmid m\ \text{or}\ Q\equiv1\pmod4.
\end{equation}
Each realizable original return target has exactly $(d/\delta)C$ preimages, where
$\delta=[\F_q(\nu):\F_q]$. With $\ord_1(q)=1$, this is also $\delta=\ord_e(q)$.
\end{theorem}
\begin{proof}
Theorem~\ref{v148:thm-relative-fibers} applies to the generator $\tau^n$, so each compatible irreducible subtarget has $C$ preimages. A root of order $\ell$ has norm-image order $\ell/\gcd(\ell,S_r)$, leaving $\varphi(\ell)/\varphi(e)$ roots over the specified $\nu$. Full degree is the condition $\ord_\ell(q)=rd$; the remaining conditions are precisely the binomial criterion of Theorem~\ref{thm:binomial-full}. Divide by the $r$ roots of each subtarget.

If an allowed $\ell$ exists, comparison at each prime gives $\ell\mid S_re$, proving necessity. Conversely, take $\ell=S_re$. It divides $Q-1$, has norm-image order $e$, and satisfies the binomial criterion. Since $r\ge2$,
\[
 S_r>q^{(r-1)d}\ge q^{rd/2},
\]
this order cannot divide the multiplicative order of any proper subfield. Its degree is therefore $rd$, and cyclic-group norm surjectivity supplies a root over the specified $\nu$.

Finally, one original target gives $d$ distinct subtargets with root norms $\nu,\nu^q,\ldots,\nu^{q^{d-1}}$. Exactly $d/\delta$ satisfy the fixed-$c$ condition; each contributes $C$.
\end{proof}

\begin{corollary}[Exact return orders and a primitivity obstruction]\label{v148:cor-sparse-order}
Under the preceding theorem, the number of members of return order $m\ell$, for $\ell\in\mathcal T_e$, is
$C\varphi(\ell)/(r\varphi(e))$; no other return order occurs.

Regardless of coprimality, norm-primitivity forces $n\mid k$. If $n\mid k$ and $\gcd(r,k/n)=1$, put $P=q^n$ and $H=k/n$. The norm-primitive count is
\[
 C_{r,H}(P)\frac{\varphi(P^r-1)}{r\varphi(P-1)}
 \quad\text{if }(-1)^{rH}\Nm_{K/\F_P}(c)\text{ is primitive},
\]
and zero otherwise. In particular, the family with $c=(-1)^r$ can have norm-primitive members only when $q=2,n=1$.
\end{corollary}
\begin{proof}
If $\beta^m=\alpha$, $\ord(\alpha)=\ell$, and $\rad(m)\mid\ell$, comparison of prime valuations in
$\ell=\ord(\beta)/\gcd(\ord(\beta),m)$ gives $\ord(\beta)=m\ell$.

If $m>1$, every original return root satisfies $\beta^m\in\F_Q^\times$, so
\[
 \ord(\beta)\le m(Q-1)<Q^m-1=q^{rn}-1,
\]
precluding primitivity. Thus $m=1$, or $n\mid k$. The original and subtarget roots then coincide; apply the fixed-primitive-norm count. For $c=(-1)^r$, the subtarget root has norm one, which can generate $\F_{q^n}^\times$ only if $q^n=2$.
\end{proof}

\subsection{Norm-one families and noncoprime weights}
If $c=(-1)^r$, then $\nu=1$. In the coprime case, existence reduces to
\[
 \rad(m)\mid S_r,\qquad\gcd(m,P-1)=1.
\]
These conditions already force $m$ odd. Orders in $\mathcal T_1$ divide $S_r$, and each original target has weight $dC$.

\begin{corollary}[Trinomials with arbitrary spacing]\label{v148:cor-trinomial-spacing}
For $F=t^{2n}+at^n+1$, no assumption $\gcd(2,H)=1$ is needed. Put
\[
 D_H(P)=
 \begin{cases}
 (P^H+1)/(P+1),&H\text{ odd},\\
 (P^H+1-2(-P)^{H/2})/(P+1),&H\text{ even}.
 \end{cases}
\]
The irreducible count is zero if $m$ is even or $\rad(m)\nmid P+1$. Otherwise, let
$L=\prod_{\rho\mid m\text{ prime}}\rho^{v_\rho(P+1)}$. The count is
\begin{equation}\label{v148:eq-general-trinomial}
 \frac{D_H(P)}2
 \sum_{\substack{\ell\mid P+1,\ \ell>2\\\ord_\ell(q)=2d,\ L\mid\ell}}
 \varphi(\ell).
\end{equation}
Writing $R_m=\rad(m)$ and $g_j=\gcd(P+1,q^j-1)$, the target count is also
\[
 \frac{\varphi(R_m)}{2R_m}
 \sum_{\substack{j\mid2d\\L\mid g_j}}\mu(2d/j)g_j.
\]
\end{corollary}
\begin{proof}
Subtargets are irreducible norm-one quadratics. Theorem~\ref{v145:thm-return-uniform} and its generator invariance give their common weight $D_H(P)$. The binomial conditions require the root order to contain the full prime power in $P+1$ for each prime $\rho\mid m$, exactly the condition $L\mid\ell$. Norm one excludes even $m$, proving the first formula.

The intersection of the norm-one group with $\F_{q^j}$ has order $g_j$. If $L\mid g_j$, it contains $g_j\varphi(R_m)/R_m$ elements with the required full prime powers in their orders; otherwise it contains none. Möbius inversion over subfields selects degree $2d$; division by two gives the other formula. For $m=1$, take $L=R_m=1$.
\end{proof}

For $n\mid k$, we have $m=1,d=n$, recovering Theorem~\ref{v145:thm-sparse-count}. In higher noncoprime degrees one can retain the weights in~\eqref{v148:eq-weighted-sparse}. For example, when $r=4,H=2,c=1$, the allowed order set $\mathcal T_1$ gives
\begin{equation}\label{v148:eq-sparse-quartic-weight}
 \#\{\mathbf a:F\text{ skew irreducible}\}
 =\frac14\sum_{\ell\in\mathcal T_1}\varphi(\ell)
       \left((P-1)^2(P+1)+2P^2\mathbf1_{\ell\mid P^2+1}\right).
\end{equation}
Indeed, the constant of a quadratic factor over $\F_{P^2}$ is $\alpha^{P^2+1}$. At $c=1$, the two corrections in the three-value formula occur together, exactly when $\ell\mid P^2+1$. The formula retains the two fiber weights determined by root order.

\begin{example}[Changing the constant changes realizable targets]\label{v148:ex-sparse-constants}
If $q,k$ are odd, then
\[
 \#\{a\in\F_{q^k}:t^4+at^2+c\text{ skew irreducible}\}
 =\begin{cases}(q^k+1)/2,&\Nm(c)\text{ nonsquare},\\0,&\Nm(c)\text{ square}.
 \end{cases}
\]
Here $r=n=2,d=1,m=2$. An element $\alpha\in\F_{q^2}$ is nonsquare exactly when its norm is nonsquare. All $q+1$ elements of such a norm fiber lie outside the base field, giving $(q+1)/2$ subtargets, each of weight $(q^k+1)/(q+1)$. Thus over $\F_{27}$, constant $1$ gives no members and constant $-1$ gives $14$; over $\F_{125}$, constant $2$ gives $63$.

Next take $q=2,k=4,r=3,n=2$. Then $P=4,H=2,m=1$, and $C_{3,2}(4)=13$. For
$t^6+at^4+bt^2+c$, relative norm $\Nm_{\F_{16}/\F_4}(c)=1$ allows root order $21$, giving four subtargets and two original targets: $52$ members, with weight $26$ per original target. A nontrivial relative norm in $\F_4^\times$ allows orders $9,63$, giving seven subtargets and seven original targets: $91$ members, with weight $13$ each. The distinction is exactly $\delta=1$ versus $\delta=2$.
\end{example}

\section{Generating Polynomials and Linear Recurrences for Extension-Field Character Sums}
\label{sec:character-sum-recurrence}

The Hasse--Davenport formula expresses lifted Gauss sums of fixed base-field characters over successive extensions as signed powers of a single base-field Gauss sum.
More generally, if a weighted orbit product is a polynomial, the corresponding extension-field character sums are power sums of finitely many numbers and hence satisfy a linear recurrence of fixed order.
We first prove the polynomial-phase case, then give a complete calculation for Kloosterman sums, which include an inverse term \cite[Chapter 8]{Feng}.

\subsection{Polynomial Phases and Degree Cancellation}
Let \(k=\F_q\) have characteristic \(p\), let \(\psi:k\to\C^\times\) be a nontrivial additive character,
and let \(E_r=\F_{q^r}\). For \(f\in k[X]\), put
\begin{equation}\label{def:character-L-series}
S_r(f)=\sum_{x\in E_r}\psi\bigl(\Tr_{E_r/k}(f(x))\bigr),\qquad
L_f(u)=\exp\left(\sum_{r\ge1}S_r(f)\frac{u^r}{r}\right).
\end{equation}
Here \(L_f(u)\) is called the associated \textbf{\(L\)-series}, defined in the formal power series ring \(\C[[u]]\).
Exponentials, logarithms, and the derivatives below are formal operations, with no convergence issue.

\begin{theorem}[Phase reduction with the character parameter]\label{v70:thm-phase-reduction}
Let \(q=p^e\), fix a primitive complex \(p\)th root of unity \(\zeta_p\), and write the nontrivial additive character as
\[
\psi(z)=\zeta_p^{\Tr_{k/\F_p}(az)},\qquad a\in k^\times.
\]
For every \(f\in k[X]\), there are unique \(R,H\in k[X]\) such that
\begin{equation}\label{v70:eq-phase-normal-form}
af=R+H^p-H,\qquad H(0)=0,
\end{equation}
and every positive exponent with nonzero coefficient in \(R\) is prime to \(p\). If \(af=\sum_{j\ge0}c_jX^j\), then
\begin{equation}\label{v70:eq-phase-reduced-coefficients}
R=c_0+\sum_{\substack{j\ge1\\p\nmid j}}
 \left(\sum_{s\ge0}c_{jp^s}^{p^{-s}}\right)X^j.
\end{equation}
Here \(c^{p^{-s}}\) denotes the unique \(b\in k\) satisfying \(b^{p^s}=c\), and all sums are finite. For every \(r\ge1\) and \(x\in E_r\),
\begin{equation}\label{v70:eq-phase-all-extensions}
\psi\bigl(\Tr_{E_r/k}(f(x))\bigr)
=\zeta_p^{\Tr_{E_r/\F_p}(R(x))}.
\end{equation}
Replacing \(af\) by this \(R\) is called \textbf{phase reduction}.
\end{theorem}
\begin{proof}
Every positive exponent has a unique expression \(jp^s\), where \(p\nmid j\). For the term \(c_{jp^s}X^{jp^s}\), put \(b=c_{jp^s}^{p^{-s}}\). If \(s\ge1\), take
\[
A=\sum_{t=0}^{s-1}(bX^j)^{p^t}.
\]
Cancellation of successive terms gives
\[
(bX^j)^{p^s}-bX^j=A^p-A.
\]
For \(s=0\), the term needs no change. Summing over all positive-degree terms and retaining the constant term \(c_0\) gives~\eqref{v70:eq-phase-normal-form} and~\eqref{v70:eq-phase-reduced-coefficients}; the resulting \(H\) has zero constant term.

Suppose that another decomposition is \(af=R'+(H')^p-H'\), and put \(A=H-H'\). Then \(R'-R=A^p-A\) and \(A(0)=0\). If \(A\ne0\), the positive degree \(\deg(A^p-A)=p\deg A\) is divisible by \(p\), whereas every positive exponent with nonzero coefficient in \(R'-R\) is prime to \(p\), a contradiction. Hence \(A=0\) and \(R=R'\).

Transitivity of the trace and \(a\in k\) give
\[
\Tr_{k/\F_p}\bigl(a\Tr_{E_r/k}(f(x))\bigr)
=\Tr_{E_r/\F_p}(af(x)).
\]
The absolute trace is invariant under the \(p\)th power, so \(\Tr_{E_r/\F_p}(H(x)^p-H(x))=0\), proving~\eqref{v70:eq-phase-all-extensions}.
\end{proof}

\begin{example}[The character parameter must be retained in phase reduction]\label{v70:ex-phase-parameter}
Let \(k=\F_4=\F_2(\omega)\), with \(\omega^2+\omega+1=0\), and take \(f=X^2+X\). For the canonical character \(\psi_1(z)=(-1)^{\Tr_{k/\F_2}(z)}\), the reduced polynomial is zero. Thus the phase is identically one on every \(E_r=\F_{4^r}\), and the character sum is \(4^r\).
For \(\psi_\omega(z)=(-1)^{\Tr_{k/\F_2}(\omega z)}\), however,
\[
\omega(X^2+X)=X+(\omega^2X)^2+\omega^2X.
\]
The reduced polynomial is now \(X\), so additive-character orthogonality makes the sum zero over every extension. Reducing \(f\) without including the character parameter would conflate these two cases.
\end{example}

Phase reduction preserves the trace phase over every extension, rather than the original polynomial function. Uniqueness in the theorem concerns the polynomial decomposition, not a representation of the trace function on one fixed field. In Example~\ref{v68:ex-trace-x5}, \(X^5\) already has no even exponent, yet \(\Tr(x^5)\equiv0\) on \(\F_{16}\). Evaluation on a fixed field can introduce further degeneracy, and absolute traces also fail to distinguish all constants.

\begin{theorem}[Generating polynomials for additive character sums]\label{thm:additive-L-polynomial}
If \(f\) has degree \(d\ge1\) and \(p\nmid d\), then \(L_f(u)\) is a polynomial with constant term one and degree at most \(d-1\).
It is then called the associated \textbf{\(L\)-polynomial}.
\end{theorem}
\begin{proof}
For a monic polynomial \(F\in k[X]\), list its roots in a splitting field with multiplicity as \(\beta_1,\ldots,\beta_N\) and define
\[
w_f(F)=\psi\left(\sum_{i=1}^N f(\beta_i)\right),\qquad w_f(1)=1.
\]
The root multiset is Frobenius-invariant, so the sum in parentheses belongs to \(k\).
The root multiset of a product is the multiset union of the two root multisets, hence \(w_f(FH)=w_f(F)w_f(H)\).
Unique factorization gives
\begin{equation}\label{eq:additive-weight-euler}
\sum_{F\text{ monic}}w_f(F)u^{\deg F}
=\prod_{P\text{ monic irreducible}}\bigl(1-w_f(P)u^{\deg P}\bigr)^{-1}.
\end{equation}
If \(P\) has degree \(e\mid r\) and \(\alpha\) is any of its roots, then
\[
w_f(P)^{r/e}
=\psi\left(\frac re\Tr_{\F_{q^e}/k}(f(\alpha))\right)
=\psi\bigl(\Tr_{E_r/k}(f(\alpha))\bigr).
\]
Each orbit has \(e\) elements. Taking the formal logarithm of the right-hand side of~\eqref{eq:additive-weight-euler} and collecting by powers of \(u\) gives
\[
\sum_P\sum_{j\ge1}\frac{w_f(P)^j u^{j\deg P}}j
=\sum_{r\ge1}\left(\sum_{\deg P\mid r}\deg P\,w_f(P)^{r/\deg P}\right)\frac{u^r}r
=\sum_{r\ge1}S_r(f)\frac{u^r}r.
\]
The last equality partitions \(E_r\) into Frobenius orbits. Thus both sides of~\eqref{eq:additive-weight-euler} equal \(L_f(u)\).

Now suppose that \(N\ge d\) and write
\[
F(X)=X^N+c_1X^{N-1}+\cdots+c_N,\qquad
s_j=\sum_{i=1}^N\beta_i^j.
\]
Applying Lemma~\ref{v58:lem-newton} to these \(N\) roots shows that \(s_j\) depends only on \(c_1,\ldots,c_j\), and that the coefficient of \(c_d\) in \(s_d\) is \(-d\).
Writing \(f(X)=a_dX^d+\cdots+a_0\) gives
\[
\sum_i f(\beta_i)=Na_0+\sum_{j=1}^d a_js_j
=-da_dc_d+B(c_1,\ldots,c_{d-1}),
\]
where \(B\) may also depend on the fixed integer \(N\), but not on \(c_d,\ldots,c_N\).
Fixing all coefficients except \(c_d\), the condition \(da_d\ne0\) and additive character orthogonality give
\[
\sum_{c_d\in k}\psi(-da_dc_d+B)=0.
\]
Thus \(\sum_{\deg F=N}w_f(F)=0\) for every \(N\ge d\).
The left-hand side of~\eqref{eq:additive-weight-euler} therefore terminates at degree \(d-1\), proving the theorem.
\end{proof}

\begin{corollary}[Three cases for the \(L\)-series of a reduced phase]\label{v70:cor-reduced-phase-L-series}
Let \(R\) be as in Theorem~\ref{v70:thm-phase-reduction}. The original sums \(S_r(f)\) and their \(L\)-series fall into the following three cases.
\begin{enumerate}
\item If \(R=c\) is constant, including zero, put \(\xi=\zeta_p^{\Tr_{k/\F_p}(c)}\). Then
\[
S_r(f)=q^r\xi^r,\qquad L_f(u)=(1-q\xi u)^{-1}.
\]
\item If \(\deg R=1\), then \(S_r(f)=0\) for all \(r\ge1\), and \(L_f(u)=1\).
\item If \(\deg R=d\ge2\), then \(p\nmid d\), and \(L_f(u)\) is a polynomial with constant term one and degree at most \(d-1\).
\end{enumerate}
\end{corollary}
\begin{proof}
Equation~\eqref{v70:eq-phase-all-extensions} expresses the original sum as the sum for \(R\) with the canonical character.
If \(R=c\), then \(\Tr_{E_r/\F_p}(c)=r\Tr_{k/\F_p}(c)\), so each summand is \(\xi^r\); substitution into the formal exponential gives the first case.
If \(R=c+vX\), with \(v\ne0\), then \(x\mapsto\zeta_p^{\Tr_{E_r/\F_p}(vx)}\) is a nontrivial additive character and sums to zero.
Finally, every positive exponent with nonzero coefficient in \(R\) is prime to \(p\), so \(p\nmid d\); apply Theorem~\ref{thm:additive-L-polynomial} to \(R\) and the canonical character of \(k\).
\end{proof}

Thus the effective degree of the phase is determined after including the character parameter and reducing. The third case leads to the finite-order recurrence below; degree reduction alone gives no bound on the absolute values of the characteristic roots.

\subsection{Characteristic Roots and Recurrences}
\begin{corollary}[A finite-order recurrence for extension-field character sums]\label{cor:character-sum-recurrence}
Under the hypotheses of Theorem~\ref{thm:additive-L-polynomial}, suppose that
\[
L_f(u)=1+\ell_1u+\cdots+\ell_mu^m
=\prod_{i=1}^m(1-\omega_i u),\qquad m=\deg L_f.
\]
Then, for all \(r\ge1\),
\begin{equation}\label{eq:character-power-sum}
S_r(f)=-\sum_{i=1}^m\omega_i^r.
\end{equation}
For \(m\ge1\), the initial values are determined by
\begin{equation}\label{eq:character-initial-newton}
S_r(f)+\ell_1S_{r-1}(f)+\cdots+\ell_{r-1}S_1(f)=r\ell_r
\qquad(1\le r\le m)
\end{equation}
and subsequent values satisfy
\begin{equation}\label{eq:character-sum-recurrence}
S_r(f)+\ell_1S_{r-1}(f)+\cdots+\ell_mS_{r-m}(f)=0
\qquad(r>m).
\end{equation}
If \(m=0\), all \(S_r(f)\) vanish. To use \(d-1\) characteristic roots uniformly, append zeros after \(\omega_1,\ldots,\omega_m\);
this does not change any positive power sum.
\end{corollary}
\begin{proof}
Take formal logarithms of the factorization of \(L_f\) and compare with~\eqref{def:character-L-series} to obtain~\eqref{eq:character-power-sum}.
On the other hand, the definition directly gives
\[
uL_f'(u)=L_f(u)\sum_{r\ge1}S_r(f)u^r.
\]
Compare coefficients of \(u^r\): the left-hand side is \(r\ell_r\) when \(1\le r\le m\) and zero when \(r>m\),
giving the initial relations and recurrence, respectively. The constant polynomial \(L_f=1\) has formal logarithm zero, proving the last assertion as well.
\end{proof}

Computation does not require first finding \(\omega_i\): given the first \(d-1\) character sums, set \(\ell_0=1\)
and compare coefficients successively in the identity from the proof to recover
\(L_f\) from \(j\ell_j=\sum_{i=1}^jS_i(f)\ell_{j-i}\),
then use~\eqref{eq:character-sum-recurrence} to compute all subsequent values.
The integer divisions needed for coefficient recovery take place in \(\C\).

The recurrence index here is the extension degree \(r\), and its coefficients lie in \(\C\).
The cyclic-module method of Chapter~\ref{chap:sequences} still applies to this constant-coefficient recurrence:
one annihilating polynomial of the sequence is \(T^m+\ell_1T^{m-1}+\cdots+\ell_m\),
whose roots are exactly \(\omega_1,\ldots,\omega_m\). However, eventual periodicity for finite-field recurrences relies on finiteness of the state space;
the complex-valued recurrences here do not meet that hypothesis.
In the Hasse--Davenport formula, the weighted generating polynomial has degree one, so only one characteristic root occurs;
for a general phase, several characteristic roots jointly control the sums over all extensions.
Polynomiality determines the recurrence order, whereas absolute values of characteristic roots determine magnitude bounds; the former alone does not imply \(\lvert\omega_i\rvert=\sqrt q\).
The diagonal-equation estimate of Theorem~\ref{thm:diagonal-weil-bound} comes from the explicit absolute-value formula for Gauss sums,
not as a general consequence of this formal recurrence.

\subsection{A Second-Order Recurrence for Kloosterman Sums}
\begin{definition}\label{def:kloosterman-sum}
Let \(a,b\in k^\times\), and define the extension-field \textbf{Kloosterman sum} by
\[
K_r(a,b)=\sum_{x\in E_r^\times}
\psi\bigl(\Tr_{E_r/k}(ax+bx^{-1})\bigr)\qquad(r\ge1).
\]
For fixed \(a,b,\psi\), abbreviate it to \(K_r\).
\end{definition}

\begin{theorem}[The lifting recurrence for Kloosterman sums]\label{thm:kloosterman-recurrence}
With the preceding notation,
\begin{equation}\label{eq:kloosterman-L}
\exp\left(\sum_{r\ge1}K_r\frac{u^r}{r}\right)=1+K_1u+qu^2.
\end{equation}
With the convention \(K_0=-2\),
\begin{equation}\label{eq:kloosterman-recurrence}
\boxed{\quad K_r=-K_1K_{r-1}-qK_{r-2}\qquad(r\ge2).\quad}
\end{equation}
In particular, \(K_2=2q-K_1^2\).
\end{theorem}
\begin{proof}
For a monic polynomial \(F\) with nonzero constant term, list its roots with multiplicity as \(\beta_i\) and put
\[
w(F)=\psi\left(\sum_i(a\beta_i+b\beta_i^{-1})\right).
\]
Also put \(w(1)=1\), and set \(w(F)=0\) when the constant term is zero.
Multiset union of root multisets gives \(w(FH)=w(F)w(H)\).
The same unique-factorization expansion as in~\eqref{eq:additive-weight-euler} gives an Euler product;
since \(w(X)=0\), only orbits of nonzero roots remain.
For each orbit whose degree divides \(r\), the \(r/\deg P\)th power of its weight equals the lifted additive character value over \(E_r/k\).
Thus the formal logarithm of this product is exactly \(\sum_{r\ge1}K_ru^r/r\).

Compute the coefficients on the monic-polynomial side. A linear polynomial \(X-x\) contributes
\(\psi(ax+bx^{-1})\), so the linear coefficient is \(K_1\).
For \(F=X^N+c_1X^{N-1}+\cdots+c_N\) and \(c_N\ne0\), when \(N\ge2\) we have
\[
\sum_i\beta_i=-c_1,\qquad
\sum_i\beta_i^{-1}=-c_{N-1}/c_N,
\]
so \(w(F)=\psi(-ac_1-bc_{N-1}/c_N)\).
When \(N\ge3\), the coefficients \(c_1\) and \(c_{N-1}\) are distinct free coefficients. Fixing all other coefficients and summing over \(c_1\),
the condition \(a\ne0\) makes the total sum zero.
When \(N=2\), these two coefficients coincide, so the quadratic coefficient is
\[
\sum_{c_2\in k^\times}\sum_{c_1\in k}
\psi\bigl(-(a+b/c_2)c_1\bigr)=q.
\]
In the last step, the inner sum is \(q\) only when \(c_2=-b/a\) and zero otherwise.
This proves~\eqref{eq:kloosterman-L}.

Write \(1+K_1u+qu^2=(1-\omega_1u)(1-\omega_2u)\). Then
\(\omega_1+\omega_2=-K_1\) and \(\omega_1\omega_2=q\), and taking formal logarithms gives
\(K_r=-(\omega_1^r+\omega_2^r)\).
Each \(\omega_i\) satisfies \(\omega_i^2+K_1\omega_i+q=0\), so their power sums satisfy~\eqref{eq:kloosterman-recurrence}.
Negating the zeroth power sum gives \(K_0=-2\); this is an initial-value convention for the recurrence, not a sum over a degree-zero extension.
\end{proof}

\begin{corollary}[A Dickson representation of Kloosterman sums]\label{cor:kloosterman-dickson}
In \(\Z[Z,A]\), define the Dickson polynomials of the first kind with parameter \(A\) by
\[
\begin{gathered}
D_0(Z,A)=2,\qquad D_1(Z,A)=Z,\\
D_r(Z,A)=ZD_{r-1}(Z,A)-AD_{r-2}(Z,A)\quad(r\ge2).
\end{gathered}
\]
Then, for \(r\ge0\),
\[
K_r=-D_r(-K_1,q).
\]
\end{corollary}
\begin{proof}
For \(r=0,1\) the right-hand side has values \(-2,K_1\), respectively, and satisfies~\eqref{eq:kloosterman-recurrence}; it therefore agrees term by term with \(K_r\).
\end{proof}

For parameter \(A=1\), these are the polynomials \(D_r(Z)\) of Definition~\ref{def:dickson-polynomial}.
A general parameter retains the product of the two roots: if \(z_1+z_2=Z\) and \(z_1z_2=A\), then
\(D_r(Z,A)=z_1^r+z_2^r\), as the same initial values and recurrence verify.
For permutation polynomials, this representation descends exponentiation of root parameters to a polynomial;
here it expresses the power sum of two complex characteristic roots as a polynomial in \(K_1\) and \(q\).
In this formula \(q\) is a positive integer in \(\C\), not an element of the base field of characteristic \(p\).

\begin{example}[Recovering extension Kloosterman sums from the binary field]\label{exm:binary-kloosterman}
Take \(k=\F_2\), \(a=b=1\), and \(\psi(t)=(-1)^t\). The only nonzero element contributes \(\psi(1+1)=1\), so \(K_1=1\).
Thus
\[
L(u)=1+u+2u^2,\qquad
K_2=3,\quad K_3=-5,\quad K_4=-1,\quad K_5=11.
\]
For example, in \(\F_4=\F_2(\omega)\) with \(\omega^2+\omega+1=0\),
\(x=1\) gives \(x+x^{-1}=0\), and the other two nonzero elements give \(x+x^{-1}=1\).
Both values have trace zero to \(\F_2\), so all three summands are one, directly verifying \(K_2=3\).
No new enumeration of extension-field elements is needed for other extension degrees; the same second-order recurrence determines the sums.
\end{example}

\section{From Kloosterman Sums to Curve Point Counts and Zeta Functions}
\label{sec:kloosterman-curve-points}

Kloosterman sums can also record solution counts in fibers of two-variable equations.
We apply the Artin--Schreier trace criterion to a family of curves in characteristic two, then recover their extension-field point counts from the second-order recurrence of the preceding section.

\subsection{The Curve Model and Trace Fibers}

Fix \(k=\F_q\), \(q=2^m\), and \(b\in k^\times\), and write \(k_n=\F_{q^n}\).
To include both affine solutions and points at infinity, first specify projective coordinates.
\leavevmode\phantomsection\label{def:projective-plane-points}\ignorespaces
The \textbf{projective plane} \(\mathbf P^2(K)\) over a field \(K\) consists of equivalence classes of nonzero triples \((X,Y,Z)\) under multiplication by a common nonzero scalar, denoted by \([X:Y:Z]\).
For \(Z\ne0\), take the unique representative \([x:y:1]\), where \(x=X/Z\) and \(y=Y/Z\); points with \(Z=0\) are called \textbf{points at infinity}.

Consider the homogeneous cubic equation
\begin{equation}\label{eq:binary-curve-model}
 E_b:\qquad Y^2Z+XYZ=X^3+b^2Z^3.
\end{equation}
\leavevmode\phantomsection\label{def:rational-curve-point}\ignorespaces
For an extension field \(K/k\), a \(K\)-\textbf{rational point} is a point of \(\mathbf P^2(K)\) satisfying the equation; their set is denoted by \(E_b(K)\).
For \(Z=0\), the equation gives \(X=0\), so the unique point at infinity is \(O=[0:1:0]\). Hence
\[
 E_b(K)=\{(x,y)\in K^2:y^2+xy=x^3+b^2\}\ \sqcup\ \{O\},
\]
where \((x,y)\) abbreviates the projective point \([x:y:1]\).

\leavevmode\phantomsection\label{def:nonsingular-plane-curve}\ignorespaces
A plane curve given by a homogeneous equation \(H=0\) is \textbf{nonsingular} if, at every curve point over every extension field, its three formal partial derivatives do not vanish simultaneously.
Here take \(H=Y^2Z+XYZ+X^3+b^2Z^3\). Then
\[
 H_X=YZ+X^2,\qquad H_Y=XZ,\qquad
 H_Z=Y^2+XY+b^2Z^2.
\]
At a point with \(Z=1\), the condition \(H_Y=H_X=0\) successively gives \(X=Y=0\), but then \(H=b^2\ne0\).
At \(O\) we have \(H_Z=1\). Thus \(E_b\) is nonsingular.
\leavevmode\phantomsection\label{def:binary-elliptic-curve}\ignorespaces
A nonsingular projective cubic with a specified rational point \(O\) is a plane model of an \textbf{elliptic curve}; see \cite[Chapter III]{SilvermanAEC} for the general theory.
The nonzero condition \(b\ne0\) is essential: if \(b=0\), the point \((0,0)\) is singular.

Put
\[
 \psi_n(t)=(-1)^{\Tr_{k_n/\F_2}(t)}.
\]
Transitivity of the trace gives \(\psi_n=\psi_1\circ\Tr_{k_n/k}\), so this is exactly the lifted character of Definition~\ref{def:kloosterman-sum}.
\leavevmode\phantomsection\label{def:curve-point-count}\ignorespaces
Write \(N_n=\#E_b(k_n)\); the count always includes \(O\).

\Needspace{5\baselineskip}
\begin{proposition}[Curve point counts and Kloosterman sums]\label{prop:binary-curve-point-count}
For every \(n\ge1\),
\begin{equation}\label{eq:binary-curve-kloosterman}
 \boxed{\qquad N_n=q^n+1+K_n(1,b).\qquad}
\end{equation}
\end{proposition}
\begin{proof}
When \(x=0\), bijectivity of squaring implies that \(y^2=b^2\) has a unique solution \(y=b\).
When \(x\ne0\), put \(z=y/x\) to transform the equation into
\[
 z^2+z=x+\frac{b^2}{x^2}.
\]
By the Artin--Schreier trace criterion, the equation \(z^2+z=t\) has \(1+\psi_n(t)\) solutions: two when the trace is zero and none when it is one.
Since \(\Tr_{k_n/\F_2}(u^2)=\Tr_{k_n/\F_2}(u)\), the solution count for fixed \(x\ne0\) is
\[
 1+\psi_n\left(x+\frac{b^2}{x^2}\right)
 =1+\psi_n\left(x+\frac b x\right).
\]
Summing and including the one point with \(x=0\) and the point \(O\) gives
\[
 N_n=2+\sum_{x\in k_n^\times}\left(1+\psi_n(x+b/x)\right)
 =q^n+1+K_n(1,b).
\]
\end{proof}

The trace condition on Artin--Schreier fibers thus converts curve point counting into a one-variable Kloosterman sum.

\subsection{Recovering All Extension-Field Point Counts from One Count}

\begin{corollary}[A second-order recurrence for extension-field point counts]\label{cor:binary-curve-recurrence}
Put \(a=q+1-N_1\), and for \(n\ge1\) put \(A_n=q^n+1-N_n\), with the additional convention \(A_0=2\). Then
\begin{equation}\label{eq:binary-curve-recurrence}
 A_1=a,\qquad A_n=aA_{n-1}-qA_{n-2}\quad(n\ge2).
\end{equation}
In particular, \(N_1\) and \(q\) determine all \(N_n\), and
\[
 A_2=a^2-2q,\qquad N_2=(q+1)^2-a^2.
\]
\end{corollary}
\begin{proof}
By Proposition~\ref{prop:binary-curve-point-count}, we have \(A_n=-K_n(1,b)\) and \(a=-K_1(1,b)\).
Substitute these into Theorem~\ref{thm:kloosterman-recurrence} to obtain the recurrence; \(K_0=-2\) corresponds to \(A_0=2\).
Taking \(n=2\) gives the last two formulas.
\end{proof}

The value \(A_0\) is only an initial-value convention for the recurrence and involves no degree-zero extension.
All \(A_n,N_n\) are integers; the \(q\) in the recurrence is also an integer, not a field element of characteristic two.

\begin{example}[Recovering extension-field point counts from four base-field points]\label{exm:binary-curve-counts}
Take \(q=2\) and \(b=1\). Over \(k\), the curve has points
\[
 O,\qquad(0,1),\qquad(1,0),\qquad(1,1).
\]
Thus \(N_1=4\) and \(a=-1\), and the recurrence gives
\[
\begin{array}{c|rrrrrr}
 n&1&2&3&4&5&6\\ \hline
 A_n&-1&-3&5&1&-11&9\\
 N_n&4&8&4&16&44&56
\end{array}
\]
These values also follow from the character sums in Example~\ref{exm:binary-kloosterman}.
Since neither \(\F_4\) nor \(\F_8\) contains the other, \(N_3<N_2\) is not constrained by a subfield inclusion; finite-field inclusions are determined by divisibility of extension degrees.
\end{example}

\subsection{The Point-Count Generating Function and Frobenius Fixed Points}

\leavevmode\phantomsection\label{def:curve-zeta}\ignorespaces
Define the \textbf{zeta function} of \(E_b/k\) as the formal power series
\begin{equation}\label{eq:curve-zeta-definition}
 Z(E_b/k,u)=\exp\left(\sum_{n\ge1}N_n\frac{u^n}{n}\right)\in\Q[[u]].
\end{equation}
The exponential and logarithm here are formal operations; no convergence of a complex variable is involved.

\begin{corollary}[A rational generating expression for curve point counts]\label{cor:binary-curve-zeta}
With the notation above,
\begin{equation}\label{eq:binary-curve-zeta}
 \boxed{\qquad Z(E_b/k,u)=\frac{1-au+qu^2}{(1-u)(1-qu)}.\qquad}
\end{equation}
\end{corollary}
\begin{proof}
Substitute \(N_n=q^n+1+K_n(1,b)\) into the definition.
The contributions from \(q^n\) and \(1\) give \((1-qu)^{-1}\) and \((1-u)^{-1}\), respectively;
by~\eqref{eq:kloosterman-L}, the last contribution gives \(1+K_1(1,b)u+qu^2=1-au+qu^2\).
Multiplying the three expressions proves the result.
\end{proof}

Let \(\alpha,\beta\in\C\) be the two roots of \(T^2-aT+q\), counted with multiplicity. Then
\[
 \alpha+\beta=a,\qquad \alpha\beta=q,\qquad
 A_n=\alpha^n+\beta^n,\qquad
 N_n=q^n+1-\alpha^n-\beta^n.
\]
The power-sum identity follows from the same initial values and second-order recurrence.
For example, Example~\ref{exm:binary-curve-counts} gives
\[
 Z(E_1/\F_2,u)=\frac{1+u+2u^2}{(1-u)(1-2u)}.
\]

Coordinatewise Frobenius returns this point count to the language of fixed points.
For a finite extension \(k_s/k\), put
\[
 \pi_b(x,y)=(x^q,y^q),\qquad \pi_b(O)=O.
\]
Since the curve coefficients belong to \(k\), the map \(\pi_b\) permutes \(E_b(k_s)\). When \(n\mid s\),
\[
 \Fix\bigl(\pi_b^n:E_b(k_s)\to E_b(k_s)\bigr)=E_b(k_n),
\]
since an affine point is fixed exactly when both coordinates belong to \(k_n\), and \(O\) is included on both sides.
Thus \(Z(E_b/k,u)\) records the fixed-point counts of successive Frobenius powers.
Section~\ref{sec:orbit-euler} connects these fixed-point counts with Euler products over orbits.

\leavevmode\phantomsection\label{def:curve-frobenius-polynomial}\ignorespaces
In elliptic curve theory, \(P_b(T)=T^2-aT+q\) is the \textbf{characteristic polynomial of Frobenius}; see \cite[Chapter V]{SilvermanAEC}.
The character-sum calculation above gives the point-count difference sequence \((A_n)\) with annihilating polynomial \(P_b\).
Here \(\alpha,\beta\) are complex characteristic roots of a recurrence, not conjugates of an element of a finite field;
the permutation \(\pi_b\) of the curve points is also distinct from the earlier \(k\)-linear operators on extension fields.
The recurrence and rationality of the generating function alone do not imply \(\lvert\alpha\rvert=\lvert\beta\rvert=\sqrt q\); this absolute-value statement requires a separate proof.

\begin{remark}[The Hasse bound and recurrence roots]\label{v104:rem-hasse-bound}
The \textbf{Hasse bound} for elliptic curves states that every elliptic curve \(E/\F_q\) satisfies
\[
 \bigl|\#E(\F_q)-(q+1)\bigr|\le2\sqrt q;
\]
see \cite[Chapter V, Theorem 1.1]{SilvermanAEC} for a proof. Applying it to \(E_b\), whose nonsingularity was checked above, gives
\[
 |a|=|q+1-N_1|\le2\sqrt q,\qquad |K_1(1,b)|\le2\sqrt q.
\]
Since \(a\) is real, the discriminant of \(T^2-aT+q\) is nonpositive. Its roots are complex conjugates, or a repeated real root in the equality case.
Together with \(\alpha\beta=q\), this yields \(\lvert\alpha\rvert=\lvert\beta\rvert=\sqrt q\), and hence, for every \(n\ge1\),
\[
 |N_n-(q^n+1)|=|K_n(1,b)|\le2q^{n/2}.
\]
The additional input for the absolute-value bound is Hasse's theorem. The phase \(x+b/x\) remains rational, so the later theorem for polynomial phases does not apply directly.
\end{remark}

\subsection{\texorpdfstring{$p$}{p}-adic slopes of point-count recurrences and Newton polygons}
\label{v132:subsec-newton-curves}
The characteristic polynomial recovered from extension-field point counts also carries $p$-adic information about its roots. For $P(T)=\sum_{i=0}^d c_iT^i$ with $c_0c_d\ne0$, define its \textbf{Newton polygon} to be the lower convex hull of the points $(i,v_p(c_i))$; a zero coefficient is a point of infinite height. The negatives of its slopes, counted with their horizontal lengths, give the valuations of the roots in $\overline{\mathbb Q}_p$. The quadratic examples below can be verified directly from the roots and coefficients without invoking the general Newton-polygon theorem.

\begin{example}[Equal residue-field sizes but different Frobenius slopes]
\label{v132:ex-elliptic-slopes}
Consider the smooth elliptic curves over $\F_5$
\[
 E_1:\ y^2=x^3+x+1,\qquad E_2:\ y^2=x^3+1.
\]
Direct point counts, including the point at infinity, give $\#E_1(\F_5)=9$ and $\#E_2(\F_5)=6$. Their Frobenius traces are $a_1=5+1-9=-3$ and $a_2=5+1-6=0$, giving
\[
 P_1(T)=T^2+3T+5,\qquad P_2(T)=T^2+5.
\]
The two roots $\alpha,\beta$ of $P_1$ are algebraic integers. Their sum $-3$ is a $5$-adic unit and their product is $5$. Thus at least one root is a unit; because the valuations sum to one,
\[
 \{v_5(\alpha),v_5(\beta)\}=\{0,1\}.
\]
For $P_2$ each root satisfies $\alpha^2=-5$ and has valuation $1/2$. The Newton polygons pass through
\[
 (0,1),(1,0),(2,0),\qquad\text{and}\qquad(0,1),(2,0),
\]
respectively. These are the ordinary and supersingular slope patterns for elliptic curves. All complex Frobenius eigenvalues still have absolute value $\sqrt5$; the distinction is $5$-adic.
\end{example}

For general curves, lifting Frobenius to a $p$-adic setting provides an approach to zeta-function computations. Dwork proved the rationality of zeta functions of finite-field varieties by $p$-adic methods, while Kedlaya used Monsky--Washnitzer cohomology to compute point counts on hyperelliptic curves\cite{Dwork1960Zeta,Kedlaya2001Counting}. The elementary polygon above follows from the point-count recurrence and quadratic algebra; general cohomology and point-count algorithms require substantially more theory.

If the coefficients of the Frobenius characteristic polynomial are known to be integers with an explicit bound, sufficiently high $p$-adic precision determines them exactly. For instance, if $|c|\le B$ and $c\bmod p^s$ is known with $p^s>2B$, there is at most one integer in $[-B,B]$ having that residue. Precision can thus be converted into exact point counts, provided the modular coefficients can be computed by a suitable Frobenius-lifting algorithm.

\section{Weil bounds for polynomial character sums}
\label{v91:sec-weil-bounds}

Throughout this section, write \(k=\F_q\) and \(p=\Char k\).
The generating polynomial expresses the character sums over all extension fields as power sums of finitely many complex numbers, but does not yet bound their absolute values.
This section supplies the missing step. The argument again starts from root sets and Frobenius:
an auxiliary polynomial bounds trace and power-value fibers over extension fields, and power sums then bound the characteristic roots.
The auxiliary construction only needs a square-root estimate whose constant is independent of the extension degree;
the degree of the generating polynomial subsequently supplies the precise constant in the final bound.
The construction of high-multiplicity zeros by counting coefficient parameters is called \textbf{Stepanov's method}.
We use a form adapted to polynomial phases, proving nonvanishing and the multiplicity conditions separately\cite{CochranePinner}.

\par
\begin{figure}[!htbp]
\centering
\begin{tikzpicture}[every node/.style={font=\small\linespread{1.05}\selectfont},
 bx/.style={draw=Jade,fill=CoverMint,rounded corners,align=center,text width=5.5cm,inner sep=5pt},
 ar/.style={-{Stealth[length=2mm]},semithick,Jade}]
\node[bx] (ta) at (-3.3,0) {Trace-fiber construction\\Nonvanishing; Hasse conditions\\Lemma~\ref{v91:lem-trace-fiber-stepanov}};
\node[bx] (ma) at (3.3,0) {Multiplicative-fiber construction\\Nonvanishing and multiple zeros\\Lemma~\ref{v91:lem-multiplicative-nonzero}};
\node[bx] (tr) at (-3.3,-2.5) {Uniform rough bound\\in even extensions\\Corollary~\ref{v91:cor-additive-rough-bound}};
\node[bx] (mr) at (3.3,-2.5) {Uniform rough bound\\over extensions\\Lemma~\ref{v91:lem-multiplicative-rough}};
\node[bx,text width=9cm] (radius) at (0,-4.8) {Power-sum radius lemma~\ref{v91:lem-power-sum-radius}\\Every characteristic root has modulus \(\le\sqrt q\)};
\node[bx] (tl) at (-3.3,-7.1) {Additive generating polynomial\\\(D\le d-1\)\\Theorem~\ref{thm:additive-L-polynomial}};
\node[bx] (ml) at (3.3,-7.1) {Multiplicative generating polynomial\\\(D\le r-1\)\\Lemma~\ref{v91:lem-multiplicative-L}};
\node[bx,text width=9cm] (end) at (0,-9.6) {Combine degree and root-modulus bounds\\Additive: \((d-1)\sqrt q\); multiplicative: \((r-1)\sqrt q\)};

\draw[ar] (ta)--(tr); \draw[ar] (ma)--(mr);
\draw[ar] (tr)--(radius); \draw[ar] (mr)--(radius);
\draw[ar] (tl.north)--(radius.south west); \draw[ar] (ml.north)--(radius.south east);
\draw[ar] (tl)--(end); \draw[ar] (ml)--(end);
\draw[ar] (radius.south)--(end.north);
\end{tikzpicture}
\caption{The two proof routes to the Weil bounds. The upward arrows supply the power-sum representations. The rough constant may be large, but it must be independent of the extension degree.}\label{v103:fig-weil-dependencies}
\end{figure}
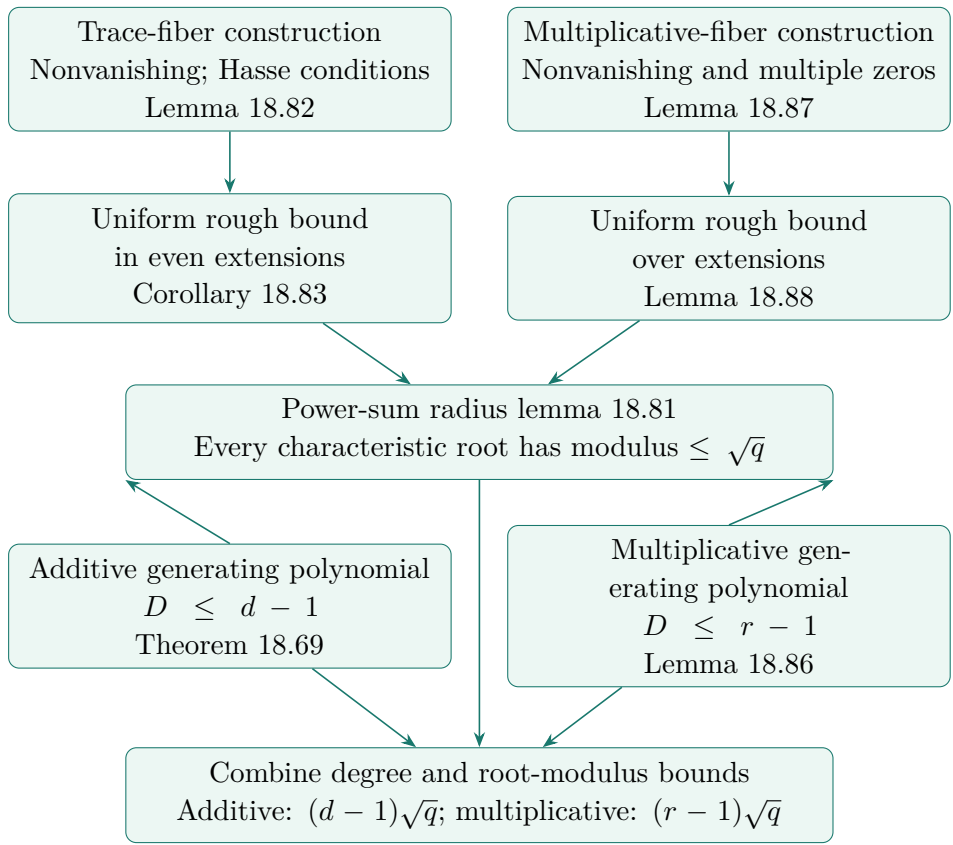

\Needspace{14\baselineskip}
Figure~\ref{v103:fig-weil-dependencies} separates two necessary inputs: the finite recurrence bounds the number of characteristic roots, and a rough bound uniform over extensions controls their moduli. The final constant follows only after both inputs are combined.
\FloatBarrier

\subsection{Bounding characteristic roots by extension-field power sums}

\begin{lemma}[Power sums and maximal modulus]\label{v91:lem-power-sum-radius}
Let \(\alpha_1,\ldots,\alpha_t\in\C\) and \(R>0\). If a constant \(C\) satisfies
\[
\left|\sum_{i=1}^t\alpha_i^n\right|\le CR^n
\]
for every sufficiently large integer \(n\), then \(|\alpha_i|\le R\) for every \(i\).
More generally, the same conclusion follows if, for some fixed positive integer \(h\),
\[
\left|\sum_{i=1}^t\alpha_i^{hn}\right|\le CR^{hn}
\]
for every sufficiently large \(n\).
In fact, for every positive integer \(h\), one has the exact formula
\begin{equation}\label{v130:eq-power-sum-radius}
\limsup_{n\to\infty}\left|\sum_{i=1}^t\alpha_i^{hn}\right|^{1/(hn)}
=\max_{1\le i\le t}|\alpha_i|.
\end{equation}
Both the empty sum and the maximum of the empty set are taken to be zero.
\end{lemma}
\begin{proof}
It suffices to prove the first assertion and then apply it to \(\alpha_i^h\) and \(R^h\).
Suppose that \(\rho=\max_i|\alpha_i|>R\). Write the distinct numbers of absolute value \(\rho\) as
\(\rho\beta_1,\ldots,\rho\beta_s\), with positive integer multiplicities \(m_1,\ldots,m_s\) and \(|\beta_j|=1\).
After division by \(\rho^n\), all remaining terms tend to zero, so the hypothesis implies
\(u_n=\sum_jm_j\beta_j^n\to0\). However, finite geometric sums give
\[
\lim_{N\to\infty}\frac1N\sum_{n=1}^N|u_n|^2
=\sum_{j=1}^s m_j^2>0.
\]
Indeed, the diagonal terms contribute \(m_j^2\), whereas every off-diagonal ratio
\(\beta_i\overline{\beta_j}\ne1\) has mean tending to zero. This contradicts \(u_n\to0\).
The upper bound in the exact formula follows from the triangle inequality. If its left-hand side were strictly less than \(\max_i|\alpha_i|\), choose \(R>0\) between them. The bound in the lemma would then hold for all sufficiently large \(n\), contradicting the conclusion just proved. If every \(\alpha_i\) is zero, or the list is empty, both sides are zero directly.
\end{proof}

Thus uniform estimates over all extensions of even degree already suffice.
We use this observation below to divide the Frobenius powers into two blocks of equal length.

\subsection{An auxiliary polynomial for trace fibers}

\begin{lemma}[A coarse square-root bound for trace fibers]\label{v91:lem-trace-fiber-stepanov}
Let \(f\in k[X]\), where \(k=\F_q\), \(d=\deg f\ge1\), and \(p\nmid d\).
Put \(E=\F_{q^{2u}}\), \(Q=q^{2u}\), and \(A=q^u\). If \(A\ge8dq^2\), then, for every \(b\in k\),
\[
N_b=\#\{x\in E:\Tr_{E/k}(f(x))=b\}
\le\frac Qq+8d(d+1)q\sqrt Q.
\]
\end{lemma}
\begin{proof}
Define
\[
H(X)=\sum_{t=u}^{2u-1}f(X)^{q^t},\qquad
B_b(X)=b-\sum_{t=0}^{u-1}f(X)^{q^t}.
\]
Every point \(x\) in the specified fiber satisfies \(H(x)=B_b(x)\). Moreover,
\[
\deg H=dQ/q,\qquad \deg B_b\le dA/q,
\]
and every exponent occurring in \(H\) is divisible by \(A\).
Choose the integers
\[
S=Q/q^2,\qquad
K=\left\lfloor\frac{A}{4dq^2}\right\rfloor,\qquad M=qK.
\]
The hypothesis ensures that \(K\ge1\) and \(M<A\). Consider
\begin{equation}\label{v91:eq-trace-auxiliary}
P(X)=\sum_{i=0}^{q-1}\sum_{j=0}^{K}
e_{ij}(X)H(X)^iX^{Qj},\qquad \deg e_{ij}<S,
\end{equation}
with \(q(K+1)S\) undetermined coefficients in the polynomials \(e_{ij}\in k[X]\).

First we show that a nonzero coefficient vector produces a nonzero polynomial.
The basic terms \(X^tH(X)^iX^{Qj}\), with \(0\le t<S\), have degrees
\[
t+idQ/q+jQ.
\]
If two such degrees agree, reduction modulo \(Q/q\), together with \(S<Q/q\), gives \(t=t'\).
It follows that \(d(i-i')+q(j-j')=0\). Since \(\gcd(d,q)=1\) and \(0\le i,i'<q\),
we obtain \(i=i'\) and \(j=j'\). The degrees are therefore distinct, so the basic terms are linearly independent.
Consequently, \(P\ne0\) whenever the coefficients are not all zero.

Next we arrange for every point of the fiber to be a zero of multiplicity at least \(M\).
Using the Hasse derivatives of Definition~\ref{ffs:def-hasse}, impose the following \(M\) polynomial identities:
\begin{equation}\label{v91:eq-trace-linear-conditions}
\sum_{i=0}^{q-1}\sum_{j=0}^{K}
\partial^{[r]}e_{ij}(X)B_b(X)^iX^j=0
\qquad(0\le r<M).
\end{equation}
Each left-hand side has degree less than \(S+K+dA\), so these identities impose at most
\(M(S+K+dA)\) homogeneous linear conditions. The parameter choices give
\[
K(K+dA)
\le \frac{A^2}{16d^2q^4}+\frac{A^2}{4q^2}
=S\left(\frac{1}{16d^2q^2}+\frac14\right)<S.
\]
Hence
\[
M(S+K+dA)=qK(S+K+dA)<q(K+1)S.
\]
There are fewer conditions than unknowns, so a nonzero solution exists; the preceding argument ensures that its polynomial \(P\) is nonzero.

If \(x\) belongs to the fiber, the Frobenius identities give, in a formal variable \(Y\),
\[
H(x+Y)\equiv H(x)\pmod{Y^A},\qquad
(x+Y)^{Qj}\equiv x^{Qj}=x^j\pmod{Y^Q}.
\]
For \(0\le r<M<A\), the coefficient of \(Y^r\) in \(P(x+Y)\) is therefore precisely
the value at \(x\) of the left-hand side of~\eqref{v91:eq-trace-linear-conditions}, because \(H(x)=B_b(x)\).
These coefficients vanish. The Hasse multiplicity criterion shows that \(x\) is a zero of multiplicity at least \(M\).
Adding the multiplicities at distinct points gives \(MN_b\le\deg P\). On the other hand,
\[
\deg P<S+dQ+KQ,
\]
and therefore
\[
N_b\le\frac Qq+\frac{dQ+S}{qK}.
\]
Since \(A\ge8dq^2\), we have \(K\ge A/(8dq^2)\), so the last term is at most
\[
8dq\bigl(d+q^{-2}\bigr)A\le8d(d+1)q\sqrt Q.
\]
\end{proof}

The three constraints in the proof have distinct roles: the different degrees prove nonvanishing,
the count of linear conditions proves existence, and Hasse derivatives establish the required multiplicities.
Together they reduce the fiber count to the root bound for a polynomial in one variable.

\begin{corollary}[A uniform estimate over extensions of even degree]\label{v91:cor-additive-rough-bound}
In the notation of Lemma~\ref{v91:lem-trace-fiber-stepanov}, let \(\psi\) be a nontrivial additive character of \(k\). Then
\[
\left|\sum_{x\in E}\psi\bigl(\Tr_{E/k}(f(x))\bigr)\right|
\le16d(d+1)q^2\sqrt Q
\]
for every sufficiently large \(u\).
\end{corollary}
\begin{proof}
Put \(B=8d(d+1)q\) and \(\delta_b=N_b-Q/q\). The fibers partition \(E\), so
\(\sum_b\delta_b=0\), while the lemma gives \(\delta_b\le B\sqrt Q\). Thus
\[
\sum_{b\in k}|\delta_b|
=2\sum_{\delta_b>0}\delta_b\le2qB\sqrt Q.
\]
Additive character orthogonality removes the \(Q/q\) terms, and the triangle inequality finishes the proof.
\end{proof}

\par\Needspace{14\baselineskip}
\subsection{The Weil bound for additive character sums}

\begin{theorem}[The Weil bound for polynomial additive character sums]\label{v91:thm-additive-weil}
Let \(\psi\) be a nontrivial additive character of \(k=\F_q\), and let \(f\in k[X]\) have degree \(d\ge1\), with \(p\nmid d\). Then
\[
\boxed{\left|\sum_{x\in k}\psi(f(x))\right|\le(d-1)\sqrt q.}
\]
More generally, for every \(r\ge1\),
\[
|S_r(f)|\le(d-1)q^{r/2}.
\]
\end{theorem}
\begin{proof}
By Theorem~\ref{thm:additive-L-polynomial} and Corollary~\ref{cor:character-sum-recurrence},
\[
L_f(u)=\prod_{i=1}^{D}(1-\omega_i u),\qquad
D\le d-1,\qquad S_r(f)=-\sum_{i=1}^{D}\omega_i^r.
\]
If \(D=0\), all the character sums vanish. Otherwise, Corollary~\ref{v91:cor-additive-rough-bound} gives
\[
\left|\sum_{i=1}^{D}(\omega_i^2)^u\right|\le Cq^u
\]
for every sufficiently large \(u\), with \(C\) independent of \(u\).
Applying Lemma~\ref{v91:lem-power-sum-radius} to \(\omega_i^2\) gives
\(|\omega_i^2|\le q\), and hence \(|\omega_i|\le\sqrt q\). Therefore
\[
|S_r(f)|\le\sum_{i=1}^{D}|\omega_i|^r\le Dq^{r/2}\le(d-1)q^{r/2}.
\]
\end{proof}

For an arbitrary polynomial \(f\), first include the additive character parameter and obtain the reduced phase \(R\) as in Theorem~\ref{v70:thm-phase-reduction}.
If \(R\) is nonconstant, its degree \(d\) is prime to \(p\), and the theorem gives the bound in terms of this effective degree.
If \(R\) is constant, the character sum can have absolute value \(q^r\), and the exact formula of Corollary~\ref{v70:cor-reduced-phase-L-series} applies instead.
Thus phase reduction identifies the nondegeneracy condition, while counting high-multiplicity zeros establishes square-root cancellation.

\subsection{Multiplicative character sums: finite recurrences from residue-class algebras}

We now prove the general bound for multiplicative character sums. Compared with the additive case, norms take the place of traces,
and summation over residue classes in the Chinese remainder theorem replaces coefficient cancellation in Newton's identities.
All multiplicative characters are still assigned the value zero at zero, including those that are trivial on the multiplicative group.

\begin{theorem}[The Weil bound for multiplicative character sums]\label{v91:thm-multiplicative-weil}
Let $\chi$ be a multiplicative character of $\F_q^\times$ of order $n>1$,
and let $f\in\F_q[X]$ be a nonconstant polynomial that cannot be written as
$f=a h^n$, where $a\in\F_q^\times$ and $h\in\F_q[X]$.
Let $r$ be the number of distinct roots of $f$ in a splitting field, so that
$r=\deg\operatorname{rad}(f)$. Then
\[
\left|\sum_{x\in\F_q}\chi(f(x))\right|\le(r-1)\sqrt q.
\]
\end{theorem}

We first encode the sums over all extension fields in a single polynomial. Write
\[
\chi_m=\chi\circ\Nm_{\F_{q^m}/\F_q},\qquad
S_m=\sum_{x\in\F_{q^m}}\chi_m(f(x)).
\]
Surjectivity of the norm ensures that $\chi_m$ still has order $n$.

\begin{lemma}[The generating polynomial for multiplicative character sums]\label{v91:lem-multiplicative-L}
Under the hypotheses of Theorem~\ref{v91:thm-multiplicative-weil}, there is a polynomial $L(T)$ with complex coefficients,
degree $D\le r-1$, and constant term one such that
\[
L(T)=\exp\left(\sum_{m\ge1}S_m\frac{T^m}{m}\right).
\]
Consequently, there are $\alpha_1,\ldots,\alpha_D\in\C$ such that
$S_m=-\sum_{j=1}^D\alpha_j^m$ for every $m\ge1$.
\end{lemma}
\begin{proof}
Write
\[
f=a\prod_{i=1}^sP_i^{e_i},\qquad
B=\prod_{i=1}^sP_i,\qquad \deg B=r,
\]
where the $P_i$ are distinct monic irreducible polynomials.
For a monic polynomial $G$ with roots $\beta_1,\ldots,\beta_k$, counted with multiplicity, define
\[
\lambda(G)=\chi\left(\prod_{j=1}^k f(\beta_j)\right),
\qquad \lambda(1)=1.
\]
The product in parentheses is fixed by Frobenius and hence belongs to $\F_q$.
Taking unions of root multisets gives $\lambda(GH)=\lambda(G)\lambda(H)$.

Choose a root $\gamma_i$ of each $P_i$, and put $d_i=\deg P_i$.
Interchanging the two products over roots gives
\[
\lambda(G)=c^{\deg G}\eta(G),\qquad
c=\chi\bigl(a(-1)^{\deg f}\bigr),
\]
where, if $(G,B)=1$,
\[
\eta(G)=\prod_{i=1}^s
\chi\left(\Nm_{\F_{q^{d_i}}/\F_q}(G(\gamma_i))\right)^{e_i},
\]
and otherwise $\eta(G)=0$.
The Chinese remainder theorem decomposes $(\F_q[X]/(B))^\times$ as the direct product of the groups
$\F_{q^{d_i}}^\times$, so $\eta$ is a multiplicative character of this unit group.
At least one $e_i$ is not divisible by $n$. Surjectivity of the norm makes the character on that factor nontrivial,
and hence $\eta$ is nontrivial.

For $k\ge r$, each residue class modulo $B$ has exactly $q^{k-r}$ monic representatives of degree $k$.
Nontrivial character orthogonality therefore gives
\[
\sum_{\substack{G\text{ monic}\\\deg G=k}}\lambda(G)
=c^kq^{k-r}\sum_{u\in\F_q[X]/(B)}\eta(u)=0.
\]
Thus the formal power series
\[
L(T)=\sum_{G\text{ monic}}\lambda(G)T^{\deg G}
\]
is a polynomial of degree at most $r-1$. Unique factorization gives the Euler product
\[
L(T)=\prod_{P\text{ monic irreducible}}
\bigl(1-\lambda(P)T^{\deg P}\bigr)^{-1}.
\]
If the degree $d$ of $P$ divides $m$, then any root $\beta$ of $P$ in $\F_{q^m}$ satisfies
\[
\chi_m(f(\beta))=\lambda(P)^{m/d}.
\]
Each such $P$ contributes $d$ conjugate roots. Taking the formal logarithmic derivative of the Euler product gives
\[
T\frac{L'(T)}{L(T)}=\sum_{m\ge1}S_mT^m.
\]
The constant term one then determines the stated exponential formula.
Finally, factor $L$ over $\C$ as $\prod_{j=1}^D(1-\alpha_jT)$ and compare coefficients in its logarithmic derivative.
\end{proof}

\subsection{An auxiliary polynomial for multiplicative fibers}

The degree of the generating polynomial has already determined the constant in the final bound.
It remains to establish a square-root bound uniform in the extension degree.
We use only polynomials over finite fields, Hasse derivatives, norms, and linear equations.

\begin{lemma}[Nonvanishing of the auxiliary polynomial]\label{v91:lem-multiplicative-nonzero}
Let $K=\F_Q$, let $n>1$ divide $Q-1$, and let $g\in K[X]$ have degree $l\ge1$ and satisfy $g(0)\ne0$.
Suppose that the multiplicities of the roots of $g$ in a splitting field have greatest common divisor one.
Put $\rho=(Q-1)/n$. If a positive integer $S$ satisfies
\[
n(S-1)+(n-1)l<Q,
\]
then, for any finite family of polynomials $e_{jk}$ of degree less than $S$,
\[
F(X)=\sum_{j=0}^{n-1}\sum_{k=0}^{K_0}
 e_{jk}(X)g(X)^{\rho j}X^{Qk}
\]
is identically zero if and only if all $e_{jk}$ are zero.
\end{lemma}
\begin{proof}
Only the forward implication needs proof. Suppose that $F=0$ but some $e_{jk}\ne0$.
After removing the smallest common factor $X^{Qk}$ from all the terms, we may assume that some $e_{j0}\ne0$.
Write
\[
E_j(X)=\sum_ke_{jk}(X)X^{Qk},\qquad
U(X,Y)=\sum_{j=0}^{n-1}E_j(X)Y^j.
\]
Take the product over all $n$th roots of unity, which belong to $K$.
Invariance under $Y\mapsto\zeta Y$ gives
\[
\prod_{\zeta^n=1}U(X,\zeta Y)
=\sum_{s=0}^{n-1}C_s(X)Y^{ns}.
\]
Each $C_s$ is a homogeneous expression of degree $n$ in $E_0,\ldots,E_{n-1}$, and therefore has the form
\[
C_s(X)=\sum_{k\ge0}A_{sk}(X)X^{Qk},
\qquad \deg A_{sk}\le n(S-1).
\]
Here $A_{s0}$ is precisely the coefficient obtained by retaining only the polynomials $e_{j0}$.
After substituting $Y=g(X)^\rho$, one factor on the left is $F=0$.
Multiply by $g^{n-1}$ and use $g(X)^Q\equiv g(0)^Q\pmod{X^Q}$ to obtain
\[
P(X):=\sum_{s=0}^{n-1}
A_{s0}(X)g(X)^{n-1-s}g(0)^{Qs}
\equiv0\pmod{X^Q}.
\]
Since $\deg P\le n(S-1)+(n-1)l<Q$, it follows that $P=0$.

Let $L/K$ be a finite splitting field of $g$, put $E=L(X)$, and choose $y$ in an extension of $E$ satisfying
\[
y^n=\frac{g(0)^Q}{g(X)}.
\]
This element has degree exactly $n$ over $E$. Indeed, write its degree as $t\le n$.
Using the determinant of multiplication as the norm of a finite extension gives
\[
\Nm_{E(y)/E}(y)^n
=\left(\frac{g(0)^Q}{g(X)}\right)^t.
\]
If $\gamma$ is a root of $g$ of multiplicity $m_\gamma$, comparison of the exponents of $X-\gamma$ in these rational functions gives
$n\mid t m_\gamma$.
The greatest common divisor of the $m_\gamma$ is one, so $n\mid t$, and hence $t=n$.

On the other hand, the $k=0$ part of the preceding product identity gives
\[
0=\frac{P(X)}{g(X)^{n-1}}
=\prod_{\zeta^n=1}
\left(\sum_{j=0}^{n-1}e_{j0}(X)(\zeta y)^j\right).
\]
One factor must vanish, contradicting the fact that the minimal polynomial of $y$ has degree $n$.
\end{proof}

\begin{lemma}[A uniform coarse square-root bound]\label{v91:lem-multiplicative-rough}
Let $K=\F_Q$, let $\chi$ be a multiplicative character of order $n>1$,
and let $g\in K[X]$ have degree $l\ge1$.
Suppose that the multiplicities of its roots in a splitting field have greatest common divisor one. Then
\[
\left|\sum_{x\in K}\chi(g(x))\right|
\le8n^2(l+1)^2\sqrt Q.
\]
The constant on the right depends only on $n,l$ and is independent of the size of the field.
\end{lemma}
\begin{proof}
Put $A=\sqrt Q$. If $A<8n(l+1)$, the result follows from the trivial bound $Q$.
We may therefore assume that $A\ge8n(l+1)$.
Since $Q>l$, a translation of the variable makes $g(0)\ne0$.
Put
\[
\rho=\frac{Q-1}{n},\qquad
N_a=\#\{x\in K:g(x)\ne0,\ g(x)^\rho=a\}
\quad(a\in\mu_n(K)).
\]
Choose
\[
S=\left\lfloor\frac{Q}{2n}\right\rfloor,
\quad K_0=\left\lfloor\frac{A}{4n(l+1)}\right\rfloor,
\quad R=nK_0.
\]
These parameters satisfy
\[
S\ge\frac{Q}{4n},\qquad
K_0\ge\frac{A}{8n(l+1)},\qquad R<Q,
\qquad n(S-1)+(n-1)l<Q.
\]
Fix $a\in\mu_n(K)$, and regard the coefficients of the $e_{jk}$ in the preceding lemma's auxiliary polynomial as unknowns.

We first record a degree bound for Hasse derivatives. For any integer $b\ge0$, polynomial $e$, and $r\ge0$, the expression
\[
H_r(X)=g(X)^r g(X)^{-b}\partial^{[r]}(e(X)g(X)^b)
\]
is a polynomial and satisfies
\[
\deg H_r\le\deg e+r(l-1)
\]
whenever it is nonzero. To verify this, expand
\[
e(X+T)g(X+T)^b
=e(X+T)g(X)^b
\left(1+\frac{g(X+T)-g(X)}{g(X)}\right)^b
\]
in powers of $T$. The coefficient of $T^r$ contains at most $r$ factors $g(X)$ in its denominator.
Suppose that a term uses $t$ positive-degree Taylor terms of $g$, with total $T$-degree $s$.
After multiplication by $g^r$, its degree is at most
\[
\deg e-(r-s)+(r-t)l+tl-s
=\deg e+r(l-1).
\]
This proves both polynomiality and the degree bound.

Define $H_{jkr}$ by taking $e=e_{jk}$ and $b=\rho j$.
Since $r<R<Q$, the positive-order Hasse derivatives of $X^{Qk}$ vanish at these orders.
At a point $x$ in the chosen fiber, we have $x^Q=x$ and $g(x)^\rho=a$, so
\[
g(x)^r\partial^{[r]}F(x)=
\sum_{j,k}a^jH_{jkr}(x)x^k.
\]
Thus every point of the fiber is a zero of $F$ of multiplicity at least $R$ provided that
\[
\sum_{j,k}a^jH_{jkr}(X)X^k
\]
vanishes identically for every $0\le r<R$.
These conditions are homogeneous linear equations in the unknown coefficients.
There are $n(K_0+1)S$ unknowns, whereas the number of equations is at most
\[
R(S+K_0)+(l-1)\frac{R(R-1)}2.
\]
The latter number is strictly smaller than the former, because
\[
K_0^2\left(1+\frac{n(l-1)}2\right)
\le\frac{Q}{32n(l+1)}<S.
\]
There is therefore a solution whose coefficients are not all zero, and Lemma~\ref{v91:lem-multiplicative-nonzero}
ensures that the corresponding polynomial $F$ is nonzero.

Comparing the number of zeros, counted with multiplicity, with the degree gives
\[
N_a\le\frac{\deg F}{R}
\le\frac{QK_0+(n-1)l\rho+S-1}{nK_0}
\le\frac Qn+8(l+1)^2\sqrt Q.
\]
If $z\le l$ is the number of roots of $g$ in $K$, then
$\sum_{a\in\mu_n(K)}N_a=Q-z$.
Applying the upper bound to the other $n-1$ fibers gives
\[
N_a\ge\frac Qn-8(n-1)(l+1)^2\sqrt Q-l.
\]
The map $t\mapsto t^\rho$ and the character $\chi$ have the same kernel, so there is a nontrivial complex character $\vartheta$ of
$\mu_n(K)$ such that $\chi(t)=\vartheta(t^\rho)$.
Since $\sum_a\vartheta(a)=0$,
\begin{align*}
\left|\sum_{x\in K}\chi(g(x))\right|
&=\left|\sum_{a\in\mu_n(K)}\vartheta(a)(N_a-Q/n)\right|\\
&\le\bigl(8n(n-1)(l+1)^2+nl\bigr)\sqrt Q\\
&\le8n^2(l+1)^2\sqrt Q.
\end{align*}
\end{proof}

\begin{proof}[Proof of Theorem~\ref{v91:thm-multiplicative-weil}]
Again write $f=a\prod_iP_i^{e_i}$, and put
\[
e=\gcd_i e_i,\qquad h=\prod_iP_i^{e_i/e}.
\]
Finite fields are perfect, so each $P_i$ has distinct roots.
The multiplicities of the roots of $h$ in a splitting field therefore have greatest common divisor one.
Moreover, $f=ah^e$ and $\chi^e$ is nontrivial: otherwise $n\mid e$, so $n$ would divide every $e_i$, contrary to the hypothesis.

After norm lifting, $\chi_m^e$ still has the same order as $\chi^e$.
Apply Lemma~\ref{v91:lem-multiplicative-rough} over each $\F_{q^m}$, and note that $|\chi_m(a)|=1$. We obtain
\[
|S_m|\le Cq^{m/2}\qquad(m\ge1),
\]
where $C$ depends only on the original $f,\chi$ and not on $m$.
Lemma~\ref{v91:lem-multiplicative-L} expresses $S_m$ as a power sum of at most $r-1$ complex numbers.
Lemma~\ref{v91:lem-power-sum-radius} then bounds the absolute value of each characteristic root by $\sqrt q$,
and hence $|S_1|\le(r-1)\sqrt q$.
\end{proof}

\begin{remark}
The exclusion of $f=ah^n$ is necessary. In that case,
\[
\sum_{x\in\F_q}\chi(f(x))
=\chi(a)\bigl(q-\#\{x\in\F_q:h(x)=0\}\bigr),
\]
which is generally of order $q$.
The auxiliary polynomial is only needed to establish the square-root scale; its coarse constant disappears from the final estimate.
The precise constant $r-1$ comes from the degree of the generating polynomial, determined by the residue-class algebra.
This identifies the separate roles of degree control and size control in the argument.
\end{remark}

This section proves bounds for polynomial additive phases and polynomial multiplicative character sums. Kloosterman sums involve a rational phase; their earlier recurrence does not acquire an absolute-value bound merely from these polynomial estimates.

\par\Needspace{11\baselineskip}
\section*{Exercises}
\phantomsection
\addcontentsline{toc}{section}{Exercises}

\begin{exercise}[Low-degree power-value rigidity and the square-root boundary]\label{v111:ex-low-degree-power}
Let $v\ge2$ divide $q-1$, and let a nonconstant $F\in\F_q[X]$ satisfy
\[
F(\F_q)\subseteq\{0\}\cup(\F_q^\times)^v.
\]
\begin{enumerate}
\item Use Theorem~\ref{v91:thm-multiplicative-weil} to prove that $\deg F<\sqrt q$ implies $F=G^v$ for some $G\in\F_q[X]$.
Treat $F=cH^v$ separately, explaining how the power-value condition determines the constant $c$.
\item Over $\F_9$, take $F=X^3+X$ and $v=4$. Prove that $F$ is a minimal value set polynomial with three values but is not a fourth power. Thus the strict degree bound cannot be replaced by $\deg F\le\sqrt q$.
\item Extend the boundary example to $q=s^2$, with $s$ any prime power: take $F=X^s+X$ and $v=s+1$, and determine every nonempty fiber size. Explain why the assertion concerns powers as polynomials, rather than a pointwise choice of roots.

\item Strengthen the first assertion by replacing \(\deg F<\sqrt q\) with \(\rho<\sqrt q\), where \(\rho\) is the number of distinct roots of \(F\) in a splitting field. No upper bound on \(\deg F\) is imposed in this version. Show more precisely that if \(F\) is not a \(v\)th power, then
\[
q-z\le(\rho-1)\sqrt q,
\]
where \(z\) is its number of zeros in \(\F_q\). Explain why the trace-polynomial examples also attain the boundary \(\rho=\sqrt q\).
\end{enumerate}
\end{exercise}

\begin{exercise}[The lower bound for primitive elements of nonzero trace]\label{v88:ex-primitive-nonzero-bound}
For \((q,n)=(2,3)\) and \((3,2)\), compute
\(Q,M,R,\vartheta(M),W(R)\), and use
\eqref{v88:eq-primitive-nonzero-trace-lower} to verify the existence
of primitive elements of nonzero trace. Explain why the bound does not
require separate counts for each nonzero trace value.
\end{exercise}

\begin{exercise}[The sixth moment of cubic-phase character sums]\label{v58:ex-cubic-sixth-moment}
Let \(K=\F_Q\) have characteristic \(p>3\), and let \(\psi\) be a nontrivial additive character. Prove
\[
\sum_{a,b,c\in K}\left|\sum_{t\in K}\psi(at+bt^2+ct^3)\right|^6
=6Q^6-9Q^5+4Q^4.
\]
\end{exercise}

\begin{exercise}\label{ex:fourier-delta}
Find the Fourier transforms of the point indicator \(1_{\{b\}}\) and the constant function one, and recover orthogonality from Fourier inversion.
\end{exercise}
\begin{exercise}\label{ex:warning-quadratic}
Let \(q\) be odd and \(a,b\in\F_q^\times\). Use a count of squares to prove that \(aX^2+bY^2=1\) is solvable,
and compare with the conclusion of the Chevalley--Warning theorem for the ternary homogeneous equation \(aX^2+bY^2-Z^2=0\).
\end{exercise}

\begin{exercise}[Norm fibers and quadratic-form counts]\label{ex:quadratic5}
Without enumerating pairs, count the solutions of \(x^2+y^2=1\) and \(x^2+y^2=0\) in \(\F_5\).
Now let \(q\) be odd and \(a\in\F_q^\times\) be a nonsquare. Use both the norm fibers of \(\F_q(\sqrt a)/\F_q\) and the value distribution of rank-two quadratic forms to count the solutions of \(x^2-ay^2=b\) for every \(b\in\F_q\), and explain how \(q+1\) corresponds in the two counts.
\end{exercise}
\begin{exercise}\label{ex:correlation7}
Replace the binary period-seven sequence of Chapter~\ref{chap:sequences} by \(a_t=(-1)^{s_t}\).
Prove that its periodic autocorrelation \(\sum_{t=0}^6a_{t+\tau}a_t\) is seven when \(7\mid\tau\) and minus one otherwise.
\end{exercise}

\begin{exercise}[The Fourier transform of a subspace indicator]\label{ex:wu-fourier-subspace}
Let \(U\) be an \(\F_p\)-subspace of \(K=\F_{p^m}\), and let
\(U^\perp=\{a:\Tr_{K/\F_p}(au)=0\ \text{for all }u\in U\}\).
Prove that, with the Fourier convention of this chapter,
\[
\widehat{1_U}(a)=|U|\,1_{U^\perp}(a).
\]
Use this to find the Fourier transform of the indicator of the affine coset \(b+U\).
\end{exercise}
\par

\begin{exercise}[The frequency action of an invertible linear substitution]\label{ex:wu2-fourier-change}
Let \(A:K\to K\) be an invertible \(\F_p\)-linear map. Prove \(\widehat{F\circ A}(a)=\widehat F((A^*)^{-1}(a))\).
Apply this to \(A(x)=x^p\) and write the transformation of the frequency parameter.
\end{exercise}

\begin{exercise}[The adjoint formula for a quotient phase]\label{v70:ex-quotient-phase}
Let \(K=\F_{p^m}\), let \(\psi(z)=\zeta_p^{\Tr_{K/\F_p}(z)}\), and let \(L:K\to K\) be \(\F_p\)-linear, with adjoint \(L^*\) relative to the absolute trace pairing. Put \(U=\im L\) and \(W=\ker L^*\), and take any function \(g:K\to K\). For \(b\in K\), write
\[
S(b)=\sum_{x\in K}\psi(g(Lx)+bx).
\]
Prove that \(S(b)=0\) if \(b\notin\im L^*\), whereas if \(b=L^*c\), then
\[
S(b)=|\ker L|\sum_{y\in U}\psi(g(y)+cy)
=\sum_{t\in W}\sum_{y\in K}\psi(g(y)+(c+t)y).
\]
Show that the right-hand side is independent of the choice of \(c\).
\end{exercise}

\begin{exercise}[The discriminant of a trace quadratic form]\label{ex:wu3-trace-square-odd}
Let \(q\) be odd, with \(E=\F_{q^n}\) and \(Q(x)=\Tr_{E/\F_q}(x^2)\).
Prove that \(Q\) is nondegenerate and that the determinant \(\Delta\) of its quadratic coefficient matrix satisfies \(\eta(\Delta)=(-1)^{n-1}\).
Use this to count the solutions of \(Q(x)=b\) for every value. Compute explicitly the case \(q=3,n=2\).
\end{exercise}
\begin{exercise}[Zero polar form does not mean zero quadratic form]\label{ex:wu3-trace-square-even}
Let \(q\) be even and \(E=\F_{q^n}\). Prove that \(Q(x)=\Tr_{E/\F_q}(x^2)\) has zero polar form,
but \(Q\) is surjective and every value has \(q^{n-1}\) preimages.
When \(q=2\), find the Fourier support of the phase function \((-1)^{Q(x)}\),
and explain why it differs from the orthogonal complement of the polar radical.
\end{exercise}
\begin{exercise}[Adding a linear term to a trace quadratic form]\leavevmode\label{ex:wu3-quadratic-linear-shift}
In Example~\ref{exm:f64-quadratic}, let \(Q(x)=\Tr(x^3)\) and \(A=\sigma+\sigma^{-1}\).
For \(b\in\F_{64}\) and \(c\in\F_2\), count the solutions of \(Q(x)+\Tr(bx)=c\).
Distinguish the cases \(b\in\im A\) and \(b\notin\im A\), and express this condition using a relative trace.
\end{exercise}
\begin{exercise}[The boundary cases of the Jacobi lifting formula]\label{ex:wu3-jacobi-boundary}
Prove that the identity in Corollary~\ref{cor:jacobi-lifting} holds whenever at least one of \(\chi,\lambda\) is nontrivial,
including the case \(\chi\lambda=1\). Explain why the case of two trivial characters cannot be included.
\end{exercise}
\begin{exercise}[A cubic diagonal equation]\label{ex:wu3-diagonal-seven}
Over \(\F_7\), count the solutions of \(x^3+y^3=c\) for every \(c\).
Take a cubic character \(\chi(3)=\zeta\), verify \(J(\chi,\chi)=2+3\zeta^2\), and recover the counts using Proposition~\ref{prop:diagonal-jacobi-count}.
\end{exercise}

\begin{exercise}[The degree condition for generating polynomials and trace equations]\label{ex:character-L-boundary}
Take \(k=\F_2\) and \(\psi(t)=(-1)^t\).
\begin{enumerate}
\item For \(f=X^3\), compute only \(S_1(f)\) and \(S_2(f)\) to determine \(L_f(u)\) and all \(S_r(f)\),
and count the solutions of \(\Tr_{\F_{2^r}/\F_2}(x^3)=0\).
Check that the result for \(r=6\) agrees with Example~\ref{exm:f64-quadratic}.
\item For \(g=X^2+X\), compute all \(S_r(g)\) and \(L_g(u)\).
Explain why \(p\nmid\deg f\) cannot simply be omitted from Theorem~\ref{thm:additive-L-polynomial}.
\item For \(h=X^6+X^3+X\), find the reduced phase, all \(S_r(h)\), and \(L_h(u)\). Explain why a polynomial whose original degree is divisible by the characteristic may still have a polynomial \(L\)-series after reduction.
\end{enumerate}
\end{exercise}

\begin{exercise}[The zero fiber and the limits of square-root bounds]\leavevmode\phantomsection\label{ex:diagonal-zero-fiber}
Retaining the notation of Theorem~\ref{thm:diagonal-weil-bound}, prove
\[
N_0=q^{n-1}-(q-1)
\sum_{\substack{\chi_i\in\mathcal X_{d_i}\setminus\{\varepsilon\}\\
                 \chi_1\cdots\chi_n=\varepsilon}}
\left(\prod_i\chi_i(a_i^{-1})\right)J(\chi_1,\ldots,\chi_n),
\]
and hence \(\lvert N_0-q^{n-1}\rvert\le(q-1)M_0q^{(n-2)/2}\).
Let \(q\) be odd. Use \(x^2-y^2=0\) to show why the estimate for a nonzero right-hand side cannot be applied directly to a zero right-hand side.
Then take all \(d_i=2\) and use the parity of the product character to explain why the nondegenerate quadratic forms in Theorem~\ref{thm:quadratic-odd-count}
have different error scales in odd and even dimensions, and recover their exact solution counts.
\end{exercise}

\begin{exercise}[The quadratic obstruction to primitive elements of zero trace]\label{ex:primitive-zero-trace-boundary}
Let \(E=\F_{q^2}\) and \(k=\F_q\). Prove that if \(x\in E^\times\) and \(\Tr_{E/k}(x)=0\),
then \(\ord(x)\mid2(q-1)\), so \(x\) is not a primitive element of \(E\).
In characteristic two, further determine the zero-trace fiber.
Then prove that~\eqref{eq:primitive-trace-zero-criterion} cannot hold when \(n=2\),
showing that this sufficient condition is consistent with the structural obstruction above.
\end{exercise}

\begin{exercise}[Extension-field point counts after adding a quadratic term]\leavevmode\phantomsection\label{ex:curve-quadratic-term}
Retain the notation of Section~\ref{sec:kloosterman-curve-points}, take \(c\in k\), and put
\[
 E_{b,c}(k_n)=\{(x,y)\in k_n^2:y^2+xy=x^3+cx^2+b^2\}\ \sqcup\ \{O\},
 \qquad \varepsilon=(-1)^{\Tr_{k/\F_2}(c)}.
\]
Prove
\[
 \#E_{b,c}(k_n)=q^n+1+\varepsilon^nK_n(1,b).
\]
If \(a=q+1-\#E_b(k)\), define \(Z(E_{b,c}/k,u)\) by the exponential generating expression for point counts and prove
\[
 Z(E_{b,c}/k,u)=\frac{1-\varepsilon au+qu^2}{(1-u)(1-qu)}.
\]
Compare the point counts of \(E_{b,c}\) and \(E_b\) over even-degree extensions, and find the first four extension-field point counts when \(q=2\) and \(b=c=1\).
\end{exercise}

\begin{exercise}[Root power sums, unions of orbits, and characteristic restrictions]
\label{v40:ex:root-power-fourier}\leavevmode\par\nobreak
\begin{enumerate}
\item Over \(\F_2\), take \(N=15\) and a primitive 15th root of unity \(\zeta\). For \(T=\{0,5,10\}\), compute all \(s_j=\sum_{a\in T}\zeta^{aj}\) for \(0\le j<15\), and use Proposition~\ref{v40:prop:root-power-composition} to recover the additive stabilizer and the largest exponent in a monomial composition of \(P_T\). Compare with the constituent cosets \(C_0\) and \(C_5\) separately.
\item In \(\F_9\), take an element \(\xi\) of order four. For \(T'=\{0,1,2\}\subseteq\Z/4\Z\), compute \(s'_j=\sum_{a\in T'}\xi^{aj}\) for \(0\le j<4\). Explain why \(s'_0=0\) does not affect determination of the additive stabilizer from the frequencies with nonzero coefficients.
\item Determine whether \(T'\) is stable under multiplication by three and whether \(P_{T'}(X)=\prod_{a\in T'}(X-\xi^a)\) belongs to \(\F_3[X]\). Explain the respective hypotheses needed for the cyclic Fourier stabilizer criterion and descent of the root product.
\end{enumerate}
\end{exercise}

\begin{exercise}[Distinct-element sums in a multiplicative subgroup and Jacobi sums]\label{v85:ex-restricted-sums}
Let \(H\le\F_q^\times\) have index \(d\), and let \(H^\perp\) be the group of multiplicative characters trivial on \(H\).
Extend every multiplicative character, including the trivial one, by zero at zero. Put
\[
 S_H(s)=\#\{(x,y)\in H^2:x\ne y,\ x+y=s\}.
\]
\begin{enumerate}
\item For \(s\ne0\), prove
\[
 S_H(s)=\frac1{d^2}\sum_{\chi,\lambda\in H^\perp}
 (\chi\lambda)(s)J(\chi,\lambda)-\epsilon_s,
\]
where \(\epsilon_s=1\) if the characteristic is not two and \(s/2\in H\), and \(\epsilon_s=0\) otherwise.
\item Using Proposition~\ref{prop:jacobi}, prove
\[
 S_H(s)\ge
 \frac{q+1-3d-(d-1)(d-2)\sqrt q}{d^2}-\epsilon_s\qquad(s\ne0).
\]
Separate the cases of two trivial characters, exactly one trivial character, inverse nontrivial characters, and the remaining pairs.
\item Determine \(S_H(0)\), and prove \(S_H(hs)=S_H(s)\) for \(h\in H\). If \(q=p^r\), prove also \(S_H(s^p)=S_H(s)\).
\end{enumerate}
\end{exercise}

\begin{exercise}[Products of distinct elements and constant terms]\label{v85:ex-distinct-products}
Let \(H\) be a multiplicative subgroup of a finite field, of order \(n\ge2\).
\begin{enumerate}
\item For every \(1\le t<n\), prove \(\{\prod_{a\in S}a:S\in\binom Ht\}=H\).
Write \(H=\langle\gamma\rangle\), put \(g=\gcd(t,n)\), and consider the exponent sets
\(\{0,1,\ldots,t-2,t-1+j\}\), \(0\le j<g\), and their translates. For \(t=1\), the initial consecutive part is empty.
\item Let \(n\ge3\), \(\beta\in H\), and put \(\nu_{H,3}(\beta)=\#\{x\in H:x^3=\beta\}\). Prove
\[
 \#\{(x,y,z)\in H^3:x,y,z\text{ pairwise distinct},\ xyz=\beta\}
 =n^2-3n+2\nu_{H,3}(\beta),
\]
and characterize \(\nu_{H,3}(\beta)\) using \(H^3=\{h^3:h\in H\}\).
\item Explain why \textup{(1)} says exactly that the constant terms of monic degree-\(t\) polynomials with distinct roots in \(H\) range over \((-1)^tH\).
Why must the case \(t=n\) be excluded from \textup{(1)}?
\end{enumerate}
\end{exercise}

\begin{exercise}[Short relations versus multiplicative orders]\label{v95:ex-short-order}
Let $H\le\F_{q^k}^\times$ have order $n$, with $k=\ord_n(q)\ge2$, and use the notation of \eqref{v95:eq-lambda-definition}.
\begin{enumerate}
\item Prove directly that $\lambda_q(H)=2\iff\gcd(n,q-1)>1$.
\item If $\lambda_q(H)=k+1$, prove $n\mid\Phi_k(q)$, $n\equiv1\pmod k$, and $n\ge(k+1)^{\Omega(n)}$.
\item For $q=2,n=7$, explain why all nonzero cyclotomic cosets have length three but $\lambda_2(H)=3$.
\end{enumerate}
\end{exercise}

\begin{exercise}[Cyclotomic counting of dependent triples]\label{v95:ex-triples}
Suppose $\gcd(n,q-1)=1$, $k=\ord_n(q)\ge3$, and put $U=\F_q^\times H$.
\begin{enumerate}
\item Prove $\mathcal C(U)\ge q-2$ and explain the origin of the $q-2$ scalar solutions.
\item Without invoking Theorem~\ref{v95:thm-triple-cyclotomic}, double-count ordered triples to derive $6T_3(H)=n(\mathcal C(U)-q+2)$.
\item Show that for $q=2$ this gives ordinary sum-freeness. Why cannot $U$ simply be replaced by $H$ when $q>2$?
\end{enumerate}
\end{exercise}

\begin{exercise}[Three moments of Gaussian periods]\label{v95:ex-period-three}
Use the notation of Section~\ref{v95:subsec-period-moments}: $Q=q^k$, $M=|U|$, and Gaussian periods $\eta_1,\ldots,\eta_e$.
\begin{enumerate}
\item Evaluate $\sum_i\eta_i$ and $\sum_i\eta_i^2$ by additive-character orthogonality.
\item Normalize the equation $u_1+u_2+u_3=0$ multiplicatively to show $\sum_i\eta_i^3=Q\mathcal C(U)-M^2$.
\item When $g=1$ and $k\ge3$, deduce a necessary and sufficient criterion for $\lambda_q(H)\ge4$ using only the third moment.
\end{enumerate}
\end{exercise}

\begin{exercise}[Hyperplane intersections in dimension three and order seven]\label{v95:ex-seven-periods}
Assume $k=3$, $g=\gcd(n,q-1)=1$, and $\lambda_q(H)=4$. Let $A_j=\#\{[a]:N_H(a)=j\}$ for $j=0,1,2$.
\begin{enumerate}
\item Using the number of projective hyperplanes and their first two intersection moments, prove
\[
A_2=\binom n2,\qquad A_1=n(q-n+2),\qquad
A_0=q^2+q+1-A_1-A_2.
\]
\item For $q=11,n=7$, determine all nineteen Gaussian period values and multiplicities, and verify the third moment.
\item Explain why $q=2,n=7$ has the same nonzero Frobenius cosets but cannot satisfy this distribution.
\end{enumerate}
\end{exercise}

\begin{exercise}[Exceptional characteristics for Fourier minors]\label{v95:ex-fourier-minor}
Let $\ell\ge3$ be prime, $p\ne\ell$, $q=p^f$, and $k=\ord_\ell(q)\ge2$. Put $\Gamma=\langle q\bmod\ell\rangle$. You may use the theorem that every square minor of the complex Fourier matrix of prime order is nonzero.
\begin{enumerate}
\item For each $k$-element subset $S\subset\Z/\ell\Z$, show that $D_{\Gamma,S}(\zeta_\ell)$ is a nonzero algebraic integer whose absolute integer norm is at most $(k!)^{\ell-1}$.
\item Show that for $p>(k!)^{\ell-1}$ every such minor remains nonzero after reduction to characteristic $p$.
\item Apply the Moore determinant criterion to deduce $\lambda_q(H)=k+1$. Identify precisely which external theorem has been used.
\end{enumerate}
\end{exercise}

\begin{exercise}[Local rank of the power-sum map]\label{v109:ex-power-jacobian}
Let $N\ge m\ge1$ and $\Char\F_q=p>m$. Define $\Phi:\F_q^N\to\F_q^m$ by taking its $j$-th coordinate to be $\sum_i x_i^j$.
Write its Jacobian matrix and prove that the determinant of any $m$ columns is $m!\prod_{i<j}(x_j-x_i)$, with indices ordered along the selected columns. Show that
\[
\rank J_\Phi(x)=\min\{m,\#\{x_1,\ldots,x_N\}\}.
\]
Explain why this local rank of the formal differential does not assert injectivity on the finite point set. Give a counterexample to the rank formula when $p\le m$.
\end{exercise}

\begin{exercise}[Two equal moments and different empty fibers]\label{v115:ex-two-moments}
On a cyclic group of order three, give two counting functions the values \((0,3,3)\) and \((1,1,4)\), in a fixed order.
Find their means, \(V\), and numbers of empty fibers, and compare the bound in Proposition~\ref{v115:prop-empty-cosets}.
Compute the nontrivial Fourier coefficients and explain why their squared magnitudes do not recover whether each fiber is empty.
These are counting functions on a finite group; no realization as exact-order trace distributions is assumed.
\end{exercise}

\begin{exercise}\label{v145:ex-quadratic-fibers}
Let $P$ be a prime power, $a\in\F_{P^2}$, and $A_a=\left(\begin{smallmatrix}0&-1\\1&-a\end{smallmatrix}\right)$. Determine all fiber sizes of $a\mapsto-\operatorname{tr}(A_aA_{a^P})$, including reducible quadratic targets.
Further determine which repeated-root returns are scalar matrices and which have a nontrivial Jordan block.
\end{exercise}

\begin{exercise}\label{v146:ex-generator-values}
Compare \(\Theta_{k,r}\) and \(\Theta_{k,1}\): do equal fiber sizes imply equal functions? First determine whether \(\Theta_{k,r}\) is constant on conjugate orbits of \(\F_{P^k}/\F_P\). Then take \(P=2,k=5\), write \(\F_{32}=\F_2(\alpha)\) with \(\alpha^5+\alpha^3+1=0\), and compare \(\Theta_{5,1}(1+\alpha)\) and \(\Theta_{5,2}(1+\alpha)\).
\end{exercise}

\begin{exercise}[Alternating deviations among three constants]\label{v148:ex-repeated-profile}
Take $K=\F_4=\F_2(\omega)$ and $m\ge1$. Find all fiber sizes of target $h=(Z+1)^{2m}$ at constants $1,\omega,\omega^2$, and determine which are largest. For target $h=(Z^2+Z+1)^e$, determine when the three compatible constants have equal fibers.
\end{exercise}

\chapter{Exact-Period Layers, Affine Compression, and Orbit Counts}
\label{chap:period-layers}
Exact-period layers separate orbits of prescribed length from the full fixed-point sets. Polynomial divisibility, coefficient gaps, and affine invariants refine this decomposition, while Euler products reconnect local orbits with global point counts.

Parity cancellation follows by comparing multiplicities in finite products, and invariants are determined by compatibility between translations and dilations.
The former describes the coefficient structure of exact-period layers; the latter compresses affine orbits in a layer into new root coordinates, whose factorization is recovered by return maps.

\section{Fixed-Point Sets and Exact-Period Layers}
Fix \(k=\F_q\), and write
\[
A_{q,m}(X)=X^{q^m}-X,\qquad
C_{q,m}(X)=1-X^{q^m-1}\qquad(m\ge1).
\]
\leavevmode\phantomsection\label{def:exact-period-layer}\ignorespaces In \(\F_{q^n}\), let \(\mathcal P_{q,n}\) be the set of elements whose least positive period under the relative Frobenius
\(\sigma:x\mapsto x^q\) is exactly \(n\), and define
\[
E_{q,n}(X)=\prod_{\alpha\in\mathcal P_{q,n}}(X-\alpha),
\qquad D_{q,n}=|\mathcal P_{q,n}|.
\]
An empty product is taken to be one. The set \(\mathcal P_{q,n}\) is called an \textbf{exact-period layer},
and \(E_{q,n}\) its \textbf{exact-period polynomial}. The period here is always the Frobenius period;
it does not require the element to be primitive in the multiplicative group.

\begin{theorem}\leavevmode\phantomsection\label{thm:exact-period-layer}
For every \(n\ge1\), the polynomial \(E_{q,n}\) is monic and has no repeated roots, and
\[
E_{q,n}(X)=\prod_{f\in\Irr_q(n)}f(X)\in\F_p[X],
\qquad D_{q,n}=nI_q(n).
\]
In particular, \(E_{q,1}=X^q-X\). For every \(m\ge1\),
\begin{equation}\label{eq:period-layer-product}
A_{q,m}=\prod_{d\mid m}E_{q,d}.
\end{equation}
Exact-period polynomials with distinct indices are pairwise coprime.
\end{theorem}
\begin{proof}
The orbit theorem shows that the exact period of \(\alpha\) equals its degree over \(k\).
Thus \(\mathcal P_{q,n}\) is precisely the disjoint union of the root sets of all monic irreducible polynomials of degree \(n\);
each such polynomial has \(n\) distinct roots, giving the product and degree formulas.
For \(n=1\), this set is exactly \(k\).

Every root of \(A_{q,m}\) has period dividing \(m\), and every element of period \(d\mid m\) is a root of \(A_{q,m}\).
These layers are disjoint, and \(A_{q,m}\) has no repeated roots, giving~\eqref{eq:period-layer-product}.
Two distinct indices can be compared in a common extension field of order \(q^{\lcm(m,n)}\).
An element cannot have two different least positive periods, so the layers have no common root and their polynomials are coprime.

Finally, the absolute Frobenius \(\phi:x\mapsto x^p\) commutes with \(\sigma\) and is invertible on \(\F_{q^n}\).
If \(\sigma^t(\phi(\alpha))=\phi(\alpha)\), applying \(\phi^{-1}\) gives \(\sigma^t(\alpha)=\alpha\); the converse also holds.
Thus \(\phi\) permutes \(\mathcal P_{q,n}\), and all coefficients of the exact-period polynomial are fixed by the \(p\)-th power map, hence lie in \(\F_p\).
\end{proof}

When \(n>1\), the set \(\mathcal P_{q,n}\) contains neither zero nor one, so it is not a subfield.
The subfield \(\F_{q^n}\) is the union of the exact-period layers for all \(d\mid n\).
Here subfields correspond to divisibility conditions on periods, while individual layers correspond to equality conditions on least periods.

\par\Needspace{12\baselineskip}
\section{First Occurrence and Parity Cancellation}
\begin{lemma}[First occurrence of an irreducible factor]\label{lem:period-first-return}
Let \(f\in k[X]\) be monic irreducible, with \(f\ne X\).
There is a least positive integer \(\rho\) such that \(f\mid C_{q,\rho}\), and
\[
\rho=\deg f,\qquad
f\mid C_{q,m}\quad\Longleftrightarrow\quad \rho\mid m.
\]
\end{lemma}
\begin{proof}
Since \(A_{q,m}=-XC_{q,m}\) and \(f\ne X\), we have \(f\mid C_{q,m}\) if and only if \(f\mid A_{q,m}\).
Theorem~\ref{thm:divisibility-main} therefore gives \(f\mid C_{q,m}\Longleftrightarrow\deg f\mid m\), whose least positive index is \(\deg f\).
\end{proof}

This first occurrence is precisely the least return period of a root in Theorem~\ref{thm:orbit}; we now use this period law for parity cancellation.

\begin{lemma}[Strong divisibility]\label{lem:C-strong-divisibility}
For positive integers \(a,b\), the polynomials \(\gcd(C_{q,a},C_{q,b})\) and \(C_{q,\gcd(a,b)}\) differ by a nonzero constant factor.
Every \(C_{q,m}\) has no repeated roots.
\end{lemma}
\begin{proof}
Since \(A_{q,m}=-XC_{q,m}\) and \(C_{q,m}(0)=1\), the zero-root factor in each \(A_{q,m}\) is exactly \(X\).
Canceling this common factor in the greatest-common-divisor identity of Corollary~\ref{cor:gcd-poly} gives the stated relation.
Also, \(C_{q,m}\) divides the polynomial \(A_{q,m}\), which has no repeated roots, so \(C_{q,m}\) has none either.
\end{proof}

\leavevmode\phantomsection\label{def:distinct-prime-count}\ignorespaces Let \(n>1\), with distinct prime divisors \(\ell_1,\ldots,\ell_r\).
Write
\[
\rad(n)=\ell_1\cdots\ell_r,\qquad \omega(n)=r,\qquad
n_S=\frac{n}{\prod_{i\in S}\ell_i}
\quad(S\subseteq\{1,\ldots,r\}).
\]
Here \(\rad(n)\) retains the notation for the square-free kernel used in Section~\ref{sec:cyclotomic-med}.
Define
\[
P_+=\prod_{|S|\text{ even}}C_{q,n_S},\qquad
P_-=\prod_{|S|\text{ odd}}C_{q,n_S}.
\]

\begin{theorem}[Parity cancellation]\leavevmode\phantomsection\label{thm:parity-period}
We have \(P_-\mid P_+\), and
\begin{equation}\label{eq:period-mobius}
\frac{P_+}{P_-}=E_{q,n}
=\prod_{d\mid n}C_{q,d}^{\,\mu(n/d)}
=\prod_{d\mid n}A_{q,d}^{\,\mu(n/d)}.
\end{equation}
The last two products are initially interpreted in the rational function field; they are in fact polynomials.
For \(n>1\), we have \(E_{q,n}(0)=1\).
\end{theorem}
\begin{proof}
Compare the monic irreducible factors of \(P_+\) and \(P_-\) one at a time.
The zero-root factor \(X\) occurs in none of the \(C_{q,m}\).
Suppose that another factor \(f\) occurs, with first-occurrence index \(\rho\).
By Lemma~\ref{lem:period-first-return}, it occurs in \(C_{q,n_S}\) if and only if \(\rho\mid n_S\);
as \(C_{q,n_S}\) has no repeated roots, each occurrence contributes multiplicity one.
Any occurrence also implies \(\rho\mid n\).

If \(\rho<n\), there is some \(\ell_j\mid n\) with \(v_{\ell_j}(\rho)<v_{\ell_j}(n)\).
Add or remove the index \(j\) from the subset \(S\), giving the involution \(S\mapsto S\triangle\{j\}\).
This operation only changes the \(\ell_j\)-valuation of \(n_S\) between \(v_{\ell_j}(n)\) and \(v_{\ell_j}(n)-1\);
both are at least \(v_{\ell_j}(\rho)\), and all other prime valuations are unchanged.
Thus
\[
\rho\mid n_S\quad\Longleftrightarrow\quad
\rho\mid n_{S\triangle\{j\}}.
\]
This involution has no fixed point and reverses the parity of the subset size, so \(f\) has equal multiplicities in \(P_+\) and \(P_-\).

If \(\rho=n\), only \(S=\varnothing\) satisfies \(\rho\mid n_S\),
so \(f\) occurs once more in \(P_+\).
By unique factorization, every denominator factor cancels, and the irreducible factors of the quotient are exactly those with first-occurrence index \(n\), each occurring once.
Since the first-occurrence index equals the degree, these are precisely the polynomials in \(\Irr_q(n)\).
There are \(2^{r-1}\) even subsets and \(2^{r-1}\) odd subsets, and every \(C_{q,m}\) has leading coefficient \(-1\),
so the quotient is monic and equals \(E_{q,n}\).

The condition \(\mu(n/d)\ne0\) corresponds exactly to \(d=n_S\), in which case \(\mu(n/d)=(-1)^{|S|}\),
giving the product formula in terms of \(C\).
Since \(A_{q,d}=-XC_{q,d}\) and \(\sum_{d\mid n}\mu(n/d)=0\), the total exponent of \(-X\) is zero, giving the product formula in terms of \(A\).
Every \(C_{q,d}(0)=1\), so the quotient has constant term one.
\end{proof}

\section{The First Nonconstant Term and Existence of Irreducible Polynomials}
For a nonzero polynomial \(H(X)\), let \(\operatorname{ord}_X(H)\) denote the degree of its lowest nonzero term.
The following congruences with negative exponents are interpreted in the truncated ring \(k[X]/(X^N)\):
a polynomial with constant term one is invertible there. For example,
\[
(1-X^a)^{-1}\equiv1+X^a+\cdots+X^{ja}\pmod{X^N},
\qquad j=\left\lfloor\frac{N-1}{a}\right\rfloor.
\]
Thus these congruence calculations involve only finitely many coefficients.

\begin{theorem}[The first nonconstant term and nonconstancy]\label{thm:period-gap}
Let \(n>1\), \(s=n/\rad(n)\), and \(M=q^s-1\). Then
\[
E_{q,n}(X)\equiv
1+(-1)^{\omega(n)+1}X^M\pmod{X^{M+1}}.
\]
In particular,
\[
\operatorname{ord}_X(E_{q,n}-1)=q^{n/\rad(n)}-1,
\]
so \(E_{q,n}\) is nonconstant.
\end{theorem}
\begin{proof}
Among all indices \(n_S\) in the parity products, the smallest is \(s\), attained only for \(S=\{1,\ldots,r\}\).
Thus \(C_{q,s}=1-X^M\) is the only factor that can contribute a nonzero term of degree \(M\);
all other factors and their inverses are one modulo \(X^{M+1}\).
The exponent of this factor is \((-1)^r\).
If \(r\) is even, we retain \(1-X^M\);
if \(r\) is odd, we retain \((1-X^M)^{-1}\equiv1+X^M\pmod{X^{M+1}}\), using \(2M\ge M+1\).
The two cases give the stated formula.
The coefficient \((-1)^{r+1}\) is nonzero in every characteristic, so \(M\) is the least positive degree with a nonzero coefficient in \(E_{q,n}\).
\end{proof}

\begin{corollary}[Existence from cancellation and coefficients]\label{cor:euclid-frobenius-existence}
For every positive integer \(n\), there is an irreducible polynomial of degree \(n\) in \(\F_q[X]\).
\end{corollary}
\begin{proof}
For \(n=1\), take \(X-a\), with \(a\in k\).
If \(n>1\), the parity cancellation theorem gives a polynomial all of whose irreducible factors have degree \(n\);
Theorem~\ref{thm:period-gap} ensures that this polynomial is nonconstant.
By the lemma on existence of irreducible factors, it has at least one irreducible factor, so at least one irreducible polynomial of degree \(n\) exists.
\end{proof}

This argument has two parts: cancellation by divisibility removes shorter periods, and a local coefficient ensures that the remaining polynomial is nonconstant.
Comparing degrees on the two sides of~\eqref{eq:period-mobius} also gives
\[
D_{q,n}=\sum_{d\mid n}\mu(n/d)q^d.
\]
This agrees with the count of irreducible polynomials in Chapter~\ref{chap:frobenius}; the first nonconstant term supplies local structural information beyond the degree formula.

\begin{proposition}[Self-reciprocity]\label{prop:layer-reciprocal}
If \(n>1\) and \(D=D_{q,n}\), then
\[
X^DE_{q,n}(X^{-1})=E_{q,n}(X).
\]
Consequently, if \(E_{q,n}=\sum_{i=0}^D a_iX^i\), then \(a_i=a_{D-i}\);
the low-degree coefficient gap in Theorem~\ref{thm:period-gap} also occurs at the high-degree end.
\end{proposition}
\begin{proof}
For a nonzero element \(\alpha\), we have \((\alpha^{-1})^{q^t}=\alpha^{-1}\) if and only if \(\alpha^{q^t}=\alpha\).
Thus inversion permutes \(\mathcal P_{q,n}\).
The exact-period polynomial has constant term one, so \(X^DE_{q,n}(X^{-1})\) is monic;
it has exactly the same simple roots as \(E_{q,n}\), hence the two polynomials are equal.
Comparison of coefficients gives the final assertion.
\end{proof}

\begin{example}[Geometric series for prime-power layers]
If \(n=\ell^a\), where \(\ell\) is prime, only two M\"obius exponents are nonzero, so
\[
E_{q,\ell^a}(X)=
\frac{1-X^{q^{\ell^a}-1}}{1-X^{q^{\ell^{a-1}}-1}}.
\]
Put \(Q=q^{\ell^{a-1}}\), \(M=Q-1\), and \(L=1+Q+\cdots+Q^{\ell-1}\). Then \(Q^\ell-1=ML\), giving
\[
E_{q,\ell^a}(X)=\sum_{j=0}^{L-1}X^{jM}.
\]
In particular,
\[
E_{2,4}=1+X^3+X^6+X^9+X^{12}.
\]
Its three irreducible factors are \(X^4+X+1\), \(X^4+X^3+1\), and \(X^4+X^3+X^2+X+1\).
The sparsity of this product reflects the prime-power structure of the index four.
\end{example}

\section{The Square-Free Kernel and Base-Field Lifting}
\begin{theorem}\label{thm:period-radical}
Let \(n>1\), \(r=\rad(n)\), and \(s=n/r\). Then
\[
E_{q,n}(X)=E_{q^s,r}(X).
\]
\end{theorem}
\begin{proof}
In the M\"obius product, the nonzero exponents correspond to \(n/d=t\), where \(t\mid r\).
Substituting \(d=n/t=s(r/t)\) gives
\[
\begin{aligned}
E_{q,n}(X)
&=\prod_{t\mid r}\bigl(X^{q^{s(r/t)}}-X\bigr)^{\mu(t)}\\
&=\prod_{t\mid r}\bigl(X^{(q^s)^{r/t}}-X\bigr)^{\mu(t)}
=E_{q^s,r}(X).
\end{aligned}
\]
All expressions are polynomials by the preceding theorems, so this is also a polynomial identity.
\end{proof}

The formula also has a direct orbit interpretation. If an element has \(q\)-Frobenius period \(d\),
then its \(q^s\)-Frobenius period is \(d/\gcd(d,s)\).
In the fixed field of the \(q^{sr}=q^n\)-th power map, we have \(d\mid n\).
For each \(\ell\mid n\), write \(v_\ell(n)=a\); then \(v_\ell(s)=a-1\).
The condition
\[
\frac{d}{\gcd(d,s)}=r
\]
requires \(v_\ell(d)-\min\{v_\ell(d),a-1\}=1\), which is equivalent to \(v_\ell(d)=a\).
Comparing all prime divisors gives \(d=n\).
Thus the polynomials on the two sides record the same root set, using different Frobenius steps.

\begin{proposition}[A truncation formula for low-degree coefficients]\label{prop:period-filtration}
Put \(Q=q^s\), and list the positive divisors of \(r=\rad(n)\) as \(1=e_0<e_1<\cdots<e_t=r\).
For \(0\le j<t\), we have
\[
E_{q,n}(X)\equiv
\prod_{i=0}^{j}(1-X^{Q^{e_i}-1})^{\mu(r/e_i)}
\pmod{X^{Q^{e_{j+1}}-1}}.
\]
In particular, if \(\ell\) is the smallest prime divisor of \(n\), then
\[
E_{q,n}(X)\equiv
(1-X^{Q-1})^{(-1)^{\omega(n)}}
\pmod{X^{Q^\ell-1}}.
\]
\end{proposition}
\begin{proof}
The square-free-kernel formula and the M\"obius product give \(E_{q,n}=\prod_{e\mid r}(1-X^{Q^e-1})^{\mu(r/e)}\).
For \(i\ge j+1\), we have \(Q^{e_i}-1\ge Q^{e_{j+1}}-1\);
the corresponding factors and their inverses are all one in the specified truncated ring, so they may be omitted.
Taking \(j=0\) and noting that \(e_1=\ell\) and \(\mu(r)=(-1)^{\omega(n)}\) gives the second formula.
\end{proof}

\begin{example}
Take \(q=2,n=12\). Then \(r=6,s=2,Q=4\), so
\[
E_{2,12}=E_{4,6},\qquad
E_{2,12}(X)\equiv1+X^3\pmod{X^{15}}.
\]
The second formula identifies the first nonconstant term, of degree three, and shows that all terms of degrees four through \(14\) vanish.
Self-reciprocity gives the corresponding coefficient gap at the high-degree end.
\end{example}

\par\Needspace{12\baselineskip}
\section{Connection with Cyclotomic Polynomials}\label{v48:sec-period-cyclotomic}
\begin{proposition}\label{prop:period-cyclotomic}
For \(n>1\), in \(\F_q[X]\) we have
\[
E_{q,n}(X)=
\prod_{\substack{e\mid q^n-1\\ \ord_e(q)=n}}\overline{\Phi}_e(X),
\]
where \(\overline{\Phi}_e\) denotes the reduction of the integer cyclotomic polynomial in characteristic \(p\).
\end{proposition}
\begin{proof}
Every element of \(\mathcal P_{q,n}\) is nonzero, and its multiplicative order \(e\) divides \(q^n-1\), hence is coprime to \(p\).
Proposition~\ref{prop:degree-from-order} shows that its Frobenius period is \(\ord_e(q)\).
Thus its exact period is \(n\) if and only if its multiplicative order \(e\) satisfies the conditions in the product.
For a fixed such \(e\), all elements of order \(e\) are precisely the distinct roots of \(\overline{\Phi}_e\).
Different multiplicative orders give disjoint root sets, and these roots exhaust \(\mathcal P_{q,n}\), proving the equality.
\end{proof}

For example, \(E_{2,4}=\overline{\Phi}_5\,\overline{\Phi}_{15}\):
the roots of one quartic factor have multiplicative order five, and the roots of the other two quartic factors have multiplicative order \(15\).
All three have Frobenius period four, showing that an exact-period layer may contain elements of different multiplicative orders.
\subsection{Decomposition by Multiplicative Order and Counting}
\begin{proposition}[Decomposition of an exact-period layer by multiplicative order]\label{v48:prop-period-order-count}
Let \(n>1\). For each \(e\) satisfying \(e\mid q^n-1\) and \(\ord_e(q)=n\),
the set \(\mathcal P_{q,n}\) contains exactly \(\varphi(e)\) elements of multiplicative order \(e\);
they form \(\varphi(e)/n\) Frobenius orbits. Consequently,
\begin{equation}\label{v48:eq-period-order-count}
D_{q,n}=nI_q(n)
=\sum_{\substack{e\mid q^n-1\\ \ord_e(q)=n}}\varphi(e)
=\sum_{d\mid n}\mu(n/d)q^d.
\end{equation}
In particular, the number of irreducible polynomials of degree \(n\) whose roots are primitive elements of \(\F_{q^n}\) is
\(\varphi(q^n-1)/n\).
\end{proposition}
\begin{proof}
The cyclic group \(\F_{q^n}^{\times}\) contains exactly \(\varphi(e)\) elements of order \(e\).
Their Frobenius periods are determined by their multiplicative orders, so the listed order classes partition \(\mathcal P_{q,n}\).
Every orbit contains \(n\) elements, so dividing the count in each order class by \(n\) gives the corresponding number of irreducible factors.
The last equality is the degree formula from~\eqref{eq:period-mobius}.
Also, \(\ord_{q^n-1}(q)=n\): we have \(q^n\equiv1\pmod{q^n-1}\),
whereas for \(0<d<n\), we have \(0<q^d-1<q^n-1\), ruling out a shorter period.
A primitive element has multiplicative order \(q^n-1\); applying the preceding count to this order class gives the last assertion.
\end{proof}

The two counts here use, respectively, inclusions among subfields and element orders in the multiplicative group.
The former excludes elements in proper subfields; the latter determines which multiplicative orders give the prescribed conjugate-orbit length.
An element of order \(e\) is a primitive \(e\)-th root of unity; only when \(e=q^n-1\)
is it a primitive element of the whole field \(\F_{q^n}\).
For \(n=1\), zero must be counted separately:
\[
E_{q,1}(X)=X\prod_{e\mid q-1}\overline\Phi_e(X),
\qquad D_{q,1}=1+\sum_{e\mid q-1}\varphi(e)=q.
\]
Here \(X\) is the minimal polynomial of zero, which cannot be obtained from the decomposition by multiplicative order.

\subsection{Base-Field Extension and Reciprocal Pairing}
For an order class satisfying \(\ord_e(q)=n\), extend the base field to \(\F_{q^r}\) and put \(g=\gcd(n,r)\). Then
\[
\ord_e(q^r)=\frac{n}{g}.
\]
By Theorem~\ref{thm:factor-base-extension}, each original orbit splits into \(g\) orbits of length \(n/g\),
so \(\overline\Phi_e\) has exactly \(g\varphi(e)/n\) irreducible factors over the new base field.
Correspondingly, the whole polynomial \(E_{q,n}\) splits into \(gI_q(n)\) irreducible factors of degree \(n/g\).
All roots of the original layer lie in the exact-period layer of period \(n/g\) over the new base field, but the latter may also contain roots of other original periods;
the two exact-period polynomials therefore cannot in general be identified.

Inversion preserves multiplicative order, so reciprocal pairing occurs separately within each \(\overline\Phi_e\).
If \(\ord_e(q)=n>1\), Corollary~\ref{v48:cor-reciprocal-order} shows that
the irreducible factors in this order class are self-reciprocal if and only if
\begin{equation}\label{v48:eq-period-reciprocal-order}
n\text{ is even},\qquad q^{n/2}\equiv-1\pmod e.
\end{equation}
This condition depends only on \(e\), so all irreducible factors in an order class are either self-reciprocal or paired with distinct reciprocals.
In particular, factors in an exact-period layer of odd degree occur in pairs.
Self-reciprocity of the exact-period polynomial can therefore arise either from self-reciprocal individual factors or from pairings of distinct factors.

\Needspace{5\baselineskip}
\begin{example}[Six multiplicative orders in one period]\label{v48:ex-period-order-34}
Retain the extension \(\F_{81}/\F_3\) of Example~\ref{v48:ex-power-subfields}, taking \(q=3,n=4\).
Among the positive divisors of \(3^4-1=80\), those satisfying \(\ord_e(3)=4\) are exactly
\(5,10,16,20,40,80\), so
\[
E_{3,4}=\overline\Phi_5\,\overline\Phi_{10}\,\overline\Phi_{16}\,
\overline\Phi_{20}\,\overline\Phi_{40}\,\overline\Phi_{80}.
\]
The numbers of roots and irreducible quartic factors in each order class are as follows:
\begin{center}
\begin{tabular}{c|rrrrrr}
\(e\)&\(5\)&\(10\)&\(16\)&\(20\)&\(40\)&\(80\)\\\hline
\(\varphi(e)\)&\(4\)&\(4\)&\(8\)&\(8\)&\(16\)&\(32\)\\
\(\varphi(e)/4\)&\(1\)&\(1\)&\(2\)&\(2\)&\(4\)&\(8\)
\end{tabular}
\end{center}
Thus \(D_{3,4}=72=3^4-3^2\) and \(I_3(4)=18\). Eight factors have primitive elements as roots;
the roots of the other ten still generate \(\F_{81}\) over \(\F_3\), but do not generate its multiplicative group.
Since \(3^2=9\), condition~\eqref{v48:eq-period-reciprocal-order} holds only for \(e=5,10\).
The factor of root order \(10\) appeared in Example~\ref{v48:ex-reciprocal-ten}; the other self-reciprocal irreducible factor is
\[
\overline\Phi_5=X^4+X^3+X^2+X+1.
\]
The remaining \(16\) quartic factors form eight reciprocal pairs.
After extending the base field to \(\F_9\), each quartic factor splits into two quadratics, giving \(36\) in total;
after extending to \(\F_{27}\), all \(18\) quartic factors remain irreducible, since \(\gcd(4,3)=1\).
Thus a single original period can accommodate different multiplicative orders and give different factorizations as the base-field step changes.
\end{example}

\section{Affine Invariants and Compression Coordinates}
\label{sec:affine-compression}
\subsection{The Ring of Affine Invariants}
\label{sec:affine-invariant-ring}
\leavevmode\phantomsection\label{def:affine-compression-coordinate}\ignorespaces Put
\[
Z_q(X)=X^q-X,\qquad
\Theta_q(X)=(X^q-X)^{q-1}.
\]
We consider the one-dimensional affine group \(\operatorname{AGL}_1(k)\) defined in Chapter~\ref{chap:root-symmetry},
consisting of all base-field transformations \(X\mapsto aX+b\), where \(a\in k^\times,b\in k\). They act on every extension field of \(k\).
We call \(\Theta_q(\alpha)\) the compression coordinate of \(\alpha\). Proposition~\ref{prop:affine-orbit-separation} will show that, outside the base field, two elements have the same compression coordinate if and only if they lie in the same affine-group orbit.
\leavevmode\phantomsection\label{def:affine-invariant}\ignorespaces A polynomial \(F\) satisfying \(F(aX+b)=F(X)\) for all such \(a,b\)
is called an affine invariant. This requires equality as ordinary polynomials, not merely equality of values at the points of \(k\).

\begin{theorem}[The ring of affine invariants in one variable]\label{thm:affine-invariants}
We have
\[
k[X]^{\operatorname{AGL}_1(k)}=k[\Theta_q(X)].
\]
That is, every affine invariant has a unique expression as a polynomial in \(\Theta_q\).
\end{theorem}
\begin{proof}
First consider translations \(X\mapsto X+b\), where \(b\in k\). Taking \(W=k\) gives \(L_W=X^q-X=Z_q\).
By Theorem~\ref{ffs:thm-additive-composition}, the translation-invariant polynomials are exactly the elements of \(k[Z_q]\).

Now include dilations. Write \(F=H(Z_q)\). Since \(a^q=a\), we have \(Z_q(aX)=aZ_q(X)\).
Thus dilation invariance is equivalent to \(H(aZ_q)=H(Z_q)\).
Substitution of the nonconstant polynomial \(Z_q\) defines an injective map \(k[T]\to k[X]\):
if \(H\ne0\), then \(H(Z_q)\) is nonzero and has degree \(q\deg H\).
Hence the condition is equivalent to \(H(aT)=H(T)\) for every \(a\in k^\times\).
Writing \(H=\sum_j h_jT^j\) and comparing coefficients gives \(h_j(a^j-1)=0\).
Since \(k^\times\) is cyclic, for \(h_j\ne0\) this holds if and only if \(q-1\mid j\).
Thus \(H(T)=G(T^{q-1})\), or \(F=G(\Theta_q)\).
Finally, \(\Theta_q\) is nonconstant, so its substitution map is also injective, ensuring uniqueness of \(G\).
\end{proof}

\subsection{Compressed Root Sets and Orbit Quotients}
\label{sec:affine-orbits-return}
Orbit quotients are defined on page~\pageref{def:orbit-quotient-free}. A free action ensures that a unique group element sends a specified point to another point of the same orbit; see page~\pageref{ffs:def-point-stabilizer}.

\begin{proposition}[Affine invariants separate orbits]\leavevmode\phantomsection\label{prop:affine-orbit-separation}
Let \(\alpha,\beta\notin k\) lie in a common finite extension. Then
\[
\Theta_q(\alpha)=\Theta_q(\beta)
\quad\Longleftrightarrow\quad
\beta=a\alpha+b\quad\text{for some }a\in k^\times,\ b\in k.
\]
Thus every nonempty fiber of \(\Theta_q\) outside the base field is a single \(\operatorname{AGL}_1(k)\)-orbit,
and each such orbit has \(q(q-1)\) elements.
\end{proposition}
\begin{proof}
Since \(\alpha,\beta\notin k\), both \(Z_q(\alpha)\) and \(Z_q(\beta)\) are nonzero.
If their compression values are equal, then
\[
a=\frac{Z_q(\beta)}{Z_q(\alpha)}
\quad\text{satisfies}\quad a^{q-1}=1.
\]
Since the root set of \(X^q-X\) is \(k\), we have \(a\in k^\times\). Put \(b=\beta-a\alpha\); then \(b^q-b=Z_q(\beta)-aZ_q(\alpha)=0\), so \(b\in k\).
Conversely, \(Z_q(a\alpha+b)=aZ_q(\alpha)\) and \(a^{q-1}=1\) directly give equality of the compression values.

If \(a\alpha+b=\alpha\), then \((a-1)\alpha=-b\). Since \(\alpha\notin k\), we must have \(a=1,b=0\),
so the affine group acts freely on these elements, and each orbit has size equal to the group order \(q(q-1)\).
\end{proof}

\Needspace{11\baselineskip}
\begin{theorem}[Affine compression of exact-period layers]\leavevmode\phantomsection\label{thm:period-compression}
For every \(n>1\), there is a unique \(G_{q,n}(T)\in\F_p[T]\) such that
\[
E_{q,n}(X)=G_{q,n}(\Theta_q(X)).
\]
It satisfies
\[
G_{q,n}(0)=1,\qquad
\deg G_{q,n}=\frac{D_{q,n}}{q(q-1)}.
\]
\end{theorem}
\begin{proof}
For \(a\in k^\times,b\in k\), we have
\[
A_{q,d}(aX+b)=aA_{q,d}(X).
\]
Substitution in~\eqref{eq:period-mobius} gives total exponent \(\sum_{d\mid n}\mu(n/d)=0\) for \(a\), so \(E_{q,n}(aX+b)=E_{q,n}(X)\).
The affine-invariant ring theorem first gives a unique \(G_{q,n}\) in \(k[T]\).

By Theorem~\ref{thm:exact-period-layer} and \(\Theta_q=(X^q-X)^{q-1}\), we have \(E_{q,n},\Theta_q\in\F_p[X]\).
Let \(G^{[p]}\) denote the polynomial obtained by raising the coefficients of \(G\) to the \(p\)-th power.
Applying the absolute Frobenius coefficientwise to the identity gives \(E_{q,n}=G_{q,n}^{[p]}(\Theta_q)\).
Uniqueness of the representation gives \(G_{q,n}^{[p]}=G_{q,n}\), hence \(G_{q,n}\in\F_p[T]\).
Setting \(X=0\) gives \(G_{q,n}(0)=E_{q,n}(0)=1\).
Since \(\deg\Theta_q=q(q-1)\) and \(E_{q,n}\) is nonconstant, comparison of degrees gives the formula.
\end{proof}

\Needspace{11\baselineskip}
\begin{corollary}[Frobenius on a compressed layer]\label{cor:compressed-orbit-period}
Let \(n>1\). The map \([\alpha]\mapsto\Theta_q(\alpha)\) gives a bijection
\[
\mathcal P_{q,n}/\operatorname{AGL}_1(k)
\longrightarrow\{c:G_{q,n}(c)=0\},
\]
compatible with Frobenius. The polynomial \(G_{q,n}\) is monic and has no repeated roots, and its irreducible factors in \(k[T]\) correspond to Frobenius orbits on this orbit quotient.
For \(\alpha\in\mathcal P_{q,n}\), we have
\begin{equation}\label{eq:compressed-root-degree}
\begin{split}
[k(\Theta_q(\alpha)):k]
=\min\{t\ge1:\ &\alpha^{q^t}=a\alpha+b\\
&\text{for some }a\in k^\times,\ b\in k\}.
\end{split}
\end{equation}
\end{corollary}
\begin{proof}
Since \(k(a\alpha+b)=k(\alpha)\), the affine group preserves \(\mathcal P_{q,n}\).
It also commutes with Frobenius, so \(\sigma\) induces a permutation on the orbit quotient.
The identity \(E_{q,n}=G_{q,n}(\Theta_q)\) shows that compression values are roots of \(G_{q,n}\);
Proposition~\ref{prop:affine-orbit-separation} shows that equality of two values corresponds exactly to membership in the same affine orbit.
Conversely, if \(G_{q,n}(c)=0\), choose a root \(\alpha\) of \(\Theta_q(X)-c\) in a finite extension.
Substitution gives \(E_{q,n}(\alpha)=0\), so \(\alpha\in\mathcal P_{q,n}\). This proves the bijection.

Each affine orbit has \(q(q-1)\) elements, so the number of distinct compression values is \(D_{q,n}/(q(q-1))=\deg G_{q,n}\), showing that \(G_{q,n}\) has no repeated roots.
It is monic because both \(E_{q,n}\) and \(\Theta_q\) are monic.
Also, \(\Theta_q(\alpha^q)=\Theta_q(\alpha)^q\), so the bijection is compatible with \(\sigma\).
The root-orbit theorem therefore identifies the irreducible factors of \(G_{q,n}\) in \(k[T]\) with the Frobenius orbits on the orbit quotient.
Finally, \(\Theta_q(\alpha)^{q^t}=\Theta_q(\alpha)\) is equivalent to \(\Theta_q(\alpha^{q^t})=\Theta_q(\alpha)\); applying orbit separation gives~\eqref{eq:compressed-root-degree}.
\end{proof}

The period of an original element requires Frobenius to return to that element; the period of its compression value only requires a return to the same affine orbit.
This explains why root degrees may drop after compression; for example, when \(n=2\), all compression values equal \(-1\); see Exercise~\ref{ex:wu-compressed-not-period}.
Section~\ref{sec:compression-return-lifting} identifies the period-reduction factor \(n/t\) with the internal return period of the fiber
and uses this period to determine the factorization obtained by lifting a compressed factor.

\section{Compressed Periods, Return Orders, and Factor Lifting}
\label{sec:compression-return-lifting}

The Frobenius period of a compression value determines the first return time of its fiber.
The return permutation within the fiber further determines the degrees of the original roots and the degrees and number of lifted factors.

\subsection{Compressed Periods and Orders of Affine Returns}
\label{sec:affine-return-order}
\begin{proposition}[Compressed periods and orders of affine returns]\label{prop:affine-return-order}
Let \(n>1\), let \(\alpha\) have degree \(n\) over \(k=\F_q\), and put
\[
c=\Theta_q(\alpha),\qquad t=[k(c):k].
\]
Then \(t\mid n\), there is a unique \(g(X)=aX+b\in\operatorname{AGL}_1(k)\) satisfying \(\sigma^t(\alpha)=g(\alpha)\), and
\begin{equation}\label{eq:affine-period-drop}
\boxed{\quad\frac nt=\ord(g),\qquad
\frac nt=p\ \text{or}\ \frac nt\mid q-1.\quad}
\end{equation}
\end{proposition}
\begin{proof}
Since \(k(c)\subseteq k(\alpha)\), the tower formula gives \(t\mid n\).
Corollary~\ref{cor:compressed-orbit-period} gives the required affine transformation,
and the free action in Proposition~\ref{prop:affine-orbit-separation} ensures uniqueness.
Since the coefficients of \(g\) lie in \(k\), we have \(\sigma g=g\sigma\), and induction gives
\[
\sigma^{jt}(\alpha)=g^j(\alpha)\qquad(j\ge1).
\]
The left-hand side first equals \(\alpha\) at \(j=n/t\). The right-hand side equals \(\alpha\) if and only if \(g^j\) is the identity,
because \(\alpha\notin k\) has trivial point stabilizer in \(\operatorname{AGL}_1(k)\). Thus \(n/t=\ord(g)\).

If \(a=1,b=0\), then \(g\) has order one; if \(a=1,b\ne0\),
then \(g^j(X)=X+jb\), so its order is exactly \(p\).
If \(a\ne1\), put \(z=b/(1-a)\in k\). Then
\[
g(X)-z=a(X-z).
\]
Thus a translation conjugates \(g\) to \(X\mapsto aX\), whose order \(\ord(a)\) divides \(q-1\).
This exhausts the cases and proves~\eqref{eq:affine-period-drop}.
\end{proof}

Every permitted nontrivial period-reduction factor occurs. Fix \(k=\F_q\). For every integer \(e>1\) satisfying \(e=p\) or \(e\mid q-1\),
there exists \(\alpha\notin k\) with original period \(n=e\) and compressed period \(t=1\), hence period-reduction factor \(e\).
Specifically, a root of \(\Theta_q(\alpha)=1\) gives a nonidentity translation and period-reduction factor \(p\);
if \(c\in k^\times\) has order \(e>1\), a root of \(\Theta_q(\alpha)=c\) gives a return transformation conjugate to a dilation, with period-reduction factor \(e\).
These roots all lie outside \(k\); see Exercise~\ref{ex:wu7-compressed-linear-factors} for their degrees and return transformations.

\noindent\textbf{A fiber interpretation of the period drop.}
\phantomsection\label{sec:affine-return-meaning}
Retain the assumption \(n>1\). Thus \(\alpha\notin k\) and \(c=\Theta_q(\alpha)\ne0\), and the affine group acts freely at \(\alpha\). Put \(e=\ord(g)=n/t\), and write the full fiber and original orbit in \(k(\alpha)\) as
\[
F_c=\{\beta:\Theta_q(\beta)=c\},\qquad
\mathcal O_\alpha=\{\sigma^j(\alpha):0\le j<n\}.
\]
Then
\begin{equation}\label{eq:full-fiber-orbit-intersection}
\begin{aligned}
F_c&=\operatorname{AGL}_1(k)\alpha,
&|F_c|&=q(q-1),\\
\mathcal O_\alpha\cap F_c&=\langle g\rangle\alpha,
&|\mathcal O_\alpha\cap F_c|&=e.
\end{aligned}
\end{equation}
The first line follows from Proposition~\ref{prop:affine-orbit-separation}; all these roots in fact lie in \(k(\alpha)\). Equivariance gives
\(\Theta_q(\sigma^j(\alpha))=\sigma^j(c)\), so \(\sigma^j(\alpha)\in F_c\) if and only if \(t\mid j\).
These intersection points are exactly \(\sigma^{jt}(\alpha)=g^j(\alpha)\), for \(0\le j<e\), and the free action ensures that they are distinct, proving the second line.

Thus \(n=te\): \(\Theta_q\) retains the period \(t\) between fibers, while the internal return period \(e=[k(\alpha):k(c)]\) is the factor removed by compression.

For the whole fiber, use the free action to identify \(F_c\) with the affine group by \(h\mapsto h(\alpha)\). Since every \(h\) has coefficients in \(k\),
\[
\sigma^t(h(\alpha))=h(\sigma^t(\alpha))=(hg)(\alpha).
\]
The first return corresponds to right multiplication \(h\mapsto hg\), so the full fiber splits into \(q(q-1)/e\) return orbits of length \(e\). The intersection \(\mathcal O_\alpha\cap F_c\) is just the one containing \(\alpha\), and need not be the full fiber.
This also explains the number of factors in Corollary~\ref{cor:compressed-fiber-lifting}.
Changing the base point to \(h(\alpha)\) changes the return element to \(hgh^{-1}\), preserving its order, so the internal period is independent of the base point.
Right multiplication must be distinguished from applying the fixed affine expression \(g\) to every point of the fiber: the latter corresponds to left multiplication \(h\mapsto gh\), which is generally different.
The classification of orders of affine transformations restricts the period-reduction factor to \(p\) or a divisor of \(q-1\), strengthening the divisibility relation obtained from the fiber size alone.

\subsection{Computing the Affine Return Order from the Compression Value}
\label{sec:explicit-affine-return}

For a given compression value \(c\), the dilation coefficient of the return transformation is determined by a norm; when that coefficient is one,
a Frobenius orbit sum further distinguishes the identity return from a nontrivial translation.
\leavevmode\phantomsection\label{def:affine-return-sum}\ignorespaces To express this test, write
\[
[j]_q=\frac{q^j-1}{q-1}\quad(j\ge0),\qquad
R_{q,m}(T)=\sum_{j=0}^{m-1}T^{[j]_q}\quad(m\ge1).
\]
In particular, \([0]_q=0\) and \(R_{q,1}=1\).

\begin{proposition}[Determining the return order from the compression value]\label{prop:explicit-affine-return}
Let \(c\ne0\) have degree \(t\) over \(k=\F_q\), choose any root \(\alpha\) of \(\Theta_q(X)-c\), and put
\[
a=\Nm_{\F_{q^t}/k}(c)=c^{[t]_q},
\qquad z=\alpha^q-\alpha.
\]
Then the unique affine transformation satisfying \(\sigma^t(\alpha)=g(\alpha)\) is
\[
g(X)=aX+b,\qquad b=(1-a)\alpha+zR_{q,t}(c)\in k.
\]
Its order depends only on \(c\), and is given by
\begin{equation}\label{eq:explicit-affine-return-order}
\boxed{\quad e(c)=\ord(g)=
\begin{cases}
\ord(a),&a\ne1,\\
p,&a=1,\ R_{q,t}(c)\ne0,\\
1,&a=1,\ R_{q,t}(c)=0.
\end{cases}\quad}
\end{equation}
The first line is the multiplicative order of \(a\) in \(k^\times\), and
\[
[k(\alpha):k]=t\,e(c).
\]
\end{proposition}
\begin{proof}
Since \(c\ne0\), we have \(\alpha\notin k\) and \(z\ne0\). The equality \(\Theta_q(\alpha)=c\) gives
\(z^{q-1}=c\), so \(z^q=cz\). Successive \(q\)-th powers and induction give
\[
z^{q^j}=c^{[j]_q}z\qquad(j\ge0).
\]
On the other hand, \(c^{q^t}=c\) ensures that \(\sigma^t(\alpha)\) remains in the same compressed fiber.
Proposition~\ref{prop:affine-orbit-separation} gives a unique \(g(X)=a_0X+b\), with \(a_0\in k^\times,b\in k\).
Taking the \(q\)-th power of \(\alpha^{q^t}=a_0\alpha+b\) and subtracting gives \(z^{q^t}=a_0z\).
Comparison with the preceding formula yields \(a_0=c^{[t]_q}=a\).

Next sum the differences of consecutive Frobenius iterates, canceling the intermediate terms:
\begin{equation}\label{eq:affine-return-telescoping}
\alpha^{q^t}-\alpha
=\sum_{j=0}^{t-1}(\alpha^q-\alpha)^{q^j}
=zR_{q,t}(c).
\end{equation}
Thus \(b=\alpha^{q^t}-a\alpha=(1-a)\alpha+zR_{q,t}(c)\); its membership in \(k\) is already ensured by affine orbit separation.

If \(a\ne1\), a translation conjugates \(g\) to \(X\mapsto aX\), so its order is \(\ord(a)\).
If \(a=1\), then \(b=zR_{q,t}(c)\). Since \(z\ne0\), we have \(b=0\) if and only if \(R_{q,t}(c)=0\);
in this case the return is the identity, and otherwise it is a nontrivial translation of order \(p\).
Finally, Proposition~\ref{prop:affine-return-order} gives \([k(\alpha):k]=t\ord(g)\).
The expression for \(e(c)\) does not involve the chosen root \(\alpha\), so all roots in the fiber have the same total period.
\end{proof}

\begin{corollary}[A trace test when the norm is one]\label{cor:affine-return-trace}
Retain the notation above. If \(\Nm_{\F_{q^t}/k}(c)=1\), we may choose \(w\in\F_{q^t}^\times\) with \(w^{q-1}=c\), and
\[
e(c)=
\begin{cases}
1,&\Tr_{\F_{q^t}/k}(w)=0,\\
p,&\Tr_{\F_{q^t}/k}(w)\ne0.
\end{cases}
\]
This test is independent of the choice of \(w\). In particular, when \(q=2\),
\[
cR_{2,t}(c)=\Tr_{\F_{2^t}/\F_2}(c),
\]
so \(e(c)=1\) or two according as the absolute trace of \(c\) is zero or one.
\end{corollary}
\begin{proof}
Multiplicative Hilbert 90 gives a nonzero solution of \(w^{q-1}=c\). Iterating \(w^q=cw\) gives
\[
\Tr_{\F_{q^t}/k}(w)=wR_{q,t}(c).
\]
Since \(w\ne0\), vanishing of the trace is equivalent to \(R_{q,t}(c)=0\).
The ratio of any two such \(w\) satisfies \(u^{q-1}=1\), hence lies in \(k^\times\), and does not change whether the trace vanishes.
For \(q=2\), we may take \(w=c\) directly, or expand \(cR_{2,t}(c)=\sum_{j=0}^{t-1}c^{2^j}\) to obtain the last formula.
\end{proof}

Thus compression retains the period \(t\), and the internal period removed by compression can be computed directly from the norm of \(c\) and the orbit sum above.
The trace and norm of Chapter~\ref{chap:trace} control translation returns and dilation returns, respectively.
The polynomial recurrence for \(R_{q,t}\) and its connection with exact-period layers will be developed in Section~\ref{sec:period-R-recursion}.

\Needspace{12\baselineskip}
\subsection{Lifting Irreducible Factors}
\label{sec:compressed-factor-lifting}
\begin{corollary}[Lifting compressed factors and decomposing fibers]\label{cor:compressed-fiber-lifting}
Retain the notation of Proposition~\ref{prop:affine-return-order}, put \(e=n/t=\ord(g)\),
and let \(H\in k[T]\) be the monic minimal polynomial of \(c\). Then:
\begin{enumerate}
\item Over \(k(c)\), the polynomial \(\Theta_q(X)-c\) splits into \(q(q-1)/e\) distinct monic irreducible factors, each of degree \(e\).
\item Over \(k\), the polynomial \(H(\Theta_q(X))\) splits into \(q(q-1)/e\) distinct monic irreducible factors, each of degree \(n\).
\end{enumerate}
In particular, if \(H\mid G_{q,n}\) and \(\deg H=t\), the number of lifted factors is \(q(q-1)t/n\),
and every lifted factor has degree \(n\).
\end{corollary}
\begin{proof}
Since \(\alpha\notin k\), we have \(c\ne0\). Proposition~\ref{prop:affine-orbit-separation} gives
\(q(q-1)\) distinct roots of \(\Theta_q(X)-c\), namely \(\operatorname{AGL}_1(k)\alpha\).
This polynomial also has degree \(q(q-1)\), so these are all its roots, each simple.
The Frobenius of the base field \(k(c)=\F_{q^t}\) is \(\sigma^t\). Under the fiber identification above,
it corresponds to right multiplication \(h\mapsto hg\) on the group, every cycle having length \(e\).
The correspondence between root orbits and irreducible factors immediately gives the first assertion.

The roots of \(H(\Theta_q)\) lie in disjoint fibers over \(c,c^q,\ldots,c^{q^{t-1}}\).
Every root in each fiber is simple. By Lemma~\ref{lem:fiber-return},
each return orbit of length \(e\) extends to a total orbit of length \(te=n\), and return orbits correspond bijectively to total orbits.
Thus the number of factors is unchanged, and every degree is multiplied by \(t\), proving the second assertion.
Starting instead from \(H\mid G_{q,n}\), choose a root \(c\) and adjoin a root \(\alpha\) of \(\Theta_q(X)-c\).
The identity \(E_{q,n}=G_{q,n}(\Theta_q)\) shows that \(\alpha\) has degree exactly \(n\), so the conclusions above still apply.
\end{proof}

\begin{corollary}[Compressed factors and affine orbits of irreducible polynomials]\label{v129:cor-affine-factor-orbits}
Let \(n>1\) and \(G=\operatorname{AGL}_1(k)\). For \(u(X)=aX+b\in G\) and \(f\in\Irr_q(n)\), define
\[
(u\cdot f)(X)=a^n f\bigl((X-b)/a\bigr).
\]
Choose a root \(\alpha\) of \(f\), put \(c=\Theta_q(\alpha)\), \(H=m_c\), and \(t=\deg H\), and let \(g\in G\) be the unique element with \(\alpha^{q^t}=g(\alpha)\). Then
\begin{equation}\label{v129:eq-affine-factor-orbits}
\operatorname{Stab}_G(f)=\langle g\rangle,
\qquad
H(\Theta_q(X))=\prod_{u\langle g\rangle\in G/\langle g\rangle}(u\cdot f)(X).
\end{equation}
The factors are distinct, and \(|\operatorname{Stab}_G(f)|=n/t\). The map \(f\mapsto H\) induces a bijection
\[
\Irr_q(n)/G\ \longrightarrow\
\{H\in k[T]:H\text{ monic irreducible and }H\mid G_{q,n}\}.
\]
\end{corollary}
\begin{proof}
The root set of \(u\cdot f\) is \(u(\mathcal O_\alpha)\). Since \(u\) commutes with Frobenius, \(u\cdot f=f\) holds exactly when \(u(\alpha)=\sigma^j(\alpha)\) for some \(j\).
Applying \(\Theta_q\) gives \(c=c^{q^j}\), so \(t\mid j\). Write \(j=t\ell\); then
\(u(\alpha)=g^\ell(\alpha)\). The affine action is free at \(\alpha\notin k\), so \(u=g^\ell\). The reverse inclusion holds as well, proving the stabilizer formula.
Distinct left cosets \(u\langle g\rangle\) therefore give distinct irreducible factors. Their root sets have union
\(\bigcup_{i=0}^{t-1}G\sigma^i(\alpha)\), exactly the root set of \(H(\Theta_q)\).
All these roots are simple and both sides are monic, proving the product formula.

Replacing the root of \(f\) replaces \(c\) by a conjugate, while applying an element of \(G\) leaves \(c\) unchanged, so the quotient map is well defined.
Every root of a factor \(H\mid G_{q,n}\) lifts to a root of \(E_{q,n}\), proving surjectivity. The product formula places all irreducible factors above the same \(H\) in one \(G\)-orbit, proving injectivity.
\end{proof}

Thus the compressed-factor degree \(t\), the original-factor degree \(n\), and the order \(e\) of its affine stabilizer over the base field satisfy \(n=te\).
The number \(q(q-1)/e\) in Corollary~\ref{cor:compressed-fiber-lifting} counts both return cycles and the polynomials in an affine orbit.
This stabilizer uses affine transformations with coefficients in \(k\); allowing coefficients in an extension can enlarge it, as in Chapter~\ref{chap:root-symmetry}.

The order of the return transformation thus gives both the relative extension degree and the factor degrees of the fiber polynomial.
In affine compression, the general fiber factorization of Theorem~\ref{thm:fiber-factorization} becomes a case in which all return orbits have equal length;
the same integer \(e\) therefore determines both the number and the degrees of the factors.

\begin{corollary}[Computing the lifted factorization type from a compressed factor]\leavevmode\phantomsection\label{cor:explicit-compressed-lifting}
Let \(H(T)=T^t+h_{t-1}T^{t-1}+\cdots+h_0\in k[T]\) be monic irreducible, with \(H\ne T\). Put
\[
a=(-1)^t h_0\in k^\times,
\qquad r(T)=R_{q,t}(T)\bmod H(T),\quad\deg r<t.
\]
Use \(a\) and whether \(r=0\) in~\eqref{eq:explicit-affine-return-order} to determine the integer \(e\).
Then \(H(\Theta_q(X))\) has exactly \(q(q-1)/e\) distinct monic irreducible factors, each of degree \(te\).
For \(q>2\), this composition is always reducible; for \(q=2\), it is irreducible if and only if \(h_{t-1}=1\).
\end{corollary}
\begin{proof}
Choose a root \(c\) of \(H\). It has degree \(t\), and \(H\ne T\) ensures \(c\ne0\).
The product of the roots gives \(\Nm_{k(c)/k}(c)=(-1)^th_0=a\).
Since \(H\) is the minimal polynomial of \(c\), we have \(R_{q,t}(c)=0\) if and only if \(H\mid R_{q,t}\), or \(r=0\).
Proposition~\ref{prop:explicit-affine-return} therefore gives degree \(te\) for each lifted root;
Corollary~\ref{cor:compressed-fiber-lifting} then gives the number of factors and the absence of repeated roots.

If \(q>2\), the return order is \(p\) or a divisor of \(q-1\), both strictly smaller than \(q(q-1)\), so there is more than one factor.
If \(q=2\), the fiber has size two, so irreducibility of the composition is equivalent to \(e=2\).
Corollary~\ref{cor:affine-return-trace} reduces this to \(\Tr_{k(c)/\F_2}(c)=1\).
The sum of the roots shows that this trace equals \(h_{t-1}\), giving the final criterion.
\end{proof}

This corollary determines a base-field element from the constant term, computes its multiplicative order, and, when needed, takes one polynomial remainder to determine the degrees and number of all lifted factors.
In the binary case, it agrees with the Artin--Schreier composition criterion of Chapter~\ref{chap:root-fibers}.

\begin{example}[A period four compressed to period two]\label{exm:compression-return-four}
Choose \(\alpha\in\F_{16}\) satisfying \(\alpha^4+\alpha+1=0\), and put \(\beta=\alpha^2+\alpha\).
Example~\ref{exm:quartic-model} already showed that \(\alpha\) has degree four over \(\F_2\).
From \(\alpha^4=\alpha+1\), we obtain \(\beta^2=\beta+1\), so \(\beta\in\F_4\setminus\F_2\) and the compressed period is two.
The original orbit passes successively through \(\alpha,\alpha^2,\alpha+1,\alpha^2+1\), and the two fibers are
\[
F_\beta=\{\alpha,\alpha+1\},\qquad F_{\beta^2}=\{\alpha^2,\alpha^2+1\}.
\]
The first-return map \(\sigma^2\) is the translation \(x\mapsto x+1\) of order two on each fiber.
Thus the total period is \(4=2\cdot2\), and \(X^2+X-\beta\) is irreducible over \(\F_4\).
For \(H(T)=T^2+T+1\), Corollary~\ref{cor:compressed-fiber-lifting} gives
the irreducibility of \(H(X^2+X)=X^4+X+1\) over \(\F_2\).

Now take the quotient by the full translation space \(\F_4\) from Example~\ref{ffs:ex-quartic-center}.
Since \(L_{\F_4}(\alpha)=\alpha^4+\alpha=1\), all four roots map to the fixed point \(1\) in the quotient, and the fiber is \(\alpha+\F_4\).
Still taking \(\sigma(x)=x^2\), and using \(\alpha\) as base point, the coordinate expression of \(\sigma\) on this fiber is
\[
R(v)=v^2+\beta,\qquad v\in\F_4.
\]
Since \(\beta^2+\beta=1\), we have
\[
R^2(v)=v+1,\qquad R^4(v)=v.
\]
The map \(R^2\) has no fixed point, so the four points form a single return cycle of length four.
The two translation quotients satisfy the composition relation
\[
L_{\F_4}(X)=X^4+X
=\bigl(L_{\F_2}(X)\bigr)^2+L_{\F_2}(X).
\]
Thus taking the quotient by \(\F_2\)-translations and then substituting \(Y\mapsto Y^2+Y\) gives the same result as taking the quotient by \(\F_4\)-translations.
The two period decompositions of the same original orbit are
\begin{center}\small
\renewcommand{\arraystretch}{1.15}
\begin{tabular}{@{}lccc@{}}
\toprule
Translation group&Quotient period&Internal return period&Original period\\
\midrule
\(\F_2\)&\(2\)&\(2\)&\(4=2\cdot2\)\\
\(\F_4\)&\(1\)&\(4\)&\(4=1\cdot4\)\\
\bottomrule
\end{tabular}
\end{center}
The quotient period and the internal return period depend on the chosen translation group, while their product is determined by the original orbit length.
The full affine stabilizer corresponds to the quotient polynomial \(U^3=(X^4+X+1)^3\); see Proposition~\ref{ffs:prop-invariant-ring}.
\end{example}

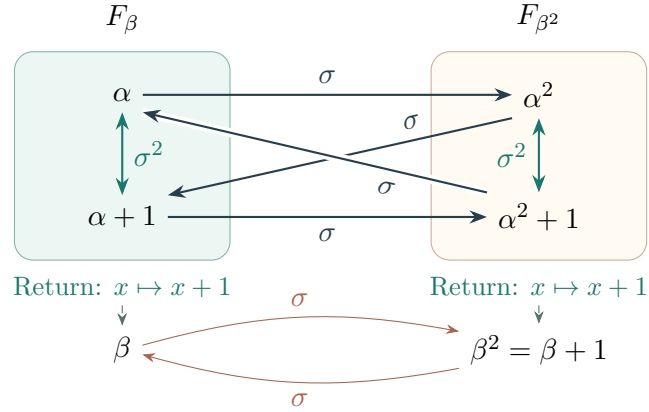
\begin{figure}[htbp]
\centering
\begin{tikzpicture}[>=Stealth,scale=0.95,every node/.style={font=\small}]
\filldraw[fill=CoverMint,draw=Jade!55,rounded corners=7pt] (-4.4,-1.45) rectangle (-1.4,1.45);
\filldraw[fill=CoverPaper,draw=Clay!55,rounded corners=7pt] (1.4,-1.45) rectangle (4.4,1.45);
\node at (-2.9,1.9) {\(F_\beta\)};
\node at (2.9,1.9) {\(F_{\beta^2}\)};
\node (a) at (-2.9,0.85) {\(\alpha\)};
\node (b) at (-2.9,-0.85) {\(\alpha+1\)};
\node (d) at (2.9,0.85) {\(\alpha^2\)};
\node (e) at (2.9,-0.85) {\(\alpha^2+1\)};
\draw[->,InkSlate,thick] (a) -- node[above] {\(\sigma\)} (d);
\draw[->,InkSlate,thick] (b) -- node[below] {\(\sigma\)} (e);
\draw[->,InkSlate,thick] (d.south west) -- node[pos=0.28,above,sloped] {\(\sigma\)} (b.north east);
\draw[->,InkSlate,thick,preaction={draw=white,line width=3.4pt,arrows=-}]
 (e.north west) -- node[pos=0.28,below,sloped] {\(\sigma\)} (a.south east);
\draw[<->,Jade,thick] (a) -- node[right] {\(\sigma^2\)} (b);
\draw[<->,Jade,thick] (d) -- node[left] {\(\sigma^2\)} (e);
\node[font=\footnotesize,text=Jade] at (-2.9,-1.82) {Return: \(x\mapsto x+1\)};
\node[font=\footnotesize,text=Jade] at (2.9,-1.82) {Return: \(x\mapsto x+1\)};
\node (c) at (-2.9,-2.7) {\(\beta\)};
\node (c2) at (2.9,-2.7) {\(\beta^2=\beta+1\)};
\draw[->,MutedSlate,dashed] (-2.9,-2.08) -- (c);
\draw[->,MutedSlate,dashed] (2.9,-2.08) -- (c2);
\draw[->,Clay] (c) to[bend left=13] node[above] {\(\sigma\)} (c2);
\draw[->,Clay] (c2) to[bend left=13] node[below] {\(\sigma\)} (c);
\end{tikzpicture}
\caption{The four upper \(\sigma\)-arrows form a single orbit of period four; \(\sigma^2\) interchanges the two points in each fiber. Compression retains the period two below, and the internal return period is two, so the original period is \(4=2\cdot2\).}
\label{fig:compression-return}
\end{figure}

\Needspace{11\baselineskip}
\begin{corollary}[Restrictions on compressed factor degrees]\label{cor:compressed-factor-degrees}
Let \(n>1\). If \(\gcd(n,p(q-1))=1\), every irreducible factor of \(G_{q,n}\) in \(k[T]\) has degree \(n\).
When \(q=2\), each monic irreducible factor \(H\mid G_{2,n}\) satisfies
\[
\boxed{\quad \deg H=n\quad\text{or}\quad\deg H=\frac n2.\quad}
\]
When \(n\) is odd, only \(\deg H=n\) is possible.
\end{corollary}
\begin{proof}
Choose a root \(c\) of an irreducible factor of \(G_{q,n}\) in \(k[T]\), and write \(t=[k(c):k]\).
In a finite extension, choose \(\Theta_q(\alpha)=c\). The identity \(E_{q,n}=G_{q,n}(\Theta_q)\) gives \([k(\alpha):k]=n\).
Proposition~\ref{prop:affine-return-order} therefore applies.
If \(\gcd(n,p(q-1))=1\), then \(n/t\) divides \(n\) and is either \(p\) or a divisor of \(q-1\), so it must be one.
If \(q=2\), then \(p=2\) and \(q-1=1\), so \(n/t\) can only be one or two.
Here \(t=\deg H\); if \(n\) is odd, \(n/2\) is not an integer, so only \(t=n\) is possible.
\end{proof}

\Needspace{30\baselineskip}
\begin{example}[Dilation and translation returns over the ternary field]\label{exm:ternary-affine-return}
Over \(k=\F_3\), we have \(\Theta_3(X)=(X^3-X)^2\). Choose
\[
\alpha^2+1=0,\qquad \beta^3-\beta-1=0.
\]
Neither polynomial has a root in \(\F_3\), so they are irreducible of degrees two and three, respectively.
Thus \([k(\alpha):k]=2\) and \([k(\beta):k]=3\).

From \(\alpha^2=-1\), we obtain \(\alpha^3=-\alpha\), so
\[
\Theta_3(\alpha)=(-2\alpha)^2=\alpha^2=2.
\]
The compression value lies in the base field and has period one. The return element determined by the base point \(\alpha\) is the dilation \(g(X)=-X\),
of order two; hence the period-reduction factor is two, dividing \(q-1\).
Take \(h(X)=X+1\). Then \(h(\alpha)=\alpha+1\) lies in the same full fiber, but
\[
(\alpha+1)^3=-\alpha+1\ne-\alpha-1=g(\alpha+1).
\]
Here the Frobenius return gives \((hg)(\alpha)\), whereas applying \(g\) gives \((gh)(\alpha)\), explicitly distinguishing right from left multiplication.
On the other hand, \(\beta^3=\beta+1\) gives \(\Theta_3(\beta)=1\).
Here the return element is the translation \(g(X)=X+1\), of order three, and the period-reduction factor is \(3=p\).
The two cases are compared below.
\begin{center}\small
\begin{tabular}{@{}ccccc@{}}
\toprule
Root&Original period \(n\)&Value \(c\)&Return \(g(X)\)&Order \(e=n/t\)\\
\midrule
\(\alpha\)&\(2\)&\(2\)&\(-X\)&\(2\)\\
\(\beta\)&\(3\)&\(1\)&\(X+1\)&\(3\)\\
\bottomrule
\end{tabular}
\end{center}

The fiber of each nonzero compression value has \(q(q-1)=6\) roots, and \(t=1\).
Corollary~\ref{cor:compressed-fiber-lifting} therefore gives three irreducible quadratic factors for \(\Theta_3(X)-2\)
and two irreducible cubic factors for \(\Theta_3(X)-1\). The complete factorizations are
\begin{align*}
\Theta_3(X)-2
&=(X^2+1)(X^2+X+2)(X^2+2X+2),\\
\Theta_3(X)-1
&=(X^3-X-1)(X^3-X+1).
\end{align*}
Expansion of the right-hand sides verifies the identities; none of the listed quadratic or cubic factors has a root in \(\F_3\), so all are irreducible.
Both fibers have size six, but the return orders are two and three, giving three and two orbits, respectively.
This illustrates both sources of affine return orders and how the return order controls the degrees and number of factors.
\end{example}

\subsection{Return Mechanisms for Power and Affine Compression}
\label{sec:power-affine-comparison}

The power map of Chapter~\ref{chap:cyclotomic} and the affine invariant of this chapter both send one Frobenius orbit to another.
The image period records motion between fibers, while the ratio of the original period to the image period records the first-return period within a fiber.
Lemma~\ref{lem:fiber-return} applies in both cases; the specific return maps on the fibers must be determined separately.

\noindent\textbf{Power compression.}
In the notation of Section~\ref{sec:cyclotomic-med}, let \(C=C_i\) and \(\alpha=\zeta^i\),
take \(r=\kappa_C\), \(t=\omega_C\), and put \(\beta=\alpha^r\).
Theorem~\ref{prop:med-binomial-descent} gives
\[
[k(\alpha):k]=\tau=rt,\qquad [k(\beta):k]=t.
\]
On the conjugate root set of \(\alpha\), the map \(x\mapsto x^r\) sends each complete equal-difference root block of size \(r\) to a single value.
The fiber over \(\beta\) is precisely \(\alpha\mu_r\), and \(\sigma^t\) forms a single cycle of length \(r\) on this fiber.
Use \(\varepsilon\mapsto\alpha\varepsilon\) as fiber coordinates, and put
\[
Q=q^t,\qquad \delta=\alpha^{Q-1}.
\]
Since \(\beta^Q=\beta\), we have \(\delta^r=1\), so \(\delta\in\mu_r\). Direct calculation gives the first-return map
\begin{equation}\label{eq:power-compression-return}
\sigma^t(\alpha\varepsilon)=\alpha\delta\varepsilon^Q,
\qquad \mathcal R(\varepsilon)=\delta\varepsilon^Q
\quad(\varepsilon\in\mu_r).
\end{equation}
Since the original root set is a single orbit of length \(rt\), the least period of \(\mathcal R\) starting from \(1\) is exactly \(r\).
Thus the period-reduction factor here equals the size of the entire root block.

\noindent\textbf{Affine compression.}
For \(c=\Theta_q(\alpha)\ne0\), put \(n=[k(\alpha):k]\) and \(t=[k(c):k]\), and write \(\sigma^t(\alpha)=g(\alpha)\), with \(e=\ord(g)=n/t\). Equation~\eqref{eq:full-fiber-orbit-intersection} shows that the full fiber has \(q(q-1)\) elements, whereas its intersection with one original orbit has only \(e\); the first return is right multiplication \(h\mapsto hg\) in group coordinates.
Thus the period-reduction factor is the return period of the chosen root and need not equal the size of the full fiber. In power compression above, the whole of \(\alpha\mu_r\) lies in one original orbit, so the two sizes agree.

Both cases satisfy
\[
[k(\alpha):k]
=[k(\pi(\alpha)):k]\,[k(\alpha):k(\pi(\alpha))],
\]
where \(\pi\) is the power map or \(\Theta_q\), respectively, and the relative degree is precisely the return period of the chosen root within the fiber.
The internal period is \(\kappa_C\) for power compression and \(\ord(g)\) for affine compression.

The forms of the two return maps differ.
The map \(\varepsilon\mapsto\varepsilon^Q\) in~\eqref{eq:power-compression-return} is the identity only when \(r\mid Q-1\);
in general, \(\mathcal R\) cannot simply be written as multiplication by \(\delta\). Its iterates are in fact
\[
\mathcal R^j(\varepsilon)=\delta^{1+Q+\cdots+Q^{j-1}}\varepsilon^{Q^j}
\qquad(j\ge1),
\]
as follows by successive composition.
For example, when \(q=5,N=32,C=C_1\), we have \(\tau=8\), \(\kappa_C=8\), and \(\omega_C=1\).
Take \(\alpha\) to be a primitive 32nd root of unity, and put \(\xi=\alpha^4\). Then \(\mu_8=\langle\xi\rangle\),
and the return map is \(\varepsilon\mapsto\xi\varepsilon^5\). In coordinates \(\varepsilon=\xi^j\), it becomes
\[
j\longmapsto 1+5j\pmod8.
\]
Exercise~\ref{ex:wu-compressed-separable} gives the full return calculation and compares this permutation with multiplication by \(\xi\) alone.
In the affine case, the right-multiplication formula follows from the commutativity of the group elements with Frobenius.
Thus the common period multiplication law comes from equivariance of the compression map with Frobenius, while the specific return-order test also depends on the group action on the fiber.

\Needspace{17\baselineskip}
\section{Sparse Recurrences and Base-Field Lifting}
\label{sec:period-R-recursion}
Retain \([j]_q\) and \(R_{q,m}\) from Section~\ref{sec:explicit-affine-return}.
The orbit sum in the return calculation also describes the elements outside the base field in a fixed-point set: first factor out \(Z_q\) to remove the base-field roots, and then use \(R_{q,m}\) to record their compression coordinates.
\begin{proposition}\label{prop:period-R}
For \(R_{q,m}(T)=\sum_{j=0}^{m-1}T^{[j]_q}\), we have
\[
A_{q,m}(X)=Z_q(X)R_{q,m}(\Theta_q(X)),
\]
and
\[
R_{q,1}=1,\qquad R_{q,m+1}=1+T R_{q,m}^q.
\]
For \(n>1\), we also have
\[
G_{q,n}=\prod_{d\mid n}R_{q,d}^{\,\mu(n/d)},\qquad
R_{q,m}=\prod_{\substack{d\mid m\\d>1}}G_{q,d}.
\]
\end{proposition}
\begin{proof}
The Frobenius identity and telescoping give
\[
A_{q,m}=\sum_{j=0}^{m-1}(X^q-X)^{q^j}
=Z_q\sum_{j=0}^{m-1}(Z_q^{q-1})^{[j]_q},
\]
which proves the first formula.
Also, \([j+1]_q=1+q[j]_q\), so
\[
1+T\left(\sum_{j=0}^{m-1}T^{[j]_q}\right)^q
=1+\sum_{j=0}^{m-1}T^{[j+1]_q}=R_{q,m+1}.
\]
Substitute \(A_{q,d}=Z_qR_{q,d}(\Theta_q)\) in the M\"obius product.
The total exponent of \(Z_q\) is zero, giving the product expression for \(G\).
On the other hand, canceling \(E_{q,1}=Z_q\) in~\eqref{eq:period-layer-product} gives
\[
R_{q,m}(\Theta_q)
=\prod_{\substack{d\mid m\\d>1}}G_{q,d}(\Theta_q).
\]
Injectivity of substitution gives the factorization formula for \(R\).
Thus the quotients defined formally above belong to \(k[T]\).
\end{proof}

\begin{corollary}[The first nonconstant term after compression]
Let \(n>1\) and \(s=n/\rad(n)\). Then
\[
\operatorname{ord}_T(G_{q,n}-1)=[s]_q,\qquad
[T^{[s]_q}]G_{q,n}=(-1)^{\omega(n)+1}.
\]
\end{corollary}
\begin{proof}
The lowest nonzero term of \(\Theta_q=(X^q-X)^{q-1}\) is \(X^{q-1}\), with coefficient one.
Thus, for nonzero \(H\) with \(H(0)=0\), we have \(\operatorname{ord}_XH(\Theta_q)=(q-1)\operatorname{ord}_T H\),
and the first nonzero coefficient is unchanged.
Take \(H=G_{q,n}-1\) and apply Theorem~\ref{thm:period-gap}.
\end{proof}

\begin{example}[The compressed degree-four layer]
For \(q=2\), we have
\[
R_{2,2}=1+T,\qquad R_{2,4}=1+T+T^3+T^7.
\]
Hence
\[
G_{2,4}(T)=\frac{1+T+T^3+T^7}{1+T}
=1+T^3+T^4+T^5+T^6.
\]
Substituting \(\Theta_2=X^2+X\) gives \(E_{2,4}=1+X^3+X^6+X^9+X^{12}\).
The degree of the first nonconstant term is determined by \(s=2\) and \([2]_2=3\).
Furthermore, in \(\F_2[T]\),
\[
G_{2,4}(T)=(T^2+T+1)(T^4+T+1).
\]
Denote the two factors by \(H_2=T^2+T+1\) and \(H_4=T^4+T+1\).
For a root \(c\) of \(H_2\), we have \(R_{2,2}(c)=1+c\ne0\), so the return order is two and the total degree is \(2\cdot2=4\).
For a root \(d\) of \(H_4\), the identities \(d^4=d+1\) and \(d^7=d^3+d+1\) give
\[
R_{2,4}(d)=1+d+d^3+d^7=0.
\]
Thus the return order is one, and the total degree is again four. Corollary~\ref{cor:explicit-compressed-lifting} gives one lifted quartic factor and two lifted quartic factors, respectively.
Both factors are irreducible: the quadratic has no base-field root, and irreducibility of the quartic was shown in Example~\ref{exm:quartic-model}. Their degrees two and four are exactly the two compressed periods allowed by Corollary~\ref{cor:compressed-factor-degrees}.

Lifting the two factors separately along \(\Theta_2=X^2+X\) gives
\[
\begin{aligned}
(\Theta_2^2+\Theta_2+1)&=X^4+X+1,\\
(\Theta_2^4+\Theta_2+1)&=
(X^4+X^3+1)(X^4+X^3+X^2+X+1).
\end{aligned}
\]
The first identity corresponds to the period-four orbit in Example~\ref{exm:compression-return-four}:
the compressed period is two, and the return transformation \(X\mapsto X+1\) has order two, giving total period \(4=2\cdot2\).
In the second identity, the compression value still has degree four and the internal return period is one; the two original period-four orbits are interchanged by the translation \(X\mapsto X+1\).
The roots of \(X^4+X^3+1\) have multiplicative order \(15\), while those of \(X^4+X^3+X^2+X+1\) have multiplicative order five.
Thus original elements in the same compressed fiber have the same field degree but need not have the same multiplicative order;
the degree of the compression value itself is determined by the return period.
Accordingly, the three irreducible quartics form one singleton orbit and one two-element orbit under \(\operatorname{AGL}_1(\F_2)\). Their stabilizers have orders two and one, corresponding to compressed-factor degrees two and four.
\end{example}

\Needspace{12\baselineskip}
\begin{proposition}[Compatibility of base-field lifting]\label{prop:period-base-change}
For \(s\ge1\), define
\[
\tau_{q,s}(T)=T^{[s]_q}R_{q,s}(T)^{q^s-1}.
\]
Then, for \(s,t\ge1\),
\[
\begin{aligned}
\Theta_{q^s}(X)&=\tau_{q,s}(\Theta_q(X)),\\
R_{q,st}(T)&=R_{q,s}(T)R_{q^s,t}(\tau_{q,s}(T)),\\
\tau_{q,st}&=\tau_{q^s,t}\circ\tau_{q,s}.
\end{aligned}
\]
If \(n>1\), \(r=\rad(n)\), and \(s=n/r\), then
\[
G_{q,n}=G_{q^s,r}\circ\tau_{q,s}.
\]
\end{proposition}
\begin{proof}
From \(A_{q,s}=Z_qR_{q,s}(\Theta_q)\), we obtain
\[
\Theta_{q^s}
=A_{q,s}^{q^s-1}
=\Theta_q^{[s]_q}R_{q,s}(\Theta_q)^{q^s-1},
\]
which gives the first identity.
Next use
\[
\frac{A_{q,st}}{A_{q,1}}
=\frac{A_{q,s}}{A_{q,1}}\,
  \frac{A_{q^s,t}}{A_{q^s,1}},
\]
apply the expression for \(R\) over the two base fields, and substitute the first identity.
Injectivity of \(T\mapsto\Theta_q\) gives the second identity.
For the third, express \(\Theta_{q^{st}}\) in terms of \(\Theta_q\) either directly
or through \(\Theta_{q^s}\). The results agree, and injectivity of substitution gives the claim.
Finally, the square-free-kernel formula gives
\[
G_{q,n}(\Theta_q)=E_{q,n}
=E_{q^s,r}=G_{q^s,r}(\tau_{q,s}(\Theta_q)).
\]
Injectivity of substitution again yields the last identity.
\end{proof}

\section{Binary Moore Forms and Homogeneous Period Layers}
Here \emph{binary} refers to forms in two variables. The affine transformations considered above are represented in the ratio coordinate \(X/Y\) by matrices
\(\bigl(\begin{smallmatrix}a&b\\0&1\end{smallmatrix}\bigr)\).
Allowing any \(g=\bigl(\begin{smallmatrix}a&b\\c&d\end{smallmatrix}\bigr)\in\operatorname{GL}_2(k)\)
transforms the ratio coordinate to \((aX+bY)/(cX+dY)\).
For \(\alpha\notin k\), the denominator \(c\alpha+d\) is nonzero; the transformation and its inverse are defined over \(k\), so
\[
k\!\left(\frac{a\alpha+b}{c\alpha+d}\right)=k(\alpha).
\]
Thus each exact-period layer \(\mathcal P_{q,n}\) with \(n>1\) is stable under these fractional linear transformations.
Using two variables expresses this larger symmetry through polynomial identities without denominators and turns the preceding univariate recurrence into compatibility relations under the linear group. Put
\[
\Lambda_{q,m}(X,Y)=X^{q^m}Y-XY^{q^m}\quad(m\ge0),
\qquad
\Delta_q=\Lambda_{q,1}^{\,q-1},\quad
\Gamma_q=\frac{\Lambda_{q,2}}{\Lambda_{q,1}}.
\]
The expression \(\Gamma_q\) is indeed a polynomial: we have \(\Lambda_{q,1}=Y\prod_{a\in k}(X-aY)\),
and each distinct linear factor divides \(\Lambda_{q,2}\).
For any invertible matrix \(g\in\operatorname{GL}_2(k)\), expansion gives
\[
\Lambda_{q,m}(g(X,Y))=\det(g)\Lambda_{q,m}(X,Y).
\]
Thus \(\Delta_q,\Gamma_q\) are invariant under this linear group, with degrees \(q^2-1,q^2-q\), respectively;
they are called the binary Dickson invariants.

\begin{proposition}[A binary Frobenius recurrence]\label{prop:binary-moore-recursion}
Define
\[
Q_{q,0}=0,\qquad Q_{q,1}=1,\qquad
Q_{q,m+2}(U,V)=V^{q^m}Q_{q,m+1}(U,V)-U^{q^m}Q_{q,m}(U,V).
\]
Then, for \(m\ge0\),
\[
\Lambda_{q,m}=\Lambda_{q,1}Q_{q,m}(\Delta_q,\Gamma_q),
\]
and, for \(m\ge1\),
\[
R_{q,m}(T)=Q_{q,m}(T,1+T).
\]
\end{proposition}
\begin{proof}
For \(Z=X\) and \(Z=Y\), direct expansion and cancellation of equal terms give
\[
\Lambda_{q,1}Z^{q^2}-\Lambda_{q,2}Z^q+\Lambda_{q,1}^{\,q}Z=0.
\]
Dividing by the nonzero polynomial \(\Lambda_{q,1}\) in the fraction field gives
\[
Z^{q^2}-\Gamma_qZ^q+\Delta_q Z=0.
\]
After taking \(q^m\)-th powers, multiply the identities for \(Z=X,Y\) by \(Y,X\), respectively, and subtract to obtain
\[
\Lambda_{q,m+2}
=\Gamma_q^{q^m}\Lambda_{q,m+1}
-\Delta_q^{q^m}\Lambda_{q,m}.
\]
This recurrence and the initial values \(\Lambda_{q,0}=0,\Lambda_{q,1}\) agree with the definition of \(Q_{q,m}\),
so induction gives the first identity.
When \(Y=1\), we have \(\Delta_q(X,1)=\Theta_q(X)\)
and \(\Gamma_q(X,1)=1+\Theta_q(X)\).
The first identity becomes
\[
A_{q,m}=Z_qQ_{q,m}(\Theta_q,1+\Theta_q).
\]
Compare this with \(A_{q,m}=Z_qR_{q,m}(\Theta_q)\), cancel the nonzero factor \(Z_q\),
and apply injectivity of substitution to obtain the second identity.
\end{proof}

\leavevmode\phantomsection\label{def:polynomial-homogenization}\ignorespaces For \(n>1\), write \(E_{q,n}(X)=\sum_{i=0}^{D}a_iX^i\), where \(D=D_{q,n}\). Its \textbf{homogenization} is the homogeneous polynomial of total degree \(D\)
\[
\mathcal E_{q,n}(X,Y)=\sum_{i=0}^{D}a_iX^iY^{D-i}
=Y^{D}E_{q,n}(X/Y).
\]
The expression on the right has no denominator after expansion.
From \(\Lambda_{q,d}=Y^{q^d+1}A_{q,d}(X/Y)\),
\(\sum_{d\mid n}\mu(n/d)=0\), and the degree formula, we obtain
\[
\mathcal E_{q,n}
=\prod_{d\mid n}\Lambda_{q,d}^{\,\mu(n/d)}.
\]
The left-hand side ensures that the right-hand side is a polynomial, and it is invariant under \(\operatorname{GL}_2(k)\):
the total exponent of \(\det(g)\) is zero.
The recurrence directly gives the first few layers:
\[
\begin{aligned}
\mathcal E_{q,2}&=\Gamma_q,\\
\mathcal E_{q,3}&=\Gamma_q^{q+1}-\Delta_q^q,\\
\mathcal E_{q,4}&=
\Gamma_q^{q^2-1}(\Gamma_q^{q+1}-\Delta_q^q)-\Delta_q^{q^2}.
\end{aligned}
\]
For example, the third identity follows from \(\mathcal E_{q,4}=\Lambda_{q,4}/\Lambda_{q,2}\)
and \(Q_{q,4}=V^{q^2}(V^{q+1}-U^q)-U^{q^2}V\).
Setting \(Y=1\) recovers the corresponding univariate period-layer and affine-compression formulas.

Exact-period layers record degree conditions on roots; the M\"obius product records cancellation on the divisor lattice of the index;
affine compression records symmetry of the root set under translations and dilations;
the binary recurrence expresses the same Frobenius relation as invariants under linear transformations.
Root sets and the corresponding group actions connect these descriptions.

\section{Period Counts and Euler Products}
\label{sec:orbit-euler}

Exact-period layers can also be encoded together in a formal power series. The series in this section have coefficients in \(\Q\),
so \(1/r\), formal exponentials, and formal logarithms have their usual meanings.

\begin{proposition}\leavevmode\phantomsection\label{prop:orbit-euler}
Let \(\sigma\) permute a finite set \(S\), put \(N_r=|\Fix(\sigma^r)|\), and let \(c_d\) be the number of orbits of length \(d\).
Then
\[
Z_S(u):=\exp\left(\sum_{r\ge1}\frac{N_r}{r}u^r\right)
=\prod_{d\ge1}(1-u^d)^{-c_d}
=\det(I-uP_\sigma)^{-1},
\]
where \(P_\sigma\) is the permutation matrix of \(\sigma\) on the \(\Q\)-vector space with basis \(S\).
\end{proposition}
\begin{proof}
An orbit of length \(d\) contributes \(d\) to \(N_r\) when \(d\mid r\), and zero otherwise.
Thus \(N_r=\sum_{d\mid r}dc_d\), giving
\[
\sum_{r\ge1}\frac{N_r}{r}u^r
=\sum_{d\ge1}c_d\sum_{j\ge1}\frac{u^{dj}}j
=-\sum_{d\ge1}c_d\log(1-u^d).
\]
Taking formal exponentials gives the product formula. The permutation matrix of a \(d\)-cycle satisfies \(\det(I-uP)=1-u^d\), as follows from its characteristic polynomial \(T^d-1\).
Multiplying the determinants of the cycle blocks gives the last expression.
\end{proof}

If \(\pi:S\to O\) satisfies Lemma~\ref{lem:fiber-return}, with \(|O|=n\) and \(\tau=\sigma^n|_{\pi^{-1}(a)}\), then
\begin{equation}\label{v129:eq-return-euler}
Z_{S,\sigma}(u)=Z_{\pi^{-1}(a),\tau}(u^n).
\end{equation}
The subscript specifies the permutation used for counting. Each return cycle of length \(j\) contributes \((1-u^{nj})^{-1}\), proving the identity.
This expresses the factor-degree counts of Chapter~\ref{chap:root-fibers} as one product. After a change of base field, Corollary~\ref{v129:cor-return-base-change} describes the splitting of each product factor. Root multiplicities do not enter these permutation data and must still be retained separately.

Proposition~\ref{prop:finite-field-euler} already gave the unweighted finite-field identity~\eqref{eq:finite-field-euler}.
Only finitely many Frobenius orbits occur in each fixed degree, so the calculation above, regrouped by periods, still applies.
The determinant expression for a finite permutation set further shows that periodic-point counts can be encoded together by a linear operator.
Here the permutation matrix acts on a vector space with evaluation points as a basis and a coefficient field of characteristic zero; it should be distinguished from the \(k\)-linear Frobenius operator on the extension field itself.

Assigning a complex weight \(w_{\mathcal O}\) to each orbit \(\mathcal O\), the same calculation gives
\[
\prod_{\mathcal O}(1-w_{\mathcal O}u^{|\mathcal O|})^{-1}
=\exp\left(\sum_{r\ge1}\frac{u^r}{r}
\sum_{\substack{\mathcal O\\|\mathcal O|\mid r}}
|\mathcal O|w_{\mathcal O}^{r/|\mathcal O|}\right).
\]
Setting all weights to one recovers ordinary point counts; applying a multiplicative character to orbit products and an additive character to orbit sums
gives the weighted product in Theorem~\ref{thm:hasse-davenport}.
The proof also applies to infinitely many orbits, provided that only finitely many have any given length.
Thus irreducible-polynomial counts and Gauss-sum lifting use the same principle of regrouping sums by orbits.

\par\Needspace{11\baselineskip}
\section*{Exercises}
\phantomsection
\addcontentsline{toc}{section}{Exercises}
\begin{exercise}[The degree-six exact-period layer]\label{ex:wu-E26}
Write the M\"obius quotient for \(E_{2,6}\), and find its degree and the number of irreducible sextic factors it contains.
Express it as a product of reduced cyclotomic polynomials, and find its remainder modulo \(X^3\).
\end{exercise}
\begin{exercise}[Changing the period step]\label{ex:wu-layer-radical}
Prove that an element has degree \(12\) over \(\F_q\) if and only if it has degree six over \(\F_{q^2}\).
When comparing the two sides, place the element in a common finite extension and regard the base fields as its specified subfields.
Explain why the converse fails if two is replaced by three and six by four.
\end{exercise}
\begin{exercise}[Root degrees of compressed polynomials]\label{ex:wu-compressed-not-period}
Prove that \(G_{q,2}(T)=1+T\) and \(G_{q,3}(T)=1+T+T^{q+1}\).
Use the quadratic case to show that roots of a compressed polynomial need not retain the original exact period.
\end{exercise}
\begin{exercise}[Invariants versus polynomial functions]\label{ex:wu-polynomial-invariance}
Prove that \(\Theta_q(X)\) is invariant under the affine group but represents the zero function on \(\F_q\).
Show that this does not imply that it is the zero polynomial.
Further explain why checking only \(X\in\F_q\) is insufficient when testing \(F(aX+b)=F(X)\).
\end{exercise}
\begin{exercise}[Separation of compressed layers and power-compression return]\label{ex:wu-compressed-separable}\leavevmode\par\nobreak
\begin{enumerate}
\item Prove that the \(G_{q,n}\) for distinct indices \(n>1\) are pairwise coprime.
\item Take \(k=\F_5\), let \(\alpha\) be a primitive 32nd root of unity, and put \(\beta=\alpha^8\), \(\xi=\alpha^4\).
Prove that \(\beta\in k\) and \([k(\alpha):k]=8\), and determine the root fiber of the power map \(x\mapsto x^8\) over \(\beta\) and its first-return map.
Writing elements of \(\mu_8\) as \(\varepsilon=\xi^j\), prove that the return map in exponent coordinates is \(j\mapsto1+5j\pmod8\), and write the full cycle starting from exponent \(0\).
Compare the first two iterates at \(1\) of the actual return map and multiplication by \(\xi\) alone. Explain why the fifth-power action cannot be omitted,
and why equal return periods do not imply equal return permutations.
\end{enumerate}
\end{exercise}
\begin{exercise}[Affine specialization of the binary recurrence]\leavevmode\phantomsection\label{ex:wu-binary-specialization}
Starting from \(Q_{q,0}=0,Q_{q,1}=1\), compute \(Q_{q,2},Q_{q,3},Q_{q,4}\).
Verify that substituting \(U=T,V=1+T\) turns \(Q_{q,3}\) into \(1+T+T^{q+1}\), and use this to explain \(\mathcal E_{q,3}(X,1)=E_{q,3}(X)\).
\end{exercise}

\begin{exercise}[Point counts and base-field lifting]\label{ex:wu2-zeta-basechange}
For \(S=\F_{q^n}\) with \(\sigma:x\mapsto x^q\), prove that \(N_r=q^{\gcd(n,r)}\), and write the finite Euler product for \(Z_S(u)\).
Then explain how one original orbit of length \(d\) contributes new product factors when \(\sigma\) is replaced by \(\sigma^s\).
\end{exercise}

\begin{exercise}[Linear factors of compressed layers and realization of period drops]\leavevmode\phantomsection\label{ex:wu7-compressed-linear-factors}
Let \(n>1\), and let \(J_{q,n}(T)\) denote the product of all monic linear factors of \(G_{q,n}(T)\).
Prove that
\[
J_{q,n}(T)=
\begin{cases}
T-1,&n=p,\\
\overline\Phi_n(T),&n\mid q-1,\\
1,&\text{otherwise},
\end{cases}
\]
where the second case still requires \(n>1\), and \(\overline\Phi_n\) is the reduction in \(k\) of the integer cyclotomic polynomial.
Deduce that every nontrivial period-reduction factor permitted by Proposition~\ref{prop:affine-return-order} can be realized.
\end{exercise}

\chapter{From Unramified \texorpdfstring{$p$}{p}-adic Extensions to Finite Fields}
\label{revadic:chapter}

The equation $X^{p^n}-X=0$ gives a finite field in characteristic $p$, and zero together with a group of roots of unity in characteristic zero. If these roots are adjoined to $\mathbb Q_p$, can reduction of the integer ring recover the finite field and its Frobenius? How much information about multiplication, linear relations among conjugates, and coefficient cancellation survives this passage?

An unramified cyclotomic extension supplies the residue field, and Hensel uniqueness matches roots and factors on both sides. Normal integral bases lift the cyclic modules of Chapter~9 to integer rings. For coefficients indistinguishable after one reduction, integral subset counts and cyclotomic remainders first detect exact zeros; local factors and Witt carries then determine divisibility levels. In the character sums of Chapter~18, this precision information also connects point-count recurrences with Newton slopes.

The first carry also constrains linear maps, digit subspaces, and code permutations. For low-degree subfield representations, the Smith exponents of a Frobenius difference record the precision to which they persist.

Index constraints on phase differences further restrict higher cancellation. For Fermat curves, Jacobi sums supply the orbit return scalars, so carries among exponents determine slopes directly.

\section{Compatible congruences, completeness, and Hensel lifting}
\label{v135:sec-foundations}

The residues of an integer modulo $p,p^2,\ldots$ are compatible. Retaining every precision gives the ring
\[
 \mathbb Z_p=\varprojlim_s\mathbb Z/p^s\mathbb Z
 =\{(x_s)_{s\ge1}:x_{s+1}\equiv x_s\pmod{p^s}\},
\]
with componentwise addition and multiplication. Every element has a unique expansion
$x=\sum_{j\ge0}a_jp^j$, where $0\le a_j<p$; its finite truncations specify the corresponding congruence classes.
For $x\ne0$, its \emph{$p$-adic valuation} is
$v_p(x)=\min\{j:a_j\ne0\}$; put $v_p(0)=+\infty$.
Thus $x=p^{v_p(x)}u$ with $u\bmod p\ne0$. The inverses of $u$ at finite precision are compatible and define $u^{-1}\in\mathbb Z_p$.
Consequently $v_p(xy)=v_p(x)+v_p(y)$, and $\mathbb Z_p$ is an integral domain.
Its fraction field $\mathbb Q_p$ is the \emph{field of $p$-adic numbers}. The extended valuation
$v_p(x/y)=v_p(x)-v_p(y)$ satisfies
\[
 v_p(x+y)\ge\min\{v_p(x),v_p(y)\},
\]
with equality when the two valuations differ. The absolute value $|x|_p=p^{-v_p(x)}$ defines a metric. A Cauchy sequence eventually stabilizes at each fixed precision and therefore has a unique limit; both $\mathbb Z_p$ and $\mathbb Q_p$ are complete.

Let $A$ be an integral domain and $\pi\in A$ a nonzero nonunit. If every nonzero element of $A$ has a unique expression $\pi^ju$, with $j\ge0$ and $u\in A^\times$, then $A$ is a \emph{discrete valuation ring}; $\pi$ is a \emph{uniformizer}. Its nonzero ideals are $(\pi^j)$, its unique maximal ideal is $(\pi)$, and $A/(\pi)$ is its \emph{residue field}. It is \emph{complete} when $A\simeq\varprojlim_s A/(\pi^s)$. For instance, $\mathbb Z_p$ has uniformizer $p$ and residue field $\mathbb F_p$.

\begin{lemma}[Hensel lifting over a complete discrete valuation ring]
\label{v135:lem-hensel}
Let $A$ be a complete discrete valuation ring with uniformizer $\pi$ and residue field $k=A/(\pi)$.
\begin{enumerate}
\item If $F\in A[X]$ and $\bar a\in k$ satisfies $\bar F(\bar a)=0$ and $\bar F'(\bar a)\ne0$, there is a unique $a\in A$ with $F(a)=0$ and $a\bmod\pi=\bar a$.
\item If a monic polynomial $F$ satisfies $\bar F=\bar G\bar H$ for relatively prime monic polynomials $\bar G,\bar H\in k[X]$, there is a unique monic factorization $F=GH$ with the prescribed reductions.
\end{enumerate}
\end{lemma}
\begin{proof}
If $F(a_s)\equiv0\pmod{\pi^s}$, put $a_{s+1}=a_s+\pi^st$. Taylor expansion gives
\[
 F(a_{s+1})\equiv F(a_s)+\pi^stF'(a_s)\pmod{\pi^{s+1}}.
\]
The derivative is a unit, so $t\bmod\pi$ is unique. Successive corrections converge to a root; uniqueness at every precision proves uniqueness of the limit.

For factors with $F-G_sH_s\in\pi^sA[X]$, put
$G_{s+1}=G_s+\pi^sU$ and $H_{s+1}=H_s+\pi^sV$, imposing
$\deg U<\deg\bar G$ and $\deg V<\deg\bar H$. The correction equation is
\[
 \bar U\bar H+\bar G\bar V
 =\overline{(F-G_sH_s)/\pi^s}.
\]
The linear map on the left is injective: a zero relation and coprimality force $\bar G\mid\bar U$ and $\bar H\mid\bar V$. Both spaces have dimension $\deg F$, so the map is bijective. Unique corrections at each precision converge coefficientwise to the required unique factorization.
\end{proof}

The $p$-adic valuation extends uniquely to every finite extension $L/\mathbb Q_p$. Its \emph{ring of integers}
$\mathcal O_L=\{x:v_p(x)\ge0\}$ is a complete discrete valuation ring and is the integral closure of $\mathbb Z_p$ in $L$. If $\pi$ is a uniformizer and $p=u\pi^e$, then $e$ is the \emph{ramification index}. The extension is \emph{unramified} when $e=1$, in which case $p$ itself is a uniformizer. This valuation-extension theorem is given in\cite[Theorem 10.4]{Sutherland2021Local}. The cyclotomic construction below also supplies a free basis, completeness, and a uniformizer for its integer ring directly.

\section{Roots of unity and residue fields}
\label{revadic:sec-roots}

\subsection{From an orbit modulo \texorpdfstring{$p$}{p} to a \texorpdfstring{$p$}{p}-adic root}

When $p\nmid m$, the polynomial $T^m-1$ is squarefree modulo $p$, so its coprime factors and their roots have unique Hensel lifts.

\begin{lemma}[Reduction of prime-to-$p$ torsion]
\label{revadic:lem-torsion}
Let $m\ge1$ be prime to $p$ and let $\zeta\in\overline{\mathbb Q}_p$ be a primitive $m$th root of unity. Then $\zeta$ is integral and its reduction still has multiplicative order $m$.
\end{lemma}
\begin{proof}
Integrality follows from the monic equation $T^m-1=0$. If $\bar\zeta^d=1$ for some $1\le d<m$, then $\zeta^d$ and $1$ are two $m$th roots of unity reducing to the same residue class. Since $m$ is a $p$-adic unit, the derivative $mT^{m-1}$ is nonzero modulo the maximal ideal at a unit root. Uniqueness in Hensel's lemma forces $\zeta^d=1$, a contradiction. Reduction here takes place in the integer ring of a finite extension containing $\zeta$.
\end{proof}

\begin{lemma}[Lifting cyclotomic factors]
\label{revadic:lem-factor}
Suppose $m\ge2$, $p\nmid m$, and $r=\ord_m(p)$. Every monic irreducible factor of $\Phi_m(T)$ over $\mathbb F_p$ has degree $r$. Its corresponding monic Hensel factor $F(T)\in\mathbb Z_p[T]$ is irreducible over $\mathbb Q_p$, and
\[
 A=\mathbb Z_p[T]/(F),\qquad
 A/pA\simeq\mathbb F_p[T]/(\bar F)
\]
is a complete discrete valuation ring with uniformizer $p$.
\end{lemma}
\begin{proof}
Choose a primitive $m$th root $\alpha$ in $\overline{\mathbb F}_p$. Its Frobenius orbit has length equal to the least positive $r$ for which $p^r\equiv1\pmod m$. Thus every irreducible factor of $\bar\Phi_m$ has degree $r$. Since $T^m-1$ is squarefree modulo $p$, these factors are pairwise coprime. Hensel's lemma lifts each of them uniquely to a monic factor $F$ of $\Phi_m$.

Irreducibility of $\bar F$ implies irreducibility of $F$ over $\mathbb Q_p$. Indeed, a factorization over $\mathbb Q_p$ could, by Gauss's lemma for the discrete valuation ring $\mathbb Z_p$, be chosen with positive-degree monic factors in $\mathbb Z_p[T]$; reduction would contradict irreducibility of $\bar F$. The ring $A$ is a finite free $\mathbb Z_p$-module and hence is complete. If $u\notin pA$, choose $v\in A$ with $uv=1-pz$. The convergent series $\sum_{j\ge0}(pz)^j$ inverts $1-pz$, so $u$ is a unit. Hence $pA$ is the unique maximal ideal. Every nonzero $x\in A$ can, by $p$-adic separatedness, be written uniquely as $p^ju$ with $u\notin pA$, hence with $u$ a unit. Thus $A$ is a complete discrete valuation ring. In particular it is integrally closed and equals the integer ring of its fraction field.
\end{proof}

\subsection{The unramified cyclotomic tower}
\begin{theorem}[Constructing a field of $p^n$ elements from $\mathbb Q_p$]
\label{revadic:thm-main}
For each $n\ge1$, set $m_n=p^n-1$. If $m_n>1$, choose a primitive $m_n$th root $\zeta_n\in\overline{\mathbb Q}_p$ and put
\[
 K_n=\mathbb Q_p(\zeta_n).
\]
If $m_n=1$ (equivalently, $p=2,n=1$), set $\zeta_n=1$ and $K_n=\mathbb Q_2$. Then
\begin{equation}\label{revadic:eq-main}
 [K_n:\mathbb Q_p]=n,\qquad
 \mathcal O_n=\mathbb Z_p[\zeta_n],\qquad
 \mathcal O_n/p\mathcal O_n\simeq\mathbb F_{p^n}.
\end{equation}
The reduction $\bar\zeta_n$ generates the multiplicative group of the residue field. This construction therefore gives both the cardinality of the field and the cyclicity of its multiplicative group.
\end{theorem}
\begin{proof}
For $n=1$, $m_n=p-1$, and $\mathbb Q_p$ already contains all $(p-1)$st roots of unity, by Hensel lifting from $\mathbb F_p^\times$. The case $p=2$ uses the stated convention. If $n>1$, then $m_n>1$ and
\[
 p^n\equiv1\pmod{m_n},\qquad
 0<p^d-1<m_n\quad(0<d<n).
\]
Consequently $\ord_{m_n}(p)=n$. Lemma~\ref{revadic:lem-factor} gives $[K_n:\mathbb Q_p]=n$ and a residue field with $p^n$ elements. Lemma~\ref{revadic:lem-torsion} shows that $\bar\zeta_n$ has order $m_n=p^n-1$, exactly the order of the nonzero multiplicative group.
\end{proof}

\begin{theorem}[Uniqueness of the unramified extension]
\label{revadic:thm-unique}
Inside a fixed $\overline{\mathbb Q}_p$, the field $K_n$ is the unique unramified extension of degree $n$. More precisely, if $L/\mathbb Q_p$ is finite, $p$ remains a uniformizer of its integer ring $\mathcal O_L$, and its residue field $k_L$ has $p^n$ elements, then $L=K_n$.
\end{theorem}
\begin{proof}
Cyclicity of the finite-field multiplicative group provides $a\in k_L^\times$ of order $p^n-1$. The derivative of $T^{p^n-1}-1$ is a unit at $a$ modulo the maximal ideal, so Hensel's lemma in the complete integer ring gives a unique lift $u\in\mathcal O_L$ with $u^{p^n-1}=1$ and $\bar u=a$. Reduction preserves prime-to-$p$ orders, so $u$ is a primitive $(p^n-1)$st root of unity and $K_n=\mathbb Q_p(u)\subseteq L$.

Choose a basis of $k_L/\mathbb F_p$ and arbitrary lifts $e_1,\ldots,e_n\in\mathcal O_L$. For $x\in\mathcal O_L$, expand its residue in this basis, subtract the lifted expansion, and divide by $p$. Repeating and using completeness gives $x=\sum_i a_i e_i$ with $a_i\in\mathbb Z_p$. A nonzero relation, divided by the smallest common power of $p$ and reduced modulo $p$, would contradict independence of the residue basis. Thus $\mathcal O_L$ has $\mathbb Z_p$-rank $n$, and $[L:\mathbb Q_p]=n$. Comparing this with $[K_n:\mathbb Q_p]=n$ gives $L=K_n$.\end{proof}

\begin{corollary}[Finite precision and Witt rings]
\label{revadic:cor-witt}
Write $k_n=\mathcal O_n/p\mathcal O_n$. For all $s\ge1$,
\[
 \mathcal O_n/p^s\mathcal O_n\simeq\operatorname{GR}(p^s,n),
 \qquad \mathcal O_n\simeq\varprojlim_s\operatorname{GR}(p^s,n)
 \simeq W(k_n).
\]
After an identification of residue fields is fixed, these are compatible isomorphisms of complete local rings. The last isomorphism is not induced by an embedding of rings $k_n\hookrightarrow\mathcal O_n$.
\end{corollary}
\begin{proof}
By Lemma~\ref{revadic:lem-factor}, $\mathcal O_n=\mathbb Z_p[T]/(F)$ with $\bar F$ irreducible of degree $n$. Hence
\[
 \mathcal O_n/p^s\mathcal O_n=(\mathbb Z/p^s\mathbb Z)[T]/(F\bmod p^s)
\]
is precisely the Galois-ring model of Chapter~5. After fixing the residue-field identification, uniqueness of Hensel roots gives compatible isomorphisms at all precisions. Completeness gives the inverse limit, and Theorem~\ref{v132:thm-unramified-limit} and its Witt representation identify it with $W(k_n)$.
\end{proof}

\begin{example}[A cubic unramified extension from a seventh root]
\label{revadic:ex-seven}
Take $p=2,n=3$. Since $\operatorname{ord}_7(2)=3$, the extension $\mathbb Q_2(\zeta_7)$ is unramified of degree three, with integer ring $\mathbb Z_2[\zeta_7]$ and residue field of size eight. The two orbits in Example~\ref{v132:ex-seven-period} give the two canonical cubic factors, whose full coefficients appear in Example~\ref{v134:ex-seven-full-profile}. The same orbit now also determines the local extension degree and residue-field degree.
\end{example}

\begin{example}[A cyclotomic factorization over $\mathbb Q_3$]
\label{revadic:ex-eight}
Take $p=3,n=2$, so that $m=8$ and $\ord_8(3)=2$. A primitive eighth root $\zeta$ gives an unramified quadratic extension, although $\Phi_8(T)=T^4+1$ is reducible over $\mathbb Q_3$. Put $u=\zeta+\zeta^3$. Then $u^2=-2$, and since $-2\equiv1\pmod3$, Hensel lifting puts both roots of $X^2+2$ in $\mathbb Z_3$, hence $u\in\mathbb Z_3$. Exchanging the two exponent orbits if necessary, take $u\equiv1\pmod3$. Consequently
\[
 T^4+1=(T^2-uT-1)(T^2+uT-1).
\]
Both quadratic factors are irreducible over $\mathbb Q_3$, and their reductions are $T^2\mp T+2$. Thus the degree of an unramified extension is controlled by a \emph{single Frobenius orbit}, not necessarily by the full degree $\varphi(m)$ of a cyclotomic polynomial.
\end{example}

\begin{theorem}[Arithmetic Frobenius on the integer ring]
\label{revadic:thm-frob}
The extension $K_n/\mathbb Q_p$ is cyclic Galois. Its Galois group is generated by a unique automorphism $\phi$ satisfying
\[
 \phi(\zeta_n)=\zeta_n^p,\qquad
 \overline{\phi(x)}=\bar x^{\,p}\quad(x\in\mathcal O_n).
\]
Reduction induces an isomorphism
\begin{equation}\label{revadic:eq-gal}
 \operatorname{Gal}(K_n/\mathbb Q_p)
 \xrightarrow{\ \sim\ }
 \operatorname{Gal}(k_n/\mathbb F_p)
 \simeq\mathbb Z/n\mathbb Z.
\end{equation}
In general, $\phi(x)\ne x^p$ for $x\in\mathcal O_n$.
\end{theorem}
\begin{proof}
Let $F$ be the minimal polynomial of $\zeta_n$ over $\mathbb Q_p$. The roots of its reduction are $\bar\zeta_n^{p^j}$, $0\le j<n$. Each lifts uniquely to a root of $F$ in $\mathcal O_n$. Since $F\mid X^{p^n-1}-1$, these lifts are roots of unity; uniqueness in each residue class identifies them with $\zeta_n^{p^j}$. Thus $K_n$ is a splitting field and $\zeta_n\mapsto\zeta_n^p$ defines an automorphism $\phi$ of order $n$. On $\mathbb Z_p[\zeta_n]$ it reduces to the $p$th-power map. Both Galois groups are cyclic of order $n$, and $\phi$ reduces to a generator, so reduction is an isomorphism. For $n=1$ it is the identity. Finally, $\phi(p)=p\ne p^p$.
\end{proof}

\section{Two realizations of the same root equation}
\label{revadic:sec-two-roots}

Consider the equation in two characteristics:
\[
                       X^{p^n}-X=0.\tag{\(*\)}
\]
In characteristic $p$ its roots form a finite field. In characteristic zero they consist of zero and a group of roots of unity. The second set is not closed under ordinary addition, but its roots correspond individually to those in characteristic $p$. This is where the constructions meet.

\begin{theorem}[Root realizations, operations, and extension compatibility]
\label{revadic:thm-two-roots}
Let $n\ge1$, put $K_n=\mathbb Q_p(\zeta_{p^n-1})$ (with $K_n=\mathbb Q_p$ if $p^n-1=1$), and write $\mathcal O_n$ for its integer ring and $k_n=\mathcal O_n/p\mathcal O_n$. Set
\[
 R_n=\{x\in\mathcal O_n:x^{p^n}=x\},\qquad
 S_n=\{a\in\overline{\mathbb F}_p:a^{p^n}=a\}.
\]
Then
\[
 R_n=\{0\}\cup\mu_{p^n-1},\qquad S_n\simeq k_n\simeq\mathbb F_{p^n},
\]
after fixing an identification of $k_n$ with $S_n$, reduction gives a bijection $R_n\xrightarrow{\sim}S_n$. Its inverse $a\mapsto[a]$ preserves zero, unity, multiplication, inverses on nonzero elements, and Frobenius:
\[
 [ab]=[a][b],\qquad [a^p]=\phi([a])=[a]^p.
\]
The set $R_n$ is not closed under ordinary addition, however; this bijection is not a ring isomorphism.

Inside a common algebraic closure, choose compatible residue-field identifications. For $d\mid n$, reduction commutes with $R_d\subseteq R_n$ and $S_d\subseteq S_n$, and these inclusions also commute with Frobenius.
\end{theorem}
\begin{proof}
In characteristic zero, $X^{p^n}-X=X(X^{p^n-1}-1)$. Thus the roots are zero and the $(p^n-1)$st roots of unity. Since $p\nmid p^n-1$, distinct such roots remain distinct on reduction and retain their orders. They give $p^n$ distinct roots in $k_n$, hence all its elements. A chosen field isomorphism then identifies them individually with $S_n$. In characteristic $p$, the identities
\[
 (a+b)^{p^n}=a^{p^n}+b^{p^n},\qquad
 (ab)^{p^n}=a^{p^n}b^{p^n}
\]
show closure under addition and multiplication, and nonzero roots are closed under inverses. The derivative is $-1$, so there are $p^n$ distinct roots in an algebraic closure, forming a field. Reduction and uniqueness in Hensel lifting establish the asserted multiplicative and Frobenius compatibility. If $R_n$ were additively closed, the sum of $p$ copies of $1$, namely $p$, would be a root; but $p$ is neither zero nor a root of unity.

If $d\mid n$, then $p^d-1\mid p^n-1$, giving inclusions of the root groups and their generated integer rings. Reduction restricts to the same map on these subrings, proving compatibility.
\end{proof}

One must distinguish two uses of Frobenius. On $R_n$, arithmetic Frobenius agrees with $x\mapsto x^p$; on all of $\mathcal O_n$ it is the ring automorphism $\phi$, not pointwise exponentiation: $\phi(p)=p\ne p^p$.

\begin{proposition}[Unique roots at every finite precision]
\label{revadic:prop-precision-roots}
For $s\ge1$ set
\[
 R_{n,s}=\{x\in\mathcal O_n/p^s\mathcal O_n:
                     x^{p^n}=x\}.
\]
Each $a\in k_n$ has exactly one lift in $R_{n,s}$, namely $[a]\bmod p^s$. Thus $|R_{n,s}|=p^n$, reduction maps between successive precisions are compatible, and
\[
 R_n\simeq\varprojlim_s R_{n,s}
 \qquad\text{as sets equipped with multiplication and Frobenius}.
\]
\end{proposition}
\begin{proof}
For $F(X)=X^{p^n}-X$ we have $F'(X)=p^nX^{p^n-1}-1\equiv-1\pmod p$. Hensel's lemma lifts each residual root uniquely, including at each finite precision. Compatibility follows from uniqueness.
\end{proof}

\subsection{Recovering addition on the roots by \texorpdfstring{$p$}{p}-adic projection}
\label{revadic:subsec-transport}

The preceding correspondence preserves ordinary multiplication, but this does not mean that the characteristic-zero root set contains no information about finite-field addition. Addition can be recovered from the complete integer ring by a $p$-adic limit. For $u\in\mathcal O_n$, put
\[
 \tau(u)=\lim_{r\to\infty}u^{p^{nr}}.
\]
The existence of this limit is part of the following statement.

\begin{proposition}[Teichm\"uller projection and transported addition]
\label{revadic:prop-teich-projection}
The limit exists for every $u\in\mathcal O_n$, and
\[
 \tau(u)=[\bar u],\qquad \tau(u)\equiv u\pmod p,
 \qquad\tau(uv)=\tau(u)\tau(v).
\]
For $x,y\in R_n$ define
\begin{equation}\label{revadic:eq-transport-add}
 x\boxplus y=\tau(x+y)=\lim_{r\to\infty}(x+y)^{p^{nr}}.
\end{equation}
Then $(R_n,\boxplus,\cdot)$ is a field of $p^n$ elements, and reduction gives a field isomorphism
\[
 (R_n,\boxplus,\cdot)\xrightarrow{\ \sim\ }k_n.
\]
Moreover $\phi(x\boxplus y)=\phi(x)\boxplus\phi(y)$; if $x+y$ already belongs to $R_n$, then $x\boxplus y=x+y$.

More precisely, let $q=p^n$, $s\ge1$, and choose an integer $r\ge0$ with $nr\ge s-1$. For every $u\in\mathcal O_n$,
\begin{equation}\label{v135:eq-finite-projection}
 \tau(u)\equiv u^{q^r}\pmod{p^s\mathcal O_n}.
\end{equation}
Thus transported addition on $R_{n,s}$ has the finite formula
\[
 x\boxplus_s y=(x+y)^{q^r}\quad\text{in }\mathcal O_n/p^s\mathcal O_n.
\]
The resulting field $(R_{n,s},\boxplus_s,\cdot)$ is isomorphic to $k_n$, and reduction between precisions is a field isomorphism.
\end{proposition}
\begin{proof}
If $u\in p\mathcal O_n$, then $v_p(u^{p^{nr}})\to\infty$, and the limit is zero. If $u$ is a unit, write $u=[\bar u](1+pz)$, with $z\in\mathcal O_n$. Since $[\bar u]^{p^n}=[\bar u]$,
\[
 u^{p^{nr}}=[\bar u](1+pz)^{p^{nr}}.
\]
The binomial theorem gives $w^p\in1+p^{t+1}\mathcal O_n$ whenever $w\in1+p^t\mathcal O_n$ and $t\ge1$. Iteration implies $(1+pz)^{p^{nr}}\to1$, also for $p=2$. Thus $\tau(u)=[\bar u]$. Multiplicativity follows from uniqueness of the lift. The operation in \eqref{revadic:eq-transport-add} is exactly the addition transported from $k_n$ through the bijection $R_n\simeq k_n$; all field axioms and Frobenius compatibility follow.

For a unit $u$, the binomial estimate gives $u^{q^r}\equiv[\bar u]\pmod{p^{nr+1}}$. If $u\in p\mathcal O_n$, then $v_p(u^{q^r})\ge q^r\ge nr+1$. This proves~\eqref{v135:eq-finite-projection} whenever $nr\ge s-1$. Roots modulo $p^s$ correspond individually to $k_n$; the displayed operation transports its addition, and uniqueness makes reduction compatible with it.
\end{proof}

\begin{remark}[What the two realizations of the root equation mean]
\label{revadic:rem-master-root}
The roots of $X^{p^n}-X$ correspond individually in both characteristics. In characteristic $p$ their ordinary sum is again a root; in characteristic zero one first adds inside $\mathcal O_n$ and then projects back to the root set. The difference
\[
 x+y-(x\boxplus y)\in p\mathcal O_n
\]
is the starting point of the Witt carry: multiplication remains within the root set, while higher digits record the difference between ordinary and transported addition.
\end{remark}

\subsection{Residual roots and corrections of integral equations}
\label{revadic:subsec-root-correction}

Transported addition restores closure of the root set. When the defining polynomial itself changes, however, the Teichmüller representative need not solve the new equation. The discrepancy is already visible modulo $p^2$.

\begin{proposition}[Error valuations and Hensel correction of a simple root]
\label{revadic:prop-root-correction}
Let $F(X)\in\mathcal O_n[X]$, let $a\in k_n$ be a simple root of $\bar F$, and let $z\in\mathcal O_n$ be its unique root lift. For every $t\in\mathcal O_n$ with $\bar t=a$,
\begin{equation}\label{v143:eq-root-error}
 v_p(z-t)=v_p(F(t)).
\end{equation}
If $F(t)\ne0$ and $e=v_p(F(t))\ge1$, then
\begin{equation}\label{v143:eq-newton-precision}
 z\equiv t-\frac{F(t)}{F'(t)}\pmod{p^{2e}\mathcal O_n}.
\end{equation}
In particular, taking $t=[a]$ gives
\begin{equation}\label{revadic:eq-root-correction}
 z\equiv t-\frac{F(t)}{F'(t)}\pmod{p^2\mathcal O_n}.
\end{equation}
The canonical representative $[a]$ agrees with the root modulo $p^s$ if and only if $F([a])\equiv0\pmod{p^s}$; they agree exactly if and only if $F([a])=0$.

In the equal-characteristic ring $k_n[[u]]$, if $F_0(X,u)=f_0(X)+uf_1(X)+O(u^2)$ and $a$ is a simple root of $f_0$, the unique deformed root satisfies
\[
 z(u)=a-u\frac{f_1(a)}{f_0'(a)}+O(u^2).
\]
Both corrections are governed by the simple-root condition: the first measures the error in an integral equation, while the second measures a formal deformation of its coefficients.
\end{proposition}
\begin{proof}
Write $F(X)=(X-z)H(X)$. Division by the monic polynomial $X-z$ gives $H\in\mathcal O_n[X]$, and
$H(t)\bmod p=\bar F'(a)\ne0$. Thus the second factor in $F(t)=(t-z)H(t)$ is a unit, proving~\eqref{v143:eq-root-error}, including the infinite valuation when $t=z$.
If $F(t)\ne0$, then $v_p(z-t)=e$. Taylor expansion gives
\[
 0=F(z)=F(t)+F'(t)(z-t)+(z-t)^2U,
 \qquad U\in\mathcal O_n.
\]
Since $F'(t)$ is a unit, division proves~\eqref{v143:eq-newton-precision}. The equal-characteristic formula follows by comparing the coefficient of $u$ in $F_0(a+ub,u)$.
\end{proof}

\begin{example}[One residual root and two integral equations]
\label{revadic:ex-root-correction}
Take $p=2$, $k_2=\mathbb F_4$, and $a^2+a+1=0$, with $t=[a]=\zeta_3$. For $F_1=X^2+X+1$ one has $F_1(t)=0$. The polynomial $F_2=X^2+X-1$ has the same reduction modulo two, but $F_2(t)=-2$. Its Hensel root $z$ therefore satisfies
\[
 z\equiv\zeta_3+\frac{2}{2\zeta_3+1}\pmod{4\mathcal O_2},
 \qquad v_2(z-\zeta_3)=1.
\]
Thus the same residual root does not distinguish the two integral equations; the distinction first appears in their lifting correction.
\end{example}

\subsection{Four constructions at fixed parameters}
\label{v137:subsec-construction-comparison}
A common irreducible defining polynomial allows root sets, quotient rings, exponent orbits, and unramified lifts to be compared in the same coordinates.

\begin{example}[Comparing constructions of the fields with eight and nine elements]
\label{v137:ex-construction-comparison}
First take $f_8=X^3+X+1\in\mathbb F_2[X]$. It has no binary root and is therefore irreducible. A root $\alpha$ satisfies $\alpha^3=\alpha+1$ and has order seven, so
\[
 \{x\in\overline{\mathbb F}_2:x^8=x\}
 =\mathbb F_2(\alpha)\simeq\mathbb F_2[X]/(f_8).
\]
The nonzero binary cyclotomic cosets are $\{1,2,4\}$ and $\{3,5,6\}$. Choose a seventh root of unity $\zeta$ reducing to $\alpha$. Example~\ref{v134:ex-seven-full-profile} gives the canonical factor
\[
 F_8=X^3-\eta X^2+\eta'X-1
 \equiv X^3+2X^2+X+3\pmod4.
\]
This follows directly from~\eqref{v137:eq-seven-periods} and $\eta\equiv2\pmod4$, without repeating factor lifting.

Next take $f_9=X^2+X+2\in\mathbb F_3[X]$, which has no root in $\mathbb F_3$. If $f_9(\beta)=0$, then
\[
 \beta^2=2\beta+1,\quad\beta^3=2-\beta,\quad\beta^4=-1.
\]
Hence $\beta$ has order eight, and
\[
 \{x\in\overline{\mathbb F}_3:x^9=x\}
 =\mathbb F_3(\beta)\simeq\mathbb F_3[X]/(f_9).
\]
Choose an eighth root of unity $\xi$ reducing to $\beta$; its ternary cyclotomic coset is $\{1,3\}$. With the root $u$ of $u^2=-2$, $u\equiv1\pmod3$, from Example~\ref{revadic:ex-eight}, the canonical factor of $\xi$ is
\[
 F_9=X^2+uX-1\equiv X^2+4X+8\pmod9,
\]
since $u\equiv4\pmod9$. In the two cases, $\mathcal O_n=\mathbb Z_p[\zeta]$ or $\mathbb Z_p[\xi]$, respectively, and reduction gives the corresponding finite field.
\end{example}

\begin{center}\small
\begin{tabular}{p{0.15\textwidth}p{0.36\textwidth}p{0.38\textwidth}}
\toprule
Construction&Information retained&Data and operations required\\\midrule
Root set&The field in a common algebraic closure, Frobenius, and subfields&Closure of the roots of $X^{p^n}-X$; arithmetic coordinates are not yet chosen\\
Quotient ring&A power basis and its multiplication relations&A verified irreducible polynomial of degree $n$; naive product reduction costs $O(n^2)$ base-field operations\\
Cyclotomic cosets&Root orders, conjugacy orbits, and factor degrees&A modulus and a chosen root; an orbit of length $n$ takes $n$ exponent shifts, while coefficients still require root products or period relations\\
Unramified lift&Canonical roots of unity, coefficient divisibility, and compatibility across precisions&A residue factor and coprime lifting; its $n$ coefficients modulo $p^s$ record $ns$ base-$p$ digits\\\bottomrule
\end{tabular}
\end{center}
These operation counts assume that the defining polynomial or orbit data are already known. They exclude finding an irreducible polynomial, factoring the modulus, and computing multiplicative orders; $ns$ is a representation size, not an operation bound for lifting.
Viewing the coefficients of $f_8,f_9$ as integers also defines the corresponding unramified extensions, but their chosen roots are generally not canonical roots of unity. For example, their constant coefficients are $1,2$, respectively, whereas both canonical factors above have constant coefficient $-1$. Having the same residue field does not require the chosen generators to retain the same relations in characteristic zero.

\section{Frobenius, subfields, and cyclotomic cosets}
\label{revadic:sec-frob}

\begin{theorem}[The unramified tower and the subfield tower]
\label{revadic:thm-tower}
Inside a fixed $\overline{\mathbb Q}_p$, for positive integers $d,n$,
\[
 K_d\subseteq K_n\quad\Longleftrightarrow\quad d\mid n.
\]
If $d\mid n$, then $K_n/K_d$ is a cyclic unramified extension of degree $n/d$, with $\operatorname{Gal}(K_n/K_d)=\langle\phi^d\rangle$. On residue fields, $k_d\subseteq k_n$ and $k_d$ is exactly the fixed field of $\phi^d$. In particular, the unique intermediate field of degree $d$ over $\mathbb Q_p$ is $K_d$.
\end{theorem}
\begin{proof}
If $d\mid n$, then $p^d-1\mid p^n-1$, so $\mu_{p^d-1}\subseteq\mu_{p^n-1}$ and hence $K_d\subseteq K_n$. Conversely, inclusion implies $d\mid n$ by the tower law. By Theorem~\ref{revadic:thm-frob} and the subgroup structure of a cyclic Galois group, there is exactly one intermediate field of degree $d$, corresponding to $\langle\phi^d\rangle$. Its residue field consists of the solutions of $x^{p^d}=x$ and has $p^d$ elements.
\end{proof}

\begin{corollary}[Compatibility of the root equation across the unramified tower]
\label{revadic:cor-root-tower}
If $d\mid n$, then $R_d\subseteq R_n$ and
\[
 R_d=\{x\in R_n:\phi^d(x)=x\}
     =\{x\in R_n:x^{p^d}=x\}.
\]
Upon reduction this becomes $k_d=\{a\in k_n:a^{p^d}=a\}$. The same Teichm\"uller representatives thus describe subfield inclusions, Galois fixed points, and inclusions of root sets. Under transported addition, these inclusions are field embeddings as well.
\end{corollary}
\begin{proof}
Since $d\mid n$, we have $p^d-1\mid p^n-1$ and therefore $\mu_{p^d-1}\subseteq\mu_{p^n-1}$. On $R_n$ one has $\phi(x)=x^p$, whence $\phi^d(x)=x^{p^d}$. Its fixed points are precisely zero and the roots of $x^{p^d-1}=1$. Reduction yields the usual finite-field subfield criterion. Transported addition is obtained from the same reduction map, and so it is compatible with these inclusions.
\end{proof}

\begin{corollary}[Two descriptions of irreducible-factor orbits]
\label{revadic:cor-orbits}
Let $m\ge2$, $p\nmid m$, and $n=\ord_m(p)$. Over $\mathbb F_p$, irreducible factors of $\Phi_m$ are indexed by $p$-cyclotomic cosets of unit exponents modulo $m$; over $\mathbb F_{p^d}$ they are indexed by $p^d$-cosets. If $d\mid n$, the canonical integral lift of each factor is the corresponding factor of $\Phi_m$ over $\mathcal O_d[T]$, retaining the exponent set of its roots. For all roots of $T^m-1$, nonunit exponents are allowed; their lifted factors need not divide $\Phi_m$.
\end{corollary}
\begin{proof}
We have $\phi^d(\zeta_m^i)=\zeta_m^{p^di}$, and the fixed coefficient field is $K_d$. Uniqueness in Hensel's lemma associates a unique monic integral factor with every squarefree factor after reduction.
\end{proof}

\begin{remark}
For a fixed $m$, take $n=\ord_m(p)$. If $d$ does not divide $n$, one can still group the $p^d$-orbits in the compositum $K_{\operatorname{lcm}(d,n)}$. Enlarging the coefficient field and increasing the ring precision are different operations: the former changes the step size of Frobenius; the latter refines the coefficients attached to the same roots.
\end{remark}

\begin{theorem}[Embeddings and compatible root sets at every precision]
\label{v135:thm-hom-lifting}
For positive integers $d,n,s$, reduction gives natural bijections
\begin{align*}
 \operatorname{Hom}_{\mathbb Z_p}(\mathcal O_d,\mathcal O_n)
 &\simeq\operatorname{Hom}_{\mathbb F_p}(k_d,k_n),\\
 \operatorname{Hom}_{\mathbb Z/p^s\mathbb Z}
 (\mathcal O_d/p^s\mathcal O_d,\mathcal O_n/p^s\mathcal O_n)
 &\simeq\operatorname{Hom}_{\mathbb F_p}(k_d,k_n).
\end{align*}
All homomorphisms preserve unity. These sets are empty if $d\nmid n$ and contain $d$ elements if $d\mid n$. The unique lift of $\iota:k_d\to k_n$ satisfies
\begin{equation}\label{v135:eq-hom-digits}
 \widetilde\iota\!\left(\sum_{j\ge0}p^j[a_j]\right)
 =\sum_{j\ge0}p^j[\iota(a_j)].
\end{equation}
Modulo $p^s$, truncate the sums at $j=s-1$. These bijections commute with composition, Frobenius, and reduction of precision. On $R_d$ and $R_n$ they preserve multiplication and transported addition.
\end{theorem}
\begin{proof}
Write $\mathcal O_d=\mathbb Z_p[T]/(F)$ with $\bar F$ irreducible of degree $d$. A $\mathbb Z_p$-algebra homomorphism is determined by the image of $T$, hence by a root of $F$ in $\mathcal O_n$. A residue-field homomorphism is likewise determined by a root of $\bar F$ in $k_n$. Every residue root is simple, so Hensel uniqueness gives the bijection; successive corrections give the same result modulo $p^s$. The finite-field embedding criterion gives existence and the number of maps.

Uniqueness for the root equation gives $\widetilde\iota([a])=[\iota(a)]$. Repeatedly taking the residue, subtracting its Teichmüller representative, and dividing by $p$ gives a unique digit expansion. Algebra homomorphisms preserve $p^s$ ideals and are therefore continuous, so termwise application yields~\eqref{v135:eq-hom-digits}. Uniqueness also gives compatibility with composition, reduction, Frobenius, and transported addition.
\end{proof}

\section{Trace, norm, normal integral bases, and duality}
\label{revadic:sec-trace}

\subsection{Reduction of integral operations}
\begin{proposition}[Compatibility of trace and norm]
\label{revadic:prop-tr-nm}
Suppose $d\mid n$, $r=n/d$, and $x\in\mathcal O_n$. Then
\begin{align}
 \Tr_{K_n/K_d}(x)&=\sum_{j=0}^{r-1}\phi^{dj}(x)\in\mathcal O_d,\\
 \Nm_{K_n/K_d}(x)&=\prod_{j=0}^{r-1}\phi^{dj}(x)\in\mathcal O_d,
\end{align}
and
\begin{equation}\label{revadic:eq-red-tr-nm}
 \overline{\Tr_{K_n/K_d}(x)}=\Tr_{k_n/k_d}(\bar x),\qquad
 \overline{\Nm_{K_n/K_d}(x)}=\Nm_{k_n/k_d}(\bar x).
\end{equation}
\end{proposition}
\begin{proof}
The first formulas are those of the cyclic Galois group. The sum and product are fixed by $\phi^d$ and hence belong to $K_d$; they are integral, so they lie in $\mathcal O_d$. Reducing term by term, and using $\overline{\phi^{dj}(x)}=\bar x^{p^{dj}}$, proves the second formulas.
\end{proof}

\begin{proposition}[Exact compatibility of Teichmüller norms]
\label{revadic:prop-teichnorm}
For $a\in k_n$, let $[a]\in\mathcal O_n$ denote the unique Teichmüller representative satisfying $[a]^{p^n}=[a]$. If $d\mid n$, then
\[
 \phi([a])=[a^p],\qquad
 \Nm_{K_n/K_d}([a])=[\Nm_{k_n/k_d}(a)].
\]
For trace, in general only the congruence
\[
 \Tr_{K_n/K_d}([a])\equiv[\Tr_{k_n/k_d}(a)]\pmod p
\]
holds, not equality.
\end{proposition}
\begin{proof}
Hensel's lemma gives the unique lift of each simple residue root of $T^{p^n}-T$. A product of Teichmüller representatives is again such a representative and reduces to the finite-field norm, proving exact equality. A sum is compatible after reduction, but Teichmüller representatives are not generally closed under addition.
\end{proof}

\begin{example}[The first lost digit of a trace]
\label{revadic:ex-trace}
Take $K_2=\mathbb Q_2(\zeta_3)$, where $\zeta_3^2+\zeta_3+1=0$. Then
\[
 \Tr_{K_2/\mathbb Q_2}(\zeta_3)=\zeta_3+\zeta_3^2=-1,
 \qquad\Tr_{\mathbb F_4/\mathbb F_2}(\bar\zeta_3)=1.
\]
Since $[1]=1$, the integral trace differs from the Teichmüller lift of the finite trace by $-2$. By contrast, $\Nm(\zeta_3)=1$ agrees exactly on both sides. Thus trace and norm both commute with reduction, but only norm automatically commutes exactly with Teichmüller lifting.
\end{example}

For $d\mid n$, arbitrary lifts of a basis of $k_n/k_d$ form a basis of $\mathcal O_n$ over $\mathcal O_d$. Successive residue expansions prove spanning, and reduction after removal of the smallest common power of $p$ proves independence, as in Theorem~\ref{revadic:thm-unique}. Thus this module is free of rank $n/d$.

\subsection{Normal integral bases and trace duality}
\begin{theorem}[An integral lifting criterion for normal bases]
\label{revadic:thm-normal}
Let $\theta\in\mathcal O_n$. The following conditions are equivalent:
\begin{enumerate}
\item $\theta,\phi(\theta),\ldots,\phi^{n-1}(\theta)$ form a $\mathbb Z_p$-basis of $\mathcal O_n$;
\item $\bar\theta,\bar\theta^p,\ldots,\bar\theta^{p^{n-1}}$ form a normal basis of $k_n/\mathbb F_p$.
\end{enumerate}
Consequently, \emph{every integral lift} of a finite-field normal element generates a normal integral basis of the unramified extension. Conversely, reduction of a normal integral basis is a normal basis.
\end{theorem}
\begin{proof}
Choose any $\mathbb Z_p$-basis $e_1,\ldots,e_n$ of $\mathcal O_n$, whose reductions form an $\mathbb F_p$-basis of $k_n$. Express the $n$ Galois conjugates of $\theta$ as column vectors with respect to this basis, giving $M\in M_n(\mathbb Z_p)$. Its reduction $\bar M$ is precisely the coordinate matrix of the Frobenius conjugates in $k_n$. Hence
\[
 M\in\operatorname{GL}_n(\mathbb Z_p)
 \iff \det M\not\equiv0\pmod p
 \iff\bar M\in\operatorname{GL}_n(\mathbb F_p).
\]
These are exactly the two basis conditions. The normal basis theorem of Chapter~\ref{chap:normal} supplies $\bar\theta$ and hence a normal integral basis. Replacing $\mathbb Z_p,\phi,\mathbb F_p$ by $\mathcal O_d,\phi^d,k_d$ gives the same criterion for every relative extension with $d\mid n$.
\end{proof}

\begin{proposition}[Trace duality commutes with reduction]
\label{revadic:prop-dual}
Let $e_1,\ldots,e_n$ be an integral basis of $\mathcal O_n/\mathbb Z_p$. The trace Gram matrix
\[
 G=\bigl(\Tr_{K_n/\mathbb Q_p}(e_ie_j)\bigr)_{i,j}
\]
lies in $\operatorname{GL}_n(\mathbb Z_p)$. There is therefore a unique dual integral basis $e_1^*,\ldots,e_n^*\in\mathcal O_n$ with $\Tr(e_ie_j^*)=\delta_{ij}$; the reductions of these bases are trace-dual bases of $k_n/\mathbb F_p$. In particular, the integral trace pairing of an unramified extension is perfect.
\end{proposition}
\begin{proof}
Modulo $p$, $G$ is the trace-pairing Gram matrix in the finite field with respect to $\bar e_i$. Finite-field extensions are separable, so this trace pairing is nondegenerate and $\det\bar G\ne0$. Thus $\det G$ is a $\mathbb Z_p$-unit, and all entries of $G^{-1}$ lie in $\mathbb Z_p$. The columns of $G^{-1}$ define the dual integral basis, whose reduction gives the asserted finite-field duality.
\end{proof}

Proposition~\ref{v141:prop-period-duality} applies the same trace-duality principle to a period phase algebra, recovering exceptional-phase idempotents from trace moments even when their ordinary trace is zero.

\begin{example}[A normal integral basis and its trace dual in degree two]
\label{revadic:ex-basis}
Continue with $K_2=\mathbb Q_2(\zeta_3)$. Since $\zeta_3^2=-1-\zeta_3$, the coordinate matrix of $\zeta_3,\zeta_3^2$ relative to $(1,\zeta_3)$ is
\[
 \begin{pmatrix}0&-1\\1&-1\end{pmatrix},
 \qquad\det=1.
\]
Thus they form a normal integral basis over $\mathbb Z_2$, reducing to a normal basis of $\mathbb F_4/\mathbb F_2$. The trace Gram matrix of this normal basis is
\[
 \begin{pmatrix}-1&2\\2&-1\end{pmatrix},\qquad
 \det=-3\in\mathbb Z_2^\times.
\]
Its inverse has $2$-adic integral entries and provides the dual integral basis. The two unit determinants show separately why normality and trace duality survive reduction.
\end{example}

\begin{proposition}[The cyclic-module normality criterion on both sides]
\label{revadic:prop-cyclic-normal}
Let $\gamma\in k_n$ be normal and let $\Theta\in\mathcal O_n$ be any integral lift. For
\[
 f(T)=\sum_{j=0}^{n-1}a_jT^j\in\mathbb Z_p[T],\qquad
 \theta=f(\phi)\Theta,
\]
one has
\[
 \theta\text{ generates a normal integral basis of }\mathcal O_n/\mathbb Z_p
 \iff \gcd(\bar f(T),T^n-1)=1\quad\text{in }\mathbb F_p[T].
\]
This remains true when $p\mid n$, so that $T^n-1$ has repeated factors. Moreover, there is an isomorphism of $\mathbb Z_p[T]$-modules
\[
 \mathcal O_n\simeq\mathbb Z_p[T]/(T^n-1),\qquad
 T\cdot x=\phi(x),
\]
whose reduction is the Frobenius cyclic-module description of the finite field.
\end{proposition}
\begin{proof}
By Theorem~\ref{revadic:thm-normal}, $\Theta,\phi(\Theta),\ldots,\phi^{n-1}(\Theta)$ form an integral basis. Thus the cyclic map $\mathbb Z_p[T]/(T^n-1)\to\mathcal O_n$ generated by $\Theta$ is an isomorphism. The conjugates of $f(\phi)\Theta$ form an integral basis precisely when multiplication by $f$ is invertible in this quotient. By reduction and the unit-determinant criterion, this is equivalent to $\bar f$ being a unit in $\mathbb F_p[T]/(T^n-1)$, or $\gcd(\bar f,T^n-1)=1$. No squarefreeness of $T^n-1$ is needed.
\end{proof}

\begin{proposition}[Taking invariants commutes with increasing precision]
\label{revadic:prop-fixed-precision}
If $d\mid n$, inclusion induces, for every $s\ge1$, an isomorphism
\begin{equation}\label{revadic:eq-fixed-precision}
 \mathcal O_d/p^s\mathcal O_d
 \ \xrightarrow{\ \sim\ }\
 \bigl(\mathcal O_n/p^s\mathcal O_n\bigr)^{\langle\phi^d\rangle}.
\end{equation}
Thus taking fixed rings along the unramified tower and truncating precision commute. The assertion concerns entire quotient rings, not only the fixed Teichmüller roots.
\end{proposition}
\begin{proof}
Every class modulo $p^s$ has a unique finite expansion
\[
 x=\sum_{i=0}^{s-1}p^i[a_i],\qquad a_i\in k_n,
 \qquad\phi^d(x)=\sum_{i=0}^{s-1}p^i[a_i^{p^d}].
\]
By uniqueness of digits, $\phi^d(x)=x$ if and only if every $a_i$ belongs to $k_d$. These are exactly the digit expansions of $\mathcal O_d/p^s\mathcal O_d$. Also $p^s\mathcal O_n\cap\mathcal O_d=p^s\mathcal O_d$, so the map induced by inclusion is injective.
\end{proof}

This commutation extends the finite-field fixed-field criterion digit by digit. A relative normal integral basis expresses the same condition as equality of cyclic coordinates and also describes the image of Frobenius difference.

\begin{corollary}[Additive Hilbert 90 over integers and at every finite precision]
\label{v136:cor-integral-hilbert90}
Let $d\mid n$, put $\sigma=\phi^d$, $r=n/d$, and write
$A_{e,s}=\mathcal O_e/p^s\mathcal O_e$.
For every $s\ge1$, there is an exact sequence of $A_{d,s}$-modules
\begin{equation}\label{v136:eq-integral-hilbert90}
 0\longrightarrow A_{d,s}\longrightarrow A_{n,s}
 \xrightarrow{\ \sigma-1\ }A_{n,s}
 \xrightarrow{\ \Tr\ }A_{d,s}\longrightarrow0,
 \qquad \Tr=\sum_{i=0}^{r-1}\sigma^i.
\end{equation}
Replacing $A_{d,s},A_{n,s}$ by $\mathcal O_d,\mathcal O_n$ also gives an exact sequence.
Thus $\sigma(x)-x=a$ is soluble in $A_{n,s}$ if and only if $\Tr(a)=0$; each nonempty solution fiber has exactly $p^{ds}$ elements.
\end{corollary}
\begin{proof}
Choose a relative normal integral basis $\theta,\sigma\theta,\ldots,\sigma^{r-1}\theta$.
Its sum $t=\Tr(\theta)\in\mathcal O_d$ is a unit, since otherwise reduction gives a nontrivial linear relation among normal-basis elements.
For $x=\sum_i a_i\sigma^i\theta$, one has $\Tr(x)=t\sum_i a_i$; hence trace is surjective and its kernel consists of coordinate vectors with sum zero.
These are precisely the image of cyclic difference: the equations $b_i=a_{i-1}-a_i$, with indices modulo $r$, can be solved successively, with sole compatibility condition $\sum_i b_i=0$.
The difference kernel consists of constant coordinate vectors, namely $\mathcal O_d t=\mathcal O_d$.
All arguments remain valid modulo $p^s$, and the fiber size is $|A_{d,s}|=p^{ds}$.
\end{proof}

This connects additive Hilbert 90 from Chapter~\ref{chap:trace} with the cyclic modules of Chapter~\ref{chap:normal}.
The proof uses cyclic coordinates and a unit trace without dividing by $r$, so it also applies when $p\mid r$.

\subsection{Seventh roots: primitivity, normality, and self-duality}
\label{revadic:subsec-seven-unified}
Continue with $p=2,n=3$, and choose $\alpha=\bar\zeta_7$ satisfying $\alpha^3+\alpha+1=0$. Since $\alpha$ has order seven, it generates $\mathbb F_8^\times$, but
\[
 \alpha+\alpha^2+\alpha^4=0.
\]
Thus $\alpha$ is not normal. A primitive root of unity gives a multiplicative generator, not necessarily a normal-basis generator.

Take $\theta=\zeta_7^3$. Its reductions of conjugates are $\alpha^3,\alpha^6,\alpha^5$, with
\[
 \alpha^3=1+\alpha,\quad
 \alpha^6=1+\alpha^2,\quad
 \alpha^5=1+\alpha+\alpha^2.
\]
Their coordinate matrix relative to $(1,\alpha,\alpha^2)$ is
\[
 M=\begin{pmatrix}1&1&1\\1&0&1\\0&1&1\end{pmatrix},
 \qquad\det M=1\in\mathbb F_2.
\]
Hence these elements are a normal basis of $\mathbb F_8/\mathbb F_2$, and $\theta,\phi(\theta),\phi^2(\theta)$ form a normal integral basis of $\mathcal O_3/\mathbb Z_2$.

Use $\eta,\eta'$ from Example~\ref{v132:ex-seven-period}.
The trace Gram matrix of the integral conjugates is
\[
 G=\begin{pmatrix}
 \eta'&\eta&\eta\\
 \eta&\eta'&\eta\\
 \eta&\eta&\eta'
 \end{pmatrix}.
\]
Since $\bar\eta=0$ and $\bar\eta'=1$, we obtain $G\equiv I_3\pmod2$. Thus the reduced normal basis is also a \emph{self-dual normal basis}. Over the integer ring, the trace-dual integral basis is given by the columns of $G^{-1}$. Since $G\ne I_3$, this integral basis itself is not self-dual.

The full canonical-factor profiles are given in Example~\ref{v134:ex-seven-full-profile}; the same periods also determine the precision at which the trace pairing departs from self-duality.

\begin{example}[Trace-pairing precision for seventh roots]
\label{revadic:ex-eta}
Equation~\eqref{v137:eq-seven-periods} gives $\eta'-1=\eta^2$. For the trace Gram matrix $G$ of the normal integral basis above, let $J_3$ denote the all-ones matrix. Then
\[
 G-I_3=\eta(J_3-I_3)+\eta^2I_3.
\]
The off-diagonal entries have valuation exactly one, while the diagonal deviations have valuation exactly two. This lift is self-dual modulo two but is no longer self-dual modulo four.
On the other hand,
\[
 \det G=(\eta'-\eta)^2(\eta'+2\eta)
 =-7(\eta-1)\in\mathbb Z_2^\times.
\]
The trace pairing therefore remains nondegenerate at every precision. Loss of self-duality and degeneracy of the pairing are distinct, and the original period identities decide both.
\end{example}

\section{Mixed- and equal-characteristic lifting}
\label{revadic:sec-carries}

\begin{proposition}[No ring section in mixed characteristic]
\label{revadic:prop-section}
The reduction map $\mathcal O_n\twoheadrightarrow k_n$ has no unital ring-homomorphic section. The Teichmüller map $a\mapsto[a]$ preserves multiplication but not, in general, addition; every $x\in\mathcal O_n$ has a unique expansion
\begin{equation}\label{revadic:eq-digits}
 x=\sum_{j\ge0}p^j[a_j],\qquad a_j\in k_n.
\end{equation}
\end{proposition}
\begin{proof}
A unital ring homomorphism $s:k_n\to\mathcal O_n$ would carry $p\cdot1_{k_n}=0$ to $p\cdot1_{\mathcal O_n}=0$, contradicting characteristic zero. To construct the expansion, choose $a_0=\bar x$ and divide $x-[a_0]$ by $p$; continue recursively. Completeness gives convergence and $p$-adic separatedness gives uniqueness. Multiplicativity of Teichmüller representatives follows from uniqueness of Hensel lifts of prime-to-$p$ roots.
\end{proof}

Proposition~\ref{v133:prop-first-carry} expresses the first Witt carry as a finite-field polynomial. For the present purpose, write
\[
 [a]+[b]=[a+b]+p[\kappa(a,b)]\pmod{p^2}.
\]
Thus $\kappa(a,b)$ measures the first deviation of the finite-field relation $a+b$ when lifted to mixed characteristic. If $\sum_i a_i=0$, one may only infer $\sum_i[a_i]\in p\mathcal O_n$; it might vanish identically or become nonzero at some higher precision. The finite-field relation alone does not decide which.

\begin{proposition}[Equal-characteristic lifting has an actual ring section]
\label{revadic:prop-equal}
Let $k=\mathbb F_{p^n}$. The formal power-series ring $k[[t]]$ is a complete discrete valuation ring with residue field $k$, and the constant embedding
\[
 k\hookrightarrow k[[t]]\longrightarrow k[[t]]/(t)\simeq k
\]
is a unital ring section preserving both addition and multiplication. Although $k[[t]]$ and $\mathcal O_n$ have residue fields of the same size, they have different characteristics. If $f_0(X)\in k[X]$ has a simple root $a$, then every $F(X,t)\in k[[t]][X]$ with $F(X,0)=f_0(X)$ admits a unique root $x(t)\in k[[t]]$ with $x(0)=a$.
\end{proposition}
\begin{proof}
The first statements follow from coefficientwise arithmetic on formal power series. The last statement follows from $\partial_XF(a,0)\ne0$ and the $t$-adic Hensel lemma.
\end{proof}

\begin{example}[Two lifts of the same root]
\label{revadic:ex-two-lifts}
Let $k=\mathbb F_4=\mathbb F_2(a)$, with $a^2+a+1=0$. In the equal-characteristic ring $k[[t]]$, the constant $a$ satisfies the very same equation. In $\mathcal O_2=\mathbb Z_2[\zeta_3]$, its Teichmüller lift $[a]=\zeta_3$ also satisfies the corresponding integral polynomial, but the relation $a+a^2=1$ becomes $\zeta_3+\zeta_3^2=-1\ne1$. The first construction retains all ring identities of the constant field; the second retains specified integral root equations and multiplication, while its $2$-adic carry records the discrepancy.
\end{example}

\subsection{The first carry, multiplicative rigidity, and translation stabilizers}
\label{v140:sec-carry-rigidity}
The first carry records more than the discrepancy in one addition. A linear map or subspace that accommodates every first carry is constrained by field multiplication. Throughout, let $q=p^f$ and $E=\mathbb F_{q^n}$, and write $x^{1/p}$ for the unique $p$th root of $x$ in $E$. By Proposition~\ref{v133:prop-first-carry}, with
\[
 a_i=-\frac1p\binom pi\bmod p\in\mathbb F_p^\times
 \quad(1\le i<p),
\]
we have, for $t\in\mathbb F_p$,
\begin{equation}\label{v140:eq-carry-polarization}
 \kappa(x,ty)=\sum_{i=1}^{p-1}a_i t^{p-i}(x^iy^{p-i})^{1/p}.
\end{equation}
A polynomial of degree at most $p-1$ is determined by its values on $\mathbb F_p$. Hence each coefficient in~\eqref{v140:eq-carry-polarization} is an $\mathbb F_p$-linear combination of these carry values.

\begin{theorem}[Linear symmetries of the Teichmüller set]
\label{v140:thm-teich-rigidity}
Let $s\ge2$, $R_s=\mathcal O_f/p^s\mathcal O_f$, and $A_s=\mathcal O_{fn}/p^s\mathcal O_{fn}$, and let
\[
 \mathcal T_s=\{[x]:x\in E\}\subseteq A_s
\]
be the Teichmüller set. An $R_s$-linear automorphism $L$ satisfies $L(\mathcal T_s)=\mathcal T_s$ if and only if there are unique $a\in E^\times$ and $0\le j<n$ such that
\begin{equation}\label{v140:eq-teich-rigidity}
 L(z)=[a]\phi^{fj}(z)\qquad(z\in A_s).
\end{equation}
The same conclusion holds for $\mathcal O_f$-linear automorphisms of the complete integer ring $\mathcal O_{fn}$.
\end{theorem}
\begin{proof}
Reduction gives an $\mathbb F_q$-linear bijection $g:E\to E$. Uniqueness of canonical representatives gives $L([x])=[g(x)]$. Dividing $L$ by the unit $[g(1)]$, we may suppose $g(1)=1$. Apply $L$ to the first-carry identity modulo $p^2$ and use $g(x+y)=g(x)+g(y)$ to obtain
\[
 g(\kappa(x,y))=\kappa(g(x),g(y)).
\]
Replace $y$ by $ty$ and compare coefficients in~\eqref{v140:eq-carry-polarization}. Since $a_i\ne0$, this gives
\begin{equation}\label{v140:eq-mixed-monomials}
 g\bigl((x^iy^{p-i})^{1/p}\bigr)
       =\bigl(g(x)^ig(y)^{p-i}\bigr)^{1/p}
       \quad(1\le i<p).
\end{equation}
Taking $i=1$ and $y=1$ shows that $g$ commutes with inverse Frobenius, hence also with Frobenius. If $p=2$, equation~\eqref{v140:eq-mixed-monomials} now gives $g(xy)=g(x)g(y)$. If $p>2$, take $i=2$ and $y=1$ to obtain $g(x^2)=g(x)^2$; expansion of $(x+y)^2$ and division by two again give multiplicativity. Thus $g$ is a field automorphism fixing $\mathbb F_q$, so $g(x)=x^{q^j}$. Restoring the scalar, $L$ agrees with $[a]\phi^{fj}$ on all canonical representatives. The finite digit expansion then gives~\eqref{v140:eq-teich-rigidity} on $A_s$.

Conversely, Frobenius and Teichmüller multiplication preserve the representative set. The value $L(1)$ determines $a$, and the residue-field automorphism determines $j$, proving uniqueness. On the complete integer ring, linear maps are continuous, so the same argument applies to the convergent digit expansion.
\end{proof}

At $s=1$ the stabilizer is the full general linear group $\operatorname{GL}_{\mathbb F_q}(E)$. The carry condition in the second digit forces the normalized linear map to preserve multiplication.

\begin{theorem}[Additive closure with the same subspace in every digit]
\label{v140:thm-digit-subspaces}
Let $0\ne V\subseteq E$ be an $\mathbb F_q$-linear subspace and $s\ge2$, and define
\[
 \mathcal T_s(V)=\left\{\sum_{i=0}^{s-1}p^i[x_i]:x_i\in V\right\}\subseteq A_s.
\]
The following conditions are equivalent:
\begin{enumerate}
\item $\mathcal T_s(V)$ is an additive subgroup;
\item $\kappa(V,V)\subseteq V$;
\item $V=a\mathbb F_{q^d}$ for some $a\in E^\times$ and $d\mid n$.
\end{enumerate}
Under the third condition, at every precision,
\begin{equation}\label{v140:eq-subspace-ring}
 \mathcal T_s(V)=[a]\bigl(\mathcal O_{fd}/p^s\mathcal O_{fd}\bigr).
\end{equation}
\end{theorem}
\begin{proof}
The first condition remains true after reduction modulo $p^2$; comparison of the second digit of $[x]+[y]$ gives the second condition. Choose $a\in V\setminus\{0\}$. Since $\kappa(ax,ay)=a\kappa(x,y)$, multiplication by $[a]^{-1}$ allows us to assume $1\in V$.

Under the second condition, coefficient extraction in~\eqref{v140:eq-carry-polarization} gives $(x^iy^{p-i})^{1/p}\in V$ for $x,y\in V$ and $1\le i<p$. With $i=1,y=1$, we obtain $V^{1/p}\subseteq V$. Inverse Frobenius is injective on this finite set, so $V^{1/p}=V$ and $V^p=V$. If $p=2$, this immediately yields $xy\in V$. If $p>2$, taking $i=2,y=1$ first gives $x^2\in V$, and polarization gives multiplicative closure. A finite subring containing one is a field; since it contains $\mathbb F_q$, it is $\mathbb F_{q^d}$ for some $d\mid n$. Restoring the scalar gives the third condition.

Finally, the digit expansion in an unramified subring gives~\eqref{v140:eq-subspace-ring}, so the third condition implies the first.
\end{proof}

When full closure fails, the surviving translations can still be determined exactly. For a subset $B$ of an additive group, write $H_+(B)=\{u:B+u=B\}$.
\Needspace{8\baselineskip}
\begin{proposition}[The translation stabilizer of a two-digit set]
\label{v140:prop-digit-stabilizer}
Put
\[
 U(V)=\{a\in V:\kappa(a,x)\in V\text{ for every }x\in V\}.
\]
Then $U(V)$ is an $\mathbb F_q$-linear subspace, and
\begin{equation}\label{v140:eq-digit-stabilizer}
\begin{aligned}
 H_+(\mathcal T_2(V))&=\{[a]+p[b]:a\in U(V),\ b\in V\},
 \\ |H_+(\mathcal T_2(V))|&=|U(V)|\,|V|.
\end{aligned}
\end{equation}
If $p=2$, then also $U(V)=\{a\in V:aV\subseteq V^{[2]}\}$, where $V^{[2]}=\{x^2:x\in V\}$ is the pointwise square image, distinct from the product of subspaces in Definition~\ref{v70:def-subspace-product}.
\end{proposition}
\begin{proof}
A stabilizing translation must lie in $\mathcal T_2(V)$, since this set contains zero. Adding $[a]+p[b]$ to an arbitrary $[x]+p[y]$ gives digits $a+x$ and $b+y+\kappa(a,x)$, proving the stated criterion. Inclusion under translation is equality because the set is finite. The stabilizer is an additive subgroup whose reduction is $U(V)$. Moreover, $[\lambda]\mathcal T_2(V)=\mathcal T_2(V)$ for every $\lambda\in\mathbb F_q^\times$, so this image is closed under $\mathbb F_q$-scalars. The characteristic-two formula follows from $\kappa(a,x)=\sqrt{ax}$.
\end{proof}

\begin{theorem}[Nonzero translation directions of the first carry]\label{v147:thm-carry-directions}
Let $0\ne V\subseteq E=\mathbb F_{q^n}$ be an $\mathbb F_q$-subspace, retain $U(V)$ from Proposition~\ref{v140:prop-digit-stabilizer}, and write
\[
 H(V)=\{c\in E:cV\subseteq V\}.
\]
This is the maximal field of scalars in Proposition~\ref{v70:prop-scalar-support}.
\begin{enumerate}
\item If $p>2$, then $U(V)\ne0$ if and only if $V=a\mathbb F_{q^d}$ for some $a\ne0$ and $d\mid n$; in that case $U(V)=V$. Thus, if $V$ is not a scalar multiple of a subfield, then
\begin{equation}\label{v147:eq-odd-carry-stabilizer}
 H_+(\mathcal T_2(V))=p[V],\qquad
 |H_+(\mathcal T_2(V))|=|V|.
\end{equation}
\item If $p=2$, then $U(V)\ne0$ if and only if $V^{[2]}=aV$ for some $a\in V\setminus\{0\}$. For every such $a$,
\begin{equation}\label{v147:eq-binary-carry-directions}
 U(V)=aH(V),\qquad
 |H_+(\mathcal T_2(V))|=|H(V)|\,|V|.
\end{equation}
Equivalently, $a^{-1}V$ contains one and is stable under the squaring Frobenius; it need not be closed under multiplication.
\end{enumerate}
Here $p[V]=\{p[x]:x\in V\}\subseteq A_2$.
\end{theorem}
\begin{proof}
Take $a\in U(V)\setminus\{0\}$ and put $W=a^{-1}V$. Homogeneity of the first carry gives $1\in W$ and $\kappa(1,x)\in W$ for every $x\in W$.

If $p>2$, applying~\eqref{v140:eq-carry-polarization} for $t\in\mathbb F_p$ gives
\[
 \kappa(1,tx)=\sum_{j=1}^{p-1}a_{p-j}t^j(x^j)^{1/p}\in W.
\]
Extracting the coefficients of $t$ and $t^2$ yields $x^{1/p},(x^2)^{1/p}\in W$. The first inclusion and finiteness imply that inverse Frobenius maps $W$ bijectively to itself. Hence $W$ is also closed under $p$th powers, so $x^2\in W$. The identity
\[
 xy=\frac{(x+y)^2-x^2-y^2}{2}
\]
then gives multiplicative closure. Thus $W$ is a finite subfield containing $\mathbb F_q$, and $W=\mathbb F_{q^d}$. Conversely, a scalar multiple of a subfield is closed under the first carry, so $U(V)=V$. Formula~\eqref{v147:eq-odd-carry-stabilizer} now follows from Proposition~\ref{v140:prop-digit-stabilizer}.

If $p=2$, the same proposition gives $aV\subseteq V^{[2]}$. The two spaces have equal dimension, so they are equal. Conversely, if $a\in V\setminus\{0\}$ satisfies this equality, then for $b\in V$,
\[
 b\in U(V)\iff bV\subseteq aV\iff b/a\in H(V).
\]
Since $a\in V$, we have $aH(V)\subseteq V$, proving $U(V)=aH(V)$. Finally, $V^{[2]}=aV$ is equivalent to $(a^{-1}V)^{[2]}=a^{-1}V$, while $a\in V$ is equivalent to $1\in a^{-1}V$.
\end{proof}

\begin{example}[The first carry of a plane in the eight-element field]
\label{v140:ex-plane-carry}
Let $E=\mathbb F_2(\alpha)$ with $\alpha^3=\alpha+1$, and put $V=\operatorname{span}_{\mathbb F_2}\{1,\alpha\}$. Then
\[
 \kappa(1,\alpha)=\alpha^4=\alpha^2+\alpha\notin V.
\]
Thus the sixteen elements of $\mathcal T_2(V)$ do not form an additive subgroup. More precisely, $V^{[2]}=\operatorname{span}\{1,\alpha^2\}$ and $V\cap V^{[2]}=\mathbb F_2$. If $a\in U(V)$, then $a\cdot1\in V^{[2]}$ forces $a\in\mathbb F_2$, while $1\notin U(V)$. Hence $U(V)=0$, and the full translation stabilizer consists of the four elements of $2[V]$.
\end{example}

\begin{corollary}[Recovering a subfield from the carry closure of a subspace]
\label{v144:cor-carry-closure}
Let $0\ne V_0\subseteq E$ be an $\mathbb F_q$-linear subspace, and define
\[
 V_{i+1}=V_i+\operatorname{span}_{\mathbb F_q}
                 \{\kappa(x,y):x,y\in V_i\}.
\]
For $a\in V_0\setminus\{0\}$, put $K=\mathbb F_q(a^{-1}V_0)$ and $d=[K:\mathbb F_q]$. The chain stabilizes at $aK$, with at most $d-\dim_{\mathbb F_q}V_0$ strict inclusions. The stable term is called the \emph{carry closure} of $V_0$.
\end{corollary}
\begin{proof}
The space $aK$ is closed under the first carry and contains every $V_i$. Bounded dimension forces stabilization. By Theorem~\ref{v140:thm-digit-subspaces}, the stable term is a scaled subfield containing $a$, hence has the form $aK'$. Containment of $V_0$ forces $K\subseteq K'$, and the reverse inclusion is already known. Each strict inclusion increases dimension by at least one.
\end{proof}

The plane in Example~\ref{v140:ex-plane-carry} becomes the whole field of eight elements after adjoining $\kappa(1,\alpha)=\alpha^2+\alpha$. This construction adjoins carries and takes linear spans; a map preserving only the original partial set of canonical representatives need not preserve the subsequent closures.

\subsection{First-carry obstructions on partial root sets}
\label{v144:sec-partial-carry}
Preserving all canonical representatives makes the first carry recover field multiplication. On a multiplicative subgroup, mixed exponents may coincide modulo the subgroup order; these coincidences determine which multiplicative relations the first carry can still distinguish.

Let $q=p^f$, $p\nmid N$, $n=\ord_N(q)\ge2$, and put $E=\mathbb F_{q^n}$ and $H=\mu_N\subset E^\times$. Retain
$R_s=\mathcal O_f/p^s\mathcal O_f$ and $A_s=\mathcal O_{fn}/p^s\mathcal O_{fn}$, and write $[H]$ for the canonical representatives of $H$. Take an $\mathbb F_q$-linear bijection
\[
 g(X)=\sum_{i=0}^{n-1}c_iX^{q^i},\qquad g(H)=H.
\]
For a multi-index $\mathbf e=(e_0,\ldots,e_{n-1})$, write $|\mathbf e|=\sum_i e_i$ and $\operatorname{supp}\mathbf e=\{i:e_i>0\}$, and set
\begin{equation}\label{v144:eq-carry-obstruction}
\begin{split}
 W_g(X)&=\sum_{\substack{|\mathbf e|=p\\|\operatorname{supp}\mathbf e|\ge2}}
 \gamma_{\mathbf e}\prod_i c_i^{e_i}\,
 X^{\sum_i e_iq^i},\\
 \gamma_{\mathbf e}&=\frac1p\binom p{e_0,\ldots,e_{n-1}}\bmod p\in\mathbb F_p^\times,
 \qquad P_N=\{pq^i\bmod N:0\le i<n\}.
\end{split}
\end{equation}
The displayed multinomial coefficients are divisible by $p$ exactly once, so each $\gamma_{\mathbf e}$ is well defined and nonzero. After reduction modulo $X^N-1$, the coefficients outside $P_N$ form the \emph{first-carry obstruction} to preserving $[H]$.

\begin{proposition}[A two-digit lifting criterion for subgroup representatives]
\label{v144:prop-partial-lift}
There is an $R_2$-linear automorphism $L:A_2\to A_2$ satisfying $L([x])=[g(x)]$ for $x\in H$ if and only if
\begin{equation}\label{v144:eq-partial-lift}
 \operatorname{supp}(W_g\bmod(X^N-1))\subseteq P_N.
\end{equation}
When this holds, write the unique remainder as
\[
 W_g(X)\equiv\sum_{i=0}^{n-1}b_iX^{pq^i}\pmod{X^N-1}.
\]
The lift is unique and equals
\begin{equation}\label{v144:eq-partial-lift-map}
 L(z)=\sum_{i=0}^{n-1}\bigl([c_i]+p[b_i^{1/p}]\bigr)\phi^{fi}(z).
\end{equation}
The exponents in the remainder are read modulo $N$; the $n$ elements of $P_N$ are distinct.
\end{proposition}
\begin{proof}
Choose an integral basis of $\mathcal O_{fn}/\mathcal O_f$. Its conjugate matrix reduces to an invertible Moore matrix, so its determinant is a unit. Every $R_2$-linear map therefore has a unique expression $\sum_i C_i\phi^{fi}$ with $C_i\in A_2$. If its reduction is $g$, write $C_i=[c_i]+p[d_i]$.

For $y=[a]+p[b]\in A_2$, one has $y^p-\phi(y)=-p[b^p]$. Thus $y$ is a canonical representative if and only if $y^p=\phi(y)$. Expanding $L([x])^p$ for $x\in H$ and subtracting $\phi(L([x]))$ gives
\[
 \frac{L([x])^p-\phi(L([x]))}{p}\bmod p
 =W_g(x)-\sum_i d_i^p x^{pq^i}.
\]
Division by $p$ followed by reduction denotes the natural identification $pA_2\simeq E$. The roots of $X^N-1$ are exactly $H$ and are simple. Vanishing on $H$ is therefore equivalent to a zero remainder modulo $X^N-1$. Since $\ord_N(q)=n$, the pure exponents are distinct, giving precisely~\eqref{v144:eq-partial-lift} and $d_i=b_i^{1/p}$. The constructed map reduces to a bijection and hence is invertible; its values on $H$ are the canonical representatives of $g(x)$.

If $\theta\in H$ has order $N$, then $1,\theta,\ldots,\theta^{n-1}$ form an $\mathbb F_q$-basis of $E$, and their canonical lifts form an integral basis. Thus $[H]$ spans $A_2$, which also proves uniqueness directly from the prescribed values.
\end{proof}

The criterion retains cancellation between mixed terms; it does not require each term to vanish separately. Separated exponents impose stronger restrictions.

\begin{corollary}[Separated mixed exponents force semilinearity]
\label{v144:cor-partial-rigidity}
Suppose that for every pair $i\ne j$ there is a multi-index $\mathbf e$ of total weight $p$ and support exactly $\{i,j\}$ such that
\[
 \sum_k e_kq^k\pmod N
\]
lies outside $P_N$ and differs from the exponent of every other mixed multi-index modulo $N$. Then, for every $s\ge2$, all $R_s$-linear automorphisms preserving $[H]$ are exactly
\begin{equation}\label{v144:eq-partial-group}
 z\longmapsto[a]\phi^{fj}(z),\qquad a\in H,\quad0\le j<n.
\end{equation}
The same holds over the complete integer ring. This group has order $Nn$; its residue maps are semilinear for the one-dimensional $E$-structure.
\end{corollary}
\begin{proof}
Reduce modulo $p^2$ and apply Proposition~\ref{v144:prop-partial-lift}. Each selected mixed exponent has a single coefficient, which must vanish, so $c_i^{e_i}c_j^{e_j}=0$. No two coefficients can be nonzero. Hence $g(x)=ax^{q^j}$, with $a\in H$ because $g(H)=H$. The original map takes $[x]$ to the unique canonical representative $[a][x]^{q^j}$ for $x\in H$. These representatives contain an integral basis, proving~\eqref{v144:eq-partial-group} at every precision and over the complete ring. Multiplication and Frobenius give the converse and the group order.
\end{proof}

\begin{example}[Multiplicative rigidity on two proper subgroups]
\label{v144:ex-eighty-five}
Take $p=q=2$ and $N=85$. Then $n=8$, and $H$ has index three in $\mathbb F_{256}^\times$. The mixed terms are $c_ic_jX^{2^i+2^j}$ for $i<j$. The cyclic distance between the indices is $1,2,3$, or $4$, giving the last four orbits in the following table:
\[
\begin{array}{c|l}
 a&C_2(a)\pmod{85}\\\hline
 1&1,2,4,8,16,32,64,43\\
 3&3,6,12,24,48,11,22,44\\
 5&5,10,20,40,80,75,65,45\\
 9&9,18,36,72,59,33,66,47\\
 17&17,34,68,51
\end{array}
\]
The last four orbits are disjoint, have $8+8+8+4=\binom82$ elements in total, and avoid the pure orbit $C_2(1)$. Thus every mixed exponent is separated. Modulo four, the linear automorphisms preserving these eighty-five canonical roots are already restricted to $[a]\phi^j$, giving $85\cdot8=680$ maps. Preservation of all 255 nonzero canonical representatives is unnecessary.

The same separation occurs in odd characteristic. Take $p=q=3$ and $N=13$, so $n=3$ and $H$ has index two in $\mathbb F_{27}^\times$. The pure exponents are $\{1,3,9\}$, while the seven mixed multi-indices have exponents
\[
 5,11,7,2,6,8,0\pmod{13}.
\]
They are distinct and avoid the pure exponents. Modulo nine, exactly $13\cdot3=39$ linear automorphisms preserve $[H]$.
\end{example}

\begin{example}[Five roots of unity admit extra permutations]
\label{v144:ex-five-permutations}
Take $p=q=2$, $N=5$, and $n=4$. A primitive fifth root $\theta$ satisfies
\[
 1+\theta+\theta^2+\theta^3+\theta^4=0,
\]
and any four roots form an integral basis of $\mathcal O_4/\mathbb Z_2$. Indeed, $\Phi_5\bmod2$ is irreducible, and deleting any one root gives a coordinate determinant that is a unit.
Every permutation of the five roots therefore extends uniquely to an integral linear automorphism and preserves the displayed relation. The full stabilizer is $S_5$, at every finite precision as well. The semilinear maps $[a]\phi^j$ give only $5\cdot4=20$ elements. A transposition fixing the other three roots is not among them, since a nonidentity affine map $i\mapsto b+2^ji$ modulo five fixes at most one index.

Here $H$ still generates the whole degree-four extension, but mixed exponents coincide. Multiplicative generation and preservation of these representatives at every precision do not replace the carry-separation condition.

Ninth roots give a block decomposition as well. For a primitive ninth root $\vartheta\in\mathcal O_6$, the equation $\Phi_9=X^6+X^3+1$ partitions the roots into
\[
 \{\vartheta^r,\vartheta^{r+3},\vartheta^{r+6}\},\qquad r=0,1,2.
\]
Each triple sums to zero, and choosing any two roots from each triple gives an integral basis of size six. Arbitrary permutations of the triples and within them extend to linear automorphisms, giving a subgroup of order $6^3\cdot6=1296$, larger than the semilinear subgroup of order $9\cdot6=54$. These maps preserve three rank-two root blocks without necessarily preserving multiplication in the degree-six field.
\end{example}

\subsection{Precision of low-degree subfield representations}
\label{v140:sec-profile}
Subfields lift canonically along the unramified tower, but their intersections with low-degree power spaces need not be preserved. Let $\alpha\in E=\mathbb F_{q^n}$ have degree $n$ over $\mathbb F_q$, put $\theta=[\alpha]\in\mathcal O_{fn}$, and take $d\mid n$. By the integral-basis reduction criterion, $1,\theta,\ldots,\theta^{n-1}$ is an integral basis over $\mathcal O_f$. For $0\le t<n$, set
\[
 L_t=\operatorname{span}_{\mathcal O_f}\{1,\theta,\ldots,\theta^t\},
 \qquad D_{t,d}=(\phi^{fd}-1)|_{L_t}:L_t\longrightarrow\mathcal O_{fn}.
\]
Let $\bar D_{t,d,s}$ denote its reduction modulo $p^s$, and define
\begin{equation}\label{v140:eq-precision-profile}
 W_{t,d}^{(s)}=\operatorname{im}\left(\ker\bar D_{t,d,s}
                   \longrightarrow L_t/pL_t\right),
 \qquad \rho_s(t,d)=\dim_{\mathbb F_q}W_{t,d}^{(s)}.
\end{equation}
We call $\rho_s(t,d)$ the \textbf{precision profile of representations}. By Proposition~\ref{revadic:prop-fixed-precision}, this space records exactly the residue classes that lift to subring elements modulo $p^s$ without exceeding degree $t$. At $s=1$ it is the degree profile of $\mathbb F_{q^d}\cap\mathcal P_t(\alpha)$ from Chapter~13. It measures liftable residue classes, rather than the sizes of kernels over different quotient rings.

\Needspace{10\baselineskip}
\begin{theorem}[Smith exponents and representation precision]
\label{v140:thm-profile-smith}
Suppose $D_{t,d}$ has rank $r$ and nonzero \textbf{Smith factors} $p^{e_1},\ldots,p^{e_r}$, with $0\le e_1\le\cdots\le e_r$. Then
\begin{align}
 \rho_s(t,d)&=t+1-r+\#\{i:e_i\ge s\},\label{v140:eq-smith-profile}\\
 |\ker\bar D_{t,d,s}|&=q^{\,s(t+1-r)+\sum_{i=1}^r\min(s,e_i)}.\label{v140:eq-smith-count}
\end{align}
Moreover,
\[
 \operatorname{rank}_{\mathcal O_f}(L_t\cap\mathcal O_{fd})=t+1-r,
 \qquad
 \rho_s(t,d)-\rho_{s+1}(t,d)=\#\{i:e_i=s\}.
\]
Once the precision exceeds every positive exponent, $W_{t,d}^{(s)}$ is the reduction of the integral intersection $L_t\cap\mathcal O_{fd}$. In particular, every low-degree representation modulo $p$ persists to arbitrary precision if and only if there are no positive Smith exponents.
\end{theorem}
\begin{proof}
Over the discrete valuation ring $\mathcal O_f$, choose a nonzero matrix entry of least valuation as pivot. It divides every other entry, so invertible row and column operations clear its row and column. Induction on the remaining matrix gives diagonal entries $p^{e_i}$ and zero columns, absorbing units into the invertible operations. This is the Smith form used here.

In these coordinates the kernel conditions are $p^{e_i}u_i\equiv0\pmod{p^s}$. If $e_i<s$, then $u_i$ is divisible by $p^{s-e_i}$ and has zero residue modulo $p$. If $e_i\ge s$, its residue is arbitrary. The remaining $t+1-r$ zero columns are unrestricted, giving~\eqref{v140:eq-smith-profile}. The $i$th nonzero diagonal condition has $q^{\min(s,e_i)}$ solutions, and each zero column has $q^s$, proving~\eqref{v140:eq-smith-count}.

The integral kernel consists of the zero-column directions and equals $L_t\cap\mathcal O_{fd}$ by the fixed-ring criterion. It is a direct summand, so its reduction has dimension equal to its rank. Subtracting adjacent profiles gives the jump formula. Once $s$ exceeds every $e_i$, only the zero-column directions survive in the image modulo $p$. All invertible coordinate changes commute with reduction, so the conclusions are independent of the chosen bases.
\end{proof}

\begin{example}[Quadratic subfield representatives modulo four]
\label{v140:ex-quartic-depth}
Use Example~\ref{sfr:ex-quartic-comparison}: $\alpha^4+\alpha+1=0$, and $\beta=\alpha^2+\alpha$ generates $\mathbb F_4/\mathbb F_2$. Take $\theta=[\alpha]\in\mathcal O_4$ and put $S=\theta+\theta^4$. Since $\alpha+\alpha^4=1$, the first carry gives
\begin{equation}\label{v140:eq-quartic-carry}
 S\equiv1+2[\alpha^{10}]\pmod4.
\end{equation}
Here $\alpha$ has order fifteen: degree four excludes orders one and three, while $\alpha^5=\alpha^2+\alpha\ne1$ excludes order five. Thus $u=\alpha^{10}$ has order three and lies in $\mathbb F_4\setminus\mathbb F_2$.

For $z=c_0+c_1\theta+c_2\theta^2$ with $c_i\in\mathbb Z_2$, the exact identity
\[
 (\phi^2-1)z=(\theta^4-\theta)(c_1+c_2S)
\]
follows from $\theta^8-\theta^2=(\theta^4-\theta)(\theta^4+\theta)$. The first factor is a unit, since its residue is one. If $z\bmod4$ belongs to $\mathcal O_2/4\mathcal O_2$ and its residue generates the four-element subfield, both $c_1$ and $c_2$ are odd and
\[
 c_1+c_2S\equiv0\pmod4.
\]
Substituting~\eqref{v140:eq-quartic-carry}, dividing by two, and reducing gives $u$ plus an element of $\mathbb F_2$, which cannot vanish.

The calculation also determines the whole precision profile. Write $S=1+2v$. Since $\bar v=u$, the residues of $1,v$ are independent and extend to an integral basis. Multiplication by the unit $\theta^4-\theta$ preserves Smith factors; the columns $(1,S)$ become $(1,2v)$. Hence the nonzero Smith exponents of $D_{2,2}$ are exactly $0,1$, and
\[
 \rho_1(2,2)=2,\qquad \rho_s(2,2)=1\ (s\ge2),
 \qquad L_2\cap\mathcal O_2=\mathbb Z_2.
\]
The full integral power basis provides a cubic representative of $[\beta]$. Thus subring elements whose residues generate $\mathbb F_4/\mathbb F_2$ require degree at least two modulo two, and degree at least three at every precision $2^s$ with $s\ge2$ and in the complete integer ring; these bounds are attained.

Low-degree elements of higher divisibility may still occur: $2(\theta^2+\theta)\bmod4$, for example, belongs to the quadratic subring but has zero residue. It does not contradict the statement about generators. On the other hand, choose the Hensel root $\vartheta$ above $\alpha$ of the integral polynomial
\[
 F(X)=(X^2+X)^2+(X^2+X)+1.
\]
Then $\vartheta^2+\vartheta$ satisfies $Y^2+Y+1=0$ exactly, so its quadratic representation persists. This lift differs from the canonical root of unity $\theta$; a precision profile must be considered relative to the chosen lift.
\end{example}

\section{Cyclotomic orbits, periodic digits, and coefficient valuations}
\label{revadic:sec-periods}

A cyclotomic coset can be recorded by the periodic digits of a rational number, or used to form polynomial coefficients through canonical roots of unity. The former turns orbit steps into digit shifts and connects with Gauss-sum valuations; the latter also requires detecting cancellation among subset sums, distinguishing exact zeros from finite divisibility.

\subsection{Cyclotomic cosets, periodic expansions, and digit sums}
\label{v138:subsec-periodic-digits}
Let $q=p^f$, $m\ge2$, and $p\nmid m$. We use the integer digits $0,1,\ldots,q-1$, listed from low to high powers; each base-$q$ digit splits into $f$ base-$p$ digits. This is the ordinary positional expansion in $\mathbb Z_p$, distinct from the Teichmüller-digit expansion in~\eqref{revadic:eq-digits}.
An expansion is \textbf{purely periodic} if a finite block repeats from its first digit. Write $\operatorname{Per}_q(x)$ for its shortest period block, again in low-to-high order.

\begin{theorem}[Cyclotomic cosets and purely periodic digit blocks]
\label{v138:thm-periodic-cosets}
Take $1\le a<m$, put $d=|C_q(a)|$, and set
\begin{equation}\label{v138:eq-periodic-block}
 \beta_a=-\frac am,\qquad
 E_a=\frac{a(q^d-1)}m
     =\sum_{i=0}^{d-1}e_iq^i,\qquad 0\le e_i<q.
\end{equation}
Then $E_a$ is an integer with $0<E_a<q^d-1$, and in $\mathbb Z_p$,
\[
 \beta_a=\sum_{r\ge0}E_aq^{dr}
        =\sum_{r\ge0}\sum_{i=0}^{d-1}e_iq^{dr+i}.
\]
Its least base-$q$ digit period is exactly $d$. If $qa=mb+r$, $0<r<m$, then
\begin{equation}\label{v138:eq-digit-rotation}
 \beta_r=q\beta_a+b,\qquad b=e_{d-1},\qquad
 \operatorname{Per}_q(\beta_r)=(e_{d-1},e_0,\ldots,e_{d-2}).
\end{equation}
For $d=1$, the rotation is the identity.
Thus the nonzero $q$-cyclotomic cosets modulo $m$ correspond bijectively to the cyclic equivalence classes of these period blocks.
\end{theorem}
\begin{proof}
Divisibility $m\mid a(q^d-1)$ makes $E_a$ integral, and $0<a<m$ gives its range. The $p$-adic geometric series gives $E_a/(1-q^d)=-a/m$. Successive blocks occupy disjoint sets of $d$ digit positions, so this is the ordinary positional expansion.
If it has period $h$, then $(1-q^h)\beta_a$ is an integer. Hence $q^ha\equiv a\pmod m$, so $d\mid h$. Since $d$ is already a period, it is the least one.

Reduction of $qE_a$ modulo $q^d-1$ gives
\[
 e_{d-1}+\sum_{i=0}^{d-2}e_iq^{i+1}.
\]
This still lies in $1,\ldots,q^d-2$ and is therefore $E_r$. Comparing with $qE_a=E_r+b(q^d-1)$ gives $b=e_{d-1}$ and~\eqref{v138:eq-digit-rotation}.
The $d$ elements of a coset give $d$ distinct rotations. Conversely, a period block determines its value $-a/m$ uniquely, so different cosets cannot give the same cyclic equivalence class.
\end{proof}

Long division also determines the direction. Put $r_0=a$ and write
\[
 qr_i=mb_i+r_{i+1},\qquad 0<r_{i+1}<m.
\]
Then $r_d=r_0$, and the real and $p$-adic expansions satisfy
\begin{equation}\label{v138:eq-real-adic-reversal}
 \frac am=0.\overline{b_0b_1\cdots b_{d-1}}{}_q,
 \qquad \operatorname{Per}_q(-a/m)=(b_{d-1},\ldots,b_0).
\end{equation}
Indeed, iterating the divisions gives $a(q^d-1)/m=\sum_{i=0}^{d-1}b_iq^{d-1-i}$. The same remainders determine both expansions, with opposite reading directions. The zero residue separately gives the all-zero block; the all-$(q-1)$ block represents $-1$ and is excluded by $1\le a<m$.

\begin{proposition}[Digit steps and Teichmüller orbits]
\label{v138:prop-period-steps}
Put $n=m/\gcd(a,m)$, $t=\operatorname{ord}_n(p)$, and $g=\gcd(t,f)$. Then
\begin{equation}\label{v138:eq-period-steps}
 d=\frac t g,\qquad fd=\operatorname{lcm}(t,f).
\end{equation}
The least base-$p$ digit period of $\beta_a$ is $t$. Splitting the $d$ base-$q$ digits in~\eqref{v138:eq-periodic-block} produces $fd$ digits, consisting of $f/g$ repetitions of the shortest base-$p$ block. The full $p$-cyclotomic coset splits into $g$ distinct $q$-cyclotomic cosets.

Choose a primitive $m$th root of unity $\zeta_m$ in a common unramified extension, and put $\xi_a=\zeta_m^a$. Then $\xi_a\in\mathcal O_{fd}$ and
\[
 \phi^f(\xi_a)=\xi_{qa},\qquad
 [K_f(\xi_a):K_f]=d.
\]
Reduction gives a finite-field Frobenius orbit of the same length. Thus in exponent, digit-block, and canonical-root coordinates, the same permutation is written as
\[
 \begin{aligned}
 a&\longmapsto qa\bmod m,\\
 (e_0,\ldots,e_{d-1})&\longmapsto(e_{d-1},e_0,\ldots,e_{d-2}),\\
 \xi_a&\longmapsto\phi^f(\xi_a).
 \end{aligned}
\]
\end{proposition}
\begin{proof}
The condition $q^ha\equiv a\pmod m$ is equivalent to $p^{fh}\equiv1\pmod n$, whose least positive solution is $t/g$. Applying Theorem~\ref{v138:thm-periodic-cosets} with base $p$ gives least digit period $t$. Taking steps of length $f$ around a cycle of length $t$ produces $g$ cycles of length $t/g$.
The order of $\xi_a$ is $n$, so $\xi_a\in K_t\subseteq K_{fd}$. Unramified Frobenius and reduction preserve the order of these prime-to-$p$ roots of unity, giving relative orbit length and extension degree $d$.
\end{proof}

Here $\beta_a$ encodes an exponent, whereas $\xi_a$ lifts a root of unity; $\beta_a$ is not an expansion of $\xi_a$. In particular, arithmetic Frobenius fixes $\mathbb Z_p$ pointwise, so $\phi(\beta_a)=\beta_a$. Equation~\eqref{v138:eq-digit-rotation} describes a map on the finite set $\{-a/m:1\le a<m\}$ induced by the exponent permutation, not the restriction of arithmetic Frobenius to $\mathbb Z_p$.

\begin{example}[Two cyclic word classes for seventh roots]
\label{v138:ex-seven-words}
Take $p=q=2,m=7$. Here $d=3$ and $E_a=a$. The blocks in the two cosets are
\[
\begin{array}{c|c|c|c}
 a&\operatorname{Per}_2(-a/7)&a&\operatorname{Per}_2(-a/7)\\\hline
 1&(1,0,0)&3&(1,1,0)\\
 2&(0,1,0)&6&(0,1,1)\\
 4&(0,0,1)&5&(1,0,1)
\end{array}
\]
Each column follows its coset. For instance,
\[
 -\frac17=1+2^3+2^6+\cdots,\quad
 -\frac27=2+2^4+2^7+\cdots,\quad
 -\frac47=2^2+2^5+2^8+\cdots.
\]
The real expansion $1/7=0.\overline{001}{}_2$ reads the same block in reverse. These two cyclic word classes correspond to the two root orbits in Example~\ref{v132:ex-seven-period}; the digit encoding does not change the canonical lift of any root.
\end{example}

\begin{example}[Grouping pairs of digits and splitting an orbit]
\label{v138:ex-four-fifteen}
Take $p=2,q=4,m=15$. The coset $C_2(1)=\{1,2,4,8\}$ splits into $C_4(1)=\{1,4\}$ and $C_4(2)=\{2,8\}$. Their base-four blocks are
\[
 (1,0)\longleftrightarrow(0,1),\qquad
 (2,0)\longleftrightarrow(0,2).
\]
Splitting them into binary digits in low-to-high order gives
\[
 (1,0\mid0,0)\longleftrightarrow(0,0\mid1,0),\qquad
 (0,1\mid0,0)\longleftrightarrow(0,0\mid0,1).
\]
Here $t=4,f=2,g=2,d=2$. Shifting one digit visits all four positions, whereas shifting two digits visits only two: one binary cyclotomic coset splits into two quaternary cosets.
\end{example}

Periodic digits also enter the Gauss-sum valuations of Chapter~\ref{chap:characters}. Fix a compatible $p$-adic embedding, put $Q=q^d$, let $\omega:\mathbb F_Q^\times\to\mu_{Q-1}$ be the Teichmüller character, and set $\psi_Q(x)=\zeta_p^{\operatorname{Tr}_{\mathbb F_Q/\mathbb F_p}(x)}$, where $\zeta_p$ is a primitive $p$th root of unity. Write
\[
 g_a(Q)=G(\omega^{-E_a},\psi_Q).
\]
Since $1\le E_a\le Q-2$, the multiplicative character is nontrivial. The required divisibility is $m\mid a(Q-1)$, not necessarily $m\mid Q-1$.

\begin{corollary}[Digit sums of period blocks and Gauss-sum valuations]
\label{v138:cor-word-gauss}
With the notation above and valuation normalized by $v_p(p)=1$,
\begin{equation}\label{v138:eq-word-gauss}
 v_p(g_a(Q))=\frac{s_p(E_a)}{p-1}
           =\frac1{p-1}\sum_{i=0}^{d-1}s_p(e_i).
\end{equation}
The Gauss sums themselves are equal on each $q$-cyclotomic coset, and even on the full $p$-cyclotomic coset.
If $E_a^{(0)}=a(p^t-1)/m$, then
\begin{equation}\label{v138:eq-normalized-digit-weight}
 \frac{v_p(g_a(Q))}{fd}
 =\frac{s_p(E_a^{(0)})}{(p-1)t}.
\end{equation}
Thus, as the grouping length $f$ varies, the valuation normalized by $\log_p Q=fd$ depends only on the shortest base-$p$ block.
\end{corollary}
\begin{proof}
The first equality in~\eqref{v138:eq-word-gauss} is the Stickelberger formula of Theorem~\ref{v132:thm-gauss-digit-valuation}. The second follows because $q=p^f$, so the blocks occupy disjoint base-$p$ positions.
If $r\equiv qa\pmod m$, then $E_r\equiv qE_a\pmod{Q-1}$. Substitution $y=x^q$ in the Gauss sum, together with invariance of the absolute trace, gives $g_r(Q)=g_a(Q)$; the same argument applies with $p$ in place of $q$.
Finally, $fd=(f/g)t$ and
\[
 E_a=E_a^{(0)}\bigl(1+p^t+\cdots+p^{(f/g-1)t}\bigr).
\]
Repeating the shortest block $f/g$ times multiplies its digit sum by $f/g$, proving~\eqref{v138:eq-normalized-digit-weight}.
\end{proof}

\begin{example}[Seventh roots: digit sums and an exact period identity]
\label{v138:ex-seven-gauss}
In Example~\ref{v138:ex-seven-words}, take $Q=8$. The digit sums of the two word classes are one and two, giving Gauss-sum valuations one and two.
The periods of Chapter~5 also give these values directly. Choose $\alpha=\bar\zeta$ with $\alpha^3+\alpha+1=0$, use $\eta,\eta'$ from Example~\ref{v132:ex-seven-period}, and take $\omega(\alpha)=\zeta$. The nonzero trace-zero elements are exactly $\alpha,\alpha^2,\alpha^4$. Since $\sum_{i=0}^6\zeta^{-ai}=0$,
\begin{equation}\label{v138:eq-seven-gauss-period}
 g_a(8)=2\sum_{i\in\{1,2,4\}}\zeta^{-ai}
 =\begin{cases}
 2\eta',&a\in\{1,2,4\},\\
 2\eta,&a\in\{3,5,6\}.
 \end{cases}
\end{equation}
The known values $v_2(\eta)=1$ and $v_2(\eta')=0$ recover the two valuations, and $g_1(8)g_3(8)=8$. Here the Fourier expansion of the trace-zero set connects Gauss sums to canonical-factor coefficients; the Gauss-sum digit formula cannot simply be assigned to the periods $\eta,\eta'$ themselves.
\end{example}

\begin{example}[The same binary period over different Gauss-sum fields]
\label{v138:ex-four-seven}
Take $p=2,q=4,m=7,a=1$. Now $t=d=3$, but $fd=6$, and
\[
 E_1=9=1+2\cdot4,
 \qquad \operatorname{Per}_4(-1/7)=(1,2,0).
\]
Its binary expansion is $(1,0\mid0,1\mid0,0)$, namely two repetitions of $(1,0,0)$. Thus $g_1(64)=G(\omega_{64}^{-9},\psi_{64})$ has valuation two, whereas $g_1(8)$ has valuation one; division by $6,3$, respectively, gives $1/3$ in both cases.
The character $\omega_{64}^{-9}$ is the norm lift of $\omega_8^{-1}$, and Hasse--Davenport also gives $g_1(64)=-g_1(8)^2$. The shortest binary period is unchanged; the finite field over which the Gauss sum is taken has changed.
\end{example}

\begin{example}[Ternary digit blocks for eighth roots]
\label{v138:ex-three-eight}
Take $p=q=3,m=8$. The coset $\{1,3\}$ has blocks $(1,0),(0,1)$, while $\{5,7\}$ has blocks $(2,1),(1,2)$. Over $Q=9$, their Gauss sums have valuations $1/2,3/2$, respectively. These use the same root orbits as Example~\ref{revadic:ex-eight}, but the additive character also introduces a cube root of unity; the fractional valuations are consistent with $v_3(\zeta_3-1)=1/2$.
For $a=4$, instead, $d=1$, $Q=3$, and $E_a=1$, while $\beta_a=-1/2$ has shortest ternary block $(1)$. Here $8\nmid Q-1$ but $8\mid4(Q-1)$; the corresponding root is $\zeta_8^4=-1$, and the Gauss sum is the quadratic Gauss sum over $\mathbb F_3$.
\end{example}

\subsection{Minimal rotations, reduction, and relative phases}
\label{v139:sec-digit-structure}
To compare different cosets modulo the same $m$, use a common block length
\[
 L=\ord_m(q),\qquad M=q^L-1,\qquad c=M/m,
 \qquad E_a=ca\quad(0\le a<m).
\]
Here $E_a$ represents $L$ digits, including leading zeros; its shortest period need not be $L$. On the multiples of $c$ modulo $M$, put
\begin{equation}\label{v139:eq-two-actions}
 R(E)=[qE]_M,\qquad A(E)=[E+c]_M.
\end{equation}
Then $R(E_a)=E_{[qa]_m}$ and $A(E_a)=E_{[a+1]_m}$. Thus $R$ rotates digits, while $A$ performs addition with carries modulo $q^L-1$. A carry past position $L-1$ returns to position zero because $q^L\equiv1\pmod M$. These permutations satisfy $RA=A^qR$.

\begin{proposition}[Coset leaders as constrained minimal rotations]
\label{v139:prop-digit-leaders}
With the notation of Definition~\ref{v50:def-coset-leader},
\begin{equation}\label{v139:eq-digit-leader}
 \lambda_m(a)=\frac1c\min_{0\le j<L}R^j(E_a).
\end{equation}
In particular, $a$ is a coset leader precisely when the digits of $E_a$, read from high to low, are lexicographically no larger than any cyclic rotation. If $\Delta(m,q)$ denotes the largest nonzero coset leader, then
\[
 \Delta(m,q)=\frac1c
 \max_{\substack{1\le E<M\\c\mid E}}\min_{0\le j<L}R^j(E).
\]
Moreover,
\begin{equation}\label{v139:eq-largest-leader}
 \Delta(m,q)\le
 \left\lfloor\frac{m((q-1)q^{L-1}-1)}{q^L-1}\right\rfloor,
 \qquad
 \Delta(q^L-1,q)=(q-1)q^{L-1}-1
\end{equation}
whenever the modulus in question is at least two.
\end{proposition}
\begin{proof}
The integers $R^j(E_a)$ are exactly $c$ times the least nonnegative representatives of the coset, proving~\eqref{v139:eq-digit-leader}. Comparing finite base-$q$ integers means comparing their digits from high to low. Every allowed nonzero block has a digit at most $q-2$, since the all-$(q-1)$ block represents $M$ and is excluded. Move such a digit to the highest position. The resulting integer is at most $(q-1)q^{L-1}-1$, giving the bound. When $c=1$, the high-to-low block $(q-2)(q-1)\cdots(q-1)$ attains it: its other rotations are no smaller.
\end{proof}
The order here is the order of the finite integers $E_a$, not an order on $\mathbb Z_p$. For a general modulus, divisibility by $c$ restricts the allowable words. To prove that a candidate is largest, it suffices to construct its minimal rotation and show that every allowable word has a rotation no larger than it.

Reduction and relative phases can also be read through the same rational encodings. Set
\[
 \mathcal D_m=\{-a/m+\mathbb Z:a\in\mathbb Z/m\mathbb Z\}
 \subseteq\mathbb Z_p/\mathbb Z.
\]
Only the additive quotient group is used. If $h\mid m$, reduction of exponents modulo $h$ corresponds to
\begin{equation}\label{v139:eq-digit-reduction}
 \mathcal D_m\longrightarrow\mathcal D_h,
 \qquad x\longmapsto(m/h)x.
\end{equation}
For $h\mid m\mid n$, these maps compose by multiplication of their integer factors, and all commute with multiplication by $q$.

\begin{proposition}[Local words and synchronized phases]
\label{v139:prop-digit-crt}
Let $m=\prod_{i=1}^s m_i$, where the $m_i>1$ are pairwise coprime. Put $n_i=m/m_i$ and choose $u_i$ with $u_in_i\equiv1\pmod{m_i}$. The maps
\[
 x\longmapsto(x_i)_i=(n_ix)_i,
 \qquad (x_i)_i\longmapsto\sum_i u_ix_i
\]
are inverse isomorphisms between $\mathcal D_m$ and $\prod_i\mathcal D_{m_i}$, commuting with multiplication by $q$.
Fix a local orbit of length $d_i$ in each factor. The global orbits lying over these local orbits have length $\lcm_i d_i$ and are parametrized by
\begin{equation}\label{v139:eq-digit-phase}
 \left(\prod_i\mathbb Z/d_i\mathbb Z\right)
 \big/\langle(1,\ldots,1)\rangle.
\end{equation}
Thus local cyclic word classes require relative phases to determine a global class.
\end{proposition}
\begin{proof}
For $x=-a/m+\mathbb Z$, its $i$th image is $-a/m_i+\mathbb Z$. The inverse formula is the Chinese remainder theorem, since $\sum_i u_in_i\equiv1\pmod m$. Multiplication by $q$ rotates every local word once. Its action on the phase tuple is addition of $(1,\ldots,1)$, whose order is $\lcm_i d_i$; the quotient gives exactly the global orbits, as in Theorem~\ref{v44:thm-local-orbit-synchronization}.
\end{proof}
Along a divisibility chain, compatible choices in~\eqref{v139:eq-digit-reduction} recover the compatible orbit systems of Chapter~5. Projection or multiplication may involve carries; it is not generally deletion of selected digit positions.

\begin{example}[Two phases modulo fifteen]\label{v139:ex-digit-phase}
For $q=2,m=15$, the nonzero orbit modulo three has length two, and the unit orbit modulo five has length four. The phase quotient has $2\cdot4/4=2$ elements. The two global unit orbits are
\[
 C_2(1)=\{1,2,4,8\},\qquad C_2(7)=\{7,14,13,11\}.
\]
Their low-to-high blocks are the rotations of $(1,0,0,0)$ and $(1,1,1,0)$, respectively. They have the same local orbit classes: the pairs of representatives $(1,1)$ and $(1,2)$ modulo $(3,5)$ differ by a relative phase. They do not give the same global word class.
\end{example}

\subsection{Carries, complete blocks, and binomials}
\label{v139:sec-digit-binomial}
Choose a primitive $m$th root $\zeta_m\in\mathcal O_{fL}$. There is an isomorphism of groups
\begin{equation}\label{v139:eq-digit-root-group}
 \Theta_m:\mathcal D_m\longrightarrow\mu_m,
 \qquad -a/m+\mathbb Z\longmapsto\zeta_m^a.
\end{equation}
Addition of encodings corresponds to multiplication of roots; multiplication by $q$ corresponds to $\phi^f$. This is not an identification of the rational number $-a/m$ with a root of unity. For a single orbit $C$, retain the notation $\widetilde M_C(X)=\prod_{x\in C}(X-\zeta_m^x)$ from Chapter~5.

\begin{proposition}[A complete fiber and its binomial]
\label{v139:prop-digit-binomial}
Let $r\mid m$, $h=m/r$, and $B\subseteq\mathbb Z/m\mathbb Z$ have $r$ elements. The following conditions are equivalent:
\begin{enumerate}
\item $B=b+h\mathbb Z/m\mathbb Z$ for some $b$;
\item the $r$ encodings $-x/m+\mathbb Z$, $x\in B$, have the same image under multiplication by $r$;
\item the canonical-root polynomial is
\begin{equation}\label{v139:eq-lifted-block-binomial}
 \prod_{x\in B}(X-\zeta_m^x)=X^r-\zeta_m^{br}
\end{equation}
for some $b$.
\end{enumerate}
The identity holds in $\mathcal O_{fL}[X]$ and at every quotient precision $p^j$, $j\ge1$; its residue-field form is Proposition~\ref{prop:ed-binomial-correspondence}.
\end{proposition}
\begin{proof}
Equality after multiplication by $r$ means equality of exponents modulo $h$. Each such fiber has exactly $r$ elements, proving the first equivalence. Under~\eqref{v139:eq-digit-root-group}, the fiber becomes $\zeta_m^b\mu_r$, which is precisely the root set of $X^r-\zeta_m^{br}$. Conversely, every root of that binomial has the same $r$th power; the roots are distinct because $p\nmid r$. Reducing an exact integral identity preserves it at every precision.
\end{proof}

For $x'=[x+h]_m$, actual rational representatives satisfy
\[
 -\frac{x'}m=-\frac xm-\frac1r+\left\lfloor\frac{x+h}{m}\right\rfloor.
\]
The integer correction keeps the representative in $(-1,0]$. A finite set in $\mathbb Z_p$ cannot be invariant under a nonzero translation; the complete-block symmetry belongs to $\mathcal D_m$, or to this normalized operation with carries.

\begin{theorem}[Exact period compression and multiple equal-difference decompositions]
\label{v139:thm-period-compression}
Let $C=C_q(a)$, $1\le a<m$, $d=|C|$, and $\kappa=|H_+(C)|$. For $r\mid m$, put $h=m/r$ and
\[
 \beta_a^{(r)}=-\frac{[a]_h}{h}
 =r\left(-\frac am\right)+\left\lfloor\frac ah\right\rfloor,
 \qquad d_r=\operatorname{per}_q(\beta_a^{(r)}),
\]
where $\operatorname{per}_q$ is the least digit period and $\operatorname{per}_q(0)=1$. Call $d=rd_r$ exact period compression. Then $d_r\mid d$, $d/d_r\le r$, and the following are equivalent:
\begin{enumerate}
\item $r\mid\kappa$, equivalently $C+h=C$;
\item $C$ is a union of complete equal-difference blocks of size $r$;
\item $d=r d_r$;
\item $\widetilde M_C(X)=G_r(X^r)$, where
\begin{equation}\label{v139:eq-compression-factor}
 G_r(Y)=\prod_{j=0}^{d_r-1}(Y-\zeta_m^{ra q^j})
\end{equation}
is the canonical irreducible orbit factor over $K_f$ for the quotient orbit modulo $h$.
\end{enumerate}
When these conditions hold,
\begin{equation}\label{v139:eq-compression-leader}
 \lambda_m(a)=\lambda_h([a]_h),\qquad
 \{R^{j d_r}(E_a):0\le j<r\}
 =E_a+\frac Mr\mathbb Z/M\mathbb Z.
\end{equation}
The largest such $r$ is $\kappa$; it is also the largest monomial composition exponent at every precision $\mathcal O_f/p^j\mathcal O_f$. In particular, the least period after multiplication by $\kappa$ and normalization is $\omega_C=d/\kappa$.
\end{theorem}
\begin{proof}
Reduction modulo $h$ maps $C$ onto a single orbit of length $d_r$. Multiplication by $q$ bijects its fibers, so each has $d/d_r$ elements. The ambient fibers have $r$ elements. Equality $d/d_r=r$ holds exactly when every nonempty fiber is complete, proving the first three equivalences. Multiplying the binomials from Proposition~\ref{v139:prop-digit-binomial} over these fibers gives~\eqref{v139:eq-compression-factor}; their $r$th powers form a single orbit of length $d_r$, so $G_r$ is irreducible over $K_f$.

Conversely, a monomial composition makes the root set invariant under multiplication by $\mu_r$, hence the exponent set invariant under addition of $m/r$. Each full fiber contains its least representative modulo $h$, proving the leader identity. The permutation $q^{d_r}$ traverses a fiber in $r$ steps, giving the second identity in~\eqref{v139:eq-compression-leader}. It describes the same set, not necessarily its order under successive additions of $M/r$.

Theorem~\ref{thm:med-coset-partitions} identifies all permissible $r$ as the divisors of $\kappa$. Since the reduction of $\widetilde M_C$ is separable and has nonzero constant term, every monomial composition exponent of that reduction is prime to $p$; its root ratios force it to divide $m$. Thus an additional composition at a finite precision would give a larger allowed exponent after reduction, which is impossible. The exact composition for $\kappa$ already holds in $\mathcal O_f$.
\end{proof}

For $a=1$, the residue-field complete-fiber criterion is Proposition~\ref{v142:prop-full-fibers}; here it also records the digit period and the canonical factor at every precision.

Consequently,
\[
 \kappa_C=\max\{r\mid m:d=r\operatorname{per}_q(\beta_a^{(r)})\}.
\]
The arithmetic formula of Theorem~\ref{thm:med-minimal-step} gives the same invariant: for $n=m/\gcd(m,a)$,
\[
 \omega_C=\ord_{b(n)}(q),\qquad \kappa_C=\ord_n(q)/\omega_C,
\]
with $b(n)$ as defined in Chapter~5. Over $K_{fs}$ all irreducible factors of $\widetilde M_C$ are binomials exactly when $\omega_C\mid s$: the $q^s$-suborbits must be complete fibers. Thus the minimum extension degree for binomial factorization is the compressed digit period. The residue-field and finite-precision versions use the same orbit partition.

\begin{example}[Modulus nine: carries, compression, and a binomial]
\label{v139:ex-nine-compression}
For $q=2,m=9$, take $C=C_2(1)$, $L=6$, $M=63$, $c=7$. Its complete blocks and low-to-high digits are
\[
\begin{array}{c|c|c|c}
 x\in B_1&E_x\text{ in binary}&x\in B_2&E_x\text{ in binary}\\\hline
 1&111000&2&011100\\
 4&001110&5&110001\\
 7&100011&8&000111
\end{array}
\]
Here $\kappa_C=3$, $\omega_C=2$, and rotation by two positions traverses each block. In the following notation a parenthesized word is repeated in low-to-high order:
\[
 -1/9=(111000)_2^\infty,\qquad
 3(-1/9)=-1/3=(10)_2^\infty.
\]
Multiplication with carries reduces the least period from six to two, while
\[
 \widetilde M_C(X)=X^6+X^3+1=G(X^3),\qquad G(Y)=Y^2+Y+1.
\]
The other nonzero binary coset is $\{3,6\}$, with least period two and trivial additive stabilizer. Its six-digit encoding is $(101010)$, and the two nonzero leaders are $1,3$. Thus $\Delta(9,2)=3$, strictly below the bound $\lfloor31/7\rfloor=4$ in~\eqref{v139:eq-largest-leader}.

Keeping $m=9$ but taking $q=4$ gives the single complete block $C_4(1)=\{1,4,7\}$. Now
\[
 -1/9=(310)_4^\infty,\qquad 3(-1/9)=(1)_4^\infty,
 \qquad F_{C_4(1)}(X)=X^3-\zeta_9^3.
\]
The coefficient $\zeta_9^3$ belongs to $\mathcal O_2$, and its reduction lies in $\mathbb F_4$. The period has compressed to one, as a binomial requires.
\end{example}

\begin{proposition}[Digit sums and orbit sums]\label{v139:prop-orbit-digit-sum}
Using the shortest block $E=a(q^d-1)/m=\sum_{i=0}^{d-1}e_iq^i$, one has
\begin{equation}\label{v139:eq-orbit-digit-sum}
 \sum_{b\in C_q(a)}b=\frac{m}{q-1}\sum_{i=0}^{d-1}e_i,
 \qquad
 \lambda_m(a)\le\frac{m}{d(q-1)}\sum_i e_i,
\end{equation}
where $b$ denotes its least nonnegative representative. If $q=p$, the Gauss sum over $\mathbb F_{p^d}$ with the conventions of Corollary~\ref{v138:cor-word-gauss} satisfies
\[
 v_p\bigl(G(\omega^{-E},\psi)\bigr)=\frac1m\sum_{b\in C_p(a)}b.
\]
\end{proposition}
\begin{proof}
In the sum of all $d$ rotations, each digit visits each position once. Their sum is $(q^d-1)\sum_i e_i/(q-1)$; it is also $(q^d-1)\sum_{b\in C_q(a)}b/m$. The inequality compares the minimum with the average. For $q=p$, apply the digit-sum valuation formula. For $q=p^f$ with $f>1$, the orbit-sum identity uses base-$q$ digit values, whereas the valuation uses their base-$p$ digit sums.
\end{proof}

\subsection{Digit coordinates and BCH bounds}
\label{v139:sec-digit-bch}
Fix the primitive $m$th root $\bar\zeta_m$ in the residue field. Let $Z\subsetneq\mathbb Z/m\mathbb Z$ be $q$-stable and let $\mathcal C_Z$ be the nonzero cyclic code over $\mathbb F_q$ with defining set $Z$. The required consecutive objects in the BCH bound are exponents in $Z$, not coset leaders. Put
\[
 S_Z=\{ca:a\in Z\}\subseteq\mathbb Z/M\mathbb Z.
\]
In these formulas the subscript of $E_b$ is reduced modulo $m$. Frobenius closure is $R(S_Z)=S_Z$. A consecutive interval $b,b+1,\ldots,b+t-1$ is contained in $Z$ precisely when
\begin{equation}\label{v139:eq-bch-carry-run}
 E_b,A(E_b),\ldots,A^{t-1}(E_b)\in S_Z.
\end{equation}
Thus rotation describes closure, while consecutive roots are detected by repeated addition with carries.

\begin{proposition}[The longest covered interval and the omitted gaps]
\label{v139:prop-digit-bch-gap}
Let $\ell(Z)$ be the greatest length of a cyclic consecutive interval contained in $Z$, with $\ell(\varnothing)=0$. List the omitted exponents as
\[
 0\le t_1<\cdots<t_k<m,\qquad t_{k+1}=t_1+m.
\]
Then the largest bound obtained by varying the starting exponent in the BCH argument for this fixed root is
\begin{equation}\label{v139:eq-bch-gap}
 B(Z)=\ell(Z)+1=\max_i(t_{i+1}-t_i)
     =\frac1c\max_i(E_{t_{i+1}}-E_{t_i}),
 \qquad d_H(\mathcal C_Z)\ge B(Z).
\end{equation}
In the last expression set $E_{t_{k+1}}=E_{t_1}+M$.
\end{proposition}
\begin{proof}
The interval strictly between two consecutive omitted exponents has length $t_{i+1}-t_i-1$. Every covered interval is contained in one of these gaps. Apply Theorem~\ref{v56:thm-bch-bound}; if $Z$ is empty, the inequality is simply $d_H\ge1$. Scaling by $c$ gives the last equality.
\end{proof}
This determines the strongest consecutive-root bound for the chosen exponent coordinates, not necessarily the minimum distance. If the primitive root is replaced by $\bar\zeta_m^u$, $u\in(\mathbb Z/m\mathbb Z)^\times$, the defining set becomes $u^{-1}Z$. Maximizing $B(u^{-1}Z)$ also allows all coprime step sizes; in digit coordinates the step is $A^u$.

For a narrow-sense BCH code with design parameter $2\le\delta\le m$, write $Z_\delta=\bigcup_{i=1}^{\delta-1}C_q(i)$. Proposition~\ref{v139:prop-digit-leaders} gives
\begin{equation}\label{v139:eq-bch-threshold}
 a\in Z_\delta\ \Longleftrightarrow\
 0<\min_jR^j(E_a)<c\delta,
 \qquad
 \dim\mathcal C_{Z_\delta}
 =m-\sum_{\substack{\lambda\in\Gamma_m(q)\\1\le\lambda<\delta}}|C_q(\lambda)|.
\end{equation}
The least omitted positive exponent is the least coset leader at least $\delta$, or $m$ if there is none. Hence the defining set stays unchanged as $\delta$ increases up to that next leader. Once $\delta>\Delta(m,q)$, the code is the repetition code. These conclusions follow from the same rotation threshold; no ordering of whole cosets is required.

\begin{example}[The length-fifteen decoding example in digit coordinates]
\label{v139:ex-digit-bch-fifteen}
For $q=2,m=15$, $L=4$ and $c=1$. The nonzero leaders, high-to-low minimal words, and orbit lengths are
\[
\begin{array}{c|c|c}
 \lambda&\text{minimal word}&|C_2(\lambda)|\\\hline
 1&0001&4\\3&0011&4\\5&0101&2\\7&0111&4
\end{array}
\]
For $Z=C_2(1)\cup C_2(3)$, the omitted exponents are $0,5,7,10,11,13,14$. Their largest cyclic gap is five, so $B(Z)=5$ and the dimension is $15-4-4=7$. Design parameters four and five give this same defining set. Example~\ref{v118:ex-bch-fifteen-decoding} gives a weight-five generator and completes its two-error decoding; the digit description determines both the root interval and the stability range of the design parameter.
\end{example}

\begin{corollary}[Complete fibers and interleaved copies of a shorter code]
\label{v139:cor-compressed-code}
Let $r\mid m$, $h=m/r$, and $Z=\pi^{-1}(Z_0)$ for a proper $q$-stable set $Z_0\subsetneq\mathbb Z/h\mathbb Z$, where $\pi$ is reduction modulo $h$. If $g_0(Y)$ generates the length-$h$ code $\mathcal C_{Z_0}$, then
\[
 g_Z(X)=g_0(X^r),\qquad
 \mathcal C_Z\simeq\mathcal C_{Z_0}^{\oplus r}
 \quad\text{by a coordinate permutation}.
\]
Consequently,
\begin{equation}\label{v139:eq-compressed-code}
 \dim\mathcal C_Z=r\dim\mathcal C_{Z_0},\qquad
 d_H(\mathcal C_Z)=d_H(\mathcal C_{Z_0}),\qquad B(Z)=B(Z_0).
\end{equation}
\end{corollary}
\begin{proof}
The complete-root fibers give the polynomial identity. Every residue modulo $X^{rh}-1$ has a unique expression
\[
 \sum_{i=0}^{r-1}X^iu_i(X^r),\qquad \deg u_i<h.
\]
 It belongs to the ideal $(g_0(X^r))$ exactly when every $u_i$ belongs to $(g_0(Y))$ modulo $Y^h-1$. Grouping coordinates by their residue modulo $r$ gives the direct sum, hence the dimension and distance. The indicator of $Z$ repeats with period $h$ and has an omitted position in each period, so its longest covered cyclic interval has the same length as for $Z_0$.
\end{proof}
For the binary modulus-nine orbit in Example~\ref{v139:ex-nine-compression}, $g_Z=X^6+X^3+1$ gives a $[9,3,3]$ code: it consists of three interleaved copies of the $[3,1,3]$ repetition code. The digit period compresses from six to two, while the longest covered exponent interval has length two at both moduli. Which other carry constraints can force long covered intervals without making the defining set a union of complete fibers?

Digit sums determine these Gauss-sum valuations; canonical-factor coefficients are sums of many products of roots of unity. Which joint distributions of periodic words can also detect local cancellation in those sums? The exterior-power traces and integer subset sums below retain combinatorial relations absent from a single digit sum. For curve eigenvalues already expressed through Gauss or Jacobi sums, equation~\eqref{v138:eq-normalized-digit-weight} connects with the normalized slopes in Section~\ref{revadic:sec-newton}.

\subsection{Canonical factors and their valuation profiles}
Let $m\ge2$ be prime to $p$ and $C$ a $p^d$-cyclotomic coset modulo $m$. In a suitable unramified extension containing the $m$th roots of unity, put
\[
 \widetilde M_C(X)=\prod_{i\in C}(X-\zeta_m^i)
 =X^r+c_1X^{r-1}+\cdots+c_r.
\]
Its coefficients lie in the integer ring fixed by the relevant $\phi^d$. Define the \emph{valuation profile}
\begin{equation}\label{revadic:eq-profile}
 \mathcal V_p(C)=\bigl(v_p(c_1),\ldots,v_p(c_r)\bigr),
 \qquad v_p(0)=+\infty.
\end{equation}
Valuation zero indicates a nonzero coefficient after reduction. A positive finite valuation indicates a coefficient that vanishes modulo $p$ but not in characteristic zero. Infinite valuation signifies an exact zero. The ordinary reduced factor distinguishes only the first case from the latter two.

\begin{proposition}[Exact multiplicative symmetry forces exact coefficient gaps]
\label{revadic:prop-exact-gaps}
Suppose a root of unity $u$ of order $e$ stabilizes the root set multiplicatively:
$u\cdot\{\zeta_m^i:i\in C\}=\{\zeta_m^i:i\in C\}$. Then $e\mid r$ and
\[
 \widetilde M_C(X)\in\mathcal O_d[X^e].
\]
Thus $c_j=0$ for every $e\nmid j$, and $v_p(c_j)=+\infty$. In contrast, $\bar c_j=0$ alone does not imply $c_j=0$.
\end{proposition}
\begin{proof}
Multiplication by $u$ acts freely on the nonzero roots, each orbit having size $e$, so $e\mid r$. Invariance of the root set gives
\[
 \widetilde M_C(uX)=u^r\widetilde M_C(X)=\widetilde M_C(X).
\]
Comparison of the coefficients of $X^{r-j}$ yields $(u^{-j}-1)c_j=0$. In a field of characteristic zero, $u^j\ne1$ forces $c_j=0$.
\end{proof}

\subsection{Lifting Lucas Atoms and Root-Difference Valuations}
\label{v146:sec-adic-lucas}
In Chapter 6, a Lucas atom in characteristic two is the square of the optimal coefficient polynomial. Lifting to an unramified integer ring separates the two branches over each repeated root and determines their exact collision depth.

The symmetric decomposition in Lemma~\ref{v103:lem-symmetric-two} uses only recurrences with integer coefficients. Thus the homogeneous cyclotomic identity defines the same \(\Lambda_n\) in \(\Z[B,C]\), before reduction to the residue field.

\begin{proposition}[Root differences and discriminant valuations of Lucas atoms]\label{v146:prop-adic-lucas}
Let \(q=p^f\), let \(\bar c\in\F_q^\times\) be primitive, and put \(c=[\bar c]\in\mathcal O_f^\times\). Regard \(\Lambda_{q+1}(B,C)\) as the polynomial in \(\Z[B,C]\) defined by the homogeneous cyclotomic identity, and put
\[
 \mathcal L(B)=\Lambda_{q+1}(B,c),\qquad d=\varphi(q+1)/2.
\]
Then \(\mathcal L\) has \(2d\) distinct roots in \(\mathcal O_f\), all of which are units. If \(p\ne2\), they remain distinct modulo \(p\). If \(p=2\), they form \(d\) pairs \(\{\beta,-\beta\}\), with
\begin{equation}\label{v146:eq-lucas-root-depth}
 v_2(\beta-(-\beta))=1,
\end{equation}
whereas differences between roots from different pairs are units. With \(v_p(p)=1\),
\begin{equation}\label{v146:eq-lucas-discriminant}
 v_p(\operatorname{disc}\mathcal L)
 =\begin{cases}0,&p\ne2,\\\varphi(q+1),&p=2.\end{cases}
\end{equation}
\end{proposition}
\begin{proof}
For each primitive quadratic \(X^2+bX+\bar c\), choose a root \(\theta\in\F_{q^2}\) and put
\[
 \beta=-([\theta]+[\theta^q]).
\]
Frobenius exchanges the two Teichmüller representatives, so \(\beta\in\mathcal O_f\). Their product is \([\bar c]=c\), and their ratio still has order \(q+1\); hence the homogeneous cyclotomic identity gives \(\mathcal L(\beta)=0\). Also \(\bar\beta=b\ne0\), since \(b=0\) would give root-ratio order at most two. Here \(\beta\) is the negative sum of the lifted roots and generally differs from \([b]\).

For \(p\ne2\), the fixed-constant-term count of Chapter 6 gives \(2d\) distinct values of \(b\). Thus the corresponding \(\beta\) already give all roots, with unit differences. For \(p=2\), the count gives \(d\) distinct values of \(b\). Since \(\mathcal L\) has only even powers of \(B\), each \(-\beta\) is also a root. The element \(\beta\) is a unit, so \(\beta\ne-\beta\) and \(v_2(2\beta)=1\). Distinct values of \(b\) make differences between different pairs units. These are all \(2d\) roots.

The discriminant of a monic polynomial is the product of squared differences over all unordered pairs of distinct roots. In odd characteristic every contribution has valuation zero; in characteristic two exactly \(d\) pairs contribute two each. The total is \(2d=\varphi(q+1)\).
\end{proof}

For example, when \(q=4\) and \(\bar c=\omega\), the polynomial \(\mathcal L=B^4-3cB^2+c^2\) reduces modulo two to \(P_{4,\omega}^2\) from Example~\ref{v146:ex-lucas-four}. It has four distinct roots in characteristic zero and discriminant \(400c^6\), of binary valuation four. Each pair coalesces modulo two and separates modulo four: the multiplicity records exactly one level of root-difference precision.

\subsection{Cyclotomic root blocks as relative norms and exterior-power traces}
\label{revadic:subsec-orbit-norm}

For a single cyclotomic coset, the canonical factor is not merely a product of roots: it is a \emph{relative norm polynomial}. Thus every coefficient becomes an explicit linear-algebraic quantity over an integer ring, without first computing the coefficients of a lifted factor.

For a free module of rank $r$, its $j$th exterior power has basis $e_{i_1}\wedge\cdots\wedge e_{i_j}$, with $i_1<\cdots<i_j$. The entries of $\bigwedge^jM$ are the corresponding $j$th-order minors of $M$; its trace is the sum of all principal minors of that order. Thus the coefficient of $X^{r-j}$ in $\det(XI-M)$ is $(-1)^j\Tr(\bigwedge^jM)$, an identity compatible with reduction of coefficients.

\begin{theorem}[Orbit factors, relative norms, and exterior-power traces]
\label{revadic:thm-orbit-exterior}
Let $m\ge2$ and $p\nmid m$, put $h=\ord_m(p)$, and choose a primitive $m$th root $\zeta\in K_h$. Given $d\ge1$ and an exponent $i\bmod m$, write
\[
 C=\{i,p^di,\ldots,p^{d(r-1)}i\}\pmod m,
 \qquad r=|C|,
\]
and inside $K_{\operatorname{lcm}(d,h)}$ put $L=K_d(\zeta^i)$ and $B=\mathcal O_L$. Then $L/K_d$ is unramified of degree $r$. Let $M_{\zeta^i}$ be the $\mathcal O_d$-linear operator on $B$ given by multiplication by $\zeta^i$. Write $Q_C=\widetilde M_C$ and $c_0=1$. Then
\begin{align}
 Q_C(X)&=\operatorname{N}_{L/K_d}(X-\zeta^i)
       =\det(XI-M_{\zeta^i}),\label{revadic:eq-orbit-norm}\\
 (-1)^jc_j&=\operatorname{Tr}_{\mathcal O_d}
       \!\bigl(\textstyle\bigwedge^j M_{\zeta^i}\bigr)
           \quad(0\le j\le r,\ c_0=1).\label{revadic:eq-exterior-trace}
\end{align}
In particular, $-c_1=\operatorname{Tr}_{L/K_d}(\zeta^i)$ and $(-1)^rc_r=\operatorname{N}_{L/K_d}(\zeta^i)$. Reducing~\eqref{revadic:eq-exterior-trace} modulo $p^s$ gives the exact exterior-power trace congruence in the associated Galois ring; modulo $p$ it is the exterior-power trace of a finite-field multiplication operator.
\end{theorem}
\begin{proof}
Inside $K_{\operatorname{lcm}(d,h)}$, the group fixing $K_d$ is generated by $\phi^d$. Hence the distinct conjugates of $\zeta^i$ over $K_d$ are precisely $\zeta^{p^{dj}i}$, $0\le j<r$. Its minimal polynomial and its relative norm polynomial are therefore both $Q_C$, and $[L:K_d]=r$. Unramifiedness is inherited from the unramified tower. Choose an integral $\mathcal O_d$-basis of $B$ and represent multiplication by $\zeta^i$ by an integral matrix. Over a splitting field in characteristic zero its eigenvalues are the $r$ conjugates just listed. The eigenvalues of its $j$th exterior power are their products over all $j$-subsets. Its trace is the $j$th elementary symmetric polynomial, proving~\eqref{revadic:eq-exterior-trace}. Both the original and the exterior-power matrices have integral entries and may be reduced modulo $p^s$. The reduced integral basis remains a basis, and the reduced operator is exactly multiplication by the residue class in the finite field.
\end{proof}

The formula also marks a precise limit: the Frobenius orbit specifies \emph{which} exterior-power trace must be computed, but the length of that orbit alone cannot determine its $p$-adic valuation. Exact zero, low-level divisibility, and higher-level cancellation correspond respectively to a zero exterior-power trace, its vanishing modulo $p^s$, and cancellation at greater precision. The subset-count formula below decomposes the same trace into orbit-period sums in the exponent group.

\subsection{An exact coefficient formula from subset sums and Gauss periods}
\label{revadic:subsec-orbit-period}

Computing successive Witt digits tests individual coefficients, but does not always display their relation to the underlying orbit. A complementary description first counts exponent sums modulo $m$ and then uses orbit-period sums (generalized Gauss periods) to return these counts to the $p$-adic integer ring. It works equally well for unit and nonunit exponent orbits.

Let $p\nmid m$, let $\zeta\in\mu_m$ be primitive, and let $C\subseteq\mathbb Z/m\mathbb Z$ be stable under multiplication by $p^d$ (in particular, $C$ may be a single $p^d$-cyclotomic coset). Set $r=|C|$ and
\[
 Q_C(X)=\prod_{i\in C}(X-\zeta^i)
       =\sum_{j=0}^r c_jX^{r-j},\qquad c_0=1.
\]
For $0\le j\le r$ and $s\bmod m$ define the integer
\begin{equation}\label{revadic:eq-subsets}
 N_{C,j}(s)=\#\left\{I\subseteq C:|I|=j,\quad
                         \sum_{i\in I}i\equiv s\pmod m\right\}.
\end{equation}
Let $\mathscr O_d$ denote the orbits of multiplication by $p^d$ on $\mathbb Z/m\mathbb Z$. For each $D\in\mathscr O_d$, put
$\eta_D=\sum_{s\in D}\zeta^s$ (including $D=\{0\}$, with $\eta_D=1$).

\begin{theorem}[The subset-count--Gauss-period formula]
\label{revadic:thm-period-profile}
For each $j$, the function $N_{C,j}(s)$ is constant on every $D\in\mathscr O_d$; denote the common value by $N_{C,j}(D)$. Then
\begin{equation}\label{revadic:eq-period-profile}
 \boxed{\displaystyle
 c_j=(-1)^j\sum_{D\in\mathscr O_d}N_{C,j}(D)\eta_D.}
\end{equation}
Every $\eta_D$ lies in $\mathcal O_d$, and so does $c_j$. Consequently, for every $s\ge1$,
\[
 v_p(c_j)\ge s\ \Longleftrightarrow\
 \sum_{D\in\mathscr O_d}(-1)^jN_{C,j}(D)\eta_D\equiv0\pmod{p^s}.
\]
If the minimum valuation among the nonzero terms $N_{C,j}(D)\eta_D$ is attained at a unique orbit, then
\begin{equation}\label{revadic:eq-period-min}
 v_p(c_j)=\min_{D\in\mathscr O_d}
 \bigl(v_p(N_{C,j}(D))+v_p(\eta_D)\bigr).
\end{equation}
If the minimum is attained more than once, it is only a lower bound; the exact valuation depends on cancellation of the terms of lowest valuation.
\end{theorem}
\begin{proof}
Multiplication by $p^d$ permutes $C$ and maps a $j$-subset with exponent sum $s$ to one with sum $p^ds$. Thus $N_{C,j}(p^ds)=N_{C,j}(s)$. Expanding the elementary symmetric function gives
\[
 (-1)^jc_j=\sum_{\substack{I\subseteq C\\|I|=j}}
 \zeta^{\sum_{i\in I}i}=\sum_{s\bmod m}N_{C,j}(s)\zeta^s.
\]
Grouping by orbits proves \eqref{revadic:eq-period-profile}. Since $\phi^d$ fixes $\eta_D$, it lies in the corresponding fixed integer ring. The congruence criterion follows from the equality. The nonarchimedean triangle inequality yields the valuation bound, with equality when the minimum is unique.
\end{proof}

\subsection{Integral coordinates for exterior-power traces: from computation to criteria}
\label{revadic:subsec-integral-resolvent}

Equation~\eqref{revadic:eq-exterior-trace} identifies a coefficient with the exterior-power trace of an integral operator, while~\eqref{revadic:eq-period-profile} writes the same trace as a sum of periods weighted by exponent-subset counts. These descriptions meet in an elementary integer polynomial. More importantly, it separates \emph{exact vanishing in characteristic zero} from \emph{high divisibility at a chosen prime} by two distinct, verifiable conditions.

Retain the notation of the preceding subsection, with $m\ge2$, $p\nmid m$, $\zeta=\zeta_m$, and $C$ stable under multiplication by $p^d$. For $0\le j\le|C|$ define
\begin{equation}\label{revadic:eq-count-poly}
 A_{C,j}(T)=\sum_{s=0}^{m-1}N_{C,j}(s)T^s\in\mathbb Z[T],
 \qquad
 A_{C,j}(\zeta)=(-1)^jc_j
 =\operatorname{Tr}\!\bigl(\textstyle\bigwedge^jM\bigr).
\end{equation}
The final equality applies when $C$ is one $p^d$-orbit, by Theorem~\ref{revadic:thm-orbit-exterior}. Thus $A_{C,j}$ is an integral representative of the exterior-power trace in the group algebra $\mathbb Z[\mathbb Z/m\mathbb Z]$. Ordinary polynomial division can now be performed before computing any $p$-adic roots of unity.

\begin{theorem}[Cyclotomic remainders: exact zeros and valuation bounds]
\label{revadic:thm-integral-resolvent}
Divide by the monic $m$th cyclotomic polynomial in $\mathbb Z[T]$:
\[
 A_{C,j}(T)=H_{C,j}(T)\Phi_m(T)+B_{C,j}(T),\qquad
 B_{C,j}(T)=\sum_{t=0}^{\varphi(m)-1}b_tT^t,
\]
where $\deg B_{C,j}<\varphi(m)$, allowing the zero remainder. Put
\[
 \nu_{C,j}=\min_{0\le t<\varphi(m)}v_p(b_t),\qquad v_p(0)=+\infty.
\]
Then the following statements hold.
\begin{enumerate}
\item $c_j=0$ if and only if $B_{C,j}=0$, equivalently $\Phi_m\mid A_{C,j}$ in $\mathbb Z[T]$. This criterion for an \emph{exact zero} does not depend on the choice of a prime above $p$.
\item One always has $v_p(c_j)\ge\nu_{C,j}$. If $B_{C,j}\ne0$, let $\bar\zeta$ be the reduction of $\zeta$ in $k_h$, where $h=\operatorname{ord}_m(p)$, and let $f(T)\in\mathbb F_p[T]$ be its minimal polynomial. Then
\begin{equation}\label{revadic:eq-leading-resolvent}
 v_p(c_j)=\nu_{C,j}
 \quad\Longleftrightarrow\quad
 f(T)\nmid\overline{p^{-\nu_{C,j}}B_{C,j}(T)}.
\end{equation}
The bar denotes coefficientwise reduction modulo $p$. If the polynomial on the right is divisible by $f$, then $v_p(c_j)>\nu_{C,j}$.
\item In particular, if $\operatorname{ord}_m(p)=\varphi(m)$, so that $\Phi_m$ is irreducible modulo $p$, then for every $j$ one has the exact formula
\begin{equation}\label{revadic:eq-sharp-resolvent}
 \boxed{\displaystyle v_p(c_j)=\min_t v_p(b_t).}
\end{equation}
In this case the entire valuation profile follows from integer subset counts and polynomial remainders modulo $\Phi_m$, \emph{without} an additional computation of period valuations.
\end{enumerate}
\end{theorem}
\begin{proof}
Equation~\eqref{revadic:eq-count-poly} gives $(-1)^jc_j=B_{C,j}(\zeta)$. The minimal polynomial of a primitive $m$th root of unity over $\mathbb Q$ is $\Phi_m$; hence $B_{C,j}(\zeta)=0$ exactly when $B_{C,j}=0$. This proves the first statement. Since $\zeta$ is a $p$-adic unit, the valuation of $B_{C,j}(\zeta)$ is at least the minimum valuation of its coefficients. If $\nu_{C,j}<\infty$, factor out $p^{\nu_{C,j}}$; the remaining sum reduces to zero if and only if $\bar\zeta$ is a root of $\overline{p^{-\nu_{C,j}}B_{C,j}}$, equivalently if $f$ divides that polynomial. This proves the second assertion. When $\Phi_m\bmod p$ is irreducible, $f=\bar\Phi_m$ has degree $\varphi(m)$, whereas the nonzero polynomial $\overline{p^{-\nu_{C,j}}B_{C,j}}$ has smaller degree. Divisibility by $f$ is then impossible, giving the claimed equality.
\end{proof}

The distinction is structural. \emph{All} rational cyclotomic conjugates determine exact vanishing, but a fixed $p$-adic embedding tests divisibility through just one irreducible factor $f$ of $\Phi_m\bmod p$. When the latter polynomial is reducible, a nonzero integral remainder may nevertheless vanish modulo $p$ in the selected local component. This is precisely the source of the cancellation in the seventh-root example.

\begin{theorem}[The minimum-valuation theorem across local components]
\label{revadic:thm-local-min}
Retain the notation of Theorem~\ref{revadic:thm-integral-resolvent}, with $B_{C,j}\ne0$ and $\nu_{C,j}<\infty$. Put $U=(\mathbb Z/m\mathbb Z)^\times$ and, for each $u\in U$, define
\[
 c_j^{(u)}=(-1)^jA_{C,j}(\zeta_m^u),
 \qquad \nu=\min_t v_p(b_t).
\]
Then
\begin{equation}\label{revadic:eq-local-min}
 \boxed{\displaystyle\min_{u\in U}v_p(c_j^{(u)})=\nu.}
\end{equation}
These valuations are constant on the $p$-cyclotomic cosets of $u$. Thus cancellation at a local factor may increase the valuation of certain components, but it \emph{cannot increase the valuations of every component simultaneously}.
\end{theorem}
\begin{proof}
Each $c_j^{(u)}$ equals $\pm B_{C,j}(\zeta_m^u)$, so its valuation is at least $\nu$. Put $D(T)=p^{-\nu}B_{C,j}(T)\in\mathbb Z_p[T]$. Its reduction modulo $p$ is nonzero and has degree less than $\varphi(m)$. Since $p\nmid m$, the reduction of $\Phi_m$ is squarefree, with roots precisely $\bar\zeta_m^u$ for $u\in U$. If $D$ vanished at all these roots, then $\bar\Phi_m\mid\bar D$, contradicting the degree bound. Hence some $u$ attains valuation $\nu$. Finally, $\phi(B(\zeta_m^u))=B(\zeta_m^{pu})$, and $\phi$ preserves valuations, proving constancy along $p$-cosets.
\end{proof}

\begin{corollary}[Multiplicative symmetry makes the bound sharp]
\label{revadic:cor-symmetry-sharp}
Assume the hypotheses of Theorem~\ref{revadic:thm-local-min}, in particular $B_{C,j}\ne0$.
Let $H=\{h\in U:hC=C\}$ be the multiplicative stabilizer of the exponent set. If
\[
 U=H\langle p\rangle,
\]
then $v_p(c_j^{(u)})=\nu_{C,j}$ for \emph{every} $u\in U$, in particular for the originally chosen $\zeta_m$. This condition is more flexible than $\operatorname{ord}_m(p)=\varphi(m)$: symmetries of the exponent set can act together with Frobenius to force equal local valuations.
\end{corollary}
\begin{proof}
For $h\in H$, multiplication by $h$ permutes the $j$-subsets of $C$, hence $A_{C,j}(\zeta_m^h)=A_{C,j}(\zeta_m)$. If $u=hp^t$, then $A_{C,j}(\zeta_m^u)=\phi^t(A_{C,j}(\zeta_m))$, preserving valuation. The hypothesis forces the valuations at all $u$ to agree, and Theorem~\ref{revadic:thm-local-min} gives their common value.
\end{proof}

\begin{proposition}[Affine exponent symmetries and complementary coefficients]
\label{v141:prop-affine-profile}
Let $C\subseteq\mathbb Z/m\mathbb Z$, $r=|C|$, and $p\nmid m$. Write $c_j^{(u)}$ for the $j$th coefficient of $\prod_{i\in C}(X-\zeta_m^{ui})$. Define
\[
 H_{\rm aff}(C)=\{a\in(\mathbb Z/m\mathbb Z)^\times:
                    aC+b=C\text{ for some }b\bmod m\}.
\]
This is a subgroup. If $aC+b=C$, then, for every unit $u$,
\begin{equation}\label{v141:eq-affine-profile}
 c_j^{(ua)}=\zeta_m^{-ubj}c_j^{(u)},\qquad
 c_{r-j}^{(u)}=c_r^{(u)}c_j^{(-u)}\quad(0\le j\le r).
\end{equation}
Consequently coefficient valuations are constant on cosets of
$H_{\rm aff}(C)\langle p\rangle$. If this group is the full unit group, every nonzero coefficient attains the integer-remainder minimum at every phase. If
$-1\in H_{\rm aff}(C)\langle p\rangle$, then the profile is symmetric:
\[
 v_p(c_{r-j}^{(u)})=v_p(c_j^{(u)})\quad(0\le j\le r).
\]
Here $c_0=1$ and the constant coefficient $c_r$ is a unit; the equalities include infinite valuations of zero coefficients.
\end{proposition}
\begin{proof}
Composition and inversion of the maps $i\mapsto ai+b$ prove the subgroup assertion. The first identity follows by scaling every root by $\zeta_m^{ub}$. Complementing a $j$-element subset of the roots gives the second identity, including its sign. Teichmüller roots are units and Frobenius preserves valuation, proving constancy on the stated cosets. Theorem~\ref{revadic:thm-local-min} gives the common minimum; the complementary-coefficient identity gives the last assertion.
\end{proof}

For $B_{C,j}\ne0$, this yields a natural invariant beyond the lowest integral valuation $\nu_{C,j}$. For each $p$-coset $D\subseteq U$, define its \emph{excess local cancellation depth}
\[
 \epsilon_{C,j}(D)=v_p\bigl(B_{C,j}(\zeta_m^u)\bigr)-\nu_{C,j}
 \quad(u\in D).
\]
It is a nonnegative integer, and it vanishes on at least one coset. The problem is therefore not to predict an unconstrained list of valuations, but a nonnegative profile restricted by Frobenius, multiplicative stabilizers, and a global minimum. If $B_{C,j}=0$, the coefficient vanishes identically in every component and no such $\epsilon_{C,j}$ is defined.

\begin{corollary}[Degree bound for the number of excess-cancellation orbits]
\label{revadic:cor-excess-degree}
Under the hypotheses of Theorem~\ref{revadic:thm-local-min}, put
\[
 D(T)=\overline{p^{-\nu_{C,j}}B_{C,j}(T)}\in\mathbb F_p[T],
 \qquad h=\operatorname{ord}_m(p).
\]
Let $t$ be the number of $p$-Frobenius orbits of primitive $m$th roots with $\epsilon_{C,j}>0$. Then
\begin{equation}\label{revadic:eq-excess-degree}
 \boxed{\displaystyle t\le\left\lfloor\frac{\deg D}{h}\right\rfloor
 \le\frac{\varphi(m)}h-1.}
\end{equation}
In particular, if $\deg D<h$, no local component has first-level excess cancellation.
\end{corollary}
\begin{proof}
Each orbit with excess cancellation corresponds to a distinct irreducible factor $f$ of degree $h$ dividing $D$. Their product divides $D$, so $ht\le\deg D<\varphi(m)$, as required.
\end{proof}

\begin{proposition}[A norm constraint on local cancellation depths]
\label{v137:prop-norm-depth}
Retain the assumptions and notation of Theorem~\ref{revadic:thm-local-min}. Suppose $B=B_{C,j}\ne0$, put $\nu=\nu_{C,j}$, and set
\[
 D_0(T)=p^{-\nu}B(T)\in\mathbb Z[T],\quad
 U=(\mathbb Z/m\mathbb Z)^\times,\quad
 H=\{u\in U:uC=C\},\quad G=H\langle p\rangle.
\]
Write $e=|G|$ and $r=[U:G]$. For $\gamma\in U/G$, define
\[
 \epsilon_\gamma=v_p(D_0(\zeta_m^u))\quad(u\in\gamma),
 \qquad R_j=\operatorname{Res}_T(\Phi_m,D_0).
\]
Then $R_j$ is a nonzero integer, $\epsilon_\gamma$ is independent of the representative, and
\begin{equation}\label{v137:eq-norm-depth}
 \min_\gamma\epsilon_\gamma=0,\qquad
 v_p(R_j)=e\sum_{\gamma\in U/G}\epsilon_\gamma.
\end{equation}
In particular, if $r=2$, the two depths are exactly
\[
 \left\{0,\frac{v_p(R_j)}e\right\}.
\]
If $v_p(R_j)>0$, the local factor test modulo $p$ identifies the positive entry.
In general, $\max_\gamma\epsilon_\gamma\le v_p(R_j)/e$; if $v_p(R_j)=e$, exactly one $G$-coset has depth one and all others have depth zero.
\end{proposition}
\begin{proof}
The group $H$ preserves subset-sum counts, and Frobenius preserves valuations, so the depth is constant on each $G$-coset. Theorem~\ref{revadic:thm-local-min} gives minimum zero.
Since $\Phi_m$ is irreducible over $\mathbb Q$, whereas $0\ne D_0$ and $\deg D_0<\varphi(m)$, these polynomials are coprime. Hence
\[
 0\ne R_j=\prod_{u\in U}D_0(\zeta_m^u)\in\mathbb Z.
\]
Taking valuations and grouping by $G$-cosets gives~\eqref{v137:eq-norm-depth}. The remaining claims follow because the depths are nonnegative integers.
\end{proof}

This is a \textbf{norm constraint} on local cancellation depths. When period relations already give the norm, it bounds the depths directly; when a high-degree resultant must still be computed, successive local factor tests may be simpler. For $r\ge3$, the sum and minimum usually do not determine all entries, so the positions of the local orbits remain relevant.

\begin{proposition}[A local cyclotomic-factor test at arbitrary precision]
\label{revadic:prop-local-factor-test}
Continue with the notation of Theorem~\ref{revadic:thm-integral-resolvent}. Let $F(T)\in\mathbb Z_p[T]$ be the monic Hensel factor of $\Phi_m$ satisfying $F\bmod p=f$ and $F(\zeta)=0$, and set $F_s=F\bmod p^s$. Then for every $s\ge1$,
\begin{equation}\label{revadic:eq-local-factor-test}
 \boxed{\displaystyle
 v_p(c_j)\ge s
 \iff A_{C,j}(T)\equiv F_s(T)U_s(T)\pmod{p^s}
 \text{ for some }U_s\in(\mathbb Z/p^s\mathbb Z)[T].}
\end{equation}
Thus one may successively test local-factor divisibility modulo $p,p^2,\ldots$ to locate the first nonzero digit of each nonzero coefficient. Unlike the Witt-digit method, this test works directly with the \emph{integer polynomial} furnished by the orbit counts.
\end{proposition}
\begin{proof}
The unramified Hensel factorization gives the evaluation isomorphism
\[
 (\mathbb Z/p^s\mathbb Z)[T]/(F_s)
 \ \simeq\ \mathcal O_h/p^s\mathcal O_h,
 \qquad h=\operatorname{ord}_m(p),
\]
with $T$ mapping to $\zeta$. Since $A_{C,j}(\zeta)=(-1)^jc_j$, vanishing modulo $p^s$ is exactly the assertion that $A_{C,j}$ is zero in this quotient. Division by the monic polynomial $F_s$ gives the displayed divisibility criterion.
\end{proof}

The local-factor criterion retains the position of each component. If only the valuations and their multiplicities are needed, Smith form can instead be applied directly to the integer counting polynomial.

\begin{theorem}[Recovering local valuation multisets from an integer multiplication matrix]
\label{v143:thm-valuation-smith}
Let $m\ge2$, $p\nmid m$, and $h=\ord_m(p)$. Let $B\in\mathbb Z[T]$ satisfy $0\ne B$ and $\deg B<\varphi(m)$. On the free module
\[
 \mathcal R_m=\mathbb Z_p[T]/(\Phi_m),
\]
write $M_B$ for multiplication by $B$, with Smith exponents $e_1\le\cdots\le e_{\varphi(m)}$. Let $D$ run through the $p$-Frobenius orbits of primitive $m$th roots, and put
$w_D=v_p(B(\zeta_m^u))$ for $u\in D$. Then
\begin{equation}\label{v143:eq-valuation-smith}
 \{e_1,\ldots,e_{\varphi(m)}\}
 =\bigsqcup_D\{\underbrace{w_D,\ldots,w_D}_{h\text{ times}}\}
 \quad\text{as multisets}.
\end{equation}
If $\Delta_i$ is the positive greatest common divisor of all $i$th-order minors of the integer power-basis matrix $M_B$, with $\Delta_0=1$, then
\[
 e_i=v_p(\Delta_i)-v_p(\Delta_{i-1}).
\]
Alternatively, set $k_0=0$ and
$k_s=\log_p|\ker(M_B\bmod p^s)|$ for $s\ge1$. Then
\begin{equation}\label{v143:eq-valuation-kernels}
 k_s=h\sum_D\min(s,w_D),\qquad
 \#\{D:w_D\ge s\}=\frac{k_s-k_{s-1}}h.
\end{equation}
\end{theorem}
\begin{proof}
Coprime Hensel factorization and the Chinese remainder theorem give
\[
 \mathcal R_m\simeq\prod_D\mathbb Z_p[T]/(F_D)
 \simeq\prod_D\mathcal O_h.
\]
On component $D$, the multiplier is $B(\zeta_m^u)=p^{w_D}v_D$, where $v_D$ is a unit. It is nonzero because $\Phi_m$ is irreducible over $\mathbb Q$ and $\deg B<\varphi(m)$. Multiplication by $v_D$ is an invertible integral operator. Thus this rank-$h$ component contributes $h$ copies of $w_D$, proving~\eqref{v143:eq-valuation-smith}.

Invertible integral row and column operations preserve the ideals generated by minors of each order. In Smith form, the least valuation of an $i$th-order minor is $e_1+\cdots+e_i$, giving the gcd formula. Modulo $p^s$, the kernel on component $D$ has $p^{h\min(s,w_D)}$ elements. Multiplication of these cardinalities and successive differences give~\eqref{v143:eq-valuation-kernels}.
\end{proof}

For $B=B_{C,j}\ne0$, the least Smith exponent is $\nu_{C,j}$, and the sum of all exponents is
$v_p\operatorname{Res}(\Phi_m,B)$; the intermediate exponents give the multiplicity of each depth. Integer minors therefore determine the multiset of all local valuations, but do not label the corresponding orbits. Local factors or phase trace duality retain those labels.

\begin{corollary}[Frobenius-conjugate cosets have identical valuation profiles]
\label{revadic:cor-profile-frob}
For each coefficient, $c_j(pC)=\phi(c_j(C))$; consequently
\[
 \mathcal V_p(pC)=\mathcal V_p(C).
\]
If $u$ is the smallest positive integer for which $p^uC=C$, then $\prod_{t=0}^{u-1}Q_{p^tC}(X)$ is the coefficient norm of $Q_C(X)$ from $K_u$ to $\mathbb Q_p$. It belongs to $\mathbb Z_p[X]$, and its root multiset is the union, with multiplicities, of the conjugate root sets.
\end{corollary}
\begin{proof}
The identity $\phi(\zeta^i)=\zeta^{pi}$ sends $Q_C$ to $Q_{pC}$, and $\phi$ preserves valuations. The $u$ factors are the distinct Galois conjugates of this polynomial over its coefficient field $K_u$; their product is its coefficient norm.
\end{proof}

\subsection{Phase differences and the index of a generating order}
\label{v141:sec-phase-index}
The norm constrains the valuations of the phase values themselves. Their pairwise differences contain additional information: two phases can remain congruent beyond the depth forced in their direction. We first keep all phases together, so that a local cancellation remains distinguishable from a common divisibility.

Let $\Gamma$ be a finite abelian group of order $R\ge2$, let $\tau\in\Gamma$ have order $f$, and put $P=\langle\tau\rangle$. With $\phi$ the Frobenius of $\mathcal O_f/\mathbb Z_p$, define the \textbf{unramified phase algebra}
\begin{equation}\label{v141:eq-phase-algebra}
 \mathcal A=\{(z_u)_{u\in\Gamma}\in\mathcal O_f^\Gamma:
 z_{\tau u}=\phi(z_u)\},\qquad
 \mathcal B=\mathcal A/p\mathcal A.
\end{equation}
Choosing one representative of each $P$-coset identifies these algebras with
$\prod_{\Gamma/P}\mathcal O_f$ and $\prod_{\Gamma/P}\mathbb F_{p^f}$, respectively. In particular, $\mathcal A$ has rank $R$ over $\mathbb Z_p$ and has a perfect trace pairing. A direction $a\in\Gamma$ acts by $(\sigma_a z)_u=z_{ua}$.

For a cyclotomic coefficient, take $x=A_{C,j}(\zeta_m)\ne0$, and use its \emph{coefficient stabilizer}
\[
 H_j=\{a\in(\mathbb Z/m\mathbb Z)^\times:
                  A_{C,j}(\zeta_m^a)=A_{C,j}(\zeta_m)\}.
\]
Then $\Gamma=(\mathbb Z/m\mathbb Z)^\times/H_j$, $\tau=pH_j$, and $x_u=A_{C,j}(\zeta_m^u)$ give an element of~\eqref{v141:eq-phase-algebra} with distinct coordinates. The subgroup $H_j$ may strictly contain the stabilizer of $C$. It is determined over the integers by
\[
 a\in H_j\iff \Phi_m(T)\mid A_{C,j}(T^a)-A_{C,j}(T),
\]
using a nonnegative representative of $a$. If $\Gamma$ is trivial, there is just one phase and no pairwise difference to consider.

Now let $x\in\mathcal A$ have distinct coordinates. Write
\[
 P_x(Y)=\prod_{u\in\Gamma}(Y-x_u)\in\mathbb Z_p[Y],\qquad
 I_x=\log_p|\mathcal A/\mathbb Z_p[x]|.
\]
The polynomial has coefficients in $\mathbb Z_p$ because Frobenius permutes its roots. Distinctness makes $1,x,\ldots,x^{R-1}$ a basis after inverting $p$, so $I_x$ is finite. For $a\ne1$, set
\begin{equation}\label{v141:eq-direction-depth}
\begin{aligned}
 h(a)&=\min_u v_p(x_u-x_{ua}),&
 \epsilon_{u,a}&=v_p(x_u-x_{ua})-h(a),\\
 y_a&=p^{-h(a)}(x-\sigma_a x),& b_a&=\bar y_a\in\mathcal B.
\end{aligned}
\end{equation}
Put $h(1)=+\infty$, and let $M_a$ be multiplication by $b_a$ on $\mathcal B$. These $\epsilon_{u,a}$ concern \emph{differences}; they are distinct from the coefficient depths $\epsilon_{C,j}$ defined above.

\begin{theorem}[Index excess and successive cancellation layers]
\label{v141:thm-index-budget}
Define
\[
 I_{\rm forced}=\frac R2\sum_{a\ne1}h(a),\qquad
 E_x=I_x-I_{\rm forced},\qquad
 g_f=\begin{cases}f,&f\text{ odd},\\f/2,&f\text{ even}.\end{cases}
\]
Then $I_{\rm forced}$ and $E_x$ are nonnegative integers, and
\begin{align}
 I_x&=\tfrac12v_p(\operatorname{disc}P_x)
      =\sum_{\{u,v\}\subseteq\Gamma}v_p(x_u-x_v),\label{v141:eq-index-discriminant}\\
 E_x&=\frac12\sum_{a\ne1}\dim_{\mathbb F_p}\ker M_a+E_{\ge2},\label{v141:eq-budget}\\
 E_{\ge2}&=\sum_{\{u,v\}\subseteq\Gamma}
             \max\{\epsilon_{u,u^{-1}v}-1,0\}.
\end{align}
The number $E_x$ is the \textbf{index excess}. Both $E_x$ and $E_{\ge2}$ belong to $g_f\mathbb Z_{\ge0}$. In particular,
\[
 E_x=\tfrac12\sum_{a\ne1}\dim\ker M_a
 \iff \epsilon_{u,a}\le1\quad(u\in\Gamma, a\ne1).
\]
If the nonnegative difference between the two sides of this equality is less than $g_f$, it is zero. If it equals $g_f$, exactly one shortest Frobenius orbit of unordered pairs has excess two; every other pair has excess at most one.
\end{theorem}
\begin{proof}
In an integral basis of $\mathcal A$, the trace Gram determinant is a unit. The change-of-basis matrix from $1,x,\ldots,x^{R-1}$ has determinant of valuation $I_x$. The Gram determinant in this power basis is $\operatorname{disc}P_x$, proving the first equality; the Vandermonde determinant proves the second.

Since $h(a^{-1})=h(a)$, summing $h(u^{-1}v)$ over unordered pairs gives $I_{\rm forced}$. Thus
$E_x=\sum_{\{u,v\}}\epsilon_{u,u^{-1}v}$ is a nonnegative integer. After extending the residue field to an algebraic closure, multiplication by $b_a$ is diagonal, with one entry for each phase. Its zero entries are exactly the $u$ for which $\epsilon_{u,a}\ge1$. Scalar extension preserves dimension, so
\[
 \dim\ker M_a=\#\{u:\epsilon_{u,a}\ge1\}.
\]
More generally, the dimension of an intersection of such kernels counts the phases satisfying all the corresponding first-layer conditions. Summing
$\epsilon=\sum_{r\ge1}\mathbf1_{\epsilon\ge r}$ over directed pairs and dividing by two proves~\eqref{v141:eq-budget}.

Translation by $\tau$ acts freely on vertices, with orbit length $f$. An unordered pair has orbit length $f$, except that for even $f$ a pair $\{u,\tau^{f/2}u\}$ has orbit length $f/2$. Indeed, a power of $\tau$ stabilizing a pair either fixes both vertices or exchanges them. Frobenius preserves $h$ and $\epsilon$, so each sum is an integer-weighted sum of these orbit lengths. This gives the divisibilities. Nonnegativity gives the equality criterion and the two final assertions.
\end{proof}

The required integer invariant is therefore one discriminant valuation; the first exceptional phases are detected by multiplication kernels over finite fields. Their comparison determines when no further excess layer is possible.

\subsection{Minimal excess and the last congruence layer}
\label{v141:sec-minimal-excess}
The global directional depths themselves form a filtration. Set
\[
 K_s=\{a\in\Gamma:h(a)\ge s\}\quad(s\ge0),\qquad
 h_{\max}=\max_{a\ne1}h(a).
\]
The identity $x-\sigma_{ab}x=(x-\sigma_ax)+\sigma_a(x-\sigma_bx)$ gives
$h(ab)\ge\min\{h(a),h(b)\}$, so each $K_s$ is a subgroup. The values at all vertices of a coset $uK_s$ are congruent modulo $p^s$.

\begin{theorem}[The last layer under minimal positive excess]
\label{v141:thm-minimal-excess}
Suppose $E_x=g_f$. There is exactly one Frobenius orbit $\mathscr E$ of unordered pairs such that
\[
 \epsilon_{u,u^{-1}v}=\begin{cases}1,&\{u,v\}\in\mathscr E,\\0,&\text{otherwise}.\end{cases}
\]
If $a=u^{-1}v$ is an exceptional direction, then $h(a)=h_{\max}$. At precision $s=h_{\max}+1$, the congruent distinct pairs are exactly $\mathscr E$. Their graph has one of the following forms:
\begin{enumerate}
\item If $f$ is even, $a=\tau^{f/2}$, and the graph is a matching of $f/2$ edges within one Frobenius vertex orbit.
\item If $f$ is odd and $a\notin P$, it is a matching of $f$ edges between two Frobenius vertex orbits.
\item If $f$ is odd and $a\in P$, then $a$ has order three, and the graph consists of disjoint triangles on one Frobenius vertex orbit, with $f$ edges in total.
\end{enumerate}
In particular,
\begin{equation}\label{v141:eq-last-depth}
 v_p(x_u-x_v)\le h_{\max}+1\quad(u\ne v),
\end{equation}
and equality holds precisely on $\mathscr E$.
\end{theorem}
\begin{proof}
The pair-orbit lengths in Theorem~\ref{v141:thm-index-budget} force exactly one shortest orbit to have excess one. For even $f$ it is the antipodal matching. For odd $f$, a direction outside $P$ gives a matching between two vertex orbits; a direction in $P$ gives disjoint odd cycles of length $\operatorname{ord}(a)$.

Put $s=h(a)+1$ and suppose $k=|K_s|\ge2$. The cosets $uK_s$ and $uaK_s$ are distinct. Every pair joining these two cosets is congruent modulo $p^s$, whereas its direction lies outside $K_s$. All $k^2$ joining edges must therefore belong to $\mathscr E$. A matching cannot contain these edges. In the odd-cycle case every vertex has degree two, so $k\le2$. If $k=2$, its unique nonidentity element has order two and cannot belong to the odd-order group $P$. Yet all endpoints of $\mathscr E$ lie in one $P$-coset, which cannot contain $uK_s$. This is again impossible. Hence $K_s=\{1\}$ and $h(a)=h_{\max}$.

At this precision the congruence graph has exactly the exceptional edges. Congruence is an equivalence relation, so every connected component of this graph is complete. An odd cycle is complete only when its length is three. The stated classification and~\eqref{v141:eq-last-depth} follow.
\end{proof}

\begin{proposition}[Finite-field kernels locate the exceptional phases]
\label{v141:prop-exceptional-kernels}
Under Theorem~\ref{v141:thm-minimal-excess}, the kernels vanish outside directions $a,a^{-1}$. If $a\ne a^{-1}$, each of these two kernels has dimension $f$: they are disjoint for a matching between vertex orbits, and equal for triangles. If $a=a^{-1}$, there is only one exceptional direction, and
\[
 \dim\ker M_a=\begin{cases}
 f,&f\text{ even},\ a=\tau^{f/2},\\
 2f,&f\text{ odd},\ a\notin P.
 \end{cases}
\]
In every case
\begin{equation}\label{v141:eq-phase-idempotent}
 e_a=1-b_a^{p^f-1}\in\mathcal B
\end{equation}
is the idempotent indicating the exceptional starting phases. It is primitive except in the second case of the displayed formula, where it is the sum of two primitive idempotents.
\end{proposition}
\begin{proof}
A kernel counts starting vertices of directed exceptional edges. A direction different from its inverse gives one starting orbit, while an involution counts both endpoints of every edge. In the antipodal case these endpoints form one orbit; in the odd-$f$ matching case they form two. The assertions about intersections follow from their supports. In each field component, $1-b_a^{p^f-1}$ is one exactly when $b_a$ is zero, proving the last assertion.
\end{proof}

\begin{corollary}[A cutoff for valuations of the phase values]
\label{v141:cor-phase-cutoff}
Under Theorem~\ref{v141:thm-minimal-excess}, put $s=h_{\max}+1$ and
$Z_t=\{u:v_p(x_u)\ge t\}$. Then $|Z_s|\le3$ and $|Z_{s+1}|\le1$. More precisely,
\[
 \begin{array}{ll}
 f>3:&Z_s=\varnothing,\\
 f=2\text{ or }3:&|Z_s|\in\{0,f\},\quad Z_{s+1}=\varnothing,\\
 f=1:&|Z_s|\le2,\quad |Z_{s+1}|\le1.
 \end{array}
\]
For $f=2$ a nonempty $Z_s$ is an antipodal edge; for $f=3$ it is a triangle. If all $x_u$ are nonzero, then
\begin{equation}\label{v141:eq-phase-norm}
 v_p\!\left(\prod_u x_u\right)=\sum_{t\ge1}|Z_t|.
\end{equation}
Thus when $f=1$, the depth of a single phase remaining after level $s$ is determined by this norm and the other known depths.
\end{corollary}
\begin{proof}
Any two vertices of $Z_t$ are congruent modulo $p^t$. At level $s$ they form a complete subgraph of a matching or a union of triangles; at level $s+1$ no distinct vertices are congruent. Also $Z_t$ is Frobenius-stable, so its size is a multiple of $f$. These observations give all the bounds. Summing the valuations of nonzero coordinates gives~\eqref{v141:eq-phase-norm}.
\end{proof}

\begin{example}[A triangle at the last layer]
\label{v141:ex-triangle}
Take $p=2$, $\Gamma=\mathbb Z/6\mathbb Z$, and $\tau=2$. Choose $\alpha\in\mathbb F_8$ of degree three, and prescribe
\[
 x_{2i}=2[\alpha^{2^i}],\qquad x_{2i+1}=[\alpha^{2^i}]
 \quad(0\le i<3).
\]
These coordinates define an element of the phase algebra with $f=3$. Differences within the even vertices have valuation one; differences within the odd vertices and between the two sets are units. Every nonzero direction has minimum zero, so $h_{\max}=0$ and $E_x=3=g_f$. The exceptional triangle consists of the three even vertices, exactly $Z_1$, while $Z_2$ is empty. Thus the three-vertex bound is attained.
\end{example}

The minimal-excess condition controls additional collisions beyond $h(a)$; the ordinary refinement of the subgroups $K_s$ may still occur. Nor do pairwise differences alone determine the valuations of the values: adding a common element of $\mathbb Z_p$ leaves all differences unchanged. The low-level sets $Z_t$ and the norm retain information about their position relative to zero.

\subsection{Period trace duality and explicit phase profiles}
\label{v141:sec-period-duality}
For prime conductor, trace duality turns a kernel support into labelled phases. Let $m=\ell\ne p$ be prime, let $H\le(\mathbb Z/\ell\mathbb Z)^\times$ have order $g$, put $\Gamma=(\mathbb Z/\ell\mathbb Z)^\times/H$, and take $\tau=pH$. Define
$\eta_b=\sum_{h\in H}\zeta_\ell^{bh}$, with $b\in\Gamma$. In the phase algebra its coordinate at $u$ is $\eta_{ub}$. Write $\Tr$ for the algebra trace, over $\mathbb Z_p$ or $\mathbb F_p$ as appropriate.

\begin{proposition}[Period bases and recovery of an idempotent support]
\label{v141:prop-period-duality}
The periods form a $\mathbb Z_p$-basis of $\mathcal A$, with trace-dual basis
\begin{equation}\label{v141:eq-period-dual}
 \Tr(\eta_a\eta_b)=\ell\mathbf1_{aH=-bH}-g,
 \qquad \eta_a^\vee=\frac{\eta_{-a}-g}{\ell}.
\end{equation}
For any $e\in\mathcal B$, its coordinate at phase $u$ is
\begin{equation}\label{v141:eq-phase-recovery}
 e(u)=\frac{g\Tr_{\mathcal B/\mathbb F_p}(e)+
       \sum_{b\in\Gamma}\bar\eta_{-ub}
                     \Tr_{\mathcal B/\mathbb F_p}(e\bar\eta_b)}{\ell}.
\end{equation}
For $e=e_a$ from~\eqref{v141:eq-phase-idempotent}, these coordinates are zero or one and locate every exceptional starting phase. The formula remains valid when $p\mid f$.
\end{proposition}
\begin{proof}
Expanding the product and summing over phases gives
\[
 \Tr(\eta_a\eta_b)=\sum_{r\in H}\sum_{t\in(\mathbb Z/\ell\mathbb Z)^\times}
                         \zeta_\ell^{t(a+br)}.
\]
The inner sum is $\ell-1$ when $a+br=0$, and $-1$ otherwise. This proves the Gram formula. Its matrix is $\ell P-gJ$, where $P$ is the permutation $aH\mapsto-aH$ and $J$ is the all-one matrix. Since $gR=\ell-1$, its determinant is $\pm\ell^{R-1}$, a $p$-adic unit. Perfect trace duality on $\mathcal A$ therefore makes the periods an integral basis. Using $\Tr(\eta_a)=-1$ verifies the displayed dual basis directly.

Expand $e$ in this dual basis and use $\sum_b\eta_b=-1$ to obtain
\[
 e=\frac{g\Tr(e)+\sum_b\Tr(e\bar\eta_b)\bar\eta_{-b}}{\ell}.
\]
Taking coordinate $u$ gives~\eqref{v141:eq-phase-recovery}. Only $\ell$ is inverted; no division by $f$ occurs.
\end{proof}

\begin{example}[Two small period algebras]
\label{v141:ex-period-indices}
Take $p=2$ and use the cosets $3^iH$ to label the periods $\eta_i$.
For $\ell=17$, let $H=\{1,4,13,16\}=C_4(1)$; indices are modulo four. Counting exponent sums gives
\[
 \eta_0^2=4+2\eta_1+\eta_2,\quad
 \eta_0\eta_1=2\eta_0+\eta_2+\eta_3,\quad
 \eta_0\eta_2=\eta_0+\eta_1+\eta_2+\eta_3=-1.
\]
Their translates give the remaining products. For $\ell=43$, take $H=\langle3^3\rangle=C_2(1)$, of order fourteen; indices are modulo three. Here
\[
 \eta_0^2=14+3\eta_0+6\eta_1+4\eta_2,\qquad
 \eta_0\eta_1=6\eta_0+4\eta_1+4\eta_2.
\]
Again translation and commutativity give the complete multiplication table. Using $1=-\sum_i\eta_i$, the characteristic polynomial of multiplication by $\eta_0$ is respectively
\[
 \begin{array}{c|c|c|c}
 \ell&f&P_{\eta_0}(Y)&\operatorname{disc}P_{\eta_0}\\\hline
 17&2&Y^4+Y^3-6Y^2-Y+1&4\cdot17^3\\
 43&1&Y^3+Y^2-14Y+8&4\cdot43^2
 \end{array}
\]
The period basis shows $h(a)=0$ for every nonidentity direction: the coordinate vector of a difference of distinct periods has unit content. Thus $I_x=E_x=1=g_f$ in both cases.

For $\ell=17$, the exceptional graph is one antipodal edge, and the three nonidentity directions have kernel dimensions $(0,2,0)$. Every period is a $2$-adic unit, since their product is one. Extra divisibility of a difference need not give positive valuation of either value.

For $\ell=43$, the residue algebra is $\mathbb F_2^3$, and the residue polynomial is $Y^2(Y+1)$. Exactly two phases are zero modulo two, and their difference has valuation one; both nonidentity direction kernels have dimension one. The product of the three periods has valuation three. The two positive valuations cannot both be one, since two elements of $2\mathbb Z_2\setminus4\mathbb Z_2$ differ by a multiple of four. Consequently the complete valuation multiset is
\[
 \{v_2(\eta_0),v_2(\eta_1),v_2(\eta_2)\}=\{0,1,2\}.
\]
The norm and the exceptional difference determine these depths; specifying which zero residue phase has depth two requires its position relative to zero as well.

Integer minors recover the same multiset. The component argument of Theorem~\ref{v143:thm-valuation-smith} applies to any finite product of unramified integer rings; here every component has degree one. Multiplication by $\eta_0$ in the period basis $(\eta_0,\eta_1,\eta_2)$ has matrix
\begin{equation}\label{v143:eq-period-smith-example}
 M_{\eta_0}=\begin{pmatrix}-11&6&4\\-8&4&4\\-10&4&6\end{pmatrix}.
\end{equation}
It has an odd entry, rank one modulo two, a second-order minor equal to $-26$, and determinant $-8$. Its Smith exponents are therefore $0,1,2$, giving the valuation multiset of the three periods directly, but still without assigning phase labels.
\end{example}

All the index and graph arguments used an abelian phase group, not a prime conductor. They therefore apply to composite $m$ as well. The explicit period dual basis above uses prime conductor; at composite conductor one instead uses a trace-dual integral basis of the phase algebra.

\subsection{Compressed orbits and recovery of labelled valuation profiles}
\label{v144:sec-profile-recovery}
Integer remainders detect exact zeros, Smith exponents give valuation multisets, and trace duality recovers phase coordinates. Before applying these tools, one can compress the root set and then pass to the fixed algebra of the symmetries retained by each coefficient.

\begin{proposition}[Full-fiber compression preserves coefficient valuations]
\label{v144:prop-profile-compression}
Let $p\nmid m$ and let $\varnothing\ne C\subseteq\mathbb Z/m\mathbb Z$ be stable under multiplication by $q$. Put
\[
 \kappa=|H_+(C)|,\qquad d=m/\kappa,\qquad
 C_0=C\bmod d.
\]
Then $C$ is the full inverse image of $C_0$, and $|C|=\kappa|C_0|$. With $\zeta_d=\zeta_m^\kappa$, for $u\in(\mathbb Z/m\mathbb Z)^\times$ one has
\begin{equation}\label{v144:eq-profile-compression}
 Q_C^{(u)}(X)=Q_{C_0}^{(\bar u)}(X^\kappa),\qquad
 c_j^{(u)}=
 \begin{cases}
 c_{0,j/\kappa}^{(\bar u)},&\kappa\mid j,\\
 0,&\kappa\nmid j.
 \end{cases}
\end{equation}
Superscripts multiply root exponents by the indicated unit, $c_{0,t}$ is the $t$th coefficient of $Q_{C_0}$, and $\bar u=u\bmod d$. The original valuation profile is the pullback of the compressed profile under reduction of unit groups, with $+\infty$ at the other positions.
\end{proposition}
\begin{proof}
The additive stabilizer is the subgroup of order $\kappa$ in $\mathbb Z/m\mathbb Z$, namely $d\mathbb Z/m\mathbb Z$. The roots in each fiber satisfy
\[
 \prod_{t=0}^{\kappa-1}(X-\zeta_m^{u(a+td)})
 =X^\kappa-\zeta_d^{ua}.
\]
Multiplication over $a\in C_0$ gives~\eqref{v144:eq-profile-compression}, the coefficient form of Theorem~\ref{v139:thm-period-compression}. Reduction of unit groups is surjective, so every compressed phase occurs. The valuation statement includes exact zeros.
\end{proof}

Assume henceforth that $d\ge2$. If $d=1$, then $C=\mathbb Z/m\mathbb Z$ and $Q_C=X^m-1$ applies directly. Write $U_d=(\mathbb Z/d\mathbb Z)^\times$.

\begin{theorem}[Valuation multisets in a coefficient fixed algebra]
\label{v144:thm-invariant-smith}
Let $d\ge2$, $p\nmid d$, $0\ne B\in\mathbb Z[T]$, $\deg B<\varphi(d)$, and choose a subgroup $H\le U_d$ such that
\begin{equation}\label{v144:eq-coefficient-invariance}
 B(T^a)\equiv B(T)\pmod{\Phi_d(T)}\qquad(a\in H).
\end{equation}
Put $\Lambda=\mathbb Z[T]/(\Phi_d)$, let $\sigma_a$ denote substitution $T\mapsto T^a$, and set
\[
 r=[U_d:H],\qquad f_H=\ord_{U_d/H}(pH),\qquad
 \mathcal A_H=\Lambda^H\otimes_{\mathbb Z}\mathbb Z_p.
\]
Then $\Lambda^H$ is a free integer module of rank $r$, and
\begin{equation}\label{v144:eq-invariant-algebra}
 \mathcal A_H=(\mathbb Z_p[T]/(\Phi_d))^H
 \simeq\prod_{U_d/(H\langle p\rangle)}\mathcal O_{f_H}.
\end{equation}
In any integer basis of $\Lambda^H$, multiplication by $x=B\bmod\Phi_d$ has an $r\times r$ integer matrix. Its $p$-Smith exponents are
\begin{equation}\label{v144:eq-invariant-smith}
 \bigsqcup_{D\in U_d/(H\langle p\rangle)}
 \{\underbrace{w_D,\ldots,w_D}_{f_H\text{ times}}\},
 \qquad w_D=v_p(B(\zeta_d^u)),\quad u\in D.
\end{equation}
Thus the matrix size drops from $\varphi(d)$ to $\varphi(d)/|H|$.
\end{theorem}
\begin{proof}
Each $\sigma_a$ has an integer power-basis matrix. The fixed module is defined by the integer linear equations $(\sigma_a-I)v=0$. Smith form of the stacked matrix shows that an integer kernel basis is also a basis of its $p$-adic kernel. Taking fixed elements therefore commutes with extension to $\mathbb Z_p$, without division by $|H|$.

In an unramified integer ring containing all $d$th roots of unity, label evaluation coordinates by $u\in U_d$. The original algebra has $z_{pu}=\phi(z_u)$; invariance additionally requires $z_{ua}=z_u$ for $a\in H$. Coordinates are therefore indexed by $U_d/H$. Multiplication by $pH$ has orbits of length $f_H$, each contributing $\mathcal O_{f_H}$. This proves~\eqref{v144:eq-invariant-algebra} and the rank statement.

Condition~\eqref{v144:eq-coefficient-invariance} ensures that multiplication by $x$ preserves the integer module. Irreducibility of $\Phi_d$ over $\mathbb Q$ gives $B(\zeta_d^u)\ne0$. In each local component the multiplier is $p^{w_D}$ times a unit, contributing $f_H$ identical Smith exponents. The component argument of Theorem~\ref{v143:thm-valuation-smith} gives~\eqref{v144:eq-invariant-smith}.
\end{proof}

One may take $H$ to be the multiplicative stabilizer of the compressed root set, or the full stabilizer of the coefficient. The latter is determined by~\eqref{v144:eq-coefficient-invariance} and may lower the rank further. The affine multiplier group in Proposition~\ref{v141:prop-affine-profile} only preserves coefficients up to unit factors; it cannot replace the fixed group here without checking the stated condition.

\begin{proposition}[Precision for labelling a known valuation multiset]
\label{v144:prop-label-precision}
Under the hypotheses of Theorem~\ref{v144:thm-invariant-smith}, put $x_u=B(\zeta_d^u)$ and let the distinct Smith exponents be $a_1<\cdots<a_t$. If $t=1$, every phase has valuation $a_1$ and no further evaluation is needed. If $t\ge2$, put
\begin{equation}\label{v144:eq-label-precision}
 S=a_{t-1}+1.
\end{equation}
The phase values modulo $p^S$ determine the complete labelled valuation profile: a nonzero reduction gives its exact valuation, and a zero reduction must have valuation $a_t$. Equivalently, it suffices to determine
\[
 Z_{a_i+1}=\{uH:x_u\equiv0\pmod{p^{a_i+1}}\},
 \qquad1\le i<t.
\]
Successive differences of these nested zero sets label the valuation layers.

For an integral basis $b_1,\ldots,b_r$ of $\mathcal A_H$ and its trace-dual basis $b_1^\vee,\ldots,b_r^\vee$, the needed phase values are recovered by
\begin{equation}\label{v144:eq-trace-labels}
 x_u=\sum_{i=1}^r\Tr_{\mathcal A_H/\mathbb Z_p}(xb_i^\vee)(b_i)_u
 \pmod{p^S}.
\end{equation}
For prime conductor, the period trace-dual basis of Proposition~\ref{v141:prop-period-duality} applies directly.
\end{proposition}
\begin{proof}
By~\eqref{v144:eq-invariant-smith}, each $w_D$ lies in the known set $\{a_1,\ldots,a_t\}$. A nonzero reduction modulo $p^S$ has a first nonzero digit, giving the exact valuation below $S$. A zero reduction has $w_D\ge S>a_{t-1}$, forcing $w_D=a_t$. The nested zero sets give the same classification.

The decomposition~\eqref{v144:eq-invariant-algebra} gives a perfect trace pairing and hence an integral trace-dual basis. Taking coordinate $u$ in the dual-basis expansion gives~\eqref{v144:eq-trace-labels}. Reduction modulo $p^S$ commutes with these operations; no hypothesis $p\nmid f_H$ is needed.
\end{proof}

This precision bound cannot be lowered in general. In $\mathbb Z_p\times\mathbb Z_p$, the elements $x=(p^a,p^b)$ and $y=(p^b,p^a)$, where $0\le a<b$, have the same Smith exponents and agree modulo $p^a$, but have different labelled profiles. Modulo $p^{a+1}$ they are distinguished. For $a=0$, precision modulo $p^0$ means no precision data. Additional relations among specified phases may permit lower precision.

These results give a finite recovery procedure. First compress $C$ to $C_0$ using the additive stabilizer, retaining only coefficient indices divisible by $\kappa$. Next read the integer counting polynomials as coefficients of powers of $Y$ in
\[
 \prod_{a\in C_0}(1+YT^a)
 \quad\text{in }\mathbb Z[Y,T]/(T^d-1),
\]
and reduce modulo $\Phi_d$; zero remainders give exact zeros. For each nonzero remainder, choose a known coefficient fixed group and compute the Smith exponents of the smaller multiplication matrix. Finally use local factors or~\eqref{v144:eq-trace-labels}, only at the required precision $S$, to attach phase labels, and pull the profile back by~\eqref{v144:eq-profile-compression}.

The required arithmetic data include the factorization of the modulus, the relevant groups and orbits, and an integer basis and multiplication table of the chosen fixed module. Phase evaluation also requires a fixed residue root and compatible lifts. For a fixed algebra of rank $r$, once the trace-dual data and phase values of the basis modulo $p^S$ are known, choose one representative of each Frobenius orbit. All evaluations in~\eqref{v144:eq-trace-labels} then require $O(r^2)$ arithmetic operations in $\mathbb Z/p^S\mathbb Z$, expressing each degree-$f_H$ component in $f_H$ base-ring coordinates. This bound excludes integer factorization, Smith computation, preparation of integral bases, and the bit lengths of coefficients.

\begin{example}[Three-dimensional recovery at the square of forty-three]
\label{v144:ex-profile-tower}
Take $p=q=2$, $m=43^2$, and $C=C_2(1)\pmod{43^2}$. Since
\[
 \ord_{43}(2)=14,\qquad 2^{14}-1=43\cdot381,
 \qquad43\nmid381,
\]
stable lifting in Chapter~5 gives $\ord_{43^2}(2)=602$. The orbit $C$ is the full inverse image of $H=C_2(1)$ modulo 43. The latter is a nonempty proper subset modulo 43 and has trivial additive stabilizer, so $|H_+(C)|=43$. Consequently
\[
 Q_C(X)=Q_H(X^{43}).
\]
Only 14 nonleading coefficients need recovery; the other 588 coefficients vanish exactly. For coefficient 43, $c_{43}=-\eta_0$, where the $\eta_i$ are the three periods of Example~\ref{v141:ex-period-indices}. The uncompressed integer multiplication matrix has size $\varphi(43^2)=1806$. Reduction to modulus 43, followed by passage to the fixed algebra of the order-fourteen multiplicative stabilizer, leaves the existing three-by-three matrix~\eqref{v143:eq-period-smith-example}.

Its valuation multiset is $\{0,1,2\}$, so attaching phase labels requires only reduction modulo four. After a phase translation, take $\eta_0$ to be a unit. Its minimal polynomial $Y^3+Y^2-14Y+8$ has the unique odd root 1 modulo four. The existing product relation and period sum give
\[
 \eta_0\equiv1,\qquad
 \eta_0\eta_1\equiv2\eta_0\pmod4,
 \qquad \eta_0+\eta_1+\eta_2=-1.
\]
Hence $(\eta_0,\eta_1,\eta_2)\equiv(1,2,0)\pmod4$, with labelled valuations $(0,1,2)$. The last exact value comes from the Smith multiset; precision modulo eight is unnecessary. For the original modulus, $c_{43}^{(u)}$ has valuation $i$ when $u\bmod43\in3^iH$, with labels based at the chosen unit phase.
\end{example}

\subsection{Prime-power moduli: subset-sum uniformity and exact zeros}
\label{revadic:subsec-prime-power-profile}

For a prime-power modulus, the cyclotomic remainder has particularly transparent integral coordinates. Exact zeros can consequently be read directly from subset-sum counts, with no need to approximate $p$-adic roots of unity.

\begin{theorem}[Uniformity criteria for prime-power cyclotomic coefficients]
\label{revadic:thm-prime-power-profile}
Let $m=\ell^a$, where $\ell\ne p$ is prime and $a\ge1$, and put $u=\ell^{a-1}$. Let $C\subseteq\mathbb Z/m\mathbb Z$ be stable under multiplication by $p^d$. For a fixed $j$, define, for $0\le s<u$ and $0\le t\le\ell-2$,
\[
 D_j(s,t)=N_{C,j}(s+tu)-N_{C,j}(s+(\ell-1)u).
\]
Then
\begin{equation}\label{revadic:eq-primepower-rem}
 B_{C,j}(T)=\sum_{s=0}^{u-1}\sum_{t=0}^{\ell-2}D_j(s,t)T^{s+tu},
\end{equation}
so that
\begin{equation}\label{revadic:eq-primepower-zero}
 \boxed{\displaystyle
 c_j=0\iff
 N_{C,j}(s)=N_{C,j}(s+u)=\cdots=N_{C,j}(s+(\ell-1)u)
 \quad(0\le s<u).}
\end{equation}
Equivalently, define the additive stabilizer of the integer-valued count function by
\[
 \operatorname{Stab}^{+}(N_{C,j})=
 \{h\bmod m:N_{C,j}(s+h)=N_{C,j}(s)\ \text{for every }s\}.
\]
Every nontrivial subgroup of a cyclic prime-power group contains its unique subgroup of order $\ell$. Consequently,~\eqref{revadic:eq-primepower-zero} is also equivalent to
\[
 \boxed{c_j=0\iff\operatorname{Stab}^{+}(N_{C,j})\ne\{0\}.}
\]
This is distinct from the multiplicative stabilizer of the original root set: a subset-sum distribution can have an additive period even when the roots have no nontrivial multiplicative stabilizer.

If $\delta_j=\min_{s,t}v_p(D_j(s,t))$, then $v_p(c_j)\ge\delta_j$. Moreover, whenever $\operatorname{ord}_{\ell^a}(p)=\varphi(\ell^a)$,
\begin{equation}\label{revadic:eq-primepower-exact}
 \boxed{\displaystyle v_p(c_j)=\delta_j.}
\end{equation}
For this explicit class of parameters, the entire valuation profile is determined by the \emph{departure from uniformity} in integral subset-sum counts. If all differences vanish, the valuation is understood as $+\infty$.
\end{theorem}
\begin{proof}
The cyclotomic polynomial for a prime power is
$\Phi_{\ell^a}(T)=1+T^u+\cdots+T^{(\ell-1)u}$.
Thus, modulo $\Phi_{\ell^a}$,
$T^{s+(\ell-1)u}\equiv-\sum_{t=0}^{\ell-2}T^{s+tu}$.
Collecting the terms of $A_{C,j}$ gives~\eqref{revadic:eq-primepower-rem}. Its degree is strictly less than $\varphi(\ell^a)=(\ell-1)u$, so it is the unique remainder. The remaining assertions are precisely the first, second, and third parts of Theorem~\ref{revadic:thm-integral-resolvent}.
\end{proof}

\begin{example}[Local cancellation for seventh roots and an exact valuation for nineteenth roots]
\label{revadic:ex-sharp-vs-cancel}
For $p=2,m=7,C=\{1,2,4\}$ and $j=1$, the counting polynomial is $A_{C,1}(T)=T+T^2+T^4$. Its degree is less than six, so $B_{C,1}=A_{C,1}$ and $\nu_{C,1}=0$. Modulo two, $A_{C,1}=T(T^3+T+1)$, so only the local component corresponding to $T^3+T+1$ has excess cancellation. Example~\ref{v132:ex-seven-period} already determines the two period valuations, giving local depths $1,0$. The lowest layer of the integer remainder is unchanged; its evaluations at the different local roots distinguish the components.

Now take $p=2$ and $m=19$. Since $\operatorname{ord}_{19}(2)=18=\varphi(19)$, the $4$-cyclotomic coset
\[
 C=\{1,4,16,7,9,17,11,6,5\}
\]
is the set of nonzero quadratic residues modulo nineteen. Every nonzero residue class occurs exactly twice as the unordered sum of two \emph{distinct} elements of $C$, while zero never occurs. Consequently,
\[
 A_{C,2}(T)=2\sum_{s=1}^{18}T^s
 =2\bigl(\Phi_{19}(T)-1\bigr),
 \qquad \boxed{c_2=-2,\quad v_2(c_2)=1.}
\]
No local periods need to be evaluated: both the exterior-power trace and its positive valuation are determined by finite subset-sum counts. The following proposition explains why this uniformity is not accidental.
\end{example}

\begin{proposition}[The second exterior-power trace of a quadratic-residue root block]
\label{revadic:prop-qr-exterior}
Let $\ell\equiv3\pmod8$ be prime, let $H\subseteq(\mathbb Z/\ell\mathbb Z)^\times$ be the set of quadratic residues, and let $\zeta$ be a primitive $\ell$th root of unity. Then
\begin{equation}\label{revadic:eq-qr-exterior}
 \boxed{\displaystyle
 \sum_{\{x,y\}\subset H\atop x\ne y}\zeta^{x+y}
 =-\frac{\ell-3}{8}.}
\end{equation}
For $\ell>3$, the coefficient of $X^{|H|-2}$ in $Q_H(X)=\prod_{x\in H}(X-\zeta^x)$ is exactly $-(\ell-3)/8$. If, moreover, $2$ is a primitive root modulo $\ell$, then $H$ is a single $4$-cyclotomic coset, and~\eqref{revadic:eq-qr-exterior} gives the second exterior-power trace of its unramified relative norm, with exact $2$-adic valuation
\[
 v_2(c_2)=v_2\!\left(\frac{\ell-3}{8}\right)
 \quad(\ell>3).
\]
For $\ell=3$ there is no second exterior-power coefficient; the identity itself is interpreted as an empty sum.
\end{proposition}
\begin{proof}
Let $\chi$ be the quadratic character modulo $\ell$, with $\chi(0)=0$. Since $\ell\equiv3\pmod8$, $\chi(-1)=\chi(2)=-1$. For a fixed $s\ne0$, use
$1_H(x)=(1_{x\ne0}+\chi(x))/2$ and $\sum_x\chi(x)\chi(s-x)=-\chi(-1)=1$ to obtain the number of \emph{ordered} pairs $x,y\in H$ with $x+y=s$:
\[
 R(s)=\frac{\ell-1-2\chi(s)}4.
\]
The diagonal $x=y$ contributes $(1-\chi(s))/2$ pairs, since $x=y=s/2$ lies in $H$ exactly when $\chi(s)\chi(2)=1$. The number of unordered pairs with distinct entries is therefore
\[
 \frac12\left(R(s)-\frac{1-\chi(s)}2\right)
 =\frac{\ell-3}{8},
\]
independent of $s$. There is no pair summing to zero because $\chi(-1)=-1$. Summing over $s\ne0$ and using $\sum_{s=1}^{\ell-1}\zeta^s=-1$ proves the identity. If $2$ is primitive modulo $\ell$, its square $4$ generates $H$, so the relative-norm and exterior-power interpretation follows from Theorem~\ref{revadic:thm-orbit-exterior}.
\end{proof}

\begin{proposition}[Exact local valuations of quadratic periods]
\label{v137:prop-quadratic-period-depth}
Let $\ell\equiv3\pmod4$ be prime, let $H$ be the set of nonzero quadratic residues modulo $\ell$, and put $N=(\mathbb Z/\ell\mathbb Z)^\times\setminus H$. Let $p\ne\ell$ be prime with $p\bmod\ell\in H$. In an unramified $p$-adic extension containing a primitive $\ell$th root of unity $\zeta$, set
\[
 \eta_H=\sum_{a\in H}\zeta^a,\qquad
 \eta_N=\sum_{a\in N}\zeta^a.
\]
Then $\eta_H,\eta_N\in\mathbb Z_p$, and
\begin{equation}\label{v137:eq-quadratic-period-depth}
 \eta_H+\eta_N=-1,\quad
 \eta_H\eta_N=\frac{\ell+1}{4},\quad
 \{v_p(\eta_H),v_p(\eta_N)\}
 =\left\{0,v_p\!\left(\frac{\ell+1}{4}\right)\right\}.
\end{equation}
If $p\mid(\ell+1)/4$, the period reducing to zero modulo $p$ has positive valuation. Thus the valuations of the coefficients of $X^{(\ell-3)/2}$ in $Q_H$ and $Q_N$ are determined without successive lifting. In particular, for $p=2$ and $\ell\equiv7\pmod8$, the positive valuation is $v_2(\ell+1)-2$.
\end{proposition}
\begin{proof}
Since $p\bmod\ell\in H$, Frobenius fixes both periods, so they belong to $\mathbb Z_p$. Their sum is $-1$.
Put $h=(\ell-1)/2$. Since $-1\in N$, exactly $h$ pairs in $H\times N$ sum to zero. Each nonzero residue occurs equally often as a sum: multiplication by a square preserves the two sets, while multiplication by a nonsquare followed by interchanging the coordinates gives a bijection. The common number is
\[
 \frac{h^2-h}{\ell-1}=\frac{\ell-3}{4}.
\]
Consequently $\eta_H\eta_N=h-(\ell-3)/4=(\ell+1)/4$. Their sum is a unit, so at least one period is a unit and the other has the valuation of their product.
\end{proof}

For $p=2$ and $\ell=7,31,127$, the positive valuations are $1,3,5$, respectively. The first is Example~\ref{v132:ex-seven-period}. The quadratic-residue root blocks for the latter two contain three and nine binary cyclotomic cosets, since $\operatorname{ord}_{31}(2)=5$ and $\operatorname{ord}_{127}(2)=7$. Thus the formula determines the trace of the entire root block; its valuation cannot be assigned directly to the traces of the individual irreducible factors. In the dyadic case, each block corresponds to one irreducible factor precisely when $\operatorname{ord}_\ell(2)=(\ell-1)/2$.

Here $A_H(T)=\sum_{a\in H}T^a$ is already a nonzero cyclotomic remainder, with lowest integer valuation zero. Writing $h=(\ell-1)/2$, we have
\[
 \operatorname{Res}_T(\Phi_\ell,A_H)
 =(\eta_H\eta_N)^h=\left(\frac{\ell+1}{4}\right)^h.
\]
The multiplicative stabilizer gives $G=H$ and $[U:G]=2$, realizing the two-component case of Proposition~\ref{v137:prop-norm-depth}. An odd-characteristic example is $p=3,\ell=11$: since $\operatorname{ord}_{11}(3)=5$, the periods of the two orbits have valuations $0,1$.

Three levels of structure can now be compared. Multiplicative symmetry of a root set sometimes forces an \emph{exact gap} in its polynomial. Integral uniformity of subset-exponent sums detects more general exact zeros. When the local cyclotomic factor introduces no additional linear relations modulo $p$, the departure from uniformity even gives an \emph{exact positive valuation}. Comparing seventh and nineteenth roots shows why a high-precision prediction must retain both integer counting data and the local Frobenius orbit, rather than relying on a single vanishing coefficient modulo $p$.

\begin{corollary}[Power-sum periods and Newton's recurrence]
\label{revadic:cor-period-newton}
Write $P_h(C)=\sum_{i\in C}\zeta^{hi}=\eta_{hC}$, with multiplicities retained in $hC$. Then, in $K_d[[T]]$,
\begin{equation}\label{revadic:eq-period-exp}
 \sum_{j=0}^r c_jT^j=\prod_{i\in C}(1-\zeta^iT)
 =\exp\left(-\sum_{h\ge1}\frac{P_h(C)}hT^h\right),
\end{equation}
so that
\begin{equation}\label{revadic:eq-period-rec}
 jc_j=-\sum_{h=1}^jP_h(C)c_{j-h}\qquad(1\le j\le r).
\end{equation}
Thus all coefficients of the canonical lift can be recovered from Gauss periods attached to multiples of the exponent set. If $p\mid j$, division by $j$ costs $p$-adic working precision, and the integral subset-count formula \eqref{revadic:eq-period-profile} is preferable for valuation tests.
\end{corollary}
\begin{proof}
Take the formal logarithm of the finite product and use $\log(1-z)=-\sum_{h\ge1}z^h/h$. Differentiating with respect to $T$ and comparing coefficients yields Newton's recurrence.
\end{proof}

\begin{example}[Newton recurrence and precision for seventh roots]
\label{revadic:ex-period-count}
Use the notation of Example~\ref{v132:ex-seven-period}, with $C=\{1,2,4\}$. Since $2C=C$, we have $P_1(C)=P_2(C)=\eta$. Newton recurrence gives
\[
 c_1=-\eta,\qquad 2c_2=\eta^2-\eta=2\eta',
\]
the last equality following directly from~\eqref{v137:eq-seven-periods}.
Knowing only $\eta\bmod2$ does not recover $c_2\bmod2$ from $2c_2\bmod2$. To obtain $c_2\bmod2^s$ by division by two, the numerator must be known modulo $2^{s+1}$.
In contrast, the exponent sums of the three two-element subsets of $C$ are exactly $\{3,5,6\}$, giving $c_2=\eta'$ without division. This explains why the second coefficient is a unit although $v_2(\eta)=1$, and compares the precision required by power-sum recurrence and the subset-sum expression.
\end{example}

Integer remainders, stabilizers, and norms constrain coefficient valuations directly. When these conditions do not yet determine the local cancellation depth, Witt carries can be computed successively in the residue field. First consider precision modulo $p^2$.

\begin{proposition}[A finite-field test for the first coefficient carry]
\label{revadic:prop-first-profile}
Use the notation of~\eqref{revadic:eq-profile}, and let $\alpha_i=\bar\zeta_m^i$ for $i\in C$. For a fixed $1\le j\le r$, list in any order the $j$-element subsets $I_1,\ldots,I_N$ of $C$, and put
\[
 b_h=\prod_{i\in I_h}\alpha_i\in k_{\operatorname{lcm}(d,\operatorname{ord}_m(p))},
 \qquad s_0=0,\quad s_h=\sum_{t=1}^h b_t.
\]
If $\bar c_j=(-1)^js_N=0$, the first Witt carry $\kappa$ computes the first nonzero candidate digit entirely in the finite field:
\begin{equation}\label{revadic:eq-first-profile}
 \overline{c_j/p}=(-1)^j\sum_{h=1}^N\kappa(s_{h-1},b_h).
\end{equation}
Consequently $v_p(c_j)=1$ precisely when the right side of~\eqref{revadic:eq-first-profile} is nonzero; if it vanishes, then $v_p(c_j)\ge2$ or $c_j=0$. The test uses only finite-field products of roots and the first carry, not a full $p$-adic factor calculation.
\end{proposition}
\begin{proof}
The canonical roots $\zeta_m^i$ are the Teichmüller representatives $[\alpha_i]$, so the $j$th elementary symmetric sum is
\[
 (-1)^jc_j=\sum_{h=1}^N[b_h].
\]
Apply $[a]+[b]\equiv[a+b]+p[\kappa(a,b)]\pmod{p^2}$ successively. If the first $h-1$ terms sum to $[s_{h-1}]+p[d_{h-1}]$, adding $[b_h]$ gives
$[s_h]+p[d_{h-1}+\kappa(s_{h-1},b_h)]$ modulo $p^2$. Since $s_N=0$, the stated identity follows. As $c_j\in\mathcal O_d$ and $c_j/p\in\mathcal O_d$, the right side belongs to the appropriate residue subfield.
\end{proof}

\begin{proposition}[Successive digits of cyclotomic coefficients and fixed fields]
\label{revadic:prop-all-digits}
Let $\widetilde M_C(X)$ be a canonical cyclotomic factor, with $C$ stable under multiplication by $p^d$, and take $d$ to be a permissible degree for its fixed coefficient field. Each coefficient $c_j\in\mathcal O_d$ has a unique expansion
\[
 c_j=\sum_{h\ge0}p^h[a_{j,h}],\qquad a_{j,h}\in k_d.
\]
Consequently,
\[
 v_p(c_j)=\min\{h:a_{j,h}\ne0\},\qquad
 v_p(c_j)=+\infty\iff a_{j,h}=0\ \text{for every }h.
\]
These digits can be computed successively by the Witt operations already introduced in Chapter~5. Define
\[
 B_j=(-1)^j\sum_{\substack{I\subseteq C\\|I|=j}}
                   \biggl[\prod_{i\in I}\alpha_i\biggr]
 \quad(\alpha_i=\bar\zeta_m^i).
\]
Then $B_j=c_j$. Repeating the operations ``take the residue, subtract its Teichmüller representative, and divide by $p$'' yields $a_{j,0},a_{j,1},\ldots$ in order. In particular, if the first $h$ digits vanish and $a_{j,h}\ne0$, then $v_p(c_j)=h$ exactly.
\end{proposition}
\begin{proof}
Stability of $C$ under $p^d$ means that $\phi^d$ permutes its roots and fixes the coefficients. The unramified Galois correspondence gives $c_j\in\mathcal O_d$. The Teichmüller expansion in a complete discrete valuation ring exists uniquely; its digits lie in $k_d$. The elementary-symmetric-function formula and multiplicativity of Teichmüller lifts give $B_j=c_j$. The valuation is the position of the first nonzero digit.
\end{proof}

There are two different kinds of vanishing: a coefficient may have several zero leading digits and eventually become nonzero, or an exact symmetry may force every digit to vanish. Finite layers of carries detect the former; exact root-set structure often explains the latter. In particular, an additive-translation stabilizer of field elements cannot, \emph{by itself}, force the corresponding characteristic-zero coefficient to vanish: Teichmüller lifting does not preserve addition.

\section{Fourier transforms and coefficient fields of Gauss sums}
\label{v135:sec-fourier}
\subsection{Unramified Fourier transforms and divisibility levels}
\label{revadic:subsec-fourier}
Let $m\ge1$ with $p\nmid m$, choose the integer ring $\mathcal O$ of an unramified extension containing $\mu_m$, and let $\zeta$ be a primitive $m$th root of unity. For $f=(f_x)_{x\in\mathbb Z/m\mathbb Z}\in\mathcal O^m$, set
\[
 \widehat f(a)=\sum_{x\bmod m}f_x\zeta^{ax},\qquad a\bmod m.
\]

\begin{theorem}[Fourier transforms commute with reduction, precision, and Frobenius]
\label{revadic:thm-fourier}
The transform is an $\mathcal O$-module automorphism, with inverse
\[
 f_x=\frac1m\sum_{a\bmod m}\widehat f(a)\zeta^{-ax}.
\]
For every $s\ge1$ it commutes with reduction modulo $p^s$. Modulo $p$ it is the cyclic-group Fourier transform with values in the corresponding finite field. In particular,
\[
 \widehat f(a)\in p^s\mathcal O
 \iff\sum_x(f_x\bmod p^s)(\zeta\bmod p^s)^{ax}=0
 \quad\text{in }\mathcal O/p^s\mathcal O.
\]
If every $f_x\in\mathbb Z_p$, then
\[
              \phi(\widehat f(a))=\widehat f(pa).
\]
Thus Fourier coordinates on a $p$-cyclotomic coset form a conjugate family in the same unramified extension. Their reductions reproduce the cyclotomic orbit description of Chapter 5, while successive divisibility levels record additional arithmetic information.
\end{theorem}
\begin{proof}
Since $p\nmid m$, $m$ is a unit in $\mathcal O$. The geometric-sum identity
\[
 \sum_{a\bmod m}\zeta^{a(x-y)}=m\delta_{xy}
\]
gives the inverse formula. All expressions use only ring operations in $\mathcal O$ and therefore commute with every reduction. Divisibility is precisely vanishing in the quotient ring. Finally $\phi$ fixes $\mathbb Z_p$ and takes $\zeta$ to $\zeta^p$; applying it termwise proves the formula.
\end{proof}

\begin{remark}[Complementary uses of the two directions]
Finite-field Fourier support or cosets first identify zeros of $\widehat f(a)\bmod p$; Witt carries can then detect the first nonzero digit. In the reverse direction, an exact Fourier vanishing or a proved higher valuation in characteristic zero immediately yields a finite-field relation on reduction. This Fourier transform lives entirely in the unramified layer. Gauss sums below also involve $p$th roots of unity, which must be treated separately.
\end{remark}

\subsection{Defining sets and permutation symmetries}
\label{v140:sec-code-symmetry}
Let $q=p^f$, $N\ge2$, and $p\nmid N$. Choose an unramified integer ring $\mathcal O\supseteq\mathcal O_f$ containing $\mu_N$, with primitive root $\zeta$. Let $S\subseteq\mathbb Z/N\mathbb Z$ be a $q$-stable spectral support, with defining set $Z=S^c$. Fourier inversion gives the matrix
\begin{equation}\label{v140:eq-integral-projector}
 (\Pi_S)_{ij}=\frac1N\sum_{a\in S}\zeta^{a(j-i)}
 \quad(0\le i,j<N).
\end{equation}
Its entries are fixed by $\phi^f$, hence belong to $\mathcal O_f$, and $\Pi_S^2=\Pi_S$. Write
\[
 R_s=\mathcal O_f/p^s\mathcal O_f,\qquad
 \mathcal C_s(S)=\Pi_S R_s^N,\qquad
 \mathcal C_\infty(S)=\Pi_S\mathcal O_f^N.
\]
These codes are free modules of rank $|S|$; their Fourier coordinates vanish precisely as required on $Z$. Freeness and rank also follow directly from the canonical coprime factors of $X^N-1$ and CRT. Permutation automorphisms below allow only coordinate permutations, excluding coordinatewise multiplication by units and automorphisms of the alphabet.

\begin{proposition}[Exact preservation depth from a spectral projector]
\label{v140:prop-permutation-depth}
For a coordinate permutation $\sigma$ with matrix $P_\sigma$, define
\[
 B_\sigma=(I-\Pi_S)P_\sigma\Pi_S,\qquad
 \nu_S(\sigma)=\min_{i,j}v_p((B_\sigma)_{ij}),
\]
with value $+\infty$ for the zero matrix. Then
\begin{equation}\label{v140:eq-permutation-depth}
 \sigma\in\operatorname{PAut}(\mathcal C_s(S))
 \iff s\le\nu_S(\sigma),\qquad
 \sigma\in\operatorname{PAut}(\mathcal C_\infty(S))
 \iff B_\sigma=0.
\end{equation}
Consequently, the permutation groups decrease with precision and eventually equal the permutation group of the integral code.
\end{proposition}
\begin{proof}
For an idempotent projector, $\ker(I-\Pi_S)=\operatorname{im}\Pi_S$. Hence $B_\sigma\equiv0\pmod{p^s}$ exactly when $P_\sigma$ maps the image module into itself. Since $P_\sigma$ has finite order, inclusion is equality; the same argument applies over the complete integer ring. The first equivalence is entrywise divisibility. There are only finitely many coordinate permutations. Above the maximum of the finite values $\nu_S(\sigma)$, only those with $B_\sigma=0$ remain. If there are no finite values, all levels already agree.
\end{proof}

\begin{proposition}[Multipliers and the common permutation group persist at every precision]
\label{v140:prop-common-symmetry}
For $u\in(\mathbb Z/N\mathbb Z)^\times$, the coordinate multiplier $i\mapsto ui$ preserves $\mathcal C_s(S)$ at any precision if and only if $uS=S$. Thus the multiplier stabilizer at every level is $H(S)$.

For each $s\ge1$, the group of permutations preserving every cyclic code over $R_s$ individually equals $\mathcal G_q(N)$ from Chapter~17. It is also the intersection over just the canonical free lifts. Its criterion remains
\begin{equation}\label{v140:eq-common-lift}
 \sigma(i)-\sigma(j)\in\langle q\rangle(i-j)
 \qquad(i,j\in\mathbb Z/N\mathbb Z).
\end{equation}
\end{proposition}
\begin{proof}
A multiplier permutes Fourier frequencies, proving sufficiency over the integer ring; reduction modulo $p$ gives necessity from the spectral support.

A permutation preserving all canonical free lifts reduces to an element of $\mathcal G_q(N)$, so Theorem~\ref{v73:thm-common-permutation-orbits} gives~\eqref{v140:eq-common-lift}. Conversely, under this condition, for each $q$-orbit $O$ the exact orbit sum
\[
 h_O(z)=\sum_{a\in O}\zeta^{-az}
\]
is constant on the $q$-orbit of $z$. Thus $h_O(\sigma(i)-\sigma(j))=h_O(i-j)$ holds in the integer ring. The permutation commutes with every $\Pi_O$ and preserves all canonical free codes at every level.

Finally, the CRT components of $R_s[X]/(X^N-1)$ are unramified local rings, and each ideal is a $p^e$ multiple of its component, with $0\le e\le s$. Indeed, an element of least valuation in a nonzero ideal is $p^e$ times a unit and generates that ideal. Every cyclic code is therefore a direct sum of modules $p^{e_O}\operatorname{im}\Pi_O$. An $R_s$-linear permutation commuting with all projectors preserves each of these codes as well.
\end{proof}

A single residue code may have additional permutations. For the following family, the first carry determines all of them that persist.

\begin{theorem}[Permutation groups with all nonzero field elements as check columns]
\label{v140:thm-code-group}
Let $n\ge2$, $E=\mathbb F_{q^n}$, and $N=q^n-1$. Choose a primitive element $\alpha\in E$ and put $\theta=[\alpha]$. For $s\ge1$, define
\begin{equation}\label{v140:eq-full-column-code}
 \mathcal H_s=\left\{c\in R_s^N:
        \sum_{i=0}^{N-1}c_i\theta^i=0
        \text{ in }\mathcal O_{fn}/p^s\mathcal O_{fn}\right\}.
\end{equation}
This is the canonical cyclic code with defining set $C_q(1)$ and generator
\[
 g_s(X)=\prod_{j=0}^{n-1}(X-\theta^{q^j})\pmod{p^s}.
\]
It is a free $R_s$-module of rank $N-n$. Label coordinate $i$ by $\alpha^i\in E^\times$. Then
\begin{align}
 \operatorname{PAut}(\mathcal H_1)&=\operatorname{GL}_{\mathbb F_q}(E),\label{v140:eq-residue-group}\\
 \operatorname{PAut}(\mathcal H_s)&=
 \{x\mapsto ax^{q^j}:a\in E^\times,\ 0\le j<n\}
 \quad(s\ge2).\label{v140:eq-lifted-group}
\end{align}
The second group is also the permutation group of the complete integral code. It has order $n(q^n-1)$ and, in exponent coordinates, is generated by the cyclic shift and the Frobenius multiplier.
\end{theorem}
\begin{proof}
The first $n$ check columns $1,\theta,\ldots,\theta^{n-1}$ form an integral basis, so the check map $h_s:R_s^N\to\mathcal O_{fn}/p^s\mathcal O_{fn}$ is surjective. The monic polynomial $g_s$ has degree $n$. Divide a code polynomial by $g_s$; its remainder vanishes at $\theta$, so independence of the integral basis forces the remainder to be zero. The kernel is therefore $(g_s)$, with rank $N-n$ by its coprime complementary factor, or directly by monic division.

If $P_\sigma$ preserves the kernel, it induces a unique $R_s$-linear automorphism $L$ of the quotient with $h_sP_\sigma=Lh_s$. Evaluation on standard coordinate vectors gives $L(\theta^i)=\theta^{\sigma(i)}$. Conversely, a linear automorphism permuting these check columns gives a coordinate permutation preserving the kernel.

At $s=1$, the columns run through all nonzero elements of $E$, giving exactly $\operatorname{GL}_{\mathbb F_q}(E)$. At $s\ge2$, they run through all nonzero Teichmüller representatives. Theorem~\ref{v140:thm-teich-rigidity} forces $L=[a]\phi^{fj}$ and shows that all such maps work. Writing $a=\alpha^b$, the coordinate permutation is $i\mapsto b+q^ji\pmod N$. Distinct pairs $a,j$ give distinct permutations, proving the order formula. Any permutation of the integral code preserves its reduction modulo $p^2$, so it belongs to this list. Conversely, every listed map acts over the complete integer ring.
\end{proof}

Perfect trace duality also gives a uniform description of the dual check kernel:
\[
 \mathcal H_s^\perp=\left\{
 \bigl(\Tr_{A_s/R_s}(b\theta^i)\bigr)_{i=0}^{N-1}:b\in A_s\right\},
 \qquad A_s=\mathcal O_{fn}/p^s\mathcal O_{fn}.
\]
Indeed, the dual of the check map sends $b$ to this trace vector, and its image is the orthogonal complement of the kernel. Surjectivity of the check map makes the dual map injective. Taking the orthogonal complement twice recovers the original code, so these trace codes have the same permutation group.

\begin{example}[Seventh roots and the disappearance of additional permutations]
\label{v140:ex-hamming-symmetry}
Take $q=2,n=3$, selecting the residue orbit of seventh roots by $\alpha^3+\alpha+1=0$. Then $\mathcal H_1$ is the binary $[7,4,3]$ Hamming code, and
\[
 |\operatorname{PAut}(\mathcal H_1)|=(8-1)(8-2)(8-4)=168,
 \qquad |\operatorname{PAut}(\mathcal H_s)|=7\cdot3=21\quad(s\ge2).
\]
Modulo four, its canonical generator is $X^3+2X^2+X-1$. It reduces to $X^3+X+1$, and its product with $X^3-X^2+2X-1$ equals $\Phi_7(X)$ modulo four. This is also the canonical factor supplied by the seventh-root periods; see Example~\ref{v132:ex-seven-period} and~\cite{CalderbankSloane1995Padic}.

An explicit permutation that does not persist modulo four is induced by the $\mathbb F_2$-linear bijection
\[
 g(1)=1,\qquad g(\alpha)=\alpha,\qquad g(\alpha^2)=\alpha^2+1.
\]
Indeed,
\[
 g(\kappa(1,\alpha))=\alpha^2+\alpha+1
 \ne\alpha^2+\alpha=\kappa(g(1),g(\alpha)).
\]
It preserves the residue code but not its canonical lift. Strict reduction also occurs in odd characteristic: for $q=3,n=2$, the group order drops from $(9-1)(9-3)=48$ to $8\cdot2=16$.
\end{example}

\begin{corollary}[Permutation groups of canonical codes with subgroup check columns]
\label{v144:cor-subgroup-code}
Let $N\mid q^n-1$ and $\ord_N(q)=n\ge2$. Choose $\alpha\in\mathbb F_{q^n}$ of order $N$, put $\theta=[\alpha]$, and define
\[
 \mathcal C_s(H)=\left\{c\in R_s^N:
                 \sum_{i=0}^{N-1}c_i\theta^i=0\text{ in }A_s\right\}.
\]
This is the canonical free cyclic code of rank $N-n$ with defining set $C_q(1)$. Under the mixed-exponent separation hypothesis of Corollary~\ref{v144:cor-partial-rigidity}, for every $s\ge2$ and over the complete integer ring,
\begin{equation}\label{v144:eq-subgroup-code-group}
 \operatorname{PAut}(\mathcal C_s(H))
 =\{i\mapsto b+q^ji\pmod N:b\in\mathbb Z/N\mathbb Z,\ 0\le j<n\}.
\end{equation}
In general, a permutation of the residue code preserves the canonical code modulo $p^2$ exactly when its induced map $g$ on the check quotient satisfies~\eqref{v144:eq-partial-lift}.
\end{corollary}
\begin{proof}
The first $n$ check columns form an integral basis, so the check map is surjective. Freeness of its kernel, the defining set, and the linear automorphism induced on the check quotient by a kernel-preserving permutation follow from the monic-division and quotient arguments in Theorem~\ref{v140:thm-code-group}. Here the quotient automorphism need only preserve $[H]$. Corollary~\ref{v144:cor-partial-rigidity} gives~\eqref{v144:eq-subgroup-code-group}; without separation, Proposition~\ref{v144:prop-partial-lift} supplies the exact two-digit criterion.
\end{proof}

Example~\ref{v144:ex-eighty-five} gives a length-85 rank-77 code modulo four and a length-13 rank-10 code modulo nine, with permutation-group orders 680 and 39. Example~\ref{v144:ex-five-permutations} shows that extra symmetries can persist at every precision: fifth roots give a repetition code, and ninth roots give a coordinate direct sum of three length-three repetition codes.

\subsection{Gauss sums and ramified extensions}
The $m$th roots above have $p\nmid m$ and therefore lie in unramified extensions. A standard additive character, however, uses a primitive $p$th root $\xi_p$. For $p>2$, the element $\xi_p-1$ is a ramification parameter, with
\[
 v_p(\xi_p-1)=\frac1{p-1}.
\]
Thus defining these Gauss sums generally requires adjoining $p$th roots of unity; one cannot assume that all summands lie in an unramified integer ring. Let $q=p^f$ and fix a Teichmüller character $\omega:\mathbb F_q^\times\to\mu_{q-1}$ and the additive character $\psi(x)=\xi_p^{\Tr_{\mathbb F_q/\mathbb F_p}(x)}$. Define
\[
 G(\omega^{-a},\psi)=\sum_{x\in\mathbb F_q^\times}\omega(x)^{-a}\psi(x),
 \qquad 1\le a\le q-2.
\]

\begin{corollary}[Complex absolute values and digit-sum valuations of Gauss sums]
\label{revadic:thm-gauss}
For a nontrivial multiplicative character $\chi$,
\[
 G(\chi,\psi)G(\chi^{-1},\psi)=\chi(-1)q,
 \qquad |G(\chi,\psi)|_{\mathbb C}=\sqrt q.
\]
For the fixed $p$-adic embedding and $1\le a\le q-2$, one also has
\begin{equation}\label{revadic:eq-stick}
 v_p\bigl(G(\omega^{-a},\psi)\bigr)=\frac{s_p(a)}{p-1},
\end{equation}
where $s_p(a)$ is the base-$p$ digit sum.
\end{corollary}
\begin{proof}
The first two identities follow from character orthogonality in Chapter~\ref{chap:characters}. Equation~\eqref{revadic:eq-stick} is the Stickelberger valuation formula quoted in Theorem~\ref{v132:thm-gauss-digit-valuation}\cite{GrossKoblitz1979,Robert2001GrossKoblitz}.
\end{proof}

\begin{proposition}[Unramified Frobenius invariance of Gauss sums]
\label{revadic:prop-gauss-invariance}
Let $q=p^f$, take the Teichmüller character $\omega$ in $K_f=\mathbb Q_p(\mu_{q-1})$, and put $L=K_f(\xi_p)$. Extend the unramified Frobenius $\phi$ to $L$ by requiring $\phi(\xi_p)=\xi_p$. Then, for every integer $a$,
\[
 \phi\bigl(G(\omega^{-a},\psi)\bigr)
       =G(\omega^{-ap},\psi)
       =G(\omega^{-a},\psi).
\]
In particular $G(\omega^{-a},\psi)\in\mathbb Q_p(\xi_p)$. For $p>2$ the extension $\mathbb Q_p(\xi_p)/\mathbb Q_p$ is totally ramified of degree $p-1$; for $p=2$, $\xi_2=-1\in\mathbb Q_2$ and this extra extension is trivial.
\end{proposition}
\begin{proof}
For $p>2$, the polynomial $\Phi_p(1+T)=((1+T)^p-1)/T$ is Eisenstein over both $\mathbb Z_p$ and $\mathcal O_f$. Hence $[K_f(\xi_p):K_f]=p-1=[\mathbb Q_p(\xi_p):\mathbb Q_p]$, so the two fields are linearly disjoint. The same polynomial gives $v_p(\xi_p-1)=1/(p-1)$. Thus $\phi$ extends uniquely while fixing $\xi_p$, and
\[
 \phi(G(\omega^{-a},\psi))
   =\sum_{x\ne0}\omega(x)^{-ap}\psi(x).
\]
Make the substitution $y=x^p$. Finite-field Frobenius is bijective, and $\Tr(y)=\Tr(x)$, hence $\psi(y)=\psi(x)$. This gives the second equality. The fixed field of $\phi$ in $L$ is $\mathbb Q_p(\xi_p)$, proving the claim. The case $p=2$ follows directly.
\end{proof}

The multiplicative character takes values in an unramified extension, but Frobenius invariance of the full sum places the Gauss sum in $\mathbb Q_p(\xi_p)$. Its valuation can still be nonzero and is given by~\eqref{revadic:eq-stick}.

Example~\ref{v138:ex-seven-gauss} writes this phenomenon as $g_1(8)=2\eta'$ and $g_3(8)=2\eta$: both complex absolute values are $\sqrt8$, but the dyadic valuations are one and two.

\section{From point-counting recurrences to Newton slopes}
\label{revadic:sec-newton}

Character orthogonality in Chapter~18 can turn curve point counts into finite-field sums, but their complex absolute values do not generally determine $p$-adic valuations of Frobenius roots. The preceding unramified Fourier transform lifts multiplicative symmetry, while Gauss sums show why additive characters also introduce ramification. Point-count recurrences determine the $L$-polynomial; the lower convex hull of its coefficient valuations determines the Newton slopes.

Let $P(T)=\sum_{i=0}^d b_iT^i$ with $b_0b_d\ne0$. Form the lower convex hull of the points $(i,v_p(b_i))$, with zero coefficients assigned infinite height. The slopes of its edges, counted with their horizontal lengths, equal the valuations of the reciprocal roots $\lambda_j$, where $P(T)=b_0\prod_j(1-\lambda_jT)$. The product formula below proves this interpretation of its slopes.

\begin{proposition}[The Newton polygon of a product]
\label{revadic:prop-valuation-product}
Let $u_1,\ldots,u_r$ be nonzero elements of a finite $p$-adic extension, and set
\[
 P(T)=\prod_{i=1}^r(1-u_iT)=\sum_{j=0}^r c_jT^j.
\]
Then
\[
 v_p(c_j)\ge\min_{|I|=j}\sum_{i\in I}v_p(u_i).
\]
Equality holds if the $j$-element subset realizing this minimum is unique. In particular, when the valuations $v_p(u_i)$ are known, the Newton slopes of $P$, with multiplicity, are precisely those valuations. 
\end{proposition}
\begin{proof}
Expand $c_j=(-1)^j\sum_{|I|=j}\prod_{i\in I}u_i$. The nonarchimedean inequality gives the bound, with equality for a unique minimal term. Order $w_i=v_p(u_i)$ so that $w_1\le\cdots\le w_r$, and set $s_j=\sum_{i=1}^jw_i$. Every point $(j,v_p(c_j))$ lies on or above $(j,s_j)$. At $j=r$ and at every index with $w_j<w_{j+1}$, the minimal subset is unique, so equality holds. These are exactly the vertices of the convex polygonal line with successive slopes $w_1,\ldots,w_r$, proving the assertion.
\end{proof}

\subsection{Arbitrary smooth projective curves: the first \texorpdfstring{$g$}{g} point counts determine all slopes}
\label{revadic:subsec-curve-general}

Rationality, integrality, and the functional equation of a curve zeta function give a finite reconstruction problem. Formal logarithms and Newton identities reduce it to an integer recurrence.

\begin{theorem}[Genus $g$ and the first $g$ extension-field counts determine the Newton polygon]
\label{revadic:thm-genus-counts}
Let $X/\mathbb F_q$ be a smooth projective geometrically connected curve of genus $g\ge1$, where $q=p^f$, and let $N_r=\#X(\mathbb F_{q^r})$. By rationality and the functional equation for curve zeta functions\cite[pp.~1--2]{Milne2015RH}, write
\[
 Z_X(T)=\exp\left(\sum_{r\ge1}\frac{N_r}{r}T^r\right)
        =\frac{P_X(T)}{(1-T)(1-qT)},\qquad
 P_X(T)=\sum_{j=0}^{2g}b_jT^j\in\mathbb Z[T].
\]
Then $b_0=1$ and, for $1\le j\le g$, the following are \emph{integer equalities}:
\begin{equation}\label{revadic:eq-curve-rec}
 \boxed{\displaystyle
 j b_j=\sum_{r=1}^{j}(N_r-1-q^r)b_{j-r}.}
\end{equation}
The functional equation determines the other coefficients:
\begin{equation}\label{revadic:eq-curve-symmetry}
 \boxed{\displaystyle b_{2g-j}=q^{g-j}b_j\quad(0\le j\le g).}
\end{equation}
Thus \emph{$q,g,N_1,\ldots,N_g$ suffice} to recover $P_X(T)$ exactly and therefore all Newton slopes from the lower convex hull of the coefficient valuations. With the normalization $v_q=v_p/f$, the $2g$ slopes, counted with multiplicities, are symmetric about $1/2$.
\end{theorem}
\begin{proof}
Take the formal logarithm of $Z_X(T)$ and subtract the logarithms of $(1-T)^{-1}$ and $(1-qT)^{-1}$. This gives
\[
 \log P_X(T)=\sum_{r\ge1}\frac{N_r-1-q^r}{r}T^r.
\]
Differentiate, multiply by $T$, and compare the coefficients of $T^j$ to obtain \eqref{revadic:eq-curve-rec}. The classical functional equation
$P_X(T)=q^gT^{2g}P_X(1/(qT))$ gives \eqref{revadic:eq-curve-symmetry}. Starting from $b_0=1$, the recurrence computes $b_1,\ldots,b_g$, and symmetry supplies the remaining coefficients. These are exact integer operations, even when $p\mid j$; only computations at fixed $p$-adic precision require extra guard digits. The Newton polygon product theorem and functional equation imply the statements about slopes and symmetry.
\end{proof}

\begin{corollary}[The first $g$ point counts and unit reciprocal roots]
\label{revadic:cor-unit-roots}
Retain the notation of Theorem~\ref{revadic:thm-genus-counts}. Let $u(X)$ denote the number, with multiplicities, of reciprocal roots of $P_X(T)$ having $p$-adic valuation zero. Then
\begin{equation}\label{revadic:eq-unit-root-count}
 \boxed{\displaystyle
 u(X)=\deg\bigl(P_X(T)\bmod p\bigr)
      =\max\bigl(\{0\}\cup\{1\le j\le g:b_j\not\equiv0\pmod p\}\bigr).}
\end{equation}
The degree of a nonzero constant polynomial is taken to be zero. Thus $q,g,N_1,\ldots,N_g$ already determine the multiplicity of slope zero, without determining all the slopes. In the theory of Jacobians this invariant is the $p$-rank of the curve.
\end{corollary}
\begin{proof}
The polynomial $Y^{2g}P_X(Y^{-1})$ is monic with integer coefficients. Its roots, the reciprocal roots of $P_X$, are therefore algebraic integers with nonnegative $p$-adic valuations. By the Newton polygon, the horizontal segment of slope zero has length equal to the largest index $j$ for which $b_j$ is a $p$-adic unit, namely $\deg(P_X\bmod p)$. Equation~\eqref{revadic:eq-curve-symmetry} shows that all $b_j$ with $j>g$ are divisible by $p$, so this largest index is at most $g$. The first $g$ coefficients come from the point-count recurrence, proving the formula. Its interpretation as the $p$-rank of the Jacobian is standard; the argument itself does not require that interpretation.
\end{proof}

Once the point-count recurrence of Chapter 18 gives $P_X$, the degree $\deg(P_X\bmod p)$ counts the unit reciprocal roots, without computing the positive slopes.

\begin{example}[Complete reconstruction in genus two]
\label{revadic:ex-genus-two}
For $g=2$, put $A_r=N_r-1-q^r$. Then
\[
 b_1=A_1,\qquad b_2=\frac{A_1^2+A_2}{2},
\]
and
\begin{equation}\label{revadic:eq-genus-two}
 P_X(T)=1+b_1T+b_2T^2+qb_1T^3+q^2T^4.
\end{equation}
Thus just two extension-field point counts determine all four slopes. In characteristic two, the division by $2$ must be performed over the integers, not interpreted as division inside the residue field.
\end{example}

The required inputs are $N_1,\ldots,N_g$. Determining these counts from a curve equation may remain the main arithmetic problem. Other formal exponential series require their own rationality, degree, and functional equation.

\subsection{Cyclotomic Frobenius orbits as blocks of Newton slopes}
\label{revadic:subsec-orbit-slopes}

There is also a purely linear-algebraic bridge from the orbit theory of Chapter~5 to the slopes considered in Chapter~18. Let $E$ be a finite $p$-adic extension containing $\mu_m$, where $p\nmid m$, and let $V$ be a finite-dimensional $E$-vector space with $S\in\operatorname{GL}(V)$, $S^m=1$. Since $m$ is a unit, we have the eigenspace decomposition $V=\bigoplus_{a\bmod m}V_a$, where $S$ acts on $V_a$ as $\zeta_m^a$. If an invertible $E$-linear operator $F$ permutes these spaces by $a\mapsto qa$ (or $a\mapsto q^{-1}a$), it has a block-cyclic structure indexed by cyclotomic orbits.

\begin{proposition}[Characteristic factors and slopes of an orbit block]
\label{revadic:prop-orbit-slope}
Let $D=(a,qa,\ldots,q^{h-1}a)$ be an orbit of length $h$, with $q$ prime to $m$. Suppose $F(V_{q^ia})=V_{q^{i+1}a}$, all these spaces having the same dimension, and put $W_D=\bigoplus_{i=0}^{h-1}V_{q^ia}$. Then
\begin{equation}\label{revadic:eq-orbit-slope}
 \boxed{\displaystyle
 \det(1-TF\mid W_D)=\det(1-T^hF^h\mid V_a).}
\end{equation}
If the eigenvalues of $F^h\mid V_a$ are $\beta_1,\ldots,\beta_t$ (with multiplicities), the valuations of the eigenvalues of $F\mid W_D$ are
\[
 \frac{v_p(\beta_1)}h,\ldots,\frac{v_p(\beta_t)}h,
 \quad\text{each repeated $h$ times}.
\]
Taking the union over all orbit blocks gives the full Newton-slope multiset of $F$.
\end{proposition}
\begin{proof}
Choose a basis of $V_a$ and transport it by $F,F^2,\ldots,F^{h-1}$ to the other spaces. Let $B$ be the matrix of $F^h|_{V_a}$ in the original basis. If $h>1$, the nonzero blocks of $F$ are identity blocks on the subdiagonal and $B$ in the upper-right corner; for $h=1$, the matrix is $B$. Eliminating the identity diagonal blocks of $I-TF$ gives
\[
 \det(I-TF)=\det(I-T^hB).
\]
Over an algebraic closure, the right side is $\prod_{j=1}^t(1-\beta_jT^h)$. Each factor has reciprocal roots satisfying $\lambda^h=\beta_j$, all $h$ of valuation $v_p(\beta_j)/h$. No diagonalizability of $B$ is required.
\end{proof}

\begin{remark}[Additional hypotheses needed for curves]
If a curve $L$-polynomial has a determinant realization and a prime-to-$p$ cyclic automorphism group supplies the stated eigenspaces, the proposition applies. Slopes in each orbit block still depend on the eigenvalues of the return operator $F^h$. The hypothesis requires $E$-linearity; a semilinear Frobenius must first be replaced by a suitable iterate with a specified linear realization.
\end{remark}

\subsection{Slopes of Fermat curves}
\label{revadic:subsec-fermat-slopes}

For general curves the orbit-block formula alone does not give explicit slopes. For curves whose point counts can be expressed through Jacobi sums, however, cyclotomic exponents from Chapter~5, character sums from Chapter~18, and unramified lifting combine into an explicit slope formula. Fermat curves provide a family in which the whole argument can be carried out.

For nontrivial multiplicative characters $\rho,\lambda$ of $\mathbb F_q$ (extended by zero at zero), put
\[
 J(\rho,\lambda)=\sum_{x\in\mathbb F_q}\rho(x)\lambda(1-x).
\]
If $\rho\lambda$ is nontrivial, character orthogonality gives
$J(\rho,\lambda)=G(\rho,\psi)G(\lambda,\psi)/G(\rho\lambda,\psi)$.

\begin{theorem}[Fermat orbit factors and return valuations]
\label{revadic:thm-fermat-slopes}
Let $m\ge3$, $p\nmid m$, and $q=p^f$. The Fermat curve
$X_m:X^m+Y^m=Z^m$ has genus $g=(m-1)(m-2)/2$. Let
\[
 \mathcal A_m=\{(a,b):1\le a,b\le m-1,\ a+b\not\equiv0\pmod m\},
\]
and let $D$ run through its orbits under $(a,b)\mapsto(qa,qb)$. For an orbit of length $h$, choose $(a,b)\in D$, set $Q=q^h$, and let $1\le c<m$ satisfy $c\equiv a+b\pmod m$. The integers
\[
 A=\frac{a(Q-1)}m,\qquad B=\frac{b(Q-1)}m,\qquad C=\frac{c(Q-1)}m
\]
are defined even when $m\nmid Q-1$. With the Teichmüller character $\omega_Q$ on $\mathbb F_Q$, put
\[
 J_D=J_{\mathbb F_Q}(\omega_Q^{-A},\omega_Q^{-B}).
\]
This value is independent of the chosen representative, and
\begin{equation}\label{v141:eq-fermat-orbit-factor}
 P_{X_m/\mathbb F_q}(T)=\prod_D(1+J_DT^{|D|}).
\end{equation}
Each orbit contributes the normalized slope
\begin{equation}\label{revadic:eq-fermat-slopes}
 \lambda_D=\frac{s_p(A)+s_p(B)-s_p(C)}{fh(p-1)},
 \qquad\text{with multiplicity }h.
\end{equation}
If $m\mid q-1$, every orbit has length one. Writing
$\chi=\omega_q^{-(q-1)/m}$ gives the split form
\begin{equation}\label{revadic:eq-fermat-factor}
 P_{X_m}(T)=\prod_{(a,b)\in\mathcal A_m}
                  (1+J(\chi^a,\chi^b)T).
\end{equation}
\end{theorem}
\begin{proof}
The equality $q^ha\equiv a$ implies $m\mid a(Q-1)$, and similarly for $b,c$. All three characters are nontrivial. Replacing $(a,b)$ by $(qa,qb)$ changes both characters by the same Frobenius power; substitution in the Jacobi sum leaves it unchanged.

Over $\mathbb F_{q^r}$ the nontrivial characters detecting $m$th powers are indexed by
$1\le a<m$ with $m\mid a(q^r-1)$, namely
$\omega_{q^r}^{-a(q^r-1)/m}$. Expanding the numbers of solutions of $x^m=u$ and $y^m=1-u$ gives the sum of the corresponding Jacobi sums. Terms with $a+b\equiv0$ equal minus the character value at $-1$; the points at infinity cancel these terms and add one. Hence
\[
 N_r=q^r+1+\sum_{\substack{(a,b)\in\mathcal A_m\\q^r(a,b)=(a,b)}}
 J_{\mathbb F_{q^r}}\!\left(
 \omega_{q^r}^{-a(q^r-1)/m},\omega_{q^r}^{-b(q^r-1)/m}\right).
\]
An orbit $D$ appears exactly when $h\mid r$. Its characters are the norm lifts of those over $\mathbb F_{q^h}$: the exponent is multiplied by $(q^r-1)/(q^h-1)$. The Hasse--Davenport formula from Theorem~\ref{thm:hasse-davenport}, applied to the quotient of Gauss sums, gives
\[
 J_{\mathbb F_{q^r}}=(-1)^{r/h-1}J_D^{r/h}.
\]
There are $h$ equal terms from $D$, so
\[
 N_r=q^r+1-\sum_{D:\,|D|\mid r}|D|(-J_D)^{r/|D|}.
\]
Substitution in the exponential definition of the zeta function gives~\eqref{v141:eq-fermat-orbit-factor}.

The reciprocal roots of $1+J_DT^h$ satisfy $\alpha^h=-J_D$. By the Gauss--Jacobi quotient and the Stickelberger valuation formula in Theorem~\ref{v132:thm-gauss-digit-valuation},
\[
 v_p(J_D)=\frac{s_p(A)+s_p(B)-s_p(C)}{p-1}.
\]
Dividing by $fh$ gives~\eqref{revadic:eq-fermat-slopes}. Thus the return scalar for this orbit factor is explicitly $-J_D$, rather than an unspecified value attached to an orbit of length $h$.
\end{proof}

\begin{corollary}[Base extension and normalized slopes]
\label{revadic:cor-fermat-nonsplit}
The multiset of normalized Newton slopes of a smooth projective geometrically connected curve is unchanged by extending the base field from $\mathbb F_q$ to $\mathbb F_{q^r}$.
In particular, for a Fermat curve with $m\ge3$ and $p\nmid m$, take $r=\ord_m(q)$ and $Q=q^r=p^{fr}$. The split formula~\eqref{revadic:eq-fermat-factor} then applies over $\mathbb F_Q$. For all $(a,b)\in\mathcal A_m$, let $1\le c<m$ satisfy $c\equiv a+b\pmod m$. Compute the Jacobi-sum valuations with exponents $a(Q-1)/m,b(Q-1)/m,c(Q-1)/m$ and divide by $fr$. The resulting multiset is the normalized slope multiset over the original field.
\end{corollary}
\begin{proof}
If the reciprocal Frobenius roots over $\mathbb F_q$ are $\alpha_1,\ldots,\alpha_{2g}$, then after base extension to $\mathbb F_{q^r}$ they are $\alpha_i^r$. Hence
\[
 v_Q(\alpha_i^r)=\frac{r\,v_p(\alpha_i)}{fr}
 =v_q(\alpha_i).
\]
The normalized slopes do not change under finite base extension. Since $m\mid Q-1$, the theorem applies over $\mathbb F_Q$.
\end{proof}

\begin{remark}
This explains another role of the intermediate fields of Chapter~\ref{chap:frobenius} and the cyclotomic cosets of Chapter~5: \emph{enlarge the ground field until the multiplicative symmetry splits, while the unramified Galois tower preserves normalized slopes}. For example, the Fermat cubic over $\mathbb F_5$ already has slopes $1/2,1/2$, as the digit-sum calculation over $\mathbb F_{25}$ shows. Extending the base field and increasing Witt precision remain fundamentally different operations.
\end{remark}

\begin{example}[A Fermat cubic in two different characteristics]
\label{revadic:ex-fermat-cubic}
Take $m=3$. For $q=p=7$, the possible values of $A$ are $2$ and $4$, with base-seven digit sums $2$ and $4$. The only ordered pairs are $(a,b)=(1,1),(2,2)$, and \eqref{revadic:eq-fermat-slopes} gives slopes $0,1$: the cubic is ordinary in characteristic seven. For $q=25$ and $p=5$, the possible values of $A$ are $8$ and $16$, each having base-five digit sum $4$. Both ordered pairs yield $v_5(J)=1$, and $f=2$, so both normalized slopes equal $1/2$: this is the supersingular case. As an independent check, direct counts give $\#X_3(\mathbb F_7)=9$ and $\#X_3(\mathbb F_{25})=36$, yielding the polynomials $1+T+7T^2$ and $1+10T+25T^2$, consistent with these slopes. Neither example requires point counts over all extensions.
\end{example}

\begin{corollary}[Slopes as carry frequencies along an exponent orbit]
\label{v141:cor-orbit-carries}
For $(a,b)\in\mathcal A_m$, let $d$ be its orbit length under multiplication by $p$. Let $a_i,b_i$ be the least positive residues of $p^ia,p^ib$ modulo $m$. The corresponding normalized slope is
\begin{equation}\label{v141:eq-carry-slope}
 \lambda(a,b)=\frac1d\sum_{i=0}^{d-1}\mathbf1_{a_i+b_i>m}.
\end{equation}
It is independent of the choice of $q=p^f$. In particular, $p\equiv1\pmod m$ gives slopes zero and one, each with multiplicity $g$. If $p^t\equiv-1\pmod m$ for some $t$, all slopes equal $1/2$.
\end{corollary}
\begin{proof}
Put $L=fh$ in Theorem~\ref{revadic:thm-fermat-slopes}. For any of $a,b,c$, write $r_i$ for the least positive residue of $p^ia$ and use
$p r_i=m d_i+r_{i+1}$. Since $r_L=r_0$, summing gives
\[
 \frac{s_p(a(p^L-1)/m)}{p-1}
       =\sum_{i=0}^{L-1}\frac{r_i}{m}.
\]
Indeed, the $d_i$ are the repeating long-division digits, read in reverse order in the integer $a(p^L-1)/m$. Subtracting the identity for $c$ from those for $a,b$ gives one exactly when $a_i+b_i>m$; equality with $m$ is excluded by $a+b\not\equiv0$. The $p$-orbit period $d$ divides $L$, proving~\eqref{v141:eq-carry-slope}.

Negation sends a carry to a noncarry, so if it occurs within the $p$-orbit the average is $1/2$. When $p\equiv1\pmod m$, the average is a single zero or one; negation pairs the two types. Since $|\mathcal A_m|=2g$, each occurs $g$ times.
\end{proof}

\begin{example}[Equal orbit lengths with different return slopes]
\label{v141:ex-fermat-orbits}
For $m=7,p=q=2$, all thirty admissible pairs lie in orbits of length three. Representatives, grouped by their number of carries in one period, are
\[
\begin{array}{c|l|c|c}
\text{carries}&\text{representatives}&\text{slope}&\text{multiplicity}\\\hline
0&(1,2),(1,4)&0&6\\
1&(1,1),(1,5),(3,2)&1/3&9\\
2&(1,3),(3,1),(3,3)&2/3&9\\
3&(3,5),(3,6)&1&6
\end{array}
\]
For instance, $(1,1)\mapsto(2,2)\mapsto(4,4)$ has one carry, whereas $(1,2)\mapsto(2,4)\mapsto(4,1)$ has none. Their orbit factors have the same degree three, but the valuations of their return scalars are one and zero. These two orbit types already have return fields $\mathbb F_8$; no further extension or higher Witt precision is needed.

Even the common orbit length two does not force slope $1/2$. For $m=8$, the orbit of $(1,1)$ in characteristic three is $(1,1),(3,3)$ and has slope zero. In characteristic seven it is $(1,1),(7,7)$ and has slope $1/2$. The difference is the carry distribution, not the orbit length.
\end{example}

The multiplicative symmetry of the Fermat equation gives a Jacobi-sum factorization, so the digit-sum formula determines all slopes. For other curves, a corresponding character-sum factorization must first be established before applying this valuation formula.

\begin{proposition}[Slope criterion for a second-order recurrence]
\label{revadic:prop-slopes}
Let $q=p^f$, $a\in\mathbb Z$, and $P(T)=1-aT+qT^2=(1-\alpha T)(1-\beta T)$. If $v_p(a)<f/2$, then
\[
 \{v_p(\alpha),v_p(\beta)\}=
 \{v_p(a),f-v_p(a)\}.
\]
If $v_p(a)\ge f/2$ (including $a=0$, with $v_p(0)=+\infty$), both roots have valuation $f/2$.
\end{proposition}
\begin{proof}
The equations $\alpha+\beta=a$ and $\alpha\beta=q$ show that the two valuations sum to $f$. If the valuation of the sum is less than $f/2$, the two valuations cannot coincide; the ultrametric inequality then identifies the smaller one with $v_p(a)$. Otherwise, if one valuation were below $f/2$, the other would exceed $f/2$, forcing the sum to have valuation equal to the smaller one, a contradiction. Both valuations must therefore equal $f/2$.
\end{proof}

\begin{example}[Different slopes over the same finite field]
\label{revadic:ex-curves}
Over $\mathbb F_5$, consider
\[
 E_1:\ y^2=x^3+x+1,\qquad E_2:\ y^2=x^3+1.
\]
Direct counting gives $\#E_1(\mathbb F_5)=9$ and $\#E_2(\mathbb F_5)=6$, hence traces $a_1=-3$ and $a_2=0$. The characteristic polynomials provided by the elliptic-curve point-counting theorem are
\[
 P_1(T)=1+3T+5T^2,\qquad P_2(T)=1+5T^2.
\]
The reciprocal roots of the first have valuations $\{0,1\}$, and those of the second $\{1/2,1/2\}$. Under a complex embedding, all Frobenius eigenvalues have absolute value $\sqrt5$, yet their $5$-adic slopes differ. 
\end{example}

Use $g_1=g_1(8)$ and $g_3=g_3(8)$ from Example~\ref{v138:ex-seven-gauss}, where $v_2(g_1)=1$ and $v_2(g_3)=2$. The polynomial
\[
                 P(T)=(1-g_1T)(1-g_3T)
\]
therefore has slopes $1,2$; its linear coefficient has valuation exactly $1$ and its quadratic coefficient valuation $3$. 

\subsection{Jacobi valuation distributions with a fixed product character}
Cyclic carries can be viewed as two-state transitions between adjacent digits\cite{EvansHollmannKrattenthalerXiang1999}. Fixing the product character lets these transitions record the valuations of all its Jacobi decompositions.

\begin{proposition}[Two-state carry transfer]\label{v150:prop-jacobi-carries}
Let $q=p^f\ge3$, $N=q-1$, and $1\le t\le N-1$, with $t=\sum_{i=0}^{f-1}t_i p^i$.
Fix a compatible $p$-adic embedding and Teichm\"uller character $\omega$, and put
\[
 P_t(z)=\sum_{\substack{1\le a,b\le N-1\\a+b\equiv t\pmod N}}
 z^{v_p(J(\omega^{-a},\omega^{-b}))}.
\]
Normalize $v_p(p)=1$, and set
\[
 M_d(z)=\begin{pmatrix}d+1&(p-1-d)z\\d&(p-d)z\end{pmatrix}
 \qquad(0\le d<p).
\]
Then
\begin{equation}\label{v150:eq-jacobi-carries}
 P_t(z)=\operatorname{tr}\bigl(M_{t_0}(z)\cdots M_{t_{f-1}}(z)\bigr)-2-2z^f.
\end{equation}
Moreover, $P_t(1)=q-3$, $P_{pt\bmod N}(z)=P_t(z)$, and
\begin{equation}\label{v150:eq-carry-reciprocity}
 P_{N-t}(z)=z^fP_t(z^{-1}).
\end{equation}
\end{proposition}
\begin{proof}
All three characters are nontrivial. The Gauss--Jacobi quotient and Theorem~\ref{v132:thm-gauss-digit-valuation} give
\[
 v_p(J(\omega^{-a},\omega^{-b}))
 =\frac{s_p(a)+s_p(b)-s_p(t)}{p-1}.
\]
Write cyclic carries as $c_i\in\{0,1\}$, with $c_f=c_0$ and
$a_i+b_i+c_i=t_i+pc_{i+1}$.
For input and output carries $(0,0),(0,1),(1,0),(1,1)$, the numbers of digit pairs are respectively
$t_i+1,p-1-t_i,t_i,p-t_i$. Weighting each output carry by $z$ gives $M_{t_i}(z)$, and the trace of their product enumerates all closed carry sequences.
Summing the digit equations shows that the valuation above is exactly $\sum_i c_i$.

The matrix model allows integers $a,b$ in $[0,N]$. Since $t\ne0\pmod N$, the excluded endpoints are precisely
$(0,t),(t,0),(N,t),(t,N)$. The first two have no carry, and the last two have $f$ carries, giving the correction $2+2z^f$.
Every remaining pair has a unique cyclic carry sequence, proving~\eqref{v150:eq-jacobi-carries}.
For fixed nonzero $t$, each $a\in\{1,\ldots,N-1\}\setminus\{t\}$ determines $b$, giving $q-3$ pairs.
Simultaneous multiplication by $p$ rotates digits and preserves valuations. Replacing $a,b,t$ by $N-a,N-b,N-t$ complements the digit sums and replaces a valuation by $f$ minus itself, proving~\eqref{v150:eq-carry-reciprocity}.
\end{proof}

\begin{example}[Equal digit sums with different valuation distributions]\label{v150:ex-carry-distribution}
Take $q=32$. The parameters $t=3,5$ both have binary digit sum two; the corresponding characters both have order 31 and Frobenius orbit length five. Yet~\eqref{v150:eq-jacobi-carries} gives
\[
\begin{array}{c|rrrrrr}
\text{valuation}&0&1&2&3&4&5\\\hline
t=3&2&2&5&10&4&6\\
t=5&2&4&5&4&8&6
\end{array}
\]
Each row sums to $29=q-3$. Relative digit positions enter the matrix product, so a digit sum alone cannot recover the joint distribution.
\end{example}

The valuation distribution does not record unit parts. In fact, for every fixed nontrivial product character $\omega^{-t}$,
\begin{equation}\label{v150:eq-jacobi-total-sum}
 \sum_{\substack{1\le a,b\le N-1\\a+b\equiv t\pmod N}}
 J(\omega^{-a},\omega^{-b})=
 \begin{cases}2,&p=2,\\N\omega(2)^t+2,&p\ne2.\end{cases}
\end{equation}
To prove this, first allow $a$ to run through all residues modulo $N$. Character orthogonality requires $x/(1-x)=1$. There is no solution in characteristic two; in odd characteristic the unique solution $x=1/2$ contributes $N\omega(2)^t$.
Delete the two terms $a=0,t$, both equal to $-1$, to obtain the formula.
In characteristic two, the summands can thus have different valuations while their sum is always two. Determining further local cancellation requires their unit information as well.

\section{Constructions and coefficient questions}
\label{revadic:sec-end}

\begin{center}
\begin{tikzpicture}[>=Stealth, every node/.style={align=center,font=\small}]
\node (finite) at (0,0) {$k_n=\mathbb F_{p^n}$};
\node (witt) at (-4,-1.5) {$W(k_n)$};
\node (cycl) at (4,-1.5) {$\mathcal O_{K_n}$};
\node (base) at (4,-3) {$\mathbb Q_p(\zeta_{p^n-1})$};
\draw[->,thick,Jade] (finite) -- node[above left]{Witt lift} (witt);
\draw[->,thick,Jade] (base) -- node[right]{integer ring} (cycl);
\draw[->,thick,Jade] (cycl) -- node[above right]{reduction} (finite);
\draw[<->,thick,InkSlate] (witt) -- node[below]{complete local ring isomorphism} (cycl);
\end{tikzpicture}
\end{center}

The two constructions meet at the level of complete unramified integer rings. Teichmüller lifting turns multiplicative relations among finite-field roots into exact relations among roots of unity. Reduction of the unramified cyclotomic ring recovers the same finite field and Frobenius. Ordinary addition of lifted roots followed by Teichmüller projection recovers residue-field addition; equation~\eqref{v135:eq-finite-projection} realizes this projection at every finite precision, and Witt carries record the difference between the ordinary sum and its projection. The equal-characteristic ring $k_n[[t]]$ preserves both operations on constants, but records a different type of deformation rather than the same mixed-characteristic divisibility levels.

The root orbits of Chapter~5 specify canonical factors and, through subset counts, exterior-power traces, and cyclotomic remainders, distinguish exact zero coefficients from divisibility after reduction. Theorem~\ref{revadic:thm-local-min} shows that excess local cancellation cannot affect every component, while Corollary~\ref{revadic:cor-symmetry-sharp} determines a symmetry condition under which the lower bound is sharp. The Frobenius cyclic modules of Chapter~9 make normal and integral normal bases subject to the same unit-determinant test, and taking invariants commutes with truncating precision. The Fourier relations of Chapter~18 lift throughout the unramified tower; the passage to additive characters, Gauss sums, and Newton slopes requires separate inputs concerning ramification and curve $L$-functions.

The mixed monomials in the first carry also recover multiplication, determining linear maps that preserve all canonical representatives and subspaces that remain additively closed with identical digit sets. Smith exponents of a Frobenius difference control the precision of the low-degree subfield representations of Chapter~13. Canonical spectral projectors detect the preservation depth of permutations for the defining sets of Chapter~17. These criteria distinguish persistence of subfields and common symmetries from the loss of individual low-degree representations or additional permutations. For partial root sets, mixed exponents modulo $X^N-1$ record the carry obstruction: separated exponents give semilinear rigidity, while complete root blocks permit larger stabilizers. Coefficient profiles are first compressed by full fibers and fixed groups; Smith exponents then determine depths, and finite-precision trace coordinates determine positions.

\paragraph{Distribution of local cancellation.}
The norm gives the sum of the valuations of phase values; the generating-order index gives the sum for their differences. At minimal positive excess, Theorem~\ref{v141:thm-minimal-excess} and Proposition~\ref{v141:prop-exceptional-kernels} determine the positions and depths of all exceptional differences, while low-level zero sets still record the position of the values relative to zero. At excess $2g_f$ or $3g_f$, can intersections of multiplication kernels, difference identities, and the congruence-subgroup chain exclude combinations of exceptional orbits that cannot occur together?

\paragraph{Constraints on slopes from orbits.}
Theorem~\ref{revadic:thm-fermat-slopes} obtains the return scalar $-J_D$ from point counts and character sums; Corollary~\ref{v141:cor-orbit-carries} then reads its valuation from exponent carries. For other extension-field sums defined by root sets or multiplicative subgroups, which phase relations express return scalars through Gauss sums, Jacobi sums, or products of known periods? Orbit lengths and coefficient-difference cancellation graphs alone do not supply such a return identity.

\section*{Exercises}
\addcontentsline{toc}{section}{Exercises}
\begin{exercise}[The degree-one boundary case]\label{revadic:ex1}
Show that $\mathbb Q_p$ contains every $(p-1)$st root of unity. Explain why $p=2,n=1$ needs a separate convention for $\operatorname{ord}_{p^n-1}(p)$.
\end{exercise}

\begin{exercise}[Another cyclotomic orbit]\label{revadic:ex2}
For $p=3,m=8$, find the ternary cyclotomic cosets modulo eight and verify irreducibility modulo three of the two quadratic factors in Example~\ref{revadic:ex-eight}.
\end{exercise}

\begin{exercise}[The unramified tower]\label{revadic:ex3}
If $d\mid n$, prove $\mathcal O_d=\mathcal O_n\cap K_d$ and identify the fixed elements of $\phi^d$ in $k_n$ with $k_d$.
\end{exercise}

\begin{exercise}[The root equation at finite precision]\label{revadic:ex11}
For $p=2,n=3$, prove that $x^8=x$ has exactly eight solutions in $\mathcal O_3/2^s\mathcal O_3$ for every $s\ge1$. Determine the precisions at which the solution set is closed under addition.
\end{exercise}

\begin{exercise}[Recovering addition from the root equation]\label{revadic:ex15}
Show that $\tau(u)=\lim_{r\to\infty}u^{p^{nr}}$ is not a ring homomorphism on $\mathcal O_n$, although $\tau(u+v)\equiv u+v\pmod p$. Compute $1\boxplus1$ for $p=2,n=1$ and explain why $\boxplus$ is associative.
\end{exercise}

\begin{exercise}[First correction of the same residual root]\label{revadic:ex24}
Take $p=2$ and $a\in\mathbb F_4$ satisfying $a^2+a+1=0$. Explain why $F_1=X^2+X+1$ and $F_2=X^2+X-1$ have the same simple root modulo two but different integral Hensel lifts, and determine the first digit of the difference.
More generally, for $r\ge1$, let $H_r=X^2+X+1-2^r$. Find the valuation of the distance between its root $z_r$ above $a$ and $[a]$, and give a root formula modulo $2^{2r}$.
\end{exercise}

\begin{exercise}[An explicit trace-dual basis]\label{revadic:ex4}
In Example~\ref{revadic:ex-basis}, compute the trace-dual integral basis of $\zeta_3,\zeta_3^2$ and reduce it modulo two.
\end{exercise}

\begin{exercise}[Lifting all normal elements]\label{revadic:ex5}
Extend Theorem~\ref{revadic:thm-normal} to a relative unramified extension $K_n/K_d$ for $d\mid n$, replacing $\phi$ by $\phi^d$. Let $\theta,z\in\mathcal O_n$, with $\bar\theta=\eta\in k_n$ normal over $k_d$. Does $\theta+pz$ still generate a normal integral basis?
\end{exercise}

\begin{exercise}[The asymmetry of trace and norm]\label{revadic:ex6}
Verify $\Nm(\zeta_3)=1$ and $\Tr(\zeta_3)=-1$ for $K_2=\mathbb Q_2(\zeta_3)$. Prove exact compatibility of Teichmüller norms in any unramified tower and explain why it need not hold for trace.
\end{exercise}

\begin{exercise}[Cyclic-module normality and repeated factors]\label{revadic:ex12}
Let $p=2,n=4$. Show that $T^4-1=(T-1)^4$ over $\mathbb F_2$, and use Proposition~\ref{revadic:prop-cyclic-normal} to explain why the normality test does not require distinct factors of $T^n-1$.
\end{exercise}

\begin{exercise}[An equal-characteristic deformation]\label{revadic:ex8}
In $\mathbb F_3[[t]]$, find through $t^3$ the root $x(t)$ satisfying $x(0)=1$ and $x(t)^2-(1+t)=0$. Explain why there is no unital ring section $\mathbb F_3\to\mathbb Z_3$.
\end{exercise}

\begin{exercise}[Which periodic words come from a given modulus?]\label{v138:ex-word-denominator}
Fix $q=p^f$ and $p\nmid m$. Given a base-$q$ block of length $D$, neither all zero nor all $q-1$, let $E$ be its integer value in low-to-high order.
Prove that its infinite repetition represents some $-a/m$, $1\le a<m$, if and only if the reduced denominator of $E/(q^D-1)$ divides $m$. When this holds, determine $a$ and relate the least period $d$ to $D$ and $C_q(a)$.
\end{exercise}

\begin{exercise}[Digit complements and extension of Gauss sums]\label{v138:ex-digit-complement}
Let $Q=q^d$, $1\le a<m$, and $m\mid a(Q-1)$, and put $E_a=a(Q-1)/m$. Over each $\mathbb F_{Q^h}$, define $g_a(Q^h)$ using compatible Teichmüller characters, the standard additive character, and exponent $a(Q^h-1)/m$.
Prove
\[
 g_a(Q)g_{m-a}(Q)=\omega^{-E_a}(-1)Q,\qquad
 g_a(Q^h)=(-1)^{h-1}g_a(Q)^h.
\]
Then interpret the sum of the valuations for inverse characters and the growth of valuations under extension through the periodic digit blocks. Why can the unnormalized valuation increase while the least digit period remains unchanged?
\end{exercise}

\begin{exercise}[A tower of ternary moduli in characteristic two]\label{v139:ex-three-power}
Let $s\ge1$, $m=3^s$, and $C=C_2(1)$ modulo $m$. Determine the least binary period of $-1/m$, the largest complete equal-difference block size in $C$, and the least period after multiplication by that size. Determine the canonical factor for $C$ and the dimension and minimum distance of the binary cyclic code with defining set $C$.
\end{exercise}

\begin{exercise}[The defining words of the dual code]\label{v139:ex-dual-words}
Let $\varnothing\ne Z\subsetneq\mathbb Z/m\mathbb Z$ be $q$-stable, and use $S_Z$, $M$, and $c$ from Section~\ref{v139:sec-digit-bch}. Determine the encoded defining set of $\mathcal C_Z^\perp$ from $S_Z$, and express its strongest consecutive-root bound for the fixed primitive root using the elements of $Z$. Apply the result to $q=2,m=15$, $Z=C_2(1)\cup C_2(3)$.
\end{exercise}

\begin{exercise}[Two-digit translations of a trace hyperplane]\label{v140:ex-trace-hyperplane}
Let $E=\mathbb F_{2^n}$, $n\ge2$, and $V=\ker\Tr_{E/\mathbb F_2}$. Determine the full translation stabilizer of $\mathcal T_2(V)$ and its order, distinguishing odd and even $n$. When is $\mathcal T_2(V)$ itself an additive subgroup?
\end{exercise}

\begin{exercise}[Reciprocal power bases and precision profiles]\label{v140:ex-reciprocal-profile}
Use Example~\ref{v140:ex-quartic-depth}, put $\gamma=\alpha^{-1}$, and define $\rho_s^\gamma(t,2)$ using the power basis of $[\gamma]=\theta^{-1}$. Determine the profile for every $s\ge1$ and $0\le t\le3$, and compare it with the profile for $\theta$. Does the difference in representation degrees between the residue power bases persist at higher precision?
\end{exercise}

\begin{exercise}[Two fixed points and return degrees]\label{v140:ex-two-point-group}
Use Theorem~\ref{v140:thm-code-group}. For distinct $a,b\in E^\times$, put $d=[\mathbb F_q(b/a):\mathbb F_q]$. Find the order of the subgroup fixing coordinates $a,b$ pointwise at $s\ge2$, and describe the orbits on ordered pairs of distinct coordinates. Compare with the residue code $\mathcal H_1$: which pair orbits split as precision increases?
\end{exercise}

\begin{exercise}[Seventy-three check columns]\label{v144:ex-seventy-three}
Let $\alpha\in\mathbb F_{512}$ have multiplicative order 73, and put $\theta=[\alpha]$. For $s\ge2$, determine all coordinate permutation automorphisms of
\[
 \mathcal C_s=\left\{c\in(\mathbb Z/2^s\mathbb Z)^{73}:\sum_{i=0}^{72}c_i\theta^i=0\right\},
\]
and its free rank. Is the group the same over the complete integer ring?
\end{exercise}

\begin{exercise}[Exact zeros in a valuation profile]\label{revadic:ex7}
If a root set is stable under multiplication by a root of unity of order $e$, show that its root polynomial belongs to $\mathcal O[X^e]$. Use $m=7,p=2$ to illustrate that a zero coefficient in the finite field need not lift to an exact zero.
\end{exercise}

\begin{exercise}[Exterior-power traces of an orbit factor]\label{revadic:ex19}
Take $p=2,m=7,C=\{1,2,4\}$. Regard multiplication by $\zeta_7$ as an endomorphism of a free $\mathbb Z_2$-module of rank three. Use Theorem~\ref{revadic:thm-orbit-exterior} to explain why $\eta$ and $\eta'$ appear respectively in its first and second exterior-power traces.
\end{exercise}

\begin{exercise}[Period counts and the first coefficient digit]\label{revadic:ex16}
For $p=2,m=7,C=\{1,2,4\}$, compute all coefficients using $N_{C,j}(s)$ for $j=1,2,3$. Explain why $c_1$ vanishes modulo two but is nonzero in $\mathbb Z_2$.
\end{exercise}

\begin{exercise}[Cyclotomic remainders and exact zero coefficients]\label{revadic:ex21}
Let $p=2,m=9,d=2,C=\{1,4,7\}$. Compute the remainders of $A_{C,1}(T)$ and $A_{C,2}(T)$ modulo $\Phi_9(T)$, and prove $Q_C(X)=X^3-\zeta_9^3$. Explain why these are exact zeros, not merely local cancellations.
\end{exercise}

\begin{exercise}[An involutory direction outside Frobenius]\label{v141:ex-involution}
In the phase algebra $\mathcal A=\mathbb Z_3^4$ take $\Gamma=\mathbb Z/4\mathbb Z$, $\tau=0$, and $x=(1,2,4,3)$, with directions acting by cyclic translation of coordinates. Determine the index excess, all directional depths and multiplication-kernel dimensions, and the exceptional idempotent. Compare the valuation profiles of $x$ and $x+2$.
\end{exercise}

\begin{exercise}[Zero trace of an exceptional idempotent]\label{v141:ex-period-idempotent}
For the conductor-seventeen algebra in Example~\ref{v141:ex-period-indices}, find the idempotent indicating the zero components of $\bar\eta_0-\bar\eta_2$. Determine its trace and its four trace moments against $\bar\eta_i$. How do these moments recover information lost by its trace alone?
\end{exercise}

\begin{exercise}[A return field smaller than a splitting field]\label{v141:ex-short-return}
For $X_{12}/\mathbb F_5$, determine the orbit factors and normalized slopes associated with $(3,3)$, $(9,9)$, and $(1,7)$, leaving Jacobi sums in exact character-sum form. Which of the return fields fail to contain all twelfth roots of unity?
\end{exercise}

\begin{exercise}[Integral remainders and local factors]\label{revadic:ex23}
Let $p=2,m=7,C=\{1,2,4\}$. Verify that $\Phi_7\bmod2$ has two distinct cubic factors, and show that $A_{C,1}(T)$ has a nonzero integral remainder but reduces to zero at the chosen local root of $T^3+T+1$. Explain why irreducibility cannot be omitted in part (3) of Theorem~\ref{revadic:thm-integral-resolvent}.
\end{exercise}

\begin{exercise}[Global minimum valuation and excess local cancellation]\label{revadic:ex25}
Let $m=7,p=2,C=\{1,2,4\}$, and $j=1$. Compute the minimum coefficient valuation $\nu$ of the integral cyclotomic remainder, compare the valuations at the two modulo-two Frobenius orbits of primitive roots, and verify Theorem~\ref{revadic:thm-local-min}. Explain why the stabilizer criterion does not rule out excess cancellation.
\end{exercise}

\begin{exercise}[A second exterior-power trace for nineteenth roots]\label{revadic:ex22}
Take $p=2,m=19,d=2,C=\langle4\rangle$. Without computing a lifted polynomial, prove $N_{C,2}(0)=0$ and $N_{C,2}(s)=2$ for $s\ne0$, and determine $c_2$ and $v_2(c_2)$.
\end{exercise}

\begin{exercise}[Detecting the seventh period by a Witt carry]\label{revadic:ex10}
Let $\mathbb F_8=\mathbb F_2(\alpha)$, with $\alpha^3=\alpha+1$. For $b_1=\alpha,b_2=\alpha^2,b_3=\alpha^4$, use the characteristic-two identity $\kappa(a,b)^2=ab$ to compute the right side of~\eqref{revadic:eq-first-profile}. Deduce $v_2(\eta)=1$ without using $\eta\eta'=2$.
\end{exercise}

\begin{exercise}[Frobenius cosets for unramified Fourier transforms]\label{revadic:ex13}
Let $p\nmid m$ and $f_x\in\mathbb Z_p$. Prove $v_p(\widehat f(pa))=v_p(\widehat f(a))$ and interpret the restriction this puts on valuations along one $p$-cyclotomic coset.
\end{exercise}

\begin{exercise}[Gauss sums and a Frobenius fixed field]\label{revadic:ex14}
Use $\Tr_{\mathbb F_{p^f}/\mathbb F_p}(x^p)=\Tr(x)$ to show directly that $G(\omega^{-pa},\psi)=G(\omega^{-a},\psi)$. Explain why this does not imply the digit-sum valuation formula.
\end{exercise}

\begin{exercise}[Quadratic Newton slopes]\label{revadic:ex9}
For $P(T)=1-aT+p^fT^2$, give its slopes when $v_p(a)=0$, when $0<v_p(a)<f/2$, and when $v_p(a)\ge f/2$. Explain why $|\alpha|_{\mathbb C}=|\beta|_{\mathbb C}=p^{f/2}$ cannot determine these slopes.
\end{exercise}

\begin{exercise}[Point-count reconstruction in genus two]\label{revadic:ex17}
Given the point counts $N_1,N_2$ of a smooth projective genus-two curve over $\mathbb F_q$, derive all five coefficients of $P_X(T)$ and describe a Newton-polygon computation requiring no further extension-field counts.
\end{exercise}

\begin{exercise}[Detecting unit roots from point counts]\label{revadic:ex20}
For a genus-two curve put $A_r=N_r-1-q^r$. Give a criterion in terms of $A_1,A_2$ for the multiplicity of slope zero in its Newton polygon. Explain why in characteristic two the integer $(A_1^2+A_2)/2$ must be computed before reduction.
\end{exercise}

\begin{exercise}[Digit sums for Fermat curves]\label{revadic:ex18}
Apply Theorem~\ref{revadic:thm-fermat-slopes} to $m=3,q=7$ and $m=3,q=25$. Compute all normalized slopes and verify their symmetry about $1/2$. Identify the two theorem-level inputs needed before applying the digit-sum formula.
\end{exercise}

\begin{exercise}[Root-set addition at finite precision]\label{v135:ex-finite-add}
In $\mathcal O_2/3^5\mathcal O_2$, put $R=\{x:x^9=x\}$. Prove that $x\star y=(x+y)^{81}$ and ordinary multiplication make $R$ a field of nine elements. Give an exponent that still works in $\mathcal O_2/3^6\mathcal O_2$.
\end{exercise}
\begin{exercise}[All homomorphisms between two quotient rings]\label{v135:ex-hom-count}
Let $p=2$. For every $s\ge1$, count the unital ring homomorphisms from $\mathcal O_2/2^s\mathcal O_2$ to $\mathcal O_6/2^s\mathcal O_6$. Describe their action on the roots of $x^4=x$ and determine whether they can be chosen compatibly at all precisions.
\end{exercise}

\appendix

\chapter{Field Extensions and Polynomial Tools}
\label{chap:polynomial-tools}

This appendix treats extensions and polynomial operations over arbitrary fields. The basic concepts of field extensions relate the base field, vector-space coordinates, and algebraic relations;
division with remainder, Bézout's identity, and interpolation provide the computational tools for adjoining roots and decomposing root sets.
The proofs of these general results do not depend on the existence or structure theorems for finite fields.

\section{Field Embeddings and a Specified Base Field}\label{sec:extension-embeddings}

Field homomorphisms preserve addition, multiplication, and the identity element. Every field homomorphism \(\iota:k\to E\) is injective:
if a nonzero \(a\in k\) satisfies \(\iota(a)=0\), then
\[
1=\iota(1)=\iota(aa^{-1})=\iota(a)\iota(a^{-1})=0,
\]
contrary to \(1\ne0\) in a field.
\leavevmode\phantomsection\label{def:field-embedding}\ignorespaces Thus \(\iota\) is an injective field homomorphism, called a \textbf{field embedding};
it gives an isomorphism from \(k\) to the subfield \(\iota(k)\).
Identifying \(k\) with \(\iota(k)\) through this isomorphism allows us to write the inclusion \(k\subseteq E\).
A field equipped with a specified embedding \(\iota:k\to E\) is therefore also regarded as an extension field of \(k\).
Scalar multiplication in its \(k\)-vector space structure is explicitly \(a\cdot x=\iota(a)x\).

\phantomsection\label{def:base-field-morphism}
Let \(\iota:k\to E\) and \(\iota':k\to E'\) be specified embeddings.
A field homomorphism \(\tau:E\to E'\) satisfying \(\tau\circ\iota=\iota'\) is called a \textbf{\(k\)-homomorphism};
if it is also bijective, it is called a \textbf{\(k\)-isomorphism}.
When both images of the base field are denoted by \(k\), this condition says that \(\tau(a)=a\) for every \(a\in k\).
Thus a \(k\)-homomorphism is automatically a \(k\)-linear map.
A \(k\)-isomorphism from a field \(E\) to itself is called a \textbf{\(k\)-automorphism}; these maps form the group \(\Aut_k(E)\) under composition.

For example, in the quotient construction obtained by adjoining a single root \(k[X]/(f)\), the specified embedding is \(a\mapsto\overline a\).
Identifying constant classes in the quotient with elements of the original field through this map gives the inclusion \(k\subseteq k[X]/(f)\) used in the main text.
Similarly, a \(k\)-isomorphism between two models obtained by adjoining a root must map the chosen root to a corresponding root and preserve the specified embedding of the base field.

\section{Bases, Extension Degrees, and Restriction of Scalars}\label{sec:extension-bases}

For a finite tuple \(e_1,\ldots,e_n\) in a field extension \(E/k\), spanning means that every \(x\in E\) is a \(k\)-linear combination of these elements;
linear independence means that \(\sum_i a_ie_i=0\) with \(a_i\in k\) holds only when all \(a_i=0\).
Together, these two conditions are precisely the existence and uniqueness of coordinate representations used to define the extension degree in the main text.
Dimension is determined by linear relations over the base field; the same elements may have different linear relations over different base fields.

\begin{example}[Finite extensions and finite fields]
Every complex number has a unique expression \(a+bi\) with \(a,b\in\mathbb R\).
Thus \(1,i\) is a basis of \(\C/\mathbb R\) and \([\C:\mathbb R]=2\), although both fields have infinitely many elements.
If the base field \(k\) has \(q\) elements, a coordinate tuple of length \(n\) has exactly \(q^n\) possible values,
so every extension \(E/k\) of degree \(n\) satisfies \(|E|=q^n\).
Conversely, if \(E\) is a finite set, the length of every \(k\)-linearly independent tuple is bounded by \(|E|\); successively extending a linearly independent tuple therefore produces a finite basis.
\end{example}

\phantomsection\label{def:restriction-scalars}
Let \(k\subseteq L\), and let \(V\) be an \(L\)-vector space.
Keeping addition in \(V\) and allowing only elements of \(k\) as scalars gives the \(k\)-vector space obtained from \(V\) by \textbf{restriction of scalars}.
If \(u_1,\ldots,u_a\) is a basis of \(L/k\) and \(v_1,\ldots,v_b\) is a basis of \(V\) over \(L\),
the two successive expansions and coefficient comparison in Lemma~\ref{lem:tower} apply equally to \(V\), giving the \(k\)-basis \(u_iv_j\). Hence
\[
\dim_k V=[L:k]\dim_L V.
\]
Taking \(V=E\) gives the tower formula for dimensions. This identity describes the further expansion of each \(L\)-coordinate into several \(k\)-coordinates.

If \(k\subseteq L\subseteq E\) and \([E:k]<\infty\), then \(L\) is a \(k\)-linear subspace of \(E\), so \([L:k]<\infty\).
Also, a \(k\)-basis of \(E\) spans \(E\) over \(L\), so \([E:L]<\infty\).
The tower formula therefore applies and gives \([L:k]\mid[E:k]\).
Thus divisibility of intermediate-field degrees follows from successive basis expansions, independently of the classification of subfields.

\section{Generated Subfields and Algebraic Extensions}\label{sec:generated-algebraic-extensions}

Chapter~\ref{chap:generators} defines \(k(S)\) as the smallest subfield containing \(k\cup S\),
and \(k[\alpha]\) as the set of all polynomial values at \(\alpha\).
For any \(\alpha\in E\), the generated subfield \(k(\alpha)\) also has the following description in terms of fractions.

\begin{proposition}[Fractions in a generated subfield]\label{prop:generated-fractions}
Let \(k\subseteq E\) and \(\alpha\in E\). Then
\[
k(\alpha)=
\left\{\frac{A(\alpha)}{B(\alpha)}:
A,B\in k[X],\ B(\alpha)\ne0\right\}.
\]
More generally, if \(S\subseteq E\), every element of \(k(S)\) belongs to a subfield generated by some finite subset of \(S\).
\end{proposition}
\begin{proof}
Sums, differences, and products of fractions on the right have the same form, and the inverse of a nonzero fraction is obtained by interchanging numerator and denominator.
Thus the right-hand side is a subfield containing \(k\) and \(\alpha\), and hence contains \(k(\alpha)\).
Conversely, every subfield containing \(k\) and \(\alpha\) contains all these fractions, proving equality.

For the second assertion, let \(U\) be the union of all \(k(S_0)\), where \(S_0\) ranges over the finite subsets of \(S\).
If \(x\in k(S_0)\) and \(y\in k(S_1)\), both belong to \(k(S_0\cup S_1)\),
so their sum, difference, product, and quotient when the denominator is nonzero all belong to \(U\).
Hence \(U\) is a subfield containing \(k\cup S\); since every \(k(S_0)\subseteq k(S)\), we obtain \(U=k(S)\).
\end{proof}

\leavevmode\phantomsection\label{def:finitely-generated-extension}\ignorespaces If \(E=k(\alpha_1,\ldots,\alpha_r)\), then \(E/k\) is called a \textbf{finitely generated extension}.
Finite generation restricts only the number of generators; it does not require them to be algebraic over the base field.

\begin{example}[One generator and infinite dimension]
Let \(t\) be an indeterminate, and let \(k(t)\) be the rational function field consisting of fractions \(A(t)/B(t)\) with \(A,B\in k[T]\) and \(B\ne0\).
Then \(k(t)/k\) is generated by the single element \(t\), but \(1,t,t^2,\ldots\) are linearly independent over \(k\), since a polynomial in \(k[t]\) is zero precisely when all its coefficients are zero.
Thus \([k(t):k]=\infty\). Moreover, \(t^{-1}\in k(t)\) whereas \(t^{-1}\notin k[t]\), so \(k[t]\ne k(t)\).
\end{example}

\phantomsection\label{def:algebraic-extension}
If every element of \(E\) is algebraic over \(k\), then \(E/k\) is called an \textbf{algebraic extension}.
An element that is not algebraic over \(k\) is called a \textbf{transcendental element}.
For an algebraic element \(\alpha\), Proposition~\ref{prop:minpoly} has already established
\[
k(\alpha)=k[\alpha],\qquad
[k(\alpha):k]=\deg m_{\alpha,k}.
\]
The key is that a relation of least degree both expresses higher powers in terms of lower powers and excludes nontrivial linear relations among those lower powers.

\begin{proposition}[Finite extensions and algebraic generators]\label{prop:finite-algebraic-generation}
Let \(E/k\) be a field extension.
\begin{enumerate}
\item If \([E:k]=n<\infty\), then \(E/k\) is algebraic, and the degree of the minimal polynomial of every \(\alpha\in E\) divides \(n\).
\item If \(\alpha_1,\ldots,\alpha_r\in E\) are all algebraic over \(k\), then \(k(\alpha_1,\ldots,\alpha_r)/k\) is finite, and
\[
[k(\alpha_1,\ldots,\alpha_r):k]
\le\prod_{i=1}^{r}\deg m_{\alpha_i,k}.
\]
\end{enumerate}
Thus finite extensions are precisely the finitely generated algebraic extensions.
\end{proposition}
\begin{proof}
If \(\dim_k E=n\), the \(n+1\) elements \(1,\alpha,\ldots,\alpha^n\) are linearly dependent,
so some nonzero \(A\in k[X]\) satisfies \(A(\alpha)=0\). This proves that \(\alpha\) is algebraic.
Applying the tower formula to \(k\subseteq k(\alpha)\subseteq E\) and using Proposition~\ref{prop:minpoly} gives the divisibility of degrees.

For the second assertion, put \(K_0=k\) and \(K_i=k(\alpha_1,\ldots,\alpha_i)\).
In \(K_{i-1}[X]\), the polynomial \(m_{\alpha_i,k}\) remains nonzero and has \(\alpha_i\) as a root,
so \(\alpha_i\) is algebraic over \(K_{i-1}\), and its minimal polynomial divides \(m_{\alpha_i,k}\).
Proposition~\ref{prop:minpoly} therefore gives
\[
[K_i:K_{i-1}]=\deg m_{\alpha_i,K_{i-1}}
\le\deg m_{\alpha_i,k}.
\]
Successive applications of the tower formula give the asserted finiteness and bound.
Finally, any finite extension is generated by a finite \(k\)-basis, whose elements are algebraic by the first assertion;
the converse follows from the second assertion.
\end{proof}

\begin{corollary}[Algebraic elements form a subfield]\label{cor:algebraic-elements-field}
In any extension field \(E/k\), the elements algebraic over \(k\) form a subfield of \(E\).
\end{corollary}
\begin{proof}
Every element of the base field is algebraic. If \(\alpha,\beta\) are algebraic over \(k\), then \(k(\alpha,\beta)/k\) is finite.
It contains \(\alpha\pm\beta\) and \(\alpha\beta\), and also \(\alpha/\beta\) when \(\beta\ne0\).
By the preceding proposition, these elements are algebraic over \(k\), proving closure under the field operations.
\end{proof}

Every element of an algebraic extension belongs to a finite extension, but the entire algebraic extension need not be finite.
The algebraic closure of a finite field constructed in Theorem~\ref{thm:closure} provides such an example: it is an increasing union of finite extensions whose degrees are unbounded.

\section{Simple Extensions, Splitting Fields, and Separability}\label{sec:extension-root-models}

If \(\alpha\) is algebraic over \(k\), Proposition~\ref{prop:minpoly} gives
\[
k[X]/(m_{\alpha,k})\simeq k(\alpha),\qquad
\overline{A}\longmapsto A(\alpha).
\]
This isomorphism identifies the field defined by the quotient with the subfield generated by a root in a specified extension field.
Different images of the root determine different embeddings; the precise statement and proof of the extension lemma appear in Lemma~\ref{lem:extend}.
Adjoining all roots of a specified polynomial gives its splitting field as in Definition~\ref{def:splitting-field}.
A splitting field is generated by finitely many algebraic elements and is therefore a finite extension; its existence is guaranteed by Lemma~\ref{lem:finite-splitting}.
Proposition~\ref{prop:splitting-unique} shows that any two splitting fields admit an isomorphism preserving the base field, although this isomorphism need not be unique.

\phantomsection\label{def:separable-extension}
An algebraic element \(\alpha\) is \textbf{separable} over \(k\) if its minimal polynomial has no repeated roots in a splitting field.
If every element of an algebraic extension \(E/k\) is separable over \(k\), then \(E/k\) is called a \textbf{separable extension}.
The criterion relating repeated roots to the formal derivative is given in Proposition~\ref{prop:derivative-roots} of this appendix.
A finite extension over a general field can be inseparable.

\begin{example}[Power relations and an inseparable extension]\label{ex:inseparable-rational-extension}
Let \(u\) be an indeterminate, with \(E=\F_p(u)\) and \(k=\F_p(u^p)\).
The element \(u\) satisfies the polynomial \(X^p-u^p\) in \(k[X]\), so \([E:k]\le p\).
On the other hand, if
\[
\sum_{i=0}^{p-1}a_i(u^p)u^i=0,\qquad a_i(u^p)\in k,
\]
clearing a common denominator of the coefficients gives \(\sum_{i=0}^{p-1}P_i(u^p)u^i=0\) with \(P_i\in\F_p[T]\).
The degrees of all monomials in the \(i\)th term are congruent to \(i\) modulo \(p\), so monomials from different terms cannot cancel.
Since \(u\) is an indeterminate, every \(P_i\) must be zero, and hence every \(a_i\) is zero.
Thus \(1,u,\ldots,u^{p-1}\) are linearly independent over \(k\), \([E:k]=p\), and \(m_{u,k}(X)=X^p-u^p\).
But in \(E[X]\),
\[
X^p-u^p=(X-u)^p,
\]
so the extension is inseparable. For every \(R(T)\in\F_p(T)\), we have \(R(u)^p=R(u^p)\),
so the image of the absolute Frobenius map is exactly \(\F_p(u^p)=k\), a proper subfield of \(E\);
the map is injective but not surjective.
\end{example}

Over finite fields, the Frobenius orbit theorem settles the matter uniformly: Corollary~\ref{cor:separable} guarantees that every irreducible polynomial is separable,
and Theorem~\ref{thm:orbit} identifies \([k(\alpha):k]\) with the least Frobenius period of \(\alpha\).
This connects the linear dimension of an extension, the algebraic relation of a root, and its period under iteration.

\section{Divisibility, Bézout's Identity, and Unique Factorization}

The root bound and its inductive proof are given in Lemma~\ref{lem:root-bound}. We first establish division with remainder and the associated algebraic results.

\begin{lemma}[Division with remainder]\label{lem:division}
Let \(k\) be a field, and let \(A,B\in k[X]\) with \(B\ne0\). There exist unique \(Q,R\in k[X]\) such that
\[
A=QB+R,\qquad R=0\ \text{or}\ \deg R<\deg B.
\]
Consequently, \(a\) is a root of \(A\) if and only if \(X-a\) divides \(A\).
\end{lemma}
\begin{proof}
Successively cancel the highest-degree term of \(A\) using the leading term of \(B\). The degree strictly decreases, giving division with remainder.
If two representations had different quotients, subtracting them would express a nonzero multiple of \(B\) as a polynomial of degree less than \(\deg B\), which is impossible. The quotient is therefore unique, and so is the remainder.
Division by \(B=X-a\) has remainder exactly \(A(a)\).
\end{proof}

\begin{definition}\label{def:irreducible-polynomial}\label{def:monic-polynomial}
A polynomial \(f\in k[X]\) of positive degree is called an \textbf{irreducible polynomial} if it cannot be expressed as a product of two polynomials of positive degree.
A polynomial with leading coefficient \(1\) is called a \textbf{monic polynomial}.
\end{definition}

\begin{lemma}[Existence of irreducible factors]\label{lem:factor}
Every nonconstant polynomial has a monic irreducible factor and can be expressed as a constant times a finite product of monic irreducible polynomials.
\end{lemma}
\begin{proof}
If the polynomial is irreducible, divide by its leading coefficient. Otherwise, factor it into two nonconstant polynomials of smaller degree and treat each in turn.
Induction on the degree shows that this process terminates.
\end{proof}

This factorization result does not prescribe factor degrees; it therefore applies to adjoining finitely many roots without assuming the existence of irreducible polynomials of specified degrees.

\begin{lemma}[Bézout's identity]\label{lem:bezout}
Let \(A,B\in k[X]\) be not both zero, and let \(D\) be their monic greatest common divisor. Then there exist \(U,V\in k[X]\) such that
\[
UA+VB=D.
\]
In particular, if \(A,B\) are coprime, there exist \(U,V\) such that \(UA+VB=1\).
\end{lemma}
\begin{proof}
Among all nonzero \(UA+VB\), choose one of least degree and make it monic; denote it by \(D\).
Divide \(A\) by \(D\): \(A=QD+R\). The remainder \(R=A-QD\) is still a polynomial linear combination of \(A,B\).
If \(R\ne0\), it has smaller degree, a contradiction. Hence \(D\mid A\), and similarly \(D\mid B\).
Every common divisor also divides \(UA+VB=D\), so \(D\) is precisely the monic greatest common divisor.
\end{proof}

\begin{example}
In \(\F_2[X]\), let \(f=X^2+X+1\). Then
\[
(X+1)X+f=1.
\]
Thus in the quotient ring modulo \(f\), the element \(X+1\) is the multiplicative inverse of \(X\).
This is the role of Bézout's identity in the root-adjoining construction.
\end{example}

\begin{corollary}[Primality of irreducible factors and unique factorization]\label{cor:unique-factorization}
If an irreducible polynomial \(f\mid AB\), then \(f\mid A\) or \(f\mid B\).
Every nonzero polynomial factors uniquely, up to the order of the factors, as a constant times a product of monic irreducible factors.
\end{corollary}
\begin{proof}
If \(f\nmid A\), then \(\gcd(f,A)=1\), so we can write \(UA+Vf=1\). Multiplying by \(B\) gives \(f\mid B\).
Compare two irreducible factorizations. Every factor in the first must divide some factor in the second;
both are monic and irreducible, so they are equal. Cancelling these factors successively proves uniqueness.
\end{proof}

\section{The Chinese Remainder Theorem and Interpolation}
Decomposing a root set divides a polynomial problem into local problems. The Chinese remainder theorem combines local congruence conditions into a single polynomial residue class.

\begin{definition}\label{v101:def-ideal}
An \textbf{ideal} of a commutative ring \(R\) is a nonempty subset \(I\) closed under addition and subtraction and satisfying \(rI\subseteq I\) for \(r\in R\).
The congruence \(a\equiv b\pmod I\) modulo an ideal means \(a-b\in I\); operations on representatives make the congruence classes into the quotient ring \(R/I\).
\end{definition}
In \(k[X]\), every nonzero ideal has the form \((g)=\{gh:h\in k[X]\}\):
choose a monic nonzero polynomial \(g\) of least degree in the ideal and divide any element by \(g\). The remainder is still in the ideal and must be zero.
The kernel of a ring homomorphism is an ideal. If a homomorphism \(\psi:R\to S\) has kernel \(I\), then \(R/I\to\im\psi,\ \bar a\mapsto\psi(a)\) is an isomorphism, since it is well defined and surjective with zero kernel.
An evaluation homomorphism therefore factors through the corresponding quotient ring.

\begin{theorem}[Chinese remainder theorem]\label{thm:crt-basic}
Let \(s\ge1\), let \(g_1,\ldots,g_s\in k[X]\) be pairwise coprime polynomials of positive degree, and put \(G=\prod_i g_i\). Then
\[
k[X]/(G)\ \longrightarrow\ \prod_{i=1}^s k[X]/(g_i),\qquad
\bar A\longmapsto(A\bmod g_1,\ldots,A\bmod g_s)
\]
is a ring isomorphism. In particular, every tuple of local remainders corresponds to exactly one polynomial of degree less than \(\deg G\).
\end{theorem}
\begin{proof}
First check that the map is well defined. If \(G\mid A-B\), then \(g_i\mid A-B\) for every index, so the corresponding remainder classes agree.
The remainder map preserves addition, multiplication, and the identity, and hence is a ring homomorphism.

Put \(G_i=G/g_i\). Pairwise coprimality gives \(\gcd(G_i,g_i)=1\);
Bézout's identity gives \(u_iG_i+v_ig_i=1\). Set \(e_i=u_iG_i\). Then
\[
e_i\equiv1\pmod{g_i},\qquad e_i\equiv0\pmod{g_j}\quad(j\ne i).
\]
For any tuple of remainder classes \(\bar a_i\), choose representatives \(a_i\) and take \(A=\sum_i a_ie_i\). Modulo \(g_j\), only the \(j\)th summand survives, so \(A\equiv a_j\pmod{g_j}\).
This proves surjectivity.

If \(A\) maps to zero, then \(g_i\mid A\) for every index. The product of coprime polynomials also divides \(A\):
for example, from \(g_1\mid A\) write \(A=g_1B\), and use \(\gcd(g_1,g_2)=1\) to obtain \(g_2\mid B\);
then proceed by induction. Thus \(G\mid A\), and the kernel of the map is zero.
Finally, every class modulo \(G\) has a unique remainder of degree less than \(\deg G\), giving the asserted uniqueness.
\end{proof}

The elements \(\bar e_i\) satisfy \(\bar e_i^2=\bar e_i\), \(\bar e_i\bar e_j=0\) for \(i\ne j\), and \(\sum_i\bar e_i=1\),
and are called \textbf{idempotents} separating the components. The Berlekamp algorithm in Chapter~\ref{chap:algorithms} uses these separating functions to distinguish root orbits.
The proof for pairwise coprime integer moduli is identical.

\begin{corollary}[Lagrange interpolation]\label{cor:interpolation-basic}
For distinct \(a_1,\ldots,a_s\in k\) and arbitrary \(b_1,\ldots,b_s\in k\), the unique polynomial of degree less than \(s\)
satisfying \(A(a_i)=b_i\) is
\[
A(X)=\sum_{i=1}^s b_i\prod_{j\ne i}\frac{X-a_j}{a_i-a_j}.
\]
\end{corollary}
\begin{proof}
Take \(g_i=X-a_i\); the displayed products are the corresponding Lagrange basis polynomials. Substitution of each \(a_i\) verifies existence; a nonzero difference of two solutions could not have \(s\) roots.
\end{proof}

\phantomsection\label{v85term:def-formal-derivative}\ignorespaces
The \textbf{formal derivative} of a polynomial \(A(X)=\sum_i a_iX^i\) is defined by \(A'(X)=\sum_{i\geq1} i a_iX^{i-1}\). It satisfies the usual sum and product rules,
which can be verified directly by comparing coefficients, without the analytic notion of a limit.

\begin{proposition}[Repeated roots and the derivative]\label{prop:derivative-roots}
Let \(P\in k[X]\) be a nonzero polynomial. In an extension field where \(P\) splits, \(P\) has a repeated root if and only if \(\gcd(P,P')\ne1\).
\end{proposition}
\begin{proof}
A nonzero constant polynomial has no roots, and its monic greatest common divisor with its zero derivative is \(1\), so the assertion holds.
Now suppose \(P\) has positive degree. For any root \(a\) of \(P\), write in the given splitting extension
\[
P(X)=(X-a)^mQ(X),\qquad Q(a)\ne0.
\]
By the product rule, if \(m=1\), then \(P'(a)=Q(a)\ne0\);
if \(m\ge2\), then \(P'(a)=0\). Thus the repeated roots are exactly the common roots of \(P,P'\).

Let \(D=\gcd(P,P')\). If \(D=1\), Bézout's identity \(UP+VP'=1\) excludes a common root. If \(D\) is nonconstant, then since \(D\mid P\) and \(P\) splits in the specified extension,
\(D\) also splits there and has a root \(a\). Since \(D\mid P'\), this is also a root of \(P'\).
Hence \(P\) has a repeated root if and only if \(D\ne1\).
\end{proof}

\begin{corollary}[The Euler--Lagrange coefficient identity]\label{v84:cor-euler-lagrange}
Let \(S\subseteq k\), with \(|S|=n\ge1\), and put \(P_S(X)=\prod_{a\in S}(X-a)\).
If \(\deg A<n\), then
\begin{equation}\label{v84:eq-euler-lagrange}
 \sum_{a\in S}\frac{A(a)}{P_S'(a)}=[X^{n-1}]A.
\end{equation}
In particular, for \(0\le j\le n-1\),
\[
 \sum_{a\in S}\frac{a^j}{P_S'(a)}=
 \begin{cases}0,&j<n-1,\\1,&j=n-1.\end{cases}
\]
\end{corollary}
\begin{proof}
Lagrange interpolation gives \(A(X)=\sum_{a\in S}A(a)P_S(X)/((X-a)P_S'(a))\).
Each \(P_S(X)/(X-a)\) is monic of degree \(n-1\); compare the coefficients of \(X^{n-1}\).
\end{proof}

\chapter{Algebraic Closures and Compatible Finite Periods}
\label{app:algebraic-closure}

An algebraic closure of a finite field can be constructed as an increasing union of finite extension fields. This extends the process of adjoining finitely many roots to a common extension field for all algebraic equations.

\begin{definition}\label{def:algebraic-closure}
A field \(\Omega\) is called an \textbf{algebraically closed field} if every nonconstant polynomial in \(\Omega[X]\) has a root in \(\Omega\).
Equivalently, every nonconstant polynomial splits completely over \(\Omega\): after finding a root, divide out the corresponding linear factor and proceed by induction.
If, in addition, \(\Omega/k\) is algebraic, meaning that every element of \(\Omega\) is algebraic over \(k\),
then \(\Omega\) is called an \textbf{algebraic closure} of \(k\).
\end{definition}

\begin{theorem}\label{thm:closure}
Every finite field \(k\) has an algebraic closure, which can be constructed as a union of finite extension fields.
\end{theorem}
\begin{proof}
Put \(K_0=k\). After constructing a finite field \(K_j\), consider all monic polynomials in \(K_j[X]\) of degrees from \(1\) to \(j+1\). There are only finitely many, since each coefficient has finitely many choices.
Their product is a nonconstant polynomial \(P_j\).
By the splitting-extension construction in Lemma~\ref{lem:finite-splitting}, there is a finite extension field \(K_{j+1}\supseteq K_j\) in which \(P_j\) splits completely.
Each of the polynomials above then splits completely in \(K_{j+1}\). This gives an increasing sequence of fields
\[
K_0\subseteq K_1\subseteq K_2\subseteq\cdots.
\]
Put \(\Omega=\bigcup_{j\geq0}K_j\). Any finite collection of elements belongs to some common \(K_j\),
so sums, differences, products, and inverses of nonzero elements remain in \(\Omega\); hence \(\Omega\) is a field.

Take any nonconstant \(A\in\Omega[X]\) of degree \(D\). It has only finitely many coefficients,
all of which therefore belong to some \(K_N\). Choose \(j\geq N\) with \(j+1\geq D\).
After division by its leading coefficient, \(A\) is one of the polynomials considered in constructing \(K_{j+1}\),
so it splits completely in \(K_{j+1}\subseteq\Omega\). Thus \(\Omega\) is algebraically closed.

Finally, \(\alpha\in\Omega\) belongs to some finite field \(K_j\). Hence \(\alpha^{|K_j|}=\alpha\), so \(\alpha\) is a root of the nonzero polynomial \(X^{|K_j|}-X\), whose coefficients lie in the prime field and hence in \(k\).
Thus \(\alpha\) is algebraic over \(k\), and \(\Omega\) is an algebraic closure of \(k\).
\end{proof}

Adjoining finitely many roots makes each stage possible; the increasing union ensures that any finite set of coefficients eventually belongs to a single stage. Together, these properties ensure that the resulting field is algebraically closed.

\section*{Common Notation in an Algebraic Closure}
Once an algebraic closure \(\overline{\F}_p\) of \(\F_p\) is fixed, we may define simultaneously in this common extension field
\[
\F_{p^n}=\{a\in\overline{\F}_p:a^{p^n}=a\}.
\]
For each \(n\), this set has \(p^n\) elements, and
\[
\F_{p^a}\subseteq\F_{p^b}\quad\Longleftrightarrow\quad a\mid b.
\]
Moreover,
\[
\overline{\F}_p=\bigcup_{n\geq1}\F_{p^n}.
\]
For the last equality, simply note that every \(\alpha\in\overline{\F}_p\) is algebraic over \(\F_p\),
so \(\F_p(\alpha)\) is a finite field, all of whose elements satisfy the corresponding power identity.

The family \(\{\F_{p^n}\}\) is not nested in the usual order of \(n\); inclusion is governed by divisibility.
To express the closure as an increasing union, we may write
\[
\overline{\F}_p=\bigcup_{s\geq1}\F_{p^{s!}},
\qquad
\F_{p^{s!}}\subseteq\F_{p^{(s+1)!}}.
\]
This holds because every positive integer \(n\) divides all sufficiently large \(s!\).

\section{Compatible Finite Periods and Absolute Automorphisms}
\label{sec:profinite-frobenius}

In \(\Omega=\overline{\F}_q\), every element belongs to some finite subfield.
Thus an \(\F_q\)-automorphism is completely determined by its restrictions to these finite subfields.
Each restriction is a power of Frobenius, and the restrictions must agree on common subfields.

\begin{definition}\leavevmode\phantomsection\label{def:profinite-integers}
Define the additive group
\[
\widehat{\Z}
=\left\{(a_n)_{n\ge1}:a_n\in\Z/n\Z,\quad
a_m\equiv a_n\pmod n\ \text{whenever }n\mid m\right\},
\]
with componentwise addition; componentwise multiplication makes it a ring. It is called the \textbf{profinite completion} of the ring of integers, and is also denoted by \(\varprojlim_n\Z/n\Z\).
\end{definition}

\begin{theorem}\label{thm:absolute-frobenius}
Let \(\sigma(x)=x^q\). Every \(\tau\in\Aut_{\F_q}(\Omega)\) corresponds uniquely to an \((a_n)\in\widehat{\Z}\) through the rule
\[
\tau|_{\F_{q^n}}=\sigma^{a_n}|_{\F_{q^n}}.
\]
This correspondence gives a group isomorphism \(\Aut_{\F_q}(\Omega)\simeq\widehat{\Z}\).
\end{theorem}
\begin{proof}
Every \(\tau\) preserves the root set of \(x^{q^n}=x\) and therefore maps the unique subfield \(\F_{q^n}\) to itself.
The classification of finite-field automorphisms gives a unique residue class \(a_n\pmod n\). If \(n\mid m\), compatibility of restrictions gives \(a_m\equiv a_n\pmod n\).

Conversely, given a compatible family, define \(\tau(x)=\sigma^{a_n}(x)\) when \(x\in\F_{q^n}\).
If \(x\) also lies in another finite subfield, take their compositum; compatibility ensures that the two definitions agree.
Any two elements belong to a common finite subfield, so this map preserves sums, products, and the identity.
The compatible family \((-a_n)\) defines its inverse, so it is an automorphism.
Composition corresponds to addition of exponents, and the finite subfields have union \(\Omega\), so the resulting correspondence is a group isomorphism.
\end{proof}

An ordinary integer \(a\) gives the family \((a\bmod n)_n\). These do not exhaust the compatible families.
For example, for each \(n\) use the Chinese remainder theorem to prescribe
\[
a_n\equiv0\pmod{2^{v_2(n)}},\qquad
a_n\equiv1\pmod{n/2^{v_2(n)}}.
\]
These congruences define a compatible family. If it came from an integer, that integer would be divisible by every power of two and hence would be zero;
yet the condition modulo three requires it to be one, a contradiction.

\leavevmode\phantomsection\label{def:procyclic-group}\ignorespaces For \(\tau\in\Aut_{\F_q}(\Omega)\) and \(n\ge1\), put
\[
U_n(\tau)=\{\eta\in\Aut_{\F_q}(\Omega):
\eta|_{\F_{q^n}}=\tau|_{\F_{q^n}}\}.
\]
These sets form a neighborhood basis at \(\tau\): any finite set of restriction conditions can be expressed simultaneously on the compositum of the corresponding finite subfields.
This defines a topology on the automorphism group, corresponding to the topology on \(\widehat{\Z}\) defined by finitely many residue-class conditions.
In this topology the subgroup of integer powers of Frobenius is dense: for any \(\tau\) and \(n\),
choose an integer \(a\) representing \(a_n\); then \(\sigma^a\in U_n(\tau)\).
\phantomsection\label{def:profinite-group}
An inverse limit of finite groups, equipped with the topology defined by finitely many coordinate conditions, is called a profinite group. Its elements are tuples satisfying the compatibility conditions imposed by the transition homomorphisms.
A profinite group with a dense cyclic subgroup is called a \textbf{procyclic group}.
Thus \(\Aut_{\F_q}(\Omega)\) is procyclic. To see that it is not cyclic as an abstract group, identify it with \(\widehat{\Z}\).
If some \(a\in\widehat{\Z}\) generated this group, then \(1=(1\bmod n)_n\) would equal \(ma\) for an integer \(m\).
The case \(m=0\) is impossible after reduction modulo two. If \(|m|>1\), reduction modulo any prime divisor of \(|m|\) would give \(1=0\).
Hence \(m=\pm1\), so \(a=\pm1\). These elements generate only the embedded copy of \(\Z\), which is a proper subgroup, as shown by the compatible family constructed above.

\section{Pointwise Periods and a Common Period}
\begin{proposition}\label{prop:pointwise-uniform-period}
Let \(K\) be a field containing \(\F_q\), and put \(\sigma(x)=x^q\). Then
\begin{enumerate}
\item \(K/\F_q\) is algebraic if and only if every \(x\in K\) satisfies some \(\sigma^{m_x}(x)=x\), where \(m_x\ge1\) may depend on \(x\);
\item \(K\) is finite if and only if there is an \(N\ge1\) such that \(\sigma^N=\id_K\).
\end{enumerate}
\end{proposition}
\begin{proof}
If \(K/\F_q\) is algebraic, then for every \(x\) the simple extension \(\F_q(x)\) is a field of \(q^{m_x}\) elements for some exponent,
so \(x^{q^{m_x}}=x\). Conversely, this identity makes \(x\) a root of the nonzero polynomial \(X^{q^{m_x}}-X\in\F_q[X]\),
so every \(x\) is algebraic over \(\F_q\). This proves the first assertion.
If \(K\) is finite, write \(|K|=q^N\); then \(\sigma^N=\id_K\).
Conversely, if the same \(N\) works for every element, all elements of \(K\) are roots of \(X^{q^N}-X\),
and the root bound gives \(|K|\le q^N\), so \(K\) is finite.
\end{proof}
Every element of an algebraic closure has finite period, but no common finite period controls all elements:
it contains finite extensions of every degree and is therefore infinite. For a finite field, a single integer is a common multiple of all pointwise periods.
The root bound converts a common period into finiteness of the whole field, connecting this with the Frobenius fixed-point construction in Chapter~\ref{chap:foundations}.

\section{Intermediate Fields Determined by Allowed Periods}
\label{sec:intermediate-periods}

Fix an algebraic closure \(\Omega\) of \(k=\F_q\), and put \(\sigma(x)=x^q\).
Inclusion and composita of finite subfields are governed by divisibility and least common multiples; these rules also characterize arbitrary intermediate fields in \(\Omega\).

\begin{proposition}[Intermediate fields and allowed periods]\leavevmode\phantomsection\label{prop:intermediate-periods}
For an intermediate field \(k\subseteq K\subseteq\Omega\), define
\[
\mathcal D(K)=\{n\ge1:\F_{q^n}\subseteq K\}.
\]
The map \(K\mapsto\mathcal D(K)\) is a bijection between intermediate fields and nonempty sets \(\mathcal D\) of positive integers satisfying the following conditions:
\begin{enumerate}
\item if \(n\in\mathcal D\) and \(d\mid n\), then \(d\in\mathcal D\);
\item if \(m,n\in\mathcal D\), then \(\lcm(m,n)\in\mathcal D\).
\end{enumerate}
The inverse map is
\[
\boxed{\qquad K_{\mathcal D}=\bigcup_{n\in\mathcal D}\F_{q^n}.\qquad}
\]
The set \(\mathcal D(K)\) is also exactly the set of exact Frobenius periods of elements of \(K\).
\end{proposition}
\begin{proof}
First, \(1\in\mathcal D(K)\). If \(d\mid n\) and \(\F_{q^n}\subseteq K\), then
\(\F_{q^d}\subseteq\F_{q^n}\subseteq K\).
If \(\F_{q^m},\F_{q^n}\subseteq K\), their compositum is also contained in \(K\);
by Corollary~\ref{cor:intersection-compositum}, this compositum is \(\F_{q^{\lcm(m,n)}}\).
Thus \(\mathcal D(K)\) satisfies both conditions.
Every \(\alpha\in K\) is algebraic over \(k\), so \(k(\alpha)=\F_{q^n}\) for some \(n\);
this \(n\) belongs to \(\mathcal D(K)\), and hence
\[
K=\bigcup_{n\in\mathcal D(K)}\F_{q^n}.
\]

Conversely, suppose a nonempty set \(\mathcal D\) satisfies the stated conditions. Closure under positive divisors gives \(1\in\mathcal D\),
so \(K_{\mathcal D}\) contains \(k\).
Take any \(x\in\F_{q^m}\) and \(y\in\F_{q^n}\) with \(m,n\in\mathcal D\).
By the second condition, they belong to the common finite field \(\F_{q^{\lcm(m,n)}}\) contained in the union.
Thus \(x-y,xy\), and \(x^{-1}\) when \(x\ne0\), remain in the union, so \(K_{\mathcal D}\) is an intermediate field.

By definition, \(\mathcal D\subseteq\mathcal D(K_{\mathcal D})\).
If \(r\in\mathcal D(K_{\mathcal D})\), choose a generator \(\alpha\) of \(\F_{q^r}/k\).
Since \(\alpha\in K_{\mathcal D}\), there exists \(n\in\mathcal D\) such that \(\alpha\in\F_{q^n}\).
Then \(\F_{q^r}=k(\alpha)\subseteq\F_{q^n}\), so \(r\mid n\), and the first condition gives \(r\in\mathcal D\).
This proves \(\mathcal D(K_{\mathcal D})=\mathcal D\), so the two constructions are inverse to each other.

Finally, by the orbit theorem, the exact period of an element \(\alpha\) is \([k(\alpha):k]\).
If this period is \(n\), then \(\F_{q^n}=k(\alpha)\subseteq K\);
conversely, if \(\F_{q^n}\subseteq K\), a \(k\)-generator of this field has exact period \(n\).
Thus the allowed degrees are precisely the exact periods that occur.
\end{proof}

\begin{example}[An infinite field allowing only power-of-two periods]\label{exm:two-power-period-field}
The set \(\mathcal D=\{2^j:j\ge0\}\) is closed under positive divisors and least common multiples, so
\[
K=\bigcup_{j\ge0}\F_{q^{2^j}}
\]
is an intermediate field in \(\Omega\), and \(\mathcal D(K)=\mathcal D\).
Its finite subfields have unbounded orders, so \(K\) is infinite; nevertheless, every element has Frobenius period a power of two.
Since \(3\notin\mathcal D(K)\), we have \(\F_{q^3}\nsubseteq K\).
No irreducible cubic over \(k\) can have a root in \(K\), since such a root would generate the subfield \(\F_{q^3}\).
Thus \(K/k\) is an infinite algebraic extension, but \(K\) is not algebraically closed.
Proposition~\ref{prop:pointwise-uniform-period} also shows that Frobenius on \(K\) has no common finite period valid for every element.
\end{example}

Inside a fixed algebraic closure, specifying a subfield containing the base field is equivalent to specifying a family of exact periods closed under positive divisors and least common multiples.
The first condition corresponds to inclusion of finite subfields; the second puts any finite collection of elements in a common finite extension field.
Together, they express closure under field operations as closure of periods.

\section{Uniqueness of the Isomorphism Type and Choice of an Isomorphism}
\label{sec:isomorphism-choices}

Finite fields of the same order have a unique isomorphism type, but choosing an isomorphism between two concrete models usually still requires specifying the image of a root.
This distinction occurs in quotient-ring representations, explicit isomorphisms, and compatible identifications of finite extension fields.

\leavevmode\phantomsection\label{def:isomorphism-set}\ignorespaces
Let \(k=\F_q\), and let \(K,L\) be two finite extensions of \(k\) of degree \(n\).
Write \(\operatorname{Iso}_k(K,L)\) for the set of all \(k\)-isomorphisms.
Uniqueness of splitting fields over the specified base field shows that this set is nonempty.

\begin{proposition}[All isomorphisms between two finite-field models]\label{prop:finite-field-isomorphism-choices}
Choose \(\iota_0\in\operatorname{Iso}_k(K,L)\), and put \(\sigma_L(y)=y^q\). Then
\begin{equation}\label{eq:finite-field-isomorphism-choices}
 \operatorname{Iso}_k(K,L)
 =\{\sigma_L^j\circ\iota_0:0\le j<n\},
\end{equation}
and the \(n\) maps on the right are distinct.
\end{proposition}
\begin{proof}
For any \(\iota\in\operatorname{Iso}_k(K,L)\), the composite \(\iota\circ\iota_0^{-1}\) belongs to \(\Aut_k(L)\).
Theorem~\ref{thm:automorphisms} says it equals \(\sigma_L^j\) for a unique \(0\le j<n\).
Conversely, every \(\sigma_L^j\circ\iota_0\) is a \(k\)-isomorphism. Uniqueness also gives distinctness.
\end{proof}

If \(K=k(\alpha)\) and \(\beta=\iota_0(\alpha)\), these isomorphisms send \(\alpha\) respectively to
\[
 \beta,\ \beta^q,\ \ldots,\ \beta^{q^{n-1}}.
\]
Thus degree data determine the isomorphism type, whereas specifying the image of a generator determines a particular isomorphism.
Finding a root of a minimal polynomial in the explicit-isomorphism construction of Chapter~\ref{chap:algorithms} makes precisely this choice.
After choosing \(\iota_0\), we may index all isomorphisms by \(\Z/n\Z\); if the starting point is changed to \(\sigma_L^t\circ\iota_0\), the index of the same map changes from \(j\) to \(j-t\).
This indexing depends on the chosen starting point.

For example, let
\[
 K=\F_2[X]/(X^2+X+1),\qquad L=\F_2[Y]/(Y^2+Y+1),
\]
and let \(\alpha,\beta\) denote the residue classes of \(X,Y\), respectively.
The two isomorphisms are determined by \(\alpha\mapsto\beta\) and \(\alpha\mapsto\beta^2=\beta+1\), respectively.
Naming \(\beta\) specifies the first isomorphism, but interchanging \(\beta\) and \(\beta+1\) does not change the isomorphism type of \(L\).
When \(n=1\), there is only one \(k\)-isomorphism, in agreement with the classification above.

Within a fixed algebraic closure \(\Omega=\overline\F_q\), identifications of finite levels must also be compatible.
Theorem~\ref{thm:absolute-frobenius} has already shown that the restrictions of a global automorphism are \(\sigma^{a_n}|_{\F_{q^n}}\),
where \((a_n)\in\widehat\Z\) and \(a_m\equiv a_n\pmod n\) whenever \(n\mid m\).
Each restriction is an integer power of Frobenius on its finite level, but a single ordinary integer need not represent them simultaneously on every level.
Thus the finite-level groups are cyclic, whereas the full group is procyclic. The passage between them is governed by compatibility of restrictions, not by a common finite period.

\chapter{Arithmetic Completion of the Primitive Normal Basis Theorem}
\label{v93:chap-primitive-normal-arithmetic}

Section~\ref{v93:sec-primitive-normal-proof} reduces the existence of primitive normal elements to a joint count for two local generation conditions.
This appendix determines exactly when the basic bound~\eqref{v93:eq-easy} holds. The eight pairs not covered by that bound are treated in the main text by structural exclusions, inclusion--exclusion among elements of full degree, and a local density sieve.
The arithmetic retains the two structures used in the main proof: Frobenius orbits control the number of normality conditions, and multiplicative orders control the possible integer prime factors.

Throughout, write
\[
 k=\F_q,\qquad E=\F_{q^n},\qquad Q=q^n,\qquad N=Q-1,
 \qquad n=p^a m,\quad p=\Char k,\quad p\nmid m.
\]
Let
\[
 M(T)=T^m-1=\prod_{i=1}^{s}P_i(T),\qquad
 \kappa=\gcd(n,q-1),\qquad
 R=\frac{Q-1}{(q-1)\kappa},
\]
where the \(P_i\) are distinct monic irreducible polynomials. For a positive integer \(u\), put \(W(u)=2^{\omega(u)}\), where \(\omega(u)\) counts its distinct prime factors.
By the scaling lemma~\ref{v93:lem-scaling} and the basic joint estimate, the required sufficient condition is
\[
 q^{n/2}>(2^{\omega(R)}-1)(2^s-1).
\]
The finite tables below give monotone bounds for infinite parameter families, the remaining boundary cases, and integer comparisons that can be checked individually.

\section{Uniform Arithmetic Bounds for the Two Degree Distributions}
\label{v93:sec-degree-controls}
\subsection{Irreducible Factor Counts: Low-Degree Roots and the Remaining Degree}
Put \(a_j=\gcd(n,q^j-1)=\gcd(m,q^j-1)\). Exactly \(a_j\) distinct roots of \(M=T^m-1\) belong to \(\F_{q^j}\). Hence the numbers of irreducible factors of degrees exactly one, two, and three are, respectively,
\[
 a_1,\qquad \frac{a_2-a_1}{2},\qquad \frac{a_3-a_1}{3}.
\]
Bounding the degrees of the remaining factors below by two, three, and four, respectively, gives
\begin{equation}\label{v93:eq-degree-budgets}
\boxed{\begin{aligned}
 s&\le \frac{m+a_1}{2},\\
 s&\le \frac m3+\frac{a_1}{2}+\frac{a_2}{6},\\
 s&\le \frac m4+\frac{5a_1}{12}+\frac{a_2}{4}+\frac{a_3}{12}.
\end{aligned}}
\end{equation}
For example, the third inequality follows from
\(s\le m/4+(3/4)a_1+(1/2)(a_2-a_1)/2+(1/4)(a_3-a_1)/3\).
These bounds remain valid when \(p\mid n\), since the characteristic power has already been removed from \(m\).

For exact values, Corollary~\ref{cor:cyclotomic-equal-degree} gives the cyclotomic factorization formula
\begin{equation}\label{v93:eq-exact-degree-count}
 s=\sum_{d\mid m}\frac{\varphi(d)}{\ord_d(q)},\qquad \ord_1(q):=1.
\end{equation}
This involves only the divisors of \(m\) and multiplicative orders modulo those divisors.

\subsection{Multiplicative Prime Factors: Retaining Only the Allowed Small Primes}
\begin{lemma}[A power bound from the small prime factors]\label{v93:lem-small-prime-bound}
Let \(h\ge1\) be an integer, and let \(\mathcal P\) be a set of primes less than \(2^h\) containing every prime factor of \(R\) below \(2^h\). Put \(r=|\mathcal P|\) and \(D=\prod_{\ell\in\mathcal P}\ell\), with the empty product equal to one. Then
\begin{equation}\label{v93:eq-small-prime-bound}
 W(R)\le 2^r(R/D)^{1/h}.
\end{equation}
\end{lemma}
\begin{proof}
Each prime factor \(\ell\ge2^h\) of \(R\) satisfies \(2\le\ell^{1/h}\). For each smaller prime factor, retain the correction factor \(2/\ell^{1/h}>1\). Including also the primes in \(\mathcal P\) that do not divide \(R\) can only increase the upper bound. Multiplying these inequalities proves the claim.
\end{proof}
Taking \(h=4\) and \(\mathcal P=\{2,3,5,7,11,13\}\), the inequality \(2^{24}<625\cdot30030\) gives
\begin{equation}\label{v93:eq-fourth-universal}
 W(R)<5R^{1/4}.
\end{equation}
If \(R\) is odd, omit two; then \(2^{20}<81\cdot15015\) gives \(W(R)<3R^{1/4}\).

The following conditions determine which small primes can occur. If \(\ell\nmid q(q-1)\), then
\begin{equation}\label{v93:eq-order-admissible}
 \ell\mid R\quad\Longleftrightarrow\quad \ord_\ell(q)\mid n.
\end{equation}
If \(\ell\) is an odd prime factor of \(q-1\), Lemma~\ref{lem:valuation-lift} gives
\begin{equation}\label{v93:eq-valuation-admissible}
 v_\ell(R)=\max\{0,v_\ell(n)-v_\ell(q-1)\}.
\end{equation}
If \(q\) is odd, then \(R\) is odd when \(n\) is odd; when \(n\) is even,
\[
 v_2(R)=v_2(q+1)+v_2(n)-1-\min\{v_2(n),v_2(q-1)\}.
\]
Indeed, for an odd prime \(\ell\mid q-1\), subtracting the valuation of the denominator from \(v_\ell(q^n-1)=v_\ell(q-1)+v_\ell(n)\) gives~\eqref{v93:eq-valuation-admissible}. If \(q\) is odd and \(n\) is even, apply Lemma~\ref{lem:valuation-lift} to \(q^2\) with exponent \(n/2\) to obtain \(v_2(q^n-1)=v_2(q-1)+v_2(q+1)+v_2(n)-1\); subtracting the valuation of the denominator again gives the displayed formula.

\begin{lemma}[A criterion for a monotone parameter family]\label{v93:lem-family}
Within a parameter family, suppose that \(\alpha,\beta\ge0\) are rational numbers, \(\delta\ge1\) is an integer, and \(s\le\alpha n+\beta\), \(\kappa\ge\delta\). Choose an integer \(h\ge3\) and \(\mathcal P\) as in Lemma~\ref{v93:lem-small-prime-bound}. Put
\[
 L=\lcm\bigl(2,\operatorname{den}(h\alpha),\operatorname{den}(h\beta)\bigr),\qquad
 K=(q-1)\delta D.
\]
Here \(\operatorname{den}(x)\) denotes the positive denominator of the rational number \(x\) in lowest terms. If a positive integer \(n_0\) satisfies
\begin{align}
 q^{L(h/2-1)}&>2^{Lh\alpha},\label{v93:eq-family-slope}\\
 q^{Ln_0(h/2-1)}K^L&>2^{Lh(r+\alpha n_0+\beta)},\label{v93:eq-family-base}
\end{align}
then every \(n\ge n_0\) in this family satisfies the basic bound~\eqref{v93:eq-easy}.
\end{lemma}
\begin{proof}
From \(R<Q/((q-1)\delta)\) and the small-prime power bound,
\[
 (W(R)-1)(2^s-1)
 <2^{r+\alpha n+\beta}(Q/K)^{1/h}.
\]
Raising the sufficient condition that the right side be less than \(Q^{1/2}\) to the power \(Lh\) gives precisely the integer comparison~\eqref{v93:eq-family-base}, with \(n_0\) replaced by \(n\). By~\eqref{v93:eq-family-slope}, the ratio of the left side to the right side strictly increases when \(n\) increases by one. It therefore suffices to check the initial value. All exponents in these comparisons are nonnegative integers.
\end{proof}

The two bounds in this section come from Frobenius orbit lengths and multiplicative orders of integer prime factors. We first treat uniform ranges, then divide the remaining parameters into monotone families.

\section{Three Uniform Ranges}
\label{v93:sec-uniform-ranges}
\subsection{Degrees That Are Powers of the Characteristic}
For every positive integer \(u\),
\begin{equation}\label{v93:eq-squarefree-elementary}
 (W(u)-1)^2<2u;\qquad u\text{ odd: }(W(u)-1)^2<u.
\end{equation}
Both inequalities were proved in Lemma~\ref{v88:lem-squarefree-divisor-inequality}.

If \(n=p^a\), then \(s=1\) and \(\kappa=1\). For \(q=2\), the integer \(R=Q-1\) is odd; for \(q\ge3\), we have \(2R<Q\). Equation~\eqref{v93:eq-squarefree-elementary} immediately yields the basic bound. Thus all characteristic-power degrees are covered.

\subsection{Degrees Three and Five}
Let \(n=\nu\) be an odd prime with \(\nu\ne p\). Every prime factor \(\ell\) of \(R=\Phi_\nu(q)/\gcd(\nu,q-1)\) satisfies
\begin{equation}\label{v93:eq-prime-degree-primes}
 \ell\equiv1\pmod{2\nu}.
\end{equation}
Indeed, if \(\ell\ne\nu\), the congruence \(\Phi_\nu(q)\equiv0\pmod\ell\) implies that \(q\) has order \(\nu\), so \(\nu\mid\ell-1\); moreover, \(\ell\) is odd. If \(\nu\mid\Phi_\nu(q)\), then \(q\equiv1\pmod\nu\), and this prime occurs with valuation exactly one, so it is removed by the denominator.

For \(n=3\), we may take \(\mathcal P=\{7,13\}\) and \(h=4\). If \(q\equiv1\pmod3\), then \(s=3,\kappa=3\), and for \(q\ge8\),
\[
 91\cdot3(q-1)q^3>256\cdot7^4.
\]
The inequality is strict at \(q=8\), and the left side increases with \(q\). The smaller values \(q=4,7\) give \(R=7,19\), respectively, and both satisfy \(7<q^{3/2}\).
If \(q\equiv-1\pmod3\), then \(s=2,\kappa=1\), and for \(q\ge5\),
\[
 91(q-1)q^3>256\cdot3^4.
\]
Only \((q,n)=(2,3)\) remains. The case \(p=3\) is already covered by characteristic-power degrees.

For \(n=5\) and \(p\ne5\), take \(\mathcal P=\{11\}\) and \(h=4\). Put \(d=\ord_5(q)\); then \(s=1+4/d\). The sufficient condition to check is
\[
 11\kappa(q-1)q^5>16(2^s-1)^4.
\]
If \(d=1\), then \(q\ge11\) and \(\kappa=5,s=5\), and the inequality holds at the initial value. If \(d=2\), then \(s=3\), and the inequality holds for \(q\ge5\); the remaining case \(q=4\) has \(R=11\cdot31\), giving \(21<32\). If \(d=4\), then \(s=2\), and the inequality holds for \(q\ge3\); the remaining case \(q=2\) has \(R=31\), giving \(3<\sqrt{32}\). Thus every extension of degree five satisfies the basic bound.

\subsection{Larger Base Fields: \texorpdfstring{\(q\ge16\)}{q at Least 16}}
First suppose that \(n\ge6\) and \(n\nmid q-1\). Since \(\kappa\) is a proper divisor of \(n\), we have \(s\le(n+\kappa)/2\le3n/4\). By~\eqref{v93:eq-fourth-universal}, a sufficient condition for the basic bound is
\[
 (q/8)^n(q-1)>625.
\]
At \(q=16,n=6\), this follows from \(64\cdot15>625\), and the left side increases with both \(q\) and \(n\).

If \(n\mid q-1\), then \(s=n,\kappa=n\). For \(q\ge23\), it suffices to have
\[
 (q/16)^n n(q-1)>625,
\]
which follows from the initial comparison \(23^6\cdot132>625\cdot16^6\) and monotonicity.
It remains to treat \(q=16,17,19\). When \(n\mid q-1\), no odd prime factor of \(q-1\) divides \(R\); if \(q\equiv1\pmod4\), then \(R\) is also odd. In these three cases, Lemma~\ref{v93:lem-small-prime-bound} gives \(W(R)^4\le c_qR\), where
\[
 c_{16}=\frac{2^{12}}{7\cdot11\cdot13}<16,\quad
 c_{17}=\frac{2^{20}}{3\cdot5\cdot7\cdot11\cdot13}<81,\quad
 c_{19}=\frac{2^{20}}{2\cdot5\cdot7\cdot11\cdot13}<105.
\]
For \(n\ge6\), the quantity \((q/16)^nn(q-1)\) is at least \(90,96,108\), respectively, so the sufficient condition holds in each case.

Degrees three and five have already been treated. Degree four is a characteristic power in even characteristic. In odd characteristic, \(q\equiv1\pmod4\) gives \(s=4,\kappa=4\), with \(R\) odd, and for \(q\ge17\) the required comparison is
\[
 4q^4(q-1)>81\cdot15^4
\]
If \(q\equiv3\pmod4\), then \(s=3,\kappa=2\), and for \(q\ge16\) the required comparison is
\[
 2q^4(q-1)>625\cdot7^4
\]
Thus the basic bound holds for all \(q\ge16,n\ge3\).

\section{33 Monotone Families over Small Base Fields}
\label{v93:sec-small-field-families}
The only remaining base-field orders are \(2,3,4,5,7,8,9,11,13\). Each row of Table~\ref{v93:tab-families} covers every \(n\ge n_0\) satisfying its conditions. Here \(\delta\) is a lower bound for \(\kappa\), and \(\mathcal P\) is chosen using~\eqref{v93:eq-order-admissible} and the valuation formulas. The integer comparison at the initial value of row \(i\) is given explicitly in Table~\ref{v93:tab-certificates}.

All bounds for the number of irreducible factors follow from~\eqref{v93:eq-degree-budgets}. For example, for odd degrees over the binary field, use the third inequality with \(a_1=1,a_2\in\{1,3\},a_3\in\{1,7\}\); the resulting constant terms are \(3/4,5/4,7/4\). For even degrees over the binary field, \(m\le n/2\), and if four divides the degree, \(m\le n/4\). For \(q=3\), group the degrees by \(\gcd(n,8)\); for \(q=4,5\), use \(\gcd(n,15),\gcd(n,24)\), respectively. For \(q=7,8,9,11\), only the first inequality and \(\gcd(n,q-1)\) are needed. Finally, if \(n\nmid12\), then \(\gcd(n,12)\le n/2\), giving the row for \(q=13\).

The sets of small primes can also be checked by computing powers modulo small integers. For example, among primes below 64, only \(7,23,31,47\) can divide \(R\) for odd degrees over the binary field: the multiplicative orders of two modulo these primes are \(3,11,5,23\), respectively, while the orders for all other odd primes below 64 are even. For \(q=3,n\equiv2\pmod4,3\nmid n\), the only allowed primes below 32 are \(2,11,23\), and the latter two have orders five and eleven. All other rows require only primes below 32; their orders appear in Table~\ref{v93:tab-orders}.
\par\medskip
\begingroup\footnotesize
\setlength{\tabcolsep}{2.5pt}\renewcommand{\arraystretch}{1.32}
\begin{longtable}{@{}rr>{\raggedright\arraybackslash}p{0.26\textwidth}cccr>{\raggedright\arraybackslash}p{0.29\textwidth}r@{}}
\caption{Infinite parameter families; each row covers all admissible $n\ge n_0$}\label{v93:tab-families}\\
\toprule No.&$q$&Conditions&$\alpha$&$\beta$&$h$&$\delta$&$\mathcal P$&$n_0$\\\midrule\endfirsthead
\toprule No.&$q$&Conditions&$\alpha$&$\beta$&$h$&$\delta$&$\mathcal P$&$n_0$\\\midrule\endhead
\bottomrule\endfoot
1&2&$2\nmid n,\ 3\nmid n,\ 7\nmid n$&$\frac{1}{4}$&$\frac{3}{4}$&6&1&$7,23,31,47$&23\\
2&2&$2\nmid n$; exactly one of $3,7$ divides $n$&$\frac{1}{4}$&$\frac{5}{4}$&6&1&$7,23,31,47$&29\\
3&2&$2\nmid n,\ 21\mid n$&$\frac{1}{4}$&$\frac{7}{4}$&6&1&$7,23,31,47$&35\\
4&2&$n\equiv2\pmod4$&$\frac{1}{8}$&$\frac{7}{4}$&4&1&$3,7,11$&26\\
5&2&$4\mid n$&$\frac{1}{8}$&$\frac{1}{2}$&4&1&$3,5,7,11,13$&20\\
6&3&$3\mid n$&$\frac{1}{6}$&$1$&4&1&$2,5,7,11,13$&12\\
7&3&$2\nmid n,\ 3\nmid n$&$\frac{1}{3}$&$\frac{2}{3}$&4&1&$11$&11\\
8&3&$n\equiv2\pmod4,\ 3\nmid n$&$\frac{1}{3}$&$\frac{4}{3}$&5&2&$2,11,23$&18\\
9&3&$n\equiv4\pmod8,\ 3\nmid n$&$\frac{1}{3}$&$\frac{5}{3}$&5&2&$2,5,11,23,29$&22\\
10&3&$8\mid n,\ 3\nmid n$&$\frac{1}{3}$&$\frac{7}{3}$&5&2&$2,5,11,17,23,29$&28\\
11&4&$2\mid n$&$\frac{1}{4}$&$\frac{3}{2}$&4&1&$3,5,7,11,13$&12\\
12&4&$2\nmid n,\ 3\nmid n,\ 5\nmid n$&$\frac{1}{3}$&$\frac{2}{3}$&4&1&$\varnothing$&2\\
13&4&$2\nmid n,\ 5\mid n,\ 3\nmid n$&$\frac{1}{3}$&$\frac{4}{3}$&4&1&$11$&7\\
14&4&$2\nmid n,\ 3\mid n,\ 5\nmid n$&$\frac{1}{3}$&$2$&4&3&$3,7$&13\\
15&4&$2\nmid n,\ 15\mid n$&$\frac{1}{3}$&$4$&4&3&$3,7,11$&26\\
16&5&$2\nmid n$&$\frac{1}{3}$&$1$&5&1&$11,19,31$&3\\
17&5&$n\equiv2\pmod4$&$\frac{1}{3}$&$2$&5&2&$3,7,11,19,23,29,31$&9\\
18&5&$4\mid n,\ 3\nmid n$&$\frac{1}{3}$&$\frac{10}{3}$&5&4&$2,3,11,13,17,23,29$&14\\
19&5&$12\mid n$&$\frac{1}{3}$&$6$&5&4&$2,3,7,11,13,17,19,23,29,31$&23\\
20&7&$3\nmid n$&$\frac{1}{2}$&$1$&4&1&$2,5,11$&9\\
21&7&$2\nmid n,\ 3\mid n$&$\frac{1}{2}$&$\frac{3}{2}$&4&3&$3$&6\\
22&7&$6\mid n$&$\frac{1}{2}$&$3$&5&6&$2,3,5,11,13,17,19,23,29,31$&15\\
23&8&$7\nmid n$&$\frac{1}{2}$&$\frac{1}{2}$&4&1&$3,5,11,13$&5\\
24&8&$7\mid n$&$\frac{1}{2}$&$\frac{7}{2}$&5&7&$3,5,7,11,13,17,19,23,29,31$&13\\
25&9&$2\nmid n$&$\frac{1}{2}$&$\frac{1}{2}$&4&1&$7,11,13$&1\\
26&9&$n\equiv2\pmod4$&$\frac{1}{2}$&$1$&4&2&$5,7,11,13$&4\\
27&9&$n\equiv4\pmod8$&$\frac{1}{2}$&$2$&4&4&$5,7,11,13$&6\\
28&9&$8\mid n$&$\frac{1}{2}$&$4$&4&8&$2,5,7,11,13$&15\\
29&11&$2\nmid n,\ 5\nmid n$&$\frac{1}{2}$&$\frac{1}{2}$&4&1&$7$&1\\
30&11&$2\mid n,\ 5\nmid n$&$\frac{1}{2}$&$1$&4&2&$2,3,7,13$&5\\
31&11&$2\nmid n,\ 5\mid n$&$\frac{1}{2}$&$\frac{5}{2}$&4&5&$5,7$&5\\
32&11&$10\mid n$&$\frac{1}{2}$&$5$&4&10&$2,3,5,7,13$&16\\
33&13&$n\nmid12$&$\frac{3}{4}$&$0$&4&1&$2,3,5,7,11$&8\\
\end{longtable}\endgroup

\paragraph{The divisibility cases for \(q=13\).}
If \(n\mid12\), then after degrees one, two, and three have been treated, only \(n=4,6,12\) remain. For \(n=6,12\), we have \(s=n,\kappa=n\), and \(R\) is coprime to six. Taking \(h=4,\mathcal P=\{5,7,11\}\), it suffices to check
\[
 13^n\cdot12n\cdot385>2^{4n+12}.
\]
For integers \(n\ge5\), the expression \(n(13/16)^n\) decreases. Hence, for \(6\le n\le12\), it suffices to check
\[
 13^{12}\cdot144\cdot385>2^{60},
\]
which follows by direct integer multiplication. Only \((13,4)\) is left for the next table.

\paragraph{The explicit boundary set.}
From the degrees \(3\le n<n_0\) in each row, remove degrees three and five and all characteristic powers. Exactly the 42 pairs in Table~\ref{v93:tab-boundary} remain. Exhaustiveness is checked by listing the integers below \(n_0\le35\) that satisfy the divisibility conditions; the table lists the remaining degrees for each base field. Together with \((2,3)\), retained in Section~\ref{v93:sec-uniform-ranges}, these give precisely the eight exceptions to the basic bound stated below.

\section{The 42 Boundary Comparisons}
\label{v93:sec-boundary-comparisons}
In Table~\ref{v93:tab-boundary}, \(t_0\) is an upper bound for \(\omega(R)\), and \(B=(2^{t_0}-1)(2^s-1)\). Except in the starred rows, the last two columns satisfy \(B^2<q^n\). Every value of \(s\) follows from~\eqref{v93:eq-exact-degree-count}; the multiset of factor degrees is also listed for independent verification.

Rows marked \(\dagger\) have prime degree \(n\ge7\). No claim that \(R\) is prime is needed: equation~\eqref{v93:eq-prime-degree-primes} and
\[
 R<(2n+1)^{t_0+1}
\]
give \(\omega(R)\le t_0\). Apart from the integers discussed below, every factor used as a prime in the table is at most 683, so its primality can be checked by trial division by primes at most 23.

For four types of larger factors, only an upper bound for the number of distinct prime factors is needed:
\begin{itemize}
\item Since \(262657=\Phi_{27}(2)\) is not divisible by three, two has order 27 modulo every prime factor, so each such prime is at least 55. The inequality \(262657<55^4\) gives at most three distinct prime factors. This yields \(t_0=5\) for \((2,27)\).
\item Neither \(547=\Phi_{14}(3)\) nor \(1093=\Phi_7(3)\) is divisible by two or seven, so every prime factor is at least 15. Both integers are less than \(15^3\), and each therefore has at most two distinct prime factors. This yields \(t_0=5\) for \((3,14)\).
\item Neither \(1181\) nor \(1201\) is divisible by two, three, or five, and both are less than \(7^4\). Each therefore has at most three distinct prime factors. This yields \(t_0=7\) for \((3,20)\) and \(t_0=5\) for \((7,8)\).
\item We have \(5\cdot3221=\Phi_5(11)\) and \(13421=\Phi_{10}(11)\). Neither \(3221\) nor \(13421\) is divisible by two or five, so every prime factor is \(1\pmod{10}\). The three smallest such primes are \(11,31,41\), whose product 13981 exceeds both integers. Thus each has at most two distinct prime factors, giving \(t_0=6\) for \((11,10)\).
\end{itemize}
The general fact used here is that if a prime \(\ell\nmid qd\) divides \(\Phi_d(q)\), then \(\ord_\ell(q)=d\). This follows directly from the squarefreeness of \(X^d-1\) modulo \(\ell\) and its cyclotomic factorization.
\par\medskip
\begingroup\footnotesize
\setlength{\tabcolsep}{3pt}\renewcommand{\arraystretch}{1.28}
\begin{longtable}{@{}rrl>{\raggedright\arraybackslash}p{0.32\textwidth}rrr@{}}
\caption{Boundary data; a degree with superscript $[c]$ occurs $c$ times}\label{v93:tab-boundary}\\
\toprule$q$&$n$&Factor degrees&$R$ or its factorization&$t_0$&$B$&$q^n$\\\midrule\endfirsthead
\toprule$q$&$n$&Factor degrees&$R$ or its factorization&$t_0$&$B$&$q^n$\\\midrule\endhead
\bottomrule\endfoot
2&6&$1,2$&$3^{2}\cdot7\ *$&2&9&64\\
2&7&$1,3^{[2]}$&$127^{\dagger}$&1&7&128\\
2&9&$1,2,6$&$7\cdot73$&2&21&512\\
2&10&$1,4$&$3\cdot11\cdot31$&3&21&1024\\
2&11&$1,10$&$2047^{\dagger}$&2&9&2048\\
2&12&$1,2$&$3^{2}\cdot5\cdot7\cdot13$&4&45&4096\\
2&13&$1,12$&$8191^{\dagger}$&2&9&8192\\
2&14&$1,3^{[2]}$&$3\cdot43\cdot127$&3&49&16384\\
2&15&$1,2,4^{[3]}$&$7\cdot31\cdot151\ *$&3&217&32768\\
2&17&$1,8^{[2]}$&$131071^{\dagger}$&3&49&131072\\
2&18&$1,2,6$&$3^{3}\cdot7\cdot19\cdot73$&4&105&262144\\
2&19&$1,18$&$524287^{\dagger}$&3&21&524288\\
2&21&$1,2,3^{[2]},6^{[2]}$&$7^{2}\cdot127\cdot337$&3&441&2097152\\
2&22&$1,10$&$3\cdot23\cdot89\cdot683$&4&45&4194304\\
2&27&$1,2,6,18$&$7\cdot73\cdot262657$&5&465&134217728\\
3&4&$1^{[2]},2$&$2^{2}\cdot5\ *$&2&21&81\\
3&6&$1^{[2]}$&$2\cdot7\cdot13$&3&21&729\\
3&7&$1,6$&$1093^{\dagger}$&2&9&2187\\
3&8&$1^{[2]},2^{[3]}$&$2^{3}\cdot5\cdot41\ *$&3&217&6561\\
3&10&$1^{[2]},4^{[2]}$&$2\cdot11^{2}\cdot61$&3&105&59049\\
3&14&$1^{[2]},6^{[2]}$&$2\cdot547\cdot1093$&5&465&4782969\\
3&16&$1^{[2]},2^{[3]},4^{[2]}$&$2^{4}\cdot5\cdot17\cdot41\cdot193$&5&3937&43046721\\
3&20&$1^{[2]},2,4^{[4]}$&$2^{2}\cdot5^{2}\cdot11^{2}\cdot61\cdot1181$&7&16129&3486784401\\
4&6&$1^{[3]}$&$5\cdot7\cdot13$&3&49&4096\\
4&9&$1^{[3]},3^{[2]}$&$3\cdot7\cdot19\cdot73$&4&465&262144\\
4&10&$1,2^{[2]}$&$5^{2}\cdot11\cdot31\cdot41$&4&105&1048576\\
4&15&$1^{[3]},2^{[6]}$&$7\cdot11\cdot31\cdot151\cdot331$&5&15841&1073741824\\
5&4&$1^{[4]}$&$3\cdot13\ *$&2&45&625\\
5&6&$1^{[2]},2^{[2]}$&$3^{2}\cdot7\cdot31$&3&105&15625\\
5&8&$1^{[4]},2^{[2]}$&$2\cdot3\cdot13\cdot313\ *$&4&945&390625\\
5&12&$1^{[4]},2^{[4]}$&$3^{2}\cdot7\cdot13\cdot31\cdot601$&5&7905&244140625\\
7&4&$1^{[2]},2$&$2^{3}\cdot5^{2}$&2&21&2401\\
7&6&$1^{[6]}$&$2^{2}\cdot19\cdot43\ *$&3&441&117649\\
7&8&$1^{[2]},2^{[3]}$&$2^{4}\cdot5^{2}\cdot1201$&5&961&5764801\\
7&12&$1^{[6]},2^{[3]}$&$2^{3}\cdot5^{2}\cdot13\cdot19\cdot43\cdot181$&6&32193&13841287201\\
8&7&$1^{[7]}$&$42799^{\dagger}$&3&889&2097152\\
9&4&$1^{[4]}$&$5\cdot41$&2&45&6561\\
9&8&$1^{[8]}$&$5\cdot17\cdot41\cdot193$&4&3825&43046721\\
11&4&$1^{[2]},2$&$2^{2}\cdot3\cdot61$&3&49&14641\\
11&10&$1^{[10]}$&$2\cdot3\cdot3221\cdot13421$&6&64449&25937424601\\
13&4&$1^{[4]}$&$5\cdot7\cdot17$&3&105&28561\\
13&7&$1,2^{[3]}$&$5229043^{\dagger}$&5&465&62748517\\
\end{longtable}\endgroup

\begin{proposition}[The exact range of the basic bound]\label{v93:prop-arithmetic-complete}
For every prime power \(q\) and every \(n\ge3\), the basic bound~\eqref{v93:eq-easy} holds except for
\begin{equation}\label{v93:eq-eight-v2}
 (2,3),(2,6),(2,15),(3,4),(3,8),(5,4),(5,8),(7,6)
\end{equation}
and each of these eight pairs fails to satisfy that bound.
\end{proposition}
\begin{proof}
Section~\ref{v93:sec-uniform-ranges} treats characteristic-power degrees, cubic extensions other than \((2,3)\), all quintic extensions, and all \(q\ge16\). Table~\ref{v93:tab-families}, together with the argument for \(q=13,n\mid12\), reduces the remaining base fields to the 42 boundary pairs. Every unstarred row of Table~\ref{v93:tab-boundary} satisfies the basic bound; the starred rows and \((2,3)\) are exactly the eight pairs listed above. For the seven starred rows, the factorizations involve only small integers, the actual values of \(\omega(R)\) agree with the tabulated bounds, and direct comparison shows that the basic bound fails. For \((2,3)\), we have \(R=7,s=2\), and the right side of the basic bound is \(3>\sqrt8\), so this pair also fails to satisfy the bound.
\end{proof}

\section{Integer Initial Values for the Monotone Comparisons}
\label{v93:sec-integer-certificates}
The row numbers refer to Table~\ref{v93:tab-families}. Define
\[
 A=Ln_0(h/2-1),\quad B=Lh(r+\alpha n_0+\beta),\quad K=(q-1)\delta\prod_{\ell\in\mathcal P}\ell.
\]
Each quadruple in the table specifies the integer inequality \(q^A K^L>2^B\). The slope comparison~\eqref{v93:eq-family-slope} is checked with the same \(q,h,\alpha,L\). For example, row 8 gives
\(3^{162}\cdot 2024^{6}>2^{310}\), with successive ratio \(3^9/2^{10}>1\). The table records one initial value for each of the 33 families; no interval of extension degrees is checked term by term.
\par\medskip
\begingroup\small
\setlength{\tabcolsep}{5pt}\renewcommand{\arraystretch}{1.15}
\begin{longtable}{@{}rrrrrr@{}}
\caption{Each row asserts $q^A K^L>2^B$}\label{v93:tab-certificates}\\
\toprule No.&$q$&$A$&$K$&$L$&$B$\\\midrule\endfirsthead
\toprule No.&$q$&$A$&$K$&$L$&$B$\\\midrule\endhead
\bottomrule\endfoot
1&2&92&234577&2&126\\
2&2&116&234577&2&150\\
3&2&140&234577&2&174\\
4&2&52&231&2&64\\
5&2&40&15015&2&64\\
6&3&72&20020&6&192\\
7&3&66&22&6&128\\
8&3&162&2024&6&310\\
9&3&198&293480&6&420\\
10&3&252&4989160&6&530\\
11&4&24&45045&2&76\\
12&4&12&3&6&32\\
13&4&42&33&6&112\\
14&4&78&189&6&200\\
15&4&156&2079&6&376\\
16&5&27&25916&6&150\\
17&5&81&726010824&6&360\\
18&5&126&155661792&6&450\\
19&5&207&641793568416&6&710\\
20&7&18&660&2&68\\
21&7&12&54&2&44\\
22&7&45&1031453949240&2&205\\
23&8&10&15015&2&56\\
24&8&39&4913732008185&2&200\\
25&9&2&8008&2&32\\
26&9&8&80080&2&56\\
27&9&12&160160&2&72\\
28&9&30&640640&2&132\\
29&11&2&70&2&16\\
30&11&10&10920&2&60\\
31&11&10&1750&2&56\\
32&11&32&273000&2&144\\
33&13&16&27720&2&88\\
\end{longtable}\endgroup

\section{Order Data for the Allowed Small Primes}
\label{v93:sec-small-prime-orders}
The table gives \(\ord_\ell(q)\); a dash means \(\ell\mid q\). Each entry is obtained by repeated multiplication or squaring modulo the prime. Combining these entries with the divisibility conditions in Table~\ref{v93:tab-families} and~\eqref{v93:eq-valuation-admissible} verifies the allowed sets row by row. When the order is one, the valuation formulas must be used: the condition \(\ord_\ell(q)\mid n\) alone does not determine whether \(\ell\) divides \(R\) after cancellation of the denominator.
\par\medskip\begin{table}[htbp]\centering\scriptsize
\caption{Multiplicative orders modulo small primes}\label{v93:tab-orders}
\setlength{\tabcolsep}{5pt}
\begin{tabular}{@{}r|rrrrrrrrrrr@{}}\toprule
$q\backslash\ell$&2&3&5&7&11&13&17&19&23&29&31\\\midrule
2&--&2&4&3&10&12&8&18&11&28&5\\
3&1&--&4&6&5&3&16&18&11&28&30\\
4&--&1&2&3&5&6&4&9&11&14&5\\
5&1&2&--&6&5&4&16&9&22&14&3\\
7&1&1&4&--&10&12&16&3&22&7&15\\
8&--&2&4&1&10&4&8&6&11&28&5\\
9&1&--&2&3&5&3&8&9&11&14&15\\
11&1&2&1&3&--&12&16&3&22&28&30\\
13&1&1&4&2&10&--&4&18&11&14&30\\
\bottomrule\end{tabular}\end{table}\par

For odd degrees over the binary field, the primes greater than 31 and less than 64 must also be checked. Their orders are
\[
\begin{array}{c|rrrrrrr}
\ell&37&41&43&47&53&59&61\\\hline
\ord_\ell(2)&36&20&14&23&52&58&60
\end{array}
\]
Thus only 47 must be added to the allowed set.

\chapter{Solutions to the Exercises}
\Needspace{7\baselineskip}
\section*{Solutions to Chapter~\ref{chap:foundations}}
\pdfbookmark[1]{Solutions to Chapter \getrefnumber{chap:foundations}}{solutions.chap:foundations}

\solutionheading{ex:orders}
\begin{solution}
Neither \(12\) nor \(36\) is a prime power, so neither is possible.
For \(16=2^4\), the possible characteristic is \(2\) and the dimension over the prime field is \(4\);
for \(27=3^3\), the possible characteristic is \(3\) and the dimension is \(3\).
The order and characteristic constraints give only necessary conditions; the existence of fields with \(16\) and \(27\) elements follows from Theorem~\ref{thm:root-construction}.
\end{solution}

\solutionheading{ex:additive}
\begin{solution}
Choose an \(\F_p\)-basis \(e_1,\ldots,e_n\). The map \((a_1,\ldots,a_n)\mapsto\sum_i a_ie_i\) is an additive-group isomorphism.
Every nonzero element has additive order \(p\). If \(n>1\), the group has \(p^n>p\) elements, so no element generates the whole group.
\end{solution}

\solutionheading{ex:tower}
\begin{solution}
A subfield must have order \(2^d\) with \(d\mid12\), so the candidates are \(2,4,8,16,64,4096\).
The tower formula gives only a necessary condition; Chapter~\ref{chap:frobenius} proves that each candidate indeed corresponds to a unique subfield.
\end{solution}

\solutionheading{ex:wu-ring-power}
\begin{solution}
The integer congruences \(a^3\equiv a\pmod2\) and \(a^3\equiv a\pmod3\) follow from the power identities in the two prime fields, respectively;
the integer Chinese remainder theorem gives \(a^3\equiv a\pmod6\).
However, \(\bar2,\bar3\) are both nonzero and their product is zero, so this ring is not a field.
For the second question, if \(R\) is nonzero, every nonzero \(x\) has inverse \(x^{m-1}\), so \(R\) is a field.
The definition of a field excludes the zero ring; if \(0=1\) is allowed, the condition also holds in the zero ring, which has no nonzero elements, so the zero ring is an exception.
\end{solution}

\solutionheading{ex:bezout}
\begin{solution}
In \(\F_3[X]\), we have \((X+1)(X+2)=X^2+2\); subtracting \(X^2+1\) gives \(1\).
Thus the inverse is \(\overline X+2\).
The values of \(X^2+1\) at \(0,1,2\) are \(1,2,2\), respectively, so it has no roots.
A reducible quadratic has a linear factor, so this quadratic is irreducible and the quotient ring is a field.
\end{solution}

\solutionheading{ex:reducible}
\begin{solution}
Write \(f=AB\) with \(0<\deg A,\deg B<\deg f\).
The classes \(\overline A,\overline B\) are nonzero, but \(\overline A\,\overline B=\overline f=0\).
Thus the quotient ring has nonzero zero divisors and is not a field.
\end{solution}

\solutionheading{ex:split-bound}
\begin{solution}
Denote the known root by \(a\in k\) and write \(P=(X-a)Q\). If \(D=1\), take \(L=k\), with the convention \(0!=1\).
If \(D>1\), apply the proof of the splitting-extension construction in Lemma~\ref{lem:finite-splitting} to the degree-\(D-1\) polynomial \(Q\) to obtain the required extension field and bound.
This bound asserts only the existence of such an extension field. If \(k=\F_q\) and \(P\) already splits in \(\F_{q^m}\),
it also splits in every \(\F_{q^{mt}}\), whose degrees over \(k\) can be arbitrarily large.
\end{solution}

\solutionheading{ex:roots-needs-splitting}
\begin{solution}
In \(\F_2\) there are only the two roots \(0,1\), so the root set is \(\F_2\) itself and is indeed a field.
But \(X^4-X\) has not yet split completely over this field.
To conclude that there are four roots, one needs both complete splitting in the chosen extension field and the absence of repeated roots.
\end{solution}

\solutionheading{ex:supp-crt}
\begin{solution}
Interpolation gives \(A=2X\). We may take \(e_1=(X+1)/2=2X+2\) and \(e_2=(1-X)/2=X+2\).
Modulo \(X^2-1\), they satisfy \(e_i^2=e_i\), \(e_1e_2=0\), \(e_1+e_2=1\), and \(2e_1+e_2=2X\).
\end{solution}

\solutionheading{ex:wu-rootfield}
\begin{solution}
A subfield is closed under multiplication. Conversely, if \(V\) is closed under multiplication, it is a finite subring with identity and has no zero divisors as a subring of a field.
For \(a\in V\setminus\{0\}\), the multiplication map \(x\mapsto ax\) is injective on \(V\), hence surjective by finiteness, so some \(x\in V\) satisfies \(ax=1\).
Thus \(V\) is a subfield, and the proposition on root polynomials of finite subfields gives the stated identity.
If that identity holds, \(V\) is the full fixed-point set of the \(p^n\)th-power map and is therefore a subfield.
\end{solution}

\section*{Solutions to Chapter~\ref{chap:generators}}
\pdfbookmark[1]{Solutions to Chapter \getrefnumber{chap:generators}}{solutions.chap:generators}

\solutionheading{ex:f8iso}
\begin{solution}
The polynomials \(f,g\) are nonzero at both \(0,1\). A reducible cubic has a linear factor, so both are irreducible.
Direct expansion in characteristic two gives
\[
g(\alpha+1)=(\alpha+1)^3+(\alpha+1)^2+1
=\alpha^3+\alpha+1=0.
\]
Thus the required isomorphism is
\[
a+b\beta+c\beta^2\longmapsto
a+b(\alpha+1)+c(\alpha+1)^2,
\qquad a,b,c\in\F_2.
\]
Its image contains \(\alpha+1\) and \(1\), hence \(\alpha\), so it is the entire target field.
\end{solution}

\solutionheading{ex:iso-condition}
\begin{solution}
If an isomorphism extends \(\tau\) and sends \(\alpha\) to \(\beta\), applying it to \(m(\alpha)=0\) necessarily gives \(m^\tau(\beta)=0\). If this condition fails, the map does not preserve the polynomial relation and cannot be a field homomorphism.
The one-step extension lemma proves that this necessary condition is also sufficient.
\end{solution}

\solutionheading{ex:finite-domain}
\begin{solution}
In a finite commutative integral domain \(R\), fix a nonzero \(a\). The map \(x\mapsto ax\) is injective:
\(ax=ay\) implies \(a(x-y)=0\) and hence \(x=y\). Finiteness makes it surjective, so there exists \(b\) with \(ab=1\).
Thus \(R\) is a field.
The root set contains \(0,1\), is closed under subtraction and multiplication, and has no zero divisors as a subring of a field;
the root bound also makes it finite, so the preceding result establishes closure under inverses at once.
\end{solution}

\solutionheading{ex:wu-embedding-compatibility}
\begin{solution}
A homomorphism preserving the relation \(\alpha^2+1=0\) must send \(\alpha\) to a root of \(X^2+1\).
Its two distinct roots are \(\alpha,-\alpha\), and the one-step extension lemma gives the corresponding homomorphisms.
The linear map \(a+b\alpha\mapsto a+b(1+\alpha)\) preserves addition, but \((1+\alpha)^2=2\alpha\ne-1\), so it does not preserve \(\alpha^2=-1\) and hence does not preserve multiplication.
\end{solution}

\solutionheading{ex:primitive-test}
\begin{solution}
Let \(h=\gcd(N,j)\), \(N=hN_0\), and \(j=hj_0\), where \(\gcd(N_0,j_0)=1\). Then
\[
(\zeta^j)^t=1\ \Longleftrightarrow\ N\mid jt
\ \Longleftrightarrow\ N_0\mid t.
\]
Thus the order is \(N_0=N/h\). The element \(\zeta^j\) generates \(G\) if and only if \(\gcd(N,j)=1\),
and there are \(\varphi(N)\) such residue classes modulo \(N\).
\end{solution}

\solutionheading{ex:power-eq}
\begin{solution}
Choose a primitive element \(\zeta\), write \(a=\zeta^u\) and \(x=\zeta^v\), and put \(N=Q-1\).
The equation becomes \(tv\equiv u\pmod N\). Let \(g=\gcd(t,N)\). This congruence is solvable if and only if \(g\mid u\):
necessity is clear; for sufficiency, the integer Bézout identity gives an inverse of \(t/g\) modulo \(N/g\).
When it is solvable, \(v\) is unique modulo \(N/g\), giving \(g\) solutions modulo \(N\).
On the other hand, \(a^{N/g}=1\) is equivalent to \(N\mid uN/g\), namely \(g\mid u\).
\end{solution}

\solutionheading{ex:cosets}
\begin{solution}
The set \(k^\times\) is a subgroup of \(E^\times\); intersecting cosets are equal, and each coset has \(q-1\) elements.
Choosing a representative \(z_i\) from each coset gives the stated disjoint union and counting identity. Hence
\[
s=\frac{q^m-1}{q-1}=1+q+\cdots+q^{m-1}.
\]
The root set of an irreducible polynomial is usually not closed under multiplication and does not even contain \(1\).
For example, the roots of \(X^2+X+1\) in the four-element field are \(\alpha,\alpha^2\); their product is \(1\), but \(1\) is not a root.
The multiplicative group of the subfield generated by these roots admits a coset decomposition; the root set of one irreducible polynomial is generally not a multiplicative subgroup.
\end{solution}

\solutionheading{ex:wu-powers-groups}
\begin{solution}
We have \(\alpha\notin\F_3\), so the element generates the quadratic extension, but its order is four.
Let \(\beta=1+\alpha\). Then \(\beta^2=2\alpha\) and \(\beta^4=2\), so \(\ord(\beta)=8\).
The primitive elements are precisely \(\beta,\beta^3,\beta^5,\beta^7\), explicitly \(1+\alpha,1-\alpha,-1-\alpha,-1+\alpha\).
The four cosets are \(\{1,-1\}\), \(\{\alpha,-\alpha\}\), \(\{1+\alpha,-1-\alpha\}\), and \(\{1-\alpha,-1+\alpha\}\).
\end{solution}

\solutionheading{ex:phi5}
\begin{solution}
We have \(h(0)=h(1)=1\), so there is no linear factor.
Modulo the unique monic irreducible quadratic \(X^2+X+1\), we have \(X^3\equiv1\) and \(X^4\equiv X\),
hence \(h\equiv X+1\ne0\). Thus \(h\) has no irreducible quadratic factor either and is an irreducible quartic.
Consequently \(|\F_2(\beta)|=2^4=16\). Also \((X-1)h=X^5-1\), so \(\beta^5=1\); since \(\beta\ne1\), its order is exactly five.
\end{solution}

\solutionheading{ex:inverse16}
\begin{solution}
From \(\alpha^4+\alpha=1\) we get \(\alpha(\alpha^3+1)=1\), so \(\alpha^{-1}=\alpha^3+1\). Direct calculation gives
\[
(\alpha+1)(\alpha^3+\alpha^2+\alpha)
=\alpha^4+\alpha=1,
\]
hence \((\alpha+1)^{-1}=\alpha^3+\alpha^2+\alpha\).
\end{solution}

\solutionheading{ex:primitive-polynomial}
\begin{solution}
Take \(E=\F_{q^n}\) and \(k=\F_q\), and choose a generator \(\zeta\) of \(E^\times\).
Since the powers of \(\zeta\) exhaust \(E^\times\), we have \(k(\zeta)=E\),
so its minimal polynomial has degree \([k(\zeta):k]=[E:k]=n\).
This polynomial has the primitive element \(\zeta\) as a root and is therefore a primitive polynomial of degree \(n\).
\end{solution}

\section*{Solutions to Chapter~\ref{chap:frobenius}}
\pdfbookmark[1]{Solutions to Chapter \getrefnumber{chap:frobenius}}{solutions.chap:frobenius}

\solutionheading{ex:subfields18}
\begin{solution}
The positive divisors of \(18\) are \(1,2,3,6,9,18\), so the subfields are exactly \(\F_{3^d}\) for these values of \(d\).
The intersection is \(\F_{3^3}\) and the compositum is \(\F_{3^{18}}\). Moreover,
\[
\Fix(\phi^{12})=\F_{3^6},\qquad
\Aut_{\F_{3^6}}(\F_{3^{18}})=\langle\phi^6\rangle
=\{\id,\phi^6,\phi^{12}\}.
\]
\end{solution}

\solutionheading{ex:embedding}
\begin{solution}
The image of an embedding is a subfield of the target containing \(\F_q\) and having \(q^a\) elements, so the subfield classification requires \(a\mid b\).
Conversely, if \(a\mid b\), the target has a unique subfield of \(q^a\) elements.
The source is isomorphic to it over \(\F_q\), so an embedding exists.
Fix an embedding \(\iota\). Every other embedding \(\jmath\) has the same image,
so \(\iota^{-1}\circ\jmath\) is an \(\F_q\)-automorphism of the source.
There are exactly \(a\) such automorphisms, hence exactly \(a\) embeddings.
\end{solution}

\solutionheading{ex:counting}
\begin{solution}
By \eqref{eq:mobius},
\[
I_2(6)=\frac{2^6-2^3-2^2+2}{6}=9.
\]
Each irreducible sextic contributes six distinct roots, so there are \(6\cdot9=54\) elements of degree six.
Alternatively, the proper divisors of six are one, two, and three, so the elements of smaller degree form \(\F_4\cup\F_8\), whose two fields intersect in \(\F_2\).
The complement has size \(64-(4+8-2)=54\). The first calculation groups elements into individual orbits; the second removes the layers with smaller periods from the fixed-point set.
The same partition by periods ensures that the two agree.
\end{solution}

\solutionheading{ex:factor-extension}
\begin{solution}
Choose a root \(\alpha\) of \(f\) in \(\F_{q^{\lcm(d,s)}}\), and denote all roots by \(\alpha_i=\alpha^{q^i}\), \(i\in\Z/d\Z\).
Frobenius relative to \(\F_{q^r}\) is \(\sigma^r\) and acts by \(i\mapsto i+r\),
giving \(g\) orbits, each of length \(m=d/g\). On one such orbit, use coordinates \(\alpha_{i+rj}\) with \(j\in\Z/m\Z\); then \(\sigma^s=(\sigma^r)^{s/r}\) corresponds to \(j\mapsto j+s/r\). Thus this orbit splits further into \(h=\gcd(m,s/r)\) orbits, each of length \(m/h=d/(gh)\).
The root-orbit theorem translates these orbit partitions into the corresponding irreducible factorizations.

Write \(d=gd_0,r=gr_0,s=rt\). Then \(\gcd(d_0,r_0)=1\), and hence
\[
\gcd(d,s)=g\gcd(d_0,r_0t)=g\gcd(d_0,t)=gh.
\]
Every \(\sigma^s\)-orbit lies in a unique \(\sigma^r\)-orbit; successive extensions therefore give not only the same number and degrees of factors,
but also the same actual partition of the roots.

When \(d=12,r=2,s=8\), the first stage consists of two six-cycles, on the even and odd indices:
\[
(0\ 2\ 4\ 6\ 8\ 10),\qquad(1\ 3\ 5\ 7\ 9\ 11).
\]
The second-stage step is eight, splitting them respectively into
\[
\begin{array}{ll}
(0\ 8\ 4),&(2\ 10\ 6),\\
(1\ 9\ 5),&(3\ 11\ 7).
\end{array}
\]
Thus there are first two sextic factors and then four cubic factors; direct factorization over \(\F_{q^8}\) gives these same four orbits.
\end{solution}

\solutionheading{ex:wu-affine-degree}
\begin{solution}
A field containing \(\alpha\) contains \(a\alpha+b\), and \(\alpha=a^{-1}((a\alpha+b)-b)\), so the two generated fields are equal.
If \(f=m_{\alpha,k}\), then
\[
g(X)=a^d f((X-b)/a)
\]
is monic of degree \(d\) and vanishes at \(a\alpha+b\).
This element also has degree \(d\), so \(g\) is its minimal polynomial.
In \(\F_9\), the order of \(\alpha\) is four and that of \(1+\alpha\) is eight, showing that the order need not be preserved.
\end{solution}

\solutionheading{ex:multivariate-normalization}
\begin{solution}
By the subfield classification, \(s\mid N\). The distinct conjugates of \(L\) are \(L,\sigma(L),\ldots,\sigma^{s-1}(L)\),
and each occurs \(N/s\) times in the product of length \(N\); the product of all constant factors is
\(\prod_{i=0}^{N-1}u^{q^i}=u^{(q^N-1)/(q-1)}\). Multiplying these two contributions gives the formula.
In the example \(M^{(2)}=\omega^2(X+Y)\), so
\[
 MM^{(2)}=\omega^3(X+Y)^2=(X+Y)^2.
\]
The coefficients of \(M\) generate \(\F_4\), but its normalized factor is \(X+Y\), whose coefficients already lie in \(\F_2\).
The two coefficient conjugates are associates, so their product has repeated factors.
Normalization removes the effect of constant multiples, making the coefficient-orbit length equal to the orbit length of distinct linear factors.
\end{solution}

\solutionheading{ex:wu2-monogenic}
\begin{solution}
If \(a=(a_{d,j})\) generates the whole product, its projection to each component generates that field, so the minimal polynomial of \(a_{d,j}\) has degree \(d\).
These minimal polynomials must be pairwise distinct. If two components had the same minimal polynomial \(P\) of degree \(d\),
the evaluation map from \(k[X]\) to their product would have kernel \((P)\), and its image would have only \(q^d\) elements, preventing surjectivity onto a product of \(q^{2d}\) elements.
Thus each degree requires \(m_d\) distinct irreducible polynomials, proving necessity.
Conversely, choose the required distinct irreducible polynomials of each degree and choose a root in each component as \(a_{d,j}\).
The Chinese remainder theorem makes \(P\mapsto(P(a_{d,j}))\) surjective, so \(a\) generates the whole algebra.
Since \(I_2(1)=2<3\), the algebra \(\F_2^3\) is not generated by one element; every quotient algebra of \(\F_2[X]\) is generated by \(\overline X\), giving the conclusion.
\end{solution}

\solutionheading{ex:wu3-algebra-hom}
\begin{solution}
Specifying a homomorphism \(A\to B\) is equivalent to independently specifying its component homomorphisms \(A\to\F_{q^{e_j}}\). The latter must select one field factor of \(A\):
the coordinate idempotents are pairwise orthogonal and sum to one, and their images can only be zero or one. Orthogonality prevents two images from both being one, while their sum ensures that at least one is one.
A factor \(\F_{q^{d_i}}\) embeds in the target if and only if \(d_i\mid e_j\), in which case there are exactly \(d_i\) embeddings.
Summing over the choice of factor and multiplying over target coordinates gives the formula; an empty sum is zero, meaning that no homomorphism exists.
Equivalently, a cycle of length \(e_j\) maps equivariantly into one of length \(d_i\) if and only if \(d_i\mid e_j\);
then the first point has \(d_i\) possible images, determining all the others.
An automorphism can only permute factors of the same degree; for degree \(d\) there are \(m_d!\) permutations,
and \(d\) field automorphisms for each factor, giving the final product.
\end{solution}

\solutionheading{ex:wu3-algebra-fixed}
\begin{solution}
The map \(\sigma_A\) raises each coordinate to the \(q\)th power; in every field factor the fixed field of \(\sigma_A^r\) is \(\F_{q^{\gcd(d_i,r)}}\), giving the stated fixed algebra and dimension.
On evaluation points, a cycle of length \(d_i\) splits into \(\gcd(d_i,r)\) cycles under step \(r\).
Its number of fixed points is \(d_i\) if \(d_i\mid r\) and zero otherwise.
For example, the \(q\)-Frobenius of a cubic field has a one-dimensional fixed field, but its three evaluation points form a three-cycle with no fixed point.
\end{solution}

\solutionheading{fra:ex-five-dimensions}
\begin{solution}
The two inversions give
\[
 (B_1,B_2,B_3,B_4,B_5)=(2,1,1,0,0),
 \qquad a_2=a_3=1,
\]
and all other \(a_j=0\). For example, \(B_3=(b_3-b_1)/2=1\) and \(B_4=(b_4-b_2)/2=0\), whereas \(a_1=B_1-B_2-B_3-B_5=0\).
Therefore \(A\simeq\F_{q^2}\times\F_{q^3}\).
\end{solution}

\solutionheading{fra:ex-two-counts}
\begin{solution}
On the root set or evaluation-point set, an orbit of \(P_j\) of length \(j\) is fixed pointwise exactly when \(j\mid r\), so
\[
 z_r=2\,\mathbf1_{2\mid r}+3\,\mathbf1_{3\mid r}.
\]
On algebra elements, \(A\simeq\F_{q^2}\times\F_{q^3}\); the fixed fields of the \(r\)th iterate on the field components have dimensions \(\gcd(2,r)\) and \(\gcd(3,r)\), respectively, so
\[
 b_r(A)=\gcd(2,r)+\gcd(3,r).
\]
In particular, \(z_1=0\) and \(b_1(A)=2\). There are no fixed evaluation points, whereas the algebra elements constant along each orbit form a two-dimensional space; these concern different Frobenius actions.
\end{solution}

\section*{Solutions to Chapter~\ref{chap:orders}}
\pdfbookmark[1]{Solutions to Chapter \getrefnumber{chap:orders}}{solutions.chap:orders}

\solutionheading{ex:order12}
\begin{solution}
The roots of \(X^2+X+1\) are the nontrivial cube roots of unity, so this factor has order three; \(X+1\) has order one.
The least common multiple of the distinct factor orders is three. The largest multiplicity is three, and the least power of two at least three is four, so the required order is \(3\cdot4=12\).
One can also see directly that \(X^{12}-1=(X^3-1)^4\) contains the required factors, while the odd part of every possible exponent must be divisible by three and its power-of-two part must be at least four.
\end{solution}

\solutionheading{ex:order-degree}
\begin{solution}
Modulo 13, \(2^6\equiv-1\) and \(2^4\equiv3\), so \(\ord_{13}(2)=12\).
For an element \(a\) of order \(d\) in any group, \(a^s\) has order \(d/\gcd(d,s)\).
Thus the three degrees are, respectively,
\[
\ord_{13}(2)=12,\qquad \ord_{13}(4)=6,\qquad \ord_{13}(8)=4.
\]
\end{solution}

\solutionheading{ex:reciprocal3}
\begin{solution}
Reversing the coefficients gives \(X^3+X^2+1\), different from the original polynomial.
Structurally, an irreducible self-reciprocal polynomial other than a linear polynomial must have even degree, whereas the original polynomial is an irreducible cubic.
The coefficient sequence of \(X^4+X^3+X^2+X+1\) is palindromic, so it is self-reciprocal; its roots are the nontrivial fifth roots of unity, and inversion stays in the same Frobenius orbit.
\end{solution}

\solutionheading{ex:primitive4count}
\begin{solution}
There are \(\varphi(15)/4=8/4=2\) primitive polynomials and \((2^4-2^2)/4=3\) irreducible quartics in total. The two primitive polynomials are \(X^4+X+1\) and \(X^4+X^3+1\);
the remaining polynomial \(X^4+X^3+X^2+X+1\) has roots of order five and is not primitive.
\end{solution}

\solutionheading{ex:supp-projective}
\begin{solution}
We have \(X^2\equiv-1\pmod f\), so the projective order is two and the ordinary order is four; a primitive polynomial of degree two would have order eight, so this polynomial is not primitive.
In general, if \(\bar X\) has order \(q^d-1\), its powers exhaust all nonzero elements of the quotient ring.
Every nonzero element is therefore invertible, so the quotient ring is a field and \(f\) must be irreducible.
\end{solution}

\solutionheading{ex:wu-matrix-polynomial}
\begin{solution}
For every \(P\in k[T]\), we have \(P(M_a)(x)=P(a)x\), so \(P(M_a)=0\) if and only if \(P(a)=0\): for necessity, take \(x=1\).
The two monic annihilating polynomials of least degree are therefore equal.
By contrast, \(P(\sigma)(a)=\sum_i c_i a^{q^i}\) describes linear relations among conjugate elements.
For example, if \(\alpha^3+\alpha+1=0\) in \(\F_8\), the element has minimal polynomial \(m_\alpha=T^3+T+1\).
From \(\alpha^4=\alpha^2+\alpha\) we obtain \((\sigma^2+\sigma+\id)(\alpha)=0\).
Also, \(\alpha\ne0\) and \(\alpha^2\ne\alpha\), so none of \(1,T,T+1\) annihilates \(\alpha\); hence its minimal annihilating polynomial with respect to \(\sigma\) is \(\mu_\alpha=T^2+T+1\).
\end{solution}

\section*{Solutions to Chapter~\ref{chap:cyclotomic}}
\pdfbookmark[1]{Solutions to Chapter \getrefnumber{chap:cyclotomic}}{solutions.chap:cyclotomic}

\solutionheading{ex:phi9}
\begin{solution}
The 2-cyclotomic cosets modulo nine are
\[
\{0\},\qquad\{1,2,4,8,7,5\},\qquad\{3,6\}.
\]
The roots of \(\Phi_9=X^6+X^3+1\) correspond to the middle six-element orbit, so it is irreducible over \(\F_2\).
Over \(\F_4\), Frobenius has step two on the original orbit, splitting it into two three-element orbits; hence the polynomial is a product of two irreducible cubics.
Over \(\F_8\) the step is three, giving a product of three irreducible quadratics.
\end{solution}

\solutionheading{ex:cyclop}
\begin{solution}
The integer polynomial is \(\Phi_{12}=X^4-X^2+1\). Thus
\[
\overline{\Phi_{12}}=
\begin{cases}
(X^2+X+1)^2,&\text{in }\F_2[X],\\
(X^2+1)^2,&\text{in }\F_3[X].
\end{cases}
\]
The first case removes the power-of-two part of 12, leaving primitive cube roots of unity; the second leaves primitive fourth roots of unity.
Proposition~\ref{prop:cyclotomic-p} gives multiplicities \(2^{2-1}(2-1)=2\) and \(3^{1-1}(3-1)=2\), respectively.
\end{solution}

\solutionheading{ex:coset-negative}
\begin{solution}
The roots of the monic reciprocal polynomial are the inverses of the original roots; since \((\zeta^i)^{-1}=\zeta^{-i}\), the orbit \(C\) becomes \(-C\).
Modulo 15, each of \(C_0,C_3,C_5\) is invariant under negation, so the corresponding polynomials \(X+1\), \(\Phi_5\), and \(\Phi_3\) are self-reciprocal.
The cosets \(C_1\) and \(C_7\) are interchanged, so the two primitive quartic factors are reciprocal to each other.
\end{solution}

\solutionheading{ex:degree5extension}
\begin{solution}
By Theorem~\ref{thm:factor-base-extension}, the numbers of factors are \(\gcd(5,6)=1\) and \(\gcd(5,10)=5\), respectively.
Thus the polynomial remains irreducible over the first field and splits into five distinct linear factors over the second.
\end{solution}

\solutionheading{ex:wu-splitting-degree}
\begin{solution}
Put \(m=\lcm_i d_i\). By \(d_i\mid m\), every \(f_i\) divides \(X^{q^m}-X\), so all roots lie in \(\F_{q^m}\).
Let \(\F_{q^s}\) be the splitting field generated by these roots within \(\F_{q^m}\). Then \(s\mid m\).
The degree \(d_i\) of each root also divides \(s\), so \(m\mid s\) and hence \(s=m\).
Multiplicities do not change the set of roots and therefore do not change the splitting field.
\end{solution}

\solutionheading{ex:wu2-factor-coefficients}
\begin{solution}
Since \(\gcd(6,4)=2\), there are two cubic factors over \(\F_{q^4}\), and the coefficients of each generate \(\F_{q^2}\).
Since \(\gcd(6,9)=3\), there are three quadratic factors over \(\F_{q^9}\), each with exact coefficient field \(\F_{q^3}\).
Over \(\F_{q^{12}}\), the polynomial \(f\) splits into six linear factors \(X-\alpha^{q^j}\), whose coefficients each generate \(\F_{q^6}\).
These conclusions use both Theorem~\ref{thm:factor-base-extension} and Proposition~\ref{prop:factor-coefficient-field}.
The coefficients jointly generate a field; each coefficient need not generate it individually. For example, in Example~\ref{exm:base-extension-four} the coefficient field of \(X^2+X+\omega\) is \(\F_4\), but the coefficient of \(X\) is one and belongs to \(\F_2\).
\end{solution}

\solutionheading{ex:med-seventytwo}
\begin{solution}
The multiplicative order of five is two modulo eight and six modulo nine, so \(\tau=\lcm(2,6)=6\).
Since \(\rad(72)=6\), we have \(\rho=\ord_6(5)=2\). Although \(8\mid72\), we have \(5^2\equiv1\pmod4\), so \(\omega_C=2\) and \(\kappa_C=3\).
Successive multiplication by five gives
\[
C_1=\{1,5,25,53,49,29\}
=\{1,25,49\}\sqcup\{5,29,53\}.
\]
Both blocks have common difference 24; the allowed block sizes are one and three.
Theorem~\ref{prop:med-binomial-descent} gives \(M_{C_1}=G(X^3)\), where \(G\) is the quadratic minimal polynomial of \(\alpha^3\).
The least extension degree is two. Over \(\F_{25}\), the factorization consists of the two irreducible cubic binomials
\(X^3-\alpha^3\) and \(X^3-\alpha^{15}\).
\end{solution}

\solutionheading{ex:med-boundaries}
\begin{solution}
If \(d\) is square-free, then \(\rad(d)=d\), so \(\rho=\ord_d(q)=\tau\) and \(8\nmid d\).
The least-step formula gives \(\omega_C=\tau\) and \(\kappa_C=1\), so only the singleton decomposition is possible.
For the degree-eight polynomial, extending the base field gives \(\gcd(8,m)\) factors, each of degree \(8/\gcd(8,m)\).
The values \(m=2,3,4\) give two quartic factors, one degree-eight factor, and four quadratic factors, respectively.
Also \(\omega_C=2\), so all factors are binomials when \(m=2,4\), but not when \(m=3\).
Thus the original orbit length determines the factor degrees, while complete equal-difference blocks further determine the binomial form of the factors.
\end{solution}

\solutionheading{ex:supp-binomial-factors}
\begin{solution}
Apply the global divisibility criterion in Corollary~\ref{cor:med-all-binomial}.
Modulo 16 we have \(\ord_{16}(5)=4\), and \(16\mid4(5-1)\), so all factors over \(\F_5\) are binomials.
But \(\ord_{16}(3)=4\) and \(16\nmid4(3-1)\), so this is not true of all factors over \(\F_3\).
\end{solution}

\solutionheading{ex:med-eighteen-five}
\begin{solution}
Successive multiplication by five gives
\[
\begin{aligned}
C_0&=\{0\},&C_1&=\{1,5,7,11,13,17\},\\
C_2&=\{2,4,8,10,14,16\},&C_3&=\{3,15\},\\
C_6&=\{6,12\},&C_9&=\{9\}.
\end{aligned}
\]
Their lengths are \(1,6,6,2,2,1\), respectively, with sum 18.
For \(C_1,C_2\), the reduced moduli are 18 and nine, respectively; both have \(\omega_C=\ord_3(5)=2\) and \(\kappa_C=3\). The coarsest decompositions are
\[
\begin{aligned}
C_1&=\{1,7,13\}\sqcup\{5,11,17\},\\
C_2&=\{2,8,14\}\sqcup\{4,10,16\}.
\end{aligned}
\]
Each block has common difference six. The reduced moduli for \(C_3,C_6\) are six and three, respectively; both have \(\tau=\omega_C=2\) and \(\kappa_C=1\), so only singleton decompositions are allowed.
All three parameters of \(C_0,C_9\) equal one. In particular, although \(C_3=\{3,15\}\) has only two elements, it is not a complete equal-difference block: such a two-element block modulo 18 must have common difference nine, whereas the difference here is 12.
\end{solution}

\solutionheading{ex:med-twentyfive}
\begin{solution}
The order of traversal under multiplication by four is
\[
1\mapsto4\mapsto16\mapsto14\mapsto6\mapsto24\mapsto21\mapsto9\mapsto11\mapsto19\mapsto1.
\]
Thus \(\tau=10\), \(\rho=\ord_5(4)=2\), \(\omega_C=2\), and \(\kappa_C=5\). The coarsest decomposition is
\[
C_1=\{1,6,11,16,21\}\sqcup\{4,9,14,19,24\}.
\]
Both sets are complete equal-difference blocks with common difference five, and multiplication by \(4^2\) traverses each block.
Let \(\gamma=\zeta^5\). It has order five and degree two over \(\F_4\), so
\[
M_{C_1}(X)=(X^5-\gamma)(X^5-\gamma^4)
\quad\text{in }\F_{16}[X],
\]
both factors are irreducible, and the least required extension degree is two.
The orbit of index five is \(C_5=\{5,20\}\), with reduced modulus five, so \(\tau=\omega_{C_5}=2\) and \(\kappa_{C_5}=1\).
Its binomial factorization over \(\F_{16}\) is already a factorization into linear factors.
\end{solution}

\solutionheading{v40:ex-unit-fixed-counts}
\begin{solution}
The group \(\F_{q^t}^{\times}\) is cyclic of order \(q^t-1\), and the equation \(x^N=1\) has exactly
\(\gcd(N,q^t-1)\) solutions. This proves the first identity.
The polynomial \(X^N-1\) has no repeated roots, and Frobenius orbits partition its roots; an irreducible factor of degree \(r\) contributes an orbit of \(r\) elements.
This orbit is fixed pointwise by \(\sigma^t\) if and only if \(r\mid t\). Hence
\[
 \gcd(N,q^t-1)=\sum_{r\mid t}ra_r.
\]
Möbius inversion over the integers gives the first counting formula.
On the other hand, \(\mu_N\) has exactly \(\varphi(d)\) roots of multiplicative order \(d\);
by Proposition~\ref{prop:degree-from-order}, their Frobenius period is \(\operatorname{ord}_d(q)\).
Summing the numbers of roots of period \(r\) and dividing by \(r\) gives the second counting formula.

When \(q=2,N=15\), we have \(\operatorname{ord}_{15}(2)=4\), so the only possible orbit lengths are one, two, and four.
The corresponding fixed-point counts are
\[
 \gcd(15,2-1)=1,\qquad
 \gcd(15,2^2-1)=3,\qquad
 \gcd(15,2^4-1)=15.
\]
Consequently,
\[
 a_1=1,\qquad a_2=\frac{3-1}{2}=1,
 \qquad a_4=\frac{15-3}{4}=3,
\]
all other \(a_r\) are zero, and \(a_1+2a_2+4a_4=15\).

The Chinese remainder theorem gives
\[
 A\simeq\prod_{r\ge1}\F_{q^r}^{\,a_r}.
\]
On each field component, \(\sigma_A\) is Frobenius relative to \(k\); the fixed field of its \(t\)th iterate is
\(\F_{q^{\gcd(r,t)}}\), of dimension \(\gcd(r,t)\) over \(k\).
Summing over components gives the required dimension formula. In the example above, when \(t=1\) the root set has one fixed point,
whereas the fixed-space dimension is \(a_1+a_2+a_4=5\).
The former counts only roots fixed pointwise, namely orbits of length one; the latter retains a copy of the base field in each field component,
and thus counts all root orbits. Fixed-point counts on a root set must not be confused with dimensions of linear fixed spaces in an algebra.
\end{solution}

\solutionheading{v40:ex-uniform-root-degrees}
\begin{solution}
Every nonidentity \(N\)-th root of unity has Frobenius period \(m\), and hence belongs to \(\F_{q^t}\) if and only if \(m\mid t\). The root one always belongs to the base field, so
\[
\mu_N\cap\F_{q^t}=\begin{cases}
\mu_N,&m\mid t,\\
\{1\},&m\nmid t.
\end{cases}
\]

When \(q=2,N=9\), the only prime divisor is three, and \(\operatorname{ord}_3(2)=2\), but \(9\nmid2^2-1\).
Explicitly,
\[
 S_9=(X^2+X+1)(X^6+X^3+1).
\]
The two factors are the reductions of \(\Phi_3\) and \(\Phi_9\) over \(\F_2\), respectively;
\(\operatorname{ord}_3(2)=2\) and \(\operatorname{ord}_9(2)=6\) show that both are irreducible.
Their degrees differ, showing that prime-power moduli can further lengthen the Frobenius period.

Finally, \(9^3-1=728=8\cdot91\), and \(9\) is not one modulo either seven or 13.
The orders modulo both numbers divide the prime three and are not one, so
\(\operatorname{ord}_7(9)=\operatorname{ord}_{13}(9)=3\).
Proposition~\ref{v124:prop-uniform-cyclotomic} gives \((91-1)/3=30\) cubic factors of \(S_{91}\).
The three layers of roots with multiplicative orders seven, 13, and 91 contribute, respectively,
\[
 \frac{\varphi(7)}3=2,\qquad
 \frac{\varphi(13)}3=4,\qquad
 \frac{\varphi(91)}3=24
\]
factors. These layers have different multiplicative orders, but every nontrivial root belongs to
\(\F_{9^3}\setminus\F_9\) and has the same Frobenius period.

Finally, take \(N=q^s-1\), with \(s\ge2\). Since \(0<q^j-1<N\) for \(1\le j<s\), one has \(\ord_N(q)=s\). If \(q>2\), then \(\F_q^\times\) contains a nonidentity \(N\)-th root of degree one, ruling out a common degree. If \(q=2\) and \(s\) is composite, choose a proper divisor \(1<t<s\) of \(s\); a root of order \(2^t-1\) has degree \(t\), again ruling out a common degree. If \(q=2\) and \(s\) is prime, every root degree divides \(s\), while \(\F_2^\times\) contains only one. Thus every nonidentity root has degree \(s\). The required pairs are exactly \(q=2\) with \(s\) prime.
\end{solution}

\solutionheading{v44:ex-local-phases-315}
\begin{solution}
Using \(\ord_3(2)=2\) and \(v_3(2^2-1)=1\),
formula~\eqref{v44:eq-prime-power-order} gives \(\ord_9(2)=6\).
The six units modulo nine form one orbit, represented by one.
Modulo five, \(\ord_5(2)=4\), so the four units also form one orbit, represented by one.
The units modulo seven split into
\[
 \{1,2,4\},\qquad\{3,6,5\},
\]
each of length three, represented by one and three, respectively.

For either unit orbit modulo seven, the local lengths are \((6,4,3)\), so
\[
 L_1=6,\qquad L_2=12,\qquad L_3=12,
 \qquad\delta_2=2,\quad\delta_3=3.
\]
The recursive phase representatives are \((0,j_2,j_3)\), where \(0\le j_2<2\) and \(0\le j_3<3\).
They give six global orbits, each of length 12.
Each of the two unit orbits modulo seven contributes six, so the primitive 315th roots of unity form 12 Frobenius orbits.
Their number also satisfies \(12\cdot12=144=\varphi(315)\).
By Corollary~\ref{v44:cor-crt-cyclotomic-factorization}, the polynomial \(\Phi_{315}\) splits into exactly 12 irreducible factors of degree 12.

The CRT coefficients \(280,126,225\) equal one in their respective components modulo nine, five, and seven, and zero in the other two components.
Thus the three local components of the stated expression are \(1,2^{j_2},b2^{j_3}\).
Proposition~\ref{v44:prop-recursive-phase-representatives} ensures that these outputs represent all unit orbits without repetition.
Ordered first by \(b=1,3\), then by \(j_2=0,1\) and \(j_3=0,1,2\), the 12 representatives are
\[
 1,226,46,127,37,172,136,181,271,262,307,82.
\]
These are orbit representatives, not necessarily the least nonnegative elements of the orbits.
\end{solution}

\solutionheading{v124:ex-local-coefficient-fields}
\begin{solution}
One has \(P_8=X^4+1\) and \(P_3=X^2+X+1\), so both have exact coefficient field \(\F_{13}\). Also, \(H=\langle5,7\rangle\) is a subgroup of the units modulo twenty-four, and its reductions modulo eight and three are the full unit groups.

However, \(13\notin H\), whereas \(13^2\equiv1\pmod{24}\). Proposition~\ref{v50:prop-multiplier-blocks} therefore gives \(\F_{13^2}\) as the exact coefficient field of \(B\). Its conjugate is
\[
B^{[13]}(X)=\prod_{a\in\{13,17,19,23\}}(X-\zeta^a).
\]
The two exponent sets partition the units modulo twenty-four, so
\[
B(X)B^{[13]}(X)=\Phi_{24}(X)=X^8-X^4+1.
\]
Membership in both local images does not ensure membership in the global subgroup: \(13\bmod8\in H_8\) and \(13\bmod3\in H_3\), but \(13\bmod24\notin H\). The difference in coefficient fields records this compatibility condition.
\end{solution}

\solutionheading{v44:ex-local-orbits-81}
\begin{solution}
Here \(\ell=3\), \(h=\ord_3(4)=1\), and \(b=v_3(4-1)=1\), so
\[
 \ord_{3^a}(4)=3^{a-1}\qquad(a\ge1).
\]
Modulo \(81\), indices of valuation \(y\in\{0,1,2,3\}\) reduce to the unit group modulo \(3^{4-y}\);
their orbit length is \(3^{3-y}\), and the number of orbits is
\[
 \frac{\varphi(3^{4-y})}{3^{3-y}}=2.
\]
Adding the singleton orbit of the zero index gives nine orbits in total.
Thus the distribution of irreducible-factor degrees is
\[
 1,1,1,3,3,9,9,27,27,
\]
with total degree \(3+6+18+54=81\).

Take a primitive 81st root of unity \(\zeta\) and put \(\eta_a=\zeta^{3^{3-a}}\). Then \(\eta_a\) has order \(3^{a+1}\).
Thus \(\eta_a^{3^a}=\zeta^{27}\) is a primitive cube root of unity, equal to \(\omega\) or \(\omega^2\).
The degree of \(\eta_a\) over \(\F_4\) is
\(\ord_{3^{a+1}}(4)=3^a\).
Thus the binomial of degree \(3^a\) with root \(\eta_a\),
\(X^{3^a}-\eta_a^{3^a}\), is its minimal polynomial.
Using \(\eta_a^{-1}\) as a root similarly gives the other irreducible binomial, with the other primitive cube root of unity as its constant.

All roots of the polynomials \(X^{3^a}-\omega\) and \(X^{3^a}-\omega^2\) have order \(3^{a+1}\):
a root \(\alpha\) satisfies \(\alpha^{3^{a+1}}=1\), while \(\alpha^{3^a}\ne1\).
Different values of \(a\) give different root orders, and the two binomials for the same \(a\) have no common root.
Including the root one, the sum of the factor degrees is 81, giving exactly the stated complete irreducible factorization.

The primitive 81st roots correspond to \(a=3\), so
\[
 \Phi_{81}(X)=(X^{27}-\omega)(X^{27}-\omega^2)
 =X^{54}+X^{27}+1.
\]
All indices in these two orbits are units modulo 81, namely the layer of valuation \(y=0\);
the other three layers correspond to roots of unity of orders 27, nine, and three, respectively.
\end{solution}

\solutionheading{v44:ex-cyclotomic-return-factors}
\begin{solution}
The parent orbit has length \(\tau=2\). On lifting modulo nine, \(r=3\), and the return map is
\[
T(k)=4k+\frac{4-1}{3}\equiv k+1\pmod3.
\]
It has a single three-cycle, so the unique lift is
\[
C_9(1)=\{1,2,4,8,7,5\},\qquad |C_9(1)|=6.
\]
The parent factor is \(M_{C_3(1),3}=X^2+X+1\), and the compatible-factor identity gives
\[
M_{C_9(1),9}(X)=M_{C_3(1),3}(X^3)=X^6+X^3+1.
\]
It corresponds to a single complete Frobenius orbit and is therefore irreducible.

On lifting modulo 15, \(r=5\), and the return map is \(T(k)=4k+1\pmod5\), with orbits
\[
\{0,1\},\qquad\{2,4\},\qquad\{3\}.
\]
The admissible digits \(0,2,3\) give compatible representatives one, seven, and ten, respectively; the child orbits are
\[
C_{15}(1)=\{1,2,4,8\},\quad
C_{15}(7)=\{7,14,13,11\},\quad
C_{15}(10)=\{10,5\}.
\]
Their lengths are four, four, and two, respectively. From the factorization for 15th roots of unity in this chapter, or by direct multiplication and irreducibility checks, we obtain
\[
X^{10}+X^5+1
=(X^2+X+1)(X^4+X+1)(X^4+X^3+1).
\]
To check irreducibility independently, neither quartic has a root in \(\F_2\); the only monic irreducible quadratic is \(X^2+X+1\),
and the remainders of \(X^4+X+1\) and \(X^4+X^3+1\) modulo it are one and \(X\), respectively. Thus both quartics are irreducible.
Their correspondence with the specific four-element exponent orbits depends on the chosen compatible primitive roots of unity; the factorization itself and the factor degrees do not.
\end{solution}

\solutionheading{v44:ex-cyclotomic-compatible-leaders}
\begin{solution}
We have \(C_9(1)=\{1,2,4,8,7,5\}\), so \(\tau=6\). Here \(r=5\), and the return map is
\[
T(k)=2^6 k+\frac{2^6-1}{9}\equiv4k+2\pmod5.
\]
Its orbits are \(\{0,2\}\), \(\{1\}\), and \(\{3,4\}\), so the admissible digits are \(0,1,3\) and the compatible representatives are
\[
1+9\cdot0=1,\qquad1+9\cdot1=10,\qquad1+9\cdot3=28.
\]
The child-orbit lengths are 12, six, and 12, respectively. Listing the orbits directly gives
\[
\begin{aligned}
C_{45}(1)&=\{1,2,4,8,16,32,19,38,31,17,34,23\},\\
C_{45}(10)&=\{10,20,40,35,25,5\},\\
C_{45}(28)&=\{28,11,22,44,43,41,37,29,13,26,7,14\}.
\end{aligned}
\]
Their least nonnegative representatives are one, five, and seven, respectively. The latter two reduce modulo nine to five and seven, not to the fixed parent representative one.
Nevertheless, five and seven still belong to \(C_9(1)\), so all three child orbits reduce to the parent orbit modulo nine, and orbit compatibility is unchanged.
The chosen representatives are no longer compatible under reduction; to continue using the original admissible-digit recursion, one must retain the compatible representatives one, ten, and 28.
\end{solution}

\solutionheading{v44:ex-binomial-six}
\begin{solution}
We have \(2^2=4\) and \(2^4=1\), so \(g=4\) and \(6g=24\).
The exponent set is
\[
 S_{6,4}=\{1,5,9,13,17,21\},
\]
with 5-cyclotomic orbits
\[
 \{1,5\},\qquad\{9,21\},\qquad\{13,17\}.
\]
Each has length two, so \(X^6-2\) has exactly three irreducible quadratic factors.

Since \(3^3=2\) in \(\F_5\), putting \(Y=X^2\) gives
\[
 X^6-2=(X^2-3)(X^4+3X^2+4)
 =(X^2-3)(X^2+X+2)(X^2-X+2).
\]
The nonzero squares in \(\F_5\) are \(1,4\), so \(3\) is not a square and \(X^2-3\) is irreducible.
Both other quadratics have discriminant \(1-8=3\), also a nonsquare, so they too are irreducible.

Frobenius on \(\F_5\) is the identity, so the unique fifth root of \(2\) is again \(2\).
By Corollary~\ref{v44:cor-binomial-multiplicities},
\[
 X^{30}-2=(X^6-2)^5
 =(X^2-3)^5(X^2+X+2)^5(X^2-X+2)^5.
\]
The three distinct irreducible factors all have degree two, determined by the orbit lengths in \(S_{6,4}\);
their common multiplicity five is determined by the characteristic-power part of the degree \(30=5\cdot6\).
\end{solution}

\solutionheading{v46:ex-cyclotomic-stable-binary}
\begin{solution}
Successive multiplication by seven gives
\[
\begin{aligned}
 C_2(1)&=\{1\},& C_4(1)&=\{1,3\},& C_8(1)&=\{1,7\},\\
 C_{16}(1)&=\{1,7\},& C_{32}(1)&=\{1,7,17,23\}.
\end{aligned}
\]
The lengths are one, two, two, two, and four, respectively.
This system first differs from the zero principal system modulo two, so \(u=0\).
Since \(Q=7\equiv3\pmod4\) and \(v_2(Q+1)=3\),
Theorem~\ref{v46:thm-cyclotomic-stable-growth} gives the stabilization index \(s=4\).
The length already doubles from modulus two to four, but splitting still occurs from four to eight and from eight to 16;
only from modulus 16 onward is the lift unique at every level, with the length doubling each time.

The nonzero squares in \(\F_7\) are \(1,2,4\), whereas the discriminant of \(f\) is
\(1-24=5\), a nonsquare, so \(f\) is irreducible.
Take a root \(\alpha\in\F_{49}\). Its conjugate is \(\alpha^7\), and the norm equals the constant coefficient, so
\[
 \alpha^8=\alpha\alpha^7=6=-1,
 \qquad \alpha^{16}=1.
\]
Thus \(\alpha\) has order 16.
Choose \(\zeta_{16}=\alpha\) and compatible primitive roots of unity along the chain of power-of-two moduli.
Then \(f\) corresponds to \(C_{16}(1)\) modulo 16, precisely the first stable level of the system above.
Corollary~\ref{v46:cor-cyclotomic-stable-factor-substitution} gives
\[
 f(X^{2^t})=X^{2^{t+1}}+X^{2^t}+6
\]
irreducible for every \(t\ge0\), of degree \(2^{t+1}\).
\end{solution}

\solutionheading{v133:ex-no-section}
\begin{solution}
For an additive map $s$, $ps(a)=s(pa)=0$, so absence of $p$-torsion
forces $s(a)=0$. A unital ring map would require $s(1)=1\ne0$, a contradiction.
The constant embedding is evidently a ring section of $\F_q[[t]]\to\F_q$.
For another section, write $s(a)=a+u$, $u\in t\F_q[[t]]$.
Since $s(a)^q=s(a)$, one has $u^q=u$. A nonzero $u$ would have
$v_t(u^q)>v_t(u)$, which is impossible.
\end{solution}

\solutionheading{v133:ex-cocycle}
\begin{solution}
Modulo $p^2$, formula~\eqref{v133:eq-first-carry} gives
\[
 ([a]+[b])+[c]\equiv[a+b+c]+p[\kappa(a,b)+\kappa(a+b,c)].
\]
Changing parentheses gives
$[a]+([b]+[c])\equiv[a+b+c]+p[\kappa(b,c)+\kappa(a,b+c)]$.
As $p[x]\equiv p[y]\pmod{p^2}$ if and only if $x=y$, this proves the
cocycle identity; commutativity proves symmetry.
Multiplying the first-carry formula by $[u]$ and using
$[u][x]=[ux]$ proves homogeneity.
For $K=\F_4$, $\kappa(a,b)=(ab)^2$, which is not always $ab$:
if $\omega^2+\omega+1=0$, then $\kappa(1,\omega)=\omega^2\ne\omega$.
\end{solution}

\solutionheading{v133:ex-zero-sum-depth}
\begin{solution}
In $W(\F_2)=\Z_2$, $[1]+[1]=2$ has valuation $1$, whereas
$[0]+[0]=0$ has valuation $+\infty$.
The Teichmüller representative $[\omega]\in W(\F_4)$ is still a
primitive cube root of unity, so
$[\omega]^2+[\omega]+1=0$, equivalently
$[1]+[\omega]+[\omega^2]=0$. These are all finite-field zero sums,
but their lifted valuations differ.
\end{solution}

\solutionheading{v133:ex-trace-norm}
\begin{solution}
The degree of $E/K$ is two, so the finite-field trace of either
$a=0$ or $a=1$ is zero. Their lifted traces are $0$ and $2$, respectively,
and their first obstructions are $0$ and $1$.
For any $a\in E$,
\[
 \Nm_{\Z_{q^d}/\Z_q}([a])=
 \prod_{i=0}^{d-1}[a^{q^i}]
 =\left[\prod_{i=0}^{d-1}a^{q^i}\right]
 =[\Nm_{E/K}(a)].
\]
Each step of this product uses multiplicativity of Teichmüller lifts;
there is no corresponding justification for an exact trace identity.
\end{solution}

\solutionheading{v133:ex-formal-deform}
\begin{solution}
Coefficient comparison, or direct substitution, gives
$B=-\sum_{j\ge0}t^{p^j}$ and $B^p-B=t$.
If $B_1,B_2\in tK[[t]]$ were distinct solutions, their difference
$D\in tK[[t]]\setminus\{0\}$ would satisfy $D^p-D=0$, contrary to
the $t$-adic valuations. Thus the solution is unique.
The roots are $\alpha(t)+j$ for $j\in\F_p$, and their product is $c+t$.
As $B(t)\in K[[t]]$, the trace is
$\Tr_{E/K}(\alpha)+pB(t)=\Tr_{E/K}(\alpha)$.
When $p=2$, the root sum of $X^2-X-c$ is $1$; when $p>2$,
the coefficient of $X^{p-1}$ in $X^p-X-c$ is zero, so the trace is $0$.

If $u\in K((t))$ is a root, $v_t(u)<0$ would imply
$v_t(u^p-u)=pv_t(u)<0$, while $v_t(c+t)=0$, a contradiction.
Thus $u\in K[[t]]$, and reduction gives $u_0^p-u_0=c$, contradicting
$\Tr_{K/\F_p}(c)\ne0$.
Hence the polynomial has no root in $K((t))$.
An Artin--Schreier polynomial $X^p-X-a$ over a field of characteristic $p$
is either split or irreducible, giving the result.
\end{solution}

\solutionheading{v134:ex-stabilizer}
\begin{solution}
For $h\in H_C$, multiplication by $\widetilde\zeta^h$ permutes the lifted roots.
Comparing $\widetilde M_C(\widetilde\zeta^hX)$ with
$(\widetilde\zeta^h)^d\widetilde M_C(X)$ shows
$a_j\ne0\Rightarrow\kappa_C\mid j$.
A further zero coefficient modulo $p$ merely implies $a_j\in p\Z_q$.
In Example~\ref{v125:ex-160-hensel},
$\widetilde F_{5,2}=X^8+3X^6+3X^4+X^2+8$.
All odd-degree coefficients vanish exactly by composition through $X^2$,
whereas the coefficients of $X^6,X^4$ are $3\bmod9$, each with valuation one.
\end{solution}

\solutionheading{v134:ex-subset-periods}
\begin{solution}
For $C=\{1,2,4\}$, the subset sums for $j=1$ are $1,2,4$;
for $j=2$ they are $3,5,6$; and for $j=3$ the only sum is $0$.
Each occurs once. Thus $(a_1,a_2,a_3)=(-\eta,\eta',-1)$,
with profile $(1,0,0)$. For $C'=\{3,5,6\}$, the corresponding
sums are $3,5,6$ for $j=1$, $1,2,4$ for $j=2$, and $0$ for $j=3$.
The coefficients are $(-\eta',\eta,-1)$ and the profile is $(0,1,0)$.
\end{solution}

\solutionheading{v134:ex-newton-bound}
\begin{solution}
Comparing the first two coefficients of
$TE'(T)/E(T)=-\sum_{k\ge1}P_kT^k$ gives
$a_1=-P_1$ and $2a_2=P_1^2-P_2$.
The two numerator terms have valuations $2u,w$, so the ultrametric
inequality gives the bound. When their valuations differ, the lower-valuation
term cannot cancel, and equality follows. When $2u=w$, the leading
unit parts may cancel and must be compared. At $p=2$, computing
$a_2\bmod2^s$ by division normally requires the numerator
$P_1^2-P_2\bmod2^{s+1}$.
\end{solution}

\solutionheading{v134:ex-base-change}
\begin{solution}
Since $8\equiv1\pmod7$, the three subcosets are $\{1\},\{2\},\{4\}$,
and $\widetilde M_C=(X-\widetilde\zeta)(X-\widetilde\zeta^2)
(X-\widetilde\zeta^4)$. Each linear factor has unit constant coefficient,
so its profile is $(0)$. Nevertheless the sum of the three units appearing
as the quadratic coefficient of the product has valuation one.
The numerical profiles do not retain their relative phases.
In general there are $\gcd(s,d)$ subcosets of length $d/\gcd(s,d)$;
$q$-power Frobenius permutes them transitively, so their corresponding
coefficient valuations agree.
\end{solution}

\solutionheading{v134:ex-first-coefficient-carry}
\begin{solution}
The sum $\sum_{|J|=j}[b_J]$ is invariant under relative $q$-power
Frobenius and thus belongs to $\Z_q$. Division by $p$ still gives an element
of $\Z_q$, whose residue therefore lies in $\F_q$.
Example~\ref{v132:ex-seven-period} gives $\eta\equiv2\pmod4$.
Both relevant Teichm\"uller sums equal $\eta$, hence
$\mathcal B_{1,1}(C_1)=\mathcal B_{1,2}(C_3)=1\in\F_2$.
\end{solution}

\solutionheading{v47:ex-galois-extensions}
\begin{solution}
Every element has a unique expression \(a+b\omega\) with \(a,b\in\Z/4\Z\), so \(|A|=16\).
The quotient ring modulo \(2\) is \(\F_4\); each residue class has \(4\) lifts, and there are \(3\) nonzero residue classes,
so the number of units is \(12\). By \(\omega^3=1\) and \(\overline\omega\ne1\),
the Teichmüller set is \(\{0,1,\omega,\omega^2\}\).

The absolute Frobenius fixes the prime subring and sends \(\omega\) to \(\omega^2=-1-\omega\), so
\[
 \phi(a+b\omega)=a+b\omega^2=(a-b)-b\omega.
\]
If this equals \(a+b\omega\), comparison of basis coordinates gives \(b=0\); hence the fixed subring is exactly \(\Z/4\Z\).
For \(z=1+\omega=-\omega^2\), we have
\[
 \phi(z)=1+\omega^2=-\omega,
 \qquad z^2=1+2\omega+\omega^2=\omega.
\]
Their difference is \(2\omega\ne0\). Also \(z^3=-1\ne1\) and \(z^6=1\),
and \(z^2=\omega\ne1\), so its order is exactly \(6\).
\end{solution}

\solutionheading{v47:ex-hensel-seven}
\begin{solution}
Write \(g=X^3+X+1\in\F_2[X]\). The initial local inverse is
\[
 7^{-1}Xg'=X^3+X=1\pmod g.
\]
Monic division modulo \(4\) gives \(X^7-1\bmod G_1=2+2X^2\),
so \(E_1=1+X^2\) and \(\Delta_1=1+X^2\), and hence
\[
 G_2=X^3+2X^2+X+3\pmod4.
\]
With this integer-coefficient representative, the remainder modulo \(8\) is \(4+4X+4X^2\),
so \(E_2=1+X+X^2\). At each level the local inverse modulo \(2\) remains \(1\), and therefore
\[
 G_3=G_2+4(1+X+X^2)=X^3+6X^2+5X+7\pmod8.
\]
Division by \(G_3\), computed modulo \(16\), gives remainder \(8\), so \(E_3=\Delta_3=1\), whence
\[
 G_4=X^3+6X^2+5X+15\pmod{16}.
\]
The degree and reduction remain unchanged at each level, so the resulting factors are basic irreducible polynomials.

Since \(15^{-1}=15\pmod{16}\), its monic reciprocal polynomial is
\[
 G_4^\vee(X)=15X^3G_4(X^{-1})
 =X^3+11X^2+10X+15.
\]
Its reduction is the other irreducible cubic factor, so uniqueness of lifting gives the complete factorization
\[
 \begin{aligned}
 X^7-1={}&(X-1)(X^3+6X^2+5X+15)\\
 &\cdot(X^3+11X^2+10X+15)
 \quad\text{in }(\Z/16\Z)[X].
 \end{aligned}
\]
Direct multiplication gives the product \(1+X+\cdots+X^6\) of the two cubic factors modulo \(16\);
multiplication by \(X-1\) then gives \(X^7-1\).
\end{solution}

\solutionheading{v85:ex-coset-orbits}
\begin{solution}
We have \(|H|=5\), and the quotient is indexed by \(\Z/3\Z\). Squaring has orbits \(\{0\}\) and \(\{1,2\}\).
Put \(\omega=\gamma^5\). Its order is three, so \(\omega\in\F_4\setminus\F_2\) and \(\omega^2+\omega+1=0\).
The three root polynomials are \(T^5+1,T^5+\omega,T^5+\omega^2\), with smallest coefficient fields \(\F_2,\F_4,\F_4\), respectively.
The last two constants have sum and product one, giving the displayed degree-ten polynomial.
Each coset has five elements, whereas its quotient orbit has length one or two. For instance, \(\gamma\) has Frobenius period four, but its coset has period two.
Return of the entire root set does not require pointwise return of its elements.
\end{solution}

\section*{Solutions to Chapter~\ref{chap:trace}}
\pdfbookmark[1]{Solutions to Chapter \getrefnumber{chap:trace}}{solutions.chap:trace}

\solutionheading{ex:trace-f9}
\begin{solution}
Since \(\Tr_{\F_9/\F_3}(1)=2\ne0\), the Artin--Schreier criterion makes \(X^3-X-1\) irreducible.
Also \(\alpha^3=-\alpha\), so \(\Tr(\alpha)=0\), whereas \(\alpha^3-\alpha=-2\alpha=\alpha\).
Thus all roots of \(X^3-X-\alpha\) are \(\alpha,\alpha+1,\alpha+2\), and all lie in \(\F_9\).
\end{solution}

\solutionheading{ex:trace-iterate}
\begin{solution}
Put \(L=\F_{q^g}\) and \(\tau=\sigma^s\). Its order on \(E\) is \(n/g\), and its fixed field is \(L\).
The maps \(\tau\) and \(\sigma^g\) generate the same cyclic group, so \(1+\tau+\cdots+\tau^{n/g-1}\) is exactly the relative trace operator.
Hence \(\Tr_{E/L}((\tau-1)x)=0\).
The kernel of \(\tau-1\) is \(L\), and its image has \(q^{n-g}\) elements; the relative trace is surjective, so the trace-zero space has the same number of elements.
They are therefore equal. Each nonempty fiber is an additive coset of the kernel \(L\), giving \(q^g\) solutions.
\end{solution}

\solutionheading{ex:trace-norm-fibers}
\begin{solution}
Every \(x\) satisfying the conditions is a root of \(P(X)=X^2-tX+u\). Put \(\Delta=t^2-4u\).
If \(\Delta\) is a nonsquare, \(P\) is irreducible and Frobenius interchanges its two roots, so both satisfy the required trace and norm conditions.
If \(\Delta=0\), the unique root \(x=t/2\in k\) satisfies \(\Tr(x)=2x=t\) and \(\Nm(x)=x^2=u\), giving exactly one solution.
If \(\Delta\) is a nonzero square, both roots lie in \(k\).
The trace of a base-field element in a quadratic extension is twice that element; if one satisfies \(\Tr(x)=t\), then \(x=t/2\),
and substitution into \(P\) forces \(\Delta=0\), a contradiction. Thus there are no solutions in this case. In all three cases the number of solutions is \(1-\chi(\Delta)\).
\end{solution}

\solutionheading{ex:q-artin-schreier}
\begin{solution}
In a finite splitting extension, every root \(\alpha\) satisfies \(\alpha^q=\alpha+a\).
Since \(a^q=a\), induction gives \(\alpha^{q^j}=\alpha+ja\).
The condition \(a\ne0\) ensures that the first return to \(\alpha\) occurs exactly at time \(p\), so every root has degree \(p\) over \(\F_q\).
The derivative is minus one, so there are \(q\) distinct roots; grouping them into orbits of length \(p\) gives exactly \(q/p\) irreducible factors.
\end{solution}

\solutionheading{ex:supp-hilbert-formula}
\begin{solution}
If \(R_1\) is also a right inverse, then \(D(R_1-R)=0\), so \((R_1-R)(a)\in\ker D=k\) for every \(a\).
Thus \(R_1-R\) is a linear functional taking values in \(k\). Conversely, for such a functional \(\ell\), we have \(D(R+\ell)=DR=\id\).
Since \(\dim_k\ker\Tr=n-1\), there are \(q^{n-1}\) functionals and hence the same number of right inverses.
Choosing a trace-one element gives a particular right inverse; the solutions obtained from any two such choices differ, at every parameter value, by a base-field element.
This also includes \(n=1\): the trace kernel is then zero, and there is only one right inverse.
\end{solution}

\solutionheading{ex:supp-trace-factor}
\begin{solution}
Let \(\omega^2+\omega+1=0\). The trace-one elements of \(\F_4\) are \(\omega,\omega^2\), so
\[
X^4+X+1=(X^2+X+\omega)(X^2+X+\omega^2).
\]
Each constant term has absolute trace one, so both quadratic factors are irreducible over \(\F_4\).
\end{solution}

\solutionheading{ex:wu-subfield-pairing}
\begin{solution}
Since \(xy\in L\), the relative trace is \(\Tr_{E/L}(xy)=sxy\).
Trace transitivity gives the formula.
If \(p\nmid s\), this is a nonzero multiple of the nondegenerate trace pairing of \(L/k\) and is therefore nondegenerate;
if \(p\mid s\), the restriction is identically zero.
For each nonzero \(x\in L\), nondegeneracy on \(E\) guarantees some \(y\in E\) with nonzero pairing;
such a \(y\) need not lie in \(L\).
\end{solution}

\solutionheading{ex:wu2-equivariant-functionals}
\begin{solution}
The commutation condition is equivalent to \(\ell\) vanishing on \(\im(\sigma-\id)=\ker\Tr\).
Choose \(u\in E\) with \(\Tr(u)=1\), possible by surjectivity of the trace. For every \(x\),
\(x-\Tr(x)u\in\ker\Tr\), so \(\ell(x)=\ell(u)\Tr(x)\). Taking \(c=\ell(u)\) gives existence;
evaluation at \(u\) gives uniqueness. Such a functional restricts to the identity on \(k\) if and only if \(cn=1\).
Thus an equivariant projection exists if and only if \(p\nmid n\); it is then unique and equals \(\Tr/n\).
\end{solution}

\solutionheading{ex:wu10-diagonalization}
\begin{solution}
Choose \(\alpha\in\F_8\) satisfying \(\alpha^3+\alpha+1=0\). This cubic has no root in \(\F_2\) and is therefore irreducible,
so \(\mathcal B=(1,\alpha,\alpha^2)\) is an \(\F_2\)-basis.
From \(\alpha^4=\alpha^2+\alpha\), the Frobenius map given by squaring has matrix
\[
[\sigma]_{\mathcal B}=\begin{pmatrix}
1&0&0\\
0&0&1\\
0&1&1
\end{pmatrix}.
\]
The lower-right two-dimensional block \(C=\left(\begin{smallmatrix}0&1\\1&1\end{smallmatrix}\right)\) satisfies
\(C^2+C+I=0\) and is not scalar, so its minimal polynomial is \(T^2+T+1\).
The minimal polynomial on the first coordinate is \(T+1\); the two are coprime, so the full matrix has minimal polynomial
\(T^3-1=(T+1)(T^2+T+1)\).
This polynomial has no repeated roots, so the matrix is diagonalizable if and only if \(T^2+T+1\) splits over the chosen scalar field.
Its roots are primitive cube roots of unity, which lie in \(\F_{2^r}\) if and only if \(3\mid 2^r-1\), namely if \(r\) is even.
Thus the matrix is not diagonalizable over \(\F_2\) or \(\F_8\), but is diagonalizable over \(\F_{64}\); the least scalar extension is \(\F_4\).
Here scalar extension applies to the matrix representation; the original map \(x\mapsto x^2\) is not \(\F_8\)-linear.
Separability of \(E/\F_2\) says that conjugates of field elements are distinct; it does not ensure that this linear operator's eigenvalues already lie in the specified scalar field.
\end{solution}

\solutionheading{v45:ex-quadratic-first-zero-f5}
\begin{solution}
The two discriminants are \(-8=2\) and \(1-8=3\), respectively, both nonsquares in \(\F_5\), so both polynomials are irreducible.
The recurrence \(u_n=-bu_{n-1}-2u_{n-2}\) gives
\[
 \begin{array}{c|rrrrrr}
 &u_1&u_2&u_3&u_4&u_5&u_6\\\hline
 f_0&1&0&3&0&4&0\\
 f_1&1&4&4&3&4&0
 \end{array}
\]
For \(f_0\), the first zero term of positive index occurs at \(r=2\), and \(\lambda=-2u_1=3\).
Thus the minimal binomial multiple is \(X^2-3=X^2+2\) and the ordinary order is \(2\ord(3)=8\), so the polynomial is not primitive.
For \(f_1\), the first zero term of positive index occurs at \(r=6\), and \(\lambda=-2u_5=2\).
The minimal binomial multiple is \(X^6-2\) and the ordinary order is \(6\ord(2)=24=5^2-1\), so \(f_1\) is primitive.

By \eqref{v45:eq-first-zero-quotient},
\[
 \frac{X^7-2X}{X^2+X+2}
 =X^5+4X^4+4X^3+3X^2+4X.
\]
The coefficients, ordered from leading to constant, are \(1,4,4,3,4,0\); the first zero coefficient is indeed in position six.
The condition \(b\ne0\) alone is insufficient. For example, when \(b=2\),
\[
 X^2+2X+2=(X-1)(X-2)
 \quad\text{in }\F_5[X],
\]
so the quadratic is reducible. If irreducibility is already known, excluding \(b=0\) within this fixed constant-term fiber does give the two primitive quadratics.
\end{solution}

\solutionheading{v45:ex-optimal-quadratics-f4}
\begin{solution}
The homogeneous cyclotomic polynomial is
\[
 \Phi_5^{\mathrm{hom}}(\alpha,\beta)
 =\alpha^4+\alpha^3\beta+\alpha^2\beta^2+\alpha\beta^3+\beta^4.
\]
In characteristic two, put \(B=\alpha+\beta\) and \(C=\alpha\beta\). Then
\(B^4=\alpha^4+\beta^4\), \(B^2C=\alpha^3\beta+\alpha\beta^3\), \(C^2=\alpha^2\beta^2\).
Thus \(\Lambda_5(B,C)=B^4+B^2C+C^2\).

In \(A\), the relation \(C^2+C+1=0\) gives \(C^3=1\) and \(C^4=C\), so
\[
 (B^2+C^2B+C)^2=B^4+C^4B^2+C^2
 =B^4+CB^2+C^2.
\]
Squaring is injective on \(A\simeq k\times k\), hence also on \(A[B]\), so this square root in \(A[B]\) is unique.
At \(C=\omega\) and \(C=\omega^2\), respectively, we have
\[
 \begin{aligned}
 P(B,\omega)&=B^2+\omega^2B+\omega=(B-1)(B-\omega),\\
 P(B,\omega^2)&=B^2+\omega B+\omega^2=(B-1)(B-\omega^2).
 \end{aligned}
\]
By Theorem~\ref{v45:thm-optimal-quadratic-criterion}, all primitive quadratics are
\[
 X^2+X+\omega,\quad X^2+\omega X+\omega,
 \qquad X^2+X+\omega^2,\quad X^2+\omega^2X+\omega^2.
\]
There are exactly two for each fixed primitive constant term, consistently with \(\varphi(5)/2=2\).
Since \(\Lambda_5(B,c)=P(B,c)^2\), the two polynomials have the same zeros and both detect primitivity.
However, each of the two roots of \(\Lambda_5(B,c)\) has multiplicity two, and its degree is four;
the optimal polynomial \(P(B,c)\) is monic, square-free, and of degree two.
\end{solution}

\solutionheading{v45:ex-quadratic-cosets-nine}
\begin{solution}
The square of \(1+\beta\) is \(-\beta\) and its fourth power is \(-1\), so it has order eight and generates \(K^\times\).
Thus the nonzero squares in \(K\) are \(1,2,\beta,-\beta\).
The discriminant of \(F_t\) is \(\delta_t=1+\beta+t\):
both \(\delta_0=1+\beta\) and \(\delta_1=2+\beta\) are nonsquares, whereas \(\delta_2=\beta\) is a square.
Hence \(F_0,F_1\) are irreducible, whereas
\[
 F_2(X)=(X-\beta)(X+1+\beta).
\]
There is one reducible member and two irreducible members, consistently with \((q-1)/2=1\) and \((q+1)/2=2\).

If \(\theta\) is a root of \(F_0\), then \(\Nm_{\F_{81}/K}(\theta)=-\beta\), and this norm has order four.
A primitive element of \(K\) must have order eight, so \(F_0\) is not primitive.

The constant term of \(F_1\) is \(c=-(1+\beta)=(1+\beta)^5\), of order eight in \(K^\times\).
Use \(K\) as the base field, of order \(Q=9\). For the Lucas recurrence of \(F_1=X^2+X+c\),
\[
 U_2(1,c)=-1\ne0,\qquad
 U_5(1,c)=1+c^2=1-\beta\ne0.
\]
If \(\theta\) is a root of \(F_1\), then \(\theta^{10}=c\in K\),
so its projective order \(r=\ord(\theta K^\times)\) divides ten.
This projective order is the least positive index at which the Lucas sequence vanishes; the two terms above exclude \(r\mid2\) and \(r\mid5\), respectively, so \(r=10\).
Also \(\theta^{10}=c\) has order eight, so \(\ord(\theta)=10\cdot8=80\).
Thus \(F_1\) is a primitive quadratic over \(K\), and both its roots generate \(\F_{81}^\times\).
\end{solution}

\solutionheading{v74:ex-norm-preservers}
\begin{solution}
If \(T(x)=0\), norm preservation gives \(\Nm(x)=0\), so \(x=0\).
Thus \(T\) is injective and, by finiteness, bijective.
Take \(H=\ker(\Nm:E^\times\to\F_q^\times)\).
The norm is onto and \(q>2\), so \(H\) is a proper subgroup of index \(q-1>1\).
Additivity transfers norm preservation to each difference: for \(x\ne y\),
\[
\Nm\!\left(\frac{T(x)-T(y)}{x-y}\right)
=\frac{\Nm(T(x-y))}{\Nm(x-y)}=1.
\]
Hence the difference quotient belongs to \(H\).
Write \(q=p^r\) and apply Carlitz--McConnel to \(E=\F_{p^{rm}}\).
It gives \(T(x)=ax^{p^s}\), with \(a\in H\) and \(0\le s<rm\); the constant term vanishes since \(T(0)=0\).
For \(c\in\F_q\), the equality \(T(c)=cT(1)\) gives \(ac^{p^s}=ac\).
The \(p^s\)-Frobenius is therefore the identity on \(\F_q\), which means \(r\mid s\).
Writing \(s=rj\) yields the required form \(T(x)=ax^{q^j}\).
Conversely, \(q^j\)-Frobenius is \(\F_q\)-linear and preserves the norm, so
\[
\Nm(ax^{q^j})=\Nm(a)\Nm(x)^{q^j}=\Nm(x).
\]
Linearity supplies only \(T(x)-T(y)=T(x-y)\); preservation of norm fibers is the additional condition that restricts the multiplicative structure.

Over \(\F_2\), the norm \(\Nm_{\F_{2^m}/\F_2}(x)\) equals one at every nonzero element and zero at zero, so every linear bijection preserves it.
The given images of the three basis vectors are linearly independent; in fact \(T^2=\id\), so \(T\) is bijective.
If \(T(x)=ax^{2^j}\), then \(T(1)=1\) forces \(a=1\).
Since \(T(\theta)=\theta\) and \(\theta\) has degree three over \(\F_2\), only \(j=0\) is possible in \(0\le j<3\).
This would make \(T\) the identity, contrary to \(T(\theta^2)=\theta^2+1\).
Thus \(q>2\) is substantive: the quotient must retain more than one nonzero norm label.
\end{solution}

\solutionheading{v82:ex-adjacent-norm-fibers}
\begin{solution}
Expanding with the two conjugates gives
\[
 \Nm(x+1)=(x+1)(x^q+1)=\Nm(x)+\Tr(x)+1.
\]
Thus the conditions are equivalent to \(\Tr(x)=t\) and \(\Nm(x)=a\). For odd \(q\), apply Exercise~\ref{ex:trace-norm-fibers}.

Let \(q=2^s\). If \(t=0\), then \(x^q=x\), so \(x\in k\), and \(\Nm(x)=x^2=a\) has exactly one solution in \(k\).
If \(t\ne0\), substituting \(z=x/t\) gives \(z^2+z=a/t^2\).
When \(\Tr_{k/\F_2}(a/t^2)=1\), this polynomial is irreducible over \(k\); its two roots lie in \(E\) and are exchanged by the \(q\)-Frobenius. Both corresponding values of \(x\) have the required trace and norm.
When that absolute trace is zero, both roots lie in \(k\), so their relative traces for \(E/k\) are zero rather than \(t\), and neither is admissible. Hence
\[
 M(a,b)=
 \begin{cases}
 1,&t=0,\\
 2,&t\ne0,\ \Tr_{k/\F_2}(a/t^2)=1,\\
 0,&t\ne0,\ \Tr_{k/\F_2}(a/t^2)=0.
 \end{cases}
\]
More generally, if \(X^2-tX+a\) splits over \(k\) into two distinct roots, Frobenius fixes each root instead of exchanging them.
In odd characteristic, a base-field root \(x\) with relative trace \(2x=t\) would force the roots to coincide; in characteristic two with \(t\ne0\), every base-field root has relative trace zero.
Thus complete splitting in \(E\) alone is not the prescribed trace-and-norm condition.
\end{solution}

\solutionheading{v82:ex-primitive-constant-term}
\begin{solution}
Take a root \(\alpha\in\F_{Q^m}\) of \(f\). Primitivity gives \(\ord(\alpha)=Q^m-1\), so
\[
 \Nm_{\F_{Q^m}/\F_Q}(\alpha)
 =\alpha^{(Q^m-1)/(Q-1)}=(-1)^m f(0)
\]
has order \(Q-1\). This proves the first assertion.

Write \(Q-1=2s\), where \(s\) is odd. If \(\lambda\in\F_Q^\times\) is primitive, then \(-1=\lambda^s\), hence \(-\lambda=\lambda^{s+1}\).
Since \(\gcd(2s,s+1)=2\), the order of \(-\lambda\) is \(s\), so it is not primitive.
If \(m\) is odd and \(f(0)=\lambda\) is primitive, the first assertion would require \(-\lambda\) to be primitive, a contradiction.

The norm condition is necessary for primitivity, but does not force a root to have maximal multiplicative order.
For example, \(X^2+2\) in Exercise~\ref{v45:ex-quadratic-first-zero-f5} is irreducible and has primitive constant term in \(\F_5\), but its roots have order eight rather than \(5^2-1=24\).
Irreducibility ensures field generation, whereas the norm records only a quotient of the multiplicative order; together these still do not replace primitivity.
\end{solution}

\solutionheading{v146:ex-lucas-characteristic}
\begin{solution}
For independent \(\alpha,\beta\), use \(\alpha^n-\beta^n=(\alpha^m-\beta^m)^{p^v}\) and injectivity of symmetric substitution to obtain
\[
 U_n=\begin{cases}
 (B^2-4C)^{(p^v-1)/2}U_m^{p^v},&p\ne2,\\
 B^{2^v-1}U_m^{2^v},&p=2.
 \end{cases}
\]
For \(v=0\) this is an identity. For two distinct nonzero roots, the discriminant factor is nonzero, so vanishing still means that the root-ratio order divides \(m\), while the multiplicities of existing factors are multiplied by \(p^v\). For \(v>0\), the additional discriminant zeros represent coincident roots, not new nonzero root-ratio orders divisible by \(p\): roots of unity in characteristic \(p\) have order prime to \(p\).
\end{solution}

\section*{Solutions to Chapter~\ref{chap:dual}}
\pdfbookmark[1]{Solutions to Chapter \getrefnumber{chap:dual}}{solutions.chap:dual}

\solutionheading{ex:dual-f4}
\begin{solution}
We have \(\Tr(1)=0\), \(\Tr(\alpha)=\Tr(\alpha^2)=1\), and \(\alpha^3=1\).
Direct verification then gives \((\alpha^2,1)\) as the dual basis of \((1,\alpha)\).
For \((\alpha,\alpha^2)\), the squares of both basis vectors have trace one and their product has trace \(\Tr(1)=0\), so this basis is self-dual.
\end{solution}

\solutionheading{ex:trace-orthogonal}
\begin{solution}
We have \(x\in k^\perp\) if and only if \(\Tr(cx)=c\Tr(x)=0\) for every \(c\in k\), equivalently \(\Tr(x)=0\).
Also, \(x\) is orthogonal to every vector in \(U+V\) if and only if it is orthogonal to all vectors in each of \(U\) and \(V\), giving the second identity.
\end{solution}

\solutionheading{ex:no-selfdual-power}
\begin{solution}
If \(n\geq3\), self-duality requires the distinct basis vectors \(1,\alpha^2\) to have trace pairing zero, and the basis vector \(\alpha\) to have trace pairing one with itself.
But both pairings equal \(\Tr(\alpha^2)\), a contradiction.
If \(n=2\), self-duality requires \(\Tr(1)=1\), whereas \(\Tr(1)=2\cdot1_k\).
The identity \(2\cdot1_k=1_k\) would imply \(1_k=0\), again a contradiction.
\end{solution}

\solutionheading{ex:supp-discriminant}
\begin{solution}
We have \(\Tr(1)=2\), \(\Tr(\alpha)=0\), and \(\Tr(\alpha^2)=1\), so the Gram matrix is \(\operatorname{diag}(2,1)\) and its discriminant is two.
We have \(f'(\alpha)=2\alpha\), with norm \((2\alpha)(-2\alpha)=1\); multiplying by \((-1)^{2(2-1)/2}=-1=2\) again gives two.
\end{solution}

\solutionheading{ex:wu-even-selfdual-sum}
\begin{solution}
Self-duality gives \(\Tr(b_i^2)=1\). In characteristic two, \(\Tr(b_i^2)=\Tr(b_i)^2\), and squaring is injective, so \(\Tr(b_i)=1\).
Expanding \(1\) by the self-dual coordinate recovery formula gives \(1=\sum_i\Tr(b_i)b_i=\sum_i b_i\).
If \(n\) is even, taking traces gives \(\Tr(1)=n\cdot1=0\), so the identities are fully consistent.
\end{solution}

\solutionheading{ex:wu2-sum-subfields}
\begin{solution}
By the intersection formula, \(L\cap M=\F_{q^g}\), where \(g=\gcd(a,b)\). Hence \(\dim_k(L+M)=a+b-g\). Proposition~\ref{prop:relative-trace-adjoint} gives
\[
(L+M)^\perp=L^\perp\cap M^\perp
=\ker\Tr_{E/L}\cap\ker\Tr_{E/M}.
\]
If \(L+M\) is a field, it contains the compositum of the two fields and is contained in that compositum, so it must equal \(\F_{q^{\lcm(a,b)}}\).
Writing \(a=gu,b=gv\), the dimension identity becomes \(u+v-1=uv\), namely \((u-1)(v-1)=0\).
Thus a necessary condition is \(a\mid b\) or \(b\mid a\); in either case one field contains the other and the sum is the larger field, proving sufficiency.
\end{solution}

\solutionheading{ex:wu7-trace-radical}
\begin{solution}
\begin{enumerate}
\item By Theorem~\ref{thm:algebra-trace-coefficient},
\[
B_{\Tr}(u,v)=B_\lambda(f'(\alpha)u,v).
\]
The coefficient pairing is nondegenerate, so \(u\in\rad B_{\Tr}\) if and only if \(f'(\alpha)u=0\).
Write \(f=hf_1\) and \(f'=hg_1\); then \(\gcd(f_1,g_1)=1\). For a polynomial representative \(U\),
\[
f\mid f'U\quad\Longleftrightarrow\quad f_1\mid g_1U
\quad\Longleftrightarrow\quad f_1\mid U.
\]
Thus the radical is the ideal \((f_1(\alpha))\). The map
\[
k[T]/(h)\longrightarrow(f_1(\alpha)),\qquad
\overline V\longmapsto f_1(\alpha)V(\alpha)
\]
is well defined and surjective; since \(f\mid f_1V\) if and only if \(h\mid V\), it is also injective.
Thus the radical has dimension \(\deg h\) and the rank is \(d-\deg h\). This includes \(f'=0\): then \(h=f\) and \(f_1=1\), and the radical is all of \(A\).
The coefficient pairing is nondegenerate; the radical of the trace pairing consists precisely of the elements annihilated by multiplication by the derivative.

\item Here \(A=k[\varepsilon]/(\varepsilon^e)\). Translation leaves the highest-degree coefficient unchanged, so
\(\lambda(\sum_{i=0}^{e-1}a_i\varepsilon^i)=a_{e-1}\).
Let \(J_e=(\delta_{i+j,e-1})_{0\le i,j<e}\) be the matrix with ones on the antidiagonal and zeros elsewhere. Then
\[
G_\lambda=J_e,
\qquad f'(\alpha)=e\varepsilon^{e-1}.
\]
Indeed, when \(i+j<e\), the value \(\lambda(\varepsilon^{i+j})\) is one only for \(i+j=e-1\); when \(i+j\ge e\), the product is zero.
The derivative formula further gives
\[
\Tr_{A/k}\left(\sum_{i=0}^{e-1}a_i\varepsilon^i\right)=ea_0,
\qquad G_{\Tr}=\operatorname{diag}(e,0,\ldots,0).
\]
Thus \(B_\lambda\) has rank \(e\) and zero radical.
If \(p\nmid e\), then \(B_{\Tr}\) has rank one and radical \((\varepsilon)\) of dimension \(e-1\);
if \(p\mid e\), then \(B_{\Tr}=0\), and the radical is all of \(A\).
In both cases the nilradical is \((\varepsilon)\): an element with nonzero constant term is a unit, whereas an element with zero constant term has zero \(e\)th power.
Thus when \(e=p\) the constant subspace also lies in the radical of the trace pairing; in odd characteristic with \(e=2\), its radical equals the one-dimensional nilradical.

\item In characteristic two, \(P'=1\), so
\[
f'=(T+1)^2P^2=T^6+1,
\qquad h=(T+1)^2P^2,
\qquad f/h=P.
\]
By the first part, \(\rad B_{\Tr}=(P(\alpha))\), with basis
\[
P(\alpha),\ \alpha P(\alpha),\ \ldots,\ \alpha^5P(\alpha).
\]
The representatives have degrees two through seven and are therefore linearly independent; the radical has dimension six, so they form a basis.
Hence the trace pairing has rank two and the coefficient pairing rank eight.
Under the Chinese remainder decomposition, multiplicity two makes the trace on the first component zero, while multiplicity three makes the trace on the second component equal to the trace of its residue field \(\F_4\);
the latter has rank two, in agreement with the derivative calculation.

To find the two components of \(e_1\), put \(\varepsilon=T+1\). Modulo \((T+1)^2\) we have
\(T=1+\varepsilon\) and \(P=1+\varepsilon\), hence \(TP^3=(1+\varepsilon)^4=1\);
modulo \(P^3\) we have \(TP^3=0\). Thus the two coordinates of \(e_1\) are \((1,0)\), making it a nonzero idempotent and, in particular, not nilpotent.
Also \(e_1\in(P(\alpha))=\rad B_{\Tr}\), so its trace pairing with every \(v\) is zero.
On the other hand, \(TP^3\) is monic of degree seven and requires no reduction modulo \(f\), so \(\lambda(e_1)=1\).
Thus this nonnilpotent element belongs to the radical of the trace pairing but has nonzero value under the coefficient functional.

\item The extension degree proved in the appendix is \(p\), so the listed powers form a basis.
When \(0\le i,j<p\), if \(i+j<p\) the product is \(u^{i+j}\); if \(i+j\ge p\), then
\[
u^{i+j}=u^p u^{i+j-p},\qquad u^p\in k,\quad i+j-p\le p-2.
\]
Thus the highest-coefficient functional still satisfies \(\lambda(u^iu^j)=\delta_{i+j,p-1}\), namely \(G_\lambda=J_p\).
But \((T^p-u^p)'=0\), so the algebra trace is identically zero and \(G_{\Tr}=0\).
The two pairings have ranks \(p\) and zero, and radicals zero and the whole of \(E\), respectively.
Since \(E\) is a field, its nilradical is zero.

These Gram matrices are identical to those in the second part when \(e=p\), but that algebra has a nonzero nilpotent ideal of dimension \((p-1)\), whereas this is an inseparable field extension.
Thus degeneracy of the trace pairing detects inseparability of the defining polynomial and, over arbitrary base fields, cannot be decided solely by reducedness.
The perfection of finite fields makes every irreducible factor separable; this is what makes reducedness equivalent to nondegeneracy of the trace pairing.
\end{enumerate}
\end{solution}

\solutionheading{ex:wu10-local-trace}
\begin{solution}
The Chinese remainder components are \(A'_1=\F_2[T]/((T+1)^3)\) and \(A_2=\F_2[T]/(P^3)\), where \(P=T^2+T+1\).
The nilradical is \((\overline{T+1})A'_1\times(\overline P)A_2\), of dimension \(2+4=6\),
and the reduced quotient is still \(\F_2\times\F_4\). Proposition~\ref{prop:algebra-periodic-reduction} therefore gives a three-dimensional periodic part.
Both components have multiplicity three, whose image in the base field is one, so the local trace formula becomes
\[
\Tr_{A'/\F_2}(u)=\bar u_1+\Tr_{\F_4/\F_2}(\bar u_2).
\]
Thus the trace pairing has rank \(1+2=3\) and radical exactly the six-dimensional nilradical above.
The highest-coefficient pairing is always nondegenerate and hence has rank \(\dim A'=9\).
Corollary~\ref{cor:algebra-fixed-points-reduced} gives \(2^{1+\gcd(2,r)}\) solutions: four when \(r\) is odd and eight when it is even.
The reduced quotient and periodic-point data are unchanged; the total algebra dimension and nilradical dimension each increase by one.
The image of the first component's multiplicity in the base field changes from zero to one, so its residue field contributes again to the algebra trace and the trace-pairing rank also increases by one.
\end{solution}

\solutionheading{v85:ex-quadratic-coordinates}
\begin{solution}
Since \(3\) is not a square in \(\F_5\), the polynomial \(Q_1\) is irreducible. The identity \(\theta^2=3\) gives \(\theta^5=-\theta\).
As \(p_S=T^2-T\), the first coordinate is
\[
 \frac{3+\theta}{3-\theta}=2+\theta,
 \qquad (2+\theta)(2-\theta)=1.
\]
The second is \(p_S(2)/p_S(3)=2/1=2\), and the third is
\(p_S'(2)/p_S(2)=3/2=4\).
The fifth-power Frobenius sends these coordinates to \(2-\theta,2,0\), respectively. The last map sends the entire additive group to zero, so it is not an automorphism: it removes the nilpotent direction in the ring of dual numbers.
\end{solution}

\solutionheading{v97:ex-node-change}
\begin{solution}
The Lagrange basis for the first node set is
\[
 \ell_0=\frac{(X-1)(X-2)}2,\qquad
 \ell_1=-X(X-2),\qquad
 \ell_2=\frac{X(X-1)}2.
\]
Evaluation at the second set gives
\[
 M=\begin{pmatrix}1&3&6\\4&6&6\\3&6&3\end{pmatrix},
 \qquad\det M=1\quad\text{in }\F_7.
\]
The prescribed polynomial has evaluation columns \(v_\alpha=(1,4,2)^{\mathsf T}\) and \(v_\beta=(2,4,1)^{\mathsf T}\), and direct multiplication verifies the transpose formula.
The two node sets are disjoint and each consists of distinct elements, so Proposition~\ref{v97:prop-cauchy-node-change} makes every square minor nonzero.
However, the same-representative map sends \(\overline X\) to \(\overline X\).
If it preserved multiplication, it would preserve the relation \(f(\overline X)=0\).
In \(k[X]/(g)\) this relation fails, since the root \(3\) of \(g\) satisfies \(f(3)=6\ne0\).
Thus a change of linear coordinates is not the same as a ring isomorphism between the quotient algebras.
\end{solution}

\section*{Solutions to Chapter~\ref{chap:linearized}}
\pdfbookmark[1]{Solutions to Chapter \getrefnumber{chap:linearized}}{solutions.chap:linearized}

\solutionheading{ex:subspace-line}
\begin{solution}
Dividing the variable by \(u\) gives
\[
\prod_{a\in k}(X-au)=u^q\bigl((X/u)^q-X/u\bigr)=X^q-u^{q-1}X.
\]
Put \(L(X)=X^q-u^{q-1}X\) and \(b=L(v)\). Linear independence ensures \(b\ne0\), and the recurrence in this chapter gives
\[
L_{ku+kv}(X)=L(X)^q-b^{q-1}L(X).
\]
Its degree is \(q^2\) and its root set is exactly the two-dimensional space \(ku+kv\).
\end{solution}

\solutionheading{ex:kernel-gcd}
\begin{solution}
The corresponding operator polynomial is \(1+T^2=(T+1)^2\), which divides \(T^6-1=(T^3+1)^2\). Hence the kernel has dimension two and the image dimension \(6-2=4\).
Alternatively, \(x+x^4=0\) is equivalent to \(x^4=x\), so the kernel is exactly \(\F_4\subseteq\F_{64}\).
\end{solution}

\solutionheading{ex:trace-adjoint}
\begin{solution}
For any \(x,y\in E\), since \(\Tr(x),\Tr(y)\in k\), the \(k\)-linearity of the trace gives
\[
\Tr(\Tr(x)y)=\Tr(x)\Tr(y)=\Tr(x\Tr(y)).
\]
Thus the trace operator is self-adjoint. Its image is \(k\), so \(k=\im\Tr=(\ker\Tr)^\perp\).
\end{solution}

\solutionheading{ex:linearized-not-ordinary}
\begin{solution}
Both are monic polynomials without repeated roots whose full root sets are the respective subspaces, so ordinary divisibility is equivalent to containment of root sets.
The ordinary least common multiple has root set \(U\cup V\), whereas \(L_{U+V}\) has root set \(U+V\), generally a different set.
For example, choose \(\alpha\in\F_4\setminus\F_2\), and take \(U=\{0,1\}\) and \(V=\{0,\alpha\}\).
The former least common multiple is \(X(X-1)(X-\alpha)\) of degree three; the latter polynomial is \(X^4-X\) of degree four.
\end{solution}

\solutionheading{ex:supp-symbolic-division}
\begin{solution}
Take \(Q(Y)=Y^4+Y^2+Y\); direct substitution and collection of terms in characteristic two give \(Q(X^2+X)=X^8+X\).
In operator notation this is simply \(T^3+1=(T^2+T+1)(T+1)\).
Since \(Q(Y)\) has no constant term, the ordinary polynomial \(X^2+X\) also divides this composite.
\end{solution}

\solutionheading{ex:wu-noncommuting}
\begin{solution}
We have \(L\circ M=\alpha X^4\), whereas \(M\circ L=\alpha^2X^4\), so they differ.
On \(\F_4\) their values are \(x\mapsto\alpha x\) and \(x\mapsto\alpha^2x\), respectively, and differ at \(x=1\).
However, \(X^4\) and \(X\) are distinct polynomials representing the same function on \(\F_4\),
since their difference is the root polynomial \(X^4-X\).
\end{solution}

\solutionheading{ex:wu-centralizer}
\begin{solution}
We have \(A\sigma^s=\sum_iM_{a_i}\sigma^{i+s}\), whereas \(\sigma^sA=\sum_iM_{a_i^{q^s}}\sigma^{i+s}\). Reducing exponents modulo \(n\) only makes a cyclic rearrangement,
so uniqueness of the linearized representation makes commutation equivalent to \(a_i^{q^s}=a_i\) for all coefficients.
The fixed field is \(\F_{q^g}\), of dimension \(g\) over \(k\); the \(n\) coefficients can be chosen independently, giving dimension \(ng\).
\end{solution}

\solutionheading{ex:ore-dickson-comparison}
\begin{solution}
When all coefficients lie in \(k\), we have \(\mathcal L_a=a(\sigma)\). By Theorem~\ref{thm:linearized-gcd},
\[
\dim_k\ker\mathcal L_a=\deg\gcd(a,T^n-1).
\]
A linear endomorphism of a finite vector space is a permutation if and only if its kernel is zero, giving the criterion in the first part.
If \(u(T)a(T)+v(T)(T^n-1)=1\), substitution of \(\sigma\) gives \(u(\sigma)a(\sigma)=\id_E\).
Thus the inverse map is \(u(\sigma)\); reducing \(u\) modulo \(T^n-1\) gives a linearized representative of degree less than \(q^n\).

When the coefficients lie in \(E\), we can still write \(L=\sum_iM_{a_i}\sigma^i\), but
\[
(M_a\sigma^i)(M_b\sigma^j)=M_{ab^{q^i}}\sigma^{i+j}.
\]
This differs from multiplication in the commutative polynomial ring \(E[T]\), so a Bézout identity there cannot directly be converted into an operator identity of the same form.
Proposition~\ref{prop:dickson-rank} gives \(\dim_k\ker L=n-\rank_E D_L\), so the permutation property is equivalent to \(\det D_L\ne0\).

In the third part, \(T^2-1=(T+1)^2\) whereas \(1+\omega\ne0\), so the displayed greatest common divisor is one.
On the other hand,
\[
D_L=\begin{pmatrix}1&\omega\\ \omega^2&1\end{pmatrix},\qquad
\det D_L=1-\omega^3=0.
\]
Solving \(x+\omega x^2=x(1+\omega x)=0\) directly gives \(\ker L=\{0,\omega^2\}\).
Thus \(L\) is not a permutation; coprimality in the commutative polynomial ring does not detect invertibility of the operator here.
\end{solution}

\solutionheading{v70:ex-two-dimensional-intersections}
\begin{solution}
Write \(H(W)=\F_{q^d}\). Since \(W\) is an \(H(W)\)-space, \(d\mid2\), so \(d=1\) or \(d=2\).
If \(d=2\), then \(2\mid n\), and \(W\) is one-dimensional over \(\F_{q^2}\); thus \(W=b\F_{q^2}\).
Conversely, the maximal field of scalars of such a space is exactly \(\F_{q^2}\): if \(aW\subseteq W\), then \(ab\in b\F_{q^2}\), whence \(a\in\F_{q^2}\).
Example~\ref{v70:ex-trace-hyperplane-intersection} shows that its intersection dimensions can only be zero or two, giving the first set of values.

In every other case, \(d=1\), so \(N_2=q-1\). The counting identity in Proposition~\ref{v70:prop-intersection-count} gives
\[
 (q-1)N_1+(q^2-1)(q-1)=(q^2-1)^2.
\]
Therefore
\[
 N_1=\frac{(q^2-1)((q^2-1)-(q-1))}{q-1}=q(q^2-1).
\]
Using \(N_0+N_1+N_2=q^n-1\) now yields \(N_0=q^n-q^3\).
If \(n=2\), the only two-dimensional subspace is \(E=\F_{q^2}\), which belongs to the first case. Thus the second case requires \(n\ge3\).
Finally, each scalar multiple is represented by \(q^2-1\) or \(q-1\) nonzero parameters, respectively. The two totals are therefore
\[
 \frac{q^n-1}{q^2-1},\qquad \frac{q^n-1}{q-1}.
\]
\end{solution}

\solutionheading{v116:ex-inverse-sum}
\begin{solution}
Put \(Q=L_U/X=\prod_{u\ne0}(X-u)\). Its logarithmic derivative at zero gives \(Q'(0)/Q(0)=\sum_{u\ne0}u^{-1}\), while \(Q(0)=c_0\) and \(Q'(0)=c_1\).
In dimension four, \(c_1=0\) says exactly that the all-one vector belongs to \(\mathcal C(U)\); it does not determine the other relations.
The cubic is irreducible, so the given space is all of \(\F_8\). Inversion permutes its nonzero elements, and the reciprocal rank is three.
The rational function model has rank six and \(\deg M=6\), whereas \(\deg\alpha=3\). Proposition~\ref{v117:prop-specialization-loss} gives \(\dim(W\cap(X^3+X+1))=6-3=3\). These three directions are precisely the additional relations after denominators are cleared: they vanish at \(\alpha\) but are nonzero polynomials.
\end{solution}

\Needspace{9\baselineskip}
\section*{Solutions to Chapter~\ref{chap:normal}}
\pdfbookmark[1]{Solutions to Chapter \getrefnumber{chap:normal}}{solutions.chap:normal}

\solutionheading{v88:ex-moore-sign}
\begin{solution}
Since \(x^4=x\), the determinant is \(x^2-x^4=x^2+x=\Tr_{4/2}(x)\).
Every \(t\in\F_2\) has \(t^2=t\), so this equals \(\Tr_{4/2}(x)^2\).
For \(n=2\), \(\lfloor(n-1)/2\rfloor=0\), agreeing with the formula.
\end{solution}

\solutionheading{v88:ex-trace-quantifier}
\begin{solution}
We have \(\Tr_{8/2}(1)=1+1+1=1\), but the three conjugates of \(1\)
are linearly dependent. Since \(\ord_3(2)=2\),
\(T^3-1=(T-1)(T^2+T+1)\), with the second factor irreducible.
An element of degree three and nonzero trace has both Frobenius components
nonzero, and is therefore normal. The full-degree hypothesis matters.
\end{solution}

\solutionheading{v88:ex-affine-nine-cosets}
\begin{solution}
For \(x=a+b\xi\), the determinant of its conjugate coordinate columns
is \(ab\). At \(b=0\), no element is normal; for either nonzero \(b\),
\(a=\pm1\) gives the two normal elements. The three counts are \(0,2,2\).
In local coordinates, \(U=k\) allows the \(T-1\) component to vary freely,
while the \(T+1\) component is fixed and is nonzero precisely when \(b\ne0\). The formula gives \(3(1-1/3)=2\).
\end{solution}

\solutionheading{v88:ex-correlation-f8}
\begin{solution}
A self-adjoint element of \(\F_2[T]/(T^3-1)\) has the form
\(a+c(T+T^2)\). Evaluation at \(T=1\) yields \(a\), so \(a=1\).
For \(c=1\), the element \(1+T+T^2\) divides \(T^3-1\) and is not a unit.
Only \(b=1\) remains. The correlation polynomial of any normal element
must be such a self-adjoint unit with value one at \(T=1\); it is therefore
\(1\), exactly the self-duality condition.
\end{solution}

\solutionheading{v88:ex-primitive-normal-nine}
\begin{solution}
The multiplicative group has order eight, and
\((1+\xi)^2=2\xi\), \((1+\xi)^4=2=-1\), so \(1+\xi\) has order eight.
Its powers with exponents \(1,3,5,7\) are the four primitive elements,
namely \(\pm1\pm\xi\). The conjugate-coordinate determinant for
\(a+b\xi\) is \(ab\), so exactly those four elements are normal.
\end{solution}

\solutionheading{v58:ex-normal-traces-binary-six}
\begin{solution}
A normal element has nonzero trace, which in \(\F_2\) must equal one, so the two target values have the same trace to the intersection field \(k\).
Put \(P=T^2+T+1\). Then
\[
 F=T^6-1=(T+1)^2P^2,\qquad
 G=\lcm(T^2-1,T^3-1)=(T+1)^2P.
\]
Thus the total number of solutions is \(2^{6-4}=4\), and the number of normal solutions is
\[
 \frac{\varphi_2((T+1)^2P^2)}{\varphi_2((T+1)^2P)}
 =\frac{(4-2)(16-4)}{(4-2)(4-1)}=4.
\]
Hence every solution is normal. If the two intermediate traces of \(x\) are normal, using them as targets shows that \(x\) is normal;
the converse follows from the proposition on descent of normality.
Here \(G\) already contains every distinct irreducible factor of \(F\), so every lift of the fixed unit residue class in the remaining multiplicity directions is still a unit.
The minimal polynomial \(F\) of Frobenius has repeated factors, so Frobenius is not diagonalizable, but the cyclic module and reduction of its units still apply.
\end{solution}

\solutionheading{ex:normal-count8}
\begin{solution}
Since \(T^3-1=(T-1)(T^2+T+1)\), the number of normal elements is
\[
2^3(1-2^{-1})(1-2^{-2})=3.
\]
The example provides three distinct normal elements \(\eta,\eta^2,\eta^4\), so these exhaust the normal elements.
\end{solution}

\solutionheading{ex:normal-p-power}
\begin{solution}
Nondegeneracy of the trace pairing makes taking orthogonal complements a bijection between hyperplanes and one-dimensional subspaces.
By \(\sigma^*=\sigma^{-1}\), this bijection preserves Frobenius invariance:
if \(H\) is invariant and \(c\in H^\perp\), then for \(x\in H\) we have \(\Tr(x\sigma c)=\Tr(\sigma^{-1}x\,c)=0\), so \(\sigma c\in H^\perp\); the converse is similar.

A one-dimensional space \(kc\) with \(c\ne0\) is invariant if and only if \(\sigma c=ac\) for some \(a\in k^\times\).
From \(\sigma^n c=c\) we obtain \(a^n=1\). Conversely, if \(a\in k^\times\) and \(a^n=1\),
then \(T-a\mid T^n-1\), and Theorem~\ref{thm:invariant-duality} gives \(\dim_k\ker(\sigma-a\id)=1\).
Thus each such \(a\) corresponds to exactly one invariant line and to the hyperplane
\[
H_a=\{x\in E:\Tr_{E/k}(cx)=0\},
\qquad 0\ne c\in\ker(\sigma-a\id).
\]
Nonzero \(c\) on the same line differ only by a nonzero base-field multiple, so the definition is independent of the choice of \(c\).
Different \(a\) give different lines and hyperplanes.
The equation \(a^n=1\) has exactly \(\gcd(n,q-1)\) solutions in the cyclic group \(k^\times\), giving the required count.

If \(n=p^a\), then since \(p\nmid q-1\) the only possible eigenvalue is one, whose eigenspace is \(k\).
Taking \(c=1\) gives precisely the hyperplane \(\ker\Tr_{E/k}\).
\end{solution}

\solutionheading{ex:normal-selfdual8}
\begin{solution}
The element \(\eta=\alpha+1\) satisfies \(\Tr(\eta)=1\), so \(\Tr(\eta^2)=1\). Also
\[
\eta\eta^2=(\alpha+1)(\alpha^2+1)=\alpha^2,
\]
so \(\Tr(\eta\eta^2)=0\). By Frobenius invariance, each basis vector has trace pairing one with itself and zero with the others; hence the basis is self-dual.
The extension \(\F_{2^4}/\F_2\) has a self-dual basis, but its degree is divisible by four, so it has no self-dual normal basis.
For \(\F_{3^4}/\F_3\) the base-field order is odd and the extension degree even, so even a self-dual basis does not exist.
\end{solution}

\solutionheading{ex:minpoly-vs-qorder}
\begin{solution}
We have \(m_\alpha(X)=X^3+X+1\). Substituting \(\alpha=\eta+1\) gives \(m_\eta(X)=X^3+X^2+1\). Both are minimal polynomials of field elements and have degree three.
Since \(\alpha+\alpha^2+\alpha^4=0\), we have \(\mu_\alpha(T)\mid T^2+T+1\).
The elements \(\alpha,\alpha^2\) are independent over \(\F_2\), so \(\mu_\alpha=T^2+T+1\).
The element \(\eta\) is normal, so \(\mu_\eta=T^3-1\). Thus the two elements have the same orbit length but different dimensions of their orbit spans.
\end{solution}

\solutionheading{ex:supp-linear-primitive}
\begin{solution}
We have \(\mathcal L_h=X^4+X^2+X\). The polynomial \(h\) is irreducible over \(\F_2\) and has order three, equal to \(2^2-1\).
Thus the three nonzero elements of its two-dimensional root space form a single Frobenius orbit, and \(X^3+X+1\) is irreducible.
\end{solution}

\solutionheading{ex:supp-nonnormal-six}
\begin{solution}
Under the module isomorphism \(k[T]/(F)\simeq E\) given by a normal basis, \(\alpha_\mu\) corresponds to the residue class of \(F/\mu\).
Every \(h\in k[T]\) satisfies
\[
h(\sigma)\alpha_\mu=0
\quad\Longleftrightarrow\quad F\mid h(F/\mu)
\quad\Longleftrightarrow\quad \mu\mid h.
\]
Thus its monic annihilating polynomial of least degree is exactly \(\mu\), and the map
\(k[T]/(\mu)\to k[\sigma]\alpha_\mu\), \(\bar h\mapsto h(\sigma)\alpha_\mu\)
is a module isomorphism, giving the dimension identity.
Moreover, \(\sigma^r(\alpha_\mu)=\alpha_\mu\) if and only if \(\mu\mid T^r-1\),
so both the orbit length and the field extension degree equal \(\ord(\mu)\).

Write \(S=1+T+\cdots+T^{n-1}=F/(T-1)\).
Since \(T^m-1\) has no repeated factors, the multiplicity of \(T-1\) in \(F\) is exactly \(p^a\).
Except for this factor, \(S\) and \(F\) have the same irreducible-factor multiplicities; thus, assuming \(\mu\mid F\),
\[
\Tr_{E/k}(\alpha_\mu)=S(\sigma)\alpha_\mu\ne0
\quad\Longleftrightarrow\quad \mu\nmid S
\quad\Longleftrightarrow\quad (T-1)^{p^a}\mid\mu.
\]

Consequently, assuming \(\mu\mid F\) so that the construction is possible, the two additional conditions in the second part say that \(\alpha_\mu\) has nonzero trace and degree \(n\), respectively;
the conclusion \(\mu=F\) is equivalent to normality of \(\alpha_\mu\).
Conversely, the monic annihilating polynomial of least degree of every nonzero element is among these factors, so the proposed implication is exactly equivalent to the property in Theorem~\ref{thm:trace-normal-classification}.
The answer is precisely \(n=p^a\), or \(n=\ell\ne p\) prime with \(\ord_\ell(q)=\ell-1\).
In every other case, the proof of that theorem already gives a proper factor \(\mu\) satisfying the two conditions; substituting it into \(\alpha_\mu\) gives a counterexample without separately choosing roots in the field.

In particular, when \(q=2,n=6\) take \(\mu=(T+1)^2(T^2+T+1)\). Then
\(F/\mu=T^2+T+1\), so \(\alpha_\mu=\eta+\eta^2+\eta^4\), exactly the construction in Example~\ref{exm:full-degree-nonnormal-f64}.
\end{solution}

\solutionheading{ex:wu-three-generators}
\begin{solution}
Example~\ref{exm:quartic-model} proves that \(X^4+X+1\) is irreducible and that its root \(\alpha\) has order 15, so \(\alpha\) is primitive; the cubic coefficient of its minimal polynomial is zero, so its absolute trace is zero.
Since the extension degree four is a power of the characteristic two, normality is equivalent to nonzero trace; thus \(\alpha\) is not normal.
The element \(\alpha^3\) has order five, degree four, minimal polynomial \(\Phi_5=X^4+X^3+X^2+X+1\), and trace one, so it is normal but not primitive.
The element \(\alpha^{-1}\) has minimal polynomial \(X^4+X^3+1\) and is again primitive with trace one, hence normal.
The two primitive quartics each have four roots, giving eight primitive elements;
there are \(16-8=8\) elements of nonzero trace, all normal.
Four are roots of \(X^4+X^3+1\) and four of \(\Phi_5\), so there are exactly four primitive normal elements.
\end{solution}

\solutionheading{ex:wu2-invariant-complements}
\begin{solution}
Invariant subspaces correspond bijectively to monic factors, and each exponent can be chosen independently from zero to \(e_i\), giving the product.
When \(p\nmid n\), all \(e_i=1\). Given \(U_h\), take \(U_{F/h}\); the intersection and sum formulas show that their intersection is zero and their sum is \(E\).
When some \(e_i>1\), consider \(U_{P_i}\). If it had an invariant complement \(U_g\), then \(\gcd(P_i,g)=1\) and \(\lcm(P_i,g)=F\). The first equality forces \(P_i\nmid g\), but the second requires the multiplicity of \(P_i\) in the least common multiple to be \(e_i>1\), a contradiction. Since \(F\) has repeated factors exactly when \(p\mid n\), the assertion follows.
\end{solution}

\solutionheading{ex:wu7-rootfield-extremes}
\begin{solution}
The polynomial \(\mathcal L_f\) is monic of degree \(q^d\) with nonzero linear coefficient, so it is exactly the root polynomial of \(W_f\).
If \(W_f\) is a subfield, its size is \(q^d\), and Proposition~\ref{prop:subfield-root-polynomial} gives \(\mathcal L_f(X)=X^{q^d}-X\). Comparing the coefficients of each \(X^{q^i}\) gives \(f=T^d-1\).
Conversely, the linearized root set of this \(f\) is exactly \(\F_{q^d}\).

Corollary~\ref{cor:linearized-splitting-field} gives splitting-field degree \(N=\ord(f)\).
The order bound in the unit group of the quotient algebra and Proposition~\ref{prop:maximal-order-primitive} give \(N\le q^d-1\), with equality if and only if \(f\) is primitive.
Then \(f\) is irreducible, and every nonzero root has Frobenius annihilator polynomial \(f\),
so every nonzero root has orbit length \(\ord(f)=q^d-1\). There are only this many nonzero roots, so they form one orbit.

In both cases the root space has dimension \(d\) and \(q^d\) elements;
in the subfield case the splitting-field degree is \(d\), whereas in the maximal-orbit case it is \(q^d-1\).
The dimension and cardinality describe the size of the root space; the last gives the number of Frobenius iterations required to fix all roots.
The two cases coincide only when \(q=2,d=1\).
\end{solution}

\solutionheading{wfr:ex-norm-classes}
\begin{solution}
Here \(K=\F_4\) and \(d=2,m=3\), and the norm is \(a^{21}\). It maps onto the group \(K^\times\) of order three, with \(21\) preimages of each norm value.
The norm value one forms one \(k\)-conjugacy class, and the other two elements \(\omega,\omega^2\) form another. Thus there are two similarity classes.
The first contains \(21\) weights and has minimal polynomial \(T^3+1\) and characteristic polynomial \((T^3+1)^2\), with no cyclic vectors.
The second contains \(42\) weights, has both minimal and characteristic polynomial \(T^6+T^3+1\), and all its operators have cyclic vectors.
\end{solution}

\solutionheading{wfr:ex-nonprimitive-norm}
\begin{solution}
The norm \(E^\times\to K^\times\) is surjective, so the required \(a\) exists.
We have \(i\notin\F_3\), with minimal polynomial \(P_i=T^2+1\), so \(r=d=2\) and \(m=2\). Theorem~\ref{wfr:thm-norm-decomposition} gives
\[
 \mu_{S_a}(T)=\chi_{S_a}(T)=T^4+1.
\]
This has degree \(\dim_k E=4\), so a cyclic vector exists. But \(i\) has multiplicative order four and \(|K^\times|=8\), so it generates \(K\) over \(k\) without generating the entire multiplicative group \(K^\times\).
\end{solution}

\solutionheading{v40:ex-necklace-models}
\begin{solution}
Write \(a=\sum_{j=0}^{m-1}a_jq^j\). Modulo \(q^m-1\) we have
\[
qa\equiv a_{m-1}+\sum_{j=0}^{m-2}a_jq^{j+1},
\]
so the digit word rotates cyclically. The all-\(q-1\) word represents \(q^m-1\) and gives the same residue class as the zero word, so it must be removed; every other word represents a unique class, whose first return time is its least cyclic period.
In the normal basis \((\eta,\sigma\eta,\ldots,\sigma^{m-1}\eta)\),
\[
\sigma\left(\sum_{j=0}^{m-1}b_j\sigma^j\eta\right)
=\sum_{j=0}^{m-1}b_j\sigma^{j+1}\eta,
\]
so the coefficients rotate in the same direction.

Let \(P_r\) be the number of primitive words of length \(r\) on an alphabet of \(q\) symbols.
Every word of length \(r\) is a repetition of a primitive word of length \(e\mid r\), so
\(q^r=\sum_{e\mid r}P_e\). Möbius inversion gives
\(P_r=\sum_{e\mid r}\mu(r/e)q^e\).
Each primitive necklace has \(r\) words; repeating \(m/r\) times gives exactly one orbit of words of length \(m\) with least period \(r\), so there are \(P_r/r=I_q(r)\) such orbits.
Removing the all-\(q-1\) word only affects period one; the two models therefore have equal counts when \(r>1\) and differ by one when \(r=1\).

For \(q=2,m=6\), the lengths and orbit counts are
\[
\begin{array}{c|rrrr}
r&1&2&3&6\\\hline
\text{Exponent orbits modulo }63&1&1&2&9\\
\text{Frobenius orbits in }\F_{64}&2&1&2&9
\end{array}
\]
The totals weighted by length are \(1+2+6+54=63\) and \(2+2+6+54=64\), respectively.
A normal basis and labels for the alphabet give word coordinates to field elements; digit words represent integer residue classes.
The same form of cyclic action does not mean that these correspondences preserve addition or multiplication; the two full sets even have different sizes.
\end{solution}

\solutionheading{v82:ex-universal-normal-translation}
\begin{solution}
First suppose \(p\nmid n\). If every full-degree element has a normal translate, take any full-degree \(\alpha\) of nonzero trace.
Some translate is normal, so all the conditions for \(P\ne T-1\) in Proposition~\ref{v82:prop-normal-translates} hold. Since \(\alpha\) itself has nonzero trace, it is normal as well.
Conversely, if every full-degree element of nonzero trace is normal, precisely one parameter in each full-degree translate gives trace zero, and all other \(q-1\) elements are normal.
Thus the two properties are equivalent. Theorem~\ref{thm:trace-normal-classification}, with \(n\ge2\), gives exactly the case in which \(n=\ell\ne p\) is prime and \(\ord_\ell(q)=\ell-1\).

Now let \(p\mid n\), and write \(n=p^a m\), with \(a\ge1\) and \(p\nmid m\). Translation preserves normality, so the property asks whether every full-degree element is normal.
Fix a normal element \(\eta\). For any nonconstant monic divisor \(\mu\mid F=T^n-1\), the element
\[
 \alpha_\mu=(F/\mu)(\sigma)\eta
\]
has Frobenius annihilator \(\mu\) and degree \(\ord(\mu)\), by Theorem~\ref{thm:root-module-count}.
If \(m>1\), use the proper divisor from the proof of Theorem~\ref{thm:trace-normal-classification}:
\[
 \mu=(T-1)^{p^a}\frac{T^m-1}{T-1}.
\]
It has order \(n\), but \(\mu\ne F\), so \(\alpha_\mu\) has full degree without being normal.
If \(m=1\) and \(p^a\ne2\), put \(b=p^{a-1}+1\). Then \(p^{a-1}<b<p^a\), and the order formula with multiplicities gives
\[
 \mu=(T-1)^b\ne F,\qquad \ord(\mu)=p^a=n,
\]
again producing the required counterexample.

Finally, if \(n=p=2\), every full-degree \(\theta\) has \(\theta^q\ne\theta\), so \(\Tr_{E/k}(\theta)=\theta+\theta^q\ne0\); it is therefore normal.
All its base-field translates are normal too. If \(n=1\), every \(\theta+k=k\) contains \(q-1\) nonzero elements, precisely the normal elements for a degree-one extension. Thus \(n=1\) is an additional trivial boundary case.
\end{solution}

\solutionheading{v82:ex-subfield-translates}
\begin{solution}
Put \(R=k[T]/(F)\). The element \(B(\sigma)\eta\) belongs to \(L\) if and only if \(F\mid DB\), equivalently, \(M\mid B\).
Thus \(L\) corresponds to the ideal \((\overline M)\), whose elements are uniquely \(\overline{MV}\) as \(V\) runs through residues modulo \(D\).
The translate \(\theta+L\) therefore corresponds to \(A+(\overline M)\), and normality is equivalent to \(A+MV\) being a unit in \(R\).

Write
\[
 F=\prod_{P\mid F}P^{e_P},\qquad
 D=\prod_{P\mid F}P^{f_P},\qquad 0\le f_P\le e_P,
\]
where the factors \(P\) are distinct, monic, and irreducible.
If \(f_P<e_P\), then \(P\mid M\), so the residue modulo \(P\) remains \(A\bmod P\).
If it is zero, the whole translate is nonnormal; if it is nonzero, all \(q^{f_P\deg P}\) parameters in this component are admissible.
If \(f_P=e_P\), then \(M\) is a unit in \(k[T]/(P^{e_P})\), and \(V\) runs freely through that component. Exactly
\(q^{e_P\deg P}(1-q^{-\deg P})\) parameters give units.
The Chinese remainder theorem makes the components independent, so
\begin{equation}\label{v82:eq-subfield-normal-count}
 \#\{u\in L:\theta+u\text{ is normal}\}
 =\begin{cases}
 0,&\gcd(A,M)\ne1,\\[1mm]
 \displaystyle q^d\prod_{\substack{P\mid F\\ f_P=e_P}}
       (1-q^{-\deg P}),&\gcd(A,M)=1.
 \end{cases}
\end{equation}
The product is over distinct irreducible factors. The nonzero branch also equals \(\varphi_q(F)/\varphi_q(M)\), with \(\varphi_q(1)=1\).
If \(M=1\), the compatibility condition is empty and the formula counts all normal elements in \(E\).
The condition \(P\nmid A\) is also equivalent to \((F/P)(\sigma)\theta\ne0\), so although the proof uses normal coordinates, the conclusion is independent of the chosen normal element \(\eta\).

If \(v_p(d)<v_p(n)\), every factor of \(D\) has strictly smaller multiplicity than in \(F\).
There is consequently no factor with \(f_P=e_P\), and \(\gcd(A,M)=1\) is equivalent to \(\gcd(A,F)=1\).
The count is therefore \(q^d\) or zero, according as \(\theta\) is normal or not.
Thus the capacity of a subfield translation to repair normality depends on primary multiplicities, not merely on the dimension of the parameter space.

Finally, take \(q=2,n=6\), and put \(P=T^2+T+1\), so \(F=(T+1)^2P^2\).
For \(d=2\), we have \(D=(T+1)^2\) and \(M=P^2\). A translate contains normal elements exactly when \(P\nmid A\), in which case there are \(2^2(1-2^{-1})=2\).
For \(d=3\), we have \(D=(T+1)P\) and \(M=(T+1)P\); the number is eight or zero, and is eight only when \(\theta\) itself is normal.
For \(\theta=\eta+\eta^2\), the coordinate polynomial is \(A=T+1\), and
\[
 \mu_\theta=F/(T+1)=(T+1)P^2,\qquad
 \ord(\mu_\theta)=6,\qquad \deg\mu_\theta=5.
\]
Thus \(\theta\) has full degree but is not normal; \(\theta+\F_4\) contains two normal elements, while \(\theta+\F_8\) contains none.
\end{solution}

\solutionheading{v145:ex-right-factors}
\begin{solution}
The linearized polynomial is $X^{Q^r}-X$, with root space $\F_{Q^r}$. Since $r\mid s$, the return $x\mapsto x^{Q^s}$ is the identity. By Theorem~\ref{v145:thm-skew-roots}, the answer is the number of $j$-dimensional subspaces of an $r$-dimensional $\F_Q$-space:
\[
 \prod_{i=0}^{j-1}\frac{Q^r-Q^i}{Q^j-Q^i}.
\]
The numerator counts ordered linearly independent $j$-tuples; the denominator counts ordered bases in each $j$-dimensional subspace. For $j=0$, the empty product is one.
\end{solution}

\solutionheading{v148:ex-factor-lattices}
Local types are partitions $\lambda\vdash e$ with $\ell(\lambda)\le k$. Type $(e)$ has a chain of submodules given by successive power ideals, and formula~\eqref{v148:eq-local-type} shows that this type always occurs. If $\ell(\lambda)\ge2$, the submodule killed by $\pi$ has dimension at least two over $\F_{q^{\deg\pi}}$ and contains distinct one-dimensional subspaces, so the right factors do not form a chain. Thus every member has the property exactly when $e=1$ or $k=1$; at least one member has it for all $e,k$.

\section*{Solutions to Chapter~\ref{chap:algorithms}}
\pdfbookmark[1]{Solutions to Chapter \getrefnumber{chap:algorithms}}{solutions.chap:algorithms}

\solutionheading{ex:rabin-prime-degree}
\begin{solution}
When the degree \(n\) is prime, Rabin's criterion becomes \(X^{q^n}\equiv X\pmod f\) and \(\gcd(f,X^q-X)=1\).
The second condition only excludes linear factors; a prime-degree counterexample showing that the first cannot be omitted appears in Exercise~\ref{ex:wu-rabin-missing}.
For the second part, \((X^2+X+1)^2\in\F_2[X]\) has no root in \(\F_2\) but is a reducible quartic.
\end{solution}

\solutionheading{ex:squarefree-example}
\begin{solution}
Write \(A=X+1\) and \(B=X^2+X+1\). Then \(f'=A^2B^2\), with initial values \(c=A^2B^2\) and \(w=A\).
The values \(c,w\) in the following table are those at the start of each step.
\begin{center}
\begin{tabular}{ccccc}
\toprule
\(i\)&\(c\)&\(w\)&\(y=\gcd(w,c)\)&\(z=w/y\)\\
\midrule
1&\(A^2B^2\)&\(A\)&\(A\)&1\\
2&\(AB^2\)&\(A\)&\(A\)&1\\
3&\(B^2\)&\(A\)&1&\(A\)\\
\bottomrule
\end{tabular}
\end{center}
The third pass outputs the factor \(A\) with multiplicity three; after updating, \(w=1\) and \(c=B^2\).
Taking the square root of the residual \(c\) gives \(B\), which is square-free; multiplying its recursively obtained multiplicity one by two gives the factor \(B\) with multiplicity two.
\end{solution}

\solutionheading{v56:ex-factorization-pipeline}
\begin{solution}
\textbf{Separating multiplicities.}
In \(\F_3[X]\),
\[
F'=2X^7+2X^6+X^3+X,\qquad
c=\gcd(F,F')=X^3-1,\qquad
w=X^5+2X^4+X^3+X^2+2.
\]
In the first pass, \(y=\gcd(w,c)=1\) and \(z=w\), so the output is the square-free polynomial of multiplicity one
\[
z_1=w.
\]
After updating, \(w=1\), with residual \(c=X^3-1=(X-1)^3\).
Taking its cube root gives \(X-1\); multiplying the recursive multiplicity one by three gives the other output \(z_3=X-1\).
Thus \(F=z_1z_3^3\). Direct calculation also gives
\(\gcd(z_1,z_1')=1\) and \(\gcd(z_1,z_3)=1\), as required for a square-free factorization.

\textbf{Grouping by degree.}
The linear polynomial \(z_3\) is already irreducible, so only \(z_1\) needs to be factored.
Starting with \(r=z_1\) and \(u=X\), the first step gives
\[
u=X^3\bmod z_1=X^3,\qquad
g_1=\gcd(z_1,X^3-X)=X+1.
\]
After removing the product of the linear factors,
\[
r=H=\frac{z_1}{X+1}=X^4+X^3+X+2,
\]
and the remainder of \(u\) modulo \(H\) is still \(X^3\).
In the second step,
\[
X^4\equiv2X^3+2X+1,\qquad X^8\equiv1\pmod H
\]
give
\[
u\leftarrow u^3=X^9\bmod H=X,\qquad
g_2=\gcd(H,u-X)=H.
\]
Distinct-degree factorization has already removed all linear factors, so \(H\) is exactly a product of two distinct irreducible quadratics.

\textbf{Separating equal-degree factors.}
Take \(h=X\). Since \(H(0)=2\), we have \(\gcd(H,h)=1\).
Each field component has order \(Q=3^2=9\), so its quadratic character is given by \(h^{(Q-1)/2}=X^4\).
Computing modulo \(H\) gives
\[
b=2X^3+2X+1,\qquad b-1=2X(X^2+1),
\]
hence
\[
g_+=\gcd(H,b-1)=X^2+1,\qquad
H/g_+=X^2+X+2.
\]
In the first quotient \(X^4=1\), and in the second \(X^4=2=-1\), so the chosen element has different quadratic-character values in the two components.
Neither quadratic has a root in \(\F_3\), so both are irreducible. The final result is
\[
\boxed{F=(X-1)^3(X+1)(X^2+1)(X^2+X+2).}
\]
Direct multiplication gives
\[
(X^2+1)(X^2+X+2)=H,\qquad
(X+1)H=z_1,\qquad (X-1)^3z_1=F,
\]
confirming the factorization and all multiplicities.
Square-free factorization preserves multiplicities, distinct-degree factorization determines root-orbit lengths, and equal-degree splitting separates root orbits of the same length; together the three outputs give the complete factorization.
\end{solution}

\solutionheading{ex:berlekamp7}
\begin{solution}
The 2-cyclotomic cosets modulo seven are \(\{0\}\), \(\{1,2,4\}\), and \(\{3,6,5\}\), so \(X^7-1\) has three distinct irreducible factors.
The Berlekamp fixed algebra is isomorphic to \(\F_2^3\) and hence has dimension three.
Each root orbit corresponds to a field component, whose Frobenius fixed part has dimension one.
\end{solution}

\solutionheading{ex:berlekamp-value-multiplicity}
\begin{solution}
The polynomial \(P\) has no root in \(\F_2\) and is therefore irreducible. Likewise \(Q\) has no base-field root, and modulo \(P\) we have \(X^3=1\) and \(X^4=X\), hence \(Q\equiv1\pmod P\).
The polynomial \(P\) is the only monic irreducible quadratic over \(\F_2\), so \(Q\) has neither a linear nor a quadratic factor and is irreducible.
Modulo \(P\) we have \(h=0\), and modulo \(Q\) we have \(h=1\); the Chinese remainder theorem therefore gives \(h^2=h\).
Its two values correspond to field components of dimensions two and four, so
\[
m_h(Y)=Y(Y+1),\qquad \chi_h(Y)=Y^2(Y+1)^4.
\]
In characteristic two, \(\chi_h'=0\), hence \(\gcd(\chi_h,\chi_h')=\chi_h\), and the required quotient is one.
The roots of the original polynomial \(f\) are simple, but a single eigenvalue of the multiplication operator is repeated throughout a root orbit;
when the characteristic divides this multiplicity, one differentiation loses the corresponding information.
Obtaining \(m_h=Y^2+Y\) directly from the first linear relation among \(1,h,h^2\) avoids this loss.

For \(F=X^3+X^2\), we have \(F'=X^2\) and \(\gcd(F,F')=X^2\), so \(F/\gcd(F,F')=X+1\), losing the original root zero.
The part not retained by the derivative may be a \(p\)th power; the residual part here is \(X^2\), whose square root \(X\) must be taken to recover this root and its multiplicity.
Thus square-free factorization must retain and recursively process the residual part with zero derivative, rather than merely take one derivative quotient.
\end{solution}

\solutionheading{ex:supp-berlekamp-repeated}
\begin{solution}
There is only one distinct irreducible factor, so the fixed algebra has dimension one and consists only of the two constants.
Repeated factors add no new root orbits; hence the fixed algebra counts one distinct factor but cannot exclude its square.
Square-free factorization must precede the use of this dimension-one condition as an irreducibility test.
\end{solution}

\solutionheading{ex:wu-rabin-missing}
\begin{solution}
Both factors are irreducible, so \(f\) has degree five and no root in \(\F_2\).
Nevertheless, \(f\) already has a nontrivial factorization.
Its two irreducible-factor degrees, two and three, do not divide five, so neither factor divides \(X^{2^5}-X\).
Thus the first congruence \(X^{32}\equiv X\pmod f\) fails.
The second condition tests only \(\gcd(f,X^2-X)=1\) and does hold here.
\end{solution}

\solutionheading{ex:wu2-idempotents}
\begin{solution}
In the decomposed algebra \(\prod_iK_i\), idempotency is the componentwise condition \(e_i^2=e_i\). In a field only zero and one satisfy this,
giving \(2^s\) choices. Components with value one are exactly the selected orbits; the equivariant-function model gives the same interpretation.
If \(e\) is represented by a polynomial, \(\gcd(f,e)\) collects the factors where its value is zero, and \(\gcd(f,e-1)\) those where its value is one.
Both greatest common divisors are taken monic; they are therefore coprime and satisfy
\[
f=\gcd(f,e)\,\gcd(f,e-1).
\]
The monicity of the left-hand side is why no extra constant factor is needed.
If the idempotent is neither zero nor one in \(k[X]/(f)\), both factors are nonconstant.
\end{solution}

\solutionheading{ex:wu3-f64-coordinates}
\begin{solution}
By the dual coordinate formula, the six coordinates of \(\alpha\) are \(\Tr(\alpha(\gamma^*)^{2^i})\) for \(0\le i<6\), respectively.
Substituting \(\gamma^*\) and its conjugates from the example gives \((0,0,0,0,1,1)\), hence \(\alpha=\gamma^{16}+\gamma^{32}\); this can also be checked by adding the two power-basis expansions given there.
The relative trace is
\[
\gamma+\gamma^4+\gamma^{16}
=\alpha^5+\alpha^4+\alpha^3+\alpha=u^2.
\]
The elements \(u,u^2\) are distinct nonzero elements of \(\F_4\) and are linearly independent over \(\F_2\), so they form a normal basis.
This is also a concrete instance of Proposition~\ref{prop:normal-descent} on relative-trace descent of normal elements.
\end{solution}

\solutionheading{ex:wu3-isomorphism-certificate}
\begin{solution}
By Theorem~\ref{thm:operator-algebra}, the ring of operators commuting with \(\sigma\) is \(R=\F_2[T]/((T+1)^4)\). Its units are precisely the classes nonzero modulo \(T+1\), namely those satisfying \(a(1)=1\),
of which there are \(2^4-2^3=8\). The corresponding operators automatically satisfy \(a(\sigma)(1)=a(1)=1\), so there are eight linear maps as stated.
Those preserving multiplication are exactly the four field automorphisms \(1,\sigma,\sigma^2,\sigma^3\) of \(E/\F_2\).
The other four do not preserve multiplication; Example~\ref{exm:linear-not-field} gives one of them.
\end{solution}

\solutionheading{v97:ex-sparse-uniqueness}
\begin{solution}
Put \(D=L_s(f)\). Degree minimality and \(f(0)\ne0\) imply that every degree-minimal multiple has nonzero coefficients at exponents \(0\) and \(D\).
If a proper subset of the support of \(P\) were dependent, there would be a nonzero multiple \(R\) supported on that subset.
If \(\deg R<D\), minimality fails immediately. If \(\deg R=D\), normalize \(R\) to be monic and subtract it from \(P\).
The result is a nonzero multiple of smaller degree with at most \(s\) terms, again a contradiction.
The same subtraction proves uniqueness of the monic relation on a fixed support.

For distinct \(P,Q\), the polynomial \(P-Q\) is nonzero of degree less than \(D\).
Their leading terms cancel, and both supports also contain the constant position, so
\[
 \operatorname{wt}(P-Q)\le\operatorname{wt}(P)+\operatorname{wt}(Q)-3.
\]
Minimality forces \(\operatorname{wt}(P-Q)\ge s+1\), giving the inequality.
For \(s=2,3\), one has \(2s<s+4\), so two distinct degree-minimal multiples cannot exist.

The given quartic has five nonzero coefficients, hence \(L_4(f)>4\), while
\[
 (X+1)f=X^5-X^2-X+1,\qquad
 (X+2)f=X^5+X^4-X^3-1.
\]
These are distinct monic four-term multiples. Thus \(L_4(f)=5\), and \(\mathcal M_4(f)\) is not a singleton.
\end{solution}

\solutionheading{v97:ex-sparse-base-change}
\begin{solution}
Fix an exponent set \(I\). The remainder matrix with columns \(X^i\bmod f\) is defined over \(k\), and scalar extension does not change its rank.
Thus a nonzero \(f\)-multiple over \(K\) supported within \(I\) exists if and only if such a multiple exists over \(k\).
Minimizing \(D\) over sets \(I\subseteq\{0,\ldots,D\}\) with \(|I|\le s\) proves the first claim.
This guarantees the existence of a base-field relation; it does not assert that every relation over \(K\) has coefficients in \(k\).

Substitution into a multiple attaining \(L_s(f)\) gives \(L_s(f(X^a))\le aL_s(f)\).
Conversely, suppose \(P=f(X^a)U\ne0\) and \(\operatorname{wt}(P)\le s\). Write uniquely
\[
 U=\sum_{r=0}^{a-1}X^rU_r(X^a),\qquad
 P=\sum_{r=0}^{a-1}X^r(fU_r)(X^a).
\]
Terms from different residues cannot cancel. Any nonzero component has at most \(s\) terms, so \(\deg(fU_r)\ge L_s(f)\), and hence
\(\deg P\ge r+a\deg(fU_r)\ge aL_s(f)\).

Finally, coefficient Frobenius is an automorphism of \(k\), preserving divisibility, degree, and number of terms.
Thus \(L_s(f^{[p^u]})=L_s(f)\).
Apply the second claim to \(f^{[p^u]}\) and use \(f(X)^{p^u}=f^{[p^u]}(X^{p^u})\) to obtain the third.
\end{solution}

\solutionheading{v97:ex-norm-cancellation}
\begin{solution}
Coefficient comparison gives \(a_j^r=a_{d-j}/a_0\), whence
\[
 t_{d+1-i}=\frac{a_{d+1-i}}{a_{d-i}}
 =\frac{a_{i-1}^r}{a_i^r}=t_i^{-r}.
\]
The prescribed norm equality is equivalent to \(t_i^{r+1}=1\), so \(t_i^{-r}=t_i\).
The adjacent ratios at the distinct positions \(i\) and \(d+1-i\) therefore coincide.
Consequently,
\[
 P(X)=(1-t_iX)F(X)
\]
has zero coefficients at both positions.
Its constant and leading coefficients are nonzero, its degree is \(d+1\), and it has at most \(d\) terms.
The monic multiplier is \(X-t_i^{-1}\), differing only by a nonzero scalar.
At \(2i=d+1\), reflection gives no second position; the norm-one condition already follows from \(t_i=t_i^{-r}\) and does not imply a repeated ratio.
The argument uses a given norm equality and does not assert that such a noncentral position always exists.
\end{solution}

\solutionheading{v97:ex-projection-sparse}
\begin{solution}
Write \(V=k[X]/(f)\) and \(\beta=X\bmod f\). Since \(f(0)\ne0\), every \(\beta^j\) is nonzero.
For \(d=2\), the space \(V\) has \(q+1\) one-dimensional subspaces.
Two of the \(q+2\) vectors \(1,\beta,\ldots,\beta^{q+1}\) must be proportional, giving the required binomial multiple.

For \(d\ge3\), put \(W=\operatorname{span}_k\{1,\beta,\ldots,\beta^{d-4}\}\), taking \(W=0\) when \(d=3\).
The power basis gives \(\dim W=d-3\), so \(V/W\) is three-dimensional.
Consider the \(q+2\) images of \(\beta^{d-3},\ldots,\beta^{d+q-2}\).
A zero image or a proportional pair already gives a relation with at most two terms.
Otherwise their one-dimensional spans form a set \(S\) of \(q+2\) distinct points.

Suppose no three points of \(S\) lie in a common two-dimensional subspace.
There are \(q+1\) two-dimensional subspaces through any given point of \(S\), and the other \(q+1\) points must occupy them separately.
Thus every two-dimensional subspace meeting \(S\) contains exactly two of its points.
There are \(q^2+q+1>q+2\) one-dimensional subspaces in total; choose one, \(P\), outside \(S\).
The two-dimensional subspaces through \(P\) partition \(S\) into groups of size zero or two, making \(|S|=q+2\) even, a contradiction.
Hence three images satisfy a nontrivial relation.
Lift it to \(V\). Its sum lies in \(W\); subtract its expression using at most \(d-3\) power-basis terms.
The result is a nonzero \(f\)-multiple of degree at most \(d+q-2\) with at most \(d\) terms.
The quotient used here is a quotient vector space: \(W\) is not assumed to be an ideal.

Finally, \(\alpha^q=\alpha^{-1}\), so the prescribed \(f\) lies in \(k[X]\) and has three distinct roots.
For \(0\le j\le q\), its root-evaluation columns are \((1,\alpha^j,\alpha^{-j})^{\mathsf T}\).
Scale the \(j\)th column by \(\alpha^j\) and reorder the rows to obtain \((1,\alpha^j,\alpha^{2j})^{\mathsf T}\).
Any three are Vandermonde columns at distinct nodes and hence independent; smaller column sets are independent as well.
There is therefore no nonzero \(f\)-multiple of degree at most \(q\) with at most three terms.
Since \(f\mid X^{q+1}-1\), equality \(L_3(f)=q+1\) follows.
\end{solution}

\solutionheading{v114:ex-quadratic-count}
\begin{solution}
A nonzero element \(y\) of a quadratic extension is normal if and only if \(\Tr(y)\ne0\) and \(y\notin k\). Indeed, dependence of \(y,y^q\) gives \(y^q=cy\) with \(c\in k^\times\) and \(c^2=1\); thus \(c=1\) or \(-1\), giving \(y\in k\) or \(\Tr(y)=0\), respectively. The converse is immediate; in characteristic two the two cases coincide.
Here \(y=(\alpha+t)^{-1}\) always has full degree, and
\[
 \Tr(y)=\frac{\Tr(\alpha)+2t}{\Nm(\alpha+t)}.
\]
The denominator is nonzero. In odd characteristic the trace vanishes exactly at \(t=-\Tr(\alpha)/2\). In characteristic two, \(\Tr(\alpha)=\alpha+\alpha^q\ne0\), so there are no exceptions.
\end{solution}

\solutionheading{v114:ex-cubic-classification}
\begin{solution}
Multiplication and reduction by \(\theta^3=d\) give
\[
\begin{aligned}
 &(w+u\theta+v\theta^2)
 \bigl((w^2-uvd)+(v^2d-uw)\theta+(u^2-vw)\theta^2\bigr)\\
 &\hspace{25mm}=w^3+u^3d+v^3d^2-3uvdw.
\end{aligned}
\]
The right-hand side is the norm: the three Frobenius conjugates are obtained by \(\theta\mapsto\zeta^i\theta\), and their product has this expansion.
Since \(\alpha\notin k\), the denominator is nonzero for every \(w\in k\). Irreducibility implies that \(d\ne0\) is not a cube.
If \(u\ne0,v=0\), the numerators are \(w^2,-uw,u^2\); if \(u=0,v\ne0\), they are \(w^2,v^2d,-vw\). In either case their root sets have union \(\{0\}\).

If \(uv\ne0\), put \(r_1=v^2d/u\) and \(r_2=u^2/v\). Equality \(r_1=r_2\) would give \(d=(u/v)^3\). Either \(r_1^2=uvd\) or \(r_2^2=uvd\) gives the same contradiction to irreducibility. Thus the three root sets are pairwise disjoint, and
\[
\begin{array}{c|c}
 \text{Condition}&b_{q,3}(\alpha)\\\hline
 uv=0&1\\
 uv\ne0,\ q\text{ odd},\ uvd\text{ a square}&4\\
 uv\ne0,\ q\text{ odd},\ uvd\text{ a nonsquare}&2\\
 uv\ne0,\ q\text{ even}&3
\end{array}
\]
In even characteristic the square map is bijective. The quadratic numerator has a repeated root but only one distinct root, giving count three.
In odd characteristic, \(u=1,v=d^{-1}\) makes \(uvd=1\), attaining four. In even characteristic any nonzero \(u,v\) attain three.
All these \(\alpha\) lie outside the base field and therefore have full degree, since three is prime. This proves the stated maxima; the odd-characteristic value two is obtained by choosing \(uvd\) to be any nonsquare.
\end{solution}

\solutionheading{v114:ex-quartic-bound}
\begin{solution}
Multiply the four numerators by \(T+\alpha\) and compare coefficients of \(1,\theta,\theta^2,\theta^3\) to obtain~\eqref{v114:eq-quartic-numerators}.
The discriminant relation and degree drops are treated in the proof of Proposition~\ref{v115:prop-quartic-bound}; the root sets attaining equality are listed in Example~\ref{v115:ex-quartic-seven}.
\end{solution}

\solutionheading{v114:ex-reciprocal-average}
\begin{solution}
An element in a subfield of degree \(d<n\) has conjugate-span dimension at most \(d\), so it is not normal. Thus every normal element belongs to \(\mathcal G_{q,n}\).
For fixed \(t\in k\), the identity \(\alpha=((\alpha+t)^{-1})^{-1}-t\) shows that translating and taking a reciprocal preserves field generation. The inverse map \(y\mapsto y^{-1}-t\) also preserves it, proving bijectivity on this set.
When \(n\mid q-1\), equation~\eqref{v114:eq-eigen-normal} gives exactly \((q-1)^n\) normal elements.
Thus each fixed \(t\) has exactly \(G_{q,n}-(q-1)^n\) generators whose reciprocal translate is not normal. Counting the same pairs first by \(\alpha\) and then by \(t\) gives
\[
 \sum_{\alpha\in\mathcal G_{q,n}}b_{q,n}(\alpha)
 =q\bigl(G_{q,n}-(q-1)^n\bigr).
\]
Divide by \(G_{q,n}\) to obtain the formula. The full-degree count follows from M\"obius inversion of subfield counts, independently of the normal-element count.
For \((q,n)=(7,3)\), we have \(G_{7,3}=343-7=336\), and the average is \(7(1-216/336)=5/2\); the two families in Example~\ref{v114:ex-cubic-comparison} have one and four exceptions in the same extension.
For \((q,n)=(5,4)\), we have \(G_{5,4}=625-25=600\), and the average is \(5(1-256/600)=43/15\), whereas the fixed family in Example~\ref{v114:ex-quartic-failure} has five exceptions.
An average count and extremal behavior can coexist; they answer different questions.
\end{solution}

\solutionheading{v115:ex-reciprocal-eight}
\begin{solution}
The cubic has no root in \(\F_2\), so it is irreducible. The test for \(Y+1\) is \(1+\sigma+\sigma^2=\Tr_{E/\F_2}\); the test for the quadratic factor is \(1+\sigma\), whose kernel is the base field.
Every reciprocal in the family has degree three, so the second test is automatically nonzero. Also,
\[
D(T)=T^3+T+1,\qquad
\Tr((\alpha+t)^{-1})=\frac{t^2+1}{t^3+t+1}.
\]
Thus \(t=1\) is the unique exception. Direct calculation gives
\[
\alpha^{-1}=\alpha^2+1,\qquad
(\alpha+1)^{-1}=\alpha^2+\alpha.
\]
Since \(\Tr(\alpha)=\Tr(\alpha^2)=0\) and \(\Tr(1)=1\), these reciprocals have traces one and zero, respectively, confirming the test.
\end{solution}

\Needspace{9\baselineskip}
\section*{Solutions to Chapter~\ref{chap:root-fibers}}
\pdfbookmark[1]{Solutions to Chapter \getrefnumber{chap:root-fibers}}{solutions.chap:root-fibers}

\solutionheading{v70:ex-octic-trace-zero}
\begin{solution}
First, \(f(0)=f(1)=1\), so \(f\) has no linear factor. Each row below gives the identity \(f=PQ+R\):
\[
\begin{array}{c|c|c}
P&Q&R\\\hline
X^2+X+1&X^6+X^4+X^3&X+1\\
X^3+X+1&X^5+X^4+X^3&X+1\\
X^3+X^2+1&X^5+X^2+X&X^2+1\\
X^4+X+1&X^4+X^3+X&X^2+1\\
X^4+X^3+1&X^4+1&X\\
X^4+X^3+X^2+X+1&X^4+X^2&X^2+X+1
\end{array}
\]
The left column lists all monic irreducible polynomials of degrees two, three, and four: in degrees two and three it suffices to exclude base-field roots;
the monic quartics without such roots are the three listed and \(X^4+X^2+1=(X^2+X+1)^2\),
and none of the three listed quartics is divisible by \(X^2+X+1\). Since every remainder is nonzero, \(f\) has no irreducible factor of degree at most four.
A reducible polynomial of degree eight must have such a factor, so \(f\) is irreducible.

The next-to-leading coefficient gives \(\Tr_{\F_{2^8}/\F_2}(\alpha)=1\). Also,
\[
f'(X)=X^6+X^2+1,\qquad f'(0)=f'(1)=1,
\]
so Proposition~\ref{v70:prop-mobius-trace-zero} excludes every fractional linear transformation.

For \(\beta=\alpha^2+\alpha\), the trace is \(1^2+1=0\).
If \(\beta\) did not generate \(\F_{2^8}\), the relation \(\alpha^2+\alpha+\beta=0\) would force
\([\F_{2^8}:\F_2(\beta)]=2\), \(\F_2(\beta)=\F_{16}\), and \(\alpha^{16}=\alpha+1\).
Trace transitivity would then give
\[
\Tr_{\F_{2^8}/\F_2}(\alpha)
=\Tr_{\F_{16}/\F_2}(\alpha+\alpha^{16})
=\Tr_{\F_{16}/\F_2}(1)=0,
\]
contrary to the trace being one. This is the case \(4\mid n\) of Proposition~\ref{v70:prop-trace-zero-generator}.

Finally, the binomial expansion gives
\[
f(X+1)=X^8+X^7+X^6+X^5+X^4+X+1,
\]
and substitution into \(g\) and multiplication give the same polynomial:
\[
\begin{aligned}
g(X^2+X)&=X^{16}+X^{10}+X^9+X^8+X^3+X^2+1,\\
f(X)f(X+1)&=X^{16}+X^{10}+X^9+X^8+X^3+X^2+1.
\end{aligned}
\]
Thus \(g(\beta)=0\). Since \(\beta\) has degree eight and \(g\) is monic of degree eight, \(m_\beta=g\), with no separate irreducibility test for \(g\) required.
\end{solution}

\solutionheading{ex:supp-binomial}
\begin{solution}
\begin{enumerate}
\item If \(n=1\), every \(X+a\) is irreducible, so \(\mathcal A_q(1)=\F_q\). If \(n\ge2\), then \(X^n\) is reducible, so \(0\notin\mathcal A_q(n)\).
Assume henceforth that \(n\ge2\) and \(a\ne0\), and put \(P=q-1\) and \(e=\ord(-a)\).
Apply Theorem~\ref{thm:binomial-full} to \(X^n-(-a)\). Its first condition is equivalent to
\[
v_\ell(e)=v_\ell(P)>0\qquad(\ell\mid n,\ \ell\text{ prime}).
\]
Thus each prime divisor of \(n\) divides \(P\), and \(4\mid n\) requires \(q\equiv1\pmod4\). This proves the necessity of \eqref{v77:eq-binomial-existence}.
Conversely, if these conditions hold, take a primitive element \(g\in\F_q^\times\) and set \(a=-g\).
Then \(e=P\), so all conditions of the criterion hold, proving \(\mathcal A_q(n)\ne\varnothing\).
Only the prime divisors of \(n\) are restricted here, not their exponents in \(n\); consequently \(n\mid q-1\) is not required.
If the characteristic \(p\) divides \(n\), the existence condition already fails. Directly, taking \(b^p=-a\) gives \(X^n+a=(X^{n/p}-b)^p\), again proving reducibility.

\item Suppose now that \eqref{v77:eq-binomial-existence} holds, and write uniquely
\[
-a=g^k,\qquad 0\le k\le P-1.
\]
Then \(e=P/\gcd(P,k)\). For each prime \(\ell\mid n\), we already have \(\ell\mid P\), and hence
\[
\begin{aligned}
v_\ell(e)=v_\ell(P)>0
&\iff \ell\nmid k\\
&\iff (g^k)^{P/\ell}\ne1.
\end{aligned}
\]
The last equivalence uses the fact that \(g\) has order \(P\).
These conditions hold for every prime divisor of \(n\) precisely when \(\gcd(k,n)=1\), proving both \eqref{v77:eq-binomial-constant-exponents} and \eqref{v77:eq-binomial-coefficient-tests}.
In the cyclic group \(\F_q^\times\), the element \(g^k\) is an \(\ell\)th power precisely when \(\ell\mid k\).
Thus the power test also says that \(-a\) is not an \(\ell\)th power for any prime \(\ell\mid n\).
The order condition can equivalently be written as \(\rad(n)\mid e\) and \(\gcd(n,P/e)=1\), as in the criterion in the text.
Changing to a primitive element \(g'=g^u\) gives \(\gcd(u,P)=1\).
Since \(\rad(n)\mid P\), this only permutes the exponents satisfying the coprimality condition and leaves the coefficient set unchanged.

\item Put \(R=\rad(n)\). Under the existence conditions, \(R\mid P\), and \(\gcd(k,n)=1\) is equivalent to \(\gcd(k,R)=1\).
The integers \(0,1,\ldots,P-1\) contain \(P/R\) complete systems of residues modulo \(R\), each contributing \(\varphi(R)\) admissible exponents.
Since \(k\mapsto-g^k\) is a bijection on these \(P\) exponents,
\[
|\mathcal A_q(n)|=\frac{P}{R}\varphi(R)
=(q-1)\prod_{\ell\mid n}\left(1-\frac1\ell\right)
=(q-1)\frac{\varphi(n)}n.
\]
Combining this with the case \(n=1\) and the nonexistence cases from part~(1) proves \eqref{v77:eq-binomial-count}.
Formula~\eqref{v77:eq-binomial-constant-exponents} also shows that \(\mathcal A_q(m)=\mathcal A_q(n)\) whenever \(m,n\ge2\) satisfy the existence conditions and have the same squarefree kernel.
Existence itself still includes the restriction modulo four: for example, \(\rad(2)=\rad(4)\) but \(\mathcal A_3(2)\ne\mathcal A_3(4)\).

\item Over \(\F_5\), we have \(5\equiv1\pmod4\), and the required condition is exactly that \(-a\) be a nonsquare. Thus \(\mathcal A_5(8)=\{2,3\}\).
Over \(\F_7\), the nonzero cubes are \(1,6\), giving \(\mathcal A_7(9)=\{2,3,4,5\}\).
Over \(\F_3\), degree four violates the restriction modulo four, so \(\mathcal A_3(4)=\varnothing\).
Consequently \(X^8-2\) is irreducible over \(\F_5\), with constant term \(a=3\); \(X^9-2\) is irreducible over \(\F_7\), with \(a=5\); and \(X^4+1\) is reducible over \(\F_3\), as factored in Example~\ref{v77:ex-binomial-constant-families}.
The same prime-by-prime conditions give \(\mathcal A_5(2^r)=\{2,3\}\) and \(\mathcal A_7(3^r)=\{2,3,4,5\}\) for every \(r\ge1\).

For \(q=2\) and \(n\ge2\), the polynomial \(X^n\) is reducible, whereas \(X^n+1\) has the root \(1\). Thus \(\mathcal A_2(n)=\varnothing\).
This does not extend to characteristic two in general: in \(\F_4=\F_2(\omega)\), where \(\omega^2+\omega+1=0\), the only nonzero cube is \(1\), and hence \(\mathcal A_4(3)=\{\omega,\omega^2\}\).

The nonsquares in \(\F_7\) are \(3,5,6\), so \(\mathcal A_7(2)=\{1,2,4\}\).
All three coefficients are themselves squares, but their negatives are nonsquares; the sign in the criterion cannot be omitted.
Finally, take \(n=6\) and \(a=5\) over \(\F_7\). The element \(-a=2\) is not a sixth power, since every nonzero sixth power equals \(1\).
Nevertheless, \(2=3^2\) is a square and
\[
X^6+5=X^6-2=(X^3-3)(X^3+3).
\]
Excluding only \(n\)th powers therefore cannot replace the separate test for each prime divisor of \(n\).
\end{enumerate}
\end{solution}

\solutionheading{ex:supp-fiber-binary}
\begin{solution}
Take \(P=X^2+X+1\). It is irreducible and \(c_1/c_0=1\) has absolute trace one, so \(X^2P(X+X^{-1})=X^4+X^3+X^2+X+1\) is irreducible.
Its fiber polynomial is an irreducible quadratic over \(\F_4\); the return map has period two, producing one orbit of length four.
After extending the base field to \(\F_4\), the four-cycle splits under \(\sigma^2\) into two two-cycles, so there are exactly two distinct irreducible quadratic factors.
\end{solution}

\solutionheading{ex:supp-affine-quartic}
\begin{solution}
Take \(u=v=b=1\). Then \(L=(X^2+X)^2+(X^2+X)=X^4+X\).
Both quadratics to be checked equal \(X^2+X+1\), so \(L+1\) is irreducible.
If a degree-eight \(2\)-linearized polynomial plus a constant were irreducible, its derivative would be nonzero, its root set would be a three-dimensional \(\F_2\)-affine space, and Frobenius \(x\mapsto x^2\) would act as an eight-cycle.
Choose affine coordinates and write this action as \(T(v)=Av+b\) with \(A\in\operatorname{GL}_3(\F_2)\) and \(b\in\F_2^3\), and define the four-dimensional augmented matrix
\[
U=\begin{pmatrix}A&b\\0&1\end{pmatrix}.
\]
This matrix representation is injective and preserves composition, so \(\ord(U)=\ord(T)=8\). Thus \(N=U-I\) satisfies \(N^8=U^8-I=0\) and is nilpotent.
A nilpotent operator on a four-dimensional space satisfies \(N^4=0\); hence \(U^4=I+N^4=I\), contrary to its having order eight.
\end{solution}

\solutionheading{ex:supp-composite-as}
\begin{solution}
The absolute trace in the criterion is the identity on \(\F_3\), and \(3b-c_2=-c_2\ne0\). Hence every \(b\) satisfies the condition, and the composite has degree nine.
\end{solution}

\solutionheading{ex:wu-quartic-hypothesis}
\begin{solution}
We have \(f(0)=f(1)=1\), so there is no linear factor.
Modulo the unique irreducible quadratic \(X^2+X+1\), the relation \(X^3=1,X^4=X\) gives \(f\equiv X\ne0\), so \(f\) is irreducible.
If the root set were \(\beta+\{0,u,v,u+v\}\), the sum of the four roots would be \(4\beta+u+v+(u+v)=0\).
But the sum of the roots of \(f\), given by its cubic coefficient, is one, a contradiction.
Thus \(L\) must be \(2\)-linearized; merely requiring \(L\) to be monic of degree four is insufficient for this criterion.
\end{solution}

\solutionheading{ex:wu2-fiber-multiplicity}
\begin{solution}
The polynomial \(P\) has no root in \(k\) and is irreducible; the base orbit of \(\alpha\) is \(\{\alpha,\alpha^2\}\), of length two.
Thus the first Frobenius iterate preserving the starting fiber is \(\sigma^2\), rather than \(\sigma\).

In the first part, \(X^2-\alpha=(X-\alpha^2)^2\), so the fiber has only the root \(\alpha^2\), with multiplicity two.
The return map fixes this root, so the total orbit length is \(2\cdot1=2\) and the multiplicity remains two; hence
\[
P(X^2)=(X^2+X+1)^2.
\]
For the second fiber, if \(\beta^2+\beta=\alpha\), then
\[
\beta^4+\beta=\alpha^2+\alpha=1,
\qquad \sigma^2(\beta)=\beta+1.
\]
The two distinct fiber roots are \(\beta,\beta+1\), and the return map interchanges them, so the total orbit length is \(2\cdot2=4\).
Since the fiber polynomial has derivative one, the multiplicity is one, and \(P(X^2+X)=X^4+X+1\) is irreducible.

In the second part, \(L'=1\), so \(S_\alpha\) has exactly eight distinct elements, and every \(\beta\in S_\alpha\) satisfies
\(S_\alpha=\beta+\F_8\). Since \(\alpha^4=\alpha\), the map \(\sigma^2\) preserves this fiber.
For \(x\in S_\alpha\), we have
\[
x^{64}+x=(x^8+x)^8+(x^8+x)=\alpha^8+\alpha
=\alpha^2+\alpha=1.
\]
Thus \(R^3(x)=x+1\) and \(R^6(x)=x\). The return-cycle lengths therefore divide six; lengths one and three are excluded because \(R^3\) has no fixed points, leaving only two and six.

If \(R^2(x)=x\), then \(x^{16}=x\). Squaring \(x^8+x=\alpha\) gives \(x^{16}+x^2=\alpha^2\), so
\[
x^2+x=\alpha^2.
\]
Conversely, this quadratic equation implies
\[
x^4+x=\alpha^4+\alpha^2=1,\qquad
x^8+x=1+\alpha^2=\alpha,
\]
and also \(x^{16}+x=(x^8+x)^2+(x^2+x)=0\).
Thus the fixed points of \(R^2\) are exactly the two distinct roots of \(X^2+X+\alpha^2\).
They form one two-cycle; the other six points must form one six-cycle.

In the third part, Lemma~\ref{lem:fiber-return} lifts these two return cycles to orbits of lengths
\[
2\cdot2=4,\qquad 2\cdot6=12
\]
under \(\sigma\). Thus \(P(L)\) has exactly one irreducible quartic factor and one irreducible factor of degree 12 over \(k\).
Also \((P(L))'=P'(L)L'=1\), so both factors have multiplicity one.
The fiber polynomial of the two-cycle is \(X^2+X+\alpha^2\); multiplying by its coefficient conjugate gives
\[
(X^2+X+\alpha^2)(X^2+X+\alpha)=X^4+X+1.
\]
Dividing \(P(L)\) by this quartic gives the listed degree-12 factor. Its irreducibility is already determined by the six-step return period and two-step base period, without further trial division.

The return period computes a root's degree over \(\F_4\); its degree over the base field \(\F_2\) must also include the base-orbit length two as a factor.
The return map acts on the set of distinct roots and cannot detect multiplicities; the inseparable fiber in the first part shows that these must be retained separately by the fiber polynomial and carried through conjugate descent.
\end{solution}

\solutionheading{ex:root-spectrum-base-extension}
\begin{solution}
By Theorem~\ref{thm:factor-base-extension} on extension of the base field, each irreducible factor of degree \(d\) splits over \(\F_{q^r}\) into irreducible factors of degree
\(d/\gcd(d,r)\). Taking the set of these degrees over all original factors gives the first assertion; multiplicities do not affect this set.

If \(d_0\mid d\), then for every prime \(\ell\),
\[
v_\ell\!\left(\frac{d_0}{\gcd(d_0,r)}\right)
=\max\{v_\ell(d_0)-v_\ell(r),0\}
\le\max\{v_\ell(d)-v_\ell(r),0\}
=v_\ell\!\left(\frac d{\gcd(d,r)}\right).
\]
Thus \(d_0/\gcd(d_0,r)\mid d/\gcd(d,r)\).
If \(d_0\) is the least member of the original root-field spectrum, its transformed degree lies in the new spectrum and divides every degree there, proving the second assertion.

For \(f=P_2P_3\), the original root-field spectrum is \(\{2,3\}\) and has no least member under divisibility;
although the two root fields intersect in \(\F_q\), the polynomial \(f\) has no base-field root, so this intersection is not a field generated by a root of \(f\).
After extending the base field to \(\F_{q^2}\), the spectrum becomes \(\{1,3\}\), and the least root field is the new base field.
Thus the greatest common divisor must itself belong to the root-field spectrum for it to determine the least field generated by a root.
\end{solution}

\solutionheading{v79:ex-rational-generation}
\begin{solution}
(1) Since \(\deg B\le\deg R<\ell\le n\), a nonzero \(B\) cannot vanish at \(\alpha\).
Lemma~\ref{v70:lem-low-degree-generation} gives
\[
[E:k(R(\alpha))]\mid n,\qquad
[E:k(R(\alpha))]\le\deg R<\ell.
\]
The only divisor of \(n\) below \(\ell\) is one, so \(k(R(\alpha))=E\).

(2) Write \(A(X)=\sum_{i=0}^r a_iX^i\) and \(B(X)=\sum_{i=0}^r b_iX^i\), with coefficients in \(k\).
The equation \(A(\alpha)-\beta B(\alpha)=0\) gives at most \(n\) homogeneous linear equations in \(2r+2>n\) unknowns, so it has a nonzero solution \((A,B)\).
If \(B=0\), then \(A(\alpha)=0\); since \(\deg A\le r<n\), this forces \(A=0\), a contradiction. Hence \(B\ne0\), and \(\deg B<n\) also gives \(B(\alpha)\ne0\).
Cancel common factors of \(A,B\). Each such factor has degree below \(n\) and is nonzero at \(\alpha\), so cancellation preserves the value \(\beta\) and does not increase either degree.
The resulting rational function is nonconstant, since a constant function in \(k(X)\) would take a value in \(k\), contrary to \(\beta\notin k\).
This is also the symmetric degree-bound case of Proposition~\ref{v70:prop-rational-reconstruction}.

(3) A monic irreducible binomial of degree \(n\) splits in \(E\), so choose one of its roots \(\beta\in E\). It has degree \(n\) over \(k\), and part~(2) represents it as \(R(\alpha)\) with \(\deg R\le r\).
Thus existence of a binomial target gives accessibility under a degree bound depending on \(n\). For \(n=2,3\), this bound is one, and for \(n=4,5\), it is two; for \(n\ge6\), the bound exceeds two and does not settle accessibility by quadratic rational functions.
\end{solution}

\solutionheading{v79:ex-binomial-defining-equation}
\begin{solution}
(1) Put \(\alpha_i=\alpha^{q^i}\). Since the coefficients of \(A,B\) lie in \(k\), all \(B(\alpha_i)\) are nonzero and \(R(\alpha_i)=\beta^{q^i}\).
The root-product formula for the resultant gives
\[
\begin{aligned}
\operatorname{Res}_X(f,BY-A)
&=\prod_{i=0}^{n-1}\bigl(B(\alpha_i)Y-A(\alpha_i)\bigr)\\
&=\left(\prod_{i=0}^{n-1}B(\alpha_i)\right)
  \prod_{i=0}^{n-1}(Y-\beta^{q^i})\\
&=\Nm_{E/k}\bigl(B(\alpha)\bigr)\,m_\beta(Y)^{n/d}.
\end{aligned}
\]
The norm factor is nonzero. Dividing by it gives a monic polynomial, which is the new degree-\(n\) minimal polynomial exactly when \(d=n\).
If \(d<n\), each of the \(d\) points of the image orbit has \(n/d\) preimages in the source orbit. This repeated occurrence accounts for the exponent \(n/d\), although the source polynomial \(f\) is separable.

(2) The discriminant of \(f\) is \(1-4=2\), a nonsquare in \(\F_5\), so \(f\) is irreducible.
From \(\alpha^2=-\alpha-1\), one obtains
\[
\beta^2=(\alpha+3)^2=\alpha^2+\alpha+4=3.
\]
Also, \(\alpha=\beta-3\), so \(k(\beta)=k(\alpha)\), and \(m_\beta(Y)=Y^2-3\). Equivalently,
\[
\operatorname{Res}_X(f,Y-X-3)=f(Y-3)=Y^2-3.
\]
The original defining polynomial has three nonzero terms, whereas the new one has two.
Since \((\alpha-1)f(\alpha)=\alpha^3-1=0\) and \(f(1)\ne0\), we have \(\ord(\alpha)=3\).
The element \(3\in\F_5^\times\) has order four, so \(\beta^2=3\) gives \(\ord(\beta)=8\).
Taking \(\beta\) as a primitive eighth root of unity, its conjugate roots are \(\beta,\beta^5\), with exponent orbit
\[
C_1=\{1,5\}=1+4\Z/8\Z.
\]
This is a complete equal-difference block and gives the binomial root polynomial by Proposition~\ref{prop:ed-binomial-correspondence}.
The change preserves the field degree but changes the multiplicative order and the exponent modulus.
\end{solution}

\solutionheading{v109:ex-binary-trace}
Take a root $\alpha\in\F_{2^n}$ of $f$. The Artin--Schreier criterion says that $Y^2+Y=\alpha$ has no solution in this field exactly when $\Tr(\alpha)=1$, which is the next-to-leading coefficient.
In the nonsplit case, a root $\beta$ has degree $2n$ over $\F_2$, so the composition is irreducible. In the split case both roots lie in $\F_{2^n}$ and generate it, since $\alpha=\beta^2+\beta$. All $2n$ roots are distinct, so there are exactly two degree-$n$ factors.

Put $B_1(e)=eA_1(e)$, the number of degree-$e$ elements with trace one to $\F_2$. On the subfield $\F_{2^e}\subseteq\F_{2^n}$, trace is multiplied by $n/e$. Hence
\[
2^{n-1}=\sum_{\substack{e\mid n\\n/e\text{ odd}}}B_1(e).
\]
Fix the 2-adic valuation of $n$ and invert over its odd divisors to obtain the formula for $A_1$.
Also $A_0(n)+A_1(n)=I_2(n)=n^{-1}\sum_{d\mid n}\mu(d)2^{n/d}$.
For odd $n$ the formulas cancel directly. For even $n$, only even divisors $d=2e$ remain, with $e$ odd. Using $\mu(2e)=-\mu(e)$ gives
\[
2A_1(n)-I_2(n)=\frac1n\sum_{\substack{e\mid n/2\\e\text{ odd}}}\mu(e)2^{n/(2e)}=A_1(n/2).
\]
The pairs $(A_1,A_0)$ are respectively $(1,1),(1,0),(1,1),(2,1),(5,4)$.

\solutionheading{v115:ex-return-observations}
\begin{solution}
The fixed-point counts are \((2,2,2,6)\) and \((2,6,2,6)\), respectively. Inversion gives
\[
a_1=r_1,\qquad a_2=(r_2-r_1)/2,\qquad
a_3=(r_3-r_1)/3,\qquad a_4=(r_4-r_2)/4.
\]
Thus the first return has two fixed points and one four-cycle; the second has two fixed points and two two-cycles.
The full factor degrees are \((n,n,4n)\) and \((n,n,2n,2n)\). Fiber size and \(r_1\) coincide, whereas \(r_2\) already distinguishes them.
\end{solution}

\Needspace{9\baselineskip}
\section*{Solutions to Chapter~\ref{chap:root-symmetry}}
\pdfbookmark[1]{Solutions to Chapter \getrefnumber{chap:root-symmetry}}{solutions.chap:root-symmetry}

\solutionheading{ffs:ex-weighted-two}
\begin{solution}
The two roots have different multiplicities, so every stabilizing transformation must fix zero and one separately; hence \(\mathcal A_f\) is trivial and \(V_f=\{0\},M_f=\{1\}\).
In \(\F_3\) we have \(f=X^3+X^2+X\), so
\[
f(X+B)-f(X)=2BX+B^3+B^2+B.
\]
The monic greatest common divisor of the two nonzero coefficients is \(B\), so \(\Psi_f=B\).
Ignoring multiplicities, \(X\mapsto1-X\) also preserves \(\{0,1\}\), so the affine stabilizer of the root set alone has two elements.
\end{solution}

\solutionheading{ffs:ex-power-four}
\begin{solution}
Expansion gives \(f=X^8+X^4+1\), so \(e=4\), and
\[
f(X+B)-f(X)=B^8+B^4=(B^2+B)^4.
\]
Thus \(\Psi_f=B^8+B^4\), \(L_{V_f}=B^2+B\), and \(V_f=\F_2\).
The positive exponents have greatest common divisor four, whose part prime to two is one, so \(M_f=\{1\}\).
There are only two distinct roots; the identity and \(X\mapsto X+1\) already give two distinct permutations, and an affine transformation is determined by the images of two points, so the full group consists exactly of these two translations.
\end{solution}

\solutionheading{ffs:ex-six-slopes}
\begin{solution}
The derivative of \(Y^6-1\) is \(Y^5\), nonzero at all six roots, so \(f\) has six distinct roots and \(e=1\).
The number \(|V_f|\) is a power of five dividing six, so \(V_f=\{0\}\), \(L_V=X\), and \(\Psi_f=B\).
Expansion gives
\[
f=X^6+3X^5+3X+3,
\]
and the nonzero positive exponents have greatest common divisor one, so \(M_f=\{1\}\).
Every transformation \(x\mapsto2+a(x-2)\) with \(a\in\mu_6\) preserves the root set, so the slope image contains at least six elements.
When the translation group is trivial, the full affine stabilizer is a cyclic dilation group with a common center; its special orbit is a singleton, and every other orbit has size \(t\).
If this point belonged to the root set, then \(t\mid5\), contradicting \(6\mid t\); thus it does not, and \(t\mid6\).
Therefore \(t=6,c=2,\delta=2,U=X-2,w=0,H(Y)=Y-1\).
A base-field center exists, and \(|\mu_6\cap\F_5^\times|=\gcd(6,4)=2\), so the base-field affine stabilizer has order two.

Here \(d=6,Q=m=1,t=6\), so \(r=0\) and \(\operatorname{sp}(f)=1\).
The powers defining \(E_\Gamma\) are \(1\) and \((X-2)^6=f+1\); the expansion above gives
\[
E_\Gamma=\{0,1,5,6\},\qquad E_0=\{0,1,2,3,4,5,6\}.
\]
Thus
\[
|\operatorname{Supp}_X(f)|=|E_\Gamma|=4<7=|E_0|
=\frac{\deg f}{\operatorname{sp}(f)}+1.
\]
The full affine symmetry therefore gives a sharper support bound than the support period alone.
\end{solution}

\solutionheading{ffs:ex-two-additive-types}
\begin{solution}
Neither quadratic has a root in \(\F_3\), so both are irreducible.
The roots of \(P_1\) have trace zero in \(\F_9/\F_3\), while those of \(P_2\) have trace \(-1=2\).
For \(L=X^3-X\), these traces are exactly the return-translation constants.
Thus \(P_1(L)\) splits into three distinct quadratic factors, whereas \(P_2(L)\) is an irreducible sextic.
Direct multiplication verifies
\[
P_1(L)=(X^2+1)(X^2+X+2)(X^2+2X+2).
\]
None of the three quadratics has a root in \(\F_3\), so this is the complete irreducible factorization.
\end{solution}

\solutionheading{ffs:ex-four-cycle}
\begin{solution}
Successive application gives
\[
(0,0)\longmapsto(0,1)\longmapsto(1,0)
\longmapsto(1,1)\longmapsto(0,0).
\]
Thus \(N_1=N_2=0,N_4=4\), and \eqref{ffs:eq-return-mobius} gives \(C_1=C_2=0,C_4=1\).
One may also calculate \(A^2=I\), \(b_2=(1,0)\ne0\), and \(b_4=0\) directly and substitute into the fixed-point formula.

Take \(c^4+c=1\) and put \(\omega=c^2+c\). Then \(\omega^2+\omega=1\), so \(1,\omega\) is a basis of \(\F_4/\F_2\).
The root fiber is \(c+\F_4\), and squaring in the coordinates \(c+v\leftrightarrow v\) is \(v\mapsto v^2+\omega\).
By \(\omega^2=\omega+1\), its linear part in the basis \(1,\omega\) is exactly the stated \(A\), and its translation part is \(b\).
Thus this four-cycle is the Frobenius root orbit of \(X^4+X+1\).
\end{solution}

\solutionheading{ffs:ex-degree-restrictions}
\begin{solution}
The number \(|V_f|\) can only be one, two, or four, and must divide \(\deg f\).
If it equals four, we must also have \(v_2(\deg f)=2\).
Thus the allowed sizes are one or two in degree six; one, two, or four in degree 12; and again only one or two in degree 24.
The last case satisfies \(4\mid24\), but \(24/4=6\) is even, so the exact condition for a four-element translation block excludes it.
\end{solution}

\solutionheading{ffs:ex-rational-centers}
\begin{solution}
If \(c_0\in C\cap k\), then
\[
C\cap k=(c_0+V)\cap k=c_0+(V\cap k).
\]
Thus its size is \(|V\cap k|\).
The condition \(U\in k[X]\) ensures that Frobenius permutes its roots, but does not imply that any individual root is fixed.
In the quartic example, \(L_V(x)=x^4+x=0\) for every \(x\in\F_2\), whereas \(\delta=1\).
Thus \(\delta\notin L_V(\F_2)\), and the center equation has no solution in the base field.
\end{solution}

\solutionheading{ffs:ex-product-new-symmetry}
\begin{solution}
Take \(\alpha^3-\alpha=1\). The two root sets are \(\alpha+\F_3\) and \(-\alpha+\F_3\), respectively, with \(\alpha\in\F_{27}\).
Their full affine stabilizers are, respectively,
\begin{align*}
\mathcal A_{f_1}&=\{X\mapsto\alpha+a(X-\alpha)+v:a\in\{1,-1\},\ v\in\F_3\},\\
\mathcal A_{f_2}&=\{X\mapsto-\alpha+a(X+\alpha)+v:a\in\{1,-1\},\ v\in\F_3\}.
\end{align*}
Each group has exactly six elements, since its translation group is \(\F_3\) and its allowed slopes must lie in \(\F_3^\times\).
If a transformation of slope \(-1\) preserved both cosets, its affine transformation on the quotient coordinate \(L\) would have to fix \(1,-1\) separately, which is impossible.
Thus the intersection of the groups consists only of translations by \(\F_3\).

The product is \(L^2-1\), with six distinct roots; its translation group is still \(\F_3\), but slope \(-1\) can interchange the two root cosets.
Consequently,
\[
\mathcal A_{f_1f_2}=\{X\mapsto aX+v:a\in\{1,-1\},\ v\in\F_3\}.
\]
In particular, \(X\mapsto-X\) preserves the product but neither individual factor, showing that the intersection group is generally strictly smaller than the stabilizer of the product.
\end{solution}

\solutionheading{ffs:ex-not-all-compositions}
\begin{solution}
The two cubic factors of \(f=R(R+1)\) are coprime.
If \(R(x)=0\), then \(x=0\) or \(x^2=-1\), and the corresponding \(R'(x)=3x^2+1\) equals one or \(-2\), respectively, both nonzero.
If \(R(x)=-1\) and \(R'(x)=0\), then \(x^2=2\), hence \(R(x)=3x=-1\), giving \(x=2\); but this forces \(x^2=4\ne2\), a contradiction.
Thus neither cubic has repeated roots, and \(f\) has six distinct roots in total.

Suppose \(f(aX+b)=a^6f(X)\). The coefficient of \(X^5\) on the right is zero and that on the left is \(6a^5b\), so \(b=0\).
Also \(m=1\), and the set of nonzero exponent differences is \(\{1,2,3,5\}\), with greatest common divisor one; hence \(M_f=\{1\}\) and \(\mathcal A_f\) is trivial.
Nevertheless, \(f=F(R)\) and \(F(Y)=Y^2+Y\) still give a nontrivial composition of degrees two and three.
Thus the compositional decompositions given by affine stabilizers do not exhaust general compositional decompositions.
\end{solution}

\solutionheading{fas:ex-semi-support}
\begin{solution}
Here \(Q=3,t=2,\delta=0,U=L\). Proposition~\ref{fas:prop-semi-invariant-basis} gives the two spaces as
\[
\operatorname{span}_{\F_3}\{1,L^2,L^4\},\qquad
\operatorname{span}_{\F_3}\{L,L^3\},
\]
of dimensions three and two, respectively. Expanding
\[
L^2=X^6+X^4+X^2,\quad L^3=X^9-X^3,\quad
L^4=X^{12}-X^{10}-X^6+X^4,
\]
shows that the common support sets are \(\{0,2,4,6,10,12\}\) and \(\{1,3,9\}\), respectively.
\end{solution}

\solutionheading{fas:ex-nonrational-line}
\begin{solution}
Choose \(i\in\F_9\) with \(i^2=-1\), and choose \(\alpha\) with \(\alpha^3+\alpha=i\). Writing \(\alpha=iy\) gives
\(y^3-y=-1\). Since \(\Tr_{\F_9/\F_3}(-1)=1\ne0\), this cubic is irreducible over \(\F_9\).
Also \(i=\alpha^3+\alpha\in\F_3(\alpha)\), so the tower formula gives \([\F_3(\alpha):\F_3]=6\), proving \(f\) irreducible.
From \(i^3=-i\) and \(\alpha^3=i-\alpha\) we obtain \(\alpha^9=\alpha+i\); hence
\[
b_3=i,\qquad b_3^3=-b_3,\qquad V_f=i\F_3,\qquad L_V=X^3+X.
\]
The translation quotient is given by \(\beta=i\), \(G=Y^2+1\), and \(D=2\), with \(h_V=2\).
The unique candidate center supplied by the trace is \(\delta=2^{-1}\Tr_{\F_9/\F_3}(i)=0\).
The element \(i\) has order four, and \(s_3(4)=2\), so \(t_0=2\). From \(L_V(b)=0\) we obtain
\[
\mathcal A_f=\{X\mapsto aX+b:a\in\{1,-1\},\ b\in i\F_3\}.
\]
Since \(i\F_3\cap\F_3=\{0\}\), the base-field group consists only of \(X\mapsto X\) and \(X\mapsto-X\).
Here every nonzero translation lies outside the base field.
\end{solution}

\solutionheading{fas:ex-prime-counts}
\begin{solution}
Over \(\F_3\) the count is \(3(3-1)/3=2\). Every nonzero \(b\) satisfies \(b^2=1\), and the trace condition requires \(a_0\ne0\),
so the two polynomials are \(X^3-X-1\) and \(X^3-X+1\). Both have full affine stabilizer \(S_3\); their stabilizer over the base field is the translation group of order three.
Over \(\F_4\) the characteristic is two and three divides \(4-1\); formula~\eqref{fas:eq-prime-binomial-count} gives the count
\(4(4-1)(3-1)/3=8\). Their full affine stabilizers are dilation groups \(\mu_3\) centered at points of the base field.
The three monic irreducible quartics over \(\F_2\) are
\[
X^4+X+1,\qquad X^4+X^3+1,\qquad X^4+X^3+X^2+X+1.
\]
None has a base-field root or is divisible by \(X^2+X+1\), so all are irreducible; the irreducible-counting formula gives a total of \((2^4-2^2)/4=3\), so this list contains all monic irreducible quartics.
The first has zero cubic and quadratic coefficients and full affine stabilizer \(A_4\); the other two have nonzero cubic coefficient and trivial full affine stabilizer.
\end{solution}

\solutionheading{fas:ex-twelve-four-translations}
\begin{solution}
The polynomial \(G\) has no root in \(\F_2\) and is irreducible. Take a root \(\beta\in\F_8\); then \(\Tr_{\F_8/\F_2}(\beta)=1\).
If \(x^4+x=\beta\), successive fourth powers give
\[
x^{64}=x+\beta+\beta^4+\beta^2=x+1.
\]
The fiber \(L(X)=\beta\) has four distinct roots, permuted by \(8\)-Frobenius.
The square of this permutation is the nonzero translation \(x\mapsto x+1\), so its order is four; a permutation of four points with order four must be a four-cycle.
Thus \(L(X)-\beta\) is irreducible over \(\F_8\). Also \(\beta=L(x)\), so the tower formula gives \([\F_2(x):\F_2]=12\), proving \(f\) irreducible.
The group \(V_f\) contains \(\F_4\), and the size bound forces \(V_f=\F_4\).
If a slope of order three existed, Corollary~\ref{fas:cor-group-types} would force \(\deg f=4\), a contradiction.
Thus \(\mathcal A_f\) consists exactly of the four translations, and is isomorphic to \(C_2\times C_2\).
One can also argue on the quotient: \(G\) has degree three, and \(3\nmid2-1\), so Theorem~\ref{fas:thm-prime-tame} makes the full affine stabilizer of \(G\) trivial; no nonidentity quotient transformation is available to lift.
\end{solution}

\section*{Solutions to Chapter~\ref{chap:subfield-representations}}
\pdfbookmark[1]{Solutions to Chapter \getrefnumber{chap:subfield-representations}}{solutions.chap:subfield-representations}

\solutionheading{v88:ex-hankel-two-nodes}
\begin{solution}
The first four terms are \((2,3,0,4)\).
The determinant of \(\begin{pmatrix}2&3\\3&0\end{pmatrix}\) is
\(-9=1\in\F_5\), also equal to \(a_1a_2(z_2-z_1)^2=1\).
Writing \(s_{r+2}=us_{r+1}+vs_r\), the first four terms give
\(0=3u+2v\) and \(4=3v\), hence \(v=3=-2\) and \(u=3\).
Thus the monic recurrence polynomial is
\(T^2-3T+2=(T-1)(T-2)\).
\end{solution}

\solutionheading{sfr:ex-quadratic-count}
\begin{solution}
Every element of \(K_2\) outside the base field generates \(K_2/k\), so the first nontrivial jump is the minimum generator-representation degree.
The jump sequence is \((0,\kappa_f(2))\), and the defect formula follows.
Before the first jump, \(U_{f,2}(t)=k\) and hence \(N_f(2,t)=0\); after it, \(U_{f,2}(t)=K_2\),
so \(N_f(2,t)=q^2-q\). Here \(t\) is still restricted to \(0\le t<n\).
\end{solution}

\solutionheading{sfr:ex-quartic-matrices}
\begin{solution}
With rows and columns both indexed by \((1,X,X^2,X^3)\), we have
\[
M_{f_1,2,3}=\begin{pmatrix}
0&1&1&1\\0&0&0&1\\0&0&0&1\\0&0&0&0
\end{pmatrix},\qquad
M_{f_2,2,3}=\begin{pmatrix}
0&1&0&1\\0&1&1&1\\0&0&0&0\\0&1&1&1
\end{pmatrix}.
\]
Their kernels are, respectively, \(\operatorname{span}_{\F_2}\{1,X^2+X\}\) and
\(\operatorname{span}_{\F_2}\{1,X^3+X\}\).
The first three columns of the first matrix have rank one and kernel dimension two; those of the second have rank two and kernel dimension one.
Thus the first polynomial has a quadratic inner component, whereas the second does not.
\end{solution}

\solutionheading{sfr:ex-cyclotomic-four}
\begin{solution}
The polynomial \(f\) has no base-field root and has remainder \(X+1\) modulo \(X^2+X+1\), so it is irreducible.
Also \((X+1)f=X^5+1\), so \(\alpha^5=1\) and \(\alpha\ne1\).
Put \(\eta=\alpha^3+\alpha^2\). Then
\[
\eta^2+\eta+1=\alpha^6+\alpha^4+\alpha^3+\alpha^2+1
=\alpha+\alpha^4+\alpha^3+\alpha^2+1=0.
\]
Thus \(\eta\) generates the four-element subfield and has representation degree at most three.
The polynomial \(f\) has a nonzero cubic term, whereas a composition of two quadratics has none, so the minimum is not two.
Hence \(\kappa_f(2)=3\), \(\Lambda_{f,2}=\{0,3\}\), \(\rho_{f,2}=(1,1,1,2)\), and \(\Delta_f(2)=1\).
\end{solution}

\solutionheading{sfr:ex-tower-eight}
\begin{solution}
Irreducibility of \(P\) is proved in Example~\ref{sfr:ex-quartic-comparison}.
Take a root \(\beta\in\F_{16}\), so \(\Tr_{\F_{16}/\F_2}(\beta)=1\).
If \(\gamma^2+\gamma=\beta\), successive squaring and addition give
\[
\gamma^{16}+\gamma=\beta+\beta^2+\beta^4+\beta^8=1.
\]
Thus \(\gamma\notin\F_{16}\), and the quadratic fiber polynomial is irreducible over \(\F_{16}\).
Also \(\beta\in\F_2(\gamma)\); the tower formula gives \([\F_2(\gamma):\F_2]=8\), so \(f\) is irreducible.

Relative to \(K_4\), the inner component of the original decomposition attains the least degree two, so
\[
\Lambda_{f,4}=\{0,2,4,6\},\qquad\kappa_f(4)=2,\qquad\Delta_f(4)=0.
\]
Relative to \(K_2\), the outer component \(P\) has jump sequence \((0,3)\), minimum generator-representation degree three, and defect one.
Proposition~\ref{sfr:prop-tower-scaling} gives
\[
\Lambda_{f,2}=\{0,6\},\qquad\kappa_f(2)=6,\qquad\Delta_f(2)=2.
\]
An inner component of degree four would require \(\kappa_f(2)=8/2=4\) relative to the quadratic subfield, contrary to the displayed identity.
\end{solution}

\solutionheading{sfr:ex-truncation-tower}
\begin{solution}
(1) By \(\beta^2=\beta+1\), we have \(\deg_\alpha(\beta^2)=2\), whereas \(2\deg_\alpha(\beta)=4=n\).
Multiplication has reached the point where reduction modulo \(f_1\) lowers degrees, so additivity of degrees is no longer preserved.

(2) Choose \(A\in k[Y]\) so that \(\deg A=\kappa_P(c)\) and \(A(\beta)\) generates \(K_c/k\).
Then \(A(h(\alpha))\) also generates \(K_c/k\). Reduction of \(A(h)\) modulo \(f\) cannot increase the degree, so
\[
\kappa_f(c)\le\deg\bigl(A(h)\bmod f\bigr)
\le\deg A\deg h=\kappa_P(c)\kappa_f(d).
\]
In general, one cannot exclude degree reduction modulo \(f\) or assert that this compositional representation is optimal.
When \(\deg h=n/d\), the condition \(\deg A<d\) ensures that each \(A(h)\) has degree less than \(n\),
and these give unique representatives of every element of \(K_d\), yielding equality in the tower relation.
\end{solution}

\solutionheading{str:ex-trace-weight}
\begin{solution}
An element \(v\) belongs to \(\mathcal P_t^{\perp_{\Tr}}\) if and only if, for every \(u\in\mathcal P_t\),
\[
0=\Tr_{E/k}(uv)=\left\langle u,f'(\alpha)v\right\rangle_\alpha.
\]
By Proposition~\ref{str:prop-filter-orthogonal} this is equivalent to \(f'(\alpha)v\in\mathcal P_{n-2-t}\), giving the formula.
In the quartic example, \(\rho_{f,2}(2)=1\), but the unweighted first term is \(\Tr_{E/K}(1)=2=0\), whose prefix rank is zero.
Here \(n-2-t=0\), and the sum of the two dimensions is one rather than \(d=2\), so the weight cannot be omitted.
\end{solution}

\solutionheading{str:ex-hankel-zero-first}
\begin{solution}
Here \(f'=1\) and \(\alpha^4=\alpha+1\), so
\[
c_0=0,\qquad c_1=1,\qquad c_2=1,\qquad c_3=\alpha^3+(\alpha+1)^3=\beta+1.
\]
Hence
\[
\mathcal H_{f,2}=\begin{pmatrix}0&1\\1&1\end{pmatrix},\qquad
\mathcal H_{f,2}\begin{pmatrix}b_0\\b_1\end{pmatrix}
=\begin{pmatrix}1\\\beta+1\end{pmatrix}.
\]
Solving gives \(b_1=1,b_0=\beta\), so \(B_{f,2}=X^2+X+\beta\).
The trace ranks are \((0,1,1,2)\), the jump positions are \((1,3)\), and the corresponding degree jumps are \((0,2)\).
The coordinate polynomial space is \(\operatorname{span}_{\F_2}\{X^2+X,1\}\) and contains the constants, so the least polynomial degree is two and the normalized inner component is \(X^2+X\).
\end{solution}

\solutionheading{v70:ex-rational-rectangle-count}
\begin{solution}
First exclude the zero numerator. Let \(C_{a,b}\) count pairs with \(A\ne0\), monic \(B\), \(\deg A\le a\), \(\deg B\le b\), and \(\gcd(A,B)=1\). Let \(S_{a,b}\) count the same pairs without the coprimality condition. Then
\[
S_{a,b}=\frac{(q^{a+1}-1)(q^{b+1}-1)}{q-1}.
\]
Every pair \((A,B)\) has a unique monic greatest common divisor \(D\). If \(\deg D=j\), there are \(q^j\) choices for \(D\) and \(C_{a-j,b-j}\) choices for \((A/D,B/D)\). Hence
\[
S_{a,b}=\sum_{j=0}^{\min(a,b)}q^jC_{a-j,b-j}.
\]
For \(a,b\ge1\), subtracting \(qS_{a-1,b-1}\) cancels all terms with \(j\ge1\), giving
\[
C_{a,b}=S_{a,b}-qS_{a-1,b-1}=q^{a+b+1}-1.
\]
If \(a=0\), every nonzero constant numerator is coprime to every denominator, and direct counting gives \(C_{0,b}=q^{b+1}-1\). If \(b=0\), the denominator is \(1\), giving \(C_{a,0}=q^{a+1}-1\).
Reduction and the monic denominator ensure that each nonzero rational function is counted once. Adding the zero function \(0/1\) proves the formula.
The count itself does not require \(a+b<n\). When \(a+b=n-1\), there are \(q^n=|E|\) functions, in agreement with the reconstruction bijection. The strict degree condition is needed for injectivity of evaluation, not for the count.
\end{solution}

\solutionheading{str:ex-count-rational}
\begin{solution}
By Corollary~\ref{str:cor-rational-uniqueness}, every such rational function is an invertible fractional linear transformation of a fixed nonconstant \(R=A/C\).
Two invertible matrices define the same fractional linear function if and only if they differ by a nonzero scalar.
Substitution of a nonconstant \(R\) preserves this uniqueness: \(R\) is transcendental over \(k\) in \(k(X)\), so the substitution map \(k(Y)\to k(X)\) is injective.
An invertible two-by-two matrix has \(q^2-1\) choices for its first column and \(q^2-q\) for its second; dividing by the \(q-1\) scalars in each equivalence class gives
\[
\frac{(q^2-1)(q^2-q)}{q-1}=q(q^2-1).
\]
If \(1\in\mathcal W_{f,d}\), take the normalized \(h\) and \(1\) as a basis. A rational function with coprime numerator and denominator can be a polynomial only if the denominator is a nonzero constant,
so all such polynomials are \(ah+b\) with \(a\in k^\times,b\in k\), giving \(q(q-1)\) in total.
If the space does not contain \(1\), Corollary~\ref{str:cor-polynomial-direction} excludes a polynomial representation.
If \(d>2\), then \(2m<n\). Clearing denominators gives a polynomial of degree at most \(2m\), so equal values force equality of the rational functions.
If \(d=2\), put \(\beta=R(\alpha)\in K_2\setminus k\). The group \(\mathrm{PGL}_2(k)\) acts everywhere on \(K_2\setminus k\): an equation \(c\beta+e=0\) with \((c,e)\ne(0,0)\) would force \(\beta\in k\).
The action is transitive. Indeed, \(1,\beta\) is a basis of \(K_2/k\), and every \(\gamma\in K_2\setminus k\) has a unique expression \(u\beta+v\) with \(u\in k^\times,v\in k\).
By the orbit--stabilizer formula, all value fibers have size
\[
\frac{|\mathrm{PGL}_2(k)|}{|K_2\setminus k|}
=\frac{q(q^2-1)}{q^2-q}=q+1.
\]
Distinct group elements give distinct functions after substitution of the nonconstant \(R\), so this is precisely the required number of functions.
\end{solution}

\solutionheading{str:ex-six-two-subfields}
\begin{solution}
Reduction using \(\alpha^6=\alpha+1\) gives \(\beta^2=\beta+1\), so \(\beta\) generates \(K_2/\F_2\).
Direct multiplication and reduction also give
\[
\beta(\alpha+1)=\alpha^3+\alpha^2+1.
\]
Thus \(X^3+X^2+\beta X+\beta+1\) vanishes at \(\alpha\). It is monic of degree three,
and \([E:K_2]=3\), so it equals the relative minimal polynomial \(B_{f,2}\).
Put
\[
A=X^3+X^2+1,\qquad C=X+1,\qquad P=Y^2+Y+1.
\]
Then \(B_{f,2}=A-\beta C\), and
\[
f=A^2+AC+C^2=C^2P(A/C).
\]
By \(A(1)=1\), the polynomials \(A,C\) are coprime, so \(\kappa_f^{\mathrm{rat}}(2)=3\).
The four-element subfield contains only \(0,1,\beta,\beta+1\); the last two have unique representatives of degree five, so
\[
\kappa_f(2)=5,\qquad\Lambda_{f,2}=\{0,5\},\qquad\Delta_f(2)=5-3=2.
\]
Relative to \(K_3\), Example~\ref{str:ex-rational-six} gives \(\eta_f(3)=3\) and
\(\kappa_f^{\mathrm{rat}}(3)=\kappa_f(3)=3>2\).
The same sextic root attains the rational representation degree bound for the quadratic subfield but not for the cubic subfield; the difference lies in the coefficient ranks of the respective root blocks.
\end{solution}

\solutionheading{str:ex-shifted-hankel}
\begin{solution}
(1) Since \(n\ge2\) and \(f\) is irreducible, \(\alpha\ne0\). The matrix \((c_{s+i+j})\) is the Gram matrix of the \(K\)-basis \(1,\alpha,\ldots,\alpha^{m-1}\) for the pairing
\[
(u,v)\longmapsto\Tr_{E/K}\left(\frac{\alpha^suv}{f'(\alpha)}\right).
\]
The weight is nonzero, so the pairing is nondegenerate and the matrix invertible.
Taking the starting index in \eqref{str:eq-relative-recurrence} to be \(s,\ldots,s+m-1\) gives an invertible linear system recovering all \(b_i\) from this consecutive window.

(2) The pairing \((u,v)\mapsto\Tr_{E/K}(wuv)\) is likewise nondegenerate.
If \(R\) gives a recurrence for \((a_j)\), then \(\Tr_{E/K}(w\alpha^jR(\alpha))=0\) for \(0\le j<m\),
so \(wR(\alpha)=0\), hence \(B_{f,d}\mid R\). Taking weighted traces of \(B_{f,d}(\alpha)=0\) directly gives the corresponding recurrence, so the minimal polynomial is exactly \(B_{f,d}\).
Invertibility of the Hankel matrix again follows from its Gram-matrix interpretation.
Degree-profile duality also uses the highest-coefficient formula, whose normalization an arbitrary \(w\) generally lacks; Exercise~\ref{str:ex-trace-weight} already gives a counterexample with \(w=1\).
\end{solution}

\solutionheading{v115:ex-profile-affine}
\begin{solution}
The expression \((c\alpha+a)^j\) is a polynomial in \(\alpha\) of degree \(j\), with leading coefficient \(c^j\ne0\). The change-of-basis matrix is triangular and invertible.
Intersection with any fixed subfield gives identical filtered spaces, so all four kinds of data are unchanged.
Inversion need not preserve them. In Example~\ref{sfr:ex-quartic-comparison}, \(\alpha\) and \(\gamma=\alpha^{-1}\) give jump sets \(\{0,2\}\) and \(\{0,3\}\) for \(\F_4\), minimum generator-representation degrees two and three, and defects zero and one.
\end{solution}

\Needspace{9\baselineskip}
\solutionheading{v116:ex-residual-types}
\begin{solution}
The degree of \(R\) is \(2d\), and every irreducible factor has degree a positive multiple of \(d\). The permitted types are one irreducible factor of degree \(2d\), two distinct irreducible factors of degree \(d\), or the square of one irreducible factor of degree \(d\).
The first two types give one two-cycle and two fixed points, respectively, on the added part of a fiber. Hence \((r_1,r_2)=(0,2)\) and \((2,2)\).
The third type gives a single distinct added root with multiplicity two, so the counts are \((1,1)\). The return permutation records this distinct root; its multiplicity must still be retained in the fiber polynomial.
These types are allowed by the degree condition; this does not prove their realization for every fixed set of parameters.
\end{solution}

\Needspace{9\baselineskip}
\section*{Solutions to Chapter~\ref{chap:root-constructions}}
\pdfbookmark[1]{Solutions to Chapter \getrefnumber{chap:root-constructions}}{solutions.chap:root-constructions}

\solutionheading{ex:supp-composed-six}
\begin{solution}
In the tensor-product basis, multiplication by \(u\otimes v\) has matrix \(M_u\otimes M_v\).
Its diagonal entries are all products of diagonal entries of the two matrices, so its matrix trace is \(\operatorname{tr}(M_u)\operatorname{tr}(M_v)\),
giving the first identity. In the product-of-fields model, the multiplication matrix is block diagonal, each block being the multiplication matrix in the corresponding field; hence the algebra trace is the sum of the field traces of the blocks.
This also shows that the two models compute the same linear functional.

Applying the first identity to tensor-product basis vectors gives
\[
\Tr_{C/k}\bigl((b_i\otimes c_j)(b_{i'}^*\otimes c_{j'}^*)\bigr)
=\delta_{ii'}\delta_{jj'}.
\]
Thus the two displayed bases are dual and the trace pairing is nondegenerate. No coprimality assumption on \(m,n\) is needed.

Finally, \(C\simeq\F_4\times\F_4\) is not a field and its identity \((1,1)\) has trace zero.
Choose \(\omega\in\F_4\) with \(\omega^2+\omega+1=0\); then \(\Tr_{C/\F_2}((\omega,0))=\Tr_{\F_4/\F_2}(\omega)=1\).
Thus although the identity pairs to zero with itself, it is not orthogonal to \((\omega,0)\);
every nonzero vector likewise has a nonorthogonal vector in some nonzero component. Nondegeneracy requires this latter condition, not nonzero self-pairings for all vectors.
\end{solution}

\solutionheading{ex:supp-recursive}
\begin{solution}
The polynomial \(F_0\) is self-reciprocal; setting \(Y=X^2+X+1\) gives \(F_1=Y^2+Y+1=X^4+X+1\).
This is the case \(P=X,a=1\) of the recursive prime-field construction, so \(\deg F_j=2^{j+1}\) and all the polynomials are irreducible.
This recursion combines reciprocal transformation with Artin--Schreier substitution; its definition differs from the quadratic-fiber recursion of the Q-transform.
\end{solution}

\solutionheading{ex:supp-trace-descent}
\begin{solution}
The choice \(s=\ord_\ell(q)\) gives \(\ell\mid q^s-1\), allowing construction of a binomial extension in the auxiliary field using a non-\(\ell\)th power.
Since \(s\mid\ell-1\), we also have \(\gcd(s,\ell)=1\). Consequently, \(q^{\ell^e}\) still has order \(s\) modulo \(\ell\), so the power exponents in the trace sum are distinct modulo \(\ell\).
If the descended element lay in the preceding subfield, these distinct powers would give a nonzero relation of degree less than \(\ell\), a contradiction.
The divisibility \(\ell\mid q^s-1\) supplies the required roots of unity in the auxiliary field; the coprimality \(\gcd(s,\ell)=1\) ensures that the distinct exponent terms cannot cancel during descent.
\end{solution}

\solutionheading{ex:wu-composed-collisions}
\begin{solution}
Let the roots be \(\alpha,\alpha^2\) with \(\alpha^3=1\).
The four composed-sum values are \(0,1,1,0\), so
\[
f\oplus f=X^2(X+1)^2=X^4+X^2.
\]
The four composed-product values are \(\alpha^2,1,1,\alpha\), so
\[
f\odot f=(X+1)^2(X^2+X+1)=X^4+X^3+X+1.
\]
The original roots are simple, but distinct pairs of roots have the same composed value, allowing repeated roots in the composed polynomial.
Under the hypotheses of Theorem~\ref{thm:composed-products}, including nonzero roots in the multiplicative case, coprime degrees ensure that each composed value has degree over the base field equal to the product of the original degrees, excluding such coincidences.
\end{solution}

\solutionheading{ex:power-residue-boundaries}
\begin{solution}
Let \(\omega^2+\omega+1=0\). The two roots are \(\omega,\omega^2\), and the stated exponent is one,
so evaluation at the two roots gives the different values \(\omega\) and \(\omega^2\).
Here \(\mu_3\) is not contained in \(\F_2\), and Frobenius permutes the symbol values.
The condition \(d\mid q-1\) ensures \(\mu_d\subseteq\F_q\), making the symbol value fixed by Frobenius and hence independent of the chosen root.

In the second part, \(2^2=-1\in\F_5\), so
\[
\left(\frac{2}{X(X-1)}\right)_2
=\left(\frac{2}{X}\right)_2\left(\frac{2}{X-1}\right)_2=(-1)^2=1.
\]
But the evaluation isomorphism \(\F_5[X]/(F)\simeq\F_5\times\F_5\) sends \(2\) to \((2,2)\),
and \(2\) is not a square in \(\F_5\), so \((2,2)\) is not a square either.
Multiplicative extension records only the product of the component symbol values, which can cancel each other;
a power-residue condition in the quotient ring must instead hold in every component.
\end{solution}

\solutionheading{v55:ex-stability-base-extension}
\begin{solution}
Level \(n\) is irreducible over the original base field and has degree \(d^n\). Theorem~\ref{thm:factor-base-extension} shows that
after extending the base field it splits into \(\gcd(d^n,r)\) irreducible factors, so it remains irreducible if and only if \(\gcd(d^n,r)=1\).
Each \(d^n\) has the same prime divisors as \(d\), so all levels remain irreducible exactly when \(\gcd(d,r)=1\).
Necessity already follows from the first level.
\end{solution}

\solutionheading{v55:ex-function-stability}
\begin{solution}
Both polynomials take value one at \(0,1\), so they represent the same function.
Stability of \(f\) is given in Example~\ref{v55:ex-binary-stable-cubic}. Direct computation gives
\[
g^{\circ2}=X^4+X+1,\qquad
g^{\circ3}=X^8+X^4+X^2+X+1.
\]
The latter is an affine linearized polynomial of degree \(2^3\).
Theorem~\ref{thm:affine-linearized-obstruction} with \(p=2\) and \(v=3\) excludes its irreducibility, so \(g\) is not stable.
Equality of polynomial functions on the base field records only iteration on that finite set.
Iterated irreducibility requires the next root fiber at each level to form a single orbit under the relative Frobenius over the minimal field of the root at that level; the function algebra modulo \(X^2-X\) does not retain those fibers or ordinary degrees.
\end{solution}

\solutionheading{v55:ex-inverse-quartic-boundary}
\begin{solution}
We have \(\psi(0)=\psi(1)=1\), so there is no linear factor.
The only monic irreducible quadratic over \(\F_2\) is \(X^2+X+1\); the remainder of \(\psi\) modulo it is one, so there is no quadratic factor either and \(\psi\) is irreducible.
Computing the rational composite gives
\[
\Psi^{\circ2}=\frac{\psi^4}{\psi^4+\psi^3+1}.
\]
The denominator is one modulo \(\psi\), so numerator and denominator are coprime. Expanding the denominator gives
\[
X^{16}+X^{12}+X^9+X^8+X^6+X^3+X^2+X+1.
\]
Multiplying the two stated degree-eight factors gives this expression, so the reduced denominator at the second level is reducible and \(\psi\) is not inversely stable.
This contrasts with inverse stability of \(X^2+X+1\): the same reciprocal form must still satisfy the fiber criterion at the next level.
\end{solution}

\section*{Solutions to Chapter~\ref{chap:functions}}
\pdfbookmark[1]{Solutions to Chapter \getrefnumber{chap:functions}}{solutions.chap:functions}

\solutionheading{v72:exer-parallel-subfield}
\begin{solution}
By Theorem~\ref{v72:thm-parallel-difference}, the map \(L=f-f(0)\) is \(\F_p\)-linear and \(s(t)=L(t)/t\).
If \(s(ct)=s(t)\) for every \(c\in K^\times\) and \(t\in\F_q^\times\), then
\[
L(ct)=ct\,s(ct)=cL(t)\qquad(t\ne0).
\]
The same equality is immediate when \(c=0\) or \(t=0\), so \(L\) is \(K\)-linear.
Apply Theorem~\ref{thm:linear-map-representation} with \(k=K\), \(E=\F_q\), and \([E:K]=n/r\).
This gives the unique representation \(L(X)=\sum_{j=0}^{n/r-1}a_jX^{p^{rj}}\); add \(b=f(0)\).
Conversely, the displayed polynomial satisfies \(L(ct)=cL(t)\) for \(c\in K\).
Hence \(s(ct)=L(ct)/(ct)=L(t)/t=s(t)\) for every nonzero \(c,t\).
\end{solution}

\solutionheading{v72:exer-carlitz-mcconnel}
\begin{solution}
The subgroups \(H_5\) and \(H_3\) have orders five and three, respectively.
For \(H_5\), the index is \(d=3\), and \(\ord_3(2)=2\).
The exponents \(0\le j<4\) satisfying \(3\mid2^j-1\) are therefore \(j=0,2\).
Theorem~\ref{v72:thm-carlitz-mcconnel} gives exactly
\[
f(x)=ax+b\quad\text{or}\quad f(x)=ax^4+b,
\qquad a\in H_5,\quad b\in\F_{16}.
\]
There are \(16\cdot5\cdot2=160\) such functions. The normalized functions are \(x\mapsto x\) and \(x\mapsto x^4\).

For \(H_3\), the index is \(d=5\). The successive powers \(2,4,8,16\) are congruent to \(2,4,3,1\) modulo five,
so \(\ord_5(2)=4\). Only \(j=0\) is allowed, giving
\[
f(x)=ax+b,\qquad a\in H_3,\quad b\in\F_{16}.
\]
There are \(16\cdot3=48\) such functions, and the identity is the only normalized one.
In each case, the listed functions satisfy the condition because \(a\in H\) and \(t^{2^j-1}\in H\) for every \(t\ne0\).
The difference between the two classifications is governed by the order of \(2\) modulo the subgroup index: it is two for \(H_5\) and four for \(H_3\).
Preserving the subgroup \(H\) itself is insufficient: the induced Frobenius power must fix every coset in the multiplicative quotient.
\end{solution}

\solutionheading{v74:ex-exponent-rigidity}
\begin{solution}
Since \(H\subseteq\F_q^\times\) and \(u\ge1\), the condition forces \(f(x)-f(y)\ne0\) whenever \(x\ne y\).
Thus \(f\) is a permutation.
Choose \(ur\equiv1\pmod d\) and pass to the cyclic quotient \(G=\F_q^\times/H\), of order \(d\):
\[
\pi(f(x)-f(y))^u\pi(x-y)^v=1.
\]
It follows that
\begin{equation}\label{v74:eq-quotient-power-action}
\pi(f(x)-f(y))=\pi(x-y)^{-rv}.
\end{equation}
For fixed \(y\), the differences \(f(x)-f(y)\), with \(x\ne y\), run through \(\F_q^\times\), so the left-hand side runs through \(G\).
The power map \(z\mapsto z^{-rv}\) is therefore onto \(G\), hence an automorphism.
Thus \(\gcd(rv,d)=1\); since \(\gcd(r,d)=1\), we obtain \(\gcd(v,d)=1\).

\textbf{(2)} If \(v\equiv-u\pmod d\), then \(t^{u+v}\in H\) for every \(t\ne0\), and
\[
(f(x)-f(y))^u(x-y)^v
=\left(\frac{f(x)-f(y)}{x-y}\right)^u(x-y)^{u+v}.
\]
Raising to the \(u\)th power is an automorphism of the quotient, so the original condition is equivalent to the difference quotient belonging to \(H\).
Carlitz--McConnel gives exactly \(ax^{p^j}+b\), with \(a\in H\), \(b\in\F_q\), \(0\le j<n\), and \(d\mid p^j-1\).
Under the present congruence, \(d\mid up^j+v\) is also equivalent to \(d\mid p^j-1\), agreeing with the general formula.

\textbf{(3)} In~\eqref{v74:eq-quotient-power-action}, the map \(\kappa(z)=z^{-rv}\) is one fixed group automorphism.
Lenstra's Theorem~3 gives \(f(x)=ax^{p^j}+b\), initially with only \(a\ne0\), \(b\in\F_q\), and \(0\le j<n\).
Substitute this into the original condition and put \(t=x-y\):
\[
a^u t^{up^j+v}\in H\qquad(t\in\F_q^\times).
\]
Taking \(t=1\) yields \(a^u\in H\).
As \(\gcd(u,d)=1\), the \(u\)th power map of the quotient is invertible, so \(a\in H\).
The remaining condition is \(t^{up^j+v}\in H\) for every nonzero \(t\).
For a primitive element \(\zeta\), we have \(H=\langle\zeta^d\rangle\); thus the condition is equivalent to \(d\mid up^j+v\).
Conversely, the stated parameters give
\[
(f(x)-f(y))^u(x-y)^v=a^u(x-y)^{up^j+v}\in H.
\]
This proves sufficiency and completes the classification.

For \((u,v)=(1,-1)\), the difference quotient belongs to \(H\), with \(d\mid p^j-1\).
For \((u,v)=(1,1)\), the product of the two differences belongs to \(H\), with \(d\mid p^j+1\).
On the quotient, these are the identity and inversion, respectively.
In \(\F_{16}\), take \(H=\langle\zeta^3\rangle\), of index three.
The residue \(2^j\bmod3\) alternates between \(1\) and \(-1\), so the difference-quotient condition permits \(j=0,2\), whereas the product condition permits \(j=1,3\).
The same affine Frobenius form appears in both, but with different Frobenius powers.

Finally, for \(u=d,v=0\), every nonzero \(z\) has \(z^d\in H\).
Thus every permutation satisfies \((f(x)-f(y))^d\in H\) for \(x\ne y\).
The \(d\)th power map of the quotient collapses every element to one and discards all coset labels.
The coprimality hypothesis is therefore substantive.

The main step is not rearranging exponents but converting the mixed difference condition into a fixed action on the quotient:
\[
\begin{gathered}
(f(x)-f(y))^u(x-y)^v\in H\\
\Downarrow\\
\pi(f(x)-f(y))=\kappa\bigl(\pi(x-y)\bigr)\\
\Downarrow\\
\text{one fixed group automorphism on difference cosets}\\
\Downarrow\\
\text{affine Frobenius rigidity of the field map}.
\end{gathered}
\]
Difference quotients are only one expression of compatibility between addition and multiplication.
Difference products and the exponent constraints here exhibit similar rigidity not because one arithmetic expression is intrinsically special,
but because all nonzero difference cosets transform by the same group automorphism.
This converts local difference constraints into restrictions on the global form of the function.
\end{solution}

\solutionheading{v76:ex-frobenius-shape-factor}
\begin{solution}
The linear and constant terms do not contribute to the second divided difference. Substituting \((x+\lambda t)^{p^j}=x^{p^j}+\lambda^{p^j}t^{p^j}\) into equation~\eqref{v76:def-fixed-shape-difference} gives
\[
K_{\lambda,f}(x,t)=a\frac{\lambda^{p^j}-\lambda}{\lambda(\lambda-1)}t^{p^j-2}.
\]
Since \(\lambda\ne0,1\), the condition \(C_j(\lambda)=0\) is equivalent to \(\lambda^{p^j}=\lambda\), that is,
\(\lambda\in\F_{p^{\gcd(j,n)}}\setminus\{0,1\}\), with the subfield taken inside \(\F_q\).

If all the required values lie in \(H\), setting \(t=1\) gives \(aC_j(\lambda)\in H\); in particular, every shape factor is nonzero.
Fix a shape and divide by its value at \(t=1\). This gives \(t^{p^j-2}\in H\) for every \(t\ne0\).
Writing \(H=\langle\zeta^d\rangle\) for a primitive element \(\zeta\) of \(\F_q\), this is equivalent to \(d\mid p^j-2\).
Conversely, this divisibility condition puts every \(t^{p^j-2}\) in \(H\); multiplying by \(aC_j(\lambda)\in H\) proves sufficiency.
If \(p^j=2\), the exponent is zero and \(C_j(\lambda)=1\), recovering the constant-value condition for a quadratic term.

Choose \(\lambda_0\in\Lambda\). A suitable \(a\) exists if and only if the divisibility condition holds, every shape factor is nonzero, and
\[
\frac{C_j(\lambda)}{C_j(\lambda_0)}\in H\qquad(\lambda\in\Lambda).
\]
In that case, one may take \(a=C_j(\lambda_0)^{-1}\), and the full set of permitted coefficients is \(C_j(\lambda_0)^{-1}H\).
This classifies only the prescribed family \(aX^{p^j}+bX+c\). Excluding its nonquadratic members does not automatically exclude arbitrary polynomial functions.
\end{solution}

\solutionheading{v76:ex-multishape-recovery}
\begin{solution}
For each \(\lambda\), put \(K_\lambda=\F_p(\lambda)\) and \(s_\lambda=|K_\lambda|>2\).
Theorem~\ref{v76:thm-fixed-shape-constant} gives a representation
\[
f(X)=a_\lambda X^2+L_\lambda(X)+f(0),
\]
where \(L_\lambda\) is the standard \(s_\lambda\)-linearized polynomial.
Its exponents belong to \(1,s_\lambda,s_\lambda^2,\ldots\), none is two, and the entire expression has degree less than \(q\).
Thus every \(a_\lambda\) is the coefficient of \(X^2\) in the same function representative. They are all equal, say to \(a\).
Put \(L=f-aX^2-f(0)\). This is an additive map commuting with scalar multiplication by every \(\lambda\in\Lambda\), hence with all scalar multiplications from the generated field \(K\).
Therefore \(L\) is \(K\)-linear. Conversely, apply the same theorem to each \(\lambda\).

If \(K=\F_q\), every \(K\)-linear map has the form \(L(X)=b_1X\), giving a polynomial of degree at most two.
If \(K\) is proper, put \(s=|K|\). Since \(\Lambda\) is nonempty, \(s>2\), and also \(s<q\).
The function \(aX^2+X^s\) satisfies the constant-value condition for every shape in \(\Lambda\), but \(2<s<q\) prevents a representative of degree at most two.
This proves the equivalence in part~(1).

For part~(2), choose any ordered pair of distinct nodes from \(x,x+t,x+\lambda t\). Take the first node as the new base point and the second minus the first as the new displacement. The relative parameter of the third point is one of the six listed values.
Symmetry of the divided difference in its nodes shows that this merely rewrites the same value. The new base point and nonzero displacement still range over all possibilities, so every transformation preserves the condition.
Conversely, any three pairwise distinct nodes can be written as \(x,x+t,x+\nu t\), with \(t\ne0\) and \(\nu\ne0,1\).
If \(\Lambda\) meets every shape orbit, the subgroup condition holds for all triples. Theorem~\ref{v76:thm-second-difference-rigidity} gives \(f(X)=aX^2+bX+c\), with \(a\in H\).
This is a sufficient orbit-covering condition, not a claim that every orbit representative is indispensable.

Finally, take \(q=16\). Since \(\gcd(3,4)=1\), the elements of \(\F_{16}\) satisfying \(\lambda^8=\lambda\) are precisely \(\F_2\). Thus the specified \(\lambda\) has \(C_3(\lambda)\ne0\).
Set
\[
f(X)=C_3(\lambda)^{-1}X^8.
\]
The preceding exercise gives \(K_{\lambda,f}(x,t)=t^6\in(\F_{16}^\times)^3\).
Yet the unique representative has degree eight. Although \(\F_2(\lambda)=\F_{16}\), the function is not quadratic.
Thus subfield generation exactly governs the constant-value condition in part~(1), but it does not alone ensure rigidity under a subgroup-value condition.
\end{solution}

\solutionheading{v58:ex-linearized-minimal-values}
\begin{solution}
Put \(r=\dim_{\F_q}\ker(L|_E)\). Since \(L\) is linear, every nonempty fiber has \(q^r\) elements,
so \(|\mathcal V_L|=q^{n-r}\). On the other hand, \(q^t\mid q^n\), so the minimum value-set bound is \(q^{n-t}\).
The two are equal if and only if \(r=t\), that is, if the kernel size equals the ordinary degree.
This is equivalent to \(L\) having \(q^t\) distinct roots, all in \(E\),
and thus to \(L\) dividing the polynomial \(X^{q^n}-X\) without repeated roots.

In \(\F_{27}\), the kernel of \(X^3-X\) is \(\F_3\), of size three, so it is a minimal value set polynomial.
If \(x^3+x=0\) and \(x\ne0\), then \(x^2=-1\), making \(x\) have multiplicative order four;
but \(4\nmid26=|\F_{27}^\times|\), so this is impossible.
Thus \(X^3+X\) has only zero in its kernel and is a permutation, with value-set size 27, greater than the degree-based bound nine.
\end{solution}

\solutionheading{v70:ex-linearized-values-all-extensions}
\begin{solution}
Corollary~\ref{v70:cor-all-extension-kernels} gives kernel dimension \(d_m\), and hence kernel size \(q^{d_m}\).
Every nonempty fiber is an additive coset of the kernel. There are therefore \(q^{m-d_m}\) nonempty fibers, each of size \(q^{d_m}\), giving
\[
C_{\mathcal L_h|_{K_m}}=q^{m-d_m}(q^{d_m})^2=q^{m+d_m}.
\]
Their equal sizes give equality in the collision bound. This argument allows \(p\mid m\) and repeated factors of \(h\),
and does not require the ordinary degree of \(\mathcal L_h\) to be less than \(q^m\); the coefficients must lie in the fixed base field \(k\) for the ordinary gcd formula to apply.
If \(h=0\) is considered separately, the zero map has value-set size one and collision count \(q^{2m}\).

For the second part, put \(P=T^2+T+1\). Its roots have order three, and the operator polynomials of \(A\) and \(A^{\circ3}\) are \(P\) and \(P^3\), respectively.
When \(3\nmid m\), both are coprime to \(T^m-1\), so both permute \(\F_{2^m}\).
When \(3\mid m\), write \(m=2^v u\) with \(u\) odd. Then
\(T^m-1=(T^u-1)^{2^v}\), and \(T^u-1\) is squarefree.
The two kernel dimensions are consequently
\[
2,\qquad 2\min\{3,2^{v_2(m)}\}.
\]
Thus both are permutations exactly when \(3\nmid m\). In particular,
\[
\begin{array}{c|cc|cc}
 &\dim\ker A&\dim\ker A^{\circ3}&|A(K_m)|&|A^{\circ3}(K_m)|\\ \hline
m=6&2&4&16&4\\
m=12&2&6&1024&64
\end{array}
\]
The kernel dimensions in the table are over \(\F_2\). These values follow from the order and multiplicity of one factor, without evaluating the maps point by point.

Finally, \(T\) acts on \(K_m\) as the invertible Frobenius \(\sigma\), so
\(\mathcal L_{T^ah}=\sigma^a\circ\mathcal L_h\) has the same kernel, value-set size, and collision count;
one may also use \(\gcd(T,T^m-1)=1\) directly. For one-sided recurrences, \(T\) acts as a shift \(S\) that need not be invertible,
and the zero-root primary component of the minimal polynomial gives exactly the transient; see Theorem~\ref{thm:transient-period}.
The meaning of an operator polynomial depends on the specified operator.
\end{solution}

\solutionheading{v58:ex-factorial-value-count}
\begin{solution}
For \(b\in\F_Q\), put \(r_b=N_b(F)\); then \(0\le r_b\le m\).
Choosing an ordered tuple of \(j\) distinct elements in a fiber of size \(r_b\) gives
\((r_b)_j=r_b(r_b-1)\cdots(r_b-j+1)\) choices. Therefore
\[
M_j(F)=\sum_b(r_b)_j=j!\sum_b\binom{r_b}{j}.
\]
These identities are interpreted in the integers, not by dividing by \(j!\) in the finite field.
For \(r>0\), the binomial theorem gives
\(\sum_{j=1}^{r}(-1)^{j-1}\binom rj=1\); for \(r=0\), the sum is zero.
Summing over fibers gives the value-set formula. Taking \(j=1,2\) gives \(M_1=Q\) and \(M_2=C_F-Q\), respectively.

The nonempty fibers in the example have sizes \(3,6,6,6,6\), so
\[
\bigl(M_1,M_2,M_3,M_4,M_5,M_6\bigr)
=\bigl(27,126,486,1440,2880,2880\bigr).
\]
Hence
\[
|\mathcal V_F|=27-63+81-60+24-4=5.
\]
Two-point collisions give a lower bound, whereas the full collection of multiple-point counts recovers the exact number of nonempty fibers.
\end{solution}

\solutionheading{ex:function-inverse}
\begin{solution}
Since \(5\cdot5\equiv1\pmod6\), the inverse of \(X^5\) is again \(X^5\): exponents of nonzero elements are computed modulo six, and zero is also fixed.
Composition gives \(X^{25}\), which is not equal to \(X\) as a polynomial, but \(X^{25}\equiv X\pmod{X^7-X}\); that is, they take the same value at every element of \(\F_7\).
\end{solution}

\solutionheading{ex:delta-idempotent}
\begin{solution}
The function \(\delta_a\) has value one at \(a\) and zero everywhere else.
Thus the three identities hold pointwise at all \(q\) points.
The theorem on polynomial functions over a finite field makes pointwise equality equivalent to equality in the quotient ring modulo \(X^q-X\), proving the assertion.
These mutually orthogonal idempotents are the coordinate idempotents in the decomposition of the function algebra into \(q\) field components.
\end{solution}

\solutionheading{ex:wu-function-frobenius}
\begin{solution}
Write \(h\) as a polynomial with base-field coefficients. Raising to the \(q\)th power fixes the coefficients, and \(X^q=X\), so \(h^q=h\).
On the other hand, \(R\) is isomorphic to \(\prod_{a\in\F_q}\F_q\);
if \(a\ne b\), the point-indicator functions \(\delta_a,\delta_b\) are nonzero but have zero product, so it is not a field.
Operations on the root set are inherited from the chosen extension field and therefore have no zero divisors;
operations on functions on the root set are defined pointwise, and nonzero functions with different supports may have zero product.
\end{solution}

\solutionheading{ex:wu2-trace-shears}
\begin{solution}
Since \(g,h\) take values in the base field and \(p\mid n\), we have \(\Tr(h(t))=0\). Thus \(\Tr(F_h(x))=\Tr(x)\),
and substitution into the definition gives \(F_g(F_h(x))=x+h(\Tr(x))+g(\Tr(x))\). The inverse map is \(F_{-g}\).
The trace is surjective, so \(F_g=\id\) if and only if \(g\) is the zero function.
The identity \(F_g^m=F_{mg}\) shows that the permutation associated with nonzero \(g\) has order exactly \(p\). This gives a permutation subgroup isomorphic to the additive function group \(k^k\).
\end{solution}

\solutionheading{v70:ex-additive-return}
\begin{solution}
For any \(x\in x_0+W\), the traces during its first \(\ell\) steps are \(t_0,\ldots,t_{\ell-1}\).
Repeated substitution into \(F(x)=ax+g(\Tr(x))\) gives
\[
F^\ell(x)=a^\ell x+B.
\]
The base cycle first returns after \(\ell\) steps, so \(\Tr(F^\ell(x_0))=t_0=\Tr(x_0)\). Hence
\(\delta=F^\ell(x_0)-x_0\in W\). Substituting \(x=x_0+v\) and subtracting \(x_0\) gives \(R(v)=a^\ell v+\delta\).
By Lemma~\ref{lem:fiber-return}, each cycle of \(R\) of length \(e\) corresponds to exactly one cycle of \(F\) of length \(\ell e\),
and these cycles cover all trace fibers over the chosen base cycle.

If \(a^\ell\ne1\), then \(v_*=-\delta/(a^\ell-1)\in W\) is the unique fixed point of \(R\).
In the coordinate \(v=v_*+w\), the return map becomes \(w\mapsto a^\ell w\).
Put \(e=\ord(a^\ell)\). For \(w\ne0\), the equality \((a^\ell)^j w=w\) is equivalent to \((a^\ell)^j=1\),
so every nonzero \(w\) has cycle length \(e\).
Thus \(F\) has one cycle of length \(\ell\) and \((q^{n-1}-1)/e\) cycles of length \(\ell e\).

If \(a^\ell=1\) and \(\delta=0\), then \(R\) is the identity, so \(F\) has \(q^{n-1}\) cycles of length \(\ell\).
If \(a^\ell=1\) and \(\delta\ne0\), then \(R^j(v)=v+j\delta\).
The nonzero element \(\delta\) has additive order \(p\), so every cycle of \(R\) has length \(p\).
Consequently, \(F\) has \(q^{n-1}/p\) cycles of length \(\ell p\).

Finally, replacing \(x_0\) by \(x_0+w_0\), with \(w_0\in W\), replaces \(\delta\) by
\(\delta+(a^\ell-1)w_0\). The two return maps are conjugate by a translation of coordinates and have the same cycle data.
When \(n=1\), we have \(W=\{0\}\) and \(\delta=0\), so the nonzero-translation case cannot occur;
in the first case, the number of cycles of length \(\ell e\) is automatically zero.
\end{solution}

\solutionheading{ex:agw-return-comparison}
\begin{solution}
We have \(\widetilde g(1)=2\), \(\widetilde g(3)=4\), and \(\widetilde g(9)=5\). Since \(5^4=1\), the fourth powers of these three values agree with those in the original example,
so the induced map remains \((1,3,9)\), and the multiplicative-fiber criterion makes \(\widetilde F\) a permutation.
Take \(z_0=x_0=1\). The return coefficient is
\[
\widetilde\delta=2^9\cdot4^3\cdot5=1\quad\text{in }\F_{13}.
\]
Relative to \(\mu_4=\langle5\rangle\), the return map becomes \(v\mapsto3v\pmod4\), with two fixed points \(0,2\) and one two-cycle \((1,3)\).
Thus \(\widetilde F\) has two three-cycles and one six-cycle on the nonzero elements, together with the fixed point \(0\).
The original example and this problem have the same base cycle but different return coefficients; changing the fiber cycles changes the total cycle decomposition.
\end{solution}

\solutionheading{ex:wu10-dickson-inverse}
\begin{solution}
Applying the root-parameter identity twice gives
\[
D_m(D_r(u+u^{-1}))=u^{mr}+u^{-mr}=D_{mr}(u+u^{-1}).
\]
Regard \(u\) as an indeterminate. Substituting \(u+u^{-1}\) into a nonzero polynomial of degree \(d\) leaves the coefficient of the highest power \(u^d\) equal to the original leading coefficient,
so this substitution is injective. Hence the ordinary polynomial identity \(D_m\circ D_r=D_{mr}\) follows.

By \(5\cdot38\equiv1\pmod{63}\), the polynomial \(D_{38}\) gives the inverse function of \(D_5\) on \(\F_8\).
Using the defining recurrence and reducing modulo \(X^8-X\) gives
\[
D_{38}(X)\equiv X^5+X^4+X^3+X^2+X\pmod{X^8-X}.
\]
This is the unique polynomial of degree less than eight representing the inverse function.
Its composite with \(D_5\) equals \(X\) only modulo \(X^8-X\);
in the ordinary polynomial ring, \(D_{38}\circ D_5=D_{190}\) is not a linear polynomial.
\end{solution}

\solutionheading{afr:ex-sixtyfour-cycles}
\begin{solution}
Here \(d=1,m=6,p=2\). By \(\ord(G)=12\) and \(v_2(\ell)=v_2(6)+1=2\), cycle lengths can only be four or 12.
Formula~\eqref{afr:eq-affine-fixed} gives \(N_4(G)=2^{\gcd(6,4)}=4\), so there is exactly one four-cycle.
The remaining \(64-4=60\) elements form five 12-cycles.
\end{solution}

\solutionheading{afr:ex-weighted-sixteen}
\begin{solution}
The norm is \(c=a^5=\omega\ne1\), so every \(b\) gives a unique fixed point; taking it as origin gives the same linear operator \(S_a\).
Thus all \(F_{a,b}\) have the same cycle type. Here \(d=m=2\), and \eqref{afr:eq-linear-fixed} gives
\[
 N_1=1,\qquad N_2=1,\qquad N_3=4,\qquad N_6=16.
\]
The last equality shows that every cycle length divides six. Successively subtracting shorter cycles gives one fixed point, one three-cycle, and two six-cycles.
\end{solution}

\solutionheading{v53:ex-cycles-root-sign}
\begin{solution}
The two irreducible factors \(T-1,T+1\) have orders one and two, respectively, and the largest multiplicity is two.
The characteristic is three, so \(\ord(A)=3\cdot2=6\) and every cycle length divides six.
By \eqref{v53:eq-cyclic-gcd-count} and \eqref{v53:eq-linear-cycle-count},
\[
\begin{array}{c|rrrr}
t&1&2&3&6\\\hline
\deg\gcd(f,T^t-1)&1&2&2&3\\
N_t(A)&3&9&9&27\\
C_t(A)&3&3&2&2
\end{array}
\]
For example, \(C_6=(27-9-9+3)/6=2\). These cycles contain \(3+2\cdot3+3\cdot2+6\cdot2=27\) points.
There are \(10\) cycles in total, so \(\operatorname{sgn}(A)=(-1)^{27-10}=-1\).
On the other hand, the characteristic polynomial of \(A\) is \(f\), so \(\det A=-1\), a nonsquare in \(\F_3\);
Theorem~\ref{v53:thm-linear-sign} gives the same sign \(-1\).

Proposition~\ref{prop:recurrence-root-correspondence} identifies this cyclic module compatibly with \(W_f=\ker\mathcal L_f\),
with multiplication by \(\overline T\) corresponding to \(\sigma(x)=x^3\). Thus \(W_f\) has three fixed points, three two-cycles, two three-cycles, and two six-cycles.
Its splitting field is \(\F_{3^{\ord(f)}}=\F_{3^6}\): every root has period dividing six, and some root has period six, excluding a smaller splitting field.
\end{solution}

\solutionheading{v97:ex-dickson-functional}
\begin{solution}
The polynomial \(H\) is squarefree and has no zero root.
Lagrange interpolation for \(\deg R<N\), followed by comparison of highest coefficients, gives
\[
 \Lambda(R)=\sum_{H(t)=0}\frac{R(t)}{H'(t)}.
\]
Substitute the derivative formula of Corollary~\ref{v97:cor-dickson-weights} to obtain the desired expression.
Taking \(R=Y^j\) yields
\[
 \sum_{H(t)=0}(t-2)J(t)t^{j-1}
 =\begin{cases}-2,&j=N-1,\\0,&0\le j<N-1.\end{cases}
\]
For \(j=0\), the inverses exist because the roots are nonzero.
Frobenius permutes the roots of \(H\), and \(R,J\) are defined over the base field, so it permutes the summands and fixes their sum.
This explains directly why the root-value expression recovers a base-field coefficient.
\end{solution}

\solutionheading{v106:ex-redei-rational}
\begin{solution}
Put \(\iota(Z)=Z^{-1}\). Coefficient Frobenius sends \(\theta\) to \(-\theta\), and hence \(\phi\) to \(\iota\circ\phi\). Since \(Z^n\) commutes with \(\iota\), the function \(R_n\) is fixed and belongs to \(\F_q(X)\). Explicitly,
\[
R_n(X)=\theta\frac{(X+\theta)^n+(X-\theta)^n}{(X+\theta)^n-(X-\theta)^n},
\]
reduced as a rational function and evaluated on the projective line.
If \(m\) is odd, put \(Q=q^m\). Then \(\theta^Q=-\theta\), and \(\phi\) maps \(\mathbb P^1(\F_Q)\) bijectively to \(\mu_{Q+1}\): finite-point images satisfy \(z^Q=z^{-1}\), infinity maps to \(1\), and \(\theta(z+1)/(z-1)\) verifies the inverse direction. Permutation is thus equivalent to \(\gcd(n,Q+1)=1\). If \(m\) is even, \(\theta\in\F_Q\) and \(\phi\) is a fractional linear permutation of \(\mathbb P^1(\F_Q)\). The power map fixes \(0,\infty\); on the remaining points it is a permutation exactly when \(\gcd(n,Q-1)=1\).
For \(q=7,\delta=3,n=3\), expansion gives
\[
R_3(X)=\frac{X^3+2X}{3X^2+3}.
\]
Since \(7\equiv1\pmod3\), odd \(m\) gives \(3\nmid7^m+1\), while even \(m\) gives \(3\mid7^m-1\). Thus permutation holds exactly in odd-degree extensions.
\end{solution}

\solutionheading{v116:ex-stability-gap}
\begin{solution}
If \(\mathcal E_H(f)\le\varepsilon q^2\), choose a nearest model and put \(r=D_H(f)\le\eta q<(q-h+1)/2\). The local bound gives
\(\mathcal E_H(f)\ge r(q-h+1)\ge(1-1/d)qr\).
If \(\mathcal E_H(f)>\varepsilon q^2\), use \(D_H(f)\le q\).
Conversely, given \(\mathcal E_H(f)\ge cqD_H(f)\), choose any \(\eta>0\) in the stated range and take \(\varepsilon=c\eta\).
\end{solution}

\Needspace{9\baselineskip}
\section*{Solutions to Chapter~\ref{chap:sequences}}
\pdfbookmark[1]{Solutions to Chapter \getrefnumber{chap:sequences}}{solutions.chap:sequences}

\solutionheading{ex:sequence-trace8}
\begin{solution}
We have \(\Tr(1)=1\) and \(\Tr(\alpha)=\Tr(\alpha^2)=0\), so the trace sequence has initial state \((1,0,0)\).
It satisfies the same recurrence and hence equals the sequence in the main text by uniqueness.
All nonzero solutions are parametrized by \(c\in\F_8^\times\); since \(\alpha\) is primitive, every \(c\) has a unique expression \(\alpha^j\) with \(0\le j<7\).
Then \(\Tr(c\alpha^t)=\Tr(\alpha^{t+j})\), exactly a cyclic shift.
\end{solution}

\solutionheading{ex:sequence-repeated}
\begin{solution}
The irreducible factor \(T^2+T+1\) has order three, the largest multiplicity is three, and the required power of two is four, so the least period is 12.
In general, put \(d=\deg h\geq1\) and \(C(Z)=Z^dh(Z^{-1})\) and define the generating function of the sequence by \(1/C(Z)\).
Since \(C(0)=1\), this formal power series exists uniquely; comparing coefficients gives the recurrence \(h(S)s=0\).
The numerator is one and cannot cancel against the denominator, so the minimal polynomial is exactly \(h\). If \(h=1\), take the zero sequence.
\end{solution}

\solutionheading{ex:sequence-generating}
\begin{solution}
The initial state and recurrence give
\[
G_s(Z)=\frac{1+Z^2}{1+Z^2+Z^3}
=\frac{1+Z^3+Z^5+Z^6}{1-Z^7}
\qquad\text{in }\F_2[[Z]].
\]
In the first expression, numerator and denominator are coprime, and the reciprocal of the denominator is \(T^3+T+1\), which is therefore the minimal polynomial.
The second expression groups terms in blocks of length seven, displaying the period directly.
\end{solution}

\solutionheading{ex:wu-transient}
\begin{solution}
The sequence satisfies \(s_{t+3}=s_{t+2}\) and is therefore annihilated by \(T^2(T+1)\).
Its minimal polynomial must divide this polynomial.
For proper factors, \((S+1)s\ne0\) and \(S(S+1)s\ne0\), whereas \(S^j s\ne0\) for every \(j\ge0\);
thus no proper factor annihilates the sequence, and the minimal polynomial is exactly \(T^2(T+1)\).
The generating function is
\[
G(Z)=1+\frac{Z^2}{1+Z}
=\frac{1+Z+Z^2}{1+Z}
=Z+\frac1{1+Z}.
\]
In \(\F_2[T]\) we have \(T^2+(T+1)^2=1\), so in \eqref{eq:transient-projectors} we may take \(U=1\) and \(V=T+1\). Hence
\[
s_{\mathrm{nil}}=(S+1)^2s=(0,1,0,0,\ldots),\qquad
s_{\mathrm{per}}=S^2s=(1,1,1,1,\ldots).
\]
The two parts have generating functions \(Z\) and \(1/(1+Z)\), respectively, in agreement with the decomposition above.
The denominator records only the eventual period one; the additional polynomial term \(Z\) records the transient part.
The factor \(T^2\) of the sequence's minimal polynomial cannot be recovered from the reduced denominator alone.
\end{solution}

\solutionheading{ex:wu10-transient-repeated}
\begin{solution}
For any \(P\in k[T]\), the correspondence between shifting and multiplication by \(\alpha\) gives
\[
\bigl(P(S)s^{(c)}\bigr)_j
=\lambda\bigl(P(\alpha)c(\alpha)\alpha^j\bigr).
\]
Since \(1,\alpha,\ldots,\alpha^8\) span \(A\), the expression vanishes for every \(j\ge0\)
if and only if \(\lambda(P(\alpha)c(\alpha)v)=0\) for every \(v\in A\).
Nondegeneracy of the coefficient pairing makes this equivalent to \(P(\alpha)c(\alpha)=0\), namely \(f\mid Pc\).
Write \(d=\gcd(f,c)\), \(f=df_1\), and \(c=dc_1\); then \(\gcd(f_1,c_1)=1\),
so \(f\mid Pc\) is equivalent to \(f_1\mid P\). Thus the minimal polynomial is exactly \(f_1\).
This reduction determines both the recurrence and its minimality.

Take \(c=T(T^2+T+1)\). Then
\[
h=T^2(T+1)^2(T^2+T+1)
=T^6+T^5+T^3+T^2.
\]
Extracting coefficients directly from \(c=T^3+T^2+T\) gives the first six terms \(0,0,0,0,0,1\),
followed by the recurrence
\[
s_{j+6}=s_{j+5}+s_{j+3}+s_{j+2}\qquad(j\ge0).
\]
Write \(g=(T+1)^2(T^2+T+1)=T^4+T^3+T+1\).
If the sequence has period \(N\) starting at index \(u\), then \(h\mid T^u(T^N-1)\).
Since \(T\nmid T^N-1\), necessarily \(u\ge2\).
Also, the roots of \(T^2+T+1\) have multiplicative order three, so \(3\mid N\);
if \(N\) is odd, \(T^N-1\) has nonzero derivative at \(T=1\) and cannot contain \((T+1)^2\), so we also require \(2\mid N\).
Conversely, \(g\mid T^6-1=(T+1)^2(T^2+T+1)^2\).
Thus the least transient length is two and the least eventual period six, and both are attained simultaneously.

In \(k[T]\) there is the explicit Bézout identity
\[
(T^3+T+1)T^2+(T+1)g=1.
\]
Thus \eqref{eq:transient-projectors} gives
\[
\begin{aligned}
s_{\mathrm{per}}&=(S^5+S^3+S^2)s,\\
s_{\mathrm{nil}}&=(S^5+S^3+S^2+1)s.
\end{aligned}
\]
The recurrence above gives \(s_6=s_7=1\) and \(s_8=s_9=s_{10}=0\).
Hence the first two terms of \(s_{\mathrm{nil}}\) are both one, and \(S^2s_{\mathrm{nil}}=0\), so
\[
s_{\mathrm{nil}}=(1,1,0,0,0,\ldots).
\]
Meanwhile, the first six terms of \(s_{\mathrm{per}}\) are \(1,1,0,0,0,1\),
and this block of length six repeats thereafter. Adding gives the first two terms \(0,0\) of the original sequence,
followed from index two by repeated copies of \(0,0,0,1,1,1\).
This also verifies directly the separation of the two components: the finitely supported part vanishes exactly after two shifts,
whereas the purely periodic part is periodic from the initial term.

Finally, put \(c'=T^3(T+1)^2(T^2+T+1)\). The formula in the first part gives
\[
\begin{array}{c|c|c|c}
\text{Parameter}&\text{Minimal polynomial}&\text{Least transient length}&\text{Least eventual period}\\\hline
1&f&3&6\\
c&T^2(T+1)^2(T^2+T+1)&2&6\\
c'&T^2+T+1&0&3
\end{array}
\]
The nonzero-root part in the first row is \((T^3+1)^2=T^6+1\), of order six;
the third row is the purely periodic sequence repeating \(0,1,1\).
Although all three are annihilated by \(f\), common factors of the parameter and \(f\) remove some factors or lower their multiplicities.
Thus an arbitrary nonzero annihilating polynomial gives only an upper bound on the transient length and a multiple of the eventual period;
the exact values must be read from the minimal polynomial of the sequence.
\end{solution}

\solutionheading{ex:wu2-recurrence-profile}
\begin{solution}
Since \(\varphi_2(f)=2^4(1-2^{-1})(1-2^{-2})=6\), six sequences have minimal polynomial exactly \(f\).
The two irreducible factors have orders one and three; multiplicity two contributes an additional factor two, so these sequences have least period six.
A sequence has least period three precisely when it satisfies \(h_s\mid T^3-1=(T+1)(T^2+T+1)\)
and does not have least period one. The possible minimal polynomials are \(T^2+T+1\) and \((T+1)(T^2+T+1)\),
giving \(3\) sequences each, six in total. Alternatively, the solution space of period dividing three has \(2^3\) elements, while period dividing one gives two; their difference is six.
\end{solution}

\solutionheading{ex:wu3-bm-short}
\begin{solution}
The Berlekamp--Massey algorithm gives, for initial segments of lengths \(1,\ldots,6\), respectively,
\[
\begin{aligned}
&(1,1+Z),\quad(1,1),\quad(1,1),\\
&(3,1+Z^3),\quad(3,1+Z^3),\quad(3,1+Z^2+Z^3).
\end{aligned}
\]
Nonzero discrepancies occur at positions \(0,1,3,5\); Theorem~\ref{thm:berlekamp-massey} proves minimality at every step.
The final output is \(T^3+T+1\). Although the connection polynomial for the first two terms is the constant one, the relation is imposed only from index one onward; it represents \(s_1=0\). Setting the length to zero would also require \(s_0=0\), contrary to the input.
\end{solution}

\solutionheading{ex:wu3-bm-sharp}
\begin{solution}
Let \(a_j\) be one only at \(j=D-1\) and zero elsewhere; let \(b_j\) be one at \(j\equiv D-1\pmod D\) and zero elsewhere.
The sequences agree at indices \(0,\ldots,2D-2\) but have values zero and one, respectively, at index \(2D-1\).
The polynomial \(T^D\) annihilates \(a\), and \(T^D-1\) annihilates \(b\).
If either had a monic annihilating polynomial of degree \(L<D\), the recurrence equation corresponding to index \(D-1\)
would have all previous terms zero and would force this term to be zero, contrary to its value one. Thus both minimal degrees are \(D\).
This proves that when only a degree bound is known, the uniform recovery guarantee from \(2D\) terms cannot be reduced in the worst case to \(2D-1\) terms.
\end{solution}

\solutionheading{ex:wu7-f64-two-polynomials}
\begin{solution}
Example~\ref{exm:f64-computation} proves \(f=T^6+T+1\) primitive, so \(m_{\alpha,\F_2}=f\) and \(\ord(\alpha)=63\), and the Frobenius orbit of \(\alpha\) has length six.
Successive squaring and reduction using \(\alpha^6=\alpha+1\) give
\[
\alpha^8=\alpha^3+\alpha^2,\qquad
\alpha^{16}=\alpha^4+\alpha+1,\qquad
\alpha^{32}=\alpha^3+1.
\]
The six conjugates sum to zero, so \(1+T+\cdots+T^5\) annihilates \(\alpha\). They all belong to \(\langle1,\alpha,\alpha^2,\alpha^3,\alpha^4\rangle_{\F_2}\).
Conversely, \(\alpha^3=\alpha^8+\alpha^2\) and \(1=\alpha^{32}+\alpha^3\),
so all five powers belong to the span of the conjugates. As part of a power basis, they are independent, so this span has dimension exactly five.
Thus the displayed monic annihilating polynomial of degree five already has least degree, giving the required \(\mu_\alpha\).

Term \(j\) of the sequence arises from the multiplicative iteration \(M_\alpha^j(1)=\alpha^j\),
whereas \(\mu_\alpha\) describes linear relations in the squaring iteration \(\sigma^j(\alpha)=\alpha^{2^j}\).
By Theorem~\ref{thm:sequence-trace}, a nonzero trace sequence has minimal polynomial \(f\) and least period 63.
Nondegeneracy of the trace pairing ensures that the sequence for this parameter \(c=1\) is nonzero.
In the realization of Example~\ref{exm:two-sixtyfour-rootsets}, this multiplicative period becomes the Frobenius period of the nonzero roots in \(W\),
not the Frobenius period of \(\alpha\) in \(E\).
\end{solution}

\section*{Solutions to Chapter~\ref{chap:codes}}
\pdfbookmark[1]{Solutions to Chapter \getrefnumber{chap:codes}}{solutions.chap:codes}

\solutionheading{v88:ex-constacyclic-bch-five}
\begin{solution}
Multiplication in \(\F_5[X]\) gives
\((X^2+2)(X^2+3)=X^4+5X^2+6=X^4+1\).
Every \(a\in\F_5^\times\) has \(a^4=1\ne-1\), so \(-1\)
has no fourth root in the base field. In a splitting field, choose
\(\omega^4=-1\) and write any root of \(f\) as \(\omega\zeta^b\).
Corollary~\ref{v88:cor-constacyclic-bch-bound} with \(\delta=2\) gives
\(d_H(C)\ge2\). The generator \(f=X^2+2\) itself is a codeword
of weight two, proving equality.
\end{solution}

\solutionheading{v83:ex-interpolated-split-fiber}
\begin{solution}
Since \(K\) has characteristic two, squaring is a bijection. Thus
\(P_V(a)=U(a)^2+aV(a)^2=0\) is equivalent to \(V(a)=U(a)/\sqrt a\).
At \(s\) distinct prescribed points, evaluation on polynomials of degree less than \(s\) gives an isomorphism
\[
 K[Y]_{<s}\longrightarrow K^A,\qquad V\longmapsto(V(a))_{a\in A}
\]
by Lagrange interpolation. Hence there is exactly one such \(V\).
Now \(L_A\mid P_V\), and \(R_A=P_V/L_A\) is monic of degree \(s\).
Since \(A\subseteq H\), we have \(L_A\mid Y^h-1\), and \(Y^h-1\) is square-free. Consequently,
\[
 P_V\mid Y^h-1
 \quad\Longleftrightarrow\quad
 R_A\mid Y^h-1\ \text{and}\ \gcd(R_A,L_A)=1.
\]
The latter condition puts the remaining \(s\) roots in \(H\) and keeps them distinct from the specified roots.
\end{solution}

\solutionheading{v83:ex-sixteen-coset-fibers}
\begin{solution}
The group \(H\) is cyclic of order \(5\).
If \(P_v\) has two distinct roots \(y,z\in H\), their product is \(yz=u^2\in H\).
Since \(|K^\times/H|=3\) is odd, we must have \(u\in H\).
Conversely, if \(u\in H\), each \(y\in H\) determines \(z=u^2/y\in H\).
The two roots coincide precisely when \(y^2=u^2\), or \(y=u\).
There are therefore \((5-1)/2=2\) unordered pairs of distinct roots. Each pair determines \(v^2=y+z\), and hence a unique \(v\).
No such pair exists when \(u\notin H\).
For \(u=1\), Proposition~\ref{v83:prop-derivative-fibers} says that \(F_v\) is square-free exactly when \(v\ne0\), yielding fifteen polynomials.
The polynomial \(G=X^3+1\) has three roots in \(K\), so its full support has length \(16-3=13\).
The two split locator polynomials yield codewords of weight seven, while Proposition~\ref{v82:prop-goppa-support} gives the matching lower bound. Thus the minimum distance is seven.
\end{solution}

\solutionheading{v44:ex-interleaving-intrinsic}
\begin{solution}
If \(f=F(X^d)\), the standard representative \(c=\sum_{j=0}^{d-1}X^j a_j(X^d)\) is divisible by \(f\)
if and only if every \(a_j\) is divisible by \(F\). Thus each component \(X^j a_j(X^d)\) still belongs to \(C\),
and every element of \(C\) decomposes uniquely into a sum of elements of the \(C\cap V_j^{(d)}\), proving that the first condition implies the third.
Under the third condition these components are precisely the coordinate projections, so the third implies the second.

Conversely, suppose the second holds. By \(\deg f<N\), the polynomial \(f\) is a nonzero standard representative in \(C\).
Also \(f(0)\ne0\), so \(\pi_0^{(d)}(f)\ne0\). It belongs to \(C\) and is therefore divisible by \(f\);
its degree is at most \(\deg f\), so it must equal \(cf\) with \(c\in k^\times\).
Comparing constant terms gives \(c=1\), so \(\pi_0^{(d)}(f)=f\), namely \(f\in k[X^d]\).
This proves equivalence of the three conditions.

By the definition of the support period, \(f\in k[X^d]\) if and only if \(d\mid\operatorname{sp}(f)\).
Proposition~\ref{v44:prop-interleaving} also gives \(\operatorname{sp}(f)\mid N\), so the largest possible parameter is exactly \(\operatorname{sp}(f)\).
\end{solution}

\solutionheading{v54:ex-effective-support}
\begin{solution}
If \(z\) is a root of \(Y^6+Y^3+1\), then \(z^9=1\) but \(z^3\ne1\):
otherwise \(z^6+z^3+1=1\). Thus \(z\) has order \(9\), belongs to \(\F_{64}\), and satisfies \(z\ne1\).
The corresponding root \(x=z+1\) of \(f\) is nonzero, so \(x^{63}=1\).
Also \((Y^6+Y^3+1)'=Y^2\) is nonzero at these roots, so \(f\) is square-free, giving \(f\mid X^{63}-1\).
In the outer variable \(Y=U^3\), the generator polynomial is \(B=Y^2+Y+1\).
Using \(1,Y\) as the remainder basis, the remainder matrix of \(1,Y,Y^2\) is
\[
M=\begin{pmatrix}1&0&1\\0&1&1\end{pmatrix}.
\]
Its rank is two and its kernel is generated by \((1,1,1)\), so \(C\cap W=\F_2f\) has dimension one.
Expansion gives \(f=X^6+X^4+X^3+X+1\), with effective support \(\{0,1,3,4,6\}\).
The common support of \(W\) is \(\{0,1,2,3,4,6\}\). The quadratic-coefficient row is \((0,1,1)\),
which belongs to the row space of \(M\) and therefore vanishes on \(\ker M\).
The unique nonzero codeword in the subcode has weight five, so \(d_H(C)\le5\);
this does not exclude smaller weights for other codewords in \(C\), and hence does not imply \(d_H(C)=5\).
\end{solution}

\solutionheading{v73:exer-independent-multipliers}
\begin{solution}
Ordering the CRT coordinates modulo three and five gives
\[
K=\{1,2,4,8\}
\longleftrightarrow\{(1,1),(2,2),(1,4),(2,3)\}.
\]
The projections are \(\{1,2\}\) and \(\{1,2,3,4\}\), the full respective unit groups.
The independent local groups, however, are \(L_3=\{1\}\) and \(L_5=\{1,4\}\).
Fixing the coordinate modulo five leaves only the identity multiplier; fixing the coordinate modulo three permits only even iterates.

The modulus is squarefree, so only the first exception in the rigidity theorem could occur, and its only possible prime is five.
But \(L_5\ne\F_5^\times\), excluding that case too.
Theorem~\ref{v73:thm-ring-rigidity} therefore gives precisely
\[
f(x)=2^jx+b,\qquad 0\le j<4,\quad b\in\Z/15\Z,
\]
numbering \(60\).
A projection records which multipliers can appear in one coordinate, but forgets the synchronization constraint with the other coordinate.
This is why it cannot replace the independent local group.
Finally, Theorem~\ref{v73:thm-common-permutation-orbits} identifies the same maps as the common permutation group of all binary length-fifteen cyclic codes.
\end{solution}

\solutionheading{v73:exer-family-versus-common}
\begin{solution}
From \((\mu_uc)_i=c_{u^{-1}i}\), the substitution \(i=uj\) gives
\[
\widehat{\mu_uc}(a)=\sum_i c_{u^{-1}i}\zeta^{ai}
=\sum_jc_j\zeta^{auj}=\widehat c(ua).
\]
Thus the spectral support becomes \(u^{-1}S(C)\). This remains \(q\)-stable, so the image is again a cyclic code.
Alternatively, \(\mu_uR_n\mu_u^{-1}=R_n^u\) directly verifies shift invariance.

The two nonzero binary cyclotomic orbits modulo seven are \(O_1=\{1,2,4\}\) and \(O_3=\{3,5,6\}\).
Since \(3^{-1}=5\pmod7\), we have \(5O_1=O_3\) and \(5O_3=O_1\).
Hence \(\mu_3\) exchanges the irreducible cyclic codes supported on these orbits instead of preserving them separately.
It does not belong to \(\mathcal G_2(7)\), consistently with \(\mu_3(0)=0\) and \(\mu_3(1)=3\notin\langle2\rangle\).
\end{solution}

\solutionheading{v82:ex-cyclic-covering}
\begin{solution}
Suppose a proper subspace \(U\) is cyclically covering. Choose nonzero \(w\in V\) with \(U\subseteq H=w^\perp\); then \(H\) is cyclically covering too.
Define
\[
 \Phi_w(x)=\bigl(\langle x,R_N^iw\rangle\bigr)_{i=0}^{N-1}.
\]
The standard coordinate pairing is invariant under cyclic shift, so \(\Phi_w(R_Nx)=R_N\Phi_w(x)\), making \(\im\Phi_w\) a cyclic code.
It is nonzero because \(x\mapsto\langle x,w\rangle\) is a nonzero linear functional.
For each \(x\), the covering property gives some \(i\) with \(x\in R_N^iH\), whence
\(\langle x,R_N^iw\rangle=\langle R_N^{-i}x,w\rangle=0\).
Every image vector therefore has a zero coordinate, so this nonzero cyclic code has no full-weight word.

Conversely, let \(C=(g)\subseteq k[X]/(X^N-1)\) be a nonzero cyclic code without a full-weight word.
Under the coefficient identification, multiplication by the generator polynomial \(g\) gives a surjection \(\Psi:V\to C\) commuting with \(R_N\).
Let \(\lambda(x)\) be the zeroth coordinate of \(\Psi(x)\).
This functional is nonzero: a nonzero coordinate of a nonzero word in \(C\) can be shifted to position zero, and surjectivity of \(\Psi\) supplies a preimage.
Thus \(H=\ker\lambda\) is a proper hyperplane.
For any \(x\in V\), the vector \(\Psi(x)\) has a zero coordinate \(j\), so
\[
 \lambda(R_N^{-j}x)
 =\bigl(R_N^{-j}\Psi(x)\bigr)_0
 =\Psi(x)_j=0.
\]
Hence \(x\in R_N^jH\), proving that \(H\) is cyclically covering.
No step assumes that \(X^N-1\) is square-free, so the argument also applies to repeated-root lengths.
\end{solution}

\solutionheading{v84:ex-projective-basis}
\begin{solution}
At \((b,1)\), the value \(F_a=1-(b-a)^{q-1}\) is one for \(b=a\) and zero otherwise, while \(F_\infty=b^q-b=0\).
At \((1,0)\), every \(F_a\) is zero and \(F_\infty=1\).
Thus the interpolant of the prescribed values \((z_a)_{a\in k},z_\infty\) is
\(\sum_{a\in k}z_aF_a+z_\infty F_\infty\).
Uniqueness follows from Proposition~\ref{v84:prop-projective-evaluation}.
Although \(F_\infty\) vanishes at every affine point, it takes the value one at infinity.
It records precisely the information retained by projective evaluation from an affine vanishing relation.
\end{solution}

\solutionheading{v84:ex-ternary-kernel}
\begin{solution}
Since \(a_i^q=a_i\) and \(b_i^q=b_i\) in \(\F_q\), we have
\[
 L^qM-LM^q=\sum_{i,j}(a_ib_j-b_i a_j)X_i^qX_j.
\]
The diagonal terms vanish. Pairing the indices \((i,j)\) and \((j,i)\) gives the stated identity.
By Theorem~\ref{v86:thm-frobenius-kernel}, every such expression belongs to the critical kernel.
After a change of linear coordinates, each new variable is a linear form, so the identity expresses every new generating binomial as a linear combination of the old ones. Applying the same argument to the inverse coordinate change proves invariance of the kernel space.
If \(L,M\) are linearly independent, some \(2\times2\) coefficient minor \(a_ib_j-a_jb_i\) is nonzero; linear independence of the generators therefore makes \(L^qM-LM^q\ne0\).
In the example over \(\F_5\), the three minors are \(1,1,2\), and hence
\[
 L^5M-LM^5
 =(X_1^5X_2-X_1X_2^5)+(X_1^5X_3-X_1X_3^5)
  +2(X_2^5X_3-X_2X_3^5).
\]
This is a nonzero homogeneous form of degree six, so it cannot be a fifth power.
\end{solution}

\solutionheading{v84:ex-dimension-two}
\begin{solution}
Two distinct elements \(x,y\in H=\langle\alpha\rangle\) are dependent exactly when
\(xy^{-1}\in k^\times\setminus\{1\}\).
In a cyclic group, \(|H\cap k^\times|=\gcd(N,q-1)\), so the two-dimensional column criterion is precisely the stated coprimality condition.
Since \(N\mid(q-1)(q+1)\), coprimality gives \(N\mid q+1\) and excludes even \(N\).
Conversely, if \(N>1\) is an odd divisor of \(q+1\), then \(\gcd(N,q-1)=1\), and \(q\equiv-1\pmod N\) has order exactly two.
The subgroup therefore gives a two-dimensional irreducible MDS code.
Nontriviality requires \(2\le N-2\), which together with oddness is \(N\ge5\).

For \(q=4,N=5\), one has \(\alpha^4=\alpha^{-1}\), hence \(\omega\in\F_4\).
Divide \(1+\alpha+\alpha^2+\alpha^3+\alpha^4=0\) by \(\alpha^2\) to obtain \(\omega^2+\omega+1=0\).
The check polynomial is \(X^2+\omega X+1\), and the code has parameters \([5,2,4]\).
The traces of nonidentity subgroup elements are \(\omega\) or \(\omega^2\), both nonzero.
For a quadratic extension in characteristic two, the established normality criterion therefore makes all these elements normal.
The same conclusion follows directly from Theorem~\ref{v84:thm-subgroup-rigidity}.
\end{solution}

\solutionheading{v84:ex-low-degree-multiplier}
\begin{solution}
Set \(a_i=0\) outside the given range. In coefficient coordinates, multiplication \(B\mapsto m_\alpha B\) has matrix
\[
T_s=(a_{i-j})_{\substack{0\le i\le K+s\\0\le j\le s}}.
\]

For \(B=\sum_{j=0}^s b_jX^j\), the coefficient vector of \(m_\alpha B\) in degrees \(0,\ldots,K+s\) is
\(T_s(b_0,\ldots,b_s)^{\mathsf T}\).
Requiring \(s+1\) selected coefficients to vanish admits nonzero parameters exactly when the corresponding square matrix is singular.
The product \(m_\alpha B\) is nonzero, and at least \(s+1\) of its \(K+s+1\) coefficient positions vanish, so its weight is at most \(K\).
Conversely, such a multiple of weight at most \(K\) has at least \(s+1\) zero coefficient positions in this range, and its parameter vector is a nonzero kernel vector for the corresponding row minor.
Every nonzero multiple of degree less than \(N\) has \(\deg B\le N-K-1\).
Applying Theorem~\ref{v84:thm-sparse-multiple} over all \(s\) gives the desired equivalent criterion.
In the example, take \(s=2\) and rows \(1,2,3\). Their submatrix and kernel vector are
\[
 \begin{pmatrix}1&2&0\\1&1&2\\0&1&1\end{pmatrix},
 \qquad
 \begin{pmatrix}1&2&0\\1&1&2\\0&1&1\end{pmatrix}
 \begin{pmatrix}2\\2\\1\end{pmatrix}=0\quad\text{over }\F_3.
\]
This corresponds to \(B=X^2+2X+2\) and \(m_\alpha B=X^4+1\).
For \(s=0\), all three coefficients of \(m_\alpha=X^2+X+2\) are nonzero. For \(s=1\), the four rows of \(T_1\) are \((2,0),(1,2),(1,1),(0,1)\), and any two are linearly independent. Thus no such multiple exists, and the least multiplier degree is two.
\end{solution}

\solutionheading{v142:ex-partial-fibers}
\begin{solution}
We have $\ord_{21}(2)=6$ and $\ord_7(2)=3$. Each fiber of
\[
 C_{21}(1)=\{1,2,4,8,16,11\}\longrightarrow C_7(1)=\{1,2,4\}
\]
has only two elements, whereas a complete fiber has size three. Thus $6<3\cdot3$, and Proposition~\ref{v142:prop-full-fibers} excludes $h=h_0(X^3)$, where $h_0=m_{\xi^3}$.
More generally, any identity $h=G(X^3)$ would give $G(\xi^3)=0$ with $\deg G=2<\deg h_0=3$, again a contradiction.

The element $\xi^3$ has order seven and generates $\F_8/\F_2$. Proposition~\ref{v142:prop-subfield-distance} gives dimension three for both the restriction image and its kernel, as well as
\[
 d_H(C_\xi^\perp)\le4,\qquad
 d_H(C_\xi)\le21-6+1-(7-3)=12<16.
\]
Hence the original code is not MDS. The restriction image is in fact a $[7,3,4]$ code: a nonzero trace functional has four zeros in $\F_8$, three of them nonzero, and therefore has weight four on the seven nonzero elements.
The corresponding seven-coordinate subcode of the dual also contains a word of weight three: for distinct $x,y\in\F_8^\times$, the elements $x,y,x+y$ give a three-term relation. Two distinct nonzero elements cannot be dependent over $\F_2$, so $d_H(C_\xi^\perp)=3$.
If the original code were a coordinate direct sum of three copies of its restriction image, it would have dimension nine rather than six. The subfield supplies short relations; complete fibers are needed to make the three blocks' trace parameters independent.
\end{solution}

\solutionheading{v84:ex-fifteen-obstruction}
\begin{solution}
Since \(\ord_{15}(2)=4\), the dimension is four and the support is \(S=\{1,2,4,8\}\).
Its three residue classes modulo three have sizes \(0,2,2\), violating the lower bound \(\lfloor4/3\rfloor=1\).
On the subgroup side, \(\gamma=\alpha^5\) has order three and belongs to \(\F_4\), whose dimension over \(\F_2\) is only two.
Thus \(1,\gamma,\gamma^2\) are dependent; explicitly, \(1+\gamma+\gamma^2=0\).
Equivalently, \(1+X^5+X^{10}\) vanishes at \(\alpha\), has degree less than fifteen and weight three, and violates the sparse-multiple criterion.
The three descriptions record the same short relation in frequency layers, subfield dimension, and polynomial coefficients.
\end{solution}

\solutionheading{v84:ex-extension-splitting}
\begin{solution}
Scalar extension preserves parameters, so \(C_{16}\) is still an \([5,2,4]\) MDS code.
The original check polynomial \(X^2+\omega X+1\) splits over \(\F_{16}\) as \((X-\alpha)(X-\alpha^{-1})\), so the code is no longer irreducible.
With evaluation at powers of \(\alpha\), the support is still \(S=\{1,4\}\subseteq\Z/5\Z\).
Multiplication by two sends it to \(\{2,3\}\), whereas multiplication by four preserves it.
Its set return period is therefore two, and its smallest field of definition is \(\F_4\).
The formula \(\gcd(\ord_5(2),4)=4\) applies to a single faithful irreducible orbit.
Here the support has split into two singleton orbits under the sixteenth-power Frobenius, so that special formula does not apply.
The general set-period criterion still applies and correctly gives two.
\end{solution}

\solutionheading{v84:ex-arithmetic-orders}
\begin{solution}
If \(q\equiv1\pmod N\), then \(N\mid q-1\), proving the first part immediately.
Otherwise, under \(N\mid q^2-1\), the order of \(q\) modulo \(N\) is two.
Put \(\nu=\gcd(q-1,N)\); the multiplier fixed-point restriction gives \(\nu\le2\).
Write \(N=\nu N_0\), \(q-1=\nu c\). Since \(\gcd(N_0,c)=1\), divisibility gives \(N_0\mid q+1\), and hence \(N\mid2(q+1)\).
The elementary bound \(N\le2q-2<2(q+1)\) makes \(N\) a proper divisor of that integer, so \(N\le q+1\).

If \(3\mid N\) or \(4\mid N\), synchronization of orders gives \(\ord_N(q)\mid2\), and the first part applies.
For the two odd prime divisors in the question, \(\ord_N(q)\) divides both \(\ell_1-1\) and \(\ell_2-1\), so again it is at most two.
Finally put \(d=\gcd(N,q^2-1)\).
If \(d>2\), the order modulo \(d\) is at most two; synchronization gives \(N\mid q^2-1\).
If \(K\) is even and \(d>1\), Corollary~\ref{v84:cor-even-length-parity} excludes \(d=2\), reducing to the preceding case.
\end{solution}

\solutionheading{v84:ex-gluing-boundaries}
\begin{solution}
At \(d=N\), the nonzero homogeneous form \(L_1\cdots L_N\) vanishes on every \(H_i\).
Even zero data therefore fail to determine a unique form.
In the second example, all three pairwise intersections are \(\{X=Y=0\}\), on which \(g_3=X\) is zero, so compatibility holds.
If a homogeneous linear \(G\) vanished on both \(H_1,H_2\), then \(X\mid G\) and \(Y\mid G\), hence \(XY\mid G\), forcing \(G=0\).
But \(X|_{H_3}\ne0\), a contradiction.
The forms \(X,Y,X-Y\) are dependent, violating the hypothesis that every three are independent.
Consequently, the intersections on a new section no longer provide distinct linear factors.
\end{solution}

\solutionheading{v116:ex-five-completion}
\begin{solution}
The polynomial \(F_0\) has five distinct roots, so \(2\deg H_A=\deg F_0\ge5\), while interpolation gives \(\deg H_A\le4\).
If the degree is three, \(F_0/R_A\) is linear and splits. Since \(F_0'=1\), its root is distinct from the prescribed roots, giving completion size six.
If the degree is four, \(S=F_0/R_A\) is a squarefree cubic coprime to \(R_A\). An eight-point completion exists exactly when
\[
\deg\gcd(S,X^q-X)=3.
\]
Otherwise the degree gap gives \(c_q(A)\ge10\). For \(q=16\), Corollary~\ref{v116:cor-completion-size} rules out completion sizes \(10,12,14\), so the full field is required and the minimum is sixteen.
Thus test \(\deg H_A=3\) first, then the cubic splitting condition; all remaining cases have value sixteen.
\end{solution}

\Needspace{9\baselineskip}
\section*{Solutions to Chapter~\ref{chap:characters}}
\pdfbookmark[1]{Solutions to Chapter \getrefnumber{chap:characters}}{solutions.chap:characters}

\solutionheading{v111:ex-low-degree-power}
\begin{solution}
Write $d=\deg F$. First suppose $F=cH^v$, with $c\in\F_q^\times$ and $H\in\F_q[X]$.
Since $d<q$, some $x_0\in\F_q$ has $F(x_0)\ne0$. Hence $c=F(x_0)/H(x_0)^v\in(\F_q^\times)^v$.
Writing $c=b^v$ gives $F=(bH)^v$.
This also handles a power over the algebraic closure: finite fields are perfect, so every root multiplicity of $F$ is divisible by $v$ if and only if every exponent in its monic irreducible factorization over $\F_q$ is divisible by $v$. In that case $F=cH^v$.

If $F$ is not a $v$th power, the preceding argument excludes every $cH^v$.
Take a multiplicative character $\chi$ of order $v$, with $\chi(0)=0$. Let $z$ be the number of zeros in $\F_q$, and $\rho$ the number of distinct roots in a splitting field.
The power-value condition and the Weil bound give
\[
q-d\le q-z=\sum_{x\in\F_q}\chi(F(x))
\le(\rho-1)\sqrt q\le(d-1)\sqrt q.
\]
Thus $q+\sqrt q\le d(\sqrt q+1)$, or $d\ge\sqrt q$, a contradiction.
No minimal-value-set hypothesis was used: the inputs are the power-value restriction and the ordinary degree bound.

For $q=9$, the map $x^3+x$ is the relative trace from $\F_9$ to $\F_3$, hence is surjective with three elements in every fiber.
The fourth powers in $\F_9^\times$ form its unique subgroup of order two, namely $\F_3^\times$. Therefore
\[
F(\F_9)=\F_3=\{0\}\cup(\F_9^\times)^4,
\qquad |\mathcal V_F|=3=\left\lceil\frac93\right\rceil.
\]
But $\deg F=3$ is not divisible by four, so $F$ is not a fourth power.

More generally, $x^s+x=\Tr_{\F_{s^2}/\F_s}(x)$.
Theorem~\ref{thm:trace-norm-surjective} gives image $\F_s$, a kernel of size $s$, and size $s$ for every fiber.
Also $(\F_{s^2}^\times)^{s+1}=\F_s^\times$, so the power-value condition holds and
$|\mathcal V_F|=s=\lceil s^2/s\rceil$. Its ordinary degree $s$ is smaller than $s+1$, excluding an $(s+1)$st power decomposition as a polynomial.
This includes characteristic two. Choosing a root of each function value does not supply a polynomial root of controlled degree, nor does it make the original polynomial a $v$th power in $\F_q[X]$.
Thus a value set that is itself a subfield cannot, merely because it is a multiplicative power image, force a polynomial power decomposition with the originally specified exponent. Compare the actual hypotheses of Theorem~\ref{v109:thm-mvsp-classification}.

For the strengthening, first suppose \(\rho<q\). Since \(z\le\rho\), a nonzero value exists; the argument with \(cH^v\) therefore applies without any bound on the ordinary degree. If \(F\) is not a \(v\)th power, the multiplicative Weil bound gives \(q-z\le(\rho-1)\sqrt q\). If instead \(\rho\ge q\), this inequality follows directly from \(q-z\le q\le(q-1)\sqrt q\); here \(v\ge2\) implies \(q\ge3\). Thus the asserted inequality holds in either case.
In particular, \(q-\rho\le(\rho-1)\sqrt q\) forces \(\rho\ge\sqrt q\). This proves the stronger assertion. The derivative of \(X^s+X\) is one, so it has exactly \(s\) distinct roots in a splitting field and attains the boundary.
The bound counts distinct roots over a splitting field, not merely zeros in the base field; it does not discard the root multiplicities when deciding whether the polynomial is a \(v\)th power.
\end{solution}

\solutionheading{v88:ex-primitive-nonzero-bound}
\begin{solution}
For \((q,n)=(2,3)\), \(Q=8\), \(M=R=7\),
\(\vartheta=6/7\), and \(W(R)=2\). The lower bound is
\(\frac37(8-\sqrt8)>0\).
For \((q,n)=(3,2)\), \(Q=9\), \(M=8\), \(R=4\),
\(\vartheta=1/2\), and \(W(R)=2\); the bound is
\(\frac13(9-3)=2>0\).
Summing the nonzero trace fibers avoids estimating their individual counts.
\end{solution}

\solutionheading{v58:ex-cubic-sixth-moment}
\begin{solution}
By \eqref{v58:eq-moment-power-systems}, the left-hand side is \(Q^3\) times the number of ordered pairs of ordered triples with equal first three power sums.
Proposition~\ref{v58:prop-power-sum-multiset} shows that the two triples must be permutations of the same multiset.
There are \(\binom{Q}{3}\) multisets of three distinct elements. Each has six permutations and contributes \(36\) ordered pairs, giving \(36\binom{Q}{3}\) pairs in total.
There are \(Q(Q-1)\) multisets with multiplicity type \(2+1\): choose the repeated element first, then the other element.
Each has three permutations and contributes \(9\) ordered pairs, giving \(9Q(Q-1)\) pairs in total.
There are \(Q\) multisets with three equal elements, each with only one permutation.
Thus the total number of pairs is
\[
36\binom{Q}{3}+9Q(Q-1)+Q=6Q^3-9Q^2+4Q,
\]
so multiplication by \(Q^3\) gives the required identity.
\end{solution}

\solutionheading{ex:fourier-delta}
\begin{solution}
The transform of a point-indicator function is \(\widehat{1_{\{b\}}}(a)=\psi(-ab)\).
The transform of the constant function one is \(Q\) at \(a=0\) and zero elsewhere.
Substitution of the first identity into the inversion formula gives \(1_{\{b\}}(x)=Q^{-1}\sum_a\psi(a(x-b))\), exactly the additive orthogonality relation.
\end{solution}

\solutionheading{ex:warning-quadratic}
\begin{solution}
Both sets \(A=\{ax^2:x\in k\}\) and \(B=\{1-by^2:y\in k\}\) have \((q+1)/2\) elements.
Their sizes sum to \(q+1>q\), so they intersect; an intersection point gives \(ax^2+by^2=1\).
The Chevalley--Warning theorem guarantees a nonzero solution of the homogeneous ternary equation \(aX^2+bY^2-Z^2=0\).
This existence statement alone does not ensure \(Z\ne0\) for the resulting solution, so one cannot directly divide by \(Z^2\) to obtain the prescribed affine equation.
The intersection argument for sets of squares gives a more specific conclusion here.
\end{solution}

\solutionheading{ex:quadratic5}
\begin{solution}
The element \(-1=4\) is a square in \(\F_5\). The rank-two quadratic-form count gives \(5-1=4\) solutions of \(x^2+y^2=1\) and \(5+4=9\) solutions in the zero fiber.

Put \(K=\F_q(\sqrt a)=\F_{q^2}\). Since \(a\) is a nonsquare, the other conjugate of \(\sqrt a\) is \(-\sqrt a\), so
\[
\Nm_{K/\F_q}(x+y\sqrt a)=x^2-ay^2.
\]
The norm is zero exactly when its argument is zero; hence when \(b=0\) the only solution is \((0,0)\).
The norm map on nonzero elements is surjective, with kernel of size \((q^2-1)/(q-1)=q+1\);
each nonzero fiber is a coset of this kernel, so every \(b\ne0\) has \(q+1\) solutions.

On the other hand, the matrix of \(Q(x,y)=x^2-ay^2\) has determinant \(-a\) and rank two,
and its type sign in Theorem~\ref{thm:quadratic-odd-count} is \(\chi((-1)(-a))=\chi(a)=-1\). Thus
\[
N_Q(b)=q-v(b)=
\begin{cases}1,&b=0,\\q+1,&b\ne0,\end{cases}
\qquad
v(b)=\begin{cases}q-1,&b=0,\\-1,&b\ne0.\end{cases}
\]
An anisotropic binary quadratic form is precisely the norm form of a quadratic field extension; hence \(q+1\) counts both the elements in a nonzero norm fiber and the representations of each nonzero value by this quadratic form.
\end{solution}

\solutionheading{ex:correlation7}
\begin{solution}
We have \(a_{t+\tau}a_t=(-1)^{s_{t+\tau}+s_t}\), and addition and subtraction coincide in \(\F_2\).
Using the trace representation, rewrite the sum as
\[
\sum_{x\in\F_8^\times}(-1)^{\Tr((\alpha^\tau-1)x)}.
\]
If seven divides \(\tau\), every term is one and the sum is seven; otherwise this is the sum of a nontrivial additive character over the nonzero elements, equal to minus one.
\end{solution}

\solutionheading{ex:wu-fourier-subspace}
\begin{solution}
By definition, \(\widehat{1_U}(a)=\sum_{u\in U}\psi(-au)\).
If \(a\in U^\perp\), every term is one and the sum is \(|U|\).
Otherwise there is a \(u_0\in U\) such that \(\psi(-au_0)\ne1\);
translating the summation index by \(u_0\) multiplies the sum by this scalar different from one without changing its value, so the sum is zero.
For a coset, the substitution \(x=b+u\) gives
\[
\widehat{1_{b+U}}(a)=\psi(-ab)|U|\,1_{U^\perp}(a).
\]
\end{solution}

\solutionheading{ex:wu2-fourier-change}
\begin{solution}
Changing variables by \(y=A(x)\) gives \(\Tr(ax)=\Tr(aA^{-1}(y))=\Tr((A^*)^{-1}(a)y)\), and substitution into the definition proves the formula.
If \(K=\F_{p^m}\) and \(A=\sigma:x\mapsto x^p\), then \(A^*=\sigma^{-1}\), so the frequency parameter becomes \(a^p\).
This also shows that the actions of a linear change of variables on functions and frequencies are linked by adjointness under the trace pairing.
\end{solution}

\solutionheading{v70:ex-quotient-phase}
\begin{solution}
The kernel-image orthogonality relation for adjoints gives \(\im L^*=(\ker L)^\perp\). If \(b\notin\im L^*\), choose \(w\in\ker L\) such that \(\psi(bw)\ne1\). Replacing \(x\) by \(x+w\) in the sum and using \(L(x+w)=Lx\) gives \(S(b)=\psi(bw)S(b)\), hence \(S(b)=0\).

If \(b=L^*c\), then \(\Tr_{K/\F_p}(bx)=\Tr_{K/\F_p}(cLx)\). Each \(y\in U\) has \(|\ker L|\) preimages, so
\[
S(b)=|\ker L|\sum_{y\in U}\psi(g(y)+cy).
\]
Also \(U=W^\perp\), and additive-character orthogonality gives
\[
1_U(y)=\frac1{|W|}\sum_{t\in W}\psi(ty).
\]
The maps \(L\) and \(L^*\) have the same rank, hence \(|W|=|\ker L|\). Substitute this indicator expansion into the preceding expression and interchange the finite sums to obtain the second identity.
If \(L^*c'=L^*c\), then \(c'-c\in W\), so \(c'+W=c+W\), and the last sum is unchanged. Since \(W\) is orthogonal to \(U\), the first sum is also unchanged.
\end{solution}

\solutionheading{ex:wu3-trace-square-odd}
\begin{solution}
In any basis \((b_i)\), the coefficient matrix of \(Q\) is \((\Tr(b_i b_j))\), invertible by nondegeneracy of the trace pairing.
Let the conjugate matrix be \(V=(b_j^{q^i})_{0\le i,j<n}\). The coefficient matrix is then \(V^{\mathsf T}V\), so \(\Delta=(\det V)^2\).
Raising every entry of \(V\) to the \(q\)th power cyclically permutes its rows, hence \((\det V)^q=(-1)^{n-1}\det V\).
Therefore \(\eta(\Delta)=\Delta^{(q-1)/2}=(\det V)^{q-1}=(-1)^{n-1}\).
If \(n\) is even, the number of solutions is
\[
q^{n-1}-\eta(-1)^{n/2}v(b)q^{n/2-1}.
\]
If \(n\) is odd, the number of solutions is
\[
q^{n-1}+\eta(-1)^{(n-1)/2}\eta(b)q^{(n-1)/2}.
\]
When \(q=3,n=2\), zero has five preimages and every nonzero value has two.
\end{solution}

\solutionheading{ex:wu3-trace-square-even}
\begin{solution}
We have \(Q(x)=\Tr(x)^2\), so the map is additive and its polar form is zero.
The trace is surjective and squaring is bijective on \(\F_q\), so \(Q\) is surjective.
The equation \(Q(x)=b\) is equivalent to \(\Tr(x)=b^{q/2}\), whose fibers have size \(q^{n-1}\).
In Theorem~\ref{thm:quadratic-even-count}, this is the case \(R=E\) with \(Q|_R\ne0\).
When \(q=2\), we have \(Q(x)=\Tr(x)\), so
\[
\widehat F_Q(b)=\sum_{x\in E}(-1)^{\Tr((1+b)x)}
=2^n1_{\{1\}}(b).
\]
The support is \(\{1\}\), whereas \(R^\perp=\{0\}\). They differ by translation by one,
exactly the affine coset determined by \(Q|_R\ne0\) in Proposition~\ref{prop:quadratic-fourier-support}.
\end{solution}

\solutionheading{ex:wu3-quadratic-linear-shift}
\begin{solution}
Example~\ref{exm:f64-quadratic} gives \(\widehat F_Q(0)=16\) and the support
\(\im A=(\F_4)^\perp=\ker\Tr_{\F_{64}/\F_4}\).
If \(b=Az\), Proposition~\ref{prop:quadratic-fourier-support} gives
\(\widehat F_Q(b)=16(-1)^{Q(z)}\), where \(Q(z)\) is independent of the chosen solution.
The number of solutions is therefore
\[
N_{Q,b}(c)=32+8(-1)^{c+Q(z)}.
\]
It is 40 when \(c=Q(z)\) and 24 otherwise.
If \(b\notin\im A\), then \(\widehat F_Q(b)=0\). By \(N_{Q,b}(c)=(64+(-1)^c\widehat F_Q(b))/2\), each of the two values has 32 preimages.
Here the linear-term parameter first determines whether the frequency lies in the support; translation of the summation variable then determines the sign within the support.
\end{solution}

\solutionheading{ex:wu3-jacobi-boundary}
\begin{solution}
The case where all three characters \(\chi,\lambda,\chi\lambda\) are nontrivial has already been proved.
If \(\lambda=\chi^{-1}\) and \(\chi\ne\varepsilon\), then
\[
J_{E_n}(\chi_n,\chi_n^{-1})=-\chi_n(-1)=-\chi(-1)^n
=(-1)^{n-1}(-\chi(-1))^n.
\]
Here we use \(\Nm_{E_n/k}(-1)=(-1)^n\).
If exactly one character is trivial, the Jacobi sums over both the base and extension fields are minus one, and \((-1)^{n-1}(-1)^n=-1\), so the identity still holds.
If both characters are trivial, the respective values are \(q-2\) and \(q^n-2\) and generally do not satisfy the lifting identity; for example, \(q=2,n=2\).
\end{solution}

\solutionheading{ex:wu3-diagonal-seven}
\begin{solution}
The only nonzero cubes are \(1,-1\), each with three preimages; zero has only one.
Thus the solution counts for \(c=0,\pm1,\pm2,\pm3\) are \(19,6,9,0\), respectively, with a total of 49 over all seven values of \(c\).
With three as a generator, the discrete logarithms \(0,1,\ldots,5\) correspond to the nonzero elements \(1,3,2,6,4,5\), respectively. In \(J(\chi,\chi)\), the values \(x=0,1\) contribute nothing;
the contributions of \(x=2,3,4,5,6\) are \(\zeta^2,1,\zeta^2,1,\zeta^2\), respectively, so the sum is \(2+3\zeta^2\).
The conjugate sum is \(2+3\zeta\), and both Jacobi sums of inverse characters are minus one, since \(\chi(-1)=1\).
For \(c\ne0\), substitution into the solution-count formula gives
\[
5+\chi(c)^2(2+3\zeta^2)+\chi(c)(2+3\zeta),
\]
which takes the values six, nine, and zero when \(\chi(c)=1,\zeta^2,\zeta\), respectively, in agreement with direct counting.
For zero, the second identity in the proposition gives \(7+6\cdot2=19\).
\end{solution}

\solutionheading{ex:character-L-boundary}
\begin{solution}
\textup{(1)} Over \(\F_2\), we have \(S_1(f)=1-1=0\). Over \(\F_4\), every nonzero element satisfies \(x^3=1\),
and \(\Tr_{\F_4/\F_2}(1)=0\), so \(S_2(f)=4\).
Theorem~\ref{thm:additive-L-polynomial} gives \(\deg L_f\le2\),
and \eqref{eq:character-initial-newton} further gives \(\ell_1=0\) and \(2\ell_2=4\), so
\[
L_f(u)=1+2u^2,\qquad
S_{2j+1}(f)=0\ (j\ge0),\qquad
S_{2j}(f)=-2(-2)^j\ (j\ge1).
\]
The trace takes only the values zero and one, so the required count is \((2^r+S_r(f))/2\). When \(r=6\), we have \(S_6(f)=16\) and the count is 40,
in agreement with Example~\ref{exm:f64-quadratic}.

\textup{(2)} For every \(x\in\F_{2^r}\), Frobenius invariance of the trace gives
\(\Tr(x^2+x)=\Tr(x)^2+\Tr(x)=0\).
Thus \(S_r(g)=2^r\), and
\[
L_g(u)=\exp\left(\sum_{r\ge1}\frac{(2u)^r}{r}\right)=\frac1{1-2u}.
\]
This formal power series is not a polynomial. The absolute trace of \(g(x)\) is zero over every extension field, so the additive-character phase is identically one.
This illustrates the role of the degree condition in the coefficient-cancellation argument: the coefficient \(-d a_d\) of the leading term must be nonzero.

\textup{(3)} We have \(h=X+(X^3)^2+X^3\), so the reduced polynomial in Theorem~\ref{v70:thm-phase-reduction} is \(R=X\), with \(H=X^3\). Thus the phase over every extension is that of \(X\), and orthogonality gives \(S_r(h)=0\) and \(L_h(u)=1\). As in \textup{(2)}, the original degree is divisible by the characteristic, but the reduced phases are respectively linear and constant, with different outcomes.
\end{solution}

\solutionheading{ex:diagonal-zero-fiber}
\begin{solution}
Expanding as in the proof of Theorem~\ref{thm:diagonal-weil-bound}, the constant term still contributes \(q^{n-1}\),
and terms mixing constants with nontrivial characters still vanish.
For tuples consisting entirely of nontrivial characters, the quantity \(S(0)\) in the proof of Proposition~\ref{prop:multiple-jacobi-gauss}
is zero when the total character is nontrivial and \(-(q-1)J\) when it is trivial, giving the stated identity and bound.
If \(n=1\), there are no character tuples satisfying the condition, and the formula gives \(N_0=1\).

In odd characteristic, \(x^2-y^2=(x-y)(x+y)\), whose zero set consists of two lines meeting at the origin,
so \(N_0=2q-1\). Hence \(N_0-q=q-1>\sqrt q\),
whereas the general nonzero-right-hand-side estimate would give only \(\sqrt q\) here, showing that \(b\ne0\) cannot be omitted.

When all \(d_i=2\), the unique character tuple is \((\eta,\ldots,\eta)\),
where \(\eta\) is the quadratic character and the total character is \(\eta^n\). Write \(\Delta=\prod_i a_i\) and \(G=G(\eta,\psi)\). Then
\[
J(\eta,\ldots,\eta)=
\begin{cases}
G^{n-1},&n\text{ odd},\\
-G^n/q,&n\text{ even}.
\end{cases}
\]
Substitution of \(G^2=\eta(-1)q\) into the solution-count formula gives, in odd dimension,
\[
N_b=q^{n-1}+\eta\bigl((-1)^{(n-1)/2}\Delta\bigr)\eta(b)q^{(n-1)/2};
\]
and in even dimension,
\[
N_b=q^{n-1}+\eta\bigl((-1)^{n/2}\Delta\bigr)v(b)q^{(n-2)/2},
\qquad
v(b)=\begin{cases}q-1,&b=0,\\-1,&b\ne0.\end{cases}
\]
This is exactly the formula in Theorem~\ref{thm:quadratic-odd-count} for \(m=r=n\).
The two scales correspond to a nontrivial and a trivial total character, respectively.
\end{solution}

\solutionheading{ex:primitive-zero-trace-boundary}
\begin{solution}
Trace zero is equivalent to \(x+x^q=0\). For \(x\ne0\), we have \(x^{q-1}=-1\),
hence \(x^{2(q-1)}=1\). Since \(q+1>2\),
we have \(2(q-1)<q^2-1\), so \(x\) cannot generate \(E^\times\).
In characteristic two, trace zero is further equivalent to \(x^q=x\), so the trace-zero fiber is exactly \(k\).

When \(n=2\), we have \(R=q+1>1\) and hence \(W(R)-1\ge1\).
Formula~\eqref{eq:primitive-trace-zero-criterion} would then require
\(q-1>(W(q+1)-1)(q-1)\), which is impossible.
Thus this sufficient condition cannot hold on the trace-zero fiber of a quadratic extension, consistently with the absence of primitive trace-zero elements.
\end{solution}

\solutionheading{ex:curve-quadratic-term}
\begin{solution}
When \(x\ne0\), put \(z=y/x\) to obtain
\[
 z^2+z=x+c+b^2/x^2.
\]
The same trace count as in Proposition~\ref{prop:binary-curve-point-count} gives
\[
 \#E_{b,c}(k_n)=q^n+1+\psi_n(c)K_n(1,b).
\]
Since \(c\in k\), trace transitivity gives
\[
 \Tr_{k_n/\F_2}(c)=n\Tr_{k/\F_2}(c),\qquad \psi_n(c)=\varepsilon^n.
\]
This proves the point-count formula. Replacing \(u\) by \(\varepsilon u\) and using \eqref{eq:kloosterman-L} gives
\[
 \exp\left(\sum_{n\ge1}\varepsilon^nK_n(1,b)\frac{u^n}{n}\right)
 =1+\varepsilon K_1(1,b)u+q\varepsilon^2u^2
 =1-\varepsilon au+qu^2.
\]
Multiplication by \((1-u)^{-1}(1-qu)^{-1}\) then gives the zeta function.
Even \(n\) satisfy \(\varepsilon^n=1\), so the two curves have equal point counts over every extension of even degree.
When \(q=2\) and \(b=c=1\), we have \(\varepsilon=-1\) and \(a=-1\); the first four point counts are
\[
 2,\qquad8,\qquad14,\qquad16.
\]
\end{solution}

\solutionheading{v40:ex:root-power-fourier}
\begin{solution}
(1) Put \(\eta=\zeta^5\). Then \(\eta\) has order three and \(\eta^2+\eta+1=0\). Hence
\[
 s_j=1+\eta^j+\eta^{2j}
 =\begin{cases}1,&3\mid j,\\0,&3\nmid j.\end{cases}
\]
The frequencies with nonzero coefficients are \(\{0,3,6,9,12\}\), whose greatest common divisor together with 15 is three.
Thus \(H_+(T)=5\Z/15\Z\), and the maximum monomial composition exponent of \(P_T=X^3-1\) is three.
For \(C_0\), all root power sums equal one, giving a trivial stabilizer; for \(C_5\), the root power sums are \(\eta^j+\eta^{2j}\), zero when \(3\mid j\) and one otherwise. The frequencies with nonzero coefficients include one, again giving a trivial stabilizer.
The corresponding polynomials \(X-1\) and \(X^2+X+1\) both have maximum monomial composition exponent one.
An exponent translation stabilizing the union can cross between its parts, corresponding to multiplication on the root set; thus the union's stabilizer cannot be recovered as the intersection of the two stabilizers.

(2) Since \(\xi^2=-1\), termwise calculation gives
\[
 (s'_0,s'_1,s'_2,s'_3)=(0,\xi,1,-\xi).
\]
The frequencies with nonzero coefficients are \(\{1,2,3\}\), whose greatest common divisor together with four is one, so \(H_+(T')=\{0\}\).
Here \(|T'|=3=0\) in \(\F_9\) makes the zero frequency vanish; nevertheless inversion remains valid, since four is invertible in \(\F_9\).
The zero frequency imposes no restriction on translations in any case, so it alone cannot determine whether the root set is empty or the stabilizer nontrivial.

(3) We have \(3T'=\{0,2,3\}\ne T'\), so this is not a Frobenius-stable defining set over \(\F_3\). Direct calculation also gives
\[
 P_{T'}(X)=(X^2-1)(X-\xi)
 =X^3-\xi X^2-X+\xi\notin\F_3[X],
\]
since \(\xi\) has order four, whereas \(\F_3^\times\) has order two.
The proof of the Fourier stabilizer criterion only requires that the value field contain a primitive \(N\)th root of unity and have characteristic not dividing \(N\); it does not require \(qT=T\).
Frobenius stability is used to place the root product and all root power sums in the specified base field.
\end{solution}

\solutionheading{v85:ex-restricted-sums}
\begin{solution}
\textup{(1)} The cyclic quotient \(\F_q^\times/H\) has \(d\) characters. Orthogonality gives
\[
 \mathbf1_H(x)=\frac1d\sum_{\chi\in H^\perp}\chi(x)
 \qquad(x\in\F_q),
\]
with both sides zero at \(x=0\). Without the distinctness condition, the count is \(\sum_x\mathbf1_H(x)\mathbf1_H(s-x)\).
Expand and substitute \(x=su\) to obtain the Jacobi-sum expression. In odd characteristic the only possible diagonal solution is \(x=y=s/2\), contributing \(\epsilon_s\); in characteristic two, \(s\ne0\) rules out diagonal solutions.

\textup{(2)} The term with two trivial characters is \(q-2\). The \(2(d-1)\) terms with exactly one trivial character and the \(d-1\) terms with inverse nontrivial characters have absolute value one.
The remaining \((d-1)(d-2)\) terms have absolute value \(\sqrt q\). Since each factor \((\chi\lambda)(s)\) has absolute value one, the triangle inequality for real parts gives
\[
 d^2(S_H(s)+\epsilon_s)\ge q-2-3(d-1)-(d-1)(d-2)\sqrt q.
\]
\textup{(3)} If \(s=0\), then \(y=-x\). Thus \(S_H(0)=|H|\) in odd characteristic when \(-1\in H\), and is zero otherwise; in characteristic two the two coordinates coincide.
Simultaneous scaling \((x,y)\mapsto(hx,hy)\), and simultaneous \(p\)th powers, are bijections preserving distinctness. They give the two invariances.
Consequently nonzero sum counts are first grouped by multiplicative cosets, which are then permuted by Frobenius.
\end{solution}

\solutionheading{v85:ex-distinct-products}
\begin{solution}
\textup{(1)} Use \(\gamma\) to turn products into sums of exponents modulo \(n\). Since \(g\mid n-t\) and \(n>t\), we have \(n-t\ge g\).
Each displayed set is therefore a \(t\)-element subset of \(\{0,\ldots,n-1\}\). Their sums differ successively by one and cover every residue modulo \(g\).
Translating each element of a set by \(r\) increases its sum by \(tr\). These increments are precisely the multiples of \(g\) modulo \(n\).
The translates of the \(g\) initial sets therefore produce every sum modulo \(n\), proving the product assertion.

\textup{(2)} Without distinctness, any \(x,y\in H\) determine \(z=\beta/(xy)\), giving \(n^2\) solutions.
Each of the three equality events \(x=y,y=z,z=x\) has \(n\) solutions. The intersection of any two, and of all three, is \(x=y=z\) with \(x^3=\beta\).
Inclusion--exclusion gives \(n^2-3n+3\nu_{H,3}(\beta)-\nu_{H,3}(\beta)\).
The cube map on a cyclic group gives
\[
 \nu_{H,3}(\beta)=
 \begin{cases}\gcd(3,n),&\beta\in H^3,\\0,&\beta\notin H^3.\end{cases}
\]
In particular, for \(n\ge4\), every \(\beta\in H\) has a representation by three distinct elements because \(n^2-3n>0\).

\textup{(3)} The constant term is \((-1)^t\prod_{a\in S}a\). For \(t=n\), the only choice is \(S=H\),
whose root polynomial is \(T^n-1\). Its constant term is only \(-1\), not a coset of size \(n\ge2\).
\end{solution}

\solutionheading{v95:ex-short-order}
\begin{solution}
(1) Two distinct subgroup elements are linearly dependent exactly when their ratio belongs to $\F_q^\times\setminus\{1\}$. The cyclic group $H\cap\F_q^\times$ has order $\gcd(n,q-1)$, proving the criterion.
(2) The minimal polynomial of an element of order $d\mid n$, $d>1$, supplies a relation of support at most $\ord_d(q)+1$. Hence the extreme condition forces $\ord_d(q)=k$. Every prime $\ell\mid n$ then satisfies $k\mid\ell-1$ and $\ell\ge k+1$. Multiplying the prime factors modulo $k$ gives $n\equiv1\pmod k$; multiplying their lower bounds gives $n\ge(k+1)^{\Omega(n)}$. Every such $\ell$ contributes only to the factor $\Phi_k(q)$ in $q^k-1$, whence $n\mid\Phi_k(q)$.
(3) The $2$-cyclotomic cosets modulo seven are $\{0\}$, $\{1,2,4\}$ and $\{3,5,6\}$. Yet the sum of any two distinct elements of $\F_8^\times$ is again in this group. Thus a three-term relation exists, whereas no two-term relation does, so $\lambda_2(H)=3$.
\end{solution}

\solutionheading{v95:ex-triples}
\begin{solution}
(1) Every $t\in\F_q^\times\setminus\{-1\}$ satisfies $t,1+t\in\F_q^\times\subset U$; there are $q-2$ such values.
(2) For $g=1$, each element of $U$ is uniquely a base-field scalar times an element of $H$. A nonscalar solution $u=tx$, $1+u=sy$ corresponds to the relation $1+tx-sy=0$ among the distinct elements $1,x,y\in H$. Conversely, given an ordered dependent triple $(h_0,h_1,h_2)$, divide by $h_0$ and normalize the constant coefficient to recover this $u$ uniquely. Therefore the number of ordered dependent triples is $n(\mathcal C(U)-q+2)$. Dividing by the six orders of each three-element subset gives the formula.
(3) If $q=2$, then $\F_q^\times=\{1\}$, $U=H$, and $q-2=0$. Thus the absence of dependent triples is exactly $(1+H)\cap H=\varnothing$. For $q>2$, the three coefficients of a linear relation need not all equal one. Their normalization requires replacing $H$ by $U$.
\end{solution}

\solutionheading{v95:ex-period-three}
\begin{solution}
For $a\in L$, write $\eta(a)=\sum_{u\in U}\Psi(au)$ and $\eta(0)=M$. Orthogonality yields
\[\sum_{a\in L}\eta(a)^r=Q\#\{(u_1,\ldots,u_r)\in U^r:\sum u_i=0\}.\]
For $r=1$, the right side vanishes. For $r=2$, it equals $QM$ because $-U=U$. For $r=3$, fix the first term and write the other two as $uv$ and $-u(1+v)$; there are exactly $M\mathcal C(U)$ such triples. Each nonzero period occurs $M$ times, so subtracting $M^r$ at the zero parameter and dividing by $M$ gives
\[\sum_i\eta_i=-1,\quad\sum_i\eta_i^2=Q-M,\quad\sum_i\eta_i^3=Q\mathcal C(U)-M^2.\]
For $g=1$ and $k\ge3$, no two-term dependency exists. By the cyclotomic triple criterion, $\lambda_q(H)\ge4$ precisely when $\mathcal C(U)=q-2$, equivalently $\sum_i\eta_i^3=Q(q-2)-M^2$.
\end{solution}

\solutionheading{v95:ex-seven-periods}
\begin{solution}
(1) Trace hyperplanes in dimension three correspond to the $q^2+q+1$ lines of $\mathrm{PG}(2,q)$. Each projective point lies on $q+1$ lines, so $\sum N_H(a)=n(q+1)$. Any two distinct projective points from $H$ determine exactly one line; general position ensures that no such line contains a third point of $H$. Hence $\sum\binom{N_H(a)}2=\binom n2$. It follows that $A_2=\binom n2$, $A_1+2A_2=n(q+1)$, and $A_0+A_1+A_2=q^2+q+1$, giving the formulas.
(2) For $q=11,n=7$, one obtains $A_2=21$, $A_1=42$, $A_0=70$. Each $U=\F_{11}^\times H$-coset contains seven projective parameters, and the period value is $\eta=11N_H-7$. Thus among the nineteen periods the value $-7$ occurs ten times, $4$ six times, and $15$ three times. Their third moment is
\[
10(-7)^3+6\cdot4^3+3\cdot15^3=7079
=1331\cdot9-70^2.
\]
(3) For $q=2,n=7$, $H=\F_8^\times$ contains relations $x+y=z$, so $N_H(a)$ can attain three and cannot have the distribution restricted by $N_H(a)\le2$. The two Frobenius multipliers modulo seven, two and four, generate the same subgroup of order three, so the cyclotomic coset partitions coincide.
\end{solution}

\solutionheading{v95:ex-fourier-minor}
\begin{solution}
(1) The determinant $D_{\Gamma,S}(\zeta_\ell)$ is a polynomial with integer coefficients in a root of unity, hence an algebraic integer. Chebotarëv's theorem on minors makes it nonzero. The determinant expansion has at most $k!$ signed terms of complex absolute value one; the same bound holds under each of the $\ell-1$ embeddings of $\Q(\zeta_\ell)$. Thus its nonzero integer norm has absolute value at most $(k!)^{\ell-1}$.
(2) If some minor vanished upon reduction modulo a prime $p$, the corresponding cyclotomic integer would lie in a prime ideal above $p$, so $p$ would divide its integer norm. This contradicts its nonzero norm having absolute value less than $p$.
(3) In characteristic $p$, the image of $\zeta_\ell$ is a primitive $\ell$th root $\beta$, and the row exponents $\Gamma=\langle q\rangle$ are exactly the Frobenius powers $\beta^{q^j}$. Thus the corresponding matrix is a Moore matrix; nonvanishing of every $k$-column minor says that all $k$ distinct subgroup elements are independent over $\F_q$, so $\lambda_q(H)=k+1$. The only external theorem used is Chebotarëv's nonvanishing theorem for prime-order Fourier minors; the remaining steps use integer norms, determinant estimates, and the Moore criterion.
\end{solution}

\solutionheading{v109:ex-power-jacobian}
The entries are $J_{ji}=jx_i^{j-1}$. Factoring $j$ from row $j$ gives the stated Vandermonde determinant.
Let $r$ be the number of distinct coordinates. Equal coordinates give equal columns, so rank is at most $\min(m,r)$.
Choose $t=\min(m,r)$ distinct coordinates and the first $t$ rows. The resulting minor is $t!$ times a nonzero Vandermonde product. It is nonzero because $p>m$, proving equality.
For $m=N=2$ in odd characteristic with $x_1\ne x_2$, the local rank is two, but swapping the coordinates preserves both power sums. Thus local rank does not imply global injectivity.
Over $\F_2$, take $m=N=2$ and $(x_1,x_2)=(0,1)$. The second row is zero, and rank is one despite two distinct coordinate values.

\solutionheading{v115:ex-two-moments}
\begin{solution}
Both functions have total six, mean two, and \(V=6\), but their numbers of empty fibers are one and zero.
The bound is \(\lfloor3\cdot6/(3\cdot4+6)\rfloor=1\), compatible with both cases.
Let \(\omega=e^{2\pi i/3}\) and \(\chi(j)=\omega^j\). With the conjugation convention of the book, the two pairs of nontrivial coefficients are
\[
(-3,-3),\qquad (3\omega,3\omega^2).
\]
Both squared-magnitude sums are eighteen, but inversion gives different values. Here the phases lost by the second moment determine whether an empty fiber occurs.
\end{solution}

\solutionheading{v145:ex-quadratic-fibers}
\begin{solution}
Matrix multiplication gives $b_2(a)=2-a^{P+1}=2-\Nm_{\F_{P^2}/\F_P}(a)$. Thus the fiber over $b=2$ consists only of $a=0$, and every fiber over $b\ne2$ has $P+1$ elements: a nonzero norm fiber is a coset of the norm-one group of order $P+1$. The target $b=2$ corresponds to $(Z+1)^2$. The fiber sizes sum to $1+(P-1)(P+1)=P^2$.
For $a=0$, the return is $-I_2$. In odd characteristic the other repeated-root target, $b=-2$, has $\Nm(a)=4$, so $a\ne0$; the return is nonscalar and has a nontrivial Jordan block with eigenvalue one. In characteristic two the targets coincide, and the only repeated-root return is the identity at $a=0$.
\end{solution}

\solutionheading{v146:ex-generator-values}
\begin{solution}
Applying the \(P\)-power Frobenius entrywise to the matrix product gives \(\Theta_{k,r}(a^P)=\Theta_{k,r}(a)^P=\Theta_{k,r}(a)\), so the function is constant on conjugate orbits.

In the given quintic extension, put \(a=1+\alpha\). The cyclic matching expansion gives
\[
 \Theta_{5,1}(a)=\Nm(a)+\Tr(a^7)+\Tr(a),\qquad
 \Theta_{5,2}(a)=\Nm(a)+\Tr(a^{21})+\Tr(a).
\]
From \(\alpha^5=\alpha^3+1\), we obtain \(a^7=1+\alpha^3\) and \(a^{21}=\alpha+\alpha^2\). Newton identities give \(\Tr(\alpha)=\Tr(\alpha^2)=\Tr(\alpha^3)=0\), whereas \(\Tr(1)=1\). Also \(\Nm(a)=1\), so the two values are one and zero, respectively. Both irreducible-target fibers have size \(D_5(2)=11\), but the functions differ.
\end{solution}

\Needspace{9\baselineskip}
\solutionheading{v148:ex-repeated-profile}
The first target has $U_2(h)=2^{2m+1}-1$ and $E(h)=1$. Its unique commutative norm preimage $(Z+1)^m$ has constant one. Put $a=(-2)^m$ and $B=(2^{2m+1}-1-a)/3$. The three fibers are $B+a,B,B$: the first is largest for even $m$, and the last two are largest for odd $m$.
The second target has $E(h)=e+1$. Its commutative norm preimages are $(Z+\omega)^j(Z+\omega^2)^{e-j}$, $0\le j\le e$, with constants $\omega^{2e-j}$. The three residue classes occur equally often exactly when $3\mid e+1$. Thus the fibers are uniform exactly when $e\equiv2\pmod3$.

\section*{Solutions to Chapter~\ref{chap:period-layers}}
\pdfbookmark[1]{Solutions to Chapter \getrefnumber{chap:period-layers}}{solutions.chap:period-layers}

\solutionheading{ex:wu-E26}
\begin{solution}
We have
\[
E_{2,6}=\frac{(1-X^{63})(1-X)}{(1-X^7)(1-X^3)}.
\]
The degree is \(63+1-7-3=54\), so there are \(54/6=9\) irreducible sextic factors.
Nonzero elements of the degree-six field have orders dividing \(63\);
the divisors satisfying \(\ord_e(2)=6\) are \(e=9,21,63\), so \(E_{2,6}=\overline\Phi_9\,\overline\Phi_{21}\,\overline\Phi_{63}\).
Modulo \(X^3\), both the denominator and \(1-X^{63}\) equal one, so the remainder is \(1+X\).
\end{solution}

\solutionheading{ex:wu-layer-radical}
\begin{solution}
An element of degree 12 over \(\F_q\) has degree \(12/\gcd(12,2)=6\) over \(\F_{q^2}\).
Conversely, an element of degree six over \(\F_{q^2}\) is fixed by the \(q^{12}\)th-power map;
its \(q\)-period \(d\) satisfies \(d\mid12\) and \(d/\gcd(d,2)=6\).
Comparing the 2-adic valuations and odd parts gives \(d=12\).
If the base field is instead \(\F_{q^3}\), an element of degree four over \(\F_q\)
still has degree four over \(\F_{q^3}\), since \(\gcd(4,3)=1\).
Thus relative degree four does not distinguish original degree four from degree 12.
\end{solution}

\solutionheading{ex:wu-compressed-not-period}
\begin{solution}
Both two and three are prime, so \(G_{q,\ell}=R_{q,\ell}\).
Substitution into the explicit formula for \(R\) gives both identities.
If \(\alpha\) has degree two over \(k\), then \(1+\Theta_q(\alpha)=0\), so its compressed value is \(\Theta_q(\alpha)=-1\in k\).
The original element has period two, but the compressed value has period one.
By Corollary~\ref{cor:compressed-orbit-period}, compression identifies elements in the same affine-group orbit;
the degree of a compressed root is determined by its Frobenius period on the orbit quotient. In the quadratic case that quotient has one element, so the period is one.
\end{solution}

\solutionheading{ex:wu-polynomial-invariance}
\begin{solution}
We have \(Z_q(aX+b)=aZ_q(X)\), so \(\Theta_q(aX+b)=a^{q-1}\Theta_q(X)=\Theta_q(X)\).
Also \(x^q=x\) for every \(x\in\F_q\), so \(\Theta_q(x)=0\) at these points.
But its leading term is \(X^{q(q-1)}\), so it is not the zero polynomial.
Pointwise equality on the base field only says that the difference is divisible by \(X^q-X\); the definition of an invariant requires the difference polynomial itself to be zero.
For example, \(F=(X^q-X)X\) vanishes everywhere on the base field, but for nonzero \(b\in k\),
\(F(X+b)-F(X)=b(X^q-X)\ne0\), so it is not translation-invariant.
\end{solution}

\solutionheading{ex:wu-compressed-separable}
\begin{solution}
\textup{(1)} If \(G_{q,m},G_{q,n}\) have a common root \(c\), adjoin a root \(\alpha\) of \(\Theta_q(X)-c\) in a finite extension field;
then \(E_{q,m}(\alpha)=E_{q,n}(\alpha)=0\).
An element has a unique least positive period, so \(m=n\).
The nonconstant polynomial \(\Theta_q(X)-c\) acquires a root in a finite extension field; thus a common root of compressed layers lifts to a common root of the original period layers.

\textup{(2)} The multiplicative order of \(\beta=\alpha^8\) is four, so \(\beta^5=\beta\) and \(\beta\in\F_5^\times\).
Also \(5^4\equiv17\not\equiv1\pmod{32}\) and \(5^8\equiv1\pmod{32}\),
so \(\ord_{32}(5)=8\) and hence \([k(\alpha):k]=8\).
The element \(\xi=\alpha^4\) has order eight, so \(\mu_8=\langle\xi\rangle\).
The eight distinct roots of \(X^8-\beta\) are exactly \(\alpha\mu_8\), all belonging to \(k(\alpha)\).
Since \(\beta\in k\), the base orbit has length one and the first-return map is the restriction of \(\sigma:x\mapsto x^5\) to this fiber.

In the fiber coordinates \(x=\alpha\varepsilon\),
\[
\sigma(\alpha\varepsilon)=\alpha^5\varepsilon^5
=\alpha\xi\varepsilon^5,
\qquad \mathcal R(\varepsilon)=\xi\varepsilon^5.
\]
Putting \(\varepsilon=\xi^j\) gives \(\mathcal R(\xi^j)=\xi^{1+5j}\). The full cycle starting at exponent \(0\) is
\[
0\longmapsto1\longmapsto6\longmapsto7\longmapsto4
\longmapsto5\longmapsto2\longmapsto3\longmapsto0.
\]
These eight residue classes are distinct, so the return map traverses the whole fiber, with internal period eight and total period \(1\cdot8=8\).
In this example the full fiber is also the entire Frobenius orbit of \(\alpha\).

The actual iteration satisfies \(\mathcal R(1)=\xi\) and \(\mathcal R^2(1)=\xi^6\).
If we instead use \(\mathcal M(\varepsilon)=\xi\varepsilon\), then \(\mathcal M(1)=\xi\) and \(\mathcal M^2(1)=\xi^2\).
By \(\ord(\xi)=8\), we have \(\xi^6\ne\xi^2\), so the two return permutations differ.
Multiplication by \(\xi\) alone also has order eight, but the same order does not determine the permutation of a fiber;
the fifth-power action in this example cannot be omitted.
\end{solution}

\solutionheading{ex:wu-binary-specialization}
\begin{solution}
The recurrence gives
\[
Q_{q,2}=V,\quad Q_{q,3}=V^{q+1}-U^q,\quad
Q_{q,4}=V^{q^2}(V^{q+1}-U^q)-U^{q^2}V.
\]
In characteristic \(p\) we have \((1+T)^q=1+T^q\), so
\[
(1+T)^{q+1}-T^q=(1+T)(1+T^q)-T^q
=1+T+T^{q+1}.
\]
After setting \(Y=1\), we have \((\Delta_q,\Gamma_q)=(\Theta_q,1+\Theta_q)\);
hence the homogeneous period-layer polynomial \(\mathcal E_{q,3}\) specializes to \(R_{q,3}(\Theta_q)=G_{q,3}(\Theta_q)=E_{q,3}\).
\end{solution}

\solutionheading{ex:wu2-zeta-basechange}
\begin{solution}
The fixed field of \(\sigma^r\) is \(\F_{q^{\gcd(n,r)}}\), proving the point-count formula.
In \(S\) only orbit lengths dividing \(n\) occur, with \(I_q(d)\) orbits of each length \(d\mid n\), so
\[
Z_S(u)=\prod_{d\mid n}(1-u^d)^{-I_q(d)}.
\]
On a \(d\)-cycle, \(\sigma^s\) has \(g=\gcd(d,s)\) cycles of length \(d/g\),
so one original factor \((1-u^d)^{-1}\) becomes \((1-u^{d/g})^{-g}\).
This is the same orbit formula used to factor irreducible polynomials after extending the base field.
\end{solution}

\solutionheading{ex:wu7-compressed-linear-factors}
\begin{solution}
Since \(G_{q,n}(0)=1\), it suffices to consider \(c\in k^\times\). In a finite extension field take \(\Theta_q(\alpha)=c\) and write \(z=\alpha^q-\alpha\). By \(c\ne0\), we have \(\alpha\notin k\) and \(z\ne0\).
Since \(\Theta_q(\alpha^q)=c^q=c\), the orbit-separation proposition gives a unique \(\alpha^q=a\alpha+b\) with \(a\in k^\times,b\in k\). Hence
\[
z^q=(\alpha^q)^q-\alpha^q=a(\alpha^q-\alpha)=az,
\qquad c=z^{q-1}=a.
\]
If \(c=1\), then \(\alpha^q=\alpha+b\) and \(b\ne0\); the return transformation is a nontrivial translation,
so \([k(\alpha):k]=p\). If \(c\ne1\), the return transformation has order \(\ord(c)\),
and hence \([k(\alpha):k]=\ord(c)\). Here the compressed value lies in the base field, so the compressed period is \(t=1\).

By \(E_{q,n}=G_{q,n}(\Theta_q)\), the element \(c\) is a root of \(G_{q,n}\)
if and only if \(\alpha\) above has degree \(n\). Thus \(c=1\) contributes only the factor \(T-1\) when \(n=p\);
for \(c\ne1\), it collects exactly the base-field units of order \(n\). Such units exist if and only if \(n\mid q-1\),
and their root polynomial is \(\overline\Phi_n\). Since \(p\nmid q-1\), the two nonempty cases do not overlap.

Taking \(c=1\) realizes period-reduction factor \(p\); for every \(e>1\) with \(e\mid q-1\),
choose \(c\) of order \(e\) in the cyclic group \(k^\times\) to realize period-reduction factor \(e\).
This proves sharpness of the restriction on nontrivial ratios in the proposition.
\end{solution}

\Needspace{7\baselineskip}
\section*{Solutions to Chapter~\ref{revadic:chapter}}
\pdfbookmark[1]{Chapter \getrefnumber{revadic:chapter} solutions}{solutions.revadic:chapter}

\solutionheading{revadic:ex1}
\begin{solution}
The polynomial $T^{p-1}-1$ has $p-1$ distinct roots in $\mathbb F_p$, each of which lifts by Hensel's lemma to $\mathbb Z_p$. For $p=2$ the modulus $p-1=1$, so the usual multiplicative order is undefined; the extension is simply taken to be trivial.
\end{solution}

\solutionheading{revadic:ex2}
\begin{solution}
The full list is $\{0\}$, $\{1,3\}$, $\{2,6\}$, $\{4\}$, and $\{5,7\}$. The reduced factors are $T^2-T+2$ and $T^2+T+2$. Their discriminant is $1-8=-7\equiv2\pmod3$, a nonsquare, so they are irreducible.
\end{solution}

\solutionheading{revadic:ex3}
\begin{solution}
The intersection $\mathcal O_n\cap K_d$ is the ring of integral elements of $K_d$, namely $\mathcal O_d$. The elements fixed by $\phi^d$ satisfy $x^{p^d}=x$; there are exactly $p^d$ such elements in $k_n$, forming $k_d$.
\end{solution}

\solutionheading{revadic:ex11}
\begin{solution}
By Proposition~\ref{revadic:prop-precision-roots}, each of the eight residue elements uniquely lifts to a solution modulo $2^s$. For $s\ge2$, closure under addition would put $1+1=2$ among the solutions, whereas $2^8\not\equiv2\pmod{2^s}$. Modulo two, the entire residue field is recovered and is additively closed.
\end{solution}

\solutionheading{revadic:ex15}
\begin{solution}
If $\tau$ were additive, then $\tau(p)=p\tau(1)=p$, whereas $\tau(p)=0$. Both sides have the same residue, proving $\tau(u+v)\equiv u+v\pmod p$. For $p=2,n=1$, $1\boxplus1=0$. Associativity follows because reduction identifies the transported addition with addition in $k_n$.
\end{solution}

\solutionheading{revadic:ex24}
\begin{solution}
Both equations reduce to $X^2+X+1$ and have derivative $1$ modulo two. For $t=[a]=\zeta_3$, $F_1(t)=0$, while $F_2(t)=-2$ and $F_2'(t)=2t+1$ is a unit. By Proposition~\ref{revadic:prop-root-correction}, $z-t\equiv2/(2t+1)\equiv2\pmod4$; the first digit of the correction is $1$.
In general $H_r(t)=-2^r$, so $v_2(z_r-t)=r$ and $z_r\equiv t+2^r/(2t+1)\pmod{2^{2r}}$.
\end{solution}

\solutionheading{revadic:ex4}
\begin{solution}
The inverse Gram matrix is $\frac13\begin{psmallmatrix}1&2\\2&1\end{psmallmatrix}$. Therefore the dual basis is $(\zeta_3+2\zeta_3^2)/3$ and $(2\zeta_3+\zeta_3^2)/3$. Modulo two, $3^{-1}=1$ and the reductions are $\bar\zeta_3$ and $\bar\zeta_3^2$. In this example the normal basis is self-dual after reduction.
\end{solution}

\solutionheading{revadic:ex5}
\begin{solution}
Replace $\mathbb Z_p$, $\phi$, and $\mathbb F_p$ in the determinant argument by $\mathcal O_d$, $\phi^d$, and $k_d$. The same unit-determinant criterion works. The reduction of $\theta+pz$ is unchanged, so it also generates a normal integral basis.
\end{solution}

\solutionheading{revadic:ex6}
\begin{solution}
We have $\zeta_3\zeta_3^2=1$ and $\zeta_3+\zeta_3^2=-1$. Teichmüller representatives are multiplicatively closed, so the norm product is exactly the unique Teichmüller lift of the reduced norm. Trace is a sum and may acquire carries.
\end{solution}

\solutionheading{revadic:ex12}
\begin{solution}
In characteristic two, $(T-1)^4=T^4-1$. Relative to an integral lift $\Theta$ of a normal element, the criterion is $\gcd(\bar f,(T-1)^4)=1$, equivalently $\bar f(1)\ne0$. Repeated factors do not obstruct the unit test in the cyclic quotient module.
\end{solution}

\solutionheading{revadic:ex8}
\begin{solution}
Write $x=1+at+bt^2+ct^3+\cdots$. Comparing coefficients modulo three gives $2a=1$, hence $a=2$; $a^2+2b=0$, hence $b=1$; and $2c+2ab=0$, hence $c=1$. Thus $x=1+2t+t^2+t^3+O(t^4)$. A unital section would imply $3\cdot1_{\mathbb Z_3}=0$, a contradiction.
\end{solution}

\solutionheading{v138:ex-word-denominator}
\begin{solution}
The repeated block has value $E/(1-q^D)$. Thus the required condition is that $a=mE/(q^D-1)$ be integral; the inequalities $0<E<q^D-1$ then give $1\le a<m$. Equivalently, the reduced denominator $n=(q^D-1)/\gcd(E,q^D-1)$ divides $m$.
Now $n=m/\gcd(a,m)$, so the least period is $d=\operatorname{ord}_n(q)=|C_q(a)|$ and $d\mid D$. The given block consists of $D/d$ repetitions of the shortest block.
\end{solution}

\solutionheading{v138:ex-digit-complement}
\begin{solution}
The exponents satisfy $E_{m-a}=Q-1-E_a$, so the characters are inverse and the Gauss-sum product formula gives the first identity. Also,
$\omega_{Q^h}^{-a(Q^h-1)/m}=\omega_Q^{-E_a}\circ\Nm_{\mathbb F_{Q^h}/\mathbb F_Q}$; the standard additive character is lifted by the trace. Hasse--Davenport gives the second identity.
Write $Q=p^L$. All $L$ digits of $Q-1$ are $p-1$, so the digits of $E_a$ and $Q-1-E_a$ are complementary. Their digit sums add to $L(p-1)$, and the valuations add to $L$. The extension exponent is
\[
 \frac{a(Q^h-1)}m=E_a(1+Q+\cdots+Q^{h-1}).
\]
Its base-$p$ block is repeated $h$ times, multiplying the valuation by $h$ but leaving its quotient by $\log_p(Q^h)=hL$ unchanged. Both blocks encode the same rational number $-a/m$, so the least digit period is unchanged.
\end{solution}

\solutionheading{v139:ex-three-power}
\begin{solution}
An odd exponent cannot give $2^t\equiv1\pmod3$. For even $t=2u$, the valuation-lifting lemma gives $v_3(2^{2u}-1)=1+v_3(u)$. Hence $\ord_{3^s}(2)=2\cdot3^{s-1}$, so $C$ consists of all units modulo $3^s$. Put $r=3^{s-1}$. The two complete blocks are $1+3\mathbb Z/m\mathbb Z$ and $2+3\mathbb Z/m\mathbb Z$, both of size $r$; their stabilizer cannot be larger since its order divides $\gcd(m,2r)=r$. Multiplication by $r$ sends $-1/m$ to $-1/3$, of least binary period two. The canonical factor is
\[
 \Phi_{3^s}(X)=X^{2r}+X^r+1.
\]
Corollary~\ref{v139:cor-compressed-code} identifies the code with $r$ interleaved copies of the $[3,1,3]$ repetition code. Its parameters are $[3^s,r,3]$. This includes $s=1$, when the block size is one.
\end{solution}

\solutionheading{v139:ex-dual-words}
\begin{solution}
The Euclidean dual defining set is $-Z^c$. Let $\mathcal E_m=\{0,c,\ldots,(m-1)c\}$. Its encoding is
\[
 S_{\mathcal C_Z^\perp}=-(\mathcal E_m\setminus S_Z)\pmod M.
\]
For a nonzero encoded integer $E$, negation is $M-E$, the digitwise $(q-1)$-complement in $L$ positions. Zero is sent to zero, represented by the all-zero block, not by the excluded all-$(q-1)$ block. The omitted exponents of the dual are $-Z$. Negation preserves the multiset of cyclic gaps, so the required bound is the largest cyclic gap between elements of $Z$ itself. For the stated example, $Z=\{1,2,3,4,6,8,9,12\}$, with gaps $1,1,1,2,2,1,3,4$. The bound is four and the dual dimension is eight.
\end{solution}

\solutionheading{v140:ex-trace-hyperplane}
Since $\Tr(x^2)=\Tr(x)$, we have $V^{[2]}=V$. The condition $aV\subseteq V$ is equivalent to $\Tr(ax)=0$ for every $x\in V$. Nondegeneracy of the trace pairing and $V^\perp=\mathbb F_2\cdot1$ give $U(V)=V\cap\mathbb F_2$. Since $\Tr(1)=n\bmod2$, for odd $n$ the stabilizer is $2[V]$, of order $2^{n-1}$. For even $n$ it is $\{[a]+2[b]:a\in\mathbb F_2,b\in V\}$, of order $2^n$. The set $\mathcal T_2(V)$ has $2^{2n-2}$ elements, so it equals its own stabilizer exactly when $n=2$. Then $V=\mathbb F_2$ and $\mathcal T_2(V)=\mathbb Z/4\mathbb Z$.
For even $n>2$, Example~\ref{v70:ex-trace-hyperplane-intersection} gives $H(V)=\mathbb F_2$, so Theorem~\ref{v147:thm-carry-directions} also yields $U(V)=H(V)\subsetneq V$.

\solutionheading{v140:ex-reciprocal-profile}
Example~\ref{sfr:ex-quartic-comparison} gives the residue profile $(1,1,1,2)$ for $\gamma$. Precision profiles decrease with $s$, and constants always persist. In the full power space $t=3$, Proposition~\ref{revadic:prop-fixed-precision} lifts every direction of the quadratic subfield. Thus $\rho_s^\gamma(t,2)=(1,1,1,2)$ for every $s\ge1$. The profile for $\theta$ is $(1,1,2,2)$ at $s=1$ and $(1,1,1,2)$ at $s\ge2$. The advantage of a quadratic representative in the residue field disappears at the second level.

\solutionheading{v140:ex-two-point-group}
Write the permutation as $x\mapsto c x^{q^j}$. Fixing $a$ forces $c=a/a^{q^j}$; fixing $b$ as well is equivalent to $(b/a)^{q^j}=b/a$, or $d\mid j$. Hence the pointwise stabilizer has order $n/d$. Acting on a pair changes its ratio to $(b/a)^{q^j}$, and the scalar can send its first coordinate to any chosen nonzero element. Thus pair orbits are parameterized by the $q$-Frobenius orbits of ratios in $E^\times\setminus\{1\}$.

The residue group $\operatorname{GL}_n(q)$ is transitive on linearly independent ordered pairs. Dependent pairs of distinct coordinates are parameterized by their ratios $\lambda\in\mathbb F_q^\times\setminus\{1\}$. The latter remain single orbits at every precision; the independent-pair orbit separates at higher precision into the Frobenius orbits of ratios in $E\setminus\mathbb F_q$. There is no split if this set already forms one orbit, as for $q=2,n=2$.

\solutionheading{revadic:ex7}
\begin{solution}
The root set decomposes into orbits of size $e$; its root polynomial is a product of factors $X^e-a$. In the seventh-root example, $\eta\equiv0\pmod2$ but $\eta\eta'=2$ with $\eta'$ a unit, so $v_2(\eta)=1$.
\end{solution}

\solutionheading{v144:ex-seventy-three}
\begin{solution}
Since $\ord_{73}(2)=9$, the check quotient has rank nine and the code has rank 64. Mixed exponents with cyclic index distances $1,2,3,4$ form $C_2(3),C_2(5),C_2(9),C_2(17)$, respectively. Successive doubling shows that these four nine-element orbits are disjoint and avoid $C_2(1)$. They contain $4\cdot9=\binom92$ terms, proving separation. Corollary~\ref{v144:cor-subgroup-code} gives all permutations
\[
 i\longmapsto b+2^ji\pmod{73},\qquad b\in\mathbb Z/73\mathbb Z,\quad0\le j<9.
\]
The group has order $73\cdot9=657$, also over the complete integer ring.
\end{solution}

\solutionheading{revadic:ex19}
\begin{solution}
The three eigenvalues of multiplication are $\zeta,\zeta^2,\zeta^4$, so the trace of the first exterior power is $\eta$, while that of the second is $\zeta^3+\zeta^5+\zeta^6=\eta'$. Thus the characteristic polynomial is $X^3-\eta X^2+\eta'X-1$, exactly the canonical factor.
\end{solution}

\solutionheading{revadic:ex16}
\begin{solution}
The one-element subset sums are $1,2,4$, the two-element sums are $3,5,6$, and the three-element sum is $0$. Thus $(c_1,c_2,c_3)=(-\eta,\eta',-1)$. Since $\eta\eta'=2$ and $\bar\eta'=1$, we have $v_2(\eta)=1$: vanishing modulo two does not imply an exact zero in characteristic zero.
\end{solution}

\solutionheading{revadic:ex21}
\begin{solution}
Here $\Phi_9(T)=1+T^3+T^6$. Since $C=\{1,4,7\}$ is an arithmetic progression with step three, $A_{C,1}(T)=T+T^4+T^7=T\Phi_9(T)$. The pair sums are $5,8,2$, so $A_{C,2}(T)=T^2+T^5+T^8=T^2\Phi_9(T)$. Both cyclotomic remainders vanish identically. Therefore $Q_C(X)=X^3-\zeta_9^{1+4+7}=X^3-\zeta_9^3$; the missing coefficients vanish in characteristic zero.
\end{solution}

\solutionheading{v141:ex-involution}
Here $f=1$. Of the six pairwise differences only $x_0-x_2=-3$ is divisible by three, and it has valuation one. Thus $I_x=1$. Each nonzero direction has at least one unit difference, so $h(1)=h(2)=h(3)=0$ in additive group notation and $E_x=1=g_f$. The residue difference for direction two vanishes at coordinates zero and two; its kernel has dimension two, and the other two kernels vanish. The exceptional idempotent is $(1,0,1,0)$, the sum of two primitive idempotents. Thus an involution outside the trivial Frobenius subgroup need not have a kernel of dimension $f$. The profiles of $x$ and $x+2=(3,4,6,5)$ are respectively $(0,0,0,1)$ and $(1,0,1,0)$, although all their differences coincide.

\solutionheading{v141:ex-period-idempotent}
Put $b=\bar\eta_0+\bar\eta_2$. Reducing the product identities modulo two gives
$\bar\eta_0^2=\bar\eta_2$ and $\bar\eta_2^2=\bar\eta_0$, so $b^2=b$. Hence
$e=1-b^3=1+b=\bar\eta_1+\bar\eta_3$. In this case $-1\in H$, $g=4$, and $\ell=17$; the trace Gram matrix reduces to the identity. Consequently
\[
 \Tr(e)=0,\qquad
 (\Tr(e\bar\eta_i))_{i=0}^3=(0,1,0,1).
\]
Formula~\eqref{v141:eq-phase-recovery} therefore recovers $e(u)=\bar\eta_{u3}+\bar\eta_{u3^3}$, where the subscripts on the right are multiplicative coset labels. The ordinary trace vanishes because the supported field component has degree two in characteristic two; the trace pairing with all basis elements remains nondegenerate.

\solutionheading{v141:ex-short-return}
The first two pairs are fixed by multiplication by five. Over $\mathbb F_5$, with Teichmüller character $\omega_5$, their factors are
\[
 1+J_{\mathbb F_5}(\omega_5^{-1},\omega_5^{-1})T,
 \qquad
 1+J_{\mathbb F_5}(\omega_5^{-3},\omega_5^{-3})T,
\]
with slopes zero and one. Here $12\nmid4$, but the individual exponents $a(5-1)/12$ are integers. The third orbit is $(1,7),(5,11)$; its return field is $\mathbb F_{25}$ and its factor is
\[
 1+J_{\mathbb F_{25}}(\omega_{25}^{-2},\omega_{25}^{-14})T^2.
\]
There is one carry in the two steps, so the slope is $1/2$ with multiplicity two. Equivalently, the base-five digit sums of $2,14,16$ are $2,6,4$, giving valuation $(2+6-4)/4=1$ for the return scalar. The first two return fields do not contain all twelfth roots; the third does, since $12\mid24$.

\solutionheading{revadic:ex23}
\begin{solution}
Modulo two, $\Phi_7(T)=(T^3+T+1)(T^3+T^2+1)$. The integral remainder is $B_{C,1}(T)=T+T^2+T^4\ne0$, but modulo two this is $T(T^3+T+1)$, which vanishes at the chosen local root. Since $\operatorname{ord}_7(2)=3<\varphi(7)=6$, part (3) of Theorem~\ref{revadic:thm-integral-resolvent} does not apply. This is local congruence vanishing, not an exact zero in characteristic zero.
\end{solution}

\solutionheading{revadic:ex25}
\begin{solution}
The remainder of $A_{C,1}(T)=T+T^2+T^4$ is itself, so $\nu=0$. At a primitive root reducing to a root of $T^3+T+1$, the period $\eta=\zeta+\zeta^2+\zeta^4$ has valuation $1$; at a root of the other cubic factor, the corresponding period $\eta'=\zeta^3+\zeta^5+\zeta^6$ has valuation $0$. Thus the minimum is indeed $\nu$. The multiplicative stabilizer is $H=\langle2\rangle$, and $H\langle2\rangle=H\ne(\mathbb Z/7\mathbb Z)^\times$, so the sufficient hypothesis of Corollary~\ref{revadic:cor-symmetry-sharp} fails.
\end{solution}

\solutionheading{revadic:ex22}
\begin{solution}
The set $C$ is the group of quadratic residues modulo nineteen and $19\equiv3\pmod8$. Proposition~\ref{revadic:prop-qr-exterior} gives $N_{C,2}(0)=0$ and $N_{C,2}(s)=(19-3)/8=2$ for $s\ne0$. Thus $A_{C,2}(T)=2(\Phi_{19}(T)-1)$, so $c_2=-2$ and $v_2(c_2)=1$.
\end{solution}

\solutionheading{revadic:ex10}
\begin{solution}
Since $\alpha^4=\alpha+\alpha^2$, the partial sums are $s_1=\alpha$, $s_2=\alpha^4$, and $s_3=0$. In characteristic two, the square-root map is $x\mapsto x^4$, so
$\kappa(0,\alpha)=0$, $\kappa(\alpha,\alpha^2)=\sqrt{\alpha^3}=\alpha^5$, and
$\kappa(\alpha^4,\alpha^4)=\sqrt{\alpha^8}=\alpha^4$.
Since $\alpha^5=\alpha^4+1$, the sum of the carries is $1$. Thus $\eta/2\equiv1\pmod2$, proving $v_2(\eta)=1$.
\end{solution}

\solutionheading{revadic:ex13}
\begin{solution}
Theorem~\ref{revadic:thm-fourier} gives $\widehat f(pa)=\phi(\widehat f(a))$. Unramified Frobenius preserves $p$-adic valuations, so the transform coordinates have equal valuations along a coset.
\end{solution}

\solutionheading{revadic:ex14}
\begin{solution}
Substitute $y=x^p$ and use $\psi(y)=\psi(x)$. This proves Frobenius symmetry, but does not determine the first nonzero digit of the sum and therefore does not imply Stickelberger's digit-sum formula.
\end{solution}

\solutionheading{revadic:ex9}
\begin{solution}
The three cases give $\{0,f\}$, $\{v_p(a),f-v_p(a)\}$, and $\{f/2,f/2\}$. Complex absolute values bound archimedean magnitudes, not valuations at the chosen prime $p$.
\end{solution}

\solutionheading{revadic:ex17}
\begin{solution}
Set $A_r=N_r-1-q^r$. The coefficients are $b_0=1$, $b_1=A_1$, $b_2=(A_1^2+A_2)/2$, $b_3=qb_1$, and $b_4=q^2$. Compute their $p$-adic valuations; the slopes of the lower convex hull of $(j,v_p(b_j))$, divided by $f$, are the normalized Newton slopes.
\end{solution}

\solutionheading{revadic:ex20}
\begin{solution}
Here $b_1=A_1$ and $b_2=(A_1^2+A_2)/2$. If $b_2\not\equiv0\pmod p$, the slope-zero multiplicity is $2$; otherwise it is $1$ if $b_1\not\equiv0\pmod p$, and $0$ if both vanish. In characteristic two the integer division must precede reduction because $2$ is not invertible in $\mathbb F_2$.
\end{solution}

\solutionheading{revadic:ex18}
\begin{solution}
For $q=7$, $s_7(2)=2,s_7(4)=4$ give slopes $0,1$. For $q=25$, the base-five expansions $8=(13)_5$, $16=(31)_5$ both have digit sum $4$, giving slopes $1/2,1/2$. Both multisets are symmetric about $1/2$. The two non-elementary inputs are Hasse--Davenport lifting and the Stickelberger digit-sum valuation; the point-count decomposition itself uses character orthogonality.
\end{solution}

\solutionheading{v135:ex-finite-add}
\begin{solution}
Take $n=2,s=5,r=2$ in~\eqref{v135:eq-finite-projection}. Since $nr=s-1$, the exponent $9^2=81$ works, and reduction gives the nine-element residue field. For $s=6$, choose $r=3$, giving the exponent $9^3=729$.
\end{solution}
\solutionheading{v135:ex-hom-count}
\begin{solution}
A unital ring map fixes $\mathbb Z/2^s\mathbb Z$. Theorem~\ref{v135:thm-hom-lifting} gives exactly two maps, lifting the inclusion $k_2\hookrightarrow k_6$ and its composite with squaring Frobenius. On Teichmüller roots they send $x$ to $x$ and $x^2$, respectively. Each residue embedding has a unique compatible lift, so each map forms a compatible system through all precisions.
\end{solution}

\backmatter
\chapter*{Structural Overview and Quick Reference}

\addcontentsline{toc}{chapter}{Structural Overview and Quick Reference}
\markboth{Structural Overview and Quick Reference}{Structural Overview and Quick Reference}

\noindent\phantomsection\label{v98:par-problem-structure}\textbf{From a problem to its structure.}
Let \(k=\F_q\subseteq E=\F_{q^n}\), with \(\sigma(x)=x^q\).
Frobenius stability and root orbits describe coefficient descent and factor degrees; linear operators and cyclic modules describe conjugate relations, normality, and kernels; cyclic groups provide the setting for multiplicative orders, power classes, and norms.
For composition and periods, a compatible projection separates the whole action into a base orbit and return on one fiber.

A construction must also distinguish satisfying a relation from attaining its least degree or period.
A prescribed Frobenius annihilator can be realized using the corresponding divisor of \(T^n-1\) and a normal element. For a prescribed base-field permutation, a trace condition determines part of the parameters, while the trace kernel supplies the remaining freedom; see \hyperref[v98:par-prescribed-annihilator]{Prescribing the Frobenius annihilator} and \hyperref[v98:par-prescribed-trace-permutation]{Choosing parameters from a base-field permutation}.
A recursive construction further requires that each operation preserve the hypotheses needed at the next step.

The table lists objects, structural correspondences, and references; the cited definitions and theorems supply the hypotheses.
Entries concerning subfield representations and weighted relative traces retain the conventions \(E=k(\alpha)\), \(K_d=\F_{q^d}\), \(1<d<n\), \(d\mid n\), and \(m=n/d\).

Passing from cyclotomic cosets over the residue field to unramified $p$-adic lifts preserves root orbits while Witt carries and coefficient valuations record the precision lost on reduction. Equal-characteristic formal deformation instead retains coefficient-field operations and follows the variation of roots, traces and norms. Complex magnitudes and $p$-adic valuations of Gauss sums, and point-count recurrences and Newton slopes of Frobenius, provide complementary arithmetic data; see Sections~\ref{v132:sec-adic-lift}, \ref{v132:subsec-adic-gauss}, and~\ref{v132:subsec-newton-curves}.

\begingroup\small
\renewcommand{\arraystretch}{1.15}
\begin{longtable}{@{}>{\raggedright\arraybackslash}p{0.27\textwidth}>{\raggedright\arraybackslash}p{0.51\textwidth}>{\raggedright\arraybackslash}p{0.14\textwidth}@{}}
\toprule
Object&Structural correspondence and interpretation&Main reference\\
\midrule
\endhead
\bottomrule
\endfoot
Residue fields of unramified cyclotomic extensions&Prime-to-$p$ roots retain their order; transported addition recovers the finite field&\ref{revadic:thm-main}, \ref{revadic:prop-teich-projection}\\
Normal elements and normal integral bases&Every lift satisfies the same unit-determinant test; trace duality and finite-precision invariants commute with reduction&\ref{revadic:thm-normal}, \ref{revadic:prop-fixed-precision}\\
Additive Hilbert 90 at every precision&Relative-trace zero is equivalent to being a Frobenius difference; each nonempty fiber has size $p^{ds}$&\ref{v136:cor-integral-hilbert90}\\
Exact zeros and local divisibility&Integral cyclotomic remainders detect exact zeros; local factors test precision, and the least local valuation equals the coefficient-content valuation&\ref{revadic:thm-integral-resolvent}, \ref{revadic:thm-local-min}\\
Curve counts and slopes&The first $g$ counts and the functional equation determine every coefficient; Fermat Jacobi factors yield digit-sum slopes&\ref{revadic:thm-genus-counts}, \ref{revadic:thm-fermat-slopes}\\
Valuations for a fixed product character&Cyclic carry matrices record every Jacobi valuation; unit parts still determine cancellation&\ref{v150:prop-jacobi-carries}\\

\multicolumn{3}{@{}l}{\textbf{I. Fixed points, orbits, and descent}}\\*
Intermediate fields of \(E/k\)&Fixed-point sets of iterates; inclusion corresponds to divisibility of degrees&\ref{thm:subfields}\\
Identifying finite extensions of the same degree&After one isomorphism is chosen, all others arise by postcomposition with Frobenius; a global identification requires compatibility at every finite level&\ref{prop:finite-field-isomorphism-choices}, \ref{thm:absolute-frobenius}\\
Frobenius orbits of multiplicative cosets&Multiplication on cyclic quotient indices; return of the root set determines the coefficient field&\ref{v85:prop-coset-frobenius}\\
Root sets of monic nonconstant separable polynomials&Stability corresponds to coefficient descent; under this condition, transitivity corresponds to irreducibility&\ref{thm:root-descent}\\
Right factors of monic skew polynomials&Factor degree equals root-subspace dimension; irreducibility is detected by the return characteristic polynomial&\ref{v145:thm-skew-roots}, \ref{v145:prop-return-matrix}\\
Discriminants of squarefree polynomials&The sign of Frobenius on roots determines the square class, recovering only the parity of the number of irreducible factors&\ref{v105:prop-discriminant-parity}\\
Normalized multivariate linear factors&A single coefficientwise Frobenius orbit gives an irreducible product of conjugates; values at base-field points are norms&\ref{thm:multivariate-linear-descent}, \eqref{eq:linear-form-norm-value}\\
Square-free quotient algebras&Each root orbit contributes one field factor&\ref{thm:equivariant-functions}\\
Base-field coefficients and Fourier spectra&Coefficients descend exactly when spectral coordinates are Frobenius-equivariant; each orbit contributes its field component&\ref{v96:prop-fourier-descent}\\
Elements fixed by Frobenius in the Berlekamp algebra&The minimal polynomial records the actual values; the characteristic polynomial of multiplication and the resultant retain the weights given by the component dimensions&\ref{prop:berlekamp-value-polynomial}\\
Value fibers of a fixed element&Indicator polynomials for base-field values give orthogonal idempotent projections; equal values may still combine primitive components&\ref{v82:prop-value-idempotents}\\
Unknown field components of equal degree&Quadratic characters or absolute traces produce two-valued elements in the fixed algebra; greatest common divisors separate the values&\ref{v48:prop-edf-success}\\
Canonical output from randomized computations&Order the full candidate set in a fixed representation; growth of multiplicative orders under successive root extraction ensures a fixed irreducible output&\ref{v58:prop-canonical-nonsquare}\\
Root-finding using elliptic curves&Greatest common divisors partition roots by square classes of cubic values; deterministic root-finding also depends on group data for a curve and its twist&\ref{v105:rem-elliptic-root-finding}\\

Finite reduced algebras and homomorphisms&Finite sets with permutations; homomorphisms correspond to equivariant maps in the reverse direction&\ref{thm:finite-orbit-duality}\\
Several algebra generators&Exact tuple-orbit counts determine generator number and the count of ordered generating tuples&\ref{v107:prop-generating-tuples}\\
Frobenius fixed-space dimensions&Two integer inversions recover field components of every degree; the linear similarity type of Frobenius determines the isomorphism type of the reduced algebra&\ref{fra:thm-recovery}, \ref{fra:cor-similarity}\\
Base-field extension&Changes the step size along orbits; the intersection degree determines the number of components&\ref{thm:factor-base-extension}, \ref{thm:tensor-fields}\\
Local orbits for coprime moduli&Common iteration forms synchronized orbits; the phase quotient distinguishes the global factors&\ref{v44:thm-local-orbit-synchronization}\\
Phases and lifting branches&The phase kernel describes global branches; prime valuations determine its group structure and whether reduction splits&\ref{v125:prop-phase-lifting}, \ref{v126:cor-phase-splitting}\\
Lifting moduli and cyclotomic systems&Affine return determines branch lengths; compatible orbits give towers of irreducible factors&\ref{v44:thm-cyclotomic-affine-return}, \ref{v44:prop-cyclotomic-factor-tower}\\
Nonprincipal local cyclotomic systems&Deviation depth and local multiplicative order determine the stabilization index; a finite level then uniquely determines the tail&\ref{v46:thm-cyclotomic-stable-growth}\\
Stable towers of irreducible factors&Each level lifts uniquely; power substitution preserves irreducibility, with degrees growing by prime powers&\ref{v46:cor-cyclotomic-stable-factor-substitution}\\
Exact-period layers&Local orbit counts and the global Euler product determine each other&\ref{prop:orbit-euler}\\
Common degree of nonidentity roots of unity&\(N\mid\Phi_m(q)\) and \(\gcd(N,m)=1\) characterize a common degree for all nonidentity roots&\ref{v124:prop-uniform-cyclotomic}\\
Proportions of irreducible polynomials&Synchronized orbit products and a prime recursion give positivity; for fixed degree \(n\ge2\), the irreducible proportion strictly increases with field size&\ref{v105:prop-irreducible-density}\\
Stratification of exact-period layers by multiplicative order&Euler totient counts in the order strata agree with subfield M\"obius counts; base-field step size and inversion determine splitting and reciprocal pairing, respectively&\ref{v48:prop-period-order-count}, \eqref{v48:eq-period-reciprocal-order}\\
Full-degree elements of the norm-one group&Inverse pairs, order strata, and odd-divisor inversion give one count; each self-reciprocal orbit has $n$ pairs&\ref{v145:prop-inverse-count}\\
Root-of-unity orbits in residue fields&Teichm\"uller roots and Hensel uniqueness determine Galois-ring factors; the quotient ring decomposes as a product indexed by orbits&\ref{v47:thm-galois-cyclotomic-factorization}\\
Unramified infinite-precision lifts&The inverse limit of Galois rings retains cyclotomic orbits; coefficient valuations distinguish exact zeros from vanishing modulo $p$&\ref{v132:thm-unramified-limit}, \ref{v132:prop-coset-adic-factor}\\
Valuation profiles of canonical factors&Orbit subset-sum counts and Gauss periods give exact coefficients; Newton partitions bound valuations, while cancellation requires extra information&\ref{v134:thm-subset-periods}, \ref{v134:prop-newton-profile}\\
Local valuation multisets&Smith exponents repeat component valuations by the local degrees; minors and kernels modulo powers recover depth multiplicities&\ref{v143:thm-valuation-smith}\\
Recovery of labelled valuation profiles&Full fibers and fixed groups reduce matrix rank; the second-largest distinct valuation plus one bounds labelling precision&\ref{v144:thm-invariant-smith}, \ref{v144:prop-label-precision}\\
Cyclotomic cosets and periodic digits&The least base-$q$ period of $-a/m$ is the coset length; changing the step splits a $p$-orbit according to a greatest common divisor&\ref{v138:thm-periodic-cosets}, \ref{v138:prop-period-steps}\\
Periodic words and Gauss-sum valuations&Stickelberger reads digit sums from period blocks; the shortest base-$p$ block determines the valuation normalized by field degree&\ref{v138:cor-word-gauss}, \ref{v138:ex-seven-gauss}\\
Minimal rotations of periodic words&Divisibility restricts allowed words; minimal rotations determine leaders and narrow-sense BCH thresholds&\ref{v139:prop-digit-leaders}, \ref{v139:eq-bch-threshold}\\
Complete fibers and digit compression&Exact period loss under multiplication with carries equals block size; canonical factors are monomial compositions&\ref{v139:prop-digit-binomial}, \ref{v139:thm-period-compression}\\
Consecutive roots and carry steps&The largest omitted gap gives the fixed-root BCH bound; complete fibers give a direct sum of shorter codes&\ref{v139:prop-digit-bch-gap}, \ref{v139:cor-compressed-code}\\
Carry and multiplicative rigidity&The second digit determines linear maps preserving canonical representatives&\ref{v140:thm-teich-rigidity}\\
Digit subspaces and translations&Additive closure characterizes scalar multiples of subfields; nonzero carry directions distinguish odd characteristic from characteristic two&\ref{v147:thm-carry-directions}\\
Precision of low-degree representations&Smith exponents determine disappearing directions and kernel sizes&\ref{v140:thm-profile-smith}\\
Permutation symmetries of canonical codes&Projector valuations determine depth; all nonzero check columns give a complete group classification&\ref{v140:thm-code-group}\\
Lifting rigidity on partial root sets&Mixed-exponent remainders detect two-digit lifts; separation determines all subgroup-code permutations&\ref{v144:prop-partial-lift}, \ref{v144:cor-subgroup-code}\\
Witt carries and zero-sum obstructions&Teichmüller lifting is multiplicative; Witt carries reconstruct addition, and the first zero-sum obstruction detects divisibility modulo $p^2$&\ref{v133:prop-first-carry}, \ref{v133:prop-zero-sum-obstruction}\\
Equal-characteristic formal deformation&Constant embedding is a ring section; simple roots lift uniquely, while Artin--Schreier deformation fixes the trace and varies the norm&\ref{v133:prop-equal-hensel}, \ref{v133:ex-artin-schreier-deformation}\\
Canonical divisor lifts of orbit unions&Hensel uniqueness preserves all monomial composition exponents; the maximal exponent is unchanged at every coefficient precision&\ref{v47:prop-galois-med}\\
Composition exponents over finite local rings&Coprime Hensel lifting over the nilpotent maximal ideal preserves the maximum composition exponent at every precision&\ref{v125:prop-chain-composition}\\
Exact local cancellation depths&Stabilizers compress the local components; a resultant norm gives the sum of depths, and quadratic periods give direct valuation formulas&\ref{v137:prop-norm-depth}, \ref{v137:prop-quadratic-period-depth}\\
Affine exponents and complementary coefficients&Unit factors identify valuations across phases; negation and complements give profile symmetry&\ref{v141:prop-affine-profile}\\
Phase differences and index excess&The discriminant and first multiplication kernels detect termination of excess cancellation&\ref{v141:thm-index-budget}\\
Minimal excess and trace recovery&The last congruence layer consists of matchings or triangles; period trace duality recovers labels&\ref{v141:thm-minimal-excess}\\
Orbit returns and carry slopes&Jacobi sums give base-field orbit factors; exponent carry frequencies determine normalized slopes&\ref{v141:eq-fermat-orbit-factor}\\
Coefficientwise permutation lifts&Permutation of the residue field separates residue classes; a nowhere-zero derivative gives injectivity within each class over finite commutative local rings with nonzero maximal ideal&\ref{v105:prop-permutation-lifting}\\
Lengths divisible by the residue characteristic&Twisting indices and power substitutions give canonical pairwise comaximal blocks, retaining noncomaximal factorizations within each block&\ref{v47:prop-galois-ramified-blocks}\\
\addlinespace[0.7ex]
\multicolumn{3}{@{}l}{\textbf{II. Linear operators, cyclic modules, and duality}}\\*
Multiplication and Frobenius&The same tensor-product evaluation isomorphism gives diagonal action and cyclic shift, respectively&\ref{thm:conjugate-coordinates}, \eqref{eq:conjugate-operator-models}\\
Two root-evaluation coordinate systems&Lagrange interpolation gives a Cauchy node-change matrix, normalized by derivative weights&\ref{v97:prop-cauchy-node-change}\\
Field matrix models in the same dimension&A simultaneous basis change identifies the regular representations; matrix order and Frobenius conjugacy period record the two notions of periodicity&\ref{v48:prop-matrix-models}, \ref{v48:prop-two-matrix-periods}\\
Conjugate coefficients and trace correlation&The identity \(g_\alpha=h_\alpha h_\alpha^*\) connects normality with cyclic trace duality&\ref{v88:prop-conjugate-correlation}\\
Absolute Frobenius in normal coordinates&Entrywise powers and a circulant matrix give the semilinear action; the iterated norm product is the cyclic shift&\ref{v97:prop-absolute-frobenius-coordinates}\\
Affine fibers with Frobenius-invariant direction&Local residue coordinates give exact normal-element counts&\ref{v88:thm-invariant-affine-normal-count}\\
All annihilator data along an invariant direction&Fixed and free primary components give constant orders and polynomial Euler distributions&\ref{v107:thm-affine-annihilator-distribution}\\
Arbitrary linear constraints&Kernel intersections and inclusion--exclusion give the exact annihilator distribution&\ref{v110:prop-affine-count}\\
Prescribed cyclic correlation&Unit norm equations with involution realize symmetric trace correlations&\ref{v88:thm-prescribed-correlation}\\
Primitive normal elements&Character expansions jointly count the local conditions in the exponent group and Frobenius module; Gauss sums and a local sieve give existence&\ref{v48:prop-primitive-normal-test}, \ref{v93:thm-PNBT}\\
Weighted trace fibers of normal elements&One nonzero trace functional from each local component gives an affine space of normal elements of maximum dimension&\ref{v93:thm-patch}\\
Normality and primitivity in translates&Base-field translation can repair only the trace condition; primary multiplicities control subfield translates, after which primitivity is screened separately&\ref{v82:prop-normal-translates}; Exercise~\ref{v82:ex-subfield-translates}\\
Several prescribed intermediate traces&Normal coordinates turn traces into polynomial reductions; compatible congruences and unit lifting give exact counts of all solutions and of normal solutions&\ref{v58:thm-prescribed-normal-traces}\\
Independent components of complete normality&Pair-orbit coordinates and fibre sums determine intermediate-field obstructions; essential conditions and joint ranks yield exact probability factors&\ref{v131:thm-fibre-normal}; \ref{v131:thm-component-product}; \ref{v131:prop-component-exact}\\
Completely normal elements in restricted translation families&Base-field translations have an exact two-value count; intermediate-field translations are controlled by affine joint ranks and image compatibility&\ref{v131:cor-base-complete-count}; \ref{v131:prop-intermediate-translate-count}\\
Normal elements among reciprocals of translates&Coordinate greatest common divisors determine all exceptions; component ranks bound their degrees, and quartic discriminant constraints give a further improvement&\ref{v115:prop-primary-reciprocal}, \ref{v117:cor-primary-rank-bound}, \ref{v115:prop-quartic-bound}\\
Normality obstructions and local spectra&Trace orthogonal complements give detecting functionals and exact Fourier support; frequency intersections count normal elements in affine fibers&\ref{prop:normal-test}, \ref{v128:prop-local-spectrum}, \ref{v128:cor-affine-spectrum}\\
Primary components and row-space intersections&One irreducible primary component forces a chain; several components give disjoint idempotent images&\ref{v150:prop-row-intersection}\\
Annihilator strata and exact spectra&Subtracting adjacent primary layers gives the full spectrum; its dimension is the minimum amount of linear information&\ref{v110:prop-annihilator-info}, \ref{v130:prop-annihilator-spectrum}\\

Frobenius differences&The image is the trace kernel; the characteristic determines whether an equivariant splitting exists&\ref{prop:trace-splitting}\\
Root-set residues in quadratic quotients&The three factorization types give norm-one products, value ratios, and reciprocal sums; a nondegenerate functional fixes one projective unit class&\ref{v85:cor-rootset-quotient-coordinates}, \ref{v85:prop-functional-unit-fiber}\\
Coefficient pairings on quotient algebras&The derivative converts the coefficient pairing into the trace pairing; separability is equivalent to nondegeneracy of the latter&\ref{thm:algebra-trace-coefficient}, \ref{prop:trace-pairing-separability}\\
Trace observations and a specified linear target&A target factors through the observations exactly when their kernel lies in its kernel&\ref{v150:ex-three-traces}\\
Realizing invertible linear operators inside fields&Realization as Frobenius on a subspace is possible exactly when a cyclic vector exists&\ref{thm:frobenius-realization}\\
Base-field linearized polynomials&Under composition, factor multiplicities control splitting-field growth; for a fixed polynomial over varying extensions, factor orders and multiplicities determine every kernel&\ref{v68:cor-linearized-iterate-splitting}, \ref{v70:cor-all-extension-kernels}\\
Weighted Frobenius&The minimal polynomial of the relative norm determines the cyclic-module decomposition; the norm generates the fixed subfield exactly when a base-field cyclic vector exists&\ref{wfr:thm-norm-decomposition}, \ref{wfr:cor-cyclic-realization}\\
Linear operators&Combinations of multiplication and Frobenius; the trace pairing gives adjoints&\ref{thm:operator-algebra}\\
Linearized remainder interpolation&Restrictions to root subspaces must preserve all relations in their sum; independent prescription is equivalent to a direct sum&\ref{v105:prop-linearized-remainder-compatibility}\\
The trace pairing on an operator image&The adjoint kernel gives the radical of the restricted pairing; a rank difference measures its dimension&\ref{v68:cor-image-trace-radical}\\
Scalars and dilated intersections of subspaces&Root-polynomial support recovers the maximal field of scalars; kernels and ranks of quotient multiplication give intersection dimensions, and trace orthogonality recovers the scalar field from products&\ref{v70:prop-scalar-support}, \ref{v70:prop-quotient-multiplication}\\
Scalar intersections of power spaces&The degree bound lifts evaluated relations to identities; reduced fractions give dimensions and full frequencies&\ref{v150:cor-power-distribution}\\
Linear operators in conjugate coordinates&Dickson-matrix rank determines the number of roots; matrix multiplication and transposition correspond to composition and passage to the trace adjoint&\ref{prop:dickson-rank}, \ref{prop:dickson-operations}\\
Cycle decompositions of invertible linear operators&Cyclic-module multiplicities determine dimensions of fixed spaces of iterates; M\"obius inversion of fixed-point counts recovers all cycles&\ref{v53:lem-cyclic-decomposition}, \ref{v53:thm-linear-cycle-count}\\
General affine permutations&Translation height decides compatibility of iterates; elementary divisors count each compatible fiber&\ref{v107:thm-affine-cycles}\\
Signs of linear permutations&In odd characteristic, the determinant's square class determines the sign; in even characteristic, odd linear permutations occur only in dimension two over \(\F_2\)&\ref{v53:thm-linear-sign}, \ref{v53:cor-affine-sign}\\
Invariant subspaces&Correspond to divisors of \(T^n-1\); duality corresponds to reciprocal complementary factors&\ref{thm:invariant-duality}\\
Products of invariant subspaces&When \(p\nmid n\), multiplying eigenvalue sets in a common splitting field determines the product space; repeated factors require more than root sets&\ref{v70:thm-semisimple-product}\\
Subfield degree profiles&Related by a dimension-complement formula to the base-field ranks of weighted relative-trace prefixes; degree jumps and trace-rank jumps correspond in reverse order&\ref{str:thm-profile-trace-duality}, \ref{str:cor-reversed-jumps}\\
Product pairings and recurrence recovery&A nondegenerate pairing makes the Hankel matrix invertible; the observation kernel determines the quotient algebra and shortest recurrence order&\ref{v127:lem-pairing-recovery}, \ref{v127:prop-observation-kernel}\\
Finite prefixes of relative traces&The first \(2m\) trace values in fixed field coordinates recover the root block and profiles of divisor subfields; M\"obius inversion gives exact generator counts and minimum generator-representation degrees&\ref{str:thm-hankel-recovery}, \ref{v70:cor-trace-prefix-generators}\\
Recurrence solution spaces&When the recurrence polynomial has nonzero constant term, these spaces and the corresponding linearized root spaces are isomorphic as cyclic modules; shift corresponds to Frobenius&\ref{prop:recurrence-root-correspondence}\\
Cyclic codes of trace sequences&Cyclic right shift corresponds to multiplication in an extension; a single root orbit retains one field component&\ref{v40:prop-cyclic-trace-model}\\
Projective roots and evaluation&Linear factors record rational directions; the full root form generates the homogeneous evaluation kernel&\ref{v84:prop-projective-evaluation}\\
Interpolation and duality&The highest-coefficient identity gives derivative weights and the coordinate at infinity&\ref{v84:prop-rs-dual}\\
Frobenius derivation&Characteristic powers give the low-degree kernel; at degree \(q+1\) the kernel equals the first homogeneous vanishing relations of projective rational points&\ref{v86:thm-frobenius-kernel}\\
Constacyclic MDS and reciprocal multiples&Full-support generator and check polynomials; reciprocal check multiples detect dual distance&\ref{v87:thm-constacyclic-mds-reciprocal}\\
Trace codes and sparse relations&MDS column independence excludes low-weight multiples and forces normality of nonidentity subgroup elements&\ref{v84:thm-sparse-multiple}\\
The relation kernel of a trace code&The inverse-root relation map has the dual code as kernel; shortest relation length is dual distance&\ref{v96:prop-relation-kernel}\\
Complete orbit fibers&The orbit degree attains the fiber bound exactly when the minimal polynomial is a power substitution&\ref{v142:prop-full-fibers}\\
Independent blocks of an irreducible code&Complete fibers make relative-trace parameters independent, giving a direct sum and equal minimum distance&\ref{v142:cor-irreducible-interleaving}\\
Subfield intersections and distance loss&The image and kernel of subgroup-coordinate restriction bound the code distance and dual distance&\ref{v142:prop-subfield-distance}\\
Orbitwise normality versus joint independence&Normal orbits and full-support minimal polynomials need not exclude a cross-orbit three-term relation&\ref{v96:ex-nineteen}\\
The projective image of a subgroup&With no scalar repetitions, maximal general position is equivalent to an arc, an MDS trace code, and a hyperplane bound&\ref{v96:cor-arc-trace-equivalence}\\
Layers of Fourier support&Proportional repetitions and duality give balance; the zero residue class inherits the MDS property&\ref{v84:thm-balanced-support}\\
Coefficient fields and root-set return&The return period under prime-field Frobenius determines the smallest field in fixed coordinates&\ref{v84:thm-code-field}\\
Multipliers and arithmetic synchronization&Orbitwise difference subcodes restrict fixed points and synchronize orders modulo divisors&\ref{v84:cor-multiplier-orders}\\

Common invariance of cyclic codes&Preserving every root component is equivalent to commuting with every Fourier spectral projection; orbit sums turn this into preservation of difference orbits&\ref{v73:thm-common-permutation-orbits}\\
The common permutation group&Containment of unit reduction kernels is equivalent to equality of multiplier indices; outside three local exceptions, shifts and Frobenius multipliers generate the common affine group&\ref{v73:thm-common-permutation-arithmetic}\\

Exponential moments and Hankel rank&Distinct nodes yield an exact Hankel determinant and minimum recurrence order&\ref{v88:thm-hankel-vandermonde}\\
Period sums of roots of unity&The first \(2d\) orbit power sums recover the irreducible factor through a Hankel system, even when the characteristic divides \(d\)&\ref{v97:cor-period-recovery}\\
Power-sum prefix identification&A nonzero second derivative gives a singleton fiber; zero second derivative gives an affine line whose points are distinguished by the next power sum&\ref{v123:thm-power-trace-pencil}\\
Consecutive roots and codeword support&Frobenius closure determines the base-field generator polynomial; Vandermonde nondegeneracy gives the BCH distance bound&\ref{v56:thm-bch-bound}\\
Equality in the BCH bound&Interpolation weights characterize extremal supports; matching power-trace prefixes of two root orbits give signed words attaining the bound&\ref{v123:prop-bch-equality}, \ref{v123:cor-paired-bch-word}\\
Syndromes and error recovery&Hankel rank gives the number of errors, the shortest recurrence gives reciprocal error locations, and Vandermonde equations recover their values&\ref{v118:prop-syndrome-recurrence}, \ref{v118:ex-bch-fifteen-decoding}\\
Determinants of two kinds of powers&Consecutive powers use distinct nodes; Frobenius powers use the Moore criterion for base-field independence&\ref{v96:rem-vandermonde-moore}\\
Column normalization of code coordinates&Nonzero column scaling preserves Hamming distance, but independent extension-field scalars may change rank distance&\ref{v105:ex-rank-column-normalization}\\
Supports in binary Goppa codes&Divisibility of the support derivative by \(G^2\) gives the distance bound; a splitting condition characterizes equality&\ref{v82:prop-goppa-support}; \ref{v82:cor-goppa-equality}\\
Derivative-constrained multiplicative-coset root sets&The polynomial Euler function counts algebraic parameters for a fixed derivative; power-subgroup splitting identifies the actual root sets&\ref{v83:prop-derivative-fibers}; \ref{v83:ex-sixtyfour-fibers}\\
Cyclic-shift coverings&A hyperplane normal vector converts covering into the absence of full-weight codewords in a nonzero cyclic code&Exercise \ref{v82:ex-cyclic-covering}\\
Arbitrary root constraints&Finite-field Fourier submatrices can be singular; consecutive exponents supply the Vandermonde nondegeneracy used by the BCH bound&\ref{v68:rem-singular-fourier-minor}\\
Support periods of generator polynomials&Exponent residue classes are the indecomposable coordinate blocks; interleaving preserves Hamming weight, and the compressed code gives a distance bound&\ref{v67:prop-coordinate-blocks}, \ref{v67:cor-core-distance}\\
Affine semi-invariants and code spaces&Outer divisibility retains special-orbit multiplicities; the remainder-matrix kernel gives the subcode dimension, and coefficient-row-space tests determine effective support&\ref{v54:thm-outer-divisibility}, \ref{v54:prop-test-subcode}\\
Outer divisibility and quotient coordinates&The coefficient pairing gives a recurrence matrix with the same kernel; Hasse derivative tests retain root multiplicities and give the same rank&\ref{v67:cor-hankel-divisibility}, \eqref{v67:eq-hasse-divisibility}\\
Initial parts and tails of general recurrences&Zero-root multiplicity gives the transient length and the nonzero-root part gives the eventual period; the polynomial part and denominator of the reduced generating function recover these parts&\ref{thm:transient-period}, \ref{v70:prop-rational-generating-annihilator}\\
Local multiplicities in quotient algebras&The trace is a multiplicity-weighted sum of residue-field traces; the determinant of multiplication is a product of residue-field norms raised to their multiplicities&\ref{prop:local-multiplicity-trace}, \ref{prop:local-multiplicity-determinant}\\
Finite algebras with nilpotents&The stable Frobenius kernel is the nilradical and its stable image realizes the reduced quotient; local dimensions over residue fields determine the powers in characteristic polynomials, the trace weights, and the norm exponents&\ref{prop:algebra-periodic-reduction}, \ref{v117:prop-local-dimension-trace}\\
Quotients of recurrence, quotient-algebra, and root-space models&The coefficient-functional identity makes \(v(S)\), natural reduction, and \(\mathcal L_v\) compatible&\ref{prop:quotient-model-compatibility}\\
Finite sequence observations&Discrepancies give complexity lower bounds; a sufficiently long prefix recovers the minimal polynomial&\ref{thm:berlekamp-massey}\\
\addlinespace[0.7ex]
\multicolumn{3}{@{}l}{\textbf{III. Fiber return, symmetry, and composition}}\\*
Changes of generator and defining equations&Equivariant orbit images record minimal polynomials and degree loss; complete equal-difference exponent orbits give binomials of minimum support size&\ref{v79:prop-generator-orbit-image}, \ref{v79:rem-med-defining-equation}\\
Sparse relations for a fixed defining polynomial&Remainder columns retain the exponent range; repeated adjacent ratios give simultaneous cancellation by a linear multiplier&\ref{v97:def-sparse-degree}, \ref{v97:prop-linear-cancellation}\\
Parallel difference quotients&Independence of the base point is equivalent to prime-field linearity after subtracting the constant term; a linearized polynomial gives the unique representative&\ref{v72:thm-parallel-difference}\\
Multiplicative cosets of field differences&A proper subgroup permits passage from low-degree function identities to polynomial identities, recovering Frobenius and the affine coefficients&\ref{v72:thm-carlitz-mcconnel}\\
Additive stabilizer of a multiplicative orbit&The complete-coset condition \(C+I=C\) is equivalent to the multiplicative containment \(1+I\subseteq H\)&\ref{v74:prop-orbit-additive-stabilizer}\\
Global relabeling of cosets&The multiplicative subgroup generates a hidden subfield \(K\); the permitted maps are \(K\)-semilinear automorphisms followed by translations&\ref{v74:rem-hidden-subfield}\\
Norm-preserving linear maps&Additivity transfers norm preservation to differences; the norm kernel recovers norm-one scalars and relative Frobenius&Exercise~\ref{v74:ex-norm-preservers}\\
Mixed exponent constraints on differences&An invertible power map on the quotient makes all difference cosets transform by one fixed group automorphism&Exercise~\ref{v74:ex-exponent-rigidity}\\
Equal-difference structure and non-affine symmetries&At odd prime-power moduli, the \(p\)-part of the multiplier order determines both maximal blocks and exceptions to local rigidity&\ref{v74:thm-med-nonaffine}\\
Cyclotomic differences over residue rings&Chinese remainders isolate independent local multipliers; prime-power rigidity and three local perturbations give an exact affine-rigidity criterion&\ref{v73:thm-ring-rigidity}\\
Multiple equal-difference structure of cyclotomic cosets&Maximal equal-difference blocks give the monomial composition exponent; block periods give the least extension degree for binomial factorization&\ref{thm:med-minimal-step}, \ref{prop:med-binomial-descent}\\
Coset leaders and order-preserving filtering&Scaling preserves orbit lengths and least representatives; maximal equal-difference blocks reduce the leader search to the set of blocks&\ref{v50:thm-leader-scaling}, \ref{v50:prop-med-leader-compression}\\
Stable growth of primitive root-of-unity layers&A common multiplicative period and valuation thresholds determine a finite base level; all higher-exponent irreducible factors arise by power substitution&\ref{v50:prop-cyclotomic-uniform-stabilization}\\
Coefficient fields and blocks of multiplier orbits&The Frobenius period of a root block gives its exact coefficient field; the group index counts irreducible factors in a base-field block&\ref{v50:prop-multiplier-blocks}\\
Three-sparse root-of-unity factorizations&At most two equal-difference root blocks give a trinomial shape; squaring the roots excludes additional prime factors in odd characteristic, while characteristic two retains an orbit-restricted family of powers of seven&\ref{v49:thm-odd-three-sparse}, \ref{v49:thm-binary-three-sparse}\\
Trace parameters of roots of unity of low order&Inverse conjugacy gives quadratic factors; trace recurrences and multiplicative-order growth recover their irreducible power substitutions&\ref{v49:prop-low-order-cyclotomic}\\
Root power sums of unions of orbits&The greatest common divisor of the modulus and the frequency indices with nonzero power sums determines the additive stabilizer of the exponent set and the maximum monomial composition exponent of its root polynomial&\ref{v40:prop:root-power-composition}\\
Root ratios and multiplicities&The least period of a common root power and the characteristic-power multiplicity together determine the minimal binomial multiple&\ref{v40:prop-minimal-binomial}\\
Binomial factor degrees and splitting fields&Norm classes constrain factor degrees; support periods constrain the scalar action remaining on the roots&\ref{v97:prop-norm-factor-degrees}, \ref{v97:prop-support-splitting}\\
Least extension for binomial factorization of a general binomial&\(\ord_{b(ng)}(q)\) controls all root orbits simultaneously, where \(g=\ord(\lambda)\)&\ref{v124:cor-binomial-extension}\\
Ratios of quadratic conjugate roots&The first zero term of positive index in a Lucas recurrence recovers the projective order and minimal binomial multiple&\ref{v45:prop-lucas-first-zero}\\
Quadratic roots with prescribed primitive norm&Cyclotomic atoms detect the exact order of the root ratio; the unique square root in characteristic two gives a square-free optimal coefficient criterion&\ref{v45:thm-optimal-quadratic-criterion}\\
Characteristic-power indices in Lucas recurrences&The Frobenius identity separates periods of separable root ratios from zeros caused by repeated roots&\ref{v151:cor-lucas-frobenius}, \ref{v151:ex-lucas-characteristic}\\
Lucas atoms and exact-order traces&Inverse pairing in the norm-one group pulls back to a coefficient criterion; the same trace polynomial prescribes return-eigenvalue order&\ref{v146:prop-lucas-pullback}, \ref{v146:cor-return-order}\\
Unramified lifts of Lucas atoms&In characteristic two, opposite roots coalesce modulo two; root differences have valuation one and the discriminant valuation equals the atom degree&\ref{v146:prop-adic-lucas}\\
Quadratic root sets in families of additive cosets&The two Frobenius eigenspaces distinguish reducible and irreducible members; norm fibers give exact counts&\ref{v45:prop-quadratic-coset-count}\\
Conjugate root pairs&A square parameter records each root pair; the relative norm descends to a function of that parameter&\ref{v45:prop-quadratic-norm-descent}\\
Root-pair phases and composed images&Image periods determine within-orbit multiplicity; equal image orbits merge factors; coprime root orders ensure injectivity&\ref{v129:prop-phase-images}, \ref{v40:prop-coprime-root-orders}\\
Equivariant root fibers&Total period is base period times return period; base extension changes both the base orbit and the internal return&\ref{lem:fiber-return}, \ref{v129:cor-return-base-change}\\
Several prescribed norms&Norm compatibility on intersections is equivalent to compatibility of integer congruences; a least common multiple determines the joint multiplicative fiber size&\ref{v58:prop-compatible-norms}\\
\(v\)-th powers of subspace polynomials, \(v\mid q-1\)&Equal-sized nonzero fibers from additive cosets and the base-field group of roots of unity give a minimal value set; the zero fiber is treated separately&\ref{v58:prop-subspace-power-values}\\
Frobenius difference after a norm&Multiplicative norm fibers and additive fibers of the Frobenius difference give a minimal value set; the zero fiber requires a separate count&\ref{v105:prop-norm-trace-values}\\

Additive and multiplicative fiber permutations&The AGW criterion separates bijectivity on the base set from injectivity within fibers; explicit return on multiplicative and trace fibers determines cycle lengths&\ref{prop:multiplicative-permutation-return}, Exercise~\ref{v70:ex-additive-return}\\
Sparse nonlinear involutions&A trace involution on a subfield lifts through a multiplicative quotient, preserving an explicit inverse&\ref{v68:ex-sparse-involution}\\
Affine envelopes of square-free polynomials&The linear span of root differences gives the least monic affine \(q\)-polynomial multiple&\ref{prop:affine-q-multiple}\\
Root-field spectra of composite polynomials&Base-orbit length times relative degrees of roots in the fiber; for additive fibers, affine averaging ensures a least root field, and an image condition determines its degree&\ref{prop:composition-root-spectrum}, \ref{prop:affine-least-period-computation}, \ref{cor:additive-least-root-field}\\
Root multiplicities and translation symmetry&The greatest common divisor of difference coefficients recovers the translation subspace and the common \(p\)-power&\ref{ffs:thm-psi}\\
First-order affine symmetries&Roots whose multiplicities are nonzero in the base field give linear constraints; the rank of the algebra trace pairing determines the dimension of the symmetry space&\ref{v67:prop-first-order-affine}\\
Full affine root orbits&Descent of the center coset gives a standard inner polynomial; the special-orbit multiplicity is retained separately&\ref{ffs:thm-normal-form}\\
Polynomials with a fixed scaling factor&Semi-invariants consist of powers of the center polynomial in a specified residue class; a degree bound gives a finite basis&\ref{fas:prop-semi-invariant-basis}\\
Common support of affine semi-invariants&For a nonzero constant term, the full affine group gives a support bound at least as sharp as the support-period bound&\ref{v67:cor-affine-support-bounds}\\
Conjugate differences in irreducible root sets&A specified Frobenius difference detects nonzero translations and identifies the four-element exception in characteristic two&\ref{fas:thm-translation-difference}\\
Affine symmetries of translation quotients&Slopes lift exactly when they preserve the translation space; base-field centers and support-index differences recover outer slopes, and root fibers give all lifts&\ref{fas:prop-compatible-lifting}, \ref{v70:prop-outer-center-recovery}\\
Nontrivial symmetries in prime degree&Artin--Schreier type when the degree equals the characteristic, and a translated prime-degree binomial otherwise&\ref{fas:thm-prime-characteristic}, \ref{fas:thm-prime-tame}\\
Root fibers of affine quotients&Frobenius conjugacy on parameters and the return transformation together determine factor degrees and multiplicities&\ref{ffs:thm-full-return}\\
Low-degree subfield representations&The minimum generator-representation degree equals the relative extension degree exactly when the corresponding polynomial composition exists; equivalently, the representation defect is zero&\ref{sfr:thm-optimal-composition}, \ref{sfr:thm-zero-defect}\\
Lower bounds for rational representation degrees&\(\kappa_f^{\mathrm{rat}}(d)=n/d\iff\eta_f(d)=2\); under this condition, the coordinate space contains constants exactly when polynomial representations also attain the lower bound&\ref{str:thm-rational-rank}, \ref{str:cor-polynomial-direction}\\
Fibers of recurrence quotients&A linearized map realizes quotient states; affine return determines the internal period of an original state&\ref{prop:recurrence-fiber-return}\\
Power and affine compression&The image period times the fiber return period gives the original period; the return map must be determined in the respective fiber coordinates&Section\allowbreak~\ref{sec:power-affine-comparison}\\
Affine compression&The compressed period is the period on the orbit quotient; the period-reduction factor is the order of the affine return transformation&\ref{prop:affine-return-order}\\
Return order for a prescribed compressed value&The norm determines the scaling coefficient; the orbit sum distinguishes identity return from nontrivial translation&\ref{prop:explicit-affine-return}\\
Compressed factors and affine orbits&Each compressed factor corresponds to an affine orbit of irreducible polynomials; return order is the base-field affine stabilizer order&\ref{cor:compressed-fiber-lifting}, \ref{v129:cor-affine-factor-orbits}\\
Affine Frobenius&The norm controls linear fixed points, and the trace decides whether translation can be removed; fixed-point counts recover all cycles&\ref{afr:thm-fixed-counts}, \eqref{afr:eq-cycle-count}\\
Periods of affine Frobenius maps with nonzero translation trace&When the norm is one and the translation trace is nonzero, cycle lengths \(\ell\) satisfy \(v_p(\ell)=v_p(\ord(S_a))+1\); fixed-point counts give the length and number of shortest cycles&\ref{afr:cor-period-jump}, \ref{afr:cor-shortest-cycles}\\
Power maps on inverse-root orbits&Dickson polynomials induce power maps on conjugacy parameters; two groups of roots of unity jointly control permutation behavior&\ref{thm:dickson-permutation}\\
Reciprocal root-of-unity orbits and a Frobenius equation&Two reciprocal-orbit root sets factor \(Y^q-Y^2+2\) and determine the derivative weights for interpolation&\ref{v97:thm-dickson-frobenius}, \ref{v97:cor-dickson-weights}\\
Reciprocal pairing of irreducible factors&Negation of exponents gives an orbit involution; for root order greater than two, self-reciprocity requires half-period inversion and the corresponding norm-one condition&\ref{v48:cor-reciprocal-order}\\
Persistent iteration of the R-transform&For odd \(q\), the initial polynomial is monic irreducible with \(P(1)P(-1)\) nonsquare, and has even degree if \(q\equiv3\pmod4\). Preservation of the square class of a product of base-field values ensures that each return swaps the two roots&\ref{thm:r-transform-iteration}, \ref{v48:ex-r-transform-tower}\\
Polynomial iteration with a fixed target&Irreducibility at every level is equivalent to irreducibility of each new root fiber; affine conjugacy must transform the target as well&\ref{v55:cor-iterated-capelli}, \ref{v55:prop-pointed-conjugacy}\\
Critical orbits of quadratic iteration&Norms of fiber parameters give base-field nonsquare conditions; orbit repetition makes the test finite&\ref{v55:thm-quadratic-stability}, \ref{v55:cor-quadratic-finite-test}\\

Normalization of general quadratics&Affine conjugacy transforms the target as well; the discriminant and the critical-orbit values weighted by the leading coefficient give the square-class criterion&\ref{v56:cor-general-quadratic-stability}\\
Cubic iteration in characteristic two&Reciprocal coordinates lift the fiber to a cube-root extension; the norm preserves noncubicity&\ref{v55:thm-char2-target-cubic}\\
Mixed compositions of odd-degree translated powers&The smallest forward-invariant subset of the base field containing the constant terms controls power-residue conditions for all compositions&\ref{v55:thm-mixed-power-stability}\\
Inversely stable Artin--Schreier families&Denominators of rational iterates agree with reciprocal additive recursion; base-field extension requires \(p\nmid e\)&\ref{v55:prop-inverse-as}\\
\addlinespace[0.7ex]
\multicolumn{3}{@{}l}{\textbf{IV. Character expansions and counting}}\\*

Restricted coefficients in translation orbits&Translation preserves irreducibility and discriminant; derivative values turn coefficient-coset conditions into joint character sums&\eqref{v105:eq-translate-coset-count}\\
Irreducible quadratics with square coefficients&Quadratic character sums give an exact count; zero discriminant supplies the correction at small fields&\ref{v105:prop-square-coefficient-quadratics}\\
Power-residue symbols&The norm equals a resultant of conjugate values; exchanging the two root sets gives polynomial reciprocity&\ref{prop:residue-norm-resultant}, \ref{thm:polynomial-reciprocity}\\
Linear equations&Trace orthogonality to the image is equivalent to the Fourier-support condition&\ref{prop:fourier-adjoint}\\
Supports of a function and its Fourier transform&The product of support sizes is at least the group order; equality holds precisely for nonzero scalar multiples of subgroup-coset indicators times characters&\ref{v90:thm-fourier-uncertainty}\\
Trace-evaluation parameters&Relative trace to the field generated by power values removes coefficient redundancy; the orbit degree gives the parameter-space dimension&\ref{v68:prop-trace-evaluation-parameters}\\
Translations of the norm-one group&Intersection bounds identify the unique exception to additive indecomposability; exact counts give additive energy and a Fourier fourth moment&\ref{v150:ex-norm-one-indecomposable}, \ref{v151:prop-normone-intersections}, \ref{v151:cor-normone-moment}\\
Short zero sums in a norm-one group&A nonzero sum determines an unordered root pair; triple and quadruple relations determine binary distance bounds&\ref{v153:prop-normone-short-sums}, \ref{v153:cor-normone-binary-code}\\
Value sets and fiber collisions&Pairs with the same image give a lower bound for the value set; orthogonality expresses the collision count as the mean squared absolute value of character sums&\ref{v58:prop-value-collision}, \ref{v58:prop-collision-fourier}\\
Exact value distributions from intersection fields&A reciprocal coordinate turns nonzero fibers into affine linearized equations; the intersection field controls their kernel sizes&\ref{v68:prop-affine-values}\\
Even moments of character sums&Counts of solutions to power-sum equations determine moments; when the characteristic exceeds the degree, Newton identities recover root multisets from power sums&\ref{v58:eq-moment-power-systems}, \ref{v58:prop-power-sum-multiset}\\

Quadratic forms&Polarization corresponds to a sum of adjoint operators; the radical and type determine the value distribution&\ref{prop:trace-quadratic-operator}\\
Quadratic phase functions&The polarization radical determines the direction and amplitude of Fourier support; the values of the quadratic form on that radical determine its translation&\ref{prop:quadratic-fourier-support}\\
A field-valued quadratic map&In characteristic two, adjoint incidence relates differential kernels to the polar radicals of trace components; double counting yields a first identity for their dimension distributions&\ref{v70:prop-dual-kernel-incidence}\\
Lifting characters to extensions&Trace and norm determine orbit weights; the Euler product gives Gauss-sum lifting&\ref{thm:hasse-davenport}\\
Fibers of irreducible quadratic returns&Normalized root bases and character sums give uniform sizes, determining cubic trace--norm fibers and sparse skew irreducible counts&\ref{v145:thm-return-uniform}, \ref{v145:thm-sparse-count}\\
Return spectra for arbitrary Frobenius generators&Semilinear normal forms and power-map fibers preserve irreducible-target counts; a linear equation counts split distinct-root targets&\ref{v146:lem-semilinear-normal}, \ref{v146:thm-split-fibers}\\
Repeated returns and constant terms&Target totals recover repeated-root fibers; a power-image condition determines constant-term normalization&\ref{v146:cor-repeated-fibers}, \ref{v146:prop-return-scaling}\\
Return norms and common subfields&Local lattices count all fibers; relative norms control constants and give uniformity for coprime degrees&\ref{v148:thm-all-fibers}, \ref{v148:thm-relative-fibers}\\
Even degrees over a quadratic coefficient field&Hankel kernel shifts remove singular phases; factor constants give the deviations, and cyclic phase divisibility decides uniformity&\ref{v148:thm-all-even-fibers}, \ref{v149:prop-uniformity}\\
Return orders and ordinary factorization&Local partitions recover orbit lengths; primitivity is transitivity on nonzero roots&\ref{v148:thm-ordinary-return}, \ref{v148:thm-even-order}\\
Primitive coefficient criteria&Lucas pullback cancels the excess zero-root multiplicity; optimal degrees vary across parameter strata&\ref{v148:thm-optimal-pullback}\\
A fixed middle coefficient&Cubic right roots give a value-set complement; trace--norm weights recover slices and moments&\ref{v148:thm-cubic-value-set}, \ref{v148:thm-cubic-primitive-slices}\\
Equally spaced powers&Norm transfer separates target selection from fiber weights; root orders give existence and primitivity obstructions&\ref{v148:thm-sparse-coprime}, \ref{v148:cor-sparse-order}\\
Products of Gauss sums&Power-map fibers and Newton identities compute the same leading coefficient&\ref{v91:thm-dh-product}\\
Square-root bounds for character sums&Linear conditions create multiple roots; extension-field power sums turn rough bounds into bounds on characteristic roots&\ref{v91:thm-additive-weil}, \ref{v91:thm-multiplicative-weil}\\
Primitive elements with prescribed trace&Expand the indicators of primitivity and of the trace condition separately; restricted characters determine the error scale&\ref{thm:primitive-prescribed-trace}\\
Recurrences for extension-field character sums&Degree cancellation in a weighted orbit product gives a finite-order recurrence; Kloosterman sums reduce to order two&\ref{cor:character-sum-recurrence}, \ref{thm:kloosterman-recurrence}\\
Polynomial phases after taking traces&Reduction modulo \(H^p-H\), with the character parameter included, preserves phases over every extension; constant, linear, and higher-degree reductions give different forms of the \(L\)-series&\ref{v70:thm-phase-reduction}, \ref{v70:cor-reduced-phase-L-series}\\
Extension-field point counts for \(E_b\) in characteristic two&Trace fibers give a Kloosterman sum; the first point count determines a second-order recurrence and the zeta function&\ref{prop:binary-curve-point-count}, \ref{cor:binary-curve-zeta}\\
Two-point coefficient fibers&Fractional-linear involutions turn coefficient conditions into subgroup intersections; fixed-point correction imposes distinctness&\ref{v85:prop-pair-involution}\\
Sums and products of distinct subgroup elements&Jacobi expansions with diagonal removal; cyclic exponent translations for prescribed products and constant terms&Exercises \ref{v85:ex-restricted-sums}, \ref{v85:ex-distinct-products}\\
Subset products and codeword symbols&Multiplicative Fourier coefficients decide uniformity; deleting one point obstructs a minimum-weight three-design&\ref{v150:thm-ers-design}\\
Solution counts and errors for diagonal equations&Character expansions of power-map fibers give exact counts; the product character determines the two scales of Jacobi sums&\ref{prop:multiple-jacobi-gauss}, \ref{thm:diagonal-weil-bound}\\
Jacobi data and field addition&A faithful character row recovers addition; a partial table retains joint coset counts&\ref{v150:prop-jacobi-recovery}\\
Reciprocal ranks of subspaces&Clearing denominators retains rational-function relations; intersection with the minimal-polynomial ideal measures the exact loss under evaluation&\ref{v117:prop-specialization-loss}, \ref{v116:prop-rational-specialization}\\
Small positive representation defects&Truncated multiplication forces arithmetic progressions of jumps; residual factors and fiber returns realize the first positive layer&\ref{v116:thm-small-defect}, \ref{v116:prop-defect-construction}\\
Realizing subfield jumps&Multiplicative closure of a reduced degree basis is equivalent to Frobenius fixed-point equations; binary quartics give a complete small classification&\ref{v137:thm-jump-realization}, \ref{v137:ex-quartic-classification}\\
Minimal constant-derivative completion&Inverse-Frobenius interpolation gives a degree gap; splitting of residual factors decides completion&\ref{v116:prop-completion-gap}, \ref{v116:prop-four-completion}\\
Local stability of difference quotients&Agreement with a model yields a lower bound on bad pairs; global recovery requires a further argument&\ref{v116:prop-local-stability}, Exercise~\ref{v116:ex-stability-gap}\\
\end{longtable}
\endgroup

\noindent\textbf{Fields, orbits, and linear subspaces.}
The fixed points of a Frobenius iterate form a subfield; a finite Frobenius-stable root set has its monic root polynomial over the base field; a Frobenius-invariant \(k\)-subspace of \(E\) corresponds to a monic divisor of \(T^n-1\).
These correspondences use, respectively, preservation of field operations, root-set stability, and the cyclic-module structure of Frobenius.
A normal element gives a module isomorphism \(k[T]/(T^n-1)\simeq E\), with \(T\) acting as \(\sigma\), generally not a ring isomorphism.
Field multiplication supplies further comparisons: the scalars preserving a subspace form a subfield, and ranks of multiplication maps into quotient spaces determine dimensions of intersections with dilates.
Trace duality expresses these conditions as orthogonality; for semisimple Frobenius, eigenvalue products in conjugate coordinates determine products of invariant subspaces.

\noindent\textbf{The root bound and generation.}
The root bound excludes low-degree polynomial relations on an entire finite field or an excessively large subspace.
This exclusion underlies the proofs of multiplicative cyclicity, trace nondegeneracy, the Moore criterion, and the normal basis theorem; see Section~\ref{sec:root-bound-rigidity}.
It constrains the possible relations; Frobenius preserves the operations that give the relevant sets their field or linear-space structure.

\noindent\textbf{Three degrees and periods.}
For \(\alpha\ne0\), set \(m=\ord(\alpha)\), and let \(\mu_\alpha\) be its Frobenius annihilator polynomial. Then
\[
[k(\alpha):k]=\min\{j\ge1:\alpha^{q^j}=\alpha\}
=\ord_m(q)=\ord(\mu_\alpha).
\]
Here \(\ord_1(q)=1\). On the other hand, \(\deg\mu_\alpha\) is the dimension of the span of its conjugates,
which need not equal the field extension degree. Normality requires \(\mu_\alpha=T^n-1\), while multiplicative primitivity requires \(m=q^n-1\).

\noindent\textbf{Subfields, representations, and relative traces.}
Fix \(E=k(\alpha)\). The subfield \(K_d\) is determined by a Frobenius fixed-point condition,
but its intersections with low-degree power-basis spaces depend on the choice of generator.
Weighted relative traces express these intersections as orthogonal complements of trace subspaces in \(K_d\); a shortest recurrence of the corresponding trace sequence recovers the root-block polynomial of \(\alpha\) over \(K_d\).
Thus base-field rank records representation degree, recurrences over the subfield give relative algebraic relations, and root-block coefficient rank determines whether rational representations attain the degree lower bound.
Intermediate traces then give the degree profiles of divisor subfields; M\"obius inversion separates counts of elements lying in a subfield into counts of elements generating it.

\Needspace{16\baselineskip}
\noindent\textbf{Pairings, observations, and error correction.}
For a monogenic quotient algebra, nondegeneracy of \(\ell(uv)\) preserves the multiplication relation of its generator in \(\ell(t^j)\); a degenerate pairing retains the quotient by the observation kernel.
Relative-trace and syndrome recovery use two coordinate descriptions: a power basis over the subfield, and evaluations at error locations over a splitting field.
Hankel equations recover the recurrence coefficients at once, while Berlekamp--Massey corrects a connection relation term by term; both need an order bound to pass from a finite prefix to the entire sequence. If ordinary trace removes a nilpotent direction, longer observations cannot recover it; a functional that detects that direction is required.

\noindent\textbf{Local obstructions, candidate families, and spectra.}
Normality requires avoidance of \(U_P=P(\sigma)E\). The operator \((F/P)(\sigma)\) detects this obstruction, while \(V_P\) parametrizes the corresponding trace conditions. Substitution of a reciprocal family and removal of denominators give coordinate greatest common divisors locating exceptional parameters. Character expansion of the local indicators gives exact support \(\bigoplus_PV_P\). The former identifies successful candidates; the latter intersects the same conditions with affine fibers or multiplicative generation. Kernels and ranks already suffice for linear fiber counts; Gauss sums provide further control of the error once multiplicative conditions are present.

\noindent\textbf{Images, returns, and repeated constructions.}
A root-set image can shorten one orbit or merge distinct orbits; the minimal polynomials of phase images and their exponents record these losses. Inverse images require the internal return above each base orbit, and recursion further requires preservation of the traces, norms, or values controlling that return. Affine compression expresses the same problem through group orbits of irreducible polynomials, whereas Euler products retain only cycle lengths and counts. Explicit factors require root coordinates, and multiplicities require the fiber polynomials.

\noindent\textbf{Frequently used formulas.}
For monic irreducible \(f\in k[X]\) with \(\deg f=d\),
\[
f\mid X^{q^n}-X\iff d\mid n,
\qquad I_q(n)=\frac1n\sum_{e\mid n}\mu(e)q^{n/e}.
\]
Over \(\F_{q^r}\), it splits into \(\gcd(d,r)\) distinct irreducible factors, each of degree \(d/\gcd(d,r)\).
If \(P\) runs over the distinct monic irreducible factors of \(T^n-1\), then
\[
\#\{\text{normal elements}\}=q^n\prod_P(1-q^{-\deg P}),
\qquad \#\{\text{primitive elements}\}=\varphi(q^n-1).
\]
Every finite extension has a normal basis. A self-dual basis exists if and only if \(q\) is even or \(n\) is odd.
A self-dual normal basis exists if and only if \(n\) is odd, or \(q\) is even and \(n\equiv2\pmod4\).
These results involve different generation and pairing conditions, whose connections are made explicit by cyclic modules and trace adjoints.

\chapter*{Index of Notation}
\addcontentsline{toc}{chapter}{Index of Notation}
\markboth{Index of Notation}{Index of Notation}
Notation is grouped by subject, with each meaning specified and appendix definitions placed at the end of the relevant group. Local uses of the same letter are distinguished by context. Page links refer to the stated meaning.
\begingroup\small
\renewcommand{\arraystretch}{1.04}
\begin{longtable}{@{}p{0.27\textwidth}p{0.55\textwidth}>{\raggedleft\arraybackslash}p{0.12\textwidth}@{}}
\toprule
Symbol&Meaning&Page\\
\midrule
\endfirsthead
\toprule
Symbol&Meaning&Page\\
\midrule
\endhead
\bottomrule
\endfoot
\addlinespace[0.7ex]\multicolumn{3}{@{}l}{\textbf{Unramified construction and arithmetic precision}}\\*
\(\mathbb Z_p,\ \mathbb Q_p\)&Compatible congruence limit and its fraction field&\pageref{v135:sec-foundations}\\
\(K_n,\ \mathcal O_n,\ k_n\)&Unramified extension, integer ring, residue field&\pageref{revadic:thm-main}\\
\(R_n,\ S_n\)&Root sets of $X^{p^n}-X$ in two characteristics&\pageref{revadic:thm-two-roots}\\
\(\tau,\ \boxplus\)&Teichmüller projection and transported addition&\pageref{revadic:prop-teich-projection}\\
\(W_g,\ P_N\)&First-carry mixed polynomial and pure-exponent set&\pageref{v144:eq-carry-obstruction}\\
\(\mathcal A_H,\ f_H\)&Coefficient fixed algebra and local component degree&\pageref{v144:eq-invariant-algebra}\\
\(S=a_{t-1}+1\)&Phase-labelling precision for a known valuation multiset&\pageref{v144:prop-label-precision}\\
\(M_B,\ \Delta_i,\ k_s\)&Integer multiplication matrix, determinantal gcds, and logarithmic kernel sizes&\pageref{v143:thm-valuation-smith}\\
\(A_{C,j},\ B_{C,j},\ \nu_{C,j}\)&Subset-count polynomial, cyclotomic remainder, coefficient valuation&\pageref{revadic:thm-integral-resolvent}\\
\(\epsilon_{C,j}(D)\)&Excess local cancellation depth&\pageref{revadic:cor-excess-degree}\\
\(\Pi_\Lambda\)&Remainder projection at jump positions&\pageref{v137:eq-jump-projection}\\
\(\operatorname{Per}_q(x)\)&Shortest period block, low-to-high order&\pageref{v138:subsec-periodic-digits}\\
\(\operatorname{per}_q(x)\)&Length of the shortest digit period&\pageref{v139:thm-period-compression}\\
\(\mathcal D_m,\Theta_m\)&Rational encoding group and its canonical-root isomorphism&\pageref{v139:eq-digit-root-group}\\
\(R,A\)&Rotation and addition with carries on periodic words&\pageref{v139:eq-two-actions}\\
\(\Delta(m,q),B(Z)\)&Largest coset leader; fixed-root consecutive-root bound&\pageref{v139:prop-digit-leaders}, \pageref{v139:prop-digit-bch-gap}\\
\(\mathcal T_s(V),U(V)\)&Repeated digit set and its allowable first translation digits&\pageref{v140:prop-digit-stabilizer}\\
\(V^{[2]}\)&Pointwise square image of a subspace&\pageref{v140:prop-digit-stabilizer}\\
\(D_{t,d},\rho_s(t,d)\)&Frobenius difference on a low-degree space and its precision profile&\pageref{v140:eq-precision-profile}\\
\(B_\sigma,\nu_S(\sigma)\)&Spectral-projector obstruction and preservation depth of a permutation&\pageref{v140:prop-permutation-depth}\\
\(H_{\rm aff}(C),H_j\)&Affine multiplier group and a single-coefficient stabilizer&\pageref{v141:prop-affine-profile}\\
\(\mathcal A,\mathcal B,\Gamma,\tau\)&Unramified phase algebra, residue algebra, phase group, and Frobenius element&\pageref{v141:eq-phase-algebra}\\
\(h(a),K_s,E_x\)&Directional depth, congruence subgroup, and index excess&\pageref{v141:thm-index-budget}\\
\(J_D,\lambda_D\)&Jacobi return value and normalized slope of a Fermat orbit&\pageref{v141:eq-fermat-orbit-factor}\\
\(P_t(z),\ M_d(z)\)&Jacobi valuation distribution for a fixed product character and its carry matrix&\pageref{v150:prop-jacobi-carries}\\
\(Z_X(T),\ P_X(T),\ N_r\)&Curve zeta function, numerator, extension-field counts&\pageref{revadic:thm-genus-counts}\\
\(u(X)\)&Multiplicity of slope zero&\pageref{revadic:cor-unit-roots}\\

\addlinespace[0.7ex]\multicolumn{3}{@{}l}{\textbf{Fields and groups}}\\*
\(p,\ q=p^r\)&Prime and prime power&\pageref{def:basic-notation}\\
\(\F_q\)&Finite field with \(q\) elements, unique up to isomorphism&\pageref{thm:unique}\\
\([E:k]\)&Extension degree, or \(k\)-vector-space dimension&\pageref{def:extension-degree}\\
\(k(S),\ k(\alpha)\)&Subfield generated by \(k\) and the specified set of elements&\pageref{def:generated-subfield}\\
\(k[\alpha]\)&Polynomial subring generated by \(k\) and \(\alpha\)&\pageref{def:polynomial-subring}\\
\(m_{\alpha,k}(X)\)&Minimal polynomial of a field element&\pageref{def:algebraic-minpoly}\\
\(E^\times\)&Multiplicative group of nonzero elements of \(E\)&\pageref{thm:cyclic}\\
\(\Aut_k(E)\)&Automorphism group fixing \(k\) pointwise&\pageref{thm:automorphisms}\\
\(\Fix(\sigma^s)\)&Fixed field of the \(s\)-th iterate&\pageref{prop:fixed-gcd}\\
\(\mathcal A(S)\)&Equivariant function algebra of a finite stable root set&\pageref{thm:equivariant-functions}\\
\(S_L(A)\)&Evaluation points of an algebra in a splitting field, with their Frobenius action&\pageref{thm:finite-orbit-duality}\\
\(E_a\otimes_kE_b\)&Tensor product of two field extensions&\pageref{thm:tensor-fields}\\
\(\GR(p^s,r)\)&Galois ring of characteristic \(p^s\) and residue field \(\F_{p^r}\)&\pageref{v47:def-galois-ring}\\
\(\mathcal T(A)\)&Teichm\"uller representatives of the Galois ring \(A\)&\pageref{v47:def-teichmuller}\\
\(W(K),W_s(K)\)&Witt vectors over the perfect field $K=\F_q$ and their length-$s$ truncation&\pageref{v133:subsec-witt-carry}\\
\([a],\ \kappa(a,b)\)&Teichmüller lift and first additive carry&\pageref{v133:prop-first-carry}\\
\(\mathcal O_1\)&First obstruction to lifting a zero-sum relation&\pageref{v133:prop-zero-sum-obstruction}\\
\(K[[t]],\ v_t\)&Equal-characteristic formal series ring and its $t$-adic valuation&\pageref{v133:subsec-equichar}\\
\(\mathcal V_p(C),\ m_j(D),\ \eta_D\)&Canonical valuation profile, subset-sum count, and orbit period&\pageref{v136:eq-valuation-profile}, \pageref{v134:thm-subset-periods}\\
\(\phi,\sigma\) (Galois rings)&Absolute and relative Frobenius automorphisms acting on Teichm\"uller digits&\pageref{v47:thm-galois-frobenius}, \pageref{v47:def-galois-relative-frobenius}\\
\(\mathfrak n,A_{\mathrm{per}}\)&Nilradical and Frobenius-periodic subalgebra of a finite algebra&\pageref{prop:algebra-periodic-reduction}\\
\(A_i,K_i,\ell_i\)&Local component, its periodic subfield, and dimension over the residue field&\pageref{v117:prop-local-dimension-trace}\\
\(R_Q,\ G_Q\)&Quadratic quotient algebra and its units modulo base-field scalars&\pageref{v85:def-projective-unit-group}\\
\(H(W)\)&Maximal field of scalars preserving \(W\)&\pageref{v70:def-scalar-field}\\
\(K_d\)&Unique intermediate field of degree \(d\) in \(E/k\)&\pageref{sfr:def-subfield-coordinate}\\
\(K=\F_p(D)\)&Hidden scalar field generated by the multiplicative subgroup \(D\)&\pageref{v74:rem-hidden-subfield}\\
\(\mathcal R_K(n)\)&Group preserving cyclotomic difference orbits&\pageref{v73:def-difference-orbit}\\
\(\operatorname{Aff}_K(n)\)&Affine permutations with multiplier in \(K\)&\pageref{v73:def-difference-orbit}\\
\(\Delta_n(d)\)&Kernel of unit reduction from modulus \(n\) to \(d\)&\pageref{v73:def-reduction-kernel}\\
\(L_i\)&Independent local multiplier group in CRT coordinates&\pageref{v73:eq-independent-local-group}\\
\(\widehat{\Z}\)&Group of compatible families of integer residue classes&\pageref{def:profinite-integers}\\
\(\mathcal D(K)\)&Extension degrees and exact Frobenius periods allowed in the intermediate field \(K\)&\pageref{prop:intermediate-periods}\\
\(\operatorname{Iso}_k(K,L)\)&Set of all \(k\)-isomorphisms between two finite extensions&\pageref{def:isomorphism-set}\\
\addlinespace[0.7ex]\multicolumn{3}{@{}l}{\textbf{Orbits and orders}}\\*
\(\phi,\sigma\) (finite fields)&Absolute and relative Frobenius: \(x^p,x^q\)&\pageref{def:frobenius-maps}\\
\(\Irr_q(d)\)&Set of monic irreducible polynomials of degree \(d\)&\pageref{cor:factor-product}\\
\(\mu(n)\)&M\"obius function on the positive integers&\pageref{def:mobius-function}\\
\(B_e(A),a_j\)&Number of components with degree divisible by \(e\), and number of degree-\(j\) field components&\pageref{fra:eq-second-inversion}\\
\(I_q(d)\)&Number of monic irreducible polynomials of degree \(d\)&\pageref{def:irreducible-count}\\
\(A_0(n),A_1(n)\)&Binary irreducible-polynomial counts by trace&\pageref{v109:ex-binary-trace}\\
\(J_{q,r}(d)\)&Number of length-\(d\) Frobenius orbits of \(r\)-tuples&\pageref{v107:prop-generating-tuples}\\
\(M_n(Z)\)&Necklace polynomial, with \(M_n(q)=I_q(n)\)&\pageref{v105:def-necklace-polynomial}\\
\(\ord_N(q)\)&Order of \(q\) in the unit group modulo \(N\)&\pageref{prop:degree-from-order}\\
\(\ord(a),\ord(f)\)&Multiplicative order of an element; order of a polynomial&\pageref{def:poly-order}\\
\(\varphi(N)\)&Euler totient function on the positive integers&\pageref{sec:euler-polynomial}\\
\(\varphi_q(F)\)&Number of units in \(\F_q[T]/(F)\)&\pageref{def:polynomial-euler}\\
\(B_f(X)\)&Minimal monic binomial multiple of a polynomial&\pageref{v40:prop-minimal-binomial}\\
\(\mu_N\)&Group of \(N\)-th roots of unity, with \(p\nmid N\)&\pageref{def:roots-of-unity}\\
\(\Phi_N(X)\)&\(N\)-th cyclotomic polynomial&\pageref{def:cyclotomic-polynomial}\\
\(C_i\)&\(q\)-cyclotomic coset modulo \(N\)&\pageref{def:cyclotomic-coset}\\
\(K^\times/H\)&Multiplicative quotient and its Frobenius coset orbits&\pageref{v85:prop-coset-frobenius}\\
\([a]_n\)&Least nonnegative residue&\pageref{v50:def-coset-leader}\\
\(\lambda_n(a),\Gamma_n(q)\)&Coset leader and the set of all coset leaders&\pageref{v50:def-coset-leader}\\
\(H_+(C)\)&Additive stabilizer of a cyclotomic coset&\pageref{def:additive-stabilizer}\\
\(\omega_C,\kappa_C\)&Least Frobenius step and largest complete equal-difference block size&\pageref{thm:med-coset-partitions}\\
\(b(d)\)&Small modulus retaining the required modulo-four information&\pageref{v124:def-small-modulus}\\
\(\rad(n)\)&Square-free kernel of a positive integer&\pageref{def:squarefree-kernel}\\
\(\omega_N\)&Least extension degree over which \(X^N-1\) factors into irreducible binomials&\pageref{cor:med-all-binomial}\\
\(\Omega_{n,\lambda}\)&Least extension degree for factorization of a general binomial into irreducible binomials&\pageref{v124:cor-binomial-extension}\\
\(c_j\) (root-of-unity traces)&\(u^j+u^{-j}\), a root-pair coefficient computed by a second-order recurrence&\pageref{v49:eq-coefficient-recurrence}\\
\(\kappa_T\)&Order of the translation stabilizer of a union of exponent orbits; equivalently, the maximum monomial composition exponent of the corresponding root polynomial&\pageref{v40:prop:root-power-composition}\\
\(s_j\)&\(j\)-th root power sum of a set of roots of unity&\pageref{v40:prop:root-power-composition}\\
\(\mathcal P\)&Synchronized phase quotient of local cycles&\pageref{v44:eq-phase-quotient}\\
\(\rho,\ K\)&Reduction of phase quotients and its fiber size&\pageref{v125:eq-phase-fiber}\\
\(B/\langle c\rangle\)&Kernel of phase reduction, with \(B=\prod_i\Z/(T_i/\tau_i)\Z\), \(c=(L/\tau_i)_i\)&\pageref{v126:eq-phase-exact}\\
\(u_\ell(r)\)&Base-\(\ell\) logarithm of the number of phase-kernel elements annihilated by \(\ell^r\)&\pageref{v126:eq-phase-kernel-torsion}\\
\(H_d\)&Reduction image of a multiplier group in the unit group modulo \(d\)&\pageref{v50:prop-multiplier-blocks}\\
\(e_{q,H}(d)\)&Exact coefficient-field degree of a multiplier-orbit polynomial&\pageref{v124:eq-coefficient-degree}\\
\(B_{H,i}\)&Monic root-polynomial block of a multiplier orbit&\pageref{v50:prop-multiplier-blocks}\\
\(C_m(\gamma),\tau_m(\gamma)\)&Cyclotomic orbit modulo \(m\) and its length&\pageref{v44:def-cyclotomic-system}\\
\(\mathcal C_m\)&Set of all cyclotomic orbits modulo \(m\)&\pageref{v44:def-cyclotomic-system}\\
\(\pi_{n,m}\)&Orbit reduction map between moduli related by divisibility&\pageref{v44:def-cyclotomic-system}\\
\(\mathsf{PC}_{q;\ell,m}\)&Set of all local cyclotomic systems in the \(\ell\)-direction with base level \(m\)&\pageref{v44:def-local-cyclotomic-system}\\
\(\mathsf{PC}_q\)&Set of all total cyclotomic systems&\pageref{v44:def-total-cyclotomic-system}\\
\(T(k)\)&Affine return permutation on a modulus-lifting fiber&\pageref{v44:eq-cyclotomic-return}\\
\((P_i)\)&Principal cyclotomic system above a fixed base-level orbit&\pageref{v46:def-principal-cyclotomic-system}\\
\(u,s\)&Deviation depth and stabilization index of a nonprincipal cyclotomic system&\pageref{v46:def-cyclotomic-deviation-degree}, \pageref{v46:def-cyclotomic-stable-degree}\\
\(F_C,B_C\)&Basic irreducible factors and canonical comaximal blocks of root-of-unity orbits&\pageref{v47:thm-galois-cyclotomic-factorization}, \pageref{v47:prop-galois-ramified-blocks}\\
\(r,\lambda\)&Projective order of a quadratic and the base-field constant at first return&\pageref{v45:prop-lucas-first-zero}\\
\(C_\Psi,\gamma_r,s_r\)&Root-pair composition polynomial, phase image, and its degree&\pageref{v129:prop-phase-images}\\
\(\left(\frac{g}{F}\right)_d\)&Polynomial \(d\)-th power-residue symbol&\pageref{def:power-residue-symbol}\\
\(F_{a,b}\)&Affine Frobenius permutation \(x\mapsto ax^{q^s}+b\)&\pageref{afr:def-affine-frobenius}\\
\(\theta\) (translation trace)&Relative trace of the scaled translation term&\pageref{afr:eq-trace-obstruction}\\
\(N_t(F)\)&Number of fixed points of the \(t\)-th iterate&\pageref{afr:thm-fixed-counts}\\
\(C_t(F)\)&Number of cycles of length \(t\)&\pageref{afr:def-cycle-count}\\
\(\beta_A(t)\)&Fixed-space dimension of the invertible linear operator \(A^t\)&\pageref{v53:eq-linear-fixed-dimension}\\
\(h_A(b)\)&Obstruction height of a translation relative to its linear part&\pageref{v107:def-translation-height}\\
\(\operatorname{sgn}(F)\)&Sign of a finite permutation: one for even and minus one for odd&\pageref{v53:def-permutation-sign}\\
\(\mathcal P_{q,n}\)&Exact-period layer: elements of least Frobenius period \(n\) over \(\F_q\)&\pageref{def:exact-period-layer}\\
\(E_{q,n}(X),D_{q,n}\)&Exact-period layer polynomial and its degree&\pageref{def:exact-period-layer}\\
\(\omega(n)\)&Number of distinct prime factors&\pageref{def:distinct-prime-count}\\
\(\Theta_q\)&Invariant for affine compression&\pageref{def:affine-compression-coordinate}\\
\(G_{q,n}\)&Compressed period-layer polynomial in \(\F_p[T]\)&\pageref{thm:period-compression}\\
\(R_{q,m}\)&Orbit-sum polynomial for affine return&\pageref{def:affine-return-sum}\\
\(e(c)\)&Affine return order determined by the compressed value \(c\)&\pageref{prop:explicit-affine-return}\\
\(u\cdot f,\operatorname{Stab}_G(f)\)&Affine action on monic polynomials and base-field stabilizer&\pageref{v129:cor-affine-factor-orbits}\\
\(\tau_{q,s}\)&Base-field lifting polynomial&\pageref{prop:period-base-change}\\
\(Z_S(u)\)&Formal power series encoding periodic-point counts&\pageref{prop:orbit-euler}\\
\addlinespace[0.7ex]\multicolumn{3}{@{}l}{\textbf{Traces and bases}}\\*
\(U_s\)&Norm-one subgroup for \(\F_{q^{2s}}/\F_{q^s}\)&\pageref{v48:def-norm-one-group}\\
\(C_n(q),\ S_q(2n)\)&Full-degree quadratic targets; monic self-reciprocal irreducibles&\pageref{v145:eq-target-count}\\
\(\Tr_{E/k},\Nm_{E/k}\)&Relative trace and norm&\pageref{def:trace-norm}\\
\(\phi_\theta(a)\)&Fractional-linear parameter for the norm-one group&\pageref{v85:prop-norm-one-parametrization}\\
\(t_q\)&Number of primitive quadratics with a fixed primitive constant term&\pageref{v45:prop-fixed-constant-primitive-count}\\
\(b_i^*\)&Dual basis vector corresponding to \(b_i\)&\pageref{thm:dual-basis}\\
\(U^\perp\)&Orthogonal complement for the trace pairing&\pageref{def:trace-orthogonal-complement}\\
\(\Tr_{A/k}\)&Matrix trace of a multiplication operator on a finite-dimensional algebra&\pageref{def:algebra-trace}\\
\(\lambda(u)\)&Coefficient of \(T^{\deg f-1}\) in the standard remainder in \(k[T]/(f)\)&\pageref{thm:algebra-trace-coefficient}\\
\(\operatorname{disc}(f)\)&Root-difference discriminant of a monic polynomial&\pageref{v105:def-polynomial-discriminant}\\
\(\rad(B|_W)\)&Radical of the trace pairing restricted to an operator image&\pageref{v68:cor-image-trace-radical}\\
\(h_\alpha,g_\alpha\)&Conjugate-coefficient and cyclic-correlation polynomials&\pageref{v88:eq-conjugate-correlation}\\
\(M_\alpha,C_\alpha\)&Moore matrix and cyclic multiplication matrix of conjugate coefficients&\pageref{v88:eq-moore-circulant}\\
\(C^{[p^j]}\)&Entrywise Frobenius on a matrix&\pageref{v97:prop-absolute-frobenius-coordinates}\\
\(F,G\) (prescribed intermediate traces)&Frobenius annihilator polynomial and least common multiple of the intermediate-trace moduli&\pageref{v58:thm-prescribed-normal-traces}\\
\(\pi_i,S(U)\)&Local Frobenius quotient and free component indices of an affine space with invariant direction&\pageref{v88:thm-invariant-affine-normal-count}\\
\(A_P,U_P,V_P\)&Local normality detector, obstruction kernel, and trace dual&\pageref{prop:normal-test}\\
\(U_i,V_i\)&Local bad subspaces for normality and their trace orthogonal complements&\pageref{v93:thm-patch}\\
\(\mathcal A_{\boldsymbol b}\)&Weighted trace affine fiber consisting of normal elements&\pageref{v93:thm-patch}\\
\(P_\theta,PN_\theta\)&Counts of primitive and primitive normal elements in a base-field translate&\pageref{v82:rem-primitive-normal-translates}\\
\(\beta_t=(\theta-t)^{-1}\)&Family of reciprocals of translates (candidates for normality)&\pageref{v68:prop-reciprocal-normal-family}\\
\(b_{q,n}(\alpha),\ B_q(n)\)&Exception count for a fixed reciprocal family and its maximum&\pageref{v114:def-bad-count}, \pageref{v114:def-max-bad}\\
\(\langle u,v\rangle_\alpha\)&Trace pairing for \(E/k\) with weight \(1/f'(\alpha)\)&\pageref{str:def-weighted-pairing}\\
\(U^{\perp_\alpha}\)&Orthogonal complement of \(U\) for the weighted trace pairing&\pageref{str:def-weighted-pairing}\\
\(c_j\) (relative trace)&\(\Tr_{E/K_d}(\alpha^j/f'(\alpha))\)&\pageref{str:def-relative-trace-sequence}\\
\(\mathcal C_{f,d}(s),r_{f,d}(s)\)&Base-field span of a relative-trace prefix and its dimension&\pageref{str:def-trace-rank}\\
\(\mu_i\) (trace rank)&Successive positions where the trace-rank function increases&\pageref{str:def-trace-jumps}\\
\(K_j,\ d_j,\ T_j\)&Field generated by power values, its degree, and the trace-evaluation parameter map&\pageref{v68:prop-trace-evaluation-parameters}\\
\(\mathcal C,\kappa,R\)&Normality-preserving scaling group, \(\gcd(n,q-1)\), and order of the multiplicative quotient&\pageref{v93:lem-scaling}\\
\addlinespace[0.7ex]\multicolumn{3}{@{}l}{\textbf{Operators and recurrences}}\\*
\(b_r(A)\)&Fixed-space dimension of the \(r\)-th Frobenius iterate&\pageref{fra:def-fixed-dimensions}\\
\(U_n(B,C),u_n\)&Lucas polynomials of the first kind and their coefficient specializations&\pageref{v45:def-lucas-polynomials}\\
\(L_V(X)\)&Monic root polynomial of the subspace \(V\)&\pageref{thm:subspace-polynomial}\\
\(I(U),r(U)\)&Reciprocal span of a binary subspace and its dimension&\pageref{v116:def-inverse-rank}\\
\(\mathcal C(U)\)&Reciprocal relation space of a four-dimensional binary subspace&\pageref{v116:def-inverse-relations}\\
\(\mathscr L_{<d}(E)\)&Space of linearized polynomials over \(E\) of \(q\)-degree less than \(d\)&\pageref{v105:prop-linearized-remainder-compatibility}\\
\(\operatorname{rem}_i L\)&Right composition remainder modulo the \(i\)-th subspace polynomial&\pageref{v105:prop-linearized-remainder-compatibility}\\
\(\mathcal L_a(X)\)&Linearization \(\sum a_iX^{q^i}\) of \(a(T)=\sum a_iT^i\)&\pageref{def:linearized-symbol}\\
\(K[t;\sigma],\ \Theta\)&Skew polynomial ring; Frobenius return fixing coefficients&\pageref{v145:thm-skew-roots}\\
\(\mathcal N(f),\ U_k(h),\ F_{h,k}(c)\)&Return norm, total fiber, and fixed-constant fiber&\pageref{v148:sec-norm-fibers}\\
\(T_k(\lambda;Q)\)&Fiber count of a local partition&\pageref{v148:prop-local-types}\\
\(C_{n,k}(q)\)&Common irreducible-target fiber for coprime degrees&\pageref{v148:eq-Cnk}\\
\(E_h(c),\ E(h),\ B_m(q)\)&Commutative norm-lift counts and even-degree baseline&\pageref{v148:thm-all-even-fibers}\\
\(\rho_\pi,\ r_\pi\)&Ratio of conjugate factor constants and its exponent in the norm-one group&\pageref{v149:prop-uniformity}\\
\(H_{k,c}(A)\)&Optimal parameter polynomial for norm-primitivity&\pageref{v148:prop-optimal-strata}\\
\(w_q(T,N),\ P_D(t,n)\)&Quadratic trace--norm weight and primitive cubic slice&\pageref{v148:thm-cubic-primitive-slices}\\
\(L^*\)&Trace adjoint of a linear map&\pageref{prop:linearized-adjoint}\\
\(UV\)&Base-field linear span of all products \(uv\)&\pageref{v70:def-subspace-product}\\
\(T_{a,W},\ I_W(a)\)&Quotient multiplication map and base-field dimension of \(W\cap aW\)&\pageref{v70:prop-quotient-multiplication}\\
\(P_r(\alpha)\)&Space spanned by the first $r$ powers&\pageref{v150:prop-power-intersection}\\
\(D_L\)&Dickson matrix of a linearized polynomial&\pageref{def:dickson-matrix}\\
\(\mu_a(T)\)&Frobenius annihilator polynomial of \(a\)&\pageref{sec:q-order}\\
\(S_a\)&Weighted Frobenius operator \(x\mapsto ax^{q^s}\)&\pageref{wfr:def-weighted}\\
\(c,P_c\) (weighted operator)&Relative norm of the weight and its base-field minimal polynomial&\pageref{wfr:def-weighted}\\
\(U_h\)&Invariant root space \(\ker h(\sigma)\)&\pageref{thm:invariant-duality}\\
\(u^*\) (quotient ring)&Involution \(u(T)\mapsto u(T^{-1})\)&\pageref{def:quotient-involution}\\
\(m_h\)&Value polynomial of a fixed element&\pageref{prop:berlekamp-value-polynomial}\\
\(e_c=1-(h-c)^{q-1}\)&Idempotent projection onto the value fiber \(h=c\) of a fixed element&\pageref{v82:eq-value-idempotents}\\
\(\mathrm{NP}\)&Class of decision problems admitting a deterministic polynomial-time verifier, with polynomial-length accepted certificates for yes-instances and no accepted certificate for no-instances&\pageref{v58:def-np-bpp}\\
\(\mathrm{BPP}\)&Class of decision problems admitting randomized polynomial-time algorithms with error probability at most \(1/3\) on every input&\pageref{v58:def-np-bpp}\\
\(\mathcal H_{f,d}\)&Hankel matrix of order \(m\) for the relative-trace sequence&\pageref{str:eq-hankel}\\
\(I_\ell\)&Common kernel of the linear observation \(\ell\) over all shifts&\pageref{v127:eq-observation-kernel}\\
\(s_j(f)\)&\(j\)-th power sum of the roots of a polynomial&\pageref{v123:sec-power-trace-fibers}\\
\(h_f,\ f_c\)&Derivative inverse representative and \(f+c h_f\)&\pageref{v123:thm-power-trace-pencil}\\
\(\mathcal P(f)\)&Derivative partner of a binary square-free polynomial&\pageref{v123:cor-binary-derivative-pair}\\
\(\delta_a\)&Polynomial indicator function of the point \(a\)&\pageref{def:point-indicator}\\
\(V\) (\(h=|H|\))&Function space spanned by the translation orbit of a monomial&\pageref{v74:eq-translation-power-space}\\
\(\mathcal A_f(x,y,z)\)&Three-point difference determinant&\pageref{v76:def-three-point-determinant}\\
\(V(x,y,z)\)&Vandermonde normalizing factor for three-point divided differences&\pageref{v76:def-second-divided-difference}\\
\(\mathcal K_f(x,y,z)\)&Second divided difference, symmetric in its three nodes&\pageref{v76:def-second-divided-difference}\\
\(K_{\lambda,f}(x,t)\)&Second divided difference of shape \(\lambda\)&\pageref{v76:def-fixed-shape-difference}\\
\(N_\lambda(H)\)&Number of nondegenerate four-term relations of a prescribed shape in the subgroup&\pageref{v76:def-shape-four-term-count}\\
\(C_j(\lambda)\)&Frobenius shape factor&\pageref{v76:def-frobenius-shape-factor}\\
\(S\)&Shift operator on sequences&\pageref{def:sequence-shift}\\
\(h_s(T)\)&Minimal polynomial of a sequence&\pageref{def:sequence-minpoly}\\
\(\mathcal R_f,W_f\)&Recurrence solution space and root space of \(\mathcal L_f\)&\pageref{prop:recurrence-root-correspondence}\\
\((L,C)\)&Connection length and connection polynomial of a finite prefix&\pageref{def:finite-linear-complexity}\\
\(\mathcal D_q\)&Frobenius derivation on a polynomial ring&\pageref{v84:def-frobenius-derivation}\\
\(R_d\) (Frobenius derivation)&Homogeneous polynomial space of degree \(d\)&\pageref{v86:thm-frobenius-kernel}\\
\(I(\mathbb P^{r-1}(k))_{q+1}\)&Critical-degree component of the homogeneous vanishing ideal&\pageref{v86:eq-vanishing-kernel}\\
\addlinespace[0.7ex]\multicolumn{3}{@{}l}{\textbf{Root sets and fibers}}\\*
\(h^\vee\)&Monic reciprocal polynomial \(T^{\deg h}h(T^{-1})/h(0)\)&\pageref{def:reciprocal-polynomial}\\
\([X^j]P\)&Coefficient of \(X^j\) in \(P\)&\pageref{def:coefficient-extraction}\\
\(k[X^r]\)&Subring of polynomials of the form \(G(X^r)\)&\pageref{def:coefficient-extraction}\\
\(\Lambda_n(B,C)\)&Lucas atom defined by a homogeneous cyclotomic polynomial&\pageref{v45:def-lucas-atom}\\
\(V_n(B,C)\)&Second-kind Lucas polynomial, recording quadratic root power sums&\pageref{v146:prop-lucas-factorization}\\
\(H_{q,s}(Y),H_q(Y)\)&Exact-order trace polynomial of the norm-one group; \(H_q=H_{q,q+1}\)&\pageref{v146:eq-order-trace}\\
\(B_{k,r}(a,c),\Theta_{k,r}(a)\)&Negative return trace and its constant-term-one specialization&\pageref{v146:eq-return-spectrum}\\
\(\mathcal I_P,\mathcal C(\Phi)\)&Irreducible quadratic targets and normalized cyclic vectors&\pageref{v146:eq-theta-definition}, \pageref{v146:lem-semilinear-normal}\\
\(\mathcal P_c\)&Coefficient set of primitive quadratics with a fixed primitive constant term&\pageref{v45:def-primitive-coefficient-set}\\
\(P_{q,c}(B)\)&Optimal coefficient-determining polynomial for a fixed primitive constant term&\pageref{v45:def-optimal-determining-polynomial}\\
\(A=k[C]/(\Phi_{q-1}(C))\)&Coefficient algebra organizing all primitive constant terms&\pageref{v45:prop-global-optimal-uniqueness}\\
\(D,\Delta,\Phi(w)\)&Difference of additive cosets, normalized difference, and quadratic root-set equation&\pageref{v45:eq-quadratic-normalization}\\
\(V_+,V_-\)&Frobenius eigenspaces for eigenvalues +1 and -1 associated with a quadratic root set&\pageref{v45:prop-quadratic-planes}\\
\(Z_\Delta\)&Nonsquare parameter set recording irreducible conjugate root pairs&\pageref{v45:prop-quadratic-norm-descent}\\
\(p_S(T)\)&Monic root polynomial of a finite set \(S\)&\pageref{v85:cor-rootset-quotient-coordinates}\\
\(\operatorname{wt}(f)\)&Number of nonzero coefficients of \(f\)&\pageref{v79:def-defining-complexity}\\
\(w_q(n),\ \theta_q(n)\)&Minimum support size and minimum highest degree of a nonleading term among monic irreducible polynomials of degree \(n\ge2\)&\pageref{v79:def-defining-complexity}\\
\(L_s(f)\)&Minimum degree of an at-most-\(s\)-term multiple of fixed \(f\)&\pageref{v97:def-sparse-degree}\\
\(\mathcal A_q(n)\)&Constant terms \(a\) for which \(X^n+a\) is irreducible&\pageref{v77:def-binomial-constants}\\
\(\mathcal Q,\mathcal R\)&Q-transform and R-transform&\pageref{thm:reciprocal-transform}, \pageref{def:r-transform}\\
\(\mathcal D_k(f)\)&Root-field spectrum of \(f\) over \(k\)&\pageref{def:root-field-spectrum}\\
\(\nu_f,S_f\)&Root-multiplicity function and set of distinct roots&\pageref{ffs:def-root-multiplicity}\\
\(V_f,M_f,\mathcal A_f\)&Additive, multiplicative, and full affine stabilizers&\pageref{ffs:def-stabilizers}\\
\(\mathcal A_f(k)\)&Affine stabilizer with base-field parameters&\pageref{ffs:def-stabilizers}\\
\(e\)&Largest common \(p\)-power divisor of the root multiplicities&\pageref{ffs:prop-common-power}\\
\(\partial^{[j]}f\)&\(j\)-th Hasse derivative&\pageref{ffs:def-hasse}\\
\(\Psi_f\)&Monic greatest common divisor of translation-difference coefficients&\pageref{ffs:def-translation-polynomial}\\
\(\mathcal A_P\)&Largest base-field subspace giving constant translation differences&\pageref{v109:prop-additive-directions}\\
\(\operatorname{sp}(f),f_{\mathrm{sp}}\)&Support period and support-primitive core&\pageref{v44:def-support-period}\\
\(I_f\)&First-order affine symmetry space&\pageref{v67:def-first-order-affine}\\
\(Q,t\)&Order of the translation group and of the slope image of the full affine stabilizer&\pageref{ffs:def-affine-data}\\
\(C,\delta,U\)&Center coset, its quotient value, and its monic root polynomial&\pageref{ffs:thm-center-descent}\\
\(\mathcal W_{\Gamma,r}(N)\)&Space of semi-invariants with fixed scaling factor \(a^r\) and degree at most \(N\)&\pageref{fas:def-semi-invariant}\\
\(\operatorname{Supp}_X(P)\)&Exponents of the nonzero coefficients of \(P\)&\pageref{fas:def-common-support}\\
\(E_\Gamma,E_0\)&Common support sets for full affine semi-invariants and for invariants under translations and scalings about the origin&\pageref{v67:cor-affine-support-bounds}\\
\(R_\beta\)&Affine return permutation on an additive root fiber&\pageref{ffs:thm-additive-return}\\
\(\tau,\mathcal R\)&Frobenius conjugacy and fiber return on group parameters&\pageref{ffs:thm-full-return}\\
\(b_p\)&Conjugate difference \(\sigma^{n/p}(\alpha)-\alpha\) of an irreducible root&\pageref{fas:eq-prime-difference}\\
\(\mathcal N_E(V),h_V\)&Scalar group preserving the translation space \(V\), and its order&\pageref{fas:def-scalar-normalizer}\\
\(\mathcal A_{G,E}^{[V]}\)&Affine stabilizer of the quotient polynomial with slopes preserving \(V\)&\pageref{fas:eq-compatible-quotient}\\
\(s_q(N)\)&Size of the multiplicative stabilizer of an irreducible polynomial with root order \(N\)&\pageref{fas:eq-order-multiplier}\\
\(N_z,t_z\)&Multiplicative order of a centered quotient root and number of compatible slopes&\pageref{fas:eq-center-test}\\
\(\deg_\alpha(u)\)&Degree of the representative of \(u\) in the power basis of the generator \(\alpha\)&\pageref{sfr:def-representative-degree}\\
\(\mathcal P_t(\alpha)\)&Space of power-basis representatives of degree at most \(t\)&\pageref{sfr:def-degree-filtration}\\
\(U_{f,d}(t)\)&Intersection of \(K_d\) with the low-degree space \(\mathcal P_t(\alpha)\)&\pageref{sfr:def-profile}\\
\(\rho_{f,d}(t)\)&Dimension of the low-degree intersection space: the subfield degree profile&\pageref{sfr:def-profile}\\
\(\kappa_f(d)\)&Minimum generator-representation degree of a subfield&\pageref{sfr:def-generator-degree}\\
\(\Lambda_{f,d},\lambda_i\)&Set of degree jumps and its ordered jump positions&\pageref{sfr:def-jumps-defect}\\
\(\Delta_f(d)\)&Representation defect of the degree-jump sequence relative to the least degree&\pageref{sfr:def-representation-defect}\\
\(T_{f,d,t},M_{f,d,t}\)&Linear map and matrix for the subfield fixed-point equation&\pageref{sfr:def-fixed-point-matrix}\\
\(N_f(d,t)\)&Number of subfield-generator representatives of degree at most \(t\)&\pageref{sfr:prop-exact-count}\\
\(B_{f,d},b_i\)&Monic minimal polynomial of \(\alpha\) over \(K_d\) and its coefficients&\pageref{str:def-block-polynomial}\\
\(\mathcal V_{f,d},\eta_f(d)\)&Base-field span of root-block coefficients and its rank&\pageref{str:def-coefficient-rank}\\
\(\mathcal W_{f,d}\)&Coordinate-polynomial space of the root block&\pageref{str:def-coordinate-space}\\
\(\kappa_f^{\mathrm{rat}}(d)\)&Minimum rational generator-representation degree of a subfield&\pageref{str:def-rational-generator-degree}\\
\(\widehat F\) (reciprocal polynomial)&\(X^{\deg F}F(X^{-1})\), without monic normalization&\pageref{sec:recursive-as}\\
\(h^{\circ n}\)&\(n\)-fold compositional iterate of \(h\)&\pageref{v55:def-dynamical-irreducibility}\\
\((h,\tau)\)&Iteration pair with specified target \(\tau\)&\pageref{v55:def-dynamical-irreducibility}\\
\(C,\tau,\gamma\) (quadratic iteration)&Normalized constant term, transformed target, and original critical point&\pageref{v56:cor-general-quadratic-stability}\\
\(\Psi,A_n,B_n\) (rational iteration)&\(\Psi=1/\psi\), numerator and monic denominator of the reduced \(n\)-th iterate, \mbox{\(n\ge1\)}&\pageref{v55:def-inverse-stability}\\
\(\Gamma_H,\mathcal E_H(f),D_H(f)\)&Exact difference-quotient models, bad ordered pairs, and minimum disagreements&\pageref{v116:def-stability}\\
\(\mathcal V_F,\ N_b(F)\)&Value set of a polynomial and the size of each value fiber&\pageref{v58:def-value-set}\\
\(W^v\)&\(v\)-th power image of a set&\pageref{v58:def-power-image}\\
\(C_F\text{ (collision count)}\)&Number of ordered pairs with the same image, including diagonal pairs&\pageref{v58:def-collision-number}\\
\(P_N(X)\)&\(X(X-1)^{N-1}\), with \(N=p^e\), whose value distribution is governed by an intersection field&\pageref{v68:prop-affine-values}\\
\(\Phi_{A,B},\ N_H(A,B,C)\)&Two-point coefficient map and ordered distinct-coordinate count&\pageref{v85:sec-pair-coefficient-fibers}\\
\(\mathbb P^1(k)\)&Projective line over a finite field&\pageref{v85:def-projective-line}\\
\(D_m(X)\)&Dickson polynomial with parameter one&\pageref{def:dickson-polynomial}\\
\(P_m(Y)\)&Reciprocal-orbit root polynomial for roots of unity of odd order&\pageref{v97:def-reciprocal-orbit-polynomial}\\
\(M_j(F)\)&Number of ordered \(j\)-tuples of pairwise distinct preimages with the same image&\pageref{v58:ex-factorial-value-count}\\
\(P_V(Y),F_V(X)\)&Parameter polynomial and multiplicative-coset root polynomial for a fixed derivative&\pageref{v83:eq-derivative-fiber}\\
\(c_q(A),H_A,F_Q\)&Minimum completion size, inverse-Frobenius interpolant, and all candidates&\pageref{v116:def-completion}\\
\(B_q(U,V)\)&Projective root form \(U^qV-UV^q\)&\pageref{v84:eq-projective-root}\\
\(R_S(U,V)\)&Normalized homogeneous form of a projective root set&\pageref{v84:prop-projective-descent}\\
\addlinespace[0.7ex]\multicolumn{3}{@{}l}{\textbf{Codes}}\\*
\(\operatorname{wt}_H,d_H\)&Hamming weight and minimum distance of a nonzero linear code&\pageref{v44:def-hamming-weight}\\
\(\operatorname{wt}_R(c)\)&Dimension of the coordinate span over the base field&\pageref{v99:def-rank-weight}\\
\(d_R(\mathcal C)\)&Minimum rank distance of a nonzero linear code&\pageref{v99:def-rank-distance}\\
\(A_N\)&Quotient algebra \(k[X]/(X^N-1)\) for cyclic shift&\pageref{v40:def-cyclic-code}\\
\(R_N\)&Cyclic right-shift operator of length \(N\)&\pageref{v40:def-cyclic-code}\\
\(T(C)\)&Defining set of a cyclic code&\pageref{v40:def-defining-set}\\
\(C_\zeta\)&Trace cyclic code of a single root-of-unity orbit&\pageref{v40:prop-cyclic-trace-model}\\
\(\widehat c(i)\)&Finite-field-valued Fourier coordinates and their Frobenius descent condition&\pageref{v96:prop-fourier-descent}\\
\(R_{N,\lambda}\)&Quotient algebra \(k[X]/(X^N-\lambda)\) for constacyclic codes&\pageref{v54:def-constacyclic-code}\\
\(\operatorname{rev}_N\)&Coefficient reversal \(F\mapsto X^{N-1}F(X^{-1})\) for length \(N\)&\pageref{v87:thm-constacyclic-mds-reciprocal}\\
\(b,\delta\) (BCH codes)&Starting exponent of consecutive zeros and designed distance&\pageref{v56:def-bch-code}\\
\(S_j,\Lambda(Z)\)&Syndromes and the error-locator polynomial&\pageref{v118:def-syndrome}, \pageref{v118:def-error-locator}\\
\(\Gamma(L,G)\)&Binary Goppa code with location set \(L\) and Goppa polynomial \(G\)&\pageref{v82:def-binary-goppa-code}\\
\(F_S(X)\)&Root polynomial of the nonzero coordinate positions&\pageref{v82:def-support-root-polynomial}\\
\(\operatorname{Supp}_X(D)\)&Effective support of a subcode in the original coefficient coordinates&\pageref{v54:def-effective-support}\\
\(W_{r,J},B_r,\rho_{r,J}\)&Power space of semi-invariants, outer divisibility factor, and rank of the remainder map&\pageref{v54:thm-outer-divisibility}\\
\(\mathcal H_J\)&Selected-column Hankel matrix for outer divisibility&\pageref{v67:eq-hankel-divisibility}\\
\(\mathcal D_J(H)\)&Hasse derivative matrix for outer divisibility&\pageref{v67:eq-hasse-divisibility}\\
\(\operatorname{ev}_d\)&Homogeneous evaluation in a fixed degree&\pageref{v84:prop-projective-evaluation}\\
\(\operatorname{RS}_K(S),\ \operatorname{ERS}_K(S)\)&Ordinary and extended Reed--Solomon evaluation codes&\pageref{v84:prop-rs}\\
\(\mathcal E_H\)&Evaluation space defined by differential interpolation&\pageref{v84:prop-differential-interpolation}\\
\(C_\alpha\)&Irreducible cyclic trace code in inverse-root coordinates&\pageref{v84:eq-trace-code-convention}\\
\(\mathcal R_\alpha\)&Inverse-root linear-relation map, with kernel the dual trace code&\pageref{v96:eq-relation-map}\\
\(\mathcal D_H\)&Binary linear-relation kernel for a norm-one group&\pageref{v153:cor-normone-binary-code}\\
\(J_d,D_d\)&Subgroup coordinate set and the trace subcode vanishing there&\pageref{v142:eq-subgroup-restriction}\\
\(S(C)\)&Fourier frequencies allowed to be nonzero in a cyclic code&\pageref{v84:def-fourier-support}\\
\(\F_{p^s}\)&Smallest field of definition in fixed cyclic coordinates&\pageref{v84:thm-code-field}\\
\(\operatorname{PAut}(C),\mathcal G_q(n)\)&Permutation group of one code and the group common to all cyclic codes&\pageref{v73:def-common-permutation-group}\\
\(h_O,\ P_O\)&Cyclotomic orbit sum and Fourier spectral projection&\pageref{v73:eq-fourier-orbit-projector}\\
\(\iota_q(d)\)&Multiplier index \(\varphi(d)/\ord_d(q)\)&\pageref{v73:def-multiplier-index}\\
\(\mathcal P(H)\)&Projective image of a multiplicative subgroup over the base field&\pageref{v96:def-projective-arc}\\
\addlinespace[0.7ex]\multicolumn{3}{@{}l}{\textbf{Character sums}}\\*
\(\operatorname{PN}_q(n)\)&Number of primitive normal elements of \(\F_{q^n}/\F_q\)&\pageref{v93:cor-exact-pn}\\
\(\psi_a,\chi_j\)&Additive and multiplicative characters&\pageref{thm:characters}\\
\(\widehat F\) (Fourier transform)&Finite Fourier transform&\pageref{def:fourier-transform}\\
\(S_d(a)\)&Family of polynomial-phase character sums whose even moments are determined by counts for power-sum equations&\pageref{v58:def-character-moment}\\
\(A_*F\)&Summation over fibers of the linear map \(A\)&\pageref{prop:fourier-adjoint}\\
\(G(\chi,\psi)\)&Gauss sum&\pageref{def:gauss-sum}\\
\(J(\chi,\lambda)\)&Jacobi sum&\pageref{def:jacobi-sum}\\
\(A_r(c),\ B_r(u,v)\)&Fiber counts for one or two prescribed subset products&\pageref{v150:prop-subset-products}\\
\(\mathcal F_r(\chi,\lambda)\)&Multiplicative Fourier coefficient of joint-product fibers&\pageref{v150:prop-joint-products}\\
\(\lambda_q(H)\)&Linear relation length of a multiplicative subgroup&\pageref{v95:eq-lambda-definition}\\
\(N_H(a)\)&Intersection size of a trace hyperplane and a subgroup&\pageref{v95:subsec-trace-hyperplanes}\\
\(r_H(\delta),\ \mathcal E_+(H)\)&Translation intersection count and additive energy&\pageref{v151:prop-normone-intersections}, \pageref{v151:def-additive-energy}\\
\(U=\F_q^\times H\)&Scalar saturation of the multiplicative subgroup&\pageref{v95:subsec-cyclotomic-triples}\\
\(\mathcal C(U)\)&Scalar-saturated cyclotomic intersection number&\pageref{v95:subsec-cyclotomic-triples}\\
\(T_3(H)\)&Number of dependent three-element subsets&\pageref{v95:thm-triple-cyclotomic}\\
\(\eta_i\)&Gaussian periods of the scalar-saturated subgroup&\pageref{v95:eq-gaussian-period}\\
\(B_r\)&Number of hyperplanes containing a prescribed subspace&\pageref{v95:subsec-factorial-moments}\\
\(\mathfrak M_s(H)\)&Factorial moment of trace-hyperplane intersections&\pageref{v95:thm-factorial-rank}\\
\(R_s(H)\)&Number of s-subsets with rank s-1&\pageref{v95:eq-rank-step}\\
\(D_{\Gamma,S}(\zeta_n)\)&Cyclotomic Fourier minor&\pageref{v95:subsec-examples-fourier}\\
\(J(\chi_1,\ldots,\chi_n)\)&Multiple Jacobi sum&\pageref{def:multiple-jacobi-sum}\\
\(S_a(\chi)\)&Multiplicative character sum over a trace fiber&\pageref{def:trace-fiber-character-sum}\\
\(N_a\)&Number of primitive elements of trace \(a\)&\pageref{def:primitive-trace-count}\\
\(P_{\ne0}\)&Number of primitive elements with nonzero trace&\pageref{v88:eq-primitive-nonzero-trace-lower}\\
\(\mathscr N(e,g)\)&Joint count for multiplicative freeness and local normality conditions&\pageref{v93:prop-basic}\\
\(\Theta(g),W_q(g)\)&Density of polynomial local conditions and squarefree divisor count&\pageref{v128:prop-local-spectrum}\\
\(V_g,c_g(b)\)&Frequency space and normalized expansion coefficients of local generation conditions&\pageref{v128:prop-local-spectrum}\\
\(J_D\)&Translation stabilizer of an annihilator stratum and orthogonal direction of its exact spectrum&\pageref{v130:prop-annihilator-spectrum}\\
\(N_q^\square\)&Number of monic irreducible quadratics with nonzero square nonleading coefficients&\pageref{v105:prop-square-coefficient-quadratics}\\
\(B_Q,\rad B_Q\)&Polar form of a quadratic form and its radical&\pageref{def:quadratic-form}\\
\(N_Q(b)\)&Number of preimages of \(b\) under a quadratic form&\pageref{def:quadratic-value-count}\\
\(\epsilon(Q)\)&Type sign of a quadratic form in characteristic two&\pageref{def:quadratic-type-sign}\\
\(F_Q,\mathcal S_Q\)&Quadratic phase function and its Fourier support&\pageref{prop:quadratic-fourier-support}\\
\(\chi_n,\psi_n\)&Characters lifted to an extension by norm and trace&\pageref{thm:hasse-davenport}\\
\(M_k(a),\ b_k(a)\)&Frobenius product of quadratic companion matrices and its negative trace&\pageref{v145:eq-quadratic-return}\\
\(D_k(P)\)&Fiber size of an irreducible quadratic return target&\pageref{v145:eq-return-fiber-size}\\
\(\mathcal N_n(q,nk)\)&Irreducible parameters of the sparse skew polynomial $t^{2n}+at^n+1$&\pageref{v145:thm-sparse-count}\\
\(S_r(f)\)&Additive character sum with polynomial phase over an extension&\pageref{def:character-L-series}\\
\(L_f(u)\)&Generating series of additive character sums&\pageref{def:character-L-series}\\
\(K_r(a,b)\)&Kloosterman sum over an extension&\pageref{def:kloosterman-sum}\\
\(\mathbf P^2(K)\)&Projective plane over \(K\)&\pageref{def:projective-plane-points}\\
\(E_b(K),O\)&Set of \(K\)-rational points of \(E_b\) and its point at infinity&\pageref{def:rational-curve-point}\\
\(N_n\) (curve point count)&\(\#E_b(\F_{q^n})\), including the point at infinity&\pageref{def:curve-point-count}\\
\(a,A_n\)&\(a=q+1-N_1\), \(A_n=q^n+1-N_n\)&\pageref{cor:binary-curve-recurrence}\\
\(Z(E_b/k,u)\)&Exponential generating expression for curve point counts&\pageref{def:curve-zeta}\\
\(P_b(T)\)&\(T^2-aT+q\), governing the recurrence for curve point counts&\pageref{def:curve-frobenius-polynomial}\\
\(S_H(s),\ H^\perp\)&Distinct-element sum count and annihilating character group&\pageref{v85:ex-restricted-sums}\\
\end{longtable}
\endgroup

\chapter*{Index}
\addcontentsline{toc}{chapter}{Index}
\markboth{Index}{Index}
Main entries are alphabetical, with indented subentries. \emph{See} points to a preferred entry in this index; \emph{see also} links related, distinct concepts. Page references lead to the principal definition or the relevant result; Roman-numbered pages refer to the thematic reading guide.
\begingroup
\linespread{1}\fontsize{11}{16}\selectfont\setlength{\parskip}{0pt}

\ifdefined\IdxPageWidth\else\newlength{\IdxPageWidth}\fi
\newcommand{\IdxLine}[3]{%
  \par\begingroup\setlength{\IdxPageWidth}{0pt}%
  \if\relax\detokenize{#3}\relax\else\settowidth{\IdxPageWidth}{#3}\addtolength{\IdxPageWidth}{0.6em}\fi
  \parindent=0pt\leftskip=#1\relax\hangindent=0.85em\hangafter=1
  \interlinepenalty=10000\relax
  \if\relax\detokenize{#3}\relax\rightskip=0pt plus 1fil\parfillskip=0pt plus 1fil\relax\else\rightskip=0pt plus 0.8em\relax\advance\rightskip by\IdxPageWidth\relax\parfillskip=-\IdxPageWidth\relax\fi
  \noindent#2\if\relax\detokenize{#3}\relax\else\nobreak\dotfill\nobreak\hbox to\IdxPageWidth{\hfil#3}\fi\par
  \endgroup}
\newcommand{\IdxHead}[3]{\par\addvspace{0.40\baselineskip}%
  \IdxLine{0pt}{\hypertarget{#1}{\textbf{#2}}}{#3}\nobreak}
\newcommand{\IdxSub}[3]{\IdxLine{1.2em}{\hypertarget{#1}{#2}}{#3}}
\newcommand{\IdxSee}[3]{\IdxLine{0pt}{#1\enspace\emph{#2}~#3}{}}
\newcommand{\IdxGroup}[1]{\par\addvspace{0.85\baselineskip}%
  \noindent\textbf{\color{Jade}#1}\par\nobreak}

\setlength{\columnsep}{1.4em}
\begin{multicols}{2}
\raggedcolumns
\IdxGroup{A}
\IdxSee{Absolute trace}{see}{\hyperlink{idx101.en.e4}{Trace: Absolute trace}}
\IdxSee{Additive character}{see}{\hyperlink{idx101.en.e7}{Character: Additive character}}
\IdxHead{idx151.en.energy}{Additive energy}{\pageref{v151:def-additive-energy}}
\IdxSub{idx151.en.normmoment}{Norm-one group and a fourth moment}{\pageref{v151:cor-normone-moment}}
\IdxSee{Additive Fourier transform}{see}{\hyperlink{idx101.en.e8}{Fourier transform: Additive transform}}
\IdxSee{Additive stabilizer of a cyclotomic coset}{see}{\hyperlink{idx101.en.e9}{Cyclotomic coset: Additive stabilizer}}
\IdxHead{idx101.en.adjoint}{Adjoint map}{\pageref{def:adjoint-map}}
\IdxSub{idx101.en.e368}{Fourier transform}{\pageref{prop:fourier-adjoint}}
\IdxSub{idx101.en.e17}{Linearized operators}{\pageref{prop:linearized-adjoint}}
\IdxSub{idx101.en.e21}{Relative trace and inclusion}{\pageref{prop:relative-trace-adjoint}}
\IdxLine{1.2em}{\emph{See also} \hyperlink{idx101.en.coordinates}{Field coordinates}; \hyperlink{idx101.en.pairing}{Trace pairing}}{}
\IdxSee{Adjoint operator}{see}{\hyperlink{idx101.en.adjoint}{Adjoint map}}
\IdxHead{idx101.en.affinefrob}{Affine Frobenius map}{\pageref{afr:def-affine-frobenius}}
\IdxSub{idx101.en.e646}{Characteristic part of periods}{\pageref{afr:cor-period-jump}}
\IdxSub{idx101.en.e183}{Cycle counts}{\pageref{afr:eq-cycle-count}}
\IdxSub{idx101.en.e346}{Fixed-point counts}{\pageref{afr:thm-fixed-counts}}
\IdxSub{idx101.en.e484}{Least dividing period}{\pageref{lem:affine-least-period}}
\IdxSub{idx101.en.e140}{Least-period computation}{\pageref{prop:affine-least-period-computation}}
\IdxSub{idx101.en.e849}{Shortest cycles}{\pageref{afr:cor-shortest-cycles}}
\IdxSub{idx101.en.e780}{Trace obstruction to translation removal}{\pageref{afr:eq-trace-obstruction}}
\IdxHead{idx141.en.affine}{Affine multiplier group}{\pageref{v141:prop-affine-profile}}
\IdxSee{Affine \(q\)-polynomial}{see}{\hyperlink{idx101.en.e28}{Linearized polynomial: Affine \(q\)-polynomial}}
\IdxHead{idx101.en.affine}{Affine stabilizer}{\pageref{ffs:def-stabilizers}}
\IdxSub{idx101.en.e27}{Affine normal form and root multiplicities}{\pageref{ffs:thm-normal-form}}
\IdxSub{idx129.en.affine}{Affine orbits of irreducible polynomials}{\pageref{v129:cor-affine-factor-orbits}}
\IdxSub{idx101.en.e38}{Affine transformations}{\pageref{ffs:def-affine-slope}}
\IdxSub{idx101.en.e61}{Base-field affine stabilizer of an irreducible polynomial}{\pageref{ffs:thm-arithmetic-cyclic}}
\IdxSub{idx101.en.e76}{Bounds on the translation stabilizer of an irreducible root set}{\pageref{ffs:thm-geometric-bound}}
\IdxSub{idx101.en.e85}{Center coset}{\pageref{ffs:prop-affine-orbits}}
\IdxSub{idx101.en.e88}{Center of a transformation}{\pageref{ffs:def-affine-center}}
\IdxSub{idx101.en.e146}{Conjugate-difference criterion for four translations}{\pageref{fas:eq-four-translation-test}}
\IdxSub{idx101.en.e36}{Defined over the base field}{\pageref{ffs:def-stabilizers}}
\IdxSub{idx101.en.e228}{Descent of the center coset to the base field}{\pageref{ffs:thm-center-descent}}
\IdxSub{idx101.en.e341}{First-order affine symmetry space}{\pageref{v67:def-first-order-affine}}
\IdxSub{idx101.en.e86}{Group structure}{\pageref{ffs:thm-affine-group}}
\IdxSub{idx101.en.e387}{Group types and the four-element exception}{\pageref{fas:cor-group-types}}
\IdxSub{idx101.en.e434}{Image of the slope map}{\pageref{ffs:def-slope-image}}
\IdxSub{idx101.en.e443}{Invariant ring}{\pageref{ffs:prop-invariant-ring}}
\IdxSub{idx101.en.e25}{Irreducible quartics}{\pageref{fas:cor-quartic-coefficients}}
\IdxSub{idx101.en.e487}{Lift}{\pageref{fas:prop-compatible-lifting}}
\IdxSub{idx101.en.e131}{Lifting compatible quotient slopes}{\pageref{fas:prop-compatible-lifting}}
\IdxSub{idx101.en.e580}{Multiplicative symmetry in prime degree}{\pageref{fas:thm-prime-tame}}
\IdxSub{idx101.en.e622}{One-dimensional affine group}{\pageref{ffs:def-stabilizers}}
\IdxSub{idx101.en.e575}{Orders and stabilizing multipliers}{\pageref{fas:prop-order-multiplier}}
\IdxSub{idx101.en.e637}{Origin-preserving scalings and full slopes}{\pageref{v54:cor-origin-slopes}}
\IdxSub{idx101.en.e760}{Recovery from coefficients}{\pageref{fas:coefficient-recovery}}
\IdxSub{idx101.en.e761}{Recovery from root values}{\pageref{fas:thm-full-recovery}}
\IdxSub{idx101.en.e241}{Same-order roots with different groups}{\pageref{fas:ex-thirteen}}
\IdxSub{idx101.en.e819}{Scalar stabilizer}{\pageref{fas:def-scalar-normalizer}}
\IdxSub{idx101.en.e830}{Semidirect product}{\pageref{def:semidirect-affine}}
\IdxSub{idx101.en.e861}{Slope map}{\pageref{ffs:def-slope-image}}
\IdxSub{idx101.en.e863}{Slope of an affine transformation}{\pageref{ffs:def-affine-slope}}
\IdxSub{idx101.en.e871}{Special orbit}{\pageref{ffs:prop-affine-orbits}}
\IdxSub{idx101.en.e942}{Trace formula for the center of the quotient}{\pageref{fas:eq-trace-center}}
\IdxSub{idx101.en.e377}{Translations from Frobenius differences}{\pageref{fas:thm-translation-difference}}
\IdxSee{Affine stabilizer of a root multiset}{see}{\hyperlink{idx101.en.affine}{Affine stabilizer}}
\IdxHead{idx101.en.agw}{AGW criterion}{}
\IdxSub{idx101.en.e571}{Multiplicative fibers}{\pageref{prop:multiplicative-agw}}
\IdxSub{idx101.en.e941}{Prescribing the base permutation}{\pageref{v98:par-prescribed-trace-permutation}}
\IdxSub{idx101.en.e39}{Statement and proof}{\pageref{lem:fiber-permutation}}
\IdxSub{idx101.en.e659}{Trace fibers}{\pageref{thm:trace-permutations}}
\IdxHead{idx101.en.algebraic}{Algebraic extensions}{\pageref{def:algebraic-extension}}
\IdxSub{idx101.en.e48}{Algebraically closed field}{\pageref{def:algebraic-closure}}
\IdxSub{idx101.en.e43}{Algebraic closure}{\pageref{def:algebraic-closure}}
\IdxSub{idx101.en.e44}{Algebraic element}{\pageref{def:algebraic-minpoly}}
\IdxSub{idx101.en.e337}{Finite extensions and algebraic generators}{\pageref{prop:finite-algebraic-generation}}
\IdxSub{idx101.en.e340}{Finitely generated extension}{\pageref{def:finitely-generated-extension}}
\IdxSub{idx101.en.e954}{Transcendental element}{\pageref{def:algebraic-extension}}
\IdxSee{Almost perfect nonlinear mapping}{See}{\hyperlink{idx150.en.apn}{Quadratic form: Almost perfect nonlinear mappings}}
\IdxHead{idx101.en.annihilator}{Annihilating polynomial}{\pageref{def:annihilating-polynomial}}
\IdxSub{idx130.en.affine}{Affine counts for arbitrary directions}{\pageref{v130:eq-stratum-affine-count}}
\IdxSub{idx106.en.annfiber}{Affine fibers}{\pageref{v106:prop-annihilator-fiber}}
\IdxSub{idx101.en.e53}{Annihilating polynomials of operators and vectors}{\pageref{def:operator-annihilator}}
\IdxSub{idx110.en.arbitrary-count}{Counts under arbitrary linear constraints}{\pageref{v110:prop-affine-count}}
\IdxSub{idx130.en.stratum}{Exact spectra of annihilator strata}{\pageref{v130:prop-annihilator-spectrum}}
\IdxSub{idx110.en.ann-info}{Exact strata and minimal linear information}{\pageref{v110:prop-annihilator-info}}
\IdxSub{idx101.en.e371}{Frobenius annihilator polynomial}{\pageref{sec:q-order}}
\IdxSub{idx107.en.ann-distribution}{Full distribution along an invariant direction}{\pageref{v107:thm-affine-annihilator-distribution}}
\IdxSub{idx101.en.e683}{Prescribed Frobenius annihilator polynomial}{\pageref{v98:par-prescribed-annihilator}}
\IdxHead{idx135.en.3}{Arithmetic Frobenius}{\pageref{revadic:thm-frob}}
\IdxHead{idx101.en.artinschreier}{Artin--Schreier equation}{}
\IdxSub{idx101.en.e56}{Equation and irreducibility}{\pageref{thm:artin-schreier}}
\IdxSub{idx101.en.e307}{Explicit solution of additive Hilbert 90}{\pageref{prop:hilbert90-explicit}}
\IdxSub{idx101.en.e57}{Recursive lifting}{\pageref{thm:recursive-as}}
\IdxGroup{B}
\IdxHead{idx101.en.bch}{BCH code}{\pageref{v56:def-bch-code}}
\IdxSub{idx101.en.e65}{BCH bound}{\pageref{v56:thm-bch-bound}}
\IdxSub{idx123.en.examples}{Binary and ternary equality examples}{\pageref{v123:ex-bch-equality-small}}
\IdxSub{idx139.en.bchfibers}{Complete root fibers and direct sums}{\pageref{v139:cor-compressed-code}}
\IdxSub{idx101.en.e153}{Constacyclic BCH bound}{\pageref{v88:cor-constacyclic-bch-bound}}
\IdxSub{idx101.en.e229}{Designed distance}{\pageref{v56:def-bch-code}}
\IdxSub{idx139.en.bchdigits}{Digit thresholds and omitted gaps}{\pageref{v139:prop-digit-bch-gap}}
\IdxSub{idx118.en.locator}{Error-locator polynomial}{\pageref{v118:def-error-locator}}
\IdxSub{idx101.en.e406}{Generator polynomial for consecutive roots of unity}{\pageref{v50:eq-leader-interval-polynomial}}
\IdxSub{idx123.en.equality}{Interpolation weights for equality}{\pageref{v123:prop-bch-equality}}
\IdxSub{idx118.en.syndrome}{Syndromes}{\pageref{v118:def-syndrome}}
\IdxSub{idx123.en.word}{Trace-prefix fibers and extremal words}{\pageref{v123:cor-paired-bch-word}}
\IdxSub{idx118.en.decode}{Two-error decoding example}{\pageref{v118:ex-bch-fifteen-decoding}}
\IdxHead{idx101.en.bm}{Berlekamp--Massey algorithm}{}
\IdxSub{idx101.en.e68}{Algorithm}{\pageref{alg:berlekamp-massey}}
\IdxSub{idx101.en.e150}{Connection length}{\pageref{def:finite-linear-complexity}}
\IdxSub{idx101.en.e151}{Connection polynomial}{\pageref{def:finite-linear-complexity}}
\IdxSub{idx101.en.e152}{Connection relation}{\pageref{def:finite-linear-complexity}}
\IdxSub{idx101.en.e245}{Discrepancy}{\pageref{def:recurrence-discrepancy}}
\IdxSub{idx101.en.e1001}{Empty and zero prefixes; connection length}{\pageref{v101:par-bm-boundaries}}
\IdxSub{idx101.en.e492}{Finite-prefix complexity}{\pageref{def:finite-linear-complexity}}
\IdxSub{idx127.en.syndrome-rank}{Syndrome rank and error count}{\pageref{v127:eq-syndrome-rank}}
\IdxSub{idx118.en.syndrome-recovery}{Syndrome recurrence and error correction}{\pageref{v118:prop-syndrome-recurrence}}
\IdxHead{idx101.en.binomial}{Binomial}{}
\IdxSub{idx101.en.e103}{Coefficient classification for irreducible binomials}{\pageref{v77:eq-binomial-constant-exponents}}
\IdxSub{idx101.en.e110}{Coefficient support and binomial splitting}{\pageref{v97:prop-support-splitting}}
\IdxSub{idx101.en.e157}{Constant terms of irreducible binomials}{\pageref{v77:def-binomial-constants}}
\IdxSub{idx101.en.e168}{Counting irreducible binomials}{\pageref{v77:eq-binomial-count}}
\IdxSub{idx101.en.e73}{Defining equation}{\pageref{v79:eq-binomial-optimum}}
\IdxSub{idx101.en.e301}{Exact splitting degree}{\pageref{v97:eq-binomial-splitting}}
\IdxSub{idx101.en.e450}{Irreducibility criterion for binomials}{\pageref{thm:binomial-full}}
\IdxSub{idx101.en.e451}{Irreducibility criterion for prime-degree binomials}{\pageref{prop:prime-binomial}}
\IdxSub{idx124.en.binomial-extension}{Least extension for binomial factorization of a general binomial}{\pageref{v124:cor-binomial-extension}}
\IdxSub{idx101.en.e545}{Minimum extension degree for binomial factorization}{\pageref{cor:med-binomial-extension}}
\IdxSub{idx101.en.e599}{Norm classes and factor degrees}{\pageref{v97:prop-norm-factor-degrees}}
\IdxSub{idx101.en.e809}{Root-of-unity polynomials with only binomial irreducible factors}{\pageref{cor:med-all-binomial}}
\IdxGroup{C}
\IdxHead{idx141.en.carryslope}{Carry frequencies and Newton slopes}{\pageref{v141:cor-orbit-carries}}
\IdxSee{Cauchy matrix}{see}{\hyperlink{idx101.en.e83}{Interpolation: Cauchy matrix}}
\IdxSee{Center of an affine transformation}{see}{\hyperlink{idx101.en.e88}{Affine stabilizer: Center of a transformation}}
\IdxHead{idx101.en.char}{Character}{}
\IdxSub{idx101.en.e7}{Additive character}{\pageref{def:characters}}
\IdxSub{idx101.en.e569}{Multiplicative character}{\pageref{def:characters}}
\IdxSub{idx101.en.e91}{Orthogonality relations}{\pageref{prop:character-orthogonality}}
\IdxHead{idx101.en.single93}{Characteristic}{\pageref{def:characteristic}}
\IdxHead{idx101.en.charsum}{Character sum}{}
\IdxSub{idx105.en.square-coefficient-count}{Exact count of quadratics with square coefficients}{\pageref{v105:prop-square-coefficient-quadratics}}
\IdxSub{idx101.en.e403}{Generating polynomial of additive character sums}{\pageref{thm:additive-L-polynomial}}
\IdxSub{idx109.en.jacobian}{Jacobian rank of the power-sum map}{\pageref{v109:ex-power-jacobian}}
\IdxSub{idx105.en.translate-character-count}{Joint character sums of translated coefficients}{\pageref{v105:eq-translate-coset-count}}
\IdxSub{idx101.en.e498}{Linear recurrence for extension-field character sums}{\pageref{cor:character-sum-recurrence}}
\IdxSub{idx101.en.e475}{L-polynomial}{\pageref{thm:additive-L-polynomial}}
\IdxSub{idx101.en.e476}{L-series}{\pageref{def:character-L-series}}
\IdxSub{idx101.en.e557}{Moment of character sums}{\pageref{v58:def-character-moment}}
\IdxSub{idx101.en.e92}{Over a trace fiber}{\pageref{def:trace-fiber-character-sum}}
\IdxSub{idx101.en.e662}{Phase reduction}{\pageref{v70:thm-phase-reduction}}
\IdxSub{idx101.en.e719}{Quadratic character sum}{\pageref{prop:quadratic-sum}}
\IdxSub{idx101.en.e803}{Root counting by character sums}{\pageref{eq:character-root-count}}
\IdxSub{idx101.en.e902}{Summation and the top coefficient}{\pageref{v91:prop-sum-top-coefficient}}
\IdxHead{idx101.en.single96}{Chevalley--Warning theorem}{\pageref{thm:chevalley-warning}}
\IdxHead{idx101.en.crt}{Chinese remainder theorem}{\pageref{thm:crt-basic}}
\IdxSub{idx101.en.e69}{Bézout identity}{\pageref{lem:bezout}}
\IdxSub{idx101.en.e112}{Comaximal ideals}{\pageref{v47:prop-comaximal-crt}}
\IdxSub{idx105.en.direct-remainders}{Direct sums and remainder reconstruction}{\pageref{v105:eq-linearized-remainder-codimension}}
\IdxSub{idx101.en.e431}{Ideal}{\pageref{v101:def-ideal}}
\IdxSub{idx105.en.pairwise-compatibility}{Pairwise right coprimality and compatibility}{\pageref{v105:ex-linearized-crt-three-lines}}
\IdxHead{idx141.en.coeffstab}{Coefficient stabilizer}{\pageref{v141:sec-phase-index}}
\IdxHead{idx106.en.collision}{Collision factors, descent}{\pageref{v106:prop-quartic-rational}}
\IdxSee{Companion matrix}{see}{\hyperlink{idx101.en.e123}{Field coordinates: Companion matrix}}
\IdxHead{idx135.en.4}{Compatible congruences}{\pageref{v135:sec-foundations}}
\IdxHead{idx135.en.5}{Complete discrete valuation ring}{\pageref{v135:sec-foundations}}
\IdxSee{Completely normal element}{see}{\hyperlink{idx101.en.e132}{Normal element: Completely normal elements}}
\IdxHead{idx101.en.composed}{Composed sums and products}{}
\IdxSub{idx101.en.e133}{Composed product}{\pageref{def:composed-products}}
\IdxSub{idx101.en.e134}{Composed product with coprime root orders}{\pageref{v40:prop-coprime-root-orders}}
\IdxSub{idx101.en.e135}{Composed sum}{\pageref{def:composed-products}}
\IdxSub{idx101.en.e472}{Kronecker product}{\pageref{prop:composed-matrices}}
\IdxSub{idx101.en.e474}{Kronecker sum}{\pageref{prop:composed-matrices}}
\IdxSub{idx101.en.e523}{Matrix computation of composed polynomials}{\pageref{prop:composed-matrices}}
\IdxSub{idx129.en.phases}{Phase images and composition multiplicities}{\pageref{v129:prop-phase-images}}
\IdxHead{idx101.en.complexity}{Computational complexity}{}
\IdxSub{idx101.en.e70}{Binary exponentiation}{\pageref{v101:def-binary-exponentiation}}
\IdxSub{idx101.en.e79}{Canonical output}{\pageref{v58:def-pseudo-deterministic}}
\IdxSub{idx101.en.e480}{Las Vegas algorithm}{\pageref{v48:def-las-vegas}}
\IdxSub{idx105.en.elliptic-preprocessing}{Preprocessing for elliptic root-finding}{\pageref{v105:rem-elliptic-root-finding}}
\IdxSub{idx101.en.e716}{Pseudo-deterministic algorithm}{\pageref{v58:def-pseudo-deterministic}}
\IdxSub{idx101.en.e868}{Sparse representation}{\pageref{v58:remark-sparse-complexity}}
\IdxHead{idx101.en.constacyclic}{Constacyclic code}{\pageref{v54:def-constacyclic-code}}
\IdxSub{idx101.en.e155}{Constacyclic MDS codes: full-support generators}{\pageref{v87:thm-constacyclic-mds-reciprocal}}
\IdxSub{idx101.en.e260}{Dual generator}{\pageref{v87:thm-constacyclic-mds-reciprocal}}
\IdxSub{idx101.en.e643}{Outer divisibility with a distinguished orbit}{\pageref{v54:thm-outer-divisibility}}
\IdxSub{idx101.en.e865}{Sparse-multiple criterion for constacyclic codes}{\pageref{v87:eq-constacyclic-dual-multiple}}
\IdxHead{idx116.en.completion-head}{Constant-derivative completion}{\pageref{v116:def-completion}}
\IdxSub{idx116.en.completion-head1}{Five-point criterion over the sixteen-element field}{\pageref{v116:ex-five-completion}}
\IdxSub{idx116.en.completion-head0}{Interpolation and the degree gap}{\pageref{v116:prop-completion-gap}}
\IdxHead{idx137.en.construction}{Construction of finite fields}{\pageref{v137:ex-construction-comparison}}
\IdxSee{Coset leader}{see}{\hyperlink{idx101.en.e167}{Cyclotomic coset: Leaders}}
\IdxHead{idx101.en.curve}{Curve point counts}{}
\IdxSub{idx101.en.e273}{Elliptic curve}{\pageref{def:binary-elliptic-curve}}
\IdxSub{idx101.en.e372}{Frobenius characteristic polynomial of an elliptic curve}{\pageref{def:curve-frobenius-polynomial}}
\IdxSub{idx101.en.e592}{Nonsingular curve}{\pageref{def:nonsingular-plane-curve}}
\IdxSub{idx101.en.e593}{Nonsingular plane curve}{\pageref{def:nonsingular-plane-curve}}
\IdxSub{idx105.en.quadratic-twist-partition}{Quadratic twists and root partitions}{\pageref{v105:eq-elliptic-root-partition}}
\IdxSub{idx101.en.e744}{Rational generating function}{\pageref{cor:binary-curve-zeta}}
\IdxSub{idx101.en.e746}{Rational point}{\pageref{def:rational-curve-point}}
\IdxSub{idx101.en.e747}{Rational point on a curve}{\pageref{def:rational-curve-point}}
\IdxSub{idx101.en.e823}{Second-order recurrence}{\pageref{cor:binary-curve-recurrence}}
\IdxSub{idx101.en.e843}{Set of rational points}{\pageref{def:rational-curve-point}}
\IdxSub{idx101.en.e990}{Zeta function}{\pageref{def:curve-zeta}}
\IdxHead{idx101.en.cycliccode}{Cyclic code}{\pageref{v40:def-cyclic-code}}
\IdxSub{idx140.en.codegroups}{Canonical lifts and permutation groups}{\pageref{v140:thm-code-group}}
\IdxSub{idx144.en.code}{Canonical subgroup check codes}{\pageref{v144:cor-subgroup-code}}
\IdxSub{idx142.en.interleaving}{Complete-fiber interleaving}{\pageref{v142:cor-irreducible-interleaving}}
\IdxSub{idx101.en.e211}{Defining set}{\pageref{v40:def-defining-set}}
\IdxSub{idx101.en.e439}{Interleaving decomposition of cyclic codes}{\pageref{v44:prop-interleaving}}
\IdxSub{idx101.en.e452}{Irreducible cyclic code}{\pageref{v40:def-irreducible-cyclic-code}}
\IdxSub{idx101.en.e857}{Simple-root cyclic code}{\pageref{v84:def-simple-root-code}}
\IdxSub{idx142.en.parity}{Single parity-check code}{\pageref{v142:ex-twenty-five-code}}
\IdxSub{idx101.en.e864}{Smallest field of definition}{\pageref{v84:def-code-field}}
\IdxSub{idx142.en.code25}{Twenty-five-coordinate example}{\pageref{v142:ex-twenty-five-code}}
\IdxHead{idx101.en.cyclicmodule}{Cyclic module}{\pageref{def:cyclic-module}}
\IdxSub{idx101.en.e191}{Cyclic decomposition}{\pageref{v53:lem-cyclic-decomposition}}
\IdxSub{idx101.en.e194}{Cyclic vector}{\pageref{def:cyclic-module}}
\IdxSub{idx101.en.e196}{Cyclic vector lemma}{\pageref{lem:cyclic-vector}}
\IdxSub{idx102.en.par-cyclic-eight}{Eight-element field: projections and components}{\pageref{v102:par-cyclic-eight}}
\IdxSub{idx101.en.e272}{Elementary divisors}{\pageref{v53:eq-elementary-divisors}}
\IdxSub{idx101.en.e347}{Fixed-point counts in cyclic modules}{\pageref{v53:eq-cyclic-gcd-count}}
\IdxSub{idx150.en.row}{Primary components and row-space intersections}{\pageref{v150:prop-row-intersection}}
\IdxSub{idx101.en.e750}{Realization of invertible linear operators by Frobenius}{\pageref{thm:frobenius-realization}}
\IdxSub{idx102.en.par-repeated-root-chain}{Repeated factors and nilpotent chains}{\pageref{v102:par-repeated-root-chain}}
\IdxSee{Cyclic vector}{see}{\hyperlink{idx101.en.e194}{Cyclic module: Cyclic vector}}
\IdxHead{idx101.en.coset}{Cyclotomic coset}{\pageref{def:cyclotomic-coset}}
\IdxSub{idx101.en.e201}{Absence of nonsingleton blocks}{\pageref{exm:med-fifteen-trivial}}
\IdxSub{idx101.en.e9}{Additive stabilizer}{\pageref{def:additive-stabilizer}}
\IdxSub{idx101.en.e10}{Additive-stabilizer criterion for multiplicative orbits}{\pageref{v74:prop-orbit-additive-stabilizer}}
\IdxSub{idx101.en.e22}{Admissible digit path}{\pageref{v44:def-admissible-digit-path}}
\IdxSub{idx101.en.e278}{Binomial descent: example}{\pageref{exm:med-binomial-twentyseven}}
\IdxSub{idx104.en.sparseboundary}{Boundary of the trinomial synchronization condition}{\pageref{v104:ex-four-sparse-boundary}}
\IdxSub{idx101.en.e99}{Coarsest equal-difference decomposition}{\pageref{thm:med-coset-partitions}}
\IdxSub{idx101.en.e100}{Coarsest partition}{\pageref{def:med-cosets}}
\IdxSub{idx142.en.fibers}{Complete fibers and irreducible composition}{\pageref{v142:prop-full-fibers}}
\IdxSub{idx101.en.e181}{CRT synchronization of local orbits}{\pageref{v44:thm-local-orbit-synchronization}}
\IdxSub{idx101.en.e254}{Divisibility filter for leaders}{\pageref{v50:eq-leader-filter}}
\IdxSub{idx101.en.e266}{Dyadic exception for equal-difference orbits}{\pageref{prop:ed-coset-criterion}}
\IdxSub{idx101.en.e276}{Equal-difference block}{\pageref{def:med-cosets}}
\IdxSub{idx101.en.e138}{Equal-difference orbit computation}{\pageref{exm:med-nine}}
\IdxSub{idx124.en.coefficient-field}{Exact coefficient fields of multiplier orbits}{\pageref{v50:prop-multiplier-blocks}}
\IdxSub{idx101.en.e532}{Leader compression}{\pageref{v50:prop-med-leader-compression}}
\IdxSub{idx101.en.e167}{Leaders}{\pageref{v50:def-coset-leader}}
\IdxSub{idx124.en.local-global}{Local and global coefficient fields}{\pageref{v124:ex-local-coefficient-fields}}
\IdxSub{idx101.en.e546}{Minimum Frobenius step}{\pageref{thm:med-minimal-step}}
\IdxSub{idx124.en.small-modulus}{Minimum step at a small modulus}{\pageref{v124:def-small-modulus}}
\IdxSub{idx101.en.e566}{Multiple equal-difference decomposition}{\pageref{def:med-cosets}}
\IdxSub{idx101.en.e582}{Multiplier orbits and factor blocks}{\pageref{v50:prop-multiplier-blocks}}
\IdxSub{idx101.en.e567}{Non-affine symmetries}{\pageref{v74:thm-med-nonaffine}}
\IdxSub{idx101.en.e594}{Nonunit indices and reduced moduli}{\pageref{exm:med-nonunit}}
\IdxSub{idx101.en.e627}{Orbit representatives and lengths modulo odd prime powers}{\pageref{v44:prop-prime-power-local-orbits}}
\IdxSub{idx101.en.e634}{Order-preserving scaling}{\pageref{v50:thm-leader-scaling}}
\IdxSub{idx139.en.cosetdigits}{Periodic digits and exact compression}{\pageref{v139:thm-period-compression}}
\IdxSub{idx126.en.phase-exact}{Phase-group exact sequence}{\pageref{v126:eq-phase-exact}}
\IdxSub{idx126.en.phase-split}{Phase kernels and splitting}{\pageref{v126:cor-phase-splitting}}
\IdxSub{idx101.en.e660}{Phase of a local orbit}{\pageref{v44:def-local-phase}}
\IdxSub{idx101.en.e661}{Phase quotient}{\pageref{v44:eq-phase-quotient}}
\IdxSub{idx125.en.phase-fiber}{Phase reduction and global branches}{\pageref{v125:prop-phase-lifting}}
\IdxSub{idx125.en.phase-invariant}{Phases and composition exponents}{\pageref{v125:cor-phase-composition}}
\IdxSub{idx101.en.e279}{Root blocks and binomials}{\pageref{prop:ed-binomial-correspondence}}
\IdxSub{idx125.en.structure}{Stabilizers, composition, and binomial descent}{\pageref{v125:thm-orbit-structure}}
\IdxSub{idx101.en.e901}{Successive equal-difference decompositions}{\pageref{exm:med-twentyseven}}
\IdxSub{idx142.en.roots25}{Twenty-fifth roots and five-element fibers}{\pageref{v142:ex-twenty-five-roots}}
\IdxLine{1.2em}{\emph{See also} \hyperlink{idx101.en.system}{Cyclotomic system}; \hyperlink{idx101.en.cyclopol}{Cyclotomic polynomial}}{}
\IdxHead{idx101.en.cyclopol}{Cyclotomic polynomial}{\pageref{def:cyclotomic-polynomial}}
\IdxSub{idx101.en.e199}{Cyclicity of unit groups modulo integers}{\pageref{prop:cyclotomic-irreducible}}
\IdxSub{idx101.en.e519}{Low-order quadratic cyclotomic lifting}{\pageref{v49:prop-low-order-cyclotomic}}
\IdxSub{idx101.en.e891}{Stable substitution of cyclotomic factors}{\pageref{v50:eq-cyclotomic-uniform-substitution}}
\IdxSub{idx101.en.e963}{Uniform stabilization of primitive-root layers}{\pageref{v50:prop-cyclotomic-uniform-stabilization}}
\IdxHead{idx135.en.6}{Cyclotomic remainder}{\pageref{revadic:thm-integral-resolvent}}
\IdxHead{idx101.en.system}{Cyclotomic system}{\pageref{v44:def-cyclotomic-system}}
\IdxSub{idx101.en.e204}{Cyclotomic factor tower}{\pageref{v44:prop-cyclotomic-factor-tower}}
\IdxSub{idx101.en.e230}{Deviation depth}{\pageref{v46:def-cyclotomic-deviation-degree}}
\IdxSub{idx101.en.e445}{Inverse limit}{\pageref{v44:def-cyclotomic-system}}
\IdxSub{idx101.en.e454}{Irreducible substitution in a stable factor tower}{\pageref{v46:cor-cyclotomic-stable-factor-substitution}}
\IdxSub{idx101.en.e511}{Local phase of a cyclotomic orbit}{\pageref{v44:def-local-phase}}
\IdxSub{idx101.en.e506}{Local systems}{\pageref{v44:def-local-cyclotomic-system}}
\IdxSub{idx101.en.e700}{Principal cyclotomic system}{\pageref{v46:def-principal-cyclotomic-system}}
\IdxSub{idx101.en.e886}{Stabilization index}{\pageref{v46:def-cyclotomic-stable-degree}}
\IdxSub{idx101.en.e889}{Stable growth of local cyclotomic systems}{\pageref{v46:thm-cyclotomic-stable-growth}}
\IdxSub{idx101.en.e925}{Three prime-lifting branches}{\pageref{v44:cor-cyclotomic-prime-lifts}}
\IdxSub{idx101.en.e926}{Total system}{\pageref{v44:def-total-cyclotomic-system}}
\IdxGroup{D}
\IdxHead{idx101.en.defining}{Defining polynomial}{\pageref{v79:sec-defining-equations}}
\IdxSub{idx101.en.e90}{Change of generator}{\pageref{v79:sec-defining-equations}}
\IdxSub{idx103.en.sixconstruction}{Complete degree-six construction}{\pageref{v103:ex-composite-six}}
\IdxSub{idx101.en.e158}{Construction by power substitution}{\pageref{cor:power-composition}}
\IdxSub{idx101.en.e159}{Construction by trace descent}{\pageref{thm:trace-descent-construction}}
\IdxSub{idx101.en.e125}{Construction hypotheses}{\pageref{tab:construction-map}}
\IdxSub{idx101.en.e369}{Fractional linear changes of roots}{\pageref{cor:mobius-transform}}
\IdxSub{idx101.en.e517}{Low-degree changes and field generation}{\pageref{v70:lem-low-degree-generation}}
\IdxSub{idx101.en.e553}{Minimum support size of defining polynomials}{\pageref{v79:def-defining-complexity}}
\IdxSub{idx101.en.e743}{Rational expressions in a field generator}{\pageref{prop:generated-fractions}}
\IdxSub{idx101.en.e745}{Rational image of a Frobenius orbit}{\pageref{v79:prop-generator-orbit-image}}
\IdxSub{idx101.en.e748}{Rational reconstruction}{\pageref{v70:prop-rational-reconstruction}}
\IdxSub{idx101.en.e911}{Symmetric degree bound for rational representation}{\pageref{v79:ex-rational-generation}}
\IdxSub{idx101.en.e914}{Tail degree}{\pageref{v79:def-defining-complexity}}
\IdxSub{idx101.en.e953}{Trace-zero generators}{\pageref{v70:prop-trace-zero-generator}}
\IdxHead{idx101.en.degreeprofile}{Degree profile}{\pageref{sfr:def-profile}}
\IdxSub{idx116.en.small}{Arithmetic progressions at small defect}{\pageref{v116:thm-small-defect}}
\IdxSub{idx137.en.quartic}{Complete binary quartic classification}{\pageref{v137:ex-quartic-classification}}
\IdxSub{idx101.en.e173}{Counting subfield generators of low representation degree}{\pageref{sfr:prop-exact-count}}
\IdxSub{idx101.en.e649}{Defect one: restrictions}{\pageref{sfr:prop-defect-one}}
\IdxSub{idx101.en.e215}{Degree of a representative}{\pageref{sfr:def-representative-degree}}
\IdxSub{idx101.en.e243}{Dimension bound}{\pageref{sfr:thm-profile-bound}}
\IdxSub{idx101.en.e213}{Jump sequence}{\pageref{sfr:def-jumps-defect}}
\IdxSub{idx115.en.defect}{Least positive defect in degree eight}{\pageref{v115:ex-defect-two}}
\IdxSub{idx101.en.e542}{Minimum degree and compositional decomposition}{\pageref{sfr:thm-optimal-composition}}
\IdxSub{idx101.en.e547}{Minimum generator-representation degree}{\pageref{sfr:def-generator-degree}}
\IdxSub{idx101.en.e676}{Power-basis filtration}{\pageref{sfr:def-degree-filtration}}
\IdxSub{idx140.en.profileprecision}{Precision and Smith exponents}{\pageref{v140:thm-profile-smith}}
\IdxSub{idx137.en.jumpsub}{Realization by multiplication and Frobenius}{\pageref{v137:thm-jump-realization}}
\IdxSub{idx101.en.e759}{Recovery of subfield representations from coefficients}{\pageref{sfr:prop-rank-recovery}}
\IdxSub{idx101.en.e784}{Representation defect}{\pageref{sfr:def-representation-defect}}
\IdxSub{idx116.en.residual}{Residual factors and defect constructions}{\pageref{v116:cor-residual-factor}}
\IdxSub{idx101.en.e219}{Transfer in subfield towers}{\pageref{sfr:prop-tower-scaling}}
\IdxSub{idx101.en.e960}{Truncated additive closure of the degree set}{\pageref{sfr:prop-truncated-closure}}
\IdxSub{idx101.en.e968}{Uniqueness of the normalized inner component}{\pageref{sfr:cor-normalized-inner}}
\IdxSub{idx101.en.e280}{Zero defect: equivalent conditions}{\pageref{sfr:thm-zero-defect}}
\IdxHead{idx101.en.descent}{Descent}{}
\IdxSub{idx101.en.e107}{Coefficient fields of factors}{\pageref{prop:factor-coefficient-field}}
\IdxSub{idx101.en.e224}{Descent, factorization, and compositional blocks}{\pageref{box:descent-factor-composition}}
\IdxSub{idx101.en.e375}{Fourier coordinates}{\pageref{v96:prop-fourier-descent}}
\IdxSub{idx101.en.e225}{Linear relation length}{\pageref{v95:eq-descent-lambda}}
\IdxSub{idx101.en.e376}{Multivariate linear factors}{\pageref{thm:multivariate-linear-descent}}
\IdxSub{idx101.en.e619}{Normalized linear factor}{\pageref{sec:multivariate-linear-descent}}
\IdxSub{idx101.en.e227}{Root sets}{\pageref{thm:root-descent}}
\IdxSub{idx101.en.e817}{Scalar extension versus change of characteristic}{\pageref{v96:rem-scalar-extension}}
\IdxLine{1.2em}{\emph{See also} \hyperlink{idx101.en.cycliccode}{Cyclic code}; \hyperlink{idx101.en.defining}{Defining polynomial}}{}
\IdxSee{Deviation depth of a cyclotomic system}{see}{\hyperlink{idx101.en.e230}{Cyclotomic system: Deviation depth}}
\IdxHead{idx101.en.dicksonmatrix}{Dickson matrix}{\pageref{def:dickson-matrix}}
\IdxSub{idx101.en.e236}{And trace adjoint}{\pageref{prop:dickson-operations}}
\IdxSub{idx101.en.e232}{Dickson determinant and permutation sign}{\pageref{v53:cor-dickson-sign}}
\IdxSub{idx101.en.e234}{Dickson matrices and cycle structure}{\pageref{v53:eq-dickson-iterate-kernel}}
\IdxSub{idx101.en.e739}{Rank and root count via the Dickson matrix}{\pageref{prop:dickson-rank}}
\IdxHead{idx101.en.dickson}{Dickson polynomial}{\pageref{def:dickson-polynomial}}
\IdxSub{idx101.en.e233}{Dickson--Frobenius equation}{\pageref{v97:thm-dickson-frobenius}}
\IdxSub{idx101.en.e657}{Permutation criterion for Dickson polynomials}{\pageref{thm:dickson-permutation}}
\IdxSub{idx101.en.e752}{Reciprocal-orbit root polynomial}{\pageref{v97:def-reciprocal-orbit-polynomial}}
\IdxHead{idx101.en.difference}{Difference rigidity}{}
\IdxSub{idx101.en.e47}{Algebraic rigidity from differences}{\pageref{v73:sec-difference-rigidity}}
\IdxSub{idx116.en.local}{Bad pairs and local stability}{\pageref{v116:prop-local-stability}}
\IdxSub{idx101.en.e82}{Carlitz--McConnel theorem}{\pageref{v72:thm-carlitz-mcconnel}}
\IdxSub{idx101.en.e202}{Cyclotomic difference-orbit preservation}{\pageref{v73:def-difference-orbit}}
\IdxSub{idx101.en.e203}{Cyclotomic difference-orbit rigidity}{\pageref{v73:thm-ring-rigidity}}
\IdxSub{idx101.en.e240}{Difference-product rigidity}{\pageref{v74:ex-exponent-rigidity}}
\IdxSub{idx101.en.e353}{Fixed-shape subfield linearity}{\pageref{v76:thm-fixed-shape-constant}}
\IdxSub{idx101.en.e385}{Frobenius shape factor}{\pageref{v76:def-frobenius-shape-factor}}
\IdxSub{idx116.en.gap}{Global gaps and linear stability}{\pageref{v116:ex-stability-gap}}
\IdxSub{idx101.en.e435}{Independent local multiplier group}{\pageref{v73:eq-independent-local-group}}
\IdxSub{idx101.en.e503}{Linearized maps characterized by parallel differences}{\pageref{v72:thm-parallel-difference}}
\IdxSub{idx101.en.e513}{Local rigidity modulo an odd prime power}{\pageref{v73:lem-odd-local-rigidity}}
\IdxSub{idx101.en.e512}{Local rigidity modulo a power of two}{\pageref{v73:lem-two-local-rigidity}}
\IdxSub{idx101.en.e554}{Mixed exponent constraints on differences}{\pageref{v74:ex-exponent-rigidity}}
\IdxSub{idx101.en.e570}{Multiplicative coset preservation}{\pageref{v72:subsec-coset-rigidity}}
\IdxSub{idx101.en.e648}{Parallel difference-quotient invariance}{\pageref{v72:thm-parallel-difference}}
\IdxSub{idx101.en.e757}{Recovery from several shapes}{\pageref{v76:ex-multishape-recovery}}
\IdxSub{idx101.en.e822}{Second-order multiplicative subgroup rigidity}{\pageref{v76:thm-second-difference-rigidity}}
\IdxSub{idx101.en.e850}{Sidon set}{\pageref{v76:rem-sidon-shape-relations}}
\IdxSub{idx101.en.e956}{Translation-power function space}{\pageref{v74:eq-translation-power-space}}
\IdxSub{idx101.en.e969}{Unit reduction kernel}{\pageref{v73:def-reduction-kernel}}
\IdxHead{idx103.en.diffset}{Difference sets}{\pageref{v103:def-difference-set}}
\IdxSub{idx101.en.e859}{Singer difference set}{\pageref{v103:prop-singer-difference-set}}
\IdxSee{Discrete logarithm}{see}{\hyperlink{idx101.en.e249}{Multiplicative group: Discrete logarithm}}
\IdxHead{idx101.en.divideddiff}{Divided difference}{}
\IdxSub{idx101.en.e352}{Fixed-shape divided difference}{\pageref{v76:def-fixed-shape-difference}}
\IdxSub{idx101.en.e821}{Second divided difference}{\pageref{v76:def-second-divided-difference}}
\IdxSub{idx101.en.e844}{Shape parameter}{\pageref{v76:def-fixed-shape-difference}}
\IdxSub{idx101.en.e845}{Shape-specific four-term relation count}{\pageref{v76:def-shape-four-term-count}}
\IdxSub{idx101.en.e921}{Three-point difference determinant}{\pageref{v76:def-three-point-determinant}}
\IdxHead{idx101.en.dualbasis}{Dual basis}{\pageref{thm:dual-basis}}
\IdxSub{idx101.en.e255}{Dual bases in a tower of fields}{\pageref{prop:dual-tower}}
\IdxSub{idx101.en.e263}{Dual of a power basis}{\pageref{thm:power-dual}}
\IdxSub{idx101.en.e824}{Self-dual basis}{\pageref{def:self-dual-basis}}
\IdxSub{idx101.en.e938}{Trace dual of a tensor product basis}{\pageref{ex:supp-composed-six}}
\IdxSee{Dual of a power basis}{see}{\hyperlink{idx101.en.e263}{Dual basis: Dual of a power basis}}
\IdxHead{idx101.en.dynamical}{Dynamical irreducibility}{}
\IdxSub{idx101.en.e177}{Critical orbit}{\pageref{def:critical-point-value}}
\IdxSub{idx101.en.e178}{Critical point}{\pageref{def:critical-point-value}}
\IdxSub{idx101.en.e180}{Critical value}{\pageref{def:critical-point-value}}
\IdxSub{idx101.en.e267}{Dynamically irreducible polynomial}{\pageref{v55:def-dynamical-irreducibility}}
\IdxSub{idx101.en.e268}{Dynamically irreducible set}{\pageref{v55:def-dynamical-set}}
\IdxSub{idx101.en.e339}{Finite test for quadratic dynamical irreducibility}{\pageref{v55:cor-quadratic-finite-test}}
\IdxSub{idx101.en.e359}{Forward-invariant-set criterion for mixed compositions}{\pageref{v55:thm-mixed-power-stability}}
\IdxSub{idx101.en.e398}{General quadratic stability criterion}{\pageref{v56:cor-general-quadratic-stability}}
\IdxSub{idx101.en.e447}{Inversely stable polynomial}{\pageref{v55:def-inverse-stability}}
\IdxSub{idx101.en.e666}{Pointed affine conjugacy}{\pageref{v55:prop-pointed-conjugacy}}
\IdxSub{idx101.en.e846}{Shifted power maps and stability}{\pageref{v55:thm-shifted-power-stability}}
\IdxSub{idx101.en.e888}{Stable cubic family in characteristic two}{\pageref{v55:cor-linwang-stable-cubic}}
\IdxSub{idx101.en.e890}{Stable polynomial}{\pageref{v55:def-dynamical-irreducibility}}
\IdxGroup{E}
\IdxHead{idx135.en.1}{Embeddings of unramified integer rings}{\pageref{v135:thm-hom-lifting}}
\IdxSee{Entrywise Frobenius}{see}{\hyperlink{idx101.en.e274}{Normal basis: Entrywise Frobenius}}
\IdxHead{idx135.en.7}{Equal-characteristic lifting}{\pageref{revadic:prop-equal}}
\IdxHead{idx101.en.evaluation}{Evaluation}{\pageref{v85term:def-evaluation}}
\IdxSub{idx101.en.e287}{At squares}{\pageref{v84:prop-differential-interpolation}}
\IdxSub{idx101.en.e288}{Codes}{\pageref{v85term:def-evaluation-code}}
\IdxSub{idx101.en.e289}{Coordinate}{\pageref{v40:def-finite-field-dft}}
\IdxSub{idx101.en.e290}{Direction}{\pageref{v84:ex-hyperoval}}
\IdxSub{idx101.en.e311}{Extended evaluation space}{\pageref{v84:prop-differential-interpolation}}
\IdxSub{idx101.en.e390}{Full evaluation set}{\pageref{v85term:def-full-evaluation-set}}
\IdxSub{idx101.en.e423}{Homogeneous evaluation}{\pageref{v84:prop-projective-evaluation}}
\IdxSub{idx101.en.e424}{Homogeneous evaluation map}{\pageref{v84:prop-projective-evaluation}}
\IdxSub{idx101.en.e438}{Injectivity of evaluation maps}{\pageref{v84:prop-projective-evaluation}}
\IdxSub{idx101.en.e477}{\(L\)-valued point}{\pageref{v85term:def-valued-point}}
\IdxSub{idx101.en.e292}{Map and homomorphism}{\pageref{v85term:def-evaluation-map}}
\IdxSub{idx101.en.e579}{Multiplicative subgroup evaluation set}{\pageref{v85term:def-subgroup-evaluation-set}}
\IdxSub{idx105.en.normalized-evaluation}{Normalized evaluation and one-dimensional remainders}{\pageref{v105:eq-linearized-point-remainder}}
\IdxSub{idx101.en.e294}{Point of an algebra}{\pageref{v85term:def-valued-point}}
\IdxSub{idx101.en.e293}{Points}{\pageref{v85term:def-evaluation-point}}
\IdxSub{idx101.en.e707}{Projective evaluation}{\pageref{v85term:def-projective-evaluation}}
\IdxSub{idx101.en.e295}{Sets}{\pageref{v85term:def-evaluation-set}}
\IdxSub{idx101.en.e296}{Spaces}{\pageref{v85term:def-evaluation-space}}
\IdxSub{idx101.en.e297}{Vector}{\pageref{v84:prop-projective-evaluation}}
\IdxSee{Evaluation homomorphism}{see}{\hyperlink{idx101.en.e292}{Evaluation: Map and homomorphism}}
\IdxSee{Evaluation point of an algebra}{see}{\hyperlink{idx101.en.e294}{Evaluation: Point of an algebra}}
\IdxHead{idx101.en.exactperiod}{Exact-period layer}{}
\IdxSub{idx101.en.e23}{Affine compression}{\pageref{thm:period-compression}}
\IdxSub{idx101.en.e26}{Affine invariant}{\pageref{def:affine-invariant}}
\IdxSub{idx101.en.e108}{Coefficient gaps}{\pageref{thm:period-gap}}
\IdxSub{idx101.en.e299}{Definition and layer polynomial}{\pageref{def:exact-period-layer}}
\IdxSub{idx101.en.e218}{Degrees of compressed factors}{\pageref{cor:compressed-factor-degrees}}
\IdxSub{idx101.en.e325}{Fiber interpretation of period reduction}{\pageref{sec:affine-return-meaning}}
\IdxSub{idx101.en.e382}{Frobenius on compressed values}{\pageref{cor:compressed-orbit-period}}
\IdxSub{idx101.en.e495}{Linear factors on compressed layers}{\pageref{ex:wu7-compressed-linear-factors}}
\IdxSub{idx101.en.e576}{Multiplicative-order strata}{\pageref{v48:prop-period-order-count}}
\IdxSub{idx101.en.e650}{Parity cancellation}{\pageref{thm:parity-period}}
\IdxSub{idx101.en.e753}{Reciprocal pairing}{\pageref{v48:eq-period-reciprocal-order}}
\IdxSub{idx101.en.e841}{Separation of orbits by affine invariants}{\pageref{prop:affine-orbit-separation}}
\IdxHead{idx141.en.phasekernels}{Exceptional phase kernels}{\pageref{v141:prop-exceptional-kernels}}
\IdxHead{idx135.en.8}{Excess local cancellation depth}{\pageref{revadic:cor-excess-degree}}
\IdxSub{idx141.en.cutoff}{Minimal-excess cutoff}{\pageref{v141:cor-phase-cutoff}}
\IdxHead{idx135.en.9}{Exterior-power trace}{\pageref{revadic:thm-orbit-exterior}}
\IdxGroup{F}
\IdxSee{Factorization after base extension}{see}{\hyperlink{idx101.en.e319}{Factorization algorithm: Factorization after base extension}}
\IdxHead{idx101.en.factor}{Factorization algorithm}{}
\IdxSub{idx101.en.e67}{Berlekamp algebra}{\pageref{thm:berlekamp}}
\IdxSub{idx101.en.e80}{Cantor--Zassenhaus algorithm}{\pageref{sec:equal-degree-factorization}}
\IdxSub{idx101.en.e212}{Degree-by-degree irreducibility test}{\pageref{prop:ben-or}}
\IdxSub{idx101.en.e251}{Distinct-degree factorization}{\pageref{sec:distinct-degree-factorization}}
\IdxSub{idx105.en.elliptic-rootfinding}{Elliptic root-finding and root partitions}{\pageref{v105:rem-elliptic-root-finding}}
\IdxSub{idx101.en.e275}{Equal-degree factorization}{\pageref{sec:equal-degree-factorization}}
\IdxSub{idx101.en.e319}{Factorization after base extension}{\pageref{thm:factor-base-extension}}
\IdxSub{idx101.en.e436}{Indicator function}{\pageref{thm:berlekamp}}
\IdxSub{idx101.en.e734}{Rabin criterion}{\pageref{thm:rabin}}
\IdxSub{idx101.en.e738}{Random roots of split polynomials}{\pageref{prop:split-values-random}}
\IdxSub{idx101.en.e881}{Square-free factorization}{\pageref{sec:squarefree-factorization}}
\IdxSub{idx101.en.e900}{Success probability}{\pageref{v48:prop-edf-success}}
\IdxSub{idx101.en.e974}{Value polynomial}{\pageref{prop:berlekamp-value-polynomial}}
\IdxSee{Family of reciprocals of translates}{see}{\hyperlink{idx101.en.e612}{Normal element: Reciprocals of translates}}
\IdxHead{idx135.en.10}{Fermat curve slopes}{\pageref{revadic:thm-fermat-slopes}}
\IdxSub{idx141.en.return}{Orbit factors and return scalars}{\pageref{v141:eq-fermat-orbit-factor}}
\IdxHead{idx101.en.fibers}{Fiber}{\pageref{sec:fibers}}
\IdxSub{idx101.en.e81}{Capelli's root-fiber criterion}{\pageref{thm:capelli}}
\IdxSub{idx101.en.e105}{Coefficient fiber}{\pageref{v85:sec-pair-coefficient-fibers}}
\IdxSub{idx101.en.e106}{Coefficient fiber of a root set}{\pageref{v85:sec-pair-coefficient-fibers}}
\IdxSub{idx101.en.e449}{Irreducibility criterion for additive compositions}{\pageref{prop:general-additive-composition}}
\IdxSub{idx101.en.e486}{Least root field}{\pageref{def:root-field-spectrum}}
\IdxSub{idx101.en.e805}{Root fiber}{\pageref{def:root-fiber}}
\IdxSub{idx101.en.e806}{Root-field spectrum}{\pageref{def:root-field-spectrum}}
\IdxSub{idx101.en.e327}{Size}{\pageref{v58:def-value-set}}
\IdxSub{idx101.en.e897}{Subnormal polynomial}{\pageref{def:root-field-spectrum}}
\IdxLine{1.2em}{\emph{See also} \hyperlink{idx101.en.return}{Return map}; \hyperlink{idx101.en.agw}{AGW criterion}}{}
\IdxHead{idx101.en.coordinates}{Field coordinates}{}
\IdxSub{idx101.en.e94}{Characteristic polynomial of a multiplication operator}{\pageref{prop:characteristic-multiplication}}
\IdxSub{idx101.en.e123}{Companion matrix}{\pageref{sec:regular-representation}}
\IdxSub{idx101.en.e145}{Conjugate coordinates}{\pageref{thm:conjugate-coordinates}}
\IdxSub{idx101.en.e165}{Coordinate representations of field elements}{\pageref{sec:field-coordinates}}
\IdxSub{idx101.en.e310}{Exponent representation}{\pageref{sec:field-coordinates}}
\IdxSub{idx101.en.e669}{Polynomial basis representation}{\pageref{sec:field-coordinates}}
\IdxSub{idx101.en.e731}{Quotient-ring representation}{\pageref{thm:quotient-model}}
\IdxSub{idx101.en.e772}{Regular representation}{\pageref{sec:regular-representation}}
\IdxSub{idx101.en.e895}{Structure constants}{\pageref{def:multiplication-constants}}
\IdxSub{idx101.en.e928}{Tower representation}{\pageref{sec:field-coordinates}}
\IdxSub{idx101.en.e987}{Worked example in the field of 64 elements}{\pageref{exm:f64-computation}}
\IdxSub{idx101.en.e988}{Zech logarithm}{\pageref{def:zech-logarithm}}
\IdxHead{idx101.en.field}{Field extension}{\pageref{def:field-extension}}
\IdxSub{idx101.en.e60}{Base field}{\pageref{def:field-extension}}
\IdxSub{idx101.en.e98}{Classification of subfields}{\pageref{thm:subfields}}
\IdxSub{idx101.en.e137}{Compositum}{\pageref{cor:intersection-compositum}}
\IdxSub{idx101.en.e160}{Construction of finite fields from root sets}{\pageref{thm:root-construction}}
\IdxSub{idx101.en.e161}{Construction of splitting extensions}{\pageref{lem:finite-splitting}}
\IdxSub{idx101.en.e313}{Degree}{\pageref{def:extension-degree}}
\IdxSub{idx101.en.e329}{Embeddings}{\pageref{def:field-embedding}}
\IdxSub{idx101.en.e315}{Extension degree and Frobenius period}{\pageref{thm:orbit}}
\IdxSub{idx101.en.e316}{Extension field}{\pageref{def:field-extension}}
\IdxSub{idx101.en.e336}{Finite extension}{\pageref{def:extension-degree}}
\IdxSub{idx101.en.e343}{Fixed field}{\pageref{prop:fixed-gcd}}
\IdxSub{idx101.en.e399}{Generated subfield}{\pageref{def:generated-subfield}}
\IdxSub{idx101.en.e440}{Intermediate fields and allowed Frobenius periods}{\pageref{prop:intermediate-periods}}
\IdxSub{idx101.en.e442}{Intersections and composita}{\pageref{cor:intersection-compositum}}
\IdxSub{idx101.en.e457}{Isomorphism of finite fields of the same order}{\pageref{thm:unique}}
\IdxSub{idx101.en.e494}{Linear disjointness}{\pageref{def:linearly-disjoint}}
\IdxSub{idx101.en.e630}{Order of a finite field}{\pageref{thm:cardinality}}
\IdxSub{idx101.en.e651}{Perfect field}{\pageref{def:perfect-field}}
\IdxSub{idx101.en.e684}{Prime field}{\pageref{def:prime-subfield}}
\IdxSub{idx101.en.e787}{Restriction of scalars}{\pageref{def:restriction-scalars}}
\IdxSub{idx101.en.e839}{Separable extension}{\pageref{def:separable-extension}}
\IdxSub{idx101.en.e856}{Simple extension}{\pageref{def:generated-subfield}}
\IdxSub{idx101.en.e875}{Splitting field}{\pageref{def:splitting-field}}
\IdxSub{idx101.en.e898}{Subring generated by one element}{\pageref{def:polynomial-subring}}
\IdxSub{idx101.en.e927}{Tower law}{\pageref{lem:tower}}
\IdxSub{idx101.en.e967}{Uniqueness of splitting fields}{\pageref{prop:splitting-unique}}
\IdxHead{idx101.en.isomorphism}{Field isomorphism}{}
\IdxSub{idx101.en.e305}{Explicit isomorphisms of finite fields}{\pageref{prop:construct-field-isomorphism}}
\IdxSub{idx101.en.e317}{Extension of isomorphisms}{\pageref{lem:extend}}
\IdxSub{idx101.en.e459}{Isomorphisms between finite-field models}{\pageref{prop:finite-field-isomorphism-choices}}
\IdxSub{idx101.en.e458}{Isomorphism type and choice of an isomorphism}{\pageref{sec:isomorphism-choices}}
\IdxSub{idx101.en.e465}{\(k\)-automorphism}{\pageref{def:base-field-morphism}}
\IdxSub{idx101.en.e466}{\(k\)-homomorphism}{\pageref{def:base-field-morphism}}
\IdxSub{idx101.en.e468}{\(k\)-isomorphism}{\pageref{def:base-field-morphism}}
\IdxHead{idx101.en.algebra}{Finite algebra}{}
\IdxSub{idx101.en.e40}{Algebra of equivariant functions}{\pageref{v102:par-four-root-components}, \pageref{thm:equivariant-functions}}
\IdxSub{idx101.en.e49}{Algebras and algebra homomorphisms}{\pageref{def:algebra-language}}
\IdxSub{idx101.en.e162}{Contravariant correspondence between algebras and orbits}{\pageref{thm:finite-orbit-duality}}
\IdxSub{idx101.en.e223}{Descending ideal filtration}{\pageref{def:ideal-filtration}}
\IdxSub{idx117.en.index}{Dimension versus nilpotency index}{\pageref{v117:ex-dimension-index}}
\IdxSub{idx101.en.e261}{Dual numbers}{\pageref{exm:dual-numbers-trace}, \pageref{v85:prop-quadratic-unit-groups}}
\IdxSub{idx101.en.e334}{Filtration}{\pageref{def:ideal-filtration}}
\IdxSub{idx101.en.e335}{Finite étale algebra}{\pageref{v102:par-component-idempotents}, \pageref{def:finite-etale-algebra}}
\IdxSub{idx101.en.e350}{Fixed points of iterates in algebras with repeated factors}{\pageref{cor:algebra-fixed-points-reduced}}
\IdxSub{idx101.en.e386}{Frobenius similarity type and algebra isomorphism type}{\pageref{fra:cor-similarity}}
\IdxSub{idx107.en.generators}{Generator number and ordered generating tuples}{\pageref{v107:prop-generating-tuples}}
\IdxSub{idx101.en.e432}{Idempotent}{\pageref{def:reduced-algebra}}
\IdxSub{idx101.en.e433}{Idempotent projections onto value fibers}{\pageref{v82:prop-value-idempotents}}
\IdxSub{idx121.en.charpoly}{Local components and multiplication characteristic polynomials}{\pageref{v121:eq-local-characteristic}}
\IdxSub{idx117.en.trace}{Local dimensions and trace weights}{\pageref{v117:prop-local-dimension-trace}}
\IdxSub{idx101.en.e515}{Local trace formula with factor multiplicities}{\pageref{prop:local-multiplicity-trace}}
\IdxSub{idx104.en.monogenic}{Monogenic criterion}{\pageref{v104:prop-monogenic-reduced}}
\IdxSub{idx101.en.e587}{Nilpotent element}{\pageref{def:reduced-algebra}}
\IdxSub{idx101.en.e588}{Nilpotent ideal}{\pageref{def:local-residue}}
\IdxSub{idx101.en.e589}{Nilradical}{\pageref{sec:algebra-periodic-part}}
\IdxSub{idx101.en.e590}{Nilradical and the periodic part of Frobenius}{\pageref{prop:algebra-periodic-reduction}}
\IdxSub{idx101.en.e639}{Orthogonal idempotents}{\pageref{def:reduced-algebra}}
\IdxSub{idx101.en.e652}{Perfect fields and traces of reduced algebras}{\pageref{rem:perfectness-trace}}
\IdxSub{idx101.en.e755}{Recovering algebra decompositions from fixed spaces}{\pageref{fra:thm-recovery}}
\IdxSub{idx101.en.e768}{Reduced algebra}{\pageref{def:reduced-algebra}}
\IdxSub{idx101.en.e770}{Reducedness and nondegeneracy of the trace pairing}{\pageref{cor:reduced-trace-pairing}}
\IdxHead{idx135.en.2}{Finite-precision root projection}{\pageref{v135:eq-finite-projection}}
\IdxSee{Fixed field}{see}{\hyperlink{idx101.en.e343}{Field extension: Fixed field}}
\IdxHead{idx135.en.11}{Fixed rings at finite precision}{\pageref{revadic:prop-fixed-precision}}
\IdxHead{idx101.en.fouriercode}{Fourier support of a cyclic code}{\pageref{v84:def-fourier-support}}
\IdxSub{idx101.en.e59}{Balanced Fourier support}{\pageref{v84:thm-balanced-support}}
\IdxSub{idx101.en.e348}{Fixed-point rigidity of multipliers}{\pageref{v84:thm-multiplier-fixed}}
\IdxSub{idx101.en.e581}{Multiplier index}{\pageref{v73:def-multiplier-index}}
\IdxSub{idx101.en.e583}{Multiplier stabilizer}{\pageref{v84:def-multiplier-stabilizer}}
\IdxSub{idx101.en.e363}{Spectral projections}{\pageref{v73:eq-fourier-orbit-projector}}
\IdxSub{idx101.en.e913}{Synchronization of orders on divisors}{\pageref{v84:cor-multiplier-orders}}
\IdxSub{idx101.en.e989}{Zero-class compression}{\pageref{v84:prop-fourier-compression}}
\IdxHead{idx101.en.fourier}{Fourier transform}{}
\IdxSub{idx101.en.e8}{Additive transform}{\pageref{def:fourier-transform}}
\IdxSub{idx101.en.e5}{Additive versus cyclic transforms}{\pageref{v96:rem-two-fourier}}
\IdxSub{idx101.en.e362}{Characteristic bound for minors}{\pageref{v95:eq-chebotarev-bound}}
\IdxSub{idx111.en.fourier-cosets}{Coset invariance and support}{\pageref{v111:prop-fourier-cosets}}
\IdxSub{idx101.en.e247}{Discrete Fourier transform}{\pageref{v40:def-finite-field-dft}}
\IdxSub{idx121.en.uncertainty}{Equality in the support bound}{\pageref{v121:par-uncertainty-equality}}
\IdxSub{idx101.en.e302}{Exceptional characteristics of minors}{\pageref{v95:eq-chebotarev-bound}}
\IdxSub{idx101.en.e338}{Finite-field values on cyclic groups}{\pageref{v40:def:cyclic-fourier}}
\IdxSub{idx102.en.par-four-character-table}{Four-element character table and trace fibers}{\pageref{v102:par-four-character-table}}
\IdxSub{idx102.en.par-fourier-four}{Four-element supports and inversion}{\pageref{v102:par-fourier-four}}
\IdxSub{idx101.en.e364}{Fourier support}{\pageref{def:fourier-support}}
\IdxSub{idx101.en.e365}{Fourier support and magnitude of quadratic phases}{\pageref{prop:quadratic-fourier-support}}
\IdxSub{idx112.en.fourier-cosets}{Full stabilizer and smallest quotient}{\pageref{v112:cor-spectral-stabilizer}}
\IdxSub{idx101.en.e360}{Inversion}{\pageref{eq:fourier}}
\IdxSub{idx111.en.quotient-fourier}{Quotient reduction}{\pageref{v111:eq-quotient-fourier}}
\IdxSub{idx101.en.e361}{Singular submatrices}{\pageref{v68:rem-singular-fourier-minor}}
\IdxSub{idx101.en.e367}{Support uncertainty}{\pageref{v90:thm-fourier-uncertainty}}
\IdxSee{Fractional linear transformation}{see}{\hyperlink{idx101.en.e369}{Defining polynomial: Fractional linear changes of roots}}
\IdxSee{Frobenius annihilator polynomial}{see}{\hyperlink{idx101.en.e371}{Annihilating polynomial: Frobenius annihilator polynomial}}
\IdxHead{idx101.en.frobderivation}{Frobenius derivation}{\pageref{v84:def-frobenius-derivation}}
\IdxSub{idx101.en.e176}{Critical degree}{\pageref{v86:eq-vanishing-kernel}}
\IdxSub{idx104.en.projectivecode}{Dimensions and gluing for projective evaluation codes}{\pageref{v104:ex-projective-gluing-code}}
\IdxSub{idx101.en.e374}{Frobenius derivation, critical kernel}{\pageref{v86:thm-frobenius-kernel}}
\IdxSub{idx101.en.e426}{Homogeneous vanishing ideal}{\pageref{v86:eq-vanishing-kernel}}
\IdxSub{idx101.en.e518}{Low-degree kernel}{\pageref{v86:thm-frobenius-kernel}}
\IdxSub{idx101.en.e708}{Projective Frobenius relations}{\pageref{v86:eq-critical-kernel}}
\IdxHead{idx101.en.frob}{Frobenius map}{\pageref{def:frobenius-first}}
\IdxSub{idx101.en.e2}{Absolute map}{\pageref{def:frobenius-maps}}
\IdxSub{idx101.en.e58}{Automorphism group}{\pageref{thm:automorphisms}}
\IdxSub{idx101.en.e122}{Commutant algebra}{\pageref{thm:operator-algebra}}
\IdxSub{idx101.en.e345}{Fixed points}{\pageref{def:frobenius-first}}
\IdxSub{idx101.en.e349}{Fixed-point sets}{\pageref{def:frobenius-first}}
\IdxSub{idx101.en.e244}{Fixed-space dimensions}{\pageref{fra:def-fixed-dimensions}}
\IdxSub{idx101.en.e383}{Frobenius on multiplicative quotient groups}{\pageref{v85:prop-coset-frobenius}}
\IdxSub{idx101.en.e633}{Order}{\pageref{prop:frobenius-order}}
\IdxSub{idx101.en.e775}{Relative map}{\pageref{def:frobenius-maps}}
\IdxSub{idx101.en.e833}{Semisimple operator}{\pageref{cor:frobenius-semisimple}}
\IdxSub{idx101.en.e835}{Semisimplicity}{\pageref{cor:frobenius-semisimple}}
\IdxGroup{G}
\IdxSee{Gabidulin code}{see}{\hyperlink{idx101.en.e394}{Rank metric: Gabidulin code}}
\IdxHead{idx101.en.galoisring}{Galois ring}{\pageref{v47:def-galois-ring}}
\IdxSub{idx101.en.e1}{Absolute Frobenius}{\pageref{v47:thm-galois-frobenius}}
\IdxSub{idx101.en.e63}{Basic irreducible polynomial}{\pageref{v47:def-basic-irreducible}}
\IdxSub{idx101.en.e77}{Canonical comaximal blocks}{\pageref{v47:prop-galois-ramified-blocks}}
\IdxSub{idx101.en.e89}{Chain ring}{\pageref{v47:def-chain-ring}}
\IdxSub{idx125.en.precision}{Coefficient precision in the modulus-160 factor tower}{\pageref{v125:ex-160-hensel}}
\IdxSub{idx105.en.coefficientwise-lift}{Coefficientwise lift}{\pageref{v105:def-coefficientwise-lift}}
\IdxSub{idx101.en.e565}{Composition of orbit unions under lifting}{\pageref{v47:prop-galois-med}}
\IdxSub{idx125.en.chain}{Composition over finite local rings}{\pageref{v125:prop-chain-composition}}
\IdxSub{idx101.en.e166}{Coprime Hensel factor lifting}{\pageref{v47:thm-hensel-coprime}}
\IdxSub{idx104.en.twist}{Coset twist in residue degree two}{\pageref{v104:ex-gr-four-twist}}
\IdxSub{idx101.en.e205}{Cyclotomic factorization}{\pageref{v47:thm-galois-cyclotomic-factorization}}
\IdxSub{idx101.en.e419}{Hensel lifting}{\pageref{v47:thm-hensel-coprime}}
\IdxSub{idx101.en.e508}{Local Hensel correction for roots of unity}{\pageref{v47:prop-galois-local-lift}}
\IdxSub{idx101.en.e514}{Local ring}{\pageref{v47:def-local-ring}}
\IdxSub{idx101.en.e644}{\(p\)-adic expansion}{\pageref{v47:eq-teichmuller-expansion}}
\IdxSub{idx105.en.permutation-lift}{Permutation lifts and derivative criterion}{\pageref{v105:prop-permutation-lifting}}
\IdxSub{idx101.en.e774}{Relative Frobenius}{\pageref{v47:def-galois-relative-frobenius}}
\IdxSub{idx101.en.e785}{Residue field}{\pageref{v47:def-local-ring}}
\IdxSub{idx101.en.e858}{Simple-root Hensel lifting}{\pageref{v47:lem-hensel-simple-root}}
\IdxSub{idx101.en.e915}{Teichmüller representative}{\pageref{v47:def-teichmuller}}
\IdxSub{idx101.en.e970}{Unramified extension}{\pageref{v47:def-unramified-extension}}
\IdxHead{idx101.en.periods}{Gaussian period}{}
\IdxSub{idx101.en.e84}{Cayley graph spectrum}{\pageref{v95:subsec-period-moments}}
\IdxSub{idx101.en.e318}{Factorial hyperplane moment}{\pageref{v95:thm-factorial-rank}}
\IdxSub{idx101.en.e397}{Of a scalar-saturated subgroup}{\pageref{v95:eq-gaussian-period}}
\IdxSub{idx101.en.e653}{Period-moment criterion for linear relation length}{\pageref{v95:cor-period-determines}}
\IdxSub{idx106.en.spectrum}{Real and integral spectra}{\pageref{v106:prop-real-integral-spectrum}}
\IdxSub{idx101.en.e818}{Scalar-saturated cyclotomic intersection number}{\pageref{v95:subsec-cyclotomic-triples}}
\IdxSub{idx101.en.e919}{Third Gaussian-period moment}{\pageref{v95:thm-period-third}}
\IdxSub{idx101.en.e933}{Trace-code weights and Gaussian periods}{\pageref{v96:eq-period-weight}}
\IdxSub{idx101.en.e943}{Trace hyperplane intersection}{\pageref{v95:prop-hyperplane-criterion}}
\IdxHead{idx101.en.gauss}{Gauss sum}{\pageref{def:gauss-sum}}
\IdxSub{idx101.en.e209}{Davenport--Hasse multiplication formula}{\pageref{v91:thm-dh-product}}
\IdxSub{idx132.en.stickelberger}{Digit-sum $p$-adic valuation}{\pageref{v132:thm-gauss-digit-valuation}}
\IdxSub{idx101.en.e416}{Hasse--Davenport lifting}{\pageref{thm:hasse-davenport}}
\IdxSub{idx145.en.normgauss}{Quadratic norm-quotient character sums}{\pageref{v145:lem-norm-one-gauss}}
\IdxHead{idx135.en.12}{Gauss sum coefficient field}{\pageref{revadic:prop-gauss-invariance}}
\IdxHead{idx101.en.generatorfunction}{Generating function}{\pageref{def:sequence-generating-function}}
\IdxSub{idx101.en.e355}{Formal exponential}{\pageref{sec:counting-euler}}
\IdxSub{idx101.en.e356}{Formal logarithm}{\pageref{sec:counting-euler}}
\IdxSub{idx101.en.e358}{Formal power series ring}{\pageref{sec:counting-euler}}
\IdxSub{idx101.en.e536}{Minimal polynomial from a reduced generating function}{\pageref{v70:prop-rational-generating-annihilator}}
\IdxHead{idx101.en.goppa}{Goppa code (binary)}{\pageref{v82:def-binary-goppa-code}}
\IdxSub{idx101.en.e544}{Derivative criterion and distance}{\pageref{v82:prop-goppa-support}}
\IdxSub{idx101.en.e907}{Support root polynomial}{\pageref{v82:def-support-root-polynomial}}
\IdxHead{idx101.en.gaction}{Group action}{}
\IdxSub{idx101.en.e370}{Free action}{\pageref{ffs:def-point-stabilizer}}
\IdxSub{idx101.en.e664}{Point stabilizer}{\pageref{ffs:def-point-stabilizer}}
\IdxGroup{H}
\IdxHead{idx101.en.hamming}{Hamming distance}{\pageref{v44:def-hamming-weight}}
\IdxSub{idx101.en.e411}{Hamming weight}{\pageref{v44:def-hamming-weight}}
\IdxSub{idx101.en.e549}{Minimum Hamming distance}{\pageref{v44:def-minimum-distance}}
\IdxSub{idx105.en.diagonal-scaling}{Nonzero column scaling}{\pageref{v105:ex-rank-column-normalization}}
\IdxHead{idx101.en.hankel}{Hankel matrix}{\pageref{str:def-hankel}, \pageref{v97:cor-period-recovery}}
\IdxSub{idx101.en.e413}{Hankel recovery from period sums}{\pageref{v97:cor-period-recovery}}
\IdxSub{idx101.en.e415}{Hankel--Vandermonde determinant}{\pageref{v88:thm-hankel-vandermonde}}
\IdxSub{idx148.en.kernel}{Kernel shifts and singular strata}{\pageref{v148:lem-hankel-shift}}
\IdxSub{idx101.en.e654}{Period sums of roots of unity}{\pageref{v97:cor-period-recovery}}
\IdxSub{idx101.en.e800}{Root-block recovery for a degree-six polynomial}{\pageref{str:ex-hankel-six}}
\IdxSub{idx101.en.e801}{Root-block recovery from a trace prefix}{\pageref{str:thm-hankel-recovery}}
\IdxHead{idx104.en.hassebound}{Hasse bound}{\pageref{v104:rem-hasse-bound}}
\IdxHead{idx101.en.hasse}{Hasse derivative}{\pageref{ffs:def-hasse}, \pageref{v97:prop-hasse-coordinates}}
\IdxSub{idx101.en.e479}{Largest common \(p\)-power divisor of root multiplicities}{\pageref{ffs:prop-common-power}}
\IdxSub{idx101.en.e418}{Matrix for divisibility}{\pageref{v67:eq-hasse-divisibility}}
\IdxSub{idx101.en.e782}{Remainder coordinates and multiplicity evaluations}{\pageref{v97:prop-hasse-coordinates}}
\IdxSub{idx101.en.e807}{Root multiplicity function}{\pageref{ffs:def-root-multiplicity}}
\IdxSee{Hensel lifting}{see}{\hyperlink{idx101.en.e419}{Galois ring: Hensel lifting}}
\IdxHead{idx135.en.13}{Hensel root correction}{\pageref{revadic:prop-root-correction}}
\IdxSub{idx143.en.error}{Error valuation and doubled precision}{\pageref{v143:eq-root-error}}
\IdxHead{idx101.en.hilbert90}{Hilbert's Theorem 90}{\pageref{thm:hilbert90}}
\IdxSub{idx136.en.hilbert}{Additive exact sequence at every precision}{\pageref{v136:cor-integral-hilbert90}}
\IdxSub{idx101.en.e308}{Explicit solution of multiplicative Hilbert 90}{\pageref{prop:multiplicative-hilbert-explicit}}
\IdxSee{Hyperoval}{see}{\hyperlink{idx101.en.e430}{Projective space: Hyperoval}}
\IdxGroup{I}
\IdxHead{idx141.en.index}{Index excess}{\pageref{v141:thm-index-budget}}
\IdxSee{Indicator function}{see}{\hyperlink{idx101.en.e436}{Factorization algorithm: Indicator function}}
\IdxHead{idx135.en.14}{Integral normal basis}{\pageref{revadic:thm-normal}}
\IdxHead{idx135.en.15}{Integral trace duality}{\pageref{revadic:prop-dual}}
\IdxSub{idx141.en.tracedual}{Period idempotent recovery}{\pageref{v141:prop-period-duality}}
\IdxHead{idx101.en.interpolation}{Interpolation}{}
\IdxSub{idx104.en.differentialcode}{A seven-coordinate MDS code from differential interpolation}{\pageref{v104:ex-differential-eight}}
\IdxSub{idx101.en.e83}{Cauchy matrix}{\pageref{v97:prop-cauchy-node-change}}
\IdxSub{idx116.en.selfdual}{Constant derivative and ordinary self-dual evaluation codes}{\pageref{v116:cor-unweighted-self-dual}}
\IdxSub{idx101.en.e242}{Differential interpolation}{\pageref{v84:prop-differential-interpolation}}
\IdxSub{idx101.en.e283}{Euler--Lagrange identity}{\pageref{v84:cor-euler-lagrange}}
\IdxSub{idx101.en.e444}{Inverse evaluation map from trace duality}{\pageref{eq:inverse-conjugate-tensor-map}}
\IdxSub{idx101.en.e478}{Lagrange interpolation}{\pageref{cor:interpolation-basic}}
\IdxSub{idx101.en.e804}{Root-node change}{\pageref{v97:prop-cauchy-node-change}}
\IdxSub{idx101.en.e952}{Trace representation of the top-degree coefficient}{\pageref{prop:coefficient-trace}}
\IdxSub{idx101.en.e441}{Weight}{\pageref{v84:prop-rs-dual}}
\IdxSee{Involution}{see}{\hyperlink{idx101.en.e448}{Permutation polynomial: Involution}}
\IdxHead{idx101.en.irreducible}{Irreducible polynomial}{}
\IdxSub{idx109.en.trace-rec}{Binary trace-count recurrence}{\pageref{v109:ex-binary-trace}}
\IdxSub{idx101.en.e169}{Counting irreducible polynomials}{\pageref{eq:mobius}}
\IdxSub{idx101.en.e284}{Euler product for counting irreducible polynomials}{\pageref{prop:finite-field-euler}}
\IdxSub{idx101.en.e304}{Existence of irreducible polynomials of every degree}{\pageref{thm:irreducible-existence}}
\IdxSub{idx101.en.e453}{Irreducibility}{\pageref{def:irreducible-polynomial}}
\IdxSub{idx105.en.irreducible-proportion}{Proportion as field size increases}{\pageref{v105:prop-irreducible-density}}
\IdxSub{idx105.en.square-coefficients}{Quadratics with square coefficients}{\pageref{v105:prop-square-coefficient-quadratics}}
\IdxSub{idx105.en.restricted-coefficients}{Restricted coefficients and translation orbits}{\pageref{v105:subsec-restricted-coefficients}}
\IdxGroup{J}
\IdxHead{idx101.en.jacobi}{Jacobi sum}{\pageref{def:jacobi-sum}}
\IdxSub{idx150.en.recovery}{Addition from a faithful character row}{\pageref{v150:prop-jacobi-recovery}}
\IdxSub{idx101.en.e298}{Exact counting for binary diagonal equations}{\pageref{prop:diagonal-jacobi-count}}
\IdxSub{idx150.en.joint}{Joint-product fibers and Newton identities}{\pageref{v150:prop-joint-products}}
\IdxSub{idx101.en.e568}{Multiple Jacobi sum}{\pageref{def:multiple-jacobi-sum}}
\IdxSub{idx150.en.partial}{Partial character tables and joint coset counts}{\pageref{v150:ex-partial-jacobi}}
\IdxSub{idx150.en.stabilizer}{Power-permutation stabilizers of a Jacobi row}{\pageref{v150:cor-jacobi-stabilizer}}
\IdxSub{idx101.en.e786}{Restricted sums and Jacobi sums}{\pageref{v85:ex-restricted-sums}}
\IdxSub{idx150.en.carries}{Valuation distributions for a fixed product character}{\pageref{v150:prop-jacobi-carries}}
\IdxSub{idx101.en.e986}{Weil-type bound for diagonal equations}{\pageref{thm:diagonal-weil-bound}}
\IdxHead{idx137.en.jump}{Jump sequence realization}{\pageref{v137:thm-jump-realization}}
\IdxGroup{K}
\IdxHead{idx101.en.kloosterman}{Kloosterman sum}{\pageref{def:kloosterman-sum}}
\IdxSub{idx101.en.e182}{Curve point counts and Kloosterman sums}{\pageref{prop:binary-curve-point-count}}
\IdxSub{idx101.en.e239}{Dickson representation of Kloosterman sums}{\pageref{cor:kloosterman-dickson}}
\IdxSub{idx101.en.e489}{Lifting recurrence for Kloosterman sums}{\pageref{thm:kloosterman-recurrence}}
\IdxGroup{L}
\IdxSee{Largest common \(p\)-power divisor of root multiplicities}{see}{\hyperlink{idx101.en.e479}{Hasse derivative: Largest common \(p\)-power divisor of root multiplicities}}
\IdxHead{idx141.en.last}{Last congruence layer}{\pageref{v141:thm-minimal-excess}}
\IdxSee{Las Vegas algorithm}{see}{\hyperlink{idx101.en.e480}{Computational complexity: Las Vegas algorithm}}
\IdxHead{idx101.en.code}{Linear code}{\pageref{v40:def-linear-code}}
\IdxSub{idx101.en.e95}{Check polynomial}{\pageref{v40:def-generator-check}}
\IdxSub{idx101.en.e101}{Codeword}{\pageref{v40:def-linear-code}}
\IdxSub{idx101.en.e163}{Coordinate direct-sum decomposition}{\pageref{v67:def-coordinate-direct-sum}}
\IdxSub{idx101.en.e164}{Coordinate-indecomposable code}{\pageref{v67:def-coordinate-direct-sum}}
\IdxSub{idx101.en.e258}{Dual code}{\pageref{v100:def-code-dual}}
\IdxSub{idx101.en.e270}{Effective support of a subcode}{\pageref{v54:def-effective-support}}
\IdxSub{idx101.en.e392}{Full-weight codeword}{\pageref{v82:def-full-weight-codeword}}
\IdxSub{idx101.en.e404}{Generator matrix}{\pageref{v100:def-generator-matrix}}
\IdxSub{idx101.en.e405}{Generator polynomial}{\pageref{v40:def-generator-check}}
\IdxSub{idx101.en.e414}{Hankel submatrix for outer divisibility}{\pageref{v67:cor-hankel-divisibility}}
\IdxSub{idx101.en.e437}{Information set}{\pageref{v100:def-information-set}}
\IdxSub{idx101.en.e655}{Permutation automorphism group}{\pageref{v73:def-common-permutation-group}}
\IdxSub{idx101.en.e656}{Permutation automorphism group of a code}{\pageref{v73:def-common-permutation-group}}
\IdxSub{idx101.en.e848}{Shortened code}{\pageref{v84:prop-mds-linear}}
\IdxSub{idx101.en.e918}{Test subcode and remainder matrix}{\pageref{v54:prop-test-subcode}}
\IdxHead{idx101.en.linearized}{Linearized polynomial}{\pageref{sec:linearized-definition}}
\IdxSub{idx101.en.e28}{Affine \(q\)-polynomial}{\pageref{def:affine-space}}
\IdxSub{idx101.en.e29}{Affine q-polynomial multiple}{\pageref{prop:affine-q-multiple}}
\IdxSub{idx105.en.compatible-remainders}{Compatibility of remainders}{\pageref{v105:prop-linearized-remainder-compatibility}}
\IdxSub{idx106.en.hidden}{Hidden linearization of binomials}{\pageref{v106:prop-hidden-linearization}}
\IdxSub{idx101.en.e469}{Kernel dimension and polynomial gcd}{\pageref{thm:linearized-gcd}}
\IdxSub{idx101.en.e470}{Kernels of a linearized polynomial over all extensions}{\pageref{v70:cor-all-extension-kernels}}
\IdxSub{idx101.en.e502}{Linearization and the Euclidean algorithm}{\pageref{eq:linearized-euclidean-gcd}}
\IdxSub{idx101.en.e505}{Linearly primitive root}{\pageref{def:linear-primitive-root}}
\IdxSub{idx101.en.e534}{Minimal linearized polynomial}{\pageref{prop:minimal-linearized}}
\IdxSub{idx101.en.e636}{Ore correspondence}{\pageref{ex:ore-dickson-comparison}}
\IdxSub{idx101.en.e674}{Polynomial representation of linear maps}{\pageref{thm:linear-map-representation}}
\IdxSub{idx101.en.e725}{Quartic polynomials with an affine root space}{\pageref{prop:affine-quartic}}
\IdxSub{idx101.en.e788}{Restrictions on affine root spaces}{\pageref{thm:affine-linearized-obstruction}}
\IdxSub{idx105.en.right-remainder}{Right composition remainder}{\pageref{v105:prop-linearized-remainder-compatibility}}
\IdxSub{idx101.en.e797}{Right division under composition}{\pageref{thm:composition-division}}
\IdxSub{idx101.en.e813}{Root sets of linearized polynomials and subfields}{\pageref{ex:wu7-rootfield-extremes}}
\IdxSub{idx101.en.e814}{Root space}{\pageref{prop:linearized-rootspace}}
\IdxSub{idx145.en.roots}{Skew right factors and return-invariant subspaces}{\pageref{v145:thm-skew-roots}}
\IdxSub{idx101.en.e876}{Splitting field}{\pageref{cor:linearized-splitting-field}}
\IdxSub{idx101.en.e878}{Splitting fields of linearized iterates}{\pageref{v68:cor-linearized-iterate-splitting}}
\IdxSub{idx101.en.e976}{Value sets of linearized polynomials over extensions}{\pageref{v70:ex-linearized-values-all-extensions}}
\IdxSee{Linearly primitive root}{see}{\hyperlink{idx101.en.e505}{Linearized polynomial: Linearly primitive root}}
\IdxHead{idx101.en.linearmap}{Linear map}{}
\IdxSub{idx107.en.affine-cycles}{Affine cycles and translation obstructions}{\pageref{v107:thm-affine-cycles}}
\IdxSub{idx101.en.e184}{Cycle counts for invertible linear maps}{\pageref{v53:thm-linear-cycle-count}}
\IdxSub{idx101.en.e185}{Cycle counts for linearized maps with base-field coefficients}{\pageref{v53:eq-base-coefficient-cycles}}
\IdxSub{idx101.en.e186}{Cycle type}{\pageref{def:permutation-cycle-type}}
\IdxSub{idx150.en.mixed}{Mixed traces and irreducible returns}{\pageref{v150:prop-mixed-return}}
\IdxSub{idx101.en.e621}{Odd linear permutations in characteristic two}{\pageref{v53:ex-F4-sign}}
\IdxSub{idx110.en.affine-order}{Order increase under translation}{\pageref{v110:cor-affine-order}}
\IdxSub{idx101.en.e853}{Permutation sign}{\pageref{v53:def-permutation-sign}}
\IdxSub{idx101.en.e832}{Semilinear map}{\pageref{def:semilinear-map}}
\IdxSub{idx101.en.e852}{Sign and determinant}{\pageref{v53:thm-linear-sign}}
\IdxSub{idx101.en.e854}{Sign of an affine permutation}{\pageref{v53:cor-affine-sign}}
\IdxHead{idx101.en.recurrence}{Linear recurrence}{}
\IdxSub{idx101.en.e128}{Compatibility of three quotient maps}{\pageref{prop:quotient-model-compatibility}}
\IdxSub{idx101.en.e129}{Compatibility of top-coefficient extraction and reduction}{\pageref{eq:coefficient-quotient-compatibility}}
\IdxSub{idx101.en.e499}{Definition and shift}{\pageref{def:sequence-shift}}
\IdxSub{idx101.en.e306}{Explicit projections onto transient and periodic components}{\pageref{eq:transient-projectors}}
\IdxSub{idx127.en.hankel-rank}{Hankel rank and recurrence order}{\pageref{v127:eq-observation-rank}}
\IdxSub{idx101.en.e483}{Least eventual period}{\pageref{def:transient-period}}
\IdxSub{idx101.en.e485}{Least period of a sequence}{\pageref{def:sequence-minpoly}}
\IdxSub{idx101.en.e529}{Maximum-length sequences}{\pageref{def:m-sequence}}
\IdxSub{idx101.en.e537}{Minimal polynomial of a relative trace sequence}{\pageref{str:thm-minimal-recurrence}}
\IdxSub{idx101.en.e538}{Minimal polynomial of a sequence}{\pageref{def:sequence-minpoly}}
\IdxSub{idx101.en.e539}{Minimal preperiod}{\pageref{def:transient-period}}
\IdxSub{idx101.en.e765}{Natural quotient maps}{\pageref{eq:coefficient-quotient-compatibility}}
\IdxSub{idx127.en.pairing}{Nondegenerate pairings and finite recovery}{\pageref{v127:lem-pairing-recovery}}
\IdxSub{idx127.en.observation}{Observation kernels and quotient algebras}{\pageref{v127:prop-observation-kernel}}
\IdxSub{idx101.en.e728}{Quotient-algebra model}{\pageref{thm:recurrence-algebra}}
\IdxSub{idx101.en.e764}{Recurrence polynomial}{\pageref{str:def-trace-recurrence}}
\IdxSub{idx101.en.e766}{Recurrence spaces and root spaces of linearized polynomials}{\pageref{prop:recurrence-root-correspondence}}
\IdxSub{idx101.en.e126}{Squaring versus multiplication}{\pageref{exm:mixed-algebra-recurrence}}
\IdxSub{idx101.en.e892}{State of a sequence}{\pageref{def:sequence-shift}}
\IdxSub{idx101.en.e893}{State space}{\pageref{def:sequence-shift}}
\IdxSub{idx101.en.e951}{Trace representation}{\pageref{thm:sequence-trace}}
\IdxLine{1.2em}{\emph{See also} \hyperlink{idx101.en.bm}{Berlekamp--Massey algorithm}; \hyperlink{idx101.en.generatorfunction}{Generating function}}{}
\IdxHead{idx135.en.16}{Local cyclotomic factor test}{\pageref{revadic:prop-local-factor-test}}
\IdxSee{Logarithmic derivative}{see}{\hyperlink{idx101.en.e516}{Polynomial: Logarithmic derivative}}
\IdxSee{L-polynomial}{see}{\hyperlink{idx101.en.e475}{Character sum: L-polynomial}}
\IdxSee{L-series}{see}{\hyperlink{idx101.en.e476}{Character sum: L-series}}
\IdxHead{idx101.en.lucas}{Lucas sequence}{\pageref{v45:def-lucas-polynomials}}
\IdxSub{idx146.en.lucasfactor}{Atom factorization of the Lucas recurrence}{\pageref{v146:prop-lucas-factorization}}
\IdxSub{idx151.en.lucaschar}{Characteristic-power indices and repeated-root zeros}{\pageref{v151:cor-lucas-frobenius}}
\IdxSub{idx101.en.e109}{Coefficient recurrence for root-of-unity trace parameters}{\pageref{v49:eq-coefficient-recurrence}}
\IdxSub{idx146.en.tracepoly}{Exact-order trace polynomial}{\pageref{v146:eq-order-trace}}
\IdxSub{idx101.en.e342}{First zero term of positive index in a Lucas sequence}{\pageref{v45:prop-lucas-first-zero}}
\IdxSub{idx101.en.e520}{Lucas atom}{\pageref{v45:def-lucas-atom}}
\IdxSub{idx101.en.e521}{Lucas polynomials of the first kind}{\pageref{v45:def-lucas-polynomials}}
\IdxSub{idx146.en.secondkind}{Lucas polynomials of the second kind}{\pageref{v146:prop-lucas-factorization}}
\IdxSub{idx148.en.optimal}{Optimal pullback through return coefficients}{\pageref{v148:thm-optimal-pullback}}
\IdxSub{idx146.en.pullback}{Quadratic pullback of a Lucas atom}{\pageref{v146:prop-lucas-pullback}}
\IdxSub{idx146.en.adic}{Root differences of lifted Lucas atoms}{\pageref{v146:prop-adic-lucas}}
\IdxSub{idx103.en.symtwo}{Two-variable symmetric polynomials}{\pageref{v103:lem-symmetric-two}}
\IdxGroup{M}
\IdxSee{Maximal field of scalars}{see}{\hyperlink{idx101.en.e525}{Subspace: Maximal field of scalars}}
\IdxSee{Maximum rank distance code}{see}{\hyperlink{idx101.en.e528}{Rank metric: Maximum rank distance code}}
\IdxHead{idx101.en.mds}{MDS code}{\pageref{v84:def-mds}}
\IdxSub{idx101.en.e860}{Singleton bound}{\pageref{v84:prop-mds-linear}}
\IdxSub{idx101.en.e866}{Sparse-multiple MDS criterion}{\pageref{v84:thm-sparse-multiple}}
\IdxSee{Minimal binomial multiple}{see}{\hyperlink{idx101.en.e533}{Order of a polynomial: Minimal binomial multiple}}
\IdxSee{Minimal polynomial}{see}{\hyperlink{idx101.en.e535}{Polynomial: Minimal polynomial}}
\IdxSee{Minimal preperiod}{see}{\hyperlink{idx101.en.e539}{Linear recurrence: Minimal preperiod}}
\IdxSee{Minimal preperiod and least eventual period}{see}{\hyperlink{idx101.en.e539}{Linear recurrence: Minimal preperiod}}
\IdxSee{Minimal value set polynomial}{see}{\hyperlink{idx101.en.e541}{Value set: Minimal value set polynomial}}
\IdxSee{Minimum generator-representation degree}{see}{\hyperlink{idx101.en.e547}{Degree profile: Minimum generator-representation degree}}
\IdxHead{idx135.en.17}{Minimum local valuation}{\pageref{revadic:thm-local-min}}
\IdxSee{Minimum rational generator-representation degree}{see}{\hyperlink{idx101.en.e552}{Rational representation: Minimum rational generator-representation degree}}
\IdxSee{Monogenic algebra}{see}{\hyperlink{idx104.en.monogenic}{Finite algebra: Monogenic criterion}}
\IdxHead{idx101.en.moore}{Moore determinant}{\pageref{thm:moore}}
\IdxSub{idx105.en.moore-vandermonde-scaling}{Column normalization and Vandermonde matrices}{\pageref{v105:eq-rank-column-normalization}}
\IdxSub{idx101.en.e562}{Moore--circulant determinant identity}{\pageref{v88:prop-moore-circulant}}
\IdxSub{idx101.en.e763}{Recurrence for two-variable Moore determinants}{\pageref{prop:binary-moore-recursion}}
\IdxSub{idx101.en.e564}{Resultant identity}{\pageref{v88:eq-moore-resultant}}
\IdxSub{idx101.en.e977}{Vandermonde and Moore determinants}{\pageref{v96:rem-vandermonde-moore}}
\IdxSee{Multiplicative character}{see}{\hyperlink{idx101.en.e569}{Character: Multiplicative character}}
\IdxHead{idx101.en.multiplicative}{Multiplicative group}{}
\IdxSub{idx101.en.e114}{Common period of nonidentity roots of unity}{\pageref{v124:prop-uniform-cyclotomic}, \pageref{v40:ex-uniform-root-degrees}}
\IdxSub{idx101.en.e198}{Cyclicity}{\pageref{thm:cyclic}}
\IdxSub{idx101.en.e249}{Discrete logarithm}{\pageref{def:discrete-logarithm}}
\IdxSub{idx101.en.e253}{Distinct-element products in a cyclic group}{\pageref{v85:ex-distinct-products}}
\IdxSub{idx101.en.e309}{Exponent of a group}{\pageref{lem:exponent}}
\IdxSub{idx101.en.e421}{Hidden subfield of a multiplicative subgroup}{\pageref{v74:rem-hidden-subfield}}
\IdxSub{idx101.en.e573}{Multiplicative order after centering}{\pageref{fas:def-centered-order}}
\IdxSub{idx101.en.e574}{Multiplicative order and degree of an element}{\pageref{prop:degree-from-order}}
\IdxSub{idx150.en.smallfrob}{Small-exponent Frobenius and additive indecomposability}{\pageref{v150:prop-small-frobenius}}
\IdxSub{idx101.en.e877}{Splitting field of roots of unity}{\pageref{prop:unity-splitting}}
\IdxSub{idx150.en.subsetproducts}{Uniform subset-product fibers}{\pageref{v150:prop-subset-products}}
\IdxHead{idx101.en.scaling}{Multiplicative stabilizer}{\pageref{ffs:def-stabilizers}}
\IdxSub{idx111.en.order-stabilizer}{Exact-order sets}{\pageref{v111:eq-order-stabilizer}}
\IdxHead{idx101.en.subgroup}{Multiplicative subgroup, short relations}{}
\IdxSub{idx101.en.e206}{Cyclotomic full-degree subgroup rigidity}{\pageref{v95:prop-normal-cyclotomic}}
\IdxSub{idx101.en.e220}{Dependent triple counting formula}{\pageref{v95:thm-triple-cyclotomic}}
\IdxSub{idx101.en.e259}{Dual distance and linear relation length}{\pageref{v96:prop-relation-kernel}}
\IdxSub{idx101.en.e389}{Full-degree subgroup rigidity}{\pageref{v84:thm-subgroup-rigidity}}
\IdxSub{idx101.en.e500}{Linear relation length}{\pageref{v95:def-linear-relation-length}}
\IdxSub{idx101.en.e526}{Maximal linear general position}{\pageref{v95:def-linear-relation-length}}
\IdxSub{idx101.en.e591}{Nineteen-element subgroup counterexample}{\pageref{v96:ex-nineteen}}
\IdxSub{idx101.en.e608}{Normal-basis decomposition of a multiplicative subgroup}{\pageref{v95:prop-normal-cyclotomic}}
\IdxSub{idx101.en.e617}{Normality in a multiplicative subgroup}{\pageref{v84:thm-subgroup-rigidity}}
\IdxSub{idx101.en.e629}{Orbitwise normality versus joint independence}{\pageref{v96:ex-nineteen}}
\IdxSub{idx101.en.e635}{Order-seven subgroup example}{\pageref{v95:ex-seven}}
\IdxSub{idx101.en.e847}{Short relations and order bound}{\pageref{v95:prop-order-lambda}}
\IdxSub{idx142.en.distance}{Subfield intersections and distance loss}{\pageref{v142:eq-subfield-distance}}
\IdxSub{idx106.en.cover}{Subfield translation covers}{\pageref{v106:prop-subfield-cover}}
\IdxSub{idx101.en.e965}{Unique weighted sums in a subgroup}{\pageref{v95:prop-additive-unique}}
\IdxGroup{N}
\IdxSee{Necklace}{see}{\hyperlink{idx101.en.e585}{Orbit: Necklaces}}
\IdxSee{Necklace polynomial}{see}{\hyperlink{idx105.en.necklace-polynomial}{Orbit: Necklace polynomial}}
\IdxHead{idx135.en.18}{Newton polygon of a product}{\pageref{revadic:prop-valuation-product}}
\IdxSee{Nonzero trace and normality}{see}{\hyperlink{idx101.en.e595}{Normal element: Nonzero-trace criterion}}
\IdxSee{Nonzero-trace criterion without a degree hypothesis}{see}{\hyperlink{idx101.en.e596}{Normal element: Nonzero trace without a degree hypothesis}}
\IdxHead{idx101.en.norm}{Norm}{\pageref{def:trace-norm}}
\IdxSub{idx152.en.normindecomp}{Additive indecomposability of the norm-one group}{\pageref{v150:ex-norm-one-indecomposable}}
\IdxSub{idx101.en.e16}{Adjacent norm fibers}{\pageref{v82:ex-adjacent-norm-fibers}}
\IdxSub{idx153.en.normcode}{Binary relation code of a norm-one group}{\pageref{v153:cor-normone-binary-code}}
\IdxSub{idx101.en.e600}{Compatibility condition}{\pageref{v58:def-norm-compatibility}}
\IdxSub{idx145.en.cubic}{Cubic trace--norm fibers}{\pageref{v145:cor-cubic-trace-norm}}
\IdxSub{idx101.en.e601}{Descent on conjugate root pairs}{\pageref{v45:def-root-pair-norm-descent}}
\IdxSub{idx151.en.norminter}{Exact intersections of norm-one translates}{\pageref{v151:prop-normone-intersections}}
\IdxSub{idx101.en.e344}{Fixed-norm primitive quadratic count}{\pageref{v45:prop-fixed-constant-primitive-count}}
\IdxSub{idx105.en.norm-minimal-values}{Frobenius difference after a norm}{\pageref{v105:prop-norm-trace-values}}
\IdxSub{idx101.en.e464}{Joint fiber}{\pageref{v58:def-joint-norm-fiber}}
\IdxSub{idx101.en.e497}{Linear functionals and projective unit fibers}{\pageref{v85:prop-functional-unit-fiber}}
\IdxSub{idx101.en.e556}{Möbius parametrization of the norm-one group}{\pageref{v85:prop-norm-one-parametrization}}
\IdxSub{idx101.en.e604}{Norm-one group: parametrization}{\pageref{v85:prop-norm-one-parametrization}}
\IdxSub{idx101.en.e605}{Norm-one subgroup: definition}{\pageref{v48:def-norm-one-group}}
\IdxSub{idx101.en.e606}{Norm-preserving linear maps}{\pageref{v74:note-norm-preservers}}
\IdxSub{idx101.en.e776}{Relative norm}{\pageref{def:trace-norm}}
\IdxSub{idx153.en.normsums}{Short zero sums in a norm-one group}{\pageref{v153:prop-normone-short-sums}}
\IdxHead{idx101.en.normalbasis}{Normal basis}{\pageref{def:normal-basis}}
\IdxSub{idx101.en.e139}{Computing normal elements and dual bases}{\pageref{sec:construct-normal-dual}}
\IdxSub{idx101.en.e144}{Conjugate-coefficient polynomial}{\pageref{v88:prop-conjugate-correlation}}
\IdxSub{idx101.en.e189}{Cyclic correlation polynomial}{\pageref{prop:normal-circulant-test}}
\IdxSub{idx101.en.e190}{Cyclic-correlation realization theorem}{\pageref{v88:thm-prescribed-correlation}}
\IdxSub{idx101.en.e226}{Descent of normal bases}{\pageref{prop:normal-descent}}
\IdxSub{idx101.en.e262}{Dual basis}{\pageref{prop:dual-normal}}
\IdxSub{idx101.en.e274}{Entrywise Frobenius}{\pageref{v97:prop-absolute-frobenius-coordinates}}
\IdxSub{idx101.en.e609}{Existence theorem}{\pageref{thm:normal-basis}}
\IdxSub{idx104.en.complexity}{Multiplication complexity lower bound}{\pageref{v104:lem-normal-complexity}}
\IdxSub{idx101.en.e610}{Normal coordinates of the absolute Frobenius}{\pageref{v97:prop-absolute-frobenius-coordinates}}
\IdxSub{idx101.en.e682}{Prescribed cyclic correlation}{\pageref{v88:thm-prescribed-correlation}}
\IdxSub{idx101.en.e825}{Self-dual normal basis}{\pageref{thm:self-dual-normal}}
\IdxSub{idx101.en.e962}{Type I optimal basis}{\pageref{v90:thm-gauss-optimal}}
\IdxLine{1.2em}{\emph{See also} \hyperlink{idx101.en.normality}{Normal element}; \hyperlink{idx101.en.pn}{Primitive normal basis}}{}
\IdxHead{idx101.en.normality}{Normal element}{\pageref{def:normal-element}}
\IdxSub{idx101.en.e14}{Additive stabilizer of the normal set}{\pageref{v93:prop-normal-linear-info}}
\IdxSub{idx128.en.affine}{Affine fiber counts from frequency intersections}{\pageref{v128:cor-affine-spectrum}}
\IdxSub{idx101.en.e378}{Affine space with Frobenius-invariant direction}{\pageref{v88:thm-invariant-affine-normal-count}}
\IdxSub{idx101.en.e33}{All-normal affine spaces}{\pageref{v93:def-normal-affine-space}}
\IdxSub{idx106.en.complete}{Base-field translation of completely normal elements}{\pageref{v106:lem-complete-translate}}
\IdxSub{idx101.en.e62}{Base-field translations and normality}{\pageref{v82:prop-normal-translates}}
\IdxSub{idx101.en.e132}{Completely normal elements}{\pageref{v82:def-completely-normal}}
\IdxSub{idx117.en.primary}{Component ranks and reciprocal-family bounds}{\pageref{v117:cor-primary-rank-bound}}
\IdxSub{idx101.en.e170}{Counting}{\pageref{thm:normal-count}}
\IdxSub{idx101.en.e616}{Criterion}{\pageref{prop:normal-test}}
\IdxSub{idx131.en.essential}{Essential conditions and their reduction}{\pageref{v131:prop-essential-reduction}}
\IdxSub{idx131.en.translate.count}{Exact count of completely normal base-field translates}{\pageref{v131:cor-base-complete-count}}
\IdxSub{idx128.en.spectrum}{Exact spectrum of the normality indicator}{\pageref{v128:prop-local-spectrum}}
\IdxSub{idx115.en.reciprocal}{Exceptional parameters in the cyclic module}{\pageref{v115:prop-primary-reciprocal}}
\IdxSub{idx131.en.fibre}{Fiber-sum criterion across intermediate fields}{\pageref{v131:thm-fibre-normal}}
\IdxSub{idx131.en.f21}{Four components in degree twenty-one}{\pageref{v131:ex-f2-21-components}}
\IdxSub{idx101.en.e379}{Frobenius-invariant hyperplane}{\pageref{ex:normal-p-power}}
\IdxSub{idx101.en.e388}{Full-degree nonnormal element}{\pageref{exm:full-degree-nonnormal-f64}}
\IdxSub{idx131.en.pair}{Independent pair-orbit coordinates}{\pageref{v131:prop-pair-coordinates}}
\IdxSub{idx131.en.translate.relative}{Intermediate-field translation families}{\pageref{v131:prop-intermediate-translate-count}}
\IdxSub{idx131.en.joint}{Joint ranks and exact component probabilities}{\pageref{v131:prop-component-exact}}
\IdxSub{idx101.en.e510}{Local normality condition}{\pageref{v93:def-local-normality}}
\IdxSub{idx128.en.local}{Local obstructions and trace duality}{\pageref{prop:normal-test}}
\IdxSub{idx101.en.e527}{Maximal affine dimension}{\pageref{v93:thm-patch}}
\IdxSub{idx101.en.e550}{Minimal linear information}{\pageref{v93:prop-normal-linear-info}}
\IdxSub{idx101.en.e595}{Nonzero-trace criterion}{\pageref{thm:trace-normal-classification}}
\IdxSub{idx101.en.e596}{Nonzero trace without a degree hypothesis}{\pageref{v88:cor-trace-universal-normal}}
\IdxSub{idx101.en.e615}{Normality boundary for joint trace fibers}{\pageref{v93:thm-ordinary-boundary}}
\IdxSub{idx101.en.e618}{Normality-preserving scalings}{\pageref{v93:lem-scaling}}
\IdxSub{idx101.en.e620}{Normalizing full-degree elements}{\pageref{v82:ex-universal-normal-translation}}
\IdxSub{idx106.en.resolvent}{Partial resolvents and recovery}{\pageref{v106:prop-partial-resolvents}}
\IdxSub{idx101.en.e614}{Prescribed intermediate traces}{\pageref{v58:thm-prescribed-normal-traces}}
\IdxSub{idx109.en.ktrace}{Prescribed trace and k-normal counts}{\pageref{v109:thm-knormal-trace}}
\IdxSub{idx101.en.e751}{Reciprocal changes of generators}{\pageref{v68:ex-reciprocal-normal-f343}}
\IdxSub{idx128.en.power}{Reciprocal exceptions in characteristic-power degree}{\pageref{v128:eq-power-degree-exceptions}}
\IdxSub{idx101.en.e612}{Reciprocals of translates}{\pageref{v68:prop-reciprocal-normal-family}, \pageref{v114:sec-reciprocal-exceptions}}
\IdxSub{idx114.en.extrema}{Reciprocals of translates: average and maximum}{\pageref{v114:par-extremal-question}, \pageref{v114:ex-reciprocal-average}}
\IdxSub{idx114.en.exceptions}{Reciprocals of translates: exceptional parameters}{\pageref{v114:prop-reciprocal-count}}
\IdxSub{idx101.en.e781}{Relative trace reduction of normality}{\pageref{v93:thm-trace}}
\IdxSub{idx112.en.resolvent}{Relative trace retained by full resolvents}{\pageref{v112:eq-resolvent-trace}}
\IdxSub{idx131.en.components}{Resolvent components and product formula}{\pageref{v131:thm-component-product}}
\IdxSub{idx101.en.e613}{Subfield translates}{\pageref{v82:ex-subfield-translates}}
\IdxLine{1.2em}{\emph{See also} \hyperlink{idx101.en.normalbasis}{Normal basis}; \hyperlink{idx101.en.trace}{Trace}}{}
\IdxHead{idx137.en.normdepth}{Norm constraint on cancellation depths}{\pageref{v137:prop-norm-depth}}
\IdxSee{Norm expression for the power-residue symbol}{see}{\hyperlink{idx101.en.e603}{Power-residue symbol: Norm expression}}
\IdxSee{Norm-one group}{see}{\hyperlink{idx101.en.e604}{Norm: Norm-one group: parametrization}}
\IdxSee{Norm-primitive polynomial}{see}{\hyperlink{idx148.en.primitive}{Skew polynomial: Norm-primitivity}}
\IdxGroup{O}
\IdxSee{Optimal determining polynomial}{see}{\hyperlink{idx101.en.e624}{Primitive element: Optimal determining polynomial}}
\IdxHead{idx101.en.orbits}{Orbit}{\pageref{def:orbit-conjugacy}}
\IdxSub{idx101.en.e143}{Conjugacy}{\pageref{def:orbit-conjugacy}}
\IdxSub{idx101.en.e149}{Conjugate roots and Frobenius orbits}{\pageref{thm:orbit}}
\IdxSub{idx101.en.e281}{Equivariant map}{\pageref{def:equivariant-map}}
\IdxSub{idx101.en.e285}{Euler product for periodic orbit counting}{\pageref{prop:orbit-euler}}
\IdxSub{idx101.en.e384}{Frobenius orbit phase}{\pageref{v40:def-composed-phase}}
\IdxSub{idx101.en.e555}{Möbius function}{\pageref{def:mobius-function}}
\IdxSub{idx105.en.necklace-polynomial}{Necklace polynomial}{\pageref{v105:def-necklace-polynomial}}
\IdxSub{idx101.en.e585}{Necklaces}{\pageref{v40:def-necklace}}
\IdxSub{idx101.en.e667}{Pointwise periodicity and a common period}{\pageref{prop:pointwise-uniform-period}}
\IdxSub{idx101.en.e698}{Primitive words}{\pageref{v40:def-primitive-word}}
\IdxSub{idx105.en.necklace-product}{Product formula for synchronized orbit counts}{\pageref{v105:eq-necklace-product}}
\IdxSub{idx101.en.e626}{Quotient}{\pageref{def:orbit-quotient-free}}
\IdxHead{idx135.en.19}{Orbit blocks and slopes}{\pageref{revadic:prop-orbit-slope}}
\IdxHead{idx101.en.order}{Order of a polynomial}{\pageref{def:poly-order}}
\IdxSub{idx101.en.e74}{Binomial order of a polynomial}{\pageref{v45:def-binomial-order}}
\IdxSub{idx101.en.e490}{Lifting valuations and multiplicative orders at primes}{\pageref{lem:valuation-lift}}
\IdxSub{idx101.en.e530}{Maximum monomial composition exponent}{\pageref{def:coefficient-extraction}}
\IdxSub{idx101.en.e533}{Minimal binomial multiple}{\pageref{v40:def-minimal-binomial}}
\IdxSub{idx101.en.e712}{Projective order}{\pageref{def:projective-order}}
\IdxSub{idx101.en.e783}{Repeated factors and polynomial order}{\pageref{thm:order-multiplicity}}
\IdxSub{idx101.en.e812}{Root ratios and projective order}{\pageref{v40:prop-minimal-binomial}}
\IdxSee{Ore correspondence}{see}{\hyperlink{idx101.en.e636}{Linearized polynomial: Ore correspondence}}
\IdxSee{Orthogonal complement}{see}{\hyperlink{idx101.en.e638}{Trace pairing: Orthogonal complement}}
\IdxGroup{P}
\IdxSee{p-adic expansion}{see}{\hyperlink{idx101.en.e644}{Galois ring: \(p\)-adic expansion}}
\IdxHead{idx132.en.padics}{$p$-adic lifts}{\pageref{v132:thm-unramified-limit}}
\IdxSub{idx140.en.subspaces}{Additive closure of digit subspaces}{\pageref{v140:thm-digit-subspaces}}
\IdxSub{idx134.en.base}{Base change and valuation profiles}{\pageref{v134:prop-base-change-profile}}
\IdxSub{idx139.en.bch}{BCH bounds and carry steps}{\pageref{v139:prop-digit-bch-gap}}
\IdxSub{idx140.en.groups}{Canonical code permutation groups}{\pageref{v140:thm-code-group}}
\IdxSub{idx144.en.closure}{Carry closure and generated subfields}{\pageref{v144:cor-carry-closure}}
\IdxSub{idx139.en.binomial}{Complete digit fibers and binomials}{\pageref{v139:prop-digit-binomial}}
\IdxSub{idx139.en.code}{Complete fibers and interleaved codes}{\pageref{v139:cor-compressed-code}}
\IdxSub{idx138.en.steps}{Digit grouping and Frobenius steps}{\pageref{v138:prop-period-steps}}
\IdxSub{idx138.en.rotation}{Digit rotation and cyclotomic cosets}{\pageref{v138:eq-digit-rotation}}
\IdxSub{idx139.en.orbit}{Digit sums and orbit sums}{\pageref{v139:prop-orbit-digit-sum}}
\IdxSub{idx133.en.equal}{Equal-characteristic ring section and formal series}{\pageref{v133:subsec-equichar}}
\IdxSub{idx139.en.compression}{Exact period compression}{\pageref{v139:thm-period-compression}}
\IdxSub{idx134.en.stabilizer}{Exact zero coefficients of canonical lifts}{\pageref{v134:prop-stabilizer-zero}}
\IdxSub{idx140.en.rigidity}{First carry and multiplicative rigidity}{\pageref{v140:thm-teich-rigidity}}
\IdxSub{idx144.en.obstruction}{First-carry obstruction on partial root sets}{\pageref{v144:prop-partial-lift}}
\IdxSub{idx134.en.carry}{First Witt carry for canonical coefficients}{\pageref{v134:prop-first-coefficient-carry}}
\IdxSub{idx133.en.carry}{First Witt carry polynomial}{\pageref{v133:prop-first-carry}}
\IdxSub{idx133.en.deform}{Formal Hensel lift and Artin--Schreier deformation}{\pageref{v133:prop-equal-hensel}, \pageref{v133:ex-artin-schreier-deformation}}
\IdxSub{idx134.en.phase}{Gauss-period phases and Fourier decomposition}{\pageref{v134:eq-phase-gauss-fourier}}
\IdxSub{idx132.en.cosets}{High-precision coefficients of cyclotomic cosets}{\pageref{v132:problem-coset-valuations}}
\IdxSub{idx132.en.seven}{Lifted seventh-root periods}{\pageref{v132:ex-seven-period}}
\IdxSub{idx133.en.obstruction}{Lifted zero-sum obstruction}{\pageref{v133:prop-zero-sum-obstruction}}
\IdxSub{idx139.en.leaders}{Minimal rotations and coset leaders}{\pageref{v139:prop-digit-leaders}}
\IdxSub{idx144.en.rigidity}{Mixed-exponent separation and semilinearity}{\pageref{v144:cor-partial-rigidity}}
\IdxSub{idx132.en.newton}{Newton polygons and elliptic-curve slopes}{\pageref{v132:subsec-newton-curves}}
\IdxSub{idx134.en.newton}{Newton power sums and valuation profiles}{\pageref{v134:prop-newton-profile}}
\IdxSub{idx134.en.subset}{Orbit subset sums and lifted coefficients}{\pageref{v134:thm-subset-periods}}
\IdxSub{idx138.en.gauss}{Periodic digit blocks and Gauss sums}{\pageref{v138:cor-word-gauss}}
\IdxSub{idx140.en.depth}{Permutation preservation depth}{\pageref{v140:prop-permutation-depth}}
\IdxSub{idx140.en.profile}{Precision profiles of representations}{\pageref{v140:thm-profile-smith}}
\IdxSub{idx138.en.periodic}{Purely periodic expansion}{\pageref{v138:thm-periodic-cosets}}
\IdxSub{idx139.en.phase}{Relative phases of periodic words}{\pageref{v139:prop-digit-crt}}
\IdxSub{idx138.en.seven}{Seventh-root Gauss sums and periods}{\pageref{v138:eq-seven-gauss-period}}
\IdxSub{idx137.en.sevennewton}{Seventh-root Newton recurrence and precision}{\pageref{revadic:ex-period-count}}
\IdxSub{idx137.en.sevenpair}{Seventh-root trace-pairing precision}{\pageref{revadic:ex-eta}}
\IdxSub{idx132.en.slopes}{Slope questions for exponential sums}{\pageref{v132:problem-gauss-slopes}}
\IdxSub{idx132.en.digits}{Teichmüller expansions and valuations}{\pageref{v132:eq-adic-expansion}}
\IdxSub{idx147.en.directions}{Translation directions and maximal fields of scalars}{\pageref{v147:thm-carry-directions}}
\IdxSub{idx140.en.translations}{Translation stabilizers of digit sets}{\pageref{v140:prop-digit-stabilizer}}
\IdxSub{idx133.en.witt}{Witt vectors and additive carries}{\pageref{v133:subsec-witt-carry}}
\IdxSee{Partial resolvent}{see}{\hyperlink{idx106.en.resolvent}{Normal element: Partial resolvents}}
\IdxSee{Period sums of roots of unity}{see}{\hyperlink{idx101.en.e654}{Hankel matrix: Period sums of roots of unity}}
\IdxHead{idx141.en.perioddual}{Period trace-dual basis}{\pageref{v141:prop-period-duality}}
\IdxHead{idx101.en.symmetry}{Permutation group of cyclic codes}{}
\IdxSub{idx101.en.e115}{Common permutation group}{\pageref{v73:def-common-permutation-group}}
\IdxSub{idx101.en.e116}{Common permutation group: arithmetic criterion}{\pageref{v73:thm-common-permutation-arithmetic}}
\IdxSub{idx101.en.e117}{Common permutation group of cyclic codes}{\pageref{v73:def-common-permutation-group}}
\IdxSub{idx101.en.e118}{Common permutation symmetries: difference-orbit criterion}{\pageref{v73:thm-common-permutation-orbits}}
\IdxSub{idx140.en.commonlift}{Persistence over unramified rings}{\pageref{v140:prop-common-symmetry}}
\IdxHead{idx101.en.permutation}{Permutation polynomial}{\pageref{def:permutation-polynomial}}
\IdxSub{idx101.en.e408}{Gold power map}{\pageref{v90:prop-gold}}
\IdxSub{idx101.en.e420}{Hermite's criterion}{\pageref{thm:hermite}}
\IdxSub{idx101.en.e448}{Involution}{\pageref{def:involution}}
\IdxSub{idx105.en.ring-lift}{Lifts to finite local rings}{\pageref{v105:prop-permutation-lifting}}
\IdxSub{idx101.en.e867}{Sparse nonlinear involutions}{\pageref{v68:ex-sparse-involution}}
\IdxHead{idx135.en.20}{Point-count reconstruction}{\pageref{revadic:thm-genus-counts}}
\IdxHead{idx101.en.polynomial}{Polynomial}{}
\IdxSub{idx123.en.pair}{Binary derivative pairing}{\pageref{v123:cor-binary-derivative-pair}}
\IdxSub{idx101.en.e670}{Division with remainder}{\pageref{lem:division}}
\IdxSub{idx101.en.e671}{Euler totient function}{\pageref{def:polynomial-euler}}
\IdxSub{idx101.en.e354}{Formal derivative}{\pageref{v85term:def-formal-derivative}}
\IdxSub{idx105.en.derivative-lift}{Formal derivatives and permutation lifts}{\pageref{v105:ex-F32-permutation-lifts}}
\IdxSub{idx101.en.e516}{Logarithmic derivative}{\pageref{v85:def-logarithmic-derivative}}
\IdxSub{idx101.en.e535}{Minimal polynomial}{\pageref{def:algebraic-minpoly}}
\IdxSub{idx101.en.e558}{Monic polynomial}{\pageref{def:monic-polynomial}}
\IdxSub{idx101.en.e586}{Newton identities}{\pageref{eq:root-power-newton}}
\IdxSub{idx123.en.odd}{Odd characteristic and the second derivative}{\pageref{v123:rem-odd-characteristic}}
\IdxSub{idx101.en.e680}{Power sum}{\pageref{v40:prop:root-power-composition}}
\IdxSub{idx123.en.prefix}{Power-sum prefix fibers}{\pageref{v123:thm-power-trace-pencil}}
\IdxSub{idx101.en.e732}{Quotient rings and adjoining a root}{\pageref{prop:adjoin}}
\IdxSub{idx101.en.e789}{Resultant}{\pageref{def:resultant-root-product}}
\IdxSub{idx101.en.e75}{Root bound}{\pageref{lem:root-bound}}
\IdxSub{idx101.en.e802}{Root bounds and nondegeneracy}{\pageref{sec:root-bound-rigidity}}
\IdxSub{idx101.en.e811}{Root power sums}{\pageref{v40:prop:root-power-composition}}
\IdxSub{idx101.en.e836}{Separability}{\pageref{def:separable-squarefree}}
\IdxSub{idx101.en.e840}{Separable polynomial}{\pageref{def:separable-squarefree}}
\IdxSub{idx101.en.e675}{Splitting into linear factors}{\pageref{sec:multivariate-linear-descent}}
\IdxSub{idx101.en.e883}{Square-free kernel}{\pageref{def:squarefree-kernel}}
\IdxSub{idx101.en.e884}{Square-free polynomial}{\pageref{def:separable-squarefree}}
\IdxSub{idx101.en.e908}{Support size}{\pageref{v79:def-defining-complexity}}
\IdxSub{idx101.en.e964}{Unique factorization}{\pageref{cor:unique-factorization}}
\IdxSee{Polynomial basis representation}{see}{\hyperlink{idx101.en.e669}{Field coordinates: Polynomial basis representation}}
\IdxHead{idx101.en.functions}{Polynomial function}{\pageref{thm:polynomial-functions}}
\IdxSub{idx105.en.representative-derivative}{Representatives and formal derivatives}{\pageref{v105:subsec-polynomial-lifts}}
\IdxSee{Polynomial representation of linear maps}{see}{\hyperlink{idx101.en.e674}{Linearized polynomial: Polynomial representation of linear maps}}
\IdxSee{Power-basis degree filtration}{see}{\hyperlink{idx101.en.e676}{Degree profile: Power-basis filtration}}
\IdxHead{idx101.en.residue}{Power-residue symbol}{\pageref{def:power-residue-symbol}}
\IdxSub{idx101.en.e603}{Norm expression}{\pageref{prop:residue-norm-resultant}}
\IdxSub{idx101.en.e673}{Polynomial reciprocity law}{\pageref{thm:polynomial-reciprocity}}
\IdxHead{idx140.en.precision}{Precision profile of representations}{\pageref{v140:eq-precision-profile}}
\IdxHead{idx140.en.preservation}{Preservation depth of a permutation}{\pageref{v140:prop-permutation-depth}}
\IdxHead{idx101.en.primitive}{Primitive element}{\pageref{def:primitive-element}}
\IdxSub{idx101.en.e695}{Constant-term norms}{\pageref{v82:ex-primitive-constant-term}}
\IdxSub{idx101.en.e172}{Counting primitive polynomials}{\pageref{cor:primitive-count}}
\IdxSub{idx101.en.e407}{Global uniqueness of optimal determining polynomials}{\pageref{v45:prop-global-optimal-uniqueness}}
\IdxSub{idx101.en.e687}{Indicator}{\pageref{prop:primitive-character-indicator}}
\IdxSub{idx101.en.e688}{Nonzero trace}{\pageref{v88:thm-primitive-nonzero-trace}}
\IdxSub{idx101.en.e623}{Optimal coefficient criterion for primitive quadratics}{\pageref{v45:thm-optimal-quadratic-criterion}}
\IdxSub{idx101.en.e624}{Optimal determining polynomial}{\pageref{v45:def-optimal-determining-polynomial}}
\IdxSub{idx101.en.e689}{Prescribed trace}{\pageref{thm:primitive-prescribed-trace}}
\IdxSub{idx101.en.e694}{Primitive polynomial}{\pageref{def:primitive-polynomial}}
\IdxSub{idx101.en.e696}{Primitive quadratic root-ratio criterion}{\pageref{v45:prop-quadratic-root-ratio}}
\IdxSub{idx101.en.e697}{Primitive-root order test by trace recurrence}{\pageref{v49:criterion-trace-order}}
\IdxSub{idx101.en.e699}{Primitivity criterion}{\pageref{prop:primitive-test-poly}}
\IdxSub{idx111.en.rprimitive}{$r$-primitive elements}{\pageref{v111:def-r-primitive}}
\IdxHead{idx101.en.pn}{Primitive normal basis}{}
\IdxSub{idx101.en.e171}{A counting example}{\pageref{v48:ex-primitive-normal-f16}}
\IdxSub{idx101.en.e55}{Arithmetic completion}{\pageref{v93:chap-primitive-normal-arithmetic}}
\IdxSub{idx101.en.e1003}{Basic bound: exact exceptional range}{\pageref{v93:prop-arithmetic-complete}}
\IdxSub{idx101.en.e691}{Complete proof}{\pageref{v93:thm-PNBT}}
\IdxSub{idx101.en.e174}{Criterion}{\pageref{v48:prop-primitive-normal-test}}
\IdxSub{idx101.en.e269}{e-free element}{\pageref{v93:def-e-free}}
\IdxSub{idx101.en.e690}{Existence theorem}{\pageref{thm:primitive-normal}}
\IdxSub{idx101.en.e324}{Fiber counts}{\pageref{v93:cor-exact-pn}}
\IdxSub{idx101.en.e393}{g-condition}{\pageref{v93:def-local-normality}}
\IdxSub{idx101.en.e463}{Joint character count}{\pageref{v93:prop-basic}}
\IdxSub{idx101.en.e507}{Local density sieve}{\pageref{v93:thm-sieve}}
\IdxSub{idx101.en.e561}{Monotone parameter families}{\pageref{v93:lem-family}}
\IdxSub{idx101.en.e685}{Primitive completely normal element}{\pageref{v88:cor-primitive-normal-special-degrees}}
\IdxSub{idx101.en.e692}{Primitive normal elements}{\pageref{def:primitive-normal-element}}
\IdxSub{idx101.en.e693}{Primitive normal elements in a translate}{\pageref{v82:rem-primitive-normal-translates}}
\IdxSub{idx101.en.e729}{Quotient generation and full degree}{\pageref{v93:lem-quotient-full-degree}}
\IdxSub{idx101.en.e677}{Small-prime power bounds}{\pageref{v93:lem-small-prime-bound}}
\IdxSub{idx101.en.e980}{Weighted trace}{\pageref{v93:def-weighted-trace}}
\IdxSub{idx101.en.e981}{Weighted trace fiber of normal elements}{\pageref{v93:thm-patch}}
\IdxSee{Primitive normal element}{see}{\hyperlink{idx101.en.e692}{Primitive normal basis: Primitive normal elements}}
\IdxSee{Primitive word}{see}{\hyperlink{idx101.en.e698}{Orbit: Primitive words}}
\IdxHead{idx101.en.profinite}{Profinite group}{\pageref{def:profinite-group}}
\IdxSub{idx101.en.e3}{Absolute Galois group of a finite field}{\pageref{thm:absolute-frobenius}}
\IdxSub{idx101.en.e701}{Procyclic group}{\pageref{def:procyclic-group}}
\IdxSub{idx101.en.e703}{Profinite completion}{\pageref{def:profinite-integers}}
\IdxSee{Projective line}{see}{\hyperlink{idx101.en.e711}{Projective space: Projective line}}
\IdxSee{Projective order}{see}{\hyperlink{idx101.en.e712}{Order of a polynomial: Projective order}}
\IdxSee{Projective plane}{see}{\hyperlink{idx101.en.e713}{Projective space: Projective plane}}
\IdxHead{idx101.en.projective}{Projective space}{}
\IdxSub{idx101.en.e425}{Homogeneous gluing}{\pageref{v84:prop-gluing}}
\IdxSub{idx101.en.e427}{Homogenization}{\pageref{def:polynomial-homogenization}}
\IdxSub{idx101.en.e430}{Hyperoval}{\pageref{v84:ex-hyperoval}}
\IdxSub{idx101.en.e663}{Point at infinity}{\pageref{def:projective-plane-points}}
\IdxSub{idx101.en.e706}{Projective arc}{\pageref{v84:def-projective-arc}, \pageref{v96:def-projective-arc}}
\IdxSub{idx101.en.e709}{Projective image of a multiplicative subgroup}{\pageref{v96:def-projective-arc}}
\IdxSub{idx101.en.e710}{Projective involution}{\pageref{v85:prop-pair-involution}}
\IdxSub{idx101.en.e711}{Projective line}{\pageref{v85:def-projective-line}}
\IdxSub{idx101.en.e713}{Projective plane}{\pageref{def:projective-plane-points}}
\IdxSub{idx101.en.e714}{Projective root form}{\pageref{v84:prop-projective-root}}
\IdxSub{idx101.en.e715}{Projective unit group}{\pageref{v85:def-projective-unit-group}}
\IdxHead{idx138.en.pure}{Purely periodic expansion}{\pageref{v138:thm-periodic-cosets}}
\IdxGroup{Q}
\IdxSee{\(q\)-order polynomial}{see}{\hyperlink{idx101.en.e371}{Annihilating polynomial: Frobenius annihilator polynomial}}
\IdxHead{idx101.en.quadratic}{Quadratic form}{\pageref{def:quadratic-form}}
\IdxSub{idx150.en.apn}{Almost perfect nonlinear mappings and difference kernels}{\pageref{v150:ex-apn-boundary}}
\IdxSub{idx101.en.e50}{Anisotropic binary quadratic form}{\pageref{v45:def-binary-quadratic-types}}
\IdxSub{idx101.en.e51}{Anisotropic form}{\pageref{def:quadratic-form}}
\IdxSub{idx101.en.e54}{Arf invariant}{\pageref{def:arf-invariant}}
\IdxSub{idx101.en.e72}{Binary quadratic form}{\pageref{v45:def-binary-quadratic-forms}}
\IdxSub{idx101.en.e78}{Canonical forms of quadratic forms}{\pageref{thm:quadratic-classification}}
\IdxSub{idx106.en.sesqui}{Diagonal recovery of sesquilinear forms}{\pageref{v106:prop-sesquilinear-recovery}}
\IdxSub{idx101.en.e333}{Field-valued quadratic map}{\pageref{v70:def-field-valued-quadratic}}
\IdxSub{idx101.en.e429}{Hyperbolic plane}{\pageref{def:orthogonal-quadratic-language}}
\IdxSub{idx101.en.e455}{Isometry}{\pageref{def:quadratic-isometry}}
\IdxSub{idx101.en.e460}{Isotropic vector}{\pageref{def:quadratic-form}}
\IdxSub{idx101.en.e640}{Orthogonal sum}{\pageref{def:orthogonal-quadratic-language}}
\IdxSub{idx101.en.e668}{Polar form}{\pageref{def:quadratic-form}}
\IdxSub{idx101.en.e723}{Quadratic phase function}{\pageref{sec:quadratic-fourier-support}}
\IdxSub{idx101.en.e874}{Split binary quadratic form}{\pageref{v45:def-binary-quadratic-types}}
\IdxSub{idx101.en.e880}{Square class}{\pageref{def:square-class}}
\IdxSub{idx101.en.e912}{Symplectic basis}{\pageref{def:orthogonal-quadratic-language}}
\IdxSub{idx101.en.e935}{Trace component}{\pageref{v70:def-field-valued-quadratic}}
\IdxSub{idx101.en.e946}{Trace quadratic forms and adjoint operators}{\pageref{prop:trace-quadratic-operator}}
\IdxSub{idx101.en.e971}{Value distribution of quadratic forms in characteristic two}{\pageref{thm:quadratic-even-count}}
\IdxSub{idx101.en.e972}{Value distribution of quadratic forms in odd characteristic}{\pageref{thm:quadratic-odd-count}}
\IdxSub{idx101.en.e982}{Weighted trace form}{\pageref{v90:subsec-weighted-trace}}
\IdxHead{idx137.en.perioddepth}{Quadratic periods, local valuations}{\pageref{v137:prop-quadratic-period-depth}}
\IdxHead{idx101.en.quadraticroot}{Quadratic root set}{}
\IdxSub{idx101.en.e720}{Quadratic factorization in additive coset families}{\pageref{v45:prop-quadratic-coset-count}}
\IdxSub{idx101.en.e724}{Quadratic quotient coordinates of root sets}{\pageref{v85:cor-rootset-quotient-coordinates}}
\IdxSub{idx101.en.e961}{Two-dimensional Frobenius eigenspaces of quadratic root sets}{\pageref{v45:prop-quadratic-planes}}
\IdxSee{Quasi-order of a polynomial}{see}{\hyperlink{idx101.en.e712}{Order of a polynomial: Projective order}}
\IdxGroup{R}
\IdxSee{Rabin's irreducibility criterion}{see}{\hyperlink{idx101.en.e734}{Factorization algorithm: Rabin criterion}}
\IdxSee{Rank distance}{see}{\hyperlink{idx101.en.e740}{Rank metric: Rank distance}}
\IdxHead{idx101.en.rankmetric}{Rank metric}{}
\IdxSub{idx105.en.column-normalization}{Column normalization}{\pageref{v105:ex-rank-column-normalization}}
\IdxSub{idx105.en.column-scaling}{Extension-field column scaling}{\pageref{v105:ex-rank-column-normalization}}
\IdxSub{idx101.en.e394}{Gabidulin code}{\pageref{v90:thm-rank-metric}}
\IdxSub{idx101.en.e528}{Maximum rank distance code}{\pageref{v99:def-mrd}}
\IdxSub{idx101.en.e551}{Minimum rank distance}{\pageref{v99:def-rank-distance}}
\IdxSub{idx101.en.e740}{Rank distance}{\pageref{v99:def-rank-distance}}
\IdxSub{idx101.en.e742}{Rank-metric Singleton bound}{\pageref{v99:par-rank-singleton}}
\IdxSub{idx103.en.rankexample}{Rank versus Hamming weight}{\pageref{v103:ex-rank-eight}}
\IdxSub{idx101.en.e741}{Rank weight}{\pageref{v99:def-rank-weight}}
\IdxSee{Rank weight}{see}{\hyperlink{idx101.en.e741}{Rank metric: Rank weight}}
\IdxHead{idx101.en.rational}{Rational representation}{}
\IdxSub{idx101.en.e156}{Constant direction and polynomial inner components}{\pageref{str:cor-polynomial-direction}}
\IdxSub{idx101.en.e175}{Criterion for rational representations to attain the degree lower bound}{\pageref{str:thm-rational-rank}}
\IdxSub{idx101.en.e214}{Degree of a rational function}{\pageref{str:def-rational-degree}}
\IdxSub{idx101.en.e552}{Minimum rational generator-representation degree}{\pageref{str:def-rational-generator-degree}}
\IdxSub{idx101.en.e749}{Rational representations attaining the degree lower bound: fractional linear classification}{\pageref{str:cor-rational-uniqueness}}
\IdxHead{idx122.en.guide}{Reading routes}{\pageref{v122:guide}}
\IdxSub{idx122.en.decoding}{Linear recurrences and codes}{\pageref{v122:route-decoding}}
\IdxSub{idx122.en.normal}{Normal bases and cyclic modules}{\pageref{v122:route-normal}}
\IdxSub{idx122.en.reciprocals}{Reciprocal families and normality obstructions}{\pageref{v122:route-reciprocals}}
\IdxSub{idx130.en.orbits}{Root exponents, phases, and coefficient lifting}{\pageref{v125:route-orbits}}
\IdxSub{idx122.en.fibers}{Root fibers and recursive irreducible constructions}{\pageref{v122:route-fibers}}
\IdxSub{idx122.en.subfields}{Subfield representations and relative trace recovery}{\pageref{v122:route-subfields}}
\IdxSub{idx122.en.fourier}{Trace duality and Fourier analysis}{\pageref{v122:route-fourier}}
\IdxHead{idx101.en.reciprocal}{Reciprocal polynomial}{\pageref{def:reciprocal-polynomial}}
\IdxSub{idx101.en.e148}{Conjugate reciprocity and norm cancellation}{\pageref{v97:ex-norm-cancellation}}
\IdxSub{idx145.en.inverse}{Full-degree inverse pairs and irreducible counts}{\pageref{v145:prop-inverse-count}}
\IdxSub{idx101.en.e461}{Iterated R-transform}{\pageref{v48:ex-r-transform-tower}}
\IdxSub{idx101.en.e718}{Q-transform}{\pageref{thm:reciprocal-transform}}
\IdxSub{idx101.en.e733}{R-transform}{\pageref{def:r-transform}}
\IdxSub{idx101.en.e826}{Self-reciprocal polynomial}{\pageref{def:self-reciprocal-polynomial}}
\IdxSub{idx101.en.e827}{Self-reciprocal transform}{\pageref{thm:reciprocal-transform}}
\IdxHead{idx116.en.inverse-head}{Reciprocal rank}{\pageref{v116:def-inverse-rank}}
\IdxSub{idx117.en.loss}{Exact loss under evaluation}{\pageref{v117:prop-specialization-loss}}
\IdxSub{idx116.en.inverse-head1}{Relations across subspace sections}{\pageref{v116:def-inverse-relations}}
\IdxSub{idx116.en.inverse-head0}{Relations preserved by evaluation}{\pageref{v116:prop-rational-specialization}}
\IdxSee{Recurrence polynomial}{see}{\hyperlink{idx101.en.e764}{Linear recurrence: Recurrence polynomial}}
\IdxHead{idx106.en.redei}{R\'edei rational function}{\pageref{v106:ex-redei-rational}}
\IdxHead{idx101.en.rs}{Reed--Solomon code}{\pageref{v84:prop-rs}}
\IdxSub{idx105.en.moore-normalization}{Column normalization of a Moore matrix}{\pageref{v105:ex-rank-column-normalization}}
\IdxSub{idx101.en.e265}{Duality of evaluation spaces}{\pageref{v84:prop-rs-dual}}
\IdxSub{idx101.en.e312}{Extended Reed--Solomon code}{\pageref{v84:prop-rs}}
\IdxSub{idx101.en.e391}{Full-length Reed--Solomon code}{\pageref{v85term:def-full-length-rs}}
\IdxSub{idx150.en.design}{Three-design obstruction for minimum-weight codewords}{\pageref{v150:thm-ers-design}}
\IdxSee{Relative Frobenius map of a finite field extension}{see}{\hyperlink{idx101.en.e775}{Frobenius map: Relative map}}
\IdxSee{Relative trace}{see}{\hyperlink{idx101.en.e778}{Trace: Relative trace}}
\IdxHead{idx101.en.tracerecovery}{Relative trace recovery}{}
\IdxSub{idx101.en.e264}{Duality between degree profiles and relative traces}{\pageref{str:thm-profile-trace-duality}}
\IdxSub{idx101.en.e548}{Minimum generator-representation degrees from trace prefixes}{\pageref{v70:cor-trace-prefix-generators}}
\IdxSub{idx150.en.threetraces}{Three traces and a specified linear target}{\pageref{v150:ex-three-traces}}
\IdxSub{idx101.en.e948}{Trace-rank criterion for polynomial composition}{\pageref{str:cor-trace-composition}}
\IdxSub{idx101.en.e949}{Trace-rank function}{\pageref{str:def-trace-rank}}
\IdxSub{idx101.en.e950}{Trace-rank jump positions}{\pageref{str:def-trace-jumps}}
\IdxSub{idx101.en.e979}{Weighted relative trace sequence}{\pageref{str:def-relative-trace-sequence}}
\IdxHead{idx135.en.21}{Residue field construction}{\pageref{revadic:thm-main}}
\IdxSee{Resultant}{see}{\hyperlink{idx101.en.e789}{Polynomial: Resultant}}
\IdxHead{idx101.en.return}{Return map}{\pageref{def:fiber-return-map}}
\IdxSub{idx101.en.e30}{Affine return for cyclotomic lifting}{\pageref{v44:thm-cyclotomic-affine-return}}
\IdxSub{idx101.en.e31}{Affine return on recurrence quotients}{\pageref{prop:recurrence-fiber-return}}
\IdxSub{idx101.en.e32}{Affine return on trace fibers}{\pageref{v70:ex-additive-return}}
\IdxSub{idx129.en.return}{Base extension and return cycles}{\pageref{v129:cor-return-base-change}}
\IdxSub{idx101.en.e124}{Comparison of dilation and translation returns}{\pageref{exm:ternary-affine-return}}
\IdxSub{idx101.en.e141}{Computing the lifted factorization type from a compressed factor}{\pageref{cor:explicit-compressed-lifting}}
\IdxSub{idx101.en.e142}{Computing the return order from a compressed value}{\pageref{prop:explicit-affine-return}}
\IdxSub{idx129.en.euler}{Euler products for return maps}{\pageref{v129:eq-return-euler}}
\IdxSub{idx101.en.e320}{Factorization from return maps on additive fibers}{\pageref{ffs:thm-additive-return}}
\IdxSub{idx101.en.e794}{Full affine fibers}{\pageref{ffs:thm-full-return}}
\IdxSub{idx101.en.e488}{Lifting compressed factors and fiber decomposition}{\pageref{cor:compressed-fiber-lifting}}
\IdxSub{idx101.en.e496}{Linear fixed-point equations and factor counts}{\pageref{ffs:prop-return-counts}}
\IdxSub{idx101.en.e795}{Multiplicative fibers}{\pageref{prop:multiplicative-permutation-return}}
\IdxSub{idx129.en.images}{Orbit images and multiplicities}{\pageref{v129:eq-orbit-image-product}}
\IdxSub{idx101.en.e632}{Order of affine return and period reduction}{\pageref{prop:affine-return-order}}
\IdxSub{idx102.en.par-two-fiber-returns}{Quartic compositions with two return types}{\pageref{v102:par-two-fiber-returns}}
\IdxSub{idx115.en.return}{Return fixed points and factor-degree counts}{\pageref{v115:cor-return-counts}}
\IdxSub{idx101.en.e796}{Return maps for power compression and affine compression}{\pageref{sec:power-affine-comparison}}
\IdxSub{idx101.en.e792}{Root-fiber interpretation}{\pageref{def:fiber-return-map}}
\IdxSub{idx101.en.e936}{Trace criterion for return when the norm is one}{\pageref{cor:affine-return-trace}}
\IdxSee{Return norm}{see}{\hyperlink{idx148.en.all}{Skew polynomial: All return-norm fibers}}
\IdxHead{idx101.en.rootblock}{Root block}{}
\IdxSub{idx101.en.e104}{Coefficient criterion for blocks of conjugate roots}{\pageref{sfr:prop-block-coefficients}}
\IdxSub{idx101.en.e642}{Orthogonality reversing the degree filtration}{\pageref{str:prop-filter-orthogonal}}
\IdxSub{idx101.en.e758}{Recovery of representation structures}{\pageref{str:cor-unified-recovery}}
\IdxSub{idx101.en.e777}{Relative root-block polynomial}{\pageref{str:def-block-polynomial}}
\IdxSub{idx101.en.e798}{Root-block coefficient rank}{\pageref{str:def-coefficient-rank}}
\IdxSub{idx101.en.e799}{Root-block coordinate-polynomial space}{\pageref{str:def-coordinate-space}}
\IdxSee{Root fiber}{see}{\hyperlink{idx101.en.e805}{Fiber: Root fiber}}
\IdxSee{Root-field spectrum}{see}{\hyperlink{idx101.en.e806}{Fiber: Root-field spectrum}}
\IdxHead{idx101.en.unity}{Root of unity}{\pageref{def:roots-of-unity}}
\IdxSub{idx101.en.e130}{Compatible primitive roots of unity}{\pageref{v44:lem-primitive-root-lift}}
\IdxSub{idx101.en.e351}{Fixed points of roots of unity and factorization type}{\pageref{v40:ex-unit-fixed-counts}}
\IdxHead{idx101.en.rootset}{Root set}{}
\IdxSub{idx101.en.e127}{Comparison of two root sets of size 64}{\pageref{exm:two-sixtyfour-rootsets}}
\IdxSub{idx101.en.e147}{Conjugate-pair root extension}{\pageref{v83:prop-conjugate-pair-extension}}
\IdxSub{idx101.en.e221}{Derivative construction from multiplicative cosets}{\pageref{v83:prop-two-mult-cosets}}
\IdxSub{idx101.en.e222}{Derivative fiber of a polynomial}{\pageref{v83:prop-derivative-fibers}}
\IdxSub{idx116.en.fourpoints}{Four-point completion and a quadratic trace obstruction}{\pageref{v116:prop-four-completion}}
\IdxSub{idx116.en.completion}{Minimal constant-derivative completion}{\pageref{v116:def-completion}}
\IdxSub{idx101.en.e584}{Multiset}{\pageref{v58:def-multiset}}
\IdxSub{idx101.en.e810}{Root polynomial of a finite subfield}{\pageref{prop:subfield-root-polynomial}}
\IdxSub{idx101.en.e879}{Splitting test for a derivative fiber}{\pageref{v83:eq-fiber-splitting-test}}
\IdxSee{Root space}{see}{\hyperlink{idx101.en.e814}{Linearized polynomial: Root space}}
\IdxSee{$r$-primitive element}{see}{\hyperlink{idx111.en.rprimitive}{Primitive element: $r$-primitive elements}}
\IdxSee{R-transform}{see}{\hyperlink{idx101.en.e733}{Reciprocal polynomial: R-transform}}
\IdxGroup{S}
\IdxSee{Self-dual basis}{see}{\hyperlink{idx101.en.e824}{Dual basis: Self-dual basis}}
\IdxSee{Self-dual normal basis}{see}{\hyperlink{idx101.en.e825}{Normal basis: Self-dual normal basis}}
\IdxSee{Self-reciprocal polynomial}{see}{\hyperlink{idx101.en.e826}{Reciprocal polynomial: Self-reciprocal polynomial}}
\IdxHead{idx101.en.semiinvariant}{Semi-invariant}{\pageref{fas:def-semi-invariant}}
\IdxSub{idx101.en.e64}{Basis of a space of semi-invariants}{\pageref{fas:prop-semi-invariant-basis}}
\IdxSee{Semilinear map}{see}{\hyperlink{idx101.en.e832}{Linear map: Semilinear map}}
\IdxSee{Sesquilinear form}{see}{\hyperlink{idx106.en.sesqui}{Quadratic form: Diagonal recovery}}
\IdxSee{Set of evaluation points}{see}{\hyperlink{idx101.en.e295}{Evaluation: Sets}}
\IdxSee{Shortened code}{see}{\hyperlink{idx101.en.e848}{Linear code: Shortened code}}
\IdxSee{Sidon set}{see}{\hyperlink{idx101.en.e850}{Difference rigidity: Sidon set}}
\IdxSee{Sidon set in characteristic two}{see}{\hyperlink{idx101.en.e850}{Difference rigidity: Sidon set}}
\IdxSee{Simple-root Hensel lifting}{see}{\hyperlink{idx101.en.e858}{Galois ring: Simple-root Hensel lifting}}
\IdxSee{Singleton bound}{see}{\hyperlink{idx101.en.e860}{MDS code: Singleton bound}}
\IdxHead{idx145.en.skew}{Skew polynomial}{\pageref{v145:sec-skew-roots}}
\IdxSub{idx148.en.all}{All return-norm fibers}{\pageref{v148:thm-all-fibers}}
\IdxSub{idx148.en.commutative}{Commutative-norm factor allocation}{\pageref{v148:prop-commutative-norm}}
\IdxSub{idx148.en.orders}{Complete return-order spectrum}{\pageref{v148:thm-even-order}}
\IdxSub{idx145.en.rightcount}{Counting right factors}{\pageref{v145:ex-right-factors}}
\IdxSub{idx148.en.cubeclasses}{Cube-class invariance of primitive slices}{\pageref{v148:prop-cubic-classes}}
\IdxSub{idx148.en.cubics}{Cubic slices with a fixed middle coefficient}{\pageref{v148:thm-cubic-value-set}}
\IdxSub{idx146.en.matching}{Cyclic matching expansion}{\pageref{v146:prop-return-matching}}
\IdxSub{idx148.en.even}{Even-degree fibers and the three-value formula}{\pageref{v148:thm-all-even-fibers}}
\IdxSub{idx146.en.order}{Exact return-order counts}{\pageref{v146:cor-return-order}}
\IdxSub{idx149.en.uniformity}{Factor phases and uniformity criterion}{\pageref{v149:prop-uniformity}}
\IdxSub{idx145.en.irreducibility}{Irreducibility and invariant subspaces}{\pageref{v145:thm-skew-roots}}
\IdxSub{idx148.en.local}{Local partitions and fiber types}{\pageref{v148:prop-local-types}}
\IdxSub{idx146.en.higher}{Nonuniform quartic return fibers}{\pageref{v146:ex-quartic-counterexample}}
\IdxSub{idx146.en.normal}{Normalized cyclic vector}{\pageref{v146:lem-semilinear-normal}}
\IdxSub{idx148.en.primitive}{Norm-primitivity}{\pageref{v148:thm-ordinary-return}}
\IdxSub{idx148.en.transfer}{Norm transfer for equally spaced powers}{\pageref{v148:prop-sparse-transfer}}
\IdxSub{idx151.en.optimal-skew}{Optimal coefficient-criterion pullback}{\pageref{v148:thm-optimal-pullback}}
\IdxSub{idx148.en.strata}{Optimal degrees and parameter strata}{\pageref{v148:prop-optimal-strata}}
\IdxSub{idx148.en.ordinary}{Ordinary and projective factorization patterns}{\pageref{v148:thm-ordinary-return}}
\IdxSub{idx146.en.quartic}{Quartic trace--norm distribution}{\pageref{v146:cor-quartic-distribution}}
\IdxSub{idx148.en.relative}{Relative norms and common subfields}{\pageref{v148:thm-relative-fibers}}
\IdxSub{idx146.en.repeated}{Repeated-root return fibers}{\pageref{v146:cor-repeated-fibers}}
\IdxSub{idx145.en.return}{Return characteristic polynomial}{\pageref{v145:prop-return-matrix}}
\IdxSub{idx146.en.spectrum}{Return spectrum and fixed field}{\pageref{v146:def-return-spectrum}}
\IdxSub{idx146.en.scaling}{Semilinear scaling of the constant term}{\pageref{v146:prop-return-scaling}}
\IdxSub{idx148.en.obstruction}{Spacing obstruction to norm-primitivity}{\pageref{v148:cor-sparse-order}}
\IdxSub{idx148.en.sparse}{Sparse-family counts and existence}{\pageref{v148:thm-sparse-coprime}}
\IdxSub{idx145.en.sparse}{Sparse irreducible counts}{\pageref{v145:thm-sparse-count}}
\IdxSub{idx146.en.split}{Split distinct-root return fibers}{\pageref{v146:thm-split-fibers}}
\IdxSub{idx145.en.fibers}{Uniform quadratic return fibers}{\pageref{v145:thm-return-uniform}}
\IdxHead{idx140.en.smith}{Smith factors and exponents}{\pageref{v140:thm-profile-smith}}
\IdxSub{idx144.en.invariants}{Coefficient fixed algebras and rank reduction}{\pageref{v144:thm-invariant-smith}}
\IdxSub{idx143.en.valuations}{Local valuation multisets}{\pageref{v143:thm-valuation-smith}}
\IdxHead{idx101.en.sparse}{Sparse multiple}{}
\IdxSub{idx101.en.e15}{Adjacent coefficient ratios}{\pageref{v97:prop-linear-cancellation}}
\IdxSub{idx101.en.e102}{Coefficient cancellation by a linear multiplier}{\pageref{v97:prop-linear-cancellation}}
\IdxSub{idx101.en.e543}{Minimum degree}{\pageref{v97:def-sparse-degree}}
\IdxSub{idx101.en.e560}{Monomial substitution in sparse multiples}{\pageref{v97:ex-sparse-base-change}}
\IdxSub{idx101.en.e816}{Scalar extension of sparse relations}{\pageref{v97:ex-sparse-base-change}}
\IdxSub{idx101.en.e920}{Three-dimensional projection for sparse relations}{\pageref{v97:ex-projection-sparse}}
\IdxSub{idx101.en.e966}{Uniqueness of degree-minimal sparse multiples}{\pageref{v97:ex-sparse-uniqueness}}
\IdxLine{1.2em}{\emph{See also} \hyperlink{idx101.en.sparsepoly}{Sparse polynomial}; \hyperlink{idx101.en.order}{Order of a polynomial}}{}
\IdxSee{Sparse nonlinear involutions}{see}{\hyperlink{idx101.en.e867}{Permutation polynomial: Sparse nonlinear involutions}}
\IdxHead{idx101.en.sparsepoly}{Sparse polynomial}{}
\IdxSub{idx125.en.zeros}{Equal stabilizers and additional coefficient zeros}{\pageref{v125:ex-extra-cancellation}}
\IdxSub{idx101.en.e303}{Exceptional seven-power family}{\pageref{v49:thm-binary-three-sparse}}
\IdxSub{idx101.en.e815}{Root-squaring obstruction to three-sparsity}{\pageref{v49:ex-square-four-terms}}
\IdxSub{idx101.en.e870}{Sparse trinomial families}{\pageref{prop:sparse-trinomial-family}}
\IdxSub{idx101.en.e885}{Square-root estimate for sparse families}{\pageref{prop:sparse-square-root}}
\IdxSub{idx101.en.e0}{Three-sparse factorization of \(X^n-1\)}{\pageref{v49:def-three-sparse}}
\IdxSub{idx101.en.e923}{Three-sparsity criterion in characteristic two}{\pageref{v49:thm-binary-three-sparse}}
\IdxSub{idx101.en.e924}{Three-sparsity criterion in odd characteristic}{\pageref{v49:thm-odd-three-sparse}}
\hypertarget{idx101.en.e872}{}\IdxSee{Spectral support}{see}{\hyperlink{idx101.en.fouriercode}{Fourier support of a cyclic code}}
\IdxSee{Splitting field}{see}{\hyperlink{idx101.en.e875}{Field extension: Splitting field}}
\IdxSee{Square-free polynomial}{see}{\hyperlink{idx101.en.e884}{Polynomial: Square-free polynomial}}
\IdxSee{Stable polynomial}{see}{\hyperlink{idx101.en.e890}{Dynamical irreducibility: Stable polynomial}}
\IdxHead{idx101.en.subspace}{Subspace}{}
\IdxSub{idx101.en.e37}{Affine subspace}{\pageref{def:affine-space}}
\IdxSub{idx101.en.e197}{Cyclically covering subspace}{\pageref{v82:def-cyclically-covering}}
\IdxSub{idx105.en.sum-remainders}{Direct sums and independent remainders}{\pageref{v105:prop-linearized-remainder-compatibility}}
\IdxSub{idx150.en.pgl}{Fractional linear equivalence of power spaces}{\pageref{v150:eq-power-pgl}}
\IdxSub{idx150.en.powerdistribution}{Intersection distributions of power spaces}{\pageref{v150:cor-power-distribution}}
\IdxSub{idx101.en.e481}{Lattice}{\pageref{sec:invariant-duality}}
\IdxSub{idx101.en.e482}{Lattice of invariant subspaces and duality}{\pageref{thm:invariant-duality}}
\IdxSub{idx101.en.e525}{Maximal field of scalars}{\pageref{v70:def-scalar-field}, \pageref{v147:thm-carry-directions}}
\IdxSub{idx101.en.e572}{Multiplicative intersection dimension}{\pageref{v70:def-quotient-multiplication}}
\IdxSub{idx101.en.e899}{Polynomial}{\pageref{thm:subspace-polynomial}}
\IdxSub{idx101.en.e702}{Product of subspaces}{\pageref{v70:def-subspace-product}}
\IdxSub{idx101.en.e730}{Quotient multiplication map}{\pageref{v70:def-quotient-multiplication}}
\IdxSub{idx116.en.inverse}{Reciprocal ranks and rational function models}{\pageref{v116:ex-inverse-models}}
\IdxSub{idx116.en.inversesum}{Reciprocal sums and relation spaces}{\pageref{v116:ex-inverse-sum}}
\IdxSub{idx150.en.powers}{Scalar intersections of power spaces}{\pageref{v150:prop-power-intersection}}
\IdxSub{idx101.en.e834}{Semisimple product formula for invariant subspaces}{\pageref{v70:thm-semisimple-product}}
\IdxSee{Subspace polynomial}{see}{\hyperlink{idx101.en.e899}{Subspace: Polynomial}}
\IdxHead{idx101.en.support}{Support of a polynomial}{}
\IdxSub{idx101.en.e119}{Common support bounds from affine symmetries}{\pageref{v67:cor-affine-support-bounds}}
\IdxSub{idx101.en.e120}{Common support set}{\pageref{fas:def-common-support}}
\IdxSub{idx101.en.e559}{Monomial composition and root power sums of orbit unions}{\pageref{v40:prop:root-power-composition}}
\IdxSub{idx101.en.e1002}{Singleton bound for compressed components}{\pageref{v67:cor-core-distance}}
\IdxSub{idx101.en.e903}{Support period}{\pageref{v44:def-support-period}, \pageref{v97:prop-support-splitting}}
\IdxSub{idx101.en.e904}{Support-primitive core}{\pageref{v44:def-support-core}}
\IdxSub{idx101.en.e905}{Support-primitive decomposition}{\pageref{v44:def-support-decomposition}}
\IdxSub{idx101.en.e906}{Support-primitive polynomial}{\pageref{v44:def-support-primitive}}
\IdxSee{Support period}{see}{\hyperlink{idx101.en.e903}{Support of a polynomial: Support period}}
\IdxSee{Support-primitive core}{see}{\hyperlink{idx101.en.e904}{Support of a polynomial: Support-primitive core}}
\IdxSee{Support-primitive decomposition}{see}{\hyperlink{idx101.en.e905}{Support of a polynomial: Support-primitive decomposition}}
\IdxSee{Support-primitive polynomial}{see}{\hyperlink{idx101.en.e906}{Support of a polynomial: Support-primitive polynomial}}
\IdxGroup{T}
\IdxHead{idx135.en.22}{Teichmüller projection}{\pageref{revadic:prop-teich-projection}}
\IdxSee{Teichmüller representative}{see}{\hyperlink{idx101.en.e915}{Galois ring: Teichmüller representative}}
\IdxSee{Teichmüller set}{see}{\hyperlink{idx101.en.e915}{Galois ring: Teichmüller representative}}
\IdxHead{idx101.en.tensor}{Tensor product}{}
\IdxSub{idx101.en.e328}{Field decomposition}{\pageref{thm:tensor-fields}}
\IdxSub{idx103.en.tensorfour}{Four-element tensor product and idempotents}{\pageref{v103:ex-tensor-four}}
\IdxSub{idx101.en.e917}{Tensor-product isomorphism by conjugate evaluation}{\pageref{eq:conjugate-tensor-map}}
\IdxHead{idx106.en.tower}{Tower matrices, compatible product bases}{\pageref{v106:par-tower-matrices}}
\IdxHead{idx101.en.trace}{Trace}{\pageref{def:trace-norm}}
\IdxSub{idx101.en.e4}{Absolute trace}{\pageref{def:trace-norm}}
\IdxSub{idx101.en.e932}{Autocorrelation}{\pageref{v88:eq-conjugate-correlation}}
\IdxSub{idx101.en.e934}{Compatibility condition}{\pageref{v58:def-trace-compatibility}}
\IdxSub{idx101.en.e282}{Equivariant splitting}{\pageref{prop:trace-splitting}}
\IdxSub{idx101.en.e939}{Evaluation, a nonzero function with zero trace}{\pageref{v68:ex-trace-x5}}
\IdxSub{idx101.en.e940}{Evaluation, coefficient parameters}{\pageref{v68:prop-trace-evaluation-parameters}}
\IdxSub{idx101.en.e321}{Factorization of trace fibers}{\pageref{prop:trace-fiber-factorization}}
\IdxSub{idx101.en.e947}{Rank and first-order affine symmetries}{\pageref{v67:prop-first-order-affine}}
\IdxSub{idx101.en.e778}{Relative trace}{\pageref{def:trace-norm}}
\IdxSub{idx101.en.e910}{Surjectivity (also for norms)}{\pageref{thm:trace-norm-surjective}}
\IdxSub{idx105.en.trace-minimal-values}{Trace kernels as minimal value sets}{\pageref{v105:ex-F81-norm-values}}
\IdxSub{idx101.en.e955}{Transitivity (also for norms)}{\pageref{prop:trace-tower}}
\IdxSub{idx123.en.trace-prefix}{Unweighted power-trace identification}{\pageref{v123:sec-power-trace-fibers}}
\IdxLine{1.2em}{\emph{See also} \hyperlink{idx101.en.pairing}{Trace pairing}; \hyperlink{idx101.en.tracerecovery}{Relative trace recovery}}{}
\IdxHead{idx101.en.tracecode}{Trace code}{}
\IdxSub{idx104.en.tracefive}{A length-five MDS example over the four-element field}{\pageref{v104:ex-trace-five}}
\IdxSub{idx101.en.e501}{Linear-relation map}{\pageref{v96:eq-relation-map}}
\IdxSub{idx142.en.relativeblocks}{Relative-trace block coordinates}{\pageref{v142:eq-relative-trace-blocks}}
\IdxSub{idx142.en.restriction}{Restriction to subgroup coordinates}{\pageref{v142:prop-subfield-distance}}
\IdxHead{idx111.en.trace-distribution}{Trace distribution}{\pageref{v111:eq-trace-count}}
\IdxSub{idx115.en.empty}{Empty trace fibers and second moments}{\pageref{v115:prop-empty-cosets}}
\IdxSub{idx111.en.trace-fourier}{Fourier quotient and uniformity}{\pageref{v111:cor-trace-fourier}}
\IdxSub{idx111.en.trace-cosets}{Multiplicative coset symmetry}{\pageref{v111:thm-order-trace-symmetry}}
\IdxSub{idx112.en.trace-fourier}{Squared deviation from uniformity}{\pageref{v112:eq-fourier-nonuniformity}}
\IdxSub{idx111.en.trace-uniform}{Uniform nonzero traces of primitive k-normal elements}{\pageref{v111:cor-primitive-trace-uniform}}
\IdxHead{idx101.en.pairing}{Trace pairing}{\pageref{def:trace-pairing}}
\IdxSub{idx101.en.e41}{Algebra trace}{\pageref{def:algebra-trace}}
\IdxSub{idx101.en.e42}{Coefficient pairing and derivative}{\pageref{thm:algebra-trace-coefficient}}
\IdxSub{idx101.en.e250}{Discriminant}{\pageref{def:discriminant}}
\IdxSub{idx105.en.discriminant-parity}{Discriminants and factorization parity}{\pageref{v105:prop-discriminant-parity}}
\IdxSub{idx101.en.e931}{Eight-dimensional quotient: trace and periodic part}{\pageref{exm:mixed-algebra-trace}}
\IdxSub{idx101.en.e409}{Gram matrix}{\pageref{def:gram-matrix}}
\IdxSub{idx101.en.e509}{Local multiplicities and multiplication determinants}{\pageref{prop:local-multiplicity-determinant}}
\IdxSub{idx101.en.e647}{Nondegeneracy}{\pageref{def:pairing-radical}}
\IdxSub{idx101.en.e945}{On an operator image}{\pageref{v68:cor-image-trace-radical}}
\IdxSub{idx101.en.e638}{Orthogonal complement}{\pageref{def:trace-orthogonal-complement}}
\IdxSub{idx105.en.polynomial-discriminant}{Polynomial discriminant}{\pageref{v105:def-polynomial-discriminant}}
\IdxSub{idx101.en.e735}{Radical of a pairing}{\pageref{def:pairing-radical}}
\IdxSub{idx101.en.e736}{Radical of a restricted trace pairing}{\pageref{v68:cor-image-trace-radical}}
\IdxSub{idx101.en.e737}{Radical of the trace pairing and the derivative}{\pageref{ex:wu7-trace-radical}}
\IdxSub{idx101.en.e837}{Separability and the algebra trace pairing}{\pageref{prop:trace-pairing-separability}}
\IdxSub{idx105.en.stickelberger}{Stickelberger's discriminant formula}{\pageref{v105:prop-discriminant-parity}}
\IdxSub{idx101.en.e983}{Weighted trace pairing}{\pageref{str:def-weighted-pairing}}
\IdxSee{Translation quotient}{see}{\hyperlink{idx101.en.e957}{Translation stabilizer: Translation quotient}}
\IdxHead{idx101.en.translation}{Translation stabilizer}{}
\IdxSub{idx101.en.e13}{Additive symmetry in prime degree}{\pageref{fas:thm-prime-characteristic}}
\IdxSub{idx109.en.additive}{Additive translation space}{\pageref{v109:prop-additive-directions}}
\IdxSub{idx101.en.e524}{Maximal additive inner component}{\pageref{ffs:cor-maximal-additive}}
\IdxSub{idx101.en.e959}{Polynomial}{\pageref{ffs:def-translation-polynomial}}
\IdxSub{idx101.en.e756}{Recovery from outer centers}{\pageref{v70:prop-outer-center-recovery}}
\IdxSub{idx101.en.e762}{Recovery of the translation stabilizer from coefficients}{\pageref{ffs:thm-psi}}
\IdxSub{idx101.en.e12}{Root multisets}{\pageref{ffs:def-stabilizers}}
\IdxSub{idx101.en.e957}{Translation quotient}{\pageref{fas:def-translation-quotient}}
\IdxSub{idx140.en.translationdigits}{Two-digit sets}{\pageref{v140:prop-digit-stabilizer}}
\IdxHead{idx135.en.23}{Transported addition}{\pageref{revadic:eq-transport-add}}
\IdxSee{Type I optimal normal basis}{see}{\hyperlink{idx101.en.e962}{Normal basis: Type I optimal basis}}
\IdxGroup{U}
\IdxHead{idx135.en.24}{Uniformizer}{\pageref{v135:sec-foundations}}
\IdxHead{idx135.en.25}{Unit reciprocal roots}{\pageref{revadic:cor-unit-roots}}
\IdxHead{idx135.en.0}{Unramified construction reading route}{\pageref{v132:route-adic}}
\IdxHead{idx135.en.26}{Unramified cyclotomic extension}{\pageref{revadic:thm-main}}
\IdxHead{idx141.en.phase}{Unramified phase algebra}{\pageref{v141:eq-phase-algebra}}
\IdxGroup{V}
\IdxHead{idx135.en.27}{Valuation profile of an orbit}{\pageref{revadic:eq-profile}}
\IdxSub{idx144.en.compression}{Full fibers and coefficient compression}{\pageref{v144:prop-profile-compression}}
\IdxSub{idx143.en.minors}{Integer minors and depth multiplicities}{\pageref{v143:thm-valuation-smith}}
\IdxSub{idx144.en.precision}{Second-largest distinct valuation and phase labels}{\pageref{v144:prop-label-precision}}
\IdxHead{idx101.en.values}{Value set}{\pageref{v58:def-value-set}}
\IdxSub{idx101.en.e111}{Collision count}{\pageref{v58:def-collision-number}}
\IdxSub{idx101.en.e541}{Minimal value set polynomial}{\pageref{v58:def-minimal-value-set}}
\IdxSub{idx105.en.norm-trace-values}{Norm fibers and trace kernels}{\pageref{v105:prop-norm-trace-values}}
\IdxSub{idx101.en.e678}{Power image}{\pageref{v58:def-power-image}}
\IdxSub{idx111.en.power-boundary}{Power-value rigidity and the square-root boundary}{\pageref{v111:ex-low-degree-power}}
\IdxSub{idx101.en.e726}{Quartic value sets over extension fields}{\pageref{v68:ex-affine-values-small}}
\IdxSub{idx101.en.e767}{Rédei root polynomial}{\pageref{v90:subsec-redei}}
\IdxSub{idx109.en.mvsp}{Structural classification and five values over F729}{\pageref{v109:thm-mvsp-classification}}
\IdxSub{idx101.en.e973}{Value distributions and intersection fields}{\pageref{v68:prop-affine-values}}
\IdxGroup{W}
\IdxHead{idx101.en.weighted}{Weighted Frobenius operator}{\pageref{wfr:def-weighted}}
\IdxSub{idx101.en.e332}{Field generation by a norm and cyclic vectors}{\pageref{wfr:cor-cyclic-realization}}
\IdxSub{idx101.en.e598}{Norm and cyclic module decomposition}{\pageref{wfr:thm-norm-decomposition}}
\IdxSub{idx101.en.e855}{Similarity classification of weighted Frobenius operators}{\pageref{wfr:cor-similarity}}
\IdxHead{idx101.en.weil}{Weil bound}{}
\IdxSub{idx112.en.weil-power}{Distinct-root refinement}{\pageref{v111:ex-low-degree-power}}
\IdxSub{idx101.en.e984}{For additive character sums}{\pageref{v91:thm-additive-weil}}
\IdxSub{idx101.en.e985}{For multiplicative character sums}{\pageref{v91:thm-multiplicative-weil}}
\IdxSub{idx111.en.weil-power}{Low-degree power-value rigidity}{\pageref{v111:ex-low-degree-power}}
\IdxSub{idx101.en.e681}{Power sums and maximal modulus}{\pageref{v91:lem-power-sum-radius}}
\IdxSub{idx101.en.e894}{Stepanov's method}{\pageref{v91:sec-weil-bounds}}
\IdxSee{Weil bound for additive character sums}{see}{\hyperlink{idx101.en.e984}{Weil bound: For additive character sums}}
\IdxSee{Weil bound for multiplicative character sums}{see}{\hyperlink{idx101.en.e985}{Weil bound: For multiplicative character sums}}
\IdxHead{idx136.en.witt-coordinates}{Witt coordinates and Teichmüller digits}{\pageref{v136:eq-witt-coordinates}}
\IdxGroup{Z}
\IdxSee{Zech logarithm}{see}{\hyperlink{idx101.en.e988}{Field coordinates: Zech logarithm}}
\IdxSee{Zeta function}{see}{\hyperlink{idx101.en.e990}{Curve point counts: Zeta function}}
\end{multicols}
\endgroup

\end{document}